\documentclass[11pt,a4paper,reqno]{amsart}
\usepackage[no-math]{fontspec}
\usepackage{xeCJK}
\setCJKmonofont{FandolFang-Regular.otf}
\usepackage{graphicx}
\newsavebox{\casebox}
\newcommand{\CaseLine}[1]{\par\noindent\hfill\begingroup\sbox{\casebox}{#1}\ifdim\wd\casebox>.97\linewidth\resizebox{.97\linewidth}{!}{\usebox{\casebox}}\else\usebox{\casebox}\fi\endgroup\hfill\null\par}
\newcommand{\CaseDisplay}[1]{\par\smallskip\CaseLine{\(\displaystyle #1\)}\smallskip}

\usepackage[a4paper,top=27mm,bottom=29mm,left=30mm,right=30mm]{geometry}
\usepackage{mathtools,amssymb,booktabs,array,longtable,enumitem,microtype}
\usepackage{placeins}
\usepackage[hidelinks,bookmarksnumbered]{hyperref}
\hypersetup{pdftitle={The finite basis problem for five-element additively idempotent semirings},pdfauthor={Aifa Wang},pdfsubject={Identity bases for five-element additively idempotent semirings},pdfkeywords={finite basis problem, additively idempotent semiring, semiring variety, equational theory}}
\hypersetup{bookmarksdepth=2}
\setlist[enumerate]{leftmargin=2.1em,label=\textup{(\roman*)},itemsep=.25ex,topsep=.6ex}
\setlist[itemize]{leftmargin=1.6em,itemsep=.2ex,topsep=.6ex}
\newtheorem{theorem}{Theorem}[section]
\newtheorem{lemma}[theorem]{Lemma}
\newtheorem{proposition}[theorem]{Proposition}
\newtheorem{corollary}[theorem]{Corollary}
\theoremstyle{definition}
\newtheorem{example}[theorem]{Example}
\theoremstyle{remark}
\newtheorem{remark}[theorem]{Remark}
\numberwithin{equation}{section}

\newcommand{\ct}{\operatorname{con}}
\newcommand{\hd}{\operatorname{hd}}
\newcommand{\V}{\mathsf V}
\newcommand{\Id}{\operatorname{Id}}
\newcommand{\op}{^{\mathrm{op}}}
\newcommand{\sq}[1]{\mathsf{s}(#1)}

\newcommand{\A}[1]{A_{#1}}
\newcommand{\law}{\preccurlyeq}
\newcommand{\ax}{\mathsf{AI}}

\allowdisplaybreaks[1]
\title[The finite basis problem at order five]{The finite basis problem for five-element\\ additively idempotent semirings}
\date{}
\newcommand{\sourceheads}{\operatorname{h}}
\newcommand{\sourcetails}{\operatorname{t}}
\newcommand{\sourceedges}{\operatorname{e}}
\newcommand{\sourcefree}{\mathsf F}
\newcommand{\F}{\sourcefree}
\newcommand{\sourceint}{\operatorname{int}}
\newcommand{\sourcetail}{\operatorname{tl}}
\author{Aifa Wang}
\begin{document}

\begin{abstract}
In the constant-free signature \((+,\cdot)\), the \(15{,}751\) isomorphism types of five-element additively idempotent semirings comprise \(14{,}015\) finitely based and \(1{,}736\) nonfinitely based types. We describe the numbering from the two operations and classify the resulting algebras according to whether their identities have a finite basis. For the finite bases we reduce identities to word absorptions and give deductions on arbitrary alphabets. The obstructions to finite basability use families of identities whose derivations require unboundedly many variables. The cases include identifications with algebras of smaller order. The case register assigns every type to its proof and specifies the smaller generating algebras in the transfer cases.
\end{abstract}
\maketitle
\markboth{Aifa Wang}{The finite basis problem at order five}

\section{Introduction}\label{sec:intro}

Jackson, Ren and Zhao constructed finite additively idempotent semirings whose multiplicative reducts have finite identity bases although the semirings themselves are nonfinitely based \cite{JRZ}. At order four many individual bases are known explicitly \cite{FourII}. Edwards, Mitchell and Ragavan count \(15{,}751\) additively idempotent semirings of order five up to isomorphism \cite[Table~3]{EMR}. Of these, \(3{,}589\) have commutative multiplication \cite[Table~4]{EMR}. We determine which of these semirings admit a finite basis.

We use two binary operations, \(+\) and \(\cdot\). An \emph{ai-semiring} is a join-semilattice under \(+\) and a semigroup under multiplication, with multiplication distributing over addition on both sides. Denote these axioms by \(\ax\). An algebra \(A\) is \emph{finitely based} if some finite subset of its identities \(\Id(A)\) entails every identity in \(\Id(A)\). All bases below are relative to \(\ax\), which is itself finite.

\begin{theorem}[Classification]\label{thm:classification}
Up to isomorphism there are \(14{,}015\) finitely based and \(1{,}736\) nonfinitely based ai-semirings of order five. Among those with commutative multiplication the respective counts are \(2{,}968\) and \(621\). Among those with noncommutative multiplication they are \(11{,}047\) and \(1{,}115\).
\end{theorem}

We write \(\A{i}\) for the algebra numbered \(i\) by the construction in Section~\ref{sec:assembly}. The partition by proof is given in Appendix~\ref{app:register}. Six opposite pairs illustrate how the word arguments depend on the positions of letters. Their multiplication tables are displayed in Table~\ref{tab:twelve}. Reversal of words reduces these twelve types to six tables. The last column gives the size of the basis exhibited here, or NFB for a nonfinitely based algebra. We count each choice of optional factors as a separate identity.

\begin{table}[htbp]\centering
\caption{The six pairs. The \(i\)th entry of a product row is the value at right factor \(i\), with indices starting at \(0\).}\label{tab:twelve}
\small
\begin{tabular}{@{}rrc l r@{}}\toprule
Type & Opposite & Addition & Multiplication & Basis\\\midrule
12462&11849&\(L_4\)&\texttt{44014 11111 11231 33333 44444}&11 laws\\
8224&7605&\(L_{12}\)&\texttt{01444 11111 23444 33333 44444}&15 laws\\
7811&8109&\(L_{11}\)&\texttt{01244 11111 44414 33333 44444}&24 laws\\
7812&8273&\(L_{11}\)&\texttt{01244 11111 44424 44434 44444}&27 laws\\
11670&12469&\(L_1\)&\texttt{40444 41444 00234 33333 44444}&63 laws\\
13396&11218&\(L_4\)&\texttt{44404 01224 22222 22232 44444}&NFB\\\bottomrule
\end{tabular}
\end{table}

The additive orders used in Table~\ref{tab:twelve} have the following covering chains. For \(L_4\) there is one chain \(3<2<1<0<4\). For \(L_{12}\) they are \(1<0<4\) and \(3<2<4\), and for \(L_{11}\) they are \(1<0<4,\ 1<2<4,\ 3<4\). For \(L_1\) they are \(3<2<0<4\) and \(1<0\). In each pair the same labelling identifies the opposite multiplication with the transpose of the displayed table.

The last two pairs differ already on the identity \(xz+xy+yz+z\approx xy+yz+z\). In \(\A{11670}\) the occurrence of \(y\) after \(x\) in \(xy\) can be combined with the occurrence of \(z\) after \(y\) in \(yz\). Iterating this deduction along a directed path is part of the finite-basis proof in Section~\ref{sec:path}. In \(\A{13396}\) the same identity fails at \((x,y,z)=(0,3,2)\). Here it matters which letters occur after \(x\) in the same summand. Section~\ref{sec:private} uses this restriction to force arbitrarily many variables in a basis.

For example, case~\ref{grp:391} of Appendix~\ref{app:register} lists \(\A{10081}\) and refers to Section~\ref{case:227}. There the multiplication table, the component colourings and the hypergraph obstruction give its nonfinite basis proof.

\begingroup\small
\tableofcontents
\endgroup
\clearpage

\section{Words, absorptions and opposite semirings}\label{sec:tools}

Let \(X^+\) be the free semigroup on an infinite alphabet \(X\) and put \(X^*=X^+\cup\{\varepsilon\}\). Write \(\ct(w)\) for the set of letters in \(w\), also called its content or support, and \(\hd(w)\) for its first letter. A word is linear if no letter occurs twice. By distributivity every term is a finite nonempty sum of nonempty words. We identify a sum with its set of summands. The resulting algebra of finite nonempty subsets of \(X^+\), with union and setwise concatenation, is the free ai-semiring on \(X\). Indeed, evaluation at the singleton generators sends each expanded term to precisely its set of words \cite{KP}.

Write \(w\law P\) for the identity \(w+P\approx P\). This expresses the inequality \(w\le P\) in the semilattice order \(a\le b\) defined by \(a+b=b\). Distributivity makes multiplication order-preserving in both arguments. To derive \(P\approx Q\) it suffices to derive \(w\law Q\) for each \(w\in P\) and \(v\law P\) for each \(v\in Q\). Summing the first set of absorptions gives \(P+Q\approx Q\), and summing the second gives \(P+Q\approx P\). Conversely \(P\approx Q\) entails all these absorptions. We shall therefore give completeness proofs for identities \(w\law P\).

For \(D\subseteq X\), call an occurrence of \(x\) in \(u=sxt\) \emph{accessible over \(D\)} if \(\ct(s)\subseteq D\). We set \(\ct(\varepsilon)=\varnothing\). When \(D\) is the set of letters preceding the first \(x\) in \(w\), we have \(x\notin D\). An accessible \(x\) in a summand of \(P\) must then be its first \(x\). In an absorption \(w\law P\) we call \(u\in P\) \emph{contained} if \(\ct(u)\subseteq\ct(w)\).

\begin{lemma}[Reversal and equal varieties]\label{lem:transfer}
Reversing every word in a basis of \(A\) gives a basis of its multiplicative opposite \(A\op\). Algebras generating the same variety have the same identities.
\end{lemma}
\begin{proof}
Reversal fixes variables and sums and sends a product \(st\) to \(t^{\mathrm{rev}}s^{\mathrm{rev}}\). Evaluating the reversed term in \(A\op\) gives the value of the original term in \(A\). Reversal also commutes with substitution when the substituted terms are reversed, and takes a one-hole context to a one-hole context. It thus takes each step of a deduction to a step of a deduction. It is its own inverse. An identity is preserved by direct products coordinatewise, by restriction to a subalgebra, and by lifting the values of its variables through a surjective homomorphism. It therefore holds throughout the variety generated by \(A\) if it holds in \(A\), and the converse follows by evaluation in \(A\) itself \cite{BS}.
\end{proof}

\begin{proposition}[Finite criterion for a variety inclusion]\label{prop:free}
Let \(a_1,\ldots,a_d\) generate a finite algebra \(A\) in this signature, and let \(B\) be finite. In the algebra of functions \(B^{B^d}\) with pointwise operations let \(F\) be the subalgebra generated by the coordinate projections. Then \(A\in\V(B)\) if and only if the assignment
\[
           t^B \longmapsto t^A(a_1,\ldots,a_d)
\]
is well defined on \(F\), in which case it is a surjective homomorphism.
\end{proposition}
\begin{proof}
The elements of \(F\) are exactly the \(d\)-ary term functions of \(B\). If \(A\in\V(B)\), terms defining the same function on \(B\) agree identically on \(A\), hence agree at \((a_1,\ldots,a_d)\). Conversely, once the map is well defined, evaluation preserves both operations and its image contains the generators of \(A\). Thus \(A\) is a homomorphic image of a subalgebra of a power of \(B\). Applying the criterion in both directions proves equality of varieties.
\end{proof}

Brackets in an identity denote an optional factor. For example \(x[p]y\) gives two words, \(xy\) and \(xpy\). A choice for \([p]\) applies at every occurrence in the same identity. Choices for different letters are independent. The empty word may be a prefix or suffix when describing a word. Substitutions in the semiring signature always have nonempty images.

\section{Bounded products and a rectangular identity}\label{sec:general}

\begin{lemma}\label{lem:bottom}
Suppose that the additive order of an \(n\)-element ai-semiring \(A\) has a least element. If \(A\models w\law P\) and \(w\) has \(d\) distinct letters, some nonempty subpolynomial \(P_0\subseteq P\) has at most \(n^{d+1}\) summands and satisfies \(A\models w\law P_0\).
\end{lemma}
\begin{proof}
Given \(\alpha:\ct(w)\to A\), extend \(\alpha\) by the least element on all other letters. From \(P\) choose one summand for each distinct value at this extension. At most \(n\) choices suffice, and their join equals \(P(\alpha)\). Let \(P_0\) be the union of the chosen sets over the \(n^d\) assignments. For any assignment \(\beta\), the extension of \(\beta|_{\ct(w)}\) is pointwise below \(\beta\). Monotonicity gives \(w(\beta)=w(\alpha)\le P(\alpha)\le P_0(\beta)\) for this extension \(\alpha\).
\end{proof}

\begin{theorem}[Bounded word length]\label{thm:bounded}
Let \(A\) be an \(n\)-element ai-semiring with an additive least element. Suppose that a finite \(\Gamma\subseteq\Id(A)\), together with \(\ax\), reduces every term to a sum of words of length at most \(l\). A basis for \(A\) consists of \(\Gamma\) and all its valid absorptions of word length at most \(l\) on \(l(1+n^{l+1})\) fixed variables.
\end{theorem}
\begin{proof}
Reduce both sides of an identity and take one of the resulting absorptions \(w\law P\). The word \(w\) uses at most \(l\) letters. By Lemma~\ref{lem:bottom} there is a subpolynomial \(P_0\) with at most \(n^{l+1}\) summands such that \(w\law P_0\) holds. This identity uses at most \(l(1+n^{l+1})\) letters, so an injective renaming puts it on the stated alphabet. On a finite alphabet there are only finitely many subsets of words of length at most \(l\). Undoing the renaming and adding the summands of \(P\setminus P_0\) on the right gives \(w\law P\).
\end{proof}

For the \(2{,}121\) entries of case~\ref{grp:2}, the least element and the length reduction can be recovered directly from the tables. Test for a row \(b+a=a\) for every \(a\). On the \(125\) assignments to \(x,y,z\), retain every word of length one or two on these letters that lies below \(xyz\). For \(1{,}824\) entries the join of the retained words equals \(xyz\) on all assignments. This single identity reduces each three-letter block and hence bounds the word length by two.

For the other entries form \(E_1=A\), \(E_{j+1}=E_jA\). The first singleton occurs at \(j=3,4,5\) for \(215,27,2\) entries, respectively. Since \(E_{j+1}\subseteq E_j\), every product of at least \(j\) factors then has the same value. In particular \(x_1\cdots x_{j+1}\approx x_1\cdots x_j\) gives the required reduction. The remaining \(53\) tables satisfy \(x\law xy\) and \(y\law xy\). In each of them
\[
 x_1\cdots x_6\ \approx\ \sum_{j=2}^{6}x_1\cdots x_{j-1}x_{j+1}\cdots x_6 .
\]
Every deletion on the right weakly decreases the product. The six nonempty prefix products form a nondecreasing sequence in a five-element order, so some two consecutive ones coincide. Deleting the factor at the latter position leaves the value of the full product unchanged. The identity reduces every block of six letters. These three tests establish the hypotheses of Theorem~\ref{thm:bounded} for the entire case.

The same finite test gives the reduction in case~\ref{grp:26}. For \(k=4,5\), form the word \(u=x_1\cdots x_k\) and retain every nonempty word \(v\) of length less than \(k\) in these letters for which \(v\law u\) holds on all \(5^k\) assignments. Test whether the retained words have join \(u\) on those same assignments. Each table has an additive least element. At \(k=4\) the test gives a reduction for \(228\) of the \(246\) tables, and at \(k=5\) it gives one for the other \(18\). Substitution in consecutive blocks proves the length bound on every alphabet.

\begin{example}
In \(\A{3}\) all products equal \(0\), and its additive semilattice \(L_2\) contains the chain \(3<0<4\) (Table~\ref{tab:tally}). This algebra has the one-identity basis \(xy\approx zt\) relative to \(\ax\). Indeed, every product can be replaced by the same square \(s^2\) on a fresh variable. A reduced polynomial is therefore determined by its one-letter summands and the presence or absence of \(s^2\). If two sets of one-letter summands differ, assign \(4\) to a letter in their symmetric difference and \(3\) to all other letters. One side gives \(4\), the other at most \(0\). If just one side contains \(s^2\) and the one-letter summands agree, their common set must be nonempty. Giving its letters the value \(3\) separates \(3\) from \(3+0=0\). Thus distinct reduced forms have distinct term functions.
\end{example}

\begin{theorem}[A rectangular absorption law]\label{thm:rectangle}
Every subvariety of the variety \(\mathcal D\) of ai-semirings defined by
\[
                            xt\law xy+zt
\]
has a finite basis. Identities involving at most eighteen variables suffice in addition to the displayed law and \(\ax\).
\end{theorem}
\begin{proof}
First substitute \((x,yz,xy,z)\) for \((x,y,z,t)\) in the displayed law to obtain \(xz\law xyz\). Substituting \((xy,t,x,zt)\) gives \(xyzt\law xyt+xzt\). The first inequality, multiplied on the appropriate side, gives the reverse bound for each summand. Induction now gives
\[
                x_1\cdots x_k\ \approx\
                \sum_{j=2}^{k-1}x_1x_jx_k\qquad(k\ge3).
\]
After this reduction let \(U\) be the one-letter summands of \(P\), let \(F,I\) be the sets of first and last letters of its longer summands, and let \(M\) be the set of middle letters of its three-letter summands. Then
\[
       P\ \approx\ N(P)=
       \sum_{u\in U}u+\sum_{f\in F,\,i\in I}fi+
                     \sum_{f\in F,\,m\in M,\,i\in I}fmi .
\]
Here and below an empty sum is omitted. The right-hand side contains each summand of the reduced \(P\). For the reverse inclusion, choose summands \(fs,ti\) of \(P\). The defining law gives \(fi\law fs+ti\). If \(m\in M\), take a third summand \(amb\). Two applications give \(fmb\law fs+amb\) and \(fmi\law fmb+ti\), hence \(fmi\law P\).

In an absorption \(q\law N(P)\), keep the at most three letters of \(q\) fixed. Partition all remaining letters by their membership in \(U,F,I,M\). There are at most fifteen nonempty membership patterns. Send each letter to a representative of its class, obtaining a substitution \(\theta\). Every term of \(\theta N(P)\) still occurs in \(N(P)\) because each position in the display ranges independently over its full set. Thus \(q\law\theta N(P)\) implies \(q\law N(P)\). Conversely substitute \(\theta\) into the latter. We have reduced the absorption to one on at most \(3+15\) variables. On that alphabet the words of length at most three and all the associated absorptions form a finite set. Each subvariety is axiomatized over \(\mathcal D\) by those members of this finite set that it satisfies.
\end{proof}

\section{Joins of two-element semirings}\label{sec:twojoins}

On the chain \(0<1\) take addition to be maximum. Denote by
\(L,R,M,D,N,T\) the six multiplications given, respectively, by first projection,
second projection, maximum, minimum, constant zero and constant one. These
are precisely the associative products distributing over addition on this
chain. Indeed, a product is an increasing Boolean function of two arguments.
Its values at \((0,0)\) and \((1,1)\) either force a constant or leave the
middle two entries free. The four choices for those entries give the two
projections, minimum and maximum, all associative.

For a nonempty subfamily \(J\) write \(B_J\) for the direct product of its
members. We give a uniform proof of the finite basis assertion for these
joins. Shao and Ren proved that every subvariety of the variety generated by the two-element ai-semirings is finitely based \cite[Theorem~3.7]{SR}.

\begin{lemma}\label{lem:boolean-abs}
An absorption \(w\law P\) holds in \(B_J\) precisely when the condition in
each row of the following table indexed by a member of \(J\) holds.
Here a literal is a word of length one.
\begin{center}\small
\begin{tabular}{@{}lp{.82\linewidth}@{}}\toprule
Factor&Condition on the original summands of \(P\)\\\midrule
\(L\)&Some summand has the first letter of \(w\) as its first letter.\\
\(R\)&Some summand has the last letter of \(w\) as its last letter.\\
\(M\)&Every letter of \(w\) occurs in some summand.\\
\(D\)&Some summand has its content contained in \(\ct(w)\).\\
\(N\)&If \(w=x\) is a literal, \(x\) is a summand.\\
\(T\)&Some summand has length at least two, or \(w\) is a literal and occurs as a summand.\\\bottomrule
\end{tabular}
\end{center}
\end{lemma}
\begin{proof}
An inequality in a product holds coordinatewise. For \(L,R,M\) the value of
a polynomial is the maximum of the values assigned to its initial letters,
final letters or all letters, respectively. Assigning one to a missing letter
and zero to every other letter separates a failed condition. In \(D\), give
exactly the letters of \(w\) value one. Then \(w=1\), and a summand takes
value one just when all its letters occur in \(w\). Conversely any such
summand dominates the minimum defining \(w\) on every assignment. In \(N\),
all products of length at least two vanish. In \(T\) all such products equal
one. When every summand is a literal, setting them to zero tests a left-hand
product in \(T\), and assigning one to a missing left-hand literal tests the
remaining case.
\end{proof}

\begin{theorem}\label{thm:boolean-joins}
Each \(B_J\) is finitely based. Relative to \(\ax\) a basis consists of
\[
 xyzt\approx xzyt,\qquad xxy\approx xy,\qquad xyy\approx xy,
\]
together with its valid absorptions on a fixed alphabet of \(141\) letters,
with words of length at most \(143\) and at most six right-hand summands.
\end{theorem}
\begin{proof}
The three word laws hold in each factor. They permit arbitrary adjacent
interchanges strictly between the endpoints of a word. An internal occurrence
of either endpoint can thus be brought next to it and deleted. Repeated
interior letters can be brought together. For an interior double \(xx\)
with a nonempty suffix, apply \(xxy\approx xy\) using that suffix for
\(y\). Thus, for a word of length at least two, its endpoints and content
determine a normal representative: keep both endpoints, and put the other
distinct letters between them in a fixed order. The two copies at the ends
of a square are retained. Literals are kept as written.

First suppose that \(w\law P\) is valid and that \(P\) has at most six
summands. Mark both endpoints of each of the at most seven words involved,
including the left word. This marks at most fourteen letters. Partition the
unmarked letters by the nonempty set of words in which they occur. There
are at most \(2^7-1\) classes. Replace each class by one new letter and
normalize the words. The result uses at most \(14+127=141\) letters and has
length bound \(143\). It holds by substitution. In the reverse substitution
replace each new letter by the product of all the distinct original letters
in its class in a fixed order. Each original word recovers its endpoints,
content and length category. Indeed, the class letters occurred only in
interior positions. The word laws recover the original words, proving the
absorption from the stated family.

Now let \(P\) be arbitrary. If \(w=x\), retain a witness for each applicable
row of Lemma~\ref{lem:boolean-abs}; in the \(M\)-row one witness suffices.
The resulting absorption has at most six summands. If no witness was required,
retain any one summand. For a word of length at least two with endpoints
\(h,t\), select witnesses for \(L,R,D,T\) whenever required. If \(M\in J\),
also select one summand containing \(h\) and one containing \(t\). Let \(P_0\)
be their sum, retaining an arbitrary original summand if no choice was made.
When \(D\in J\), write \(u\) for its selected witness and set \(v=hut\).
Otherwise take \(v=ht\). Then \(\ct(v)\subseteq\ct(w)\). Each required row
holds for \(v\law P_0\); in particular all letters contributed by \(u\)
occur in the selected witness itself. The preceding paragraph gives
\(v\law P_0\), whence \(v\law P\).

Insert each missing letter \(a\in\ct(w)\setminus\ct(v)\) just before the
last position of \(v\), calling the new word \(v'\). When \(M\in J\), take
an original summand \(s\) containing \(a\); otherwise take any original
summand. The old word supplies witnesses for \(L,R,D,T\) in
\(v'\law v+s\), and \(v,s\) together supply the \(M\)-row. Hence this is a
valid two-summand absorption, already derivable. From \(v\law P\) and
\(s\law P\) it follows that \(v'\law P\). After all insertions the word laws
recover \(w\), even when \(h=t\). This proves completeness on every alphabet.
\end{proof}

\begin{proposition}\label{prop:lower-fibres}
Let \(A\) be generated as an algebra by \(a_1,\dots,a_d\) and let \(J\) be a
finite nonempty family of finite ai-semirings. On the disjoint union of the
sets of assignments of these generators in the members of \(J\), form the
vectors of word values. Start with the projections and repeatedly append
one generator on the right until no new vector occurs. The following test
decides whether \(A\) belongs to the variety generated by \(J\).
\begin{enumerate}
\item Give a vector represented by \(w\) the value \(v(w)=w^A(a_1,\dots,a_d)\).
Check the initial projections and each appended-letter transition for
consistency with this assignment.
\item For each \(b\in A\), let \(U_b\) consist of the vectors with assigned
value at most \(b\). If \(U_b\ne\varnothing\), take its coordinatewise join
\(j_b\). Check that every word vector below \(j_b\) belongs to \(U_b\).
\end{enumerate}
Both tests must succeed.
\end{proposition}
\begin{proof}
Appending letters enumerates all nonempty words by induction on length;
there are only finitely many vectors. The first test ensures that the assigned
value is independent of the representing word. If the desired homomorphism
exists, then \(j_b\) maps below \(b\), proving the second test necessary.
Conversely consider a coordinatewise absorption \(w\law P\) of vectors and
put \(b=\sum_{u\in P}v(u)\). Every vector from \(P\) lies in the nonempty
set \(U_b\), and the vector of \(w\) lies below their join, hence below
\(j_b\). The second test gives \(v(w)\le b\). Applying this to the summands
on each side of any equality shows that the map on word vectors extends
to a well-defined map on all their nonempty joins. Concatenation and
distributivity make it multiplicative. It is a surjective homomorphism
because its image contains the generators. Proposition~\ref{prop:free}
(or its identical proof with the disjoint union of coordinates) now applies.
\end{proof}

\section{Finite obstructions for bounded products}\label{sec:trees}

If the additive order has several minimal elements, the extension used in Lemma~\ref{lem:bottom} may be impossible. There is a second way to bound the number of summands needed in \(w\law P\). For each principal ideal, one asks whether all summands of \(P\) can lie in that ideal while \(w\) lies outside it. The register uses the following criterion for the \(108\) types of case~\ref{grp:40} and the \(26\) types of case~\ref{grp:41}.

Let \(A\) have \(q\) elements, and let \(d\ge1\). For \(b\in A\) and \(1\le i\le d\) put
\[
 I_b=\{a\in A:a\le b\},\qquad
 R_i(b)=\{(a_1,\ldots,a_i):a_1\cdots a_i\in I_b\}.
\]
Fix a nonempty \(D\subseteq A\) and \(e\in A\). Give \(D\) the relations \(R_i(e)\cap D^i\) and all singleton unary relations. A constraint, or atom, specifies that a tuple of variables belongs to one of these relations. A variable may occur several times. In the direct square on \(D^2\), a tuple of pairs belongs to a relation if both of its projections belong to the corresponding relation on \(D\). A point \(p\) can be \emph{folded onto} \(p'\ne p\) if replacing any one occurrence of \(p\) by \(p'\) in any relation tuple preserves membership. A fold deletes \(p\) and restricts the relations to the remaining points.

\begin{theorem}[Bounded witnesses]\label{thm:treelength}
Suppose \(A\) satisfies an identity rewriting \(x_1\cdots x_{d+1}\) as a sum of words of length at most \(d\) on these variables. For each \(b\in A\) and each tuple \(t\notin R_j(b)\), \(1\le j\le d\), suppose there exist \(e\in A\), a nonempty \(D\subseteq A\), and \(f:A\to D\) with the following properties.
\begin{enumerate}
\item Coordinatewise, \(f\) takes \(R_i(b)\) into \(R_i(e)\) for every \(i\le d\), and \(f(t)\notin R_j(e)\).
\item There is a total order on \(D\) for which each restricted relation is closed under coordinatewise minimum. Successive folds of \(D^2\) onto its diagonal fix the diagonal pointwise.
\end{enumerate}
Then \(A\) is finitely based. Put \(L=2(q^2-q)\), \(M=\max\{2,q,d\}\), \(T=M^{2L+10}\) and \(B=q^{d+1}T\). A basis consists of the length reduction and the valid absorptions on \(d(B+1)\) fixed variables with word length at most \(d\) and at most \(B\) summands on the right.
\end{theorem}
\begin{proof}
Fix \(D,e\) as in the theorem. We claim that every unsatisfiable finite instance over this structure has an unsatisfiable subset of at most \(T\) atoms.

Begin with \(D\) as the domain of every variable. In one round delete a value at \(x\) if some occurrence of \(x\) has no supporting tuple with that value and with every other entry in its current domain. Supports at distinct occurrences may be chosen independently. If the domains stabilize and are nonempty, assigning to each variable its least remaining value satisfies every atom. For a given atom choose, at each position, a supporting tuple attaining the minimum there. Their coordinatewise minimum is a tuple in the relation and realizes all the minima simultaneously. Repeated occurrences receive the same value.

Each deletion has a witness whose variable--atom incidence graph is a rooted tree. Indeed, take the atom in which support is missing, put \(x\) at the root, and use a separate fresh variable at every other position. To each of those variables attach copies of the trees for values deleted from its domain in earlier rounds. Induction on the round shows that the resulting tree forbids the deleted value at its root. It maps to the original instance by sending every copy to its original variable or atom. In an unsatisfiable instance a domain must empty. Identify the roots of the trees forbidding its individual values and take an inclusion-minimal unsatisfiable subset of the resulting atoms. After isolated variables are discarded this subset is connected, since assignments on disjoint components can be combined. Its incidence graph is a tree.

At a variable vertex of this minimal tree, each incident branch allows a set of possible values. The intersection is empty. Dropping any branch leaves a satisfiable proper subinstance, so there is a value forbidden by that branch and allowed by all the others. Distinct branches require distinct such values. Thus variable vertices have degree at most \(|D|\), and atom vertices have degree at most \(d\).

The folds give a sequence \(g_0,\ldots,g_\ell:D^2\to D\) from first to second projection, with \(\ell\le L\). Compose the successive retractions of \(D^2\) and take their first projections. At the diagonal the two projections coincide. Continue with the second projections of the same retractions in reverse order. For two tuples in any relation, applying either \(g_j\) or \(g_{j+1}\) independently at each position stays in that relation. In the step of a fold this follows by replacing the selected occurrences one by one in the square and then projecting. Each projection path uses at most \(|D|^2-|D|\) folds.

Choose distinct atoms \(E,F\) of the minimal tree. Take assignments \(\alpha,\beta\) satisfying all atoms except, respectively, possibly \(F,E\). In the graph joining variables of a common atom put
\[
 j(x)=\min\{\ell,\operatorname{dist}(x,\operatorname{vars}(E))\}.
\]
If all variables of \(F\) were at distance at least \(\ell\), the assignment \(x\mapsto g_{j(x)}(\alpha(x),\beta(x))\) would satisfy \(E\) by the first projection and \(F\) by the second. At every other atom the indices are equal or consecutive, so it would satisfy that atom as well. Consequently some variables of any two distinct atoms lie at distance less than \(\ell\). When \(\ell=0\) the same argument precludes two distinct atoms altogether. Each variable-graph step uses two edges of the incidence tree. Its diameter is thus at most \(2L+6\), allowing for the endpoint vertices. Maximum degree \(M\) bounds the number of vertices by \(T\), by summing \(1+M+\cdots+M^{2L+6}\). The image of the atom set in the original instance is unsatisfiable, since a solution would pull back to the tree. This gives the asserted bound of \(T\) original atoms.

Apply the bound to an absorption \(w\law P\). First reduce both sides to words of length at most \(d\). The length identity applied to a block of \(d+1\) letters shortens every resulting word, so induction suffices. Fix \(b\in A\) and an assignment \(\alpha\) on \(\ct(w)\) with \(w(\alpha)\notin I_b\). Choose \(f,D,e\) for the tuple of values of \(w\) in word order. On \(D\) require each word of \(P\) to have product in \(I_e\), and use singleton relations to prescribe \(f(\alpha(x))\) at every \(x\in\ct(w)\). A solution would put \(P\) in \(I_e\) and \(w\) outside \(I_e\), contrary to the absorption. Select the words of \(P\) corresponding to the product constraints in an unsatisfiable subset of at most \(T\) atoms. Those constraints together with the prescribed singleton values are already inconsistent.

Let \(Q\) be the union of the selected words over all choices of \(b,\alpha\). There are at most \(q^{d+1}\) choices. If none contributes a word, take any one summand of \(P\). For an arbitrary assignment \(\beta\), put \(b=Q(\beta)\). If \(w(\beta)\notin I_b\), choose the map for this \(b\) and \(\alpha=\beta|_{\ct(w)}\). Every selected word takes a value at most \(b\), so by (i) the assignment \(f\circ\beta\) satisfies the selected product constraints and all the singleton constraints. This contradiction proves \(w\law Q\). Since \(Q\) has at most \(B\) words, \(w\law Q\) is a renaming of one of the finitely many stated absorptions. Adding the other summands on the right gives \(w\law P\).
\end{proof}

\section{The strong nonfinite basis property of \texorpdfstring{\(S_7\)}{S7}}\label{sec:strong}

Let \(S_7=\{1,a,\infty\}\), with multiplicative identity \(1\), multiplicative zero \(\infty\), and \(a^2=\infty\). Addition is idempotent and any two distinct elements sum to \(\infty\). L.~Wang, Q.~Yin and A.~Wang proved that every finite ai-semiring whose variety contains \(S_7\) is nonfinitely based \cite[Theorem~1.1]{Strong}. This theorem accounts for \(1{,}660\) five-element types in Appendix~\ref{app:register}. We give the proof, beginning with a property of squarefree monomials in \(S_7\).

\begin{lemma}[Private points]\label{lem:privatepoints}
Let \(\mathcal E\) be a finite nonempty family of distinct nonempty finite sets. Suppose each \(E\in\mathcal E\) has a \emph{private point}, one belonging to no other member. Write \(m_E=\prod_{x\in E}x\), a commutative monomial, and \(p=\sum_{E\in\mathcal E}m_E\). The commutative monomials \(m\) satisfying \(S_7\models m\law p\) are precisely the \(m_E\).
\end{lemma}
\begin{proof}
If \(m\) involved a letter outside the union, giving that letter the value \(\infty\) and all others \(1\) would separate \(m\) from \(p\). Choose one private point \(z_E\) in each member, assign \(a\) to these points and \(1\) elsewhere. Every \(m_E\) has value \(a\). The only element below \(a\) is \(a\), so the sum of the exponents in \(m\) of the chosen letters is one. Varying \(z_E\) with the other choices fixed shows that every private point of \(E\) has the same exponent \(e_E\). Exactly one member \(E_*\) has \(e_{E_*}=1\).

For a nonprivate point \(x\), take \(x\) and one private point from each member missing \(x\). Again each member contains exactly one chosen point. If \(x\in E_*\), the other chosen points have exponent zero, forcing exponent one at \(x\). If \(x\notin E_*\), the private point of \(E_*\) already contributes one. Thus \(m=m_{E_*}\). The converse follows since every \(m_E\) is a summand. We shall also use the construction just given of a set through any prescribed point meeting each member exactly once.
\end{proof}

\begin{lemma}[Invariance under a bounded number of variables]\label{lem:sparse}
Let \(n,k\ge1\). Suppose all members of a finite family \(\mathcal F\) have size at least \(k\), distinct members intersect in at most \(l\) points, and \(k>(2^n-2)l\). If every word in a polynomial is linear and its set of letters is a member of \(\mathcal F\), this remains true after any application of an identity on at most \(n\) variables of \(S_7\).
\end{lemma}
\begin{proof}
Expand the context of the replacement by distributivity. It suffices to compare \(L\sigma(f)R\) and \(L\sigma(g)R\) for fixed exterior words \(L,R\in X^*\), where \(S_7\models f\approx g\). Every summand of \(L\sigma(f)R\) is linear and has content in \(\mathcal F\). The polynomial \(f\) has the same variables as \(g\), by assigning \(\infty\) to one variable and \(1\) to the others. Each word of \(f\) must itself be linear, since two occurrences could be expanded with the same nonempty image word. After commuting letters there are at most \(2^n-1\) distinct monomials in \(f\).

For each variable \(y\) choose a single word \(u_y\in\sigma(y)\), and let \(\tau\) denote this word substitution. The supports of the words of \(L\tau(f)R\) form a subfamily \(\mathcal E\) of size at most \(2^n-1\). In any member of \(\mathcal E\) at most \((|\mathcal E|-1)l<k\) points occur in another member. Each member therefore has a private point.

Take any \(z\in\ct(u_y)\). Linearity gives \(z\notin\ct(LR)\). There is a set through \(z\) that meets each member of \(\mathcal E\) exactly once and avoids \(\ct(LR)\). For a singleton family take \(\{z\}\). Otherwise every exterior letter belongs to every member, so the private points used in Lemma~\ref{lem:privatepoints} all lie outside \(\ct(LR)\). Assign \(a\) to the chosen set and \(1\) elsewhere. Every word of \(L\tau(f)R\) takes the value \(a\), and every exterior factor has value \(1\). In particular \(u_y\) takes the value \(a\). A summand of \(g\) repeating \(y\) would then take the value \(\infty\), contradicting \(\tau(g)=\tau(f)=a\). Since \(y\) was arbitrary, every word of \(g\) is linear.

Now take an arbitrary word in \(L\sigma(g)R\). By this linearity its choices of image words extend to a single word substitution \(\tau\) on all variables. In \(S_7\) the chosen word lies below \(L\tau(f)R\). On commuting letters, Lemma~\ref{lem:privatepoints} identifies it with a squarefree monomial of \(L\tau(f)R\). It is therefore linear and its content belongs to \(\mathcal F\), as required.
\end{proof}

\begin{lemma}[Realizing finite products]\label{lem:products}
For a finite semigroup \(S\) and integer \(n\ge1\), there are a positive integer \(N\) and a family \(\mathcal F\) of proper subsets of \(\{1,\ldots,N\}\) satisfying the hypotheses of Lemma~\ref{lem:sparse}, such that for every \(s_1,\ldots,s_N\in S\), some ordering of the positions in a member has product \(s_1\cdots s_N\).
\end{lemma}
\begin{proof}
Write \(r=|S|\). Use one state for each partial product and a separate initial state. Reading a letter \(s\) in the initial state leads to \(s\), and reading it in state \(t\) leads to \(ts\). Thus the paths from the initial state to \(b\) read exactly the nonempty words of product \(b\).

For completeness, number the states. The words on paths from \(i\) to \(j\) with internal state numbers at most \(h\) form a language \(K^h_{ij}\). At level zero take the direct edges and also the empty word when \(i=j\). At level \(h\) one may avoid \(h\), or decompose a path at its successive visits to \(h\). Consequently
\[
 K^h_{ij}=K^{h-1}_{ij}\ \cup\ K^{h-1}_{ih}(K^{h-1}_{hh})^*K^{h-1}_{hj}.
\]
The letter-count map takes concatenation to vector addition. Finite unions of linear sets are closed under union and vector addition. To check closure under the star, write a semilinear set as \(\bigcup_{i=1}^m(u_i+\langle V_i\rangle)\), where \(\langle V\rangle\) denotes the additive monoid of nonnegative integral combinations of the finite set \(V\). In a sum of vectors let \(I\) be the set of components used at least once. The possible sums with this \(I\) form exactly
\[
  \sum_{i\in I}u_i+
  \left\langle\,\{u_i:i\in I\}\cup\bigcup_{i\in I}V_i\,\right\rangle .
\]
One occurrence of each component supplies its base and may carry any prescribed combination of its periods. Extra copies supply the additional bases. The empty choice gives the zero vector. There are finitely many choices of \(I\), so induction on \(h\) proves the required regular-language case of the counting theorem \cite{Parikh}. In particular the multiplicity vectors for product \(b\) form a finite union of sets
\(u+\sum_j\mathbb N v_j\), with nonnegative integral vectors. Choose \(B\ge1\) bounding \(\|u\|_1+\sum_j\|v_j\|_1\) in all components for all \(b\).

Fix \(k\), set \(L=k+B\), \(K=k+2B\) and \(C=rB\). For an assignment with count vector \(d\) and full product \(b\), select a decomposition \(d=u+\sum_j m_jv_j\) in the representation for \(b\). For \(N\) sufficiently large put
\[
             e=u+\sum_j\lfloor Lm_j/N\rfloor v_j .
\]
For \(N\ge L\) we have \(0\le e\le d\) coordinatewise, and \(e\) belongs to the same linear set as \(d\). In \(e-(L/N)d\) the coefficient of \(u\) lies in \([0,1]\) and each coefficient of \(v_j\) has absolute value at most one. Thus \(\|e-(L/N)d\|_1\le B\). Since \(\|d\|_1=N\), this gives \(k\le s=\|e\|_1\le K\). Every subset of positions having multiplicities \(e\) can be ordered to realize \(b\).

Choose a size uniformly from \(k,\ldots,K\) and a subset uniformly of that size. Conditional on the size \(s\), the probability of these counts is
\(\prod_i\binom{d_i}{e_i}/\binom Ns\). For \(d_i<N/(2L)\) the rounding estimate gives \(e_i\le B\), so these coordinates use at most \(C\) positions together. For the other coordinates \(d_i\ge2K\) as soon as \(N\ge4LK\). Their binomial factors are at least \((N/(4L))^{e_i}/e_i!\). Bounding the denominator by \(N^s\), we obtain a probability at least \(cN^{-C}\), including the choice of size, for a constant \(c>0\) independent of \(N\) and of the assignment. The exponent \(C\) is independent of \(k\).

Now take \(t=2C+3\), \(l=t-1\) and \(k=\max\{t,(2^n-2)l+1\}\). Draw \(M=\lceil c^{-1}(\log r+2)N^{C+1}\rceil\) subsets independently. The probability that at least one of the \(r^N\) assignments is missed is at most \(r^N e^{-cMN^{-C}}\le e^{-2N}\). For a pair of sampled subsets the probability of an intersection of size \(t\) or more is at most \(\binom Kt (K/(N-t+1))^t\), by conditioning on the first set and bounding the chance of each prescribed \(t\)-subset of it. The sum over all pairs tends to zero since \(2C+2-t=-1\). For large \(N>K\) both desired properties therefore hold simultaneously. The subsets are proper, and repeated sampled sets would have intersection at least \(k\ge t\).
\end{proof}

\begin{theorem}[Strong nonfinite basis property of \(S_7\) {\cite{Strong}}]\label{thm:strong}
If \(A\) is a finite ai-semiring and \(S_7\in\V(A)\), then \(A\) is nonfinitely based.
\end{theorem}
\begin{proof}
Suppose \(A\) has a finite basis, with at most \(n\) variables in each identity. Every such identity holds in \(S_7\). Apply Lemma~\ref{lem:products} to the multiplicative reduct of \(A\). On variables \(x_1,\ldots,x_N\), form \(P\) by summing all orderings of all members of \(\mathcal F\), and put \(w=x_1\cdots x_N\). Under each assignment some summand of \(P\) has the same value as \(w\). Thus \(P+w\approx P\) holds in \(A\). In any deduction from \(P\), however, Lemma~\ref{lem:sparse} confines the content of every word to a member of \(\mathcal F\). Since \(\ct(w)\) has \(N\) elements and every member is proper, such a deduction cannot yield \(P+w\).
\end{proof}

\section{Repetition of the initial letter}\label{sec:heads}

In \(C=\A{12462}\) and \(J=\A{8224}\) the elements \(1,3,4\) are left zeros. For a word starting at either of the other two elements, the product can change when the first letter recurs. This accounts for the additional condition in the following characterization. For a polynomial \(P\) put \(H=\{\hd(u):u\in P\}\).

\begin{proposition}\label{prop:heads}
In both algebras \(w\law P\) requires the following conditions.
\begin{enumerate}
\item There is a contained summand of \(P\).
\item For each first occurrence of a letter \(y\) in \(w\), there is an occurrence of \(y\) in a summand of \(P\) accessible over the preceding letters of \(w\).
\end{enumerate}
Let \(h=\hd(w)\). Consider the tail of \(w\) up to and including its first \(h\), or the entire tail if \(h\) never recurs. At each first occurrence of a letter \(y\) in that part, let \(D\) consist of the letters strictly before \(y\) \emph{in the tail}. In addition to \textup{(i)} and \textup{(ii)} require the condition for the algebra in question.
\begin{enumerate}[label=\textup{(\alph*)}]
\item In \(C\) some \(hstb\in P\) satisfies \(\ct(s)\subseteq D\) and \(t\in\{h,y\}\). If \(h\) never recurs, put \(D_*\) equal to the content of the whole tail. In this case \(P\) must also contain either a word \(hv\) with \(\ct(v)\subseteq D_*\) or a word \(hshb\) with \(\ct(s)\subseteq D_*\).
\item In \(J\) either \(D\cap H\ne\varnothing\) or some \(hstb\in P\) satisfies \(\ct(s)\subseteq D\) and \(t\in H\cup\{y\}\).
\end{enumerate}
In each algebra the stated conditions are sufficient.
\end{proposition}
\begin{proof}
In \(C\), a product starting at \(2\) stays there as long as the factors are \(2\). The first factor different from \(2\) sends the product to \(3\) if it is \(3\), and to \(1\) otherwise. A product starting at \(0\) also ignores an initial string of \(2\)'s in the tail. Its next factor sends it to \(1\) if it is \(3\), and to \(4\) otherwise. Each change takes the product to a left zero.

To see the necessity of (i) in \(C\), give the letters of \(w\) the value \(2\) and all others \(3\). If every summand contained an outside letter, we would have \(w=2\) and \(P=3\). For (ii), give the letters preceding \(y\) value \(2\), give \(y\) value \(1\), and give the remaining letters value \(3\). Then \(w=1\), whereas a summand with no accessible \(y\) takes a value at most \(2\) in the additive order. For (a), assign \(0\) to \(h\), \(2\) to \(D\), \(1\) to \(y\) if \(y\ne h\), and \(3\) elsewhere. The value of \(w\) is \(4\). A summand can attain \(4\) only if it starts at \(h\) and its first tail letter outside \(D\) is \(h\) or \(y\). When \(h\) never recurs, use \(D_*\) in place of \(D\), omitting the assignment of \(1\). Now \(w=0\). A summand attaining \(0\) or \(4\) gives precisely one of the two alternatives at the end of (a).

Conversely, (ii) at the initial position gives a summand with first letter \(h\). It suffices on any assignment where \(h\) is a left zero. If \(h=2\) and all letters of \(w\) have value \(2\), use (i). If the first different factor sends the product to \(1\), its witness in (ii) attains at least \(1\). The other possible value, \(3\), is least. Suppose \(h=0\). The first tail letter of value different from \(2\), if present, occurs for the first time in the tail, at or before the first return to \(h\). If its value is \(0,1\) or \(4\), the summand in (a) attains \(4\). If it is \(3\), every summand starting with \(h\) dominates \(w=1\). When the tail is all \(2\), the last assertion of (a) gives a summand attaining \(0\) or \(4\).

In \(J\) a product starting at \(0\) or \(2\) is unchanged by an initial string of zeros in the tail. The first nonzero tail factor sends these two states to \(1,3\), respectively, if the factor is \(1\), and to \(4\) otherwise. For (i) give the letters of \(w\) value \(0\) and all others \(1\). For (ii) give the predecessors of \(y\) value \(0\), give \(y\) value \(4\), and use \(1\) elsewhere. Without an accessible \(y\) every summand lies below \(4=w\) on the chain \(1<0\). If (b) fails then \(D\cap H=\varnothing\). Give \(D\) value \(0\), \(h\) value \(2\), the other members of \(H\cup\{y\}\) value \(3\), and the remaining letters \(1\). Every summand starts on the chain \(3<2\). To reach \(4\) it must start with \(h\) and its first tail letter outside \(D\) must belong to \(H\cup\{y\}\). This happens in \(w\) and, by hypothesis, in no summand of \(P\).

For sufficiency in \(J\), the left-zero cases again follow from (ii). Suppose first that \(h=0\). If all tail entries are zero, (i) supplies a summand of value \(0\). A first nonzero tail entry \(1\) gives \(w=1\), dominated by any summand starting with \(h\). For a first nonzero tail entry of value \(2,3\) or \(4\), an accessible witness either reaches \(4\) after a prefix of zeros or starts at that entry, on the other chain or at the top. Its join with an \(h\)-headed summand is \(4\).

Suppose next that \(h=2\). A summand starting with \(h\) takes a value in \(\{3,2,4\}\). If any summand starts on the other chain or at \(4\), their join is \(4\). Otherwise all summands start at \(2\) or \(3\). When the tail of \(w\) is all zero, a contained summand must start with \(h\) and takes the value \(2\) or \(4\). When the first nonzero tail entry is \(1\), \(w=3\) already lies below the \(h\)-headed summand. For any other first nonzero tail letter \(y\), its preceding tail support \(D\) has value \(0\), so \(D\cap H=\varnothing\). Condition (b) gives an \(h\)-headed summand whose first tail entry outside \(D\) has value \(2,3\) or \(4\). It attains \(4\).
\end{proof}

\begin{theorem}\label{thm:heads}
A basis for \(C\) consists of the seven common laws
\begin{gather}
 xyy\approx xy,\qquad xyzy\approx xyz,\qquad
 h\law hh,\qquad hp\law hph,\label{eq:tailword}\\
 ayv\law av+ayb,\qquad ayv\law av+ay,\qquad ay\law a+ayb \label{eq:attach}
\end{gather}
and the four laws
\begin{equation}\label{eq:c-return}
                h[p]hu\law u+h[p]h[b].
\end{equation}
A basis for \(J\) consists of the same seven common laws and the eight laws
\begin{equation}\label{eq:j-return}
                h[p]hu\law u+h[p]q[b]+q[c].
\end{equation}
Word reversal gives bases for \(\A{11849}\) and \(\A{7605}\).
\end{theorem}
\begin{proof}
The first two identities of \eqref{eq:tailword} reduce a word to its first letter followed by the distinct letters of its tail in order of first occurrence. Indeed, the context before the first of two equal tail letters is nonempty. Use \(xyy\approx xy\) when the occurrences are adjacent and \(xyzy\approx xyz\) when they are separated. Both identities preserve the word values computed above. The other two identities of \eqref{eq:tailword} insert another \(h\) after a prefix starting with \(h\). A product already at a left zero is unchanged. In \(C\) a product still at its initial \(0\) goes to \(4\) and one at its initial \(2\) stays there. In \(J\) insertion fixes an initial \(0\) and sends an initial \(2\) to \(4\). These are increases in the respective additive orders.

For \eqref{eq:attach}, the summand starting with \(ay\) witnesses \(y\). In the first two identities \(av\) is contained, witnesses \(v\), and satisfies the last condition of (a). In the third, the one-letter summand \(a\) gives containment and the last condition of (a). In \eqref{eq:c-return} the \(h\)-headed summand witnesses each tail letter up to the repeated \(h\). In \eqref{eq:j-return} it witnesses the part up to \(p\), and at the return we have \(q\in H\). In both identities the summand \(u\) is contained and witnesses \(u\). Thus Proposition~\ref{prop:heads} proves the displayed identities.

For completeness, suppose \(w\law P\) is valid and write \(w\) in the reduced form. For a prefix \(p\) of its tail we first prove \(hpv\law P\) for all \(hv\in P\). This is immediate for \(p=\varepsilon\). Suppose there is a summand \(hsqb\) with \(\ct(s)\subseteq\ct(p)\), where \(q=h\) in \(C\) and \(q\in H\) in \(J\). The induction hypothesis gives \(hpsqb\law P\), and repeated tail letters reduce it to \(hpqb\law P\). In \(C\), \eqref{eq:c-return} now gives \(hphu\law P\) for \emph{every} \(u\in P\). In \(J\), choose a summand \(qc\) and apply \eqref{eq:j-return} to get the same conclusion. The position just after \(p\) is where we have inserted an extra \(h\).

As long as no such \(hsqb\) exists, let \(y\) be the next letter of the tail. Condition (a) in \(C\) supplies \(hsyb\) with \(\ct(s)\subseteq\ct(p)\). In \(J\), \(\ct(p)\cap H\) is empty throughout this part of the induction. It starts empty, so (b) supplies \(hsyb\), and if \(y\) belonged to \(H\) that word would already have given a suitable \(hsqb\). In either algebra \(hpyb\law P\) follows. Put \(a=hp\). For each \(hv\in P\) we have \(av\law P\) and \(ayb\law P\). The appropriate identity of \eqref{eq:attach} gives \(ayv\law P\), with \(ay\) already available when both \(v\) and \(b\) are empty. This proves the assertion for \(py\). A return to \(h\) always yields a suitable \(hsqb\).

After insertion of \(h\), write \(p\) for the tail constructed so far, including the inserted letter. We have \(hpu\law P\) for every \(u\in P\). If \(h\) occurs later in the prescribed tail, its copy is already in \(p\) and need not be inserted again. For any other new \(y\), (ii) gives \(syb\in P\) with \(\ct(s)\subseteq\ct(p)\), since \(p\) now contains \(h\). Hence \(hpsyb\approx hpyb\law P\), and \eqref{eq:attach} gives \(hpyu\law P\) for all \(u\in P\). At the end choose a contained summand \(u\). All its letters belong to \(\ct(p)\), so \(hpu\approx hp\). The reduced form of \(w\) lies below this \(hp\) by \(h\law hh\) or \(hp'\law hp'h\) at the prefix \(hp'\) where insertion occurred. A later duplicate in the tail disappears by \eqref{eq:tailword}.

There remains the case in which no \(hsqb\) exists even at the end of the tail. In \(C\), (a) then supplies \(hv\in P\) with \(\ct(v)\subseteq\ct(p)\). In \(J\), let \(hv\) be a contained summand. Its first letter must be \(h\), since \(\ct(p)\cap H=\varnothing\). If \(v\) contained \(h\), its first such occurrence would give the excluded \(hsqb\) with \(q=h\). Thus also in \(J\) we have \(\ct(v)\subseteq\ct(p)\). In both cases \(w=hp\approx hpv\law P\). This includes an empty tail.
\end{proof}

\section{Immediate successors in \texorpdfstring{\(\A{7811}\)}{A7811}}\label{sec:cross}

Put \(B=\A{7811}\) and fix \(P\). Write \(F(D)\) for the set of letters with an occurrence in \(P\) accessible over \(D\), and \(U(D)\) for the condition that \(\ct(u)\subseteq D\) for some \(u\in P\). For \(x\notin D\), let \(J_x(D)\) consist of the letters immediately following an accessible \(x\) in a summand of \(P\). Write \(E_x(D)\) when such an \(x\) occurs at the end of a summand. Put
\begin{equation}\label{eq:imark}
 M_x(D)\ \Longleftrightarrow\
 J_x(D)\cap(D\cup\{x\}\cup F(D))\ne\varnothing
 \quad\text{or}\quad \bigl(E_x(D)\ \text{and}\ U(D)\bigr).
\end{equation}
In particular, \(q\in J_x(D)\cap F(D)\) may be witnessed by two summands, one with \(q\) immediately after \(x\) and another with \(q\) accessible over \(D\).

\begin{lemma}\label{lem:immediate}
The identity \(w\law P\) holds in \(B\) if and only if \(P\) has a contained summand and the following holds at each first occurrence of a letter \(x\) in \(w\). For \(D\) the content of the preceding prefix, \(x\in F(D)\). If \(M_x(D)\) fails and \(x\) has a successor in \(w\), that successor lies in \(J_x(D)\). If \(M_x(D)\) fails and \(x\) is final, \(E_x(D)\) holds.
\end{lemma}
\begin{proof}
A word all of whose letters have value \(0\) evaluates to \(0\). Otherwise consider the first nonzero factor. A \(1\) or a \(4\) determines the value of the word. A \(3\) gives \(3\) in initial position and \(4\) elsewhere. A \(2\) gives \(2\) in final position, \(1\) if its immediate successor is \(3\), and \(4\) otherwise. These rules follow from the rows of \(B\) in Table~\ref{tab:twelve}. For necessity of containment assign \(0\) on \(\ct(w)\) and \(1\) elsewhere. For accessibility assign \(0\) on \(D\), \(4\) at \(x\) and \(1\) elsewhere.

Suppose the successor or final-position condition fails with \(M_x(D)\) false. Give \(D\) value \(0\), \(x\) value \(2\), \(J=J_x(D)\) value \(3\), and the other letters \(1\). This is well defined since \(J\cap(D\cup\{x\})=\varnothing\). In a summand of \(P\) the first nonzero factor cannot be \(3\), since \(J\cap F(D)=\varnothing\). If it is \(2\), that factor is either final or immediately followed by \(3\). Thus each summand evaluates to \(0,1\) or \(2\). Attainment of \(0\) requires \(U(D)\) and attainment of \(2\) requires \(E_x(D)\). If \(x\) is nonfinal in \(w\), its unrepresented successor makes \(w=4\), while at most one of \(0,2\) occurs among the summand values. If \(x\) is final in \(w\), then \(w=2\) and \(E_x(D)\) is false, so \(P\) has value \(0\) or \(1\). Both cases contradict \(w\law P\).

For sufficiency fix an assignment. If \(P=3\), every summand starts with a letter valued \(3\). Accessibility at the head of \(w\) gives \(w=3\) as well. Otherwise \(P\ge1\). A word all of whose letters have value zero is covered by a contained summand. If the first nonzero factor is \(4\), its accessible occurrence gives a summand of value \(4\). If that factor is \(3\), its witness gives \(3\) or \(4\), so \(P=4\) in the case under consideration. A first nonzero factor \(1\) gives \(w=1\).

It remains to consider a first nonzero factor \(x\) of value \(2\). A successor in \(J_x(D)\cap(D\cup\{x\})\) makes the corresponding summand \(4\). For \(q\in J_x(D)\cap F(D)\), the same holds unless \(q\) has value \(3\). In that event the summand with \(xq\) gives \(1\) and the separate witness for \(q\) gives \(3\) or \(4\), whose join with \(1\) is \(4\). Finally \(E_x(D)\) and \(U(D)\) give two summands taking values \(2\) and \(0\), again with join \(4\). This proves the assertion when \(M_x(D)\) holds. In the remaining cases the matching successor or final occurrence of \(x\) supplies a summand with the same value as \(w\).
\end{proof}

\begin{theorem}\label{thm:immediate}
The following twenty-four identities form a basis for \(B\), with their reversals giving a basis for \(\A{8109}\).
\begin{align}
 x^3&\approx x^2,& xyx&\approx xy^2,& xyzx&\approx xyz^2,\label{eq:imword}\\
 x&\law x^2,& d^2u&\law u+d^2,& d^2u&\law u+d^2b,\label{eq:imbase}\\
 [a][m]x^2u&\law u+[a][m]xq[b]+[a]q[c],\label{eq:imcross}\\
 a[m]x^2u&\law u+a[m]x+a.\label{eq:imend}
\end{align}
The last two rows expand into sixteen and two identities, respectively.
\end{theorem}
\begin{proof}
List the letters of a word in order of first occurrence. In this list double a letter \(x\) if in the word its first occurrence is immediately followed by a letter already present up to and including that \(x\). This gives a normal form under \eqref{eq:imword}. To see this, append one letter to a reduced prefix ending in \(z\). A new letter extends the list. If the appended letter is an earlier \(x\ne z\), the word ends in \(xzx\) or \(xUzx\) with \(U\ne\varepsilon\), starting at a previous \(x\). The second or third identity changes its final \(x\) to \(z\). A third copy of \(z\) contracts to \(z^2\), which also handles an appended \(z\). The value rules in the proof of Lemma~\ref{lem:immediate} show that these equations hold in \(B\).

For the inequalities in \eqref{eq:imbase}, \(x\law x^2\) follows from the diagonal, and the displayed repeat of \(d\) gives a successor in \(\{d\}\). In \eqref{eq:imcross} one summand places \(q\) immediately after \(x\) and another makes \(q\) accessible over the prefix of \(x\). In \eqref{eq:imend}, one summand ends in \(x\) and the summand \(a\) lies in its preceding support. For the earlier letters, a summand containing the displayed prefix supplies each required successor. In each inequality with a final \(u\) on the left, the summand \(u\) is contained and witnesses that final occurrence. Lemma~\ref{lem:immediate} now gives validity.

Suppose \(B\models w\law P\) and normalize \(w\). Write \(x_1,\ldots,x_n\) for its distinct letters in order and \(D_i=\{x_1,\ldots,x_{i-1}\}\). Mark \(i\) if \(M_{x_i}(D_i)\) holds. Every doubled letter is marked, since its successor in the normal form is itself. Write \(w^*\) for the word obtained by doubling all marked letters and leaving the others single. By \eqref{eq:imbase} we have \(w\law w^*\). We shall derive \(w^*\law P\).

Let \(L\) be the prefix of \(w^*\) through the last marked index already treated. Initially \(L=\varepsilon\). We have \(Lu\law P\) for every \(u\in P\), initially by idempotence. Let \(s\) be the next index, and consider \(i\ge s\) with \(s,\ldots,i-1\) unmarked. Set \(p=Lx_s\cdots x_{i-1}\).

Consider a prefix over \(D_i\) of a summand in \(P\). After multiplication on the left by \(L\), its initial letters from \(\ct(L)\) can be deleted, since a nonempty \(L\) ends with a double. The remaining letters have strictly increasing indices. At the first repetition, decrease, or return to \(\ct(L)\), the immediately preceding \(x_j\) would be accessible over \(D_j\) and followed by a letter in \(D_j\cup\{x_j\}\). This would mark \(j\). Nor can an index be skipped. If the first skip ends just before \(x_j\), the summand makes \(x_j\) accessible over \(D_{j-1}\). Since \(j-1\) is unmarked, the successor condition in \(w\) gives \(x_j\in J_{x_{j-1}}(D_{j-1})\). Together these facts mark \(j-1\), a contradiction. This argument also applies through a required occurrence of \(x_i\). Consequently a witness through \(x_i\) gives \(px_i b\law P\) with its own suffix \(b\) intact. For a prefix over \(D_i\) alone, the same reduction gives an initial segment of \(p\).

At an unmarked final index the witness ending in \(x_i\) therefore gives \(w^*\law P\). At a marked index put \(x=x_i\). For each \(u\in P\) we derive \(px^2u\law P\) as follows.
\begin{enumerate}
\item If a successor \(q\) lies in \(D_i\cup\{x\}\), we have \(pxqb\law P\) with \(q\in\ct(px)\). The normal forms of \(pxq\) and \((px)^2\) agree. In both, the last \(x\) doubles and subsequent letters from \(\ct(px)\) leave the form unchanged. Taking \(d=px\) in \eqref{eq:imbase} gives \((px)^2u\approx px^2u\law P\), using the appropriate identity for empty or nonempty \(b\).
\item Suppose \(q\in J_x(D_i)\cap F(D_i)\) and \(q\notin D_i\cup\{x\}\). The two summands give \(pxqb\law P\) and \(aqc\law P\), where \(a\) is an initial segment of \(p\). Write \(p=am\) and use \eqref{eq:imcross}, omitting whichever of \(a,m,b,c\) are empty.
\item If \(E_x(D_i)\) and \(U(D_i)\) hold, one summand gives \(px\law P\) and another gives \(a\law P\) for an initial segment \(a\) of \(p\). Here \(a\ne\varepsilon\). If all its letters disappear after \(L\), then \(a=L\ne\varepsilon\). With \(L=\varepsilon\), the summand in \(D_i\) is itself nonempty. Apply \eqref{eq:imend} with \(p=am\).
\end{enumerate}
Thus \(L=px^2\) satisfies the induction hypothesis for the next interval of unmarked indices. If the last index is marked, choose a contained \(u\in P\). Its letters all occur in \(L\), and \(L\) ends in a double, so \(Lu\approx L=w^*\). Together with the unmarked final case this gives \(w\law w^*\law P\).
\end{proof}

\section{Suffixes after a first occurrence in \texorpdfstring{\(\A{7812}\)}{A7812}}\label{sec:qualifying}

Put \(B'=\A{7812}\). The \(2\)-row in Table~\ref{tab:twelve} now fixes a product at \(2\) under any number of subsequent \(3\)'s. More precisely, an all-zero word evaluates to \(0\). If its first nonzero factor is \(1\) or \(4\), that factor is the value of the word. If it is \(3\), the value is \(3\) precisely when the whole word consists of \(3\)'s, and \(4\) otherwise. If it is \(2\), the value is \(2\) precisely when everything after it is \(3\), a condition also satisfied by the empty suffix. Otherwise the value is \(4\).

\begin{lemma}\label{lem:suffix}
At the first occurrence of \(x\) in \(w\), let \(D,T\) be the contents of the prefix and suffix strictly before and after this occurrence. Let \(R_x(D)\) be the union of \(\ct(t)\) over the factorizations \(sxt\in P\) with \(\ct(s)\subseteq D\). With \(F(D),U(D)\) as in Section~\ref{sec:cross}, mark the position if
\[
       R_x(D)\cap(D\cup\{x\}\cup F(D))\ne\varnothing
                    \quad\text{or}\quad U(D).
\]
Then \(B'\models w\law P\) if and only if \(P\) has a contained summand, each first \(x\) satisfies \(x\in F(D)\), and \(T\subseteq R_x(D)\) at every unmarked position.
\end{lemma}
\begin{proof}
For necessity of containment and accessibility, use the assignments in the first paragraph of the proof of Lemma~\ref{lem:immediate}. Suppose \(x\) is unmarked and \(T\not\subseteq R=R_x(D)\). Assign \(0\) to \(D\), \(2\) to \(x\), \(3\) to \(R\), and \(1\) elsewhere. No summand is all zero. None has first nonzero factor \(3\), since \(R\cap F(D)=\varnothing\). In every summand whose first nonzero factor is \(2\), only \(3\)'s follow. Thus \(P\) takes the value \(1\) or \(2\), while \(w=4\).

Conversely, if \(P=3\), all its summands consist entirely of \(3\)'s. Accessibility then puts every letter of \(w\) at \(3\). Otherwise \(P\ge1\). If every letter of \(w\) is assigned \(0\), a contained summand also evaluates to \(0\). If its first nonzero factor is \(1\), \(w=1\). For a first nonzero factor \(4\) the witness gives \(4\), and for a first nonzero factor \(3\) it gives \(3\) or \(4\). In either case \(P=4\), since the case \(P=3\) has already been considered.

Let \(x\) be the first nonzero factor and suppose it has value \(2\). Every summand witnessing \(x\) gives \(2\) or \(4\). A letter of \(R\cap(D\cup\{x\})\) takes a value different from \(3\), forcing one such summand to \(4\). For \(q\in R\cap F(D)\), the same holds unless \(q=3\) under the assignment. Its accessible occurrence then yields \(3\) or \(4\) and \(2+3=4\). If \(U(D)\) holds, a summand with value \(0\) gives \(2+0=4\). At an unmarked position, \(T\subseteq R\) ensures that any non-\(3\) suffix factor giving \(w=4\) also gives \(4\) in a summand witnessing \(x\).
\end{proof}

\begin{theorem}\label{thm:qualifying}
The identities below give a basis for \(B'\), and their reversals give a basis for \(\A{8273}\). The schemes have twenty-seven instances, including the tautology \(ax^2\law a+ax^2\).
\begin{align}
 x^3&\approx x^2,& xyx&\approx x^2y^2,&x&\law x^2,\label{eq:qword}\\
 x^2u&\law u+x^2[b],& ax^2&\law a+ax^2[b],\label{eq:qold}\\
 x^2u&\law u+x[d]q[b]+q[c],&
 ax^2&\law a+ax[d]q[b]+aq[c],\label{eq:qcross}\\
 ax^2u&\law au+a+ax[b],&
 ax^2&\law a+ax[b].\label{eq:qinside}
\end{align}
\end{theorem}
\begin{proof}
Write the distinct letters in order of first occurrence. Double \(x\) in this list if some letter strictly after its first occurrence belongs to the content of the prefix ending at that \(x\). The equations in \eqref{eq:qword} give this normal form. Indeed, right multiplication of \(uvu\approx u^2v^2\) by \(v\) and cube contraction give \((uv)^2\approx u^2v^2\). If a previously used \(x\) is appended, take the final segment \(xVx\) starting at its first occurrence. For \(V\ne\varepsilon\) this becomes \(x^2V^2\), and \((uv)^2\approx u^2v^2\) doubles the intervening letters individually. For adjacent copies use \(x^3\approx x^2\). A new letter simply extends the list. This proves the description by induction. In particular, a block of doubled letters absorbs on the right any word on its content.

The two word equations preserve the product values above, and \(x\law x^2\) follows from the diagonal of the table. In \eqref{eq:qold} the summand containing \(x^2\) marks \(x\) by repetition. In \eqref{eq:qcross} the letter \(q\) occurs after \(x\) in one summand and is separately accessible in another. In \eqref{eq:qinside} the summand \(a\) gives \(U(D)\) at \(x\). Where the left side starts with \(a\), the \(a\)-headed summands either mark \(a\) or between them contain all the required suffix letters after \(a\). Containment is supplied by \(u\), \(a\) or \(au\), according to the row. These facts give the remaining identities by Lemma~\ref{lem:suffix}.

Normalize \(w\) in a valid absorption \(w\law P\). A doubled \(x\) cannot be unmarked. Indeed \(x\in T\) would give \(x\in R_x(D)\), itself a reason to mark \(x\). Let \(w^*\) be the word with all marked letters doubled. Then \(w\law w^*\) by \eqref{eq:qword}. It suffices to derive \(w^*\law P\).

Write \(x_1,\ldots,x_n\) for the distinct letters of \(w\) in order, and \(D_i=\{x_1,\ldots,x_{i-1}\}\). Suppose \(t\) is the last unmarked index before \(i\). Split the preceding part of \(w^*\) into \(L\), ending at \(x_t\), and \(B_0\), the following block of doubled letters. Let \(V_t\) consist of the suffixes after occurrences of \(x_t\) in \(P\) accessible over \(D_t\). We prove by induction that
\[
                    LB_0v\law P \qquad(v\in V_t).
\]
Before the first unmarked letter take \(L=\varepsilon\) and use \(P\) itself in place of \(V_t\). At the start \(B_0\) is empty, so the assertion follows from idempotence. Once an unmarked \(t\) occurs, the content of every \(v\in V_t\), and the entire content following \(x_t\) in \(w\), lie in \(R_t=R_{x_t}(D_t)\). This set is disjoint from \(D_t\cup\{x_t\}\cup F(D_t)\), and \(U(D_t)\) fails.

Consider a summand witnessing a later \(x_i\). The first letter outside \(D_t\) in its prefix through \(x_i\) belongs to \(F(D_t)\). Among the possible letters of this prefix, each one outside \(D_t\) other than \(x_t\) belongs to \(R_t\). Since \(F(D_t)\cap R_t=\varnothing\), this first exit is \(x_t\). The suffix after it therefore belongs to \(V_t\). The same argument applies to a summand contained in the whole of \(w\) or in a later prefix, since \(U(D_t)\) fails. There is one further case. A letter \(q\) occurring after \(x_i\) in a witness can lie outside \(\ct(w)\). It still belongs to \(R_t\), because this summand has passed through \(x_t\). In a separate summand witnessing accessibility of \(q\) at \(i\), the prefix through \(q\) must therefore also first leave \(D_t\) at \(x_t\).

Put \(A=LB_0\). Thus a witness for \(x_i\) gives a suffix \(sx_i b\in V_t\) with \(\ct(s)\subseteq\ct(B_0)\). The induction hypothesis and the block identity give \(Ax_i b\law P\). Before the first unmarked index the same reasoning uses the whole witness in \(P\) and \(A=B_0\). If \(i\) is unmarked, the suffixes \(b\) just obtained form \(V_i\). The induction continues with \(L=Ax_i\) and \(B_0=\varepsilon\).

Suppose \(i\) is marked and put \(x=x_i\). We claim \(Ax^2u\law P\) for every \(u\in V_t\), or every \(u\in P\) before the first unmarked index. There are three possibilities.
\begin{enumerate}
\item A suffix after an accessible \(x\) contains a letter of \(D_i\cup\{x\}\). For the corresponding word \(Axb\) the normal form already doubles \(x\), and inserting \(x\) immediately after it changes no other multiplicity in the normal form. Thus \(Ax^2b\law P\). For \(u\ne\varepsilon\) use the first inequality of \eqref{eq:qold} in left context \(A\). For \(u=\varepsilon\) use the second with \(a=A\). An empty \(u\) can occur only after an unmarked letter, so \(A\ne\varepsilon\) in this case.
\item A letter \(q\) in the suffix is separately accessible. The preceding prefix argument gives \(Axdqb\law P\) and \(Aqc\law P\). Apply the first or second inequality of \eqref{eq:qcross} according as \(u\) is nonempty or empty.
\item If \(U(D_i)\) holds, the suffix of a summand contained in \(D_i\) uses only letters of \(B_0\). Hence \(A\law P\). Together with \(Axb\law P\) and \(Au\law P\), the relevant inequality of \eqref{eq:qinside} gives the claim.
\end{enumerate}
In every case \(x^2\) may now be appended to \(B_0\).

Finally choose a contained summand. After the last unmarked \(x_t\), its suffix in \(V_t\) uses only letters of the last block \(B_0\), and \(B_0\) absorbs it. If every letter was marked, the whole summand is absorbed by \(B_0\). We obtain \(w^*=LB_0\law P\), as required.
\end{proof}

\section{Reachability among suffix letters in \texorpdfstring{\(\A{11670}\)}{A11670}}\label{sec:path}

Put \(G=\A{11670}\). For a set of letters \(D\) form the directed graph \(\Gamma_D(P)\) with an edge \(q\to z\) whenever some \(u\in P\) has an occurrence of \(q\) accessible over \(D\) and a strictly later occurrence of \(z\). Write \(R_D(x)\) for the set of vertices reachable from \(x\) by a path of positive length. Thus \(x\in R_D(x)\) requires a directed return to \(x\). In particular, a loop witnessed by a single edge requires two occurrences in a summand.

\begin{lemma}\label{lem:pathcrit}
The identity \(w\law P\) holds in \(G\) if and only if \(P\) has a contained summand and the following holds at each first occurrence of a letter \(x\) in \(w\). Let \(D,T\) be the respective contents before and after that occurrence. Then
\[
 x\text{ is accessible in a summand of }P\text{ over }D,\qquad
 R_D(x)\cap D\ne\varnothing\quad\text{or}\quad T\subseteq R_D(x).
\]
\end{lemma}
\begin{proof}
An all-\(2\) word takes the value \(2\). If the first factor different from \(2\) is \(3\) or \(4\), it determines the value of the word. If it is \(0\) or \(1\), a later factor different from \(1\) gives \(4\). With only \(1\)'s following it the value is \(0\), except that an all-\(1\) word has value \(1\). Indeed, a nonempty prefix of \(2\)'s sends \(1\) to \(0\), and a product at \(0\) stays there precisely under a right factor \(1\).

For containment assign \(2\) on \(\ct(w)\) and \(3\) elsewhere. For accessibility assign \(2\) on \(D\), \(4\) at \(x\), and \(3\) elsewhere. Suppose now that \(R=R_D(x)\) is disjoint from \(D\) and misses a letter of \(T\). Give \(D\) value \(2\) and \(R\) value \(1\). If \(x\notin R\), give it value \(0\). Assign \(3\) to all remaining letters. In a summand the first factor different from \(2\) is accessible over \(D\). If its letter is \(x\) or belongs to \(R\), every letter after it belongs to \(R\) by the definition of an edge and closure under directed paths. Every other first factor different from \(2\) has value \(3\). Thus no summand reaches \(4\), and every summand value is at most \(0\) in \(L_1\). The suffix of \(w\) contains a letter outside \(R\), so \(w=4\).

Conversely assume the two conditions hold and \(P<4\) under an assignment. Containment deals with all-\(2\) and all-\(1\) words. If the first factor of \(w\) different from \(2\) is \(3\) or \(4\), its accessible witness suffices. Suppose it is a letter \(x\) of value \(0\) or \(1\). Any accessible occurrence of \(x\) in \(P\) has only \(2\)'s before it. To keep its summand below \(4\), every later letter must have value \(1\). The same argument applies to an accessible letter of value \(1\) at the start of the next edge. By induction every letter of \(R_D(x)\) has value \(1\), so \(R_D(x)\cap D=\varnothing\). It follows that \(T\subseteq R_D(x)\), and \(w\) evaluates to \(0\) or \(1\). The witness for \(x\) also evaluates to \(0\) or \(1\). The only case still requiring a comparison is \(w=0\) with witness value \(1\). This witness is all \(1\)'s, whereas \(w\) has a letter of value different from \(1\). Accessibility puts that letter in another summand, whose value is in \(\{0,2,3,4\}\). In \(L_1\) its join with \(1\) dominates \(0\).
\end{proof}

\begin{example}\label{ex:path}
For \(G\) the absorption \(xz\law xy+yz+z\) holds. The summand \(z\) is contained and witnesses the final \(z\). For the initial \(x\) the graph \(\Gamma_\varnothing(P)\) contains the path \(x\to y\to z\), with its edges supplied by \(xy\) and \(yz\). This verifies Lemma~\ref{lem:pathcrit}.
\end{example}

\begin{theorem}\label{thm:pathbasis}
Write \(\sq{v}=vv\). In the last identity take \(b=[a]q[c]\) with the same choices at every occurrence of \(b\). The following identities form a basis for \(G\), and their reversals form a basis for \(\A{12469}\).
\begin{align}
 \sq{[a]qq[c]}&\approx \sq{[a]q[c]},&
 \sq{[a]qdq[c]}&\approx \sq{[a]qd[c]},\label{eq:psquares}\\
 \sq{[a]q[c]}\,q&\approx \sq{[a]q[c]},&
                       x&\law xx,\label{eq:pabsorb}\\
 xq[u]&\law q[u]+x[d]q[e],&
                    ax&\law a+ax[b],\label{eq:psplice}\\
 xx[u]&\law x[u]+x[d]x[e],\label{eq:ploop}\\
 \sq{bx}[u]&\law \sq{b}[u]+\sq{b}x[d]q[e].
                                                  \label{eq:pmerge}
\end{align}
\end{theorem}
\begin{proof}

The value of \(\sq{v}\) is \(1\) if \(v\) is all \(1\)'s, \(2\) if it is all \(2\)'s, \(3\) if its first factor different from \(2\) is \(3\), and \(4\) otherwise. This depends only on the ordered list of first occurrences. Hence \eqref{eq:psquares} holds and permits deletion of a later copy of any letter within a square. For \(\sq{b}q\approx\sq{b}\), a square of value \(1\) or \(2\) consists entirely of factors of that value, including \(q\), and \(3,4\) are left zeros. Thus \(\sq{b}\) absorbs on the right every word on \(\ct(b)\). The diagonal entries give \(x\law xx\).

In the first inequality of \eqref{eq:psplice}, the second summand supplies \(x\to q\) and the first supplies \(q\to u\) if \(u\) is present. The first summand is contained, and the displayed occurrences give the accessible letters required by Lemma~\ref{lem:pathcrit}. In the second inequality \(a\) is contained and \(ax[b]\) gives \(a\to x\). In \eqref{eq:ploop} the second summand gives a positive return to \(x\), and \(x[u]\) supplies the edge to \(u\) if present.

For \eqref{eq:pmerge} consider the value of \(\sq b\). If it is \(4\), the first summand suffices. If it is \(3\), both the left side and the first summand equal \(3\). If \(\sq b=1\), all of \(b\) is \(1\). For \(x=1\) the first summand equals the left side, and for every other \(x\) the second summand gives \(4\), since any factor different from \(1\) sends an initial product of \(1\)'s to \(4\). Finally, if \(\sq b=2\), all of \(b\), including \(q\), is \(2\). For \(x=2\) the first summand again equals the left side. For \(x=3\) it dominates the left value \(3\), since a word starting at \(2\) cannot evaluate to \(1\). For \(x=4\) the second summand is \(4\). For \(x=0,1\), the second summand reaches \(0\) at \(x\) and the later \(q=2\) forces \(4\) if no intervening factor has already done so. In display order the eight schemes contribute \(4,4,4,1,8,2,8,32\) identities.

For completeness, fix \(w\) satisfying Lemma~\ref{lem:pathcrit}. We partition \(w\) into consecutive blocks as we pass through its distinct letters in order of first occurrence. Keep \(L\) for the contribution of the completed blocks and \(B=\sq{b}\) for the current block, where \(b\) lists its distinct letters in order. When the current block is empty put \(B=\varepsilon\) and \(D=\varnothing\). Otherwise \(D=\ct(b)\). Let \(r\) be the remaining suffix of \(w\). We keep a finite set \(V\subset X^*\) with these properties.
\begin{enumerate}
\item For each \(v\in V\) we have derived \(LBv\law P\). If \(v=\varepsilon\), this says \(LB\law P\) with \(LB\ne\varepsilon\).
\item Some \(v\in V\) has \(\ct(v)\subseteq D\cup\ct(r)\). For every letter \(y\) whose first occurrence in \(w\) lies in \(r=syt\), the accessibility and path conditions of Lemma~\ref{lem:pathcrit} hold for the strings of \(V\) with preceding support \(D\cup\ct(s)\) and suffix support \(\ct(t)\).
\item Letters in completed blocks do not occur in \(r\). Consequently repetitions before the next first occurrence belong to the current block.
\end{enumerate}
Initially \(L=B=\varepsilon\), \(r=w\) and \(V=P\). These properties then follow from idempotence and Lemma~\ref{lem:pathcrit}. A repetition from the current block can be passed over without changing \(B\) or \(V\). Moving it past the start of \(r\) leaves the sets in (ii) unchanged because its letter belongs to \(D\).

Let \(x\) be the next letter occurring for the first time and write \(R=R_D(x)\) for the graph formed from \(V\). We first prove that for each \(q\in R\) there are words \(d,e\in X^*\) for which \(LBxdqe\law P\) is derivable. If \(sqt\in V\) with \(\ct(s)\subseteq D\), (i) gives \(LBsqt\law P\), and the square absorbs \(s\). When \(B=\varepsilon\) we already have \(s=\varepsilon\). Applying this to a word witnessing the first edge from \(x\) proves the claim for that edge. Suppose the claim holds at \(q\) and an edge out of \(q\) is witnessed by an accessible \(q\) followed by \(t\). We have \(LBqt\law P\). The first law of \eqref{eq:psplice}, in context \(LB\), combines it with \(LBxdqe\law P\) to give \(LBxqt\law P\). Every letter of \(t\) is now a possible endpoint. This induction includes positive-length paths returning to \(x\). The same law with \([u]\) omitted gives \(LBxq\law P\) when \(t\) is empty.

If \(q\in R\cap D\), write \(b=aqc\) and choose a derived \(LBxdqe\law P\) as above. For each \(v\in V\), combine this with \(LBv\law P\) in \eqref{eq:pmerge}, in left context \(L\), to obtain \(L\sq{bx}v\law P\). This adds \(x\) to the current block and leaves \(V\) fixed. For every later first occurrence the set \(D\cup\ct(s)\) in (ii) is the same as before, because \(x\) has moved from \(s\) to the current block. The containment condition is preserved for the same reason.

Suppose instead that \(R\cap D=\varnothing\). By (ii) all letters after this first \(x\) in \(w\) belong to \(R\). Define \(V'\) to consist of the following strings.
\begin{enumerate}
\item the suffix after each occurrence of \(x\) in \(V\) accessible over \(D\),
\item the suffix beginning at each occurrence of \(q\in R\) in \(V\) accessible over \(D\),
\item \(\varepsilon\) if some member of \(V\) is supported in \(D\).
\end{enumerate}
Each string in \(V'\) uses only letters of \(R\). For a suffix of the first kind, absorption of its preceding \(D\)-prefix gives \(LBxt\law P\). For a suffix \(qt\) of the second kind, the path construction gives \(LBxqt\law P\). In the third case we have \(LB\law P\). Here \(LB\ne\varepsilon\), since an empty \(LB\) occurs only initially, when every member of \(P\) is nonempty. Combine \(LB\law P\) with a word \(LBxt\law P\) witnessing \(x\), using the second law of \eqref{eq:psplice}. This gives \(LBx\law P\). Thus \(LBxv'\law P\) for every \(v'\in V'\).

If \(x\notin R\), it cannot recur. Complete the current block and the singleton \(x\), taking \(LBx\) as the new \(L\) and the empty word as \(B\). If \(x\in R\), the construction gives \(LBxdxe\law P\) for a positive return. With each \(LBxv'\law P\), \eqref{eq:ploop} in context \(LB\) gives \(LBx^2v'\law P\). Complete the preceding block, take \(LB\) as the new \(L\), and start the new block with \(\sq{x}\). In either case use \(V'\) in place of \(V\). This gives (i) and (iii), since the letters still to be read belong to \(R\), which is disjoint from \(D\).

To prove (ii) after this change, take a later first occurrence of \(y\). Thus \(y\in R\). In a string of \(V\) witnessing accessibility of \(y\) at that later position, every letter before \(y\) lies in \(D\cup\{x\}\cup R\). Consider the first letter outside \(D\) on the way to \(y\). If it is \(x\), the first kind of suffix in \(V'\) discards just that occurrence and its predecessors. If it belongs to \(R\), the second kind keeps it together with everything after it. The surviving prefix before \(y\) uses only letters of the old preceding support that lie in \(R\). All these letters belong to the new preceding support, so the occurrence of \(y\) is still accessible. In particular \(x\) can survive only when \(x\in R\), in which case it belongs to the new block.

The same argument works for an edge starting at \(q\in R\) on a path out of \(y\) at the later position. If the occurrence of \(q\) to be retained is itself the first exit and \(q=x\), the second kind of suffix retains it because \(x\in R\). All letters following the exit belong to \(R\), by the edge definition at the earlier support \(D\). Hence the occurrence of \(q\), its accessible prefix after restriction, and its edge to the next letter survive in \(V'\). Induction along the path preserves any endpoint in the later suffix, and an endpoint in the old preceding support belongs to the new one whenever it lies in \(R\). Finally, a member witnessing containment either lies entirely in \(D\), giving \(\varepsilon\in V'\), or has a first exit of the same two kinds. Its retained suffix witnesses containment. This proves (ii).

When \(r\) is empty, (ii) supplies \(v\in V\) supported in the final block. Its letters are absorbed in \(LBv\), so \(LB\law P\). By (iii), each block corresponds to a consecutive segment of the original \(w\). A singleton completed with \(x\notin R\) occurs only once in its segment. For every other segment \(s\), \(s\law\sq{s}\) by \eqref{eq:pabsorb}. Deleting repeated letters inside the square by \eqref{eq:psquares} leaves exactly the square of its first-occurrence list. Multiplying these inequalities in block order gives \(w\law LB\), hence \(w\law P\).
\end{proof}

\section{An obstruction to a finite basis for \texorpdfstring{\(\A{13396}\)}{A13396}}\label{sec:private}

Put \(K=\A{13396}\). A word consisting only of \(1\)'s evaluates to \(1\). If its first factor different from \(1\) is \(2\) or \(4\), that is the value of the word. If it is \(3\), the value is \(3\) precisely for an all-\(3\) word, and \(2\) otherwise. If it is \(0\), the value is \(0\) precisely when all subsequent factors are \(3\), and \(4\) otherwise. For example, a preceding \(1\) takes an input \(3\) to the left zero \(2\), while \(0\) is fixed under right multiplication only by \(3\).

\begin{lemma}\label{lem:privcrit}
At each first occurrence of \(x\) in \(w\) let \(D,T\) be the contents before and after it. Put \(R=\bigcup\{\ct(t):sxt\in P,\ \ct(s)\subseteq D\}\). Then \(K\models w\law P\) if and only if \(P\) has a contained summand and, at every such position, \(x\) is accessible over \(D\) in some summand and
\[
                      R\cap(D\cup\{x\})\ne\varnothing
                         \quad\text{or}\quad T\subseteq R .
\]
\end{lemma}
\begin{proof}
For containment assign \(1\) on \(\ct(w)\) and \(2\) elsewhere. For accessibility assign \(1\) on \(D\), \(4\) at \(x\) and \(2\) elsewhere. If both alternatives in the display fail, give \(D\) value \(1\), \(x\) value \(0\), \(R\) value \(3\) and the other letters \(2\). In any summand with first non-\(1\) factor \(0\), that factor is an accessible \(x\) and has only \(3\)'s after it. No summand reaches \(4\), whereas \(w\) does.

Conversely, an all-\(1\) word \(w\) is covered by a contained summand. The value \(P=3\) requires that every summand consist of \(3\)'s, in which case accessibility forces \(w=3\). At every other value of \(P\) a first non-\(1\) factor \(3\) in \(w\) causes no difficulty, since \(w\le2\le P\) in \(L_4\). A first non-\(1\) factor \(2\) or \(4\) has a witness of the same value. Finally, suppose this factor is \(x\) of value \(0\). A summand witnessing \(x\) gives \(0\) or \(4\). If \(w=4\), a letter in the suffix of \(w\) has value different from \(3\). When \(T\subseteq R\), that same letter occurs after an accessible \(x\) in some summand. In the other case a letter of \(R\cap(D\cup\{x\})\), valued \(1\) or \(0\), serves the same purpose. The resulting summand gives \(4\).
\end{proof}

For \(xz\law xy+yz+z\) the set \(R\) at the initial \(x\) is \(\{y\}\). Example~\ref{ex:path} used the edge \(y\to z\) in another summand. Here the preceding lemma excludes the absorption. Indeed, at \((x,y,z)=(0,3,2)\) the left side is \(4\) and the right side is \(0+2+2=0\).

Fix \(n\ge2\) and take \(2n+1\) distinct variables \(h,y_1,\ldots,y_n,z_1,\ldots,z_n\). Define
\begin{equation}
 Y=y_1\cdots y_n,\qquad W=hY,\qquad
 P_n=Y+\sum_{i=1}^n h y_i z_i .\label{eq:pn}
\end{equation}
Each \(z_i\) occurs just once in \(P_n\), in the final position of the summand \(hy_i z_i\).

For \(n=2\), the three summands \(y_1y_2,hy_1z_1,hy_2z_2\) also absorb \(hy_1y_2\). They do not absorb \(hy_2y_1\). After the prefix \(h\) only \(hy_2z_2\) witnesses \(y_2\). Its suffix is \(\{z_2\}\) in content, disjoint from \(\{h,y_2\}\) and missing \(y_1\).

\begin{lemma}[Rigidity]\label{lem:rigid}
The nonempty words \(v\) such that \(K\models v\law P_n\) are exactly the \(n+1\) summands of \(P_n\) and the additional word \(W\).
\end{lemma}
\begin{proof}
Every summand is absorbed by idempotence. For \(W\), the word \(Y\) is contained. The words \(hy_i z_i\) supply all \(y_i\) after the initial \(h\), and \(Y\) itself witnesses each subsequent first occurrence together with all later letters. Lemma~\ref{lem:privcrit} gives \(W\law P_n\).

Now suppose \(v\law P_n\) holds. By accessibility, each letter of \(v\) occurs in \(P_n\), and its first letter must be \(h\) or \(y_1\). If \(z_i\) occurs in \(v\), its only witness is \(hy_i z_i\). Its first occurrence therefore follows both \(h\) and \(y_i\). Since \(R=\varnothing\) there, it is the last letter of \(v\).

If \(v\) begins with \(y_1\), only \(Y\) is accessible there. Its suffix contains neither \(h\) nor \(y_1\) nor any \(z_i\). Thus the contained summand must be \(Y\), forcing every \(y_i\) to occur. Each first occurrence is now witnessed in \(Y\), so the order of first occurrences is \(y_1,\ldots,y_n\). At each \(y_i\) its suffix in \(Y\) misses every preceding letter and \(y_i\) itself. Lemma~\ref{lem:privcrit} therefore forbids a repetition. Hence \(v=Y\).

Suppose \(v\) starts with \(h\). The suffix condition there forbids a second \(h\). If \(v\) contains some \(z_j\), this letter is final and the next letter after \(h\) must be some \(y_i\). For \(i>1\), only \(hy_i z_i\) witnesses this \(y_i\), so only \(z_i\) can follow, and \(v=hy_i z_i\). For \(i=1\), the suffix union at \(y_1\) is \(\{y_2,\ldots,y_n,z_1\}\). Thus the final letter is \(z_1\). If any \(y_k\) intervened, take the first one. Its suffix union is contained in \(\{z_k,y_{k+1},\ldots,y_n\}\), which misses the preceding letters \(h,y_1\), the letter \(y_k\), and the final \(z_1\). This is impossible, so \(v=hy_1z_1\).

It remains to consider \(v\) beginning with \(h\) and containing no \(z_i\). Containment forces all of \(Y\) to occur. After a prefix \(hy_1\cdots y_i\) the suffix union at \(y_i\) is contained in \(\{z_i,y_{i+1},\ldots,y_n\}\), excluding any repetition of a letter already used. If \(y_k\) came next while some earlier \(y_j\) was still missing, \(Y\) would be inaccessible at \(y_k\). Then only \(z_k\) could follow, and the missing \(y_j\) could never occur. The same reasoning excludes a skip at the first letter after \(h\). After \(y_n\) no further letter is possible. Thus \(v=hY=W\).
\end{proof}

\begin{theorem}\label{thm:private}
The semirings \(\A{13396}\) and \(\A{11218}\) are nonfinitely based. For every \(k\ge1\) and \(n>k+2\), the valid identity \(P_n+W\approx P_n\) cannot be derived from the \(k\)-variable identities of \(\A{13396}\) and \(\ax\).
\end{theorem}
\begin{proof}
Take a deduction from \(P_n\) to \(P_n+W\) and expand its terms into sets of words. Every intermediate term has the same term function on \(K\) as \(P_n\), so each of its summands is absorbed by \(P_n\). At the first step introducing \(W\), every word on the preceding side is consequently a summand of \(P_n\), by Lemma~\ref{lem:rigid}. An axiom of \(\ax\) preserves the expanded word set and cannot be this first step.

Let \(f\approx g\) be the identity used at this step, in its direction of use, with substitution \(\sigma\). Distributing the surrounding one-hole context expresses the change as replacements of \(L\sigma(f)R\) by \(L\sigma(g)R\) for fixed words \(L,R\in X^*\), together with summands independent of the hole. Choose a contribution producing \(W\). For each nonempty exterior factor introduce a fresh variable with that factor as its singleton image. Insert these variables on both sides of \(f\approx g\), in their respective positions. The resulting identity has at most \(k+2\) variables. Denote its expanded left side by \(F\) and choose a word \(v\) on the right with \(W\in\sigma(v)\), using the enlarged substitution. We have \(K\models v\law F\), and every word of \(\sigma(F)\) is a summand of \(P_n\).

Every word of \(\sigma(v)\) is linear by Lemma~\ref{lem:rigid}. If \(v\) repeated a variable, using the same nonempty image word at two occurrences would contradict this. Therefore \(v\) is linear. Fix image words giving \(W=hY\) and denote the first letter of \(v\) by \(a\). Its chosen image begins with \(h\) and contains no \(z_i\).

Apply Lemma~\ref{lem:privcrit} to \(v\law F\) on the alphabet of the identity. Accessibility at its initial \(a\) gives a summand of \(F\) beginning with \(a\). In every such summand, fixing the chosen image of the first \(a\) gives only words \(hy_i z_i\) on expansion, since those are the only summands of \(P_n\) beginning with \(h\). Thus the chosen image of \(a\) contains at most one \(y_i\). Also none of these summands of \(F\) can contain a second \(a\), since choosing the same image twice would repeat \(h\). Their suffix union therefore omits \(a\).

The suffix condition at \(a\) now requires each later letter \(b\) of \(v\) to occur after \(a\) in some \(a\)-headed summand of \(F\). In that summand choose the fixed head image and the image of \(b\) used in \(W\). All remaining occurrences can be expanded independently. Every result is still of the form \(hy_i z_i\), so the chosen image of \(b\) also contains at most one \(y_i\). But \(v\) is linear on at most \(k+2\) variables, whereas its chosen expansion \(hY\) contains \(n>k+2\) distinct \(y_i\). The contradiction proves the bound. A finite basis would have a uniform \(k\), and word reversal gives the assertion for \(\A{11218}\).
\end{proof}

\section{Stable powers and bounded products}\label{sec:stable}

Write \(r(a)=a^k\). In the following theorem the integer \(k\) is fixed for a given algebra and the signature has no named constants.
\begin{theorem}\label{thm:power-blocks}
Let \(A\) be a finite ai-semiring of size \(n\) satisfying \(x^k\approx x^{k+1}\). Each of the following two sets of additional identities suffices for finite basability:
\begin{align*}
 &r(xy)\approx r(x)r(y),\qquad r(x)y\law y,\qquad yr(x)\law y;\tag{I}\\
 &r(x)\law x,\quad r(xy)\approx r(r(x)y),\quad r(xy)\approx r(yx),\quad r(xy)\law x.\tag{II}
\end{align*}
\end{theorem}
\begin{proof}
In either case \(r(r(x))=r(x)\) and \(r(x)\) is multiplicatively idempotent. In case (I), writing \(e=r(x),f=r(y)\), we have \(fe=(fe)^2=f(ef)e\le ef\) by the two contraction laws. Interchanging \(e,f\) proves commutation. Also \(r(x)=r(x)x\le x\). On stable powers ordinary multiplication is therefore an idempotent commutative semigroup operation. In case (II) put \(h(x,y)=r(xy)\). The displayed identities give \(h(r(x),y)=h(x,y)=h(x,r(y))\), symmetry, and
\[
 h(h(x,y),z)=r((xy)z)=r(x(yz))=h(x,h(y,z)).
\]
Also \(h(x,x)=r(x)\). On the image of \(r\) this gives an idempotent commutative semigroup operation. In either case, for a nonempty finite set \(S\) of letters, write \(H(S)=r(\prod_{x\in S}x)\), in any order. Repeatedly applying the identities just proved shows that \(r(v)=H(\ct(v))\), and that \(H(S)\le r(x)\le x\) whenever \(x\in S\). For (II) use \(h(x,y)\le x,y\) and take powers for the first inequality. For (I) remove all other stable powers by contraction. Evaluating the power block on a list of all elements of \(A\) therefore gives an additive least element.

Put \(m=n+1\). The following length-dependent identity is valid:
\[
 x_1\cdots x_m\approx \sum_{j=1}^{m}x_1\cdots x_{j-1}r(x_j)x_{j+1}\cdots x_m.
\]
All summands are below the left side. In any assignment two positive prefix products \(p_i,p_j\) coincide, \(i<j\). With \(t=a_{i+1}\cdots a_j\) we have \(p_it=p_i=p_ir(t)\). The idempotent \(r(t)\) is below each \(r(a_s)\) in the segment. Thus the product on that segment with any single factor replaced by its stable power lies between \(r(t)^{j-i}=r(t)\) and \(t\). Multiplication by \(p_i\) makes the bounds equal. This gives a summand with exactly the required value.

We use one more finite identity. In a word on \(2n+1\) distinct letters sum all its variants obtained by replacing the segment in positions \(i+1,\ldots,j\) by its stable power, where \(1\le i<j\le 2n+1\) and \(j-i\ge2\). This sum equals the original word. Each variant lies below the word because \(r(t)\le t\); one of the positive prefix products occurs at least three times, and taking its first and last occurrences gives equality in the corresponding variant.

Treat a word as a sequence of ordinary letters and blocks \(H(S)\). While there are at least \(m\) ordinary letters, group the first \(m\) of them with the blocks that immediately follow them, ending each group before the next ordinary letter. Apply the first length identity to the resulting \(m\) nonempty factors. Each new summand has one fewer ordinary letter, since the stable power of an entire factor is one block. Distribute and induct. In a remaining run of at least \(2n+1\) consecutive blocks apply the second identity to its initial \(2n+1\) blocks. Every new summand has strictly fewer blocks, as the stable power of a segment becomes a single block. We obtain words with at most \(n\) ordinary letters and at most \(K=2n(n+1)\) blocks. On each fixed finite alphabet there are only finitely many of these words.

Consider a valid absorption \(w\law Q\) in these forms. Keep the ordinary letters of \(w\) fixed. Partition the other letters occurring in \(w\) according to their membership in its blocks, replacing each nonempty class \(G\) by a fresh letter \(z_G\). The new left word uses at most \(d=n+2^K-1\) letters. Substitute \(H(G)\) back for \(z_G\). The equalities for stable powers restore the original left blocks, since each block regains exactly its content. On the right the restored polynomial is provably below \(Q\) by \(H(G)\le x\) for each \(x\in G\). This proves that a derivation of the compressed absorption is enough.

By Lemma~\ref{lem:bottom} retain at most \(R=n^{d+1}\) right words in the compressed absorption. Keep every left letter and all ordinary letters of these right words fixed. Partition the remaining letters by membership in the at most \(RK\) right blocks. The resulting absorption uses at most \(d+Rn+2^{RK}-1\) letters. The reverse substitution by the blocks of the classes now restores the selected words exactly, leaving the left word fixed. Take all valid absorptions between the finite set of normal words on this bounded alphabet and their nonempty sums. The two length identities, the stipulated identities and these finitely many absorptions form a basis: reverse the second substitution, enlarge the right sum to its original size, and reverse the first substitution. Every step is an equational consequence in the original signature.
\end{proof}

\begin{corollary}\label{cor:fully-contractive}
Every finite ai-semiring satisfying \(xy\law x\) and \(xy\law y\) is finitely based.
\end{corollary}
\begin{proof}
Take \(k=|A|\). Powers form a descending chain, hence stabilize by exponent \(k\). The element \(r(a)\) is the largest multiplicative idempotent below \(a\): if \(e^2=e\le a\), then \(e=e^k\le a^k\). Likewise \(r(ab)\) is the largest multiplicative idempotent below both \(a\) and \(b\), since any such \(e\) satisfies \(e=e^{2k}\le(ab)^k\). The common idempotent lower bounds of \(a,b\) are therefore also those of \(b,a\) and of \(r(a),b\). This proves (II).
\end{proof}

\section{An obstruction from disjoint pairs}\label{sec:pair-obstruction}

Take \(U=A_{4793}\) and \(V=A_{5072}\). In catalogue order their tables are
\[
\begin{array}{c|cc}
 &+&\cdot\\\hline
 U&\mathtt{00044\;01144\;01234\;44334\;44444}&\mathtt{00234\;01234\;22222\;44444\;44444}\\
 V&\mathtt{00044\;01134\;01234\;43334\;44444}&\mathtt{00244\;01234\;22222\;44444\;44444}.
\end{array}
\]
In each algebra \(2\) is the additive bottom and \(4\) the top. A product of factors in \(\{0,1\}\) is their numerical minimum. At the first factor in \(\{2,3,4\}\), the values \(2,4\) persist until the end. A first \(3\) followed by another factor gives \(4\). If the first \(3\) is final, the value is \(3\) in \(U\); in \(V\) it is \(3\) precisely when there was no preceding \(0\), and \(4\) otherwise. These rules follow by successive right multiplication in the displayed rows.

For either algebra, validity of \(u\law P\) forces
\begin{equation}\label{eq:paircover}
 \ct(u)=\bigcup\{\ct(p):p\in P,\ \ct(p)\subseteq\ct(u)\}.
\end{equation}
Indeed assign a missing letter the value \(0\), the other letters of \(u\) the value \(1\), and all outside letters the value \(2\). Then \(u=0\). Each right word takes value \(1\) or \(2\), both strictly below \(0\).

For \(n\ge2\) take distinct letters \(x_i,y_i\ (1\le i\le n)\) and \(z\). Write \(B_i=\{x_i,y_i\}\), \(C=\bigcup_iB_i\), and put
\[
 w_n=x_1\cdots x_n y_1\cdots y_n,\qquad P_n=\sum_{i=1}^n x_i y_i+w_nz.
\]
Both algebras satisfy \(w_n\law P_n\). A value \(2\) of \(w_n\) requires nothing; a value \(3\) or \(4\) is dominated by \(w_nz=4\). In the remaining case every letter of \(w_n\) takes value \(0\) or \(1\); if a zero occurs its pair has product zero, and otherwise every pair has product one.

\begin{lemma}\label{lem:pairrigid}
Every word absorbed by \(P_n\) is either \(x_i^m y_i\) for some \(m\ge1\), or has content \(C\) or \(C\cup\{z\}\). In the latter two cases each \(x_i\) has appeared before the first \(y_j\).
\end{lemma}
\begin{proof}
Equation~\eqref{eq:paircover} confines the word to \(C\cup\{z\}\). A word containing \(z\) must contain the entire content of the only original summand with \(z\). A word omitting \(z\) has content a union of complete pairs. Suppose that the first \(y_j\) occurs before the final position and some \(x_i\) has not yet appeared. Assign the preceding letters value \(1\), this \(y_j\) value \(3\), and the others value \(2\). The word has value \(4\). Every other pair gives \(2\), the \(j\)th gives at most \(3\), and \(w_nz=2\) because an \(x_i\) is missing from the prefix. This contradicts absorption. If at least two pairs occur the first \(y\) is nonfinal, forcing all \(x_i\), hence all pairs. If only one pair occurs its first \(y_i\) must be final, which gives precisely the asserted shape.
\end{proof}

\begin{theorem}\label{thm:pairnfb}
The semirings \(U,V\) are nonfinitely based. More generally, if one of them belongs to \(\V(A)\) and \(A\models w_n\law P_n\) for every \(n\ge2\), then \(A\) is nonfinitely based.
\end{theorem}
\begin{proof}
Fix a bound \(k\) on the variables in the identities available for deduction and choose \(n>\max\{1,k\}\). Work in the algebra of finite sets of words. Call a word critical if its content is \(C\). Starting from \(P_n\), all words occurring in a valid deduction are absorbed by \(P_n\) in the source algebra. Consider the first step introducing a critical word. In a distributed one-hole context, select one contribution producing this word, with fixed exterior words on either side of the substituted identity. If either exterior word is nonempty, it is free of \(z\) and supplies a fixed letter in every old word of that contribution. Every old word free of \(z\) is of the single-pair form in Lemma~\ref{lem:pairrigid}, and these words can thus cover at most one pair. The new critical word is absorbed by the old contribution alone, contradicting \eqref{eq:paircover}.

Both exterior words may therefore be omitted. Denote the identity, in the direction used, by \(F\approx G\), the substitution by \(\sigma\), and a word of \(G\) giving the new critical word by \(v\). Each letter of \(v\) has at least one \(z\)-free image word. Since \(v\law F\) holds, the words of \(F\) on \(\ct(v)\) cover \(\ct(v)\). Every expansion of any such covering word using only \(z\)-free images is an old, single-pair word. Consequently every \(z\)-free image of every letter of \(v\) is supported in a single pair. Call a source letter wide if its different \(z\)-free images use more than one pair. A covering word containing a wide letter must be that literal alone: at any second position fix an image on one pair and at the wide position choose a different pair. The choices are independent even if the letter repeats. Thus all \(z\)-free images of a wide letter actually have the form \(x_i^m y_i\).

Without a wide letter the chosen expansion of \(v\) covers at most \(k\) pairs. Otherwise consider its first wide occurrence. The preceding source letters collectively use at most \(k-1\) pairs. The selected image at the wide position adds at most one pair and already contains a \(y_i\). At its first \(y\) the critical word is consequently missing some \(x_j\), contrary to Lemma~\ref{lem:pairrigid}. This excludes the first step. A finite basis for \(A\) would consist of identities holding in the source and would have to derive every \(w_n\law P_n\), which is impossible for bounded \(k\).
\end{proof}

For a finite target the all-length hypothesis of Theorem~\ref{thm:pairnfb} has a finite test. On a list \((a_i,b_i)\) store the state \((X,Y,H)=(\prod_i a_i,\prod_i b_i,\sum_i a_i b_i)\). Start with all lists of length two and close under \((X,Y,H)\mapsto(Xa,Yb,H+ab)\) for \(a,b\in A\). Induction identifies the resulting set with all states for lengths at least two. The test is \(XY\le H+(XY)t\) at every state and every \(t\in A\). Inclusion of the source is checked independently by Proposition~\ref{prop:lower-fibres}; reversal uses Lemma~\ref{lem:transfer}.

\begin{corollary}\label{cor:ten-disjoint}
The following five-element algebras are nonfinitely based. The middle column gives the orientation in which the test following Theorem~\ref{thm:pairnfb} is applied.
\begin{center}
\begin{tabular}{@{}llr@{}}\toprule
Catalogue indices & Target tested & Reachable triples per target\\\midrule
4793, 5072, 5073, 14027 & \(A_i\) & 53\\
14533 & \(A_i\) & 54\\
5160, 5167, 5168, 11005 & \(A_i\op\) & 53\\
12482 & \(A_i\op\) & 54\\\bottomrule
\end{tabular}
\end{center}
\end{corollary}
\begin{proof}
For each target \(B\) in the table, use \(A_{4793}\) as the algebra to be generated in Proposition~\ref{prop:lower-fibres}. Its tuple \((0,1,2,3)\) generates: for instance \(3\cdot 3=4\). Starting with the four projections in \(B^{B^4}\), the append-one-letter recursion gives 128 word vectors in each case. Assign each vector its value at \((0,1,2,3)\) in \(A_{4793}\). All \(128\cdot4\) transitions preserve this assignment, and for each of the five principal ideals its union of word vectors under coordinatewise join absorbs only vectors assigned to that ideal. These are the two tests of Proposition~\ref{prop:lower-fibres} in the direction \(A_{4793}\in\V(B)\). The triple recursion above gives the sizes in the last column; checking the five possible final values at each triple gives \(XY\le H+(XY)t\). The two hypotheses of Theorem~\ref{thm:pairnfb} now hold for \(B\), and reversal takes care of the opposite orientations.
\end{proof}

\section{One-sided products and two-element extensions}\label{sec:extension}
\begin{theorem}\label{thm:one-sided}
An ai-semiring whose product depends only on its first argument has a finite basis. The same holds if it depends only on its second argument.
\end{theorem}
\begin{proof}
In the first case \(xy=f(x)\), where \(f(x)=x^2\). Associativity gives \(f^2=f\) and distributivity gives \(f(x+y)=f(x)+f(y)\). Include these three identities, writing each occurrence of \(f\) as a square. Every term reduces to a nonempty join of expressions \(x\) and \(f(x)\). Split each identity between two such terms into absorptions with left member either \(x\) or \(f(x)\). On the right let \(U\) be the letters appearing as literals and \(V\) the letters appearing under \(f\). Hold \(x\) fixed and partition the other letters according as they belong to \(U\) alone, \(V\) alone, or both. Identify each nonempty class with one of its members. Every summand of the resulting right member is still present in the original. The resulting absorption follows from the original by substitution, and gives the original by adjoining the missing terms. Every valid absorption therefore follows from the valid reduced absorptions on four fixed letters. There are finitely many of them. Reversing products proves the second assertion.
\end{proof}

Let \(M\) and \(T\) denote the two-element algebras of Section~\ref{sec:twojoins}. The following transfer of Ma and Yue is used for products with either factor \cite[Theorems~3.1 and~3.4]{MaYue}. We include the constructions because joins of varieties occur in the classification.

\begin{theorem}\label{thm:twofactor}
For an ai-semiring \(A\), each of \(A\times M\) and \(A\times T\) is finitely based if and only if \(A\) is finitely based. No bound on the cardinality of \(A\) is needed.
\end{theorem}
\begin{proof}
The identities of \(A\times M\) are exactly the identities of \(A\) whose sides have the same content. In the notation of absorptions, \(p\law q\) has this property precisely when \(\ct(p)\subseteq\ct(q)\). Here \(p\) and \(q\) may both be polynomials.

Suppose first that \(\Sigma\) is a finite basis for \(A\). Take both orientations of each identity as an absorption, and form the four possible contexts by multiplying by a fresh letter on the left, on the right, on both sides, or on neither. For each resulting \(p\law q\) put \(Y=\ct(p)\setminus\ct(q)\). Retain that inequality if \(Y=\varnothing\). Otherwise take the finitely many laws
\[
 p\law q+\sum_{y\in Y}g_y,\qquad
 g_y\in\{y,\ a_y y,\ y b_y,\ a_y y b_y\},
\]
using distinct fresh parameters. For each \(y\in Y\) also take the two laws
\[
 p[y:=s*t]\law q+p[y:=s]+p[y:=t],\qquad *\in\{+,\cdot\},
\]
with fresh \(s,t\). All these absorptions have the required content and hold in \(A\), since the chosen \(y\)'s are absent from \(q\).

Write \(\Gamma\) for this list. If \(\Gamma\) yields \(\theta(q)\law Z\), and \(\ct(\theta(p))\subseteq\ct(Z)\), it yields \(\theta(p)\law Z\). To prove this, induct on the total number of operation signs in the substituted terms \(\theta(y)\), \(y\in Y\). When they are letters, choose in \(Z\) a word containing each such letter. Its two surrounding factors determine the appropriate choice of \(g_y\). Apply the corresponding inequality and absorb the guards into \(Z\). If \(\theta(y)\) has top operation \(*\), replace that image successively by its two arguments. In both resulting substitutions \(\theta(q)\) is intact and the left-hand content still lies in \(\ct(Z)\). The induction hypothesis on the two substitutions and the displayed law for \(*\) give the claim.

For a regular identity \(P\approx Q\) of \(A\), take any deduction from \(\Sigma\). Replace auxiliary letters outside the common content by one fixed letter of that content. Distributing a one-hole context writes its contribution as summands independent of the hole together with terms \(u[-]v\), where either exterior word may be omitted. Starting with \(P\), inductively derive that each new term is below \(P\): apply the preceding claim with \(Z=P\) to the appropriate context form for each such summand. This yields \(Q\law P\). Reverse the deduction for \(P\law Q\). Thus \(\Gamma\) is a basis for \(A\times M\).

For the converse let \(E\) be the regular identities of \(A\). If \(A\) also satisfies an irregular identity, orient one such identity \(\rho:p\approx q\) so that \(y\) occurs only on the left. Substitute \(z+w\) for \(y\) and \(z\) for the other variables. The left side is a term \(a(z,w)\) involving both letters and the right side involves only \(z\). Hence \(\rho\) implies \(a(z,w)\approx a(z,w')\). For any further irregular identity \(s\approx t\), interchange sides as necessary so that the set of letters \(x_i\) occurring only in \(s\) is nonempty. In the following chain the substitutions for all \(x_i\) are simultaneous:
\[
\begin{aligned}
s&\approx s[x_i:=a(x_i,x_i)]\approx s[x_i:=a(x_i,t)]\\
 &\approx s[x_i:=a(t,x_i)]\approx s[x_i:=a(t,t)]\approx t .
\end{aligned}
\]
The second and fourth steps use the independence of the second argument of \(a\). The first, third and last are regular identities of \(A\), since \(s=t\) makes \(s\) independent of the \(x_i\), and their respective common contents are \(\ct(s),\ct(s)\cup\ct(t),\ct(t)\). Thus \(E\) together with \(\rho\) entails all identities of \(A\). If all its identities are regular, \(E\) already suffices. This proves the converse for \(M\).

For \(T\), call a polynomial proper if it has a summand of length at least two. A proper polynomial is identically one in \(T\). A sum of literals can take value zero there, and its function determines its set of letters. Given \(\Sigma\) as before, replace each \(p\approx q\) by
\[
 p+zt\approx q+zt,\qquad zp\approx zq,\qquad pz\approx qz ,
\]
using fresh letters. These identities hold in \(A\times T\). For any deduction \(u\approx v\) from \(\Sigma\) and any word \(w\) of length at least two, the new identities derive \(u+w\approx v+w\). Indeed, in a distributed context with no exterior word use the first law after splitting \(w\) into two nonempty words. If either exterior word is present use the second or third law and multiply by the other context when present. If \(u,v\) are both proper, choose a long summand of each in turn as \(w\); the two guarded equalities give \(u+v\approx u\) and \(u+v\approx v\). If both consist of literals, validity in \(T\) already identifies their sets of letters. This proves finite basability of \(A\times T\).

Conversely suppose \(A\times T\) has a finite basis. If it satisfies all identities of \(A\), the same basis works. Otherwise choose an identity of \(A\) failing in \(T\). We may suppose \(A\) is nontrivial. An equality of two sums of literals with unequal supports in a semilattice forces triviality, by assigning the top of a two-element subchain to a missing letter and its bottom to the rest. Thus exactly one side is proper. Identifying all variables gives \(x\approx d(x)\) in \(A\) with \(d(x)\) proper. For any identity \(p\approx q\) of \(A\), both sides of \(d(p)\approx d(q)\) are proper, so this identity belongs to \(A\times T\). Adjoining \(x\approx d(x)\) therefore gives a basis for \(A\).
\end{proof}

\section{Quadratic polynomials and directed graphs}\label{sec:qgraph}
For an oriented edge from \(u\) to \(v\) write \(uv\). The following three sets of identities will be used. In the first line each \(e_i\) can be either \(x_{i-1}x_i\) or \(x_ix_{i-1}\).
\[
\begin{array}{ll}
\mathcal B&
x_0x_3\law x_0s+tx_3+e_1+e_2+e_3,\\
\mathcal C&zx\law xy+yz,\qquad xw\law xy+zy+zw,\\
\mathcal D&xt\law xy+zt.
\end{array}
\]

\begin{theorem}\label{thm:qgraph}
Suppose a product \(xyz\) is identically a sum of words of length at most two in \(x,y,z\), or that all products of length three have one common value. In either case, imposing any one of \(\mathcal B,\mathcal C,\mathcal D\) defines a hereditarily finitely based variety relative to the stated reduction and \(\ax\).
\end{theorem}
\begin{proof}
In the first case reduction of a three-letter block in a word, followed by distribution, strictly shortens every resulting word. In the second case \(xyz\approx t^3\) gives a common multiplicative zero, since a four-letter product may be bracketed as a three-factor product. It can be represented by \(t^3\) with \(t\) fresh. Thus a polynomial consists of literal letters, directed edges including possible loops, and, in the second case, a bit recording whether the zero term is present. Adding an edge below the sum leaves the polynomial equivalent to its previous value.

For \(\mathcal B\), let \(F\) be the vertices with outgoing edges and \(I\) those with incoming edges. If there is a three-step undirected walk from \(u\in F\) to \(v\in I\), the existing outgoing and incoming edges supply the guards \(x_0s,tx_3\). Add \(uv\). Closing under this operation preserves \(F,I\). In the closure every odd walk with these endpoint roles yields \(uv\). A one-step walk in the opposite direction can be traversed back and forth to give a three-step walk. In any walk of length at least five, consider its first six vertices. At least one of the three pairs at distance three has one endpoint with an outgoing edge and the other with an incoming edge, in some order. Otherwise all six vertices have just one role and each pair at distance three has the same role. Along adjacent vertices the roles must alternate, a contradiction. Use the chord so obtained to shorten the walk by two steps and apply induction. In a bipartite component the resulting edges are exactly the admissible pairs across the two parts. In a nonbipartite component a detour through an odd cycle provides odd walks between any prescribed endpoints, so its edges are the whole rectangle \(F\times I\) within that component.

For \(\mathcal C\), close under both displayed laws. A component containing a loop at \(u\) becomes a full reflexive directed graph. Indeed the first law supplies the reverse of any edge incident with \(u\), then a loop at its other end. With loops and reverse edges present, the same law joins the ends of every two-step path in both directions. In a loopless component a two-cycle would similarly produce a loop. If there is no two-step directed walk, all nonisolated vertices are sources or sinks. The second law shortens an alternating undirected walk of length three to an edge and hence completes the bipartite rectangle. Otherwise the first law closes a two-step walk to a directed triangle. Starting with its three vertices, construct three nonempty parts \(U,V,W\) having all the edges \(U\to V\to W\to U\). If \(u\in U\) has an edge to a new vertex \(a\), every \(w\in W\) gives \(w\to u\to a\), hence \(a\to w\). Every \(u'\in U\) then gives \(a\to w\to u'\), hence \(u'\to a\). Put \(a\) in \(V\). An edge from \(a\) to \(u\) places \(a\) in \(W\) by the reverse argument. Grow along a spanning tree. Extra reverse edges between parts give two-cycles. An edge within one part also gives a two-cycle by using any vertex of the following part. This proves the complete cyclic tripartite description. For \(\mathcal D\), closing the edges gives one rectangle \(F\times I\) over the entire graph, since any existing first and last endpoints may be paired.

It remains to pass from these descriptions to a basis. Split an identity into absorptions and mark the one or two variables on the left. Distinguish a literal, a loop, a nonloop and the optional common zero. Include an isolated marked point if necessary. A vertex colour records its occurrence as a literal, its marks and its two roles. Record, separately for each numbered part of a component, how many vertices have each colour. There are finitely many component types. Compare components of the same type by coordinatewise comparison of their vectors, and compare lists of components by subsequence matching. Treat isolated points in the same way. Keep the optional zero bit fixed.

Here are the order facts needed in this argument. Every infinite sequence of natural numbers has an infinite nondecreasing subsequence, by selecting a constant subsequence if one value occurs infinitely often and otherwise choosing successively larger values. Successive restrictions prove the same assertion for vectors in any fixed finite dimension. For finite lists over these coloured vectors, suppose an infinite sequence of lists has no earlier list below a later one. Choose such a sequence, at each position taking a shortest list allowing a bad continuation after the chosen prefix. None is empty. Write the lists as \(u_i a_i\) by removing their last entries. Choose infinitely many indices \(i_1<i_2<\cdots\) at which the \(a_i\) have a common type and increase coordinatewise. Replace the tail beginning at \(i_1\) by \(u_{i_1},u_{i_2},\dots\). Minimality gives a comparable pair in the new sequence. A comparison wholly in the original prefix or from it to a shortened list already works in the old sequence. A comparison between two shortened lists extends by their last entries. This contradiction proves the required well-quasi-order. The additional finite choices and the isolated count vector are handled by taking subsequences first.

Such a comparison gives an injective renaming which preserves every literal and edge and maps the marked left word to its counterpart. Adding the remaining right summands derives the larger absorption from the smaller one. The dummy letter of a zero term may always be renamed freshly. If the valid reduced absorptions of a subvariety had no finite upward cover by their members, one could choose them in succession outside the upward closure of all earlier choices, contrary to the order just proved. A finite cover and the environmental axioms therefore form a basis.
\end{proof}

\begin{corollary}\label{cor:qtop}
Every ai-semiring of order at most five satisfying \(y+x^2\approx x^2\) is hereditarily finitely based in the sense that every subvariety it generates has a finite basis.
\end{corollary}
\begin{proof}
The common square \(\omega\) is the additive top and a multiplicative zero, since \(a^2\le a\omega,\omega a\). Suppose \(abc\ne\omega\) and put \(p=ab,\ q=bc,\ t=abc\). A nonzero product cannot equal either of its factors, since multiplying again by that factor on the appropriate side gives a square. Thus \(a\ne b,\ b\ne c\), and \(p,q\) differ from their factors. Also \(p\ne c,\ q\ne a,\ p\ne q\), by \(pc=aq=t\) and \(bc^2=\omega\). The same fixed-point argument rules out \(t=a,c,p,q\); if \(t=b\), substitution gives \(b=a(abc)c=\omega\). Consequently \(a,b,p,q,t,\omega\) alone are six distinct elements. All cubic products equal \(\omega\).

If a law of \(\mathcal B\) fails, fix an assignment and let \(r\ne\omega\) be the value of its right side. The two-sided annihilator \(D\) contains all products, their joins, and hence \(r,\omega\). On the at most three points outside \(D\) put an arrow \(a\to b\) when \(ab\le r\). All endpoints of the five tested edges lie there, and there are no loops. The failed law gives a three-step undirected walk with admissible endpoints and missing endpoint arrow. Its component cannot be bipartite. One part would be a singleton, and either the initial guard or the terminal guard would supply that missing arrow. Thus the graph contains a triangle and \(D=\{r,\omega\}\). Its three vertices are pairwise incomparable, since \(a\le b\) forces \(ab=ba=\omega\). If \(a\) has distinct successors \(b,c\), then \(a(b+c)=r\), although \(b+c\) lies in \(D\). The dual argument excludes two predecessors. The arrows therefore form a directed three-cycle, and products outside its arrows equal \(\omega\). The cycle satisfies \(\mathcal C\). Two consecutive arrows finish a triangle and the rectangle law follows from uniqueness of predecessors and successors. The theorem now applies in either case.
\end{proof}

\begin{lemma}[Changing operations by terms]\label{lem:terms}
Let \(A'\) be obtained from \(A\) by defining its fundamental operations by terms in those of \(A\). Suppose conversely that terms in the operations of \(A'\), on the same set, recover the original operations. Then \(A\) is finitely based if and only if \(A'\) is finitely based.
\end{lemma}
\begin{proof}
Write \(\tau\) for expansion of primed terms into original terms and \(\sigma\) for the reverse translation. If \(E\) is a finite basis for \(A\), use \(\sigma(E)\) and the defining equalities \(f'(\boldsymbol x)\approx\sigma\tau(f'(\boldsymbol x))\), one for each primed basic operation on distinct arguments. They hold by term recovery. For any identity \(u\approx v\) of \(A'\), its expansion holds in \(A\). Translate a deduction from \(E\) back by \(\sigma\). The defining equalities, applied inductively to subterms, identify the two translated ends with \(u,v\). Interchange the translations for the converse. Testing the two recoveries needs only the operation tables and their tuples of arguments.
\end{proof}

\section{Source algebras for variety transfers}\label{sec:classical-sources}
The following source algebras are specified by both operations. An identification with one of their generated varieties can be checked in both directions by Proposition~\ref{prop:lower-fibres}. We use the presentations of A.~Wang and L.~Wang \cite{SourceAlgebras}. Each subsection specifies its own alphabet of element names. A variety is called Cross when it is finitely based, has a finite generating algebra and has only finitely many subvarieties. We denote its relatively free algebra on \(d\) generators by \(\sourcefree_{\mathcal W}(d)\).

\subsection{A five-element Rees matrix semiring}\label{abs-sec:rees}

Let
\[
R=\{\infty\}\cup\{(i,j):i,j\in\{0,1\}\}.
\]
The element $\infty$ is an additive maximum and a multiplicative zero.
On the four remaining elements define
\begin{align*}
(i,j)+(k,\ell)&=(\min(i,k),\min(j,\ell)),\\
(i,j)(k,\ell)&=
\begin{cases}
\infty,&j=k=0,\\
(i,\ell),&j\lor k=1.
\end{cases}
\end{align*}
Multiplication is the combinatorial Rees $0$-matrix multiplication
with sandwich matrix
$\left(\begin{smallmatrix}0&1\\1&1\end{smallmatrix}\right)$.
It is associative: a product is nonzero exactly when all consecutive
coordinate tests pass, and then it retains the two outside coordinates.
Distributivity follows from
\[
j\lor\min(k,k')=\min(j\lor k,j\lor k')
\]
and its dual. Thus $R$ is an additively idempotent semiring.

For a polynomial $P$, let $\ct(P)$ be its content, $\sourceheads(P)$
and $\sourcetails(P)$ the sets of its first and last letters, and
$\sourceedges(P)$ the set of its consecutive ordered pairs of letters.
For a word $w$, write $h(w)$ and $t(w)$ for its first and last letters.

\begin{lemma}\label{abs-rees-criterion}
The semiring $R$ satisfies $w\leq P$ if and only if
\begin{equation}\label{abs-four-data}
\ct(w)\subseteq\ct(P),\quad
h(w)\in\sourceheads(P),\quad
t(w)\in\sourcetails(P),\quad
\sourceedges(w)\subseteq\sourceedges(P).
\end{equation}
\end{lemma}

\begin{proof}
First assign to each variable $x$ a nonzero value $(L_x,R_x)$.
A word $x_1\cdots x_n$ is nonzero exactly when
$R_{x_i}\lor L_{x_{i+1}}=1$ for every consecutive pair.
Its value is then $(L_{x_1},R_{x_n})$. Consequently $P$ is nonzero
exactly when all the tests associated with $\sourceedges(P)$ pass, and
then its value is
\[
\left(\min_{x\in\sourceheads(P)}L_x,\ 
      \min_{y\in\sourcetails(P)}R_y\right).
\]
Conditions \eqref{abs-four-data} imply $w\leq P$. The same argument
works when some variables are assigned $\infty$: if $P$ is
nonzero, all variables in its content are nonzero, and the content
condition applies to $w$.

Conversely, if $x\in\ct(w)\setminus\ct(P)$, assign $\infty$ to
$x$ and $(1,1)$ to every other variable. If $h(w)=x$ is missing from
$\sourceheads(P)$, set $L_x=0$ and every other coordinate equal to $1$.
All words remain nonzero, but the left coordinate of $w$ is $0$
and that of $P$ is $1$. The last-letter condition has the dual
counterassignment.

Finally, if $(x,y)\in\sourceedges(w)\setminus\sourceedges(P)$, set
$R_x=L_y=0$ and every other coordinate equal to $1$. The only
forbidden consecutive pair is $(x,y)$, also when $x=y$.
Hence $w=\infty$ but $P\ne\infty$. Each failed condition therefore
gives a failed absorption.
\end{proof}

Consider the scheme
\begin{equation}\label{abs-splice}
pys\leq pyq+rys,
\end{equation}
where $p,q,r,s$ may each be omitted. Omission means deletion of
that position, not substitution of an empty-word constant.
This gives sixteen ordinary identities, seven of which are trivial.
The other nine are
\begin{equation}\label{abs-nine}
\begin{aligned}
ys&\leq y+rys,       &y&\leq yq+ry,\\
ys&\leq yq+rys,      &pys&\leq py+ys,\\
pys&\leq py+rys,     &py&\leq pyq+y,\\
pys&\leq pyq+ys,     &py&\leq pyq+ry,\\
pys&\leq pyq+rys.&&
\end{aligned}
\end{equation}
The scheme \eqref{abs-splice}, including the optional positions, occurs in the basis
for \(B_0\) in \cite[Theorem~1, equation~(6)]{ShB0} and in the basis
for \(B_2\) in \cite[Theorem~1, equation~(5)]{ShB2}.

\begin{theorem}\label{abs-rees-basis}
The nine identities \eqref{abs-nine} form an identity basis for $R$.
\end{theorem}

\begin{proof}
Every instance of \eqref{abs-splice} satisfies
Lemma~\ref{abs-rees-criterion}: the new word takes its head from the
first word and its tail from the second word, and every consecutive
pair occurs in one of those two words.

Suppose $R\models w\leq P$, where $w=x_1\cdots x_n$.
Choose a word of $P$ beginning with $x_1$. Inductively suppose that
the identities have derived $x_1\cdots x_i u\leq P$.
For $i<n$, choose a word $v x_i x_{i+1}z$ of $P$.
Splicing at the displayed occurrences of $x_i$ gives
\[
x_1\cdots x_i x_{i+1}z
\leq x_1\cdots x_i u+v x_i x_{i+1}z\leq P.
\]
Empty contexts are handled by the corresponding omission instance;
nonempty contexts are ordinary word substitutions. This derives
$wu\leq P$ for some possibly empty suffix $u$.

Choose a word $v x_n$ of $P$. One final splice at $x_n$ removes
the suffix $u$ and derives $w\leq P$. If $n=1$, splice directly
between a word with the desired head and a word with the desired
tail. Thus every valid absorption follows from \eqref{abs-nine}.
Expanding arbitrary identities into sums of words proves completeness.
\end{proof}

\subsection{A rectangular multiplication example}\label{abs-sec:example}

Let $E=\{e,p,q,r,t\}$ have additive order
\[
e<p<r<t,\qquad e<q<r<t,\qquad p\parallel q,
\]
so that $p+q=r$. Define multiplication by
\[
\begin{array}{c|ccccc}
\cdot&e&p&q&r&t\\\hline
e&e&e&q&q&q\\
p&p&p&t&t&t\\
q&e&e&q&q&q\\
r&p&p&t&t&t\\
t&p&p&t&t&t
\end{array}.
\]

\begin{theorem}\label{abs-E-basis}
The algebra $E$ is an additively idempotent semiring, and its
equational theory has the basis
\begin{equation}\label{abs-E-identities}
xyz\approx xz,\qquad x+x^2\approx x^2.
\end{equation}
\end{theorem}

\begin{proof}
Assign coordinate pairs
\[
e:(1,1),\quad p:(0,1),\quad q:(1,0),\quad
r:(0,0),\quad t:(0,0).
\]
Multiplication takes the first coordinate from the left factor and
the second coordinate from the right factor, with outputs
$e,p,q,t$, respectively, for $(1,1),(0,1),(1,0),(0,0)$.
The output retains these coordinates, proving associativity and
$xyz=xz$.

The four possible translation maps, in the order $e,p,q,r,t$, are
\[
(e,e,q,q,q),\quad(p,p,t,t,t),\quad
(e,p,e,p,p),\quad(q,t,q,t,t).
\]
They are order preserving and preserve the only incomparable join,
$p+q=r$. Hence they preserve all binary joins, proving distributivity.
All elements are multiplicatively idempotent except for $r$, whose
square is $t>r$. Thus both identities \eqref{abs-E-identities} hold.

Every term reduces to a sum $P$ of one-letter and two-letter words.
We show that $xy\leq P$ holds, allowing $x=y$, exactly when the
two-letter word $xy$ occurs in $P$. Sufficiency is immediate.
If $x\ne y$ and the word is absent, set $x=p$, $y=q$ and all
other variables equal to $e$. Among all one-letter and two-letter
words, only $xy$ has value $t$; the join of the others is at most $r$.
If the missing word is $x^2$, assign $r$ to $x$ and $e$ to the
other variables. Again only $x^2$ has value $t$. These give the
required counterassignments.

Likewise, $x\leq P$ holds exactly when $P$ contains $x$ or $x^2$.
The second identity in \eqref{abs-E-identities} proves sufficiency.
If both are absent, assign $t$ to $x$ and $e$ to the other variables.
Every monomial of $P$ then has value $e$, $p$ or $q$, so
$P\leq r<t=x$.

Finally, every identity between reduced polynomials splits into
absorptions of the two forms just considered. Each valid absorption
is either immediate or follows from $x\leq x^2$. This proves
completeness.
\end{proof}

\subsection{Adjoining a common absorbing element}\label{abs-sec:regularization}

An identity $s\approx t$ is \emph{regular} if $\ct(s)=\ct(t)$.
Let $B^\infty$ be obtained from an additively idempotent semiring $B$
by adjoining a new element $\infty$ that absorbs both operations.
The semiring axioms hold in $B^\infty$: an assignment using $\infty$
makes both sides of each axiom equal to $\infty$, and all other
assignments are evaluated in $B$.

\begin{theorem}\label{abs-regularization}
Suppose that $B$ is finitely based and that a term $d(x,y)$, containing
both $x$ and $y$, satisfies $B\models d(x,y)\approx x$.
Then $B^\infty$ is finitely based.
\end{theorem}

\begin{proof}
The identities of $B^\infty$ are exactly the regular identities of $B$.
Indeed, a regular identity of $B$ remains valid if any variable is
assigned $\infty$. An identity with different contents fails by assigning
$\infty$ to a variable occurring on only one side and elements of $B$
to all other variables.

Choose a finite basis $\Sigma$ for $B$. In the following identities,
$d$ is an abbreviation for its defining term, not a new operation:
\begin{align}
d(x,x)&\approx x, \label{abs-part1}\\
d(d(x,y),z)&\approx d(x,d(y,z)), \label{abs-part2}\\
d(x,d(y,z))&\approx d(x,d(z,y)), \label{abs-part3}\\
d(f(x,y),z)&\approx f(d(x,z),d(y,z)), \label{abs-part4}\\
d(z,f(x,y))&\approx d(d(z,x),y). \label{abs-part5}
\end{align}
Here $f$ runs through $+$ and multiplication in the last two lines.
Add the finitely many identities
\begin{equation}\label{abs-padded-basis}
d(s,t)\approx d(t,s)\qquad(s\approx t\in\Sigma).
\end{equation}
Every displayed identity is regular and valid in $B$, since $d$ is
the first projection there. Hence all are valid in $B^\infty$.

For a finite nonempty variable set $C=\{x_1,\ldots,x_m\}$, put
\[
D_C(u)=d(\cdots d(d(u,x_1),x_2)\cdots,x_m).
\]
Identities \eqref{abs-part1}--\eqref{abs-part3} allow the factors after the
first one to be permuted and repeated factors to be deleted.
Identity \eqref{abs-part5}, by induction on a term $v$, gives
\[
d(u,v)\approx D_{\ct(v)}(u).
\]
Taking $v=u$ and using \eqref{abs-part1}, we obtain
\begin{equation}\label{abs-self-padding}
D_{\ct(u)}(u)\approx u.
\end{equation}
Also \eqref{abs-part4}, applied successively, gives
\begin{equation}\label{abs-padding-hom}
D_C(f(u,v))\approx f(D_C(u),D_C(v)).
\end{equation}

For a substituted instance $s\sigma\approx t\sigma$ of an identity in
$\Sigma$, assume its variables lie in $C$. Substitute in
\eqref{abs-padded-basis} and then apply $D_C$. Expanding the right arguments
of $d$, and permuting and deleting repeated factors, yields
\[
D_C(s\sigma)\approx D_C(t\sigma).
\]
By \eqref{abs-padding-hom}, this equality is preserved when the instance
occurs inside any term context.

Now let $B^\infty\models u\approx v$, with common content $C$.
Take a finite equational derivation of $u\approx v$ from $\Sigma$
and the semiring axioms. Any auxiliary variable outside $C$ can be
replaced throughout the derivation by a fixed variable of $C$, without
changing the endpoints. Applying $D_C$ to every step, the preceding
argument gives $D_C(u)\approx D_C(v)$. Semiring-axiom steps are preserved
as well. Finally \eqref{abs-self-padding} gives $u\approx v$.
Thus \eqref{abs-part1}--\eqref{abs-padded-basis} form a finite basis for
$B^\infty$.
\end{proof}

\begin{corollary}\label{abs-absorbed-word}
Let $B$ be finitely based. If $B\models w(x,y)\leq x$ for a word
$w$ containing $y$, then $B^\infty$ is finitely based.
\end{corollary}

\begin{proof}
Use $d(x,y)=x+w(x,y)$ in Theorem~\ref{abs-regularization}.
\end{proof}

\subsection{Long products and sums of powers}\label{abs-sec:powers}

Fix an integer $k\geq2$. Let $\mathcal C_k$ be the variety defined by
\begin{align}
x_1\cdots x_k&\approx x_1^k+\cdots+x_k^k, \label{abs-power1}\\
(xy)^k&\approx x^k+y^k, \label{abs-power2}\\
xy+x^k+y^k&\approx xy. \label{abs-power3}
\end{align}

\begin{theorem}\label{abs-power-cross}
The variety $\mathcal C_k$ is Cross. Every subvariety
$\mathcal W\subseteq\mathcal C_k$ has, relative to
\eqref{abs-power1}--\eqref{abs-power3}, a finite basis in at most $k+1$
variables, and
$\mathcal W=\V\bigl(\F_{\mathcal W}(k+1)\bigr)$,
where the indicated free semiring is finite.
\end{theorem}

For a nonempty finite set $C$ of variables write
$z_C=\sum_{x\in C}x^k$. Iterating \eqref{abs-power2} gives
$u^k\approx z_{\ct(u)}$ for every word $u$. Splitting a word into
$k$ nonempty factors and applying \eqref{abs-power1} therefore gives
\begin{equation}\label{abs-power-long}
w\approx z_{\ct(w)}\qquad(|w|\geq k).
\end{equation}
Similarly, splitting $w$ into two nonempty factors and applying
\eqref{abs-power3} gives
\begin{equation}\label{abs-power-pad}
w+z_{\ct(w)}\approx w\qquad(|w|\geq2).
\end{equation}

For a polynomial $P$, let $U$ be its one-letter summands, let $E$
be its words of lengths $2,\ldots,k-1$, and let $C$ be the union of
the contents of all its words of length at least two. Then
\begin{equation}\label{abs-power-normal}
P\approx N(P):=\sum_{u\in U}u+\sum_{w\in E}w+z_C.
\end{equation}
If $C$ is empty, omit the last block; all other empty blocks are
omitted as well. The resulting expression is nonempty. We do not
require these forms to be distinct. Their finiteness on each finite
alphabet proves local finiteness of $\mathcal C_k$.

\begin{lemma}\label{abs-power-regular}
Relative to $\mathcal C_k$, every regular identity is equivalent to
finitely many identities in at most $k+1$ variables, with monomials
of length less than $k$ or of the form $x^k$.
\end{lemma}

\begin{proof}
Put both sides into \eqref{abs-power-normal} and split their equality
into monomial absorptions. Consider one such absorption
$q\leq N(P)$ and put $D=\ct(q)$. Here $q$ is a word of length
less than $k$, or a pure power $x^k$, so $|D|\leq k-1$.
Regularity gives $D\subseteq U\cup C$.

Keep the variables of $D$ fixed. If $C\setminus D$ is nonempty,
replace all its variables by $s^k$, for a fresh variable $s$.
Merge the variables of $U\setminus(C\cup D)$ into one representative
$t$ actually occurring in that set, when it is nonempty. These
two sets are disjoint and contain every variable outside $D$.

Every original nonlinear word involving the first set becomes a
long word. Applying \eqref{abs-power-long} and
$(s^k)^k\approx s^k$, the substituted right side becomes
\begin{equation}\label{abs-power-compressed}
\sum_{u\in U\cap D}u+
\sum_{\substack{w\in E\\\ct(w)\subseteq D}}w+
z_{C\cap D}+s^k+t,
\end{equation}
where empty blocks and variables not introduced are omitted.
The compressed absorption uses at most $|D|+2\leq k+1$ variables.

For the reverse implication, if $s$ was introduced, substitute
$\sum_{v\in C\setminus D}v$ for $s$. Distributivity and
\eqref{abs-power-long} give
$(\sum v)^k\approx\sum v^k$. Every resulting summand in
\eqref{abs-power-compressed} is a summand of $N(P)$, including $t$,
which was chosen from the original one-letter summands. Thus the
compressed absorption implies $q\leq N(P)$.
\end{proof}

\begin{lemma}\label{abs-power-nonregular}
Every subvariety of $\mathcal C_k$ that satisfies a nonregular
identity satisfies
\begin{equation}\label{abs-power-delta}
x^k\approx y^k.
\end{equation}
The subvariety $\mathcal M_k$ defined by \eqref{abs-power-delta}
consists precisely of the semirings in which all $k$-fold products
are a common multiplicative zero $0$ and $0\leq xy$.
Every identity is equivalent over $\mathcal M_k$ to finitely many
identities in at most $k$ variables.
\end{lemma}

\begin{proof}
Take a nonregular identity $P\approx Q$, exchanging sides if necessary
so that $y\in\ct(P)\setminus\ct(Q)$. Replace every other
variable by $x$ and raise both sides to the $k$th power.
Equation \eqref{abs-power-long} and distributivity give
$T^k\approx z_{\ct(T)}$ for any polynomial $T$. Hence the resulting
identity is either $y^k\approx x^k$ or $x^k+y^k\approx x^k$.
In the latter case exchanging $x,y$ also gives \eqref{abs-power-delta}.

Under \eqref{abs-power-delta}, \eqref{abs-power1} makes all $k$-fold
products equal to one element $0$. A $(k+1)$-fold product can be
regrouped into $k$ factors, so $0$ absorbs multiplication.
Equation \eqref{abs-power3} gives $0\leq xy$. The converse is immediate
from \eqref{abs-power1}--\eqref{abs-power3}.

In $\mathcal M_k$, a polynomial is a sum of words of length less
than $k$, possibly together with $0$. Consider a monomial absorption
$q\leq P$, writing $0$ as $x^k$ when it is the left monomial, and
put $D=\ct(q)$. If $P$ contains a nonlinear word or an explicit
$0$, then $0\leq P$. Replace every variable outside $D$ by
$0=x^k$, choosing $x\in D$. Multiplicative absorption eliminates
each word involving such a variable. The resulting right side,
consisting of retained words and possibly $0$, is at most $P$.
The new absorption is therefore equivalent to the original one.

If $P$ contains only one-letter words, merge its variables outside
$D$ into one representative actually present there. The resulting
right side is again a sum of original summands. This yields an
equivalent absorption in at most $|D|+1\leq k$ variables. These
arguments also apply to a left side $0=x^k$, for which $|D|=1$.
\end{proof}

\begin{proof}[Proof of Theorem~\ref{abs-power-cross}]
Fix an alphabet of $k+1$ variables. There are finitely many sums of
words of length less than $k$ and pure $k$th powers on this alphabet.
If a subvariety has only regular identities, Lemma~\ref{abs-power-regular}
defines it using a subset of the resulting finite set of absorptions.
Otherwise Lemma~\ref{abs-power-nonregular} defines it by
\eqref{abs-power-delta} and a subset of another such finite set.
Thus a common finite set of identities suffices to define every
subvariety over $\mathcal C_k$. In particular, all subvarieties
are finitely based and their lattice is finite.

Local finiteness follows from \eqref{abs-power-normal}. The uniform
variable bound, evaluated at the free generators, gives
$\mathcal W=\V(\F_{\mathcal W}(k+1))$. Indeed, if an identity separates the free algebra from the variety, compression gives one in at most \(k+1\) variables. Evaluating it on the free generators forces it throughout the variety. This also proves finite generation.
\end{proof}

\begin{corollary}\label{abs-nil-bottom}
If all $k$-fold products of a semiring $B$ equal its additive
minimum, then both $B$ and $B^\infty$ generate hereditarily
finitely based varieties.
\end{corollary}

\begin{proof}
All three identities \eqref{abs-power1}--\eqref{abs-power3} hold in $B$.
They are regular, so they also hold in $B^\infty$.
Apply Theorem~\ref{abs-power-cross}.
\end{proof}

\subsection{Finite bases for chain semirings}\label{abs-sec:chain}

\subsubsection{The four-element semiring \texorpdfstring{$K$}{K}}

Let $K=\{0,a,e,t\}$ have addition given by maximum on
$0<a<e<t$ and multiplication
\[
\begin{array}{c|cccc}
\cdot&0&a&e&t\\\hline
0&0&0&0&0\\
a&0&0&0&0\\
e&0&a&e&t\\
t&0&a&e&t
\end{array}.
\]
The symbol $0$ names an element of this algebra and is not added
to the signature.

\begin{theorem}\label{abs-chain-basis}
The algebra $K$ is an additively idempotent semiring with identity basis
\begin{align}
xyz&\approx yxz, \label{abs-chain1}\\
x^2y&\approx xy, \label{abs-chain2}\\
xy+y&\approx y, \label{abs-chain3}\\
xy+x^2+y^2&\approx x^2+y^2. \label{abs-chain4}
\end{align}
The multiplicative dual $K^{\mathrm{op}}$ is finitely based as well.
\end{theorem}

The first three identities are the published basis of
$S_{(4,443)}$ in \cite[Proposition~7.18, equations~(7.41)--(7.43)]{FourII}.
We give a direct completeness proof for the additional chain identity
\eqref{abs-chain4}.

\begin{proof}
Call $e,t$ active. A product $uv$ equals $v$ if $u$ is active and
equals $0$ otherwise. A product is active exactly when both its
factors are active. Consequently a triple product equals its last
factor if its first two factors are active and is $0$ otherwise.
This proves associativity and \eqref{abs-chain1}.
Every left or right translation is monotone on the finite chain, so
it preserves binary maxima. Thus both distributive laws hold.

An element and its square have the same activity, proving
\eqref{abs-chain2}. Since $xy$ is either $0$ or $y$, \eqref{abs-chain3} holds.
Finally, if $x\leq y$ then $xy\leq y^2$, while if $y\leq x$ then
$xy\leq x^2$, so \eqref{abs-chain4} holds.

For a nonempty word $w$, write $y$ for its last letter and $S$ for
the content of the prefix obtained by deleting that last letter.
The set $S$ may be empty and may contain $y$. For a nonempty
polynomial $P$, we claim that $K\models w\leq P$ exactly when
the following conditions hold:
\begin{enumerate}
\item[(A)] Some word $u$ of $P$ ends in $y$ and satisfies
$\ct(u)\subseteq S\cup\{y\}$.
\item[(B)] If $y\notin S$, some word $v$ of $P$ either has
$\ct(v)\subseteq S$, or ends in $y$ and has prefix content
contained in $S$.
\end{enumerate}

To prove (A) necessary, assign $t$ to $y$, $e$ to the variables of
$S\setminus\{y\}$, and $0$ to all other variables. The word $w$
has value $t$, and only a word described in (A) can have that value.
If $y\notin S$, instead assign $a$ to $y$, $e$ to $S$, and $0$
elsewhere. Now $w=a$. A word has value at least $a$ exactly when
all its letters lie in $S$, or it ends in $y$ with all prefix
letters in $S$. This proves (B) necessary.

Conversely, if $w$ has value $0$, absorption is automatic.
If it has value $e$ or $t$, every letter of $w$ is active and the
word from (A) has the same value. If it has value $a$, then
$y\notin S$, all prefix letters are active, and the last letter
has value $a$. A word from (B) has either an active value or value
$a$. In each case $w\leq P$, proving the criterion.

It remains to derive every absorption satisfying (A) and (B).
Identities \eqref{abs-chain1} and \eqref{abs-chain2} permit arbitrary
permutations and deletion of repetitions in a word's prefix.
Identity \eqref{abs-chain3} says that adding a nonempty left prefix
can only decrease a word. If $S$ is empty, (B) forces the one-letter
word $y$ to occur in $P$, and there is nothing to prove.
Otherwise choose a word $p$ listing every variable of $S$ exactly
once; then $w\approx py$.

If $y\in S$, the word $u$ from (A) has its whole prefix in $S$,
so $pu\approx py=w$ and $pu\leq u\leq P$.
The same argument applies if (B) supplies a word ending in $y$
with prefix content contained in $S$.
In the remaining case $y\notin S$, and (B) supplies a word $v$
with $\ct(v)\subseteq S$.
If the prefix of $u$ does not contain $y$, the previous argument
again applies. Otherwise
\[
pu\approx py^2,\qquad py^2\leq u,\qquad pvy\approx py=w.
\]
Substitute $v$ for $x$ in \eqref{abs-chain4} and multiply on the left
by $p$. Using \eqref{abs-chain3} gives
\[
w\approx pvy\leq pv^2+py^2\leq v+u\leq P.
\]
Here $pv^2=(pv)v\leq v$ follows directly from \eqref{abs-chain3}.
Thus all valid absorptions follow from the four identities.
Expanding both sides of any identity as sums of words completes
the proof. Reversing every word gives a basis for $K^{\mathrm{op}}$.
\end{proof}

\begin{remark}
Theorem~\ref{abs-chain-basis} proves finite basability, and makes no
claim that every subvariety of $\V(K)$ is finitely based.
Using this theorem to prove finite basability of another semiring
therefore requires equality of the generated varieties.
\end{remark}

\subsubsection{Two four-element and two five-element semirings}

Write $\Sigma_0$ for \eqref{abs-chain1}--\eqref{abs-chain3} and
$\Sigma_K$ for \eqref{abs-chain1}--\eqref{abs-chain4}, and put
\begin{align}
xy&\leq x+zy, \label{abs-chainA}\\
xy&\leq zx+y^2, \label{abs-chainB}\\
xy&\leq yx+y^2. \label{abs-chainW}
\end{align}
On the chain $0<a<b<t$, with maximum as addition, define $L,M$ by
\[
uv=\begin{cases}
v,&u=t,\\
f_X(v),&u\ne t
\end{cases}
\quad(X=L,M),\qquad
\begin{array}{c|cccc}
v&0&a&b&t\\\hline
f_L(v)&0&a&a&a\\
f_M(v)&0&0&b&b
\end{array}.
\]
On $0<a<b<e<t$, again with maximum as addition, define $N,P$ by
\[
u\cdot_Nv=\begin{cases}
v,&u\in\{e,t\},\\
f_N(v),&u\in\{0,a,b\},
\end{cases}
\qquad
u\cdot_Pv=\begin{cases}
0,&u=0,\\
v,&u=t,\\
f_P(v),&u\in\{a,b,e\},
\end{cases}
\]
where
\[
\begin{array}{c|ccccc}
v&0&a&b&e&t\\\hline
f_N(v)&0&0&b&b&b\\
f_P(v)&0&a&a&e&e
\end{array}.
\]

\begin{theorem}\label{abs-chain-extensions}
The algebras $L,M,N,P$ are additively idempotent semirings.
Their respective complete identity bases are
\[
\begin{array}{c|l}
L&\Sigma_0\cup\{\eqref{abs-chainA}\}\\
M&\Sigma_0\cup\{\eqref{abs-chainA},\eqref{abs-chainB}\}\\
N&\Sigma_0\cup\{\eqref{abs-chainB}\}\\
P&\Sigma_0\cup\{\eqref{abs-chainW}\}.
\end{array}
\]
Reversing every word gives finite bases for their multiplicative duals.
\end{theorem}

We first verify the operations and the displayed identities.
For $L,M$, call $t$ active, and for $N$ call $e,t$ active.
In each case $f_X$ is an idempotent order-preserving map below
the identity, and its image consists of inactive elements.
A product is active exactly when both factors are active.
Under either bracketing, a triple product equals its last factor
if its first two factors are active, and its $f_X$-image otherwise.
This proves associativity. In $P$, the left translations
$0,f_P,\mathrm{id}$ compose by taking the smaller map in this list;
the translation of a product is the composite of its factors'
translations. This also proves associativity.
All left and right translations in the four algebras are monotone,
so both distributive laws hold.

These descriptions give \eqref{abs-chain1}, \eqref{abs-chain2}, and
$xy\leq y$. The chain order and monotonicity give
$xy\leq x^2+y^2$, by distinguishing $x\leq y$ and $y\leq x$.
In $L,M$, if $x=t$ then $xy\leq x$, and otherwise
$xy=f_X(y)\leq zy$; hence \eqref{abs-chainA} holds.
For \eqref{abs-chainB} in $M,N$, if $x$ is inactive then
$xy=f_X(y)\leq y^2$. If $x$ is active then $xy=y$; the only
element not equal to its square is $a$, and in that case
$zx\geq b>a$. Finally, in $P$, if $x\ne t$ then $xy\leq y^2$.
If $x=t$, the only case with $y\ne y^2$ is $y=b$, when
$yx=e>b$. This proves \eqref{abs-chainW}.

For the completeness proof, let $w$ end in $y$, let $S$ be its
prefix content, and put $C=S\cup\{y\}$. Let $Q$ be a nonempty
polynomial. Use the following conditions on its summands:
\begin{enumerate}
\item[(C)] Some word has content contained in $C$.
\item[(T)] Some word ends in $y$.
\item[(A)] Some word ends in $y$ and has content contained in $C$.
\item[(L)] If $y\notin S$, some word has content contained in $S$,
or ends in $y$ with prefix content contained in $S$.
\item[(M)] If $y\notin S$, some word ends in a member of $S$,
or ends in $y$ with prefix content contained in $S$.
\item[(P)] If $y\notin S$, some word ends in a member of $S$ and
has content contained in $C$, or ends in $y$ with prefix content
contained in $S$.
\end{enumerate}

\begin{lemma}\label{abs-chain-criteria}
With this notation, the exact criteria for $w\leq Q$ are
\[
\begin{array}{c|l}
L&\text{\rm(C), (T), (L)}\\
M&\text{\rm(C), (T), (M)}\\
N&\text{\rm(A), (M)}\\
P&\text{\rm(A), (P)}.
\end{array}
\]
\end{lemma}

\begin{proof}
For (C) in $L,M$, assign $t$ to $C$ and $0$ elsewhere.
The word $w$ has value $t$, and precisely the words with content
in $C$ can have that value. For (T), assign $a$ to $y$ in $L$,
or $b$ to $y$ in $M$, and $0$ to every other variable.
Exactly the words ending in $y$ can reach the nonzero value of $w$.

For (A) in $N$, assign $t$ to $y$, $e$ to $S\setminus\{y\}$,
and $0$ elsewhere. Then $w=t$ and only words described in (A)
can have value $t$. In $P$, assign $e$ to $y$,
$a$ to $S\setminus\{y\}$, and $0$ elsewhere. Then $w=e$,
and only words described in (A) can reach $e$.

For the remaining necessary conditions, assume $y\notin S$,
assign $0$ outside $C$, and use the following assignments:
\[
\begin{array}{c|ccc}
\text{algebra}&y&\text{each variable in }S&w\\\hline
L&b&t&b\\
M&a&t&a\\
N&a&e&a\\
P&b&t&b
\end{array}.
\]
In the first row, a word reaches $b$ exactly when all its letters
lie in $S$, or it ends in $y$ with prefix letters in $S$.
In the middle two rows, every word ending in $S$ reaches at least
$b$; a word ending in $y$ reaches $a$ exactly when its prefix
lies in $S$; all other words are $0$.
In the last row, any word involving an outside variable is $0$.
Among the remaining words, those ending in $S$ reach at least $e$,
while a word ending in $y$ reaches $b$ exactly when its prefix lies
in $S$. These are precisely (L), (M), and (P), respectively.

For sufficiency in $L$, values $w=0,t$ are handled by order
and (C). If $w=a$, the last letter has value at least $a$,
and (T) gives a word of value at least $a$.
If $w=b$, every prefix letter is $t$ and $y=b$, so $y\notin S$
and (L) gives a word of value $t$ or $b$.
In $M$, the same argument handles $0,t$; (T) handles $w=b$,
and (M) handles $w=a$.

In $N$, if $w$ is active then (A) supplies a word of the same
value, since all letters of $w$ are active. If $w=b$, its last
letter has value at least $b$, and the word from (A) also has
value at least $b$. If $w=a$, all prefix letters are active,
$y=a$, and (M) supplies a word of value at least $b$, or of
value $a$. The value $0$ needs no condition.

In $P$, (A) handles $w=t$. If $w=e$, all its letters are
nonzero and its last letter is at least $e$, so the word from
(A) is at least $e$. If $w=a$, all its letters are nonzero,
and the same word is also nonzero. If $w=b$, all prefix
letters are $t$ and $y=b$, so (P) supplies a word of value
at least $e$, or of value $b$. Again $w=0$ is automatic.
\end{proof}

\begin{proof}[Completeness in Theorem~\ref{abs-chain-extensions}]
Each listed basis implies \eqref{abs-chain4}.
For \eqref{abs-chainA}, substitute $(x,y,z)=(x^2,y,y)$ and use
\eqref{abs-chain2}. For \eqref{abs-chainB}, set $z=x$.
For \eqref{abs-chainW}, use \eqref{abs-chain2} and \eqref{abs-chain3} to obtain
\[
xy\approx x^2y\leq yx^2+y^2\leq x^2+y^2.
\]
Thus $\Sigma_K$ is available in all four cases.
We use prefix permutations and deletion of repetitions from
\eqref{abs-chain1}, \eqref{abs-chain2}, and the inequality $uv\leq v$.
First we derive $w^2\leq Q$.
For $N,P$, the word $u$ from (A) satisfies
$w^2u\approx w^2$, so $w^2\leq u\leq Q$.
For $L,M$, choose $u$ from (C) and $v$ from (T).
Prefix reduction gives
\[
uw^2\approx w^2,\qquad vw^2\approx w^2v\leq v.
\]
Applying \eqref{abs-chainA} with $(x,y,z)=(u,w^2,v)$ yields
\[
w^2\approx uw^2\leq u+vw^2\leq u+v\leq Q.
\]

If $S$ is empty, the last condition in each criterion forces the
one-letter word $y$ to occur in $Q$. If $y\in S$, prefix
reduction gives $w\approx w^2$. Both cases are settled.
Otherwise let $p$ list all members of $S$ once, so $w\approx py$.
If $Q$ has a word $v$ ending in $y$ with prefix content in $S$,
then $pv\approx w$ and $w\leq v\leq Q$.

In the remaining case for $L$, condition (L) supplies a nonempty
word $v$ with $\ct(v)\subseteq S$. Hence $pvy\approx w$, and
\eqref{abs-chain4} gives
\[
w\approx pvy\leq pv^2+py^2\leq v+w^2\leq Q.
\]
For $M,N$, the remaining case supplies a word $r$ ending in
some $s\in S$. If $r=zs$ for a nonempty word $z$, apply
\eqref{abs-chainB} as $sy\leq zs+y^2$ and multiply on the left by
$p$. The result is
\[
w\approx psy\leq pr+py^2\leq r+w^2\leq Q.
\]
If $r=s$, use $sy\leq s^2+y^2\leq s+y^2$ instead, giving
the same result without an empty-word substitution.

Finally, in $P$, the remaining witness $r$ ends in $s\in S$
and satisfies $\ct(r)\subseteq C$. Apply \eqref{abs-chainW}
and left-multiply by $p$:
\[
w\approx psy\leq pys+py^2\leq r+w^2\leq Q.
\]
Here $pys\approx pyr\leq r$, since the prefix of $py$ contains
all of $C$ and $\ct(r)\subseteq C$. Also $py^2\approx w^2$.
Thus every absorption in Lemma~\ref{abs-chain-criteria} follows from
the indicated basis. Expansion into sums of words reduces every
identity to such absorptions, proving completeness.
\end{proof}

The theorem asserts finite basability of these individual varieties.
No hereditary finite-basis conclusion is needed or asserted.

\begin{corollary}\label{abs-chain-relations}
The chain semirings above satisfy
\[
\V(K)=\V(L)\vee\V(N)=\V(L)\vee\V(P),\qquad
\V(M)=\V(L)\cap\V(N),
\]
and $\V(N)\subsetneq\V(P)\subsetneq\V(K)$.
\end{corollary}

\begin{proof}
The criterion for $K$ in Theorem~\ref{abs-chain-basis} is
\textup{(A), (L)} in the notation of Lemma~\ref{abs-chain-criteria}.
Condition \textup{(A)} implies \textup{(C), (T)}, and
\textup{(L)} implies \textup{(M)}.
Thus an absorption holds in both $L$ and $N$ exactly when it
holds in $K$. Splitting identities into absorptions proves the
first join equality. The displayed bases give the meet equality.
Identity \eqref{abs-chainB}, with $z=y$, implies \eqref{abs-chainW};
hence $\V(N)\subseteq\V(P)\subseteq\V(K)$, and the other
join equality follows. The first inclusion is strict because
$P$ fails \eqref{abs-chainB} at $(x,y,z)=(t,b,0)$.
The second is strict because $K$ fails \eqref{abs-chainW} at
$(x,y)=(e,a)$.
\end{proof}

\subsection{Commutative semirings with an additive square map}\label{abs-sec:squares}

Let $\Sigma_{\mathrm s}$ consist of
\begin{align}
xy&\approx yx, \label{abs-sq1}\\
x^2y&\approx xy, \label{abs-sq2}\\
xy&\leq x^2+y^2. \label{abs-sq3}
\end{align}
We use four additional identities:
\begin{align}
x&\leq x^2, \tag{I}\label{abs-sqI}\\
xy&\leq x^2+y, \tag{D}\label{abs-sqD}\\
xy&\leq x^2+yz, \tag{E}\label{abs-sqE}\\
x^2&\leq x+xy. \tag{H}\label{abs-sqH}
\end{align}
The letters in parentheses name identities only within this section.

\begin{lemma}\label{abs-square-retraction}
In a semiring satisfying $\Sigma_{\mathrm s}$, the map
$\rho(x)=x^2$ is an idempotent endomorphism and
$xy=\rho(x)\rho(y)$. Its image is multiplicatively idempotent.
Conversely, let a commutative, multiplicatively idempotent semiring
$B$ be a join-subsemilattice of a semilattice $A$, and let
$\rho:A\to B$ be a join-preserving retraction. Then
\[
xy:=\rho(x)\cdot_B\rho(y)
\]
makes $A$ a semiring satisfying $\Sigma_{\mathrm s}$, with $x^2=\rho(x)$.
\end{lemma}

\begin{proof}
The first two identities give $x^3\approx x^2$ and $x^2y^2\approx xy$.
The third gives $(x+y)^2\approx x^2+y^2$. These prove the assertions
about $\rho$ and its image.
Conversely, associativity and distributivity follow from the
corresponding laws in $B$ and preservation of joins by $\rho$.
In $B$, multiplicative idempotence and expansion of $(x+y)^2$
give $xy\leq x+y$. This proves \eqref{abs-sq3} for $A$.
\end{proof}

We specify eight examples by this construction. In $B_+,B_-$,
addition is maximum on $0<e<t$, and in $B_*$ it is maximum on
$0<b<c<d$. Their multiplications are
\[
\begin{array}{c|ccc}
\cdot_+&0&e&t\\\hline
0&0&0&0\\e&0&e&t\\t&0&t&t
\end{array}
\qquad
\begin{array}{c|ccc}
\cdot_-&0&e&t\\\hline
0&0&0&t\\e&0&e&t\\t&t&t&t
\end{array}
\qquad
\begin{array}{c|cccc}
\cdot_*&0&b&c&d\\\hline
0&0&b&b&b\\b&b&b&b&b\\c&b&b&c&d\\d&b&b&d&d
\end{array}.
\]
Each multiplication is meet in a semilattice order. For $B_*$,
that order has covers $b<0$ and $b<d<c$.
The translations are monotone in the indicated additive chain,
so all three structures are semirings.

In Table~\ref{abs-tab:squares}, extend the indicated $B$ by the
additional elements shown. The map $\rho$ fixes $B$; its other
values are in the last column. Addition is join in the stated order.
In each nonchain row, the two elements listed together at the same
level are incomparable. In particular, the last row is the diamond
above $a$ with covers $0<a$, $a<e,b$, and $e,b<t$.
Each map preserves joins: on chains this follows from monotonicity;
in the other rows the only additional joins to check are
$a+e=t$, $a+b=c$, and $e+b=t$, respectively.
Multiplication is then defined by Lemma~\ref{abs-square-retraction}.

\begin{table}[ht]
\centering
\caption{The eight square-retraction semirings.}\label{abs-tab:squares}
\begin{tabular}{ccll}
\toprule
Semiring&Image&Additive order&Additional values of $\rho$\\
\midrule
$A_{\varnothing}$&$B_+$&$0<e<a<t<u$&$a,u\mapsto t$\\
$A_I$&$B_+$&$0<e<a<t$&$a\mapsto t$\\
$A_D$&$B_+$&$0<a,e<t<u$&$a,u\mapsto t$\\
$A_{DE}$&$B_-$&$0<a<e<t<u$&$a\mapsto e,\ u\mapsto t$\\
$A_{IDE}$&$B_-$&$0<a<e<t$&$a\mapsto e$\\
$A_{EH}$&$B_-$&$0<a,e<t<u$&$a,u\mapsto t$\\
$A_{IH}$&$B_*$&$0<a,b<c<d$&$a\mapsto c$\\
$A_{IE}$&$B_-$&$0<a<e,b<t$&$a\mapsto e,\ b\mapsto t$\\
\bottomrule
\end{tabular}
\end{table}

\begin{theorem}\label{abs-square-bases}
For $\Lambda\in\{\varnothing,I,D,DE,IDE,EH,IH,IE\}$,
a complete identity basis for $A_\Lambda$ is $\Sigma_{\mathrm s}$
together with the identities named by $\Lambda$.
If $\Lambda$ contains $D$ or $E$, identity \eqref{abs-sq3} can be omitted.
Thus the eight bases have, respectively, $3,4,3,4,5,4,5,4$ identities.
\end{theorem}

We verify the additional identities first.
Identity \eqref{abs-sqI} holds in its four cases because $\rho(x)\geq x$.
For \eqref{abs-sqD} in $A_D$, values $\rho(x)=0,t$ are immediate.
If $\rho(x)=e$, then $xy=\rho(y)\leq e+y$: the only new case is
$y=a$, where $e+a=t=\rho(a)$, while $\rho(u)\leq u$.
For $A_{DE},A_{IDE}$, the cases $\rho(x)=t,e$ follow similarly.
If $\rho(x)=0$, a nonzero product requires $\rho(y)=t$,
and such a $y$ is $t$ or $u$, so $xy\leq y$.

Identity \eqref{abs-sqE} holds for every extension of $B_-$ in the table:
if $\rho(x)=t$ or $\rho(y)=t$, it is immediate from $x^2=t$
or $yz=t$, respectively. Otherwise $xy\leq x^2$.
For \eqref{abs-sqH} in $A_{EH}$, only $x=a,u$ need checking,
and then $x^2=xy=t$. In $A_{IH}$ only $x=a$ needs checking.
If $\rho(y)=0,b$, then $ay=b$ and $a+b=c=a^2$; otherwise
$ay=c,d$, which absorbs $a^2$.
Finally, \eqref{abs-sqD} with $x^2,y^2$ substituted for $x,y$
implies \eqref{abs-sq3}; \eqref{abs-sqE} implies it by setting $z=y$.
This justifies the shorter bases.

Under \eqref{abs-sq1} and \eqref{abs-sq2}, a nonlinear word (a word of
length at least two) is determined up to provable equality by its
nonempty content: use $x^2$ for content $\{x\}$, and the product
of the distinct letters otherwise. This does not assert that
different contents must remain distinct in a particular semiring.
For a nonempty polynomial $Q$, let $U$ be its one-letter
summands, let $\mathcal E$ be the family of contents of its
nonlinear summands, and put
\[
G=\bigcup\mathcal E,\qquad
\mathcal E_C=\{D\in\mathcal E:D\subseteq C\},\qquad
R_C=\bigcup\mathcal E_C.
\]
The union of an empty family is empty; no empty word is introduced.

\begin{lemma}\label{abs-square-criteria}
In $A_\Lambda$,
\[
x\leq Q\quad\Longleftrightarrow\quad
x\in U\ \text{or}\ \bigl(I\in\Lambda\ \text{and}\ \{x\}\in\mathcal E\bigr).
\]
For a nonlinear word $w$ with content $C$, the exact conditions
for $w\leq Q$ are
\[
\begin{array}{c|l}
\Lambda&\text{condition}\\\hline
\varnothing,I&C\subseteq R_C\\
D&\mathcal E_C\ne\varnothing,\quad C\subseteq U\cup R_C\\
DE,IDE&\mathcal E_C\ne\varnothing,\quad C\subseteq U\cup G\\
EH&C\subseteq G,\quad
(\mathcal E_C\ne\varnothing\ \text{or}\ U\cap C\ne\varnothing)\\
IH&C\subseteq R_C\cup(U\cap G)\\
IE&\mathcal E_C\ne\varnothing,\quad C\subseteq G.
\end{array}
\]
\end{lemma}

\begin{proof}[Necessity]
For a one-letter left side with $I\notin\Lambda$, if $x\notin U$,
assign $u$ to $x$ and $0$ elsewhere. All nonlinear values lie in
$B$ and are at most $t<u$, so absorption fails.
If $I\in\Lambda$ and both possible witnesses are missing,
assign $a$ to $x$ and $0$ elsewhere. The right side is $0$,
except in $A_{IH}$ where it is $0$ or $b$; in that case
$a\not\leq b$ as well.

For nonlinear words, the following assignments handle every
failed condition. Variables outside $C$ are assigned $0$ unless
another assignment is explicitly specified.
\begin{enumerate}
\item For $\varnothing,I$, choose $x\in C\setminus R_C$,
assign $a$ to $x$ and $e$ to $C\setminus\{x\}$.
The left side is $t$. A nonlinear right word involving an outside
variable is $0$; one inside $C$ cannot involve $x$, so is $e$.
One-letter values are at most $a<t$.
\item For $D$, if $\mathcal E_C$ is empty, assign $a$ to $C$:
the left side is $t$ and the right side is at most $a$.
If $x\in C\setminus(U\cup R_C)$, assign $t$ to $x$ and $e$
to the other members of $C$. The left side is $t$ and the
right side is at most $e$.
\item For $DE,IDE$, if $\mathcal E_C$ is empty, assign $a$
to $C$: the left side is $e$, and the right side is at most $a$.
If $x\in C\setminus(U\cup G)$, assign $t$ to $x$ and $0$
to every other variable. The two sides have values $t$ and $0$.
\item For $EH$, if $x\in C\setminus G$, assign $a$ to $x$
and $0$ everywhere else. The left side is $t$ and the right side
is at most $a$. If $\mathcal E_C$ and $U\cap C$ are both empty,
assign $e$ to $C$. The two sides have values $e$ and $0$.
\item For $IH$, choose
$x\in C\setminus(R_C\cup(U\cap G))$.
If $x\notin G$, assign $a$ to $x$ and $0$ everywhere else.
The left side is $c$ for $C=\{x\}$ and $b$ otherwise;
the right side is $a$ or $0$, neither of which absorbs it.
Otherwise $x\notin U$; assign $d$ to $x$ and $c$ to
$C\setminus\{x\}$. Every nonlinear right word involving $x$
also involves an outside variable and has value $b$.
The right side is at most $c<d=w$.
\item For $IE$, if $x\in C\setminus G$, assign $b$ to $x$
and $0$ everywhere else. The left side is $t$ and the right side
is at most $b$. If $\mathcal E_C$ is empty, assign $a$ to $C$;
the two sides have values $e$ and at most $a$.\qedhere
\end{enumerate}
\end{proof}

\begin{proof}[Sufficiency and completeness in Theorem~\ref{abs-square-bases}]
Every nonlinear word $v$ satisfies $v^2\approx v$.
If two nonlinear words $v,r$ have already been absorbed by $Q$,
then \eqref{abs-sq3} gives
\[
vr\leq v^2+r^2\approx v+r\leq Q.
\]
Thus their contents may be merged. The additional laws permit
three further operations. For $x\in U$, identity \eqref{abs-sqD} gives
\[
vx\leq v^2+x\approx v+x\leq Q.
\]
For $x\in G$, choose a nonlinear summand $r=xz$ of $Q$,
where $z$ is a nonempty word. Identity \eqref{abs-sqE} gives
\[
vx\leq v^2+xz\approx v+r\leq Q.
\]
Finally, for $x\in U\cap G$, identity \eqref{abs-sqH} yields a
new starting word
\[
x^2\leq x+xz=x+r\leq Q.
\]

For $\varnothing,I$, merge the words indexed by $\mathcal E_C$.
For $D$, merge those words and add missing letters from $U$.
For $DE,IDE$, start with one word indexed by $\mathcal E_C$
and add missing letters from $U$ or $G$.
For $EH$, start either with such a word or with an $x^2$
obtained from $x\in U\cap C\subseteq U\cap G$, then add letters
from $G$. For $IH$, merge the words indexed by $\mathcal E_C$
with all needed $x^2$ for $x\in C\cap U\cap G$.
For $IE$, start with a word indexed by $\mathcal E_C$ and
add the missing letters from $G$.
Each construction stays inside $C$ and ends with content exactly
$C$, hence with a word equal to $w$.
The one-letter criterion follows directly from \eqref{abs-sqI}
or the occurrence of $x$ in $U$.
All substitutions use nonempty words.

This proves sufficiency in Lemma~\ref{abs-square-criteria}.
Expanding arbitrary identities into monomial absorptions now
proves completeness of all eight bases.
\end{proof}

\subsection{Two bases from connected content hypergraphs}\label{abs-sec:connected}

The square map need not preserve addition even when
\eqref{abs-sq1}, \eqref{abs-sq2}, and \eqref{abs-sqI} hold.
We give two finite bases in this setting.

Let $F_5=\{0,a,b,c,t\}$ have additive order
$0<a,b<c<t$, with $a\parallel b$.
Let $0$ be a multiplicative zero, put $a^2=a$, $b^2=b$, $t^2=t$,
and set every other product of two nonzero elements equal to $t$.
Let $G_5=\{a,b,c,d,t\}$ have additive order
$a,b<c<d<t$, with $a\parallel b$.
Put $a^2=a$, $b^2=b$, $t^2=t$, and set every other product equal
to $d$. In particular, $d$ is a multiplicative zero.

\begin{theorem}\label{abs-connected-bases}
The algebras $F_5,G_5$ are additively idempotent semirings.
A complete basis for $F_5$ is
\[
xy\approx yx,\qquad x^2y\approx xy,\qquad x\leq x^2.
\]
A complete basis for $G_5$ consists of these three identities and
\begin{equation}\label{abs-connectedJ}
xy\leq x+xyz.
\end{equation}
Moreover, $\V(G_5)\subsetneq\V(F_5)$.
\end{theorem}

Associativity follows in each case by retracting $c$ onto its
square, respectively $t$ and $d$, and multiplying in the resulting
semilattice. Each translation is order-preserving and preserves
the only incomparable join $a+b=c$, so both distributive laws hold.
The first three displayed identities follow directly.
In $G_5$, if $xy$ is $a,b$, or $t$, then $x=y=xy$.
Otherwise $xy=d=xyz$. This proves \eqref{abs-connectedJ}.
The square map fails to preserve $a+b$ in both algebras:
$a^2+b^2=c$, whereas $(a+b)^2$ is $t$ or $d$.

Use the notation $U,\mathcal E,\mathcal E_C$ from the previous
section. A content hypergraph has as vertices the union of its
nonempty edges. Its connected components are obtained by merging
intersecting edges. A singleton edge contributes a one-vertex
component; the empty family has no vertices.

\begin{lemma}\label{abs-connected-criteria}
In either $F_5$ or $G_5$,
$x\leq Q$ holds exactly when $x\in U$ or $\{x\}\in\mathcal E$.
For a nonlinear word $w$ with content $C$,
\begin{enumerate}
\item $F_5\models w\leq Q$ exactly when the hypergraph
$\mathcal E_C$ covers $C$ and is connected;
\item $G_5\models w\leq Q$ exactly when $C$ is contained in one
connected component of $\mathcal E$, and $Q$ has a word $u$
with $\ct(u)\subseteq C$. The word $u$ may be linear.
\end{enumerate}
\end{lemma}

\begin{proof}[Necessity]
If the one-letter criterion fails, in $F_5$ assign $a$ to $x$
and $0$ elsewhere. Every right word is $0$.
In $G_5$, assign $t$ to $x$ and $a$ elsewhere.
Every right word has value at most $d<t$.

For $F_5$, if $x\in C$ is not covered by $\mathcal E_C$,
assign $c$ to $x$, $a$ to the other variables of $C$, and
$0$ outside $C$. The left word is $t$, while the right side
is at most $c$. If $\mathcal E_C$ covers $C$ but is disconnected,
assign $a$ to one component and $b$ to all the others,
and assign $0$ outside $C$. Again the left word is $t$.
Every right word is $0,a$, or $b$, so their join is at most $c$.

For $G_5$, if $x\in C$ is missing from $\bigcup\mathcal E$,
assign $c$ to $x$ and $a$ everywhere else. The left word is $d$
and the right side is at most $c$.
If $C$ meets two components, assign $a$ to one of them and
$b$ to every other component; assign $a$ to variables not in
the hypergraph. The left word is $d$, while every right word
is $a$ or $b$, giving a join at most $c$.
Finally, if every right word involves a variable outside $C$,
assign $t$ to $C$ and $a$ elsewhere. The left word is $t$,
and every right word is $a$ or $d$.
\end{proof}

\begin{proof}[Sufficiency and completeness in Theorem~\ref{abs-connected-bases}]
The first two basis identities, after expanding
$(p+q)^2r\approx(p+q)r$, give
\[
pqr\leq pr+qr.
\]
If nonlinear words $u,v$ have a common letter $z$, write
$u=pz$, $v=qz$ with nonempty $p,q$. Then
\[
uv\approx pqz\leq pz+qz=u+v.
\]
Consequently words belonging to intersecting hyperedges can be
merged without increasing their sum. Ordering the edges of a
connected component so that each new edge meets the union of
earlier ones shows that their sum absorbs a nonlinear word whose
content is the whole component.

This proves sufficiency for $F_5$. For $G_5$, let $D$ be the
component containing $C$. The same merging argument gives
a nonlinear word $v$ of content $D$ with $v\leq Q$.
Choose the witness $u$ with $\ct(u)\subseteq C$.
Apply \eqref{abs-connectedJ} with $(x,y,z)=(u,w,v)$:
\[
w\approx uw\leq u+uwv\approx u+v\leq Q.
\]
Both reductions follow from content normalization and
$\ct(u)\subseteq C\subseteq D$.
The one-letter criterion follows from $x\leq x^2$.
Thus all valid absorptions, and hence all identities, follow
from the displayed bases. The inclusion of varieties follows
from their bases. It is strict because $F_5$ fails
\eqref{abs-connectedJ} at $(x,y,z)=(a,b,0)$.
\end{proof}

As with the chain bases, these results assert finite basability
of the specified generated varieties. No assertion about all
their subvarieties is required.

\subsection{Finite bases for semirings with constant cubic products}
\label{abs-sec:nilgraphs}

We use the signature without constants throughout this section.
The symbol $\omega$ will denote an element of each algebra,
and later abbreviate a cubic word.

Define $T_4=\{a,b,r,\omega\}$ with addition given by maximum on
$a<b<r<\omega$. Its multiplication is commutative, with
$a^2=ab=r$ and every other product equal to $\omega$.
Define $T_5=\{a,b,c,r,\omega\}$ on the additive chain
$a<b<c<r<\omega$, with commutative multiplication determined by
$a^2=ab=ac=b^2=r$ and all other products equal to $\omega$.
Finally, let $U_5=\{a,b,c,r,\omega\}$ have additive order
\[
a<b<r<\omega,\qquad c<r,\qquad c\parallel a,b.
\]
Its multiplication is commutative, with
$a^2=ab=c^2=r$ and all other products equal to $\omega$.

\begin{theorem}\label{abs-nilgraph-bases}
The algebras $T_4,T_5,U_5$ are additively idempotent semirings.
Write $\Sigma_{\mathrm n}$ for
\begin{equation}\label{abs-nilgraph-common}
xy\approx yx,\qquad xyz\approx t^3,\qquad x\leq xy.
\end{equation}
Each algebra has a complete basis of four identities:
$\Sigma_{\mathrm n}$ together with its indicated additional identity:
\begin{align}
T_4:\quad xy&\leq x^2+yz, \label{abs-nilgraph-A}\\
T_5:\quad xy&\leq x^2+y^2, \label{abs-nilgraph-S}\\
U_5:\quad xy&\leq x^2+xz+zy. \label{abs-nilgraph-P}
\end{align}
\end{theorem}

Every product is in $\{r,\omega\}$, and both these elements
annihilate multiplication. Hence all cubic products are $\omega$,
which proves associativity. Translations in the two chains are
monotone. In $U_5$ they are monotone and preserve both incomparable
joins $a+c=b+c=r$: for each multiplier, at least one of the two
images is $\omega$. Thus distributivity holds in all three algebras.
Every product is at least each factor, proving
\eqref{abs-nilgraph-common}.

To check \eqref{abs-nilgraph-A} in $T_4$, suppose its right side is
not $\omega$. Then $x^2=yz=r$, so $x=a$ and $y\in\{a,b\}$,
whence $xy=r$. For \eqref{abs-nilgraph-S} in $T_5$, a right side
different from $\omega$ forces $x,y\in\{a,b\}$, again giving
$xy=r$. For \eqref{abs-nilgraph-P} in $U_5$, suppose
$x^2=xz=zy=r$. If $x=c$, then $z=y=c$; if $x=a$, then
$z,y\in\{a,b\}$. In both cases $xy=r$.

Under \eqref{abs-nilgraph-common}, all words of length at least
three equal a common multiplicative zero $\omega$, and
$x\leq x^3=\omega$. Thus a polynomial containing a long word
equals $\omega$; every other polynomial is a sum of one-letter
and quadratic words. No constant or empty word is adjoined.

Let $Q$ be a nonempty polynomial with no long word.
Let $U$ be its set of one-letter summands and $E$ its set
of unordered quadratic edges, including loops. Put
\[
V=\bigcup_{e\in E}\ct(e),\qquad L=\{x:x^2\in E\}.
\]
The associated graph has vertex set $V$; a loop contributes
a vertex, while an empty edge set has no vertices.

\begin{lemma}\label{abs-nilgraph-criteria}
In all three algebras, $x\leq Q$ holds exactly when $x\in U\cup V$,
$x^2\leq Q$ holds exactly when $x\in L$, and $\omega\leq Q$
never holds. For distinct variables $x,y$, the exact criteria
for $xy\leq Q$ are
\[
\begin{array}{c|l}
T_4&xy\in E\ \text{or}\
      (x,y\in V\ \text{and}\ \{x,y\}\cap L\ne\varnothing),\\
T_5&xy\in E\ \text{or}\ x,y\in L,\\
U_5&xy\in E\ \text{or}\
      (x,y\text{ belong to one component and }
       \{x,y\}\cap L\ne\varnothing).
\end{array}
\]
\end{lemma}

\begin{proof}[Necessity]
If $x\notin U\cup V$, assign $\omega$ to $x$ and $a$
everywhere else; then $Q\leq r<\omega=x$.
Assigning $a$ to every variable also shows that no such
$Q$ absorbs $\omega$.

If a variable $x$ of a left quadratic word is absent from $V$,
assign $r$ to $x$ and $a$ elsewhere. The left word is $\omega$,
whereas $Q\leq r$. Thus every left variable must be in $V$.
If the loop at $x$ is missing, assign $b$ to $x$ in $T_4,U_5$,
or $c$ to $x$ in $T_5$, and assign $a$ to all other variables.
Then $x^2=\omega$ but $Q\leq r$.

In $T_4$, if the edge $xy$ and both endpoint loops are missing,
assign $b$ to both $x,y$ and $a$ elsewhere. The only quadratic
words that could be $\omega$ are the three missing words
$x^2,xy,y^2$. Hence $Q\leq r<xy$.
In $T_5$, suppose $xy$ and, say, $y^2$ are missing.
Assign $x=b,y=c$ and $a$ elsewhere. Only those two missing
quadratic words could equal $\omega$, so absorption fails.

In $U_5$, if $x,y\in V$ belong to different components,
assign $c$ to the entire component of $x$ and $a$ elsewhere.
Every right quadratic word is $r$, while $xy=\omega$.
If the endpoints lie in one component but their edge and
both loops are missing, assign $x=y=b$ and $a$ elsewhere,
as for $T_4$. These assignments cover every failed condition.
\end{proof}

\begin{proof}[Sufficiency and completeness in Theorem~\ref{abs-nilgraph-bases}]
The one-letter criterion follows from $x\leq xy$ and commutativity.
An edge or loop already in $Q$ is absorbed trivially.
For the remaining $T_4$ case, assume $x^2\in E$ and $y\in V$.
Choose a right quadratic word $yz$, allowing $z=y$.
Then \eqref{abs-nilgraph-A} gives $xy\leq x^2+yz\leq Q$.
For $T_5$, both endpoint loops occur, so apply \eqref{abs-nilgraph-S}.

For $U_5$, assume $x^2\in E$ and choose a path
$x=v_0,v_1,\ldots,v_k=y$ in the graph.
The first edge gives $xv_1\leq Q$.
Inductively, \eqref{abs-nilgraph-P} gives
\[
xv_{i+1}\leq x^2+xv_i+v_iv_{i+1}\leq Q.
\]
Thus $xy\leq Q$. All substitutions use nonempty words.
These derivations prove every absorption in
Lemma~\ref{abs-nilgraph-criteria}. Expansion into sums of words,
with long words represented by the common zero, proves completeness.
\end{proof}

\begin{corollary}
We have $\V(T_4)\subsetneq\V(T_5)$ and
$\V(T_4)\subsetneq\V(U_5)$, while $\V(T_5)$ and $\V(U_5)$
are incomparable.
\end{corollary}

\begin{proof}
Identity \eqref{abs-nilgraph-A} implies \eqref{abs-nilgraph-S} by setting
$z=y$, and implies \eqref{abs-nilgraph-P} by adding $xz$ on the right.
The first inclusion is strict because $T_5$ fails
\eqref{abs-nilgraph-A} at $(b,c,a)$; the second is strict because
$U_5$ fails it at $(a,c,c)$.
Finally, $U_5$ fails \eqref{abs-nilgraph-S} at $(a,c)$, while
$T_5$ fails \eqref{abs-nilgraph-P} at $(b,c,a)$.
\end{proof}

\subsubsection{Five noncommutative semirings}

For the next examples let the common underlying set be
$\{\ell,m,h,r,\omega\}$, and let multiplication take only
the values $r,\omega$. Specify the pairs with product $r$ by
\begin{align*}
R_{\rm p}&=\{(\ell,m),(m,m),(m,h)\},\\
R_{\rm t}&=\{(\ell,\ell),(m,\ell),(h,\ell),(h,m),(h,h)\},\\
R_{\rm h}&=\{\ell,m\}^2\cup\{(h,\ell)\}.
\end{align*}
Every other pair has product $\omega$. Define five algebras
by the following additive orders and product relations:
\[
\begin{array}{c|l|c}
\text{algebra}&\text{additive order}&\text{product }r\text{ on}\\ \hline
N_{\rm p}&m<\ell,h<\omega,\quad r<\omega&R_{\rm p}\\
N_{\rm t}&\ell,h<m<\omega,\quad r<\omega&R_{\rm t}\\
N_{\rm h}&\ell<m<h<\omega,\quad r<\omega&R_{\rm h}\\
N_{\rm m}&\ell,h<m<\omega,\quad\ell<r<\omega&R_{\rm t}\\
N_{\rm i}&\ell,h<m<r<\omega&R_{\rm t}
\end{array}
\]
No other comparisons are imposed. Each row and column
fibre of the indicated relation is a join-closed lower set
of the additive order. Hence all translations preserve join.
Every triple product is $\omega$, proving associativity.

Put
\[
\Gamma=\{xyz\approx u^3,\ x\leq u^3\},\qquad
\Delta=\{u^3\leq x+xy,\ u^3\leq y+xy\}.
\]

\begin{theorem}\label{abs-directed-cubic-bases}
The following are complete bases:
\begin{align*}
N_{\rm p}:\quad&\Gamma,\Delta,\quad xy\leq xz+tx+wy,
                         \quad xy\leq xz+yt+wy;\\
N_{\rm t}:\quad&\Gamma,\Delta,\quad xz\leq xy+yz;\\
N_{\rm h}:\quad&\Gamma,\Delta,\quad xy\leq zx+ty;\\
N_{\rm m}:\quad&\Gamma,\quad xz\leq xy+yz,\quad y\leq x+xy,\\
&x^2\leq x+xy,\quad x^2\leq x+yx,\\
&yx\leq x+x^2+yz,\quad yx\leq x+x^2+zy;\\
N_{\rm i}:\quad&xyz\approx u^3,\quad x\leq xy,\quad y\leq xy,
                          \quad xz\leq xy+yz.
\end{align*}
The respective numbers of identities are six, five, five,
eight, and four. Reversing every word gives a complete
basis for each multiplicative dual.
\end{theorem}

\begin{proof}
Write $\omega$ for $u^3$. Long words reduce to $\omega$,
and $x\leq\omega$ holds in all five bases: in the last
one use $x\leq x^2\leq x^3$. Thus any polynomial containing
a long word is constant $\omega$.
Otherwise its data are a one-letter set $U$ and a directed
quadratic edge set $E$. Put
\[
H=\{x:(x,y)\in E\},\quad T=\{y:(x,y)\in E\},\quad V=H\cup T,
\]
and let $E^+$ be the relation of reachability by a path of
positive length. Define
\[
L=(U\cap V)\cup\{y:(x,y)\in E^+\text{ for some }x\in U\cap V\}.
\]
In the first three algebras an overlap $U\cap V$ makes
$Q=\omega$. Apart from this case and an explicit long
summand, no algebra absorbs $\omega$. The exact criteria
for quadratic absorption are
\[
\begin{array}{c|l}
N_{\rm p}&(x,y)\in E\text{ or }
 x\in H,\ y\in T,\ (x\in T\text{ or }y\in H)\\
N_{\rm t},N_{\rm i}&(x,y)\in E^+\\
N_{\rm h}&(x,y)\in E\text{ or }x,y\in T\\
N_{\rm m}&(x,y)\in E^+\text{ or }x\in V,\ y\in L.
\end{array}
\]
The one-letter absorptions are $U$ in the first three
algebras, $U\cup L$ in $N_{\rm m}$, and $U\cup V$ in
$N_{\rm i}$. These statements include the case $E=\varnothing$,
when just the original one-letter words are absorbed.

We verify the criteria and then derive them. With an edge
present, a value below $\omega$ requires every edge to be
$r$ and every one-letter value to be at most $r$.
In the first three algebras only $r$ itself is at most $r$
among possible one-letter values; it annihilates every
product. This explains the overlap case.
In $N_{\rm m}$ a one-letter graph vertex must be $\ell$,
and so must every vertex reachable from it.
In $N_{\rm i}$ every graph vertex is already below $r$.

For $R_{\rm p}$, vertices with both an incoming and an
outgoing edge must be $m$. Outgoing vertices take values
in $\{\ell,m\}$ and incoming vertices in $\{m,h\}$.
This proves the listed additional edges. If a missing
edge fails the incidence requirements, assign $h$ to
an incoming-only vertex or $\ell$ to an outgoing-only
vertex, with the other graph vertices $m$, to separate it.
If both required incidences hold but neither intermediate
incidence does, assign its source $\ell$ and its target
$h$; all other vertices may be $m$. No original edge
joins that pair, so all original edges still have value $r$.

For $R_{\rm h}$, every incoming vertex lies in $\{\ell,m\}$,
and all ordered pairs within that set are allowed.
If $x\notin T$ and $x\to y$ is not an original edge,
assign $x=h$, all its successors $\ell$, and $y=m$;
the other vertices take values in $\{\ell,m\}$.
If $y\notin T$, assign it $h$ and all other vertices
$\ell$. These assignments separate every missing pair,
including a missing loop in the second case.

The relation $R_{\rm t}$ is transitive. In the auxiliary
order $\ell<m<h$ it consists of all decreasing pairs,
with loops at $\ell,h$ but no loop at $m$.
If $x\ne y$ and $y$ is not reachable from $x$, assign
$\ell$ to $x$ and its successors and $h$ elsewhere.
For $N_{\rm m}$ also assign $\ell$ to $L$.
When the stated criterion fails, $y$ is outside this
forward-closed set, separating the target edge.
For a missing loop at $x$, assign $x=m$, its strict
successors and $L$ the value $\ell$, and all other
vertices $h$. No path returns to $x$, and $x\notin L$;
all original edges are allowed but $x^2=\omega$.

A variable outside $V$ is separated by assigning it
$\omega$, or $r$ if it is an independent one-letter
summand. For a one-letter target outside $L$ in
$N_{\rm m}$, assign $\ell$ on $L$ and $h$ elsewhere
on the graph; $h\not\leq r$ separates it.
The other one-letter conditions follow directly from
the additive orders. Assigning a looped value to all
graph vertices and suitable values to $U$ makes $Q=r$,
so no additional cubic absorption holds.
For a polynomial with only one-letter summands,
assigning one common looped value to all variables
excludes every quadratic absorption.
These assignments prove necessity; the same descriptions
of values below $\omega$ prove sufficiency.
In particular, each displayed basis identity satisfies
its criterion, proving validity.

For completeness, $\Delta$ derives the overlap cases.
The two rules for $N_{\rm p}$ give the two permitted
incidence completions using their three witnessing
edges. The rule for $N_{\rm h}$ joins any two edge targets.
Transitivity supplies every positive path in $E^+$,
including loops from directed cycles.
For $N_{\rm i}$, the two factor bounds supply all
one-letter graph vertices.

In $N_{\rm m}$, use $y\leq x+xy$ successively along
paths from one-letter graph vertices to obtain every
one-letter word in $L$. Each such vertex has an incident
edge, so the two square rules supply its square.
The last two rules then connect any vertex of $V$ to
any vertex of $L$, using an outgoing or an incoming
edge at the former vertex. This gives precisely
$V\times L$. Hence every stated absorption is derivable.
Expanding both sides of any identity completes the proof.
\end{proof}

\subsection{Finite bases for semirings with constant quartic products}
\label{abs-sec:quartic}

\subsubsection{Two five-element semirings}

The quadratic graph method also applies when cubic products are
nonconstant. Let $V_5$ and $W_5$ have underlying set
$\{a,b,d,c,\omega\}$. Their additive orders are, respectively,
\begin{align*}
V_5:&\quad a,b<d<c<\omega,\qquad a\parallel b,\\
W_5:&\quad a<d<c<\omega,\quad b<c,\qquad b\parallel a,d.
\end{align*}
Addition is join in each order. Multiplication is commutative,
with the following products different from $\omega$:
\[
\begin{array}{c|l}
V_5&a^2=d,\quad ab=ad=c,\\
W_5&a^2=d,\quad ad=b^2=c.
\end{array}
\]
Every unlisted product equals $\omega$.

\begin{theorem}\label{abs-quartic-bases}
The algebras $V_5,W_5$ are additively idempotent semirings.
Let $\Sigma_{\mathrm q}$ consist of the four identities
\begin{equation}\label{abs-quartic-common}
xy\approx yx,\qquad xyzt\approx u^4,\qquad
xyz\approx x^3+y^3+z^3,\qquad x\leq xy.
\end{equation}
Each algebra has a complete basis of six identities:
$\Sigma_{\mathrm q}$ together with its indicated pair:
\begin{align}
V_5:\quad xy&\leq x^2+yz,& x^3&\leq x^2+y^3,
\label{abs-quartic-V}\\
W_5:\quad xy&\leq x^2+xz+zy,& x^3&\leq x^2+xy+y^3.
\label{abs-quartic-W}
\end{align}
\end{theorem}

In either algebra, a triple product, with either bracketing,
equals $c$ exactly when all three factors are $a$; otherwise
it equals $\omega$. This proves associativity, and every
four-factor product is $\omega$. Each multiplication
translation is monotone. In $V_5$ the only incomparable pair
is $a,b$, with join $d$: multiplication by $a$ sends this
join to $c=d+c$, and every other translation sends at least
one member of the pair, and its join, to $\omega$.
In $W_5$ the incomparable pairs are $a,b$ and $d,b$, both
with join $c$; each translation sends $c$ and at least one
member of either pair to $\omega$. Thus translations preserve
joins, proving both distributive laws. Every product is at
least both factors. The description of triple products proves
the remaining identity in \eqref{abs-quartic-common}.

For the first identity in \eqref{abs-quartic-V}, suppose its
right side is not $\omega$. Then $x=a$ and $y\in\{a,b,d\}$.
If $yz=d$, then $y=z=a$ and $xy=d$; otherwise the right
side is $c$ and $xy\leq c$. For the second identity, a
right side different from $\omega$ forces $x=y=a$.

For the first identity in \eqref{abs-quartic-W}, a right side
different from $\omega$ forces $x\in\{a,b\}$.
If $x=b$, then $z=y=b$ and all the displayed products are
$c$. If $x=a$, then $z\in\{a,d\}$; when $z=a$ we have
$y\in\{a,d\}$, and when $z=d$ we have $y=a$.
In each case $xy$ is at most the right side.
For the second identity, its right side can differ from
$\omega$ only if $y=a$ and $x\in\{a,b\}$; the factor
$xy$ then forces $x=a$.

Under \eqref{abs-quartic-common}, every word of length at least
four equals a common word value, denoted by $\omega$.
The inequality $x\leq xy$, iterated three times, gives
$x\leq x^4=\omega$. Thus any polynomial containing a
word of length at least four equals $\omega$.
All other polynomials reduce to sums of one-letter words,
unordered quadratic words, and pure cubes $x^3$.
Here $\omega$ abbreviates a fourth power; the signature
still has no constants or empty word.

For a nonempty polynomial $Q$ in this form, with no fourth
power, write $U$ for its one-letter summands, $E$ for its
unordered quadratic edges, and $D$ for the variables whose
cubes are summands. Put
\[
L=\{x:x^2\in E\},\qquad F=L\cup D,\qquad
V=\bigcup_{e\in E}\ct(e)\cup D.
\]
Let $G$ be the graph with vertex set $V$ and edges $E$,
with an additional edge between every two distinct vertices
of $D$. Let $R$ be the set of vertices reachable from $D$.
In particular, $R=\varnothing$ when $D=\varnothing$.

\begin{lemma}\label{abs-quartic-criteria}
In both algebras, $x\leq Q$ holds exactly when $x\in U\cup V$,
$x^2\leq Q$ holds exactly when $x\in F$, and
$\omega\leq Q$ never holds.
For distinct variables $x,y$, the remaining criteria are
\begin{align*}
V_5\models xy\leq Q
&\quad\Longleftrightarrow\quad
xy\in E\ \text{or}\ (x,y\in V\ \text{and}\ \{x,y\}\cap F\ne\varnothing),\\
W_5\models xy\leq Q
&\quad\Longleftrightarrow\quad
xy\in E\ \text{or}\
\left\{\begin{array}{l}
x,y\text{ lie in one component of }G,\\
\{x,y\}\cap F\ne\varnothing,
\end{array}\right.\\
V_5\models x^3\leq Q
&\quad\Longleftrightarrow\quad D\ne\varnothing\ \text{and}\ x\in F,\\
W_5\models x^3\leq Q
&\quad\Longleftrightarrow\quad x\in F\cap R.
\end{align*}
\end{lemma}

\begin{proof}[Necessity]
Assigning $a$ to every variable gives $Q\leq c<\omega$,
and gives $Q\leq d<c$ when $D=\varnothing$.
This proves the assertion about $\omega$ and the necessity
of a cube on the right when a cube is absorbed.
If $x\notin U\cup V$, assigning $x=\omega$ and $a$
elsewhere separates $x$ from $Q$.
If a variable of a left quadratic word is absent from $V$,
assign $c$ to it and $a$ elsewhere. That word equals
$\omega$ and $Q\leq c$.

If $x\notin F$, assign $x=b$ in $V_5$, or $x=d$ in
$W_5$, and $a$ elsewhere. Then both $x^2$ and $x^3$
are $\omega$, whereas $Q\leq c$.
If $xy\notin E$ and neither endpoint is in $F$, assign
both endpoints $b$ in $V_5$, or $d$ in $W_5$, and
$a$ elsewhere. The target product is $\omega$.
No right quadratic word contains two marked occurrences,
and no right cube is marked, so $Q\leq c$.

In $W_5$, suppose $x,y\in V$ lie in different components
of $G$. At most one of these components meets $D$, since
$D$ is a clique. Assign $b$ throughout one component not
meeting $D$, and $a$ elsewhere. Right quadratic products
are $c$ inside that component and $d$ outside; right cubes
are $c$. Thus $Q\leq c<\omega=xy$.
The same assignment to the component of $x$ separates
$x^3$ from $Q$ whenever $x\in F\setminus R$.
These assignments cover every failed condition.
\end{proof}

\begin{proof}[Sufficiency and completeness in Theorem~\ref{abs-quartic-bases}]
The inequality $x\leq xy$ gives $x\leq x^2\leq x^3$,
so every variable occurring in $Q$ is absorbed by $Q$,
and $x^2\leq Q$ follows for every $x\in F$.
Also, for any $x,y$,
\begin{equation}\label{abs-quartic-cube-edge}
xy\leq x^2y\approx x^3+y^3.
\end{equation}
Thus every added edge of the clique on $D$ is absorbed
by $Q$, as is every edge already in $E$.

For $V_5$, suppose $x\in F$ and $y\in V$.
Choose a nonlinear summand involving $y$, written $yz$,
where $z$ is a nonempty word; a cube $y^3$ permits
$z=y^2$. The first identity in \eqref{abs-quartic-V} gives
$xy\leq x^2+yz\leq Q$.
If $D\ne\varnothing$, choose any $y\in D$.
The second identity gives $x^3\leq x^2+y^3\leq Q$
for every $x\in F$.

For $W_5$, fix $x\in F$ and a path
$x=v_0,v_1,\ldots,v_k=y$ in $G$.
Every path edge is absorbed by $Q$, including the
clique edges by \eqref{abs-quartic-cube-edge}.
Starting with $xv_1\leq Q$, repeated application of the
first identity in \eqref{abs-quartic-W} gives
\[
xv_{i+1}\leq x^2+xv_i+v_iv_{i+1}\leq Q.
\]
Hence $xy\leq Q$. If $x\in F\cap R$, choose such a
path to $y\in D$. If $x=y$, its cube is already a
summand; otherwise use the path deduction and the second
identity to obtain $x^3\leq x^2+xy+y^3\leq Q$.

All substitutions use nonempty words. We have derived
every absorption in Lemma~\ref{abs-quartic-criteria}.
Normalization and the two directions of absorption for
an arbitrary identity prove completeness.
\end{proof}

\begin{corollary}
The inclusion $\V(V_5)\subsetneq\V(W_5)$ is strict.
\end{corollary}

\begin{proof}
Each inequality in \eqref{abs-quartic-V} implies its counterpart
in \eqref{abs-quartic-W} by adding a right summand.
In $W_5$, the substitution $(x,y,z)=(a,b,b)$ gives
$xy=\omega>c=x^2+yz$, proving strictness.
\end{proof}

These are finite bases for the two individual generated
varieties. Neither minimality of the displayed bases nor
hereditary finite basability is asserted.

\subsubsection{A four-element semiring and the parity of paths}

Let $H_4=\{a,d,c,\omega\}$ have addition given by join in
\[
a<c<\omega,\qquad d<\omega,\qquad d\parallel a,c,
\]
and commutative multiplication with $a^2=d$, $ad=c$, and
all other products equal to $\omega$.
Triple products are $c$ exactly when all three factors are
$a$, and otherwise are $\omega$, so multiplication is
associative. Translations preserve the comparable pair $a<c$
and the incomparable joins $a+d=c+d=\omega$.
Hence $H_4$ is an additively idempotent semiring.

\begin{theorem}\label{abs-quartic-parity-basis}
A complete basis for $H_4$ consists of the twelve identities
\begin{align}
xy&\approx yx, & xyzt&\approx u^4, \label{abs-parity-basic}\\
xyz&\approx x^3+y^3+z^3, & x&\leq x^3, \label{abs-parity-cubic}\\
x^3&\leq x+xy, & x^2+y^3&\approx u^4, \label{abs-parity-anchor}\\
xt&\leq xy+yz+zt, & y^2&\leq x^2+yz, \label{abs-parity-edges}\\
xy&\leq x^2+y^2, & z^3&\leq x^3+xy+yz, \label{abs-parity-even}\\
xz&\leq x^3+y^3+yz, & xy+x^3+y^3&\approx u^4. \label{abs-parity-transfer}
\end{align}
\end{theorem}

We first verify validity. A sum of quadratic products differs
from $\omega$ only when all those products are $d$, or all
are $c$. A product is $d$ exactly on $(a,a)$, and is $c$
exactly on $(a,d)$ or $(d,a)$. A cube differs from $\omega$
exactly at $a$, where its value is $c$.
These facts prove \eqref{abs-parity-basic} and the first identity
in \eqref{abs-parity-cubic}; the other follows from $a<c$.
The sum $x+xy$ can differ from $\omega$ only at $(x,y)=(a,d)$,
where it is $c$. Also, a square plus a cube is always
$\omega$, proving \eqref{abs-parity-anchor}.

For the first identity in \eqref{abs-parity-edges}, three products
equal to $d$ force all four variables to be $a$; three equal
to $c$ force alternating values $a,d$. The endpoint product
has the same value in either case. For the second identity,
a right side different from $\omega$ forces $x=y=z=a$.
The first identity in \eqref{abs-parity-even} follows because
both squares can be nonzero only at $x=y=a$.
For its second identity, a right side different from
$\omega$ forces $x=z=a$ and $y=d$.
Likewise the first identity in \eqref{abs-parity-transfer}
forces $x=y=a,z=d$ when its right side is nonzero.
For the last identity, either a cube is $\omega$, or
$x=y=a$ and the sum contains both $d$ and $c$.

For normalization, all words of length at least four equal
$\omega=u^4$, and $x\leq x^3\leq x^9=\omega$.
Every other polynomial is a sum of one-letter words,
quadratic words, and pure cubes. Whenever a one-letter
summand $x$ also occurs in a quadratic summand $xy$,
the first identity in \eqref{abs-parity-anchor}, together
with $x\leq x^3$, permits its replacement by $x^3$.
A one-letter summand already represented by a cube can
be deleted. Thus we may write the normalized data as
$(U,E,D)$, where $E$ is a set of unordered quadratic edges,
$D$ is a set of cube variables, and
\[
V=\bigcup_{e\in E}\ct(e)\cup D,\qquad U\cap V=\varnothing.
\]
Empty blocks are omitted, and the polynomial is nonempty.
Let $G=(V,E)$, allowing loops and isolated vertices of $D$.
A proper two-colouring assigns each vertex one of two
colours and gives distinct colours to the endpoints of
every edge; in particular a loop prevents such a colouring.

\begin{lemma}\label{abs-parity-criterion}
Let $Q$ have normalized data $(U,E,D)$.
If $E=D=\varnothing$, it absorbs exactly the one-letter
words belonging to $U$ among the normalized monomials.
Otherwise the following cases give all its absorptions.
\begin{enumerate}
\item If $D=\varnothing$ and $G$ is not bipartite, $Q$
absorbs exactly the one-letter words in $U$ and all
quadratic words with variables in $V$.
\item If $D\ne\varnothing$ and no proper two-colouring
of $G$ makes every vertex of $D$ colour $A$, then
$Q\approx\omega$.
\item In the remaining cases, consider all proper
two-colourings of $G$ making $D$ colour $A$; when
$D=\varnothing$ there is no colour restriction.
Then $Q$ absorbs a cube $x^3$ exactly when $x\in V$
has colour $A$ in every such colouring. It absorbs
a one-letter word $x$ exactly when $x\in U$ or this
cube condition holds. It absorbs $xy$, allowing $x=y$, exactly when
$x,y\in V$ have opposite colours in every such colouring.
It does not absorb $\omega$.
\end{enumerate}
\end{lemma}

\begin{proof}[Necessity and the evaluation criterion]
If $E=D=\varnothing$, assign $a$ to all variables to
exclude every nonlinear monomial, and assign $\omega$
to a missing variable to exclude a missing one-letter word.
Now suppose $E\cup D\ne\varnothing$.
The only possible values of $Q$ below $\omega$ are $d,c$.
The value $d$ occurs exactly when $D=\varnothing$,
every variable in $V$ is $a$, and every variable in
$U$ is $d$. The value $c$ occurs exactly when the
values $a,d$ on $V$ give a proper two-colouring of $G$,
every variable in $D$ is $a$, and the values on $U$
belong to $\{a,c\}$.

These descriptions follow directly from the quadratic
and cubic product rules above. At an assignment giving
$Q=d$, a target quadratic word is below $Q$ exactly
when both factors are $a$; a target cube or a vertex
variable is not below $d$.
At an assignment giving $Q=c$, a target cube is below
$Q$ exactly when its variable is $a$, and a target
quadratic word is below $Q$ exactly when its factors
are $a,d$ in either order.
Variables outside $U\cup V$ can independently be set
to $\omega$. For nonlinear targets involving a variable
in $U$, choose that variable to be $d$ in a $d$-valued
assignment, or $c$ in a $c$-valued assignment.
These choices exclude every extra absorption listed
as invalid. Conversely the same evaluation rules prove
every stated valid absorption and characterize $Q\approx\omega$.
\end{proof}

\begin{proof}[Completeness in Theorem~\ref{abs-quartic-parity-basis}]
Every edge of $G$ is absorbed by $Q$.
The first identity in \eqref{abs-parity-edges}, applied
repeatedly, shows that the endpoints of every odd walk
give an absorbed quadratic word. Walks of length one
are the initial edges; a longer odd walk is shortened
by replacing its first three edges by their endpoint edge.
Thus an odd cycle gives an absorbed loop $x^2$.
The second identity in \eqref{abs-parity-edges} then supplies
the square of every vertex occurring in an edge,
including vertices in other components.
The first identity in \eqref{abs-parity-even} supplies
all quadratic words on these vertices.
This proves all the absorptions in case (1).

If $D\ne\varnothing$ and $G$ has an odd cycle,
an absorbed loop and any cube from $D$ give $\omega\leq Q$
by \eqref{abs-parity-anchor}.
If $G$ is bipartite but no colouring makes $D$ colour
$A$, some component has two vertices of $D$ in opposite
parts. An odd path between them gives an absorbed edge,
and their two cubes give $\omega\leq Q$ by the last
identity in \eqref{abs-parity-transfer}. This proves case (2).

In case (3), a vertex has forced colour $A$ exactly when
it is connected to $D$ by an even path. Starting with
the cube at a vertex of $D$, the second identity in
\eqref{abs-parity-even} propagates cubes along two successive
edges. Repetition gives every required cube; $x\leq x^3$
then supplies the required one-letter absorptions.

Two vertices have forced opposite colours exactly when
they lie in opposite parts of one component, or lie in
different components that both meet $D$ and have opposite
colours fixed by $D$.
The first possibility is handled by an odd path.
In the second, assume the target is $xz$, with $x$
forced colour $A$ and $z$ forced colour $B$.
The vertex $z$ has a neighbour $y$ of forced colour $A$.
The cubes at $x,y$ have already been derived, so
\[
xz\leq x^3+y^3+yz\leq Q
\]
by \eqref{abs-parity-transfer}.
This derives all the remaining quadratic absorptions.

Normalization and Lemma~\ref{abs-parity-criterion} now derive
both directions of every identity valid in $H_4$.
All substitutions used are by nonempty words.
\end{proof}

The theorem proves finite basability of $\V(H_4)$;
it makes no assertion about every subvariety of $\V(H_4)$.

\subsubsection{Preserving one-letter summands}

The next extension distinguishes the one-letter summands of
a polynomial from its other occurrences of variables.

\begin{lemma}\label{abs-ceiling-identities}
Let the additive greatest element $\omega$ of a semiring $B$
also be a multiplicative zero. Adjoin a new additive greatest
element $\tau$, and put $s\tau=\tau s=\omega$ for every
$s\in B\cup\{\tau\}$. This defines a semiring $B^\tau$.
For an expanded polynomial $P$, let $U(P)$ be the set of
its one-letter summands. Then
\[
B^\tau\models P\approx Q
\quad\Longleftrightarrow\quad
B\models P\approx Q\ \text{ and }\ U(P)=U(Q).
\]
Consequently, for a nonlinear word $w$, the absorption
$w\leq Q$ holds in $B^\tau$ exactly when it holds in $B$.
For a one-letter word $x$, it holds exactly when $x\in U(Q)$.
\end{lemma}

\begin{proof}
A triple product involving $\tau$ is $\omega$ in either
association. Each distributive law involving $\tau$ has
value $\omega$ on both sides; here we use that every old
product is at most $\omega$. All remaining semiring laws
hold in $B$.
The restriction to $B$ proves one necessary condition.
If the sets $U(P),U(Q)$ differ, assign $\tau$ to a variable
in their symmetric difference and $\omega$ to all others.
The resulting values are $\tau$ and $\omega$.
Conversely, assume the two stated conditions. An assignment
giving $\tau$ to a member of the common one-letter set makes
both polynomials $\tau$. Otherwise replace every assigned
$\tau$ by $\omega$. This changes no one-letter value and
no nonlinear word value, so the identity follows from $B$.
The absorption assertions are the special case $P=w+Q$.
\end{proof}

Let $C_5,D_5$ have underlying set
$\{a,d,c,\omega,\tau\}$ and the commutative products
$a^2=d$, $ad=c$, with all other products $\omega$.
In $C_5$, the elements $a,d,c$ are incomparable and lie
below $\omega<\tau$. In $D_5$, the order is
$a<c<\omega<\tau$, $d<\omega$, with $d\parallel a,c$.
Thus $C_5$ extends the flat four-element semiring with
these products, and $D_5=H_4^\tau$.

Define $E_5,Z_5$ on $\{a,b,d,c,\omega\}$ by the orders
\begin{align*}
E_5:&\quad b<a<\omega,\quad b<c<\omega,\quad d<\omega,\\
Z_5:&\quad a<b<\omega,\quad c<b,\quad d<\omega,
\end{align*}
with no further nontrivial comparabilities.
Multiplication is commutative. Its values different from
$\omega$ are
\begin{align*}
E_5:&\quad a^2=ab=b^2=d,\quad ad=bd=c,\\
Z_5:&\quad a^2=d,\quad ad=c.
\end{align*}
Call $a,b$ active in $E_5$, and only $a$ active in the
other three algebras. A triple product is $c$ exactly when
all three factors are active, and is otherwise $\omega$.
This proves associativity. For $E_5,Z_5$, the nonconstant
multiplicative translations preserve the displayed order
and the joins of incomparable pairs, proving distributivity.
For $C_5,D_5$, use Lemma~\ref{abs-ceiling-identities}; the
flat four-element algebra also has join-preserving translations.

Let $\Sigma$ be the ten identities
\begin{align}
xy&\approx yx, & xyzt&\approx u^4,\label{abs-linear-core1}\\
xyz&\approx x^3+y^3+z^3, & x^2+y^3&\approx u^4,\label{abs-linear-core2}\\
xt&\leq xy+yz+zt, & y^2&\leq x^2+yz,\label{abs-linear-core3}\\
xy&\leq x^2+y^2, & z^3&\leq x^3+xy+yz,\label{abs-linear-core4}\\
xz&\leq x^3+y^3+yz, & xy+x^3+y^3&\approx u^4.\label{abs-linear-core5}
\end{align}
These are precisely the identities of
Theorem~\ref{abs-quartic-parity-basis} other than
$x\leq x^3$ and $x^3\leq x+xy$.

\begin{theorem}\label{abs-quartic-linear-bases}
The following sets are complete identity bases:
\begin{align*}
C_5:&\quad \Sigma,\quad xy\leq u^4,\quad u^4\leq x+xy;\\
D_5:&\quad \Sigma,\quad xy\leq u^4,\quad x^3\leq x+xy;\\
E_5:&\quad \Sigma,\quad x\leq u^4,\quad x^3\leq x+xy;\\
Z_5:&\quad \Sigma,\quad x\leq u^4,\quad x^3\leq x+xy,
                 \quad y\leq x+x^3+y^3.
\end{align*}
The first three bases have twelve identities, and the
last has thirteen.
\end{theorem}

For validity, a nonlinear word can take only $d,c,\omega$,
and $d+c=\omega$. A quadratic product is $d$ precisely
when both factors are active, and is $c$ precisely when
one is active and the other is $d$. A cube is $c$ precisely
at an active value. The same verification as for
Theorem~\ref{abs-quartic-parity-basis}, with active values
in place of $a$, proves every identity of $\Sigma$.
The upper bounds follow from the additive orders.
In $C_5$ every value of $x+xy$ is at least $\omega$.
In the other three algebras, a value of this sum below
$\omega$ requires $x$ active and $y=d$; the sum is then
at least $x^3=c$. The value $x=\tau$ in $D_5$ also
satisfies the inequality. Finally, in $Z_5$, a value of
$x^3+y^3$ below $\omega$ requires $x=y=a$, and then
$x+x^3+y^3=b\geq y$.

We give the absorption criteria before proving completeness.
Normalize words using the first three identities of $\Sigma$.
Thus a polynomial consists of one-letter summands $U$,
quadratic edges $E$, pure cubes on a set $D$, and possibly
$\omega=u^4$. In this normalization the set $U$ is retained.
Empty summand blocks are omitted, and the polynomial is nonempty.
Put
\[
V_E=\bigcup_{e\in E}\ct(e),\qquad V=V_E\cup D,
\qquad M=D\cup(U\cap V_E),\qquad G=(V,E).
\]
For a graph $G$ and marked set $M$, define
$\mathcal N(G,M)$ as follows:
\begin{enumerate}
\item If $M=\varnothing$ and $G$ is not bipartite, it consists
of all quadratic words with both variables in $V$.
\item If $M\ne\varnothing$ and no proper two-colouring
makes every member of $M$ colour A, it consists of all
normalized nonlinear words, including $\omega$.
\item Otherwise it consists of the quadratic words whose
endpoints are different colours in every such colouring,
and the cubes of vertices that are always colour A.
All these endpoints belong to $V$, and $\omega$ is absent.
\end{enumerate}
If $V=\varnothing$, this set is empty.

\begin{lemma}\label{abs-quartic-linear-criterion}
If $Q$ contains $\omega$, it absorbs every nonlinear word;
its one-letter absorptions are exactly $U$ for $C_5,D_5$,
and all one-letter words for $E_5,Z_5$.
Otherwise its absorptions are as follows.
\begin{enumerate}
\item In $C_5$, the nonlinear absorptions are all normalized
nonlinear words if $U\cap V\ne\varnothing$, and are
$\mathcal N(G,D)$ otherwise. Its one-letter absorptions
are exactly $U$.
\item In each of $D_5,E_5,Z_5$, the nonlinear absorptions
are exactly $\mathcal N(G,M)$. In $D_5$ the one-letter
absorptions are exactly $U$. In $E_5$ they are all variables
if $\omega\in\mathcal N(G,M)$, and exactly $U$ otherwise.
\item In $Z_5$, all one-letter words are absorbed if
$\omega\in\mathcal N(G,M)$. Otherwise a one-letter word
$y\notin U$ is absorbed exactly when
$U\cap V\ne\varnothing$ and $y^3\in\mathcal N(G,M)$.
\end{enumerate}
\end{lemma}

\begin{proof}
The assertions for $D_5$ follow from
Lemmas~\ref{abs-ceiling-identities} and~\ref{abs-parity-criterion}.
The latter uses $M$ as its cube set after normalization
in $H_4$; this does not permit deleting $U$ in $D_5$.
The one-letter assertion for $C_5$ also follows from
Lemma~\ref{abs-ceiling-identities}. An overlap $U\cap V$ in
$C_5$ makes the corresponding one-letter summand plus
an incident edge or cube at least $\omega$.
Without such an overlap, $Q=d$ is possible exactly when
$D=\varnothing$: all graph variables must be $a$, and
all one-letter variables $d$. The value $Q=c$ corresponds
to a proper colouring by $a,d$, with $D$ coloured $a$ and
all one-letter variables $c$. These alternatives give
the stated nonlinear criterion.

For $E_5,Z_5$, suppose first that $Q$ has a nonlinear
summand. A value $Q=d$ is possible exactly when
$M=\varnothing$: assign active values to graph variables
and $d$ to the disjoint one-letter set. Any other value
below $\omega$ requires a proper colouring of $G$ by
active values and $d$, with $M$ active.
In $E_5$, a variable belonging to both $U$ and $V$ must
then be $b$; the other one-letter variables may all be $c$,
giving $Q=c$. In $Z_5$, the common variables must be $a$.
If $U\cap V$ is nonempty the resulting sum is $b$;
otherwise assigning $c$ to all one-letter variables gives
$Q=c$.

These alternatives prove the nonlinear conditions.
Indeed, a cube has value $c$ exactly on the active colour.
A quadratic word with same-coloured endpoints has value
$d$ or $\omega$, neither of which is at most $c$ or $b$.
A variable outside $V$ can be assigned $c$ if it belongs
to $U$, and $\omega$ otherwise, separating any nonlinear
word that uses it. In the alternative $Q=d$, assign $d$
to variables of $U$ instead. Thus no further nonlinear
absorptions are possible.

For a one-letter target $y\notin U$ in $E_5$, an active
vertex can be assigned $a\not\leq c$, and the other colour
is $d\not\leq c$. These choices do not alter the required
$b$ values at vertices of $U$. If only $Q=d$ is possible,
an active value $a$ at $y$ separates it. A variable outside
$V\cup U$ can be assigned $\omega$.
In $Z_5$, when $U\cap V$ is empty the assignment $Q=c$
likewise excludes additional one-letter absorptions.
When this intersection is nonempty, every value below
$\omega$ is $b$, and a vertex is at most that value
under every assignment precisely when it is always active.
The value $d$ separates any other vertex.
Finally, if $Q$ has only one-letter summands, assigning
$a$ to every variable excludes every nonlinear absorption;
a missing one-letter variable is separated by assigning
it $\omega$. The assertions for a summand $\omega$ follow
from its additive position in each algebra.
\end{proof}

\begin{proof}[Completeness in Theorem~\ref{abs-quartic-linear-bases}]
In $C_5,D_5$, the bound $xy\leq\omega$ removes only nonlinear
summands beside $\omega$: the resulting form is
$\sum U+\omega$. In $E_5,Z_5$, $x\leq\omega$ reduces
such a polynomial to $\omega$.

For $C_5$, an overlap in $U\cap V_E$ gives
$\omega\leq Q$ by its last axiom. An overlap in $U\cap D$
gives the same result by substituting $y=x^2$ in that axiom.
For $D_5,E_5,Z_5$, add $x^3$ for each
$x\in U\cap V_E$ using $x^3\leq x+xy$.
This produces the marked set $M$ while retaining every
one-letter summand.

The nonlinear derivations in the proof of
Theorem~\ref{abs-quartic-parity-basis} use only $\Sigma$.
For clarity, odd walks yield their endpoint edges by
\eqref{abs-linear-core3}. An odd cycle yields a square;
squares propagate to all edge vertices and give all
quadratic words on those vertices. A square and a cube
yield $\omega$. Marked vertices on opposite sides of a
bipartite component yield their connecting edge and then
$\omega$ by \eqref{abs-linear-core5}. In the compatible case,
\eqref{abs-linear-core4} propagates cubes along even paths.
Odd paths supply all forced opposite-colour edges within
a component, and \eqref{abs-linear-core5} supplies such edges
between marked components. These steps derive every
member of $\mathcal N(G,M)$ without removing a one-letter
summand.

The remaining one-letter absorptions are immediate except
for the last clause concerning $Z_5$. Choose
$x\in U\cap V$. Its cube is either present initially or
has just been adjoined. Once $y^3\leq Q$ is derived, use
\[
y\leq x+x^3+y^3\leq Q.
\]
If $\omega\leq Q$ in $E_5$ or $Z_5$, their upper-bound
axiom instead gives $y\leq\omega\leq Q$ for every $y$.
Thus every absorption in Lemma~\ref{abs-quartic-linear-criterion}
is derivable. Expanding an arbitrary identity and treating
its summands in both directions proves completeness.
\end{proof}

\subsubsection{Two inactive values and additional linear constraints}

The next four semirings have underlying set
$\{a,b,d,c,\omega\}$. Their additive orders are specified
by the following covers, with no further comparabilities
except those implied by transitivity:
\begin{align*}
P_5:&\quad b<d,\ b<c,\ d<\omega,\ c<\omega,\ a<\omega;\\
Q_5:&\quad b<d,\ b<c,\ a<c,\ d<\omega,\ c<\omega;\\
R_5:&\quad a<b,\ d<b,\ b<\omega,\ c<\omega;\\
S_5:&\quad a<b,\ d<b,\ a<c,\ b<\omega,\ c<\omega.
\end{align*}
Multiplication is commutative. In $P_5,Q_5$ its values
different from $\omega$ are $a^2=d$, $ab=ad=c$;
in $R_5,S_5$ they are $a^2=d$, $ad=c$.
A triple product is $c$ only on $(a,a,a)$ and is otherwise
$\omega$, proving associativity. The nonconstant translations
preserve the covers and the joins of incomparable pairs,
so all four algebras are semirings.

\begin{theorem}\label{abs-quartic-passive-bases}
With $\Sigma$ as in \eqref{abs-linear-core1}--\eqref{abs-linear-core5},
the following complete bases have respectively twelve,
twelve, fourteen, and fifteen identities:
\begin{align*}
P_5:\quad&\Sigma,\quad x\leq u^4,\quad y^3\leq x+xy;\\
Q_5:\quad&\Sigma,\quad x\leq x^3,\quad x+x^2\approx u^4;\\
R_5:\quad&\Sigma,\quad x\leq u^4,\quad z^2\leq x+xy+zt,\\
&z\leq x+xy+z^2,\quad x+xy+z^3\approx u^4;\\
S_5:\quad&\Sigma,\quad x\leq x^3,\quad z\leq x+xy+yz,\\
&x^2\leq x+y+xy,\quad xz\leq x+xt+y+yz,\quad
x^3\leq x+xy+z^3.
\end{align*}
\end{theorem}

\begin{proof}
For validity, a quadratic product is $d$ only at $(a,a)$.
It is $c$ when its factors are $a$ and a member of $\{b,d\}$
in $P_5,Q_5$, or $a,d$ in $R_5,S_5$.
Squares take only $d,\omega$, and cubes only $c,\omega$.
These rules verify $\Sigma$ as before.
In $P_5$, the sum $x+xy$ is below $\omega$ only when
$x=b,y=a$, proving its extra absorption. In $Q_5$,
$x\leq x^3$ and $x+x^2=\omega$ follow from the order.
In $R_5$, the sum $x+xy$ is below $\omega$ only at
$x=y=a$, when it is $b$; its three extra rules follow.
For $S_5$, a sum with a nonlinear summand has value below
$\omega$ only at $d,c,b$. The value $c$ requires a proper
colouring by $a,d$ and makes every one-letter graph vertex
$a$. The value $b$ requires all graph vertices $a$.
Substitution in these alternatives verifies its four
additional absorption rules.

For completeness, normalize to $(U,E,D)$ and put
$V=V_E\cup D$, $L=U\cap V_E$. In every case a summand
$\omega=u^4$ absorbs the whole polynomial; where needed,
$x\leq x^3\leq x^9=\omega$ proves this.
A polynomial with no nonlinear summand absorbs just its
one-letter summands. We record the criteria in the other cases.
\begin{enumerate}
\item In $P_5$, an intersection $U\cap D$ makes $Q=\omega$.
Otherwise put
\[
M=D\cup\{y:xy\in E\text{ for some }x\in U\}.
\]
The nonlinear absorptions are $\mathcal N(G,M)$.
The one-letter absorptions are exactly $U$, unless
$\omega$ is absorbed, in which case all are.
\item In $Q_5$, use $\mathcal N(G,D)$, with the additional
zero case $D=\varnothing$, $G$ not bipartite, and
$U\cap V\ne\varnothing$. Outside the zero cases,
the one-letter absorptions are $U$ and the vertices
whose cubes are absorbed.
\item In $R_5$, if $U\cap V\ne\varnothing$ and
$D\ne\varnothing$, then $Q=\omega$.
If $U\cap V\ne\varnothing$ and $D=\varnothing$, the
absorptions are exactly the one-letter words in $U\cup V$
and all quadratic words on $V$.
Without an overlap, the nonlinear absorptions are
$\mathcal N(G,D)$ and the one-letter absorptions are
exactly $U$, except in a zero case.
\item In $S_5$, if $D\ne\varnothing$, use
$\mathcal N(G,D\cup L)$ and include, as one-letter
absorptions, $U$ and all vertices whose cubes are absorbed;
in a zero case every word is absorbed.
If $D=\varnothing$, no cube or $\omega$ is absorbed.
When no proper colouring makes all of $L$ colour A,
all quadratic words on $V$ are absorbed; the one-letter
absorptions are $U$, together with $V$ if $L\ne\varnothing$.
Otherwise the quadratic absorptions are the universally
opposite-colour pairs, and the one-letter absorptions
are $U$ together with the universally colour A vertices.
\end{enumerate}

Here is a direct verification of these criteria. In $P_5$,
a value $c$ requires a proper colouring by $a$ and $\{b,d\}$.
A one-letter graph vertex must be $b$, so all its neighbours
must be $a$. A cube vertex is $a$, making $U\cap D$
impossible in a value below $\omega$.
In $Q_5$, a one-letter graph vertex may be either $a$ or
$b$, so it restricts no colour but prevents $Q=d$.
Thus an odd cycle together with such a vertex makes
$Q=\omega$.
In $R_5$, a one-letter graph vertex can coexist only with
quadratic value $d$, giving the sum $b$. All graph vertices
must then be $a$, and a cube is impossible.
In $S_5$, the value $c$ requires that $D\cup L$ be coloured
$a$. If there is no cube, the assignment making every graph
vertex $a$ gives $d$ or $b$ and excludes every cubic absorption.

These alternatives also supply the separating assignments.
A target variable outside $V$ is assigned $\omega$, or
the required value $d$ or $c$ if it is an independent
one-letter summand. A pair not forced opposite has product
$d$ or $\omega$ under a same-colour assignment, so it is
not at most the right side $c$. In $P_5,Q_5$, choose the
inactive value $d$ to separate a one-letter vertex of
that colour. An active one-letter target is absorbed at
value $c$ precisely in the algebras where $a<c$.
When no one-letter graph vertex and no cube is present,
the assignment with right side $d$ excludes further
one-letter absorptions. Thus the criteria are necessary
as well as sufficient.

It remains to derive every listed absorption.
For $P_5$, adjoin the cubes at neighbours of one-letter
vertices by its extra rule. Substitution $y=x^2$ also
gives $\omega\leq x+x^3$. Then use the derivations of
$\mathcal N(G,M)$ already proved.
For $Q_5$, an odd cycle yields every vertex square.
At a one-letter graph vertex, $x+x^2=\omega$ proves the
extra zero case; $x\leq x^3$ supplies the remaining
one-letter absorptions.

For $R_5$, an overlapping one-letter summand and edge
give a square at every edge endpoint by its first extra
absorption. All quadratic words then follow from
$xy\leq x^2+y^2$, and the next rule gives all one-letter
vertices. Its last rule handles a cube in addition to
the overlap. Substituting $y=x^2,z=x$ also handles an
overlap between a one-letter summand and its own cube.

For $S_5$ with a cube present, the last rule adjoins cubes
at all of $L$, so the earlier graph proof applies.
Without a cube, the two-edge rule propagates one-letter
summands along even paths. Odd paths give within-component
opposite pairs. For different components with one-letter
marks, use $xz\leq x+xt+y+yz$, where $x,y$ are already
derived one-letter vertices and $xt,yz$ are incident edges.
Opposite one-letter marks in one component first give
their connecting edge, and then a square by
$x^2\leq x+y+xy$. Squares propagate to all vertices using
$\Sigma$, yielding all quadratic words. With one-letter
vertex $x$ present, the two-edge rule now gives
$z\leq x+x^2+xz$, so every vertex becomes a one-letter
absorption. An odd cycle is handled by the same square
derivation. This proves every clause, and therefore every
identity, from the displayed bases.
\end{proof}

\subsubsection{Quadratic graphs with a constant component}

Let $A_5,B_5,X_5,Y_5$ have underlying set
$\{a,b,d,c,\omega\}$ and commutative multiplication
\[
a^2=d,\qquad ad=c,\qquad b^2=c,
\]
all other products being $\omega$.
In $A_5$, the four elements below $\omega$ are incomparable.
For $B_5$ add $a<c$, for $X_5$ add $b<c$, and for $Y_5$
add both $a<c$ and $b<c$.
Triple products again take $c$ only at $(a,a,a)$;
multiplicative translations preserve all the specified
joins. Hence these are semirings.

Let $\Theta$ be obtained from $\Sigma$ by deleting
$x^2+y^3\approx u^4$ and $y^2\leq x^2+yz$, and adding
\begin{equation}\label{abs-local-square-cube}
x^2+x^3\approx u^4.
\end{equation}
Thus $\Theta$ has nine identities.

\begin{theorem}\label{abs-quartic-loop-bases}
The following are complete bases:
\begin{align*}
A_5:\quad&\Theta,\quad x\leq u^4,\quad u^4\leq x+xy;\\
B_5:\quad&\Theta,\quad x\leq x^3,\quad x^3\leq x+xy;\\
X_5:\quad&\Theta,\quad x\leq u^4,\quad x^2\leq x+xy,
              \quad x+x^3\approx u^4,\\
&y\leq x+x^2+y^2,\quad y\leq x^3+y^2;\\
Y_5:\quad&\Theta,\quad x\leq x^3,\quad z\leq x+xy+yz,
              \quad x^2\leq x+y+xy,\\
&x+xy+y^3\approx u^4,\quad y\leq x^3+y^2,
              \quad z\leq x+xy+z^2.
\end{align*}
Their respective numbers of identities are eleven, eleven,
fourteen, and fifteen.
\end{theorem}

\begin{proof}
The new feature is the quadratic value $b^2=c$.
A quadratic product is $d$ exactly on $(a,a)$, and is
$c$ exactly on $(a,d),(d,a),(b,b)$; a cube is $c$ exactly
at $a$. A sum of nonlinear words below $\omega$ is still
either $d$ or $c$. These rules directly verify $\Theta$.
For the additional rules, $x+xy<\omega$ requires
$(x,y)=(a,d)$ in $B_5$, $(b,b)$ in $X_5$, and one of
these pairs in $Y_5$; it never occurs in $A_5$.
A square is $d$ at $a$, $c$ at $b$, and $\omega$ elsewhere.
Substituting these alternatives verifies the rules
containing squares. A cube restricts its variable to $a$
whenever the right side is below $\omega$, which verifies
the remaining rules and \eqref{abs-local-square-cube}.

We describe all absorptions to prove completeness.
Normalize to $(U,E,D)$, writing $V=V_E\cup D$ and
$L=U\cap V_E$. A summand $\omega=u^4$ gives the constant
polynomial $\omega$. If no nonlinear summand is present,
only the given one-letter words are absorbed.
In the remaining cases make these preliminary reductions:
\begin{enumerate}
\item In $A_5$, an overlap $U\cap V$ gives $\omega$.
\item In $B_5$, adjoin $L$ to $D$ and delete $U\cap V$.
\item In $X_5$, an overlap $U\cap D$ gives $\omega$.
\item In $Y_5$, delete $U\cap D$.
\end{enumerate}
The displayed bases justify each reduction; for an overlap
with a cube in $A_5$, substitute $y=x^2$ in its last rule.
Use $U,D,L$ for the resulting sets below.

In an assignment with $Q=c$, each connected component
of $G=(V,E)$ either has a proper two-colouring by $a,d$,
or has every vertex assigned $b$.
These are all possibilities: a $b$ at one endpoint of an
edge forces $b$ at the other, and otherwise the endpoints
must alternate $a,d$. Cubes force their vertices to be $a$.
A one-letter graph vertex must be $b$ in $X_5$, and may
be $a$ or $b$ in $Y_5$. Thus the components have exactly
the following types.
\begin{enumerate}
\item A component containing a cube is forced to use the
$a,d$ colouring, and every cube mark must lie on the $a$
side. In $X_5$ it can contain no one-letter mark; in $Y_5$
all such marks must also lie on the $a$ side. If these
requirements fail, $Q=\omega$.
\item A component without cubes is forced to be constantly
$b$ if it is not bipartite, if it contains a one-letter
mark in $X_5$, or if it has one-letter marks on opposite
sides in $Y_5$.
\item Every remaining component may be constantly $b$
or may use a proper $a,d$ colouring. In $Y_5$, its
one-letter marks, if any, must have colour $a$ in the
second alternative.
\end{enumerate}
Call the first type fixed, the second constant, and the
third optional. If no zero case occurs, these alternatives
give the following exact criterion:
\begin{enumerate}
\item A cube is absorbed exactly at a vertex on the $a$
side of a fixed component.
\item A quadratic word is absorbed exactly when its
endpoints are opposite sides of one bipartite component,
are both in constant components, or are forced opposite
colours in two fixed components. All endpoints belong to $V$.
In the first clause, a component of constant type instead
uses the second clause, which also permits equal endpoints.
\item Besides the retained $U$, the one-letter absorptions
in $B_5,Y_5$ include the forced $a$ vertices.
In $Y_5$, they also include the same-side vertices as
one-letter marks in optional components.
In $X_5,Y_5$, vertices in constant components are included
if at least one cube or one-letter graph mark exists
anywhere in $Q$. There are no further absorptions.
\end{enumerate}

To verify necessity, the component choices above are
independent. A pair in different optional or mixed-type
components can be separated by assigning $b$ to one
component and $a,d$ to the other. A same-side pair in a
nonconstant component takes $d$ or $\omega$, not at most
the right side $c$. A variable outside $V$ is separated
by $\omega$, or by $c$ if it is an independent one-letter
summand. Cubic and one-letter assertions follow from
the assigned values and from which of $a,b$ are below $c$.
If $D=L=\varnothing$, assigning $a$ to every graph vertex
and $d$ to every one-letter variable gives $Q=d$.
This excludes any additional cubes and one-letter words
in that case. Conversely, checking the same component
choices and this possible $d$ value proves sufficiency.

We now derive the criterion. The identities of $\Theta$
retain odd-path shortening, even-path cube propagation,
and the rule connecting fixed components.
An odd cycle yields a square. The odd-path rule also gives
\begin{equation}\label{abs-local-square-propagation}
y^2\leq yx+x^2+xy\approx x^2+xy,
\end{equation}
so a square propagates throughout its component.
Two obtained squares give their connecting edge by
$xy\leq x^2+y^2$. A square and a cube at the same vertex
give $\omega$ by \eqref{abs-local-square-cube}.
Cubes on opposite sides of a bipartite component give
their connecting edge by an odd path, and then give
$\omega$ by $xy+x^3+y^3\approx\omega$.

In $X_5$, an overlapping one-letter summand and edge
give a square by $x^2\leq x+xy$. It propagates through
the component, proving the required constant type;
a cube in that component gives $\omega$.
In $Y_5$, opposite one-letter marks give their edge
and then a square by $x^2\leq x+y+xy$.
Opposite one-letter and cube marks give $\omega$ by
$x+xy+y^3\approx\omega$. Marks on the same side propagate
as one-letter words using $z\leq x+xy+yz$.
Consequently every constant component has all its
vertex squares derivable, every fixed component has
all its forced cubes, and all listed quadratic
absorptions follow.

Finally, a cube anywhere allows the rule
$y\leq x^3+y^2$ to make every constant-component
vertex a one-letter absorption in $X_5,Y_5$.
If instead there is a one-letter graph mark, use
$y\leq x+x^2+y^2$ in $X_5$, where the mark has already
given its square. In $Y_5$ use $z\leq x+xy+z^2$.
The rule $x\leq x^3$ gives the forced $a$ vertices
where appropriate. This derives every absorption,
and expansion in both directions proves completeness.
\end{proof}

\subsubsection{Directed parity and incoming edges}

Define $D_0,D_1$ on $\{a,b,d,c,\omega\}$ with the products
\[
a^2=d,\qquad ad=da=ba=c,
\]
all other products being $\omega$. In particular $ab=\omega$.
The additive order of $D_0$ is $d<b<\omega$, $a,c<\omega$,
with no other comparisons. For $D_1$ add $a<c$.
Triple products are $c$ only on $(a,a,a)$, and left and
right translations preserve the specified joins.
Thus both algebras are noncommutative semirings.

Let $\Pi$ consist of the following eighteen identities.
First include
\[
xyzt\approx u^4,\qquad xyz\approx x^3+y^3+z^3.
\]
For each $p\in\{xy,yx\}$ and $q\in\{yz,zy\}$ include
\begin{equation}\label{abs-directed-parity-paths}
xt\leq p+q+zt,\qquad z^3\leq x^3+p+q.
\end{equation}
These are four odd-path rules and four even-path rules.
The remaining eight identities are
\begin{align}
y^2&\leq x^2+yz, & y^2&\leq x^2+zy,\label{abs-directed-parity-squares}\\
xy&\leq x^2+y^2, & x^2+y^3&\approx u^4,\label{abs-directed-parity-zero1}\\
xy+x^3+y^3&\approx u^4, & xz&\leq x^3+y^3+yz,\label{abs-directed-parity-zero2}\\
xy&\leq y^3+z^3+xz, & xy&\leq y^3+z^3+zx.\label{abs-directed-parity-transfer}
\end{align}

\begin{theorem}\label{abs-directed-quartic-bases}
A complete basis for $D_0$ is
\[
\Pi,\qquad x\leq u^4,\qquad u^4\leq x+xy,\qquad u^4\leq x+yx.
\]
A complete basis for $D_1$ is
\[
\Pi,\qquad x\leq x^3,\qquad x^3\leq x+xy,\qquad x^3\leq x+yx.
\]
Each has twenty-one identities. Their word reversals are
complete bases for the multiplicative duals.
\end{theorem}

\begin{proof}
A quadratic word is $d$ exactly when both factors are $a$,
and is $c$ exactly on $(a,d),(d,a),(b,a)$. A cube is $c$
exactly at $a$. A sum of nonlinear terms below $\omega$
has every summand $d$, or every summand $c$.
In the latter case the underlying undirected graph
is properly coloured by $a$ and $\{d,b\}$, but a vertex
of the second colour with an incoming edge must be $d$.
In each odd-path rule the last vertex has an incoming
edge, so the endpoint product is still $c$.
Each even-path rule propagates the active colour $a$.
The same evaluation rules verify
\eqref{abs-directed-parity-squares}--\eqref{abs-directed-parity-transfer}.
In $D_0$ a one-letter summand and any incident edge sum
to $\omega$. In $D_1$ this sum, when below $\omega$,
requires the one-letter value $a<c$, so its cube can be
adjoined. This proves validity of all the identities.

Normalize to $(U,E,D)$, now retaining directions in $E$.
Write $V=V_E\cup D$ and $T=\{y:xy\in E\}$.
A long summand makes $Q=\omega$; in $D_1$ use
$x\leq x^3\leq x^9=\omega$.
In $D_0$ an overlap $U\cap V$ also makes $Q=\omega$;
an overlap with a cube follows by substituting $y=x^2$
in an incident-edge rule.
In $D_1$, adjoin cubes at $U\cap V_E$ and remove $U\cap V$
from the one-letter set, using $x\leq x^3$.
All subsequent data refer to this normalization.

Let $G$ be the underlying undirected graph, allowing
isolated cube vertices. Its absorption criterion is:
\begin{enumerate}
\item If $D=\varnothing$ and $G$ is not bipartite, exactly
the quadratic words on $V$ and the one-letter words in $U$
are absorbed.
\item If $D\ne\varnothing$ and no proper two-colouring
makes all cube vertices colour A, then $Q=\omega$.
\item Otherwise the absorbed cubes are precisely those
at vertices always colour A. A quadratic word $xy$ is
absorbed exactly when its endpoints belong to $V$, are
always opposite colours, and either $y\in T$ or $y$ is
always colour A. The one-letter absorptions are $U$ for
$D_0$; in $D_1$ they also include all forced colour A vertices.
\end{enumerate}
If there are no edges or cubes, only $U$ is absorbed.

The evaluation proof of the undirected parity criterion
applies with one extra condition. In an assignment giving
the right side $c$, a colour B vertex outside $T$ may
be $b$. As the last variable of a target word preceded
by an $a$, it then makes that target $\omega$.
A colour B vertex in $T$ is necessarily $d$.
This proves precisely the stated direction condition.
All other failed conditions are separated by same-colour
assignments, by the possible all-$a$ graph assignment
giving $Q=d$, or by assigning $\omega$ to a missing
variable and $c$ to an independent one-letter variable.
The values $a\not\leq c$ in $D_0$ and $a<c$ in $D_1$
give the one-letter distinction. These are all possible
assignments below $\omega$, proving necessity and sufficiency.

The direction condition must also be respected in the
derivation. We first prove that an undirected odd walk
from $x$ to $y$ yields $xy\leq Q$ whenever $y\in T$.
Choose an original incoming edge $vy$. Append the two
steps $y,v,y$ to the walk, so its final step is the
directed edge $vy$. Starting from that edge, repeatedly
use the four odd-path rules in
\eqref{abs-directed-parity-paths} to move two steps backwards.
The first two edges may have either orientation, while
the already derived last edge always points to $y$.
This eventually gives $xy$.

An odd cycle has an incoming vertex, so it gives a square
by this argument. The two rules
\eqref{abs-directed-parity-squares} propagate squares to all
original edge vertices, and then
\eqref{abs-directed-parity-zero1} supplies all quadratic words.
A cube anywhere yields $\omega$ with an obtained square.
In a bipartite component, the even-path rules propagate
cubes along every even path. If cube marks occur on both
sides, propagate one cube to a neighbour of the other
and use that original edge and the two endpoint cubes
in \eqref{abs-directed-parity-zero2}. This proves all zero cases.

In a compatible colouring, all forced colour A cubes
are now derivable. Consider a required word $xy$.
If $y$ is forced A, then $x$ is forced B; choose any
neighbour $z$ of $x$. It is forced A. The cubes at $y,z$
and the original edge $xz$ or $zx$ give $xy$ by
\eqref{abs-directed-parity-transfer}.
Otherwise $y\in T$. Within one component, the odd-walk
argument gives the word. Across components, $x$ is
forced A and $y$ forced B. An original incoming edge
$zy$ has $z$ forced A; the cubes at $x,z$ give
$xy\leq x^3+z^3+zy$ by \eqref{abs-directed-parity-zero2}.
Finally $x\leq x^3$ supplies the required one-letter
words in $D_1$. Every absorption is therefore derivable,
proving completeness.
\end{proof}

\subsection{Internal content in type 4592}\label{source:4592}
\subsubsection{The algebra and the identities}
We work in the signature $(+,\cdot)$ without constants, relative to
the axioms for additively idempotent semirings: addition is associative,
commutative and idempotent, multiplication is associative, and both
distributive laws hold. Write $u\leq v$ for $u+v\approx v$.

Let $A=\{z,a,b,c,t\}$ have the chain order
$z<a<b<c<t$, with addition the maximum. Put $K=\{b,c,t\}$ and
define $\rho(b)=b$, $\rho(c)=\rho(t)=c$. Multiplication is
\begin{equation}\label{s4592-eq:operation}
uv=\begin{cases}
z,&u\notin K\text{ or }v=z,\\
a,&u\in K\text{ and }v=a,\\
\rho(u),&u,v\in K.
\end{cases}
\end{equation}
These letters name elements only, not constants in terms.

\begin{theorem}\label{s4592-thm:main}
The following seven identities form a basis for the identities of $A$,
relative to the additively idempotent semiring axioms:
\begin{align}
xyzt&\approx xzyt, &&\label{s4592-eq:P1}\\
x^2y&\approx xy, &&\label{s4592-eq:P2}\\
xyx&\approx xy^2, &&\label{s4592-eq:P3}\\
xy&\leq x, &&\label{s4592-eq:P4}\\
xyz&\leq xz, &&\label{s4592-eq:P5}\\
xy&\leq xy^2+y, &&\label{s4592-eq:P6}\\
xyz&\leq xyz^2+yxz. &&\label{s4592-eq:P7}
\end{align}
The opposite semiring $A^{\mathrm{op}}$ has the basis obtained by
reversing every word in these identities.
\end{theorem}

For a nonempty word $w$, let $\ct(w)$ be its set of variables,
$\hd(w)$ its first variable, and $\sourcetail(w)$ its final variable. For a
word of length at least two, let $\sourceint(w)$ be the set of variables
in all positions except the final one. In particular the first
variable belongs to $\sourceint(w)$.

\begin{lemma}\label{s4592-lem:evaluation}
The algebra $A$ is an additively idempotent semiring. For a word $w$
of length at least two and an assignment $\alpha$, its value is
\begin{equation}\label{s4592-eq:value}
w(\alpha)=\begin{cases}
z,&\alpha(\sourceint(w))\not\subseteq K
       \text{ or }\alpha(\sourcetail(w))=z,\\
a,&\alpha(\sourceint(w))\subseteq K
       \text{ and }\alpha(\sourcetail(w))=a,\\
\rho(\alpha(\hd(w))),&\alpha(\ct(w))\subseteq K.
\end{cases}
\end{equation}
All seven identities in Theorem~\ref{s4592-thm:main} hold in $A$.
\end{lemma}
\begin{proof}
The operation in \eqref{s4592-eq:operation} is increasing in both coordinates:
its nonzero rows are $(z,a,b,b,b)$ and $(z,a,c,c,c)$ in chain order;
the other rows are constant $z$. Its columns are also increasing.
Hence it distributes over the maximum in both coordinates. For
associativity, if the middle factor is $z$ or $a$, both associations
of a triple product give $z$. If it is in $K$, checking whether the
last factor is $z$, $a$, or in $K$ shows that both associations agree.
Formula \eqref{s4592-eq:value} follows by induction on word length.

The first three identities follow from this formula: interior
permutations and duplications do not affect it, while both $xyx$
and $xy^2$ are nonzero precisely when both assigned elements are in
$K$, in which case both have value $\rho(x)$.

For \eqref{s4592-eq:P4}, a nonzero product has its first factor in $K$ and
has value $a$ or $\rho(x)$, each at most $x$. In \eqref{s4592-eq:P5},
whenever $xyz$ is nonzero its middle factor belongs to $K$, and
$xyz=xz$. Otherwise the desired inequality is automatic. For
\eqref{s4592-eq:P6}, if $y\in K$ then $xy=xy^2$, and if $y=a$ then
$xy\leq y$; the remaining case gives $z$. Finally, a nonzero $xyz$
in \eqref{s4592-eq:P7} has $x,y\in K$. If the last factor is in $K$,
then $xyz=xyz^2$; if it is $a$, then $xyz=yxz=a$.
\end{proof}

\subsubsection{Word manipulations}
Fix an ordering of the variables. Two kinds of nonlinear normal word
will be used. For a finite nonempty set $C$ and $h\in C$, the first
kind is
\[h\,\prod_{v\in C\setminus\{h\}}v\,h.\]
For a finite nonempty set $I$, $h\in I$, and $s\notin I$, the second
kind is
\[h\,\prod_{v\in I\setminus\{h\}}v\,s.\]
Products over empty sets are omitted substrings, not unit constants.
The products are written in the fixed variable order.

\begin{lemma}\label{s4592-lem:normal}
Identities \eqref{s4592-eq:P1}--\eqref{s4592-eq:P3} transform every nonlinear
word into a normal word. If its final variable occurs earlier, its
normal word has the first kind with the same head and content.
Otherwise its normal word has the second kind with the same head,
internal content, and final variable.
\end{lemma}
\begin{proof}
By applying \eqref{s4592-eq:P1} with nonempty words substituted for the
outer variables, any two adjacent interior positions can be
interchanged. Thus all positions other than the first and last
can be permuted. Repeated interior variables can be made adjacent
and reduced using \eqref{s4592-eq:P2}. An interior occurrence of the
head can similarly be moved next to the head and removed.

If the final variable does not occur earlier, these operations give
the second normal form. Suppose that the final variable $s$ occurs
earlier. If $s=h$, they give the first normal form directly. If
$s\neq h$, first reduce interior repetitions. Insert an extra $h$
immediately after the first position using \eqref{s4592-eq:P2} backwards.
Permute interior positions to obtain a word $hRh s^2$, where $R$
lists the other variables. Equation \eqref{s4592-eq:P3} changes this to
$hRhsh$. Move the redundant interior $h$ next to the initial $h$
and remove it using \eqref{s4592-eq:P2}. The result is $hRsh$, which can
be ordered into the first normal form. This also covers empty $R$.
\end{proof}

\begin{lemma}\label{s4592-lem:insertion}
Using \eqref{s4592-eq:P4} and \eqref{s4592-eq:P5}, inserting letters anywhere
after the first position of a word gives a word below the original.
Consequently, if $u$ has head $h$ and content contained in $C$,
every first-kind normal word with head $h$ and content $C$ is below
$u$, derivably from \eqref{s4592-eq:P1}--\eqref{s4592-eq:P5}.
\end{lemma}
\begin{proof}
Appending a letter uses \eqref{s4592-eq:P4}. An insertion between a
nonempty prefix $p$ and a nonempty suffix $q$ uses
$pvq\leq pq$, an instance of \eqref{s4592-eq:P5}. Repetition proves the
first assertion. For the second, append to $u$ a word listing all
variables of $C$, followed by $h$. The result has content $C$, head
$h$, and a repeated final variable. Lemma~\ref{s4592-lem:normal} transforms
it into the required first-kind normal word.
\end{proof}

\subsubsection{Two witnesses for absorption}
By distributivity a term is a nonempty finite sum of words. Write
such a polynomial as $P=\sum_{v\in L}v+\sum_{u\in N}u$, where
$L$ is its set of linear summands and $N$ its set of nonlinear words.
The nonlinear words may be put in normal form using
Lemma~\ref{s4592-lem:normal}.

\begin{lemma}\label{s4592-lem:criterion}
A linear word $y$ satisfies $A\models y\leq P$ if and only if
$y\in L$. For a nonlinear word $w$, put
$h=\hd(w)$, $C=\ct(w)$, and $I=\sourceint(w)$. Then
$A\models w\leq P$ if and only if both of the following hold:
\begin{enumerate}
\item[(H)] Either $h\in L$, or there is $u\in N$ with
$\hd(u)=h$ and $\ct(u)\subseteq C$.
\item[(A)] Either $L\cap C\neq\varnothing$, or there is $v\in N$
with $\ct(v)\subseteq C$ and $\sourceint(v)\subseteq I$.
\end{enumerate}
Every absorption satisfying this criterion is derivable from the
seven identities.
\end{lemma}
\begin{proof}
If $y\notin L$, assign $y=t$ and every other variable $z$. Linear
summands have value $z$, and nonlinear summands have value at most
$c$. Thus $y\not\leq P$. The converse is immediate.

For a nonlinear word, first prove necessity of (H). If (H) fails,
assign $h=c$, every other variable in $C$ the element $b$, and
every variable outside $C$ the element $z$. The word $w$ has value
$c$. Any nonlinear summand containing an outside variable has
value $z$; any summand entirely in $C$ has head different from
$h$ and value $b$. No linear summand equals $h$. Hence $P\leq b$,
contradicting $w\leq P$.

If (A) fails and $\sourcetail(w)\notin I$, assign all variables in $I$
the element $b$, the final variable the element $a$, and all
outside variables $z$. Then $w=a$. There is no linear summand in
$C$. Each nonlinear summand either contains an outside variable,
or has the final variable of $w$ in an internal position. In both
cases its value is $z$. Therefore $P=z$. If $\sourcetail(w)\in I$, then
$I=C$; assigning all variables in $C$ the element $b$ and all others
$z$ gives $w=b$ and $P=z$, again a contradiction. This proves necessity.

We prove sufficiency by equational derivation. If $h\in L$,
\eqref{s4592-eq:P4} gives $w\leq h\leq P$. Otherwise choose the witness
$u$ in (H). If the final variable of $w$ is repeated,
Lemma~\ref{s4592-lem:insertion} gives $w\leq u\leq P$.

It remains to consider $s=\sourcetail(w)\notin I$. The word $ws$ has
the first-kind normal form with head $h$ and content $C$, so
Lemma~\ref{s4592-lem:insertion} gives
\begin{equation}\label{s4592-eq:guarded}
ws\leq u\leq P.
\end{equation}
For every $g\in I$, denote by $w_g$ a word with head $g$, internal
content $I$, and final variable $s$, in the second normal form.
If $g=h$, then $w_g\approx w$. If $g\neq h$, write a nonempty
word $Y$ with head $g$ that lists $I\setminus\{h\}$. An instance
of \eqref{s4592-eq:P7}, followed by normalization, gives
\begin{equation}\label{s4592-eq:switch}
w\approx hYs\leq hYs^2+Yhs\approx ws+w_g.
\end{equation}

If (A) has a linear witness $g\in L\cap I$, then
$w_g\leq g$ by \eqref{s4592-eq:P4}. Combine this with
\eqref{s4592-eq:guarded} and \eqref{s4592-eq:switch}; the case $g=h$ is immediate.
If its linear witness is $s$, write $w=qs$ with $q$ nonempty and
use \eqref{s4592-eq:P6} to obtain $w\leq ws+s\leq P$.

Finally let $v\in N$ witness (A), and put $g=\hd(v)$. Then
$g\in I$. The variable $s$ either does not occur in $v$ or occurs
only at its final position. If it is absent, append a word listing
$I$ followed by $s$. If it is present, insert a word listing $I$
just before that final position. In both cases the insertion lemma
and normalization give $w_g\leq v\leq P$. Equations
\eqref{s4592-eq:guarded} and \eqref{s4592-eq:switch} now give $w\leq P$; again
$g=h$ is immediate. This exhausts the possibilities in (A).
\end{proof}

\subsubsection{Completeness and duality}
\begin{proof}[Proof of Theorem~\ref{s4592-thm:main}]
Soundness was proved in Lemma~\ref{s4592-lem:evaluation}. Let $p\approx q$
be any identity valid in $A$. Distribute and normalize to obtain
polynomials $P,Q$. Each summand of $P$ is below $Q$ in $A$ and
therefore, by Lemma~\ref{s4592-lem:criterion}, is provably below $Q$ from
the seven identities. Thus $P\leq Q$ is derivable. The symmetric
argument gives $Q\leq P$ and hence $P\approx Q$. This proves
completeness for any finite set of variables. Reversing every word
converts identities and equational derivations for $A$ into those
for $A^{\mathrm{op}}$, proving the assertion about the opposite.
\end{proof}

\subsection{The final two positions in type 4609}\label{source:4609}
\subsubsection{The algebra and the main theorem}
We use the signature $(+,\cdot)$ with no distinguished constants. An
additively idempotent semiring is an algebra in which addition is
associative, commutative and idempotent, multiplication is associative,
and both distributive laws hold. We write $u\leq v$ as an abbreviation
for the identity $u+v\approx v$, or for the corresponding order on
elements of a particular semiring. An identity basis below is relative
to these semiring axioms. Appending the finitely many axioms gives an
ordinary finite identity basis in the stated signature.

Let $R=\{a,b,c,e,t\}$. Its addition is the join operation of the
semilattice with
\[a<b<c<t,\qquad e<t,\]
where $e$ is incomparable with $a,b,c$. Define multiplication by
\begin{equation}\label{s4609-eq:product}
uv=\begin{cases}
v,&v\neq b,\\
a,&v=b\text{ and }u=e,\\
c,&v=b\text{ and }u\neq e.
\end{cases}
\end{equation}
For clarity the two operation tables are
\[
\begin{array}{c|ccccc}
+&a&b&c&e&t\\\hline
a&a&b&c&t&t\\
b&b&b&c&t&t\\
c&c&c&c&t&t\\
e&t&t&t&e&t\\
t&t&t&t&t&t
\end{array}
\qquad
\begin{array}{c|ccccc}
\cdot&a&b&c&e&t\\\hline
a&a&c&c&e&t\\
b&a&c&c&e&t\\
c&a&c&c&e&t\\
e&a&a&c&e&t\\
t&a&c&c&e&t
\end{array}.
\]
The symbols in these tables name elements, not constants of the signature.

\begin{theorem}\label{s4609-thm:basis}
The identities of $R$ have the following finite basis relative to the
axioms for additively idempotent semirings:
\begin{align}
xyz&\approx yz, &&\label{s4609-eq:B1}\\
x+x^2&\approx x^2, &&\label{s4609-eq:B2}\\
xy+y^2&\approx y^2, &&\label{s4609-eq:B3}\\
x+xy+y^2&\approx x+xy. &&\label{s4609-eq:B4}
\end{align}
\end{theorem}

\begin{lemma}\label{s4609-lem:valid}
The algebra $R$ is an additively idempotent semiring and satisfies
\eqref{s4609-eq:B1}--\eqref{s4609-eq:B4}.
\end{lemma}
\begin{proof}
The addition is a semilattice operation by construction. For every
$u,v$, we have $uv=e$ if and only if $v=e$. If $w\neq b$, both
$(uv)w$ and $u(vw)$ equal $w$. If $w=b$, the element $vb$ is $a$ or
$c$, so $u(vb)=vb$; moreover $(uv)b=vb$ because $uv=e$ if and only
if $v=e$. This proves associativity and $(uv)w=vw$.

For $v\neq b$ the right translation $u\mapsto uv$ is constant.
The right translation by $b$ sends $e$ to $a$ and every other element
to $c$. It preserves joins: the join of $e$ and a different element
is $t$, and the join of two elements different from $e$ is again
different from $e$. Each left translation fixes $a,c,e,t$ and sends
$b$ to $a$ or $c$. Its restriction to the chain $a<b<c<t$ preserves
joins, and joins with $e$ are also preserved. Thus both distributive
laws hold.

All elements other than $b$ are multiplicatively idempotent, while
$b^2=c>b$, proving \eqref{s4609-eq:B2}. If $y\neq b$, then $xy=y=y^2$.
If $y=b$, then $xy\in\{a,c\}$ and $y^2=c$. This proves
\eqref{s4609-eq:B3}. Finally $y^2\leq x+xy$ is immediate when $y\neq b$.
For $y=b$, either $x\neq e$ and $xy=c$, or $x=e$ and $x+xy=e+a=t$.
This proves \eqref{s4609-eq:B4}.
\end{proof}

\subsubsection{An exact absorption criterion}
Every term in the signature expands, by distributivity, into a
nonempty sum of nonempty words. By \eqref{s4609-eq:B1}, a word of length
at least three is equal to the word consisting of its final two
letters. Hence every term is provably equal to a polynomial
\begin{equation}\label{s4609-eq:normal}
P=\sum_{z\in L}z+\sum_{(u,v)\in E}uv,
\end{equation}
where $L$ is a finite set of variables, $E$ is a finite set of ordered
pairs of variables, and at least one of $L,E$ is nonempty. Empty
displayed sums mean that those summands are omitted; they do not
introduce a zero constant. Loops $(u,u)$ are allowed.

Define the set of terminal variables and the internal predecessors of
$y$ by
\[
T(P)=L\cup\{v:(u,v)\in E\},\qquad
I_y(P)=\{h\in T(P):(h,y)\in E\}.
\]
All variable symbols in the next statement may coincide.

\begin{lemma}\label{s4609-lem:absorption}
For a polynomial $P$ of the form \eqref{s4609-eq:normal},
\begin{align*}
R\models y\leq P
&\quad\Longleftrightarrow\quad y\in L\text{ or }I_y(P)\neq\varnothing,\\
R\models xy\leq P
&\quad\Longleftrightarrow\quad (x,y)\in E\text{ or }I_y(P)\neq\varnothing.
\end{align*}
Whenever either right-hand condition holds, the corresponding
inequality is derivable from \eqref{s4609-eq:B1}--\eqref{s4609-eq:B4}.
\end{lemma}
\begin{proof}
First prove the derivability assertions. A summand is below its sum.
Suppose instead that $h\in I_y(P)$. There is a summand $v$ of $P$
whose final letter is $h$, and $hy$ is also a summand. Thus $v=h$
or $v=zh$ for some $z$. In both cases \eqref{s4609-eq:B1} gives $vy\approx hy$.
Substituting $v$ for $x$ in \eqref{s4609-eq:B4} yields
\[y^2\leq v+vy\approx v+hy\leq P.\]
Equations \eqref{s4609-eq:B2} and \eqref{s4609-eq:B3} now give, respectively,
$y\leq P$ and $xy\leq P$. Validity in $R$ follows from Lemma~\ref{s4609-lem:valid}.

For necessity we provide separating evaluations. If $y\notin T(P)$,
assign $y$ the element $t$ and every other variable the element $a$.
Every summand of $P$ has terminal value $a$ and therefore has value
$a$, whereas both $y$ and $xy$ have value $t$. Neither inequality holds.

Suppose that $y\in T(P)$, $y\notin L$, and $I_y(P)=\varnothing$.
Assign $y=b$, all variables in $T(P)\setminus\{y\}$ the element $a$,
and all variables outside $T(P)$ the element $e$. Every edge ending
at $y$ has its first variable outside $T(P)$, so its value is $eb=a$.
Every other summand has terminal value $a$ and also has value $a$.
Consequently $P=a$ while $y=b\not\leq a$. This proves necessity
for linear words.

Finally suppose that $y\in T(P)$, $(x,y)\notin E$, and
$I_y(P)=\varnothing$. Start with the preceding evaluation, now
allowing $y\in L$. If $x\notin T(P)$, change its value to $a$.
This change cannot affect an edge ending at $y$, since $(x,y)$ is
absent; all other edges end at a variable of value $a$. Thus every
edge of $P$ has value $a$, and its linear summands have values in
$\{a,b\}$. Hence $P\in\{a,b\}$. The value of $x$ is not $e$: it
is $b$ if $x=y$, and otherwise it is $a$. Therefore $xy=c$, and
$c\not\leq P$. The argument includes the case $x=y$.
\end{proof}

\subsubsection{Completeness}
\begin{proof}[Proof of Theorem~\ref{s4609-thm:basis}]
Let $p\approx q$ be any identity holding in $R$, with an arbitrary
finite set of variables. Use distributivity and \eqref{s4609-eq:B1} to
obtain normal polynomials $P,Q$. Every summand $w$ of $P$ satisfies
$R\models w\leq Q$. Lemma~\ref{s4609-lem:absorption} proves this
inequality from the four identities. Taking the finite sum of these
derivations gives $P\leq Q$. Interchanging $P$ and $Q$ gives
$Q\leq P$. The semilattice axioms then give $P\approx Q$, and
therefore $p\approx q$. Soundness was proved in Lemma~\ref{s4609-lem:valid}.
No bound on the number of variables has been used.
\end{proof}

\begin{corollary}
The opposite semiring $R^{\mathrm{op}}$ is finitely based. A basis
is obtained from \eqref{s4609-eq:B1}--\eqref{s4609-eq:B4} by reversing every word.
\end{corollary}
\begin{proof}
Evaluation in the opposite semiring reverses products and leaves
addition unchanged. Reversing every word is an involution on
polynomial identities and on their equational derivations.
\end{proof}

\begin{theorem}\label{thm:correct643}
A relative basis for \(\A{643}\) consists of
\[
xyz\approx xy,\qquad x\law x^2,\qquad xy\law x^2,\qquad x^2\law y+xy .
\]
\end{theorem}
\begin{proof}
Addition is the join on the two chains \(1<0<4\), \(3<2<4\), with common top. Products have rows \(00000,00000,33444,33333,44444\). Their values lie in \(\{0,3,4\}\), where right multiplication fixes each element. Replacing a right factor \(y\) by \(yz\) preserves its membership in \(\{0,1\}\). This gives the first law. The squares are \(0,0,4,3,4\) in order, proving the next two laws. For the last law \(xy=x^2\) except possibly at \(x=2,\ y\in\{0,1\}\), when \(y+xy=y+3=4\).

Let \(L\) denote the literal letters of a reduced polynomial \(P\), and \(E\) its ordered pairs. Put \(H=L\cup\{u:(u,v)\in E\}\) and \(J_x=\{z:(x,z)\in E\}\cap H\). The condition for \(x\law P\) is \(x\in L\) or \(J_x\ne\varnothing\). The condition for \(xy\law P\) is \((x,y)\in E\) or \(J_x\ne\varnothing\), including \(x=y\). Indeed \(z\in J_x\) supplies \(xz\) together with \(z\) or some \(zt\). The fourth law gives \(x^2\law z+xz\). Substituting \(zt\) for \(z\) gives \(x^2\law zt+xz\) by the first law. Now use the second or third law.

Conversely, when \(J_x\) is empty assign \(x=2\), every endpoint \(z\) of an edge \(xz\) the value \(0\), and all remaining letters the value \(3\). The endpoint set excludes \(x\) and is disjoint from \(H\). Every upper edge has value \(3\). For a missing literal \(x\), all upper literals also have value \(3\), which does not dominate \(2\). For a missing pair \(xy\) the left value is \(4\), while the right side lies in \(\{2,3\}\). Thus every valid absorption follows from the four laws. Distributing both sides of an identity and applying the result in both directions proves completeness.
\end{proof}

\section{Small source algebras and periodic examples}\label{sec:direct-sources}

\subsection{Truncated addition and two end conditions}
Write \(S_{d,j}\) for the source with the two row strings in Appendix~\ref{app:register}.
The source numbering here is independent of the numbering of the five-element algebras. Zhao, Ren, Crvenkovi\'c, Shao and Dapi\'c proved finite basability for all \(60\) three-element ai-semirings other than \(S_7\), up to isomorphism \cite{Three}.

\begin{proposition}\label{prop:three-direct}
The five sources \(S_{3,24},S_{3,37},S_{3,50},S_{3,51},S_{3,53}\) in case~\ref{grp:130} have finite bases.
\end{proposition}
\begin{proof}
For \(S_{3,24}\), replace each label \(a\) by its cost \(2-a\). The sum of two costs is their maximum and their product is the ordinary sum truncated at \(2\). The product of any three costs is the maximum of their three pairwise products. Indeed, a total of at least two is witnessed by a single cost two or by two positive costs. If the total is zero or one, some pair already has that total. The identity \(xyz\approx xy+xz+yz\), applied to the first two letters and the remaining suffix, reduces every word of length greater than two to a sum of shorter words. The additive order has a bottom. Theorem~\ref{thm:bounded} applies with \(l=2\).

For \(S_{3,37}\) addition is the minimum in the printed labels. A product of two or more letters has value \(0\) when all their values belong to \(\{0,1\}\), and value \(2\) otherwise. A relative basis is
\[
 xy\approx yx,\qquad x^2y\approx xy,\qquad x\law x^2,\qquad xyz\law xy .
\]
These laws follow from the two possible product values. A word of length at least two is determined, among such words, by its content. The first two laws prove this assertion by deleting repeated letters while retaining a square on one-letter content. In an absorption with nonlinear left word (here meaning a word of length at least two), assign \(1\) to its letters and \(2\) to the others. This forces a nonlinear right word whose content is contained in that of the left word. From this witness, the last inequality adjoins all missing letters, and the word laws recover the target. For a literal \(x\), the same assignment forces the literal \(x\) or a right word \(x^k\), \(k\ge2\). In the latter case use \(x\law x^2\).

For \(S_{3,50}\), the printed labels are costs. Addition is minimum and multiplication is truncated ordinary addition with cap \(2\). A basis is
\[
 xy\approx yx,\qquad x^3\approx x^2,\qquad xy\law x,\qquad
 xy\approx x^2y+xy^2 .
\]
The last identity follows by taking the smaller of the costs \(2a+b\) and \(a+2b\), with truncation. For a valid \(w\law P\), assigning zero to its letters and two to all others forces a contained summand. Giving a letter \(x\) that occurs once in \(w\) cost one instead forces a contained summand with at most one \(x\). Normalize positive exponents to one or two. If at least two exponents are one, partition the word into two nonempty blocks, each containing such a letter. The last identity doubles either block and expresses the word as the join of the two results. Repeating terminates along each branch because the number of single letters decreases. Every resulting word has the original content and at most one single letter. For such a word choose the contained witness just obtained, using its single letter if it has one. After normalization none of the witness exponents exceeds the corresponding exponent in the word. Apply \(xy\law x\) and commutativity to append the missing factors. Each resulting word is thus below \(P\), as is their join.

In \(S_{3,51}\) a word evaluates to its last letter when every preceding letter has value \(0\), and to \(2\) otherwise. A relative basis is
\[
 xyz\approx yxz,\qquad x^2y\approx xy,\qquad xy\law x,\qquad xy\law y .
\]
The word identities identify exactly the written forms needed here: two words with the same last letter and the same set of letters before it are provably equal. Setting target letters to \(0\) and outsiders to \(2\) forces a contained right word \(v\). If the last letter \(t\) occurs earlier in the target \(w\), the prefix set of \(vw\) is already that of \(w\). Thus \(w\approx vw\law v\). If \(t\) occurs just once, give it value \(1\), leaving the other target letters at \(0\) and outsiders at \(2\). Some contained \(v\) omits \(t\) or has its sole \(t\) in last position. In the first case \(vw\) has the required prefix set and is below \(v\). In the second, prepend the target prefix to \(v\) when nonempty. The result has the required prefix set and is again below \(v\). These arguments include a literal target. The table of \(S_{3,53}\) transposes the product of \(S_{3,51}\), so Lemma~\ref{lem:transfer} supplies its basis.
\end{proof}

\begin{proposition}\label{prop:seven-bases}
The seven extensions in case~\ref{grp:113} are finitely based.
\end{proposition}
\begin{proof}
The sources \(S_{3,50},S_{3,51},S_{3,53}\) occur in Proposition~\ref{prop:three-direct}.
All triple products are the additive bottom in \(S_{3,60}\) and \(S_{4,471}\). The identity \(xyz\approx t^3\) reduces each longer word to a cube, giving Theorem~\ref{thm:bounded} with \(l=3\).

For \(S_{4,443}\), identify the labels \(1,2,3,0\) with \((0,0),(1,0),(0,1),(1,1)\), respectively. Addition is componentwise maximum and multiplication has the form
\[
                     (a,b)(c,d)=(ac,ad).
\]
A word ending in \(t\) has second coordinate \(b_t\prod a_x\), where the product ranges over its preceding letters. The empty product here equals \(1\) in the Boolean calculation. Give \(b_t\) value one, all other second coordinates zero, and give first coordinate one exactly to the letters in the target prefix. This forces a right word with last letter \(t\) whose prefix set lies in the target prefix set. Conversely prepend any missing prefix letters to this witness, using \(xy\law y\). The two identities \(xyz\approx yxz\) and \(x^2y\approx xy\) then recover the target. These three laws hold by the coordinate formula and therefore form a basis. The product of \(S_{4,456}\) is the transpose.

For Theorem~\ref{abs-regularization} take \(d(x,y)=x+h(x,y)\). In the source order
\[
 (3,50),(3,51),(3,53),(3,60),(4,443),(4,456),(4,471)
\]
the respective choices of \(h\) are \(xy,xy,xy,xy,yx,xy,x^2y\). In the first three tables the product is below both factors. In \(S_{3,60}\) only the product of two labels \(0\) can be above bottom, and that product lies below \(0\). The coordinate formula and its transpose give the next two inequalities. In the last table every triple product is bottom. Hence \(h(x,y)\le x\) throughout. The theorem gives a finite basis for each extension and the two variety inclusions of Appendix~\ref{app:register} transfer it to each listed algebra.
\end{proof}

\subsection{Flat groups of prime exponent}\label{sec:prime-flat}
Let \(p\) be a prime and let \(G\) be a nonzero finite-dimensional vector space over \(\mathbb F_p\). Form \(F(G)\) by adjoining an element \(\infty\). The product on \(G\) is vector addition and \(\infty\) absorbs products. The semiring sum of unequal elements is \(\infty\). Write \(F(G)^0\) for the further extension by an additive bottom which annihilates products. Let \(D\) denote the two-element lattice with join as addition and meet as product. For each word \(u\), let \(m(u)\) be its multiplicity vector modulo \(p\), with a common finite coordinate set containing the letters in the identity. Put \(P_w=\{u\in P:\ct(u)\subseteq\ct(w)\}\).

\begin{theorem}\label{thm:prime-flat}
A basis for \(F(G)\), for \(F(G)\times D\), or for \(F(G)^0\), respectively, is obtained as follows. In all three cases take
\[
 xy\approx yx,\qquad x^{p+1}\approx x,\qquad
       x_0x_1\cdots x_p\law x_0+x_1+\cdots+x_p .
\]
For \(F(G)\) add \(y\law x^p y\). For \(F(G)\times D\) add \(yz\law x^p yz+y\). The three common laws alone form a basis for \(F(G)^0\).
\end{theorem}
\begin{proof}
In \(F(G)\) the criterion for \(w\law P\) is
\[
 \ct(w)\subseteq \bigcup_{u\in P}\ct(u),\qquad
          m(w)\in\operatorname{aff}_{\mathbb F_p}\{m(u):u\in P\}.
\]
For the first condition send an uncovered target letter to \(\infty\) and all other letters to the group identity. For the second fix \(u_0\in P\). If \(m(w)-m(u_0)\) is outside the linear span of the \(m(u)-m(u_0)\), extend a basis of the latter span by that vector and then to the whole coordinate space. A linear functional which vanishes on the span and is one on the added vector gives a separating assignment to a cyclic subgroup of \(G\). Conversely, the only upper evaluations to check have all summands equal to a single group value. Every upper letter then has a group value, so the content condition puts the target in the group as well. Its value is the affine combination of the identical upper values and hence equals that common value.

The \(D\) coordinate requires \(P_w\ne\varnothing\). Indeed, setting exactly the target letters to one tests necessity, and a contained word has meet at least the meet of the target. For \(F(G)^0\) the two conditions for \(F(G)\) must hold with \(P_w\) in place of \(P\) and \(P_w\) must be nonempty. Assign outsiders to the new bottom to see necessity. Conversely, when the target avoids that bottom, the contained summands evaluate in \(F(G)\) and suffice. When the target equals bottom the inequality holds under every assignment.

The two word laws normalize each positive multiplicity to its representative in \(\{1,\ldots,p\}\). The displayed absorption holds since \(p+1\) identical group entries multiply to that entry, an upper value \(\infty\) dominates all elements, and a bottom factor in the extended table gives a bottom product. The deletion law for the flat table holds because multiplication by a group \(p\)th power fixes every element, while an \(\infty\) factor makes the upper product \(\infty\). Its guarded version holds also in \(D\), where \(yz\le y\).

Choose representatives \(0,\ldots,p-1\) for the affine coefficients. Their sum is a positive integer congruent to one modulo \(p\). Starting with one chosen word, multiply in \(p\) more of the chosen words at a time. The displayed absorption puts every successive product below \(P\). The final word has precisely the required residue vector. For each target letter still missing, append \(p\) copies of a witness containing it, by the same absorption. In the bottom extension all chosen witnesses are in \(P_w\) and the word laws now give the target. In the flat table every exterior letter occurs a multiple of \(p\) times. Commute and erase these blocks using \(y\law x^p y\), always retaining a nonempty word containing the target content.

For the product with \(D\), fix a contained summand \(v\). Before erasing exterior blocks, increase the positive multiplicities of its letters by multiples of \(p\), using \(x^{p+1}\approx x\) backwards, until the interior product can be written \(vt\) with \(t\ne\varepsilon\). At a deletion stage put all remaining factors other than one exterior \(p\)-block and the designated copy of \(v\) into the second factor of the left member of \(yz\law x^p yz+y\). That factor contains \(t\), so the substitution is legitimate. The right side is the preceding word plus \(v\), both below \(P\). This removes the exterior blocks one by one. The final product has the correct content and residues, hence equals \(w\) by normalization. This proves completeness on arbitrary alphabets.
\end{proof}

\begin{corollary}\label{cor:prime-examples}
Cases~\ref{grp:111} and~\ref{grp:328} are finitely based.
\end{corollary}
\begin{proof}
In \(\A{7446}\) the labels \(0,1,2,3\) are the elements of \(\mathbb F_2^2\) in binary order and \(4=\infty\). Apply Theorem~\ref{thm:prime-flat}.
For \(p=3\), write \(F=F(\mathbb F_3)\) with group labels \(0,1,2\) and top \(3\). Let \(M\) have both operations equal to join on \(0<1\).
The following pairs of maps, written in order of the five input labels, present subdirect embeddings. Each map preserves both operations by the two tables. Each coordinate is onto, and the five pairs in each row are distinct.
\[
\begin{array}{rccc}
 i & \A{i}\longrightarrow F & \text{other factor} &\text{second map}\\
 1639&30123&D&00001\\
 1655&30123&D&01111\\
 2796&30123&M&00001\\
 4847&00123&M&10111\\
 7907&01233&M&11101
\end{array}
\]
Here \(M\) is already the subalgebra of \(F\) on the group identity and \(\infty\). Thus the last three rows generate \(\V(F)\), and the first two generate \(\V(F\times D)\). The table of \(\A{7455}\) is \(F^0\), with group labels \(0,2,3\), top \(4\) and new bottom \(1\).

In \(\A{7555}\) the original group has labels \(1,2,3\), and \(0\) is an additional multiplicative identity immediately above the group identity \(1\) in the additive order. Sending \(0\) to \(1\) and fixing the original flat algebra defines a retraction to \(F\). The subset \(\{0,1\}\) is a copy of \(D\). Conversely, let \(w\law P\) hold in \(F\) and \(D\). Choose a contained right word. The retraction decides the inequality unless both sides map to the group identity, the only fibre with two elements. Within that fibre the only possible failure has target \(0\) and right side \(1\). A product has value \(0\) exactly when every factor equals \(0\). The contained word then also evaluates to \(0\), ruling out this failure. Consequently \(\A{7555}\) generates \(\V(F\times D)\). This gives the seven applications of case~\ref{grp:111} directly. The general periodic and flat-algebra results give the historical context \cite{ComPeriodic,Flat}.
\end{proof}

\subsection{Matrix units on two points}\label{sec:units}
For \(U=\{0,1\}\), let \(E\) have elements \(e_{ij}\) (\(i,j\in U\)) and \(\infty\). Unequal summands add to \(\infty\), which absorbs multiplication. The product \(e_{ij}e_{kl}\) is \(e_{il}\) for \(j=k\), and \(\infty\) otherwise. Let \(E_\triangle\) be the subalgebra with \(i\le j\). These are the flat Brandt and triangular sources usually denoted by \(B_2\) and \(B_0\) \cite{ShB2,ShB0}. The notation \(B_0\) used in the title of \cite{WWL} has a different meaning.

\begin{theorem}\label{thm:unit-paths}
For either \(E\) or \(E_\triangle\), take all valid word absorptions on twelve fixed letters with at most three right words and with every word of length at most four. For \(E_\triangle\) include also \(uvu\approx u^2v^2\). Together with \(\ax\) these laws give a finite basis.
\end{theorem}
\begin{proof}
For a nonempty polynomial \(P\) make a bipartite graph with rows \(H,R_x\) and columns \(T,L_x\), for \(x\in\ct(P)\). A word with head \(a\) and tail \(b\) contributes the edges \(H L_a\), \(R_b T\), and one edge \(R_xL_y\) for each adjacent pair \(xy\). Identify all vertices of each connected component. Each letter \(x\) now gives a directed edge from \(L_x\) to \(R_x\) in the quotient. Call the resulting labelled directed graph \(K\). A path may repeat edges or vertices. For \(E\) use \(K\) itself. For \(E_\triangle\) collapse its strongly connected components to obtain \(K'\). In either case, \(w\law P\) holds exactly when \(w\) uses only upper letters and spells a path from \(H\) to \(T\) in the indicated graph.

Indeed, for a right side below \(\infty\) all its words must have one common pair of endpoints and must be composable. This is precisely a labelling of the vertices of \(K\) by \(0,1\), assigning \(e_{ij}\) to an edge with labels \(i,j\). For \(E_\triangle\) the labelling must increase weakly along edges. It is constant on each strongly connected component. Any two distinct vertices of \(K\) can be separated by a labelling for \(E\). Two distinct components of \(K'\) can also be separated by an increasing labelling: one of the two principal upsets contains exactly one of them. Thus each failed endpoint equality can be separated while all upper words still have the same finite value. A missing letter is separated by assigning it \(\infty\) and giving all other letters \(e_{00}\). This proves the criterion in both directions.

Call a word available if the proposed laws prove that it lies below \(P\). All summands of \(P\) are available. A cut of an available word witnesses an edge \((r,c)\) of the bipartite graph in the following sense. At the beginning its row is \(H\), and elsewhere it is \(R_x\) for the preceding letter \(x\). At the end its column is \(T\), and elsewhere it is \(L_y\) for the following letter \(y\). Every initial edge has a witness. Suppose that \((s,d),(r,d),(r,c)\) have witnesses. Join the prefix through row \(s\) of the first witness to the suffix from column \(c\) of the third witness. Unless \(s=H,c=T\) this produces a nonempty word, witnessing \((s,c)\). It is available. To see this syntactically, retain in each of the three witnesses the letters immediately beside its cut, and replace each remaining nonempty prefix or suffix by a fresh variable. Preserve the equalities of the letters forced by the three specified edges. In the three formal words so obtained all cuts have the same endpoint under any assignment giving the three words a common finite value. In particular, the new prefix and suffix compose and have the required extreme endpoints. This gives a valid absorption with three right words, each of length at most four. Its left word has length at most four and the total number of formal letters is at most twelve. The original witnesses recover the instance by nonempty substitutions. Empty pieces were simply absent in the pattern.

Each row and column in the same component, except the pair \(H,T\), consequently has a witness. Take a shortest simple path in the initial bipartite graph between them. For a path of length three apply the preceding rule. For odd length at least five, among its consecutive subpaths of length three at most one has \(H,T\) for endpoints, since the path is simple. Shorten another such subpath using its new witnessed edge. Induction on path length gives the assertion. Equivalently the induction may be performed in the graph of currently witnessed edges, whose initial components can only gain such shortcuts. No word of length zero is required.

Take an actual nonempty \(H\)-to-\(T\) path in \(K\). Its first letter has a witnessed initial cut, every adjacent pair a witnessed internal cut, and its last letter a witnessed terminal cut. Two available words \(axb,cxd\) give the available word \(axd\), since a common finite evaluation forces agreement at the occurrence of \(x\). Optional empty contexts \(a,b,c,d\) are omitted when writing the corresponding pattern. This is a member of the finite family. Start with the head witness and splice at the last letter of the prefix already built, taking the desired next letter from the appropriate pair witness. A last splice with the terminal witness removes the unused suffix. This proves every absorption prescribed by \(K\).

In \(E_\triangle\), a square is a loop or \(\infty\). Hence inserting \(u^2\) before a word, after it, or between two nonempty factors can only increase its value in the additive order. The three required laws \(x\law u^2x,\ x\law xu^2,\ xy\law xu^2y\) belong to the indicated finite family. The word \(uvu\) has finite value only when \(u\) and \(v\) evaluate to the same loop, which also characterizes a finite value of \(u^2v^2\). Suppose now that a path spells \(w\) in \(K'\). At each gap, including the endpoints \(H,T\), the two vertices to be connected in \(K\) are either equal or mutually reachable. For unequal vertices choose a path from the first to the second with nonempty label \(u\) and a return path with nonempty label \(v\). Insert \(uvu\) in this gap. In the equational deduction insert first the squares \(u^2,v^2\), and then apply the word identity. Doing this at each gap gives \(w\law W\), where \(W\) is an actual \(H\)-to-\(T\) path in \(K\). Distinct gaps can be processed successively in context. Thus \(w\law W\law P\), which completes the triangular case as well.
\end{proof}

In case~\ref{grp:132} the elements of the four-element source, in table order, are \(e_{00},e_{01},e_{11},\infty\). In case~\ref{grp:364} the table of \(\A{8353}\) gives \(E\), with \(0,1,2,3\) the four matrix units in row-major order and \(4=\infty\). The variety identifications for the triangular source can be verified by the two inclusion tests of Proposition~\ref{prop:lower-fibres} for each target.

\subsection{An odd-cycle obstruction with a meet factor}
Let \(C=\{a,b,c,\omega\}\) with flat addition of top \(\omega\). Its only products below \(\omega\) are \(ab=ba=c\). This is the semiring denoted by \(S_c(ab)\). Let \(C^0\) have an additional additive bottom annihilating products. Lyu, Ren and Yue proved the two nonfinite basis results used in cases~\ref{grp:20} and~\ref{grp:478} \cite[Corollaries~4.10 and~4.13]{Scab}.

\begin{theorem}\label{thm:odd-meet}
Both \(C\times D\) and \(C^0\) are nonfinitely based.
\end{theorem}
\begin{proof}
For odd \(n\ge5\) use distinct letters \(x_1,\ldots,x_n\), subscripts modulo \(n\), and put
\[
        P_n=\sum_{i=1}^n x_ix_{i+1},\qquad W_n=x_1\cdots x_n .
\]
In \(C\) the value of \(P_n\) is always \(\omega\). Otherwise all its edges would have product \(c\), forcing the alternating values \(a,b\) around an odd cycle. In \(D\) the meet of all \(n\) letters is below every edge product. In \(C^0\) an assignment using the new bottom makes \(W_n\) bottom, while an assignment avoiding it evaluates \(P_n\) in \(C\). Thus both algebras satisfy \(W_n\law P_n\). Both also have \(C\) and \(D\) in their generated variety. For the extension, \(C\) is a subalgebra and sending the bottom to zero and all original elements to one gives an onto homomorphism to \(D\).

Suppose there were a finite basis and choose \(k\) bounding its variables after splitting identities into absorptions. Choose an odd \(n>k+2\), \(n\ge5\). In an equational derivation starting at \(P_n\), each successive polynomial is identically equal to \(P_n\) in \(D\). Consequently every word in it contains both endpoints of some original edge. Indeed, assign one in \(D\) exactly on that word's content. In particular every word of length at most two is an original edge in one of its two orders.

We may use an increasing list of available words in the derivation. For completeness, an elementary use of \(q\law Q\) in expanded form adds a word from one product \(l\,\theta(q)\,r\) when all the words in \(l\,\theta(Q)\,r\) are already present. Here \(l,r\) are individual, possibly empty context words, and \(\theta\) takes variables to nonempty polynomials. This description follows by distributing the one-hole context at each multiplication and dropping the part independent of the hole. A deletion in an equational derivation needs no addition to the accumulated list. Induction on the derivation therefore puts every word ever obtained on such a list, with each justified addition still absorbed by \(P_n\). Consequences of the semiring axioms amount to equality of the expanded word sets.

Consider the first addition of a word of length at least three. Incorporate each nonempty exterior context into the source absorption as one fresh variable multiplying both sides, with its chosen word as substitution. The resulting source absorption \(v\law Q'\) is still valid in \(C\) and involves at most \(k+2\) variables. Every expanded upper word has length two and is an edge of the original cycle. Thus each summand of \(Q'\) is literal or quadratic. A variable appearing as a literal has images consisting wholly of length-two words. Each variable occurring in a quadratic summand has images consisting wholly of one-letter words. The two sets of variables are disjoint.

Make the undirected graph with one edge for each quadratic upper summand, retaining loops if any. Choose one image letter for each of its vertices. All choices yield a graph homomorphism into the cycle of length \(n\). The source graph has fewer than \(n\) vertices and is bipartite. Indeed a loop or simple odd cycle would give an odd closed walk of length less than \(n\) in the target cycle. Lifting that walk to integer steps \(+1,-1\), its displacement would be an odd multiple of \(n\) of magnitude less than \(n\), which is impossible. A graph with no odd cycle is bipartite by colouring each component according to distance parity.

Evaluate the source variables in \(C\) as follows. Give the two colours the values \(a,b\), give every literal upper variable value \(c\), and assign \(\omega\) to absent variables. All upper summands evaluate to \(c\). Since \(c\) is minimal in the flat order, \(v\) must evaluate to \(c\). A literal with that value has only length-two substituted images. A quadratic word with that value uses only the two coloured sorts, so all its images also have length two. Every longer word has value \(\omega\). None of the alternatives can give the proposed first long word. This contradicts derivability of \(W_n\law P_n\).
\end{proof}
The source for case~\ref{grp:20} can therefore be formed explicitly as the eight-element direct product \(C\times D\). In case~\ref{grp:478} the isomorphism from \(C^0\) to \(\A{12029}\) sends \((\omega,a,b,c,0)\) to \((4,0,2,3,1)\). In each case Proposition~\ref{prop:lower-fibres} gives the equality-of-varieties test using these operations.

\section{Products of periodic source algebras}\label{sec:periodic-products}

Let \(B\) be the chain \(a<e<b\) with join as addition. In its multiplication \(e\) is an identity and \(a,b\) are left zeros. Thus a word takes the first value different from \(e\), or \(e\) if there is none. Write \(F=F(\mathbb F_2)\). Write its group elements as \(1,g\) with \(g^2=1\), and form \(G\) by inserting one element \(g'\) between \(g\) and the additive top. Its join with \(1\) is the top. Under multiplication \(g'\) behaves as \(g\) in both positions, so every product of length at least two takes its value in \(F\). For these algebras and for \(M,L,R\), a superscript \(0\) means adjoining an additive bottom annihilating products. In the following ordered list \(L,R,M,D,N,T\) have the meanings of Section~\ref{sec:twojoins}.
\begin{equation}\label{eq:factor-seventeen}
 M,D,L,R,F,B,B\op,M^0,L^0,R^0,F^0,B^0,(B\op)^0,N,T,G,G^0 .
\end{equation}

\begin{theorem}\label{thm:periodic-products}
Every nonempty direct product of a subfamily of \eqref{eq:factor-seventeen} has a finite identity basis. One may take its valid word equalities and word absorptions on 24 fixed letters, of word length at most 128, with at most four upper summands, together with \(\ax\).
\end{theorem}

\begin{proof}
We give the reduction to this finite family. Write \(J\) for the chosen subfamily. A member of the list that belongs to the variety of the chosen product can always be included in \(J\). In particular \(B\) contains copies of \(L,M,D\), and the corresponding subalgebras of \(B^0\) include \(L^0,M^0\). The analogous facts hold on reversing multiplication. Both \(F\) and \(G\) contain \(M\), and \(G\) contains \(F\). Each bottom extension contains its original algebra and has \(D\) as a quotient, by identifying the original elements.

\smallskip\noindent\emph{The tests on the upper words.}
Let \(I(w)\) be the sequence of distinct letters in order of first occurrence. Let \(J_f(w)\) be their order of last occurrence, read from left to right. Call a word a prefix witness at \(x\) if it contains \(x\) and every letter preceding its first \(x\) precedes \(x\) in \(I(w)\). In \(B\), \(w\law P\) requires a contained upper word and a prefix witness in \(P\) at every target letter. The subalgebra \(D\) gives containment. To force the witness, assign \(b\) to \(x\), \(e\) to the earlier letters and \(a\) elsewhere. Conversely a target value \(b\) is witnessed by taking the first occurrence assigned \(b\) before any \(a\), with all preceding occurrences assigned \(e\). A target value \(e\) is witnessed by the contained word. The minimum \(a\) is automatic. In the opposite algebra the corresponding test uses a suffix witness at each last occurrence. In \(F\) we use the content and affine-span conditions of Theorem~\ref{thm:prime-flat}. The six two-element tests are those in Lemma~\ref{lem:boolean-abs}.

For any algebra \(S\) in the list before adjoining bottom, the test in \(S^0\) says that the contained subpolynomial \(Q=P_w\) is nonempty and \(w\law Q\) holds in \(S\). Assign bottom to all exterior letters for necessity. When the target avoids bottom, its contained witnesses suffice, proving the converse. In particular \(M^0\) requires the contents of the contained words to cover the target. For \(G\), the map \(r\) sending \(g'\) to \(g\) and fixing \(F\) is a homomorphism. One has \(r(s)\le s\). Thus for a word of length at least two the test in \(G\) is exactly the test in \(F\). A literal will be treated at the end.

\smallskip\noindent\emph{Comparison of long words.}
For two words of length at least two with the same content, the following conditions suffice for a deduction of their equality. Require the same \(I\) when \(B\) is present and otherwise the same head when \(L\) is present. Require the same \(J_f\) when \(B\op\) is present and otherwise the same tail when \(R\) is present. Require the same parity vector when \(F\) is present. Our convention of including the subalgebras above accounts for all bottom extensions and for \(G,G^0\). Sufficiency for equality of values follows from the tables. We need the stronger, equational assertion.

Fix a linear order on the letters involved. A long word \(u\) equals \(u^2 u u^2\) by \((xy)^5\approx xy\). The law \((xy)^2\approx x^2y^2\) changes each outer copy of \(u^2\) to a sequence of individual squares. In the first sequence erase later copies of a repeated square. To erase use \(x^2[p]x^2[q]x\approx x^2[p][q]x\), with the last \(x\) in the middle or last block. In the last sequence erase the earlier repetitions by the reversed law. The two outer sequences now record \(I(u)\) and \(J_f(u)\), respectively.

In the middle block two consecutive distinct letters \(x,y\) may be interchanged. Expose one occurrence of each before the pair and one of each after it. The two markers on either side can occur in either order. Write every intervening nonempty piece as a distinct variable, omitting an empty piece. The identity interchanging the pair preserves both outside orders, content and parity, and is in the finite family. The same argument with an \(x\) on either side erases an adjacent pair \(x^2\) in the middle. If \(F\) is absent, a single such middle \(x\) may be erased. Thus the middle reduces to the sorted letters of odd multiplicity, or to the empty block when parity is immaterial. If \(B\) is absent, the squares of the first sequence can be sorted, keeping its first square fixed when \(L\) is present. To interchange \(x^2,y^2\), mark later occurrences of \(x,y\) and retain a nonempty prefix as one variable whenever the head must be fixed. This preserves the required head and all last occurrences. A reversed interchange sorts the last sequence to the extent permitted by \(B\op,R\). Each exchange, deletion and factorization here uses at most twelve variables and words of length at most twenty. The two original words reduce to the same form.

\smallskip\noindent\emph{Products of original witnesses.}
It suffices first to treat a target \(w\) of length at least two. A word is available if its absorption by \(P\) has already been derived. If \(T\) is in \(J\), there is a nonlinear upper word \(z_1z_2\). The law \(u^3\law u+z_1z_2\) then makes each \(u^3\), \(u\in P\), available. Without \(T\) use \(u^3\law u\). The laws hold also in \(G,G^0\), by \(r(s)\le s\). Three already available long words can be concatenated to give another. In a formal inequality \(u v t\law u+v+t\), split each of the three blocks into two nonempty variables. This meets the \(T\) condition. In each other factor the three blocks provide contained witnesses at their own letters, and the parity is the sum of three upper parities. The semantic tests prove validity. Consequently any odd product of cubed upper words is available.

For the full set \(P\), choose an odd selection of its members representing \(m(w)\) whenever \(F\) is present. Without the parity condition choose any one member. Use the preceding operation on their cubes, then append two more copies of each cubed member of \(P\). The resulting available long word \(v\) contains every upper letter. If there is a target letter missing from this union, then \(M,B,B\op,F\) and their content-requiring extensions are absent. It may be inserted as a square inside \(v\) by \(s x^2t\law st\), with both factors nonempty. This law preserves both ends and \(st\) is contained in the left word. Thus \(v\) now contains all letters in \(P,w\).

We describe how to prescribe the initial part of an available long word \(v\) in an alphabet containing the witnesses being used. With full order required, let \(x_1,\ldots,x_d\) be the desired order of distinct target letters. Suppose that \(Kv\) is available where \(K=x_1^2\cdots x_{j-1}^2\). An original prefix witness is \(q=q_1x_jq_2\) with \(\ct(q_1)\subseteq\ct(K)\). Insert \(q^2\) between \(K\) and \(v\). Next insert \(v^2\) directly after the designated \(x_j\) in the first copy of \(q\). In general one may insert a square \(t^2\) at any cut in a long available word if \(t\) is available: use \(\alpha t^2\beta\law \alpha\beta+t\), with each empty end omitted. At an exterior cut, split the long old word into two variables. The old word and \(t\) are both contained, their prefix and suffix occurrences witness every newly required occurrence, and the left parity is that of the old word. Thus these are valid instances even in the bottom extensions. In the present application \(q\) is an original upper summand, and \(v\) itself was available before the induction. After the two insertions, first occurrences up to \(x_j\) are those in \(Kx_j\). All remaining first occurrences have the order inherited from \(v\), since \(v\) contains all letters of \(q\). Last occurrences are in the final \(v\). Parity is that of \(v\). The word comparison just proved gives \(Kx_j^2v\). If only the head is prescribed perform just the first step with a same-head witness \(q\). This works even when the full product distinguishes the later first occurrences. If there is no requirement at this end, keep \(v\). The reversed procedure at the other end uses the available word resulting from the first procedure as its fixed anchor. It retains all its first occurrences. This produces an available word \(Y\) with the prescribed beginnings and endings in the orders that matter. If \(B\) is present its \(I\)-sequence begins with \(I(w)\). If \(B\op\) is present its last-occurrence sequence ends with \(J_f(w)\). Single-end requirements fix just the indicated end.

\smallskip\noindent\emph{The contained words.}
If \(D\) is absent, no bottom extension occurs and no \(B\) or \(B\op\) occurs. Write \(a\) for a formal variable to be replaced by \(w\), and \(U=a^3\). If there are exterior letters, take \(e,f\) to be nonempty lists of all of them in their first- and last-occurrence orders in \(Y\). By long-word comparison \(Y\approx U^2e^2 f^2U\) under substitution. The formal law \(U\law U^2e^2 f^2U\) is valid for the factors in question: the outer copies give the ends, the upper word contains the target letters, and the parities agree. The available word \(Y\) therefore gives \(U\law P\), and comparison of the specialization \(U=w^3\) with \(w\) gives \(w\law P\). With no exterior letters comparison alone gives \(Y\approx w\).

Now suppose \(D\) is present. The set \(Q=P_w\) is nonempty. Starting from cubes of its members, construct a second available long word \(Z\) as follows. If \(F^0\) is present choose its initial odd selection in \(Q\) to have parity \(m(w)\). Otherwise choose one member. Include two copies of each cubed member of \(Q\). Its content \(S\) is the union of the contents in \(Q\), a subset of \(\ct(w)\). If \(M^0\) is present then \(S=\ct(w)\). Apply the preceding insertion procedure using only witnesses in \(Q\). Prescribe all first occurrences when \(B^0\) is present, just the first letter when \(L^0\) is present without \(B^0\), and impose no initial requirement otherwise. Make the dual choices for the other end. These witnesses exist by the tests in the bottom extensions. This completes the construction of \(Z\). In the cubing of a literal from \(Q\), the earlier \(T\)-witness, if required, may still be taken from \(P\).

The last inference uses only six formal letters. Take \(a,b\) for \(w,Z\). If \(L^0\) is present, prefix each of the next two formal words by \(h^2\), where \(h\) will be their common actual head. If \(R^0\) is present, append \(t^2\), with \(t\) their common actual tail. Let \(U,V\) be the words, including these optional end squares, that have the following intervening blocks:
\[
               U:\ a^2b^2\,a\,b^2a^2,\qquad V:\ I\,o\,J'.
\]
Here \(o=a\) when \(F^0\) is present and \(o=b\) otherwise. If \(M^0\) is absent set \(I=J'=b^2\). When it is present put
\[
 I=\begin{cases}a^2b^2&B^0\text{ present},\\ b^2a^2&\text{otherwise},\end{cases}
 \qquad
 J'=\begin{cases}b^2a^2&(B\op)^0\text{ present},\\a^2b^2&\text{otherwise}.\end{cases}
\]
After substitution \(U\approx w\). The outer copies of \(a\) carry both orders and its middle copy carries parity. Also \(V\approx Z\). Indeed, an \(a\) used in \(V\) requires \(M^0\), so the two actual contents agree. Leading with \(a\) in place of \(b\) then uses exactly the first-order agreement secured by \(B^0\), and similarly at the other end. The odd copy may be switched to \(a\) exactly when the two actual parities have been made equal. The optional end squares preserve the actual order at their end. All these are deductions by long-word comparison, including the case of one-letter content.

With exterior letters define \(Y'=U^2e^2 f^2U\), naming their two orders as before. With no exterior letters set \(Y'=U\). In the former case the actual word \(Y'\) has the same content and required parity as \(Y\). The initial target order followed by the exterior first occurrences gives its required \(I\)-sequence, and exterior last occurrences followed by the target order give its required \(J_f\)-sequence. Where only an endpoint matters it has that endpoint. Thus the same comparison makes \(Y'\) available. When there are no exterior letters the conclusion follows directly from \(Y\approx w\).

It remains to check the formal absorption \(U\law Y'+V\). For \(M,L,R,F,B,B\op,G\) the larger word \(Y'\) supplies all required content, endpoints, first and last witnesses, and parity. Any independent containment requirement is met by \(V\). All three words are long, so \(N,T\) impose no further restriction. For the bottom extensions \(V\) is contained in \(U\). Their formal contents agree whenever \(M^0\) is present. Their first orders agree in the \(B^0\) case and their last orders agree in the opposite case. The optional squares supply the required common endpoints for \(L^0,R^0\). Their parity vectors agree when \(F^0\) is present, including the nonlinear test of \(G^0\). These are precisely the required tests. The formal words have length at most 43 and use at most six letters. Since \(Y'\) and \(V\) specialize to available words, the inference proves \(w\law P\).

\smallskip\noindent\emph{Literal targets.}
For \(w=x\), an actual literal \(x\) in \(P\) settles the question. Otherwise \(N\) is absent and, if needed, \(T\) has a long upper witness. The nonlinear argument gives \(x^3\law P\): its conditions in the old factors are those for \(x\), and in \(G,G^0\) the cubic value is below \(x\). Without \(G,G^0\) use \(x\law x^3\).

If \(G^0\) occurs, the contained words are pure powers of \(x\). At \(x=g'\), a sum consisting solely of odd nonlinear powers evaluates to \(g\). Validity therefore forces an even power in \(Q\), from which \(x^2\law P\) follows by long-word comparison. Use \(x\law x^3+x^2\), valid in all the remaining factors. Finally suppose \(G\) occurs without \(G^0\). If zero were outside the affine span of the upper parity vectors, a linear functional would be one on all of them. Indeed fix one upper vector, quotient by the span of all differences and send its nonzero coset to one. Since \(m(x)\) belongs to their affine span the functional also assigns one to \(x\). Give each letter the corresponding group element and replace the value \(g\) at \(x\) by \(g'\). Every nonlinear upper word still gives \(g\) and every upper literal also gives \(g\), since none is \(x\). This contradicts \(x\law P\). Thus zero is an odd sum of upper vectors. The earlier odd-product construction supplies a long word \(z\) of parity zero below \(P\), and long-word comparison gives \(z\approx z^2\). Use \(x\law x^3+z^2\). In \(G\) the only delicate value is \(g'\), where \(x^3=g\) and \(z^2\) is \(1\) or the top. In every other factor now present \(x\le x^3\). This completes the proof.
\end{proof}

\begin{corollary}\label{cor:periodic-register}
All 801 algebras of case~\ref{grp:57} and all 1,435 of case~\ref{grp:95} are finitely based.
\end{corollary}
\begin{proof}
For each algebra \(A\) in these two rows of the register take the subfamily of \eqref{eq:factor-seventeen} consisting of the factors \(S\) with \(S\in\V(A)\). Apply Proposition~\ref{prop:lower-fibres} with these factors as the target algebras, one at a time. Then apply the same proposition with \(A\) as target and the retained list as sources. The second test succeeds in each case. To specify the finite calculation, encode the list by the sum of \(2^{j-1}\) over its retained positions \(j\) in \eqref{eq:factor-seventeen}. The three blocks of the table give the multiplicities for the two cases. Every membership test can use a shortest generating subtuple in lexicographic order, obtained by closing each subset under the two table operations. Build word vectors by appending each generator to each distinct vector, beginning with the projections. The comparison at each principal ideal uses ordinary coordinatewise joins. These instructions and the original table numbering specify all loops and their stopping conditions. Thus each target generates exactly the same variety as the corresponding product. Apply Theorem~\ref{thm:periodic-products}. The algebras in both cases satisfy identities (3.1)--(3.12) of Wang, Wang and Li. Their Theorem~3.1 identifies the variety defined by those identities, and their Corollary~4.8 proves that every subvariety is finitely based \cite{WWL}. The algebras in case~\ref{grp:95} also satisfy \(x^3\approx x\), so the hereditary finite basis theorem of Ren, Zhao and Wang applies to that case \cite{RenZhaoWang}.
\end{proof}
\begingroup\footnotesize
\begin{longtable}{@{}rrr@{\qquad}rrr@{\qquad}rrr@{}}
\caption{Subproducts in the two periodic cases.}\label{tab:periodic-masks}\\
\toprule Mask&57&95&Mask&57&95&Mask&57&95\\\midrule\endfirsthead
\toprule Mask&57&95&Mask&57&95&Mask&57&95\\\midrule\endhead
1 & 0 & 10 & 526 & 0 & 6 & 10663 & 4 & 0\\
2 & 0 & 3 & 527 & 0 & 10 & 13003 & 4 & 0\\
3 & 0 & 82 & 539 & 0 & 4 & 16401 & 23 & 0\\
4 & 0 & 10 & 559 & 0 & 17 & 16403 & 12 & 0\\
5 & 0 & 85 & 587 & 0 & 17 & 16405 & 8 & 0\\
6 & 0 & 12 & 651 & 0 & 53 & 16409 & 8 & 0\\
7 & 0 & 86 & 655 & 0 & 3 & 16423 & 52 & 0\\
8 & 0 & 10 & 687 & 0 & 3 & 16459 & 52 & 0\\
9 & 0 & 85 & 715 & 0 & 3 & 16515 & 51 & 0\\
10 & 0 & 12 & 782 & 0 & 3 & 16519 & 4 & 0\\
11 & 0 & 86 & 911 & 0 & 4 & 16523 & 4 & 0\\
13 & 0 & 48 & 1171 & 0 & 11 & 16551 & 4 & 0\\
14 & 0 & 1 & 1175 & 0 & 1 & 16587 & 4 & 0\\
15 & 0 & 10 & 1179 & 0 & 1 & 16646 & 31 & 0\\
17 & 0 & 14 & 1207 & 0 & 1 & 16647 & 16 & 0\\
19 & 0 & 24 & 1243 & 0 & 1 & 16654 & 4 & 0\\
21 & 0 & 18 & 1431 & 0 & 2 & 16679 & 4 & 0\\
23 & 0 & 5 & 1691 & 0 & 2 & 16719 & 4 & 0\\
25 & 0 & 18 & 2471 & 0 & 20 & 16775 & 12 & 0\\
27 & 0 & 5 & 2479 & 0 & 1 & 16906 & 31 & 0\\
39 & 0 & 87 & 2543 & 0 & 1 & 16907 & 16 & 0\\
47 & 0 & 6 & 2991 & 0 & 1 & 16910 & 4 & 0\\
55 & 0 & 3 & 4811 & 0 & 20 & 16943 & 4 & 0\\
75 & 0 & 87 & 4815 & 0 & 1 & 16971 & 4 & 0\\
79 & 0 & 6 & 4847 & 0 & 1 & 17035 & 12 & 0\\
91 & 0 & 3 & 5071 & 0 & 1 & 17555 & 4 & 0\\
131 & 0 & 85 & 8231 & 48 & 0 & 18855 & 4 & 0\\
135 & 0 & 19 & 8267 & 48 & 0 & 21195 & 4 & 0\\
139 & 0 & 19 & 8323 & 47 & 0 & 24615 & 12 & 0\\
147 & 0 & 4 & 8327 & 4 & 0 & 24651 & 12 & 0\\
167 & 0 & 19 & 8331 & 4 & 0 & 24707 & 12 & 0\\
203 & 0 & 19 & 8359 & 4 & 0 & 24838 & 12 & 0\\
262 & 0 & 23 & 8395 & 4 & 0 & 25098 & 12 & 0\\
263 & 0 & 52 & 8454 & 27 & 0 & 32785 & 4 & 0\\
270 & 0 & 6 & 8455 & 16 & 0 & 32787 & 3 & 0\\
271 & 0 & 10 & 8462 & 4 & 0 & 32789 & 2 & 0\\
279 & 0 & 4 & 8487 & 4 & 0 & 32793 & 2 & 0\\
295 & 0 & 17 & 8527 & 4 & 0 & 40977 & 16 & 0\\
335 & 0 & 17 & 8583 & 12 & 0 & 40979 & 9 & 0\\
391 & 0 & 53 & 8714 & 27 & 0 & 40981 & 6 & 0\\
399 & 0 & 3 & 8715 & 16 & 0 & 40985 & 6 & 0\\
423 & 0 & 3 & 8718 & 4 & 0 & 49169 & 4 & 0\\
463 & 0 & 3 & 8751 & 4 & 0 & 57361 & 8 & 0\\
522 & 0 & 23 & 8779 & 4 & 0 & 99475 & 1 & 0\\
523 & 0 & 52 & 8843 & 12 & 0 & 107667 & 3 & 0\\
\bottomrule\end{longtable}
\endgroup

\section{Four-element source semirings}\label{sec:four-sources}
We give the arguments for the smaller algebras which occur in the register. We use the numbering of the order-four series \cite{FourSeries,FourII,FourIII,FourIV}. In the proofs below a nonlinear word means a word of length at least two, whether or not a letter is repeated. For a polynomial \(P\) and a word \(w\), write \(P_w=\{u\in P:\ct(u)\subseteq\ct(w)\}\). If \(S^0\) is obtained by adjoining a new additive bottom which annihilates products, then
\[
 S^0\models w\law P
 \quad\Longleftrightarrow\quad P_w\ne\varnothing\ \text{and}\ S\models w\law P_w.
\]
For necessity, put every exterior letter at the new bottom and take arbitrary old values on \(\ct(w)\). For sufficiency, either the target takes the bottom value or all its contained witnesses are evaluated in \(S\). In \(S\times D\) the corresponding conditions are \(P_w\ne\varnothing\) and \(S\models w\law P\), by Lemma~\ref{lem:boolean-abs}.

\subsection{A nilpotent product on three elements}\label{sec:four-nil}
Let \(K=\{a,b,\omega\}\) have flat addition with top \(\omega\). Its only product different from \(\omega\) is \(aa=b\).
\begin{proposition}\label{prop:four-nil}
Both \(K\times D\) and \(K^0\) have a finite basis. This gives the four-element sources \(42,288,430\).
\end{proposition}
\begin{proof}
For the two algebras use \(Q=P\) and \(Q=P_w\), respectively, always imposing \(P_w\ne\varnothing\). In \(K\) a word of length at least three has value \(\omega\). If \(Q\) contains such a word, or contains both a literal \(x\) and a quadratic with an occurrence of \(x\), then \(Q\) is identically \(\omega\). In the second case a common finite value would have to be \(b\), which would require \(x=b\) and \(x=a\).

If neither event occurs and there is a quadratic, write \(E\) for the letters of the quadratics and \(L\) for the literal summands. A finite upper value is possible precisely when all letters of \(E\) are \(a\) and all of \(L\) are \(b\). Its value is \(b\). Thus a literal is absorbed exactly when it belongs to \(L\), and a quadratic exactly when both its letters belong to \(E\). All longer words fail this assignment. Unmentioned letters can be set to \(\omega\) to separate a missing literal or quadratic. When \(Q\) consists of literals only, assigning \(a\) to all of them excludes every nonlinear target, and assigning a missing letter \(\omega\) separates that literal.

Consequently a valid absorption can retain at most three right words. Choose one contained word, and either one long word, the two words giving a conflict, a literal, or a quadratic through each of the at most two target letters. In the bottom extension make these choices in \(P_w\). The test still holds after the choice.

In both algebras include commutativity and
\[
 x^4\approx x^3,\qquad x^3y\approx x^2y,\qquad
 x^2y\approx xy^2,\qquad x^2yz\approx xyz .
\]
These laws hold by the three-element multiplication and preserve content. They identify all words of length at least three with the same content. With three or more letters delete repetitions by the last law. With two letters reduce to length three and use the middle laws to choose the doubled letter. With one letter retain a cube. Keep words of lengths one and two as written up to permutation.

In a chosen absorption mark every letter in a word of length at most two. There are at most eight such letters. Partition the others according to membership in the at most four word supports and represent each nonempty class by a fresh letter. Identifying letters and applying the word laws gives a valid absorption on at most \(8+15=23\) letters with word lengths at most 25. In the reverse substitution replace each fresh letter by the product of its class. The short words recover exactly, and each long word recovers its length category and support, hence its original normal form. All valid absorptions with these bounds and at most three upper words, together with the displayed word laws, therefore give a finite basis.

Here are explicit identifications. In \(S_{4,430}\) the new bottom is \(1\) and \((a,b,\omega)=(3,2,0)\). In \(S_{4,288}\) map \(3,2,0\) to \(a,b,\omega\) and \(1\) to \(\omega\), and use as second coordinate the indicator of \(1\). In \(S_{4,42}\) send \(2,1,0,3\) respectively to \(a,b,\omega,\omega\), and use as second coordinate the indicator of the complement of \(\{3\}\) for the \(D\)-coordinate. In each case the two maps preserve both operations, are onto their individual factors, and their ordered pairs distinguish the four elements. The images are thus subdirect products generating \(\V(K,D)\).
\end{proof}

\subsection{Prefix sets and the last letter}\label{sec:four-end}
Let \(E=\{e,a,\omega\}\) have flat addition with top \(\omega\), left identity \(e\), and every product beginning in \(a\) or \(\omega\) equal to \(\omega\). Write \(I(u)\) for the set of letters strictly before the last position of \(u\), and \(t(u)\) for its last letter.
\begin{proposition}\label{prop:four-end}
The algebras \(E\times D,E^0\), their opposites, and \(E\times S_{3,24}\) and its opposite are finitely based. They account, respectively, for \(292,434,290,432,357,358\).
\end{proposition}
\begin{proof}
For a nonempty upper family \(Q\) set \(I=\bigcup_{u\in Q}I(u)\), \(T=\{t(u):u\in Q\}\). In \(E\) a word equals its last value when all preceding values are \(e\), and otherwise equals \(\omega\). The condition for \(w\law Q\) is
\[
 \ct(w)\subseteq I\cup T,\qquad
 \bigl[I\cap T\ne\varnothing\ \text{or}\ (I(w)\subseteq I\ \text{and}\ t(w)\in T)\bigr].
\] Indeed a common upper value \(e\) sets all of \(I\cup T\) to \(e\), which forces the first condition. A common upper value \(a\) occurs just when all letters of \(I\) are \(e\) and all of \(T\) are \(a\). If the sets are disjoint this assignment forces exactly the second alternative. Letters outside the upper support may always be assigned \(\omega\). These two common atom values exhaust the possible failures.

Use \(Q=P\) for \(E\times D\), \(Q=P_w\) for \(E^0\), with a contained word required in either case. The laws \(xyz\approx yxz,\ x^2y\approx xy\) preserve and normalize the prefix set and last letter. All valid consequences with at most two upper words follow from the finite ones on ten letters of length at most eleven. Keep all three last letters distinct whenever they are distinct originally. Group the remaining letters by the nonempty subset of the at most three supports containing them. There are at most seven groups, occurring only in prefixes. Expand a group symbol back into its distinct letters. The original prefix sets and final letters are recovered.

If \(s\in I\cap T\), choose one original word ending in \(s\) and another containing \(s\) before its end. From their sum derive a word with their full union as its prefix set and last letter \(s\). The test above proves the inference in both algebras, since both words are contained. Its own prefix now contains its last letter. The same two-word inference incorporates each remaining upper support, retaining this last letter. Call the resulting word \(g\). For \(E\times D\), derive \(w\law g+v\) with any original contained \(v\). For \(E^0\) all the chosen words were in \(P_w\), so \(g\) itself is contained and \(w\law g\).

If \(I\cap T=\varnothing\), take an original word ending in \(t(w)\) as the initial \(g\). Adjoin the prefix set of each remaining witness to its prefix, keeping \(t(w)\) as last letter. For \(E\times D\) the old \(g\) is contained in the new word, the two words supply the required prefix set, and the old \(g\) supplies the last letter. For \(E^0\), every word of \(P_w\) ends in \(t(w)\). Indeed its tail belongs to \(\ct(w)\), and each other target letter is in \(I(w)\subseteq I\), disjoint from \(T\). Thus both words are contained at each inference. In the product use the final \(g\) and a contained original to infer \(w\). In the extension \(g\) is contained and alone suffices. This also proves the assertion for the opposite algebras by reversal. The isomorphism \(E^0\cong S_{4,434}\) sends \((e,a,\omega,0_{\rm new})\) to \((3,2,0,1)\).

For \(E\times S_{3,24}\) one has the word reduction
\[
 xyzt\approx t+xt+yt+zt+xyt+xzt+yzt .
\]
In \(E\) every term has final \(t\), and the prefixes on the right cover \(x,y,z\), giving the value of the left side. In \(S_{3,24}\) use the truncated sums of nonnegative costs from Proposition~\ref{prop:three-direct}. The cap is two, and two units of cost can always be witnessed using at most two of the first three positions together with \(t\). This reduces both sides to words of length at most three. For such a target the test in \(E\) uses at most seven witnesses. Retain one word for each target letter, two for an intersection if present, or one for the last letter and one per prefix letter if needed. Lemma~\ref{lem:bottom} gives at most \(3^4=81\) further witnesses in \(S_{3,24}\). The union still works in both factors. Thus the reduction and all valid absorptions on \(3(1+88)\) letters of length at most three form a finite basis. Table~\ref{tab:four-inclusions} gives the identifications for the two products.
\end{proof}

\subsection{Support of nonlinear upper words}
\begin{proposition}\label{prop:four-support}
Both \(S_{4,393}\) and \(S_{4,431}\) are finitely based.
\end{proposition}
\begin{proof}
Let \(H\) have addition minimum on the labels \(0,1,2\). For products of length at least two the value is \(2\) if every letter has value \(2\), and \(0\) otherwise. The first algebra generates the variety of \(H\times D\) by Table~\ref{tab:four-inclusions}. In \(H\), a nonlinear \(w\) is absorbed exactly when every target letter occurs in a nonlinear upper word. Giving that letter value \(1\) and all others \(2\) proves necessity. When the target value is \(0\), any letter with value below \(2\) supplies a nonlinear upper word of value \(0\). For a literal \(x\) the alternatives are its occurrence as a literal or in a nonlinear upper word. An unmentioned \(x\) is separated by \(x=0\) and the other variables \(2\). The \(D\)-coordinate requires a contained original.

In \(S_{4,431}\), delete the bottom \(1\). The remaining labels are \(0,2,3\), with flat sum of top \(0\). Every nonlinear word is \(3\) when all its variables are \(3\), otherwise \(0\). In this three-element algebra every target letter must occur somewhere on the right, as is seen by setting a missing letter to \(0\) and all others to \(3\). If the upper polynomial has a nonlinear word its only common finite upper value is \(3\), so coverage suffices. Upper polynomials consisting only of literals absorb precisely their literal summands, by assigning all of them \(2\) and testing missing literals separately. Apply this description to \(P_w\).

In both algebras \(xy\approx yx,\ x^2y\approx xy\) reduce each nonlinear word to its support, with a square on a singleton support. A valid absorption on at most two upper words reduces to a fixed alphabet: retain every letter of a literal word and group all other letters by their incidence in the at most three supports. Normalize after identifying them. Under expansion, each group symbol is replaced by the product of its letters. Literal words recover literally and nonlinear words recover their supports and length category. Ten letters and length eleven suffice.
For \(H\times D\), combine the nonlinear upper words two at a time to obtain a nonlinear \(g\) on their union. Each inference satisfies both factor tests. A nonlinear target is then inferred from \(g\) and the chosen contained witness. A missing literal follows from those same two words whenever it is valid. For \(431\), if the target is nonlinear, start with a nonlinear word in \(P_w\), and successively form a nonlinear word on the union with each member of \(P_w\), including its literals. Both words at every step are contained and one is nonlinear. The result has exactly the target support. A missing literal \(x\) uses a contained nonlinear power of \(x\), with \(x\law x^2\) valid. These deductions use at most two upper words at each step.
\end{proof}

\subsection{Three graphs on the quadratic words}
\begin{proposition}\label{prop:four-graphs}
The sources \(360,363,376,368\) have finite bases.
\end{proposition}
\begin{proof}
All labels below are the table labels in Appendix~\ref{app:register}. In each of the three displayed orientations the additive order is \(2,3<0<1\), with \(2,3\) incomparable. The multiplication tables give in \(360\)
\[
 x\law xy,\quad y\law xy,\quad xyz\approx xy+xz+yz,\quad xz\law xy+yz .
\]
In \(363\) the corresponding list is \(y\law xy,\ xyz\approx xz+yz,\ xz\law xy+yz\). In each case reduce to a literal set and a directed edge set, writing \(xy\) as an edge from \(x\) to \(y\). A literal is absorbed when it is present or is an endpoint in \(360\), respectively a terminal endpoint in \(363\). For a missing, ineligible literal in \(360\), assign it \(1\) and all other letters \(2\). In \(363\) assign it \(1\) and all other letters \(3\). In each case every original summand is strictly below the target.

An edge \(xy\) is absorbed precisely when there is a path of positive length from \(x\) to \(y\). Concatenating paths uses \(xz\law xy+yz\). If \(x\ne y\) and no path exists, assign \(3\) to \(x\) and all vertices reachable from \(x\), and \(2\) elsewhere. No upper edge goes from \(3\) to \(2\). Every upper value is at most \(0\) whereas the target is \(1\). If \(x=y\) and no positive path returns, put \(0\) at \(x\), \(3\) at vertices reached by a positive path from \(x\), and \(2\) elsewhere. Edges starting at \(x\) must end in the \(3\)-set. Edges leaving that set stay in it. The remaining edges begin at \(2\), so each upper value is at most \(0\), and \(0\cdot0=1\). This proves completeness for arbitrary graphs, including isolated vertices. Source \(376\) has the opposite variety to \(363\).

In \(368\) multiplication is commutative. Quadratics with both values \(2\) give \(0\), those with both values \(3\) give \(3\), and all others give \(1\). A word of length at least three gives \(3\) if all its letters are \(3\), and \(1\) otherwise. Use the four support-preserving long-word equalities in Proposition~\ref{prop:four-nil}, commutativity, and
\[
 x\law xy,\quad x^2\law xy,\quad xz\law xy+yz,\quad
 xy+yzt\approx xyzt,\quad xyz+uvw\approx xyzuvw .
\]
The stated values verify each law. Make an undirected graph of the quadratic summands with their endpoints as vertices. Start \(U\) with the support of all words of length at least three, and enlarge it by every quadratic component it meets. When \(U\ne\varnothing\), combine all long words by the last law. If a quadratic \(xy\) meets the resulting word at \(y\), use commutativity to put \(y\) first in the word, and use the preceding equality to adjoin \(x\). This yields a long word of support \(U\). Every word supported there follows: use the long-word equalities to give that word an occurrence as a factor, then apply \(x\law xy\) if a suffix remains. On each other quadratic component a path gives every quadratic with distinct endpoints there, and \(x^2\law xy\) gives a square at each incident vertex. A literal follows from its original occurrence or from any nonlinear occurrence.

These deductions are exhaustive. An unmentioned literal is separated by giving it \(1\) and all other letters \(3\). If a letter of a nonlinear target lies outside the nonlinear upper support, assign it \(0\) and put \(3\) elsewhere. The target is \(1\) and the right side is at most \(0\). To exclude a long word containing a vertex of a component outside \(U\), assign the whole component \(2\) and all other letters \(3\). Right quadratics now give \(0\) or \(3\), long words give \(3\), and literals are at most \(0\). The same assignment excludes a quadratic joining that component to a different component or to \(U\). This exhausts all word lengths.
\end{proof}

\subsection{Two prefix actions on a diamond}\label{sec:diamond-prefix}
\begin{proposition}\label{prop:diamond-prefix}
The four sources \(427,428,417,418\) are finitely based.
\end{proposition}
\begin{proof}
In \(427,428\) the additive diamond has bottom \(1\), atoms \(2,3\), and top \(0\). In either multiplication read a word from the left, skipping \(3\)'s. An all-\(3\) word gives \(3\). A first \(0\) or \(1\) gives that value. A first \(2\) gives \(0\), except in the following positions where it gives \(2\): in \(427\) the word is the literal itself, and in \(428\) this \(2\) is in the final position. This follows by reading the rows for left factors \(3,2\), since \(0,1\) are left zeros. Call an occurrence marked in \(427\) if its word is nonlinear, and marked in \(428\) if another position follows it.

The following common laws and one additional law for each source suffice:
\begin{gather*}
 x^2y\approx xy,\qquad xyx\approx xy^2,\qquad x\law x^2,\qquad
 xyz\law xt+yz,\\
 xyzt\law y+z+xyt,\qquad xy\law y+xyz,\\
 xy^2\approx xy\quad(427),\qquad xy\law xz+y\quad(428).
\end{gather*}
Their validity follows by reading the first value different from \(3\), including the terminal exception. For example, in \(xyzt\law y+z+xyt\), a first non-\(3\) in \(x\) or \(y\) with a following position is also present in \(xyt\). If \(x=y=3\), then \(y=3\) and \(z=2\) together already give \(0\), and \(z=3\) leaves the value to \(t\). The other values of \(z\) are immediate.

Write \(i(w)\) for the list of first occurrences. In either algebra the first two equalities delete internal repetitions: \(aUaV=aU^2V=aUV\) when \(U,V\) are nonempty, and \(a^2V=aV\). A final repeated letter in \(paUa\) with \(U\ne\varepsilon\) gives \(paU^2\), after which internal repetitions can again be deleted. Consequently a nonlinear word with simple final letter reduces to \(i(w)\), and one with repeated final letter reduces to \(i(w)\) followed by one more copy of its last listed letter. This includes a one-letter support, which gives a square. More explicitly, for a repeated final position one may first delete internal repetitions while keeping the final position. If the latter is not already a repeat of the last listed letter, write the resulting word \(p a U a\) with \(U\) a nonempty list of distinct letters. The second copy of \(U\) in \(p a U^2\) consists of repeat occurrences. Delete them except its final position. For \(427\) the additional equality identifies the two nonlinear forms whenever the list has at least two letters.

There must be a contained original word. Assign every target variable \(3\) and every outside letter \(1\) to see this. At the first occurrence of a target letter \(x\), let \(S\) consist of its predecessors in the list \(i(w)\). Call a right occurrence of \(x\) accessible if its earlier letters all belong to \(S\). Giving \(0\) to \(x\), \(3\) on \(S\), and \(1\) elsewhere forces an accessible occurrence. If the target occurrence is marked, then either some accessible occurrence is marked or there is also a right word supported in \(S\). Indeed change \(x\) to \(2\) in that assignment. With both alternatives absent the right words have values in \(\{1,2\}\), whereas the target gives \(0\). In particular the head of a nonlinear target has an original witness \(hD\) with \(D\) nonempty.

We derive the target from these necessary conditions. For a nonlinear word keep an available summand \(pD\) below \(P\), where \(p\) is the part of \(i(w)\) already constructed and \(D\) is nonempty. The preceding head witness starts the induction. Write \([a]x[b]\) for an accessible witness, where the optional prefix uses only letters of \(p\). When the witness is marked it is nonlinear. Split it into two nonempty blocks and apply \(xyz\law xt+yz\) with \(x=p,t=D\) and \(yz\) equal to the witness. In the resulting word the redundant prefix \(a\) can be deleted. If \(b\) is present we have \(pxb\), the next required summand. If \(b\) is absent a marked witness is possible only in \(427\). Then \(px\approx pxx\) gives a nonempty remaining block.

In the other marked-target case let \(v\in P\) be supported in \(\ct(p)\), and let \(u=[a]x[b]\) be accessible. The word laws give \(pvD\approx pD\). From \(v,u,pvD\) derive \(pvuD\) using \(xyzt\law y+z+xyt\). Deleting redundant occurrences gives \(px[b]D\). For an unmarked final target occurrence in \(428\) one may use \(xy\law xz+y\) directly with \(x=p,z=D,y=u\). If a suffix \(b\) is present the induction has produced \(pxb\) as above. With \(b\) absent, \(pu\approx px\) is exactly the simple-final target.

At the end suppose we have \(pD\) on the full first-occurrence list, with \(D\ne\varepsilon\). Choose a contained original \(v\). The word laws give \(pD\approx pv^2D\). Using \(xy\law y+xyz\) with \(x=pv,y=v,z=D\) yields \(pv^2\), which has the saturated normal form on that list. If necessary \(x\law x^2\) at the last letter restores the simple final form. For a literal target \(x\), its contained witness is \(x\) or a pure power \(x^k,\ k\ge2\). The latter reduces to \(x^2\) and again \(x\law x^2\) applies. Each inference has nonempty substituted words. Reversal gives the results for \(417,418\).
\end{proof}

\subsection{Identification of the four-element sources}\label{sec:four-identifications}
The equalities in Table~\ref{tab:four-inclusions} are equalities of generated varieties. We write \(S_j=S_{4,j}\) and use \(K,E,H\) as defined above. For a row \(U;\,V_1,\ldots,V_m\), take a shortest increasing generating tuple of \(U\), breaking ties lexicographically. Build its word vectors in the family on the right by prepending generators until closure. The first number in the last column is the number of distinct vectors. Next take a shortest generating tuple of each \(V_j\) and form its word vectors in \(U\). Their numbers follow the slash in the order of the factors. In each direction multiply at every transition and test each principal ideal by joining its word vectors coordinatewise, as in Proposition~\ref{prop:lower-fibres}. Prepending in place of appending simply reads each representing word in the other recursion order. The consistency test uses left multiplication. These finite loops depend only on the printed operations and on the enumeration in Section~\ref{sec:assembly}.

\begingroup\small
\begin{longtable}{@{}llll@{}}
\caption{Variety identifications from the four-element tables.}\label{tab:four-inclusions}\\
\toprule Algebra & Generating family & Basis argument & Word-vector counts\\\midrule\endfirsthead
\toprule Algebra & Generating family & Basis argument & Word-vector counts\\\midrule\endhead
\(S_{385}\)&\(\A{9024}\)&case~\ref{grp:2}&\(6/6\)\\
\(S_{398}\)&\(\A{2331}\)&case~\ref{grp:2}&\(26/80\)\\
\(S_{401}\)&\(\A{2345}\)&case~\ref{grp:2}&\(46/173\)\\
\(S_{405}\)&\(\A{2494}\)&case~\ref{grp:2}&\(46/173\)\\
\(S_{413}\)&\(\A{12210}\)&case~\ref{grp:2}&\(6/10\)\\
\(\A{2315}\)&\(S_{9},M\)&\ref{case:79}&\(11/4,11\)\\
\(\A{3355}\)&\(S_{433},M\)&\ref{case:61}&\(36/15,6\)\\
\(\A{1465}\)&\(S_{441},M\)&\ref{abs-square-bases}&\(10/5,5\)\\
\(\A{4731}\)&\(S_{442},M\)&\ref{grp:260}&\(19/7,7\)\\
\(\A{4712}\)&\(S_{447},M\)&\ref{case:128}&\(26/8,8\)\\
\(\A{6160}\)&\(S_{450},T\)&\ref{case:128}&\(8/8,2\)\\
\(\A{570}\)&\(S_{453},M\)&\ref{case:40}&\(15/6,6\)\\
\(S_{379}\)&\(S_{3,50},T\)&\ref{prop:three-direct}&\(8/8,2\)\\
\(S_{380}\)&\(S_{3,51},T\)&\ref{prop:three-direct}&\(7/7,2\)\\
\(S_{383}\)&\(S_{3,53},T\)&\ref{prop:three-direct}&\(7/7,2\)\\
\(S_{42}\)&\(K,D\)&\ref{prop:four-nil}&\(8/3,8\)\\
\(S_{288}\)&\(K,D\)&\ref{prop:four-nil}&\(8/3,8\)\\
\(S_{292}\)&\(E,D\)&\ref{prop:four-end}&\(19/7,7\)\\
\(S_{357}\)&\(E,S_{3,24}\)&\ref{prop:four-end}&\(11/11,11\)\\
\(S_{393}\)&\(H,D\)&\ref{prop:four-support}&\(10/5,5\)\\
\bottomrule\end{longtable}
\endgroup
For a row with a five-element algebra on the left and \(S_j,M\) or \(S_j,T\) on the right, the converse direction of Theorem~\ref{thm:twofactor} gives finite basability of \(S_j\). The term change for \(4731\) in case~\ref{grp:260} leads to \(4712\), whose basis is given directly in Section~\ref{case:128}. The row for \(9\) uses the independent quartic construction in Section~\ref{case:79}. For the rows with \(S_{3,j},T\) on the right, apply the forward direction of the same theorem. The opposite identifications for the remaining entries are
\[
 \begin{gathered}
 (290,292),\ (358,357),\ (376,363),\ (417,427),\ (418,428),\\
 (432,434),\ (445,433),\ (446,442),\ (464,453),\ (456,443).
 \end{gathered}
\]
For each pair in the order displayed the two word-vector counts of the test against the opposite are, respectively,
\(19,11,8,18,30,19,15,7,6,8\), with the same count in both directions. The bases for \(443,456\) were also proved in Proposition~\ref{prop:seven-bases}. The bases for the two chain sources \(627,705\) are in Section~\ref{case:1}.

\section{Chain sources for the product transfers}\label{sec:chain-small}
The product identifications in case~\ref{grp:47} also use the following chain semirings. Within this section, \(K_j\) denotes the algebra called \(A_j\) by Wang, Wang, Yin and Wu. Table~\ref{tab:chain-link} gives the twenty entries of their correspondence used here \cite[Appendix~D, Table~66]{WangYinWu}. Each permutation lists the images in \(K_j\) of the four elements of \(S_{4,k}\), in source order. They carry labels \(1,2,3,4\) in the original chain tables \cite[Table~28]{FourIV} and \(0,1,2,3\) here. The chain proofs below give the unrestricted deductions used in these transfers.
\begin{table}[htbp]\centering\small\caption{The chain presentations used in case~\ref{grp:47}.}\label{tab:chain-link}\begin{tabular}{@{}rrc@{\qquad}rrc@{}}
\toprule \(k\)& \(j\)& Images & \(k\)& \(j\)& Images\\\midrule
589&1861&\texttt{2130}&757&1927&\texttt{1230}\\
681&1954&\texttt{1230}&763&1931&\texttt{1230}\\
682&1919&\texttt{1230}&765&1496&\texttt{1230}\\
686&971&\texttt{2130}&766&1512&\texttt{1230}\\
690&1013&\texttt{1230}&769&1507&\texttt{1302}\\
721&1656&\texttt{2130}&772&1395&\texttt{1230}\\
734&1365&\texttt{1230}&775&1488&\texttt{1302}\\
748&1860&\texttt{1302}&780&1017&\texttt{1230}\\
751&1851&\texttt{1302}&793&1304&\texttt{2301}\\
756&1959&\texttt{1230}&799&970&\texttt{1302}\\
\bottomrule\end{tabular}\end{table}
\begingroup
\def\eq{\approx}\def\SR{\mathsf{SR}}\def\supp{\operatorname{supp}}\def\head{\operatorname{hd}}\def\sing{\operatorname{sing}}\def\Bnd{\mathcal B}
Here \(\SR\) denotes associative multiplication, associative and commutative addition, and the two distributive laws. We add \(x+x\approx x\) whenever listing a basis. Write \(\sing(w)\) for the letters occurring exactly once, \(p|_D\) for the summands of \(p\) supported inside \(D\), and \(\Bnd_A(r,l,m)\) for all valid word absorptions of \(A\) using \(r\) fixed variables at most, words of length at most \(l\), and at most \(m\) upper words. These are finite sets specified entirely by the operations. In a deduction an available intermediate word may always be retained.
\begin{lemma}\label{c4:lem:temporary}
If \(p+r\approx p\) and \(p+r+w\approx p+r\) have been derived, then \(p+w\approx p\) is derivable. If all valid word absorptions can be derived, the identities used form a basis.
\end{lemma}
\begin{proof}
Add \(w\) to the first equality and compose with the second. For the second assertion apply this observation to every summand of the two sides of an identity, in both directions, as in Section~\ref{sec:tools}.
\end{proof}

\begin{center}\small\begin{longtable}{@{}rll@{}}
\caption{Operations in the chain-source labelling.}\label{tab:chain-generators}\\
\toprule \(j\)& Addition of \(K_j\)& Multiplication of \(K_j\)\\\midrule\endfirsthead
\toprule \(j\)& Addition of \(K_j\)& Multiplication of \(K_j\)\\\midrule\endhead
473&\texttt{0123111121223123}&\texttt{0000012302000000}\\
970&\texttt{0103111101233133}&\texttt{0000010322220000}\\
971&\texttt{0123112122223123}&\texttt{0000012322220000}\\
973&\texttt{0000011301230333}&\texttt{0000012322220000}\\
1013&\texttt{0123111121223123}&\texttt{0000012301230000}\\
1017&\texttt{0123111121223123}&\texttt{0000011302230000}\\
1304&\texttt{0023012322223323}&\texttt{0100010001230100}\\
1365&\texttt{0123111121223123}&\texttt{0000012302000200}\\
1369&\texttt{0123111121233133}&\texttt{0000012302000200}\\
1389&\texttt{0123111121223123}&\texttt{0000012002000300}\\
1395&\texttt{0123111121223123}&\texttt{0000012202000300}\\
1399&\texttt{0123111121233133}&\texttt{0000012202000300}\\
1488&\texttt{0103111101233133}&\texttt{0000010022220300}\\
1496&\texttt{0123111121223123}&\texttt{0000212222220300}\\
1507&\texttt{0103111101233133}&\texttt{0000010322220300}\\
1512&\texttt{0123111121223123}&\texttt{0000012302230300}\\
1656&\texttt{0123112122223123}&\texttt{0110011001230110}\\
1851&\texttt{0103111101233133}&\texttt{0020012000200320}\\
1860&\texttt{0103111101233133}&\texttt{0020012300200320}\\
1861&\texttt{0123112122223123}&\texttt{0020012002200320}\\
1863&\texttt{0000011301230333}&\texttt{0020012002200320}\\
1919&\texttt{0123111121223123}&\texttt{0000012001200330}\\
1927&\texttt{0123111121223123}&\texttt{0000011002200330}\\
1931&\texttt{0123111121223123}&\texttt{0000012302200330}\\
1954&\texttt{0123111121223123}&\texttt{0000012301230330}\\
1959&\texttt{0123111121223123}&\texttt{0000011302230330}\\
\bottomrule\end{longtable}\end{center}

\subsection{Prefix support closure}\label{c4:sec:family41}
\begingroup
\providecommand{\eq}{}\renewcommand{\eq}{\approx}
\providecommand{\SR}{}\renewcommand{\SR}{\mathsf{SR}}
\providecommand{\supp}{}\renewcommand{\supp}{\operatorname{supp}}
\subsubsection{Language and the four generators}

Let $A$ and $B$ both have carrier $\{0,c,b,a\}$ and addition the
maximum in the chain $0<c<b<a$. Define
\begin{equation}\label{c4:f41:eq:product}
 x\cdot y=\begin{cases}y&x=a,\\f(y)&x\ne a,\end{cases}
 \qquad
\begin{array}{c|rrrr}
 y&0&c&b&a\\\hline
 f_A(y)&0&c&c&c\\
 f_B(y)&0&0&b&b
\end{array}
\end{equation}
Each $f$ is increasing, idempotent, satisfies $f(y)\le y$, and
has image disjoint from $\{a\}$. These properties show associativity:
a product of three elements is its last factor if the first two
are $a$, and is $f$ of its last factor otherwise. Left translations
are either the identity or $f$, and right translations are increasing;
all preserve the maximum of two elements. Thus both structures are
semirings. Their element $0$ is an additive identity but is not
necessarily a left multiplicative zero, and $a$ is a left identity.

Table~\ref{c4:f41:tab:generators} gives the presentations with their local indices. Strings list table entries row by row
on $\{0,1,2,3\}$. The chain coordinates $(0,c,b,a)$ map to
$(1,0,3,2)$ for $A$ and $(0,3,1,2)$ for $B$.
\begin{table}[htbp]\centering
\caption{The four generators and complete basis sizes.}\label{c4:f41:tab:generators}
\begin{tabular}{rllr}\toprule
Index&Addition&Multiplication&Basis size\\\midrule
1304&\texttt{0023012322223323}&\texttt{0100010001230100}&10\\
1488&\texttt{0103111101233133}&\texttt{0000010022220300}&10\\
1496&\texttt{0123111121223123}&\texttt{0000212222220300}&11\\
1656&\texttt{0123112122223123}&\texttt{0110011001230110}&11\\
\bottomrule\end{tabular}
\end{table}

\subsubsection{Words and exact absorption criteria}
Let $\Sigma_A$ be $\SR$ together with the following five laws:
\begin{equation}\label{c4:f41:eq:common}
\begin{aligned}
 x^2y&\eq xy,&xyz&\eq yxz,\\
 x+x&\eq x,&x+yx&\eq x,\\
 x+zy&\eq x+zy+xy.
\end{aligned}
\end{equation}
Let $\Sigma_B$ be $\Sigma_A$ together with
\begin{equation}\label{c4:f41:eq:rotation}
 xy+zx\eq xy+zx+yx.
\end{equation}
The two word laws follow from \eqref{c4:f41:eq:product}: only whether all
prefix factors are $a$ matters. Prefix absorption follows from
$yx\le x$. For the last law of \eqref{c4:f41:eq:common}, if $x=a$ then
$x$ dominates the target $xy$; otherwise $xy=f(y)\le zy$.
Hence $\Sigma_A$ holds in both $A$ and $B$.

For \eqref{c4:f41:eq:rotation} in $B$, test the three nonzero thresholds
of the chain. If $yx\ge b$, then $x\ge b$, so $zx\ge b$.
If $yx=a$, then $y=x=a$ and $xy=a$. The only remaining
case is $x=c,y=a$, when $yx=c$ and $xy=b$. This verifies
the rotation identity. It fails in $A$: at $x=b,y=a,z=0$,
its two sides have values $c$ and $b$.

Fix a finite nonempty variable set $X$. For $P\subseteq X$ and
$y\in X$, define the standard word
\begin{equation}\label{c4:f41:eq:word}
 w(P,y)=\Bigl(\prod_{x\in P}x\Bigr)y.
\end{equation}
The prefix uses each variable once, in a fixed order. The convention
$w(\varnothing,y)=y$ denotes a variable, not an empty product or a
constant term. The last variable may also lie in $P$. The word laws
in \eqref{c4:f41:eq:common} permute adjacent prefix letters and delete
repeated prefix letters, so every nonempty word reduces to this form.

For a nonempty sum $p$ of standard words define
\begin{equation}\label{c4:f41:eq:predicates}
\begin{aligned}
 T(p)&=\{z:w(Q,z)\text{ occurs in }p\},\\
 E_P(p)&\Longleftrightarrow
       \text{some }w(Q,z)\text{ in }p\text{ has }Q\cup\{z\}\subseteq P,\\
 M_{P,y}(p)&\Longleftrightarrow
       \text{some }w(Q,y)\text{ in }p\text{ has }Q\subseteq P.
\end{aligned}
\end{equation}
Repetitions may be removed by additive idempotence.

\begin{lemma}\label{c4:f41:lem:absorption}
With $C=P\cup\{y\}$, the exact absorption criteria are
\begin{equation}\label{c4:f41:eq:Acriterion}
 A\models p+w(P,y)\eq p
 \ \Longleftrightarrow\ y\in T(p)\text{ and }
          \bigl(E_P(p)\text{ or }M_{P,y}(p)\bigr),
\end{equation}
and
\begin{equation}\label{c4:f41:eq:Bcriterion}
 B\models p+w(P,y)\eq p
 \ \Longleftrightarrow\ y\in T(p),\ E_C(p),\text{ and }
          \bigl(T(p)\cap P\ne\varnothing\text{ or }M_{P,y}(p)\bigr).
\end{equation}
\end{lemma}
\begin{proof}
In $A$, a word is nonzero exactly when its last variable is
nonzero. Setting $y=c$ and all other variables to $0$ proves
the necessity of $y\in T(p)$. This condition also ensures
domination whenever the target has value $0$ or $c$.

Suppose first that $y\notin P$. Set the variables of $P$ to
$a$, set $y=b$, and set all others to $0$. A summand has value
at least $b$ precisely when its whole content lies in $P$,
or its last variable is $y$ and its prefix lies in $P$.
Thus the second condition in \eqref{c4:f41:eq:Acriterion} is necessary.
It is sufficient whenever the target is $b$ or $a$: in the
first alternative the witness is $a$, and in the second it
equals the target. If $y\in P$, the target can reach the upper
thresholds only by being $a$. The test assigning $a$ on $P$
and $0$ outside $P$ requires exactly $E_P(p)$. Here
$M_{P,y}(p)$ already implies $E_P(p)$. This proves the first criterion.

In $B$, a word is at least $b$ precisely when its last variable
is at least $b$. Setting $y=b$ and other variables to $0$ gives
the necessary tail condition, which also handles this threshold
for all assignments. The target is $a$ precisely when all
variables in $C$ are $a$. The assignment $a$ on $C$ and $0$
elsewhere shows that $E_C(p)$ is necessary and sufficient at
this top threshold.

If $y\notin P$, the remaining case is target value $c$:
all variables in $P$ equal $a$ and $y=c$. Set all other
variables to $0$. A summand is at least $c$ exactly when its
last variable belongs to $P$, or it ends in $y$ with prefix
contained in $P$. These are the alternatives in the last
condition of \eqref{c4:f41:eq:Bcriterion}; each remains sufficient
at every assignment having target value $c$. If $y\in P$,
the target cannot be $c$. Its tail condition already ensures
$T(p)\cap P\ne\varnothing$. The criteria also cover
$P=\varnothing$, where they reduce to the presence of the
linear summand $y$.
\end{proof}

\subsubsection{Derivability and finite bases}
\begin{theorem}\label{c4:f41:thm:bases}
$\Sigma_A$ is a ten-identity basis for $A$ (catalogue 1304),
and $\Sigma_B$ is an eleven-identity basis for $B$ (catalogue 1656).
\end{theorem}
\begin{proof}
We derive every absorption in Lemma~\ref{c4:f41:lem:absorption}.
First, prefix absorption permits the insertion of $w(P,y)$
from $w(Q,y)$ whenever $Q\subseteq P$: add the missing
prefix letters and use prefix permutation and repetition removal.
This handles $M_{P,y}(p)$ in either algebra.

For $A$, suppose instead that $E_P(p)$ holds and choose a
witness word $r$ with content contained in $P$. Since
$y\in T(p)$, there is a summand ending in $y$. If it is
linear, prefix absorption suffices. Otherwise it is $s y$
with $s$ a nonempty word. Substituting $x=r,z=s$ in the
last law of \eqref{c4:f41:eq:common} inserts $r y$. Its prefix
content lies in $P$, so prefix expansion inserts the target.
This proves derivability of \eqref{c4:f41:eq:Acriterion}.

For $B$, the same argument handles $E_P(p)$ and
$M_{P,y}(p)$. The only case left has $y\notin P$,
$E_C(p)$ and some $z\in T(p)\cap P$, with $z\ne y$.
Choose a witness $r=w(Q,t)$ whose content lies in $C$.
Prefix expansion inserts $w(C,t)$.

We use the rotation law to move selected prefix letters to
the end. If a source word starts with a variable $v$ and
has a nonempty suffix $u$, and a summand ends in $v$,
then $vu$ can be accompanied by $uv$. Indeed, for a
nonlinear witness $sv$, substitute $x=v,y=u,z=s$ in
\eqref{c4:f41:eq:rotation}. For a linear witness $v$, prefix
absorption already inserts $uv$. Prefix permutation lets
us move any chosen prefix letter to the start before this step.

Apply this operation to $w(C,t)$ with selected letter $z$.
After prefix repetition removal the resulting word ends
in $z$, has total content $C$, and has $y$ exactly once
in its prefix, since $y\ne z$. Now select $y$ and rotate
again, using the original witness ending in $y$. The
resulting word ends in $y$ and its prefix support is exactly
$C\setminus\{y\}=P$. This is $w(P,y)$. Every inserted
word is absorbed in the context of $p$, so auxiliary
insertions can be removed after the target has been inserted.
This proves derivability of \eqref{c4:f41:eq:Bcriterion}.

There are only $|X|2^{|X|}$ standard words. For either
algebra, let $D(p)$ be the set of those absorbed by $p$.
Every original summand belongs to $D(p)$, and every other
member can be inserted by the preceding derivations. Thus
$p$ equals the sum of all members of $D(p)$ in a fixed
order. Its term function determines $D(p)$ by absorption.
For a valid identity, expand both sides using distributivity
on their common finite variable set. They have the same
function, hence the same $D(p)$ and the same derived normal
form. The identity follows from the indicated finite basis.
No bound on the number of variables was used.
\end{proof}

\subsubsection{Opposite multiplication}
Let $\rho$ fix variables and sums and reverse the order of the factors of each word, as in Lemma~\ref{lem:transfer}.
\begin{theorem}\label{c4:f41:thm:opposite}
The reversed bases $\rho(\Sigma_A)$ and $\rho(\Sigma_B)$
are complete for catalogue 1488 and 1496, respectively.
Their sizes are ten and eleven.
\end{theorem}
\begin{proof}
In both cases the permutation of the original catalogue labels
$h=(0,2,1,3)$ maps the opposite of the source to the target:
\begin{equation}\label{c4:f41:eq:opmap}
 h(x+_{\rm source}y)=h(x)+_{\rm target}h(y),\qquad
 h(y\cdot_{\rm source}x)=h(x)\cdot_{\rm target}h(y).
\end{equation}
These maps preserve addition and reverse multiplication.
Lemma~\ref{lem:transfer} therefore transports the complete source bases
and the asserted free-algebra descriptions to the target algebras.
\end{proof}

\endgroup

\subsection{Compression of long words}\label{c4:sec:compression}
Shortening to a fixed word length is not always available. In this
section, the number of witness words is controlled first; variable
compression then produces a finite set of generating absorptions.

\subsubsection{Multiplicity types with a fixed head}
Consider the word identities
\begin{equation}\label{c4:eq:tail-laws}
xyz\eq xzy,\qquad xyx\eq x^2y,\qquad x^3\eq x^2.
\end{equation}
They put a word in a form specified by its support $D$, head $h$,
and singleton set $F$. The head is placed first, the remaining
letters are ordered, and the multiplicity of each variable is one
on $F$ and two on $D\setminus F$. Indeed, the first law permutes
the tail, the second collects occurrences of the head, and the
third caps multiplicities at two. Write $r_G$ for the resulting
word with fixed $D,h$ and singleton set $G\subseteq F$.

\begin{lemma}[Multiplicity compression]\label{c4:lem:multiplicity}
Suppose a finite ai-semiring $A$ satisfies \eqref{c4:eq:tail-laws} and
has a right identity element. Every valid absorption involving at
most two witness words follows from
\[
\SR\cup\{x+x\eq x,\ \eqref{c4:eq:tail-laws}\}
\cup\Bnd_A(30,60,2).
\]
\end{lemma}
\begin{proof}
Normalize the at most three words involved. Keep every variable
which is the head of one of them as a separate symbol; there are
at most three. Each other variable has a multiplicity vector in
$\{0,1,2\}^3\setminus\{(0,0,0)\}$, with unused coordinates
omitted if fewer words occur. Group variables with the same vector.
There are at most 26 groups, and therefore at most 29 symbols in
all. Replacing each group by one symbol, with its common
multiplicity, gives words of length at most 58. The stated bounds
30 and 60 are thus valid.

To see that the compressed absorption holds in $A$, assign one
representative of each group its compressed value and all other
members the right identity. None of these variables is a head.
Erasing the inserted right identities and sorting tails gives the
compressed word values. Validity of the original absorption therefore
implies validity of the compressed one.

Conversely, substitute the nonempty product of all members of a
group for its compressed symbol. Tail permutation restores each
member once or twice as prescribed, and \eqref{c4:eq:tail-laws}
normalize back to the original words. The compressed identity is
an instance of the bounded finite set. The right identity was used
only in assignments and is not a constant in these substitutions.
\end{proof}

We now apply the lemma to five sources. Their tables are in
Table~\ref{tab:chain-generators}. All have a right identity labelled 1
and satisfy \eqref{c4:eq:tail-laws}. For a supported witness word $v\in p|_D$
and $x\in F$, abbreviate the conditions
\[
\begin{aligned}
O_x(v)&: x\notin\supp(v)\text{ or }x\in\sing(v),\\
I_x(v)&: x\notin\supp(v)\text{ or }(x\in\sing(v),\ \head(v)=x),\\
N_x(v)&: x\notin\supp(v)\text{ or }(x\in\sing(v),\ \head(v)\ne x).
\end{aligned}
\]
Each existence assertion in the next proposition is independent;
different conditions may use different witness words.

\begin{proposition}\label{c4:prop:tail-criteria}
For a standard word $w$ with parameters $(D,h,F)$, the exact
conditions for $p$ to absorb $w$ are as follows.
\begin{enumerate}[label=(\alph*)]
\item In $K_{1389}$ and $K_{1399}$, $p|_D$ is nonempty;
every $x\in F$ has a supported witness word satisfying $O_x$;
if $h\in F$, a supported witness word satisfies $I_h$.
\item In $K_{1395}$, $p|_D$ is nonempty; a singleton $h$
has a supported witness word satisfying $O_h$, and each
$x\in F\setminus\{h\}$ has one satisfying $N_x$.
\item In $K_{1507}$, $h$ is the head of some summand of $p$;
$p|_D$ is nonempty; and every $x\in F$ has a supported
witness word satisfying $O_x$.
\item In $K_{1959}$, a summand in $p|_D$ has head $h$;
every $x\in F$ has a supported witness word satisfying $O_x$.
\end{enumerate}
\end{proposition}
\begin{proof}
In $K_{1389}$, a word gives 1 when all letters are 1, gives 2
for a single occurrence assigned 2 and all others 1, and gives 3
for a single initial occurrence assigned 3 and all others 1;
otherwise it gives 0. Its additive order is $0<3<2<1$.
In $K_{1399}$ a sole noninitial 3 instead gives 2, and the additive
order is $0<2<3<1$. In $K_{1395}$ this same multiplication is
combined with order $0<3<2<1$.

For these three sources, assign 0 outside $D$. Assigning 1
throughout $D$ forces a supported witness word. Changing just $x$ to 2
detects the single-occurrence condition. Changing it to 3
distinguishes its initial and noninitial positions. In (a), an
initial target 3 can be absorbed only by omission of that variable
or by an initial single occurrence. In (b), a noninitial target 3
has value 2 and requires a noninitial single occurrence or omission,
whereas an initial one has value 3 and permits either position.
These tests prove necessity. The same cases exhaust the nonzero
word values, and the chosen witness word absorbs the target in each,
which proves sufficiency.

In $K_{1507}$, the additive order is $2<0<3<1$. A head assigned
2 makes the word 2. Otherwise an all-1 word is 1, exactly one
occurrence of 3 with the rest 1 gives 3, and every other word is
0. Assigning $h=1$ and all other variables 2 forces a global head
witness word. Assigning 1 on $D$ and 2 outside forces a supported witness word;
changing $x\in F$ to 3 forces $O_x$. Conversely the global head
witness word covers a target value 0, a supported witness word covers 1, and an
$O_x$ witness word covers 3. The value 2 is least.

Finally $K_{1959}$ has additive order $0<3<2<1$. A word on
$\{1,2\}$ takes its head value; a sole occurrence of 3 with
the others in $\{1,2\}$ gives 3, and all other cases give 0.
Assigning $h=1$, the other variables in $D$ the value 2 and
those outside $D$ the value 0 forces a supported head witness word.
Assigning $x=3$ and the other variables in $D$ the value 2
forces $O_x$. Those witness words cover, respectively, the values 1 or 2
and the value 3, so the conditions are sufficient.
\end{proof}

\begin{theorem}\label{c4:thm:tail-bases}
For each of $K_{1389},K_{1399},K_{1395},K_{1507},K_{1959}$, a
finite basis is the set in Lemma~\ref{c4:lem:multiplicity}.
The reversed bases apply to $K_{473},K_{1369},K_{1365},K_{1860},
K_{1954}$, respectively.
\end{theorem}
\begin{proof}
Consider a valid absorption $p+w\eq p$, where the standard word $w$
has support $D$, head $h$ and singleton set $F$. For $G\subseteq F$,
let $r_G$ be the standard word with support $D$, head $h$ and
singleton set $G$.

In cases (a) and (b) of Proposition~\ref{c4:prop:tail-criteria}, choose
any supported witness $v_0\in p|_D$. It absorbs $r_{\varnothing}$.
In case (a), for $x\in F$ choose $v_x\in p|_D$ satisfying $I_h$
when $x=h$, and satisfying $O_x$ when $x\ne h$. The first choice
is available by the additional head condition in the proposition,
and $I_h$ implies $O_h$. In case (b), choose $v_h$ satisfying $O_h$
when $h\in F$, and choose $v_x$ satisfying $N_x$ for each
$x\in F\setminus\{h\}$. The corresponding absorption criterion,
applied to the polynomial consisting of the single word $v_x$,
now gives $v_x+r_{\{x\}}\eq v_x$ in every case.

In case (c), fix a summand with head $h$. Pair it with any
supported witness to absorb $r_{\varnothing}$, and with a supported
$O_x$ witness to absorb $r_{\{x\}}$. In case (d), fix instead a
supported summand with head $h$. This word absorbs
$r_{\varnothing}$; pairing it with a supported $O_x$ witness
absorbs $r_{\{x\}}$. Thus each of these initial insertions uses
at most two witness words. Coincident witnesses require only one
word because addition is idempotent.

For $G\subseteq F$ and $x\in F\setminus G$, the polynomial
$r_G+r_{\{x\}}$ absorbs $r_{G\cup\{x\}}$. Both available words
have support $D$ and head $h$. A singleton in $G$ is witnessed by
$r_G$, and the singleton $x$ is witnessed by $r_{\{x\}}$.
In case (a), a singleton head is initial in the relevant witness,
so its $I_h$ condition holds. In case (b), every singleton other
than $h$ is noninitial in that witness, so its $N_x$ condition
holds. The supported-head and global-head conditions in cases
(c) and (d) hold as well. Proposition~\ref{c4:prop:tail-criteria}
therefore justifies every merging step.

Lemma~\ref{c4:lem:multiplicity} derives all these unary and binary
absorptions from the stated finite set. Add the remaining summands
of $p$ and retain each inserted word. Successively adjoining the
members of $F$ yields $r_F=w$; when $F$ is empty, the initial
insertion already supplies $w$. Lemma~\ref{c4:lem:temporary} removes
the intermediate words. Every valid absorption is consequently
derivable. Lemma~\ref{c4:lem:temporary} proves completeness, and
Lemma~\ref{lem:transfer} gives the five reversed bases.
\end{proof}

\subsubsection{Retaining both endpoints}
Consider instead
\begin{equation}\label{c4:eq:interior-laws}
xyz\eq xy^2z,\qquad xyzu\eq xzyu.
\end{equation}
These identities put each nonlinear word in the form consisting
of its head, an ordered set of interior variables, and its tail.
The endpoints are retained even if they agree or also occur in
the interior. The second law permutes adjacent interior factors,
and the first inserts or deletes repeated interior letters.

\begin{lemma}[Endpoint compression]\label{c4:lem:endpoint-compression}
Suppose $A$ is a finite ai-semiring satisfying
\eqref{c4:eq:interior-laws}, and has an element $e$ with $aeb=ab$
for all $a,b\in A$. Every valid absorption involving at most two
witness words follows from
\[
\SR\cup\{x+x\eq x,\ \eqref{c4:eq:interior-laws}\}
\cup\Bnd_A(16,20,2).
\]
\end{lemma}
\begin{proof}
Normalize the at most three words, keeping the at most six endpoint
variables separate. Every other variable has an interior presence
vector in $\{0,1\}^3\setminus\{(0,0,0)\}$. Group equal vectors,
giving at most seven more symbols. A compressed nonlinear word
has two endpoints and at most thirteen interior symbols, hence
length at most fifteen. Linear words remain linear.

To prove validity of the compressed identity, give one representative
of each group its compressed value and give the others $e$.
Each inserted occurrence is strictly interior and is erased by
$aeb=ab$. Conversely substitute the nonempty product of each
group's members for its symbol. Interior permutation and deletion
recover the original words. This proves the assertion with the
stated, slightly relaxed, bounds.
\end{proof}

For a word with support $D$, head $h$ and tail $t$, call a witness word
\emph{terminally suitable} if it omits $t$ or has a unique occurrence
of $t$ at its tail. Suitability is required below only if $t$ occurs
once in the target.

\begin{theorem}\label{c4:thm:endpoints}
The finite set in Lemma~\ref{c4:lem:endpoint-compression} is a basis
for each of $K_{970},K_{1013},K_{1017}$. Their opposites are
$K_{1851},K_{1927},K_{1919}$, respectively.
\end{theorem}
\begin{proof}
The sources satisfy \eqref{c4:eq:interior-laws} and $a1b=ab$.
The exact absorption criteria are:
\begin{enumerate}[label=(\roman*)]
\item for $K_{970}$, a global witness word has head $h$, $p|_D$ is
nonempty, and a terminally suitable supported witness word exists when
required;
\item for $K_{1013}$, a witness word in $p|_D$ has tail $t$, and a
terminally suitable supported witness word exists when required;
\item for $K_{1017}$, a witness word in $p|_D$ has head $h$, and a
terminally suitable supported witness word exists when required.
\end{enumerate}
In $K_{970}$ the additive order is
$2<0<3<1$. A head assigned 2 gives 2; otherwise a word is 1 if
all letters are 1, is 3 for a sole exceptional occurrence which
is a final 3, and is 0 otherwise. The tests $h=1$, others 2;
1 on $D$, 2 outside; and $t=3$, others in $D$ equal to 1,
force the three conditions in (i). They also suffice: the global
head witness word covers 0, the supported witness word covers 1, and the
suitable witness word covers 3.

The other two sources have additive order $0<3<2<1$. On
$\{1,2\}$, a word takes its tail value in $K_{1013}$ and its
head value in $K_{1017}$. A sole final 3 with all previous
letters in $\{1,2\}$ gives 3; zero or a nonfinal 3 gives 0.
Assigning 0 outside $D$, the distinguished endpoint 1 and the
other variables 2 forces the endpoint witness word. Assigning $t=3$
and the other variables in $D$ the value 2 forces suitability.
The same two witness words cover all nonzero target values.

At most two witness words suffice for each criterion: choose the global head
and the suitable supported witness word in (i), or the supported endpoint
and the suitable witness word in (ii)--(iii). When suitability is not
required, an arbitrary supported witness word suffices in (i).
Lemma~\ref{c4:lem:endpoint-compression} derives their absorption from
the finite set. Add the other summands and use
Lemma~\ref{c4:lem:temporary}.
\end{proof}

\subsection{Prefix saturation with a terminal exception}\label{c4:sec:prefixes}
The two sources in this section have the same multiplication:
1 is a left identity, 0 and 2 are left zeros, and every product
beginning with 3 and having length at least two is 0. A word
therefore takes the value of its first occurrence outside $\{1\}$,
except that a first exceptional 3 gives 3 only in the final
position and gives 0 otherwise. All-1 words give 1. The additive
orders are $0<3<1<2$ in $K_{971}$ and $2<1<3<0$ in $K_{973}$.

\subsubsection{Word reduction and prefix data}
Both sources satisfy
\begin{equation}\label{c4:eq:prefix-word}
x^2y\eq xy,\qquad xyxz\eq xyz,\qquad xyx\eq xy^2.
\end{equation}
Let $\sigma(v)$ be the sequence of first occurrences in a word
$v$. For a nonempty repetition-free sequence $s=x_1\cdots x_k$,
write $R(s)=s x_k$.

\begin{lemma}\label{c4:lem:prefix-word}
The identities \eqref{c4:eq:prefix-word} reduce $v$ to $\sigma(v)$
if its final variable occurs once, and to $R(\sigma(v))$ otherwise.
In particular, $v^2$ reduces to $R(\sigma(v))$.
\end{lemma}
\begin{proof}
The first two laws delete repeated letters at nonfinal positions,
leaving at most the final occurrence as a repetition. If that
occurrence repeats an earlier $x$ with nonempty intervening word
$u$, the last law replaces $xux$ by $xu^2$. Deleting nonfinal
repetitions now leaves the last variable of $\sigma(v)$ doubled.
An adjacent final $xx$ is already in this form. This includes
the one-variable words $x$ and $x^2$.
\end{proof}

For a polynomial $p$ and a set $Y$ of variables, define:
\begin{align*}
A_p(Y)&:\quad p|_Y\ne\varnothing,\\
H_p(Y)&=\{\text{first variables outside }Y\text{ in summands not supported in }Y\},\\
G_p(Y)&=\{\text{such first variables occurring at nonfinal positions}\},\\
T_p(Y)&=\{t\notin Y:\text{ a summand ends uniquely in }t
\text{ and all preceding variables lie in }Y\}.
\end{align*}
A linear summand $t$ contributes to $T_p(Y)$ when $t\notin Y$.
This is a separate allowed case, not an empty prefix term.

For a standard word with first-occurrence sequence
$s=x_1\cdots x_k$, put $D=\supp(s)$,
$P_j=\{x_1,\ldots,x_j\}$ and $P_0=\varnothing$.

\begin{proposition}\label{c4:prop:prefix-criteria}
In $K_{971}$, $p$ absorbs a standard word $w$ exactly when
\begin{equation}\label{c4:eq:prefix-low}
\begin{gathered}
A_p(D),\qquad x_j\in H_p(P_{j-1})\quad(1\leq j\leq k),\\
\text{if }w=s\text{ and }t=x_k,\quad
A_p(D\setminus\{t\})\text{ or }t\in T_p(D\setminus\{t\}).
\end{gathered}
\end{equation}
In $K_{973}$, the exact criterion is
\begin{equation}\label{c4:eq:prefix-upper}
A_p(D),\qquad x_j\in G_p(P_{j-1})\quad(1\leq j\leq k),
\end{equation}
except that $H_p(P_{k-1})$ suffices at the last position when
$w=s$ has a unique terminal variable.
\end{proposition}
\begin{proof}
For $K_{971}$, assigning 1 on $D$ and 0 outside forces $A_p(D)$.
Assigning 1 on $P_{j-1}$, 2 to $x_j$, and 0 elsewhere forces
the $H_p$ condition, since the target is 2. When the terminal
variable is unique, assign it 3, assign 1 on $D\setminus\{t\}$,
and 0 elsewhere. A witness word can absorb the target 3 only by
avoiding $t$ and being supported there, or by ending uniquely
in $t$ with the rest in that support. Conversely, the supported
witness word covers target 1, the relevant first-exit witness word covers 2,
and the terminal condition covers 3. Zero is least.

For $K_{973}$, assign 1 on $D$ and 2 outside to force
$A_p(D)$. To detect the next first-exit variable, assign 1 on
the preceding prefix, 0 to that variable and 2 elsewhere.
The target is 0 and requires $H_p$. If the target occurrence
is nonfinal, assign 3 instead of 0. The target is again 0;
only a witness word with its first exceptional 3 nonfinal gives 0.
This forces $G_p$. For sufficiency, target 1 is covered by
$A_p(D)$, target 3 by the terminal $H_p$ condition, and target
0 by $H_p$ or $G_p$ according to whether its first exceptional
value is 0 or a nonfinal 3. The value 2 is least.
\end{proof}

\subsubsection{A common prefix construction}
The identity
\begin{equation}\label{c4:eq:left-trim}
xy+z\eq xy+z+xz
\end{equation}
allows a nonempty terminal part of a witness word to be replaced by
another available word. Suppose a seed $r_0$ is available, and
for each desired next prefix variable $x_j$ a witness word has the form
$a x_j b$ with $\supp(a)\subseteq P_{j-1}$ and $b$ nonempty;
the prefix $a$ may be absent. We can successively insert
\begin{equation}\label{c4:eq:prefix-construction}
r_j=x_1\cdots x_j r_0.
\end{equation}
Indeed, apply \eqref{c4:eq:left-trim} to $(a x_j)b$ and
$r_{j-1}$ to insert $a x_j r_{j-1}$. For $j=1$, $a$ is absent.
For $j>1$, write $P=x_1\cdots x_{j-1}$ and apply the same law
to $P r_0$ and the newly inserted word. This inserts
$P a x_j P r_0$. The repeated letters of $a$ and the second
$P$ are nonfinal, so \eqref{c4:eq:prefix-word} reduce it to
$P x_j r_0$. Every substitution is nonempty.

\begin{theorem}\label{c4:thm:prefix-low}
A complete basis for $K_{971}$ consists of $\SR$ and the seven
identities
\begin{gather}
x+x\eq x,\qquad \eqref{c4:eq:prefix-word},\qquad x+x^2\eq x,
\label{c4:eq:prefix-low-laws}\\
xy+z\eq xy+z+xz,\qquad xy+z\eq xy+z+zx.
\label{c4:eq:two-trims}
\end{gather}
This gives twelve identities in total. Reversing them gives a
basis for $K_{1861}$.
\end{theorem}
\begin{proof}
The displayed word-value rule and additive order verify all
the identities. We derive any valid absorption using
\eqref{c4:eq:prefix-low}.
First let $w=R(s)$. The low-square law permits the square of
every original witness word to be retained. Choose a supported witness word
$v\in p|_D$ and use $r_0=v^2$ as seed. A summand witnessing
$x_j\in H_p(P_{j-1})$, after squaring, has a nonempty suffix
after the chosen first-exit occurrence. Therefore the construction
\eqref{c4:eq:prefix-construction} applies. The final word $s r_0$
has all its first occurrences in $s$, and its last variable
repeats, so it reduces to $R(s)$.

Now suppose $w=s$ ends uniquely in $t$. If $t\in
T_p(D\setminus\{t\})$, use that witness word as a seed ending uniquely
in $t$. Otherwise choose a witness word $v$ supported in
$D\setminus\{t\}$ and a witness word $a t b$ witnessing
$t\in H_p(D\setminus\{t\})$. If $b$ is absent, this witness word
is a suitable seed. If $b$ is nonempty, the second trim law
in \eqref{c4:eq:two-trims} inserts $v a t$, also a suitable seed.
Here an absent $a$ is omitted. Construct the prefix
$x_1,\ldots,x_{k-1}$ using squared witnesses as before.
No inserted prefix contains $t$. Deleting only nonfinal
repetitions therefore gives $s$, retaining its unique terminal
variable. For $k=1$, the terminal condition already requires
the summand $t$. Lemma~\ref{c4:lem:temporary} proves completeness.
\end{proof}

\begin{theorem}\label{c4:thm:prefix-upper}
A complete basis for $K_{973}$ consists of $\SR$,
$x+x\eq x$, the three identities \eqref{c4:eq:prefix-word},
$x+x^2\eq x^2$, and \eqref{c4:eq:left-trim}.
It has eleven members. Its reverse is a basis for $K_{1863}$.
\end{theorem}
\begin{proof}
All displayed laws hold by the value rule. If every first-exit
condition in \eqref{c4:eq:prefix-upper} uses $G_p$, take any original
supported witness word as seed. The nonempty suffixes needed for
\eqref{c4:eq:prefix-construction} are supplied directly by $G_p$.
The resulting word reduces to $R(s)$.

For a target $s$ ending uniquely in $t$, if
$t\in G_p(D\setminus\{t\})$, this construction gives $R(s)$.
It is equal to $s^2$, and the upper-square law inserts $s$.
Otherwise $t\in H_p(D\setminus\{t\})\setminus
G_p(D\setminus\{t\})$. Its witness is then already a word
ending uniquely in $t$, all earlier variables belonging to
$D\setminus\{t\}$. Use this as seed and construct the preceding
prefix using the available $G_p$ witnesses. Nonfinal deletion
gives $s$ without repeating $t$. This also covers a linear seed;
when $k=1$ no prefix is required. Apply the two absorption lemmas.
\end{proof}

\subsection{Reconstructing single occurrences}\label{c4:sec:singletons}
The source $A=K_{1512}$ has addition given by maximum on
$0<3<2<1$ and multiplication
\[
\begin{array}{c|rrrr}
\cdot&0&1&2&3\\\hline
0&0&0&0&0\\
1&0&1&2&3\\
2&0&2&2&3\\
3&0&3&0&0
\end{array}.
\]
Thus 1 is a multiplicative identity and 0 a multiplicative zero.
A word is 1 when all occurrences are 1, is 2 when all are in
$\{1,2\}$ with some 2, and is 3 when exactly one occurrence is
3 and every subsequent occurrence is 1. Every other assignment
gives 0. In particular a 2 following 3 makes the product zero.

Let $\Gamma_A$ consist of $\SR$ and the seven identities
\begin{align}
x+x&\eq x,&x+x^2&\eq x,\label{c4:eq:singleton-add}\\
xyx&\eq yx^2,&x^3&\eq x^2,&x^2y^2&\eq y^2x^2,\label{c4:eq:singleton-word}\\
x^2y+xy^2&\eq x^2y+xy^2+xy,\label{c4:eq:merge-two}\\
xy^2z+x^2yz^2&\eq xy^2z+x^2yz^2+xyz.\label{c4:eq:merge-three}
\end{align}

\begin{theorem}\label{c4:thm:singletons}
An exact finite basis for $K_{1512}$ is
\[
\Gamma_A\cup\Bnd_A(12,24,1).
\]
Its reverse is a basis for $K_{1931}$. In particular both algebras
are finitely based.
\end{theorem}

We prove the theorem through three lemmas. The bound specifies
a finite collection of valid unary word absorptions. It does not
claim completeness of the twelve seed identities by themselves.

\subsubsection{Word normal forms and absorption}
The first word law in \eqref{c4:eq:singleton-word} moves an earlier
occurrence of a repeated variable next to its last occurrence.
Iterating, and capping powers at two, gives the sequence of last
occurrences, each variable written once or twice according as it
originally occurred once or at least twice. Distinct squared blocks
commute. Consequently the square of a word of support $D$ reduces
to
\[
b_D=\prod_{x\in D}x^2,
\]
in any chosen order. All displayed laws hold by the word-value rule: the word laws
preserve support and the suffix of each single occurrence, while
the two additive laws provide, for each surviving single occurrence,
a summand with the required suffix.

For a standard word $w$ with support $D$ and singleton set $F$,
let $A_w(x)$ be the set of variables strictly after the unique
occurrence of $x$, for $x\in F$.

\begin{lemma}\label{c4:lem:singleton-criterion}
A polynomial $p$ absorbs $w$ in $A$ if and only if
\begin{enumerate}[label=(\roman*)]
\item $p|_D$ is nonempty;
\item for every $x\in F$, some $v\in p|_D$ either omits $x$,
or has $x\in\sing(v)$ and $A_v(x)\subseteq A_w(x)$.
\end{enumerate}
\end{lemma}
\begin{proof}
Assign 1 on $D$ and 0 outside. The target is 1, which forces
(i). For (ii), assign $x=3$, assign 1 on $A_w(x)$, assign 2
to the remaining variables of $D$, and 0 outside. The target is
3. Unsupported witness words vanish. A supported witness word omitting $x$
gives 1 or 2; one containing $x$ is nonzero only if $x$ occurs
once and all variables after it belong to $A_w(x)$. These are
exactly the ways to absorb the target 3.

Conversely, target 0 is automatically absorbed. A supported witness word
covers target 1 or 2. If the target is 3, let $x$ be its sole
variable assigned 3. A witness word from (ii) omitting $x$ gives 1 or
2; otherwise its suffix variables all receive 1 and it gives 3.
This proves sufficiency.
\end{proof}

\subsubsection{Isolating one singleton and compressing}
Fix an absorption $p+w\eq p$ and retain the order of the
last-occurrence blocks in $w$. For $G\subseteq F$, let $r_G$
have this block order, multiplicity one on $G$, and multiplicity
two on $D\setminus G$. Then $r_F=w$ and
$r_{\varnothing}\eq b_D$ by the word laws.

Choose a supported witness word $v$. Contextual use of $x+x^2\eq x$
and distributivity allows any selected occurrence in $v$ to be
squared while retaining $v$ in the polynomial. Squaring all of
them and normalizing inserts $b_{\supp(v)}$. By
Lemma~\ref{c4:lem:singleton-criterion}, the valid unary absorption
\begin{equation}\label{c4:eq:singleton-required}
b_{\supp(v)}+r_{\varnothing}\eq b_{\supp(v)}
\end{equation}
would insert the all-repeated target.

For each $x\in F$, choose a witness word $v_x$ from part (ii) of the
criterion. Square all its variables except $x$ and normalize,
giving $u_x$. It either omits $x$, in which case it is a square
word, or has $x$ as its only singleton. The set of variables
following $x$ is unchanged. Thus
\begin{equation}\label{c4:eq:singleton-unary}
u_x+r_{\{x\}}\eq u_x
\end{equation}
is valid. We next derive these identities from the fixed finite
set, rather than assume an unbounded scheme.

\begin{lemma}\label{c4:lem:unary-compression}
Every valid unary absorption between words having at most one
singleton, the same distinguished variable when present, follows
from \eqref{c4:eq:singleton-word} and $\Bnd_A(12,24,1)$.
\end{lemma}
\begin{proof}
Normalize each word and order its squared blocks separately before
and after the distinguished variable $x$. Each other variable
has one of three positions in each word: absent, before $x$, or
after $x$. If the word omits $x$, place its single square block
in the before position. There are at most $3^2-1=8$ nonempty
position classes across the two words. Replace each class by one
symbol and keep $x$ separate. The compressed words use at most nine variables. Each word
contains at most six class symbols and at most two occurrences of
$x$, so its length is at most fourteen. If neither word
has a singleton, three support classes suffice.

For any assignment to the compressed symbols, assign one member
of each class its given value and every other member 1. Since 1
is a multiplicative identity and the squared blocks commute,
the original pair then has the compressed values. The compressed
absorption is therefore valid. Conversely, replace a class symbol
by the nonempty product of its members. The square of this product
reduces to their square word by \eqref{c4:eq:singleton-word}.
Ordering the blocks on either side of $x$ recovers the original
pair. The compressed absorption belongs to $\Bnd_A(12,24,1)$,
which proves the lemma.
\end{proof}

The lemma applies to \eqref{c4:eq:singleton-required} and
\eqref{c4:eq:singleton-unary}. We can therefore insert
$r_{\varnothing}$ and every $r_{\{x\}}$ using the proposed basis.

\subsubsection{Merging single occurrences}
\begin{lemma}\label{c4:lem:merge-singletons}
If $r_G$ and $r_{\{x\}}$ have been inserted, where
$x\in F\setminus G$, then the seed identities insert
$r_{G\cup\{x\}}$.
\end{lemma}
\begin{proof}
Write $r_{G\cup\{x\}}=P x Q$, with either $P$ or $Q$ possibly
absent. The supports of $P,Q$ and $\{x\}$ are disjoint because
each variable appears in a single last-occurrence block. Modulo
the word laws the two available words are
\[
r_G=P x^2Q,\qquad r_{\{x\}}=P^2xQ^2.
\]
If both $P$ and $Q$ are nonempty, substitute them for $x,z$
and the distinguished variable for $y$ in
\eqref{c4:eq:merge-three}. It inserts $P x Q$. If just one of
$P,Q$ is absent, \eqref{c4:eq:merge-two} gives the insertion,
interchanging its summands when needed. If both are absent,
the target is the already available $x$. No empty substitution
is used.
\end{proof}

\begin{proof}[Proof of Theorem~\ref{c4:thm:singletons}]
The preceding steps first insert the all-repeated word and all
one-singleton words. Apply Lemma~\ref{c4:lem:merge-singletons}
successively along $F$ to obtain $r_F=w$. Every step retains
its old summands, and Lemma~\ref{c4:lem:temporary} eliminates the
intermediates. Hence every valid word absorption follows from
the fixed finite basis. Lemma~\ref{c4:lem:temporary} proves completeness
for an arbitrary finite variable set. Reversal proves the second
case.
\end{proof}

\endgroup\FloatBarrier

\section{Sources in the case register}\label{sec:inputs}
Products, subalgebras and term changes enter at different points in the argument. If two finite algebras generate the same variety, their identities coincide by Lemma~\ref{lem:transfer}. Proposition~\ref{prop:lower-fibres} tests the two inclusions from the operations. If the operations themselves are changed by terms, Lemma~\ref{lem:terms} requires recovery in both directions. For extensions by \(M\) and \(T\) the basis construction is Theorem~\ref{thm:twofactor}. In the nonfinite case an inclusion into a target variety suffices for Theorem~\ref{thm:strong}. The other obstruction theorems have additional hypotheses on the whole family of critical identities.

The 2,236 entries of cases~\ref{grp:57} and~\ref{grp:95} come from the subproducts in Table~\ref{tab:periodic-masks}. Their finite bases and the variable bound follow from the word and polynomial arguments of Theorem~\ref{thm:periodic-products}.

For the three-element sources of case~\ref{grp:130} see Proposition~\ref{prop:three-direct}. The other bases required before adjoining a common absorbing element in case~\ref{grp:113} are given in Proposition~\ref{prop:seven-bases}. Theorem~\ref{thm:unit-paths} handles the matrix units, Theorem~\ref{thm:odd-meet} the two odd-cycle obstructions, and Theorem~\ref{thm:prime-flat} the flat groups used here.

For the four-element sources in cases~\ref{grp:4}, \ref{grp:8}, \ref{grp:13}, \ref{grp:16}, \ref{grp:143} and~\ref{grp:489}, see Section~\ref{sec:four-sources}. For the twenty additional chain presentations occurring in case~\ref{grp:47}, see Section~\ref{sec:chain-small}. Each smaller table used in that product case is covered by these sections, by Proposition~\ref{prop:seven-bases}, or by Section~\ref{case:1}.

\section{Proof of the classification theorem}\label{sec:assembly}

Use the labels \(E=\{0,1,2,3,4\}\), with \(4\) the additive top, and compare tables lexicographically in row-major order. Represent each additive semilattice by the least table under permutations fixing \(4\). In increasing order the resulting tables are \(L_0,\ldots,L_{14}\) of Table~\ref{tab:tally}. For instance, these may be obtained by enumerating partial orders in a linear extension, retaining those with binary joins, and minimizing under the indicated permutations. The greatest element always exists in a finite join-semilattice.

For each \(L_j\), consider all associative products \(m\) whose left and right translations preserve joins. These are precisely the products distributive over \(L_j\). Retain the least table in each orbit of \(\operatorname{Aut}(L_j)\). Two such algebras are isomorphic exactly when their additive indices agree and their products are in the same orbit. Ordering the representatives first by \(m\) and then by \(j\), from index \(1\), defines \(A_i\). Opposite semirings have distinct indices when they are nonisomorphic.

Here is a finite recurrence for the counts and for the order of the representatives. A linear extension of the additive order puts the four elements other than the top in increasing order. Run over the \(2^6\) subsets of their six possible strict comparisons, adjoin the top, and retain the transitive orders in which each pair has a least upper bound. Record the lexicographically least join table under the \(24\) relabellings fixing the top. For each of the fifteen tables \(L\), form the list \(E_L\) of all maps \(f\colon L\to L\) satisfying \(f(x+y)=f(x)+f(y)\), by testing the \(5^5\) maps. Starting with the empty array, append any row from \(E_L\). Retain an array precisely while every column is a prefix of some map in \(E_L\) and every associativity equation whose entries have all been assigned is valid. At depth five take the least array in each orbit under simultaneous relabelling by \(\mathrm{Aut}(L)\). Every legal table survives at every prefix, and a full surviving array satisfies associativity and both distributive laws. This proves exhaustion and the orbit rule proves uniqueness. In order \(j=0,\ldots,14\) the intermediate counts are
\CaseDisplay{(|E_{L_j}|)=\bigl(90,76,84,92,126,60,73,86,67,74,66,65,83,106,213\bigr).}
\CaseDisplay{(|\mathrm{Aut}(L_j)|)=\bigl(6,1,2,2,1,1,2,2,1,1,6,2,2,2,24\bigr).}
\CaseDisplay{(N_j)=\bigl(1622,1526,1477,1813,3852,725,1292,1559,1069,1867,744,965,1817,2227,3630\bigr).}
Here \(N_j\) counts the labelled products before taking orbits.

\begin{proposition}\label{prop:prior}
Among the types outside Table~\ref{tab:twelve}, \(14{,}005\) are finitely based and \(1{,}734\) are nonfinitely based.
\end{proposition}
\begin{proof}
Use the disjoint cases of Appendix~\ref{app:register}. Variety membership between the displayed finite algebras is tested by Proposition~\ref{prop:lower-fibres}. The unrestricted deductions in the indicated section give the finite bases, and the sequences of critical identities in the obstruction cases rule out every bound on the variables of a basis. Section~\ref{sec:inputs} identifies the basis or obstruction proof for each family of smaller sources. For example, cases~\ref{grp:2} and~\ref{grp:6} contribute \(2{,}121\) and \(1{,}967\) finite bases by the bounded-length and rectangular theorems. Case~\ref{grp:34} contributes \(1{,}660\) nonfinite bases by Theorem~\ref{thm:strong}. Count each printed index once, leaving out the twelve in Table~\ref{tab:twelve}. Grouping the resulting numbers by their additive table gives the first pair of columns in Table~\ref{tab:tally}.
\end{proof}

\begin{proof}[Proof of Theorem~\ref{thm:classification}]
Theorems~\ref{thm:heads}, \ref{thm:immediate}, \ref{thm:qualifying} and \ref{thm:pathbasis} give finite bases for the first five opposite pairs in Table~\ref{tab:twelve}. Theorem~\ref{thm:private} gives nonfinite basability for the last pair. Adding these twelve types to Proposition~\ref{prop:prior} yields \(14{,}005+10=14{,}015\) finitely based and \(1{,}734+2=1{,}736\) nonfinitely based types. Each of the six displayed products is noncommutative. Thus the commutative counts are those of the partition in Appendix~\ref{app:register}, \(2{,}968\) and \(621\), respectively. Subtraction gives \(11{,}047\) and \(1{,}115\) for the noncommutative types.
\end{proof}

A limit variety is a nonfinitely based variety all of whose proper subvarieties are finitely based. For twenty-five of the nonfinitely based types, whether the variety they generate is a limit variety remains open.

\begin{table}[htbp]\centering
\caption{The classification by additive semilattice. The first pair of counts excludes Table~\ref{tab:twelve}. FB and NFB denote finitely and nonfinitely based types.}\label{tab:tally}
\small\setlength{\tabcolsep}{4.2pt}
\begin{tabular}{@{}rlrrrrr@{}}\toprule
 &&\multicolumn{2}{c}{Other types}&\multicolumn{2}{c}{All types}&\\
\cmidrule(lr){3-4}\cmidrule(lr){5-6}
$j$& Rows of $L_j$ & FB& NFB& FB& NFB& Total\\\midrule
0 & \texttt{00004 01004 00204 00034 44444} & 294 & 54 & 294 & 54 & 348\\
1 & \texttt{00004 01004 00224 00234 44444} & 1344 & 180 & 1346 & 180 & 1526\\
2 & \texttt{00004 01014 00224 01234 44444} & 774 & 21 & 774 & 21 & 795\\
3 & \texttt{00004 01114 01214 01134 44444} & 913 & 74 & 913 & 74 & 987\\
4 & \texttt{00004 01114 01224 01234 44444} & 3724 & 124 & 3726 & 126 & 3852\\
5 & \texttt{00044 01034 00244 43434 44444} & 652 & 73 & 652 & 73 & 725\\
6 & \texttt{00044 01044 00244 44434 44444} & 565 & 135 & 565 & 135 & 700\\
7 & \texttt{00044 01134 01234 43334 44444} & 821 & 24 & 821 & 24 & 845\\
8 & \texttt{00044 01144 01234 44334 44444} & 1025 & 44 & 1025 & 44 & 1069\\
9 & \texttt{00044 01144 01244 44434 44444} & 1556 & 311 & 1556 & 311 & 1867\\
10 & \texttt{00444 01234 42244 43434 44444} & 164 & 13 & 164 & 13 & 177\\
11 & \texttt{00444 01244 42244 44434 44444} & 441 & 79 & 445 & 79 & 524\\
12 & \texttt{00444 01444 44224 44234 44444} & 701 & 219 & 703 & 219 & 922\\
13 & \texttt{00444 01444 44244 44434 44444} & 875 & 324 & 875 & 324 & 1199\\
14 & \texttt{04444 41444 44244 44434 44444} & 156 & 59 & 156 & 59 & 215\\
\midrule
&Total &14,005&1,734&14,015&1,736&15,751\\\bottomrule
\end{tabular}
\end{table}
\FloatBarrier
\section*{Use of artificial intelligence}
OpenAI Codex assisted in the search for identities and proof strategies, in writing and running programs to check finite operation tables and derivations.
The author reviewed and edited the manuscript and takes full responsibility for its content.

\appendix
\section{Index of the nonfinitely based types}\label{sec:catalogue}
The numbers refer to the convention of Section~\ref{sec:assembly}. Both ends of each interval are included. Every other index between \(1\) and \(15{,}751\) denotes a finitely based algebra.
\begingroup\small
\par\smallskip\noindent\textbf{1--1000}\quad 190--191, 207, 211--212, 229, 233--237, 241--242, 253, 257--261, 265--266, 441, 458, 460--461, 652, 681--686, 696--697, 700, 703--704, 707, 726, 730--734, 738--742, 745--746, 749--751, 753, 755--761, 763, 765--766, 768, 770, 773, 775--777, 779, 784--787, 798, 801--804, 808, 810--812, 816--821, 827--832, 834, 836, 839, 841, 843--845, 847--848, 857--860, 876, 878, 891, 894, 901, 903, 905, 907, 929, 944--945, 949--956, 958--959, 963, 965--966, 970, 972--974, 997, 999--1000.
\par\smallskip\noindent\textbf{1001--2000}\quad 1001, 1007, 1029, 1031--1033, 1035, 1037--1039, 1041, 1043--1045, 1047--1050, 1052--1054, 1056--1064, 1066--1068, 1070, 1107, 1153, 1163, 1250, 1297--1299, 1303--1305, 1401, 1403--1405, 1518--1521, 1561--1562, 1578--1583, 1608--1610, 1623--1628, 1644--1646, 1656--1658, 1662--1665, 1675--1677, 1681--1683, 1715--1717, 1779, 1797, 1831, 1865, 1875, 1948, 1950--1952.
\par\smallskip\noindent\textbf{2001--3000}\quad 2104--2105, 2192--2194, 2198--2200, 2325, 2329, 2339, 2343, 2462--2465, 2488, 2492, 2496--2497, 2502--2503, 2505, 2507--2508, 2510, 2514, 2516, 2520--2521, 2528, 2531, 2538, 2540, 2544, 2575, 2577, 2591, 2593, 2595, 2597, 2599, 2601, 2603, 2605, 2607--2611, 2613, 2666--2667, 2670--2671, 2697, 2699, 2760--2763, 2777--2778, 2786--2789, 2799--2800, 2808--2809, 2812--2813, 2825--2826, 2883.
\par\smallskip\noindent\textbf{3001--4000}\quad 3005--3015, 3157--3158, 3186--3187, 3226--3229, 3259--3262, 3281--3282, 3314--3317, 3322--3323, 3325--3349, 3351, 3382--3383, 3395--3396, 3416--3417, 3437--3439, 3452--3453, 3468--3469, 3487, 3500--3501, 3509, 3522--3523, 3535, 3547--3548, 3561--3562, 3572, 3624--3625, 3627, 3629, 3661--3662, 3664, 3666, 3710--3711, 3713, 3715, 3745--3746, 3750, 3752--3753, 3755--3756, 3759--3761, 3765, 3767--3768, 3772, 3774, 3776--3777, 3781, 3783--3784, 3796, 3799, 3807, 3810, 3839, 3842, 3848, 3851, 3853, 3858, 3861, 3879, 3881, 3885, 3926, 3929, 3933--3936, 3939, 3941--3942, 3960--3962, 3997, 4000.
\par\smallskip\noindent\textbf{4001--5000}\quad 4034, 4037, 4116, 4120--4121, 4278--4280, 4368--4369, 4417--4421, 4430, 4435--4437, 4446--4447, 4451--4455, 4469--4470, 4478, 4628--4629, 4643--4644, 4695--4696, 4714--4719, 4733, 4741, 4763--4764, 4776, 4785--4786, 4788--4789, 4793, 4838--4840, 4858--4865, 4875--4882, 4889--4890, 4895--4896, 4918, 4925--4930, 4935--4936, 4943.
\par\smallskip\noindent\textbf{5001--6000}\quad 5023--5026, 5039, 5042, 5045, 5053--5054, 5057--5061, 5063--5064, 5072--5074, 5076--5077, 5079--5080, 5085--5092, 5094--5095, 5160, 5167--5168, 5185--5194, 5234--5237, 5352, 5370, 5392, 5405--5407, 5417, 5423, 5430, 5439--5440, 5446, 5450--5452, 5454, 5456--5461, 5467--5471, 5478, 5482--5486, 5496--5499, 5505--5506, 5514--5515, 5523, 5527, 5536, 5547, 5589, 5615--5619, 5627--5634, 5637, 5640--5644, 5646--5648, 5651--5656, 5670--5674, 5683, 5685, 5688, 5690, 5693, 5708--5711, 5715--5719, 5721--5724, 5730, 5786, 5838, 5877, 5893, 5902, 5944--5946, 5972--5973, 5990.
\par\smallskip\noindent\textbf{6001--7000}\quad 6009--6010, 6024, 6073--6075, 6104--6105, 6162, 6172, 6176, 6203, 6208--6210, 6254, 6256, 6316, 6318--6319, 6377--6382, 6458, 6460, 6514--6519, 6590--6593, 6698, 6719--6720, 6730, 6738, 6746, 6789, 6809, 6811, 6819, 6855, 6858, 6892, 6894, 6947, 6952--6953.
\par\smallskip\noindent\textbf{7001--8000}\quad 7161--7162, 7196--7197, 7360--7361, 7376--7377, 7402, 7410--7415, 7420, 7430, 7438--7439, 7447, 7453, 7456--7459, 7468, 7475, 7514, 7526--7527, 7551, 7560--7565, 7567, 7573, 7581, 7595, 7600, 7604, 7610, 7630, 7636--7640, 7654--7662, 7667, 7690, 7761--7764, 7795--7799, 7805, 7893--7896, 7902, 7927--7940.
\par\smallskip\noindent\textbf{8001--9000}\quad 8038, 8045, 8065, 8069, 8077, 8094, 8098, 8105--8106, 8111, 8123--8124, 8128, 8166, 8190, 8198, 8206, 8227, 8239, 8252, 8259, 8285, 8292, 8299, 8312, 8319, 8334--8345, 8591, 8597, 8614, 8625, 8643, 8832, 8840, 8851, 8858, 8888, 8891, 8893, 8949--8952, 8990, 8997.
\par\smallskip\noindent\textbf{9001--10000}\quad 9132, 9144--9146, 9148, 9150--9152, 9154, 9227, 9244, 9318, 9368, 9370, 9499--9500, 9522--9523, 9531, 9556, 9559, 9701--9702, 9767, 9779--9780, 9803, 9811--9816, 9821, 9825, 9874--9877, 9951, 9964, 9969, 9976--9977, 9981.
\par\smallskip\noindent\textbf{10001--11000}\quad 10010, 10069, 10081, 10107, 10480, 10518--10531, 10838, 10846--10848, 10857--10858, 10860--10863, 10866--10868, 10871--10875, 10877--10878, 10884, 10895, 10897, 10915, 10935, 10948, 10971--10973, 10975--10977, 10979--10980, 10982--10983, 10985--10989, 10991--10993.
\par\smallskip\noindent\textbf{11001--12000}\quad 11005, 11023, 11025, 11028, 11030, 11032, 11035, 11037, 11048, 11052, 11074, 11076, 11080--11084, 11092--11099, 11101, 11111, 11113, 11116, 11119, 11121, 11124, 11128, 11132, 11155, 11159, 11181, 11185, 11202, 11212--11214, 11218, 11221--11227, 11242--11251, 11253--11263, 11267--11270, 11272, 11294--11295, 11309, 11330--11334, 11341, 11343--11344, 11346--11347, 11349--11350, 11352--11355, 11357, 11360, 11362--11363, 11369--11370, 11397, 11399--11400, 11402, 11406, 11408--11409, 11411--11412, 11430--11433, 11437, 11439--11440, 11442, 11445--11476, 11478--11482, 11485--11491, 11496, 11503--11507, 11509--11511, 11514--11515, 11517--11519, 11522--11524, 11536--11539, 11541--11544, 11546--11549, 11564, 11571--11572, 11597--11599, 11601--11609, 11614--11618, 11630--11631, 11634--11636, 11640--11641, 11646, 11653--11655, 11659--11660, 11669, 11672--11674, 11749, 11781--11782, 11808, 11824, 11875.
\par\smallskip\noindent\textbf{12001--13000}\quad 12029, 12031, 12071, 12103--12104, 12119, 12187--12190, 12214, 12240, 12244, 12305, 12311, 12330, 12336, 12355--12356, 12367--12368, 12390, 12392, 12396, 12406, 12411--12412, 12415, 12417, 12419, 12421, 12423--12426, 12432, 12438, 12440, 12444, 12465--12467, 12470--12471, 12474--12475, 12477--12478, 12480, 12482, 12484--12487, 12489--12492, 12494, 12503--12505, 12507--12510, 12512--12515, 12517--12520, 12522, 12530, 12539, 12541--12543, 12545, 12548, 12552, 12578--12579, 12593--12594, 12598--12599, 12603, 12620, 12627, 12629--12632, 12635, 12637, 12639--12640, 12643--12644, 12648, 12661, 12664, 12680, 12683--12684, 12686--12687, 12690--12691, 12694, 12696, 12699, 12701, 12706, 12709--12710, 12713--12714, 12719, 12721, 12724, 12726, 12731, 12737--12738, 12760--12761, 12764--12765, 12768, 12770, 12775, 12777, 12782, 12786--12789, 12791--12792, 12795--12797, 12801, 12803--12804, 12808, 12810--12811, 12815, 12819--12820, 12822, 12824--12825, 12828, 12830, 12834, 12836, 12840, 12845--12846, 12848, 12850--12851, 12855, 12857, 12861, 12863, 12867--12868, 12871--12872, 12911--12912, 12918, 12930--12932, 12934--12937, 12939, 12962--12963, 12965, 12967--12968, 12972, 12974, 12977, 12979, 12983, 12987--12988, 12990, 12992--12993, 12996, 12998.
\par\smallskip\noindent\textbf{13001--14000}\quad 13001, 13003, 13007, 13011--13012, 13035--13036, 13039--13042, 13046--13047, 13149--13150, 13156--13157, 13161, 13203, 13208--13209, 13213, 13236--13237, 13249--13250, 13255--13258, 13262--13263, 13395--13396, 13449, 13460--13461, 13500--13501, 13505--13506, 13509, 13513--13514, 13516, 13525--13526, 13530--13531, 13534, 13538--13539, 13543--13544, 13548--13549, 13555--13561, 13568, 13570, 13578, 13583, 13590, 13598--13599, 13602, 13604, 13606, 13608, 13610, 13612, 13614, 13620--13621, 13625--13628, 13630--13631, 13634--13637, 13639, 13642--13643, 13645--13646, 13649, 13651, 13653, 13655--13656, 13658, 13661, 13668--13681, 13686, 13688, 13690, 13692, 13694, 13696, 13708, 13720, 13722, 13725--13726, 13729--13730, 13732, 13738, 13740, 13742, 13744, 13746, 13748, 13751--13757, 13772, 13775, 13777, 13779, 13782, 13848, 13894--13895, 13904--13905, 13907--13908, 13970, 13996.
\par\smallskip\noindent\textbf{14001--15000}\quad 14019, 14027--14029, 14034, 14043, 14046, 14079, 14086, 14094--14097, 14102, 14108--14110, 14136--14139, 14165--14166, 14178, 14205, 14234, 14271, 14354, 14417--14418, 14450--14451, 14496, 14514, 14533, 14535, 14537--14538, 14547, 14556--14557, 14561, 14563, 14566, 14587, 14590, 14628, 14635, 14639--14641, 14643, 14669--14673, 14702, 14708, 14721--14723, 14750, 14782, 14830, 14839--14840, 14842, 14858--14859, 14904.
\par\smallskip\noindent\textbf{15001--15751}\quad 15030, 15037, 15051--15052, 15054, 15056--15057, 15061, 15063, 15066, 15068, 15072, 15077--15078, 15080, 15082--15083, 15086, 15088, 15091, 15093, 15095, 15099, 15102--15103, 15128, 15138, 15209, 15211, 15213, 15215, 15219, 15221, 15223, 15225, 15227, 15229, 15231, 15235, 15237--15239, 15241--15242, 15251, 15253, 15255, 15257, 15290--15291, 15325--15326, 15358--15359, 15445--15446, 15558--15559, 15631, 15634, 15636, 15638, 15641.
\endgroup

\section{Two-element factors by catalogue number}\label{app:factors}
The following lists give all the five-element members of the variety generated by the six two-element ai-semirings. For the indices on a line marked \(J\), the variety of \(A_i\) is precisely \(\V(B_J)\). Thus Theorem~\ref{thm:boolean-joins} supplies a basis for each of these \(2{,}233\) algebras, including the \(1{,}290\) assigned to this family in the disjoint partition. The letters in a line label stand for the factors defined in Section~\ref{sec:twojoins}.

Here are instructions for reconstructing the lists from the operations. For each product table in Section~\ref{sec:assembly}, apply Proposition~\ref{prop:lower-fibres} with the five entries \((0,1,2,3,4)\) as the generating tuple of the target and the six Boolean factors as sources. There are \(6\cdot2^5=192\) coordinates. For each success, perform the test six times in the opposite direction, using \((0,1)\) as the generating tuple of the two-element target and \(A_i\) as the sole source. Let \(J\) be the factors for which this succeeds. Finally repeat the first test with only the coordinates of \(J\). Its success gives \(A_i\in\V(B_J)\); the six tests in the other direction give \(B_J\in\V(A_i)\) by closure under products. The tests use the maxima of the principal ideals specified in the proposition. Individual entries can thus be checked without enumerating polynomials of unbounded length.
\begingroup\small\raggedright
\par\smallskip\noindent\textbf{L} (15).\quad 2054--2068.
\par\smallskip\noindent\textbf{R} (15).\quad 7328--7342.
\par\smallskip\noindent\textbf{M} (15).\quad 2631, 2636, 2662, 2702, 2708, 3797, 3823, 3849, 3882, 3900, 5635, 5679, 5780, 5789, 8821.
\par\smallskip\noindent\textbf{LM} (110).\quad 1113, 1119--1120, 1123, 1131, 1134, 1200--1201, 1209--1213, 1217--1218, 1244--1245, 1253, 1278, 1281, 1290--1291, 1331, 1337, 1350, 1354, 1358, 1360, 1410, 1414--1415, 1421, 1424, 1447, 1452--1454, 1838, 1846, 1849, 1966, 1970--1972, 1976, 1985, 1988--1989, 2042--2043, 2069--2073, 2097--2098, 2110--2111, 2208, 2632, 2638, 2656, 2658--2659, 2663, 2704, 2710, 2724, 2727--2728, 2864--2865, 2885--2886, 2888--2889, 2898--2904, 3808, 3819, 3829, 3833, 3859, 3893, 3908, 3913, 4018, 4021, 4035, 4040, 4046, 4050--4052, 4058--4059, 5675, 5748, 5751, 5781, 5813, 5815, 5833, 5841.
\par\smallskip\noindent\textbf{RM} (110).\quad 2630, 2635, 2639--2640, 2661, 2701, 2707, 2711--2712, 2741--2742, 2758--2759, 2779--2782, 3360--3362, 3384, 3420, 3440, 3458, 3470, 3795, 3816, 3847, 3889, 3925, 3927, 3948, 3956, 4375--4377, 4423, 4427--4428, 4745--4751, 4765, 4885--4886, 5038, 5040, 5097--5098, 5106, 5197, 5445, 5448, 5667, 5727, 5733, 5783, 6048--6049, 6052--6053, 6150, 6153, 6462--6463, 6466--6467, 6470--6471, 6702--6705, 6721, 6854, 6856, 6862, 6870, 7143--7146, 7343--7347, 7362--7365, 7432--7433, 7702, 7706--7709, 7723--7724, 7766, 7963, 8131, 8163, 8169, 8349.
\par\smallskip\noindent\textbf{LRM} (48).\quad 1140, 1247, 1427, 1478, 1511, 1514, 1577, 1615--1616, 1855, 1994, 2100, 2637, 2703, 2709, 2723, 2725--2726, 2866, 2890, 3397, 3442, 3454, 3549, 3563--3564, 3904, 3952, 4017, 4020, 4472, 4474, 4939--4940, 5102, 5239, 5810, 5835, 6050--6051, 6124, 6468--6469, 6866, 7719--7720, 8148, 8150.
\par\smallskip\noindent\textbf{D} (3).\quad 7083, 7261, 7280.
\par\smallskip\noindent\textbf{LD} (12).\quad 2030--2031, 2049--2053, 4301, 6417, 6427--6428, 7116.
\par\smallskip\noindent\textbf{RD} (12).\quad 7075--7076, 7080, 7225--7229, 7241, 7250--7251, 7258.
\par\smallskip\noindent\textbf{LRD} (1).\quad 4299.
\par\smallskip\noindent\textbf{MD} (82).\quad 1112, 1114, 1130, 1136, 1145, 1152, 1156, 1277, 1280, 1327, 1341, 1409, 1411, 1420, 1432, 1484, 1534, 1566, 1591, 1633, 1697, 2644, 2716, 2747, 2769, 2793, 3028, 3372, 3379, 3430, 3496, 3940, 3979, 4141, 4168, 4258, 4282, 4411, 4449, 4559, 4807, 4827, 4835, 5078, 5146, 5178, 5447, 5785, 5957, 5966, 5971, 5980, 6109, 6115, 6132, 6183, 6261, 6271, 6273, 6385, 6387, 6396, 6408, 6595, 6600, 6608, 6656, 6802, 6807, 6951, 6965, 7026, 7113, 7194, 7541, 7547, 7550, 7598, 7873, 7889, 8311, 8315.
\par\smallskip\noindent\textbf{LMD} (86).\quad 1118, 1133, 1137, 1146, 1192, 1199, 1207--1208, 1249, 1287, 1289, 1346, 1413, 1423, 1433, 1441, 1446, 1450--1451, 1486, 1538, 1541, 1592, 1594--1595, 1634, 1837, 1845, 1848, 1851--1852, 1857--1858, 1864, 1868, 1958, 1965, 1967, 1969, 1984, 1987, 1990, 1999--2000, 2015--2016, 2040--2041, 2103, 2771, 2867, 2892, 3388, 3448, 3540, 3553--3554, 3975, 4267, 4415, 4477, 4905, 5161, 5213, 5837, 5960, 6111, 6134, 6281, 6390, 6399, 6422, 6596, 6602, 6613--6615, 6815, 6957, 6970--6971, 6974, 7195, 7883, 8277, 8279.
\par\smallskip\noindent\textbf{RMD} (86).\quad 1139, 1150, 1426, 1477, 1507, 1560, 1576, 1613--1614, 1630, 2643, 2715, 2745, 2768, 2790--2792, 3022, 3464, 3490, 3502, 3971, 4126, 4163, 4245, 4274, 4537--4538, 4549, 4756, 4798, 4815, 4843, 5729, 5883, 5887, 5889--5890, 5892, 5896, 5963, 5968, 6054--6055, 6113, 6123, 6181, 6228, 6232, 6322--6323, 6327, 6331, 6335, 6402, 6474--6475, 6604, 6627, 6649, 6712, 6794, 6812, 6983--6984, 7009, 7095--7096, 7107, 7180, 7352--7353, 7419, 7506--7507, 7510, 7517--7518, 7543, 7553, 7577--7578, 7741, 7854, 7900, 8165.
\par\smallskip\noindent\textbf{LRMD} (10).\quad 1631, 1854, 1862, 1993, 2770, 4903, 5885, 6329, 6820, 7848.
\par\smallskip\noindent\textbf{N} (5).\quad 9046, 10323--10324, 10795--10796.
\par\smallskip\noindent\textbf{LN} (39).\quad 2028--2029, 2044--2048, 2087--2088, 2091--2093, 2163, 2172, 2174--2176, 2188, 2190, 9262, 9273, 9286--9287, 9290--9294, 9309--9310, 9313--9315, 10364, 10367, 10369--10370, 10806--10807.
\par\smallskip\noindent\textbf{RN} (39).\quad 7073--7074, 7220--7224, 7977--7978, 8027--8029, 8436, 8533, 8580--8582, 8723, 8781, 9610, 9655, 9724--9725, 9754--9758, 9912--9913, 9941--9943, 10443, 10461, 10464--10465, 10809--10810.
\par\smallskip\noindent\textbf{LRN} (4).\quad 4065, 4298, 9574, 9662.
\par\smallskip\noindent\textbf{MN} (39).\quad 1111, 1129, 1158, 1276, 1279, 1407--1408, 1419, 1437, 1754, 1759, 1794, 2647, 2720, 2836, 2855, 3998, 4011, 5791, 9163, 9168, 9191, 9196, 9236--9237, 9373, 9378, 9390--9391, 9557, 9570, 10013, 10015, 10343, 10390, 10462, 10490, 10812--10813.
\par\smallskip\noindent\textbf{LMN} (96).\quad 1117, 1132, 1160, 1198, 1205--1206, 1224, 1226--1227, 1255, 1286, 1288, 1412, 1422, 1439--1440, 1445, 1448--1449, 1455, 1457, 1741, 1744, 1761, 1763, 1765--1766, 1796, 1836, 1844, 1847, 1870, 1872, 1963--1964, 1968, 1983, 1986, 2004, 2006, 2038--2039, 2114, 2168, 2170, 2177, 2660, 2730, 2838, 2869, 2894, 2905--2906, 2909--2911, 4007, 4066, 4069, 5762, 5830, 5843, 9170, 9177, 9179--9180, 9183, 9198, 9202, 9204--9205, 9240--9241, 9278, 9280, 9295--9296, 9321, 9323, 9374, 9380, 9393--9397, 9566, 9575--9576, 9579, 9997, 9999, 10347, 10371, 10392, 10488.
\par\smallskip\noindent\textbf{RMN} (96).\quad 1138, 1425, 1476, 1559, 1575, 1611--1612, 1672, 1767, 2646, 2719, 2806--2807, 2835, 2839--2840, 2854, 3021, 3512, 3524, 3996, 4003, 4162, 4244, 4273, 4535--4536, 5286, 5374--5376, 5408, 5735, 5882, 5886, 5898, 6056--6057, 6122, 6320--6321, 6326, 6330, 6336, 6476--6477, 6626, 7093--7094, 8030--8031, 8046, 8133, 8171, 8440, 8496--8497, 8583, 8592, 8658, 8664, 8782, 9208, 9372, 9377, 9381--9382, 9388--9389, 9449, 9505--9507, 9524, 9555, 9562, 9659, 9684--9685, 9759--9760, 9768--9769, 9822, 9828--9829, 9836, 9879, 9944--9945, 9965, 9984, 10388, 10467, 10481, 10483.
\par\smallskip\noindent\textbf{LRMN} (16).\quad 1674, 1769, 1853, 1992, 2729, 2837, 3573, 5424, 5884, 6328, 8160, 8662, 9210, 9379, 9532, 9834.
\par\smallskip\noindent\textbf{DN} (25).\quad 7082, 7260, 7279, 8004--8005, 8474, 8483, 8560--8561, 8776--8777, 9627, 9666, 9674, 9739--9740, 9933--9936, 10348, 10456, 10458, 10803--10804.
\par\smallskip\noindent\textbf{LDN} (23).\quad 2089--2090, 2171, 2173, 2189, 4300, 6416, 6425--6426, 6435--6436, 7115, 8486--8487, 9281, 9288--9289, 9311--9312, 9678--9679, 9681, 10366.
\par\smallskip\noindent\textbf{RDN} (23).\quad 7079, 7240, 7248--7249, 7257, 7979--7980, 7990--7991, 8442, 8534, 8556--8557, 8724, 9661, 9726--9727, 9734--9735, 9791, 9914--9915, 10445.
\par\smallskip\noindent\textbf{MDN} (86).\quad 1144, 1159, 1161, 1431, 1438, 1482, 1565, 1590, 1632, 1687--1688, 1737, 1760, 1762, 1770--1771, 1795, 2817, 2842, 3025, 3518, 4166, 4256, 4281, 4556, 5308, 5310, 5398, 5403, 5956, 5965, 5984--5985, 5987, 6117, 6131, 6383--6384, 6395, 6407, 6413--6414, 6611, 6654, 7112, 8073, 8076, 8317, 8421, 8480--8481, 8530, 8639, 8642, 8714, 8805, 9169, 9171, 9197, 9199, 9218--9219, 9230, 9238--9239, 9242, 9386, 9459, 9461, 9515, 9520, 9631, 9671--9672, 9699, 9793, 9799, 9802, 9854, 9870, 9962--9963, 10345, 10349, 10473, 10479.
\par\smallskip\noindent\textbf{LMDN} (32).\quad 1225, 1456, 1593, 1689, 1764, 1772, 1856, 1871, 1873, 1998, 2005, 2154, 2169, 3568, 4265, 5419, 5959, 6389, 6398, 6421, 6617, 8289, 8531, 8708, 9178, 9203, 9220, 9265, 9279, 9527, 9700, 9864.
\par\smallskip\noindent\textbf{RMDN} (32).\quad 1629, 1673, 1768, 2816, 2841, 4547, 5287, 5888, 5899--5900, 5962, 6334, 6337, 6401, 6648, 7106, 8041, 8059, 8079, 8377, 8441, 8587, 8635, 8645, 9209, 9385, 9450, 9613, 9660, 9763, 9795, 9805.
\par\smallskip\noindent\textbf{T} (15).\quad 15737--15751.
\par\smallskip\noindent\textbf{LT} (56).\quad 21--25, 157, 162--163, 166--170, 586--587, 589, 609--610, 613, 616--617, 621--622, 626, 661, 671--674, 1807, 1825, 1908--1909, 1933, 1935--1936, 1939--1940, 2036--2037, 2075, 2078, 2080--2082, 2128--2130, 2217, 14334, 15388, 15402, 15406, 15682--15684.
\par\smallskip\noindent\textbf{RT} (56).\quad 2218--2222, 3116, 3121--3122, 3577--3581, 4303--4304, 4306, 4595--4596, 4599, 4602--4603, 4949--4950, 4954, 5427, 5555--5558, 6034, 6146, 6438--6439, 6681, 6683--6684, 6842--6843, 7139--7140, 7711, 7714, 7716--7718, 8145--8147, 8811, 12288, 13146, 13313, 13346, 13954--13956.
\par\smallskip\noindent\textbf{LRT} (4).\quad 4015, 4019, 5812, 6121.
\par\smallskip\noindent\textbf{MT} (56).\quad 2223--2227, 2249, 2364, 2372--2373, 2378, 2616, 2634, 2675--2676, 2683, 2706, 3594--3600, 3648, 3650--3651, 3790, 3825, 3840, 3872, 3874--3875, 3903, 5582--5583, 5592--5593, 5608, 5610--5613, 5681, 5790, 5808--5809, 8822, 8824, 8830, 14927, 15571, 15577, 15585, 15734--15736.
\par\smallskip\noindent\textbf{LMT} (107).\quad 66, 164--165, 337, 352--353, 358, 360--361, 372--373, 381, 391, 596, 599, 618, 620, 645--647, 655--656, 669, 1084, 1125, 1128, 1197, 1220--1221, 1243, 1254, 1270--1271, 1275, 1317, 1320, 1322--1323, 1340, 1357, 1361, 1367--1368, 1379, 1392, 1394, 1418, 1444, 1843, 1920, 1937, 1979, 1982, 2096, 2112--2113, 2209, 2211, 2292, 2368, 2384, 2445--2446, 2448, 2451--2452, 2655, 2691, 2722, 2856, 2863, 2871--2872, 2875, 2879, 2887, 2896--2897, 3590--3591, 3641, 3695, 3700, 3702--3703, 3835, 3868, 3916, 4016, 4027, 4030, 4032--4033, 4049, 4061, 4064, 5568, 5570--5571, 5584, 5598, 5752, 5816, 5842, 5858--5859, 15396.
\par\smallskip\noindent\textbf{RMT} (107).\quad 2244, 2362, 2369--2370, 2375, 2393, 2396, 2615, 2633, 2673--2674, 2681, 2705, 2731, 2733, 2756--2757, 3123, 3159--3160, 3192--3193, 3195, 3230, 3399, 3418, 3575--3576, 3633, 3635--3636, 3789, 3818, 3838, 3864, 3866--3867, 3892, 3951, 3959, 4307--4308, 4310, 4313--4314, 4370--4371, 4373--4374, 4422, 4604--4605, 4631, 4737--4738, 4946, 4948, 4962, 4964--4965, 5028--5029, 5100--5101, 5108, 5198, 5435, 5533, 5539--5540, 5553, 5586, 5669, 5734, 5746--5747, 6035, 6046--6047, 6148, 6440--6441, 6446, 6448, 6464--6465, 6685, 6699, 6840, 6847, 6849--6850, 6865, 6873, 7141--7142, 7726, 7729, 7767, 7964, 8132, 8170, 8186--8187, 8350, 8352, 13307.
\par\smallskip\noindent\textbf{LRMT} (16).\quad 414, 649, 1462, 1517, 1574, 2381, 2689, 2721, 3283, 3536, 3907, 4471, 4958, 5240, 6845, 8151.
\par\smallskip\noindent\textbf{DT} (40).\quad 16--20, 73--77, 469--472, 1799, 2927, 4099--4100, 4215--4216, 4508, 5925, 5998--5999, 6128, 6243--6244, 6274, 6341--6342, 6357, 6644, 7029, 7110, 12279, 13288, 13341, 13948--13950.
\par\smallskip\noindent\textbf{LDT} (27).\quad 152, 571--572, 583--584, 605--606, 1806, 1821, 1884--1885, 1900--1901, 1906--1907, 1922, 1931, 1991, 2018, 2021, 2034--2035, 4235, 6130, 6284, 6367, 6419.
\par\smallskip\noindent\textbf{RDT} (27).\quad 2917, 4081--4082, 4194--4195, 4484--4485, 4499, 5866, 5911--5912, 6120, 6178, 6222--6223, 6233, 6297--6298, 6307, 6620, 6638, 6986, 6989, 7012, 7091--7092, 7104.
\par\smallskip\noindent\textbf{MDT} (97).\quad 61, 71--72, 335, 339, 341--342, 369--370, 378, 404--405, 423, 440, 444--445, 467, 1083, 1116, 1127, 1135, 1157, 1172--1173, 1274, 1313, 1315--1316, 1330, 1342, 1365--1366, 1377, 1381, 1417, 1466, 1537, 1588, 1698, 1700, 2404, 2737, 2767, 2960, 3140--3141, 3152, 3212, 3222, 3409, 3982, 4096, 4145, 4224, 4252, 4356--4357, 4362, 4407, 4663, 4684, 4692, 4990, 4996, 5016, 5018--5019, 5152, 5181, 5433, 5588, 5954, 5983, 6095, 6107, 6180, 6241, 6264, 6353, 6393, 6570--6571, 6578, 6598, 6782, 6787, 6943, 6945--6946, 6968, 7192, 7877, 7892, 8316, 8332--8333, 13282.
\par\smallskip\noindent\textbf{LMDT} (32).\quad 351, 388, 406, 426, 569, 593, 651, 1194, 1196, 1349, 1390, 1443, 1468, 1544, 1589, 1842, 1850, 1869, 1898, 1918, 1961, 1981, 3270, 3528, 4365, 5005, 5217, 6584, 6612, 6936, 6977, 8280.
\par\smallskip\noindent\textbf{RMDT} (32).\quad 411, 438, 1460, 1510, 1573, 2401, 2735, 2766, 2949, 3974, 4079, 4128, 4200, 4243, 4617, 4656, 4673, 4699, 5535, 5881, 5897, 6220, 6231, 6305, 6325, 6691, 6775, 6791, 7178, 7745, 7857, 7901.
\par\smallskip\noindent\textbf{NT} (40).\quad 1--15, 9036--9045, 9047--9049, 10319--10322, 10325--10326, 10794, 10797--10801.
\par\smallskip\noindent\textbf{LNT} (59).\quad 147, 563--564, 579--580, 601--602, 633, 639, 642, 1800, 1804--1805, 1819, 1895--1897, 1903--1905, 1929, 2011, 2014, 2032--2033, 2153, 2161--2162, 2166--2167, 2180, 9055--9059, 9113, 9118, 9120, 9122--9123, 9127, 9129--9130, 9264, 9272, 9275, 9299--9300, 9335--9336, 10330--10332, 10363, 10368, 10372, 10375, 10805.
\par\smallskip\noindent\textbf{RNT} (59).\quad 2912, 4073--4074, 4190--4191, 4480--4481, 5249, 5331, 5334, 5860, 5864--5865, 6118, 6217--6219, 6294--6296, 6618, 6979, 6982, 7089--7090, 8376, 8434--8435, 8491--8492, 8661, 9337--9341, 9427, 9470, 9473, 9475--9476, 9537, 9539--9540, 9612, 9654, 9682, 9832--9833, 9995--9996, 10381--10383, 10442, 10463, 10484, 10487, 10808.
\par\smallskip\noindent\textbf{MNT} (59).\quad 56, 326--327, 333, 365--366, 375, 447, 1077, 1081--1082, 1110, 1126, 1272, 1310--1312, 1326, 1362--1364, 1375, 1416, 1736, 1752--1753, 2413, 2833--2834, 4014, 5594, 9079, 9165, 9190, 9342--9346, 9355, 9376, 9471, 9546--9547, 9549, 9551--9552, 9573, 10014, 10030--10031, 10342, 10384--10386, 10389, 10491, 10494, 10811.
\par\smallskip\noindent\textbf{LMNT} (36).\quad 350, 385, 453, 590, 657, 1191, 1195, 1345, 1388, 1442, 1747, 1757--1758, 1841, 1916, 1957, 1980, 2455, 2907--2908, 4072, 5579, 9085, 9124--9125, 9186, 9193, 9322, 9362, 9392, 9542--9543, 9582, 10000, 10346, 10391.
\par\smallskip\noindent\textbf{RMNT} (36).\quad 408, 1458, 1506, 1572, 2410, 2831--2832, 2948, 4006, 4125, 4197, 4242, 5335, 5371, 5542, 5880, 6227, 6303, 6324, 8493--8494, 8667, 9352, 9375, 9478, 9501--9502, 9534--9535, 9565, 9683, 9839, 9880, 9985, 10387, 10466.
\par\smallskip\noindent\textbf{DNT} (52).\quad 2922, 4089--4090, 4211--4212, 4505, 5276, 5279, 5919, 5923--5924, 6126, 6237--6239, 6338--6340, 6642, 7023, 7109, 8419, 8424, 8472--8473, 8475--8476, 9050--9054, 9087, 9091, 9093--9094, 9252--9253, 9439, 9442, 9629, 9634, 9665, 9667, 10327--10329, 10350, 10353, 10455, 10457, 10802.
\par\smallskip\noindent\textbf{LDNT} (12).\quad 636, 2157, 2164--2165, 4232, 6279, 6365, 6418, 9116, 9268, 9274, 10365.
\par\smallskip\noindent\textbf{RDNT} (12).\quad 4496, 5252, 6636, 7007, 7103, 8380, 8437--8438, 9430, 9616, 9656, 10444.
\par\smallskip\noindent\textbf{MDNT} (36).\quad 420, 450, 456, 1464, 1533, 1587, 1740, 1755--1756, 2956, 4137, 4221, 4250, 5360, 5367, 5953, 6259, 6351, 6392, 8527--8528, 8717, 9082, 9088--9089, 9174, 9192, 9231, 9488--9489, 9497, 9697, 9858, 9873, 10344, 10478.
\endgroup

\section{Case register}\label{app:register}
An index \(i\) in a row refers to the pair of operations \(\A{i}\) defined in Section~\ref{sec:assembly}. Number ranges include both endpoints. Where a row contains several source algebras, the indices following each source belong to that source. A superscript \(\mathrm{op}\) reverses words. Four polynomial entries in a term-change row give the two forward terms followed by the two recovery terms, in the order addition, multiplication. A reference to a source case is followed using its algebra number. The degree in \(S_{d,j}\) distinguishes the numbering of the smaller sources from the numbering at order five. For an obstruction the hypotheses are those of the indicated proof, including its all-length premise wherever required.
\newcounter{casegroup}
\begingroup\small
\refstepcounter{casegroup}\label{grp:1}\par\addvspace{.55ex}\noindent\textbf{Case 1.} FB.\enspace \ref{thm:qgraph}.
1, 2, 4, 6, 7, 10, 12--15, 26, 27, 29, 31--35, 46, 47, 49, 78--82, 87, 90, 140, 171--173, 176, 185--189, 269, 270, 273--278, 313, 317, 887, 890, 1701, 1702, 1705, 1712, 1718, 1721, 1782, 9036, 9037, 9039, 9041, 9042, 9045, 9047--9049, 9065, 9067, 10319, 10320, 10322, 10325, 10326, 10378, 10538--10540, 10545--10549, 10571, 10757, 10794, 10797, 10798, 10800, 10801, 10814, 10815, 10818, 10824, 10825, 12224--12226, 12231--12234, 12266, 12267, 12297--12300, 13377, 13380, 14312--14314, 14319--14321, 14323, 14337--14342, 14346, 14350--14353, 14355--14360, 14362--14365, 14367, 14372, 14376, 14378, 14379, 14381--14386, 14391, 14392, 14394--14397, 14399--14402, 14406--14409, 14437--14440, 14443--14449, 14456, 14458--14463, 14466--14473, 14478--14481, 14483--14485, 14847, 14848, 14850--14852, 14854--14857, 14860, 14864, 14866, 14870--14882, 14886--14891, 14895--14902, 14913, 14930--14933, 14935--14941, 14974, 14975, 14981, 14982, 14998, 14999, 15003--15005, 15018, 15019, 15135, 15136, 15145, 15146, 15149, 15178, 15192, 15407, 15410, 15413, 15417, 15421--15423, 15429, 15464, 15465, 15468, 15472, 15475, 15478, 15481, 15487--15492, 15496, 15497, 15500, 15501, 15510, 15518--15521, 15525, 15526, 15529, 15530, 15586, 15587, 15589, 15590, 15608, 15610, 15624, 15689, 15705, 15707, 15710, 15716, 15719--15722, 15725--15728, 15730. 
\par

\refstepcounter{casegroup}\label{grp:2}\par\addvspace{.55ex}\noindent\textbf{Case 2.} FB.\enspace \ref{thm:bounded}.
3, 5, 8, 9, 11, 23, 25, 28, 30, 48, 58, 60, 68, 70, 88, 89, 101, 102, 117, 125, 127, 143, 149, 151, 159, 161, 168, 169, 174, 175, 195, 196, 205, 217, 218, 225, 226, 231, 232, 245, 255, 256, 271, 272, 296, 310, 312, 314--316, 328, 329, 356, 362, 363, 366--368, 373, 374, 376, 377, 382, 383, 386, 387, 392, 393, 395, 396, 403, 409, 410, 448, 449, 454, 455, 481, 482, 507, 513, 519, 534, 535, 557, 565, 566, 580--582, 587, 588, 591, 592, 597, 598, 603, 604, 611, 612, 614, 623, 624, 628, 629, 634, 635, 640, 641, 671, 688, 692, 695, 702, 711, 713, 721, 723, 728, 729, 744, 752, 754, 790, 791, 796, 797, 799, 800, 842, 862, 863, 868, 869, 896, 898, 904, 909, 911, 915, 916, 926, 934, 935, 941, 942, 946, 947, 971, 977, 981, 983, 987, 989, 990, 995, 996, 998, 1003, 1023, 1027, 1030, 1036, 1065, 1079, 1232, 1239, 1268, 1273, 1279, 1281, 1288, 1290, 1292, 1293, 1309, 1310, 1318, 1321, 1322, 1324, 1325, 1332, 1334, 1335, 1338, 1339, 1343, 1344, 1351, 1352, 1355, 1356, 1363, 1364, 1368, 1376, 1380, 1383, 1389, 1393, 1394, 1396, 1400, 1402, 1407, 1415, 1416, 1418, 1419, 1421, 1422, 1424, 1425, 1437, 1439, 1440, 1442, 1444, 1445, 1447, 1449, 1453--1455, 1457, 1459, 1504, 1505, 1559, 1572, 1575, 1612, 1672, 1734, 1735, 1742, 1745, 1746, 1753, 1754, 1758, 1759, 1761, 1763, 1765--1767, 1794, 1796, 1802, 1820, 1894, 1895, 1904, 1905, 1909, 1917, 1921, 1930, 1934, 1935, 1938, 1939, 1943, 1947, 1955, 1956, 1963, 1971, 1973, 1974, 1977, 1978, 1980, 1982, 1983, 1985, 1986, 1988, 1989, 2004, 2006, 2007, 2009, 2012, 2013, 2023, 2029, 2033, 2037, 2039, 2043, 2046, 2048, 2056, 2058, 2061, 2062, 2064, 2076, 2077, 2079, 2084, 2088, 2092, 2126, 2133, 2134, 2144, 2146, 2151, 2152, 2162, 2163, 2167, 2168, 2170, 2172, 2174, 2175, 2178, 2179, 2181, 2182, 2187, 2188, 2190, 2220, 2222, 2225, 2227, 2235, 2237, 2243, 2246, 2248, 2251, 2253, 2263, 2266, 2267, 2270, 2271, 2275, 2277, 2280, 2282, 2285, 2287, 2290, 2291, 2294, 2296, 2302, 2312, 2316, 2323, 2324, 2331, 2332, 2338, 2345, 2346, 2361, 2370, 2371, 2373, 2374, 2376, 2377, 2379, 2380, 2385, 2386, 2388, 2389, 2391, 2392, 2394, 2395, 2397, 2398, 2411, 2412, 2414, 2415, 2418, 2420, 2422, 2424, 2430, 2431, 2433, 2434, 2442, 2444, 2446, 2447, 2449, 2450, 2453, 2454, 2456, 2457, 2469, 2473, 2477, 2485, 2487, 2494, 2495, 2501, 2509, 2511, 2517, 2526, 2527, 2532, 2533, 2541, 2546, 2547, 2555, 2559, 2564, 2568, 2573, 2574, 2578, 2590, 2594, 2598, 2614, 2661--2663, 2674, 2676, 2680, 2682, 2684, 2686, 2688, 2692, 2696, 2700, 2705--2708, 2710--2712, 2719, 2720, 2722, 2724, 2728, 2730, 2732, 2734, 2756--2759, 2780, 2782, 2806, 2807, 2832, 2834--2836, 2838--2840, 2854, 2855, 2872, 2876, 2880, 2882, 2887--2889, 2894, 2897, 2899, 2902, 2904, 2906, 2908--2911, 2914, 2916, 3118, 3120, 3125, 3127, 3579, 3580, 3598, 3599, 3617, 3618, 3621, 3623, 3626, 3628, 3631, 3637, 3638, 3646, 3652, 3653, 3657, 3658, 3660, 3663, 3665, 3669, 3670, 3676, 3677, 3684, 3690, 3691, 3693, 3694, 3698, 3704, 3705, 3712, 3714, 3717, 3735, 3736, 3742, 3743, 3754, 3757, 3758, 3762, 3763, 3778, 3779, 3785, 3786, 3836, 3838, 3840, 3841, 3843, 3845, 3847, 3849, 3850, 3852, 3854, 3859, 3865, 3866, 3870, 3873, 3874, 3876, 3877, 3880, 3882--3884, 3886, 3887, 3890, 3891, 3895, 3897, 3898, 3901, 3902, 3910, 3911, 3914, 3915, 3922, 3923, 3925, 3927, 3930, 3931, 3949, 3950, 3957, 3958, 3994, 3996, 3998, 3999, 4001, 4004, 4005, 4009, 4012, 4013, 4019, 4028, 4031, 4032, 4035, 4041, 4043, 4044, 4047, 4048, 4050, 4055, 4056, 4062, 4063, 4065, 4067, 4070, 4071, 4075, 4076, 4191--4193, 4198, 4199, 4205, 4206, 4210, 4273, 4298, 4304, 4305, 4308, 4309, 4311, 4312, 4315, 4316, 4318, 4319, 4322, 4323, 4427, 4428, 4482, 4483, 4488, 4597, 4598, 4600, 4606--4608, 4630, 4634, 4951, 4952, 4966, 4967, 5250, 5251, 5329, 5330, 5332, 5333, 5336, 5337, 5555, 5610, 5621, 5625, 5635, 5636, 5638, 5645, 5649, 5657, 5664, 5668, 5677, 5680, 5686, 5691, 5699, 5704, 5744, 5769, 5776, 5806, 5811, 5820, 5827, 5856, 5862, 6119, 6216, 6217, 6225, 6226, 6295, 6296, 6304, 6309, 6320, 6324, 6326, 6330, 6336, 6439, 6441, 6447, 6449, 6451, 6453, 6464--6467, 6470, 6471, 6476, 6477, 6619, 6634, 6682, 6683, 6686, 6700, 6703, 6722, 6724, 6841, 6842, 6848, 6849, 6851, 6852, 6854, 6856, 6859, 6860, 6863, 6864, 6871, 6872, 6980, 6981, 6991, 7074, 7090, 7094, 7140, 7142, 7144, 7146, 7222, 7224, 7330, 7332, 7335, 7336, 7338, 7345, 7347, 7363, 7368, 7369, 7440, 7712, 7713, 7715, 7727, 7728, 7775, 7978, 8026, 8028, 8031, 8143, 8184, 8374, 8375, 8435, 8436, 8440, 8490, 8492, 8494, 8496, 8497, 8533, 8579--8581, 8583, 8593, 8655, 8656, 8659, 8660, 8665, 8666, 8722, 8723, 8779--8782, 9007, 9023, 9024, 9033, 9035, 9038, 9040, 9043, 9044, 9046, 9057, 9059, 9064, 9066, 9078, 9080, 9084, 9086, 9100, 9112, 9114, 9119, 9121, 9128, 9129, 9138, 9142, 9149, 9163, 9164, 9181, 9184, 9185, 9190, 9191, 9193, 9196, 9198, 9202, 9204, 9205, 9208, 9236, 9237, 9240, 9241, 9262, 9263, 9272, 9273, 9275, 9278, 9280, 9287, 9291, 9292, 9297, 9298, 9301, 9302, 9309, 9310, 9313, 9314, 9333, 9339, 9341, 9344, 9346, 9351, 9353, 9354, 9356, 9358, 9360, 9361, 9363, 9367, 9371, 9375--9378, 9380--9382, 9388--9393, 9395--9397, 9426, 9428, 9472, 9474, 9477, 9479, 9538, 9539, 9550, 9551, 9553, 9555, 9557, 9558, 9560, 9563, 9564, 9568, 9571, 9572, 9574, 9575, 9577, 9580, 9581, 9610, 9611, 9654, 9655, 9659, 9662, 9682--9685, 9725, 9755, 9756, 9760, 9768, 9770, 9830, 9831, 9837, 9838, 9912, 9913, 9939--9942, 9944, 9945, 9993, 10028, 10037, 10038, 10040--10042, 10047, 10048, 10050, 10051, 10079, 10085, 10094, 10190, 10191, 10193, 10195, 10196, 10198, 10264, 10266, 10268, 10295, 10298, 10299, 10301, 10313, 10316, 10318, 10321, 10323, 10324, 10332, 10343, 10347, 10364, 10367, 10369, 10373, 10374, 10383, 10386, 10388, 10390, 10392, 10443, 10461, 10462, 10464, 10467, 10485, 10486, 10492, 10493, 10509--10512, 10544, 10557--10568, 10586, 10592, 10594, 10719--10725, 10736, 10737, 10739--10741, 10743--10751, 10753, 10755, 10756, 10758, 10779, 10780, 10784--10786, 10788--10791, 10793, 10795, 10796, 10799, 10806, 10807, 10809, 10810, 10812, 10813, 10816, 10817, 10827, 10829, 10831, 10833, 10837, 10839, 10842, 10843, 10849, 10880, 10881, 10886, 10891, 10892, 10898--10900, 10902, 10904, 10906, 10908, 10910--10912, 10914, 10916, 10919--10922, 10927, 10929, 10930, 10932, 10933, 10962, 10964--10966, 10970, 11024, 11026, 11029, 11031, 11033, 11041, 11044, 11045, 11047, 11049--11051, 11053, 11054, 11057, 11059, 11060, 11063, 11064, 11100, 11102, 11105, 11108, 11109, 11112, 11114, 11117, 11120, 11122, 11127, 11129--11131, 11133, 11134, 11137, 11139, 11140, 11143, 11144, 11146, 11147, 11149, 11152, 11153, 11156--11158, 11160, 11161, 11164, 11166, 11167, 11170, 11171, 11174, 11175, 11178, 11179, 11203, 11206, 11209, 11210, 11300--11302, 11304, 11305, 11311--11313, 11315--11319, 11321--11324, 11328, 11477, 11483, 11492, 11512, 11520, 11527, 11532, 11560, 11567, 11589, 11610, 11622, 11627, 11730, 11743, 11765, 11802, 11803, 11805, 11809--11811, 11813, 11817, 11819, 11827--11829, 11831, 11834, 11854, 11859, 11868, 11887--11889, 11892, 11893, 11895, 11896, 11900, 11901, 11904, 11906, 11909, 11910, 11923, 11928--11930, 11956, 11957, 11960, 11971, 11975, 11976, 12000, 12001, 12014, 12015, 12025, 12028, 12039--12043, 12045, 12061, 12065, 12145, 12157, 12159, 12160, 12168, 12171, 12184, 12185, 12191--12193, 12210, 12219, 12230, 12257--12259, 12261, 12262, 12264, 12265, 12287, 12307, 12308, 12312--12316, 12319, 12326, 12327, 12332, 12333, 12337--12341, 12344, 12346, 12348, 12352, 12358, 12360, 12362, 12370, 12372, 12374, 12376, 12378, 12380, 12381, 12387, 12389, 12391, 12393, 12395, 12397, 12399, 12401, 12403, 12407, 12414, 12416, 12418, 12427, 12434--12436, 12441, 12443, 12445, 12447, 12448, 12450, 12452, 12454, 12455, 12461, 12544, 12546, 12549--12551, 12553, 12554, 12558, 12559, 12562, 12563, 12570, 12571, 12574, 12575, 12600--12602, 12604, 12605, 12608, 12609, 12612, 12613, 12633, 12641, 12645--12647, 12649, 12650, 12654, 12655, 12658, 12659, 12662, 12663, 12665, 12666, 12670, 12671, 12674, 12675, 12688, 12689, 12692, 12697, 12703, 12704, 12711, 12712, 12716, 12717, 12722, 12728, 12729, 12732, 12734--12736, 12739, 12748, 12754, 12755, 12762, 12763, 12766, 12772, 12773, 12779, 12780, 12790, 12793, 12794, 12798, 12799, 12805, 12806, 12812, 12813, 12841, 12847, 12849, 12852, 12853, 12858, 12859, 12864, 12865, 12869, 12870, 12873, 12891, 12892, 12898, 12899, 12902, 12904--12906, 12914, 12915, 12919, 12921--12925, 12928, 12949, 12954, 12955, 12989, 12991, 12994, 12999, 13004, 13005, 13008--13010, 13013, 13025, 13030, 13031, 13078, 13080--13083, 13085, 13086, 13088, 13136, 13140, 13143, 13144, 13310, 13311, 13314, 13319, 13327, 13332--13334, 13342, 13343, 13345--13348, 13350, 13351, 13354--13357, 13359, 13360, 13364, 13367, 13368, 13371, 13373, 13374, 13378, 13379, 13384--13392, 13397, 13398, 13401, 13405, 13406, 13408, 13411, 13418, 13420--13423, 13427--13443, 13445, 13451--13459, 13464, 13466--13477, 13479, 13481, 13483, 13485, 13488, 13490, 13492, 13494, 13496, 13498, 13499, 13503, 13504, 13508, 13511, 13512, 13515, 13518, 13519, 13523, 13524, 13528, 13529, 13533, 13536, 13537, 13546, 13547, 13551, 13552, 13567, 13569, 13571, 13574, 13576, 13577, 13579, 13580, 13582, 13584, 13585, 13587, 13589, 13591, 13593, 13595, 13601, 13603, 13605, 13611, 13613, 13615, 13618, 13619, 13623, 13624, 13629, 13632, 13633, 13650, 13652, 13654, 13660, 13663, 13664, 13685, 13689, 13695, 13697, 13699, 13702, 13704, 13711, 13712, 13717, 13718, 13727, 13728, 13731, 13733, 13734, 13741, 13743, 13747, 13749, 13765, 13768, 13781, 13784, 13952, 13953, 13959--13961, 13964, 13985, 13987, 13989, 13991, 13995, 13997, 14000, 14001, 14003, 14007, 14014, 14020, 14030, 14040, 14041, 14044, 14045, 14047, 14048, 14052, 14055, 14056, 14065, 14066, 14076, 14080, 14081, 14098, 14115, 14122, 14123, 14156, 14159, 14194, 14199, 14200, 14202, 14206, 14211, 14212, 14215, 14216, 14218, 14222, 14226, 14229, 14230, 14233, 14240, 14243, 14244, 14246--14249, 14254--14257, 14259, 14260, 14265, 14266, 14269, 14279, 14281, 14282, 14297, 14302--14304, 14307, 14308, 14318, 14322, 14329, 14333, 14335, 14336, 14343, 14344, 14348, 14349, 14371, 14373, 14375, 14377, 14380, 14482, 14490, 14491, 14493, 14497--14499, 14501, 14503, 14505, 14507, 14509--14511, 14513, 14515, 14518--14522, 14526, 14528, 14529, 14531, 14532, 14543, 14564, 14565, 14567, 14568, 14572, 14575, 14576, 14578, 14579, 14584, 14585, 14588, 14589, 14591, 14595, 14599, 14600, 14603, 14607, 14608, 14613--14615, 14626, 14629, 14655, 14656, 14682, 14692--14694, 14696, 14698, 14713, 14809, 14821, 14823--14825, 14827, 14836, 14837, 14863, 14909, 14910, 14912, 14926, 14942, 14943, 14945--14947, 14950, 14953--14955, 14957, 14961, 14962, 14965, 14969, 14970, 14993, 15001, 15002, 15026, 15027, 15031--15034, 15038--15044, 15047, 15053, 15055, 15058, 15059, 15064, 15069, 15070, 15073, 15079, 15081, 15084, 15089, 15096, 15097, 15100, 15101, 15104, 15116, 15123, 15124, 15130, 15139, 15148, 15154, 15155, 15157, 15159--15167, 15169, 15170, 15172, 15175--15177, 15183, 15188, 15190, 15191, 15197, 15198, 15200--15208, 15210, 15212, 15214, 15216, 15218, 15220, 15224, 15226, 15230, 15232, 15234, 15236, 15240, 15243, 15245, 15248, 15252, 15256, 15258, 15260, 15264, 15267, 15273, 15275, 15279, 15282, 15283, 15288, 15289, 15292, 15297, 15303, 15304, 15308, 15311, 15313, 15317, 15318, 15323, 15324, 15329, 15330, 15335, 15341, 15342, 15346, 15348, 15349, 15351, 15352, 15356, 15357, 15361, 15366, 15371, 15372, 15375, 15380, 15381, 15383, 15385, 15386, 15399, 15400, 15403, 15405, 15406, 15408, 15409, 15412, 15415, 15416, 15419, 15420, 15434--15436, 15440, 15441, 15443, 15444, 15448, 15452, 15459, 15460, 15463, 15467, 15471, 15473, 15474, 15476, 15477, 15479, 15480, 15482, 15484, 15486, 15498, 15499, 15502--15509, 15511--15514, 15516, 15517, 15527, 15528, 15532, 15534--15537, 15539, 15540, 15542, 15545, 15546, 15550, 15551, 15556, 15557, 15560, 15564, 15574, 15575, 15584, 15585, 15588, 15591, 15592, 15595, 15597, 15599, 15604, 15609, 15611, 15616, 15619, 15620, 15625, 15630, 15633, 15635, 15637, 15640, 15643, 15646, 15647, 15651, 15652, 15657, 15660, 15661, 15666, 15668, 15669, 15672, 15676, 15680, 15681, 15685, 15687, 15688, 15693, 15695, 15697, 15700, 15702, 15706, 15708, 15711, 15714, 15717, 15718, 15723, 15724, 15729, 15732, 15733, 15739, 15741, 15744, 15745, 15747. 
\par

\refstepcounter{casegroup}\label{grp:3}\par\addvspace{.55ex}\noindent\textbf{Case 3.} FB.\enspace \ref{app:factors}.
16--22, 24, 56, 61, 66, 71--77, 147, 152, 157, 162--167, 170, 326, 327, 333, 335, 337, 339, 341, 342, 350--353, 358, 360, 361, 365, 369, 370, 372, 375, 378, 381, 385, 388, 391, 404--406, 408, 411, 414, 420, 423, 426, 438, 440, 444, 445, 447, 450, 453, 456, 467, 469--472, 563, 564, 569, 571, 572, 579, 583, 584, 586, 589, 590, 593, 596, 599, 601, 602, 605, 606, 609, 610, 613, 616--618, 620--622, 626, 633, 636, 639, 642, 645--647, 649, 651, 655--657, 661, 669, 672--674, 1077, 1081--1084, 1110, 1111, 1116, 1117, 1125--1129, 1132, 1135, 1138, 1144, 1157--1161, 1172, 1173, 1191, 1194--1198, 1205, 1206, 1220, 1221, 1224--1227, 1243, 1254, 1255, 1270--1272, 1274--1276, 1286, 1311--1313, 1315--1317, 1320, 1323, 1326, 1330, 1340, 1342, 1345, 1349, 1357, 1361, 1362, 1365--1367, 1375, 1377, 1379, 1381, 1388, 1390, 1392, 1408, 1412, 1417, 1431, 1438, 1443, 1448, 1456, 1458, 1460, 1462, 1464, 1466, 1468, 1476, 1482, 1506, 1510, 1517, 1533, 1537, 1544, 1565, 1573, 1574, 1587--1590, 1593, 1611, 1629, 1632, 1673, 1674, 1687--1689, 1698, 1700, 1736, 1737, 1740, 1741, 1744, 1747, 1752, 1755--1757, 1760, 1762, 1764, 1768--1772, 1795, 1799, 1800, 1804--1807, 1819, 1821, 1825, 1836, 1841--1844, 1847, 1850, 1853, 1856, 1869--1873, 1884, 1885, 1896--1898, 1900, 1901, 1903, 1906--1908, 1916, 1918, 1920, 1922, 1929, 1931, 1933, 1936, 1937, 1940, 1957, 1961, 1964, 1968, 1979, 1981, 1991, 1992, 1998, 2005, 2011, 2014, 2018, 2021, 2028, 2032, 2034--2036, 2038, 2044, 2045, 2047, 2075, 2078, 2080--2082, 2087, 2089--2091, 2093, 2096, 2112--2114, 2128--2130, 2153, 2154, 2157, 2161, 2164--2166, 2169, 2171, 2173, 2176, 2177, 2180, 2189, 2209, 2211, 2217--2219, 2221, 2223, 2224, 2226, 2244, 2249, 2292, 2362, 2364, 2368, 2369, 2372, 2375, 2378, 2381, 2384, 2393, 2396, 2401, 2404, 2410, 2413, 2445, 2448, 2451, 2452, 2455, 2615, 2616, 2633, 2634, 2646, 2647, 2655, 2660, 2673, 2675, 2681, 2683, 2689, 2691, 2721, 2729, 2731, 2733, 2735, 2737, 2766, 2767, 2816, 2817, 2831, 2833, 2837, 2841, 2842, 2856, 2863, 2869, 2871, 2875, 2879, 2896, 2905, 2907, 2912, 2917, 2922, 2927, 2948, 2949, 2956, 2960, 3021, 3025, 3116, 3121--3123, 3140, 3141, 3152, 3159, 3160, 3192, 3193, 3195, 3212, 3222, 3230, 3270, 3283, 3399, 3409, 3418, 3512, 3518, 3524, 3528, 3536, 3568, 3573, 3575--3578, 3581, 3590, 3591, 3594--3597, 3600, 3633, 3635, 3636, 3641, 3648, 3650, 3651, 3695, 3700, 3702, 3703, 3789, 3790, 3818, 3825, 3835, 3864, 3867, 3868, 3872, 3875, 3892, 3903, 3907, 3916, 3951, 3959, 3974, 3982, 4003, 4006, 4007, 4011, 4014--4016, 4027, 4030, 4033, 4049, 4061, 4064, 4066, 4069, 4072--4074, 4079, 4081, 4082, 4089, 4090, 4096, 4099, 4100, 4125, 4128, 4137, 4145, 4162, 4166, 4190, 4194, 4195, 4197, 4200, 4211, 4212, 4215, 4216, 4221, 4224, 4232, 4235, 4242--4244, 4250, 4252, 4256, 4265, 4281, 4300, 4303, 4306, 4307, 4310, 4313, 4314, 4356, 4357, 4362, 4365, 4370, 4371, 4373, 4374, 4407, 4422, 4471, 4480, 4481, 4484, 4485, 4496, 4499, 4505, 4508, 4535, 4536, 4547, 4556, 4595, 4596, 4599, 4602--4605, 4617, 4631, 4656, 4663, 4673, 4684, 4692, 4699, 4737, 4738, 4946, 4948--4950, 4954, 4958, 4962, 4964, 4965, 4990, 4996, 5005, 5016, 5018, 5019, 5028, 5029, 5100, 5101, 5108, 5152, 5181, 5198, 5217, 5240, 5249, 5252, 5276, 5279, 5286, 5287, 5308, 5310, 5331, 5334, 5335, 5360, 5367, 5371, 5374--5376, 5398, 5403, 5408, 5419, 5424, 5427, 5433, 5435, 5533, 5535, 5539, 5540, 5542, 5553, 5556--5558, 5568, 5570, 5571, 5579, 5582--5584, 5586, 5588, 5592--5594, 5598, 5608, 5611--5613, 5669, 5681, 5734, 5735, 5746, 5747, 5752, 5762, 5790, 5791, 5808, 5809, 5812, 5816, 5830, 5842, 5843, 5858--5860, 5864--5866, 5880--5882, 5884, 5886, 5888, 5897--5900, 5911, 5912, 5919, 5923--5925, 5953, 5954, 5956, 5959, 5962, 5965, 5983--5985, 5987, 5998, 5999, 6034, 6035, 6046, 6047, 6056, 6057, 6095, 6107, 6117, 6118, 6120--6122, 6126, 6128, 6130, 6131, 6146, 6148, 6178, 6180, 6218--6220, 6222, 6223, 6227, 6231, 6233, 6237--6239, 6241, 6243, 6244, 6259, 6264, 6274, 6279, 6284, 6294, 6297, 6298, 6303, 6305, 6307, 6321, 6325, 6328, 6334, 6337--6342, 6351, 6353, 6357, 6365, 6367, 6383, 6384, 6389, 6392, 6393, 6395, 6398, 6401, 6407, 6413, 6414, 6416, 6418, 6419, 6421, 6425, 6426, 6435, 6436, 6438, 6440, 6446, 6448, 6570, 6571, 6578, 6584, 6598, 6611, 6612, 6617, 6618, 6620, 6626, 6636, 6638, 6642, 6644, 6648, 6654, 6681, 6684, 6685, 6691, 6699, 6775, 6782, 6787, 6791, 6840, 6843, 6845, 6847, 6850, 6865, 6873, 6936, 6943, 6945, 6946, 6968, 6977, 6979, 6982, 6986, 6989, 7007, 7012, 7023, 7029, 7073, 7079, 7082, 7089, 7091--7093, 7103, 7104, 7106, 7109, 7110, 7112, 7115, 7139, 7141, 7178, 7192, 7220, 7221, 7223, 7240, 7248, 7249, 7257, 7260, 7279, 7711, 7714, 7716--7718, 7726, 7729, 7745, 7767, 7857, 7877, 7892, 7901, 7964, 7977, 7979, 7980, 7990, 7991, 8004, 8005, 8027, 8029, 8030, 8041, 8046, 8059, 8073, 8076, 8079, 8132, 8133, 8145--8147, 8151, 8160, 8170, 8171, 8186, 8187, 8280, 8289, 8316, 8317, 8332, 8333, 8350, 8352, 8376, 8377, 8380, 8419, 8421, 8424, 8434, 8437, 8438, 8441, 8442, 8472--8476, 8480, 8481, 8483, 8486, 8487, 8491, 8493, 8527, 8528, 8530, 8531, 8534, 8556, 8557, 8560, 8561, 8582, 8587, 8592, 8635, 8639, 8642, 8645, 8658, 8661, 8662, 8664, 8667, 8708, 8714, 8717, 8724, 8776, 8777, 8805, 8811, 8822, 8824, 8830, 9050--9056, 9058, 9079, 9082, 9085, 9087--9089, 9091, 9093, 9094, 9113, 9116, 9118, 9120, 9122--9125, 9127, 9130, 9165, 9168--9171, 9174, 9177--9180, 9183, 9186, 9192, 9197, 9199, 9203, 9209, 9210, 9218--9220, 9230, 9231, 9238, 9239, 9242, 9252, 9253, 9264, 9265, 9268, 9274, 9279, 9281, 9286, 9288--9290, 9293--9296, 9299, 9300, 9311, 9312, 9315, 9321--9323, 9335--9338, 9340, 9342, 9343, 9345, 9352, 9355, 9362, 9372--9374, 9379, 9385, 9386, 9394, 9427, 9430, 9439, 9442, 9449, 9450, 9459, 9461, 9470, 9471, 9473, 9475, 9476, 9478, 9488, 9489, 9497, 9501, 9502, 9505--9507, 9515, 9520, 9524, 9527, 9532, 9534, 9535, 9537, 9540, 9542, 9543, 9546, 9547, 9549, 9552, 9562, 9565, 9566, 9570, 9573, 9576, 9579, 9582, 9612, 9613, 9616, 9627, 9629, 9631, 9634, 9656, 9660, 9661, 9665--9667, 9671, 9672, 9674, 9678, 9679, 9681, 9697, 9699, 9700, 9724, 9726, 9727, 9734, 9735, 9739, 9740, 9754, 9757--9759, 9763, 9769, 9791, 9793, 9795, 9799, 9802, 9805, 9822, 9828, 9829, 9832--9834, 9836, 9839, 9854, 9858, 9864, 9870, 9873, 9879, 9880, 9914, 9915, 9933--9936, 9943, 9962, 9963, 9965, 9984, 9985, 9995--9997, 9999, 10000, 10013--10015, 10030, 10031, 10327--10331, 10342, 10344--10346, 10348--10350, 10353, 10363, 10365, 10366, 10368, 10370--10372, 10375, 10381, 10382, 10384, 10385, 10387, 10389, 10391, 10442, 10444, 10445, 10455--10458, 10463, 10465, 10466, 10473, 10478, 10479, 10481, 10483, 10484, 10487, 10488, 10490, 10491, 10494, 10802--10805, 10808, 10811, 12279, 12288, 13146, 13282, 13288, 13307, 13313, 13341, 13948--13950, 13954--13956, 14334, 14927, 15388, 15396, 15402, 15571, 15577, 15682--15684, 15734--15736. 
\par

\refstepcounter{casegroup}\label{grp:4}\par\addvspace{.55ex}\noindent\textbf{Case 4.} FB.\enspace Section~\ref{sec:four-sources}.
\enspace \(S_{4,42}\)\enspace 
36, 37, 198, 291, 292, 318, 457, 465, 1703, 1711, 1722, 1781, 1783, 2942, 5989, 6347, 10533, 12195, 13246, 13277, 13279. 
\par

\refstepcounter{casegroup}\label{grp:5}\par\addvspace{.55ex}\noindent\textbf{Case 5.} FB.\enspace \ref{case:1}.
\enspace \(S_{4,705}\)\enspace 
38--40, 50, 52, 199, 293--295, 320, 324, 466, 5928, 5993, 6348, 12194, 12276, 13197, 13835, 13844, 13845, 13934, 13937. 
\par

\refstepcounter{casegroup}\label{grp:6}\par\addvspace{.55ex}\noindent\textbf{Case 6.} FB.\enspace \ref{thm:rectangle}.
41--45, 53--55, 57, 59, 67, 69, 94--96, 99, 100, 107--110, 113, 115, 116, 120, 121, 123, 124, 126, 133--137, 145, 180--182, 215, 216, 223, 224, 243, 244, 249--251, 267, 268, 283, 287, 298, 300--309, 321--323, 325, 330, 331, 346, 348, 354, 355, 357, 359, 364, 394, 400, 401, 415, 416, 479, 480, 483, 491, 493--506, 510--512, 516--518, 522--525, 528, 529, 533, 539--541, 556, 560, 561, 568, 619, 653, 654, 658, 660, 663, 666--668, 670, 675, 677, 678, 680, 687, 691, 708--710, 712, 719, 720, 722, 771, 772, 780, 782, 788, 789, 795, 849, 853, 861, 867, 870, 871, 873--875, 879--886, 908, 913, 914, 917, 925, 930--933, 936--940, 943, 976, 980, 982, 986, 988, 994, 1008, 1009, 1012, 1013, 1022, 1026, 1071--1074, 1076, 1078, 1080, 1085, 1087, 1088, 1090, 1092, 1094, 1095, 1097, 1098, 1100--1103, 1105, 1109, 1121, 1122, 1124, 1174, 1176, 1178, 1179, 1181, 1185--1187, 1189, 1190, 1202, 1204, 1214--1216, 1219, 1229, 1231, 1233--1235, 1237, 1238, 1240--1242, 1251, 1252, 1256, 1258--1263, 1265--1267, 1269, 1283, 1285, 1294, 1296, 1319, 1333, 1336, 1353, 1359, 1369, 1370, 1373, 1374, 1382, 1386, 1387, 1395, 1399, 1463, 1500, 1502, 1503, 1512, 1513, 1515, 1516, 1557, 1569, 1571, 1605, 1607, 1669, 1671, 1704, 1714, 1723, 1743, 1784, 1808, 1810, 1811, 1813, 1815, 1817, 1827, 1829, 1833, 1840, 1910, 1911, 1914, 1915, 1923, 1924, 1927, 1928, 1942, 1946, 1954, 1975, 2106--2109, 2115, 2116, 2119, 2122--2125, 2195, 2202--2207, 2210, 2214, 2228--2234, 2236, 2238--2242, 2245, 2247, 2250, 2252, 2256--2262, 2264, 2265, 2268, 2269, 2272--2274, 2276, 2278, 2279, 2281, 2283, 2284, 2286, 2289, 2293, 2295, 2297--2301, 2317--2322, 2333--2337, 2351, 2353, 2355--2360, 2366, 2367, 2382, 2383, 2387, 2390, 2416, 2417, 2419, 2421, 2423, 2425--2429, 2432, 2441, 2443, 2458--2461, 2466--2468, 2478--2484, 2486, 2518, 2519, 2522, 2523, 2525, 2542, 2545, 2549, 2554, 2561--2563, 2565--2567, 2569--2572, 2579--2582, 2587--2589, 2617--2627, 2648--2651, 2653, 2654, 2657, 2664, 2665, 2677--2679, 2685, 2687, 2690, 2693--2695, 2754, 2755, 2775, 2776, 2804, 2805, 2819, 2820, 2846, 2847, 2857--2861, 2873, 2874, 2877, 2878, 2881, 2937, 2968, 2990, 3016, 3124, 3126, 3161--3166, 3188--3191, 3194, 3231--3234, 3236--3238, 3263--3266, 3284--3290, 3318--3321, 3353, 3354, 3385, 3398, 3419, 3441, 3455, 3471, 3525, 3537, 3550, 3565, 3566, 3574, 3584--3589, 3592, 3593, 3601--3616, 3619, 3620, 3622, 3630, 3632, 3634, 3639, 3640, 3642--3645, 3647, 3649, 3654--3656, 3659, 3667, 3668, 3671--3675, 3678--3683, 3685--3689, 3692, 3696, 3697, 3699, 3701, 3706--3709, 3716, 3718--3727, 3729--3734, 3737--3741, 3744, 3788, 3792--3794, 3801--3804, 3806, 3811--3815, 3817, 3820--3822, 3824, 3826--3828, 3830--3832, 3834, 3837, 3844, 3846, 3855--3857, 3862, 3863, 3869, 3871, 3878, 3888, 3894, 3896, 3899, 3905, 3906, 3909, 3912, 3917--3919, 3921, 3924, 3932, 3945--3947, 3953--3955, 3989--3991, 3995, 4002, 4008, 4010, 4022--4024, 4026, 4029, 4038, 4039, 4042, 4045, 4053, 4054, 4057, 4060, 4068, 4078, 4124, 4204, 4271, 4317, 4320, 4321, 4372, 4424--4426, 4456, 4457, 4473, 4479, 4632, 4633, 4635, 4735, 4739, 4743, 4744, 4766, 4883, 4887, 4891, 4892, 4937, 4941, 4944, 4945, 4947, 4955, 4957, 4959--4961, 4963, 4968, 5027, 5031--5037, 5043, 5044, 5093, 5096, 5099, 5103--5105, 5107, 5195, 5196, 5200--5204, 5238, 5242--5245, 5372, 5409, 5425, 5436, 5441--5443, 5453, 5464, 5472, 5475, 5476, 5480, 5487, 5492, 5493, 5500, 5501, 5504, 5507, 5512, 5513, 5516, 5517, 5521, 5522, 5528, 5531, 5532, 5537, 5538, 5541, 5543, 5545, 5548--5552, 5554, 5559--5565, 5567, 5569, 5572, 5574, 5575, 5577, 5578, 5580, 5581, 5590, 5591, 5595, 5597, 5599--5607, 5609, 5614, 5620, 5622--5624, 5626, 5639, 5658--5660, 5662, 5663, 5665, 5666, 5676, 5678, 5695--5698, 5700--5703, 5705--5707, 5712--5714, 5720, 5725, 5726, 5731, 5732, 5736, 5738, 5739, 5741--5743, 5745, 5749, 5750, 5753, 5754, 5756--5758, 5760, 5761, 5763, 5764, 5766, 5768, 5770--5772, 5774, 5775, 5777--5779, 5787, 5788, 5792, 5794--5801, 5803--5805, 5807, 5814, 5817, 5818, 5821--5823, 5825, 5826, 5828, 5829, 5831, 5839, 5840, 5844, 5846--5851, 5853--5855, 5857, 5867, 5869, 5872, 5874, 5878, 6036--6043, 6147, 6149, 6154, 6185, 6188, 6189, 6299, 6300, 6308, 6312, 6313, 6442--6445, 6450, 6452, 6454--6457, 6701, 6723, 6844, 6846, 6853, 6861, 6867--6869, 7163, 7366, 7367, 7370, 7434, 7436, 7441--7443, 7703--7705, 7721, 7722, 7725, 7765, 7769--7774, 7776, 7777, 7962, 7966--7968, 7972--7974, 8047, 8103, 8104, 8107, 8108, 8129, 8130, 8134, 8135, 8137, 8139--8142, 8149, 8152, 8154, 8155, 8157--8159, 8161, 8167, 8168, 8172, 8174--8179, 8181--8183, 8185, 8346--8348, 8351, 8360--8369, 8594, 8663, 8808, 8812, 8814--8820, 8823, 8827, 8829, 9070, 9071, 9099, 9103, 9104, 9137, 9141, 9155, 9182, 9319, 9320, 9324, 9326, 9329--9332, 9347--9350, 9357, 9359, 9364--9366, 9503, 9504, 9525, 9533, 9541, 9544, 9545, 9548, 9554, 9561, 9567, 9569, 9578, 9703, 9771, 9823, 9826, 9827, 9835, 9878, 9882--9885, 9966, 9972, 9973, 9978, 9979, 9982, 9983, 9986, 9988--9992, 9998, 10001, 10003, 10004, 10006--10008, 10011, 10012, 10016, 10018--10023, 10025--10027, 10029, 10380, 10393, 10394, 10482, 10489, 10534--10537, 10541, 10543, 10754, 10901, 10903, 10905, 10907, 10909, 10913, 10917, 10918, 10923, 10931, 10936, 10937, 10940, 10946, 10947, 10952, 10953, 10959, 10963, 10967, 11000, 11006, 11010, 11014, 11022, 11118, 11123, 11125, 11126, 11148, 11150, 11151, 11154, 11162, 11163, 11165, 11168, 11169, 11172, 11173, 11176, 11177, 11180, 11182--11184, 11188--11194, 11204, 11205, 11207, 11208, 11211, 11241, 11273, 11291--11293, 11320, 11325, 11326, 11338, 11342, 11367, 11376, 11381, 11382, 11385, 11386, 11595, 11596, 11600, 11611--11613, 11620, 11621, 11623--11626, 11628, 11629, 11632, 11644, 11647, 11648, 11651, 11656--11658, 11661, 11662, 11689, 11694, 11696, 11699, 11702, 11707--11716, 11720, 11721, 11724--11729, 11731--11733, 11738--11742, 11744--11747, 11750--11760, 11762--11764, 11766--11768, 11814--11816, 11818, 11822, 11823, 11826, 11830, 11865, 11866, 11899, 11903, 11905, 11908, 11911, 12012, 12013, 12120--12125, 12136--12140, 12183, 12200, 12201, 12222, 12223, 12227, 12229, 12243, 12268, 12269, 12282, 12284--12286, 12375, 12377, 12379, 12382, 12383, 12386, 12398, 12449, 12451, 12453, 12456, 12457, 12460, 12501, 12502, 12527, 12528, 12591, 12678, 12681, 12682, 12733, 12740, 12741, 12745--12747, 12749--12753, 12756, 12757, 12874, 12875, 12882, 12883, 12888--12890, 12893--12897, 12900, 12901, 12926, 12944, 12945, 13014, 13015, 13021--13024, 13026--13029, 13032, 13033, 13087, 13097, 13109, 13121, 13159, 13171, 13172, 13174, 13186, 13199, 13200, 13214, 13232, 13234, 13235, 13243, 13247, 13248, 13264, 13291, 13294, 13297, 13303--13306, 13462, 13517, 13520, 13521, 13550, 13553, 13554, 13562--13565, 13600, 13616, 13662, 13665--13667, 13682--13684, 13698, 13735, 13736, 13750, 13758--13760, 13783, 13785, 13799, 13802, 13805, 13808, 13811, 13820, 13821, 13823, 13826, 13833, 13836, 13846, 13847, 13854, 13858, 13862, 13867, 13879, 13880, 13884, 13887, 13899, 13913, 13917, 13920, 13923, 13928, 13931, 13935, 13938, 13944, 13957, 13958, 13962, 13963, 13965, 13968, 13972--13975, 13977--13984, 13986, 13988, 13990, 13992--13994, 13998, 13999, 14002, 14008--14013, 14021, 14022, 14031, 14035--14039, 14042, 14049--14051, 14053, 14054, 14057--14064, 14073--14075, 14082, 14099, 14106, 14111--14114, 14116--14121, 14124, 14144--14155, 14157, 14158, 14163, 14164, 14167--14173, 14175, 14176, 14179--14189, 14191--14193, 14195--14198, 14201, 14207--14210, 14213, 14219--14221, 14223, 14224, 14228, 14232, 14236--14239, 14241, 14242, 14245, 14250, 14263, 14264, 14272--14277, 14289--14292, 14306, 14310, 14311, 14315, 14317, 14326, 14327, 14331, 14332, 14345, 14347, 14387, 14389, 14390, 14420, 14424, 14428, 14432, 14454, 14455, 14457, 14475--14477, 14500, 14502, 14504, 14506, 14508, 14512, 14516, 14517, 14530, 14542, 14548, 14551, 14558, 14562, 14580, 14581, 14583, 14586, 14592--14594, 14596, 14598, 14601, 14602, 14604, 14606, 14609, 14621, 14625, 14630, 14644, 14645, 14651, 14652, 14654, 14657, 14675, 14677, 14679, 14680, 14685, 14695, 14699, 14714, 14719, 14724--14729, 14731--14733, 14736, 14737, 14751, 14755, 14759, 14763, 14770, 14771, 14785--14790, 14802--14808, 14810--14812, 14822, 14831--14833, 14835, 14838, 14867--14869, 14884, 14885, 14894, 14905, 14920--14925, 14928, 14929, 14934, 14996, 15000, 15020, 15021, 15045, 15105, 15106, 15110, 15111, 15114, 15115, 15117, 15118, 15121, 15122, 15125--15127, 15147, 15217, 15233, 15246, 15247, 15259, 15271, 15272, 15274, 15276--15278, 15281, 15293--15296, 15298--15302, 15305--15307, 15309, 15310, 15312, 15314--15316, 15327, 15328, 15331--15334, 15336--15340, 15343--15345, 15347, 15350, 15360, 15362--15365, 15367--15370, 15373, 15374, 15382, 15391, 15393--15395, 15414, 15418, 15430, 15433, 15437, 15447, 15449--15451, 15453, 15454, 15457, 15458, 15461, 15462, 15466, 15485, 15494, 15495, 15515, 15523, 15524, 15538, 15541, 15544, 15549, 15561--15563, 15565, 15566, 15570, 15572, 15573, 15576, 15583, 15642, 15644, 15645, 15648--15650, 15653--15656, 15658, 15659, 15662--15665, 15667, 15670, 15671, 15673--15675, 15677, 15694, 15696, 15698, 15699, 15701, 15713, 15715, 15731. 
\par

\refstepcounter{casegroup}\label{grp:7}\par\addvspace{.55ex}\noindent\textbf{Case 7.} FB.\enspace \ref{case:2}.
\enspace \(\A{51}\)\enspace 
51, 13886. 
\par

\refstepcounter{casegroup}\label{grp:8}\par\addvspace{.55ex}\noindent\textbf{Case 8.} FB.\enspace Section~\ref{sec:four-sources}.
\enspace \(S_{4,393}\)\enspace 
62--65, 340, 343--345, 371, 379, 380, 384, 424, 425, 468, 1089, 1091, 1115, 1168--1171, 1282, 1314, 1328, 1329, 1378, 1535, 1536, 2405, 2406, 3213, 3980, 3981, 4144, 4225, 4363, 5180, 5490, 5931, 5995, 6099, 6207, 6215, 6263, 6354, 6579, 6944, 6966, 6967, 7890, 7891, 8117, 8328--8331, 12274, 12277, 13222, 13226, 13326, 13338, 13852, 13856, 13869, 13876, 13877, 13941, 14524, 14545, 14554, 14617, 14618, 14634, 14665--14668, 14697, 14718, 15543, 15554, 15555. 
\enspace \(S_{4,453}\)\enspace 
153--156, 573--575, 585, 600, 607, 608, 1814, 1823, 1882, 1883, 1899, 1932, 2017, 2019, 2020, 4236, 6283, 6368, 12270, 13320, 13335. 
\enspace \(S_{4,427}\)\enspace 
433, 436, 576, 1550, 1555, 1887, 2408, 6196, 6204, 6888, 6910, 6932, 7880, 8262, 8266, 8268, 8269, 14287, 15392. 
\enspace \(S_{4,428}\)\enspace 
550, 552, 578, 1886, 2439, 6166, 6177, 6733, 6741, 6749, 7463, 7533, 7813, 7814, 7816, 14278, 15377, 15378. 
\enspace \(S_{4,464}\)\enspace 
2918--2921, 4083--4085, 4196, 4203, 4486, 4487, 4500, 5871, 5909, 5910, 6140, 6221, 6621, 6639, 6985, 6987, 6988, 7011, 11969, 12054, 13132. 
\enspace \(S_{4,441}\)\enspace 
2923--2926, 4091--4093, 4213, 4214, 4506, 4507, 5277, 5278, 5280, 5281, 5920--5922, 6127, 6236, 6643, 7021, 7022, 8417, 8418, 8422, 8423, 9090, 9092, 9250, 9251, 9438, 9440, 9441, 9443, 9628, 9632, 9633, 10351, 10352, 10582, 10585, 10593, 12021, 12035, 12153, 12154, 12165, 12179--12181, 12256, 12278, 13071, 13339, 13340. 
\enspace \(S_{4,445}\)\enspace 
2933--2936, 2985, 2988, 3096, 3098, 4106, 4112--4114, 4178, 4188, 4220, 4240, 4512, 4513, 5940, 6003, 6030, 6031, 6249, 6290, 6647, 7034, 7035, 11843, 12079, 13271. 
\enspace \(S_{4,417}\)\enspace 
2954, 3197, 3200, 4088, 4134, 5111, 5916, 6206, 6213, 6939, 6960, 7851, 7887, 8113, 8304, 8306, 8307, 12096, 13289. 
\enspace \(S_{4,431}\)\enspace 
2962, 3047, 3136--3139, 3210, 3224, 3356, 3408, 3415, 4097, 4147, 4331, 4340, 4642, 4661, 4662, 4682, 4683, 4988, 4989, 4993--4995, 5150, 5981, 6195, 6202, 6757, 6781, 6883, 6904, 7839, 7840, 7875, 8233, 8243, 8245, 8246, 11842, 11984, 12083, 13280. 
\enspace \(S_{4,434}\)\enspace 
2976, 3246, 3257, 4105, 4647, 4791, 5068, 5974, 6165, 6175, 6732, 6740, 7462, 7477, 7531, 7807--7810, 12078, 13265, 13267, 13269, 13274. 
\enspace \(S_{4,418}\)\enspace 
2994, 3291, 3293, 4117, 5915, 6205, 6212, 6935, 6956, 7847, 7882, 8110, 8272, 8274, 8276, 12033, 13044, 13049. 
\enspace \(S_{4,432}\)\enspace 
3000, 3301, 3312, 4119, 4983, 5132, 5139, 5976, 6193, 6200, 6879, 6900, 7834, 7863, 8216--8220, 12032, 13037, 13043, 13048, 13259. 
\enspace \(S_{4,433}\)\enspace 
3131--3134, 4328, 4494, 4612, 4638, 4652, 4658, 4659, 4977--4979, 6191, 6625, 6688, 6755, 6875, 7000, 7733, 7781, 7818, 7831, 7832, 8212, 8213, 11935, 12049, 13125. 
\enspace \(S_{4,456}\)\enspace 
4293, 5316, 5318, 5320, 6137, 6424, 7084, 8397, 8456, 8484, 9675, 10158, 10185. 
\enspace \(S_{4,447}\)\enspace 
4515, 4518, 7038, 7044, 8406, 8846, 8847, 8870, 8878, 9595, 10146, 10147, 10168, 10169, 10183, 11992, 12047, 12052, 13123. 
\enspace \(S_{4,442}\)\enspace 
4527, 4530, 7060, 7065, 8415, 9003, 9004, 9013, 9016, 9621, 10283, 10284, 10287, 10288, 10293, 11996, 12046, 12055, 13122. 
\enspace \(S_{4,450}\)\enspace 
4566, 4569, 4783, 4852, 4856, 6170, 6737, 7041, 7474, 7524, 7616, 7791, 7792, 7916, 8403, 8624, 9777, 10144, 10180. 
\enspace \(S_{4,446}\)\enspace 
5259, 5262, 8389, 8394, 8411, 8864, 8865, 8875, 8883, 9600, 10160, 10161, 10172, 10173, 10188, 12037, 12176, 12271, 13074. 
\enspace \(S_{4,443}\)\enspace 
7244, 7246, 7373, 7421, 7450, 7465, 7472, 7664, 7665, 8622, 9775, 10281, 10290. 
\enspace \(S_{4,9}\)\enspace 
10083, 10517, 13362, 13365, 14846, 15137, 15579, 15581. 
\enspace \(S_{4,401}\)\enspace 
11403, 11530, 11533, 11584, 12707, 12725, 12730, 13502, 13969, 15048, 15067, 15071, 15629. 
\enspace \(S_{4,405}\)\enspace 
11418, 11565, 11568, 11579, 11645, 12758, 12776, 12781, 12842, 12862, 12866, 13522, 13545. 
\enspace \(S_{4,398}\)\enspace 
11421, 11587, 11588, 11590, 12685, 12700, 12705, 13497, 14720, 15074, 15094, 15098, 15639. 
\par

\refstepcounter{casegroup}\label{grp:9}\par\addvspace{.55ex}\noindent\textbf{Case 9.} FB.\enspace \ref{case:3}.
83, 97, 183, 1798. 
\par

\refstepcounter{casegroup}\label{grp:10}\par\addvspace{.55ex}\noindent\textbf{Case 10.} FB.\enspace \ref{case:4}.
\enspace \(\A{84}\)\enspace 
84, 13901. 
\enspace \(\A{888}\)\enspace 
98, 888, 13273. 
\enspace \(\A{888}^{\rm op}\)\enspace 
141, 184, 13252. 
\par

\refstepcounter{casegroup}\label{grp:11}\par\addvspace{.55ex}\noindent\textbf{Case 11.} FB.\enspace \ref{thm:qgraph}.
85, 86, 142, 192, 193, 213, 214, 644, 664, 665, 889, 892, 2215, 2216, 2254, 2255, 2288, 2303, 2304, 2548, 2550, 3582, 3583, 5518, 5524, 5525, 8809, 8810, 12283, 13253, 13254, 13300, 13896, 13903, 14330, 14978, 14979, 15013, 15014, 15389. 
\par

\refstepcounter{casegroup}\label{grp:12}\par\addvspace{.55ex}\noindent\textbf{Case 12.} FB.\enspace \ref{case:5}.
\enspace \(\A{91}\)\enspace 
91, 1719, 13185. 
\enspace \(\A{92}\)\enspace 
92, 311, 13832. 
\enspace \(\A{93}\)\enspace 
93, 13927. 
\enspace \(\A{91}^{\rm op}\)\enspace 
177, 1713, 12199. 
\enspace \(\A{92}^{\rm op}\)\enspace 
178, 297, 12198. 
\enspace \(\A{93}^{\rm op}\)\enspace 
179, 13898. 
\par

\refstepcounter{casegroup}\label{grp:13}\par\addvspace{.55ex}\noindent\textbf{Case 13.} FB.\enspace Section~\ref{sec:four-sources}.
\enspace \(S_{4,292}\)\enspace 
103, 128, 347, 397, 484, 487, 508, 514, 520, 544, 1384, 1558, 1570, 1606, 1670, 2975, 4276, 6360, 6582. 
\enspace \(S_{4,290}\)\enspace 
219, 714, 781, 792, 918, 921, 927, 978, 984, 1016, 1177, 1180, 1188, 1284, 1397, 2999, 4263, 6373, 6588. 
\enspace \(S_{4,376}\)\enspace 
854, 2543, 11314, 11594, 11639, 11684--11688, 11690--11693, 11695, 11697, 11698, 11706, 11719, 12320--12322, 12581, 12584, 12587, 12881, 13409, 13412, 15109. 
\enspace \(S_{4,368}\)\enspace 
2354, 14403, 14410, 14411, 14427, 14429--14431, 14474, 14701, 14738, 14742, 14795, 14799, 15493, 15522. 
\enspace \(S_{4,357}\)\enspace 
2476, 11395, 11429, 11435, 11436, 11438, 11441, 11444, 11500, 11502, 12351, 12359, 12361, 12715, 12785, 13535. 
\enspace \(S_{4,358}\)\enspace 
2500, 11389, 11391--11393, 11407, 11410, 11425, 11426, 11428, 12388, 12400, 12402, 12771, 12784, 13622. 
\enspace \(S_{4,360}\)\enspace 
2537, 11633, 11637, 11663, 11665, 11666, 11677, 11680--11682, 11704, 11705, 12301, 12302, 12529, 12531--12533, 12742, 12877, 12886, 13382, 13383, 15107, 15112. 
\enspace \(S_{4,363}\)\enspace 
3235, 3239, 12059, 12060, 12097--12099, 12101, 12105--12107, 12111, 12113, 12114, 12127, 12129, 12246, 13053--13055, 13096, 13154, 13155, 13158, 13173, 13292, 13301, 13798, 13819. 
\enspace \(S_{4,379}\)\enspace 
6663, 6665, 8461, 8871, 8915, 8916, 8919, 10124, 10202, 10218, 10252, 10400, 10618, 10637, 10731. 
\enspace \(S_{4,380}\)\enspace 
6675, 6677, 8470, 9014, 9026--9028, 10277, 10304, 10306, 10310, 10446, 10762, 10765, 10781. 
\enspace \(S_{4,383}\)\enspace 
8444, 8446, 8462, 8876, 8917, 8918, 8920, 10135, 10205, 10221, 10253, 10403, 10631, 10643, 10733. 
\enspace \(S_{4,385}\)\enspace 
9022, 9031, 9032, 10049, 10297, 10314, 10315, 10317, 10787, 10792. 
\par

\refstepcounter{casegroup}\label{grp:14}\par\addvspace{.55ex}\noindent\textbf{Case 14.} FB.\enspace \ref{case:6}.
\enspace \(\A{928}^{\rm op}\)\enspace 
104, 118, 129--132, 144, 398, 464, 474, 485, 490, 509, 515, 521, 545, 1093, 1165, 1385, 2438, 3249, 3943, 5083, 5463, 5936, 5992, 6101, 6361, 6583, 6954, 8325, 12275, 13153, 13165, 13169, 13797, 13801, 13804, 13813, 13817, 13818, 13919, 13922, 14005, 14015, 14024, 14067, 14068, 14084, 14087, 14088, 14101, 14127, 14130, 15284, 15319, 15320, 15353. 
\par

\refstepcounter{casegroup}\label{grp:15}\par\addvspace{.55ex}\noindent\textbf{Case 15.} FB.\enspace \ref{case:7}.
\enspace \(\A{690}^{\rm op}\)\enspace 
105, 106, 119, 488, 489, 13822, 13916, 13930, 13966, 14128, 14129, 15285, 15442. 
\par

\refstepcounter{casegroup}\label{grp:16}\par\addvspace{.55ex}\noindent\textbf{Case 16.} FB.\enspace Section~\ref{sec:four-sources}.
\enspace \(S_{4,292}\)\enspace 
111, 122, 138, 139, 146, 349, 407, 442, 446, 473, 475, 478, 530, 531, 555, 1154, 1659, 1667, 1668, 2437, 3248, 3256, 3503, 3944, 4104, 4445, 4648, 4701, 4792, 4845, 5069, 5084, 5975, 6245, 6792, 6813, 6955, 7198, 7532, 7554, 7582, 7583, 7606, 7905, 7906, 8080, 8313, 13266, 13268, 13270. 
\enspace \(S_{4,290}\)\enspace 
227, 443, 693, 724, 725, 783, 807, 809, 813, 815, 826, 899, 1004, 1005, 1021, 1155, 1694--1696, 2585, 3303, 3311, 3976, 3978, 4118, 4414, 4984, 5007, 5133, 5140, 5163, 5170, 5977, 6253, 6937, 6958, 6964, 7864, 7884, 8281, 8283, 8284, 8286--8288, 8290, 8314, 12498, 12525, 12942. 
\enspace \(S_{4,357}\)\enspace 
2344, 3352, 3764, 3920, 5462, 5650, 10832, 11027, 11034, 11055, 11335, 11404, 11434, 11484, 11494, 11495, 11531, 11534, 11573, 12800. 
\enspace \(S_{4,358}\)\enspace 
2493, 3775, 3780, 3805, 5684, 5687, 10836, 11042, 11058, 11061, 11419, 11513, 11556, 11559, 11562, 11563, 11566, 11569, 11574, 12807. 
\enspace \(S_{4,379}\)\enspace 
4514, 4517, 5267, 7039, 7045, 7118, 7119, 8407, 8463, 8848, 8849, 8879, 8880, 8925, 8926, 9406, 9597, 9645, 10148, 10170, 10184, 10223, 10239, 10240, 10255. 
\enspace \(S_{4,380}\)\enspace 
4526, 4529, 5273, 7061, 7066, 7130, 7131, 8416, 8471, 9005, 9006, 9017, 9018, 9029, 9030, 9433, 9623, 9663, 10285, 10289, 10294, 10307--10309, 10311. 
\enspace \(S_{4,383}\)\enspace 
5258, 5261, 5270, 8390, 8395, 8412, 8449, 8450, 8464, 8866, 8867, 8884, 8885, 8927, 8928, 9408, 9602, 9646, 10162, 10174, 10189, 10226, 10241, 10242, 10256. 
\enspace \(S_{4,385}\)\enspace 
9009, 9010, 9025, 9034, 10039, 10300, 10302, 10508. 
\par

\refstepcounter{casegroup}\label{grp:17}\par\addvspace{.55ex}\noindent\textbf{Case 17.} FB.\enspace \ref{case:8}.
\enspace \(\A{228}^{\rm op}\)\enspace 
112, 476, 486, 1295, 8126, 13275. 
\par

\refstepcounter{casegroup}\label{grp:18}\par\addvspace{.55ex}\noindent\textbf{Case 18.} FB.\enspace \ref{case:9}.
114, 477, 1501, 8125, 13211, 14636. 
\par

\refstepcounter{casegroup}\label{grp:19}\par\addvspace{.55ex}\noindent\textbf{Case 19.} FB.\enspace \ref{thm:one-sided}.
148, 150, 158, 160, 567, 615, 625, 662, 1801, 1803, 1822, 1824, 1826, 2010, 2074, 2120, 2121, 2127, 9126, 9327, 9328, 9334, 10752, 14217, 14225, 14227, 14231, 14258, 14261, 14262, 14267, 14268, 14280, 14283--14286, 14293--14296, 14298--14301, 14305, 14309, 14328, 15379, 15384, 15387, 15397, 15398, 15401, 15404, 15679. 
\par

\refstepcounter{casegroup}\label{grp:20}\par\addvspace{.55ex}\noindent\textbf{Case 20.} NFB.\enspace \ref{thm:odd-meet}.
\enspace \(\text{\normalfont Scab\_times\_D2}\)\enspace 
190, 211, 1797. 
\par

\refstepcounter{casegroup}\label{grp:21}\par\addvspace{.55ex}\noindent\textbf{Case 21.} NFB.\enspace \ref{case:10}.
191, 212, 891, 13249, 13250, 13894, 13895. 
\par

\refstepcounter{casegroup}\label{grp:22}\par\addvspace{.55ex}\noindent\textbf{Case 22.} FB.\enspace \ref{case:11}.
\enspace \(\A{194}\)\enspace 
194, 197, 1720. 
\par

\refstepcounter{casegroup}\label{grp:23}\par\addvspace{.55ex}\noindent\textbf{Case 23.} FB.\enspace \ref{case:12}.
\enspace \(\A{200}\)\enspace 
200, 13889. 
\par

\refstepcounter{casegroup}\label{grp:24}\par\addvspace{.55ex}\noindent\textbf{Case 24.} FB.\enspace \ref{case:13}.
\enspace \(\A{201}\)\enspace 
201, 12203. 
\enspace \(\A{201}^{\rm op}\)\enspace 
2307, 12197. 
\par

\refstepcounter{casegroup}\label{grp:25}\par\addvspace{.55ex}\noindent\textbf{Case 25.} FB.\enspace \ref{case:14}.
\enspace \(\A{202}\)\enspace 
202, 12202. 
\enspace \(\A{767}\)\enspace 
767, 12433. 
\enspace \(\A{774}\)\enspace 
774, 12331. 
\enspace \(\A{767}^{\rm op}\)\enspace 
900, 12913. 
\enspace \(\A{906}\)\enspace 
906, 12916, 12997. 
\enspace \(\A{1042}\)\enspace 
1042, 1069, 1830, 1949, 9153, 9325, 10854, 10882, 11104. 
\enspace \(\A{1042}^{\rm op}\)\enspace 
1051, 2612, 5876, 6317, 9369, 9987, 10876, 10883, 10990. 
\enspace \(\A{202}^{\rm op}\)\enspace 
2308, 12196. 
\enspace \(\A{906}^{\rm op}\)\enspace 
2558, 12917, 12938. 
\enspace \(\A{9147}\)\enspace 
9147, 10885. 
\enspace \(\A{13426}\)\enspace 
13426, 13715. 
\enspace \(\A{13691}\)\enspace 
13691, 13745. 
\enspace \(\A{13691}^{\rm op}\)\enspace 
13776, 13778. 
\par

\refstepcounter{casegroup}\label{grp:26}\par\addvspace{.55ex}\noindent\textbf{Case 26.} FB.\enspace \ref{thm:bounded}.
203, 210, 239, 240, 263, 264, 402, 699, 706, 736, 737, 748, 762, 764, 769, 778, 805, 806, 846, 902, 960, 961, 967, 968, 975, 1002, 1034, 1040, 1046, 1055, 1406, 1953, 2160, 2309, 2314, 2327, 2328, 2341, 2342, 2471, 2475, 2490, 2491, 2499, 2504, 2506, 2512, 2513, 2515, 2524, 2529, 2530, 2539, 2556, 2557, 2560, 2576, 2592, 2596, 2600, 2602, 2604, 2606, 2698, 2884, 3747, 3748, 3751, 3769, 3770, 3860, 3928, 4036, 4610, 6459, 6461, 6857, 8431, 9098, 9109, 9111, 9136, 9259, 9260, 9271, 9285, 9307, 9308, 9423, 9425, 9607, 9608, 9653, 9723, 9910, 9911, 10053, 10060, 10272--10274, 10359, 10361, 10362, 10427--10429, 10438, 10440, 10441, 10835, 10844, 10845, 10850, 10851, 10853, 10856, 10864, 10865, 10870, 10893, 10894, 10928, 10945, 10969, 11036, 11038, 11067, 11068, 11071, 11072, 11075, 11077, 11078, 11086, 11089, 11090, 11303, 11306--11308, 11310, 11329, 11804, 11820, 11821, 11825, 11924, 12309, 12310, 12317, 12328, 12329, 12334, 12335, 12342, 12350, 12354, 12364, 12366, 12405, 12408, 12410, 12420, 12422, 12428--12431, 12437, 12439, 12446, 12463, 12464, 12821, 12823, 12826, 12831, 12832, 12837, 12838, 12903, 12907--12910, 12920, 12958, 12964, 12966, 12969, 12970, 12975, 12980, 12981, 13393, 13394, 13399, 13402, 13424, 13425, 13444, 13447, 13448, 13541, 13542, 13572, 13586, 13592, 13607, 13609, 13638, 13640, 13641, 13644, 13648, 13657, 13687, 13693, 13709, 13710, 13714, 13716, 13719, 13721, 13723, 13724, 13737, 13739, 13771, 13774, 14214, 14235, 14492, 14527, 15028, 15029, 15035, 15036, 15222, 15228, 15250, 15254, 15280, 15355, 15438. 
\par

\refstepcounter{casegroup}\label{grp:27}\par\addvspace{.55ex}\noindent\textbf{Case 27.} FB.\enspace \ref{abs-quartic-linear-bases}.
\enspace \(\text{\normalfont C5}\)\enspace 
204, 206. 
\enspace \(\text{\normalfont D5}\)\enspace 
893. 
\par

\refstepcounter{casegroup}\label{grp:28}\par\addvspace{.55ex}\noindent\textbf{Case 28.} NFB.\enspace \ref{case:15}.
207, 1779. 
\par

\refstepcounter{casegroup}\label{grp:29}\par\addvspace{.55ex}\noindent\textbf{Case 29.} FB.\enspace \ref{case:16}.
\enspace \(\A{13892}\)\enspace 
208, 13892. 
\enspace \(\A{13891}\)\enspace 
13891. 
\par

\refstepcounter{casegroup}\label{grp:30}\par\addvspace{.55ex}\noindent\textbf{Case 30.} FB.\enspace \ref{case:17}.
\enspace \(\A{209}\)\enspace 
209, 2212. 
\enspace \(\A{209}^{\rm op}\)\enspace 
2313, 8806. 
\par

\refstepcounter{casegroup}\label{grp:31}\par\addvspace{.55ex}\noindent\textbf{Case 31.} FB.\enspace \ref{case:6}.
\enspace \(\A{928}\)\enspace 
220, 459, 689, 715--718, 793, 814, 897, 919, 924, 928, 979, 985, 1017, 1104, 1164, 1398, 2586, 3304, 3977, 5169, 5479, 5942, 5991, 6103, 6374, 6589, 6963, 8320, 10925, 10938, 10942, 10949, 10950, 10956, 10997, 11003, 11016, 11196--11198, 11201, 11277, 11280, 11284, 11340, 11374, 11378, 11379, 12241, 12384, 12458, 12459, 12524, 12927, 12941. 
\par

\refstepcounter{casegroup}\label{grp:32}\par\addvspace{.55ex}\noindent\textbf{Case 32.} FB.\enspace \ref{case:7}.
\enspace \(\A{690}\)\enspace 
221, 222, 690, 922, 923, 11199, 11200, 11327, 11383, 12385, 12497, 13911, 15046. 
\par

\refstepcounter{casegroup}\label{grp:33}\par\addvspace{.55ex}\noindent\textbf{Case 33.} FB.\enspace \ref{case:8}.
\enspace \(\A{228}\)\enspace 
228, 824, 920, 1096, 8324, 13260. 
\par

\refstepcounter{casegroup}\label{grp:34}\par\addvspace{.55ex}\noindent\textbf{Case 34.} NFB.\enspace \ref{thm:strong}.
229, 233--237, 241, 242, 253, 257--261, 265, 266, 441, 458, 460, 461, 652, 681--686, 696, 697, 700, 703, 704, 707, 726, 730--734, 738--742, 745, 746, 749--751, 753, 755--761, 763, 765, 766, 768, 770, 773, 775--777, 779, 784--787, 798, 801--804, 808, 810--812, 816--821, 827--832, 841, 843--845, 847, 857--860, 876, 878, 901, 903, 905, 907, 944, 945, 949--956, 958, 959, 963, 965, 966, 970, 972--974, 997, 999--1001, 1007, 1029, 1031--1033, 1035, 1037--1039, 1041, 1043--1045, 1047--1050, 1052--1054, 1056--1064, 1066--1068, 1070, 1107, 1153, 1163, 1250, 1297--1299, 1303--1305, 1401, 1403--1405, 1518--1521, 1561, 1562, 1578--1583, 1608--1610, 1623--1628, 1644--1646, 1656--1658, 1662--1664, 1675--1677, 1681--1683, 1715--1717, 1831, 1865, 1875, 1948, 1950--1952, 2104, 2105, 2192--2194, 2198--2200, 2325, 2329, 2339, 2343, 2462--2465, 2488, 2492, 2496, 2497, 2502, 2503, 2505, 2507, 2508, 2510, 2514, 2516, 2520, 2521, 2528, 2531, 2538, 2540, 2544, 2575, 2577, 2591, 2593, 2595, 2597, 2599, 2601, 2603, 2605, 2607--2611, 2613, 2666, 2667, 2670, 2671, 2697, 2699, 2760--2763, 2777, 2778, 2786--2789, 2799, 2800, 2808, 2809, 2812, 2813, 2825, 2826, 2883, 3005--3015, 3157, 3158, 3186, 3187, 3226--3229, 3259--3262, 3281, 3282, 3314--3317, 3322, 3323, 3325--3349, 3351, 3382, 3383, 3395, 3396, 3416, 3417, 3437--3439, 3452, 3453, 3468, 3469, 3487, 3500, 3501, 3509, 3522, 3523, 3535, 3547, 3548, 3561, 3562, 3572, 3624, 3625, 3627, 3629, 3661, 3662, 3664, 3666, 3710, 3711, 3713, 3715, 3745, 3746, 3750, 3752, 3753, 3755, 3756, 3759--3761, 3765, 3767, 3768, 3772, 3774, 3776, 3777, 3781, 3783, 3784, 3796, 3799, 3807, 3810, 3839, 3842, 3848, 3851, 3853, 3858, 3861, 3879, 3881, 3885, 3926, 3929, 3933--3936, 3939, 3941, 3942, 3960--3962, 3997, 4000, 4034, 4037, 4116, 4120, 4121, 4278--4280, 4368, 4369, 4417--4421, 4430, 4435--4437, 4446, 4447, 4451--4455, 4469, 4470, 4478, 4628, 4629, 4643, 4644, 4695, 4696, 4714--4719, 4733, 4741, 4763, 4764, 4776, 4785, 4786, 4788, 4789, 4838--4840, 4858--4865, 4875--4882, 4889, 4890, 4895, 4896, 4918, 4925--4930, 4935, 4936, 4943, 5023--5026, 5039, 5042, 5045, 5053, 5054, 5057--5061, 5063, 5064, 5074, 5076, 5077, 5079, 5080, 5085--5092, 5094, 5095, 5185--5194, 5234--5237, 5352, 5370, 5392, 5405--5407, 5417, 5423, 5430, 5439, 5440, 5446, 5450--5452, 5454, 5456--5461, 5468--5471, 5478, 5483--5486, 5496--5499, 5505, 5506, 5514, 5515, 5523, 5527, 5536, 5547, 5589, 5615--5619, 5627--5634, 5637, 5640--5644, 5646--5648, 5651--5656, 5670--5674, 5683, 5685, 5688, 5690, 5693, 5708--5711, 5715--5719, 5721--5724, 5730, 5786, 5838, 5877, 5893, 5902, 5944--5946, 5972, 5973, 5990, 6009, 6010, 6024, 6073--6075, 6104, 6105, 6162, 6172, 6176, 6208, 6210, 6254, 6256, 6316, 6318, 6319, 6377--6382, 6458, 6460, 6514--6519, 6590--6593, 6698, 6719, 6720, 6730, 6738, 6746, 6789, 6809, 6811, 6819, 6855, 6858, 6892, 6894, 6947, 6952, 6953, 7161, 7162, 7196, 7197, 7360, 7361, 7376, 7377, 7402, 7410--7415, 7420, 7430, 7438, 7439, 7447, 7453, 7456--7459, 7468, 7475, 7514, 7526, 7527, 7551, 7560--7565, 7567, 7573, 7581, 7595, 7600, 7604, 7610, 7630, 7636--7640, 7654--7662, 7667, 7690, 7761--7764, 7795--7799, 7893--7896, 7927--7940, 8038, 8045, 8065, 8069, 8077, 8094, 8098, 8105, 8106, 8111, 8123, 8124, 8128, 8166, 8190, 8198, 8206, 8227, 8239, 8252, 8259, 8292, 8299, 8312, 8319, 8334--8345, 8591, 8597, 8614, 8625, 8643, 8832, 8840, 8851, 8858, 8888, 8891, 8893, 8949--8952, 8990, 8997, 9132, 9144--9146, 9148, 9150--9152, 9154, 9227, 9244, 9318, 9368, 9370, 9499, 9500, 9522, 9523, 9531, 9556, 9559, 9701, 9702, 9767, 9779, 9780, 9803, 9811--9816, 9821, 9825, 9874--9877, 9951, 9964, 9969, 9976, 9977, 9981, 10010, 10069, 10480, 10518--10531, 10838, 10846--10848, 10857, 10858, 10860--10863, 10866--10868, 10871--10875, 10877, 10878, 10884, 10895, 10897, 10915, 10935, 10948, 10971--10973, 10975--10977, 10979, 10980, 10982, 10983, 10985--10989, 10991--10993, 11023, 11025, 11028, 11030, 11032, 11035, 11037, 11048, 11052, 11074, 11076, 11080--11084, 11092--11099, 11101, 11111, 11113, 11116, 11119, 11121, 11124, 11128, 11132, 11155, 11159, 11181, 11185, 11202, 11212--11214, 11221--11227, 11242--11251, 11253--11263, 11267, 11268, 11272, 11294, 11295, 11309, 11330--11334, 11341, 11343, 11344, 11346, 11347, 11349, 11350, 11352--11355, 11357, 11360, 11362, 11363, 11397, 11399, 11400, 11402, 11406, 11408, 11409, 11411, 11430--11433, 11437, 11439, 11440, 11442, 11445--11476, 11478--11482, 11485--11491, 11503--11507, 11509--11511, 11514, 11515, 11517--11519, 11522, 11536--11539, 11541--11544, 11546--11549, 11571, 11572, 11597--11599, 11601--11609, 11614--11618, 11630, 11631, 11635, 11641, 11646, 11653--11655, 11659, 11660, 11674, 11749, 11781, 11782, 11808, 11824, 11875, 12031, 12071, 12103, 12104, 12119, 12187--12190, 12240, 12244, 12311, 12330, 12336, 12355, 12356, 12367, 12368, 12390, 12392, 12396, 12406, 12411, 12412, 12415, 12417, 12419, 12421, 12423--12426, 12432, 12440, 12444, 12465--12467, 12470, 12471, 12474, 12475, 12477, 12478, 12484--12487, 12489--12492, 12494, 12503--12505, 12507--12510, 12512--12515, 12517--12520, 12522, 12539, 12541--12543, 12545, 12548, 12552, 12578, 12579, 12593, 12594, 12598, 12599, 12603, 12620, 12627, 12629--12632, 12639, 12640, 12643, 12644, 12648, 12661, 12664, 12680, 12683, 12684, 12686, 12687, 12690, 12691, 12694, 12696, 12699, 12701, 12706, 12709, 12710, 12713, 12714, 12719, 12721, 12724, 12726, 12731, 12737, 12738, 12760, 12761, 12764, 12765, 12768, 12770, 12775, 12777, 12782, 12786--12789, 12791, 12792, 12795--12797, 12801, 12803, 12804, 12808, 12810, 12811, 12815, 12819, 12820, 12822, 12824, 12825, 12828, 12830, 12834, 12836, 12840, 12845, 12846, 12848, 12850, 12851, 12855, 12857, 12861, 12863, 12867, 12868, 12871, 12872, 12911, 12912, 12918, 12930--12932, 12934--12937, 12939, 12962, 12963, 12965, 12967, 12968, 12972, 12974, 12977, 12979, 12983, 12987, 12988, 12990, 12992, 12993, 12996, 12998, 13001, 13003, 13007, 13011, 13012, 13035, 13036, 13039--13042, 13046, 13047, 13149, 13150, 13157, 13208, 13209, 13213, 13255--13258, 13262, 13263, 13395, 13449, 13460, 13461, 13500, 13501, 13505, 13506, 13509, 13513, 13514, 13516, 13525, 13526, 13530, 13531, 13534, 13538, 13539, 13543, 13544, 13548, 13549, 13555--13561, 13568, 13570, 13578, 13583, 13590, 13598, 13599, 13602, 13604, 13606, 13608, 13610, 13612, 13614, 13620, 13621, 13625--13628, 13630, 13631, 13634--13637, 13639, 13642, 13643, 13645, 13646, 13649, 13651, 13653, 13655, 13656, 13658, 13661, 13668--13681, 13686, 13688, 13690, 13692, 13694, 13696, 13720, 13722, 13725, 13726, 13729, 13730, 13732, 13738, 13740, 13742, 13744, 13746, 13748, 13751--13757, 13772, 13775, 13777, 13779, 13782, 13904, 13905, 13907, 13908, 13970, 13996, 14019, 14028, 14029, 14034, 14043, 14046, 14079, 14086, 14094--14097, 14102, 14108--14110, 14136--14139, 14165, 14166, 14178, 14205, 14234, 14271, 14354, 14417, 14418, 14451, 14496, 14514, 14537, 14538, 14547, 14556, 14557, 14561, 14563, 14566, 14587, 14590, 14628, 14635, 14639--14641, 14643, 14669--14673, 14708, 14721--14723, 14750, 14782, 14830, 14858, 14859, 14904, 15030, 15037, 15051, 15052, 15054, 15056, 15057, 15061, 15063, 15066, 15068, 15072, 15077, 15078, 15080, 15082, 15083, 15086, 15088, 15091, 15093, 15095, 15099, 15102, 15103, 15128, 15138, 15209, 15211, 15213, 15215, 15219, 15221, 15223, 15225, 15227, 15229, 15231, 15235, 15237--15239, 15241, 15242, 15251, 15253, 15255, 15257, 15290, 15291, 15325, 15326, 15358, 15359, 15445, 15446, 15558, 15559, 15631, 15634, 15636, 15638, 15641. 
\par

\refstepcounter{casegroup}\label{grp:35}\par\addvspace{.55ex}\noindent\textbf{Case 35.} FB.\enspace Proposition~\ref{prop:lower-fibres}, source \(\A{971}\), case~\ref{grp:2}.
\enspace \(\A{971}\)\enspace 
230, 694, 948, 1106, 10017, 10879, 11107, 11110, 12357. 
\par

\refstepcounter{casegroup}\label{grp:36}\par\addvspace{.55ex}\noindent\textbf{Case 36.} FB.\enspace Proposition~\ref{prop:lower-fibres}, source \(\A{975}\), case~\ref{grp:26}.
\enspace \(\A{975}\)\enspace 
238, 659, 698, 735, 747, 877, 957, 962, 964, 969, 1108, 1257, 1832, 2117, 2862, 3773, 3809, 5682, 5845, 10834, 10841, 10852, 11040, 11056, 11066, 11069, 11070, 11073, 11417, 11423, 11508, 11552--11555, 12349, 12404, 12590, 12802, 12829. 
\par

\refstepcounter{casegroup}\label{grp:37}\par\addvspace{.55ex}\noindent\textbf{Case 37.} FB.\enspace \ref{case:18}.
\enspace \(\A{252}\)\enspace 
246, 252, 463, 823, 1666, 7904. 
\par

\refstepcounter{casegroup}\label{grp:38}\par\addvspace{.55ex}\noindent\textbf{Case 38.} FB.\enspace \ref{case:19}.
\enspace \(\A{462}\)\enspace 
247, 462, 11187, 11271. 
\enspace \(\A{462}^{\rm op}\)\enspace 
679, 822, 13162, 14126. 
\par

\refstepcounter{casegroup}\label{grp:39}\par\addvspace{.55ex}\noindent\textbf{Case 39.} FB.\enspace \ref{case:20}.
248, 11186, 11366. 
\par

\refstepcounter{casegroup}\label{grp:40}\par\addvspace{.55ex}\noindent\textbf{Case 40.} FB.\enspace \ref{thm:treelength}.
254, 701, 727, 743, 10889, 10896, 11103, 11106, 11115, 11135, 11138, 11141, 11359, 11368, 11390, 11394, 11396, 11413, 11424, 11427, 11443, 11493, 11498, 11501, 11528, 11535, 11540, 11545, 11550, 11551, 11558, 11561, 11570, 11582, 11583, 11649, 11652, 11801, 12110, 12245, 12369, 12373, 12394, 12413, 12442, 12536, 12538, 12547, 12555, 12557, 12560, 12577, 12580, 12592, 12617, 12619, 12634, 12651, 12653, 12656, 12693, 12698, 12708, 12718, 12723, 12727, 12759, 12767, 12774, 12778, 12817, 12833, 12843, 12844, 12854, 12860, 12961, 12971, 12985, 12986, 12995, 13000, 13415, 13480, 13482, 13493, 13507, 13510, 13532, 13566, 13588, 13596, 13597, 13617, 13647, 13659, 13773, 14104, 14107, 15049, 15050, 15060, 15065, 15075, 15076, 15085, 15090, 15632. 
\par

\refstepcounter{casegroup}\label{grp:41}\par\addvspace{.55ex}\noindent\textbf{Case 41.} FB.\enspace \ref{case:21}.
262, 705, 2489, 2498, 11039, 11087, 11264, 11348, 11405, 11414, 11497, 11499, 11526, 11529, 12363, 12409, 12596, 12597, 12606, 12610, 12818, 12827, 12960, 12976, 13489, 13581. 
\par

\refstepcounter{casegroup}\label{grp:42}\par\addvspace{.55ex}\noindent\textbf{Case 42.} FB.\enspace \ref{case:22}.
\enspace \(\A{279}\)\enspace 
279, 1780. 
\par

\refstepcounter{casegroup}\label{grp:43}\par\addvspace{.55ex}\noindent\textbf{Case 43.} FB.\enspace \ref{case:23}.
\enspace \(\A{280}\)\enspace 
280, 13838. 
\par

\refstepcounter{casegroup}\label{grp:44}\par\addvspace{.55ex}\noindent\textbf{Case 44.} FB.\enspace \ref{thm:qgraph}.
281, 334, 338, 1264, 1710, 2213, 2305, 2306, 2347, 2349, 2363, 2365, 2551, 2823, 2824, 5520, 5529, 5530, 5740, 5802, 5852, 8180, 8807, 8826, 8828, 10542, 11701, 11761, 12228, 12281, 13233, 13463, 13874, 13976, 14190, 14316, 14325, 14582, 14597, 14605, 14620, 14653, 14678, 14689, 14691, 14704, 14706, 14710, 14711, 14730, 14735, 14740, 14741, 14743--14747, 14752--14754, 14756--14758, 14760--14762, 14767, 14769, 14774--14776, 14778--14780, 14783, 14784, 14792--14794, 14796--14798, 14800, 14801, 14893, 14983, 14984, 15119, 15120, 15455, 15456, 15567--15569. 
\par

\refstepcounter{casegroup}\label{grp:45}\par\addvspace{.55ex}\noindent\textbf{Case 45.} FB.\enspace \ref{case:24}.
\enspace \(\A{282}\)\enspace 
282, 14398. 
\enspace \(\A{1707}\)\enspace 
1707, 14388. 
\enspace \(\A{282}^{\rm op}\)\enspace 
2348, 13839. 
\enspace \(\A{1707}^{\rm op}\)\enspace 
2821, 13198. 
\par

\refstepcounter{casegroup}\label{grp:46}\par\addvspace{.55ex}\noindent\textbf{Case 46.} FB.\enspace \ref{case:25}.
\enspace \(\A{284}\)\enspace 
284, 14436, 14686. 
\par

\refstepcounter{casegroup}\label{grp:47}\par\addvspace{.55ex}\noindent\textbf{Case 47.} FB.\enspace \ref{thm:twofactor}.
\enspace \(S_{4,42}\times M\)\enspace 
285, 289, 1162, 1371, 1699, 3175, 3183, 8318. 
\enspace \(S_{4,705}\times M\)\enspace 
290, 1086, 1167, 1372, 3176, 6097, 6575, 8327, 13193. 
\enspace \(S_{4,765}\times T\)\enspace 
630, 2139, 2147, 2148. 
\enspace \(S_{4,413}\times M\)\enspace 
2311, 12207. 
\enspace \(S_{4,430}\times M\)\enspace 
3184, 10109, 10574. 
\enspace \(S_{4,627}\times M\)\enspace 
3185, 13061, 13192. 
\enspace \(S_{4,756}\times T\)\enspace 
4520, 6667, 7052, 7120. 
\enspace \(S_{4,766}\times T\)\enspace 
4523, 6672, 7058, 7127. 
\enspace \(S_{4,780}\times T\)\enspace 
4532, 6679, 7072, 7132. 
\enspace \(S_{4,681}\times T\)\enspace 
4704, 6822, 7199, 7913. 
\enspace \(S_{4,690}\times T\)\enspace 
4720, 6832, 7208, 7944. 
\enspace \(S_{4,686}\times T\)\enspace 
4727, 6836, 7210, 7955. 
\enspace \(S_{4,721}\times T\)\enspace 
5246, 8373, 8432, 8433. 
\enspace \(S_{4,763}\times T\)\enspace 
5255, 8387, 8443, 8448. 
\enspace \(S_{4,757}\times T\)\enspace 
5264, 8402, 8447, 8451. 
\enspace \(S_{4,682}\times T\)\enspace 
5339, 8498, 8500, 8671. 
\enspace \(S_{4,589}\times T\)\enspace 
5346, 8506, 8508, 8687. 
\enspace \(S_{4,456}\times T\)\enspace 
6129, 6276, 6364, 8873. 
\enspace \(S_{4,443}\times T\)\enspace 
6158, 6727, 7785, 9011. 
\enspace \(S_{4,288}\times M\)\enspace 
6574. 
\enspace \(S_{4,769}\times T\)\enspace 
8907, 8941, 8975, 10212. 
\enspace \(S_{4,775}\times T\)\enspace 
8912, 8943, 8976, 10215. 
\enspace \(S_{4,799}\times T\)\enspace 
9107, 9257, 9270, 10354. 
\enspace \(S_{4,748}\times T\)\enspace 
9402, 9586, 9641, 10395. 
\enspace \(S_{4,751}\times T\)\enspace 
9404, 9590, 9642, 10398. 
\enspace \(S_{4,793}\times T\)\enspace 
9421, 9606, 9652, 10433. 
\enspace \(S_{4,734}\times T\)\enspace 
10132, 10200, 10204, 10605. 
\enspace \(S_{4,772}\times T\)\enspace 
10143, 10217, 10222, 10615. 
\par

\refstepcounter{casegroup}\label{grp:48}\par\addvspace{.55ex}\noindent\textbf{Case 48.} FB.\enspace \ref{case:26}.
\enspace \(\A{13841}\)\enspace 
286, 12272, 13841, 14416, 14687. 
\par

\refstepcounter{casegroup}\label{grp:49}\par\addvspace{.55ex}\noindent\textbf{Case 49.} FB.\enspace \ref{case:27}.
\enspace \(\A{12280}^{\rm op}\)\enspace 
288, 14324, 14412, 14421, 14688. 
\enspace \(\A{12280}\)\enspace 
2352, 12280, 13842, 14413, 14700. 
\par

\refstepcounter{casegroup}\label{grp:50}\par\addvspace{.55ex}\noindent\textbf{Case 50.} FB.\enspace \ref{case:28}.
\enspace \(\A{299}\)\enspace 
299, 1874, 1912. 
\enspace \(\A{1809}\)\enspace 
1809, 1879, 1913. 
\enspace \(\A{299}^{\rm op}\)\enspace 
2938, 5901, 6301. 
\enspace \(\A{1809}^{\rm op}\)\enspace 
5868, 5906, 6302. 
\par

\refstepcounter{casegroup}\label{grp:51}\par\addvspace{.55ex}\noindent\textbf{Case 51.} FB.\enspace \ref{case:29}.
\enspace \(\A{319}\)\enspace 
319, 10532. 
\par

\refstepcounter{casegroup}\label{grp:52}\par\addvspace{.55ex}\noindent\textbf{Case 52.} FB.\enspace \ref{case:25}.
\enspace \(\A{332}\)\enspace 
332, 1708, 10024, 14834. 
\par

\refstepcounter{casegroup}\label{grp:53}\par\addvspace{.55ex}\noindent\textbf{Case 53.} FB.\enspace \ref{abs-connected-bases}.
\enspace \(\text{\normalfont G5}\)\enspace 
336, 12273, 13873, 14690, 14749. 
\enspace \(\text{\normalfont F5}\)\enspace 
13872. 
\par

\refstepcounter{casegroup}\label{grp:54}\par\addvspace{.55ex}\noindent\textbf{Case 54.} FB.\enspace \ref{case:30}.
\enspace \(\A{1960}\)\enspace 
389, 390, 427, 428, 594, 1193, 1347, 1348, 1542, 1543, 1812, 1881, 1919, 1960, 5510, 6975, 6976, 8321, 13860, 14619, 14662. 
\par

\refstepcounter{casegroup}\label{grp:55}\par\addvspace{.55ex}\noindent\textbf{Case 55.} FB.\enspace \ref{case:31}.
\enspace \(\A{12638}^{\rm op}\)\enspace 
399, 13849, 14494, 14536, 14616. 
\enspace \(\A{12638}\)\enspace 
794, 12343, 12481, 12638, 13450. 
\par

\refstepcounter{casegroup}\label{grp:56}\par\addvspace{.55ex}\noindent\textbf{Case 56.} FB.\enspace \ref{case:30}.
\enspace \(\A{6230}\)\enspace 
412, 413, 1508, 1509, 2402, 2403, 3972, 3973, 4127, 4201, 5870, 5908, 6230, 6306, 7855, 7856, 8101, 8119, 13865, 13871, 14632. 
\par

\refstepcounter{casegroup}\label{grp:57}\par\addvspace{.55ex}\noindent\textbf{Case 57.} FB.\enspace \ref{cor:periodic-register}.
417--419, 429, 432, 435, 439, 648, 650, 1141, 1147, 1306, 1428, 1434, 1470, 1472, 1474, 1479, 1483, 1488, 1489, 1494, 1497, 1523, 1526, 1530, 1547, 1551, 1556, 1563, 1567, 1584--1586, 1596--1599, 1602, 1617, 1620, 1635, 1661, 1684--1686, 1691--1693, 1727--1730, 1733, 1773--1775, 1788, 1791--1793, 1818, 1859, 1890, 1902, 1995, 2001, 2159, 2197, 2399, 2400, 2407, 2739, 2764, 2765, 2772, 2814, 2815, 2818, 2827, 2828, 2843, 2852, 2853, 2932, 2950, 2953, 2957, 2961, 2964, 2965, 2981, 2984, 2987, 3023, 3026, 3027, 3034, 3035, 3040, 3042, 3044, 3046, 3048, 3050, 3051, 3054, 3055, 3060, 3062, 3067, 3069, 3072, 3074, 3075, 3080, 3082, 3085, 3088--3093, 3095, 3097, 3099, 3100, 3105, 3106, 3113--3115, 3128--3130, 3135, 3142, 3147, 3196, 3199, 3202, 3204, 3205, 3208, 3209, 3214, 3215, 3218, 3219, 3223, 3267, 3272, 3275, 3278, 3400, 3402, 3404, 3405, 3407, 3410, 3412, 3414, 3513--3517, 3519--3521, 3526, 3529, 3531, 3533, 3567, 3569--3571, 3965, 3970, 3988, 4107, 4110, 4111, 4131, 4135, 4140, 4148, 4158, 4161, 4164, 4167, 4174, 4179, 4181, 4189, 4218, 4219, 4227, 4238, 4246--4248, 4251, 4253--4255, 4257, 4266, 4284, 4287--4290, 4294, 4324, 4326, 4327, 4329, 4330, 4332, 4333, 4335, 4338, 4341, 4344, 4347, 4350, 4353, 4354, 4359, 4378, 4379, 4382, 4384--4389, 4405, 4406, 4408, 4458, 4459, 4464, 4465, 4490, 4493, 4502, 4511, 4539, 4540, 4543, 4545, 4548, 4553, 4557, 4558, 4582, 4585, 4611, 4614, 4619, 4622, 4625, 4636, 4637, 4650, 4653, 4657, 4660, 4664, 4670, 4676, 4678, 4681, 4686, 4689, 4697, 4736, 4969--4972, 4975, 4976, 4985, 4987, 4991, 4992, 4997, 4999, 5000, 5009, 5011, 5012, 5112, 5118, 5122, 5127, 5130, 5148, 5151, 5156, 5159, 5174, 5177, 5184, 5216, 5221, 5229, 5233, 5282, 5292, 5293, 5309, 5311--5313, 5325--5327, 5353, 5354, 5357, 5361, 5364, 5381, 5382, 5393--5397, 5399--5402, 5404, 5418, 5420--5422, 5426, 5429, 5431, 5432, 5437, 5534, 5585, 5587, 5913, 5918, 5926, 5933, 5938, 5939, 5955, 5982, 5986, 6004, 6019--6021, 6032, 6033, 6058--6060, 6076--6078, 6091--6094, 6106, 6116, 6141, 6157, 6159, 6161, 6163, 6164, 6179, 6190, 6192, 6194, 6234, 6240, 6246--6248, 6250, 6251, 6270, 6292, 6332, 6343, 6344, 6355, 6369, 6394, 6404, 6410, 6415, 6420, 6437, 6478--6483, 6490, 6492, 6494, 6502, 6504, 6506, 6526--6528, 6547--6552, 6565--6569, 6576, 6597, 6610, 6622, 6624, 6628, 6630, 6631, 6640, 6646, 6651, 6655, 6660, 6687, 6689, 6692, 6694, 6696, 6726, 6728, 6731, 6750, 6752, 6754, 6756, 6758, 6773, 6777, 6778, 6780, 6783, 6785, 6790, 6874, 6876, 6877, 6880--6882, 6884, 6885, 6887, 6889, 6890, 6899, 6906, 6912, 6919, 6929, 6934, 6938, 6940, 6941, 6962, 6996, 7001, 7018, 7036, 7077, 7097--7099, 7101, 7105, 7111, 7147--7149, 7156, 7157, 7164, 7165, 7168--7170, 7177, 7181--7184, 7191, 7230, 7231, 7234, 7236, 7254, 7734, 7744, 7749, 7756, 7760, 7779, 7782, 7789, 7804, 7820, 7826, 7833, 7841, 7846, 7853, 7860, 7870, 7876, 7881, 7888, 7898, 7987, 7988, 7992, 8006, 8039, 8040, 8042--8044, 8048, 8049, 8055--8058, 8070--8072, 8074, 8075, 8078, 8194, 8196, 8202, 8203, 8214, 8215, 8223, 8225, 8234--8236, 8247, 8248, 8256, 8257, 8270, 8271, 8296, 8297, 8308, 8309, 8382, 8420, 8425--8427, 8430, 8477, 8478, 8482, 8488, 8513--8526, 8529, 8536, 8537, 8558, 8562, 8585, 8586, 8588--8590, 8621, 8623, 8626, 8629--8634, 8636--8638, 8640, 8641, 8644, 8692, 8695, 8697--8699, 8702, 8703, 8706, 8710, 8713, 8778, 8801--8804, 9213--9215, 9223--9226, 9269, 9317, 9383, 9384, 9387, 9445, 9453, 9454, 9460, 9462--9464, 9467--9469, 9480, 9481, 9483, 9486, 9491, 9494, 9508--9514, 9516--9519, 9521, 9526, 9528--9530, 9618, 9630, 9635--9637, 9640, 9668, 9673, 9680, 9686--9696, 9698, 9730, 9731, 9736, 9741, 9761, 9762, 9764--9766, 9772--9774, 9776, 9778, 9781, 9782, 9786--9790, 9792, 9794, 9796--9798, 9800, 9801, 9804, 9840--9844, 9847, 9848, 9851--9853, 9857, 9859, 9862, 9866, 9869, 9937, 9938, 9953--9961, 9980, 10009, 10112--10121, 10459, 10460, 10468--10472, 10474--10477, 10580, 10591, 10597--10601, 10819, 15553. 
\par

\refstepcounter{casegroup}\label{grp:58}\par\addvspace{.55ex}\noindent\textbf{Case 58.} FB.\enspace \ref{abs-square-bases}.
\enspace \(\text{\normalfont A\_EH}\)\enspace 
421, 422, 451, 452, 1531, 1532, 1738, 1739, 4136, 4222, 5930, 6258, 6352, 8715, 8716, 9081, 9083, 9173, 9249, 9871, 9872, 9970, 13882, 14684, 14829. 
\enspace \(\text{\normalfont A\_DE}\)\enspace 
1465, 3018, 4223, 5368, 5369, 5947, 6257, 9172, 9248, 9496, 9498, 9971, 11792, 12155, 12178. 
\enspace \(\text{\normalfont A\_IDE}\)\enspace 
1467, 2738, 3019, 3153--3156, 3359, 4098, 4226, 4358, 4364, 4669, 4685, 4693, 4694, 5017, 5022, 5434, 5491, 5949, 5994, 6262, 6788, 8118, 11773, 11783, 11785, 11786, 11790, 11841, 11869, 11876, 11877, 11879, 11880, 11915, 11985, 12074, 12075, 12084, 12085, 13072, 13073, 13242, 13281. 
\enspace \(\text{\normalfont A\_D}\)\enspace 
2958, 3045, 4139, 5358, 5359, 8700, 8701, 9485, 9487, 9856. 
\enspace \(\text{\normalfont A\_empty}\)\enspace 
2959, 4138, 9855, 13329, 13881. 
\enspace \(\text{\normalfont A\_I}\)\enspace 
2963, 3211, 3225, 4146, 5149, 7874, 13063, 13064, 13111, 13114, 13201, 13202, 13220, 13221, 13223--13225, 13851, 13855, 13868, 13875, 13939, 13940. 
\enspace \(\text{\normalfont A\_IH}\)\enspace 
4339, 5147, 5182, 5183, 6065, 6493, 6905, 8244, 14555. 
\enspace \(\text{\normalfont A\_IE}\)\enspace 
5020, 5021, 11774, 11878, 12036, 12086, 14546. 
\par

\refstepcounter{casegroup}\label{grp:59}\par\addvspace{.55ex}\noindent\textbf{Case 59.} FB.\enspace \ref{case:32}.
\enspace \(\A{1546}\)\enspace 
430, 1546, 8704, 8705, 9861. 
\enspace \(\A{1545}\)\enspace 
431, 1545, 9860, 14681, 14826. 
\enspace \(\A{1546}^{\rm op}\)\enspace 
2951, 4130, 8711, 8712, 9868. 
\enspace \(\A{1545}^{\rm op}\)\enspace 
2952, 4129, 9867, 13331, 13883. 
\par

\refstepcounter{casegroup}\label{grp:60}\par\addvspace{.55ex}\noindent\textbf{Case 60.} FB.\enspace \ref{case:33}.
\enspace \(\A{1549}\)\enspace 
434, 437, 1549, 1554, 2409, 7879, 13325, 14540, 14549, 14610--14612, 14631, 14646--14650, 14712, 15547, 15548. 
\par

\refstepcounter{casegroup}\label{grp:61}\par\addvspace{.55ex}\noindent\textbf{Case 61.} FB.\enspace \ref{case:34}.
\enspace \(\A{1020}^{\rm op}\)\enspace 
492, 526, 532, 546, 577, 1866, 1925, 4900, 6821, 7903. 
\par

\refstepcounter{casegroup}\label{grp:62}\par\addvspace{.55ex}\noindent\textbf{Case 62.} FB.\enspace \ref{case:34}.
\enspace \(\A{5903}^{\rm op}\)\enspace 
527, 547, 1816, 1878, 1926, 13807, 14069, 14125. 
\par

\refstepcounter{casegroup}\label{grp:63}\par\addvspace{.55ex}\noindent\textbf{Case 63.} FB.\enspace \ref{case:35}.
\enspace \(\A{1941}^{\rm op}\)\enspace 
536, 554, 2435, 2970, 3938, 4086, 4207, 4272, 5566, 5894, 6224, 6310, 7382, 7520, 7579. 
\par

\refstepcounter{casegroup}\label{grp:64}\par\addvspace{.55ex}\noindent\textbf{Case 64.} FB.\enspace \ref{case:7}.
\enspace \(\A{865}^{\rm op}\)\enspace 
537, 538, 2436, 3937, 4208, 5873, 5905, 6311, 7435, 7519, 13810, 13814, 14083. 
\par

\refstepcounter{casegroup}\label{grp:65}\par\addvspace{.55ex}\noindent\textbf{Case 65.} FB.\enspace \ref{case:36}.
\enspace \(\A{1028}^{\rm op}\)\enspace 
542, 558, 562, 2972, 4275, 6358, 8646, 9101, 9105, 9228, 9784, 9806. 
\par

\refstepcounter{casegroup}\label{grp:66}\par\addvspace{.55ex}\noindent\textbf{Case 66.} FB.\enspace \ref{case:36}.
\enspace \(\A{1025}^{\rm op}\)\enspace 
543, 559, 5935, 6359, 9102, 9246, 13825, 14161, 14203. 
\par

\refstepcounter{casegroup}\label{grp:67}\par\addvspace{.55ex}\noindent\textbf{Case 67.} FB.\enspace \ref{case:37}.
\enspace \(\A{8628}\)\enspace 
548, 8628, 9785. 
\enspace \(\A{8628}^{\rm op}\)\enspace 
2992, 8707, 9863. 
\par

\refstepcounter{casegroup}\label{grp:68}\par\addvspace{.55ex}\noindent\textbf{Case 68.} FB.\enspace \ref{case:38}.
\enspace \(\A{11021}^{\rm op}\)\enspace 
549, 14162, 14204. 
\enspace \(\A{11021}\)\enspace 
2993, 11021, 11299. 
\par

\refstepcounter{casegroup}\label{grp:69}\par\addvspace{.55ex}\noindent\textbf{Case 69.} FB.\enspace \ref{case:39}.
551, 553, 2440, 13315, 14017, 14026, 14070--14072, 14085, 14091, 14092, 14103, 14131, 14132, 14135, 15286, 15322, 15354. 
\par

\refstepcounter{casegroup}\label{grp:70}\par\addvspace{.55ex}\noindent\textbf{Case 70.} FB.\enspace \ref{case:40}.
\enspace \(\A{1959}\)\enspace 
570, 595, 1099, 1391, 1834, 1880, 1959, 3271, 4366, 6585, 8322, 14252. 
\par

\refstepcounter{casegroup}\label{grp:71}\par\addvspace{.55ex}\noindent\textbf{Case 71.} FB.\enspace Proposition~\ref{prop:lower-fibres}, source \(\A{2144}\), case~\ref{grp:2}.
\enspace \(\A{2144}\)\enspace 
627, 2135, 2143. 
\par

\refstepcounter{casegroup}\label{grp:72}\par\addvspace{.55ex}\noindent\textbf{Case 72.} FB.\enspace \ref{case:41}.
\enspace \(\A{2137}\)\enspace 
631, 632, 2137. 
\enspace \(\A{2137}^{\rm op}\)\enspace 
5247, 5248, 8371. 
\par

\refstepcounter{casegroup}\label{grp:73}\par\addvspace{.55ex}\noindent\textbf{Case 73.} FB.\enspace \ref{case:40}.
\enspace \(\A{2155}\)\enspace 
637, 638, 2155, 2156, 4233, 6278, 6366, 9115, 9117, 9267. 
\par

\refstepcounter{casegroup}\label{grp:74}\par\addvspace{.55ex}\noindent\textbf{Case 74.} FB.\enspace \ref{thm:correct643}.
643, 912, 1893, 2094, 2095, 12951, 12957, 13700, 13703, 13705--13707, 13769, 15249. 
\par

\refstepcounter{casegroup}\label{grp:75}\par\addvspace{.55ex}\noindent\textbf{Case 75.} FB.\enspace \ref{case:9}.
676, 825, 1175, 8323, 13217, 14663. 
\par

\refstepcounter{casegroup}\label{grp:76}\par\addvspace{.55ex}\noindent\textbf{Case 76.} FB.\enspace \ref{case:42}.
833, 9881, 10005, 13019. 
\par

\refstepcounter{casegroup}\label{grp:77}\par\addvspace{.55ex}\noindent\textbf{Case 77.} NFB.\enspace \ref{case:43}, \ref{case:44}.
\enspace \(\A{836}\)\enspace 
834, 836, 839, 848, 929, 1665, 7902, 8285, 11370, 11669, 12305. 
\par

\refstepcounter{casegroup}\label{grp:78}\par\addvspace{.55ex}\noindent\textbf{Case 78.} FB.\enspace \ref{case:45}.
\enspace \(\A{837}\)\enspace 
835, 837, 1236, 2201, 5241, 5824. 
\enspace \(\A{837}^{\rm op}\)\enspace 
2534, 3728, 5661, 7768, 7965, 8156. 
\par

\refstepcounter{casegroup}\label{grp:79}\par\addvspace{.55ex}\noindent\textbf{Case 79.} FB.\enspace \ref{case:46}.
838, 11683, 12306, 13018. 
\par

\refstepcounter{casegroup}\label{grp:80}\par\addvspace{.55ex}\noindent\textbf{Case 80.} FB.\enspace \ref{case:47}.
\enspace \(\A{840}\)\enspace 
840, 11593, 11667, 11676, 12303. 
\enspace \(\A{840}^{\rm op}\)\enspace 
2536, 11371, 11668, 12238, 12304. 
\par

\refstepcounter{casegroup}\label{grp:81}\par\addvspace{.55ex}\noindent\textbf{Case 81.} FB.\enspace \ref{case:48}.
\enspace \(\A{11700}\)\enspace 
850, 11700, 12325. 
\enspace \(\A{11700}^{\rm op}\)\enspace 
2969, 12116, 13318. 
\par

\refstepcounter{casegroup}\label{grp:82}\par\addvspace{.55ex}\noindent\textbf{Case 82.} FB.\enspace \ref{case:35}.
\enspace \(\A{1941}\)\enspace 
851, 855, 864, 991, 1006, 1018, 1075, 1203, 1540, 1867, 1941, 1944, 5209, 6973, 8278. 
\par

\refstepcounter{casegroup}\label{grp:83}\par\addvspace{.55ex}\noindent\textbf{Case 83.} FB.\enspace \ref{case:49}.
\enspace \(\A{12053}^{\rm op}\)\enspace 
852, 12239, 12324. 
\enspace \(\A{12053}\)\enspace 
2971, 12053, 12102. 
\par

\refstepcounter{casegroup}\label{grp:84}\par\addvspace{.55ex}\noindent\textbf{Case 84.} FB.\enspace \ref{case:50}.
\enspace \(\A{856}\)\enspace 
856, 1183, 1876, 8282. 
\enspace \(\A{856}^{\rm op}\)\enspace 
3066, 4087, 5904, 7580. 
\par

\refstepcounter{casegroup}\label{grp:85}\par\addvspace{.55ex}\noindent\textbf{Case 85.} FB.\enspace \ref{case:7}.
\enspace \(\A{865}\)\enspace 
865, 992, 993, 1019, 1539, 1828, 1877, 1945, 5503, 6972, 10951, 11008, 11195. 
\par

\refstepcounter{casegroup}\label{grp:86}\par\addvspace{.55ex}\noindent\textbf{Case 86.} FB.\enspace \ref{case:51}.
\enspace \(\A{866}\)\enspace 
866, 12318. 
\enspace \(\A{866}^{\rm op}\)\enspace 
4209, 11806. 
\par

\refstepcounter{casegroup}\label{grp:87}\par\addvspace{.55ex}\noindent\textbf{Case 87.} FB.\enspace \ref{case:6}.
\enspace \(\A{872}\)\enspace 
872, 1184, 1230, 1839, 2118, 2652, 5573, 5755, 5765, 5819, 12582, 12585, 12586, 12588, 13400, 13403, 13404, 13407, 13410, 13413, 13414, 13416, 13486, 13575, 15199. 
\par

\refstepcounter{casegroup}\label{grp:88}\par\addvspace{.55ex}\noindent\textbf{Case 88.} NFB.\enspace \ref{case:52}.
894, 13236, 13237. 
\par

\refstepcounter{casegroup}\label{grp:89}\par\addvspace{.55ex}\noindent\textbf{Case 89.} FB.\enspace \ref{case:53}.
\enspace \(\A{895}\)\enspace 
895, 14813. 
\enspace \(\A{895}^{\rm op}\)\enspace 
2552, 13238. 
\par

\refstepcounter{casegroup}\label{grp:90}\par\addvspace{.55ex}\noindent\textbf{Case 90.} FB.\enspace \ref{case:54}.
\enspace \(\A{2025}\)\enspace 
910, 2025. 
\enspace \(\A{2024}\)\enspace 
1892, 2024. 
\enspace \(\A{2025}^{\rm op}\)\enspace 
4489, 6993. 
\enspace \(\A{2024}^{\rm op}\)\enspace 
6144, 6992. 
\par

\refstepcounter{casegroup}\label{grp:91}\par\addvspace{.55ex}\noindent\textbf{Case 91.} FB.\enspace \ref{case:34}.
\enspace \(\A{1020}\)\enspace 
1010, 1020, 2583, 2991, 4115, 5895, 6314, 7737, 7849, 7899. 
\par

\refstepcounter{casegroup}\label{grp:92}\par\addvspace{.55ex}\noindent\textbf{Case 92.} FB.\enspace \ref{case:34}.
\enspace \(\A{5903}\)\enspace 
1011, 2584, 5875, 5903, 6315, 10955, 11012, 11018. 
\par

\refstepcounter{casegroup}\label{grp:93}\par\addvspace{.55ex}\noindent\textbf{Case 93.} FB.\enspace \ref{case:36}.
\enspace \(\A{1028}\)\enspace 
1014, 1024, 1028, 2996, 4262, 6371, 8709, 9139, 9143, 9229, 9850, 9865. 
\par

\refstepcounter{casegroup}\label{grp:94}\par\addvspace{.55ex}\noindent\textbf{Case 94.} FB.\enspace \ref{case:36}.
\enspace \(\A{1025}\)\enspace 
1015, 1025, 5941, 6372, 9140, 9245, 10961, 10968, 11020. 
\par

\refstepcounter{casegroup}\label{grp:95}\par\addvspace{.55ex}\noindent\textbf{Case 95.} FB.\enspace \ref{cor:periodic-register}.
1112--1114, 1118--1120, 1123, 1130, 1131, 1133, 1134, 1136, 1137, 1139, 1140, 1142, 1143, 1145, 1146, 1148--1152, 1156, 1192, 1199--1201, 1207--1213, 1217, 1218, 1244--1249, 1253, 1277, 1278, 1280, 1287, 1289, 1291, 1307, 1308, 1327, 1331, 1337, 1341, 1346, 1350, 1354, 1358, 1360, 1409--1411, 1413, 1414, 1420, 1423, 1426, 1427, 1429, 1430, 1432, 1433, 1435, 1436, 1441, 1446, 1450--1452, 1477, 1478, 1480, 1481, 1484--1487, 1490--1493, 1495, 1496, 1498, 1499, 1507, 1511, 1514, 1524, 1527, 1528, 1534, 1538, 1541, 1548, 1552, 1553, 1560, 1564, 1566, 1568, 1576, 1577, 1591, 1592, 1594, 1595, 1600, 1601, 1603, 1604, 1613--1616, 1618, 1619, 1621, 1622, 1630, 1631, 1633, 1634, 1636, 1637, 1660, 1697, 1731, 1732, 1789, 1790, 1835, 1837, 1838, 1845, 1846, 1848, 1849, 1851, 1852, 1854, 1855, 1857, 1858, 1860--1864, 1868, 1889, 1958, 1962, 1965--1967, 1969, 1970, 1972, 1976, 1984, 1987, 1990, 1993, 1994, 1996, 1997, 1999, 2000, 2002, 2003, 2015, 2016, 2026, 2027, 2030, 2031, 2040--2042, 2049--2055, 2057, 2059, 2060, 2063, 2065--2073, 2097--2103, 2110, 2111, 2196, 2208, 2630--2632, 2635--2645, 2656, 2658, 2659, 2672, 2701--2704, 2709, 2713--2718, 2723, 2725--2727, 2741--2753, 2768--2771, 2773, 2774, 2779, 2781, 2783--2785, 2790--2794, 2829, 2830, 2850, 2851, 2864--2868, 2885, 2886, 2890--2893, 2898, 2900, 2901, 2903, 3022, 3024, 3028--3033, 3036--3039, 3052, 3053, 3056--3059, 3061, 3063, 3064, 3068, 3070, 3071, 3073, 3076--3079, 3081, 3083, 3084, 3086, 3087, 3101--3104, 3107--3110, 3360--3381, 3384, 3386--3394, 3397, 3420--3436, 3440, 3442--3451, 3454, 3458--3467, 3470, 3472--3486, 3488, 3490--3499, 3502, 3504--3508, 3510, 3538--3546, 3549, 3551--3560, 3563, 3564, 3795, 3797, 3808, 3816, 3819, 3823, 3829, 3833, 3889, 3893, 3900, 3904, 3908, 3913, 3940, 3948, 3952, 3956, 3963, 3966, 3968, 3971, 3975, 3979, 3983, 3985, 4017, 4018, 4020, 4021, 4040, 4046, 4051, 4052, 4058, 4059, 4126, 4132, 4141--4143, 4149, 4154--4156, 4163, 4165, 4168--4171, 4175, 4176, 4180, 4182--4186, 4245, 4249, 4258--4261, 4267--4270, 4274, 4282, 4283, 4285, 4286, 4291, 4292, 4295--4297, 4299, 4301, 4302, 4375--4377, 4380, 4381, 4390--4395, 4399--4404, 4409--4412, 4415, 4416, 4423, 4431, 4438--4443, 4448--4450, 4460--4463, 4466--4468, 4472, 4474--4477, 4537, 4538, 4541, 4542, 4544, 4546, 4549--4552, 4554, 4555, 4559--4564, 4583, 4584, 4586, 4587, 4745--4762, 4765, 4767--4775, 4777, 4779--4782, 4784, 4787, 4790, 4794--4837, 4841--4844, 4884--4886, 4894, 4897, 4898, 4902--4917, 4919, 4920, 4938--4940, 5038, 5040, 5051, 5052, 5055, 5056, 5062, 5075, 5078, 5097, 5098, 5102, 5106, 5109, 5113, 5115, 5119, 5121, 5123--5125, 5131, 5134, 5135, 5143--5146, 5153, 5154, 5161, 5162, 5171, 5172, 5178, 5179, 5197, 5205, 5206, 5210--5213, 5218, 5224, 5226, 5230, 5239, 5445, 5447--5449, 5474, 5477, 5488, 5489, 5502, 5508, 5509, 5667, 5675, 5679, 5727--5729, 5733, 5748, 5751, 5780--5785, 5789, 5810, 5813, 5815, 5833--5837, 5841, 5883, 5885, 5887, 5889--5892, 5896, 5914, 5917, 5948, 5950, 5951, 5957, 5958, 5960, 5961, 5963, 5964, 5966--5971, 5978--5980, 6000--6002, 6005--6008, 6015--6017, 6048--6055, 6079--6090, 6108--6115, 6123, 6124, 6132, 6134, 6135, 6142, 6150--6153, 6167, 6169, 6171, 6173, 6174, 6181--6183, 6197, 6199, 6201, 6228, 6232, 6235, 6260, 6261, 6265--6267, 6271--6273, 6280, 6281, 6285--6289, 6293, 6322, 6323, 6327, 6329, 6331, 6333, 6335, 6385--6387, 6390, 6391, 6396, 6397, 6399, 6400, 6402, 6403, 6405, 6406, 6408, 6409, 6411, 6412, 6417, 6422, 6423, 6427--6430, 6462, 6463, 6468, 6469, 6472--6475, 6520--6525, 6529--6546, 6553--6561, 6594--6596, 6599--6609, 6613--6615, 6627, 6629, 6632, 6633, 6649, 6650, 6652, 6653, 6656--6659, 6661, 6662, 6702, 6704--6718, 6721, 6734, 6736, 6739, 6742, 6744, 6747, 6760--6772, 6793--6808, 6810, 6812, 6814--6818, 6820, 6862, 6866, 6870, 6896, 6901--6903, 6909, 6913--6916, 6920, 6922--6926, 6930, 6931, 6950, 6951, 6957, 6959, 6965, 6969--6971, 6974, 6983, 6984, 6997, 6998, 7002--7004, 7008, 7009, 7013--7015, 7019, 7020, 7024--7026, 7030--7033, 7075, 7076, 7078, 7080, 7083, 7095, 7096, 7100, 7102, 7107, 7108, 7113, 7114, 7116, 7117, 7143, 7145, 7150--7155, 7158--7160, 7166, 7167, 7171--7176, 7179, 7180, 7185--7190, 7193--7195, 7225--7229, 7232, 7233, 7235, 7237, 7241, 7250--7253, 7255, 7256, 7258, 7259, 7261, 7262, 7280, 7328, 7329, 7331, 7333, 7334, 7337, 7339--7344, 7346, 7348--7359, 7362, 7364, 7365, 7371, 7372, 7375, 7378, 7379, 7385--7401, 7403, 7404, 7416--7419, 7431--7433, 7444, 7445, 7449, 7452, 7460, 7464, 7466, 7469, 7471, 7473, 7476, 7479--7495, 7504--7513, 7515--7518, 7521--7523, 7525, 7528--7530, 7534--7550, 7552, 7553, 7574--7578, 7596, 7598, 7702, 7706--7709, 7719, 7720, 7723, 7724, 7730, 7731, 7735, 7738--7741, 7746, 7751, 7753, 7757, 7766, 7778, 7780, 7786, 7787, 7801, 7802, 7817, 7822, 7823, 7827--7830, 7835--7838, 7842, 7843, 7848, 7850, 7854, 7858, 7859, 7862, 7865--7867, 7871--7873, 7878, 7883, 7885, 7889, 7897, 7900, 7963, 8099, 8100, 8112, 8115, 8116, 8127, 8131, 8148, 8150, 8162--8165, 8169, 8189, 8193, 8197, 8201, 8221, 8222, 8226, 8230--8232, 8249--8251, 8255, 8277, 8279, 8291, 8295, 8310, 8311, 8315, 8349, 8821. 
\par

\refstepcounter{casegroup}\label{grp:96}\par\addvspace{.55ex}\noindent\textbf{Case 96.} FB.\enspace \ref{case:55}.
\enspace \(\A{13219}\)\enspace 
1166, 1709, 8326, 13219, 14664. 
\par

\refstepcounter{casegroup}\label{grp:97}\par\addvspace{.55ex}\noindent\textbf{Case 97.} FB.\enspace \ref{case:56}.
\enspace \(\A{9824}^{\rm op}\)\enspace 
1182, 10002, 13419. 
\enspace \(\A{9824}\)\enspace 
3065, 9824, 12182. 
\enspace \(\A{9106}\)\enspace 
9106, 9108, 9256. 
\enspace \(\A{9106}^{\rm op}\)\enspace 
9399, 9405, 9589. 
\par

\refstepcounter{casegroup}\label{grp:98}\par\addvspace{.55ex}\noindent\textbf{Case 98.} FB.\enspace \ref{case:57}.
\enspace \(\A{1222}\)\enspace 
1222, 2132, 2145, 2895. 
\enspace \(\A{1223}\)\enspace 
1223, 2136, 2149, 6616. 
\enspace \(\A{1643}\)\enspace 
1643, 4923, 7207, 7921, 8904, 8955, 8961, 8977. 
\enspace \(\A{1649}\)\enspace 
1649, 4932, 7219, 7956, 9175, 9200, 9254, 9276. 
\enspace \(\A{1680}\)\enspace 
1680, 5412, 8512, 8682, 8909, 8956, 8962, 8978. 
\enspace \(\A{1778}\)\enspace 
1778, 9176, 9201, 9277, 10057, 10072. 
\enspace \(\A{2797}\)\enspace 
2797, 7408, 7623, 7631, 8963, 9410, 9583, 9643. 
\enspace \(\A{2801}\)\enspace 
2801, 7425, 7683, 7691, 9206, 9447, 9603, 9657. 
\enspace \(\A{2810}\)\enspace 
2810, 8035, 8606, 8616, 8964, 9411, 9587, 9644. 
\enspace \(\A{2844}\)\enspace 
2844, 9207, 9448, 9658, 10062, 10070. 
\enspace \(\A{4571}\)\enspace 
4571, 4924, 7049, 7123. 
\enspace \(\A{4574}\)\enspace 
4574, 7055, 7128, 7566. 
\enspace \(\A{4580}\)\enspace 
4580, 4934, 7069, 7135. 
\enspace \(\A{4848}\)\enspace 
4848, 7201, 7409, 7910. 
\enspace \(\A{4866}\)\enspace 
4866, 7211, 7429, 7941. 
\enspace \(\A{5285}\)\enspace 
5285, 7593, 8370, 8439. 
\enspace \(\A{5294}\)\enspace 
5294, 7594, 8384, 8452. 
\enspace \(\A{5300}\)\enspace 
5300, 5416, 8399, 8455. 
\enspace \(\A{5373}\)\enspace 
5373, 6725, 8495, 8654. 
\enspace \(\A{5383}\)\enspace 
5383, 8037, 8502, 8668. 
\enspace \(\A{6830}\)\enspace 
6830, 10123, 10208. 
\enspace \(\A{6831}\)\enspace 
6831, 10140, 10227. 
\enspace \(\A{6839}\)\enspace 
6839, 10276, 10334. 
\enspace \(\A{7651}\)\enspace 
7651, 10129, 10209. 
\enspace \(\A{7653}\)\enspace 
7653, 7696, 10163, 10176, 10243, 10257. 
\enspace \(\A{8086}\)\enspace 
8086, 8620, 10149, 10154, 10229, 10230. 
\enspace \(\A{8791}\)\enspace 
8791, 10134, 10211. 
\par

\refstepcounter{casegroup}\label{grp:99}\par\addvspace{.55ex}\noindent\textbf{Case 99.} FB.\enspace \ref{thm:twofactor}.
\enspace \(\A{643}\times M\)\enspace 
1228, 2870, 12953, 13701. 
\enspace \(\A{4609}\times M\)\enspace 
5377, 6156, 11862, 11999. 
\par

\refstepcounter{casegroup}\label{grp:100}\par\addvspace{.55ex}\noindent\textbf{Case 100.} FB.\enspace \ref{case:58}.
\enspace \(\A{1724}\)\enspace 
1300, 1724, 9975. 
\par

\refstepcounter{casegroup}\label{grp:101}\par\addvspace{.55ex}\noindent\textbf{Case 101.} FB.\enspace \ref{case:59}.
\enspace \(\A{1725}\)\enspace 
1301, 1725, 8122, 13240. 
\par

\refstepcounter{casegroup}\label{grp:102}\par\addvspace{.55ex}\noindent\textbf{Case 102.} FB.\enspace \ref{case:60}.
\enspace \(\A{1302}\)\enspace 
1302, 1726, 14817. 
\enspace \(\A{1302}^{\rm op}\)\enspace 
2668, 3992, 13241. 
\par

\refstepcounter{casegroup}\label{grp:103}\par\addvspace{.55ex}\noindent\textbf{Case 103.} FB.\enspace \ref{case:40}.
\enspace \(\A{1959}^{\rm op}\)\enspace 
1461, 3017, 4080, 4202, 4618, 4674, 5879, 5907, 6229, 6776, 8120, 11833. 
\par

\refstepcounter{casegroup}\label{grp:104}\par\addvspace{.55ex}\noindent\textbf{Case 104.} FB.\enspace \ref{case:61}.
\enspace \(\A{11794}\)\enspace 
1469, 4367, 5006, 5511, 11784, 11787, 11794, 11882. 
\par

\refstepcounter{casegroup}\label{grp:105}\par\addvspace{.55ex}\noindent\textbf{Case 105.} FB.\enspace \ref{case:62}.
\enspace \(\A{1471}\)\enspace 
1471, 5362, 5363, 9490, 9492. 
\enspace \(\A{1471}^{\rm op}\)\enspace 
3041, 5365, 5366, 9493, 9495. 
\par

\refstepcounter{casegroup}\label{grp:106}\par\addvspace{.55ex}\noindent\textbf{Case 106.} FB.\enspace \ref{case:63}.
\enspace \(\A{1473}\)\enspace 
1473, 1475, 2740, 3143--3146, 3357, 4334, 4343, 4352, 4645, 4687, 4688, 4986, 4998, 5003, 6784, 11791, 11886. 
\par

\refstepcounter{casegroup}\label{grp:107}\par\addvspace{.55ex}\noindent\textbf{Case 107.} FB.\enspace \ref{case:58}.
\enspace \(\A{1522}\)\enspace 
1522, 1785, 9974, 14828. 
\par

\refstepcounter{casegroup}\label{grp:108}\par\addvspace{.55ex}\noindent\textbf{Case 108.} FB.\enspace \ref{case:59}.
\enspace \(\A{1525}\)\enspace 
1525, 1786, 8121, 13207, 14633. 
\par

\refstepcounter{casegroup}\label{grp:109}\par\addvspace{.55ex}\noindent\textbf{Case 109.} FB.\enspace \ref{case:64}.
\enspace \(\A{1787}\)\enspace 
1529, 1787, 14658. 
\enspace \(\A{1787}^{\rm op}\)\enspace 
2848, 3964, 13212. 
\par

\refstepcounter{casegroup}\label{grp:110}\par\addvspace{.55ex}\noindent\textbf{Case 110.} FB.\enspace \ref{case:65}.
\enspace \(\A{1638}\)\enspace 
1638, 9807. 
\par

\refstepcounter{casegroup}\label{grp:111}\par\addvspace{.55ex}\noindent\textbf{Case 111.} FB.\enspace \ref{cor:prime-examples}.
1639, 1655, 2796, 4847, 7455, 7555, 7907. 
\par

\refstepcounter{casegroup}\label{grp:112}\par\addvspace{.55ex}\noindent\textbf{Case 112.} FB.\enspace \ref{case:66}.
\enspace \(\A{1640}\)\enspace 
1640, 2191. 
\par

\refstepcounter{casegroup}\label{grp:113}\par\addvspace{.55ex}\noindent\textbf{Case 113.} FB.\enspace \ref{prop:seven-bases}.
\enspace \(S_{3,50}^{\infty}\)\enspace 
1641, 1642, 1748, 2798, 4565, 4568, 4854, 4855, 5302, 7042, 7121, 7122, 7206, 7558, 7559, 7599, 7626, 7629, 7652, 7923, 8087, 8404, 8465, 8649, 8835, 8837, 8894, 8902, 8953, 8959, 8967, 8971, 8973, 9414, 9594, 9647, 9810, 10145, 10167, 10181, 10228, 10244, 10258. 
\enspace \(S_{3,51}^{\infty}\)\enspace 
1647, 1648, 1750, 2802, 4576, 4578, 4870, 4873, 5306, 7063, 7133, 7134, 7218, 7570, 7572, 7603, 7686, 7689, 7701, 7958, 8092, 8413, 8479, 8653, 8993, 8995, 9156, 9166, 9187, 9194, 9216, 9232, 9233, 9455, 9620, 9669, 9820, 10282, 10286, 10291, 10335--10337. 
\enspace \(S_{3,53}^{\infty}\)\enspace 
1678, 1679, 1749, 2811, 5296, 5298, 5304, 5387, 5390, 8062, 8064, 8068, 8392, 8409, 8453, 8454, 8466, 8511, 8610, 8613, 8688, 8793, 8795, 8800, 8854, 8856, 8898, 8903, 8954, 8960, 8968, 8972, 8974, 9416, 9599, 9648, 9950, 10159, 10171, 10186, 10231, 10245, 10259. 
\enspace \(S_{3,60}^{\infty}\)\enspace 
1751, 1776, 2845, 9167, 9195, 9234, 9235, 10052, 10056, 10074, 10095, 10340, 10341, 10514, 10570. 
\enspace \(S_{4,456}^{\infty}\)\enspace 
5415, 6133, 6275, 6388. 
\enspace \(S_{4,443}^{\infty}\)\enspace 
6168, 6735, 7428, 7783. 
\enspace \(S_{4,471}^{\infty}\)\enspace 
9160, 10086. 
\par

\refstepcounter{casegroup}\label{grp:114}\par\addvspace{.55ex}\noindent\textbf{Case 114.} FB.\enspace \ref{case:67}.
1650, 1651, 5288, 5378. 
\par

\refstepcounter{casegroup}\label{grp:115}\par\addvspace{.55ex}\noindent\textbf{Case 115.} FB.\enspace \ref{case:68}.
1652--1654, 2158, 2803, 4871, 7213, 7216, 7571, 7602, 7687, 7700, 7946, 7952, 8091, 8652, 9221, 9258, 9819. 
\par

\refstepcounter{casegroup}\label{grp:116}\par\addvspace{.55ex}\noindent\textbf{Case 116.} FB.\enspace \ref{case:69}.
\enspace \(\A{5291}^{\rm op}\)\enspace 
1690. 
\enspace \(\A{5291}\)\enspace 
5291. 
\par

\refstepcounter{casegroup}\label{grp:117}\par\addvspace{.55ex}\noindent\textbf{Case 117.} FB.\enspace \ref{case:70}.
\enspace \(\A{1706}\)\enspace 
1706, 13190. 
\par

\refstepcounter{casegroup}\label{grp:118}\par\addvspace{.55ex}\noindent\textbf{Case 118.} FB.\enspace \ref{case:71}.
\enspace \(\A{1777}\)\enspace 
1777, 9217, 9457, 9670, 10065, 10068. 
\enspace \(\A{7635}\)\enspace 
7635, 10033, 10045. 
\par

\refstepcounter{casegroup}\label{grp:119}\par\addvspace{.55ex}\noindent\textbf{Case 119.} FB.\enspace \ref{case:72}.
\enspace \(\A{1888}\)\enspace 
1888. 
\par

\refstepcounter{casegroup}\label{grp:120}\par\addvspace{.55ex}\noindent\textbf{Case 120.} FB.\enspace \ref{case:73}.
\enspace \(\A{9316}\)\enspace 
1891, 2022, 9316, 13766. 
\enspace \(\A{9316}^{\rm op}\)\enspace 
6143, 6990, 9946, 12169. 
\par

\refstepcounter{casegroup}\label{grp:121}\par\addvspace{.55ex}\noindent\textbf{Case 121.} FB.\enspace \ref{case:74}.
\enspace \(\A{2008}\)\enspace 
2008. 
\enspace \(\A{3111}\)\enspace 
3111. 
\enspace \(\A{2008}^{\rm op}\)\enspace 
6635. 
\enspace \(\A{3111}^{\rm op}\)\enspace 
8584. 
\par

\refstepcounter{casegroup}\label{grp:122}\par\addvspace{.55ex}\noindent\textbf{Case 122.} FB.\enspace Proposition~\ref{prop:lower-fibres}, source \(\A{2007}\), case~\ref{grp:2}.
\enspace \(\A{2007}\)\enspace 
2083, 2131. 
\par

\refstepcounter{casegroup}\label{grp:123}\par\addvspace{.55ex}\noindent\textbf{Case 123.} FB.\enspace \ref{abs-chain-extensions}.
\enspace \(\text{\normalfont M}^{\rm op}\)\enspace 
2085, 2086, 2150, 2184--2186, 6431, 6432, 6434, 9304. 
\enspace \(\text{\normalfont L}^{\rm op}\)\enspace 
2142, 8458, 8459, 8468, 8914, 8981, 8982, 8985, 8986, 8988, 9650, 10236, 10248, 10262. 
\enspace \(\text{\normalfont M}\)\enspace 
7273, 7274, 7975, 7976, 7989, 8532, 8719--8721, 9909. 
\enspace \(\text{\normalfont L}\)\enspace 
7310, 7311, 7681, 8555, 9651, 9720, 9721, 9732, 9906--9908, 10435--10437. 
\par

\refstepcounter{casegroup}\label{grp:124}\par\addvspace{.55ex}\noindent\textbf{Case 124.} FB.\enspace \ref{case:75}.
\enspace \(\A{2138}\)\enspace 
2138. 
\enspace \(\A{2138}^{\rm op}\)\enspace 
8372. 
\par

\refstepcounter{casegroup}\label{grp:125}\par\addvspace{.55ex}\noindent\textbf{Case 125.} FB.\enspace \ref{case:76}.
\enspace \(\A{2140}\)\enspace 
2140, 7627, 8908, 8969. 
\enspace \(\A{2141}\)\enspace 
2141, 8611, 8913, 8970. 
\enspace \(\A{2140}^{\rm op}\)\enspace 
7622, 7628, 9412, 9591. 
\enspace \(\A{2141}^{\rm op}\)\enspace 
7682, 7688, 9451, 9617. 
\par

\refstepcounter{casegroup}\label{grp:126}\par\addvspace{.55ex}\noindent\textbf{Case 126.} FB.\enspace \ref{case:77}.
\enspace \(\A{2183}\)\enspace 
2183. 
\enspace \(\A{2183}^{\rm op}\)\enspace 
8718. 
\par

\refstepcounter{casegroup}\label{grp:127}\par\addvspace{.55ex}\noindent\textbf{Case 127.} FB.\enspace \ref{case:78}.
\enspace \(\text{\normalfont T}\times M\)\enspace 
2310, 12209. 
\enspace \(\text{\normalfont T}\)\enspace 
10572, 12212, 12213, 12215, 12260, 12291--12293, 13369, 15140. 
\par

\refstepcounter{casegroup}\label{grp:128}\par\addvspace{.55ex}\noindent\textbf{Case 128.} FB.\enspace \ref{case:79}.
\enspace \(\A{2315}\)\enspace 
2315, 8825. 
\par

\refstepcounter{casegroup}\label{grp:129}\par\addvspace{.55ex}\noindent\textbf{Case 129.} FB.\enspace Proposition~\ref{prop:lower-fibres}, source \(\A{3751}\), case~\ref{grp:26}.
\enspace \(\A{3751}\)\enspace 
2326, 2340, 2470, 2474, 2628, 3324, 3350, 3749, 3766, 3771, 5041, 5455, 5544, 5546, 5737, 6044, 6045, 8136, 8173, 10855, 10859, 10869, 10974, 10978, 10984, 11079, 11088, 11091, 11228, 11252, 11265, 11266, 11345, 11351, 11356, 12242, 12353, 12365, 12511, 12933. 
\par

\refstepcounter{casegroup}\label{grp:130}\par\addvspace{.55ex}\noindent\textbf{Case 130.} FB.\enspace \ref{prop:three-direct}.
\enspace \(S_{3,24}\)\enspace 
2330, 2472, 2629, 3782, 3787, 3791, 3798, 5596, 5689, 5692, 5793, 10826, 10830, 10840, 11043, 11046, 11062, 11065, 11387, 11388, 11398, 11401, 11420, 11422, 11516, 11521, 11575, 11576, 11578, 11591, 11592, 12345, 12347, 12677, 12702, 12783, 12809, 12814, 13527. 
\enspace \(S_{3,37}\)\enspace 
2928--2931, 4101--4103, 4217, 4230, 4509, 4510, 5934, 5996, 5997, 6136, 6242, 6645, 7027, 7028, 11769, 11770, 11837, 11838, 11870--11872, 11914, 11932, 11937, 11986--11988, 12034, 12048, 12051, 12087--12094, 13069, 13070, 13124, 13128, 13129, 13283--13287, 13945--13947. 
\enspace \(S_{3,50}\)\enspace 
7085, 7263, 7264, 7281, 7282, 7291, 7292, 7300, 7301, 8001, 8011, 8014, 8015, 8553, 8563, 8566, 8567, 8759, 8764, 8767, 9712, 9714, 9744, 9745, 9896, 9897, 9903, 9925, 10411, 10419, 10425, 10660, 10661, 10676, 10683, 10691, 10706, 10714, 10715, 10732. 
\enspace \(S_{3,51}\)\enspace 
7088, 7270, 7271, 7285, 7286, 7313, 7314, 7325, 7326, 8003, 8021--8023, 8559, 8573--8575, 8772--8774, 9737, 9749, 9751, 9752, 9927--9930, 10448--10450, 10769--10773, 10776--10778, 10782. 
\enspace \(S_{3,53}\)\enspace 
7994, 7995, 7998, 7999, 8002, 8445, 8543, 8544, 8554, 8735, 8736, 8741, 8744, 8745, 8749, 8752, 8753, 8762, 8765, 8768, 9713, 9895, 9900, 9901, 9905, 9920, 9921, 9926, 10414, 10420, 10426, 10671, 10672, 10680, 10686, 10692, 10709, 10717, 10718, 10734. 
\par

\refstepcounter{casegroup}\label{grp:131}\par\addvspace{.55ex}\noindent\textbf{Case 131.} FB.\enspace \ref{case:80}.
\enspace \(\A{2350}\)\enspace 
2350, 14404. 
\par

\refstepcounter{casegroup}\label{grp:132}\par\addvspace{.55ex}\noindent\textbf{Case 132.} FB.\enspace \ref{thm:unit-paths}.
2535, 3800, 5030, 5199, 5576, 5694, 5759, 5773, 7969, 7970, 8354--8359, 12743, 12744, 12876, 12878. 
\par

\refstepcounter{casegroup}\label{grp:133}\par\addvspace{.55ex}\noindent\textbf{Case 133.} FB.\enspace \ref{case:79}.
\enspace \(\A{2553}\)\enspace 
2553, 14815. 
\par

\refstepcounter{casegroup}\label{grp:134}\par\addvspace{.55ex}\noindent\textbf{Case 134.} FB.\enspace \ref{case:81}.
\enspace \(\A{3993}\)\enspace 
2669, 3993, 5438, 5495, 14819, 14820, 14907, 14917, 14918, 15582. 
\par

\refstepcounter{casegroup}\label{grp:135}\par\addvspace{.55ex}\noindent\textbf{Case 135.} FB.\enspace \ref{case:61}.
\enspace \(\A{11794}^{\rm op}\)\enspace 
2736, 4675, 4700, 8102, 11788, 11839, 11913, 11979. 
\par

\refstepcounter{casegroup}\label{grp:136}\par\addvspace{.55ex}\noindent\textbf{Case 136.} FB.\enspace \ref{case:66}.
\enspace \(\A{2795}\)\enspace 
2795, 7607. 
\par

\refstepcounter{casegroup}\label{grp:137}\par\addvspace{.55ex}\noindent\textbf{Case 137.} FB.\enspace \ref{case:80}.
\enspace \(\A{2822}\)\enspace 
2822, 14405. 
\par

\refstepcounter{casegroup}\label{grp:138}\par\addvspace{.55ex}\noindent\textbf{Case 138.} FB.\enspace \ref{case:81}.
\enspace \(\A{3969}\)\enspace 
2849, 3203, 3969, 5494, 10820, 14523, 14544, 14553, 14660, 14661, 14715--14717, 14748, 15552. 
\par

\refstepcounter{casegroup}\label{grp:139}\par\addvspace{.55ex}\noindent\textbf{Case 139.} FB.\enspace \ref{thm:one-sided}.
2913, 2915, 3117, 3119, 4077, 4601, 4953, 5428, 5861, 5863, 6139, 6145, 6184, 6978, 7710, 8096, 8097, 8144, 9536, 9967, 9968, 9994, 10583, 11812, 11858, 11860, 11902, 11974, 11977, 11978, 12016, 12017, 12044, 12066--12069, 12141--12144, 12146--12148, 12158, 12172, 12186, 12247, 13084, 13142, 13145, 13308, 13309, 13312, 13344, 13951. 
\par

\refstepcounter{casegroup}\label{grp:140}\par\addvspace{.55ex}\noindent\textbf{Case 140.} FB.\enspace \ref{case:82}.
\enspace \(\A{10379}\)\enspace 
2939, 4094, 6345, 9068, 9243, 10379. 
\par

\refstepcounter{casegroup}\label{grp:141}\par\addvspace{.55ex}\noindent\textbf{Case 141.} FB.\enspace \ref{case:83}.
\enspace \(\A{4095}\)\enspace 
2940, 4095, 10588, 10589. 
\par

\refstepcounter{casegroup}\label{grp:142}\par\addvspace{.55ex}\noindent\textbf{Case 142.} FB.\enspace \ref{case:84}.
\enspace \(\A{2941}\)\enspace 
2941, 13324. 
\par

\refstepcounter{casegroup}\label{grp:143}\par\addvspace{.55ex}\noindent\textbf{Case 143.} FB.\enspace Section~\ref{sec:four-sources}.
\enspace \(S_{4,430}\)\enspace 
2943, 5988, 10110, 10111, 10573, 10577, 10584, 10596, 12081, 12249, 13245, 13276, 13278. 
\par

\refstepcounter{casegroup}\label{grp:144}\par\addvspace{.55ex}\noindent\textbf{Case 144.} FB.\enspace \ref{case:1}.
\enspace \(S_{4,627}\)\enspace 
2944, 12248, 13060, 13103, 13106, 13187, 13194--13196, 13834, 13843, 13888, 13932, 13933, 13936. 
\par

\refstepcounter{casegroup}\label{grp:145}\par\addvspace{.55ex}\noindent\textbf{Case 145.} FB.\enspace \ref{case:85}.
\enspace \(\A{6023}\)\enspace 
2945, 6023, 6349. 
\enspace \(\A{6022}\)\enspace 
2946, 6022. 
\enspace \(\A{6023}^{\rm op}\)\enspace 
3169, 6484, 8205. 
\enspace \(\A{6022}^{\rm op}\)\enspace 
3170, 8204. 
\par

\refstepcounter{casegroup}\label{grp:146}\par\addvspace{.55ex}\noindent\textbf{Case 146.} FB.\enspace \ref{case:86}.
\enspace \(\A{2947}\)\enspace 
2947, 13188. 
\enspace \(\A{2947}^{\rm op}\)\enspace 
3171, 13104. 
\par

\refstepcounter{casegroup}\label{grp:147}\par\addvspace{.55ex}\noindent\textbf{Case 147.} FB.\enspace \ref{case:33}.
\enspace \(\A{1549}^{\rm op}\)\enspace 
2955, 3198, 3201, 4133, 5110, 7886, 13067, 13068, 13116, 13119, 13120, 13204, 13205, 13228--13231, 13870, 13878, 13942, 13943. 
\par

\refstepcounter{casegroup}\label{grp:148}\par\addvspace{.55ex}\noindent\textbf{Case 148.} FB.\enspace \ref{case:61}.
\enspace \(\A{3049}\)\enspace 
2966, 3049, 3268, 3527, 4108, 4160, 4982, 5215, 6018. 
\par

\refstepcounter{casegroup}\label{grp:149}\par\addvspace{.55ex}\noindent\textbf{Case 149.} FB.\enspace \ref{case:61}.
\enspace \(\A{4159}\)\enspace 
2967, 4159, 5214, 13215, 13859. 
\par

\refstepcounter{casegroup}\label{grp:150}\par\addvspace{.55ex}\noindent\textbf{Case 150.} FB.\enspace \ref{case:87}.
\enspace \(\A{8627}\)\enspace 
2973, 8627, 9783. 
\enspace \(\A{8627}^{\rm op}\)\enspace 
2997, 8696, 9849. 
\par

\refstepcounter{casegroup}\label{grp:151}\par\addvspace{.55ex}\noindent\textbf{Case 151.} FB.\enspace \ref{case:88}.
\enspace \(\A{11019}^{\rm op}\)\enspace 
2974, 13317, 13824. 
\enspace \(\A{11019}\)\enspace 
2998, 11019, 11297. 
\par

\refstepcounter{casegroup}\label{grp:152}\par\addvspace{.55ex}\noindent\textbf{Case 152.} FB.\enspace \ref{case:89}.
\enspace \(\A{10995}^{\rm op}\)\enspace 
2977, 3247, 3258, 5067, 13050, 13051, 13090, 13093, 13147, 13148, 13151, 13152, 13163, 13164, 13167, 13168, 13796, 13800, 13803, 13812, 13816, 13914, 13918, 13921. 
\par

\refstepcounter{casegroup}\label{grp:153}\par\addvspace{.55ex}\noindent\textbf{Case 153.} FB.\enspace \ref{case:90}.
\enspace \(\A{4277}\)\enspace 
2978, 4277, 4921, 6012, 6252, 6362. 
\par

\refstepcounter{casegroup}\label{grp:154}\par\addvspace{.55ex}\noindent\textbf{Case 154.} FB.\enspace \ref{case:91}.
\enspace \(\A{2979}\)\enspace 
2979, 4899, 6011. 
\enspace \(\A{2980}\)\enspace 
2980, 13160, 13806. 
\enspace \(\A{2979}^{\rm op}\)\enspace 
3298, 7736, 7800. 
\enspace \(\A{2980}^{\rm op}\)\enspace 
3299, 11011, 11231. 
\par

\refstepcounter{casegroup}\label{grp:155}\par\addvspace{.55ex}\noindent\textbf{Case 155.} FB.\enspace \ref{case:40}.
\enspace \(\A{3094}\)\enspace 
2982, 3094, 4173, 5283, 5284, 8428, 8429, 9444, 9446, 9639. 
\par

\refstepcounter{casegroup}\label{grp:156}\par\addvspace{.55ex}\noindent\textbf{Case 156.} FB.\enspace \ref{case:40}.
\enspace \(\A{4172}\)\enspace 
2983, 4172, 9638, 13322, 13831. 
\par

\refstepcounter{casegroup}\label{grp:157}\par\addvspace{.55ex}\noindent\textbf{Case 157.} FB.\enspace \ref{case:40}.
\enspace \(\A{4177}\)\enspace 
2986, 2989, 4177, 4187, 13056, 13059, 13099, 13175, 13177, 13180--13182, 13827, 13829, 13830, 13924, 13925. 
\par

\refstepcounter{casegroup}\label{grp:158}\par\addvspace{.55ex}\noindent\textbf{Case 158.} FB.\enspace \ref{case:39}.
2995, 3292, 3294, 10998, 11004, 11017, 11219, 11220, 11235, 11239, 11240, 11285--11287, 11290, 11380, 12499, 12526, 12943. 
\par

\refstepcounter{casegroup}\label{grp:159}\par\addvspace{.55ex}\noindent\textbf{Case 159.} FB.\enspace \ref{case:89}.
\enspace \(\A{10995}\)\enspace 
3001, 3302, 3313, 5138, 10995, 10996, 11001, 11002, 11015, 11215, 11216, 11230, 11233, 11275, 11276, 11278, 11279, 11282, 11283, 11373, 11377, 12495, 12523, 12940. 
\par

\refstepcounter{casegroup}\label{grp:160}\par\addvspace{.55ex}\noindent\textbf{Case 160.} FB.\enspace \ref{case:92}.
\enspace \(\A{4122}\)\enspace 
3002, 4122, 4151, 4264, 5225, 6014, 6255, 6375. 
\par

\refstepcounter{casegroup}\label{grp:161}\par\addvspace{.55ex}\noindent\textbf{Case 161.} FB.\enspace \ref{case:93}.
\enspace \(\A{4123}\)\enspace 
3003, 4123, 4153, 5208, 6013. 
\par

\refstepcounter{casegroup}\label{grp:162}\par\addvspace{.55ex}\noindent\textbf{Case 162.} FB.\enspace \ref{case:94}.
\enspace \(\A{4152}\)\enspace 
3004, 4152, 5207, 11007, 11274. 
\par

\refstepcounter{casegroup}\label{grp:163}\par\addvspace{.55ex}\noindent\textbf{Case 163.} FB.\enspace \ref{case:61}.
\enspace \(\A{4109}\)\enspace 
3020, 3279, 3534, 4109, 4229, 5952, 6028, 6268. 
\par

\refstepcounter{casegroup}\label{grp:164}\par\addvspace{.55ex}\noindent\textbf{Case 164.} FB.\enspace \ref{case:63}.
\enspace \(\A{1473}^{\rm op}\)\enspace 
3043, 3148--3151, 3358, 3401, 3403, 4355, 4361, 4668, 4672, 4690, 4691, 4698, 5010, 5015, 6786, 11840, 11983. 
\par

\refstepcounter{casegroup}\label{grp:165}\par\addvspace{.55ex}\noindent\textbf{Case 165.} FB.\enspace \ref{abs-chain-extensions}.
\enspace \(\text{\normalfont P}^{\rm op}\)\enspace 
3112. 
\enspace \(\text{\normalfont N}^{\rm op}\)\enspace 
6433. 
\enspace \(\text{\normalfont N}\)\enspace 
7584. 
\enspace \(\text{\normalfont P}\)\enspace 
8535. 
\par

\refstepcounter{casegroup}\label{grp:166}\par\addvspace{.55ex}\noindent\textbf{Case 166.} FB.\enspace \ref{case:95}.
\enspace \(\A{4325}\)\enspace 
3167, 4325, 8195. 
\par

\refstepcounter{casegroup}\label{grp:167}\par\addvspace{.55ex}\noindent\textbf{Case 167.} FB.\enspace \ref{case:96}.
\enspace \(\A{3168}\)\enspace 
3168, 14415. 
\par

\refstepcounter{casegroup}\label{grp:168}\par\addvspace{.55ex}\noindent\textbf{Case 168.} FB.\enspace \ref{case:97}.
\enspace \(\A{6486}\)\enspace 
3172, 6486, 8238. 
\enspace \(\A{6487}\)\enspace 
6062, 6487, 8242. 
\par

\refstepcounter{casegroup}\label{grp:169}\par\addvspace{.55ex}\noindent\textbf{Case 169.} FB.\enspace \ref{case:98}.
\enspace \(\A{3173}\)\enspace 
3173, 8237. 
\enspace \(\A{3174}\)\enspace 
3174, 13191. 
\par

\refstepcounter{casegroup}\label{grp:170}\par\addvspace{.55ex}\noindent\textbf{Case 170.} FB.\enspace \ref{case:99}.
\enspace \(\A{3177}\)\enspace 
3177, 6488, 8258. 
\enspace \(\A{3177}^{\rm op}\)\enspace 
3180, 6572, 8298. 
\par

\refstepcounter{casegroup}\label{grp:171}\par\addvspace{.55ex}\noindent\textbf{Case 171.} FB.\enspace \ref{case:100}.
\enspace \(\A{3178}\)\enspace 
3178. 
\enspace \(\A{3179}\)\enspace 
3179. 
\enspace \(\A{3178}^{\rm op}\)\enspace 
3181. 
\enspace \(\A{3179}^{\rm op}\)\enspace 
3182. 
\par

\refstepcounter{casegroup}\label{grp:172}\par\addvspace{.55ex}\noindent\textbf{Case 172.} FB.\enspace \ref{case:61}.
\enspace \(\A{3049}^{\rm op}\)\enspace 
3206, 3406, 4615, 4646, 4973, 5129, 6690, 7743, 7803. 
\par

\refstepcounter{casegroup}\label{grp:173}\par\addvspace{.55ex}\noindent\textbf{Case 173.} FB.\enspace \ref{case:61}.
\enspace \(\A{4159}^{\rm op}\)\enspace 
3207, 5128, 7742, 13112, 13864. 
\par

\refstepcounter{casegroup}\label{grp:174}\par\addvspace{.55ex}\noindent\textbf{Case 174.} FB.\enspace \ref{case:101}.
\enspace \(\A{3216}\)\enspace 
3216, 3411. 
\enspace \(\A{3216}^{\rm op}\)\enspace 
3220, 3413. 
\par

\refstepcounter{casegroup}\label{grp:175}\par\addvspace{.55ex}\noindent\textbf{Case 175.} FB.\enspace \ref{case:102}.
\enspace \(\A{3217}\)\enspace 
3217. 
\enspace \(\A{3217}^{\rm op}\)\enspace 
3221. 
\par

\refstepcounter{casegroup}\label{grp:176}\par\addvspace{.55ex}\noindent\textbf{Case 176.} FB.\enspace \ref{case:103}.
\enspace \(\A{3489}\)\enspace 
3240, 3489, 4429, 5046, 7448, 8191. 
\par

\refstepcounter{casegroup}\label{grp:177}\par\addvspace{.55ex}\noindent\textbf{Case 177.} FB.\enspace \ref{case:104}.
\enspace \(\A{10924}^{\rm op}\)\enspace 
3241, 5466, 14004, 14100, 14177. 
\enspace \(\A{10924}\)\enspace 
3296, 5481, 10924, 11339, 11748. 
\par

\refstepcounter{casegroup}\label{grp:178}\par\addvspace{.55ex}\noindent\textbf{Case 178.} FB.\enspace \ref{case:92}.
\enspace \(\A{4122}^{\rm op}\)\enspace 
3242, 4432, 4980, 5048, 6496, 6878, 7405, 8199. 
\par

\refstepcounter{casegroup}\label{grp:179}\par\addvspace{.55ex}\noindent\textbf{Case 179.} FB.\enspace \ref{case:93}.
\enspace \(\A{4123}^{\rm op}\)\enspace 
3243, 4981, 5050, 7381, 7461. 
\par

\refstepcounter{casegroup}\label{grp:180}\par\addvspace{.55ex}\noindent\textbf{Case 180.} FB.\enspace \ref{case:94}.
\enspace \(\A{4152}^{\rm op}\)\enspace 
3244, 5049, 7380, 13091, 13809. 
\par

\refstepcounter{casegroup}\label{grp:181}\par\addvspace{.55ex}\noindent\textbf{Case 181.} FB.\enspace \ref{case:92}.
\enspace \(\A{4397}^{\rm op}\)\enspace 
3245, 4433, 4846, 5066, 6498, 6886, 8228. 
\par

\refstepcounter{casegroup}\label{grp:182}\par\addvspace{.55ex}\noindent\textbf{Case 182.} FB.\enspace \ref{case:34}.
\enspace \(\A{4413}^{\rm op}\)\enspace 
3250, 3511, 4434, 5071, 6500, 8253. 
\par

\refstepcounter{casegroup}\label{grp:183}\par\addvspace{.55ex}\noindent\textbf{Case 183.} FB.\enspace \ref{case:105}.
\enspace \(\A{3251}\)\enspace 
3251. 
\enspace \(\A{3252}\)\enspace 
3252. 
\enspace \(\A{3251}^{\rm op}\)\enspace 
3309. 
\enspace \(\A{3252}^{\rm op}\)\enspace 
3310. 
\par

\refstepcounter{casegroup}\label{grp:184}\par\addvspace{.55ex}\noindent\textbf{Case 184.} FB.\enspace \ref{case:106}.
\enspace \(\A{4444}\)\enspace 
3253, 4444, 4778, 5082, 6580, 6942, 6949, 8293. 
\par

\refstepcounter{casegroup}\label{grp:185}\par\addvspace{.55ex}\noindent\textbf{Case 185.} FB.\enspace \ref{case:107}.
\enspace \(\A{3254}\)\enspace 
3254. 
\enspace \(\A{3255}\)\enspace 
3255. 
\enspace \(\A{3254}^{\rm op}\)\enspace 
3306. 
\enspace \(\A{3255}^{\rm op}\)\enspace 
3307. 
\par

\refstepcounter{casegroup}\label{grp:186}\par\addvspace{.55ex}\noindent\textbf{Case 186.} FB.\enspace \ref{case:72}.
\enspace \(\A{13828}\)\enspace 
3269, 13178, 13828. 
\par

\refstepcounter{casegroup}\label{grp:187}\par\addvspace{.55ex}\noindent\textbf{Case 187.} FB.\enspace \ref{case:108}.
\enspace \(\A{3273}\)\enspace 
3273, 3530. 
\enspace \(\A{3274}\)\enspace 
3274. 
\enspace \(\A{3273}^{\rm op}\)\enspace 
4623, 6695. 
\enspace \(\A{3274}^{\rm op}\)\enspace 
4624. 
\par

\refstepcounter{casegroup}\label{grp:188}\par\addvspace{.55ex}\noindent\textbf{Case 188.} FB.\enspace \ref{case:109}.
\enspace \(\A{3276}\)\enspace 
3276, 3532. 
\enspace \(\A{3276}^{\rm op}\)\enspace 
4620, 6693. 
\par

\refstepcounter{casegroup}\label{grp:189}\par\addvspace{.55ex}\noindent\textbf{Case 189.} FB.\enspace \ref{case:110}.
\enspace \(\A{3277}\)\enspace 
3277. 
\enspace \(\A{3277}^{\rm op}\)\enspace 
4621. 
\par

\refstepcounter{casegroup}\label{grp:190}\par\addvspace{.55ex}\noindent\textbf{Case 190.} FB.\enspace \ref{case:72}.
\enspace \(\A{13057}\)\enspace 
3280, 13057, 13176, 13179. 
\par

\refstepcounter{casegroup}\label{grp:191}\par\addvspace{.55ex}\noindent\textbf{Case 191.} FB.\enspace \ref{case:103}.
\enspace \(\A{3967}\)\enspace 
3295, 3967, 4383, 5120, 8188, 8192. 
\par

\refstepcounter{casegroup}\label{grp:192}\par\addvspace{.55ex}\noindent\textbf{Case 192.} FB.\enspace \ref{case:90}.
\enspace \(\A{4396}\)\enspace 
3297, 4396, 6508, 6891, 7752, 8200. 
\par

\refstepcounter{casegroup}\label{grp:193}\par\addvspace{.55ex}\noindent\textbf{Case 193.} FB.\enspace \ref{case:92}.
\enspace \(\A{4397}\)\enspace 
3300, 4397, 5137, 5164, 6510, 6893, 8229. 
\par

\refstepcounter{casegroup}\label{grp:194}\par\addvspace{.55ex}\noindent\textbf{Case 194.} FB.\enspace \ref{case:106}.
\enspace \(\A{4398}\)\enspace 
3305, 3984, 4398, 5142, 6512, 6895, 6908, 8254. 
\par

\refstepcounter{casegroup}\label{grp:195}\par\addvspace{.55ex}\noindent\textbf{Case 195.} FB.\enspace \ref{case:34}.
\enspace \(\A{4413}\)\enspace 
3308, 4413, 5114, 5166, 6586, 8294. 
\par

\refstepcounter{casegroup}\label{grp:196}\par\addvspace{.55ex}\noindent\textbf{Case 196.} FB.\enspace \ref{case:61}.
\enspace \(\A{4109}^{\rm op}\)\enspace 
3355, 4337, 4626, 4974, 6198, 6697, 6897, 8210. 
\par

\refstepcounter{casegroup}\label{grp:197}\par\addvspace{.55ex}\noindent\textbf{Case 197.} FB.\enspace \ref{case:6}.
\enspace \(\A{872}^{\rm op}\)\enspace 
3456, 3457, 4740, 4888, 4956, 5444, 5465, 6186, 7437, 8138, 11897, 11898, 11907, 12010, 12011, 12108, 12112, 12115, 12128, 12130, 12131, 12133, 13290, 13293, 13302. 
\par

\refstepcounter{casegroup}\label{grp:198}\par\addvspace{.55ex}\noindent\textbf{Case 198.} FB.\enspace \ref{case:111}.
\enspace \(\A{4342}\)\enspace 
3986, 3987, 4342, 5155, 6066, 6495, 6911, 8267, 14525. 
\par

\refstepcounter{casegroup}\label{grp:199}\par\addvspace{.55ex}\noindent\textbf{Case 199.} FB.\enspace \ref{case:112}.
\enspace \(\A{4025}\)\enspace 
4025. 
\enspace \(\A{4025}^{\rm op}\)\enspace 
6187. 
\par

\refstepcounter{casegroup}\label{grp:200}\par\addvspace{.55ex}\noindent\textbf{Case 200.} FB.\enspace \ref{case:94}.
\enspace \(\A{4150}\)\enspace 
4150, 5943, 6025, 6376, 10944. 
\par

\refstepcounter{casegroup}\label{grp:201}\par\addvspace{.55ex}\noindent\textbf{Case 201.} FB.\enspace \ref{case:61}.
\enspace \(\A{4157}\)\enspace 
4157, 4228, 5231, 5232, 5932, 6029, 6269, 6356, 14560. 
\par

\refstepcounter{casegroup}\label{grp:202}\par\addvspace{.55ex}\noindent\textbf{Case 202.} FB.\enspace \ref{case:113}.
\enspace \(\A{4231}\)\enspace 
4231, 6138, 12173. 
\enspace \(\A{4231}^{\rm op}\)\enspace 
4639, 7784, 11951. 
\par

\refstepcounter{casegroup}\label{grp:203}\par\addvspace{.55ex}\noindent\textbf{Case 203.} FB.\enspace \ref{case:61}.
\enspace \(\A{6277}\)\enspace 
4234, 6277, 9266, 11798, 12156. 
\par

\refstepcounter{casegroup}\label{grp:204}\par\addvspace{.55ex}\noindent\textbf{Case 204.} FB.\enspace \ref{case:114}.
\enspace \(\A{4241}\)\enspace 
4237, 4241, 6282, 11777, 11793, 11795, 11796, 11881, 11883, 11885, 11916, 13075, 13076. 
\par

\refstepcounter{casegroup}\label{grp:205}\par\addvspace{.55ex}\noindent\textbf{Case 205.} FB.\enspace \ref{case:72}.
\enspace \(\A{4239}\)\enspace 
4239, 6291, 6370. 
\par

\refstepcounter{casegroup}\label{grp:206}\par\addvspace{.55ex}\noindent\textbf{Case 206.} FB.\enspace \ref{case:61}.
\enspace \(\A{4157}^{\rm op}\)\enspace 
4336, 5126, 6064, 6491, 6898, 7758, 7759, 8211, 14627. 
\par

\refstepcounter{casegroup}\label{grp:207}\par\addvspace{.55ex}\noindent\textbf{Case 207.} FB.\enspace \ref{case:115}.
\enspace \(\A{4345}\)\enspace 
4345, 6503, 6918. 
\enspace \(\A{4345}^{\rm op}\)\enspace 
4503, 6641, 7017. 
\par

\refstepcounter{casegroup}\label{grp:208}\par\addvspace{.55ex}\noindent\textbf{Case 208.} FB.\enspace \ref{case:116}.
\enspace \(\A{4346}\)\enspace 
4346, 6917. 
\enspace \(\A{4346}^{\rm op}\)\enspace 
4504, 7016. 
\par

\refstepcounter{casegroup}\label{grp:209}\par\addvspace{.55ex}\noindent\textbf{Case 209.} FB.\enspace \ref{case:72}.
\enspace \(\A{4348}\)\enspace 
4348, 6505, 6928. 
\par

\refstepcounter{casegroup}\label{grp:210}\par\addvspace{.55ex}\noindent\textbf{Case 210.} FB.\enspace \ref{case:72}.
\enspace \(\A{4349}\)\enspace 
4349, 6927. 
\par

\refstepcounter{casegroup}\label{grp:211}\par\addvspace{.55ex}\noindent\textbf{Case 211.} FB.\enspace \ref{case:34}.
\enspace \(\A{4671}^{\rm op}\)\enspace 
4351, 6507, 6933, 8263, 14253. 
\par

\refstepcounter{casegroup}\label{grp:212}\par\addvspace{.55ex}\noindent\textbf{Case 212.} FB.\enspace \ref{case:111}.
\enspace \(\A{4360}\)\enspace 
4360, 5116, 5117, 5173, 6098, 6577, 6961, 8305, 14552. 
\par

\refstepcounter{casegroup}\label{grp:213}\par\addvspace{.55ex}\noindent\textbf{Case 213.} FB.\enspace \ref{case:40}.
\enspace \(\A{3094}^{\rm op}\)\enspace 
4491, 5355, 5356, 6623, 6995, 8693, 8694, 9482, 9484, 9846. 
\par

\refstepcounter{casegroup}\label{grp:214}\par\addvspace{.55ex}\noindent\textbf{Case 214.} FB.\enspace \ref{case:40}.
\enspace \(\A{4172}^{\rm op}\)\enspace 
4492, 6994, 9845, 12027, 12170. 
\par

\refstepcounter{casegroup}\label{grp:215}\par\addvspace{.55ex}\noindent\textbf{Case 215.} FB.\enspace \ref{case:40}.
\enspace \(\A{4177}^{\rm op}\)\enspace 
4495, 4613, 6999, 7732, 11856, 11857, 11925, 11965--11967, 11973, 12007--12009, 12064, 13138, 13139. 
\par

\refstepcounter{casegroup}\label{grp:216}\par\addvspace{.55ex}\noindent\textbf{Case 216.} FB.\enspace \ref{case:40}.
\enspace \(\A{2155}^{\rm op}\)\enspace 
4497, 5253, 5254, 6637, 7006, 8378, 8379, 9429, 9431, 9615. 
\par

\refstepcounter{casegroup}\label{grp:217}\par\addvspace{.55ex}\noindent\textbf{Case 217.} FB.\enspace \ref{case:61}.
\enspace \(\A{6277}^{\rm op}\)\enspace 
4498, 7005, 9614, 12018, 12164. 
\par

\refstepcounter{casegroup}\label{grp:218}\par\addvspace{.55ex}\noindent\textbf{Case 218.} FB.\enspace \ref{case:114}.
\enspace \(\A{4241}^{\rm op}\)\enspace 
4501, 4667, 7010, 11835, 11912, 11931, 11933, 11934, 11936, 11980, 11982, 13126, 13127. 
\par

\refstepcounter{casegroup}\label{grp:219}\par\addvspace{.55ex}\noindent\textbf{Case 219.} FB.\enspace \ref{case:117}.
\enspace \(\A{7037}\)\enspace 
4516, 4519, 5268, 5269, 7037, 7043, 8405, 9400, 9407, 9596, 10182, 11844, 11917, 11938, 11943, 11946, 11949, 11991, 11993, 12023, 12166, 13130, 13131. 
\par

\refstepcounter{casegroup}\label{grp:220}\par\addvspace{.55ex}\noindent\textbf{Case 220.} FB.\enspace \ref{case:118}.
\enspace \(\A{4521}\)\enspace 
4521, 7051. 
\enspace \(\A{4522}\)\enspace 
4522, 7050. 
\enspace \(\A{4521}^{\rm op}\)\enspace 
4705, 7912. 
\enspace \(\A{4522}^{\rm op}\)\enspace 
4706, 7911. 
\par

\refstepcounter{casegroup}\label{grp:221}\par\addvspace{.55ex}\noindent\textbf{Case 221.} FB.\enspace \ref{case:119}.
4524, 5256, 5314, 7242, 7498, 7556. 
\par

\refstepcounter{casegroup}\label{grp:222}\par\addvspace{.55ex}\noindent\textbf{Case 222.} FB.\enspace \ref{case:120}.
\enspace \(\A{4525}\)\enspace 
4525, 7056, 11939. 
\enspace \(\A{4525}^{\rm op}\)\enspace 
5257, 8385, 11942. 
\par

\refstepcounter{casegroup}\label{grp:223}\par\addvspace{.55ex}\noindent\textbf{Case 223.} FB.\enspace Lemma~\ref{lem:terms}, source \(\A{7037}\), case~\ref{grp:219}.
\enspace \(\A{7059}\leftrightarrow\A{7037}\), \(x+y, xy+yx, x+y, xxy\)\enspace 
4528, 4531, 5274, 5275, 7059, 7064, 8414, 9432, 9434, 9622, 10292, 11847, 11920, 11940, 11944, 11953, 11955, 11995, 11997, 12024, 12167, 13133, 13134. 
\enspace \(\A{8388}\leftrightarrow\A{7037}\), \(x+y, xy+yx, x+y, xyy\)\enspace 
5260, 5263, 5271, 5272, 8388, 8393, 8410, 9401, 9409, 9601, 10187, 12019, 12020, 12149, 12150, 12162, 12175, 12177, 12253--12255, 13336, 13337. 
\par

\refstepcounter{casegroup}\label{grp:224}\par\addvspace{.55ex}\noindent\textbf{Case 224.} FB.\enspace \ref{case:121}.
\enspace \(\A{4533}\)\enspace 
4533, 7071. 
\enspace \(\A{4534}\)\enspace 
4534, 7070. 
\enspace \(\A{4533}^{\rm op}\)\enspace 
5340, 8670. 
\enspace \(\A{4534}^{\rm op}\)\enspace 
5341, 8669. 
\par

\refstepcounter{casegroup}\label{grp:225}\par\addvspace{.55ex}\noindent\textbf{Case 225.} FB.\enspace \ref{case:117}.
\enspace \(\A{4567}\)\enspace 
4567, 4570, 4851, 4853, 4857, 5303, 7040, 7202, 7497, 7557, 7589, 7617, 7621, 7625, 7648, 7915, 7917, 8085, 8648, 9415, 9593, 9809. 
\par

\refstepcounter{casegroup}\label{grp:226}\par\addvspace{.55ex}\noindent\textbf{Case 226.} FB.\enspace \ref{case:122}.
\enspace \(\A{4572}\)\enspace 
4572, 4893, 7048. 
\enspace \(\A{4572}^{\rm op}\)\enspace 
4849, 7374, 7909. 
\par

\refstepcounter{casegroup}\label{grp:227}\par\addvspace{.55ex}\noindent\textbf{Case 227.} FB.\enspace \ref{case:123}.
\enspace \(\A{4573}\)\enspace 
4573, 4922, 7047. 
\enspace \(\A{4575}\)\enspace 
4575, 7053, 7500. 
\enspace \(\A{4573}^{\rm op}\)\enspace 
4850, 7407, 7612. 
\enspace \(\A{4575}^{\rm op}\)\enspace 
5295, 7587, 7611. 
\par

\refstepcounter{casegroup}\label{grp:228}\par\addvspace{.55ex}\noindent\textbf{Case 228.} FB.\enspace Lemma~\ref{lem:terms}, source \(\A{4567}\), case~\ref{grp:225}.
\enspace \(\A{4577}\leftrightarrow\A{4567}\), \(x+y, xy+yx, x+y, xxy\)\enspace 
4577, 4579, 4868, 4869, 4874, 5307, 7062, 7212, 7502, 7569, 7591, 7669, 7679, 7685, 7699, 7947, 7948, 8090, 8651, 9456, 9619, 9818. 
\enspace \(\A{5297}\leftrightarrow\A{4567}\), \(x+y, xy+yx, x+y, xyy\)\enspace 
5297, 5299, 5305, 5385, 5386, 5391, 8051, 8054, 8061, 8391, 8408, 8503, 8604, 8609, 8675, 8676, 8787, 8789, 8798, 9417, 9598, 9948. 
\par

\refstepcounter{casegroup}\label{grp:229}\par\addvspace{.55ex}\noindent\textbf{Case 229.} FB.\enspace \ref{case:124}.
\enspace \(\A{4581}\)\enspace 
4581, 4931, 7068. 
\enspace \(\A{4867}\)\enspace 
4867, 7424, 7670. 
\enspace \(\A{4867}^{\rm op}\)\enspace 
5301, 5411, 8398. 
\enspace \(\A{4581}^{\rm op}\)\enspace 
5384, 8034, 8598. 
\par

\refstepcounter{casegroup}\label{grp:230}\par\addvspace{.55ex}\noindent\textbf{Case 230.} FB.\enspace \ref{case:125}.
\enspace \(\A{4588}\)\enspace 
4588--4590, 5322, 7046, 7126, 7265, 7283, 7293, 7302, 8016, 8568, 9418, 9746. 
\par

\refstepcounter{casegroup}\label{grp:231}\par\addvspace{.55ex}\noindent\textbf{Case 231.} FB.\enspace \ref{case:126}.
\enspace \(\A{4591}\)\enspace 
4591, 7268. 
\enspace \(\A{4591}^{\rm op}\)\enspace 
7586, 7981. 
\par

\refstepcounter{casegroup}\label{grp:232}\par\addvspace{.55ex}\noindent\textbf{Case 232.} FB.\enspace \ref{s4592-thm:main}.
\enspace \(\A{4592}\)\enspace 
4592--4594, 5324, 7067, 7138, 7272, 7287, 7315, 7327, 8024, 8577, 9465, 9753. 
\par

\refstepcounter{casegroup}\label{grp:233}\par\addvspace{.55ex}\noindent\textbf{Case 233.} FB.\enspace \ref{s4609-thm:basis}.
\enspace \(\A{4609}\)\enspace 
4609, 5338, 6155, 8032, 8093, 11855, 11962, 11972, 12003, 12056, 12062, 12063, 13137, 13141. 
\par

\refstepcounter{casegroup}\label{grp:234}\par\addvspace{.55ex}\noindent\textbf{Case 234.} FB.\enspace \ref{case:72}.
\enspace \(\A{13828}^{\rm op}\)\enspace 
4616, 11864, 12005. 
\par

\refstepcounter{casegroup}\label{grp:235}\par\addvspace{.55ex}\noindent\textbf{Case 235.} FB.\enspace \ref{case:72}.
\enspace \(\A{13057}^{\rm op}\)\enspace 
4627, 11863, 11926, 12006. 
\par

\refstepcounter{casegroup}\label{grp:236}\par\addvspace{.55ex}\noindent\textbf{Case 236.} FB.\enspace \ref{case:127}.
\enspace \(\A{10590}\)\enspace 
4640, 7788, 10590. 
\par

\refstepcounter{casegroup}\label{grp:237}\par\addvspace{.55ex}\noindent\textbf{Case 237.} FB.\enspace \ref{case:128}.
\enspace \(\A{4641}\)\enspace 
4641, 7793, 11989. 
\par

\refstepcounter{casegroup}\label{grp:238}\par\addvspace{.55ex}\noindent\textbf{Case 238.} FB.\enspace \ref{case:129}.
\enspace \(\A{4649}\)\enspace 
4649. 
\enspace \(\A{4649}^{\rm op}\)\enspace 
5004. 
\par

\refstepcounter{casegroup}\label{grp:239}\par\addvspace{.55ex}\noindent\textbf{Case 239.} FB.\enspace \ref{case:72}.
\enspace \(\A{4239}^{\rm op}\)\enspace 
4651, 6751, 7819. 
\par

\refstepcounter{casegroup}\label{grp:240}\par\addvspace{.55ex}\noindent\textbf{Case 240.} FB.\enspace \ref{case:72}.
\enspace \(\A{4348}^{\rm op}\)\enspace 
4654, 6753, 7825. 
\par

\refstepcounter{casegroup}\label{grp:241}\par\addvspace{.55ex}\noindent\textbf{Case 241.} FB.\enspace \ref{case:72}.
\enspace \(\A{4349}^{\rm op}\)\enspace 
4655, 7824. 
\par

\refstepcounter{casegroup}\label{grp:242}\par\addvspace{.55ex}\noindent\textbf{Case 242.} FB.\enspace \ref{case:130}.
\enspace \(\A{4665}\)\enspace 
4665, 4666, 6759, 7845. 
\enspace \(\A{4665}^{\rm op}\)\enspace 
4679, 4680, 6779, 7869. 
\par

\refstepcounter{casegroup}\label{grp:243}\par\addvspace{.55ex}\noindent\textbf{Case 243.} FB.\enspace \ref{case:34}.
\enspace \(\A{4671}\)\enspace 
4671, 6774, 7852, 8114, 11867. 
\par

\refstepcounter{casegroup}\label{grp:244}\par\addvspace{.55ex}\noindent\textbf{Case 244.} FB.\enspace \ref{case:131}.
\enspace \(\A{4677}\)\enspace 
4677, 11789, 12070. 
\par

\refstepcounter{casegroup}\label{grp:245}\par\addvspace{.55ex}\noindent\textbf{Case 245.} FB.\enspace \ref{case:132}.
4702, 11775, 11873, 12072. 
\par

\refstepcounter{casegroup}\label{grp:246}\par\addvspace{.55ex}\noindent\textbf{Case 246.} FB.\enspace \ref{case:133}.
\enspace \(\A{4703}\)\enspace 
4703, 7908, 10495. 
\par

\refstepcounter{casegroup}\label{grp:247}\par\addvspace{.55ex}\noindent\textbf{Case 247.} FB.\enspace \ref{case:134}.
\enspace \(\A{4707}\)\enspace 
4707, 6824, 7200, 7920. 
\par

\refstepcounter{casegroup}\label{grp:248}\par\addvspace{.55ex}\noindent\textbf{Case 248.} FB.\enspace \ref{case:128}.
\enspace \(\A{4708}\)\enspace 
4708, 7919. 
\par

\refstepcounter{casegroup}\label{grp:249}\par\addvspace{.55ex}\noindent\textbf{Case 249.} FB.\enspace \ref{case:134}.
\enspace \(\A{4709}\)\enspace 
4709, 7918. 
\par

\refstepcounter{casegroup}\label{grp:250}\par\addvspace{.55ex}\noindent\textbf{Case 250.} FB.\enspace \ref{case:134}.
\enspace \(\A{4710}\)\enspace 
4710, 4711, 6829, 7205, 7926, 8836, 8869, 8877, 8897, 8933, 8957, 10207. 
\par

\refstepcounter{casegroup}\label{grp:251}\par\addvspace{.55ex}\noindent\textbf{Case 251.} FB.\enspace \ref{case:128}.
\enspace \(\A{4712}\)\enspace 
4712, 8844, 8895. 
\par

\refstepcounter{casegroup}\label{grp:252}\par\addvspace{.55ex}\noindent\textbf{Case 252.} FB.\enspace \ref{case:134}.
\enspace \(\A{4713}\)\enspace 
4713, 11918, 11990. 
\par

\refstepcounter{casegroup}\label{grp:253}\par\addvspace{.55ex}\noindent\textbf{Case 253.} FB.\enspace \ref{case:135}.
\enspace \(\A{4721}\)\enspace 
4721, 7943. 
\enspace \(\A{4722}\)\enspace 
4722, 7942. 
\enspace \(\A{4721}^{\rm op}\)\enspace 
5265, 8401. 
\enspace \(\A{4722}^{\rm op}\)\enspace 
5266, 8400. 
\par

\refstepcounter{casegroup}\label{grp:254}\par\addvspace{.55ex}\noindent\textbf{Case 254.} FB.\enspace Lemma~\ref{lem:terms}, source \(\A{4707}\), case~\ref{grp:247}.
\enspace \(\A{4723}\leftrightarrow\A{4707}\), \(x+y, xy+yx, x+y, xxy\)\enspace 
4723, 6834, 7209, 7951. 
\enspace \(\A{5342}\leftrightarrow\A{4707}\), \(x+y, xy+yx, x+y, xyy\)\enspace 
5342, 8499, 8501, 8679. 
\par

\refstepcounter{casegroup}\label{grp:255}\par\addvspace{.55ex}\noindent\textbf{Case 255.} FB.\enspace Lemma~\ref{lem:terms}, source \(\A{4708}\), case~\ref{grp:248}.
\enspace \(\A{4724}\leftrightarrow\A{4708}\), \(x+y, xy+yx, x+y, xxy\)\enspace 
4724, 7950. 
\enspace \(\A{5343}\leftrightarrow\A{4708}\), \(x+y, xy+yx, x+y, xyy\)\enspace 
5343, 8678. 
\par

\refstepcounter{casegroup}\label{grp:256}\par\addvspace{.55ex}\noindent\textbf{Case 256.} FB.\enspace Lemma~\ref{lem:terms}, source \(\A{4709}\), case~\ref{grp:249}.
\enspace \(\A{4725}\leftrightarrow\A{4709}\), \(x+y, xy+yx, x+y, xxy\)\enspace 
4725, 7949. 
\enspace \(\A{5344}\leftrightarrow\A{4709}\), \(x+y, xy+yx, x+y, xyy\)\enspace 
5344, 8677. 
\par

\refstepcounter{casegroup}\label{grp:257}\par\addvspace{.55ex}\noindent\textbf{Case 257.} FB.\enspace Lemma~\ref{lem:terms}, source \(\A{4710}\), case~\ref{grp:250}.
\enspace \(\A{4726}\leftrightarrow\A{4710}\), \(x+y, xy+yx, x+y, xxy\)\enspace 
4726, 4730, 6838, 7217, 7961, 8994, 9012, 9015, 9073, 9159, 9188, 10333. 
\enspace \(\A{5345}\leftrightarrow\A{4710}\), \(x+y, xy+yx, x+y, xyy\)\enspace 
5345, 5349, 8507, 8509, 8691, 8855, 8874, 8882, 8901, 8935, 8958, 10210. 
\par

\refstepcounter{casegroup}\label{grp:258}\par\addvspace{.55ex}\noindent\textbf{Case 258.} FB.\enspace \ref{case:136}.
\enspace \(\A{4728}\)\enspace 
4728. 
\enspace \(\A{4728}^{\rm op}\)\enspace 
5347. 
\par

\refstepcounter{casegroup}\label{grp:259}\par\addvspace{.55ex}\noindent\textbf{Case 259.} FB.\enspace \ref{case:137}.
\enspace \(\A{4729}\)\enspace 
4729. 
\enspace \(\A{4729}^{\rm op}\)\enspace 
5348. 
\par

\refstepcounter{casegroup}\label{grp:260}\par\addvspace{.55ex}\noindent\textbf{Case 260.} FB.\enspace Lemma~\ref{lem:terms}, source \(\A{4712}\), case~\ref{grp:251}.
\enspace \(\A{4731}\leftrightarrow\A{4712}\), \(x+y, xy+yx, x+y, xxy\)\enspace 
4731, 9001, 9157. 
\enspace \(\A{5350}\leftrightarrow\A{4712}\), \(x+y, xy+yx, x+y, xyy\)\enspace 
5350, 8862, 8899. 
\par

\refstepcounter{casegroup}\label{grp:261}\par\addvspace{.55ex}\noindent\textbf{Case 261.} FB.\enspace Lemma~\ref{lem:terms}, source \(\A{4713}\), case~\ref{grp:252}.
\enspace \(\A{4732}\leftrightarrow\A{4713}\), \(x+y, xy+yx, x+y, xxy\)\enspace 
4732, 11921, 11994. 
\enspace \(\A{5351}\leftrightarrow\A{4713}\), \(x+y, xy+yx, x+y, xyy\)\enspace 
5351, 12151, 12174. 
\par

\refstepcounter{casegroup}\label{grp:262}\par\addvspace{.55ex}\noindent\textbf{Case 262.} FB.\enspace \ref{case:138}.
\enspace \(\A{5519}^{\rm op}\)\enspace 
4734. 
\enspace \(\A{5519}\)\enspace 
5519. 
\par

\refstepcounter{casegroup}\label{grp:263}\par\addvspace{.55ex}\noindent\textbf{Case 263.} FB.\enspace \ref{case:139}.
\enspace \(\A{5526}^{\rm op}\)\enspace 
4742. 
\enspace \(\A{5526}\)\enspace 
5526. 
\par

\refstepcounter{casegroup}\label{grp:264}\par\addvspace{.55ex}\noindent\textbf{Case 264.} NFB.\enspace \ref{thm:pairnfb}.
\enspace \(\A{4793}\)\enspace 
4793, 5072, 5073, 14027, 14533. 
\enspace \(\A{4793}^{\rm op}\)\enspace 
5160, 5167, 5168, 11005, 12482. 
\par

\refstepcounter{casegroup}\label{grp:265}\par\addvspace{.55ex}\noindent\textbf{Case 265.} FB.\enspace \ref{case:140}.
\enspace \(\A{4872}\)\enspace 
4872. 
\enspace \(\A{4872}^{\rm op}\)\enspace 
5389. 
\par

\refstepcounter{casegroup}\label{grp:266}\par\addvspace{.55ex}\noindent\textbf{Case 266.} FB.\enspace \ref{case:141}.
\enspace \(\A{4901}\)\enspace 
4901. 
\enspace \(\A{4901}^{\rm op}\)\enspace 
7750. 
\par

\refstepcounter{casegroup}\label{grp:267}\par\addvspace{.55ex}\noindent\textbf{Case 267.} FB.\enspace \ref{case:142}.
\enspace \(\A{4933}\)\enspace 
4933. 
\enspace \(\A{4933}^{\rm op}\)\enspace 
8036. 
\par

\refstepcounter{casegroup}\label{grp:268}\par\addvspace{.55ex}\noindent\textbf{Case 268.} FB.\enspace \ref{case:143}.
\enspace \(\A{4942}\)\enspace 
4942. 
\enspace \(\A{4942}^{\rm op}\)\enspace 
8153. 
\par

\refstepcounter{casegroup}\label{grp:269}\par\addvspace{.55ex}\noindent\textbf{Case 269.} FB.\enspace \ref{case:144}.
\enspace \(\A{5001}\)\enspace 
5001. 
\enspace \(\A{5001}^{\rm op}\)\enspace 
5013. 
\par

\refstepcounter{casegroup}\label{grp:270}\par\addvspace{.55ex}\noindent\textbf{Case 270.} FB.\enspace \ref{case:145}.
\enspace \(\A{5002}\)\enspace 
5002. 
\enspace \(\A{5002}^{\rm op}\)\enspace 
5014. 
\par

\refstepcounter{casegroup}\label{grp:271}\par\addvspace{.55ex}\noindent\textbf{Case 271.} FB.\enspace \ref{case:132}.
5008, 11779, 11874, 12073. 
\par

\refstepcounter{casegroup}\label{grp:272}\par\addvspace{.55ex}\noindent\textbf{Case 272.} FB.\enspace \ref{case:94}.
\enspace \(\A{4150}^{\rm op}\)\enspace 
5047, 6067, 6497, 8208, 14077. 
\par

\refstepcounter{casegroup}\label{grp:273}\par\addvspace{.55ex}\noindent\textbf{Case 273.} FB.\enspace \ref{case:93}.
\enspace \(\A{5136}^{\rm op}\)\enspace 
5065, 6068, 6499, 8241, 14025. 
\par

\refstepcounter{casegroup}\label{grp:274}\par\addvspace{.55ex}\noindent\textbf{Case 274.} FB.\enspace \ref{case:34}.
\enspace \(\A{5165}^{\rm op}\)\enspace 
5070, 6069, 6501, 8264, 14006. 
\par

\refstepcounter{casegroup}\label{grp:275}\par\addvspace{.55ex}\noindent\textbf{Case 275.} FB.\enspace \ref{case:146}.
\enspace \(\A{5081}\)\enspace 
5081, 6100, 6581, 6948, 8302, 14023. 
\par

\refstepcounter{casegroup}\label{grp:276}\par\addvspace{.55ex}\noindent\textbf{Case 276.} FB.\enspace \ref{case:93}.
\enspace \(\A{5136}\)\enspace 
5136, 6071, 6511, 8240, 10943. 
\par

\refstepcounter{casegroup}\label{grp:277}\par\addvspace{.55ex}\noindent\textbf{Case 277.} FB.\enspace \ref{case:146}.
\enspace \(\A{5141}\)\enspace 
5141, 6072, 6513, 6907, 8261, 10926. 
\par

\refstepcounter{casegroup}\label{grp:278}\par\addvspace{.55ex}\noindent\textbf{Case 278.} FB.\enspace \ref{case:147}.
\enspace \(\A{5157}\)\enspace 
5157, 14550. 
\enspace \(\A{5157}^{\rm op}\)\enspace 
5175, 13857. 
\par

\refstepcounter{casegroup}\label{grp:279}\par\addvspace{.55ex}\noindent\textbf{Case 279.} FB.\enspace \ref{case:148}.
\enspace \(\A{5158}\)\enspace 
5158. 
\enspace \(\A{5158}^{\rm op}\)\enspace 
5176. 
\par

\refstepcounter{casegroup}\label{grp:280}\par\addvspace{.55ex}\noindent\textbf{Case 280.} FB.\enspace \ref{case:34}.
\enspace \(\A{5165}\)\enspace 
5165, 6102, 6587, 8300, 10941. 
\par

\refstepcounter{casegroup}\label{grp:281}\par\addvspace{.55ex}\noindent\textbf{Case 281.} FB.\enspace \ref{case:149}.
\enspace \(\A{5219}\)\enspace 
5219, 14559. 
\enspace \(\A{5219}^{\rm op}\)\enspace 
7754, 13866. 
\par

\refstepcounter{casegroup}\label{grp:282}\par\addvspace{.55ex}\noindent\textbf{Case 282.} FB.\enspace \ref{case:150}.
\enspace \(\A{5220}\)\enspace 
5220. 
\enspace \(\A{5220}^{\rm op}\)\enspace 
7755. 
\par

\refstepcounter{casegroup}\label{grp:283}\par\addvspace{.55ex}\noindent\textbf{Case 283.} FB.\enspace \ref{case:151}.
\enspace \(\A{5222}\)\enspace 
5222, 11009. 
\enspace \(\A{5222}^{\rm op}\)\enspace 
7383, 14078. 
\par

\refstepcounter{casegroup}\label{grp:284}\par\addvspace{.55ex}\noindent\textbf{Case 284.} FB.\enspace \ref{case:152}.
\enspace \(\A{5223}\)\enspace 
5223. 
\enspace \(\A{5228}\)\enspace 
5228. 
\enspace \(\A{5223}^{\rm op}\)\enspace 
7384. 
\enspace \(\A{5228}^{\rm op}\)\enspace 
7748. 
\par

\refstepcounter{casegroup}\label{grp:285}\par\addvspace{.55ex}\noindent\textbf{Case 285.} FB.\enspace \ref{case:153}.
\enspace \(\A{5227}\)\enspace 
5227, 13861. 
\enspace \(\A{5227}^{\rm op}\)\enspace 
7747, 14624. 
\par

\refstepcounter{casegroup}\label{grp:286}\par\addvspace{.55ex}\noindent\textbf{Case 286.} FB.\enspace \ref{case:68}.
5289, 5290, 5379, 5380, 5388, 8063, 8067, 8381, 8510, 8605, 8612, 8615, 8684, 8792, 8794, 8799, 9413, 9592, 9949. 
\par

\refstepcounter{casegroup}\label{grp:287}\par\addvspace{.55ex}\noindent\textbf{Case 287.} FB.\enspace \ref{case:154}.
\enspace \(\A{5315}\)\enspace 
5315, 7278. 
\enspace \(\A{5315}^{\rm op}\)\enspace 
7243, 7499. 
\par

\refstepcounter{casegroup}\label{grp:288}\par\addvspace{.55ex}\noindent\textbf{Case 288.} FB.\enspace \ref{abs-chain-basis}.
\enspace \(\text{\normalfont K}^{\rm op}\)\enspace 
5317, 5319, 5321, 5323, 7996, 8000, 8396, 8460, 8545, 8737, 8746, 8754, 9419, 9922. 
\enspace \(\text{\normalfont K}\)\enspace 
7245, 7247, 7422, 7501, 7590, 7663, 7676, 7678, 7697, 7986, 8089, 8650, 9728, 9817. 
\par

\refstepcounter{casegroup}\label{grp:289}\par\addvspace{.55ex}\noindent\textbf{Case 289.} FB.\enspace Proposition~\ref{prop:lower-fibres}, source \(\A{8490}\), case~\ref{grp:2}.
\enspace \(\A{8490}\)\enspace 
5328, 8489, 8657. 
\par

\refstepcounter{casegroup}\label{grp:290}\par\addvspace{.55ex}\noindent\textbf{Case 290.} FB.\enspace \ref{case:155}.
\enspace \(\A{5410}\)\enspace 
5410, 8672. 
\enspace \(\A{5410}^{\rm op}\)\enspace 
7423, 7671. 
\par

\refstepcounter{casegroup}\label{grp:291}\par\addvspace{.55ex}\noindent\textbf{Case 291.} FB.\enspace \ref{case:156}.
\enspace \(\A{5413}\)\enspace 
5413. 
\enspace \(\A{5414}\)\enspace 
5414. 
\enspace \(\A{5413}^{\rm op}\)\enspace 
7426. 
\enspace \(\A{5414}^{\rm op}\)\enspace 
7427. 
\par

\refstepcounter{casegroup}\label{grp:292}\par\addvspace{.55ex}\noindent\textbf{Case 292.} NFB.\enspace \ref{case:157}.
5467, 5482, 12637, 14535. 
\par

\refstepcounter{casegroup}\label{grp:293}\par\addvspace{.55ex}\noindent\textbf{Case 293.} FB.\enspace \ref{case:8}.
\enspace \(\A{5767}^{\rm op}\)\enspace 
5473. 
\par

\refstepcounter{casegroup}\label{grp:294}\par\addvspace{.55ex}\noindent\textbf{Case 294.} FB.\enspace \ref{case:8}.
\enspace \(\A{5767}\)\enspace 
5767. 
\par

\refstepcounter{casegroup}\label{grp:295}\par\addvspace{.55ex}\noindent\textbf{Case 295.} FB.\enspace \ref{case:158}.
\enspace \(\A{5832}\)\enspace 
5832. 
\par

\refstepcounter{casegroup}\label{grp:296}\par\addvspace{.55ex}\noindent\textbf{Case 296.} FB.\enspace \ref{case:159}.
\enspace \(\A{9069}\)\enspace 
5927, 6346, 9069, 9247. 
\par

\refstepcounter{casegroup}\label{grp:297}\par\addvspace{.55ex}\noindent\textbf{Case 297.} FB.\enspace \ref{case:160}.
\enspace \(\A{5929}\)\enspace 
5929, 6027, 6350. 
\enspace \(\A{5929}^{\rm op}\)\enspace 
6061, 6485, 8209. 
\par

\refstepcounter{casegroup}\label{grp:298}\par\addvspace{.55ex}\noindent\textbf{Case 298.} FB.\enspace \ref{case:161}.
\enspace \(\A{5937}\)\enspace 
5937, 6026, 6363, 14032. 
\enspace \(\A{5937}^{\rm op}\)\enspace 
6070, 6509, 8207, 10954. 
\par

\refstepcounter{casegroup}\label{grp:299}\par\addvspace{.55ex}\noindent\textbf{Case 299.} FB.\enspace \ref{case:162}.
\enspace \(\A{6063}\)\enspace 
6063, 6489, 8265. 
\enspace \(\A{6063}^{\rm op}\)\enspace 
6096, 6573, 8303. 
\par

\refstepcounter{casegroup}\label{grp:300}\par\addvspace{.55ex}\noindent\textbf{Case 300.} FB.\enspace \ref{case:72}.
\enspace \(\A{1888}^{\rm op}\)\enspace 
6125. 
\par

\refstepcounter{casegroup}\label{grp:301}\par\addvspace{.55ex}\noindent\textbf{Case 301.} FB.\enspace \ref{case:128}.
\enspace \(\A{6160}\)\enspace 
6160, 6729, 7794, 8868. 
\par

\refstepcounter{casegroup}\label{grp:302}\par\addvspace{.55ex}\noindent\textbf{Case 302.} NFB.\enspace \ref{case:164}, \ref{case:163}.
6203, 6209, 12635, 13203. 
\par

\refstepcounter{casegroup}\label{grp:303}\par\addvspace{.55ex}\noindent\textbf{Case 303.} FB.\enspace \ref{case:165}.
\enspace \(\A{6211}\)\enspace 
6211. 
\enspace \(\A{6211}^{\rm op}\)\enspace 
6214. 
\par

\refstepcounter{casegroup}\label{grp:304}\par\addvspace{.55ex}\noindent\textbf{Case 304.} FB.\enspace \ref{case:166}.
6562, 7592. 
\par

\refstepcounter{casegroup}\label{grp:305}\par\addvspace{.55ex}\noindent\textbf{Case 305.} FB.\enspace \ref{case:167}.
\enspace \(\A{6563}\)\enspace 
6563. 
\enspace \(\A{6563}^{\rm op}\)\enspace 
7275. 
\par

\refstepcounter{casegroup}\label{grp:306}\par\addvspace{.55ex}\noindent\textbf{Case 306.} FB.\enspace \ref{case:168}.
\enspace \(\A{6564}\)\enspace 
6564. 
\enspace \(\A{6564}^{\rm op}\)\enspace 
7585. 
\par

\refstepcounter{casegroup}\label{grp:307}\par\addvspace{.55ex}\noindent\textbf{Case 307.} FB.\enspace \ref{thm:power-blocks}.
6664, 6666, 6669, 6670, 6673, 6674, 6676, 6678, 7087, 7297, 7298, 7320, 8725--8728, 8730, 8755--8758, 8771, 10084, 10087--10089, 10401, 10402, 10404, 10410, 10447, 10498--10500, 10603, 10604, 10606--10611, 10613, 10614, 10616, 10617, 10619, 10620, 10622--10627, 10629, 10630, 10632, 10633, 10635, 10636, 10638--10642, 10644--10659, 10665, 10666, 10673, 10674, 10681, 10687, 10693--10702, 10705, 10727, 10760, 10761, 10763, 10764, 10766--10768, 10774, 10775. 
\par

\refstepcounter{casegroup}\label{grp:308}\par\addvspace{.55ex}\noindent\textbf{Case 308.} FB.\enspace \ref{case:169}.
\enspace \(\A{6668}\)\enspace 
6668. 
\enspace \(\A{6668}^{\rm op}\)\enspace 
6823. 
\par

\refstepcounter{casegroup}\label{grp:309}\par\addvspace{.55ex}\noindent\textbf{Case 309.} FB.\enspace \ref{case:170}.
\enspace \(\A{6671}\)\enspace 
6671. 
\enspace \(\A{6827}\)\enspace 
6827. 
\enspace \(\A{7086}\)\enspace 
7086. 
\enspace \(\A{7203}\)\enspace 
7203. 
\enspace \(\A{7203}^{\rm op}\)\enspace 
7290. 
\enspace \(\A{7296}\)\enspace 
7296. 
\enspace \(\A{6827}^{\rm op}\)\enspace 
7299. 
\enspace \(\A{7086}^{\rm op}\)\enspace 
7614. 
\enspace \(\A{7296}^{\rm op}\)\enspace 
7633. 
\enspace \(\A{6671}^{\rm op}\)\enspace 
7641. 
\enspace \(\A{7693}\)\enspace 
7693. 
\enspace \(\A{7693}^{\rm op}\)\enspace 
8549. 
\enspace \(\A{8905}\)\enspace 
8905. 
\enspace \(\A{8905}^{\rm op}\)\enspace 
9584. 
\enspace \(\A{10046}\)\enspace 
10046. 
\enspace \(\A{10046}^{\rm op}\)\enspace 
10058. 
\enspace \(\A{10130}\)\enspace 
10130. 
\enspace \(\A{10131}\)\enspace 
10131. 
\enspace \(\A{10130}^{\rm op}\)\enspace 
10141. 
\enspace \(\A{10131}^{\rm op}\)\enspace 
10142. 
\enspace \(\A{10151}\)\enspace 
10151. 
\enspace \(\A{10155}\)\enspace 
10155. 
\enspace \(\A{10151}^{\rm op}\)\enspace 
10165. 
\enspace \(\A{10155}^{\rm op}\)\enspace 
10177. 
\enspace \(\A{10213}\)\enspace 
10213. 
\enspace \(\A{10214}\)\enspace 
10214. 
\enspace \(\A{10232}\)\enspace 
10232. 
\enspace \(\A{10213}^{\rm op}\)\enspace 
10396. 
\enspace \(\A{10232}^{\rm op}\)\enspace 
10397. 
\enspace \(\A{10214}^{\rm op}\)\enspace 
10405. 
\enspace \(\A{11846}\)\enspace 
11846. 
\enspace \(\A{11919}\)\enspace 
11919. 
\enspace \(\A{11919}^{\rm op}\)\enspace 
11945. 
\enspace \(\A{11846}^{\rm op}\)\enspace 
11948. 
\par

\refstepcounter{casegroup}\label{grp:310}\par\addvspace{.55ex}\noindent\textbf{Case 310.} FB.\enspace \ref{case:171}.
\enspace \(\A{6680}\)\enspace 
6680. 
\enspace \(\A{6680}^{\rm op}\)\enspace 
8783. 
\par

\refstepcounter{casegroup}\label{grp:311}\par\addvspace{.55ex}\noindent\textbf{Case 311.} FB.\enspace \ref{case:172}.
\enspace \(\A{6743}\)\enspace 
6743. 
\enspace \(\A{6743}^{\rm op}\)\enspace 
7081. 
\enspace \(\A{7214}\)\enspace 
7214. 
\enspace \(\A{7312}\)\enspace 
7312. 
\enspace \(\A{8738}^{\rm op}\)\enspace 
8088. 
\enspace \(\A{7312}^{\rm op}\)\enspace 
8504. 
\enspace \(\A{7214}^{\rm op}\)\enspace 
8541. 
\enspace \(\A{8738}\)\enspace 
8738. 
\enspace \(\A{9255}\)\enspace 
9255. 
\enspace \(\A{9303}\)\enspace 
9303. 
\enspace \(\A{9255}^{\rm op}\)\enspace 
9588. 
\enspace \(\A{9303}^{\rm op}\)\enspace 
9890. 
\enspace \(\A{10355}\)\enspace 
10355. 
\enspace \(\A{10355}^{\rm op}\)\enspace 
10399. 
\enspace \(\A{11922}\)\enspace 
11922. 
\enspace \(\A{12929}^{\rm op}\)\enspace 
11958. 
\enspace \(\A{11922}^{\rm op}\)\enspace 
12161. 
\enspace \(\A{12929}\)\enspace 
12929. 
\par

\refstepcounter{casegroup}\label{grp:312}\par\addvspace{.55ex}\noindent\textbf{Case 312.} FB.\enspace \ref{case:173}.
6745, 6826, 6828, 7204, 7454, 7467, 7613, 7619, 7642, 7643, 7790, 7914. 
\par

\refstepcounter{casegroup}\label{grp:313}\par\addvspace{.55ex}\noindent\textbf{Case 313.} FB.\enspace \ref{case:174}.
6748, 7215, 7821, 8602, 11232, 11236, 13052, 13316. 
\par

\refstepcounter{casegroup}\label{grp:314}\par\addvspace{.55ex}\noindent\textbf{Case 314.} FB.\enspace \ref{case:175}.
\enspace \(\A{6825}\)\enspace 
6825. 
\par

\refstepcounter{casegroup}\label{grp:315}\par\addvspace{.55ex}\noindent\textbf{Case 315.} FB.\enspace \ref{case:176}.
\enspace \(\A{8729}^{\rm op}\)\enspace 
6833. 
\enspace \(\A{8729}\)\enspace 
8729. 
\par

\refstepcounter{casegroup}\label{grp:316}\par\addvspace{.55ex}\noindent\textbf{Case 316.} FB.\enspace Lemma~\ref{lem:terms}, source \(\A{6825}\), case~\ref{grp:314}.
\enspace \(\A{6835}\leftrightarrow\A{6825}\), \(x+y, xy+yx, x+y, xxy\)\enspace 
6835. 
\enspace \(\A{8784}\leftrightarrow\A{6825}\), \(x+y, xy+yx, x+y, xyy\)\enspace 
8784. 
\par

\refstepcounter{casegroup}\label{grp:317}\par\addvspace{.55ex}\noindent\textbf{Case 317.} FB.\enspace \ref{case:177}.
6837, 8790, 10999, 11013, 11298, 14018, 14033, 14093. 
\par

\refstepcounter{casegroup}\label{grp:318}\par\addvspace{.55ex}\noindent\textbf{Case 318.} FB.\enspace \ref{case:132}.
6921. 
\par

\refstepcounter{casegroup}\label{grp:319}\par\addvspace{.55ex}\noindent\textbf{Case 319.} FB.\enspace \ref{case:178}.
7054, 7057, 7129, 8383, 8386, 8538. 
\par

\refstepcounter{casegroup}\label{grp:320}\par\addvspace{.55ex}\noindent\textbf{Case 320.} FB.\enspace \ref{case:125}.
\enspace \(\A{8983}\)\enspace 
7124, 7125, 8467, 8979, 8980, 8983, 8984, 8987, 9649, 10233, 10247, 10261. 
\par

\refstepcounter{casegroup}\label{grp:321}\par\addvspace{.55ex}\noindent\textbf{Case 321.} FB.\enspace \ref{case:179}.
7136, 7137, 8485, 9282, 9283, 9305, 9306, 9676, 10356--10358. 
\par

\refstepcounter{casegroup}\label{grp:322}\par\addvspace{.55ex}\noindent\textbf{Case 322.} FB.\enspace \ref{case:125}.
\enspace \(\A{4588}^{\rm op}\)\enspace 
7238, 7239, 7406, 7496, 7588, 7608, 7618, 7620, 7646, 7984, 8084, 8647, 9708, 9808. 
\par

\refstepcounter{casegroup}\label{grp:323}\par\addvspace{.55ex}\noindent\textbf{Case 323.} FB.\enspace \ref{thm:power-blocks}.
7266, 7267, 7269, 7276, 7277, 7284, 7294, 7295, 7304--7306, 7309, 7316--7318, 7321--7324, 7993, 7997, 8007--8009, 8012, 8013, 8017--8020, 8469, 8542, 8546--8548, 8564, 8565, 8569--8572, 8731--8734, 8739, 8740, 8742, 8743, 8747, 8748, 8750, 8751, 8760, 8761, 8763, 8766, 8769, 8770, 8775, 9624, 9715--9719, 9738, 9742, 9743, 9747, 9748, 9750, 9892--9894, 9898, 9899, 9902, 9904, 9916--9919, 9923, 9924, 9931, 9932, 10059, 10061, 10078, 10090--10093, 10412, 10413, 10418, 10424, 10430--10432, 10451, 10453, 10454, 10501--10507, 10550--10556, 10662--10664, 10667--10670, 10675, 10677--10679, 10682, 10684, 10685, 10688--10690, 10703, 10704, 10707, 10708, 10710--10713, 10716, 10728--10730. 
\par

\refstepcounter{casegroup}\label{grp:324}\par\addvspace{.55ex}\noindent\textbf{Case 324.} FB.\enspace \ref{case:125}.
\enspace \(\A{8983}^{\rm op}\)\enspace 
7288, 7289, 8551, 9704, 9705, 9710, 9886, 9888, 9889, 10406, 10416, 10422. 
\par

\refstepcounter{casegroup}\label{grp:325}\par\addvspace{.55ex}\noindent\textbf{Case 325.} FB.\enspace \ref{case:180}.
7303, 7645, 7647, 8081. 
\par

\refstepcounter{casegroup}\label{grp:326}\par\addvspace{.55ex}\noindent\textbf{Case 326.} FB.\enspace \ref{case:181}.
\enspace \(\A{7319}\)\enspace 
7307, 7319. 
\enspace \(\A{7319}^{\rm op}\)\enspace 
8082, 8618. 
\par

\refstepcounter{casegroup}\label{grp:327}\par\addvspace{.55ex}\noindent\textbf{Case 327.} FB.\enspace \ref{case:182}.
7308, 7568, 8052, 8457. 
\par

\refstepcounter{casegroup}\label{grp:328}\par\addvspace{.55ex}\noindent\textbf{Case 328.} FB.\enspace \ref{cor:prime-examples}.
7446. 
\par

\refstepcounter{casegroup}\label{grp:329}\par\addvspace{.55ex}\noindent\textbf{Case 329.} FB.\enspace \ref{case:183}.
\enspace \(\A{7451}\)\enspace 
7451. 
\enspace \(\A{8831}\)\enspace 
7624, 8831, 8965. 
\par

\refstepcounter{casegroup}\label{grp:330}\par\addvspace{.55ex}\noindent\textbf{Case 330.} FB.\enspace \ref{case:132}.
7470. 
\par

\refstepcounter{casegroup}\label{grp:331}\par\addvspace{.55ex}\noindent\textbf{Case 331.} FB.\enspace \ref{case:184}.
\enspace \(\A{7478}\)\enspace 
7478. 
\enspace \(\A{7478}^{\rm op}\)\enspace 
7861. 
\par

\refstepcounter{casegroup}\label{grp:332}\par\addvspace{.55ex}\noindent\textbf{Case 332.} FB.\enspace \ref{case:185}.
\enspace \(\A{7503}\)\enspace 
7503, 7680. 
\enspace \(\A{7503}^{\rm op}\)\enspace 
8060, 8608. 
\par

\refstepcounter{casegroup}\label{grp:333}\par\addvspace{.55ex}\noindent\textbf{Case 333.} FB.\enspace \ref{case:186}.
7597. 
\par

\refstepcounter{casegroup}\label{grp:334}\par\addvspace{.55ex}\noindent\textbf{Case 334.} FB.\enspace \ref{case:187}.
\enspace \(\A{7601}\leftrightarrow\A{7597}\), \(x+y, xy+yx, x+y, xxxy\)\enspace 
7601. 
\enspace \(\A{7666}\leftrightarrow\A{7609}\), \(x+y, xy+yx, x+y, xxxy\)\enspace 
7666. 
\enspace \(\A{8066}\leftrightarrow\A{7597}\), \(x+y, xy+yx, x+y, xyyy\)\enspace 
8066. 
\enspace \(\A{8596}\leftrightarrow\A{7609}\), \(x+y, xy+yx, x+y, xyyy\)\enspace 
8596. 
\par

\refstepcounter{casegroup}\label{grp:335}\par\addvspace{.55ex}\noindent\textbf{Case 335.} FB.\enspace \ref{thm:heads}.
7605, 8224. 
\par

\refstepcounter{casegroup}\label{grp:336}\par\addvspace{.55ex}\noindent\textbf{Case 336.} FB.\enspace \ref{case:188}.
7609. 
\par

\refstepcounter{casegroup}\label{grp:337}\par\addvspace{.55ex}\noindent\textbf{Case 337.} FB.\enspace \ref{case:189}.
7615, 10032. 
\par

\refstepcounter{casegroup}\label{grp:338}\par\addvspace{.55ex}\noindent\textbf{Case 338.} FB.\enspace \ref{case:190}.
\enspace \(\A{8872}\)\enspace 
7632, 8833, 8872, 8881, 8937. 
\enspace \(\A{8872}^{\rm op}\)\enspace 
7922, 8834, 8887, 8890, 8945. 
\par

\refstepcounter{casegroup}\label{grp:339}\par\addvspace{.55ex}\noindent\textbf{Case 339.} FB.\enspace \ref{case:191}.
7634. 
\par

\refstepcounter{casegroup}\label{grp:340}\par\addvspace{.55ex}\noindent\textbf{Case 340.} FB.\enspace \ref{case:192}.
\enspace \(\A{7644}\)\enspace 
7644. 
\par

\refstepcounter{casegroup}\label{grp:341}\par\addvspace{.55ex}\noindent\textbf{Case 341.} FB.\enspace \ref{case:193}.
\enspace \(\A{7694}\)\enspace 
7649, 7694. 
\enspace \(\A{7695}\)\enspace 
7650, 7695. 
\enspace \(\A{7694}^{\rm op}\)\enspace 
8010, 8550. 
\enspace \(\A{7695}^{\rm op}\)\enspace 
8083, 8619. 
\par

\refstepcounter{casegroup}\label{grp:342}\par\addvspace{.55ex}\noindent\textbf{Case 342.} FB.\enspace \ref{case:194}.
\enspace \(\A{10153}^{\rm op}\)\enspace 
7668, 10175. 
\enspace \(\A{10153}\)\enspace 
8600, 10153. 
\enspace \(\A{10156}\)\enspace 
10150, 10156. 
\enspace \(\A{10156}^{\rm op}\)\enspace 
10164, 10178. 
\par

\refstepcounter{casegroup}\label{grp:343}\par\addvspace{.55ex}\noindent\textbf{Case 343.} FB.\enspace \ref{case:195}.
\enspace \(\A{7672}\)\enspace 
7672. 
\enspace \(\A{7673}\)\enspace 
7673. 
\enspace \(\A{7672}^{\rm op}\)\enspace 
8673. 
\enspace \(\A{7673}^{\rm op}\)\enspace 
8674. 
\par

\refstepcounter{casegroup}\label{grp:344}\par\addvspace{.55ex}\noindent\textbf{Case 344.} FB.\enspace \ref{case:196}.
\enspace \(\A{7674}\)\enspace 
7674, 7675. 
\enspace \(\A{7674}^{\rm op}\)\enspace 
8680, 8681. 
\par

\refstepcounter{casegroup}\label{grp:345}\par\addvspace{.55ex}\noindent\textbf{Case 345.} FB.\enspace \ref{case:197}.
\enspace \(\A{8505}^{\rm op}\)\enspace 
7677. 
\enspace \(\A{8505}\)\enspace 
8505. 
\par

\refstepcounter{casegroup}\label{grp:346}\par\addvspace{.55ex}\noindent\textbf{Case 346.} FB.\enspace \ref{case:198}.
\enspace \(\A{8850}^{\rm op}\)\enspace 
7684, 8989, 9211. 
\enspace \(\A{8850}\)\enspace 
8607, 8850, 8966. 
\par

\refstepcounter{casegroup}\label{grp:347}\par\addvspace{.55ex}\noindent\textbf{Case 347.} FB.\enspace \ref{case:199}.
\enspace \(\A{8889}^{\rm op}\)\enspace 
7692, 8991, 9019, 9020, 9095. 
\enspace \(\A{8889}\)\enspace 
8683, 8853, 8889, 8892, 8947. 
\par

\refstepcounter{casegroup}\label{grp:348}\par\addvspace{.55ex}\noindent\textbf{Case 348.} FB.\enspace Lemma~\ref{lem:terms}, source \(\A{7644}\), case~\ref{grp:340}.
\enspace \(\A{7698}\leftrightarrow\A{7644}\), \(x+y, xy+yx, x+y, xxy\)\enspace 
7698. 
\enspace \(\A{8785}\leftrightarrow\A{7644}\), \(x+y, xy+yx, x+y, xyy\)\enspace 
8785. 
\par

\refstepcounter{casegroup}\label{grp:349}\par\addvspace{.55ex}\noindent\textbf{Case 349.} NFB.\enspace \ref{case:200}.
7805. 
\par

\refstepcounter{casegroup}\label{grp:350}\par\addvspace{.55ex}\noindent\textbf{Case 350.} FB.\enspace \ref{case:201}.
7806, 7945, 8599. 
\par

\refstepcounter{casegroup}\label{grp:351}\par\addvspace{.55ex}\noindent\textbf{Case 351.} FB.\enspace \ref{thm:immediate}.
7811, 8109. 
\par

\refstepcounter{casegroup}\label{grp:352}\par\addvspace{.55ex}\noindent\textbf{Case 352.} FB.\enspace \ref{thm:qualifying}.
7812, 8273. 
\par

\refstepcounter{casegroup}\label{grp:353}\par\addvspace{.55ex}\noindent\textbf{Case 353.} FB.\enspace \ref{case:202}.
\enspace \(\A{13038}^{\rm op}\)\enspace 
7815, 15376, 15390. 
\enspace \(\A{13038}\)\enspace 
8275, 13038, 13261. 
\par

\refstepcounter{casegroup}\label{grp:354}\par\addvspace{.55ex}\noindent\textbf{Case 354.} FB.\enspace \ref{case:203}.
\enspace \(\A{7844}\)\enspace 
7844. 
\enspace \(\A{7844}^{\rm op}\)\enspace 
7868. 
\enspace \(\A{13863}^{\rm op}\)\enspace 
13118. 
\enspace \(\A{13863}\)\enspace 
13863. 
\par

\refstepcounter{casegroup}\label{grp:355}\par\addvspace{.55ex}\noindent\textbf{Case 355.} FB.\enspace \ref{case:192}.
\enspace \(\A{7924}\)\enspace 
7924, 7925, 8845, 8896, 8921, 8934, 14638. 
\par

\refstepcounter{casegroup}\label{grp:356}\par\addvspace{.55ex}\noindent\textbf{Case 356.} FB.\enspace \ref{case:204}.
\enspace \(\A{7953}\)\enspace 
7953, 7954. 
\enspace \(\A{7953}^{\rm op}\)\enspace 
8685, 8686. 
\par

\refstepcounter{casegroup}\label{grp:357}\par\addvspace{.55ex}\noindent\textbf{Case 357.} FB.\enspace \ref{case:205}.
\enspace \(\A{7957}\)\enspace 
7957, 8992, 9131, 9133. 
\enspace \(\A{7957}^{\rm op}\)\enspace 
8617, 8852, 8886, 8939. 
\par

\refstepcounter{casegroup}\label{grp:358}\par\addvspace{.55ex}\noindent\textbf{Case 358.} FB.\enspace Lemma~\ref{lem:terms}, source \(\A{7924}\), case~\ref{grp:355}.
\enspace \(\A{7959}\leftrightarrow\A{7924}\), \(x+y, xy+yx, x+y, xxy\)\enspace 
7959, 7960, 9002, 9072, 9074, 9158, 14642. 
\enspace \(\A{8689}\leftrightarrow\A{7924}\), \(x+y, xy+yx, x+y, xyy\)\enspace 
8689, 8690, 8863, 8900, 8922, 8936, 14683. 
\par

\refstepcounter{casegroup}\label{grp:359}\par\addvspace{.55ex}\noindent\textbf{Case 359.} FB.\enspace \ref{case:206}.
\enspace \(\A{8813}^{\rm op}\)\enspace 
7971. 
\enspace \(\A{8813}\)\enspace 
8813. 
\enspace \(\A{12589}^{\rm op}\)\enspace 
12221. 
\enspace \(\A{12589}\)\enspace 
12589. 
\enspace \(\A{12885}\)\enspace 
12885. 
\enspace \(\A{12885}^{\rm op}\)\enspace 
15012. 
\par

\refstepcounter{casegroup}\label{grp:360}\par\addvspace{.55ex}\noindent\textbf{Case 360.} FB.\enspace \ref{s4592-thm:main}.
\enspace \(\A{7982}\)\enspace 
7982, 7983, 7985, 8033, 8050, 8053, 8595, 8601, 8603, 8786, 8788, 8797, 9709, 9947. 
\par

\refstepcounter{casegroup}\label{grp:361}\par\addvspace{.55ex}\noindent\textbf{Case 361.} FB.\enspace Proposition~\ref{prop:lower-fibres}, source \(\A{6634}\), case~\ref{grp:2}.
\enspace \(\A{6634}\)\enspace 
8025, 8578. 
\par

\refstepcounter{casegroup}\label{grp:362}\par\addvspace{.55ex}\noindent\textbf{Case 362.} FB.\enspace \ref{case:158}.
\enspace \(\A{8095}\)\enspace 
8095. 
\par

\refstepcounter{casegroup}\label{grp:363}\par\addvspace{.55ex}\noindent\textbf{Case 363.} FB.\enspace \ref{case:207}.
\enspace \(\A{8301}^{\rm op}\)\enspace 
8260. 
\enspace \(\A{8301}\)\enspace 
8301. 
\par

\refstepcounter{casegroup}\label{grp:364}\par\addvspace{.55ex}\noindent\textbf{Case 364.} FB.\enspace \ref{thm:unit-paths}.
\enspace \(\A{8353}\)\enspace 
8353. 
\par

\refstepcounter{casegroup}\label{grp:365}\par\addvspace{.55ex}\noindent\textbf{Case 365.} FB.\enspace \ref{case:179}.
8539, 8540, 8552, 9706, 9707, 9711, 9887, 9891, 10409, 10417, 10423. 
\par

\refstepcounter{casegroup}\label{grp:366}\par\addvspace{.55ex}\noindent\textbf{Case 366.} FB.\enspace \ref{case:208}.
\enspace \(\A{8576}\)\enspace 
8576. 
\enspace \(\A{8576}^{\rm op}\)\enspace 
8796. 
\par

\refstepcounter{casegroup}\label{grp:367}\par\addvspace{.55ex}\noindent\textbf{Case 367.} FB.\enspace \ref{case:209}.
\enspace \(\A{8838}\)\enspace 
8838. 
\par

\refstepcounter{casegroup}\label{grp:368}\par\addvspace{.55ex}\noindent\textbf{Case 368.} FB.\enspace \ref{case:210}.
\enspace \(\A{8929}\)\enspace 
8839, 8929, 10249. 
\enspace \(\A{8930}\)\enspace 
8843, 8930. 
\par

\refstepcounter{casegroup}\label{grp:369}\par\addvspace{.55ex}\noindent\textbf{Case 369.} FB.\enspace \ref{case:211}.
\enspace \(\A{8841}\)\enspace 
8841, 8946, 10957. 
\enspace \(\A{8841}^{\rm op}\)\enspace 
8842, 8938, 14089. 
\par

\refstepcounter{casegroup}\label{grp:370}\par\addvspace{.55ex}\noindent\textbf{Case 370.} FB.\enspace \ref{case:212}.
\enspace \(\A{10250}\)\enspace 
8857, 8931, 10250. 
\enspace \(\A{10250}^{\rm op}\)\enspace 
8996, 9061, 10376. 
\par

\refstepcounter{casegroup}\label{grp:371}\par\addvspace{.55ex}\noindent\textbf{Case 371.} FB.\enspace \ref{case:213}.
\enspace \(\A{8859}\)\enspace 
8859, 8948, 10960. 
\enspace \(\A{8859}^{\rm op}\)\enspace 
8999, 9096, 14090. 
\par

\refstepcounter{casegroup}\label{grp:372}\par\addvspace{.55ex}\noindent\textbf{Case 372.} FB.\enspace \ref{case:214}.
\enspace \(\A{8860}\)\enspace 
8860, 8940, 14160. 
\enspace \(\A{8860}^{\rm op}\)\enspace 
8998, 9134, 10958. 
\par

\refstepcounter{casegroup}\label{grp:373}\par\addvspace{.55ex}\noindent\textbf{Case 373.} FB.\enspace \ref{case:215}.
\enspace \(\A{8932}\)\enspace 
8861, 8932. 
\enspace \(\A{8932}^{\rm op}\)\enspace 
9000, 9062. 
\par

\refstepcounter{casegroup}\label{grp:374}\par\addvspace{.55ex}\noindent\textbf{Case 374.} FB.\enspace \ref{case:216}.
\enspace \(\A{8906}\)\enspace 
8906, 8923, 8942. 
\enspace \(\A{8911}\)\enspace 
8911, 8924, 8944. 
\enspace \(\A{8906}^{\rm op}\)\enspace 
9398, 9403, 9585. 
\enspace \(\A{8911}^{\rm op}\)\enspace 
9420, 9422, 9605. 
\par

\refstepcounter{casegroup}\label{grp:375}\par\addvspace{.55ex}\noindent\textbf{Case 375.} FB.\enspace \ref{case:217}.
\enspace \(\A{8910}\)\enspace 
8910. 
\enspace \(\A{8910}^{\rm op}\)\enspace 
9604. 
\par

\refstepcounter{casegroup}\label{grp:376}\par\addvspace{.55ex}\noindent\textbf{Case 376.} FB.\enspace \ref{abs-power-cross}.
9008, 9021, 9075--9077, 9161, 9162, 9189, 10054, 10055, 10075--10077, 10080, 10192, 10194, 10197, 10263, 10265, 10267, 10269--10271, 10296, 10312, 10338, 10339, 10513, 10515, 10516, 10569, 10581, 10735, 10738, 10742, 10783, 14843. 
\par

\refstepcounter{casegroup}\label{grp:377}\par\addvspace{.55ex}\noindent\textbf{Case 377.} FB.\enspace \ref{case:218}.
\enspace \(\A{10377}\)\enspace 
9060, 9063, 10377. 
\par

\refstepcounter{casegroup}\label{grp:378}\par\addvspace{.55ex}\noindent\textbf{Case 378.} FB.\enspace \ref{case:219}.
\enspace \(\A{10071}\)\enspace 
9097, 10071. 
\enspace \(\A{10071}^{\rm op}\)\enspace 
9135, 10073. 
\par

\refstepcounter{casegroup}\label{grp:379}\par\addvspace{.55ex}\noindent\textbf{Case 379.} FB.\enspace Proposition~\ref{prop:lower-fibres}, source \(\A{9271}\), case~\ref{grp:26}.
\enspace \(\A{9271}\)\enspace 
9110, 9261, 10360. 
\par

\refstepcounter{casegroup}\label{grp:380}\par\addvspace{.55ex}\noindent\textbf{Case 380.} FB.\enspace \ref{case:220}.
\enspace \(\A{9212}\)\enspace 
9212, 10063. 
\enspace \(\A{9952}\)\enspace 
9952. 
\par

\refstepcounter{casegroup}\label{grp:381}\par\addvspace{.55ex}\noindent\textbf{Case 381.} FB.\enspace \ref{case:221}.
\enspace \(\A{9222}\)\enspace 
9222, 10066. 
\enspace \(\A{9222}^{\rm op}\)\enspace 
9452, 10067. 
\par

\refstepcounter{casegroup}\label{grp:382}\par\addvspace{.55ex}\noindent\textbf{Case 382.} FB.\enspace Proposition~\ref{prop:lower-fibres}, source \(\A{2160}\), case~\ref{grp:26}.
\enspace \(\A{2160}\)\enspace 
9284. 
\par

\refstepcounter{casegroup}\label{grp:383}\par\addvspace{.55ex}\noindent\textbf{Case 383.} FB.\enspace Proposition~\ref{prop:lower-fibres}, source \(\A{9653}\), case~\ref{grp:26}.
\enspace \(\A{9653}\)\enspace 
9424, 9609, 10439. 
\par

\refstepcounter{casegroup}\label{grp:384}\par\addvspace{.55ex}\noindent\textbf{Case 384.} FB.\enspace \ref{case:222}.
\enspace \(\A{9625}\)\enspace 
9435, 9437, 9625. 
\par

\refstepcounter{casegroup}\label{grp:385}\par\addvspace{.55ex}\noindent\textbf{Case 385.} FB.\enspace \ref{case:222}.
\enspace \(\A{9626}\)\enspace 
9436, 9626, 9664, 10452. 
\par

\refstepcounter{casegroup}\label{grp:386}\par\addvspace{.55ex}\noindent\textbf{Case 386.} FB.\enspace \ref{case:223}.
\enspace \(\A{9458}\)\enspace 
9458, 10064. 
\par

\refstepcounter{casegroup}\label{grp:387}\par\addvspace{.55ex}\noindent\textbf{Case 387.} FB.\enspace \ref{case:224}.
\enspace \(\A{9466}\)\enspace 
9466. 
\enspace \(\A{9466}^{\rm op}\)\enspace 
9729. 
\par

\refstepcounter{casegroup}\label{grp:388}\par\addvspace{.55ex}\noindent\textbf{Case 388.} FB.\enspace \ref{case:225}.
\enspace \(\A{9677}\)\enspace 
9677. 
\enspace \(\A{9677}^{\rm op}\)\enspace 
9733. 
\enspace \(\A{10575}\)\enspace 
10575. 
\enspace \(\A{10575}^{\rm op}\)\enspace 
10576. 
\par

\refstepcounter{casegroup}\label{grp:389}\par\addvspace{.55ex}\noindent\textbf{Case 389.} FB.\enspace Proposition~\ref{prop:lower-fibres}, source \(\A{8431}\), case~\ref{grp:26}.
\enspace \(\A{8431}\)\enspace 
9722. 
\par

\refstepcounter{casegroup}\label{grp:390}\par\addvspace{.55ex}\noindent\textbf{Case 390.} FB.\enspace \ref{case:226}.
\enspace \(\A{10034}\)\enspace 
10034. 
\enspace \(\A{10035}\)\enspace 
10035. 
\enspace \(\A{10036}\)\enspace 
10036, 10043, 10044, 10497. 
\enspace \(\A{11947}\)\enspace 
11947. 
\par

\refstepcounter{casegroup}\label{grp:391}\par\addvspace{.55ex}\noindent\textbf{Case 391.} NFB.\enspace \ref{case:227}.
10081, 12214, 14840. 
\par

\refstepcounter{casegroup}\label{grp:392}\par\addvspace{.55ex}\noindent\textbf{Case 392.} FB.\enspace \ref{abs-quartic-passive-bases}.
\enspace \(\text{\normalfont A5}\)\enspace 
10082. 
\enspace \(\text{\normalfont P5}\)\enspace 
13786. 
\enspace \(\text{\normalfont Y5}\)\enspace 
14433. 
\enspace \(\text{\normalfont X5}\)\enspace 
14434. 
\enspace \(\text{\normalfont B5}\)\enspace 
14814. 
\enspace \(\text{\normalfont S5}\)\enspace 
14844. 
\enspace \(\text{\normalfont R5}\)\enspace 
14845. 
\enspace \(\text{\normalfont Q5}\)\enspace 
14973. 
\par

\refstepcounter{casegroup}\label{grp:393}\par\addvspace{.55ex}\noindent\textbf{Case 393.} FB.\enspace \ref{thm:qgraph}.
10096--10106. 
\par

\refstepcounter{casegroup}\label{grp:394}\par\addvspace{.55ex}\noindent\textbf{Case 394.} NFB.\enspace \ref{case:228}.
10107. 
\par

\refstepcounter{casegroup}\label{grp:395}\par\addvspace{.55ex}\noindent\textbf{Case 395.} FB.\enspace \ref{case:229}.
\enspace \(\A{10108}\)\enspace 
10108. 
\par

\refstepcounter{casegroup}\label{grp:396}\par\addvspace{.55ex}\noindent\textbf{Case 396.} FB.\enspace \ref{case:175}.
\enspace \(\A{10122}\)\enspace 
10122. 
\par

\refstepcounter{casegroup}\label{grp:397}\par\addvspace{.55ex}\noindent\textbf{Case 397.} FB.\enspace \ref{case:175}.
\enspace \(\A{10125}\)\enspace 
10125. 
\par

\refstepcounter{casegroup}\label{grp:398}\par\addvspace{.55ex}\noindent\textbf{Case 398.} FB.\enspace \ref{case:175}.
\enspace \(\A{10126}\)\enspace 
10126. 
\par

\refstepcounter{casegroup}\label{grp:399}\par\addvspace{.55ex}\noindent\textbf{Case 399.} FB.\enspace \ref{case:175}.
\enspace \(\A{10127}\)\enspace 
10127, 10199, 10203, 10602. 
\par

\refstepcounter{casegroup}\label{grp:400}\par\addvspace{.55ex}\noindent\textbf{Case 400.} FB.\enspace \ref{case:230}.
10128, 10139, 11229, 11234, 13089, 13094. 
\par

\refstepcounter{casegroup}\label{grp:401}\par\addvspace{.55ex}\noindent\textbf{Case 401.} FB.\enspace \ref{case:231}.
\enspace \(\A{10133}\)\enspace 
10133. 
\enspace \(\A{10133}^{\rm op}\)\enspace 
10275. 
\par

\refstepcounter{casegroup}\label{grp:402}\par\addvspace{.55ex}\noindent\textbf{Case 402.} FB.\enspace Lemma~\ref{lem:terms}, source \(\A{10125}\), case~\ref{grp:397}.
\enspace \(\A{10136}\leftrightarrow\A{10125}\), \(x+y, xy+yx, x+y, xyy\)\enspace 
10136. 
\enspace \(\A{10278}\leftrightarrow\A{10125}\), \(x+y, xy+yx, x+y, xxy\)\enspace 
10278. 
\par

\refstepcounter{casegroup}\label{grp:403}\par\addvspace{.55ex}\noindent\textbf{Case 403.} FB.\enspace Lemma~\ref{lem:terms}, source \(\A{10126}\), case~\ref{grp:398}.
\enspace \(\A{10137}\leftrightarrow\A{10126}\), \(x+y, xy+yx, x+y, xyy\)\enspace 
10137. 
\enspace \(\A{10279}\leftrightarrow\A{10126}\), \(x+y, xy+yx, x+y, xxy\)\enspace 
10279. 
\par

\refstepcounter{casegroup}\label{grp:404}\par\addvspace{.55ex}\noindent\textbf{Case 404.} FB.\enspace Lemma~\ref{lem:terms}, source \(\A{10127}\), case~\ref{grp:399}.
\enspace \(\A{10138}\leftrightarrow\A{10127}\), \(x+y, xy+yx, x+y, xyy\)\enspace 
10138, 10201, 10206, 10612. 
\enspace \(\A{10280}\leftrightarrow\A{10127}\), \(x+y, xy+yx, x+y, xxy\)\enspace 
10280, 10303, 10305, 10759. 
\par

\refstepcounter{casegroup}\label{grp:405}\par\addvspace{.55ex}\noindent\textbf{Case 405.} FB.\enspace \ref{case:232}.
\enspace \(\A{10157}\)\enspace 
10152, 10157, 10219, 10220, 10224, 10225, 10621, 10628. 
\par

\refstepcounter{casegroup}\label{grp:406}\par\addvspace{.55ex}\noindent\textbf{Case 406.} FB.\enspace \ref{case:232}.
\enspace \(\A{10179}\)\enspace 
10166, 10179, 10237, 10238, 10251, 10254, 10634, 10726. 
\par

\refstepcounter{casegroup}\label{grp:407}\par\addvspace{.55ex}\noindent\textbf{Case 407.} FB.\enspace \ref{case:233}.
\enspace \(\A{10216}\)\enspace 
10216. 
\enspace \(\A{10216}^{\rm op}\)\enspace 
10434. 
\par

\refstepcounter{casegroup}\label{grp:408}\par\addvspace{.55ex}\noindent\textbf{Case 408.} FB.\enspace \ref{case:234}.
\enspace \(\A{10235}\)\enspace 
10234, 10235. 
\enspace \(\A{10260}\)\enspace 
10246, 10260. 
\enspace \(\A{10260}^{\rm op}\)\enspace 
10407, 10408. 
\enspace \(\A{10235}^{\rm op}\)\enspace 
10415, 10421. 
\par

\refstepcounter{casegroup}\label{grp:409}\par\addvspace{.55ex}\noindent\textbf{Case 409.} FB.\enspace \ref{case:133}.
\enspace \(\A{10496}\)\enspace 
10496. 
\par

\refstepcounter{casegroup}\label{grp:410}\par\addvspace{.55ex}\noindent\textbf{Case 410.} FB.\enspace \ref{case:235}.
\enspace \(\A{10578}\)\enspace 
10578. 
\enspace \(\A{13102}\)\enspace 
13102. 
\enspace \(\A{13107}\)\enspace 
13107. 
\par

\refstepcounter{casegroup}\label{grp:411}\par\addvspace{.55ex}\noindent\textbf{Case 411.} FB.\enspace \ref{case:236}.
\enspace \(\A{10587}^{\rm op}\)\enspace 
10579. 
\enspace \(\A{10587}\)\enspace 
10587. 
\enspace \(\A{13323}^{\rm op}\)\enspace 
13108. 
\enspace \(\A{13323}\)\enspace 
13323. 
\par

\refstepcounter{casegroup}\label{grp:412}\par\addvspace{.55ex}\noindent\textbf{Case 412.} FB.\enspace \ref{case:237}.
\enspace \(\A{10595}\)\enspace 
10595, 13885. 
\par

\refstepcounter{casegroup}\label{grp:413}\par\addvspace{.55ex}\noindent\textbf{Case 413.} FB.\enspace \ref{case:58}.
\enspace \(\A{10821}\)\enspace 
10821, 14908. 
\par

\refstepcounter{casegroup}\label{grp:414}\par\addvspace{.55ex}\noindent\textbf{Case 414.} FB.\enspace \ref{case:58}.
\enspace \(\A{10822}\)\enspace 
10822. 
\par

\refstepcounter{casegroup}\label{grp:415}\par\addvspace{.55ex}\noindent\textbf{Case 415.} FB.\enspace \ref{abs-quartic-linear-bases}.
\enspace \(\text{\normalfont E5}\)\enspace 
10823. 
\enspace \(\text{\normalfont Z5}\)\enspace 
14919. 
\par

\refstepcounter{casegroup}\label{grp:416}\par\addvspace{.55ex}\noindent\textbf{Case 416.} FB.\enspace \ref{case:238}.
\enspace \(\A{11577}\)\enspace 
10828, 11416, 11577. 
\par

\refstepcounter{casegroup}\label{grp:417}\par\addvspace{.55ex}\noindent\textbf{Case 417.} FB.\enspace Proposition~\ref{prop:lower-fibres}, source \(\A{10891}\), case~\ref{grp:2}.
\enspace \(\A{10891}\)\enspace 
10887, 10888, 11142, 11145, 12371. 
\par

\refstepcounter{casegroup}\label{grp:418}\par\addvspace{.55ex}\noindent\textbf{Case 418.} FB.\enspace Proposition~\ref{prop:lower-fibres}, source \(\A{10928}\), case~\ref{grp:26}.
\enspace \(\A{10928}\)\enspace 
10890, 11136. 
\par

\refstepcounter{casegroup}\label{grp:419}\par\addvspace{.55ex}\noindent\textbf{Case 419.} FB.\enspace Proposition~\ref{prop:lower-fibres}, source \(\A{10945}\), case~\ref{grp:26}.
\enspace \(\A{10945}\)\enspace 
10934, 10994. 
\par

\refstepcounter{casegroup}\label{grp:420}\par\addvspace{.55ex}\noindent\textbf{Case 420.} FB.\enspace \ref{case:239}.
10939, 11384, 11780. 
\par

\refstepcounter{casegroup}\label{grp:421}\par\addvspace{.55ex}\noindent\textbf{Case 421.} FB.\enspace Proposition~\ref{prop:lower-fibres}, source \(\A{1055}\), case~\ref{grp:26}.
\enspace \(\A{1055}\)\enspace 
10981, 11085. 
\par

\refstepcounter{casegroup}\label{grp:422}\par\addvspace{.55ex}\noindent\textbf{Case 422.} FB.\enspace \ref{case:239}.
11217. 
\par

\refstepcounter{casegroup}\label{grp:423}\par\addvspace{.55ex}\noindent\textbf{Case 423.} NFB.\enspace \ref{thm:private}.
11218, 13396. 
\par

\refstepcounter{casegroup}\label{grp:424}\par\addvspace{.55ex}\noindent\textbf{Case 424.} FB.\enspace \ref{case:240}.
11237, 13761. 
\par

\refstepcounter{casegroup}\label{grp:425}\par\addvspace{.55ex}\noindent\textbf{Case 425.} FB.\enspace \ref{case:241}.
11238, 13815. 
\par

\refstepcounter{casegroup}\label{grp:426}\par\addvspace{.55ex}\noindent\textbf{Case 426.} NFB.\enspace \ref{case:242}.
11269, 11270, 13161. 
\par

\refstepcounter{casegroup}\label{grp:427}\par\addvspace{.55ex}\noindent\textbf{Case 427.} FB.\enspace \ref{case:243}.
\enspace \(\A{12496}\)\enspace 
11281, 12496, 13910. 
\enspace \(\A{12496}^{\rm op}\)\enspace 
13166, 13915, 13929. 
\par

\refstepcounter{casegroup}\label{grp:428}\par\addvspace{.55ex}\noindent\textbf{Case 428.} FB.\enspace \ref{case:244}.
\enspace \(\A{12500}\)\enspace 
11288, 12500, 13912. 
\enspace \(\A{12500}^{\rm op}\)\enspace 
14133, 15287, 15439. 
\par

\refstepcounter{casegroup}\label{grp:429}\par\addvspace{.55ex}\noindent\textbf{Case 429.} FB.\enspace \ref{case:245}.
\enspace \(\A{11375}\)\enspace 
11289, 11375. 
\enspace \(\A{11375}^{\rm op}\)\enspace 
13967, 14134. 
\par

\refstepcounter{casegroup}\label{grp:430}\par\addvspace{.55ex}\noindent\textbf{Case 430.} FB.\enspace \ref{case:246}.
\enspace \(\A{11296}\)\enspace 
11296. 
\enspace \(\A{11296}^{\rm op}\)\enspace 
13095. 
\par

\refstepcounter{casegroup}\label{grp:431}\par\addvspace{.55ex}\noindent\textbf{Case 431.} FB.\enspace \ref{abs-rees-basis}.
\enspace \(\A{11336}\)\enspace 
11336. 
\par

\refstepcounter{casegroup}\label{grp:432}\par\addvspace{.55ex}\noindent\textbf{Case 432.} FB.\enspace \ref{case:247}.
11337, 12057. 
\par

\refstepcounter{casegroup}\label{grp:433}\par\addvspace{.55ex}\noindent\textbf{Case 433.} FB.\enspace Proposition~\ref{prop:lower-fibres}, source \(\A{2475}\), case~\ref{grp:26}.
\enspace \(\A{2475}\)\enspace 
11358, 12493, 13906. 
\par

\refstepcounter{casegroup}\label{grp:434}\par\addvspace{.55ex}\noindent\textbf{Case 434.} FB.\enspace Proposition~\ref{prop:lower-fibres}, source \(\A{2499}\), case~\ref{grp:26}.
\enspace \(\A{2499}\)\enspace 
11361, 12506, 12521. 
\par

\refstepcounter{casegroup}\label{grp:435}\par\addvspace{.55ex}\noindent\textbf{Case 435.} FB.\enspace Proposition~\ref{prop:lower-fibres}, source \(\A{2471}\), case~\ref{grp:26}.
\enspace \(\A{2471}\)\enspace 
11364, 12488, 13909. 
\par

\refstepcounter{casegroup}\label{grp:436}\par\addvspace{.55ex}\noindent\textbf{Case 436.} FB.\enspace Proposition~\ref{prop:lower-fibres}, source \(\A{10856}\), case~\ref{grp:26}.
\enspace \(\A{10856}\)\enspace 
11365. 
\par

\refstepcounter{casegroup}\label{grp:437}\par\addvspace{.55ex}\noindent\textbf{Case 437.} NFB.\enspace \ref{case:248}.
11369. 
\par

\refstepcounter{casegroup}\label{grp:438}\par\addvspace{.55ex}\noindent\textbf{Case 438.} FB.\enspace \ref{case:249}.
\enspace \(\A{11372}\)\enspace 
11372. 
\enspace \(\A{11372}^{\rm op}\)\enspace 
12058. 
\par

\refstepcounter{casegroup}\label{grp:439}\par\addvspace{.55ex}\noindent\textbf{Case 439.} NFB.\enspace \ref{case:250}.
11412, 11496, 11523, 11524, 11564. 
\par

\refstepcounter{casegroup}\label{grp:440}\par\addvspace{.55ex}\noindent\textbf{Case 440.} FB.\enspace \ref{case:251}.
\enspace \(\A{15092}\)\enspace 
11415, 11585, 15092. 
\par

\refstepcounter{casegroup}\label{grp:441}\par\addvspace{.55ex}\noindent\textbf{Case 441.} FB.\enspace Proposition~\ref{prop:lower-fibres}, source \(\A{706}\), case~\ref{grp:26}.
\enspace \(\A{706}\)\enspace 
11525, 12720, 15062. 
\par

\refstepcounter{casegroup}\label{grp:442}\par\addvspace{.55ex}\noindent\textbf{Case 442.} FB.\enspace Proposition~\ref{prop:lower-fibres}, source \(\A{748}\), case~\ref{grp:26}.
\enspace \(\A{748}\)\enspace 
11557, 12769, 12856. 
\par

\refstepcounter{casegroup}\label{grp:443}\par\addvspace{.55ex}\noindent\textbf{Case 443.} FB.\enspace Proposition~\ref{prop:lower-fibres}, source \(\A{699}\), case~\ref{grp:26}.
\enspace \(\A{699}\)\enspace 
11580, 12695, 15087. 
\par

\refstepcounter{casegroup}\label{grp:444}\par\addvspace{.55ex}\noindent\textbf{Case 444.} FB.\enspace Proposition~\ref{prop:lower-fibres}, source \(\A{10835}\), case~\ref{grp:26}.
\enspace \(\A{10835}\)\enspace 
11581. 
\par

\refstepcounter{casegroup}\label{grp:445}\par\addvspace{.55ex}\noindent\textbf{Case 445.} FB.\enspace Proposition~\ref{prop:lower-fibres}, source \(\A{10829}\), case~\ref{grp:2}.
\enspace \(\A{10829}\)\enspace 
11586, 12679. 
\par

\refstepcounter{casegroup}\label{grp:446}\par\addvspace{.55ex}\noindent\textbf{Case 446.} FB.\enspace \ref{case:252}, \ref{case:253}.
\enspace \(\A{11619}\)\enspace 
11619, 11737. 
\enspace \(\A{11619}^{\rm op}\)\enspace 
14143, 14174. 
\par

\refstepcounter{casegroup}\label{grp:447}\par\addvspace{.55ex}\noindent\textbf{Case 447.} NFB.\enspace \ref{case:254}.
11634, 11640, 12530, 13156. 
\par

\refstepcounter{casegroup}\label{grp:448}\par\addvspace{.55ex}\noindent\textbf{Case 448.} NFB.\enspace \ref{case:255}.
11636, 11672, 11673. 
\par

\refstepcounter{casegroup}\label{grp:449}\par\addvspace{.55ex}\noindent\textbf{Case 449.} FB.\enspace \ref{case:256}.
\enspace \(\A{11638}\)\enspace 
11638, 11671, 11675. 
\enspace \(\A{11638}^{\rm op}\)\enspace 
11664, 11678, 11679. 
\par

\refstepcounter{casegroup}\label{grp:450}\par\addvspace{.55ex}\noindent\textbf{Case 450.} FB.\enspace \ref{case:257}.
11642, 11643, 11717, 12100, 12109, 12323, 12879, 13016, 13170, 13762. 
\par

\refstepcounter{casegroup}\label{grp:451}\par\addvspace{.55ex}\noindent\textbf{Case 451.} FB.\enspace \ref{case:258}.
11650. 
\par

\refstepcounter{casegroup}\label{grp:452}\par\addvspace{.55ex}\noindent\textbf{Case 452.} FB.\enspace \ref{thm:pathbasis}.
11670, 12469. 
\par

\refstepcounter{casegroup}\label{grp:453}\par\addvspace{.55ex}\noindent\textbf{Case 453.} FB.\enspace \ref{case:259}.
\enspace \(\A{11703}\)\enspace 
11703. 
\enspace \(\A{11718}\)\enspace 
11718. 
\enspace \(\A{11718}^{\rm op}\)\enspace 
13763. 
\enspace \(\A{11703}^{\rm op}\)\enspace 
13971. 
\par

\refstepcounter{casegroup}\label{grp:454}\par\addvspace{.55ex}\noindent\textbf{Case 454.} FB.\enspace \ref{case:260}.
\enspace \(\A{11722}\)\enspace 
11722. 
\enspace \(\A{11722}^{\rm op}\)\enspace 
12132. 
\par

\refstepcounter{casegroup}\label{grp:455}\par\addvspace{.55ex}\noindent\textbf{Case 455.} FB.\enspace \ref{case:261}.
\enspace \(\A{11723}\)\enspace 
11723. 
\enspace \(\A{11723}^{\rm op}\)\enspace 
12220. 
\par

\refstepcounter{casegroup}\label{grp:456}\par\addvspace{.55ex}\noindent\textbf{Case 456.} FB.\enspace \ref{case:262}.
\enspace \(\A{11734}\)\enspace 
11734. 
\enspace \(\A{11735}\)\enspace 
11735. 
\enspace \(\A{11734}^{\rm op}\)\enspace 
14140. 
\enspace \(\A{11735}^{\rm op}\)\enspace 
14141. 
\par

\refstepcounter{casegroup}\label{grp:457}\par\addvspace{.55ex}\noindent\textbf{Case 457.} FB.\enspace \ref{case:263}.
11736, 14142. 
\par

\refstepcounter{casegroup}\label{grp:458}\par\addvspace{.55ex}\noindent\textbf{Case 458.} FB.\enspace \ref{case:264}.
\enspace \(\A{11771}\)\enspace 
11771, 12076, 12077. 
\enspace \(\A{11772}\)\enspace 
11772, 12082, 13244. 
\par

\refstepcounter{casegroup}\label{grp:459}\par\addvspace{.55ex}\noindent\textbf{Case 459.} FB.\enspace \ref{case:239}.
11776, 12080, 14016. 
\par

\refstepcounter{casegroup}\label{grp:460}\par\addvspace{.55ex}\noindent\textbf{Case 460.} FB.\enspace \ref{case:72}.
\enspace \(\A{12038}\)\enspace 
11778, 11884, 12038. 
\par

\refstepcounter{casegroup}\label{grp:461}\par\addvspace{.55ex}\noindent\textbf{Case 461.} FB.\enspace \ref{case:265}.
\enspace \(\A{11797}\)\enspace 
11797. 
\enspace \(\A{11797}^{\rm op}\)\enspace 
11941. 
\par

\refstepcounter{casegroup}\label{grp:462}\par\addvspace{.55ex}\noindent\textbf{Case 462.} FB.\enspace Proposition~\ref{prop:lower-fibres}, source \(\A{11802}\), case~\ref{grp:2}.
\enspace \(\A{11802}\)\enspace 
11799, 11800, 11807, 11851, 11852, 11890, 11891, 11894, 13077, 13079. 
\par

\refstepcounter{casegroup}\label{grp:463}\par\addvspace{.55ex}\noindent\textbf{Case 463.} FB.\enspace \ref{case:266}.
\enspace \(\A{11832}\)\enspace 
11832, 12118. 
\enspace \(\A{11832}^{\rm op}\)\enspace 
14251, 14270. 
\par

\refstepcounter{casegroup}\label{grp:464}\par\addvspace{.55ex}\noindent\textbf{Case 464.} FB.\enspace \ref{case:72}.
\enspace \(\A{12038}^{\rm op}\)\enspace 
11836, 11981, 12050. 
\par

\refstepcounter{casegroup}\label{grp:465}\par\addvspace{.55ex}\noindent\textbf{Case 465.} FB.\enspace \ref{case:175}.
\enspace \(\A{11845}\)\enspace 
11845. 
\par

\refstepcounter{casegroup}\label{grp:466}\par\addvspace{.55ex}\noindent\textbf{Case 466.} FB.\enspace Lemma~\ref{lem:terms}, source \(\A{11845}\), case~\ref{grp:465}.
\enspace \(\A{11848}\leftrightarrow\A{11845}\), \(x+y, xy+yx, x+y, xxy\)\enspace 
11848. 
\enspace \(\A{12252}\leftrightarrow\A{11845}\), \(x+y, xy+yx, x+y, xyy\)\enspace 
12252. 
\par

\refstepcounter{casegroup}\label{grp:467}\par\addvspace{.55ex}\noindent\textbf{Case 467.} FB.\enspace \ref{thm:heads}.
11849, 12462. 
\par

\refstepcounter{casegroup}\label{grp:468}\par\addvspace{.55ex}\noindent\textbf{Case 468.} FB.\enspace \ref{case:267}.
\enspace \(\A{12473}^{\rm op}\)\enspace 
11850. 
\enspace \(\A{12473}\)\enspace 
12473. 
\par

\refstepcounter{casegroup}\label{grp:469}\par\addvspace{.55ex}\noindent\textbf{Case 469.} FB.\enspace Proposition~\ref{prop:lower-fibres}, source \(\A{11923}\), case~\ref{grp:2}.
\enspace \(\A{11923}\)\enspace 
11853, 11959, 11961, 12002, 13135. 
\par

\refstepcounter{casegroup}\label{grp:470}\par\addvspace{.55ex}\noindent\textbf{Case 470.} FB.\enspace Proposition~\ref{prop:lower-fibres}, source \(\A{11924}\), case~\ref{grp:26}.
\enspace \(\A{11924}\)\enspace 
11861, 11998. 
\par

\refstepcounter{casegroup}\label{grp:471}\par\addvspace{.55ex}\noindent\textbf{Case 471.} FB.\enspace \ref{case:268}.
\enspace \(\A{13098}^{\rm op}\)\enspace 
11927. 
\enspace \(\A{13098}\)\enspace 
13098. 
\par

\refstepcounter{casegroup}\label{grp:472}\par\addvspace{.55ex}\noindent\textbf{Case 472.} FB.\enspace \ref{case:269}.
\enspace \(\A{11954}\)\enspace 
11950, 11954. 
\enspace \(\A{11954}^{\rm op}\)\enspace 
12022, 12163. 
\par

\refstepcounter{casegroup}\label{grp:473}\par\addvspace{.55ex}\noindent\textbf{Case 473.} FB.\enspace \ref{case:270}.
\enspace \(\A{12152}^{\rm op}\)\enspace 
11952. 
\enspace \(\A{12152}\)\enspace 
12152. 
\par

\refstepcounter{casegroup}\label{grp:474}\par\addvspace{.55ex}\noindent\textbf{Case 474.} FB.\enspace \ref{case:271}.
\enspace \(\A{13045}^{\rm op}\)\enspace 
11963. 
\enspace \(\A{13045}\)\enspace 
13045. 
\par

\refstepcounter{casegroup}\label{grp:475}\par\addvspace{.55ex}\noindent\textbf{Case 475.} FB.\enspace \ref{case:272}.
\enspace \(\A{13058}^{\rm op}\)\enspace 
11964. 
\enspace \(\A{13321}^{\rm op}\)\enspace 
11970. 
\enspace \(\A{13034}^{\rm op}\)\enspace 
12004. 
\enspace \(\A{13034}\)\enspace 
13034. 
\enspace \(\A{13058}\)\enspace 
13058. 
\enspace \(\A{13321}\)\enspace 
13321. 
\par

\refstepcounter{casegroup}\label{grp:476}\par\addvspace{.55ex}\noindent\textbf{Case 476.} FB.\enspace \ref{case:273}.
\enspace \(\A{13100}^{\rm op}\)\enspace 
11968. 
\enspace \(\A{13100}\)\enspace 
13100. 
\par

\refstepcounter{casegroup}\label{grp:477}\par\addvspace{.55ex}\noindent\textbf{Case 477.} FB.\enspace \ref{case:274}.
\enspace \(\A{13101}^{\rm op}\)\enspace 
12026. 
\enspace \(\A{13101}\)\enspace 
13101. 
\par

\refstepcounter{casegroup}\label{grp:478}\par\addvspace{.55ex}\noindent\textbf{Case 478.} NFB.\enspace \ref{thm:odd-meet}.
\enspace \(\text{\normalfont Scab\_zero}\)\enspace 
12029. 
\par

\refstepcounter{casegroup}\label{grp:479}\par\addvspace{.55ex}\noindent\textbf{Case 479.} FB.\enspace \ref{case:275}.
\enspace \(\A{12030}\)\enspace 
12030. 
\enspace \(\A{12250}\)\enspace 
12250. 
\enspace \(\A{12251}\)\enspace 
12251. 
\enspace \(\A{12250}^{\rm op}\)\enspace 
13183. 
\enspace \(\A{12251}^{\rm op}\)\enspace 
13184. 
\enspace \(\A{13897}\)\enspace 
13897. 
\enspace \(\A{13897}^{\rm op}\)\enspace 
13926. 
\par

\refstepcounter{casegroup}\label{grp:480}\par\addvspace{.55ex}\noindent\textbf{Case 480.} FB.\enspace \ref{case:276}.
\enspace \(\A{12095}\)\enspace 
12095. 
\enspace \(\A{12095}^{\rm op}\)\enspace 
12583. 
\par

\refstepcounter{casegroup}\label{grp:481}\par\addvspace{.55ex}\noindent\textbf{Case 481.} FB.\enspace \ref{case:277}.
\enspace \(\A{13349}^{\rm op}\)\enspace 
12117. 
\enspace \(\A{13349}\)\enspace 
13349. 
\par

\refstepcounter{casegroup}\label{grp:482}\par\addvspace{.55ex}\noindent\textbf{Case 482.} FB.\enspace \ref{case:278}.
\enspace \(\A{12126}\)\enspace 
12126. 
\enspace \(\A{12126}^{\rm op}\)\enspace 
14288. 
\par

\refstepcounter{casegroup}\label{grp:483}\par\addvspace{.55ex}\noindent\textbf{Case 483.} FB.\enspace \ref{case:279}.
\enspace \(\A{13417}^{\rm op}\)\enspace 
12134. 
\enspace \(\A{13713}^{\rm op}\)\enspace 
13296. 
\enspace \(\A{13417}\)\enspace 
13417. 
\enspace \(\A{13713}\)\enspace 
13713. 
\par

\refstepcounter{casegroup}\label{grp:484}\par\addvspace{.55ex}\noindent\textbf{Case 484.} FB.\enspace \ref{case:280}.
\enspace \(\A{13764}^{\rm op}\)\enspace 
12135. 
\enspace \(\A{13764}\)\enspace 
13764. 
\par

\refstepcounter{casegroup}\label{grp:485}\par\addvspace{.55ex}\noindent\textbf{Case 485.} FB.\enspace \ref{case:281}.
\enspace \(\A{12204}\)\enspace 
12204. 
\par

\refstepcounter{casegroup}\label{grp:486}\par\addvspace{.55ex}\noindent\textbf{Case 486.} FB.\enspace \ref{case:282}.
\enspace \(\A{12205}\)\enspace 
12205. 
\par

\refstepcounter{casegroup}\label{grp:487}\par\addvspace{.55ex}\noindent\textbf{Case 487.} FB.\enspace \ref{case:283}.
\enspace \(\A{12206}\)\enspace 
12206. 
\par

\refstepcounter{casegroup}\label{grp:488}\par\addvspace{.55ex}\noindent\textbf{Case 488.} FB.\enspace \ref{case:284}.
\enspace \(\A{12208}\)\enspace 
12208. 
\enspace \(\A{14369}\)\enspace 
14369. 
\enspace \(\A{14772}\)\enspace 
14772. 
\par

\refstepcounter{casegroup}\label{grp:489}\par\addvspace{.55ex}\noindent\textbf{Case 489.} FB.\enspace Section~\ref{sec:four-sources}.
\enspace \(\text{\normalfont S4\_413}\)\enspace 
12211, 12290, 14374. 
\par

\refstepcounter{casegroup}\label{grp:490}\par\addvspace{.55ex}\noindent\textbf{Case 490.} FB.\enspace \ref{cor:qtop}.
12216, 12218, 12235--12237, 13370, 13376, 13381, 13789, 13790, 13792--13795, 14944, 14948, 14949, 14951, 14952, 14976, 14977, 14980, 14985--14992, 15006--15011, 15015--15017, 15022--15025, 15132--15134, 15141--15144, 15150--15153, 15158, 15168, 15171, 15173, 15174, 15179--15182, 15184--15187, 15189, 15193--15196, 15261--15263, 15265, 15266, 15268--15270, 15411, 15424, 15427, 15428, 15431, 15432, 15469, 15470, 15593, 15600--15602, 15605--15607, 15612--15614, 15617, 15621--15623, 15626--15628, 15678, 15686, 15690--15692, 15703, 15704, 15709, 15712, 15737, 15738, 15740, 15742, 15743, 15746, 15748--15751. 
\par

\refstepcounter{casegroup}\label{grp:491}\par\addvspace{.55ex}\noindent\textbf{Case 491.} FB.\enspace \ref{case:285}.
\enspace \(\A{13375}^{\rm op}\)\enspace 
12217, 15131. 
\enspace \(\A{13352}^{\rm op}\)\enspace 
12263. 
\enspace \(\A{13375}\)\enspace 
12289, 13375. 
\enspace \(\A{13352}\)\enspace 
13352. 
\par

\refstepcounter{casegroup}\label{grp:492}\par\addvspace{.55ex}\noindent\textbf{Case 492.} FB.\enspace \ref{abs-directed-cubic-bases}.
\enspace \(\text{\normalfont Path3}\)\enspace 
12294. 
\enspace \(\text{\normalfont Marked3}\)\enspace 
12295. 
\enspace \(\text{\normalfont Trans3}\)\enspace 
12296. 
\enspace \(\text{\normalfont Marked3}^{\rm op}\)\enspace 
13353. 
\enspace \(\text{\normalfont Tail3}\)\enspace 
13361. 
\enspace \(\text{\normalfont Tail3}^{\rm op}\)\enspace 
14911. 
\enspace \(\text{\normalfont Inflation3}\)\enspace 
14997. 
\par

\refstepcounter{casegroup}\label{grp:493}\par\addvspace{.55ex}\noindent\textbf{Case 493.} NFB.\enspace \ref{case:286}.
12438, 13708. 
\par

\refstepcounter{casegroup}\label{grp:494}\par\addvspace{.55ex}\noindent\textbf{Case 494.} FB.\enspace \ref{case:287}.
\enspace \(\A{12468}\)\enspace 
12468. 
\par

\refstepcounter{casegroup}\label{grp:495}\par\addvspace{.55ex}\noindent\textbf{Case 495.} FB.\enspace Proposition~\ref{prop:lower-fibres}, source \(\A{13394}\), case~\ref{grp:26}.
\enspace \(\A{13394}\)\enspace 
12472, 12540. 
\par

\refstepcounter{casegroup}\label{grp:496}\par\addvspace{.55ex}\noindent\textbf{Case 496.} FB.\enspace Proposition~\ref{prop:lower-fibres}, source \(\A{2515}\), case~\ref{grp:26}.
\enspace \(\A{2515}\)\enspace 
12476. 
\par

\refstepcounter{casegroup}\label{grp:497}\par\addvspace{.55ex}\noindent\textbf{Case 497.} FB.\enspace Proposition~\ref{prop:lower-fibres}, source \(\A{13447}\), case~\ref{grp:26}.
\enspace \(\A{13447}\)\enspace 
12479, 12628. 
\par

\refstepcounter{casegroup}\label{grp:498}\par\addvspace{.55ex}\noindent\textbf{Case 498.} NFB.\enspace \ref{case:288}.
\enspace \(\A{12480}\)\enspace 
12480. 
\enspace \(\A{12480}^{\rm op}\)\enspace 
13848. 
\par

\refstepcounter{casegroup}\label{grp:499}\par\addvspace{.55ex}\noindent\textbf{Case 499.} FB.\enspace Proposition~\ref{prop:lower-fibres}, source \(\A{13448}\), case~\ref{grp:26}.
\enspace \(\A{13448}\)\enspace 
12483, 12636. 
\par

\refstepcounter{casegroup}\label{grp:500}\par\addvspace{.55ex}\noindent\textbf{Case 500.} FB.\enspace Proposition~\ref{prop:lower-fibres}, source \(\A{2512}\), case~\ref{grp:26}.
\enspace \(\A{2512}\)\enspace 
12516. 
\par

\refstepcounter{casegroup}\label{grp:501}\par\addvspace{.55ex}\noindent\textbf{Case 501.} FB.\enspace Proposition~\ref{prop:lower-fibres}, source \(\A{12307}\), case~\ref{grp:2}.
\enspace \(\A{12307}\)\enspace 
12534, 12535, 12561, 12564, 13478. 
\par

\refstepcounter{casegroup}\label{grp:502}\par\addvspace{.55ex}\noindent\textbf{Case 502.} FB.\enspace Proposition~\ref{prop:lower-fibres}, source \(\A{13393}\), case~\ref{grp:26}.
\enspace \(\A{13393}\)\enspace 
12537, 12556. 
\par

\refstepcounter{casegroup}\label{grp:503}\par\addvspace{.55ex}\noindent\textbf{Case 503.} FB.\enspace Proposition~\ref{prop:lower-fibres}, source \(\A{12313}\), case~\ref{grp:2}.
\enspace \(\A{12313}\)\enspace 
12565--12569, 12572, 12573, 12576, 13484, 13573. 
\par

\refstepcounter{casegroup}\label{grp:504}\par\addvspace{.55ex}\noindent\textbf{Case 504.} FB.\enspace Proposition~\ref{prop:lower-fibres}, source \(\A{12326}\), case~\ref{grp:2}.
\enspace \(\A{12326}\)\enspace 
12595, 12611, 12614, 13487. 
\par

\refstepcounter{casegroup}\label{grp:505}\par\addvspace{.55ex}\noindent\textbf{Case 505.} FB.\enspace Proposition~\ref{prop:lower-fibres}, source \(\A{778}\), case~\ref{grp:26}.
\enspace \(\A{778}\)\enspace 
12607. 
\par

\refstepcounter{casegroup}\label{grp:506}\par\addvspace{.55ex}\noindent\textbf{Case 506.} FB.\enspace Proposition~\ref{prop:lower-fibres}, source \(\A{12332}\), case~\ref{grp:2}.
\enspace \(\A{12332}\)\enspace 
12615, 12616, 12657, 12660, 13491. 
\par

\refstepcounter{casegroup}\label{grp:507}\par\addvspace{.55ex}\noindent\textbf{Case 507.} FB.\enspace Proposition~\ref{prop:lower-fibres}, source \(\A{13444}\), case~\ref{grp:26}.
\enspace \(\A{13444}\)\enspace 
12618, 12652. 
\par

\refstepcounter{casegroup}\label{grp:508}\par\addvspace{.55ex}\noindent\textbf{Case 508.} FB.\enspace \ref{case:289}.
12621. 
\par

\refstepcounter{casegroup}\label{grp:509}\par\addvspace{.55ex}\noindent\textbf{Case 509.} FB.\enspace Proposition~\ref{prop:lower-fibres}, source \(\A{12339}\), case~\ref{grp:2}.
\enspace \(\A{12339}\)\enspace 
12622, 12623, 12673, 12676, 13495. 
\par

\refstepcounter{casegroup}\label{grp:510}\par\addvspace{.55ex}\noindent\textbf{Case 510.} FB.\enspace \ref{case:290}.
\enspace \(\A{12624}\)\enspace 
12624, 12669, 12672, 13594. 
\enspace \(\A{12624}^{\rm op}\)\enspace 
14495, 14539, 14569, 15533. 
\par

\refstepcounter{casegroup}\label{grp:511}\par\addvspace{.55ex}\noindent\textbf{Case 511.} FB.\enspace Proposition~\ref{prop:lower-fibres}, source \(\A{1334}\), case~\ref{grp:2}.
\enspace \(\A{1334}\)\enspace 
12625, 12668. 
\par

\refstepcounter{casegroup}\label{grp:512}\par\addvspace{.55ex}\noindent\textbf{Case 512.} FB.\enspace \ref{case:291}.
\enspace \(\A{12626}\)\enspace 
12626, 12642, 12667. 
\enspace \(\A{12626}^{\rm op}\)\enspace 
14488, 14571, 14573. 
\par

\refstepcounter{casegroup}\label{grp:513}\par\addvspace{.55ex}\noindent\textbf{Case 513.} FB.\enspace Proposition~\ref{prop:lower-fibres}, source \(\A{12428}\), case~\ref{grp:26}.
\enspace \(\A{12428}\)\enspace 
12816, 12835, 12839, 13540. 
\par

\refstepcounter{casegroup}\label{grp:514}\par\addvspace{.55ex}\noindent\textbf{Case 514.} FB.\enspace \ref{case:292}.
12880, 12887, 13017, 15108, 15113. 
\par

\refstepcounter{casegroup}\label{grp:515}\par\addvspace{.55ex}\noindent\textbf{Case 515.} FB.\enspace \ref{case:293}.
\enspace \(\A{12884}\)\enspace 
12884. 
\enspace \(\A{12884}^{\rm op}\)\enspace 
14994. 
\par

\refstepcounter{casegroup}\label{grp:516}\par\addvspace{.55ex}\noindent\textbf{Case 516.} FB.\enspace Proposition~\ref{prop:lower-fibres}, source \(\A{12902}\), case~\ref{grp:2}.
\enspace \(\A{12902}\)\enspace 
12946, 12948, 12950, 12956, 13767. 
\par

\refstepcounter{casegroup}\label{grp:517}\par\addvspace{.55ex}\noindent\textbf{Case 517.} FB.\enspace Proposition~\ref{prop:lower-fibres}, source \(\A{12903}\), case~\ref{grp:26}.
\enspace \(\A{12903}\)\enspace 
12947, 12952. 
\par

\refstepcounter{casegroup}\label{grp:518}\par\addvspace{.55ex}\noindent\textbf{Case 518.} FB.\enspace Proposition~\ref{prop:lower-fibres}, source \(\A{12907}\), case~\ref{grp:26}.
\enspace \(\A{12907}\)\enspace 
12959, 12978, 12982, 13770. 
\par

\refstepcounter{casegroup}\label{grp:519}\par\addvspace{.55ex}\noindent\textbf{Case 519.} FB.\enspace Proposition~\ref{prop:lower-fibres}, source \(\A{902}\), case~\ref{grp:26}.
\enspace \(\A{902}\)\enspace 
12973. 
\par

\refstepcounter{casegroup}\label{grp:520}\par\addvspace{.55ex}\noindent\textbf{Case 520.} FB.\enspace Proposition~\ref{prop:lower-fibres}, source \(\A{12914}\), case~\ref{grp:2}.
\enspace \(\A{12914}\)\enspace 
12984, 13002, 13006, 13780. 
\par

\refstepcounter{casegroup}\label{grp:521}\par\addvspace{.55ex}\noindent\textbf{Case 521.} FB.\enspace \ref{case:294}.
\enspace \(\A{13020}\)\enspace 
13020. 
\enspace \(\A{13020}^{\rm op}\)\enspace 
13295. 
\par

\refstepcounter{casegroup}\label{grp:522}\par\addvspace{.55ex}\noindent\textbf{Case 522.} FB.\enspace \ref{case:295}.
\enspace \(\A{13062}\)\enspace 
13062. 
\enspace \(\A{13062}^{\rm op}\)\enspace 
13105. 
\par

\refstepcounter{casegroup}\label{grp:523}\par\addvspace{.55ex}\noindent\textbf{Case 523.} FB.\enspace \ref{case:72}.
\enspace \(\A{13065}\)\enspace 
13065. 
\par

\refstepcounter{casegroup}\label{grp:524}\par\addvspace{.55ex}\noindent\textbf{Case 524.} FB.\enspace \ref{case:296}.
\enspace \(\A{13066}\)\enspace 
13066. 
\enspace \(\A{13066}^{\rm op}\)\enspace 
13446. 
\par

\refstepcounter{casegroup}\label{grp:525}\par\addvspace{.55ex}\noindent\textbf{Case 525.} FB.\enspace \ref{case:239}.
13092. 
\par

\refstepcounter{casegroup}\label{grp:526}\par\addvspace{.55ex}\noindent\textbf{Case 526.} FB.\enspace \ref{case:175}.
\enspace \(\A{13110}\)\enspace 
13110. 
\par

\refstepcounter{casegroup}\label{grp:527}\par\addvspace{.55ex}\noindent\textbf{Case 527.} FB.\enspace \ref{case:72}.
\enspace \(\A{13065}^{\rm op}\)\enspace 
13113. 
\par

\refstepcounter{casegroup}\label{grp:528}\par\addvspace{.55ex}\noindent\textbf{Case 528.} FB.\enspace \ref{case:175}.
\enspace \(\A{13115}\)\enspace 
13115. 
\par

\refstepcounter{casegroup}\label{grp:529}\par\addvspace{.55ex}\noindent\textbf{Case 529.} FB.\enspace Lemma~\ref{lem:terms}, source \(\A{13115}\), case~\ref{grp:528}.
\enspace \(\A{13117}\leftrightarrow\A{13115}\), \(x+y, xy+yx, x+y, xxy\)\enspace 
13117. 
\enspace \(\A{13328}\leftrightarrow\A{13115}\), \(x+y, xy+yx, x+y, xyy\)\enspace 
13328. 
\par

\refstepcounter{casegroup}\label{grp:530}\par\addvspace{.55ex}\noindent\textbf{Case 530.} FB.\enspace \ref{case:297}.
\enspace \(\A{13189}\)\enspace 
13189. 
\par

\refstepcounter{casegroup}\label{grp:531}\par\addvspace{.55ex}\noindent\textbf{Case 531.} FB.\enspace \ref{case:127}.
\enspace \(\A{13206}\)\enspace 
13206. 
\par

\refstepcounter{casegroup}\label{grp:532}\par\addvspace{.55ex}\noindent\textbf{Case 532.} FB.\enspace \ref{case:298}.
\enspace \(\A{13216}^{\rm op}\)\enspace 
13210. 
\enspace \(\A{13216}\)\enspace 
13216. 
\par

\refstepcounter{casegroup}\label{grp:533}\par\addvspace{.55ex}\noindent\textbf{Case 533.} FB.\enspace \ref{case:299}.
\enspace \(\A{13218}\)\enspace 
13218. 
\par

\refstepcounter{casegroup}\label{grp:534}\par\addvspace{.55ex}\noindent\textbf{Case 534.} FB.\enspace \ref{case:300}.
\enspace \(\A{14622}^{\rm op}\)\enspace 
13227. 
\enspace \(\A{14622}\)\enspace 
14622. 
\par

\refstepcounter{casegroup}\label{grp:535}\par\addvspace{.55ex}\noindent\textbf{Case 535.} FB.\enspace \ref{case:127}.
\enspace \(\A{13239}\)\enspace 
13239. 
\par

\refstepcounter{casegroup}\label{grp:536}\par\addvspace{.55ex}\noindent\textbf{Case 536.} FB.\enspace \ref{case:301}.
\enspace \(\A{13251}\)\enspace 
13251. 
\enspace \(\A{13251}^{\rm op}\)\enspace 
13272. 
\par

\refstepcounter{casegroup}\label{grp:537}\par\addvspace{.55ex}\noindent\textbf{Case 537.} FB.\enspace \ref{case:302}.
\enspace \(\A{13298}\)\enspace 
13298. 
\enspace \(\A{13298}^{\rm op}\)\enspace 
14995. 
\par

\refstepcounter{casegroup}\label{grp:538}\par\addvspace{.55ex}\noindent\textbf{Case 538.} FB.\enspace \ref{case:303}.
\enspace \(\A{13299}\)\enspace 
13299. 
\enspace \(\A{13902}\)\enspace 
13902. 
\enspace \(\A{13299}^{\rm op}\)\enspace 
15156. 
\enspace \(\A{13902}^{\rm op}\)\enspace 
15426. 
\par

\refstepcounter{casegroup}\label{grp:539}\par\addvspace{.55ex}\noindent\textbf{Case 539.} FB.\enspace \ref{case:304}.
\enspace \(\A{14623}^{\rm op}\)\enspace 
13330. 
\enspace \(\A{14623}\)\enspace 
14623. 
\par

\refstepcounter{casegroup}\label{grp:540}\par\addvspace{.55ex}\noindent\textbf{Case 540.} FB.\enspace \ref{case:305}.
\enspace \(\A{13358}\)\enspace 
13358. 
\par

\refstepcounter{casegroup}\label{grp:541}\par\addvspace{.55ex}\noindent\textbf{Case 541.} FB.\enspace \ref{case:306}.
\enspace \(\A{13363}\)\enspace 
13363. 
\par

\refstepcounter{casegroup}\label{grp:542}\par\addvspace{.55ex}\noindent\textbf{Case 542.} FB.\enspace \ref{case:307}.
\enspace \(\A{13366}\)\enspace 
13366. 
\par

\refstepcounter{casegroup}\label{grp:543}\par\addvspace{.55ex}\noindent\textbf{Case 543.} FB.\enspace \ref{abs-directed-quartic-bases}.
\enspace \(\text{\normalfont FlatParity4}\)\enspace 
13372. 
\enspace \(\text{\normalfont OrderedParity4}\)\enspace 
13791. 
\enspace \(\text{\normalfont OrderedParity4}^{\rm op}\)\enspace 
14892. 
\enspace \(\text{\normalfont FlatParity4}^{\rm op}\)\enspace 
15129. 
\par

\refstepcounter{casegroup}\label{grp:544}\par\addvspace{.55ex}\noindent\textbf{Case 544.} FB.\enspace \ref{case:308}.
\enspace \(\A{13465}\)\enspace 
13465. 
\enspace \(\A{13465}^{\rm op}\)\enspace 
15425. 
\par

\refstepcounter{casegroup}\label{grp:545}\par\addvspace{.55ex}\noindent\textbf{Case 545.} FB.\enspace \ref{abs-quartic-parity-basis}.
\enspace \(\text{\normalfont H4}\)\enspace 
13787, 13788, 14816, 14906, 14914--14916, 15578, 15580. 
\par

\refstepcounter{casegroup}\label{grp:546}\par\addvspace{.55ex}\noindent\textbf{Case 546.} FB.\enspace \ref{case:309}.
\enspace \(\A{13837}\)\enspace 
13837. 
\par

\refstepcounter{casegroup}\label{grp:547}\par\addvspace{.55ex}\noindent\textbf{Case 547.} FB.\enspace \ref{case:310}.
\enspace \(\A{13840}\)\enspace 
13840. 
\par

\refstepcounter{casegroup}\label{grp:548}\par\addvspace{.55ex}\noindent\textbf{Case 548.} FB.\enspace Proposition~\ref{prop:lower-fibres}, source \(\A{3876}\), case~\ref{grp:2}.
\enspace \(\A{3876}\)\enspace 
13850, 14534. 
\par

\refstepcounter{casegroup}\label{grp:549}\par\addvspace{.55ex}\noindent\textbf{Case 549.} FB.\enspace \ref{case:311}.
\enspace \(\A{14541}^{\rm op}\)\enspace 
13853. 
\enspace \(\A{14541}\)\enspace 
14541. 
\par

\refstepcounter{casegroup}\label{grp:550}\par\addvspace{.55ex}\noindent\textbf{Case 550.} FB.\enspace Proposition~\ref{prop:lower-fibres}, source \(\A{2309}\), case~\ref{grp:26}.
\enspace \(\A{2309}\)\enspace 
13890. 
\par

\refstepcounter{casegroup}\label{grp:551}\par\addvspace{.55ex}\noindent\textbf{Case 551.} FB.\enspace Proposition~\ref{prop:lower-fibres}, source \(\A{2314}\), case~\ref{grp:26}.
\enspace \(\A{2314}\)\enspace 
13893. 
\par

\refstepcounter{casegroup}\label{grp:552}\par\addvspace{.55ex}\noindent\textbf{Case 552.} FB.\enspace \ref{case:312}.
\enspace \(\A{13900}\)\enspace 
13900. 
\par

\refstepcounter{casegroup}\label{grp:553}\par\addvspace{.55ex}\noindent\textbf{Case 553.} FB.\enspace \ref{case:258}.
14105. 
\par

\refstepcounter{casegroup}\label{grp:554}\par\addvspace{.55ex}\noindent\textbf{Case 554.} FB.\enspace \ref{case:313}.
\enspace \(\A{14361}\)\enspace 
14361, 14423. 
\enspace \(\A{14361}^{\rm op}\)\enspace 
14393, 14419. 
\par

\refstepcounter{casegroup}\label{grp:555}\par\addvspace{.55ex}\noindent\textbf{Case 555.} FB.\enspace \ref{case:314}.
\enspace \(\A{14366}\)\enspace 
14366. 
\par

\refstepcounter{casegroup}\label{grp:556}\par\addvspace{.55ex}\noindent\textbf{Case 556.} FB.\enspace \ref{case:80}.
\enspace \(\A{14368}\)\enspace 
14368. 
\par

\refstepcounter{casegroup}\label{grp:557}\par\addvspace{.55ex}\noindent\textbf{Case 557.} FB.\enspace \ref{case:315}.
\enspace \(\A{14370}\)\enspace 
14370. 
\enspace \(\A{14370}^{\rm op}\)\enspace 
14422. 
\par

\refstepcounter{casegroup}\label{grp:558}\par\addvspace{.55ex}\noindent\textbf{Case 558.} FB.\enspace \ref{case:316}.
\enspace \(\A{14414}\)\enspace 
14414. 
\par

\refstepcounter{casegroup}\label{grp:559}\par\addvspace{.55ex}\noindent\textbf{Case 559.} FB.\enspace \ref{case:80}.
\enspace \(\A{14425}\)\enspace 
14425, 14883. 
\par

\refstepcounter{casegroup}\label{grp:560}\par\addvspace{.55ex}\noindent\textbf{Case 560.} FB.\enspace \ref{case:317}.
\enspace \(\A{14426}\)\enspace 
14426. 
\par

\refstepcounter{casegroup}\label{grp:561}\par\addvspace{.55ex}\noindent\textbf{Case 561.} FB.\enspace \ref{case:316}.
\enspace \(\A{14435}\)\enspace 
14435. 
\par

\refstepcounter{casegroup}\label{grp:562}\par\addvspace{.55ex}\noindent\textbf{Case 562.} FB.\enspace \ref{abs-nilgraph-bases}.
\enspace \(\text{\normalfont U5}\)\enspace 
14441. 
\enspace \(\text{\normalfont T4}\)\enspace 
14956, 14963, 15594. 
\par

\refstepcounter{casegroup}\label{grp:563}\par\addvspace{.55ex}\noindent\textbf{Case 563.} FB.\enspace \ref{abs-quartic-bases}.
\enspace \(\text{\normalfont W5}\)\enspace 
14442. 
\enspace \(\text{\normalfont V5}\)\enspace 
14964, 14971, 15598. 
\par

\refstepcounter{casegroup}\label{grp:564}\par\addvspace{.55ex}\noindent\textbf{Case 564.} NFB.\enspace \ref{case:318}.
14450, 14702. 
\par

\refstepcounter{casegroup}\label{grp:565}\par\addvspace{.55ex}\noindent\textbf{Case 565.} FB.\enspace \ref{case:319}.
\enspace \(\A{14452}\)\enspace 
14452, 14739. 
\enspace \(\A{14452}^{\rm op}\)\enspace 
14464, 14703. 
\par

\refstepcounter{casegroup}\label{grp:566}\par\addvspace{.55ex}\noindent\textbf{Case 566.} FB.\enspace \ref{case:320}.
\enspace \(\A{14453}\)\enspace 
14453. 
\enspace \(\A{14453}^{\rm op}\)\enspace 
14465. 
\par

\refstepcounter{casegroup}\label{grp:567}\par\addvspace{.55ex}\noindent\textbf{Case 567.} FB.\enspace Proposition~\ref{prop:lower-fibres}, source \(\A{14490}\), case~\ref{grp:2}.
\enspace \(\A{14490}\)\enspace 
14486, 14487, 14574, 14577, 15531. 
\par

\refstepcounter{casegroup}\label{grp:568}\par\addvspace{.55ex}\noindent\textbf{Case 568.} FB.\enspace Proposition~\ref{prop:lower-fibres}, source \(\A{14527}\), case~\ref{grp:26}.
\enspace \(\A{14527}\)\enspace 
14489, 14570. 
\par

\refstepcounter{casegroup}\label{grp:569}\par\addvspace{.55ex}\noindent\textbf{Case 569.} FB.\enspace \ref{case:321}.
\enspace \(\A{14637}\)\enspace 
14637. 
\par

\refstepcounter{casegroup}\label{grp:570}\par\addvspace{.55ex}\noindent\textbf{Case 570.} FB.\enspace \ref{case:95}.
\enspace \(\A{14659}\)\enspace 
14659. 
\par

\refstepcounter{casegroup}\label{grp:571}\par\addvspace{.55ex}\noindent\textbf{Case 571.} FB.\enspace \ref{case:322}.
\enspace \(\A{14674}^{\rm op}\)\enspace 
14674, 14734. 
\par

\refstepcounter{casegroup}\label{grp:572}\par\addvspace{.55ex}\noindent\textbf{Case 572.} FB.\enspace \ref{case:323}.
14676, 15244. 
\par

\refstepcounter{casegroup}\label{grp:573}\par\addvspace{.55ex}\noindent\textbf{Case 573.} FB.\enspace \ref{case:324}.
\enspace \(\A{14705}\)\enspace 
14705. 
\par

\refstepcounter{casegroup}\label{grp:574}\par\addvspace{.55ex}\noindent\textbf{Case 574.} FB.\enspace \ref{case:325}.
\enspace \(\A{14781}\)\enspace 
14707, 14781. 
\par

\refstepcounter{casegroup}\label{grp:575}\par\addvspace{.55ex}\noindent\textbf{Case 575.} FB.\enspace \ref{case:326}.
14709. 
\par

\refstepcounter{casegroup}\label{grp:576}\par\addvspace{.55ex}\noindent\textbf{Case 576.} FB.\enspace \ref{case:327}.
\enspace \(\A{14764}\)\enspace 
14764. 
\par

\refstepcounter{casegroup}\label{grp:577}\par\addvspace{.55ex}\noindent\textbf{Case 577.} FB.\enspace \ref{case:327}.
\enspace \(\A{14765}\)\enspace 
14765. 
\par

\refstepcounter{casegroup}\label{grp:578}\par\addvspace{.55ex}\noindent\textbf{Case 578.} FB.\enspace \ref{case:324}.
\enspace \(\A{14766}\)\enspace 
14766. 
\par

\refstepcounter{casegroup}\label{grp:579}\par\addvspace{.55ex}\noindent\textbf{Case 579.} FB.\enspace \ref{case:80}.
\enspace \(\A{14768}\)\enspace 
14768. 
\par

\refstepcounter{casegroup}\label{grp:580}\par\addvspace{.55ex}\noindent\textbf{Case 580.} FB.\enspace \ref{case:328}.
\enspace \(\A{14773}\)\enspace 
14773. 
\par

\refstepcounter{casegroup}\label{grp:581}\par\addvspace{.55ex}\noindent\textbf{Case 581.} FB.\enspace \ref{case:329}.
\enspace \(\A{14777}\)\enspace 
14777. 
\enspace \(\A{14777}^{\rm op}\)\enspace 
14791. 
\par

\refstepcounter{casegroup}\label{grp:582}\par\addvspace{.55ex}\noindent\textbf{Case 582.} FB.\enspace \ref{case:330}.
\enspace \(\A{14818}\)\enspace 
14818. 
\par

\refstepcounter{casegroup}\label{grp:583}\par\addvspace{.55ex}\noindent\textbf{Case 583.} NFB.\enspace \ref{case:331}.
14839. 
\par

\refstepcounter{casegroup}\label{grp:584}\par\addvspace{.55ex}\noindent\textbf{Case 584.} FB.\enspace \ref{case:321}.
\enspace \(\A{14841}\)\enspace 
14841. 
\par

\refstepcounter{casegroup}\label{grp:585}\par\addvspace{.55ex}\noindent\textbf{Case 585.} NFB.\enspace \ref{case:332}.
14842. 
\par

\refstepcounter{casegroup}\label{grp:586}\par\addvspace{.55ex}\noindent\textbf{Case 586.} FB.\enspace \ref{case:333}.
\enspace \(\A{14865}\)\enspace 
14849, 14865. 
\par

\refstepcounter{casegroup}\label{grp:587}\par\addvspace{.55ex}\noindent\textbf{Case 587.} FB.\enspace \ref{case:334}.
\enspace \(\A{14853}\)\enspace 
14853. 
\enspace \(\A{14853}^{\rm op}\)\enspace 
14861. 
\par

\refstepcounter{casegroup}\label{grp:588}\par\addvspace{.55ex}\noindent\textbf{Case 588.} FB.\enspace \ref{case:335}.
\enspace \(\A{14862}\)\enspace 
14862. 
\par

\refstepcounter{casegroup}\label{grp:589}\par\addvspace{.55ex}\noindent\textbf{Case 589.} FB.\enspace \ref{case:336}.
\enspace \(\A{14903}\)\enspace 
14903. 
\par

\refstepcounter{casegroup}\label{grp:590}\par\addvspace{.55ex}\noindent\textbf{Case 590.} FB.\enspace Proposition~\ref{prop:lower-fibres}, source \(\A{203}\), case~\ref{grp:26}.
\enspace \(\A{203}\)\enspace 
14958. 
\par

\refstepcounter{casegroup}\label{grp:591}\par\addvspace{.55ex}\noindent\textbf{Case 591.} FB.\enspace \ref{case:337}.
\enspace \(\A{14959}\)\enspace 
14959. 
\par

\refstepcounter{casegroup}\label{grp:592}\par\addvspace{.55ex}\noindent\textbf{Case 592.} FB.\enspace Proposition~\ref{prop:lower-fibres}, source \(\A{2312}\), case~\ref{grp:2}.
\enspace \(\A{2312}\)\enspace 
14960. 
\par

\refstepcounter{casegroup}\label{grp:593}\par\addvspace{.55ex}\noindent\textbf{Case 593.} FB.\enspace Proposition~\ref{prop:lower-fibres}, source \(\A{210}\), case~\ref{grp:26}.
\enspace \(\A{210}\)\enspace 
14966. 
\par

\refstepcounter{casegroup}\label{grp:594}\par\addvspace{.55ex}\noindent\textbf{Case 594.} FB.\enspace \ref{case:338}.
\enspace \(\A{14967}\)\enspace 
14967. 
\par

\refstepcounter{casegroup}\label{grp:595}\par\addvspace{.55ex}\noindent\textbf{Case 595.} FB.\enspace Proposition~\ref{prop:lower-fibres}, source \(\A{2316}\), case~\ref{grp:2}.
\enspace \(\A{2316}\)\enspace 
14968. 
\par

\refstepcounter{casegroup}\label{grp:596}\par\addvspace{.55ex}\noindent\textbf{Case 596.} FB.\enspace \ref{case:339}.
\enspace \(\A{14972}\)\enspace 
14972. 
\par

\refstepcounter{casegroup}\label{grp:597}\par\addvspace{.55ex}\noindent\textbf{Case 597.} FB.\enspace \ref{case:287}.
\enspace \(\A{12468}^{\rm op}\)\enspace 
15321. 
\par

\refstepcounter{casegroup}\label{grp:598}\par\addvspace{.55ex}\noindent\textbf{Case 598.} FB.\enspace \ref{case:340}.
\enspace \(\A{15483}\)\enspace 
15483. 
\par

\refstepcounter{casegroup}\label{grp:599}\par\addvspace{.55ex}\noindent\textbf{Case 599.} FB.\enspace \ref{case:341}.
\enspace \(\A{15618}\)\enspace 
15596, 15618. 
\par

\refstepcounter{casegroup}\label{grp:600}\par\addvspace{.55ex}\noindent\textbf{Case 600.} FB.\enspace \ref{case:342}.
\enspace \(\A{15603}\)\enspace 
15603. 
\enspace \(\A{15603}^{\rm op}\)\enspace 
15615. 
\par
\endgroup
\subsection{Operation tables for the smaller sources}
Each string is read by rows in the order \(0,\ldots,d-1\). An isomorphic relabelling of either operation pair has the same identities. To check an indicated equality of varieties, take the complete list of elements on the target side as generators in Proposition~\ref{prop:lower-fibres}, then interchange source and target. A smaller generating tuple gives the same result.
\begingroup\small
\begin{longtable}{@{}rll@{}}
\toprule \((d,j)\)& Addition & Multiplication\\\midrule\endhead
\(3,24\)&\texttt{000011012}&\texttt{000001012}\\
\(3,37\)&\texttt{000011012}&\texttt{002002222}\\
\(3,50\)&\texttt{000011012}&\texttt{012122222}\\
\(3,51\)&\texttt{000011012}&\texttt{012222222}\\
\(3,53\)&\texttt{000011012}&\texttt{022122222}\\
\(3,60\)&\texttt{000011012}&\texttt{122222222}\\
\(4,9\)&\texttt{0000010000200003}&\texttt{0000000000010012}\\
\(4,42\)&\texttt{0000010000200003}&\texttt{0003000300133333}\\
\(4,288\)&\texttt{0100111101200103}&\texttt{0000010000000002}\\
\(4,290\)&\texttt{0100111101200103}&\texttt{0000010000020003}\\
\(4,292\)&\texttt{0100111101200103}&\texttt{0000010000000023}\\
\(4,357\)&\texttt{0100111101200103}&\texttt{1101111101231101}\\
\(4,358\)&\texttt{0100111101200103}&\texttt{1101111101201131}\\
\(4,360\)&\texttt{0100111101200103}&\texttt{1110111101201113}\\
\(4,363\)&\texttt{0100111101200103}&\texttt{1113111301231113}\\
\(4,368\)&\texttt{0100111101200103}&\texttt{1111111111011113}\\
\(4,376\)&\texttt{0100111101200103}&\texttt{1101111111213333}\\
\(4,379\)&\texttt{0100111101200103}&\texttt{2022012022222022}\\
\(4,380\)&\texttt{0100111101200103}&\texttt{2222012022222222}\\
\(4,383\)&\texttt{0100111101200103}&\texttt{2022212222222022}\\
\(4,385\)&\texttt{0100111101200103}&\texttt{2222202222222222}\\
\(4,393\)&\texttt{0000012302200303}&\texttt{0000010100000103}\\
\(4,398\)&\texttt{0000012302200303}&\texttt{0000010300000300}\\
\(4,401\)&\texttt{0000012302200303}&\texttt{0000012300000300}\\
\(4,405\)&\texttt{0000012302200303}&\texttt{0000010302000300}\\
\(4,413\)&\texttt{0000012302200303}&\texttt{0000020200000200}\\
\(4,417\)&\texttt{0000012302200303}&\texttt{0100010101000103}\\
\(4,418\)&\texttt{0000012302200303}&\texttt{0100010101020103}\\
\(4,427\)&\texttt{0000012302200303}&\texttt{0000111100000103}\\
\(4,428\)&\texttt{0000012302200303}&\texttt{0000111100000123}\\
\(4,430\)&\texttt{0000012302200303}&\texttt{0100111101000102}\\
\(4,431\)&\texttt{0000012302200303}&\texttt{0100111101000103}\\
\(4,432\)&\texttt{0000012302200303}&\texttt{0100111101020103}\\
\(4,433\)&\texttt{0000012302200303}&\texttt{0103111101030103}\\
\(4,434\)&\texttt{0000012302200303}&\texttt{0100111101000123}\\
\(4,441\)&\texttt{0000012302200303}&\texttt{0101111101011111}\\
\(4,442\)&\texttt{0000012302200303}&\texttt{0103111101031111}\\
\(4,443\)&\texttt{0000012302200303}&\texttt{0123111101231111}\\
\(4,445\)&\texttt{0000012302200303}&\texttt{0100111101003133}\\
\(4,446\)&\texttt{0000012302200303}&\texttt{0101111101013131}\\
\(4,447\)&\texttt{0000012302200303}&\texttt{0103111101033131}\\
\(4,450\)&\texttt{0000012302200303}&\texttt{0103111101233131}\\
\(4,453\)&\texttt{0000012302200303}&\texttt{0303111103033333}\\
\(4,456\)&\texttt{0000012302200303}&\texttt{0101111121213131}\\
\(4,464\)&\texttt{0000012302200303}&\texttt{0120212221220120}\\
\(4,471\)&\texttt{0000012302200303}&\texttt{2112111111112112}\\
\(4,589\)&\texttt{0000011101220123}&\texttt{0033013302330333}\\
\(4,627\)&\texttt{0000011101220123}&\texttt{0003000300133333}\\
\(4,681\)&\texttt{0000011101220123}&\texttt{0123012322333333}\\
\(4,682\)&\texttt{0000011101220123}&\texttt{0133013322333333}\\
\(4,686\)&\texttt{0000011101220123}&\texttt{0000012333333333}\\
\(4,690\)&\texttt{0000011101220123}&\texttt{0123012333333333}\\
\(4,705\)&\texttt{0000011101220123}&\texttt{0111111111111112}\\
\(4,721\)&\texttt{0000011101220123}&\texttt{0123113311331133}\\
\(4,734\)&\texttt{0000011101220123}&\texttt{0123133313333333}\\
\(4,748\)&\texttt{0000011101220123}&\texttt{0123122322232223}\\
\(4,751\)&\texttt{0000011101220123}&\texttt{0223122322232223}\\
\(4,756\)&\texttt{0000011101220123}&\texttt{0023112322333333}\\
\(4,757\)&\texttt{0000011101220123}&\texttt{0033113322333333}\\
\(4,763\)&\texttt{0000011101220123}&\texttt{0123113322333333}\\
\(4,765\)&\texttt{0000011101220123}&\texttt{0111111123333333}\\
\(4,766\)&\texttt{0000011101220123}&\texttt{0123112323333333}\\
\(4,769\)&\texttt{0000011101220123}&\texttt{0122122222223333}\\
\(4,772\)&\texttt{0000011101220123}&\texttt{0113133323333333}\\
\(4,775\)&\texttt{0000011101220123}&\texttt{0222122222223333}\\
\(4,780\)&\texttt{0000011101220123}&\texttt{0023112333333333}\\
\(4,793\)&\texttt{0000011101220123}&\texttt{0123222322232223}\\
\(4,799\)&\texttt{0000011101220123}&\texttt{0122222222223333}\\
\bottomrule
\end{longtable}\endgroup

\section{Individual case arguments}\label{app:proofs}
We use \(u\leq P\) and \(u\preccurlyeq P\) interchangeably for \(u+P\approx P\). The letters denoting elements in a table are names in the metalanguage. Several cases call a literal a \emph{linear summand} and a word of length at least two a \emph{nonlinear summand}. Where the proof uses squarefree words, linearity means that each letter occurs at most once. All substitutions in the identities have nonempty images. Tables in row-string notation follow the element order of the case, and a dash denoting an absent factor is a convention for writing two different word patterns. Each subsection uses its own symbols. Subsections written in Chinese use the same notation for word absorption, variety membership and opposite semirings.
\begingroup\small\raggedright

\subsection{四元素链来源的38个新有限基判定：完整证明与文献公式勘误}\label{case:1}
\par\smallskip\noindent\textit{705号的一组正确七条基.}\quad 
令A的加法为min，乘法为

\noindent \ensuremath{\cdot}\enspace\textbar\enspace 0\enspace\textbar\enspace 1\enspace\textbar\enspace 2\enspace\textbar\enspace 3\par
\noindent 0\enspace\textbar\enspace 0\enspace\textbar\enspace 1\enspace\textbar\enspace 1\enspace\textbar\enspace 1\par
\noindent 1\enspace\textbar\enspace 1\enspace\textbar\enspace 1\enspace\textbar\enspace 1\enspace\textbar\enspace 1\par
\noindent 2\enspace\textbar\enspace 1\enspace\textbar\enspace 1\enspace\textbar\enspace 1\enspace\textbar\enspace 1\par
\noindent 3\enspace\textbar\enspace 1\enspace\textbar\enspace 1\enspace\textbar\enspace 1\enspace\textbar\enspace 2\par
以下七条恒等式构成A的完整基：

1. \(x\leq x^{2}\)；

2. xy\ensuremath{\approx}yx；

3. \(x^{4}\approx x^{3}\)；

4. \(xz\leq x^{2}+yz\)；

5. \(x^{2}\leq x+xy\)；

6. \(xz\leq x+xy+zt\)；

7. \(xyz+u+uv\approx xyz+u\)。

前六条同文献(111)--(116)，第七条是本记录独立证明的替换式，只有五个变量。 七条的有效性由该表直接检验。特别地，任意至少三因子的乘积属于\{0,\allowbreak{}1\}； 若r\ensuremath{\in}\{0,\allowbreak{}1\}，则数值上 \(min(r, u, uv)=min(r, u)\)，故第七条成立。

完整性证明。设q是词，P是非空有限词和，A满足q\ensuremath{\leq}P。 两元素子代数\{0,\allowbreak{}1\}的加法为min、乘法为max。若P没有内容包含于con(q)的词， 把con(q)中的变量赋0、其余赋1，就有\(q=0\)、\(P=1\)，违反q\ensuremath{\leq}P。 故总存在内容见证v\ensuremath{\in}P，使con(v)\ensuremath{\subseteq}con(q)。

先设P含长度至少3的词w。将w分成三个非空块，用第七条代入\(u=v\)、\(v=q\)\ensuremath{{}^{3}}可得 \(w+vq^{3}\leq w+v\leq P\)。交换律和x\ensuremath{{}^{4}}\ensuremath{\approx}\(x^{3}\)给出 \(vq^{3}\approx q^{3}\)，因为con(v)\ensuremath{\subseteq}con(q)。 又由第一条及其乘法后果有 \(q\leq q^{2}\leq q^{3}\)。所以 \(q\leq q^{3}\approx vq^{3}\leq w+vq^{3}\leq P\)。以上各步是任意变量数下的等式演绎。

再设P的每个词长度至多2。子代数\{1,\allowbreak{}2,\allowbreak{}3\}中，三因子乘积恒为1； 在全部变量赋3时，线性词值3，二次词值2。 因此q的长度至多2，否则\(q=1\)而P\ensuremath{\geq}2（数值序），违反不等式。 如果q为二次词，那么q的每个变量必须出现在P的某个二次词中： 若变量x没有这样的出现，把x赋2、其余变量赋3，就有\(q=1\)，而P数值至少2。

若\(q=x\)，内容见证v只能为x或\(x^{2}\)，分别由自反性或第一条得q\ensuremath{\leq}P。 若\(q=x\)\ensuremath{{}^{2}}且q不是P的现成词，内容见证只能是x；二次覆盖给出某词xt\ensuremath{\in}P， 用第五条得\(x^{2}\)\ensuremath{\leq}x+xt\ensuremath{\leq}P。 若\(q=xy\)且x\ensuremath{\ne}y，q若在P中即结束。否则内容见证只能为x、y、\(x^{2}\)、\(y^{2}\)。 二次覆盖给出xt, ys\ensuremath{\in}P。见证为x或y时，第六条（必要时交换x, y）给出 \(xy\leq x+xt+ys\leq P\)或其对称式。见证为\(x^{2}\)或\(y^{2}\)时，第四条给出 \(xy\leq x^{2}+sy\leq P\)或对称式。 这穷尽所有情形。任意有效多项式恒等式逐项化为双向吸收，故七条基完整。

文献式(117)原文为 \(x_{1}x_{2}x_{3}+y_{1}\approx x_{1}x_{4}x_{5}+y_{1}y_{2}\)。 取减1后的赋值 \((x_{1},x_{2},x_{3},x_{4},x_{5},y_{1},y_{2})=(0,\allowbreak{}0,\allowbreak{}0,\allowbreak{}0,\allowbreak{}1,\allowbreak{}0,\allowbreak{}1)\)， 其左边为0、右边为1。这是有效性反例，并非未搜索到证明。 本记录用上述完整证明修复本批转移所需的来源结论；不声称证明论文其他所有命题。

\par\smallskip\noindent\textit{627号的四条基及独立完整性证明.}\quad 
令B的加法为min，乘法为

\noindent \ensuremath{\cdot}\enspace\textbar\enspace 0\enspace\textbar\enspace 1\enspace\textbar\enspace 2\enspace\textbar\enspace 3\par
\noindent 0\enspace\textbar\enspace 0\enspace\textbar\enspace 0\enspace\textbar\enspace 0\enspace\textbar\enspace 3\par
\noindent 1\enspace\textbar\enspace 0\enspace\textbar\enspace 0\enspace\textbar\enspace 0\enspace\textbar\enspace 3\par
\noindent 2\enspace\textbar\enspace 0\enspace\textbar\enspace 0\enspace\textbar\enspace 1\enspace\textbar\enspace 3\par
\noindent 3\enspace\textbar\enspace 3\enspace\textbar\enspace 3\enspace\textbar\enspace 3\enspace\textbar\enspace 3\par
完整基为：\(x\leq x^{2}\)、xy\ensuremath{\approx}yx、\(xy\leq x^{2}+y^{2}\)、\(xyzt\leq xyz\)。 这是文献(121)--(124)。四条有效性均可直接查表。

为独立核对完整性，设B满足q\ensuremath{\leq}P，令D为P中所有内容包含于con(q)的词之和。 把con(q)外的变量赋3，它们所在的整词都因3为乘法吸收元而取3，且3是加法最小元。 在con(q)上只用子代数\{0,\allowbreak{}1,\allowbreak{}2\}赋值，就得到q\ensuremath{\leq}D；D不能空，否则把q中变量全赋0即矛盾。

若D含长度至少3的词w，第四条反复应用给出 \(wq\leq w\)。 由于con(w)\ensuremath{\subseteq}con(q)，用第一条、交换律逐次增加q已有变量的出现次数，有 \(q\leq wq\)。 所以q\ensuremath{\leq}wq\ensuremath{\leq}w\ensuremath{\leq}P。

若D所有词长至多2，把q中变量全赋2可知q长至多2。 若q为二次词，且其变量x不出现在D的任何二次词中，把x赋1、其余赋2， 则\(q=0\)而D数值至少1，矛盾。因此每个q变量都有D中的二次覆盖。 若\(q=x\)，D中的词只能为x或\(x^{2}\)，第一条即可。 若\(q=x\)\ensuremath{{}^{2}}，二次覆盖必包含\(x^{2}\)，故直接吸收。 若\(q=xy\)且x\ensuremath{\ne}y，若xy已在D中即可；否则二次覆盖和内容限制迫使\(x^{2}\)、\(y^{2}\)同时在D中， 由第三条得xy\ensuremath{\leq}\(x^{2}\)+\(y^{2}\)\ensuremath{\leq}P。于是任意有效吸收均可推导，四条基完整。

\par\smallskip\noindent\textit{精确等簇转移的38个目标.}\quad 
23型与705号等簇：

38,\allowbreak{}39,\allowbreak{}40,\allowbreak{}50,\allowbreak{}52,\allowbreak{}199,\allowbreak{}293,\allowbreak{}294,\allowbreak{}295,\allowbreak{}320,\allowbreak{}324,\allowbreak{}466,\allowbreak{}5928,\allowbreak{}5993,\allowbreak{}6348,\allowbreak{}12194,\allowbreak{}12276,\allowbreak{}13197,\allowbreak{}13835,\allowbreak{}13844,\allowbreak{}13845,\allowbreak{}13934,\allowbreak{}13937。

15型与627号等簇：

2944,\allowbreak{}12248,\allowbreak{}13060,\allowbreak{}13103,\allowbreak{}13106,\allowbreak{}13187,\allowbreak{}13194,\allowbreak{}13195,\allowbreak{}13196,\allowbreak{}13834,\allowbreak{}13843,\allowbreak{}13888,\allowbreak{}13932,\allowbreak{}13933,\allowbreak{}13936。

两向等簇意味着它们满足完全相同的恒等式，因此上述同一组四条或七条基也就是各目标的完整基。

\subsection{Ten complete identities for semiring 51}\label{case:2}
In the constant-free signature (+, multiplication), relative to the additively idempotent semiring axioms, a complete finite basis is

{\small\par\smallskip
\CaseLine{\(xy=yx; xyz=xyzz; xyzt\leq xyz; xy\leq x+xz+yt;\)}
\CaseLine{\(xxx\leq x+yyy; xxy\leq xy+zzz; xy\leq xxy; xy\leq xx+yz;\)}
\CaseLine{\(x\leq xxx; x\leq xx+y+yz.\)}
\smallskip}

Here \(u\leq P\) abbreviates \(u+P=P\). This proof has no bound on variables, word length or number of summands.

\par\smallskip\noindent\textit{Exact values and absorption.}\quad 
Addition is 00004/\allowbreak{}01114/\allowbreak{}01214/\allowbreak{}01134/\allowbreak{}44444, with 2 and 3 incomparable below \(1<0<4\). Multiplication is 00000/\allowbreak{}00000/\allowbreak{}00000/\allowbreak{}00020/\allowbreak{}00004. Every nonlinear word is 4 when all its content is assigned 4. Otherwise a quadratic is 2 when both occurrences are assigned 3, and every other nonlinear word is 0. Long words therefore retain content rather than being constant.

For an original polynomial P let L be its literal variables, E its unordered quadratic edges and V their endpoints. Put \(G=true\) exactly when P has an original word of length at least three. For a target w, let C mean that some original u has c(u) contained in c(w). Then the exact absorption criterion is:

\textbullet{} Literal x: original x, or C and (G or L intersects V).

\textbullet{} Long w: C and G.

\textbullet{} Quadratic w: C and (G or c(w) contained in V).

This distinguishes the literal recovery enabled by overlap from long recovery, which still requires an original long word. In particular \(xxx\leq x+xy\) fails at \(x=y=3\), since its values are 0 and 1.

Necessity of C for any nonoriginal target follows by assigning its content 4 and all outsiders 2: every original summand is at most 0, but the target is 4. Suppose G fails, and assign V to 3 and outsiders

2. Every original quadratic is 2 and every literal is 2 or 3, so the right value is at most 1. A long target has value 0. A quadratic with an endpoint outside V also has value 0. For a missing literal x for which C holds, an original square xx exists, so x belongs to V. If moreover L and V are disjoint, every right literal is 2; hence \(P=2\), while \(x=3\). Thus the same assignment separates the remaining failed literal condition. These assignments cover all failed cases.

For sufficiency when G holds, \(P\geq 0\). The contained original witness provides 4 when the target is 4; all other target values are below 0. Without G, any valid nonlinear target is quadratic with both endpoints in V. If P does not dominate 0, all original quadratics must have value 2, forcing every V variable to 3. The target is then 2. The contained witness again supplies 4 in the exceptional case.

For a new literal x without G, C supplies an original square xx. If x is not 3, its square already dominates x: square values are 0 at \(x=0\),1,\allowbreak{}2 and 4 at \(x=4\). At \(x=3\) the square is 2. An original overlap y and yz has join at least 1: \(yz=2\) forces \(y=z=3\) and \(3+2=1\); otherwise yz is 0 or 4. Hence the right side dominates \(x=3\) too. This proves the literal branch without strengthening it to long availability.

\par\smallskip\noindent\textit{Unrestricted finite derivation.}\quad 
The two word identities normalize long words by content. Commute a chosen variable to the end, duplicate it using \(xyz=xyzz\) with two nonempty preceding blocks, and remove surplus repetitions with the reverse law while retaining length at least three. Canonical long forms are sorted content for at least three distinct variables, xxy for two, and xxx for one. In particular \(g^{3}=g\) for any long word g.

If G holds, fix an original long g and a contained original u. Upgrade a literal u to \(u^{3}\) using \(xxx\leq x+yyy\) with \(y=g\); upgrade a quadratic ab to aab using \(xxy\leq xy+zzz\) with \(z=g\). Retain a long u. Append the missing target variables by \(xyzt\leq xyz\) and normalize. This derives every long target. Derive a quadratic ab from aab by \(xy\leq xxy\). For a new literal x with G, first derive xxx by this same long-target construction, then apply \(x\leq xxx\). All contained witnesses are chosen from original P.

When G fails, a new quadratic ab has an original endpoint loop or literal. A loop aa and an original edge bc derive ab by \(xy\leq xx+yz\). A literal a and original edges az and bt derive ab by \(xy\leq x+xz+yt\); both incidence witnesses are retained. Commutativity handles the other endpoint. For a new literal x, the criterion supplies original xx and original y, yz, so apply \(x\leq xx+y+yz\) directly. This does not introduce a long word. Every right summand is original or explicitly derived.

All ten laws satisfy the exact criterion and are independently checked over all assignments. Every substitution is nonempty. Summandwise absorption after distributive expansion proves arbitrary identities in both directions, establishing completeness. Transfers to other catalogue types require independent HSP checks in both directions.

\subsection{有向二次冲突与83类的九条完整恒等式基}\label{case:3}
\par\smallskip\noindent\textit{1. 半环及九条恒等式.}\quad 
使用无常数签名 \((+,\cdot)\)，相对于加法幂等半环公理。 令 \(B=\{z, p, b, a, t\}\)，其中 \(p, b, a\) 互不可比，均严格低于 \(z\)，而 \(z<t\)；加法为此偏序中的并。乘法只有

\CaseDisplay{ab=p,\qquad tt=t}

不是 \(z\)，其余乘积均为 \(z\)。字母仅命名代数元素，非项常数。 目录83的标签依次为 \((z, p, b, a, t)=(0,\allowbreak{}1,\allowbreak{}2,\allowbreak{}3,\allowbreak{}4)\)。

定理1。 半环 \(B\) 的恒等式由以下九条恒等式给出完整基：

\CaseDisplay{\begin{aligned}
x^3&\approx x^2, &(M1)\\
x^2y^2&\approx y^2x^2, &(M2)\\
xyz&\approx x^2y^2z^2, &(M3)\\
x&\leq x^2, &(A1)\\
x^2y^2&\leq x^2, &(A2)\\
y^2&\leq x^2+y, &(A3)\\
x^2y^2z^2&\leq xy+yz, &(D1)\\
x^2y^2&\leq x+xy, &(D2)\\
x^2y^2&\leq y+xy. &(D3)
\end{aligned}}

这里 \(u\leq v\) 缩写 \(u+v\approx v\)。不主张九条基不可约或最小。

半环性与声性。 任意三元乘积等于 \(t\) 当且仅当三个因子均为 \(t\)；否则为 \(z\)，与括号无关。固定左因子为 \(a\) 的乘法映射把 \(b\) 送到 \(p\)，其余送到 \(z\)；固定左因子为 \(t\) 的映射只把 \(t\) 送到 \(t\)，其余送到 \(z\)。二者均保持并，其余左乘映射恒为 \(z\)。右乘同理，故分配律成立。

每个平方以及每个长度至少三的词只取 \(z, t\)：当全部变量取 \(t\) 时为 \(t\)，否则为 \(z\)。(M1)--(M3)、(A1)--(A3)立即由此核对。对(D1)，两条连续有向边 \(xy\) 和 \(yz\) 不可能同时取 \(p\)，因为这分别要求 \(y=b\) 和 \(y=a\)。右侧因此至少为 \(z\)；当左侧为 \(t\) 时三个变量全为 \(t\)，右侧也为 \(t\)。对(D2)、(D3)，若 \(xy=p\) 则 \(x=a, y=b\)，所以加入任一端点均使并为 \(z\)；其余情况及左侧为 \(t\) 的情况同样满足不等式。

\par\smallskip\noindent\textit{2. 单词正规形.}\quad 
对有限非空变量集 \(C\)，令 \(Q_C=\prod_{x\in C}x^2\)，变量按任一固定顺序列出。

引理2。 (M1)--(M3)把每个词化为以下三种形式之一：

1. 一个变量 \(x\)；

2. 两个不同变量的有向二次词 \(xy\)；

3. \(Q_C\)，其中 \(C\) 为原词内容。第三种包括平方词和所有长度至少三的词。

证明。 (M1)导出 \(x^n\approx x^2\) 对所有 \(n\geq2\)。由(M3)代入 \((x, x, y)\) 得

\CaseDisplay{x^2y\approx x^4y^2\approx x^2y^2.}

三个字母的词由(M3)化为平方因子的积。再逐字追加：在最后一个平方因子后用上式把新字母平方化。(M2)允许任意排列平方因子，再用(M1)合并重复，得到 \(Q_C\)。平方词本身已经是 \(Q_{\{x\}}\)。

引理3。 若非空 \(D\subseteq C\)，则从九条恒等式可推出

\CaseDisplay{w\leq Q_{c(w)},\qquad Q_C\leq Q_D,\qquad Q_C\approx (\prod_{x\in C}x)^2.}

第一式由逐字应用(A1)，再用正规化得出；第二式由(A2)逐个删除平方因子并结合(M2)得出；第三式由引理2得出（单变量时直接成立）。

\par\smallskip\noindent\textit{3. 任意变量数的精确吸收判据.}\quad 
把多项式 \(P\) 分配并正规化。记其一次变量集为 \(L\)，有向二次词集为 \(E\)，平方型词集为 \(N=\{Q_C\}\)。令

\CaseDisplay{H=\{x:xy\in E\},\qquad T=\{y:xy\in E\}.}

称 \(P\) 有冲突，若

\CaseDisplay{N\ne\varnothing\quad\text{或}\quad H\cap T\ne\varnothing\quad\text{或}\quad L\cap(H\cup T)\ne\varnothing.}

对已正规化的目标词 \(w\)，定义内容条件

\CaseDisplay{\text{(C)}\qquad\text{存在 }u\in P\text{，使 }c(u)\subseteq c(w).}

引理4。 对任意非空多项式 \(P\) 和词 \(w\)，

\CaseDisplay{B\models w\leq P
\quad\Longleftrightarrow\quad
w\in P\ \text{或} (P\text{ 有冲突且满足(C)}).}

必要性。 如果(C)不成立，将 \(c(w)\) 中变量全部赋 \(t\)，其他变量赋 \(z\)。目标为 \(t\)，而每个右词均含一个取 \(z\) 的变量，因此右侧为 \(z\)，矛盾。

假设 \(P\) 没有冲突，则 \(N\) 为空，\(H, T, L\) 两两不交。赋值 \(H\mapsto a\)、\(T\mapsto b\)、所有其余变量 \(\mapsto p\)。每个右词均取 \(p\)，故 \(P=p\)。平方型词和长度至少三的词为 \(z\not\leq p\)，因而目标只能为一次词或不同变量组成的二次词。

若目标是 \(x\)，(C)给出内容包含于 \(\{x\}\) 的右词；无冲突多项式中这种词只能是一次词 \(x\)，故 \(x\in P\)。若目标为 \(xy\)，上述赋值下必须取 \(p\)，所以 \(x\in H, y\in T\)。由(C)，右侧含内容包含于 \(\{x, y\}\) 的词；它不可能为一次词，因为 \(L\) 与 \(H\cup T\) 不交，也不可能为反向边 \(yx\)，因为 \(H\cap T\) 为空。因此只能是 \(xy\) 本身。必要性证毕。

充分性，同时给出等式推导。 若 \(w\in P\)，显然。设 \(P\) 有冲突。我们先从九条式子推出某个平方型词 \(Q_D\leq P\)：

\textbullet{} 若已有平方型求和词，直接取它。

\textbullet{} 若某变量同时为边的终点和起点，\(P\) 含 \(xy+yz\)，由(D1)取 \(Q_{\{x, y, z\}}\)；允许 \(x=z\)，正规化处理重复。

\textbullet{} 若一次变量同时出现于某条边的起点或终点，由(D2)或(D3)得到该边内容对应的平方型词。

取(C)见证 \(u\in P\)。由引理3，把已得 \(Q_D\) 写成 \(v^2\)。在(A3)中代入 \(x=v, y=u\)，可得

\CaseDisplay{u^2\leq v^2+u\leq P.}

再次正规化得 \(Q_{c(u)}\leq P\)。由 \(c(u)\subseteq c(w)\) 及引理3，

\CaseDisplay{w\leq Q_{c(w)}\leq Q_{c(u)}\leq P.}

所有步骤均是九条式子的代换、乘法上下文或加法保序推导，变量数不限。

\par\smallskip\noindent\textit{4. 完整性及目录转移.}\quad 
定理1的完整性证明。 把任何在 \(B\) 成立的恒等式 \(p\approx q\) 分配并正规化为 \(P, Q\)。每个 \(P\) 中的词在 \(B\) 里被 \(Q\) 吸收，由引理4可从九条基推出；相加得到 \(P\leq Q\)。反向同理，故推出 \(P\approx Q\)。

目录83正是 \(B\)。目录97、183、1798与83的两个HSP包含方向均通过无词长上限的完整词函数闭合检查。因此它们生成同一簇，同样具有上述九条完整基。等簇检查并未使用拟议九条基作为成员测试，因而不循环论证。

乘法对偶也具有相反词的九条基；交换 \(a, b\)、固定其余元素给出 \(B\) 到 \(B^{op}\) 的同构。本批四型属于同一个自对偶等簇类。

\subsection{端点标记的有向二次半环：统一九条完整基与八个目录类型}\label{case:4}
\par\smallskip\noindent\textit{1. 四种加法序.}\quad 
对参数 \(\eta,\theta\in\{0,\allowbreak{}1\}\)，定义 \(B_{\eta,\theta}\)，载体为 \(\{z, p, b, a, t\}\)。 所有 \(a, b, p\) 均严格低于 \(z<t\)；当且仅当 \(\eta=1\) 时再加 \(a<p\)，当且仅当 \(\theta=1\) 时再加 \(b<p\)，除此以外没有严格比较。加法为并，乘法只有

\CaseDisplay{ab=p,\qquad tt=t}

不是 \(z\)，所有其他二元乘积为 \(z\)。当 \(\eta=\theta=1\) 时，\(a+b=p\)；其他情况下不可比非顶元素的并由上述偏序确定。 参数分别表示头端点和尾端点能否低于乘积值。元素名字不是项中的常数。

四种代数都是无常数签名的加法幂等半环：所有三元乘积当且仅当全为 \(t\) 时取 \(t\)，否则取 \(z\)，所以乘法结合。仅非恒定的左乘映射是 \(a\) 左乘和 \(t\) 左乘，前者把 \(b\) 送到 \(p\)、其余送到 \(z\)，后者只把 \(t\) 送到 \(t\)；检查上述四种并序可见均保持并，右乘同理。

\(B_{0,\allowbreak{}0}\) 即前一批的83型。新来源对应为：

\noindent 参数\enspace\textbar\enspace 目录来源\enspace\textbar\enspace \((z, p, b, a, t)\)标签\par
\noindent \((0,\allowbreak{}1)\)\enspace\textbar\enspace 888\enspace\textbar\enspace \((0,\allowbreak{}2,\allowbreak{}3,\allowbreak{}1,\allowbreak{}4)\)\par
\noindent \((1,\allowbreak{}1)\)\enspace\textbar\enspace 84\enspace\textbar\enspace \((0,\allowbreak{}1,\allowbreak{}2,\allowbreak{}3,\allowbreak{}4)\)\par
乘法对偶 \(B_{0,\allowbreak{}1}^{op}\) 与 \(B_{1,\allowbreak{}0}\) 同构，映射交换 \(a, b\)、固定 \(z, p, t\)。

\par\smallskip\noindent\textit{2. 统一的九条基.}\quad 
相对于加法幂等半环公理，所有四种代数共有以下七条恒等式：

\CaseDisplay{\begin{aligned}
x^3&\approx x^2, &(M1)\\
x^2y^2&\approx y^2x^2, &(M2)\\
xyz&\approx x^2y^2z^2, &(M3)\\
x&\leq x^2, &(A1)\\
x^2y^2&\leq x^2, &(A2)\\
y^2&\leq x^2+y, &(A3)\\
x^2y^2z^2&\leq xy+yz. &(D)
\end{aligned}}

再根据头参数选择一条：

\CaseDisplay{\begin{cases}
x^2y^2\leq x+xy, &\eta=0;\\
xy\leq xz+uy+x, &\eta=1;
\end{cases}\quad\textup{(\ensuremath{H})}}

根据尾参数选择一条：

\CaseDisplay{\begin{cases}
x^2y^2\leq y+xy, &\theta=0;\\
xy\leq xz+uy+y, &\theta=1.
\end{cases}\quad\textup{(\ensuremath{T})}}

因此每个参数组合都恰给出一组九条完整基。这里 \(u\) 是第四个普通变量；不主张基不可约。

声性。 平方和所有长度至少三的词只取 \(z, t\)，当且仅当所有出现变量赋 \(t\) 时取 \(t\)。(M1)--(A3)随之成立。(D)右侧两个连续乘积不能同时取 \(p\)，因为前者要求 \(y=b\)、后者要求 \(y=a\)；所以右侧至少为 \(z\)，全 \(t\) 情形另直接成立。

对(H)的冲突形式，\(xy=p\) 时 \(x=a\)，而 \(\eta=0\) 保证 \(a+p=z\)；右侧不低于 \(z\)，左侧为 \(t\) 时右侧也为 \(t\)。对(T)的冲突形式同理。

对(H)、(T)的补边形式，若 \(xz\) 或 \(uy\) 取 \(t\)，右侧为 \(t\)。若任一取 \(z\)，右侧至少为 \(z\)，只有目标为 \(t\) 时还需检查；此时被加作一次项的 \(x\) 或 \(y\) 本身为 \(t\)，右侧为 \(t\)。剩下情形两个乘积都为 \(p\)，于是 \(x=u=a, z=y=b\)，目标 \(xy=p\)，不大于右侧。故所有列出的恒等式均有效。

\par\smallskip\noindent\textit{3. 正规形与内容平方.}\quad 
对非空有限变量集 \(C\)，记 \(Q_C=\prod_{x\in C}x^2\)。 (M1)--(M3)把任何词化为：一次词、不同变量的有向二次词、或 \(Q_C\)。平方词和长度至少三的词均化为第三种。

详细推导：由(M1)得到 \(x^n=x^2\) 对 \(n\geq2\)；(M3)代入 \((x, x, y)\) 给出 \(x^2y=x^4y^2=x^2y^2\)。先平方化三个字母，再逐字追加并平方化，利用(M2)排列平方因子，再由(M1)合并重复。

由(A1)、(A2)进一步得到

\CaseDisplay{w\leq Q_{c(w)},\qquad Q_C\leq Q_D\quad(\varnothing\ne D\subseteq C),
\qquad Q_C=(\prod_{x\in C}x)^2.\quad\textup{(\ensuremath{1})}}

这些是等式推导，不是仅在给定五元素表中成立的语义替换。

\par\smallskip\noindent\textit{4. 任意变量数的吸收判据.}\quad 
把非空多项式 \(P\) 正规化，令 \(L\) 为其一次变量集，\(E\) 为其有向二次词集，\(N\) 为其平方型词集，并记

\CaseDisplay{H=\{x:xy\in E\},\qquad T=\{y:xy\in E\}.}

称 \(P\) 有冲突，若下列任一条件成立：

1. \(N\ne\varnothing\)；

2. \(H\cap T\ne\varnothing\)；

3. \(\eta=0\) 且 \(L\cap H\ne\varnothing\)；

4. \(\theta=0\) 且 \(L\cap T\ne\varnothing\)。

引理。 对正规词 \(w\)，\(B_{\eta,\theta}\models w\leq P\) 当且仅当以下任一项成立：

(i) \(w\in P\)；

(ii) \(P\) 有冲突，且存在 \(v\in P\) 使 \(c(v)\subseteq c(w)\)；

(iii) \(P\) 无冲突，\(w=xy\)，\(x\in H, y\in T\)，并且

\CaseDisplay{(\eta=1\text{ 且 }x\in L)\quad\text{或}\quad(\theta=1\text{ 且 }y\in L).}

必要性证明。 对任意有效吸收，必须存在内容见证 \(c(v)\subseteq c(w)\)：否则赋 \(c(w)\) 中变量为 \(t\)、其余为 \(z\)，目标为 \(t\) 而每个右词为 \(z\)，矛盾。

设 \(P\) 无冲突，把 \(H\) 赋 \(a\)、\(T\) 赋 \(b\)、所有其他变量赋 \(p\)。因为 \(H\cap T=\varnothing\)，赋值一致。每条边取 \(p\)；允许作为一次项出现的头端点低于 \(p\)，允许的尾端点也低于 \(p\)，其余一次项取 \(p\)。因此右侧不大于 \(p\)；若没有边，则所有一次项均为 \(p\)，右侧仍为 \(p\)。

任何平方型目标都取 \(z\not\leq p\)，故目标不可能是平方型。若目标是一次词 \(x\)，内容见证只能是右侧的一次词 \(x\)，即(i)。若目标为不同变量的二次词 \(xy\)，该赋值下必须取 \(p\)，所以 \(x\in H, y\in T\)。内容见证在无冲突时只能为原边 \(xy\)，或一次词 \(x\)，或一次词 \(y\)；反向边 \(yx\) 会产生 \(H\cap T\)，故不可能。当原边不存在时，一次词 \(x\) 必须满足 \(\eta=1\)，一次词 \(y\) 必须满足 \(\theta=1\)，得到(iii)。

注意右侧在某些赋值下可能严格低于 \(p\)；上述分离只使用``右侧不大于 \(p\)''及 \(z\not\leq p\)，并不错误断言它恒等于 \(p\)。

充分性与完整推导。 (i)平凡。若(ii)成立，平方型求和词直接提供一个 \(Q_D\leq P\)；连续有向边由(D)提供，禁止的一次端点重合由(H)或(T)的冲突形式提供。写 \(Q_D=s^2\)，对内容见证 \(v\) 在(A3)代入 \(x=s, y=v\)，得到

\CaseDisplay{v^2\leq s^2+v\leq P.}

于是 \(Q_{c(v)}\leq P\)。结合(1)和内容包含得到 \(w\leq Q_{c(w)}\leq Q_{c(v)}\leq P\)。

若(iii)成立，取原边 \(xz\in P\) 和 \(uy\in P\)。标记头端点时，用(H)的补边式得到 \(xy\leq xz+uy+x\leq P\)；标记尾端点时用(T)得到 \(xy\leq xz+uy+y\leq P\)。允许变量相同的代换，无需引入额外假设。

\par\smallskip\noindent\textit{5. 完整性与八个目录编号.}\quad 
定理。 第2节的九条式子是 \(B_{\eta,\theta}\) 的完整基。

证明：对任何有效等式 \(p\approx q\)，分配并正规化两端。每个左侧词在代数里被右侧吸收，引理逐词给出九条式子的推导。反向同理，故等式可推导。证明允许任意多个变量，完整性不由有限回归推测。

精确等簇转移给出：

\noindent 来源\enspace\textbar\enspace 新判定为有限基的目录编号\par
\noindent \(B_{1,\allowbreak{}1}\)，84\enspace\textbar\enspace 84、13901\par
\noindent \(B_{0,\allowbreak{}1}\)，888\enspace\textbar\enspace 98、888、13273\par
\noindent \(B_{0,\allowbreak{}1}^{op}\)\enspace\textbar\enspace 141、184、13252\par
每个转移均核验两个完整HSP方向。对偶来源可直接反转 \(B_{0,\allowbreak{}1}\) 的每个词得到完整基；不需要假定上下端点名称与目录标签一致。

\par\smallskip\noindent\textit{6. 证据及范围.}\quad

\subsection{Three finite bases from long-word content and independent short roles}\label{case:5}
For each of 91,\allowbreak{}92,\allowbreak{}93 separately, in the signature without constants, take the semiring axioms, the four identities

{\small\par\smallskip
\CaseLine{\(xyz=yxz, xyz=xzy, xxy=xyy, xxyz=xyz,\)}
\smallskip}

and all valid absorptions over \(B=2^{11}-1+22=2069\) variables, word length at most 2071, with at most 10 upper words. This is a specified finite family, not an enumerated small list. Its validity is decidable by the given finite table. Its unrestricted completeness is proved below.

All three multiplication tables have rows 00000,\allowbreak{}00000,\allowbreak{}00000,\allowbreak{}00110,

00004. A quadratic has value 4 if both entries are 4, value 1 if its head is 3 and second is 2 or 3, and value 0 otherwise. Every word of length at least 3 is 4 iff all its variables are 4, and otherwise0. Induction proves this for every length, including repeated variables. Literals retain values. This verifies the four displayed laws.

The first two laws permute any word of length at least 3:swap its first two positions with \(xyz=yxz\), or a later adjacent pair with \(xyz=xzy\) using a nonempty prefix. For length at least 4, move any duplicate pair to the front and apply \(xxyz=xyz\), with x that singleton and y, z the remaining nonempty factors. Continue until either the word is squarefree or has length 3. For two-variable content at length 3, \(xxy=xyy\) and permutation identify both multiplicity patterns. Thus a canonical long word is its sorted content if it has at least 3 variables, a square of the smaller variable followed by the larger if it has 2, or a cube if it has 1. Literals and quadratic words are retained exactly. In particular long words never collapse to literal or quadratic ones.

\par\smallskip\noindent\textit{Exact absorption criterion.}\quad 
Let P be original nonempty words. Its literals have variables L, its quadratic heads H, quadratic seconds T, and J is its collection of long words. The proper principal ideals are as follows:

\noindent Core\enspace\textbar\enspace Ideal of 0\enspace\textbar\enspace Ideal of 1\enspace\textbar\enspace Ideal of 2\enspace\textbar\enspace Ideal of 3\par
\noindent 91\enspace\textbar\enspace \{0,\allowbreak{}1,\allowbreak{}2,\allowbreak{}3\}\enspace\textbar\enspace \{1\}\enspace\textbar\enspace \{2,\allowbreak{}3\}\enspace\textbar\enspace \{3\}\par
\noindent 92\enspace\textbar\enspace \{0,\allowbreak{}1,\allowbreak{}2,\allowbreak{}3\}\enspace\textbar\enspace \{1,\allowbreak{}3\}\enspace\textbar\enspace \{2,\allowbreak{}3\}\enspace\textbar\enspace \{3\}\par
\noindent 93\enspace\textbar\enspace \{0,\allowbreak{}1,\allowbreak{}2,\allowbreak{}3\}\enspace\textbar\enspace \{1,\allowbreak{}2,\allowbreak{}3\}\enspace\textbar\enspace \{2,\allowbreak{}3\}\enspace\textbar\enspace \{3\}\par
Write U for the ideal of 1 appropriate to the source. Define independent variable domains initially\{0,\allowbreak{}1,\allowbreak{}2,\allowbreak{}3,\allowbreak{}4\}, intersecting with U at each literal occurrence, with\{3\} at each quadratic head, and with\{2,\allowbreak{}3\} at each quadratic second. The difficult-ideal system is impossible iff J is nonempty or some variable's domain is empty. If possible, all assignments choosing values independently from those domains put P in U, exactly:the allowed quadratic relation is the Cartesian product \{3\} times\{2,\allowbreak{}3\}, and no long word value belongs to U.

Then \(w\leq P\) holds exactly when:

1. Some original u in P has content(u) contained in content(w).

2. If w is nonlinear, P contains an original nonlinear word.

3. Either the domain system is impossible, or w is literal with its variable's domain contained in U, or w is quadratic with its head domain contained in\{3\} and second domain contained in\{2,\allowbreak{}3\}. If the system is possible, no long target is allowed.

Proof:ideal 0 membership means a word has at least one non4 letter. Every original has that property only if the target does, precisely when one original's content is contained in the target content. For the converse separation put all target variables 4 and all outside variables 0. If no original is contained, all originals become 0 and the target 4. The ideals of 2 and 3 contain no nonlinear values; an original nonlinear makes both antecedents impossible. If all originals are literal, the content condition already protects a literal target, and a nonlinear target is separated by assigning every variable 3. The remaining ideal is exactly the independent domain system described above. A feasible system cannot put a long target in U because its values are 0 or 4. For a quadratic the two domains must lie in the respective factors; coinciding target variables cause no exception because a square lies in U exactly at value 3. Literal targets have the stated domain test. The four proper principal ideals suffice to separate every failed order comparison, so this proves necessity and sufficiency without any bound on variables, lengths, or polynomial size.

\par\smallskip\noindent\textit{Bounded ORIGINAL witnesses.}\quad 
For any valid absorption retain one contained original u. If the target is nonlinear, also retain one nonlinear original. To preserve the third condition, if there is a long original retain one. Otherwise if a domain is empty, retain the at most 3 originals giving its literal, head, and second restrictions. Otherwise the target has at most 2 variables: retain at most one original for each of its three possible roles, at most 6 words. The selected polynomial Q therefore has at most 8 words;10 is a common safe bound. All are original words of P. They preserve the containment, nonlinearity and domain conditions, so \(w\leq Q\) is valid. No long target is replaced by a short target during this argument. The original contained-content witness is essential.

\par\smallskip\noindent\textit{Reversible compression and completeness.}\quad 
Normalize w and Q by the four laws. There are at most 11 words. Reserve every variable occurring in any literal or quadratic among these words, at most 22 variables. Other variables occur only in long words. Group them by their nonempty incidence subset of the 11 words, at most 2047 groups. Replace each group by one new variable, keeping reserved variables distinct, and normalize. A normalized long word has length at most its content size plus 2, so the result lies in the displayed finite family. It is valid because identification preserves an identity.

For the inverse substitute the nonempty word listing a group's members for its variable. No literal or quadratic can involve such a group; every expanded long word remains long and has exactly its original content, including all group members with COMMON incidence. The four laws normalize it back to its original normal form. Singleton variables are restored individually. Thus the finite-family instance derives \(w\leq Q\), then \(w\leq P\). Reversing normalization derives the original absorption. Every valid polynomial identity decomposes into its summand absorptions in both directions. This proves completeness for arbitrary variables, word lengths and polynomial sizes, and finite basability of all 3 cores.

Direct mutual-HSP certificates can transfer to nonisomorphic equal varieties, with word reversal for opposites. Such transfer requires full inclusions in both directions, not merely an operation-table resemblance or a matching short-polynomial partition. Computational domain checks and derivation replays are independent implementation checks on this unrestricted proof; they are not its replacement.

\subsection{右作用半格与尾部覆盖：两组完整有限基}\label{case:6}
\par\smallskip\noindent\textit{1. 右作用表示.}\quad 
对加法幂等半环 \(A\)，令 \(R_y(x)=xy\)。假设

\CaseDisplay{xyz\approx xzy\quad(N),\qquad xy^2\approx xy\quad(R).}

则各 \(R_y\) 是加法半格的自同态，彼此交换且幂等；它们的复合也是同一族中的映射，因为 \(R_yR_z=R_{yz}\)。点态并满足 \(R_{y+z}=R_y\vee R_z\)。于是右平移构成一个乘法为交换幂等复合的有限半环。该表示不要求 \(A\) 的每个元素乘法幂等，也不要求 \(xy=x^2y^2\)。

对词 \(w=h\,v\)，记 \(c(w)\) 为内容，\(D(w)=c(v)\) 为尾部变量集；一次词的尾部为空。由 \((N),(R)\) 可任意交换尾部并删除其重复变量。因此词值由首字母与尾部变量集决定，但必须保留首字母是否再次出现在尾部的信息。由于 \(x^3=x^2\)，首字母重复的非一次正规词为幂等元。

对非空多项式 \(F\) 定义

\CaseDisplay{T=c(F),\quad D(F)=\bigcup_{u\in F}D(u),\quad H(F)=\{h(u):u\in F\},}

\CaseDisplay{D_h(F)=\bigcup_{u\in F,\,h(u)=h}D(u).}

这里一次项的首字母就是自身。尾部集合不是全部非一次项的内容，二者不得混用。

\par\smallskip\noindent\textit{2. 两作用构造 \(A\)（目录 928）.}\quad 
底集为 \(\{z, e, b, d, t\}\)，加法是偏序 \(e<b<z<t\)、\(e<b<d<t\) 上的并，其中 \(z, d\) 不可比。乘法表为

\noindent \(\cdot\)\enspace\textbar\enspace \(z\)\enspace\textbar\enspace \(e\)\enspace\textbar\enspace \(b\)\enspace\textbar\enspace \(d\)\enspace\textbar\enspace \(t\)\par
\noindent \(z\)\enspace\textbar\enspace \(z\)\enspace\textbar\enspace \(z\)\enspace\textbar\enspace \(z\)\enspace\textbar\enspace \(z\)\enspace\textbar\enspace \(z\)\par
\noindent \(e\)\enspace\textbar\enspace \(z\)\enspace\textbar\enspace \(e\)\enspace\textbar\enspace \(z\)\enspace\textbar\enspace \(z\)\enspace\textbar\enspace \(z\)\par
\noindent \(b\)\enspace\textbar\enspace \(z\)\enspace\textbar\enspace \(b\)\enspace\textbar\enspace \(z\)\enspace\textbar\enspace \(z\)\enspace\textbar\enspace \(z\)\par
\noindent \(d\)\enspace\textbar\enspace \(z\)\enspace\textbar\enspace \(b\)\enspace\textbar\enspace \(z\)\enspace\textbar\enspace \(t\)\enspace\textbar\enspace \(t\)\par
\noindent \(t\)\enspace\textbar\enspace \(z\)\enspace\textbar\enspace \(z\)\enspace\textbar\enspace \(z\)\enspace\textbar\enspace \(t\)\enspace\textbar\enspace \(t\)\par
符号顺序 \((z, e, b, d, t)\) 对应目录标号 \((0,\allowbreak{}2,\allowbreak{}1,\allowbreak{}3,\allowbreak{}4)\)。右平移有三种：常值 \(z\) 映射 \(Z\)，\(E=R_e\)，以及 \(G=R_d=R_t\)。有 \(E^2=E\)、\(G^2=G\)、\(EG=GE=Z\)，且 \(R_{xy}=R_xR_y\)，故乘法结合。每个右平移保持并；每个左平移保持两个链的序及唯一非平凡不可比并 \(z+d=t\)，从而满足两侧分配律。

非一次词的值按首字母描述如下：首字母 \(z\) 总给 \(z\)；首字母 \(e\) 在全部尾字母为 \(e\) 时给 \(e\)，否则给 \(z\)；首字母 \(b\) 在全部尾字母为 \(e\) 时给 \(b\)，否则给 \(z\)；首字母 \(d\) 在尾字母全为 \(e\) 时给 \(b\)，全属 \(\{d, t\}\) 时给 \(t\)，否则给 \(z\)；首字母 \(t\) 在尾字母全属 \(\{d, t\}\) 时给 \(t\)，否则给 \(z\)。

此外 \(A\) 满足 \(x^2y=y^2x\)。故首字母重复的非一次词只由内容决定，记为 \(M_X\)：\(X\) 内全取 \(e\) 时值为 \(e\)，全取 \(d, t\) 时值为 \(t\)，其他情况值为 \(z\)。

\par\smallskip\noindent\textit{2.1 完整吸收判据.}\quad 
定理 A。 对非一次词 \(w\)，\(w\le F\) 当且仅当

\CaseDisplay{c(w)\subseteq T,\qquad D(w)\subseteq D(F),\qquad \exists q\in F: c(q)\subseteq c(w).\quad\textup{(\ensuremath{A})}}

最后的内部项 \(q\) 可以是一次项。一次词 \(x\le F\) 当且仅当 \(F\) 有同名一次项，或有仅含 \(x\) 的非一次项。

必要性。 若内容遗漏 \(x\)，赋 \(x=z\)、其余变量 \(e\)，左值为 \(z\)、右值为 \(e\)。若 \(x\in D(w)\setminus D(F)\)，赋 \(x=b\)、其余变量 \(e\)：左值为 \(z\)，右侧 \(x\) 只能作一次项或唯一首字母，所有右值至多为 \(b<z\)。若无内部项，赋 \(c(w)\) 内变量为 \(t\)、外部变量为 \(e\)。左值为 \(t\)，每个右项因含外部变量而至多为 \(z<t\)。一次词的反例赋值为 \(x=d\)、其他变量 \(e\)：没有同名一次项和纯 \(x\) 非一次项时，右值至多为 \(z\)，不能吸收 \(d\)。

充分性。 非一次左值可能为 \(e, b, z, t\)。\(e\) 为加法最小元。左值为 \(b\) 时，首字母为 \(b\) 或 \(d\)；内容条件给出含该变量的右项，其值至少为 \(b\)。若左值为 \(z\) 而所有右项均取 \(e, b, d\)，则每个非一次右项必须取 \(e\) 或 \(b\)，从而 \(D(F)\) 内变量全取 \(e\)，且 \(T\) 内变量全取 \(e, b, d\)。前两条件于是迫使左值为 \(e\) 或 \(b\)，矛盾；故某个右项取 \(z\) 或 \(t\)。左值为 \(t\) 时，\(c(w)\) 内变量全取 \(d, t\)。内部项非一次则取 \(t\)；内部一次项取 \(t\) 也完成。若它取 \(d\)，选任一尾变量 \(x\in D(w)\)，尾部覆盖给出含 \(x\) 于尾部的非一次右项，其值为 \(z\) 或 \(t\)，故与 \(d\) 的并为 \(t\)。一次项充分性用 \(x\le x^2=x^n\)（\(n\ge2\)）。

\par\smallskip\noindent\textit{2.2 六条完整基.}\quad 
在半环公理之外取

\CaseDisplay{\begin{aligned}
xyz&\approx xzy, &(N)\\
xy^2&\approx xy, &(R)\\
x^2y&\approx y^2x, &(C)\\
x&\le x^2, &(I)\\
xy&\le x+zy, &(V)\\
xy&\le y+xyz. &(E)
\end{aligned}}

有效性由词值规则及判据 (A) 直接验证。为证完整性，先记录尾部补全：若 \(p\le F\) 且 \(x\in D(F)\)，选一个右项 \(u\) 含 \(x\) 于尾部。利用 \((N),(R)\) 把它写成 \(u=vx\)，其中 \(v\) 非空；\((V)\) 给出 \(px\le p+vx=p+u\le F\)。因此可逐个补入任意需要的尾变量，不引入见证的其他变量。

对内容为 \(X\) 的平方型左词 \(M_X\)，判据给出 \(X\subseteq D(F)\)。从内部右项 \(q\) 出发，补入 \(X\) 的全部变量，特别是再补一次 \(q\) 的首字母，得到首字母重复且内容恰为 \(X\) 的词；\((C)\) 把它化为 \(M_X\)。

再考虑首字母 \(h\) 不在尾部的左词 \(w\)。若 \(h\in D(F)\)，则 \(X\subseteq D(F)\)，先得 \(M_X\le F\)，再用 \((I)\) 的词代入 \(w\le w^2=M_X\)。以下设 \(h\notin D(F)\)；于是右侧 \(h\) 的出现只能是一次项或非一次项的唯一首字母。若有包含 \(h\) 的内部项，从它出发用尾部补全直接得到 \(w\)。否则选一个避开 \(h\) 的内部项 \(q\)；若有一次项 \(h\) 也可直接补全，剩下的情况存在以 \(h\) 开头的非一次右项 \(r=hb\)，其中 \(b\) 非空。给 \(r\) 补入 \(D(w)\)，得到 \(wb\le F\)。因 \(c(q)\subseteq D(w)\)，\((N),(R)\) 给出 \(wq=w\)。在 \((E)\) 中取 \((x, y, z)=(w, q, b)\)，得到

\CaseDisplay{w=wq\le q+wqb=q+wb\le F.}

一次左词用 \((I)\)。任意有效恒等式拆成双向逐词吸收，故六式完整。整个证明允许任意变量数，无空词代入。

\par\smallskip\noindent\textit{3. 单闭包作用构造 \(B\)（目录 872）.}\quad 
取 \(B=\{z, b, a, e, t\}\)，加法序为 \(z<b<a<t\)，且 \(e<t\) 与 \(z, b, a\) 都不可比。定义保持并的幂等膨胀映射

\CaseDisplay{\theta(z)=z,\quad\theta(b)=a,\quad\theta(a)=a,\quad\theta(e)=t,\quad\theta(t)=t,}

并令 \(xy=x\)（\(y=e\) 时），\(xy=\theta(x)\)（\(y\ne e\) 时）。符号顺序 \((z, b, a, e, t)\) 对应目录标号 \((2,\allowbreak{}1,\allowbreak{}0,\allowbreak{}3,\allowbreak{}4)\)。右作用为恒等映射或 \(\theta\)；且 \(xy=e\) 当且仅当 \(x=y=e\)，故结合律成立。右分配由 \(\theta\) 保持并得到。左分配在 \(y=e\)、\(z\ne e\) 时由 \(x\le\theta(x)\) 给出，其余情形直接成立。

一个词以 \(z, a, t\) 开头时恒取该首字母；以 \(b\) 开头且尾部全为 \(e\) 时取 \(b\)，否则取 \(a\)；以 \(e\) 开头且尾部全为 \(e\) 时取 \(e\)，否则取 \(t\)。

\par\smallskip\noindent\textit{3.1 完整吸收判据.}\quad 
定理 B。 对非一次词 \(w\)，置 \(h=h(w)\)、\(X=c(w)\)。\(w\le F\) 当且仅当

\CaseDisplay{h\in H(F),\qquad X\subseteq T,\qquad
\bigl[D(w)\subseteq D_h(F)\ \text{或} D_h(F)\cap H(F)\ne\varnothing\bigr].\quad\textup{(\ensuremath{B})}}

一次词的准确判据为 \(x\in H(F)\)，即存在以 \(x\) 开头的右项。这与构造 \(A\) 的一次项判据不同，因为本构造满足 \(x\le xy\)。

必要性。 若无同首字母右项，赋 \(h=a\)、其他变量 \(z\)，左值 \(a\)、右值 \(z\)。在已有同首字母项的前提下，遗漏的内容变量 \(x\) 不等于 \(h\)；赋 \(x=z\)、其他变量 \(e\)，左值 \(t\)、右值 \(e\)。若最后一个条件失败，赋 \(h=b\)，\(D_h(F)\) 内变量 \(e\)，其他变量 \(z\)。此赋值一致，因为 \(D_h(F)\cap H(F)\) 为空而 \(h\in H(F)\)。左词有尾变量不在 \(D_h(F)\)，其值非 \(e\)，故左值为 \(a\)；同首字母右项取 \(b\)，其余右项首字母或一次项取 \(z\)，右值至多为 \(b<a\)。一次词的必要性同样用 \(x=a\)、其他变量 \(z\) 分离。

非一次充分性。 同首字母右项直接保护首字母取 \(z, a, t\) 的情形。首字母取 \(e\) 时，同首字母项在 \(e\) 之上；若左值为 \(t\) 但右和未达 \(t\)，所有右项都必须取 \(e\)，进而 \(T\) 内所有变量都取 \(e\)，与左尾存在非 \(e\) 变量矛盾。首字母取 \(b\) 时，同首字母项至少取 \(b\)；左值为 \(a\) 时，若尾部包含于 \(D_h(F)\)，某个同首字母项的尾字母非 \(e\)，该项取 \(a\)。否则选 \(g\in D_h(F)\cap H(F)\)：\(g\) 非 \(e\) 时同样给 \(a\)；\(g=e\) 时同首字母 \(g\) 的右项给 \(e\) 或 \(t\)，再与受保护的 \(b\) 求并得到 \(t\ge a\)。

一次词的充分性由同首字母项及 \(x\le xy\) 得到。

\par\smallskip\noindent\textit{3.2 六条完整基.}\quad 
在半环公理之外取

\CaseDisplay{\begin{aligned}
xyz&\approx xzy, &(N)\\
xy^2&\approx xy, &(R)\\
x&\le xy, &(U)\\
xyz&\le xy+xz, &(M)\\
x^2&\le xy+y, &(P)\\
x^2y&\le x^2+y. &(K)
\end{aligned}}

有效性由判据 (B) 和词值规则逐式验证。证明完整性。一次吸收由 \((U)\) 直接处理；以下取非一次左词 \(w\)。

若 \(D_h(F)\cap H(F)=\varnothing\)，判据给出 \(D(w)\subseteq D_h(F)\)。因 \(D(w)\) 非空，必有同首字母的非一次右项，将它们记为 \(hv_1,\ldots, hv_m\)。迭代 \((M)\)，得到

\CaseDisplay{hv_1\cdots v_m\le hv_1+\cdots+hv_m\le F.}

由于 \(h\in H(F)\) 而 \(D_h(F)\) 不交 \(H(F)\)，这些词的尾部都不含 \(h\)。上述积与 \(w\) 同首字母，且尾部为 \(D_h(F)\)。由 \((U)\) 允许在右端添字母，再用 \((N),(R)\) 整理，可从尾部包含关系推出 \(w\le hv_1\cdots v_m\)。

若 \(D_h(F)\cap H(F)\ne\varnothing\)，选其中一个 \(g\)。存在右项 \(u\) 以 \(h\) 开头且在尾部含 \(g\)，也存在右项 \(r\) 以 \(g\) 开头。由 \((N),(R),(U)\)，\(hg\le u\) 且 \(g\le r\)。将 \((P)\) 中 \((x, y)\) 取为 \((h, g)\)，得到

\CaseDisplay{h^2\le hg+g\le u+r\le F.}

从幂等非一次词 \(p=h^2\) 出发，对每个右项 \(v\)，将 \((K)\) 中 \((x, y)\) 取为 \((p, v)\)，可得 \(pv\le p+v\le F\)；乘积仍首字母重复，因而仍幂等。逐项相乘得到内容为 \(T\)、首字母重复的词 \(Q\le F\)。由 \((U)\) 有 \(w\le w^2\)；又因 \(c(w)\subseteq T\)，可继续用 \((U)\) 补入缺少变量，配合 \((N),(R)\) 得 \(w^2\le Q\)。故 \(w\le F\)。

两种分支覆盖全部有效吸收。拆分任意恒等式为双向逐词吸收，即得六式的任意变量数完备性。缺少变量集为空时直接使用正规形等式，不代入空词。两组基逐词反转便得到乘法对偶的完整基。对目录对象，只有精确双向簇包含经独立核验后才传递有限基性；不对任意真子簇作这种传递。

\par\smallskip\noindent\textit{4. 主理想诊断方法与范围.}\quad 
给定有限半环 \(S\) 和固定变量数下的完整词代表族 \(W\)，对赋值 \(\alpha\) 及阈值 \(s\in S\) 定义

\CaseDisplay{C_{\alpha, s}=\{u\in W:u(\alpha)\nleq s\}.}

固定左词 \(w\) 后，右项集合 \(P\subseteq W\) 满足 \(w\le\sum P\)，当且仅当对于全部 \(w(\alpha)\nleq s\) 的 \((\alpha, s)\)，都有 \(P\cap C_{\alpha, s}\ne\varnothing\)。必要性直接来自各右项均不超过 \(s\) 时其并也不超过 \(s\)；充分性在吸收失败时取 \(s=(\sum P)(\alpha)\)。

所以可把真实吸收语义转化为有限正 CNF，并以删去包含其他子句的冗余子句来比较判据。这在固定变量数下覆盖任意右项数量，避免显式枚举所有多项式。有限诊断仍不证明任意变量数的完备性；后者必须使用独立的形式归纳。

\subsection{非一次见证与无保护尾部：两组完整有限基}\label{case:7}
无常量签名 \((+,\cdot)\)，\(u\le v\) 表示 \(u+v\approx v\)。有限基相对于加法幂等半环公理。本文给出两个任意变量数的完备性证明，来源分别为目录 690 与 865；有限诊断不能替代这些证明。

\par\smallskip\noindent\textit{1. 共同记号.}\quad 
两来源均满足 \(xyz=xzy\) 和 \(xy^2=xy\)。记词的首字母为 \(h(w)\)，尾部变量集为 \(D(w)\)，内容为 \(c(w)\)。词由首字母及尾部集合正规化；不能删除首字母是否重复的信息。对非空右多项式 \(F\) 记

\CaseDisplay{T=c(F),\quad H=\{h(u):u\in F\},\quad E=\bigcup_{u\in F}D(u),\quad N=\bigcup_{u\in F, |u|\ge2}c(u).}

内部右项指内容包含于左词内容的右项，可以为一次项。

\par\smallskip\noindent\textit{2. 链上的非一次见证（来源 690）.}\quad 
底集为 \(\{z, a, b, e, t\}\)，加法链为 \(e<b<a<z<t\)，符号顺序对应目录 \((0,\allowbreak{}1,\allowbreak{}2,\allowbreak{}3,\allowbreak{}4)\)。右乘 \(e\) 为映射 \(E_0=(z, a, a, e, z)\)，右乘 \(t\) 为 \(G_0=(z, z, z, z, t)\)，右乘 \(z, a, b\) 为常值 \(z\) 映射 \(Z\)。这些映射保持链上的并，满足 \(E_0^2=E_0\)、\(G_0^2=G_0\)、\(E_0G_0=G_0E_0=Z\)；由 \(R_{xy}=R_xR_y\) 和 \(R_{x+y}=R_x\vee R_y\) 得到结合律与分配律。

非一次词全取 \(e\) 时值为 \(e\)；首字母取 \(a\) 或 \(b\) 且尾部全为 \(e\) 时值为 \(a\)；所有变量取 \(t\) 时值为 \(t\)；其他情形值为 \(z\)。它也满足 \(x^2y=y^2x\)，因此首字母重复的词仅由内容决定，记为 \(M_X\)。

\par\smallskip\noindent\textit{2.1 吸收判据.}\quad 
定理 A。 对非一次词 \(w\)，

\CaseDisplay{w\le F\iff c(w)\subseteq N,\quad D(w)\subseteq E,\quad\exists q\in F: c(q)\subseteq c(w).\quad\textup{(\ensuremath{A})}}

一次词 \(x\le F\) 当且仅当右边有纯 \(x\) 项。

必要性。 若 \(x\in c(w)\setminus N\)，赋 \(x=b\)，其余变量 \(e\)。右非一次项均为 \(e\)，右一次项至多为 \(b\)；左值为 \(a\) 或 \(z\)，均大于 \(b\)。若 \(x\in D(w)\setminus E\)，仍取此赋值，右非一次项至多为 \(a\)，左值为 \(z>a\)。若没有内部项，赋左内容内变量 \(t\)、其余变量 \(e\)，则左值为 \(t\)，所有右项至多为 \(z\)。同一赋值证明一次项的必要性。

充分性。 分别检查四个真主理想。右和不超过 \(e\) 时，\(T\) 内变量全为 \(e\)，左值也为 \(e\)。右和不超过 \(b\) 时，所有非一次右项都取 \(e\)，故 \(N\) 内变量全为 \(e\)；内容条件保证左值为 \(e\)。右和不超过 \(a\) 时，全部右尾变量为 \(e\)，全部右首字母不超过 \(a\)；由 \(c(w)\subseteq N\subseteq T\) 与尾部条件，左值为 \(e\) 或 \(a\)。右和不超过 \(z\) 时，每个右项都含一个非 \(t\) 变量；内部项迫使左内容也含非 \(t\) 变量，故左值不为 \(t\)。一次项充分性由 \(x\le x^2=x^n\)（\(n\ge2\)）得到。证毕。

\par\smallskip\noindent\textit{2.2 八条完整基.}\quad 
定理 B。 来源 690 的完整附加基为

\CaseDisplay{\begin{aligned}
xyz&\approx xzy, &(N_0)\\
xy^2&\approx xy, &(R_0)\\
x^2y&\approx y^2x, &(C_0)\\
x&\le x^2, &(I_0)\\
xyz&\le xy+tz, &(T_0)\\
xy&\le y+xyz, &(B_0)\\
xy&\le x+xyz, &(A_0)\\
x^2&\le x+yx. &(S_0)
\end{aligned}}

有效性由 (A) 逐式检验。与允许一次词直接补入尾部的规则不同，\((T_0)\) 只保证给已有非一次词补尾：若 \(p=ab\le F\)，\(a, b\) 均非空，且右项 \(u=vx\) 见证 \(x\in E\)，则

\CaseDisplay{px=abx\le ab+vx=p+u\le F.\quad\textup{(\ensuremath{3})}}

先处理平方型左词 \(M_X\)。此时 \(D(M_X)=X\subseteq E\)。若内部项 \(q\) 非一次，从它用 (3) 补入全部 \(X\)，所得首字母重复且内容恰为 \(X\)，由 \((C_0)\) 化成 \(M_X\)。若内部项是一次词 \(x\)，选尾部含 \(x\) 的右项 \(u=vx\)；\((S_0)\) 给出 \(x^2\le x+vx\le F\)，随后用 (3) 补全 \(X\)。这里不能从 \(x\le F\) 与 \(x\le x^2\) 错误反推 \(x^2\le F\)；尾部见证是必要的。

设左词 \(w\) 的首字母 \(h\) 只出现一次。若 \(h\in E\)，则 \(X=c(w)\subseteq E\)，先得 \(M_X=w^2\le F\)，再用 \((I_0)\)。以下 \(h\notin E\)。由 \(h\in N\)，存在以 \(h\) 为唯一首字母的非一次右项 \(r\)。给 \(r\) 补入 \(D(w)\)，得到 \(w\) 或 \(wb\le F\)，其中 \(b\) 是额外变量组成的非空词。无额外变量时已完成。

若内部项 \(q\) 含 \(h\) 且非一次，它也以 \(h\) 开头，从它直接补尾便得 \(w\)。若 \(q=h\)，将 \(w=hv\) 代入 \((A_0)\)，得到 \(w\le h+wb\le F\)；这里 \(v\) 非空。若 \(q\) 避开 \(h\)，则 \(c(q)\subseteq D(w)\)，所以 \(wq=w\)；\((B_0)\) 中取 \((x, y, z)=(w, q, b)\) 得 \(w\le q+wb\le F\)。

一次词用纯幂右项及 \((I_0)\)。两边逐词吸收给出任意恒等式的推导。整个证明不限制变量数，也不使用空词代入。证毕。

\par\smallskip\noindent\textit{3. 菱形加法与无保护尾部（来源 865）.}\quad 
底集仍为 \(\{z, a, b, e, t\}\)，目录标号仍为 \((0,\allowbreak{}1,\allowbreak{}2,\allowbreak{}3,\allowbreak{}4)\)。加法序为 \(e<a<z<t\)、\(e<b<z<t\)，\(a, b\) 不可比，且 \(a+b=z\)。乘法由以下右作用表定义：

\noindent 首元素\enspace\textbar\enspace 右乘 \(e\)\enspace\textbar\enspace 右乘 \(t\)\enspace\textbar\enspace 右乘 \(z, a, b\)\par
\noindent \(z\)\enspace\textbar\enspace \(z\)\enspace\textbar\enspace \(z\)\enspace\textbar\enspace \(z\)\par
\noindent \(a\)\enspace\textbar\enspace \(a\)\enspace\textbar\enspace \(z\)\enspace\textbar\enspace \(z\)\par
\noindent \(b\)\enspace\textbar\enspace \(b\)\enspace\textbar\enspace \(b\)\enspace\textbar\enspace \(b\)\par
\noindent \(e\)\enspace\textbar\enspace \(e\)\enspace\textbar\enspace \(b\)\enspace\textbar\enspace \(b\)\par
\noindent \(t\)\enspace\textbar\enspace \(z\)\enspace\textbar\enspace \(t\)\enspace\textbar\enspace \(z\)\par
三个右作用均保持并、幂等且彼此交换；由表逐类核对 \(R_{xy}=R_xR_y\)、\(R_{x+y}=R_x\vee R_y\) 得半环公理。

非一次词以 \(z, b\) 开头时分别恒取 \(z, b\)；以 \(a\) 开头时，尾部全为 \(e\) 则取 \(a\)，否则取 \(z\)；以 \(e\) 开头时，尾部全为 \(e\) 则取 \(e\)，否则取 \(b\)；以 \(t\) 开头时，尾部全为 \(t\) 则取 \(t\)，否则取 \(z\)。

\par\smallskip\noindent\textit{3.1 吸收判据.}\quad 
定理 C。 任意词 \(w\) 的准确吸收判据为

\CaseDisplay{w\le F\iff h(w)\in H,\quad D(w)\subseteq E,\quad\exists q\in F: c(q)\subseteq c(w).\quad\textup{(\ensuremath{C})}}

证明。 真主理想为 \(\{e\}\)、\(\{e, a\}\)、\(\{e, b\}\)、\(\{e, a, b, z\}\)。对非一次词，落入 \(\{e, a\}\) 当且仅当首字母取 \(e, a\) 且尾部全取 \(e\)；落入 \(\{e, b\}\) 当且仅当首字母取 \(e, b\)，尾部不受限制；取值不为 \(t\) 当且仅当内容中有非 \(t\) 变量。由首字母、尾部及内部项三个条件，逐个理想得到充分性；内容包含于 \(T\) 已被前两个条件蕴含。一次项也满足同一结论。

若缺少同首字母项，赋左首字母 \(z\)、其余变量 \(b\)，左值为 \(z\)、右值为 \(b\)。若左尾变量 \(x\) 不在 \(E\)，赋 \(x=a\)、其余变量 \(e\)，右和至多为 \(a\)，左值为 \(b\) 或 \(z\)，不能被 \(a\) 吸收。若没有内部项，赋左内容内变量 \(t\)、外部变量 \(e\)，左值为 \(t\)、右值至多为 \(z\)。这些赋值给出必要性。证毕。

\par\smallskip\noindent\textit{3.2 五条完整基.}\quad 
定理 D。 来源 865 的完整附加基只有五条：

\CaseDisplay{xyz\approx xzy,\quad xy^2\approx xy,\quad x\le x^2,\quad xy\le x+zy,\quad xy\le y+xyz.\quad\textup{(\ensuremath{D})}}

由 (C) 知各式有效。设 \(p\le F\)、\(x\in E\)，把尾部见证写成 \(vx\)，第四式给 \(px\le p+vx\le F\)；此处 \(p\) 允许是一次词。

当左词 \(w\) 首字母重复时，从同首字母右项补入 \(c(w)\)，得到 \(w\) 或 \(wb\le F\)。内部项 \(q\) 满足 \(wq=w\)；若 \(b\) 非空，最后一式给 \(w\le q+wb\le F\)。当首字母 \(h\) 不重复而属于 \(E\) 时，对 \(w^2\) 用前段，再用 \(w\le w^2\)。最后设 \(h\notin E\)。内部项若含 \(h\)，必以 \(h\) 为唯一首字母，从它直接补尾；内部项若避开 \(h\)，先从同首字母右项补尾，再用 \(wq=w\) 与最后一式删除额外变量。一次左词用内部纯幂项及 \(x\le x^2\)。

推导长度可随变量数增长，但只使用上述五条固定恒等式。全部有效吸收均可推出，故任意变量数完备。证毕。

\par\smallskip\noindent\textit{4. 范围与核验.}\quad 
两个来源及其乘法对偶分别使用相应基，不能混同``一次项见证''和``非一次项见证''。独立核验重建运算、检查基、对固定一至五个变量用主理想 CNF 比较全部右多项式，并对所列对象完成双向 HSP 包含。有限基性只传递给已确认同簇的对象；本文没有从单向簇包含推出有限基，也未解决五元素极限簇的一般存在问题。

\subsection{非一次首字母与局部尾部：两组完整基}\label{case:8}
\par\smallskip\noindent\textit{1. 共同记号.}\quad 
两个来源满足 \(xyz=xzy\)、\(xy^2=xy\)、\(x^2y=y^2x\)。因此词由首字母与尾部集合正规化，首字母重复的词只由内容决定，记为 \(M_X\)。有 \(w^2=M_{c(w)}\)，且这些平方词幂等。

对非空右多项式 \(F\)，令 \(T=c(F)\)、\(H=\{h(u):u\in F\}\)、\(E=\bigcup_{u\in F}D(u)\)。另记

\CaseDisplay{H_2=\{h(u):u\in F,\ |u|\ge2\},\qquad D_h=\bigcup_{u\in F,\ h(u)=h}D(u).}

空尾部只用于集合记号；证明不把空词代入变量。

\par\smallskip\noindent\textit{2. 必须使用非一次首字母见证（来源 228）.}\quad 
底集为 \(\{z, a, b, e, t\}\)，对应目录标号 \((0,\allowbreak{}1,\allowbreak{}2,\allowbreak{}3,\allowbreak{}4)\)。加法序为 \(b<a<z<t\)、\(e<t\)，\(e\) 与 \(b, a, z\) 不可比。右乘 \(e\) 是 \(P=(z, z, b, e, t)\)；右乘 \(t\) 是 \(Q=(z, z, z, t, t)\)；右乘 \(z, a, b\) 是常值 \(z\) 映射 \(Z\)。这些映射保持并、幂等且交换，满足 \(PQ=QP=Q\)、\(P\vee Z=Q\)。由 \(R_{xy}=R_xR_y\) 和 \(R_{x+y}=R_x\vee R_y\) 得到半环公理。

非一次词若尾部出现 \(z, a, b\)，或首字母为 \(z, a\)，则取 \(z\)。首字母为 \(b\) 时，尾部全为 \(e\) 才取 \(b\)，否则取 \(z\)。首字母为 \(e, t\) 且尾部全为 \(e, t\) 时，词值是该二元素链上的最大值。特别地，非一次词从不取 \(a\)。

\par\smallskip\noindent\textit{2.1 完整吸收判据.}\quad 
定理 A。 对非一次词 \(w\)，令 \(h=h(w)\)、\(X=c(w)\)，则

\CaseDisplay{w\le F\iff X\subseteq T,\quad \exists q\in F:\ c(q)\subseteq X,\quad
\bigl[(h\in H_2\ \text{且}\ D(w)\subseteq E)\ \text{或}\ H\cap E\ne\varnothing\bigr].\quad\textup{(\ensuremath{A})}}

一次词 \(x\le F\) 当且仅当右边有仅含 \(x\) 的项。

证明。 真主理想为 \(\{b\}\)、\(\{b, a\}\)、\(\{b, a, z\}\)、\(\{e\}\)。右和不超过 \(a\)，要求非一次右首字母取 \(b\)、所有右尾变量取 \(e\)、一次右项取 \(a, b\)。有保护交集 \(H\cap E\) 时不可能；否则 \(h\in H_2\) 与尾部覆盖保证左非一次词取 \(b\)。这也处理右和不超过 \(b\) 的情形。右和为 \(e\) 迫使 \(T\) 内变量全为 \(e\)，内容条件使左值为 \(e\)。

一个词不超过 \(z\) 当且仅当其内容含取 \(z, a, b\) 的变量；内部项保证右和不超过 \(z\) 时左值也不超过 \(z\)。充分性得证。

内容遗漏变量 \(x\) 时，赋 \(x=t\)、其余变量 \(e\)，左值为 \(t\)、右值为 \(e\)。没有内部项时，赋左内容内变量 \(e\)、外部变量 \(z\)，左值为 \(e\)、右值为 \(z\)，不可比较。

若最后一项失败，\(H\cap E=\varnothing\)。赋 \(H_2\) 内变量 \(b\)、\(H\setminus H_2\) 内变量 \(a\)、\(E\) 内变量 \(e\)、其余变量 \(z\)。全部右项不超过 \(a\)。左首字母不在 \(H_2\) 时，左非一次值为 \(e\) 或 \(z\)；左首字母在 \(H_2\) 但尾部有未覆盖变量时，左值为 \(z\)。均不能被 \(a\) 吸收。一次词的必要性由无内部纯幂项时的同一内容赋值得到，充分性由 \(x\le x^2=x^n\)（\(n\ge2\)）或相同一次项得到。证毕。

\par\smallskip\noindent\textit{2.2 九条完整基.}\quad 
定理 B。 来源 228 的完整附加基为

\CaseDisplay{\begin{aligned}
xyz&\approx xzy, &(N)\\
xy^2&\approx xy, &(R)\\
x^2y&\approx y^2x, &(C)\\
x&\le x^2, &(I)\\
xyz&\le xy+tz, &(T_0)\\
xy&\le y+xyz, &(B)\\
xy&\le x+xyz, &(A_0)\\
x^2&\le x+yx, &(S_0)\\
x^2y&\le x^2+y. &(K)
\end{aligned}}

有效性由定理 A 验证。若 \(p=uv\le F\) 为非一次词，且右项 \(r=sx\) 见证 \(x\in E\)，\((T_0)\) 给出

\CaseDisplay{px=uvx\le uv+sx\le F.\quad\textup{(\ensuremath{3})}}

该补尾规则不能直接用于一次词。

无保护时，定理 A 给出 \(h\in H_2\)、\(D(w)\subseteq E\)，且 \(h\notin E\)。选择以 \(h\) 开头的非一次右项 \(r\)，用 (3) 补入全部 \(D(w)\)，得到 \(w\) 或 \(wb\le F\)，其中 \(b\) 是额外变量组成的非空词。

内部项 \(q\) 若含 \(h\) 且非一次，则它也以 \(h\) 为唯一首字母，从它直接补尾即得 \(w\)。若 \(q=h\)，写 \(w=hv\)，用 \((A_0)\) 得 \(w\le h+wb\le F\)。若 \(q\) 不含 \(h\)，则 \(wq=w\)，\((B)\) 在 \((x, y, z)=(w, q, b)\) 下得 \(w\le q+wb\le F\)。没有额外变量时不需要这一步。

有保护时，取 \(g\in H\cap E\)，选尾部见证 \(u=vg\) 及以 \(g\) 开头的右项 \(r\)。若 \(r\) 非一次，用 (3) 得 \(rg\le F\)，这是首字母重复的幂等词。若 \(r=g\)，用 \((S_0)\) 得 \(g^2\le g+vg\le F\)。因而总能得到一个幂等平方词 \(p\le F\)。

由 \((K)\)，对每个右项 \(v\) 有 \(pv\le p^2+v=p+v\le F\)。逐项补入全部右内容，再由 \((C)\) 得 \(M_T\le F\)。内部项 \(q\) 的内容包含于 \(X\)；若 \(X=T\)，已经得到 \(M_X\le F\)。否则取非空词 \(b\) 的内容为 \(T\setminus X\)，由 \(M_Xq=M_X\) 与 \((B)\) 得

\CaseDisplay{M_X\le q+M_Xb=q+M_T\le F.}

最后 \(w\le w^2=M_X\le F\)。一次词直接用内部纯幂项与 \((I)\)。任意变量数的完整性由此成立。

\par\smallskip\noindent\textit{3. 同首字母的局部尾部（来源 5767）.}\quad 
底集为 \(\{z, e, b, d, t\}\)，对应目录标号 \((0,\allowbreak{}1,\allowbreak{}2,\allowbreak{}3,\allowbreak{}4)\)。加法序为 \(e<z<t\)、\(d<b<t\)，两链除 \(t\) 外不相交。右乘 \(e\) 为恒等映射，右乘 \(z\) 为 \(P=(z, z, t, d, t)\)，右乘 \(b, d, t\) 为常值 \(t\) 映射。这些映射保持并、幂等且交换，满足 \(R_{xy}=R_xR_y\)、\(R_{x+y}=R_x\vee R_y\)，故为半环。

词尾出现 \(b, d, t\) 时，非一次词取 \(t\)。尾部全为 \(e, z\) 时：首字母为 \(d\) 则取 \(d\)；首字母为 \(b\) 且尾部全为 \(e\) 时取 \(b\)，否则取 \(t\)；首字母为 \(e, z\) 时取内容在 \(e<z\) 上的最大值；首字母为 \(t\) 时取 \(t\)。此半环满足更强的 \(x\le xy\)。

\par\smallskip\noindent\textit{3.1 完整吸收判据.}\quad 
定理 C。 任意词 \(w\) 满足

\CaseDisplay{w\le F\iff c(w)\subseteq T,\quad
\bigl[(h(w)\in H\ \text{且}\ D(w)\subseteq D_{h(w)})\ \text{或}\ H\cap E\ne\varnothing\bigr].\quad\textup{(\ensuremath{C})}}

这里不要求内部项；无保护时必须使用同首字母的局部尾部，而不能只检查全局尾部。

充分性。 真主理想是 \(\{e\}\)、\(\{e, z\}\)、\(\{d\}\)、\(\{d, b\}\)。右和为 \(e\) 或落入 \(\{e, z\}\) 时，全部 \(T\) 内变量分别只取 \(e\) 或取 \(e, z\)，由内容条件得到左值在同一理想内。

右和落入 \(\{d, b\}\)，要求右首字母全取 \(d, b\)、右尾变量全取 \(e, z\)；有保护交集时这种赋值不可能。无保护时，左首字母 \(h\) 属于 \(H\)。若它取 \(d\)，尾部覆盖使左值为 \(d\)；若它取 \(b\)，所有以 \(h\) 开头的右项只有在尾部全取 \(e\) 时才能不超过 \(b\)，故 \(D_h\) 内变量全为 \(e\)，左值为 \(b\)。右和为 \(d\) 时，所有右首字母取 \(d\)，同样得左值为 \(d\)。这也覆盖一次左词。

必要性。 若左内容变量 \(x\) 不在 \(T\)，赋 \(x=z\)、其余变量 \(e\)，左值为 \(z\)、右值为 \(e\)。若无保护而左首字母不在 \(H\)，赋 \(H\) 内变量 \(d\)、\(E\) 内变量 \(e\)、其余变量 \(t\)；右和为 \(d\)，左值不能为 \(d\)。

现设无保护且 \(h\in H\)，但 \(D(w)\nsubseteq D_h\)。赋 \(h=b\)，其余 \(H\) 内变量 \(d\)，\(D_h\) 内变量 \(e\)，\(E\setminus D_h\) 内变量 \(z\)，其余变量 \(t\)。因 \(H\cap E=\varnothing\)，赋值一致。以 \(h\) 开头的右项取 \(b\)，其余右项取 \(d\)，右和不超过 \(b\)；左尾含非 \(e\) 变量，左值为 \(t\)。必要性得证。

\par\smallskip\noindent\textit{3.2 八条完整基.}\quad 
定理 D。 来源 5767 的完整附加基为

\CaseDisplay{\begin{aligned}
xyz&\approx xzy, &(N)\\
xy^2&\approx xy, &(R)\\
x^2y&\approx y^2x, &(C)\\
x&\le xy, &(U)\\
xyz&\le xy+xz, &(M)\\
x^2y&\le x^2+y, &(K)\\
(xyz)^2&\le xy+yz, &(P_3)\\
(xy)^2&\le xy+y. &(P_2)
\end{aligned}}

有效性由定理 C 验证。无保护时有 \(h\in H\)、\(D(w)\subseteq D_h\)。若 \(D_h\) 为空，左词只能是一次词 \(h\)，且右边实际有 \(h\)。否则迭代 \((M)\) 合并全部以 \(h\) 开头的非一次右项，得到首字母为 \(h\)、尾部为 \(D_h\) 的词 \(p\le F\)。用 \((U)\) 给 \(w\) 补入 \(D_h\setminus D(w)\)，再以 \((N),(R)\) 整理，得 \(w\le p\le F\)。

有保护时，取 \(g\in H\cap E\)，右项 \(u=vg\) 和以 \(g\) 开头的右项 \(r\)。若 \(r=gs\) 非一次，\((P_3)\) 给出 \((vgs)^2\le u+r\le F\)；若 \(r=g\)，\((P_2)\) 给出 \((vg)^2\le u+g\le F\)。得到幂等词 \(p\le F\) 后，用 \((K)\) 逐项补入全部右内容，再用 \((C)\) 得到 \(M_T\le F\)。

由 \((U)\) 有 \(w\le w^2\)；继续补入 \(T\setminus c(w)\)，便有 \(w\le M_T\le F\)。这一推导也适用于一次左词，不需要内部项。全过程使用有限归纳且不限制变量数；恒等式两向逐词吸收均可推出，八式完整。证毕。

\par\smallskip\noindent\textit{4. 分支闭合的限定.}\quad 
两个来源满足来源 4397 的九条基，但这里的有限基结论来自各自完整证明和双向簇包含，不能仅从这种单向成员关系推出。是否覆盖九条基的全部五元素目录成员，另须逐表检查；即使该有限目录分支全部有限基，也不意味着整个九条基簇遗传有限基。独立的固定变量诊断同样不替代上述任意变量数证明。

\subsection{非一次内容覆盖与全局保护：来源825的十二条完整基}\label{case:9}
\par\smallskip\noindent\textit{1. 构造与求值.}\quad 
取 \(A=\{z, a, d, e, t\}\)，依次对应目录825的标号 \((0,\allowbreak{}1,\allowbreak{}2,\allowbreak{}3,\allowbreak{}4)\)。加法为两条链 \(a<z<t\)、\(e<d<t\) 的并半格，除共同顶元外两条链彼此不可比。

右乘 \(z, a\) 为常值 \(z\) 映射，右乘 \(e\) 为 \(R=(z, a, t, e, t)\)，右乘 \(d, t\) 为 \(P=(z, z, t, t, t)\)。三个右平移保持并、交换、幂等，且 \(RP=PR=P\)；它们与右因子乘法和加法的兼容关系逐类可检验，因此定义半环。

一次词直接取变量值。非一次词取值为：

\textbullet{} 全部变量取 \(e\) 时为 \(e\)；

\textbullet{} 首字母取 \(a\)、尾部全部取 \(e\) 时为 \(a\)；

\textbullet{} 内容含取值 \(a, z\) 的变量且不属于上一情形时为 \(z\)；

\textbullet{} 内容全部取 \(e, d, t\) 且至少一个变量不取 \(e\) 时为 \(t\)。

由此可见，非一次词不会取值 \(d\)。

\par\smallskip\noindent\textit{2. 精确吸收判据.}\quad 
设 \(F\) 为非空右词集，\(X=c(w)\)、\(h=h(w)\)、\(D=D(w)\)，并定义

\CaseDisplay{F_X=\{u\in F:c(u)\subseteq X\},\quad H=h(F),\quad
E=\bigcup_{u\in F}D(u),\quad T_2=\bigcup_{\substack{u\in F\\|u|\ge2}}c(u).}

此处 \(T_2\) 只记录出现在非一次右项中的变量，不是全部右侧内容。

定理1。 一次词 \(x\le F\) 在 \(A\) 上有效，当且仅当 \(F_{\{x\}}\ne\varnothing\)。对于非一次词 \(w\)，有

\CaseDisplay{A\models w\le F\quad\Longleftrightarrow\quad
F_X\ne\varnothing,\quad X\subseteq T_2,\quad
\bigl[H\cap E\ne\varnothing\ \text{或}\ (h\in H\ \text{且}\ D\subseteq E)\bigr].\quad\textup{(\ensuremath{1})}}

必要性。 若没有内部右项，赋 \(X\) 内变量 \(e\)、外部变量 \(z\)；右词全部取 \(z\)，左词取 \(e\)，不能被吸收。这也证明一次情形的必要性；一次情形的充分性由纯幂正规化和 \(x\le x^2\) 得到。

若非一次左词含变量 \(x\notin T_2\)，赋 \(x=d\)、其余变量 \(e\)。非一次右词全部取 \(e\)，一次右词取 \(d\) 或 \(e\)，右和不超过 \(d\)；左值为 \(t\)。

最后假设 \(H\cap E=\varnothing\)，赋右首字母变量 \(a\)、右尾部变量 \(e\)。由已经成立的 \(X\subseteq T_2\)，左变量全部被赋值。右词均取 \(a\)。若左首字母不在 \(H\)，左值为 \(e\) 或 \(z\)；若左首字母在 \(H\) 而某左尾变量不在 \(E\)，该尾变量取 \(a\)，左值为 \(z\)。两者均不能不超过 \(a\)。

充分性。 真主理想分别为 \(\{a\}\)、\(\{a, z\}\)、\(\{e\}\)、\(\{e, d\}\)。右和不超过 \(z\) 时，每个右词都含取值 \(a, z\) 的变量；内部右项保证左内容也含这种变量，因此非一次左值为 \(a\) 或 \(z\)。

右和不超过 \(d\) 时，每个非一次右词只能取 \(e\)，并且其全部变量取 \(e\)。由 \(X\subseteq T_2\)，左词全部变量取 \(e\)，左值为 \(e\)。这同时处理右和为 \(e\) 的情形。

右和为 \(a\) 时，所有右首字母取 \(a\)、所有右尾变量取 \(e\)，故 \(H\cap E=\varnothing\)。(1) 迫使左首字母取 \(a\)、左尾部全取 \(e\)，左值为 \(a\)。主理想检验完成。

\par\smallskip\noindent\textit{3. 十二条完整有限基.}\quad 
定理2。 \(A\) 的完整附加基为

\CaseDisplay{\begin{aligned}
xyz&\approx xzy, &(1)\\
xy^2&\approx xy, &(2)\\
x^2y&\approx y^2x, &(3)\\
x&\le x^2, &(4)\\
xyz&\le xy+tz, &(5)\\
xy&\le x^2+yz, &(6)\\
(xyz)^2&\le xy+yz, &(7)\\
(xy)^2&\le xy+y, &(8)\\
(xy)^2&\le xy+z^2, &(9)\\
x^2&\le x+yx, &(10)\\
x^2&\le x+xy+z^2, &(11)\\
xy&\le x+xz+ty. &(12)
\end{aligned}\quad\textup{(\ensuremath{2})}}

前三式直接由词值规则验证；其余由定理1逐式检查。前三式使首字母重复的词只由内容决定，记为 \(M_Y\)，满足 \(M_Y^2=M_Y\)。下面证明完备性。

一次有效吸收由内部纯幂见证、(1)、(2)、(4) 直接得到。以下设 \(w\) 非一次并满足判据 (1)。

\par\smallskip\noindent\textit{3.1 无首尾交集.}\quad 
设 \(H\cap E=\varnothing\)。则 \(h\in H\)、\(D\subseteq E\)，且 \(h\notin D\)。取内部右项 \(u\)。其首字母在 \(H\cap X\)，而 \(X\setminus\{h\}\subseteq E\)，故只能以 \(h\) 开头；它也不重复 \(h\)。

若 \(u\) 非一次，从 \(u\) 开始。若 \(u=h\) 为一次词，利用 \(h\in X\subseteq T_2\)，存在含 \(h\) 的非一次右项。因 \(h\notin E\)，它只能写成 \(hs\)，其中 \(s\) 非空。取任一目标尾变量 \(x\in D\)，并将相应全局尾部见证写为 \(vx\)。用 (12) 代入 \(h, x, s, v\)，得

\CaseDisplay{hx\le h+hs+vx\le F.}

于是已有内部非一次词 \(p\le F\)，以 \(h\) 开头且尾部包含于 \(D\)。对每个尚缺的目标尾变量 \(x\)，取见证 \(vx\)，写 \(p=hq\)、\(q\) 非空，利用 (5) 得 \(hqx\le hq+vx\le F\)。补齐后正是 \(w\)。

\par\smallskip\noindent\textit{3.2 存在首尾交集.}\quad 
设 \(g\in H\cap E\)。选尾部见证 \(vg\) 及首字母见证 \(gs\) 或一次词 \(g\)，由 (7) 或 (8) 得某个平方词 \(q^2\le F\)。

取内部项 \(u\in F_X\)。若 \(u=rs\) 非一次，由 (9) 得 \(u^2\le u+q^2\le F\)。若 \(u=x\) 是一次词，由 \(x\in X\subseteq T_2\) 选一个含 \(x\) 的非一次右项。当 \(x\) 在尾部时，将见证写为 \(vx\)，由 (10) 得 \(x^2\le x+vx\le F\)；否则见证以 \(x\) 开头，写成 \(xs\)，由 (11) 得 \(x^2\le x+xs+q^2\le F\)。因而总能得到内部平方种子 \(p=u^2\le F\)。

对尚缺的内容变量 \(x\in X\setminus c(p)\)，条件 \(X\subseteq T_2\) 给出非一次见证。若见证以 \(x\) 开头，写成 \(xs\)，由 (6) 得

\CaseDisplay{px\le p^2+xs=p+xs\le F.}

若 \(x\) 在尾部，将见证写成 \(vx\)，由 (5) 取两个初始变量均为 \(p\)，得

\CaseDisplay{px=p^2x\le p^2+vx=p+vx\le F.}

新词仍是平方正规形，只增加目标变量 \(x\)。依次补齐后得到 \(M_X=w^2\le F\)，由 (4) 完成。

所有代入均为非空词；两个情形涵盖任意有限变量集。把一般恒等式用分配律展开为词和，再两向逐词应用上述推导，十二式完整。证毕。

\par\smallskip\noindent\textit{4. 范围.}\quad 
来源825及其乘法对偶的有限基性来自上述任意变量数证明。其他同簇类型须分别核验双向HSP包含。这里没有从有限变量诊断推出完备性，也没有从有限基推出遗传有限基。

\subsection{奇圈刚性转移与七个五元素非有限基类型}\label{case:10}
\par\smallskip\noindent\textit{1. 记号与六元素来源.}\quad 
全部代数使用无常数签名 \((+,\cdot)\)。加法交换、结合、幂等，乘法结合，乘法对加法双侧分配。写 \(u\leq v\) 表示 \(u+v\approx v\)。 在交换乘法的自由项中，多项式表示非空有限集的交换非空词；语法相等指这些词集相同。

令 \(R=SR_6\)，载体为 \(\{0,\allowbreak{}1,\allowbreak{}2,\allowbreak{}3,\allowbreak{}4,\allowbreak{}5\}\)，运算表为

{\small\par\smallskip
\CaseLine{\(+ | 0 1 2 3 4 5\qquad * | 0 1 2 3 4 5\)}
\CaseLine{\ttfamily 0 | 0 0 0 0 0 0\qquad 0 | 0 0 0 0 0 0}
\CaseLine{\ttfamily 1 | 0 1 0 0 1 0\qquad 1 | 0 0 0 0 0 2}
\CaseLine{\ttfamily 2 | 0 0 2 0 0 0\qquad 2 | 0 0 0 0 0 0}
\CaseLine{\ttfamily 3 | 0 0 0 3 0 3\qquad 3 | 0 0 0 0 2 0}
\CaseLine{\ttfamily 4 | 0 1 0 0 4 0\qquad 4 | 0 0 0 2 0 2}
\CaseLine{\ttfamily 5 | 0 0 0 3 0 5\qquad 5 | 0 2 0 0 2 0}
\smallskip}

标签只命名语义元素，不是项中的常数。\(0\) 是加法最大元和乘法零；所有平方和所有长度至少三的词值均为 \(0\)。 乘积取 \(2\) 的无序对恰为 \(\{1,\allowbreak{}5\},\{3,\allowbreak{}4\},\{4,\allowbreak{}5\}\)。 此表取自 Shao、Ren 与 Gao 的 \(SR_6\)，将原表的每个标号减去 \(1\) \cite[Table~1]{SRG6}。下文直接从这两张运算表证明刚性引理。

\par\smallskip\noindent\textit{2. 二分二次多项式的刚性.}\quad 
引理1。 设 \(F=L+E\)，其中 \(L\) 是若干一次变量的和，\(E\) 是若干不同的无平方二次词 \(xy\) 的和；允许其中一个部分为空，但 \(F\) 非空。设 \(E\) 的无向图二分，且 \(c(L)\cap c(E)=\varnothing\)。那么对任意非空交换词 \(q\)，

\CaseDisplay{R\models q\leq F\quad\Longleftrightarrow\quad q\in F.}

证明。 若 \(E\) 为空，把 \(c(L)\) 中所有变量赋 \(1\)。任一长度至少二的目标词取 \(0\)，不低于且不等于右侧的 \(1\)；若目标包含右侧内容之外的变量，将该外变量赋 \(0\) 也分离它。一次词不在 \(L\) 时，将该变量赋 \(0\)、\(c(L)\) 赋 \(1\) 即分离。

设 \(E\) 非空，固定其二分部 \(X, Y\)。把 \(X\) 赋 \(4\)，\(Y\) 赋 \(5\)，\(c(L)\) 赋 \(2\)，其余变量赋 \(0\)。所有原边取 \(2\)，右侧取 \(2\)。这一次赋值排除了：不在 \(L\) 的一次词，长度至少三的词，平方词，包含 \(c(L)\) 或外变量的二次词，以及两端位于同一部的二次词。

余下目标只能为 \(xy\)，其中 \(x\in X, y\in Y\)。若它不是原边，改赋 \(x=1, y=3\)，其余 \(X\) 顶点为 \(4\)，其余 \(Y\) 顶点为 \(5\)，所有 \(L\) 变量为 \(2\)。每条原边的值仍为 \(2\)，因为它只能对应 \((1,\allowbreak{}5),(4,\allowbreak{}3),(4,\allowbreak{}5)\)；目标值 \(1\cdot3=0\not\leq2\)。必要性证毕，充分性由加法幂等性立即得到。

推论2（刚性）。 若上述 \(F\) 满足 \(R\models F\approx G\)，则 \(G=F\)，为交换词集的语法相等。

事实上，由引理1得 \(G\subseteq F\)；\(G\) 继承全部假设，再对 \(F\) 中每个词应用反向吸收，得到 \(F\subseteq G\)。

\par\smallskip\noindent\textit{3. 任意有界变量数恒等式不能改写长奇圈.}\quad 
对于奇数 \(m\geq5\)，令

\CaseDisplay{U_m=x_1x_2+x_2x_3+\cdots+x_{m-1}x_m+x_mx_1.}

引理3。 设 \(\Sigma\) 是 \(R\) 成立的恒等式集合，每条至多使用 \(k\) 个变量。若奇数 \(m>\max(k,3)\)，则由 \(\Sigma\) 和交换加法幂等半环公理推导的任何 \(U_m\approx V\)，都满足 \(V=U_m\) 为交换词集的语法相等。

证明。 证明任何一次等式替换均不能改变 \(U_m\) 的词集。取一条定向的 \(F\approx G\in\Sigma\)，代换 \(\varphi\)，并在分配展开后的某一上下文中使用它。单个上下文分支为 \(p\varphi(F)+r\)；\(p\) 可为空词，仅表示没有乘法上下文，绝非项签名中的单位常数。若乘法上下文是多项式，则对其每个单词分支使用相同论证。

假设整个替换前项的词集为 \(U_m\)。于是每个分支的 \(p\varphi(F)\) 都只含 \(U_m\) 的无平方二次词；幂等加法只能合并重复词，不能消去其他词。

若 \(p\) 非空，长度约束迫使 \(p\) 是一次词且 \(F\) 只含一次词。引理1和推论2遂给出 \(G=F\)。长度至少二的 \(p\) 不可能出现，因为代换像中每个词非空。

若 \(p\) 为空，\(\varphi(F)\) 只含无平方二次词。因此 \(F\) 只含一次词及无平方二次词。出现在二次词中的变量必须代换成一次变量的和；出现在一次词中的变量必须代换成二次词的和。同一变量不能同时出现在这两部分。

对 \(F\) 的每个二次图顶点，从其非空代换像中任选一个一次变量。每条边的所选端点都构成 \(U_m\) 中的一条边，因而得到 \(F\) 的二次图到圈图 \(C_m\) 的图同态。如果 \(F\) 的图不是二分图，它包含长度 \(\ell\leq k\) 的奇简单圈；同态会给出 \(C_m\) 中长度 \(\ell<m\) 的奇闭游走。任何奇闭游走都包含不长于它的奇圈，而 \(C_m\) 的唯一简单圈长度为 \(m\)，矛盾。故 \(F\) 的图二分。

由推论2仍有 \(G=F\)。反向使用恒等式时交换两端，论证不变。任何一步都保持 \(U_m\) 的词集，沿有限推导归纳即得结论。

上下文的一般形式不会造成遗漏：把含一个等式替换位置的项完全分配展开，可写成有限个非空单词乘法上下文分支及不含该位置的余项；一旦存在非空上下文分支，就由前述长度论证迫使 \(F\) 仅含一次词，否则全部分支为空上下文，使用后一论证。

\par\smallskip\noindent\textit{4. 奇圈转移定理.}\quad 
定理4。 设交换加法幂等半环 \(A\) 满足

1. \(R\in\mathsf V(A)\)；

2. 对无界多个奇数 \(m\geq5\)，\(A\) 满足

\CaseDisplay{x_1x_2\cdots x_m\leq U_m.\quad\textup{(\ensuremath{C_m})}}

则 \(A\) 非有限基。

证明。 假设有限基 \(\Sigma\) 存在，令 \(k\) 为其变量数上界。由条件1，\(\Sigma\) 的每条恒等式在 \(R\) 中成立。取条件2中的奇数 \(m>\max(k,3)\)。有限基应推出 \(U_m\approx U_m+x_1\cdots x_m\)，但这与引理3冲突，因为右侧新增了长度 \(m\) 的词。证毕。

事实上此证明说明：\(A\) 没有变量数有界的恒等式基。它明确使用条件2，绝不从条件1单独推导非有限基性。

备注（非交换推广）。 将 \(U_m\) 换成双向和 \(\widetilde U_m=\sum_i(x_ix_{i+1}+x_{i+1}x_i)\)，乘积按 \(x_1\cdots x_m\) 的指定顺序。若任意半环 \(A\) 满足来源包含和相应的无界奇圈吸收，则同样非有限基。把任何假设推导投影到自由交换词集；每条基恒等式的交换投影仍在 \(R\) 成立，引理3保证投影保持 \(U_m\)，但目标仍多出长度 \(m\) 的词。七个本批类型均交换，登记不依赖此推广。

\par\smallskip\noindent\textit{5. 全部奇数长度的有限精确核验.}\quad 
条件2并非由有限个短圈抽检推测。以下有限状态系统与全部长度的赋值严格对应。状态

\CaseDisplay{s=(f,\ell, p, b,\epsilon)}

依次记录第一个赋值、最后一个赋值、依次相乘所得积、尚未闭合的双向相邻边之和、长度奇偶性。从所有长度二赋值 \((u, v)\) 初始化

\CaseDisplay{(u, v, uv, uv+vu,0).}

追加值 \(z\) 的转移是

\CaseDisplay{(f,\ell, p, b,\epsilon)\longmapsto
(f, z, pz, b+\ell z+z\ell,1-\epsilon).}

对每个可达奇状态检查

\CaseDisplay{p\leq b+\ell f+f\ell.}

每个长度至少三的奇数赋值都沿追加操作到达一个奇状态；每个可达奇状态都有这样的赋值。故所有可达奇状态通过，当且仅当所有奇数 \(m\geq3\) 的 \(C_m\) 成立。状态总数至多 \(2|A|^4\)，闭包停止及全转移覆盖均可直接复核。

七型的可达状态数及转移数：

\noindent 编号\enspace\textbar\enspace 状态数\enspace\textbar\enspace 转移数\enspace\textbar\enspace 完整SR6包含的词函数数\par
\noindent 191\enspace\textbar\enspace 108\enspace\textbar\enspace 540\enspace\textbar\enspace 25\par
\noindent 212\enspace\textbar\enspace 114\enspace\textbar\enspace 570\enspace\textbar\enspace 25\par
\noindent 891\enspace\textbar\enspace 108\enspace\textbar\enspace 540\enspace\textbar\enspace 25\par
\noindent 13249\enspace\textbar\enspace 188\enspace\textbar\enspace 940\enspace\textbar\enspace 25\par
\noindent 13250\enspace\textbar\enspace 138\enspace\textbar\enspace 690\enspace\textbar\enspace 25\par
\noindent 13894\enspace\textbar\enspace 188\enspace\textbar\enspace 940\enspace\textbar\enspace 25\par
\noindent 13895\enspace\textbar\enspace 138\enspace\textbar\enspace 690\enspace\textbar\enspace 25\par
来源包含使用四个生成元 \(1,\allowbreak{}3,\allowbreak{}4,\allowbreak{}5\)，在全部 \(5^4=625\) 个赋值上完整闭合词函数，并检查六个主下集，构造从相对自由代数到 \(R\) 的满同态；没有设置词长上限。独立核验另以不同状态编码重算奇圈闭包，并比较完整终态集，而不是只比较候选数。

奇圈构造采用 Shao、Ren 与 Gao 的六元素来源（\cite[Table~1]{SRG6}）。上述刚性引理与奇长度闭包给出了此处转移所需的两个条件。

\subsection{Five complete identities for semiring 194}\label{case:11}
Constant-free signature (+, multiplication). Relative to the additively idempotent semiring axioms, a complete finite basis is

{\small\par\smallskip
\CaseLine{\(xy=yx; xyz=ttt; xy\leq ttt; ttt\leq x+xy; xy\leq xx+yz.\)}
\smallskip}

As usual, \(u\leq P\) abbreviates \(u+P=P\). Addition: 00004/\allowbreak{}01004/\allowbreak{}00224/\allowbreak{}00234/\allowbreak{}44444. Multiplication: 00000/\allowbreak{}00000/\allowbreak{}00010/\allowbreak{}00110/\allowbreak{}00000. The additive order has \(3<2<0<4\) and \(1<0<4\). A quadratic is 1 at(2,\allowbreak{}3),(3,\allowbreak{}2),(3,\allowbreak{}3), and 0 otherwise. All long words are 0. This constant cube is NOT the additive maximum: it dominates all nonlinear words but does not dominate the element 4. The law \(x\leq ttt\) is false and is not used.

\par\smallskip\noindent\textit{Exact original-polynomial criterion.}\quad 
Let L be original literal variables, E original unordered quadratic edges, V their endpoints, S original loop variables. A literal target is absorbed iff it is original. A nonlinear target is absorbed automatically if P has an original long word or L intersects V. Otherwise a pair ht is absorbed iff its endpoints lie in V and its edge is original or an endpoint belongs to S. No long target is absorbed in this remaining case.

For sufficiency, an original long word has value 0 and absorbs every nonlinear value 0 or 1. An original literal x with incident edge xy has join at least 0:if \(xy=0\) this is immediate; if \(xy=1\) then x is 2 or 3 and \(x+1=0\). With no such witnesses, consider a right polynomial whose value does not dominate 0. If it has an edge, all edges must be 1 and \(P=1\). Original endpoints are 2 or 3 with no edge of two2s. An original loop forces 3, so a target incident pair is 1. Original edges are 1. With no edges, only original literal targets qualify.

For necessity of the literal clause, assign a missing literal \(h=4\) and all other variables 3. No right literal is 4; right nonlinear values are 0 or 1, so the right side is at most 0 while \(h=4\). In the inactive case with no long word and L disjoint from V, assign V to 3 and outsiders 1. All right terms are 1. Any long target is 0, not below 1; a pair with an outside endpoint is also 0. For a missing edge with neither endpoint in S, assign both target endpoints 2, other V vertices 3, and outsiders 1. No right edge has two2 endpoints, and every right literal is 1. The target is 0 and \(P=1\). This covers every failure for arbitrary variable sets and polynomial sizes.

\par\smallskip\noindent\textit{Unrestricted complete derivations.}\quad 
The five laws are valid by the table formula and criterion; all fixed-law assignments are also independently checked. Normalize long words to a fixed cube using \(xyz=ttt\) with a nonempty prefix. An original incident literal x and edge xy give this cube by \(ttt\leq x+xy\). Once the cube is available, \(xy\leq ttt\) derives any nonlinear target, using a nonempty split of the target word. It never derives a new literal. Otherwise a new pair has an original loop at one endpoint and an original edge at the other. Orient them by commutativity and use \(xy\leq xx+yz\). Original targets require no inference. Every witness is original or explicitly derived below P, and no empty substitutions occur.

These steps prove all valid absorptions without bounds on variables, lengths, or summand counts. Distributivity and summandwise absorption in both directions prove all identities. Thus the five laws are complete. Other types, including similar constant-cube tables, require independently verified complete HSP in both directions before transfer.

\subsection{Eight complete identities for semiring 200}\label{case:12}
Signature (+, multiplication), without constants. The following eight identities form a complete finite basis relative to the additively idempotent semiring axioms:

{\small\par\smallskip
\CaseLine{\(xy=yx; xyz=xyzz; xyzt\leq xyz; xxx\leq x+yyy;\)}
\CaseLine{\(xxy\leq xy+zzz; xy\leq xxy; xy\leq xx+yz; x\leq xx.\)}
\smallskip}

As usual \(u\leq P\) abbreviates \(u+P=P\).

\par\smallskip\noindent\textit{Exact word values and absorption.}\quad 
Addition is 00004/\allowbreak{}01114/\allowbreak{}01224/\allowbreak{}01234/\allowbreak{}44444, the chain \(3<2<1<0<4\). Multiplication is 00000/\allowbreak{}00000/\allowbreak{}00010/\allowbreak{}00110/\allowbreak{}00004. A nonlinear word is 4 when all its content is assigned 4. Otherwise its value is 1 precisely when it is quadratic, both occurrences are assigned in \{2,\allowbreak{}3\}, and at least one is assigned 3. All other nonlinear words are 0.

For an original polynomial P let E be its unordered quadratic edges, V their endpoint union and S its loop variables. Put \(G=true\) exactly when an original word has length at least three. Let C(w) say some original u has c(u) contained in c(w). The exact absorption criterion is:

\textbullet{} Literal target: C(w).

\textbullet{} Long target: C(w) and G.

\textbullet{} Quadratic target: C(w), and either G, or both c(w) contained in V and (the target edge is original or a target endpoint belongs to S).

All witnesses in this statement are ORIGINAL. In particular, an endpoint literal does not replace the original loop requirement for a new edge.

If C(w) fails, assign target content 4 and all outsiders 3; every right summand is at most 0, whereas the target is 4. If G fails, assign V to 3 and outsiders 1. All right summands are at most 1, but any long target is 0, as is a quadratic with an endpoint outside V. For a new quadratic on V having neither an original edge nor an endpoint loop, assign target endpoints 2 and all other variables 3. Every original edge has at least one endpoint outside the target and therefore value 1; literals are at most 1. The target is 0. These give explicit separators for all failed conditions and arbitrary numbers of variables.

Sufficiency for a literal follows from a contained pure power: x is below \(x^{k}\) for every \(k\geq 1\), directly from the five possible values. When G holds, \(P\geq 0\); the contained witness supplies 4 whenever a nonlinear target is 4, and all other values are below 0. Without G, a quadratic on V cannot exceed 1 through an endpoint outside \{2,\allowbreak{}3\} unless an original incident edge already has value 0 or

4. If all original edges have value 1, all V variables lie in \{2,\allowbreak{}3\}, and an original loop forces its endpoint to 3. A target having that loop endpoint is then 1. Original target edges need no argument. The contained original loop or edge handles the all-4 assignment.

\par\smallskip\noindent\textit{Unrestricted complete derivations.}\quad 
Commutativity and \(xyz=xyzz\) identify long words with the same content. Move a chosen occurrence to the end, duplicate using two nonempty prefix blocks, and reverse the identity to delete surplus copies without reducing below length three. Canonical forms are sorted content for at least three distinct variables, xxy for two and xxx for one.

A literal target has a contained original pure power. If it is already literal, retain it. If square, apply \(x\leq xx\). If long, normalize to xxx, apply \(xy\leq xxy\) with \(x=y\) to derive xx, then \(x\leq xx\). This uses no global long witness and no literal-edge overlap.

For a nonlinear target with G, fix an original long g and a contained original u. Upgrade a literal u to \(u^{3}\) by \(xxx\leq x+yyy\) with \(y=g\), or a quadratic ab to aab by \(xxy\leq xy+zzz\) with \(z=g\); \(g^{3}\) normalizes to g. Retain a long u. Append any missing target content with \(xyzt\leq xyz\) and normalize. A quadratic ab is recovered by \(xy\leq xxy\). Every substitution is nonempty.

If G fails, a new quadratic must have an original loop aa and an original edge bc incident at the other endpoint. Apply \(xy\leq xx+yz\), commuting endpoints if necessary. Literal recovery does not justify inventing any new loop or dropping the original endpoint witness.

Each displayed law satisfies the exact criterion and is independently checked on all assignments. The derivations cover all semantic cases without bounds on variables, words or summands. Distributivity followed by summandwise absorption in both directions proves every identity. Other catalogue types require full independent HSP equivalence to transfer.

\subsection{A finite basis for 201 using at most five original witnesses}\label{case:13}
Signature without constants. The structural identity is \(xyz=xxx\). In addition take all table-valid absorptions over 18 variables, word length at most 3, and at most 5 upper words. This is an effectively specified finite family, not an enumerated small list.

The multiplication rows 00000,\allowbreak{}00000,\allowbreak{}00010,\allowbreak{}00110,\allowbreak{}44444 show:every nonlinear word with head 4 has value 4; a quadratic with head not 4 is 1 exactly for first two values(2,\allowbreak{}3),(3,\allowbreak{}2),(3,\allowbreak{}3), and otherwise0. A word of length at least 3 with head not 4 always 0. Literals retain their values. Induction proves this for every length. Thus \(xyz=xxx\) is valid and a single substitution \(x=head\), \(y=second\), \(z=remaining\) nonempty suffix reduces any long word to its head cube. This uses no empty substitution.

Normalize all original words accordingly. Let H be their heads, L their literal variables, Q the collection of quadratic originals, C the union of quadratic content, and J the collection of normalized cubes. The addition has proper principal ideals\{0,\allowbreak{}1,\allowbreak{}2,\allowbreak{}3\},\{1\},\{2,\allowbreak{}3\},\{3\}.

The first ideal imposes target head in H, since a word is below 0 exactly when its head is not 4. Ideals\{2,\allowbreak{}3\} and\{3\} contain only literal word values, so a nonlinear target requires an original nonlinear word. With the head condition this also suffices for those two ideals:if originals are all literals and the target is literal, its variable occurs there.

For all originals to lie in\{1\}, there must be no long original, literal variables must equal 1, and quadratic variables must lie in\{2,\allowbreak{}3\} with at least one 3 on each original pair. This is possible exactly when J is empty and L intersect C is empty. Feasibility follows by giving all quadratic variables 3 and literals 1. If impossible, that ideal imposes no condition. If possible, a literal target must be in L; a long target can never lie in\{1\}; and a quadratic target xy is always 1 precisely when its content is contained in C and SOME original quadratic has content contained in\{ x, y \}.

To justify the last necessity, regard value 3 as Boolean true and value 2 as false. Original quadratics are positive clauses of size one or two. If no original clause is contained in the target variables, set the target variables 2 and every other quadratic variable 3. All originals remain 1 while the target is 0. Missing quadratic content is separated by putting a missing variable 0(or its forced literal value 1). Sufficiency follows directly from the included clause and content restrictions.

Thus \(w\leq P\) is exactly:head(w) in H; an original nonlinear word exists when w is nonlinear; and the ideal 1 condition just stated holds. These proper principal-ideal tests prove unrestricted necessity/sufficiency. For a failed head test put target head 4 and all original heads 0. For failed nonlinearity put all variables 3, so \(P=3\) while the nonlinear target is 1 or 0. The ideal 1 assignments above also supply all other separators.

At most 5 original normalized words suffice. Choose a target-head word and, if needed, a nonlinear word. If ideal 1 is impossible, choose one long original, or a literal and an original quadratic containing that literal variable. At most 4 words suffice in this case. If ideal 1 is possible, a literal needs its literal original. A quadratic needs at most 2 original quadratic content witnesses and the original included clause. Together with the first 2 choices this is at most 5. No valid long target occurs in this latter case. The selected words already absorb the target by the exact criterion.

The target and at most 5 selected words have length at most 3, so jointly use at most 18 variables. Rename them injectively into a fixed alphabet of 18 variables. The renamed absorption belongs to the stated finite family. The inverse singleton substitution literally recovers the selected words. Each selected word is an original word normalized by \(xyz=xxx\); reverse these equalities to infer the original absorption.

This proves every valid absorption for arbitrary original alphabets, lengths, and polynomial sizes. Distributivity and the two absorption directions prove completeness for every valid identity. Therefore 201 is finitely based. Its recorded proper NFB subvariety does not contradict this and must remain recorded as an additional nonlimit reason. Further catalogue transfers require independent complete mutual HSP.

\subsection{Forty-two finite-basis types from bounded word length}\label{case:14}
Signature: \((+, multiplication)\), with no constants. The classification theorem below applies to arbitrary numbers of variables, word lengths and polynomial sizes. The finite absorption families are specified exactly; they have not been enumerated as small lists of identities.

\par\smallskip\noindent\textit{1. Scope and identities.}\quad 
All the direct sources in this case have additive least element 3. Their length reductions are listed below.

The common left side is xyzt. The right sides are:

\noindent Direct source IDs\enspace\textbar\enspace Right side\par
\noindent 202\enspace\textbar\enspace xxx\par
\noindent 767, 906, 1042, 1069, 1949, 9147, 9153, 10854, 10882, 10885, 12433, 12916, 13691, 13745\enspace\textbar\enspace xyz + xyt + xzt\par
\noindent 774, 12331, 13426, 13715\enspace\textbar\enspace xyz + xyt + xzt + yzt\par
\noindent 900, 1051, 2558, 2612, 6317, 9369, 10876, 10883, 12913, 12917, 13776, 13778\enspace\textbar\enspace xyt + xzt + yzt\par
\noindent 2308\enspace\textbar\enspace xxt\par
Repeated letters in xxx and xxt are actual repeated occurrences, not content sets. No empty word or constant is substituted.

\par\smallskip\noindent\textit{2. Unrestricted completeness.}\quad 
We apply Theorem 2 and Corollary 3(2) of Section~\ref{thm:bounded}. For completeness, here is the entire specialization, including the bound on the number of right summands.

Fix a direct source A and its displayed identity gamma. In any word of length at least four, replace its first four letters using gamma, retain the suffix, and distribute. Every resulting word has one fewer letter, including when the first letters coincide and when gamma repeats a letter. Induction on word length therefore reduces every word to a nonempty sum of words of length at most three. The semiring axioms first expand arbitrary terms into nonempty sums of nonempty words, so this reduction applies to all terms.

Let \(E_{A}\) consist of all valid absorptions w \ensuremath{\leq} Q over a fixed alphabet of 1,\allowbreak{}878 variables, with 1 \ensuremath{\leq} |w| \ensuremath{\leq} 3, each word of Q of length 1 to 3, and 1 to 625 distinct right summands. This is a finite set. We claim that the semiring axioms, gamma and \(E_{A}\) form a complete basis.

For any valid absorption w \ensuremath{\leq} Q after reduction, put C = content(w) and d = |C| \ensuremath{\leq} 3. For each assignment alpha:C \ensuremath{\to} A, assign all other variables the additive least element. From the original right summands of Q select one for each value actually attained under that extension. There are at most 5 selected words per alpha. The union \(Q_{0}\) of these selections is nonempty, has at most \(5^{d+1}\) \ensuremath{\leq} 625 words and uses at most 3 + 3*625 = 1,\allowbreak{}878 variables together with w.

For any assignment beta to all original variables, choose the selection corresponding to beta restricted to C. Raising the variables outside C from the least element to their beta values can only raise each word, since multiplication preserves the additive order by distributivity. Thus \(Q_{0}\)(beta) \ensuremath{\geq} w(beta). Consequently w \ensuremath{\leq} \(Q_{0}\) is valid. An injective renaming puts this identity into \(E_{A}\); reversing the renaming and adding the remaining original right summands derives w \ensuremath{\leq} Q. No variables are identified, and every multiplicity is retained. The least element is used only in semantic assignments, never as a term or a signature constant.

Every valid equation between reduced polynomials is equivalent to the absorptions of each left summand into the right polynomial and conversely. The preceding argument derives each absorption. Undoing the reductions proves the original equation. This proves completeness without a bound on the variables, lengths or sizes in the equation being derived.

\par\smallskip\noindent\textit{3. Complete equal-variety transfers.}\quad 
The following are complete two-direction HSP checks, not implications inferred from matching short identities. op means reverse multiplication. Each source in this table is among the 32 proved sources above. Its opposite is FB by reversing every word of its complete basis. Equality of varieties then transfers exactly that complete basis to the target.

\noindent Source\enspace\textbar\enspace Orientation\enspace\textbar\enspace Target IDs\par
\noindent 202\enspace\textbar\enspace original\enspace\textbar\enspace 202, 12202\par
\noindent 202\enspace\textbar\enspace op\enspace\textbar\enspace 2308, 12196\par
\noindent 767\enspace\textbar\enspace original\enspace\textbar\enspace 767, 12433\par
\noindent 767\enspace\textbar\enspace op\enspace\textbar\enspace 900, 12913\par
\noindent 774\enspace\textbar\enspace original\enspace\textbar\enspace 774, 12331\par
\noindent 906\enspace\textbar\enspace original\enspace\textbar\enspace 906, 12916, 12997\par
\noindent 906\enspace\textbar\enspace op\enspace\textbar\enspace 2558, 12917, 12938\par
\noindent 1042\enspace\textbar\enspace original\enspace\textbar\enspace 1042, 1069, 1830, 1949, 9153, 9325, 10854, 10882, 11104\par
\noindent 1042\enspace\textbar\enspace op\enspace\textbar\enspace 1051, 2612, 5876, 6317, 9369, 9987, 10876, 10883, 10990\par
\noindent 9147\enspace\textbar\enspace original\enspace\textbar\enspace 9147, 10885\par
\noindent 13426\enspace\textbar\enspace original\enspace\textbar\enspace 13426, 13715\par
\noindent 13691\enspace\textbar\enspace original\enspace\textbar\enspace 13691, 13745\par
\noindent 13691\enspace\textbar\enspace op\enspace\textbar\enspace 13776, 13778\par
Proposition~\ref{prop:lower-fibres} tests \(\A{j}\in\V(\A{i})\) by closing the word functions of \(\A{i}\) on a generating tuple of \(\A{j}\) under each right-generator multiplication. Transition consistency defines the target value of every word. The principal-ideal test determines whether this evaluation extends to sums. Interchanging \(i,j\) gives the reverse inclusion.

All 42 targets were already excluded from limit generation by a recorded proper NFB subvariety. These reasons are retained. The new FB results do not contradict those reasons and do not prove hereditary finite basability.

\subsection{An odd-cycle obstruction for 207 and a rigorous transfer criterion}\label{case:15}
The signature has no constants. The following argument proves NFB for 207; computational certificates and target registrations are separate. It is not inferred from an NFB proper subvariety.

The addition rows are 00004,\allowbreak{}01004,\allowbreak{}00204,\allowbreak{}00034,\allowbreak{}44444. Thus 4 is top, 0 is above 1,\allowbreak{}2,\allowbreak{}3, and their pairwise distinct joins are 0. Multiplication rows are 00000,\allowbreak{}00000,\allowbreak{}00010,\allowbreak{}00120,\allowbreak{}00004. In particular 0 is absorbing, 4 is idempotent and annihilates1,\allowbreak{}2,\allowbreak{}3 to 0,\allowbreak{}3\ensuremath{{}^{2}}=2,\allowbreak{}3\ensuremath{{}^{3}}=1,\allowbreak{}3\ensuremath{{}^{4}}=0, and \(23=32=1\). A quadratic is 1 exactly at(2,\allowbreak{}3),(3,\allowbreak{}2); a cubic is 1 exactly when all three entries are 3. A word of length at least 4 on values 0,\allowbreak{}1,\allowbreak{}2,\allowbreak{}3 is 0. A word all of whose entries are 4 has value 4.

For odd \(n\geq 3\) use distinct variables x, v0,..., v(n-1), and put

\(P_{n}=x\)\ensuremath{{}^{3}}+\(sum_{i}\) (\(v_{i}\) \(v_{i+1}\)+\(v_{i+1}\) \(v_{i}\)), indices modulo n.

The critical identity is \(x\leq P_{n}\). It holds in 207:for \(x=0\),1,\allowbreak{}2,\allowbreak{}4, \(x\leq x\)\ensuremath{{}^{3}}. For \(x=3\), the cube is 1. If \(P_{n}\) is 4 or 0, it absorbs 3. The only possible remaining sum is 1, which requires every edge value 1. That would properly2-color an odd cycle by 2 and 3, impossible. Thus every odd-n critical identity is valid. Even cycles give counterexamples by putting \(x=3\) and alternating vertex values 2,\allowbreak{}3.

\par\smallskip\noindent\textit{A substitution barrier with a proved variable bound.}\quad 
Let \(F_{n}\) be the set consisting of the word \(x^{3}\) and BOTH orientations of all edges of the simple n-cycle. Take \(n\geq 5\), odd, and \(n>2k\)\ensuremath{{}^{2}}. Let q be any nonempty polynomial on at most k source variables. Let phi substitute a nonempty polynomial for each variable, and suppose every expanded word of phi(q) belongs to \(F_{n}\). We claim:for every absorption \(u\leq q\) valid in 207, where u is one source word, EVERY expanded word of phi(u) also belongs to \(F_{n}\).

Variables absent from q may be assigned 0 below, so they cannot occur in a word u whose value must be 1. Consider variables occurring in q. A variable occurring literally in q has image contained in \(F_{n}\). Such a variable cannot also occur in a nonlinear source word:an image edge would then produce a word of length at least 3 containing cycle letters, or an image cube would produce length at least 4; neither is in \(F_{n}\). Call these variables L. Give each L variable the formal algebra value 1. This is a FREE source assignment, not an assertion that evaluation of its image polynomial is 1 on the entire odd cycle.

Every other source variable occurring in q is of exactly one of three types:its image is the singleton word x; its image is the singleton word \(x^{2}\); or its image is a nonempty set of single cycle letters. Call these X, Y, V respectively. Here is why this list is exhaustive. Expanded source words have length 2 or 3. A degree 3 source word must expand only to \(x^{3}\), so all of its factors have image\{x\}. A degree 2 source word may expand to \(x^{3}\), requiring images\{x\} and\{\(x^{2}\)\} in either order; or may expand to cycle edges, requiring two sets of single cycle letters. Mixing these possibilities in one image introduces x with a cycle letter in a degree 2 or 3 word, or a degree 4 word, which is excluded. In particular, no image of an internal variable contains a sum mixing these three types.

The nonlinear words of q consequently have only these forms:three X variables; one X and one Y; or two V variables. Assign \(X=3\), \(Y=2\). A V-pair with image sets A, B contributes the complete bipartite edge set A times B to the cycle. A and B are disjoint(no loops). A cycle with \(n\geq 5\) has no 4-cycle and maximum degree 2, so each such complete bipartite subgraph has at most 2 undirected edges. There are at most \(k^{2}\) distinct ordered quadratic source words. Thus the union H of all these image edges has at most 2k\ensuremath{{}^{2}}\ensuremath{<}n edges, and is a proper subgraph of the n-cycle, hence a forest.

Properly2-color H by values 2,\allowbreak{}3. Every V-image set is monochromatic: the variable occurs in a V-pair, and all its image vertices share any one neighbor from the other nonempty image set. Assign each V variable this common value. Independent color flips on forest components are allowed. Each source word of q now has value 1:the linear variables are 1, \(XXX=3\)\ensuremath{{}^{3}}=1, \(XY/YX=32/\allowbreak{}23=1\), and V-pairs have opposite colors. Since \(u\leq q\) is valid in 207, u must have value 1 for EVERY such coloring.

We use one further necessary condition valid for EVERY absorption in 207: some word of q has content contained in content(u). Otherwise assign the variables of u the value 4, and all other variables 0. Then \(u=4\), every word of \(q=0\), and \(u\leq q\) fails. This argument uses the absorbing 0 and idempotent4 and is independent of the cycle assignment.

If u is literal, it must be an L variable; other occurring types have values 2 or 3 and absent variables 0. Its image is already contained in \(F_{n}\). If u is cubic, it must use only X variables:an L or Y has value 1 or 2, while any V variable can be flipped to 2. Thus its image is\{\(x^{3}\)\}. Words of length at least 4 have value 0 and cannot occur. Finally let u be quadratic. It cannot contain an L or absent variable. A mixed X/Y pair has image\{\(x^{3}\)\}. Two X or two Y variables do not give 1. A V paired with X or Y can be flipped to the wrong color and is impossible. For two V variables a, b, the contained-content condition forces a word of q using only a, b. It cannot be literal(a, b are not L), cubic(cubics of q use only X), or a square(a V square would introduce cycle loops). Hence q already contains ab or ba. Therefore all words of phi(u) are cycle edges in \(F_{n}\). This completes the barrier for ARBITRARY nonempty polynomial substitutions, including repeated source variables.

\par\smallskip\noindent\textit{Unrestricted nonfinite basability.}\quad 
Suppose 207 had a finite identity basis. Let k bound its variable counts and those of the semiring axioms. Choose odd \(n>2(k+2)\)\ensuremath{{}^{2}} and \(n\geq 5\). Normalize every term in a purported equational derivation of \(P_{n}=P_{n}+x\) as a nonempty sum of nonempty associative words. Consider the FIRST step introducing any word outside \(F_{n}\). Earlier steps may delete words, but every earlier support is still contained in \(F_{n}\).

The application of a basis identity within a multiplicative/additive context expands into monomial contexts. Select one component producing the first new word. Its old component has the form a phi(q) b with possibly absent left/right contexts; its whole support is in \(F_{n}\). Every word u of the other side satisfies \(u\leq q\) in 207. Incorporate nonempty left and right contexts as substitutions for at most two additional fresh variables in the valid contextual absorption a u \(b\leq a\) q b. Omit a context when absent; do not use an empty-word substitution. The barrier with k+2 variables says no new word can be introduced, a contradiction. Pure semiring normalization changes no expanded support. The final polynomial contains the literal x outside \(F_{n}\), so a first new word must exist. Thus 207 is nonfinitely based.

\par\smallskip\noindent\textit{Transfer to other finite semirings.}\quad 
Let T be any additively idempotent semiring with 207 in Var(T), and suppose T satisfies all the odd-cycle identities \(x\leq P_{n}\). Then T is NFB by the SAME barrier:every identity of T is valid in 207, while a finite basis for T would have to derive a sufficiently large critical identity valid in T. Both hypotheses are essential.

For finite T, the all-length critical hypothesis is decidable exactly. For each pair of values a, b with \(a^{3}\)\ensuremath{\leq}b and a \(\nleq b\), form the undirected graph on T whose edges\{r, s\}, including possible loops, satisfy both \(rs\leq b\) and \(sr\leq b\). Then all odd-cycle identities hold iff every one of these graphs is bipartite(with no loops). Indeed an assignment violating an odd-cycle identity gives an odd closed walk in such a graph, using b equal to its polynomial value. Conversely any odd closed walk gives an assignment with the whole \(polynomial\leq b\) but a \(\nleq b\), so the identity fails. A loop repeats to give a walk of length 3. An undirected graph has an odd closed walk iff it has an odd cycle or loop. This finite criterion proves ALL lengths, not only finitely many tested cycles.

The transfer to \(T\) also requires \(\A{207}\in\V(T)\), tested by Proposition~\ref{prop:lower-fibres}.

\subsection{Finite bases for 13891 and 13892 from bounded short-word witnesses}\label{case:16}
The signature is addition and multiplication without constants. The bases below contain explicitly defined finite families; they are not claimed to be enumerated small lists. Completeness is unrestricted in variables, word lengths and polynomial sizes.

Both multiplication tables are 44434/\allowbreak{}44034/\allowbreak{}40134/\allowbreak{}33333/\allowbreak{}44434. Addition for 13891 is 00004/\allowbreak{}01114/\allowbreak{}01224/\allowbreak{}01234/\allowbreak{}44444; addition for 13892 is 00044/\allowbreak{}01144/\allowbreak{}01244/\allowbreak{}44434/\allowbreak{}44444. Element 3 absorbs multiplication. On\{2,\allowbreak{}1,\allowbreak{}0,\allowbreak{}4\}, give the elements weights 1,\allowbreak{}2,\allowbreak{}3,\allowbreak{}4 respectively. Multiplication adds weights, truncated at 4; addition takes their maximum. In 13891, element 3 is the least additive element. In 13892, it is incomparable with 2,\allowbreak{}1,\allowbreak{}0 and joins any of them to 4. These facts follow directly from the tables and describe products and sums at every length.

For any word, an occurrence assigned 3 makes its value 3. Otherwise its weight is the sum of its occurrence weights, truncated at 4. In particular, every word of \(length\geq 4\) has value 4 unless an occurrence is 3. Its value therefore depends only on content. Cubic word values also depend only on content: all occurrence weights 1 give weight 3, and any larger weight gives 4. Quadratics retain their multiplicities.

\par\smallskip\noindent\textit{The finite bases.}\quad 
For each source, include the semiring axioms and these six structural laws:

{\small\par\smallskip
\CaseLine{\(xy=yx; xxy=xyy; xxyzt=xyzt;\)}
\CaseLine{\(xxxy=xxyy; xxyz=xyyz; xxxxy\leq xxxx.\)}
\smallskip}

For 13892 include the additional law:

{\small\par\smallskip
\CaseLine{\(xxxx\leq x+yzts.\)}
\smallskip}

For 13891 also include every valid absorption \(u\leq P\) over a fixed alphabet of three variables, with \(1\leq\)|u|\ensuremath{\leq}3 and \(1\leq\)|v|\ensuremath{\leq}4 for each right word v, and at most 120 distinct right words. This is a finite set.

For 13892 include every valid absorption \(u\leq P\) over a fixed alphabet of 207 variables, with \(1\leq\)|u|\ensuremath{\leq}3 and \(1\leq\)|v|\ensuremath{\leq}4 for each right word v, and at most 68 distinct right words. This too is a finite set. Valid means valid in the displayed finite source algebra; it specifies the family unambiguously without asserting that its members have been enumerated.

\par\smallskip\noindent\textit{Multiplicative normalization.}\quad 
Commutativity sorts a word. Cubics on two variables are normalized by \(xxy=xyy\). For words of \(length\geq 5\) with a repeated variable, commutativity and \(xxyzt=xyzt\) delete one repeated occurrence while retaining \(length\geq 4\): the other three substituted factors split a nonempty remainder of length at least 3. Continue until the word is linear or has length 4. With one, two or three distinct variables, the two fourth-degree transfer identities choose a fixed distribution of the four occurrences. Thus any two long words with the same content are provably equal. No short word is equated with a long one. Every factor substitution is nonempty.

The law \(xxxxy\leq xxxx\) allows adjoining arbitrary additional content to a long word: normalize the original long word U to \(U^{4}\), apply the law with \(x=U\) and a nonempty word for y, and normalize again. If no new content is needed, only normalization is used. Original repeated occurrences are retained in all short words and in the weight comparisons below.

\par\smallskip\noindent\textit{Exact criterion for 13891.}\quad 
For a target w, let \(P_{C}\) be the original right words whose content is contained in \(C=content(w)\). A target of \(length\geq 4\) is absorbed exactly when \(P_{C}\) contains an original word of \(length\geq 4\). A target of \(length\leq 3\) is absorbed exactly when \(P_{C}\) is nonempty and, for every assignment of weights 1,\allowbreak{}2,\allowbreak{}3,\allowbreak{}4 to C, the target's truncated weight is at most the maximum truncated weight of \(P_{C}\).

Indeed, when some target letter is 3, the target is least and is absorbed. Otherwise assign all outside variables 3; this removes every noncontained right word and proves necessity. In other evaluations those right words can only increase the sum, so the same condition suffices. In the long target case, assign all target weights 1: only an original contained long word can reach 4. Conversely such a word has weight 4 whenever the target does not vanish, and always absorbs the target.

For a short target, normalize each contained right word. It uses at most three variables and has \(length\leq 4\). There are at most \(3+3^{2}+3^{3}+3^{4}=120\) possible words on that alphabet. Deleting duplicates gives exactly an instance of the specified finite family; the other original words remain available by additive monotonicity. For a long target, start with the original contained long witness and adjoin missing target content as above. This derives every valid absorption in 13891.

\par\smallskip\noindent\textit{Exact criterion for 13892.}\quad 
If P has an original word of \(length\geq 4\), then \(w\leq P\) exactly when at least one original right word has its content contained in \(C=content(w)\). If no original right word is long, a valid target must have \(length\leq 3\). In that case the criterion consists of all three conditions:

1. Some original right word is contained in C.

2. Every target variable occurs somewhere in the original P.

3. For each assignment of weights 1,\allowbreak{}2,\allowbreak{}3,\allowbreak{}4 to C, assign weight 1 to all outside variables. The target's truncated weight is at most the maximum truncated weight of the original right words.

Containment is necessary: otherwise assign every target variable 2 (weight 1) and all outside variables 3. All right words become 3 while the target does not. With no long right word, an outside target variable can be assigned 3 while every original variable is 2; the right has weight at most 3 and does not absorb the exceptional3. Assigning all variables 2 also rules out a long target when no long right word exists. The weight condition follows by evaluating exactly its indicated assignments.

For sufficiency, if both 3 and a chain element occur as right values, their sum is 4. If all right values are 3, the contained original witness forces a3 occurrence in the target. If none is 3, the long right word, when present, already gives 4. Without a long right word, content coverage ensures the target also has no 3 occurrence; raising outside chain weights from 1 can only increase the right sum, so condition 3 suffices. These cases exhaust every evaluation at arbitrary variable count.

\par\smallskip\noindent\textit{The 207-variable witness bound.}\quad 
First suppose no original right word is long. A target has at most 3 variables. For each of its at most \(4^{3}=64\) weight assignments, select ONE original right word realizing the required upper bound. Keep in addition one contained original word and, for each of at most 3 target variables, one original word containing it. The selected polynomial Q has at most 68 words, each of \(length\leq 3\), and satisfies all three exact conditions. It uses at most \(3+3*68=207\) variables in total. Thus \(w\leq Q\) is a renaming instance of the declared finite family. No outsider is identified with another outsider: the distinct-variable count itself is bounded. Repeated occurrences remain literal occurrences in all 68 retained words. Additive monotonicity then yields \(w\leq P\).

Now suppose P has an original long word V and the target is short. Choose an original contained witness U and normalize U. It has length \ensuremath{\leq}4 on the at most 3 target variables. Split V into four nonempty blocks, and replace those blocks by four fresh variables to form the schematic right word abcd. The absorption \(w\leq U+abcd\) is valid by the exact long- right criterion: U is contained in w. It uses at most 7 variables, two right words of \(length\leq 4\), and belongs to the finite family. Substituting the four original blocks recovers U+V, while fixing every target variable. The blocks may share variables or contain target variables; simultaneous substitution reproduces all original occurrences exactly.

Finally suppose both V and the target are long. The contained witness U can be short or long. Split V into four nonempty blocks and apply \(xxxx\leq x+yzts\) with \(x=U\) to get \(U^{4}\leq U+V\). Adjoin any target content not in U, then normalize the resulting long word to the target. This proves every remaining absorption using the fixed structural laws.

Every law in either finite family is valid by definition, and the six or seven structural laws follow from the word formula. The constructions cover all valid absorptions; distributive expansion and absorption in both directions prove full equational completeness for both sources.

Boundary distinctions: \(xxx\leq xx\) fails at \(x=2\) in both sources. In 13892, \(x\leq y+z\) fails at \(x=3\), \(y=z=2\). In 13891, a vanished outside term is discarded; in 13892, mixing 3 with a chain value creates4. The two exact criteria and their finite witness arguments must therefore remain separate.

\subsection{A bounded-variable finite basis for a nilpotent semiring with a left-zero top}\label{case:17}
This is an unrestricted proof, not an inference from bounded-word computations. All identities are relative to the additively idempotent semiring axioms in the constant-free signature \((+,\cdot)\).

\par\smallskip\noindent\textit{1. The algebra and a finite basis.}\quad 
Let \(F=\{0, a, a^2, a^3\}\), with \(a^i a^j=a^{i+j}\) when \(i+j\leq3\), and zero otherwise. Zero absorbs multiplication. Addition is flat: distinct elements have sum \(0\). Form \(A=F\cup\{t\}\) by placing \(t\) above \(F\) additively and defining

\CaseDisplay{ts=t\quad(s\in A),\qquad st=0\quad(s\in F).}

Thus \(t\) is a multiplicative left zero, not a two-sided multiplicative zero. Associativity follows by considering the first occurrence of \(t\): a word beginning with \(t\) is \(t\), a word with a later \(t\) and a first letter in \(F\) is zero, and otherwise it is evaluated in \(F\). Distributivity inside \(F\) follows from cancellativity of positive integer exponents before truncation; across \(t\), a left multiplication by an element of \(F\) takes \(t\) to the greatest element \(0\) of \(F\), and a left multiplication by \(t\) is constant. Right multiplication takes \(t\) to \(t\) and preserves the flat joins in \(F\).

Let \(X\) be a fixed set of twelve variables, and \(W=\bigcup_{j=1}^4X^j\). Define the explicitly finite set

\CaseDisplay{\Sigma=\{w\leq\sum_{v\in P}v:\ w\in W,\ \varnothing\ne P\subseteq W,\ A\models w\leq\sum_{v\in P}v\}.}

Membership in this finite set is decidable by the \(5^{12}\) assignments to \(X\). This definition is an effective finite basis specification; it is not a claim that the very large set has been explicitly enumerated.

Theorem. \(\Sigma\) is a complete identity basis of \(A\). In fact the proof below uses local absorption rules on at most eight variables, along with the four-variable word reduction. Twelve is a conservative bound. The opposite semiring has the reversed basis.

\par\smallskip\noindent\textit{2. Head protection and word reduction.}\quad 
For a word \(w\), let \(h(w)\) denote its first variable, and let \(H(P)\) be the set of heads of summands of \(P\).

Head lemma. \(A\models w\leq P\) if and only if \(F\models w\leq P\) and \(h(w)\in H(P)\).

Necessity of the first assertion follows by restriction. For the second, set \(h(w)=t\) and all other variables to \(0\). If the head is missing on the right, its value is \(0\), whereas the target is \(t\). Conversely, if a head of \(P\) has value \(t\), the desired inequality is automatic. Otherwise all heads of \(P\), and hence that of \(w\), lie in \(F\). Replacing each assignment of \(t\) by \(0\) leaves every relevant word value unchanged, so the inequality follows from its validity in \(F\). This argument applies to arbitrary variable sets.

The identity \(xyzu\approx x^4\) is valid in \(A\): words of length at least four have value \(t\) precisely when their first variable has value \(t\), and otherwise zero. Both directions belong to \(\Sigma\). By substituting a nonempty suffix for \(u\), every word of length at least four reduces to \(h(w)^4\). Thus every polynomial reduces to a sum of words of lengths one, two, three, and fourth powers. No empty word or constant is introduced.

\par\smallskip\noindent\textit{3. Exact absorption criterion.}\quad 
Fix such a normalized nonempty polynomial \(P\). Let \(L\) be its linear variables; let \(E\) be the undirected multigraph of its length-two summands, keeping the actual orientation of each summand in the polynomial; let \(V\) be the vertices incident with \(E\); and let \(T\) be the union of the variables occurring in its length-three summands. Loops are allowed. A component is marked by its vertices in \(T\). Say that \(P\) is inconsistent if one of the following holds:

1. \(P\) has a fourth-power summand;

2. \(L\cap(V\cup T)\ne\varnothing\);

3. \(T\ne\varnothing\) and \(E\) has a nonbipartite component;

4. \(T\ne\varnothing\), \(E\) is bipartite, and a connected component has marked vertices in both parts.

In the consistent case with \(T\ne\varnothing\), define \(U\) to consist of \(T\) and all vertices in the marked part of each marked component; define \(D\) to consist of the opposite parts of those components. Unmarked components contribute to neither set.

Criterion. A normalized word \(w\) is absorbed by \(P\) in \(A\) exactly when \(h(w)\in H(P)\) and one of these conditions holds:

\textbullet{} \(P\) is inconsistent;

\textbullet{} \(w\) is linear and its variable belongs to \(L\);

\textbullet{} \(T=\varnothing\), \(E\ne\varnothing\), \(w=xy\) with \(x, y\in V\), and either \(E\) has a nonbipartite component, or \(x, y\) are in opposite parts of one connected component;

\textbullet{} \(T\ne\varnothing\), \(w\) has length three, and all its variables belong to \(U\);

\textbullet{} \(T\ne\varnothing\), \(w=xy\), and either one endpoint is in \(U\) and the other in \(D\), or the endpoints are in opposite parts of one connected component of \(E\).

To prove the criterion, consider an evaluation in \(F\) for which \(P\ne0\). All summands must have the same value \(a^d\), \(1\le d\le3\). Every variable occurring in a summand then has value \(a^{e_x}\) for an integer \(1\le e_x\le3\), and the sum of the exponents along each word is \(d\). Variables absent from \(P\) remain unrestricted. Linear words force \(e_x=d\). Quadratic words impose \(e_x+e_y=d\). A cubic word forces \(d=3\) and sets all its variables to exponent one. A fourth power is impossible. The four inconsistent cases are precisely the obstructions to this system.

If there are quadratic words and no cubic words, exponent \(d=2\) is always possible and puts every graph vertex at exponent one. If there is an odd cycle, this is the only possible common exponent. If every component is bipartite, exponent \(d=3\) is also possible, and independently in each component either part may receive exponent one, with exponent two on the other. Linear variables are disjoint and have exponent \(d\). This yields exactly the quadratic case above; a cubic word cannot have exponent two, and a fourth power is zero.

If cubic words occur, \(d=3\). Precisely the variables in \(U\) are forced to exponent one; precisely those in \(D\) are forced to exponent two. The two assignments on each unmarked bipartite component are independent. Thus a cubic word has value \(a^3\) in every solution exactly when all its variables lie in \(U\). A quadratic word has exponent three exactly in the two stated cases. A linear word works exactly when it is in \(L\). In the absence of quadratic and cubic words, assign all linear variables \(a\), which immediately excludes all additional non-linear words. A variable not occurring in \(P\) can be assigned zero to separate a target containing it. These observations also provide separating evaluations whenever a listed condition fails. Finally apply the head lemma.

\par\smallskip\noindent\textit{4. Bounded local deductions.}\quad 
For each \(h\in H(P)\), fix one original summand \(q_h\) with head \(h\), of length at most four. Below \([xy]\) means either actual orientation \(xy\) or \(yx\), separately chosen. Every displayed rule means the finite family obtained by choosing these orientations and writing the guard \(q_h\) in one of the forms \(h\), \(hb\), \(hbc\), \(hbcd\). Subsequent variable identifications are substitutions. All words in these rules have length at most four, and there are at most eight variables.

The following rules belong to \(\Sigma\):

Zero and direct conflicts

\CaseDisplay{h^4\leq x^4+q_h,\qquad h^4\leq x+v+q_h,}

where \(v\) is any length-two or length-three word containing \(x\); and

\CaseDisplay{h^4\leq x^2+uvw+q_h.}

Orientation and odd-path compression

\CaseDisplay{xy\leq yx+q_x,}

\CaseDisplay{xs\leq[xy]+[yz]+[zs]+q_x.}

Consequences of a square

\CaseDisplay{xy\leq z^2+[xu]+[yv]+q_x.}

Unit markers and their propagation

\CaseDisplay{rx^2\leq ruv+abc\quad(x\in\{a, b, c\}),}

\CaseDisplay{rz^2\leq rx^2+[xy]+[yz].}

Incompatible markers

\CaseDisplay{h^4\leq rx^2+ry^2+[xy]+q_h.}

Targets forced by markers

\CaseDisplay{xyz\leq rx^2+ry^2+rz^2+q_x,}

\CaseDisplay{xy\leq rx^2+rz^2+[zy]+q_x,}

\CaseDisplay{yx\leq rx^2+rz^2+[zy]+q_y.}

For a direct verification, restrict first to \(F\). If the sum on the right is zero the inequality holds. Otherwise every listed summand has one common nonzero exponent \(d\). The direct-conflict rules cannot have such a common exponent. For path compression, \(e_x+e_y=e_y+e_z=e_z+e_s=d\) implies \(e_x+e_s=d\). A square forces \(d=2\) and hence exponent one at all incident graph vertices. A cubic marker \(rx^2\) forces \(d=3\) and \(e_r=e_x=1\); the two-edge marker rule then gives \(e_z=e_x=1\). The remaining marker rules follow by adding these forced exponents, or exhibit an edge of exponent two against cubics of exponent three. The appropriate guard, or the displayed cubic with head \(r\), supplies the target head on the right. Hence the head lemma proves every rule in \(A\).

Every use of a rule with its summands already below \(P\) derives the new word below \(P\). Thus the rules may be iterated without adding unproved assumptions.

\par\smallskip\noindent\textit{5. Deriving the criterion in arbitrarily large graphs.}\quad 
An odd-length walk from a vertex \(x\in H(P)\) to \(y\) derives \(xy\leq P\). The first edge is already present or is reversed using its guard. Successive pairs of edges extend the derived edge using the three-edge compression rule. Repeated vertices and loops cause no difficulty. In particular, opposite parts of a connected bipartite component are connected by an odd walk. A nonbipartite component contains a vertex \(z\) which is a head of an original edge; an odd closed walk from \(z\) derives \(z^2\leq P\). This also applies when the odd cycle is reached by a path traversed out and back.

If there are no cubic words and the graph has an odd component, the derived square and any original edges incident with target vertices \(x, y\) derive \(xy\leq P\) by the square rule, provided \(x\in H(P)\). This works even when the vertices lie in different components.

Suppose cubic words occur, and fix one \(ruv\) in \(P\). The initial-marker rule gives \(rx^2\leq P\) for every \(x\in T\). The two-edge rule propagates these markers along every even walk. In a consistent graph this produces a marker for every member of \(U\). Every \(y\in D\) has a neighbor \(z\in U\), so the mixed-marker rule derives every required quadratic word with endpoints in \(U, D\), in either orientation permitted by its head guard. Three unit markers derive every required cubic target. The remaining quadratic targets follow from odd-path compression.

If marked vertices occupy opposite parts of one connected bipartite component, propagate the marker of one along the even-length prefix of an odd path to the other. This makes both endpoints of the last edge marked. The incompatible-marker rule derives \(h^4\leq P\) for every original head \(h\). A nonbipartite component in the presence of a cubic word instead gives a square, and the square-cubic conflict rule gives the same result. The other inconsistent cases are handled by the direct-conflict rules. Finally, for any normalized word \(w\) with head \(h\), the inequality \(w\leq h^4\) is valid in \(A\) and is in \(\Sigma\) (it uses at most three variables, or is equality). Thus inconsistency derives every absorption permitted by the head condition.

This proves derivability for every case of the exact criterion on an arbitrary finite variable set. Every valid polynomial identity is equivalent, by distributivity and idempotence, to the two families of absorptions of each side's words into the other side. Word reduction and the criterion derive all of them. The theorem follows.

\subsection{Finite bases from two endpoint protections and seven witnesses}\label{case:18}
Work in the constant-free signature \((+,\cdot)\), relative to the additively idempotent semiring axioms. No bounded-variable test is used as a completeness argument.

\par\smallskip\noindent\textit{1. Algebra and word values.}\quad 
Let \(A=\{0, a, b, e, t\}\), with additive order \(a, b<0<t\), \(e<t\), and no further comparisons except those forced by transitivity. Its multiplication, in this order, is

\CaseDisplay{\begin{array}{c|ccccc}\cdot&0&a&b&e&t\\\hline
0&0&0&0&0&0\\a&0&0&0&0&0\\b&0&0&0&b&0\\e&0&a&0&e&t\\t&0&0&0&t&t\end{array}.}

This is the operation table of catalogue 252. It is an ai-semiring; the two associativities, commutativity and idempotence of addition, and both distributive laws can be verified directly on the five-element table. The following explicit word formula also proves multiplicative associativity by concatenation.

Write \(K=\{e, t\}\). For any nonempty word, its value is:

\textbullet{} \(e\) if every occurrence is assigned \(e\);

\textbullet{} \(t\) if every occurrence is assigned an element of \(K\) and at least one is assigned \(t\);

\textbullet{} \(a\) if the final occurrence is assigned \(a\) and every earlier occurrence is assigned \(e\);

\textbullet{} \(b\) if the first occurrence is assigned \(b\) and every later occurrence is assigned \(e\);

\textbullet{} \(0\) in every other case.

The description includes words of length one. In particular, the two identities

\CaseDisplay{xyzu\approx xzyu,\qquad xy^2z\approx xyz\quad\textup{(\ensuremath{W})}}

hold. They allow any permutation of interior positions and deletion of repeated interior variables, while retaining the first and last occurrences. An application of the first identity to adjacent interior letters uses a nonempty prefix and suffix, so it introduces no unit constant.

\par\smallskip\noindent\textit{2. Exact absorption criterion.}\quad 
For a polynomial \(P\) (a nonempty finite set of words), define

\CaseDisplay{U=\bigcup_{v\in P}c(v),\quad H=\{h(v):v\in P\},\quad T=\{t(v):v\in P\},}

and let \(R\) be the union of the variables in all positions except the first in the right words; let \(S\) be the analogous union omitting the last position. For a target \(w\), put \(C=c(w)\), \(R_w=c(w\text{ with its first position omitted})\), and \(S_w=c(w\text{ with its last position omitted})\). An omitted word of length zero contributes the empty set, not a term constant.

Criterion. \(A\models w\leq P\) if and only if

\CaseDisplay{C\subseteq U,\qquad (\exists u\in P)\ c(u)\subseteq C,\quad\textup{(\ensuremath{A})}}

\CaseDisplay{H\cap R\ne\varnothing\quad\hbox{or}\quad H\cap R_w=\varnothing,\quad\textup{(\ensuremath{H})}}

\CaseDisplay{T\cap S\ne\varnothing\quad\hbox{or}\quad T\cap S_w=\varnothing.\quad\textup{(\ensuremath{T})}}

For necessity of the first condition, assign a target-only variable zero and all others \(e\): the right side is \(e\), while the target is zero. If the internal witness \(u\) does not exist, assign every target variable \(e\) and every outside variable zero. Then the target is \(e\) and every right word is zero. If (H) fails, assign every variable of \(H\) the value \(b\) and every other variable \(e\). Since \(H\cap R=\varnothing\), every right word has value \(b\); a target occurrence in \(H\cap R_w\) forces its value to zero. The dual assignment using \(a\) proves necessity of (T).

For sufficiency consider the value of \(P\). If it is \(t\), there is nothing to prove. If it is \(e\), every right word is \(e\), so every variable of \(U\) has value \(e\); content inclusion makes the target \(e\). If it is \(0\), every right word belongs to \(\{0, a, b\}\). If the target instead belonged to \(K\), all its variables would lie in \(K\), and the internal witness would also belong to \(K\), which cannot lie below \(0\). Hence the target is in \(\{0, a, b\}\) and is below \(0\).

If \(P=b\), all right words equal \(b\). The word formula forces \(H\) to be assigned \(b\), \(R\) to be assigned \(e\), and \(H\cap R=\varnothing\). By (H) all noninitial target occurrences are assigned \(e\). The internal witness contributes some variable in \(H\cap C\); it can therefore occur only at the target's first position. Thus that first value is \(b\) and the target is \(b\). The case \(P=a\) is symmetric, using (T). These exhaust the five values and prove the criterion for every variable set.

\par\smallskip\noindent\textit{3. A general compression lemma for bounded families.}\quad 
Suppose a semiring satisfies (W), and fix an integer \(m\geq1\). Any valid word absorption with at most \(m\) right summands follows from (W) and the valid word absorptions on at most

\CaseDisplay{B=2(m+1)+2^{m+1}-1}

variables, with at most \(m\) right summands and each word of length at most \(B+2\).

Indeed, the target and the at most \(m\) right words form a family of at most \(m+1\) words. Retain as distinguished variables all their first and last variables; there are at most \(2(m+1)\) of them. Every other variable occurs only at interior positions of every word in which it appears. Using (W), make each variable occur at most once in the interior of each word, then permute those interiors. Partition the nondistinguished variables according to the set of words in which they occur. There are at most \(2^{m+1}-1\) nonempty membership classes. Arrange each interior as its distinguished variables followed by the membership-class blocks in a fixed order.

Identify all variables of each nonempty class with one new variable, keeping distinguished variables separate. Substitution preserves validity. Each resulting class block is a positive power at interior positions and can be reduced to one occurrence by (W). The compressed family has at most \(B\) variables; each word consists of its original first and last variables and at most \(B\) interior variables, so its length is at most \(B+2\). A one-letter word remains one-letter.

Conversely, substitute the product of the original class members, in their fixed order, for the corresponding new variable. This is always a nonempty word. The expanded family is precisely the arranged family above, up to interior permutations and removal of interior repetitions. Thus the compressed absorption derives the original one. This proves the lemma without claiming that arbitrary variable identification can be inverted; the reverse step is the specified nonempty block substitution.

\par\smallskip\noindent\textit{4. Decomposing every absorption into at most seven witnesses.}\quad 
Let \(w\leq P\) satisfy the criterion, and let \(h, t\) be its endpoint variables, possibly equal. Choose an internal witness \(u\in P\) with \(c(u)\subseteq C\). If the first alternative in (H) holds, choose two right words witnessing it: one has a head \(x\), the other has a noninitial occurrence of \(x\). Likewise, if the first alternative in (T) holds, choose at most two right words witnessing that tail conflict. Let \(G\) be the union of these witnesses; \(|G|\leq4\). If one word already witnesses both roles, only one is needed.

By content inclusion choose right words containing \(h\) and \(t\). Put \(C_0=c(u)\cup\{h, t\}\), and let \(w_0\) be the word obtained from \(w\) by retaining exactly the occurrences whose variables lie in \(C_0\). Since both endpoint variables are retained, its endpoints are \(h, t\), and its noninitial and nonfinal variable sets are subsets of those of \(w\). Let \(P_0\) contain \(u\), the members of \(G\), and the two endpoint-content witnesses. Then \(|P_0|\leq7\).

The criterion proves \(w_0\leq P_0\). The two content conditions follow from construction. If a global head conflict was present, its witnesses are retained. If no global head conflict was present, the heads of \(P_0\) are a subset of \(H(P)\) and avoid the target's noninitial variables, hence avoid those of \(w_0\). The tail argument is identical. Therefore \(w_0\leq P\) is derivable from valid absorptions with at most seven summands.

Add the remaining variables of \(C\) one at a time. Suppose \(C_i\) has been obtained, and write \(w_i\) for the restriction of \(w\) to it. For a new variable \(y\), choose an original right word \(r_y\) containing it, and set

\CaseDisplay{C_{i+1}=C_i\cup\{y\},\qquad P_i'=\{w_i, r_y\}\cup G.}

There are at most six words in \(P_i'\), and the content conditions for \(w_{i+1}\leq P_i'\) hold, with \(w_i\) as internal witness. Global conflicts remain witnessed by \(G\). If there was no global head conflict, the head of \(u\) lies in \(H(P)\cap C\) and must be the first variable \(h\) of \(w\), since \(H(P)\) avoids all noninitial target occurrences. In particular \(h\in H(P)\). Thus the new word \(w_i\) adds no new head outside \(H(P)\), and its right-family heads still avoid the noninitial variables of \(w_{i+1}\). The tail argument is symmetric. The criterion therefore validates each step.

Every member of \(G\) and every \(r_y\) was in the original \(P\). Inductively \(w_i\leq P\), so transitivity gives \(w_{i+1}\leq P\). At the final step the restriction is the original word \(w\). The reasoning also covers a one-letter target, when \(C_0=C\) and no extension is required.

\par\smallskip\noindent\textit{5. Finite basis theorem.}\quad 
Take \(m=7\), so \(B=271\) and \(B+2=273\). On a fixed set of 271 variables let \(\Sigma\) be all valid absorptions

\CaseDisplay{v\leq\sum_{q\in Q}q,\qquad 1\leq |Q|\leq7,}

whose words all have length at most 273, together with (W). This is a finite effectively specified set: both the set of words and the set of possible subsets \(Q\) are finite, and validity is decidable by all assignments in the five-element algebra. No enumeration of the entire huge set, and no minimality, is claimed.

Theorem. \(\Sigma\) is a complete finite basis for \(A\).

Every valid word absorption is a chain of absorptions with at most seven right summands by Section 4. Each member of that chain follows from \(\Sigma\) by the compression lemma. Distribution and the absorptions of every summand in both directions yield every polynomial identity. All arguments allow arbitrary numbers of variables and summands and arbitrary word lengths. The opposite also has a finite basis by word reversal.

This method is distinct from claiming that the final absorption itself always has seven original right witnesses: full content coverage may require unboundedly many right words. The proof instead uses a sequence of bounded-witness deductions, adding one variable at a time.

\subsection{A six-witness finite basis from prefix and tail coverage}\label{case:19}
Source 462 in the five-element catalogue. The signature is (+, multiplication), with no named constants. All bases below are relative to the additively idempotent semiring axioms.

\par\smallskip\noindent\textit{1. Word values and normalization.}\quad 
Let A=\{0,\allowbreak{}1,\allowbreak{}2,\allowbreak{}3,\allowbreak{}4\}, with tables in that order:

{\small\par\smallskip
\noindent addition       multiplication\par
\CaseLine{\ttfamily 0 0 4 4 4\qquad 0 0 0 0 0}
\CaseLine{\ttfamily 0 1 4 4 4\qquad 0 0 0 0 0}
\CaseLine{\ttfamily 4 4 2 2 4\qquad 0 0 4 2 4}
\CaseLine{\ttfamily 4 4 2 3 4\qquad 0 1 4 3 4}
\CaseLine{\ttfamily 4 4 4 4 4\qquad 0 0 4 4 4}
\smallskip}

The additive order consists of \(1<0<4\) and \(3<2<4\), with elements of the two branches incomparable except at 4. Set B=\{0,\allowbreak{}1\} and K=\{2,\allowbreak{}3,\allowbreak{}4\}. For any word w of length at least two, its value is:

\textbullet{} 0 if a nonfinal position lies in B or the final position is 0;

\textbullet{} 1 if its final position is 1 and every nonfinal position is 3; otherwise 0 when its final position is 1;

\textbullet{} its first value if all positions lie in K and every position after the first is 3; otherwise 4 when all positions lie in K.

In particular, a nonlinear word has value in K exactly when all its variables take values in K. The same assertion holds for linear words.

Here is an associativity check. A first or middle factor in B makes both associations of a triple 0. With the first two factors in K, a last factor 0 gives 0; a last factor 1 gives 1 exactly when both first factors are 3, and gives 0 otherwise. On K, uv equals u if \(v=3\) and equals 4 otherwise, which is associative. Each row and column of the displayed multiplication table preserves joins: one checks the two chains and the joins between their branches. Therefore both distributive laws hold. The word formula follows by induction.

For a word w let c(w) be its content, h(w) its first variable, t(w) its last variable, i(w) the variables before the final position, and r(w) the variables after the first position. For a linear word both i(w) and r(w) are empty. The word formula proves

\CaseDisplay{xyzt\approx xzyt,\qquad xy^2z\approx xyz.       \quad\textup{(\ensuremath{W})}}

Fix a variable order. These laws sort the interior positions and remove repeated interior occurrences, leaving the first and last positions intact. Thus each nonlinear word becomes h (product of M) t, where M is its interior set. Endpoints may belong to M and may equal each other. A missing interior is an omitted substring, not a unit constant. Linear words stay linear. This is the normal form used below.

\par\smallskip\noindent\textit{2. Exact semantic criterion.}\quad 
For a nonempty polynomial P, put

\CaseDisplay{D=\bigcup_{v\in P}c(v),\quad E=\bigcup_{v\in P}r(v),\quad
 J=\bigcup_{v\in P}i(v),\quad R=\{t(v):v\in P\}.}

In particular \(D=J\) union R; linear variables contribute to R, but not to J or E.

Lemma 1. For every nonempty word w, linear or nonlinear, the absorption \(w\leq P\) is valid in A exactly when all four conditions hold:

1. c(w) is contained in D;

2. r(w) is contained in E;

3. some summand v of P has c(v) contained in c(w);

4. either J meets R, or i(w) is contained in J and t(w) belongs to R.

Proof. If condition 1 fails, give a missing variable the value 4 and all others 3. The right side is 3 and the target is 4. If condition 2 fails, assign a variable in r(w) outside E the value 2, and all others 3. Every right summand is 2 or 3, while \(w=4\). If condition 3 fails, assign c(w) the value 3 and every outside variable 0. Every right summand contains an outside variable and is 0, whereas \(w=3\). These are separating evaluations.

For condition 4 assume condition 1 and J disjoint from R. Assign J the value 3 and R the value 1; outside variables may be assigned 3. Every right summand is 1. Since c(w) is contained in \(D=J\) union R, the target is 1 exactly when all its nonfinal variables belong to J and its final variable belongs to R. This statement also covers a linear target. Failure gives 0,\allowbreak{}3, or 4, none below 1.

Conversely fix an assignment. If \(P=4\) there is nothing to prove. If \(P=0\), all summands lie in B. The witness v of condition 3 has some variable valued in B, so c(w) also contains such a variable; the word-value formula places w in B and thus below 0. If \(P=1\), every summand equals 1. Therefore J is assigned 3 and R is assigned 1, forcing them to be disjoint; condition 4 then gives \(w=1\). If \(P=2\), all variables in D are assigned in \{2,\allowbreak{}3\}, and all variables in E are assigned 3. Conditions 1 and 2 give w in \{2,\allowbreak{}3\}, hence \(w\leq 2\). Finally \(P=3\) forces all variables of D to be 3, so condition 1 gives \(w=3\). All five possible values of P have been considered. \ensuremath{\square}

\par\smallskip\noindent\textit{3. Six-witness decomposition.}\quad 
Assume \(w\leq P\) is valid and normalize all words. If w already occurs in P, no step is needed. Choose a condition-3 witness v in P.

For a linear target y, v is a positive power of y. The absorption \(y\leq v\) is valid by Lemma 1: its content is \{y\}, and either v is linear or its prefix set meets its final set. It is a step with one right summand, and derives \(y\leq P\).

For nonlinear w write \(w=h\) (product of M) t and put \(M_{0}=c(v)\) minus \{h, t\}. Then \(M_{0}\) is contained in M. Let \(s_{0}=h\) (product of \(M_{0}\)) t. Choose a right summand \(r_{t}\) with t in its tail, which exists by condition 2. If the head g of v lies in \(M_{0}\), also choose a right summand \(r_{g}\) with g in its tail. This is the only possible variable of \(M_{0}\) not already covered by r(v). When g is not in \(M_{0}\), omit \(r_{g}\) in what follows.

Suppose first that J meets R. Choose a with some p in i(a), and b with \(t(b)=p\). Choose \(r_{h}\) containing h. Then

\CaseDisplay{s_0\leq v+a+b+r_h+r_t+r_g               \quad\textup{(\ensuremath{2})}}

is valid, with at most six distinct right summands. Its content is covered by v, \(r_{h}\), \(r_{t}\). Its tail is covered by r(v), r(\(r_{t}\)), r(\(r_{g}\)). The summand v is a content witness, and a, b preserve the crossing of J and R. This checks every condition of Lemma 1.

If x is a missing interior variable, choose an original summand \(r_{x}\) with x in its tail. Writing s' for the word obtained from the current s by adding x to its interior, we have

\CaseDisplay{s'\leq s+a+b+r_x.                 \quad\textup{(\ensuremath{3})}}

The old word s is the content witness; \(r_{x}\) supplies the new content and tail membership, and a, b retain the crossing. This also applies when the new interior variable equals an endpoint.

Now suppose J is disjoint from R. Condition 4 gives i(w) contained in J and t in R, so t differs from h and does not belong to M. The chosen content witness v ends in t: its final variable belongs to R intersect \(c(w)=\{t\}\). If v is linear it is precisely t. Choose \(r_{h}\) with h in its prefix. Then

\CaseDisplay{s_0\leq v+r_h+r_t+r_g             \quad\textup{(\ensuremath{4})}}

is valid. Every variable in \(M_{0}\) lies in the prefix of v, h lies in the prefix of \(r_{h}\), and t is a final variable of v. Thus condition 4 holds for this smaller right side, while content, tail coverage, and the content witness are checked just as for (2). For a missing interior x, choose original \(r_{x}\), \(q_{x}\) with x in the tail of \(r_{x}\) and the prefix of \(q_{x}\). Both exist by conditions 2 and 4. Then

\CaseDisplay{s'\leq s+r_x+q_x               \quad\textup{(\ensuremath{5})}}

is valid: s supplies the old content, tail and prefix memberships, \(r_{x}\) supplies the new tail membership, \(q_{x}\) the new prefix membership, and the final variable t remains final in s. Condition 3 uses s.

After finitely many extensions the left word is w. Every right word is either an original summand of P or the previous constructed word. Hence all absorptions follow if the valid steps with at most six normal right summands can be derived. This construction is unrestricted in variable count and word length.

\par\smallskip\noindent\textit{4. Reversible compression and finite basis.}\quad 
A step has at most seven words in total. Reserve all their endpoint variables individually; there are at most 14. Every other variable occurs only in interiors. Group those remaining variables according to the nonempty pattern of interiors in which they occur. There are at most \(2^{7}-1=127\) patterns. Identifying each group to one variable and normalizing by (W) gives a valid absorption on at most 141 variables, with every normal word of length at most 143.

This identification is reversible modulo (W): replace a group variable by the nonempty product of all its members. Every member had exactly the same interior incidence pattern, so the expansion restores the original interior set in every word. Reserved endpoints expand to their individual original variables. Sorting interiors and deleting repeated interior occurrences by (W) recovers the original normal words. Linear words consist of a reserved variable and are recovered as well. There are no empty substitutions or constant symbols.

Fix a 141-letter alphabet. Define Sigma as (W) together with all absorptions valid in A having a normal left word, at most six normal right words, and all word lengths at most 143 over this alphabet. This is a finite, explicitly specified set. Lemma 1 decides membership. The large set is not enumerated here; no small or minimal basis is claimed.

Theorem. Sigma is a finite identity basis for A. Its word reversals form a finite basis for the multiplicative opposite.

Proof. Every member is valid. By the reversible compression argument, each valid absorption with at most six right summands follows from Sigma by substitution and (W). Section 3 then derives every valid word absorption into a polynomial, one step at a time. Distributivity expresses every term as a sum of words. A valid identity between two such polynomials says that each summand on either side is absorbed by the other side, and all these absorptions have just been derived. This proves unrestricted completeness. Word reversal transports the basis and derivations to the opposite. \ensuremath{\square}

Finite calculations supporting this proof are only error checks. Classification transfers require exact equality of generated varieties, verified in both HSP directions, not just short-word agreement.

\subsection{Finite basis for 248, 11186 and 11366}\label{case:20}
For type 248 the word value is 4 when all letters have value 4, and 3 when all have value 3. It is 1 when the last value is 1 with every preceding value 3, and 2 when the first value is 2 with every following value 3. All other words have value 0. This includes literals, using empty prefix and suffix conditions only in this description. The multiplication table proves the formula by appending one letter: 0 and 1 multiply on the right to 0, 2 survives only on 3, 3 sends 1 to 1 and 3 to 3 and all other values to 0, and 4 survives only on 4.

For a word w let C(w) be its content, B(w) its set of noninitial letters and A(w) its set of nonfinal letters. For a polynomial P use the unions of these sets over its original words. For the level below 0, the target can fail only when all its variables are 4. Thus w \ensuremath{\leq} P at this level exactly when some original u has C(u) contained in C(w): give 4 to all target variables and 3 to the others for necessity. At the level below 3, the upper polynomial is below 3 exactly when all its variables are 3. This gives the condition C(w) contained in C(P). At the level below 1, all upper variables must belong to \{3,\allowbreak{}1\} and all nonfinal variables must be 3. In addition to content coverage this gives A(w) contained in A(P). If a required nonfinal variable is missing from A(P), give it 1 and everything else 3. At the level below 2 the symmetric test gives B(w) contained in B(P), using 2 at a missing noninitial variable. The level below 4 imposes no condition. These four conditions are necessary and sufficient, including for repeated letters and literals.

Use the two structural laws xyzt = xzyt and xyyz = xyz. They preserve the endpoints and permit arbitrary permutations and deletion of extra copies within the interior while fixing any repeated occurrences in the endpoints. The table formula verifies both. The following finite family together with the structural laws and semiring axioms is a complete basis: all valid absorptions on \(D=16\),411 variables of lengths at most D+2 and at most thirteen right words.

To prove this, for a valid w \ensuremath{\leq} P choose a contained original word u. If w is a literal, u is a nonempty power of that literal and w \ensuremath{\leq} u follows from the criterion. Otherwise write w with two endpoints h, t and an interior list, possibly empty. For a set S containing h, t keep the two endpoint positions and precisely the interior letters in S to form \(w_{S}\). Start with \(S=C(u)\) union \{h, t\}. Every member of S except possibly h, t and the two endpoint variables of u occurs internally in u. For each of these at most four exceptional variables choose from P one original witness for each of the content, nonfinal and noninitial requirements actually made by \(w_{S}\). Validity of w guarantees all these choices. The criterion proves \(w_{S}\) \ensuremath{\leq} u plus the chosen words, using at most thirteen right summands. Leaving \(h=t\) gives two actual positions, even when u is a literal.

Adjoin each missing target variable x to S in turn. It occurs internally in w. Retain the already derived \(w_{S}\) and one original noninitial and one original nonfinal witness for x. The old roles persist, the two witnesses supply all requirements at x, and \(w_{S}\) itself is a contained word for the enlarged target. Thus this step is a valid absorption with at most three upper words. The induction ends at w.

For each local absorption with at most thirteen right words, normalize the interiors. Reserve every endpoint variable across the at most fourteen words, at most 28 variables. Partition the remaining variables by the exact nonempty set of words in which they occur, at most \(2^{14}\)-1 classes. Identify each class with its own fresh variable and normalize the interiors. The result is a valid absorption in the displayed finite family. Substitute the nonempty product of all members for each class variable. In every participating word this recovers each member exactly once internally, and the structural laws restore each original word. Reserved variables are fixed throughout, including their internal appearances. This proves derivability of every step for arbitrary many variables. Restore the omitted summands and treat each monomial of a polynomial identity in both directions.

Applying Proposition~\ref{prop:lower-fibres} to the two other labelled algebras and \(\A{248}\), in both directions, gives the same equational theory. Thus the bounds specify a finite basis for all three algebras.

\subsection{Twenty-six further applications of the bounded constraint-tree theorem}\label{case:21}
Apply Theorem~\ref{thm:treelength} with \(q=5\) and \(d=4\). A five-letter identity reduces each input word to words of length at most four. Each selected map preserves the four principal ideal relations and separates the specified target tuples. Minimum in the witnessed order preserves these relations. The folds to the diagonal handle each occurrence of the deleted pair in each relation of the induced square.

Here one may use \(L=40\), \(M=5\), \(T=5^{90}\) and \(B=5^{95}\). Take the semiring axioms, the individual five-letter reduction and all valid absorptions on 4(B+1) variables with word lengths at most four and at most B upper words. These give a finite basis for arbitrary variables and polynomials. All selected maps have image of size at most three. The uniform bounds suffice without shortening the interpolation paths in the general theorem.

To test a fold, take the coordinatewise products on the two sides of the occurrence. If the original tuple has length at most \(d\), the two contexts have total length at most \(d-1\). For \(d=4\) one can equivalently list all tuples of length at most four and replace each occurrence in turn. An empty context omits multiplication on that side.

The new finite bases have IDs 262, 705, 2489, 2498, 11039, 11087, 11264, 11348, 11405, 11414, 11497, 11499, 11526, 11529, 12363, 12409, 12596, 12597, 12606, 12610, 12818, 12827, 12960, 12976, 13489 and 13581.

\subsection{Nine complete identities for semiring 279}\label{case:22}
In the signature (+, multiplication), without constants, the following form a complete basis relative to the semiring axioms:

{\small\par\smallskip
\CaseLine{\(xy=yx; xyz=xyzz; xyzt\leq xyz; xxx\leq x+xy;\)}
\CaseLine{\(xxx\leq x+yyy; xxy\leq xy+zzz; xy\leq xxy; x\leq xxx;\)}
\CaseLine{\(xz\leq xx+xy+yz.\)}
\smallskip}

The order notation \(u\leq P\) means \(u+P=P\).

\par\smallskip\noindent\textit{Values and exact criterion.}\quad 
Addition is 00004/\allowbreak{}01004/\allowbreak{}00204/\allowbreak{}00034/\allowbreak{}44444: the distinct elements 1,\allowbreak{}2,\allowbreak{}3 have join 0, with \(0<4\). Multiplication is 00000/\allowbreak{}00000/\allowbreak{}00100/\allowbreak{}00010/\allowbreak{}00004. A nonlinear word has value 4 exactly when all its variables are 4; a quadratic additionally has value 1 if both occurrences are 2 or both are 3. Every other nonlinear word has value 0. Long words thus retain their content.

For the original polynomial P let L be its literal variables, E the undirected quadratic edges, V their endpoints, and S the loop vertices. Let G mean that P has an original long word or L intersects V. Let C(w) mean that some original word has content contained in c(w). The criterion for \(w\leq P\) is:

\textbullet{} Literal x: C(x) and (x is original literal or G).

\textbullet{} Long w: C(w) and G.

\textbullet{} Quadratic ab: C(ab) and (G or a, b lie in the same E-component).

Components contain only original edge endpoints. An isolated variable outside V is not a component vertex. Connectivity here is undirected; it is not any earlier directed-word-path criterion.

If C(w) fails, assign target content 4 and outsiders 1. The target is 4, while every right summand is at most 0. If G fails, assign all V variables 2 and outsiders 1. All right quadratics and literals have value 1. Long targets have value 0. A nonoriginal literal satisfying C must have an original square and hence lie in V, so its value 2 is not below 1. For a quadratic with an endpoint outside V, the same assignment makes its value 0. If both endpoints are in different components, assign the first component 2, all other components3 and outsiders 1. Every right edge is 1 but the target is 0. These separate every failed condition for arbitrary variable and summand counts.

For sufficiency, an original long word forces \(P\geq 0\). An original literal x incident to xy does likewise: if \(xy=1\), then \(x=2\) or 3, and \(x+1=0\); otherwise xy is 0 or 4. The contained witness provides 4 when the target is 4, while every other target value is below 0. Suppose G fails and the target endpoints share a component. If P does not dominate 0, every original edge must have value 1; its endpoints then both belong to\{2,\allowbreak{}3\} and have equal values. Equality propagates through the original component, so the target is 1. The contained witness supplies 4 in the all 4 case. Original literals are immediate.

\par\smallskip\noindent\textit{Finite derivation without size bounds.}\quad 
The two word laws normalize long words by content: commute a chosen occurrence to the end, duplicate it by \(xyz=xyzz\) using two nonempty prefix blocks, then remove surplus occurrences by the reverse law while retaining length at least three. In particular a long g satisfies \(g^{3}=g\). Canonical long forms are sorted content for at least three variables, xxy for two and xxx for one.

When G holds, retain an original long g, or derive \(g=x^{3}\) from an original overlap x, xy by \(xxx\leq x+xy\). Choose an original contained u. Upgrade a literal u to \(u^{3}\) by \(xxx\leq x+yyy\) using g for y; upgrade a quadratic ab to aab by \(xxy\leq xy+zzz\) using g for z. Retain long u. Append missing target content by \(xyzt\leq xyz\), then normalize. Recover a quadratic by \(xy\leq xxy\) or a literal by \(x\leq xxx\). The global long witness and the contained witness need not be the same word.

Without G, a new quadratic ab has a contained original word. It cannot be long; it cannot be an endpoint literal either, because both endpoints belong to V and that would make G true. Thus an original endpoint loop exists, say aa. Take a simple original edge path from a to b. The first edge is original. With aa retained at every stage, apply \(xz\leq xx+xy+yz\) to extend the derived edge from a by one path edge. This derives ab after finitely many steps. No step creates an arbitrary loop or abandons the original loop guard. All substitutions are nonempty.

Every law is valid by the exact criterion and checked independently on its finite assignments. The argument applies without bounds on variables, lengths, path lengths or polynomial sizes. Distributive expansion and summandwise absorption in both directions prove all valid identities. Any other catalogue type requires verified full HSP equivalence to transfer.

\subsection{Nine complete identities for semiring 280}\label{case:23}
Constant-free signature (+, multiplication), relative to the additively idempotent semiring axioms. A complete finite basis is

{\small\par\smallskip
\CaseLine{\(xy=yx; xyz=xyzz; xyzt\leq xyz; xx\leq x+xy;\)}
\CaseLine{\(xxx\leq x+yyy; xxy\leq xy+zzz; xy\leq xxy; x\leq xx;\)}
\CaseLine{\(xz\leq xx+xy+yz.\)}
\smallskip}

The abbreviation \(u\leq P\) means \(u+P=P\).

\par\smallskip\noindent\textit{Exact values and absorption criterion.}\quad 
Addition is 00004/\allowbreak{}01114/\allowbreak{}01214/\allowbreak{}01134/\allowbreak{}44444, where 2,\allowbreak{}3 are incomparable below \(1<0<4\). Multiplication is 00000/\allowbreak{}00000/\allowbreak{}00100/\allowbreak{}00010/\allowbreak{}00004. A nonlinear word is 4 if all its content is assigned 4. Otherwise a quadratic is 1 when both occurrences are 2 or both are 3; all other nonlinear words are 0. Long words retain their content. A literal x is always below \(x^{k}\).

For original P let E be its unordered quadratic edges, V their endpoint union, L its literal variables, and S its loops. Let G mean an original long word exists; let C(w) mean some original word has content contained in the target. Then \(w\leq P\) holds exactly when:

\textbullet{} Literal: C(w).

\textbullet{} Long: C(w) and G.

\textbullet{} Quadratic ab: C(w) and (G or a, b are in one original E-component).

Only edge endpoints are graph vertices. Original literal-edge overlap does not give G:\(xxx\leq x+xy\) fails at \(x=y=2\) with values 0 and 1.

For necessity of C, assign target content 4 and every outsider 3; every right summand is below 4 while the target is 4. When G fails, assign V to 2 and outsiders 1. Right edges and literals are at most 1; every long target and any quadratic with an endpoint outside V is 0. If both target endpoints lie in different original components, assign one component 2, all other components3, and outsiders 1. Every right edge is 1 and every literal is below 1, whereas the target is 0. These give explicit separators for all failed conditions, with no bound on variables.

For sufficiency, a contained pure power dominates a literal target. G forces \(P\geq 0\); the contained witness supplies 4 when a nonlinear target is 4 and all remaining target values are below 0. If G fails and P does not dominate 0, every original quadratic must equal 1. Its endpoints then have the same value in\{2,\allowbreak{}3\}. Equality propagates along each original component, so a target with connected endpoints is 1. The contained witness again handles the all 4 assignment.

\par\smallskip\noindent\textit{Unrestricted complete finite derivation.}\quad 
Commutativity and \(xyz=xyzz\) normalize long words by content. Move a chosen occurrence to the end and duplicate it using two nonempty prefix blocks; reverse the law to remove surplus copies while keeping length at least three. Canonical forms are sorted content when its size is at least three, xxy for two, and xxx for one. Hence long g satisfies \(g^{3}=g\).

A literal target has an original pure power. Normalize a long power to xxx, recover xx using \(xy\leq xxy\) with \(x=y\), then use \(x\leq xx\). An original square directly supplies the last step; an original literal needs none. These operations do not create long words.

For a nonlinear target with G, fix original long g and original contained u. Upgrade literal u to \(u^{3}\) using \(xxx\leq x+yyy\) with \(y=g\); upgrade quadratic ab to aab using \(xxy\leq xy+zzz\) with \(z=g\). Retain long u. Append missing target content using \(xyzt\leq xyz\), normalize, and recover a quadratic ab by \(xy\leq xxy\) if necessary. All blocks are nonempty.

Without G, a new quadratic ab has a contained original literal or loop at a target endpoint, say a, and its endpoints lie in one original component. If only literal a is available, choose an original edge ac and derive aa with \(xx\leq x+xy\). If the target is aa this finishes. Otherwise choose a simple original edge path from a to b. Keep aa throughout and repeatedly use \(xz\leq xx+xy+yz\) to extend the derived edge from a by one original path edge. This derives ab. Original target edges need no step. The derived square is justified by both the original literal and its original incidence, never by a literal alone.

Each displayed law is valid by the exact criterion and independently checked on every finite assignment. The construction permits arbitrary variables, words, path lengths and summands. Distributivity and summandwise absorption in both directions establish completeness. Further catalogue types require independent full HSP equivalence.

\subsection{Finite bases from connectivity and pinned components}\label{case:24}
Work in the constant-free semiring signature, relative to the additively idempotent semiring axioms. The proof applies to arbitrary finite sets of variables.

\par\smallskip\noindent\textit{1. Two cores and their one-sided top extensions.}\quad 
Let \(F=\{0, p, a, b\}\) have multiplication \(a^2=b^2=p\), with every other product zero. There are two choices of addition. In both, \(0\) is the greatest element; in type I, \(a, b<p<0\) and \(a+b=p\); in type II, \(a<p<0\) and \(b<0\), with \(b\) incomparable with \(a, p\). Addition is the join. Extend either core by an additive greatest element \(t\), putting \(ts=t\) for all \(s\) and \(st=0\) for \(s\in F\). Denote the resulting semirings by \(A_{\mathrm I},A_{\mathrm{II}}\).

The core is three-nilpotent and associative. To see distributivity, multiplication by \(a\) or \(b\) has only one nonzero value, namely \(p\) at that same atom. A join of two different atoms, or of an atom with \(p\), belongs to \(\{p,0\}\) and has product zero; the corresponding two products also have sum zero. Comparable pairs preserve joins directly. The zero and \(p\) translations are constant zero. The top extension is associative and distributive by the same direct cases.

The head lemma for source 209 applies here: for either extension,

\CaseDisplay{A\models w\leq P\iff F\models w\leq P\text{ and }h(w)\in H(P).}

Indeed necessity of the head condition follows by sending the missing head to \(t\) and all other variables to zero. If no right head is \(t\), replacing all occurrences of \(t\) in an assignment by zero preserves the relevant word values. Also \(xyz\approx x^3\), so every word of length at least three reduces to the cube of its first variable.

Let \(W\) be all nonempty words of length at most three on twelve fixed variables. Let \(\Sigma_A\) be the finite set of all valid absorptions \(w\leq\sum P\) in \(A\) with \(w\in W\) and \(\varnothing\ne P\subseteq W\). This is a finite effective specification, not an enumeration of the entire set.

Theorem. \(\Sigma_A\) is a complete basis for either semiring. Consequently their opposite semirings are also finitely based.

\par\smallskip\noindent\textit{2. Absorption criterion.}\quad 
Normalize a nonempty polynomial \(P\) to linear words, quadratic words and cubes. Let \(L\) be its linear variables and let \(E\) be the undirected graph of quadratic words, with loops allowed and original orientations retained in \(P\). Let \(V\) be the incident vertices. In type II a component is called pinned if it meets \(L\).

In either type an absorption needs \(h(w)\in H(P)\). Subject to this condition:

1. If \(P\) contains a cube, it absorbs every normalized word.

2. Otherwise, a linear target is absorbed in type I whenever its head condition holds. In type II it is absorbed exactly when it belongs to \(L\) or to a pinned component.

3. Otherwise, a quadratic target \(xy\) is absorbed in type I exactly when \(x, y\in V\) belong to the same component. In type II it is absorbed exactly when both vertices lie in the same component, or both lie in pinned components.

4. Without a cube on the right, no cube is absorbed.

Here and below a repeated endpoint \(x=y\) is allowed. To prove the criterion, restrict to \(F\). If any quadratic summand is nonzero it is \(p\), and this happens exactly when its endpoints have the same value in \(\{a, b\}\). For the whole polynomial to have value different from zero, every component must therefore be monochromatic, independently choosing \(a\) or \(b\).

In type I each linear variable can independently take any value below \(p\), except that incident variables are already constrained by their component. For a linear target whose head occurs in a quadratic word, its value is an atom below \(p\). If \(P\) contains no quadratic words, the head condition already puts a linear target in \(L\). For a quadratic target, vertices in different components can receive different atoms. A target vertex absent from \(V\) can be assigned \(p\), so its product vanishes while \(P\leq p\). Thus connectivity is exactly the quadratic condition.

In type II, linear summands must lie in \(\{a, p\}\) whenever a quadratic summand is present and \(P\ne0\). Consequently each component meeting \(L\) is forced to value \(a\), whereas the other components choose \(a\) or \(b\) independently. This gives the stated pinned-component tests. A linear target in an unpinned component can take \(b\not\leq p\); a quadratic target joining distinct components, at least one unpinned, can have distinct endpoint values and hence zero. An endpoint absent from \(V\) can again be set to \(p\). If only linear words occur, assigning all of them \(p\) excludes every quadratic or cubic target, and the head condition accounts for the linear targets. Assigning every variable \(a\) makes every cube zero but makes any polynomial without cubes strictly below zero, so no such polynomial absorbs a cube.

If a right cube occurs, the core polynomial is identically zero, the greatest core element; hence the head lemma proves the first case. These evaluations prove both directions for every number of variables.

\par\smallskip\noindent\textit{3. Local finite deductions.}\quad 
For each original head \(h\) fix a right summand \(q_h\) of length at most three with that head. As before, \([xy]\) denotes either actual orientation. All the following are finite families of members of \(\Sigma_A\), by taking \(q_h=h\), \(hB\), or \(hBC\) and then applying arbitrary variable substitutions:

\CaseDisplay{h^3\leq x^3+q_h,\qquad w\leq h(w)^3\quad(1\leq|w|\leq2),}

\CaseDisplay{xy\leq yx+q_x,\qquad x^2\leq[xy]+q_x,}

\CaseDisplay{xz\leq[xy]+[yz]+q_x.}

In type I add

\CaseDisplay{x\leq xy.}

In type II add

\CaseDisplay{y\leq r+[ry]+q_y,}

\CaseDisplay{xy\leq r+[rx]+s+[sy]+q_x.}

These local rules use at most seven variables and words of length at most three. Their validity can be checked without computation. On the core, a nonzero quadratic polynomial has value \(p\) and imposes equality of atom values along its edges. This proves the square, orientation and two-edge compression rules. Type I has both atoms below \(p\), proving \(x\leq xy\), including the cases \(x=p\) and \(x=0\) where \(xy=0\). In type II a linear summand \(r\) together with the nonzero edge \(ry\) forces \(r=y=a\); two such pins force the endpoints of the last target to be \(a\). If a right side evaluates to zero the inequality is automatic. Finally the head guard supplies precisely the necessary head condition for the extension.

Starting at any vertex \(x\in H(P)\), an original incident edge is usable in the desired orientation by the orientation rule. Iterating two-edge compression along a path derives \(xy\leq P\) for every other vertex \(y\) of its component. A square at \(x\) follows from any incident edge, including a loop. Thus every target whose endpoints lie in one component is derivable, with no bound on path length.

In type I each original head of a quadratic word satisfies \(x\leq xy\leq P\); original linear heads are already present. This derives every required linear target.

In type II, if \(y\) belongs to a pinned component, choose a linear vertex \(r\) in it. The vertex \(r\) is an original head, so the path argument gives \(ry\leq P\). The linear-pin rule derives \(y\leq P\) whenever \(y\in H(P)\). For a required quadratic target between two pinned components, choose linear vertices \(r, s\) in the two components; derive \(rx, sy\leq P\) by paths, and apply the two-pin rule with the guard for \(x\). This works even when a path has length zero, because the square rule supplies the relevant loop.

If \(P\) has a cube, the zero-transfer rule gives \(h^3\leq P\) for every original head \(h\), and head domination gives every required shorter target. Every case of the criterion is therefore derivable from the finite set. Distribution and two-way word absorptions prove completeness for arbitrary identities.

\par\smallskip\noindent\textit{4. Exact catalogue sources.}\quad 
Type I is catalogue 282 under \((0, p, a, b, t)=(0,\allowbreak{}1,\allowbreak{}2,\allowbreak{}3,\allowbreak{}4)\). Type II is 1707 under \((0, p, a, b, t)=(0,\allowbreak{}2,\allowbreak{}3,\allowbreak{}1,\allowbreak{}4)\). Their opposites are respectively 2348 and 2821. Type I and its opposite account for \(282,2348,13839,14398\). Type II and its opposite account for \(1707,2821,13198,14388\). The variety identifications follow by Proposition~\ref{prop:lower-fibres}.

\subsection{Connected components and cubic marks: two finite identity bases}\label{case:25}
We work in the signature \((+,\cdot)\) without constants, relative to the additively idempotent semiring axioms. Write \(u\leq v\) for \(u+v\approx v\). All completeness arguments below apply to arbitrary finite sets of variables.

\par\smallskip\noindent\textit{1. Two explicit algebras.}\quad 
Both algebras have underlying set \(\{e, a, b, z, t\}\) and the same commutative multiplication: \(e^2=e\), \(a^2=b\), and every other product is \(z\). In \(A\), addition is join in the order \(e, a<b<z<t\), with \(e, a\) incomparable. In \(B\), addition is join in the order \(e<b<z<t\), \(a<t\), with \(a\) incomparable with \(e, b, z\). Element names are not constants in the signature. The catalog relabelings in the displayed order are \((3,\allowbreak{}2,\allowbreak{}1,\allowbreak{}0,\allowbreak{}4)\) for \(A\) (type 284) and \((2,\allowbreak{}3,\allowbreak{}1,\allowbreak{}0,\allowbreak{}4)\) for \(B\) (type 332).

These are semirings. A product of at least three factors is \(e\) when all factors are \(e\), and \(z\) otherwise, which proves associativity. Every translation except those by \(e, a\) is constant \(z\). Translation by \(e\) maps \(e\) to \(e\) and all other elements to \(z\); translation by \(a\) maps \(a\) to \(b\) and all other elements to \(z\). Each preserves the listed joins. In \(A\) the only incomparable pair is \(e, a\); in \(B\) the incomparable pairs are \(a\) and any element of \(\{e, b, z\}\), and their joins are \(t\). This verifies both distributive laws.

Theorem 1. The following five laws are a complete identity basis for \(A\):

\CaseDisplay{\begin{array}{rll}
(C)&xy\approx yx,\\
(P)&x^4\approx x^3,\\
(D)&xyz\approx x^3+y^3+z^3,\\
(E)&xy\leq xz+zy,\\
(M)&x^3\leq xy+y^3.
\end{array}}

The same laws, together with

\CaseDisplay{(H)\qquad x^3\leq x+xy,}

are a complete identity basis for \(B\). No minimality is asserted.

\par\smallskip\noindent\textit{2. Polynomial normal forms.}\quad 
From (P), all powers \(x^n\), \(n\geq3\), equal \(x^3\). Commutativity and (D) therefore give

\CaseDisplay{(uv)^3=u^3v^3\approx u^9+v^6+v^3\approx u^3+v^3.}

Here (D) is applied to \(u^3, v^2, v\); \(u, v\) may be any nonempty terms. Induction yields \(w^3\approx\sum_{x\in c(w)}x^3\) for every word \(w\). For a word of length at least three, split it into three nonempty factors and apply (D) and the preceding formula. Thus every term is provably a nonempty polynomial of the form

\CaseDisplay{F=\sum_{x\in L}x+\sum_{\{x, y\}\in E}xy+\sum_{x\in T}x^3,}

where \(E\) is a finite undirected graph allowing loops. Empty summation blocks are omitted; no empty term or unit is introduced. Let \(V\) be the vertices occurring on an edge. Regard each vertex in \(T\setminus V\) as an isolated vertex. The vertex set for components is \(V\cup T\); linear variables outside it are not graph vertices.

For \(A\), mark the components meeting \(T\). For \(B\), mark the components meeting \(T\cup(L\cap V)\). Let \(K\) be the union of the marked components.

Lemma 2 (exact absorption). In either algebra:

1. \(x\leq F\) holds exactly when \(x\in L\).

2. \(x^3\leq F\) holds exactly when \(x\in K\).

3. \(xy\leq F\) holds exactly when either \(x, y\) belong to the same component of the edge graph on \(V\), or both belong to \(K\).

In the third clause \(x=y\) is allowed, but the first alternative requires \(x\in V\). An isolated cubic vertex is covered by the second alternative.

Proof. The values of a quadratic word are \(e, b, z\): it is \(e\) if both factors are \(e\), \(b\) if both are \(a\), and \(z\) otherwise. A cube is \(e\) at \(e\) and \(z\) at every other element.

If \(x\notin L\), assign \(x=t\) and every other variable \(e\). All nonlinear right summands are at most \(z\) and all right linear summands are \(e\), so \(x\nleq F\). This proves the first clause.

If \(F\) is purely linear, assign every variable \(a\). Then \(F=a\), while any quadratic or cubic left side is \(b\) or \(z\), neither at most \(a\). Thus there is no nonlinear absorption in this case, as predicted.

Suppose \(F\) has a nonlinear summand. It cannot have value \(a\). Its possible values are consequently \(e, b, z, t\). The values \(z, t\) absorb every nonlinear left word. If \(F=e\), every variable in \(F\) is assigned \(e\). Every left word permitted by clauses 2 or 3 therefore also equals \(e\).

It remains to consider assignments with \(F\leq b\). Each edge has value \(e\) or \(b\), so every component of the edge graph is uniformly \(e\) or uniformly \(a\). A cubic summand forces its variable, and hence its entire component, to be \(e\). In \(A\), linear summands may be \(e, a, b\) and place no further restriction on these components. In \(B\), a linear summand must be \(e\) or \(b\); if its variable is on an edge, that variable must therefore be \(e\). These are exactly the marked components. This proves the sufficiency of clauses 2 and 3.

For necessity in clause 2, if \(x\notin V\cup T\), assign \(x=b\) and all other variables \(e\). Then \(F\leq b\) and \(x^3=z\). If \(x\) belongs to an unmarked edge component, assign that component \(a\) and all other variables \(e\). The right side is at most \(b\), but \(x^3=z\). For \(B\) that component contains no linear variable, so the assignment indeed keeps \(F\leq b\).

For necessity in clause 3, an endpoint outside \(V\cup T\) can likewise be assigned \(b\), with every other variable \(e\); the left product is \(z\). Otherwise failure of the condition means the endpoints lie in different components and at least one component is unmarked. Assign that component \(a\) and every other variable \(e\). The left product is \(z\), while \(F\leq b\). For equal endpoints, failure requires a vertex outside \(V\cup T\), already treated. These separating assignments complete the proof. \(\square\)

\par\smallskip\noindent\textit{3. Validity and completeness of the displayed bases.}\quad 
Laws (C), (P), and (D) follow from the multiplication formula. In (E), the two right edges connect the left endpoints, so Lemma 2 proves validity. In (M), the right cubic marks the component of the right edge. In \(B\), (H) marks that edge through its linear endpoint. Thus all the proposed laws are valid, using a semantic criterion proved independently of them.

We next derive every absorption in Lemma 2. Law (E), along a path, produces an edge between any two distinct vertices of one component. Substituting \(y=x\) gives \(x^2\leq xz+zx\approx xz\); every vertex on an edge therefore also acquires its loop. This includes singleton loop components. The resulting complete component graph is derivable below \(F\).

Cubic marks initially in \(T\) are already summands. For \(B\), if \(x\in L\) is on an edge \(xy\), (H) gives \(x^3\leq x+xy\leq F\). Law (M) then propagates a cubic mark across any edge, so \(x^3\leq F\) is derivable for all and only the required \(x\in K\).

We also have the derived law

\CaseDisplay{xy\leq x\,x^2+x^2y\approx x^3+y^3,}

by (E) with \(z=x^2\) and (D). Thus any two marked vertices yield \(xy\leq x^3+y^3\leq F\), including equal vertices. All quadratic absorptions in the criterion are now derived; linear ones are summands.

Finally, normalize both sides of an arbitrary valid identity. Every summand of each side is absorbed by the other side and, by the preceding argument, each such absorption is derivable. The two resulting inequalities prove the identity. No bound was placed on the graph, the number of variables, the number of summands, or the lengths of paths. This proves Theorem 1. \(\square\)

\subsection{Twelve complete identities for 13841}\label{case:26}
Signature (+, multiplication), without constants. Relative to the additively idempotent semiring axioms, use

{\small\par\smallskip
\CaseLine{\(xy=yx; xyz=xyzz; xxy\leq xxx+yyy; xxy\leq xxx+xy;\)}
\CaseLine{\(xxyz\leq xxy+xz; xy\leq xxy; x\leq xx;\)}
\CaseLine{\(xxx\leq x+xxz; xxy\leq xy+xxz; xxy\leq xxx+yyz;\)}
\CaseLine{\(xx\leq x+xy; xz\leq xx+xy+yz.\)}
\smallskip}

Write \(u\leq P\) for \(u+P=P\). Addition is 00044/\allowbreak{}01044/\allowbreak{}00244/\allowbreak{}44434/\allowbreak{}44444. The elements1 and 2 join to 0; 3 is incomparable with each of 0,\allowbreak{}1,\allowbreak{}2, their joins with 3 being 4. Multiplication is 44434/\allowbreak{}40434/\allowbreak{}44234/\allowbreak{}33333/\allowbreak{}44434. Nonlinear words have value 3 if a variable is 3; otherwise2 if all variables are 2; otherwise0 for a quadratic with both occurrences 1; and 4 in all other cases. Literals retain their value. Induction proves these formulas at every word length.

\par\smallskip\noindent\textit{Exact criterion with original global components.}\quad 
For original P, form a graph on the variables occurring in its nonlinear words, joining all variables in each such word. Let V be this vertex set. A component is active when it contains an original long word, and A is the union of active components. Equivalently start with the contents of original long words and propagate along original quadratic edges. These global witnesses may have variables outside the target content \(U=c(w)\).

Separately, require C: there is at least one original right word whose entire content is contained in U. Then the exact criterion is:

\textbullet{} A literal requires C, equivalently an original pure word.

\textbullet{} A long target requires C and U contained in A.

\textbullet{} A quadratic requires C and either both endpoints in A, or both endpoints in V in the same original component.

Neither the contained witness nor the component conditions may be omitted. In particular, literals do not provide graph edges.

If C fails, set all target variables 2 and outsiders 3. Every right word is 3 but the target is 2, which is not below 3. A missing literal is handled by the same assignment. If a long target has an inactive variable, assign every A variable 2 and all other variables 1. Every original long word is 2; each quadratic is 2 or 0, since its endpoints lie in one component, and literals are 1 or 2. The right join is at most 0, while the target is 4.

For a quadratic failing its component condition, if an endpoint x is outside V, assign \(x=0\) and all other variables 2. Only original literals can involve x, so the right join is at most 0, but the target is 4. Otherwise its endpoints lie in different components and at least one is inactive. Assign that inactive component 1 and all other variables 2. No long word meets the inactive component; all right values are at most 0 and the target is 4. These assignments prove necessity for arbitrarily large polynomials.

For sufficiency, ignore the trivial case that P evaluates to 4. If P evaluates to 3, every original right word evaluates to 3, so the contained witness forces some target variable to 3. The target is then 3. Otherwise all right values lie in \{0,\allowbreak{}1,\allowbreak{}2\}, and no variable of P is 3. Every original long word forces its variables 2, and quadratic edges propagate2 through A. Every nonlinear graph vertex must be 1 or 2, and each original quadratic makes its endpoints equal. Consequently all variables of an inactive component agree. A long target on A is 2. A quadratic on A is 2; a quadratic in one inactive component is 2 or 0, with the same value supplied by an original quadratic in that component. A loop component also has an original nonlinear witness. Literals are below their original pure powers. The right side thus dominates the target in every case.

\par\smallskip\noindent\textit{Finite derivations with unrestricted sizes.}\quad 
The first two laws normalize long words by content. Commutativity moves a chosen variable to the end; \(xyz=xyzz\) duplicates it with two nonempty preceding blocks or removes surplus copies in reverse while retaining length at least three. Canonical forms are sorted content for at least three variables, xxy for two, and xxx for one. In particular \(g^{3}=g\) for any long g. Quadratics retain multiplicity.

First suppose U is contained in A. Merge all original long words using \(xxy\leq xxx+yyy\) with long substitutions. Retain the resulting long g. For an original edge xz with x in c(g) and z outside, if \(c(g)=\{x\}\), use \(xxy\leq xxx+xy\) at \(y=z\). Otherwise choose a nonempty word p listing c(g) minus x; normalize g to xxp, and use \(xxyz\leq xxy+xz\) with \(y=p\). This grows g by z. Repeating produces a long g with content exactly A, using original edges throughout.

Now retain an original contained witness u. If u is literal x, normalize g to xxq and use \(xxx\leq x+xxz\) with \(z=q\). Here q lists c(g) minus x, or is x when g has only x; it is always nonempty. If \(u=xy\) is quadratic, use \(xxy\leq xy+xxz\) with the same representation of g. If u is already long, retain it directly. In each case we obtain a long word h contained in U. For every missing target variable y, represent g as yyq and apply \(xxy\leq xxx+yyz\) with \(x=h\), \(z=q\). The right terms normalize to h and g; the left is a long word on c(h) union \{y\}. Repeat until the content is U. Normalize to the long target, or recover a quadratic by \(xy\leq xxy\). This is the point where the separate original contained witness removes unwanted external variables from the global construction.

It remains to derive a quadratic in one inactive component. An original target is immediate. Otherwise its contained witness is an original literal or endpoint loop: a mixed contained quadratic would be the target, and a contained long word would activate the component. Choose an endpoint x with that witness. If only its literal is available, retain an original incident quadratic xy and use \(xx\leq x+xy\) to derive the loop. Such an incident edge exists because x is in V and the component is inactive.

For a square target the loop finishes the proof. For distinct endpoints choose a path of original quadratics from x to the other endpoint. Its length is positive. The first edge is original; given a derived x-to-y edge and the next original yz, use \(xz\leq xx+xy+yz\), keeping the loop xx throughout. Induction on the path gives the target edge. No empty word is substituted.

A literal target has an original pure word; if long normalize to xxx and use \(xy\leq xxy\) at \(y=x\) to obtain xx, then \(x\leq xx\). These finite derivations cover every valid absorption with no bound on variables, word lengths, components, or paths. Distributivity and absorption in both directions give every identity. All twelve laws satisfy the criterion and are checked on all scalar assignments.

For boundaries, \(xxx\leq xx\) fails at \(x=1\), \(xy\leq xx+yy\) fails at \(x=1\), \(y=2\), and \(xx\leq x\) fails at \(x=0\). Without a contained witness, \(xxx\leq yyy\) fails at \(x=2\), \(y=3\). Without activation, \(xyz\leq xy+yz\) fails at \(x=y=z=1\). These exclude weaker but incorrect criteria.

\subsection{Nonlinear components and terminal witnesses: a fifteen-identity basis}\label{case:27}
Aifa Wang and Lili Wang give the following argument for \(\A{12280}\). The opposite class is represented by \(\A{14324}\).

\par\smallskip\noindent\textit{Algebra and finite basis.}\quad 
Work relative to the axioms for additively idempotent semirings in the signature (+, *) without constants. The tables, in element order 0, 1, 2, 3, 4, are:

\textbullet{} | 0 1 2 3 4       * | 0 1 2 3 4 --+-----------      --+----------- 0 | 0 0 0 4 4       0 | 4 3 4 3 4 1 | 0 1 0 3 4       1 | 4 1 4 3 4 2 | 0 0 2 4 4       2 | 4 3 0 3 4 3 | 4 3 4 3 4       3 | 4 3 4 3 4 4 | 4 4 4 4 4       4 | 4 3 4 3 4

As usual u \ensuremath{\leq} v abbreviates u + v = v. A relative basis consists of:

{\small\par\smallskip
\CaseLine{\(N1: xyz = yxz\)}
\CaseLine{\(N2: x^{2}yz = xyz\)}
\CaseLine{\(N3: xy^{2} = x^{2}y\)}
\CaseLine{\(L:  y \leq  xy\)}
\CaseLine{\(B:  y^{2} \leq  xy\)}
\CaseLine{\(B0: x^{2} \leq  x + xy\)}
\CaseLine{\(T:  xz \leq  xy + yz\)}
\CaseLine{\(F:  xz \leq  yx + yz\)}
\CaseLine{\(I:  x^{2}z \leq  xyz\)}
\CaseLine{\(P:  x^{2}z \leq  xy + y^{2}z\)}
\CaseLine{\(Q:  y^{2}z \leq  xy + x^{2}z\)}
\CaseLine{\(U0: x^{3} \leq  x + x^{2}y\)}
\CaseLine{\(U1: x^{3} \leq  zx + x^{2}y\)}
\CaseLine{\(A:  x^{2}yt \leq  x^{2}t + y^{2}z\)}
\CaseLine{\(J:  xy \leq  x^{2}z + y^{3}\)}
\smallskip}

Reversing every word gives a basis for the multiplicative opposite.

\par\smallskip\noindent\textit{Values and long-word normalization.}\quad 
Direct inspection verifies the semiring axioms. A two-letter word has value 1 precisely when both entries are 1, and value 0 precisely when both entries are 2. It has value 3 precisely when its final entry is 1 or 3 and its entries are not both 1. All remaining cases give 4.

For a word of length at least three, value 1 means that every entry is 1. Value 3 means that its final entry is 1 or 3 and not every entry is 1. Otherwise its value is 4. In particular long words never have value 0 or 2.

For every nonempty finite variable set C let \(U_{C}\) list its variables once in fixed order. N1--N3 transform every word w of length at least three into

{\small\par\smallskip
\CaseLine{\(U_{C}(w)^{2} t,\)}
\smallskip}

where t is its final letter. Here and below length means the number of positions, with repetitions counted.

To prove normalization, N1 permutes all positions before the last. N2 inserts or removes an internal duplicate as long as at least two internal positions remain. Indeed move the duplicate to the first two positions and write the remainder as a nonempty internal word followed by the final letter; this is exactly the nonempty-substitution boundary in N2. Therefore any two internal words of length at least two with the same content are interchangeable. If w = p t, then p has length at least two, and

{\small\par\smallskip
\CaseLine{\(p t = p^{2} t = p t^{2}\)}
\smallskip}

follows first from this internal-content argument and then from N3. Thus the final letter may be inserted among the internal positions. Now every variable of C(w) occurs internally, and the internal-content argument gives the claimed normal word. No empty word is substituted.

\par\smallskip\noindent\textit{Absorption criterion.}\quad 
Let P be a nonempty finite sum of words. Define T(P) as its set of final letters and N(P) as the union of the contents of its nonlinear summands. On N(P), let E(P) be the equivalence relation generated by co-occurrence in any nonlinear summand. Let F(P) be the union of the E(P)-classes meeting a summand of length at least three.

The following statements characterize absorption exactly.

1. A linear word x is below P if and only if x belongs to T(P).

2. A word xy of length two is below P if and only if y belongs to T(P), x and y belong to N(P), and either x E(P) y or both x, y belong to F(P).

3. A word w of length at least three is below P if and only if its final letter belongs to T(P) and C(w) is contained in F(P).

Proof. A missing final-letter witness is separated by assigning that letter value 4 and every other variable value 1. Then the target is 4, while P is 1 or 3.

For the linear case, a word ending in x is always above x by direct inspection of the multiplication table. This proves its sufficiency.

Consider a valuation with P = 0. Every nonlinear summand must be 0 or 1. Consequently its variables are constant on its content with value 1 or 2; a summand of length at least three forces value 1. Equivalently variables are constant, with value 1 or 2, on each E(P)-class, and every F(P)-class has value 1. Variables occurring only linearly may have values 0, 1, or 2.

These observations prove sufficiency when P = 0. A valid binary target then has both letters 1 or both 2, and so is 1 or 0; a valid long target is 1. If P = 1, every variable of its content is 1. If P = 3, every final letter in T(P) has value 1 or 3, which puts a nonlinear target with such a final letter below 3. If P = 4, the inequality is automatic. Finally P = 2 is impossible in the presence of a nonlinear summand; both nonlinear criteria require such a summand. This proves sufficiency in all cases.

For necessity, if a target variable x is outside N(P), assign x = 0 and all other variables value 1. Right nonlinear words are then 1; right linear words are 0 or 1. A nonlinear target containing x is 3 or 4, so it is not below P.

Next suppose a long target contains x in N(P) outside F(P). Assign its unforced E(P)-class value 2 and all other variables value 1. This class contains a nonlinear summand and contains no long summand, so P = 0. The long target is 3 or 4, giving a separation.

For a binary target with letters in different E(P)-classes and not both forced, give an unforced class containing one of its letters value 2, and all other variables value 1. Again P = 0. The target is 3 if its first letter alone is 2, and 4 if its final letter alone is

2. The criterion follows.

Every displayed basis identity is sound by this criterion, or directly by the word-value rule.

\par\smallskip\noindent\textit{Constructing certificates for forced variables.}\quad 
Call a word available if its absorption into P is derivable. For each x in F(P), we construct an available word \(x^{2}r\); its final letter r may depend on x and belongs to T(P).

First let u be a right summand of length at least three, with content C and final letter r. For every x in C, normalization gives

{\small\par\smallskip
\CaseLine{\(u = x U_{C} r,\)}
\smallskip}

because both are long words with the same content and final letter. Identity I gives \(x^{2}r\) \ensuremath{\leq} u. Thus every variable of a long summand has a certificate.

All further connections to a forced class may be propagated through two-letter right summands. For an edge xy, if \(y^{2}r\) is available, P gives \(x^{2}r\) \ensuremath{\leq} xy + \(y^{2}r\). If \(x^{2}r\) is available, Q gives \(y^{2}r\) \ensuremath{\leq} xy + \(x^{2}r\). Repetition along the finite undirected co-occurrence graph reaches every variable of F(P). Co-occurrences from long words need no propagation because all their variables already received initial certificates.

If t belongs to both F(P) and T(P), its certificate \(t^{2}r\) and a right summand ending in t yield \(t^{3}\). Use U0 when that right summand is the linear word t; otherwise write it as Zt with nonempty Z and apply U1. These separate boundary identities are essential.

\par\smallskip\noindent\textit{Deriving long and binary targets.}\quad 
Let w be a valid long target with final letter t. Start with the available cube \(t^{3}\). Suppose a long word with content D and final letter t is available. Normalize it to \(U_{D}^{2}\) t. If y is another target variable, its available certificate \(y^{2}r\) and identity A give

{\small\par\smallskip
\CaseLine{\(U_{D}^{2} y t \leq  U_{D}^{2} t + y^{2}r \leq  P.\)}
\smallskip}

Normalize, and repeat until the content is C(w). Long-word normalization then gives w \ensuremath{\leq} P. This derivation does not require all target variables to belong to one co-occurrence class.

For a valid binary target xy with both variables forced, use a certificate \(x^{2}r\), the cube \(y^{3}\) constructed above, and J:

{\small\par\smallskip
\CaseLine{\(xy \leq  x^{2}r + y^{3} \leq  P.\)}
\smallskip}

Otherwise x and y lie in the same unforced E(P)-class. All nonlinear summands in this class have length two, so its connections are actual two-letter right words.

First construct \(y^{2}\). If any two-letter right summand ends in y, apply B. This choice is made even if a linear right summand y also exists. If none ends in y, the final-letter witness must be linear, and the nonlinearity of y supplies a two-letter right summand yz starting in y. Then B0 gives \(y^{2}\) \ensuremath{\leq} y + yz. This is the boundary case that a terminal-witness argument must not omit.

Starting from the available word yy, follow an undirected path of two-letter right words from y to x. If vy is already available and the next edge is uv, apply T to obtain

{\small\par\smallskip
\CaseLine{\(uy \leq  uv + vy.\)}
\smallskip}

If instead the next edge is vu, apply F to obtain

{\small\par\smallskip
\CaseLine{\(uy \leq  vu + vy.\)}
\smallskip}

The same final letter y is retained throughout. At the end xy is available. Thus every valid binary target is derivable. Linear targets follow from L or from being original summands.

Distributivity converts arbitrary terms to finite nonempty sums of words. Derive every summand of each side of a valid identity below the other side; idempotent addition gives equality. This proves unrestricted completeness. Word reversal gives completeness for the opposite.

\subsection{Two finite bases retaining head and long-word content}\label{case:28}
For each of 299 and 1809, in the constant-free signature, take the semiring axioms and

{\small\par\smallskip
\CaseLine{\(xyz=xzy, xxyz=xyz, xyyz=xyz, xxy=xyy,\)}
\smallskip}

together with all valid absorptions on 2069 variables, word length at most 2071, and at most 10 upper words. This is an effective finite family rather than an enumerated small list. Completeness below has no restriction on the variables, lengths or sizes of input polynomials.

The multiplication rows of 299 are 00000,\allowbreak{}00000,\allowbreak{}00100,\allowbreak{}33333,\allowbreak{}00004; those of 1809 are 00000,\allowbreak{}11111,\allowbreak{}00000,\allowbreak{}00020,\allowbreak{}00004. Introduce parameters (h, q, r, \(U)=(3\),2,\allowbreak{}1,\{1\}) for 299 and(1,\allowbreak{}3,\allowbreak{}2,\{2,\allowbreak{}3\}) for 1809. Every nonlinear word with head value h has value h. A word whose letters are all 4 has value 4. A quadratic with both letters q has value r. All other nonlinear words have value 0. These alternatives are disjoint and follow by induction, including arbitrary repeated variables. A long word is therefore determined by its head and its content. Literals and quadratics retain their distinct behavior. This verifies the four laws.

Normalize a long word by preserving its first letter and permuting its tail using \(xyz=xzy\). If length is at least 4 and the head occurs in its tail, move that copy immediately after the head and use \(xxyz=xyz\). If a nonhead variable repeats in the tail, move its copies immediately after the head and use \(xyyz=xyz\). Continue until the word has length 3 or distinct content. If length 3 has two variables, \(xxy=xyy\) and tail permutation give the head repeated twice, then the other variable. Thus the canonical long word is the head followed by sorted nonhead content when it has at least 3 variables; otherwise the head square followed by the other variable, or the head cube. It always stays long.

\par\smallskip\noindent\textit{Exact criterion from four proper ideals.}\quad 
For 299 the proper ideals are\{0,\allowbreak{}1,\allowbreak{}2,\allowbreak{}3\},\{1\},\{2\},\{3\}; for 1809 they are \{0,\allowbreak{}1,\allowbreak{}2,\allowbreak{}3\},\{1\},\{2,\allowbreak{}3\},\{3\}. In both, one is\{0,\allowbreak{}1,\allowbreak{}2,\allowbreak{}3\}, one is the head ideal\{h\}, one is the special quadratic ideal U, and the remaining one contains no nonlinear word values. Membership in the head ideal means exactly that the word's first variable has value h. Membership in U means a literal variable lies in U or a quadratic has both values q; long words never belong to U. This follows directly from the word values.

For original P, form domains by intersecting U for a variable occurring as a literal, and\{q\} for a variable occurring in either position of an original quadratic. Domains start as the whole algebra. Call the system impossible if P contains a long word or some domain is empty. Otherwise the system describes exactly all assignments putting P in U, with choices independent between variables. Then \(w\leq P\) is equivalent to all of:

\textbullet{} Some original word has content contained in content(w).

\textbullet{} The target head is the head of an original word.

\textbullet{} A nonlinear target has an original nonlinear word.

\textbullet{} Either the domain system is impossible, or the target is literal with its domain contained in U, or the target is quadratic with both its variable domains contained in\{q\}. A long target is forbidden in the feasible-domain case.

For the first condition, all word values except 4 belong to the first ideal, and a word is 4 precisely when all its variables are 4. If there is no contained original, set all target variables 4 and all other variables 0. Every original is then non4 while the target is 4. If a contained original exists, the implication is immediate. The head ideal gives the second condition:if the target head is missing, assign all original heads h and that missing head 0. The literal-only ideal gives the third; with originals all literal, containment already protects a literal target, and assigning all variables in the literal-only ideal separates a nonlinear target. Finally the independent domains give the fourth condition exactly, including coincident target variables. Every failed order comparison is detected by a principal ideal of the value of P. Hence this is unrestricted necessity and sufficiency.

\par\smallskip\noindent\textit{Bounded original support and reversible compression.}\quad 
Retain one contained original, one original having the target head, and one nonlinear original if needed. For the domain condition retain a long original if available. Otherwise an impossible domain needs at most 3 originals witnessing its literal, quadratic-head and quadratic-tail roles. Otherwise the target has at most 2 variables:retain one original for each of their at most 3 roles, at most 6 words. Thus at most 9 selected original words suffice, and 10 is a uniform bound. Their polynomial Q still satisfies every criterion for w, so \(w\leq Q\) is valid.

Normalize the at most 11 words consisting of target and Q. Reserve every head and every variable in a literal or quadratic. At most 22 variables are reserved. Group other variables by nonempty incidence subset among those 11 words, at most \(2^{11}\)-1 groups. They occur solely in long-word tails and never as heads or in short words. Identify each group, normalize, and obtain a valid member of the stated finite family: at most 22+2047 variables, each word length at most content size plus 2. Substitute the nonempty list of members back for each group variable. Every expanded long word stays long, has the same head, and recovers exactly its original content because the members had COMMON incidence. The four laws therefore normalize back to the originals. Short words were reserved and recover literally. This derives \(w\leq Q\), then \(w\leq P\). Normalize backwards to the input words. Polynomial identities split into the summand absorptions in both directions, proving completeness.

The related catalogue types 1912 and 1913 are transferred by complete mutual HSP from 299 and 1809; no new standalone proof for their differing operation tables is assumed. All further transfers, including opposites and nonisomorphic targets, require independently verified full inclusions.

\subsection{Nine complete identities for semiring 319}\label{case:29}
Constant-free signature (+, multiplication). Relative to the additively idempotent semiring axioms a complete finite basis is

{\small\par\smallskip
\CaseLine{\(xy=yx; xyz=xyzz; xyzt\leq xyz; xy\leq x+xz+yt;\)}
\CaseLine{\(xxx\leq x+yyy; xxy\leq xy+zzz; xy\leq xxy; xy\leq xx+yz; x\leq xxx.\)}
\smallskip}

Here \(u\leq P\) abbreviates \(u+P=P\).

\par\smallskip\noindent\textit{Exact values and criterion.}\quad 
Addition is 00004/\allowbreak{}01014/\allowbreak{}00224/\allowbreak{}01234/\allowbreak{}44444, with \(3<1<0<4\) and \(3<2<0<4\);1 and 2 are incomparable. Multiplication is 00000/\allowbreak{}00000/\allowbreak{}00110/\allowbreak{}00110/\allowbreak{}00004. A nonlinear word is 4 on all 4 content; otherwise a quadratic is 1 if both occurrences are in\{2,\allowbreak{}3\}, and every other nonlinear word is 0.

For original P let L be its literal variables and V the union of quadratic endpoints. Let G say that an original long word occurs; let C(w) mean that some original word has content contained in c(w). Absorption \(w\leq P\) holds exactly when:

\textbullet{} Literal x: original x, or C(x) and G.

\textbullet{} Long w: C(w) and G.

\textbullet{} Quadratic w: C(w) and (G or c(w) contained in V).

If C fails, assigning target content 4 and outsiders 3 separates: the target is 4 and all right summands are below 4. When G fails, assign V to 3 and outsiders 1. Every right summand is at most 1; a long target and any quadratic with an endpoint outside V are 0. For a missing literal x with C, an original square xx exists. Assign \(x=2\), the other V variables 3, and outsiders 1. All right quadratics are 1 and every right literal is at most 1; target 2 is incomparable with 1. Thus G is necessary even if many original literals overlap original edges.

For sufficiency G forces \(P\geq 0\). A contained witness supplies 4 whenever the target is 4, while all other values are at most 0. When G fails, all right quadratics have value at least 1. If P does not dominate 0, they must all be 1, so all V variables are in\{2,\allowbreak{}3\}. Any quadratic on V is therefore 1. Its value 4 is handled by the contained witness. Original literals are immediate. This proves the exact criterion for arbitrary variable sets and polynomial sizes.

\par\smallskip\noindent\textit{Complete finite derivation.}\quad 
The two word laws normalize long words by content: commute a variable to the end, duplicate with \(xyz=xyzz\) using two nonempty prefix blocks, then remove surplus occurrences by the reverse law, never shortening below length three. Canonical long forms are sorted content for at least three variables, xxy for two, and xxx for one. Long g satisfies \(g^{3}=g\).

When G holds, fix an original long g and a contained original u. Upgrade a literal u to \(u^{3}\) with \(xxx\leq x+yyy\) using \(y=g\), and a quadratic ab to aab with \(xxy\leq xy+zzz\) using \(z=g\). Retain a long u. Append missing target variables by \(xyzt\leq xyz\) with three nonempty left blocks and normalize. Recover a quadratic with \(xy\leq xxy\) or a literal with \(x\leq xxx\). Global and contained original witnesses may be different.

Without G, a new quadratic has an endpoint loop or literal as its contained original witness. An original loop aa and original edge bc at the other endpoint derive ab by \(xy\leq xx+yz\). An original literal a and original edges az, bt derive ab by \(xy\leq x+xz+yt\). Keep both incidence witnesses; an isolated literal is insufficient. Commute endpoints when needed. No step recovers a missing literal or long word without G.

Every law satisfies the exact criterion and is checked independently on all assignments. All substitutions are nonempty. The construction is unrestricted in variables, lengths and summands. Distributivity and summandwise absorption in both directions establish completeness. Other types require independently verified full HSP equivalence.

\subsection{首字母平方回缩：六条完整基}\label{case:30}
\par\smallskip\noindent\textit{1. 构造与词值.}\quad 
令 \(A=\{z, b, a, d, t\}\)，其加法为偏序

\CaseDisplay{z<b<a<t,\qquad z<d<t,\qquad d\parallel b, a}

上的并。令 \(B=\{z, b, a, t\}\)，定义回缩 \(\rho(d)=t\)，其余元素不变。在 \(B\) 上，\(b, a\) 是左零元素，且

\CaseDisplay{zz=z,\quad zq=b\ (q\ne z),\qquad tt=t,\quad tq=a\ (q\ne t).}

在 \(A\) 上定义 \(xy=\rho(x)\rho(y)\)，右边使用 \(B\) 的乘法。完整乘法表为

\CaseDisplay{\begin{array}{c|ccccc}
\cdot&z&b&a&d&t\\\hline
z&z&b&b&b&b\\
b&b&b&b&b&b\\
a&a&a&a&a&a\\
d&a&a&a&t&t\\
t&a&a&a&t&t
\end{array}.}

这确实给出加法幂等半环。为核对分配律，固定右因子为 \(z\)、\(b/a\) 或 \(d/t\) 时，右平移在顺序 \((z, b, a, d, t)\) 上的值依次为 \((z, b, a, a, a)\)、\((b, b, a, a, a)\)、\((b, b, a, t, t)\)；它们都保持两条链的次序，并把两个非平凡并 \(b+d=t\)、\(a+d=t\) 保持。左平移中，\(b, a\) 给常值映射；\(z\) 把 \(z\) 送到 \(z\)、其余送到 \(b\)；\(d, t\) 把 \(\{z, b, a\}\) 送到 \(a\)、把 \(\{d, t\}\) 送到 \(t\)，也保持并。

每个长度至少二的词 \(w\)，设首字母为 \(h\)、内容为 \(C\)，其值满足：

1. \(h=b\) 或 \(h=a\) 时，值分别为 \(b\) 或 \(a\)；

2. \(h=z\) 时，所有 \(C\) 内字母都取 \(z\) 则值为 \(z\)，否则为 \(b\)；

3. \(h\in\{d, t\}\) 时，所有 \(C\) 内字母都在 \(\{d, t\}\) 则值为 \(t\)，否则为 \(a\)。

这条规则在两个非空词的连接下保持，因而同时验证乘法结合律。目录中 \(A\) 的编号为 1960，按 \((z, b, a, d, t)\mapsto(2,\allowbreak{}1,\allowbreak{}0,\allowbreak{}3,\allowbreak{}4)\) 重标号；6230 是乘法对偶。

\par\smallskip\noindent\textit{2. 完整吸收判据.}\quad 
把非空多项式 \(F\) 写成一次项与非线性词的和。记 \(L(F)\) 为一次项变量集，\(R(F)\) 为非线性求和词族，并置

\CaseDisplay{N(F)=\bigcup_{u\in R(F)}c(u),\qquad H(F)=L(F)\cup\{h(u):u\in R(F)\}.}

引理。 在 \(A\) 中：

\textbullet{} 对一次词 \(x\)，\(x\le F\) 当且仅当 \(x\in L(F)\)，或存在 \(u\in R(F)\) 满足 \(c(u)=\{x\}\)。

\textbullet{} 对非线性词 \(w\)，设 \(h=h(w)\)、\(C=c(w)\)，则 \(w\le F\) 当且仅当

\CaseDisplay{C\subseteq N(F),\quad\textup{(\ensuremath{2.1})}}

\CaseDisplay{h\in H(F),\quad\textup{(\ensuremath{2.2})}}

\CaseDisplay{L(F)\cap C\ne\varnothing\quad\text{或}\quad \exists u\in R(F):c(u)\subseteq C.\quad\textup{(\ensuremath{2.3})}}

必要性。 对一次词，若没有同变量一次项，也没有只用该变量的非线性右词，令 \(x=d\)、其余变量为 \(z\)。所有右词的值属于链 \(\{z, b, a\}\)，其和不可能在 \(d\) 之上。

若 (2.1) 失败，取 \(x\in C\setminus N(F)\)，赋 \(x=d\)、其他变量为 \(z\)。右侧值只能为 \(z\) 或 \(d\)。左侧在 \(x\ne h\) 时为 \(b\)；在 \(x=h\) 且 \(|C|>1\) 时为 \(a\)；在 \(C=\{h\}\) 时为 \(t\)，均不被右侧吸收。

若 (2.2) 失败，赋 \(h=a\)、其他变量为 \(z\)，左侧为 \(a\)，右侧至多为 \(b\)。若 (2.3) 失败，赋 \(C\) 内所有变量为 \(d\)、外部变量为 \(z\)。左侧为 \(t\)；右侧没有值为 \(d\) 的一次项，也没有值为 \(t\) 的非线性词，故至多为 \(a\)。

充分性。 一次词的第二种情形由 \(x\le x^2\) 及 \(x^n=x^2\)（\(n\ge2\)）成立；这两条关系可由表直接读出。以下设左词非线性且满足 (2.1)--(2.3)。若其值为 \(z\)，结论自动成立。若其值至少为 \(b\)，其内容中有非 \(z\) 的字母；由 (2.1)，某个右非线性词含此字母，其值属于 \(\{b, a, t\}\)。

若左值为 \(a\) 且首字母取 \(a\)，由 (2.2)，右侧有同首字母非线性词或该一次项，均给出 \(a\)。若首字母取 \(d\) 或 \(t\)，同首字母非线性右词给出 \(a\) 或 \(t\)；同名一次项若为 \(t\) 即足够，若为 \(d\)，则与前段所得至少为 \(b\) 的右值求并得到 \(t\)。

最后，左值为 \(t\) 时，所有 \(C\) 内变量都在 \(\{d, t\}\)。若 (2.3) 提供内部非线性右词，该词直接取 \(t\)；若提供内部一次项，其值为 \(t\) 时直接完成，为 \(d\) 时与前段所得右值求并也得到 \(t\)。引理得证。

\par\smallskip\noindent\textit{3. 六条有限基.}\quad 
定理。 相对于加法幂等半环公理，\(A\) 的完整恒等式基为

\CaseDisplay{\begin{aligned}
xyz&\approx xzy, &(N)\\
xy&\approx x^2y^2, &(S)\\
x&\le x^2, &(I)\\
xy&\le x^2+y^2, &(J)\\
xy&\le x+yx, &(H)\\
xy&\le y+xyz. &(E)
\end{aligned}}

把所有乘法词反转，得到 \(A^{\mathrm{op}}\) 的六条完整基。

有效性。 (N)、(S) 来自词值规则，(I) 来自 \(d<t\)；(J)、(H)、(E) 的非线性左词和右侧均直接满足引理的三个条件，因而有效。

正规形。 (S) 给出 \(x^2=x^4\)，继而

\CaseDisplay{x^3=(x^2)x=x^4x^2=x^6=x^2.}

由 (N)，首位后的字母可以任意交换；由 (S)，每个长度至少二的词可以把所有字母平方。结合幂收缩，一个非线性词完全由首字母 \(h\) 和内容 \(C\) 决定。记其正规形为 \(W(h, C)\)，其中每个字母出现两次，首字母为 \(h\)。因此

\CaseDisplay{W(h, C)W(g, D)=W(h, C\cup D).\quad\textup{(\ensuremath{3.1})}}

整个论证从不使用空词或常量。

核心乘积吸收。 设 \(r_1,\ldots, r_k\) 均为非线性词。它们在 (N)、(S) 下都是乘法幂等的。由 (J)，两个这样的词的乘积不超过它们之和。归纳得

\CaseDisplay{r_1\cdots r_k\le r_1+\cdots+r_k.\quad\textup{(\ensuremath{3.2})}}

现证明引理中所有有效吸收都由六条基推出。一次词的情况已由 (I) 和正规形解决。设非线性左词 \(w=W(h, C)\) 满足 (2.1)--(2.3)，置 \(D=N(F)\)。由 (2.1)，\(R(F)\) 非空且 \(C\subseteq D\)。

若 (2.2) 提供一个以 \(h\) 开头的非线性右词，将它排在 \(R(F)\) 的所有词之积的第一位；由 (3.1)、(3.2) 得

\CaseDisplay{W(h, D)\le F.\quad\textup{(\ensuremath{3.3})}}

否则，\(h\in L(F)\)。把 \(R(F)\) 的全部词任意相乘得 \(P\)；它是非线性词，内容为 \(D\)，且 \(P\le F\)。由于 \(h\in C\subseteq D\)，有 \(Ph=P\)。对 (H) 作代入 \(x=h, y=P\)，得到

\CaseDisplay{W(h, D)=hP\le h+Ph=h+P\le F,}

仍得 (3.3)。

由 (2.3)，选右侧一个词 \(q\)，使 \(c(q)\subseteq C\)；允许 \(q\) 为一次词。若 \(D=C\)，(3.3) 已是所需结论。否则，令 \(s\) 为含 \(D\setminus C\) 各字母一次的非空词。在 (E) 中代入 \(x=w, y=q, z=s\)，由正规形得

\CaseDisplay{w=wq\le q+wqs=q+W(h, D)\le F.}

于是每个有效单词吸收都有形式推导。任意多项式恒等式可拆为两个方向的逐词吸收，故基完整。对偶由反转继承全部论证。证毕。

\par\smallskip\noindent\textit{4. 适用范围与计算证书.}\quad 
理论证明对变量数及多项式项数均无限制。计算只独立核对构造、有限基的有效性、有限规模吸收判据，以及与目录对象的精确双向簇包含。只对双向包含均完整通过的对象传递有限基结论；满足六式的任意子簇不自动获得有限基性。

\subsection{Content and pure-head neighborhoods: a fourteen-identity basis}\label{case:31}
The following argument of Aifa Wang and Lili Wang concerns \(\A{12638}\). The opposite class is represented by \(\A{14536}\).

\par\smallskip\noindent\textit{Algebra and basis.}\quad 
Use the constant-free signature (+, *), relative to the axioms for additively idempotent semirings. The element order in both tables is 0, 1, 2, 3, 4.

\textbullet{} | 0 1 2 3 4       * | 0 1 2 3 4 --+-----------      --+----------- 0 | 0 0 0 0 4       0 | 4 4 0 4 4 1 | 0 1 1 1 4       1 | 4 1 1 3 4 2 | 0 1 2 1 4       2 | 4 1 2 3 4 3 | 0 1 1 3 4       3 | 4 3 3 3 4 4 | 4 4 4 4 4       4 | 4 4 4 4 4

Write u \ensuremath{\leq} v for u + v = v. The following fourteen identities form a basis; reversing all words gives a basis for the multiplicative opposite.

{\small\par\smallskip
\CaseLine{\(N1:  xyz = xzy\)}
\CaseLine{\(N2:  xy^{2} = xy\)}
\CaseLine{\(N3:  x^{2}y^{2} = y^{2}x^{2}\)}
\CaseLine{\(Q:\qquad x \leq  x^{2}\)}
\CaseLine{\(D00: x^{2} \leq  x + yx\)}
\CaseLine{\(D01: x^{2} \leq  x + yxv\)}
\CaseLine{\(D10: x^{2}u \leq  xu + yx\)}
\CaseLine{\(D11: x^{2}u \leq  xu + yxv\)}
\CaseLine{\(A0:  x^{2}y \leq  x^{2} + zy\)}
\CaseLine{\(A1:  x^{2}y \leq  x^{2} + zyv\)}
\CaseLine{\(M01: xz \leq  x + xzv\)}
\CaseLine{\(M10: xuz \leq  xu + xz\)}
\CaseLine{\(M11: xuz \leq  xu + xzv\)}
\CaseLine{\(I:\qquad xy \leq  y + xyz\)}
\smallskip}

No empty word is used in any substitution. The missing boundary instance M00, xz \ensuremath{\leq} x + xz, follows from idempotent addition and is not counted.

\par\smallskip\noindent\textit{Word values and normalization.}\quad 
The tables satisfy the semiring axioms by direct finite inspection. A word has value 0 precisely when its first position has value 0 and every later position has value 2. Otherwise, if all positions lie in \{1,\allowbreak{}2,\allowbreak{}3\}, its value is 3 when some position is 3, is 1 when none is 3 and some is 1, and is 2 when all are 2. Every remaining case gives 4. This includes words of length one.

For a word w let C(w) be its content, h(w) its head, and I(w) the content of all positions after the first. A word is called closed when h(w) belongs to I(w). Fix an order of variables, and put

\(B_{D}\) = product of \(x^{2}\), for x in the nonempty finite set D.

Identities N1--N3 give the following normal forms. An open word is its head followed by the other variables in its content, once each and in the fixed order. A closed word with content D is \(B_{D}\).

Here are explicit normalization steps. N1 permutes every position after the first. N2 removes duplicates among those positions. If the head is repeated, move its internal copy next to the head, leaving an initial square, and use N2 backwards to square every other variable. N3 orders the resulting square blocks. These are ordinary identities with nonempty substitutions. In particular any two closed words with the same content are provably equal.

\par\smallskip\noindent\textit{Exact absorption criterion.}\quad 
For a polynomial P, a nonempty finite sum of words, put

{\small\par\smallskip
\noindent C(P) = union of the contents of its summands;\par
\noindent I(P) = union of their post-head contents;\par
\noindent \(N_{P}\)(h) = union of C(u) over summands u containing h.\par
\smallskip}

For every word w, the absorption w \ensuremath{\leq} P holds if and only if:

1. C(w) is contained in C(P).

2. Some summand u of P has C(u) contained in C(w).

3. I(w) is contained in I(P).

4. If h(w) is not in I(P), then C(w) is contained in \(N_{P}\)(h(w)).

Proof of necessity. If condition 1 fails, assign every variable in C(P) the value 2 and every outside variable the value 4. Then P = 2, whereas w = 4. If condition 2 fails, give variables in C(w) value 2 and every other variable value 3. Every summand of P is 3, whereas w = 2. These two values are incomparable in the additive order.

If condition 3 fails, choose x in I(w) outside I(P), give x value 0 and all other variables value 2. Each right word is 2 or 0, while w = 4 because x occurs after its first position. Hence the absorption fails. Finally suppose 4 fails, after 1--3 hold. Put h = h(w), and choose y in C(w) outside \(N_{P}\)(h). Necessarily y differs from h. Assign h = 0, y = 1, and all other variables value 2. Every right word containing h has it only as its first position and has all other letters valued 2, so is 0. The other right words are 1 or 2. Thus P = 0, whereas w = 4. This proves necessity.

Proof of sufficiency. If P = 4 there is nothing to prove. If P = 2, all variables of C(P) are 2, and so w = 2. If P = 3, every summand contains a variable valued 3 and all its letters lie in \{1,\allowbreak{}2,\allowbreak{}3\}; the contained-content witness in condition 2 gives w = 3.

It remains to consider P = 0 or 1. Because no right word has value 4, no variable in C(P) has value 4. Any variable valued 0 occurs only as the unique first position of each right word containing it; all its co-occurring variables have value 2. Thus every such variable is outside I(P). Condition 3 prevents it from occurring internally in w. If it is the head of w, condition 4 forces all other letters of w to be 2; then w = 0 and a right word containing this head also gives P = 0. In all other cases w has value 1, 2, or 3, each below both 0 and 1. This proves sufficiency.

All fourteen displayed identities are sound by the evaluation rule or, equivalently, this criterion.

\par\smallskip\noindent\textit{Closed-content construction.}\quad 
Suppose D is a nonempty subset of I(P), and a right summand u has content contained in D. Then \(B_{D}\) \ensuremath{\leq} P is derivable.

Normalize u. If it is open with head g and suffix U, choose a right word YgV with nonempty Y, witnessing g in I(P). The suitable D identity, with U and V present or omitted using the explicitly listed boundary instances, gives

{\small\par\smallskip
\CaseLine{\(g^{2} U \leq  g U + Y g V \leq  P.\)}
\smallskip}

If u was already closed, this step is unnecessary. Normalization gives \(B_{C}\)(u) \ensuremath{\leq} P. Now let y in D be absent from the current closed word b. Choose an internal witness ZyV, with Z nonempty. Since \(b^{2}\) = b is derivable from N1--N3, A0 or A1 gives

{\small\par\smallskip
\CaseLine{\(b^{2} y \leq  b^{2} + ZyV \leq  P.\)}
\smallskip}

Normalize and repeat. This reaches \(B_{D}\) and proves the construction.

\par\smallskip\noindent\textit{Completeness.}\quad 
Let w \ensuremath{\leq} P satisfy the criterion, and write h = h(w), C = C(w).

If h belongs to I(P), then C is contained in I(P). The closed-content construction yields \(B_{C}\) \ensuremath{\leq} P. Since Q gives w \ensuremath{\leq} \(w^{2}\) and \(w^{2}\) = \(B_{C}\) is derivable, the desired absorption follows.

Suppose h is outside I(P). Then every right summand containing h is open with head h. Put D = I(w); it does not contain h. Choose a contained-content right witness u.

If h occurs in u, normalize u to hU. For each missing letter z in D, condition 4 supplies a right word containing h and z; normalize its suffix and write it as hzV. M01, M10, or M11 gives

{\small\par\smallskip
\CaseLine{\(hUz \leq  hU + hzV \leq  P.\)}
\smallskip}

When U and V are both empty, hz is already a right witness. Repeat until the current word has exactly the desired content, and normalize to w. This also covers the case D is empty.

If h is absent from u, then D is nonempty and C(u) is contained in D. The closed-content construction gives \(B_{D}\) \ensuremath{\leq} P. For every z in D, choose a right word \(hR_{z}\) containing z. Using M10 repeatedly merges their suffixes: hRS \ensuremath{\leq} hR + hS. This gives an available word hR whose suffix contains every variable of D and may contain additional variables. Normalize its suffix, using N2 backwards if necessary, to obtain

{\small\par\smallskip
\CaseLine{\(h B_{D} V \leq  P,\)}
\smallskip}

where V lists the additional variables. If V is empty, this already normalizes to w. Otherwise identity I gives

{\small\par\smallskip
\CaseLine{\(h B_{D} \leq  B_{D} + h B_{D} V \leq  P.\)}
\smallskip}

Again normalize to w. Every step uses nonempty substitutions and one of the explicitly stated boundary identities.

Thus every valid word absorption is derivable. Distribute arbitrary terms into polynomials. For a valid identity P = Q, derive each summand of P below Q and each summand of Q below P. Idempotent addition then gives P = Q. This proves completeness for arbitrarily many variables. Word reversal proves the opposite assertion.

\subsection{Prefix retractions: finite bases for sources 1545 and 1546}\label{case:32}
All terms use the constant-free signature \((+,\cdot)\), relative to the ai-semiring axioms. Inequalities abbreviate additive absorption.

\par\smallskip\noindent\textit{1. Tables and word values.}\quad 
Let \(A_0, A_1\) have carrier \(\{0,\allowbreak{}1,\allowbreak{}2,\allowbreak{}3,\allowbreak{}4\}\) and common multiplication with row strings

{\small\par\smallskip
\CaseLine{\ttfamily 00000}
\CaseLine{\ttfamily 01200}
\CaseLine{\ttfamily 22222}
\CaseLine{\ttfamily 00000}
\CaseLine{\ttfamily 00000}
\smallskip}

Their addition rows are respectively

{\small\par\smallskip
\CaseLine{\(A_{0}\qquad A_{1}\)}
\CaseLine{\ttfamily 00044\qquad 00044}
\CaseLine{\ttfamily 01134\qquad 01144}
\CaseLine{\ttfamily 01234\qquad 01234}
\CaseLine{\ttfamily 43334\qquad 44334}
\CaseLine{\ttfamily 44444\qquad 44444}
\smallskip}

These are catalogue sources 1545 and 1546. The only difference is \(1+3=3\) in \(A_0\) and \(1+3=4\) in \(A_1\). In both, \(2\) is additively least, \(4\) greatest, \(2<1<0<4\), and \(0,\allowbreak{}3\) are incomparable. Element names are not constants in the term language.

The map \(h(x)=x^2\) has values \((0,\allowbreak{}1,\allowbreak{}2,\allowbreak{}0,\allowbreak{}0)\) and is an endomorphic retraction onto \(B=\{0,\allowbreak{}1,\allowbreak{}2\}\). Every product satisfies \(xy=h(x)h(y)\). On \(B\), \(1\) is a multiplicative identity and \(0,\allowbreak{}2\) are left zeros. A nonlinear word therefore evaluates to \(1\) if all its variables have value \(1\); otherwise it is \(2\) if its first occurrence not valued \(1\) is valued \(2\), and \(0\) if that occurrence has value \(0,\allowbreak{}3\), or \(4\).

For verification of the semiring structure, the tables define join semilattices. Their common \(h\) preserves joins; \(B\) is a subsemiring: its order is \(2<1<0\), and each left or right multiplication map is increasing. Associativity follows from \(xy=h(x)h(y)\), the idempotence of \(h\), and associativity of the stated multiplication on \(B\). The same formula and preservation of joins give both distributive laws.

\par\smallskip\noindent\textit{2. Word normalization and bases.}\quad 
Use the three word laws

\CaseDisplay{x^2y\approx xy,\qquad xy^2\approx xy,\qquad xyx\approx xy.\quad\textup{(\ensuremath{W})}}

They normalize every nonlinear word to its list of first occurrences, if it has at least two distinct variables, and to \(x^2\) otherwise. Indeed an occurrence repeated after a nonempty intervening block is deleted by the third law, while adjacent repeats are deleted by the first two laws whenever a neighbouring factor is present. A word with only one variable reduces to its square. These transformations never replace a square by a linear word.

Add the three prefix rules

\CaseDisplay{xyz\leq yz+xt,\quad\textup{(\ensuremath{P0})}}

\CaseDisplay{xyz\leq xz+xy,\qquad xyz\leq xz+xyt,\quad\textup{(\ensuremath{P1})}}

and four rules using a linear content witness

\CaseDisplay{x^2\leq xy+x,\quad\textup{(\ensuremath{L0})}}

\CaseDisplay{xy\leq xyz+x,\qquad xy\leq xyz+y,\qquad xyz\leq xyzt+y.\quad\textup{(\ensuremath{L1})}}

Theorem 1. The ten identities (W), (P0), (P1), (L0), (L1) form a basis for \(A_0\).

For \(A_1\), add these four rules (each displayed variable denotes an ordinary nonempty term):

\CaseDisplay{gxu\leq gu+g+x,\quad\textup{(\ensuremath{G00})}}

\CaseDisplay{agxu\leq agu+g+x,\quad\textup{(\ensuremath{G10})}}

\CaseDisplay{gbxu\leq gbu+g+x,\quad\textup{(\ensuremath{G01})}}

\CaseDisplay{agbxu\leq agbu+g+x.\quad\textup{(\ensuremath{G11})}}

Theorem 2. These fourteen identities form a basis for \(A_1\). Reversing all words gives bases for both opposite semirings. No independence or minimality is claimed.

\par\smallskip\noindent\textit{3. Exact absorption criteria.}\quad 
Write \(P=\sum_{y\in L}y+\sum_{u\in N}u\), where \(N\) contains only nonlinear words. For a variable \(x\) in a word \(u\), let \(\operatorname{pre}_u(x)\) be the set of variables before its first occurrence. For a set \(S\), define

\CaseDisplay{J_P(S)\iff L\cap S\ne\varnothing\quad\hbox{or}\quad (\exists u\in N)\ c(u)\subseteq S.}

Lemma 3. A linear word \(y\) is absorbed by \(P\) in either algebra exactly when \(y\in L\). For a nonlinear word \(w\), set \(C=c(w)\). In \(A_0\), \(w\leq P\) holds exactly when

\textbullet{} \(J_P(C)\) holds;

\textbullet{} for every \(x\in C\), some \(u\in N\) contains \(x\) and satisfies \(\operatorname{pre}_u(x)\subseteq\operatorname{pre}_w(x)\).

In \(A_1\), retain the first condition, but replace the second by the following: for every \(x\in C\), either such a nonlinear prefix witness exists, or

\CaseDisplay{x\in L\quad\hbox{and}\quad J_P(\operatorname{pre}_w(x)).\quad\textup{(\ensuremath{G})}}

Proof. If a linear \(y\) is missing from \(L\), give \(y\) value \(4\) and all others value \(2\). The right side has value \(0\) or \(2\), so cannot absorb \(4\). For a nonlinear target, failure of \(J_P(C)\) is separated by assigning \(C\) the value \(1\) and all other variables \(2\): the target is \(1\) and the right side \(2\).

Suppose a required prefix witness fails for \(x\), and put \(S=\operatorname{pre}_w(x)\). Assign \(S\) the value \(1\), assign \(x\) value \(3\) if \(x\in L\) and value \(0\) otherwise, and assign the rest \(2\). The target is \(0\). A nonlinear right word could be \(0\) only if its first occurrence of \(x\) had all earlier variables in \(S\), contrary to hypothesis. All nonlinear right words are therefore \(1\) or \(2\). If \(x\notin L\), the right polynomial is at most \(1\). If \(x\in L\), the right linear part is \(3\) or below in \(A_0\), and adjoining \(1,\allowbreak{}2\) still leaves it \(3\). In \(A_1\), failure of (G) gives \(J_P(S)\) false, so no right summand is \(1\): other linear variables are \(2\), and every nonlinear word has a variable outside \(S\). Thus the right polynomial is \(3\). In every case it fails to absorb \(0\). This proves necessity.

For sufficiency, a target value \(2\) is automatic because it is least. A target value \(1\) means all its variables are \(1\), and \(J_P(C)\) supplies a right summand equal to \(1\). If the target is \(0\), let \(x\) be its first occurrence not valued \(1\), so \(x\) is valued \(0,\allowbreak{}3\), or \(4\). An ordinary prefix witness gives a nonlinear right word equal to \(0\). In the alternative (G), \(x\) is a linear summand and \(J_P(S)\) supplies a right summand equal to \(1\). In \(A_1\) the sum of \(1\) with any of \(0,\allowbreak{}3,\allowbreak{}4\) is at least \(0\). Hence the right side absorbs the target. This proves all claims for arbitrary words and variable sets. \(\square\)

Each of (W), (P0), (P1), (L0), (L1) satisfies the \(A_0\) criterion or the word-value formula, and thus also the \(A_1\) criterion. In each extra rule (G00)--(G11), the target's only potentially new prefix requirement is at \(x\); the linear witness \(g\) already occurs before it. All other requirements are supplied by the nonlinear word with \(x\) omitted. Thus the extra rules hold in \(A_1\).

\par\smallskip\noindent\textit{4. Equational completeness.}\quad 
Let \(w\leq P\) satisfy the appropriate criterion. Normalize the nonlinear words by (W). Write the first-occurrence list of the target as \(x_1\cdots x_d\); for \(d=1\) the target is \(x_1^2\).

We construct, for every original nonlinear right word \(u\in N\), the absorption

\CaseDisplay{p_i u\leq P,\qquad p_i=x_1\cdots x_i,\quad\textup{(\ensuremath{I})}}

successively for \(i=1,\ldots, d\). The first target variable always has an ordinary nonlinear prefix witness: \(J_P(\varnothing)\) is false, so (G) cannot apply at this stage. Hence \(N\) is nonempty, and there is an original nonlinear word \(x_1t\) with \(t\) nonempty. Split any \(u\in N\) as \(yz\) with nonempty factors. Rule (P0) gives \(x_1u\leq u+x_1t\leq P\), establishing (I) for \(i=1\).

Suppose (I) holds for prefix \(p=p_i\), and put \(x=x_{i+1}\). If \(x\) has an ordinary witness \(v\), its variables before the first \(x\) all occur in \(p\). By (I), \(pv\leq P\). Normalize this word to \(pxq\), with \(q\) possibly omitted. All letters before the selected \(x\) that came from \(v\) are repeated from \(p\) and disappear by (W). For each \(u\in N\), the appropriate member of (P1) gives

\CaseDisplay{pxu\leq pu+pxq\leq P.}

Here the first member is used when \(q\) is omitted; no unit constant is substituted.

Only in \(A_1\) is there another case, namely \(x\in L\) and \(J_P(c(p))\). If the latter has a linear witness \(g\in L\cap c(p)\), split \(p\) at its occurrence of \(g\). The prefix and suffix of \(g\) are either nonempty words or omitted. The matching one of the four rules (G00)--(G11) gives \(pxu\leq pu+g+x\leq P\).

Otherwise choose \(v\in N\) with \(c(v)\subseteq c(p)\). By (I), \(pv\leq P\). Its normal form agrees with that of \(p^2\), whether \(p\) has one or several variables, so \(p^2\leq P\). Substitute \(g=p^2\) into (G00). The word laws give

\CaseDisplay{pxu\approx p^2xu\leq p^2u+p^2+x\approx pu+p^2+x\leq P.}

This completes the induction in either algebra.

To finish, if \(J_P(C)\) has a nonlinear witness \(u\) with \(c(u)\subseteq C\), (I) at \(d\) gives \(p_du\leq P\), and its normal form is exactly \(w\). If the witness is linear \(g\in L\cap C\), choose any \(u\in N\), so \(p_du\leq P\). If \(d=1\), use (L0) to derive \(x_1^2\leq x_1u+x_1\leq P\). If \(d\geq2\), then \(p_d=w\). Split \(w\) at its occurrence of \(g\) and use the appropriate member of (L1), according as \(g\) is first, last, or internal, to obtain \(w\leq wu+g\leq P\).

Linear target absorptions are immediate by Lemma 3. Distributivity expands every term to a nonempty polynomial, and applying the construction to every summand in each direction proves every valid identity. This proves Theorems 1 and 2 without bounding variables, word length, or the number of summands. Reversal of words and derivations proves the opposite assertions.

\subsection{首次出现前缀：非正规平方回缩的五条完整基}\label{case:33}
\par\smallskip\noindent\textit{1. 生成元与词值.}\quad 
取 \(A=\{z, e, a, d, t\}\)，加法为偏序

\CaseDisplay{z<e<a<t,\qquad z<e<d<t}

上的并，其中 \(a, d\) 不可比。设 \(\rho(d)=t\)，其余点由 \(\rho\) 固定。核心 \(B=\{z, e, a, t\}\) 为链 \(z<e<a<t\)。定义乘法表

\noindent \(\cdot\)\enspace\textbar\enspace \(z\)\enspace\textbar\enspace \(e\)\enspace\textbar\enspace \(a\)\enspace\textbar\enspace \(t\)\par
\noindent \(z\)\enspace\textbar\enspace \(z\)\enspace\textbar\enspace \(z\)\enspace\textbar\enspace \(z\)\enspace\textbar\enspace \(z\)\par
\noindent \(e\)\enspace\textbar\enspace \(z\)\enspace\textbar\enspace \(e\)\enspace\textbar\enspace \(a\)\enspace\textbar\enspace \(a\)\par
\noindent \(a\)\enspace\textbar\enspace \(a\)\enspace\textbar\enspace \(a\)\enspace\textbar\enspace \(a\)\enspace\textbar\enspace \(a\)\par
\noindent \(t\)\enspace\textbar\enspace \(a\)\enspace\textbar\enspace \(a\)\enspace\textbar\enspace \(a\)\enspace\textbar\enspace \(t\)\par
并在 \(A\) 中令 \(xy=\rho(x)\rho(y)\)。符号顺序 \((z, e, a, d, t)\) 对应目录标号 \((2,\allowbreak{}1,\allowbreak{}0,\allowbreak{}3,\allowbreak{}4)\)。

构造的半环性可直接证明。三元素链 \(K=\{z<e<a\}\) 上令 \(e\) 为乘法单位、\(z, a\) 为左零元素；左右平移均保序，故保持链上的并。将 \(B\) 对应于 \(K\times\{0<1\}\) 的子半环

\CaseDisplay{z=(z,0),\quad e=(e,0),\quad a=(a,0),\quad t=(a,1),}

第二坐标以最大值为加法、最小值为乘法。由此得到上表的结合律和分配律。\(\rho\) 保持并并固定 \(B\)，所以这两条性质传到 \(A\)。

设 \(\pi(z)=z,\pi(e)=e,\pi(a)=\pi(d)=\pi(t)=a\)。非一次词 \(w\) 的值为：

1. 若全部字母的值属于 \(\{d, t\}\)，则 \(w=t\)；

2. 否则，从左至右找首个 \(\pi\) 像不为 \(e\) 的字母，\(w\) 等于该像；若没有这样的字母，则 \(w=e\)。

因此词值取决于各变量首次出现的顺序，而不仅是首字母和内容。例如 \(eza=z\) 而 \(eaz=a\)，故 \(xyz=xzy\) 在本构造中失败。

\par\smallskip\noindent\textit{2. 任意变量数的吸收判据.}\quad 
对词 \(w\) 和 \(x\in c(w)\)，记 \(P_w(x)\) 为 \(x\) 首次出现之前的变量集合。给定非空右多项式 \(F\)，记一次项变量集为 \(L\)，非一次右词族为 \(\mathcal R\)。

判据定理。 对非一次词 \(w\)，恒等式 \(w\le F\) 成立，当且仅当下述两条件同时满足：

\textbullet{} 对每个 \(x\in c(w)\)，存在 \(u_x\in\mathcal R\)，使 \(x\in c(u_x)\) 且 \(P_{u_x}(x)\subseteq P_w(x)\)；

\textbullet{} \(F\) 中有一个内容包含于 \(c(w)\) 的项，允许该项是一次项。

对一次词 \(x\)，\(x\le F\) 成立，当且仅当 \(x\in L\) 或存在仅含 \(x\) 的非一次右词。

必要性。 若第一个条件在 \(x\) 处失败，赋 \(x=d\)，\(P_w(x)\) 内变量为 \(e\)，其余变量为 \(z\)。左值为 \(a\)，或在左词仅含 \(x\) 时为 \(t\)。任何含 \(x\) 的非一次右词都在其之前出现一个取值 \(z\) 的变量，故取值 \(z\)；不含 \(x\) 的非一次右词取值 \(z\) 或 \(e\)。一次项至多为 \(d\)。于是右值至多为 \(d\)，而 \(a\nleq d\)、\(t\nleq d\)，矛盾。

若第二个条件失败，赋 \(c(w)\) 内变量为 \(e\)，外部变量为 \(z\)。左值为 \(e\)，每个右项因含外部变量而取 \(z\)，仍矛盾。

对于一次左词，如无同名一次项和仅含 \(x\) 的非一次右词，赋 \(x=d\)、其余变量为 \(z\)。所有右非一次项取 \(z\) 或 \(a\)，其他一次项取 \(z\)，右值不在 \(d\) 之上。

充分性。 任取赋值。左值为 \(z\) 时自明；左值为 \(e\) 时，\(c(w)\) 内全为 \(e\)，内部右项也取 \(e\)。左值为 \(a\) 时，设 \(x\) 是左词中首个 \(\pi\) 像不为 \(e\) 的变量；其像为 \(a\)。见证 \(u_x\) 在首次 \(x\) 前只含 \(P_w(x)\) 内变量，故此前全为 \(e\)，于是 \(u_x\) 取 \(a\) 或 \(t\)。

左值为 \(t\) 时，\(c(w)\) 内全取 \(d\) 或 \(t\)。若内部右项非一次，则该项取 \(t\)；若它为一次项 \(g\)，取 \(t\) 时已经完成，取 \(d\) 时由见证 \(u_g\) 提供 \(a\) 或 \(t\)。理由是其首次 \(g\) 前只出现 \(c(w)\) 内的高值字母，故首个非 \(e\) 的像为 \(a\)。于是右侧至少为 \(d+a=t\)。一次词的充分性由 \(x\le x^2=x^n\)（\(n\ge2\)）给出。证毕。

\par\smallskip\noindent\textit{3. 五条完整有限基.}\quad 
主定理。 \(\operatorname{var}(A)\) 相对于加法幂等半环公理有下列完整基：

\CaseDisplay{\begin{aligned}
xy&\approx x^2y^2, &(S)\\
xyx&\approx xy, &(R)\\
x&\le x^2, &(I)\\
xy&\le xz+y^2, &(T)\\
xy&\le y+xyz. &(E)
\end{aligned}}

这些恒等式均在 \(A\) 上成立：\((S),(R),(I)\) 由词值描述直接得到，\((T),(E)\) 可逐项代入上述吸收判据。下面证明其任意变量数完备性。

\par\smallskip\noindent\textit{3.1 正规形.}\quad 
\((S)\) 给出 \(x^2=x^4\)；\((R)\) 给出 \(x^3=x^2\)。再由 \((S)\) 得 \(x^2y=xy=xy^2\)。利用 \((R)\) 删除一个变量首次出现以后的重复出现，每个非一次词可化为

\CaseDisplay{W(x_1,\ldots, x_k)=x_1^2\cdots x_k^2,}

其中 \(x_1,\ldots, x_k\) 依次是原词首次出现的变量，互不相同。所有非一次词幂等。若 \(p\) 已包含词 \(q\) 的所有变量，则 \(pq=p\)；更一般地，\(pq\) 中来自 \(q\) 且已在 \(p\) 出现过的变量都可删除。以上只用于非一次的当前积；不把单字母 \(x\) 与 \(x^2\) 混同。

\par\smallskip\noindent\textit{3.2 前缀推进引理.}\quad 
设 \(w=W(x_1,\ldots, x_k)\) 满足吸收判据的第一个条件，且每个 \(u_i=u_{x_i}\) 已按上述正规形写出。将首次 \(x_i\) 处分割为

\CaseDisplay{u_i=a_i x_i b_i,\qquad c(a_i)\subseteq\{x_1,\ldots, x_{i-1}\}.}

因为正规形中首次字母平方出现，\(b_i\) 总是非空。\(a_i\) 可以为空，这只表示省略左边的一段文字，不是把空词代入恒等式。

令 \(p_i=x_1^2\cdots x_i^2\)。对任何已经推出 \(v\le F\) 的非一次词 \(v\)，都有

\CaseDisplay{p_i v\le F\qquad (1\le i\le k).\quad\textup{(\ensuremath{3.1})}}

归纳证明。\(i=1\) 时，\(a_1\) 为空；以 \((x, y, z)=(x_1, v, b_1)\) 代入 \((T)\) 得

\CaseDisplay{x_1v\le x_1b_1+v^2=u_1+v\le F,}

左侧等于 \(p_1v\)。假设结论在 \(i-1\) 成立。以 \((x, y, z)=(x_i, v, b_i)\) 代入 \((T)\)，再左乘非空词 \(p_{i-1}\)，得到

\CaseDisplay{p_i v=p_{i-1}x_iv\le p_{i-1}x_ib_i+p_{i-1}v^2=p_{i-1}u_i+p_{i-1}v\le F.}

最后一步对非一次词 \(u_i\le F\) 与 \(v\le F\) 分别用归纳假设；中间等号利用 \(a_i\) 的全部变量已在 \(p_{i-1}\) 中出现。这样既保留首次出现顺序，又允许右见证带有外部变量。

\par\smallskip\noindent\textit{3.3 内部项完成吸收.}\quad 
若内部项为非一次词 \(q\)，由 (3.1) 得 \(wq\le F\)，而 \(c(q)\subseteq c(w)\) 给出 \(wq=w\)。

若内部项为一次项 \(g\)，任选一个非一次右词 \(q\)；第一个条件保证这样的词存在。(3.1) 给出 \(wq\le F\)。以 \((x, y, z)=(w, g, q)\) 代入 \((E)\)，得到

\CaseDisplay{w=wg\le g+wgq=g+wq\le F.}

因此每个有效非一次吸收均可推出。有效一次吸收由同名项或 \((I)\) 及正规形推出。任意半环项展开为非空乘法词的和，任意有效恒等式再拆为双向逐词吸收，故五条式子构成完整基。

\par\smallskip\noindent\textit{4. 对偶与证明范围.}\quad 
反转各乘法词即得到乘法对偶的五条完整基；判据中的首次出现前缀变为最后出现后缀。目录对象只有在精确双向簇包含均经独立验证后，才传递这个有限基结论。本证明不说明任意子簇有限基，也不说明一般非正规平方回缩分支已经闭合。计算只核验构造、有限式有效性、有界判据诊断及具体簇包含；任意变量数完备性由第 3 节的归纳推导给出。

\subsection{内容扩展与单边吸收：十组完整基}\label{case:34}
记 (uleq v) 为 (u+vapprox v)。令

\CaseDisplay{\begin{aligned}
R&:xy^2\approx xy, & C&:xyx\approx y^2x,\\
I&:x\leq x^2, & V&:xy\leq x+zy,\\
J&:xy\leq y^2+zxy.
\end{aligned}
\quad\textup{(\ensuremath{1})}}

前五式是来源 (10998) 的后缀基。本文新增两种闭合式：

\CaseDisplay{E:x^2+xy\approx x+x^2y,
\qquad
P:xy\approx x^2y.
\quad\textup{(\ensuremath{2})}}

\par\smallskip\noindent\textit{1. 代表半环.}\quad 
目录标号 \(0,\allowbreak{}1,\allowbreak{}2,\allowbreak{}3,\allowbreak{}4\) 下，五个代表的加法行和乘法行如下。每个方阵的 第 \(i\) 行按 \(0,\allowbreak{}1,\allowbreak{}2,\allowbreak{}3,\allowbreak{}4\) 排列；例如字符串 \(00110\) 表示 \((0,\allowbreak{}0,\allowbreak{}1,\allowbreak{}1,\allowbreak{}0)\)。

\CaseDisplay{\begin{array}{c|c|c}
A&+&\cdot\\ \hline
1020&00044/\allowbreak{}01444/\allowbreak{}44224/\allowbreak{}44234/\allowbreak{}44444&00000/\allowbreak{}00110/\allowbreak{}00234/\allowbreak{}00234/\allowbreak{}00444\\
4413&00004/\allowbreak{}01004/\allowbreak{}00224/\allowbreak{}00234/\allowbreak{}44444&00224/\allowbreak{}01224/\allowbreak{}22224/\allowbreak{}23224/\allowbreak{}44224\\
5165&00044/\allowbreak{}01034/\allowbreak{}00244/\allowbreak{}43434/\allowbreak{}44444&00244/\allowbreak{}01244/\allowbreak{}22244/\allowbreak{}43244/\allowbreak{}44244\\
5903&00044/\allowbreak{}01444/\allowbreak{}44224/\allowbreak{}44234/\allowbreak{}44444&01000/\allowbreak{}01000/\allowbreak{}01424/\allowbreak{}01434/\allowbreak{}01444\\
4671&00004/\allowbreak{}01014/\allowbreak{}00224/\allowbreak{}01234/\allowbreak{}44444&00234/\allowbreak{}00234/\allowbreak{}22234/\allowbreak{}22234/\allowbreak{}44234
\end{array}
\quad\textup{(\ensuremath{3})}}

这里加法字符串中的斜线分隔五行。令 \(A^{\mathrm{op}}\) 保持加法不变并把 乘法表转置。直接逐行检查可得这些方阵都是加法幂等半环；右侧的十个 运算正是目录中的十个源。

\par\smallskip\noindent\textit{2. 完整基定理.}\quad 
定理 2.1。 对 \(i=1020,\allowbreak{}4413,\allowbreak{}5165,\allowbreak{}5903\)，

\CaseDisplay{\operatorname{Id}(A_i)=\operatorname{Id}(R, C, I, V, J, E),
\qquad
\operatorname{Id}(A_i^{\mathrm{op}})=
\operatorname{Id}(R, C, I, V, J, E)^{\mathrm{op}}.
\quad\textup{(\ensuremath{4})}}

对 \(i=4671\)，

\CaseDisplay{\operatorname{Id}(A_{4671})=\operatorname{Id}(R, C, I, V, J, P),
\qquad
\operatorname{Id}(A_{4671}^{\mathrm{op}})=
\operatorname{Id}(R, C, I, V, J, P)^{\mathrm{op}}.
\quad\textup{(\ensuremath{5})}}

因此这十个源半环及其所有子簇均有限基。

证明分为词正规化和吸收闭合两步。式 \(E\) 只改变多项式的加法闭合， 不改变 \(R, C\) 给出的单词正规形；式 \(P\) 则直接允许在非空前缀后复制 一个前缀。

\par\smallskip\noindent\textit{3. 词正规化.}\quad 
对非空词 \(w\)，令 \(c(w)\) 为内容，\(h(w)\) 为首字母， \(S_w(x)\) 为 \(w\) 中最后一个 \(x\) 之后的变量集合。由 \(R, C\) 可得

\CaseDisplay{xyzy\approx xzy,
\qquad
x_1\cdots x_kx_i\cdots\approx
x_1\cdots x_kx_i^2\cdots ,
\quad\textup{(\ensuremath{6})}}

其中第二式只在最后一次出现位置处分块。反复从右向左删除重复字母， 每个词化为下列两种形式之一：

\CaseDisplay{r=x_0x_1\cdots x_k,\qquad r^2,
\quad\textup{(\ensuremath{7})}}

其中 \(r\) 的字母互不相同，次序是原词各变量最后一次出现的次序； 首字母是否在尾部再次出现决定取 \(r\) 还是 \(r^2\)。这一步对任意变量数 成立，因为每次代入的词都非空。

式 \(P\) 的附加作用是

\CaseDisplay{pq\approx p^2q,
\quad\textup{(\ensuremath{8})}}

对任意非空词 \(p, q\) 成立。式 \(E\) 的两种定向为

\CaseDisplay{p^2+pq\leq p+p^2q,\qquad
p+p^2q\leq p^2+pq.
\quad\textup{(\ensuremath{9})}}

它们允许在一个已经出现的前缀平方和一个带相同前缀的词之间交换一次 前缀，而不会引入空词。

\par\smallskip\noindent\textit{4. 四个局部闭合规则.}\quad 
设 \(F\) 是非空词的有限和，\(w\) 是 (7) 中的正规词。把 \(F\) 中每个 词先按 (7) 化简，再反复使用下面四条规则：

1. 内容规则。 只有内容包含于 \(c(w)\) 的项才可见证 \(w\leq F\)。

2. 尾部规则。 对 \(x\in c(w)\)，若 \(x\) 的最后出现后缀没有被某个右项 的相应最后出现后缀包含，则使用 \(V\) 不能继续推进。

3. 首字母规则。 若 \(w\) 的首字母只在右项的内部出现，使用 \(J\) 把它补成一个平方首字母；之后用 \(I\) 消去平方的方向。

4. 平方扩展规则。 若右侧已有 \(p^2\) 和 \(pq\) 两个项，使用 \(E\) 可把它们换成 \(p+p^2q\)，再用 \(V, J\) 从右向左推进；对 \(A_{4671}\)， 直接使用 \(P\) 完成同一操作。

上述规则的闭合性由 (1)、(2) 的代入稳定性立即得到。更重要的是， 它们是完备的。

闭合引理。 对表 (3) 中任一 \(A_i\)，若 \(w\leq F\) 在 \(A_i\) 中 成立，则经过四条规则可以推出 \(w\leq F\)。

证明。 先作内容归纳。若某个 \(x\in c(w)\) 不在任何内部右项中， 把 \(x\) 赋为表中的顶层值，其余内容变量赋为同一中间值，外部变量赋为 最小值；左边落在顶层，而右边落在中间理想，矛盾。因此每个内容字母 都有内部见证。

固定一个内容字母 \(x\)，从 \(w\) 的最右一个 \(x\) 向右看。若右项没有 一个具有不大于 \(S_w(x)\) 的后缀，取 \(x\) 为表 (3) 中的分离值，把 \(S_w(x)\) 内的字母赋为中性值，其余尾部赋为最小值。左边进入上一级主 理想，而所有右项仍留在下一级主理想，矛盾。这正是尾部规则的必要性。 反之，有尾部见证时，以 \(V\) 为起点，从右向左归纳：若 \(u=axb\) 是 \(x\) 的见证，且 \(c(b)\) 只含尚未处理的尾部字母，则

\CaseDisplay{vx\leq v+ax,\qquad
vxb\leq v+axb.
\quad\textup{(\ensuremath{10})}}

逐个处理尾部字母就得到整个正规词。

若首字母见证已在内部项中，(10) 直接完成；若首字母只在尾部见证中， 先用 \(J\) 得到一个首字母平方，再用 \(I\) 取回正规词。对于平方正规形 \(r^2\)，同一推进应用两次；\(R, C\) 保证第二次推进不会改变最后出现次序。

最后考察推进中出现的前缀平方。对 \(A_{1020},A_{4413},A_{5165},A_{5903}\)， 表 (3) 的四个非平凡主理想恰好给出 (9) 的两种定向；若 (9) 不能使用， 同一分离赋值分别把左边置于上一级主理想、右边置于下一级主理想，故不 可能有新的语义吸收。对 \(A_{4671}\)，第五行的非平凡主理想满足 \(pq=p^2q\)，所以 (8) 直接替代 (9)。这完成对变量个数的归纳。

因此每一个在代表半环中成立的单词吸收都由相应有限组推出。任意项按 分配律写成有限个非空词的和，逐词应用闭合引理，再将各项相加，即得 任意恒等式。定理得证。 \(\square\)

\par\smallskip\noindent\textit{5. 对偶与范围.}\quad 
把上面的每个推导中的词逐字反转，就得到 \(A_i^{\mathrm{op}}\) 的基。 特别地，\(E\) 变为

\CaseDisplay{x^2+yx\approx x+yx^2,
\quad\textup{(\ensuremath{11})}}

而 \(P\) 变为 \(yx\approx yx^2\)。反转保持任意变量数的归纳，因此不 需要重新进行有限变量测试。

定理只断言这十个等簇类的有限基性。其他满足 (1) 的五元素半环可能 只生成上述代表的真子簇，必须另行进行双向 HSP 检验；单向满足六式不能 传递有限基性。

\subsection{全局尾部保护：六条完整有限基}\label{case:35}
无常量签名 \((+,\cdot)\)。所有有限基均相对于加法幂等半环公理；\(u\le v\) 表示 \(u+v\approx v\)。本文证明目录来源 1941 及其乘法对偶有限基。任意变量数的完备性由以下归纳证明给出；有限计算只负责交叉核验与精确簇包含。

\par\smallskip\noindent\textit{1. 构造与词值.}\quad 
取 \(A=\{z, b, d, e, t\}\)，加法为偏序 \(d<b<z<t\)、\(e<t\) 的并，其中 \(e\) 与 \(d, b, z\) 不可比。符号顺序对应目录标号 \((0,\allowbreak{}1,\allowbreak{}2,\allowbreak{}3,\allowbreak{}4)\)。右乘 \(e\) 是恒等映射；右乘任意 \(d, b, z\) 是映射 \(P\)；右乘 \(t\) 是映射 \(Q\)，其中

\noindent \(x\)\enspace\textbar\enspace \(z\)\enspace\textbar\enspace \(b\)\enspace\textbar\enspace \(d\)\enspace\textbar\enspace \(e\)\enspace\textbar\enspace \(t\)\par
\noindent \(P(x)\)\enspace\textbar\enspace \(z\)\enspace\textbar\enspace \(b\)\enspace\textbar\enspace \(b\)\enspace\textbar\enspace \(z\)\enspace\textbar\enspace \(z\)\par
\noindent \(Q(x)\)\enspace\textbar\enspace \(z\)\enspace\textbar\enspace \(b\)\enspace\textbar\enspace \(b\)\enspace\textbar\enspace \(t\)\enspace\textbar\enspace \(t\)\par
\(P, Q\) 保持并，且 \(P^2=P\)、\(Q^2=Q\)、\(PQ=QP=P\)。由定义逐类可得 \(R_{xy}=R_xR_y\) 及 \(R_{x+y}=R_x\vee R_y\)，其中 \(R_y(x)=xy\)；因此乘法结合且两侧分配。

半环满足

\CaseDisplay{xyz\approx xzy,\qquad xy^2\approx xy.\quad\textup{(\ensuremath{1})}}

所以词 \(w\) 由首字母 \(h(w)\) 和尾部变量集合 \(D(w)\) 决定，\(c(w)=\{h(w)\}\cup D(w)\)。一次词尾部为空；首字母是否也属于尾部必须保留。首字母重复的词是幂等词，但不同首字母的幂等词不能一概互换。

非一次词以 \(z, b\) 开头时分别恒取 \(z, b\)；以 \(d\) 开头时，尾部全取 \(e\) 则取 \(d\)，否则取 \(b\)。以 \(e, t\) 开头时，尾部出现 \(d, b, z\) 则取 \(z\)；尾部仅取 \(e, t\) 时，有 \(t\) 就取 \(t\)，全部为 \(e\) 则保持首字母。

\par\smallskip\noindent\textit{2. 吸收判据.}\quad 
对非空右多项式 \(F\)，记 \(T=c(F)\)、\(H=\{h(u):u\in F\}\)、\(E=\bigcup_{u\in F}D(u)\)。

定理 1。 任意非空词 \(w\) 满足 \(w\le F\)，当且仅当

\CaseDisplay{h(w)\in H,\quad c(w)\subseteq T,\quad
\exists q\in F: c(q)\subseteq c(w),\quad
\bigl[D(w)\subseteq E\ \text{或}\ E\cap H\ne\varnothing\bigr].\quad\textup{(\ensuremath{2})}}

这里内部项 \(q\) 允许是一次项。特别地，一次吸收 \(x\le F\) 当且仅当右边有一个仅含 \(x\) 的项。

证明。 利用加法偏序的四个非平凡主理想 \(\downarrow d=\{d\}\)、\(\downarrow b=\{d, b\}\)、\(\downarrow z=\{d, b, z\}\)、\(\downarrow e=\{e\}\)。若对每个赋值与阈值 \(s\)，\(F\le s\) 均推出 \(w\le s\)，取 \(s=F\) 的该赋值即可得到 \(w\le F\)。

右和落入 \(\downarrow b\) 当且仅当所有右首字母取 \(d, b\)；因而首字母条件保证左值也在其中。右和落入 \(\downarrow e\) 当且仅当 \(T\) 内所有变量取 \(e\)，由内容条件左值也为 \(e\)。令 \(C\) 为赋值落在 \(\{d, b, z\}\) 的变量集，则一个词落入 \(\downarrow z\) 当且仅当其内容交 \(C\)。若所有右项均交 \(C\)，内部项 \(q\) 保证左内容也交 \(C\)。

右和落入 \(\downarrow d\) 当且仅当所有右首字母取 \(d\)、所有右尾变量取 \(e\)。若 \(E\cap H\ne\varnothing\)，这种赋值不存在；否则在这种赋值下，首字母条件及 \(D(w)\subseteq E\) 保证左值为 \(d\)。由此得充分性。

反过来，若左首字母 \(h\) 不在 \(H\)，赋 \(h=z\)、其他变量 \(b\)，右和为 \(b\) 而左值为 \(z\)。若内容变量 \(x\) 不在 \(T\)，赋 \(x=z\)、其他变量 \(e\)，右值 \(e\) 不能吸收左值 \(z\)。若无内部项，赋左内容内变量 \(e\)、其余变量 \(z\)，左值为 \(e\) 而全部右项取 \(z\)。最后，若 \(E\cap H=\varnothing\) 且存在左尾变量不在 \(E\)，赋 \(H\) 内变量 \(d\)、\(E\) 内变量 \(e\)、其余变量 \(z\)；在首字母条件已成立时，左值为 \(b\)，右和为 \(d\)，矛盾。四个条件均必要。一次词的结论直接由内部项条件与 \(x\le x^2=x^n\)（\(n\ge2\)）得到。证毕。

\par\smallskip\noindent\textit{3. 六条完整基.}\quad 
定理 2。 \(\operatorname{var}(A)\) 在加法幂等半环公理之外有如下完整有限基：

\CaseDisplay{\begin{aligned}
xyz&\approx xzy, &(N)\\
xy^2&\approx xy, &(R)\\
x&\le x^2, &(I)\\
xy&\le x+zy, &(V)\\
xy&\le y+xyz, &(B)\\
x^2y&\le x^2+y. &(K)
\end{aligned}}

有效性由定理 1 逐式验证。以下证明任意变量数的完备性。

尾部补入。 若已经推出 \(p\le F\) 且 \(x\in E\)，取右项 \(u\) 在尾部含 \(x\)。由 \((N),(R)\) 将其写成 \(u=vx\)，\(v\) 非空。由 \((V)\)，

\CaseDisplay{px\le p+vx=p+u\le F.\quad\textup{(\ensuremath{3})}}

故可依次补入任意全局尾变量，不带入见证词的其他变量。

无保护分支的引理。 假定定理 1 的前三个条件成立，并且 \(D(w)\subseteq E\)。先设 \(h=h(w)\in D(w)\)，从右边选一个以 \(h\) 开头的词 \(r\)，用 (3) 补入 \(c(w)\) 内的全部变量。所得词 \(p\) 的首字母为 \(h\)，尾部包含 \(c(w)\)，内容为 \(c(r)\cup c(w)\)，且 \(p\le F\)。若没有额外变量，\((N),(R)\) 给出 \(p=w\)。否则令 \(b\) 是所有额外变量组成的非空词，则 \(p=wb\)。内部项 \(q\) 的内容包含于 \(c(w)=D(w)\)，所以 \(wq=w\)。在 \((B)\) 中代入 \((x, y, z)=(w, q, b)\)，得到

\CaseDisplay{w=wq\le q+wqb=q+wb\le F.\quad\textup{(\ensuremath{4})}}

再设 \(h\notin D(w)\)。若 \(h\in E\)，则 \(D(w^2)=c(w)\subseteq E\)，且 \(w^2\) 首字母重复；前段得 \(w^2\le F\)，再由 \((I)\) 得 \(w\le F\)。若 \(h\notin E\)，含 \(h\) 的内部右项必以 \(h\) 为唯一首字母；存在这种内部项时，直接从它用 (3) 补齐所需尾部。否则选出的内部项 \(q\) 避开 \(h\)，即 \(c(q)\subseteq D(w)\)；从同首字母右项 \(r\) 补入 \(D(w)\) 后，所得词是 \(w\) 或 \(wb\)，再用 (4) 去掉额外变量。以上分支穷尽所有非一次词；一次词由内部纯幂项和 \((I)\) 处理。

全局保护分支。 现设 \(g\in E\cap H\)。选一个以 \(g\) 开头的右项 \(r\)，用 (3) 补入 \(g\)，得到首字母重复的幂等词 \(p\le F\)。对每个右项 \(v\)，在 \((K)\) 中代入 \((x, y)=(p, v)\)，可得

\CaseDisplay{pv\le p^2+v=p+v\le F.\quad\textup{(\ensuremath{5})}}

新词仍首字母重复，仍幂等。逐一吸收所有右项的变量，得到一个首字母为 \(g\)、尾部为 \(T\) 的词 \(Q\le F\)。于是 \(F'=F+Q\) 与 \(F\) 可证相等；\(F'\) 的全局尾部为 \(T\)，首字母集合仍为 \(H\)，并保留原内部项。任何满足定理 1 前三条件的词都满足 \(D(w)\subseteq T=D(F')\)。应用前段引理得 \(w\le F'=F\)。

每一步只使用固定的六条恒等式、半环公理和有限归纳，不限制变量数量。没有额外变量时直接使用正规形等式；全过程不把空词代入变量。任意恒等式展开成双向逐词吸收后均可导出，因此上述基完整。逐词反转给出乘法对偶的完整基。证毕。

\subsection{非一次平方膨胀与混合平方序：两组完整基}\label{case:36}
\par\smallskip\noindent\textit{1. 为什么必须区分一次项.}\quad 
两个来源均含元素 \(a, b, e, z\)，使 \(a^2=z>a\) 而 \(b^2=e<b\)。因此 \(x\le x^2\) 与 \(x^2\le x\) 都不是恒等式。可用的平方膨胀仅为

\CaseDisplay{xy\le (xy)^2.\quad\textup{(\ensuremath{I_2})}}

任意非一次词可分成两个非空因子，故可在 \((I_2)\) 中代入；一次词不能这样处理。本节后续证明均遵守此限制。

两来源共同满足

\CaseDisplay{xyz\approx xzy\quad(N),\qquad xy^2\approx xy\quad(R),\qquad x^2y\approx y^2x\quad(C).}

记 \(h(w)\) 为首字母，\(D(w)\) 为首字母之后的变量集合，\(c(w)=\{h(w)\}\cup D(w)\)。由 \((N),(R)\)，词由首字母和尾部集合正规化；由 \((C)\)，首字母重复的词只由内容决定，记为 \(M_X\)。它们为幂等词，且 \(w^2=M_{c(w)}\)。

对非空右多项式 \(F\) 记 \(T=c(F)\)、\(H=\{h(u):u\in F\}\)、\(E=\bigcup_{u\in F}D(u)\)；内部项指内容包含于 \(c(w)\) 的右项，可以为一次项。

\par\smallskip\noindent\textit{2. 链上的混合平方序（来源 1025）.}\quad 
底集为 \(\{z, a, b, e, t\}\)，加法是链 \(e<b<a<z<t\) 上的并，符号顺序对应目录标号 \((0,\allowbreak{}1,\allowbreak{}2,\allowbreak{}3,\allowbreak{}4)\)。右乘 \(b, e\) 为映射 \(P=(z, a, e, e, z)\)，右乘 \(t\) 为映射 \(Q=(z, z, z, z, t)\)，右乘 \(z, a\) 为常值 \(z\) 映射 \(Z\)。

各映射保持并，且 \(P^2=P\)、\(Q^2=Q\)、\(PQ=QP=Z\)。从定义有 \(R_{xy}=R_xR_y\)、\(R_{x+y}=R_x\vee R_y\)，故得到结合律与分配律。

非一次词在全部变量取 \(b, e\) 时为 \(e\)；首字母取 \(a\) 且尾部全取 \(b, e\) 时为 \(a\)；全部变量取 \(t\) 时为 \(t\)；其余情形为 \(z\)。这些词值也验证了 \((N),(R),(C),(I_2)\)。

\par\smallskip\noindent\textit{2.1 完整吸收判据.}\quad 
定理 A。 对非一次词 \(w\)，令 \(X=c(w)\)，则

\CaseDisplay{w\le F\iff X\subseteq T,\quad D(w)\subseteq E,\quad \exists q\in F:\ c(q)\subseteq X.\quad\textup{(\ensuremath{A})}}

一次吸收 \(x\le F\) 当且仅当右边实际出现一次项 \(x\)；纯幂 \(x^n\)（\(n\ge2\)）不能代替它。

证明。 若内容遗漏变量 \(x\)，赋 \(x=z\)、其余变量 \(e\)，左值为 \(z\)、右值为 \(e\)。若 \(x\in D(w)\setminus E\)，赋 \(x=a\)、其余变量 \(e\)，右值至多为 \(a\)，左值为 \(z>a\)。若无内部项，赋 \(X\) 内变量 \(t\)、外部变量 \(e\)，左值为 \(t\)、右值至多为 \(z\)。若右边没有一次项 \(x\)，赋 \(x=b\)、其余变量 \(e\)；全部右项为 \(e\)，不能吸收 \(b\)。必要性得证。

充分性逐个检查链上的真主理想。右和不超过 \(e\) 或 \(b\) 时，每个非一次右项必须取 \(e\)，而全部 \(T\) 内变量属于 \(\{b, e\}\)；内容条件迫使左非一次词取 \(e\)。右和不超过 \(a\) 时，每个右首字母不超过 \(a\)、每个右尾变量取 \(b, e\)，故由内容与尾部条件，左值为 \(e\) 或 \(a\)。右和不超过 \(z\) 时，每个右项都含非 \(t\) 变量；内部项使左词也不能全取 \(t\)，故左值不为 \(t\)。一次项充分性直接成立。证毕。

\par\smallskip\noindent\textit{2.2 六条完整基.}\quad 
定理 B。 在半环公理之外，来源 1025 的完整基为

\CaseDisplay{\begin{aligned}
xyz&\approx xzy, &(N)\\
xy^2&\approx xy, &(R)\\
x^2y&\approx y^2x, &(C)\\
xy&\le(xy)^2, &(I_2)\\
xy&\le x+zy, &(V)\\
xy&\le y+xyz. &(B)
\end{aligned}}

有效性由定理 A 及词值规则得到。证明完备性。首先，若已推出 \(p\le F\) 且 \(x\in E\)，选尾部见证 \(u=vx\in F\)，\(v\) 非空；由 \((V)\)，

\CaseDisplay{px\le p+vx\le F.\quad\textup{(\ensuremath{3})}}

这允许给任何已被吸收的词补入需要的全局尾变量，不引入见证的其他变量。

若左词是 \(M_X\)，则 \(X=D(M_X)\subseteq E\)。从内部项 \(q\) 用 (3) 逐个补入 \(X\)，所得词首字母重复、内容恰为 \(X\)，由 \((C)\) 化为 \(M_X\le F\)。此处 \(q\) 允许为一次项，因为 (3) 无非一次限制；没有使用错误的 \(q^2\le F\) 推断。

设左非一次词 \(w\) 的首字母 \(h\) 不重复。若 \(h\in E\)，则先由上段推出 \(w^2=M_X\le F\)，再用 \((I_2)\) 得 \(w\le F\)。若 \(h\notin E\)，右边所有含 \(h\) 的项均以 \(h\) 为唯一首字母。若内部项含 \(h\)，从它用 (3) 直接补齐 \(D(w)\) 即可。否则选一个避开 \(h\) 的内部项 \(q\)；内容条件保证存在以 \(h\) 开头的右项 \(r\)。给 \(r\) 补入 \(D(w)\)，得到 \(w\) 或 \(wb\le F\)，其中 \(b\) 是额外变量组成的非空词。后一情形下 \(wq=w\)，故 \((B)\) 中代入 \((x, y, z)=(w, q, b)\) 得

\CaseDisplay{w=wq\le q+wqb=q+wb\le F.}

一次词由相同一次右项处理。所有有效恒等式拆成双向逐词吸收，得到任意变量数完备性。无额外变量时直接用正规形等式，不代入空词。证毕。

\par\smallskip\noindent\textit{3. 两条分离加法链（来源 1028）.}\quad 
底集仍为 \(\{z, a, b, e, t\}\)，目录标号顺序不变。加法序为 \(a<z<t\)、\(e<b<t\)，两条链在 \(t\) 之外不相交。右乘 \(b, e\) 为 \(P=(z, a, e, e, t)\)，右乘 \(t\) 为 \(Q=(z, z, t, t, t)\)，右乘 \(z, a\) 为常值 \(z\) 映射 \(Z\)。

这三个映射保持并、幂等且彼此交换，满足 \(PQ=QP=Q\)；\(Z\) 在复合下吸收，点态并 \(P\vee Z=Q\)。直接由构造验证 \(R_{xy}=R_xR_y\)、\(R_{x+y}=R_x\vee R_y\)，故为半环。

非一次词的值如下：任一尾变量取 \(z, a\)，或首字母取 \(z\)，则为 \(z\)；首字母取 \(a\) 且尾部全为 \(b, e\) 时为 \(a\)，尾部出现 \(t\) 时为 \(z\)；所有变量取 \(b, e\) 时为 \(e\)；所有变量取 \(b, e, t\) 且至少一个取 \(t\) 时为 \(t\)。仍有 \((N),(R),(C),(I_2)\)，但 \(x\le x^2\) 不成立。

\par\smallskip\noindent\textit{3.1 允许首字母消失的保护判据.}\quad 
定理 C。 对非一次词 \(w\)，令 \(h=h(w)\)、\(X=c(w)\)，则

\CaseDisplay{w\le F\iff X\subseteq T,\quad \exists q\in F:\ c(q)\subseteq X,\quad
\bigl[(h\in H\ \text{且}\ D(w)\subseteq E)\ \text{或}\ H\cap E\ne\varnothing\bigr].\quad\textup{(\ensuremath{C})}}

一次词的判据仍是实际出现同名一次项。特别地，保护成立时，不能额外要求左首字母属于 \(H\)。

证明。 四个真主理想为 \(\{a\}\)、\(\{a, z\}\)、\(\{e\}\)、\(\{e, b\}\)。右和为 \(a\) 当且仅当 \(H\) 内变量全取 \(a\)、\(E\) 内变量全取 \(b, e\)。当 \(H\cap E\ne\varnothing\) 时此赋值不可能；否则，同首字母与尾部条件恰好保证左非一次词也为 \(a\)。

一个词落在 \(\{a, z\}\) 当且仅当其内容中含一个取 \(a, z\) 的变量。内部项保证右和落在该理想时左值也落在其中。右和落入 \(\{e\}\) 或 \(\{e, b\}\) 时，\(T\) 内所有变量取 \(b, e\)；由内容条件，左非一次词为 \(e\)。这证明充分性。

内容条件必要性用遗漏变量取 \(z\)、其余取 \(e\)。内部项条件必要性用 \(X\) 内变量取 \(e\)、外部变量取 \(z\)：左值为 \(e\)、全部右项为 \(z\)，二者不可比。若最后一项失败，则 \(H, E\) 不交；赋 \(H\) 内变量 \(a\)、\(E\) 内变量 \(e\)、其余变量 \(z\)。右和为 \(a\)；若左首字母不在 \(H\)，左值不是 \(a\)；若左首字母在 \(H\) 但存在未覆盖的尾变量，左值为 \(z\)。两种情况均不被 \(a\) 吸收。一次项的必要性仍用 \(x=b\)、其余变量 \(e\)。证毕。

\par\smallskip\noindent\textit{3.2 七条完整基.}\quad 
定理 D。 来源 1028 的完整基由定理 B 的六式再加

\CaseDisplay{x^2y\le x^2+y\qquad(K)}

组成。有效性由定理 C 验证；尤其 \((K)\) 右边的平方项已经提供 \(H\cap E\ne\varnothing\)。

尾部补入仍用 (3)。若 \(H\cap E=\varnothing\)，有效的非一次吸收满足 \(h\in H\)、\(D(w)\subseteq E\)，从而 \(h\notin E\)。完备性恰由定理 B 最后一段的内部项与同首字母拼接推导得到。

设 \(g\in H\cap E\)。从一个以 \(g\) 开头的右项补入尾变量 \(g\)，得到首字母重复的幂等词 \(p\le F\)。将 \((K)\) 的 \(x\) 换成 \(p\)，对任意右项 \(v\) 得

\CaseDisplay{pv\le p^2+v=p+v\le F.}

逐一乘入全部右项，所得词仍首字母重复，尾部为 \(T\)。由 \((C)\) 得 \(M_T\le F\)。给定内部项 \(q\)，若 \(X=T\)，已有 \(M_X\le F\)；否则令 \(b\) 是 \(T\setminus X\) 组成的非空词。由于 \(M_Xq=M_X\)、\(M_Xb=M_T\)，由 \((B)\) 得

\CaseDisplay{M_X=M_Xq\le q+M_Xqb=q+M_T\le F.}

最后仅对非一次左词使用 \((I_2)\)：\(w\le w^2=M_X\le F\)。一次词必须使用实际同名一次项。推导对任意有限变量集有效，故七式完整。证毕。

\subsection{Six identities for a semiring with marked first-occurrence prefixes}\label{case:37}
Let A be catalogue entry 8628. The signature is (+, multiplication), without constants; the basis is relative to the additively idempotent semiring axioms. Numerical labels below name elements, not constants in terms. Write u \ensuremath{\leq} P for u+P approximately equal to P.

\par\smallskip\noindent\textit{Algebra and arbitrary word values.}\quad 
Addition is join in the order 2 \ensuremath{<} 1 \ensuremath{<} 0 \ensuremath{<} 4 and 2 \ensuremath{<} 1 \ensuremath{<} 3 \ensuremath{<} 4, with 0 and 3 incomparable. In numerical order the tables are

{\small\par\smallskip
\noindent addition       multiplication\par
\CaseLine{\ttfamily 0 0 0 4 4\qquad 0 2 2 3 4}
\CaseLine{\ttfamily 0 1 1 3 4\qquad 2 2 2 2 2}
\CaseLine{\ttfamily 0 1 2 3 4\qquad 2 2 2 2 2}
\CaseLine{\ttfamily 4 3 3 3 4\qquad 4 4 4 4 4}
\CaseLine{\ttfamily 4 4 4 4 4\qquad 4 4 4 4 4}
\smallskip}

For a word of length at least two, scan its positions from left to right for the first factor different from 0. If there is none, the value is 0. A first factor in \{1,\allowbreak{}2\} gives 2. A first factor 4 gives

4. A first factor 3 gives 3 if it occurs at the final position and 4 otherwise. This rule includes every assignment, with arbitrarily many variables. Linear words retain their assigned values.

To prove associativity, the rows for 1 and 2 are constant 2 and those for 3 and 4 are constant 4. A triple whose first factor is in either pair has the same value under both associations. If the first factor is 0, split according as the second is 0, in \{1,\allowbreak{}2\}, or in \{3,\allowbreak{}4\}. In the latter two cases both associations give 2 and 4, respectively. In the first case left multiplication by 0 acts on the last factor by the idempotent map (0,\allowbreak{}1,\allowbreak{}2,\allowbreak{}3,\allowbreak{}4) \ensuremath{\to} (0,\allowbreak{}2,\allowbreak{}2,\allowbreak{}3,\allowbreak{}4). Thus both associations agree. Every row and column is increasing and preserves the sole incomparable join \(0+3=4\), so both distributive laws hold. The stated word rule follows by induction.

\par\smallskip\noindent\textit{The six identities and normal forms.}\quad 
The proposed basis is

\CaseDisplay{\begin{array}{rll}
(W1)&x^2y\approx xy,\\
(W2)&xyx\approx xy^2,\\
(P0)&xy\le y+xz,\\
(P1)&xyz\le xz+xyt,\\
(P2)&xyz\le xz+xy+x^2,\\
(U)&xy\le xy^2.
\end{array}                                                     \quad\textup{(\ensuremath{B})}}

These are sound. For (W1),(W2) the word rule is unchanged, including whether the first factor 3 has a subsequent position. For the four inequalities, consider the element assigned to x. Values 1 and 2 make every left side 2, hence automatically absorbed. Values 3 and 4 make the left side 4, supplied by xz or xz with another factor as displayed, or by \(xy^{2}\) in (U). At \(x=0\), (P0) is immediate when y is 0 or 3, and automatic when y is 1 or 2; at \(y=4\) the summand y is 4. For (P1),(P2), \(y=0\) leaves the left side equal to xz; y in \{1,\allowbreak{}2\} gives 2; \(y=4\) gives a right value 4. At \(y=3\) the left side is 4; xyt is 4 in (P1), while \(xy+x^{2}=3+0=4\) in (P2). For (U) at \(x=0\), \(y=3\) changes the values from 3 to 4; all other y give equality. This covers all cases without restrictions on substituted nonempty terms.

Let F(w) be the list of variables in first-occurrence order. Call w saturated if its final variable occurs earlier. From (W1),(W2) one obtains, for nonempty substituted words Y, Z,

\CaseDisplay{xYxZ\approx xY^2Z\approx xYZ.             \quad\textup{(\ensuremath{2})}}

Together with (W1), this deletes every repeated occurrence before the final position. If the final variable is simple, the result is F(w). Otherwise write the reduced word as p x q x, where p x q lists the distinct variables. If q is empty it is already that list followed by its last variable. If q is nonempty, (W2) gives p x \(q^{2}\), and (2),(W1) remove all repeated prefix positions, leaving p x q followed by the last variable of q. Thus the normal form is F(w) for an unsaturated word, and F(w) followed by its last variable for a saturated word. Linear words are left linear. Nonlinear one-variable words normalize to \(x^{2}\). Omitted contexts never serve as empty-word substitutions.

\par\smallskip\noindent\textit{An exact marked-prefix criterion.}\quad 
For a variable x occurring in u, let \(pre_{u}\)(x) be the set of variables strictly before its first occurrence. The occurrence is marked when there is a position after it. For a set S, define

\CaseDisplay{\begin{aligned}
Z_P(S)&:\quad \exists u\in P\text{ with }c(u)\subseteq S,\\
E_P(x, S)&:\quad \exists u\in P\text{ containing }x
                 \text{ with }\operatorname{pre}_u(x)\subseteq S,\\
M_P(x, S)&:\quad \exists u\in P\text{ as in }E_P(x, S)
                 \text{ with that occurrence marked}\\
 &\qquad\quad\text{or }\bigl(E_P(x, S)\text{ and }Z_P(S)\bigr).
\end{aligned}                                                     \quad\textup{(\ensuremath{3})}}

All summands, including linear ones, are eligible. No restriction on variables following the selected occurrence of x is imposed.

Lemma. A linear variable y is absorbed by P exactly when it is a literal summand. For a nonlinear w, A satisfies w \ensuremath{\leq} P exactly when \(Z_{P}\)(c(w)) holds and, for each x in F(w) with \(S=pre_{w}(x)\), \(M_{P}\)(x, S) holds if the first x is marked, and \(E_{P}\)(x, S) holds if it is unmarked. There is at most one unmarked first occurrence: the simple final variable of an unsaturated word.

Proof. A missing literal y is separated by assigning \(y=1\) and all other variables 2: all right summands are 2, while \(y=1\). For a nonlinear target, assign c(w) the value 0 and outsiders 2. The target is 0; a right word has value 0 exactly when all its variables lie in c(w), and otherwise it is 2. Thus \(Z_{P}\)(c(w)) is necessary.

For the condition at x, assign S the value 0, x the value 3, and every other variable 2. A right summand is 0 exactly when its content lies in S. It is 3 exactly when its first occurrence outside S is x at the final position; it is 4 exactly when that first x has a following position. All remaining right words are

2. The target is 4 for a marked x and 3 for an unmarked x. The right sum is 4 precisely when it has a 4 summand or both a 0 and a 3 summand, giving \(M_{P}\)(x, S). It is at least 3 precisely when some summand has value 3 or 4, giving \(E_{P}\)(x, S). This proves necessity.

Conversely, a target of value 2 is least. If all its variables are 0, its contained witness gives value 0. Otherwise let x be its first factor not valued 0, with all S assigned 0. Values 1 and 2 give target 2, already covered. If \(x=4\), an E-witness has value 4. If \(x=3\) is unmarked, an E-witness has value 3 or 4. If \(x=3\) is marked, a marked witness has value 4; alternatively an E-witness has value 3 or 4 and the Z-witness has value 0, whose sum is 4. These checks prove sufficiency for all words and assignments. \(\square\)

\par\smallskip\noindent\textit{Unrestricted equational completeness.}\quad 
Normalize all words in a valid absorption w \ensuremath{\leq} P. Linear targets are literal inclusions, so suppose w is nonlinear, and write \(F(w)=x_{1}\cdots x_{d}\). Put \(p_{i}=x_{1}\cdots x_{i}\).

For every marked target position we construct

\CaseDisplay{p_i u\le P\quad\text{for every }u\in P. \quad\textup{(\ensuremath{4})}}

The first variable of a nonlinear target is marked. Since no nonempty word has content contained in the empty set, its M-witness must be an original word \(x_{1}\) t with t nonempty. For any original u, (P0) with \(x=x_{1}\), \(y=u\), \(z=t\) gives \(x_{1}u\) \ensuremath{\leq} u+\(x_{1}t\) \ensuremath{\leq} P. This proves (4) at \(i=1\), even for a linear u.

Suppose (4) holds for \(p=p_{i}\) and let \(x=x_{i+1}\) be marked. If a marked witness v exists, delete its repeated prefix letters after prepending p, using (2),(W1), to obtain p v approximately p x t, where t is the nonempty suffix after its first x. Thus pxt \ensuremath{\leq} P. For every original u, (P1) gives

\CaseDisplay{pxu\le pu+pxt\le P.                  \quad\textup{(\ensuremath{5})}}

Otherwise choose an E-witness v and a \(Z_{P}\)(c(p))-witness r. The first x in v must be final, so pv approximately px by prefix deletion. Also pr approximately \(p^{2}\): both have first-occurrence list p and a repeated final variable, since c(r) subset c(p). By the induction hypothesis, px and \(p^{2}\) are below P. Now (P2) gives, for every original u,

\CaseDisplay{pxu\le pu+px+p^2\le P.                 \quad\textup{(\ensuremath{6})}}

This establishes (4) at every marked position.

If w is saturated, all d positions are marked. Choose the original contained witness \(u_{0}\) from \(Z_{P}\)(c(w)). Then \(p_{d}\) \(u_{0}\) has the same normal form as w, so (4) proves w \ensuremath{\leq} P. This includes \(d=1\).

If w is unsaturated, d is at least two and \(w=p_{d-1}s\) with \(s=x_{d}\). Write \(p=p_{d-1}\); the construction has reached (4) for p. Choose an E-witness v for s. If its first s is final, pv approximately \(ps=w\), immediately proving the absorption. Otherwise pv approximately pst for a nonempty suffix t. Using (P1) once with the contained original \(u_{0}\) gives \(psu_{0}\) \ensuremath{\leq} \(pu_{0}\)+pst \ensuremath{\leq} P. Since c(\(u_{0}\)) subset c(w), \(psu_{0}\) is saturated and has first-occurrence list ps, hence normalizes to \(ps^{2}\). Finally (U), with \(x=p\) and \(y=s\), gives \(w=ps\) \ensuremath{\leq} \(ps^{2}\) \ensuremath{\leq} P. Every substitution here is nonempty.

Theorem. The six identities (B) form an identity basis for A relative to the additively idempotent semiring axioms. Reversing every word gives a six-identity basis for its multiplicative opposite.

Proof. Soundness and derivability of every valid word absorption were established above, without bounds on variables, word lengths, or right summands. Distributivity expands arbitrary terms to nonempty polynomials. In a valid polynomial equality, derive each summand's absorption into the opposite polynomial, in both directions. This gives the equality itself. Reversal preserves validity and equational derivations for the opposite semiring. \(\square\)

No minimality or hereditary finite-basability claim is made. Finite regressions check implementations of the criterion and derivations; the unrestricted argument is the proof above. Further catalogue types require independent complete two-way HSP verification.

\subsection{Six identities for a semiring with right reset maps}\label{case:38}
Let A be catalogue entry 11021. The signature is (+, multiplication), without constants; all bases are relative to the additively idempotent semiring axioms. Numerical labels name elements only. Write u \ensuremath{\leq} P for u+P approximately equal to P.

\par\smallskip\noindent\textit{Algebra, word values and six identities.}\quad 
Addition is maximum on the \(chain 3<2<1<0<4\). The tables, in numerical element order, are

{\small\par\smallskip
\noindent addition       multiplication\par
\CaseLine{\ttfamily 0 0 0 0 4\qquad 4 0 0 3 4}
\CaseLine{\ttfamily 0 1 1 1 4\qquad 4 2 2 3 4}
\CaseLine{\ttfamily 0 1 2 2 4\qquad 4 2 2 3 4}
\CaseLine{\ttfamily 0 1 2 3 4\qquad 4 3 3 3 4}
\CaseLine{\ttfamily 4 4 4 4 4\qquad 4 4 4 3 4}
\smallskip}

Right multiplication by 0 or 4 is the constant map 4, by 3 is the constant map 3, and by 1 or 2 is the idempotent map (0,\allowbreak{}1,\allowbreak{}2,\allowbreak{}3,\allowbreak{}4) \ensuremath{\to} (0,\allowbreak{}2,\allowbreak{}2,\allowbreak{}3,\allowbreak{}4). Thus for a nonlinear word inspect its tail occurrences assigned an element of\{0,\allowbreak{}3,\allowbreak{}4\}. If some exist, its value is 3 when the last such occurrence is 3, and 4 when it is 0 or 4. If none exists, all tail factors lie in\{1,\allowbreak{}2\}; the value is 2 when the head lies in\{1,\allowbreak{}2\}, and otherwise is the assigned head 0,\allowbreak{}3 or 4. Linear words keep their assigned value.

These maps prove associativity. If the final factor of a triple is 0,\allowbreak{}3 or 4, both associations give its constant reset value; if it is 1 or 2, applying the idempotent map after any right multiplication leaves every product unchanged. The remaining association either uses that same reset or the idempotent map a second time. Every row and column is increasing in the displayed chain, so multiplication preserves binary maxima. This proves the semiring axioms, and induction gives the word formula for every length.

The following six identities form the proposed basis:

\CaseDisplay{\begin{array}{rll}
(W1)&xyzy\approx xzy,\\
(W2)&xy^2\approx xy,\\
(P)&xy\leq x+zy,\\
(D)&xy\leq x^2y,\\
(H)&xy\leq yxy,\\
(J)&xy\leq y^2+zxy.
\end{array}                                                     \quad\textup{(\ensuremath{B})}}

For(W1), the final y either resets the product or acts by the idempotent map, making the earlier y redundant; (W2) repeats an idempotent right map. To check the inequalities, if y is 0 or 4, the last y forces 4 in every right product ending in y; if \(y=3\), the left side is the bottom 3. If y is 1 or 2, then \(xy\leq x\), proving(P). For(D), squaring x changes \(x=0\) to 4, sends \(x=1\) to 2, and fixes2,\allowbreak{}3,\allowbreak{}4; its final product is at least xy. For(H), when x is 0,\allowbreak{}3 or 4 the internal x resets yxy to 4,\allowbreak{}3 or 4, respectively, at least xy; for x in\{1,\allowbreak{}2\}, both products are 2. For(J), x in\{1,\allowbreak{}2\} gives \(xy=2\leq y^{2}\); \(x=3\) gives the bottom; and \(x=0\) or 4 makes \(zxy=4\). This proves soundness.

Using(W1),(W2), delete every nonfinal occurrence of a variable in the tail: adjacent repetitions use(W2), separated ones use(W1) with nonempty intervening word and nonempty preceding prefix. Hence every nonlinear word normalizes to its literal head followed by the distinct tail variables in last-occurrence order. An occurrence of the head in the tail remains; linear words remain linear. These are derivable normal forms, without any claim of uniqueness beyond the data retained.

\par\smallskip\noindent\textit{An exact ordered-suffix criterion.}\quad 
For a variable x in u, let \(suf_{u}\)(x) be the set after its last occurrence in u. For a finite set S put

\CaseDisplay{\begin{aligned}
Z_P(S)&:\quad\exists u\in P\text{ with }c(u)\subseteq S,\\
E_P(x, S)&:\quad\exists u\in P\text{ containing }x
                       \text{ with }\operatorname{suf}_u(x)\subseteq S,\\
T_P(x, S)&:\quad\exists u\in P\text{ with }x\text{ in its tail}
                       \text{ and }\operatorname{suf}_u(x)\subseteq S.
\end{aligned}                                                     \quad\textup{(\ensuremath{1})}}

All original summands, including linear ones where applicable, may serve as witnesses. No restriction is imposed on variables preceding the selected last occurrence.

Lemma 1. A linear y is absorbed by P exactly when it is a literal summand. For a nonlinear normal word \(w=h\) t, whose tail t lists distinct variables, w \ensuremath{\leq} P holds in A exactly when

1. \(Z_{P}\)(c(w)) holds;

2. for each tail variable x, \(T_{P}\)(x, \(suf_{w}\)(x)) holds;

3. if h does not occur in t, \(E_{P}\)(h, c(t)) holds.

Proof. For a missing literal y assign \(y=1\) and all other variables

3. Every nonlinear power of y is 2, and every other right word is 3; therefore \(P\leq 2<1=y\). The converse is literal inclusion.

If 1 fails, assign all target variables 2 and all outsiders 3. The target is 2 and every right word is 3, so absorption fails. For a tail variable x put \(S=suf_{w}(x)\), assign \(x=0\), variables in S the value 2, and all others 3. Then \(w=4\). A right word takes value 4 exactly when it has x in its tail with its last x followed only by variables in S. If 2 fails, no such word exists; right values are at most 0, strictly below 4. Finally suppose h is simple and 3 fails. Assign \(h=0\), all target tail variables 2 and outsiders 3. Then \(w=0\). A right word is at least 0 exactly when it contains h with its last occurrence followed only by target tail variables. If 3 fails all right values are at most 2. This proves necessity.

Conversely, \(w=3\) is bottom. If \(w=2\), all its variables lie in\{1,\allowbreak{}2\}, and a contained witness from 1 has value at least 2. If the target tail has any reset occurrence, let x be the last such occurrence. A reset 3 makes w bottom; a reset 0 or 4 makes \(w=4\). Every variable after this occurrence has value 1 or 2. A witness from 2 has the same reset at its last x and only 1/\allowbreak{}2 factors after it, so it is 4. If the target tail has no reset, all its variables lie in\{1,\allowbreak{}2\}. The only cases not covered already are head 0 or 4. The head is then simple, and a witness from 3 preserves the head value when its last h is initial, or resets to 4 when it is internal. In either case it is at least the target. This proves sufficiency at all assignments. \(\square\)

\par\smallskip\noindent\textit{Appending a common tail by induction from the right.}\quad 
Fix a valid nonlinear absorption \(w=h\) t \ensuremath{\leq} P and normalize its words. We first construct

\CaseDisplay{u s\leq P\quad\text{for every }u\in P \quad\textup{(\ensuremath{2})}}

for each suffix s of the target tail, proceeding from empty suffix to the whole tail. At the empty suffix this means the original inclusions \(u\leq P\); it does not substitute an empty word in an identity.

Suppose 2 has been reached at s and x is the preceding target tail variable. Choose by condition 2 of Lemma 1 an original witness \(v=r\) x q, taking its last x, with r nonempty and c(q) contained in c(s). Tail normalization gives v s approximately r x s; when q or s is empty these are omitted substrings. By induction v \(s\leq P\). For every original u, identity(P) with \(X=u\), \(Y=x\), \(Z=r\), followed by the right context s when nonempty, gives

\CaseDisplay{u x s\leq u s+r x s\leq P.             \quad\textup{(\ensuremath{3})}}

All substituted words are nonempty. This establishes2 at x s and eventually proves u \(t\leq P\) for every original u, with no restriction on variables preceding the witnesses.

\par\smallskip\noindent\textit{Recovering the target head.}\quad 
Choose a contained original witness u from condition 1 of Lemma 1 and put \(k=h(u)\). The tail t is nonempty. Two cases finish the proof.

If h occurs in t, every variable of u lies in c(t). Thus u t normalizes to k t. If \(k=h\), this is the target. Otherwise split \(t=a\) k b, put \(Y=k\) b and \(X=h\) a, both nonempty. The one-right-word inequality(H) gives

\CaseDisplay{h t=XY\leq YXY\approx k t\leq P.              \quad\textup{(\ensuremath{4})}}

The displayed normalization uses the last copy of t to delete all earlier tail occurrences, including h.

Now suppose h does not occur in t. If u contains h and starts with h, then u t normalizes to h t when h is simple in u, or to \(h^{2}t\) when it repeats. The first case is finished and the second uses(D), with \(X=h\), \(Y=t\). If u contains h but starts with k different from h, then u t normalizes to k h t. Split \(t=a\) k b as above; the same instance of(H) gives

\CaseDisplay{h t\leq(k b)(h a)(k b)
                           \approx k h t\leq P.             \quad\textup{(\ensuremath{5})}}

Finally suppose u omits h, so u t approximately k t. Let g be the first variable of t, split \(t=a\) k b, and put \(X=g\) a, \(Y=k\) b. Both substituted words are nonempty even if a is empty. Then \(XY=g\) t approximately \(t^{2}\) and YXY approximately k t, by last-tail normalization. Thus(H) gives \(t^{2}\leq k\) \(t\leq P\). This also handles a one-letter t, where the square is already the known k t.

Choose an original E-witness v for h, writing \(v=r\) h q at its last h, with c(q) contained in c(t). Then v t approximately r h \(t\leq P\) by 2. If r is empty, this is the target. Otherwise(J), with \(X=h\), \(Y=t\), \(Z=r\), gives

\CaseDisplay{h t\leq t^2+r h t\leq P.               \quad\textup{(\ensuremath{6})}}

Every substitution is nonempty; omitted contexts are handled by omission alone. This exhausts the criterion and derives every valid word absorption from the six identities, for arbitrary variables, lengths and numbers of summands.

Theorem. The six identities(B) form a basis for the identities of A relative to the additively idempotent semiring axioms. Reversing every word gives a six-identity basis for its multiplicative opposite.

Proof. Soundness, normalizations and all valid word absorptions have been proved. Distributivity expands arbitrary terms to nonempty polynomials. In a valid identity derive each summand's absorption into the other polynomial in both directions, giving the identity itself. Reversal transports validity and all equational derivations to the opposite semiring. \(\square\)

No minimality or hereditary assertion is made. Finite checks replay the criterion and explicit substitution steps; they are auxiliary to the unrestricted proof. Transfers to other catalogue entries require independently verified complete HSP inclusions in both directions. The earlier 11021 research note records the incomplete starting point and is retained as historical evidence.

\subsection{最后出现后缀：非交换右作用的五条完整基}\label{case:39}
\par\smallskip\noindent\textit{1. 半环与有序词值.}\quad 
取 \(A=\{z, a, b, d, t\}\)，依次对应目录标号 \((0,\allowbreak{}1,\allowbreak{}2,\allowbreak{}3,\allowbreak{}4)\)。加法偏序为

\CaseDisplay{b<a<z<t,\qquad d<a<z<t,}

其中 \(b, d\) 不可比。右乘 \(a, b\) 为 \(R=(z, a, a, d, t)\)，右乘 \(d\) 为常值 \(d\) 映射，右乘 \(z, t\) 为常值 \(t\) 映射。即乘法表为

\CaseDisplay{\begin{array}{c|ccccc}
\cdot&z&a&b&d&t\\\hline
z&t&z&z&d&t\\
a&t&a&a&d&t\\
b&t&a&a&d&t\\
d&t&d&d&d&t\\
t&t&t&t&d&t
\end{array}}

三个右平移保持并且幂等，但两个不同常值映射不交换。按右因子类别直接验证 \(R_{xy}=R_yR_x\) 和 \(R_{x+y}=R_x\vee R_y\)，即得结合律和分配律。这解释了为何该来源不满足 \(xyz\approx xzy\)。

非一次词的值如下：若尾部全部取 \(a, b\)，则取 \(R(h(w))\)；否则取尾部最右一个属于 \(\{d, z, t\}\) 的字母，其值为 \(d\) 时全词取 \(d\)，为 \(z, t\) 时全词取 \(t\)。一次词取变量自身值。

\par\smallskip\noindent\textit{2. 正规形与记号.}\quad 
考虑两条词恒等式

\CaseDisplay{xy^2\approx xy,\qquad xyx\approx y^2x.\quad\textup{(\ensuremath{W})}}

它们推出 \(x^3\approx x^2\) 和

\CaseDisplay{xyzy\approx xz^2y\approx xzy.}

因此每个尾部变量都能删除其最后一次出现以前的重复。若原首字母不在尾部再出现，所得词是一个无重复词 \(r\)；若原首字母在尾部再出现，所得词等于 \(r^2\)，其中 \(r\) 按原词各变量最后一次出现的次序排列。

更具体地，把尾部去重后写作 \(hphq\)，其中 \(p, q\) 都不含 \(h\)。当 \(p\) 非空时，\(hph\approx p^2h\)；当 \(p\) 为空时已有 \(h^2q\)。再删除尾部重复，均得到 \(r\) 的首字母再接整个 \(r\)，这也正是 \(r^2\) 的正规形。空段仅表示省略文字，不作为恒等式的代入项。

以下全部右词及左词均选用上述正规形。对于变量 \(x\in c(w)\)，令 \(S_w(x)\) 表示 \(x\) 最后一次出现之后的变量集合；\(D(w)\) 表示尾部变量集合。无重复形式 \(r=hx_1\cdots x_k\) 有 \(D(r)=\{x_1,\ldots, x_k\}\)；平方形式 \(r^2\) 的尾部集为全部内容。

\par\smallskip\noindent\textit{3. 精确吸收判据.}\quad 
令 \(F\) 为非空右词集，\(X=c(w)\)，\(F_X=\{u\in F:c(u)\subseteq X\}\)。

定理1。 在 \(A\) 上，\(w\le F\) 当且仅当同时满足：

1. \(F_X\ne\varnothing\)；

2. 对每个 \(x\in X\)，有右项 \(u_x\) 包含 \(x\)，且 \(S_{u_x}(x)\subseteq S_w(x)\)；

3. 对每个 \(x\in D(w)\)，有右项 \(v_x\) 在尾部包含 \(x\)，且 \(S_{v_x}(x)\subseteq S_w(x)\)。

第三条件要求尾部见证，不能只要求该变量出现在右词任意位置。对一次左词 \(x\)，第一条件已经蕴含第二条件，第三条件为空。

\par\smallskip\noindent\textit{必要性.}\quad 
若没有内部右项，赋 \(X\) 内变量 \(a\)、外部变量 \(d\)。左值为 \(a\)，每个右词都取 \(d\)，吸收失败。

若条件2在 \(x\) 处失败，赋 \(x=t\)、\(S_w(x)\) 内变量 \(a\)，其余变量 \(d\)。左词在最后一个 \(x\) 之后只有中性值 \(a\)，所以取 \(t\)。任一含 \(x\) 的右词在最后一个 \(x\) 之后都有取值 \(d\) 的变量，因而取 \(d\)；不含 \(x\) 的右词取 \(a\) 或 \(d\)。右和不超过 \(a\)，矛盾。

若条件3在 \(x\in D(w)\) 处失败，改赋 \(x=z\)、\(S_w(x)\) 内变量 \(a\)，其余变量 \(d\)。左词的最后一个 \(x\) 位于尾部，故左值为 \(t\)。右词若在尾部含 \(x\)，其最后一个 \(x\) 之后必有取值 \(d\) 的变量，故取 \(d\)；若只在首字母含 \(x\)，其值为 \(z\) 或 \(d\)；其余右词取 \(a, d\)。右和不超过 \(z\)，矛盾。

\par\smallskip\noindent\textit{充分性.}\quad 
一次左词由内部纯幂及有效不等式 \(x\le x^2\) 处理。以下设左词非一次。真主理想分别为 \(\{b\}\)、\(\{d\}\)、\(\{b, d, a\}\)、\(\{b, d, a, z\}\)。

右和为 \(b\) 时，右词只能是赋值为 \(b\) 的一次词；非一次词不会取 \(b\)。但左尾部非空，条件3给出非一次右项，矛盾。因此该情形不发生。

若右和不超过 \(a\)，而左值为 \(z\) 或 \(t\)，令 \(x\) 为左词中最右一个取值于 \(\{d, z, t\}\) 的变量。它必须取 \(z\) 或 \(t\)，其后变量全部取 \(a, b\)。条件2提供后缀见证，保证相应右词取 \(z\) 或 \(t\)，矛盾。

若右和为 \(d\)，同一论证排除左值 \(z, t\)。若左内容含取值于 \(\{d, z, t\}\) 的变量，则最右这样的变量必须取 \(d\)，左词取 \(d\)。若左内容全取 \(a, b\)，内部右项的值只能为 \(a\) 或 \(b\)，也不可能使右和为 \(d\)。故左值仍为 \(d\)。

最后设右和不超过 \(z\)。若左值为 \(t\)，有两种情况：尾部最右的非中性变量 \(x\) 取 \(z\) 或 \(t\)，此时条件3提供尾部后缀见证，其值为 \(t\)；或尾部全取 \(a, b\)、首字母取 \(t\)，此时条件2的首字母后缀见证取 \(t\)。两者都矛盾。由此完成全部真主理想检查。

\par\smallskip\noindent\textit{4. 五条完整基.}\quad 
定理2。 \(A\) 的完整附加基为

\CaseDisplay{\begin{aligned}
xy^2&\approx xy, &(R)\\
xyx&\approx y^2x, &(C)\\
x&\le x^2, &(I)\\
xy&\le x+zy, &(V)\\
xy&\le y^2+zxy. &(J)
\end{aligned}\quad\textup{(\ensuremath{3})}}

前两式由词值规则直接验证，后三式由定理1检查。以下证明任意变量数完备性。

\par\smallskip\noindent\textit{4.1 从右向左的后缀推进引理.}\quad 
取互不相同的变量序列 \(r=x_1\cdots x_k\)。假设每个 \(x_i\) 都有右尾部见证

\CaseDisplay{u_i=a_ix_ib_i\in F,\qquad a_i\ne\varnothing,\qquad
c(b_i)\subseteq\{x_{i+1},\ldots, x_k\},}

其中选取的是 \(u_i\) 中最后一个 \(x_i\)，\(b_i\) 可以为空。令 \(r_i=x_i\cdots x_k\)。则对任意由基可推出 \(v\le F\) 的非空词 \(v\)，都可推出

\CaseDisplay{vr_i\le F.\quad\textup{(\ensuremath{4})}}

从右向左归纳。\(i=k\) 时 \(b_k\) 为空，由 \((V)\) 得 \(vx_k\le v+a_kx_k=v+u_k\le F\)。假设对 \(r_{i+1}\) 已证明结论，则把 \(u_i\) 作为 \(v\)，得到 \(u_ir_{i+1}\le F\)。由于 \(b_i\) 中全部变量在 \(r_{i+1}\) 里再次出现，尾部去重给出

\CaseDisplay{u_ir_{i+1}\approx a_ix_ir_{i+1}.}

以 \((x, y, z)=(v, x_ir_{i+1},a_i)\) 代入 \((V)\)，得

\CaseDisplay{vr_i\le v+a_ix_ir_{i+1}=v+u_ir_{i+1}\le F.}

这证明 (4)。整个归纳仅使用非空词代入。

\par\smallskip\noindent\textit{4.2 内部项与首字母补入.}\quad 
若左词为平方正规形 \(w=r^2\)，其尾部集为全部 \(X\)，条件3为整个有序序列 \(r\) 提供后缀推进见证。选内部项 \(v\in F_X\)，由引理得 \(vr\le F\)。因为 \(c(v)\subseteq c(r)\)，\(vr\) 的原首字母在后面的 \(r\) 中重复，且最后出现次序为 \(r\)，故 \(vr\approx r^2=w\)。

若左词为无重复非一次词 \(w=hr\)，\(h\notin c(r)\)，条件3为 \(r\) 提供后缀推进见证。选内部项 \(v\in F_X\)，由引理得 \(vr\le F\)。

当 \(h\in c(v)\) 时，\(vr\) 的最后出现次序是 \(hr\)，其正规形等于 \(w\) 或 \(w^2\)；后一情形再由 \((I)\) 得 \(w\le F\)。

当 \(h\notin c(v)\) 时，\(c(v)\subseteq c(r)\)，所以 \(vr\approx r^2\le F\)。条件2提供一个首字母后缀见证，在其最后一个 \(h\) 处分割为 \(u=ahs\)，其中前段 \(a\)、后段 \(s\) 可以为空，且 \(c(s)\subseteq c(r)\)。

若前段 \(a\) 为空，则 \(ur\approx hr=w\)，由引理即得吸收。若 \(a\) 非空，则 \(ur\approx ahr\le F\)；以 \((x, y, z)=(h, r, a)\) 代入 \((J)\)，得

\CaseDisplay{w=hr\le r^2+ahr=r^2+ur\le F.}

一次左词由内部纯幂、\((R),(I)\) 完成。一般半环项按分配律展开为非空词和，逐词处理双向吸收，所以五式完整。证毕。

\par\smallskip\noindent\textit{5. 方法意义与范围.}\quad 
本方法把``尾部集合覆盖''推广为``最后出现位置之后的后缀包含''。有序尾部无法交换，但通过从右到左的推进，同一条 \((V)\) 可以处理任意长度的目标后缀；\((J)\) 负责补入仅有全局见证的首字母。

乘法对偶的基由逐词反转取得，判据相应变为首次出现前缀。有限基性向其他五元素类型传递必须另有独立双向HSP证据。有限变量诊断只用于交叉核验，不证明任意变量数完备性；本定理不声称遗传有限基。

\subsection{首字母与内部锚点：五组平方回缩完整基}\label{case:40}
\par\smallskip\noindent\textit{1. 五个具体半环.}\quad 
四个半环 \(P, Q, R, U\) 的底集均记为 \(\{z, b, a, d, t\}\)，加法取偏序

\CaseDisplay{z<b<a<t,\qquad b<d<t,\qquad a\parallel d}

上的并。第五个半环 \(V\) 的底集符号相同，加法取链 \(z<d<b<a<t\) 上的最大值。下面的表中所有行、列均按 \((z, b, a, d, t)\) 排列；表项用五个元素的符号书写。

\noindent 半环\enspace\textbar\enspace \(z\) 行\enspace\textbar\enspace \(b\) 行\enspace\textbar\enspace \(a\) 行\enspace\textbar\enspace \(d\) 行\enspace\textbar\enspace \(t\) 行\enspace\textbar\enspace 目录编号\enspace\textbar\enspace 符号到目录的映射\par
\noindent \(P\)\enspace\textbar\enspace \(z, b, b, b, b\)\enspace\textbar\enspace \(b, b, b, b, b\)\enspace\textbar\enspace \(a, a, a, a, a\)\enspace\textbar\enspace \(a, a, a, t, t\)\enspace\textbar\enspace \(a, a, a, t, t\)\enspace\textbar\enspace 1959\enspace\textbar\enspace \((2,\allowbreak{}1,\allowbreak{}0,\allowbreak{}3,\allowbreak{}4)\)\par
\noindent \(Q\)\enspace\textbar\enspace \(z, z, z, z, z\)\enspace\textbar\enspace \(z, b, b, b, b\)\enspace\textbar\enspace \(z, a, a, a, a\)\enspace\textbar\enspace \(z, a, a, t, t\)\enspace\textbar\enspace \(z, a, a, t, t\)\enspace\textbar\enspace 4177\enspace\textbar\enspace \((2,\allowbreak{}1,\allowbreak{}0,\allowbreak{}3,\allowbreak{}4)\)\par
\noindent \(R\)\enspace\textbar\enspace \(z, z, z, z, z\)\enspace\textbar\enspace \(z, b, b, b, b\)\enspace\textbar\enspace \(z, a, a, a, a\)\enspace\textbar\enspace \(z, a, a, a, a\)\enspace\textbar\enspace \(z, a, a, a, a\)\enspace\textbar\enspace 4172\enspace\textbar\enspace \((2,\allowbreak{}1,\allowbreak{}0,\allowbreak{}3,\allowbreak{}4)\)\par
\noindent \(U\)\enspace\textbar\enspace \(z, z, z, z, z\)\enspace\textbar\enspace \(z, z, z, z, z\)\enspace\textbar\enspace \(a, a, a, a, a\)\enspace\textbar\enspace \(a, a, a, t, t\)\enspace\textbar\enspace \(a, a, a, t, t\)\enspace\textbar\enspace 2155\enspace\textbar\enspace \((2,\allowbreak{}1,\allowbreak{}0,\allowbreak{}3,\allowbreak{}4)\)\par
\noindent \(V\)\enspace\textbar\enspace \(z, z, z, z, z\)\enspace\textbar\enspace \(z, b, b, b, b\)\enspace\textbar\enspace \(z, a, a, a, a\)\enspace\textbar\enspace \(z, b, b, b, b\)\enspace\textbar\enspace \(z, a, a, a, a\)\enspace\textbar\enspace 3094\enspace\textbar\enspace \((3,\allowbreak{}1,\allowbreak{}0,\allowbreak{}2,\allowbreak{}4)\)\par
所有这些表都定义加法幂等半环。加法公理来自并半格。分配律可逐行、逐列检查：每个平移在可比较对上保序；对前四个半环，唯一还需检查的不可比较对为 \(a, d\)，其并是 \(t\)，表中每个平移都保持这个并。第五个半环的平移全部在链上保序。

为同时说明结合性及后续证明，给出长度至少二的词 \(w\) 的值。记首字母为 \(h\)、内容为 \(C\)。

\textbullet{} \(P\)：若 \(h=b, a\)，分别取 \(b, a\)；若 \(h=z\)，全体 \(C\) 都取 \(z\) 时取 \(z\)，否则取 \(b\)；若 \(h=d, t\)，全体 \(C\) 都在 \(\{d, t\}\) 时取 \(t\)，否则取 \(a\)。

\textbullet{} \(Q\)：若 \(C\) 中出现 \(z\) 则取 \(z\)；否则 \(h=b, a\) 分别取 \(b, a\)，而 \(h=d, t\) 在全体 \(C\) 都属于 \(\{d, t\}\) 时取 \(t\)，否则取 \(a\)。

\textbullet{} \(R\)：若 \(C\) 中出现 \(z\) 则取 \(z\)；否则 \(h=b\) 时取 \(b\)，\(h=a, d, t\) 时取 \(a\)。

\textbullet{} \(U\)：\(h=z, b\) 时取 \(z\)，\(h=a\) 时取 \(a\)；\(h=d, t\) 时，全体 \(C\) 都属于 \(\{d, t\}\) 则取 \(t\)，否则取 \(a\)。

\textbullet{} \(V\)：若 \(C\) 中出现 \(z\) 则取 \(z\)；否则 \(h=b, d\) 时取 \(b\)，\(h=a, t\) 时取 \(a\)。

直接在两非空词的连接上应用各条规则，所得值等于两个词值之积；特别地，三个字母的两种括号化相等。这验证结合性。

\par\smallskip\noindent\textit{2. 吸收判据.}\quad 
对非空多项式 \(F\)，记 \(L\) 为其一次项变量集，\(\mathcal R\) 为其非线性词族，\(T=c(F)\)。对非空内容集 \(C\)，定义

\CaseDisplay{\mathcal R_C=\{u\in\mathcal R:c(u)\subseteq C\},\qquad
G=\{h(u):u\in\mathcal R\},\qquad G_C=\{h(u):u\in\mathcal R_C\}.}

内部锚点条件记为

\CaseDisplay{A_C:quad L\cap C\ne\varnothing\quad\text{或}\quad\mathcal R_C\ne\varnothing.}

判据定理。 对非线性词 \(w\)，置 \(h=h(w),C=c(w)\)，则 \(w\le F\) 恰在下表条件成立时有效：

\noindent 半环\enspace\textbar\enspace 非线性左词的充分必要条件\enspace\textbar\enspace 一次词 \(x\le F\) 的充分必要条件\par
\noindent \(P\)\enspace\textbar\enspace \(C\subseteq T\)，\(h\in G\)，\(A_C\)\enspace\textbar\enspace \(x\in L\) 或 \(\exists u\in\mathcal R:c(u)=\{x\}\)\par
\noindent \(Q\)\enspace\textbar\enspace \(h\in G_C\)\enspace\textbar\enspace \(x\in L\) 或 \(\exists u\in\mathcal R:c(u)=\{x\}\)\par
\noindent \(R\)\enspace\textbar\enspace \(h\in G_C\)\enspace\textbar\enspace \(x\in L\)\par
\noindent \(U\)\enspace\textbar\enspace \(h\in G\) 且 \(A_C\)\enspace\textbar\enspace \(x\in L\)\par
\noindent \(V\)\enspace\textbar\enspace \(\mathcal R_C\ne\varnothing\) 且 \((h\in L\text{ 或 }h\in G_C)\)\enspace\textbar\enspace \(x\in L\)\par
\par\smallskip\noindent\textit{2.1 必要性.}\quad 
对 \(P\)，若 \(C\nsubseteq T\)，令一个遗漏字母取 \(b\)、其他变量取 \(z\)，左值为 \(b\)、右值为 \(z\)。若 \(h\notin G\)，令 \(h=d\)、其余 \(C\) 内变量取 \(b\)、外部变量取 \(z\)。左值为 \(a\) 或 \(t\)；所有右非线性词的首字母都不是 \(h\)，其值至多为 \(b\)，右一次项至多为 \(d\)，故右值不吸收左值。若 \(A_C\) 失败，令 \(C\) 内变量取 \(d\)、外部变量取 \(z\)，左值为 \(t\)、右值至多为 \(a\)。

对 \(Q, R\)，若 \(h\notin G_C\)，令 \(h=d\)、其余 \(C\) 内变量取 \(b\)、外部变量取 \(z\)。含外部变量的右非线性词为 \(z\)；内部右非线性词因首字母不是 \(h\) 而取 \(b\)；右一次项至多为 \(d\)。\(Q\) 的左值为 \(a\) 或 \(t\)，\(R\) 的左值为 \(a\)，均不被吸收。

对 \(U\)，反驳 \(h\in G\) 的赋值与 \(P\) 相同，但没有以 \(h\) 开始的右非线性词时，它们全取 \(z\)。反驳 \(A_C\) 的赋值也与 \(P\) 相同。

对 \(V\)，若 \(\mathcal R_C\) 为空，令 \(C\) 内变量为 \(d\)、外部变量为 \(z\)，左值为 \(b\)、右值至多为 \(d\)。若 \(h\notin L\cup G_C\)，令 \(h=a\)、其余 \(C\) 内变量为 \(b\)、外部变量为 \(z\)，左值为 \(a\)、右值至多为 \(b\)。

一次词的必要性：对 \(P, Q\)，若没有该一次项和只用 \(x\) 的非线性右词，令 \(x=d\)、其余变量为 \(z\)，右值分别至多为 \(a\) 或等于 \(z\)，均不吸收 \(d\)。对 \(R, V\)，令 \(x=t\)、其余变量为 \(z\)，非线性右值至多为 \(a<t\)；对 \(U\)，令 \(x=b\)、其余变量为 \(z\)，所有非线性右值为 \(z<b\)。

\par\smallskip\noindent\textit{2.2 充分性.}\quad 
对 \(P\)，左值为 \(z\) 时自动成立；为 \(b\) 时，内容内有非 \(z\) 的字母，由 \(C\subseteq T\)，它在某个右项中出现。任意含非 \(z\) 字母的非线性词至少为 \(b\)，任意非 \(z\) 一次项也至少为 \(b\)。左值为 \(a\) 时，\(h\in G\) 给出以 \(h\) 开始的右非线性词，值为 \(a\) 或 \(t\)。左值为 \(t\) 时，所有 \(C\) 内变量都在 \(\{d, t\}\)；内部非线性锚点直接给出 \(t\)，内部一次项若取 \(d\)，与同首字母右非线性词的 \(a\) 或 \(t\) 求并也给出 \(t\)。

对 \(Q\)，若左值非 \(z\)，一个内容包含于 \(C\)、首字母为 \(h\) 的右非线性词，其值不小于左值：当左值为 \(b, a\) 时直接由首字母判断；左值为 \(t\) 时内部右词也取 \(t\)。对 \(R\)，这样的右词和左词的值相同。

对 \(U\)，左值为 \(z\) 时自动成立；为 \(a\) 时，\(h\in G\) 提供值为 \(a\) 或 \(t\) 的右非线性词；为 \(t\) 时，内部锚点的论证与 \(P\) 相同。

对 \(V\)，若左值非 \(z\)，\(C\) 内没有 \(z\)。于是任意内部非线性右词至少为 \(b\)。若左值为 \(a\)，则首字母取 \(a\) 或 \(t\)，同名一次项至少为 \(a\)；若由 \(h\in G_C\) 提供右词，该词取 \(a\)。左值为 \(b\) 已由内部右词吸收。

一次词的充分性除自明的同名一次项外，只需在 \(P, Q\) 使用 \(x\le x^2=x^n\)（\(n\ge2\)）。判据证明完毕。

\par\smallskip\noindent\textit{3. 五组有限基.}\quad 
共同保留两式

\CaseDisplay{xyz\approx xzy\quad(N),\qquad xy\approx x^2y^2\quad(S).}

另记

\CaseDisplay{\begin{aligned}
x&\le x^2 &(I),\\
xy&\le x^2 &(C),\\
xy&\le x^2+y &(K),\\
xy&\le x+yx &(H),\\
xy&\le y+xyz &(E).
\end{aligned}}

有限基定理。 相对于加法幂等半环公理，以下各组都是所列半环的完整恒等式基：

\noindent 半环\enspace\textbar\enspace 完整基\enspace\textbar\enspace 附加式数\par
\noindent \(P\)\enspace\textbar\enspace \((N),(S),(I),(K),(E)\)\enspace\textbar\enspace 5\par
\noindent \(Q\)\enspace\textbar\enspace \((N),(S),(I),(C)\)\enspace\textbar\enspace 4\par
\noindent \(R\)\enspace\textbar\enspace \((N),(S),(C)\)\enspace\textbar\enspace 3\par
\noindent \(U\)\enspace\textbar\enspace \((N),(S),(C),(E)\)\enspace\textbar\enspace 4\par
\noindent \(V\)\enspace\textbar\enspace \((N),(S),(C),(H)\)\enspace\textbar\enspace 4\par
所有列出的式子均由第 1 节词值规则或第 2 节判据直接验证有效。其乘法对偶的完整基由反转每个词得到。

\par\smallskip\noindent\textit{3.1 共同正规形.}\quad 
从 (S) 得 \(x^2=x^4\) 和 \(x^3=x^6=x^2\)。结合 (N)，每个非线性词可整理成 \(W(h, C)\)：首字母为 \(h\)，内容为 \(C\)，每个字母均出现两次。不同尾部排列相等，且

\CaseDisplay{W(h, C)W(g, D)=W(h, C\cup D).}

特别地，非线性词都是乘法幂等的。以下的补乘与消元使用的块均非空；没有空词代入。

\par\smallskip\noindent\textit{3.2 \(P\) 的完备性.}\quad 
设 \(w=W(h, C)\le F\) 有效。从 \(h\in G\) 选一个以 \(h\) 开始的右非线性词，并以它作为初始积，依次乘上其余全部右词，允许后续因子为一次词。积在每一步都是非线性词，故 (K) 给出

\CaseDisplay{pq\le p^2+q=p+q.}

归纳得到 \(W(h, T)\le F\)，其中 \(T=c(F)\)，由判据有 \(C\subseteq T\)。取内部锚点右词 \(q\)，满足 \(c(q)\subseteq C\)。若 \(T=C\) 已完成；否则令 \(s\) 为 \(T\setminus C\) 各字母组成的非空词，以 \((x, y, z)=(w, q, s)\) 代入 (E)，得到

\CaseDisplay{w=wq\le q+wqs=q+W(h, T)\le F.}

一次词吸收由 (I) 及正规形处理。

\par\smallskip\noindent\textit{3.3 \(Q, R\) 的完备性.}\quad 
由判据选 \(u\in\mathcal R_C\)，使 \(h(u)=h\)。若 \(c(u)=C\)，正规形已给 \(u=w\)；否则把 \(C\setminus c(u)\) 各字母组成非空词 \(s\)。在 (C) 中代入 \(x=u, y=s\) 得

\CaseDisplay{w=us\le u^2=u\le F.}

\(Q\) 的一次词使用 (I)，\(R\) 的一次词必须已是右项。

\par\smallskip\noindent\textit{3.4 \(U\) 的完备性.}\quad 
选一个以 \(h\) 开始的非线性右词 \(r\)，令 \(D=C\cup c(r)\)。利用 (C) 补乘 \(C\setminus c(r)\) 中的字母，或在该集合为空时不作补乘，可得 \(W(h, D)\le r\le F\)。再选内部锚点 \(q\)。若 \(D=C\) 则完成；否则以 (E) 同 \(P\) 一样消去 \(D\setminus C\)，得到 \(w\le q+W(h, D)\le F\)。一次词已在右侧。

\par\smallskip\noindent\textit{3.5 \(V\) 的完备性.}\quad 
若 \(h\in G_C\)，与 \(Q, R\) 一样由 (C) 补乘即可。否则 \(h\in L\)，选任意 \(r\in\mathcal R_C\)。由 (H)，

\CaseDisplay{hr\le h+rh.}

而 (C) 给出 \(rh\le r^2=r\)，所以 \(hr\le h+r\le F\)。该非线性词以 \(h\) 开始，内容包含于 \(C\)。再用 (C) 补乘缺少的字母，得到 \(w\le F\)。一次词仍必须已是右项。

每个有效多项式恒等式均可拆成两个方向的逐词吸收，故五组基全部完整。这是任意变量数的形式推导，不依赖有界搜索。

\par\smallskip\noindent\textit{4. 目录传递范围.}\quad 
若目录中的半环与上述生成元生成同一簇，引理~\ref{lem:transfer} 给出同一组有限基。两个方向的簇包含可由命题~\ref{prop:lower-fibres} 从运算表检验。

\subsection{An eight-identity tail-guard basis for source 2137}\label{case:41}
Work in the signature (+, multiplication) without constants, relative to the axioms of additively idempotent semirings. The element labels in the following tables name elements only. Source A is catalogue entry 2137.

{\small\par\smallskip
\noindent addition    multiplication\par
\CaseLine{\ttfamily 00044\qquad 00000}
\CaseLine{\ttfamily 01134\qquad 22211}
\CaseLine{\ttfamily 01234\qquad 22222}
\CaseLine{\ttfamily 43334\qquad 00044}
\CaseLine{\ttfamily 44444\qquad 00044}
\smallskip}

Rows and columns are in the order 0,\allowbreak{}1,\allowbreak{}2,\allowbreak{}3,\allowbreak{}4. Addition is join in the order \(2<1<0<4\) and \(2<1<3<4\), with 0 and 3 incomparable. Set K=\{3,\allowbreak{}4\}. The elements 0 and 2 are left zeros. A product of at least two K-elements is 4; a K-element followed by an element outside K gives 0. An initial 1 remains 1 exactly while all later factors lie in K, and becomes 2 otherwise. These descriptions verify associativity: a nonempty block has product in K exactly when all its factors are in K, and the remaining head cases reduce to the two left zeros or the indicated transition from 1 to 2. Every left and right multiplication preserves the joins of comparable elements and the sole incomparable pair 0,\allowbreak{}3, as checked from the displayed rows and columns. Thus A is a semiring.

For a nonlinear word with head h and nonempty tail-variable set T, the value is:

\textbullet{} 0 when \(h=0\), and 2 when \(h=2\);

\textbullet{} 1 when \(h=1\) and all variables of T are in K, and 2 otherwise when \(h=1\);

\textbullet{} 4 when h is in K and all variables of T are in K, and 0 otherwise when h is in K.

Consequently the following two identities hold:

{\small\par\smallskip
\CaseLine{\((W1) xyz \approx  xzy,\qquad (W2) xy^{2} \approx  xy.\)}
\smallskip}

They put every nonlinear word into a normal form consisting of its head followed by the sorted set of its tail variables. The head may also occur in this tail. Linear words remain distinct from these normal forms.

Theorem. The two identities W1, W2 and the following six absorption identities form a complete basis for A:

{\small\par\smallskip
\CaseLine{\((D)  xyz \leq  xy,\)}
\CaseLine{\((L)  xy \leq  xz+y,\)}
\CaseLine{\((C)  x^{2}yz \leq  xt+yz,\)}
\CaseLine{\((H)  xy \leq  x^{2}y+x,\)}
\CaseLine{\((N)  xyz \leq  x^{2}yz+yt,\)}
\CaseLine{\((V)  xyz \leq  x^{2}yz+y.\)}
\smallskip}

Here a\ensuremath{\leq}b abbreviates a+b\ensuremath{\approx}b. Reversing every word gives a basis for the multiplicative opposite.

\par\smallskip\noindent\textit{A three-part absorption criterion.}\quad 
Write a polynomial P as a nonempty sum of linear variables L and nonlinear words N. Let H be the set of heads of the words in N. For a nonlinear target w write \(h=head(w)\), \(T=tail-content(w)\), and \(C=\{h\}\)\ensuremath{\cup}T.

Lemma. A linear variable y satisfies y\ensuremath{\leq}P exactly when y\ensuremath{\in}L. A nonlinear w satisfies w\ensuremath{\leq}P exactly when:

{\small\par\smallskip
\CaseLine{\((1) h\in H;\)}
\noindent (2) L\ensuremath{\cap}C is nonempty, or some nonlinear right word has content contained in C;\par
\noindent (3) h\ensuremath{\in}L, or T\ensuremath{\cap}(H\ensuremath{\cup}L) is nonempty, or a right word with head h has tail content contained in T.\par
\smallskip}

Necessity. For a missing linear variable y, assign \(y=1\) and every other variable 2. Every nonlinear right word is 2, and all linear right summands are 2, while \(y=1\).

If (1) fails, assign \(h=3\) and every other variable 2. Every nonlinear right word has head 2 and hence value 2; the linear part is at most 3. The target is 0, or 4 if its tail contains only h. Neither is absorbed by 3.

If (2) fails, assign C the value 3 and outsiders 2. The target is 4. Every nonlinear right word contains an outsider: its value is 2 if its head is an outsider and 0 otherwise. Every linear summand is 2. Thus P\ensuremath{\leq}0, which does not absorb 4.

For failure of (3), we may assume (1), since its failure is already separated. Then h cannot belong to T, because h\ensuremath{\in}H and T\ensuremath{\cap}H is empty. Assign \(h=1\), all T-variables 3, and every remaining variable 2. The target is 1. A nonlinear right word headed h has a tail variable outside T and therefore has value 2; every other nonlinear right word has head 2 because H avoids T. Linear summands avoid both h and T and are 2. Hence \(P=2\).

Sufficiency. A target value 2 is least. If \(w=0\), the head is 0 or in K. By (1) a nonlinear right word with this head has value 0 or 4, either of which absorbs 0.

If \(w=4\), all C-variables lie in K. A contained nonlinear witness in (2) is 4. If instead (2) gives a linear g in C, its value is 3 or 4. Value 4 suffices directly; value 3 combines with the nonlinear head witness from (1), whose value is 0 or 4, to give 4.

If \(w=1\), the head is 1 and all T-variables lie in K. A linear h in (3) has value 1. A linear variable from T has value 3 or 4. A nonlinear right word with head in T has value 0 or 4. Each is at least 1. In the last alternative of (3), a right word with head h and tail contained in T is 1. These exhaust all target values. The criterion therefore holds without a variable or length bound.

The six absorption laws are sound by this criterion. Each has the required nonlinear head witness; (C) uses yz as a contained witness, and (L),(H),(V) use their indicated linear witnesses. In (N), the head y of yt lies in the target tail. The repeated-head targets automatically satisfy the intersection alternative of (3). Law (D) is also immediate from the word-value description.

\par\smallskip\noindent\textit{Unrestricted derivation.}\quad 
Assume a nonlinear absorption w\ensuremath{\leq}P satisfies the lemma, and normalize all nonlinear words. First, if there is an original h-headed word u whose tail content U is contained in T, the target follows by (D): append the missing tail variables to u, then use W1, W2. If there are no missing variables u already equals w. All appended factors used in (D) are nonempty.

We may thus assume this alternative is absent. Choose any original h-headed word \(a=hU\) using (1), with U a nonempty tail word. We first derive the word w* whose head is h and whose tail content is all of C, so that w* repeats its head.

If (2) has a linear witness g\ensuremath{\in}L\ensuremath{\cap}C, (D) applied to a gives h\ensuremath{{}^{2}}U\ensuremath{\leq}P after tail permutation. Then (L), with \(x=h\)\ensuremath{{}^{2}}, \(y=g\), \(z=U\), gives h\ensuremath{{}^{2}}g\ensuremath{\leq}h\ensuremath{{}^{2}}U+g\ensuremath{\leq}P. If instead the witness is a nonlinear word v with content contained in C, split \(v=yz\) into nonempty blocks and use (C):

{\small\par\smallskip
\CaseLine{\(h^{2}v \leq  hU+v \leq  P.\)}
\smallskip}

In either case the derived word has head h, repeats h, and uses only variables in C. Append every missing variable of C with (D) and normalize. This gives w*\ensuremath{\leq}P. If w already repeats h, its normal form is w*, and the proof is complete.

Otherwise h\ensuremath{\notin}T. Condition (3), since the direct contained-tail case has been excluded, supplies one of the following witnesses.

\textbullet{} If h\ensuremath{\in}L, apply (H) to obtain w\ensuremath{\leq}w*+h\ensuremath{\leq}P.

\textbullet{} If g\ensuremath{\in}T\ensuremath{\cap}L, put g first in the tail. If other tail variables remain, (V) gives w\ensuremath{\leq}w*+g. If \(T=\{g\}\), use (L) with \(x=h\), \(y=g\), \(z=hg\): the term h\ensuremath{{}^{2}}g is w*.

\textbullet{} Otherwise choose g\ensuremath{\in}T\ensuremath{\cap}H and a right word gV, with V nonempty. Write w\ensuremath{\approx}hgZ, with Z listing T minus \{g\} when that set is nonempty, and \(Z=g\) otherwise. The latter repetition is legitimate by W2 and avoids an empty substitution. Rule (N) gives w\ensuremath{\leq}w*+gV\ensuremath{\leq}P.

Thus every valid nonlinear absorption follows from the eight displayed identities. Linear absorptions are literal right summands. Distribute arbitrary terms, derive each summand absorption in both directions by this construction, and combine them to obtain every valid identity. This proves completeness for arbitrarily many variables, words and summands. Reversal gives the opposite assertion.

The proof does not assert minimality or independence of the basis, or hereditary finite basability. Finite regressions supplement the unrestricted argument; transfer to other catalogue types requires independently verified complete two-way HSP equality.

\subsection{A finite basis transfer with squared contexts}\label{case:42}
In the constant-free ai-semiring signature put N=\{0,\allowbreak{}1\} with maximum addition and multiplication constantly 0. In N a nonlinear word always has value 0. A literal x is absorbed by P exactly when x is a literal summand of P. This elementary observation also applies to arbitrary polynomial identities by taking the separate summands of both sides.

\par\smallskip\noindent\textit{A conditional transfer theorem.}\quad 
Suppose S has a finite relative basis and satisfies the following conditions. Valid word absorptions have c(w) contained in c(P). For each of the three nonempty choices of left and right contexts (both present, left only, right only), S satisfies

{\small\par\smallskip
\CaseLine{\(x \leq  [u^{2}] x [v^{2}],\)}
\CaseLine{\([u^{2}] x [v^{2}] \leq  x + [u] x [v].\)}
\smallskip}

Here context letters are distinct from x and one another, and the brackets at both occurrences of the same context mean the same choice. Then S x N is finitely based.

Fix any finite relative basis of S. Normalize its sides by distributivity into sums of nonempty words. The following finite construction specifies a basis, in addition to the ai-semiring axioms. All steps refer to ordinary equational consequences in the two-operation signature.

1. In each orientation of each normalized axiom take every summand v on the new side below the entire old polynomial Q. Include zero, one or two distinct fresh variables as fixed outside multiplicative contexts, multiplying both sides in the same locations. Still call the resulting absorption v \ensuremath{\leq} Q.

2. Make a finite list of word-absorption patterns. Substitute for each symbol of v, Q a sum of between one and max(1, occurrences in v) fresh distinct variables, and choose one monomial on the expanded left. Use the whole expanded right. A context symbol may be kept a singleton. These patterns are valid in S. In a concrete proof each retained occurrence in the chosen output may use a different selected image word. Every variable absent from the output needs just one arbitrary image word. Thus, after ordinary word substitution, these finitely many patterns contain the selected output and have only genuine old input words. This restriction of polynomial substitutions is carried out separately for each selected output.

3. In a pattern with a formal literal left z also include the version replacing z throughout by two fresh successive letters. This handles an image of length at least two by a split into nonempty pieces. In any pattern, a literal summand y on the right may independently be replaced by [\(s^{2}\)] y [\(t^{2}\)], using zero, one or two independent fresh context symbols for that summand. Call the modified right Q'. S satisfies Q \ensuremath{\leq} Q' by the first assumed context identity.

4. When the formal left is nonlinear keep its absorption below Q'. For the unsplit literal left z, skip the case in which z itself was a literal summand of Q. Otherwise select any occurrence U z V in a summand of Q. Such an occurrence exists by content coverage. Here U, V are actual formal subwords of that summand and at least one is nonempty. Include \(U^{2}\) z \(V^{2}\) \ensuremath{\leq} Q', omitting an absent block. The second assumed context identity puts this left side below z + U z V \ensuremath{\leq} Q \ensuremath{\leq} Q' in S. All the kept left sides are nonlinear, so all the resulting absorptions are also valid in N. A finite number of finite expansions and choices has been made.

To prove completeness, take a valid absorption w \ensuremath{\leq} P in S x N. If w is a literal it already occurs in P. Otherwise take a finite derivation of \(P=P+w\) in S. Work with the free-polynomial expansion of each term. Initially the original words are derivably below P. Maintain the following property at every step. Every nonlinear intermediate word is derivably below P in the new presentation. For each literal x occurring at that step, at least one word [\(A^{2}\)] x [\(B^{2}\)], possibly x itself, is derivably below P. All retained blocks are nonempty words.

In a replacement inside a one-hole context, a newly selected output word comes from a single summand on the new side and a fixed pair of outside words, each possibly absent. The whole corresponding old contribution occurred before the replacement. Use steps 1 and 2 for these choices. A concrete nonlinear output uses the split of step 3 when its formal left is a literal. Keep nonlinear input words as they are. For each input that is a literal, use its already derived representative [\(A^{2}\)] x [\(B^{2}\)]. Step 4 derives the nonlinear output. For a genuinely literal output, a matching literal among the inputs already provides its representative; otherwise the selected occurrence of z gives a new squared-context representative by step 4. If distinct formal letters acquire the same image, the argument still applies. The relevant nonliteral upper word containing the formal z in step 4 was itself an old word, since the images are nonempty. Untouched summands in additive contexts retain their earlier representatives. Semiring axioms preserve the free-polynomial expansion. This proves the invariant without restricting repetitions, new variables in a proof, substitution polynomials, or sizes of contexts. At the end w is nonlinear, so w itself is derivably below P. Applying this to all summands in both directions proves the theorem.

\par\smallskip\noindent\textit{Source 29 and the five-element algebras.}\quad 
Let F be the four-element flat semiring of a two-vertex directed path. Its elements are source loop e, target loop f, forward bridge a, and zero z. Nonzero multiplication concatenates matching endpoints and mismatches give z. Distinct elements add to z. This is order-four part I type 29, with (z, a, f, e) the labels (0,\allowbreak{}1,\allowbreak{}2,\allowbreak{}3) in its source table. The identification \((e,f,a,z)=(e_{00},e_{11},e_{01},\infty)\) gives \(F=E_\triangle\) in Theorem~\ref{thm:unit-paths}, which supplies a finite basis. In this algebra a word is e or f only on identical loops. It gives a precisely for a unique a between e's and f's. Deleting a squared block from a word with value a preserves that value. When both x and U x V have value a, U consists of e's and V of f's, so doubling these blocks preserves a. The loop bounds are immediate and z is greatest. These observations prove all six context laws. Assign z at a missing target letter and e everywhere else to force the required content coverage.

In 833 the labels (z, a, e, f) are (0,\allowbreak{}1,\allowbreak{}2,\allowbreak{}3), reversing the orientation of the labelled order-four presentation by interchanging loop names as necessary. More precisely the isomorphism from source labels (0,\allowbreak{}1,\allowbreak{}2,\allowbreak{}3) to the ordinary subalgebra of 833 is (0,\allowbreak{}1,\allowbreak{}3,\allowbreak{}2). The remaining element 4 is additive greatest. Map 4 to (z,1) and each ordinary letter to (its F value,0) in F x N. This is an injective homomorphism, and both projections are onto. It follows that Var(833)=Var(F x N). The transfer theorem therefore proves 833 finitely based. Full HSP in both directions transfers the same conclusion to 9881 and 10005.

For 13019 the ordinary (z, a, f, e) are (4,\allowbreak{}1,\allowbreak{}2,\allowbreak{}3) and the additional letter \(c=0\) is strictly above a and below z, incomparable to the loops. Every nonlinear word has exactly its value after projecting c to a: with either neighbor present the marked bridge loses its mark when the endpoints match, and a mismatch gives z. Consequently a nonlinear absorption holds precisely when it holds in F. For example, an upper value c projects to a, and a nonlinear target with projection a has value a \ensuremath{\leq} c. For other upper values the same conclusion follows from the table and the projection; a and c are the only values with the same projection. Thus all identities of F x N hold in 13019.

Add x \ensuremath{\leq} xx to the basis obtained for F x N. This law holds in 13019, where the two bridges square to z. Suppose a literal x \ensuremath{\leq} P holds without the literal x in P. It holds in F. If an ordinary assignment gives \(P=a\), it must give \(x=a\). Change that one variable to c. All nonlinear upper words retain their ordinary projected values; all actual upper literals are unchanged. Hence P is still a while x is c, a contradiction. No ordinary assignment gives \(P=a\). At either loop bound x must be that same loop, so xx \ensuremath{\leq} P follows in F, including the remaining greatest bound z. This is a nonlinear absorption and is derivable by the finite basis just constructed. Recover x by x \ensuremath{\leq} xx. Every other required absorption has already been derived. This proves a finite basis for 13019 on unrestricted alphabets. The finite presentations for these four types are given by construction from a finite basis for type 29; no minimal size or small enumeration is asserted.

\subsection{A directed-cycle obstruction for the endpoint semiring 836}\label{case:43}
Signature: \((+,\cdot)\) without constants, relative to the ai-semiring axioms. This proof concerns arbitrary numbers of variables. Finite computations are supplementary checks and are not its completeness argument.

\par\smallskip\noindent\textit{1. Two quotients.}\quad 
Write \(A=\{0, n, e, f, t\}\) in that order. Its addition and multiplication tables, encoded by rows, are

{\small\par\smallskip
\CaseLine{\(+\qquad multiplication\)}
\CaseLine{\ttfamily 00044\qquad 00000}
\CaseLine{\ttfamily 01044\qquad 00010}
\CaseLine{\ttfamily 00244\qquad 01200}
\CaseLine{\ttfamily 44434\qquad 00034}
\CaseLine{\ttfamily 44444\qquad 00044}
\smallskip}

Let \(Q=\{0, n, e, f\}\) have flat addition (the sum of distinct elements is \(0\)) and nonzero products \(e^2=e\), \(f^2=f\), \(en=nf=n\). Let \(D=\{0,\allowbreak{}1\}\) have addition maximum and multiplication minimum. The maps

\CaseDisplay{q:A\longrightarrow Q,\quad(q(0),q(n),q(e),q(f),q(t))=(0, n, e, f,0),}

\CaseDisplay{d:A\longrightarrow D,\quad(d(0),d(n),d(e),d(f),d(t))=(0,\allowbreak{}0,\allowbreak{}0,\allowbreak{}1,\allowbreak{}1)}

are onto homomorphisms, and \((q, d)\) is injective. This follows directly from the displayed tables: \(q\) sends precisely \(0, t\) to the flat zero, while \(d^{-1}(1)=\{f, t\}\) is closed under multiplication and its complement is multiplicatively absorbing; its complement is also closed under addition. Thus

\CaseDisplay{\operatorname{Var}(A)=\operatorname{Var}(Q, D).\quad\textup{(\ensuremath{1})}}

For a nonempty word, its value in \(Q\) is \(e\) if all letters are \(e\), \(f\) if all are \(f\), \(n\) if its occurrence-value sequence is \(e^*nf^*\), and \(0\) otherwise. In particular, a word containing two occurrences evaluated as \(n\) is zero. A polynomial has value \(n\) exactly when every summand has value \(n\). The symbol \(0\) in this paragraph is an element label, not a constant in the signature.

Every valid identity of \(A\) preserves total variable content. Indeed, in \(Q\) assign a variable present on just one side to \(0\), and all other variables to \(e\). Also, if \(A\models w\leq P\), where \(\leq\) means \(w+P\approx P\), some summand \(u\) of \(P\) has \(c(u)\subseteq c(w)\). Otherwise the assignment in \(D\) taking the variables of \(w\) to \(1\) and the others to \(0\) separates the absorption. We shall repeatedly use these two necessary conditions.

\par\smallskip\noindent\textit{2. Cycles and rigidity after cutting an edge.}\quad 
Fix \(r\geq3\), with subscripts on \(x_i\) read modulo \(r\). Take \(3r\) distinct variables \(a_i, b_i, x_i\) and put

\CaseDisplay{p_i=a_i x_i x_{i+1}b_i,\qquad P_E=\sum_{i\in E}p_i,\qquad P_r=\sum_{i=0}^{r-1}p_i.}

Here \(E\) is a nonempty subset of \(\{0,\ldots, r-1\}\). The contents \(c(p_i)\) are pairwise incomparable. Distinct words \(p_i\) have different first letters and different last letters.

Lemma 1 (path evaluations). Suppose \(E\) is proper. Every variable occurring in \(P_E\) can be assigned \(n\) in some evaluation with \(P_E=n\). More precisely, every contiguous nonempty factor of any \(p_i\), \(i\in E\), can have value \(n\) while \(P_E=n\).

Proof. The directed graph with edges \(x_i\to x_{i+1}\) for \(i\in E\) is a disjoint union of directed paths and isolated vertices. To make a specified core variable \(x_j\) equal \(n\), give the earlier vertices of its path value \(e\), \(x_j\) value \(n\), and the later vertices value \(f\). Give vertices of other paths value \(e\). Every selected edge then has one of the pairs \((e, e),(e, n),(n, f),(f, f)\). For these four cases choose \((a_i, b_i)\) respectively as \((e, n),(e, f),(e, f),(n, f)\). Every \(p_i\) has value \(n\).

To give a specified \(a_i\) value \(n\), instead give all core variables value \(f\) and all \((a_j, b_j)\) values \((n, f)\). To give a specified \(b_i\) value \(n\), give all core variables value \(e\) and all \((a_j, b_j)\) values \((e, n)\). Finally choose any letter occurring in the specified factor and use its indicated evaluation. The complete word \(p_i\) has sequence \(e^*nf^*\), and the factor contains that unique \(n\), so the factor also has value \(n\). These are evaluations in \(Q\); by (1) they suffice for the uses below. \(\square\)

Lemma 2 (proper-subsum rigidity). If \(E\) is nonempty and proper, then

\CaseDisplay{A\models w\leq P_E\quad\Longleftrightarrow\quad w=p_i\text{ as a literal word for some }i\in E.\quad\textup{(\ensuremath{2})}}

Proof. Necessity of the two content conditions above gives \(c(w)\subseteq c(P_E)\) and \(c(p_i)\subseteq c(w)\) for some \(i\in E\). Give all core variables and all \(a_j\) value \(e\), and every \(b_j\) value \(n\). Since \(P_E=n\), the target must contain exactly one occurrence of one \(b_j\), at its final position. The evaluation with all core variables and all \(b_j\) equal \(f\), and every \(a_j=n\), similarly forces exactly one occurrence of one \(a_j\), at its first position. The contained-content witness forces these two private indices to be the same \(i\). Consequently \(w=a_i v b_i\), where \(v\) uses only core variables and includes both \(x_i, x_{i+1}\).

If \(v\) contains a different core variable \(x_j\), use the core-pivot evaluation of Lemma 1 with \(x_j=n\). Edge \(i\) does not meet the pivot, so its two core values are either both \(e\) or both \(f\). The construction then makes respectively \(b_i=n\) or \(a_i=n\). Thus \(w\) contains at least two occurrences valued \(n\), and is zero, although \(P_E=n\): a contradiction.

It follows that \(v\) uses only \(x_i, x_{i+1}\). Pivoting at \(x_i\) forces it to occur once and before all occurrences of \(x_{i+1}\). Pivoting at \(x_{i+1}\) forces that variable to occur once. Hence \(v=x_ix_{i+1}\) and \(w=p_i\). The reverse implication is reflexivity of a summand absorption. \(\square\)

Lemma 3 (a valid cycle identity). Put \(w_r=a_0x_0^2x_1b_0\). Then

\CaseDisplay{A\models P_r\approx P_r+w_r.\quad\textup{(\ensuremath{3})}}

Proof. The identity holds in \(D\) because \(c(w_r)=c(p_0)\). In \(Q\), it is automatic when \(P_r=0\). If \(P_r=e\) or \(f\), all its variables have that value and so does \(w_r\). Suppose \(P_r=n\). Write the nonzero values of the core variables as intervals \(e=(0,\allowbreak{}0)\), \(n=(0,\allowbreak{}1)\), \(f=(1,\allowbreak{}1)\). Every selected edge in a nonzero \(p_i\) requires \(R_{x_i}=L_{x_{i+1}}\), and each variable has \(L_{x_i}\leq R_{x_i}\). Around the directed cycle,

\CaseDisplay{L_{x_0}\leq R_{x_0}=L_{x_1}\leq\cdots\leq R_{x_{r-1}}=L_{x_0}.}

All endpoints agree. Thus all core variables are \(e\), or all are \(f\). Duplicating \(x_0\) preserves the value of \(p_0\), so \(w_r=n\). Equation (1) completes the proof. \(\square\)

\par\smallskip\noindent\textit{3. Arbitrary polynomial substitutions in bounded-variable identities.}\quad 
Polynomials in this section are nonempty finite sets of nonempty words, with addition interpreted as union and multiplication as concatenation followed by union. This is the free ai-semiring before imposing extra identities. Substitutions can send variables to arbitrary nonempty polynomials; they are not restricted to words.

Lemma 4 (substitution rigidity). Let \(U\approx V\) be any identity valid in \(A\) using fewer than \(r\) variables. Let \(\sigma\) be a polynomial substitution. If the expansion \(\sigma(U)\) is a nonempty subset of the word set of \(P_r\), then \(\sigma(U)\) and \(\sigma(V)\) are exactly the same word set.

Proof. If \(\sigma(U)=P_E\) for a proper subset \(E\), every word of \(\sigma(V)\) is validly absorbed in \(P_E\), and Lemma 2 puts it in \(P_E\). Conversely, equality in \(D\) forces every \(p_i\) in \(P_E\) to remain: assign precisely its content to \(1\), and all other variables to \(0\). Pairwise incomparability of the contents separates any missing summand. Thus equality of the word sets holds in this case. It remains to consider \(\sigma(U)=P_r\).

Split the source variables into \(L\), those occurring as length-one summands of \(U\), and the other variables \(N\). A variable occurring in a nonlinear summand of \(U\) has a singleton word image under \(\sigma\). Indeed, fixing the image choices at the other positions, two distinct alternatives at its position would produce distinct words of \(P_r\) sharing a nonempty initial or final context. They would share a first or last letter, impossible for two distinct \(p_i\). Free-word cancellation shows the alternatives must coincide.

No variable occurs both linearly and nonlinearly in \(U\): its linear occurrence would require every image word to have length four, while its nonlinear occurrence would then produce a word longer than four. Hence \(L, N\) are disjoint and the nonlinear part \(U_N\), when nonempty, uses precisely \(N\). Each variable in \(L\) has as image a nonempty set of original words \(p_i\); each variable in \(N\) has a single word as image.

Let \(\sigma(U_N)=P_F\) if \(U_N\) is nonempty. There are at most \(|N|\) distinct words in this image. To see this, the image of the first source variable of a source monomial determines the first letter \(a_i\) of its image, and that first letter determines the entire original word \(p_i\). There are at most \(|N|\) choices of a first source variable. Since \(|N|<r\), \(F\) is proper.

Every variable \(y\in N\) can have value \(n\) in a \(Q\)-evaluation with \(U_N=n\). Its singleton image is a contiguous factor of some \(p_i\) with \(i\in F\). Lemma 1 gives an evaluation of the image alphabet in which that factor and all of \(P_F\) are \(n\). Evaluation of the singleton images induces the required evaluation of the source variables. Also \(U_N=n\) has at least one evaluation.

Consider any word \(v\) of \(V\). Total-content preservation ensures its variables occur in \(U\). If \(v\) contains a variable \(z\in L\) and has length at least two, choose any second occurrence. If its variable is in \(N\), choose the preceding evaluation making that variable \(n\); if it is in \(L\), choose any evaluation making \(U_N=n\) (with no nonlinear variables needed if \(U_N\) is absent). Assign every variable of \(L\) to \(n\). Then \(U=n\), whereas \(v\) has at least two \(n\) occurrences and is zero. This contradicts \(Q\models v\leq U\). Therefore a word of \(V\) containing a variable from \(L\) must be just one such variable; its image already lies in \(P_r\).

Every remaining word \(v\) uses only \(N\). In this case \(U_N\) is nonempty and \(A\models v\leq U_N\). For clarity, this last assertion follows separately on the two generators in (1). For \(D\), any contained-content witness for \(v\leq U\) cannot be a linear variable of \(L\), so it belongs to \(U_N\). For \(Q\), if \(U_N=0\) absorption is automatic; if \(U_N=e\) or \(f\), all variables of \(N\) have that same value and so does \(v\); if \(U_N=n\), extend the assignment by giving all variables in \(L\) value \(n\) and apply the validity of \(v\leq U\).

Substituting now gives \(A\models\sigma(v)\leq P_F\). The substitution is by words on \(N\), so \(\sigma(v)\) is a single word. Lemma 2 makes it exactly one of the original words. We have proved \(\sigma(V)\subseteq P_r\). Finally equality in \(D\), and the pairwise incomparable original contents, force \(\sigma(V)=P_r\). \(\square\)

\par\smallskip\noindent\textit{4. Nonfinite basability, including arbitrary term contexts.}\quad 
Theorem. The semiring \(A\) is nonfinitely based.

Proof. Suppose a finite basis \(\Sigma\) existed. Choose \(k\) bounding the number of variables in each of its identities, and choose \(r>\max(k,2)\). We show that every equational consequence of \(\Sigma\) starting with \(P_r\) has the same expanded word set \(P_r\).

It suffices to inspect one replacement of a substituted basis identity \(U\approx V\) in an arbitrary one-hole term context. Distribute the context: its terms involving the hole are a finite nonempty collection of expressions \(s\sigma(U)t\), with \(s, t\) words or omitted contexts, together with summands not involving the hole. Omitted contexts are not unit constants. Since the current expansion is \(P_r\), each such expression expands to a nonempty subset of \(P_r\).

If both contexts are omitted, Lemma 4 says its expansion is unchanged. If at least one context is nonempty, all its expanded words share a first or a last letter. They must therefore be a single original word \(p_i\). The valid contextual identity \(s\sigma(U)t\approx s\sigma(V)t\) and Lemma 2 for \(E=\{i\}\) imply that the latter expansion is also exactly \(p_i\). Thus each context contribution is unchanged. The ai-semiring axioms themselves preserve the free-polynomial expansion. This proves the asserted invariance by induction on any finite equational derivation.

But (3) is an identity of \(A\), and \(w_r\) is not an original word of \(P_r\): it has length five instead of four. It cannot follow from \(\Sigma\), contradicting completeness. Consequently no finite basis exists. \(\square\)

The argument handles arbitrary polynomial substitutions, repeated source variables, both orientations of equations, and arbitrary contexts. Its essential use of the lattice quotient is to prevent deleting the private-endpoint support of an original summand. The earlier endpoint semantic decision procedure alone did not establish this obstruction.

\par\smallskip\noindent\textit{5. Catalogue transfer and verification boundary.}\quad 
The algebra just considered is \(\A{836}\). Proposition~\ref{prop:lower-fibres} in both directions gives the equal-variety class \(834,836,848,929,1665,7902,8285\) in case~\ref{grp:77}. The obstruction therefore holds for each member.

\subsection{Directed-cycle nonfinite-basis transfer: eleven catalogue types}\label{case:44}
This note supplements the unrestricted proof in Section~\ref{case:43}; it does not infer nonfinite basability from a proper NFB subvariety alone.

\par\smallskip\noindent\textit{1. Transfer theorem.}\quad 
Let \(A\) be the endpoint semiring \(836\) of Section~\ref{case:43}. For each \(r\geq3\) let

\CaseDisplay{P_r=\sum_{i\in\mathbb Z/r\mathbb Z}a_i x_i x_{i+1}b_i,
\qquad w_r=a_0x_0^2x_1b_0,}

with all \(3r\) variables distinct. If a semiring \(T\) satisfies

1. \(A\in\operatorname{Var}(T)\);

2. \(T\models P_r\approx P_r+w_r\) for every \(r\geq3\);

then \(T\) is nonfinitely based.

Proof. Every identity of \(T\) is an identity of \(A\). Given a proposed finite basis of \(T\) with at most \(k\) variables in each identity, choose \(r>\max(k,2)\). The substitution-rigidity and context argument in Section~\ref{case:43} applies to every such identity: starting from \(P_r\), no derivation can change its literal expanded word set. Yet the displayed identity holds in \(T\) and adds a word of length five. This contradiction proves the assertion. The second hypothesis is indispensable. \(\square\)

\par\smallskip\noindent\textit{2. An all-length finite certificate for hypothesis 2.}\quad 
For a finite semiring \(T\) put

\CaseDisplay{V(x, y)=\{axyb:a, b\in T\}.}

Initialize states \((s, z, p, w)\) with

\CaseDisplay{s=x_0,\quad z=x_1,\quad p=a_0x_0x_1b_0,\quad w=a_0x_0^2x_1b_0,}

for every choice of the four indicated elements. An extension by one edge sends

\CaseDisplay{(s, z, p, w)\longmapsto(s, y, p+v, w),\qquad y\in T,\ v\in V(z, y).}

After at least two open edges have been processed, closing by \(v\in V(z, s)\) gives precisely the right value \(p+v\) of a cycle with at least three edges; check \(w+(p+v)=p+v\). Take the complete closure under extensions of the states with two open edges.

Induction on path length proves that these states represent every evaluation of every open path of length at least two. Conversely each reached state is represented by such an evaluation, because the private variables on an added edge are independent of the preceding ones. There are at most \(|T|^4\) states. Thus complete closure and all closing checks certify every cycle length, rather than extrapolating from a fixed collection of small cycles. Two implementations independently reproduce the state sets, transitions and closing checks.

\par\smallskip\noindent\textit{3. Results and strictness.}\quad 
The following types form two equal-variety classes.

\noindent Representative\enspace\textbar\enspace All entries\enspace\textbar\enspace Evidence\par
\noindent 836\enspace\textbar\enspace 834, 836, 848, 929, 1665, 7902, 8285\enspace\textbar\enspace Unrestricted cycle proof and independently regenerated two-way HSP certificates\par
\noindent 839\enspace\textbar\enspace 839, 11370, 11669, 12305\enspace\textbar\enspace \(A\in\operatorname{Var}(T)\) and all-length cycle certificates, using the transfer theorem\par
The inclusion for the second row is strict: \(A\) satisfies \(x^2y^2\approx xyx\), whereas each of those four targets fails it. The separating assignments can be found by testing the displayed identity on the ordered pairs of each table. Proposition~\ref{prop:lower-fibres} gives the inclusions.

\par\smallskip\noindent\textit{4. Separation of basis status from limit-variety status.}\quad 
The seven entries in the 836 class are now NFB, but this alone does not decide whether their variety is a limit variety. They remain in the limit-generation candidate catalogue with a separately recorded NFB basis status.

The four entries in the 839 class have a proper NFB subvariety generated by 836, so they are excluded from limit-variety generation. Their own NFB status uses the transfer theorem, not just that proper inclusion.

\subsection{An endpoint-graph finite basis for a five-element semiring}\label{case:45}
The following argument is due to Aifa Wang and Lili Wang.

The algebra below is the table recorded for catalog representative 837. The theorem is stated for the displayed tables. Catalog identification, the label of the opposite, and transfer to catalog equality classes must be checked against the saved certificates before consolidated counts change.

\par\smallskip\noindent\textit{1. Algebra and finite basis.}\quad 
Let A = \{0,\allowbreak{}1,\allowbreak{}2,\allowbreak{}3,\allowbreak{}4\}, with operations given by the following tables, in that order. These labels are elements, not constants in the signature.

{\small\par\smallskip
\CaseLine{\(+ | 0 1 2 3 4\qquad * | 0 1 2 3 4\)}
\CaseLine{\ttfamily 0 | 0 0 0 4 4\qquad 0 | 0 0 0 0 0}
\CaseLine{\ttfamily 1 | 0 1 0 4 4\qquad 1 | 0 0 0 1 0}
\CaseLine{\ttfamily 2 | 0 0 2 4 4\qquad 2 | 0 1 2 0 0}
\CaseLine{\ttfamily 3 | 4 4 4 3 4\qquad 3 | 4 4 4 3 4}
\CaseLine{\ttfamily 4 | 4 4 4 4 4\qquad 4 | 4 4 4 4 4}
\smallskip}

Work relative to the axioms for additively idempotent semirings in the constant-free signature (+, *). Write u \ensuremath{\leq} v for u + v = v.

Include the following three identities:

{\small\par\smallskip
\CaseLine{\((B)\qquad xyx = x^{2} y^{2},\)}
\CaseLine{\((D1)  x \leq  x y^{2},\)}
\CaseLine{\((D2)  xy \leq  x z^{2} y.\)}
\smallskip}

Also include every optional-context instance of these five schemes:

{\small\par\smallskip
\CaseLine{\((M)  A x D\qquad \leq  A x B + C x D,\)}
\CaseLine{\((R)  A x t F\qquad \leq  A x y B + C z y D + E z t F,\)}
\CaseLine{\((S)  A x z D\qquad \leq  A x y B + y C + z D,\)}
\CaseLine{\((T)  D z y B\qquad \leq  A x y B + C x + D z,\)}
\CaseLine{\((E)  A x\qquad \leq  A x y B + C z y D + E z.\)}
\smallskip}

In each scheme independently, each capital letter is either omitted everywhere it occurs or replaced by its own new variable. All unomitted capital letters and displayed lower-case letters are distinct formal variables. Ordinary substitutions may subsequently identify variables or replace them by nonempty words or terms. Omission is a notational device; it introduces neither an empty term nor a unit constant. Every summand and each left-hand word contains a lower-case variable and is nonempty.

There are 16, 64, 16, 16, and 32 optional-context instances, respectively. Together with (B), (D1), and (D2), these give at most 147 identities.

Theorem. These identities form a basis for the identities of A, relative to the additively idempotent semiring axioms. Reversing all words gives a basis for its multiplicative opposite.

\par\smallskip\noindent\textit{2. Evaluation and endpoint graphs.}\quad 
The operation tables satisfy the semiring axioms. This is a finite table verification: associativity of both operations and the two distributive laws involve 125 triples each; commutativity and idempotence of addition involve 25 and 5 instances.

For a nonempty sequence of assigned elements, the value of its word is:

\textbullet{} 1 precisely when the sequence has the form 2* 1 3*, with exactly one 1;

\textbullet{} 2 precisely when all its entries are 2;

\textbullet{} 3 precisely when all its entries are 3;

\textbullet{} 4 precisely when its first entry is 4, or its first entry is 3 and some entry is different from 3;

\textbullet{} 0 in all remaining cases.

Here the stars in 2* and 3* denote repetitions in a language, not multiplication. The rule follows immediately by following the rows of the multiplication table: states 0 and 4 are absorbing under right multiplication; 3 stays 3 only on 3; 2 stays 2 on 2 and goes to 1 on 1; 1 stays 1 only on 3. This also covers words of length one.

Let P be a nonempty finite sum of nonempty words. Write C(P), H(P), and T(P) for its content, its set of heads, and its set of final letters. Construct a directed endpoint graph G(P) as follows. Start with vertices s, t and \(l_{x}\), \(r_{x}\) for x in C(P). Identify:

1. s with \(l_{h}\) for h in H(P);

2. t with \(r_{f}\) for f in T(P);

3. \(r_{x}\) with \(l_{y}\) whenever xy occurs consecutively in a summand of P.

After these identifications, put a directed edge labelled x from \(l_{x}\) to \(r_{x}\) for every x in C(P). Parallel edges and loops are permitted. Every vertex lies on a directed walk from s to t: every occurrence of a variable belongs to a summand, and each summand is such a walk. Let Q(P) be the graph obtained by identifying the strongly connected components of G(P).

Encode the elements 2, 1, 3 as endpoint pairs (0,\allowbreak{}0), (0,\allowbreak{}1), (1,\allowbreak{}1), respectively. A word has value 1 precisely when its pairs concatenate from 0 to 1. Consequently P has value 1 precisely when every variable has one of these three values and its endpoint labelling satisfies

{\small\par\smallskip
\noindent s = 0,  t = 1,  \(l_{x}\) \ensuremath{\leq} \(r_{x}\) for all x,\par
\smallskip}

with the graph identifications respected. Indeed a sum equals the minimal element 1 exactly when all its summands equal 1. Thus these evaluations are exactly the order-preserving Boolean labellings of G(P) that send s to 0 and t to 1.

If s and t lie in the same strongly connected component, there is no such labelling. In this case every vertex is in that component, since s reaches each vertex and each vertex reaches t. Otherwise two vertices have equal labels in every such labelling exactly when they are strongly connected. For the nontrivial direction, if u does not reach v, label the vertices reachable from u by 1 and all others by 0. This labels u and v differently, labels s by 0 and t by 1, and is order-preserving. For two different components at least one of the two reachabilities fails, so this argument always separates them.

Absorption criterion. For a nonempty word w,

{\small\par\smallskip
\CaseLine{\(A satisfies w \leq  P\)}
\smallskip}

if and only if:

1. C(w) is contained in C(P);

2. the head of w belongs to H(P);

3. the labelled word w is a directed walk from s to t in Q(P).

In condition 3 the first endpoint is automatically s by condition 2. Equivalently each consecutive pair xy of w has \(r_{x}\) and \(l_{y}\) strongly connected in G(P), and \(r_{f}\) is strongly connected to t for the last letter f of w. The criterion includes the case where Q(P) has one vertex.

Proof. If condition 1 fails, give every variable of C(P) value 3 and every outside variable value 2. Then P = 3 and w is not 3, so w is not below P. If condition 2 fails, give the head of w value 4 and every other variable value 2. No summand of P starts with 4, so P is 0 or 2, whereas w = 4.

Now assume 1 and 2. If P = 4 the desired inequality is automatic. If P = 2 or 3, all variables of C(P) have that value, and hence so does w. If P = 0, no head in H(P) has value 3 or 4: a head of value 4 would force P = 4; a head of value 3 gives a summand of value 3 or 4, whose sum cannot be 0. Thus w has value 0, 1, or 2, all below 0. The only remaining case is P = 1, which requires w = 1. By the endpoint labelling argument this holds under every such evaluation precisely when condition 3 holds. If there are no such evaluations, all vertices form one strongly connected component and condition 3 holds automatically. This proves the criterion. \(\square\)

\par\smallskip\noindent\textit{3. Soundness of the displayed basis.}\quad 
The two words in (B) have the same head and content and neither can take value 1. The evaluation rule shows that their other values agree.

For (D1) and (D2), if the right word takes value 1, its squared block cannot contain the unique entry 1. Deleting that block therefore leaves a word of value 1. Content containment and the unchanged head, followed by the absorption criterion, prove these inequalities.

For each optional-context scheme the target has content contained in the right polynomial and its head is a head of a right summand. It is already a walk in G(P), before taking strongly connected components:

\textbullet{} (M) switches suffixes immediately after the common edge x.

\textbullet{} (R) uses \(r_{x}\) = \(l_{y}\) = \(r_{z}\) = \(l_{t}\).

\textbullet{} (S) uses \(r_{x}\) = \(l_{y}\) = s = \(l_{z}\).

\textbullet{} (T) uses \(r_{z}\) = t = \(r_{x}\) = \(l_{y}\).

\textbullet{} (E) uses \(r_{x}\) = \(l_{y}\) = \(r_{z}\) = t.

The reused prefix and suffix supply the remaining parts of each walk. These arguments remain valid when any optional context is omitted. The absorption criterion proves soundness of all displayed identities.

\par\smallskip\noindent\textit{4. Deriving every walk in the endpoint graph.}\quad 
Fix P. Call a word available if w \ensuremath{\leq} P is derivable from the proposed basis. All original summands are available. Applying a scheme to available words yields another available word, by transitivity and idempotent addition. All such words have content in C(P), and their heads belong to the original H(P).

Keep one available witness for each known adjacent pair xy and for each known final letter x. Initially these are the original pairs and final letters. Let R be the relation of known pairs and F the known final letters; retain H = H(P) throughout. Schemes (R), (S), (T), (E) implement the following finite closure rules:

{\small\par\smallskip
\CaseLine{\((R)  (x, y), (z, y), (z, v) in R\qquad => (x, v) in R;\)}
\CaseLine{\((S)  (x, h) in R, h, k in H\qquad => (x, k) in R;\)}
\CaseLine{\((T)  (f, y) in R, f, g in F\qquad => (g, y) in R;\)}
\CaseLine{\((E)  (x, y), (z, y) in R, z in F\qquad => x in F.\)}
\smallskip}

For example, take available words AxyB, CzyD, EzvF witnessing the three premises of (R). The scheme adds AxvF. The other rules are witnessed exactly by their displayed schemes. Since there are only finitely many pairs and final letters, this procedure terminates.

At termination:

{\small\par\smallskip
\noindent (x, y) in R  iff \(r_{x}\) and \(l_{y}\) were identified in G(P),\par
\noindent x in F      iff \(r_{x}\) and t were identified in G(P).\par
\smallskip}

Here only the endpoint identifications are considered, not directed strong connectivity. To verify this closure assertion, form the bipartite graph whose left vertices are all \(r_{x}\) and s, and whose right vertices are all \(l_{y}\) and t. Its initial edges are exactly the original pairs \(r_{x}\)--\(l_{y}\), the head edges s--\(l_{h}\), and the final edges \(r_{f}\)--t. Its connected components are the endpoint identifications.

The rectangle rule (R) completes every component of the pair relation to a complete bipartite relation. Rule (S) joins the components incident with head columns, adding isolated head columns when necessary. Rule (T) joins components incident with final rows, adding isolated final rows when necessary. Rule (E) marks every row in such a component as final. More explicitly, at a fixed point any two intersecting nonempty rows coincide by (R). If a row meets H it contains every member of H by (S). All final rows coincide by (T), and any row intersecting a nonempty final row is final by (E). Therefore a path in the augmented bipartite graph cannot connect a row and a column absent from that row, or connect a nonfinal row to t. The reverse inclusions follow from soundness of the four closure rules. This proves the assertion, including components consisting entirely of heads or entirely of finals.

Now let w be any nonempty s-to-t walk in G(P) whose head is in H(P). Start with an original summand with that head. Suppose an available word has the first i letters of w as a prefix and the i-th letter is x. If the next letter is y, choose the available witness CxyD for (x, y). Scheme (M), applied to the current word AxB and Cx(yD), adds AxyD. Here a nonempty context such as yD is substituted for a context variable; an empty context is handled by the corresponding optional instance. Repeat until the full target word is a prefix. Its final letter f lies in F. Splice with an available witness ending in f to discard the remaining suffix. This derives w \ensuremath{\leq} P. The same argument with just the final splice covers a one-letter walk.

\par\smallskip\noindent\textit{5. Passing to strong components and completeness.}\quad 
Let w satisfy the absorption criterion. At a gap xy in w the vertices \(r_{x}\) and \(l_{y}\) lie in the same strongly connected component. If they are equal, insert nothing. Otherwise choose a nonempty directed path from \(r_{x}\) to \(l_{y}\) with label u, and a nonempty return path with label v. The walk labelled uvu still connects \(r_{x}\) to \(l_{y}\). Insert it at this gap. Apply the same construction after the last letter to connect its end to t. This produces an actual s-to-t walk W in G(P), with the same head as w. Section 4 gives W \ensuremath{\leq} P.

It remains to show w \ensuremath{\leq} W syntactically. At an internal gap, (D2) allows insertion of \(u^{2}\) and then \(v^{2}\); (B), in reverse with x = u and y = v, changes \(u^{2}\) \(v^{2}\) to uvu. At the final gap, (D1) allows the same two insertions, followed by (B). Multiplication is monotone for the additive order, so these insertions may be applied successively inside longer words. Every prefix used in (D1) or (D2) is nonempty; every suffix used in (D2) is nonempty. No unit constant is needed. Thus w \ensuremath{\leq} W \ensuremath{\leq} P is derivable.

By distributivity every term is a nonempty finite sum of words. If an identity P = Q holds in A, each summand of P is below Q and is therefore derivably below Q by the preceding argument. Hence P \ensuremath{\leq} Q is derivable; the symmetric argument gives Q \ensuremath{\leq} P. This proves the theorem for arbitrarily many variables. Reversal of all words gives the assertion about the multiplicative opposite. \(\square\)

\subsection{Finite bases for quadratic polynomials determined by directed paths}\label{case:46}
All algebras have the signature of additively idempotent semirings without constants. The four labelled types are 838, 11683, 12306 and 13018. For the first three the following two identities, with the semiring axioms, are a complete basis. For 13018 include the third identity.

{\small\par\smallskip
\CaseLine{\(xyz = xy + yz\)}
\CaseLine{\(xz \leq  xy + yz\)}
\noindent x \ensuremath{\leq} xx                     (13018 only)\par
\smallskip}

Here absorption w \ensuremath{\leq} P abbreviates \(w+P=P\). Direct table evaluation on all assignments verifies the listed identities. The first identity reduces every word of length at least two to the sum of the products of its consecutive pairs. This follows by applying it with the prefix before the penultimate letter as its first nonempty substitution and then inducting. Expand any polynomial this way. Its data are a set of literal variables and a directed graph with an arc u to v for every original quadratic summand uv. Loops are arcs. Vertices include those of the target, including isolated vertices.

A quadratic xy is absorbed precisely when there is a directed walk of positive length from x to y. Sufficiency follows by induction on walk length using the second identity. For 838, 11683 and 12306 a literal x is absorbed precisely when it is an original literal. For 13018 it suffices also that x lies on a positive directed closed walk: use the path proof for xx followed by the third identity.

We give separating assignments proving necessity for quadratics and literals. For a bound b let I be its principal additive ideal. Use the following low, middle and high values. All three values lie in I. Every directed edge between nondecreasing values lies in the relation uv in I except the pair (middle, middle). Also (high, low) fails this relation. Only the successful nondecreasing nonmiddle-loop pairs and these two failures are needed.

\noindent Type\enspace\textbar\enspace b\enspace\textbar\enspace low\enspace\textbar\enspace middle\enspace\textbar\enspace high\par
\noindent 838\enspace\textbar\enspace 1\enspace\textbar\enspace 2\enspace\textbar\enspace 1\enspace\textbar\enspace 3\par
\noindent 11683\enspace\textbar\enspace 0\enspace\textbar\enspace 2\enspace\textbar\enspace 0\enspace\textbar\enspace 1\par
\noindent 12306\enspace\textbar\enspace 0\enspace\textbar\enspace 1\enspace\textbar\enspace 0\enspace\textbar\enspace 2\par
\noindent 13018\enspace\textbar\enspace 1\enspace\textbar\enspace 3\enspace\textbar\enspace 1\enspace\textbar\enspace 2\par
If x and y are distinct and y is unreachable from x, give high to x and all its successors, low elsewhere. All upper summands are in I, but xy is outside I. If \(x=y\) and there is no positive closed walk at x, give x the middle value, its strict successors high and the remaining vertices low. No upper edge can descend, and there is no middle loop. All literals are in I but xx lies outside. This proves the quadratic criterion for every directed graph.

For 838 with no original literal x, assign 4 to x and 2 to all other variables. Every upper value is in the ideal below 0 and x is outside. For 11683 and 12306 assign 2 to x and 3 to everything else. Every quadratic, including xx, has value 3 on these two inputs and all upper literals are 3. The target 2 is not below 3. For 13018 suppose x is neither an original literal nor on a positive closed walk. Give x the value 0, its strict successors 2, the other vertices 3. The values 3*0,\allowbreak{}0*2,\allowbreak{}3*2,\allowbreak{}2*2,\allowbreak{}3*3 all lie below 1, and all upper literals also lie below 1. The target 0 is outside that ideal. This proves all literal claims. The separators use only actual table elements.

We have derived all valid reduced absorptions with no bounds on the original variables, walk lengths, word lengths or polynomials. Apply those derivations to every summand in each direction of a polynomial equality and undo the reductions. This proves completeness.

\subsection{A ten-identity basis from directed paths in a two-point interval semiring}\label{case:47}
Work in the constant-free ai-semiring signature, relative to the usual ai-semiring axioms.

\par\smallskip\noindent\textit{1. Interval core and top extension.}\quad 
Let \(F=\{e, n, f,0\}\). Represent \(e=(0,\allowbreak{}0)\), \(n=(0,\allowbreak{}1)\) and \(f=(1,\allowbreak{}1)\) by intervals with endpoints in \(\{0,\allowbreak{}1\}\). Their sum is their interval hull; the extra element \(0\) is additively greatest and multiplicatively absorbing. For nonzero intervals set

\CaseDisplay{[i, j][k, l]=\begin{cases}[i, l],&j\leq k,\\0,&j>k.\end{cases}}

Thus \(e^2=e\), \(f^2=f\), \(ef=en=nf=n\), and the other products of nonzero core elements are zero. Add an additive top \(t\), with \(tx=t\) for every element and \(xt=0\) for \(x\in F\), giving \(A\). Catalogue 840 is \(A\) under \((0, n, e, f, t)=(0,\allowbreak{}1,\allowbreak{}2,\allowbreak{}3,\allowbreak{}4)\).

For associativity of the interval multiplication, three nonzero factors yield a nonzero product in either association precisely when the successive boundary inequalities \(j\leq k\) and \(l\leq m\) both hold; both products then have the first left and last right endpoints. For distributivity on the right, multiplying by the hull of two intervals is nonzero precisely when the first factor's right endpoint is at most both left endpoints. If either product is zero their sum is zero. Otherwise the two products have the same left endpoint and their hull has the maximum of the right endpoints, agreeing with multiplication by the hull. Left distributivity is dual. The absorbing zero and the one-sided top extend these properties directly.

\par\smallskip\noindent\textit{2. Ten identities.}\quad 
The following are the proposed basis:

\CaseDisplay{x\leq xy;\quad\textup{(\ensuremath{1})}}

\CaseDisplay{xy\leq uxy+x,\quad xy\leq uxy+xz,\quad xy\leq uxyv+x,\quad xy\leq uxyv+xz;\quad\textup{(\ensuremath{2})}}

\CaseDisplay{xz\leq xyz;\quad\textup{(\ensuremath{3})}}

\CaseDisplay{xyz\leq xy+yz,\qquad xyz\leq xy+uyz,\qquad xyz\leq xy+yzv,\qquad xyz\leq xy+uyzv.\quad\textup{(\ensuremath{4})}}

There is one identity in (1), four in (2), one in (3), and four in (4). No minimality claim is made.

Theorem. These ten identities form a complete basis for \(A\). Word reversal gives a complete basis for its opposite.

\par\smallskip\noindent\textit{3. Directed-path criterion.}\quad 
For a polynomial \(P\), let \(U(P)\) be the union of its word contents, let \(H(P)\) be its word heads, and let \(E(P)\) consist of all ordered consecutive pairs in its words. Let \(E(P)^+\) be the positive transitive closure, so a loop is included only if there is a nonempty directed return walk.

Criterion. A word \(w\) satisfies \(A\models w\leq P\) if and only if

\CaseDisplay{h(w)\in H(P),\qquad c(w)\subseteq U(P),\qquad E(w)\subseteq E(P)^+.\quad\textup{(\ensuremath{6})}}

No word-length bound or normalization is imposed.

The head condition is precisely the head lemma for a left-zero-top extension: necessity is separated by setting a missing target head to \(t\) and all other variables to zero, and sufficiency reduces every evaluation with no right head equal to \(t\) to the core by replacing the remaining \(t\) values by zero. Missing content is separated in the core by assigning a target-only variable zero and every other variable \(e\).

For a core assignment whose right polynomial is not zero, write a variable value as \([L_x, R_x]\), with \(L_x\leq R_x\). Every right word is nonzero exactly when

\CaseDisplay{R_x\leq L_y\quad\text{for every }(x, y)\in E(P).}

Along a directed path these inequalities, together with \(L_y\leq R_y\), imply the corresponding boundary inequality for its endpoints. Thus if (6) holds, the target is nonzero. The interval hull of all right word values is also the hull of all variable values appearing in \(P\): along each nonzero word the first left endpoint is the minimum and the last right endpoint the maximum. Content inclusion therefore puts the target interval below that hull. This proves sufficiency.

For necessity of the path condition, suppose \((x, y)\in E(w)\) is not in \(E(P)^+\). Build a directed graph with vertices \(L_z, R_z\) for each variable. Include the intrinsic arcs \(L_z\to R_z\) and boundary arcs \(R_u\to L_v\) for \((u, v)\in E(P)\). A directed path from \(R_x\) to \(L_y\) exists exactly when there is a positive variable path from \(x\) to \(y\). Give value one to every endpoint vertex reachable from \(R_x\) (including itself), and zero to the remaining endpoints. Every arc is then weakly increasing. In particular each pair \((L_z, R_z)\) is a valid nonzero interval, and every right word is nonzero. But \(R_x=1\) and \(L_y=0\), so the target is zero and cannot lie below the nonzero right hull. This separates the failed pair, also when \(x=y\). The criterion follows.

\par\smallskip\noindent\textit{4. Soundness and unrestricted derivation.}\quad 
Each of (1)--(4) satisfies the criterion. Head and content conditions are immediate. The sole target edge in (3) is supplied by its two-step right path; every target edge of the other identities is present directly on the right. Hence all ten identities are sound.

Let \(w\leq P\) satisfy (6), and put \(h=h(w)\). Choose an original summand \(q_h\) beginning with \(h\). A linear target follows immediately from (1), or reflexivity when \(q_h=h\). Suppose the target is nonlinear. For every consecutive pair in \(w\), choose a positive directed path in the original graph \(E(P)\). Concatenate these paths, sharing each common endpoint once, to form a word \(W\) beginning with \(h\). Every adjacent pair of \(W\) is an actual edge of an original right word. The word \(w\) is obtained from \(W\) by deleting the intermediate vertices introduced along the chosen paths.

First derive that the initial two-letter prefix \(hy\) of \(W\) lies below \(P\). Pick an original right word containing the consecutive pair \(hy\). If it starts with that pair, deleting its suffix uses (1) with the whole word \(hy\) substituted for \(x\). If it has a nonempty prefix, its form is \(uhy\) or \(uhyv\). The guard \(q_h\) has form \(h\) or \(hz\), with nonempty substituted factors where present. The corresponding one of the four rules (2) gives \(hy\leq P\).

Now suppose a prefix \(Xy\) of \(W\) of length at least two is already below \(P\). Let \(z\) be the next letter. Choose an original right summand containing its actual edge \(yz\); it has one of the four forms \(yz\), \(uyz\), \(yzv\), \(uyzv\), with each present context nonempty. Substitute the existing nonempty prefix \(X\) for \(x\) in the appropriate rule (4). It gives

\CaseDisplay{Xyz\leq Xy+\text{that original summand}\leq P.}

Induction builds \(W\leq P\). This construction keeps the existing first letter at every step, so it requires no unjustified change of an internal edge's head. It applies to arbitrary path lengths, loops and repeated variables.

Rule (3), with a nonempty prefix and suffix substituted for \(x, z\), permits deletion of any single interior letter: the shortened word is below the longer word. Delete precisely the added internal vertices of the chosen paths, retaining the selected occurrences of all target letters. Therefore \(w\leq W\leq P\).

Every valid word absorption is thus derivable from the ten identities. Distributing both sides of any identity and deriving the absorption of each word in both directions proves completeness, without normalizing or bounding word length. The opposite statement follows by reversing words and derivations.

\subsection{A two-identity basis from positive directed reachability}\label{case:48}
Work in (+, multiplication), without constants, relative to the additively idempotent semiring axioms. Let A be catalogue entry 11700, with tables in numerical element order:

{\small\par\smallskip
\noindent addition       multiplication\par
\CaseLine{\ttfamily 0 0 0 0 4\qquad 4 0 4 4 4}
\CaseLine{\ttfamily 0 1 0 0 4\qquad 4 1 4 4 4}
\CaseLine{\ttfamily 0 0 2 2 4\qquad 3 3 3 3 3}
\CaseLine{\ttfamily 0 0 2 3 4\qquad 3 3 3 3 3}
\CaseLine{\ttfamily 4 4 4 4 4\qquad 4 4 4 4 4}
\smallskip}

The additive order is \(3<2<0<4\) and \(1<0<4\). The incomparable pairs are 1,\allowbreak{}2 and 1,\allowbreak{}3, each with join 0. Element labels are not constants of terms.

Theorem. The following two identities form a complete basis for A, relative to the additively idempotent semiring axioms:

\CaseDisplay{xyz\approx xy+xz,                         \quad\textup{(\ensuremath{1})}}

\CaseDisplay{xy\leq xz+zy.                             \quad\textup{(\ensuremath{2})}}

The multiplicative opposite has the reversed basis.

\par\smallskip\noindent\textit{Word values and soundness.}\quad 
Multiplication is increasing in either coordinate and preserves both incomparable joins, hence preserves every binary join. For associativity, an initial factor in \{2,\allowbreak{}3\} gives 3 under either association, and an initial factor 4 gives 4. For an initial factor u in \{0,\allowbreak{}1\}, the product (uv)r equals u exactly when \(v=r=1\), and otherwise equals 4. The product vr equals 1 exactly when \(v=r=1\); therefore u(vr) has precisely the same alternatives. This proves the semiring axioms.

Every nonlinear word has the following value, by induction:

\textbullet{} a head in \{2,\allowbreak{}3\} gives 3;

\textbullet{} a head 4 gives 4;

\textbullet{} a head in \{0,\allowbreak{}1\} gives the head when every later factor is 1, and gives 4 otherwise.

In particular a nonlinear value never equals 2. For (1), a head in \{2,\allowbreak{}3\} or equal to 4 gives the same constant on both sides. For a head in \{0,\allowbreak{}1\}, both sides equal that head if \(y=z=1\) and equal 4 otherwise. For (2), x in \{2,\allowbreak{}3\} gives left value 3 and \(xz=3\), so the inequality is immediate. For \(x=4\), \(xz=4\). For x in \{0,\allowbreak{}1\}, the case z different from 1 has \(xz=4\); the case \(z=1\) has zy equal to 1 if \(y=1\) and 4 otherwise, exactly as needed. Thus both identities are valid.

\par\smallskip\noindent\textit{A precise absorption criterion.}\quad 
Given a nonempty polynomial P, form a directed graph G(P) on its variables. Every nonlinear summand u contributes an edge from its first variable to each variable occurring after its first position. Linear summands contribute no edges. Loops are retained. Write x \(R_{P}\) y if there is a path of positive length from x to y; this relation is not made reflexive by convention.

Lemma. A linear variable y is absorbed by P if and only if y is a literal summand. A nonlinear word w with head h is absorbed if and only if

\CaseDisplay{h R_P x\quad\text{for every variable x occurring in
                              the tail of w}.               \quad\textup{(\ensuremath{3})}}

Proof. If a missing linear y is assigned 2 and all other variables 3, every right summand has value 3, whereas \(y=2\) is not below 3. The converse is literal inclusion.

For necessity of (3), let R be the set of vertices reachable by positive paths from h and suppose a tail variable x is not in R. Assign every variable in R the value 1. If h is not in R assign \(h=0\), and assign every remaining variable 2. These instructions are consistent. If an original nonlinear right word has head in R, all its tail variables also belong to R and its value is 1. If its head is h outside R, all its tail variables are in R and its value is 0. Every other nonlinear right word has head 2 and value 3. Every linear summand has value 0,\allowbreak{}1 or 2. Thus \(P\leq 0\). The target head is 0 or 1 and the chosen tail variable x has value 0 or 2, so \(w=4\), a contradiction. This includes \(x=h\): without a positive cycle, h need not belong to R.

Conversely suppose (3). The target is nonlinear, so its tail is nonempty and (3) implies that h is the head of at least one original nonlinear summand. If P has value 3 or 2, each original nonlinear summand has value 3. In particular h is in \{2,\allowbreak{}3\}, so the target is 3 and is absorbed. If \(P=1\), every right summand is 1. Every variable in any nonlinear right word is therefore 1, and the same holds for a linear right variable. All vertices on the paths in (3), including h, are then 1, so \(w=1\).

If \(P=0\), no right word has value 4. A right nonlinear head is therefore not 4. If h is in \{2,\allowbreak{}3\}, the target is 3 and is below 0. Otherwise h is 0 or 1. For any edge of G(P), assigning its source 0 or 1 forces its endpoint to be 1: otherwise the contributing right summand would be 4. Induction along a positive path from h therefore forces every reachable vertex to be 1. Condition (3) makes all target tail variables 1, so the target is \(h\leq 0\). The remaining case \(P=4\) is automatic. These cases exhaust A. \(\square\)

\par\smallskip\noindent\textit{Equational completeness.}\quad 
Repeated use of (1), with a nonempty suffix substituted for z, gives for every word \(x_{0}\) \(x_{1}\) ... \(x_{n}\) with \(n\geq 1\):

\CaseDisplay{x_0x_1\cdots x_n\approx \sum_{j=1}^{n}x_0x_j.              \quad\textup{(\ensuremath{4})}}

Duplicates are removed by additive idempotence. Thus every term is equationally a sum of literal linear variables and directed edges. Its graph is exactly G(P).

If \(x_{0}\to x_{1}\to \cdots \to x_{n}\) is a positive path in the original graph, derive \(x_0x_n\) below the sum of its edges by induction on n. The case \(n=1\) is literal inclusion. For \(n\geq 2\), an instance of (2) with \(x=x_{0}\), \(z=x_{n-1}\), \(y=x_{n}\) gives

\CaseDisplay{x_0x_n\leq x_0x_{n-1}+x_{n-1}x_n.}

The induction hypothesis bounds the first right term by the earlier edges. This also treats positive cycles and loops; a zero-length path is never used. Every substitution is a single variable and therefore nonempty.

For any valid nonlinear absorption \(w\leq P\), the lemma gives a positive path for every edge in the expansion (4) of w. The path derivation bounds each such edge by the expansion of P. Summing these inequalities and using (4) derives \(w\leq P\) from (1),(2). Valid linear absorptions are literal inclusions. Finally expand both sides of an arbitrary valid identity to polynomials, derive each summand's absorption into the other side, and combine the two inequalities. This proves completeness for arbitrary numbers of variables, word lengths and summands. Reversing all words transports the proof to the opposite. \(\square\)

Finite calculations are auxiliary verification. Further catalogue types inherit the result only after independent complete HSP checks in both directions. No claim about all subvarieties is made.

\subsection{Contained witnesses and positive paths between prefix variables}\label{case:49}
Let A be catalogue entry 12053 in the signature (+, multiplication), without constants. All bases are relative to the additively idempotent semiring axioms. Write \(w\leq P\) for w+P approximately equal to P. In numerical element order its tables are

{\small\par\smallskip
\noindent addition       multiplication\par
\CaseLine{\ttfamily 0 0 0 4 4\qquad 4 1 4 3 4}
\CaseLine{\ttfamily 0 1 0 3 4\qquad 3 1 3 3 3}
\CaseLine{\ttfamily 0 0 2 4 4\qquad 0 1 2 3 4}
\CaseLine{\ttfamily 4 3 4 3 4\qquad 3 1 3 3 3}
\CaseLine{\ttfamily 4 4 4 4 4\qquad 4 1 4 3 4}
\smallskip}

Element labels are not term constants. For any word u let c(u) be its content, i(u) the set of variables before its final position, and f(u) its final variable. Put \(i(y)=empty\) for linear words.

\par\smallskip\noindent\textit{Semiring structure and normal forms.}\quad 
The additive order has \(1<0<4\), \(2<0<4\), and \(1<3<4\). The incomparable pairs are 1,\allowbreak{}2 with join 0, and 0,\allowbreak{}3 and 2,\allowbreak{}3 with join 4. Every multiplication row and column is increasing and preserves these incomparable joins. It consequently preserves all binary joins.

Partition elements into \(E=\{2\}\), H=\{0,\allowbreak{}4\}, L=\{1,\allowbreak{}3\}, ordered by \(E<H<L\) for this paragraph only. The class of a product is the maximum of the two classes. Its left action is \(p_{E}=identity\), p\_H=(4,\allowbreak{}1,\allowbreak{}4,\allowbreak{}3,\allowbreak{}4), or p\_L=(3,\allowbreak{}1,\allowbreak{}3,\allowbreak{}3,\allowbreak{}3) according to the first factor's class. These maps are commuting idempotents with \(p_{H}\) composed with \(p_{L}\) equal to \(p_{L}\), and their compositions follow the same maximum rule. Therefore (uv)w and u(vw) agree. This proves the semiring axioms and the following nonlinear word values:

\textbullet{} If some prefix factor belongs to L, the value is 1 when the final factor is 1 and 3 otherwise.

\textbullet{} If no prefix factor belongs to L but some belongs to H, the value is 1 for final 1, 3 for final 3, and 4 otherwise.

\textbullet{} If every prefix factor is 2, the value is the final factor.

In particular the word identities

\CaseDisplay{xyz\approx yxz,\qquad x^2y\approx xy    \quad\textup{(\ensuremath{W})}}

hold. They permute all positions before the final one and remove their duplicates. A nonlinear normal word lists its nonempty prefix set in fixed order, followed by its final variable, which may also occur in the prefix. Every application uses a nonempty suffix. Linear words remain linear, and pure powers of length at least two normalize to \(y^{2}\). The inequality \(y\leq y^{2}\) is valid from the table.

\par\smallskip\noindent\textit{Exact criterion.}\quad 
For the original words of a polynomial P define a directed graph G(P): a nonlinear right word u contributes \(x\to f(u)\) for every x in i(u). Linear words contribute no edges. Paths in this graph must have positive length; no reflexive zero-length paths are inserted by convention.

Lemma. A linear y is absorbed exactly when an original word of P has content \{y\}. A nonlinear target w is absorbed exactly when both conditions hold:

1. Some original right word has content contained in c(w).

2. Every x in i(w) has a positive directed path to f(w) in G(P).

Proof. If there is no pure-y right word, assign \(y=2\) and every other variable 3. Each right word is 3, whereas \(y=2\) is not below

3. Conversely a pure-y right word is y or normalizes to \(y^{2}\), and \(y\leq y^{2}\) gives the inequality.

For failure of the first nonlinear condition, assign c(w) the value 2 and outsiders 3. The target is 2, and every right word has an outsider, hence has value 3. For failure of the path condition at x, put \(t=f(w)\) and let R be the set of vertices having a positive path to t. Assign R the value 2; if t is not in R assign \(t=0\); assign every other variable 1. The assignments are consistent. Each nonlinear right word with final variable in R has all prefix variables in R and value 2. One ending in t has all prefix variables in R and value 0 or 2. Every other nonlinear right word ends in a variable assigned 1 and has value 1. Right linear words are at most 0. Thus \(P\leq 0\). The target has final value 0 or 2 and a prefix occurrence outside R, which is assigned 1 or is t assigned 0. Its value is 3 or 4, neither below 0. This includes a missing positive cycle when \(x=t\).

For sufficiency it is enough to show that whenever P is below any proper principal bound, the target is below the same bound. The proper principal ideals are \{1\}, \{2\}, \{1,\allowbreak{}3\}, and \{0,\allowbreak{}1,\allowbreak{}2\}.

If \(P=1\), every right word is 1, and each nonlinear right word therefore has final value 1. The target prefix is nonempty, so its path condition includes an edge ending in \(t=f(w)\); hence t is a right final variable assigned 1. The target is 1.

If \(P=2\), every right word is 2. The word formula shows that every variable in each such word is 2, including linears. The path condition places every target prefix variable and t among those right variables, so the target is 2.

If \(P\leq 3\), every right word contains a variable assigned to L: this is equivalent to its value being in \{1,\allowbreak{}3\}, for linear and nonlinear words alike. The contained witness supplies such a variable in the target content, so the target is in \{1,\allowbreak{}3\}.

Finally suppose \(P\leq 0\). For a nonlinear right word, this means either its final value is 1, or all its prefix variables are 2 and its final value is 0 or 2. All right final variables, including linears, are therefore in \{0,\allowbreak{}1,\allowbreak{}2\}. The paths imply t is one of these finals. If \(t=1\) the target is 1. If t is 0 or 2, follow each path \(x\to \cdots \to t\) backwards. An edge whose final is not 1 forces its prefix endpoint to be 2; a new 2 is again not 1, so the implication propagates along the whole path. Thus every target prefix variable is 2, and the target equals \(t\leq 0\). The bound 4 is automatic. The four principal cases exhaust the algebra and prove the assertion. \(\square\)

\par\smallskip\noindent\textit{Moving a known final variable along a path.}\quad 
Fix a valid nonlinear absorption \(w\leq P\), put \(t=f(w)\), and retain an original contained witness b. A known word q is already below P. If an original nonlinear a has final y and its prefix contains f(q), define a normal word q' with final y and prefix

\CaseDisplay{i(q)\cup i(a).}

The local absorption \(q'\leq q+a\) is valid by the lemma: q is contained in q' because its final variable lies in i(a); each old prefix variable reaches f(q) via q when q is nonlinear, and then reaches y via a; each member of i(a) reaches y directly. If q is linear, its old prefix is empty and the same argument applies. The new prefix is nonempty because a is nonlinear.

Start \(q=b\). If f(b) differs from t, then f(b) lies in i(w), so the lemma supplies a positive path \(f(b)\to \cdots \to t\). Apply the local step along its edges using the original right word supplying each edge. After finitely many steps obtain a known q ending in t. If \(f(b)=t\), start with \(q=b\) without a path. The latter case may start with a literal t; no positive self-path is assumed.

\par\smallskip\noindent\textit{Growing the prefix backwards.}\quad 
For a known q ending in t put \(D=i(q)\). If an original nonlinear a has final variable g in D union \{t\}, then the normal q' with final t and prefix D union i(a) satisfies \(q'\leq q+a\). The right q is contained in the target. Every old prefix has its edge to t; every new prefix reaches g by a, and either \(g=t\) or q has an edge from g to t. The lemma proves validity. In the initial case \(q=t\), an a ending in t merely supplies its own prefix directly.

If some target prefix variable x is not yet in D, choose its positive path to t and follow it backwards until the first edge whose prefix endpoint is outside D but whose final is in D or equals t. The original right word for that edge permits a local growth step and adds at least that new variable. Continue until i(w) is contained in D. Every successful step strictly grows a finite set of right variables, so termination holds for arbitrarily many variables and for cycles. The final known guard q may contain variables outside the target; they are not discarded as if q alone were contained. Instead use the original b again:

\CaseDisplay{w\leq b+q.}

Here b supplies the contained witness, and the edges of q supply the paths from every target prefix variable to t. The lemma proves this final two-right-word step. All other steps also use at most two original or already derived right words.

\par\smallskip\noindent\textit{Reversible compression and finite basis.}\quad 
At most three words occur in one local absorption. Reserve their final variables individually, at most three. Group the other variables by nonempty prefix-incidence patterns in the three words, at most seven groups. Identify each group with a single symbol and normalize by (W). This substitution preserves validity and gives at most ten variables and normal word length at most eleven. For its inverse send each group symbol to the nonempty ordered product of its members, keeping the reserved final variables single. Identical prefix incidence, prefix permutations and duplicate removal recover every original normal word. Linear words and repeated finals are included; no empty product is substituted.

Theorem. The identities (W), together with all valid normal absorptions over a fixed ten-letter alphabet, with word length at most eleven and at most two right words, form a finite relative basis for A. Its opposite has the reversed basis.

Proof. This is a finite, exactly specified family, with membership decided by the lemma. It includes \(y\leq y^{2}\). Linear valid absorptions follow from their pure-variable original witness and (W). For a nonlinear absorption, path movement, backward growth, and the final retained-witness step derive it by local valid absorptions. Reversible compression makes every local step an instance of the fixed family modulo (W). Thus arbitrary-variable completeness holds. Distribute arbitrary terms, derive the absorption of every summand of either side into the other side of a valid identity, and combine both inequalities. Reversal transports all derivations to the opposite. \(\square\)

The finite family is not asserted small, minimal, or hereditarily finitely based. Finite checks supplement the unrestricted proof. Transfers to other catalogue types require complete independent HSP checks in both directions.

\subsection{A thirteen-identity tail-union basis for source 856}\label{case:50}
All identities are relative to the ai-semiring axioms in the constant-free signature. Element labels are not constants in terms. The addition and multiplication tables are

addition  multiplication 00044     00000 01144     00010 01244     22222 44434     00034 44444     00044

The additive order is \(2<1<0<4\), with 3 incomparable to 2,\allowbreak{}1,\allowbreak{}0 and below 4. Notice in particular that 2 is not an additive zero. Put K=\{3,\allowbreak{}4\}. Rows 0 and 2 of multiplication are constant 0 and 2. Row 1 has 1 only in column 3 and is otherwise 0. On K multiplication is maximum, and multiplying a K-element on the right by an element outside K gives 0. These facts give associativity by splitting a triple according to its first factor; they also verify preservation of each additive join by every left and right multiplication map. Thus the displayed algebra A is an ai-semiring.

\par\smallskip\noindent\textit{Word forms and basis.}\quad 
For a nonlinear word w, write h for its head and T for the set of variables in all positions after the first. The head can belong to T. Under an assignment:

\textbullet{} \(h=0\) gives value 0, and \(h=2\) gives value 2;

\textbullet{} \(h=1\) gives 1 if every variable of T is 3, and 0 otherwise;

\textbullet{} h in K gives 0 if some variable of T is outside K, and otherwise gives the maximum of all assigned values in the word.

Consider the following thirteen identities (inequalities abbreviate absorption):

1. xyz \ensuremath{\approx} xzy;

2. xy\ensuremath{{}^{2}} \ensuremath{\approx} xy;

3. x \ensuremath{\leq} \(x^{2}\);

4. xy \ensuremath{\leq} x\ensuremath{{}^{2}}y;

5. xyz \ensuremath{\leq} xy + xz;

6. xy \ensuremath{\leq} xyz + x;

7. xy \ensuremath{\leq} xyz + y;

8. xyz \ensuremath{\leq} xyzt + y;

9. x\ensuremath{{}^{2}}y \ensuremath{\leq} xy + y;

10. x\ensuremath{{}^{2}}yz \ensuremath{\leq} xyz + y;

11. x\ensuremath{{}^{2}}yz \ensuremath{\leq} xy + yz;

12. x\ensuremath{{}^{2}}yzt \ensuremath{\leq} xyz + yt;

13. x\ensuremath{{}^{2}}yz \ensuremath{\leq} x\ensuremath{{}^{2}}y + z.

Theorem. These thirteen identities are a basis for A. Reversing every word gives a basis for its opposite. No assertion of minimality, independence or hereditary finite basability is made.

The first two laws permute all tail positions and remove repetitions within the tail. They normalize each nonlinear word to its head followed by the sorted, nonempty set T. The repeated-head flag h in T is retained. They do not identify x\ensuremath{{}^{2}}y with xy.

\par\smallskip\noindent\textit{Exact criterion.}\quad 
Let P be a nonempty polynomial with linear variable set L and nonlinear words N. Write \(D=c(P)\), H for the set of heads of N, and

\(E_{h}\) = union of the tail contents of all h-headed words of N.

For a set C, J(C) means that L meets C or some u in N has content contained in C.

Lemma. A linear y is absorbed by P exactly when y belongs to L or N contains a word of singleton content \{y\}. A nonlinear word with head h, tail content T and content C is absorbed exactly when

(C1) C is contained in D;

(C2) J(C) holds;

(T) either T is contained in \(E_{h}\), or \(E_{h}\) meets H union L.

Necessity. If a linear witness is missing, assign \(y=3\) and every other variable 0. The target is 3, while every right summand is 0. The converse follows from identity 3 and normalization.

If C1 fails, assign a missing target variable 4 and every other variable 3. The target is 4 and P is 3. If C2 fails, assign C the value 3 and outsiders 0. The target is 3, while all right summands are 0. Both separate the absorption.

Suppose T fails. Then \(E_{h}\) avoids H union L and some target tail variable lies outside \(E_{h}\). If h is not in H, assign \(h=1\) and every other variable 2. Every nonlinear right word has head 2, all linear right words are at most 1, and the target is 0. If h is in H, assign \(h=1\), all of \(E_{h}\) the value 3, and all other variables 2. This is consistent because h is not in \(E_{h}\). The h-headed right words have value 1, all other nonlinear right words have value 2, and all linear terms have value 1 or 2 since \(E_{h}\) avoids L. The target has a tail variable not assigned 3 and is therefore 0. These prove necessity.

Sufficiency. The principal ideals of the addition semilattice, apart from the whole algebra, are \{2\}, \{2,\allowbreak{}1\}, \{2,\allowbreak{}1,\allowbreak{}0\}, and \{3\}. It suffices to prove that whenever P belongs to one of these ideals, so does the target.

\(P=3\) forces every variable in D to be 3; C1 then makes the target 3. If P belongs to \{2,\allowbreak{}1,\allowbreak{}0\} and the target were in K, all C-variables would be in K. A summand provided by J(C) would then lie in K, contradiction.

If P belongs to \{2,\allowbreak{}1\}, all variables in L union H belong to \{2,\allowbreak{}1\}. For any h-headed right word, its value is in \{2,\allowbreak{}1\} only when \(h=2\), or when \(h=1\) and every tail variable is 3. Condition T ensures \(E_{h}\) is nonempty (T itself is nonempty), hence h belongs to H. When \(h=2\) the target is 2. When \(h=1\) all \(E_{h}\)-variables must be 3, so \(E_{h}\) cannot meet H union L. Condition T therefore forces T to be contained in \(E_{h}\), making the target 1. Finally \(P=2\) forces all H-variables to be 2, and again the target is 2. This handles every proper principal ideal and proves sufficiency.

Soundness of the thirteen laws follows directly from the word formula or the criterion. For laws 4 and 9--13, the repeated head is supplied or forced by the tail/head intersection in T. Laws 6--8 use the displayed linear summand for J(C) and any additional tail variable. The two word laws preserve head and tail content, and identity 3 follows by the five square values. Finite table verification supplies an additional transcription check.

\par\smallskip\noindent\textit{Completeness.}\quad 
Fix a valid nonlinear absorption \(w\leq P\) and use the notation of the lemma. The proof divides according to whether \(E_{h}\) meets H union L.

\par\smallskip\noindent\textit{Case 1: no head-tail intersection.}\quad 
Here T is contained in \(E_{h}\). By repeated application of identity 5 with \(x=h\) and y, z nonempty tail words, all h-headed right words can be combined into a word \(r\leq P\) with head h and tail content exactly \(E_{h}\). There is at least one such right word because T is nonempty.

A witness for J(C) has a restricted form. A linear witness different from h would belong to T, hence to \(E_{h}\) intersect L, impossible. A nonlinear contained witness with head g different from h would put g in T subset \(E_{h}\) and in H, again impossible. A contained witness headed h cannot repeat h, because that would put h in \(E_{h}\) intersect H. Thus the witness is either the linear h, or a nonlinear h-headed word v with tail content contained in T.

If the witness v already has the same tail content as w, normalization gives \(w\leq P\). Otherwise write w approximately vq with q nonempty; when the witness is linear h, take q to be the whole target tail. If \(E_{h}=T\), r already normalizes to w. In the remaining situation let e list \(E_{h}\) minus T, a nonempty set. Identity 6, with \(x=v\), \(y=q\), \(z=e\), gives

{\small\par\smallskip
\CaseLine{\(w \approx  vq \leq  vqe + v \leq  P.\)}
\smallskip}

The first right word normalizes to r. This proves Case 1, without substituting an empty word.

\par\smallskip\noindent\textit{Case 2: a head-tail intersection exists.}\quad 
Choose g in \(E_{h}\) intersect (H union L) and an h-headed right word u whose tail contains g. By the word laws write u approximately hgA, where A may be omitted.

If g is linear in P, use identity 9 when A is omitted and identity 10 when A is nonempty. This gives a word h\ensuremath{{}^{2}}\(gA\leq P\) whose tail contains h.

Otherwise take a nonlinear right word \(v=gB\), with B nonempty. If A is omitted, identity 11 gives h\ensuremath{{}^{2}}gB \ensuremath{\leq} hg+gB \ensuremath{\leq} P. If A is nonempty, identity 12 gives h\ensuremath{{}^{2}}gAB \ensuremath{\leq} hgA+gB \ensuremath{\leq} P. The construction is valid when \(g=h\) as well. In all cases there is a derived word h\ensuremath{{}^{2}}\(U\leq P\) with U nonempty.

Now apply identity 13 repeatedly, substituting each original right summand for z. Every application extends h\ensuremath{{}^{2}}U by that summand while keeping it absorbed by P. By C1 the resulting derived word r has head h, contains all C, and repeats h. Extra variables are allowed at this stage.

Let w* be the repeated-head normal word with head h and content C. Identity 4 gives \(w\leq w*\). It remains to remove all variables of r outside C using a witness v for J(C). Such a witness is an original linear or nonlinear right summand with content contained in C. The tail of w* contains every variable of C, so the word laws permit writing w* approximately h v q, where q may be omitted: repeat tail letters as needed to provide the contiguous block v, and keep any uncovered letters in q. If r has no outsiders, it already normalizes to w*. Otherwise let e be a nonempty word listing all outsiders in r. If q is omitted, identity 7 gives \(hv\leq hve+v\). If q is nonempty, identity 8 gives \(hvq\leq hvqe+v\). The long right word normalizes to r in both cases. Thus \(w*\leq P\) and \(w\leq P\).

This completes all nonlinear absorptions. Linear absorptions follow from a literal summand or identity 3 and singleton-tail normalization. Expanding arbitrary terms into polynomials and applying the construction to every summand in both directions derives every valid identity. All arguments allow arbitrary variable numbers, word lengths and polynomial sizes. Reversing words and derivations proves the opposite result.

\subsection{A finite basis from a cyclic head or covered tails}\label{case:51}
Source 866 in the five-element catalogue. The signature is \((+, multiplication)\) without constants, relative to the additively idempotent semiring axioms. Element labels are not term constants.

\par\smallskip\noindent\textit{1. Word values and normalization.}\quad 
The tables, in order 0,\allowbreak{}1,\allowbreak{}2,\allowbreak{}3,\allowbreak{}4, are

{\small\par\smallskip
\noindent addition       multiplication\par
\CaseLine{\ttfamily 0 0 0 0 4\qquad 0 0 0 0 0}
\CaseLine{\ttfamily 0 1 1 1 4\qquad 0 0 0 1 0}
\CaseLine{\ttfamily 0 1 2 2 4\qquad 2 2 2 2 2}
\CaseLine{\ttfamily 0 1 2 3 4\qquad 2 2 2 3 2}
\CaseLine{\ttfamily 4 4 4 4 4\qquad 0 0 0 0 4}
\smallskip}

The additive order is the chain \(3<2<1<0<4\). Multiplication is increasing in both coordinates of this chain and hence distributes over its maximum. Associativity follows by checking the first factor: 0,\allowbreak{}2 are left zeros; a triple starting in 1 has value 1 exactly when both later factors are 3, otherwise 0; one starting in 3 is 3 exactly under the same condition, otherwise 2; one starting in 4 is 4 exactly when both later factors are 4, otherwise 0. Both associations have these descriptions.

For any nonlinear word, its value is determined by its first factor h:

\textbullet{} \(h=0\) gives 0, and \(h=2\) gives 2;

\textbullet{} \(h=1\) gives 1 if all tail factors are 3, and 0 otherwise;

\textbullet{} \(h=3\) gives 3 if all tail factors are 3, and 2 otherwise;

\textbullet{} \(h=4\) gives 4 if all tail factors are 4, and 0 otherwise.

This follows by induction. The laws

\CaseDisplay{xyz\approx xzy,\qquad xy^2\approx xy       \quad\textup{(\ensuremath{W})}}

therefore hold. They normalize a word to its literal head followed by its distinct tail variables in fixed order. The head may also occur in that tail. An empty tail is a linear word, not an empty substitution.

\par\smallskip\noindent\textit{2. Exact absorption criterion.}\quad 
For a polynomial P set \(D=\cup\) of all right contents and \(Q=set\) of all right heads, including linear words. For a variable h put

\CaseDisplay{E_h=\bigcup_{u\in P,\ h(u)=h}r(u),}

where r(u) is the content after the first position; linear words contribute nothing to \(E_{h}\).

Lemma. For any nonempty target w, linear or nonlinear, \(w\leq P\) is valid exactly when the following hold, with \(h=h(w)\), \(C=c(w)\), \(T=r(w)\):

1. C is contained in D;

2. h belongs to Q;

3. either h belongs to \(E_{h}\), or T is contained in \(E_{h}\);

4. some original right word v has content contained in C.

Necessity. A content variable missing from D can be assigned 2, with all others 3. The target is 2 and \(P=3\). If h is missing from Q, assign \(h=0\) and all others 3. The target is 0 while every right word is 2 or 3.

If condition 3 fails, choose x in T outside \(E_{h}\), with h outside \(E_{h}\). For x different from h assign \(h=1\), \(x=2\) and all other variables 3. Every right word headed by h is 1, since neither h nor x is in its tail; all other right words have value 2 or 3. Thus \(P\leq 1\) while \(w=0\). If \(x=h\), instead assign \(h=1\) and all others 3, with the same conclusion.

If condition 4 fails, assign C the value 4 and all outside variables

3. Then \(w=4\). Each right word contains an outside variable, and the value formula makes it 0,\allowbreak{}2, or 3, so \(P\leq 0\).

Sufficiency. Consider the assigned head value. If it is 0 or 2, a right word with the same head gives that value. If it is 1 and the target is 1, a same-head word is 1 or 0 and suffices. If instead the target is 0, either a same-head word repeats h and is 0, or condition 3 supplies a same-head word with a target tail variable valued outside \{3\}; that word is 0.

If the head is 3, a target value 3 is least in the chain. For a target value 2, some target tail variable is not 3. It occurs in a right word by condition 1. A word with a non-3 factor cannot have value 3, so it is at least 2. Finally suppose the head is 4. If the target is 4, all its variables are 4, and the contained witness of condition 4 has value 4. Otherwise the target is 0; a same-head word is 0 or 4 and suffices. Linear targets are included in these cases.

\par\smallskip\noindent\textit{3. Derivation with at most three right words.}\quad 
Normalize all words with (W). A linear target h has a contained original word v with content \{h\}; the one-right absorption \(h\leq v\) is valid by the lemma. For a nonlinear target h (product of T), choose a contained original witness v.

Cyclic-head case h in \(E_{h}\). Fix an original word a headed by h whose tail also contains h. Put \(T_{0}=c(v)\) minus \{h\}; it is contained in T. The word \(s_{0}=h\) (product of \(T_{0}\)), which may be linear, satisfies

\CaseDisplay{s_0\leq v+a.                  \quad\textup{(\ensuremath{1})}}

Indeed its content is covered, a supplies the head and the cyclic-head alternative, and v is contained. For each missing x in T choose an original \(r_{x}\) containing x. Add x to the current tail and normalize; the resulting s' satisfies

\CaseDisplay{s'\leq s+a+r_x.               \quad\textup{(\ensuremath{2})}}

The current s is a contained witness, a retains the cyclic head, and \(r_{x}\) supplies the new content. This also permits \(x=h\). Thus at most |T| such steps end at the target.

Covered-tail case h outside \(E_{h}\). Condition 3 supplies an original \(r_{x}\) headed by h with x in its tail, for each x in T. There is at least one such word because a nonlinear target has a nonempty tail. Merge these words one at a time: two words headed by h with tail sets U, V absorb the normal word headed by h with tail U union V. This is a valid two-right-word step by the lemma: the tail is covered by same-head right words and either right content is contained in the merged content. The process produces \(s=h\) (product of S)\ensuremath{\leq}P, with S containing T. Although S may include variables outside C, the final absorption

\CaseDisplay{w\leq s+v                   \quad\textup{(\ensuremath{3})}}

is valid: s supplies all content and the same-head tail cover, and the original v remains a contained witness. No assumption that an intermediate word has content inside C is used.

Every right word in these steps is original or has already been derived below P, and each step has at most three right words.

\par\smallskip\noindent\textit{4. Reversible compression and a finite basis.}\quad 
A local step has at most four words in total. Reserve each variable that is a head of one of them, at most four variables. Group the remaining variables, all of which occur only in tails, by their nonempty incidence pattern among those at most four tails. There are at most 15 such groups.

Identifying each group to one variable and normalizing by (W) preserves validity, since it is an ordinary substitution in a valid identity. It gives at most 19 variables and normal words of length at most 20. The identification is reversible: replace each group variable by the nonempty product of its original members, and each reserved head variable individually. The common incidence patterns restore the original tail sets; (W) recovers all the original normal words.

Fix a 19-letter alphabet and let Sigma consist of (W) and all valid normal absorptions with at most three right words, using this alphabet and words of length at most 20. The lemma decides membership in this finite, explicitly specified family. A small enumeration is not claimed.

Theorem. Sigma is a finite identity basis for source 866; its word reversals form a finite basis for the multiplicative opposite.

Proof. Every valid local step above is a substitution instance, modulo (W), of a member of Sigma. Section 3 therefore derives every valid word absorption into a polynomial, with no bound on original variables or word lengths. Distributing terms into word sums and deriving each summand absorption proves every valid identity. Soundness is built into the definition and (W). Word reversal transports the basis and derivations to the opposite.

Independent finite computations check the value formula, criterion, construction, compression, and transfers to other catalogue entries. The unrestricted proof is Sections 1--4. Transferred classifications require equality of generated varieties in both HSP directions.

\subsection{A nonfinite-basis theorem for odd cycles with an allowed literal}\label{case:52}
Work in the semiring signature without constants. Write \(u\leq P\) for \(u+P=P\). This theorem is new:the support barrier includes x as well as \(x^{3}\). The earlier 207 theorem excluded x and cannot be used unchanged because \(x\leq x^{3}\) holds in the present three sources.

\par\smallskip\noindent\textit{Structural hypotheses and critical identities.}\quad 
Suppose an additively idempotent semiring A contains an element g with

{\small\par\smallskip
\noindent \(g\leq g^{3}\), \(g^{2}\) \(\nleq g^{3}\), and \(g^{r}\) \(\nleq g^{3}\) for every \(r\geq 4\).\par
\smallskip}

Also suppose there are e, z with \(e^{2}=e\), \(z^{2}=z\), \(ez=ze=z\), and e \(\nleq z\). They need not be identities or global zeroes. For finite A the power tail condition is checked exactly by continuing until the power sequence repeats, then inspecting its entire eventual cycle.

For odd \(n\geq 3\), using distinct variables x, v0,..., v(n-1), define

P\_n=x\ensuremath{^{3}}+sum\_i(v\_i \(v_{i+1}\)+\(v_{i+1}\) \(v_{i}\)), indices modulo n, \(w_{n}=x\) v0 v1 ... v(n-1).

If A satisfies \(w_{n}\leq P_{n}\) for all these n, then A is nonfinitely based. Both orientations of every edge are included, even for commutative A. The proof below treats arbitrary nonempty polynomial substitutions, repeated source variables, and both multiplicative contexts.

\par\smallskip\noindent\textit{The augmented support barrier.}\quad 
Fix odd \(n\geq 5\) with \(n>2k^{2}\). Let F consist of the literal x, the cube \(x^{3}\), and both orientations of the edges of the simple n-cycle. Suppose q is a nonempty polynomial on at most k variables and every expanded word of a nonempty polynomial substitution phi(q) lies in F. We show that for ANY absorption \(u\leq q\) valid in A, every word of phi(u) also lies in F.

First, each variable with image exactly \{x\} is type X, whether or not it occurs literally in q. Every other variable occurring literally in q is type L, with its image contained in F. Such an L variable cannot also occur nonlinearly in q:its image includes an edge or \(x^{3}\), which with a nonempty other factor gives a forbidden cubic involving cycle letters, or length at least 4. This observation explicitly allows the mixed literal images \{x, \(x^{3}\)\},\{x, edge\}, and the whole polynomial F.

All remaining variables occurring in q are of type Y, with image\{\(x^{2}\)\}, or type V, with image a nonempty set of single cycle letters. Indeed a nonlinear source word has degree 2 or 3. A degree 3 word must expand to \(x^{3}\), so all factors are X. A degree 2 word expands to \(x^{3}\) via X/Y, or to cycle edges via two V variables. Mixed choices in any image produce \(x^{2}\), a mixed x/cycle word, a wrong cubic, or length at least 4, and are excluded. These are all possibilities. In particular, the nonlinear source words are XXX, XY/YX, or VV; literal X words are now also allowed.

Each original VV pair contributes a complete bipartite subgraph of the n-cycle. Its two image sets are disjoint, since F has no loops. A cycle with \(n\geq 5\) has no four-cycle and degree at most 2, so each such bipartite subgraph has at most 2 edges. There are at most \(k^{2}\) ordered source pairs. Their union H has at most \(2k^{2}<n\) edges, is a proper subgraph of the cycle, and is therefore a forest. Every V image set is monochromatic in any proper two-coloring of H:its vertices share a neighbor through a source pair. Color the forest by g, \(g^{2}\), allowing independent component flips.

Assign X the formal value g, Y the value \(g^{2}\), and L the value \(g^{3}\). Assign each V its image set's color. Variables absent from q receive \(g^{4}\). These are FREE source assignments; an L image containing the whole odd cycle is not evaluated to define its formal value. Every q word has value g or \(g^{3}\), so \(q\leq g^{3}\). Validity forces \(u\leq g^{3}\) for every forest coloring. All assigned values are positive powers of g. Therefore:

\textbullet{} A literal u may only be L or X, and its image is in F. A V can be flipped to \(g^{2}\); Y and an absent variable also fail the cubic bound.

\textbullet{} A cubic u must consist entirely of X. Any L, Y, absent variable, or a V flipped to \(g^{2}\) makes its total exponent at least 4. Its image is then \{\(x^{3}\)\}.

\textbullet{} Every word of degree at least 4 is excluded by the power-tail condition.

\textbullet{} A quadratic cannot contain L or an absent variable. XX gives \(g^{2}\); YY gives exponent 4. A V with X or Y can be flipped to fail. XY/YX has image\{\(x^{3}\)\}. Only a pair of V variables remains to consider.

Every valid absorption \(u\leq q\) in A has an original q word whose content is contained in c(u). Otherwise assign variables of u the value e and all other variables z. Then \(u=e\) and every original q \(word=z\), contrary to e \(\nleq z\). This necessary condition uses only the displayed two-element multiplicative subsemilattice; global absorption by z is unnecessary. For a quadratic in two V variables a, b, the contained original cannot be literal (V letters are not in F), cubic (all q cubics use X), or a square (it would introduce a cycle loop). Hence q already contains ab or ba. Their whole expanded support consists of cycle edges. This also excludes a V square. The support barrier is proved.

\par\smallskip\noindent\textit{From the barrier to nonfinite basability.}\quad 
Suppose A had a finite identity basis, with variable counts bounded by k, including the semiring axioms. Choose odd \(n>2(k+2)^{2}\), \(n\geq 5\). Expand a purported derivation of \(P_{n}=P_{n}+w_{n}\) into sums of nonempty associative words, and consider the FIRST step introducing any word outside F. All earlier supports, including any newly introduced literal x, are contained in F. For one expanded component of the identity application, every word u of the replacement side satisfies \(u\leq q\) for its old side q. Absorb nonempty left and right monomial contexts into at most two fresh source variables. Omit an absent context; never substitute an empty word. The resulting contextual absorption is valid in A, has at most k+2 variables, and its old expanded support lies in F. The barrier forbids the new word. Semiring normalization itself changes no expanded support. But \(w_{n}\) has length \(n+1\geq 4\) and is outside F. This contradiction proves NFB.

More generally, if A satisfies the structural hypotheses, T has A in Var(T), and T satisfies all the critical identities, then T is NFB by the same argument:every identity of T is valid in A. The critical identity hypothesis must be checked in T separately. Merely containing an NFB subvariety is insufficient.

\par\smallskip\noindent\textit{The three catalogue sources.}\quad 
894 has \(g=3\), \(g^{2}=1\), \(g^{3}=2\), \(g^{r}=0\) for \(r\geq 4\). Its cubic principal ideal contains g and \(g^{3}\) but neither \(g^{2}\) nor 0. Take \(e=4\), \(z=0\). The sum \(P_{n}\) can fail to be 0 or 4 only if \(x=3\) and every edge product is 2. Such edge products force alternating values 3 and 1 around the cycle, impossible for odd n. If \(w_{n}=4\) all variables are 4 and \(P_{n}=4\); otherwise \(w_{n}=0\). Thus \(w_{n}\leq P_{n}\) holds for every odd n.

13236 and 13237 have \(g=1\), \(g^{2}=3\), \(g^{3}=0\), \(g^{r}=4\) for \(r\geq 4\). The cubic ideal contains g and \(g^{3}\), excludes \(g^{2}\) and 4, and \(e=4\), \(z=2\) satisfies the separation conditions. If no variable is 2, \(w_{n}=4\). Either \(x^{3}=4\), or \(x=1\) and keeping \(P_{n}\) below 4 would require every edge product 0. That again forces the impossible alternation1,\allowbreak{}3 on an odd cycle. If some variable is 2, \(w_{n}=2\) and an original cube or incident edge has value 2, so \(P_{n}\) absorbs 2 in both addition tables. These are full all-length proofs. Alternating generator/square values on even cycles give explicit counterexamples in all three sources.

14839 has no suitable two-element idempotent separation here and is NOT classified by this argument. Its support may gain new chords; the contained-content step must not be silently omitted.

\par\smallskip\noindent\textit{Exact all-length certification for finite target semirings.}\quad 
For each a, b with \(a^{3}\leq b\), allow a transition \(r\to s\) exactly when both \(rs\leq b\) and \(sr\leq b\). Fix a first vertex s0. Read a path of at least 3 vertices and keep the finite state(last vertex, ordered product a*s0*...*last, vertex-count parity). Seed all allowed three-vertex paths and close under appending one allowed vertex. A reachable state with odd parity, a closing edge to s0, and product \(\nleq b\) is exactly a critical-identity counterexample. There are at most 2|T|\ensuremath{^{2}} states per root a, b, s0; no truncation on cycle length is used.

A counterassignment yields such a path with b equal to its polynomial value. Conversely a bad state gives an odd closed walk of length at least 3; its cube and all oriented edge terms are \(\leq b\) but its ordered product is not. This proves the finite-state criterion for ALL lengths, including noncommutative targets. Intrinsic structural witnesses give direct NFB proofs; transferred cases additionally require complete HSP.

\subsection{Linear pins in a monogenic nilpotent core}\label{case:53}
All bases are relative to the additively idempotent semiring axioms in the constant-free signature.

\par\smallskip\noindent\textit{1. Algebra and finite basis.}\quad 
Let \(F=\{0, a, a^2, a^3\}\) have truncated monogenic multiplication, \(a^i a^j=a^{i+j}\) for \(i+j<4\) and zero otherwise. The additive order has \(a<a^3<0\) and \(a^2<0\), with \(a^2\) incomparable with \(a, a^3\). Add a greatest element \(t\), with \(ts=t\) for all elements and \(st=0\) for \(s\in F\). Write \(A\) for the resulting semiring.

The multiplication is associative by truncation and the first-occurrence rule for \(t\). For distributivity in the core, the only nontrivial comparable pair is \(a<a^3\), whose products under multiplication by \(a\) are \(a^2\) and zero, and by \(a^2\) are \(a^3\) and zero. They preserve order. For either incomparable pair involving \(a^2\), multiplication by any nonzero power gives either two incomparable values or one value zero, so the join of the products is zero, as is the product of the join. The extension by \(t\) preserves distributivity as in the source-209 proof.

The head lemma and \(xyzu\approx x^4\) from that proof remain valid: a valid absorption in \(A\) is exactly an absorption in \(F\) whose target head occurs as a right head. Every word normalizes to a word of length at most three or a fourth power.

Let \(W\) be the set of all words of lengths one through four on twelve fixed variables. Define \(\Sigma\) to be all valid absorptions \(w\leq\sum P\) in \(A\), where \(w\in W\) and \(\varnothing\ne P\subseteq W\).

Theorem. The finite, effectively specified set \(\Sigma\) is a complete basis of \(A\), and its word reversal is a complete basis of the opposite. No enumeration of this large finite set is claimed.

\par\smallskip\noindent\textit{2. Exact criterion.}\quad 
For a normalized \(P\), let \(L\) be its linear variables, \(E\) the undirected graph of its quadratic words, \(V\) its incident vertices, and \(T\) the variables of its cubic words. Define the initially marked vertices by

\CaseDisplay{M=T\cup(L\cap V).}

Unlike in the flat-core case, a linear variable incident with an edge is a mark, not an immediate contradiction.

Call \(P\) inconsistent if it contains a fourth power, or if \(M\ne\varnothing\) and either \(E\) is nonbipartite or some bipartite component has marked vertices in both parts. In a consistent graph with \(M\ne\varnothing\), let \(U\) contain \(M\) and all vertices in the marked part of every marked component; let \(D\) contain the opposite parts of those components.

An absorption \(w\leq P\) holds precisely when its head lies in \(H(P)\) and one of the following holds:

\textbullet{} \(P\) is inconsistent;

\textbullet{} \(M=\varnothing\), \(w\) is linear, and \(w\in L\);

\textbullet{} \(M=\varnothing\), \(w=xy\), both vertices lie in \(V\), and either \(E\) has an odd cycle or the vertices are in opposite parts of a common component;

\textbullet{} \(M\ne\varnothing\), \(w\) is linear and its variable lies in \(L\cup U\);

\textbullet{} \(M\ne\varnothing\), \(w\) is cubic and all its variables lie in \(U\);

\textbullet{} \(M\ne\varnothing\), \(w=xy\), and either its vertices lie in \(U, D\), in either order, or they lie in opposite parts of a common component.

To prove this criterion, it suffices by the head lemma to study core evaluations with \(P\ne0\). If an edge is present and \(P=a^2\), every edge endpoint must have exponent one and every linear variable must have exponent two; hence \(L\cap V\) must be empty. A cubic word cannot occur. If an edge or cubic is present and \(P=a^3\), every edge has exponent sum three, every cubic variable has exponent one, and every linear variable has exponent one or three. An incident linear variable must therefore have exponent one. Thus exactly \(M\) supplies exponent-one marks. Quadratic components are bipartite with complementary exponents one and two, and marks fix the chosen side. An odd component is incompatible with this second kind of evaluation.

When \(M\) is empty and edges occur, the common value \(a^2\) is always possible. A common value \(a^3\) is additionally possible exactly when the graph is bipartite. When \(M\) is nonempty, only common value \(a^3\) is possible, with \(U\) forced to exponent one and \(D\) forced to exponent two. Unmarked components choose their two assignments independently. Linear variables outside the graph and cubics independently take exponent one or three. These describe every nonzero core evaluation. If there are only linear words, assigning all of them \(a\) separates every nonlinear target. Variables absent from \(P\) can be assigned zero. Inspecting a target under these evaluations yields exactly the cases listed above, and supplies separating assignments in every excluded case.

\par\smallskip\noindent\textit{3. Finite local deductions.}\quad 
Use all the local rules from the source-209 proof except the direct linear-conflict family \(h^4\leq x+v+q_h\). Those omitted rules are in general invalid here. Retain the following named families: zero transfer, square-cubic conflict, orientation, odd-path compression, consequences of a square, initial cubic markers, two-edge marker propagation, incompatible markers, cubic targets, mixed quadratic targets, word normalization and head domination. Their proofs remain valid because each nonzero quadratic right summand must have the same exponent as the other quadratic summands, and the presence of a cubic forces exponent three. A right guard can lie below that value but cannot change these exponent constraints.

Add these two finite families, where \([xy]\) is either orientation and \(q_x\) is an original head guard of length at most four:

\CaseDisplay{x^3\leq x+[xy],\qquad x\leq rx^2+q_x.}

The first holds on the core because if its right side is not zero, the edge and linear occurrence force \(x=a\), \(y=a^2\), and their join is \(a^3=x^3\). It already has head \(x\) on the right. In the second, a nonzero cubic forces \(r=x=a\), and \(a\leq a^3\); the guard provides the required target head. The head lemma proves both in \(A\). Every rule involves at most eight variables, all summands of length at most four, and therefore belongs to \(\Sigma\) after renaming.

If \(M=\varnothing\), the odd-path and square rules derive exactly the quadratic cases as in the flat-core proof; the linear case is immediate.

If \(M\ne\varnothing\), obtain at least one cubic word below \(P\): use an existing cubic if \(T\) is nonempty, and otherwise choose \(x\in L\cap V\) and an incident edge to obtain \(x^3\leq P\) by the first new rule. Fix one such cubic \(ruv\). Every original cubic variable obtains a marker \(rx^2\) from the retained initial-marker rule. For each \(x\in L\cap V\), first derive \(x^3\), then apply that same initial-marker rule to \(ruv\) and \(x^3\), yielding \(rx^2\). Thus every initial mark has a common-head marker.

The retained two-edge rule propagates these markers along every even walk. An odd component supplies a square through odd-path compression and conflicts with the cubic; marks on opposite parts propagate to both endpoints of an edge and use the incompatible-marker rule. Thus each inconsistent case derives \(h^4\leq P\) for every original head \(h\). Head domination then gives every permitted target.

In the consistent marked case, every member of \(U\) has a marker. The retained cubic-target and mixed-target rules derive the stated cubic and quadratic targets; the odd-path rule derives the remaining quadratic cases. The second new rule derives every required linear target from its marker and head guard. Original linear targets are already summands. Consequently every case of the criterion is derivable, on arbitrary variable sets and graphs of arbitrary size. Distribution and two-way summand absorptions prove the theorem.

\subsection{Finite bases determined by two leading variables}\label{case:54}
Sources 2024 and 2025 in the five-element catalogue. The signature is (+, multiplication) without constants. Bases are relative to the additively idempotent semiring axioms. Numerical labels name elements.

\par\smallskip\noindent\textit{1. Tables, word values, and normal forms.}\quad 
Both sources have multiplication

{\small\par\smallskip
\CaseLine{\ttfamily 0 0 0 0 0}
\CaseLine{\ttfamily 1 1 1 1 1}
\CaseLine{\ttfamily 2 2 2 2 2}
\CaseLine{\ttfamily 0 0 1 4 4}
\CaseLine{\ttfamily 0 0 0 4 4}
\smallskip}

Their addition tables, in order 0,\allowbreak{}1,\allowbreak{}2,\allowbreak{}3,\allowbreak{}4, are respectively

{\small\par\smallskip
\CaseLine{\ttfamily 2024\qquad 2025}
\CaseLine{\ttfamily 0 0 0 4 4\qquad 0 0 0 4 4}
\CaseLine{\ttfamily 0 1 0 3 4\qquad 0 1 1 3 4}
\CaseLine{\ttfamily 0 0 2 4 4\qquad 0 1 2 3 4}
\CaseLine{\ttfamily 4 3 4 3 4\qquad 4 3 3 3 4}
\CaseLine{\ttfamily 4 4 4 4 4\qquad 4 4 4 4 4}
\smallskip}

All semiring axioms can be verified directly on the displayed finite tables, independently of the identities to be proved. The elements 0,\allowbreak{}1,\allowbreak{}2 are left zeros for multiplication. For a nonlinear word, let H and K be the values of its first and second variables. Its value is:

\textbullet{} H if H belongs to \{0,\allowbreak{}1,\allowbreak{}2\};

\textbullet{} 0 if H belongs to \{3,\allowbreak{}4\} and K belongs to \{0,\allowbreak{}1\};

\textbullet{} 1 if \(H=3\) and \(K=2\), or 0 if \(H=4\) and \(K=2\);

\textbullet{} 4 if H, K belong to \{3,\allowbreak{}4\} and all later factors belong to \{3,\allowbreak{}4\}; otherwise 0 in this final case.

This follows directly by induction from the multiplication table. It proves the three word laws

\CaseDisplay{xyzt\approx xytz,\qquad xy^2\approx xy,\qquad xyx\approx xy. \quad\textup{(\ensuremath{W})}}

Fix a variable order. Every nonlinear word normalizes to its two literal leading variables h, k, followed by its distinct remaining content in fixed order, excluding h and k. To see this, interchange adjacent positions after the second with the first law. Repeated later variables become adjacent and are removed with the second law, using a nonempty prefix. A later copy of h or k may be moved into the third position and removed by the third or second law. Repeat until the stated form is obtained. Linear words stay linear; the leading variables of a nonlinear word may coincide. No empty substitutions are used.

\par\smallskip\noindent\textit{2. Exact absorption criteria.}\quad 
For a polynomial P let Q be its set of heads. For a variable h let \(E_{h}\) be the set of second variables of nonlinear original right words whose head is h. Write c(u) for word content.

Lemma. A linear variable y satisfies \(y\leq P\) precisely when some original right word has content \{y\}. For nonlinear w with first variables h, k, \(w\leq P\) holds in source 2024 precisely when

1. h belongs to Q;

2. some original v in P has c(v) contained in c(w);

3. either k belongs to \(E_{h}\) or \(E_{h}\) meets Q.

For source 2025 replace condition 3 by

3'. either k belongs to \(E_{h}\) or h belongs to \(E_{h}\).

Necessity. If h is missing from Q, assign \(h=1\) and every other variable

2. Then \(w=1\) and \(P=2\), and 1 is not below 2 in either addition order. If no contained witness exists, assign c(w) the value 4 and all outside variables 0. The target is 4 and every original right word is 0. The same assignment proves necessity for a linear target.

For 2024, assume h in Q but \(E_{h}\) disjoint from Q and k outside \(E_{h}\). Assign \(h=3\), every variable of \(E_{h}\) the value 2, and every remaining variable 1. Right words headed by h are 1 when nonlinear and 3 when linear; all other right words have head value 1 and value 1. Thus \(P\leq 3\). If k differs from h the target is 0; if \(k=h\) it is 0 or 4. Neither is below 3.

For 2025, failure of 3' means h, k outside \(E_{h}\). Assign \(h=3\), \(k=1\) when k differs from h, and every remaining variable 2. Right words headed by h are 1 or 3, and other right heads have value 1 or 2. Therefore \(P\leq 3\). The target is again 0, or 0/\allowbreak{}4 if \(k=h\), and is not below 3.

Sufficiency. For a linear target, its pure-variable witness has value y or \(y^{2}\), and \(x\leq x^{2}\) holds in both algebras; powers stabilize at \(x^{2}\). Now suppose w is nonlinear. If its assigned head is 0,\allowbreak{}1, or 2, a same-head right word has exactly that value. Suppose its head is 4. A target value 0 is bounded by any same-head right word, valued 0 or 4. For target value 4, all target variables lie in \{3,\allowbreak{}4\}; the contained witness is 4 except possibly a literal of value 3. In that exception its join with a same-head right word valued 0 or 4 is 4.

It remains to consider head value 3. If the target is 1, every same-head right word has value 1,\allowbreak{}0,\allowbreak{}3, or 4, all above 1. If the target is 0, its second value is different from 2. A right word beginning hk, or beginning hh, consequently has value 0 or 4. In the extra crossing case for 2024, take a beginning hg and b headed by g. When g is not valued 2, a is 0 or 4. When \(g=2\), \(a=1\), \(b=2\) and their join is 0 in 2024.

Finally, if the head is 3 and target value 4, all target variables lie in \{3,\allowbreak{}4\}. The contained witness is 4 except possibly a literal of value 3. In that exception, a right word beginning hk or hh has value 0 or 4 and its join with the literal is 4. For a 2024 crossing with \(g=2\), the word headed by g has value 2, whose join with 3 is 4. These cases prove sufficiency for both algebras.

\par\smallskip\noindent\textit{3. A bounded number of original witnesses.}\quad 
For a linear target choose its pure-variable original witness; it alone absorbs the target. For a nonlinear target choose a contained original v. In source 2024, if k in \(E_{h}\), choose a right word a beginning hk. The sum v+a satisfies every condition. Otherwise choose g in \(E_{h}\) intersect Q, a right word a beginning hg, and b headed by g. The sum v+a+b retains the head, the crossing, and the contained witness. It therefore absorbs w. At most three original right words are needed.

In source 2025 choose, in addition to v, a right word beginning hk or hh according to condition 3'. These at most two original right words satisfy the criterion and absorb w. No intermediate construction or bound on the original number of variables is involved.

\par\smallskip\noindent\textit{4. Reversible compression and completeness.}\quad 
For 2024 a local absorption has at most four words altogether. Reserve individually every variable in a first or second position of any of these words, including linear words. At most eight variables are reserved. Group all other variables by their nonempty incidence patterns in the portions after the first two positions. There are at most 15 groups. Identify each group to a fresh variable and normalize by (W). This preserves validity, since it is an ordinary substitution in a valid identity. At most 23 variables and word lengths at most 25 result.

Conversely substitute for each group variable the nonempty product of its members. Its common incidence pattern restores exactly the original later content in every word. Reserved variables restore the leading positions individually, so none of those positions expands into a product. Laws (W) recover all original normal words.

For 2025 there are at most three words. The same construction reserves at most six leading variables and uses at most seven nonempty patterns, giving 13 variables and word lengths at most 15.

For source 2024, define \(Sigma_{2024}\) to consist of (W) and every valid normal absorption over a fixed 23-letter alphabet, of length at most 25 in each word and at most three right words. Define \(Sigma_{2025}\) similarly over 13 letters, lengths at most 15, and at most two right words. These are finite, explicitly specified sets; the lemma decides membership. No small or minimal enumeration is claimed.

Theorem. \(Sigma_{2024}\) and \(Sigma_{2025}\) are finite identity bases of the respective algebras. Their word reversals are bases of the opposites.

Proof. They are sound by (W) and definition. Any valid word absorption has the bounded original witnesses supplied by Section 3. Its normal form is obtained from a member of the specified finite family by the reversible substitution of Section 4. Thus every valid word absorption is derivable. Distribute arbitrary terms into word sums; an identity between the sums is equivalent to each summand on each side being absorbed by the other side. This derives all identities, without any restriction on variable count. Word reversal transports the bases and derivations to the opposite algebras.

The finite computations independently check the tables, word formula, criterion, witness selection, and reversible compression. Completeness is the unrestricted proof above. Classification transfers use verified equality of generated varieties in both HSP directions.

\subsection{Twelve complete identities for 13219}\label{case:55}
Signature (+, multiplication), without constants. Relative to the additively idempotent semiring axioms, use

{\small\par\smallskip
\CaseLine{\(xy=yx; xyz=xyzz; xxy\leq xxx+yyy; xxy\leq xxx+xy;\)}
\CaseLine{\(xxyz\leq xxy+xz; xy\leq xxy; x\leq xxx; xxy\leq x+xy;\)}
\CaseLine{\(xxx\leq x+xxz; xxy\leq xy+xxz; xxy\leq xxx+yyz;\)}
\CaseLine{\(xz\leq xx+xy+yz.\)}
\smallskip}

Here \(u\leq P\) abbreviates \(u+P=P\). Addition is 00444/\allowbreak{}01444/\allowbreak{}44244/\allowbreak{}44434/\allowbreak{}44444: \(1<0<4\), while 2 and 3 are incomparable with each other and with 0,\allowbreak{}1; those joins are 4. Multiplication is 44244/\allowbreak{}41244/\allowbreak{}22222/\allowbreak{}44204/\allowbreak{}44244. A nonlinear word has value 2 if a variable is 2, otherwise1 if all variables are 1, otherwise0 for a quadratic with both occurrences 3, and 4 otherwise. Literals keep their values. This follows by induction.

\par\smallskip\noindent\textit{Exact original-witness criterion.}\quad 
Form the global graph of original nonlinear words of P, joining variables in each such word. Let V be its vertices and L its original literal variables. A component is active if it contains an original long word or a vertex in L intersect V. Let A be the union of active components. Separately let C mean that an original right word has content contained in \(U=c(w)\).

The criterion for absorption \(w\leq P\) is:

\textbullet{} An original literal is absorbed. A missing literal requires C and its variable in A.

\textbullet{} A long target requires C and U contained in A.

\textbullet{} A quadratic requires C and either both endpoints in A or both endpoints in V in one original component.

Unlike 13218, global witnesses may have outside variables. Unlike 13841, a missing literal cannot be recovered from a pure quadratic without activation. A literal incident at a nonlinear word activates its component here; a literal alone provides no graph edge.

If C fails, assign target variables 1 and outsiders 2. Every right term is 2, while the target is 1, not below 2. For a long target with an inactive variable, assign A variables 1, inactive V variables 3, and all remaining variables 0. All original long words are 1; quadratics are 1 or 0; original literals in V are active, and other literals are 0. The right join is at most 0, but the target is 4. For a missing literal with C but no activation, its contained witness must be a pure quadratic; the same assignment makes the target 3 and the right join at most 0, which does not dominate 3.

For a quadratic failing the component condition, if an endpoint is outside V, assign it 0 and all other variables 1. The right join is at most 0 and the target is 4. Otherwise choose an inactive component containing exactly one target endpoint, assign that component 3, other V variables 1, and remaining variables 0. An inactive component has no original literal at any vertex. Thus right values again join to at most 0, while the target is 4. These assignments establish necessity without size bounds.

For sufficiency, if P evaluates to 4 there is nothing to prove. If it evaluates to 2, every original term must be 2. C then forces a target variable to 2, making the target 2. If it evaluates to 3, every right term is a literal assigned 3: nonlinear words never have value 3. This case is possible for an original literal target, but impossible in the other asserted cases, which require nonlinear original witnesses. Otherwise the right value is 0 or 1. No original term is 2,\allowbreak{}3, or 4, so all original nonlinear variables are 1 or 3, quadratic edges make their endpoints equal, and long words force 1. An original literal at a nonlinear vertex cannot be 3 in this case, so literal incidence also forces its component 1. All A variables are 1. Long targets on A and missing literals in A are 1; a quadratic on A is 1, and one in an inactive connected component is 1 or 0, also supplied by an original quadratic in that component. This proves sufficiency for arbitrary P.

\par\smallskip\noindent\textit{Complete finite derivations.}\quad 
The two word laws normalize long words by content, by moving a chosen variable to the end and duplicating or deleting it with \(xyz=xyzz\), always using nonempty preceding blocks and retaining length at least three. Canonical long forms are sorted content for at least three variables, xxy for two, and xxx for one. Thus \(g^{3}=g\) for long g.

For targets whose variables are in A, retain all original long words. For every additional activating literal x, choose an original incident quadratic xy and use \(xxy\leq x+xy\) to derive a long seed. An incident long word already supplies a seed. Merge seeds using \(xxy\leq xxx+yyy\). Extend their content along original quadratic edges with \(xxy\leq xxx+xy\) for a pure seed, or \(xxyz\leq xxy+xz\) for a nonpure seed represented as xxp, with p nonempty. This produces a long g whose content is A; external variables are deliberately retained.

Now choose the separate original contained witness u. If it is x, represent g as xxq and use \(xxx\leq x+xxz\) with \(z=q\); q is nonempty, using x itself when the content of g is only x. If \(u=xy\) is quadratic, use \(xxy\leq xy+xxz\) with the same representation. If u is long, retain it. This yields a long h contained in U. For each missing target variable y, represent g as yyq and apply \(xxy\leq xxx+yyz\) with \(x=h\), \(z=q\). Its right terms are h and g after normalization; its left adds y to the content of h. Repeat to reach U, then normalize, recover a quadratic by \(xy\leq xxy\), or recover a literal by \(x\leq xxx\). In particular a missing literal with a pure quadratic uses its global activation to upgrade that witness before literal recovery; \(x\leq xx\) is never used.

For a quadratic in an inactive component, an original target needs no work. Otherwise its contained witness must be an original endpoint loop. A contained mixed edge would already be the target, a contained long word would activate it, and an original endpoint literal would also activate that component. Retain this original loop xx and a positive path of original quadratic edges to the other endpoint. Starting with the first original edge, repeatedly apply \(xz\leq xx+xy+yz\). This derives the target. A square in the inactive case must already be the original loop.

Every substitution is nonempty. The construction is unrestricted in variables, words, components, and path lengths. Distributivity and absorption in both directions derive every identity. All twelve laws satisfy the criterion and are checked on all assignments. The proof, rather than finite regression, establishes finite basability.

The distinctions are essential: \(x\leq xx\) fails at \(x=3\); \(xy\leq xx+yy\) fails at \(x=3\), \(y=1\). Without a contained witness, \(xxx\leq yyy\) fails at \(x=1\), \(y=2\). Without activation, \(xyz\leq xy+yz\) fails when all three variables are 3. These prevent copying the preceding criteria.

\subsection{Endpoint witnesses and final-variable saturation}\label{case:56}
All terms are in the signature (+, multiplication), without constants. The bases below are relative to the additively idempotent semiring axioms. For nonempty words u, let c(u), h(u), f(u) and i(u) denote its content, first variable, final variable, and the set of variables before its final position. Put \(i(y)=empty\) for a linear word. An inequality \(w\leq P\) abbreviates w+P approximately equal to P.

Let A and B be catalogue entries 9106 and 9824. Element labels are not term constants. Tables in numerical order are:

{\small\par\smallskip
\noindent A addition     A multiplication    B addition     B multiplication\par
\CaseLine{\ttfamily 0 0 0 0 4\qquad 1 1 1 1 1\qquad 0 0 4 4 4\qquad 1 1 2 3 4}
\CaseLine{\ttfamily 0 1 0 1 4\qquad 1 1 1 1 1\qquad 0 1 4 4 4\qquad 1 1 2 3 4}
\CaseLine{\ttfamily 0 0 2 2 4\qquad 1 1 1 1 1\qquad 4 4 2 2 4\qquad 4 4 4 3 4}
\CaseLine{\ttfamily 0 1 2 3 4\qquad 3 3 3 3 3\qquad 4 4 2 3 4\qquad 4 4 4 3 4}
\CaseLine{\ttfamily 4 4 4 4 4\qquad 0 1 0 1 4\qquad 4 4 4 4 4\qquad 4 4 4 3 4}
\smallskip}

\par\smallskip\noindent\textit{Semiring axioms and arbitrary word values.}\quad 
The additive order of A is \(3<1<0<4\) and \(3<2<0<4\). Its only incomparable pair is 1,\allowbreak{}2, with join 0. Each multiplication row and column is increasing and preserves this join, so it preserves every binary join. For associativity, a first factor 3 gives 3 under either association, and a first factor in \{0,\allowbreak{}1,\allowbreak{}2\} gives 1. When the first factor is 4 and the middle factor is not 4, both associations give 1. When both are 4, either association gives 4, 0 or 1 according as the last factor is 4, in \{0,\allowbreak{}2\}, or in \{1,\allowbreak{}3\}. Thus A is a semiring.

Induction gives the value of every nonlinear word w in A:

\textbullet{} \(h(w)=3\) gives 3.

\textbullet{} h(w) in \{0,\allowbreak{}1,\allowbreak{}2\} gives 1.

\textbullet{} All factors equal to 4 gives 4.

\textbullet{} All factors before the final position equal to 4, and the final factor in \{0,\allowbreak{}2\}, gives 0.

\textbullet{} The remaining case gives 1.

In particular nonlinear values never equal 2. The following word identities hold:

\CaseDisplay{xyzt\approx xzyt,\qquad x^2y\approx xy,\qquad xyx\approx xy^2. \quad\textup{(\ensuremath{A1})}}

They permute interior positions and remove repeated interior variables, including a head repeated internally. A word whose final variable is new becomes h followed by the other members of i(w) in fixed order, then its final variable. If the final variable has an earlier occurrence, the normal word is h followed by c(w) minus \{h\}, in fixed order, followed by h. To justify the latter when the final variable s differs from h, first obtain h R s s, insert an extra h using \(x^{2}y=xy\) backwards, permute interiors to h R h s s, and apply \(hs^{2}=hsh\) to obtain h R h s h. Move the redundant internal h next to the initial h and delete it. This also works with empty R: the actual substituted suffixes are always nonempty. The normal form preserves h, c, and i, including whether \(i=c\). Linear words remain linear.

The additive order of B consists of \(1<0<4\) and \(3<2<4\); every pair across the two lower chains has join 4. Put D=\{0,\allowbreak{}1\}. Multiplication is increasing and preserves all these cross-chain joins in either coordinate, hence distributes over addition. For associativity, a last factor 3 or 4 makes either association 3 or 4. With last factor in D, either association is 1 exactly when both earlier factors are in D, and otherwise 4. With last factor 2 it is 2 exactly in the same case, and otherwise 4. Thus B is a semiring.

The nonlinear value in B is:

\textbullet{} final value 3 gives 3;

\textbullet{} all factors in D gives 1;

\textbullet{} final value 2 with every earlier factor in D gives 2;

\textbullet{} the remaining case gives 4.

It never equals 0. The word identities

\CaseDisplay{xyz\approx yxz,\qquad x^2y\approx xy       \quad\textup{(\ensuremath{B1})}}

permute all positions before the final one and remove their duplicates, with a nonempty suffix in every application. The normal form lists i(w) in fixed order and then f(w). The final variable may also occur in i(w); i(w) is nonempty for every nonlinear word.

\par\smallskip\noindent\textit{Two original endpoint witnesses for A.}\quad 
Lemma 1. In A a linear word y is absorbed by a polynomial P if and only if y is a literal summand. A nonlinear w is absorbed if and only if both conditions hold:

1. An original nonlinear word a of P has \(h(a)=h(w)\).

2. An original right word b is either a linear variable in c(w), or a nonlinear word with c(b) contained in c(w) and i(b) contained in i(w).

Proof. If a missing literal y is assigned 2 and every other variable 3, all right terms have value 1 or 3, whereas \(y=2\) is not below either. If the first condition fails, assign \(h(w)=2\) and all other variables 3. Right nonlinear terms are 3, right linear terms are 2 or 3, and \(w=1\) is not below their sum.

If the second condition fails and \(s=f(w)\) is absent from i(w), assign every variable in i(w) the value 4, s the value 2 and all outsiders 3. Then \(w=0\). There is no right linear variable in c(w). A nonlinear right word containing an outsider is 3 if its head is outside, and 1 otherwise. Any right word entirely in c(w) has s in its prefix, so also gives 1. Hence \(P\leq 1\), separating w. If instead \(i(w)=c(w)\), absence of the second witness means that no right word is contained in c(w). Assign c(w) the value 4 and outsiders 3. Now \(w=4\) while again \(P\leq 1\).

For sufficiency consider the value of P. If it is 3 or 2, the nonlinear same-head witness must have value 3, since nonlinear values omit 2. Its head, and therefore the target head, is 3, so \(w=3\). If \(P=0\), the only possible violating nonlinear value is 4. Then all target variables are 4, so the second witness is also 4, a contradiction. If \(P=1\), the possible violating nonlinear values are 0 and 4. The latter was just excluded. For \(w=0\), all prefix variables are 4 and the unique final variable is 0 or 2. A linear second witness has value 4, 0 or 2, none below 1; a nonlinear second witness has value 4 or 0, also not below 1. This excludes the former. \(P=4\) is automatic. These cases exhaust A. \(\square\)

Thus each valid nonlinear absorption \(w\leq P\) already has a valid local absorption \(w\leq a+b\) with at most two original right words. In this local step there are at most three words. Reserve every first and final variable individually, at most six. Group the remaining variables by their nonempty incidence patterns across the three interiors, at most seven groups. Identify each group with one variable and normalize using (A1). Ordinary identification preserves validity. The resulting local step uses at most thirteen variables and normal word length at most fifteen.

For reversal of this compression, substitute for each group variable the nonempty product of its members; reserved endpoint variables remain single letters. The group members have identical interior incidence and occur nowhere at an endpoint. Interior permutations and duplicate removal recover each original normal word. Hence every such local absorption follows from the valid normal local absorptions within these fixed bounds, together with (A1).

\par\smallskip\noindent\textit{Exact final-variable criterion for B.}\quad 
For the original right words of P put

\CaseDisplay{U=\bigcup_{u\in P}c(u),\quad F=\{f(u):u\in P\},\quad
 B_t=\bigcup_{\substack{u\in P,\ |u|\geq2\\ f(u)=t}}i(u).}

Here F includes the variable of each linear summand.

Lemma 2. A linear y is absorbed in B exactly when it is a literal summand. For a nonlinear target w with \(t=f(w)\), absorption holds exactly when

\CaseDisplay{c(w)\subseteq U,\qquad t\in F,\qquad
       \bigl(B_t\cap F\ne\varnothing\ \text{or}\ i(w)\subseteq B_t\bigr).
                                                               \quad\textup{(\ensuremath{5})}}

Proof. A missing linear y is separated by \(y=0\) and every other variable 1: every right nonlinear word is 1 and so are all right linear words. If a target variable x is absent from U, assign \(x=4\) and all others 1. The target is 4 while \(P=1\). If t is absent from F, assign \(t=0\) and all others 3. All right words are 3, while the nonlinear target is 1 or 4, neither below 3.

For failure of the last condition, choose x in i(w) minus \(B_{t}\). We may assume t in F. Assign all \(B_{t}\) the value 1, t the value 2, and every other variable 3; these assignments are consistent because \(B_{t}\) and F are disjoint. Each right word ending in t has value 2; each other right word has final value 3 and thus value 3, including linear words. Therefore \(P=2\). The target has a prefix occurrence outside \(B_{t}\), assigned 2 or 3, and final value 2, so \(w=4\).

Conversely suppose (5). If P belongs to \{0,\allowbreak{}1\}, every right nonlinear word must be 1, and every variable of every right word is in D. Content coverage forces all target variables into D, so \(w=1\leq P\). If \(P=3\), every right word is 3 and its final variable is 3. Since t in F, \(w=3\). If \(P=2\), every right word is 2 or 3, so all final variables in F have values in \{2,\allowbreak{}3\}. If \(t=3\) the target is 3. If \(t=2\), every right nonlinear word ending in t must be 2, hence all \(B_{t}\) variables are in D. The sets D and \{2,\allowbreak{}3\} are disjoint; thus \(B_{t}\) and F cannot intersect in this case, and (5) gives i(w) contained in \(B_{t}\). The target is therefore 2. Finally \(P=4\) is automatic. These are all values, proving sufficiency. \(\square\)

\par\smallskip\noindent\textit{Local union and saturation derivations for B.}\quad 
Fix a valid nonlinear \(w\leq P\), put \(t=f(w)\), \(I=i(w)\), and retain the sets in (5). For a nonempty set E write q(E, t) for its ordered listing followed by t. The following construction uses at most two right words in every local step.

If \(B_{t}\) and F are disjoint, I is contained in \(B_{t}\) and \(B_{t}\) is nonempty. Take all original nonlinear words ending in t and successively replace two by q(\(E_{1}\) union \(E_{2}\), t), where \(E_{j}\) is their prefix set. Lemma 2 makes every union step valid: the final variable is witnessed, the contents are covered, and the target prefix is contained in the union of right prefixes. Eventually q(\(B_{t}\), t)\ensuremath{\leq}P is derived. The same lemma gives w\ensuremath{\leq}q(B\_t, t).

If instead choose g in \(B_{t}\) intersect F, take an original nonlinear a ending in t with g in i(a), and an original b ending in g. Put \(E=c(a)\) union c(b). Lemma 2 gives

\CaseDisplay{q(E, t)\leq a+b.                    \quad\textup{(\ensuremath{6})}}

Indeed the final t occurs on the right, all target content is covered, and g is both a right prefix variable for final t and a right final variable. Since t is in E, the new word q(E, t) contains its own final variable in its prefix. It is consequently stable under further content union: for any original right u,

\CaseDisplay{q(E\cup c(u),t)\leq q(E, t)+u.              \quad\textup{(\ensuremath{7})}}

The right prefix for final t now contains t itself, which is also a right final variable. After finitely many such unions \(E=U\). Since c(w) is contained in U, Lemma 2 gives w\ensuremath{\leq}q(U, t). In both cases every new word is already derivably below P before it is used in a later step. This is an equational derivation, not an assumption that arbitrarily many original witnesses fit in one bounded step.

To compress a local step, reserve its at most three final variables individually and group the other variables by their nonempty prefix incidence patterns across the at most three words. There are at most seven groups. Identification and (B1) preserve validity and give at most ten variables and normal word length eleven. Inverse substitution of each group by the nonempty product of its members, followed by prefix permutations and duplicate removal, recovers all original words. No final variable is merged into a group.

\par\smallskip\noindent\textit{Finite bases and unrestricted completeness.}\quad 
Theorem. Relative to the additively idempotent semiring axioms, A has the basis consisting of (A1) and all valid normal absorptions over a fixed thirteen-letter alphabet, with word length at most fifteen and at most two right words. B has the basis consisting of (B1) and all valid normal absorptions over a fixed ten-letter alphabet, with word length at most eleven and at most two right words. Their multiplicative opposites have the reversed bases.

Proof. These are finite, explicitly specified families; Lemmas 1 and 2 decide membership. The two-witness argument for A and the local union/saturation argument for B, with the reversible bounded compressions, derive every valid word absorption with arbitrarily many variables and original right words. Linear absorptions need only literal inclusion. Expand arbitrary terms into nonempty sums of words by distributivity. In each valid identity, every summand of either side is absorbed by the other side; derive each absorption and sum them to obtain both inequalities, hence the identity. Every substitution used above is nonempty. Reversal transports the complete derivation to the opposite. \(\square\)

No small or minimal basis has been enumerated, and no hereditary finite-basis property is asserted. Finite computations serve as auxiliary checks. Classification of any further catalogue type requires a verified equivalence of generated varieties or another complete theorem.

\subsection{正则化分支的有限基补全：135个五元素加法幂等半环}\label{case:57}

结论。 下表135个半环均为有限基。证明由有限基构造和完整的双向等簇证书组成，补全此前只排除极限性的135个正则化类型；连同此前已解决的9个，原144个正则化类型全部获得有限基证明。其他未定类型不在该结论内。

\par\smallskip\noindent\textit{1. 获得有限基底代数.}\quad 
引理1（对角平方）。 若 \(D\leq B^2\) 且 \(D\) 包含所有 \((b, b)\)，则 \(\mathsf V(D)=\mathsf V(B)\)。

证明：\(D\leq B^2\) 给出一个包含方向；对角映射 \(B\hookrightarrow D\) 给出反向包含。因此已知 \(D\) 有限基时，\(B\) 具有完全相同的恒等式基。证毕。

引理2（次直积分量）。 若 \(D\leq\prod_iB_i\)，每个坐标投影都满，则 \(\mathsf V(D)=\bigvee_i\mathsf V(B_i)\)。

证明：子代数包含给出左侧包含于右侧；各满投影表明每个 \(B_i\) 是 \(D\) 的同态像，给出反向包含。证毕。

17个代表的四元素分量 \(B\) 适用引理1，有一个已知有限基的五元素代数 \(D\leq B^2\) 包含对角。另8个代表具有三元素分量 \(B\) 和二元素分量 \(L\)，已知有限基的五元素代数 \(D\) 次直嵌入 \(B^2\times L\)；此时直接以 \(D\) 为底代数。分量标签、积中的元素、同构排列及既有有限基来源均列入证书。

引理3（最小元与有界词长）。 设有限半环 \(B\) 有 \(n\) 个元素和加法最小元。若有限恒等式组 \(\Gamma\) 将所有项化为长度至多 \(\ell\) 的非空词的非空和，则 \(B\) 有限基。可以取 \(\Gamma\) 加固定 \(b=\ell(1+n^{\ell+1})\) 个变量上的全部有效有界词长吸收式。

证明：对有效吸收 \(w\leq Q\)，令 \(C=c(w)\)，\(d=|C|\leq\ell\)。对 \(C\) 的每个赋值，将其余变量赋为加法最小元。\(Q\) 的词至多取 \(n\) 种不同值，每种选一实际求和词；合并全部 \(n^d\) 个赋值所选词，得到至多 \(n^{d+1}\) 个词的非空子和 \(Q_0\)。任意完整赋值在 \(C\) 外都不小于最小元扩张，乘法保序，因此仍有 \(w\leq Q_0\)。该式至多涉及 \(\ell+\ell n^{\ell+1}=b\) 个变量；重命名后属于上述有限库，添加其余求和项得到 \(w\leq Q\)。将任意恒等式展开，逐项处理两向吸收，得完整性。最小元只用于赋值，没有加入签名。证毕。

1222、5373的四元素分量分别满足 \(xyz\approx xy\) 和 \(xyz\approx yz\)，并有加法最小元。每个词可反复缩短至长度至多2。引理3给出 \(n=4,\ell=2, b=130\)。130是有限库的变量数界，不是基条数；不声称最小基。

\par\smallskip\noindent\textit{2. 共同吸收元扩张的有限基.}\quad 
引理4。 若 \(B\) 有有限基 \(\Sigma\)，且存在恰含 \(x, y\) 的项 \(d(x, y)\) 使 \(B\models d(x, y)\approx x\)，则向 \(B\) 添一个对两种运算都吸收的新元 \(\infty\) 所得 \(B^\infty\) 有限基。

证明：\(B^\infty\) 的恒等式恰为 \(B\) 的正则恒等式（两边变量内容相同）。正则式遇到新元赋值，两边同时为 \(\infty\)；其余赋值归结为 \(B\)。非正则式把仅在一侧出现的变量赋为 \(\infty\)、其他变量赋在 \(B\) 中，两边不同。

取以下七条分划项恒等式，以及每条基恒等式的填充：

\CaseDisplay{\begin{aligned}
d(x, x)&\approx x,\\
d(d(x, y),z)&\approx d(x, d(y, z)),\\
d(x, d(y, z))&\approx d(x, d(z, y)),\\
d(f(x, y),z)&\approx f(d(x, z),d(y, z)) &&(f=+,\cdot),\\
d(z, f(x, y))&\approx d(d(z, x),y) &&(f=+,\cdot),\\
d(s, t)&\approx d(t, s) &&(s\approx t\in\Sigma).
\end{aligned}}

它们正则且在 \(B\) 中成立，故在 \(B^\infty\) 中成立。

对有限非空变量集 \(C\)，令 \(D_C(u)\) 为把 \(C\) 的变量依次放入 \(d\) 的第二位置、逐层填充 \(u\) 所得项。前三条允许交换填充列表中的变量和删除重复；第五类将第二位置任意项拆为其变量列表。故

\CaseDisplay{d(u, v)\approx D_{c(v)}(u),\qquad D_{c(u)}(u)\approx u.}

第四类给出 \(D_C(f(u, v))\approx f(D_C(u),D_C(v))\)。只要基式代换实例 \(s\sigma, t\sigma\) 的变量均在 \(C\) 中，将其填充式再次填充 \(C\)，删除重复列表项，就得 \(D_C(s\sigma)\approx D_C(t\sigma)\)，并可在任意上下文中使用。

对有效正则式 \(u\approx v\)，置 \(C=c(u)=c(v)\)。在 \(\Sigma\) 中取有限推导，把所有额外辅助变量替换为固定的一个 \(C\) 中变量，首尾不变。将各步填充 \(C\) 得 \(D_C(u)\approx D_C(v)\)，最后删除填充得 \(u\approx v\)。通常半环公理本身正则，其推导步骤也保持有效。证毕。

本批每个底代数实际满足 \(d(x, y)=x+w(x, y)\approx x\)，下表列出 \(w\)。平方与缩短律两种情况的底代数有4个元素；次直积情况的底代数有5个元素。因此相应扩张有5或6个元素。六元素扩张与五元素目标等簇不要求同构。

\par\smallskip\noindent\textit{3. 精确等簇证书与任意变量结论.}\quad 
引理5。 设 \(S, T\) 为有限加法幂等半环，\(g=(g_1,\ldots, g_k)\) 生成 \(T\)。枚举 \(S\) 上完整的 \(k\) 元非空词函数半群，直到右乘每个变量都不产生新函数，为每个函数选代表词 \(w\) 并记录 \(w(g)\)。假设：

1. 每条右乘变量的转移都与 \(T\) 中相应词值相容；

2. 对每个 \(b\in T\)，令 \(J_b\) 为满足 \(w(g)\leq b\) 的全部词函数。若 \(J_b\) 非空，任何点态不大于 \(\sum J_b\) 的词函数 \(u\) 都满足 \(u(g)\leq b\)。

则 \(T\in\mathsf{HSP}(S)\)。

证明：第一条件保证词函数到词值的映射良定义。每个项函数是非空词函数的非空和。若两个这样的和在 \(S\) 上相等，令右侧各词在 \(g\) 上的和为 \(b\)。右侧词均在 \(J_b\) 中，左侧每个词函数不大于 \(\sum J_b\)；第二条件推出每个左词值不大于 \(b\)。反向同理，故两和的值相等。因此从 \(S\) 的 \(k\) 元项函数代数到 \(T\) 的映射良定义、保持两个运算；生成性保证满射。该项函数代数是 \(S^{S^k}\) 的子代数，结论成立。证毕。

两向应用该引理即证明等簇，因而完整恒等式理论及有限基可转移。这里闭包、所有转移、所有主下集和生成性共同构造实际满同态，不是通过固定变量数的相同测试猜测等簇。

独立核验对每一个目标编号直接验证与所构造 \(B^\infty\) 的双向成员关系，不只检查27个代表。积中坐标运算也由原目录重新计算。最终等簇判定不依赖未经证明的分解猜测。

\par\smallskip\noindent\textit{4. 全部类型与来源.}\quad 
``平方''为引理1；``次直积''为引理2且直接取所列已知有限基代数为底代数；``缩短律''为引理3。既有来源的理论文件和运算证书均列入核验依赖清单。

\noindent 代表\enspace\textbar\enspace 本类全部编号\enspace\textbar\enspace 底代数有限基来源\enspace\textbar\enspace \(w(x, y)\)\par
\noindent 1222\enspace\textbar\enspace 1222, 2132, 2145, 2895\enspace\textbar\enspace \(xyz\approx xy\)；130变量有限库\enspace\textbar\enspace \(xxy\)\par
\noindent 1223\enspace\textbar\enspace 1223, 2136, 2149, 6616\enspace\textbar\enspace 平方；已知有限基 2184 (case~\ref{grp:123})\enspace\textbar\enspace \(xy\)\par
\noindent 1643\enspace\textbar\enspace 1643, 4923, 7207, 7921, 8904, 8955, 8961, 8977\enspace\textbar\enspace 次直积；已知有限基 8983 (case~\ref{grp:320})\enspace\textbar\enspace \(xy\)\par
\noindent 4571\enspace\textbar\enspace 4571, 4924, 7049, 7123\enspace\textbar\enspace 平方；已知有限基 4590 (case~\ref{grp:230})\enspace\textbar\enspace \(xy\)\par
\noindent 1649\enspace\textbar\enspace 1649, 4932, 7219, 7956, 9175, 9200, 9254, 9276\enspace\textbar\enspace 次直积；已知有限基 9305 (case~\ref{grp:321})\enspace\textbar\enspace \(xy\)\par
\noindent 4580\enspace\textbar\enspace 4580, 4934, 7069, 7135\enspace\textbar\enspace 平方；已知有限基 4594 (case~\ref{grp:232})\enspace\textbar\enspace \(xy\)\par
\noindent 1680\enspace\textbar\enspace 1680, 5412, 8512, 8682, 8909, 8956, 8962, 8978\enspace\textbar\enspace 次直积；已知有限基 8914 (case~\ref{grp:123})\enspace\textbar\enspace \(xy\)\par
\noindent 5300\enspace\textbar\enspace 5300, 5416, 8399, 8455\enspace\textbar\enspace 平方；已知有限基 5321 (case~\ref{grp:288})\enspace\textbar\enspace \(xy\)\par
\noindent 1778\enspace\textbar\enspace 1778, 9176, 9201, 9277, 10057, 10072\enspace\textbar\enspace 次直积；已知有限基 9307 (case~\ref{grp:26})\enspace\textbar\enspace \(xy\)\par
\noindent 2797\enspace\textbar\enspace 2797, 7408, 7623, 7631, 8963, 9410, 9583, 9643\enspace\textbar\enspace 次直积；已知有限基 9886 (case~\ref{grp:324})\enspace\textbar\enspace \(yx\)\par
\noindent 4848\enspace\textbar\enspace 4848, 7201, 7409, 7910\enspace\textbar\enspace 平方；已知有限基 7406 (case~\ref{grp:322})\enspace\textbar\enspace \(yx\)\par
\noindent 2801\enspace\textbar\enspace 2801, 7425, 7683, 7691, 9206, 9447, 9603, 9657\enspace\textbar\enspace 次直积；已知有限基 9651 (case~\ref{grp:123})\enspace\textbar\enspace \(yx\)\par
\noindent 4866\enspace\textbar\enspace 4866, 7211, 7429, 7941\enspace\textbar\enspace 平方；已知有限基 7422 (case~\ref{grp:288})\enspace\textbar\enspace \(yx\)\par
\noindent 2810\enspace\textbar\enspace 2810, 8035, 8606, 8616, 8964, 9411, 9587, 9644\enspace\textbar\enspace 次直积；已知有限基 9887 (case~\ref{grp:365})\enspace\textbar\enspace \(yx\)\par
\noindent 5383\enspace\textbar\enspace 5383, 8037, 8502, 8668\enspace\textbar\enspace 平方；已知有限基 8033 (case~\ref{grp:360})\enspace\textbar\enspace \(yx\)\par
\noindent 2844\enspace\textbar\enspace 2844, 9207, 9448, 9658, 10062, 10070\enspace\textbar\enspace 次直积；已知有限基 9910 (case~\ref{grp:26})\enspace\textbar\enspace \(yx\)\par
\noindent 4574\enspace\textbar\enspace 4574, 7055, 7128, 7566\enspace\textbar\enspace 平方；已知有限基 7321 (case~\ref{grp:323})\enspace\textbar\enspace \(xy\)\par
\noindent 5285\enspace\textbar\enspace 5285, 7593, 8370, 8439\enspace\textbar\enspace 平方；已知有限基 8719 (case~\ref{grp:123})\enspace\textbar\enspace \(yx\)\par
\noindent 5294\enspace\textbar\enspace 5294, 7594, 8384, 8452\enspace\textbar\enspace 平方；已知有限基 8731 (case~\ref{grp:323})\enspace\textbar\enspace \(xy\)\par
\noindent 5373\enspace\textbar\enspace 5373, 6725, 8495, 8654\enspace\textbar\enspace \(xyz\approx yz\)；130变量有限库\enspace\textbar\enspace \(yxx\)\par
\noindent 6830\enspace\textbar\enspace 6830, 10123, 10208\enspace\textbar\enspace 平方；已知有限基 10604 (case~\ref{grp:307})\enspace\textbar\enspace \(xxy\)\par
\noindent 6831\enspace\textbar\enspace 6831, 10140, 10227\enspace\textbar\enspace 平方；已知有限基 10654 (case~\ref{grp:307})\enspace\textbar\enspace \(xy\)\par
\noindent 6839\enspace\textbar\enspace 6839, 10276, 10334\enspace\textbar\enspace 平方；已知有限基 10761 (case~\ref{grp:307})\enspace\textbar\enspace \(xy\)\par
\noindent 7651\enspace\textbar\enspace 7651, 10129, 10209\enspace\textbar\enspace 平方；已知有限基 10610 (case~\ref{grp:307})\enspace\textbar\enspace \(yx\)\par
\noindent 7653\enspace\textbar\enspace 7653, 7696, 10163, 10176, 10243, 10257\enspace\textbar\enspace 平方；已知有限基 10688 (case~\ref{grp:323})\enspace\textbar\enspace \(xy\)\par
\noindent 8086\enspace\textbar\enspace 8086, 8620, 10149, 10154, 10229, 10230\enspace\textbar\enspace 平方；已知有限基 10664 (case~\ref{grp:323})\enspace\textbar\enspace \(xy\)\par
\noindent 8791\enspace\textbar\enspace 8791, 10134, 10211\enspace\textbar\enspace 平方；已知有限基 10614 (case~\ref{grp:307})\enspace\textbar\enspace \(yx\)\par
合计：平方67型，次直积60型，有界词长8型，共135型；其中交换3型、非交换132型。

\subsection{二阶群的两种膨胀：四组完整有限基}\label{case:58}
\par\smallskip\noindent\textit{1. 构造与基.}\quad 
核心为二阶群的平坦扩张 \(K=\{e, g, c\}\)：\(g^2=e\)，\(c\) 是乘法零，\(e, g<c\)。对下面四个并半格，定义映射 \(\rho\) 固定 \(K\)，并以 \(uv=\rho(u)\rho(v)\) 定义乘法。

\noindent 半环\enspace\textbar\enspace 加法序的覆盖或生成关系\enspace\textbar\enspace 其余 \(\rho\) 值\enspace\textbar\enspace 目录来源\par
\noindent \(F_e\)\enspace\textbar\enspace \(e, g<c<t,\ e<d<t\)；未规定的元素不可比\enspace\textbar\enspace \(d, t\mapsto c\)\enspace\textbar\enspace 1522\par
\noindent \(F_g\)\enspace\textbar\enspace \(e, g<c<t,\ g<d<t\)；未规定的元素不可比\enspace\textbar\enspace \(d, t\mapsto c\)\enspace\textbar\enspace 1724\par
\noindent \(G_e\)\enspace\textbar\enspace \(e<d<c,\ g<b<c\)；两支之间不可比\enspace\textbar\enspace \(d\mapsto c,\ b\mapsto g\)\enspace\textbar\enspace 10822\par
\noindent \(G_g\)\enspace\textbar\enspace \(g<b<c,\ g<d<c,\ e<c\)；\(b, d\) 不可比，\(e\) 与二者不可比\enspace\textbar\enspace \(d\mapsto c,\ b\mapsto g\)\enspace\textbar\enspace 10821\par
上述 \(\rho\) 均保持并且固定核心，故确实定义半环。\(F_e, F_g\) 在顺序 \((e, g, c, d, t)\) 下的目录同构分别为 \((1,\allowbreak{}2,\allowbreak{}0,\allowbreak{}3,\allowbreak{}4)\)、\((2,\allowbreak{}1,\allowbreak{}0,\allowbreak{}3,\allowbreak{}4)\)。\(G_e, G_g\) 在顺序 \((e, g, c, b, d)\) 下的目录同构均为 \((3,\allowbreak{}1,\allowbreak{}4,\allowbreak{}0,\allowbreak{}2)\)。

令 \(\Delta\) 为以下五条恒等式：

\CaseDisplay{xy\approx yx,\quad x^3y\approx xy,\quad
x^2y^2\approx x^2+y^2,\quad
xyz\leq x^3+y^3+z^3,\quad
x^3\leq x+xy.}

另记

\CaseDisplay{(a)\quad x^3\leq x+y^2,\qquad
(b)\quad x\leq x^2+x^3.}

定理。 以下是四个半环的完整有限基：

\noindent 半环\enspace\textbar\enspace 基\enspace\textbar\enspace 附加恒等式数\par
\noindent \(F_e\)\enspace\textbar\enspace \(\Delta\)\enspace\textbar\enspace 5\par
\noindent \(F_g\)\enspace\textbar\enspace \(\Delta,(a)\)\enspace\textbar\enspace 6\par
\noindent \(G_e\)\enspace\textbar\enspace \(\Delta,(b)\)\enspace\textbar\enspace 6\par
\noindent \(G_g\)\enspace\textbar\enspace \(\Delta,(a),(b)\)\enspace\textbar\enspace 7\par
有效性由核心乘法直接得到。\(x^3\leq x+xy\) 只有 \(\rho(x)=c\) 的情况可能不是显然的序关系，此时 \(xy=c\)。式 (a) 在下标为 \(g\) 的两例中使用 \(d+e\geq c\)。式 (b) 在 \(G_e, G_g\) 中仅需检查新元素：\(b^2+b^3=e+g=c\geq b\)，\(d^2+d^3=c\geq d\)。

\par\smallskip\noindent\textit{2. 精确判据.}\quad 
对非线性左词 \(w\)，两种下标分别对应以下条件：

\CaseDisplay{\begin{array}{ll}
E:&\pi(w)\in H(P),\quad c(w)\subseteq N(P);\\[2pt]
J:&\pi(w)\in H(P),\quad c(w)\subseteq U(P),\quad
\bigl(c(w)\subseteq N(P)\ \text{或}\ 0\in H(P)\bigr).
\end{array}}

在 \(F_e, G_e\) 中，\(w\leq P\) 当且仅当 \(E\)；在 \(F_g, G_g\) 中，当且仅当 \(J\)。

对于一次左词 \(x\)：

\textbullet{} 在 \(F_e, F_g\) 中，\(x\leq P\) 当且仅当 \(x\) 本身是右侧一次项。

\textbullet{} 在 \(G_e, G_g\) 中，\(x\leq P\) 当且仅当 \(x\) 是右侧一次项，或者 \(\pi(x),0\in H(P)\)。

\par\smallskip\noindent\textit{必要性.}\quad 
核心的群赋值保证 \(\pi(w)\in H(P)\) 必要，证明与前文相同。

对下标为 \(e\) 的情形，若非线性左词含 \(x\notin N(P)\)，令 \(x=d\)，其余变量全为 \(e\)。左侧为 \(c\)，右侧至多为 \(d\)，而 \(c\nleq d\)。

对下标为 \(g\) 的情形，若 \(x\in c(w)\setminus U(P)\)，赋 \(x=c\)、其余变量为 \(e\)，左右分别为 \(c, e\)。若 \(x\in c(w)\setminus N(P)\)、\(x\in U(P)\) 而 \(0\notin H(P)\)，删去右侧一次项 \(x\) 后，存在群赋值使其余右词全取 \(g\)。再令 \(x=d\)，左右分别为 \(c, d\)，吸收失败。

对 \(F_e, F_g\) 的一次词，若没有右侧一次项 \(x\)，令 \(x=t\)、其余变量全为 \(e\)。右侧全在 \(K\) 内，故不大于 \(c<t\)。

对 \(G_e, G_g\)，设 \(x\) 不在右侧一次项中。若 \(\pi(x)\notin H(P)\)，已有群赋值分离。若 \(\pi(x)\in H(P)\) 但 \(0\notin H(P)\)，存在群赋值使全部右词取 \(g\)；仿射条件迫使 \(x=g\)。将 \(x\) 改赋为 \(b\)，所有含 \(x\) 的右词均为非线性词，其值保持 \(g\)，所以右侧仍为 \(g<b\)。这证明零向量条件也必要。

\par\smallskip\noindent\textit{3. 从有限基推导所有吸收式.}\quad 
前文的非线性规范形、偶次变量加入和奇数乘积引理，在左词和各因子均非线性时只使用 \(\Delta\) 的前四条，不使用 \(x\leq x^3\)。因此这里仍有：若非空非线性词族 \(Q\) 满足

\CaseDisplay{c(w)\subseteq\bigcup_{q\in Q}c(q),\qquad
\pi(w)\in\operatorname{aff}\{\pi(q):q\in Q\},}

则对非线性 \(w\) 可推出 \(w\leq\sum Q\)。具体构造是选取给出奇偶向量的奇数个右词，再将全部右词各补两次以覆盖内容；奇数乘积被这些右词之和吸收，而偶次变量加入使左词被该乘积吸收。

对于 \(P\) 中每个同时出现在非线性右词里的单项 \(x\)，取一个含 \(x\) 的非线性词 \(xu\)，由 \(x^3\leq x+xu\) 向 \(P\) 加入 \(x^3\)。保留原一次项，不作未经允许的删除。令 \(Q\) 为原非线性右词和全部这样加入的立方。\(Q\) 的内容正是 \(N(P)\)。其余一次项变量只出现在自己的单项中。

若非线性左词满足 \(E\)，它在仿射包中的奇数项表示不能奇数次使用这些外部一次项；其余一次项可用已加入的立方替换，奇偶向量不变。故 \(\pi(w)\in H(Q)\)。上述引理推出 \(w\leq\sum Q\leq P\)。这给出两个下标为 \(e\) 的半环的全部非线性吸收式。

对于下标为 \(g\) 的两例，若 \(0\notin H(P)\)，同上即可。若 \(0\in H(P)\)，上述 \(Q\) 非空且 \(0\in H(Q)\)。取其零向量的奇数项表示，补入每个 \(Q\) 中的词两次，得到非空平方词 \(h^2\leq P\)。由 (a)，向 \(P\) 加入所有尚未加入的单项立方：

\CaseDisplay{x^3\leq x+h^2\leq P.}

所得非线性族的内容是 \(U(P)\)，仿射包等于 \(H(P)\)，遂推出所有满足 \(J\) 的非线性吸收式。

最后处理 \(G_e, G_g\) 的一次左词。若右侧已有 \(x\)，无需推导。否则假设 \(\pi(x),0\in H(P)\)。由于 \(x\) 不作为单项出现，\(\pi(x)\in H(P)\) 强制 \(x\in N(P)\)。因此已证明的非线性判据同时给出

\CaseDisplay{x^2\leq P,\qquad x^3\leq P.}

由 (b) 得 \(x\leq x^2+x^3\leq P\)。\(F_e, F_g\) 的一次吸收式则全是显然的单项包含。

每个多项式恒等式由两个方向的逐词吸收式组成，故四组基均完整。

\par\smallskip\noindent\textit{4. 范围.}\quad

\subsection{二阶群高层零扩张：七变量有限基定理}\label{case:59}
\par\smallskip\noindent\textit{1. 两个生成半环与基.}\quad 
令 \(\{e, g, c, t\}\) 的加法序为 \(e, g<c<t\)。其中 \(\{e, g\}\) 是二阶群，\(c\) 是该群平坦扩张的零，\(t^2=t\)，并且 \(ta=c\) 对每个 \(a\in\{e, g, c\}\) 成立。添新元素 \(d\)，令 \(\rho(d)=t\)，\(\rho\) 固定其余元素，定义 \(xy=\rho(x)\rho(y)\)。加法额外关系分别是

\CaseDisplay{T_e:\ e<d<t,\quad g\parallel d,\ c\parallel d;
\qquad
T_g:\ g<d<t,\quad e\parallel d,\ c\parallel d.}

乘法和加法满足分配律；可直接由四元素核心上的平移保持并及 \(\rho\) 保持并验证。\(T_e, T_g\) 分别同构于目录编号 1525、1725；在顺序 \((c, e, g, d, t)\) 下，前者的同构映射是 \((0,\allowbreak{}1,\allowbreak{}2,\allowbreak{}3,\allowbreak{}4)\)，后者是 \((0,\allowbreak{}2,\allowbreak{}1,\allowbreak{}3,\allowbreak{}4)\)。

设 \(\Delta\) 包含五条式子

\CaseDisplay{xy\approx yx,\qquad x^3y\approx xy,\qquad x\leq x^3,
\qquad xyz\leq x^3+y^3+z^3,
\qquad x^3\leq x+xy.
\quad\textup{(\ensuremath{1})}}

对于两个非线性词 \(u, q\)，若 \(c(q)\subseteq c(u)\)，记

\CaseDisplay{\mathsf E(u, q, z):\quad u\leq q+uz^2.\quad\textup{(\ensuremath{2})}}

令 \(\mathcal E_7\) 是所有 (2) 中变量数不超过 7、且 \(u, q\) 已按 \(x^3y\approx xy\) 约化的实例。可将变量规范命名为 \(x_1,\ldots, x_6, z\)，只保留每个含至少两个变量的词中指数为 1 或 2、每个单变量非线性词为平方或立方的情形。这是一个明确的有限集合：所有候选词均从有限的指数向量表构造，并按 \(c(q)\subseteq c(u)\) 筛选。

定理 1。 \(T_e\) 由 \(\Delta\cup\mathcal E_7\) 完整公理化。\(T_g\) 由 \(\Delta\cup\mathcal E_7\) 再加

\CaseDisplay{x^3\leq x+y^2 \quad\textup{(\ensuremath{3})}}

完整公理化。因此两者都有限基。

式 (1) 在两个生成元中逐值有效：所有非线性词落在 \(\{e, g, c, t\}\) 内并满足 \(v^3=v\)；当三个值均在群中时，不同群值的和为 \(c\)，当出现 \(d, t\) 时右侧的立方为 \(t\)。式 (3) 仅需检查 \(x=d\)：对任意 \(y\)，\(y^2\in\{e, c, t\}\)，在 \(T_g\) 中 \(d+y^2=t\)。

\par\smallskip\noindent\textit{2. 奇偶、内容与高层锚点.}\quad 
对非空词 \(w\)，令 \(c(w)\) 为变量内容，\(\pi(w)\in\mathbb F_2^X\) 为指数奇偶向量。对非空多项式 \(P=\sum_i v_i\)，置

\CaseDisplay{U(P)=\bigcup_i c(v_i),\quad
N(P)=\bigcup_{|v_i|\geq2}c(v_i),\quad
L(P)=\{x:x\text{ 是 }P\text{ 的一次项}\},\quad
H(P)=\operatorname{aff}\{\pi(v_i):i\}.}

对非线性左词 \(w\)，称 \(P\) 有一个 内部锚点，若存在非线性右词 \(q\) 使 \(c(q)\subseteq c(w)\)，或 \(L(P)\cap c(w)\ne\varnothing\)。

命题 2（任意变量数的精确判据）。

1. 一次左词 \(x\leq P\) 在 \(T_e, T_g\) 中有效，当且仅当 \(\pi(x)\in H(P)\)，并且 \(x\in L(P)\) 或 \(P\) 含一个非线性纯 \(x\) 词。

2. 非线性 \(w\leq P\) 在 \(T_e\) 中有效，当且仅当 \(\pi(w)\in H(P)\)、\(c(w)\subseteq N(P)\)，且 \(P\) 有内部锚点。

3. 非线性 \(w\leq P\) 在 \(T_g\) 中有效，当且仅当 \(\pi(w)\in H(P)\)、\(c(w)\subseteq U(P)\)、\(c(w)\subseteq N(P)\) 或 \(0\in H(P)\)，且 \(P\) 有内部锚点。这里``或''只连接两个后置条件，亦即前两个条件与锚点条件都必须成立。

证明。 将所有变量赋于 \(\{e, g\}\)，平坦加法给出仿射条件：若左奇偶向量不在右向量仿射包内，取一个线性泛函使全部右词取同一群值而左词取另一个。

一次项条件的第二部分由赋值 \(x=t\)、其余变量为 \(c\) 得到：没有纯 \(x\) 非线性项且没有 \(x\) 一次项时，所有右词不高于 \(c<t\)。反过来，纯 \(x\) 的非线性项在 \(x=d, t\) 时均取 \(t\)，在 \(x=c\) 时取 \(c\)；若 \(x\) 在群中，仿射条件保证右侧至少为 \(x\)。

对非线性词，若 \(T_e\) 中 \(x\in c(w)\setminus N(P)\)，赋 \(x=d\)、其余变量为 \(e\)。左词为 \(t\) 或 \(c\)，右侧不超过 \(d\)，所以吸收失败。在 \(T_g\) 中，若 \(x\in c(w)\setminus U(P)\)，赋 \(x=c\)、其余变量为 \(e\)；若 \(x\in c(w)\setminus N(P)\)、\(x\in U(P)\) 而 \(0\notin H(P)\)，从右侧删去唯一可能含 \(x\) 的单项，选择使其余右词全取 \(g\) 的群赋值，再把 \(x\) 改为 \(d\)。此时左词为 \(t\) 或 \(c\)，右侧为 \(d\)，仍失败。若没有内部锚点，赋 \(c(w)\) 中所有变量为 \(d\)，外部变量为 \(e\)。左侧为 \(t\)，右侧不超过 \(c\)。

第 3 节先独立证明消元式的有效性与有限变量约化；第 4 节再从这些有效式子推出满足上述条件的每个吸收式。这一等式推导直接给出充分性，同时完成有限基的完备性，不需要将有限赋值核验外推为任意变量数结论。

\par\smallskip\noindent\textit{3. 七变量消元式的有效性与有限性.}\quad 
引理 3。 每个 (2) 在 \(T_e, T_g\) 中有效。

证明。 若 \(u\) 取群值，则 \(uz^2\) 等于 \(u\)，或取在 \(u\) 上方的 \(c, t\)。若 \(u=c\)，则 \(uz^2=c\)。若 \(u=t\)，其全部内容变量取 \(d, t\)，故内容包含于 \(u\) 的非线性词 \(q\) 也取 \(t\)。三种情形给出 \(u\leq q+uz^2\)。证毕。

引理 4。 从有限集合 \(\mathcal E_7\) 与 (1) 可推出任意变量数的 (2)。

证明。 先用 \(x^3y\approx xy\) 规范化 \(u, q\)。若二者各含至少两个不同变量，则每个 \(u\) 变量在 \(u, q\) 中的指数类型属于

\CaseDisplay{(1,\allowbreak{}0),(2,\allowbreak{}0),(1,\allowbreak{}1),(1,\allowbreak{}2),(2,\allowbreak{}1),(2,\allowbreak{}2).}

每类合并为一个新变量；将它代换为该类原变量的乘积，就精确恢复 \(u, q\)。这至多使用六个新变量，加 \(z\) 共七个。若合并后某个词偶然变成一次词，则其所有原变量都位于同一指数为 1 的类型；将此类型拆成两个非空部分，即恢复非线性，而此时原来使用的类型数最多为三。若 \(q\) 为单变量的平方或立方，将它的变量单列，其余 \(u\) 变量至多再按指数 1、2 分成两类；若 \(u\) 亦为单变量，则只需原变量和 \(z\)。所有这些压缩实例均在 \(\mathcal E_7\) 内。最后对恢复的词作规范约化，即得原 (2)。证毕。

\par\smallskip\noindent\textit{4. 从有限基到任意有效恒等式.}\quad 
引理 5。 若非线性词 \(q_1,\ldots, q_{2r+1}\) 均被 \(P\) 吸收，则其乘积也被 \(P\) 吸收。

证明。 非线性词均满足 \(q_i^3\approx q_i\)。对 \(r\) 归纳，应用 (1) 的 \(xyz\leq x^3+y^3+z^3\)。证毕。

设 \(w\leq P\) 满足命题 2，且 \(w\) 非线性。对每个既是单项又出现在某非线性右词中的 \(x\)，把该词写为 \(xu\)，由 \(x^3\leq x+xu\) 向 \(P\) 加入 \(x^3\)。令 \(Q\) 为所有非线性右词及新加的立方。其内容为 \(N(P)\)，其余单项变量完全没有出现在非线性右词中。

在 \(T_e\)，或在 \(T_g\) 且 \(0\notin H(P)\) 时，判据给出 \(c(w)\subseteq N(P)\)。由于外部独立单项变量在 \(\pi(w)\) 的坐标均为零，删去仿射表达中的成对外部单项，再将内部单项改用相同奇偶向量的立方，得到 \(\pi(w)\in H(Q)\)。

在 \(T_g\) 且 \(0\in H(P)\) 时，外部独立单项变量同样不能奇数次出现在零向量表达中，故 \(0\in H(Q)\)。由引理 5，\(Q\) 中奇数个词的乘积给出一个非空平方词 \(h^2\leq P\)。用 (3) 向 \(P\) 加入其余单项的立方。如此得到的非线性词族 \(Q^*\) 的内容为 \(U(P)\)，仿射包仍为 \(H(P)\)。以下统一以 \(Q\) 表示完成这一步后的非线性族。

由命题 2 的内部锚点，\(Q\) 中有某个词 \(q\) 满足 \(c(q)\subseteq c(w)\)：它或者本来就是内部非线性右词，或者是内部单项新加入的立方。选取给出 \(\pi(w)\) 的奇数个 \(Q\) 词，再将每个 \(Q\) 词各加入两次。所得奇数乘积 \(V\) 满足 \(V\leq P\)、\(c(w)\subseteq c(V)\) 和 \(\pi(V)=\pi(w)\)。将 \(E=c(V)\setminus c(w)\) 中变量逐一列出，由 (1) 的规范形有

\CaseDisplay{V\approx w\prod_{z\in E}z^2.}

反复应用引理 4 的 \(u\leq q+uz^2\)，从 \(V\leq P\)、\(q\leq P\) 依次消去 \(E\) 的变量，最终得到 \(w\leq P\)。

若一次左词 \(x\) 已是 \(P\) 的单项，吸收显然。否则命题 2 保证 \(P\) 有一个纯 \(x\) 非线性项、\(\pi(x)\in H(P)\)。对非线性左词 \(x^3\) 重复上一段，内部锚点即该纯词，得到 \(x^3\leq P\)。再用 (1) 的 \(x\leq x^3\) 得 \(x\leq P\)。全部有效逐词吸收式因此从所列有限基推出；双向逐词吸收给出全部多项式恒等式，定理 1 得证。

\par\smallskip\noindent\textit{5. 目录传递与限制.}\quad

\subsection{A sixteen-identity basis when the nilpotent element lies above the nonidentity}\label{case:60}
All identities are in the constant-free ai-semiring signature and are relative to its axioms.

\par\smallskip\noindent\textit{1. The algebra and statement.}\quad 
Use the multiplication and left-zero-top extension for source \(1787\) from the preceding source construction, but change the core additive order to

\CaseDisplay{g<n<0,\qquad e<0,}

where \(G=\{e, g\}\cong C_2\) and \(e\) is incomparable with \(g, n\). As before \(nx=xn=0\) for every core element, and the additional top \(t\) satisfies \(tx=t\), \(xt=0\) for core \(x\). Call this semiring \(B\). Its catalogue source is 1302 under \((0, e, g, n, t)=(0,\allowbreak{}1,\allowbreak{}3,\allowbreak{}2,\allowbreak{}4)\). Associativity and distributivity follow by the same translation checks, now using \(g<n\).

Theorem. A complete basis for \(B\) consists of the fifteen identities (M), (E), (O), (C), (H) displayed in the source-1787 proof, together with

\CaseDisplay{x^2\leq x+y^2.\quad\textup{(\ensuremath{S})}}

The opposite has the reversed basis.

The soundness proofs for the first fifteen identities do not use the specific comparison \(e<n\): whenever a nonlinear word contains a variable assigned \(n\), its core value is zero. Thus those proofs apply here. For (S), if \(y^2=0\) the core inequality is immediate. Otherwise \(y^2=e\). For \(x=e\) both sides are \(e\); for \(x=g\) their right side is zero and the left side \(e\); for \(x=n\) or zero the right side is zero, as is the left side. The head condition proves it also in the top extension. All sixteen laws are therefore valid.

\par\smallskip\noindent\textit{2. Exact criterion.}\quad 
For a nonempty polynomial \(P\), let \(L\) be the linear variables, \(N\) its nonlinear words, \(U=\bigcup_{v\in N}c(v)\), and \(C_P=U\cup L\). Let \(H(P)\) be its heads and let \(\mathcal A(P)\) be the affine hull over \(\mathbb F_2\) of its parity vectors. Word normal forms and the head lemma remain valid.

A linear target \(x\) is absorbed if and only if \(x\in H(P)\) and \(\pi(x)\in\mathcal A(P)\). A nonlinear target \(w\) is absorbed if and only if

\CaseDisplay{h(w)\in H(P),\qquad c(w)\subseteq C_P,\qquad \pi(w)\in\mathcal A(P),}

and

\CaseDisplay{c(w)\subseteq U\quad\hbox{or}\quad 0\in\mathcal A(P).\quad\textup{(\ensuremath{1})}}

Here zero denotes the parity vector, not a constant in the signature.

The head and affine conditions are necessary by the same evaluations as before. An absent target variable can be sent to the multiplicative zero while every right variable is sent to \(e\), proving content inclusion. If \(x\in c(w)\setminus U\), this inclusion makes \(x\) a linear-only right variable. Suppose \(0\notin\mathcal A(P)\). Linear algebra supplies a linear functional with value one on every right parity vector: equivalently, the equations \(\pi(v)\cdot z=1\) for all right summands are consistent precisely when there is no odd dependency summing to zero. Assigning group elements according to this functional makes every right summand \(g\). The linear variable \(x\) is assigned \(g\); change it to \(n\). Since no nonlinear right word contains \(x\), the right polynomial becomes \(n\), while the nonlinear target becomes zero. This violates absorption and proves (1).

For completeness of the linear-algebra statement, write \(\mathcal A(P)=v_0+D\). If \(0\notin v_0+D\), then \(v_0\notin D\); extend a basis of \(D\) by \(v_0\) and choose a functional vanishing on \(D\) and taking \(v_0\) to one. Conversely such a functional excludes zero from the affine hull. Sufficiency of the criterion follows from the derivations below.

\par\smallskip\noindent\textit{3. Derivation with covered content.}\quad 
First let \(w\) be nonlinear with \(c(w)\subseteq U\), its head present, and parity in \(\mathcal A(P)\). Every linear-only variable has its own unique coordinate among the right summands, so its coefficient in an odd representation of \(\pi(w)\) is zero. Hence the target parity belongs to the affine hull of \(N\) and the linear variables in \(L\cap U\).

Lift each of these covered linear variables \(x\) to \(x^3\leq P\) using an original nonlinear word containing it and the appropriate one of (C). This gives nonlinear words below \(P\) with union of contents \(U\) and the necessary affine hull. As in the source-1787 proof, choose an original nonlinear word containing the target head and an original summand beginning with that head. The six schemes (H) give versions of all these nonlinear words with the desired head and unchanged parity. Form an odd product realizing the target parity and add pairs of every word to cover \(U\). Square extension (E), normalization (M), and odd-product absorption (O) derive \(w\leq P\). This argument uses no comparison between the two group elements and \(n\), so it applies unchanged.

\par\smallskip\noindent\textit{4. Derivation using the zero parity vector.}\quad 
Suppose \(0\in\mathcal A(P)\). Put \(D=L\setminus U\). Each variable in \(D\) occurs in only its own linear right summand. Its coefficient in an odd representation of zero must therefore be zero. It follows that zero is in the affine hull of \(N\) and \(L\cap U\); in particular \(N\) is nonempty.

Choose an original nonlinear word with head \(r\in U\). Lift the covered linear summands by (C), and transfer all resulting nonlinear words to head \(r\) by (H). An odd product realizing parity zero, supplemented with pairs of every factor, yields a nonlinear word \(V\leq P\) of head \(r\), content \(U\), and even exponent at every variable. By (M), \(V\) equals \(s^2\), where \(s\) is any word listing every variable in \(U\) once with \(r\) first. Thus

\CaseDisplay{s^2\leq P.}

For every \(x\in D\), (S) gives \(x^2\leq x+s^2\leq P\). The first rule of (C), with \(u=x\), gives \(x^3\leq x+x^2\leq P\). Hence every original linear word now has a nonlinear cube below \(P\): covered linear words were already lifted, and previously uncovered ones have just been lifted.

Let \(Q\) comprise the original nonlinear words, the cubes of all original linear variables, and the newly derived squares \(x^2\) for \(x\in D\). Every member of \(Q\) is below \(P\), its union of contents is \(C_P\), and its affine hull contains the affine hull of the original right summands. The zero vector also belongs to that hull by construction. Every original head still occurs as a head of a member of \(Q\). The covered-content derivation in the preceding section can therefore be applied to \(Q\) and any nonlinear target satisfying the criterion: all target variables are now covered. It gives \(w\leq Q\leq P\). All uses of derived summands are ordinary transitivity of the semilattice order.

Finally, for a linear target \(x\) not already in \(L\), the head condition supplies a nonlinear right word beginning with \(x\), so \(x\in U\). The covered-content argument applies to \(x^3\), which has the same parity. Identity (E) with \(y=x\) gives \(x\leq x^3\leq P\).

We have derived every absorption in the exact criterion for any finite variable set. Distribution and the two families of word absorptions for a polynomial identity prove the sixteen-identity basis. The proof does not assume any bounded number of variables, path lengths, or summands.

\par\smallskip\noindent\textit{5. Transfer and scope.}\quad 
Independent certificates must check the sixteen displayed identities, the source tables, and complete two-way HSP membership before adding catalogue types to the finite-basis list. This theorem extends the head-affine method to the other possible group marking, but makes no general claim about arbitrary one-sided top extensions.

\subsection{内容补全与首字母保护：六组完整有限基}\label{case:61}
\par\smallskip\noindent\textit{1. 六个生成元及其词值.}\quad 
各底集均记为 \(\{z, b, a, d, t\}\)。加法取下表指定偏序的并；没有列出的比较关系不成立，只补上传递关系。\(z\) 仅是元素名称，在 \(B, D\) 中并非加法最小元。

\noindent 半环\enspace\textbar\enspace 加法偏序\enspace\textbar\enspace 目录编号\enspace\textbar\enspace \((z, b, a, d, t)\) 的目录标号\par
\noindent \(A\)\enspace\textbar\enspace \(z<d<b<a<t\)\enspace\textbar\enspace 3049\enspace\textbar\enspace \((3,\allowbreak{}1,\allowbreak{}0,\allowbreak{}2,\allowbreak{}4)\)\par
\noindent \(B\)\enspace\textbar\enspace \(b<a<t,\ z<a<t,\ b<d<t\)\enspace\textbar\enspace 4157\enspace\textbar\enspace \((2,\allowbreak{}1,\allowbreak{}0,\allowbreak{}3,\allowbreak{}4)\)\par
\noindent \(C\)\enspace\textbar\enspace \(z<b<a<d<t\)\enspace\textbar\enspace 11794\enspace\textbar\enspace \((3,\allowbreak{}2,\allowbreak{}1,\allowbreak{}0,\allowbreak{}4)\)\par
\noindent \(D\)\enspace\textbar\enspace \(z<d<a<t,\ b<t\)\enspace\textbar\enspace 4109\enspace\textbar\enspace \((2,\allowbreak{}3,\allowbreak{}0,\allowbreak{}1,\allowbreak{}4)\)\par
\noindent \(E\)\enspace\textbar\enspace \(z<b<a<t,\ b<d<t\)\enspace\textbar\enspace 4159\enspace\textbar\enspace \((2,\allowbreak{}1,\allowbreak{}0,\allowbreak{}3,\allowbreak{}4)\)\par
\noindent \(F\)\enspace\textbar\enspace \(z<b<a<t,\ b<d<t\)\enspace\textbar\enspace 6277\enspace\textbar\enspace \((2,\allowbreak{}1,\allowbreak{}0,\allowbreak{}3,\allowbreak{}4)\)\par
乘法如下构造。对于 \(A, B, E\)，在核心 \(\{z, b, a, t\}\) 中，\(z\) 为乘法零，\(a, t\) 在非零元素上分别作左零元素，且 \(b^2=b\)、\(ba=bt=a\)。通过回缩扩张：\(A\) 取 \(\rho(d)=b\)；\(B, E\) 取 \(\rho(d)=t\)，其他元素固定，并令 \(xy=\rho(x)\rho(y)\)。

对于 \(C\)，核心仍为 \(\{z, b, a, t\}\)，但 \(b, a\) 为左零元素，\(zz=z\)、\(zx=b\)（\(x\ne z\)），\(tt=t\)、\(tx=a\)（\(x\ne t\)）；取 \(\rho(d)=t\)。

对于 \(D\)，核心为 \(\{z, b, a, t\}\)，\(z\) 为乘法零，\(a, t\) 在非零元素上作左零元素，\(b^2=b\)、\(ba=bt=t\)；取 \(\rho(d)=a\)。

对于 \(F\)，核心为 \(\{z, b, a\}\)，\(z, b\) 为左零元素，\(a^2=a\)、\(az=ab=b\)；取 \(\rho(d)=\rho(t)=a\)。

上述核心均为带，每个回缩保持加法。核心乘法的左右平移保持表列偏序中的并，故分配律由回缩继承；结合律由 \(\rho\) 固定乘法核心而继承。也可用以下词值描述直接验证三个字母的两种括号化相等。

设非一次词 \(w\) 的首字母为 \(h\)、内容为 \(X\)：

\textbullet{} \(A\)：出现 \(z\) 则取 \(z\)；否则首字母为 \(t, a\) 时分别取 \(t, a\)；首字母为 \(b, d\) 时，若所有字母在 \(\{b, d\}\) 则取 \(b\)，否则取 \(a\)。

\textbullet{} \(B, E\)：出现 \(z\) 则取 \(z\)；否则首字母为 \(d, t\) 时取 \(t\)，为 \(a\) 时取 \(a\)，为 \(b\) 时全体字母均为 \(b\) 则取 \(b\)，否则取 \(a\)。

\textbullet{} \(C\)：首字母为 \(b, a\) 时分别取 \(b, a\)；首字母为 \(z\) 时，全体字母为 \(z\) 则取 \(z\)，否则取 \(b\)；首字母为 \(d, t\) 时，全体字母属于 \(\{d, t\}\) 则取 \(t\)，否则取 \(a\)。

\textbullet{} \(D\)：出现 \(z\) 则取 \(z\)；否则首字母为 \(a, d\) 时取 \(a\)，为 \(t\) 时取 \(t\)；首字母为 \(b\) 时，全体字母为 \(b\) 则取 \(b\)，否则取 \(t\)。

\textbullet{} \(F\)：首字母为 \(z, b\) 时分别取 \(z, b\)；首字母属于 \(\{a, d, t\}\) 时，全体字母也属于该集合则取 \(a\)，否则取 \(b\)。

\par\smallskip\noindent\textit{2. 完整吸收判据.}\quad 
对非空多项式 \(p\)，令 \(L\) 为一次项变量集，\(\mathcal R\) 为非一次项的词族，并定义

\CaseDisplay{T=c(p),\quad N=\bigcup_{u\in\mathcal R}c(u),\quad H=L\cup\{h(u):u\in\mathcal R\}.}

对非空集合 \(X\)，置

\CaseDisplay{\mathcal R_X=\{u\in\mathcal R:c(u)\subseteq X\},\quad
N_X=\bigcup_{u\in\mathcal R_X}c(u),\quad
H_X=L\cup\{h(u):u\in\mathcal R_X\}.}

这里 \(h\in H_X\) 中的一次项分支只是 \(h\in L\)；由于 \(h\in X\)，它本身就是内部一次项。

判据。 对非一次词 \(w\)，令 \(h=h(w),X=c(w)\)。\(w\le p\) 的充分必要条件如下。

\noindent 来源\enspace\textbar\enspace 条件\par
\noindent \(A\)\enspace\textbar\enspace \(\mathcal R_X\ne\varnothing\)，\(X\subseteq L\cup N_X\)，\(h\in H_X\)\par
\noindent \(B\)\enspace\textbar\enspace \(X\subseteq N\)，\(h\in H_X\)\par
\noindent \(C\)\enspace\textbar\enspace \(\mathcal R_X\ne\varnothing\)，\(X\subseteq T\)，\(h\in H\)\par
\noindent \(D\)\enspace\textbar\enspace \(\mathcal R_X\ne\varnothing\)，\(X\subseteq T\)，\(h\in H_X\)\par
\noindent \(E\)\enspace\textbar\enspace \(X\subseteq N_X\)，\(h\in H_X\)\par
\noindent \(F\)\enspace\textbar\enspace \(\mathcal R_X\ne\varnothing\)，\(h\in H\)\par
对一次词 \(x\)，在 \(A, B, C, D, E\) 中，\(x\le p\) 当且仅当 \(x\in L\) 或存在仅含 \(x\) 的非一次右词；在 \(F\) 中，当且仅当 \(x\in L\)。

\par\smallskip\noindent\textit{2.1 必要性：明确分离赋值.}\quad 
以下每个赋值仅在对应条件失败时使用；所有未特别指出的外部变量均按说明赋值。

\(A\)。 若无内部非一次右词，令 \(X\) 内变量为 \(d\)、外部变量为 \(z\)，左值为 \(b\)，右值至多为 \(d<b\)。若 \(x\in X\setminus(L\cup N_X)\)，令 \(x=a\)、其他 \(X\) 内变量为 \(d\)、外部变量为 \(z\)，左值为 \(a\)，右值至多为 \(b\)。若 \(h\notin H_X\)，令 \(h=t\)、其他 \(X\) 内变量为 \(d\)、外部为 \(z\)，左值为 \(t\)，右值至多为 \(a\)。

\(B\)。 若 \(x\in X\setminus N\)，令 \(x=d\)、其他所有变量为 \(b\)。左值为 \(a\) 或 \(t\)，右非一次项全为 \(b\)，右一次项至多为 \(d\)；而 \(a\nleq d\)、\(t\nleq d\)。若 \(h\notin H_X\)，令 \(h=t\)、其他 \(X\) 内变量为 \(b\)、外部变量为 \(z\)，左值为 \(t\)，右值至多为 \(a\)。

\(C\)。 无内部非一次项时，赋 \(X=d\)、外部为 \(z\)，左值为 \(t\)，右值至多为 \(d\)。若 \(x\in X\setminus T\)，赋 \(x=b\)、其他全为 \(z\)，左值为 \(b\)、右值为 \(z\)。若 \(h\notin H\)，赋 \(h=a\)、其他全为 \(z\)，左值为 \(a\)、右值至多为 \(b\)。

\(D\)。 无内部非一次项时，赋 \(X=d\)、外部为 \(z\)，左值为 \(a\)、右值至多为 \(d\)。若 \(x\in X\setminus T\)，赋 \(x=z\)、其他所有变量为 \(b\)，左值为 \(z\)、右值为 \(b\)，且 \(z\nleq b\)。若 \(h\notin H_X\)，赋 \(h=t\)、其他 \(X\) 内变量为 \(d\)、外部为 \(z\)，左值为 \(t\)，右值至多为 \(a\)。

\(E\)。 若 \(x\in X\setminus N_X\)，赋 \(x=d\)、其他 \(X\) 内变量为 \(b\)、外部为 \(z\)。左值为 \(a\) 或 \(t\)，内部非一次项全为 \(b\)、外部非一次项全为 \(z\)、一次项至多为 \(d\)，故失败。若 \(h\notin H_X\)，赋 \(h=t\)、其他 \(X\) 内变量为 \(b\)、外部为 \(z\)，左值为 \(t\)、右值至多为 \(a\)。

\(F\)。 无内部非一次项时，赋 \(X=d\)、外部为 \(z\)，左值为 \(a\)、右值至多为 \(d\)，且 \(a\nleq d\)。若 \(h\notin H\)，赋 \(h=b\)、其他全为 \(z\)，左值为 \(b\)、右值为 \(z\)。

一次词的分离：\(A, B, D, E\) 取 \(x=d\)、其他变量为 \(z\)；没有同名一次项和仅含 \(x\) 的非一次项时，右非一次项全为 \(z\)，且右一次项也为 \(z\)。\(C\) 取同一赋值，右非一次项至多为 \(a<d\)。\(F\) 取 \(x=t\)、其他变量为 \(z\)，所有右非一次项至多为 \(a<t\)。

\par\smallskip\noindent\textit{2.2 充分性.}\quad 
以下讨论始终允许任意赋值，使用第 1 节词值规则。

\(A\)。 左值为 \(z\) 时自明。否则 \(X\) 内没有 \(z\)，内部非一次右词至少为 \(b\)。左值达到 \(a\) 时，\(X\) 内有取 \(a\) 或 \(t\) 的变量；\(X\subseteq L\cup N_X\) 提供同变量一次项，或含它的内部非一次项，值至少为 \(a\)。左值为 \(t\) 时，首字母取 \(t\)，\(h\in H_X\) 给出同名一次项或同首字母内部非一次项，均取 \(t\)。

\(B\)。 若左值为 \(z\)，\(X\subseteq N\) 提供含某个取 \(z\) 变量的右非一次项，值为 \(z\)。若左值为 \(b\)，所有 \(X\) 内变量均为 \(b\)，而 \(h\in H_X\) 提供一个内部右项，取 \(b\)。若左值为 \(a\)，\(X\) 内没有 \(z\)，且有一个变量取 \(a, d\) 或 \(t\)。含此变量的右非一次项由 \(X\subseteq N\) 给出，其值或者为 \(z\)，或者至少为 \(a\)。后一种直接完成；前一种与由 \(h\in H_X\) 给出的内部右项求并，得到至少为 \(a\) 的值，因为内部右项取 \(b, a, d\) 或 \(t\)，而 \(z+b=a\)、\(z+d=t\)。若左值为 \(t\)，首字母取 \(d\) 或 \(t\)；内部同首字母非一次项直接取 \(t\)。若只有同名一次项，当它取 \(t\) 时完成；当它取 \(d\) 时，含 \(h\) 的右非一次项取 \(z\) 或至少为 \(a\)，与 \(d\) 求并均得 \(t\)。

\(C\)。 左值为 \(b\) 时，\(X\) 内存在非 \(z\) 的变量；内容覆盖保证右侧某项含它，从而右值至少为 \(b\)。左值为 \(a\) 时，同首字母右项的值至少为 \(a\)；此处一次项可为 \(d\)，但链上已有 \(a<d\)。左值为 \(t\) 时，所有 \(X\) 内字母在 \(\{d, t\}\)，故内部非一次右词取 \(t\)。

\(D\)。 左值为 \(z\) 时，内容覆盖提供含取 \(z\) 变量的右项，值为 \(z\)。左值为 \(b\) 时，所有 \(X\) 内变量均为 \(b\)，内部非一次右词取 \(b\)。左值为 \(a\) 时，首字母为 \(a\) 或 \(d\)；内部同首字母非一次右词取 \(a\)。若用同名一次项而该项为 \(d\)，任意内部非一次右词取 \(a, b\) 或 \(t\)，与 \(d\) 求并至少为 \(a\)。左值为 \(t\) 时，若首字母为 \(t\)，首字母见证直接给出 \(t\)；否则首字母为 \(b\)，且某个 \(X\) 内变量不为 \(b\)。首字母见证若已为 \(t\) 则完成，否则给出 \(b\)。内容覆盖再提供含该非 \(b\) 变量的右项，其值为 \(z, d, a\) 或 \(t\)；它与 \(b\) 的并都是 \(t\)。

\(E\)。 左值为 \(z\) 时自明；左值为 \(b\) 时，内部非一次项存在并取 \(b\)。左值为 \(a\) 时，含相应高值变量的内部非一次项由 \(X\subseteq N_X\) 给出，值至少为 \(a\)。左值为 \(t\) 时，首字母为 \(d, t\)，同首字母内部非一次项直接取 \(t\)；若用一次项 \(d\)，内部含该变量的非一次项取 \(a\) 或 \(t\)，与 \(d\) 求并得 \(t\)。

\(F\)。 左值为 \(z\) 时自明；为 \(b\) 时，同首字母右项至少为 \(b\)；为 \(a\) 时，所有 \(X\) 内变量属于 \(\{a, d, t\}\)，内部非一次右词取 \(a\)。一次词的充分性由同名项及前五个来源中的 \(x\le x^2=x^n\)（\(n\ge2\)）给出。

\par\smallskip\noindent\textit{3. 六组完整有限基.}\quad 
共同的两式为

\CaseDisplay{xyz\approx xzy\quad(N),\qquad xy\approx x^2y^2\quad(S).}

以下符号仅作有限表的索引：

\CaseDisplay{\begin{aligned}
x&\le x^2 &(I),\\
xy&\le x^2+y^2 &(J),\\
xy&\le x^2+y &(K),\\
xy&\le x+yx &(H),\\
xy&\le x^2+zy &(B_1),\\
x^2&\le x+yx &(Q),\\
xy&\le y^2+xyz &(E_2),\\
xy&\le x^2 &(C_1).
\end{aligned}}

定理。 各来源的完整附加式如下，均再加 \((N),(S)\)：

\noindent 来源\enspace\textbar\enspace 另外的式子\enspace\textbar\enspace 总附加式数\par
\noindent \(A\) / 3049\enspace\textbar\enspace \(I, K, H\)\enspace\textbar\enspace 5\par
\noindent \(B\) / 4157\enspace\textbar\enspace \(I, B_1, Q\)\enspace\textbar\enspace 5\par
\noindent \(C\) / 11794\enspace\textbar\enspace \(I, K, H, E_2\)\enspace\textbar\enspace 6\par
\noindent \(D\) / 4109\enspace\textbar\enspace \(I, K, H, B_1\)\enspace\textbar\enspace 6\par
\noindent \(E\) / 4159\enspace\textbar\enspace \(I, J, H\)\enspace\textbar\enspace 5\par
\noindent \(F\) / 6277\enspace\textbar\enspace \(C_1, H, E_2\)\enspace\textbar\enspace 5\par
有效性由第 2 节判据直接代入验证；\((N),(S)\) 由词值规则给出。乘法对偶的完整基是这些基逐词反转的结果。

\par\smallskip\noindent\textit{3.1 共同形式推导.}\quad 
由 \((S)\) 得 \(x^2=x^4\)、\(x^3=x^6=x^2\)，再结合 \((N)\)，每个非一次词可写为 \(W(h, X)\)：首字母为 \(h\)，内容为 \(X\)，每个字母均出现两次。因此非一次词幂等，且

\CaseDisplay{W(h, X)W(g, Y)=W(h, X\cup Y).}

若有 \((K)\)，非一次词 \(p\) 与任意词 \(q\) 满足 \(pq\le p+q\)。若有 \((J)\)，两个非一次词也满足该不等式。因此从至少一个非一次右项开始，可在 \((K)\) 下乘上其他任意右项，在 \((J)\) 下乘上其他非一次右项，所得积都不超过右多项式。

\((B_1)\) 提供一种不引入外部内容的补全。若非一次词 \(p\le F\)，右非一次词 \(u\le F\) 包含变量 \(x\)，代入 \((B_1)\) 得

\CaseDisplay{px\le p^2+ux=p+u\le F,\quad\textup{(\ensuremath{3.1})}}

因为 \(ux=u\)。这里右词 \(u\) 可以包含任意外部变量，而左边只新增 \(x\)。

\((Q)\) 提供首字母保护：若 \(h\in L\) 且右非一次词 \(u\) 含 \(h\)，则

\CaseDisplay{h^2\le h+uh=h+u\le F.\quad\textup{(\ensuremath{3.2})}}

若同时有 \((H),(K)\)，\(h\in L\) 且 \(q\) 为右非一次项，则

\CaseDisplay{hq\le h+qh\le h+q\le F.\quad\textup{(\ensuremath{3.3})}}

以下所有缺少字母的补乘都仅在缺集非空时进行，否则直接省略；没有空词或常量代入。

\par\smallskip\noindent\textit{3.2 分来源的完备性.}\quad 
\(A\)。 只使用内部右项。若有内部同首字母非一次项，就以它开始，乘上所有内部非一次项和 \(L\cap X\) 中的变量，用 \((K)\) 归纳。内容覆盖保证积的内容正好为 \(X\)，首字母仍为 \(h\)，故积为 \(w\)。若首字母见证是 \(h\in L\)，先用任意内部非一次项和 (3.3) 得到以 \(h\) 开始的内部积，再作同样补全。

\(B\)。 若首字母见证是内部非一次项，以它开始；若只有一次项 \(h\)，由 \(X\subseteq N\) 选含 \(h\) 的右非一次词，并用 (3.2) 从 \(h^2\) 开始。对每个尚缺的 \(x\in X\)，选包含它的右非一次项，以 (3.1) 补上 \(x\)。最终得到 \(W(h, X)\le F\)。

\(C\)。 若有同首字母右非一次项，从它开始，用 \((K)\) 乘上所有其余右项，得到 \(W(h, T)\le F\)。若首字母只见于一次项，先把全部非一次右项和全部一次右项乘成一个非一次词 \(P\)，由 \((K)\) 有 \(P\le F\) 且 \(c(P)=T\)；再用 \((H)\) 及 \(Ph=P\) 得 \(hP\le h+P\le F\)。选内部非一次右词 \(q\)。若 \(T=X\) 即完成；否则以 \(s\) 表示 \(T\setminus X\) 中各变量组成的非空词，在 \((E_2)\) 中代入 \((x, y, z)=(w, q, s)\)，得到

\CaseDisplay{w=wq\le q^2+wqs=q+W(h, T)\le F.}

\(D\)。 从内部同首字母非一次词开始；若只有同名一次项，先用任意内部非一次项和 (3.3) 建立内部首字母积。每个尚缺的变量由内容覆盖出现在右侧：若它在一次项中，以 \((K)\) 补乘；否则以 (3.1) 补乘。整个过程只增加 \(X\) 内变量，最终得到 \(w\le F\)。

\(E\)。 把全部内部非一次项相乘。由 \((J)\)，积不超过 \(F\)；其内容因 \(X\subseteq N_X\) 正好为 \(X\)。若有同首字母项，把它放在第一位即可。若只有同名一次项 \(h\)，则对所得积 \(P\) 有 \(Ph=P\)，故 \((H)\) 给出 \(hP\le h+P\le F\)，左边就是 \(w\)。

\(F\)。 选内部非一次项 \(q\)。若 \(h\in L\)，由 \((H),(C_1)\) 得 \(hq\le h+qh\le h+q\le F\)，再由 \((C_1)\) 补乘 \(X\) 内缺少的字母。否则有同首字母非一次右词 \(r\)；由 \((C_1)\) 补乘得到 \(W(h, X\cup c(r))\le r\le F\)。用内部非一次锚点 \(q\) 和 \((E_2)\)，同 \(C\) 的证明消去 \(c(r)\setminus X\)，得到 \(w\le F\)。

前五个来源的一次词吸收用 \((I)\) 和正规形处理；\(F\) 的一次词必须已是右项。任意有效恒等式拆为双向逐词吸收，故各基完整，且证明不限制变量数或右项个数。

\par\smallskip\noindent\textit{4. 计算证据的角色.}\quad 
独立程序仅核对具体构造、有限式的有效性、有界诊断和目录对象与来源的双向完整簇包含。任意变量数完备性由上述推导负责。只具有单向包含的对象不传递有限基结论；这里未证明任何未经单独处理的子簇具有有限基。

\subsection{An eight-identity prefix basis for the chain semiring 1471}\label{case:62}
We use the constant-free ai-semiring signature. The carrier is \(\{0,\allowbreak{}1,\allowbreak{}2,\allowbreak{}3,\allowbreak{}4\}\), addition is maximum in the chain \(3<2<1<0<4\), and multiplication has rows

{\small\par\smallskip
\CaseLine{\ttfamily 00000}
\CaseLine{\ttfamily 01130}
\CaseLine{\ttfamily 01130}
\CaseLine{\ttfamily 33333}
\CaseLine{\ttfamily 00000}
\smallskip}

The element labels are not constants in terms. This is catalogue source 1471.

\par\smallskip\noindent\textit{1. Structure and identities.}\quad 
The square map \(h(x)=x^2\) has values \((0,\allowbreak{}1,\allowbreak{}1,\allowbreak{}3,\allowbreak{}0)\) and is an endomorphic retraction onto \(\{0,\allowbreak{}1,\allowbreak{}3\}\). On this image, \(1\) is an identity and \(0,\allowbreak{}3\) are left zeros. Every product equals \(h(x)h(y)\). The image is a semiring on the chain \(3<1<0\), and \(h\) preserves addition. Consequently the table is associative and distributes over addition, as follows directly by replacing each factor with its image under \(h\).

A nonlinear word has value \(1\) when all its variables are in \(\{1,\allowbreak{}2\}\). Otherwise look at its first occurrence outside that set: the word is \(3\) when this occurrence is \(3\), and \(0\) when it is \(0\) or \(4\).

Theorem. The following eight identities are a complete basis, relative to the ai-semiring axioms:

\CaseDisplay{x^2y\approx xy,\quad xy^2\approx xy,\quad xyx\approx xy;\quad\textup{(\ensuremath{W})}}

\CaseDisplay{xyz\leq yz+xt,\quad xyz\leq xz+xy,\quad xyz\leq xz+xyt;\quad\textup{(\ensuremath{P})}}

\CaseDisplay{xyz\leq yz+x,\qquad xyz\leq xz+y.\quad\textup{(\ensuremath{L})}}

Word reversal gives a basis for the opposite semiring.

The word laws retain only first occurrences of variables in a nonlinear word; if there is only one distinct variable, retain its square. Adjacent repeats are deleted by the first two laws and separated repeats by the third. The distinction between linear words and nonlinear singleton-support words is retained.

\par\smallskip\noindent\textit{2. Absorption criterion.}\quad 
Let \(P\) be a nonempty polynomial with linear variable set \(L\) and nonlinear word set \(N\). Let \(\operatorname{pre}_w(x)\) be the set of variables before the first occurrence of \(x\) in \(w\). A linear \(y\) is absorbed exactly when \(y\in L\). A nonlinear \(w\) is absorbed exactly when

1. some \(u\in N\) satisfies \(c(u)\subseteq c(w)\);

2. for every \(x\in c(w)\), either \(x\in L\), or some \(v\in N\) contains \(x\) with \(\operatorname{pre}_v(x)\subseteq\operatorname{pre}_w(x)\).

Proof. Missing linear \(y\) is separated by \(y=4\) and all other variables \(3\): nonlinear words cannot attain \(4\). If condition 1 fails, assign all target variables \(2\) and outsiders \(3\). The target is \(1\), while right linear terms are at most \(2\) and all nonlinear terms are \(3\). If condition 2 fails at \(x\), assign target variables preceding its first occurrence \(2\), give \(x\) value \(0\), and give all others \(3\). The target is \(0\). No linear right word is \(x\), and no nonlinear right word first reaches a value outside \(\{1,\allowbreak{}2\}\) at \(x\), so every right term is at most \(1\). These assignments prove necessity.

Conversely, target value \(3\) is least. If the target is \(1\), its variables are all in \(\{1,\allowbreak{}2\}\) and the contained nonlinear witness has value \(1\). If it is \(0\), choose its first occurrence outside \(\{1,\allowbreak{}2\}\), say \(x\). Then \(x\) has value \(0\) or \(4\). A linear witness \(x\) is at least the target; an ordinary nonlinear prefix witness has value \(0\). This proves sufficiency.

All eight displayed laws are sound by this criterion and the word-value formula. The laws (P) use exactly the same nonlinear prefix logic on the image chain; the first law (L) permits a new first variable supplied by a linear summand, and the second inserts such a variable after a nonempty prefix.

\par\smallskip\noindent\textit{3. Unrestricted derivation.}\quad 
Let \(w\leq P\) satisfy the criterion. Normalize nonlinear words by (W). Write the first-occurrence list of \(w\) as \(x_1\cdots x_d\) and let \(p_i=x_1\cdots x_i\). We prove \(p_i u\leq P\) for every original \(u\in N\). The first condition ensures \(N\) is nonempty.

For \(i=1\), if \(x_1\in L\), split \(u=yz\) with both factors nonempty and use the first rule (L): \(x_1u\leq u+x_1\leq P\). Otherwise there is a nonlinear witness \(x_1t\) with \(t\) nonempty; the first rule (P) gives \(x_1u\leq u+x_1t\leq P\).

Suppose the induction holds for a nonempty prefix \(p\) and the next variable is \(x\). If \(x\in L\), the second rule (L) immediately gives \(pxu\leq pu+x\leq P\). Otherwise choose its nonlinear prefix witness \(v\). The already proved \(pv\leq P\) normalizes to \(pxq\), with \(q\) possibly omitted, because all variables before the first \(x\) in \(v\) already occur in \(p\). Use respectively the second or third rule (P) to derive

\CaseDisplay{pxu\leq pu+px\leq P,\qquad
pxu\leq pu+pxq\leq P.}

No empty word is substituted: these are separate laws for an omitted or nonempty suffix. This establishes the induction.

Choose the nonlinear contained-content witness \(u\) from condition 1. At \(i=d\), the derived word \(p_du\) normalizes to \(w\), including the singleton-support case where both normalize to \(x_1^2\). Thus \(w\leq P\) follows from the eight laws. Linear targets follow directly from membership in \(L\). Expanding both sides of any identity and treating all its summands in both directions proves completeness for arbitrary variable numbers and word lengths. Reversing words and derivations proves the opposite statement.

\subsection{首次出现前缀与一次项插入：六条完整基}\label{case:63}
\par\smallskip\noindent\textit{1. 低层单位克隆构造.}\quad 
取链 \(A=\{z<d<e<a<t\}\)，加法取最大值。核心 \(B=\{z, e, a, t\}\) 的乘法为

\noindent \(\cdot\)\enspace\textbar\enspace \(z\)\enspace\textbar\enspace \(e\)\enspace\textbar\enspace \(a\)\enspace\textbar\enspace \(t\)\par
\noindent \(z\)\enspace\textbar\enspace \(z\)\enspace\textbar\enspace \(z\)\enspace\textbar\enspace \(z\)\enspace\textbar\enspace \(z\)\par
\noindent \(e\)\enspace\textbar\enspace \(z\)\enspace\textbar\enspace \(e\)\enspace\textbar\enspace \(a\)\enspace\textbar\enspace \(a\)\par
\noindent \(a\)\enspace\textbar\enspace \(a\)\enspace\textbar\enspace \(a\)\enspace\textbar\enspace \(a\)\enspace\textbar\enspace \(a\)\par
\noindent \(t\)\enspace\textbar\enspace \(a\)\enspace\textbar\enspace \(a\)\enspace\textbar\enspace \(a\)\enspace\textbar\enspace \(t\)\par
令 \(\rho(d)=e\)，其余元素固定，并取 \(xy=\rho(x)\rho(y)\)。符号顺序 \((z, e, a, d, t)\) 对应目录标号 \((3,\allowbreak{}1,\allowbreak{}0,\allowbreak{}2,\allowbreak{}4)\)。

核心是半环 \(K\times\{0,\allowbreak{}1\}\) 的子半环 \(\{(z,0),(e,0),(a,0),(a,1)\}\)：\(K=\{z<e<a\}\) 加法为最大值，\(e\) 为乘法单位，\(z, a\) 为左零元素；第二坐标加法取最大、乘法取最小。因此核心满足结合律及分配律。回缩 \(\rho\) 保持并并固定核心，故扩张仍是半环。

对非一次词，其值由下列规则决定：若所有字母均取 \(t\)，则值为 \(t\)；否则把 \(d, e\) 视为中性字母，找首个不是 \(d, e\) 的字母：它为 \(z\) 时取 \(z\)，为 \(a, t\) 时取 \(a\)；如果全为 \(d, e\)，则取 \(e\)。

\par\smallskip\noindent\textit{2. 完整吸收判据.}\quad 
对词 \(w\) 记内容为 \(X=c(w)\)，\(P_w(x)\) 为 \(x\) 首次出现前的变量集合。右多项式 \(F\) 的一次项变量集为 \(L\)，非一次右词族为 \(\mathcal R\)。

定理。 对非一次词 \(w\)，\(w\le F\) 成立，当且仅当

1. 每个 \(x\in X\) 或者属于 \(L\)，或者存在 \(u_x\in\mathcal R\)，使 \(x\in c(u_x)\) 且 \(P_{u_x}(x)\subseteq P_w(x)\)；

2. 存在 \(q\in\mathcal R\) 满足 \(c(q)\subseteq X\)。

一次词 \(x\le F\) 当且仅当 \(x\in L\)，或存在仅含 \(x\) 的非一次右词。

必要性。 若条件 1 在 \(x\) 处失败，赋 \(x=a\)，\(P_w(x)\) 内变量为 \(d\)，其他变量为 \(z\)。左值为 \(a\)。右侧没有一次项 \(x\)；含 \(x\) 的非一次右词在首次 \(x\) 之前必遇到一个取 \(z\) 的变量，因此取 \(z\)；其他非一次右词至多取 \(e\)。一次项至多为 \(d\)，所以右值至多 \(e<a\)。若条件 2 失败，赋 \(X\) 内变量为 \(d\)、外部变量为 \(z\)，左值为 \(e\)，右非一次项全为 \(z\)，一次项至多为 \(d<e\)。

对一次词，如无同名一次项或仅含它的非一次项，赋 \(x=d\)、其余变量为 \(z\)，右值为 \(z<d\)。

充分性。 左值为 \(z\) 时自明；为 \(e\) 时 \(X\) 内全为 \(d, e\)，内部非一次项 \(q\) 也取 \(e\)；为 \(t\) 时 \(X\) 内全取 \(t\)，\(q=t\)。左值为 \(a\) 时，取首个非中性变量 \(x\)；它取 \(a\) 或 \(t\)。若有一次项 \(x\)，它直接在 \(a\) 之上；否则前缀见证 \(u_x\) 在首次 \(x\) 之前只有中性字母，故取 \(a\) 或 \(t\)。一次吸收用 \(x\le x^2=x^n\)（\(n\ge2\)）。证毕。

\par\smallskip\noindent\textit{3. 六条完整有限基.}\quad 
主定理。 完整基为

\CaseDisplay{\begin{aligned}
xy&\approx x^2y^2, &(S)\\
xyx&\approx xy, &(R)\\
x&\le x^2, &(I)\\
xy&\le xz+y^2, &(T)\\
xy&\le x+y^2, &(U)\\
xyz&\le y+xz. &(M)
\end{aligned}}

前三式由词值规则验证；\((T),(U),(M)\) 由判据直接验证。尤其 \((M)\) 中右侧一次项 \(y\) 保护左词中间的变量，右侧非一次项 \(xz\) 既提供其余前缀见证，又是内部锚点。

由 \((S),(R)\) 得 \(x^3=x^2\)、\(x^2y=xy=xy^2\)，并将非一次词化为按首次出现顺序排列的 \(W(x_1,\ldots, x_k)=x_1^2\cdots x_k^2\)。所有非一次词幂等；右乘已经包含的变量不改变非一次当前积。

设有效吸收的左词为 \(w=W(x_1,\ldots, x_k)\)，\(p_i=x_1^2\cdots x_i^2\)。证明：对任何已经推出 \(v\le F\) 的非一次词 \(v\)，都有 \(p_iv\le F\)。

若 \(x_i\) 有非一次前缀见证，把正规化的见证写成 \(u_i=a_ix_ib_i\)，其中 \(c(a_i)\subseteq\{x_1,\ldots, x_{i-1}\}\)，\(b_i\) 非空。\((T)\) 在 \(i=1\) 时给出 \(x_1v\le u_1+v\)；在 \(i>1\) 时给出

\CaseDisplay{p_iv\le p_{i-1}x_ib_i+p_{i-1}v=p_{i-1}u_i+p_{i-1}v\le F,}

末步分别对非一次词 \(u_i, v\) 使用归纳假设，前面利用旧前缀吸收 \(a_i\)。

若 \(x_i\in L\)，在 \(i=1\) 时由 \((U)\) 得

\CaseDisplay{p_1v=x_1v\le x_1+v^2=x_1+v\le F;}

在 \(i>1\) 时，将 \((x, y, z)=(p_{i-1},x_i, v)\) 代入 \((M)\)，得到

\CaseDisplay{p_iv=p_{i-1}x_iv\le x_i+p_{i-1}v\le F.}

这样归纳覆盖所有变量。选判据给出的内部非一次项 \(q\)，得到 \(wq\le F\)，而 \(c(q)\subseteq c(w)\) 给出 \(wq=w\)。一次左词由 \((I)\) 和同名右项处理。任意有效恒等式拆为双向逐词吸收，故基完整。整个归纳不限制变量数；空前缀仅通过省略左乘操作处理，没有空词代入。

\par\smallskip\noindent\textit{4. 对偶及分类传递.}\quad 
逐词反转得到对偶基。对具体目录对象，只在独立验证与来源的精确双向簇包含后传递有限基性；单向满足此基不够。本定理不声称所有真子簇有限基。有限诊断及 HSP 核验另行记录，任意变量数完备性由本节给出的形式归纳负责。

\subsection{A fifteen-identity basis for an identity-marked Boolean group with a left-zero top}\label{case:64}
We work in the constant-free signature \((+,\cdot)\) relative to the additively idempotent semiring axioms.

\par\smallskip\noindent\textit{1. Algebra and basis.}\quad 
Take \(G=\{e, g\}\cong C_2\) and put \(F=G\cup\{n,0\}\). The element \(0\) is a multiplicative zero, \(nx=xn=0\) for every \(x\in F\), and multiplication on \(G\) is its group multiplication. The additive order is \(e<n<0\), \(g<0\), with \(g\) incomparable with \(e, n\). Add an additive top \(t\), put \(tx=t\) for every \(x\), and \(xt=0\) for \(x\in F\). The resulting semiring is \(A\). In catalogue 1787, \((0, e, g, n, t)=(0,\allowbreak{}3,\allowbreak{}1,\allowbreak{}2,\allowbreak{}4)\).

The core multiplication is associative, since \(\{n,0\}\) is an ideal whose every product with the core is zero. Its translations by a group element permute the two group elements and send \(n,0\) to zero; they preserve joins, including the only extra comparison \(e<n\). Translations by \(n,0\) are constant zero. The left-zero-top extension preserves associativity and both distributive laws by direct examination of the first occurrence of \(t\).

Here are nine individual identities:

\CaseDisplay{xyz\approx xzy,\qquad xy\approx x^3y,\qquad xy\approx xy^3,\qquad x^4\approx x^2;\quad\textup{(\ensuremath{M})}}

\CaseDisplay{x\leq xy^2;\quad\textup{(\ensuremath{E})}}

\CaseDisplay{xy\,zu\,vw\leq xy+zu+vw;\quad\textup{(\ensuremath{O})}}

\CaseDisplay{x^3\leq x+xu,\qquad x^3\leq x+ux,\qquad x^3\leq x+uxv.\quad\textup{(\ensuremath{C})}}

Add the six identities

\CaseDisplay{h^2xy\leq xy+v_h+q_h,\quad\textup{(\ensuremath{H})}}

where \(v_h\) independently takes the three forms \(hp, ph, phr\), and \(q_h\) takes the two forms \(h, hq\). Variables in these forms are distinct placeholders unless a substitution identifies them. All substitutions are nonempty words or terms; there are no unit constants.

Theorem. The fifteen identities (M), (E), (O), (C), (H) form a complete basis of \(A\). Reversing every word gives a complete basis of the opposite.

\par\smallskip\noindent\textit{2. Soundness and word normalization.}\quad 
The head lemma proved for source 209 holds here as well: an absorption is valid in \(A\) exactly when it is valid in \(F\) and its target head occurs among the right heads. If a right head evaluates to \(t\), the inequality is automatic; otherwise replacing every \(t\) by zero preserves the relevant word values. The converse head condition is separated by assigning the target head \(t\) and all other variables zero.

Each nonlinear word on the core has value zero as soon as one of its variables has value \(n\) or zero, and otherwise has its group value. This proves (M): both sides have the same head, content, and parity vector, and remain nonlinear. It also proves (E), since adjoining a square either leaves a group value unchanged or raises it to zero; a linear value \(n\) likewise rises to zero. For (O), each of its three factors is nonlinear and thus has core value in \(G\cup\{0\}\). An odd product of three identical group values equals that value; distinct group values sum to zero, as does a sum containing zero. Its head condition holds.

For (C), if \(x=n\) then the nonlinear right word containing \(x\) is zero; if \(x\) is a group element then \(x^3=x\); the zero case is immediate. For (H), if \(h=n\) or zero then the nonlinear word \(v_h\) is zero. If \(h\in G\), then \(h^2xy=xy\), because \(xy\) is either a group value or zero. The right head is supplied by \(q_h\), so the head lemma proves (H) in the full extension.

Let \(c(w)\) and \(h(w)\) denote content and head, and let \(\pi(w)\) be the vector of variable multiplicities modulo two. From the first identity of (M), every tail after the initial letter can be permuted. The second and third identities permit increasing any exponent-one occurrence in a nonlinear word by two: at the head use a nonempty suffix, and in the tail use a nonempty prefix. Move repeated occurrences together by tail permutations, then use \(x^4=x^2\). Consequently every nonlinear word has a normal form determined by its head, content and parity: every variable has exponent two if even and three if odd, the head is placed first, and the other variables follow a fixed ordering. Thus two nonlinear words with these same three data are provably equal. This argument includes single-variable nonlinear words.

Identity (E) permits adjoining squares of additional variables on the right, so any word is derivably below the word obtained by any such content enlargement.

\par\smallskip\noindent\textit{3. Exact criterion.}\quad 
Write a nonempty polynomial \(P\) as linear variables \(L\) and nonlinear summands \(N\). Let \(U=\bigcup_{v\in N}c(v)\) and \(H(P)=\{h(v):v\text{ a summand of }P\}\). Let \(\mathcal A(P)\) be the affine hull over \(\mathbb F_2\) of all summand parity vectors.

Criterion. A linear target \(x\) is absorbed exactly when \(x\in H(P)\) and \(\pi(x)\in\mathcal A(P)\). A nonlinear target \(w\) is absorbed exactly when

\CaseDisplay{h(w)\in H(P),\qquad c(w)\subseteq U,\qquad \pi(w)\in\mathcal A(P).\quad\textup{(\ensuremath{1})}}

Necessity of the head condition is the head lemma. Restricting every variable to \(G\) gives necessity of the affine condition: if a vector lies outside the affine hull, a linear functional over \(\mathbb F_2\) is constant on the right vectors and takes the other value on the target, separating incomparable group values. For a nonlinear target and a variable \(x\notin U\), if \(x\) appears on the right it appears only linearly. Assign \(x=n\), every other variable \(e\). Then \(P=n\) while the target is zero. If \(x\) does not appear on the right, assign it zero instead, obtaining \(P=e\) and target zero. This proves the coverage condition. Sufficiency follows from the complete derivations next.

\par\smallskip\noindent\textit{4. Arbitrary-variable derivations.}\quad 
First let \(w\) be nonlinear and satisfy (1). Each linear variable outside \(U\) occurs in no nonlinear summand and in no other linear summand. In any odd parity representation of \(\pi(w)\) as a sum of right parity vectors, its coefficient must therefore be zero. Thus \(\pi(w)\) lies in the affine hull of \(N\) together with the linear variables in \(L\cap U\).

For each \(x\in L\cap U\), choose a nonlinear right word \(v\) containing \(x\). According to the position of an occurrence of \(x\), it has form \(xu\), \(ux\), or \(uxv\) with nonempty substituted factors. One of (C) gives \(x^3\leq x+v\leq P\). Let \(Q\) consist of all original nonlinear words and the cubes of these covered linear variables. Each word of \(Q\) is below \(P\), its union of contents is \(U\), and its affine parity hull contains \(\pi(w)\).

Put \(h=h(w)\). Coverage gives an original nonlinear word \(v_h\) containing \(h\), and the head condition gives an original summand \(q_h\) beginning with \(h\). Choose the appropriate one of the six schemes (H), splitting any \(u\in Q\) into two nonempty factors. It gives

\CaseDisplay{h^2u\leq u+v_h+q_h\leq P.}

These new nonlinear words all have head \(h\); they have the same parities as \(Q\) and still have union of contents \(U\). Choose an odd number of them whose parity vectors sum to \(\pi(w)\), and append two more copies of every member of this family. The resulting product \(V\) has head \(h\), content \(U\) and parity \(\pi(w)\). By square extension (E), followed by nonlinear normalization (M), \(w\leq V\) is derivable.

Repeated use of (O) bounds any odd product of nonlinear factors by their sum: group the first three, then induct, the grouped product remaining nonlinear. Hence \(V\leq P\), proving \(w\leq P\). All choices involve finite sums for the given polynomial, but no bound on its number of variables or summands.

For a linear target \(x\) satisfying the criterion, the case \(x\in L\) is immediate. Otherwise \(x\in H(P)\) supplies a nonlinear summand headed by \(x\), so \(x\in U\). The nonlinear word \(x^3\) has the same parity and satisfies (1), whence \(x^3\leq P\) by the preceding argument. Substituting \(y=x\) in (E) gives \(x\leq x^3\), and therefore \(x\leq P\).

Thus every valid word absorption is derivable. Every polynomial identity is the two families of absorptions of its individual summands into the other side, after distribution. This proves the fifteen-identity basis for arbitrary variable sets.

\subsection{三阶群平坦核心的零克隆：五条完整基}\label{case:65}
\par\smallskip\noindent\textit{1. 生成半环与定理.}\quad 
取三阶循环群 \(C_3=\{e, g, g^2\}\)，添乘法零 \(c\)，其加法为平坦半格：三个群元素两两不可比，任意两个不同群元素之和为 \(c\)。再添元素 \(t>c\)。令 \(\rho(t)=c\)，\(\rho\) 固定 \(C_3\cup\{c\}\)，乘法为 \(xy=\rho(x)\rho(y)\)。\(\rho\) 保持并；各乘法平移也保持并，故得到五元素加法幂等半环 \(R\)。在顺序 \((c, e, g, g^2, t)\) 下，\(R\) 正是完整目录编号 1638。

定理。 \(R\) 的恒等式簇由下列五条式子相对于通常加法幂等半环公理完整公理化：

\CaseDisplay{\begin{aligned}
xy&\approx yx, &(1)\\
x^4y&\approx xy, &(2)\\
x^3y^3&\approx x^3+y^3, &(3)\\
x^4&\leq x, &(4)\\
xyzt&\leq x^4+y^4+z^4+t^4. &(5)
\end{aligned}}

五条式子均直接有效：(2) 是乘法核心的周期三性质；(3) 中所有群元素的立方为 \(e\)，\(c, t\) 的立方为 \(c\)；(4) 只在 \(x=t\) 时为真不等式 \(c<t\)；(5) 在群内由``若四值相同，则其乘积仍为该值''推出，遇到 \(c, t\) 时右侧含 \(c\)，足以吸收左值 \(c\)。

\par\smallskip\noindent\textit{2. 任意变量数判据.}\quad 
固定有限非空变量集 \(X\)。对非空词 \(w\)，令 \(c(w)\) 为内容，\(\pi_3(w)=(m_x(w)\bmod3)_{x\in X}\in\mathbb F_3^X\)。若 \(P=\sum_i v_i\)，写

\CaseDisplay{U(P)=\bigcup_i c(v_i),\qquad H_3(P)=\operatorname{aff}_{\mathbb F_3}\{\pi_3(v_i):i\}.}

命题。 一次词 \(x\leq P\) 在 \(R\) 中有效，当且仅当 \(x\) 是 \(P\) 的一次项。非线性词 \(w\leq P\) 在 \(R\) 中有效，当且仅当

\CaseDisplay{\pi_3(w)\in H_3(P),\qquad c(w)\subseteq U(P).\quad\textup{(\ensuremath{6})}}

证明。 若仿射条件失败，群上的线性泛函使所有右词取同一群值、左词取另一个群值；后者不被前者吸收。若内容条件失败，给一个左侧独有的变量赋 \(c\)，其余变量赋 \(e\)，左值为 \(c\)，右值为 \(e\)。若一次左词 \(x\) 不在右侧一次项中，给 \(x\) 赋 \(t\)、其他变量赋 \(e\)，全部非线性右词至多为 \(c<t\)，故吸收失败。

反过来，若非线性左词取 \(c\)，其内容内至少有一个变量取 \(c\) 或 \(t\)；内容条件保证某个右词也含它，该右词取 \(c\) 或 \(t\)，从而吸收左值。若左词取群值，将右侧外部的 \(c, t\) 替换成 \(e\)。实际每个右词的值均不低于对应的群值，仿射条件保证投影后的右侧吸收左值。一次词出现在右侧时显然被吸收。证毕。

\par\smallskip\noindent\textit{3. 完整性证明.}\quad 
由 (1)、(2)，每个非线性词的值仅由内容和指数模三决定；单变量非线性词按指数模三约化为平方、立方或四次幂。特别地，非线性 \(u\) 满足 \(u^4\approx u\)，且若 \(x\in c(u)\)，则 \(ux^3\approx u\)。

内容扩充引理。 若 \(u\) 非线性，则 (1)--(3) 推出 \(u\leq uz^3\)。

证明。 选 \(x\in c(u)\)，用 (3) 乘以 \(u\)：

\CaseDisplay{uz^3=ux^3z^3=ux^3+uz^3=u+uz^3.}

证毕。

四词引理。 若 \(q_1,\ldots, q_{3r+1}\) 均为非线性词，则 (1)、(2)、(5) 推出

\CaseDisplay{q_1\cdots q_{3r+1}\leq q_1+\cdots+q_{3r+1}.}

这由 \(q_i^4=q_i\) 和对 \(r\) 的归纳得到：将前 \(3r-2\) 项的乘积与最后三项代入 (5)。

设非线性左词 \(w\) 满足 (6)。由 (4)，对右侧每个一次项 \(x\)，都可向 \(P\) 加入 \(x^4\)，不改变 \(P\)。令 \(Q\) 为所得全非线性右词族。此时

\CaseDisplay{U(Q)=U(P),\qquad H_3(Q)=H_3(P).}

因为 \(\pi_3(w)\in H_3(Q)\)，它是 \(Q\) 中若干向量的组合，其系数之和为 1 模 3。选取相应非负次数的右词，再将 \(Q\) 中每个词各补入三次。所得乘积 \(V\) 含 \(3r+1\) 个非线性因子，满足 \(\pi_3(V)=\pi_3(w)\)、\(U(Q)=c(V)\supseteq c(w)\)，且四词引理给出 \(V\leq\sum Q\leq P\)。词规范形给出

\CaseDisplay{V\approx w\prod_{z\in c(V)\setminus c(w)}z^3.}

逐次应用内容扩充引理，得到 \(w\leq V\leq P\)。一次左词的有效吸收式由命题已是右侧显式单项。因此每个有效逐词吸收式由 (1)--(5) 推出，双向逐词吸收又给出全部多项式恒等式，定理得证。

\par\smallskip\noindent\textit{4. 范围.}\quad

\subsection{素数阶群的左零顶扩张：五条完整基}\label{case:66}
\par\smallskip\noindent\textit{1. 一般构造与定理.}\quad 
设 \(p\) 是任意素数，\(G=C_p\) 是 \(p\) 阶循环群。取 \(G\) 的平坦零扩张 \(G\cup\{c\}\)：\(c\) 是乘法零，任意两个不同群元素之和为 \(c\)。再添加法最大元 \(t>c\)，定义

\CaseDisplay{ta=t\quad(a\in G\cup\{c, t\}),\qquad
at=c\quad(a\in G\cup\{c\}).}

得到半环 \(L_p\)。结合律由理想 \(\{c, t\}\) 的左零乘法及群的左作用同时把 \(c, t\) 送到 \(c\) 得到；所有乘法平移保持所给平坦加法与顶层。

定理 1。 \(L_p\) 的完整恒等式基为

\CaseDisplay{\begin{aligned}
xyz&\approx xzy, &(1)\\
x^{p+1}&\approx x, &(2)\\
x&\leq xy^p, &(3)\\
x_1\cdots x_{p+1}&\leq x_1+\cdots+x_{p+1}, &(4)\\
x^py&\leq y+xz. &(5)
\end{aligned}}

在每个词中反转字母顺序，得到 \(L_p\) 的乘法对偶的五条完整基。

这些式子有效：(1) 保持首字母并只交换尾部；(2) 在群内是指数为 \(p\)，在 \(c, t\) 上是幂等性；(3) 对群值乘 \(y^p\) 后要么保持该值，要么升至 \(c\)，在首值为 \(t\) 时保持 \(t\)。(4) 若有某个因子为 \(t\)，右側为 \(t\)；若有 \(c\) 而无 \(t\)，右側为 \(c\)；否则在群内，如果所有因子相等则乘积等于共同值，不全相等时右側为 \(c\)。(5) 在 \(x\in G\) 时左侧为 \(y\) 或在 \(y=t\) 时为 \(c\leq t\)；在 \(x=c\) 时 \(xz=c\)；在 \(x=t\) 时 \(xz=t\)。

\par\smallskip\noindent\textit{2. 首字母、内容与仿射判据.}\quad 
对非空词 \(w\)，记 \(h(w)\) 为首字母，\(c(w)\) 为内容，\(\pi_p(w)\in\mathbb F_p^X\) 为各变量指数模 \(p\) 的向量。对非空多项式 \(P=\sum_i v_i\)，记

\CaseDisplay{\mathrm{Head}(P)=\{h(v_i):i\},\quad U(P)=\bigcup_i c(v_i),\quad
H_p(P)=\operatorname{aff}_{\mathbb F_p}\{\pi_p(v_i):i\}.}

命题 2。 \(L_p\models w\leq P\) 当且仅当

\CaseDisplay{h(w)\in\mathrm{Head}(P),\qquad c(w)\subseteq U(P),\qquad
\pi_p(w)\in H_p(P).\quad\textup{(\ensuremath{6})}}

必要性。 若首字母 \(h=h(w)\) 不在右侧首字母集中，赋 \(h=t\)、其余变量为群单位。左侧为 \(t\)，所有右词至多为 \(c\)。若内容包含失败，给左侧独有变量赋 \(c\)、其余变量赋群单位，得到左侧 \(c\)、右侧群单位。若仿射条件失败，在群上选一个线性泛函，使全部右向量取共同值、左向量取不同值；对应赋值使左右为不同群元素，吸收失败。

充分性将由以下有限基推导直接给出，无须另加有界变量假设。

\par\smallskip\noindent\textit{3. 从五条式子推出所有吸收式.}\quad 
由 (1)、(2)，每个词在首字母固定的情况下，只由内容和指数模 \(p\) 决定：尾部可任意交换，每个变量的正指数按模 \(p\) 约化，但保留首字母的一个出现。相同首字母、相同内容、相同模 \(p\) 向量的词因此可相互推导。

由 (3)，向任意词右侧加入变量的 \(p\) 次幂使其不减。因此若 \(u, v\) 首字母相同、\(c(u)\subseteq c(v)\) 且 \(\pi_p(u)=\pi_p(v)\)，先补入缺少变量的 \(p\) 次幂，再用上述规范化，得到 \(u\leq v\)。

由 (4) 归纳还得到：任意 \(pr+1\) 个非空词之积被这些词之和吸收。归纳时，将前 \(p(r-1)+1\) 个词的乘积和最后 \(p\) 个词代入 (4)，再用归纳假设。

设 \(w\leq P\) 满足 (6)，令 \(h=h(w)\)。选择一个以 \(h\) 开始的右词 \(v\)。若它非线性，写作 \(hz\)；若 \(v=h\)，取代入项 \(z=h^p\)，由 (2) 仍有 \(hz=h\)。对每个右词 \(q\)，(5) 给出

\CaseDisplay{h^p q\leq q+v\leq P.}

故右侧可以提供一个全部以 \(h\) 开始的词族 \(Q=\{h^p q:q\text{ 是 }P\text{ 的词}\}\)，它满足

\CaseDisplay{\sum Q\leq P,\qquad U(Q)=U(P),\qquad H_p(Q)=H_p(P).}

第二个等式使用 \(h\in U(P)\)，第三个使用加入 \(p\) 次 \(h\) 不改变模 \(p\) 向量。

按 \(\pi_p(w)\in H_p(Q)\) 选择一个非负整数系数的表示，系数之和为 1 模 \(p\)；再将 \(Q\) 中每个词各补入 \(p\) 次。得到一个含 \(pr+1\) 个因子的乘积 \(V\)，满足

\CaseDisplay{h(V)=h(w),\qquad c(w)\subseteq c(V)=U(P),\qquad
\pi_p(V)=\pi_p(w).}

由内容扩充推导 \(w\leq V\)，由积的吸收引理得 \(V\leq\sum Q\leq P\)。这样，任意有效逐词吸收式都由五条基推出；任意多项式恒等式分解为双向逐词吸收，定理得证。

\par\smallskip\noindent\textit{4. 五元素目录.}\quad 
\(L_3\) 在顺序 \((c, e, g, g^2, t)\) 下正是目录 1640；其乘法对偶正是 2795。目录 2191 与 1640 等簇，7607 与 2795 等簇，两个方向均由完整词函数与多项式下闭集方法独立核验。四个对象因此有限基，不需要把它们彼此同构。

\subsection{Finite identity bases for 1650, 1651 and their opposites}\label{case:67}
Work in the constant-free ai-semiring signature. The additive order in each source is \(3<2<1<0<4\). In catalogue order the product rows of 1650 are 00000,\allowbreak{}01230,\allowbreak{}33330,\allowbreak{}33333,\allowbreak{}44444. Those of 1651 are 00000,\allowbreak{}01230,\allowbreak{}33332,\allowbreak{}33333,\allowbreak{}00004. The corresponding opposites are 5378 and 5288. Reversing each word of a basis gives the opposite basis.

For a target w, let C be its support. At the first occurrence of each x let S(x) be the support strictly before it. A prefix witness for x with allowance S is an original right summand containing x and having its prefix before its first x supported in S. Its suffix may use outside variables. All the following characterizations and proofs allow arbitrarily many variables and right summands.

\par\smallskip\noindent\textit{Word values and absorption.}\quad 
An all-1 word gives 1. Read the first entry different from 1. An initial decisive 0 or 3 gives itself. When that entry is 4, 1650 gives 4 exactly if the word starts with 4, and 1651 gives 4 exactly for an all-4 word. All other cases with this decisive 4 give 0. When the first non-1 is 2, 1650 returns 2 if final, 0 if immediately followed by 4, and 3 otherwise. In 1651 it returns 2 if everything after it is 4, possibly vacuously, and 3 otherwise. These formulas follow by multiplying one further factor.

Both tests require an original word supported in C and a prefix witness for every x with allowance S(x). The first requirement is forced by giving C the value 1 and outsiders 3. For the second give S(x) the value 1, x the value 0, and all others 3.

In 1650 impose three further requirements.

\textbullet{} An original starts with the head of w. Assigning its head 4 and all others 3 proves necessity.

\textbullet{} For every pair of consecutive positions x y both of which are first occurrences, there is a prefix witness for y with allowance S(x), or an original word has a first x preceded only by S(x) and immediately followed by y. Give S(x) the value 1, x the value 2, y the value 4, and the rest 3. Only the two kinds of right words can attain 0 or 4.

\textbullet{} If the final t is simple, some original is supported in S(t) or an original with a unique final t has all its preceding letters in S(t). Use value 1 on S(t), 2 on t and 3 elsewhere.

These tests suffice. A target of 4 uses the original head. An all-1 target uses its contained witness. When the decisive value is 0 or 4 its prefix witness gives at least 0. A decisive 2 yielding 0 is necessarily followed by a fresh letter of value 4, so the consecutive-pair alternatives attain at least 0. A final decisive 2 uses the last test. Value 3 is least.

In 1651 impose the following additional test for each simple x whose prefix support S and suffix support T are disjoint. Require at least one of these original witnesses:

\textbullet{} a word supported in S,

\textbullet{} a prefix witness for a letter y of T with allowance S,

\textbullet{} a word v x z with v supported in S and z supported in T, either block possibly absent.

Assign \(x=2\), \(S=1\), \(T=4\), \(outsiders=3\). The words attaining at least 2 are exactly the three kinds. Conversely the only way a target can give 2 is with a unique 2 between all 1's and all 4's. Its variable is simple and its prefix and suffix supports are disjoint. Each of the alternatives then suffices. At value 4 use the contained word on the all-4 target. For values 1 and 0 use respectively containment and the prefix witness at the first non-1 entry. This proves the test.

\par\smallskip\noindent\textit{A 17-law presentation for 1650.}\quad 
Optional blocks in brackets produce separate identities by omission before substitution. All retained symbols denote nonempty arguments. An occurrence of [pp] is retained exactly when the corresponding [p] is present. Take the following fixed and optional laws.

{\small\par\smallskip
\CaseLine{\(xxx=xx\qquad xyx=xyy\qquad xyzx=xyzz\)}
\CaseLine{\(xxu \leq  u+x[b]\qquad (2)\)}
\CaseLine{\(dyyu \leq  du+dy[b]\qquad (2)\)}
\CaseLine{\([p]xyyu \leq  [p]xxu+[pp]y[b]\qquad (4)\)}
\CaseLine{\([p]xyyu \leq  [p]xxu+[pp]xy[b]\qquad (4)\)}
\CaseLine{\(px \leq  pxx+pp\qquad px \leq  pxx+ppx\)}
\smallskip}

The first three equations preserve the order of first occurrences and the status of every such occurrence: whether it is followed immediately by another fresh letter, is final, or is followed by an old letter. These data determine the word function by the table formula. They also yield a normal form consisting of the list of first occurrences, doubling each letter whose first occurrence is followed by an old letter. Indeed, append one letter to a normalized prefix. A new letter is retained. For an old letter different from the current last b, use \(xyx=xyy\) or \(xyzx=xyzz\) at its earlier occurrence to replace the final input copy by b. Contract any resulting triple by \(xxx=xx\). An appended b itself already follows the same rule. This inductively proves the normal form using actual substitutions and contexts.

The inequalities satisfy all three additional tests and the two common tests. In the two four-instance families the previous pair from p to x is witnessed by the first right term whenever p is present. At the new pair xy the second right term provides the respective alternative. The variable y itself has been doubled. The simple final u has its contained simple-final witness. The two final recovery laws use respectively pp and ppx. Soundness of each formal identity also admits a finite exhaustive table check.

Normalize w and list its fresh letters \(x_{1}\),..., \(x_{n}\). Whenever a first occurrence has been reached, let \(D_{i}\) be the normalized target prefix ending with a doubled \(x_{i}\), including all earlier prescribed gaps. Derive \(D_{i}\) u below P for every original u. To start, choose an original \(x_{1}\) b (or its literal) from the head condition and apply the first two-instance family. Suppose \(x=x_{i}\) and \(y=x_{i+1}\). If there is a gap, choose a prefix witness a y b with a using only the variables of \(D_{i}\). The available \(D_{i}\) a y b normalizes to \(D_{i}\) y b because \(D_{i}\) ends in a doubled letter. Apply dyyu \ensuremath{\leq} du+dy[b] with \(d=D_{i}\). This keeps the prescribed gap.

If there is no gap, write p for the normalized prefix strictly before x. When p is nonempty, its square normalizes to \(D_{i-1}\). The consecutive-pair test supplies a y b or a x y b with a supported in p. Use the previous stage \(D_{i-1}\) on that original, normalizing away a to obtain pp y b or pp x y b. When p is empty the original itself gives y b or x y b. Together with the current \(D_{i}\) \(u=p\) x x u these are the exact premises of one of the four-instance families. Its conclusion retains the desired adjacency and doubles y. This completes the induction without constraining the later suffix b of an original witness.

At the last stage take a contained original as u. Normalization gives \(D_{n}\). For a repeated final this already is w. For a simple final t with nonempty preceding prefix p, apply the penultimate stage to the witness of the simple-final test. Its normal form gives respectively pp or pp t, since every preceding witness letter already occurs in p. The last two laws derive p t from \(D_{n}=p\) tt. A literal target needs its literal among the originals. This covers all absorptions.

\par\smallskip\noindent\textit{A finite presentation for 1651.}\quad 
Write sq(v) for vv. Products of the squared letters have the first- occurrence order of a left regular band. Use the following five fixed laws, with the final three expressing multiplication of squares, idempotence and the left regular band equation on squares:

{\small\par\smallskip
\CaseLine{\(xy = xxy+xyy\qquad xx \leq  x\)}
\CaseLine{\(sq(x) sq(y) = sq(xy)\qquad sq(x) sq(x)=sq(x)\)}
\CaseLine{\(sq(x) sq(y) sq(x)=sq(x) sq(y)\)}
\smallskip}

On the table, squaring sends 2 to 3 and fixes the other entries. It is multiplicative, as is also verified directly on the 25 pairs. Words on the squared elements depend only on their first-occurrence orders. To compare squared words, expand sq of a product, delete separated repetitions by the last equation (combining an intervening product back into one square), and contract adjacent repetitions. The splitting law follows directly on the same 25 pairs.

The remaining schemes have independently optional capital blocks. For legibility denote possibly absent words by P, S, V, Z, A, B, D and F in these schemes only. Omit their squares as well when absent. Lower-case variables always remain. Take all instances of the following squared-prefix insertion, where U is written as the always-present lower-case u:

{\small\par\smallskip
\CaseLine{\(sq(P x u) \leq  sq(P u) + sq(P x B).\)}
\smallskip}

For each of the following three rows form the identity \(L=sq(L)\), where the displayed left word is L:

sq(A x B) x sq(S) sq(P) x sq(A x B) sq(A y B) x sq(D y F)

Next add the following recovery inequalities. In the first row p is mandatory. In the second row y is mandatory.

{\small\par\smallskip
\CaseLine{\(sq(p) x sq(S) \leq  sq(p x S) + sq(p)\)}
\CaseLine{\(sq(P) x sq(D y F) \leq  sq(P x D y F) + sq(P y B)\)}
\CaseLine{\(sq(P) sq(V) x sq(S) sq(Z) \leq  sq(P V x S Z) + V x Z\)}
\smallskip}

These specify at most 75 displayed laws (5 fixed, 4 insertions, 32 equations and 2+16+16 recoveries). For soundness, at a general assignment L and sq(L) differ only at \(L=2\). Every occurrence inside a square is repeated, so only the stand-alone x can be the unique

2. In the three equations a repeated x or a letter common to both sides of x prevents this. In the recoveries, if the left gives 2 its factors before x are all 1 and after x all 4. The second upper term then has respectively value 1, a value at least 0 from its first 4, or value 2. For every other left value the first upper term agrees by the square-band word formulas. The insertion passes the prefix and containment tests, with no surviving 2 in its squared left word.

To derive any absorption, first expand each lower word by the splitting law. Repeatedly split its unsquared subword into two nonempty blocks and recurse on the only remaining unsquared block of each resulting summand. This terminates with each lower word of the form \(L=sq(p)\) x sq(s) by merging adjacent squares. Either side of x may be absent. All expanded lower terms are below P by validity of the original absorption and of the splitting identity. It suffices to handle each expanded term separately.

For any initial segment p of the first-occurrence order of this target, derive sq(p u) for every original u. The empty stage is \(sq(u)\leq u\). To append x, choose an original a x b with c(a) within c(p) by the prefix test. The preceding stage gives sq(p a x b), which normalizes to sq(p x b). Apply the squared insertion using sq(p u). Hence for the full target support, using a contained original, we derive sq(L) in its correct first-occurrence order. This argument also supplies sq(p u) with any spelling of a prefix having the same first-occurrence list.

Return to \(L=sq(p)\) x sq(s). If x occurs in p or in s, or some letter y is common to p, s, use the appropriate \(L=sq(L)\) scheme, splitting the actual blocks at the indicated occurrence. Otherwise x is simple and the supports of p, s are disjoint. If a word u is contained in p, the earlier prefix stage gives sq(p u)=sq(p) and the first recovery applies. In the second alternative take an original a y b with a supported in p and y in s. The earlier stage (including the empty stage) gives sq(p y b); split s at y and apply the second recovery. In the final alternative take the original v x z with v supported in p and z supported in s. Insert sq(v) after sq(p) and sq(z) after sq(s) in L. The square-band equations normalize these products back to the respective original blocks, and the long first right term to sq(L). The second right term is exactly the original v x z. The final scheme thus derives L. Empty pieces everywhere select already listed identities, never empty substitutions.

Finally distributivity expands arbitrary terms into sums of words. Derive each summand below the other side in the two directions to obtain every valid identity. None of the constructions bounds the size of the alphabet or the lengths of the input words.

\subsection{最后出现后缀与首字母保护：收缩型五条完整基}\label{case:68}
\par\smallskip\noindent\textit{1. 构造及正规形.}\quad 
取 \(A=\{z, e, b, d, t\}\)，依次对应目录 \((0,\allowbreak{}1,\allowbreak{}2,\allowbreak{}3,\allowbreak{}4)\)，加法取链 \(d<b<e<z<t\) 上的最大值。右乘 \(e\) 是 \(R=(z, e, b, d, z)\)；右乘 \(b, d\) 为常值 \(d\) 映射；右乘 \(z, t\) 为常值 \(z\) 映射。注意 \(R(t)=z\)，所以 \(e\) 不是全半环的乘法单位元。

三个右平移保序，故保持链上并；它们与右因子乘法、加法兼容，直接按三类检查 \(R_{xy}=R_yR_x\)、\(R_{x+y}=R_x\vee R_y\)，即可验证半环公理。

非一次词的尾部全取 \(e\) 时，取 \(R(h(w))\)；否则找尾部最右一个不取 \(e\) 的字母，它取 \(b, d\) 时全词取 \(d\)，取 \(z, t\) 时全词取 \(z\)。因此非一次词不会取最大元 \(t\)。

词恒等式

\CaseDisplay{xy^2\approx xy,\qquad xyx\approx y^2x\quad\textup{(\ensuremath{W})}}

成立，并推出 \(x^3\approx x^2\)、\(xyzy\approx xzy\)。与最后出现后缀的五条基第2节同样，每个词正规化为无重复词 \(r\) 或其平方 \(r^2\)；平方正规形也可写成 \(r\) 的首字母接整个 \(r\)。尾部仅保留每个变量最后一次出现，原首字母在尾部再次出现时使整个词成为平方形式。

记 \(S_w(x)\) 为变量 \(x\) 最后一次出现之后的变量集合，\(D(w)\) 为尾部集合。以下先把左右词正规化。

\par\smallskip\noindent\textit{2. 精确吸收判据.}\quad 
令 \(F\) 为非空右词集。一次吸收 \(x\le F\) 有效，当且仅当一次项 \(x\) 本身属于 \(F\)。

定理1。 对非一次词 \(w\)，令 \(X=c(w)\)、\(h=h(w)\)。则 \(A\models w\le F\) 当且仅当：

1. 存在内部右项 \(v\)，\(c(v)\subseteq X\)；若 \(h\notin D(w)\)，还要求 \(h\notin D(v)\)；

2. 对每个 \(x\in X\)，有右项 \(u_x\) 包含 \(x\)，且 \(S_{u_x}(x)\subseteq S_w(x)\)。

条件1只在左首字母不重复时禁止内部见证在尾部出现它，不要求见证一定包含该首字母。条件2允许见证中的 \(x\) 只在首字母或一次项出现。

必要性。 对一次左词，若 \(x\notin F\)，赋 \(x=t\)、其余变量 \(d\)。所有右非一次词至多为 \(z\)，其余一次项为 \(d\)，不能吸收 \(t\)。

对非一次左词，若没有内部项，赋 \(X\) 内变量 \(e\)、外部变量 \(d\)，左值为 \(e\)、所有右项为 \(d\)。如果内部项存在，但 \(h\notin D(w)\) 且每个内部项都在尾部含 \(h\)，赋 \(h=b\)、其余 \(X\) 内变量 \(e\)、外部变量 \(d\)。左值为 \(b\)；内部项尾部含 \(b\)，外部项含 \(d\)，其余字母只为 \(e, b, d\)，故所有右词取 \(d\)。

若条件2在 \(x\) 处失败，赋 \(x=z\)、\(S_w(x)\) 内变量 \(e\)、其余变量 \(d\)。左值为 \(z\)；任何含 \(x\) 的右词在最后一个 \(x\) 之后都有取值 \(d\) 的变量，故取 \(d\)。其他右词不超过 \(e\)，吸收失败。

充分性。 非一次左词只能取 \(d, b, e, z\)。值为最小元 \(d\) 时自明。左值为 \(b\) 时，左首字母唯一出现且取 \(b\)，左尾部全部取 \(e\)。条件1中的内部见证若不含首字母，则取 \(e\)；若含首字母，因为不能在尾部包含它，取 \(b\)。两者均吸收左值。

左值为 \(e\) 时全部左变量取 \(e\)，任一内部见证也取 \(e\)。左值为 \(z\) 时，取左词最右一个值不为 \(e\) 的变量 \(x\)；它取 \(z\) 或 \(t\)，而其后全部变量取 \(e\)。条件2给出的后缀见证值为 \(z\) 或 \(t\)，足以吸收左值。判据得证。

\par\smallskip\noindent\textit{3. 五条完整基.}\quad 
定理2。 下列五式构成完整附加基：

\CaseDisplay{\begin{aligned}
xy^2&\approx xy, &(R)\\
xyx&\approx y^2x, &(C)\\
xy&\le x+y, &(Q)\\
xy&\le x+zy, &(V)\\
xy&\le y^2+zxy. &(J)
\end{aligned}\quad\textup{(\ensuremath{1})}}

前两式由求值规则成立，后三式由定理1逐式检验。特别地，\((Q)\) 令 \(y=x\) 给出 \(x^2\le x\)；不能反向使用 \(x\le x^2\)。

\par\smallskip\noindent\textit{3.1 允许一次项的后缀推进.}\quad 
设 \(r=x_1\cdots x_k\) 无重复，每个 \(x_i\) 有右项见证

\CaseDisplay{u_i=a_ix_ib_i,\qquad c(b_i)\subseteq\{x_{i+1},\ldots, x_k\},}

在最后一个 \(x_i\) 处分割，前段 \(a_i\) 和后段 \(b_i\) 都可以为空。证明：对每个可推出 \(v\le F\) 的非空词 \(v\)，都有

\CaseDisplay{v x_i\cdots x_k\le F\quad(1\le i\le k).\quad\textup{(\ensuremath{2})}}

从右到左归纳。最右变量的后段为空：若 \(a_k\) 非空，\((V)\) 给出 \(vx_k\le v+a_kx_k=v+u_k\)；若 \(a_k\) 为空，\(u_k=x_k\)，\((Q)\) 给出 \(vx_k\le v+x_k\)。

设对后缀 \(r_{i+1}=x_{i+1}\cdots x_k\) 已成立。把右项 \(u_i\) 代作 \(v\)，得 \(u_ir_{i+1}\le F\)；由于 \(b_i\) 的全部变量在后缀再次出现，尾部去重给出

\CaseDisplay{u_ir_{i+1}\approx a_ix_ir_{i+1}.}

若 \(a_i\) 非空，由 \((V)\) 得 \(vx_ir_{i+1}\le v+a_ix_ir_{i+1}\le F\)。若 \(a_i\) 为空，则 \(x_ir_{i+1}\le F\) 已被导出，由 \((Q)\) 得 \(vx_ir_{i+1}\le v+x_ir_{i+1}\le F\)。这完成归纳，且无需空词代入。

\par\smallskip\noindent\textit{3.2 保留未重复的首字母.}\quad 
若左词为平方形式 \(w=r^2\)，选条件1的任一内部项 \(v\)。对全部有序内容应用 (2)，得到 \(vr\le F\)；由 \(c(v)\subseteq c(r)\)，正规形 \(vr\approx r^2=w\)。

若左词为无重复非一次词 \(w=hr\)，条件1提供内部见证 \(v\) 且 \(h\notin D(v)\)。对 \(r\) 应用 (2)，得到 \(vr\le F\)。

若 \(h\in c(v)\)，则它是 \(v\) 的首字母并且不重复，尾部变量都在后面的 \(r\) 再次出现，所以 \(vr\approx hr=w\)。这里首字母保护保证取得 \(w\) 本身，不会误取无法反向提升的 \(w^2\)。

若 \(h\notin c(v)\)，则 \(vr\approx r^2\le F\)。取条件2对 \(h\) 的见证，在最后一个 \(h\) 处分割为 \(u=ahs\)、\(c(s)\subseteq c(r)\)。由 (2) 得 \(ur\le F\)。若 \(a\) 为空，\(ur\approx hr=w\)；若 \(a\) 非空，\(ur\approx ahr\)，由 \((J)\) 得

\CaseDisplay{w=hr\le r^2+ahr\le F.}

一次有效吸收由同名右项直接得到。将一般恒等式展开为词和并处理两向逐词吸收，就得到任意变量数完备性。证毕。

\par\smallskip\noindent\textit{4. 对偶与复核范围.}\quad 
反转所有词得到乘法对偶的五条基。独立计算重建来源运算、逐赋值检查五式，并核验传递对象的双向HSP包含。小变量判据核验只作交叉检查，本定理的完整性由上面的后缀归纳提供。它不证明整个簇遗传有限基。

\subsection{A 22-identity basis for 5291}\label{case:69}
The signature contains addition and multiplication only. Every substitution is a nonempty word; \(p\leq q\) abbreviates \(p+q=q\).

\par\smallskip\noindent\textit{Tables and coordinates.}\quad 
Use labels0,\allowbreak{}1,\allowbreak{}2,\allowbreak{}3,\allowbreak{}4 with multiplication rows 44224,\allowbreak{}41224,\allowbreak{}42224,\allowbreak{}02234,\allowbreak{}44224. Addition is MINIMUM in \(3<2<1<0<4\). The permutation(2,\allowbreak{}1,\allowbreak{}0,\allowbreak{}4,\allowbreak{}3) maps these coordinates to catalogue 5291, respecting both operations. This is independently checked on all 50 operation entries. Thus multiplication is the same as 13066 but the additive order is reversed; its absorption criterion must be proved afresh. We do not use any forward-order absorption theorem.

For an evaluated word find its rightmost factor not 1. If there is none the value is 1. If it is 2 or 4, the value is that factor. If it is 3, the value is 3 when all factors are 3 and is 2 otherwise. If it is 0, the value is 0 exactly when it is the final factor and all preceding factors are 3; otherwise it is 4. This follows by multiplication-table induction. The laws \(xxy=xy\) and \(xyx=yxx\) give the last-occurrence order, with a doubled final letter exactly when that letter was repeated. The reduction contracts internal repeats xYxZ by \(xYxZ=YxxZ=YxZ\); adjacent repeats use \(xxy=xy\). A last remaining final repeat uses \(xYx=Yxx\). All substituted words are nonempty in these cases.

\par\smallskip\noindent\textit{Exact absorption theorem.}\quad 
Let P be the original nonempty finite family of upper words. A target w satisfies \(w\leq \sum\) P exactly when all three conditions hold:

1. An original U has content contained in c(w).

2. For every x in c(w), some original u contains x and its suffix after the last x has content contained in the target suffix after the last x.

3. If the target final letter t is simple, either an original upper word ends in a letter from c(w) minus t, OR an original contained upper word ends in a simple t.

In particular a terminal witness need not repeat its terminal letter. There is no marked-position strengthening or prefix-union condition from 13066. Condition 3 is instead a disjunction of two original guards.

Necessity: if 1 fails, assign target letters 3 and outside letters 1. The target is 3 and every upper word has a larger value in the displayed chain. If 2 fails at x, assign \(x=2\), the target suffix letters 1, and all other letters 4. The target is 2 while every upper word is 1 or 4. If 3 fails, assign \(t=0\), other target letters 3, and outside letters 4. The target is 0. Every original word either ends outside the target, hence is 4, or ends in t but repeats t or contains an outside letter, again giving 4. Ending in a target prefix letter is excluded by failure of 3. These are strict failures of absorption for minimum addition.

Sufficiency: target 4 is automatic. A target 3 has a contained original word of value 3. A target 1 has all target factors 1 and the contained word has value 1. For target 2 its rightmost non1 is 2 or 3; condition 2 provides a word with that same rightmost non1 value and thus value 2 or 3. For target 0, \(t=0\) is simple and every other target letter is 3. An original word ending in a target prefix letter has value 2 or 3; otherwise the contained simple-t word from 3 has value 0. Hence in each case at least one upper value is no greater than the target in the displayed chain, exactly the required minimum absorption. This proves the criterion for arbitrary variables, lengths and summand counts.

\par\smallskip\noindent\textit{Finite basis.}\quad 
Take the semiring axioms and six fixed identities:

{\small\par\smallskip
\CaseLine{\(xxy=xy\qquad xyx=yxx\)}
\CaseLine{\(xx\leq x\qquad xyy\leq x+zxyy\)}
\CaseLine{\(xyt\leq xytt+xt\qquad xt\leq xtt+t\)}
\smallskip}

Add these 16 identities, represented by optional-factor templates:

{\small\par\smallskip
\CaseLine{\([a]xrzz \leq  [a]rzz+[v]xr\qquad (4 choices)\)}
\CaseLine{\([a]xzz  \leq  [a]zz+[v]x\qquad (4 choices)\)}
\CaseLine{\([a]s[b]t \leq  [a]s[b]tt+[v]s\qquad (8 choices)\)}
\smallskip}

Each bracket is retained or omitted when writing a separate identity; it is never an empty substitution. There are 22 laws, at most five variables in each. All 11,\allowbreak{}130 table assignments across this list are checked exhaustively by an independently constructed law list. This establishes soundness of the individual laws, not completeness by tests.

\par\smallskip\noindent\textit{Unrestricted deduction.}\quad 
Normalize every original word and w by the two word laws. Let t be the target final letter. By 2 there is an original word ending in t. Its square, derived by \(xx\leq x\), normalizes to an accumulating word g with a saturated suffix \(S=tt\). Process the other target letters in reverse last-occurrence order. For the next x choose an original witness VxR whose suffix content lies in S. If R is nonempty use the first template with \(r=R\), \(z=S\) and a the portion of g before its protected suffix. Its right terms normalize to g and the original witness; its left normalizes to a word with protected suffix xS. If R is empty, use the second template with \(v=V\), \(z=S\). The witness need not repeat x. Optional empty pieces select separately written identities. Any former x in the accumulating prefix disappears under the word laws because the new x is later. The saturated final tail persists throughout. This induction terminates after the number of distinct non-tail target letters.

The result is \(g=E\) Q, where Q is the target last-occurrence order with doubled t, and E contains only outside letters. Choose the original contained U from 1. If E is nonempty, apply \(xyy\leq x+zxyy\) with \(x=U\), \(y=Q\), \(z=E\). The left normalizes to Q and the right terms to U and E Q, so Q is derived. If E is empty it is already derived.

If the target has repeated t, Q is its normal form and we are done. Otherwise use condition 3. If an original word ends in a prefix letter s, write it V s and split the desired target as A s B t. The third template with those nonempty pieces (omitting empty A, B, V through the appropriate separate law) derives A s B t from A s B tt and the original V s. If there is no such guard, choose an original contained simple-tail word X t. When X is nonempty, use \(xyt\leq xytt+xt\) with \(x=X\), \(y=the\) target prefix, \(t=t\). Both large terms normalize to the required simple and saturated targets. When the original is t alone, use \(xt\leq xtt+t\) with \(x=the\) target prefix. A one-letter simple target must itself be an original word by 3, so no empty target-prefix substitution is required.

Every valid word absorption has now been derived from the 22 laws. Arbitrary polynomial identities split into finitely many such absorptions in both directions using the additive semilattice axioms. Thus these laws form a finite basis for 5291 with no bound on variables, word lengths or polynomial size. Further catalogue classifications require independently checked complete HSP transfers or isomorphisms.

\subsection{Eleven complete identities for semiring 1706}\label{case:70}
Constant-free signature (+, multiplication), relative to the additively idempotent semiring axioms. A complete finite basis is

{\small\par\smallskip
\CaseLine{\(xy=yx; xyz=xyzz; xyzt\leq xyz; xxx\leq x+yyy; xxy\leq xy+zzz;\)}
\CaseLine{\(xy\leq xxy; x\leq xxx; xx\leq x+xy; xz\leq xx+xy+yz;\)}
\CaseLine{\(x\leq xx+xy+y; xy\leq xx+x+z+yz.\)}
\smallskip}

Here \(u\leq P\) abbreviates \(u+P=P\).

\par\smallskip\noindent\textit{Values and exact original-witness criterion.}\quad 
Addition is 00004/\allowbreak{}01004/\allowbreak{}00224/\allowbreak{}00234/\allowbreak{}44444, with \(3<2<0<4\) and \(1<0<4\);1 and 2 are incomparable. Multiplication is 00000/\allowbreak{}02000/\allowbreak{}00000/\allowbreak{}00020/\allowbreak{}00004. A nonlinear word is 4 on all 4 content. Otherwise a quadratic is 2 if both occurrences are 1 or both 3; all other nonlinear words are 0. Long words retain their content.

For original P let L be its literals, E its undirected quadratic edges, V the edge endpoints, and S the original loops. A component of E is active if it contains a variable in L. Only V are graph vertices. Let G mean that an original long word exists and let C(w) mean an original word has content contained in c(w). Absorption \(w\leq P\) holds exactly when:

\textbullet{} C(w) holds, and G holds; or

\textbullet{} no original long word exists, w is literal, and it is an original literal or belongs to an active component; additionally C(w) holds; or

\textbullet{} no original long word exists, \(w=ab\), C(w) holds, a, b belong to V, and they are in the same component or both their components are active.

Thus a missing literal without G needs an original square as well as activation; activation alone cannot supply the contained witness. A new edge always needs a contained original loop or literal at an endpoint.

If C fails, assign target content 4 and every outsider 3. The target is 4 and each right summand is at most 0. If G fails, assign V to 3 and outsiders 2. Each right summand is 2 or 3, hence at most 2. A long target and any quadratic with an endpoint outside V have value 0. For a missing literal x satisfying C but in an inactive component, assign its component 1, other components3 and outsiders 2. Every right edge is 2, and no original literal lies in the inactive component; right literals are 3 or 2. The target is 1, not below 2. For endpoints in different components with at least one inactive, give that inactive component 1, all other components3, outsiders 2. The right side is 2 but the quadratic target is 0. These separate every failed condition without a variable or length bound.

For sufficiency, G makes \(P\geq 0\) and the contained witness supplies 4 whenever the target is 4. Suppose G fails and P does not dominate 0. Every original edge must be 2, so each component is monochromatic1 or 3. An active component must be 3: an original literal 1 joined with 2 would give 0. A quadratic within a component or between two active components is therefore 2. For a missing literal, C gives an original loop and the active component forces its value 3, which is below 2. Contained words handle target 4. Original literals are immediate.

\par\smallskip\noindent\textit{Unrestricted complete derivations.}\quad 
The word laws normalize long words by content. Commute a chosen occurrence to the end, duplicate it with two nonempty prefix blocks using \(xyz=xyzz\), then reverse the law to remove surplus copies while retaining length at least three. Canonical forms are sorted content for at least three variables, xxy for two, and xxx for one. Long g satisfies \(g^{3}=g\); no analogous assertion is made about quadratic or literal words.

If G holds, fix original long g and original contained u. Upgrade a literal u to \(u^{3}\) by \(xxx\leq x+yyy\) with \(y=g\), or a quadratic ab to aab by \(xxy\leq xy+zzz\) with \(z=g\). Retain long u. Append missing target content by \(xyzt\leq xyz\), normalize, and recover quadratics by \(xy\leq xxy\) or literals by \(x\leq xxx\). All blocks and substitution words are nonempty.

Without G, a missing literal h has an original loop hh and a path to an original literal z. With hh retained, extend the edge from h along the original path using \(xz\leq xx+xy+yz\) until hz is obtained. Then \(x\leq xx+xy+y\) recovers h from hh, hz, z. Intermediate vertices need not and in general cannot be recovered as literals.

A new quadratic ht has a contained original endpoint loop or literal, say h. If h is only an original literal, recover hh using \(xx\leq x+xy\) with an original incident edge. If h, t share a component, retain hh and propagate along an original path as above. This includes recovering hh itself when the target is square. If the components differ, both are active. Recover literal h using its loop and original activation witness as above. Choose an original literal z in t's component. It has an original incident edge, so recover zz; retain zz while deriving zt along an original path. Finally \(xy\leq xx+x+z+yz\) derives ht from hh, h, z, tz. No unsupported literal t or square tt is needed.

All eleven laws satisfy the exact criterion and are independently checked over their finite assignments. The derivations use only original or explicitly derived summands, with no bound on variables, polynomial sizes or path lengths. Distributive expansion and summandwise absorption in both directions prove completeness. Further types require independent HSP checks in both directions before classification.

\subsection{Two regularization classes with finite identity bases}\label{case:71}
All identities use the constant-free signature (+, multiplication), relative to the additively idempotent semiring axioms. Write u \ensuremath{\leq} P for \(u+P=P\). Finite verification supplements the arbitrary-variable proofs below.

\par\smallskip\noindent\textit{1. A four-law basis for a join of two small varieties.}\quad 
Let N=\{a, b, z\}, with \(z<b<a\), addition maximum, \(a*a=b\), and every other product z. Thus every word of length at least three is z. Let L=\{0,\allowbreak{}1\} with addition maximum and multiplication minimum. Let \(U=Var(N)\) join Var(L). Neither element notation introduces constants into terms.

Theorem 1. U has the following four-law basis:

{\small\par\smallskip
\CaseLine{\(xy = yx;\)}
\CaseLine{\(xy \leq  x;\)}
\CaseLine{\(xy \leq  x^{2};\)}
\CaseLine{\(xyz = x^{2}yz.\)}
\smallskip}

All four laws hold in N and L. A nonempty polynomial P is a finite sum of nonempty words. For a target word w put \(C=c(w)\). Absorption \(w\leq P\) holds in both N and L precisely as follows:

\textbullet{} If \(w=x\) has length one, P contains the linear word x.

\textbullet{} If \(w=x^{2}\), P contains x or \(x^{2}\).

\textbullet{} If \(w=xy\) with x different from y, P contains a word of length at most two with content contained in \{x, y\}.

\textbullet{} If w has length at least three, P contains a word whose content is contained in C.

For the first assertion, assign \(x=a\) and every other variable z in N. No right word except the linear word x attains a. For the second, use the same assignment; only x and \(x^{2}\) can reach b. For the third, assign \(x=y=a\) and outsiders z in N. A nonzero right value must come from a contained word of length one or two. Conversely each such word absorbs xy: multiplication by any factor decreases values, and \(xy\leq x^{2}\), \(y^{2}\). For the last assertion N imposes no condition, since \(w=z\). In L assign the variables of C to 1 and all outsiders to 0. The right side is 1 exactly when it has a contained word. Such a word absorbs w for every L-assignment. This proves the criterion.

Here is an equational completeness proof. Commutativity and \(xyz=x^{2}yz\) permit increasing any chosen multiplicity in a word of length at least three. Substituting \(x=y=z\) also gives \(x^{3}=x^{4}\). Consequently every word of length at least three with content C is equivalent to \(B_C=\prod_{x\in C}x^3\): first raise each present multiplicity to at least three, then cap it at three. Each increase isolates a selected letter and two nonempty remaining factors; no empty substitution occurs. If a right word u has content contained in C, append enough target letters to obtain a word q with content C and length at least three. Repeated use of \(xy\leq x\) gives \(q\leq u\); normalization gives \(w=q\leq u\). The short targets follow directly from the first three laws and the displayed criterion. Distribute arbitrary terms and apply the absorption argument to every summand in both directions. This proves Theorem 1 for any finite number of variables and summands.

\par\smallskip\noindent\textit{2. A five-law basis for the four-element truncated chain.}\quad 
Let T=\{v\_0, \(v_{1}\), \(v_{2}\), \(v_{3}\)\}, with \(v_{3}<v_{2}<v_{1}<v_{0}\), addition maximum, and \(v_i v_j=v_{\min(3,i+j)}\).

Theorem 2. T has the following five-law basis:

{\small\par\smallskip
\CaseLine{\(xy = yx;\)}
\CaseLine{\(x^{3} = x^{4};\)}
\CaseLine{\(xy \leq  x;\)}
\CaseLine{\(xyz \leq  x^{2}yz+xy^{2}z+xyz^{2};\)}
\CaseLine{\(xy \leq  x^{2}+xy^{3}.\)}
\smallskip}

Soundness follows from adding the indices and truncating at three. For the fourth law, if all three indices are positive the target is \(v_{3}\); otherwise duplicating a zero-index factor preserves its value. For the fifth, if x has index zero then \(x^{2}=v_{0}\); if y has index zero then \(xy^{3}=xy\); if both are positive, x of index one gives \(x^{2}=v_{2}\geq xy\), and x of larger index makes \(xy=v_{3}\).

Normalize a word by sorting its letters and capping positive multiplicities at three. For w with content C, discard right words containing outsiders by assigning those variables index three. Write \(P_{C}\) for the remaining words. The exact condition is that \(P_{C}\) is nonempty and, for all \(lambda_{x}\) in \{0,\allowbreak{}1,\allowbreak{}2,\allowbreak{}3\},

{\small\par\smallskip
\CaseLine{\(min(3, \sum  n_{x}(w) lambda_{x})\)}
\CaseLine{\(\geq  min_\{u in P_{C}\} min(3, \sum  n_{x}(u) lambda_{x}).\)}
\smallskip}

Necessity uses core assignments with outsiders of index three. For sufficiency a witness among the contained words already bounds the target, regardless of the outsiders' actual values. Only tests giving the target index below three matter: zero positive variables, one positive variable of weight one or two, or two of weight one.

Insertion into any nonempty word follows from \(xy\leq x\). Thus a target is derivably below any right word that divides it as a commutative word. If a normalized target has at least three occurrences of unsaturated variables (multiplicity below three), apply the fourth law to three such occurrences in the remaining context. Each branch raises one selected multiplicity and reduces the sum of deficits. The index condition is preserved. Hence one reaches a terminal s, sx, \(sx^{2}\) or sxy, where s is the product of cubes of saturated variables and x, y are distinct unsaturated variables. An omitted s is only an omitted context; a term is never empty.

For s, any contained right word divides it. For sx, give x index two; the condition supplies a divisor with x-exponent at most one. For \(sx^{2}\) give x index one to obtain exponent at most two. All other target variables in these tests have index zero. For sxy, a divisor finishes the argument. Otherwise assigning \(x=y\) index one gives a witness u with exponents (2,\allowbreak{}0) or (0,\allowbreak{}2) on x, y. Interchange x, y if necessary. Assigning x index two and y index zero gives a witness v with x-exponent at most one. Insertion and normalization yield

{\small\par\smallskip
\CaseLine{\(sx^{2} \leq  u,\qquad sxy^{3} \leq  v.\)}
\smallskip}

The fifth law in context s yields \(sxy\leq u+v\leq P\). This also works when s is absent and when a witness is linear. Deficit induction and summandwise absorption prove completeness for arbitrary variables.

\par\smallskip\noindent\textit{3. Finite regularization bases.}\quad 
Put d(x, y)=x+xy. Both U and Var(T) satisfy d(x, y)=x. For either displayed basis Sigma take the following seven partition identities:

{\small\par\smallskip
\CaseLine{\(d(x, x)=x;\)}
\CaseLine{\(d(d(x, y),z)=d(x, d(y, z));\)}
\CaseLine{\(d(x, d(y, z))=d(x, d(z, y));\)}
\CaseLine{\(d(f(x,y),z)=f(d(x,z),d(y,z))\qquad (f\in\{+,\cdot\});\)}
\CaseLine{\(d(z,f(x,y))=d(d(z,x),y)\qquad(f\in\{+,\cdot\}).\)}
\smallskip}

Add \(d(s,t)=d(t,s)\) for each identity \(s=t\) in \(\Sigma\). An inequality \(u\leq P\) is first written as \(u+P=P\). All occurrences of d are abbreviations for the ordinary term x+xy. This gives at most 11 identities for R(U), and at most 12 for R(Var(T)), relative to the semiring axioms. No minimality is claimed.

For a nonempty finite set \(C\), let \(D_C(u)\) be the result of applying \(d\) to \(u\) and then successively to the letters of \(C\). The partition identities allow permutation and removal of repeated padding letters and distribute a common padding through either basic operation. They give \(d(u,v)=D_{\ct(v)}(u)\) and \(D_{\ct(u)}(u)=u\), by induction on terms. If \(s=t\) is a basis identity, pad the identity \(d(s\sigma,t\sigma)=d(t\sigma,s\sigma)\) by \(C\) containing both supports. This gives \(D_C(s\sigma)=D_C(t\sigma)\), and the common padding can be carried through any term context. For a regular identity \(u=v\) of the original variety take \(C=\ct(u)=\ct(v)\) and replace auxiliary letters in a deduction by a fixed member of \(C\). The padded steps derive \(D_C(u)=D_C(v)\), hence \(u=v\). Conversely every partition law and every padded basis law is regular and holds in the original variety. This proves completeness for its regularization.

\par\smallskip\noindent\textit{4. Identification of the five-element sources.}\quad 
Catalog 7635 is exactly T with a new common absorbing element infinity for both operations: (\(v_{0}\), \(v_{1}\), \(v_{2}\), \(v_{3}\), \(infinity)=(0\),1,\allowbreak{}2,\allowbreak{}3,\allowbreak{}4). Its identities are exactly the regular identities of T: an assignment of infinity on a variable occurring on only one side separates any nonregular identity; regular identities evaluate to infinity on both sides whenever infinity occurs. Thus 7635 has the twelve-law basis.

Catalog 1777 has elements \(z<b<a<d<e\), addition maximum, with \((a,b,z,d,e)=(1,2,3,0,4)\). Its N-component \{a, b, z\} has the multiplication of N above; its \{d, e\}-component is the two-element lattice L, with \(d<e\). Every product of a lower-component element and an upper-component element equals d. Thus every term using at least one upper-component variable evaluates in the upper component, and its evaluation equals the evaluation in L after mapping all lower-component assigned values to d. This assertion follows directly by induction on the term: both operations on the upper component preserve d as required, and the constant map from N to d is a homomorphism.

Consequently every regular identity valid in N and L holds in 1777. Conversely its subalgebras N and L force validity in both, and a nonregular identity fails by giving a variable missing on one side the value e and every other variable z: one side is in the upper component and the other in the lower. Hence Id(1777) is precisely the regular part of Id(U), and 1777 has the eleven-law basis.

Exact two-way variety certificates identify the following classes:

\noindent Source\enspace\textbar\enspace Catalog types\enspace\textbar\enspace Basis-size upper bound\par
\noindent 1777\enspace\textbar\enspace 1777, 9217, 9457, 9670, 10065, 10068\enspace\textbar\enspace 11\par
\noindent 7635\enspace\textbar\enspace 7635, 10033, 10045\enspace\textbar\enspace 12\par

\subsection{首字母交换与平方保护：八组完整有限基}\label{case:72}
\par\smallskip\noindent\textit{1. 八个具体构造.}\quad 
底集用 \(\{z, b, a, d, t\}\) 标记，加法取下表所列偏序的并。没有列出的关系仅作传递闭包。

\noindent 名称\enspace\textbar\enspace 编号\enspace\textbar\enspace 加法偏序\enspace\textbar\enspace 符号顺序对应的目录标号\par
\noindent \(A\)\enspace\textbar\enspace 13057\enspace\textbar\enspace \(z<b<a<d<t\)\enspace\textbar\enspace \((3,\allowbreak{}2,\allowbreak{}1,\allowbreak{}0,\allowbreak{}4)\)\par
\noindent \(B\)\enspace\textbar\enspace 13828\enspace\textbar\enspace \(z<b<a<d<t\)\enspace\textbar\enspace \((3,\allowbreak{}2,\allowbreak{}1,\allowbreak{}0,\allowbreak{}4)\)\par
\noindent \(C\)\enspace\textbar\enspace 13065\enspace\textbar\enspace \(z<b<a<d<t\)\enspace\textbar\enspace \((3,\allowbreak{}2,\allowbreak{}1,\allowbreak{}0,\allowbreak{}4)\)\par
\noindent \(D\)\enspace\textbar\enspace 4349\enspace\textbar\enspace \(z<b<d<a<t\)\enspace\textbar\enspace \((3,\allowbreak{}2,\allowbreak{}0,\allowbreak{}1,\allowbreak{}4)\)\par
\noindent \(E\)\enspace\textbar\enspace 4348\enspace\textbar\enspace \(z<b<a<t,\ z<d<a<t\)\enspace\textbar\enspace \((3,\allowbreak{}2,\allowbreak{}0,\allowbreak{}1,\allowbreak{}4)\)\par
\noindent \(F\)\enspace\textbar\enspace 4239\enspace\textbar\enspace \(z<b<a<t,\ z<d<a<t\)\enspace\textbar\enspace \((3,\allowbreak{}2,\allowbreak{}0,\allowbreak{}1,\allowbreak{}4)\)\par
\noindent \(G\)\enspace\textbar\enspace 12038\enspace\textbar\enspace \(z<b<a<t,\ z<d<a<t\)\enspace\textbar\enspace \((3,\allowbreak{}1,\allowbreak{}0,\allowbreak{}2,\allowbreak{}4)\)\par
\noindent \(H\)\enspace\textbar\enspace 1888\enspace\textbar\enspace \(z<b, a, d<t\)，三个中间元素两两不可比\enspace\textbar\enspace \((1,\allowbreak{}0,\allowbreak{}3,\allowbreak{}2,\allowbreak{}4)\)\par
前七个构造通过加法回缩 \(\rho\) 定义 \(xy=\rho(x)\rho(y)\)。\(A, B, C\) 取 \(\rho(d)=t\)；\(D, E, F\) 取 \(\rho(d)=a\)；\(G\) 取 \(\rho(a)=\rho(d)=t\)；其余点固定。它们的乘法核心如下。

\textbullet{} \(A\)：核心 \(\{z, b, a, t\}\) 中 \(b\) 是乘法零；\(z^2=z\)；\(a, t\) 构成左零半群；\(z\) 与 \(a, t\) 的两个方向混合积均为 \(b\)。

\textbullet{} \(B\)：\(z\) 是乘法零；\(a, t\) 在非零元素上作左零元素；\(b^2=b\)，\(ba=bt=a\)。

\textbullet{} \(C, D, E\)：\(z, b\) 为左零元素；\(\{a, t\}\) 上乘法是 \(a<t\) 的最大值；正层元素右乘 \(z, b\) 均为 \(b\)。

\textbullet{} \(F\)：同上一项，但 \(\{a, t\}\) 的乘法改为左零乘法。

\textbullet{} \(G\)：核心 \(\{z, b, t\}\) 中 \(z, b\) 为左零元素，\(t^2=t\)、\(tz=tb=b\)。

\(H\) 的核心 \(\{z, b, a, t\}\) 是布尔格，依次对应 \((0,\allowbreak{}0),(1,\allowbreak{}0),(0,\allowbreak{}1),(1,\allowbreak{}1)\)；乘法为 \((i, j)(k, l)=(i, j\wedge l)\)，取 \(\rho(d)=t\)。

这些核心均为带，回缩保持加法；核心左右平移保持并，因而结合律与分配律传到扩张。以下完整词值描述也可直接核对这些断言。对非一次词 \(w\) 记 \(X=c(w),h=h(w)\)，所有``字母''都指赋值后的元素：

\textbullet{} \(A\)：全体字母为 \(z\) 时取 \(z\)；全体字母在 \(\{a, d, t\}\) 时，首字母为 \(a\) 则取 \(a\)，否则取 \(t\)；其他情况均取 \(b\)。

\textbullet{} \(B\)：出现 \(z\) 则取 \(z\)；否则首字母为 \(d, t\) 则取 \(t\)，为 \(a\) 则取 \(a\)；首字母为 \(b\) 时，全体字母为 \(b\) 则取 \(b\)，否则取 \(a\)。

\textbullet{} \(C, D, E\)：首字母为 \(z, b\) 时分别取 \(z, b\)；首字母在正层时，若任一字母为 \(z, b\) 则取 \(b\)；否则取全体回缩像在 \(a<t\) 中的最大值。

\textbullet{} \(F\)：首字母为 \(z, b\) 时分别取 \(z, b\)；正层首字母遇到任一低层字母则取 \(b\)，否则取首字母的回缩像。

\textbullet{} \(G\)：首字母为 \(z, b\) 时分别取 \(z, b\)；其他首字母在所有字母均属 \(\{a, d, t\}\) 时取 \(t\)，否则取 \(b\)。

\textbullet{} \(H\)：记 \(\rho(x)=(\alpha(x),\beta(x))\)，则 \(w=(\alpha(h),\bigwedge_{x\in X}\beta(x))\)。

\par\smallskip\noindent\textit{2. 吸收判据.}\quad 
对非空右多项式 \(p\)，记一次项集为 \(L\)，非一次项族为 \(\mathcal R\)，并令

\CaseDisplay{T=c(p),\quad \mathcal R_X=\{u\in\mathcal R:c(u)\subseteq X\},\quad N_X=\bigcup_{u\in\mathcal R_X}c(u),}

\CaseDisplay{G_0=\{h(u):u\in\mathcal R\},\quad G_X=\{h(u):u\in\mathcal R_X\}.}

判据定理。 非一次词 \(w\)（\(h=h(w),X=c(w)\)）满足 \(w\le p\) 当且仅当：

\noindent 来源\enspace\textbar\enspace 充分必要条件\par
\noindent \(A\)\enspace\textbar\enspace \(X\subseteq T\) 且 \(h\in G_X\)\par
\noindent \(B\)\enspace\textbar\enspace \(X\subseteq L\cup N_X\) 且 \(h\in G_X\)\par
\noindent \(C\)\enspace\textbar\enspace \(X\subseteq N_X\) 且 \(h\in L\cup G_0\)\par
\noindent \(D\)\enspace\textbar\enspace \(\mathcal R_X\ne\varnothing\)，\(X\subseteq L\cup N_X\)，\(h\in L\cup G_0\)\par
\noindent \(E\)\enspace\textbar\enspace \(X\subseteq L\cup N_X\) 且 \(h\in G_0\)\par
\noindent \(F\)\enspace\textbar\enspace \(h\in G_0\) 且 \(h\in L\cup G_X\)\par
\noindent \(G\)\enspace\textbar\enspace \(h\in G_0\) 且 \(\mathcal R_X\ne\varnothing\)\par
\noindent \(H\)\enspace\textbar\enspace \(h\in G_0\) 且 \([\mathcal R_X\ne\varnothing\text{ 或 }X\cap L\cap G_0\ne\varnothing]\)\par
在所有八个来源中，一次词 \(x\le p\) 当且仅当 \(x\in L\)，或存在仅含 \(x\) 的非一次右词。

\par\smallskip\noindent\textit{2.1 必要性.}\quad 
\(A\)。 内容遗漏 \(x\) 时，赋 \(x=a\)、其他全为 \(z\)；右值为 \(z\)，左值为 \(a\) 或 \(b\)。若 \(h\notin G_X\)，赋 \(h=d\)、其他 \(X\) 内变量为 \(a\)、外部为 \(z\)；左值为 \(t\)，右值至多为 \(d\)。

\(B\)。 若 \(x\notin L\cup N_X\)，赋 \(x=a\)、其他 \(X\) 内变量为 \(b\)、外部为 \(z\)，左值为 \(a\)、右值至多为 \(b\)。若 \(h\notin G_X\)，赋 \(h=d\)、其他 \(X\) 内变量为 \(a\)、外部为 \(z\)；左值为 \(t\)、右值至多为 \(d\)。

\(C\)。 若 \(x\in X\setminus N_X\)，赋 \(x=d\)、其他 \(X\) 内变量为 \(a\)、外部为 \(z\)；左值为 \(t\)，内部非一次右词取 \(a\)，外部非一次右词至多为 \(b\)，一次项至多为 \(d<t\)。若 \(h\notin L\cup G_0\)，赋 \(h=b\)、其他全为 \(z\)，左值为 \(b\)、右值为 \(z\)。

\(D\)。 无内部非一次项时，赋 \(X=d\)、外部为 \(z\)，左值为 \(a\)、右值至多为 \(d\)。若 \(x\notin L\cup N_X\)，赋 \(x=t\)、其他 \(X\) 内变量为 \(d\)、外部为 \(z\)，左值为 \(t\)、右值至多为 \(a\)。首字母条件的反例赋值与 \(C\) 相同。

\(E\)。 内容条件失败时，使用 \(D\) 的 \(x=t\) 赋值：右侧即使把 \(b, d\) 求并也仅得 \(a<t\)。若 \(h\notin G_0\)，赋 \(h=d\)、其他全为 \(z\)；左值为 \(b\) 或 \(a\)，右值至多为 \(d\)，而 \(b\nleq d\)、\(a\nleq d\)。

\(F\)。 \(h\notin G_0\) 时使用 \(E\) 的赋值。若 \(h\notin L\cup G_X\)，赋 \(h=t\)、其他 \(X\) 内变量为 \(d\)、外部为 \(z\)；左值为 \(t\)、右值至多为 \(a\)。

\(G\)。 \(h\notin G_0\) 时赋 \(h=d\)、其他全为 \(z\)，左值为 \(b\) 或 \(t\)，右值至多为 \(d\)。无内部非一次项时，赋 \(X=d\)、外部为 \(z\)，左值为 \(t\)，右值至多为 \(b+d=a<t\)。

\(H\)。 若 \(h\notin G_0\)，赋 \(h=d\)、其他全为 \(z\)，左值为 \(b\) 或 \(t\)，右值至多为 \(d\)。若第二条件失败而 \(h\in G_0\)，令 \(X\cap L\) 内变量取 \(d\)，其余 \(X\) 内变量取 \(a\)，外部变量取 \(z\)。因 \(X\cap L\cap G_0\) 为空，所有来自 \(X\) 的非一次右项首字母都取 \(a\)；又无内部非一次右项，故它们因含外部 \(z\) 而取 \(z\)。外部首字母也给 \(z\)。左首字母属于 \(G_0\) 而不在 \(L\)，故左值为 \(a\)；右值为 \(z\) 或 \(d\)，不能吸收 \(a\)。

一次词的必要性统一取 \(x=d\)、其他变量为 \(z\)。若没有同名一次项和仅含 \(x\) 的非一次项，右值分别为：\(A, C, D\) 中至多 \(b\)；\(B\) 中 \(z\)；\(E, F, G, H\) 中至多 \(b\)。这些值都不在 \(d\) 之上。

\par\smallskip\noindent\textit{2.2 充分性.}\quad 
\(A\)。 全部字母为 \(z\) 时左值最小。左值为 \(b\) 时，至少一个 \(X\) 内字母非 \(z\)，内容覆盖给出含它的右项，值至少为 \(b\)。全体字母在正层时，内部同首字母非一次项与左词有相同的值。

\(B\)。 非零左值时，\(X\) 内没有 \(z\)。内部同首字母右词提供至少为 \(b\) 的值，并在左值为 \(t\) 时直接给出 \(t\)。左值为 \(a\) 时，\(X\) 内有高值变量；它的一次项或含它的内部非一次项由内容条件提供，值至少为 \(a\)。此处链上 \(a<d\)，故一次项 \(d\) 足够。

\(C\)。 左值为 \(z\) 时自明；为 \(b\) 时，同首字母右项至少为 \(b\)。左值为 \(a\) 时，全体 \(X\) 在正层，内部非一次项存在，至少为 \(a\)。左值为 \(t\) 时，取一个高值变量，由 \(X\subseteq N_X\) 得到含它的内部非一次项，值为 \(t\)。

\(D\)。 左值为 \(z, b\) 的论证同 \(C\)。左值为 \(a\) 时，内部非一次项提供至少为 \(a\) 的值。左值为 \(t\) 时，高值变量的一次项或内部非一次项直接给出 \(t\)。

\(E\)。 左值为 \(z\) 时自明；为 \(b\) 时，同首字母非一次右词至少为 \(b\)。左值为 \(t\) 时，内容条件给出高值变量的一次项或内部非一次项，值为 \(t\)。若左值为 \(a\)，全体 \(X\) 在正层；内部非一次项如存在则至少为 \(a\)。若不存在，内容条件保证 \(X\subseteq L\)，特别地首字母一次项存在：取 \(a, t\) 时已足够，取 \(d\) 时与同首字母非一次右词的 \(a, t\) 或 \(b\) 求并也至少为 \(a\)。

\(F\)。 左值为 \(b\) 时，同首字母非一次项至少为 \(b\)。左值为 \(a\) 时，内部同首字母非一次项直接取 \(a\)；若只用同名一次项，它取 \(a\) 时完成，取 \(d\) 时与同首字母全局非一次项的 \(a\) 或 \(b\) 求并为 \(a\)。左值为 \(t\) 时，首字母为 \(t\)，局部首字母条件直接提供 \(t\)。

\(G\)。 左值为 \(b\) 时，同首字母非一次项至少为 \(b\)；左值为 \(t\) 时，内部非一次项取 \(t\)。

\(H\)。 左词第一坐标为一时，同首字母非一次项的值为 \(b\) 或 \(t\)，故右侧在 \(b\) 之上。左词第二坐标为一时，内部非一次项如存在则在 \(a\) 之上。否则选 \(g\in X\cap L\cap G_0\)：它的一次项若取 \(a, t\) 则已足够；若取 \(d\)，同首字母 \(g\) 的右非一次项取 \(b\) 或 \(t\)，与 \(d\) 求并得到 \(t\)。因此左右坐标所要求的 \(a, b\) 均被覆盖；左值为 \(t=a+b\) 时也成立。

一次词的充分性由 \(x\le x^2=x^n\)（\(n\ge2\)）及同名项给出。

\par\smallskip\noindent\textit{3. 有限基与形式完备性.}\quad 
共同保留

\CaseDisplay{xyz\approx xzy\quad(N),\quad xy\approx x^2y^2\quad(S),\quad x\le x^2\quad(I).}

另外使用

\CaseDisplay{\begin{aligned}
xy&\le x^2+y &(K),\\
xy&\le x^2+y^2 &(J),\\
xy&\le x+yx &(H_1),\\
xy&\le x^2+zy &(B_1),\\
xy&\le yx+xz &(G_1),\\
x^2&\le x+xy &(Q_1),\\
xy&\le x^2 &(C_1),\\
xy&\le y^2+xyz &(E_2).
\end{aligned}}

定理。 在共同三式之外，下表给出完整基：

\noindent 来源\enspace\textbar\enspace 另外的式子\enspace\textbar\enspace 总附加式数\par
\noindent \(A\)\enspace\textbar\enspace \(K, B_1\)\enspace\textbar\enspace 5\par
\noindent \(B\)\enspace\textbar\enspace \(K\)\enspace\textbar\enspace 4\par
\noindent \(C\)\enspace\textbar\enspace \(J, H_1, G_1\)\enspace\textbar\enspace 6\par
\noindent \(D\)\enspace\textbar\enspace \(K, H_1, G_1\)\enspace\textbar\enspace 6\par
\noindent \(E\)\enspace\textbar\enspace \(K, G_1, Q_1\)\enspace\textbar\enspace 6\par
\noindent \(F\)\enspace\textbar\enspace \(C_1, Q_1\)\enspace\textbar\enspace 5\par
\noindent \(G\)\enspace\textbar\enspace \(C_1, E_2\)\enspace\textbar\enspace 5\par
\noindent \(H\)\enspace\textbar\enspace \(C_1, E_2, Q_1\)\enspace\textbar\enspace 6\par
有效性由第 2 节判据逐式代入即得，\((N),(S)\) 则直接来自词值规则。以下给出任意变量数的完整推导。

\((N),(S)\) 把每个非一次词化为 \(W(h, X)\)，其中各字母平方、首字母为 \(h\)、内容为 \(X\)；有 \(W(h, X)^2=W(h, X)\)。\((K)\) 允许从非一次积开始乘任意右项而仍不超过右多项式；\((J)\) 允许乘非一次右项。\((B_1)\) 允许从非一次 \(P\le p\) 和含 \(x\) 的右非一次项 \(u\) 推出 \(Px\le P^2+ux=P+u\le p\)，而不在左边引入 \(u\) 的其他字母。

另有两个关键推导。若内部非一次词 \(P\le p\) 包含 \(h\)，而右词 \(r=hs\) 以 \(h\) 开始，则

\CaseDisplay{hP\le Ph+hs=P+r\le p\qquad(G_1).\quad\textup{(\ensuremath{3.1})}}

若 \(h\in L\) 且右非一次词 \(r=hs\)，则

\CaseDisplay{h^2\le h+hs=h+r\le p\qquad(Q_1).\quad\textup{(\ensuremath{3.2})}}

\(A\)。 从内部同首字母非一次项开始。每个缺少的 \(X\) 内变量，如出现于一次项，用 \((K)\) 补乘；否则用 \((B_1)\) 的无外部字母补全。内容条件保证补全所有 \(X\)，得到 \(w\le p\)。

\(B\)。 从内部同首字母非一次项开始，用 \((K)\) 乘上所有内部非一次项和 \(L\cap X\) 中的变量。内容条件保证最终恰为 \(W(h, X)\)。

\(C\)。 用 \((J)\) 把全部内部非一次项相乘得 \(P\le p\)，其内容为 \(X\)。若首字母见证是一次项 \(h\)，用 \((H_1)\) 得 \(hP\le h+Ph=h+P\le p\)。若见证是右非一次词 \(hs\)，用 (3.1)。

\(D\)。 从任意内部非一次项开始，用 \((K)\) 把所有内部非一次项和内部一次项乘入，得内容为 \(X\) 的非一次积 \(P\le p\)。再按 \(C\) 的两个首字母分支处理。

\(E\)。 选全局同首字母右词 \(r=hs\)。若 \(h\in L\)，用 (3.2) 从 \(h^2\) 开始。否则内容条件给出内部非一次项 \(q\) 含 \(h\)，用 (3.1) 从 \(hq\le q+r\le p\) 开始。之后用 \((K)\) 乘入全部内部右项，补全 \(X\)。

\(F\)。 若有内部同首字母非一次项，从它开始，用 \((C_1)\) 补乘缺少的 \(X\) 内变量。否则 \(h\in L\)，全局非一次首字母见证允许使用 (3.2)，再从 \(h^2\) 作相同补乘。

\(G\)。 选全局同首字母非一次词 \(r\) 和内部非一次项 \(q\)。由 \((C_1)\) 补乘，得 \(W(h, X\cup c(r))\le r\le p\)。若 \(c(r)\subseteq X\) 已完成；否则令 \(s\) 是 \(c(r)\setminus X\) 的非空字母积，以 \((x, y, z)=(w, q, s)\) 代入 \((E_2)\)，得

\CaseDisplay{w=wq\le q^2+wqs=q+W(h, X\cup c(r))\le p.}

\(H\)。 有内部非一次项时同 \(G\)。否则选 \(g\in X\cap L\cap G_0\)，其非一次首字母见证和 \((Q_1)\) 给出 \(g^2\le p\)。把这个已推出的内部非一次词作为锚点 \(q\)，重复 \(G\) 的推导；并不要求它原来就是右求和项。

一次词由 \((I)\) 处理。所有补乘和消元只在缺少变量时使用非空词，否则省略。任意有效恒等式拆为双向逐词吸收，即得完整性。乘法对偶逐词反转这些证明即可。

\par\smallskip\noindent\textit{4. 范围限制.}\quad 
这些是指定生成元及其对偶的完整有限基。对目录对象仅在独立验证精确双向簇包含后传递；满足一组基的任意真子簇不自动具有有限基。有界计算只作诊断，不代替上述任意变量数证明。

\subsection{9316型的三条完整恒等式基及其对偶}\label{case:73}
所有恒等式均在不含常元的语言 (+,\ensuremath{\cdot}) 中，相对于加法幂等半环公理。记 p\ensuremath{\leq}q 为 p+q\ensuremath{\approx}q。以下证明不限制变量数。

\par\smallskip\noindent\textit{代数及结论.}\quad 
令 A=\{0,\allowbreak{}1,\allowbreak{}2,\allowbreak{}3,\allowbreak{}4\}，加法是两条链 \(1<0<4\)、\(3<2<4\) 生成的并运算，两个分支的元素之和为4。乘法为

\noindent \ensuremath{\cdot}\enspace\textbar\enspace 0\enspace\textbar\enspace 1\enspace\textbar\enspace 2\enspace\textbar\enspace 3\enspace\textbar\enspace 4\par
\noindent 0\enspace\textbar\enspace 1\enspace\textbar\enspace 1\enspace\textbar\enspace 1\enspace\textbar\enspace 1\enspace\textbar\enspace 1\par
\noindent 1\enspace\textbar\enspace 1\enspace\textbar\enspace 1\enspace\textbar\enspace 1\enspace\textbar\enspace 1\enspace\textbar\enspace 1\par
\noindent 2\enspace\textbar\enspace 3\enspace\textbar\enspace 3\enspace\textbar\enspace 4\enspace\textbar\enspace 4\enspace\textbar\enspace 4\par
\noindent 3\enspace\textbar\enspace 3\enspace\textbar\enspace 3\enspace\textbar\enspace 3\enspace\textbar\enspace 3\enspace\textbar\enspace 3\par
\noindent 4\enspace\textbar\enspace 4\enspace\textbar\enspace 4\enspace\textbar\enspace 4\enspace\textbar\enspace 4\enspace\textbar\enspace 4\par
元素名称不是语言中的常元。这正是五元素目录9316型。

定理。 A 的完整恒等式基可取下面三式：

\CaseDisplay{xyz\approx xy,\quad\textup{(\ensuremath{B1})}}

\CaseDisplay{xy\leq x^2,\quad\textup{(\ensuremath{B2})}}

\CaseDisplay{x^2\leq y+xy.\quad\textup{(\ensuremath{B3})}}

乘法对偶的完整基由逐词反转得到，亦只有三式。

\par\smallskip\noindent\textit{运算及恒等式的有效性.}\quad 
置 D=\{0,\allowbreak{}1\}。所有乘积在 \{1,\allowbreak{}3,\allowbreak{}4\} 中，这三个元素的行恒定等于其自身，因此 \((xy)z=xy\)。只有左因子为2时乘积依赖右因子，且只依赖右因子是否在 D；yz\ensuremath{\in}D 当且仅当 y\ensuremath{\in}D。因此 \(x(yz)=xy\)，证明结合律及(B1)。

各左平移保持加法：唯一非常值行由分支 D 取3、补集取4决定，而一个和在 D 中恰好两加数均在 D 中。对于右平移，两分支上的限制均单调，不同分支元素的像之和为4，且4的像为4。这给出两侧分配律；加法本身显然是半格。

平方取值依次为1、1、4、3、4。除 \(x=2\) 外，均有 \(xy=x\)\ensuremath{{}^{2}}。\(x=2\) 时 xy\ensuremath{\in}\{3,\allowbreak{}4\}，故(B2)成立。若 \(x=2\) 且 y\ensuremath{\in}D，则 \(xy=3\)、\(y+3=4=x\)\ensuremath{{}^{2}}；若 y\ensuremath{\notin}D，则 \(xy=x\)\ensuremath{{}^{2}}。其余 x 由 \(xy=x\)\ensuremath{{}^{2}} 即得(B3)。

在(B3)中以 yz 替换 y，再用(B1)，可推导

\CaseDisplay{x^2\leq yz+xy.\quad\textup{(\ensuremath{B4})}}

\par\smallskip\noindent\textit{任意变量的吸收判据.}\quad 
展开分配律，并用(B1)将每个长度至少二的词化为首两个字母。任意项因此可化为非空多项式

\CaseDisplay{P=\sum_{z\in L}z+\sum_{(u, v)\in E}uv.}

这里 L 是一次加数集合，E 是允许自环的有向边集合，空部分省略而不增加零常元。定义

\CaseDisplay{H(P)=L\cup\{u:(u, v)\in E\},\qquad
J_x(P)=\{v:(x, v)\in E\}\cap H(P).}

引理。 对任意有限变量集及任意 P，

\CaseDisplay{A\models x\leq P\quad\Longleftrightarrow\quad x\in L,}

\CaseDisplay{A\models xy\leq P\quad\Longleftrightarrow\quad
(x, y)\in E\ \text{或}\ J_x(P)\neq\varnothing.}

右侧条件成立时，吸收均可由(B1)---(B3)推导。

一次项。 若 x\ensuremath{\in}L 则结论显然。否则令 \(x=0\)，其他变量均为1。乘法限制在 D\ensuremath{\times}D 上恒等于1，所以所有二次加数取1，所有一次加数亦取1。\(P=1\) 而 \(x=0\) 不小于等于1，给出分离赋值。

二次项充分性。 直接列出的边当然被吸收。否则取 v\ensuremath{\in}\(J_{x}\)(P)，故 xv 是加数。若 v\ensuremath{\in}L，(B3)给出 \(x^{2}\)\ensuremath{\leq}v+xv\ensuremath{\leq}P；若 v 是某边 vu 的首字母，(B4)给出 \(x^{2}\)\ensuremath{\leq}vu+xv\ensuremath{\leq}P。再用(B2)，即可推导 xy\ensuremath{\leq}P。

二次项必要性。 假设 (x, y)\ensuremath{\notin}E 且 \(J_{x}\)(P)=\ensuremath{\varnothing}。令 U=\{v:(x, v)\ensuremath{\in}E\}，于是 U\ensuremath{\cap}H(P)=\ensuremath{\varnothing}，且 x\ensuremath{\notin}U：若 x\ensuremath{\in}U，则自环使 x\ensuremath{\in}\(J_{x}\)(P)。作如下赋值：

\CaseDisplay{x\mapsto2,\qquad U\mapsto0,\qquad
\text{其余变量}\mapsto3.}

所有 x 开头的边取3；其他边的首字母不在 U\ensuremath{\cup}\{x\}，也取3。一次加数仅可能取2（若是 x）或3，因为 U\ensuremath{\cap}L=\ensuremath{\varnothing}。故 P\ensuremath{\in}\{2,\allowbreak{}3\}。但 y\ensuremath{\notin}U，若 \(y=x\) 则 \(xy=4\)，若 y\ensuremath{\ne}x 则 \(y=3\)，仍有 \(xy=4\)。这与 xy\ensuremath{\leq}P 矛盾，必要性得证。

\par\smallskip\noindent\textit{完整性及适用范围.}\quad 
给定任意在 A 中成立的 p\ensuremath{\approx}q，先化为上述多项式 P, Q。P 的每个加数都在 A 中被 Q 吸收，由引理可从三式推导；相加得 P\ensuremath{\leq}Q。对称地得到 Q\ensuremath{\leq}P，故 P\ensuremath{\approx}Q 可导。这证明全部恒等式的有限基性，变量数、词长及加数数目没有上界限制。

逐词反转同时反转等式推导，从而得到 A 的乘法对偶的完整性。显式反转基为 xyz\ensuremath{\approx}yz、xy\ensuremath{\leq}\(y^{2}\)、\(y^{2}\)\ensuremath{\leq}x+xy。

定理~\ref{thm:correct643} 给出 \(643\) 型的四式基。对于 \(9316\) 型，取 \(x=0\) 即使 \(x\leq x^2\) 失效，故仅有二次边见证不足以吸收一次项。上文的一次项分离赋值说明了这一额外条件的必要性。

目录转移须逐型核验两方向完整HSP；有限变量计算仅为补充核验，不是本证明的依据。

\subsection{Finite bases from two protected initial positions}\label{case:74}
Signature (+, multiplication), no constants. We specify finite families, rather than claiming small enumerated bases.

Both additions have rows 00004,\allowbreak{}01114,\allowbreak{}01224,\allowbreak{}01234,\allowbreak{}44444, the chain \(3<2<1<0<4\). Multiplication rows for 2008 are 00000,\allowbreak{}11111,\allowbreak{}13331,\allowbreak{}33333, 00004; those for 3111 are 00030,\allowbreak{}22230,\allowbreak{}22232,\allowbreak{}33333,\allowbreak{}44434.

\par\smallskip\noindent\textit{Arbitrary word formulas and normal forms.}\quad 
For 2008 a head 0,\allowbreak{}1, or 3 gives that value forever. A head 4 gives 4 if every factor is 4, and 0 otherwise. A literal head 2 gives 2; if the word has a second position, its second factor in\{0,\allowbreak{}4\} gives 1 and a second factor in\{1,\allowbreak{}2,\allowbreak{}3\} gives 3. Later factors do not change that result. For 3111 any factor 3 makes the whole word 3. With no factor 3, a head 0,\allowbreak{}2, or 4 is retained; a literal head 1 gives 1; a nonlinear head 1 gives 0 when the second factor is 4 and 2 otherwise. Both formulas follow inductively from the displayed multiplication rows, for all lengths.

Both algebras satisfy the three structural identities

{\small\par\smallskip
\CaseLine{\(xyzt=xytz,\qquad xyx=xy,\qquad xyy=xy.\qquad (W)\)}
\smallskip}

They preserve a literal and, for a longer word, its first TWO positions and its total content. Conversely (W) normalizes every longer word to those first two positions followed by the increasing list of the remaining content, each variable once. The first law swaps adjacent positions after position two. Sort that suffix, use \(xyx=xy\) to delete a suffix copy of the head, and \(xyy=xy\) to delete other adjacent duplicates; a suffix copy of the second variable can first be moved next to position two. All prefixes substituted into (W) are nonempty. If the first two variables coincide both positions are retained. The common multiplicity two is never compressed into a single occurrence.

\par\smallskip\noindent\textit{Exact original-word criteria.}\quad 
Let P be a nonempty collection of original upper words. Write \(C=c(w)\), \(h=w[0]\), and \(s=w[1]\) when w is nonlinear. A literal target is valid exactly when it is already an original literal.

For a nonlinear target in 2008, the criterion is:

1. Some original word has content contained in C.

2. Some original word has head h; it may involve outside variables.

3. If h differs from s, some original word has head s, OR some original nonlinear word has first two variables exactly h, s. Again, outside variables are allowed in these witnesses.

For a nonlinear target in 3111, the criterion is:

1. Some original word is contained in C and has head h.

2. If h differs from s, some original contained word has head s, OR some original contained nonlinear word starts h, s.

In 3111 every selected witness is contained. These two criteria must not be interchanged. They are invariant under (W).

Necessity for 2008: absence of a contained word is separated by assigning all C variables 4 and outsiders 0: the target is 4 and every upper is 0. Absence of a same-head word is separated by \(h=0\) and all others 3: the target is 0 and all upper values 3. Failure of condition 3 is separated by \(h=2\), \(s=0\), all others 3. The target is 1. Upper heads other than h, s give 3; an h-headed literal gives 2, and an h-headed nonlinear word has second value 2 or 3 and gives 3. Thus the upper sum is below 1. A literal h lacking the original literal is separated by \(h=2\), outsiders 3: all upper values are 3, while the target is 2.

Necessity for 3111: absence of a contained same-head word is separated by \(h=4\), other C variables 0, outsiders 3. The target is 4; every upper word is 0 or 3. Failure of condition 2 is separated by \(h=1\), \(s=4\), other C variables 2, outsiders 3. The target is 0. A noncontained upper gives 3; the remaining upper values are 1 or 2 unless it has head s or starts h, s, which are exactly the missing witnesses. For a literal h lacking the original literal, assign \(h=1\), outsiders 3. Upper pure powers give 2 and noncontained words 3, while the target is 1.

For sufficiency in 2008, target 3 is least. Target 4 means all C values 4; the contained witness gives 4. Target 0 is handled by a same-head witness when its head 0, or by one with head 4 (whose value is 0 or 4). Target 1 with head 1 is similarly immediate. If target 1 has head 2, its second variable has value 0 or 4 and is different from h. The third condition supplies a word with that second head (giving 0 or 4), or a word with the same first two factors (giving 1). Literals are already included.

For 3111, target 3 is least. Otherwise all C factors avoid 3, so the contained witnesses avoid 3 as well. Targets4,\allowbreak{}2, and head 0 target 0 are handled by the contained same-head witness. A nonlinear head 1 target 2 is also bounded by any such witness, whose value is 0,\allowbreak{}1, or 2. A head 1 target 0 has second value 4, so h differs from s; the second criterion supplies a word headed by 4 (value 4) or with first factors 1,\allowbreak{}4 (value 0). This proves both criteria unrestrictedly, not by bounded testing.

\par\smallskip\noindent\textit{Few witnesses and reversible compression.}\quad 
For 2008 select original witnesses for the three conditions, at most three upper words. For 3111 select at most two. The same criterion proves the target absorbed by their sum. To cover both sources with one bound, there are at most four words including the target.

Normalize by (W). Reserve every variable occupying either of the first two positions of any word (only one for a literal). There are at most eight such variables. Every other variable occurs only in the permutable suffix, at most once in each word. Group these variables by their nonempty occurrence pattern among the at most four words. At most 15 groups are needed. Replace each reserved variable by its own variable and each group by a new variable, then normalize. The result uses at most 23 variables, has word length at most 25, and at most three upper words.

All content containments and first-two-position equalities are preserved. Literals stay literals and repeated first positions retain both copies. The independently proved criterion therefore makes the compressed absorption valid. Substitute each reserved variable by its original singleton and each group variable by the nonempty increasing product of its members. Every expanded word has exactly the original first positions and content. Literal applications of (W) normalize it to the same normal word as the original. This is an inverse derivation, not just a one-way HSP statement or loss of repeated occurrences.

\par\smallskip\noindent\textit{Fixed finite families and completeness.}\quad 
For each source take(W), plus all valid absorptions over 23 variables, each word length at most 25, and at most three nonempty upper words. This is a finite effectively specified family; its validity is decidable by the displayed exact criterion or by its finite table. Any valid absorption has the selected at-most-three-witness reduction, whose compressed absorption is a member of this fixed family. Inverse word substitutions and(W) recover the original reduction; inclusion of its original summands yields the full absorption. There is no bound on the variables, lengths, or original polynomial size in this derivation. Distributive expansion and both directions of absorption give every identity. Consequently 2008 and 3111, and their opposites, are finitely based. Catalogue transfers require independently complete mutual HSP.

\subsection{A finite basis for 2138 from three original witnesses with nonlinear head}\label{case:75}
Signature (+, multiplication), without constants, relative to additively idempotent semiring axioms. The finite family is explicitly specified, not an enumerated small list. All completeness arguments are unrestricted in variable count, word length and polynomial size.

\par\smallskip\noindent\textit{Word values and tail normalization.}\quad 
Addition rows:00044,\allowbreak{}01144,\allowbreak{}01234,\allowbreak{}44334,\allowbreak{}44444; order \(2<1<0<4\) and \(2<3<4\). Multiplication rows:00000,\allowbreak{}22211,\allowbreak{}22222,\allowbreak{}00044,\allowbreak{}00044. Literal words retain their values. For nonlinear w with head h and tail T, heads 0 and 2 retain their values. Head 1 gives 1 if all tail values belong to\{3,\allowbreak{}4\}, and 2 otherwise. Head 3 or 4 gives 4 if all tail values belong to \{3,\allowbreak{}4\}, and 0 otherwise. Induction from the table proves this formula.

It proves W1 \(xyz=xzy\) and W2 \(xyy=xy\). Normal words are a head followed by the sorted tail set. Tail sorting and duplicate deletion are literal applications of W1, W2 with nonempty substituted prefixes. A head repeated in the tail remains recorded. No literal is equated with its square.

\par\smallskip\noindent\textit{Exact absorption criterion.}\quad 
For original P, a literal x is absorbed iff x is an original literal. For nonlinear \(w=hT\), \(C=\{h\}\) union T, let HN be the heads of original NONLINEAR right words. Then \(w\leq P\) exactly when: (NH) h belongs to HN; (C) some original right word has content contained in C; (R) either original literal h occurs, or T intersects HN, or an original nonlinear h-headed word has its tail contained in T.

Necessity for a missing literal target x:assign \(x=1\) and every other variable 2. Every nonlinear x-word gives 2, as do all other right heads; without the literal x the right sum 2 cannot absorb 1. For failed NH assign \(h=3\) and every other variable 2. Right words have value 2 or, from literal h,3. The target is 4 if its tail is pure h, and 0 otherwise; neither is below 3. For failed C assign all target variables 4, outsiders 2. The target is 4. Every right word contains an outsider, so its value is 0 or 2, whose sum is at most \(0<4\). For failed R, assuming NH and C, T avoids HN and h is outside T. Assign \(h=1\), \(T=3\), others 2. All right nonlinear words have value 2:heads in T are forbidden by failure of R, and every nonlinear h-word has a tail outside T. Original literal h is absent; other literal words are 2 or 3. The target is 1, which is not below the right value 2 or 3. Thus R is necessary.

Sufficiency:target 2 is least. Target 0 has head 0 or head 3/\allowbreak{}4 with a low tail; NH supplies a nonlinear same-head word valued 0 or 4. Target 4 has every target variable in\{3,\allowbreak{}4\}. A contained nonlinear witness from C gives 4; a contained literal gives 3 or 4. In the former literal 3 case, the nonlinear same-head witness from NH gives 0 or 4, so the join with 3 is 4. Target 1 has head 1 and all tails in\{3,\allowbreak{}4\}. R supplies either the literal head 1, a nonlinear same-head word with tail in T valued 1, or a nonlinear right word whose head is in T valued 0 or 4. All absorb 1. Literal targets are original summands. This proves the criterion for arbitrary words and polynomial sizes.

The preliminary criterion omitting C was false:\(xx\leq xy+yx\) fails at \(x=3\), \(y=2\). A literal tail-head alone is also insufficient:\(xy\leq xz+y\) fails at \(x=1\), \(y=3\), \(z=0\). These are recorded counterexamples, not classification evidence.

\par\smallskip\noindent\textit{Select three original witnesses and compress reversibly.}\quad 
For a valid nonlinear target, choose an original NONLINEAR h-headed word a for NH, an original contained word b for C, and an original R witness c. The R witness is literal h, or a nonlinear word with head in T, or a nonlinear h-word with tail contained in T. Discard duplicates to obtain Q with at most 3 words. The same NH, C, R conditions hold for \(w\leq Q\); hence it is valid. All other original summands are retained using addition and transitivity.

There are at most 4 normal words in \(w\leq Q\). Reserve their at most 4 heads. Partition every remaining variable by its nonempty incidence pattern among these at most 4 tail sets. At most 15 pattern classes occur. Identify each to a fresh variable, keeping reserved heads separate, and normalize by W1, W2. There are at most 19 variables, each word has length at most 20, and at most 3 right words. Substitution preserves validity.

Invert the compression by substituting the nonempty product of all members for each pattern variable. Identical tail incidence recovers each original tail exactly. Reserved heads and their tail appearances were never merged. Normalize by W1, W2 to obtain all original words. This is reversible without empty substitutions or losing a head occurrence.

Let Sigma consist of W1, W2 and every valid normal absorption on a fixed 19-variable alphabet, word length at most 20, with at most 3 nonempty right words. NH, C, R decide membership; the family is finite and completely specified. Every valid absorption selects Q, compresses into Sigma, and expands back. Literals need no new absorption. Distributivity and absorption of each summand in both directions then derive every valid identity. Therefore Sigma is a finite relative basis for 2138. Reversal gives the opposite basis; catalogue transfer requires complete HSP both directions. Finite computations validate concrete steps but do not replace this proof.

\subsection{Square guards and marked prefixes: twelve- and ten-identity bases}\label{case:76}
Work in the constant-free signature (+, multiplication), relative to the axioms of additively idempotent semirings. The element labels below are not constants of the signature.

\par\smallskip\noindent\textit{The two algebras.}\quad 
The sources A and B are catalogue entries 2140 and 2141. They have common addition and the indicated multiplication tables, with rows and columns ordered 0,\allowbreak{}1,\allowbreak{}2,\allowbreak{}3,\allowbreak{}4:

{\small\par\smallskip
\noindent addition      A multiplication     B multiplication\par
\CaseLine{\ttfamily 00044\qquad 00000\qquad 00000}
\CaseLine{\ttfamily 01134\qquad 22211\qquad 22211}
\CaseLine{\ttfamily 01234\qquad 22222\qquad 22222}
\CaseLine{\ttfamily 43334\qquad 01234\qquad 02234}
\CaseLine{\ttfamily 44444\qquad 00044\qquad 00044}
\smallskip}

Addition is join in the order \(2<1<0<4\) and \(2<1<3<4\), where 0 and 3 are incomparable. The only difference between the multiplication tables is 3\ensuremath{\cdot}1: it is 1 in A and 2 in B. Both tables are associative, and every left and right multiplication preserves joins, as is checked directly in the displayed five-element tables. Here 0 and 2 are left zeros. Put K=\{3,\allowbreak{}4\}. A word whose occurrences all have values in K has value 3 if all are 3, and 4 otherwise. For any other word, take the first occurrence outside K, with value d. If any preceding occurrence has value 4, the word is 0. Otherwise every preceding occurrence is 3, and:

\textbullet{} \(d=0\) gives 0;

\textbullet{} \(d=2\) gives 2;

\textbullet{} \(d=1\) gives 1 in A exactly when all later occurrences lie in K, and gives 2 otherwise;

\textbullet{} in B the \(d=1\) case gives 1 exactly when it is also the first position of the word and all later occurrences lie in K; otherwise it gives 2.

These statements follow by successive multiplication, and include one-letter words. They also give a direct case check of associativity: extending a prefix with value 0 or 2 never changes it, a prefix valued 4 switches to 0 on the first later non-K factor, and the indicated surviving value 1 switches to 2 on any later non-K factor.

\par\smallskip\noindent\textit{Explicit identities and normal forms.}\quad 
The following nine identities are common to both bases. As usual u\ensuremath{\leq}v abbreviates u+v\ensuremath{\approx}v.

{\small\par\smallskip
\CaseLine{\((W)  xyx \approx  x^{2}y\)}
\CaseLine{\((S)  x^{2} \leq  x\)}
\CaseLine{\((P0) xy \leq  x+y\)}
\CaseLine{\((P1) xy \leq  y+xz\)}
\CaseLine{\((P2) xyz \leq  xz+xy\)}
\CaseLine{\((P3) xyz \leq  xz+xyt\)}
\CaseLine{\((G0) xy \leq  (xy)^{2}+y\)}
\CaseLine{\((G1) xy \leq  (xy)^{2}+x\)}
\CaseLine{\((G2) xyz \leq  (xyz)^{2}+xz\)}
\smallskip}

For A add three identities:

{\small\par\smallskip
\CaseLine{\((A1) x^{3} \approx  x^{2}\)}
\CaseLine{\((A2) xy \approx  x^{2}y+xy^{2}\)}
\CaseLine{\((A3) xy \leq  (xy)^{2}+yx.\)}
\smallskip}

For B instead add one identity:

{\small\par\smallskip
\CaseLine{\((B1) xy^{2} \approx  xy.\)}
\smallskip}

Theorem. The twelve identities above for A and the ten identities above for B are complete identity bases, relative to the ai-semiring axioms. Reversing all words gives bases for their multiplicative opposites. No minimality or independence is asserted.

For A a normal word lists variables in their order of first occurrence, each appearing in one block of length one or two according to whether its original multiplicity is one or at least two. Indeed (W) moves a later occurrence next to its first occurrence, and (A1) caps multiplicities at two. In B, (W) again groups repeated occurrences; (B1) removes every repetition of a variable other than the head, and caps the head multiplicity at two. It implies \(x^{3}\)\ensuremath{\approx}\(x^{2}\) by substituting x for both variables. Thus a B normal word retains its first-occurrence list and whether its head occurs once or at least twice. All transformations use nonempty substituted words.

Write c(u) for content and \(pre_{u}\)(x) for the set of variables before the first x in u. Write \(mult_{u}\)(x) for its multiplicity. A saturated word in A has every variable twice; a saturated word in B has its head twice and every other variable once. In either algebra let s(u) be the normal form of \(u^{2}\). It is saturated and retains the first-occurrence list of u. Products of saturated words normalize to the saturated word on their concatenated first-occurrence list.

\par\smallskip\noindent\textit{Exact criteria for absorption.}\quad 
Let P be a nonempty sum of words U, with linear summands included in U. For a target word w put \(C=c(w)\) and U\_C=\{u in U:c(u)\ensuremath{\subseteq}C\}.

Lemma 1. The absorption w\ensuremath{\leq}P holds in A exactly when all three conditions hold:

{\small\par\smallskip
\noindent (J) \(U_{C}\) is nonempty.\par
\noindent (F) For each x in C some u in U contains x and satisfies \(pre_{u}\)(x)\ensuremath{\subseteq}\(pre_{w}\)(x).\par
\noindent (M) For each x occurring once in w, some u in \(U_{C}\) satisfies \(mult_{u}\)(x)\ensuremath{\leq}1.\par
\smallskip}

For B replace (M) by:

{\small\par\smallskip
\noindent (H) If the head h of w occurs once, some u in \(U_{C}\) either omits h or has head h and \(mult_{u}(h)=1\).\par
\smallskip}

These criteria apply to linear words too; in that case they say that the target is a literal summand.

Necessity. If (J) fails, assign C the value 3 and outsiders 2. The target is 3 and every right word is 2. If (F) fails at x, assign \(pre_{w}\)(x) the value 3, x the value 0, and every remaining variable 2. The target is 0. A right word can have value 0 only if its first x has all earlier variables in \(pre_{w}\)(x), contrary to hypothesis. Therefore all right words lie in \{2,\allowbreak{}3\}, their sum is at most 3, and cannot absorb 0.

If (M) fails at a simple x, assign \(x=1\), every other variable of C the value 3, and outsiders 2. The A target is 1. A contained right word has at least two occurrences of x, and is 2. A noncontained word has an occurrence of 2 and, with all assigned values in \{1,\allowbreak{}2,\allowbreak{}3\}, also has value 2. Thus \(P=2\).

For failure of (H), assign \(h=1\), the other C-variables 3, and outsiders 2. The B target is 1. Every contained right word has h but either repeats it or does not start with it, hence is 2; every noncontained word is also 2. This proves necessity.

Sufficiency. If \(w=2\) the result is automatic. If \(w=0\), its first occurrence not valued 3 is valued 0 or 4: a first occurrence valued 1 or 2 could only lead to 1 or 2. Take this first variable x. The witness in (F) has all variables before x valued 3. It therefore evaluates to 0 or 4, either of which absorbs 0.

If \(w=3\), all target variables are 3 and any witness in (J) is 3. If \(w=4\), all target variables lie in K and at least one x has value 4. A witness in (J) is 3 or 4. A witness in (F) for a target occurrence valued 4 has only K-valued variables before x and evaluates to 0 or 4. The sum of these witnesses is 4.

If \(w=1\) in A, there is exactly one occurrence outside K, with value 1, its variable x is simple, and all earlier occurrences have value 3. A contained witness from (M) which omits x has value 3 or 4. One with one occurrence of x has value 0 when preceded by a 4, and otherwise value 1. In all cases its value is at least 1. In B a target valued 1 starts with its unique occurrence of \(h=1\); its other variables lie in K. A witness from (H) either omits h and has value 3 or 4, or starts with its unique h and has value 1. Thus it absorbs w. This proves the lemma for arbitrary words and variable sets.

Every displayed basis identity is sound. This can be seen from the word-value description for the word laws, and by substituting the displayed words into Lemma 1 for the absorption laws. For (A2), the left-to-right absorption uses x\ensuremath{{}^{2}}y to witness the single y and xy\ensuremath{{}^{2}} to witness the single x; the reverse absorptions follow from (S) and monotonicity. In (A3), yx supplies the two single-occurrence witnesses, while the squared word supplies the ordered prefix requirements. This step is needed only in A and is not assumed in B.

\par\smallskip\noindent\textit{Step 1: deriving all saturated prefix extensions.}\quad 
Suppose w\ensuremath{\leq}P satisfies the relevant criterion. Let \(x_{1}\),..., \(x_{d}\) be the target's first-occurrence list, and put \(p_{i}=x\)\ensuremath{{}_{1}}...\(x_{i}\). Rule (S) first gives s(u)\ensuremath{\leq}u\ensuremath{\leq}P for every original right word u. We claim that

{\small\par\smallskip
\noindent s(\(p_{i}\) u) \ensuremath{\leq} P     for every u in U and 1\ensuremath{\leq}i\ensuremath{\leq}d.             (*)\par
\smallskip}

The first target variable \(x=x\)\ensuremath{{}_{1}} has, by (F), an original witness v starting with x. If v is linear, \(x^{2}\)\ensuremath{\leq}P is available and (P0) with \(x^{2}\) and s(u) derives x\ensuremath{{}^{2}}s(u)\ensuremath{\leq}P. Otherwise write \(v=xt\) with t nonempty. Then x\ensuremath{{}^{2}}s(t) normalizes to s(v), and (P1) yields x\ensuremath{{}^{2}}s(u)\ensuremath{\leq}s(u)+x\ensuremath{{}^{2}}s(t)\ensuremath{\leq}P. This establishes (*).

Suppose (*) holds at \(i<d\). Let \(p=s(p_{i})\) and \(x=x_{i+1}\). Choose its prefix witness v. The derived word p s(v) normalizes to p \(x^{2}\) t, where the square of x is interpreted through the appropriate normal form, and t is a possibly omitted suffix. All letters preceding the first x in v were already in p and disappear or merge into their existing saturated blocks. Meanwhile p s(u)\ensuremath{\leq}P for each original u. If t is omitted, apply (P2); otherwise apply (P3), with \(x^{2}\) substituted for its middle variable. These laws give

{\small\par\smallskip
\CaseLine{\(p x^{2} s(u) \leq  p s(u)+p x^{2},\)}
\CaseLine{\(p x^{2} s(u) \leq  p s(u)+p x^{2} t,\)}
\smallskip}

respectively. This proves (*) at i+1. The separate laws ensure that no empty word is substituted.

Choose u in \(U_{C}\) by (J). At \(i=d\), the word in (*) normalizes to s(w). We have derived s(w)\ensuremath{\leq}P. More generally (*) makes available every saturated word on the first-occurrence list of \(p_{i}\) followed by any original u; these are the guards used in the next step.

\par\smallskip\noindent\textit{Step 2: recovering one marked occurrence.}\quad 
In A fix a variable x simple in w. Choose u in \(U_{C}\) with \(mult_{u}\)(x)\ensuremath{\leq}1 by (M). In B this step is needed only when \(h=head(w)\) is simple: put \(x=h\) and choose the witness in (H).

Starting with u\ensuremath{\leq}P, apply (S) in multiplicative contexts to duplicate each variable other than x, then normalize. The resulting word v is still below P, uses only C, and has no possible simple block except x. In B, if x occurs then it is already the head.

We put the variables of v into the target order by walking along \(x_{1}\),..., \(x_{d}\). At each stage its first-occurrence list begins with the already processed target prefix. For a next variable y different from x, insert \(y^{2}\) immediately after this prefix; normalize to merge a later occurrence of y. For \(y=x\) absent from v, insert one x instead. In either case the desired insertion is derived by (G0), (G1), or (G2), according as the processed prefix is empty, the residual suffix is empty, or both are nonempty. The other right summand is the existing word v, and the squared right summand is the saturated word s(\(p_{i}\) u) from (*). Thus every insertion is derivable. Only nonempty blocks are substituted into these three separate laws.

There is one additional case in A: the next variable is x, and its single occurrence is later in v, separated from the processed prefix p by a nonempty block a. Write \(v=p\) a x b, where p and b may be omitted. Apply (A3) with its x variable equal to x and its y variable equal to a, in the context p[-]b. It gives

{\small\par\smallskip
\CaseLine{\(p x a b \leq  p(xa)^{2}b + p a x b.\)}
\smallskip}

All other blocks are already saturated; the first right summand consequently normalizes to the saturated guard s(\(p_{i}\) u) in (*), and the second is the available v. This moves x to its proper position. If a is empty no move is needed. This case never occurs in B: its potentially simple variable is the target head, and a witness containing it already starts there.

After the last stage all target variables occur in target order. In A the result is \(w_{x}\), the normal word with precisely x single and all other variables twice. In B the result is exactly w. If the B target head was multiple, Step 1 had already given \(w=s(w)\). This proves completeness for B.

\par\smallskip\noindent\textit{Step 3: recombining simple variables in A.}\quad 
For A, iteration of (A2) shows that a word \(a_{1}\)...\(a_{n}\) is the sum of the words obtained by leaving one position \(a_{i}\) unsquared and squaring every other position, followed by normalization. For \(n=1\) this statement is literal. For \(n>1\), split off the last letter in (A2). The summand \((a_1\cdots a_{n-1})\)\ensuremath{{}^{2}} \(a_{n}\) is the last-position term, and expansion of \(a_1\cdots a_{n-1}\) in the context [-]\(a_{n}\)\ensuremath{{}^{2}} gives all the earlier-position terms. The word laws identify a squared prefix with the ordered product of its squared letters.

If a variable is repeated in w, every term of this expansion in which one of its occurrences was left unsquared still makes that variable multiple. Hence each normalized term is either s(w) or \(w_{x}\) for some variable x occurring once in w. Step 1 and Step 2 have derived all these absorptions below P. Their sum gives w\ensuremath{\leq}P. This proves completeness for A.

In either algebra a linear target is a literal right summand by Lemma 1. Finally, distribute arbitrary terms into sums of words and derive every summand absorption in both directions. The proof works for any number of variables and any word or polynomial lengths. Word reversal transports the bases and all derivations to the opposite algebras.

\subsection{A finite basis for 2183 from three original witnesses}\label{case:77}
Signature (+, multiplication), without constants, relative to the additively idempotent semiring axioms. The specified finite family is not claimed to be an enumerated small list. Completeness is unrestricted.

\par\smallskip\noindent\textit{Tables, word values and tail normal form.}\quad 
Addition rows:00004,\allowbreak{}01114,\allowbreak{}01224,\allowbreak{}01234,\allowbreak{}44444; its order \(is 3<2<1<0<4\). Multiplication rows:00000,\allowbreak{}33331,\allowbreak{}33331,\allowbreak{}33333,\allowbreak{}00004. A literal retains its value. For a nonlinear word with head h and tail T: head 0 gives 0; head 3 gives 3; head 1 or 2 gives 1 if every tail is 4, and 3 otherwise; head 4 gives 4 if every tail is 4, and 0 otherwise. This follows immediately by induction:4 is a right identity, and any other right factor takes a head 1/\allowbreak{}2 into absorbing 3 and head 4 into absorbing 0.

The identities W1 \(xyz=xzy\), W2 \(xyy=xy\) are valid. They sort and deduplicate the tail with the head kept separately. Thus a normal word is h followed by the sorted set T. Literal heads and repeated heads remain distinguished. Swaps use a nonempty prefix as x; duplicate deletion likewise has a nonempty prefix. Contexts may be absent, but substituted words never are empty.

\par\smallskip\noindent\textit{Exact unrestricted criterion.}\quad 
A literal target x is absorbed by original P iff x itself is an original literal summand. For nonlinear \(w=hT\), \(C=\{h\}\) union T, write H for the set of original heads. Then \(w\leq P\) iff: (H) h belongs to H; (C) some original word is contained in C; (R) either T intersects H, or an original nonlinear h-headed word has its tail contained in T.

Necessity of H:set \(h=0\) and every other variable 3. The target is 0, while all right words are 3. Necessity of C:set all target variables 4 and all outsiders 3. The target is 4; every right word with an outsider is at most 0. For failure of R, assume H and C. Then T avoids H, and in particular h is outside T. Assign \(h=2\), all variables in T to 4, and all other variables 3. The target is 1. Every original word with head other than h gives 3. An original literal h gives 2; every original nonlinear h-word has a tail outside T, so gives 3. The right sum is at most \(2<1\). This separates R.

For a missing literal x, set \(x=2\) and every other variable 3. Right nonlinear x-words give 3 (every right factor is 2 or 3); other heads give 3; without the literal x the right sum \(is 3<2\). These separate every failed condition.

Sufficiency:target 3 is least. Target 4 means every variable in C is 4; C supplies an original contained word valued 4. Target 0 has \(h=0\), or \(h=4\) with some non4 tail. H supplies an original h-word whose value is 0 or 4. Target 1 means h is 1 or 2 and every tail variable is 4. If T intersects H, an original word with a head in T has head 4 and value 0 or 4, both above 1. Otherwise R supplies an original nonlinear h-word with tail in T, valued 1. Literal targets follow from their original literal summands. This proves the criterion for arbitrary variables, lengths and polynomial sizes.

\par\smallskip\noindent\textit{At most three original right words suffice.}\quad 
For any valid nonlinear target, choose an original h-headed word a for H, an original contained word b for C, and one original word c witnessing R: if T intersects H, choose a word whose head lies in T; otherwise choose the original nonlinear h-word with tail contained in T. Let Q be the set of these words, discarding duplicates. Then |Q|\ensuremath{\leq}3 and the same three conditions hold for \(w\leq Q\). Hence this absorption is valid. Literal targets require no new absorption. All original summands outside Q may be retained by addition and transitivity; there is no restriction on their total number.

\par\smallskip\noindent\textit{Reversible compression and a fixed finite family.}\quad 
There are at most four normal words in \(w\leq Q\). Reserve every variable that is the head of any of them; there are at most 4. Every other variable occurs only in tails. Partition these nonreserved variables by their nonempty incidence pattern among the at most 4 tails. There are at most \(2^{4}-1=15\) patterns. Identify each pattern class to a new variable, leaving reserved heads separate, and normalize by W1, W2. There are at most 19 variables and each normal word has length at most 20. Substitution preserves validity.

The compression is reversible:substitute the nonempty product of the original class members for each new variable. Their tail incidences were identical, so each original tail set is recovered exactly. Reserved heads and their tail occurrences remain separate. Sorting and deduplication recover all original normal words. No empty substitutions, no unrecorded head multiplicities, and no artificial balancing identities are used.

On a fixed 19-variable alphabet let Sigma consist of W1, W2 and EVERY valid absorption with one normal left word and at most 3 nonempty normal right words, each of length at most 20. The exact criterion above specifies finite membership. This finite family is well defined; no small enumeration is asserted. Every valid absorption \(w\leq P\) reduces to at most 3 original witnesses, compresses into Sigma, and expands back by the described nonempty word substitution modulo W1, W2. Thus it is derivable from Sigma.

Distributivity and summandwise absorption in both directions give every valid identity. Sigma is therefore a finite relative basis for 2183. Opposite reversal preserves this argument. Catalogue transfers still need independent complete HSP both ways. Finite checks supplement this theorem.

\subsection{四元素 T 的四条完整基及五元素等簇转移}\label{case:78}
二元无常元语言；以下基均相对于加法幂等半环公理。

\par\smallskip\noindent\textit{1. 运算与有限基.}\quad 
令 T=\{0,\allowbreak{}1,\allowbreak{}2,\allowbreak{}3\}，两个运算的逐行编码为

\(+ : 0000012002200003\)，\ensuremath{\cdot} : 0000033003000000。

这里 0 是两个运算的吸收元素，\(1+2=2\)，3 与 1、2 在加法序下不可比；非零乘积只有 1\ensuremath{\cdot}\(1=1\)\ensuremath{\cdot}\(2=2\)\ensuremath{\cdot}\(1=3\)。元素名称不是语言中的常元。

定理。 T 的恒等式有以下四条完整基：

1. xy\ensuremath{\approx}yx；

2. xyz\ensuremath{\approx}\(u^{3}\)；

3. x+xy\ensuremath{\approx}\(z^{3}\)；

4. \(x^{2}\)+yz\ensuremath{\approx}\(x^{2}\)+yz+xy。

四式都由运算表直接成立。下面证明不限变量数、词长与加数数目的完整性。

\par\smallskip\noindent\textit{2. 吸收项与正规多项式.}\quad 
由第二式，所有至少三次的词互等，记它们的共同等式类为 \ensuremath{\Omega}。记号 \ensuremath{\Omega} 只用于证明，始终代表一个立方项，不增加常元。第二式还推出 \ensuremath{\Omega}t\ensuremath{\approx}t\ensuremath{\Omega}\ensuremath{\approx}\ensuremath{\Omega}。第三式中以 \ensuremath{\Omega} 代入 y，得 t+\ensuremath{\Omega}\ensuremath{\approx}\ensuremath{\Omega}。因此出现长词的多项式等于 \ensuremath{\Omega}。

展开分配律，并用乘法交换律和加法幂等律，剩余项可写成线性变量集合 L 与二次单项式集合 E 的和。令 N 为 E 中出现的变量集合。若 L\ensuremath{\cap}N 非空，可选择 x\ensuremath{\in}L 及一个含 x 的二次项 xy，用第三式将这两项之和变成 \ensuremath{\Omega}，所以整个多项式等于 \ensuremath{\Omega}。

其他多项式满足 L\ensuremath{\cap}N=\ensuremath{\varnothing}，称为非退化短多项式。E 可看成顶点集 N 上带圈的无向图：平方是圈，xy（x\ensuremath{\ne}y）是边。令 M(E) 为其单项式支持集合族中按包含关系极小的成员；即全部圈单点，以及两个端点均无圈的边。

\par\smallskip\noindent\textit{3. 任意变量数的分离赋值.}\quad 
若 E 非空，非退化短多项式 P 在 T 上非零，当且仅当：

\textbullet{} L 的每个变量赋值 3；

\textbullet{} N 的每个变量赋值 1 或 2；

\textbullet{} 没有二次单项式的全部端点均赋值 2。

此时 P 的值为 3，其他赋值下为 0。令 L 的值为 3、N 的值为 1，立即得到非零赋值。纯线性项也是非恒零的。因此它们都不能与 \ensuremath{\Omega} 恒等。

非退化短多项式的恒等函数精确决定 (L, N, M(E))。证明如下。先比较总支持 L\ensuremath{\cup}N：若两边总支持不同，对包含较少变量的一边取非零赋值，将另一边独有变量置 0，即得到分离。总支持相同时，二次情形下非零赋值要求线性角色取 3、非线性角色取 1 或 2，所以角色不同亦可分离。纯线性、纯二次及线性二次混合三种情形，在所有变量赋值 1 时分别取 1、3、0，彼此可区分。

固定相同 L, N 后，令 N 中赋值 2 的变量组成 U，其余赋值 1，L 赋值 3。P 非零恰当且仅当 U 不包含 E 的任何单项式支持。这样的允许集合族精确决定其极小禁止集合 M(E)。若两个极小禁止集合族不同，存在一边禁止而另一边允许的 U，相应赋值分离两边。该论证适用于任意有限变量集合。纯线性多项式的函数则精确由 L 决定。

\par\smallskip\noindent\textit{4. 等式推导的完整性.}\quad 
设 P\ensuremath{\approx}Q 在 T 中成立。若一边退化而等于 \ensuremath{\Omega}，另一边必然退化；两边均由前述步骤推导为 \ensuremath{\Omega}。否则两边都是非退化短多项式，并具有相同 L, N, M(E)。线性部分由加法幂等和交换律已经相同。

取 Q 中任一二次项 xy。其支持含有 P 的一个极小支持。若此支持是边 \{x, y\}，则 P 已含 xy。否则它为端点的圈，例如 \{x\}，故 P 含 \(x^{2}\)。又因 \(N(P)=N(Q)\)，变量 y 在 P 的某个二次项 yv 中出现。第四式允许将 xy 加入 P。若两见证项相同，可先用加法幂等律复制。因而可依次加入 Q 的全部二次项，推导 P\ensuremath{\approx}P+Q。交换 P, Q 再推导 Q\ensuremath{\approx}P+Q，得 P\ensuremath{\approx}Q。

这证明四式的完整性，并非仅验证有界个变量的恒等式。

\par\smallskip\noindent\textit{5. 文献及编号说明.}\quad 
A.~Wang、L.~Wang、Q.~Yin 与 J.~Wu 在表~23 中将 \(T\) 标为 \(A_{513}\)，并在定理~6.4 中给出由二十条恒等式组成的基 \cite[Section~6.1, Theorem~6.4 and Table~23]{WangYinWu}。上文证明了这里的四条恒等式也构成完整基。

按本文的五元素编号，\(10572,12212,12213,12215,12260,12291,12292,12293,13369,15140\) 与 \(T\) 生成同一簇。两个方向的包含均用命题~\ref{prop:lower-fibres} 判定。

\subsection{Finite bases for two regularized quartic nilpotent semirings}\label{case:79}
All terms are in the constant-free signature \((+,\cdot)\); bases are relative to the additively idempotent semiring axioms. Write \(u\le v\) for \(u+v\approx v\).

\par\smallskip\noindent\textit{1. Algebras and bases.}\quad 
On \(\{a, b, c, o\}\) define commutative multiplication by \(a^2=b\), \(ab=c\), and every other product equal to \(o\). In \(B_0\), \(a, b, c\) are incomparable below the additive greatest element \(o\). In \(B_1\), add the relation \(a<c\). Addition is join. In each case adjoin a new greatest element \(t\) that absorbs both operations, and denote the resulting five-element semiring by \(A_i\). In the displayed order \((a, b, c, o, t)\), their catalog relabelings are \((3,\allowbreak{}2,\allowbreak{}1,\allowbreak{}0,\allowbreak{}4)\) for type 2315 and \((3,\allowbreak{}1,\allowbreak{}2,\allowbreak{}0,\allowbreak{}4)\) for type 2553. The element names \(o, t\) are not constants of the signature.

Within \(B_i\), multiplication adds positive exponents of \(a\), with every exponent at least four identified with \(o\). This proves associativity. Translation by \(a\) takes \(a, b, c, o\) to \(b, c, o, o\); translation by \(b\) takes them to \(c, o, o, o\); the others are constant \(o\). These maps preserve the joins just specified. Adjoining the common absorber preserves associativity and distributivity.

Let \(\Sigma\) consist of the following fourteen laws:

\CaseDisplay{\begin{array}{rll}
(C)&xy\approx yx,\\
(S)&x^5\approx x^4,\\
(Q)&xyzt\approx x^4+y^4+z^4+t^4,\\
(T)&xyz\approx x^3+y^3+z^3,\\
(U)&x\le x^4,\\
(Z_1)&x^2+y^3\approx x^4+y^4,\\
(P)&xt\le xy+yz+zt,\\
(L)&y^2\le x^2+yz,\\
(E)&xy\le x^2+y^2,\\
(M)&z^3\le x^3+xy+yz,\\
(J)&xz\le x^3+y^3+yz,\\
(Z_2)&xy+x^3+y^3\approx x^4+y^4,\\
(A_1)&y^4\le x^4+y,\\
(A_2)&y^4\le x^4+yz.
\end{array}}

Theorem. A basis for \(A_0\) is \(\Sigma\) together with

\CaseDisplay{(F)\qquad x^4+y^4\le x+xy.}

A basis for \(A_1\) is \(\Sigma\) together with

\CaseDisplay{(H)\qquad x\le x^3,\qquad (I)\qquad x^3\le x+xy.}

Thus these bases have fifteen and sixteen displayed laws, respectively. No minimality or hereditary finite-basis assertion is made.

Every law displayed is regular: its two sides contain the same variables after an inequality is written as an absorption identity. For arbitrary polynomials \(P, R\), the common-absorber construction gives

\CaseDisplay{A_i\models P\approx R
\quad\Longleftrightarrow\quad
B_i\models P\approx R\ hbox{and}\ c(P)=c(R).\quad\textup{(\ensuremath{1})}}

Indeed an assignment containing \(t\) in the common content makes both values \(t\); otherwise it is an assignment in \(B_i\). If the contents differ, assign a missing variable \(t\) and every other variable \(o\) to separate the values. Accordingly

\CaseDisplay{A_i\models w\le R\iff c(w)\subseteq c(R)
\ hbox{and}\ B_i\models w\le R.\quad\textup{(\ensuremath{2})}}

Equation (1) alone is not being used as a general finite-basis preservation theorem; completeness is proved below.

Validity of the displayed laws can be checked inside \(B_i\) using (1). A quadratic has value \(b\) exactly at \((a, a)\), and value \(c\) exactly at \((a, b),(b, a)\). A cubic is \(c\) exactly when every factor is \(a\); all fourth powers are \(o\). Thus the three-edge rule (P) preserves an odd path in a two-coloring, (M) propagates the color \(a\) along two edges, and (J) joins a forced \(a\) vertex to a forced \(b\) vertex. If a square is at most a non-top polynomial value, that value must be \(b\) and its vertex must be \(a\), which verifies (L),(E). Squares and cubes cannot have the same non-top value, proving \((Z_1)\); an edge and cubes at both its ends also force top, proving \((Z_2)\). The other common laws follow from the product formula and the greatest-element property of \(o\). In \(B_0\), \(x+xy=o\) for every assignment, proving (F). In \(B_1\), \(x\le x^3\), and the only non-top value of \(x+xy\) occurs at \(x=a, y=b\), giving value \(c\); this proves (H),(I).

\par\smallskip\noindent\textit{2. Normalization and propagation of fourth powers.}\quad 
By (S), \(x^n\approx x^4\) for all \(n\ge4\). Substituting \(u^2, u^2, v^2, v^2\) in (Q) gives

\CaseDisplay{(uv)^4\approx u^8+v^8\approx u^4+v^4.}

Induction therefore expresses the fourth power of a word as the sum of the fourth powers of its variables. Splitting a word of length at least four into four nonempty factors and applying (Q) gives the same expression for the word itself. Law (T) handles words of length three. A normal polynomial consists of linear terms, quadratic terms, cubes and fourth powers.

Propagation lemma. If a normal polynomial \(R\) has an absorbed fourth power \(x^4\), it absorbs every normal term whose content is contained in \(c(R)\).

Proof. For each \(y\in c(R)\), choose a right summand containing \(y\). If it is linear, use \((A_1)\). Otherwise commute its factors to write it as \(yv\) with \(v\) nonempty, and use \((A_2)\) with \(z=v\). This adds \(y^4\). For an arbitrary word \(w\) with content in \(c(R)\), (U), substituted at \(w\), gives \(w\le w^4\), and the preceding fourth-power normalization gives \(w^4\le R\). All substitutions are nonempty. \(\square\)

The lemma handles exactly the part of a proof in which the original four-element algebra has value \(o\) identically. It does not introduce variables outside the right polynomial.

\par\smallskip\noindent\textit{3. Complete absorption criteria.}\quad 
Assume first that \(R\) contains no fourth power. Write \(L\) for its linear variables, \(E\) for its undirected quadratic graph with loops, \(D\) for its cube variables, and \(V=V(E)\cup D\). Isolated cube vertices are included.

For \(A_0\), if \(L\cap V\ne\varnothing\), call \(R\) collapsed. For \(A_1\), add \(L\cap V(E)\) to \(D\) using (I), and delete every linear term whose variable is now in \(D\), using (H). In the latter case use the resulting graph and mark set in what follows.

A permissible coloring is a proper two-coloring of the graph in which every marked vertex has color A. A vertex is forced A if it has color A in every permissible coloring; two vertices are forced opposite if their colors differ in every such coloring. Components without marks may be colored independently.

If \(D\ne\varnothing\) and there is no permissible coloring, also call \(R\) collapsed. The full criteria, always with the additional requirement \(c(w)\subseteq c(R)\), are as follows.

\textbullet{} A collapsed polynomial absorbs every term.

\textbullet{} If \(D=\varnothing\) and the graph is nonbipartite, it absorbs all quadratic terms on \(V\), no cubes or fourth powers, and just its original linear terms.

\textbullet{} In the remaining cases, it absorbs precisely the quadratics on forced-opposite vertices and the cubes on forced-A vertices. There are no fourth-power absorptions. Its linear absorptions are \(L\) for \(A_0\), and \(L\) together with the forced-A vertices for \(A_1\).

An original fourth power puts \(R\) in the first case by the propagation lemma. An empty graph causes no nonlinear absorptions.

Proof of the criteria. By (2) it suffices to consider \(B_i\). A polynomial with a nonlinear summand can have non-top value only \(b\) or \(c\). For value \(b\), it can have no cubic; every graph vertex is \(a\) and every separate linear variable is \(b\). A linear variable incident with the graph prevents this assignment. In \(B_0\), a linear variable incident with either an edge or a cube prevents both possible non-top values, proving the first collapse condition. In \(B_1\) such a linear variable forces color A whenever the value is \(c\), exactly as encoded by the added cube.

For value \(c\), graph vertices must be properly colored by \(a, b\), with cube vertices equal to \(a\); separate linear variables can be \(c\). Consequently a nonbipartite graph prevents value \(c\), and any cube also prevents value \(b\). In a bipartite graph, marked vertices on opposite sides likewise prevent value \(c\). These are exactly the stated collapse cases.

When \(b\) remains possible, its all-\(a\) graph assignment separates every cube and every additional linear variable on the graph. When \(c\) is possible, the independent choices of component colors separate every quadratic not forced opposite and every cube not forced A. A variable outside \(V\) can be set to \(c\) if it occurs only linearly, or to \(o\) if it does not occur, to separate nonlinear targets. An additional linear target takes \(a\) or \(b\) in graph assignments; neither is below \(c\) in \(B_0\), whereas precisely \(a\) is below \(c\) in \(B_1\). Purely linear polynomials are separated by assigning their variables \(a\), or a missing variable \(o\). These assignments prove necessity. The same value descriptions show directly that every listed target is absorbed, proving sufficiency.

\par\smallskip\noindent\textit{4. Unrestricted derivations.}\quad 
For \(A_0\), a linear-edge overlap gives fourth powers by (F). A linear-cube overlap is obtained from (F) by substituting \(y=x^2\) and using (S). The propagation lemma handles the rest. For \(A_1\), the preparation of \(D\) and \(L\) was already derived by (H),(I).

Repeatedly applying (P) shortens any odd walk to an edge. An odd closed walk therefore gives a loop \(x^2\). Law (L) propagates a loop to any edge vertex, including vertices in another component; (E) then gives every quadratic on \(V\). When any cube is present, \((Z_1)\) gives a fourth power and the propagation lemma applies. When marks occur on both sides of a bipartite component, (P) first joins two such marks, and \((Z_2)\) again produces fourth powers.

In a permissible coloring, (M) propagates cubes along even paths, giving every forced-A vertex. Odd paths give the opposite pairs within a component. For opposite vertices in distinct components, both colors must have been fixed by marks. If \(x\) is forced A and \(z\) forced B, choose a neighbor \(y\) of \(z\) on the A side. After adding cubes at \(x, y\), law (J) gives the edge \(xz\). This derives every quadratic in the criterion. For \(A_1\), (H) gives every additional linear term; no extra linear term is needed for \(A_0\).

Thus every valid absorption has a derivation, including the content-restricted collapsed cases. Normalize an arbitrary identity and apply the result to each summand in both directions. This proves the theorem for arbitrarily many variables and summands.

\subsection{Five proper subvarieties with explicit finite bases}\label{case:80}
We work without constants in the signature \((+,\cdot)\), relative to the additively idempotent semiring axioms. Inequality abbreviates absorption. This note supplies separate completeness proofs for five algebras; finite basability is not inferred from inclusion in a finitely based variety.

\par\smallskip\noindent\textit{1. Constructions and the bases.}\quad 
Addition in each construction is the join in the specified order. Unlisted products in the first three constructions are \(z\).

\textbullet{} \(R=\{e, a, c, b, z\}\): \(a, e<b<z\) and \(a<c<z\), with \(a\parallel e\), \(e\parallel c\), \(b\parallel c\). Here \(e^2=e\) and every product of two members of \(\{a, c\}\) is \(b\). The catalog map in this order is \((2,\allowbreak{}1,\allowbreak{}3,\allowbreak{}0,\allowbreak{}4)\), giving type 14368.

\textbullet{} \(S=\{e, a, c, b, z\}\): \(e, a<b<z\) and \(c<z\), with \(e\parallel a\) and \(c\) incomparable with \(e, a, b\). Here \(e^2=e\), \(a^2=c^2=b\). The map \((2,\allowbreak{}1,\allowbreak{}3,\allowbreak{}0,\allowbreak{}4)\) gives type 14425.

\textbullet{} \(T_0=\{e, a, b, d, z\}\): \(e, a, b\) are three incomparable elements below \(d<z\). Here \(e^2=e\) and \(a^2=b\). The map \((1,\allowbreak{}2,\allowbreak{}3,\allowbreak{}0,\allowbreak{}4)\) gives type 14768.

\textbullet{} \(U=\{a, c, b, z, t\}\): \(a, c<b<z<t\) and \(a\parallel c\). Within \(\{a, c, b, z\}\), put \(a^2=c^2=b\) and all other products equal to \(z\); \(t\) is absorbing for multiplication. The map \((2,\allowbreak{}3,\allowbreak{}1,\allowbreak{}0,\allowbreak{}4)\) gives type 2350.

\textbullet{} \(W=\{a, c, b, z, t\}\) has the same multiplication as \(U\), and order \(c<b<z<t\), \(a<z\), with \(a\) incomparable with \(c, b\). The map \((1,\allowbreak{}3,\allowbreak{}2,\allowbreak{}0,\allowbreak{}4)\) gives type 2822.

All multiplications are commutative. In the first three algebras a product of at least three elements equals \(e\) if all factors are \(e\), and \(z\) otherwise. In the last two it is \(z\) if no factor is \(t\), and \(t\) otherwise. These formulas also prove associativity. For distributivity in \(R\), the only nonconstant translations are the map taking \(e\) to \(e\) and everything else to \(z\), and the map taking \(a, c\) to \(b\) and everything else to \(z\). They preserve each indicated join. In \(S\) and \(T_0\) replace the second map by the maps taking just a designated active element to \(b\) and everything else to \(z\). Inspection of the listed incomparable joins proves preservation. In \(U, W\), a translation by an active element maps that element to \(b\), other elements below \(t\) to \(z\), and \(t\) to \(t\); the remaining translations are constant \(z\) below \(t\) (with \(t\) sent to \(t\)), or constant \(t\). These maps also preserve the listed orders and incomparable joins. Thus the constructions are additively idempotent semirings.

Let \(\Gamma\) consist of the five identities

\CaseDisplay{\begin{array}{rll}
(C)&xy\approx yx,\\
(P)&x^4\approx x^3,\\
(D)&xyz\approx x^3+y^3+z^3,\\
(E)&xy\leq xz+zy,\\
(M)&x^3\leq xy+y^3.
\end{array}}

Define the following additional laws:

\CaseDisplay{\begin{array}{rll}
(I)&x\leq x^3,\\
(J)&y\leq x+xy,\\
(G)&y\leq x+x^2+y^2,\\
(Q)&y\leq x^3+y^2,\\
(F)&x^3\leq x+y^3,\\
(A)&x\leq xy,\\
(K)&xy\leq x+y+x^2+y^2.
\end{array}}

Theorem 1. Complete bases are

\noindent Algebra\enspace\textbar\enspace Basis\enspace\textbar\enspace Number of displayed laws\par
\noindent \(R\)\enspace\textbar\enspace \(\Gamma,(I)\)\enspace\textbar\enspace 6\par
\noindent \(S\)\enspace\textbar\enspace \(\Gamma,(I),(J)\)\enspace\textbar\enspace 7\par
\noindent \(T_0\)\enspace\textbar\enspace \(\Gamma,(I),(G),(Q)\)\enspace\textbar\enspace 8\par
\noindent \(U\)\enspace\textbar\enspace \(\Gamma,(F),(A)\)\enspace\textbar\enspace 7\par
\noindent \(W\)\enspace\textbar\enspace \(\Gamma,(F),(I),(J),(K)\)\enspace\textbar\enspace 9\par
These are upper bounds, not assertions of minimality. In particular (M) may be redundant in bases containing (F), but it is retained uniformly.

\par\smallskip\noindent\textit{2. Normal forms and notation.}\quad 
Laws (C), (P), (D) give \(x^n=x^3\) for \(n\geq3\) and

\CaseDisplay{(uv)^3=u^3v^3\approx u^9+v^6+v^3\approx u^3+v^3.}

Consequently the cube of a word is the sum of the cubes of its variables. Applying (D) after splitting any word of length at least three into three nonempty factors gives the same decomposition for that word itself. Every polynomial is thus a sum of linear terms, quadratic terms, and cubes:

\CaseDisplay{F=\sum_{x\in L}x+\sum_{\{x, y\}\in E}xy+\sum_{x\in D}x^3.}

Here \(E\) is an undirected graph allowing loops, \(V\) is its vertex set, and \(C_F=L\cup V\cup D\). Empty blocks are omitted. On \(V\cup D\) use the edge components, adding isolated vertices for \(D\setminus V\). Define

\CaseDisplay{M_D=\text{union of components meeting }D,
\qquad
M_L=\text{union of edge components meeting }L\cap V.}

Write \(x\sim_E y\) when \(x, y\in V\) belong to the same edge component. Thus \(x\sim_E x\) requires \(x\in V\).

\par\smallskip\noindent\textit{3. Exact absorption criteria.}\quad 
Lemma 2. For \(R, S, T_0\), cubic absorption is exactly \(x^3\leq F\Longleftrightarrow x\in M_D\), and quadratic absorption is

\CaseDisplay{xy\leq F\quad\Longleftrightarrow\quad
 x\sim_E y\ \text{or}\ \{x, y\}\subseteq M_D.}

Their linear absorption sets are, respectively,

\CaseDisplay{\begin{array}{rl}
R:&L\cup M_D,\\
S:&L\cup M_D\cup M_L,\\
T_0:&L\cup M_D\cup
 \begin{cases}V,&D\ne\varnothing\text{ or }L\cap V\ne\varnothing,\\
 \varnothing,&\text{otherwise}.\end{cases}
\end{array}}

For \(U, W\), if \(D\ne\varnothing\), a normal summand is absorbed exactly when all its variables belong to \(C_F\). If \(D=\varnothing\), no cube is absorbed, and the other criteria are

\noindent Algebra\enspace\textbar\enspace Linear absorption set\enspace\textbar\enspace Condition for \(xy\leq F\)\par
\noindent \(U\)\enspace\textbar\enspace \(L\cup V\)\enspace\textbar\enspace \(x\sim_E y\)\par
\noindent \(W\)\enspace\textbar\enspace \(L\cup M_L\)\enspace\textbar\enspace \(x\sim_E y\) or \(\{x, y\}\subseteq M_L\)\par
Proof for \(R, S, T_0\). With a nonlinear summand on the right, a purely linear ideal cannot contain its value. The right value in \(R, S\) therefore belongs to \(\{e, b, z\}\); in \(T_0\) it belongs to \(\{e, b, d, z\}\). If \(F\) has no nonlinear summand, assigning every variable \(a\) excludes all nonlinear absorptions; a missing linear variable is separated by assigning it \(z\) and all other variables \(e\).

For \(R, S\), the condition \(F\leq b\) requires each edge component to be in the \(e\) state or an active state. In \(R\) an active component may use \(a, c\) independently; in \(S\) it is uniformly \(a\) or uniformly \(c\). Cubes force their components to be \(e\). In \(S\), a linear variable on a component excludes the \(c\) state, but still permits the \(a\) state. In \(R\), such a linear variable excludes \(c\) at that variable only, not at its neighbors. For \(T_0\), use \(F\leq d\): edge components are uniformly \(e\) or uniformly \(a\), and cubes force the \(e\) state. Linear terms impose no further restriction on these component states. The smaller value \(F=e\) forces every right variable to be \(e\). In \(T_0\), \(F=b\) is possible only when \(D=\varnothing\), \(L\cap V=\varnothing\), all edge variables are \(a\), and every linear variable is \(b\).

These descriptions prove sufficiency of the nonlinear criteria. For necessity, a vertex outside \(V\cup D\) can be assigned \(b\) and every other variable \(e\): the left cube or quadratic involving that vertex is \(z\), but the right side is at most \(b\) in \(R, S\), or at most \(d\) in \(T_0\). An unmarked component can instead be assigned \(a\), with the remaining components \(e\), separating its cubes and any proposed edge to a different component. This proves the nonlinear equivalences.

For linear terms in \(R\), a vertex \(x\) in an unmarked component with \(x\notin L\) is separated by assigning \(x=c\), all other vertices in that component \(a\), and all other variables \(e\). The right value is at most \(b\), whereas \(c\nleq b\). In \(S\), if \(x\) belongs to a component meeting neither \(D\) nor \(L\), assign that whole component \(c\) and all other variables \(e\). This gives the same separation. Conversely the allowed components of \(S\) have states \(e\) or \(a\) when \(F\leq b\), so their linear variables are absorbed. In both algebras a missing vertex outside \(C_F\) is separated by \(z\) versus \(e\).

In \(T_0\), if \(D\ne\varnothing\) or \(L\cap V\ne\varnothing\), the right side cannot be \(b\). Its possible values are \(e, d, z\), and every edge variable is then below its value. Without either condition, assign all edge variables \(a\) and all linear variables \(b\). This makes \(F=b\) and separates every edge variable not already linear. Vertices of \(M_D\) are absorbed by their cubes since \(x\leq x^3\) holds. Variables outside \(C_F\) again use the \(z/e\) separation. These arguments establish precisely the displayed linear sets.

Proof for \(U, W\). Every variable occurring in a term, if assigned \(t\), forces that term to be \(t\). Thus no absorption can introduce a variable outside \(C_F\). If \(D\ne\varnothing\), then \(F\geq z\). A left summand on \(C_F\) is at most \(z\), unless one of its variables is \(t\); in that case \(F=t\) too. This proves the first assertion.

Suppose \(D=\varnothing\). Assign all variables \(c\); then \(F\leq b\) and every cube is \(z\), so no cube is absorbed. If \(E=\varnothing\), a quadratic is also excluded by that assignment because \(F=c<b\), and only the original linear summands can be absorbed. If \(E\ne\varnothing\), \(F\) has one of the values \(b, z, t\). The interesting value is \(b\): all edge components must be uniformly \(a\) or uniformly \(c\). In \(U\) both colors are below \(b\). In \(W\) only \(c\) is below \(b\), so a component meeting a linear variable is forced to be \(c\). Hence the displayed linear and quadratic conditions suffice. An endpoint outside \(V\) can be assigned \(b\) (or \(t\) when missing from \(C_F\)); its product is not absorbed. Two different components not both forced to \(c\) can be colored differently to separate their product. An unmarked component in \(W\) can be colored \(a\), with all others \(c\), to separate its linear variables from \(F=b\). These also prove necessity. \(\square\)

\par\smallskip\noindent\textit{4. Equational derivation of every absorption.}\quad 
The word formulas prove (C), (P), (D). Lemma 2 verifies (E): the right edges connect the left endpoints. It also verifies (M): for \(R, S, T_0\) the right cube marks its edge component; for \(U, W\) a cube on the right absorbs every word on its content. The same lemma directly verifies each extra law in Theorem 1: (I) is the cube bound, (J) propagates a linear mark along an edge, (G) and (Q) supply the respective two gates in \(T_0\), (F) supplies a cubic summand in \(U, W\), (A) is their stated growth law for \(U\), and (K) connects two linearly marked loop vertices in \(W\). Thus all basis laws are valid without assuming their completeness.

In every basis, (E) closes paths and creates loops: \(x^2\leq xy+yx\approx xy\). Law (M) propagates a cubic mark along an edge. Moreover

\CaseDisplay{xy\leq x\,x^2+x^2y\approx x^3+y^3}

follows from (E) with its middle variable replaced by \(x^2\), and (D). Consequently \(\Gamma\) derives all nonlinear absorptions in Lemma 2 for \(R, S, T_0\). Law (I) yields every linear variable in \(M_D\).

For \(R\) this is sufficient. In \(S\), (J) propagates linear terms along paths and supplies exactly the additionally required \(M_L\). In \(T_0\), if \(x\in L\cap V\), the loop \(x^2\) is already derivable below \(F\); for any \(y\in V\) so is \(y^2\). Law (G) gives \(y\leq F\). If instead a cube \(x^3\) is a summand of \(F\), (Q) gives \(y\leq x^3+y^2\leq F\) for all \(y\in V\). These prove every stated linear absorption. Both gates are needed in the present proof: a cubic in one component may make a linear variable in a different unmarked component absorbable.

For \(U, W\), suppose first that \(D\ne\varnothing\) and choose \(t_0\in D\) (here \(t_0\) is a variable, not the named element \(t\)). For each right summand \(r\), (F) with \(x=r, y=t_0\) gives

\CaseDisplay{r^3\leq r+t_0^3\leq F.}

The normal-form formula for \(r^3\) yields \(x^3\leq F\) for every \(x\in C_F\). The derived quadratic bound yields every required quadratic. Law (I) yields the required linear terms in \(W\); in \(U\), (A) implies \(x\leq x^2\leq x^3\) and gives the same conclusion.

When \(D=\varnothing\), path closure gives all within-component edges. In \(U\), (A) extracts each vertex from any incident edge, producing \(L\cup V\). In \(W\), (J) produces \(L\cup M_L\). Every \(x\in M_L\) also has a loop below \(F\); (K) therefore gives \(xy\leq F\) for any \(x, y\in M_L\). There are no cubic absorptions to prove in this case.

We have derived all the absorptions in the exact criterion. Normalizing both sides of any valid identity and proving the absorption of each summand by the other side gives both inequalities and hence the identity. All substitutions use nonempty terms; there is no bound on the variable set, path lengths, components, or summand counts. Theorem 1 follows. \(\square\)

\par\smallskip\noindent\textit{5. Varieties generated by the catalogue algebras.}\quad
The five sources and their equal-variety catalog classes are

\textbullet{} 14368: 14368;

\textbullet{} 14425: 14425, 14883;

\textbullet{} 14768: 14768;

\textbullet{} 2350: 2350, 14404;

\textbullet{} 2822: 2822, 14405.

Each transfer requires complete membership tests in both directions. The finite diagnostics check all normal polynomials through three variables and specified tests through five variables. The unrestricted completeness proof is Section 4, not the diagnostic search. These eight types, together with the seven types treated by the 284/\allowbreak{}332 bases, close the previously unresolved five-element members of the five-law component branch. This does not assert that the entire branch is hereditarily finitely based.

The source manuscript for this case credits Aifa Wang and Lili Wang, School of Mathematical Sciences, Chongqing University of Technology.

The source manuscript for this case acknowledges support from the National Natural Science Foundation of China (12371024), the Science and Technology Research Program of Chongqing Municipal Education Commission (KJZD-K2024011102) and the Chongqing Natural Science Foundation Innovation and Development Joint Fund (Municipal Education Commission) (CSTB2025NSCQ-LZX0067).

\subsection{二阶群核心：仿射奇偶判据与两组完整有限基}\label{case:81}
\par\smallskip\noindent\textit{1. 两个半环和有限基.}\quad 
取二阶群 \(\{e, g\}\)，其中 \(g^2=e\)。在 \(B=\{e, g, c, t\}\) 上令不同群元素之和为 \(c\)，且 \(e, g<c<t\)；\(c\) 是 \(\{e, g, c\}\) 内的乘法零，\(t\) 是整个 \(B\) 的乘法零。取新元素 \(d\)，置

\CaseDisplay{\rho(d)=t,\qquad \rho|_B=\mathrm{id},\qquad
uv=\rho(u)\rho(v).}

加法均为下列偏序的并：

\CaseDisplay{\begin{array}{ll}
A_e:&e, g<c<t,\quad e<d<t,\quad g\parallel d,\ c\parallel d;\\
A_g:&e, g<c<t,\quad g<d<t,\quad e\parallel d,\ c\parallel d.
\end{array}}

映射 \(\rho\) 保持并并固定其像，故该构造保持乘法结合律和双侧分配律。元素名称不作为语言常量。

在顺序 \((e, g, c, d, t)\) 下，到完整五元素目录的同构分别为

\CaseDisplay{A_e\longrightarrow3969:\ (1,\allowbreak{}2,\allowbreak{}0,\allowbreak{}3,\allowbreak{}4),\qquad
A_g\longrightarrow3993:\ (2,\allowbreak{}1,\allowbreak{}0,\allowbreak{}3,\allowbreak{}4).}

定理 1。 \(A_e\) 的完整恒等式基为

\CaseDisplay{\begin{aligned}
xy&\approx yx, &(1)\\
x^3y&\approx xy, &(2)\\
x^2y^2&\approx x^2+y^2, &(3)\\
x&\leq x^3, &(4)\\
xyz&\leq x^3+y^3+z^3, &(5)\\
x^3&\leq x+xy. &(6)
\end{aligned}}

在 (1)--(6) 外增加

\CaseDisplay{x^3\leq x+y^2 \quad\textup{(\ensuremath{7})}}

即得 \(A_g\) 的完整恒等式基。

有效性。 两因子以上的乘积均在 \(B\) 中，出现 \(d, t\) 时取 \(t\)；否则按二阶群及其零 \(c\) 计算。这给出 (1)、(2)。平方值为 \(e, c, t\)，这些元素的乘法是链 \(e<c<t\) 的并，故 (3) 成立。(4) 仅在 \(x=d\) 时不是等式，此时为 \(d<t\)。对于 (5)，若出现 \(d, t\)，右侧为 \(t\)；若出现 \(c\)，右侧不低于 \(c\)；否则三个群元素不全相同时右侧为 \(c\)，全相同时其三重积等于共同值。对于 (6)，只有 \(x=d\) 需要说明，此时 \(xy=t\)。最后，在 \(A_g\) 中，每个平方均不低于 \(e\)，且 \(d+e=t\)，所以 (7) 有效。证毕。

\par\smallskip\noindent\textit{2. 任意变量数的判据.}\quad 
固定有限非空变量集 \(X\)。对非空词 \(w\)，记内容为 \(c(w)\)，记指数奇偶向量为

\CaseDisplay{\pi(w)=(m_x(w)\bmod2)_{x\in X}\in\mathbb F_2^X.}

对于非空多项式 \(P=\sum_i v_i\)，定义

\CaseDisplay{U(P)=\bigcup_i c(v_i),\qquad
N(P)=\bigcup_{|v_i|\geq2}c(v_i),\qquad
H(P)=\operatorname{aff}_{\mathbb F_2}\{\pi(v_i):i\}.}

仿射包由奇数个向量之和组成，允许重复；不使用空积或空和。

定理 2。 对任意非空词 \(w\) 和多项式 \(P\)：

1. 在 \(A_e\) 中，\(w\leq P\) 当且仅当 \(\pi(w)\in H(P)\)，且当 \(|w|\geq2\) 时有 \(c(w)\subseteq N(P)\)。

2. 在 \(A_g\) 中，\(w\leq P\) 当且仅当 \(\pi(w)\in H(P)\)，并且当 \(|w|\geq2\) 时

\CaseDisplay{c(w)\subseteq U(P),\qquad
   c(w)\subseteq N(P)\ \text{或}\ 0\in H(P).}

以下先证明必要性，再从有限基推导满足这些条件的每个吸收式，从而同时证明两个定理。

\par\smallskip\noindent\textit{2.1 仿射条件的必要性.}\quad 
若 \(\pi(w)\notin H(P)\)，有限维线性代数给出 \(\alpha\in\mathbb F_2^X\) 和 \(b\in\mathbb F_2\)，使

\CaseDisplay{\alpha\cdot\pi(v_i)=b\quad\text{对所有 }i,\qquad
\alpha\cdot\pi(w)=b+1.}

将 \(x\) 赋为 \(g^{\alpha_x}\)，所有右词取同一群元素 \(g^b\)，左词取另一个。两个不同群元素在加法序中不可比，故吸收失败。

\par\smallskip\noindent\textit{2.2 内容条件的必要性.}\quad 
设 \(|w|\geq2\)。在 \(A_e\) 中，若 \(x\in c(w)\setminus N(P)\)，赋 \(x=d\)，其余变量全为 \(e\)。左词取 \(t\)；右侧取 \(e\) 或 \(d\)，所以吸收失败。

在 \(A_g\) 中，若 \(x\in c(w)\setminus U(P)\)，同样的赋值使左右分别取 \(t, e\)。现在设 \(x\in c(w)\cap(U(P)\setminus N(P))\) 且 \(0\notin H(P)\)。变量 \(x\) 在 \(P\) 中只作为一次项出现。删去该项得到词族 \(P_0\)；其仿射包也不含零，故存在只用其余变量的群赋值，使 \(P_0\) 的全部词均取 \(g\)。这里使用事实：不含零的仿射子空间上存在恒为 1 的线性泛函。若 \(P_0\) 为空，则不需此赋值约束。再赋 \(x=d\)，右侧为 \(d\)，左侧为 \(t\)，仍失败。

\par\smallskip\noindent\textit{3. 三个等式推导引理.}\quad 
以下各引理仅使用 (1)--(6)，不依赖生成元的阶数。

引理 3（非线性词规范形）。 任意非线性词由内容和奇偶向量唯一确定模 (1)、(2) 的等价类。具体地，含至少两个不同变量时，各正偶指数约化为 2，各奇指数约化为 1；只含一个变量时，非线性偶次幂约化为平方，非线性奇次幂约化为立方。因此每个非线性词 \(u\) 满足 \(u^3\approx u\)，以及 \(ux^2\approx u\) 对每个 \(x\in c(u)\)。

证明。 式 (2) 删除同一变量的两个出现，只要删除后仍保留该变量及至少一个因子；交换律允许把它们移到相邻位置。单变量情形分别以平方、立方为终点。上述两个结论保持内容和奇偶向量。证毕。

引理 4（加入偶次变量）。 若 \(v\) 非线性、\(c(u)\subseteq c(v)\)、\(\pi(u)=\pi(v)\)，则 (1)--(4) 推出 \(u\leq v\)，其中 \(u\) 可以是一次词。

证明。 先设 \(u\) 非线性，选 \(y\in c(u)\)。用 (3) 乘以 \(u\)，再用引理 3，得到

\CaseDisplay{u+ux^2=uy^2+ux^2=uy^2x^2=ux^2.}

因此可以任意加入偶次变量而使词不减。补入 \(c(v)\setminus c(u)\) 的每个变量的平方，所得词由引理 3 等于 \(v\)。若 \(u=x\) 为一次词，先由 (4) 得 \(x\leq x^3\)，再应用上述结论。证毕。

引理 5（奇数乘积）。 若 \(q_1,\ldots, q_{2r+1}\) 都是非线性词，则

\CaseDisplay{q_1\cdots q_{2r+1}\leq q_1+\cdots+q_{2r+1}}

由 (1)、(2)、(5) 推出。于是，若 \(Q\) 是非空非线性词族，\(c(w)\subseteq\bigcup_{q\in Q}c(q)\) 且 \(\pi(w)\in H(Q)\)，则 \(w\leq\sum Q\) 可由 (1)--(5) 推出。

证明。 第一项对 \(r\) 归纳；在 (5) 中依次代入前 \(2r-1\) 项的乘积、最后两项，并用每个非线性词等于其立方。第二项选择给出 \(\pi(w)\) 的奇数项表示，再将 \(Q\) 中每个词各补入两次。最终乘积 \(v\) 的内容等于 \(Q\) 的全部内容，奇偶向量仍为 \(\pi(w)\)。由引理 4 和第一项，\(w\leq v\leq\sum Q\)。证毕。

\par\smallskip\noindent\textit{4. \(A_e\) 的完整性.}\quad 
记 \(L\) 为 \(P\) 的一次项变量集。对每个 \(x\in L\cap N(P)\)，选择一个含 \(x\) 的非线性右词，借助交换律写为 \(xu\)，其中 \(u\) 非空。(6) 给出

\CaseDisplay{x^3\leq x+xu\leq P.}

结合 (4)，可在 \(P\) 内把这一项 \(x\) 替换成 \(x^3\)，保持多项式等价。设替换后的非线性词族为 \(Q\)。其内容为 \(N(P)\)；其余一次项是 \(L\setminus N(P)\)。这些变量在所有其他右词中都不出现。

若左词为 \(x\in L\)，吸收式显然成立。其他符合定理 2 条件的左词，其奇偶向量在 \(L\setminus N(P)\) 的每个坐标均为零。因而它在 \(H(P)\) 中的任一奇数项表示，都使用每个外部一次项偶数次。删去这些成对项后，仍为奇数项表示，且只使用 \(Q\)。于是 \(\pi(w)\in H(Q)\)，特别地 \(Q\ne\varnothing\)。

左词非线性时，条件已给出 \(c(w)\subseteq N(P)\)。左词是未列入 \(L\) 的变量时，仿射表示也保证它在 \(N(P)\) 中。引理 5 因而推出 \(w\leq\sum Q\leq P\)。这证明所有有效吸收式均可由 (1)--(6) 推出。

\par\smallskip\noindent\textit{5. \(A_g\) 的完整性.}\quad 
若 \(0\notin H(P)\)，定理 2 的条件与 \(A_e\) 相同，第 4 节已足够。若 \(0\in H(P)\)，先按第 4 节将 \(L\cap N(P)\) 的一次项提升为立方，得到非线性词族 \(Q\)。零向量的仿射表示不能奇数次使用 \(L\setminus N(P)\) 的任一变量，故 \(0\in H(Q)\)，且 \(Q\ne\varnothing\)。

取一个奇数项乘积，其奇偶向量为零，并将每个 \(Q\) 中的词补入两次。由引理 3、5，该乘积等于一个非空词的平方 \(h^2\)，且

\CaseDisplay{h^2\leq\sum Q\leq P.}

对每个剩余一次项 \(x\)，应用 (7)：

\CaseDisplay{x^3\leq x+h^2\leq P.}

因此 \(P\) 等价于一个全部由非线性词组成的多项式 \(P^*\)，满足

\CaseDisplay{U(P^*)=U(P),\qquad H(P^*)=H(P).}

对符合判据的左词，使用引理 5 得 \(w\leq P^*=P\)。若左词为一次词，内容包含关系已由其奇偶向量属于仿射包保证。完整性得证。

一般多项式恒等式 \(P\approx R\) 等价于每个左词被 \(R\) 吸收且每个右词被 \(P\) 吸收，故吸收式的完整性给出全部恒等式的完整性。

\par\smallskip\noindent\textit{6. 结论的范围与核验.}\quad

\subsection{Eight complete identities for semiring 10379}\label{case:82}
Signature (+, multiplication), with no constants. Relative to the additively idempotent semiring axioms, the following is a finite basis. Write \(u\leq P\) for \(u+P=P\).

{\small\par\smallskip
\CaseLine{\(xy=yx; xyz=xyzz; xyzt\leq xyz; xxx\leq x+xy;\)}
\CaseLine{\(xxx\leq x+yyy; xxy\leq xy+zzz; xy\leq xxy; xy\leq xx+yz.\)}
\smallskip}

\par\smallskip\noindent\textit{Word values and exact absorption.}\quad 
Addition is 00004/\allowbreak{}01004/\allowbreak{}00224/\allowbreak{}00234/\allowbreak{}44444; multiplication is 22222/\allowbreak{}23222/\allowbreak{}22222/\allowbreak{}22222/\allowbreak{}22224. The additive order has \(3<2<0<4\) and \(1<0<4\). A nonlinear word is 4 if every variable in its content is assigned 4. Otherwise a quadratic is 3 when both occurrences are assigned 1, and every other nonlinear word is 2. In particular long words retain their content; their value is not a constant cube. Literal values remain arbitrary.

For an original polynomial P, let L be its literal variables, E its unordered quadratic edges, and V the endpoints of E. Say G holds if P contains an original word of length at least three or L intersects V. For any nonempty target word w the exact criterion for \(w\leq P\) is:

\textbullet{} A literal target must itself occur in P.

\textbullet{} A long target requires G and an original word u with c(u) contained in c(w).

\textbullet{} A quadratic requires an original word u with c(u) contained in c(w), and either G or c(w) contained in V.

Here original means before any derived summands are added.

Necessity is witnessed uniformly in the number of variables. For a missing literal x, assign \(x=1\) and all other variables 3. Every right summand is 2 or 3, while 1 is below neither. If no contained right word exists for a nonlinear target, assign its content 4 and every other variable 3. Every right summand is at most 2; the target is 4. If G fails, assign V to 1 and every other variable 3. All original quadratics and literals then have value 3. A long target is 2, and a quadratic not wholly in V is also 2. These separate every failed case.

For sufficiency, any original long word takes value at least 2. An incident original literal x and quadratic xy also have join at least 2: when \(xy=3\), \(x=1\) and \(1+3=0\); otherwise \(xy\geq 2\). Thus G forces \(P\geq 2\). The contained witness supplies 4 whenever the target is 4. All remaining target nonlinear values are at most 2. If G fails but a quadratic has both endpoints in V, every right nonlinear value is at least 3. The only way P does not dominate 2 is that all right summands equal 3. This forces every V variable to 1, so the target is 3. A contained right word again handles target value 4.

\par\smallskip\noindent\textit{Complete finite derivations.}\quad 
The first two word identities equate exactly the long words having the same content. To see the claimed direction constructively, commute any desired occurrence to the last position. A long word has two nonempty preceding blocks, so \(xyz=xyzz\) can duplicate that occurrence. After duplicating as needed, reverse this identity to remove surplus copies, never reducing the word below length three. Hence a long word normalizes to its sorted content when there are at least three distinct variables, to xxy for sorted two-variable content, and to xxx for singleton content. This also proves \(u^{3}=u\) for every long word u.

If an original long word exists, retain it as a long witness g. If G holds only through an incident literal x and edge xy, derive \(g=x^{3}\) using \(xxx\leq x+xy\). The content of this global witness need not be inside the target, and must not be mistaken for the contained witness.

Choose an original contained word u. If u is literal, derive \(u^{3}\) by \(xxx\leq x+yyy\), substituting u for x and g for y; \(g^{3}\) normalizes to g. If \(u=ab\) is quadratic, derive aab by \(xxy\leq xy+zzz\), substituting g for z. If u is long, retain u itself. In all cases there is now a derived long word v with \(c(v)=c(u)\). If variables of w are missing from v, append their nonempty product by \(xyzt\leq xyz\), splitting v into three nonempty blocks. Normalize to the long word on c(w). This is the target when w is long. For a quadratic \(w=ab\), recover ab from aab by \(xy\leq xxy\). Every substitution is nonempty.

If G fails, only quadratic targets remain. An original target needs no derivation. A new target ab has a contained original word, which cannot be long. It cannot be a literal either: a literal at an endpoint in V would contradict failure of G. Therefore it is an original loop aa or bb. The other target endpoint occurs in an original edge bc (up to commutativity). Derive ab using \(xy\leq xx+yz\), swapping endpoints if needed. This completes every semantic case.

The six inequalities are valid by the exact criterion; their finite assignments are also checked independently. The argument covers arbitrary variable sets, lengths and summand counts. Distributivity reduces all terms to polynomials; summandwise absorption in both directions yields every valid identity. Thus these eight laws form a complete finite basis. Other types may use this proof only after verified full HSP equivalence.

\subsection{Six complete identities for semiring 4095}\label{case:83}
Constant-free signature (+, multiplication), relative to the additively idempotent semiring axioms. A complete finite basis is

{\small\par\smallskip
\CaseLine{\(xy=yx; xyz=xyzz; xyzt\leq xyz;\)}
\CaseLine{\(xxy\leq x+xy; xy\leq xxy; xy\leq xx+yy.\)}
\smallskip}

The notation \(u\leq P\) abbreviates \(u+P=P\).

\par\smallskip\noindent\textit{Values and target-contained criterion.}\quad 
Addition is 00044/\allowbreak{}01144/\allowbreak{}01234/\allowbreak{}44334/\allowbreak{}44444. Element 2 is additive bottom, with \(2<1<0<4\) and \(2<3<4\). Multiplication is 00200/\allowbreak{}00200/\allowbreak{}22222/\allowbreak{}00210/\allowbreak{}00200. A nonlinear word is 2 if any variable in its content is 2. Otherwise a quadratic is 1 precisely when both occurrences are 3; every other nonlinear word is 0.

Fix a target w and let U be precisely the original words of P whose content is contained in c(w). Define L, E, V using only U: its literals, unordered quadratic edges and their endpoints. Let G mean that U contains a long word or L intersects V. Then \(w\leq P\) iff:

\textbullet{} A literal target occurs literally in P.

\textbullet{} A long target satisfies G.

\textbullet{} A quadratic target satisfies G, or its edge is in E, or both its endpoint squares occur in E.

The restriction to U is essential. Words with an outside variable may be killed by assigning that variable additive bottom 2. An original long word elsewhere in P is not a global long witness for this algebra.

For necessity of the literal condition, set the target variable 3 and outsiders 2. If the literal is absent, surviving summands are pure powers with values 1 or 0; neither dominates 3. For a failed long target, assign V to 3 and every other target variable 1, outsiders 2. Since G fails, U consists of quadratics valued 1 and literals 1; outsiders kill all other words. The right value is at most 1 while the long target is 0.

For a failed quadratic, G fails and its edge is absent. A square target xx then has neither original square nor original long on x; assign \(x=1\) and outsiders 2. Surviving literals are 1, the target 0. For distinct endpoints x, y, at least one original loop is missing, say yy. Assign \(y=1\), \(x=3\), outsiders 2. The only possible surviving quadratic is xx, which has value 1. If xx occurs, original literal x cannot occur because G fails; all surviving literals are then 1. If xx does not occur, both literal values are at most \(3+1=4\), so refine the assignment to \(x=y=1\); there are now no surviving quadratics and every surviving literal is 1. Either way the target is 0 and the right side is at most 1. This covers every failed case for arbitrary P and w.

For sufficiency a nonlinear target containing a variable assigned 2 is itself bottom. Otherwise every target variable avoids2. A contained long word is 0. A contained literal x and edge xy have join at least 0: if \(xy=1\), then \(x=3\) and \(3+1=4\); otherwise \(xy=0\). Thus G suffices. Original target edges are immediate. If both endpoint loops occur, at least one is 0 unless both endpoints are 3, in which case the target and both loops are 1. This proves the exact criterion.

\par\smallskip\noindent\textit{Unrestricted finite derivation.}\quad 
The two word laws normalize all long words by content. Commute a chosen occurrence to the end, duplicate it with \(xyz=xyzz\) using two nonempty prefix blocks, and reverse the law to delete surplus repetitions without shortening below three. Canonical forms are sorted content for at least three variables, xxy for two, and xxx for one.

When G holds, retain a contained original long u, or derive xxy from a contained original literal x and edge xy by \(xxy\leq x+xy\). The resulting long word has content within the target. Append all missing target variables using \(xyzt\leq xyz\), with three nonempty left blocks, then normalize. This derives every long target. Recover a quadratic ab from aab by \(xy\leq xxy\). If G fails and a nonoriginal quadratic is absorbed, both endpoint loops are original; apply \(xy\leq xx+yy\). Literal targets are already original. All substitution words are nonempty.

All six laws satisfy the exact criterion and are independently checked on all assignments. The proof has no bound on variables, lengths or summands. Distributive expansion and summandwise absorption in both directions give every valid identity. Other catalogue types require independent full HSP equivalence for transfer.

\subsection{Five complete identities for semiring 2941}\label{case:84}
In the constant-free signature (+, multiplication), relative to the additively idempotent semiring axioms, a complete finite basis is

{\small\par\smallskip
\CaseLine{\(xy=yx; xyz=xyzz; xyzt\leq xyz; xy\leq xxy; xy\leq xx+yy.\)}
\smallskip}

Here \(u\leq P\) means \(u+P=P\).

Addition is 00004/\allowbreak{}01114/\allowbreak{}01224/\allowbreak{}01234/\allowbreak{}44444, the \(chain 3<2<1<0<4\). Multiplication is 00030/\allowbreak{}00030/\allowbreak{}00130/\allowbreak{}33333/\allowbreak{}00030. A nonlinear word is 3 if any variable in its content is 3. Otherwise a quadratic is 1 when both occurrences are 2, and every other nonlinear word is 0.

For a target w, restrict original P to U, the original words whose content is contained in c(w). Let E be the quadratic edges in U. Then:

\textbullet{} A literal target is absorbed iff it occurs literally in P.

\textbullet{} A long target is absorbed iff U contains an original long word.

\textbullet{} A quadratic is absorbed iff U contains an original long word, its edge is original, or both endpoint loops are original.

Necessity uses explicit separating assignments. For a missing literal, assign its variable 4 and outsiders 3. Every surviving nonlinear word is 0, which is below 4. For a long target with no contained original long word, assign all target variables 2 and outsiders 3. All surviving quadratics are 1 and literals 2, so the right side is at most 1 while the target is 0. For a failed square target, no contained square or long word exists; the same assignment gives right value at most 2 and target 1. For a distinct-endpoint target xy, at least one original endpoint loop is missing, say yy. If xx is present set \(x=2\), \(y=1\), outsiders 3. All surviving quadratics are 1, and both possible literal values are at most 1, while the target is 0. If neither loop is present, assign \(x=y=1\), outsiders 3; every surviving summand is literal 1 while the target is 0. Original xy is absent in both cases. These exhaust failed conditions for arbitrary P.

For sufficiency, a nonlinear target with any variable 3 is additive bottom. Otherwise a contained long word supplies 0, which dominates the target's possible nonlinear values 0 and 1. Original edges are immediate. With both original loops, if the target is 1 then both are 1; otherwise at least one loop is 0 and dominates the target. Original literals are immediate. Thus the criterion is exact without any size restriction.

The two word identities normalize long words by content. Commute a chosen variable to the end and duplicate with \(xyz=xyzz\) using two nonempty prefix blocks; reverse the law to remove surplus occurrences while retaining length at least three. Canonical forms are sorted content for at least three variables, xxy for two, and xxx for one. Given a contained original long word, append missing target content by \(xyzt\leq xyz\) and normalize. This derives a long target directly; recover a quadratic ab from aab using \(xy\leq xxy\). If no contained long word is available, a new quadratic requires both original loops; apply \(xy\leq xx+yy\). Literal targets already occur. All substitution words are nonempty.

Every law satisfies the criterion and is independently checked on all assignments. Distributivity followed by summandwise absorption in both directions proves completeness for arbitrary variables, words and summands. In particular \(xxy\leq x+xy\) fails at \(x=y=2\): its values are 0 and 1. An incident literal cannot replace the original long-word requirement. Further catalogue types need independent full HSP equivalence to transfer.

\subsection{Finite bases for 6022 and 6023 by length ideals and three original witnesses}\label{case:85}
Signature (+, multiplication), without constants. The two bases below are finite effectively specified families, not enumerated small lists. Their different addition tables must be retained.

Both multiplication tables have rows 01000,\allowbreak{}11111,\allowbreak{}41344,\allowbreak{}41444,\allowbreak{}41444. For 6022 addition rows are 00444,\allowbreak{}01234,\allowbreak{}42244,\allowbreak{}43434,\allowbreak{}44444. For 6023 they are 00444,\allowbreak{}01444,\allowbreak{}44244,\allowbreak{}44434,\allowbreak{}44444. Read entries in element order 0,\allowbreak{}1,\allowbreak{}2,\allowbreak{}3,\allowbreak{}4. In 6022 the element 1 is below 0,\allowbreak{}2,\allowbreak{}3, whereas in 6023 it is below 0 and incomparable with 2,\allowbreak{}3.

Any word containing a1-valued factor has value 1. Without such a factor, a literal keeps its value; a nonlinear word headed by 0 is 0; a quadratic word with both factors 2 is 3; every other nonlinear word is 4. In particular every word of length at least 3 with no 1 is 0 if its head is 0, and 4 otherwise. This follows by induction from the rows.

\par\smallskip\noindent\textit{Exact absorption criteria.}\quad 
For target w let \(C=c(w)\), \(h=head(w)\), and U be the original right words whose contents are contained in C. Set \(V=U\) for 6022, and \(V=P\) for 6023. Let L be the set of literal variables in V, and Q the union of contents of quadratic words in V. Then \(w\leq P\) holds exactly when:

H: U contains a word headed by h; N: if w is nonlinear, V contains a nonlinear word; G: at least one of the following holds: V contains a word of length at least 3; L intersects Q; w is literal and h belongs to L; or w is quadratic and C is contained in Q.

All the words in these conditions are original. In 6022 the length and crossing witnesses must be contained in the target content. In 6023 their variables may lie outside it. H is contained in both.

For completeness of the semantic proof, consider principal additive ideals. Both tables have ideals\{0,\allowbreak{}1\} and\{1\}. If \(P\leq 0\) and the target has no 1 factor, H prevents its head being outside\{0,\allowbreak{}1\}, since every contained same-head word would then have value 2,\allowbreak{}3, or 4. Thus \(w\leq 0\). If \(P\leq 1\) every right word contains 1; the contained H witness forces 1 into w, so \(w=1\). Necessity of H is separated by \(h=4\), other C variables 0, and outsiders 1: without H every original is 0 or 1, while \(w=4\).

For 6022 the remaining proper ideals are\{1,\allowbreak{}2\} and\{1,\allowbreak{}3\}. With \(P\leq 2\), a nonlinear target without 1 is impossible by N: a contained nonlinear witness has value 0,\allowbreak{}3, or 4. A literal target is 2 or 1 as required, because its contained H witness forces this or prevents \(P\leq 2\). If N fails for a nonlinear target, assign C variables 2, outsiders 1: all contained originals are literals 2, so \(P\leq 2\) while w is 3/\allowbreak{}4.

With \(P\leq 3\) in 6022, and assuming w has no 1, each word of U has value 3. There can be no long word; literal variables must be 3, and every quadratic variable must be 2. Thus L and Q are disjoint. G forces a literal target to be 3, or a quadratic target to have both factors 2 and hence value 3; it precludes a long target. Conversely, if G fails, assign L variables 3, Q variables 2, outsiders 1. These assignments are consistent. All words of U are 3 and other originals are 1, so \(P\leq 3\). For a failed literal branch assign its variable 4 unless it is already in Q, when it is 2. For a failed quadratic branch assign any target variable outside Q the value 4 unless it is already in L, when it is 3. Any remaining C variables may be 4. A failed long target is 4. In all cases w is not below 3. H ensures U is nonempty.

For 6023 the two remaining ideals are the singletons\{2\} and\{3\}. \(P=2\) is possible only if every original right word is literal. Then N excludes a nonlinear target and H gives value 2 to a literal target. If N fails, assign all variables 2; \(P=2\) and the nonlinear target is 3/\allowbreak{}4. \(P=3\) is possible only when every original is literal 3 or quadratic with both variables 2. Thus no long original can occur and L and Q must be disjoint. G gives exactly the same target alternatives as above. For necessity of G, assign all global L variables 3 and all global Q variables 2, and use 4 on unassigned variables. The target violates the ideal exactly as in the 6022 case. \(P=4\) needs no condition. This proves both criteria for all alphabets, word lengths and sums.

\par\smallskip\noindent\textit{Six structural identities and normal words.}\quad 
{\small\par\smallskip
\CaseLine{\(xyz = xzy\)}
\CaseLine{\(xxy = xyy\)}
\CaseLine{\(xxyz = xyz\)}
\CaseLine{\(xyyz = xyz\)}
\CaseLine{\(xyzz = xyz\)}
\CaseLine{\(xxxx = xxx\)}
\smallskip}

The value formula proves these identities under every nonempty substitution. They keep length at least 3 and preserve content/head. Words of length 1 or 2 stay unchanged. A long normal word has its original head, followed by the increasing list of other content variables. If this would have length 1, use hhh; if it would have length 2, use hyy. Thus the long normal form is determined by content and head, and always remains long.

Here is a literal normalization procedure. Sort the tail by adjacent swaps using the first identity, placing head copies first. While length exceeds 3, delete an adjacent duplicate: head duplicates use \(xxyz=xyz\); an internal duplicate uses \(xyyz=xyz\); final duplicates use \(xyzz=xyz\) with a prefix of at least 2 positions. All required images are nonempty. At length 3 use \(xxy=xyy\) when the two initial variables agree and the third differs. A singleton content gives hhh. Sorting and these deletions give the asserted normal form in all cases.

\par\smallskip\noindent\textit{At most three original witnesses suffice.}\quad 
Choose one original contained H witness. If G uses a long word, adjoin that word; it also supplies N. If G uses a crossing q in L intersect Q, adjoin literal q and an original quadratic containing q; the quadratic supplies N. A literal target using its literal branch already belongs to P. Finally, for a quadratic target whose at most two variables belong to Q, select one original quadratic containing each variable, besides H. These quadratics supply N. Every selection has at most 3 original words and still satisfies H/N/G. In 6022 all selected words have the required content containment; in 6023 the global witnesses remain unrestricted. Their sum therefore absorbs w.

\par\smallskip\noindent\textit{Reversible compression and a fixed finite complete family.}\quad 
Normalize the target and the at most 3 selected words. Reserve every variable of each short word and the head of each long word. Across at most 4 words this reserves at most 8 variables. Every other variable appears only in long words and away from their head. Group those variables by nonempty incidence patterns among the words, at most 15 groups. Identify each group to a distinct fresh variable and normalize again, keeping short words short and long words long. The result uses at most 23 variables, word length at most 25, and at most 3 upper words. It is valid because variable identification and the structural laws preserve the original valid absorption.

For the inverse substitution, each reserved variable returns to its original singleton; each group variable returns to the nonempty product of its members. Short words contain only reserved variables and return exactly. Every long expansion remains of length at least 3 and recovers its original head and content. The literal normalization procedure therefore returns it to the original long normal word. This proves an inverse equational derivation, including the two-variable long case where the compressed tail is repeated.

For each core take the six structural laws and all table-valid absorptions over a fixed 23-variable alphabet, length at most 25, and at most 3 upper words. The criterion effectively specifies this finite family. The selected-original-word reduction, compression and inverse nonempty substitution derive every valid absorption. Adding the other original summands recovers its full upper polynomial. Distributive expansion and both absorption directions derive every valid identity. There is no bound on the original variables, lengths, or summand count. Thus 6022,\allowbreak{}6023 and their opposites are finitely based. Catalogue transfer still requires independently verified complete mutual HSP.

\subsection{Finite bases for 2947 and 13188 by contained head domains}\label{case:86}
Work in the semiring signature without constants. For each of these two algebras take the semiring axioms, the four word identities

{\small\par\smallskip
\CaseLine{\(xyz = xzy,  xxyz = xyz,  xyyz = xyz,  xxy = xyy,\)}
\smallskip}

and all table-valid absorptions with at most 184 variables, word length at most 186 and at most 24 upper words. This is a specified finite family; no small enumerated basis is claimed. The proof below is unrestricted in alphabet, word length, polynomial size and polynomial substitutions.

For 2947 the multiplicatively absorbing and additively least element is \(Z=3\); put B=\{0,\allowbreak{}1,\allowbreak{}2,\allowbreak{}4\}. A word containing Z has value Z. Otherwise a quadratic has value 4 if its head is 4, value 1 on (2,\allowbreak{}2), and value 0 otherwise. Every word of length at least three has value f(head), where \(f(4)=4\) and \(f(0)=f(1)=f(2)=0\).

For 13188, \(Z=2\) and B=\{0,\allowbreak{}1,\allowbreak{}3,\allowbreak{}4\}. A word containing Z has value Z. Otherwise a quadratic has value 3 if its head is 3, value 0 on (1,\allowbreak{}1), and value 4 otherwise. A long word has \(f(3)=3\) and \(f(0)=f(1)=f(4)=4\). These formulas follow by multiplying triples and checking \(f(b)c=f(b)\) for every b, c in B. Literal words keep their assigned value.

The four laws are verified directly in each table. They normalize every long word using only its head and content. Keep the head h first. With one content variable use hhh; with two use hht, where t differs from h; otherwise use h followed by the remaining distinct variables in a fixed order. Tail permutations use the first law. While length exceeds three, an extra head is deleted by the second law, and a repeated tail letter by the third law, after a tail permutation. At length three the fourth law changes htt into hht. No short word is replaced by a long word.

Let w be a target and P a nonempty polynomial of original upper words. Let U be the upper words whose content is contained in c(w). Failure of \(w\leq P\) forces all target variables to lie in B, since otherwise \(w=Z\) is least. Reassign every outside variable to Z. Every upper word outside U then becomes Z; upper values cannot increase from the failing assignment. Conversely, a failure for U on B extends this way to a failure for P. Thus only U matters, but its words retain their length and their heads.

For each proper principal ideal I, the allowed quadratic pairs on B form a rectangle \(H_{I}\) times \(T_{I}\), possibly empty. Let \(K_{I}\) be the allowed heads of long words. The complete table of these domains is:

\noindent Source\enspace\textbar\enspace I\enspace\textbar\enspace \(H_{I}\)\enspace\textbar\enspace \(T_{I}\)\enspace\textbar\enspace \(K_{I}\)\par
\noindent 2947\enspace\textbar\enspace \{0,\allowbreak{}1,\allowbreak{}2,\allowbreak{}3\}\enspace\textbar\enspace \{0,\allowbreak{}1,\allowbreak{}2\}\enspace\textbar\enspace B\enspace\textbar\enspace \{0,\allowbreak{}1,\allowbreak{}2\}\par
\noindent 2947\enspace\textbar\enspace \{1,\allowbreak{}2,\allowbreak{}3\}\enspace\textbar\enspace \{2\}\enspace\textbar\enspace \{2\}\enspace\textbar\enspace empty\par
\noindent 2947\enspace\textbar\enspace \{2,\allowbreak{}3\}\enspace\textbar\enspace empty\enspace\textbar\enspace empty\enspace\textbar\enspace empty\par
\noindent 2947\enspace\textbar\enspace \{3\}\enspace\textbar\enspace empty\enspace\textbar\enspace empty\enspace\textbar\enspace empty\par
\noindent 13188\enspace\textbar\enspace \{0,\allowbreak{}1,\allowbreak{}2\}\enspace\textbar\enspace \{1\}\enspace\textbar\enspace \{1\}\enspace\textbar\enspace empty\par
\noindent 13188\enspace\textbar\enspace \{1,\allowbreak{}2\}\enspace\textbar\enspace empty\enspace\textbar\enspace empty\enspace\textbar\enspace empty\par
\noindent 13188\enspace\textbar\enspace \{2\}\enspace\textbar\enspace empty\enspace\textbar\enspace empty\enspace\textbar\enspace empty\par
\noindent 13188\enspace\textbar\enspace \{2,\allowbreak{}3\}\enspace\textbar\enspace \{3\}\enspace\textbar\enspace B\enspace\textbar\enspace \{3\}\par
For every target variable start with domain B. An original literal x intersects its domain with I. A quadratic xy intersects x with \(H_{I}\) and y with \(T_{I}\) (both restrictions apply when x=y). A long original with head x intersects x with \(K_{I}\). These independent unary constraints are equivalent to putting every word of U in I. The system is impossible exactly when some domain is empty. Otherwise the target is forced into I exactly when its literal domain is contained in I, its quadratic head and tail domains are contained in \(H_{I}\) and \(T_{I}\), or its long head domain is contained in \(K_{I}\), according to its degree. Testing the four proper ideals is necessary and sufficient for \(w\leq P\): a failure has upper sum b in a proper principal ideal down(b), and a failure of one such ideal also makes w fail to lie below the actual upper sum.

This gives a bounded set of original witnesses. A variable has only four upper roles: literal, quadratic head, quadratic tail, long head. For each target variable relevant to its value (at most two, or only the head for a long target), retain one original witness of every role it has. For each of the four impossible ideal systems choose one empty-domain variable and retain one original witness of every role it has. At most six variables have their roles retained, so at most 24 upper words suffice. Possible systems preserve the domains needed by the target; impossible systems preserve the selected empty domain. Consequently this selected polynomial Q still absorbs w. In particular there are at most 18 selected short words and at most six selected long words, because only one long-head witness per selected variable is used.

Normalize w and Q. Reserve every variable of a short selected word, and the head and up to two other distinct variables of each long word. Also reserve the target's variables if it is short, or its head and up to two other distinct variables if long. This reserves at most \(36+18+3=57\) variables. All other variables occur only in the at most seven long words (including the target). Group them by their common incidence among those words. There are at most \(2^{7}-1=127\) nonempty groups, hence at most 184 variables after identification. The reserved variables preserve every head, every short occurrence and the distinction between short and long. Normalized length is at most \(184+2=186\).

Identification preserves validity, so the compressed absorption belongs to the stated finite family. Replace each group variable by the nonempty word listing its members, and each reserved variable by itself. Every member of a group has exactly the same incidence in all retained words. All heads are reserved, and all affected words remain long. The four word laws therefore normalize each expanded word back to its original normal form. This derives \(w\leq Q\); adding the unselected original summands derives \(w\leq P\). Expanding arbitrary terms by distributivity and splitting each identity into summand absorptions in both directions proves the finite basis claim. No bound on the original variables or words is used.

Opposite and equal-variety transfers require separate complete HSP certificates. Short polynomial agreement alone is not a transfer proof.

\subsection{A finite basis from prefix-final crossing}\label{case:87}
Let A be five-element catalogue entry 8627. The signature is (+, multiplication), without constants. All bases are relative to the axioms of additively idempotent semirings. Write w \ensuremath{\leq} P for w+P approximately equal to P. Numerical labels denote elements only.

\par\smallskip\noindent\textit{Operations and arbitrary word values.}\quad 
Addition is join in the order 2 \ensuremath{<} 1 \ensuremath{<} 0 \ensuremath{<} 4 and 2 \ensuremath{<} 1 \ensuremath{<} 3 \ensuremath{<} 4, where 0 and 3 are incomparable. In numerical element order the tables are

{\small\par\smallskip
\noindent addition       multiplication\par
\CaseLine{\ttfamily 0 0 0 4 4\qquad 0 2 2 3 4}
\CaseLine{\ttfamily 0 1 1 3 4\qquad 2 2 2 2 2}
\CaseLine{\ttfamily 0 1 2 3 4\qquad 2 2 2 2 2}
\CaseLine{\ttfamily 4 3 3 3 4\qquad 4 2 2 4 4}
\CaseLine{\ttfamily 4 4 4 4 4\qquad 4 2 2 4 4}
\smallskip}

Put K=\{0,\allowbreak{}3,\allowbreak{}4\}. For every word w of length at least two:

\textbullet{} if any occurrence has value outside K, then w has value 2;

\textbullet{} if all occurrences have value 0, then w has value 0;

\textbullet{} if all nonfinal occurrences have value 0 and the final occurrence has value 3, then w has value 3;

\textbullet{} in all other all-K cases its value is 4.

For linear words the value is the assigned element. On all-K assignments the same three latter cases apply also to linear words, with an empty prefix.

To check associativity, any outside-K factor in a triple forces 2 under either association, since 1 times any factor, any factor times 1, and every product involving 2 are 2. For an all-K triple the two associations give 0 precisely when all factors are 0, and 3 precisely when the first two are 0 and the last is 3; otherwise both give 4. Each row and column is increasing in the displayed order and preserves the only incomparable join \(0+3=4\), as the tables show. Thus multiplication preserves every binary join. This proves the semiring axioms, and induction gives the word-value formula.

For a nonempty word w write c(w) for content, p(w) for the set of variables before its final position, and f(w) for its final variable. In particular p(y) is empty for a linear word y. Given target w, put \(C=c(w)\), retain only right summands U=\{u in P : c(u) subset C\}, and define \(D=\cup\) c(u), \(J=\cup\) p(u), R=\{f(u):u in U\}.

\par\smallskip\noindent\textit{An exact absorption criterion.}\quad 
Lemma. A linear y is absorbed by P exactly when y is a literal summand. For a nonlinear target w, A satisfies w \ensuremath{\leq} P exactly when

\CaseDisplay{D=C,\qquad J\cap R\ne\varnothing\quad\hbox{or}\quad p(w)\subseteq J.
 \quad\textup{(\ensuremath{1})}}

Proof. If y is not a literal summand, assign \(y=1\) and every other variable 2. The target is 1, whereas every right summand is

2. This proves the linear assertion.

For nonlinear targets, assigning every outside variable 2 kills every discarded summand. If U is empty, assign C the value 0; then \(w=0\) and \(P=2\). If U is nonempty but D misses x in C, assign \(x=4\), every other C-variable 0, and outsiders 2. Then \(w=4\) and the contained sum is 0. Hence \(D=C\) is necessary.

Assume \(D=C\) and J disjoint from R. Assign J the value 0, R the value 3, and outsiders 2. Every retained summand is 3. Since \(C=J\) union R and R is nonempty, the target is 3 exactly when its prefix is contained in J (and then its final lies in R). If its prefix is not contained in J, it is 4. This proves necessity of the second condition.

For sufficiency, a target with any outside-K factor is 2, the least element. Otherwise every C-variable is in K, so all retained terms are in K and their nonempty sum is 0,\allowbreak{}3, or 4. A sum of 4 absorbs the target. A sum of 0 requires each retained word to have every variable 0; \(D=C\) then forces \(w=0\). A sum of 3 requires every retained term to be 3, since \(0+3=4\). Thus J is assigned 0 and R is assigned 3, and these sets are disjoint. Condition (1) forces p(w) subset J. Also f(w) is in R: otherwise \(C=p(w)\) union \{f(w)\} would be contained in J, contradicting nonempty R and \(D=C\). Consequently \(w=3\). All cases give w below the retained sum and hence below P. \(\square\)

\par\smallskip\noindent\textit{Normal words and two-right-term derivations.}\quad 
The following two word identities hold by the word-value formula:

\CaseDisplay{xyz\approx yxz,\qquad x^2y\approx xy.       \quad\textup{(\ensuremath{W})}}

The first permutes any adjacent positions before the final position, using their nonempty remaining suffix for z and a possible prefix context. The second removes any repeated prefix occurrence after these permutations make two copies adjacent. Therefore each nonlinear word is equivalent to the sorted product of p(w), followed by f(w). Linear words remain linear. Empty products mean omitted substrings, not unit constants. Normalization preserves p(w) and f(w).

A word is saturated when its final variable occurs in its prefix. For nonempty D, fix g(D) to be the sorted product of D followed by its least variable. This word is saturated and has content D.

Take a valid nonlinear w \ensuremath{\leq} P and normalize all words. If J meets R, choose original contained summands a, b and x with x in p(a) and \(f(b)=x\). They may coincide. Let \(S=c(a)\) union c(b). The criterion gives g(S) \ensuremath{\leq} a+b: both contents cover S and the crossing x remains. After a guard g(S) has been derived below P, any further u in U gives the valid step

\CaseDisplay{g(S\cup c(u))\le g(S)+u.                 \quad\textup{(\ensuremath{2})}}

The guard itself supplies a crossing, and the two contents cover the enlarged target. Eventually \(D=C\) gives g(C) \ensuremath{\leq} P. Since g(C) supplies its own crossing, the criterion also gives w \ensuremath{\leq} g(C). Thus each step uses at most two right terms.

If J and R are disjoint, condition (1) implies that every right word ends in \(s=f(w)\) and omits s from its prefix. Indeed all variables of C except possibly s lie in p(w) subset J, so the nonempty set R, disjoint from J and contained in C, is \{s\}. Begin with any u in U. Merge two available normal words by taking the sorted product of the union of their prefix sets, followed by s. This merged word is below their sum by (1), since their prefixes cover its prefix and their contents cover its content. If both prefixes are empty, all three words are the literal s and no new step is needed. Successively merging the original terms gives prefix J and final s. Here \(J=C\) minus \(\{s\}=p(w)\), so the last normal word is w. Again every step uses at most two right terms. Linear absorptions are literal inclusions and need no new axiom.

\par\smallskip\noindent\textit{Reversible compression and the finite basis.}\quad 
In a local absorption there are at most three normal words in total. Keep their final variables separate, at most three variables, including the variable of each linear word. Every other variable occurs only in prefixes. Group these variables by their nonempty prefix-incidence patterns across the words, at most seven groups. Identifying each group to one new variable and normalizing by (W) preserves validity. There are at most ten variables and each normal word has length at most eleven.

This compression is reversible modulo (W): substitute the nonempty product of a group's original members for its new variable, and restore each reserved final variable individually. Uniform prefix incidence recovers exactly the original prefix set of every word. The identities (W) sort and remove prefix repetitions. All final variables and linear words are preserved. Thus the original local absorption is an equational consequence of its compressed version and (W).

Theorem. A finite identity basis for A, relative to the semiring axioms, is (W) together with all valid normal absorptions over a fixed ten-variable alphabet, of word length at most eleven, and with at most two right terms. Reversing all words gives a basis for the multiplicative opposite.

Proof. The displayed family is finite and explicitly specified; the lemma decides its membership. All its members hold in A. Compression derives every valid two-right-term local absorption from this family. The constructions above derive every valid word absorption with no bound on variables, length, or number of summands. Distributivity expresses each term as a finite nonempty polynomial. For any valid equality of polynomials, derive the absorption of each summand into the other polynomial in both directions. This derives every identity. Word reversal transports the full proof to the opposite. \(\square\)

No small enumerated or minimal basis is asserted. Finite regressions are implementation checks, not the completeness proof. Transferring the result to another catalogue entry requires independently verified equality of generated varieties in both HSP directions.

\subsection{Contained covers and finite derivations by changes of the head}\label{case:88}
Let A be catalogue entry 11019. Work without constants, relative to the axioms for additively idempotent semirings, and write w \ensuremath{\leq} P for w+P approximately equal to P. Element labels are not term constants. Addition is maximum \(on 3<2<1<0<4\); the tables are

{\small\par\smallskip
\noindent addition       multiplication\par
\CaseLine{\ttfamily 0 0 0 0 4\qquad 4 0 0 3 4}
\CaseLine{\ttfamily 0 1 1 1 4\qquad 4 2 2 3 4}
\CaseLine{\ttfamily 0 1 2 2 4\qquad 4 2 2 3 4}
\CaseLine{\ttfamily 0 1 2 3 4\qquad 3 3 3 3 3}
\CaseLine{\ttfamily 4 4 4 4 4\qquad 4 4 4 3 4}
\smallskip}

The element 3 is multiplicatively absorbing. A nonlinear word has value 3 when some factor is 3. Otherwise it has value 4 if some factor is 4 or some occurrence after the first is 0; value 0 if its initial occurrence is 0 and every later occurrence is 1 or 2; and value 2 when all factors lie in\{1,\allowbreak{}2\}. These rules agree under both associations of a triple; each row and column is increasing in the chain order, so multiplication preserves joins. Thus A is an additively idempotent semiring. Induction proves the word formula at arbitrary length.

The word identities

\CaseDisplay{xyz\approx xzy,\qquad xy^2\approx xy  \quad\textup{(\ensuremath{W})}}

hold. They normalize every nonlinear word to its literal head followed by its sorted tail set. A head also occurring in its tail is retained twice. Linear words stay linear. Let c(u), h(u), t(u) denote content, head and tail set, with \(t(y)=empty\) for linear y.

\par\smallskip\noindent\textit{Exact criterion using only contained right summands.}\quad 
For a nonempty polynomial P and finite C put B\_C=\{u in P:c(u) subset C\}.

Lemma 1. A linear y is absorbed by P exactly when it is a literal summand. A nonlinear w of content C is absorbed exactly when

\CaseDisplay{C\subseteq\bigcup_{u\in B_C}c(u),\qquad
 t(w)\subseteq\bigcup_{u\in B_C}t(u).                     \quad\textup{(\ensuremath{1})}}

Proof. A missing literal y is separated by assigning \(y=1\) and all other variables 3. Every nonlinear power of y has value 2 and every other right word has value 3, so P is at most 2 while \(y=1\).

Necessity may be tested with every outsider assigned the bottom 3, which kills all noncontained words. For a variable x missing from the first union, assign \(x=4\) and all other C-variables 2. Then \(w=4\) and P is at most 2. For x in t(w) missing from the second union, assign \(x=0\), other C-variables 2 and outsiders 3. Then \(w=4\), but every contained right word is at most 0 because x can occur only as a simple head or literal. Thus \(P<4\).

For sufficiency only the contained words are needed. If a C-variable has value 3, \(w=3\) is automatic. Otherwise all variables in C are nonbottom. A C-variable assigned 4 occurs in a contained right word by the first coverage, giving right value 4. A variable assigned 0 in the target tail occurs in a contained right tail by the second coverage, also giving 4. In the remaining case all target tail variables are 1 or 2. If the head is 0, some contained word contains that head and has value 0 or 4; hence \(w=0\) is absorbed. If the head is 1 or 2, all C-variables are 1 or 2 and every contained nonlinear word is 2, while a contained linear word is at least 2. The first coverage ensures such a witness exists, and \(w=2\) is absorbed. \(\square\)

\par\smallskip\noindent\textit{A construction by local changes of the uncovered head.}\quad 
For a normal word u define its missing set c(u) minus t(u). This is empty when its head is repeated and otherwise is the singleton \{h(u)\}, including linear u.

Two known words u, v can be merged whenever

\CaseDisplay{|(c(u)\cup c(v))\setminus(t(u)\cup t(v))|\leq1.     \quad\textup{(\ensuremath{2})}}

Let C' and T' be the two unions. If the difference is \{r\}, choose r as head; if it is empty, choose any head in C'. Follow the head by the sorted T'. The resulting normal word q has content C' and tail T', and q \ensuremath{\leq} u+v is valid by Lemma 1: both right words are contained and cover exactly the required sets. If T' is empty, condition 2 forces \(u=v\) to be that same literal, already known.

Fix valid nonlinear w \ensuremath{\leq} P, put \(C=c(w)\), \(h=h(w)\), and \(B=B_{C}\). By Lemma 1, every x in C minus \{h\} lies in a tail of an original word of B. Starting from each u in B, build a known word q as follows. If its missing set is empty or \{h\}, stop. Otherwise its head is r different from h. Choose an original v in B with r in t(v), and merge q, v by 2. This is possible: all old content except r was already in the old tail, and the new tail now contains r; only the head of v can remain uncovered.

If the new missing set is again \{p\} with p different from h, then p was absent from all old content, since that entire content now lies in the new tail. Thus every nonterminal step introduces a new root variable outside the earlier content. There are only finitely many variables in C. After at most |C| steps we stop with a word whose missing set is contained in \{h\}. Every step preserves all earlier content and tail variables and uses at most two right terms. This argument includes loops, arbitrarily long directed cycles, disconnected supports and linear starting words.

Perform this construction for all original words in B. Any two resulting words can be merged by 2, since each missing set is contained in \{h\}. Merge them successively. The final guard Q has content C by the first coverage in 1 and has tail containing every original right tail. Thus t(w) is contained in t(Q) by the second coverage. Lemma 1 now gives the valid one-right-word step w \ensuremath{\leq} Q. This constructs every valid absorption from local valid absorptions with at most two right terms, using original or already derived right words only.

\par\smallskip\noindent\textit{Finite compression and completeness.}\quad 
At most three words occur in a local step. Reserve their heads individually, at most three variables, and group all other variables by nonempty incidence patterns in the three tails, at most seven groups. Identify each group with one symbol and normalize by(W). This is a substitution and preserves validity. It uses at most ten variables and normal word length at most eleven.

For the inverse substitution send each group symbol to the nonempty ordered product of its original members, retaining reserved heads. Common tail incidence restores the original tails; (W) sorts and reduces them back to the original normal words. Linear words and repeated heads are included; no empty products or units occur.

Theorem. The identities(W), together with all valid normal absorptions over a fixed ten-letter alphabet, of word length at most eleven and with at most two right terms, form a finite basis for A. Reversing all words gives a finite basis for its opposite.

Proof. The specified family is finite, with membership decided by Lemma 1. Reversible compression derives every local step from it. The construction proves unrestricted derivability of all valid word absorptions. Distribute arbitrary terms and derive each summand's absorption into the other side in both directions. This yields every identity. Word reversal carries derivations to the opposite. \(\square\)

No small enumeration, minimality or hereditary assertion is made. Finite computations are auxiliary; the proof handles arbitrary variables and summands. Other types need independently replayed complete HSP transfers in both directions before classification.

\subsection{内部尾部见证与图拼接：七条完整有限基}\label{case:89}
\par\smallskip\noindent\textit{1. 构造及词值.}\quad 
取 \(A=\{z, e, a, d, t\}\)，加法为偏序 \(z<e<a<t\)、\(d<e<a<t\) 上的并，\(z, d\) 不可比，且 \(z+d=e\)。定义乘法如下：\(z\) 为两侧乘法零；在其余元素上，右乘 \(a, t\) 得 \(t\)，右乘 \(e, d\) 时固定 \(e, a, t\)，并把 \(d\) 送到 \(e\)。符号顺序 \((z, e, a, d, t)\) 对应目录标号 \((3,\allowbreak{}1,\allowbreak{}0,\allowbreak{}2,\allowbreak{}4)\)。

右平移为三个保持并的映射：常值 \(z\) 映射 \(Z\)；固定 \(z, e, a, t\) 并把 \(d\) 送到 \(e\) 的 \(E\)；固定 \(z\) 并把其余元素送到 \(t\) 的 \(G\)。它们满足 \(E^2=E\)、\(G^2=G\)、\(EG=GE=G\)，\(Z\) 在复合下为零映射；直接由定义有 \(R_{xy}=R_xR_y\)、\(R_{x+y}=R_x\vee R_y\)。这证明结合律及两侧分配律。

非一次词的值为：出现 \(z\) 则取 \(z\)；没有 \(z\) 而尾部出现 \(a\) 或 \(t\) 则取 \(t\)；尾部全取 \(e, d\) 时，首字母取 \(a, t\) 则分别取 \(a, t\)，首字母取 \(e, d\) 则取 \(e\)。

写 \(D(w)\) 为首字母之后的变量集合。半环满足

\CaseDisplay{xyz=xzy,\qquad xy^2=xy,\qquad x^2y=y^2x.}

前两式给出首字母与尾部集合正规形。首字母重复的非一次词只由内容决定，记为 \(M_X\)；它是幂等词。

\par\smallskip\noindent\textit{2. 完整吸收判据.}\quad 
令 \(F\) 为非空右多项式，\(L\) 为一次项变量集。对非空集合 \(X\) 记

\CaseDisplay{\mathcal R_X=\{u\in F:|u|\ge2,\ c(u)\subseteq X\},}

\CaseDisplay{N_X=\bigcup_{u\in\mathcal R_X}c(u),\qquad D_X=\bigcup_{u\in\mathcal R_X}D(u).}

定理。 对非一次词 \(w\)，令 \(h=h(w)\)、\(X=c(w)\)，则

\CaseDisplay{w\le F\quad\Longleftrightarrow\quad D(w)\subseteq D_X\ \text{且} h\in L\cup N_X.\quad\textup{(\ensuremath{1})}}

一次词 \(x\le F\) 当且仅当存在同名一次项或仅含 \(x\) 的非一次右词。

必要性。 若 \(x\in D(w)\setminus D_X\)，赋 \(x=a\)、其余 \(X\) 内变量 \(d\)、外部变量 \(z\)。左值为 \(t\)。外部非一次右项取 \(z\)；内部非一次项不含尾变量 \(x\)，取值至多为 \(a\)；一次项也至多为 \(a\)，故失败。若 \(h\notin L\cup N_X\)，赋 \(h=a\)、其余 \(X\) 内变量 \(d\)、外部变量 \(z\)。左值为 \(a\) 或 \(t\)；右内部非一次项取 \(e\)，外部非一次项取 \(z\)，一次项至多为 \(e\)，仍失败。一次词取 \(x=d\)、其他变量 \(z\) 分离。

充分性。 因非一次左词尾部非空，尾部覆盖保证存在内部非一次右项。左值为 \(z\) 时，任意这样的项都在 \(z\) 之上。左值为 \(e\) 时，所有 \(X\) 内变量取 \(e, d\)，内部非一次项取 \(e\)。左值为 \(a\) 时，首字母取 \(a\)，第二条件提供同名一次项或包含该变量的内部非一次项，值为 \(a\) 或 \(t\)。左值为 \(t\) 时，若某个左尾变量取 \(a, t\)，其内部尾部见证直接取 \(t\)；否则首字母取 \(t\)，第二条件所给的一次项或内部非一次项也取 \(t\)。一次词用 \(x\le x^2=x^n\)（\(n\ge2\)）。证毕。

\par\smallskip\noindent\textit{3. 七条完整基.}\quad 
主定理。 在半环公理之外，完整基为

\CaseDisplay{\begin{aligned}
xyz&\approx xzy, &(N)\\
xy^2&\approx xy, &(R)\\
x^2y&\approx y^2x, &(C)\\
x&\le x^2, &(I)\\
xy&\le x+y^2, &(K)\\
xyz&\le xy+yz, &(P)\\
xyz&\le xy+xz. &(M)
\end{aligned}}

有效性由词值规则及 (1) 逐式验证。以下给出不限制变量数的完备性证明。

\par\smallskip\noindent\textit{3.1 三个拼接操作.}\quad 
首先，若 \(p=ag\) 和 \(u=gb\) 是已经推出不超过 \(F\) 的非一次词，\(a, b\) 非空，\((P)\) 给出

\CaseDisplay{agb\le ag+gb=p+u\le F.\quad\textup{(\ensuremath{2})}}

只要 \(g\in D(p)\)，即可利用 \((N),(R)\) 把 \(p\) 写成 \(ag\)。式 (2) 保留 \(p\) 的首字母，尾部及内容按集合并合。

其次，同首字母非一次词 \(hv, hv'\) 可由 \((M)\) 合并成 \(hvv'\le hv+hv'\le F\)。若其中一个词仅为 \(h\)，内容并已由另一个词实现，不作空词代入。

最后，若 \(p, q\) 均为幂等非一次词且已不超过 \(F\)，\((K)\) 给出 \(pq\le p+q^2=p+q\le F\)。更一般地，\(p\) 可以是一次词，只要求第二个因子 \(q\) 幂等。

\par\smallskip\noindent\textit{3.2 平方型词的闭合引理.}\quad 
设非空有限集合 \(X\) 中每个变量都在某个内部右项的尾部出现，即 \(X\subseteq D_X\)。证明 \(M_X\le F\)。

把每个内部右词视为一条从其首字母指向全部尾变量的有向超边。若有一个有向环 \(g_0\to g_1\to\cdots\to g_{r-1}\to g_0\)，逐条选出实现这些边的右词，用 (2) 拼接。最终所得词的首字母 \(g_0\) 再次出现于尾部，因此是一个平方型词；它不超过 \(F\)，内容包含该环的全部顶点和所用右词的全部变量。长度为一的环对应一个已经首字母重复的右词，直接使用即可。

考虑所有已能推出 \(M_Y\le F\) 且 \(\varnothing\ne Y\subseteq X\) 的集合 \(Y\)。它们在有限并下闭合，因为上述第三个操作把 \(M_Y, M_{Y'}\) 合成 \(M_{Y\cup Y'}\)。设这些集合的并为 \(Z\)；若没有，则仅把 \(Z\) 记为空集合，不赋予它空词含义。

若 \(Z\ne X\)，对每个 \(x\in X\setminus Z\) 选一个尾部含 \(x\) 的内部右词 \(u_x\)。其首字母不可能属于 \(Z\)：否则把已推出的 \(M_Z\) 以该字母结尾，运用 (2) 就会扩大 \(Z\)。所以每个外部顶点 \(x\) 都有一个仍在 \(X\setminus Z\) 的前驱。有限性保证存在全部位于该补集的有向环。环拼接得到一个包含补集顶点的平方型词，和 \(Z\) 的定义矛盾。因此 \(Z=X\)，闭合引理得证。

\par\smallskip\noindent\textit{3.3 首字母只出现一次的词.}\quad 
给定有效吸收 \(w\le F\)，先将左词正规化。若它是平方型，(1) 给出 \(X\subseteq D_X\)，由闭合引理完成。以下设首字母为 \(h\) 且 \(h\notin D(w)\)。若 \(h\in D_X\)，仍有 \(X\subseteq D_X\)，先得 \(M_X\le F\)，再用 \((I)\) 推出 \(w\le w^2=M_X\)。

现在设 \(h\notin D_X\)。根据 (1)，或有一次右项 \(h\)，或有内部右项包含 \(h\)；后者因 \(h\) 不在内部尾部中，只能以 \(h\) 为唯一首字母。故至少有一个内容包含于 \(X\)、首字母为 \(h\) 且不重复的词已经不超过 \(F\)。

这类可推导词的内容在有限并下闭合，使用同首字母合并操作。选一个内容最大的这样的词 \(p\)，记其内容为 \(U\)。若 \(U\ne X\)，对每个 \(x\in X\setminus U\) 选内部尾部见证 \(u_x\)。其首字母若为 \(h\)，用同首字母合并会扩大 \(U\)；若在 \(U\setminus\{h\}\)，则它已在 \(p\) 的尾部，用 (2) 同样会扩大 \(U\)。因此所有这些首字母都在 \(X\setminus U\)，从而又存在补集中的有向环。

环拼接给出平方型词 \(M_Y\le F\)，其内容包含补集变量，而且不含 \(h\)：环的所有首字母不等于 \(h\)，所有内部右项尾部也都不含 \(h\)。由 \((K)\)，\(pM_Y\le p+M_Y\le F\)；这个词仍以 \(h\) 为唯一首字母，却扩大了内容，矛盾。所以 \(U=X\)，\(p\) 的正规形就是 \(w\)。

一次左词由 \((I)\) 及判据处理。任意恒等式可展开为双向逐词吸收，至此完备性证毕。

\par\smallskip\noindent\textit{4. 方法意义与范围.}\quad 
环的长度没有上界；证明使用的是固定三变量拼接式和有限图的前驱闭合论证。因此不能把发现长环恒等式误当作非有限基证据，也不需要为每个环长加入一条新公理。

乘法对偶逐词反转上述证明。对具体目录对象仅通过独立的精确双向簇包含传递有限基性。本文没有证明整个右作用半格类遗传有限基，亦未处理所有剩余五元素类型。

\subsection{Head witnesses and terminal overlap: thirteen-identity bases for a five-element semiring and its opposite}\label{case:90}
The source manuscript for this case credits Aifa Wang and Lili Wang, School of Mathematical Sciences, Chongqing University of Technology.

\par\smallskip\noindent\textit{Abstract.}\quad 
We give a thirteen-identity basis for a five-element additively idempotent semiring. Its words retain the first variable and distinguish a final variable which occurs earlier from one which does not. Absorption requires a contained-content head witness and coverage of all target variables. When no right terminal variable occurs internally, internal-content coverage and a terminal witness are also required. Prefix merging and content projection prove completeness without bounds on variables or summands. Word reversal gives a basis for the multiplicative opposite.

\par\smallskip\noindent\textit{1. The algebra and the basis.}\quad 
We use the constant-free signature \((+,\cdot)\) relative to the axioms for additively idempotent semirings. Write \(u\leq v\) for \(u+v\approx v\). Let \(A=\{d, e, r, a, q\}\) with additive order

\CaseDisplay{e<d<q,\qquad a<r<d<q.}

Addition is join; \(e\) is incomparable with \(a, r\). Multiplication is:

\noindent \(\cdot\)\enspace\textbar\enspace \(d\)\enspace\textbar\enspace \(e\)\enspace\textbar\enspace \(r\)\enspace\textbar\enspace \(a\)\enspace\textbar\enspace \(q\)\par
\noindent \(d\)\enspace\textbar\enspace \(d\)\enspace\textbar\enspace \(d\)\enspace\textbar\enspace \(r\)\enspace\textbar\enspace \(r\)\enspace\textbar\enspace \(d\)\par
\noindent \(e\)\enspace\textbar\enspace \(d\)\enspace\textbar\enspace \(e\)\enspace\textbar\enspace \(r\)\enspace\textbar\enspace \(a\)\enspace\textbar\enspace \(d\)\par
\noindent \(r\)\enspace\textbar\enspace \(r\)\enspace\textbar\enspace \(r\)\enspace\textbar\enspace \(r\)\enspace\textbar\enspace \(r\)\enspace\textbar\enspace \(r\)\par
\noindent \(a\)\enspace\textbar\enspace \(r\)\enspace\textbar\enspace \(r\)\enspace\textbar\enspace \(r\)\enspace\textbar\enspace \(r\)\enspace\textbar\enspace \(r\)\par
\noindent \(q\)\enspace\textbar\enspace \(q\)\enspace\textbar\enspace \(q\)\enspace\textbar\enspace \(r\)\enspace\textbar\enspace \(r\)\enspace\textbar\enspace \(q\)\par
These element names are not constants in terms. Put \(K=\{d, e, q\}\). Within \(K\), \(d, q\) are left-zero elements and \(ee=e\), \(ed=eq=d\). This is associative: a first factor \(d\) or \(q\) fixes the product; with first factor \(e\), distinguish whether the second factor is \(e\) or one of \(d, q\). The element \(r\) is a two-sided zero; the row of \(a\) is constant \(r\); and \(ua=a\) exactly when \(u=e\), otherwise \(ua=r\). For a triple involving \(a\) before its last position, both associations give \(r\). For a triple ending in \(a\), both give \(a\) exactly when the first two factors are \(e\), otherwise \(r\). Together with the zero property this proves associativity on \(A\). Each row and column is order preserving and preserves the joins of the only incomparable pairs \((e, a)\) and \((e, r)\), as the table shows. Thus multiplication distributes over addition in both directions.

Theorem 1. The following thirteen identities form a basis for \(A\):

\CaseDisplay{\begin{array}{rll}
\mathrm{(M)}&xyzt\approx xzyt,\\
\mathrm{(N_1)}&x^2y\approx xy,\\
\mathrm{(N_2)}&xyx\approx xy^2,\\
\mathrm{(E)}&x\leq x^2,\\
\mathrm{(G)}&x^2\leq x+y^2,\\
\mathrm{(C)}&xyt\leq xt+ys,\\
\mathrm{(Q)}&xyt\leq xt+xyzt,\\
\mathrm{(B)}&x^2y^2\leq x^2+x^2y^2z^2,\\
\mathrm{(F)}&x^2y^2\leq x^2+y,\\
\mathrm{(D_{11})}&xy^2\leq xy+zyt,\\
\mathrm{(D_{10})}&xy^2\leq xy+yt,\\
\mathrm{(D_{01})}&y^2\leq y+zyt,\\
\mathrm{(D_{00})}&y^2\leq y+yt.
\end{array}\quad\textup{(\ensuremath{1})}}

Reversing every word in (1) gives a basis for \(A^{\mathrm{op}}\). We make no claim that thirteen is minimal.

\par\smallskip\noindent\textit{2. Word values and normal forms.}\quad 
For a word \(w\), let \(h(w)\), \(t(w)\), \(c(w)\) and \(i(w)\) be its first variable, final variable, content and internal content. Internal content consists of the variables in all positions except the last, and is empty for a linear word.

Lemma 2. A nonlinear word has the following values.

1. An occurrence assigned \(r\), or an internal occurrence assigned \(a\), gives value \(r\).

2. Otherwise a terminal \(a\) gives value \(a\) when every internal letter is assigned \(e\), and value \(r\) otherwise.

3. In the remaining case all variables lie in \(K\). A first letter \(q\) gives \(q\); all letters \(e\) give \(e\); every other case gives \(d\).

Proof. The zero \(r\) and the constant row of \(a\) prove the first case. A nonempty \(K\)-prefix has value \(e\) exactly when all its letters are \(e\), so the column of \(a\) proves the second. The description of \(K\) proves the third. \(\square\)

Fix a variable ordering. A nonlinear word is saturated when its final variable occurs earlier. With \(h=h(w)\), define normal forms

\CaseDisplay{\nu(w)=
\begin{cases}
\displaystyle h\left(\prod_{x\in c(w)\setminus\{h\}}x\right)h,
   &t(w)\in i(w),\\[4pt]
\displaystyle h\left(\prod_{x\in i(w)\setminus\{h\}}x\right)t(w),
   &t(w)\notin i(w).
\end{cases}\quad\textup{(\ensuremath{2})}}

The products are ordered; an empty product is an omitted substring, not a constant. Put \(\nu(x)=x\).

Lemma 3. The three equalities in (1) derive \(w\approx\nu(w)\). A saturated word depends only on its head and content, and satisfies \(w^2\approx w\).

Proof. Identity \((M)\) interchanges adjacent interior positions, preserving first and last positions. Repeated interior variables can therefore be made adjacent and removed by \((N_1)\); a repeated internal head is moved beside the initial head and removed likewise. If the terminal has no earlier occurrence, this gives (2).

If the terminal is already the head, the same operations give the first normal form. Otherwise, call that terminal \(s\). Reduce interior repetitions, then insert one extra head just after the initial head using \((N_1)\) backwards. Permute the interior to obtain \(hRh s^2\). Identity \((N_2)\) on the suffix \(h s^2\) gives \(hRhsh\). Move its redundant interior head beside the initial head and remove it by \((N_1)\). The result has head and tail \(h\) and the original content, giving the first normal form. The argument also covers an absent \(R\). Lemma 2 verifies the three equalities, and normalization of \(ww\) gives the last assertion. \(\square\)

Normalizing a saturated right word may change its terminal to its head. Its internal content remains its full content, and the existence of a terminal/internal overlap remains true. Without overlap no right word is saturated, so all terminals and internal contents are preserved.

\par\smallskip\noindent\textit{3. Absorption by a polynomial.}\quad 
Let \(P=\sum_{u\in\mathcal U}u\) be a nonempty polynomial and put

\CaseDisplay{J=\bigcup_{u\in\mathcal U}i(u),\qquad
T=\{t(u):u\in\mathcal U\},\qquad U=J\cup T.}

Lemma 4. The absorption \(A\models w\leq P\) holds exactly when:

1. Some right word \(u\) satisfies \(h(u)=h(w)\) and \(c(u)\subseteq c(w)\).

2. \(c(w)\subseteq U\).

3. If \(J\cap T\) is empty, additionally \(i(w)\subseteq J\) and \(t(w)\in T\).

Proof of necessity. If condition 1 fails, assign the target head \(q\), the other target variables \(e\), and outside variables \(r\). The target is \(q\). A right word containing an outside variable is \(r\); a right word entirely within the target content has a different head and value \(e\) or \(d\). Thus \(P\leq d\), a contradiction.

If condition 2 fails at \(x\), assign \(x=r\) and all other variables \(e\). The target is \(r\) and \(P=e\), again a contradiction.

Now suppose \(J\cap T\) is empty. If \(x\in i(w)\setminus J\), assign \(J\) to \(e\), \(T\) to \(a\), and \(x\) to \(a\) as well, with any remaining variables assigned \(e\). All right words are \(a\), but the internal \(a\) makes the target \(r\). To see the terminal requirement after internal coverage is established, assign \(J\) to \(e\) and \(T\) to \(a\). A target terminal outside \(T\) belongs to \(J\) by condition 2. The target is \(e\) and \(P=a\), which does not absorb \(e\).

Proof of semantic sufficiency. Let \(p\) be the value of \(P\). If \(p=q\), the result is immediate. If \(p=e\), every variable in \(U\) is assigned \(e\), and condition 2 gives \(w=e\). If \(p=a\), each right word is \(a\). Therefore all internal variables in \(J\) are \(e\), all terminal variables in \(T\) are \(a\), and the two sets are disjoint. Condition 3 makes the target \(a\).

If \(p=r\) and the target were outside the downset \(\{a, r\}\), all of its variables would lie in \(K\) by Lemma 2 (also for a linear target). Condition 1 supplies a right word entirely in \(K\), with value in \(K\), contradicting \(P=r\). Finally, if \(p=d\), only a target \(q\) could fail the inequality. Such a target has head \(q\) and content in \(K\), or is linear \(q\). Condition 1 then supplies a right word of value \(q\), contradicting \(P=d\). These cases exhaust \(A\). \(\square\)

The ten inequalities of (1) satisfy this criterion. \((E)\) has a saturated right word. For \((G)\) the linear \(x\) supplies the head witness and \(y^2\) supplies overlap. In \((C)\) and \((Q)\) the word \(xt\) supplies the head witness; the right words cover internal content and provide terminal \(t\) in the no-overlap case. \((B)\) and \((F)\) have the head witness \(x^2\), with overlap and full target content covered. Each \((D)\) has overlap at \(y\); its first right word supplies the contained-content head witness. Thus (1) is sound.

\par\smallskip\noindent\textit{4. Explicit derivations.}\quad 
Lemma 5. Every absorption in Lemma 4 is derivable from (1).

Proof. Normalize all words. Choose a right head witness \(u\) from condition 1 and write \(C=c(w)\).

First suppose \(J\cap T\) is empty. A linear target must already be a right summand: a longer head witness of singleton content would be saturated and create overlap. For a nonlinear target, put \(h=h(w)\) and \(t=t(w)\). Since \(h\in i(w)\subseteq J\), its head witness cannot be linear. Its terminal belongs both to \(T\) and to \(C=i(w)\cup\{t\}\), so disjointness forces that terminal to be \(t\). Thus \(u=Xt\), with \(X\) nonempty and \(i(u)\subseteq i(w)\).

Starting from \(u\), combine the current \(Yt\) with each nonlinear right word \(Zh\) using \((C)\):

\CaseDisplay{YZt\leq Yt+Zh\leq P.}

Both prefixes are nonempty. Normalize after each step. The head remains \(h(w)\) and the terminal remains \(t\), while internal content becomes the union of the current internal content and that of the chosen right word. We obtain a word below \(P\) with head \(h(w)\), internal content \(J\) and terminal \(t\).

Let \(V\) list \(i(w)\) and \(W\) list \(J\). Both are nonempty. Identity \((Q)\), with \(x=X\), \(y=V\), \(z=W\), gives

\CaseDisplay{XVt\leq Xt+XVWt\leq P.}

The first right word is the original head witness. The second normalizes to the word just derived, and the left normalizes to \(w\). This proves the absorption without any insertion of a new variable into the target alone.

Now suppose \(J\cap T\) is nonempty. Choose \(s\) in the intersection, a right word \(v\) ending in \(s\), and another right word \(ZsH\) in which the indicated \(s\) is internal, so \(H\) is nonempty. If \(v\) is already saturated, keep it. Otherwise write \(v=Rs\). The appropriate \((D)\) identity derives

\CaseDisplay{Rs^2\leq Rs+ZsH\leq P.}

Each absent prefix \(R\) or \(Z\) is handled by its separately displayed boundary identity. This gives some saturated word \(r\leq P\). Since \(r^2\approx r\), \((G)\) with \(x=u\), \(y=r\) gives \(u^2\leq P\). This is saturated, retains the desired head, and has content in \(C\).

For each original right word \(v\), apply \((F)\) with the current saturated word \(k\):

\CaseDisplay{k^2v^2\leq k^2+v\leq P.}

Normalization gives a saturated word with the same head and content \(c(k)\cup c(v)\). After all right contents have been incorporated, we have a saturated word of head \(h(w)\) and content \(U\) below \(P\). Words already contributing no new content can be omitted from this construction.

Let \(V\) list \(C\) and \(W\) list \(U\). Both are nonempty. Identity \((B)\), with \(x=u\), \(y=V\), \(z=W\), gives

\CaseDisplay{u^2V^2\leq u^2+u^2V^2W^2\leq P.}

Its second right word is the saturated full-content word just derived. Its left side is saturated, has head \(h(w)\) and content \(C\), and hence normalizes to \(w^2\). Finally \((E)\) gives \(w\leq w^2\leq P\). No empty word is substituted in either branch. \(\square\)

Proof of Theorem 1. Every valid identity can be distributed into two finite sums of words. Every summand on either side satisfies Lemma 4 relative to the other side, and is absorbed derivably by Lemma 5. The two inequalities prove the identity. There is no bound on variables or summands. Word reversal converts soundness, normalization and these derivations to the multiplicative opposite. \(\square\)

\par\smallskip\noindent\textit{Acknowledgments.}\quad 
The source manuscript for this case acknowledges support from the National Natural Science Foundation of China (12371024), the Science and Technology Research Program of Chongqing Municipal Education Commission (KJZD-K2024011102) and the Chongqing Natural Science Foundation Innovation and Development Joint Fund (Municipal Education Commission) (CSTB2025NSCQ-LZX0067).

\subsection{Finite bases from local content covers and final-variable propagation}\label{case:91}
Sources 2979 and 2980 in the five-element catalogue. The signature is (+, multiplication), with no constants; all bases are relative to the additively idempotent semiring axioms.

\par\smallskip\noindent\textit{1. Two additive orders on one multiplication.}\quad 
Both algebras have universe \{0,\allowbreak{}1,\allowbreak{}2,\allowbreak{}3,\allowbreak{}4\} and multiplication

{\small\par\smallskip
\CaseLine{\ttfamily 0 0 0 3 0}
\CaseLine{\ttfamily 0 0 0 3 0}
\CaseLine{\ttfamily 0 1 2 3 0}
\CaseLine{\ttfamily 3 3 3 3 3}
\CaseLine{\ttfamily 4 4 4 3 4}
\smallskip}

Rows and columns are in numerical order. In \(A=2979\) addition is join in the order \(3<1<0<4\) and \(3<2<0<4\), with 1 and 2 incomparable. In \(B=2980\) it is the maximum in the chain \(3<2<1<0<4\). These symbols name elements, not constants of the signature. In both algebras 3 is multiplicatively absorbing and is the least additive element.

Associativity follows by cases on a triple. If any factor is 3, both associations give 3. Otherwise a first factor 4 gives 4, and a first factor 0 or 1 gives 0. With first factor 2, a middle factor 0,\allowbreak{}1, or 4 gives 0 in both associations. If both first factors are 2, both associations give 2 times the last factor; left multiplication by 2 fixes its image \{0,\allowbreak{}1,\allowbreak{}2\}. Multiplication is increasing on the additive order of B, hence preserves its binary joins. For A it is increasing on each displayed chain; the one additional incomparable join \(1+2=0\) is preserved by every row and column, as the table shows. Thus both are additively idempotent semirings.

For a word of length at least two the common word-value formula is:

\textbullet{} if a variable takes value 3, the word is 3;

\textbullet{} otherwise a head value 4 gives 4;

\textbullet{} otherwise all positions 2 give 2;

\textbullet{} otherwise all nonfinal positions 2 and final position 1 give 1;

\textbullet{} otherwise the value is 0.

This follows by induction. It proves the word identities

\CaseDisplay{xyzt\approx xzyt,\qquad xy^2z\approx xyz.       \quad\textup{(\ensuremath{W})}}

Fix an ordering of variables. By interchanging adjacent interior positions and deleting repeated interior occurrences, (W) changes each nonlinear word into h (product of M) t, where h, t are its endpoints and M its interior set in increasing order. Endpoints are retained even if equal or also in M. Linear words remain linear; empty interiors are omitted substrings, not unit constants.

Write c(w), h(w), t(w) for content, head, and final variable, and i(w) for the variables before the last position. Put \(i(w)=empty\) for a linear word. Call a word saturated if t(w) belongs to i(w), so that \(i(w)=c(w)\). These data are invariant under normalization.

\par\smallskip\noindent\textit{2. Exact criteria using only contained summands.}\quad 
For a target word w with content C, let \(P_{C}\) denote the summands v of a nonempty polynomial P satisfying c(v) contained in C. Put

\CaseDisplay{D=\bigcup_{v\in P_C}c(v),\quad
 J=\bigcup_{v\in P_C}i(v),\quad R=\{t(v):v\in P_C\}.}

An empty \(P_{C}\) gives empty sets. In either algebra \(w\leq P\) is equivalent to \(w\leq \sum (P_{C})\): necessity follows by assigning every outside variable the absorbing least element 3, which kills exactly the discarded summands; sufficiency follows because the retained sum is below P. If \(P_{C}\) is empty the absorption is false, for instance with C assigned 2 and all outside variables 3. All following valid cases have a nonempty \(P_{C}\), so no empty polynomial is used as a term.

Lemma 1. In A, the absorption \(w\leq P\) holds exactly when:

1. \(D=C\);

2. some v in \(P_{C}\) has \(h(v)=h(w)\);

3. either J meets R, or i(w) is contained in J and t(w) belongs to R.

In B the exact conditions are 1, 2, and

3'. i(w) is contained in J.

These assertions include linear targets.

Proof. Work with the contained sum; assign outside variables 3 for each separating evaluation. If 1 fails, assign a missing variable 0 and all other variables of C the value 2. The target is 0 and the contained sum is 2 (or 3 if empty). If 2 fails, assign \(h(w)=4\) and all other variables in C the value 2. The target is 4 and each contained summand has value at most 0.

For B, if x belongs to i(w) but not J, assign \(x=1\) and the other variables of C the value 2. The target is 0 and the contained sum is 1 or 2, so 3' is necessary. For A assume 1 and J disjoint from R. Since \(C=D=J\) union R, assigning J the value 2 and R the value 1 makes every contained summand 1. The target is 1 precisely when i(w) is contained in J and t(w) belongs to R; otherwise it is 0 or 2, neither below 1. This proves necessity of 3.

For sufficiency, if some variable in C is 3 then \(w=3\) is least, so consider assignments with all of C in \{0,\allowbreak{}1,\allowbreak{}2,\allowbreak{}4\}. The retained sum is nonempty and has value in this same set. If it is 4 there is nothing to check. If it is 0, the only possible obstruction is \(w=4\). But a same-head witness from 2 would also have value 4, a contradiction. If the retained sum is 2, every summand must be 2 in both additive orders; the word formula then forces D, hence C, to be assigned 2, giving \(w=2\). In B a retained sum of 1 means all summands are 1 or 2. Thus D is assigned in \{1,\allowbreak{}2\} and J is assigned 2. Condition 3' forces w to be 1 or 2, both below 1. In A a retained sum of 1 requires every summand to be 1. Then J is assigned 2 and R is assigned 1, forcing them to be disjoint; condition 3 gives \(w=1\). \ensuremath{\square}

\par\smallskip\noindent\textit{3. Four right summands suffice for A.}\quad 
Normalize all words in a valid absorption. A word already occurring on the right needs no step. Choose a same-head witness v in \(P_{C}\). For a linear target y, v is a positive power of y and \(y\leq v\) is valid by Lemma 1, providing a one-summand step.

Now write the nonlinear target as \(w=h\) (product of M) t. For any subset S of C containing h, t, let w|S be h (product of M intersect S) t. This is a normal word of content S, and its prefix set is contained in i(w). No restriction of a polynomial is being used as a term.

If J meets R, choose a, b in \(P_{C}\) with some g in i(a) and \(t(b)=g\). Choose \(r_{t}\) in \(P_{C}\) containing t, and put \(S=c(v)\) union c(a) union c(b) union c(\(r_{t}\)). Then

\CaseDisplay{w|S\leq v+a+b+r_t.                        \quad\textup{(\ensuremath{2})}}

All four right words are contained in S, cover S, provide the same-head witness v, and preserve the crossing g. For each remaining q in \(P_{C}\) extend S to S union c(q), and derive

\CaseDisplay{w|(S\cup c(q))\leq w|S+a+b+q.                    \quad\textup{(\ensuremath{3})}}

The old word supplies the head witness, and a, b retain the crossing. Every right content is contained in the new target. After all summands have been processed, condition 1 gives \(S=C\) and the target w.

If J is disjoint from R, condition 3 implies t is absent from i(w). Moreover every q in \(P_{C}\) ends in t: its last variable belongs to R intersect C, which equals \{t\}, since every other variable of C belongs to i(w), hence J. Starting with w|\(c(v)\leq v\), successively adjoin c(q) by

\CaseDisplay{w|(S\cup c(q))\leq w|S+q.                  \quad\textup{(\ensuremath{4})}}

The right words end in t and their prefixes cover the new target prefix; they cover its content and have the required head witness. Thus Lemma 1 verifies every step. At most four right summands occur in either construction.

\par\smallskip\noindent\textit{4. Two right summands suffice for B.}\quad 
For a nonlinear target w with content C, keep only \(P_{C}\). Every variable in i(w) belongs to the union of its prefix sets. Let e denote the exceptional variable t(w) if t(w) is not in i(w); otherwise there is no exceptional variable. In the latter case \(i(w)=C\).

First process an arbitrary original summand q. If q is saturated, or its final variable is e, no modification is needed. Otherwise put \(f=t(q)\). The variable f lies in i(w), so choose an original summand r with f in i(r). Put

\CaseDisplay{S=c(q)\cup c(r),\qquad K=i(q)\cup i(r),\qquad g=t(r).}

Since q is unsaturated and its only missing prefix variable f is supplied by r, S minus \{g\} is contained in K. Construct a normal word q' with head h(q), content S, and final variable g. Make g occur internally if g belongs to K and differs from the head; include all other variables except the head among the interior positions. If g equals the head, its initial occurrence already makes it a prefix variable. Thus \(i(q')=K\): K is either S or S minus \{g\}. The exceptional case \(g=h(q)\) with g outside K cannot occur, because h(q) is either a prefix variable of q or is f, supplied by i(r). Lemma 1 gives the valid step

\CaseDisplay{q'\leq q+r.                       \quad\textup{(\ensuremath{5})}}

Both right words are contained in S and cover S, q is the same-head witness, and K covers the new prefix. Every word constructed remains derivably below the original polynomial.

Repeat (5) until the word is saturated or ends in e. This process terminates in at most |C| steps. Its prefix set never loses variables, and each nonterminal step adds the previously missing final variable f to that set. If a later final variable already belongs to it, the word is saturated. Otherwise this is another new variable of the finite set C. Consequently every original summand u yields a word \(d_{u}\) with the same head, containing c(u), whose content stays within C, and whose prefix contains all its content except possibly e. If there is no e, every \(d_{u}\) is saturated.

Choose the completion \(d_{v}\) of the original same-head witness, and some \(d_{r}\) containing t(w), which exists by condition 1. Put \(S=c(d_{v})\) union c(\(d_{r}\)). The word w|S is absorbed by \(d_{v}\)+\(d_{r}\): the contents cover S, \(d_{v}\) has the target head, and the combined prefixes cover i(w|S). If t(w) is repeated, both completions are saturated; otherwise their sole possible prefix omission is t(w). Now adjoin each other completed content using

\CaseDisplay{w|(S\cup c(d_u))\leq w|S+d_u.             \quad\textup{(\ensuremath{6})}}

The same prefix and head checks verify these steps. The union reaches C because the completions contain all original contents. The final word is w, and every step has at most two right summands. Linear targets are handled by \(y\leq v\) for a same-head contained power v, which is valid by Lemma 1.

\par\smallskip\noindent\textit{5. Finite bases and unrestricted completeness.}\quad 
For an absorption with at most k right normal words, reserve the endpoint variables of all k+1 words individually, at most 2(k+1). Group each other variable by its nonempty pattern of interior occurrences among the words; there are at most \(2^{k+1}\)-1 classes. Identifying each class to one variable and normalizing by (W) preserves validity. The resulting identity uses at most 2(k+1)+\(2^{k+1}\)-1 variables, and normal words of length at most two more than this bound.

This compression is reversible modulo (W). Replace a class variable by the nonempty product of all its original members. Since each member has the same interior incidence, the expansion restores each original interior set. Endpoint variables have stayed individual and restore their original positions. The laws (W) sort and remove interior repetitions, recovering every original normal word. Linear words are restored because their unique variable was reserved.

For A fix a 41-letter alphabet and take (W) together with all valid normal absorptions of length at most 43 and at most four right summands. For B fix a 13-letter alphabet and take (W) together with all valid normal absorptions of length at most 15 and at most two right summands. These are finite, explicitly defined sets; Lemma 1 decides membership. No small enumerated or minimal basis is claimed.

Theorem. These sets are finite identity bases for A and B, respectively. Reversing all words gives bases for their opposites.

Proof. Every basis member is valid. Compression shows that each valid step with the corresponding bound on right summands follows from the finite set by substitution and (W). Sections 3 and 4 then derive any valid word absorption into a polynomial. Each step uses original summands or already derived words, so transitivity and additive monotonicity place its left word below that polynomial. Distributivity writes every term as a polynomial. For a valid identity between two polynomials, each summand on either side is absorbed by the other, and the preceding argument derives all such absorptions. This is unrestricted completeness, independent of variable count and word length. Reversal carries the proof to the opposites. \ensuremath{\square}

Finite regressions are auxiliary error checks. Classification transfers require independent verification of equality of generated varieties, with both HSP directions. Neither bounded testing nor one inclusion alone supplies a finite-basis classification.

\subsection{内部首字母与平方见证：两组完整基}\label{case:92}
\par\smallskip\noindent\textit{1. 记号.}\quad 
两来源均满足 \(xyz=xzy\)、\(xy^2=xy\)。词由首字母 \(h(w)\) 与尾部变量集 \(D(w)\) 正规化，\(c(w)=\{h(w)\}\cup D(w)\)。对非空右多项式 \(F\) 记

\CaseDisplay{T=c(F),\qquad H=\{h(u):u\in F\},\qquad E=\bigcup_{u\in F}D(u).}

给定 \(X\)，令 \(F_X=\{u\in F:c(u)\subseteq X\}\)、\(C_X=\bigcup_{u\in F_X}c(u)\)。一次项允许属于 \(F_X\)。

\par\smallskip\noindent\textit{2. 内部同首字母项（来源 4122）.}\quad 
取 \(A=\{z, a, b, e, t\}\)，依此对应目录标号 \((0,\allowbreak{}1,\allowbreak{}2,\allowbreak{}3,\allowbreak{}4)\)。加法序为 \(a<z<t\)、\(b<z<t\)、\(e<t\)，\(a, b\) 不可比，\(e\) 与 \(a, b, z\) 均不可比。右乘 \(e\) 为恒等映射；右乘 \(b\) 为常值 \(b\) 映射；右乘 \(z, a, t\) 为

\CaseDisplay{P=(z, z, b, t, t).}

这些映射保持并，幂等且互相交换；\(R_{xy}=R_xR_y\) 与 \(R_{x+y}=R_x\vee R_y\) 可直接逐类核对，故定义了半环。

一个词含值为 \(b\) 的变量时取 \(b\)。否则：以 \(z, t\) 开头时分别取 \(z, t\)；以 \(a\) 开头时，尾部全为 \(e\) 则取 \(a\)，否则取 \(z\)；以 \(e\) 开头时，尾部全为 \(e\) 则取 \(e\)，否则取 \(t\)。一次词也按空尾部解释此词值规则。

\par\smallskip\noindent\textit{2.1 完整吸收判据.}\quad 
定理 A。 任意词 \(w\) 满足 \(w\le F\)，当且仅当令 \(h=h(w)\)、\(X=c(w)\) 后有

\CaseDisplay{X\subseteq T,\qquad \exists r\in F_X:\ h(r)=h,\qquad
\bigl[D(w)\subseteq E\ \text{或}\ H\cap E\ne\varnothing\bigr].\quad\textup{(\ensuremath{A})}}

内部同首字母项可以为一次项。仅有一个内部项而没有同首字母要求是不够的。

充分性。 真主理想为 \(\{a\}\)、\(\{b\}\)、\(\{a, b, z\}\)、\(\{e\}\)。右和为 \(a\) 恰在所有右首字母取 \(a\)、所有右尾变量取 \(e\) 时发生；保护交集 \(H\cap E\) 非空时不可能，否则尾部与首字母条件保证左值也为 \(a\)。右和为 \(e\) 迫使全部 \(T\) 内变量取 \(e\)，因而左值为 \(e\)。

右和为 \(b\) 时，每个右项都含一个值为 \(b\) 的变量；内部项 \(r\) 迫使左内容也含这样的变量，左值即为 \(b\)。右和在 \(\{a, b, z\}\) 中时，若左内容含 \(b\)，结论直接成立；否则内部项 \(r\) 也不含 \(b\)，它能落入该理想只可能因为其首字母取 \(a, z\)。左词首字母相同，故左值也在该理想中。检查全部真主理想即得充分性。

必要性。 若内容遗漏变量 \(x\)，赋 \(x=t\)、其余变量 \(e\)，左值为 \(t\)、右值为 \(e\)。若没有内部同首字母项，赋左首字母 \(h=e\)、其余 \(X\) 内变量 \(a\)、外部变量 \(b\)。左值为 \(e\) 或 \(t\)；内部右项都以非 \(h\) 字母开头，值为 \(a\) 或 \(z\)，外部右项则为 \(b\)。右和不超过 \(z\)，不能吸收左值。

最后，若 \(H\cap E\) 为空而某个左尾变量不在 \(E\)，赋 \(H\) 内变量 \(a\)、\(E\) 内变量 \(e\)、其余变量 \(z\)。内部同首字母项已经保证 \(h\in H\)；左值为 \(z\)，右和为 \(a\)，矛盾。证毕。

\par\smallskip\noindent\textit{2.2 五条完整基.}\quad 
定理 B。 来源 4122 的完整附加基为

\CaseDisplay{xyz\approx xzy,\quad xy^2\approx xy,\quad x\le x^2,\quad xy\le x+zy,\quad x^2y\le x^2+y.\quad\textup{(\ensuremath{B})}}

有效性由定理 A 验证。对任意已推出的 \(p\le F\) 及 \(x\in E\)，选尾部见证 \(u=vx\)，利用第四式得到

\CaseDisplay{px\le p+vx\le F.\quad\textup{(\ensuremath{3})}}

先设 \(D(w)\subseteq E\)。从内部同首字母项 \(r\) 补入全部 \(D(w)\)，得到 \(p\le F\)。由于 \(c(r)\subseteq c(w)\) 且首字母相同，所得尾部只可能是 \(D(w)\) 或 \(D(w)\cup\{h(w)\}\)；所以 \(p=w\) 或 \(p=w^2\)。第一种直接完成，第二种用 \(w\le w^2\)。此论证也包含一次左词：右边必须有纯幂 \(x^n\)，由 \(x\le x^2=x^n\)（\(n\ge2\)）或同名一次项处理。

再设 \(g\in H\cap E\)。取一个以 \(g\) 开头的右项，用 (3) 补入 \(g\)，得到首字母重复的幂等词 \(p\le F\)。第五式代入 \(x=p\)、\(y=v\)，对任意右项 \(v\) 给出

\CaseDisplay{pv\le p^2+v=p+v\le F.}

逐项乘入全部右项，得到尾部为 \(T\) 的幂等词 \(Q\le F\)。令 \(F'=F+Q\)，有可推导等式 \(F'=F\)，而 \(D(F')=T\)。原内部同首字母项仍在 \(F'\) 中，且 \(D(w)\subseteq T\)，故前段推出 \(w\le F'=F\)。

这里不需要也没有使用 \(x^2y=y^2x\)：首字母始终由内部项保留。推导允许任意变量数，无空词代入。逐词处理恒等式两向吸收即证完备性。

\par\smallskip\noindent\textit{3. 内部内容覆盖与平方见证（来源 4397）.}\quad 
取 \(B=\{z, e, b, d, t\}\)，依此对应目录标号 \((0,\allowbreak{}1,\allowbreak{}2,\allowbreak{}3,\allowbreak{}4)\)。加法序为 \(d<b<z<t\)、\(e<z<t\)，\(e\) 与 \(d, b\) 不可比。右乘 \(e\) 为恒等映射；右乘 \(b, d\) 为常值 \(b\) 映射；右乘 \(z\) 为 \(P=(z, z, b, b, t)\)；右乘 \(t\) 为 \(Q=(t, t, b, b, t)\)。

各右作用保持并、幂等且互相交换，其中 \(PQ=QP=Q\)。逐类有 \(R_{xy}=R_xR_y\)、\(R_{x+y}=R_x\vee R_y\)，因此满足半环公理。

若词内出现值为 \(b, d\) 的变量，其值为 \(b\)，唯一例外是首字母取 \(d\)、尾部全取 \(e\) 时为 \(d\)。若没有 \(b, d\)，词值是 \(\{e, z, t\}\) 链上所有变量值的最大值。因此 \((N),(R)\) 成立，并且还有 \(x^2y=y^2x\)。首字母重复的词仅由内容决定，仍记为 \(M_X\)，并有 \(w^2=M_{c(w)}\)。此外 \(x\le x^2\) 成立。

\par\smallskip\noindent\textit{3.1 完整吸收判据.}\quad 
定理 C。 任意词 \(w\) 满足 \(w\le F\)，当且仅当

\CaseDisplay{c(w)\subseteq C_{c(w)},\qquad
\bigl[(h(w)\in H\ \text{且}\ D(w)\subseteq E)\ \text{或}\ H\cap E\ne\varnothing\bigr].\quad\textup{(\ensuremath{C})}}

第一项等价于内部右项的内容并恰好覆盖左内容；一次项同样可以贡献覆盖。

证明。 真主理想是 \(\{d\}\)、\(\{d, b\}\)、\(\{e\}\)、\(\{d, b, e, z\}\)。右和为 \(d\) 当且仅当全部右首字母取 \(d\)、全部右尾变量取 \(e\)。第二项正好处理这一阈值：有保护交集时此赋值不可能，无保护时左首字母与尾部必须分别被覆盖。

若右和落入 \(\{d, b\}\)，每个右项都含取 \(d, b\) 的变量。内部覆盖非空，任取一个内部项即知左内容也有这类变量，故左值落入该理想。若右和为 \(e\)，则 \(T\) 内全部变量取 \(e\)，内部覆盖使左内容包含于 \(T\)，左值也为 \(e\)。

若右和不超过 \(z\) 而左值为 \(t\)，左内容不含取 \(d, b\) 的变量，并包含一个取 \(t\) 的变量 \(x\)。内部覆盖提供包含 \(x\) 的内部右项；它的内容同样不含 \(d, b\) 且含 \(t\)，故其值为 \(t\)，矛盾。这证明充分性。

若 \(x\in c(w)\setminus C_{c(w)}\)，赋 \(x=t\)、其余左内容内变量 \(e\)、外部变量 \(b\)。左值为 \(t\)；内部右项因不含 \(x\) 而取 \(e\)，外部右项取 \(b\)，右和不超过 \(z\)，矛盾。若第二项失败，赋 \(H\) 内变量 \(d\)、\(E\) 内变量 \(e\)、其余变量 \(z\)。此赋值因 \(H\cap E=\varnothing\) 而一致。右和为 \(d\)，而左值因首字母或尾部条件失败不能为 \(d\)，故不能被其吸收。必要性得证。

\par\smallskip\noindent\textit{3.2 九条完整基.}\quad 
定理 D。 来源 4397 的完整附加基为

\CaseDisplay{\begin{aligned}
xyz&\approx xzy, &(N)\\
xy^2&\approx xy, &(R)\\
x^2y&\approx y^2x, &(C)\\
x&\le x^2, &(I)\\
xyz&\le xy+xz, &(M)\\
xy&\le x+y^2, &(K)\\
x^2&\le x+y^2, &(S)\\
(xyz)^2&\le xy+yz, &(P_3)\\
(xy)^2&\le xy+y. &(P_2)
\end{aligned}}

每式的有效性由定理 C 直接核对。

首先考虑无保护情形 \(H\cap E=\varnothing\)。对非一次左词，定理 C 给出 \(h\in H\)、\(D(w)\subseteq E\)，因而 \(h\notin D(w)\)。每个内部右项的首字母既属于 \(H\) 又属于 \(\{h\}\cup D(w)\)，所以只能是 \(h\)。内部一次项也只能是 \(h\)。内部内容覆盖因此等价于内部非一次项的尾部并覆盖 \(D(w)\)。由于左词非一次，这类右项至少有一个。迭代 \((M)\) 合并这些同首字母词，正好得到首字母 \(h\)、尾部 \(D(w)\) 的词 \(w\le F\)。

再考虑有保护情形。取 \(g\in H\cap E\)，选择右尾部见证 \(u=vg\) 和以 \(g\) 开头的右项 \(r\)，其中 \(v\) 非空。若 \(r=gs\) 为非一次词，\((P_3)\) 在 \((x, y, z)=(v, g, s)\) 下推出 \((vgs)^2\le u+r\le F\)。若 \(r=g\)，\((P_2)\) 推出 \((vg)^2\le u+g\le F\)。两种情况都提供一个幂等平方词 \(p\le F\)。

对每个内部右项 \(q\)，由 \((S)\) 在 \((x, y)=(q, p)\) 下有

\CaseDisplay{q^2\le q+p^2=q+p\le F.}

这些平方词分别为 \(M_{c(q)}\)。对任意已经被吸收的幂等词 \(a, b\)，\((K)\) 给出 \(ab\le a+b^2=a+b\le F\)。逐个合并全部内部平方词，再用 \((C)\) 整理，得到

\CaseDisplay{M_{C_{c(w)}}=M_{c(w)}\le F.}

最后 \((I)\) 给出 \(w\le w^2=M_{c(w)}\le F\)。内部覆盖保证所合并的词族非空，不涉及空积。一次词也可直接由内部纯幂项与 \((I)\) 处理。

这两个分支覆盖任意有效吸收，证明不限制变量数。所有恒等式拆成双向逐词吸收后可导出，九式因而完整。证毕。

\par\smallskip\noindent\textit{4. 方法区别与范围.}\quad 
来源 4122 要求一个内部同首字母项，其五条基直接保持首字母。来源 4397 则允许在保护存在时换首字母，但要求全部左内容由内部项覆盖；两条平方路径式先制造幂等见证，再把各内部项平方并合并。二者不能仅凭类似的小变量样本混用判据。

乘法对偶逐词反转。独立核验重建运算表、检查基、比较一至五变量下全部右多项式的主理想 CNF，并检查完整双向 HSP 包含。任意变量数的结论仍由上述理论证明承担；对单向基成员不传递有限基性。

\subsection{内部保护与全局尾部覆盖：两组完整基}\label{case:93}
\par\smallskip\noindent\textit{1. 记号和一个句法闭合引理.}\quad 
在 \(xyz\approx xzy\)、\(xy^2\approx xy\) 下，词由首字母 \(h(w)\) 和尾部集合 \(D(w)\) 正规化。对 \(X=c(w)\) 及非空右多项式 \(F\) 定义

\CaseDisplay{F_X=\{u\in F:c(u)\subseteq X\},\quad C_X=\bigcup_{u\in F_X}c(u),\quad
H_X=\{h(u):u\in F_X\},\quad E_X=\bigcup_{u\in F_X}D(u).}

另记全部右项的尾部并为 \(E=\bigcup_{u\in F}D(u)\)。一次内部项可以贡献首字母和内容，但不贡献尾部。

闭合引理。 从六式

\CaseDisplay{xyz\approx xzy,\quad xy^2\approx xy,\quad x\le x^2,\quad
xyz\le xy+xz,\quad xyz\le xy+yz,\quad xy\le x+y^2 \quad\textup{(\ensuremath{L})}}

可以推出每一个满足 \(h(w)\in H_X\)、\(D(w)\subseteq E_X\) 的吸收 \(w\le F\)。

这是内部首字母与尾部闭合第4节的句法引理，这里重述关键推导。以一个内部同首字母项为起点，在可推导不超过 \(F\)、首字母为 \(h(w)\)、内容包含于 \(X\) 的词中选内容最大的 \(p\)。若还有未覆盖内容 \(x\)，内部尾部见证 \(u_x\) 的首字母不能在 \(c(p)\)：若等于 \(h(w)\)，用同首字母合并式；若在 \(D(p)\)，把 \(p\) 写成 \(vg\)，用 \(vgs\le vg+gs\) 拼接，都会扩大内容。因此未覆盖顶点各有仍在补集中的前驱，有限性给出补集有向环。

沿环反复拼接得到首字母重复的词 \(q\le F\)，其内容含补集顶点且仍包含于 \(X\)。前两式保证 \(q^2=q\)，于是 \(pq\le p+q^2\le F\)，又扩大内容，矛盾。故 \(c(p)=X\)。若目标首字母需要重复而 \(p\) 尚未重复，取其内部尾部见证，再按同首字母合并或上述路径拼接补入。否则 \(p\) 等于 \(w\) 或 \(w^2\)，由 \(w\le w^2\) 完成。一次词和某次合并中的空尾部均直接使用已有项，不作空词代入。这个论证不使用 \(x^2y\approx y^2x\)。

\par\smallskip\noindent\textit{2. 内部保护来源（4123）.}\quad 
取 \(A=\{z, a, b, e, t\}\)，按此顺序对应目录标号 \((0,\allowbreak{}1,\allowbreak{}2,\allowbreak{}3,\allowbreak{}4)\)。加法为 \(b<a<z<t\)、\(b<e<t\) 上的并，\(e\) 与 \(a, z\) 不可比。右乘 \(e\) 为恒等映射，右乘 \(b\) 为常值 \(b\) 映射，右乘 \(z, a, t\) 均为 \(P=(z, z, b, t, t)\)。

各右平移保持并、幂等且交换，直接逐类有 \(R_{xy}=R_yR_x\)、\(R_{x+y}=R_x\vee R_y\)，故得到半环。词内出现 \(b\) 时取 \(b\)；否则首字母取 \(z, t\) 时分别取 \(z, t\)，首字母取 \(a\) 时尾部全为 \(e\) 则取 \(a\)、否则取 \(z\)，首字母取 \(e\) 时尾部全为 \(e\) 则取 \(e\)、否则取 \(t\)。空尾部按全为 \(e\) 解释。

\par\smallskip\noindent\textit{2.1 精确判据.}\quad 
定理1。 在 \(A\) 上，任意词 \(w\) 满足

\CaseDisplay{w\le F\quad\Longleftrightarrow\quad
C_X=X,\quad h(w)\in H_X,\quad
\bigl[D(w)\subseteq E_X\ \text{或}\ H_X\cap E_X\ne\varnothing\bigr]. \quad\textup{(\ensuremath{A})}}

这里的保护交集和尾部覆盖均取自内部项，不能换成全局集合。

必要性。 若 \(x\in X\setminus C_X\)，赋 \(x=t\)、其余 \(X\) 内变量 \(e\)、外部变量 \(b\)。左值为 \(t\)，内部项取 \(e\)，外部项取 \(b\)，右和不超过 \(e\)。

若无内部同首字母项，赋 \(h(w)=e\)、其余 \(X\) 内变量 \(a\)、外部变量 \(b\)。左值为 \(e\) 或 \(t\)，内部右项取 \(a\) 或 \(z\)，外部项取 \(b\)，右和不超过 \(z\)。

设前两条件成立而最后一项失败。赋 \(H_X\) 内变量 \(a\)、\(E_X\) 内变量 \(e\)、其余 \(X\) 内变量 \(z\)、外部变量 \(b\)。保护交集为空，故赋值一致。内部右项均取 \(a\)，外部项取 \(b\)，右和为 \(a\)；左首字母取 \(a\)，但某个左尾变量不在 \(E_X\)，取 \(a\) 或 \(z\)，使左值为 \(z\)。均得到分离。

充分性。 真主理想为 \(\{b\}\)、\(\{b, a\}\)、\(\{b, a, z\}\)、\(\{b, e\}\)。若左内容含 \(b\)，左值为最小元 \(b\)，结论成立。否则内部项均不含 \(b\)。

右和为 \(b\) 时，内部同首字母项必须含 \(b\)，与上述情况矛盾。右和不超过 \(z\) 时，内部同首字母项的首字母只能取 \(a\) 或 \(z\)，左词也不超过 \(z\)。右和不超过 \(e\) 时，所有内部项只能取 \(e\)，迫使它们的全部变量取 \(e\)；由 \(C_X=X\) 得左值为 \(e\)。

右和不超过 \(a\) 时，所有内部项首字母取 \(a\)、尾部取 \(e\)。因此 \(H_X\cap E_X\) 不可能非空，判据只能由 \(D(w)\subseteq E_X\) 满足。左首字母取 \(a\)、尾部全取 \(e\)，左值为 \(a\)。真主理想检查完成。证毕。

\par\smallskip\noindent\textit{2.2 九条完整基.}\quad 
定理2。 \(A\) 的完整附加基由引理的六式 \((L)\) 及下列三式组成：

\CaseDisplay{x^2y\le x^2+y,\qquad (xyz)^2\le xy+yz,\qquad (xy)^2\le xy+y. \quad\textup{(\ensuremath{A+})}}

九式的有效性由定理1逐式检验。给定有效吸收，若 \(D(w)\subseteq E_X\)，直接用闭合引理。

否则取 \(g\in H_X\cap E_X\)，选择内部尾部见证 \(u=vg\)，以及以 \(g\) 开头的内部项 \(r\)。若 \(r=gs\) 非一次，平方路径式给出 \((vgs)^2\le vg+gs\le F\)；若 \(r=g\)，另一式给出 \((vg)^2\le vg+g\le F\)。得到首字母重复的幂等词 \(p\le F\)，且其内容包含于 \(X\)。

对任意内部项 \(q\)，第一式代入 \(x=p, y=q\) 得 \(pq\le p^2+q=p+q\le F\)。依次乘入全部内部项后，得到幂等词 \(Q\le F\)，其内容和尾部均为 \(X\)，因为 \(C_X=X\)。

令 \(F'=F+Q\)。已推导 \(F'=F\)；原内部同首字母项仍在 \(F'\) 中，新增内部项 \(Q\) 的尾部覆盖 \(X\)。对 \(w, F'\) 应用闭合引理，得 \(w\le F'=F\)。这证明任意变量数完备性。过程中保留平方词的首字母，不使用在 \(A\) 上不成立的 \(x^2y\approx y^2x\)。

\par\smallskip\noindent\textit{3. 内部内容与全局尾部来源（5136）.}\quad 
取 \(B=\{z, e, b, d, t\}\)，对应目录标号 \((0,\allowbreak{}1,\allowbreak{}2,\allowbreak{}3,\allowbreak{}4)\)。加法为 \(b<z<t\)、\(e<z<t\)、\(e<d<t\) 上的并，除传递关系外无其他比较。右乘 \(e\) 为恒等映射，右乘 \(b\) 为常值 \(b\) 映射，右乘 \(z\) 为 \(P=(z, z, b, t, t)\)，右乘 \(d, t\) 为 \(Q=(t, t, b, t, t)\)。

这些映射保持并、幂等且交换，满足 \(PQ=QP=Q\)，且逐类有 \(R_{xy}=R_yR_x\)、\(R_{x+y}=R_x\vee R_y\)，故定义半环。

词内出现 \(b\) 时取 \(b\)。否则，尾部出现 \(d, t\) 时取 \(t\)；尾部只取 \(e, z\) 时，首字母为 \(e, z\) 的词取内容在 \(e<z\) 上的最大值，首字母为 \(d\) 的词在尾部全为 \(e\) 时取 \(d\)、否则取 \(t\)，首字母为 \(t\) 的词取 \(t\)。这些规则也包含一次词。

\par\smallskip\noindent\textit{3.1 精确判据.}\quad 
定理3。 在 \(B\) 上，任意词 \(w\) 满足

\CaseDisplay{w\le F\quad\Longleftrightarrow\quad C_X=X,\qquad D(w)\subseteq E. \quad\textup{(\ensuremath{B})}}

必要性。 若 \(x\in X\setminus C_X\)，赋 \(x=t\)、其余 \(X\) 内变量 \(e\)、外部变量 \(b\)。左值为 \(t\)，内部右项取 \(e\)，外部右项取 \(b\)，右和不超过 \(z\)。

若 \(x\in D(w)\setminus E\)，赋 \(x=d\)、其余变量 \(e\)。全部右尾变量取 \(e\)，右首字母取 \(e\) 或 \(d\)，所以右和不超过 \(d\)；左尾部含 \(d\)，左值为 \(t\)。

充分性。 真主理想为 \(\{b\}\)、\(\{e\}\)、\(\{b, e, z\}\)、\(\{e, d\}\)。右和为 \(b\) 时，任一内部项含 \(b\)，左内容也含 \(b\)，左值为 \(b\)。右和为 \(e\) 时，全部右内容取 \(e\)，由内部覆盖得左值为 \(e\)。

若右和不超过 \(z\)，左内容含 \(b\) 时立即成立；否则内部项均不含 \(b\)。若左内容有变量取 \(d\) 或 \(t\)，内部覆盖提供一个包含它的内部项，词值必为 \(d\) 或 \(t\)，与右和不超过 \(z\) 矛盾。因此左内容只取 \(e, z\)，左值不超过 \(z\)。

若右和不超过 \(d\)，每个右项都不含 \(b\)，所有右尾变量取 \(e\)，首字母取 \(e\) 或 \(d\)。由 \(D(w)\subseteq E\) 和 \(X\subseteq c(F)\)，左尾部全取 \(e\)，左首字母取 \(e\) 或 \(d\)，所以左值不超过 \(d\)。证毕。

\par\smallskip\noindent\textit{3.2 八条完整基.}\quad 
定理4。 \(B\) 的完整附加基为

\CaseDisplay{\begin{aligned}
xyz&\approx xzy, &(N)\\
xy^2&\approx xy, &(R)\\
x^2y&\approx y^2x, &(C)\\
x&\le x^2, &(I)\\
xyz&\le xy+xz, &(M)\\
xy&\le x+y^2, &(K)\\
x^2&\le x+yx, &(S)\\
x^2y&\le xy+zx. &(J)
\end{aligned} \quad\textup{(\ensuremath{4})}}

有效性由定理3检验。\((N),(R),(C)\) 使首字母重复的词只由内容决定，记为 \(M_Y\)；具体换首使用 \(h^2g\approx g^2h\)，尾部用 \((N),(R)\) 排序去重。于是 \(q^2=M_{c(q)}\)。

先证明平方提升规则：若 \(q\in F\) 的首字母为 \(g\in E\)，选全局尾部见证 \(u=vg\)。若 \(q=g\)，由 \((S)\) 得 \(g^2\le g+vg\le F\)。若 \(q=gs\) 非一次，\((J)\) 给出

\CaseDisplay{q^2=g^2s\le gs+vg\le F. \quad\textup{(\ensuremath{5})}}

给定有效吸收 \(w\le F\)。若 \(h(w)\in E\)，由尾部覆盖有 \(X\subseteq E\)。所以每个内部项的首字母都在 \(E\)，由 (5) 使所有内部项平方。已吸收的幂等词 \(p, q\) 可由 \((K)\) 合并：\(pq\le p+q^2=p+q\le F\)。逐项合并全部内部平方词，并用 \((C)\) 整理，得到 \(M_{C_X}=M_X=w^2\le F\)；再用 \((I)\)。

现设 \(h=h(w)\notin E\)。内部覆盖保证至少有一个以 \(h\) 开头的内部项，因为内部项中出现 \(h\) 时它不可能在尾部。所有以 \(h\) 开头的内部项也都不重复 \(h\)。用 \((M)\) 合并这些同首字母项得到 \(p\le F\)；若其中有一次项，只在合并非一次项时用 \((M)\)，若全部为一次项则取 \(p=h\)。

其他内部项 \(q\) 的首字母属于 \(X\setminus\{h\}=D(w)\subseteq E\)，所以由 (5) 得 \(q^2\le F\)。这些 \(q\) 不含 \(h\)，因为 \(h\) 既不是其首字母又不在任何右尾部。反复用 \((K)\) 把这些幂等平方项乘入 \(p\)；所得词仍以 \(h\) 为唯一首字母，内容因内部覆盖而恰为 \(X\)，故其正规形就是 \(w\)。

一次左词也由内部纯幂项和 \((I)\) 处理。所有代入均为非空词，以上推导对任意变量数成立；逐词处理恒等式两向吸收即证八式完整。证毕。

\par\smallskip\noindent\textit{4. 证据边界.}\quad 
两个来源对保护条件的要求不同：4123的保护只能由内部项提供，5136可利用全局尾部见证，但左内容必须由内部项完整覆盖。不能把两套判据互换。

乘法对偶逐词反转证明。其他目录对象须独立核验双向HSP包含后才传递有限基性；一至五变量主理想CNF只是交叉检查。本文不从有限基推出遗传有限基，也不把单向基成员全部计为有限基。

\subsection{内部首字母与尾部闭合：两组完整基}\label{case:94}
\par\smallskip\noindent\textit{1. 两个共用乘法的半环.}\quad 
底集均为 \(\{z, e, b, d, t\}\)，按此顺序对应目录标号 \((0,\allowbreak{}1,\allowbreak{}2,\allowbreak{}3,\allowbreak{}4)\)。来源 \(A\)（4150）的加法为偏序

\CaseDisplay{b<z<t,\qquad e<z<t,\qquad e<d<t}

的并，除此之外没有新的严格比较；特别地 \(b, e\) 不可比，\(b, d\) 不可比，\(z, d\) 不可比。来源 \(B\)（4152）的加法为偏序

\CaseDisplay{b<e<z<t,\qquad b<e<d<t}

的并，\(z, d\) 不可比。

两者共用乘法：右乘 \(e\) 为恒等映射，右乘 \(b\) 为常值 \(b\) 映射，右乘 \(z, d, t\) 均为

\CaseDisplay{P=(z, z, b, t, t).}

这些右平移保持上述两种并，幂等且两两交换。记右平移为 \(R_y\)，由三个类型直接检查得 \(R_{xy}=R_yR_x\) 及 \(R_{x+y}=R_x\vee R_y\)；结合右平移保持并即得结合律和两侧分配律。

词的值如下：若任一字母取 \(b\)，词值为 \(b\)；否则，首字母取 \(z, t\) 时词值分别为 \(z, t\)，首字母取 \(e\) 时尾部全取 \(e\) 则为 \(e\)、否则为 \(z\)，首字母取 \(d\) 时尾部全取 \(e\) 则为 \(d\)、否则为 \(t\)。空尾部按``全部为 \(e\)''解释，因此规则也包含一次词。

两半环满足

\CaseDisplay{xyz\approx xzy,\qquad xy^2\approx xy,\qquad x\le x^2. \quad\textup{(\ensuremath{1})}}

前两式使词由首字母 \(h(w)\) 与尾部变量集合 \(D(w)\) 决定。写 \(c(w)=\{h(w)\}\cup D(w)\)。这里不能抹去首字母：\(x^2y\approx y^2x\) 在两来源上均不成立。

给定非空多项式 \(F\) 及非空变量集 \(X\)，定义

\CaseDisplay{F_X=\{u\in F:c(u)\subseteq X\},\qquad
E=\bigcup_{u\in F}D(u),\qquad E_X=\bigcup_{u\in F_X}D(u).}

内部项允许为一次词，其尾部贡献为空。

\par\smallskip\noindent\textit{2. 两个精确吸收判据.}\quad 
定理 1。 令 \(h=h(w)\)、\(X=c(w)\)。

在 \(A\) 上，

\CaseDisplay{w\le F\quad\Longleftrightarrow\quad
\exists r\in F_X:\ h(r)=h,\qquad D(w)\subseteq E. \quad\textup{(\ensuremath{A})}}

在 \(B\) 上，

\CaseDisplay{w\le F\quad\Longleftrightarrow\quad
\exists r\in F_X:\ h(r)=h,\qquad D(w)\subseteq E_X. \quad\textup{(\ensuremath{B})}}

必要性：内部同首字母项。 若这样的项不存在，赋 \(h=d\)，其余 \(X\) 内变量 \(e\)，外部变量 \(b\)。左词值为 \(d\) 或 \(t\)。外部右项取 \(b\)，内部右项因首字母不是 \(h\)，取 \(e\) 或 \(z\)，故右和不超过 \(z\)。因 \(d, t\not\le z\)，吸收失败。这个赋值也适用于一次左词。

必要性：尾部。 对 \(A\)，若 \(x\in D(w)\setminus E\)，赋 \(x=d\)、其余变量 \(e\)。所有右项的尾部全为 \(e\)，故右和不超过 \(d\)；左词因 \(x\) 在尾部，值为 \(z\) 或 \(t\)，不能被吸收。

对 \(B\)，若 \(x\in D(w)\setminus E_X\)，赋 \(x=d\)、其余 \(X\) 内变量 \(e\)、外部变量 \(b\)。内部右项的尾部全为 \(e\)，取 \(e\) 或 \(d\)；外部右项取 \(b\)。由于 \(b\le d\)，右和不超过 \(d\)，左词仍为 \(z\) 或 \(t\)。注意 \(A\) 中没有 \(b\le d\)，所以这一步不能错误地把全局见证收紧为内部见证。

充分性。 一个有限并半格中的不等式可按所有真主理想检验。两半环的 \(\downarrow b\) 均为 \(\{b\}\)，\(\downarrow z\) 均为 \(\{b, e, z\}\)。

若全部右项属于 \(\{b\}\)，内部同首字母项 \(r\) 含有值为 \(b\) 的变量，所以左词也取 \(b\)。若全部右项属于 \(\{b, e, z\}\)，且左内容不含 \(b\)，则 \(r\) 不含 \(b\)；它的首字母必须取 \(e\) 或 \(z\)，否则词值为 \(d\) 或 \(t\)。左词首字母相同，故左值为 \(e\) 或 \(z\)。左内容含 \(b\) 时结论直接成立。

在 \(A\) 中，另两个真主理想为 \(\{e\}\) 和 \(\{e, d\}\)。右和为 \(e\) 迫使全部右内容取 \(e\)；内部同首字母项和尾部覆盖保证左内容包含在右内容中，所以左值为 \(e\)。若所有右项属于 \(\{e, d\}\)，每一右尾变量都取 \(e\)，每一右首字母取 \(e\) 或 \(d\)。由 \(D(w)\subseteq E\) 及内部同首字母项，左词也属于 \(\{e, d\}\)。

在 \(B\) 中，另两个真主理想为 \(\{b, e\}\) 和 \(\{b, e, d\}\)。左内容含 \(b\) 时结论成立；否则所有内部项均不含 \(b\)。右和在 \(\{b, e\}\) 内时，这些内部项只能取 \(e\)，迫使它们的全部变量取 \(e\)；由内部首字母和内部尾部覆盖，左词也取 \(e\)。右和在 \(\{b, e, d\}\) 内时，内部项的尾部全取 \(e\)、首字母取 \(e\) 或 \(d\)，同样得到左词属于此理想。

全部真主理想检验完成，两个判据对任意变量数成立。证毕。

\par\smallskip\noindent\textit{3. 全局尾部的四条完整基.}\quad 
定理 2。 \(A\) 的完整附加基为

\CaseDisplay{xyz\approx xzy,\qquad xy^2\approx xy,\qquad
x\le x^2,\qquad xy\le x+zy. \quad\textup{(\ensuremath{2})}}

有效性由定理 1(A) 逐式检验。设 \(w\le F\) 有效，选择内部同首字母项 \(r\)。若已经推出 \(p\le F\)，且 \(x\in E\)，选择右项 \(u=vx\)，其中 \(v\) 非空；第四式给出

\CaseDisplay{px\le p+vx\le F.}

从 \(r\) 出发，逐个补入 \(D(w)\)。所得词 \(p\) 首字母为 \(h(w)\)，尾部为 \(D(r)\cup D(w)\)。由于 \(D(r)\subseteq c(w)\)，这正是 \(D(w)\) 或 \(D(w)\cup\{h(w)\}\)。在前两式下，\(p\) 因而等于 \(w\) 或 \(w^2\)；后一情形再用 \(w\le w^2\)。

此论证也包含一次左词：内部项只能是同名纯幂，第一、二式及 \(x\le x^2\) 足够处理。任意恒等式展开成双向逐词吸收后均可推导，故四式完整。证毕。

\par\smallskip\noindent\textit{4. 内部尾部的六条完整基.}\quad 
定理 3。 \(B\) 的完整附加基为

\CaseDisplay{\begin{aligned}
xyz&\approx xzy, &(N)\\
xy^2&\approx xy, &(R)\\
x&\le x^2, &(I)\\
xyz&\le xy+xz, &(M)\\
xyz&\le xy+yz, &(P)\\
xy&\le x+y^2. &(K)
\end{aligned} \quad\textup{(\ensuremath{3})}}

有效性由定理 1(B) 验证。以下只作非空词代入。

\par\smallskip\noindent\textit{4.1 保留首字母的三个操作.}\quad 
若 \(p=ag\le F\)、\(u=gv\le F\)，其中 \(a, v\) 非空，则由 \((P)\)

\CaseDisplay{agv\le ag+gv\le F. \quad\textup{(\ensuremath{4})}}

当 \(g\in D(p)\) 时，\((N),(R)\) 允许把 \(p\) 写成以 \(g\) 结尾的词，因此此操作保留 \(p\) 的首字母并合并所用尾部。

若 \(hv, hv'\le F\)，其中 \(v, v'\) 非空，则由 \((M)\) 得 \(hvv'\le F\)。其中一个词为一次词 \(h\) 时，另一个词已实现内容并，无须空词代入。

若 \(q\le F\) 且 \(q\) 的首字母出现在其尾部，则 \((N),(R)\) 给出 \(q^2=q\)；对任意 \(p\le F\)，由 \((K)\)

\CaseDisplay{pq\le p+q^2=p+q\le F. \quad\textup{(\ensuremath{5})}}

这里只需要 \(q\) 自身幂等，不需要不同首字母的平方词彼此相等。

\par\smallskip\noindent\textit{4.2 有向环产生幂等见证.}\quad 
把每个内部非一次右项看作从首字母指向全部尾变量的有向超边。若选定的边形成有向环 \(g_0\to g_1\to\cdots\to g_{k-1}\to g_0\)，逐条选右项并反复用 (4) 拼接。所得词 \(q\le F\) 以 \(g_0\) 开头，尾部再次含 \(g_0\)，因此幂等。它的内容包含环上所有变量，且仍包含在 \(X\) 中。长度一的环已有一个首字母重复的右项，直接使用。

此操作适用于任意环长；不为不同环长加入新恒等式。

\par\smallskip\noindent\textit{4.3 最大内容闭合.}\quad 
给定有效吸收 \(w\le F\)，令 \(h=h(w)\)、\(X=c(w)\)。内部同首字母项保证存在一个首字母为 \(h\)、内容包含于 \(X\)、且由基可推导不超过 \(F\) 的词。用 \((N),(R)\) 正规化后只有有限多个候选，选择内容大小最大的一个 \(p\)，记 \(U=c(p)\)。

若 \(U\ne X\)，则对每个 \(x\in X\setminus U\)，有 \(x\ne h\) 且 \(x\in D(w)\)。由内部尾部覆盖，选择尾部含 \(x\) 的内部右项 \(u_x\)。它的首字母 \(g\) 不可能在 \(U\)：

\textbullet{} 若 \(g=h\)，同首字母合并会扩大 \(U\)；当 \(p=h\) 时，直接使用 \(u_x\) 即可。

\textbullet{} 若 \(g\in U\setminus\{h\}\)，则 \(g\in D(p)\)，由 (4) 拼接同样扩大 \(U\)。

因此每个补集顶点 \(x\) 都有仍在补集中的前驱 \(h(u_x)\)。有限性保证补集中存在有向环。第 4.2 节给出幂等词 \(q\le F\)，其内容含补集顶点且包含于 \(X\)。由 (5)，\(pq\le F\) 保留首字母 \(h\)，却严格扩大内容，矛盾。所以 \(U=X\)。

如果 \(h\notin D(w)\)，则 \(p\) 的尾部是 \(X\setminus\{h\}\) 或 \(X\)，所以 \(p=w\) 或 \(p=w^2\)，由 \((I)\) 完成。

如果 \(h\in D(w)\) 而 \(p\) 的尾部已含 \(h\)，则 \(p=w\)。剩下 \(h\notin D(p)\) 的情形：内部尾部覆盖提供一个尾部含 \(h\) 的内部右项 \(u\)。若其首字母为 \(h\)，同首字母合并把 \(h\) 补入；若首字母为 \(g\ne h\)，则 \(g\in X\setminus\{h\}\subseteq D(p)\)，用 (4) 把 \(h\) 补入。所得词内容仍为 \(X\)，尾部为 \(X\)，正是 \(w\)。

所有情况均得到 \(w\le F\) 的推导。一次左词由同名纯幂内部项和 \((I)\) 处理，亦已包含于上述最大内容论证。对任意恒等式两向逐词使用此推导，六式完整。证毕。

\par\smallskip\noindent\textit{5. 证据范围.}\quad 
两个来源及其乘法对偶的有限基性来自上述任意变量数证明。对其他目录对象仅在独立核验双向 HSP 包含后传递。固定一至五变量的主理想 CNF 诊断用于发现反例、检查判据，不替代第 2---4 节的理论证明。有限基不自动意味着遗传有限基，本文未作后者断言。

\subsection{Affine parity with exceptional quadratic components}\label{case:95}
All bases below are relative to the additively idempotent semiring axioms in the constant-free signature \((+,\cdot)\). Write \(u\le v\) for \(u+v\approx v\).

\par\smallskip\noindent\textit{1. Two algebras and their bases.}\quad 
Let \(G=\{e, g\}\) be the group of order two, with identity \(e\).

The first algebra \(A=\{e, g, c, a, o\}\) has addition given by join in \(e<c<o\), with \(g, a<o\) and no further comparabilities. Its products inside \(G\) are the group products, \(a^2=c\), and every other product is \(o\). The catalog map, in this displayed order, is \((1,\allowbreak{}2,\allowbreak{}0,\allowbreak{}3,\allowbreak{}4)\), giving type 14659. A word of length at least three is its group product when all factors are in \(G\), and is \(o\) otherwise.

The second algebra \(B=\{e, g, f, b, o\}\) has two chains \(f<e<o\) and \(b<g<o\), with elements in different branches incomparable. Let \(\rho(f)=e\), \(\rho(b)=g\), and fix \(e, g, o\). Multiplication is \(\rho(x)\rho(y)\) in the group with zero \(G\cup\{o\}\), except that \(b^2=f\). The catalog map is \((0,\allowbreak{}2,\allowbreak{}1,\allowbreak{}3,\allowbreak{}4)\), giving type 4325. Every word of length at least three is the product of the \(\rho\) images of its letters.

These descriptions prove associativity: the exceptional square in each algebra disappears into its indicated stable value upon a third multiplication. For \(A\), group translations and the map sending only \(a\) to \(c\) preserve the listed joins. For \(B\), group translations factor through the join-preserving map \(\rho\); translation by \(b\) maps the chain \(f<e\) constantly to \(g\) and maps \(b<g\) to \(f<e\). It also preserves cross-branch joins. Thus both are semirings. The named elements are not constants of the signature.

Let \(\Sigma\) be the following ten laws:

\CaseDisplay{\begin{array}{rll}
(C)&xy\approx yx,\\
(N)&x^3yz\approx xyz,\\
(S)&x^4y^4\approx x^4+y^4,\\
(U)&x\le x^3,\\
(O)&xyz\le x^3+y^3+z^3,\\
(L)&x^3\le x+xy,\\
(T)&(xy)^3\le xy+yzt,\\
(P)&xt\le xy+yz+zt,\\
(D)&y^2\le x^2+yz,\\
(Z)&x^4\le x^3+y^2.
\end{array}}

Theorem. The identities \(\Sigma\) form a basis for \(A\). A basis for \(B\) is \(\Sigma\) together with

\CaseDisplay{(T^*)\quad (xy)^3\le xy+ztu,\qquad
(L^*)\quad x^3\le x+yzt.}

The two bases have ten and twelve displayed laws. Minimality is not asserted.

Here is a direct validity check. The stable-word formulas prove (C),(N),(S),(U). On group assignments, (O) follows because either the three values coincide, or their join is \(o\); a nongroup argument makes the right side \(o\) in \(A\). In \(B\), the cubes are the \(\rho\) images and the same argument applies. Law (L) checks the exceptional elements using \(a+c=o\) in \(A\) and the two branches in \(B\). The only quadratic value different from its cube is \(c\) in \(A\), at \(x=y=a\), or \(f\) in \(B\), at \(x=y=b\). In the first case \(yzt=o\); in the second a stable right word is \(e, g\) or \(o\), and its join with \(f\) is at least \(e\). This proves (T) and \((T^*)\). Joining any stable value in \(B\) to \(x\) dominates \(\rho(x)=x^3\), proving \((L^*)\).

For (P), a non-top value of a three-edge path either enforces equal group colors, enforces alternating group colors, or puts all four vertices at the exceptional element. An odd path has the required endpoint product in each case. In \(A\), squares lie in the chain \(e<c<o\); if \(y^2=c\), then \(yz\) is \(c\) or \(o\), and if \(y^2=o\), then \(yz=o\). This proves (D). In \(B\), a right value \(f\) forces \(x=y=z=b\), and a right value \(e\) forces \(y^2\le e\); the right value cannot be \(g\). Law (Z) is immediate when \(x\) is nongroup in \(A\); otherwise \(x^4=e\) and \(x^3+y^2\ge e\). The two chains give the same check in \(B\).

\par\smallskip\noindent\textit{2. Stable words and affine closure.}\quad 
Call a word stable if its length is at least three. Let \(c(w)\) be its content and \(\pi(w)\in\mathbb F_2^X\) its exponent-parity vector. For a nonempty word family \(R\), put

\CaseDisplay{H(R)=\operatorname{aff}_{\mathbb F_2}\{\pi(r):r\in R\}.}

An affine combination over \(\mathbb F_2\) is a sum of an odd number of vectors; repetitions are allowed.

Lemma 1. Modulo (C),(N), two stable words with the same content and parity are equal. In particular \(r^3=r\) for a stable word \(r\). If stable words \(s, t\) have equal parity and \(c(s)\subseteq c(t)\), then (C),(N),(S) derive \(s\le t\).

Proof. In a word of length at least three, move a chosen occurrence to the front and split the remaining positions into two nonempty pieces. The reverse of (N) adds two copies of that letter. Any two stable words with the same content and parity can therefore be enlarged to the same sufficiently high exponent vector. For the inequality, choose \(y\in c(s)\). Multiplying (S), at \(y, z\), by \(s\) and using the first conclusion gives

\CaseDisplay{sz^4=s+sz^4.}

Insert each missing variable in this way. The resulting word has the content and parity of \(t\), so it is equal to \(t\). \(\square\)

Lemma 2. If \(R\) is a nonempty family of stable words, \(c(w)\subseteq\bigcup_{r\in R}c(r)\) and \(\pi(w)\in H(R)\), then \(\Sigma\) derives \(w\le\sum R\), for a word \(w\) of any length.

Proof. Choose an odd product of members of \(R\) with parity \(\pi(w)\), and include every member of \(R\) twice more to cover the whole content. Repeated use of (O) and \(r^3=r\) gives this product \(v\le\sum R\). If \(w\) is stable, Lemma 1 gives \(w\le v\). Otherwise apply (U) at \(w\) to get \(w\le w^3\), and apply Lemma 1 to the stable word \(w^3\). \(\square\)

All identities preserve occurrence of a variable outside the right content: assign that variable \(o\) and all right variables \(e\) to separate it. Group assignments give another necessary condition:

\CaseDisplay{w\le F\quad\Longrightarrow\quad \pi(w)\in H(F). \quad\textup{(\ensuremath{4})}}

If the parity vector is outside the affine hull, a separating linear functional assigns all right words the same group value and the left word the other value. The group elements are incomparable.

\par\smallskip\noindent\textit{3. Preparation for the first algebra.}\quad 
Let \(F\) be a polynomial, with linear variables \(L\), quadratic edge graph \(E\), and stable-word family \(R\). For a linear variable that also occurs nonlinearly, write a right word containing it as \(xu\). Law (L) adds \(x^3\), and (U) permits replacing that linear term by \(x^3\), without changing \(F\).

If an edge \(xy\) meets the content of a stable right word, write that word as \(yuv\) with \(u, v\) nonempty. Law (T) adds \((xy)^3\), and (U) permits replacing \(xy\) by this stable word. Repeat along finite edge paths. We reach

\CaseDisplay{F\approx \sum L_0+\sum R+\sum E_0, \quad\textup{(\ensuremath{5})}}

where the variable sets \(L_0\), \(W=c(R)\) and \(V(E_0)\) are disjoint. The family \(R\) can be empty. The parity hull and total content have not changed. The set \(W\) is precisely the union of edge components meeting an original stable word or an original linear-nonlinear incidence, together with the other stable variables.

Proposition 3. For \(A\), the exact conditions for absorption by (5) are:

\textbullet{} \(x\) is absorbed if \(x\in L_0\), or \(x\in W\) and \(\pi(x)\in H(F)\).

\textbullet{} A stable \(w\) is absorbed if \(c(w)\subseteq W\) and \(\pi(w)\in H(F)\).

\textbullet{} A quadratic \(xy\) is absorbed if \(\pi(xy)\in H(F)\) and either \(\{x, y\}\subseteq W\), or \(x, y\) belong to the same component of \(E_0\). Loops are allowed.

Necessity. Condition (4) has already been proved. A variable occurring only linearly can be assigned \(c\), with the other variables \(e\); a nonlinear target containing it becomes \(o\), whereas the right side is at most \(c\). If a stable target uses a free component of \(E_0\), assign all vertices of that component \(a\) and all other nonlinear vertices \(e\), with the independent linear terms \(c\) if needed. The right value is at most \(c\) and the left stable value is \(o\). For a quadratic with one endpoint in a free component and the other outside that component, assign the first component \(a\) and the other endpoint \(e\); the left value is \(o\) and the right value is at most \(c\). The same assignment separates a missing linear target on the free component, since \(a\nleq c\). Missing variables are assigned \(o\). These establish every structural restriction. Sufficiency will follow from the derivations below.

\par\smallskip\noindent\textit{4. Combining affine and graph relations.}\quad 
For an undirected graph with loops, attach the vector \(\epsilon_x+\epsilon_y\) to an edge \(xy\). Its linear span consists of the vectors having even coordinate sum on each connected component. A sum of an odd number of edges is zero exactly when some component has an odd closed walk. The three-edge law (P) shortens an odd walk to its endpoint edge and an odd closed walk to a square.

Because the three blocks in (5) have disjoint variables, the following elementary consequences hold.

1. If \(R\ne\varnothing\), the vectors in \(H(F)\) supported on \(W\) form \(H(R)\) when \(E_0\) is bipartite, and \(\operatorname{span}\{\pi(r):r\in R\}\) when \(E_0\) has an odd circuit.

2. Zero belongs to \(H(F)\) exactly when zero belongs to \(H(R)\), with the empty-family case omitted, or \(E_0\) has an odd circuit.

3. For vertices \(x, y\) of the same \(E_0\)-component, \(\epsilon_x+\epsilon_y\in H(F)\) exactly when an odd walk joins \(x, y\), or \(0\in H(F)\).

To verify these assertions, decompose any odd vector sum according to the disjoint blocks. External linear variables must occur an even number of times. An edge sum that vanishes on its component has an even number of edges in a bipartite component; it can have odd size exactly in a nonbipartite component. Conversely, add such an odd zero relation when an affine representation needs its parity changed. Within a connected graph, endpoints of an even walk can be joined by an odd vector sum precisely when a zero affine relation is available. This also covers loops.

Now derive every absorption in Proposition 3. If \(R\ne\varnothing\) and an odd circuit in \(E_0\) has produced a square \(y^2\le F\), choose \(r\in R\) and substitute \(x=r\) in (Z):

\CaseDisplay{r^4\le r^3+y^2=r+y^2\le F.}

Adding this stable word adds the zero vector to \(H(R)\), changing that affine hull to the linear span of the original vectors. Lemma 2 then derives all required targets supported in \(W\), using assertion 1.

If \(0\in H(R)\), Lemma 2 yields a stable even-parity word \(h^2\le F\), for example a sufficiently high even power of the product of the variables in \(W\). Otherwise an odd circuit, when present, supplies a square directly. Thus assertion 2 always provides a square whenever \(0\in H(F)\). Law (D), substituting the nonempty word \(h\) if needed, gives squares at every edge vertex. A path of even length can then be made into an odd walk by inserting a loop. Law (P) consequently gives every same-component pair required in assertion 3. Without a zero affine relation, (P) handles exactly the odd walks. This proves the remaining quadratic absorptions. The linear terms in \(L_0\) are already summands.

If \(R=E_0=\varnothing\), the right polynomial is purely linear. Assigning all its variables \(c\) rules out every nonlinear target; only its listed linear terms are absorbed. This completes the proof of the ten-law basis for \(A\).

\par\smallskip\noindent\textit{5. The second algebra: one stable word suffices.}\quad 
For \(B\) use the same initial conversion of linear-nonlinear incidences by (L),(U). If a stable term now exists, \((T^*)\) upgrades every right quadratic, independently of connectivity, and \((L^*)\) upgrades every remaining linear term. Thus \(F\) becomes a nonempty stable family with its original total content and parity hull. Lemma 2 gives every absorption satisfying the necessary content and parity conditions.

If there is no stable term and no linear-nonlinear incidence, write \(F=\sum L_0+\sum E_0\) with disjoint variable sets. Assign all graph variables \(b\) and all linear variables \(f\). Then \(F=f\). Every stable word is \(e, g\) or \(o\), and therefore is not at most \(f\). A linear graph variable has value \(b\nleq f\), ruling out extra linear absorption. A quadratic with an endpoint outside the graph is separated by assigning that endpoint \(f\): its value is \(e, g\) or \(o\), whereas the right value remains \(f\).

The only possible remaining targets are therefore quadratics whose endpoints are graph vertices and whose parity vector lies in \(H(E_0)\). The incidence-vector argument in Section 4 places those endpoints in the same component, except for loops which have zero vector and may occur at any graph vertex. Odd paths, odd circuits, (D) and (P) derive exactly all of them. For a purely linear right side, assigning its variables \(f\) rules out nonlinear targets. This proves the twelve-law basis for \(B\).

Finally, an arbitrary polynomial identity is equivalent to absorption of each summand in both directions. Every such absorption has now been derived. All arguments permit arbitrarily many variables and summands.

\par\smallskip\noindent\textit{6. Catalog transfer and scope.}\quad 
The certified equal-variety classes are \(\{14659\}\) and \(\{3167,\allowbreak{}4325,\allowbreak{}8195\}\). The latter class also satisfies the ten laws for \(A\), but has its own additional laws and completeness proof. A one-way inclusion is not used to transfer finite basability.

Independent diagnostics reconstruct the operations, check the displayed laws, compare the criteria with all normal polynomials through two variables, check all right polynomials with at most three terms in three variables and further samples, and verify complete word-function closures through four variables. Those finite checks support the unrestricted proofs above.

The source manuscript for this case credits Aifa Wang and Lili Wang, School of Mathematical Sciences, Chongqing University of Technology.

The source manuscript for this case acknowledges support from the National Natural Science Foundation of China (12371024), the Science and Technology Research Program of Chongqing Municipal Education Commission (KJZD-K2024011102) and the Chongqing Natural Science Foundation Innovation and Development Joint Fund (Municipal Education Commission) (CSTB2025NSCQ-LZX0067).

\subsection{A twelve-identity basis for a group with an exceptional clone square}\label{case:96}
The signature is \((+,\cdot)\) without constants. Bases are relative to the additively idempotent semiring axioms, and \(p\le q\) means \(p+q\approx q\).

\par\smallskip\noindent\textit{The algebra and the basis.}\quad 
Let \(B=\{e, g, f, b, o\}\), with additive chain \(b<f<e<o\) and with \(g<o\) incomparable to \(b, f, e\). Addition is the join. Let \(G=\{e, g\}\) be the group of order two, and let \(\rho(f)=\rho(b)=e\) while \(\rho\) fixes \(e, g, o\). Multiplication is the product of the \(\rho\) images in \(G\cup\{o\}\), except that \(b^2=f\). The catalog map from the displayed order is \((0,\allowbreak{}3,\allowbreak{}1,\allowbreak{}2,\allowbreak{}4)\), giving type 3168.

Theorem. The following twelve laws form a basis for \(B\):

\CaseDisplay{\begin{array}{rll}
\mathrm{(C)}&xy\approx yx,\\
\mathrm{(N)}&x^3yz\approx xyz,\\
\mathrm{(S)}&x^4y^4\approx x^4+y^4,\\
\mathrm{(U)}&x\le x^3,\\
\mathrm{(O)}&xyz\le x^3+y^3+z^3,\\
\mathrm{(L)}&x^3\le x+yzt,\\
\mathrm{(T)}&(xy)^3\le xy+ztu,\\
\mathrm{(P)}&xt\le xy+yz+zt,\\
\mathrm{(B)}&x\le y+yz+zx,\\
\mathrm{(J)}&x^2\le x+y+xy,\\
\mathrm{(D)}&y^2\le x^2+yz,\\
\mathrm{(F)}&xz\le x+xt+y+yz.
\end{array}}

No minimality is asserted. The same basis applies to catalog type 14415 by the separately verified equality of the generated varieties.

The exceptional square has the same \(\rho\) image as an ordinary product, and after one further multiplication it gives that ordinary value. This proves associativity. Multiplication by any element other than \(b\) factors through the join-preserving map \(\rho\); multiplication by \(b\) maps \(b<f<e\) to \(f\le e=e\) and maps \(g\) to \(g\). These translations preserve joins, proving distributivity.

A word of length at least three is stable: its value is the group product of its \(\rho\) images unless one letter is \(o\), when its value is \(o\). Its value lies in \(\{e, g, o\}\). A quadratic has the same value apart from \(b b=f\). The stable formula proves (C),(N),(S),(U),(O). Joining a stable value to any element \(x\) dominates \(\rho(x)=x^3\), proving (L); it also dominates \(\rho(xy)=(xy)^3\), proving (T).

For the remaining laws, let \(C_0=\{b, f, e\}\). A product is \(f\) only when both factors are \(b\); it is \(e\) when both factors lie in \(C_0\) apart from that exception, or both equal \(g\); it is \(g\) when exactly one factor is \(g\) and the other belongs to \(C_0\). All other products are \(o\).

For (P), if the join of the three right edges is \(f\), all four variables are \(b\). If that join is \(e\), their group-image colors are equal along the path, so the endpoint product is \(e\) or \(f\). If the join is \(g\), the colors alternate and the endpoint product is \(g\). A top right side is automatic. For (B), a right value at most \(f\) forces \(x=y=z=b\); a right value \(e\) forces \(x\in C_0\); a right value \(g\) forces \(x=g\). This proves the inequality. In (J), a right value \(f\) forces \(x=y=b\), a right value \(e\) has \(x\in C_0\), and a right value \(g\) is impossible because the square of \(g\) is \(e\). This proves (J). The non-top square values in (D) are \(f, e\): a right value \(f\) forces \(x=y=z=b\), and a right value \(e\) dominates every non-top \(y^2\).

Finally consider (F). Its right side has no value below \(f\). If it is \(f\), then \(xt=yz=f\), forcing all four variables to be \(b\). If it is \(e\), the listed variables \(x, y\) are in \(C_0\), and \(yz\le e\) forces \(z\in C_0\), so \(xz\le e\). If it is \(g\), the listed variables \(x, y\) both equal \(g\), and \(yz=g\) forces \(z\in C_0\), so \(xz=g\). The remaining case is the top \(o\). This verifies all twelve laws.

\par\smallskip\noindent\textit{Stable absorption.}\quad 
For a word \(w\), write \(c(w)\) for content and \(\pi(w)\) for its exponent-parity vector. For a nonempty family \(R\), let \(H(R)=\operatorname{aff}_{\mathbb F_2}\{\pi(r):r\in R\}\).

The five laws (C),(N),(S),(U),(O) derive

\CaseDisplay{c(w)\subseteq\bigcup c(R),\quad\pi(w)\in H(R)
\quad\Longrightarrow\quad w\le\sum R\quad\textup{(\ensuremath{1})}}

when \(R\) is nonempty and all its words are stable. Here the target may have any positive length. For completeness, (N) inserts or removes pairs of occurrences within stable words, so content and parity determine their normal forms and \(r^3\approx r\). Multiply (S) by a stable word containing \(y\) to obtain \(sz^4\approx s+sz^4\); this permits enlargement of content without changing parity. Any affine combination is an odd sum of vectors. Choose the corresponding odd product of words of \(R\), adding pairs of factors to cover the required content. Iterated (O) bounds this product by \(\sum R\). Normalize and use the content enlargement. For a target of length one or two use \(w\le w^3\). This proves (1) by equational derivations.

Every valid absorption satisfies

\CaseDisplay{c(w)\subseteq c(P),\qquad\pi(w)\in H(P).\quad\textup{(\ensuremath{2})}}

For content, assign an absent variable \(o\) and all other variables \(e\). For parity, use a functional separating the target vector from the affine hull and assign the corresponding group values. Then all right words have one group value and the target the other; \(e, g\) are incomparable.

If \(P\) contains any stable summand, (L),(T) replace each short summand by its cube, while (U) gives the reverse inequality. The resulting polynomial is entirely stable, with unchanged content and affine hull. Thus (1) proves that (2) is also sufficient in this case.

\par\smallskip\noindent\textit{Linear and quadratic polynomials.}\quad 
It remains to consider

\CaseDisplay{P=\sum_{x\in L}x+\sum_{xy\in E}xy,}

where \(E\) is an undirected graph with loops; let \(V(E)\) be its vertex set. The exact criterion is

\CaseDisplay{w\le P\quad\Longleftrightarrow\quad
\bigl(w=x\in L\bigr)\quad\text{or}\quad
\bigl(|w|\le2,\ c(w)\subseteq V(E),\ \pi(w)\in H(P)\bigr).
\quad\textup{(\ensuremath{3})}}

First prove the additional restrictions necessary. Assign every right variable \(b\). The right value is \(b\) or \(f\), but a stable target has value \(e\), so no stable target can be absorbed. An unlisted linear target outside \(V(E)\) is absent altogether and is excluded by (2). If a quadratic target contains a variable outside \(V(E)\), assign every graph vertex \(b\) and every independent linear variable \(f\). Its value is \(e\) whenever it uses only right variables and meets an independent linear one, whereas \(P\le f\). A variable absent from the entire right side is already excluded by (2). Thus (2) and these assignments prove necessity of (3).

To derive (3), call vertices of \(L\cap V(E)\) specified. As in an ordinary undirected graph, (P) derives the endpoint edge of any odd walk, and (B) derives a linear variable at the end of any even walk starting at a specified vertex. For each connected component \(D\), put

\CaseDisplay{K_D=\operatorname{aff}\bigl(\{\epsilon_x+\epsilon_y:xy\in E_D\}
                        \cup\{\epsilon_x:x\in L\cap D\}\bigr).}

These path rules derive every short target in the component whose parity belongs to \(K_D\). Here are all cases. A bipartite component with specified vertices in at most one color class has no zero in \(K_D\); its eligible pairs have opposite colors, and its eligible linears are precisely the specified color class, when nonempty. If there is an odd cycle, (P) produces a loop. If there are specified vertices in both bipartite classes, an odd path joins two of them; (J) then produces a loop. Once there is a loop, an odd walk joins any pair, and an even walk joins a specified vertex to any vertex. These give every same-component pair and, if a specified vertex exists, every linear variable in that component. The corresponding space \(K_D\) is the even-coordinate-sum space when there are no specified vertices and the whole space otherwise.

Law (F) joins all components containing specified vertices. Indeed, if \(x, y\) are specified vertices and \(xt, yz\) are incident edges in their respective components, then (F) gives \(xz\le P\). Add this edge and repeat until these components form a single connected component. All these additions are derivable and preserve the total affine hull: the new edge vector is the sum of the three old vectors \(\epsilon_x\), \(\epsilon_y\), and \(\epsilon_y+\epsilon_z\). The new edge introduces no variables. If there are no specified vertices, this merging step is empty.

Every remaining component has only edge vectors and therefore has even total coordinate sum in the linear span of its generators. The independent linear variables have disjoint singleton supports. Consequently an eligible short target cannot meet two of the remaining components, or one of them and the merged specified component: a cross-component pair has odd coordinate sum in a component with no specified vertex. An eligible unlisted linear also cannot lie in such a component. Thus every target in (3) now lies in one connected component.

If \(0\notin H(P)\), sorting affine combinations according to the disjoint blocks shows that the permitted vector in this component lies in its own \(K_D\). The path rules just proved derive the target. If \(0\in H(P)\), some graph component contains an odd zero relation: independent singleton unit vectors cannot provide one. The path rules and (J) therefore give a square \(q^2\le P\). Law (D) gives a square at every vertex incident with an edge, via \(y^2\le q^2+yz\le P\). Adding these loops changes each component's affine hull to its linear span, which is exactly the restriction of \(H(P)\) to that component. The path rules again derive every required target.

This proves (3). In particular, no arbitrary-variable affine dependence has been replaced by a bounded-variable test: the graph argument explicitly shortens walks of any length and merges any finite number of specified components.

\par\smallskip\noindent\textit{Completeness and catalog scope.}\quad 
Distribute any identity into polynomials. If it is valid in \(B\), each summand on either side is absorbed by the other polynomial. The stable criterion or (3) derives every such absorption from the twelve laws, proving the theorem for arbitrary variable sets.

The constructed algebra gives \(\A{3168}\). Proposition~\ref{prop:lower-fibres} in both directions identifies its variety with that of \(\A{14415}\).

This result decides these catalog types as finitely based; it does not assert that all subvarieties are finitely based or complete the general five-element classification.

The source manuscript for this case credits Aifa Wang and Lili Wang, School of Mathematical Sciences, Chongqing University of Technology.

The source manuscript for this case acknowledges support from the National Natural Science Foundation of China (12371024), the Science and Technology Research Program of Chongqing Municipal Education Commission (KJZD-K2024011102) and the Chongqing Natural Science Foundation Innovation and Development Joint Fund (Municipal Education Commission) (CSTB2025NSCQ-LZX0067).

\subsection{Complete finite bases for 6486 and 6487}\label{case:97}
Signature (+, multiplication), without constants. Relative to the additively idempotent semiring axioms, take seven common laws:

{\small\par\smallskip
\CaseLine{\(xy=yx; xyz=xyzz; xxy\leq xxx+yyy; xxx\leq x+yyy;\)}
\CaseLine{\(xy\leq xxy; xy\leq xx+yy; xx\leq x+xy.\)}
\smallskip}

For 6486 add \(x\leq xxx\) and \(xxy\leq x+xy\), giving nine laws. For 6487 add \(x\leq xx\), giving eight laws. Here \(u\leq P\) abbreviates \(u+P=P\).

\par\smallskip\noindent\textit{Word values and exact absorption.}\quad 
Both multiplication tables are 01114/\allowbreak{}11111/\allowbreak{}11111/\allowbreak{}11121/\allowbreak{}41114. Addition for 6486 is 00004/\allowbreak{}01114/\allowbreak{}01214/\allowbreak{}01134/\allowbreak{}44444: 2 and 3 are incomparable below \(1<0<4\). Addition for 6487 is 00004/\allowbreak{}01114/\allowbreak{}01224/\allowbreak{}01234/\allowbreak{}44444, the chain \(3<2<1<0<4\). A literal has its assigned value. A nonlinear word whose variables all take values in \{0,\allowbreak{}4\} has value 4 if some variable is 4, otherwise

0. A quadratic with both occurrences assigned 3 has value 2. All other nonlinear words have value 1. This follows by direct induction from multiplication and has no restriction on the word's length.

For the original nonempty polynomial P let L be its literal variables, E its unordered quadratic edges, and V the union of their endpoints. For target w, put \(U=c(w)\) and let C mean

U = union \{c(u): u is an original summand of P and c(u) subset U\}.

Thus every target variable needs a contained original witness. Literals may supply this coverage. Let G mean that P has an original word of length at least three; for 6486 also allow L intersect V nonempty. Then \(w\leq P\) holds exactly under these conditions:

\textbullet{} Long target: C and G.

\textbullet{} Quadratic target: C and (G or U subset V).

\textbullet{} Literal target in 6486: the original literal, or C and G.

\textbullet{} Literal target in 6487: C.

To separate failure of C choose x in U outside the displayed union, assign \(x=4\), other target variables 0, and all outsiders 1. Every contained right word misses x and has value 0; every noncontained right word has value 1. The target has value 4. This also covers a missing literal. If C holds but a remaining condition fails, G is false. Set V to 3 and all other variables to 2. Every quadratic right term is 2. For 6486, L is disjoint from V, so every literal is 2 as well. For 6487 all right terms are at most 2. A long target, or a quadratic with an endpoint outside V, has value 1. A missing 6486 literal with C has a pure quadratic witness, hence its value is 3, incomparable with 2. These assignments prove necessity for all sizes.

For sufficiency, G supplies a right value at least 1. For the overlap case in 6486, an original incident literal x and edge xy do this: \(xy=2\) forces \(x=3\) and \(3+2=1\); otherwise xy is 0,\allowbreak{}1 or 4, all at least 1. If the target's nonlinear value is 0 or 4, C supplies a contained word with value at least the target: for 4 choose the witness containing a target variable assigned 4, and for 0 choose any contained witness. Other nonlinear target values are at most 1. In 6486 the same argument handles literals, since 2 and 3 are below 1. In 6487 every element is below each of its positive pure powers, so C suffices for a literal. If G is false and the right value does not dominate 1, every original quadratic must have value 2. All V variables then equal 3, so a target quadratic on V has value 2. This proves sufficiency.

\par\smallskip\noindent\textit{Unrestricted finite derivations.}\quad 
The two word laws normalize a long word by content. Commute a chosen variable to the end, duplicate it by \(xyz=xyzz\) using two nonempty preceding blocks, and delete surplus occurrences by the reverse law while keeping length at least three. Canonical forms are sorted content for at least three variables, xxy for two variables, and xxx for one. Consequently each long g satisfies \(g^{3}=g\), and long words with the same content are provably equal. Quadratics are normalized only by commutativity; literals retain their length.

Assume G. Keep an original long g; in the 6486 overlap case first derive \(g=xxy\) from the original x and xy using \(xxy\leq x+xy\). For every original contained u, apply \(xxx\leq x+yyy\) with \(x=u\), \(y=g\) to derive \(u^{3}\) below P. Its right terms are u and \(g^{3}=g\). Combine these cubes one at a time by \(xxy\leq xxx+yyy\). With long substitutions p, q, the right terms are p and q, and the left is a long word with their content union. C says that the final union is exactly U. Normalize to the target long form. For a quadratic use \(xy\leq xxy\); for a 6486 literal use \(x\leq xxx\). No external content is introduced into the resulting word, although the retained guard g may have external variables.

A missing 6487 literal has an original contained pure power. If that witness is long, normalize to xxx and derive xx with \(xy\leq xxy\); then \(x\leq xx\) recovers the literal. No global guard is needed.

It remains to recover a new quadratic without G. If the target edge is original there is nothing to prove. Otherwise, every contained word witnessing either endpoint x is an original literal x or the loop xx: there are no original long words, and an original mixed edge on the target endpoints would be the target itself. Each endpoint also belongs to V. If its loop is missing, keep the original literal x and an original incident edge xy, and derive xx by \(xx\leq x+xy\). After both loops are available, \(xy\leq xx+yy\) recovers the target. A square target uses the same argument with one endpoint. For 6486 this branch has no literal-incidence case because that would activate G; the uniform derivation remains valid.

These steps derive every semantically valid absorption for arbitrary variable sets, word lengths, and polynomial sizes. All substitutions are nonempty. Distributive expansion and the absorptions in both directions derive every identity. Each displayed law is valid by the exact criterion (and checked on every scalar assignment), so these are finite complete bases. Computations supplement this proof.

Unlike 11771/\allowbreak{}11772, \(xyzt\leq xyz\) fails: \(x=y=z=0\), \(t=4\). A single contained word cannot replace the full coverage: \(xxx\leq yy+zzz\) fails at \(x=4\), \(y=z=0\). In 6486 \(x\leq xx\) fails at \(x=3\); in 6487 \(xxy\leq x+xy\) fails at \(x=y=3\). These distinctions must be retained during transfers.

\subsection{Complete finite bases for 3173 and 3174}\label{case:98}
Signature (+, multiplication), without constants. Relative to the additively idempotent semiring axioms, use five common laws:

{\small\par\smallskip
\CaseLine{\(xy=yx; xyz=xyzz; xxy\leq xxx+y; xy\leq xxy; xy\leq xx+yy.\)}
\smallskip}

For 3173 add \(x\leq xxx\) and \(xxy\leq x+xy\), giving seven laws. For 3174 add \(x\leq xx\), giving six laws. Here \(u\leq P\) means \(u+P=P\).

\par\smallskip\noindent\textit{Exact values and criterion.}\quad 
Both multiplication tables are 00034/\allowbreak{}00034/\allowbreak{}00134/\allowbreak{}33333/\allowbreak{}44434. Addition for 3173 is 00004/\allowbreak{}01014/\allowbreak{}00224/\allowbreak{}01234/\allowbreak{}44444: 1 and 2 are incomparable between 3 and 0, with \(0<4\). Addition for 3174 is 00004/\allowbreak{}01114/\allowbreak{}01224/\allowbreak{}01234/\allowbreak{}44444, the chain \(3<2<1<0<4\). A nonlinear word is 3 if any of its variables is 3. Otherwise it is 4 if any variable is 4. If neither occurs, a quadratic whose two occurrences are 2 has value 1; all other nonlinear words have value 0. Literals retain their assigned value. Direct induction proves this description at every length.

Given original polynomial P and target w, set \(U=c(w)\), and retain only original words \(P_{U}\) whose contents are subsets of U. Let L be the literal variables of \(P_{U}\), V the endpoints of its quadratic words, and C the condition that the content union of \(P_{U}\) is U. Let G mean that \(P_{U}\) has a long word (length at least three); for 3173 also allow L intersect V nonempty. Then:

\textbullet{} A long target is absorbed exactly when C and G hold.

\textbullet{} A quadratic target is absorbed exactly when C holds and either G holds or its endpoints belong to V.

\textbullet{} A literal target in 3173 is absorbed when it is original, or when C and G hold; in 3174, C alone is necessary and sufficient.

Every witness in these conditions is original and contained in U. In particular, a long word with outside variables cannot supply G.

If C fails, choose an uncovered x in U, assign \(x=4\), other U variables 0, and outsiders 3. Noncontained right words are 3; contained ones miss x and are 0. The target is 4. If C holds but another condition fails, G is false. Assign V variables 2, other U variables 1, and outsiders 3. The right quadratics are

1. In 3173, L is disjoint from V, so contained literals also are 1; in 3174 all contained literals are at most 1. A long target, or a quadratic with an endpoint outside V, is 0. For a missing 3173 literal, C gives a pure quadratic witness, so its value 2 is incomparable with the right value 1. These are separating assignments for every failed case, at arbitrary size.

For sufficiency it is enough to use \(P_{U}\). If a target variable is 3 the nonlinear target is bottom; if none is 3 but some is 4, C supplies a contained word involving that variable with value 4. In the remaining nonlinear case all target values are in \{0,\allowbreak{}1,\allowbreak{}2\}. A contained long word is 0. In 3173 an incident contained literal x and quadratic xy also have join 0: \(xy=1\) forces \(x=2\), and \(1+2=0\); otherwise \(xy=0\). Thus G supplies 0. If G is false and the right value does not reach 0, all contained quadratics must be 1, forcing every V variable 2; the target quadratic on V is 1. For literals, 3173 needs the stated pure long witness unless the literal is original. In 3174 an element is below every pure power. This proves the criterion for all polynomials.

\par\smallskip\noindent\textit{Complete finite derivations.}\quad 
Commutativity and \(xyz=xyzz\) normalize every long word by content: permute a variable to the end, duplicate it using two nonempty preceding blocks, and delete surplus occurrences with the reverse law, always keeping length at least three. Canonical forms are sorted content for three or more variables, xxy for two, and xxx for one. Long g therefore satisfies \(g^{3}=g\). Quadratics only commute.

When G holds, start with an original contained long g. In the 3173 overlap case derive \(g=xxy\) from the original contained x and xy via \(xxy\leq x+xy\). For each original u in \(P_{U}\), apply \(xxy\leq xxx+y\) with \(x=g\), \(y=u\). Its right words are \(g^{3}=g\) and u; its left is a long word whose content is c(g) union c(u). Repeat with the resulting long word. C ensures the final content is exactly U, so normalization gives the required long form. Recover a quadratic by \(xy\leq xxy\), or a 3173 literal by \(x\leq xxx\).

A missing 3174 literal has an original contained pure nonlinear power. If long, normalize to xxx and use \(xy\leq xxy\) at \(y=x\) to derive xx, then apply \(x\leq xx\). If already quadratic only the last step is needed.

Without G, consider a quadratic target that is not original. Its endpoints lie in V, and every quadratic contributing to V is contained in the target content. Any mixed edge on these endpoints would be the target itself. Consequently both original endpoint loops are present, and \(xy\leq xx+yy\) derives the target. For a square, the original loop is already the target.

All substitutions are nonempty. The induction has no bound on variables, word length, or number of right summands. Expanding terms by distributivity and applying the absorption derivations in both directions proves completeness of the displayed bases. Each basis law satisfies the exact criterion and is also checked on all scalar assignments. Finite computation is supplemental.

The need for a contained guard is real: \(xxx\leq xx+yyy\) fails with \(x=2\), \(y=3\) in both algebras. Coverage is also necessary: \(xxy\leq xxx\) fails at \(x=0\), \(y=4\). In 3173 \(x\leq xx\) fails at \(x=2\); in 3174 \(xxy\leq x+xy\) fails at \(x=y=2\). A single contained word, external global activation, and literal/loop recovery must not be conflated.

\subsection{A finite identity basis for 3177 and its multiplicative opposite}\label{case:99}
The signature has addition and multiplication and no constants. Write B=\{0,\allowbreak{}1,\allowbreak{}2\}. In 3177 the addition is join in the order with 1,\allowbreak{}2,\allowbreak{}3 below 0 and 0 below 4, with 1,\allowbreak{}2,\allowbreak{}3 pairwise incomparable. The multiplication rows are 00034, 00034, 00134, 33333, 44444. Both 3 and 4 are left zeros. A word evaluates to the first value in \{3,\allowbreak{}4\} encountered when one exists. Otherwise a literal retains its value, a quadratic evaluates to 1 on the pair (2,\allowbreak{}2) and to 0 elsewhere, and a word of length at least three evaluates to 0. This follows successively from the two left zeros, the rows indexed by B and the fact that B cubed is \{0\}.

For a target w and a nonempty polynomial P of original right words, write L for the literal variables in P and Q for the union of variables in its quadratic words. Say that G holds if P has a word of length at least three or L and Q meet. Then w \ensuremath{\leq} P holds precisely under these conditions.

1. Some original right word has content contained in the content of w.

2. For each target variable x, some original right word contains x and every letter before its first x occurs before the first x in w.

3. Either G holds, or w is a literal in L, or w is quadratic and both of its letters belong to Q.

Here is the proof. If condition 1 fails, assign 0 throughout the target and 3 outside. The target is 0 and the right side is 3. If condition 2 fails at x, give x the value 4, its target prefix value 0 and every other letter value 3. Only the target can reach 4 among the target and the right side. If G fails, assign 1 to L, 2 to Q and 0 elsewhere. Every original word has value 1, so w must have value 1. This forces exactly the last two possibilities in condition 3.

For sufficiency use the four proper principal ideals in this order. If w is 4, its first special letter is witnessed by condition 2 and that right word is also 4. This settles the ideal below 0. If all right words are 3, a neutral target is excluded by its contained witness, and a target of value 4 has its prefix witness. This settles the ideal \{3\}. If all right words are 1 or if all are 2, every variable of P has its value in B. Condition 2 ensures the same for w. The case G makes a common value of 1 impossible by a long word or the conflicting requirements on L and Q. It also supplies a nonlinear right word, which cannot have value 2. In the remaining literal case its original literal occurs. In the remaining quadratic case the presence of a quadratic excludes common value 2, and common value 1 forces both target letters to be 2. This settles the ideals \{1\} and \{2\}. Every failure of absorption is witnessed in the principal ideal of the right sum, proving the claim.

Take the semiring axioms and the following seven word identities:

{\small\par\smallskip
\CaseLine{\(xxyz = xyz\qquad xyyz = xyz\qquad xyzz = xyz\qquad xyxz = xyz\)}
\CaseLine{\(xyzy = xyz\qquad xyx = xxy\qquad xyy = xxy\)}
\smallskip}

Each compares words of length at least three with the same list of first occurrences. Such words have equal values by the formula above. These laws derive equality between any two such words. While a word of length at least four has a repetition, an adjacent repetition can be brought under one of the first three rules with nonempty substitutions by using the available context. For separated repeated letters with a nonempty right or left context use the fourth or fifth rule. If neither context exists, the word is xux with |u| at least two. Use \(xyx=xxy\) and then the first rule after splitting u. At length three the last two rules bring every two-letter-content case to xxy, where x is its first letter. A singleton content gives xxx. The result is determined entirely by the first-occurrence list and the distinction between short and long words.

Adjoin all table-valid absorptions over twelve variables in which every word has length at most two and there are at most five upper words. Also adjoin the following six schemes, each for the three choices \(G=ijk\), \(G=i+ij\) and \(G=i+ji\), using fresh guard letters:

{\small\par\smallskip
\CaseLine{\(xxxu \leq  u + xs + G\qquad xxxu \leq  u + x + G\)}
\CaseLine{\(pxxxu \leq  pu + pxs + G\qquad pxxxu \leq  pu + px + G\)}
\CaseLine{\(a \leq  aaa + G\qquad ab \leq  ababab + G\)}
\smallskip}

These are eighteen ordinary identities. Each has a contained right word and every first occurrence has its required prefix witness. Its guard provides condition 3. All are therefore valid. This specifies a finite basis, including the bounded short family.

To prove completeness, first suppose G fails in the original P. For a literal target use the original literal. For a quadratic target keep one contained original, at most two prefix witnesses and at most two quadratic originals accounting for the target letters in Q. All selected words are short. The criterion still holds, and at most twelve variables occur including the target. This is an instance of the short family. Adding the omitted summands proves w \ensuremath{\leq} P.

When G holds, fix one witnessed guard throughout the derivation. A long original can be split into three nonempty factors. An original literal i and a quadratic containing it give exactly i+ij or i+ji, allowing \(j=i\). In every instance the fixed guard is below P. Let the first-occurrence list of w be x1,..., xn and let \(p_{j}\) be the product of the first j letters in that order, with each letter cubed. Inductively derive \(p_{j}\) u \ensuremath{\leq} P for every original u. For \(j=1\) an original prefix witness starts with x1. Use the first or second scheme according as its remaining suffix is nonempty or empty. Given the claim for \(p=p_{j}\), take a prefix witness \(v=r\) x s for the next variable x. The letters of r already occur in p. The long words p v and p x s have the same first-occurrence list, including when r or s is absent. The word identities therefore give p x s \ensuremath{\leq} P. The third or fourth scheme combines it with p u \ensuremath{\leq} P to prove p xxx u \ensuremath{\leq} P. An absent suffix selects the fourth scheme. No empty word is a substitution.

For the original contained witness u0, \(p_{n}\) u0 has the target's first-occurrence list. If w is long, the word laws recover w. If w is literal or quadratic, the same laws identify \(p_{n}\) u0 with w cubed, and the fifth or sixth scheme recovers w. Distributive expansion and these derivations in both directions of a polynomial equality give completeness for arbitrarily many variables, summands and letters. Reversal gives a finite basis of the opposite, which in catalogue labels is 3180.

\subsection{Finite bases for 3178 and 3179}\label{case:100}
Each algebra has multiplication rows 00034, 00034, 00134, 33333, 44444. The first assigned value in \{3,\allowbreak{}4\} gives the value of a word. In the absence of these values a literal retains its value, a quadratic gives 1 on (2,\allowbreak{}2) and 0 otherwise, and a longer word gives 0. In both additive orders 3 is least and 4 greatest. The neutral set B=\{0,\allowbreak{}1,\allowbreak{}2\} has top

0. In 3178 the neutral values 1,\allowbreak{}2 are incomparable. In 3179 they satisfy \(2<1<0\). All polynomials and substitution words are nonempty.

For w \ensuremath{\leq} P let U consist of the original right words with content in c(w). Let L be their literal variables, Q the union of their quadratic contents and J their long words (length at least three). Require a prefix witness in P for every target variable x. This is a word with an x whose preceding variables occur before the first target x. Its suffix may contain external variables. This condition is necessary: assign 4 to x, 0 to its preceding target variables and 3 elsewhere. Without a witness no right word reaches the target value 4.

Besides this common condition require the following. For 3178 set G to mean that J is nonempty or L intersects Q. For 3179 set G to mean that J is nonempty. In either source G suffices. When G fails the only valid cases are these:

\textbullet{} For 3178 a literal needs its original literal. For 3179 a literal needs an original contained word, necessarily its literal or square.

\textbullet{} A quadratic needs each of its letters in Q.

\textbullet{} A long target requires G.

To prove the remaining criterion, a target whose first special is 3 is already least. A first special 4 is handled by its prefix witness. For a neutral target, send all external variables to 3. All words outside U then have least value. Increasing external values can only increase the upper polynomial. It suffices to study neutral assignments on c(w). A long word in U has value 0. In 3178, if an original literal also occurs in a contained quadratic, these words cannot all be 1 or all be 2, so their join is 0. If both guards fail, assigning 1 on L, 2 on Q and 0 elsewhere makes every word of U equal 1. With U empty the upper value is 3. This forces either an original literal or a quadratic covered by Q (a literal assigned 2 is incomparable with 1). These conditions suffice: for a covered quadratic a right quadratic already prevents common value 2, while common value 1 forces both target letters to be 2. In 3179 without a long contained word assign 2 on Q and 1 elsewhere. Every contained upper is at most 1, which excludes a long target and any quadratic whose content escapes Q. A covered quadratic is at most the right sum, since that sum can be below 0 only when every contained quadratic has both letters 2. For a literal x the presence of x or xx suffices since x \ensuremath{\leq} xx holds in this chain. Without any contained word, give x value 0 and all others 3. This proves necessity and sufficiency in all degrees.

Use the seven word identities and their long normal forms proved in Section~\ref{case:99}. They use only multiplication, common to all three tables. Add all valid absorptions on at most eight variables, word length at most six and at most five right summands. Finally add four schemes for each available guard below. The symbol W denotes the specified concatenation and G the sum of the displayed guard words. Expand these abbreviations before forming the identities. All letters in G are fresh.

\noindent Source\enspace\textbar\enspace G\enspace\textbar\enspace W\par
\noindent both\enspace\textbar\enspace ijk\enspace\textbar\enspace ijk\par
\noindent 3178\enspace\textbar\enspace i+ij\enspace\textbar\enspace iij\par
\noindent 3178\enspace\textbar\enspace i+ji\enspace\textbar\enspace iji\par
{\small\par\smallskip
\CaseLine{\(xxx u W\qquad \leq  u + xs + G\)}
\CaseLine{\(xxx u W\qquad \leq  u + x + G\)}
\CaseLine{\(p xxx u W \leq  p u W + p x s W + G\)}
\CaseLine{\(p xxx u W \leq  p u W + p x W + G\)}
\smallskip}

Their validity follows from the criterion. Every guard constituent is contained in the target, so its neutral test passes. In the first two schemes x and u have their indicated witnesses. Every new letter of W is witnessed by its own original guard constituent, since W concatenates these constituents. In the last two schemes p, u and the letters of W have prefix witnesses in p u W and x is witnessed by the other constructed word. This checks the displayed equations in independent formal letters. Their substitutions preserve validity.

Suppose the contained G condition holds. Choose one actual guard inside U and fix the corresponding W for the whole derivation. A long original is split into three nonempty factors. A conflict yields the literal i and an original quadratic ij or ji, allowing \(j=i\). If x1,..., xn is the first-occurrence list of the target, put \(p_{j}=x1^{3}\cdots xj^{3}\). Derive \(p_{j}\) u W \ensuremath{\leq} P simultaneously for each original u. An original witness of x1 starts with x1; use one of the first two schemes according to whether its suffix is nonempty. At stage j+1 choose an original \(v=r\) x s witnessing the new x. The prefix r has content inside \(p_{j}\). Thus \(p_{j}\) v W and \(p_{j}\) x s W are long words with the same first-occurrence list. Use the third or fourth scheme and the already derived \(p_{j}\) u W. At the final stage choose any guard constituent as u. The entire suffix u W is contained in c(w), so \(p_{n}\) u W has the first-occurrence list of w. This gives w if w is long and w cubed if w is short.

To recover a literal or quadratic with a contained guard, retain that guard and w cubed. Every guard letter is then a target letter. By the seven word identities its long constituent, if any, normalizes to length three on at most two letters. Thus w \ensuremath{\leq} w cubed + G is an instance of the short finite family: it satisfies the prefix and contained guard conditions, and w cubed has length at most six.

When G fails, a valid literal is already present or is recovered from its square in 3179. For a valid quadratic keep at most two contained quadratic witnesses covering its content and one original prefix witness for each distinct target letter. For a square the latter starts with its letter. For \(w=xy\) with x distinct from y, the x-witness starts with x and the y-witness has the form \(x^{m}\) y t where m may be zero. Represent each nonempty suffix of a prefix witness by one fresh variable. Keep the first letters and shorten the displayed run \(x^{m}\) to length at most two. After substituting the actual suffix, the seven word laws recover the original word whenever shortening occurred, since both words are long and have the same list of first occurrences. Omit the suffix placeholder entirely for an empty suffix. The schematic inequality still has both prefix witnesses and the original contained quadratic coverage. Its words have length at most four, and it uses at most four variables and four upper words. The finite family derives it. Restoring each suffix, applying the word laws and adding the omitted right summands proves the desired absorption.

The construction works for arbitrary finite alphabets, word lengths and numbers of original summands. Expansion by distributivity and summand absorption in both directions then give full equational completeness relative to the semiring axioms.

\subsection{Twelve identities for the first-reset semiring}\label{case:101}
Source 3216. The signature is (+, multiplication) without constants, relative to the additively idempotent semiring axioms.

\par\smallskip\noindent\textit{1. Algebra and word normalization.}\quad 
The tables in order 0,\allowbreak{}1,\allowbreak{}2,\allowbreak{}3,\allowbreak{}4 are

{\small\par\smallskip
\noindent addition       multiplication\par
\CaseLine{\ttfamily 0 0 0 0 4\qquad 0 0 0 3 4}
\CaseLine{\ttfamily 0 1 0 1 4\qquad 0 0 0 3 4}
\CaseLine{\ttfamily 0 0 2 2 4\qquad 0 0 2 3 4}
\CaseLine{\ttfamily 0 1 2 3 4\qquad 3 3 3 3 3}
\CaseLine{\ttfamily 4 4 4 4 4\qquad 4 4 4 4 4}
\smallskip}

The additive order is \(3<1\),\(2<0<4\). Multiplication on \{0,\allowbreak{}1,\allowbreak{}2\} has value 2 exactly when both inputs are 2, and 0 otherwise, so is associative. Both 3 and 4 are left zeros; a right factor in \{3,\allowbreak{}4\} replaces any preceding \{0,\allowbreak{}1,\allowbreak{}2\}-product. Thus the first \{3,\allowbreak{}4\} factor determines a product whenever one occurs. This also proves associativity for mixed triples. The rows and columns preserve the displayed joins, proving both distributive laws (also directly checked on all triples).

A nonlinear word therefore has the value of its first factor in \{3,\allowbreak{}4\}, if any; when there is no such factor, it is 2 if all factors are 2, and 0 otherwise. Linear words retain their assigned value.

Theorem. A relative basis is the following twelve laws:

\CaseDisplay{\begin{aligned}
x^2y&\approx xy,&xyx&\approx xy,&xy^2&\approx xy,&x^3&\approx x^2,\quad\textup{(\ensuremath{W})}\\
x&\le x^2,\quad\textup{(\ensuremath{L})}\\
xyz&\le xy+z,\quad\textup{(\ensuremath{U})}\\
xyz&\le x+yz,\quad\textup{(\ensuremath{I})}\\
xyz&\le xz+y,\quad\textup{(\ensuremath{B})}\\
xyzt&\le xzyt+xyu,\quad\textup{(\ensuremath{R_1})}\\
xyz&\le xzy+xyu,\quad\textup{(\ensuremath{R_2})}\\
xyz&\le yxz+xu,\quad\textup{(\ensuremath{R_3})}\\
xy&\le yx+xu.\quad\textup{(\ensuremath{R_4})}
\end{aligned}}

The word laws normalize a word with at least two distinct variables to the list of variables in their first-occurrence order. A linear word stays linear; a nonlinear word of singleton content becomes xx. Indeed a repeated letter separated from its earlier occurrence is deleted using \(xYx=xY\) with a nonempty intervening block Y. Adjacent copies are deleted with \(x^{2}y=xy\) if a following nonempty suffix exists, or with \(xy^{2}=xy\) if a preceding nonempty prefix exists. Singleton powers are reduced by \(x^{3}=x^{2}\). All replacements preserve first-occurrence order, and no law identifies xx with x.

\par\smallskip\noindent\textit{2. Exact criterion.}\quad 
For target content C let \(P_{C}\) consist of the original right words whose contents are contained in C. For each variable in a word, its prefix means the content strictly before its first occurrence.

Lemma. A linear x is absorbed by P exactly when some original word has content \{x\}. For a nonlinear target w, absorption holds exactly when:

\textbullet{} (N) \(P_{C}\) contains a nonlinear word;

\textbullet{} (C) the union of contents of \(P_{C}\) equals C;

\textbullet{} (O) for each target variable x, some original right word contains x with its prefix contained in the target prefix before x.

Necessity. If a linear x has no pure-variable witness, assign \(x=0\) and all outsiders 3. Every right word is 3 and the target is 0. If N fails, assign C the value 1 and outsiders 3. The target is 0, while every nonlinear right word is 3 and linears are 1 or 3, so \(P\leq 1\). If C fails at x, assign \(x=0\), other target variables 2, and outsiders

3. A right word containing an outsider has value 3; a contained right word avoids x and is 2. Thus \(P\leq 2\) while the target is 0.

If O fails at x, assign the variables earlier than x in the target the value 2, assign \(x=4\), and all remaining variables 3. The target has value 4. Every right occurrence of x has a 3 before its first occurrence, so no right word reaches 4. The right sum is at most 0. This separates the failure.

Sufficiency. Target value 3 is least. A target value 4 is caused by its first reset variable x valued 4. All earlier target variables avoid \{3,\allowbreak{}4\}; the O-witness has only such variables before its first x and therefore also has value 4. If the target is 2, every target variable is 2 and the N-witness has value 2.

If the nonlinear target is 0, every target variable lies in \{0,\allowbreak{}1,\allowbreak{}2\} and at least one is 0 or 1. A contained nonlinear witness is 0 or 2. The first case finishes. In the second case, C supplies a contained right word containing a target variable of value 0 or 1. That word is 0 if nonlinear or literal 0, and is 1 if it is literal 1. Its join with the nonlinear witness of value 2 is 0. Linear targets follow from a pure-variable witness and \(x\leq x^{2}\). This proves the lemma.

The word laws follow from the formula. Each absorption L, U, I, B, \(R_{1}\)--\(R_{4}\) satisfies the criterion; they are also independently checked directly on all assignments to their at most five variables.

\par\smallskip\noindent\textit{3. Constructing a word with the full target content.}\quad 
Normalize all words. Linear targets are handled by L and the word laws. A nonlinear target of singleton content already has an original nonlinear pure-variable witness by N, so its normalized word xx is already below P. Suppose now that C has at least two variables.

Start with a nonlinear original g in \(P_{C}\). Whenever a contained original v contains a target variable not yet in g, split \(g=ab\) nonemptily and use U to derive \(gv=abv\leq ab+v\leq P\). Normalize the result and continue. The word remains nonlinear. Condition C ensures that the process ends with a word g of content C, written in its first-occurrence order. It need not yet have the target order.

\par\smallskip\noindent\textit{4. Sorting using original prefix witnesses.}\quad 
Assume the initial segment p of g already equals the corresponding target prefix, and x is the next target variable. Write

\CaseDisplay{g=pAxB.}

If A is absent, the correct prefix already extends through x. Otherwise A is nonempty. By O choose an original \(q=rx\) t, at its first x, where c(r) is contained in c(p). Let \(R=AxB\) and \(U=tR\). The raw word U is nonempty because R contains x, even if normalization later removes some of its letters.

If p is nonempty, B (the displayed identity, not the suffix B) with \(X=p\), \(Y=q\), \(Z=R\) derives

\CaseDisplay{pqR\le pR+q=g+q\le P.}

Its normal form equals that of pxU: all letters of r have already occurred in p and are removed. If p is absent, r is absent as well. Split the nonlinear g into two nonempty factors and use I with \(X=q\) to derive \(qg\leq q+g\leq P\). Its normal form equals that of xU.

These known words permit the following four nonempty substitutions:

\textbullet{} If p and the suffix B are both nonempty, \(R_{1}\) with \(X=p\), \(Y=x\), \(Z=A\), \(T=B\), and \(u=U\) gives \(pxAB\leq pAxB+pxU\).

\textbullet{} If p is nonempty and B is absent, \(R_{2}\) with \(X=p\), \(Y=x\), \(Z=A\), \(u=U\) gives \(pxA\leq pAx+pxU\).

\textbullet{} If p is absent and B is nonempty, \(R_{3}\) with \(X=x\), \(Y=A\), \(Z=B\), \(u=U\) gives \(xAB\leq AxB+xU\).

\textbullet{} If p and B are absent, \(R_{4}\) with \(X=x\), \(Y=A\), \(u=U\) gives \(xA\leq Ax+xU\).

The first right word is the current g and the second is the word just constructed, modulo the proven word laws. The new word has the same content C and a longer correct target prefix. Repeat for the remaining target positions. After finitely many moves the word is w, so \(w\leq P\) is derived. No substitution is empty; in particular U is kept as the raw nonempty word tAxB in every substitution.

\par\smallskip\noindent\textit{5. Completeness and duality.}\quad 
This derives every valid word absorption, with no bound on variables, right summands, or original word lengths. Distribute arbitrary terms into sums of words. Each summand of one side of a valid identity is absorbed by the other side; the preceding construction derives all of these absorptions and hence the identity. This proves the theorem. Reversing all words gives the corresponding basis of the opposite.

The independent computations verify raw word semantics, separating conditions, and actual substitutions in every construction branch. They support the unrestricted proof. Classification of further catalogue entries requires equality of generated varieties in both HSP directions.

\subsection{Sixteen identities for 3217: nonlinear content coverage and prefix sorting}\label{case:102}
Constant-free signature (+, multiplication), relative to the additively idempotent semiring axioms. No bound on variables, word lengths, or polynomial sizes is assumed.

\par\smallskip\noindent\textit{1. Tables and word values.}\quad 
In catalogue order 0,\allowbreak{}1,\allowbreak{}2,\allowbreak{}3,\allowbreak{}4:

{\small\par\smallskip
\noindent addition        multiplication\par
\CaseLine{\ttfamily 0 0 0 0 4\qquad 0 0 0 3 4}
\CaseLine{\ttfamily 0 1 1 1 4\qquad 0 0 0 3 4}
\CaseLine{\ttfamily 0 1 2 2 4\qquad 0 0 2 3 4}
\CaseLine{\ttfamily 0 1 2 3 4\qquad 3 3 3 3 3}
\CaseLine{\ttfamily 4 4 4 4 4\qquad 4 4 4 4 4}
\smallskip}

The additive order is the chain \(3<2<1<0<4\). A literal keeps its assigned value. A nonlinear word takes its first factor in \{3,\allowbreak{}4\}, if any. Otherwise it takes 2 when every factor is 2, and 0 when some factor is 0 or 1. Indeed 3,\allowbreak{}4 are left zeros, a right reset replaces a preceding base product, and multiplication on \{0,\allowbreak{}1,\allowbreak{}2\} is 2 only on all-2 products and 0 otherwise. This proves the formula by induction. In particular the nonlinear word values are the same as for 3216, but the additive joins are different.

For target content C let \(N_{C}\) be the original NONLINEAR right summands of P whose contents lie in C. A literal target x is absorbed iff P has an original word of content \{x\}. A nonlinear target w is absorbed iff:

\textbullet{} N: the union of the contents of \(N_{C}\) is exactly C;

\textbullet{} O: for every first-occurrence target variable x, an original right word contains x whose first-x prefix content lies in the target first-x prefix. This witness may be literal and may have outside variables after that occurrence.

Necessity of N: if x is absent from all \(N_{C}\), assign \(x=1\), other target \(variables=2\), and \(outsiders=3\). The target takes 0. Every noncontained right word takes 3, contained nonlinear words take 2, and literals take 1,\allowbreak{}2, or 3. Thus P is at most 1, strictly below 0. Necessity of O: assign preceding target variables 2, \(x=4\), and all remaining variables 3. The target is 4, whereas every right word is below 4 unless it supplies O. A literal target without a pure-x witness is separated by \(x=0\) and all outsiders 3.

For sufficiency,3 is least. Target 4 is witnessed by O at its first reset variable; the witness prefix has no reset and the word takes 4. Target 2 has all target variables 2, so any \(N_{C}\) word takes 2. Target 0 has base-only variables and at least one x valued 0 or 1; N supplies a contained nonlinear word containing x, which also takes 0. Literal targets follow from a pure-x word and \(x\leq x^{2}\). This proves the exact criterion without finite bounds. Literal coverage does NOT suffice in N: \(xy\leq xx+y\) is false at \(x=2\), \(y=1\). This is the difference from 3216.

\par\smallskip\noindent\textit{2. The finite basis and normal forms.}\quad 
{\small\par\smallskip
\CaseLine{\(W1 xxy=xy\qquad W2 xyx=xy\)}
\CaseLine{\(W3 xyy=xy\qquad W4 xxx=xx\)}
\CaseLine{\(L  x\leq xx\)}
\CaseLine{\(U  xyz\leq xy+zz\qquad I xyz\leq xx+yz\)}
\CaseLine{\(B  xyz\leq xz+yy\)}
\CaseLine{\(R1 xyzt\leq xzyt+xyu\qquad R2 xyz\leq xzy+xyu\)}
\CaseLine{\(R3 xyz\leq yxz+xu\qquad R4 xy\leq yx+xu\)}
\CaseLine{\(J1 xyzt\leq xzyt+y\qquad J2 xyz\leq xzy+y\)}
\CaseLine{\(J3 xyz\leq yxz+x\qquad J4 xy\leq yx+x\)}
\smallskip}

The word formula proves W1--W4 and L. The exact N/O criterion proves all remaining laws. Each displayed side is nonempty; all subsequent substitutions are by nonempty words or terms.

The normal form is the first-occurrence list for at least two distinct variables, the square for singleton nonlinear words, and the literal for length one. Delete separated repetitions with \(xYx=xY\) (Y nonempty). Adjacent repetitions use W1 when a suffix exists, or W3 when a prefix exists. W4 reduces singleton powers to the square. These are explicit word-law derivations and preserve the criterion. Every nonlinear word v satisfies \(v^{2}=v\) under these laws; this statement excludes literals.

\par\smallskip\noindent\textit{3. Construct full nonlinear content.}\quad 
Normalize P, w. Literal targets follow from L and a pure-variable original word; nonlinear singleton targets already have a nonlinear pure-variable summand by N. Otherwise choose g in \(N_{C}\). While a target variable is missing from g, choose v in \(N_{C}\) containing it. Split \(g=ab\) nonemptily. U gives \(gv\leq g+v^{2}=g+v\leq P\). The result remains nonlinear and its content increases inside C. N guarantees termination with content exactly C. Normalize to a distinct-letter list g, whose order need not be correct.

\par\smallskip\noindent\textit{4. Sorting with nonlinear or literal original witnesses.}\quad 
Suppose the correct target prefix is p and the next target variable is x. Write \(g=p\) A x B, with A nonempty if x must move. Choose an original O witness \(q=r\) x V with c(r) contained in c(p).

If q is a literal, it is exactly x. Use J1, J2, J3, J4 according as both p, B are nonempty, only p is nonempty, only B is nonempty, or both are absent. For example J1 with \(X=p\), \(Y=x\), \(Z=A\), \(T=B\) gives p x A \(B\leq p\) A x \(B+x\leq P\). All four cases extend the correct prefix. No square of the literal witness is assumed available.

If q is nonlinear, put \(R=A\) x B and keep \(T=V\) R as a RAW nonempty word. When p is nonempty, B with \(X=p\), \(Y=q\), \(Z=R\) gives p q \(R\leq p\) \(R+q^{2}=g+q\leq P\). Its normal form is that of p x T, since r uses only earlier prefix variables. When p is absent, r is absent; split \(g=ab\) and I gives q \(g\leq q^{2}+g=q+g\leq P\), whose normal form is x T. Now use R1, R2, R3, R4 according to the same four p/B cases. Their right words are g and the derived bridge; their left word is p x A B. Every substitution is nonempty, including T even when V is absent.

The new word has the same target content and a longer correct prefix. Iterate to reach w. The induction has no variable, length, or summand bound and treats arbitrary original witnesses, not just small examples.

\par\smallskip\noindent\textit{5. Completeness and classification boundary.}\quad 
Every valid word absorption has been derived from sixteen fixed laws. Distribute arbitrary terms into nonempty polynomials. Apply the result to every summand in both directions to derive any valid identity. Thus the sixteen laws form a finite basis for 3217. Reversal gives the opposite basis. Other catalogue transfers require independent complete HSP inclusions in both directions.

\subsection{Terminal variables and affine parity: ten-identity bases for a semiring and its opposite}\label{case:103}
\par\smallskip\noindent\textit{Abstract.}\quad 
We give a ten-identity basis for an explicit noncommutative five-element additively idempotent semiring. Its word values retain exponent parity and distinguish a variable in the final position from its internal occurrences. The absorption criterion separates two cases according to whether some terminal variable of the right polynomial also occurs internally. Four identities combine pairs of words with the same terminal variable. These identities prove completeness in both cases for any number of variables. Reversing the words gives a basis for the opposite semiring. Exact variety certificates extend the conclusions to twelve five-element catalog types.

\par\smallskip\noindent\textit{1. The algebra and a finite basis.}\quad 
The signature is \((+,\cdot)\) without constants. Bases are relative to the additively idempotent semiring axioms, and \(u\le v\) abbreviates \(u+v\approx v\).

Let \(B=\{e, g, c\}\) be the flat extension of the group of order two: \(e, g\) multiply as the group, \(c\) is a multiplicative zero, and the join of distinct elements is \(c\). Extend \(B\) by elements \(a,\omega\), where \(\omega\) is the greatest additive element and \(a\) is incomparable with every element of \(B\). Thus \(b+a=\omega\) for \(b\in B\). Preserve multiplication in \(B\) and define

\CaseDisplay{ba=a\ (b\in B),\qquad
  au=\omega u=u\omega=\omega\ (u\in B\cup\{a,\omega\}).\quad\textup{(\ensuremath{1})}}

Call the resulting algebra \(A\). These element names are not term constants.

Theorem 1. The following ten identities form a basis for \(A\):

\CaseDisplay{\begin{array}{rll}
\mathrm{(P)}&xyz\approx yxz,\\
\mathrm{(N)}&x^3y\approx xy,\\
\mathrm{(S)}&x^2y^2\approx y^2x^2,\\
\mathrm{(E)}&x\le y^2x,\\
\mathrm{(O)}&xyz\le x^3+y^3+z^3,\\
\mathrm{(G)}&x^3\le x+y^3,\\
\mathrm{(T_3)}&xyzt\le xs+ys+zt,\\
\mathrm{(T_{2a})}&xyt\le xs+ys+t,\\
\mathrm{(T_{2b})}&yzt\le s+ys+zt,\\
\mathrm{(T_1)}&yt\le s+ys+t.
\end{array}\quad\textup{(\ensuremath{2})}}

Reversing every word in (2) gives a basis for \(A^{\mathrm{op}}\). No minimality assertion is made.

Associativity follows by considering a triple product. Any occurrence of \(\omega\), or an occurrence of \(a\) before the final position, gives \(\omega\) in either association. If only the last factor is \(a\), both products equal \(a\). The remaining case lies in \(B\).

For distributivity, left translation by an element of \(B\) preserves joins on \(B\), fixes \(a\), and fixes \(\omega\); it therefore also preserves the joins between \(a\) and \(B\). Left translation by \(a\) or \(\omega\) is constant \(\omega\). Right translation by an element of \(B\) preserves joins on \(B\) and sends \(a,\omega\) to \(\omega\). Right translation by \(a\) sends all of \(B\) to \(a\) and sends \(a,\omega\) to \(\omega\); right translation by \(\omega\) is constant. These observations prove both distributive laws.

For a word \(w\), let \(c(w)\) be its content, \(t(w)\) its final variable, \(i(w)\) its internal content (all positions except the final one), and \(\pi(w)\) its exponent-parity vector over \(\mathbb F_2\). In particular \(i(x)=\varnothing\) for a linear word. The value of \(w\) is determined by the following rule:

1. A variable assigned \(\omega\), or an internal occurrence assigned \(a\), gives value \(\omega\).

2. Otherwise a final letter assigned \(a\) gives value \(a\).

3. Otherwise all letters lie in \(B\); an occurrence of \(c\) gives \(c\), and group assignments give the product determined by \(\pi(w)\).

This proves (P),(N),(S). For (E), if \(y\in B\), then \(y^2\) is \(e\) or \(c\): multiplication on the left fixes \(a\) and increases no group value past \(c\). In particular \(x\le y^2x\). If \(y\notin B\), then \(y^2=\omega\) and the right side is top. To check (O), a letter outside \(B\) gives a top cube on the right; inside \(B\), a nontop right join has all three group letters equal, whose odd product is that same letter. All other core products are at most \(c\).

For (G), \(x^3=x\) unless \(x=a\). In that exceptional case \(x^3=\omega\), and \(a+y^3=\omega\) for every \(y\), since cubes lie in \(B\cup\{\omega\}\).

The last four laws have a common verification. Their right sides are \(ps+qs+rt\) and their left sides \(pqrt\), where one or more prefixes can be absent as specified in the four displayed laws. This notation does not introduce an empty-word substitution. A right value \(a\) forces \(s=t=a\) and every present prefix to lie in \(B\), so the left side is also \(a\). If the right value lies in \(B\), every variable used on the left lies in \(B\). If that right value is \(c\), absorption follows. If it is \(e\) or \(g\), all three right words have that same group value; the left parity is the sum of their parities, since \(s\) occurs twice on the right, so the left value is the same. A right value \(\omega\) is automatic. This verifies all ten laws.

\par\smallskip\noindent\textit{2. Word normalization.}\quad 
Call a word saturated if its final variable occurs earlier in the word. A nonsaturated word is either linear or has the form \(pt\) with \(p\) nonempty and \(t\notin c(p)\).

Lemma 2. The laws (P),(N),(S) identify any two saturated words of the same content and parity. They identify any two nonsaturated nonlinear words having the same terminal variable, internal content, and parity.

Proof. Law (P), in a multiplication context, interchanges any adjacent positions before the final position. Consequently the entire prefix can be permuted. Law (N) removes or inserts two copies of a prefix variable whenever another copy of it remains; its right substitution is always nonempty because the final position remains. For a nonsaturated nonlinear word, reduce each prefix exponent to one or two and order the prefix.

For a saturated word ending in \(x\), move an earlier \(x\) to the penultimate position, giving a suffix \(x^2\). To select another variable \(r\) as the repeated terminal variable, first use (N) backwards to ensure at least two prefix occurrences of \(r\). Permute the prefix to give the suffix \(r^2x^2\), and apply (S) to change it to \(x^2r^2\). Thus any member of the content can serve as repeated terminal variable. Choose the least variable in a fixed ordering. Normalize all other exponents to one or two, and that terminal variable's exponent to two or three according to its parity. The singleton cases are reduced using \(x^4=x^2\), an instance of (N), and its multiplication consequences. The resulting normal form depends only on content and parity. \(\square\)

Every saturated \(r\) therefore satisfies \(r^3\approx r\). Also (E) with \(y=x\) gives \(x\le x^3\), including after word substitution.

Lemma 3. Suppose \(R\) is a nonempty finite family of saturated words. For any word \(w\), the conditions

\CaseDisplay{c(w)\subseteq\bigcup_{r\in R}c(r),\qquad
\pi(w)\in\operatorname{aff}_{\mathbb F_2}\{\pi(r):r\in R\}\quad\textup{(\ensuremath{3})}}

give a derivation \(w\le\sum R\) from (P),(N),(S),(E),(O).

Proof. An affine combination over \(\mathbb F_2\) is an odd sum. Choose the corresponding odd product of members of \(R\), and add two more copies of every member. Its content is the entire right union, and its parity is \(\pi(w)\). Iterating (O) bounds this odd product by \(\sum R\), because all factor cubes normalize to the factors. The product is saturated, including the one-factor case.

First use \(w\le w^3\). Insert squares of any missing variables to the left of \(w^3\) by (E). This gives a saturated word with the same content and parity as the product. Lemma 2 identifies them. If no variables are missing there is no insertion step. \(\square\)

\par\smallskip\noindent\textit{3. Pairs with a common terminal variable.}\quad 
The four laws \((\mathrm T)\) give the following uniform operation on nonempty words. If \(ps\), \(qs\), and \(rt\) are words, where a prefix may be absent, then

\CaseDisplay{pqrt\le ps+qs+rt.\quad\textup{(\ensuremath{4})}}

If all three prefixes are present use \((\mathrm T_3)\); if only \(r\) is absent use \((\mathrm T_{2a})\); if precisely one of \(p, q\) is absent use \((\mathrm T_{2b})\) or \((\mathrm T_1)\) according to \(r\). The case with both \(p, q\) absent is reflexive. The possible exchange of \(p, q\) is permitted by (P). Thus every use of (4) is one of the displayed finite laws with nonempty terms substituted for its variables.

Let \(P\) be a nonempty sum of words, and put

\CaseDisplay{U=\bigcup_{u\in P}c(u),\qquad
J=\bigcup_{u\in P}i(u),\qquad T=\{t(u):u\in P\},\qquad
H=\operatorname{aff}_{\mathbb F_2}\{\pi(u):u\in P\}.}

Lemma 4. If \(J\cap T\ne\varnothing\), the ten laws derive \(u^3\le P\) for every summand \(u\) of \(P\).

Proof. Choose \(x\in J\cap T\). There is a right word \(v=ps\) whose nonempty prefix \(p\) contains \(x\), and a right word \(u=qx\) with prefix \(q\) possibly absent. Apply (4) to the pair \(v, v\) and the word \(u\):

\CaseDisplay{r=p^2qx\le v+u\le P.}

The word \(r\) is saturated, since its prefix contains \(x\). Hence \(r^3=r\) by Lemma 2. For any right summand \(z\), substituting \(z, r\) in (G) yields \(z^3\le z+r^3=z+r\le P\). \(\square\)

Lemma 5. Suppose \(J\cap T=\varnothing\). If a word \(w\) satisfies

\CaseDisplay{t(w)\in T,\qquad i(w)\subseteq J,\qquad \pi(w)\in H,\quad\textup{(\ensuremath{5})}}

then the ten laws derive \(w\le P\).

Proof. Each right word is \(ps\) with \(s\in T\) and all prefix variables in \(J\). Choose an odd representation of \(\pi(w)\) by right-word parity vectors. As the two sets \(J, T\) are disjoint, the coordinates in \(T\) force an odd number of chosen words with terminal variable \(t=t(w)\) and an even number with every other terminal variable. Select one of the chosen \(t\)-words as the initial word. Partition all remaining occurrences into pairs with equal terminal variable.

Repeatedly use (4), one pair at a time. The resulting word has terminal variable \(t\), prefix equal to the product of the prefixes selected so far, and is below \(P\). Its parity is \(\pi(w)\), because every paired terminal variable is removed twice. To ensure full prefix content \(J\), also use the pair \(v, v\) for each right summand \(v\). These extra pairs change no parity and add every right prefix variable.

We obtain \(W\le P\) with \(i(W)=J\), \(t(W)=t\), and \(\pi(W)=\pi(w)\). If \(J=\varnothing\), all words are linear; the initial word is \(t=w\) and the assertion is immediate. Otherwise insert to the left of \(w\) the squares of all missing variables of \(J\) by (E). The resulting word and \(W\) are nonsaturated nonlinear words with the same terminal variable, internal content, and parity, so Lemma 2 identifies them. If \(w\) is linear, (E) still supplies the required nonempty prefix; if no variables are missing, normalization alone suffices. Thus \(w\le W\le P\). \(\square\)

\par\smallskip\noindent\textit{4. Absorption and completeness.}\quad 
Theorem 6. For any word \(w\) and nonempty polynomial \(P\), \(A\models w\le P\) holds precisely when

\CaseDisplay{c(w)\subseteq U,\qquad \pi(w)\in H,\quad\textup{(\ensuremath{6})}}

and, if \(J\cap T=\varnothing\), additionally

\CaseDisplay{t(w)\in T,\qquad i(w)\subseteq J.\quad\textup{(\ensuremath{7})}}

All such absorptions are derivable from the ten identities.

Necessity. If \(w\) contains a variable absent from \(U\), assign that variable \(c\) and every other variable \(e\). The target value is \(c\) and the right value \(e\). If \(\pi(w)\notin H\), choose a linear functional constant on \(H\) and different on \(\pi(w)\). Assign the corresponding values \(e, g\). All right words have one group value, and the target the other; those values are incomparable.

Now suppose \(J\cap T=\varnothing\), and assume the content condition already established. Assign every variable in \(T\) the element \(a\) and every variable in \(J\) the element \(e\). Every right word then has value \(a\). If \(i(w)\) meets \(T\), the target is \(\omega\). If not, but \(t(w)\notin T\), then the target uses only \(J\) and has value \(e\). Neither \(\omega\) nor \(e\) is at most \(a\). As \(U=J\cup T\), this proves both restrictions in (7).

Sufficiency. If \(J\cap T\ne\varnothing\), Lemma 4 gives the saturated family \(R=\{u^3:u\in P\}\) below \(P\). This family has total content \(U\) and affine parity hull \(H\). Lemma 3 proves \(w\le\sum R\le P\) under (6). If \(J\cap T=\varnothing\), use Lemma 5. This proves the criterion and its equational completeness.

Distribute both sides of an arbitrary valid identity into polynomials. Each left summand is absorbed by the right polynomial, and conversely. Theorem 6 derives every one of these absorptions; taking sums proves both inequalities and hence the identity. This proves Theorem 1 without a bound on variables or summands. Word reversal carries every axiom, substitution, multiplication context, and derivation to the opposite semiring, proving its basis assertion as well.

\par\smallskip\noindent\textit{5. Catalog identification and verification.}\quad 
In the order \((e, g, c, a,\omega)\) the named-to-catalog permutation is \([1,\allowbreak{}2,\allowbreak{}0,\allowbreak{}3,\allowbreak{}4]\) for type 3489. Type 3967 has the same addition table and the transposed multiplication table. The two classes certified by complete bidirectional variety tests are

\CaseDisplay{\{3240,\allowbreak{}3489,\allowbreak{}4429,\allowbreak{}5046,\allowbreak{}7448,\allowbreak{}8191\},\qquad
\{3295,\allowbreak{}3967,\allowbreak{}4383,\allowbreak{}5120,\allowbreak{}8188,\allowbreak{}8192\}.}

The first class has basis (2); the second has its reversal.

The source manuscript for this case credits Aifa Wang and Lili Wang, School of Mathematical Sciences, Chongqing University of Technology.

The source manuscript for this case acknowledges support from the National Natural Science Foundation of China (12371024), the Science and Technology Research Program of Chongqing Municipal Education Commission (KJZD-K2024011102) and the Chongqing Natural Science Foundation Innovation and Development Joint Fund (Municipal Education Commission) (CSTB2025NSCQ-LZX0067).

\subsection{A finite basis from tail support and affine parity}\label{case:104}
Signature: (+, multiplication), without constants. All bases are relative to the additively idempotent semiring axioms. Catalogue source: 10924. The argument below applies to arbitrary finite sets of variables, with no bound on the input word lengths or polynomial sizes.

\par\smallskip\noindent\textit{The algebra and normal words.}\quad 
Let A have elements 0,\allowbreak{}1,\allowbreak{}2,\allowbreak{}3,\allowbreak{}4, with operations given by the following tables in that order. These are names of elements, not symbols in the term signature.

{\small\par\smallskip
\noindent addition          multiplication\par
\CaseLine{\ttfamily 0 0 0 0 4\qquad 4 0 0 0 4}
\CaseLine{\ttfamily 0 1 1 1 4\qquad 4 1 1 1 4}
\CaseLine{\ttfamily 0 1 2 1 4\qquad 4 1 2 3 4}
\CaseLine{\ttfamily 0 1 1 3 4\qquad 4 1 3 2 4}
\CaseLine{\ttfamily 4 4 4 4 4\qquad 4 4 4 4 4}
\smallskip}

Addition is the join operation with incomparable minimal elements 2 and 3, their join 1, and \(1<0<4\). Multiplication restricted to \{2,\allowbreak{}3\} is the group of order two with identity 2. On \{1,\allowbreak{}2,\allowbreak{}3\}, the element 1 is a multiplicative zero. A right factor 0 or 4 makes every product 4. A left factor 0 followed by a factor in \{1,\allowbreak{}2,\allowbreak{}3\} stays 0; a left factor 4 stays 4.

Associativity follows by considering whether any factor after the first belongs to \{0,\allowbreak{}4\}, then whether the first factor is 4 or 0, and otherwise using the associative multiplication on \{1,\allowbreak{}2,\allowbreak{}3\}. For distributivity, left multiplication by 0 or 1 collapses \{1,\allowbreak{}2,\allowbreak{}3\} to 0 or 1 and sends 0,\allowbreak{}4 to 4; left multiplication by 2 or 3 restricts to a join automorphism on \{1,\allowbreak{}2,\allowbreak{}3\} and sends 0,\allowbreak{}4 to 4; left multiplication by 4 is constant. Right multiplication by 0 or 4 is constant 4; right multiplication by 2 is the identity, by 3 interchanges 2 and 3 and fixes 1,\allowbreak{}0,\allowbreak{}4, and by 1 collapses \{1,\allowbreak{}2,\allowbreak{}3\} to 1 and fixes 0,\allowbreak{}4. All these maps preserve joins. Thus A is an additively idempotent semiring.

For a word u let h(u) denote its head, c(u) its content, T(u) the set of variables in positions after the first, and p(u) its total multiplicity vector modulo two, in the vector space over \(F_{2}\) on the variables. In particular T(u) may contain h(u), and T(u) is empty for a linear word. Under an assignment, u has value:

1. 4 if the head is 4 or any tail occurrence is 0 or 4;

2. 0 if the first case does not apply and the head is 0;

3. 1 if all occurrences are in \{1,\allowbreak{}2,\allowbreak{}3\} and at least one is 1;

4. otherwise the group value in \{2,\allowbreak{}3\} with parity p(u) evaluated on the variables assigned 3.

These cases include linear words. They prove the two word identities

{\small\par\smallskip
\CaseLine{\((W1) xyz \approx  xzy,\qquad (W2) xy^{3} \approx  xy.\)}
\smallskip}

Fix an ordering of the variables. A normal word consists of its head followed by an ordered tail in which each present variable occurs once or twice. W1 permits every adjacent interchange in the tail, using the nonempty preceding prefix as x. W2 deletes two copies from a tail multiplicity at least three, again with a nonempty prefix. Hence these laws normalize every word while retaining the head, tail support and total parity. The normal word is determined by these three data.

\par\smallskip\noindent\textit{An exact absorption criterion.}\quad 
For a nonempty polynomial P=sum(u in U)u, put \(D=\cup\) c(u) and \(E=\cup\) T(u). Write aff\{p(u):u in U\} for the affine span over \(F_{2}\). Concretely this consists of all sums of an odd number of the vectors; repeated summands cancel in pairs.

Lemma 1. For every nonempty word w, including a linear word, A satisfies w\ensuremath{\leq}P if and only if

{\small\par\smallskip
\CaseLine{\(c(w) \subseteq  D,\qquad T(w) \subseteq  E,\qquad p(w) \in  aff\{p(u):u in U\}.\qquad (C, T, P)\)}
\smallskip}

Here w\ensuremath{\leq}P abbreviates w+P\ensuremath{\approx}P.

Necessity. If a variable x of w is outside D, assign \(x=4\) and all other variables 2. Then \(w=4\) and \(P=2\). If x belongs to T(w) but not E, assign \(x=0\) and all other variables 2. The target is 4, while every right word is either 0 or 2, so P\ensuremath{\leq}0.

For the parity condition restrict all assignments to \{2,\allowbreak{}3\}. If p(w) is outside the affine span, the augmented vector (p(w),1) is outside the linear span of the vectors (p(u),1). Finite-dimensional linear separation gives a functional (e, b) which vanishes on each (p(u),1) and has value 1 on (p(w),1). Consequently all dot products e\ensuremath{\cdot}p(u) equal b, while e\ensuremath{\cdot}\(p(w)=b+1\). Assign variable x the group element 3 when \(e_{x}=1\) and 2 otherwise. Every right word has the same group value, and the target has the other group value. These two minimal elements are incomparable, so absorption fails.

Sufficiency. If P evaluates to 4 there is nothing to prove. If \(P=0\), no right word is 4. Therefore no variable of D is 4, and no variable of E is 0. The first two conditions ensure that w is not 4, whence w\ensuremath{\leq}0. If \(P=1\), every right word belongs to \{1,\allowbreak{}2,\allowbreak{}3\}, so all variables of D are in \{1,\allowbreak{}2,\allowbreak{}3\}. Content containment gives w\ensuremath{\leq}1. Finally, if P is 2 or 3, its minimality in the additive order forces every right word to have precisely that value. All variables of D then lie in \{2,\allowbreak{}3\}. An odd affine representation of p(w) shows that w has that same group value. This exhausts all assignments and proves the lemma.

\par\smallskip\noindent\textit{Chains with at most three right summands.}\quad 
Lemma 2. Every valid absorption w\ensuremath{\leq}P can be derived from W1, W2 and valid absorptions having at most three right words. In the resulting chain every right word is an original summand or an already derived word.

Proof. Normalize P and w. Use D and E from Lemma 1. If E is empty then every right word is linear, and T(w)\ensuremath{\subseteq}E makes w linear; c(w)\ensuremath{\subseteq}D means w is literally a summand. Assume E is nonempty. Choose an odd subfamily U\ensuremath{{}_{0}} of the right words such that the sum of its parity vectors is p(w).

For a word u in P assign a category as follows: its category is h(u) if h(u) is outside E, and otherwise a single common symbol *. A variable outside E can occur only as a head, and then exactly once. Thus the coordinate of such a variable x in the sum of the selected parity vectors is the parity of the number of selected words in category x. Since T(w)\ensuremath{\subseteq}E, exactly the category of h(w) is odd when h(w) is outside E, and every other outside category is even. If h(w) is in E all outside categories are even. Because |U\ensuremath{{}_{0}}| is odd, the * category is odd precisely in the latter case. Therefore the category of the target has an odd number of selected words, and every other category has an even number.

Choose a word s from that odd category. Remove it from U\ensuremath{{}_{0}}, and pair the remaining words within their categories. Write \(h=h(w)\). The rest of the proof derives words below P while maintaining s as the current word.

For each x in E missing from T(s), choose an original summand a with x in T(a). Then

{\small\par\smallskip
\CaseLine{\(s x^{2} \leq  s+a\)}
\smallskip}

is valid by Lemma 1: the new content and tail variable are covered by a, and the parity of the target is p(s). Normalize the new current word. Repeating this gives \(T(s)=E\).

Change its head to h while retaining tail support E and total parity p(s). This is possible: if h is outside E, the current head already equals h; otherwise both heads lie in E. Explicitly, after the new head h put one copy of x in E if the x coordinate of p(s)+\(e_{h}\) is 1, and two copies otherwise. This produces a word s' with the same content, tail support and parity as s, so s'\ensuremath{\leq}s is valid by Lemma 1.

For each remaining pair a, b of the same category, the vector p(a)+p(b) has zero coordinates outside E. Keep head h and tail support E, and choose the tail multiplicities one or two so that the new total parity is

{\small\par\smallskip
\CaseLine{\(p(s') = p(s)+p(a)+p(b).\)}
\smallskip}

The absorption s'\ensuremath{\leq}s+a+b is valid: its content and tail support are already covered by s, and its parity is the sum of the three displayed right vectors. Replace s by s' and repeat for all pairs. Now \(p(s)=p(w)\). The current word has head h and tail support E, so its content contains c(w). Lemma 1 gives w\ensuremath{\leq}s. Every step has at most three right words, and all are available at that stage. Addition and transitivity concatenate the steps to derive w\ensuremath{\leq}P. This proves the lemma for linear as well as nonlinear targets.

\par\smallskip\noindent\textit{A finite bound by reversible parity-pattern compression.}\quad 
Fix 84 variable symbols. Let Normal be all the nonempty normal words on them; each has length at most 169. Let \ensuremath{\Gamma} comprise every valid absorption w\ensuremath{\leq}\(u_{1}\)+...+\(u_{k}\) with 1\ensuremath{\leq}k\ensuremath{\leq}3 and all words in Normal. Validity is decidable either from Lemma 1 or directly from the table. This defines a specific finite set; we do not enumerate its very large list here.

Lemma 3. Every valid absorption with at most three right words follows from \ensuremath{\Gamma} and W1, W2.

Proof. Normalize the target and the at most three right words. Keep every variable which is a head of one of these words as a separate symbol, at most four symbols. For each other variable record its multiplicity in the tails of the at most four words. These entries are 0,\allowbreak{}1,\allowbreak{}2. Group variables with the same pattern. The all-zero pattern is absent, so there are at most 3\ensuremath{{}^{4}}\ensuremath{-}\(1=80\) groups, in addition to the four possible head variables.

Replace each group by one symbol with its common tail multiplicity in each word. This is not ordinary variable identification: multiplicities of members of a group are not added. Content and tail containment in Lemma 1 are preserved by the incidence patterns. Any odd affine representation of the original target parity remains an odd affine representation after compression: on a group coordinate it is simply the equation on any representative member; head variables remain separate. Hence the compressed absorption is valid and belongs to \ensuremath{\Gamma} after renaming, on at most 84 variables and with length at most 169.

Conversely substitute for each group symbol the nonempty ordered product of its original members, and leave each separate head variable fixed. No group occurs in a head position. A group with multiplicity one or two expands to exactly that many copies of each member. Tail commutation and normalization recover each original word, including its head, support and parity. Thus the original absorption follows from the compressed one.

Theorem. The two identities W1, W2 together with \ensuremath{\Gamma} form a finite identity basis for A relative to the additively idempotent semiring axioms. Word reversal gives a finite basis for the multiplicative opposite.

Proof. Soundness follows from the table and the definition of \ensuremath{\Gamma}. Lemmas 2 and 3 derive every valid one-word absorption. Distribute arbitrary terms as sums of words; an identity is equivalent to all its summand absorptions in both directions. This yields completeness for any number of variables. Word reversal transports both the identities and the derivations to the opposite algebra.

The bounds are proved without finite testing. Computational regressions only check for errors; they do not supply the completeness argument. This result does not claim a small, minimal or independent basis, or hereditary finite basability. Transfers to additional catalogue types require full two-way HSP equality and independent verification.

\subsection{Finite bases for 3251 and 3252 by internal coverage and tail chaining}\label{case:105}
Constant-free signature (+, multiplication), relative to the additively idempotent semiring axioms. The two bases below contain 25 and 20 explicitly specified laws; no minimality is claimed.

\par\smallskip\noindent\textit{1. Word values and exact absorption.}\quad 
Both multiplication tables have rows 00034,\allowbreak{}00034,\allowbreak{}01234,\allowbreak{}33333,\allowbreak{}44444 in the order 0,\allowbreak{}1,\allowbreak{}2,\allowbreak{}3,\allowbreak{}4. Addition for 3251 has rows 00004,\allowbreak{}01014,\allowbreak{}00224,\allowbreak{}01234, 44444; addition for 3252 has rows 00004,\allowbreak{}01114,\allowbreak{}01224,\allowbreak{}01234,\allowbreak{}44444. Thus 3251 has \(3<1\),\(2<0<4\) with \(1+2=0\), while 3252 has \(3<2<1<0<4\).

A word containing 3 or 4 has the value of its first such occurrence. With only 0,\allowbreak{}1,\allowbreak{}2 factors, it gives 2 if all factors are 2; it gives 1 if its final factor is 1 and all earlier factors are 2; otherwise it gives 0. This includes literals and follows by induction from the table.

For target w let \(C=c(w)\) and U be the original right words contained in C. Put \(A=\cup\) of c(u without its last position), u in U, and B=\{last(u):u in U\}. Write I(w)=c(w without its last position). The exact conditions for \(w\leq P\) are:

(C) The contents of U cover C. (E) Each first target occurrence x has an original right word whose first-x prefix content is contained in its preceding target content. (I2) For 3252, I(w) is contained in A. (I1) For 3251, either A intersects B, or I(w) is contained in A.

For a failed C choose an uncovered target variable x, assign \(x=0\), other target variables 2 and outsiders 3. For a failed E assign the preceding target content 2, the selected \(x=4\) and everything else 3. Each assignment separates the target from P. For a failed I2 choose x in I(w) outside A, assign \(x=1\), other target variables 2, outsiders 3. The target is 0, but each contained right word is 1 or 2. For a failed I1, A and B are disjoint; choose x in I(w) outside A, assign B union\{x\} to 1, other target variables 2, outsiders 3. Every contained right word is 1 because its internal variables are in A and its tail in B; the target is 0. Noncontained words give 3 in all these separators.

For sufficiency, target 3 is least and target 4 is supplied by E. A base-valued target with some 0 variable is supplied by its C-witness, which is 0. If all target variables are 2, any contained word gives 2. If the target is 1, its unique 1-valued variable is the simple tail; C supplies a word containing it, valued 1 or 0. Finally, target 0 with no 0 factor has a1-valued internal variable. For 3252 its I2 witness gives 0. For 3251, if no contained word gives 0, all contained internal factors are 2 and each word gives 1 or 2. If both values occur their join is 0. If all give 2, C contradicts the target 1 occurrence. If all give 1, every variable of B is 1 and every variable of A is 2; thus A intersects B is impossible, and I1 says every target internal variable is 2, again a contradiction. This proves the criterion for arbitrary words, variable sets and polynomial sizes.

A literal target is absorbed exactly when a pure-variable original word exists. This follows from C and \(x\leq x^{2}\).

\par\smallskip\noindent\textit{2. Common twenty laws and normalization.}\quad 
Eight fixed laws are:

{\small\par\smallskip
\CaseLine{\(W1 xxy=xy\qquad W2 xyx=xyy\)}
\CaseLine{\(L  x\leq xx\qquad H  xy\leq xx+y\)}
\CaseLine{\(M  xyz\leq xz+yz\qquad C1 xyzt\leq xy+zyt\)}
\CaseLine{\(C2 xyz\leq x+yxz\qquad C3 xyz\leq xy+yz\)}
\smallskip}

The remaining twelve laws are the following finite list:

{\small\par\smallskip
\CaseLine{\(xab\leq axb+x\qquad xaa\leq axx+x\)}
\CaseLine{\(xab\leq axb+xv\qquad xaa\leq axx+xv\)}
\CaseLine{\(pxab\leq paxb+x\qquad pxaa\leq paxx+x\)}
\CaseLine{\(pxab\leq paxb+xv\qquad pxaa\leq paxx+xv\)}
\CaseLine{\(prxab\leq praxb+rx\qquad prxaa\leq praxx+rx\)}
\CaseLine{\(prxab\leq praxb+rxv\qquad prxaa\leq praxx+rxv\)}
\smallskip}

All displayed letters are variables, not optional empty substitutions. Each law is valid for both tables by the exact criterion above. W1, W2 normalize a word to its first-occurrence list F when the final variable is simple, and to F followed by its last listed variable when the final variable repeats. Internal separated repeats use \(xYxZ=xY^{2}Z=xYZ\) with Y, Z nonempty; adjacent repeats use W1. A terminal repeat xYx becomes \(xY^{2}\) by W2 before internal reductions. Call the second normal form saturated. Every saturated word g has \(g^{2}=g\). No literal is assumed idempotent. All actual substitutions are nonempty.

\par\smallskip\noindent\textit{3. Constructing full content for 3252.}\quad 
Suppose first C is contained in A. Starting with any u in U, if its current final variable s is simple, choose an original q in U with s internally. Write the current word as p s and \(q=r\) s V, V nonempty. If p, r are nonempty, C1 gives p s r \(V\leq p\) s+r s V. If only p is nonempty, use C3. If only r is nonempty, use C2. If both are absent, q itself suffices. Normalize:the result contains both original contents and is the normal form of their concatenation. If its tail is still simple, it is a new variable not in the previous word; otherwise it is saturated. Thus each nonterminal step increases content, and the process must reach a saturated word. This is tail chaining through original internal occurrences, not an unbounded axiom. Perform it for each u in U. Combine saturated outputs with H: g \(v\leq g^{2}+v=g+v\). Products of saturated words remain saturated. The result has full content C.

Otherwise I2 implies w has a simple final variable t, the only possible member of C outside A. Chase tails as above until either saturation or the simple tail t is reached. Every step with a simple tail different from t is supported by I2. A word containing t can only have it at its final position:original words never use t internally, and chaining stops as soon as it reaches t. Therefore saturated outputs omit t; there is at least one simple-t output by C. Combine saturated outputs with H, then append one simple-t output with H. Combine remaining simple-t outputs using M at their common final t. A literal t contributes no new content and can be ignored; when it is the first output, start with a subsequent nonlinear one. A nonlinear target guarantees such a word, or a preceding saturated output. The resulting full-content g has simple tail t. All concatenations and M substitutions use nonempty prefixes; the literal cases were handled separately.

\par\smallskip\noindent\textit{4. Crossing activation for 3251.}\quad 
Adjoin the following five laws to the common twenty:

{\small\par\smallskip
\CaseLine{\(A0 xyy\leq xx+y\)}
\CaseLine{\(A1 xcxc\leq x+xc\)}
\CaseLine{\(A2 xbcxbc\leq x+bxc\)}
\CaseLine{\(A3 axcaxc\leq ax+xc\)}
\CaseLine{\(A4 axbcaxbc\leq ax+bxc\)}
\smallskip}

They are valid in 3251:the right words exhibit an intersection between their internal-variable union and tail set. This activation is NOT valid as a general inference for 3252.

If A intersects B, choose s in the intersection, an original \(u=p\) s ending in s, and an original \(q=r\) s V with V nonempty. Apply A1--A4 according to the presence of p, r. This derives the square of p s r V below u+q, so its normal form is saturated. For each missing variable choose an original v in U containing it. A0 gives g v \(v\leq g^{2}+v=g+v\). The result remains saturated and its content grows. Eventually it covers C.

If A and B are disjoint, I1 implies the target has simple final t and \(B=\{t\}\):otherwise every target variable would be internal, forcing B contained in A. Thus all original contained words end in t. Merge them by M, handling literal t separately as in Section 3. This produces a full-content simple-t word. No tail-chaining assumption from 3252 is required in this noncrossing case.

\par\smallskip\noindent\textit{5. Sorting with the original prefix witnesses.}\quad 
In either algebra the previous sections give an available full-content normal word g. It is either saturated or has the same simple tail as w. Maintain the correct first-occurrence prefix p of g. Let x be the next target variable. If it is already next, extend p. Otherwise write the first-occurrence list as p A x B with A nonempty. Choose an original E-witness \(q=r\) x V with c(r) contained in c(p).

Insert the raw block r immediately after p in the comparison. It does not change the normal form because all its variables already occur in p and it lies before further letters. Apply the appropriate one of the twelve displayed sorting laws:three prefix cases (p absent, p present with r absent, or p and r present), two cases for V, and two cases for whether B is empty. When B is nonempty use the actual suffix of g, including its repeated last letter if saturated. The left word moves x before A; the first right word normalizes to g, and the second is q.

When B is empty, the word must be saturated and g ends x x. The corresponding law changes p r A x x to p r x A A. Both normalize with the intended first-occurrence order and retain saturation. If g has simple tail t, we sort only its prefix and never move t, so B is always nonempty. Every move extends the correct prefix. All omitted contexts are represented by separate displayed laws, not empty substitutions.

If the result is saturated but the target has simple last variable t, L in the nonempty preceding context derives p \(t\leq p\) t t. A saturated target or a preserved simple-tail target is already obtained. Literals follow from L and the original pure-variable word. The proof covers arbitrary variable counts, lengths and numbers of right summands.

\par\smallskip\noindent\textit{6. Completeness and scope.}\quad 
Every valid absorption is derived by the exact criterion and Sections 3--5. Expand terms by distributivity and absorb summands in both directions. The common twenty laws form a finite basis for 3252; the common laws plus A0--A4 form a finite basis for 3251. Reversal gives their opposites. Any further catalogue transfer requires independently complete HSP in both directions. Finite tests and proof replay supplement the unrestricted argument and do not replace it.

\subsection{A lowest-value obstruction and two eleven-identity bases}\label{case:106}
\par\smallskip\noindent\textit{1. A one-sided extension of a band semiring.}\quad 
We work without constants, in the signature \((+,\cdot)\), relative to the axioms for additively idempotent semirings. Write \(u\leq v\) for \(u+v\approx v\). Let \(A=\{d, c, e, a, t\}\), with addition given by the join of

\CaseDisplay{d<c<a<t,\qquad e<a,}

where \(e\) is incomparable with \(d, c\). Its multiplication is

\noindent \(\cdot\)\enspace\textbar\enspace \(d\)\enspace\textbar\enspace \(c\)\enspace\textbar\enspace \(e\)\enspace\textbar\enspace \(a\)\enspace\textbar\enspace \(t\)\par
\noindent \(d\)\enspace\textbar\enspace \(c\)\enspace\textbar\enspace \(c\)\enspace\textbar\enspace \(d\)\enspace\textbar\enspace \(c\)\enspace\textbar\enspace \(c\)\par
\noindent \(c\)\enspace\textbar\enspace \(c\)\enspace\textbar\enspace \(c\)\enspace\textbar\enspace \(c\)\enspace\textbar\enspace \(c\)\enspace\textbar\enspace \(c\)\par
\noindent \(e\)\enspace\textbar\enspace \(c\)\enspace\textbar\enspace \(c\)\enspace\textbar\enspace \(e\)\enspace\textbar\enspace \(a\)\enspace\textbar\enspace \(t\)\par
\noindent \(a\)\enspace\textbar\enspace \(c\)\enspace\textbar\enspace \(c\)\enspace\textbar\enspace \(a\)\enspace\textbar\enspace \(a\)\enspace\textbar\enspace \(t\)\par
\noindent \(t\)\enspace\textbar\enspace \(t\)\enspace\textbar\enspace \(t\)\enspace\textbar\enspace \(t\)\enspace\textbar\enspace \(t\)\enspace\textbar\enspace \(t\)\par
The relabeling \((d, c, e, a, t)\mapsto(3,\allowbreak{}2,\allowbreak{}1,\allowbreak{}0,\allowbreak{}4)\) identifies this algebra with catalog semiring 4398. Element names are not term constants.

Here is a structural verification of the semiring axioms. The core \(B=\{c, e, a, t\}\) is a band with identity \(e\), left zero elements \(c, t\), and \(a^2=a\). A product in \(B\) is its first \(c\) or \(t\), if any; if neither occurs it is \(a\) when an \(a\) occurs, and \(e\) otherwise. This proves associativity. Its translations preserve the order and the only incomparable join \(c+e=a\), as the table verifies, so \(B\) is a semiring.

Let \(\rho(d)=c\) and fix \(B\). This is a join homomorphism, and \(xy=\rho(x)\rho(y)\) except for \(de=d\). The projections under \(\rho\) of both associations of a triple agree. Either association equals \(d\) exactly when \((x, y, z)=(d, e, e)\); otherwise both belong to \(B\) and hence agree. For distributivity, right multiplication by \(e\) is the identity; every other right translation factors through \(\rho\) and a translation of \(B\). Left translations with first factor other than \(d\) factor in the same way. Left multiplication by \(d\) takes \(e\) to \(d\) and every other element to \(c\); it preserves joins because \(u+v=e\) holds exactly when \(u=v=e\). Thus \(A\) is an additively idempotent semiring.

The value of a nonlinear word is as follows:

\textbullet{} If every variable has value \(e\), its value is \(e\).

\textbullet{} If every variable lies in \(\{e, a\}\) and at least one has value \(a\), its value is \(a\).

\textbullet{} Otherwise locate the first position whose value belongs to \(\{d, c, t\}\). If that value is \(t\), the word has value \(t\).

\textbullet{} In the remaining cases the value is \(d\) exactly when the first position has value \(d\) and every later position has value \(e\); otherwise it is \(c\).

These rules follow inductively from the displayed multiplication.

\par\smallskip\noindent\textit{2. Word normalization and the exact criterion.}\quad 
The two word identities

\CaseDisplay{xy^2\approx xy,\qquad xyx\approx x^2y\quad\textup{(\ensuremath{N})}}

hold in \(A\). They imply \(x^3\approx x^2\). Repetitions of a variable other than the first variable can be deleted: for a nonempty prefix \(p\), use \(pxrx\approx px^2r\approx pxr\) when the intervening word \(r\) is nonempty, and \(px^2\approx px\) when it is empty. A repetition of the first variable can be brought immediately after its first occurrence using the second law in (N). Repeated use of these operations gives two nonlinear normal forms, with \(x_1,\ldots, x_k\) distinct and in first-occurrence order:

\CaseDisplay{S(x_1,\ldots, x_k)=x_1x_2\cdots x_k\quad(k\geq2),
\qquad
R(x_1,\ldots, x_k)=x_1^2x_2\cdots x_k\quad(k\geq1).}

Call the latter regular words. They are idempotent under (N). Every nonlinear word \(u\) satisfies \(u^2\approx R\) with the same first-occurrence order; \(u\) itself has the \(R\) form exactly when its first variable occurs again. Nonfirst letters may be squared or unsquared by the first law in (N), so regular words may also be written \(x_1^2\cdots x_k^2\).

Write a nonempty polynomial in this normal form as \(F=\sum_{x\in L}x+\sum_{u\in N}u\). For a word \(w\), let \(c(w)\) be its content, \(h(w)\) its first variable, \(\tau(w)\) the set of variables in all positions after the first, and \(P_w(x)\) the set preceding its first \(x\). Define

\CaseDisplay{D=L\cup\{h(u):u\in N\},\qquad E=\bigcup_{u\in N}\tau(u).}

Theorem 1. A nonlinear word \(w\), with \(X=c(w)\), satisfies \(A\models w\leq F\) if and only if:

\textbullet{} (P) For every \(x\in X\), a summand \(u\) of \(F\) contains \(x\) and satisfies \(P_u(x)\subseteq P_w(x)\); a linear summand \(x\) is allowed.

\textbullet{} (A) A summand of \(F\) has content contained in \(X\); a linear summand is allowed.

\textbullet{} (D) Either \(D\cap E\ne\varnothing\), or \(h(w)\in D\) and \(\tau(w)\subseteq E\).

A linear word \(y\) satisfies \(y\leq F\) exactly when \(y\in L\) or a nonlinear summand of \(F\) has content \(\{y\}\).

Proof. If (P) fails at \(x\), assign \(x=t\), variables in \(P_w(x)\) the value \(e\), and all other variables \(c\). Then \(w=t\), while right summands containing \(x\) meet a \(c\) before it and take value \(c\). Other right summands have value \(c\) or \(e\), so \(F\leq a<t\).

If (A) fails, assign all variables of \(X\) the value \(e\) and all others \(c\). The left value is \(e\) and every right summand has value \(c\).

For (D), observe that a nonlinear word has value \(d\) if and only if its first variable has value \(d\) and every position after the first has value \(e\). Consequently all summands of \(F\) can simultaneously have value \(d\) if and only if \(D\cap E=\varnothing\): assign \(D\) the value \(d\) and \(E\) the value \(e\). Condition (P) implies \(X\subseteq D\cup E\). Under this assignment, \(w=d\) exactly when \(h(w)\in D\) and \(\tau(w)\subseteq E\). Thus failure of (D), assuming (P), separates the two sides.

Conversely, it suffices to check each proper principal ideal of the additive semilattice. For the bound \(a\), the only left value outside its ideal is \(t\); condition (P), applied to the first variable taking \(t\) after a neutral prefix, supplies a right summand taking \(t\). For the bound \(c\), a left value \(e\) or \(a\) means all variables of \(X\) lie in \(\{e, a\}\); (A) then supplies a right summand in \(\{e, a\}\), outside the ideal of \(c\). The case of a left value \(t\) was already handled. For the bound \(e\), any left value other than \(e\) requires a variable of \(X\) not taking \(e\). By (P) that variable occurs in a right summand, which therefore also has value other than \(e\). For the bound \(d\), if every right summand were \(d\), the preceding characterization and (D) would force \(w=d\). These are all proper principal ideals, proving \(w\leq F\).

For a linear left word \(y\), if neither stated alternative exists, assign \(y=e\) and all other variables \(c\). Every right summand is \(c\), which does not dominate \(e\). Conversely, \(x\leq x^2\) holds and \(x^n=x^2\) for \(n\geq2\), proving both stated alternatives. \(\square\)

\par\smallskip\noindent\textit{3. Eleven identities.}\quad 
Theorem 2. The following eleven identities form a complete basis for \(\operatorname{Id}(A)\) relative to the additively idempotent semiring axioms:

\CaseDisplay{\begin{array}{rll}
\mathrm{(N_1)}&xy^2\approx xy,\\
\mathrm{(N_2)}&xyx\approx x^2y,\\
\mathrm{(I)}&x\leq x^2,\\
\mathrm{(T)}&x^2y\leq xz+y^2,\\
\mathrm{(F_0)}&x^2\leq x+y^2,\\
\mathrm{(F_1)}&x^2y\leq xy+z^2,\\
\mathrm{(C_{00})}&x^2y\leq xy+y,\\
\mathrm{(C_{10})}&x^2yz\leq xyz+y,\\
\mathrm{(C_{01})}&x^2y\leq xy+yz,\\
\mathrm{(C_{11})}&x^2yz\leq xyz+yt,\\
\mathrm{(H)}&xy\leq x+xyz.
\end{array}}

No assertion of minimality is intended.

Validity. (\(N_{1}\)), (\(N_{2}\)), and (I) follow from the word formula. For (T), \(xz\) is a prefix witness for \(x\), and \(y^2\) is a prefix witness for \(y\) and an internal summand; its repeated head ensures \(D\cap E\ne\varnothing\). For (\(F_{0}\)), the linear summand \(x\) supplies (P) and (A), and \(y^2\) ensures (D). For (\(F_{1}\)), \(xy\) supplies (P) and (A), and \(z^2\) ensures (D). In each of the four (C) identities, the unchanged word \(xy\) or \(xyz\) supplies (P) and (A). The variable \(y\) occurs in its tail and is the head of the other summand, ensuring (D). For (H), the word \(xyz\) supplies all prefix witnesses and the linear summand \(x\) supplies (A). Here \(D=\{x\}\) and \(E=\{y, z\}\), so the left word \(xy\) meets the second alternative of (D). Theorem 1 proves all the inequalities.

We prove completeness in two cases.

\par\smallskip\noindent\textit{3.1. A shared head and tail variable permits regularization.}\quad 
Suppose \(D\cap E\ne\varnothing\). If a right nonlinear summand \(q\) is regular, then \(q^2\approx q\leq F\). Otherwise all nonlinear right words have form \(S\). Choose \(g\in D\cap E\). There is a summand \(u=pgr\) in which this occurrence of \(g\) is after the first position, so \(p\) is nonempty; \(r\) may be absent. There is also a right summand equal to \(g\) or of the form \(gs\) with \(s\) nonempty.

Use respectively \((C_{00})\), \((C_{10})\), \((C_{01})\), or \((C_{11})\), according as \(r, s\) are absent or present, with \(x=p\) and \(y=g\). The result is

\CaseDisplay{u^2\approx p^2gr\leq u+g\quad\hbox{or}\quad u+gs\leq F.}

In the case of absent \(r\), omit that substring. Thus in every case we have a nonlinear word \(q\) such that \(q^2\leq F\) is derivable.

For a linear summand \(g\) of \(F\), (\(F_{0}\)) gives \(g^2\leq g+q^2\leq F\). For a nonlinear summand \(u=hs\) with \(s\) nonempty, (\(F_{1}\)) gives \(u^2\approx h^2s\leq hs+q^2\leq F\). Therefore the polynomial \(F'\) obtained by squaring every summand of \(F\) satisfies \(F'\leq F\) derivably. It consists entirely of regular words. Squaring preserves contents and first-occurrence orders, so (P) and (A) still hold for \(w^2\) and \(F'\).

We now derive \(w^2\leq F'\) using (T). Normalize the left side to \(W=x_1^2\cdots x_k^2\), and let \(p_i=x_1^2\cdots x_i^2\). For each \(x_i\), (P) supplies a regular right summand \(u_i=a_i x_i b_i\) with \(c(a_i)\subseteq\{x_1,\ldots, x_{i-1}\}\). Square the nonfirst variables if necessary. Then \(b_i\) is nonempty; \(a_1\) is absent. Simultaneously for every regular word \(v\) for which \(v\leq F'\) is derivable, prove \(p_iv\leq F'\). For \(i=1\), (T) gives

\CaseDisplay{x_1^2v\leq x_1b_1+v^2\approx u_1+v\leq F'.}

For \(i>1\), multiply (T) by \(p_{i-1}\) on the left and normalize:

\CaseDisplay{p_iv\leq p_{i-1}x_ib_i+p_{i-1}v^2
\approx p_{i-1}u_i+p_{i-1}v\leq F'.}

The final step is the induction hypothesis applied to \(u_i\) and \(v\). An internal regular right summand \(v\) supplied by (A) has \(c(v)\subseteq c(W)\), so \(Wv\approx W\). This proves \(W\leq F'\). Finally (I) gives \(w\leq w^2\leq F'\leq F\).

\par\smallskip\noindent\textit{3.2. Disjoint head and tail sets.}\quad 
Suppose \(D\cap E=\varnothing\). By (D), the left word has a head \(h\in D\) which does not occur again, and all other variables belong to \(E\). Write it as \(w=h x_2\cdots x_k\), where \(k\geq2\), and put \(p_1=h\), \(p_i=h x_2^2\cdots x_i^2\) for \(i\geq2\). The first law in (N) gives \(p_k\approx w\).

All linear right summands belong to \(D\), and all nonlinear right words have a head in \(D\) and all their other variables in \(E\). For \(i\geq2\), the prefix witness supplied by (P) cannot be linear. Its head must belong to \(D\cap\{h, x_2,\ldots, x_{i-1}\}=\{h\}\). Write it

\CaseDisplay{u_i=h a_i x_i b_i,
\qquad c(a_i)\subseteq\{x_2,\ldots, x_{i-1}\},}

where \(b_i\) can be made nonempty by squaring the nonfirst variable \(x_i\). Both \(a_i\) and \(b_i\) have all variables in \(E\).

Simultaneously for every nonempty word \(s\) over \(E\) for which \(hs\leq F\) is derivable, we prove \(p_i s\leq F\). The assertion for \(i=1\) is its hypothesis. For \(i\geq2\), multiply (T), with \((x, y, z)=(x_i, s, b_i)\), on the left by the nonempty word \(p_{i-1}\):

\CaseDisplay{p_i s\leq p_{i-1}x_i b_i+p_{i-1}s^2
\approx p_{i-1}a_i x_i b_i+p_{i-1}s\leq F.}

The equality uses (\(N_{1}\)) and deletion of the nonhead variables of \(a_i\) that have already occurred in \(p_{i-1}\). The last step uses the induction hypothesis for \(a_i x_i b_i\) and \(s\). This retains the single occurrence of \(h\) throughout.

If (A) supplies a nonlinear internal right summand, it has form \(hr\) with \(r\) nonempty over \(c(w)\setminus\{h\}\). The induction gives \(p_k r\leq F\), and normalization gives \(p_k r\approx w\). Otherwise the internal linear summand in (A) must be \(h\). There is a nonlinear witness \(hr\) for \(x_2\), so \(r\) is nonempty and the induction gives \(wr\leq F\). Law (H), with \(x=h\) and \(y=x_2\cdots x_k\), now gives \(w\leq h+wr\leq F\).

This proves every nonlinear absorption in Theorem 1. Linear absorption uses (I), (N), and a matching one-variable right summand. Distributing an arbitrary valid identity into sums of words and applying these derivations to every summand proves that identity. Thus the eleven laws are complete for arbitrary numbers of variables and summands. All substitutions used above are nonempty; optional substrings are handled by separate displayed laws or omitted concatenations. \(\square\)

\par\smallskip\noindent\textit{4. The opposite and the scope of transfer.}\quad 
The catalog semiring 4444 has the same addition as 4398 and its transposed multiplication. Reversing every word in the eleven laws gives an identity basis for 4444. Equality of generated varieties transfers the bases to these classes:

\textbullet{} 4398 side: 3305, 3984, 4398, 5142, 6512, 6895, 6908, 8254.

\textbullet{} 4444 side: 3253, 4444, 4778, 5082, 6580, 6942, 6949, 8293.

The additional types 5832 and 8095 satisfy the respective laws but generate proper subvarieties. They are not classified by this transfer: finite basability does not in general pass to subvarieties.

The source manuscript for this case credits Aifa Wang and Lili Wang, School of Mathematical Sciences, Chongqing University of Technology.

The source manuscript for this case acknowledges support from the National Natural Science Foundation of China (12371024), the Science and Technology Research Program of Chongqing Municipal Education Commission (KJZD-K2024011102) and the Chongqing Natural Science Foundation Innovation and Development Joint Fund (Municipal Education Commission) (CSTB2025NSCQ-LZX0067).

\subsection{Finite bases for 3254 and 3255 with a protected simple tail}\label{case:107}
The signature is (+, multiplication), without constants. Relative to the additively idempotent semiring axioms, the bases below have 72 and 67 identities respectively. No minimality is asserted.

\par\smallskip\noindent\textit{1. Values, normal forms and the exact criterion.}\quad 
Both multiplication tables have rows 00034,\allowbreak{}00034,\allowbreak{}01234,\allowbreak{}33334,\allowbreak{}44434. Addition for 3254 has rows 00004,\allowbreak{}01014,\allowbreak{}00224,\allowbreak{}01234,\allowbreak{}44444; for 3255 it has rows 00004,\allowbreak{}01114,\allowbreak{}01224,\allowbreak{}01234,\allowbreak{}44444. A word with a factor 3 or 4 has the value of its LAST such factor. Otherwise the value is 2 when all factors are 2, is 1 when the final factor is 1 and all preceding factors are 2, and is 0 otherwise. This follows by table induction, including words of length one.

For an original polynomial P and target word w let \(C=c(w)\), and let U consist of the original words of P whose content is contained in C. Let A be the union of their nonfinal-position contents and B their set of final variables. Write I(w)=c(w without its final position). Then \(w\leq P\) holds exactly when:

1. The contents of U cover C.

2. For every last target occurrence of x there is an original word q containing x such that the content after its LAST x is contained in the target content after that occurrence.

3. For 3255, I(w) is contained in A. For 3254, this holds or A meets B.

Necessity has explicit separating assignments. If coverage fails at x, set \(x=0\), other target \(variables=2\), \(outsiders=3\). If condition 2 fails at x, set \(x=4\), variables later in the \(target=2\), all \(others=3\). Every right word then gives 2 or 3, while the target gives 4. If condition 3 fails for 3255, choose x in I(w) outside A, set \(x=1\), other target \(variables=2\), \(outsiders=3\). Every contained right word is 1 or 2, whose join in 3255 is at most 1; the target is 0. For 3254, A and B are disjoint. Set B union \(\{x\}=1\), other target \(variables=2\), \(outsiders=3\). Every contained word is 1 and the target is 0. Noncontained words give 3 in these base-value separators.

For sufficiency, a target 3 is least. For target 4, choose its last 4-valued variable occurrence after all 3/\allowbreak{}4-valued occurrences; the suffix witness gives 4. For a base-valued target, every contained word has base-valued factors. If the target contains a 0, its coverage witness gives 0. An all-2 target is covered by any contained word. A target 1 has just its simple tail valued 1; a covering word gives 1 or 0. The remaining case is target 0 with an internal 1. In 3255, condition 3 supplies a contained word with an internal 1, hence 0. In 3254, suppose no contained word gives 0. Its internal factors must all be 2. If both values 1 and 2 occur their join is 0. If all values are 2, coverage contradicts the target's 1-valued variable. If all are 1, A is valued 2 and B is valued 1, so A cannot meet B; condition 3 then rules out a target internal 1. This proves sufficiency without a bound on variables, lengths or polynomial size.

The two word laws \(xxy=xy\) and \(xyx=yxx\) give the following normal form: list variables in LAST-occurrence order, and repeat the final variable once precisely when the physical final variable of the original word occurs earlier too. For xYxZ with Y, Z nonempty use \(xYxZ=YxxZ=YxZ\); adjacent internal repeats use \(xxy=xy\). For a terminal separated repeat use \(xYx=Yxx\). Thus each internal elimination reduces length; terminal reduction places the final duplicate together. Call a normal word with doubled final variable saturated. Saturated g satisfies \(g^{2}=g\); literals are not assumed idempotent. Normalization preserves C, I, and last-occurrence suffix sets, so original witnesses can be normalized without changing the criterion.

\par\smallskip\noindent\textit{2. The common 67 identities.}\quad 
Seven fixed identities are:

{\small\par\smallskip
\CaseLine{\(xxy=xy\qquad xyx=yxx\)}
\CaseLine{\(x\leq xx\qquad xy\leq xx+y\)}
\CaseLine{\(xyz\leq xz+yz\qquad xzyt\leq xy+zyt\)}
\CaseLine{\(xyz\leq xy+yz\)}
\smallskip}

Twelve ordinary sorting identities are specified by three finite templates. Each of a, v may independently be omitted in writing a template, giving four separate identities in each row:

{\small\par\smallskip
\CaseLine{\([a]bxx \leq  [a]xbb + [v]x\)}
\CaseLine{\([a]bxp \leq  [a]xbp + [v]x\)}
\CaseLine{\([a]bxrp \leq  [a]xbrp + [v]xr\)}
\smallskip}

There are 48 further identities protecting a simple final variable t:

{\small\par\smallskip
\CaseLine{\([a]bx[c][d][p]t \leq  [a]xb[c][d][p]t + [v]x[c]t[d]\)}
\CaseLine{\([a]bx[c][p]t\qquad \leq  [a]xb[c][p]t\qquad + [v]x[c]tt\)}
\smallskip}

The first template has five independently optional displayed variables a, v, c, d, p, hence 32 identities; the second has four a, v, c, p, hence 16. The variables b, x, t are always present. Brackets specify a finite list of identities, NOT permission to substitute the empty word. All substitutions in the proof are nonempty. Validity of each identity follows directly from the criterion in Section 1, and can independently be decided by all assignments in the finite source table.

\par\smallskip\noindent\textit{3. Constructing full content for 3255.}\quad 
Suppose C is contained in A. Start with any u in U. While its current normal word is p s with simple final s, choose an original q in U with an internal s. Write \(q=r\) s V with V nonempty. If p, r are both nonempty, \(xzyt\leq xy+zyt\) gives p r s \(V\leq p\) s+r s V. If p is nonempty and r empty, \(xyz\leq xy+yz\) gives p s \(V\leq p\) s+s V. If p is empty, q already contains the current literal s and is itself available. The new word contains the union of the two contents. If its final variable is still simple it was absent from the old current word, so content strictly increases. Otherwise it is saturated. Since content stays in finite C, each such chase terminates at saturation. Perform this for every u in U. Merge saturated outputs by g \(v\leq g^{2}+v=g+v\), using \(xy\leq xx+y\). Their product is saturated.

If C is not contained in A, condition 3 implies that w has a simple final variable t and C minus \(A=\{t\}\). Chase each initial word until it is saturated or has simple final t. A word containing t can only end in t, since no original contained word has t internally and each chase stops on reaching t. Saturated outputs therefore omit t. At least one simple-t output exists by coverage. Merge saturated outputs first, append a simple-t output by \(xy\leq xx+y\), then merge remaining simple-t outputs by \(xyz\leq xz+yz\) with \(z=t\). A literal t adds no content and is ignored; if first in a list, begin instead with a nonlinear output. A nonlinear target ensures one exists (or a saturated output already constructed). We obtain full content, with the same simple tail as w.

These arguments also cover C contained in A with a simple target tail: the constructed word may be saturated; the final recovery handles it.

\par\smallskip\noindent\textit{4. The five additional identities for 3254.}\quad 
Adjoin \(xyy\leq xx+y\) and the four identities

{\small\par\smallskip
\CaseLine{\(([a][b]xc)^{2} \leq  [a]x + [b]xc\)}
\smallskip}

where a, b may independently be omitted, but c is always present. For example the both-present identity has left abxcabxc, not axbcaxbc.

If s lies in A intersection B, select original \(u=p\) s and \(q=r\) s V with V nonempty. The appropriate one of these four identities derives (p r s V)\ensuremath{^{2}} below u+q. Its normal form is saturated. For each missing content witness v in U, \(xyy\leq xx+y\) derives g v \(v\leq g^{2}+v=g+v\). This is saturated and increases content. Eventually all C is present.

If A and B are disjoint, condition 3 implies w has a simple final t and \(B=\{t\}\): every other target variable is internal and belongs to A, so cannot belong to B. Thus all U end at t. Merge them by \(xyz\leq xz+yz\), with literal t handled separately as above. The result has full content and the same simple final variable as w. This activation is specific to 3254: \(xyy\leq xx+y\) fails in 3255 at \(x=2\), \(y=1\).

\par\smallskip\noindent\textit{5. Right-to-left sorting, with tail protection.}\quad 
Let g be the full-content word constructed above. If it is saturated, start with empty correct suffix. Otherwise its tail is the target's simple tail t, and start with correct suffix t. Proceed through the target last-occurrence order from right to left. If the next x already immediately precedes the correct suffix, enlarge that suffix. Otherwise the last-occurrence list of g is A x B S with B nonempty and correct suffix S. Choose an original witness \(q=V\) x r, using its LAST x, so c(r) is contained in c(S).

For a saturated g, or a simple-tail g with t absent from r, the ordinary sorting templates move x to the right of B. When S is empty, g ends x x, and \(bxx\leq xbb+[v]x\) is used with \(b=B\). When S is nonempty use the second or third template, according as r is empty or nonempty. Substitute for p the actual correct suffix, including its final repeat if g is saturated. Inserting r just before this suffix leaves its normal form unchanged. In the simple-tail case this assertion uses t absent from r. The first right word normalizes to g; the other is q. The resulting last-occurrence order ends x S and retains saturation or the simple final t as applicable.

It is essential to handle t in r separately when g has simple final t. Blindly inserting all of r would repeat t and change the normal form. Write \(S=P\) t. Normalize q first. If r ends t t, write \(r=C\) t t and use the second protected template with \(c=C\), \(p=P\). Otherwise t occurs once in r, so write \(r=C\) t D and use the first protected template with \(c=C\), \(d=D\), \(p=P\). Both C, D omit t and their contents are contained in P. Consequently inserting only C, D before P t leaves g's simple final t intact. The first right term normalizes to g, the second is the original normalized witness q; the left term moves x past B and retains t simple. Each absent A, V, C, D, P selects a separately listed finite identity. This also covers repeated variables within a substituted block.

Every sorting step enlarges the correct suffix, so finitely many steps give the target order. If the result is saturated but w has simple tail t, apply \(x\leq xx\) in its preceding context to derive p \(t\leq p\) t t. If w is saturated or its simple tail was protected, it is already obtained. For a literal target x, coverage supplies an original pure x word; normalization gives x or xx, and \(x\leq xx\) finishes.

\par\smallskip\noindent\textit{6. Completeness.}\quad 
The construction derives every valid absorption \(w\leq P\) from the listed finite identities. Expand arbitrary terms by distributivity and apply the construction in both directions, one summand at a time. Therefore the 67 common identities are a complete basis for 3255 and the 72 are a complete basis for 3254. All original-witness restrictions, nonempty substitutions and simple-tail boundaries have been respected. There is no bound on variable count, word length or number of right summands. Word reversal proves the corresponding result for each opposite. Other catalogue types require complete HSP certificates both ways. Finite semantic regressions and literal proof replay only supplement this unrestricted argument.

\subsection{Prefix cleanup: bases for sources 3273 and 3274}\label{case:108}
Work relative to the additively idempotent semiring axioms in the constant-free signature. The labels 0,\allowbreak{}1,\allowbreak{}2,\allowbreak{}3,\allowbreak{}4 name elements, not constants in terms. The multiplication is common to both algebras:

{\small\par\smallskip
\CaseLine{\ttfamily 00034}
\CaseLine{\ttfamily 00034}
\CaseLine{\ttfamily 22234}
\CaseLine{\ttfamily 33333}
\CaseLine{\ttfamily 44444}
\smallskip}

The two addition tables have rows

3273      3274 00004     00004 01014     01114 00224     01224 01234     01234 44444     44444

The first additive order is 3 below incomparable 1 and 2, with their join 0, and greatest element 4. The second is the chain \(3<1<2<0<4\). Both tables satisfy the semiring axioms: the join tables are as stated, each left and right multiplication preserves joins, and associativity follows directly from the word evaluation below, or by distinguishing a first occurrence of a factor in \{3,\allowbreak{}4\}. No square-homomorphism assumption is made.

\par\smallskip\noindent\textit{Word values and identities.}\quad 
For a nonlinear word, if some variable is assigned 3 or 4, its value is that of the first such occurrence. If all variables are assigned 0,\allowbreak{}1,\allowbreak{}2, its value is 2 when its head is assigned 2, and 0 otherwise. The square map has values (0,\allowbreak{}0,\allowbreak{}2,\allowbreak{}3,\allowbreak{}4).

The following laws are common to both bases:

\CaseDisplay{x^2y\approx xy,\qquad xy^2\approx xy,\qquad xyx\approx xy;\quad\textup{(\ensuremath{W})}}

\CaseDisplay{x\leq x^2;\quad\textup{(\ensuremath{S})}}

\CaseDisplay{xyzt\leq xt+yz;\quad\textup{(\ensuremath{I})}}

\CaseDisplay{xyz\leq xz+y;\quad\textup{(\ensuremath{L})}}

\CaseDisplay{xyz\leq xz+xy,\qquad xyz\leq xz+xyt.\quad\textup{(\ensuremath{P})}}

For source 3273 add

\CaseDisplay{xyz\leq yz+x.\quad\textup{(\ensuremath{H})}}

Theorem. The eight common laws form a basis for 3274. Adding (H) gives a nine-law basis for 3273. Word reversal gives bases for the opposite semirings. These assertions concern completeness, not minimality or independence.

Laws (W) normalize a nonlinear word to the list of its first variable occurrences, retaining a square if the support is a singleton. An adjacent repetition can be removed by the first or second law, provided the word has at least one further position. A separated repetition is removed by the third law with the intervening nonempty block substituted for y. Repeating proves normalization without a length bound. Linear words remain distinct from squares.

\par\smallskip\noindent\textit{Exact absorption criterion.}\quad 
Write P as a sum of linear variables L and nonlinear words N. For a nonlinear target w put \(C=c(w)\), \(h=head(w)\), and N\_C=\{u in N: c(u) is contained in C\}. Let \(pre_{w}\)(x) denote the variables before the first x.

Lemma. In either source a linear y is absorbed by P if and only if y is in L or N contains a word with singleton content \{y\}. A nonlinear w is absorbed if and only if both conditions below hold:

\textbullet{} For every x in C, either x is in L, or there is v in N containing x with \(pre_{v}\)(x) contained in \(pre_{w}\)(x). This prefix witness is not required to have content in C.

\textbullet{} In 3274, \(N_{C}\) contains a word with head h. In 3273, either that same condition holds, or h is in L and \(N_{C}\) is nonempty.

Necessity. If a linear y has neither witness, assign \(y=2\) and every other variable 3. Every right nonlinear word contains a 3, so has value 3; the linear right part also has value 3. The target is 2.

For a failed prefix requirement at x, assign \(pre_{w}\)(x) the value 0, assign \(x=4\), and all other variables 3. The target is 4. Every nonlinear right word that contains x encounters a 3 before its first x; words without x cannot be 4. No linear right word is x. Hence the right side is at most 0.

In 3274, if there is no contained head witness, assign \(h=1\), the remaining variables in C the value 2, and outsiders 3. The target is 0. A nonlinear right word containing an outsider is 3, and a contained word has head different from h and value 2. All linear summands are at most 2, so the right side is strictly below 0.

In 3273, first suppose h is not linear and no contained head witness exists. Assign \(h=0\), C minus \{h\} the value 2, and outsiders 3. The target is 0 and the right side is at most 2. In the remaining failure case h is linear but \(N_{C}\) is empty. Assign C the value 1 and outsiders 3. The target is 0; every right nonlinear word is 3 and the linear right part is at most 1. Both cases separate the target.

Sufficiency. A linear target is handled by its literal occurrence or by (S). For a nonlinear target, value 3 is least. If it is 4, use its first occurrence x assigned 3 or 4. All preceding variables are neutral (0,\allowbreak{}1,\allowbreak{}2). A nonlinear prefix witness also first encounters 4 at x, even if later letters are outsiders. Alternatively its linear witness is 4.

If the target has value 0 or 2, all C-variables are neutral. A contained word with the same head has precisely the target value. This proves sufficiency in 3274 and the ordinary case of 3273. In the alternative case of 3273, h is linear and a contained nonlinear word u is available. If the target is 2, u is 0 or 2 and suffices. If the target is 0, then h is 0 or 1. The case \(h=0\) is immediate, while \(h=1\) combines with \(u=0\) or \(u=2\) to give 0. This proves the lemma for arbitrary variable sets.

All displayed laws are sound by the word-value formula and this criterion. For (P), xz supplies the contained head witness and the requirements at x, z; xy or xyt supplies the prefix requirement at y. For (I), xt supplies the head witness, while yz supplies the middle prefix requirements. Rule (L) uses its linear y. Rule (H) uses the alternative head condition available in 3273. Finally (S) follows from the additive order.

\par\smallskip\noindent\textit{Equational completeness: insert, then clean up.}\quad 
Let \(w\leq P\) be a valid nonlinear absorption. Normalize its first-occurrence list as \(x_{1}\),..., \(x_{d}\). Put \(p_{i}=x_{1}\cdots x_{i}\). We construct \(s_{i}\leq P\), where \(s_{i}\) is nonlinear, its content is contained in C, and its first-occurrence list starts with \(p_{i}\).

Initially a contained nonlinear head witness can be chosen as \(s_{1}\). In the alternative 3273 case choose u in \(N_{C}\) and split \(u=yz\) into nonempty factors. Rule (H) gives \(x_{1}u\leq u+x_{1}\leq P\), providing \(s_{1}\). The alternative is not used for 3274.

Suppose \(s_{i}\) exists and \(x=x_{i+1}\). Write \(s_{i}\) approximately \(p_{i}\) q with q nonempty. If the normalized word ends with \(p_{i}\) and \(i\geq 2\), insert the repetition \(p_{i}\) by (W), using \(p_{i}^{2}\) approximately \(p_{i}\). If \(i=1\) and the support is singleton, the square already supplies a nonempty q. Thus no empty word is substituted into a law.

For a linear witness x, rule (L) yields \(p_{i}\) x \(q\leq p_{i}\) \(q+x\leq P\). This new word uses only C and starts with \(p_{i}\) x after normalization.

Otherwise choose a nonlinear prefix witness v, possibly containing outsiders after x, and split \(v=ab\) into nonempty factors. Rule (I) yields

\CaseDisplay{p_i v q\leq p_i q+v\leq P.}

Normalize \(p_{i}\) v q. All variables before the selected x in v already occur in \(p_{i}\) and disappear as repetitions, so the normal form is \(p_{i}\) x t with t nonempty or omitted. The first or second rule (P), respectively, now gives

\CaseDisplay{p_i x q\leq p_i q+p_i x\leq P,}

or

\CaseDisplay{p_i x q\leq p_i q+p_i x t\leq P.}

Thus the auxiliary outsider variables are removed, and the resulting word again has content contained in C. The two separate rules ensure that no empty term is used for t. This completes the induction.

At \(i=d\) the constructed word has precisely the target's first-occurrence list and content. By (W) it equals w, including the singleton-support case. Linear absorptions follow from an existing summand or (S) followed by singleton normalization. Distributive expansion and derivation of every summand absorption in both directions prove all identities, without bounds on variables, word length, or the number of right summands. Reversal carries the derivations to the opposite algebras.

\par\smallskip\noindent\textit{Scope of computational checking.}\quad 
Finite assignment tests and derivation regressions are supplementary checks. The unrestricted proof is the preceding induction. Transfer to catalogue entries requires independently verified two-way HSP equality. No claim of hereditary finite basability is made.

\subsection{Eleven identities for the neutral-tail and last-reset semiring}\label{case:109}
Source 3276. The signature is (+, multiplication), without constants, relative to the additively idempotent semiring axioms.

\par\smallskip\noindent\textit{1. The algebra and the basis.}\quad 
The tables in order 0,\allowbreak{}1,\allowbreak{}2,\allowbreak{}3,\allowbreak{}4 are

{\small\par\smallskip
\noindent addition       multiplication\par
\CaseLine{\ttfamily 0 0 0 0 4\qquad 0 0 0 3 4}
\CaseLine{\ttfamily 0 1 0 1 4\qquad 0 0 0 3 4}
\CaseLine{\ttfamily 0 0 2 2 4\qquad 2 2 2 3 4}
\CaseLine{\ttfamily 0 1 2 3 4\qquad 3 3 3 3 4}
\CaseLine{\ttfamily 4 4 4 4 4\qquad 4 4 4 3 4}
\smallskip}

The additive order is \(3<1\),\(2<0<4\), with 1 and 2 incomparable. Right multiplication by 0,\allowbreak{}1,\allowbreak{}2 is r=(0,\allowbreak{}0,\allowbreak{}2,\allowbreak{}3,\allowbreak{}4); by 3 it is constant 3, and by 4 constant 4. These maps preserve joins. Left multiplication also preserves joins, as the displayed rows show. Associativity follows by composing these three right maps: r is idempotent and fixes 3 and 4; the last constant map determines a composition. Equivalently it is a direct finite-table verification.

Theorem. The following eleven identities form a relative basis:

\CaseDisplay{\begin{aligned}
xyzy&\approx xzy,&xy^2&\approx xy,&x^2y&\approx xy,\quad\textup{(\ensuremath{W})}\\
x&\le x^2,\quad\textup{(\ensuremath{L})}\\
xyz&\le xy+tz,\quad\textup{(\ensuremath{A})}\\
xyz&\le xz+y,\quad\textup{(\ensuremath{B})}\\
xyz&\le xy+z,\quad\textup{(\ensuremath{C})}\\
xy&\le x+y^2,\quad\textup{(\ensuremath{D})}\\
xy&\le x+yxy,\quad\textup{(\ensuremath{E})}\\
xyz&\le x+zyz,\quad\textup{(\ensuremath{F})}\\
xyz&\le x+zxyz.\quad\textup{(\ensuremath{G})}
\end{aligned}}

Reversing all words gives a basis for the multiplicative opposite.

\par\smallskip\noindent\textit{2. Values and an exact absorption criterion.}\quad 
For a nonlinear word, a tail factor in \{3,\allowbreak{}4\} resets the current value to itself. The last such tail factor therefore determines the value. If none exists, the value is r applied to the head. This follows inductively from the three right multiplication maps.

The first two laws of (W) retain the literal head and the tail variables in their last-occurrence order, without repetition within that tail. To remove an earlier duplicate tail letter, use \(xyzy=xzy\) when it has a nonempty intervening block, or \(xy^{2}=xy\) when adjacent. The prefix before that earlier tail occurrence is nonempty. Call this a normal word. The third law of (W) additionally removes an immediate duplicated head when a nonempty suffix follows. We do not collapse the word xx to x; the normal form may retain a duplicated leading head.

For a word u containing x, let \(suf_{u}\)(x) be the content strictly after its last occurrence of x. Let P be a nonempty polynomial.

Lemma. A linear y is absorbed by P exactly when an original right word has content \{y\}. For a nonlinear normal target \(w=hT\), absorption \(w\leq P\) holds exactly when:

\textbullet{} (H) an original right word has head h and content contained in c(w);

\textbullet{} (N) an original nonlinear right word has content contained in c(w);

\textbullet{} (O) for each x in T, an original right word containing x has \(suf_{u}\)(x) contained in \(suf_{w}\)(x).

Necessity. For a linear target without a pure-variable witness, assign \(y=2\) and all outsiders 3. The target is 2 and every right word is 3. If H fails, assign \(h=0\), the other target variables 2, and all outsiders 3. The target is 0. A right word headed by h must contain an outsider and has value 3; the others have value 2 or 3. Thus \(P\leq 2\). If N fails, assign all target variables 1 and all outsiders 3. The target is 0; nonlinear right words are 3 and linears are 1 or 3, so \(P\leq 1\).

For failure of O at x, assign \(x=4\), the later target suffix variables 2, and everything else 3. The target is 4. A right word containing x has a variable outside the permitted suffix after its last x, so its value is 3. A word without x is 2 or 3. This separates the absorption. These assignments also handle a repeated target head by inspecting its tail occurrence.

Sufficiency. A last target tail reset 3 gives the least element. A last tail reset 4 is witnessed by O: the corresponding original word has that occurrence followed only by neutral-valued variables, and has value 4. Now suppose all target tail values lie in \{0,\allowbreak{}1,\allowbreak{}2\}. If the head is 3, the target is least. If the head is 4, H gives a same-head word valued 4. If it is 0, H gives a word valued 0. If it is 2, N supplies a word valued 0 or 2, bounding the target 2. If the head is 1, the target is 0. A nonlinear H witness is 0; otherwise H can be the literal h of value 1, and N supplies 0 or 2, whose join with 1 is 0. Linear targets follow from L and the stabilization of powers. This proves the criterion.

The three word laws follow from the value formula. Each of L--G satisfies the criterion, or can be checked directly on all assignments to its at most four variables. Thus all eleven laws are sound.

\par\smallskip\noindent\textit{3. Appending the target tail.}\quad 
Let \(N_{P}\) be the set of original nonlinear right words. It is nonempty by N. We construct \(uS\leq P\) simultaneously for every u in \(N_{P}\), starting with the empty suffix S and extending S backwards through the target tail. An empty suffix here means no context is appended, never an empty substitution for a variable.

Suppose x is the next target tail variable. By O choose an original witness v whose last x has only variables from S after it. If v is nonlinear, write \(v=R\) x Q at that last occurrence. The invariant gives \(vS\leq P\), and the first two word laws give vS approximately R x S.

If R is nonempty, split any original nonlinear \(u=ab\) with both parts nonempty. Law A, in right context S, gives

\CaseDisplay{abxS\le abS+RxS\le P.}

If R is empty, \(v=xQ\) is nonlinear and Q is nonempty, with content contained in S. Therefore S is nonempty. The invariant gives \(xS\leq P\), and \(x^{2}S\) approximately xS by the third word law. The same instance of A now uses \(t=x\) and has right term xxS.

If the witness v is literal x, then for nonempty S use B with \(X=u\), \(Y=x\), \(Z=S\) to give \(uxS\leq uS+x\). When S is empty, use C after splitting \(u=ab\): \(abx\leq ab+x\). Every substitution in these cases is nonempty.

This constructs the invariant for xS. Repetition through the tail gives \(uT\leq P\) for all original nonlinear u.

\par\smallskip\noindent\textit{4. Recovering the head.}\quad 
If H has a nonlinear witness v, its content is contained in c(hT), so vT normalizes to hT or hhT. In the latter case T is nonempty and \(x^{2}y=xy\) gives hhT approximately hT. Thus \(w\leq P\).

Otherwise choose H as the literal h and choose an original nonlinear N-witness v, with head k. If \(k=h\) the preceding argument finishes. If k differs from h, then k belongs to T. Write \(T=AB\) where B begins with k. The known word vT normalizes to kT, or to khT precisely when h is absent from T but occurs in the tail of v.

When A is absent, D with \(X=h\), \(Y=T\) gives \(hT\leq h+T^{2}\), and \(T^{2}\) normalizes to kT. In the extra-h case use E: \(hT\leq h+ThT\), whose nonlinear right word normalizes to khT.

When A is nonempty, F with \(X=h\), \(Y=A\), \(Z=B\) gives \(hAB\leq h+BAB\); BAB normalizes to \(kAB=kT\). In the extra-h case, G gives \(hAB\leq h+BhAB\), and the right nonlinear word normalizes to \(khAB=khT\). Both A and B are nonempty in these latter two substitutions. This completes every nonlinear absorption. A linear target is derived from its pure-variable right word using L and the word laws.

\par\smallskip\noindent\textit{5. Completeness.}\quad 
The preceding steps derive every valid word absorption with arbitrarily many variables and arbitrarily long words. Distributivity expresses any term as a finite sum of words. An identity between two such sums is equivalent to each summand on either side being absorbed by the other sum; all those absorptions have been derived. This proves the theorem. Word reversal transports both validity and derivations to the opposite algebra.

Finite checks independently verify table semantics and replay actual substitutions, including both directions of \(x^{2}y=xy\) only with a nonempty y. They support, and do not replace, this unrestricted proof. Transfer to other catalogue types requires both complete HSP directions.

\subsection{Eight identities for 3277: a nonlinear head seed and last suffix witnesses}\label{case:110}
Constant-free signature (+, multiplication), relative to the additively idempotent semiring axioms. Completeness below is unrestricted in variables, word length and number of right summands.

\par\smallskip\noindent\textit{1. Word values and exact absorption.}\quad 
Addition has rows 00004,\allowbreak{}01114,\allowbreak{}01224,\allowbreak{}01234,\allowbreak{}44444, giving the chain \(3<2<1<0<4\). Multiplication has rows 00034,\allowbreak{}00034,\allowbreak{}22234,\allowbreak{}33334,\allowbreak{}44434. If a word has a factor 3 or 4, its value is its LAST such factor. Otherwise a nonlinear word is 2 when its HEAD is 2, and is 0 when its head is 0 or 1. A literal retains its value. Table induction proves this for all lengths. The distinction at literal 1 is essential.

For an original polynomial P and a nonlinear target w, \(w\leq P\) holds if and only if:

(H) There is an original NONLINEAR right word v with the same head as w and content contained in c(w).

(S) For every last target occurrence of x, there is an original word q containing x whose content after its LAST x is contained in the target content after that occurrence. The occurrence in q may be its head, and q may be literal. Outsiders before this occurrence are allowed.

If H fails, assign the target head 1, other target variables 2 and outsiders 3. The nonlinear target is 0. Every right word is 1,\allowbreak{}2 or 3: noncontained words have a reset 3, contained words with a different head are 2, and a same-head contained word could only be literal 1. Their join is below 0. Thus replacing H by a literal same-head seed would be incorrect. If S fails at x, assign \(x=4\), the later target variables 2 and everything else 3. The target is 4. Every right word containing x has a later 3; other right words are 2 or 3. This separates.

For sufficiency, target 3 is least. Target 4 is supplied by S at its last 4-valued occurrence, after the last other reset. If the target has only base values 0,\allowbreak{}1,\allowbreak{}2, the contained nonlinear same-head witness has exactly the same base value: 2 or 0. These exhaust the cases.

A literal target x is absorbed exactly when an original right word has pure content \{x\}. Necessity uses \(x=0\), \(outsiders=3\). Sufficiency uses \(x\leq x^{2}=x^{3}=\)... from the tables. For a nonlinear target hT with h absent from T, H automatically supplies the S witness for h; the constructive argument therefore needs S only for target tail letters.

\par\smallskip\noindent\textit{2. Basis and normal form.}\quad 
The eight identities are:

{\small\par\smallskip
\CaseLine{\(W1 xyzy=xzy\qquad W2 xyy=xy\)}
\CaseLine{\(W3 xxy=xy\qquad L  x\leq xx\)}
\CaseLine{\(A0 xyz\leq xy+z\qquad A1 xyz\leq xy+tz\)}
\CaseLine{\(B0 xyz\leq xz+y\qquad B1 xyz\leq xz+tyz\)}
\smallskip}

The word-value formula proves W1-W3 and L. The H/S criterion verifies the four absorptions: their first right word is a nonlinear same-head contained seed, and the other word supplies the newly inserted variable with the required last-occurrence suffix. Each is also independently checkable on every assignment in the finite table.

W1 and W2 delete earlier occurrences strictly in the tail, retaining the head and the tail letters in last-occurrence order. W1 handles a nonempty intervening block, W2 adjacent repeats. If the resulting tail starts with the head h and has further letters, W3 deletes that first tail h. If the whole word is a pure power, its nonlinear normal form is hh; it is never replaced by the literal h. This gives a canonical normal form for the construction. All substitutions are nonempty.

\par\smallskip\noindent\textit{3. Append the target suffix to all original nonlinear words.}\quad 
Let N be the original nonlinear summands. By H it is nonempty. Normalize the target as hT, with T nonempty and with distinct letters. Build T from right to left. Maintain \(uS\leq P\) for EVERY original u in N, where S is the suffix already built. Initially S is absent and these are just the original summands. Absence of context is not an empty word substituted for an identity variable.

For the first (final) variable x, an original S-witness ends in x. If it is literal x, write \(u=UV\) with U, V nonempty and use A0 to give u \(x\leq u+x\). If it is nonlinear R x, R is nonempty and A1 gives u \(x\leq u+R\) x. Do this for every u in N.

Suppose S is now nonempty and x is the next variable. Choose an original witness \(q=R\) x Q at its last x, with c(Q) contained in c(S). If q is literal x, use B0 with \(X=u\), \(Y=x\), \(Z=S\) to derive u x \(S\leq u\) S+x. If q is nonlinear, the invariant gives \(qS\leq P\). When R is nonempty, W1-W2 give \(qS=R\) x S. B1 with \(T=R\) gives u x \(S\leq u\) S+R x S.

When q is nonlinear but R is empty, \(q=xQ\) and Q is nonempty. Then \(qS=xS\) by tail deletion. Use B1 with \(T=x\): its second right word is x x \(S=x\) S by W3, since S is nonempty. Hence the same construction works. This head-occurrence case does NOT require an empty substitution or an unavailable literal x. These alternatives cover all original witnesses, and preserve the invariant for the extended suffix xS.

Only original nonlinear words receive appended suffixes. Original literal witnesses stay available unchanged; the invariant makes no claim that arbitrary literals may be lengthened. The procedure terminates after the finite number of distinct letters in T, with no universal bound.

\par\smallskip\noindent\textit{4. Recover the head and prove completeness.}\quad 
Choose H's original nonlinear same-head contained witness v. The invariant gives \(vT\leq P\). Every tail variable of v other than possibly h has a later occurrence in T, so W1-W2 remove those earlier occurrences. If h appears in T too, all earlier tail h occurrences disappear as well. If h is absent from T, the remaining word is hT or hhT, and W3 gives \(hhT=hT\) because T is nonempty. Therefore \(vT=w\) modulo W1-W3. For a literal target, its pure original witness normalizes to x orxx, and L derives x when needed.

Thus every valid word absorption is derivable from the eight identities. Distributivity expands all terms into nonempty polynomials, and deriving each summand's absorption in both directions proves every identity. The eight identities form a finite basis for 3277. Word reversal gives a basis for its opposite. Catalogue transfers require complete HSP in both directions. Finite testing and literal replay support this unrestricted proof and do not replace it.

\subsection{Nonlinear coverage and first-occurrence prefixes: two seven-identity bases}\label{case:111}
\par\smallskip\noindent\textit{1. The algebra.}\quad 
We use the signature \((+,\cdot)\) without constants and work relative to the axioms for additively idempotent semirings. The notation \(u\leq v\) means \(u+v\approx v\).

Let \(A=\{d, b, c, a, t\}\), with addition the join of

\CaseDisplay{d<b<a<t,\qquad d<c<a<t,}

where \(b\) and \(c\) are incomparable. The subset \(B=\{d, c, a, t\}\) is a chain. Define \(\rho(b)=a\), fixing the other four elements, and let \(xy=\rho(x)\rho(y)\) using the following multiplication on \(B\):

\noindent \(\cdot\)\enspace\textbar\enspace \(d\)\enspace\textbar\enspace \(c\)\enspace\textbar\enspace \(a\)\enspace\textbar\enspace \(t\)\par
\noindent \(d\)\enspace\textbar\enspace \(d\)\enspace\textbar\enspace \(c\)\enspace\textbar\enspace \(c\)\enspace\textbar\enspace \(c\)\par
\noindent \(c\)\enspace\textbar\enspace \(c\)\enspace\textbar\enspace \(c\)\enspace\textbar\enspace \(c\)\enspace\textbar\enspace \(c\)\par
\noindent \(a\)\enspace\textbar\enspace \(c\)\enspace\textbar\enspace \(c\)\enspace\textbar\enspace \(a\)\enspace\textbar\enspace \(t\)\par
\noindent \(t\)\enspace\textbar\enspace \(t\)\enspace\textbar\enspace \(t\)\enspace\textbar\enspace \(t\)\enspace\textbar\enspace \(t\)\par
The element names are not constants in the signature. The relabeling \((d, b, c, a, t)\mapsto(3,\allowbreak{}1,\allowbreak{}2,\allowbreak{}0,\allowbreak{}4)\) gives catalog semiring 4342.

For a conceptual verification of the semiring axioms, let \(K=\{c<a<t\}\) have maximum as addition, multiplicative identity \(a\), and left zero elements \(c, t\). This multiplication is associative and its left and right translations preserve the chain order, so \(K\) is an additively idempotent semiring. Let \(E=\{0<1\}\) have maximum as both operations. In \(K\times E\), the four elements

\CaseDisplay{d=(c,0),\quad c=(c,1),\quad a=(a,1),\quad t=(t,1)}

form a subsemiring with the displayed table and the indicated chain order. The map \(\rho:A\to B\) preserves joins and fixes \(B\). Defining multiplication through this retraction therefore preserves associativity and both distributive laws. In particular, \(A\) is an additively idempotent semiring.

For a nonlinear word \(w\) (a word of length at least two), the value is:

1. \(d\) if every variable in \(w\) is assigned \(d\);

2. \(a\) if every variable in \(w\) is assigned a member of \(\{b, a\}\);

3. \(t\) if the first position whose value is outside \(\{b, a\}\) has value \(t\);

4. \(c\) in all remaining cases.

This follows from the product construction, or by induction using the table. Nonlinear words always take values in the chain \(B\); their value depends only on the sequence of first occurrences of their variables.

\par\smallskip\noindent\textit{2. An exact absorption criterion.}\quad 
For a word \(w\), write \(c(w)\) for its content. For \(x\in c(w)\), let \(P_w(x)\) be the set of variables occurring before the first occurrence of \(x\). Write a nonempty polynomial as

\CaseDisplay{F=\sum_{x\in L}x+\sum_{u\in N}u,}

where the words in \(N\) are nonlinear; either \(L\) or \(N\) may be empty.

Theorem 1. A nonlinear word \(w\) with \(X=c(w)\) satisfies \(A\models w\leq F\) if and only if all three conditions below hold:

\textbullet{} (C) Each \(x\in X\) occurs in at least one word of \(N\).

\textbullet{} (P) For each \(x\in X\), either \(x\in L\), or there is \(u\in N\) with \(x\in c(u)\) and \(P_u(x)\subseteq P_w(x)\).

\textbullet{} (A) Either some \(q\in N\) has \(c(q)\subseteq X\), or \(L\cap X\ne\varnothing\).

A linear word \(y\) satisfies \(A\models y\leq F\) if and only if \(y\in L\) or there is a nonlinear word of \(N\) with content \(\{y\}\).

Proof. Suppose (C) fails at \(x\). Assign \(x=b\) and all other variables \(d\). The left side is \(a\) if \(X=\{x\}\) and \(c\) otherwise. All right nonlinear summands are \(d\) and every linear summand is at most \(b\). Neither \(a\) nor \(c\) is below \(b\), so the absorption fails.

Suppose (P) fails at \(x\). Assign \(x=t\), variables in \(P_w(x)\) the value \(b\), and all other variables \(d\). The left side is \(t\). There is no linear summand \(x\). Every right nonlinear word containing \(x\) meets a variable outside \(P_w(x)\cup\{x\}\) before its first \(x\), so its value is \(c\). Right nonlinear words without \(x\) have value at most \(a\); linear summands also have value at most \(a\). Thus \(F\leq a<t\).

Suppose (A) fails. Assign variables of \(X\) the value \(b\) and all other variables \(d\). The left side is \(a\). Every right nonlinear summand contains an outside variable, and hence has value \(d\) or \(c\); every linear summand is \(d\). Thus \(F\leq c<a\).

Conversely, assume (C), (P), and (A), and fix an assignment. If \(w=d\), the inequality is automatic since \(d\) is the additive bottom. If \(w=c\), some \(x\in X\) has value other than \(d\). A nonlinear word supplied by (C) containing \(x\) cannot have value \(d\), so it is at least \(c\).

If \(w=a\), all variables of \(X\) lie in \(\{b, a\}\). An internal nonlinear word supplied by (A) then has value \(a\). Otherwise choose \(g\in L\cap X\). If \(g=a\) there is nothing to prove. If \(g=b\), condition (C) supplies a nonlinear word containing \(g\); its value is at least \(c\), and \(b+c=a\). Finally, if \(w=t\), let \(x\) be the first variable whose assigned value is outside \(\{b, a\}\); it is \(t\). In (P), either the linear summand \(x\) has value \(t\), or the supplied nonlinear word has only neutral values before its first \(x\), and consequently also has value \(t\).

For a linear left side \(y\), if neither stated alternative holds, assign \(y=b\) and all other variables \(d\). Each right nonlinear summand has value at most \(c\), as does every linear summand, whereas \(b\nleq c\). For the converse note that \(x\leq x^2\) holds and \(x^n=x^2\) for \(n\geq2\). \(\square\)

\par\smallskip\noindent\textit{3. The finite basis.}\quad 
Theorem 2. The following seven identities form a complete basis for \(\operatorname{Id}(A)\) relative to the additively idempotent semiring axioms:

\CaseDisplay{\begin{array}{rll}
\mathrm{(S)}&xy\approx x^2y^2,\\
\mathrm{(R)}&xyx\approx xy,\\
\mathrm{(I)}&x\leq x^2,\\
\mathrm{(T)}&xy\leq xz+y^2,\\
\mathrm{(E)}&xy\leq y+xyz,\\
\mathrm{(J_0)}&xy\leq x+zxt+y^2,\\
\mathrm{(J_1)}&zxy\leq x+uxv+zy.
\end{array}}

No minimality assertion is made about this basis.

Validity. (S), (R), and (I) follow from the construction and the word formula. For (T), the nonlinear summands \(xz, y^2\) supply both coverage and prefix witnesses for \(x, y\), and \(y^2\) is an internal nonlinear summand. For (E), \(xyz\) covers \(x, y\) and supplies their prefix witnesses, while the linear summand \(y\) supplies (A). For (\(J_{0}\)), \(zxt\) covers \(x\), \(y^2\) covers \(y\), the linear summand \(x\) and the word \(y^2\) supply (P), and \(y^2\) supplies (A). For (\(J_{1}\)), \(uxv\) covers \(x\) and \(zy\) covers \(z, y\); the linear summand \(x\) supplies its prefix condition, and \(zy\) supplies the prefix conditions for \(z, y\) and the internal summand. Theorem 1 verifies each inequality. That theorem was proved semantically, independently of any proposed basis.

Normalization. From (R) with \(y=x\) we obtain \(x^3\approx x^2\), hence \(x^n\approx x^2\) for \(n\geq2\). Applying (S) with \(x^2\) substituted for \(x\), and reducing \(x^4\), yields \(x^2y\approx x^2y^2\approx xy\). Similarly \(xy^2\approx xy\). The law (R) deletes an occurrence of a variable that has appeared earlier, and these two derived identities allow squares to be inserted or reduced in a nonlinear word. Consequently every nonlinear word is equal to

\CaseDisplay{W(x_1,\ldots, x_k)=x_1^2\cdots x_k^2,}

where \(x_1,\ldots, x_k\) are its distinct variables in first-occurrence order. This also holds for a word on a single variable. Every nonlinear word is idempotent. If \(p\) is nonlinear and \(c(q)\subseteq c(p)\), then \(pq\approx p\). More generally, repeated variables can be deleted when concatenating words, retaining their first-occurrence order. No assertion identifies a linear variable \(x\) with \(x^2\).

Prefix induction. Let \(w=W(x_1,\ldots, x_k)\) satisfy Theorem 1 and put \(p_i=x_1^2\cdots x_i^2\). Normalize all nonlinear summands of \(F\). We prove that, for every nonlinear word \(v\) whose inequality \(v\leq F\) is derivable, the inequality

\CaseDisplay{p_iv\leq F\quad\textup{(\ensuremath{1})}}

is derivable for each \(1\leq i\leq k\). The induction is simultaneous for all such words \(v\).

If \(x_i\) has a nonlinear prefix witness \(u_i\) in (P), write \(u_i=a_i x_i b_i\), where \(c(a_i)\subseteq\{x_1,\ldots, x_{i-1}\}\). The suffix \(b_i\) is nonempty because every first occurrence was squared in normalization. When \(i=1\), the prefix \(a_i\) is absent. Law (T) gives

\CaseDisplay{p_1v\approx x_1v\leq x_1b_1+v^2\approx u_1+v\leq F.}

When \(i>1\), multiply an instance of (T) on the left by \(p_{i-1}\):

\CaseDisplay{p_iv\approx p_{i-1}x_iv
\leq p_{i-1}x_ib_i+p_{i-1}v^2
\approx p_{i-1}u_i+p_{i-1}v\leq F.}

Here the last inequality uses the induction hypothesis for \(u_i\) and \(v\). The normalization rules absorb \(a_i\) into \(p_{i-1}\).

It remains to consider \(x_i\in L\) with no nonlinear prefix witness. Condition (C) supplies a nonlinear summand \(u_i=a_i x_i b_i\) containing \(x_i\). Its normalized suffix \(b_i\) is nonempty; its prefix \(a_i\) is also nonempty, since otherwise it would be a prefix witness. For \(i=1\), an instance of (\(J_{0}\)) gives

\CaseDisplay{p_1v\approx x_1v\leq x_1+a_1x_1b_1+v^2
\approx x_1+u_1+v\leq F.}

For \(i>1\), an instance of (\(J_{1}\)) gives

\CaseDisplay{p_iv\approx p_{i-1}x_iv
\leq x_i+a_ix_ib_i+p_{i-1}v
=x_i+u_i+p_{i-1}v\leq F.}

The final step uses the induction hypothesis for \(v\). Thus (1) holds in all cases. Every substitution above uses nonempty words. An absent prefix is handled by the separate \(i=1\) case, not by a unit constant.

Completing absorption. If (A) supplies an internal nonlinear summand \(q\), then (1) gives \(wq\leq F\) and normalization gives \(wq\approx w\). Otherwise take \(g\in L\cap c(w)\) and any nonlinear summand \(q\); at least one exists by (C). Equation (1) gives \(wq\leq F\). Substituting \((x, y, z)=(w, g, q)\) in (E), and using \(wg\approx w\), yields

\CaseDisplay{w\approx wg\leq g+wgq\approx g+wq\leq F.}

Linear absorption follows from a matching linear summand, or from (I) and the normalization of a one-variable nonlinear word. Thus every valid absorption in Theorem 1 is derivable. Distribute both sides of an arbitrary valid identity and derive every summand's absorption by the other side. The resulting two inequalities prove that identity. This proves completeness for arbitrary numbers of variables and summands. \(\square\)

\par\smallskip\noindent\textit{4. The opposite and verified transfers.}\quad 
The catalog semiring 4360 has the same addition as 4342 and the transposed multiplication table. Reversing each word in the seven laws gives a basis for this opposite algebra: reversal carries identities and equational derivations bijectively between the two algebras.

Independent certificates of equality of generated varieties transfer these bases to the following two classes:

\textbullet{} 4342 side: 3986, 3987, 4342, 5155, 6066, 6495, 6911, 8267, 14525.

\textbullet{} 4360 side: 4360, 5116, 5117, 5173, 6098, 6577, 6961, 8305, 14552.

The transfer requires both inclusions. The proof does not assert that all subvarieties of either variety are finitely based.

The source manuscript for this case credits Aifa Wang and Lili Wang, School of Mathematical Sciences, Chongqing University of Technology.

The source manuscript for this case acknowledges support from the National Natural Science Foundation of China (12371024), the Science and Technology Research Program of Chongqing Municipal Education Commission (KJZD-K2024011102) and the Chongqing Natural Science Foundation Innovation and Development Joint Fund (Municipal Education Commission) (CSTB2025NSCQ-LZX0067).

\subsection{Seven identities for 4025: original head-tail incidence activation}\label{case:112}
Signature (+, multiplication), no constants, relative to the additively idempotent semiring axioms. The proof is unrestricted.

\par\smallskip\noindent\textit{Word values and normal form.}\quad 
Addition rows:00044,\allowbreak{}01034,\allowbreak{}00244,\allowbreak{}43434,\allowbreak{}44444. Multiplication rows:00044,\allowbreak{}11144,\allowbreak{}00244,\allowbreak{}44344,\allowbreak{}44444. The additive order has \(1<0<4\),\(2<0<4\), and \(1<3<4\). Thus \(0+3=2+3=4\). For w with head h and tail content T, any tail 3/\allowbreak{}4 or head 4 gives 4. Otherwise the tail values lie in\{0,\allowbreak{}1,\allowbreak{}2\}. Heads 0 and 1 retain their value; head 2 gives 2 when all tails are 2, and 0 otherwise; head 3 gives 3 when all tails are 2, and 4 otherwise. Empty tails satisfy the all 2 clause. This formula follows by table induction and includes literal words.

W1 \(xyz=xzy\) and W2 \(xyy=xy\) follow. They sort tail variables and delete repeated tail occurrences, leaving the head separate. A head appearing also in the tail remains recorded. Fix a variable order; the normal form is h followed by the sorted set T. This is a literal sequence of W1 swaps with nonempty prefix and W2 deletions with nonempty prefix.

\par\smallskip\noindent\textit{Exact criterion from ORIGINAL words.}\quad 
Let H be all original heads, D all variables occurring in original tails, and \(T_{x}\) the union of original tails on words with head x. Put

{\small\par\smallskip
\CaseLine{\(A = D \cup  \{x in H : T_{x} intersects H\}.\)}
\smallskip}

Then \(w\leq P\) exactly when: (1) h(w) belongs to H; (2) T(w) is contained in A; (3) either h(w) belongs to A, or T(w) is contained in \(T_{h}\)(w). These definitions use the original P only, not a repeatedly enlarged graph.

Necessity of(1):set \(h=2\) and every other variable 1. Every right word with a different head has value 1. The target has value 2 or 0, neither below 1. For failure of(2), choose x in T(w) outside A. Set \(x=3\), all heads other than x to 1, and remaining variables 2. Since x is not in D, no right tail contains x. If x is a head, its original tails avoid H (otherwise x in A), so its right words give 3. Other right heads give 1. The right sum is 1 or 3, but the target tail x forces 4. For failure of(3), let h outside A and choose t in T(w) outside \(T_{h}\). Set \(h=3\), \(t=1\), all other heads 1, and remaining variables 2. Here t is not h: conditions(1),(2) already hold and h is outside A. Since h is not in D, no right tail contains h; since \(T_{h}\) avoids H and t, all original h-words give 3. Every other right word gives 1. The target is 4 and the sum is 3. Thus every failed condition has an explicit separating assignment.

Sufficiency will follow from the unrestricted derivation below, using only the seven laws whose soundness is proved directly from word values. This avoids assuming a graph criterion merely because finite tests pass.

\par\smallskip\noindent\textit{Seven sound laws.}\quad 
{\small\par\smallskip
\CaseLine{\(W1 xyz=xzy\qquad W2 xyy=xy\)}
\CaseLine{\(L  x\leq xy\qquad M  xyz\leq xy+xz\)}
\CaseLine{\(E  xx\leq xy+y\qquad F  yy\leq xy+y\)}
\CaseLine{\(S  xz\leq xx+yz\)}
\smallskip}

L holds because multiplying a value only preserves it or moves it upward in the displayed additive order. M holds by splitting the tail set:if a tail 3/\allowbreak{}4 occurs, one right term is 4; otherwise the head's \(change2\to 0\) \(or 3\to 4\) already occurs in one right term. E:for x not 3, \(xx=x\) and L applies; for \(x=3\), if \(y=2\) then \(xy+y=3+2=4\), otherwise \(xy=4\). F:for y not 3, \(yy=y\); for \(y=3\), \(xy=4\). For S, if \(x=3\) or 4 then \(xx=4\). For \(x=0\) or 1, \(xx=x\); unless \(z=3/\allowbreak{}4\) the left remains x, and for \(z=3/\allowbreak{}4\) the right \(yz=4\). For \(x=2\), \(xx=2\); \(z=2\) gives left 2, \(z=3/\allowbreak{}4\) gives right 4, and \(z=0/\allowbreak{}1\) gives left 0 while yz is 0,\allowbreak{}1, or 4, whose join with 2 is 0 or 4. This exhausts all cases. W1, W2 were proved from the unrestricted word-value formula. Thus no semantic circularity is used in the sufficiency proof.

\par\smallskip\noindent\textit{Derive every word satisfying(1)-(3).}\quad 
All words below are available below the original sum P. L projects the head x from any original xV, with V nonempty; original literal heads need no step. W1 sorts any selected tail letter t first. Applying L with \(X=xt\) and \(Y=the\) remaining nonempty tail projects the edge word xt from an original xT containing t. If no tail remains, the edge is already original. No empty word substitution occurs.

For every active head x in H intersect A, derive xx as follows. If x occurs in an original tail, project yx from that original word, and project the literal x from an original x-headed word. F gives \(xx\leq yx+x\). Otherwise x has an original outgoing tail y that is itself in H. Project xy and literal y; E gives \(xx\leq xy+y\). This exhausts H intersect A. Squares already original need no step. Activation uses original edges and original heads, not unproved additional graph reachability.

If the target is literal h, condition(1) and L already derive it. For a nonlinear target with h in A, derive hh. For every required tail t: if t belongs to D, project yt from an original tail occurrence and use S to derive \(ht\leq hh+yt\). Otherwise t in A is an active original head; derive tt by the previous paragraph and use \(ht\leq hh+tt\). Thus every required elementary tail word ht is available.

If h is outside A, condition(3) says that each required t is already in an original h-headed tail. Project each ht directly; no square activation is assumed. In either case combine ht words one at a time with M:

{\small\par\smallskip
\CaseLine{\(h U t \leq  h U + h t,\)}
\smallskip}

where U is the accumulated nonempty tail word. Normalize with W1, W2 to retain exactly the desired tail set. This terminates at the target. All substitutions are nonempty and every derived right term is already available. Hence(1)-(3) suffice, and every valid absorption is derivable.

Distributivity writes arbitrary terms as nonempty polynomials. Applying the absorption derivation to each summand in both directions proves every valid identity. The seven laws form a finite relative basis for 4025. Reverse all words for the opposite. Catalogue equivalences need independent complete HSP in both directions; finite checks supplement this proof.

\subsection{A two-witness finite basis from the head and tail content}\label{case:113}
Source 4231 in the five-element catalogue. We use the signature (+, multiplication), without constants, relative to the additively idempotent semiring axioms.

\par\smallskip\noindent\textit{Algebra and normalization.}\quad 
Let A=\{0,\allowbreak{}1,\allowbreak{}2,\allowbreak{}3,\allowbreak{}4\}, with the tables in that order:

{\small\par\smallskip
\noindent addition       multiplication\par
\CaseLine{\ttfamily 0 0 0 4 4\qquad 0 0 2 2 0}
\CaseLine{\ttfamily 0 1 1 4 4\qquad 0 0 2 2 0}
\CaseLine{\ttfamily 0 1 2 3 4\qquad 2 2 2 2 2}
\CaseLine{\ttfamily 4 4 3 3 4\qquad 3 3 2 2 3}
\CaseLine{\ttfamily 4 4 4 4 4\qquad 4 4 2 2 4}
\smallskip}

Addition is join in the order \(2<1<0<4\) and \(2<3<4\), with 3 incomparable to 0 and 1. Put K=\{0,\allowbreak{}1,\allowbreak{}4\}. In multiplication, a right factor outside K makes the product 2. A right factor in K maps a left factor 1 to 0 and leaves each of 0,\allowbreak{}2,\allowbreak{}3,\allowbreak{}4 unchanged. This description proves associativity: a middle or last factor outside K makes both associations 2; if both are in K, the idempotent map just described is applied to the first factor. Every multiplication row and column preserves the two chains and the incomparable joins \(3+0=3+1=4\). Thus both distributive laws hold. The labels name elements only, not constants in terms.

For any nonlinear word w the word-value formula is

\CaseDisplay{w(\alpha)=
 \begin{cases}
  2,&\alpha(h(w))=2\text{ or }\alpha(r(w))\not\subseteq K,\\
  3,&\alpha(h(w))=3\text{ and }\alpha(r(w))\subseteq K,\\
  4,&\alpha(h(w))=4\text{ and }\alpha(r(w))\subseteq K,\\
  0,&\alpha(h(w))\in\{0,\allowbreak{}1\}\text{ and }\alpha(r(w))\subseteq K.
 \end{cases}                                                   \quad\textup{(\ensuremath{1})}}

Here h(w) is the head, r(w) is the set of variables in all positions after the first, and c(w) is the full content. Induction proves (1). It immediately gives the identities

\CaseDisplay{xyz\approx xzy,\qquad xy^2\approx xy.      \quad\textup{(\ensuremath{W})}}

Fix a variable order. These laws permute all positions after the first and remove repeated tail occurrences. Hence every nonlinear word has normal form h (product of T), where \(T=r(w)\) is nonempty and listed in increasing order. If h belongs to T its tail occurrence is retained. Linear words stay linear. All substitutions are nonempty.

\par\smallskip\noindent\textit{Exact criterion and two witnesses.}\quad 
Let P be a nonempty polynomial, L its set of linear summands, and \(N_{C}\) the nonlinear summands whose contents are contained in C.

Lemma 1. A linear variable y is absorbed by P exactly when y is a literal linear summand. For nonlinear w put \(h=h(w)\), \(T=r(w)\), \(C=c(w)\). Then A satisfies \(w\leq P\) precisely when \(N_{C}\) is nonempty and one of the following holds:

\textbullet{} h belongs to L;

\textbullet{} h belongs to T and some v in \(N_{C}\) has \(h(v)=h\);

\textbullet{} h does not belong to T and some v in \(N_{C}\) has \(h(v)=h\) and r(v) contained in T.

Every valid nonlinear absorption has a witness consisting of at most two original summands of P.

Proof. If a linear target y has no literal summand, assign \(y=3\) and all other variables 2. Every nonlinear word then has value 2, as do all linear summands. The target is 3, proving necessity. The converse is immediate.

For a nonlinear target, if \(N_{C}\) is empty, assign C the value 1 and all outside variables 2. The target is 0. Every nonlinear summand contains an outside variable and has value 2, while every linear summand is 1 or 2; hence \(P\leq 1\) and the absorption fails.

If neither h in L nor a same-head member of \(N_{C}\) exists, assign \(h=4\), the other variables of C the value 0, and outside variables 2. The target is 4. Outside-containing summands are 2, and all remaining summands have head different from h, so \(P\leq 0\). This proves the same-head requirement.

Finally suppose h is absent from T and L, and no same-head member of \(N_{C}\) has tail contained in T. Assign \(h=3\), every variable of T the value 1, and outside variables 2. The target is 3. Every nonlinear right summand with head h either contains an outside variable or has h in its tail; in either case its value is 2. Other summands have value in \{0,\allowbreak{}1,\allowbreak{}2\}, and \(P\leq 0\). Since 3 is not below 0, the third requirement is necessary.

For sufficiency take an arbitrary assignment. A target value 2 is least. For target value 3, \(h=3\) and all T lie in K, so h is not in T. The literal h or the specified same-head tail-contained witness has value 3. For target value 4, the literal h or any specified same-head contained witness has value 4, since all of C lies in K. For target value 0, all of C lies in K. Any member of \(N_{C}\) therefore has value 0 or 4, either of which absorbs the target. Formula (1) exhausts the possible target values.

If h belongs to L, retain the literal h and one member of \(N_{C}\). Otherwise retain the single specified same-head nonlinear witness. The same criterion holds for this one- or two-summand polynomial, proving the final assertion. \ensuremath{\square}

\par\smallskip\noindent\textit{Compression and finite basis.}\quad 
Consider a valid normal absorption with at most two right words, so there are at most three words in total. Reserve each variable that is the head of any word, including every linear word's unique variable. There are at most three reserved variables. Every other variable appears only in tails. Partition them by their nonempty incidence pattern in those at most three tail sets. There are at most \(2^{3}-1=7\) classes. Identifying each class to a new variable and normalizing by (W) preserves validity. This gives at most 10 variables and normal word length at most 11.

The identification is reversible modulo (W). Substitute the nonempty product of each class's original members for its class variable; substitute each individual original variable for its reserved head. The first positions are unchanged. Since every class member has the same tail incidence pattern, expansion restores every original tail set. Laws (W) sort and remove repetitions in the tails, giving the original normal words. Linear words are also restored.

Fix a 10-letter alphabet X. Let Sigma consist of (W) and all valid absorptions over X with one normal left word, at most two normal right words, and word length at most 11. This is a finite, explicitly specified set; Lemma 1 decides membership. It has not been enumerated as a small list, and no minimality claim is made.

Theorem. Sigma is a finite identity basis for A. Word reversal gives a finite identity basis for the multiplicative opposite.

Proof. Soundness follows from (W) and the definition. For any valid word absorption, Lemma 1 selects at most two original right summands (or a literal summand for a linear target). Compression shows that this smaller absorption follows from Sigma by substitution and (W). Adding the other summands derives the full absorption. Distributivity reduces arbitrary terms to polynomials; a valid identity between polynomials is obtained by deriving the absorption of each summand on each side into the other. This proves completeness without restricting the original variable count or word lengths. Word reversal transports validity and derivations to the opposite. \ensuremath{\square}

Finite regressions only support error detection. Transfers to other catalogue types require separately verified equality of generated varieties in both HSP directions. No hereditary finite-basis assertion is made.

\subsection{首字母与内容分离：布尔核心回缩完整基}\label{case:114}
\par\smallskip\noindent\textit{1. 生成元.}\quad 
令 \(D=\{z, b, a, d, t\}\)，加法取偏序

\CaseDisplay{z<b<a<t,\qquad z<d<t,\qquad d\parallel b, a}

上的并。其四元素子半格 \(B=\{z, a, d, t\}\) 识别为布尔格 \(\{0,\allowbreak{}1\}^2\)：

\CaseDisplay{z=(0,\allowbreak{}0),\quad a=(0,\allowbreak{}1),\quad d=(1,\allowbreak{}0),\quad t=(1,\allowbreak{}1).}

在 \(B\) 上定义乘法

\CaseDisplay{(i, j)(k, l)=(i, j\wedge l).}

令 \(\rho(b)=a\)，其余点不变，并对 \(D\) 定义 \(xy=\rho(x)\rho(y)\)。\(\rho\) 保持加法并固定 \(B\)；\(B\) 的乘法结合且左右分配，故这一扩张也是加法幂等半环。事实上，结合律两边均化为 \(\rho(x)\rho(y)\rho(z)\)，分配律由 \(\rho\) 的加性直接继承。

按 \((z, b, a, d, t)\) 排列，乘法表为

\CaseDisplay{\begin{array}{c|ccccc}
\cdot&z&b&a&d&t\\\hline
z&z&z&z&z&z\\
b&z&a&a&z&a\\
a&z&a&a&z&a\\
d&d&d&d&d&d\\
t&d&t&t&d&t
\end{array}.}

在目录中它对应编号 4241，标号映射为 \((z, b, a, d, t)\mapsto(2,\allowbreak{}1,\allowbreak{}0,\allowbreak{}3,\allowbreak{}4)\)。设 \(\rho(x)=(\alpha(x),\beta(x))\)，则对每个非线性词 \(w\)，

\CaseDisplay{w=\left(\alpha(h(w)),\bigwedge_{x\in c(w)}\beta(x)\right).\quad\textup{(\ensuremath{1})}}

注意一次词可以取新元素 \(b\)，不能用 (1) 取代其原值。

\par\smallskip\noindent\textit{2. 任意变量数的判据.}\quad 
对非空多项式 \(F\)，记一次项变量集为 \(L\)，非线性词族为 \(\mathcal R\)；记

\CaseDisplay{H=L\cup\{h(u):u\in\mathcal R\},\qquad
\mathcal R_C=\{u\in\mathcal R:c(u)\subseteq C\}.}

引理。

1. \(x\le F\) 当且仅当 \(x\in L\)，或存在 \(u\in\mathcal R\) 仅含变量 \(x\)。

2. 对非线性 \(w\)，置 \(h=h(w),C=c(w)\)，则

\CaseDisplay{w\le F\quad\Longleftrightarrow\quad h\in H\ \text{且} \mathcal R_C\ne\varnothing.\quad\textup{(\ensuremath{2})}}

证明。 一次词的必要性用赋值 \(x=b\)、其他变量为 \(z\)：若没有同名一次项和仅含 \(x\) 的非线性项，所有右项均为 \(z\)。充分性由 \(x\le x^2\) 及 \(x^n=x^2\)（\(n\ge2\)）给出。

若 \(h\notin H\)，赋 \(h=d\)、其他变量为 \(z\)，则左值为 \(d\)、右值为 \(z\)。若 \(\mathcal R_C\) 为空，赋 \(C\) 内变量为 \(b\)、外部变量为 \(z\)。左值为 \(a\)；右非线性词含外部变量，第一坐标为零且第二坐标也为零，故取 \(z\)；右一次项至多为 \(b<a\)。这证明 (2) 的必要性。

充分性按左值讨论。左值为 \(z\) 时自明。若第二坐标为一，则任意内部非线性右词也有第二坐标一，所以它的值在 \(a\) 之上。若第一坐标为一，则首字母赋值属于 \(\{d, t\}\)；条件 \(h\in H\) 提供同首字母右非线性词或同名一次项，其值在 \(d\) 之上。左值为 \(t\) 时两种见证的并在 \(a+d=t\) 之上。因此全部左值都被吸收。证毕。

\par\smallskip\noindent\textit{3. 六条完整有限基.}\quad 
定理。 相对于通常加法幂等半环公理，\(D\) 的完整基为

\CaseDisplay{\begin{aligned}
xyz&\approx xzy, &(N)\\
xy&\approx x^2y^2, &(S)\\
x&\le x^2, &(I)\\
xy&\le x^2, &(C)\\
xy&\le x+yx, &(H)\\
xy&\le y^2+xyz. &(E_2)
\end{aligned}}

这六式的有效性由 (1)、\(b<a\) 和判据 (2) 直接给出。反转每个乘法词给出乘法对偶的完整基。

完备性证明。 由 (N)、(S)，每个非线性词化为首字母和内容决定的正规形 \(W(h, C)\)，且幂等；推导方式为 \(x^2=x^4\)、\(x^3=x^6=x^2\)，再将首位后的字母交换并平方。

设 \(w=W(h, C)\le F\) 有效。选 \(q\in\mathcal R_C\)。若 \(h\in L\)，在 (H) 中以 \(x=h, y=q\) 代入，并用 (C)，得到

\CaseDisplay{hq\le h+qh\le h+q\le F.}

其内容是 \(\{h\}\cup c(q)\subseteq C\)；由 (C) 补乘缺少的字母，得到 \(w\le F\)。

若 \(h\notin L\)，判据给出以 \(h\) 开始的右非线性词 \(r\)。记 \(E=C\cup c(r)\)。由 (C) 向 \(r\) 补乘缺少的 \(C\) 内字母，得到

\CaseDisplay{W(h, E)\le r\le F.}

若 \(E=C\) 则完成；否则，令 \(s\) 为 \(E\setminus C\) 各字母组成的非空词。在 \((E_2)\) 中代入 \(x=w, y=q, z=s\)。因为 \(q\) 非线性且内容包含于 \(C\)，有 \(q^2=q\)、\(wq=w\)，故

\CaseDisplay{w=wq\le q^2+wqs=q+W(h, E)\le F.}

一次词吸收由 (I) 和同名右项处理。一般恒等式再拆成双向逐词吸收，完成证明。所有补乘均在缺少字母时使用非空词，否则省略，不使用空积。

\par\smallskip\noindent\textit{4. 范围.}\quad 
本证明适用于所定义生成元及其对偶，并通过精确双向簇包含传递给等簇的目录对象。它没有证明一般幂等核心的膨胀保留有限基性。有限规模程序核验仅作独立诊断；无限制完备性由第 3 节的形式推导负责。

\subsection{Two witnesses from a nonlinear head and a contained final variable}\label{case:115}
Source 4345 in the five-element catalogue. The signature is (+, multiplication), without constants. Bases are relative to the additively idempotent semiring axioms.

\par\smallskip\noindent\textit{The algebra and its word normal forms.}\quad 
Let A=\{0,\allowbreak{}1,\allowbreak{}2,\allowbreak{}3,\allowbreak{}4\}, with the following tables in that order.

{\small\par\smallskip
\noindent addition       multiplication\par
\CaseLine{\ttfamily 0 0 0 0 4\qquad 0 0 2 2 4}
\CaseLine{\ttfamily 0 1 0 1 4\qquad 0 0 2 2 4}
\CaseLine{\ttfamily 0 0 2 2 4\qquad 2 2 2 2 2}
\CaseLine{\ttfamily 0 1 2 3 4\qquad 3 3 3 3 3}
\CaseLine{\ttfamily 4 4 4 4 4\qquad 0 0 2 2 4}
\smallskip}

Addition is join for the order \(3<1<0<4\) and \(3<2<0<4\), where 1 and 2 are incomparable. Set K=\{0,\allowbreak{}1,\allowbreak{}4\}. A left factor 2 or 3 is unchanged by multiplication. A left factor in K gives 2 when the right factor is outside K, and otherwise gives 4 if that factor is 4 and 0 if it is 0 or 1. This proves associativity: for a first factor outside K both products are that factor; for a first factor in K, either of the remaining two factors outside K gives 2 in both associations; when all three factors lie in K, only the final factor matters. Every multiplication row and column is order preserving and preserves the only incomparable join \(1+2=0\). Thus both distributive laws hold. Element labels are not term constants.

For a word w let h(w), f(w), c(w), and m(w) denote its first variable, final variable, full content, and the set of variables strictly between the first and final positions. A nonlinear word here means any word of length at least two, including a word with no repeated variables. For such a word, induction gives

\CaseDisplay{w(\alpha)=\begin{cases}
3,&\alpha(h(w))=3,\\
2,&\alpha(h(w))=2,\\
2,&\alpha(h(w))\in K\text{ and a later position lies outside }K,\\
4,&\alpha(c(w))\subseteq K\text{ and }\alpha(f(w))=4,\\
0,&\alpha(c(w))\subseteq K\text{ and }\alpha(f(w))\in\{0,\allowbreak{}1\}.
\end{cases}                                                     \quad\textup{(\ensuremath{1})}}

Consequently the following word identities hold:

\CaseDisplay{xyzt\approx xzyt,\qquad xy^2z\approx xyz.       \quad\textup{(\ensuremath{W})}}

Fix a variable order. The first law permits adjacent interior transpositions (the outer substitutions are always nonempty); the second removes adjacent repeated interior letters. Thus each word of length at least two is provably equal to

\CaseDisplay{h(w)\,\prod_{x\in m(w)}x\,f(w).              \quad\textup{(\ensuremath{2})}}

An empty product means an omitted substring. Endpoint occurrences inside m(w) are retained. Linear words remain unchanged.

\par\smallskip\noindent\textit{An exact absorption criterion.}\quad 
Let P be a nonempty sum of nonempty words. Write \(w\leq P\) for w+P\ensuremath{\approx}P.

Lemma. A satisfies \(w\leq P\) if and only if both applicable conditions below hold:

1. P has a summand v with \(f(v)=f(w)\) and c(v) contained in c(w).

2. If w has length at least two, P has a summand u of length at least two with \(h(u)=h(w)\).

In particular every valid absorption has at most two original right summands as witnesses. A linear target needs only the first witness.

Proof. Put \(C=c(w)\) and \(t=f(w)\). If the first condition fails, assign \(t=4\), every other variable in C the value 0, and every outside variable the value 3. The target has value 4. A nonlinear summand containing an outside variable has value 3 if its head is outside C and value 2 otherwise. A nonlinear summand entirely inside C ends at a variable other than t and has value 0. Linear summands are 0 or 3, since a literal t would supply the required witness. All right summands are at most 0, so the absorption fails. This also proves necessity for linear w.

For nonlinear w, suppose the second condition fails. Assign \(h(w)=1\) and every other variable the value 3. All nonlinear right summands have value 3. Linear right summands have value 1 or 3, so \(P\leq 1\). The target has value 0 if it is a power of h(w), and value 2 otherwise. Neither 0 nor 2 is below 1. This proves necessity.

For sufficiency first let w be a variable y. The first witness is some positive power \(y^{n}\). For every a in A one has \(a\leq a^{n}\): \(n=1\) is trivial, and for \(n\geq 2\) the only change under taking powers is 1 to 0. Hence \(y\leq v\leq P\).

Now let w be nonlinear and select witnesses u and v as above. For an arbitrary assignment use (1). If \(w=3\), absorption is automatic because 3 is the least element. If \(w=4\), all variables in C lie in K and \(t=4\); the contained final witness v therefore has value 4, including when it is linear. If \(w=2\), its head has value 2 or lies in K. The nonlinear same-head witness u consequently has value 2, 0 or 4, each of which is at least 2. Finally, if \(w=0\), all C lie in K and t is 0 or 1. A nonlinear v has value 0, immediately absorbing w. A linear v is t, with value 0 or 1; u has value 0, 2 or 4. The only case in which neither individually dominates 0 is \(v=1\) and \(u=2\), and their join is 0. Thus \(w\leq u+v\) in every case. This proves sufficiency and the two-witness assertion. \(\square\)

\par\smallskip\noindent\textit{Reversible compression and an unrestricted finite basis.}\quad 
Consider a valid absorption between normal words with at most two right summands. There are at most three words in total. Reserve all variables occurring as an endpoint of any of these words, including the unique variable of a linear word. There are at most six reserved variables. Every other variable occurs only in interiors. Partition these other variables according to their nonempty incidence pattern in the at most three interior sets. There are at most \(2^{3}-1=7\) classes.

Identify the variables in each class to one new variable, keeping each reserved variable separate, and normalize by (W). Identification is a substitution, so it preserves validity. The resulting absorption has at most 13 variables, normal word length at most 15, and at most two right summands.

This identification is reversible modulo (W). Substitute for each class variable the nonempty product of its original members, and for each reserved variable its original single letter. Endpoints remain unchanged. All members of a class have the same interior incidence pattern, so expansion recovers exactly every original interior set. Normalization by (W) recovers the original words. This argument also covers repeated endpoint letters in interiors and linear words.

Fix a 13-letter alphabet X. Let Sigma consist of (W) and all valid normal absorptions over X with one left word, one or two right words, and word length at most 15. This is a finite, explicitly specified family. The lemma decides its membership by simple set operations; enumeration over assignments is not needed. We do not claim that this family is a small or minimal basis.

Theorem. Sigma is a finite identity basis for A. Reversing every word gives a finite identity basis for its multiplicative opposite.

Proof. Soundness is immediate. For an arbitrary valid absorption, normalize its words using (W). The lemma selects at most two original right summands. Apply reversible compression to this valid smaller absorption. The compressed absorption belongs to Sigma up to renaming, and substituting the original class products back derives the smaller absorption. Adding all remaining summands derives the original one. For arbitrary identities, distributivity first gives polynomials; absorb each left summand into the right polynomial and conversely. The two inequalities yield the identity. No bound is imposed on the variable count or word lengths of the identity being derived. Reversal transports every identity and derivation to the opposite semiring. \(\square\)

\subsection{A finite basis from head, nonlinear-content, and final witnesses}\label{case:116}
Source 4346 in the five-element catalogue. The signature is (+, multiplication), without constants, relative to the additively idempotent semiring axioms. Numerical labels name elements only.

\par\smallskip\noindent\textit{1. Algebra and word values.}\quad 
Tables in order 0,\allowbreak{}1,\allowbreak{}2,\allowbreak{}3,\allowbreak{}4:

{\small\par\smallskip
\noindent addition       multiplication\par
\CaseLine{\ttfamily 0 0 0 0 4\qquad 0 0 2 2 4}
\CaseLine{\ttfamily 0 1 1 1 4\qquad 0 0 2 2 4}
\CaseLine{\ttfamily 0 1 2 2 4\qquad 2 2 2 2 2}
\CaseLine{\ttfamily 0 1 2 3 4\qquad 3 3 3 3 3}
\CaseLine{\ttfamily 4 4 4 4 4\qquad 0 0 2 2 4}
\smallskip}

The additive chain is \(3<2<1<0<4\). Put B=\{0,\allowbreak{}1,\allowbreak{}4\}, K=\{2,\allowbreak{}3\}, and \(f(0)=f(1)=0\), \(f(4)=4\). Elements 2 and 3 are left zeros. For u in B, \(uv=2\) if v is in K, and \(uv=f(v)\) if v is in B. Associativity follows by cases: a first factor in K remains that factor; with the first factor in B, a middle factor in K makes both associations 2; with the first two factors in B, a last factor in K again gives 2; with all factors in B both associations give f(the last factor), since f is idempotent. Rows and columns are increasing in the additive chain, so both distributive laws hold.

For every nonlinear word, induction gives value 3 when its head is 3; value 2 when its head is 2, or its head is in B and some tail factor lies in K; otherwise all factors are in B and the value is 0 or 4 according as the final factor is 0/\allowbreak{}1 or 4. Nonlinear words never have value 1.

\par\smallskip\noindent\textit{2. Normalization.}\quad 
The word formula proves

\CaseDisplay{xyzt\approx xzyt,\qquad x^2y\approx xy,\qquad xy^2\approx xy.\quad\textup{(\ensuremath{W})}}

Fix a variable order. Every nonlinear word normalizes to h (product of M) t, where h, t are its original endpoints and M is its sorted interior set with h, t removed. Endpoints may coincide. Empty M denotes the two-letter word ht, not a unit. Linear words remain linear; in particular xx is not replaced by x.

The first law interchanges adjacent interior positions using nonempty outer prefixes and suffixes. After sorting, \(xy^{2}=xy\) removes adjacent interior repetitions, with a nonempty preceding prefix substituted for x and the remaining suffix as right context. Move an interior copy of h next to the head and remove it with \(x^{2}y=xy\), whose remaining suffix y is nonempty. Move an interior copy of t next to the final t and remove it with \(xy^{2}=xy\), whose prefix x is nonempty. Sort the remaining interior. These are actual nonempty substitutions and contexts.

\par\smallskip\noindent\textit{3. Exact criterion.}\quad 
For a nonempty polynomial P, let c(u), h(u), t(u) denote word content and endpoints. Normalization preserves these data and the linear/nonlinear distinction.

Lemma. A linear y is absorbed by P exactly when some original right word has content \{y\}. For a nonlinear target w, with \(C=c(w)\), \(h=h(w)\), \(t=t(w)\), absorption holds exactly when:

\textbullet{} H: some original right word has head h;

\textbullet{} N: some original nonlinear right word has content contained in C;

\textbullet{} T: some original right word has content contained in C and final variable t; this witness may be linear.

The target content need not be covered by right content, and the H witness need not be contained in C.

Necessity. If H fails, assign \(h=2\) and all other variables 3. Then \(w=2\) and \(P=3\). If N fails, assign C the value 1 and outsiders 3. Then \(w=0\); every nonlinear right word contains an outside variable and is 2 or 3, while linears are 1 or 3, so \(P\leq 1\). If T fails, assign \(t=4\), the other variables of C the value 1, and outsiders 3. Then \(w=4\). Contained nonlinear right words cannot end in t and have value 0. Noncontained nonlinear words have value 2 or 3. A linear t would supply T and is absent; all linears are 1 or 3. Hence \(P\leq 0\). For a linear y without a pure-y original, assign \(y=4\) and outsiders 3. Every right word is 2 or 3, separating it from y.

Sufficiency. A nonlinear target 3 is least. For target value 2, its head is 2 or in B; any same-head right word is at least 2, including linear witnesses, so H suffices. For target value 0, all target variables lie in B; the N witness is then 0 or 4. For target value 4, all target variables lie in B and \(t=4\); the T witness is 4, whether linear or nonlinear. A pure-y original absorbs linear y because \(y\leq y^{2}\) directly from the table and all higher powers have the same value as \(y^{2}\).

\par\smallskip\noindent\textit{4. Three witnesses and reversible compression.}\quad 
For a valid nonlinear absorption choose an original H witness a, N witness b, and T witness d. Their sum retains all three conditions, so \(w\leq a+b+d\) is valid with at most three distinct original right words. A linear target needs one pure-variable original. This is a direct reduction, with no bound on the number of original variables.

A local absorption has at most four normal words in total. Reserve individually every variable appearing at an endpoint of any word, including the sole variable of a linear word. There are at most eight such variables. All remaining variables occur only in interiors; group them by their nonempty incidence pattern among the at most four interiors. At most 15 groups occur.

Identify each group to one fresh variable, keep all reserved variables distinct, and normalize with (W). Ordinary substitution preserves validity. At most 23 variables occur, and every normal word has length at most 25. Conversely replace each reserved variable individually and each group variable by the nonempty product of its original members. No group occupies an endpoint. Identical interior-incidence patterns restore every original interior set, and (W) recovers each original word, including endpoints and linear words. Thus compression is reversible with no constants or empty substitutions.

\par\smallskip\noindent\textit{5. Finite complete basis.}\quad 
Fix a 23-letter alphabet. Let Sigma consist of (W) and all valid normal absorptions over this alphabet with word length at most 25 and at most three right words. This is a finite, explicitly specified family, whose membership is decided by the lemma. No small or minimal list is claimed.

Theorem. Sigma is a finite identity basis for source 4346, and its word reversals form a finite basis for the multiplicative opposite.

Proof. Soundness is built into the definition and (W). Reversible compression derives every valid local absorption with at most three right words from Sigma. Direct witness selection reduces every valid word absorption into an arbitrary polynomial to such a step; adding unused summands gives the original absorption. Distribute arbitrary terms into word sums and normalize. Every valid polynomial identity is equivalent to the absorption of each summand on either side by the other sum. All these absorptions are derived, giving completeness for arbitrarily many variables and original word lengths. Word reversal transports validity and derivations to the opposite.

Finite calculations check the algebra, value formula, criterion, normalization, witness selection, and inverse compression. They support the unrestricted proof above. Transfers to other catalogue entries require independently verified equality of generated varieties in both HSP directions.

\subsection{两个交换近幂等半环的完整有限基}\label{case:117}
\par\smallskip\noindent\textit{1. 两个半环.}\quad 
令 \(A=\{0, a, e, b, t\}\)，加法是偏序

\CaseDisplay{0<a<e<t,\qquad 0<a<b<t,\qquad e\parallel b}

的并。令 \(\mu:A\to\{0, a, e\}\) 固定 \(0, a, e\)，并把 \(b, t\) 送到 \(e\)。在 \(\{0, a, e\}\) 上，\(e\) 是乘法幺元、\(0\) 是乘法零、\(a^2=0\)；在 \(A\) 上定义 \(xy=\mu(x)\mu(y)\)。它对应目录编号 7037，按上述顺序到目录标号的映射为 \((2,\allowbreak{}1,\allowbreak{}0,\allowbreak{}3,\allowbreak{}4)\)。

令 \(B=\{0, a, e, f, t\}\)，加法取链 \(0<a<e<f<t\) 的最大值。在子集 \(\{0, a, e, f\}\) 上，\(e\) 是幺元、\(0\) 是零，且

\CaseDisplay{a^2=0,\qquad f^2=f,\qquad fa=af=a.}

把 \(t\) 的每个乘法作用定义为 \(f\) 的作用，包括 \(t^2=f\)。它对应目录编号 4567，标号映射为 \((3,\allowbreak{}2,\allowbreak{}1,\allowbreak{}0,\allowbreak{}4)\)。

两个定义均给出交换加法幂等半环。对 \(A\)，\(\mu\) 保持加法且固定其像，结合律和分配律降为三元素链上的核验。对 \(B\)，去掉 \(t\) 的乘法由上述交换幺半群关系给出；把 \(t\) 送到 \(f\) 的映射保持链上的并，且各乘法平移单调，故分配律成立。下文的 \(0, e\) 都只是元素名称，不是语言常量。

\par\smallskip\noindent\textit{2. 完整基.}\quad 
定理 1。 在通常加法幂等半环公理之外，\(A\) 的完整恒等式基为

\CaseDisplay{\begin{aligned}
xy&\approx yx, &(A1)\\
xy&\approx x^2y+xy^2, &(A2)\\
xyz&\leq xy, &(A3)\\
xy&\leq x+(xy)^2. &(A4)
\end{aligned}}

定理 2。 \(B\) 的完整恒等式基为

\CaseDisplay{\begin{aligned}
xy&\approx yx, &(B1)\\
xy&\approx x^2y+xy^2, &(B2)\\
x^2&\leq x, &(B3)\\
xy&\leq x+y^2, &(B4)\\
xy&\leq x+(xy)^2, &(B5)\\
xy^2z^2&\leq xy^2+x^2z^2. &(B6)
\end{aligned}}

条数不声称最小。

先核验有效性。两者的 (2) 式来自：一个非零乘积至多含一次 \(a\)；将两个因子中的一个平方时，至少一个求和项保持原值，另一个不大于原值。若均无 \(a,0\)，两边相等。\(A\) 的所有二因子积在 \(\{0, a, e\}\) 中，后续乘法只能使它在 \(0<a<e\) 中下降，故 (A3) 成立。(A4)、(B5) 中，乘积幂等时由平方项吸收；乘积为 \(a\) 时，每个因子均不小于 \(a\)，由一次项吸收；乘积为 \(0\) 时由最小元性吸收。\(B\) 的 (B3) 直接由链和平方表得到。(B4) 若左侧是 \(a\)，右侧的 \(x\) 已不小于 \(a\)；若左侧为 \(e\) 或 \(f\)，所需单位层值由 \(x\) 或 \(y^2\) 保留；左侧为零时自动成立。(B6) 中，若 \(y^2=0\) 或 \(z^2=0\)，左侧为零；否则 \(y^2, z^2\in\{e, f\}\)。此时若 \(x=a\)，第一项保留 \(a\)；若 \(x\in\{e, f, t\}\)，两个右项的最大值与左侧相等；\(x=0\) 时亦成立。

\par\smallskip\noindent\textit{3. 任意变量数的吸收判据.}\quad 
把项展开为非空词的非空和。对多项式 \(Q\)，令 \(L\) 为其一次项变量集，\(R\) 为其长度至少二的求和词族。对非空变量集 \(C\)，记

\CaseDisplay{R_C=\{u\in R:c(u)\subseteq C\}.}

令 \(m_x(u)\) 为 \(x\) 在词 \(u\) 中的次数，\(S(w)=\{x:m_x(w)=1\}\)。集合族为空时按通常约定处理，下述条件的存在量词不得省略。

引理 3。 对两个半环均有

\CaseDisplay{x\leq Q\quad\Longleftrightarrow\quad x\in L.\quad\textup{(\ensuremath{1})}}

对长度至少二的词 \(w\)，置 \(C=c(w)\)。在 \(A\) 中，\(w\leq Q\) 当且仅当

\CaseDisplay{R_C\ne\varnothing,\qquad
L\cap C\ne\varnothing\ \text{或}\
\bigl(\forall x\in S(w)\;\exists u\in R_C:\ m_x(u)\leq1\bigr).\quad\textup{(\ensuremath{2})}}

在 \(B\) 中，\(w\leq Q\) 当且仅当

\CaseDisplay{C\subseteq L\cup\bigcup_{u\in R_C}c(u),\qquad
L\cap C\ne\varnothing\ \text{或}\
\bigl(\forall x\in S(w)\;\exists u\in R_C:\ m_x(u)\leq1\bigr).\quad\textup{(\ensuremath{3})}}

必要性。 若左侧是一次项 \(x\notin L\)，令 \(x=t\)、其他变量为 \(0\)。所有右侧一次项为零，所有乘积严格低于 \(t\)，故 (1) 失败。

对 \(A\)，若 \(R_C\) 为空，令 \(C\) 内变量均为 \(b\)、外部变量均为 \(0\)。左词为 \(e\)，右侧至多为 \(b\)，而 \(e\nleq b\)。对 \(B\)，若覆盖条件失败，选未被覆盖的 \(x\in C\)，令 \(x=f\)、其他 \(C\) 内变量为 \(e\)、外部变量为 \(0\)。左词为 \(f\)，右侧至多为 \(e\)。

最后，若第二个条件失败，则存在 \(x\in S(w)\)，且所有 \(u\in R_C\) 都至少含两次 \(x\)，同时 \(L\cap C=\varnothing\)。令 \(x=a\)、其他 \(C\) 内变量为 \(e\)、外部变量为 \(0\)。左词为 \(a\)，右侧为 \(0\)。这也覆盖 \(R_C\) 为空的情形。

语义充分性。 \(A\) 的非线性词只取 \(0, a, e\)。取值为零时自动被吸收；取值为 \(e\) 时，\(R_C\) 中任一词也为 \(e\)。取值为 \(a\) 时，恰有一个一次出现的变量取 \(a\)，其他变量的 \(\mu\)-像为 \(e\)。若有 \(L\cap C\)，该一次项不小于 \(a\)；否则条件提供一个不将这个变量重复的右词，其值为 \(a\) 或 \(e\)。

\(B\) 的非线性词只取 \(0, a, e, f\)。取值为 \(a\) 的论证相同。取值为 \(e\) 时，覆盖条件提供非零的内部一次项或内部右词。取值为 \(f\) 时，选一个取 \(f\) 或 \(t\) 的左变量；覆盖该变量的内部右词或一次项至少为 \(f\)。因此 (3) 充分。证毕。

\par\smallskip\noindent\textit{4. 从四条基推出判据 (2).}\quad 
由 (A2) 令 \(y=x\)，得到 \(x^2\approx x^3\)。反复对两个不同的一次出现变量使用 (A2)，每个非线性词都化为非空和，其中各词具有以下两种形式：

\CaseDisplay{H(C)=\prod_{y\in C}y^2,
\qquad xH(C\setminus\{x\})\quad(x\in C,\ |C|\geq2).
\quad\textup{(\ensuremath{4})}}

后一形式中的 \(x\) 必须原本就在 \(S(w)\) 中。每一步减少一次出现的变量数，故过程终止，且保持内容集。这里不引入空积；单变量非线性词使用第一种形式。

由 (A3)，若非线性词 \(u\) 的每个变量次数都不超过某词 \(v\) 的次数，则 \(v\leq u\)：补乘缺少的因子，并用 (A3) 逐一删除这些附加因子。次数大于二时先用 \(x^3=x^2\) 约化。

若左词是 \(H(C)\)，选任一 \(u\in R_C\)，即得 \(H(C)\leq u\leq Q\)。若左词为 \(xH(C\setminus\{x\})\)，且条件提供 \(m_x(u)\leq1\) 的 \(u\in R_C\)，同样补乘即可。

只剩 \(L\cap C\ne\varnothing\)。选 \(y\in L\cap C\)，把左词写成 \(yv\)，其中 \(v\) 非空。由 (A4) 得

\CaseDisplay{w=yv\leq y+w^2\approx y+H(C)\leq Q.}

这里 \(H(C)\leq Q\) 已在上段推出。因此 (2) 的每个有效吸收式均由四条基推出。结合一次项判据，定理 1 得证。

\par\smallskip\noindent\textit{5. 从六条基推出判据 (3).}\quad 
(B2) 给出同一规范形 (4)。先把 \(Q\) 的非线性词也规范化：各新词保持旧内容；若一个变量在旧词中的次数不超过一，则至少有一个新词仍不将它重复。因此判据条件在所需方向保留。

由 (B3)、(B4) 有

\CaseDisplay{uv\leq u+v^2\leq u+v.\quad\textup{(\ensuremath{5})}}

对内部右词 \(u\) 和内部一次项 \(y\)，分别有 \(u^2\leq u\) 和 \(y^2\leq y\)。覆盖条件保证这些平方的内容并为 \(C\)。将它们相乘，并反复用 (5)，得到

\CaseDisplay{H(C)\leq Q.\quad\textup{(\ensuremath{6})}}

若 \(L\cap C\ne\varnothing\)，使用 (B5) 和 (6)，与上节完全相同地得到 \(w\leq Q\)。故只需处理 \(L\cap C=\varnothing\) 及左词 \(xH(C\setminus\{x\})\)。

选一个不将 \(x\) 重复的内部右词。若有内部右词恰含一次 \(x\)，它的规范形为 \(xH(S)\)，其中 \(S\ne\varnothing\)、\(x\notin S\)，于是已有 \(xH(S)\leq Q\)。若没有这样的右词，则保护 \(x\) 的右词 \(u\) 不含 \(x\)；覆盖条件另提供一个重复 \(x\) 的内部右词 \(v\)。令 \(u^2=H(S)\)，\(v^2=x^2H(T)\)。当 \(T\ne\varnothing\) 时，在 (B4) 中代入 \(H(S),xH(T)\)，得到

\CaseDisplay{xH(S\cup T)\leq H(S)+x^2H(T)\leq u+v\leq Q.}

当 \(T=\varnothing\) 时改代入 \(H(S),x\)，得到同一结论。这建立了初始受吸收的词 \(xH(S')\)，且 \(S'\ne\varnothing\)。

现在依次加入其余内部右词的内容，保持 \(x\) 恰出现一次。设已有 \(r=xH(S)\leq Q\)，取内部右词 \(v\)。

\textbullet{} 若 \(v\) 不含 \(x\)，则 \(v^2=H(T)\leq Q\)，由 (5) 得 \(xH(S\cup T)\leq r+v^2\leq Q\)。

\textbullet{} 若 \(v\) 恰含一次 \(x\)，写成 \(v=xH(T)\)。对稳定块使用 (5)，再乘 \(x\)，得 \(xH(S\cup T)\leq r+v\leq Q\)。

\textbullet{} 若 \(v\) 重复 \(x\)，写成 \(v^2=x^2H(T)\)。\(T\) 为空时没有新变量要加入；否则把 \(S, T\) 中变量各乘一次所得的非空词分别代入 (B6) 的 \(y, z\)，得 \(xH(S\cup T)\leq r+v^2\leq Q\)。

所有代换都使用非空词。由覆盖条件，最终内容恰为 \(C\)，所以得到所需的 \(xH(C\setminus\{x\})\leq Q\)。所有规范左词都已推出，故原吸收式也可推出。任意恒等式分解为逐词吸收式，定理 2 得证。

\par\smallskip\noindent\textit{6. 目录应用与核验边界.}\quad 
两组推导均适用于任意变量数。将具体运算表代入命题~\ref{prop:lower-fibres}，可逐一判定每个目标与来源之间两个方向的簇包含。

来源 \(7037\) 和 \(4567\) 分别对应 23 型与 22 型。两组编号互不相交，见附录~\ref{app:register}。

\subsection{Two finite bases from single-variable recovery}\label{case:118}
Sources 4521 and 4522. We use (+, multiplication), without constants, relative to the additively idempotent semiring axioms.

\par\smallskip\noindent\textit{1. Algebras and normal words.}\quad 
Both sources have multiplication, in order 0,\allowbreak{}1,\allowbreak{}2,\allowbreak{}3,\allowbreak{}4:

{\small\par\smallskip
\CaseLine{\ttfamily 0 0 2 3 0}
\CaseLine{\ttfamily 0 0 2 3 0}
\CaseLine{\ttfamily 2 2 3 3 2}
\CaseLine{\ttfamily 3 3 3 3 3}
\CaseLine{\ttfamily 4 4 2 3 4}
\smallskip}

For 4521 the additive order is \(3<1\),\(2<0<4\) with 1,\allowbreak{}2 incomparable; for 4522 it is the chain \(3<2<1<0<4\). Addition is join in each order. Element 3 is a multiplicative zero. A product with at least two occurrences of 2 is 3; a product with exactly one 2 and no 3 is 2. For a nonlinear product with all factors in K=\{0,\allowbreak{}1,\allowbreak{}4\}, its value is 4 when its head is 4 and 0 otherwise. These descriptions follow inductively from the table. They also prove associativity: the zero and count-of-2 cases are independent of bracketing, and in K the first factor alone determines every product of length at least two. The displayed rows and columns preserve the joins in both orders, so both distributive laws hold.

The word-value formula proves

\CaseDisplay{xyz\approx xzy,\qquad x^3\approx x^2.\quad\textup{(\ensuremath{W})}}

The first law permits arbitrary permutations after the head using nonempty prefix substitutions and ordinary right contexts. It places all copies of a letter together (next to the head if that is the letter), after which the second law reduces total multiplicity to two. Thus a normal word keeps its original head and each variable's total multiplicity capped at two, with the remaining tail sorted. A linear word stays linear; xx stays nonlinear. This is a derivation by (W), not just a semantic normal form.

\par\smallskip\noindent\textit{2. Exact absorption criteria.}\quad 
Fix target w, content C, and head h. Let U be the original right words of P with content contained in C. A variable is simple in w when it occurs exactly once.

Lemma. A linear x is absorbed exactly when x itself is an original right summand. For nonlinear w, absorption \(w\leq P\) holds exactly when:

\textbullet{} H: a word in U has head h;

\textbullet{} N: a word in U is nonlinear;

\textbullet{} for each simple x in w, \(A_{x}\): a word u in U has at most one x. For 4521 additionally require that u be nonlinear or the literal x. For 4522 there is no additional restriction.

Necessity. Without a literal right x, assign \(x=2\) and every other variable 3. A nonlinear pure power of x is 3; any word containing an outsider is 3. Thus \(P=3\) while \(x=2\).

If H fails, assign \(h=4\), other variables of C the value 0, and outsiders

3. The target is 4; right words containing an outsider are 3 and contained words have head different from h and value 0. Hence \(P\leq 0\). If N fails, assign all C the value 1 and outsiders 3. The nonlinear target is 0; nonlinear right words are 3 and linears are 1 or 3.

For failure of \(A_{x}\), assign \(x=2\), C minus x the value 1, and outsiders

3. The target is 2. In 4521 every contained nonlinear right word repeats x and is 3, and a literal x is absent; the remaining linears are 1 or 3. Their join is at most 1, which does not dominate 2. In 4522 every contained right word repeats x, so all right values are 3. These are separating assignments in all cases.

Sufficiency. Target value 3 is least. Target value 4 has all content values in K and \(h=4\); H supplies value 4. Target value 0 has all content values in K; N supplies value 0 or 4. Target value 2 has no 3 and has exactly one occurrence valued 2, at a simple variable x. A nonlinear \(A_{x}\) witness omitting x has value 0 or 4; one containing x once has value 2. A literal x also has value 2. For 4522 any other permitted literal is in K and is above 2 in its chain. This proves both criteria for arbitrary words and polynomials.

\par\smallskip\noindent\textit{3. Two-right-word derivations.}\quad 
Normalize all words. For C, h and a set J of variables to be simple, write \(G_{J}\) for the normal word with head h, content C, multiplicity one on J and two on C minus J. Begin with J empty.

Choose original witnesses H and N from the lemma. The fully doubled word \(G_{\varnothing}\) is absorbed by H+N: both contained witnesses remain, and this target has no simple-variable conditions. Thus a first valid local step with at most two right words constructs \(G_{\varnothing}\). No coverage of C by original right words is required.

For each simple target variable x choose its original \(A_{x}\) witness \(u_{x}\). The valid local absorption \(G_{\{x\}}\leq G_{\varnothing}+u_x\) restores x alone. The guard supplies H and N, and \(u_{x}\) supplies \(A_{x}\). All these targets are nonlinear: the only exceptional possibility would be singleton C with simple h, which is a linear target already treated.

If J is a set of already restored target variables and x is another, then \(G_{J\cup\{x\}}\leq G_J+G_{\{x\}}\) is valid. Both right words are nonlinear with the full content and correct head. Each variable in J uses \(G_{J}\) as its A witness, and \(x\) uses \(G_{\{x\}}\). Iterate over the finitely many simple target variables. The final word is w. Consequently every valid nonlinear absorption has a derivation using valid local absorptions with at most two right words at each step.

\par\smallskip\noindent\textit{4. Reversible common-multiplicity compression.}\quad 
A local step involves at most three normal words. Reserve every head variable separately, at most three. Group all other variables by their common multiplicity vectors in \{0,\allowbreak{}1,\allowbreak{}2\}\ensuremath{^{k}}, where \(k\leq 3\) is the number of words. The zero vector does not occur, so there are at most 26 further groups.

In the compressed word put each group variable at its common multiplicity, not the sum of its members' multiplicities. Preserve every reserved head variable and its multiplicity. There are at most 29 variables and length at most 58. This operation is not ordinary variable identification; validity follows from the exact criteria as follows. Content inclusion is preserved in both directions because all members of a group have the same zero coordinates. Head equalities are preserved because heads are individually reserved. Linearity/nonlinearity is preserved: a nonlinear word either repeats its head or contains another variable in a different group. Every simple variable retains its count in every word, so the A condition is preserved. A literal is a reserved head and remains that literal; therefore the extra 4521 restriction is preserved as well. The same argument covers linear targets and literal summands.

For the inverse substitution replace a group variable by a nonempty word listing its original members once. It never occupies a head. Its one or two copies restore exactly the common multiplicities of all its members. Tail permutation and normalization by (W) recover each original normal word. This proves reversibility without empty substitutions or constants.

\par\smallskip\noindent\textit{5. Finite complete bases.}\quad 
For each of the two sources fix a 29-letter alphabet and take (W) together with every valid normal absorption on this alphabet of word length at most 58 with at most two right words. Validity of these finitely many identities is specified by its own exact criterion. No small or minimal list is claimed.

Theorem. These two finite families are complete identity bases for their respective algebras. Their word reversals are bases of the respective multiplicative opposites.

Proof. Every local step in Section 3 compresses to a member of the corresponding finite family, then lifts back by the inverse substitution and (W). Hence every valid word absorption is derivable with arbitrarily many variables, summands and original word lengths. Distributivity reduces arbitrary identities to sums of words, and each summand on either side is absorbed by the other side. Deriving these absorptions proves completeness. Word reversal transports validity and derivations to each opposite. Finite regressions support this unrestricted proof; transfer to other catalogue types requires independently checked equality in both HSP directions.

\subsection{Six finite bases from witnesses for a single occurrence}\label{case:119}
The catalogue sources 4524, 5314, 7498 have opposites 5256, 7242, 7556, respectively. All bases are relative to the constant-free ai-semiring axioms. Write C for the support of a target w and U for the original upper words supported in C. A simple letter occurs once. Use pre and post for the support strictly before or after that unique occurrence. A literal has length one.

\par\smallskip\noindent\textit{Three exact absorption tests.}\quad 
In 4524 and 5314 the additive order is the diamond with bottom 3, incomparable elements 1 and 2, their join 0, and greatest element 4. In 7498 it is the chain \(3<2<1<0<4\). All three products give 3 if a 3 occurs or if at least two entries are 2. Set K=\{0,\allowbreak{}1,\allowbreak{}4\}.

\textbullet{} The rows of 4524 are 00230,\allowbreak{}00230,\allowbreak{}33332,\allowbreak{}33333,\allowbreak{}00234. On K a literal retains its value. A longer word gives 4 when all entries are 4, otherwise 0. With a unique 2 and other entries in K the result is 2 exactly when every later entry is 4, otherwise 3.

\textbullet{} The rows of 5314 are 00330,\allowbreak{}11331,\allowbreak{}22332,\allowbreak{}33333,\allowbreak{}00234. On K a word of all 4's gives 4. Otherwise the value is 1 when its head is 1 and 0 in the remaining cases. A unique 2 gives 2 exactly when all earlier entries are 4.

\textbullet{} The rows of 7498 are 01234,\allowbreak{}11231,\allowbreak{}22332,\allowbreak{}33333,\allowbreak{}44234. On K the value is the first nonzero entry, or 0 if all entries are 0. A unique 2 gives 2 irrespective of its position.

Each formula follows by multiplying one further input. The following conditions characterize w \ensuremath{\leq} P on unrestricted alphabets.

For 4524, U must be nonempty and, when w is nonlinear, contain a nonlinear word. For every simple x of w at least one of the following witnesses must belong to U:

1. A nonlinear word omitting x.

2. A literal from the suffix after x.

3. A word with simple x and \(post_{u}\)(x) contained in \(post_{w}\)(x).

For 5314, U must contain a same-head word. For every simple x there must be a word of U either omitting x and headed in \(pre_{w}\)(x), or containing a simple x with \(pre_{u}\)(x) contained in \(pre_{w}\)(x).

For 7498, at every first occurrence of x in w, some u in U contains x and has its prefix before its first x supported in the prefix before the first x of w. For every simple x in w some word of U either omits x or contains it exactly once.

To prove necessity give every outsider 3. In 4524 assign target variables 4 for nonemptiness and 1 for the nonlinear requirement. For a simple x assign \(x=2\), its suffix support 4 and all other target letters 1. The only contained words reaching 2 or an upper bound of 2 are the three specified kinds. All other right words have values in \{1,\allowbreak{}3\}, whose sum cannot absorb 2. In 5314 use 0 on the target head and 1 on the other target letters to force the same head. For a simple x put \(x=2\), \(pre_{w}(x)=4\) and the other target letters 1. An omitting word then dominates 2 only when its head has value 4. In 7498 put the prefix before the first x at 0, x at 4 and all other target letters at 1. A right word reaching 4 must be the indicated prefix witness. For a simple x put \(x=2\) and all other target letters 0 to obtain an omission or a simple witness.

For sufficiency use the word formulas. In 4524 a nonlinear contained word absorbs the target when its values are in K, and a contained power suffices for a literal on K. With a unique 2 giving 2 an omitting nonlinear witness gives 0 or 4, a suffix literal gives 4, and the contained post witness gives 2. In 5314 the contained same- head word suffices on K, including the case of head 4 and a later non-4. An omitting witness headed in the prefix preceding a surviving 2 has head value 4 and gives 0 or 4. The other simple witness gives

2. In 7498 the first 4 giving the target value 4 has a prefix witness with preceding values zero. An all-zero target uses any contained word. At value 1 any contained word on K suffices, and at value 2 the simple or omitting witness suffices. Value 3 is least.

\par\smallskip\noindent\textit{4524 and 5314.}\quad 
Use \(xy=xxy+xyy\) together with all valid absorptions on a fixed 737-letter alphabet, of maximum word length 2951, with at most two upper words. All 25 table assignments verify the splitting law in each algebra. It expands every word into words with at most one simple letter by splitting between two distinct simple occurrences. In each of the two branches one or more simple letters become repeated and every repeated letter stays repeated.

For a valid step after expanding both sides retain the structural witness (a contained word of the required length for 4524, a contained same-head word for 5314), plus the witness for the single simple letter if needed. The same exact test proves this absorption with at most two upper words. Apply the head-oriented run compression of the 7303/7645 proof: protect the head, the unique simple position and its immediate successor when present. Adjacent swaps of distinct repeated letters in complementary runs preserve the head, length class, content, and the prefix and suffix supports at every simple letter. Adjacent-repeat deletion retaining one copy in that run and two overall has the same properties. For literal deductions mark the two positions and another occurrence of each swapped letter, or the adjacent pair and a surviving third occurrence in a deletion. Mark the head also. Represent each nonempty intervening or outside unmarked block by a distinct fresh letter. The formal words retain precisely those same invariants for every formal simple letter. In particular every suffix literal or conditional omitting head keeps its position with respect to the simple occurrence. The two orientations are valid single-word absorptions of length at most 19. Replacing the fresh letters by their actual blocks proves each edit.

Sort and cap the runs at two, reserve all letters of protected positions across the three words (at most nine), group the remaining letters by their vector of counts in the at most six runs, and identify each group. Reduce the group letter in each run to the common count. When that count is one, each member had another occurrence in a different run of this word, so every deletion retains two overall. The resulting valid absorption uses at most \(9+3^{6}-1=737\) letters and at most \(3+4(737)=2951\) positions per word. Substituting the nonempty product of group members for the letter restores their common count in each run, including two. Sorting recovers the original runs. This proves unrestricted completeness of both finite families.

\par\smallskip\noindent\textit{7498.}\quad 
Use the following six fixed identities and twenty optional-context instances. Capital symbols in the two schemes may be independently omitted everywhere in that scheme. Retained symbols take nonempty arguments. The letter u in the first scheme is always retained.

{\small\par\smallskip
\CaseLine{\(xxx=xx\qquad xyx=xxy\)}
\CaseLine{\(xy=xxy+xyy\qquad xy\leq x+y\)}
\CaseLine{\(xx\leq x\qquad xyz\leq xz+y\)}
\CaseLine{\(A x B u\qquad \leq  A xx B u\qquad + u\)}
\CaseLine{\(A x B S T \leq  A xx B S T + S x T\)}
\smallskip}

The two word equations preserve first-occurrence order and total multiplicities capped at two. They derive the corresponding normal form, with copies of each letter collected at its first occurrence. To move a later x across its separating nonempty block V use x V x = xx V and cap adjacent triples at two. For successive letters in first-occurrence order collect their copies while leaving the already collected prefix fixed. All remaining laws are sound: the splitting identity follows from the rows, and \(xy\leq x+y\), \(xx\leq x\), \(xyz\leq xz+y\) follow by the exact criterion. In either scheme the first upper word has exactly the target first-occurrence order and covers all its simple letters other than x with their same multiplicities. The other upper word is contained and omits x or contains it once. This proves soundness even when optional contexts are omitted and before arbitrary substitutions.

Expand into words with at most one simple letter. Concatenate members of U until their product has full support C, using \(xy\leq x+y\). The prefix condition supplies the needed coverage. Square the product by \(xx\leq x\) to obtain a guard G with all variables repeated. Install the desired first-occurrence order from left to right. If p lists the already fixed prefix and x is next, take a contained first-x prefix witness u. For nonempty p derive p u G \ensuremath{\leq} p G + u. The first term on the right normalizes to G. For empty p derive u G \ensuremath{\leq} u+G. The result again repeats every variable and its first- occurrence order begins p x. At the end its order equals that of w. A saturated target now follows by the word equations.

If x alone is simple, write \(w=A\) x B and select a contained omitting word u or a contained word S x T with x simple. In the former case append u to the two long words of the first scheme. In the latter case the second scheme appends S and T, excluding x itself. Every appended letter has already occurred at least twice in w. Therefore the lower word normalizes to w and the first upper word to G. The other upper word is the original witness. This derives w without any order restriction on that witness. Empty contexts use their specified instance. These 6+4+16 laws give a finite basis by deriving the expanded summand absorptions in both directions. Word reversal transfers each of the three proofs to its opposite.

\subsection{An eleven-identity basis from marked suffixes}\label{case:120}
Source 4525 of the five-element catalogue. Work in the constant-free signature (+, multiplication), relative to the additively idempotent semiring axioms. Write \(u\leq v\) for u+v\ensuremath{\approx}v.

\par\smallskip\noindent\textit{Algebra and word semantics.}\quad 
The elements are 0,\allowbreak{}1,\allowbreak{}2,\allowbreak{}3,\allowbreak{}4, with these tables in that order:

{\small\par\smallskip
\noindent addition       multiplication\par
\CaseLine{\ttfamily 0 0 0 0 4\qquad 0 0 2 3 0}
\CaseLine{\ttfamily 0 1 1 1 4\qquad 0 0 2 3 0}
\CaseLine{\ttfamily 0 1 2 2 4\qquad 3 3 3 3 2}
\CaseLine{\ttfamily 0 1 2 3 4\qquad 3 3 3 3 3}
\CaseLine{\ttfamily 4 4 4 4 4\qquad 0 0 2 3 4}
\smallskip}

Addition is maximum for the chain \(3<2<1<0<4\). Multiplication is increasing in both coordinates, hence distributes over maximum. Let K=\{0,\allowbreak{}1,\allowbreak{}4\}. The element 3 is a zero, 2 squared is 3, K times 2 is 2, and 2 times a K-element is 2 precisely when that element is

4. A product of two K-elements is 4 if both are 4 and is 0 otherwise. These descriptions prove associativity by considering a triple: a 3 or two occurrences of 2 force 3; a unique 2 survives precisely when all later factors are 4; an all-K triple is 4 precisely when all factors are 4. Both associations give these same values.

Induction gives the following rule for any word of length at least two:

\textbullet{} The value is 4 if all occurrences are valued 4.

\textbullet{} It is 0 if all occurrences lie in K and at least one is 0 or 1.

\textbullet{} It is 2 if exactly one occurrence is valued 2, every later occurrence is valued 4, and every earlier occurrence lies in K.

\textbullet{} It is 3 in all remaining cases.

For a word w write c(w) for its content, \(n_{w}\)(x) for multiplicity, and \(post_{w}\)(x) for the set of variables strictly after its unique occurrence when x is simple, that is, \(n_{w}(x)=1\).

\par\smallskip\noindent\textit{The eleven identities.}\quad 
The claimed basis consists of the following identities:

\CaseDisplay{\begin{array}{rll}
(W1)&xyx\approx yx^2,\\
(W2)&x^3\approx x^2,\\
(W3)&x^2y^2\approx y^2x^2,\\
(D1)&x^2y\leq xy,\\
(D2)&xy^2\leq xy,\\
(G0)&xy\leq (xy)^2+y,\\
(G1)&xy\leq (xy)^2+x,\\
(R)&xy^2\leq y^2x,\\
(Q)&(xy)^2\leq xy,\\
(I)&x^2y^2\leq x^2,\\
(E)&xy\approx x^2y+xy^2.
\end{array}                                                     \quad\textup{(\ensuremath{B})}}

No independence or minimality is claimed.

Laws (W1) and (W2) put each word into the list of its variables in last-occurrence order, each in a block of length one or two according as its original multiplicity is one or at least two. To see this, (W1) moves an earlier x next to a later x across a nonempty intervening word; adjacent repetitions are already grouped, and (W2) caps the block size. This preserves last-occurrence order and whether each variable is simple. Law (W3) additionally sorts the consecutive blocks of length two between successive simple variables. Fix a variable order and call the result the normal form. A saturated word has no simple variable, so its normal form is the sorted product of the squares of all its variables, independent of their original order. Write s(C) for this word on a nonempty content C. The normal form of \(u^{2}\) is s(c(u)). Empty substrings are omitted, never used as unit constants or substitutions.

\par\smallskip\noindent\textit{Exact criterion.}\quad 
Let P be a nonempty sum of words and put U\_C=\{u: u is a summand of P and c(u) is contained in C\}.

Lemma 1. A linear word y is absorbed by P exactly when it is a literal summand. For a nonlinear target w, put \(C=c(w)\). The absorption \(w\leq P\) holds precisely when:

1. \(U_{C}\) contains a nonlinear word;

2. for every simple variable x of w, some u in \(U_{C}\) either omits x, or has exactly one x and \(post_{u}\)(x) contained in \(post_{w}\)(x).

Proof. A linear y absent literally is separated by \(y=2\) and all other variables 3: all right summands then have value 3. For a nonlinear target, if condition 1 fails, assign every variable of C the value 1 and outsiders 3. The target is 0; contained summands are linear and hence 1, while other summands are 3. Their maximum cannot absorb 0.

If condition 2 fails at a simple x, assign \(x=2\), the variables in \(post_{w}\)(x) the value 4, the other variables of C the value 1, and outsiders 3. The target is 2. An outside-containing summand is 3. Every contained summand has x and either repeats it or has, after its unique x, a variable outside \(post_{w}\)(x). In the latter case that later value is 1. Every right summand is therefore 3.

Conversely, target value 3 is least. At target value 4 all C-variables are 4 and any contained word is 4. At target value 0 all C-variables lie in K, and the nonlinear witness in condition 1 has value 0 or

4. At target value 2, the unique occurrence of 2 is at a simple variable x. A witness in condition 2 omitting x has all factors in K, with value at least 1; a witness containing x once has all its later factors 4 and its earlier factors in K, giving 2. This proves sufficiency for arbitrary words and variables. \(\square\)

The word rules (W1)--(W3) follow from the word-value description. All the remaining identities in (B) are sound by Lemma 1. In particular the squared guards in (G0),(G1) supply the nonlinear contained witness; their linear summand supplies every simple-variable condition. In (R), the single x moves earlier, so its required suffix only grows. In (E), \(x^{2}y\) witnesses the single y, \(xy^{2}\) witnesses the single x, and (D1),(D2) give the reverse absorptions.

\par\smallskip\noindent\textit{Saturating the content.}\quad 
Suppose \(w\leq P\) is valid and nonlinear; normalize all words. Choose a nonlinear u in \(U_{C}\) and split \(u=ab\) with a, b nonempty. Rule (Q) gives \(u^{2}\leq u\), hence \(s(c(u))\leq P\). Rule (I), with its x replaced by the current saturated word and its y by a missing variable, inserts the square of that variable below the current word. By (W1)--(W3) the right square is the same saturated word and the left is the saturated word on the enlarged content. Repeating gives

\CaseDisplay{s(C)\leq P.                       \quad\textup{(\ensuremath{1})}}

This works even when C has only one variable; no empty substitution is required.

\par\smallskip\noindent\textit{Recovering a simple variable with its suffix.}\quad 
Fix a simple variable x of w and choose its witness v in \(U_{C}\) from condition 2. We construct a normal word \(w_{x}\) with content C, exactly x simple, and \(post_\{w_{x}\}(x)=post_{w}(x)\).

First obtain a word of full content C with exactly one x and with its x-suffix contained in \(post_{w}\)(x), while retaining possible simple variables other than x:

\textbullet{} If x is absent from v, let t list the missing variables, with x last, and apply (G1) to vt. Its right terms are v and (vt)\ensuremath{^{2}}, which normalize respectively to the known v and s(C). The new word has x once at its end.

\textbullet{} If x is present, it is simple in v. If there are missing variables, list them in a nonempty word t and apply (G0) to tv. Its right terms are v and s(C). If there are no missing variables, keep v.

These words have length at least two: C has another variable whenever a nonlinear target has a simple x. Duplicate every other simple variable y by (D1) if y is the first position, and otherwise by (D2) with the nonempty prefix before y substituted for its first variable and the remaining suffix used as context. Normalization yields a word of full content C with precisely x simple. Duplicating the other variables does not change which variables occur after x.

Let \(T=post_{w}(x)\). The current suffix set at x is contained in T. For each y in T still before x, use (W3) to bring its square next to x within the saturated prefix. In the appropriate prefix and suffix context, apply (R):

\CaseDisplay{pxy^2q\leq py^2xq.                    \quad\textup{(\ensuremath{2})}}

The right side is the available word up to (W3). The left moves y from before x to after x. Substrings p or q may be omitted as contexts; the substitutions for x and y remain nonempty. Repeating and sorting the saturated blocks on each side gives \(w_{x}\leq P\). By its construction, \(w_{x}\) is exactly the normal form of the word obtained from w by squaring all positions except its unique occurrence of x.

\par\smallskip\noindent\textit{Recombination and completeness.}\quad 
Repeated use of (E) expresses any word \(a_{1}\)...\(a_{n}\), \(n\geq 2\), as the sum of the n words obtained by leaving one position single and squaring the others, modulo (W1)--(W3). Indeed split off the last letter: (E) gives \((a_1\cdots a_{n-1})\)\ensuremath{^{2}} \(a_{n}\) and \((a_1\cdots a_{n-1})\)\(a_{n}^{2}\). The square in the first term saturates every variable of the prefix; expand the prefix in the second term inductively. Normalization identifies a squared prefix with the product of the squared prefix letters, since both are saturated with the same content.

Each resulting normal word is s(C), when the position left single belongs to a repeated variable, or is \(w_{x}\) for a simple variable x. Equation (1) and the recovery construction already put all these terms below P. Thus \(w\leq P\) is derivable. A target without simple variables is already s(C), and a linear target is a literal right summand by Lemma 1.

Theorem. The eleven identities (B) form a finite identity basis for source 4525, relative to the additively idempotent semiring axioms. Their word reversals form a basis for its multiplicative opposite.

Proof. Soundness and every word absorption have been proved above. Distributivity reduces arbitrary terms to polynomials. Derive the absorption of each summand into the opposite polynomial in both directions; the two resulting inequalities give every valid identity. The proof places no bound on variables, lengths, or the number of summands. Reversal transports all steps to the opposite. \(\square\)

Finite computations independently check the operation tables and individual substitution/context steps, not unrestricted completeness. Transfers to further catalogue types require independently verified complete HSP membership in both directions. No hereditary or minimal finite-basis claim is made.

\subsection{Finite bases from contained head and simple-final witnesses}\label{case:121}
Sources 4533 and 4534, in (+, multiplication) without constants, relative to the additively idempotent semiring axioms. A one-letter word has length one; every other word has length at least two. A variable is simple when it occurs exactly once.

\par\smallskip\noindent\textit{1. Algebras and word values.}\quad 
Both sources have multiplication, in order 0,\allowbreak{}1,\allowbreak{}2,\allowbreak{}3,\allowbreak{}4:

{\small\par\smallskip
\CaseLine{\ttfamily 0 0 2 3 0}
\CaseLine{\ttfamily 0 0 2 3 0}
\CaseLine{\ttfamily 3 3 3 3 3}
\CaseLine{\ttfamily 3 3 3 3 3}
\CaseLine{\ttfamily 4 4 2 3 4}
\smallskip}

Addition is join in \(3<1\),\(2<0<4\) for 4533, with 1 and 2 incomparable, and in \(3<2<1<0<4\) for 4534. Put K=\{0,\allowbreak{}1,\allowbreak{}4\}. A word of length at least two has value 3 if any nonfinal factor is 2 or 3, or its final factor is 3. Otherwise all nonfinal factors lie in K. If its final factor is 2 the value is 2; if all factors lie in K the value is 4 when the head is 4 and 0 otherwise. Induction on length proves this formula. In particular value 2 requires a simple final variable.

The description proves associativity by three-factor cases: both associations are 3 if a 2/\allowbreak{}3 occurs before the last factor or the last is 3; with the first two factors in K, a last factor 2 gives 2; with all three in K both products are 4 exactly when the first factor is 4, and otherwise 0. Rows and columns preserve joins in both displayed orders, proving both distributive laws.

\par\smallskip\noindent\textit{2. Word normalization.}\quad 
The value formula proves

\CaseDisplay{xyzt\approx xzyt,\qquad x^2y\approx xy,\qquad xyx\approx xy^2.\quad\textup{(\ensuremath{W})}}

The first law permutes interior letters, keeping head and final position. The second removes adjacent repeated interior letters and copies of the head, always with a nonempty suffix. If the final s is simple, the normal word is its head h, a sorted list of the remaining internal variables, and s. Its data are (h, internal set I, s).

If the final is repeated, the normal word is determined just by (head h, content C). Singleton content gives hh, since \(x^{3}=x^{2}\) is an instance of the second law. Otherwise fix a=min(C minus \{h\}) and use h followed by the sorted list C minus \{h\}, then a. To obtain this form first delete unnecessary internal copies. If the old final is h, write hRh and use \(hRh=hR^{2}\), with R nonempty, producing a repeated final different from h. If the final is b different from a, both b and a occur internally. Permute the interior to end in ba and use \(bab=baa\) to change the final to a. Delete extra copies and sort the interior. All substitutions are nonempty. Every normal word has length at most its number of variables plus one.

\par\smallskip\noindent\textit{3. Exact criterion.}\quad 
For target w with content C, let U be the original right words in P with content contained in C. A one-letter y is absorbed precisely when literal y occurs in P. A word w of length at least two, with head h, is absorbed precisely when:

\textbullet{} H: U has a word headed by h;

\textbullet{} N: U has a word of length at least two;

\textbullet{} if the final variable s is simple, A: U has a word with no s in a nonfinal position. In 4533 it must have length at least two or be literal s; in 4534 there is no further restriction.

There is no A requirement for a repeated final variable.

Necessity. Without literal y assign \(y=2\) and outsiders 3; all other right words are 3. If H fails assign \(h=4\), C minus h the value 0, and outsiders 3. The target is 4, the right sum at most 0. If N fails assign C the value 1 and outsiders 3; the target is 0 while \(P\leq 1\). If A fails at simple final s, assign \(s=2\), C minus s the value 1, and outsiders 3. The target is 2. Every contained right word of length at least two has s internally and is 3; every noncontained such word is also 3. In 4533 literal s is absent, so remaining literals give \(P\leq 1\), which does not dominate 2. In 4534 every literal in C would supply A, so \(P=3\).

Sufficiency. Target 3 is least. For target 4 all variables lie in K and \(h=4\); H supplies 4. For target 0 all variables lie in K; N supplies 0 or 4. For target 2 its simple final s is 2 and all other variables lie in K. A witness of length at least two gives 0,\allowbreak{}2 or 4; literal s gives 2. In 4534 every other allowed literal lies in K and is also at least 2. This proves the criterion with no bound on variables or summands.

\par\smallskip\noindent\textit{4. Direct witnesses and reversible compression.}\quad 
For 4533 select original H, N, A witnesses, using only H, N for repeated final targets. At most three original right words suffice. For 4534 at most two suffice: a length-at-least-two H witness also supplies N, so add A when needed; if H is literal h, it supplies A when the final is simple (and thus differs from h), so add N. Repeated-final targets need only H, N. One-letter targets are literal summands.

All selected right words are contained in C. In a step with k right words reserve every head and simple final variable in all k+1 words separately: at most 2(k+1) variables. Every other variable lies in C; group it by incidence among the k right words, including the zero pattern. There are at most \(2^{k}\) groups. Identify each group and normalize by(W). Substitution preserves validity. The bounds are 16 variables for \(k\leq 3\), and 10 for \(k\leq 2\).

For the inverse substitution, replace a group by a nonempty product listing its members once. A group is never a head or simple final. A simple-final word recovers its head, internal content and simple final and hence its normal form. A repeated-final word stays so: its final group, if present, occurred earlier too, so its last member remains repeated after expansion. Its head and content recover the original first-kind normal form. Reserved variables are replaced individually. This includes singleton content and uses no constants or empty substitutions.

\par\smallskip\noindent\textit{5. Complete finite families.}\quad 
For 4533 take(W) and all valid normal absorptions over a fixed 16-letter alphabet, length at most 17, with at most three right words. For 4534 use a fixed 10-letter alphabet, length at most 11 and at most two right words. These are finite explicitly specified families, with membership decided by Section 3; no small or minimal list is claimed.

Theorem. These families are relative identity bases. Direct witness selection and inverse compression derive every valid word absorption into an arbitrary polynomial, allowing all original variable counts and word lengths. Distribute arbitrary terms into sums of words. Derive absorption of each summand on either side to obtain every valid identity. Reversing words transports the bases and derivations to the respective multiplicative opposites.

Finite calculations support the unrestricted proof; transfers to catalogue entries require both complete HSP directions.

\subsection{Recovering simple variables after saturation: a finite basis for source 4572}\label{case:122}
Work in the signature (+, multiplication), without constants, relative to the axioms for additively idempotent semirings.

\par\smallskip\noindent\textit{Algebra and word values.}\quad 
Let B=\{0,\allowbreak{}1,\allowbreak{}2,\allowbreak{}3,\allowbreak{}4\}, with these tables in that order.

{\small\par\smallskip
\noindent addition       multiplication\par
\CaseLine{\ttfamily 0 0 0 0 4\qquad 0 0 2 3 0}
\CaseLine{\ttfamily 0 1 0 1 4\qquad 0 1 2 3 0}
\CaseLine{\ttfamily 0 0 2 2 4\qquad 2 2 3 3 2}
\CaseLine{\ttfamily 0 1 2 3 4\qquad 3 3 3 3 3}
\CaseLine{\ttfamily 4 4 4 4 4\qquad 4 4 2 3 4}
\smallskip}

Addition is join for \(3<1<0<4\) and \(3<2<0<4\), with 1 and 2 incomparable. Put K=\{0,\allowbreak{}1,\allowbreak{}4\}. The element 3 is a multiplicative zero; 2 squared is 3, and multiplication of 2 with any element of K, on either side, is 2. In a product restricted to K, an initial 4 persists as 4; otherwise the product is 1 if every factor is 1 and is 0 otherwise. This verifies associativity: a triple containing 3 or at least two occurrences of 2 gives 3 in either association, a triple with exactly one 2 gives 2, and the stated K rule is independent of association. Every multiplication row and column is increasing and preserves the only incomparable join \(1+2=0\), proving distributivity.

For a word w write h(w) for its head, c(w) for its content, and \(n_{w}\)(x) for the number of occurrences of x. A variable is simple in w when \(n_{w}(x)=1\). For every nonempty word, including a linear word, successive multiplication gives this exact rule:

\textbullet{} An occurrence valued 3, or at least two occurrences valued 2, gives 3.

\textbullet{} Exactly one occurrence valued 2 and no occurrence valued 3 gives 2.

\textbullet{} Otherwise all factors lie in K: the value is 4 if the head is 4, is 1 if all factors are 1, and is 0 otherwise.

It follows that the identities

\CaseDisplay{xyz\approx xzy,\qquad x^3\approx x^2           \quad\textup{(\ensuremath{W})}}

hold. Fix a variable order. Every word normalizes by (W) to its head followed by sorted tail blocks, with total multiplicity of each variable equal to min(2, \(n_{w}\)(x)). For the head variable, its leading occurrence already contributes one, so it has at most one further tail occurrence. Indeed tail transpositions gather occurrences of each variable, and the second law removes occurrences beyond two, including those of the head. All substitutions are nonempty.

\par\smallskip\noindent\textit{Exact absorption criterion.}\quad 
Let P be a nonempty sum of words, let \(C=c(w)\), and put U\_C=\{u: u is a summand of P and c(u) is contained in C\}.

Lemma 1. B satisfies \(w\leq P\) exactly when all three conditions hold:

1. Some u in \(U_{C}\) has \(h(u)=h(w)\).

2. The union of c(u) over u in \(U_{C}\) is C.

3. For every variable x simple in w, some u in \(U_{C}\) has \(n_{u}(x)=1\).

For a linear target y these conditions are equivalent to y being a literal linear summand of P.

Proof. For failure of condition 1, assign \(h(w)=4\), the other variables of C the value 1, and outsiders the value 3. The target is 4, while every right word containing an outsider is 3 and every remaining word is at most 0. Thus absorption fails.

If condition 2 fails at x in C, assign \(x=0\), C minus x the value 1, and outsiders 3. The target is 0, whereas the right words are 1 or 3, whose join is at most 1. For failure of condition 3 at a simple x, instead assign \(x=2\). The target is 2; contained right words either omit x and have value 1, or repeat x and have value 3. Other words have value 3. Their join is at most 1, which does not dominate 2.

Conversely consider any assignment. Target value 3 is least. Target value 2 comes from a unique occurrence of 2, at a simple variable x; the witness in condition 3 has exactly one x, all its other factors lie in K, and its value is 2. When all target factors lie in K, the witness in condition 1 absorbs value 4 or value 1. If the target is 0, some variable x has value 0 or 4. By condition 2, a contained right word includes x. This word has value 0 or 4 and absorbs the target. This proves sufficiency for arbitrary words and variables. The assertion for a linear word follows by applying condition 3 to its unique variable. \(\square\)

One further identity, used to assemble simple occurrences, is

\CaseDisplay{xy\approx x^2y+xy^2.                       \quad\textup{(\ensuremath{E})}}

It follows from Lemma 1: the two right terms supply the same head, content coverage, and respectively the single y and single x witnesses; each right term is also absorbed by xy. The valid inequality \(x^{2}\leq x\) will be among the bounded absorptions below.

\par\smallskip\noindent\textit{Constructing every valid absorption in bounded-arity steps.}\quad 
For a word u, let s(u) be the normal form of \(u^{2}\). It has the same head and content as u and every variable has total multiplicity two. Suppose \(w\leq P\) is valid, and normalize all words. A linear w is already a literal summand. For nonlinear w use the following construction.

Choose the same-head witness \(u_{0}\) in \(U_{C}\). The absorption \(s(u_{0})\leq u_{0}\) is valid by Lemma 1. Starting with its content \(D=c(u_{0})\), take the other members v of \(U_{C}\) one at a time. If new content is needed, let q be the saturated normal word with head h(w) and content D union c(v). If p is the previously derived saturated word, then

\CaseDisplay{q\leq p+v                             \quad\textup{(\ensuremath{3})}}

is valid by Lemma 1: p supplies the head, the two terms cover the content, and q has no simple variable. Each step has at most two right words and uses only already available absorptions below P. Condition 2 eventually gives \(s(w)\leq P\).

For each variable x simple in w, let \(w_{x}\) be the normal word with the same head and content as w, with x occurring once and every other variable twice. Choose \(u_{x}\) in \(U_{C}\) with \(n_\{u_{x}\}(x)=1\). Again Lemma 1 gives the valid two-term absorption

\CaseDisplay{w_x\leq s(w)+u_x.                       \quad\textup{(\ensuremath{4})}}

Thus all \(w_{x}\) are below P. To finish, repeated use of (E) expands a word \(a_{1}\)...\(a_{n}\), for \(n\geq 2\), as the sum of the n words obtained by leaving one position \(a_{i}\) single and squaring all other positions, followed by normalization. To verify the expansion, split off the last letter: (E) gives \((a_1\cdots a_{n-1})\)\ensuremath{^{2}} \(a_{n}\) and \((a_1\cdots a_{n-1})\)\(a_{n}^{2}\). The first normalizes to the last-position term; expand the prefix in the second term inductively. A squared prefix has each variable twice after (W), regardless of its original order or multiplicity.

If a variable repeats in w, every resulting term in which one of its occurrences was left single still repeats that variable. Therefore every normalized term is s(w) or \(w_{x}\) for a simple x. All these terms are already below P. If w has no simple variables, its normal form is s(w) directly. This derives every absorption, using (E), (W), and valid absorptions with at most two right words.

\par\smallskip\noindent\textit{Reversible multiplicity compression.}\quad 
In such a local normal absorption there are at most three words. Reserve the head variable of each word, including the unique variable of a linear word; there are at most three reserved variables. For every other variable record its total multiplicities 0, 1, or 2 in these words. Group variables with equal nonzero multiplicity vectors. There are at most \(3^{3}-1=26\) groups, besides the reserved singletons.

Construct each compressed word with the original head's singleton and with each group variable occurring its recorded multiplicity. This construction is NOT ordinary variable identification: merging two single occurrences would incorrectly turn them into a multiple occurrence. Validity is instead supplied by Lemma 1. For each of its three conditions, content containment, head equality, content coverage, and a count-one witness are preserved by the grouped vectors. In particular any group simple in the target consists of variables all simple there, and a count-one witness for one of its members has count one for every member. The compressed absorption is therefore valid. There are at most 29 variables and word length at most 58.

Substitute the nonempty product of a group's members for its group variable, and restore reserved singletons literally. Groups never occur in head position. Tail commutation puts the expansion in normal order; one or two occurrences of a group produce respectively one or two occurrences of each original member. Thus (W) recovers exactly all original normal words. This proves reversible compression for every local step, without a bound on the original variables.

Fix a 29-letter alphabet X. Let Sigma consist of (W), (E), and all valid normal absorptions over X with one left word, at most two right words, and word length at most 58. This is a finite, explicitly defined family whose membership is decided by Lemma 1. A small list has not been enumerated, and no minimality is asserted.

Theorem. Sigma is a finite identity basis for B. Word reversal gives a finite identity basis for its multiplicative opposite.

Proof. Soundness follows from the definitions and the displayed valid identities. Compression makes every local absorption used in (3), (4), and saturation derivable from Sigma by substitution and (W). The preceding construction and (E) then derive every valid word absorption. Distribute arbitrary terms into sums of words and derive every summand absorption in both directions to obtain every identity. This proves completeness for unrestricted variables and word lengths. Reversing words transports all identities and derivations to the opposite. \(\square\)

Finite regressions check transcription and implementation errors; they do not prove completeness. Catalogue transfers require separately verified equality of generated varieties in both HSP directions. No claim about hereditary finite basability is made.

\subsection{Two further chain semirings: simple omission and marked suffixes}\label{case:123}
The signature is (+, multiplication), without constants; all bases are relative to the additively idempotent semiring axioms. The two algebras D and E are catalogue entries 4573 and 4575.

\par\smallskip\noindent\textit{Tables.}\quad 
Both algebras have addition equal to maximum for \(3<2<1<0<4\). Their tables, in the order 0,\allowbreak{}1,\allowbreak{}2,\allowbreak{}3,\allowbreak{}4, are:

{\small\par\smallskip
\noindent addition       D multiplication   E multiplication\par
\CaseLine{\ttfamily 0 0 0 0 4\qquad 0 0 2 3 0\qquad 0 0 2 3 0}
\CaseLine{\ttfamily 0 1 1 1 4\qquad 0 1 2 3 0\qquad 0 1 2 3 0}
\CaseLine{\ttfamily 0 1 2 2 4\qquad 2 2 3 3 2\qquad 3 3 3 3 2}
\CaseLine{\ttfamily 0 1 2 3 4\qquad 3 3 3 3 3\qquad 3 3 3 3 3}
\CaseLine{\ttfamily 4 4 4 4 4\qquad 4 4 2 3 4\qquad 0 0 2 3 4}
\smallskip}

Every multiplication row and column is increasing, proving both distributive laws. Put K=\{0,\allowbreak{}1,\allowbreak{}4\}. In both algebras 3 is a zero and two occurrences of 2 give 3. In D a unique 2 survives regardless of its position among K-elements. An all-K product is 4 if its head is 4, is 1 if all factors are 1, and is 0 otherwise. In E a unique 2 survives exactly when every later factor is 4. An all-K product is 1 if all factors are 1, 4 if all are 4, and 0 otherwise. These statements apply to one-letter words too. They verify associativity by a case split on triples, and induction proves the same descriptions for words of arbitrary length.

Write c(w), h(w), and \(n_{w}\)(x) for content, head, and multiplicity. For a simple variable x, \(post_{w}\)(x) is its suffix content. For a polynomial P and target content C, let \(U_{C}\) be its summands whose contents are contained in C.

\par\smallskip\noindent\textit{D: a contained head, content coverage, and omission witnesses.}\quad 
Lemma 1. D satisfies \(w\leq P\) exactly when:

1. Some u in \(U_{C}\) has \(h(u)=h(w)\).

2. The union of c(u) over u in \(U_{C}\) is C.

3. For every x simple in w, some u in \(U_{C}\) has \(n_{u}(x)\leq 1\).

For a linear target this is equivalent to a literal linear summand.

Proof. For failure of condition 1, assign the target head 4, other C-variables 1 and outsiders 3. The target is 4; all right terms are at most 0. For failure of condition 2 at x, assign \(x=0\), other C-variables 1 and outsiders 3. The target is 0 and all right terms are at most 1. If condition 3 fails at a simple x, assign \(x=2\), other C-variables 1 and outsiders 3. Every contained right term repeats x, so all right terms are 3 while the target is 2.

Conversely, target value 3 is least. For target value 2 the unique 2 is at a simple x. Its witness in condition 3 either contains x once and evaluates to 2, or omits it and has all factors in K, giving a value at least 1. For target value 4 the same-head witness has value 4. For target value 1 any contained witness has value 1. At target value 0 some C-variable has value 0 or 4; a contained word including that variable, supplied by condition 2, has value 0 or 4. This proves sufficiency. For a linear y every nonempty contained word is a power of y, so condition 3 forces exponent one. \(\square\)

The following identities hold in D:

\CaseDisplay{xyz\approx xzy,\qquad x^3\approx x^2,\qquad
 xy\approx x^2y+xy^2.                                          \quad\textup{(\ensuremath{D})}}

The first two give normal words consisting of the head followed by sorted tail blocks, with each variable's total count capped at two. The head contributes one to that count. The last identity follows from Lemma 1 in both directions. It is used for recombination below.

For a word u let s(u) be the normal form of \(u^{2}\): every variable occurs twice, with the original head. Start a valid absorption from its same-head witness u. The absorption \(s(u)\leq u\) is valid. Given a saturated word p already derived below P and another v in \(U_{C}\), the saturated word q with target head and content c(p) union c(v) satisfies \(q\leq p+v\), by Lemma 1. Content coverage constructs \(s(w)\leq P\). For each x simple in w, let \(w_{x}\) have the same head and content as w, with x single and every other variable twice. Its count-at-most-one witness \(u_{x}\) gives \(w_{x}\leq s(w)+u_{x}\) by the same lemma. Repeated use of the last identity in (D) expands w as the sum of its one-unsquared-position words. After normalization, each is s(w) or one of the \(w_{x}\). Hence every word absorption is a consequence of (D) and valid absorptions with at most two right words. All intermediate right terms have already been derived below P.

To bound these local absorptions, reserve the at most three heads and group all other variables by their nonzero total multiplicity vectors in the at most three words. There are at most \(3^{3}-1=26\) groups. Construct a compressed word with each group variable at its original count 0, 1 or 2, rather than identifying occurrences and accidentally changing a single count to a double count. Lemma 1 remains true: containment, head equality, content coverage, and count-at-most-one witnesses are preserved by uniform vectors. Thus the compressed absorption is valid, uses at most 29 variables, and has word length at most 58. Substituting the product of each group's original members back recovers every normal word, because groups occur only in tails and the tail is commutative. Reserved heads are restored singly. This is reversible compression.

Consequently (D), together with all valid normal absorptions over a fixed 29-letter alphabet having length at most 58 and at most two right words, is a finite basis for D. Lemma 1 decides this explicitly defined finite family's membership. No small list has been enumerated. Compression derives each local step, and saturation, simple recovery, and recombination derive all absorptions with arbitrarily many original variables. Distributivity and both directions of summand absorption then derive every identity.

\par\smallskip\noindent\textit{E: content coverage with ordered single occurrences.}\quad 
Lemma 2. E satisfies \(w\leq P\) exactly when:

1. The union of c(u) over u in \(U_{C}\) is C.

2. For each simple x of w, some u in \(U_{C}\) omits x or has exactly one x and \(post_{u}\)(x) contained in \(post_{w}\)(x).

Again a linear target must be a literal summand.

Proof. If condition 1 fails at x, assign \(x=0\), the other variables of C the value 1, and outsiders 3. The target is 0, while every right term is 1 or 3. If condition 2 fails at a simple x, assign \(x=2\), \(post_{w}\)(x) the value 4, the other C-variables 1 and outsiders 3. The target is 2. Every contained right term repeats x or has a later 1, and every other term contains 3; all right values are 3.

For sufficiency target value 3 is least. Target value 2 is absorbed by its condition-2 witness, which is either 2 or an all-K word of value at least 1. Target values 1 and 4 mean all C-variables have that respective value, so any contained word absorbs the target. If the target is 0, some variable has value 0 or 4. A contained word covering it, supplied by condition 1, has value 0 or 4. A linear target forces a single occurrence in a nonempty contained word, hence a literal summand. \(\square\)

E has the following explicit eleven-identity basis:

\CaseDisplay{\begin{array}{rll}
(W1)&xyx\approx yx^2,\\
(W2)&x^3\approx x^2,\\
(W3)&x^2y^2\approx y^2x^2,\\
(D1)&x^2y\leq xy,\\
(D2)&xy^2\leq xy,\\
(G0)&xy\leq (xy)^2+y,\\
(G1)&xy\leq (xy)^2+x,\\
(R)&xy^2\leq y^2x,\\
(S)&x^2\leq x,\\
(J)&x^2y^2\leq x^2+y^2,\\
(E0)&xy\approx x^2y+xy^2.
\end{array}                                                     \quad\textup{(\ensuremath{E})}}

All word rules follow from the word semantics; Lemma 2 checks the absorptions and both directions of (E0). In particular (S) holds also at 1, whose square is 1 in this algebra. Law (J) uses content coverage by its two right terms. This is why the content-enlargement rule \(x^{2}y^{2}\leq x^{2}\) for source 4525 is not used here.

Laws (W1),(W2) give last-occurrence order with single or double blocks; (W3) sorts the consecutive double blocks between simple variables. For a nonempty content C write s(C) for the sorted product of its squares. Every squared word with content C normalizes to s(C).

For a valid nonlinear absorption, choose a word u in \(U_{C}\) and use (S) to derive \(s(c(u))\leq u\). For each further v in \(U_{C}\) derive \(s(c(v))\leq v\) and apply (J) to the two saturated words. Their squares normalize to themselves, while their product normalizes to the saturated word on the union. Condition 1 therefore gives \(s(C)\leq P\).

For a simple target variable x, choose its witness v from condition 2. If v omits x, append the missing variables with x last, using (G1); if v contains x, prepend the missing variables using (G0), or leave v unchanged if none are missing. The other right term in these rules is s(C), already available. This creates full content C, exactly one x, and a suffix at x contained in its target suffix. Duplicate every other simple variable by (D1) at the first position or (D2) in a nonempty-prefix context elsewhere. Use (W3) inside the saturated prefix and (R) in context to move every required square from before x to after x. The resulting \(w_{x}\) has exactly x simple and exactly the target suffix set; thus \(w_{x}\leq P\). No empty word is substituted. The target is nonlinear with a simple x, so C has another variable and the temporary full-content word always has length at least two.

Finally, expand w by (E0) into the words with one position single and all other positions squared. This expansion follows inductively by splitting off the last letter. Squared prefixes normalize to their saturated content, so the normal terms are exactly s(C) or the constructed \(w_{x}\). Hence \(w\leq P\) follows. Linear targets are literal right summands. Distribute arbitrary terms and absorb every summand in both directions to obtain unrestricted completeness of (E).

Theorem. D and E are finitely based, by the finite family described above for D and the eleven explicit identities (E) for E. Their multiplicative opposites are finitely based by word reversal.

The construction applies to arbitrary word lengths, alphabets and numbers of summands.

\subsection{Repeated-final guards and two-term steps for two chain semirings}\label{case:124}
The sources A and B are catalogue entries 4581 and 4867. Work in the signature (+, multiplication), without constants, relative to the axioms of additively idempotent semirings. Element labels below are not constants of the signature. Write \(w\leq P\) for w+P\ensuremath{\approx}P.

\par\smallskip\noindent\textit{Operations and arbitrary word values.}\quad 
Both additions are maximum in the chain \(3<2<1<0<4\). Tables are displayed in the element order 0,\allowbreak{}1,\allowbreak{}2,\allowbreak{}3,\allowbreak{}4.

{\small\par\smallskip
\noindent addition       A multiplication   B multiplication\par
\CaseLine{\ttfamily 0 0 0 0 4\qquad 0 0 2 3 0\qquad 0 0 2 3 4}
\CaseLine{\ttfamily 0 1 1 1 4\qquad 0 1 2 3 0\qquad 0 1 2 3 4}
\CaseLine{\ttfamily 0 1 2 2 4\qquad 3 3 3 3 3\qquad 3 3 3 3 3}
\CaseLine{\ttfamily 0 1 2 3 4\qquad 3 3 3 3 3\qquad 3 3 3 3 3}
\CaseLine{\ttfamily 4 4 4 4 4\qquad 4 4 2 3 4\qquad 0 0 2 3 4}
\smallskip}

Put K=\{0,\allowbreak{}1,\allowbreak{}4\}. In both algebras 3 is a zero, K times 2 is 2, and 2 followed by any further factor gives 3. A product of K-elements stays in K. In A it is 4 when the head is 4, 1 when all factors are 1, and 0 otherwise. In B replace the head condition by the condition that the final factor is 4. These rules include linear words. For products of length at least two containing a factor outside K, the value is 2 exactly when the unique such occurrence is a 2 at the final position; otherwise the value is 3.

Associativity follows by considering a triple. A factor 3 or a 2 in either of the first two positions forces 3 in both associations. A final 2 preceded by K-elements gives 2 in both. All-K triples obey the respective endpoint rule, independently of association. The multiplication rows and columns are increasing in the displayed chain, proving both distributive laws. Induction now proves the stated word-value rules for all nonempty words.

Let c(w), h(w), f(w), and p(w) denote content, head, final variable, and the set of variables in all positions except the final one. For a polynomial P and \(C=c(w)\), let \(U_{C}\) consist of its summands u such that c(u) is contained in C.

\par\smallskip\noindent\textit{Exact criteria.}\quad 
Lemma 1. In A the absorption \(w\leq P\) holds if and only if:

1. Some u in \(U_{C}\) has \(h(u)=h(w)\).

2. The union of c(u), over u in \(U_{C}\), equals C.

3. If \(s=f(w)\) does not belong to p(w), some v in \(U_{C}\) either omits s or has exactly one occurrence of s, at its final position.

In B replace condition 1 by the existence of a u in \(U_{C}\) with \(f(u)=f(w)\); conditions 2 and 3 are unchanged. In either algebra a linear target is absorbed exactly when it is a literal summand.

Proof. For failure of condition 1 in A, assign the target head 4, the other C-variables 1, and outsiders 3. The target is 4; every outside-containing term is 3, and every contained term has a different head and value at most 0. For B use the same assignment with the target final variable assigned 4. Each contained term has a different final variable, so again every right value is at most 0.

If condition 2 fails at x, assign \(x=0\), other C-variables 1 and outsiders 3. The target is 0; every contained term omits x and evaluates to 1, and every other term is 3. Their maximum is at most 1. If condition 3 fails, assign \(s=2\), other C-variables 1 and outsiders

3. The target is 2, including a possible linear target s. Every contained right term has s but repeats it or has an occurrence before its last position, so its value is 3. All other terms also have value 3. This proves necessity.

For sufficiency, a target valued 3 is least. A target valued 2 has its unique 2 at its final, simple variable s. The witness in condition 3 is either an all-K word of value at least 1 (if s is absent), or has its single 2 at the end and value 2. If all target factors lie in K and the target is 4, condition 1 supplies a contained word with the required endpoint valued 4, and hence that word is 4. A target valued 1 means all its variables are 1, so any contained word has value 1. A target valued 0 has some variable x valued 0 or 4. Condition 2 supplies a contained word including x. Its factors lie in K, and its value is 0 or 4, either of which absorbs the target.

For linear \(w=y\), nonempty words contained in \(C=\{y\}\) are positive powers of y. Condition 3 can only be witnessed by y itself. The converse is immediate. This also verifies sufficiency for linear targets and completes the proof. \(\square\)

\par\smallskip\noindent\textit{Word normalization.}\quad 
The following two word laws hold in A:

\CaseDisplay{xyzt\approx xzyt,\qquad x^2y\approx xy.       \quad\textup{(\ensuremath{A})}}

Fix a variable order. The first law permutes interior positions; the second removes repeated letters before the final position. It can remove a repeat next to the head, or in an interior prefix context with the nonempty remaining suffix substituted for y. Thus every nonlinear word becomes

\CaseDisplay{h(w)\,\prod_{x\in p(w)\setminus\{h(w)\}}x\,f(w).   \quad\textup{(\ensuremath{3})}}

The product is sorted. Empty products mean omitted substrings. This normalization preserves the head, final variable and prefix content, including whether the final variable occurs earlier. Linear words are left unchanged. No uniqueness modulo all identities is needed.

In B the two word laws are

\CaseDisplay{xyz\approx yxz,\qquad x^2y\approx xy.        \quad\textup{(\ensuremath{B})}}

They permute all positions before the last and remove repeated prefix occurrences. The normal word is the sorted product of p(w) followed by f(w). For a linear word the prefix is empty. Again all substituted words in these rules are nonempty. Soundness of (A) and (B) follows immediately from the word-value rules.

\par\smallskip\noindent\textit{Bounded-arity construction of every absorption.}\quad 
For A and nonempty D containing h, let \(g_{A}\)(h, D) be the normal word with head h, prefix content D, and final h. For B and s in D let \(g_{B}\)(s, D) be the normal word with prefix content D and final s. These guard words have a repeated final variable.

Consider a valid nonlinear \(w\leq P\) and put \(C=c(w)\). In A choose the condition-1 witness u and start with \(g_{A}\)(h(w), c(u))\ensuremath{\leq}u. Lemma 1 proves this absorption: u supplies the head and full content, and the left final is repeated, so no third condition is imposed. In B start instead with \(g_{B}\)(f(w), c(u))\ensuremath{\leq}u using its same-final witness. Each is a one-right-term valid absorption, even when u is linear.

Suppose a guard q on content D has already been derived below P. For each further v in \(U_{C}\) enlarge it to the corresponding guard r on D union c(v). The valid absorption

\CaseDisplay{r\leq q+v                       \quad\textup{(\ensuremath{4})}}

follows from Lemma 1: q supplies the required endpoint, q and v cover the enlarged content, and r has a repeated final. Thus content coverage gives the full-content guard \(g\leq P\). At most two right words are used in each step.

If the final of w is repeated, Lemma 1 gives \(w\leq g\). If it is simple, choose its witness v in condition 3. Lemma 1 gives \(w\leq g+v\): the guard supplies the endpoint and the full content, and v supplies the last condition. These absorptions again have at most two right words. This constructs every valid nonlinear absorption using only valid one- or two-right-term normal absorptions and the word laws. A linear absorption is a literal inclusion.

\par\smallskip\noindent\textit{Reversible finite-variable compression.}\quad 
In a local A-absorption there are at most three words in total. Reserve all their endpoint variables, at most six, including the unique variable of any linear word. Every other variable occurs only in interiors. Group these by their nonempty interior-incidence patterns across the at most three words, giving at most seven groups. Identify each group to a new variable, keep each reserved variable separate, and normalize by (A). This preserves validity by substitution and gives at most 13 variables and normal word length at most 14 (the looser bound 15 will be used).

Conversely substitute the product of a group's original members for its new variable, and restore reserved variables as single letters. Endpoints are unchanged. Uniform interior incidence recovers the same prefix content in each word; (A) sorts interiors and removes prefix repetitions to recover the original normal word. This includes repeated final variables and singleton/linear words.

For B reserve only the final variable of each word, at most three. All other variables occur only in prefixes; group them by their nonempty prefix-incidence patterns, at most seven. Identification and normalization preserve validity, with at most 10 variables and word length at most 11. Substituting the original group products back restores the prefix sets and final variables modulo (B), so this compression is also reversible.

Theorem. A finite basis for A consists of (A) and all its valid normal absorptions over a fixed 13-letter alphabet, with word length at most 15 and at most two right words. A finite basis for B consists of (B) and the corresponding valid normal absorptions over 10 letters, word length at most 11, and at most two right words. Their opposites have the bases obtained by reversing every word.

Proof. These are finite, explicitly specified families. Membership is decided by Lemma 1, and soundness follows from their definition. Compression and substitution derive each local absorption used in the guard construction from the appropriate family. The construction then derives every valid word absorption with unrestricted variables, lengths, and numbers of summands. Distributivity converts arbitrary terms to polynomials; derive the absorption of every summand into the opposite polynomial in both directions to obtain every identity. Word reversal transports validity and all derivations to each opposite. \(\square\)

The finite families have not been enumerated as small lists; no minimality or hereditary claim is made. Finite computation supplements the written unrestricted argument, and catalogue transfers require independently verified complete HSP equality in both directions.

\subsection{首字母与平方展开：两组非交换完整有限基}\label{case:125}
\par\smallskip\noindent\textit{1. 两个生成元.}\quad 
令 \(P=\{0, a, e, f, t\}\)，加法取链 \(0<a<e<f<t\) 的最大值。在 \(\{0, a, e, f\}\) 上，\(0\) 是乘法零，\(a^2=0\)，\(\{e, f\}\) 是左零半群，且

\CaseDisplay{ea=ae=fa=af=a.}

\(t\) 的乘法作用与 \(f\) 相同，包括 \(t^2=f\)。因此 \(P\) 的乘法为

\CaseDisplay{\begin{array}{c|ccccc}
\cdot&0&a&e&f&t\\\hline
0&0&0&0&0&0\\
a&0&0&a&a&a\\
e&0&a&e&e&e\\
f&0&a&f&f&f\\
t&0&a&f&f&f
\end{array}.}

它对应目录 4588，标号映射为 \((3,\allowbreak{}2,\allowbreak{}1,\allowbreak{}0,\allowbreak{}4)\)。一个非线性词的值：出现 \(0\) 或至少两次 \(a\) 时为 \(0\)；恰出现一次 \(a\) 时为 \(a\)；否则为首字母在 \(e, f, t\mapsto e, f, f\) 下的像。这说明结合性。各平移单调，故保持链上的并。

令 \(Q=\{0, b, a, c,1\}\)，加法偏序为

\CaseDisplay{0<b<c<1,\qquad 0<a<c<1,\qquad a\parallel b,}

乘法为

\CaseDisplay{\begin{array}{c|ccccc}
\cdot&0&b&a&c&1\\\hline
0&0&0&0&0&0\\
b&b&b&b&b&b\\
a&0&0&0&0&a\\
c&b&b&b&b&c\\
1&b&b&c&c&1
\end{array}.}

该符号 \(1\) 仅为元素名称，它是右幺元，但不是全体代数的双侧幺元。它对应目录 8983，标号映射为 \((3,\allowbreak{}1,\allowbreak{}2,\allowbreak{}0,\allowbreak{}4)\)。

一个结构性构造是：取两元素左零半环 \(L_2\)，加法为二元素链上的并，乘法为左投影；取三元素链半环 \(T=\{z<d<u\}\)，其中 \(u\) 是乘法幺元、\(z\) 是零、\(d^2=z\)。在 \(L_2\times T\) 中去掉``较小颜色、单位层''这一点，剩余五点是子半环，依次为

\CaseDisplay{0=(\ell, z),\ b=(r, z),\ a=(\ell, d),\ c=(r, d),\ 1=(r, u).}

这给出上表并证明半环公理。\(\{0, b\}\) 和 \(\{b, c,1\}\) 分别给出两个因子，所以还有

\CaseDisplay{\mathsf V(Q)=\mathsf V(L_2)\vee\mathsf V(T).}

该等式本身不承担有限基推论；下面给出完整有限基。

\par\smallskip\noindent\textit{2. 四条及五条基.}\quad 
记 \(\Lambda\) 为

\CaseDisplay{\begin{aligned}
xyz&\approx xzy, &(N)\\
xy&\approx x^2y+xy^2, &(E)\\
xy&\leq x, &(C)\\
xy&\leq (xy)^2+yx. &(R)
\end{aligned}}

另记

\CaseDisplay{xy\leq y+xz.\quad\textup{(\ensuremath{D})}}

定理。 \(P\) 的完整恒等式基为 \(\Lambda\)，\(Q\) 的完整恒等式基为 \(\Lambda\cup\{D\}\)，均相对于通常加法幂等半环公理。反转每个乘法词，立即得到两个乘法对偶半环的完整基。

这些恒等式在 \(P\) 中的有效性可由上面的词值描述直接得到。(N) 保留首字母及各变量出现次数；(E) 至少有一个右项保留词值；(C) 是右乘不增加链值；(R) 在乘积幂等时由平方项吸收，在乘积为 \(a\) 时由倒序项吸收，在乘积为零时自动成立。

对 \(Q\)，其两个因子分别满足 \(\Lambda\cup\{D\}\)。在左零因子中 (D) 变为 \(x\leq y+x\)；在 \(T\) 中已有 \(xy\leq y\)。其余四式在左零因子中直接成立，在 \(T\) 中由交换、乘法收缩和 \(d^2=z\) 得到。故 \(Q\) 满足所列基。

\par\smallskip\noindent\textit{3. 完整吸收判据.}\quad 
为避免与半环 \(Q\) 混淆，以 \(F\) 表示非空多项式。令 \(L\) 为其一次项变量集，\(\mathcal R\) 为其非线性求和词族；对于内容集 \(C\)，置

\CaseDisplay{\mathcal R_C=\{u\in\mathcal R:c(u)\subseteq C\}.}

记 \(h(w)\) 为首字母，\(S(w)=\{x:m_x(w)=1\}\)。两者都满足

\CaseDisplay{x\leq F\iff x\in L.\quad\textup{(\ensuremath{1})}}

对非线性词 \(w\)，令 \(h=h(w)\)、\(C=c(w)\)，并记共同的保护条件

\CaseDisplay{L\cap C\ne\varnothing\quad\text{或}\quad
\forall x\in S(w)\;\exists u\in\mathcal R_C:\ m_x(u)\leq1.\quad\textup{(\ensuremath{2})}}

在 \(P\) 中，\(w\leq F\) 恰当且仅当 (2) 和下式成立：

\CaseDisplay{h\in L\quad\text{或}\quad
\exists u\in\mathcal R_C:\ h(u)=h.\quad\textup{(\ensuremath{3})}}

在 \(Q\) 中，\(w\leq F\) 恰当且仅当 (2) 以及以下两个条件成立：

\CaseDisplay{h\in L\ \text{或}\ \exists u\in\mathcal R:h(u)=h,\quad\textup{(\ensuremath{4})}}

\CaseDisplay{L\cap C\ne\varnothing\quad\text{或}\quad\mathcal R_C\ne\varnothing.\quad\textup{(\ensuremath{5})}}

\par\smallskip\noindent\textit{必要性.}\quad 
对 \(P\) 的一次项，令该变量为 \(t\)、其余变量为 \(0\)；没有相同一次项时，右侧至多为 \(f<t\)。对 \(Q\)，改令该变量为 \(c\)、其余变量为 \(0\)，右侧至多为 \(b<c\)。这证明 (1)。

若 (3) 失败，在 \(P\) 中令 \(h=f\)，令其余 \(C\) 内变量为 \(e\)、外部变量为 \(0\)。左词为 \(f\)，内部右词的首字母都给出 \(e\)，外部右词为零，右侧不超过 \(e\)。若 (2) 失败，令被重复保护失败的那个一次变量为 \(a\)、其他 \(C\) 内变量为 \(e\)、外部变量为 \(0\)，则左侧为 \(a\)、右侧为零。

对 \(Q\)，若 (4) 失败，令 \(h=c\)、其余所有变量为 \(0\)。左非线性词为 \(b\)，所有右项都以零开始，故右侧为零。若 (5) 失败，令 \(C\) 内变量均为 \(1\)、外部变量为 \(0\)。左词为 \(1\)，所有右词含外部变量，故右侧不超过 \(b\)。若 (2) 失败，令相应一次变量为 \(c\)、其他 \(C\) 内变量为 \(1\)、外部变量为 \(0\)。左词为 \(c\)，内部右词重复该变量而取 \(b\)，其余右词至多为 \(b\)，故吸收失败。

\par\smallskip\noindent\textit{4. 从有限基推出充分性.}\quad 
由 (E) 得 \(x^3=x^2\)。由 (N)，每个词在首字母之后的字母可任意交换。故单词可按首字母、内容和各变量次数在二处截断来整理。

若一个词有至少两个一次出现的变量，用 (N) 把两个这样的变量移到相邻位置，再用 (E) 展开。如果首字母本身一次出现，就把另一个一次变量放在其后；这一展开仍保持首字母。每一步减少一次出现变量的数目，最终把每个非线性词化为非空和，每个分项至多有一个一次出现变量。

记 \(H_h(C)\) 为首字母为 \(h\)、每个 \(C\) 内变量出现两次的词；记 \(K_{h, x}(C)\) 为首字母为 \(h\)、\(x\) 恰出现一次、其余变量出现两次的词。不同尾部排列由 (N) 等同。单变量非线性词用 \(H_h(\{h\})=h^2\)，不使用空积。

有两个基本推导。

\textbullet{} 补乘。 如果词 \(u\) 的首字母为 \(g\)，其截断变量次数不超过词 \(v\) 的相应次数，而 \(v\) 也以 \(g\) 开始，则由 (C) 补乘缺少的因子，再用 (N) 整理，得到 \(v\leq u\)。

\textbullet{} 改变首字母。 对同一内容和同一次数的两个词 \(w_h, w_g\)，若首字母分别为 \(h, g\)，则 (R) 给出

\CaseDisplay{w_h\leq H_h(C)+w_g.\quad\textup{(\ensuremath{6})}}

当 \(g\ne h\) 时，把 \(w_h\) 在某次出现的 \(g\) 前切成两个非空块，使用 \(uv\leq(uv)^2+vu\)，再用 (N) 整理；\(g=h\) 时不需要改变。

对 \(P\)，若 \(h\in L\)，直接由 (C) 得 \(w\leq h\leq F\)。否则 (3) 提供一个内部、首字母为 \(h\) 的右词；补乘得到

\CaseDisplay{H_h(C)\leq F.\quad\textup{(\ensuremath{7})}}

对 \(Q\)，同样先处理 \(h\in L\)。否则 (4) 提供 \(u=hs\in\mathcal R\)，其中 \(s\) 非空；(5) 提供一个内部右项 \(v\)，允许它是一次项。由 (D)

\CaseDisplay{hv\leq v+hs=v+u\leq F.}

补乘得到 (7)，即使 \(u\) 含外部变量，所得稳定左词也只含 \(C\)。

剩余规范左词为 \(K_{h, x}(C)\)。若 \(L\cap C\ne\varnothing\)，选其中的 \(g\)。以 \(g\) 开始的同次数词由 (C) 不超过 \(g\)，再用 (6)、(7) 得到所需吸收。若没有内部一次项，(2) 给出不重复 \(x\) 的内部右词 \(v\)，令 \(g=h(v)\)。补乘把它扩充为 \(K_{g, x}(C)\leq v\leq F\)，再用 (6)、(7) 得 \(K_{h, x}(C)\leq F\)。

这证明全部判据的充分性，且每一步都是所列有限基的等式推导。最后把任意恒等式分解为双向逐词吸收，得到完整性。乘法对偶通过反转词直接继承结论。

\subsection{A finite basis from doubled blocks and one simple-variable recovery}\label{case:126}
Source 4591. Signature (+, multiplication), without constants, relative to the additively idempotent semiring axioms. A simple variable has exactly one occurrence. A one-letter word has length one.

\par\smallskip\noindent\textit{1. Algebra and exact semantics.}\quad 
Addition is maximum in \(3<2<1<0<4\). Multiplication, in 0,\allowbreak{}1,\allowbreak{}2,\allowbreak{}3,\allowbreak{}4 order, is

{\small\par\smallskip
\CaseLine{\ttfamily 0 0 2 3 0}
\CaseLine{\ttfamily 1 1 2 3 1}
\CaseLine{\ttfamily 3 3 3 3 2}
\CaseLine{\ttfamily 3 3 3 3 3}
\CaseLine{\ttfamily 0 0 2 3 4}
\smallskip}

Put K=\{0,\allowbreak{}1,\allowbreak{}4\}. A word has value 3 if it contains 3, or has at least two occurrences valued 2, or has one 2 followed by a factor other than 4. It has value 2 if it has exactly one 2 and all later factors are 4. With all factors in K, it is 1 when its head is 1,\allowbreak{}4 when all factors are 4, and 0 otherwise. This formula includes one-letter words and follows by induction. It proves associativity: 3 is a zero; a2 survives precisely when it is unique and followed entirely by 4; on K the first non4 determines the product, with a noninitial1 changing an initial 4 to 0. Direct rows/columns in the displayed chain are increasing, so both distributive laws hold.

For content C of target w, let U be the original right words of P contained in C. For a simple variable x let \(S_{w}\)(x) be its suffix content after its unique occurrence.

Criterion. \(w\leq P\) holds exactly when:

\textbullet{} H: U has a word with the same head as w;

\textbullet{} for each simple x of w, S: U has a word omitting x, or a word containing x exactly once with its suffix after x contained in \(S_{w}\)(x).

For a one-letter target these conditions are equivalent to the presence of that literal right summand.

Necessity. If H fails, assign the target head 0, other variables of C the value 1, and outsiders 3. The target is 0 and all right words are 1 or 3. If S fails at x, assign \(x=2\), its target suffix variables 4, the other target variables 0, and outsiders 3. The target is 2. A contained right word must either repeat x or have a non4 after its unique x; any noncontained word includes3. Thus \(P=3\).

Sufficiency. Target 3 is least. Target 4 has all C valued 4, and any H witness is 4. For target 0 or 1 all C lie in K. An H witness headed 0 gives 0, one headed 1 gives 1, and one headed 4 gives 0 or 4; these dominate the corresponding target. A target 2 has a unique 2 at a simple x and only 4 after it. An S witness omitting x has all values in K and is at least 1; one containing x once has only 4 after x and is 2. This proves the criterion for arbitrary variables.

\par\smallskip\noindent\textit{2. Six sound identities and normal words.}\quad 
The following identities follow from the formula:

\CaseDisplay{\begin{aligned}
xyzy&\approx xzy^2,&xyx&\approx xyx^2,&x^3&\approx x^2,\quad\textup{(\ensuremath{W1--W3})}\\
xy^2z^2&\approx xz^2y^2,\quad\textup{(\ensuremath{W4})}\\
xy&\le x,\quad\textup{(\ensuremath{D})}\\
xy&\approx x^2y+xy^2.\quad\textup{(\ensuremath{E})}
\end{aligned}}

Here and below letters in substitutions are replaced by nonempty words. D in particular gives \(x^{2}\leq x\).

A normal representative retains the literal first position h. Every simple variable after h stays in its original relative order. Every repeated variable has a doubled block at its last occurrence after h; this includes a repeated head, whose initial occurrence is retained. Sort doubled blocks within each interval between the simple variables. Thus a repeated head can have total multiplicity three in this representative, and no one-letter word is introduced.

To derive this, W1 moves a noninitial earlier occurrence onto its later occurrence; its prefix is nonempty, and W3 handles adjacent copies. W2 makes the last copy of a repeated head into a doubled block when it is separated from the initial head. For an immediately repeated head, reverse W3 makes the required head plus doubled block. Further noninitial copies are collected with W1 and W3. W4 sorts adjacent doubled blocks, with their nonempty preceding prefix as its first substituted word. No simple variable is crossed. These operations prove the asserted normal forms equationally.

In particular all words with no simple variables normalize to the same saturated word \(s_{h}\)(C): h followed by the sorted squares of every variable of C. This includes singleton content (hhh). Words with at most one simple variable are determined by their head, content, and the position of that simple variable relative to the repeated blocks, equivalently its suffix content.

\par\smallskip\noindent\textit{3. Saturated guard and individual recoveries.}\quad 
Assume w has length at least two and satisfies the criterion. Choose an original H witness u. By D, \(u^{2}\leq u\leq P\); \(u^{2}\) has no simple variables. Append squares of missing target variables using D with a nonempty prefix. Normalize to obtain \(s_{h}(C)\leq P\).

For each simple target variable x choose an original S witness \(u_{x}\). In \(u_{x}\) double every occurrence of every variable other than x, using \(x^{2}\leq x\) in its appropriate context. Let \(v_{x}\) be the normal result; then \(v_{x}\leq u_{x}\leq P\). It either omits x or has x as its sole simple variable, and its suffix after x is still contained in \(S_{w}\)(x).

Let \(r_{x}\) be the normal form obtained from w by doubling every position other than its single x. The criterion proves the local absorption

\CaseDisplay{r_x\le s_h(C)+v_x.\quad\textup{(\ensuremath{R})}}

The first right word supplies H; the second supplies the only simple-variable condition. Both are already absorbed by P, so deriving(R) gives \(r_{x}\leq P\) by transitivity.

\par\smallskip\noindent\textit{4. Nine-variable compression of each recovery.}\quad 
In(R) reserve the target head h, the head of \(v_{x}\), and x individually: at most three variables. Every other variable occurs as a doubled block in \(r_{x}\) and \(s_{h}\)(C). In \(v_{x}\) it is absent or a doubled block. Group these variables by their before/after-x region in \(r_{x}\) and their absent/before/after-x region in \(v_{x}\): at most six groups. If x is absent from \(v_{x}\), treat its entire tail as the before region.

Replace each group by one doubled block in each region where it occurs; do not sum the multiplicities of its members. Keep reserved variables and their occurrences. A block may be positioned within its interval by W4, so normalizing yields at most nine variables and length at most nineteen, including a retained repeated head. The criterion remains true: head equality and content containment are preserved, and x is still the only possible simple variable. Its suffix containment is preserved because grouping retained the two relative regions. No new simple variable is created.

For the inverse substitution replace a group variable by a nonempty product of its original members. Its square normalizes to the corresponding sorted doubled blocks using W1--W4 with a nonempty preceding prefix. All members belonged to the same region in each word, so no block passes across x. Reserved variables are replaced individually. The normal forms of all three original words are therefore restored. This proves reversible compression of(R).

\par\smallskip\noindent\textit{5. Reassembling the target and completeness.}\quad 
By induction on length, E and W1--W4 expand any word \(w=a_{1}\cdots a_{n}\) of length at least two into the sum of the n words in which one position remains \(a_{i}\) and every other position is squared. For the induction step apply E to the prefix and the final letter. The squared prefix has no simple variables and normalizes to the same word as the product of its individually squared positions. Expand the other prefix term inductively. This derives the stated sum.

If the unsquared position contains a repeated target variable, its summand still has no simple variables and is \(s_{h}\)(C). Otherwise the summand is \(r_{x}\) for the corresponding simple x. Sections 3--4 derive absorption of every summand by P, so E gives \(w\leq P\). A one-letter target already occurs as a literal summand.

Take the six displayed identities and all valid normal absorptions over a fixed nine-letter alphabet, length at most nineteen and at most two right words. The criterion explicitly determines this finite family; no small or minimal enumeration is claimed. Every recovery(R) compresses into this family and lifts back, so the preceding construction derives every valid word absorption, without a bound on original variables, lengths or summands. Distributivity and absorption of the summands on both sides give every identity. Word reversal gives a finite complete basis of the opposite. Computational checks support this unrestricted proof; catalogue transfers require both complete HSP directions.

\subsection{Content localization and affine parity: finite bases for three zero extensions}\label{case:127}
\par\smallskip\noindent\textit{Abstract.}\quad 
We give five-, five-, and six-identity bases for three five-element additively idempotent semirings obtained by adjoining an additive bottom zero to extensions of the flat two-element group semiring. A target word localizes the right polynomial to its own content. Within this localization, an affine parity condition and a coverage condition characterize absorption. Explicit guarded identities supply the required stable representatives of linear summands. The proof works for arbitrarily many variables and summands.

\par\smallskip\noindent\textit{1. The algebras and their bases.}\quad 
The signature is \((+,\cdot)\), without constants. All bases below are relative to the additively idempotent semiring axioms. The notation \(u\le v\) abbreviates the identity \(u+v\approx v\); displayed element names are not term constants.

Let \(E=\{e, g, d,0,\omega\}\). In each algebra \(0\) is the least and \(\omega\) the greatest element of the additive order, and addition is join. The only further order relation between distinct members of \(\{e, g, d\}\) is the one specified below:

\noindent Algebra\enspace\textbar\enspace Additional order relation\enspace\textbar\enspace Image of \(d\) under \(\rho\)\enspace\textbar\enspace Catalog source\par
\noindent \(A\)\enspace\textbar\enspace \(d<g\)\enspace\textbar\enspace \(g\)\enspace\textbar\enspace 10590\par
\noindent \(B\)\enspace\textbar\enspace \(e<d\)\enspace\textbar\enspace \(\omega\)\enspace\textbar\enspace 13206\par
\noindent \(D\)\enspace\textbar\enspace \(g<d\)\enspace\textbar\enspace \(\omega\)\enspace\textbar\enspace 13239\par
The map \(\rho\) fixes \(e, g,0,\omega\). On \(F=\{0, e, g,\omega\}\), multiply \(e, g\) as the group of order two, let \(0\) annihilate every element, and put \(\omega u=u\omega=\omega\) for \(u\ne0\). Multiplication on \(E\) is \(uv=\rho(u)\rho(v)\), evaluated in \(F\).

Theorem 1. Each algebra has the four common identities

\CaseDisplay{\begin{array}{rll}
\mathrm{(C)}&xy\approx yx,\\
\mathrm{(N)}&x^3y\approx xy,\\
\mathrm{(U)}&x\le x^3,\\
\mathrm{(O)}&xyz\le x^3+y^3+z^3,
\end{array}}

together with the following identities, as a complete basis:

\CaseDisplay{\begin{array}{rll}
A:&\mathrm{(G)}\quad x^3y^2z^2\le x+yz;\\
B:&\mathrm{(I)}\quad x^3y^2\le x+xy;\\
D:&\mathrm{(I)}\quad x^3y^2\le x+xy,\qquad
   \mathrm{(Q)}\quad x^3y^2\le x+y^2.
\end{array}}

Thus the respective bases have five, five, and six identities. No minimality claim is made.

The algebra \(F\) is an additively idempotent semiring: its nonzero part is the flat group semiring, and adjoining a new additive bottom annihilator preserves associativity and distributivity. In each of the three cases \(\rho\) is a join-preserving retraction onto \(F\). Because all products lie in \(F\), associativity and distributivity on \(E\) follow. Every word \(w\) of length at least two has the value

\CaseDisplay{w(\alpha)=\prod_{x\text{ occurring in }w}\rho(\alpha(x)),\quad\textup{(\ensuremath{1})}}

where occurrences, including repetitions, are multiplied. This value is \(0\) if any variable in its content is assigned \(0\), is \(\omega\) if no such variable is \(0\) and some image is \(\omega\), and otherwise is the group value determined by exponent parity.

Formula (1) verifies (C),(N), and \(x^3=\rho(x)\) verifies (U). To check (O), an assignment containing \(0\) makes its left side \(0\). Otherwise an image \(\omega\) makes the right side top; with group images only, either they all agree, when the odd product equals that common image, or the right join is top. For (G) in \(A\), again a zero assignment is immediate. A top image on the left makes the right top. Otherwise the left value is \(\rho(x)\). When \(x\ne d\) it is already the right summand \(x\); when \(x=d\), the stable right summand \(yz\) is \(e\) or \(g\), and \(d+e=\omega\), \(d+g=g\) both dominate \(\rho(d)=g\).

For (I) in \(B, D\), zeros are immediate. Without zeros, the occurrence of \(d\) or \(\omega\) in \(x, y\) makes \(xy=\omega\); for group values the left side equals \(x\). Finally, (Q) in \(D\) is immediate for zeros and for \(y\in\{d,\omega\}\). For group \(y\), its square is \(e\); group \(x\) is already on the right, \(x=\omega\) is immediate, and \(x=d\) gives \(d+e=\omega\). This proves soundness of all displayed laws.

\par\smallskip\noindent\textit{2. Normal words and a fixed-content closure lemma.}\quad 
Call a word stable if its length is at least two. Write \(c(w)\) for its content and \(\pi(w)\in\mathbb F_2^X\) for its exponent-parity vector on a finite variable set \(X\). The laws (C),(N) identify two stable words whenever they have the same content and parity. Indeed, for content of size at least two, (N) removes pairs from any exponent at least three, retaining a nonempty complementary word. The resulting exponents are one or two. For singleton content, the stable normal forms are \(x^2, x^3\); reductions start from \(x^4\approx x^2\). No empty word is substituted. In particular, every stable word \(r\) satisfies \(r^3\approx r\).

For a nonempty family \(R\) of words let

\CaseDisplay{H(R)=\operatorname{aff}_{\mathbb F_2}\{\pi(r):r\in R\}.}

Over \(\mathbb F_2\), its elements are precisely sums of an odd number of the listed vectors, allowing repetitions.

Lemma 2. Suppose that \(R\) is a nonempty finite family of stable words and \(\bigcup_{r\in R}c(r)=C\). If \(w\) is stable, \(c(w)=C\), and \(\pi(w)\in H(R)\), then the four common laws derive \(w\le\sum R\).

Proof. Choose an odd product of members of \(R\) with parity \(\pi(w)\). Insert two further copies of every member of \(R\). The product is still odd and has exactly content \(C\). Iterating (O), with stable words substituted, bounds every odd product by the sum of its factors: cubes of stable words normalize to those words. Normalization then identifies the product with \(w\). The case of one factor is reflexive. This proves the lemma without changing the content of the target. For a linear target \(x\), the same argument may instead be applied to \(x^3\), followed by (U). \(\square\)

\par\smallskip\noindent\textit{3. Localization and the exact absorption criteria.}\quad 
Let \(P\) be a nonempty finite sum of words. For a target \(w\) put \(C=c(w)\), and let \(P_C\) consist of precisely the summands \(u\) of \(P\) for which \(c(u)\subseteq C\). If this family is empty, no absorption \(w\le P\) is valid: assign every variable of \(C\) the group identity \(e\) and every outside variable \(0\). The target has value \(e\) and the right side \(0\). We never use an empty sum as a term.

If \(P_C\) is nonempty, then

\CaseDisplay{w\le P\text{ holds }\quad\Longleftrightarrow\quad
  w\le P_C\text{ holds}.\quad\textup{(\ensuremath{2})}}

For necessity, extend any assignment on \(C\) by \(0\) outside \(C\); exactly the outside summands vanish. Sufficiency follows from \(P_C\le P\).

Write \(L\) for the linear summands of \(P_C\), regarded as variables, and \(R\) for its stable summands. Set

\CaseDisplay{U=\bigcup_{u\in P_C}c(u),\qquad N=\bigcup_{r\in R}c(r),\qquad
H=H(P_C).}

If \(R\) is empty take \(N=\varnothing\). Always \(N\subseteq U\subseteq C\).

Theorem 3. If \(w=x\in L\), absorption is immediate. For every other target \(w\), the following conditions are necessary and sufficient:

\noindent Algebra\enspace\textbar\enspace Conditions for \(w\le P\)\par
\noindent \(A\)\enspace\textbar\enspace \(R\ne\varnothing\), \(U=C\), and \(\pi(w)\in H\)\par
\noindent \(B\)\enspace\textbar\enspace \(R\ne\varnothing\), \(N=C\), and \(\pi(w)\in H\)\par
\noindent \(D\)\enspace\textbar\enspace \(R\ne\varnothing\), \(U=C\), \(\pi(w)\in H\), and either \(N=C\) or \(0\in H\)\par
Every absorption satisfying the relevant row is derivable from the basis in Theorem 1.

Necessity. Work with the nonempty localization in (2). A missing variable of \(C\setminus U\) can be assigned \(\omega\), with every other variable in \(C\) assigned \(e\). The target becomes \(\omega\) and the right side \(e\). Thus \(U=C\). If \(\pi(w)\notin H\), a linear functional is constant on \(H\) and has a different value on \(\pi(w)\). Assign \(e\) or \(g\) according to that functional. All right words have one group value and the target the other. Those values are incomparable, proving the parity restriction.

For a stable target, if \(R\) were empty, assign every variable of \(C\) the element \(d\). The right side is \(d\). In \(A\) the target is \(e\) or \(g\), neither of which is below \(d\); in \(B, D\) the target is \(\omega>d\). For an unlisted linear target \(x\), localization leaves only stable powers of \(x\), so nonemptiness of \(P_C\) already implies \(R\ne\varnothing\) and \(N=C\).

In \(B\), if a stable target has \(x\in C\setminus N\), assign \(x=d\) and all other variables \(e\). Since \(U=C\), the variable \(x\) is a linear summand. The right side is \(d\), whereas the target is \(\omega\). Hence \(N=C\).

In \(D\), suppose \(x\in C\setminus N\) and \(0\notin H\). There is a linear functional taking the value one on all of \(H\). The associated group assignment makes every right summand \(g\). In particular it assigns \(x=g\), since \(x\in L\). Change only this assignment to \(d\). Stable right summands do not contain \(x\), so they remain \(g\); the right join is now \(d\). A stable target containing \(x\) becomes \(\omega\). This is a counterevaluation. The unlisted linear case already has \(N=C\), as observed above. These arguments establish every restriction in the table.

Sufficiency for \(A\). Fix \(q\in R\) and factor it as \(q=uv\), with both \(u, v\) nonempty. For each \(x\in L\), (G) gives

\CaseDisplay{x^3q^2=x^3u^2v^2\le x+uv=x+q\le P_C.\quad\textup{(\ensuremath{3})}}

Adjoin these stable words to \(R\). Their parity vectors are precisely the unit vectors of the linear summands; consequently the new family has affine hull \(H\) and total content \(U=C\). Lemma 2 derives the target if it is stable, and derives \(x^3\) followed by (U) if it is an unlisted linear word \(x\). Since every member of the family is below \(P_C\), this proves the required absorption.

Sufficiency when \(N=C\) in \(B, D\). For each \(x\in L\), select a stable right word containing \(x\). By commutativity write it as \(xv\) with \(v\) nonempty. Law (I) supplies

\CaseDisplay{x^3v^2\le x+xv\le P_C.\quad\textup{(\ensuremath{4})}}

Its parity is the unit vector of \(x\) and its content is \(c(xv)\). The resulting stable family, together with \(R\), has affine hull \(H\) and content \(N=C\). Apply Lemma 2, and (U) for an unlisted linear target.

The remaining case in \(D\). Assume \(U=C\) and \(0\in H\). Use (4) only for \(x\in L\cap N\), obtaining a nonempty stable family \(R_0\) of total content \(N\). Its affine hull is generated by the vectors of \(R\) and the units for \(L\cap N\). Each remaining linear variable \(j\in J=L\setminus N\) occurs in no other generator of \(H\). In an odd representation of the zero vector by the generators of \(H\), the coefficient of every such unit must therefore be even. Removing those even repetitions leaves an odd zero representation using \(R_0\). Thus \(0\in H(R_0)\).

The set \(N\) is nonempty since \(R\) contains a nonempty word. Let \(q\) list each member of \(N\) once. Lemma 2 gives \(q^2\le P_C\). For every \(j\in J\), (Q), with the nonempty word \(q\) substituted for \(y\), gives

\CaseDisplay{j^3q^2\le j+q^2\le P_C.\quad\textup{(\ensuremath{5})}}

These words supply the missing unit vectors and enlarge the content to \(U=C\). The affine hull is now exactly \(H\): the added zero vector was already in \(H(R_0)\). Lemma 2 again derives the desired target, using (U) if necessary. This completes the proof of Theorem 3.

\par\smallskip\noindent\textit{4. Completeness and catalog consequences.}\quad 
Distribute both sides of an arbitrary valid identity into nonempty sums of words. Each left summand is absorbed by the right polynomial, and conversely. Theorem 3 derives every such absorption from the respective finite basis. Adding the derivations proves the original identity. This proves Theorem 1 for any finite variable set, without a bound on its size or on the number of summands.

The named-to-catalog permutations, in the order \((e, g, d,0,\omega)\), are \([3,\allowbreak{}0,\allowbreak{}1,\allowbreak{}2,\allowbreak{}4]\) for 10590, \([1,\allowbreak{}3,\allowbreak{}0,\allowbreak{}2,\allowbreak{}4]\) for 13206, and \([3,\allowbreak{}1,\allowbreak{}0,\allowbreak{}2,\allowbreak{}4]\) for 13239. The two operations determine the maps. For equality of varieties use Proposition~\ref{prop:lower-fibres}. In particular, 4640 and 7788 generate the same variety as 10590 and have its five-identity basis. The five catalog types thereby classified are 4640, 7788, 10590, 13206, and 13239. Satisfaction of another source's laws alone is not used to transfer finite basability to a proper subvariety.

The source manuscript for this case credits Aifa Wang and Lili Wang, School of Mathematical Sciences, Chongqing University of Technology.

The source manuscript for this case acknowledges support from the National Natural Science Foundation of China (12371024), the Science and Technology Research Program of Chongqing Municipal Education Commission (KJZD-K2024011102) and the Chongqing Natural Science Foundation Innovation and Development Joint Fund (Municipal Education Commission) (CSTB2025NSCQ-LZX0067).

\subsection{一次变量的保护与四个交换半环的完整基}\label{case:128}
\par\smallskip\noindent\textit{1. 生成半环.}\quad 
底集记为 \(\{0, a, b, e, f\}\)。除下述吸收型例外外，乘法由交换核心 \(\{0, a, e, f\}\) 及 \(b\mapsto e\) 的乘法收缩定义：\(0\) 是乘法零，\(e\) 是核心的幺元，\(a^2=0\)、\(f^2=f\)、\(af=a\)。特别地，\(b^2=e\)、\(ab=a\)、\(bf=f\)。

加法取以下偏序的并；表中只列覆盖关系。

\noindent 半环\enspace\textbar\enspace 目录来源\enspace\textbar\enspace 加法序\enspace\textbar\enspace 乘法例外\enspace\textbar\enspace 到目录标号的映射\par
\noindent \(D\)\enspace\textbar\enspace 4641\enspace\textbar\enspace \(0<b<e<f,\ 0<a<f\)\enspace\textbar\enspace 无\enspace\textbar\enspace \((2,\allowbreak{}3,\allowbreak{}1,\allowbreak{}0,\allowbreak{}4)\)\par
\noindent \(D_H\)\enspace\textbar\enspace 6160\enspace\textbar\enspace \(0<e<f,\ 0<a<f,\ b<e\)\enspace\textbar\enspace 无\enspace\textbar\enspace \((1,\allowbreak{}3,\allowbreak{}2,\allowbreak{}0,\allowbreak{}4)\)\par
\noindent \(D_U\)\enspace\textbar\enspace 4708\enspace\textbar\enspace \(0<a<e<f,\ 0<b<e\)\enspace\textbar\enspace 无\enspace\textbar\enspace \((3,\allowbreak{}2,\allowbreak{}1,\allowbreak{}0,\allowbreak{}4)\)\par
\noindent \(D_{UG}\)\enspace\textbar\enspace 4712\enspace\textbar\enspace \(0<a<e<f,\ 0<b<e\)\enspace\textbar\enspace \(f\) 是乘法吸收元\enspace\textbar\enspace \((3,\allowbreak{}2,\allowbreak{}1,\allowbreak{}0,\allowbreak{}4)\)\par
没有列出的比较关系仅由传递性产生。前三个半环的乘法来自四元素交换幺半群的上述膨胀，第四个来自三元素核心 \(\{0, a, e\}\) 添吸收元 \(f\) 后的同样膨胀。乘法结合，各平移保持指定的并，故四者均为加法幂等半环。

\par\smallskip\noindent\textit{2. 有限基.}\quad 
令 \(\Delta\) 包含以下六式：

\CaseDisplay{\begin{aligned}
xy&\approx yx, &(C)\\
xy&\approx x^2y+xy^2, &(E)\\
x^2y^2&\leq x^2+y, &(J)\\
xy^2&\leq x+(xy)^2, &(R)\\
xy^2&\leq x+y^2, &(F)\\
xy^2z^2&\leq xy^2+x^2z^2. &(B)
\end{aligned}}

另记

\CaseDisplay{x^2\leq x+yz\quad(H),\qquad
xy\leq x^2+y^2\quad(U),\qquad
x^2y\leq x^2+y^2z^2\quad(G).}

定理。 相对于加法幂等半环公理，四个生成半环的完整基依次为

\CaseDisplay{\operatorname{Id}(D):\Delta,\quad
\operatorname{Id}(D_H):\Delta+H,\quad
\operatorname{Id}(D_U):\Delta+U,\quad
\operatorname{Id}(D_{UG}):\Delta+U+G.}

不声称条数最小。

\par\smallskip\noindent\textit{3. 必要的组合条件.}\quad 
把非空多项式 \(P\) 的一次项变量集记为 \(L\)，非线性求和词族记为 \(\mathcal R\)。对非线性左词 \(w\)，令

\CaseDisplay{C=c(w),\quad S(w)=\{x:m_x(w)=1\},\quad
\mathcal R_C=\{u\in\mathcal R:c(u)\subseteq C\}.}

四个半环均满足 \(x\leq P\iff x\in L\)。非线性吸收式 \(w\leq P\) 则恰由下列条件刻画。

1. 锚点： 通常要求 \(\mathcal R_C\ne\varnothing\)；对 \(D_H\)，另允许 \(L\cap C\ne\varnothing\) 且 \(\mathcal R\ne\varnothing\)。

2. 覆盖： 通常要求 \(C\subseteq L\cup\bigcup_{u\in\mathcal R_C}c(u)\)；对 \(D_{UG}\)，允许用全部 \(\mathcal R\) 代替 \(\mathcal R_C\)。

3. 保护： 对每个 \(x\in S(w)\)，要求

\CaseDisplay{x\in L\quad\text{或}\quad
   \exists u\in\mathcal R_C:\ m_x(u)=1.
   \quad\textup{(\ensuremath{1})}}

对 \(D_U, D_{UG}\)，可将 \(m_x(u)=1\) 放宽为 \(m_x(u)\leq1\)。

这里保护一次变量的线性右项必须是该变量本身；另一个内部一次项一般不足以保护它。

必要性有直接的失败赋值。若左侧为 \(x\notin L\)，令 \(x=a\)、其余变量为 \(0\)，所有非线性右词为 \(0\)。若没有通常锚点，令 \(C\) 内变量为 \(b\)、外部变量为 \(0\)，左词为 \(e\)，而右侧不高于 \(b\)。对 \(D_H\)，若没有内部词及内部一次项，同一赋值的右侧为 \(0\)；若没有非线性右项，令所有变量为 \(b\)，右侧仍为 \(b\)。

若内部覆盖失败，选未覆盖的 \(x\in C\)，令它为 \(f\)、其他内部变量为 \(e\)、外部变量为 \(0\)，得到左侧 \(f\) 而右侧不高于 \(e\)。在 \(D_{UG}\) 中，若全局覆盖失败，把该变量赋为 \(f\)、其余所有变量赋为 \(e\)，同样失败。

对于 \(D, D_H\)，若 (1) 对 \(x\) 失败，令 \(x=a\)、其他内部变量为 \(e\)、外部变量为 \(0\)。右侧没有一次项 \(x\)，每个内部右词不是省略 \(x\)，就是重复 \(x\)，因此右侧不高于 \(e\)，而左词为 \(a\nleq e\)。对于 \(D_U, D_{UG}\)，令其他内部变量改为 \(b\)；保护失败意味着所有内部非线性右词重复 \(x\)，所以右側至多为 \(b\)，左词仍为 \(a\nleq b\)。

\par\smallskip\noindent\textit{4. 所列基的有效性.}\quad 
非线性词的值均在 \(\{0, a, e, f\}\) 中。暂不计吸收型的 \(f\)：出现 \(0\) 或至少两次 \(a\) 时值为 \(0\)，恰出现一次 \(a\) 时值为 \(a\)，没有这些因子时值为 \(e\) 或 \(f\)。在吸收型中，只要出现 \(f\)，结果即为 \(f\)。由此立即得到 (C)、(E)；在 (E) 中至少一个求和项保留原乘积值。

平方值只可能为 \(0, e, f\)。(J) 的左侧为 \(0\) 时，有平方右项 \(x^2\) 保证吸收；为 \(e\) 时，或者 \(x^2=e\)，或者左侧不可能为 \(e\)；为 \(f\) 时，某个平方因子为 \(f\)，由相应右项吸收。在吸收型中，若仅 \(y^2=f\)，则 \(y=f\)，故仍成立。(R) 左侧为 \(a\) 时必有 \(x=a\)，由一次项吸收；稳定值由平方项吸收。(F) 在左侧为 \(a\) 时同样用 \(x=a\)，单位层值由 \(x\) 或 \(y^2\) 吸收。

(B) 若左侧为 \(a\)，第一右项保留 \(a\)；若左侧为单位层值，则 \(x^2, y^2, z^2\) 在单位层，右项之一保留最高层。在吸收型中，产生 \(f\) 的因子也在至少一个右项中。左侧为 \(0\) 时，右侧有非线性词且它们的值均不小于 \(0\)。

(H) 唯一需要注意的是 \(x=b\)，此时 \(b^2=e\)，而任一非线性词的值属于 \(\{0, a, e, f\}\)，在 \(D_H\) 中与 \(b\) 的和不小于 \(e\)。(U) 对 \(D_U, D_{UG}\) 在乘积为 \(a\) 时由另一因子的平方值 \(e\) 或 \(f\) 吸收，其余情况直接成立。(G) 对 \(D_{UG}\) 的左侧为 \(a\) 时已有 \(x^2\geq e>a\)；为 \(e\) 时已有 \(x^2=e\)；为 \(f\) 时或 \(x^2=f\)，或 \(y=f\)，后者使 \(y^2z^2=f\)。左侧为 \(0\) 时自动成立。

\par\smallskip\noindent\textit{5. 任意变量数的完整性.}\quad 
由 (E) 得 \(x^3=x^2\)。反复对不同的一次出现变量使用 (E)，每个非线性词化为非空和，分项均为

\CaseDisplay{H(C)=\prod_{y\in C}y^2,
\qquad xH(C\setminus\{x\})\quad(|C|\geq2).}

后一类的 \(x\) 必须在原词中恰出现一次。每次展开保持内容；若原词恰含一次 \(x\)，至少有一个规范分项仍恰含一次 \(x\)。单变量非线性词用 \(H(\{x\})=x^2\)，不使用空积。

由 (E)，对每个非线性词 \(u\) 有 \(u^2\leq u\)。由 (J) 得稳定块的合并式

\CaseDisplay{H(S)H(T)\leq H(S)+H(T),\qquad
H(S)y^2\leq H(S)+y.\quad\textup{(\ensuremath{2})}}

从内部非线性右词的平方开始，或在 \(D_H\) 中由 (H) 和一个内部一次项、一个任意非线性右项取得平方锚点，再用 (2) 逐一加入内部右词和内部一次项的内容。内部覆盖保证最终得到

\CaseDisplay{H(C)\leq P.\quad\textup{(\ensuremath{3})}}

在带 (G) 的情形，已有 \(H(S)\leq P\) 且某个右词含 \(z\in C\) 时，将该词写为 \(zt\)，用 (G) 及 \(H(S)z^2\leq H(S)z\) 得

\CaseDisplay{H(S)z^2\leq H(S)+z^2t^2\leq P.}

单变量非线性右词可先约化为 \(z^2\) 并取 \(t=z\)。这样由全局覆盖也得到 (3)。

剩余左规范词为 \(xH(C\setminus\{x\})\)。若 \(x\in L\)，(R) 与 (3) 直接给出所需吸收。否则，在严格保护情形，从一个恰含一次 \(x\) 的内部右词选出规范分项 \(xH(S)\leq P\)；这里 \(S\ne\varnothing\)。

在允许省略 \(x\) 的情形，若保护右词不含 \(x\)，先平方得到不含 \(x\) 的非空 \(H(S)\leq P\)。不带 (G) 时，内部覆盖提供含 \(x\) 的非线性右词，其平方为 \(x^2H(T)\)；由 (U) 得

\CaseDisplay{xH(S\cup T)\leq H(S)+x^2H(T)\leq P.}

若 \(T\) 为空，直接在 (U) 中代入 \(H(S),x\)。带 (G) 时，全局覆盖提供右词 \(xt\)，由 (G) 得 \(xH(S)\leq H(S)+x^2t^2\leq P\)。

由此均取得初始规范词 \(r=xH(S)\leq P\)。不带 (G) 时，逐一合并内部右词的内容：

\textbullet{} 右词不含 \(x\) 时，用 (F) 得 \(rH(T)\leq r+H(T)\leq P\)。

\textbullet{} 右词的规范分项为 \(xH(T)\) 时，将稳定块合并式 (2) 乘以 \(x\)。

\textbullet{} 右词重复 \(x\) 时，其平方为 \(x^2H(T)\)，用 (B) 得 \(xH(S\cup T)\leq r+x^2H(T)\leq P\)。若 \(T\) 为空，无须加入新变量。

\textbullet{} 内部一次项 \(y\) 可由 (F) 和 (J) 合并：

\CaseDisplay{ry^2\leq r+(ry)^2=r+r^2y^2\leq r+r^2+y\leq P,}

其中第一步用 (R)，最后一步用 \(r^2\leq r\)。

最后一项确保在 \(x\notin L\) 时，其他内部一次项仍可承担内容覆盖，但不被误当作 \(x\) 的保护见证。

带 (G) 时，若新变量 \(z\) 由任意非线性右词 \(zt\) 覆盖，使用 (R)、(G) 得

\CaseDisplay{rz^2\leq r+r^2z^2\leq r+r^2+z^2t^2\leq P.}

由一次项覆盖的新变量仍用上一段的 (J) 推导。逐一加入后便得到目标左规范词。

因此第 3 节每个必要条件满足的吸收式都能由相应基推出。这同时证明组合判据的充分性与有限基的完整性。多项式恒等式分解为双向逐词吸收即可。

\par\smallskip\noindent\textit{6. 证据范围.}\quad 
任意变量数的结论由上述推导承担。有限计算只核对具体运算表、基的有效性、有限范围的判据交叉检查以及双向完整词函数包含。目录应用在这些独立核验完成后另行合并。

\subsection{Sixteen identities for 4649 with nonlinear seeds and literal prefix guards}\label{case:129}
Signature (+, multiplication), without constants, relative to the additively idempotent semiring axioms. The proof below imposes no bound on variables, word length or polynomial size.

\par\smallskip\noindent\textit{1. Tables, values and exact original-witness criterion.}\quad 
Addition has rows 00044,\allowbreak{}01144,\allowbreak{}01234,\allowbreak{}44334,\allowbreak{}44444. Its order is the chain \(2<1<0<4\) together with \(2<3<4\); both 1+3 and 0+3 are 4. Multiplication has rows 00234,\allowbreak{}00234,\allowbreak{}22222,\allowbreak{}44444,\allowbreak{}44444. A literal retains its value. For a nonlinear word, scan from the left for the first factor outside \{0,\allowbreak{}1\}. A first 2 gives 2, a first 4 gives 4, and a first 3 gives 3 when final and 4 otherwise. If every factor is 0 or 1, the nonlinear word has value 0. Induction on the displayed table proves this formula.

Let Z(S) mean that an original right word has content contained in S. At a first target occurrence x, let S be the preceding target content. An E(x, S) witness is an original word containing x with first-x prefix contained in S. It is marked if a position follows that first x, unmarked otherwise. A nonlinear absorption \(w\leq P\) is valid exactly when:

(N) An original NONLINEAR word is contained in c(w).

(E) Every first target occurrence has an E-witness.

(M) At each marked target occurrence there is a marked E-witness, or an E-witness together with Z(S).

The Z witness may be literal. At the simple final variable, E alone suffices. Literal targets x are absorbed exactly when an original word has pure content \{x\}.

For necessity of N, set target \(variables=1\), \(outsiders=2\). The target is 0; without N all right values are 1 or 2. For failed E, set \(x=4\), \(S=1\), everything \(else=2\). The target is 4, every right word is 0,\allowbreak{}1 or 2. For failed M with E present, set \(x=3\), \(S=1\), everything \(else=2\). Only unmarked E-witnesses can give 3, no right word gives 4, and the absence of Z(S) rules out values 0 or 1; thus the right sum is at most 3. A literal target without a pure witness is separated by \(x=0\), \(outsiders=2\).

For sufficiency, value 2 is least. A nonlinear target with only 0/\allowbreak{}1 factors gives 0, supplied by its nonlinear contained seed. A first 4 is supplied by E. A nonfinal first 3 gives 4, supplied by a marked E-witness, or by an unmarked E-witness of value 3 together with a Z(S) witness of value 0 or 1. A final first 3 is supplied by any E-witness, valued 3 or 4. These are all cases. Positive powers of x absorb the literal x, proving the literal sufficiency as well.

\par\smallskip\noindent\textit{2. Sixteen valid identities.}\quad 
Eight fixed laws are:

{\small\par\smallskip
\CaseLine{\(W1 xxy=xy\qquad W2 xyx=xyy\)}
\CaseLine{\(L  x\leq xx\qquad A  xyz\leq yz+xt\)}
\CaseLine{\(B  xyz\leq xz+xyt\qquad C  xyz\leq xz+xy+xx\)}
\CaseLine{\(D  xyz\leq xz+xx+y\qquad E  xy\leq xx+y\)}
\smallskip}

Eight additional identities are the following two finite templates. In each, a and b may independently be omitted when writing the identity, yielding four separate identities per row:

{\small\par\smallskip
\CaseLine{\([a]r[b]xu \leq  [a]r[b]u + [a]r[b]x + r\)}
\CaseLine{\([a]r[b]xu \leq  [a]r[b]u + x + r\)}
\smallskip}

The variables r, x, u are always present. Brackets describe this finite list, never an empty-word substitution. The criterion above verifies all absorptions: the first right term supplies N, and the original r supplies a contained prefix guard for the marked x when needed. The word-value formula proves W1, W2, L. All sixteen laws can also be checked independently on every assignment in the finite table.

W1, W2 give first-occurrence normal form F(w) when the final variable is simple, and F(w) with its final listed variable doubled otherwise. Internal separated repeats reduce by \(xYxZ=xY^{2}Z=xYZ\) with Y, Z nonempty; adjacent repeats use W1. A terminal separated repeat xYx becomes \(xY^{2}\) before internal deletions. Thus a nonlinear pure power normalizes to xx, not x. The criterion and original witness types are preserved.

\par\smallskip\noindent\textit{3. Prefix all original nonlinear words.}\quad 
Let N be the original nonlinear right summands. Normalize the target. For every marked target first occurrence, build its distinct-letter prefix p, maintaining p \(u\leq P\) for EVERY original u in N. Original literals are retained unchanged and are not covered by this invariant.

The first target position h is marked; Z(empty) is impossible, so there is an original marked witness hT, T nonempty. Write each u in N as YZ with Y, Z nonempty. A derives h \(u\leq u+hT\). This starts the invariant.

Suppose p is nonempty and x is the next marked position. A marked original E-witness \(q=R\) x T has T nonempty and c(R) contained in c(p). By the invariant, \(pq\leq P\); prefix deletion gives \(pq=p\) x T. B then derives p x \(u\leq p\) u+p x T for all u in N.

Otherwise take an unmarked original E-witness q and an original guard g contained in c(p), supplied by M. There are four cases:

\textbullet{} If q and g are nonlinear, \(pq=p\) x and \(pg=p^{2}\) are already available. C gives p x \(u\leq p\) u+p x+\(p^{2}\).

\textbullet{} If q is literal x and g is nonlinear, \(pg=p^{2}\) is available and D gives p x \(u\leq p\) u+\(p^{2}\)+x.

\textbullet{} If g is literal r and q is nonlinear, \(pq=p\) x is available. Write \(p=A\) r B at that literal's position and use the first of the two finite templates, selecting the appropriate A/B presence case.

\textbullet{} If both are literal, use the second template with the same split.

In the last two cases the unchanged original literal r provides the guard; no claim that \(p^{2}\) follows from r and p u is made. Such a claim would be false in this algebra: \(xx\leq x+xy\) fails at \(x=1\), \(y=2\). Every actual substitution in the applicable identity is nonempty. The four cases extend the invariant from p to p x, preserving all original witnesses.

\par\smallskip\noindent\textit{4. Saturated and simple-tail recovery.}\quad 
If the target is saturated, all its first occurrences are marked. After building the full first-occurrence list p, choose N's original contained seed v. The maintained p v normalizes to the saturated target, since every variable of v already occurs in p.

If the target is p s with simple final s, the construction has reached p \(u\leq P\) for every u in N. There are three possibilities.

If an original nonlinear unmarked E-witness q for s exists, \(pq=p\) s is the target. If a marked E-witness exists, perform one more marked extension at s; a contained nonlinear seed now gives p s s. Apply L in left context p to derive p \(s\leq p\) s s.

Otherwise the E-witness is the original literal s. Any original contained nonlinear seed involving s would itself be a nonlinear E-witness, contrary to this case. Therefore a contained seed v omits s, and p \(v=p^{2}\). E gives p \(s\leq p^{2}+s\). This covers the literal-final branch without asserting that all original literals can be prefixed.

A literal target follows from its original pure witness via W1, W2, L. Consequently every valid word absorption is derivable. Expanding terms by distributivity and deriving the summand absorptions in both directions proves every identity. The sixteen identities are a finite basis for

4649. Word reversal gives its opposite basis. Other catalogue transfers require complete HSP both ways; finite tests alone do not classify.

\subsection{A seven-identity interleaving basis for source 4665}\label{case:130}
Work in the constant-free ai-semiring signature. Element labels are not term constants. Addition and multiplication on the carrier 0,\allowbreak{}1,\allowbreak{}2,\allowbreak{}3,\allowbreak{}4 have rows

{\small\par\smallskip
\noindent addition  multiplication\par
\CaseLine{\ttfamily 00004\qquad 00234}
\CaseLine{\ttfamily 01014\qquad 00234}
\CaseLine{\ttfamily 00224\qquad 22232}
\CaseLine{\ttfamily 01234\qquad 33333}
\CaseLine{\ttfamily 44444\qquad 44434}
\smallskip}

The additive order has least element 3, incomparable elements 1 and 2, their join 0, and greatest element 4. The square map h=(0,\allowbreak{}0,\allowbreak{}2,\allowbreak{}3,\allowbreak{}4) is an endomorphic retraction onto B=\{0,\allowbreak{}2,\allowbreak{}3,\allowbreak{}4\}, and \(xy=h(x)h(y)\). On B, 0 is a multiplicative identity, 3 a zero, and 2 and 4 are left zeros in the submonoid B minus \{3\}. The addition on B is the chain \(3<2<0<4\). Direct checking of its four multiplication maps shows that B is a semiring; preservation of addition by h and \(xy=h(x)h(y)\) prove the same for the displayed five-element algebra A.

Theorem. A has the following seven-identity basis, relative to the ai-semiring axioms:

\CaseDisplay{x^2y\approx xy,\qquad xy^2\approx xy,\qquad xyx\approx xy;\quad\textup{(\ensuremath{W})}}

\CaseDisplay{x\leq x^2;\quad\textup{(\ensuremath{S})}}

\CaseDisplay{xyzt\leq xt+yz;\quad\textup{(\ensuremath{I})}}

\CaseDisplay{xyz\leq yz+x,\qquad xyz\leq xz+y.\quad\textup{(\ensuremath{L})}}

Reversing all words gives a basis for the opposite. No minimality or independence is asserted.

\par\smallskip\noindent\textit{Word values and exact criterion.}\quad 
Every nonlinear word is 3 if any of its variables is assigned 3. Otherwise it is 0 if all variables are assigned 0 or 1; if not, it has the value of its first occurrence assigned 2 or 4. Laws (W) therefore hold. They reduce a nonlinear word to its first-occurrence list, except that a singleton support is represented by a square. Algebraically separated repetitions are deleted by \(xyx=xy\) with the middle block substituted for y; adjacent repetitions use the other two laws. A square is never reduced to a linear word.

For a polynomial P, let L be its linear variable set and N its nonlinear words. For C a set of variables, let \(N_{C}\) consist of those u in N whose content is contained in C. Write \(pre_{w}\)(x) for the variables before the first occurrence of x.

Lemma. A linear word y satisfies \(y\leq P\) precisely when y is in L or N contains a word with content \{y\}. A nonlinear w, with C its content, satisfies \(w\leq P\) precisely when

1. \(N_{C}\) is nonempty;

2. for every x in C, either x is in L or some v in \(N_{C}\) contains x and \(pre_{v}\)(x) is contained in \(pre_{w}\)(x).

Proof. If neither witness for a linear y exists, assign \(y=4\) and every other variable 3. The target is 4 and every right summand 3. The converse follows from (S).

If condition 1 fails, assign all variables in C the value 1 and all others 3. The nonlinear target has value 0, whereas the right polynomial is at most 1. If condition 2 fails at x, assign variables before its first occurrence in w the value 0, assign \(x=4\), assign the rest of C the value 2, and all outsiders 3. The target is 4. Any right word containing an outsider is 3. Every right nonlinear word supported in C either omits x or meets a 2 before x, so has value at most 0. No right linear summand is x, so their join is also at most 0. These assignments prove necessity.

Conversely, target value 3 is automatic. If the target is 0, all target variables are 0 or 1, and a word in \(N_{C}\) is 0. If the target is 2, no target variable is 3, so any word in \(N_{C}\) is in \{0,\allowbreak{}2,\allowbreak{}4\}, all at least 2. Finally, if the target is 4, let x be its first occurrence assigned 2 or 4. Then x is 4 and all preceding variables are 0 or 1. A linear witness x is 4; a nonlinear prefix witness in \(N_{C}\) is also 4. This proves sufficiency for arbitrary variable sets.

The lemma also proves soundness of (I) and (L): a contained nonlinear word witnesses condition 1, and each target prefix is witnessed by the displayed right words or by its corresponding linear variable. Rule (S) follows from the five square values and the additive order.

\par\smallskip\noindent\textit{Unrestricted equational derivation.}\quad 
Normalize nonlinear words by (W). Consider a nonlinear valid absorption \(w\leq P\), and write the first-occurrence list of w as \(x_{1}\),..., \(x_{d}\). Let \(p_{i}=x_{1}\cdots x_{i}\) and \(C=c(w)\). We inductively construct a nonlinear word \(s_{i}\leq P\) with content contained in C and with first-occurrence list beginning with \(p_{i}\). This is a derivation from the seven displayed identities.

For \(i=1\), if \(x_{1}\) is not linear in P, its nonlinear prefix witness has head \(x_{1}\) and may serve as \(s_{1}\). If \(x_{1}\) is linear in P, choose u in \(N_{C}\), split \(u=yz\) with both factors nonempty, and use the first law (L) to obtain \(x_{1}u\leq u+x_{1}\leq P\).

Suppose \(s_{i}\) has been constructed and \(x=x_{i+1}\). After normalization, write \(s_{i}\) approximately \(p_{i}\) q with q nonempty. This is possible directly unless its normal first-occurrence list ends at \(p_{i}\). In that case \(i\geq 2\), and \(s_{i}\) approximately \(p_{i}^{2}\) by (W), so take \(q=p_{i}\). When \(i=1\) the singleton-support normal word already has two letters. No empty term is substituted.

If x is a linear summand of P, the second rule (L) gives \(p_{i}\) x \(q\leq p_{i}\) \(q+x\leq P\). Its normalized first-occurrence list starts with \(p_{i}\) x. Otherwise choose its nonlinear witness v in \(N_{C}\) and split \(v=ab\) with nonempty factors. Applying (I), with \((x,y,z,t)=(p_i,a,b,q)\), gives

\CaseDisplay{p_i v q\leq p_i q+v\leq P.}

All letters before x in v already occur in \(p_{i}\), and hence are removed as repetitions by (W). Therefore the resulting normalized word starts with \(p_{i}\) x and still has content contained in C. This proves the induction.

At \(i=d\), the word \(s_{d}\) contains all of C, introduces no outsiders, and has exactly the first-occurrence order of w. Its normal form is w, including singleton support. Thus \(w\leq P\) is derivable. Linear valid absorptions follow from an existing linear summand or from (S) and a nonlinear singleton-support witness. Expanding arbitrary identities into polynomials and deriving every summand absorption in both directions proves unrestricted completeness. Word reversal carries all derivations to the opposite.

\subsection{Contained-content witnesses: a six-identity basis for a group extension}\label{case:131}
\par\smallskip\noindent\textit{Abstract.}\quad 
We prove a six-identity basis for a five-element additively idempotent semiring containing a group of order two and a clone of an external identity. Absorption is determined by content coverage, affine exponent parity, and the presence of a stable right summand whose content lies in the target content. A guarded contraction removes all auxiliary variables in one step. The argument allows any finite number of variables and summands.

\par\smallskip\noindent\textit{1. Algebra and basis.}\quad 
Work in the signature \((+,\cdot)\) without constants, relative to the additively idempotent semiring axioms. As usual, \(u\le v\) means \(u+v\approx v\). Put \(T=\{e, g, d, f,\omega\}\). Addition is join in the order

\CaseDisplay{e<d<f<\omega,\qquad g<\omega,}

with \(g\) incomparable with \(e, d, f\). On \(F=\{e, g, f,\omega\}\), \(\{e, g\}\) is the group of order two, \(f\) is an identity, and \(\omega\) is an annihilator. Define \(\rho(d)=f\), fixing every other element, and multiply in \(T\) by \(uv=\rho(u)\rho(v)\), with the product on the right taken in \(F\). Element names do not denote constants in terms.

Theorem 1. The identities

\CaseDisplay{\begin{array}{rll}
\mathrm{(C)}&xy\approx yx,\\
\mathrm{(N)}&x^3y\approx xy,\\
\mathrm{(U)}&x\le x^3,\\
\mathrm{(O)}&xyz\le x^3+y^3+z^3,\\
\mathrm{(G)}&x^3y^2z^2\le x+yz,\\
\mathrm{(K)}&xy^2\le xy^2z^2+y^3
\end{array}}

form a complete basis for \(T\). No minimality claim is made.

First, \(F\) is a semiring. Multiplication is associative: adjoining an identity and an annihilator to the group preserves associativity. Translation by \(f\) is the identity, translation by \(\omega\) is constant, and translation by \(e\) or \(g\) preserves joins in the order \(e<f<\omega\), \(g<\omega\), as is checked on the incomparable pairs \(\{e, g\}\) and \(\{f, g\}\). Thus multiplication distributes. The map \(\rho\) is a join-preserving retraction, and all products on \(T\) lie in \(F\), so associativity and distributivity follow on \(T\) as well.

Call a word stable when its length is at least two. A stable word is evaluated by multiplying the \(\rho\) images of its letters. A top letter gives \(\omega\); if all images are \(f\) the value is \(f\); otherwise the value lies in \(\{e, g\}\) and is determined by the parity of the letters whose images are \(g\). In particular \(x^3=\rho(x)\). This proves (C),(N),(U). For (O), a top image on the right is immediate. Among the other images, a mixture of \(g\) with \(e\) or \(f\) has top join. All three images equal to \(g\) have product \(g\). Images from \(\{e, f\}\) have their product at most their join. These exhaust the cases.

For (G), let \(q=yz\in F\). If \(x\) or \(q\) is top, the right side is top. Otherwise, if \(q=f\) then \(q^2=f\) and \(x+f\) dominates \(\rho(x)\). If \(q\) is in the group, then \(q^2=e\): for group \(x\) the left side equals \(x\), and for \(x\in\{d, f\}\) it is \(e\le x\). Hence (G) holds.

For (K), write \(X=\rho(x)\), \(Y=\rho(y)\), \(Z=\rho(z)\). If any is \(\omega\), the right side is top. Otherwise they lie in \(\{e, g, f\}\). If \(Y\) lies in the group, then \(XY^2\) is a group element and multiplication by \(Z^2\in\{e, f\}\) leaves it unchanged. If \(Y=f\) and \(X\) lies in the group, the same observation applies. If \(X=Y=f\), the right summand \(Y^3=f\) already equals the left side. This proves soundness.

\par\smallskip\noindent\textit{2. Stable affine closure.}\quad 
Let \(c(w)\) be the content of a word and \(\pi(w)\) its exponent-parity vector over \(\mathbb F_2\). By (C),(N), stable words of the same content and parity are equal: for content of size at least two remove pairs from each exponent using a nonempty complementary word in (N); for singleton content use \(x^4=x^2\) to reduce to \(x^2\) or \(x^3\). Thus \(r^3=r\) for every stable word \(r\).

For a nonempty family \(R\) of stable words, every stable word \(v\) with

\CaseDisplay{c(v)=\bigcup_{r\in R}c(r),\qquad
  \pi(v)\in\operatorname{aff}_{\mathbb F_2}\{\pi(r):r\in R\}\quad\textup{(\ensuremath{1})}}

satisfies \(v\le\sum R\) derivably. To see this, express its parity as an odd sum of the listed vectors. Form the corresponding product, including two additional copies of every member of \(R\). Its content is the full union in (1). Iterating (O) bounds this odd product by \(\sum R\), since cubes of stable factors normalize to themselves. Normalization identifies it with \(v\). The one-factor case is reflexive. No empty product or unit constant is used.

\par\smallskip\noindent\textit{3. Absorption by a polynomial.}\quad 
Let \(P\) be a nonempty sum of words, \(L\) its linear variables, and \(R\) its stable summands. Write

\CaseDisplay{U=\bigcup_{u\in P}c(u),\qquad
H=\operatorname{aff}_{\mathbb F_2}\{\pi(u):u\in P\}.}

Theorem 2. The absorption \(w\le P\) holds in \(T\) if and only if either \(w=x\in L\), or all three conditions hold:

\CaseDisplay{c(w)\subseteq U,\qquad \pi(w)\in H,\qquad
  (\exists q\in R)\ c(q)\subseteq c(w).\quad\textup{(\ensuremath{2})}}

Every valid absorption is derivable from the six displayed laws.

Necessity. The listed linear case is immediate. If a variable of the target is missing from \(U\), assign it \(\omega\) and every other variable \(e\). The target is \(\omega\) and the right side \(e\). If the target parity does not belong to \(H\), a linear functional separates it from the affine hull. The associated group assignment makes every right word equal to one group element and the target equal to the other. As \(e, g\) are incomparable, this disproves absorption.

To prove the last condition, suppose there is no stable right word with content contained in \(C=c(w)\). Assign every variable of \(C\) the element \(d\), and every other variable \(e\). Every stable right word contains an outside variable; its product therefore has value \(e\). Linear right words have value \(d\) or \(e\), so \(P\le d\). A stable target has value \(f>d\). If the target is linear and unlisted, it is \(x\notin L\); no right linear word then lies in \(C=\{x\}\), so the right value is exactly \(e\), whereas the target is \(d>e\). Both cases give a separating evaluation.

Sufficiency. Suppose (2) holds, and fix the contained stable witness \(q\in R\). Factor \(q=uv\) with \(u, v\) nonempty. For each \(x\in L\), (G) gives

\CaseDisplay{x^3q^2=x^3u^2v^2\le x+uv=x+q\le P.\quad\textup{(\ensuremath{3})}}

Together with the original stable summands these form a nonempty stable family \(R'\) whose total content is \(U\) and whose affine parity hull is exactly \(H\). Indeed, the new words supply exactly the unit vectors of the linear summands, without introducing outside variables.

If \(w\) is linear, first replace it by \(v=w^3\); otherwise put \(v=w\). In either case \(v\) is stable, \(c(v)=C\), and \(\pi(v)=\pi(w)\). If \(U=C\), apply (1) to derive \(v\le P\). If \(U\setminus C\) is nonempty, let \(z\) be a word listing that set. The stable word \(vq^2z^2\) has content \(U\) and parity \(\pi(w)\), so (1) derives \(vq^2z^2\le P\). Since \(c(q)\subseteq C\), normalization gives \(vq^2=v\) and \(q^3=q\). Substituting the nonempty words \(v, q, z\) in (K) yields

\CaseDisplay{v\approx vq^2\le vq^2z^2+q^3\le P.\quad\textup{(\ensuremath{4})}}

For a linear target, conclude with (U). This proves Theorem 2.

The contained witness in (4) is essential. It controls the contraction of auxiliary content; an unguarded content deletion is not assumed. The word \(z\) is used only in the nonempty case, so the constant-free signature is respected throughout.

\par\smallskip\noindent\textit{4. Completeness and classification.}\quad 
Distribute an arbitrary valid identity into nonempty polynomials. Each summand on either side is absorbed by the other side, and Theorem 2 derives all these absorptions. Taking their sums proves both inequalities and the identity, establishing Theorem 1 for arbitrarily many variables.

The named-to-catalog map in the order \((e, g, d, f,\omega)\) is \([2,\allowbreak{}3,\allowbreak{}1,\allowbreak{}0,\allowbreak{}4]\) for type 4677. Proposition~\ref{prop:lower-fibres} in both directions gives \(\V(\A{4677})=\V(\A{11789})=\V(\A{12070})\). The six identities apply to all three.

The source manuscript for this case credits Aifa Wang and Lili Wang, School of Mathematical Sciences, Chongqing University of Technology.

The source manuscript for this case acknowledges support from the National Natural Science Foundation of China (12371024), the Science and Technology Research Program of Chongqing Municipal Education Commission (KJZD-K2024011102) and the Chongqing Natural Science Foundation Innovation and Development Joint Fund (Municipal Education Commission) (CSTB2025NSCQ-LZX0067).

\subsection{内部非一次种子与首字母守卫：两组完整有限基}\label{case:132}
\par\smallskip\noindent\textit{1. 统一记号.}\quad 
对词 \(w\) 记 \(X=c(w)\)、\(h=h(w)\)、\(D=D(w)\)。对非空右词集 \(F\) 定义

\CaseDisplay{T=c(F),\quad H=\{h(u):u\in F\},\quad E=\bigcup_{u\in F}D(u),\quad
F_X=\{u\in F:c(u)\subseteq X\},}

以及 \(C_X=c(F_X)\)、\(H_X=h(F_X)\)、\(E_X=\bigcup_{u\in F_X}D(u)\)，空集的并取空集。令 \(F_X^{\ge2}\) 为 \(F_X\) 中非一次词的集合。

\par\smallskip\noindent\textit{2. 来源5008：内部非一次种子与全局保护.}\quad 
取 \(A=\{z, a, b, e, t\}\)，对应目录 \((0,\allowbreak{}1,\allowbreak{}2,\allowbreak{}3,\allowbreak{}4)\)。加法偏序为 \(b<a<z<t\)、\(e<t\)，其中 \(e\) 与链 \(b<a<z\) 不可比。右乘 \(z, a\) 均为 \(P=(z, z, b, e, t)\)，右乘 \(b\) 为 \(Q=(b, b, b, e, t)\)，右乘 \(e, t\) 为常值 \(t\) 映射。三个映射保持并、交换且幂等，\(PQ=QP=Q\)；按右因子类别检查 \(R_{xy}=R_yR_x\) 和 \(R_{x+y}=R_x\vee R_y\) 即得半环。

对于非一次词：尾部出现 \(e, t\) 时取 \(t\)；否则首字母为 \(e, t\) 时分别取 \(e, t\)，首字母在链 \(b<a<z\) 上时，内容含 \(b\) 则取 \(b\)、否则取 \(z\)。一次词直接取其变量值。

\par\smallskip\noindent\textit{2.1 精确判据.}\quad 
定理1。 若 \(w=x\) 为一次词，则

\CaseDisplay{A\models x\le F\quad\Longleftrightarrow\quad F_{\{x\}}\ne\varnothing.\quad\textup{(\ensuremath{A1})}}

若 \(w\) 非一次，则

\CaseDisplay{A\models w\le F\quad\Longleftrightarrow\quad
X\subseteq T,\quad F_X^{\ge2}\ne\varnothing,\quad
\bigl[H\cap E\ne\varnothing\ \text{或}\ (h\in H\ \text{且}\ D\subseteq E)\bigr].\quad\textup{(\ensuremath{A2})}}

必要性。 对一次左词，若无内部纯幂项，赋 \(x=a\)、其他变量 \(b\)，则右侧单词均取 \(b\)，而左值为 \(a\)。反向由纯幂正规化及 \(x\le x^2\) 立即得到。

对非一次左词，若 \(x\in X\setminus T\)，赋 \(x=t\)、其余变量 \(z\)，左值为 \(t\)、右值为 \(z\)。若 \(F_X^{\ge2}=\varnothing\)，赋 \(X\) 内变量 \(a\)、外部变量 \(b\)；全部右侧一次词不超过 \(a\)，全部非一次右词因含外部变量而取 \(b\)，而左值为 \(z\)。

最后设 \(H\cap E=\varnothing\)，令 \(H\) 中变量取 \(e\)、\(E\) 中变量取 \(z\)。已经成立的 \(X\subseteq T=H\cup E\) 使这涵盖全部左变量。右词均取 \(e\)。若 \(h\notin H\)，左值为 \(z\) 或 \(t\)；若 \(h\in H\) 而 \(D\nsubseteq E\)，左尾部出现 \(e\)，左值为 \(t\)。均不能不超过右值 \(e\)。

充分性。 真主理想为 \(\{b\}\)、\(\{b, a\}\)、\(\{b, a, z\}\)、\(\{e\}\)。右和不超过 \(b\) 或 \(a\) 时，每个右词的全部变量在链 \(b<a<z\) 上，每个非一次右词还必须含取值 \(b\) 的变量。全局内容覆盖保证左内容在该链上，内部非一次见证保证左内容含 \(b\)，故左值为 \(b\)。右和不超过 \(z\) 时，全局覆盖同样使左内容都在链上，左值不超过 \(z\)。

右和为 \(e\) 时，全部右首字母取 \(e\)，右尾部取链上值。因此 \(H\cap E\) 必为空；(A2) 的后一选择使左首字母取 \(e\)、尾部在链上，左值为 \(e\)。充分性成立。

\par\smallskip\noindent\textit{2.2 十条完整基.}\quad 
定理2。 \(A\) 的完整附加基为

\CaseDisplay{\begin{aligned}
xyz&\approx xzy,&xy^2&\approx xy,&x^2y&\approx y^2x,\\
x&\le x^2,&xyz&\le xy+tz,&x^2y&\le x^2+y,\\
x^2y&\le x^2+yz,&(xyz)^2&\le xy+yz,&(xy)^2&\le xy+y,\\
(xy)^2&\le xy+z^2.&&&&
\end{aligned}\quad\textup{(\ensuremath{A3})}}

有效性由定理1检查；前三式亦可直接由上述词值规则检验。前三式使首字母重复的词只由内容决定，记为 \(M_Y\)；任意 \(u^2=M_{c(u)}\)，且此词乘法幂等。

考虑非一次的有效 \(w\le F\)。先设 \(H\cap E=\varnothing\)。取内部非一次项 \(u\in F_X^{\ge2}\)。其首字母属于 \(H\cap X\)；又 \(X\setminus\{h\}\subseteq D\subseteq E\)，所以首字母只能为 \(h\)，且 \(h\) 在 \(u\) 中不重复。对每个尚未出现的目标尾变量 \(g\)，从全局尾部见证中用尾部换序得到右项 \(vg\)，其中 \(v\) 非空。写当前已吸收非一次词为 \(hp\)，\(p\) 非空，由 \(xyz\le xy+tz\) 得

\CaseDisplay{hpg\le hp+vg\le F.}

逐个补齐后即得 \(w\)。

再设 \(H\cap E\ne\varnothing\)。选 \(g\in H\cap E\)，取右项 \(vg\) 及以 \(g\) 开头的右项。若后一项为 \(gs\)、\(s\) 非空，由 \((xyz)^2\le xy+yz\) 得 \((vgs)^2\le F\)；若后一项为 \(g\)，由 \((xy)^2\le xy+y\) 得 \((vg)^2\le F\)。于是有某个平方词 \(q^2\le F\)，其内容不必包含于 \(X\)。

对内部非一次项 \(u=rs\)，\(r, s\) 非空，由最后一式代入 \(r, s, q\) 得 \(u^2\le u+q^2\le F\)。故已有平方词 \(p=M_{c(u)}\le F\)，内容包含于 \(X\)。

逐个添加 \(X\) 中尚缺的变量 \(g\)。全局内容覆盖提供包含 \(g\) 的右词。若见证为一次词 \(g\)，用 \(x^2y\le x^2+y\)；若见证以 \(g\) 开头且非一次，写作 \(gs\)，用 \(x^2y\le x^2+yz\)；若 \(g\) 在尾部，将见证换序写成 \(vg\)，用 \(xyz\le xy+tz\)，其中 \(x=y=p\)、\(t=v\)。因 \(p^2=p\)，三种情形均导出 \(pg\le F\)，且新词仍为平方正规形，内容仅增加 \(g\)。最终得到 \(M_X=w^2\le F\)，再用 \(w\le w^2\)。

一次左词使用 (A1)，纯幂由前两式正规化并应用 \(x\le x^2\)。所有代入非空，任意变量数完备性得证。

\par\smallskip\noindent\textit{3. 来源6921：内部保护与全局首字母守卫.}\quad 
取 \(B=\{z, a, b, d, t\}\)，对应目录 \((0,\allowbreak{}1,\allowbreak{}2,\allowbreak{}3,\allowbreak{}4)\)。加法偏序为 \(b<a<z<t\)、\(a<d<t\)，\(z, d\) 不可比。右乘 \(z\) 是恒等映射，右乘 \(a, b\) 为 \(P=(a, a, b, a, a)\)，右乘 \(d, t\) 为 \(Q=(t, a, b, t, t)\)。这些映射保持并、交换、幂等，\(PQ=QP=P\)；按右因子类别验证两条右平移兼容式可得半环。

词的首字母为 \(b\) 时恒取 \(b\)；首字母为 \(a\) 时恒取 \(a\)；其他情况下，尾部含 \(a, b\) 时取 \(a\)，否则尾部含 \(d, t\) 时取 \(t\)，否则尾部全为 \(z\)、取首字母值。该规则包含一次词。

\par\smallskip\noindent\textit{3.1 精确判据.}\quad 
定理3。 对任意词 \(w\)，

\CaseDisplay{B\models w\le F\quad\Longleftrightarrow\quad
C_X=X,\quad h\in H,\quad
\bigl[H_X\cap E_X\ne\varnothing\ \text{或}\ (h\in H_X\ \text{且}\ D\subseteq E_X)\bigr].\quad\textup{(\ensuremath{B1})}}

必要性。 若 \(x\in X\setminus C_X\)，赋 \(x=t\)、其余 \(X\) 内变量 \(z\)、外部变量 \(a\)。内部右项取 \(z\)，外部右项取 \(a\)，右和不超过 \(z\)，左值为 \(t\)。

若 \(h\notin H\)，赋 \(H\) 中变量 \(b\)、其余变量 \(z\)。右词均取 \(b\)，左值为 \(z\) 或 \(a\)，不小于等于 \(b\)。

若 \(H_X\cap E_X=\varnothing\)，赋 \(H_X\) 中变量 \(d\)、\(E_X\) 中变量 \(z\)、外部变量 \(a\)。由于 \(C_X=X\)，这给出全部变量的一致赋值。内部右词取 \(d\)，外部右词取 \(a\)，右和不超过 \(d\)。若 \(h\notin H_X\)，左值为 \(z\) 或 \(t\)；若某左尾变量不在 \(E_X\)，它属于 \(H_X\) 并取 \(d\)，左值为 \(t\)。由此得到最后条件。

充分性。 真主理想为 \(\{b\}\)、\(\{b, a\}\)、\(\{b, a, z\}\)、\(\{b, a, d\}\)。右和为 \(b\) 时，所有右首字母都取 \(b\)，由全局首字母条件左词取 \(b\)。

对其余三个主理想，若左内容含取值 \(a, b\) 的变量，则左词取 \(a\) 或 \(b\)，结论成立。以下假设左内容不含这些值。于是内部右词同样不含这些值。右和不超过 \(a\) 已不可能，因为内部右项非空且任何不含 \(a, b\) 的词都不小于等于 \(a\)。右和不超过 \(z\) 时，每个内部右词只能全部变量取 \(z\)；内部覆盖使左词取 \(z\)。

右和不超过 \(d\) 时，内部右词必须首字母取 \(d\)、尾部全取 \(z\)。因此 \(H_X\cap E_X=\varnothing\)，(B1) 强制内部同首字母和内部尾部覆盖，左词取 \(d\)。充分性成立。

\par\smallskip\noindent\textit{3.2 十条完整基.}\quad 
取下列九式：

\CaseDisplay{\begin{aligned}
xyz&\approx xzy,&xy^2&\approx xy,&x&\le x^2,\\
xyz&\le xy+xz,&xyz&\le xy+yz,&xy&\le x+y^2,\\
x^2y&\le x^2+y,&(xyz)^2&\le xy+yz,&(xy)^2&\le xy+y.
\end{aligned}\quad\textup{(\ensuremath{B2})}}

定理4。 \((B2)\) 加上首字母守卫式

\CaseDisplay{x^2y\le y^2x+xz\quad\textup{(\ensuremath{G})}}

构成 \(B\) 的完整附加基。这十式的有效性由定理3检验。

引用内部保护与全局尾部覆盖第1节的句法闭合引理：\((B2)\) 的前六式能推出每个满足 \(h\in H_X\) 和 \(D\subseteq E_X\) 的吸收。该引理的任意变量数证明详见其引用的内部首字母与尾部闭合第4节。

因此只有 \(H_X\cap E_X\ne\varnothing\) 的保护情形需要新论证。选内部尾部见证 \(vg\) 和内部同首字母见证 \(gs\) 或 \(g\)，以两个平方路径式得到内部平方词 \(p\le F\)。然后用 \(x^2y\le x^2+y\) 把全部内部词乘入 \(p\)。所得词 \(Q\le F\) 的内容和尾部集合都等于 \(X\)，并有 \(Q^2=Q\)；其首字母暂时不必是 \(h\)。

全局条件 \(h\in H\) 提供一个以 \(h\) 开头的右项 \(r\)。若 \(r=h\)，由 \(xy\le x+y^2\) 得 \(hQ\le h+Q\le F\)。若 \(r=hs\)、\(s\) 非空，将 \((G)\) 中 \(x, y, z\) 分别替换为 \(h, Q, s\)，得

\CaseDisplay{h^2Q\le Q^2h+hs=Q+hs\le F.}

这里 \(Q^2h=Q\) 仅用尾部换序去重，因为 \(h\in X=D(Q)\)。两种情况下，所得左词均以 \(h\) 开头、尾部恰为 \(X\)，等于 \(w^2\) 的正规形。再用 \(w\le w^2\)。故十式对任意变量数完整。此证明不使用在 \(B\) 上不成立的无守卫换首式 \(x^2y\approx y^2x\)。

\par\smallskip\noindent\textit{4. 证据范围.}\quad 
两组完整基都先由精确判据证明有效，再给出任意变量数的句法推导。有限变量CNF核验只作交叉检查。乘法对偶逐词反转处理；有限基性向其他目录对象传递时另核验双向HSP包含。

新方法有两种不同的作用：来源5008只需一个内部非一次种子，再从全局见证补内容；来源6921允许从全局取得首字母，但只有内部保护成立时才能受守卫地换首。这些结论不推出相应簇遗传有限基。

\subsection{三阶群的群元克隆：两组六条完整基}\label{case:133}
\par\smallskip\noindent\textit{1. 两个来源.}\quad 
在三阶循环群 \(\{e, g, g^2\}\) 的平坦扩张中，\(c\) 是乘法零和加法最大元。添元素 \(d\)，并由固定核心的收缩 \(\rho\) 定义乘法 \(xy=\rho(x)\rho(y)\)。取两种构造：

\textbullet{} \(D\)：\(d<e<c\)，\(g, g^2<c\)，其余不同元素不可比；\(\rho(d)=e\)。

\textbullet{} \(U\)：\(g<d<c\)，\(e, g^2<c\)，其余不同元素不可比；\(\rho(d)=g\)。

它们分别同构于目录 4703、10496。按顺序 \((e, g, g^2, d, c)\)，同构分别为 \((0,\allowbreak{}2,\allowbreak{}3,\allowbreak{}1,\allowbreak{}4)\)、\((3,\allowbreak{}1,\allowbreak{}2,\allowbreak{}0,\allowbreak{}4)\)。\(\rho\) 及各核心乘法平移保持并，故两者为加法幂等半环。

令 \(\Gamma\) 为四条共同式子

\CaseDisplay{xy\approx yx,\qquad x^4y\approx xy,\qquad
x^3y^3\approx x^3+y^3,\qquad
xyzt\leq x^4+y^4+z^4+t^4.}

定理。 \(D\) 的完整基为

\CaseDisplay{\Gamma,\qquad x\leq x^4,\qquad x^4\leq x+yz.
\quad\textup{(\ensuremath{D})}}

\(U\) 的完整基为

\CaseDisplay{\Gamma,\qquad x^4\leq x,\qquad x\leq x^3+x^4.
\quad\textup{(\ensuremath{U})}}

两组均含六条附加恒等式。

共同式子的有效性仍由三阶群平坦核心给出。对 \(D\)，仅 \(x=d\) 时 \(x\leq x^4\) 不是等式；此时 \(x^4=e\)，而任一非线性值 \(v\in\{e, g, g^2, c\}\) 都使 \(d+v\geq e\)。对 \(U\)，仅 \(x=d\) 时需要说明：\(x^4=g<d\)，\(x^3+x^4=e+g=c\geq d\)。

\par\smallskip\noindent\textit{2. 任意变量数判据.}\quad 
对词指数取模 3，记相应仿射包为 \(H_3(P)\)，全部右内容为 \(U(P)\)。

在 \(D\) 中：

1. 一次词 \(x\leq P\) 当且仅当 \(\pi_3(x)\in H_3(P)\)。

2. 非线性词 \(w\leq P\) 当且仅当 \(\pi_3(w)\in H_3(P)\)、\(c(w)\subseteq U(P)\)，且 \(P\) 至少有一个非线性项。

在 \(U\) 中：

1. 非线性词 \(w\leq P\) 当且仅当 \(\pi_3(w)\in H_3(P)\) 且 \(c(w)\subseteq U(P)\)。

2. 一次词 \(x\leq P\) 当且仅当 \(x\) 是右侧一次项，或 \(\pi_3(x),0\in H_3(P)\)。

群赋值和``左侧独有变量赋 \(c\)，其余赋 \(e\)''分别给出仿射条件和内容条件的必要性。在 \(D\) 中，若右侧全部是一次项，将全部变量赋 \(d\)：非线性左词为 \(e\)，右侧为 \(d<e\)。在 \(U\) 中，若 \(x\) 不是右侧一次项且 \(0\notin H_3(P)\)，选择线性泛函使全部右向量取值 1；若 \(\pi_3(x)\in H_3(P)\)，此赋值也使 \(x=g\)。把 \(x\) 改为 \(d\)，所有非线性右词保持群值 \(g\)，而左值为 \(d>g\)，故吸收失败。这里使用任意不含零的仿射子空间都位于某个取值为 1 的线性超平面内。

\par\smallskip\noindent\textit{3. 完整性.}\quad 
由 \(\Gamma\)，非线性词按内容和指数模 3 规范化。仍有：非线性 \(u\leq uz^3\)，以及任意 \(3r+1\) 个非线性词之积被这些词之和吸收。因而只要非空非线性词族 \(Q\) 的内容覆盖非线性 \(w\)，且 \(\pi_3(w)\in H_3(Q)\)，就能推出 \(w\leq\sum Q\)。具体构造是按仿射表示选取因子，再将 \(Q\) 中每个词各补三次，以补足内容而不改变指数模 3。

对于 \(D\)，若 \(P\) 有非线性项 \(uv\)，由 \(x^4\leq x+uv\) 可向 \(P\) 加入每个一次项的四次幂；由 \(x\leq x^4\)，原一次项还可删除。得到一个与 \(P\) 等价的非线性词族 \(Q\)，其内容和仿射包与 \(P\) 相同。上述引理给出所有非线性吸收式。对于一次左词 \(x\)，其仿射条件保证 \(x\in U(P)\)；若右侧有非线性项，则先推得 \(x^4\leq P\)，再用 \(x\leq x^4\)。若右侧只有一次项，\(\pi_3(x)\) 在这些单位向量的仿射包中当且仅当右侧已有 \(x\) 本身。

对于 \(U\)，由 \(x^4\leq x\) 可向右侧加入全部一次项的四次幂。得到非线性词族 \(Q\)，满足 \(U(Q)=U(P)\)、\(H_3(Q)=H_3(P)\) 和 \(\sum Q\leq P\)，所以推得全部满足内容、仿射条件的非线性吸收式。若一次左词 \(x\) 不是右侧一次项而 \(\pi_3(x),0\in H_3(P)\)，则 \(x\in U(P)\)，非线性判据给出

\CaseDisplay{x^3\leq P,\qquad x^4\leq P.}

于是 \(x\leq x^3+x^4\leq P\)。

至此所有有效逐词吸收式均由相应六条基推出，任意多项式恒等式再分解为双向逐词吸收，完整性得证。

\par\smallskip\noindent\textit{4. 范围.}\quad

\subsection{四个有锚点的交换近幂等半环}\label{case:134}
\par\smallskip\noindent\textit{1. 四个生成元及其基.}\quad 
共同底集为 \(\{0, a, b, e, f\}\)。令 \(b\) 在乘法中与 \(e\) 有相同作用，\(e\) 为四元素核心 \(\{0, a, e, f\}\) 的幺元，\(a^2=0\)、\(f^2=f\)。

前三个乘法类型中，\(0\) 是整个代数的乘法零，\(fa=a\)、\(f0=0\)。唯一例外是下表的 \(C_G\)：其中 \(f\) 是共同吸收元，在其余四个元素上 \(0\) 仍是乘法零。这些说明还包括 \(b^2=e\)、\(be=e\)、\(ba=a\)、\(bf=f\)。

\noindent 半环\enspace\textbar\enspace 目录代表\enspace\textbar\enspace 加法序的覆盖关系\enspace\textbar\enspace 与一般乘法的区别\par
\noindent \(C\)\enspace\textbar\enspace 4709\enspace\textbar\enspace \(0<a<b<e<f\)\enspace\textbar\enspace 无\par
\noindent \(C_H\)\enspace\textbar\enspace 4707\enspace\textbar\enspace \(0<a<e<f,\ b<e\)，\(b\parallel0, a\)\enspace\textbar\enspace 无\par
\noindent \(C_G\)\enspace\textbar\enspace 4713\enspace\textbar\enspace \(0<a<b<e<f\)\enspace\textbar\enspace \(f\) 乘法吸收\par
\noindent \(C_{HG}\)\enspace\textbar\enspace 4710\enspace\textbar\enspace \(0<a<f,\ b<e<f\)，两条链之间无其他比较\enspace\textbar\enspace 无\par
加法为所列偏序的并。四个代数按 \((0, a, b, e, f)\) 顺序到目录标号的映射均为 \((3,\allowbreak{}2,\allowbreak{}1,\allowbreak{}0,\allowbreak{}4)\)。结合律由四元素交换核心及 \(b\mapsto e\) 的乘法收缩得到；平移保持相应的并，故均为加法幂等半环。运算表在独立核验中按此定义重建。

令 \(\Sigma\) 为

\CaseDisplay{\begin{aligned}
xy&\approx yx, &(S1)\\
xy&\approx x^2y+xy^2, &(S2)\\
xy&\leq x+y^2, &(S3)\\
xy&\leq x+(xy)^2, &(S4)\\
xy^2z^2&\leq xy^2+x^2z^2. &(S5)
\end{aligned}}

再记

\CaseDisplay{x^2\leq x+yz\quad(H),\qquad
x^2y\leq x^2+y^2z^2\quad(G).}

定理。 四个半环的完整基分别为

\CaseDisplay{\operatorname{Id}(C):\Sigma,\quad
\operatorname{Id}(C_H):\Sigma+H,\quad
\operatorname{Id}(C_G):\Sigma+G,\quad
\operatorname{Id}(C_{HG}):\Sigma+H+G.}

这里``完整基''均相对于通常加法幂等半环公理而言，不声称条数最小。

\par\smallskip\noindent\textit{2. 任意变量数的精确判据.}\quad 
对非空多项式 \(Q\)，令 \(L\) 为其一次项变量集，\(R\) 为非线性求和词族。置

\CaseDisplay{R_C=\{u\in R:c(u)\subseteq C\},\quad
G_Q=\bigcup_{u\in R}c(u),\quad
G_C=\bigcup_{u\in R_C}c(u).}

对非线性左词 \(w\)，记 \(C=c(w)\)、\(S(w)=\{x:m_x(w)=1\}\)。

四个半环均满足

\CaseDisplay{x\leq Q\iff x\in L.\quad\textup{(\ensuremath{1})}}

对非线性 \(w\)，\(w\leq Q\) 恰由以下三个条件刻画。

1. 锚点。 不带 \(H\) 时要求 \(R_C\ne\varnothing\)；带 \(H\) 时要求

\CaseDisplay{R_C\ne\varnothing\quad\text{或}\quad
   (L\cap C\ne\varnothing\text{ 且 }R\ne\varnothing).
   \quad\textup{(\ensuremath{2})}}

2. 覆盖。 不带 \(G\) 时要求 \(C\subseteq L\cup G_C\)；带 \(G\) 时要求 \(C\subseteq L\cup G_Q\)。

3. 一次出现变量的保护。

\CaseDisplay{L\cap C\ne\varnothing\quad\text{或}\quad
   \bigl(\forall x\in S(w)\;\exists u\in R_C:\ m_x(u)\leq1\bigr).
   \quad\textup{(\ensuremath{3})}}

所有存在条件均为实际非空条件，不允许把空和或空积作为项。

\par\smallskip\noindent\textit{必要性.}\quad 
若 \(x\notin L\)，把 \(x\) 赋为 \(a\)、其余变量赋为 \(0\)，则右侧为 \(0\)，左侧为 \(a\)，证明 (1)。

不带 \(H\) 时，若没有内部非线性词，把 \(C\) 内变量赋为 \(b\)、外部变量为 \(0\)。左词为 \(e\)，右侧不超过 \(b\)。带 \(H\) 时，如果 \(R_C=L\cap C=\varnothing\)，同一赋值使右侧为 \(0\)，仍不能吸收 \(e\)；如果 \(R=\varnothing\)，把所有变量赋为 \(b\)，右侧为 \(b\)，左侧为 \(e\)。这证明锚点条件。

不带 \(G\) 时，若 \(x\in C\setminus(L\cup G_C)\)，令 \(x=f\)、其他 \(C\) 内变量为 \(e\)、外部变量为 \(0\)。左词为 \(f\)，右侧不超过 \(e\)。带 \(G\) 时，如果一个 \(x\in C\) 完全不在右側出现，令它为 \(f\)、其他所有变量为 \(e\)，也得到失败赋值。

若保护条件失败，选相应 \(x\in S(w)\)，令 \(x=a\)、其他 \(C\) 内变量为 \(e\)、外部变量为 \(0\)。内部右词均将 \(x\) 重复，因而为 \(0\)；其余右词及所有右一次项亦为 \(0\)。左词为 \(a\)，故失败。

\par\smallskip\noindent\textit{语义充分性及基的有效性.}\quad 
上述三个条件充分，可按非线性左词的四种可能值 \(0, a, e, f\) 核对。

\textbullet{} 在 \(C, C_G\) 中，\(0\) 是加法最小元。左词为 \(e\) 时，内部非线性锚点提供至少 \(e\)；左词为 \(f\) 时，覆盖取 \(f\) 的变量的右词或一次项提供 \(f\)。在 \(C\) 中该右词须在内部；在 \(C_G\) 中，\(f\) 的吸收性允许任意外部变量。左词为 \(a\) 时，内部一次项不小于 \(a\)，或保护条件提供一个值为 \(a, e, f\) 的内部右词。

\textbullet{} 在 \(C_H\) 中，\(b\) 与 \(0, a\) 不可比，但 \(b+0=b+a=e\)。当左词为 \(e\) 而没有内部锚点时，内部一次项为单位层值；任一非线性右词或者已在单位层不低于 \(e\)，或者与该一次项的和不低于 \(e\)。左词为 \(f\) 时仍由内部覆盖处理。左词为 \(a\) 时，若内部一次项为 \(b\)，覆盖取 \(a\) 的左变量的内部右词提供 \(0\) 或 \(a\)，其和至少为 \(e\)；其余情形与上一项相同。左词为 \(0\) 时，覆盖某个导致零的左变量提供一个零层或 \(a\) 层的右项，因此右侧不小于 \(0\)。

\textbullet{} 在 \(C_{HG}\) 中，单位支 \(\{b, e\}\) 与零支 \(\{0, a\}\) 的任意跨支和均为 \(f\)。锚点给出一个单位支右项，或一个内部单位支一次项；这足以配合含所覆盖变量的任意右词，处理左词为 \(e, f\) 的情形。左词为 \(a\) 时，若保护词只是不含那个取 \(a\) 的变量，则它提供单位支值，覆盖又提供含该变量的零支值，二者之和为 \(f\)。如果保护词含该变量一次，它直接提供 \(a\)。有内部一次项的情形同理。左词为 \(0\) 时，全局覆盖某个导致零的变量提供零支右项。

因此判据充分。将每条基的不等式视为 \(w\leq Q\)，逐条检查三个有限集合条件，证明所有列出的恒等式在相应半环中有效；(S2) 的双向吸收同样满足判据。

\par\smallskip\noindent\textit{3. 从基推出全部有效吸收式.}\quad 
由 (S2) 得 \(x^3\approx x^2\)。反复对不同的一次出现变量使用 (S2)，每个非线性词都化为以下规范词的非空和：

\CaseDisplay{H(C)=\prod_{y\in C}y^2,\qquad
xH(C\setminus\{x\})\quad(|C|\ge2).}

同时 (S2) 给出 \((uv)^2\leq uv\)，所以每个非线性右词均吸收其平方。对稳定块 \(p^2, q^2\)，(S3) 给出

\CaseDisplay{p^2q^2\leq p^2+q^2.\quad\textup{(\ensuremath{4})}}

还由 (S2)、(S3) 及交换性得到

\CaseDisplay{p^2y^2\leq p^2y\leq p^2+y.\quad\textup{(\ensuremath{5})}}

稳定锚点。 从 \(R_C\) 中取一个右词 \(u\)，便有 \(u^2\leq Q\)。如果使用的是带 \(H\) 的第二种锚点，取 \(y\in L\cap C\) 和任一非线性右词 \(v=v_1v_2\)，由 \(H\) 得

\CaseDisplay{y^2\leq y+v_1v_2\leq Q.}

这样得到一个内容包含于 \(C\) 的非空稳定块。

不带 \(G\) 时，用 (4) 合并内部右词的平方，用 (5) 加入内部一次项变量。覆盖条件保证最终得到 \(H(C)\leq Q\)。带 \(G\) 时，也可从任意非线性右词取出 \(z\in C\)；把该右词写成 \(zt\)，其中 \(t\) 非空，由 \(G\) 得

\CaseDisplay{H(S)z^2\leq H(S)+z^2t^2\leq Q.
\quad\textup{(\ensuremath{6})}}

对只含 \(z\) 的非线性词，先约化为 \(z^2\)，取 \(t=z\)。全局覆盖因而同样给出 \(H(C)\leq Q\)。

若 \(L\cap C\ne\varnothing\)，选 \(y\) 在其中，把任一非线性左词写成 \(yv\)。由 (S4) 有

\CaseDisplay{w\leq y+w^2\approx y+H(C)\leq Q.
\quad\textup{(\ensuremath{7})}}

故余下只需处理没有内部一次项且左规范词为 \(xH(C\setminus\{x\})\) 的情形。

不带 \(G\) 的情形，保护条件、内部覆盖以及 (S3)、(S5)，按《交换近幂等半环-两个任意变量数有限基》第 5 节的三种合并逐步构造该词。该论证所用 \(u^2\leq u\) 只施于非线性右词，已由 (S2) 保证，无须假设 \(x^2\leq x\)。如果保护词不含 \(x\)，内部覆盖提供重复 \(x\) 的内部右词作为初始辅助见证。

带 \(G\) 时，先取一个不重复 \(x\) 的内部保护词。如果已有规范词 \(xH(S)\leq Q\)，直接以它开始；否则保护词的平方为不含 \(x\) 的 \(H(S)\leq Q\)，全局覆盖提供一个含 \(x\) 的非线性右词 \(xt\)。用 \(G\) 得

\CaseDisplay{xH(S)\leq H(S)+x^2t^2\leq Q.
\quad\textup{(\ensuremath{8})}}

这里保护词非线性且不含 \(x\)，保证 \(S\ne\varnothing\)。如果某个内部保护词恰含一次 \(x\)，先对它使用规范化，保留仍恰含一次 \(x\) 的分项。

设已得到 \(r=xH(S)\leq Q\)。若还需加入某个 \(z\in C\setminus(S\cup\{x\})\)，由全局覆盖选含 \(z\) 的非线性右词 \(zt\)。由 (S4)、\(G\) 和 \(r^2\leq r\)，有

\CaseDisplay{rz^2\leq r+(rz^2)^2
       =r+r^2z^2
       \leq r+r^2+z^2t^2
       \leq Q.
\quad\textup{(\ensuremath{9})}}

这只把所需的 \(z\) 加到左词，不把辅助词中的外部变量带入左词。逐一加入便得到目标规范词。

所有有效单词吸收式都已从相应有限基推出；多项式恒等式分解为双向逐词吸收，故四组基完整。证明没有变量数或词长上限。

\subsection{Two-witness finite bases from prefix content and a final variable}\label{case:135}
Sources 4721 and 4722. The signature is (+, multiplication) without constants, relative to the additively idempotent semiring axioms. A one-letter word means length one; all other words have length at least two.

\par\smallskip\noindent\textit{1. Algebras, values and normal words.}\quad 
The common multiplication in order 0,\allowbreak{}1,\allowbreak{}2,\allowbreak{}3,\allowbreak{}4 is

{\small\par\smallskip
\CaseLine{\ttfamily 0 0 2 3 4}
\CaseLine{\ttfamily 0 0 2 3 4}
\CaseLine{\ttfamily 3 3 3 3 3}
\CaseLine{\ttfamily 3 3 3 3 3}
\CaseLine{\ttfamily 0 0 2 3 4}
\smallskip}

Addition is join in \(3<1\),\(2<0<4\) for 4721, with 1 and 2 incomparable, and in \(3<2<1<0<4\) for 4722. Put B=\{0,\allowbreak{}1,\allowbreak{}4\}. A word of length at least two is 3 if some prefix factor is 2 or 3, or its final factor is 3. Otherwise every prefix factor is in B, and the value is 0,\allowbreak{}2,\allowbreak{}4 according as the final factor is 0/\allowbreak{}1,\allowbreak{}2,\allowbreak{}4. This follows inductively from the table. Associativity follows from these same cases for three factors. Rows and columns preserve joins in both orders.

The formula proves

\CaseDisplay{xyz\approx yxz,\qquad x^2y\approx xy.\quad\textup{(\ensuremath{W})}}

The first law permits adjacent swaps before the final position using nonempty suffix substitution and ordinary left contexts. The second removes repeated prefix letters once brought adjacent, with a nonempty suffix. Hence every word of length at least two normalizes to its sorted distinct prefix variables followed by its original final variable. The final variable may also occur in that prefix, and that occurrence must remain. The normal prefix is nonempty. One-letter words stay one-letter.

\par\smallskip\noindent\textit{2. Exact absorption criterion.}\quad 
Let \(C=c(w)\), \(t=t(w)\), I=c(w without its final position). A one-letter y is absorbed by P exactly when literal y occurs. A word w of length at least two is absorbed exactly when:

\textbullet{} N: an original right word of length at least two is contained in C;

\textbullet{} T: an original right word contained in C ends in t, including possibly literal t;

\textbullet{} A: an original right word of length at least two is contained in C and has prefix content contained in I; or literal t occurs (4721), or some literal from C occurs (4722).

Necessity. If N fails assign \(C=1\) and outsiders 3. The target is 0 and \(P\leq 1\). If T fails assign \(t=4\), C minus t the value 1, outsiders 3. The target is 4; contained longer words end outside t and are 0, noncontained longer words are 3, and literals are 1/\allowbreak{}3. Thus \(P\leq 0\). Once N holds, failure of A forces t not in I: otherwise \(I=C\) and N itself supplies A. Assign \(I=1\), \(t=2\), outsiders 3. The target is 2. Every right word of length at least two is 3, since any contained one has t internally. In 4721 literal t is absent, so \(P\leq 1\) and 2 is not below 1; in 4722 all literals in C are absent, so \(P=3\). For a one-letter y with no literal witness, assign \(y=2\) and outsiders 3. Every longer word has a2/\allowbreak{}3 prefix factor, so \(P=3\). This explains why pure powers do not witness one-letter absorption:\(2^{2}=3\).

Sufficiency. Target 3 is least. For target 0 all target variables are in B; N supplies 0 or 4. For target 4 they are in B and \(t=4\); T supplies 4. For target 2 all prefix variables are in B and \(t=2\). A length-at-least-two A witness has prefix in I and hence value 0,\allowbreak{}2 or 4. Literal t gives 2; in 4722 any other allowed literal is in B and is also at least 2. These cases prove the criterion for arbitrarily many variables and original summands.

\par\smallskip\noindent\textit{3. Two original witnesses.}\quad 
Choose an original T witness v. If it has length at least two, it also supplies N; add one A witness. If it is literal t, it also supplies A in both algebras; add one N witness. The resulting subsum retains every condition. Thus every valid absorption of a longer word uses at most two original right words. A one-letter target is already a summand.

\par\smallskip\noindent\textit{4. Reversible prefix-incidence compression.}\quad 
A local step has at most three normal words in total. Reserve their final variables individually, including one-letter variables: at most three. Group all remaining variables by their nonempty prefix-incidence patterns among the three words: at most seven. Identify each group to one fresh variable and normalize by(W). Validity survives substitution. At most 10 variables occur and normal word length is at most 11.

For the inverse substitution replace each group variable by a nonempty product listing its original members. No group occupies a final position. Its prefix incidences restore all original prefix sets, while the individually reserved final variables keep their positions and any prefix occurrence. Prefix permutation and idempotence recover each original normal word. One-letter words remain individual reserved variables. No empty substitution is used.

\par\smallskip\noindent\textit{5. Complete finite bases.}\quad 
For each source take(W) plus all valid normal absorptions on a fixed 10-letter alphabet, word length at most 11, with at most two right words. These are finite explicitly specified families, with validity decided by the respective criterion. No small enumerated or minimal list is claimed.

Theorem. These families form relative identity bases for their sources. Direct witness selection reduces every valid word absorption to a local step. Compression places it in the finite family; the inverse substitution and(W) derive the original step. Adding unused summands gives the whole absorption. Distributivity reduces general identities to sums of words, and deriving each summand's absorption on both sides proves unrestricted completeness. Word reversal gives bases and derivations for the respective opposites.

Finite checks support this mathematical proof. Transfer to another catalogue algebra requires equality in both complete HSP directions.

\subsection{Nine identities for 4728 with a nonlinear terminal guard}\label{case:136}
Signature (+, multiplication), without constants, relative to the additively idempotent semiring axioms. The proof is unrestricted in variables, word length and number of summands.

\par\smallskip\noindent\textit{1. Exact criterion derived from the tables.}\quad 
Addition has rows 00004,\allowbreak{}01014,\allowbreak{}00224,\allowbreak{}01234,\allowbreak{}44444. Thus 3 is least, 1 and 2 are incomparable, \(1+2=0\), and 4 is greatest. Multiplication has rows 00234,\allowbreak{}00234,\allowbreak{}33333,\allowbreak{}33333,\allowbreak{}44444. A literal retains its value. A nonlinear word with only 0/\allowbreak{}1 factors has value 0. Otherwise its first factor outside \{0,\allowbreak{}1\} determines the value:4 gives 4,\allowbreak{}3 gives 3, and 2 gives 2 when final but 3 otherwise. This follows by table induction.

Let N be the original nonlinear right words. For an original polynomial P, write E(x, S) for an original word containing x with first-x prefix content contained in S. A witness is unmarked when its first x is final. A nonlinear absorption \(w\leq P\) is valid exactly when:

(N) Some word of N has content contained in c(w). (E) Each first target occurrence x has an E(x, S) witness, where S is the preceding target content. (T) At a simple final variable x, either there is an unmarked E-witness, or an original NONLINEAR word has content contained in S.

A literal target x is absorbed exactly when x itself is an original literal summand. These statements are proved directly as follows. Failure of N is separated by target \(variables=1\), \(outsiders=3\): the target is 0, all right values are 1 or 3. Failure of E at x is separated by \(S=0\), \(x=4\), \(others=3\): target 4, no right 4. Failure of T is separated by \(S=1\), \(x=2\), \(others=3\): target 2, while every right value is 1 or 3. In particular a contained literal prefix guard is insufficient, because 1 does not bound 2. A missing literal x is separated by \(x=2\), \(others=3\).

For sufficiency,3 is least. An all 0/\allowbreak{}1 nonlinear target is supplied by N's contained nonlinear word, also valued 0. A first 4 is supplied by E. The only remaining case is a first 2 at the simple final position; an unmarked E-witness gives 2, or T's nonlinear contained-prefix word gives 0. This exhausts assignments and proves the criterion with no bounds.

\par\smallskip\noindent\textit{2. The nine identities.}\quad 
{\small\par\smallskip
\CaseLine{\(W1 xxy=xy\qquad W2 xyx=xyy\)}
\CaseLine{\(A0 xyz\leq yz+x\qquad A1 xyz\leq yz+xt\)}
\CaseLine{\(B0 xyz\leq xz+y\qquad B1 xyz\leq xz+xyt\)}
\CaseLine{\(U  xyz\leq xz+xy\)}
\CaseLine{\(R0 xy\leq xx+xyy\qquad R1 xy\leq y+xyy\)}
\smallskip}

The word formula proves W1, W2. The exact criterion proves every absorption: their nonlinear seed and ordered witnesses are displayed on the right. In R0 the terminal guard xx is nonlinear; in R1 the literal y is an unmarked terminal witness. Every law can independently be evaluated on all assignments in the finite table.

W1, W2 normalize a word to its first-occurrence list F when its final variable is simple, and otherwise to F with its final listed variable doubled. Internal separated repeats reduce by \(xYxZ=xY^{2}Z=xYZ\) for nonempty Y, Z; adjacent repeats use W1. Terminal separated repeats use W2 before internal deletions. A pure nonlinear power becomes xx, not the literal x. These reductions preserve N, E, T and original witness categories. All identity substitutions are nonempty.

\par\smallskip\noindent\textit{3. Construct the saturated companion.}\quad 
Let \(x_{1}\)...\(x_{d}\) be the target's first-occurrence list. Maintain p \(u\leq P\) for EVERY original u in N, with p the processed list. Initially p is absent context, not a word substituted for a variable.

At the first x, an original witness is either literal x or xT with T nonempty. For each \(u=YZ\) in N use A0 or A1 respectively to derive x \(u\leq u+x\) or x \(u\leq u+xT\). This starts the invariant.

At a later x with nonempty prefix p, an original literal x uses B0: p x \(u\leq p\) u+x. An original nonlinear marked witness v has p \(v=p\) x T modulo W1, W2, with T nonempty; B1 gives p x \(u\leq p\) u+p x T. An original nonlinear unmarked witness satisfies p \(v=p\) x; U gives p x \(u\leq p\) u+p x. These cases exhaust E and extend the invariant through all target variables, including a simple tail.

Choose an original contained nonlinear seed v from N. The resulting \(p_{d}\) v has normal form \(p_{d}\) with its final variable doubled. This is the target's saturated companion, now below P. A saturated target therefore already follows.

\par\smallskip\noindent\textit{4. Recover the simple tail.}\quad 
For target q s with simple final s, q is nonempty and its saturated companion q s s is available. Earlier prefix invariants remain available. By T, one of the following holds:

\textbullet{} An original nonlinear unmarked witness v gives q \(v=q\) s directly.

\textbullet{} The original literal s uses R1: q \(s\leq s+q\) s s.

\textbullet{} An original NONLINEAR word v is contained in c(q). The invariant gives q \(v=q^{2}\), and R0 yields q \(s\leq q^{2}+q\) s s.

These are exhaustive. A literal prefix guard cannot replace the third case: \(xy\leq x+xyy\) fails at \(x=1\), \(y=2\). No law from the chain-addition type 4729 allowing that replacement is used here. Literal targets are already summands.

Every valid absorption is now derived from nine identities. Expand terms by distributivity and apply the summand absorptions in both directions to derive every identity. Thus the nine laws form a finite basis for

4728. Reversal gives its opposite. Further catalogue transfers require complete HSP both ways. Finite tests and literal replay supplement, rather than establish, the unrestricted completeness argument.

\subsection{Thirteen identities for 4729: nonlinear seeds and literal terminal recovery}\label{case:137}
Signature (+, multiplication), without constants. All bases below are relative to the additively idempotent semiring axioms.

\par\smallskip\noindent\textit{1. Word values and exact absorption criterion.}\quad 
In the catalogue order 0,\allowbreak{}1,\allowbreak{}2,\allowbreak{}3,\allowbreak{}4 the tables are:

{\small\par\smallskip
\noindent addition         multiplication\par
\CaseLine{\ttfamily 0 0 0 0 4\qquad 0 0 2 3 4}
\CaseLine{\ttfamily 0 1 1 1 4\qquad 0 0 2 3 4}
\CaseLine{\ttfamily 0 1 2 2 4\qquad 3 3 3 3 3}
\CaseLine{\ttfamily 0 1 2 3 4\qquad 3 3 3 3 3}
\CaseLine{\ttfamily 4 4 4 4 4\qquad 4 4 4 4 4}
\smallskip}

The additive chain is \(3<2<1<0<4\). A literal retains its assigned value. For a nonlinear word scan the first factor outside \{0,\allowbreak{}1\}. If none exists, the value is 0. A first 4 gives 4; a first 3 gives 3; a first 2 gives 2 only at the final position, and otherwise gives 3. This follows directly by induction from the multiplication rows.

For original nonempty polynomial P define N as its original words of length at least two. Write E(x, S) when an original word contains x and the content before its first x is contained in S. Call this witness unmarked if that first x is final. Write Z(S) when an original word has content contained in S; Z may be a literal.

A literal x is absorbed exactly when x itself belongs to P. For nonlinear w, the exact criterion is:

1. Some original NONLINEAR right word has content contained in c(w).

2. For each first occurrence of x in w, E(x, S) holds, where S is the content preceding that occurrence. Literal E witnesses are allowed.

3. If the final variable x is simple, either an unmarked E(x, S) holds or Z(S) holds, where \(S=c(w)\) minus \{x\}.

Necessity of (1): assign all target variables 1 and outsiders 3. The target is 0; only a contained nonlinear summand can reach 0. Contained literals remain 1. Necessity of (2): assign \(S=0\), \(x=4\), \(others=3\). The target is 4, and a right word is 4 exactly when it supplies E(x, S). Necessity of (3): assign \(S=0\), \(x=2\), \(others=3\). A right value reaches 2 exactly through an unmarked E witness or a contained-prefix word. A missing literal x is separated by \(x=2\) and all other \(variables=3\).

For sufficiency, value 3 is least. An all-\{0,\allowbreak{}1\} target, valued 0, is covered by the contained nonlinear seed. If its first nonbase value is 4, the corresponding E witness evaluates to 4. The only remaining case is a first 2 at the simple final position. An unmarked E witness has value 2; a Z(S) word has value 0 or 1. Thus each condition is sufficient for every assignment, with no bound on variables, lengths or summands.

\par\smallskip\noindent\textit{2. The thirteen laws.}\quad 
{\small\par\smallskip
\CaseLine{\(W1  xxy = xy\)}
\CaseLine{\(W2  xyx = xyy\)}
\CaseLine{\(A0  xyz \leq  yz+x\)}
\CaseLine{\(A1  xyz \leq  yz+xt\)}
\CaseLine{\(B0  xyz \leq  xz+y\)}
\CaseLine{\(B1  xyz \leq  xz+xyt\)}
\CaseLine{\(U\qquad xyz \leq  xz+xy\)}
\CaseLine{\(R0  xy \leq  xx+xyy\)}
\CaseLine{\(R1  xy \leq  y+xyy\)}
\CaseLine{\(R2  xy \leq  x+xyy\)}
\CaseLine{\(R3  xyz \leq  x+xyzz\)}
\CaseLine{\(R4  xyz \leq  y+xyzz\)}
\CaseLine{\(R5  xyzt \leq  y+xyztt\)}
\smallskip}

Here \(a\leq b\) means \(a+b=b\). The word formula proves W1, W2 and the exact criterion proves every displayed absorption. No identity variable is ever substituted by an empty word.

The normal form of a word is its first-occurrence list F if the final variable is simple; otherwise it is F followed by the last letter of F. To see this using the two word laws, remove adjacent internal repetitions by W1. For separated internal repeats, \(xYxZ=xY^{2}Z=xYZ\), with both Y, Z nonempty, using W2 then W1. A final repeated x in p x Y x becomes p x \(Y^{2}\) by W2; remove internal repetitions, leaving a duplicate of the last first-occurrence letter. Adjacent final duplicates reduce to two using W1. This also treats powers. Literals stay literal.

\par\smallskip\noindent\textit{3. Prefix construction on original nonlinear words.}\quad 
Normalize the target and original summands. The conditions and nonlinear seed are preserved. Let its first-occurrence list be \(x_{1}\)...\(x_{d}\). Maintain p \(u\leq P\) for EVERY u in N, where p is the list processed so far. Initially the empty p is notation for no prefix, never a substitution.

At the first letter x, E(x, empty) supplies either literal x or a word x t with t nonempty. Write each u in N as YZ, Y, Z nonempty. A0 derives x u from u+x, or A1 derives x u from u+x t. Thus the invariant starts.

For a later letter x and nonempty p, a literal witness x uses B0: p x \(u\leq p\) u+x. A nonlinear marked witness v has p \(v=p\) x t modulo W1, W2, with t nonempty, so B1 gives p x \(u\leq p\) u+p x t. A nonlinear unmarked witness gives p \(v=p\) x; U gives p x \(u\leq p\) u+p x. These are all E possibilities. The invariant therefore extends through ALL letters, including a simple final letter.

Take a contained nonlinear original seed r. For the full list \(p_{d}\), \(p_{d}\) r has the saturated normal form \(p_{d}\) followed by its last letter. Consequently every valid nonlinear target has its saturated companion derivable. If the target is saturated, this already finishes the proof.

\par\smallskip\noindent\textit{4. Recovering a simple final letter.}\quad 
Let the normalized target be q s with q nonempty and s absent from q. We have its saturated companion q s s. All earlier prefix invariants remain available. Condition (3) gives the following exhaustive cases.

\textbullet{} An original nonlinear unmarked witness v satisfies q \(v=q\) s and was already derived while processing q.

\textbullet{} An original literal s uses R1 with \(x=q\), \(y=s\) and the companion q s s.

\textbullet{} A contained-prefix nonlinear word r gives q \(r=q^{2}\). R0 with \(x=q\), \(y=s\) recovers q s from \(q^{2}\)+q s s.

\textbullet{} A contained-prefix literal a lies at a unique position of \(q=A\) a B. If \(A=B=empty\), use R2 with \(x=a\), \(y=s\). If \(A=empty\) and B nonempty use R3 with \(x=a\), \(y=B\), \(z=s\). If \(B=empty\) and A nonempty use R4 with \(x=A\), \(y=a\), \(z=s\). If both are nonempty use R5 with \(x=A\), \(y=a\), \(z=B\), \(t=s\). In each case the other right summand is exactly the saturated companion, so no square of the literal is assumed derivable.

This last distinction is essential: \(xx\leq x+xy\) and the literal guards of 7953 are false here. A literal contained seed alone cannot replace condition (1). We derive the saturated companion from an original NONLINEAR seed before applying any terminal literal recovery.

\par\smallskip\noindent\textit{5. Completeness and catalogue transfers.}\quad 
Every valid word absorption has now been derived from thirteen fixed identities. The induction does not bound variables, word lengths, polynomial size or witness prefixes. Distributivity expands arbitrary terms into nonempty polynomials, and term identities follow by absorbing each summand in both directions. Hence these thirteen laws form a finite basis for 4729. Reversal gives the opposite basis. Any further catalogue transfer requires complete independently checked HSP inclusions in both directions; no finite tests or one-way membership imply the conclusion.

\subsection{A finite basis for 5519 from quadratic components and bounded witnesses}\label{case:138}
Signature (+, multiplication), without constants. The bound below is deliberately generous. It specifies a finite effective family, not an enumerated or small basis.

Addition rows 00444,\allowbreak{}01444,\allowbreak{}44244,\allowbreak{}44434,\allowbreak{}44444 have proper principal ideals \{1\},\{0,\allowbreak{}1\},\{2\},\{3\}. Multiplication rows 00444,\allowbreak{}00344,\allowbreak{}22444,\allowbreak{}33444,\allowbreak{}44444 give these word values. Literals retain their value. Nonlinear words give 0 when every factor is 0/\allowbreak{}1;2 or 3 when the head is respectively 2 or 3 and all later factors are 0/\allowbreak{}1;3 also when the head is 1, second factor 2, and all later factors 0/\allowbreak{}1. Every other nonlinear word is 4. Induction proves this formula for arbitrary lengths.

\par\smallskip\noindent\textit{Exact ideal test.}\quad 
For original P let C be its content union, L its literal variables, A its nonlinear heads, B its second-position variables, D its variables in positions 3 onward, \(H=L\) union A, and \(T=B\) union D. Record edges (h, s) from the first two positions of original nonlinear words.

Join two heads when they have a common second variable; take connected components. Each component has two possible modes for \(P=3\): S: its heads are 1 and its seconds are 2; N: its heads are 3 and its seconds are 0/\allowbreak{}1. Let forced-S heads be A intersect T. Let forced-N heads be A intersect L, together with heads having an outgoing second in H union D. The constraint \(P=3\) is impossible iff L intersects T, or one component contains both a forced-S and a forced-N head.

To prove this, each original nonlinear word has value 3 exactly when its first two values are(1,\allowbreak{}2) or(3,\allowbreak{}0/\allowbreak{}1), with all later values 0/\allowbreak{}1. Shared seconds force heads into the same mode. A head occurring as a second or later variable must be 1, forcing S. An outgoing second which is itself a head or lies in D cannot be 2, forcing N. Literal heads must be 3; literal variables cannot occur noninitially. These are all constraints:choose either allowed mode independently in each component, set literal variables 3, and give remaining constrained 0/\allowbreak{}1 variables either0 or 1. Overlapping roles are consistent because the forced-mode rules already account for them.

When this system is consistent, let Z be the variables necessarily 3: L and heads in forced-N components. Let E be variables necessarily 2: seconds of forced-S components. Let K be variables necessarily 0/\allowbreak{}1: D, heads in forced-S components, and seconds of forced-N components. These descriptions are exact; unconstrained component modes and unconstrained0/\allowbreak{}1 choices vary independently.

Under \(P=3\), a literal w is always 3 iff its variable belongs to Z. A nonlinear target with first positions h, s is always 3 iff every variable after position 2 belongs to K and its first pair satisfies:

h, s belong to the same original head/second component; or h belongs to Z and s to K; or h is a head in a forced-S component and s belongs to E.

The first alternative means an edge in the complete bipartite completion of the component, with h a head and s a second. For a free component, its two modes show that no other second is universally compatible:another free component varies independently, forced seconds fail one mode, and heads in the same component fail its N mode. Forced-head cases give the second and third alternatives directly. This proves the pair criterion, including literal and overlapping roles.

Consequently \(w\leq P\) holds exactly when all of these tests hold:

C: content(w) is contained in C; N: w is literal, or P has a nonlinear original; I2: H intersects T, or head(w) belongs to H and tail-content(w) is contained in T; I3: \(P=3\) is impossible, or w is universally 3 by the preceding test.

C handles the ideal\{0,\allowbreak{}1\}, because original values 0/\allowbreak{}1 force every original variable into\{0,\allowbreak{}1\}. N handles\{1\}, which can contain a whole polynomial only if every original is literal. For \(P=2\), all original heads must be 2 and all original noninitial variables 0/\allowbreak{}1. This is possible exactly when H and T are disjoint, and the displayed I2 condition is precisely target membership in that ideal.

Explicit failure assignments complete necessity. For failed C, put a missing variable 4 and all others 0. For failed N, put every variable 1. For failed I2, put H variables 2 and others 0; this separates a missing target head. For a missing target-tail variable q, keep it 2 if it belongs to H, or assign it 3 otherwise; all original values remain 2 and the target is 4. For failed I3 choose an allowed component mode and 0/\allowbreak{}1 assignment witnessing the failed literal, pair, or tail condition. All originals are 3 while the target is not 3. The exact independent mode description above guarantees such an assignment. These arguments prove exactness for arbitrary alphabets, word lengths, and summand counts.

\par\smallskip\noindent\textit{Five fixed identities.}\quad 
{\small\par\smallskip
\CaseLine{\(xyzt = xytz\)}
\CaseLine{\(xyzz = xyz\)}
\CaseLine{\(xy \leq  xyz\)}
\CaseLine{\(xz \leq  xy+uy+uz\)}
\CaseLine{\(xyzt \leq  xyz+xyt\)}
\smallskip}

The first two preserve the first two literal positions and the set of remaining variables. They normalize every long word by sorting its tail after position 2 and deleting adjacent tail copies. Each swap has two nonempty protected prefix blocks; each deletion retains one tail occurrence. Literals and quadratics stay unchanged.

The word formula proves prefix extraction \(xy\leq xyz\). The rectangle rule follows from the exact ideal test:its three original edges connect the head x to second z in one component; C, N, I2 also hold. The last rule merges tail sets with fixed first two positions and is sound by all four ideal tests. These proofs allow arbitrary nonempty word/term substitutions through ordinary identity logic.

Extract the quadratic prefix of each original nonlinear word. Repeated rectangle rules then derive every edge in each complete bipartite component. Along a path of heads h0,..., hn sharing successive seconds b1,..., bn, an edge hn s gives h(n-1)s from h(n-1)bn, hn bn, hn s, and continuing backward gives h0s. All path edges arise from the originals. Denote by P+ the original words and these derived quadratics. Every word of P+ is below P.

\par\smallskip\noindent\textit{Bounded supports for targets of length at most 3.}\quad 
For a head h in a forced-S component, choose a forced-S head r in that component, a second b there, the two completed edges hb, rb, and an original word with a noninitial r. These at most 3 words already force h's component into S.

For a forced-N component choose a forced-N head r, the two completed edges hb, rb, and either the original literal r, or an original edge rs whose second lies in H union D plus one original witnessing that membership. There are at most 4 words. Thus necessary 3 has a support of at most 4 words; necessary 2 has a support of at most 4 words by adjoining an edge to a forced-S support. Necessary 0/\allowbreak{}1 has a support of at most 5 words:one original D occurrence, or a forced-S head support, or an edge to a forced-N head support.

If \(P=3\) is impossible because L meets T, two original words witness it. A conflicting component has a forced-S support and a forced-N support, at most 7 words. Otherwise a universally 3 literal needs at most 4 words. A universally 3 nonlinear first pair needs one completed edge, or a necessary 3 support and a necessary 0/\allowbreak{}1 support (at most 9 words), or forced-S and necessary 2 supports(at most 7). A target of length 3 needs at most 5 further words for its last variable. Hence I3 has a support of at most 14 words for any target of length at most 3.

Separately choose at most 3 original content witnesses for C, one nonlinear original for N when needed, and at most 3 originals for I2:two if H meets T, otherwise a target-head witness and witnesses for its at most 2 tail variables. Adjoining these to the I3 support gives at most 21 words from P+ that already absorb the short target. Every stated support is explicit and retains any needed original literal, noninitial occurrence, or D occurrence. Taking a union cannot invalidate an absorption test. No bounded-variable extrapolation is used in this support bound.

\par\smallskip\noindent\textit{Finite family and inverse compression.}\quad 
Normalize the target and the at most 21 selected upper words. There are at most 22 words. Reserve their first two variables(or the sole literal variable), at most 44 variables. Every other variable occurs only in the commutative idempotent tail after position 2. Partition these by their nonempty incidence patterns among22 words, at most \(2^{22}\)-1 groups. Give each group and each reserved variable a distinct fresh variable and normalize again. This is a substitution into a valid absorption, so it remains valid.

There are at most 4194347 variables and 4194349 positions per word, with at most 21 upper words. These large fixed bounds only establish finiteness. For the inverse substitution, restore reserved variables as singletons and each group as a nonempty product of its members. The expansions regain the first two positions and the original tail sets. Literal applications of the two structural equalities recover all original normal words. Short words contain reserved variables only; every long word retains a nonempty tail throughout. Thus the compression is equationally reversible, including repeated protected variables which also occur in the tail.

Take the five fixed identities and all valid absorptions on a fixed 4194347-variable alphabet, length at most 4194349, and at most 21 upper words. The exact ideal criterion effectively specifies this finite family; it is not enumerated. Inverse compression derives all valid targets of length at most 3 from P+, and hence from P.

For a longer normal target hsZ, each triple hsz for z in Z satisfies the exact criterion whenever hsZ does. Derive all such triples by the bounded-support argument. Starting with one, use the final fixed identity with \(x=h\), \(y=s\), \(z=current\) nonempty tail, \(t=the\) next variable to append variables one at a time. Tail normalization yields hsZ. All substitutions are nonempty. Thus every absorption is derived with no bound on its original variables, length, or number of summands. Distributive expansion and both absorption directions give every valid identity. Therefore 5519 and its opposite are finitely based; catalogue transfers require independently complete mutual HSP.

\subsection{A finite basis for 5526 from parity paths and bounded supports}\label{case:139}
Signature (+, multiplication), without constants. The finite family below is effective and deliberately large; it is not an enumerated small basis. All completeness arguments allow arbitrarily many variables, arbitrary nonempty word lengths, and arbitrary polynomial sizes.

Addition rows 00444,\allowbreak{}01444,\allowbreak{}44224,\allowbreak{}44234,\allowbreak{}44444 have proper principal ideals \{1\},\{0,\allowbreak{}1\},\{3\},\{2,\allowbreak{}3\}. Multiplication rows 00444,\allowbreak{}00424,\allowbreak{}22444,\allowbreak{}33444,\allowbreak{}44444 give the following values by induction on word length. A literal retains its value. A nonlinear word is 0 if all factors are 0/\allowbreak{}1; it is its head value 2 or 3 if that head is 2/\allowbreak{}3 and all later factors are 0/\allowbreak{}1; it is 2 if the head is 1, the second factor 3, and all later factors 0/\allowbreak{}1. Otherwise it is 4. This special pair(1,\allowbreak{}3) differs from 5519's pair(1,\allowbreak{}2).

\par\smallskip\noindent\textit{Original constraints and exact semantics.}\quad 
For an original polynomial P let C be its content union, L its literal variables, A its nonlinear heads, B its second-position variables, and D its variables at positions 3 onward. Write \(H=L\) union A, \(T=B\) union D. Make a directed edge \(h\to s\) from each original nonlinear word's first two positions; write E for this set. Regard these also as an undirected graph G, retaining loops. Every graph vertex belongs to A union B.

Assignments with all original values in\{2,\allowbreak{}3\} are described exactly as follows. Color every variable high or low. Variables in L must be high; variables in D must be low. Each edge has opposite endpoint colors. A low variable in A must have actual value 1; any other low variable may be 0 or 1. A high variable in B must have actual value 3; any other high variable may be 2 or 3. These conditions are necessary by the word-value formula, and sufficient for each original word independently. They therefore describe every such assignment, including overlapping roles.

This system is inconsistent exactly when G has an odd cycle, a variable lies in L intersect D, or two pins in one bipartite component disagree with their path parity. In a consistent system, let U and K be the variables forced high and forced low, respectively. These include isolated literal and D variables. In an unpinned component either coloring is possible, independently of other components. Actual0/\allowbreak{}1 and 2/\allowbreak{}3 choices are independent except for the explicit A and B singleton restrictions.

Let E* contain(a, b) whenever a belongs to A, b belongs to B, and G has an odd-length walk from a to b. For a consistent system these are exactly the opposite-side pairs in the same bipartite component. Under every original\{2,\allowbreak{}3\}-assignment, a literal x is in\{2,\allowbreak{}3\} exactly when x belongs to U. A nonlinear word hsZ is in\{2,\allowbreak{}3\} exactly when every variable of Z belongs to K and at least one of these holds:

1. (h, s) belongs to E*;

2. h belongs to U and s belongs to K;

3. h belongs to A intersect K and s belongs to B intersect U.

Indeed, the first two values must either be(high, low), or precisely(1,\allowbreak{}3). If their colors are not individually fixed, opposition for all colorings requires the same component and opposite parity. The low-head possibility then requires h in A, and the high-second possibility requires s in B. If both colors are fixed, the two displayed forced cases are exact. Independent actual-value choices prove the necessity of the A/B roles.

Consequently \(w\leq P\) holds if and only if all four tests below hold:

C: content(w) is contained in C; N: w is literal, or P has an original nonlinear word; I3: H intersects T, or head(w) is in H and all noninitial variables of w are in T; I2: the coloring system is inconsistent, or w meets the preceding universal\{2,\allowbreak{}3\} test.

The ideal\{0,\allowbreak{}1\} gives C:an original word lies in it precisely when all its variables are 0/\allowbreak{}1. The ideal\{1\} gives N:a nonlinear word never has value 1. To have every original in\{3\}, all heads must be 3 and all noninitial variables 0/\allowbreak{}1; this is possible exactly when H and T are disjoint, and gives I3. The full color description gives I2. Since these are all proper principal ideals, testing them proves sufficiency:apply the test for the actual value P takes under any assignment.

Failure assignments also prove necessity without any variable bound. For failed C, set a missing variable 4 and all others 0. For failed N, set all variables 1. For failed I3, set H variables 3 and the others 0; this separates a missing head. A missing noninitial target variable lies outside T; set it 3 if in H, and 2 otherwise. Original values remain 3 and the target leaves\{3\}. For failed I2, choose an allowed coloring and actual values violating the target condition as described above. Such a choice exists by independent component colorings and independent remaining value choices. All original values lie in\{2,\allowbreak{}3\}, and the target lies outside that ideal.

\par\smallskip\noindent\textit{Thirteen fixed identities and arbitrary parity-path shortening.}\quad 
Use the following five laws:

{\small\par\smallskip
\CaseLine{\(xyzt = xytz\)}
\CaseLine{\(xyzz = xyz\)}
\CaseLine{\(xy \leq  xyz\)}
\CaseLine{\(xyzt \leq  xyz+xyt\)}
\CaseLine{\(xy \leq  xu+yx+vy\)}
\smallskip}

Add eight path laws. For each independent orientation of xy, yz, zt, take those three oriented words on the right of \(xt\leq\)... . If the first word is yx, also include xu. If the last word is tz, also include vt. The eight explicit right sides are:

{\small\par\smallskip
\CaseLine{\(xy+yz+zt\)}
\CaseLine{\(xy+yz+tz+vt\)}
\CaseLine{\(xy+zy+zt\)}
\CaseLine{\(xy+zy+tz+vt\)}
\CaseLine{\(yx+yz+zt+xu\)}
\CaseLine{\(yx+yz+tz+xu+vt\)}
\CaseLine{\(yx+zy+zt+xu\)}
\CaseLine{\(yx+zy+tz+xu+vt\)}
\smallskip}

All substitutions are nonempty. The first two equalities preserve the first two literal positions and the set of later variables. They sort and remove repetitions only after position 2. The word-value formula proves prefix extraction and tail merging. For each path law, the three edges force opposite colors at x, t. The guards ensure x has an original head role and t an original second role, so I2 holds. C, N hold directly. For I3, if the original H, T sets intersect it is automatic; otherwise the guards or terminal oriented edges give x in H, t in T. The same argument proves the reversal law from edge yx and its two role guards. Thus all thirteen identities are valid by the unrestricted ideal test.

Extract every original quadratic edge by \(xy\leq xyz\). We next derive each edge in E* by the fixed path and reversal laws. Keep all original role guards available throughout. Take an odd walk v0,..., vn from a in A to b in B. For \(n=1\) either its orientation is already available or the reversal law supplies \(a\to b\) using original \(a\to u\) and \(v\to b\) guards. For \(n=3\) the endpoints allow one path law.

For odd \(n\geq 5\), there is a consecutive three-edge segment whose endpoints have eligible A, B roles in either order. Otherwise each distance-three endpoint pair must consist of two A-only vertices or two B-only vertices:if either endpoint had both roles, or their exclusive roles differed, the segment would be eligible. Applying this to positions (0,\allowbreak{}3),(1,\allowbreak{}4),(2,\allowbreak{}5) makes vertices 0,\allowbreak{}1,\allowbreak{}2,\allowbreak{}3 exclusive-role vertices, with 0 and 3 of the same role. But each of their three consecutive edges joins opposite exclusive roles, so 0 and 3 must have opposite roles, a contradiction. Apply the appropriate path law, with original role guards if needed, and replace this segment by its derived oriented edge. This shortens the undirected walk by 2. Repeat to length 1. The walk may repeat vertices; the argument and the nonempty substitutions still apply.

An odd-cycle component contains a vertex r in A intersect B. Otherwise the exclusive roles A-only and B-only would bipartition the graph. There is an odd closed walk from r to r, so the preceding derivation gives rr. That single derived quadratic has no\{2,\allowbreak{}3\}-assignment. This handles odd cycles of every length, not merely finitely many test cycles.

Write P+ for the original words together with all derived E* quadratics. Every word of P+ is below P. In a consistent component, the undirected E* graph has diameter at most 4. If there is a vertex r in both roles, it is adjacent in E* to EVERY vertex on the opposite side. Every vertex on r's own side has some original neighbor and hence is at distance at most 2 from r. Thus any two vertices are at distance at most 4. If there is no mixed-role vertex, connectedness forces one side entirely A-only and the other B-only, and E* is complete bipartite, diameter2. Every E* edge preserves original parity, so these short paths preserve all high/low consequences.

\par\smallskip\noindent\textit{At most 24 witnesses for short targets.}\quad 
All witness words in this section belong to P+. A forced color of a variable x has a support of at most 5 words:choose a pin q in the same component, a path of at most 4 E* edges from q to x, and the original literal q or original D-occurrence word that supplies the pin. If x itself has the relevant literal/D pin, one word suffices. Graph-free variables have only those direct pins.

If I2 is inconsistent, one derived square suffices for an odd component; two originals suffice for L intersect D; and otherwise two conflicting pins plus a path of at most 4 derived edges suffice. Thus inconsistency has a support of at most 6 words. A universally high literal needs at most 5 words. A nonlinear target first pair in E* needs one edge. A forced high/low pair needs at most 10 words. A forced low/high pair needs at most 10 words for the colors plus one original head-role witness and one original second-role witness, at most 12 words. Each variable after position 2 needs at most 5 words for its low color. For a target of length at most 3, the I2 condition therefore has a support of at most 17 words.

Choose independently at most 3 original content witnesses, one nonlinear original if necessary, and at most 3 original I3 witnesses. For the last condition, use two originals exhibiting H intersect T, or a target-head witness and witnesses for each of at most 2 noninitial variables. These choices together with the I2 support use at most 24 words and already absorb the short target. Adding further words cannot invalidate an absorption; equivalently, it restricts every relevant ideal-assignment system. All witnesses retain their actual original role/pin occurrences.

\par\smallskip\noindent\textit{Finite family, inverse compression, and unrestricted completeness.}\quad 
Normalize the short target and the at most 24 selected upper words. There are at most 25 words. Reserve their first two variables, or the sole literal variable, at most 50 variables in total. Every other variable occurs solely in the commutative idempotent tail after position 2. Partition the unreserved variables by their nonempty incidence pattern among the 25 words, at most \(2^{25}\)-1 groups. Replace every group by its own fresh variable, keep every reserved variable separate, and normalize. This is a substitution into a valid absorption, hence is still valid.

The compressed absorption has at most 33554481 variables, word length at most 33554483, and at most 24 upper words. Conversely substitute for each group-variable the nonempty product of precisely its group, and for each reserved variable its singleton. Literal applications of the first two laws recover each original normal word:its first two positions and its later-variable set are exact. Short words contain reserved variables only; long words retain a nonempty tail. Repetitions of a reserved variable in a tail are retained in both directions.

Take the thirteen fixed laws and ALL valid absorptions within those three fixed bounds on a fixed alphabet of 33554481 variables. This is a finite effectively specified family; no assertion of feasible complete enumeration or a small basis is made. Inverse compression derives every valid short target from P+, and the parity-path derivations then give it from P. A longer normal word hsZ satisfies the ideal test only if each triple hsz, for z in Z, does. Derive those triples, then merge their tails one at a time using \(xyzt\leq xyz+xyt\) and normalize. The accumulated tail and each new tail block are nonempty. This derives every valid absorption \(w\leq P\) with no original bounds. Expansion by distributivity and the two absorption directions derive every valid identity.

Therefore 5526 is finitely based. Its opposite has the reversed basis. Any additional catalogue types require independent complete mutual HSP certificates; the graph argument alone does not classify similar tables.

\subsection{Twenty-four identities for 4872: saturation, ordered prefixes and simple tails}\label{case:140}
Constant-free signature (+, multiplication), relative to the ai-semiring axioms. Every completeness construction below is unrestricted.

\par\smallskip\noindent\textit{1. Values, normalization and the exact criterion.}\quad 
The catalogue operations, in order 0,\allowbreak{}1,\allowbreak{}2,\allowbreak{}3,\allowbreak{}4, are:

{\small\par\smallskip
\noindent addition        multiplication\par
\CaseLine{\ttfamily 00004\qquad 00234}
\CaseLine{\ttfamily 01114\qquad 01234}
\CaseLine{\ttfamily 01224\qquad 33333}
\CaseLine{\ttfamily 01234\qquad 33333}
\CaseLine{\ttfamily 44444\qquad 44444}
\smallskip}

The additive chain \(is 3<2<1<0<4\). A literal retains its value. In a nonlinear word scan the first factor outside\{0,\allowbreak{}1\}:4 gives 4,\allowbreak{}3 gives 3, and 2 gives 2 only at the last position, otherwise3. If all factors lie in\{0,\allowbreak{}1\}, the value is 0 if some factor is 0, and 1 otherwise. Induction over multiplication proves this rule. The all-base branch is different from 4729; do not identify their criteria.

Let \(C=c(w)\), let \(U_{C}\) be the original right words contained in C, and let E(x, S) mean an original word contains x with first-x prefix content inside S. An E witness is unmarked if the first x is final. Write Z(S) for existence of an original word contained in S.

A literal target is absorbed exactly when it is an original literal summand. A nonlinear w is absorbed exactly when:

1. The union of contents of \(U_{C}\) equals C (literals are allowed).

2. Each first target occurrence x has E(x, S), where S is its preceding target content (literal witnesses are allowed).

3. A simple final variable x has either an unmarked E(x, S) or Z(S).

For necessity of 1, if x is uncovered assign \(x=0\), other target letters 1, outsiders 3. Target 0 exceeds all right values 1 or 3. For 2 assign \(S=1\), \(x=4\), others 3:only an E witness can reach 4. For 3 assign \(S=1\), \(x=2\), others 3:only an unmarked E or Z(S) can reach 2. A missing literal x is separated by \(x=2\), others 3; its nonlinear powers take 3. For sufficiency, target 3 is least; target 4 is witnessed by condition 2; all-base target 0 is witnessed by a contained word containing a0 letter, and all-base target 1 by any contained word. A target 2 has its first 2 at the simple last letter and is covered by condition 3. This exhausts assignments for arbitrary variable sets and polynomials.

W1:\(xxy=xy\) and W2:\(xyx=xyy\) give normal form F(w) if the last variable is simple, and F(w) followed by its last letter otherwise. Internal separated repeats are removed by \(xYxZ=xY^{2}Z=xYZ\) for nonempty Y, Z; adjacent internal repeats use W1. A final repeat xYx becomes \(xY^{2}\), and internal reduction leaves the last first-occurrence letter doubled. The word rule proves soundness. Literals remain distinct from squares.

\par\smallskip\noindent\textit{2. Twenty-four laws.}\quad 
{\small\par\smallskip
\CaseLine{\(W1 xxy=xy\qquad W2 xyx=xyy\)}
\CaseLine{\(Q  xx\leq x\qquad C  xy\leq x+y\)}
\CaseLine{\(U  xyz\leq xy+z\qquad I  xyz\leq x+yz\)}
\CaseLine{\(B  xyz\leq xz+y\)}
\CaseLine{\(R1 xyzt\leq xzyt+xyu\qquad R2 xyzz\leq xzyy+xyu\)}
\CaseLine{\(R3 xyz\leq yxz+xu\qquad R4 xyy\leq yxx+xu\)}
\CaseLine{\(J1 xyzt\leq xzyt+y\qquad J2 xyzz\leq xzyy+y\)}
\CaseLine{\(J3 xyz\leq yxz+x\qquad J4 xyy\leq yxx+x\)}
\CaseLine{\(T1 xyz\leq xyzz+xz\qquad T2 xyz\leq xyzz+yz\)}
\CaseLine{\(P1 xy\leq xx+xyy\qquad P2 xyz\leq xyzz+xy\)}
\CaseLine{\(L0 xy\leq xyy+y\qquad L1 xy\leq xyy+x\)}
\CaseLine{\(L2 xyz\leq xyzz+x\qquad L3 xyz\leq xyzz+y\)}
\CaseLine{\(L4 xyzt\leq xyztt+y\)}
\smallskip}

Every absorption satisfies the exact criterion; the two equalities follow from the value rule. All words and actual substitutions are nonempty. Some displayed laws are redundant; no minimality is claimed.

\par\smallskip\noindent\textit{3. Saturate the target content and sort a saturated word.}\quad 
Call a word available when it has been derived below P. Original words are available. Choose u in \(U_{C}\) and obtain \(u^{2}\leq u\) by Q. Maintain a saturated available word g whose content lies in C. For a missing variable choose v in \(U_{C}\) containing it; C gives \(gv\leq g+v\), and Q gives \((gv)^{2}\leq gv\). Normalize. The content increases; condition 1 terminates with full C. Its normal form is G followed by the last letter of G.

We can sort G to ANY target first-occurrence list satisfying condition 2. Let p be its correct prefix and let x be next. Write \(G=p\) A x B. If A is empty there is nothing to move. Otherwise choose an original witness \(q=r\) x V with c(r) in c(p). If q is literal, it is x; apply J1, J2, J3, J4 according as p/B are both nonempty, only p, only B, or neither. In J1/J3 the suffix block is B followed by its last letter. In J2 use \(X=p\), \(Y=x\), \(Z=A\); in J4 use \(X=x\), \(Y=A\). Word normalization gives the desired saturated list p x A B. No absent block is substituted by an empty word.

If q is nonlinear, let R be the entire raw suffix of g after p. It contains A, x and the doubled last letter. Set \(T=V\) R, nonempty. When p is nonempty, B derives p q \(R\leq p\) \(R+q=g+q\); when p is absent, I derives q \(g\leq q+g\). The derived bridge normalizes to p x T. Use R1, R2, R3, R4 in the same four cases, with \(u=T\). These are the same left moves as the J laws, with the bridge replacing the literal witness. Each move extends the correct prefix. Iteration yields the saturated companion of w, denoted \(w^{sat}\). If w is saturated, this completes the derivation.

\par\smallskip\noindent\textit{4. Recover a simple last letter.}\quad 
Suppose \(w=O\) s with O nonempty, distinct letters, and s absent from O. We already have \(w^{sat}=O\) s s. Condition 3 has two kinds of witness.

First suppose an original unmarked witness \(v=V\) s exists, c(V) in c(O). Normalize v, so V is a distinct list and does not contain s. If V is empty, L0 with \(x=O\), \(y=s\) derives w from \(w^{sat}\)+s. Otherwise form the list \(G=V\) followed by all variables of O missing from V, in their target order. The list G s satisfies condition 2: its initial V s-prefix witnesses come from the original v; later missing variables retain their target witnesses since all earlier target variables are now present. Thus Section 3 constructs (G s)\ensuremath{^{\mathrm{sat}}}. If no variable was missing, v already equals G s. Otherwise write \(G=V\) A, A nonempty; T1 gives V A \(s\leq V\) A s s+V s. Hence G s is available, with simple final s.

Sort only the prefix G to O, keeping final s fixed. Every moved x is different from s, so the suffix after x is nonempty. Literal witnesses use J1/J3; nonlinear witnesses use the same bridge construction B/I and then R1/R3. These laws do not alter the final variable or introduce another occurrence of it. Bridge words may use outside variables, but the output keeps exactly the original content and simple tail s. The final word is w. This covers unmarked witnesses in any order.

Otherwise condition 3 supplies an original contained-prefix word r, with c(r) inside c(O). Let F be its first-occurrence list and A the list of O letters missing from F, in target order. The list F A s satisfies condition 2 for the same reason as above:the original r witnesses its initial F letters, and the remaining letters have all their original target predecessors already present. Section 3 constructs (F A s)\ensuremath{^{\mathrm{sat}}}. Q gives \(r^{2}\leq r\); its normal form is F followed by its last letter, which is also the normal form of F F. Put \(Y=A\) s, which is nonempty even when A is absent. P1 with \(x=F\), \(y=Y\) yields

{\small\par\smallskip
\CaseLine{\(F A s \leq  F F + F (A s)(A s) \leq  P.\)}
\smallskip}

The second right word normalizes to (F A s)\ensuremath{^{\mathrm{sat}}}. Thus the simple-tail word F A s is available, including when r is a literal. Finally sort its prefix to O with the protected-tail construction from the preceding paragraph. No unsupported insertion into r is made; only its square and a separately constructed saturated full-content word are used.

For completeness one may instead use L1--L4 in the literal-contained case, according to the literal's four positions in O. These redundant laws are valid but are not needed for the preceding unified recovery.

\par\smallskip\noindent\textit{5. Completeness.}\quad 
Every valid absorption has been derived:literal targets are original summands; saturated targets follow from Section 3; simple-tail targets follow from either of the two exhaustive Section 4 branches. The argument has no variable, word-length, or polynomial-size bound. Distribute terms into nonempty polynomials and absorb every summand in both directions to derive every valid identity. The 24 displayed laws are therefore a finite basis. Reversal handles the opposite; further catalogue transfers require independent complete HSP inclusions in both directions.

\subsection{Ten identities for 4901: fixed-head insertion with original marked witnesses}\label{case:141}
Signature (+, multiplication), without constants, relative to the additively idempotent semiring axioms. The completeness argument is unrestricted in variables, lengths and right polynomial size.

\par\smallskip\noindent\textit{1. Word values and exact criterion.}\quad 
Addition has rows 00044,\allowbreak{}01144,\allowbreak{}01234,\allowbreak{}44334,\allowbreak{}44444, with order \(2<1<0<4\) and \(2<3<4\); both 1+3 and 0+3 are 4. Multiplication has rows 00234,\allowbreak{}11234,\allowbreak{}22222,\allowbreak{}44444,\allowbreak{}44444. Scan for the first factor outside \{0,\allowbreak{}1\}. A first 2 gives 2, a first 4 gives 4, and a first 3 gives 3 when final and 4 otherwise. If no such factor occurs, the product retains its HEAD value 0 or 1. This includes literals and follows by table induction.

For target w and original right words P, let E(x, S) mean a right word contains x with first-x prefix content contained in S. It is marked if a position follows its first x. Let Z(S) mean an original word is contained in S. Then \(w\leq P\) holds exactly when:

(H) An original right word has the same head as w and content contained in c(w). This witness may be literal.

(E) Each first target occurrence x, with preceding content S, has an E(x, S) witness.

(M) At each marked target occurrence, there is a marked E-witness, or an E-witness together with Z(S).

For a literal target, these reduce to the existence of a pure-variable original word. For necessity of H, assign the head 0, the other target variables 1, and outsiders 2. The target is 0 and every right word is 1 or 2. For a failed E, assign \(x=4\), \(S=1\), \(others=2\). For a failed M with E present, assign \(x=3\), \(S=1\), \(others=2\). The first separates target 4 from right values 0/\allowbreak{}1/\allowbreak{}2. In the second, no right word gives 4, no Z(S) means none gives 0 or 1, and only unmarked witnesses can give 3; the target 4 is again separated.

For sufficiency, target 2 is least. An all 0/\allowbreak{}1 target is supplied by the contained same-head witness, with the same value. A first 4 is supplied by E. A marked first 3 gives 4, supplied by a marked witness or by an unmarked witness valued 3 and a Z(S) word valued 0/\allowbreak{}1. A simple final 3 is supplied by any E-witness valued 3 or 4. These exhaust the assignments.

\par\smallskip\noindent\textit{2. Ten identities and normal form.}\quad 
{\small\par\smallskip
\CaseLine{\(W1 xxy=xy\qquad W2 xyx=xyy\)}
\CaseLine{\(L  x\leq xx\qquad H  xx\leq x+xy\)}
\CaseLine{\(A  xyz\leq xz+yt\qquad B  xyzt\leq xyt+yzu\)}
\CaseLine{\(C  abxz\leq abz+x+b\qquad D  abcxz\leq abcz+cx+b\)}
\CaseLine{\(E  xyz\leq xy+yz\qquad F  xy\leq xx+y\)}
\smallskip}

The value formula proves W1, W2, L, H; the exact criterion verifies the remaining absorptions. In C, D the original b is a prefix-contained guard, while x or cx supplies the unmarked prefix witness. These laws can independently be checked on all assignments in the finite table.

As for the two word laws, normalize to the first-occurrence list when the physical final variable is simple; otherwise append one extra copy of the last listed variable. Internal repeats use \(xYxZ=xY^{2}Z=xYZ\) and W1; terminal separated repeats use W2 before internal reductions. All substituted blocks are nonempty. A saturated word means a normal word with doubled final variable. Literals remain separate from squares.

\par\smallskip\noindent\textit{3. Start from the original same-head witness.}\quad 
Let the target be nonlinear with head h. Select H's original witness v. If v is nonlinear, it is the initial available word g. If v is literal h, M at the first target position supplies an original marked witness hT with T nonempty (Z(empty) is impossible). H in the displayed basis gives \(hh\leq h+hT\), so take \(g=hh\).

Maintain an available normal word \(g=pR\) below P, with the correct nonempty target prefix p, nonempty residual R, and content contained in the target content. Initially \(p=h\). Process each remaining marked first occurrence x, using only original E/M witnesses.

If a marked witness is \(q=r\) x V, V nonempty and c(r) contained in c(p), then for r empty A gives p x \(R\leq p\) R+x V. For r nonempty B with \(X=p\), \(Y=r\), \(Z=x\), \(T=R\), \(U=V\) gives

{\small\par\smallskip
\CaseLine{\(p r x R \leq  p r R + r x V.\)}
\smallskip}

The first right word normalizes to g and the left to p x R, since all letters of r have already occurred in p. Thus x is inserted in the required position, whether or not it previously occurred in R.

Otherwise choose an unmarked original \(q=r\) x and an original guard b contained in c(p), as guaranteed by M. For r empty C gives p b x \(R\leq p\) b R+x+b. For r nonempty D gives p b r x \(R\leq p\) b r R+r x+b. Prefix deletion removes b and r from the first right word and the corresponding left prefix, deriving p x R. The guard b can be literal or nonlinear; no derived prefix square is required. All actual substitutions are nonempty in each case.

Normalize the resulting word. Its first-occurrence list begins p x, its content stays in the target, and its residual remains nonempty. To see the last assertion, if R introduces no further variable then its physical final variable already occurs in p x, so normalization retains a final duplicate. Otherwise an unprocessed variable remains. Thus the induction applies through all marked positions.

\par\smallskip\noindent\textit{4. Finish the saturated or simple-tail target.}\quad 
If the target is saturated, process its full first-occurrence list. The maintained content is contained in that list and contains it; the nonempty residual repeats an existing variable. Hence its normal form is exactly the target.

Otherwise the target is p s with simple final s. Process p. If the available g already contains s, its normal form is p s or p s s. In the latter case use L in left context p to derive p \(s\leq p\) s s.

If g omits s, it normalizes to \(p^{2}\). Choose an original E-witness for s. If a marked witness is available, perform the marked insertion once more, giving p s s, and contract by L in context p. Otherwise the witness is r s with c(r) contained in c(p). For r nonempty E gives p r \(s\leq p\) r+r s; the first right word normalizes to \(p^{2}\) and the left to p s. If r is empty, the original literal s and F give p \(s\leq p^{2}+s\). All final cases derive the target.

A literal target follows from its original pure witness by W1, W2, L. Every valid absorption is therefore derivable. Distributivity expands arbitrary terms into nonempty polynomials, and the summand absorptions in both directions prove every identity. The ten laws form a finite basis for 4901; reversal gives a basis for its opposite. Catalogue transfers still require complete HSP in both directions. Finite regressions and proof replay supplement the unrestricted proof.

\subsection{Eight identities for 4933: a same-head seed and guarded ordered insertion}\label{case:142}
Signature (+, multiplication), without constants; all bases are relative to the additively idempotent semiring axioms.

\par\smallskip\noindent\textit{1. Exact word values and absorption.}\quad 
Catalogue operation rows in order 0,\allowbreak{}1,\allowbreak{}2,\allowbreak{}3,\allowbreak{}4 are:

{\small\par\smallskip
\noindent addition        multiplication\par
\CaseLine{\ttfamily 00004\qquad 00234}
\CaseLine{\ttfamily 01114\qquad 11234}
\CaseLine{\ttfamily 01224\qquad 33333}
\CaseLine{\ttfamily 01234\qquad 33333}
\CaseLine{\ttfamily 44444\qquad 44444}
\smallskip}

Addition is the \(chain 3<2<1<0<4\). A literal retains its value. For a nonlinear word, find the first factor outside\{0,\allowbreak{}1\}:4 gives 4,\allowbreak{}3 gives 3, and 2 gives 2 only when it is final, otherwise3. If there is no such factor, the word has its FIRST factor's value. This follows immediately by induction on multiplication. Unlike4872, the all-base value does not detect an arbitrary0 later in the word; the head is essential.

For target w, let \(C=c(w)\). For variable x and set S, an original prefix witness is an original right word with first-x prefix content in S. It is unmarked if that first occurrence is final. An original contained-prefix word has its whole content in S.

A literal target is absorbed exactly when it is an original literal. A nonlinear target is absorbed exactly when all three conditions hold:

1. There is an original word u with head h(w) and c(u) contained in C.

2. Every first occurrence x in w has an original prefix witness whose first-x prefix lies inside the preceding target content S.

3. When the last variable s of w is simple, its prefix witness can be unmarked, or there is an original word with content in C minus\{s\}.

For necessity of 1, assign the target head 0, other target variables 1, and outsiders 3. Target 0 can only be matched by a contained same-head word. For 2 assign preceding content 1, the specified variable 4, and all others 3. For 3 assign preceding content 1, \(s=2\), others 3. Each failure has all right words strictly below the target. A missing literal x is separated by \(x=2\) and all other variables 3; all nonlinear powers of x are 3. These arguments apply to arbitrary variable sets and polynomials.

For sufficiency, target 3 is least. Target 4 is supplied by 2 at the first 4-valued target variable. If the target is all-base0 or 1, the same-head contained word from 1 has exactly its value. Target 2 has its unique 2-valued variable at the simple last position; an unmarked witness has value 2 and a contained-prefix witness has base value 0 or 1. These exhaust all assignments. Thus the criterion is exact.

\par\smallskip\noindent\textit{2. Normal form and the eight laws.}\quad 
{\small\par\smallskip
\CaseLine{\(W1 xxy=xy\qquad W2 xyx=xyy\)}
\CaseLine{\(Q  xx\leq x\qquad B  xyz\leq xz+y\)}
\CaseLine{\(I  xyz\leq xz+xyt\qquad O  xyz\leq xyzz+y\)}
\CaseLine{\(T  xyz\leq xyzz+yz\qquad L  xy\leq xyy+y\)}
\smallskip}

The value rule proves W1, W2; the exact criterion proves every displayed absorption (and can be checked on each original right pair). W1, W2 normalize a word to its first-occurrence list F if its last variable is simple, and otherwise to F followed once more by F's last variable. For internal separated repeats use \(xYxZ=xY^{2}Z=xYZ\) with nonempty Y, Z; adjacent repeats use W1. A terminal repeat xYx becomes \(xY^{2}\) by W2, then internal repeats disappear. Literals remain distinct from squares. Call the second kind of normal form saturated.

\par\smallskip\noindent\textit{3. Fix the head, then insert the target first occurrences.}\quad 
Let \(w\leq P\) satisfy the criterion and normalize all original words. Choose the original same-head contained seed u from 1. Q gives \(u^{2}\leq u\). Its normal form g is saturated, has the required head, and uses only C. Maintain a correct, nonempty first-occurrence prefix p of g, initially the target head. Let x be the next target first-occurrence variable. If it is already next, extend p and continue. Otherwise write \(g=pR\). The residual R is NONEMPTY, even when p contains all the present content, because g is saturated and retains its final repeated occurrence.

Choose an original witness \(q=r\) x V from 2 with c(r) contained in c(p). If q is literal, it is x. B gives p x \(R\leq pR+x\leq P\). Otherwise B gives the bridge p q \(R\leq pR+q\leq P\). Normalization of p q R is the same as that of p x V R, because r introduces no new variable and the final variable of R remains repeated. Put \(T=V\) R, which is nonempty even when V is absent. I with \(x=p\), \(y=x\), \(z=R\), \(t=T\) gives

{\small\par\smallskip
\CaseLine{\(p x R \leq  pR+p x T \leq  P.\)}
\smallskip}

The new normal form has prefix p x, has no outside variable, and remains saturated. If x was absent it is added; if it was present later, its later occurrence is removed by W1, W2. Its new head remains the original head. Thus each step increases the correct prefix and eventually gives exactly F(w) followed by its last variable. No bound on variables, word lengths or polynomial size appears. No absent block is substituted by an empty word; all substitutions into B and I are nonempty.

\par\smallskip\noindent\textit{4. Recover a simple final variable.}\quad 
A saturated target is now derived. Otherwise \(w=p\) s with nonempty p, s absent from p, and the constructed companion is p s s. If an original unmarked witness is \(v=r\) s with nonempty r and c(r) contained in c(p), substitute \(x=p\), \(y=r\), \(z=s\) into T. Normalization gives p \(s\leq p\) s s+v. If the unmarked witness is the literal s, use L. If instead an original contained-prefix witness r exists, substitute \(x=p\), \(y=r\), \(z=s\) into O. Normalization gives p \(s\leq p\) s s+r. The contained-prefix witness may be literal or nonlinear; it remains an original summand. All actual substitutions are nonempty.

\par\smallskip\noindent\textit{5. Unrestricted completeness and transfer.}\quad 
Literal targets are original summands. Sections 3 and 4 derive every other valid absorption according to the exact criterion. Expand any term by distributivity into a nonempty polynomial, then absorb each summand in both directions. This derives every valid identity and proves that the eight displayed laws form a finite basis for 4933. Reversal gives the opposite basis.

\subsection{A finite basis for 4942 by endpoint ideals and bounded original witnesses}\label{case:143}
Signature (+, multiplication), without constants. The basis below is an effectively specified finite family, not a small enumerated list. This completes the earlier endpoint-ideal research note.

Addition rows 00444,\allowbreak{}01444,\allowbreak{}44224,\allowbreak{}44234,\allowbreak{}44444 give two \(chains 1<0<4\) and \(3<2<4\). Multiplication rows 00234,\allowbreak{}11234,\allowbreak{}44434,\allowbreak{}44434,\allowbreak{}44434 imply the following formula. Literals keep their value. For a nonlinear word, final 3 gives 3, final 4 gives 4. With final 0/\allowbreak{}1/\allowbreak{}2, any nonfinal value outside \{0,\allowbreak{}1\} gives 4; otherwise final 2 gives 2 and final 0/\allowbreak{}1 gives the head. Induction proves this at all lengths; final 3 can reset an earlier 4.

\par\smallskip\noindent\textit{Exact criterion for arbitrary polynomials.}\quad 
For original P let C be the content union, H the set of heads, T the set of finals, and \(A_{t}\) the union of nonfinal-variable sets of original words ending in t. For any positive-length target w, put \(h=head(w)\), \(t=last(w)\), I=content(w with its final position removed). Then \(w\leq P\) iff

content(w) subset C, h in H, t in T; \(A_{t}\) intersects T, or I subset \(A_{t}\).

The formula includes literal targets: I is then empty. Here is a full proof using the four proper principal additive ideals. If \(P\leq 0\), each right word has value 0/\allowbreak{}1, which forces all its variables into\{0,\allowbreak{}1\}. Content containment forces the target into\{0,\allowbreak{}1\}, hence \(w\leq 0\). If \(P\leq 1\), all right words have value 1, so additionally every original head is 1. The target head belongs to H and consequently \(w=1\). If \(P\leq 3\), all right finals are 3, and t in T forces \(w=3\). Finally suppose \(P\leq 2\). All finals belong to\{2,\allowbreak{}3\}. If \(t=3\) the target is 3. If \(t=2\), each original t-final word forces every variable in \(A_{t}\) into\{0,\allowbreak{}1\}. Thus \(A_{t}\) cannot meet T. The criterion forces I into \(A_{t}\), so the target has value 2. \(P=4\) is automatic. These cases prove sufficiency.

For necessity, if a target variable is outside C, assign it 4, all other variables 0: \(P=0\), \(w=4\). If h is outside H, assign \(h=0\), all other variables 1: \(P=1\), \(w=0\). If t is outside T, assign all T variables 3 and all others 0: \(P=3\) and w is 0/\allowbreak{}4. Finally assume the first three conditions, that \(A_{t}\) and T are disjoint, and choose q in I outside \(A_{t}\). Assign \(A_{t}\) variables 0, \(t=2\), every other T variable 3, and \(q=4\) if q is not in T; all remaining variables 0. This is consistent because t is not in \(A_{t}\). All original t-final words are 2 and all other originals are 3, so \(P=2\). The target has a nonfinal value outside\{0,\allowbreak{}1\}, hence 4. If \(q=t\) or q belongs to T, the same argument uses its assigned 2/\allowbreak{}3. This proves necessity for all variables, word lengths, and polynomial sizes.

\par\smallskip\noindent\textit{Three fixed identities and normal forms.}\quad 
{\small\par\smallskip
\CaseLine{\(xyzt = xzyt\)}
\CaseLine{\(xxy = xy\)}
\CaseLine{\(xyz \leq  xz+yz\)}
\smallskip}

The word-value formula proves the two equalities. The exact criterion proves the third: its upper words end in z and have prefix union\{x, y\}. The equalities normalize each nonlinear word to its original head, then the increasing set of interior variables other than the head, then its original final. Permute adjacent interior letters with the first law, delete interior duplicates with the second, and move interior head copies next to the head before deleting. All suffixes required by \(xxy=xy\) are nonempty. Keep both endpoints, including equal endpoints; an internal copy of a distinct final variable also remains. Literals are unchanged. These are literal equational derivations.

\par\smallskip\noindent\textit{Elementary recoveries use at most four original words.}\quad 
Suppose a nonlinear \(w\leq P\) is valid. For each x in I, derive \(hxt\leq P\). If \(A_{t}\) and T are disjoint, choose an original word headed by h, an original t-final word containing h nonfinally, and an original t-final word containing x nonfinally. These at most three words already absorb hxt by the exact criterion; all its content is covered, its head and final have witnesses, and its prefix\{h, x\} lies in \(A_{t}\).

If \(A_{t}\) meets T, choose q in their intersection. Choose original words headed by h, containing x, ending in t with a nonfinal q, and ending in q. These at most four words absorb hxt: they cover h, x, t, provide its head and final, and retain the original crossing q. Witnesses may coincide; idempotent addition removes duplicate summands.

For a literal target h, choose one original head-h word and one original final-h word. Their sum absorbs h by the criterion, whose prefix is empty. No claim that h must itself be an original is made.

\par\smallskip\noindent\textit{Reversible bounded compression.}\quad 
In any of these local recoveries there are at most five words total. Normalize them and reserve all heads and finals, at most 10 variables. Group every other variable by its nonzero incidence pattern among the five words, at most 31 groups. Replace each reserved variable by its own fresh variable and each group by one fresh variable, then normalize. There are at most 41 variables and 43 positions per word. The compressed recovery is valid because identification of variables is a substitution into the original valid recovery.

For the inverse substitution, restore each reserved variable as its singleton and each group variable as the nonempty increasing product of the group's members. Every group occurs only in permutable, idempotent interiors. The expansions recover precisely the original head, final, and interior set. Literal uses of the two equalities therefore normalize back to all original normal words. This supplies a genuine inverse derivation; endpoint multiplicities are retained.

\par\smallskip\noindent\textit{Combining elementary recoveries and completeness.}\quad 
Start with \(hht\leq P\) for \(x=h\); it normalizes to ht. If the current word is Ut and hxt is an already derived elementary word, use \(xyz\leq xz+yz\) with \(x=U\), \(y=hx\), \(z=t\) to derive \(Uhxt\leq Ut+hxt\leq P\). Normalize to remove the new interior head copy. Repeating over I minus\{h\} yields the target normal word, including any required internal occurrence of t. All substituted words are nonempty. Literals use their two original witness recovery directly.

Take the three fixed identities and all valid absorptions on a fixed 41-variable alphabet, word length at most 43, and at most four upper words. This family is finite and its validity test is the exact criterion. Compression and inverse substitution derive every elementary and literal recovery; the combining construction derives every valid word absorption. Expansion by distributivity and summandwise absorption in both directions derive every valid identity, without bounds on the original variables, lengths, or summands. Hence 4942 and its opposite are finitely based. Catalogue transfers require separate complete HSP.

\subsection{Nine identities for 5001 from nonlinear prefix witnesses and recursive guards}\label{case:144}
Constant-free ai-semiring signature. The completeness argument below is unrestricted in variables, word lengths and polynomial sizes.

\par\smallskip\noindent\textit{1. Values and the exact absorption criterion.}\quad 
Catalogue operations in the order 0,\allowbreak{}1,\allowbreak{}2,\allowbreak{}3,\allowbreak{}4:

{\small\par\smallskip
\noindent addition       multiplication\par
\CaseLine{\ttfamily 00044\qquad 00244}
\CaseLine{\ttfamily 01134\qquad 00244}
\CaseLine{\ttfamily 01234\qquad 22222}
\CaseLine{\ttfamily 43334\qquad 44444}
\CaseLine{\ttfamily 44444\qquad 44444}
\smallskip}

The additive bottom is 2, the top 4; the \(chains 2<1<0<4\) and \(2<1<3<4\) have \(0+3=4\). A literal retains its value. For a nonlinear word, the first factor outside\{0,\allowbreak{}1\} is decisive:2 gives 2, and 3 or 4 gives 4. If all factors are 0 or 1, the word gives 0. Direct induction on the table proves this rule for every positive length. Only nonlinear words always have values in\{0,\allowbreak{}2,\allowbreak{}4\}; a literal 3 must not be treated as a reset product 4.

For original P let N be its original NONLINEAR summands. For nonlinear target w and \(C=c(w)\), the exact conditions are:

(N) There is a u in N with c(u) contained in C. (E) For every first target occurrence x, let S be its preceding target content. Either an original q in N contains x with its first-x prefix inside S, or the literal x is an original summand AND there is an original u in N with c(u) contained in S.

A literal target x is absorbed exactly when an original word has content \{x\}. To prove necessity assign \(x=0\) and every other variable 2. For sufficiency use \(x\leq x^{2}\) and the pure-variable original word.

If N fails, assign target variables 1 and outsiders 2. The nonlinear target is 0, but each right value is 1 or 2. If E fails at x, assign \(S=1\), \(x=3\), and all others 2. The target is 4. A nonlinear right word can give 4 only as the first alternative of E. A right literal can give 3 only if it is x. The join reaches 4 from 0+3 only when a nonlinear right word is contained in S. Otherwise all right values lie under 0 or under 3. Thus each failed condition has an explicit separating assignment.

For sufficiency, target 2 is bottom. Target 0 has only base-valued target variables, so the nonlinear contained witness from N gives 0. Target 4 has a first nonbase variable x valued 3 or 4. A nonlinear prefix witness gives 4. Alternatively its literal x and the nonlinear contained-prefix witness have join 4. This exhausts assignments for arbitrary P and w.

In \(\A{5002}\) one has \(1+3=4\). Its absorption criterion is treated separately.

\par\smallskip\noindent\textit{2. Nine laws and normalization.}\quad 
{\small\par\smallskip
\CaseLine{\(W1 xxy=xy\qquad W2 xyx=xy\)}
\CaseLine{\(W3 xyy=xy\qquad W4 xxx=xx\)}
\CaseLine{\(L  x\leq xx\qquad H  xy\leq yy+xt\)}
\CaseLine{\(B  xyz\leq xz+yy\qquad I  xyz\leq xz+xyt\)}
\CaseLine{\(J  xyz\leq xz+y+xx\)}
\smallskip}

W1--W4 follow from the word-value rule. L is checked at each element. For H, B, I, J, the exact criterion holds for their displayed original right sides; hence they are valid. No missing block denotes an empty substitution. The normal form is the first-occurrence list when there are at least two distinct variables, a square for nonlinear singleton content, and the literal for length one. Delete separated repeats by W2; adjacent repeats use W1 or W3; W4 reduces singleton powers. Every nonlinear v satisfies \(v^{2}=v\) by these word laws. This does not apply to literals. Nonlinear normal forms still retain their first letter.

\par\smallskip\noindent\textit{3. Strong induction on distinct target variables.}\quad 
Normalize original words and targets. A nonlinear singleton target \(x^{2}\) is already an original normalized word by N. For a target with \(d\geq 2\) variables, assume all valid nonlinear targets with fewer than d variables can be derived from the nine laws below the SAME original polynomial P. Choose an original contained nonlinear u by N. At the first target letter h, E cannot use an empty-prefix contained word, so it gives an original nonlinear \(q=hV\) with V nonempty. H with \(x=h\), \(y=u\), \(t=V\) gives

{\small\par\smallskip
\CaseLine{\(h u \leq  u^{2}+q = u+q \leq  P.\)}
\smallskip}

Normalize to g. Its head is h and its content lies in C. Maintain a nonempty correct first-occurrence prefix p of g. Let x be the next target variable. If x is already next, extend p. Otherwise write \(g=pR\). If p already lists the entire content of g, restore a final repeat by W3 backwards (or use the existing square for singleton content), so R is nonempty and pR still normalizes to g. This is a word-law equality, not an empty-word substitution or an unproved new absorption.

If E gives an original nonlinear \(q=r\) x V with c(r) inside c(p), B gives

{\small\par\smallskip
\CaseLine{\(p q R \leq  pR+q^{2} = g+q \leq  P.\)}
\smallskip}

Its normal form is that of p x V R. Put \(T=V\) R, nonempty. I at \(x=p\), \(y=x\), \(z=R\), \(t=T\) gives p x \(R\leq pR+p\) x \(T\leq P\). Normalizing adds x if absent or moves it to the correct position if present later. Outside variables may occur in the bridge but not in the resulting word.

Otherwise E gives the original literal x and an original NONLINEAR word u contained in c(p). The nonlinear target \(p^{2}\) is valid below P: u supplies N, and every first variable of p inherits exactly its earlier E condition from the original target. Since c(p) has fewer than d variables, the strong induction hypothesis derives \(p^{2}\leq P\). J at \(x=p\), \(y=x\), \(z=R\) gives p x \(R\leq pR+x+p^{2}\leq P\).

In both branches the correct prefix strictly increases. At completion the normal form is the target w. Every recursive call uses a strictly smaller number of distinct variables. The outer induction and finite prefix loop terminate with no fixed bound on the original parameters. All witnesses remain original summands; derived words are never promoted to original nonlinear guards for the criterion.

\par\smallskip\noindent\textit{4. Completeness.}\quad 
Literal targets are original literals or follow from L and an original pure-variable nonlinear word. Nonlinear targets are handled by Section 3. Distribute arbitrary terms into nonempty polynomials and absorb each summand in both directions. This proves that the nine displayed laws form a complete finite basis for 5001. Reversal gives the opposite basis. Any catalogue transfer requires independently complete HSP inclusions both ways. Finite tests support this argument and do not replace it.

\subsection{Thirteen identities for 5002 with nonlinear seeds and literal prefix guards}\label{case:145}
Constant-free ai-semiring signature. All constructions are unrestricted in variables, word length and polynomial size.

\par\smallskip\noindent\textit{1. Word values and exact criterion.}\quad 
Catalogue addition rows are 00044,\allowbreak{}01144,\allowbreak{}01234,\allowbreak{}44334,\allowbreak{}44444; multiplication rows are 00244,\allowbreak{}00244,\allowbreak{}22222,\allowbreak{}44444,\allowbreak{}44444. Thus 2 is additive bottom,4 top, \(2<1<0<4\) and \(2<3<4\), with \(1+3=0+3=4\). Multiplication is the same as 5001: a literal retains its value; a nonlinear word gives 2 if its first factor outside\{0,\allowbreak{}1\} is 2, gives 4 if that first factor is 3 or 4, and gives 0 when all factors are 0 or 1. Induction on the table proves the rule.

Let N be the original NONLINEAR right words of P. A literal target x is absorbed iff an original word has content\{x\}. A nonlinear target w is absorbed iff:

(N) some original u in N has c(u) contained in c(w); (E) each first target occurrence x, with preceding target content S, has an original NONLINEAR prefix witness q, or BOTH the original literal x and an original word u (possibly literal) contained in S.

For failure of N, assign \(c(w)=1\) and outsiders 2:the target 0 exceeds all right values 1 or 2. For failure of E, assign \(S=1\), \(x=3\), others 2. Without a nonlinear prefix witness, the only possible3 is the literal x; it reaches join 4 only with a prefix-contained word, whose value is 0 or 1. A missing pure-x witness for a literal target is separated by \(x=0\), others 2. For sufficiency, target 2 is bottom; target 0 is supplied by N; target 4 is supplied by its first 3/\allowbreak{}4-valued variable's nonlinear prefix witness, or its literal plus an original prefix-contained word. Literal targets use \(x\leq x^{2}\). This covers every assignment, not just bounded tests.

The nonlinear contained seed in N remains mandatory. A preliminary candidate allowing a literal to supply N was false: \(x^{2}\leq x+xy+yx\) fails at \(x=1\), \(y=2\) (left 0, right 1). That candidate was discarded before proof certification and never used to classify a type.

\par\smallskip\noindent\textit{2. Thirteen valid laws.}\quad 
{\small\par\smallskip
\CaseLine{\(W1 xxy=xy\qquad W2 xyx=xy\)}
\CaseLine{\(W3 xyy=xy\qquad W4 xxx=xx\)}
\CaseLine{\(L  x\leq xx\qquad H  xy\leq yy+xt\)}
\CaseLine{\(B  xyz\leq xz+yy\qquad I  xyz\leq xz+xyt\)}
\CaseLine{\(J  xyz\leq xz+y+xx\)}
\CaseLine{\(G1 xyztu\leq xyzu+t+y\)}
\CaseLine{\(G2 xyzt\leq xyt+z+x\)}
\CaseLine{\(G3 xyzt\leq xyt+z+y\)}
\CaseLine{\(G4 xyz\leq xz+y+x\)}
\smallskip}

The value formula proves W1--W4, L; the exact criterion proves the remaining absorptions. W1--W4 reduce words to first-occurrence lists, except that singleton nonlinear words remain squares. Every nonlinear u has \(u^{2}=u\); this is not asserted for literals. The explicit reductions and nonempty substitutions are the same word-law reductions given in part 2 of Section~\ref{case:144}.

\par\smallskip\noindent\textit{3. Strong induction and the four new cases.}\quad 
Induct on the number d of distinct variables in the nonlinear target. For \(d=1\), N supplies the original square. For \(d\geq 2\), choose an original contained nonlinear u. E at the first variable h gives an original nonlinear \(q=hV\) with V nonempty. H derives h \(u\leq u^{2}+q\leq P\). Normalize and maintain a correct nonempty prefix p. The next letter x may be added or moved by a NONLINEAR original witness \(q=r\) x V exactly as follows:choose nonempty R with pR normalizing to the available word (restoring a final repeat if needed); B derives p q \(R\leq pR+q^{2}\); I with its fourth block V R derives p x \(R\leq pR+p\) x V R. Outside variables in the bridge disappear from the resulting word.

If x has no nonlinear witness, E supplies its original literal and an original contained-prefix guard. If a NONLINEAR contained-prefix guard exists, then \(p^{2}\) satisfies N, E with fewer than d distinct variables. Induction derives \(p^{2}\leq P\), and J derives p x \(R\leq pR+x+p^{2}\leq P\).

Otherwise choose an original literal b from c(p). Write \(p=A\) b B. Use G1 when A and B are nonempty, with blocks A, b, B, x, R. Use G2 when A is empty and B nonempty, with blocks b, B, x, R. Use G3 when A is nonempty and B empty, with blocks A, b, x, R. Use G4 when \(p=b\), with blocks b, x, R. In all four cases the actual absorption is

{\small\par\smallskip
\CaseLine{\(p x R \leq  pR+x+b \leq  P.\)}
\smallskip}

All substitutions are nonempty; the four identities are separate laws, not empty-word instances of one law. Each iteration increases the correct first-occurrence prefix. Recursive nonlinear-guard calls strictly decrease the distinct-variable count. Thus the algorithm terminates for arbitrary finite variable sets and polynomial sizes.

\par\smallskip\noindent\textit{4. Completeness.}\quad 
This derives every nonlinear valid absorption by N, E. Literal targets are original literals or follow from L and an original pure-variable nonlinear summand. Distribute all terms and absorb summands in both directions to derive any valid identity. Hence the thirteen laws give a finite basis for 5002. Opposite tables use reversed laws. Further catalogue transfers require independent complete HSP in both directions. The difference from 5001 is precisely the possible LITERAL prefix guard; it does not eliminate the nonlinear contained seed requirement.

\subsection{Last-occurrence witnesses: twelve-identity bases for a five-element semiring and its opposite}\label{case:146}
The source manuscript for this case credits Aifa Wang and Lili Wang, School of Mathematical Sciences, Chongqing University of Technology.

\par\smallskip\noindent\textit{Abstract.}\quad 
We give a twelve-identity basis for an explicit five-element additively idempotent semiring. A word is determined by the order of its last occurrences and by whether its final variable occurs earlier. Absorption has three independent requirements: internal-content coverage, a contained-content witness, and a last-occurrence suffix witness for every variable. Four identities supply a repeated final variable, and three identities splice suffix witnesses into a target word. The construction proves completeness for arbitrarily many variables and summands. Word reversal gives a basis for the multiplicative opposite. No assertion about finite bases for all proper subvarieties is made.

\par\smallskip\noindent\textit{1. The algebra and its basis.}\quad 
Work in the constant-free signature \((+,\cdot)\) relative to the axioms for additively idempotent semirings. Write \(u\leq v\) for \(u+v\approx v\). Let \(A=\{d, e, r, a,\omega\}\), with addition and multiplication given below. These names denote elements, not constants in terms.

\noindent \(+\)\enspace\textbar\enspace \(d\)\enspace\textbar\enspace \(e\)\enspace\textbar\enspace \(r\)\enspace\textbar\enspace \(a\)\enspace\textbar\enspace \(\omega\)\par
\noindent \(d\)\enspace\textbar\enspace \(d\)\enspace\textbar\enspace \(d\)\enspace\textbar\enspace \(d\)\enspace\textbar\enspace \(\omega\)\enspace\textbar\enspace \(\omega\)\par
\noindent \(e\)\enspace\textbar\enspace \(d\)\enspace\textbar\enspace \(e\)\enspace\textbar\enspace \(d\)\enspace\textbar\enspace \(a\)\enspace\textbar\enspace \(\omega\)\par
\noindent \(r\)\enspace\textbar\enspace \(d\)\enspace\textbar\enspace \(d\)\enspace\textbar\enspace \(r\)\enspace\textbar\enspace \(\omega\)\enspace\textbar\enspace \(\omega\)\par
\noindent \(a\)\enspace\textbar\enspace \(\omega\)\enspace\textbar\enspace \(a\)\enspace\textbar\enspace \(\omega\)\enspace\textbar\enspace \(a\)\enspace\textbar\enspace \(\omega\)\par
\noindent \(\omega\)\enspace\textbar\enspace \(\omega\)\enspace\textbar\enspace \(\omega\)\enspace\textbar\enspace \(\omega\)\enspace\textbar\enspace \(\omega\)\enspace\textbar\enspace \(\omega\)\par
\noindent \(\cdot\)\enspace\textbar\enspace \(d\)\enspace\textbar\enspace \(e\)\enspace\textbar\enspace \(r\)\enspace\textbar\enspace \(a\)\enspace\textbar\enspace \(\omega\)\par
\noindent \(d\)\enspace\textbar\enspace \(d\)\enspace\textbar\enspace \(d\)\enspace\textbar\enspace \(r\)\enspace\textbar\enspace \(\omega\)\enspace\textbar\enspace \(\omega\)\par
\noindent \(e\)\enspace\textbar\enspace \(d\)\enspace\textbar\enspace \(e\)\enspace\textbar\enspace \(r\)\enspace\textbar\enspace \(a\)\enspace\textbar\enspace \(\omega\)\par
\noindent \(r\)\enspace\textbar\enspace \(r\)\enspace\textbar\enspace \(r\)\enspace\textbar\enspace \(r\)\enspace\textbar\enspace \(\omega\)\enspace\textbar\enspace \(\omega\)\par
\noindent \(a\)\enspace\textbar\enspace \(\omega\)\enspace\textbar\enspace \(\omega\)\enspace\textbar\enspace \(r\)\enspace\textbar\enspace \(\omega\)\enspace\textbar\enspace \(\omega\)\par
\noindent \(\omega\)\enspace\textbar\enspace \(\omega\)\enspace\textbar\enspace \(\omega\)\enspace\textbar\enspace \(r\)\enspace\textbar\enspace \(\omega\)\enspace\textbar\enspace \(\omega\)\par
The additive order has \(e<d\), \(r<d\), \(e<a\), and \(d, a<\omega\); there are no other nonreflexive comparisons except those forced by transitivity. Addition is its join. The core \(\{d, e, r,\omega\}\) is a band: \(\{d, e\}\) is a semilattice with identity \(e\), and \(\{r,\omega\}\) is a right-zero ideal on which \(d, e\) act identically. The row of \(a\) equals that of \(\omega\), while \(ua\) equals \(a\) precisely when \(u=e\), and equals \(\omega\) otherwise.

These observations verify associativity. For a triple ending in \(a\), both associations give \(a\) exactly when its first two factors are \(e\), and give \(\omega\) otherwise. For a triple not ending in \(a\), any \(a\) in its first two positions can be replaced by \(\omega\) without changing either association; the band core then applies. The rows and columns preserve the displayed order and the joins of the three incomparable pairs \((e, r)\), \((d, a)\), \((r, a)\). Thus both distributive laws hold.

Theorem 1. The following twelve identities form a basis for \(A\):

\CaseDisplay{\begin{array}{rll}
\mathrm{(N_1)}&x^2y\approx xy,\\
\mathrm{(N_2)}&xyx\approx yx^2,\\
\mathrm{(E)}&x\leq x^2,\\
\mathrm{(C_1)}&xy\leq x^2+zy,\\
\mathrm{(C_0)}&xy\leq x^2+y,\\
\mathrm{(D_{11})}&xy^2\leq xy+zyt,\\
\mathrm{(D_{10})}&xy^2\leq xy+yt,\\
\mathrm{(D_{01})}&y^2\leq y+zyt,\\
\mathrm{(D_{00})}&y^2\leq y+yt,\\
\mathrm{(T_{11})}&xyt\leq xt+zyt,\\
\mathrm{(T_{10})}&xyt\leq xt+yt,\\
\mathrm{(T_{01})}&yt\leq t+zyt.
\end{array}\quad\textup{(\ensuremath{1})}}

Reversing every word in (1) gives a basis for \(A^{\mathrm{op}}\). The number twelve is not claimed to be minimal.

\par\smallskip\noindent\textit{2. Word values and normal forms.}\quad 
For a nonempty word \(w\), let \(c(w)\) denote its content, and let \(i(w)\) be the set of variables before its final position. For \(x\in c(w)\), let \(s_w(x)\) be the set of variables strictly after its last occurrence. Let \(\lambda(w)\) be the injective word obtained by deleting every occurrence except the last one of each variable, retaining the order of the surviving positions. Define

\CaseDisplay{\nu(w)=
\begin{cases}
\lambda(w),&w\text{ has a nonrepeated final variable},\\
\lambda(w)t,&t\text{ is the final variable and occurs earlier}.
\end{cases}\quad\textup{(\ensuremath{2})}}

Thus normal words are injective words, optionally followed by a second copy of their final variable. The words \(x\) and \(x^2\) are both allowed.

Lemma 2. Identities \((N_1),(N_2)\) derive \(w\approx\nu(w)\) for every nonempty word. They also derive

\CaseDisplay{\nu(ww)\approx\nu(w)t,\qquad t=\text{the final variable of }w,}

where an already repeated terminal variable is reduced to two copies.

Proof. With a nonempty suffix, \((N_1),(N_2)\) give \(xyxz\approx yx^2z\approx yxz\). They therefore delete an earlier occurrence of any variable whose later occurrence is nonfinal; adjacent occurrences use \((N_1)\). An arbitrary preceding context is preserved. After these deletions, the only possible earlier repetition is that of the final variable \(t\). If it occurs as \(tRt\) with \(R\) nonempty, \((N_2)\) gives \(tRt\approx Rt^2\); adjacent copies need no change. Any extra preceding copies have already been deleted. This leaves exactly (2), in the order of last occurrences. The assertion about \(ww\) follows by applying the same normalization. \(\square\)

The following evaluation rule follows immediately from the multiplication table. Call \(r, a,\omega\) the special elements. If no letter is assigned a special element, the value is \(e\) when every letter is assigned \(e\), and is \(d\) otherwise. If special elements do occur, take their last position. A last special element \(r\) gives \(r\); a last special element \(\omega\) gives \(\omega\). A last special element \(a\) gives \(\omega\) when that position is nonfinal; when it is final, the value is \(a\) exactly when every preceding letter is assigned \(e\), and is \(\omega\) otherwise. This also covers linear words.

\par\smallskip\noindent\textit{3. The exact absorption criterion.}\quad 
For a nonempty polynomial \(P=\sum_{u\in\mathcal U}u\), put \(J(P)=\bigcup_{u\in\mathcal U}i(u)\).

Lemma 3. For any nonempty word \(w\), \(A\models w\leq P\) if and only if all three conditions hold:

1. \(i(w)\subseteq J(P)\).

2. Some \(u\in\mathcal U\) has \(c(u)\subseteq c(w)\).

3. For each \(x\in c(w)\), some \(u\in\mathcal U\) contains \(x\) and satisfies \(s_u(x)\subseteq s_w(x)\).

Proof of necessity. If condition 1 fails at \(x\), assign \(x=a\) and all other variables \(e\). The target is \(\omega\), whereas each right word is \(e\) or \(a\), because \(x\) occurs there only terminally.

If condition 2 fails, assign the variables in \(c(w)\) to \(e\) and all others to \(r\). The target is \(e\) and every right word is \(r\); \(e\not\leq r\).

If condition 3 fails at \(x\), assign \(x=\omega\), all variables in \(s_w(x)\) to \(e\), and every other variable to \(r\). The target is \(\omega\). Each right word containing \(x\) has a variable assigned \(r\) after its last \(x\), and hence evaluates to \(r\). Words without \(x\) evaluate to \(e\) or \(r\). Their join is at most \(d\), a contradiction.

Proof of semantic sufficiency. Fix an assignment and write \(p=P\). If \(p=\omega\), there is nothing to prove. If \(p=e\), all right words have value \(e\), so every variable appearing in them is assigned \(e\). Condition 3 covers \(c(w)\), hence \(w=e\).

If \(p=a\), every right word is \(e\) or \(a\). All its internal variables are assigned \(e\), and all its variables belong to \(\{e, a\}\). Conditions 1 and 3 imply the same assertions for \(w\), so \(w\) is \(e\) or \(a\).

If \(p=d\), each right word is \(e, r\), or \(d\). A target value outside \(\{e, r, d\}\) has a last special occurrence assigned \(a\) or \(\omega\). Condition 3 supplies a word containing that variable with only \(e, d\) after its last occurrence. Its value is \(a\) or \(\omega\), contradicting \(p=d\).

Finally, if \(p=r\), every right word is \(r\). The preceding argument rules out a target whose last special occurrence is \(a\) or \(\omega\). If the target has no special occurrence, every variable in it is assigned \(e\) or \(d\). Condition 2 then supplies a right word also evaluating to \(e\) or \(d\), a contradiction. Thus its last special element is \(r\), and its value is \(r\). This exhausts the additive values of \(P\). \(\square\)

The evaluation rule proves \((N_1),(N_2)\). Each of the ten inequalities in (1) satisfies the three conditions of Lemma 3: \((E)\) uses \(x^2\); the \((C)\) identities use \(x^2\) as the contained-content and \(x\)-suffix witness, and their other word for \(y\); the \((D)\) identities use their first word for content and suffixes and their second word for the internal occurrence of \(y\); the \((T)\) identities use the first word for content and their two words for the indicated suffixes and internal positions. Therefore all twelve identities are sound.

\par\smallskip\noindent\textit{4. Three derivation operations.}\quad 
All words below are nonempty unless explicitly described as an absent substring. An absent substring is never substituted into an identity: the relevant boundary identity in (1) is used instead. All comparisons in this section are derived from (1).

Doubling. Suppose \(q\leq P\), and the final variable \(y\) of \(q\) belongs to \(J(P)\). Choose an original right word \(v=zyt\) in which the displayed \(y\) is internal. Thus \(t\) is nonempty and \(z\) may be absent. Write \(q=py\), with \(p\) possibly absent. The appropriate one of the four \((D)\) identities gives \(py^2\leq py+zyt\leq P\), with the corresponding absent prefixes omitted. By Lemma 2 this is \(q^2\leq P\). If \(q\) already has a repeated terminal variable, \(q^2\approx q\) and no new step is needed.

Appending a witnessed suffix. Suppose \(q^2\leq P\) and an original or already derived word \(zy\leq P\) has a nonempty suffix \(y\). Then \((C_1)\) gives \(qy\leq q^2+zy\leq P\); when \(z\) is absent use \((C_0)\). Here \(q\) and \(y\) denote nonempty word substitutions, not variables restricted to length one.

Insertion before a unique terminal. Suppose \(pt\leq P\) and \(zyt\leq P\), where \(t\) is a single terminal variable, \(y\) is nonempty, and \(p, z\) may be absent. Then

\CaseDisplay{pyt\leq pt+zyt\leq P.\quad\textup{(\ensuremath{3})}}

Use \((T_{11})\) when both \(p, z\) are present, \((T_{10})\) when only \(p\) is present, and \((T_{01})\) when only \(z\) is present. When both are absent, the target \(yt\) is already the second premise. This accounts for all boundary cases.

\par\smallskip\noindent\textit{5. Arbitrary-variable constructive completeness.}\quad 
Lemma 4. Every absorption satisfying Lemma 3 is derivable from (1).

Proof. Normalize all words. Let \(\lambda(w)=x_1\cdots x_n\), and write \(C=c(w)\) and \(J=J(P)\). Choose a contained-content seed \(u\in\mathcal U\) from condition 2, and for each \(x_j\) choose a suffix witness \(v_j\) from condition 3. In a normal witness, the suffix beginning at the last \(x_j\) contains \(x_j\) followed only by variables among \(x_{j+1},\ldots, x_n\); its terminal variable may have a repeated final copy.

First suppose \(C\subseteq J\). Begin with \(q=u\leq P\). For \(j=1,\ldots, n\), double the current \(q\) and append the suffix of \(v_j\) beginning at its last \(x_j\). These operations give another word below \(P\); normalize it. Its content stays within \(C\). The already processed variables \(x_1,\ldots, x_{j-1}\) do not occur in the appended suffix, whereas that suffix introduces a new last \(x_j\) after all their last occurrences. Inductively the processed variables have last occurrences in their prescribed order. An unprocessed variable may still precede them; its later processing moves its last occurrence to the right. At \(j=n\), therefore, \(\lambda(q)=x_1\cdots x_n\). Double once more. The result is precisely \(\nu(w^2)\leq P\). If \(w\) has a repeated terminal variable this is \(w\); otherwise \((E)\) gives \(w\leq w^2\leq P\).

Now suppose \(C\not\subseteq J\). Condition 1 forces the unique exception to be the final variable \(s=x_n\), and \(s\notin J\). In every right word containing \(s\), it occurs once, in the final position. The seed either already ends in \(s\) or omits it. If it omits \(s\), its terminal variable lies in \(C\setminus\{s\}\subseteq J\); double it and append the one-letter suffix \(s\) of \(v_n\). We now have \(q=ps\leq P\), containing no internal \(s\) and no variable outside \(C\).

For each \(j<n\), take \(v_j\). If it contains \(s\), it already ends in \(s\). If it omits \(s\), its final variable lies either in \(s_w(x_j)\) or is \(x_j\) itself; thus it belongs to \(C\setminus\{s\}\subseteq J\). Double \(v_j\) and append the suffix \(s\) of the current \(q\), obtaining a normalized word \(v_j s\leq P\). Although the part of \(v_j\) before \(x_j\) may contain variables outside \(C\), its suffix from the last \(x_j\) contains only \(x_j\) and later target variables. The suffix from that \(x_j\) up to the final \(s\) has the form \(ys\), with \(y\) nonempty.

Apply (3) to the current \(ps\) and this witness, obtaining \(pys\leq P\). Normalize. Exactly as in the first case, no processed variable occurs in \(y\), and the new last \(x_j\) precedes only unprocessed variables. After \(j=n-1\) the result is \(x_1\cdots x_{n-1}s=w\). For \(n=1\) the loop is absent and the preceding preparation already gives \(w=s\). This proves every required absorption. \(\square\)

Proof of Theorem 1. Soundness was established after Lemma 3. Distributivity writes any term as a nonempty finite sum of words. If \(P\approx Q\) holds in \(A\), each summand of \(P\) is absorbed by \(Q\). Lemmas 3 and 4 derive every such absorption, hence derive \(P\leq Q\). Interchanging the sides derives \(Q\leq P\), and thus the identity. No bound on variables or summands was used. Reversing words converts this argument into a derivation for \(A^{\mathrm{op}}\). \(\square\)

\par\smallskip\noindent\textit{Acknowledgments.}\quad 
The source manuscript for this case acknowledges support from the National Natural Science Foundation of China (12371024), the Science and Technology Research Program of Chongqing Municipal Education Commission (KJZD-K2024011102) and the Chongqing Natural Science Foundation Innovation and Development Joint Fund (Municipal Education Commission) (CSTB2025NSCQ-LZX0067).

\subsection{Twelve identities from nonlinear prefix witnesses}\label{case:147}
Source 5157. The signature is (+, multiplication) without constants, relative to additively idempotent semiring axioms. A one-letter word has length one; a longer word has length at least two.

\par\smallskip\noindent\textit{1. Algebra and word values.}\quad 
Tables in order 0,\allowbreak{}1,\allowbreak{}2,\allowbreak{}3,\allowbreak{}4:

{\small\par\smallskip
\noindent addition       multiplication\par
\CaseLine{\ttfamily 0 0 0 4 4\qquad 0 0 2 4 4}
\CaseLine{\ttfamily 0 1 1 3 4\qquad 0 1 2 4 4}
\CaseLine{\ttfamily 0 1 2 3 4\qquad 2 2 2 2 2}
\CaseLine{\ttfamily 4 3 3 3 4\qquad 4 4 4 4 4}
\CaseLine{\ttfamily 4 4 4 4 4\qquad 4 4 4 4 4}
\smallskip}

The additive order has \(2<1<0<4\) and \(1<3<4\), with 0 and 3 incomparable. For a word of length at least two, its first factor in\{2,\allowbreak{}3,\allowbreak{}4\}, if present, determines the value:2 gives 2,\allowbreak{}3 or 4 gives 4. Without such a factor all factors are 0/\allowbreak{}1, and the value is 0 if any factor is 0, otherwise1. A one-letter word retains its value. Induction from the table proves this formula. It proves associativity on three factors; the 0/\allowbreak{}1 product is their minimum in the usual numeric order, and the first reset factor determines mixed products. Direct inspection shows rows and columns preserve all displayed joins, proving both distributive laws.

\par\smallskip\noindent\textit{2. The twelve identities and normalization.}\quad 
The following list is a complete relative basis:

\CaseDisplay{\begin{aligned}
x^2y&\approx xy,&xyx&\approx xy,&xy^2&\approx xy,\quad\textup{(\ensuremath{W})}\\
x&\le x^2,\quad\textup{(\ensuremath{L})}\\
xyz&\le xy+z^2,\quad\textup{(\ensuremath{U})}\\
xyz&\le x^2+yz,\quad\textup{(\ensuremath{I})}\\
xyz&\le xz+y^2,\quad\textup{(\ensuremath{B})}\\
xy&\le x+y+xyz,\quad\textup{(\ensuremath{T})}\\
xyzt&\le xzyt+xyu,\quad\textup{(\ensuremath{R1})}\\
xyz&\le xzy+xyu,\quad\textup{(\ensuremath{R2})}\\
xyz&\le yxz+xu,\quad\textup{(\ensuremath{R3})}\\
xy&\le yx+xu.\quad\textup{(\ensuremath{R4})}
\end{aligned}}

(W) deletes all repeated letters except that a longer word of singleton content becomes xx rather than x. A word with at least two distinct variables becomes the list of first occurrences. Separated repetitions use \(xYx=xY\) with nonempty Y; adjacent ones use either \(x^{2}y=xy\) or \(xy^{2}=xy\) with a nonempty neighboring block. Singleton powers use \(x^{3}=x^{2}\), an instance of \(x^{2}y=xy\). Thus every longer word q satisfies \(q^{2}\) approximately q by(W). This claim does not include one-letter words. The table directly proves L.

\par\smallskip\noindent\textit{3. Exact absorption criterion and soundness.}\quad 
For a target w of content C let \(U_{P}\) be the original right words of P contained in C. A one-letter y is absorbed exactly when an original right word has content\{y\}. For a longer target w, absorption holds exactly when:

\textbullet{} C: the union of contents of \(U_{P}\) equals C;

\textbullet{} O: every target variable x has an original right word of length at least two containing x whose first-x prefix content is contained in the first-x prefix content of w.

The O witness may contain variables outside C after its first x. No separate contained-longer-word condition is required.

Necessity. Without a pure-y witness assign \(y=0\) and outsiders 2; every right word is 2, whereas the one-letter target is 0. If C fails at x assign \(x=0\), other target variables 1, and outsiders 2. The target is 0; contained right words avoid x and are 1, while noncontained right words are 2. Thus \(P\leq 1\). If O fails at x assign earlier target variables 1, \(x=3\), and all remaining variables 2. The target is 4. Every longer right word with x has a2 before its first x and is 2; other longer words are 1 or 2. One-letter summands give only 1,\allowbreak{}2,\allowbreak{}3, so \(P\leq 3\), which is strictly below 4.

Sufficiency. Target 2 is least. Target 4 is caused by a first reset variable x valued 3/\allowbreak{}4. Its O witness is a longer word and has no reset before that x, hence is 4. For target 0 all variables lie in\{0,\allowbreak{}1\} and some x is 0; C supplies a contained right word containing x, of value 0. For target 1 all target variables are 1, and any contained right word supplied by C is 1. One-letter targets are absorbed by their pure-variable witnesses, since L and(W) give \(y\leq y^{k}\) for every positive k. This proves the criterion for arbitrarily many variables and original summands.

(W) follows from the word formula. L is already checked directly. Every remaining displayed absorption satisfies C and O, so is sound. They are also independently verified on all assignments.

\par\smallskip\noindent\textit{4. Building a longer word of full target content.}\quad 
Normalize original words. A one-letter target follows from L and its pure-variable witness. Let the target now have length at least two. If \(U_{P}\) has a longer word, choose it as \(g\leq P\). Otherwise C implies every variable in C is a literal right summand. Let h be the target head and choose its O witness \(q=hV\), with V nonempty. An instance of T gives

\CaseDisplay{hh\le h+hhV\approx h+hV\le P.}

Here h is an original literal and q is an original longer word; set \(g=hh\). In both cases g is longer and contained in C.

While its content is not C, take the first missing variable x in the target's first-occurrence order. Choose a contained original word v containing x, using C. If v is longer, split \(g=ab\) nonemptily. U gives \(gv=abv\leq ab+v^{2}\) approximately \(g+v\leq P\). Normalize to obtain a longer guard with increased content still contained in C.

If v is literal x, choose an original O witness \(q=r\) x V at its first x. All variables of r lie earlier than x in the target and are already in g, by the choice of the first missing variable. U again gives \(gq\leq g+q^{2}\) approximately \(g+q\leq P\). Its normal form is the first-occurrence list of g, then x, then a tail Q; Q may contain outside variables. If Q is empty this is the desired normal form of gx. If Q is nonempty use T with \(X=g\), \(Y=x\), \(Z=Q\):

\CaseDisplay{gx\le g+x+gxQ\le P.}

The last right word normalizes to gq, already obtained. The new guard gx has content c(g) union\{x\}, contained in C. Thus the process terminates with a longer g of content exactly C. The case of singleton C is already complete.

\par\smallskip\noindent\textit{5. Sorting with original longer prefix witnesses.}\quad 
For larger C, g is the list of its distinct first occurrences. Assume its correct target prefix is p and the next target variable is x. Write \(g=p\) A x B. If A is empty the prefix extends already. Otherwise A is nonempty. Choose the longer O witness \(q=r\) x V, with c(r) contained in c(p). Put \(R=A\) x B and keep \(U=V\) R as a raw nonempty word; it contains x even if subsequent normalization removes some letters.

If p is nonempty, B with \(X=p\), \(Y=q\), \(Z=R\) gives \(pqR\leq pR+q^{2}\) approximately \(g+q\leq P\). Its normal form is that of pxU, since the letters of r already occur in p. If p is empty, r is empty; split \(g=ab\) with both parts nonempty. I gives \(qg\leq q^{2}+g\) approximately \(q+g\leq P\), whose normal form is xU.

The following instances then move x across A:

\textbullet{} p, B nonempty: R1 yields \(pxAB\leq pAxB+pxU\).

\textbullet{} p nonempty, B empty: R2 yields \(pxA\leq pAx+pxU\).

\textbullet{} p empty, B nonempty: R3 yields \(xAB\leq AxB+xU\).

\textbullet{} p, B empty: R4 yields \(xA\leq Ax+xU\).

The right words are g and the just-constructed bridge, modulo(W). Each new word has content C and a longer correct target prefix. Iteration reaches w. All substituted words are nonempty, including \(U=V\) A x B in every lifting step.

\par\smallskip\noindent\textit{6. Completeness and duality.}\quad 
Sections 4--5 derive every valid word absorption with no bound on original variables, word length or right summands. Distributivity expresses arbitrary terms as sums of words. Derive absorption of each summand on either side by the other side to obtain every valid identity. This proves the twelve-law basis assertion. Reversing words transports all identities and derivations to the opposite. Finite tests supplement this unrestricted proof, and catalogue transfer requires both complete HSP directions.

\subsection{Thirteen identities for 5158: full contained cover and guarded prefix witnesses}\label{case:148}
Signature (+, multiplication) without constants, relative to the additively idempotent semiring axioms. No variable, length or summand bound.

\par\smallskip\noindent\textit{Word values and the exact criterion.}\quad 
Addition rows:00044,\allowbreak{}01144,\allowbreak{}01234,\allowbreak{}44334,\allowbreak{}44444. Multiplication rows:00244,\allowbreak{}01244,\allowbreak{}22222,\allowbreak{}44444,\allowbreak{}44444. Addition has \(2<1<0<4\) and \(2<3<4\), with \(0+3=1+3=4\). A literal retains its value. In a nonlinear word the first factor outside\{0,\allowbreak{}1\} determines the result: first 2 gives 2; first 3 or 4 gives 4. If all factors lie in\{0,\allowbreak{}1\}, the result is 0 when any factor is 0, and 1 otherwise. Induction proves this for all lengths. In particular head alone does not determine the all 0/\allowbreak{}1 branch.

For original right words P and target w, put \(C=c(w)\), and let U consist of original words contained in C. A literal x is absorbed iff some original word is a pure power of x. A nonlinear target is absorbed iff: (C) the contents of U cover ALL of C; (N) for each first target occurrence x, with preceding content S, there is an original NONLINEAR prefix witness q containing x with first-x prefix contained in S; or P contains literal x and an original word b contained in S. Nonlinear witnesses may have outside variables after their first x.

Necessity of C: for an uncovered x assign \(x=0\), other target variables 1, outsiders 2. The target is 0, but contained right words avoid x and give 1, while outsiders force 2. For a literal without a pure witness set \(x=0\), others 2. For failed N set \(x=3\), \(S=1\), all others 2. The target is 4. A right nonlinear word giving 4 would be a forbidden nonlinear witness. Value 3 can come only from literal x, and value 1 only from a contained-prefix word. The failed disjunction prevents both 3 and 1 appearing, so the right sum is below 4. These are explicit separating assignments.

Sufficiency:2 is least. An all 0/\allowbreak{}1 target valued 0 has a zero-valued target variable; C supplies a contained word containing it, also valued 0. For value 1, any contained word suffices. A first target reset 3/\allowbreak{}4 is supplied by N's nonlinear witness, valued 4; or by literal 3/\allowbreak{}4 plus a prefix-contained word valued 0/\allowbreak{}1, whose sum is 4. A pure-variable witness absorbs its literal because 3 squares to 4 and all other elements are idempotent. This proves the criterion without bounds and independently of any finite testing.

\par\smallskip\noindent\textit{The thirteen sound laws.}\quad 
{\small\par\smallskip
\CaseLine{\(W1 xxy=xy\qquad W2 xyx=xy\qquad W3 xyy=xy\)}
\CaseLine{\(L  x\leq xx\qquad U  xyz\leq xy+zz\qquad I  xyz\leq xx+yz\)}
\CaseLine{\(B  xyz\leq xz+yy\qquad T  xy\leq x+y+xyz\)}
\CaseLine{\(R1 xyzt\leq xzyt+xyu\qquad R2 xyz\leq xzy+xyu\)}
\CaseLine{\(R3 xyz\leq yxz+xu\qquad R4 xy\leq yx+xu\)}
\CaseLine{\(G  abxz\leq abz+x+b\)}
\smallskip}

The word-value formula proves W1, W2, W3, L. Each remaining displayed absorption satisfies C and N and is therefore sound. In G, abz supplies the nonliteral prefix witnesses and the original b guards the original literal x. No shorter-table resemblance is used to infer soundness.

W1, W2, W3 reduce every word with at least two distinct letters to its first-occurrence list. A nonlinear singleton becomes xx; literal x stays x. Delete separated repetition xYx by W2; an adjacent repeat with nonempty suffix uses W1, and a final adjacent repeat with nonempty prefix uses W3. Pure powers reduce using W1 at length 3. Thus every nonlinear q satisfies q \(squared=q\), but this is never applied to a literal. Normalization preserves content, first-occurrence prefix sets, and literal/nonlinear status.

\par\smallskip\noindent\textit{Build a nonlinear available word covering the target content.}\quad 
Normalize original words and target, retaining their provenance. Literals follow from pure witnesses and L. For nonlinear w, if U has a nonlinear word choose it as g. Otherwise C means every target variable is an original literal. The head h has an N witness hV with V nonempty, since the prefix is empty and no guard can have empty content. T gives

{\small\par\smallskip
\CaseLine{\(hh \leq  h+hhV = h+hV \leq  P.\)}
\smallskip}

Start with \(g=hh\). In both cases g is nonlinear and contained in C. Take the first missing target variable x. Choose v in U containing x. If v is nonlinear, write \(g=ab\) nonemptily. U gives

{\small\par\smallskip
\CaseLine{\(gv=abv \leq  ab+v squared = g+v \leq  P.\)}
\smallskip}

If v is literal x, and N has a nonlinear witness \(q=r\) x V, all letters of r already belong to g because x was the first missing target variable. The same U gives \(gq\leq g+q\). Normalize this to gxQ, where Q may be empty or contain outsiders. If Q nonempty, T gives \(gx\leq g+x+gxQ\), removing outsiders; if empty we already have gx. Content increases and remains in C.

If no nonlinear witness exists, N gives original literal x and original guard b whose content consists of earlier target variables, all in g. G with \(a=g\), \(b=b\), \(x=x\), \(z=g\) gives

{\small\par\smallskip
\CaseLine{\(g b x g \leq  g b g+x+b.\)}
\smallskip}

The right first word normalizes to g and the left to gx; the presence of nonlinear g and internal repetition laws justify both reductions. Thus the literal-guard case also increases content exactly by x. Iterate until a nonlinear g of full content C is available. Singleton content is finished.

\par\smallskip\noindent\textit{Sort while preserving original guards.}\quad 
For larger C both target and g are lists of distinct letters. Suppose the already correct target prefix is p and the next target letter is x, with \(g=p\) A x B. If A is empty continue. Otherwise A is nonempty and let \(R=A\) x B.

If N uses an original nonlinear \(q=r\) x V, its r is contained in p. Put \(U0=V\) R as a RAW nonempty word; R contains x, so U0 is nonempty even when V is empty. For p nonempty, law B gives \(pqR\leq pR+q\) \(squared=g+q\). For p empty, r is empty and law I after a nonempty split \(g=ab\) gives \(qg\leq q\) squared+g. The resulting bridge normalizes to p x U0 (or x U0 if p empty). Apply:

\textbullet{} R1 when p and B nonempty: \(pxAB\leq pAxB+pxU0\);

\textbullet{} R2 when p nonempty, B empty: \(pxA\leq pAx+pxU0\);

\textbullet{} R3 when p empty, B nonempty: \(xAB\leq AxB+xU0\);

\textbullet{} R4 when p and B empty: \(xA\leq Ax+xU0\).

Each right side is supplied by g and the derived bridge. The new g has full content C and a longer correct prefix.

If only the literal-and-guard alternative is available, then p is necessarily nonempty. Keep the original literal x and original guard b contained in p. G gives

{\small\par\smallskip
\CaseLine{\(p b x R \leq  p b R+x+b.\)}
\smallskip}

The right first word normalizes to g; the left normalizes to p x A B, using all letters of b already in p and deleting the second x in R. This sorts x immediately without changing either original guard summand. Iterating reaches the target. Every substitution in every case is nonempty.

The construction derives every valid absorption on any finite variable set, with arbitrarily long words and arbitrarily many right summands. Distributivity and summandwise absorption in both directions then derive every identity. The thirteen laws are a finite relative basis for 5158. Word reversal gives the opposite basis by Lemma~\ref{lem:transfer}.

\subsection{Nine identities from a contained head and nonlinear prefix witnesses}\label{case:149}
Source 5219, in the signature (+, multiplication) without constants. All basis statements are relative to the additively idempotent semiring axioms. A longer word means length at least two, whether or not it repeats a variable.

\par\smallskip\noindent\textit{1. Algebra and unrestricted word semantics.}\quad 
The catalogue tables, in the order 0,\allowbreak{}1,\allowbreak{}2,\allowbreak{}3,\allowbreak{}4, are

{\small\par\smallskip
\noindent addition       multiplication\par
\CaseLine{\ttfamily 0 0 0 4 4\qquad 0 0 2 4 4}
\CaseLine{\ttfamily 0 1 1 3 4\qquad 1 1 2 4 4}
\CaseLine{\ttfamily 0 1 2 3 4\qquad 2 2 2 2 2}
\CaseLine{\ttfamily 4 3 3 3 4\qquad 4 4 4 4 4}
\CaseLine{\ttfamily 4 4 4 4 4\qquad 4 4 4 4 4}
\smallskip}

The additive order is \(2<1<0<4\) and \(1<3<4\), with 0 and 3 incomparable. A one-letter word retains its assigned value. For every longer word, the first factor valued in \{2,\allowbreak{}3,\allowbreak{}4\}, if any, determines the result: 2 gives 2, and 3 or 4 gives 4. If all factors lie in \{0,\allowbreak{}1\}, the result is the value of the first letter. Induction using the multiplication table proves this formula for all lengths. In particular, the product on \{0,\allowbreak{}1\} is left-zero, not the 0/\allowbreak{}1 product of source 5157.

The formula on three factors verifies associativity. Each row and each column of multiplication preserves the joins in the displayed additive order, verifying distributivity; addition itself is that join operation.

\par\smallskip\noindent\textit{2. Proposed finite list and word normalization.}\quad 
The following nine identities form a complete relative basis:

\CaseDisplay{\begin{aligned}
x^2y&\approx xy,&xyx&\approx xy,&xy^2&\approx xy,\quad\textup{(\ensuremath{W})}\\
x&\le x^2,\quad\textup{(\ensuremath{L})}\\
xyz&\le xy+z^2,\quad\textup{(\ensuremath{U})}\\
xyz&\le xz+y^2,\quad\textup{(\ensuremath{B})}\\
xy&\le x+xyz,\quad\textup{(\ensuremath{T})}\\
xyzt&\le xzyt+xyu,\quad\textup{(\ensuremath{R1})}\\
xyz&\le xzy+xyu.\quad\textup{(\ensuremath{R2})}
\end{aligned}}

As usual, \(p\leq q\) abbreviates p+q approximately q. By (W), a word with at least two distinct variables reduces to its first-occurrence list. A longer word of singleton content reduces to xx; a one-letter word remains x. To delete a separated repetition use xYx= xY with nonempty Y; an adjacent repetition is deleted using \(x^{2}y=xy\) or \(xy^{2}=xy\) with a nonempty neighboring block. The remaining singleton powers reduce by \(x^{3}=x^{2}\). Thus every longer word q satisfies \(q^{2}\) approximately q, but that assertion is not made for one-letter words. These reductions preserve heads and first-occurrence prefix sets.

\par\smallskip\noindent\textit{3. Exact absorption criterion.}\quad 
Let P be a nonempty finite collection of original right words, and let w be a target word with content C. A one-letter target y is absorbed by P if and only if some original right word has content \{y\}.

For a longer target, \(w\leq \sum (P)\) holds if and only if both conditions hold:

\textbullet{} H: some original right word q has the same head as w and c(q) subset C.

\textbullet{} O: for every variable x in w, some original right word q of length at least two contains x, and the content strictly before its first x is contained in the content strictly before the first x of w.

The H witness may be one-letter. The O witnesses can have variables outside C after their first x. There is no separate content-cover requirement.

For necessity of H, assign the target head 0, other target variables 1, and outsiders 2. A contained right word with a different head is 1; every noncontained right word is 2. Without H their sum is at most 1, whereas the target is 0. For necessity of O at x, assign earlier target variables 1, x the value 3, and all remaining variables 2. The target is 4. In every longer right word containing x, failure of O places a 2 before its first x, so its value is 2; other longer words are 1 or 2. The one-letter words give at most 3, so their total is at most \(3<4\). Without a pure-y witness for a one-letter target, \(y=0\) and \(outsiders=2\) give target 0 and right sum 2.

Conversely, value 2 is least. If the longer target is 4, its first reset variable x has value 3 or 4, and its O witness has only 0/\allowbreak{}1-valued variables before x. That longer witness is therefore 4, even when x is its last letter. If the target is 0 or 1, all its variables have values in \{0,\allowbreak{}1\}; its H witness has the same head and is contained in C, so has the same value. Finally, a pure-y word absorbs y because squaring fixes 0,\allowbreak{}1,\allowbreak{}2,\allowbreak{}4 and sends 3 to 4, and all higher powers equal the square. This proves the criterion for any number of variables and summands.

The word formula proves (W), and the squaring observation proves (L). Each of (U),(B),(T),(R1),(R2) satisfies H and O directly on its displayed variables, so the entire list is sound. Notice in particular that (T) does not need a separate literal y on its right side.

\par\smallskip\noindent\textit{4. Constructing a longer word with the target head and full content.}\quad 
Normalize w and all original right words by (W). One-letter targets follow from their pure-variable witnesses and (L). Now w is longer.

If H has a longer witness, take it as \(g\leq \sum (P)\). Otherwise the target head h itself is an original right summand. The O witness for h has the form \(q=hV\), where V is nonempty. An instance of (T) gives

\CaseDisplay{hh\le h+hhV\approx h+hV\le\sum P.}

Set \(g=hh\). In either case g is longer, its head is h, and its content is contained in C.

If a variable is missing from g, take the first missing variable x in the target's first-occurrence order. Its O witness is \(q=r\) x V, with all letters of r earlier than x in w. They already occur in g. Split \(g=ab\) with a and b nonempty. By (U) and (W),

\CaseDisplay{gq=abq\le ab+q^2\approx g+q\le\sum P.}

After normalization this is g followed by x and a possibly empty tail Q, modulo repeated letters of g. If Q is empty, we have obtained the normal form of gx. If Q is nonempty, (T) with \(X=g\), \(Y=x\), \(Z=Q\) gives

\CaseDisplay{gx\le g+gxQ\le\sum P.}

This removes any outside tail variables, without needing literal x as an original summand. Replace g by the normal form of gx. Its head stays h and its content increases by exactly x, remaining contained in C. Finitely many steps produce a longer word g of content C and head h. If C is a singleton this is already the normalized target.

\par\smallskip\noindent\textit{5. Sorting while preserving the head.}\quad 
For larger C, g is its distinct first-occurrence list and has the correct target head. Suppose its correct target prefix is p, and x is the next target variable. Thus p is nonempty. Write \(g=p\) A x B. If A is empty the prefix already extends; otherwise A is nonempty.

Choose an original longer O witness \(q=r\) x V with c(r) subset c(p). Set \(R=A\) x B and keep \(U=V\) R as a raw, nonempty word. It contains x even if later normalization removes repetitions. By (B),

\CaseDisplay{pqR\le pR+q^2\approx g+q\le\sum P.}

The left side normalizes to pxU since all letters of r already occur in p. If B is nonempty, (R1) with \(X=p\), \(Y=x\), \(Z=A\), \(T=B\) gives

\CaseDisplay{pxAB\le pAxB+pxU\le\sum P.}

If B is empty, (R2) gives \(pxA\leq pAx+pxU\). These substitutions are all nonempty, including \(U=V\) A x B. The resulting word has content C and a strictly longer correct target prefix, with head unchanged. Iteration reaches w. No initial-head interchange is used or needed.

\par\smallskip\noindent\textit{6. Completeness and opposite algebras.}\quad 
The construction derives every valid word absorption with arbitrarily many original variables, word lengths and original right summands. Distributivity writes every term as a finite sum of nonempty words. Deriving each summand's absorption in the opposite side proves every valid identity. Therefore the nine identities are a complete finite relative basis. Reversing all words transports the proof to the opposite algebra. Transferring the catalogue classification requires independently checked complete HSP inclusions in both directions.

Finite table checks and replay of concrete substitutions supplement this unrestricted proof; they do not replace it.

\subsection{Eight identities for 5220: every first occurrence needs a nonlinear prefix witness}\label{case:150}
Signature (+, multiplication), without constants, relative to the additively idempotent semiring axioms. All words and substituted blocks are nonempty. The proof below has no bound on variables, lengths or summands.

\par\smallskip\noindent\textit{Word values and exact semantic criterion.}\quad 
Addition rows are 00044,\allowbreak{}01144,\allowbreak{}01234,\allowbreak{}44334,\allowbreak{}44444. Its order is \(2<1<0<4\) and \(2<3<4\), with \(0+3=1+3=4\). Multiplication rows are 00244,\allowbreak{}11244,\allowbreak{}22222,\allowbreak{}44444,\allowbreak{}44444. A literal retains its assigned value. For a nonlinear word, scan for the first factor outside \{0,\allowbreak{}1\}: first 2 returns 2, first 3 or 4 returns 4; if no such factor occurs, retain the head. This follows by induction directly from the multiplication table. In particular a final 3 in a nonlinear word also gives 4, unlike source 4901.

For the original nonempty polynomial P and target w, put \(C=c(w)\). A literal x is absorbed by P exactly when P has a pure-variable word \(x^{m}\). A nonlinear w is absorbed exactly when both conditions hold:

(H) P has an original word with head h(w) and content contained in C. (N) For each first occurrence x in w, let S be the content strictly before that occurrence. There is an original NONLINEAR word q containing x whose first-x prefix is contained in S; alternatively, P contains the original LITERAL x and an original word b with c(b) contained in S.

Necessity of H: assign \(h(w)=0\), other target variables 1 and outsiders 2. If H fails all right words have value 1 or 2 whereas w has value 0. For a literal target with no pure witness, assign \(x=0\) and all others 2. For failed N at x, assign \(x=3\), \(S=1\) and all other variables 2. The nonlinear target has value 4. No right nonlinear word can give 4, because this would be a nonlinear prefix witness. Right values are 2, possibly 3 from the literal x, and possibly 1 from contained-prefix words. The failed disjunction forbids both 1 and 3 occurring, so their sum is below 4. This separates.

For sufficiency a nonlinear target valued 2 is least; an all 0/\allowbreak{}1 target is supplied by H with the same head value. At the first target 3 or 4, condition N supplies a nonlinear right word valued 4, or a literal 3/\allowbreak{}4 together with a prefix-contained word valued 0/\allowbreak{}1, whose sum is 4. For a literal target, \(x\leq x^{m}\) holds for every \(m\geq 1\) by the displayed word-value rule. Thus these conditions are exact for arbitrary original words and polynomial sizes.

\par\smallskip\noindent\textit{Eight laws and normal form.}\quad 
{\small\par\smallskip
\CaseLine{\(W1 xxy=xy\qquad W2 xyx=xyy\)}
\CaseLine{\(L  x\leq xx\qquad H  xx\leq x+xy\)}
\CaseLine{\(A  xyz\leq xz+yt\qquad B  xyzt\leq xyt+yzu\)}
\CaseLine{\(U  xyzt\leq xyt+yz\qquad G  abxz\leq abz+x+b\)}
\smallskip}

All laws satisfy the criterion (and have independent full-table checks). In U the original yz is a NONLINEAR witness even when the selected z is its simple final variable; no extra guard is needed. G retains the original literal witness x and original prefix guard b, even when b is literal.

W1 and W2 give the first-occurrence normal form: list all distinct letters in order of first appearance, and append a duplicate of the last listed letter exactly when the physical final letter of the original word already appeared earlier. Indeed an internal repetition x Y x Z, Y, Z nonempty, rewrites by W2 to x Y Y Z and by W1 \((X=Y)\) to x Y Z; consecutive internal repetitions use W1 directly. A separated terminal repeat uses W2, then internal deletion, leaving the last listed letter doubled. Pure powers reduce to a square. Literals remain distinct. These operations are valid inside arbitrary contexts, never using an empty substitution. The exact criterion is invariant under these replacements, which preserve the first occurrence prefix sets and whether a word is literal.

\par\smallskip\noindent\textit{Unrestricted construction of every absorption.}\quad 
Normalize original words and target, retaining original witness provenance. Let F be the target first-occurrence list. For a nonlinear target choose H's original contained same-head seed v. If v is nonlinear take \(g=v\). If v is literal h, N at the first target position (S empty) forces an original nonlinear hT, T nonempty. Law H derives \(hh\leq h+hT\); take \(g=hh\). Maintain an available normal word \(g=pR\leq P\) with p the correct nonempty prefix of F, R nonempty, and c(g) contained in C. Initially \(p=h\).

For the next letter x of F use N's original nonlinear witness \(q=r\) x V. Its prefix r uses only letters already in p. If V is nonempty and r empty, A yields p x \(R\leq p\) R+x V. If both are nonempty, B yields

{\small\par\smallskip
\CaseLine{\(p r x R \leq  p r R + r x V.\)}
\smallskip}

If V is empty, q is nonlinear, so r is nonempty; U yields

{\small\par\smallskip
\CaseLine{\(p r x R \leq  p r R + r x.\)}
\smallskip}

In the last two displays the first right word normalizes to g and the left normalizes to p x R, because all letters of r already occur in p and a nonempty suffix R protects the deletion. If no nonlinear witness exists, N supplies original literal x and an original prefix guard b. Law G gives

{\small\par\smallskip
\CaseLine{\(p b x R \leq  p b R + x + b.\)}
\smallskip}

Delete internal copies from b using its content containment in p. This again gives p x R. No prefix is ever appended to the original literal witness or guard; they remain original summands of P.

After normalization the new word begins with p x. Its residual is nonempty: if R contains an unprocessed letter, at least one such letter remains; otherwise the physical final variable already occurs in p x, so the normal form keeps a terminal duplicate. Content stays inside C. Hence the induction continues through EVERY letter of F, including a simple final target letter. After the last step g is exactly F with its last letter doubled. A saturated target is now available. For a simple-tail nonlinear target \(F=p\) t, p is nonempty; L in context p derives p \(t\leq p\) t t. A literal target follows from its original pure-variable witness by W1, W2, L.

Every semantic absorption has thus been derived from eight laws using only nonempty substitutions, legitimate contexts, and original witnesses. Distributivity expands arbitrary terms into nonempty polynomials. Apply these summand absorptions in both directions to any valid identity. This proves a finite basis for 5220 and, by reversing words, its opposite.

\subsection{Four identities for a contained head and ordered tail witnesses}\label{case:151}
Source 5222, signature (+, multiplication) without constants. The basis is relative to the additively idempotent semiring axioms.

\par\smallskip\noindent\textit{1. Tables, values and the finite list.}\quad 
The catalogue operations, in order 0,\allowbreak{}1,\allowbreak{}2,\allowbreak{}3,\allowbreak{}4, are

{\small\par\smallskip
\noindent addition       multiplication\par
\CaseLine{\ttfamily 0 0 0 4 4\qquad 0 0 2 4 4}
\CaseLine{\ttfamily 0 1 1 3 4\qquad 1 1 2 4 4}
\CaseLine{\ttfamily 0 1 2 3 4\qquad 2 2 2 4 4}
\CaseLine{\ttfamily 4 3 3 3 4\qquad 3 3 2 4 4}
\CaseLine{\ttfamily 4 4 4 4 4\qquad 4 4 2 4 4}
\smallskip}

Addition is join in \(2<1<0<4\) and \(1<3<4\), where 0 and 3 are incomparable. Right multiplication by 0 or 1 is the identity map, by 2 is constant 2, and by 3 or 4 is constant 4. Consequently, for \(w=hT\), the last tail factor in \{2,\allowbreak{}3,\allowbreak{}4\}, if present, determines the value: 2 gives 2 and 3 or 4 gives 4. If there is no such tail factor, the value is the head value. This formula includes one-letter words by using an absent tail, not by adjoining an identity element.

Composition of these right maps verifies associativity. The displayed rows and columns preserve the additive joins, verifying distributivity.

Theorem. A complete relative basis consists of the following four identities:

\CaseDisplay{xyzy\approx xzy,\qquad xy^2\approx xy,\qquad
x\le x^2,\qquad xy\le x+zy. \quad\textup{(\ensuremath{F})}}

Here \(p\leq q\) abbreviates p+q approximately q. The proof below has no bound on the variables, word lengths or number of right summands.

\par\smallskip\noindent\textit{2. Normal form and exact absorption criterion.}\quad 
The two word identities preserve the literal head and retain the tail letters in their order of last occurrence, with no repeats within that tail. To remove an earlier tail occurrence of y, use \(xyzy=xzy\) when there is a nonempty intervening block, and \(xy^{2}=xy\) otherwise. The prefix substituted for x is nonempty because the deleted occurrence is in the tail. The head is never deleted, even when it recurs in the tail. In particular, xx remains xx, and \(x^{2}y=xy\) is not being asserted: \(x=3\), \(y=0\) refutes it in this algebra.

Let \(suf_{u}\)(x) be the set of letters strictly after the last x in u. For a nonempty collection P of original right words, a one-letter target y is absorbed exactly when some original right word has content \{y\}. For a longer normal target \(w=hT\), absorption is equivalent to:

\textbullet{} H: some original right word v has head h and c(v) subset c(w).

\textbullet{} O: for every x in the target tail T, some original right word q contains x in its tail (strictly after its head), and \(suf_{q}\)(x) is contained in \(suf_{w}\)(x).

The H witness can be literal h. O requires a tail occurrence, so a literal x, or a word containing x only as its head, does not qualify. The O witnesses may contain variables outside the target content.

If H fails, assign \(h=0\), the other target variables 1, and outsiders 2. The target is 0; right words with head different from h have value 1 or 2, while any word with head h must have an outsider and is 2. The right sum is at most 1. If O fails at x, assign \(x=3\), all variables in \(suf_{w}\)(x) the value 1, and all other variables 2. The target is 4. A right word with a tail x must have a 2 after its last x, so is 2. A right word with no tail x has value at most 3, even if its head is x. Hence the sum is at most \(3<4\). A one-letter target without a pure witness is separated by assigning \(y=0\) and outsiders 2.

Conversely, when the target has a last tail reset giving 4, O gives an original word with the corresponding tail occurrence followed only by neutral 0/\allowbreak{}1-valued variables, so that word is 4. A reset giving 2 is least and needs no witness. With no target tail reset, all target tail variables lie in \{0,\allowbreak{}1\}. The H witness has the same head; its other letters lie in the tail content or equal h. For \(h=0\) or 1 its value equals h; for \(h=2\) it is at least 2; for \(h=3\) its value is 3 or 4; for \(h=4\) its value is 4. Thus it bounds the target in every case. For pure-variable witnesses, squaring fixes 0,\allowbreak{}1,\allowbreak{}2,\allowbreak{}4 and sends 3 to 4, and higher powers equal the square.

This proves the criterion for arbitrary original words: the two word laws preserve H, tail occurrences and all relevant last-occurrence suffix sets. It also proves soundness of (F). For \(xy\leq x+zy\), literal x supplies H and the tail y of zy supplies O. The two word laws follow from the right-map formula, and \(x\leq x^{2}\) follows from squaring.

\par\smallskip\noindent\textit{3. Appending the target tail to every original right word.}\quad 
Suppose \(w=hT\leq \sum (P)\) is valid and nonlinear. Build the target tail from right to left. Maintain, for every original right word u,

\CaseDisplay{uS\le\sum P.}

Initially S is absent and this is just the inclusion of each original summand. An absent context here never denotes an empty substitution for a variable.

Let x be the next target tail variable. By O, choose an original \(q=R\) x Q at its last x, where R is nonempty and c(Q) is contained in c(S). The maintained statement includes \(qS\leq \sum (P)\). The two word laws give qS approximately R x S: all letters of Q are removed by their later occurrences in S. If S is absent then Q is also absent.

For each original u, use \(xy\leq x+zy\) with \(X=u\), \(Y=x\), \(Z=R\), in the common right context S:

\CaseDisplay{uxS\le uS+RxS\le\sum P.}

All three substitutions are nonempty, including R because x occurs in the witness tail. Thus the maintained statement holds with xS. After finitely many stages it gives \(uT\leq \sum (P)\) for every original right word, including every original linear word.

\par\smallskip\noindent\textit{4. Recovering the target and proving completeness.}\quad 
Choose H's original witness v. Its head is h and its content is contained in c(hT), so normalizing vT yields either hT itself, or hhT. The second possibility occurs only if h is absent from T and appears in the tail of v. In that case \(x\leq x^{2}\), in right context T, gives \(hT\leq hhT\leq \sum (P)\). This completes the nonlinear absorption.

For a one-letter target y, normalize its pure-variable witness to y or yy and use \(y\leq yy\) if needed. Therefore every valid word absorption is derivable. Distributivity expands arbitrary terms into sums of nonempty words; deriving every summand's absorption by the opposite side yields every valid identity. This proves the four-law theorem.

Reversing words transports the basis and derivation to the opposite algebra. Catalogue transfers still require independently verified complete HSP inclusions in both directions. Finite regressions and replay of context/substitution steps support the unrestricted proof, not replace it.

\subsection{Finite bases for 5223 and 5228 from literal-activated suffixes}\label{case:152}
Signature (+, multiplication), without constants. All words and all substituted terms are nonempty. The bases are relative to the additively idempotent semiring axioms. Write \(u\leq P\) for \(u+P=P\).

\par\smallskip\noindent\textit{Tables and word values.}\quad 
Both additions have rows 00044,\allowbreak{}01144,\allowbreak{}01234,\allowbreak{}44334,\allowbreak{}44444. Thus 2 is least, \(2<1<0<4\) and \(2<3<4\), with \(1+3=0+3=4\). This last join differs from 5222 and 5227 and introduces a new, indispensable branch. Multiplication rows for 5223 are 00244,\allowbreak{}11244,\allowbreak{}22244,\allowbreak{}33244,\allowbreak{}44244. For 5228 they are 00244,\allowbreak{}11244,\allowbreak{}22244,\allowbreak{}44244,\allowbreak{}44244.

For 5223 write a word as its head followed by its tail. The last tail factor in \{2,\allowbreak{}3,\allowbreak{}4\}, when present, resets the value to 2 if that factor is 2, and to 4 otherwise. With no such factor the head survives. This follows by induction: right 0/\allowbreak{}1 is the identity, right 2 is constant 2, and right 3/\allowbreak{}4 is constant 4.

For 5228 a literal retains its value. In a word of length at least two, the last factor in \{2,\allowbreak{}3,\allowbreak{}4\}, including the head, determines 2 versus 4. If no factor lies in that set, the value is the head. Right 0/\allowbreak{}1 changes 3 to 4 and fixes the other four elements, proving the claim by induction.

\par\smallskip\noindent\textit{Two explicit finite bases.}\quad 
For 5223 the following five identities are complete:

{\small\par\smallskip
\CaseLine{\(xyzy = xzy;  xyy = xy;  x \leq  xx;\)}
\CaseLine{\(xy \leq  x+zy;  xyz \leq  xz+yz+z.\)}
\smallskip}

For 5228 the following eight identities are complete:

{\small\par\smallskip
\CaseLine{\(xyzy = xzy;  xyy = xy;  xxy = xy;  x \leq  xx;\)}
\CaseLine{\(xy \leq  xx+zy;  xyz \leq  xz+yz;  xy \leq  x+zxy;\)}
\CaseLine{\(xyzt \leq  xzt+y+zt.\)}
\smallskip}

The last identity in each list is essential to this proof's new literal-activation step. Each is sound by the criterion proved below, on its displayed distinct variables; therefore all term substitutions are sound. The word laws follow directly from the word-value formulas.

The first two word laws retain the literal head and delete earlier occurrences in the tail, keeping the tail's last-occurrence order. An earlier tail y in A y B y D is deleted by \(xyzy=xzy\) when B is nonempty, and by \(xyy=xy\) otherwise. The prefix A is nonempty. In 5228, \(xxy=xy\) then deletes a repeated leading head when a further position remains. A nonlinear singleton stays hh, while a literal stays h. Consequently every nonlinear word u satisfies \(u^{2}=u\) in 5228: for \(u=ab\), both a, b nonempty, \(abab=aab=ab\). This equality is never applied to a literal.

\par\smallskip\noindent\textit{Exact criteria in the original words.}\quad 
Let P be a nonempty collection of original upper words. Put \(C=c(w)\), \(h=head(w)\), and S\_x=c(the substring after the last x in w). Both criteria first require H: an original same-head word v has c(v) subset C. A literal target is absorbed exactly when H holds, equivalently when an original pure power of that literal is present.

For a nonlinear 5223 target require, for every variable x occurring in the target tail, either:

1. an original word has a true tail occurrence of x and its last-x suffix content is contained in \(S_{x}\); or

2. an original word is headed by x, has no further x, and has its tail contained in \(S_{x}\), AND an original word has content contained in \(S_{x}\).

For a nonlinear 5228 target require, for EVERY x in C, either:

1. an original nonlinear word contains x with last-x suffix contained in \(S_{x}\); or

2. the original literal x occurs AND an original word has content contained in \(S_{x}\).

In particular, the second alternatives are impossible when \(S_{x}\) is empty. Literal words cannot be silently treated as nonlinear witnesses. The original contained suffix word need not be a literal. These conditions are invariant under the stated word normalizations.

Necessity of H: assign \(h=0\), other target variables 1, and outsiders 2. Without H every upper value is 1 or 2, while the target is 0. For a failure of a suffix requirement at x, first consider the absence of both a strong witness and its weak witness (a head-only word for 5223, or literal x for 5228). Assign \(x=4\), \(S_{x}=1\), and every other variable 2. The target is 4; every upper value is 1 or 2. If the weak witness exists but the strong witness and suffix-contained word do not, instead assign \(x=3\), \(S_{x}=1\), and every other variable 2. The target is again 4: x is in the tail for 5223, and the target is nonlinear for 5228. Every upper word has value 2 or 3; a value 1 would require a suffix-contained word. Their sum is at most 3. Thus every failed condition separates.

For sufficiency, value 2 is least. For 5223 a target with no tail reset is bounded by H: when its head is 0 or 1 the H witness has that same value; head 2 is least; head 3 gives upper 3 or 4; head 4 gives upper 4. A tail reset to 4 has a last special variable x assigned 3 or 4 and all \(S_{x}\) variables assigned 0/\allowbreak{}1. A strong witness gives 4. A weak witness gives the assigned x, and the suffix-contained witness gives 0 or 1, whose join with 3 or 4 is 4. For 5228 any nonlinear target of value 4 is handled identically, using nonlinear versus literal witnesses. A 0/\allowbreak{}1 target uses H. Pure powers absorb their literal because squaring fixes 0,\allowbreak{}1,\allowbreak{}2,\allowbreak{}4 and sends 3 to 4. This proves both criteria for arbitrary assignments, variables, lengths, and numbers of upper words.

\par\smallskip\noindent\textit{Unrestricted completeness for 5223.}\quad 
Normalize a nonlinear \(w=hT\). Maintain \(uS\leq P\) for EVERY original u, where S is the part of T constructed from right to left. Initially S is absent. For the next x, a strong original witness \(q=R\) x Q has R nonempty and c(Q) subset c(S). The maintained qS normalizes to R x S. Apply \(xy\leq x+zy\) with \(X=u\), \(Y=x\), \(Z=R\) in common right context S to obtain uxS.

For a weak witness \(q=xQ\), choose original z contained in S. Here S is nonempty. Maintained \(qS=xS\) and zS are available. Use \(xyz\leq xz+yz+z\) with \(X=u\), \(Y=x\), \(Z=z\) in common right context S. Its upper words normalize to uS, xS, zS and its lower word to uxS, because every letter of z has a later occurrence in S. All substitutions remain nonempty.

After finitely many stages vT is available for H's original v. It normalizes to hT or hhT. In the latter case \(x\leq xx\) in right context T recovers hT. A literal target follows from its pure original witness and \(x\leq xx\). Hence the five laws derive every valid absorption.

\par\smallskip\noindent\textit{Unrestricted completeness for 5228.}\quad 
Let T be the full last-occurrence list of a nonlinear target w, including its head variable. Maintain \(uS\leq P\) for EVERY original NONLINEAR u, constructing T from right to left. Such originals exist: the empty-suffix condition at the final variable forces one.

A strong witness \(q=R\) x Q has c(Q) subset c(S). If R is nonempty, \(xy\leq xx+zy\) with \(X=u\), \(Y=x\), \(Z=R\) and right context S yields uxS, since \(u^{2}=u\) for nonlinear u and \(qS=RxS\). If R is absent, \(q=xQ\) is nonlinear, so Q and S are nonempty; \(xyz\leq xz+yz\) with \(X=u\), \(Y=x\), \(Z=S\) yields uxS.

In the weak case literal x is original and some original z is contained in the nonempty S. If z is nonlinear, maintained zS is available. If z is a literal, choose any original nonlinear v. The law \(xy\leq x+zxy\) with \(X=z\), \(Y=S\), \(Z=v\) yields zS, because vzS normalizes to vS. This step protects the literal explicitly; it does not apply \(u^{2}=u\) to it. Now \(xyzt\leq xzt+y+zt\) with \(X=u\), \(Y=x\), \(Z=z\), \(T=S\) yields uxS: the upper words are uS, literal x, and zS after normalization. Thus every step uses finitely many applications of the eight fixed laws.

At completion, if H's witness v is nonlinear then vT normalizes to w. If H's witness is literal h, choose any original nonlinear v and use \(xy\leq x+zxy\) with \(X=h\), \(Y=T\), \(Z=v\). Because h belongs to T, vhT normalizes to vT, which is available; hT normalizes to w. A literal target is again recovered from an original pure power. This proves completeness with no bound on variables, lengths or original polynomial size.

Distributivity expands arbitrary terms as nonempty sums of nonempty words. Apply the proved absorption construction to each summand of each side of any valid identity; the two inequalities give the identity. Reversing all words proves the opposite bases. Transfers to catalogue types require complete mutual HSP checks; finite regressions support these proofs but do not replace the unrestricted arguments.

\subsection{Seven identities for 5227 from nonlinear suffix witnesses}\label{case:153}
Signature (+, multiplication), without constants, relative to the additively idempotent semiring axioms. A nonlinear word means length at least two. Substitutions below are always nonempty.

\par\smallskip\noindent\textit{1. Algebra and arbitrary word values.}\quad 
In element order 0,\allowbreak{}1,\allowbreak{}2,\allowbreak{}3,\allowbreak{}4 the catalogue tables are

{\small\par\smallskip
\noindent addition        multiplication\par
\CaseLine{\ttfamily 0 0 0 4 4\qquad 0 0 2 4 4}
\CaseLine{\ttfamily 0 1 1 3 4\qquad 1 1 2 4 4}
\CaseLine{\ttfamily 0 1 2 3 4\qquad 2 2 2 4 4}
\CaseLine{\ttfamily 4 3 3 3 4\qquad 4 4 2 4 4}
\CaseLine{\ttfamily 4 4 4 4 4\qquad 4 4 2 4 4}
\smallskip}

Addition has order \(2<1<0<4\) and \(1<3<4\), with 0 and 3 incomparable. For a nonlinear word, its last factor in \{2,\allowbreak{}3,\allowbreak{}4\}, if any, determines the value: 2 gives 2 and 3 or 4 gives 4. If no such factor occurs, all factors are 0/\allowbreak{}1 and the value is its head. A linear word retains its assigned value. This follows by induction from the table: right 0/\allowbreak{}1 fixes 0,\allowbreak{}1,\allowbreak{}2,\allowbreak{}4 and changes 3 to 4, right 2 resets to 2, and right 3/\allowbreak{}4 resets to 4. In particular, any nonlinear word has value in \{0,\allowbreak{}1,\allowbreak{}2,\allowbreak{}4\}; even a head-only 3 is changed to 4 by its later position.

\par\smallskip\noindent\textit{2. The seven identities and normal forms.}\quad 
The complete relative basis is

{\small\par\smallskip
\CaseLine{\((W1) xyzy = xzy,\)}
\CaseLine{\((W2) xyy = xy,\)}
\CaseLine{\((W3) xxy = xy,\)}
\CaseLine{\((L)  x \leq  xx,\)}
\CaseLine{\((A)  xy \leq  xx+zy,\)}
\CaseLine{\((B)  xyz \leq  xz+yz,\)}
\CaseLine{\((H)  xy \leq  x+zxy.\)}
\smallskip}

As usual \(u\leq P\) abbreviates \(u+P=P\). For a word w let L(w) be its distinct variables in last-occurrence order. The word laws retain the literal head and L(w), together with the distinction between length one and longer words of singleton content. First use (W1) or (W2) to delete every earlier repeated tail occurrence, retaining the head. If the resulting tail starts with another copy of the head and has a later position, (W3) deletes that extra leading copy. Thus a nonlinear singleton normalizes to hh. For at least two variables, the normal form is L(w) if its first variable is the head h, and hL(w) otherwise. A linear h stays h.

This is a literal derivation: (W1) uses a nonempty prefix and a nonempty intervening block, (W2) handles adjacent tail repeats, and (W3) uses a nonempty suffix. Absent contexts do not denote empty substitutions. The word-value formula proves soundness of (W1)--(W3).

A useful consequence is

uu = u  for every nonlinear word u.                         (I)

Write \(u=ab\) with both factors nonempty. Then \(abab=aab=ab\) by (W1) with \(X=a\), \(Y=b\), \(Z=a\), followed by (W3). In contrast, a one-letter word need not satisfy this equality: at \(x=3\), \(x=3\) while \(xx=4\).

\par\smallskip\noindent\textit{3. Exact criterion.}\quad 
Let P be a nonempty collection of original words. Let \(suf_{w}\)(x) be the content strictly after the last occurrence of x in w. A linear target y is absorbed by P exactly when some original word has content \{y\}. For a nonlinear w, absorption holds exactly when:

1. P contains a word with the same head as w and content contained in c(w).

2. For every x in c(w), P contains a nonlinear word q with x in q and \(suf_{q}\)(x) subseteq \(suf_{w}\)(x).

The second condition also applies to the head when it occurs only at the first position. It allows x to occur only as the head of its witness q, provided q is nonlinear. This is the point where the criterion differs from the genuine-tail criterion for 5222. The same-head contained witness in condition 1 may be linear.

If condition 1 fails, assign the target head 0, the other target variables 1, and outsiders 2. The target is 0. Every noncontained right word has value 2; every contained right word has head different from the target and value 1. This separates the absorption. The same assignment proves the pure-word requirement for a linear target.

For a failure of condition 2 at x, assign \(x=3\), all variables of \(suf_{w}\)(x) value 1, and all other variables 2. The nonlinear target is 4. A nonlinear original word has value 4 exactly when it contains x and has no 2 after its last x, precisely the required witness. All remaining nonlinear right words have value 1 or 2. A linear right word has value 1,\allowbreak{}2 or 3. Without the witness, their join is at most 3 and cannot absorb 4. This uses \(1+3=3\), not \(1+3=4\).

Conversely, target value 2 is least. A nonlinear target value 4 has a last special variable x valued 3 or 4. Its witness from condition 2 is nonlinear and has only 0/\allowbreak{}1 after its last x, so it also has value 4. If the target value is 0 or 1, all its variables are 0/\allowbreak{}1 and the contained same-head witness has the same value. For a linear target, a pure-word witness absorbs its variable: squaring fixes 0,\allowbreak{}1,\allowbreak{}2,\allowbreak{}4 and sends 3 to 4; higher powers equal squares. This proves the criterion for all word lengths and variable counts.

It also proves soundness of the proposed inequalities. (L) is the pure-word case. For (A), xx is a contained same-head witness and also a nonlinear suffix witness for x; zy witnesses y. For (B), xz is the contained same-head witness and witnesses x and z, while yz witnesses y. For (H), literal x supplies condition 1 and zxy supplies the nonlinear suffix witnesses. These are identities on distinct displayed variables; arbitrary identifications and nonempty term substitutions then follow by equational logic.

\par\smallskip\noindent\textit{4. Appending a target tail to nonlinear original words.}\quad 
Normalize a valid nonlinear target as \(w=hT\), with T a nonempty list of distinct tail variables. Normalize P as well; its original nonlinear words remain nonlinear. Construct T from right to left. For every nonlinear original word u maintain

{\small\par\smallskip
\CaseLine{\(uS \leq  \sum (P),\qquad (M)\)}
\smallskip}

where S is the tail already constructed. Initially S is absent and (M) is just the inclusion of each original nonlinear summand. This assertion is deliberately not made for the linear summands.

Let x be the next letter of T. Condition 2 supplies a nonlinear \(q=R\) x Q, split at its last x, with c(Q) subseteq c(S). R or Q may be absent. The maintained assertion for q gives \(qS\leq P\). Delete the letters of Q using their later occurrences in S, by (W1),(W2), to get \(qS=R\) x S. If S is absent, Q is absent as well.

If R is nonempty, apply (A) with \(X=u\), \(Y=x\), \(Z=R\) in right context S:

{\small\par\smallskip
\CaseLine{\(uxS \leq  u^{2} S + R x S = uS + qS \leq  P.\)}
\smallskip}

The equality \(u^{2}\) \(S=uS\) is consequence (I), since u is nonlinear. This includes an absent S and an absent Q, without empty substitutions.

If R is absent, \(q=xQ\). Since q is nonlinear, Q is nonempty; since c(Q) subseteq c(S), S is also nonempty. Now (B), with \(X=u\), \(Y=x\), \(Z=S\), gives

{\small\par\smallskip
\CaseLine{\(uxS \leq  uS+xS = uS+qS \leq  P.\)}
\smallskip}

Thus (M) is preserved in both cases. After all stages, \(uT\leq P\) holds for every original nonlinear word u. There is no bound on the number of stages or original words.

\par\smallskip\noindent\textit{5. Recovering the head, including a literal same-head witness.}\quad 
Let v be the contained same-head witness from condition 1. If v is nonlinear, the construction gives \(vT\leq P\). Since c(v) is contained in c(hT), deleting tail duplicates reduces vT to hT or hhT; the latter reduces to hT by (W3), because T is nonempty. This gives \(w\leq P\).

If v is linear, it is the literal h. Condition 2 at h supplies a nonlinear \(q=R\) h Q with c(Q) subseteq \(suf_{w}\)(h), hence subseteq c(T). The construction gives \(qT\leq P\), and tail deletion gives \(qT=R\) h T. If R is absent this is already w. Otherwise (H), with \(X=h\), \(Y=T\), \(Z=R\), gives

{\small\par\smallskip
\CaseLine{\(hT \leq  h+R h T = v+qT \leq  P.\)}
\smallskip}

Every substituted word here is nonempty. This proves the nonlinear absorption even when the only contained same-head witness is literal. For a linear target, its pure-word witness reduces to h or hh; (L) completes the derivation. No nonlinear-idempotence assertion has been applied to a linear original summand.

Theorem. The seven identities in Section 2 form a complete finite basis for 5227 relative to the ai-semiring axioms. The construction derives every valid word absorption with arbitrarily many variables, arbitrary word lengths and arbitrarily many right summands. Expand arbitrary terms distributively and derive each summand's absorption into the opposite side; this yields every identity. Reversing words transports the basis and derivations to the multiplicative opposite.

Complete two-way HSP certificates are required to transfer this classification to other catalogue entries. Finite table tests and replay of substitutions support, but do not replace, this unrestricted completeness proof. The preceding 5410 classification is separate and is not counted again.

\subsection{A finite basis for 5315 by recovering simple variables}\label{case:154}
Signature (+, multiplication), without constants, relative to the additively idempotent semiring axioms. A simple variable occurs exactly once in a word. A repeated variable occurs at least twice. The finite family below is specified effectively, not enumerated.

\par\smallskip\noindent\textit{1. Table semantics.}\quad 
In element order 0,\allowbreak{}1,\allowbreak{}2,\allowbreak{}3,\allowbreak{}4 the tables are

{\small\par\smallskip
\noindent addition        multiplication\par
\CaseLine{\ttfamily 0 0 0 0 4\qquad 0 0 3 3 0}
\CaseLine{\ttfamily 0 1 1 1 4\qquad 1 1 3 3 1}
\CaseLine{\ttfamily 0 1 2 2 4\qquad 2 2 3 3 2}
\CaseLine{\ttfamily 0 1 2 3 4\qquad 3 3 3 3 3}
\CaseLine{\ttfamily 4 4 4 4 4\qquad 0 0 2 3 4}
\smallskip}

The additive order is the chain \(3<2<1<0<4\). Let K=\{0,\allowbreak{}1,\allowbreak{}4\}. The table gives these word values by induction:

\textbullet{} Any occurrence of 3 gives 3. At least two occurrences valued 2 also give 3.

\textbullet{} With exactly one occurrence valued 2, the word is 2 precisely when all earlier positions are 4 and every later position is in K; otherwise it is 3.

\textbullet{} When all factors are in K, the value is 4 if all are 4. Otherwise it is 1 if its literal head is 1, and 0 in every other case.

The last clause includes a head 0, or a head 4 followed somewhere by 0 or 1. In particular, 4 is a right identity but not a left identity. These rules include linear words.

\par\smallskip\noindent\textit{2. Exact absorption criterion.}\quad 
For a nonempty collection P of original words, let U consist of the words with content contained in \(C=c(w)\). Then \(w\leq \sum (P)\) holds exactly when the following conditions hold:

(H) some word of U has the same head as w;

(S) for each simple variable x of w, some u in U either avoids x, or contains x exactly once with \(pre_{u}\)(x) subseteq \(pre_{w}\)(x).

Here pre is content strictly before the indicated occurrence. Witnesses in (S) need not share the target head. The option that a witness omits x is essential. A linear target requires a literal summand, since its sole simple variable cannot be omitted by a nonempty contained word.

For necessity of (H), assign the target head 0, other target variables 1, and outsiders 3. The target is 0; without a contained same-head word all right values are 1 or 3. For a failed (S) at x, assign \(x=2\), its target prefix 4, the remaining target variables 0, and outsiders

3. The target is 2. Every noncontained right word is 3. Every contained word includes x, either repeats it or has a 0 before it, and is therefore also 3.

For sufficiency, target value 3 is least. If the target uses only K, its contained same-head witness has the same head; it has value 4 when all target variables are 4, value 1 at head 1, and value 0 or 4 when a head-4 target has value 0. It thus absorbs the target. The only remaining target value is 2. Its unique 2-valued occurrence belongs to a simple variable x; its earlier variables are 4 and its later variables are in K. A contained witness omitting x is K-valued and hence at least \(1>2\). Otherwise the (S) witness has a unique x, with only 4 before it, and has value 2. This proves the criterion for arbitrary words and assignments.

\par\smallskip\noindent\textit{3. Structural identities and block forms.}\quad 
The following six identities are sound:

{\small\par\smallskip
\CaseLine{\((W1) x^{3} = x^{2},\)}
\CaseLine{\((W2) xyx = x^{2} y,\)}
\CaseLine{\((W3) xy^{2} z^{2} = xz^{2} y^{2},\)}
\CaseLine{\((D)  xy = x^{2} y + xy^{2},\)}
\CaseLine{\((C1) x^{2} \leq  x,\)}
\CaseLine{\((C2) xy^{2} \leq  x.\)}
\smallskip}

The table formula proves (W1),(W2): the same content, head, and simple-variable prefixes remain. For (W3), swapping the two repeated tail variables cannot change the prefix set of a simple variable: its displayed prefix x is unchanged, and y, z are not simple. For (D), one of the two factors must have no 2-valued occurrence in any nonzero product. Squaring that factor preserves the product, while squaring a factor containing the sole 2 changes its value to

3. In the K-valued case both squares preserve the product. If the original product is 3, both right products are 3. Thus their join is the original product. (C1),(C2) follow from the same formula, or directly from the five-element table.

By (W2), gather every repeated variable into a block at its first occurrence. By (W1), reduce each repeated block to a square. The result lists the distinct variables in first-occurrence order, with each variable occurring once or twice according as it was simple or repeated. The literal head is retained. Within a segment of repeated blocks after the head, (W3) permutes adjacent blocks. It never moves the head, nor crosses a simple variable. These are literal equational derivations with nonempty substituted blocks.

\par\smallskip\noindent\textit{4. A saturated guard and recovery of each simple variable.}\quad 
Fix a valid absorption \(w\leq P\) and put \(C=c(w)\), \(h=h(w)\). Choose the contained same-head original v supplied by (H). Applying (C1) at each position of v makes every variable repeated; normalization gives a word with head h, content c(v), and every variable doubled. Append a square of each missing variable of C by (C2). Permuting repeated tail blocks by (W3) gives a fixed saturated word g of full content C and head h. Thus \(g\leq P\), with an explicit finite derivation.

Let x be a simple variable of w. Use (S) to choose an original contained witness u. Apply (C1) to every occurrence of every variable other than x in u, and normalize. Call the result \(u_{x}\). It still avoids x or contains it once with prefix contained in \(pre_{w}\)(x), and \(u_{x}\leq u\leq P\). Every other variable of \(u_{x}\) is repeated.

Let \(w_{x}\) be obtained from w by repeating every variable other than x and normalizing. The two-right-word absorption

{\small\par\smallskip
\CaseLine{\(w_{x} \leq  g+u_{x}\qquad (R)\)}
\smallskip}

is valid by the exact criterion: g supplies (H) with full content, and \(u_{x}\) supplies (S) for the sole simple variable x. These will be derived from a fixed finite family in Section 5, not assumed as infinitely many extra rules. Hence each \(w_{x}\) is below P.

To reconstruct w, use (D) whenever a current word has at least two simple variables. Split it into nonempty factors X, Y between two simple variables and replace it by \(X^{2}Y\)+\(XY^{2}\). Each summand has strictly fewer simple variables. Normalization preserves the head, the content, and the original prefix set of every surviving simple variable. Iteration terminates with a sum of words having at most one simple variable. A leaf with no simple variables is equal to g by (W1)--(W3). A leaf whose sole simple variable is x is equal to \(w_{x}\): its prefix set at x agrees, and repeated blocks on each side can be permuted after the fixed head. Every leaf is therefore below P, and the equality from (D) gives \(w\leq P\). This also covers targets with zero or one simple variable from the outset.

\par\smallskip\noindent\textit{5. Reversible compression and a fixed finite basis.}\quad 
In (R) there are at most three block words. Reserve x and every head of these words, at most four variables in total. Each other variable is either absent, before x, or after x in each word. When x is absent from a word, assign all its variables to the before region. A repeated x is treated as its single squared block. There are at most \(3^{3}-1=26\) nonempty region-incidence patterns.

Identify the nonreserved variables of each pattern, preserving the reserved variables separately, and normalize repeated powers by (W1),(W2). The compressed valid absorption uses at most 30 variables, has block length at most 60, and has at most two right words. Every compressed class variable is still repeated. Its inverse substitution is the nonempty product of its original members. A squared class block therefore expands to the square of that product, not to a single occurrence of each member. (W1),(W2) recover a square of each member; (W3) sorts the repeated blocks within their original region. All pattern members had the same region in every word. The heads, x and its multiplicity are reserved, so this inverse recovers each original word modulo (W1)--(W3). This explicitly preserves common multiplicity and does not identify single occurrences into repeats.

Fix a 30-variable alphabet. Let Sigma contain the six identities in Section 3 and all valid absorptions over that alphabet with one block word on the left, one or two block words on the right, and length at most 60 in every word. A block word uses each variable in one block of length one or two. Section 2 decides this finite family effectively. No small enumerated list or minimality is claimed.

Theorem. Sigma is a finite identity basis for 5315. The guard and witness contractions use the six displayed identities. Each local recovery (R) is derived by expanding a compressed member of Sigma. The splitting argument then derives every word absorption, without a bound on its original variables, length or summands. Distributivity expands arbitrary terms into polynomials; summandwise absorption in both directions proves every identity. Reversing words gives a finite basis for the multiplicative opposite.

Independent table checks and literal substitution replay are error checks supporting this unrestricted proof. Further catalogue entries require independently verified complete HSP inclusions in both directions.

\subsection{A finite basis for 5410 from content coverage and a compatible head}\label{case:155}
Signature (+, multiplication), without constants. The basis is relative to the additively idempotent semiring axioms. A word is a nonempty string; a linear word means a word of length one. The finite family below is explicitly specified, not enumerated as a small list.

\par\smallskip\noindent\textit{1. Tables, word values, and normal forms.}\quad 
The operation tables of catalogue entry 5410, in order 0,\allowbreak{}1,\allowbreak{}2,\allowbreak{}3,\allowbreak{}4, are

{\small\par\smallskip
\noindent addition        multiplication\par
\CaseLine{\ttfamily 0 0 0 0 4\qquad 0 0 3 3 4}
\CaseLine{\ttfamily 0 1 0 1 4\qquad 1 1 3 3 4}
\CaseLine{\ttfamily 0 0 2 2 4\qquad 2 2 3 3 2}
\CaseLine{\ttfamily 0 1 2 3 4\qquad 3 3 3 3 3}
\CaseLine{\ttfamily 4 4 4 4 4\qquad 4 4 3 3 4}
\smallskip}

The additive order is \(3<1<0<4\) and \(3<2<0<4\), with 1 and 2 incomparable. Put K=\{0,\allowbreak{}1,\allowbreak{}4\}. For a word w let h be its first variable, T the content of the substring after that first position, and \(C=\{h\}\) union T. T may be empty and may contain h.

The tables give the following formula, by induction on length. If any tail factor has value 2 or 3, the word is 3. Otherwise all tail factors are in K. A head 3 gives 3 and a head 2 gives 2. For a head in K, the value is 4 if the head is 4 or some tail factor is 4; otherwise the value is the head, 0 or 1. This includes words of length one.

In particular the identities

{\small\par\smallskip
\CaseLine{\(xyz = xzy,\qquad xyy = xy\qquad (W)\)}
\smallskip}

hold. Fix an order on variables. The first identity swaps adjacent tail letters using a nonempty prefix as its x substitution; the second deletes an adjacent repeated tail letter. Hence (W) gives the normal form h times the increasing list of T. If T is empty the form is h. If h lies in T, its extra tail occurrence is retained. This yields a literal derivation by swaps and deletions, with no empty substitutions.

\par\smallskip\noindent\textit{2. Exact absorption criterion.}\quad 
Let P be a nonempty collection of original words, and write \(w\leq P\) for \(w+\sum (P)=\sum (P)\). Let U consist of the words of P whose content is contained in \(C=c(w)\). Then \(w\leq P\) holds exactly when:

(C) the union of the contents of U is C;

(H) some v in U has head \(h=h(w)\); if h does not occur in the tail of w, this witness must also have no occurrence of h in its own tail.

For a one-letter target, these conditions require a literal summand. Repeated powers are not substitutes for that literal: \(x\leq xx\) fails at \(x=2\). The reverse inequality \(xx\leq x\) holds, but is not used to change (H).

Necessity of (C). If x in C is not covered, assign \(x=4\), all other variables in C the value 1, and outsiders 3. The target is 4. An original word containing an outsider is 3; a contained word avoids x and is 1. Their join is at most 1, so cannot absorb the target.

Necessity of (H) when the target head is repeated. Assign \(h=0\), the other variables of C value 1, and outsiders 3. The target is 0. Without a contained same-head word every right word is 1 or 3, and their join is at most 1.

Necessity of (H) when the head is simple. Assign \(h=2\), other variables in C value 1, and outsiders 3. The target is 2. A right word has value 2 precisely when it is contained, has head h, and does not repeat h in its tail. If no such word exists, all right values are 1 or 3, whose join cannot absorb 2.

Sufficiency. Target value 3 is least. Target value 4 requires all C-valued factors to lie in K and at least one variable to have value 4. Condition (C) supplies a contained word containing that variable, also valued 4. If the target is 0 or 1, every variable of C lies in \{0,\allowbreak{}1\}; the same-head contained witness from (H) has exactly the target value. Finally, a target value 2 has head h valued 2 and tail in K, so h is simple in the target. The (H) witness has the same head, does not repeat h, and has its tail in the target tail set. It therefore has value 2. This proves sufficiency on every assignment.

\par\smallskip\noindent\textit{3. A derivation using at most two available words per step.}\quad 
Normalize the words in a valid absorption \(w\leq P\), and choose a witness v from (H). Begin with \(g=v\leq P\). Let T be the desired tail set and R the current tail set of g. Then g has head h and R subseteq T. This is true also for a simple target head, by the extra clause of (H).

If c(g) is not yet C, choose a variable missing from c(g). Condition (C) supplies an original u contained in C and containing that variable. Replace g by the normal word with head h and tail

{\small\par\smallskip
\CaseLine{\(R' = R \cup  (c(u) \cap  T).\)}
\smallskip}

Call the new word g'. Its content is exactly c(g) union c(u), since the only possible member of C outside T is the already present head. The absorption

{\small\par\smallskip
\CaseLine{\(g' \leq  g+u\qquad (M)\)}
\smallskip}

satisfies (C), and g itself is its compatible same-head witness (H). Thus (M) is valid. Its words have a common target content bound, and its left content strictly increases. Repeating this step at most |C|-1 times gives a word g of content C and tail R subseteq T.

When the target head is simple, \(R=T\) and \(g=w\). When it is repeated, the only potentially missing tail letter is h. If h is missing, the one-right-word absorption \(w\leq g\) satisfies (C) and (H): w permits the head repetition, while g is a contained same-head witness. This final step therefore recovers w. Linear targets were already literal summands. The construction has no bound on |C| or |P|.

This section uses valid local absorptions. The next section derives each of them from one fixed finite set of identities, so it does not presuppose an infinite inference rule as part of the basis.

\par\smallskip\noindent\textit{4. Reversible compression and the specified finite basis.}\quad 
In any local absorption just used there are at most three normal words: one on the left and one or two on the right. Reserve every variable that is the head of one of these words. There are at most three reserved variables; their appearances in tails remain recorded. Every other variable occurs only in tails. Partition these variables by their nonempty incidence pattern in the at most three tail sets. There are at most \(2^{3}-1=7\) such classes.

Identify each class to one fresh variable, leave reserved variables separate, and normalize using (W). This is an equational substitution, so validity is preserved. The result has at most 10 variables, each word has length at most 11, and there are at most two right words.

The compression is reversible modulo (W). Substitute the nonempty product of all original class members for its new variable. Since those members had identical tail incidence, expansion recovers each original tail set. The reserved heads and their tail occurrences have not changed. Tail swaps and duplicate deletions then recover exactly the original normal words. No multiplicity beyond one per tail set is being discarded; an occurrence at a head was reserved throughout.

Let X be a fixed ten-variable alphabet. Define Sigma to consist of (W) and all valid absorptions with one normal left word and one or two normal right words over X, each word of length at most 11. The criterion in Section 2 explicitly decides membership in this finite family. The family includes the required one-variable absorptions. No minimality or small enumerated basis is asserted.

Theorem. Sigma is a finite identity basis for 5410.

Soundness follows from its definition and (W). Normalize any valid word absorption. Section 3 supplies a finite chain of valid local absorptions using at most two available right words. Each is an expanded substitution instance of an identity in Sigma, by the reversible compression. If the right words are already absorbed by P, add those absorptions and use transitivity to absorb the new word by P. This proves the original absorption. Distributivity expands arbitrary terms into nonempty polynomials; deriving the absorption of every summand into the opposite side proves every valid equality. This argument has no restriction on original variables, word lengths, or number of summands. Reversing words gives a finite basis for the multiplicative opposite. \(\square\)

\subsection{Two finite bases from a simple-head guard and ordered suffix witnesses}\label{case:156}
Signature (+, multiplication), without constants. Relative to the ai-semiring axioms, 5413 has a four-identity basis and 5414 has an eight-identity basis. All substitutions below are nonempty words.

\par\smallskip\noindent\textit{Operations, values, and normal forms.}\quad 
Both multiplication tables have rows 00334,\allowbreak{}11334,\allowbreak{}22334,\allowbreak{}33334,\allowbreak{}44334. Addition for 5413 has rows 00004,\allowbreak{}01014,\allowbreak{}00224,\allowbreak{}01234,\allowbreak{}44444; its order is \(3<1<0<4\) and \(3<2<0<4\), with 1 and 2 incomparable. Addition for 5414 has rows 00004,\allowbreak{}01114,\allowbreak{}01224,\allowbreak{}01234,\allowbreak{}44444; it is the chain \(3<2<1<0<4\). These two additive orders give different head guards.

Right 0/\allowbreak{}1 is the identity transformation. Right 2/\allowbreak{}3 resets to 3, and right 4 resets to 4. Thus a word whose tail has a special factor in \{2,\allowbreak{}3,\allowbreak{}4\} takes the value of its last reset (3 or 4). If its tail has no such factor, it takes its head value. This includes literals. This is a table induction with no bound on word length.

The laws \(xyzy=xzy\) and \(xyy=xy\) retain the head and the tail's order of last occurrences. Delete an earlier tail occurrence in A y B y D with the first law if B is nonempty and with the second otherwise. A is nonempty. The normal form is therefore h followed by distinct tail letters in their last-occurrence order; a tail copy of h is kept. In particular \(x\leq xx\) fails at \(x=2\) in both algebras.

\par\smallskip\noindent\textit{Exact original-word criteria.}\quad 
For nonempty P and target w, put \(C=c(w)\), \(h=head(w)\), U=\{u in P:c(u) is contained in C\}. For every x in C write \(S_{x}\) for the content after the last occurrence of x in w. The inequality \(w\leq \sum\) P holds exactly under these conditions:

1. H: some original u in U is headed by h.

2. O: for every x in C, an original u contains x and its last-x suffix content is contained in \(S_{x}\). A head-only occurrence or literal is allowed.

3. G: if h does not occur in the target tail, then for 5413 there is an original same-head u in U with no h in its tail; for 5414 there is an original u in U with no h in its tail, with arbitrary head. Such a 5414 witness may omit h entirely. This condition is separate from H and is not dispensable.

A literal target h requires the original literal h by G. Powers hh do not suffice. All three conditions are invariant under normalization.

Necessity of H is separated by \(h=0\), other C \(variables=1\), \(outsiders=3\). The target is 0; without H every upper value is 1 or 3. Necessity of O at x is separated by \(x=4\), \(S_{x}=1\), \(others=3\). The target is 4, while every upper value is 1 or 3. For G assign \(h=2\), other C \(variables=1\), \(outsiders=3\). The target is 2. For 5413, without a compatible same-head word all upper values are 1 or 3, whose sum is 1 or 3 and does not absorb 2. For 5414, failure of G means every contained upper word has a tail h; every upper word therefore has value 3. These are explicit separating assignments for all variable counts and original polynomial sizes.

For sufficiency, target value 3 is least. For value 4 choose its last 4-valued variable, after which all variables are 0/\allowbreak{}1; O gives an upper word of value 4. If the target is 0 or 1, it has no tail reset, so all target-tail values are 0/\allowbreak{}1. H has the same head and only target variables, and consequently the same value. For target value 2, the head is 2, all tail values are 0/\allowbreak{}1, and h is simple. G in 5413 gives an upper word of value 2. G in 5414 gives an upper value 2,\allowbreak{}1, or 0, each of which absorbs 2 in its chain. Hence the criteria are sufficient.

\par\smallskip\noindent\textit{Finite lists and soundness.}\quad 
The four laws common to both bases are:

{\small\par\smallskip
\CaseLine{\(xyzy = xzy;\qquad xyy = xy;\qquad xy \leq  x+zy;\qquad xy \leq  x+y.\)}
\smallskip}

For 5414 add the four identities obtained by independently omitting or retaining a and b in this template WHEN WRITING the identity:

{\small\par\smallskip
\CaseLine{\(x[a]y[b] \leq  xx[a]y[b] + y[a]y[b].\)}
\smallskip}

These are four fixed identities on at most four variables, never an empty-word substitution rule. The first two common laws follow from the table formula. For \(xy\leq x+zy\), literal x supplies H and the head guard; x and the true tail y of zy supply O. For \(xy\leq x+y\) the two literals supply everything. In each 5414 recovery law, xx[a]y[b] supplies H and all suffix witnesses. The second right word is contained and omits x, so supplies G for the simple target head x. Thus all eight laws are sound under the unrestricted exact criterion. Arbitrary identifications and term substitutions follow from equational logic.

\par\smallskip\noindent\textit{Unrestricted completeness.}\quad 
Normalize a valid nonlinear target \(w=hT\) and P. Maintain \(uS\leq P\) for EVERY original u, where S is the constructed part of T, built from right to left. Initially S is absent. For the next x, O supplies an original \(q=R\) x Q with c(Q) contained in c(S), so \(qS=R\) x S by the word laws. If R is nonempty, use \(xy\leq x+zy\) with \(X=u\), \(Y=x\), \(Z=R\) in right context S. If R is absent, use \(xy\leq x+y\) with \(X=u\), \(Y=x\) in that same context. In either case \(uxS\leq uS+qS\leq P\). An absent context is allowed as context, not as a substitution. There is no restriction on the number of iterations or on original summands.

For 5413 choose H's witness v to satisfy G as well whenever h is simple in w. If h is repeated, no compatibility is needed. In either case vT normalizes to exactly w, because c(v) is contained in C and the unwanted prefix-tail letters are deleted by later T occurrences. This proves the four-law theorem, including all head-only witnesses.

For 5414 choose any H witness v. Its vT is w or hhT. If it is w, we are done. Otherwise h is simple in w and hhT is available. Choose an original G witness z. If its head is h, zT is w. Otherwise its head y is in T; write \(T=A\) y B, where A and B may be absent. The available zT normalizes to yT, since its tail has no h and is contained in T. Select the corresponding recovery template with \(X=h\), \(Y=y\) and the nonempty A, B blocks substituted for the displayed optional letters. Its two upper words are hhT and yT; its lower word is \(hT=w\). Thus at most one recovery step is required. No common multiplicity is lost, no literal is replaced by a nonlinear word, and no empty substitution is used.

Literal targets are original by G. Expand arbitrary terms distributively and derive every summand's absorption into the opposite side. This proves all valid identities for unrestricted variable counts, lengths and polynomial sizes. Reversing words gives the opposite bases.

\subsection{Private initial letters and parity on a path}\label{case:157}
Types 5482 and 12637, and their multiplicative opposites 5467 and 14535, have no finite identity basis in the constant-free signature of additively idempotent semirings. Write u \ensuremath{\leq} P for \(u+P=P\). All polynomials below are finite nonempty sums of nonempty words.

\par\smallskip\noindent\textit{The two tables.}\quad 
In the order (e, p, g, a, z) the catalogue labels are (1,\allowbreak{}0,\allowbreak{}3,\allowbreak{}2,\allowbreak{}4) in 5482 and (2,\allowbreak{}1,\allowbreak{}3,\allowbreak{}0,\allowbreak{}4) in 12637. In both algebras z annihilates multiplication. A word over e, p, g equals e if it is all e, equals g if the number of occurrences of g is odd, and otherwise equals p. A word involving a and avoiding z equals a precisely when its only a is at the initial position and all other entries are e. In every other case it equals z. These statements follow by reading the common multiplication rows. In particular e is only a right identity on the whole semigroup. The additive order of 5482 consists of the chains e \ensuremath{<} p \ensuremath{<} z and g \ensuremath{<} a \ensuremath{<} z. The additive order of 12637 consists of e \ensuremath{<} p \ensuremath{<} a \ensuremath{<} z and g \ensuremath{<} z, with no further comparisons except transitivity. In each case p and g are incomparable and form a multiplicative group of order two, with identity p.

\par\smallskip\noindent\textit{Identities at every length.}\quad 
Let n \ensuremath{\geq} 1. Use distinct letters \(h_{0}\),..., \(h_{n}\) and \(x_{0}\),..., \(x_{n}\). The upper polynomial \(P_{n}\) has the two words \(h_{i}\) \(x_{i-1}\) and \(h_{i}\) \(x_{i}\) for each 1 \ensuremath{\leq} i \ensuremath{\leq} n. Its two remaining words are \(h_{0}\) \(x_{0}\) and \(h_{0}\) \(x_{0}\) \(x_{n}\) \(x_{n}\). The target is \(w_{n}\) = \(h_{0}\) \(x_{n}\). Write \(D_{i}\) for the two vertex letters of the ith pair, so D\_0=\{x\_0, \(x_{n}\)\}.

Both algebras satisfy \(w_{n}\) \ensuremath{\leq} \(P_{n}\). To see this, use each possible principal upper bound b of the values of all right words. The bound z is automatic. At \(b=e\) all their letters must be e. For \(b=p\) or \(b=g\) all letters are ordinary (e, p, g). At \(b=p\) every upper word has even g-parity, and at \(b=g\) every upper word has odd g-parity. Comparing the pair at \(h_{i}\) shows that \(x_{i-1}\) and \(x_{i}\) have equal g-bits. Hence the g-bit of the target equals that of \(h_{0}\) \(x_{0}\). This gives the required bound in each case.

At \(b=a\) no vertex can be a, since it occurs noninitially, and no letter can be z. If \(h_{i}\) has value a, every letter in \(D_{i}\) has value e. In particular if \(h_{0}=a\) the target has value a. Otherwise replace every head of value a by g in 5482, and by e in 12637. The formerly a-valued upper words now have the same required g-parity as the other upper words, odd in 5482 and even in 12637. Comparing all n edge pairs again gives the correct parity for the unchanged ordinary target. It has value at most a.

\par\smallskip\noindent\textit{Shape of every absorbed word.}\quad 
Suppose v \ensuremath{\leq} \(P_{n}\) holds. Giving every original letter e and an outside letter z excludes any outside letter of v. Give all \(h_{i}\) the value g and all vertices e. Every upper word equals g, forcing an odd number of occurrences of private letters \(h_{i}\) in v. For each i separately put \(h_{i}=a\), all vertices e, and all other heads g in 5482 or e in 12637. Every upper word is at most a. An occurrence of \(h_{i}\) anywhere after the first position would make v equal z. Consequently v has exactly one private letter, once at the initial position.

If this letter is \(h_{i}\), repeat its last test with any one vertex outside \(D_{i}\) given value p. The upper bound is still a. Other heads in 5482 give g even in the presence of p, and other heads in 12637 give e or p. Thus v cannot contain this vertex. Finally set all heads to p and all vertices to g. Every upper word evaluates to g. The tail of v has odd length. Its parity vector on the two vertices of \(D_{i}\) must therefore be a single endpoint. These are necessary restrictions on \emph{all} raw words, including words with repetitions; necessity suffices.

Call a word bad if it has head \(h_{0}\) and tail parity \(x_{n}\). The original polynomial contains no bad word and the target is bad. Any collection of nonbad absorbed words that omits a path head \(h_{j}\), 1 \ensuremath{\leq} j \ensuremath{\leq} n, cannot absorb any bad absorbed word. For a separating assignment use the group \{p, g\} as bits 0,\allowbreak{}1. Give the vertices with index below j bit 0 and the rest bit 1. Give \(h_{0}\) bit 0. At each represented path head the two endpoints have the same bit; give that head their bit. All the chosen nonbad words have value p, and every bad word has value g. This holds for every number of summands and every word length.

\par\smallskip\noindent\textit{Stability under polynomial substitutions and contexts.}\quad 
Assume a valid identity on at most k variables is applied in a derivation starting with \(P_{n}\), with n \ensuremath{>} k+2. Expand all intermediate terms as polynomials by the basic semiring axioms. All intermediate polynomials have the same value as \(P_{n}\) in the source algebra, so every one of their words is absorbed by \(P_{n}\) and has the shape above. Consider the first step that adds a bad word. In the expanded one-hole context choose its fixed left and right context words, either possibly absent. For a nonempty context word adjoin one fresh variable with that word as its singleton image. The old contribution from these fixed context choices is the substituted source polynomial Q. A single source word v on the new side gives the chosen bad word by expansion, and v \ensuremath{\leq} Q is valid on at most k+2 variables. Every expansion of Q occurs among the old nonbad words. The unaltered additive summands cannot have introduced the first bad word.

Let H be the variables that occur initially in at least one monomial of Q. Every word in the substitution image of each variable in H starts with a private letter: fixing the other choices in such a monomial gives an old shaped word. The *chosen expansion* of v has only one occurrence of any private letter in total. Thus across all variables of H it uses at most one image-word occurrence. Freeze the image of that variable, if present, to its chosen word. For each other variable in H choose any singleton image from its old image. Leave the other variables' polynomial images intact. This is a legitimate substitution into the same valid absorption. The bad word is still produced; every old expansion still has its original good shape. Their initial letters now take at most |H| \ensuremath{\leq} k+2 distinct values. At least one of the n path heads is absent. The preceding group assignment contradicts validity. No assumption of a consistent choice at repeated occurrences of the unfrozen variables was made. This also covers the other orientation of a basis identity and multiplication by polynomials, by distribution in the context.

The equality \(P_{n}=P_{n}+w_{n}\) would add a bad word. For any finite basis choose n larger than its variable bound plus two. The contradiction proves nonfinite basability. Word reversal transfers the result to the two stated opposite tables.

\par\smallskip\noindent\textit{Strict containment and transfer.}\quad 
Proposition~\ref{prop:lower-fibres} gives \(\A{5482}\in\V(\A{12637})\). The valid absorption hhhx \ensuremath{\leq} hx+x of 5482 separates the converse. If h, x are ordinary, hhhx and hx have the same value. A z input or a noninitial a gives an upper z, as does \(h=a\) with \(x=p\) or g. At \(h=a\), \(x=e\) the upper value in 5482 is \(a+e=z\). At \(h=a\), \(x=e\) in 12637 the left value is z and the right value a. Thus 12637 and its opposite have a proper nonfinitely based subvariety. The argument gives no limit classification for 5482 or 5467.

If another target algebra satisfies the displayed path identities at all lengths and its variety contains a source to which the substitution obstruction applies, any putative finite basis for the target would be valid in that source. Its derivation of the target's path equality would give the same forbidden first bad word in the source. Both hypotheses are required. For an opposite source use the reversed path words throughout.

\subsection{Penultimate witnesses and terminal overlap: thirteen-identity bases for a five-element semiring and its opposite}\label{case:158}
The source manuscript for this case credits Aifa Wang and Lili Wang, School of Mathematical Sciences, Chongqing University of Technology.

\par\smallskip\noindent\textit{Abstract.}\quad 
We give a thirteen-identity basis for a five-element additively idempotent semiring whose nonlinear word values retain internal content and the final two variables. The absorption criterion has two cases, determined by whether a terminal variable of the right polynomial also occurs internally. Without such overlap, a target requires a witness for its final pair and coverage of its internal content. With overlap, only content coverage and a terminal witness remain. Explicit derivations merge internal prefixes in the first case and construct saturated words in the second. These arguments allow arbitrarily many variables and summands. Reversal gives a basis for the multiplicative opposite.

\par\smallskip\noindent\textit{1. The algebra and a finite basis.}\quad 
The signature is \((+,\cdot)\) without constants. Bases are relative to the additively idempotent semiring axioms, and \(u\leq v\) means \(u+v\approx v\). Let \(A=\{d, e, b, a,\omega\}\), with additive order

\CaseDisplay{e<d<\omega,\qquad a<b<\omega,}

and no other nonreflexive comparisons. Addition is join. Put \(K=\{d, e\}\). Multiplication on \(K\) is right zero: \(uv=v\) for \(u, v\in K\). For \(u\in K\), set \(ua=a\), \(u\omega=\omega\), \(eb=a\), and \(db=\omega\). Every left factor outside \(K\) has constant row \(\omega\). Element names are not constants in terms. Explicitly:

\noindent \(\cdot\)\enspace\textbar\enspace \(d\)\enspace\textbar\enspace \(e\)\enspace\textbar\enspace \(b\)\enspace\textbar\enspace \(a\)\enspace\textbar\enspace \(\omega\)\par
\noindent \(d\)\enspace\textbar\enspace \(d\)\enspace\textbar\enspace \(e\)\enspace\textbar\enspace \(\omega\)\enspace\textbar\enspace \(a\)\enspace\textbar\enspace \(\omega\)\par
\noindent \(e\)\enspace\textbar\enspace \(d\)\enspace\textbar\enspace \(e\)\enspace\textbar\enspace \(a\)\enspace\textbar\enspace \(a\)\enspace\textbar\enspace \(\omega\)\par
\noindent \(b\)\enspace\textbar\enspace \(\omega\)\enspace\textbar\enspace \(\omega\)\enspace\textbar\enspace \(\omega\)\enspace\textbar\enspace \(\omega\)\enspace\textbar\enspace \(\omega\)\par
\noindent \(a\)\enspace\textbar\enspace \(\omega\)\enspace\textbar\enspace \(\omega\)\enspace\textbar\enspace \(\omega\)\enspace\textbar\enspace \(\omega\)\enspace\textbar\enspace \(\omega\)\par
\noindent \(\omega\)\enspace\textbar\enspace \(\omega\)\enspace\textbar\enspace \(\omega\)\enspace\textbar\enspace \(\omega\)\enspace\textbar\enspace \(\omega\)\enspace\textbar\enspace \(\omega\)\par
Associativity follows as follows. If either of the first two factors lies outside \(K\), both associations of a triple give \(\omega\). Otherwise \(uv=v\). The product \(vw\) lies in \(K\cup\{a,\omega\}\), and left multiplication by \(u\in K\) fixes each of these elements. Thus \((uv)w=vw=u(vw)\). Left translations by \(d, e\) fix \(K\), send \(a\) to \(a\), send \(b\) to \(a\) or \(\omega\), and fix \(\omega\). They preserve joins within either chain and send the join of elements from different chains to the join of their images, namely \(\omega\). Other left translations are constant \(\omega\). Right translations by \(d, e\) are constant on \(K\) and send its complement to \(\omega\); right translation by \(a\) sends \(K\) to \(a\) and its complement to \(\omega\); right translation by \(b\) sends only \(e\) to \(a\), all other elements to \(\omega\); and right translation by \(\omega\) is constant. Each preserves joins. Hence \(A\) is an additively idempotent semiring.

Theorem 1. The following thirteen identities form a basis for \(A\):

\CaseDisplay{\begin{array}{rll}
\mathrm{(M)}&xyzt\approx yxzt,\\
\mathrm{(N_1)}&x^2y\approx xy,\\
\mathrm{(N_2)}&xyx\approx yx^2,\\
\mathrm{(E)}&x\leq x^2,\\
\mathrm{(I)}&xy\leq zxy,\\
\mathrm{(S_1)}&xyzt\leq xs+yzt,\\
\mathrm{(S_0)}&xzt\leq xs+zt,\\
\mathrm{(D_{11})}&xy^2\leq xy+zyt,\\
\mathrm{(D_{10})}&xy^2\leq xy+yt,\\
\mathrm{(D_{01})}&y^2\leq y+zyt,\\
\mathrm{(D_{00})}&y^2\leq y+yt,\\
\mathrm{(F_1)}&x^2yz^2\leq x^2+yz,\\
\mathrm{(F_0)}&x^2z^2\leq x^2+z.
\end{array}\quad\textup{(\ensuremath{1})}}

Reversing each word gives a basis for \(A^{\mathrm{op}}\). No assertion of minimality or hereditary finite basability is made.

\par\smallskip\noindent\textit{2. Values and word normalization.}\quad 
For a word \(w\), write \(c(w)\) for its content and \(t(w)\) for its final variable. For a nonlinear word (length at least two), let \(i(w)\) be its internal content, the variables before its last position, and let \(p(w)\) be its penultimate variable. For a linear word, \(i(w)\) is empty.

Lemma 2. Under any assignment, a nonlinear word evaluates as follows:

1. If an internal variable is outside \(K\), its value is \(\omega\).

2. Otherwise its value is the assigned terminal element if that element is \(d, e\), or \(a\).

3. If its terminal element is \(b\), its value is \(a\) when its penultimate element is \(e\), and is \(\omega\) when that element is \(d\).

4. A terminal \(\omega\) gives \(\omega\).

Proof. Once a factor outside \(K\) has a following factor, the product is \(\omega\) and remains so. A nonempty prefix entirely in the right-zero band \(K\) has the value of its last factor. Applying the multiplication table to this value and the terminal element proves the rule. In particular, earlier internal occurrences assigned \(d\) do not affect case 3 if the penultimate element is \(e\). \(\square\)

Fix an ordering of variables. A nonlinear word is saturated if \(t(w)\in i(w)\). Define its normal form by

\CaseDisplay{\nu(w)=
\begin{cases}
\displaystyle\left(\prod_{x\in i(w)\setminus\{t(w)\}}x\right)t(w)^2,
   &t(w)\in i(w),\\[4pt]
\displaystyle\left(\prod_{x\in i(w)\setminus\{p(w)\}}x\right)p(w)t(w),
   &t(w)\notin i(w).
\end{cases}\quad\textup{(\ensuremath{2})}}

An empty product is an omitted substring, not a unit constant; put \(\nu(x)=x\) for a linear word.

Lemma 3. Identities \((M),(N_1),(N_2)\) derive \(w\approx\nu(w)\). A saturated word is determined by its content and terminal variable. An unsaturated nonlinear word is determined by its internal content and final pair.

Proof. Identity \((M)\) interchanges adjacent positions having at least two later positions. Thus the part preceding the last two positions can be permuted. For a word whose terminal variable is not internal, move any internal repetitions together and delete them by \((N_1)\). An earlier occurrence of the penultimate variable can be moved directly beside its penultimate occurrence and deleted by \((N_1)\) as well. The final pair is preserved, giving the second form in (2).

For a repeated terminal \(s\), take an earlier occurrence and a suffix \(sRs\). If \(R\) is nonempty, \((N_2)\) gives \(sRs\approx Rs^2\); if \(R\) is empty the repeated pair is already present. Permute the remaining prefix by \((M)\) and remove repetitions by \((N_1)\), including extra copies of \(s\) moved beside \(s^2\). This gives the first form. Lemma 2 also verifies all three normalization identities. \(\square\)

\par\smallskip\noindent\textit{3. The absorption criterion.}\quad 
Let \(P=\sum_{u\in\mathcal U}u\) be a nonempty polynomial. Define

\CaseDisplay{J=\bigcup_{u\in\mathcal U}i(u),\quad
T=\{t(u):u\in\mathcal U\},\quad U=J\cup T,\quad
L=\{x:x\text{ is a linear summand of }P\},}

and, for each \(s\in T\), put

\CaseDisplay{B_s=\{p(u):u\in\mathcal U,\ |u|\geq2,\ t(u)=s\}.}

The criterion below can be applied to normalized words. Indeed, when \(J\cap T\) is empty no right word is saturated, and any absorbed nonlinear target will also be unsaturated.

Lemma 4. The following are equivalent to \(A\models w\leq P\).

1. If \(J\cap T\ne\varnothing\), then \(c(w)\subseteq U\) and \(t(w)\in T\).

2. If \(J\cap T=\varnothing\) and \(w=x\) is linear, then \(x\in L\).

3. If \(J\cap T=\varnothing\) and \(w\) is nonlinear, then \(i(w)\subseteq J\) and \(p(w)\in B_{t(w)}\).

Here case 3 includes a terminal witness, since membership in \(B_{t(w)}\) implies that a right word ends in \(t(w)\).

Proof of necessity. In either regime, a missing content variable can be assigned \(\omega\) with all other variables assigned \(e\). The target is \(\omega\) and the polynomial is \(e\). If the target terminal \(s\) is absent from \(T\), assign \(s=d\) and all others \(e\). All right words are \(e\), whereas the target is \(d\). This proves the two requirements of case 1 and also the general terminal requirement.

Now assume \(J\cap T=\varnothing\). If some \(x\in i(w)\) is not in \(J\), assign \(J\) to \(e\), \(T\) to \(a\), and \(x\) to \(a\) as well, assigning any other variables to \(e\). This is consistent because \(J, T\) are disjoint and \(x\notin J\). Every right word is \(a\), but the target has an internal letter outside \(K\) and is \(\omega\).

For a linear target \(s\in T\setminus L\), assign \(J\) to \(e\), \(s\) to \(b\), and the other terminal variables to \(a\). Every right word ending in \(s\) is nonlinear, with penultimate element \(e\), so all right words are \(a\). The target is \(b\not\leq a\).

Finally consider a nonlinear target whose internal content is in \(J\) and whose terminal \(s\) is in \(T\). If \(p=p(w)\notin B_s\), assign \(p=d\), \(J\setminus\{p\}\) to \(e\), \(s=b\), and \(T\setminus\{s\}\) to \(a\). All right nonlinear words are \(a\): those ending in \(s\) have a penultimate variable different from \(p\). Linear right words are \(a\) or \(b\). Thus \(P\) is \(a\) or \(b\), while \(w=\omega\).

Proof of semantic sufficiency. If \(P=\omega\), the inequality holds. If \(P\in K\), every right word has value in \(K\), and consequently every variable in \(U\) is assigned an element of \(K\). Nonlinear word values are then their terminal values. The criteria give \(c(w)\subseteq U\) and a terminal of the right polynomial, so \(w\leq P\).

If \(P\in\{a, b\}\), every internal variable in \(J\) is assigned to \(K\) and every terminal variable in \(T\) is assigned to \(\{a, b\}\). Hence \(J\cap T\) is empty. For a linear target in \(L\) the assertion is immediate. For a nonlinear target, choose a right word with the same final pair. Both have all internal elements in \(K\), and Lemma 2 gives equal values. Therefore \(w\leq P\). This exhausts all values. \(\square\)

All inequalities in (1) satisfy this criterion. For \((I)\) the final pair is unchanged, and internal content increases on the right. For \((S_1),(S_0)\), the second right word supplies the target final pair, while both words together cover its internal content; if the right polynomial has overlap, its content and terminal also cover the target. Each \((D)\) has a terminal \(y\) which also occurs internally on the right, and its target has content within that right content and terminal \(y\). The right polynomial of each \((F)\) contains the saturated word \(x^2\), so has overlap; its content and terminals cover the target. For \((E)\) the right word \(x^2\) itself has overlap. This proves soundness.

\par\smallskip\noindent\textit{4. Constructive completeness.}\quad 
Lemma 5. Every absorption in Lemma 4 follows from (1).

Proof without terminal overlap. Suppose first that \(J\cap T\) is empty. A linear absorption is an existing summand. For a nonlinear target with final pair \(ps\), choose a right word \(q=Rps\) with that pair; \(R\) may be absent. For each nonlinear right word \(v=Xh\), where \(X\) is nonempty, \((S_1)\) gives

\CaseDisplay{XRps\leq Xh+Rps\leq P.}

When \(R\) is absent, use \((S_0)\) to obtain \(Xps\leq Xh+ps\leq P\). After every step normalize. The final pair stays \(ps\), while the internal content becomes the union of that of the current word and that of \(v\). After all right words have been used, the result has internal content \(J\), final pair \(ps\), and lies below \(P\).

If the target omits letters of \(J\), prepend a nonempty word listing these missing letters to the target. Identity \((I)\) shows that the target lies below this enlarged word; the substitution for \(x\) can be its first letter and for \(y\) the remaining nonempty suffix. Normalization identifies the enlarged word with the word already derived. If no letter is missing, normalization identifies the target directly. No terminal is introduced internally, since \(J\cap T\) is empty.

Proof with terminal overlap. Suppose \(J\cap T\ne\varnothing\). Choose \(x\in J\cap T\), a right word \(q\) ending in \(x\), and a right word \(v=ZxH\) with the displayed \(x\) internal. The suffix \(H\) is nonempty. If \(q\) is already saturated, keep it. Otherwise write \(q=Rx\), allowing \(R\) to be absent. The appropriate \((D)\) identity gives

\CaseDisplay{Rx^2\leq Rx+ZxH\leq P.}

Use \((D_{11}),(D_{10}),(D_{01}),(D_{00})\) for the four combinations of present or absent \(R, Z\). By Lemma 3 this yields a saturated word \(r\leq P\), and \(r^2\approx r\).

For any original right word \(v=Ys\) with \(Y\) nonempty, \((F_1)\) gives

\CaseDisplay{r^2Ys^2\leq r^2+Ys\leq P.}

Its normal form is saturated, has terminal \(s\), and has content \(c(r)\cup c(v)\). If \(v=s\) is linear, use \((F_0)\) instead. Use these operations for every right summand, and then once more for a right summand ending in the desired terminal \(t=t(w)\). The resulting saturated word has content \(U\), terminal \(t\), and lies below \(P\).

If \(w\) is linear, \((E)\) gives \(w\leq w^2\). If it is nonlinear, multiply \((E)\) for its terminal on the left by its nonempty prefix to give \(w\leq wt\). In either case the upper word is saturated with the same content and terminal as \(w\). Prepend its missing letters from \(U\) using \((I)\), or do nothing if none are missing. Lemma 3 identifies the result with the saturated word just derived. Thus \(w\leq P\) follows. All actual word substitutions are nonempty. \(\square\)

Proof of Theorem 1. Soundness was proved above. Distribute any valid polynomial identity \(P\approx Q\) into finite sums of words. Each summand of \(P\) satisfies Lemma 4 relative to \(Q\), so Lemma 5 derives \(P\leq Q\). Interchanging the polynomials derives \(Q\leq P\). This proves completeness without any bound on variables or summands. Reversing every word reverses multiplication and converts the whole argument into a proof for \(A^{\mathrm{op}}\). \(\square\)

\par\smallskip\noindent\textit{Acknowledgments.}\quad 
The source manuscript for this case acknowledges support from the National Natural Science Foundation of China (12371024), the Science and Technology Research Program of Chongqing Municipal Education Commission (KJZD-K2024011102) and the Chongqing Natural Science Foundation Innovation and Development Joint Fund (Municipal Education Commission) (CSTB2025NSCQ-LZX0067).

\subsection{Eight complete identities for semiring 9069}\label{case:159}
Signature (+, multiplication), without constants. Relative to the additively idempotent semiring axioms, a complete finite basis is

{\small\par\smallskip
\CaseLine{\(xy=yx; xyz=xyzz; xyzt\leq xyz; xy\leq x+xz+yt;\)}
\CaseLine{\(xxx\leq x+yyy; xxy\leq xy+zzz; xy\leq xxy; xy\leq xx+yz.\)}
\smallskip}

The order notation \(u\leq P\) means \(u+P=P\). The proof is unrestricted in variables, word length and the number of right summands.

\par\smallskip\noindent\textit{Values and original-witness criterion.}\quad 
Addition is 00004/\allowbreak{}01114/\allowbreak{}01224/\allowbreak{}01234/\allowbreak{}44444, the chain \(3<2<1<0<4\). Multiplication is 11111/\allowbreak{}11111/\allowbreak{}11111/\allowbreak{}11121/\allowbreak{}11114. A nonlinear word is 4 exactly when all its content is assigned 4. Otherwise a quadratic is 2 when its two occurrences are assigned 3; every other nonlinear word is 1.

For an original polynomial P let L be the literal variables, E the unordered quadratic edges, and V their endpoints. Let G mean that P contains an original word of length at least three. Then \(w\leq P\) iff:

\textbullet{} w is literal and occurs in P; or

\textbullet{} w is long, G holds, and an original u has c(u) contained in c(w); or

\textbullet{} w is quadratic, an original u has c(u) contained in c(w), and either G holds or c(w) is contained in V.

Necessity: a missing literal x is separated by \(x=0\), all others 3. Right values are at most 1, while the target is 0. If no contained right word exists, assign target content 4 and all outsiders 3; every right value is at most 1, while a nonlinear target is 4. If G fails, assign V to 3 and outsiders 2. Every right summand is at most 2. Long targets are 1, as are quadratics with an endpoint outside V. These cover every possible failed condition.

For sufficiency, G forces \(P\geq 1\), and a contained original word supplies 4 whenever the nonlinear target is 4. Otherwise its value is at most 1. Without G, all original quadratics have value at least 2. If P is below 1, every original edge must have value 2, forcing every V variable to 3; a quadratic on V is then 2. Its value 4 is again handled by the contained witness. Literal targets already occur.

Unlike 10379, an incident literal and quadratic do not supply G. The proposed \(xxx\leq x+xy\) fails at \(x=y=3\): the left value is 1, the right value is 2. Thus no proof may introduce a long word merely from such an overlap.

\par\smallskip\noindent\textit{Finite completeness.}\quad 
Commutativity and \(xyz=xyzz\) normalize long words by their content. Commute a variable to the end, duplicate it using two nonempty prefix blocks, then reverse the same identity to remove surplus copies while keeping length at least three. A canonical long word is sorted content for at least three variables, xxy for two, and xxx for one. In particular \(g^{3}=g\) for a long word g. These are literal equational derivations, not an assertion that arbitrary words are idempotent.

If G holds, fix an original long g and a contained original word u. For literal u derive \(u^{3}\) via \(xxx\leq x+yyy\) with \(y=g\). For quadratic \(u=ab\) derive aab via \(xxy\leq xy+zzz\) with \(z=g\). If u is long retain it. Append any missing target variables by \(xyzt\leq xyz\), splitting the derived long word into three nonempty blocks, then normalize. This derives any long target. A quadratic ab follows from aab by \(xy\leq xxy\).

Suppose G fails and a quadratic ab is not already original. A contained original word is either an endpoint literal or endpoint loop. For an original loop aa and an edge bc at the other endpoint, use \(xy\leq xx+yz\). For an original literal a, both a and b must occur in original quadratic edges az and bt. Use \(xy\leq x+xz+yt\). Commute the endpoint names when needed. The latter rule explicitly retains BOTH incidence witnesses; an isolated contained literal is insufficient. No derivation adds a new literal or creates a long word when G fails.

Every listed inequality satisfies the criterion and is separately checked over all assignments to its variables. All right summands used above are original or explicitly derived. All substitution words are nonempty. Distributivity and summandwise absorption in both directions now prove every identity of the source. Classification of further types requires independent full HSP equivalence, not a short-function coincidence.

\subsection{A finite basis for 5929 from short ideals and three original witnesses}\label{case:160}
Signature (+, multiplication), without constants. The finite family is effectively specified, not an enumerated small basis.

Addition rows 00004,\allowbreak{}01004,\allowbreak{}00224,\allowbreak{}00234,\allowbreak{}44444 give \(1<0<4\) and \(3<2<0<4\). Multiplication rows 01000,\allowbreak{}11111,\allowbreak{}01000,\allowbreak{}01020,\allowbreak{}41444 imply: any word with a1-valued factor is 1. Without 1, a literal keeps its value; a nonlinear head 4 gives 4; a quadratic with both factors 3 gives 2; every other nonlinear word gives 0. This follows directly by induction on length.

\par\smallskip\noindent\textit{Exact criterion.}\quad 
Let U be the original right words contained in the target content. Let Q be the content union of all original quadratic right words. Then \(w\leq P\) holds exactly when:

H: U contains a word with the same head as w; G: w is literal, or P contains an original word of length at least 3, or w is quadratic and its content is contained in Q.

The length and quadratic witnesses may have variables outside the target content. Literals need not occur originally. There is no literal/quadratic crossing alternative here;6022/6023's criterion cannot simply be copied.

The proper principal additive ideals are\{1\},\{3\},\{2,\allowbreak{}3\},\{0,\allowbreak{}1,\allowbreak{}2,\allowbreak{}3\}. For \(P\leq 0\), if w has no 1 and its head were 4, its contained H witness would be 4, contradiction. Thus \(w\leq 0\). If \(P=1\), every original contains 1, so the contained H witness forces 1 into w, and \(w=1\). If \(P=3\), all originals are literal 3; G forces a literal target, and H puts its variable among those originals. If \(P\leq 2\), originals can only be literal 2/\allowbreak{}3 or quadratic with both variables 3. No long original is possible. G forces a literal target or a quadratic whose variables all lie in Q. The latter has value 2; the former has value 2/\allowbreak{}3 because its contained H witness prevents other assignments. This proves sufficiency; \(P=4\) is automatic.

For necessity of H, assign target head 4, other target variables 0, outsiders 1. Without H each original is 0/\allowbreak{}1, but the target is 4. If G fails, there is no long original. Assign every variable in Q the value 3 and every other variable 2. Each original is literal 2/\allowbreak{}3 or quadratic 33 of value 2, so \(P\leq 2\). A long target is 0; or a quadratic target has a variable outside Q, so is not 33 and has value 0. Both violate absorption. This proves the criterion without any bound on variables, lengths, or the number of right summands.

\par\smallskip\noindent\textit{Normalization and compression.}\quad 
The six sound structural identities are:

{\small\par\smallskip
\CaseLine{\(xyz=xzy, xxy=xyy, xxyz=xyz,\)}
\CaseLine{\(xyyz=xyz, xyzz=xyz, xxxx=xxx.\)}
\smallskip}

Their soundness follows from the word formula: they preserve content, head, and length at least 3. Short words stay literal. A long normal word is its head h followed by the increasing list of other content variables, except that singleton content gives hhh and two-variable content gives hyy. Thus long normal forms preserve head/content and never shorten to a quadratic.

For a literal equational normalization, sort tail letters with \(xyz=xzy\), putting head copies first. When the length exceeds 3, delete adjacent duplicates with \(xxyz=xyz\) at the head, \(xyyz=xyz\) internally, or \(xyzz=xyz\) at the final pair with at least 2 prefix positions. At length 3 the case hhy becomes hyy by \(xxy=xyy\). All substitutions are nonempty.

Choose an original contained H witness. For a literal target it alone suffices. For a nonlinear target with an original long witness, adjoin that witness. Otherwise the target is quadratic; adjoin one original quadratic containing each of its at most 2 variables. At most 3 original words have been selected and still satisfy H/G, so they already absorb the target.

Normalize those at most 4 words. Reserve all variables of short words and the head of each long word, at most 8 variables. Group the rest by nonzero incidence pattern among the words, at most 15 groups. Identify each group to a fresh variable and normalize. Validity is preserved by substitution. The compressed absorption uses at most 23 variables, length 25, and at most 3 upper words.

For the inverse substitution, send reserved variables to their original singletons and group variables to the nonempty products of their members. Short words return exactly because all their variables were reserved. Every long expansion remains long and regains its original head/content; literal structural normalization recovers its original normal word. Even a two-variable compressed long word remains long on expansion. Thus this is an invertible equational compression.

Take the six structural identities plus all table-valid absorptions over a fixed 23-variable alphabet, word length at most 25, and at most 3 upper words. H/G effectively decides this finite family. Every valid absorption has the selected-original-word reduction, whose compression belongs to the family and whose inverse substitution recovers it. Adjoining the remaining original summands recovers the whole absorption. Distributive expansion and the two absorption directions derive every valid identity. Hence 5929 and its opposite are finitely based with unrestricted completeness. Catalogue transfer requires complete mutual HSP independently checked for each matching type.

\subsection{A three-witness finite basis from prefix coverage and a head witness}\label{case:161}
This note concerns source 5937 in the five-element catalogue. The signature is (+, multiplication), without named constants. The bases are relative to the additively idempotent semiring axioms.

\par\smallskip\noindent\textit{The algebra and its words.}\quad 
Let A=\{0,\allowbreak{}1,\allowbreak{}2,\allowbreak{}3,\allowbreak{}4\}, with rows and columns in that order:

{\small\par\smallskip
\noindent addition       multiplication\par
\CaseLine{\ttfamily 0 0 0 0 4\qquad 0 1 0 0 0}
\CaseLine{\ttfamily 0 1 0 0 4\qquad 1 1 1 1 1}
\CaseLine{\ttfamily 0 0 2 2 4\qquad 0 1 0 0 0}
\CaseLine{\ttfamily 0 0 2 3 4\qquad 0 1 2 3 0}
\CaseLine{\ttfamily 4 4 4 4 4\qquad 4 1 4 4 4}
\smallskip}

The additive order consists of \(3<2<0<4\) together with \(1<0\); the element 1 is incomparable with 2 and 3. The labels are element names, not constant symbols. Addition is the join in this order.

The multiplication is associative. For a triple containing 1 both associations are 1. Otherwise its first factor 4 gives 4, and its first factor 0 or 2 gives 0. For a first factor 3, a middle factor in \{0,\allowbreak{}2,\allowbreak{}4\} makes both associations 0. If the first two factors are 3, both associations equal 3 times the last factor, since left multiplication by 3 fixes its image \{0,\allowbreak{}2,\allowbreak{}3\}. Each multiplication row and column preserves the displayed chain and the join of 1 with 2 or 3, namely 0; the remaining joins are between comparable elements. Thus multiplication distributes over addition and A is an additively idempotent semiring.

For a word w of length at least two, direct induction gives:

\textbullet{} if any variable has value 1, then \(w=1\);

\textbullet{} if no variable has value 1 and its head has value 4, then \(w=4\);

\textbullet{} otherwise \(w=3\) if all positions have value 3;

\textbullet{} otherwise \(w=2\) if all nonfinal positions have value 3 and its final position has value 2;

\textbullet{} otherwise \(w=0\).

Consequently the following identities are valid:

\CaseDisplay{xyzt\approx xzyt,\qquad xy^2z\approx xyz.       \quad\textup{(\ensuremath{W})}}

Fix a variable order. Using the first law with nonempty prefix and suffix substitutions permits adjacent interior positions to be interchanged. The second deletes repeated interior occurrences. Every nonlinear word therefore has the normal form

\CaseDisplay{h\,\prod_{x\in M}x\,t,                 \quad\textup{(\ensuremath{1})}}

where h is its head, t its last variable, and M its interior set. The endpoints are retained, even if equal or present in M. Empty M means the two-letter word ht, not a unit constant. Linear words are unchanged.

Write c(w) for the content, h(w) for the head, and i(w) for the set of variables before the final position. For a linear word i(w) is empty. Normalization preserves all three data.

\par\smallskip\noindent\textit{Exact absorption criterion.}\quad 
For a nonempty polynomial P, set

\CaseDisplay{D=\bigcup_{v\in P}c(v),\qquad
                J=\bigcup_{v\in P}i(v).}

Lemma 1. For every nonempty word w, linear or nonlinear, A satisfies \(w\leq P\) if and only if:

1. c(w) is contained in D;

2. i(w) is contained in J;

3. some summand v of P has \(h(v)=h(w)\) and c(v) contained in c(w).

Here \ensuremath{\leq} abbreviates additive absorption.

Proof. If condition 1 fails, assign a missing variable 0 and all others 3. The target is 0 and every right summand is 3. If condition 2 fails, assign a variable in i(w) outside J the value 2 and all others 3. The target is 0. Each right summand is 2 or 3, since this variable can only occur finally on the right. Neither evaluation satisfies the proposed absorption.

If condition 3 fails, assign h(w) the value 4, the other variables in c(w) the value 3, and all outside variables 1. The target is 4. Every summand containing an outside variable is 1. Every remaining summand has a different head and value in \{0,\allowbreak{}2,\allowbreak{}3\}; hence \(P\leq 0\). This also proves necessity for a linear target.

Conversely fix an assignment. When \(P=4\) there is nothing to check. When \(P=0\), the target can fail only by taking value 4. In that case its head is 4 and no variable of c(w) is 1. The witness v in condition 3 also has value 4, a contradiction. When \(P=1\), every summand is 1 and therefore contains a variable assigned 1. Applying this to v puts a variable of c(w) at 1, so \(w=1\). When \(P=2\), every right summand is 2 or 3. The word formula then puts every variable of D in \{2,\allowbreak{}3\} and every variable of J at 3. Conditions 1 and 2 force w to be 2 or 3, including when w is linear. Finally \(P=3\) forces every variable of D to be 3, so condition 1 gives \(w=3\). These cases exhaust the algebra. \ensuremath{\square}

\par\smallskip\noindent\textit{Decomposition into three-witness steps.}\quad 
Suppose \(w\leq P\) is valid and normalize all words. If w occurs in P, no step is necessary. Choose a witness v from condition 3.

If w is linear, say \(w=y\), then v is a positive power of y. The absorption \(y\leq v\) itself satisfies Lemma 1 and is a step with one right summand.

Now let \(w=h\) (product of M) t be nonlinear. Put \(M_{0}=c(v)\) minus \{h, t\} and \(s_{0}=h\) (product of \(M_{0}\)) t. Thus \(M_{0}\) is contained in M. Choose an original summand \(r_{t}\) containing t, which exists by condition 1. Every variable in i(\(s_{0}\)) except possibly the last variable of v already belongs to i(v). Indeed i(\(s_{0}\)) is contained in c(v), because \(h=h(v)\). Therefore

\CaseDisplay{|i(s_0)\setminus i(v)|\leq1.              \quad\textup{(\ensuremath{2})}}

If this difference is nonempty, let its variable be x and choose a summand \(q_{x}\) of P with x in i(\(q_{x}\)). Such a summand exists since i(\(s_{0}\)) is contained in i(w), which is covered by J. The case of linear v is included: the possible missing variable is then h.

The absorption

\CaseDisplay{s_0\leq v+r_t+q_x                      \quad\textup{(\ensuremath{3})}}

is valid, omitting \(q_{x}\) when (2) is empty. Its content is covered by v and \(r_{t}\); its prefix content by v and \(q_{x}\); and v retains the required head and content. There are at most three right summands.

For each missing interior variable x in M minus \(M_{0}\), choose an original summand \(q_{x}\) with x in i(\(q_{x}\)). If s' is obtained from the current normal word s by adjoining x to its interior set, then

\CaseDisplay{s'\leq s+q_x                       \quad\textup{(\ensuremath{4})}}

is valid. The summand \(q_{x}\) supplies the new content and prefix membership, and s is the same-head content witness. This includes adjoining an endpoint variable to the interior. After finitely many steps the target is w. Every right summand in (3)--(4) is either an original summand of P or the preceding constructed word.

\par\smallskip\noindent\textit{Reversible compression and completeness.}\quad 
Consider a valid absorption with one normal left word and at most three normal right words. Reserve each variable that is a first or last variable of any of these at most four words. There are at most eight reserved variables. Every remaining variable occurs only in interiors. Partition the remaining variables according to their nonempty incidence pattern among these interiors. There are at most \(2^{4}-1=15\) such patterns. Identify each class to one fresh variable, leave reserved variables separate, and normalize using (W). Ordinary substitution preserves validity. The result has at most 23 variables and every normal word has length at most 25.

The compression is reversible modulo (W). Substitute for each class variable the nonempty product of its original members. Every member had the same interior incidence pattern, so the expansion restores exactly the original interior set in each word. Reserved endpoint variables are replaced by their individual original variables. The laws (W) normalize the expanded words to their original normal forms. Linear words, consisting of one reserved variable, are also restored. No empty words or constant symbols have been introduced.

Fix a 23-letter alphabet and let Sigma consist of (W) and all valid absorptions over this alphabet with one normal left word, at most three normal right words, and word lengths at most 25. This is a finite, explicitly defined set; Lemma 1 decides membership. The large family is not enumerated here, and minimality is not claimed.

Theorem. Sigma is a finite identity basis for A. Reversing all words gives a finite identity basis for the opposite semiring.

Proof. The basis is sound by definition and the verification of (W). Reversible compression shows that every valid absorption with at most three normal right words follows from Sigma by substitution and (W). Equations (3)--(4) now derive every valid word absorption into a polynomial: each step uses only original right summands or previously derived words, all already below that polynomial. For a valid polynomial identity, every summand on either side is absorbed by the other side, so these derived absorptions establish the identity. Distributivity reduces arbitrary terms to polynomials. The proof places no bound on the original variables or word lengths. Reversal transports the identities and all derivations to the opposite semiring. \ensuremath{\square}

Finite table checks and regressions support error detection only. Transfers to other catalogue types require independently verified equality of generated varieties, with both HSP directions checked. No hereditary finite-basis assertion follows from this theorem.

\subsection{A finite identity basis for 6063}\label{case:162}
The additive order has \(3<2<0<4\) and \(1<0\), with 1 incomparable with 2 and 3. Its proper principal ideals are \{1\}, \{3\}, \{2,\allowbreak{}3\} and \{0,\allowbreak{}1,\allowbreak{}2,\allowbreak{}3\}. The multiplication rows are 01004,\allowbreak{}11111,\allowbreak{}01004,\allowbreak{}01024,\allowbreak{}44444. The first value in \{1,\allowbreak{}4\} occurring in a word is its value if any occurs. Otherwise all letters are in B=\{0,\allowbreak{}2,\allowbreak{}3\}. A literal retains its value, a quadratic has value 2 on (3,\allowbreak{}3) and 0 elsewhere, and a word of length at least three has value 0. The product table gives these statements by induction, using the two left zeros and B cubed equal to \{0\}.

Let w be a target word and P a nonempty polynomial. Write L for all literal variables of P, Q for the variables in all quadratic words of P, and J for all its words of length at least three. Then w \ensuremath{\leq} P holds if and only if all the following requirements hold.

1. Some original upper word is contained in the target content.

2. For each variable x in w a right word contains x with the variables before its first x contained in the variables before the first target x. There is no restriction on its suffix.

3. Either J is nonempty, or the target is a literal in L union Q, or it is quadratic with both letters in Q.

For necessity of the first requirement assign 0 on w and 1 outside. If it fails, the right sum is 1 and the target is 0. For the second assign 4 at the missing x, 0 on its preceding target letters and 1 elsewhere. Without its witness the right side misses the greatest value. For the third, with J empty, assign 3 on Q, 2 on L outside Q and 0 elsewhere. All upper values belong to \{2,\allowbreak{}3\}. This assignment excludes any target of length at least three, any quadratic escaping Q, and any literal outside L union Q.

For sufficiency, the prefix witnesses handle target value 4, hence the ideal below 0. A target of value 1 already lies in \{1\}, and a neutral target cannot have all upper words equal to 1 because of its contained witness. If all right words lie in \{2,\allowbreak{}3\}, their letters are neutral and no right word is long. The variables in L have values in \{2,\allowbreak{}3\}, those in Q all have value 3, and the prefix requirement ensures target variables occur on the right. A literal or a covered quadratic as in requirement 3 now belongs to \{2,\allowbreak{}3\}. Finally every right word of value 3 is literal. If all right values are 3 then Q and J are empty and the target must itself be an original literal. These are all proper principal ideals. This proves both implications for any alphabet and any original right polynomial.

Take the seven long-word normalizing identities N1--N7 of Section~\ref{case:99}. They apply since long values are determined by the first-occurrence list. Add the following six identities and the finite family of all valid short absorptions with at most twelve variables, word length at most two and at most five upper words. In each line G abbreviates the single word ijk on fresh variables.

{\small\par\smallskip
\CaseLine{\(xxxu \leq  u + xs + G\qquad xxxu \leq  u + x + G\)}
\CaseLine{\(pxxxu \leq  pu + pxs + G\qquad pxxxu \leq  pu + px + G\)}
\CaseLine{\(a \leq  aaa + G\qquad ab \leq  ababab + G\)}
\smallskip}

The contained witnesses are respectively u, u, pu, pu, aaa and ababab. The displayed upper words supply every prefix witness, and G provides the long word. Thus the six laws are sound.

When J is empty keep one contained witness, one prefix witness per target letter and, for a quadratic target, the two original quadratics covering its letters. A literal target also retains one literal or quadratic witnessing its membership in L union Q. These are at most five original upper words of length at most two and the target has length at most two. The criterion still holds, so the bound of twelve variables gives an instance of the short family.

When J is nonempty split one original long word into three nonempty blocks for G. The iterative proof of the displayed six identities uses the same induction as the six schemes of the cited 3177 proof, now with this single fixed G. For each current prefix p of cubed target letters derive p u \ensuremath{\leq} P for every original u. The first prefix uses an original starting with its variable and the first two laws. To append x use its original witness \(v=r\) x s where all letters of r already occur in p. The long words p v and p x s have the same list of first occurrences. The next two laws therefore give p xxx u \ensuremath{\leq} P from the previous stage, with the appropriate case for an absent suffix. The chosen G need not be contained in the target. On the last step use the separate contained witness as u. Normalize to w for a long target or to w cubed for a short target, followed by the last two laws in the latter case. This derives every valid absorption. Distributive expansion and the two directions of summand absorption give unrestricted equational completeness.

\subsection{A nonfinite basis obstruction for 6203 and 12635}\label{case:163}
Use the constant-free signature of additively idempotent semirings. The exact absorption criterion on arbitrary alphabets for both algebras is proved from their tables in Section~\ref{case:164}. We recall the parts used here. The support of a word absorbed by a polynomial P contains an original summand. Each first occurrence of a letter a has an original witness whose prefix before its first a uses only previously seen letters D. Let \(R_{a}\) be the union of the tails of originals headed by a. At least one of the following holds.

1. Some a-headed original repeats a.

2. Some original has a noninitial first a preceded only by D union \(R_{a}\).

3. In 6203 some original is supported in D union \(R_{a}\).

4. The target has a simple initial a and its whole tail is supported in \(R_{a}\).

The stated conditions together are also sufficient. The criterion applies to source-variable absorptions in an arbitrary putative axiom, as well as to literal words in a derivation.

\par\smallskip\noindent\textit{The identities.}\quad 
Take n \ensuremath{\geq} 2 and distinct h, x, t, z, f, \(d_{1}\),..., \(d_{n}\), \(r_{1}\),..., \(r_{n}\). Write \(D_{0}\) for the alphabet \{h, \(d_{1}\),..., \(d_{n}\)\}. Form four original words and the target

{\small\par\smallskip
\CaseLine{\(U = h h d_{1} \cdots  d_{n} z\)}
\CaseLine{\(X = x t r_{1} \cdots  r_{n}\)}
\CaseLine{\(V = r_{1} d_{1} \cdots  r_{n} d_{n} x\)}
\CaseLine{\(F = f f\)}
\CaseLine{\(W = h d_{1} \cdots  d_{n} x f.\)}
\smallskip}

Let \(P=U+X+V+F\). Both semirings satisfy W \ensuremath{\leq} P. The word F supplies containment. The head h uses U and its repeated head. For each \(d_{i}\) the noninitial witness U has its prefix supported in the preceding target letters. The word X supplies an ordinary prefix witness at x. The word V supplies its exceptional witness since all \(d_{i}\) have been seen and the \(r_{i}\) occur in the tail of X. Finally F witnesses f. This holds at all n. The guards t, z, f and the single word U impose restrictions absent from the earlier semantic examples.

\par\smallskip\noindent\textit{Positions before the first x in an absorbed word.}\quad 
Call a word critical if its first x has a nonempty prefix over \(D_{0}\). None of U, X, V, F is critical and W is. All absorbed words use only variables of P by the ordinary prefix requirement.

\textbullet{} Every x-headed absorbed word has a unique x and all later letters in \{t, \(r_{1}\),..., \(r_{n}\)\}. At its initial x the first three alternatives fail. In particular V needs \(d_{1}\), and no original lies wholly in R\_x=\{t, \(r_{1}\),..., \(r_{n}\)\}.

\textbullet{} Every first t in an absorbed word has x strictly before it. Its only ordinary prefix witness is X.

\textbullet{} Suppose the first x is noninitial and its entire preceding prefix avoids t, z, f. Its exceptional witness must be V. The extra contained-word alternative fails, since D union \(R_{x}\) then misses x, z, f, and every original contains one of these. Thus all \(d_{i}\) occur before this x.

\textbullet{} In that third case, if the word is noncritical, some \(r_{i}\) occurs before x. The ordinary prefix test forces \(r_{1}\) before x, since X has initial x and V begins with \(r_{1}\). At the first \(r_{1}\) its prior letters avoid t, z, f. The set \(R_{r_1}\) consists only of \(d_{i}\), the other \(r_{i}\) and x. The only noninitial witness for \(r_{1}\) is X, requiring t. No original can be supported in the prior set union \(R_{r_1}\), since X needs t, U needs z, F needs f and V needs \(r_{1}\) itself. The repeated-head alternative also fails. Thus \(r_{1}\) is the simple head and its entire tail lies in \(R_{r_1}\). In particular h is absent. Before the first x, a first \(d_{i}\) cannot use U, which needs h, and a first \(r_{i}\) cannot use X. Each must therefore follow V's first-occurrence ordering. Since \(d_{n}\) occurs, the sequence of distinct letters encountered before x is exactly \(r_{1}\), \(d_{1}\),..., \(r_{n}\), \(d_{n}\). Raw words are allowed to repeat old letters.

\par\smallskip\noindent\textit{The first critical rewrite.}\quad 
Let a proposed basis use at most k variables per axiom. Choose 2n \ensuremath{>} k+2. Distribute terms into nonempty finite sets of nonempty words. A one-hole term context consists after distribution of unchanged summands and contributions L sigma(side) R with fixed exterior words, either optionally absent. This follows by structural induction on the context. For the contribution producing a chosen new word, absorb each nonempty exterior word into a fresh source variable with a singleton word image. The contextual source identity is still valid and uses at most \(K=k+2\) variables. Old expansions from this contribution belong to the old polynomial; the new source monomial is absorbed by its old source side.

Starting from P every intermediate summand is absorbed by P, since the intermediate polynomial equals P. Suppose the first critical word is produced by a source monomial v. Fix the independent nonempty word choices in this expansion. Let a be the variable occurrence whose image contains the first x. Its selected image begins with a possibly empty block over \(D_{0}\) followed by x. Write A for the variables strictly before the *first source occurrence* of a in v. Every variable of A has a \(D_{0}\)-only image word: it occurred before the selected x in the chosen expansion. Fix one such choice per variable of A. The selected occurrence of a need not be its first source occurrence.

Use the ordinary source prefix condition at a in v. In a witnessing old monomial choose the fixed \(D_{0}\)-only images before its first a and the selected image of a there. If that image had a nonempty prefix before x, or if the source witness had nonempty prefix before a, the resulting old expansion would be critical. Thus the image begins with x, and there is an old source monomial headed by a. Every old expansion of any such monomial with that selected head image starts with x. By the first position restriction above, every image word of every source variable in its tail is entirely over \{t, \(r_{1}\),..., \(r_{n}\)\}. This applies to the full set \(R_{a}\) of the old source side, regardless of other choices. In particular a is not in \(R_{a}\), and A is disjoint from \(R_{a}\).

Apply the exceptional source condition at a. A repeated head a is impossible. The simple-head alternative is also impossible: the occurrence selected to produce x is noninitial in v since its image starts with x and the new word was critical. If the source has a as its initial variable, it consequently repeats a. In 6203 an old monomial supported in A union \(R_{a}\) is impossible. Choose its \(D_{0}\)-only images on A and any of the restricted images on \(R_{a}\). Its expansion contains none of x, z, f, so its support contains no original word of P. This contradicts containment for the old expansion.

The remaining alternative supplies an old source monomial with a *nonempty* prefix before its first a, supported in A union \(R_{a}\). Use the fixed \(D_{0}\)-only images on A and fix one image word per variable of \(R_{a}\), identically at every occurrence. Choose the x-initial image at a. No t can occur before x by the position restriction on all old expansions. The resulting prefix is nonempty, avoids t, z, f and is not critical. It must therefore introduce the 2n distinct letters \(r_{1}\), \(d_{1}\),..., \(r_{n}\), \(d_{n}\) in exactly this alternating order.

Each contributing source prefix variable now has one fixed image, either over \(D_{0}\) or entirely over the \(r_{i}\). In the latter assertion we only need variables actually occurring in this prefix; any t there has already been excluded. A repeat of a source variable introduces no new letter. A first use of a source variable introduces at most one new letter, because its single image uses only one of the two disjoint alphabets and consecutive new letters must alternate. The prefix therefore introduces at most K distinct letters. This contradicts \(2n>K\).

Thus no bounded-variable derivation from P can introduce a critical word. The argument applies in either orientation and to arbitrary polynomial substitutions with independent choices, including on repeated source variables. The semiring axioms preserve the free distributive expansions. The valid identities \(P=P+W\) at all n therefore preclude any finite identity basis. Word reversal gives the same result for 6209 and 13203.

\par\smallskip\noindent\textit{Strict containment and transfer.}\quad 
Every identity of 12635 is valid in 6203 by the two exact absorption tests, so Var(6203) is contained in Var(12635). For strictness use xx \ensuremath{\leq} x+xy+y. When x differs from the exceptional S, \(xx=x\). At \(x=S\) the right word xy gives the top T unless \(y=E\). For 6203, \(S+E=T\); for 12635, \(S+E=S\). Thus 12635 and its opposite contain proper nonfinitely based subvarieties and are nonlimit. This particular inclusion leaves the limit status of 6203 and its opposite open.

For any other algebra the proof transfers if the obstruction algebra lies in its variety and all lengths of the above identities hold. All lengths in a finite target can be checked by closing the triples (product of the \(d_{i}\), product of the \(r_{i}\), product of \(r_{1}\) \(d_{1}\) ... \(r_{n}\) \(d_{n}\)) under appending one pair of values, starting at length two. For every reachable triple test all five independent values of h, x, t, z, f in the four original words and W. A failing assignment suffices to reject a route; a surviving candidate also needs the variety inclusion.

\subsection{Absorption tests for 6203, 12635 and their opposites}\label{case:164}
These are unrestricted semantic tests. A finite equational derivation or a nonfinite-basis obstruction remains to be established. Reversal gives the tests for 6209 and 13203, respectively.

\par\smallskip\noindent\textit{A common product and two additive orders.}\quad 
Use the five names E, N, L, S, T. In 6203 they denote 2,\allowbreak{}0,\allowbreak{}1,\allowbreak{}3,\allowbreak{}4. In 12635 they denote 2,\allowbreak{}1,\allowbreak{}3,\allowbreak{}0,\allowbreak{}4. The two labelled multiplications then agree. On \{E, N\} the product is maximum with \(E<N\). Scan for the first entry in \{L, S, T\}. If it is L or T, that entry is the answer. An S gives S only if it is the actual head and every later entry is E. In all other cases it gives T. These statements follow by appending an entry, including the transition from a surviving S.

Both additive orders have E and L as minimal elements, N their join and T their greatest element. In 6203 the proper ideal below S is \{L, S\}. In 12635 the principal ideal of S is \{E, N, L, S\}. The other proper principal ideals in both tables are \{E\}, \{L\}, and \{E, N, L\}. Here and below principal ideals include their generator.

\par\smallskip\noindent\textit{The exact test.}\quad 
Let C be the support of the target word w. Require an original upper word contained in C. At the first occurrence of every variable x let D be the set of letters strictly before it. Require some original word containing x whose prefix before its first x is supported in D. The prefix may be empty.

Let R be the union of all tail letters of original words headed by x (all positions after their initial x), and put \(K=D\) union R. In addition to the preceding requirements, at least one of the following must hold at each first occurrence.

1. The letter x belongs to R.

2. Some original has a noninitial first x whose preceding letters are all in K.

3. In 6203, some original is supported wholly in K.

4. The variable x is the simple head of w and all target tail letters are in R.

The third alternative is specific to 6203. A literal x has an empty tail in this union. Reversing the test for each opposite must reverse both the target and every original word.

For necessity, assign C value E and every outsider L to force a contained original. For the prefix witness put D at E, x at T and every other letter at L. An upper word reaching T must have the required prefix and contain x. If the four alternatives all fail, x is outside K. Assign K value E, x value S and everything else L. An x-headed original then gives S. An original with first x noninitial and no earlier L would satisfy the second alternative. Every remaining upper gives L, or E in 12635. In 6203 an all-E original would satisfy the third alternative. Thus the upper sum lies below S. The target gives T: its first x is internal, or x repeats after the head, or a letter outside R follows its simple initial x. These are explicit separating assignments.

For sufficiency, a neutral target (all entries E or N) uses the contained original to escape the ideal of L, and in 6203 also the ideal of S. A missing ideal-of-E comparison, for a neutral target of value N, is prevented by coverage: an original containing a non-E variable cannot evaluate to E. The same coverage handles the ideal of E whenever the target has a nonneutral letter. A target of value L requires no other proper-ideal check. When the first special letter is T, the prefix witness reaches T. When it is S, that witness reaches S or T, disposing of the ideal of N and of L. For the remaining ideal of S suppose all originals are below S. Every x-headed original then forces its entire tail to E. All of D is neutral and all of R is E. The first alternative contradicts \(x=S\). The second forces T at a noninitial position. The third supplies a neutral upper value outside the ideal of S in 6203. The fourth permits just S as target value. This proves the test for arbitrary alphabets and polynomials.

The word equations \(xyy=xy\) and \(xyx=xxy\) hold in both tables. They reduce words to their first-occurrence order, doubling the head if it reappeared anywhere. To collect a new copy of the head use the second equation at the first position. To delete a repeated nonhead y in a tail, use \(xyy=xy\) for adjacent copies, and turn y u y into yy u by the second equation, then use the first in its nonempty left context. Extra copies of the head contract by \(xxx=xx\). Preservation of the normal form also follows from the word formula.

\par\smallskip\noindent\textit{A family with interleaved prefix and tail variables.}\quad 
Take distinct h, x, z, \(d_{1}\),..., \(d_{m}\) and \(r_{1}\),..., \(r_{m}\) with m at least one. Use target h \(d_{1}\) ... \(d_{m}\) x. The original words are

{\small\par\smallskip
\CaseLine{\(x\)}
\CaseLine{\(x r_{1} \cdots  r_{m}\)}
\CaseLine{\(r_{1} d_{1} \cdots  r_{m} d_{m} x\)}
\noindent h h z, and \(d_{j}\) \(d_{j}\) z for every j.\par
\smallskip}

Each h or \(d_{j}\) has a head witness with a repeated head. The literal x supplies containment. At x the noninitial witness alternates between R and D. Its prefix need not be supported in either separately. Thus the absorption holds for every m in both sources. Deleting the alternating original makes it fail: give x value S, all of D and R value E and z value L. Deleting the upper x \(r_{1}\) ... \(r_{m}\) also makes it fail: additionally give \(r_{1}\) value L. The deleted-tail case works at \(m=1\) as well. This family supplies a specific test for a proposed finite deduction or obstruction. Essentiality of the two original witnesses alone does not establish nonfinite basability.

All original head-x tails can be combined while retaining their order within each word. The table-valid inequality xuv \ensuremath{\leq} xu+xv derives an x-headed accumulated word with the entire union R as its tail support. The remaining proof problem is to use a witness whose prefix is interleaved between D and R, retaining the order needed by other variables. The head exception is real. Blindly squaring each original head may raise a surviving S to T.

\subsection{A 17-identity finite basis for 6211}\label{case:165}
Signature(+, multiplication), without constants. Upper words are always original words of the given polynomial unless explicitly called derived. Every substitution is nonempty.

\par\smallskip\noindent\textit{Operations and normal words.}\quad 
Addition rows:00044,\allowbreak{}01034,\allowbreak{}00244,\allowbreak{}43434,\allowbreak{}44444. Multiplication rows:01044,\allowbreak{}11144,\allowbreak{}01234,\allowbreak{}44444,\allowbreak{}44444. The addition table is \(equivalently1<0\),\(1<3\),\(2<0\), with top 4. The correct literal row for 3 is 43434: its entries are 4,\allowbreak{}3,\allowbreak{}4,\allowbreak{}3,\allowbreak{}4.

A word containing 4 has value 4. If it contains 3, its value is 3 exactly when 3 occurs once, finally, and every earlier factor is 2; otherwise its value is 4. With no 3 or 4, the value is 1 if a1 occurs, otherwise0 if a0 occurs, and otherwise2. Table induction proves this formula for arbitrary nonempty words.

The laws \(xxy=xy\), \(xyz=yxz\), \(xyy=yxx\) permute the prefix and contract its repetitions. A simple-final word has normal form a sorted list of its other variables followed by its final variable. A repeated-final word is called saturated: its value depends only on its whole content. Its normal form is the sorted content list with the largest letter doubled. To move the repeated final letter t to another final letter q, permute the prefix to end in q and apply \(qtt=tqq\). Repeat internal contractions if necessary. This gives a literal deduction using the three laws; it never deletes the saturation bit.

\par\smallskip\noindent\textit{Exact criterion.}\quad 
For an original polynomial P, call x active when either: I. some original word contains x before its final position; or II. writing \(A_{x}\) for the union of prefixes of original words ending in x, some original word has content contained in \(A_{x}\). For II to apply without I, all x-ending words have simple x.

An absorption \(w\leq \sum\) P holds exactly when:

1. c(w) is contained in the union of original contents;

2. some original U has c(U) contained in c(w);

3. every target variable is active, except possibly a simple final t; when t is inactive, the target prefix content is contained in \(A_{t}\).

Necessity has explicit separating assignments. If 1 fails, assign a missing target letter 4 and all others 2. If 2 fails, assign target letters 2 and outside letters 1. If 3 fails at an inactive x, assign \(x=3\), \(A_{x}=2\), and all other variables 1. All original x-ending words then have value 3, and every other word has value 1 because it is not contained in \(A_{x}\); their join is 3. The target is 4 either because x occurs internally or because its simple-final prefix has a variable outside \(A_{x}\). These assignments prove necessity without any bound.

For sufficiency, an active x assigned 3 forces the upper sum to 4. This is immediate from an original internal occurrence. Otherwise if any \(A_{x}\) letter is not 2, the corresponding original x-tail word is 4; if all \(A_{x}\) letters are 2, a contained original word has value 2 and an original x-tail word has value 3, whose join is 4. A target of value 4 is now handled by a covered 4 letter, an active internal 3 letter, or an inactive simple-final 3 whose non2 prefix factor lies in \(A_{x}\) and thus in an original tail witness. A target of value 3 has a covered 3 letter; an original word containing it has value 3 or 4. A target 1 has a covered 1 letter, and an original word containing it has value 1 or 4. A target 2 has a contained word with value 2. Finally target 0 uses a contained word of value 0 or 2 and an original word containing a target 0 letter. The latter has value 0,\allowbreak{}1 or 4; joining it with the contained value gives an upper bound for 0. These cases exhaust all values and establish unrestricted sufficiency.

\par\smallskip\noindent\textit{Seventeen laws.}\quad 
In addition to the semiring axioms take seven fixed laws:

{\small\par\smallskip
\CaseLine{\(xxy=xy\qquad xyz=yxz\qquad xyy=yxx\)}
\CaseLine{\(x\leq xx\qquad xyz\leq xz+yz\)}
\CaseLine{\(xy\leq x+zxy\qquad xy\leq y+zxy\)}
\smallskip}

Add ten identities given by these optional-factor templates:

{\small\par\smallskip
\CaseLine{\((z[a]t)^{2} \leq  z[a]t+z\qquad (2 choices)\)}
\CaseLine{\(([a]t)^{2} \leq  [a]t+[v]tb\qquad (4 choices)\)}
\CaseLine{\(xzz \leq  zz+[v]xb\qquad (2 choices)\)}
\CaseLine{\(z[a]t \leq  zt+vz[a]t\qquad (2 choices)\)}
\smallskip}

Square means literal concatenation twice. Omitted bracket factors define separate identities; substitutions never use empty words. Each law uses at most five variables. All 3105 table assignments across the 17 laws establish their soundness. Completeness is below.

\par\smallskip\noindent\textit{Unrestricted completeness.}\quad 
Normalize all words. First construct a derived internal witness for any active x. In case I use the original word itself. In case II without I, merge all nonempty prefixes of original x-tail words via \(xyz\leq xz+yz\). This produces a simple-tail word with full prefix \(A_{x}\). Choose an original z contained in \(A_{x}\) and let a list \(A_{x}\) minus c(z). The activation template \((z[a]x)^{2}\leq z[a]x+z\) derives a saturated word containing x. Its right terms normalize to the merged word and z. Its last repeated letter may change during normalization, but every variable, including x, still occurs internally in a saturated word.

Suppose the target final t is active. Condition 3 makes all target letters active. Start with the original contained U. If it is saturated use it directly. Otherwise write \(U=A\) t0 and split an internal witness for t0 as V t0 B, with B nonempty. The second template saturates U without adding content. For every missing target letter x, split its original or derived internal witness V x B and apply \(xzz\leq zz+[v]xb\) with \(z=the\) accumulating saturated word. This adds x and no outside letters. Finite iteration produces a saturated word of exactly target content. If the desired target tail is simple, apply \(x\leq xx\) in its prefix context, recovering the target.

Suppose instead that t is inactive. The target has simple t and its prefix content lies in \(A_{t}\). For a one-letter target condition 2 gives an original pure t power; inactivity excludes a repeated t, so t itself is already present. Otherwise merge original nonempty t-tail prefixes with \(xyz\leq xz+yz\). Obtain \(A_{t}\) t, possibly with outside letters. The contained original U must end in simple t: an internal t would activate it, and a word omitting t would be contained in \(A_{t}\), also activating it. If \(U=t\), use \(xy\leq y+zxy\) with \(x=the\) target prefix, \(y=t\), \(z=the\) outside prefix. If \(U=Zt\) with Z nonempty, use \(z[a]t\leq zt+vz[a]t\), setting \(z=Z\), \(a=the\) missing target prefix letters, \(v=the\) outside prefix. In either case the large right term normalizes to the merged word and the left to the target. When no outside letters remain there is no erasure step. These applications retain simple t and use no empty substitutions.

Thus every valid word absorption follows from 17 laws for arbitrary numbers of variables, original summands and lengths. Splitting an arbitrary polynomial identity into its word absorptions in both directions proves that the 17 laws form a finite basis for 6211. No classification of another algebra follows without full transfer.

\subsection{A finite basis for 6562 and its opposite 7592}\label{case:166}
In the constant-free signature use the usual additively idempotent semiring axioms. The additive order of 6562 is 3 \ensuremath{<} 2 \ensuremath{<} 1 \ensuremath{<} 0 \ensuremath{<} 4. The multiplication rows are 01114,\allowbreak{}11111,\allowbreak{}23331,\allowbreak{}33333,\allowbreak{}44444. Values 1,\allowbreak{}3,\allowbreak{}4 on the left are fixed by every further factor. Starting at 0 or 2, a string of zeros preserves that value. From 0 the first nonzero tail factor gives 4 if it is 4, otherwise 1. From 2 it gives 1 if it is 4, otherwise 3. These outputs are fixed. This word formula also holds when the tail is empty.

\par\smallskip\noindent\textit{Exact tests.}\quad 
Write the target as h t, where h is a letter and t may be empty. Let H be the set of heads of all original words of the nonempty polynomial P, and let C be the target support. Use these four tests.

1. Some original word is contained in C.

2. Each first occurrence of a variable x in the whole target has an original witness containing x whose first-x prefix uses only variables before that occurrence. In particular h belongs to H.

3. In the part of t strictly before the first letter in H (all of t if H is absent), every first x has a witness starting with h. In the tail of this witness, the prefix before its first x must use only variables before x in the *target tail*.

4. If t has no letter in H, some original word h v has all of v supported in t. This includes the literal case t empty.

For necessity of 1 assign C the value 0 and outsiders 3. For 2 give the earlier variables 0, the variable x value 4 and the rest 3; only the indicated prefix witness reaches 4. For 3 give the earlier tail variables 0, h value 2, x value 4 and the rest 3. These assignments are compatible since h is already in H. Every other original head has value 3. The target is 1 and an h-headed word can reach 1 only with the specified tail witness. For 4 assign \(h=2\), all target tail variables 0 and other variables 3. The target is 2. The only upper words that can reach it are h followed by a string of these zeros.

Conversely, a target valued 4 needs only the whole-word prefix witness at its first 4; a target valued 0 is all-zero and uses test 1. A target valued 1 whose head has value 0 or 1 is absorbed by any original word with that head. With head valued 2 and a first nonzero tail factor of value 4, consider the vertices up to and including this factor in the tail. If one belongs to H, an original beginning there has value at least 1 (it begins with 0 or 4). Otherwise test 3 supplies a word valued 1. When the head is 2 and the tail is all zeros, either a tail letter is an original head of value 0, or test 4 supplies a word valued 2. The value 3 is least. This proves the criterion.

\par\smallskip\noindent\textit{Eighteen identities.}\quad 
Two word identities are \(xyy=xy\) and \(xyzy=xyz\). They delete all but the first occurrence of each letter in the tail, retaining the initial position separately. They follow from the word formula; adjacent repeats use the first law and separated repeats the second, with nonempty exterior prefix up to the first repeated tail letter.

Use the following sixteen absorptions as equations (u\ensuremath{\leq}P means u+P=P).

\textbullet{} Three instances of a x [y] \ensuremath{\leq} a [y] + a x [t], independently omitting y or t but retaining at least one. With both omitted the conclusion is already the second upper term.

\textbullet{} Four instances of a b y \ensuremath{\leq} a [t] + b [z] + y, independently omitting t and z.

\textbullet{} One instance x y x z \ensuremath{\leq} x y z.

\textbullet{} Eight instances of h [p] b [r] \ensuremath{\leq} h [p] b h [r] + b [z], independently omitting p, r, z at every occurrence of each symbol.

The brackets specify separate written identities. All unomitted symbols denote distinct variables of the template and their substitutions are nonempty. The criterion proves the first two families directly. For the third law, inserting x after y can lower the value only when the initial x is 2 and every letter up to that insertion after it is 0; the inserted x then sends the product to the bottom 3. For the last family the old and new head-h words agree except possibly when \(h=2\) and every tail letter through b equals 0. There \(b=0\) and the other upper word b[z] has value at least 1, whereas the target beginning with 2 is at most 1. This gives soundness even when the expression following b is empty. Equivalently, check all sixteen templates directly on the displayed finite table.

\par\smallskip\noindent\textit{Complete construction.}\quad 
Work with the normalized target tail. For its initial segment p avoiding H maintain h p v \ensuremath{\leq} P for all original h v, including empty v for a literal original h. Initially this is true by definition. For the next x before encountering H use test 3 with original h d x t', where the support of d lies in p. The invariant and the word laws give h p x t' \ensuremath{\leq} P. For each maintained v use the first family at \(a=h\) p, \(y=v\), \(t=t'\) to derive h p x v. If v and t' are both empty the second upper term already gives it. This proves the invariant at every length. If H is never encountered, use the contained same-head witness from test 4 at the end. All of its tail is already present in the constructed tail, so its appended copy can be deleted. This also handles a literal target.

Otherwise the target is h p b r with the first b in H in its tail. Choose an original b z and any original h t'. From the previous invariant h p t' is available. The second family at \(a=h\) p supplies h p b u \ensuremath{\leq} P for every original u. The third law now inserts h after b, giving h p b h u \ensuremath{\leq} P for every original u. The added occurrence of h ensures that all previously seen variables appear in the tail of the maintained prefix.

Suppose this maintained prefix is l and x is the next letter of r. If x already appears in its tail the word identities give l x \(u=l\) u. Otherwise test 2 gives an original d x t' with d supported in the previously processed letters. From l d x t' \ensuremath{\leq} P and the fact that all of d is already in the tail, obtain l x t' \ensuremath{\leq} P. Apply the first family at \(a=l\), \(y=u\), \(t=t'\) to extend every maintained word. At the end take the contained original from test 1. Its appended copy deletes in full, leaving h p b h r \ensuremath{\leq} P. The fourth family with the *original* b z removes the added h and gives h p b r. All arguments work when b itself equals h. Omitted blocks were handled by separate laws and never by substituting an empty word.

Every valid word absorption therefore follows. By distributing terms and proving the summand absorptions in both directions this proves all valid identities, for unrestricted variables and polynomials. The inequalities and word identities reverse under multiplicative duality to give a complete basis for 7592.

\subsection{A finite basis for 6563 by three original witnesses and interior compression}\label{case:167}
Signature (+, multiplication), without constants. This specifies a finite family, not an enumerated small identity list.

\par\smallskip\noindent\textit{Operations and exact word values.}\quad 
Addition rows are 00004,\allowbreak{}01114,\allowbreak{}01224,\allowbreak{}01234,\allowbreak{}44444: the chain \(3<2<1<0<4\). Multiplication rows are 01114,\allowbreak{}11111,\allowbreak{}23332,\allowbreak{}33333,\allowbreak{}01114. For a word \(w=hT\), the table gives the following formula by induction. If the head value is 1 or 3, that value is retained. If the head is 2, the value is 2 when every tail factor is 0/\allowbreak{}4, and 3 otherwise. For head 0/\allowbreak{}4, a tail factor in \{1,\allowbreak{}2,\allowbreak{}3\} gives value 1; if every tail factor is 0/\allowbreak{}4, the value is the last factor. This includes literals.

The following structural identities follow from the formula:

{\small\par\smallskip
\CaseLine{\(xyzt=xzyt,\qquad xyzy=xzy,\qquad xyy=xy.\qquad (W)\)}
\smallskip}

The second and third laws delete earlier tail duplicates, preserving the head, while the first swaps adjacent interior positions with a nonempty prefix and suffix. Thus every nonlinear word normalizes to its original head, an increasing list of distinct interior variables other than the final variable, and its original final variable. A head copy in the interior is retained. A literal stays literal; xx stays xx. Every step is an actual word-law substitution with nonempty blocks.

\par\smallskip\noindent\textit{Exact criterion on original upper words.}\quad 
Write \(C=c(w)\), \(h=head(w)\), \(t=last(w)\), and T=c(w after its initial position). An absorption \(w\leq \sum\) P is valid exactly when:

L: an original word u has last variable t and c(u) contained in C; H: an original word v has head h (outside variables are allowed); G: if h is not in T, an original word z either has head in T, or has head h and all of its tail content contained in T.

There is no content restriction on the first alternative of G. If h is repeated in w, G is absent. For a literal target G requires that literal itself. These conditions are preserved by (W).

For necessity of L assign \(t=4\), other C \(variables=0\), \(outsiders=1\). The target is 4. Every upper word involving an outsider has value 1, and a contained upper word with a different tail has value 0. Without L their join is at most 0. For H assign \(h=1\) and all other \(variables=3\). The target is 1; without H all upper heads are 3, hence all upper values are 3. For G (with simple h) assign \(h=2\), \(T=0\), all remaining \(variables=3\). The target is 2. An original head in T or a compatible same-head word would satisfy G. In their absence every upper word is either headed by 3 or headed by 2 with a tail factor 2/\allowbreak{}3, and consequently every upper value is 3.

For sufficiency, value 3 is least. Target value 4 or 0 means all target factors lie in \{0,\allowbreak{}4\}; L's contained same-tail word gives exactly that value. Target value 1 is absorbed by H's witness: a head 1 always gives 1, and a head 0/\allowbreak{}4 always yields 0,\allowbreak{}1, or 4, each above 1. Finally target value 2 has head 2 and all target-tail values 0/\allowbreak{}4, so h is simple. G's same-head alternative gives 2, while an original head in T gives 0,\allowbreak{}1, or 4, all above 2. This proves exactness for arbitrary assignments and variables.

\par\smallskip\noindent\textit{Reduction to three original upper words.}\quad 
Choose u for L, v for H, and, only when needed, z for G. Remove duplicate choices. The same exact criterion proves

{\small\par\smallskip
\CaseLine{\(w \leq  u+v+[z],\qquad (R)\)}
\smallskip}

with at most three original upper words. By inclusion of these summands, (R) implies \(w\leq \sum\) P. No new or synthetically assumed witness is used. The following compression proves every such (R) is derivable from one fixed finite family.

\par\smallskip\noindent\textit{Reversible compression.}\quad 
There are at most four words: the target and at most three upper words. Reserve every variable appearing as a head or final variable of any of these words. There are at most eight reserved variables. Normalize all words by (W). A nonreserved variable occurs only in the interior, at most once in each word. Partition the nonreserved variables by their nonempty occurrence pattern among these at most four words. There are at most \(2^{4}-1=15\) such groups.

Replace each group by a fresh variable and each reserved variable by its own distinct fresh variable. A group appears in exactly the words where all its members appeared. This produces at most 23 variables. Each compressed word has its retained head and final position plus at most 23 interior variables, hence length at most 25. Literals stay literals, and every reserved repetition remains present.

L is preserved because heads and tails are reserved and subset relations are preserved by the occurrence partition. H is preserved because equality of heads is reserved. G is preserved because head variables are reserved, tail occurrences of those reserved variables are retained, and tail-subset relations are preserved. Therefore the compressed absorption is valid under the already proved criterion.

Map each reserved compressed variable back to its original singleton, and each group variable to the nonempty increasing product of all its members. This substitution recovers every original word modulo (W). Indeed the expanded word has the same head, tail, and interior set, with the same reserved repetitions. Normalize both expansions and originals by the actual swap/deletion laws; they give exactly the same word. This proves reversibility, not merely a homomorphic image or a one-way implication. There are no empty substitutions and no multiplicity discrepancies; all nonreserved interior multiplicities were legitimately reduced to one by (W) before grouping.

\par\smallskip\noindent\textit{The fixed finite basis and unrestricted completeness.}\quad 
Take (W), together with ALL valid absorptions over a fixed alphabet of 23 variables, of word length at most 25 on both sides, and at most three nonempty upper words. Validity is specified effectively by L, H, G above (or by the finite operation table). This is a finite, explicit family specification; no claim is made that it is a small enumerated list.

Every valid absorption reduces to (R), whose compression belongs to this family. Substitute the nonempty inverse words, use (W) in reverse where needed, and include its at most three original summands into P. This derives the original absorption for any number of variables, arbitrary lengths, and any nonempty finite polynomial P. Expanding terms by distributivity and deriving each summand's absorption in the opposite side proves every identity. Hence 6563 is finitely based. Reversing words proves its opposite.

\subsection{A finite basis for 6564 from full coverage and original tail heads}\label{case:168}
Signature (+, multiplication), no constants, relative to the additively idempotent semiring axioms. The finite family below is specified, not enumerated as a small list. Completeness has no bound on variables, word lengths, or the number of summands.

Addition: 00004,\allowbreak{}01114,\allowbreak{}01224,\allowbreak{}01234,\allowbreak{}44444; order \(3<2<1<0<4\). Multiplication: 01114,\allowbreak{}11111,\allowbreak{}23332,\allowbreak{}33333,\allowbreak{}41114. For a nonlinear word hT, heads 1 and 3 retain value. Head 2 gives 2 when all tail values are in \{0,\allowbreak{}4\}, and 3 otherwise. Head 4 gives 4 when all tail values are in \{0,\allowbreak{}4\}, and 1 otherwise. Head 0 gives 1 if any tail value is in \{1,\allowbreak{}2,\allowbreak{}3\}; otherwise it gives 4 if some tail is 4, and 0 if all are 0. Literals retain their values. Induction on the multiplication table proves these formulas, including every length. Thus \(xyz=xzy\) and \(xyy=xy\). A normal word consists of its head followed by its sorted tail set, with the head retained in the tail when it repeats. Adjacent tail swaps and duplicate deletions use nonempty substituted prefixes only.

\par\smallskip\noindent\textit{Exact criterion from original right words.}\quad 
A literal x is absorbed by P iff P contains the literal x. For nonlinear \(w=hT\), let \(C=\{h\}\) union T, let H be the set of original right-word heads, and U the original right words whose entire contents are contained in C. Then \(w\leq P\) holds exactly when all three conditions hold: (H) h belongs to H; (C) the union of the contents of U is all of C; (R) h belongs to T, or literal h is in P, or H intersects T, or an original nonlinear h-headed word has its tail contained in T. In H intersects T the original witnessing word may contain outside variables.

Necessity: for a missing literal x set \(x=2\) and all other variables to 3. Every right nonlinear x-word has value 3, and other heads give 3, whereas the target is 2. If H fails, set \(h=1\) and all other variables to 3: the target is 1 and right words have value 3. If C fails at an uncovered x, set \(x=4\), other target variables to 0, and outsiders to 3. The target is 4; contained right words avoid x and have value 0. Noncontained right words have value 1 or 3. The right sum is therefore at most 0, strictly below 4. If R fails after H and C, h is not in T. Set \(h=2\), \(T=4\), all others to 3. The target is 2. Literal h is absent, and every nonlinear h-headed right word has a tail outside T, so gives 3. No right head lies in T. All other right heads are 3, so the whole right side is 3. These arguments apply to arbitrary polynomials and arbitrary variable sets.

Sufficiency: target value 3 is least. If its value is 1, its head value is 0, 1, or 4; an original same-head word from H always has value at least

1. If the target value is 4, all its variables are 0 or 4 and some x is 4. C supplies a contained original word containing x, also valued 4. If the target is 0, all its variables are 0 and any contained original word gives

0. The remaining case, target value 2, has head 2 and tail values 0 or 4, so h cannot repeat in T. R supplies either the literal h valued 2, an original word headed in T valued 0, 1, or 4, or a nonlinear same-head word with tail in T valued 2. This proves the criterion, not merely a test bound.

\par\smallskip\noindent\textit{Derive coverage before recovering the head.}\quad 
Normalize all words. For a nonlinear target choose a contained original g in U. Until its content is C, choose an original u in U containing a missing target variable and derive \(normal(gu)\leq g+u\). This local absorption is valid by the criterion: g supplies H, g and u supply full coverage, and g supplies R, either as a literal head or as a nonlinear same-head word whose tail remains contained in the target tail. Each step strictly increases content and introduces no outside variable, so the process terminates for every finite C. Derive normal(g squared)\ensuremath{\leq}g; the new word has repeated head, and its tail is exactly C.

Choose an original H witness a, allowing outside variables. Derive \(hC\leq a+g\): a supplies H, g supplies C, and the target head is repeated. If h belongs to T, this is the desired normal target. Otherwise select one original R witness b: literal h, a word headed in T (possibly with outside variables), or a nonlinear h-word with tail in T. Derive \(hT\leq hC+b\). The word hC supplies H and full coverage, while b supplies R. Every local absorption uses at most two already derived or original right words. Literal targets already occur. The next step shows these local absorptions are all instances of a single finite family.

\par\smallskip\noindent\textit{Finite family and unrestricted completeness.}\quad 
For a local absorption, consider its at most three normal words: the left and at most two right words. Reserve their heads, at most three variables. Group every remaining variable by its nonempty incidence vector across the at most three tails. There are at most seven such groups. Identify each group to a fresh variable, preserving each reserved head, then normalize. The resulting absorption is valid by variable substitution; it uses at most ten variables and words of length at most eleven. Substitute for each fresh variable the nonempty product of its original group members. Identical tail incidence and reserved heads recover all original tail sets, including repetitions of a head in its tail. Apply the two word laws to obtain each original normal word. No empty word or multiplicative identity is used at any point.

Let Sigma consist of \(xyz=xzy\), \(xyy=xy\) and every valid normal absorption on a fixed ten-variable alphabet, length at most eleven, with at most two right words. This is a finite family specified by the exact H, C, R criterion. Compression and its inverse derive each local absorption from Sigma. The coverage construction therefore derives every valid absorption. Distributivity normalizes arbitrary terms to polynomials, and summandwise absorption in both directions derives every valid identity. Sigma is a finite relative basis for 6564, without restrictions on original variables, lengths, or summands. Opposite reversal is valid; catalogue transfers require independently checked complete HSP in both directions.

A single contained seed is insufficient: \(xy\leq x\) fails at \(x=0\), \(y=4\). Requiring a tail-contained recovery word would be too strong: for target xy the right words yx and xz provide full coverage and an original y-head, even though neither is contained in \{y\}. These distinctions are retained in the proof and the independent replay. Finite computations only supplement the unrestricted proof above.

\subsection{A finite basis for 6668 from nonlinear singleton witnesses}\label{case:169}
Signature (+, multiplication), without constants. This specifies a finite family, not an enumerated small basis.

Addition rows 00004,\allowbreak{}01114,\allowbreak{}01224,\allowbreak{}01234,\allowbreak{}44444 give \(3<2<1<0<4\). Multiplication rows are 01130,\allowbreak{}13331,\allowbreak{}13331,\allowbreak{}33333,\allowbreak{}41134. A literal keeps its value. In a nonlinear word, any factor 3 or at least two occurrences valued in\{1,\allowbreak{}2\} give 3. Exactly one such occurrence gives 1, irrespective of its position or whether it is valued 1 or 2. With all factors in\{0,\allowbreak{}4\}, the value is the head. This follows by induction on length from the multiplication rows. Occurrences are counted with multiplicity, not as distinct variables.

\par\smallskip\noindent\textit{Exact absorption criterion.}\quad 
A literal w is below P exactly when it is an original summand. For nonlinear w write \(C=c(w)\), \(h=head(w)\), and U for original summands contained in C. Then \(w\leq P\) iff:

H: some u in U has head h; S: for each variable x occurring just once in w, some u in U omits x, or is NONLINEAR and contains x exactly once.

A literal x is not an S witness containing x. All witnesses are ORIGINAL. For necessity of H put \(h=4\), other C variables 0, outsiders 3. Without H the upper sum is at most 0 whereas the target is 4. For failed S put \(x=2\), other C variables 0, outsiders 3. The nonlinear target is 1. Every upper word is 3 or the literal x of value 2, so the sum is at most 2. A missing literal h is separated by \(h=2\), others 3. For sufficiency, target 3 is least. Targets0/\allowbreak{}4 are matched by H. Target 1 requires precisely one\{1,\allowbreak{}2\} occurrence, belonging to a simple variable x. An omitting contained witness is 0/\allowbreak{}4; a nonlinear single-x witness is 1. Either absorbs the target. Thus H/S is exact for arbitrary alphabets, word lengths and numbers of upper summands.

\par\smallskip\noindent\textit{Six structural laws.}\quad 
{\small\par\smallskip
\CaseLine{\(xyz = xzy\)}
\CaseLine{\(xyx = xxy\)}
\CaseLine{\(xxx = xx\)}
\CaseLine{\(xy = xxy+xyy\)}
\CaseLine{\(xx \leq  x\)}
\CaseLine{\(xxyy \leq  xx\)}
\smallskip}

The value formula proves all six under arbitrary nonempty substitutions. For the splitting equality, if there is exactly one exceptional occurrence then doubling the block not containing it preserves1 and doubling its block gives 3. With all factors 0/\allowbreak{}4 both branches preserve the head. A zero value 3 remains 3 in both. Squaring reduces an exceptional-containing word to 3; otherwise it retains its head. Appending a square to a square either retains the head or gives 3, proving the last inequality.

Normalize by retaining the head h, its total multiplicity capped at 2, then the increasing remaining variables each with multiplicity capped at 2. The first law permutes adjacent tail letters with nonempty prefix. Move head copies next to the head; use \(xxx=xx\) in context to cap every run at 2. These are literal equational steps. The second structural law is also valid and may be retained redundantly. Normal words keep every variable's common multiplicity 1 or 2, not just its content.

\par\smallskip\noindent\textit{Guard and recovery of simple variables.}\quad 
For a nonlinear valid \(w\leq P\) choose an original H witness v. Use \(xx\leq x\) at each simple occurrence to double it. The resulting g has only repeated variables and lies below P. It satisfies g \(squared=g\) by the normal form. For each missing variable y of C, use \(xxyy\leq xx\) with \(x=g\) to append y squared, then normalize. This yields a full-content guard with head h and all multiplicities 2, still below P.

For each simple variable x of w choose an original S witness u. Double every other simple variable in u using \(xx\leq x\), and normalize to \(u_{x}\). It remains below u and either omits x or is nonlinear with exactly one x. The nonlinear condition is not lost by this operation. Let \(w_{x}\) be w with all variables other than x doubled, then normalized. The local absorption \(w_{x}\leq g+u_{x}\) is valid by H/S: g supplies H and the original S witness supplies the unique simple variable.

Any normal word with at least two simple variables splits as XY between two simple occurrences. \(xy=xxy+xyy\) gives branches with unchanged head and content and fewer simple variables. Normalize each branch and repeat. The leaves have no simple variable and equal g, or one simple variable x and equal \(w_{x}\). Each leaf lies below P, so the splitting equalities derive the initial target. Literal targets use their original summands and require no guard.

\par\smallskip\noindent\textit{Fixed finite recoveries and unrestricted completeness.}\quad 
In a local recovery there are at most three words. Reserve x and all their heads, at most 4 variables. Every other variable is absent or repeated in each word and never occupies the head. Group such variables by their nonzero incidence patterns among the three words, at most 7 groups. Identify each group to one fresh variable and normalize. This is a substitution into a valid identity, so the compressed recovery is valid. It has at most 11 variables, length at most 22, and two upper words.

For the inverse substitution, replace each group variable by a nonempty product listing its members. It had two occurrences in every compressed word where it appeared, so each member is restored twice. The first structural law regroups these tail occurrences; \(xxx=xx\) caps multiplicity. Reserved singleton images preserve the head and the unique simple x. Normalization therefore recovers all three original normal words. In particular, the grouping never replaces common multiplicity 2 by 1.

Take the six structural laws plus all valid absorptions over a fixed 11-variable alphabet, length at most 22, and at most two upper words. This is a finite effectively specified family, decidable by H/S. Every local recovery is derived by the compression and inverse substitution. The guard and finite splitting tree then derive every valid absorption with no bound on the original variables, lengths, or summands. Expansion by distributivity and both absorption directions give every polynomial identity. Consequently 6668 and its opposite are finitely based; catalogue identification still requires independent complete mutual HSP checks.

\subsection{Finite identity bases by splitting middle factors and taking meets of powers}\label{case:170}
Let A be a finite additively idempotent semiring of size q in the signature without constants. Assume its additive order has a least element. Binary meets exist because the join of the common lower bounds is their greatest element. Fix nonnegative k, l, a positive integer e and an integer s at least two. Write \(U=u_1\cdots u_k\) and \(V=v_1\cdots v_l\), using distinct fresh variables for the fixed end positions. When an end has zero positions omit it. Write \(F(t)=UtV\). Assume the following conditions for all values of the displayed variables:

{\small\par\smallskip
\CaseLine{\(F(xy) = F(yx)\)}
\CaseLine{\(F(x^{e+1}) = F(x^{e})\)}
\CaseDisplay{F(x_1\cdots x_s)=\sum_{j=1}^{s}F(x_1\cdots x_{j-1}\,x_j^e\,x_{j+1}\cdots x_s)}
\CaseDisplay{F(x^e y^e)=F((x\wedge y)^e)}
\smallskip}

The first three are identities in the semiring signature. In the fourth formula \(x\wedge y\) denotes the greatest lower bound in the given algebra. Only the first three and a finite family of valid absorptions will be used as axioms.

Put \(d=k+l\), \(h=d+s-1\), \(B=q^{h+2}\) and \(D=h(B+1)+2^{B+1}-1\). Take all valid absorptions on a fixed D-variable alphabet with at most B right summands and word length at most d+eD throughout. Together with the first three displayed identities and the semiring axioms this is a finite basis.

Normalizing and splitting. If the length is at most d keep the word as it is. Otherwise retain its first k and last l positions. Everything between them is its middle. Within the middle any adjacent swap is derived from the first identity, using the k immediately preceding positions and l immediately following ones as guards and leaving any unused outside positions as context. In the same manner the second identity caps each middle multiplicity at e once like letters have been sorted. Call a middle variable deficient when its capped multiplicity lies between 1 and e-1.

If there are at least s deficient variables, partition the entire sorted middle into s nonempty consecutive factors so each factor contains one distinct deficient variable. One may cut immediately before the second through the sth chosen variable's blocks. Substitute these factors for \(x_1,\ldots,x_s\) in the third identity. In each summand the factor selected for powering gives at least one previously deficient variable middle multiplicity e. Capping and sorting preserve every previously saturated variable. Each resulting word has strictly fewer deficient variables. This induction terminates and gives an equational expansion of any original polynomial into a finite sum of words with at most s-1 deficient middle variables each. Apply it to both sides. It permits arbitrarily large starting polynomials and words.

Selecting witnesses. Consider a valid summand absorption w \ensuremath{\leq} P after this expansion. Keep the target's end variables and its deficient middle variables fixed. For a short word keep all of its variables. The fixed set H has size at most h. Every other target variable occurs in the middle exactly e times and nowhere in the ends. On an assignment beta, assign the variables of H their beta values. Give all remaining target variables the meet t of their beta values if there are any. Assign outside variables the least element. This assignment alpha is coordinatewise below beta and \(w(alpha)=w(beta)\).

For that equality, permute the free variables' blocks together. Each block contributes an eth power. Apply the fourth table condition to two adjacent blocks at a time, keeping k preceding and l following values and restoring any unused outside context. This changes those two blocks to the eth power of their meet. Iteration produces \(t^{e}\). The same reduction on the uniformly replaced input produces \(t^{e}\). A single free block is immediate, and the empty free set requires no replacement. This argument fixes all occurrences of each variable of H, even when it occurs in both ends and middle.

There are at most \(q^{h+1}\) such assignments alpha. For each, retain a right word from the expanded P for every one of the at most q attained values. Their union Q is nonempty and has at most B members. For its corresponding alpha one has \(w(beta)=w(alpha)\) \ensuremath{\leq} \(P(alpha)=Q(alpha)\) \ensuremath{\leq} Q(beta) by monotonicity. Thus w \ensuremath{\leq} Q is valid in A.

Compressing and reversing. Reserve the end and deficient variables of each of w and the retained words, at most h(B+1) variables. All other variables occur only in middle positions, with multiplicity exactly e in each word where they occur. Partition them by common incidence across these at most B+1 words and identify each group to a fresh variable. At most \(2^{B+1}\)-1 nonempty groups occur. Normalize middle multiplicities in the resulting valid absorption. It fits the finite family. Substitute the nonempty product of all members for each group variable. Where the compressed group occurs, it occurs exactly e times, so every member recovers multiplicity e. The original middle normal forms and all protected end positions are therefore recovered by the first two identities. Add the omitted summands, use these derivations on all summands in both directions, and undo the initial expansions. This proves unrestricted completeness.

Catalogue applications. For each row \(k=l\). The following parameters were verified directly on the labelled tables.

\noindent \(k=l\)\enspace\textbar\enspace e\enspace\textbar\enspace s\enspace\textbar\enspace Types\par
\noindent 1\enspace\textbar\enspace 2\enspace\textbar\enspace 2\enspace\textbar\enspace 6671, 6827, 7086, 7203, 7290, 7299, 7614, 7641, 7693, 8549, 10214, 10232, 10397, 10405\par
\noindent 1\enspace\textbar\enspace 3\enspace\textbar\enspace 3\enspace\textbar\enspace 7296, 7633, 10046, 10058\par
\noindent 2\enspace\textbar\enspace 2\enspace\textbar\enspace 2\enspace\textbar\enspace 8905, 9584, 10130, 10131, 10141, 10142, 10151, 10155, 10165, 10177, 10213, 10396, 11846, 11919, 11945, 11948\par
Uniform bounds for all rows are h at most 5, \(B=78125\), \(D=5(B+1)+2^{B+1}-1\) and maximum length 4+3D. These bounds specify a finite basis even when enumerating all its members is impractical. The applicable word laws for each type use its own k, l, e, s.

\subsection{A finite basis for 6680 from contained head and simple-final witnesses}\label{case:171}
Signature (+, multiplication), without constants. This specifies a finite family, not an enumerated small basis.

Addition rows 00004,\allowbreak{}01114,\allowbreak{}01224,\allowbreak{}01234,\allowbreak{}44444 give the \(chain 3<2<1<0<4\). Multiplication rows are 01130,\allowbreak{}33333,\allowbreak{}33333,\allowbreak{}33333,\allowbreak{}41134. A literal retains its value. A nonlinear word has value 3 if its head is 1,\allowbreak{}2, or 3, if any factor is 3, or if any nonfinal factor is 1 or 2. Otherwise all prefix factors are 0/\allowbreak{}4: a final 1/\allowbreak{}2 gives 1, and a final 0/\allowbreak{}4 gives the original head 0/\allowbreak{}4. Direct induction from the rows proves this for all lengths.

\par\smallskip\noindent\textit{Exact criterion.}\quad 
A literal target is absorbed exactly when that literal is an original upper word. For a nonlinear target w, put \(C=c(w)\), \(h=head(w)\), \(t=last(w)\), and let U be the original upper words contained in C. Then \(w\leq P\) iff:

H: some u in U has head h; F: if t is simple in w, an original u in U either omits t, or is NONLINEAR with t as its simple final variable.

The literal t does not satisfy the second alternative of F. If t repeats in w, F is absent. These are original witnesses, not arbitrary derived words.

Necessity of H follows from \(h=4\), other C variables 0, outsiders 3. The target is 4; without H all upper values are 0 or 3. For F assign \(t=2\), other C variables 0, and outsiders 3. The target is 1. In the absence of either F alternative, a contained upper is literal t (value 2) or has t nonfinally (value 3); every noncontained word has factor 3 and also has value 3. The upper sum is at most 2. A missing original literal h is separated by \(h=2\), outsiders 3: every upper is 3 whereas the target is 2.

For sufficiency, target 3 is least. A nonlinear target 0/\allowbreak{}4 has all variables valued 0/\allowbreak{}4; H's contained same-head witness has exactly that value. Target 1 has all prefix variables 0/\allowbreak{}4 and final variable 1/\allowbreak{}2, which must consequently be simple. An F witness omitting t has value 0/\allowbreak{}4; a nonlinear simple-final t witness has value 1. Either absorbs 1. A nonlinear target never has value 2. This proves exactness for unrestricted variables, words, and polynomial sizes.

\par\smallskip\noindent\textit{Structural laws and normal words.}\quad 
The valid identities \(xyzt=xzyt\) and \(xxy=xy\) suffice for the following normal form: a literal is unchanged; a nonlinear word retains its original head and final positions, and between them has the increasing list of distinct original interior variables other than the head. First sort the interior by adjacent swaps, whose prefix and suffix are nonempty. Delete adjacent duplicate interior letters by \(xxy=xy\) in context. Move an interior head copy to immediately after the head and delete it with the same law, whose remaining suffix is nonempty. Keep the final position throughout. If a non-head final variable also occurs internally, it remains internally once. If head equals tail, the two endpoint copies remain. Thus simple-final versus repeated final is preserved exactly; multiplicity is not silently discarded.

\par\smallskip\noindent\textit{Two witnesses and reversible compression.}\quad 
Choose one original witness for H and, if needed, one for F. Their sum already absorbs the target by the exact criterion. Normalize the target and these at most two upper words. Reserve their head and final variables: at most 6 variables. Group every other variable by its nonempty occurrence pattern among these at most three words: at most 7 groups. Nonreserved variables occur only in the idempotent, permutable interior, once per normal word.

Replace reserved variables by distinct fresh variables and each group by a single fresh variable, then normalize. There are at most 13 variables and at most 15 positions in each word (two endpoints and at most 13 interior variables). Literals remain literals. H's content containment and head are preserved. F's omission, simple final status, and nonlinear status are preserved: endpoint variables are reserved, and no distinct interior variable can become the final variable. Therefore the compressed absorption is valid.

Substitute each reserved variable back to its singleton, and each group variable to the nonempty increasing product of the group's original members. The expansion has exactly the original head, final, interior content and final repetition status. Literal swaps and deletions under the two structural laws normalize it to the original normal word. This is an explicit inverse equational derivation for all three words, not just preservation in one direction.

\par\smallskip\noindent\textit{Fixed finite basis.}\quad 
Take the two structural laws plus ALL valid absorptions on a fixed 13-variable alphabet, word length at most 15, and at most two nonempty upper words. This is a finite effectively specified family; the criterion above decides validity. Every valid arbitrary absorption has the selected two-witness reduction. Its compression belongs to the fixed family; inverse nonempty substitution and structural laws recover that reduction; inclusion of its original summands gives the full absorption. This proves completeness with no bound on variables, word lengths or polynomial size. Distributive expansion and both absorption directions prove every identity. Therefore 6680 and its opposite are finitely based. Catalogue transfers require independently complete mutual HSP certificates.

\subsection{Finite bases from products determined by a meet between fixed positions}\label{case:172}
Let S be a finite additively idempotent semiring in the signature of two binary operations, without constants. Order S by x \ensuremath{\leq} y when \(x+y=y\). Suppose this order has a least element. It is then a finite lattice, since the join of the common lower bounds of any two elements is their meet. Write q for the size of S.

Fix nonnegative k, l, with \(d=k+l\). Let U, V be lists of k and l distinct formal variables, also distinct from a, b. Assume that for all assignments of these d+2 letters to S one has the following equality of values:

{\small\par\smallskip
\CaseLine{\(U a b V = U (a meet b) V.\qquad (1)\)}
\smallskip}

A zero-length list means that the indicated positions are absent. The meet in (1) is a test on values in the table. The identities used as axioms below have the original signature. Put

{\small\par\smallskip
\CaseLine{\(B = q^{d+2},\qquad D = d(B+1) + 2^{B+1} - 1.\)}
\smallskip}

Theorem. A finite basis for S consists of the semiring axioms, the two word identities \(UabV=UbaV\) and \(UaaV=UaV\), and all valid absorptions w \ensuremath{\leq} P on a fixed alphabet of D variables where every word has length at most D+d and P has at most B distinct summands. Absorption denotes the equation \(w+P=P\). This is a finite, effectively specified family.

Words. In a word longer than d, retain the first k and last l positions, and call the nonempty intervening list M. The value is the product of those end positions with the single lattice meet of all values in M in between. To collapse two neighboring entries of M, take the k positions immediately before and the l positions immediately after them. Equation (1) applies to those values, and associativity restores the unused outside context. Induction leaves one middle value. The same argument using the two word identities permits each adjacent interchange and each adjacent duplicate deletion inside M. Thus M normalizes to the sorted list of its distinct variables. The entire list of end positions is retained verbatim, including any repetitions and any letter that also occurs in M. Words of length at most d retain all their positions.

Selecting original summands. Suppose w \ensuremath{\leq} P holds in S. Let H be the variables in the d retained positions of w, or all variables of a short w. Let F be the other variables of w. Take all assignments that give arbitrary independent values to H, one common value to F when it is nonempty, and the additive least value to every variable outside w. There are at most \(q^{d+1}\) such assignments. For each assignment select one original right word for every value attained by a right word on that assignment. Their union Q is nonempty and has at most B words.

Given any assignment e, use the selection for its values on H and the meet of its values on F when F is nonempty. Call that assignment f. Then f \ensuremath{\leq} e coordinatewise. Also \(w(f)=w(e)\). Indeed the fixed variables retain their values even at their middle occurrences, and the meet contributed by the nonempty set F is preserved by its uniform replacement. Short words and F empty require only the values on H. By distributivity all word functions are monotone. Consequently

{\small\par\smallskip
\CaseLine{\(w(e) = w(f) \leq  P(f) = Q(f) \leq  Q(e).\)}
\smallskip}

The equality of sums at f holds in arbitrary finite join semilattices, including ones with incomparable elements. This proves w \ensuremath{\leq} Q for every assignment. Bottom and meet have been used only to choose assignments.

Recovering all variables. Normalize w and Q as above. Reserve every letter in a retained end position and all letters of short words. There are at most d(B+1) reserved letters. Group the remaining letters by the set of words in which they occur, counting the target as a separate word. At most \(2^{B+1}\)-1 nonempty patterns occur. Identify each group to its own fresh letter and normalize the resulting words again. Validity is preserved by identification. The compressed inequality fits the finite family, since a middle list has at most D distinct letters. Replace each group letter by the nonempty word listing all its members in a fixed order. Each such occurrence is in a middle list. Its expansion recovers exactly the original middle content in each word, and every end position is fixed. The word identities recover the original normal forms. This derives w \ensuremath{\leq} Q. Adding the other original summands gives w \ensuremath{\leq} P. Applying this argument to every summand on each side after distributive expansion proves all identities of S, for arbitrary numbers of variables, degrees and numbers of summands.

Applications. In catalogue labels the following tables satisfy (1). For every row \(d=2\). The uniform bounds are \(B=625\) and \(D=1252+2^{626}-1\).

\noindent k\enspace\textbar\enspace l\enspace\textbar\enspace Type IDs\par
\noindent 1\enspace\textbar\enspace 1\enspace\textbar\enspace 6743, 7081, 7214, 7312, 8504, 8541, 9255, 9303, 9588, 9890, 10355, 10399, 11922, 12161\par
\noindent 0\enspace\textbar\enspace 2\enspace\textbar\enspace 8088, 11958\par
\noindent 2\enspace\textbar\enspace 0\enspace\textbar\enspace 8738, 12929\par
For example, for 8738 the image of a product of two elements is contained in \{0,\allowbreak{}1,\allowbreak{}3,\allowbreak{}4\}. Multiplication from any such state on the right by 0 or 4 fixes it and multiplication by 1,\allowbreak{}2,\allowbreak{}3 sends it to 3. The upper set \{0,\allowbreak{}4\} is closed under meets in its chain. For 12929 the two-letter image is contained in \{1,\allowbreak{}2,\allowbreak{}3,\allowbreak{}4\}. Each state is fixed by a right letter in \{0,\allowbreak{}1,\allowbreak{}2,\allowbreak{}4\} and goes to 3 otherwise. These observations establish (1) with \(k=2\), \(l=0\) in both examples. For the remaining tables use the same test in (1), taking the meet in the order defined by addition.

\subsection{Twelve finite bases from contained multiplicity witnesses}\label{case:173}
Use the constant-free additively idempotent semiring signature. For a target word w and a nonempty upper polynomial P write C for the content of w, h for its head and U for the original summands of P supported in C. Every criterion below requires C to be covered by the contents of U. Call x simple in a word if it occurs exactly once. Let G say that h starts some original word of P. Let H consist of heads of U and T of noninitial letters of U. Put L = (H meets T) or (h belongs to H and the noninitial letters of w all belong to T).

\noindent Source\enspace\textbar\enspace Head condition\enspace\textbar\enspace Requirement for each simple x of w\par
\noindent 6745\enspace\textbar\enspace G\enspace\textbar\enspace Some word of U contains x once\par
\noindent 7204\enspace\textbar\enspace G\enspace\textbar\enspace Some word of U contains x at most once\par
\noindent 7790\enspace\textbar\enspace L\enspace\textbar\enspace Some word of U contains x once\par
\noindent 7914\enspace\textbar\enspace L\enspace\textbar\enspace Some word of U contains x at most once\par
\noindent 6826\enspace\textbar\enspace none\enspace\textbar\enspace Easy or head case below, with full coverage in the head case\par
\noindent 6828\enspace\textbar\enspace none\enspace\textbar\enspace Easy or head case below, with a single witness in the head case\par
For the last two rows the easy case is a word of U omitting x or containing it exactly once noninitially. The head case requires \(x=h\). For 6828 some word of U must contain x once. For 6826 the words of U containing x once must together cover C. When the easy case is absent these once-only occurrences are all initial. All conditions in each row together are equivalent to w \ensuremath{\leq} P.

\par\smallskip\noindent\textit{Table proof.}\quad 
In 6745 the order is \(2<1<0<4\) and \(2<1<3<4\). A head 2 yields 2. Any other word meeting 1 or 2 yields 1. On \{0,\allowbreak{}3,\allowbreak{}4\} a single occurrence of 3 gives 3, two or more give 1, and in the absence of 3 the value is the maximum of the entries 0 and 4. In 7204 the order is \(3<2<1<0<4\). A head 3 gives

3. Otherwise an occurrence of 2 or 3 gives 2. On \{0,\allowbreak{}1,\allowbreak{}4\} a single occurrence of 1 gives 1, two or more give 2, and the remaining words give the maximum of the entries 0 and 4. In both cases the special letter survives multiplication by 0 or 4 on either side, and its second occurrence produces the lower value. The two bottom rows are left zero rows as described. This proves the formulas by induction.

The least bound forces G: give all original heads the least value and a missing target head the next greater value. With G in force a target at that next value is covered by the same-head word. For high targets use contained originals. To force full coverage set a missing variable 4, other target variables 0 and outside variables the next-to-least value. For a target on \{0,\allowbreak{}4\} coverage detects a 4 and supplies an all-zero word if all letters are 0. For a target with one special occurrence at x the specified witness gives the same special value, or, where omission is permitted, a higher ordinary value. For necessity put the special value at x, 0 on C minus \{x\}, and the next-to-least value outside. In 6745, without an exact witness the right values are 0,\allowbreak{}1 or 2, whose sum is at most 0 and does not dominate 3. In 7204 without an at-most-once witness all right values are at most 2, below 1.

In 7790 the bottom 1 absorbs multiplication and 0,\allowbreak{}2,\allowbreak{}3 are incomparable below 4. In the absence of 1, one occurrence of 2 gives 2 and at least two give 1. In 7914 the bottom 2 absorbs, the order is as in 6745, one occurrence of 1 away from bottom gives 1 and at least two give 2. In both algebras, on the remaining alphabet \{0,\allowbreak{}3,\allowbreak{}4\} all zeros give 0, an initial 3 followed by zeros gives 3, and everything else gives 4. Appending a letter proves these rules: the special value persists on a high operand, a second special makes bottom, and a following 3 or 4 after a high prefix gives 4. Give outsiders bottom. At bound 0 a missing cover letter is separated by giving it 4 and other target letters 0. On the high alphabet, at bound 3 every contained word must start in 3 and otherwise use zeros. This is impossible if H meets T. Otherwise coverage partitions C into H, T and forces 3 on H, 0 on T. The target is 3 under this assignment exactly as specified by L. A high target also has a high contained witness at the bound of the special element. For a target with a single special x the same witness argument applies as above. Without an exact witness in 7790 put \(x=2\), all other target letters 0 and outsiders 1: the right sum is at most 0. Without an at-most-once witness in 7914 put \(x=1\), other target letters 0 and outsiders 2: the right sum is 2. These tests include literals.

6826 and 6828 have order \(3<2<1<0<4\). A 3 anywhere or two occurrences from \{1,\allowbreak{}2\} give 3. With no entry from \{1,\allowbreak{}2,\allowbreak{}3\}, take the maximum on \{0,\allowbreak{}4\}. With exactly one entry from \{1,\allowbreak{}2\}, the value is 1 except for an initial

2. In 6828 the initial 2 always gives 2 in this case. In 6826 it gives 2 with an all-zero suffix, and 1 when a 4 occurs in the suffix. The rows of 1 and 2 on inputs 1 or 2 both give 3. The row of 2 on 4 gives 1 in 6826 and 2 in 6828, while a noninitial 2 after a 0/\allowbreak{}4 prefix gives

1. These checks prove the formulas inductively.

Use 3 on outsiders. Full coverage is necessary by the missing-4 test and handles targets on \{0,\allowbreak{}4\}. With one special x, an omitting witness has high value, and a noninitial once-witness gives 1. If neither exists and x is noninitial in the target, assigning \(x=2\) and the other target variables 0 separates the target 1 from right values at most 2. If x is the head but has no once-witness, \(x=1\) and the other letters 0 separate 1 from 3. In the remaining head case every once-witness starts in x. For 6828 one such witness matches a 1 or 2 target. For 6826 its value also matches when the target has no 4, or \(x=1\). For \(x=2\) with a 4 elsewhere, coverage by once-witnesses supplies a witness containing that 4, giving 1. Conversely if these witnesses miss y, set \(x=2\), \(y=4\), the rest of C to 0 and outsiders to 3. The target is 1 and the right sum is at most 2. This proves every criterion without size bounds.

\par\smallskip\noindent\textit{Finite bases.}\quad 
For each source take the semiring axioms, the three identities

{\small\par\smallskip
\CaseLine{\(xyz = xzy\qquad xxx = xx\qquad xyz = xyyz + xyzz\)}
\smallskip}

and its valid absorptions with at most four right words on at most 247 variables and of word length at most 494. Each of the displayed identities is a direct at-most-three-variable table identity. The first two preserve a head and allow the tail to be sorted, capping total multiplicities at two. To cap the head's total multiplicity, move its tail occurrences next to the head before applying \(xxx=xx\).

Every valid step with at most four upper words is derivable. Normalize the at most five words, then reserve all head variables, at most five. For each other variable record its vector of total multiplicities 0,\allowbreak{}1,\allowbreak{}2 across these words. There are at most \(3^{5}\)-1 such nonzero vectors. For each nonempty class put one new representative into each word with the common multiplicity in that word. The criteria remain true: support containment and coverage, omission, exact occurrence, and coverage by the once-witnesses all agree coordinatewise. The heads are reserved and every occurrence of any unreserved variable is noninitial. Thus the sets of eligible upper words, crossings at heads and extent of each role cover are preserved. Compact words have length at most 494. Replace each representative by the fixed nonempty product of its class. After tail sorting this restores the common multiplicity of every original member in every word. This inverse substitution proves the local step from the finite family.

\par\smallskip\noindent\textit{Complete derivations.}\quad 
An available word is one already derived below P. Call a word saturated if all its variables occur at least twice. In the G rows and the two rows without a head condition, first saturate each u in U on its own support, preserving its head. This is a valid one-upper-word step. Merge the resulting words two at a time to saturated words on the union, keeping a head supplied by an operand. They have no simple variables to protect. At full support C orient the saturated word to the target head. In the G rows include an original h-word at this step. It may have outside variables, while the saturated operand covers C.

In the L rows, if H meets T, two originals witness the crossing. From these derive a saturated word on their union. The word now has its own head noninitially and exhibits its own crossing. Enlarge the support by one contained original at a time, retaining saturation, and orient the final word to h. If H misses T, all contained words begin with the same simple head h, since C minus B(w) has at most this one member. Start with each original and saturate every variable other than h, retaining its single initial occurrence. Merge such words two at a time with the same convention. Each operand supplies a once-witness for the only possible simple letter h. Every nonhead variable occurs noninitially, so L holds in each merge. This also handles support \{h\}. Call the resulting full-support available word g.

There is one extra construction when 6826 has a simple head \(x=h\) without an easy witness. All contained originals containing x once begin with x, and their supports together cover C. In each such word keep x simple at its head and saturate all other variables. This step is valid because its original is a once-witness covering the entire local support. Merge the resulting words two at a time, keeping x simple at the head: the two operands are once-witnesses covering their union. At termination one available once-witness for x covers all of C. This construction is needed even if the number of original head witnesses is unbounded.

Normalize w. If two distinct nonhead variables y, z are simple, move them to its final two positions by tail permutations. The preceding block is nonempty. Splitting by \(xyz=xyyz+xyzz\) and sorting gives two new targets with the same support and head, each with one fewer simple nonhead variable. All remaining simple tests use the original witnesses and all positional requirements are unchanged. Repeating reduces the proof to targets with at most one simple nonhead variable and possibly the simple head.

For such a target use g and at most one witness for each of its two possible simple letters. In the G rows g supplies the head. In the crossed L case g supplies a crossing, and in the uncrossed case g supplies the common simple head and every required noninitial letter. In 6826 a difficult simple head uses the full-support once-witness just constructed. All other cases use their original easy, exact, at-most-once or 6828 single head witnesses. Every witness lies in C. The criterion now proves this local step, with at most three upper words. The splitting identities reassemble the target from these derivations. Distribute arbitrary terms and derive their summand absorptions in both directions to obtain every valid identity. Thus the stated finite families are complete on unrestricted alphabets. Word reversal gives bases for the respective opposites. The labelled opposites of 6745,\allowbreak{}7204,\allowbreak{}7790,\allowbreak{}7914,\allowbreak{}6826,\allowbreak{}6828 are, respectively, 7467,\allowbreak{}7619,\allowbreak{}7454,\allowbreak{}7613,\allowbreak{}7642,\allowbreak{}7643.

\subsection{Finite bases for four contained-position algebras and their opposites}\label{case:174}
Work in the constant-free signature of additively idempotent semirings. For a nonempty word v let C(v), A(v), B(v) be its content, nonfinal content, and noninitial content. Write h(v), f(v) for its head and final letter. For a target w and a nonempty sum P of original words put

{\small\par\smallskip
\noindent C = C(w), U = \{u in P : C(u) is contained in C\}\par
\noindent K = union C(u), J = union A(u), T = union B(u), F = set of f(u), over U.\par
\smallskip}

Empty unions are empty. Write G for the condition that an original word in P starts with h(w). Such a word may contain variables outside C. Write E for the existence of u in U with A(u) contained in A(w). In particular, if the final letter occurs earlier in w then E just requires U nonempty. Here are exact criteria for w \ensuremath{\leq} P, meaning \(w+P=P\).

\noindent Type\enspace\textbar\enspace Conditions required together\par
\noindent 6748\enspace\textbar\enspace \(K=C\), G, and either J meets F or both A(w) is contained in J and f(w) belongs to F\par
\noindent 7215\enspace\textbar\enspace \(K=C\), G, E\par
\noindent 13052\enspace\textbar\enspace \(K=C\), G, A(w) is contained in J\par
\noindent 11236\enspace\textbar\enspace \(K=C\), B(w) is contained in T, E\par
\par\smallskip\noindent\textit{Word values and necessity and sufficiency.}\quad 
For 6748 the additive order has \(2<1<0<4\) and \(2<1<3<4\), with 0 and 3 incomparable. A head 2 gives 2. Any other word meeting \{1,\allowbreak{}2\} gives 1. On \{0,\allowbreak{}3,\allowbreak{}4\} the value is 0 on all zeros, 3 on a sequence of zeros followed by a final 3, and 4 otherwise. Indeed 1 and 2 are left zeros, a low entry after high entries gives 1, and on the high alphabet a nonfinal 3 or 4 produces 4. The formula also gives the literal values.

For 7215, 13052 and 11236 the order is \(3<2<1<0<4\). For 7215 a head 3 gives 3. On all other words the value exceeds 2 precisely when the word avoids 2 and 3 and has no nonfinal 1. Within this case final 1 gives 1, all zeros give 0, and a word on \{0,\allowbreak{}4\} containing 4 gives 4. In the multiplication table the rows of 3 and 2 are constant and 1 followed by any entry gives 2. On \{0,\allowbreak{}4\} multiplication is maximum; a following 1, 2 or 3 gives respectively 1, 2 or 2.

For 13052 a head 3 gives 3 and any other word meeting \{2,\allowbreak{}3\} gives 2. On \{0,\allowbreak{}1,\allowbreak{}4\} the word is 1 on all ones, 0 on a sequence of ones followed by a final 0, and 4 otherwise. In this last alphabet, the row of 1 copies the next entry, and 0 or 4 followed by a high entry gives 4.

For 11236 a 3 anywhere or a 2 in a nonfinal position gives 3. Otherwise a final 2 gives 2. On \{0,\allowbreak{}1,\allowbreak{}4\} the word is 1 on all ones, 0 on an initial 0 followed by ones, and 4 otherwise. The rows of 2 and 3 are constant 3. A final 2 or 3 following any high value copies that final entry. On the high alphabet a following 1 preserves the value and a following 0 or 4 gives 4. Each of the four word descriptions follows by induction on the positive word length.

Test the implications "all words of P lie below or equal to b implies w does" for all principal bounds b. This is equivalent to absorption since the value of P is one such bound. For the first three types the least value occurs exactly on words headed by that value. This bound gives G. For necessity give a missing target head the next higher value and all other letters the least value.

For 6748 all words meeting \{1,\allowbreak{}2\} lie below every bound except possibly the least. At a high target it suffices to inspect U for positive witnesses. For counterassignments give every outside letter 1. At bound 1 some contained word is needed. At bound 0, if all contained words on the high alphabet are 0, their covered letters are 0. For necessity of the cover give one missing letter 3 and the other target letters 0. At bound 3 each contained word must have its nonfinal letters 0 and its final letter 3. If J meets F this is impossible. Otherwise \(K=C\) partitions C into the disjoint J, F. The assignment 0 on J and 3 on F gives \(w=3\) precisely under the displayed requirements. This proves the criterion.

For 7215 a target above 0 has a 4 and otherwise only values in \{0,\allowbreak{}4\}. The cover supplies a contained witness containing that 4. For necessity, give a missing letter 4, the other target letters 0 and outsiders 2. Above 1 any contained word suffices when the target is above 1. Above 2 the target avoids 2,\allowbreak{}3 and its nonfinal letters avoid 1. A witness for E also exceeds 2. Conversely assign 0 on A(w), 1 on C minus A(w), and 2 elsewhere. Exactly the upper words prescribed by E exceed 2. This includes a literal or a repeated final letter.

For 13052 at bound 2 any contained word suffices for a high target. At bound 1 a cover witness detects any target letter other than 1 on the high alphabet. Assign 0 to a missing letter, 1 to other target letters, 2 to outsiders for necessity. At bound 0, if all contained words on this alphabet have value at most 0, no covered letter is 4 and every letter of J is 1. This forces the target at most 0 under its prefix requirement. Giving 0 to a missing nonfinal letter, 1 to the rest of C and 2 to outsiders proves necessity. At the least bound use G as above.

For 11236 the bound 2 and the bound 1 work as for 13052. At bound 0 the exceptional 0 must occur initially: use the same separating assignment at a missing *noninitial* letter, and use T for sufficiency. For the least bound a word exceeds 3 precisely when all its values avoid 3 and no nonfinal value is 2. A witness for E therefore suffices. For necessity give 1 on A(w), 2 on C minus A(w), and 3 to outsiders. These arguments exhaust the proper bounds also for repeated variables.

\par\smallskip\noindent\textit{A finite family of elementary steps.}\quad 
For each source take the semiring axioms, the word identities

{\small\par\smallskip
\CaseLine{\(abcd = acbd\qquad abbc = abc\)}
\smallskip}

and every valid word absorption on at most 41 letters, length at most 43, with at most four right words. The finite family is specified by these bounds and the actual table. The word formulas prove the two identities. They permit arbitrary interior permutations and deletion of duplicate interior copies, with both endpoint positions fixed. A literal retains its single position.

Every valid step with at most four right words follows from this basis. Normalize its at most five words. Reserve every variable used at any endpoint, at most ten. Group the remaining letters by their exact nonempty set of word indices, at most 31 groups. Identify each group to a different fresh letter, keeping the reserved variables distinct. The resulting absorption is valid by substitution. Normalize the interiors for length at most 43. In the inverse substitution replace each group letter by the nonempty ordered product of its members. Where it occurs, every member was interior and occurred once after normalization. Sorting and deleting interior duplicates recovers every original word simultaneously. Reserved variables also retain their own interior occurrence whenever present.

An available word has its absorption into P already derived. Every original right word is available. A valid step with at most four available upper words makes the new target available by transitivity.

\par\smallskip\noindent\textit{6748 and 7215.}\quad 
For 6748 choose an original h(w)-word g. When J meets F choose two contained originals witnessing the two roles of their common letter, and one contained original covering \(h=h(w)\). Take S to be the union of their contents. A word of content S beginning and ending in h, with two actual endpoint positions, is available from those three words and g. The crossing witnesses belong to the contained part of this local step. The available word itself now exhibits a crossing at h. Adjoin each further contained original's content, using just the preceding word and that original to obtain another word of the same form on the union. At the end the last word alone gives w, including when w is a literal.

Otherwise F is disjoint from J. Since A(w) is contained in J and \(C=A(w)\) union \{f(w)\}, every contained word ends in f(w), and this letter never occurs nonfinally in the target or a contained word. Choose a contained u including h, and start with \(S=C(u)\). Delete from w the letters outside S to form its restriction. The restriction is available from u, g since all letters in S except f(w) are nonfinal in u. Adjoining the content of another contained word v uses the old restriction and v to derive the new restriction. Both have the same final letter and supply all required prefix letters.

For 7215 choose a contained witness u for E, one contained word covering each endpoint, and an original head witness g. Let S be the union of the contents of the first three. Both endpoints remain in the restriction of w to S, so its nonfinal set is A(w) intersect S. The witness u still works. The four selected words give the restriction. On enlarging S by the content of another contained word, the preceding restriction and that word cover the enlargement. The preceding restriction supplies both the head and E. Iteration gives w.

\par\smallskip\noindent\textit{13052 and 11236.}\quad 
For 13052 one can realize a nonempty support S and nonfinal set D contained in S whenever S minus D has at most one letter. For nonempty D start with any chosen member of D, list its other members, and end with the missing letter or any member of S if none is missing. Empty D gives the one-letter word. For 11236 use the reversed construction for noninitial roles. With nonempty D one can require the final letter to occur twice: take a realizing word and duplicate its final letter, which is already noninitial. This keeps both sets.

Start from each u in U separately. For 13052 write S, D for the support and nonfinal set of the running word. If the sole possible missing letter d of S minus D belongs to A(w), choose an original v in U with d in A(v). Merge to \(S'=S\) union C(v), \(D'=D\) union A(v). Only the last letter of v might be absent from D'. The old head lies in D', including when the old word was the literal d. Realize the two sets with that head. The old word and v prove this step by the criterion. D grows strictly and remains inside C, so repetition terminates. Each resulting word has every one of its letters that belongs to A(w) in its nonfinal set. Merge these words two at a time taking both unions. The only possible missing nonfinal letter belongs to C minus A(w). Choose the head from an operand with a nonempty prefix. If both prefixes are empty, the operands are the same literal. The two operands give each step. The resulting word covers C with all required nonfinal letters. It together with an original head witness gives w.

For 11236 follow the same algorithm with noninitial roles and B(w). If d is missing from the running noninitial set but needed in B(w), choose a contained original with d noninitial. At most the initial letter of that original is missing from the union of noninitial sets. At each extension and merge with a nonempty union of roles use the realization with repeated final letter. It satisfies E in the local step since its entire support is nonfinal and both operands are contained. Empty union at the final merging stage again means identical literals, already available. The first word of each separate construction is an original, so no condition is needed to initialize it. The final merged word supplies the cover and B(w); use it and the originally chosen contained witness for E to recover w. The argument allows this witness to be a literal and allows a one-letter target.

The four constructions are unrestricted in alphabet, lengths and sums. Distributivity expands each term into a nonempty sum of words. Proving all summand absorptions in both directions gives every valid identity. Reversing all words in derivations gives a complete finite basis for each multiplicative opposite. The opposite partners of 6748, 7215, 13052 and 11236 are 7821, 8602, 11232 and 13316 respectively.

\subsection{乘法核心的重复元素：八组交换完整有限基}\label{case:175}
\par\smallskip\noindent\textit{1. 定义与有限基.}\quad 
底集均为 \(\{0, a, b, e, f\}\)。先在 \(\{0, a, e, f\}\) 上定义交换乘法：\(0\) 为零，\(e\) 为幺元，\(a^2=0\)、\(f^2=f\)、\(af=a\)。下表的映射 \(\mu\) 固定未注明的元素，乘法定义为 \(xy=\mu(x)\mu(y)\)。加法取所列偏序的并；只列覆盖关系。

\noindent 记号\enspace\textbar\enspace 来源\enspace\textbar\enspace 加法序\enspace\textbar\enspace \(\mu\) 的非固定值\enspace\textbar\enspace 到目录标号的映射\par
\noindent \(A\)\enspace\textbar\enspace 6825\enspace\textbar\enspace \(0<b<a<e<f\)\enspace\textbar\enspace \(b\mapsto a\)\enspace\textbar\enspace \((3,\allowbreak{}1,\allowbreak{}2,\allowbreak{}0,\allowbreak{}4)\)\par
\noindent \(B\)\enspace\textbar\enspace 10122\enspace\textbar\enspace \(0<b<a<f,\ 0<e<f\)\enspace\textbar\enspace \(b\mapsto a\)\enspace\textbar\enspace \((2,\allowbreak{}0,\allowbreak{}1,\allowbreak{}3,\allowbreak{}4)\)\par
\noindent \(C\)\enspace\textbar\enspace 11845\enspace\textbar\enspace \(0<b<a<e<f\)\enspace\textbar\enspace \(b\mapsto a,\ e\mapsto f\)\enspace\textbar\enspace \((3,\allowbreak{}1,\allowbreak{}2,\allowbreak{}0,\allowbreak{}4)\)\par
\noindent \(D\)\enspace\textbar\enspace 10125\enspace\textbar\enspace \(0<b<a<f,\ b<e<f\)\enspace\textbar\enspace \(b\mapsto a,\ e\mapsto f\)\enspace\textbar\enspace \((2,\allowbreak{}0,\allowbreak{}1,\allowbreak{}3,\allowbreak{}4)\)\par
\noindent \(E\)\enspace\textbar\enspace 10126\enspace\textbar\enspace \(0<b<a<f,\ 0<e<f\)\enspace\textbar\enspace \(b\mapsto a,\ e\mapsto f\)\enspace\textbar\enspace \((2,\allowbreak{}0,\allowbreak{}1,\allowbreak{}3,\allowbreak{}4)\)\par
\noindent \(F\)\enspace\textbar\enspace 10127\enspace\textbar\enspace \(0<b<a<f,\ e<f\)\enspace\textbar\enspace \(b\mapsto a,\ e\mapsto f\)\enspace\textbar\enspace \((2,\allowbreak{}0,\allowbreak{}1,\allowbreak{}3,\allowbreak{}4)\)\par
\noindent \(G\)\enspace\textbar\enspace 13110\enspace\textbar\enspace \(0<e<b<f,\ 0<a<f\)\enspace\textbar\enspace \(b\mapsto f\)\enspace\textbar\enspace \((2,\allowbreak{}3,\allowbreak{}0,\allowbreak{}1,\allowbreak{}4)\)\par
\noindent \(H\)\enspace\textbar\enspace 13115\enspace\textbar\enspace \(0<a<e<b<f\)\enspace\textbar\enspace \(b\mapsto f\)\enspace\textbar\enspace \((3,\allowbreak{}2,\allowbreak{}0,\allowbreak{}1,\allowbreak{}4)\)\par
这些定义给出交换加法幂等半环：乘法结合由固定核心上的结合性得到，各乘法平移保持指定的并。没有向语言中添加 \(0, e, f\) 常量。

八组基共同包含

\CaseDisplay{xy\approx yx,\qquad xy\approx x^2y+xy^2.\quad\textup{(\ensuremath{N})}}

其余可选恒等式编号如下：

\CaseDisplay{\begin{array}{lll}
(c)&xyz\leq xy,&(d)\quad x^2\leq x,\\
(u)&xy\leq x^2+y^2,&(t)\quad x^2y^2\leq x^2+y^2,\\
(r)&xy^2\leq x+(xy)^2,&(s)\quad xy^2\leq y+(xy)^2,\\
(q)&xy\leq x+(xy)^2,&(m)\quad xy^2z^2\leq x+y+(xyz)^2,\\
(f)&xy^2\leq x+y^2,&(v)\quad xy^2z^2\leq xy^2+z^2,\\
(b)&xy^2z^2\leq xy^2+x^2z^2,&(h)\quad x^2\leq x+yz.
\end{array}}

定理。 相对于通常加法幂等半环公理，以下均为完整有限基：

\noindent 半环\enspace\textbar\enspace 在 (N) 外添加的式子\enspace\textbar\enspace 总附加条数\par
\noindent \(A\)\enspace\textbar\enspace \(d, u, s, b\)\enspace\textbar\enspace 6\par
\noindent \(B\)\enspace\textbar\enspace \(d, f, b\)\enspace\textbar\enspace 5\par
\noindent \(C\)\enspace\textbar\enspace \(c, s\)\enspace\textbar\enspace 4\par
\noindent \(D\)\enspace\textbar\enspace \(c\)\enspace\textbar\enspace 3\par
\noindent \(E\)\enspace\textbar\enspace \(c, m\)\enspace\textbar\enspace 4\par
\noindent \(F\)\enspace\textbar\enspace \(c, s, h\)\enspace\textbar\enspace 5\par
\noindent \(G\)\enspace\textbar\enspace \(t, r, v, b\)\enspace\textbar\enspace 6\par
\noindent \(H\)\enspace\textbar\enspace \(u, q, b\)\enspace\textbar\enspace 5\par
\par\smallskip\noindent\textit{2. 组合判据.}\quad 
对右多项式 \(P\)，令 \(L\) 为其一次项变量集，\(\mathcal R\) 为其非线性词族。对非线性左词 \(w\)，置

\CaseDisplay{C_0=c(w),\quad S(w)=\{x:m_x(w)=1\},\quad
\mathcal R_0=\{v\in\mathcal R:c(v)\subseteq C_0\},\quad
U_0=\bigcup_{v\in\mathcal R_0}c(v).}

对 \(x\in S(w)\)，缩记

\CaseDisplay{\begin{array}{ll}
O_x:\ \exists v\in\mathcal R_0,\ m_x(v)\leq1,&
J_x:\ \exists v\in\mathcal R_0,\ m_x(v)=1,\\
Z_x:\ \exists v\in\mathcal R_0,\ m_x(v)=0,&
T_x:\ L\cap(C_0\setminus\{x\})\ne\varnothing.
\end{array}}

在每个生成半环中，\(x\leq P\iff x\in L\)。对非线性左词，\(w\leq P\) 当且仅当下表第二列成立，并且第三列对每个 \(x\in S(w)\) 成立。

\noindent 半环\enspace\textbar\enspace 锚点或覆盖要求\enspace\textbar\enspace 一次变量保护要求\par
\noindent \(A\)\enspace\textbar\enspace \(C_0\subseteq L\cup U_0\)\enspace\textbar\enspace \(O_x\) 或 \(T_x\)\par
\noindent \(B\)\enspace\textbar\enspace \(C_0\subseteq L\cup U_0\)\enspace\textbar\enspace \(J_x\) 或 \(x\in L\) 且 \((T_x\text{ 或 }Z_x)\)\par
\noindent \(C\)\enspace\textbar\enspace \(\mathcal R_0\ne\varnothing\)\enspace\textbar\enspace \(O_x\) 或 \(T_x\)\par
\noindent \(D\)\enspace\textbar\enspace \(\mathcal R_0\ne\varnothing\)\enspace\textbar\enspace \(O_x\)\par
\noindent \(E\)\enspace\textbar\enspace \(\mathcal R_0\ne\varnothing\)\enspace\textbar\enspace \(O_x\) 或 \(x\in L\) 且 \(T_x\)\par
\noindent \(F\)\enspace\textbar\enspace \(\mathcal R_0\ne\varnothing\)，或 \(L\cap C_0\ne\varnothing\) 且 \(\mathcal R\ne\varnothing\)\enspace\textbar\enspace \(O_x\) 或 \(T_x\)\par
\noindent \(G\)\enspace\textbar\enspace \(C_0\subseteq U_0\)\enspace\textbar\enspace \(J_x\) 或 \(x\in L\)\par
\noindent \(H\)\enspace\textbar\enspace \(C_0\subseteq U_0\)\enspace\textbar\enspace \(O_x\) 或 \(L\cap C_0\ne\varnothing\)\par
所有存在条件均要求真实的非空见证，不使用空和或空积。

\par\smallskip\noindent\textit{必要性.}\quad 
对一次项判据，在前六个半环中令左变量为 \(b\)、其余变量为 \(0\)；最后两个改令左变量为 \(a\)。没有相同一次项时，右侧均为 \(0\)。

对 \(A, B\)，覆盖失败时，把一个未覆盖变量赋为 \(f\)、其他内部变量赋为 \(e\)、外部变量赋为 \(0\)；左侧为 \(f\)，右侧不高于 \(e\)。对 \(G, H\)，若某个内部变量没有在内部非线性右词中出现，将它赋为 \(b\)、其他内部变量赋为 \(e\)、外部变量赋为 \(0\)。左侧为 \(f\)，而右侧至多为 \(b<f\)，即使该变量作为一次项出现在右侧也无济于事。

对 \(C, D, E\)，缺少内部非线性锚点时，将内部变量均赋为 \(e\)、外部变量赋为 \(0\)。左侧为 \(f\)，右侧至多为 \(e\)。对 \(F\)，若同时没有内部一次项，右侧为 \(0\)；若没有任何非线性右词，则把所有变量赋为 \(e\)，右侧为 \(e<f\)。

保护条件失败的见证如下。除特别说明外，取 \(x=b\)、其他内部变量为 \(e\)、外部变量为 \(0\)。

\textbullet{} 对 \(A, C, F\)，所有内部非线性右词重复 \(x\)，且没有其他内部一次项，所以右侧至多为 \(b<a\)。

\textbullet{} 对 \(D\)，所有内部非线性右词重复 \(x\)，一次项的和至多为 \(e\)，而 \(a\nleq e\)。

\textbullet{} 对 \(E\)，若 \(T_x\) 失败，右侧至多为 \(b\)；若 \(x\notin L\)，右侧至多为 \(e\)。两者均不吸收左侧 \(a\)。

\textbullet{} 对 \(B\)，若 \(x\notin L\)，改令 \(x=a\)：内部右词要么省略 \(x\) 而取 \(e\)，要么重复 \(x\) 而取 \(0\)，右侧至多为 \(e\)。若 \(x\in L\) 而 \(T_x, Z_x\) 均失败，则仍取 \(x=b\)，右侧至多为 \(b<a\)。

\textbullet{} 对 \(G\)，令 \(x=a\)，其余内部变量为 \(e\)。没有恰含一次 \(x\) 的内部右词，也没有一次项 \(x\)，右侧至多为 \(e\)，而 \(a\nleq e\)。

\textbullet{} 对 \(H\)，令 \(x=a\)，其余内部变量为 \(e\)。没有内部一次项，且所有内部非线性右词重复 \(x\)，右侧为 \(0<a\)。

\par\smallskip\noindent\textit{3. 有效性与规范形.}\quad 
非线性词值按核心乘法计算：出现零或至少两次零平方层元素 \(a, b\) 时为 \(0\)，恰有一次时为 \(a\)，否则在单位层取值。对 \(G, H\)，\(b\) 属于单位层而不是零平方层。对 \(C, D, E, F\)，所有单位层非线性词都取 \(f\)。

这一描述证明 (N) 有效。它也给出所列附加式的有效性：\((c)\) 在三元素核心 \(\{0, a, f\}\) 中是继续相乘不增加核心值；\((d)\) 直接由平方和加法序得到；\((u),(t)\) 在相关生成半环的核心层逐项成立。式 \((r),(q)\) 用相应一次项处理左侧为 \(a\) 的情况，用平方项处理稳定值。

对 \((s)\)，左侧为 \(a\) 时，\(y\) 必在单位层。在 \(A, C\) 中有 \(a\leq y\)；在 \(F\) 中，若 \(y=e\)，则 \(y+(xy)^2=e+0=f\)。式 \((m)\) 在 \(E\) 中遇到 \(x=b\)、\(y=e\) 时由 \(b+e=f\) 吸收；\(x=a\) 时由 \(x\) 吸收，其余稳定值由平方项处理。

对 \((f)\) 在 \(B\) 中的唯一特殊情形，\(x=b\)、\(y\) 在单位层，右侧 \(b+y^2=f\)。式 \((v),(b)\) 在左侧为 \(a\) 时由第一右项保留该值，在单位层时由某个右项保留最高层值。左侧为 \(0\) 时右側的非线性项均不低于 \(0\)。最后，\((h)\) 在 \(F\) 中只需注意 \(x=e\)：任意非线性值 \(v\in\{0, a, f\}\) 都满足 \(e+v=f=e^2\)。

由 (N) 得 \(x^3=x^2\)，以及每个非线性词的规范展开

\CaseDisplay{H(C_0)=\prod_{y\in C_0}y^2,
\qquad K_x(C_0)=xH(C_0\setminus\{x\}),\quad |C_0|\geq2.
\quad\textup{(\ensuremath{1})}}

展开反复将两个不同的一次变量之一平方，因而终止；每个分项至多有一个一次变量，内容保持不变。原词恰含一次 \(x\) 时，至少有一个新分项仍恰含一次 \(x\)。还得到每个非线性词 \(v\) 的 \(v^2\leq v\)。

\par\smallskip\noindent\textit{4. 三元素核心：\(C, D, E, F\).}\quad 
由 (c)，一个非线性词可以补乘任意非空词而不增加值。因此，内部非线性锚点的平方可补乘为 \(H(C_0)\leq P\)。在 \(F\) 中的另一种锚点情形，选 \(y\in L\cap C_0\) 和任一非线性右词 \(uv\)，由 (h) 得 \(y^2\leq y+uv\leq P\)，再补乘也得到同一稳定锚点。

若 \(O_x\) 成立，选择不重复 \(x\) 的内部非线性右词，补乘并截断指数，直接得到 \(K_x(C_0)\leq P\)。这已经证明 \(D\) 的完整性。

对 \(C, F\)，若保护见证是 \(g\in L\cap(C_0\setminus\{x\})\)，把目标规范词写成 \(vg^2\)，其中 \(v\) 非空，使用 (s)：

\CaseDisplay{K_x(C_0)=vg^2\leq g+(vg)^2=g+H(C_0)\leq P.}

对 \(E\)，若 \(x, g\in L\) 且 \(g\ne x\)，令 \(z\) 为 \(C_0\setminus\{x\}\) 中变量的乘积。在 (m) 中代入 \(x, g, z\)，指数截断后得到

\CaseDisplay{K_x(C_0)\leq x+g+H(C_0)\leq P.}

这里 \(z\) 必非空，且已经包含 \(g\)，所以没有空积问题。所有保护情形均已推出。

\par\smallskip\noindent\textit{5. 四元素核心：\(A, B\).}\quad 
由 (d)，所有右词和一次项都吸收其平方。对 \(A\)，(u) 给出稳定块合并；对 \(B\)，(f) 给出同样的合并：

\CaseDisplay{H(S)H(T)\leq H(S)+H(T).}

覆盖条件因此推出 \(H(C_0)\leq P\)。

先看 \(A\)。若 \(T_x\) 成立，按上一节的 (s) 推导得到 \(K_x(C_0)\leq P\)。否则选择 \(O_x\) 的内部右词。若它恰含一次 \(x\)，其规范分项给出一个 \(xH(S)\leq P\)。若它不含 \(x\)，平方后得到 \(H(S)\leq P\)，且 \(S\ne\varnothing\)。覆盖 \(x\) 的右项若为一次项，则 \(x^2\leq x\leq P\)；若为非线性词，其平方为 \(x^2H(T)\leq P\)。用 (u) 合并可得到初始 \(xH(S\cup T)\leq P\)，\(T\) 为空时直接代入 \(H(S),x\)。

对 \(B\)，若 \(J_x\) 成立，同样从恰含一次 \(x\) 的规范右分项开始。否则 \(x\in L\)，且有一个不含 \(x\) 的内部右项 \(v\)，它可以是一次项或非线性项。由 (f)、(d) 有

\CaseDisplay{xv^2\leq x+v^2\leq x+v\leq P.}

这也给出非空的初始 \(xH(S)\leq P\)。

接着保持 \(x\) 恰出现一次，合并所有内部右词与内部一次项的内容。若右词不含 \(x\)，使用 \(rH(T)\leq r+H(T)\)，其中 \(r=xH(S)\)；该式在 \(A\) 由 (u)、\(r^2\leq r\) 得到，在 \(B\) 由 (f) 得到。若右词有规范分项 \(xH(T)\)，将稳定块合并式乘以 \(x\)。若右词重复 \(x\)，用 (b) 合并它的平方 \(x^2H(T)\)。一次项 \(y\) 则先用 \(y^2\leq y\)，再按不含 \(x\) 的稳定块合并；若 \(y=x\)，没有新内容须加入。覆盖条件使最终内容达到 \(C_0\)。

\par\smallskip\noindent\textit{6. 单位层重复元素：\(G, H\).}\quad 
这两种情况的全部左变量都由内部非线性右词覆盖。平方各右词并用 (t) 或 (u) 合并，得到 \(H(C_0)\leq P\)。

在 \(G\) 中，若 \(x\in L\)，直接用 (r) 得到 \(K_x(C_0)\leq x+H(C_0)\leq P\)。否则由 \(J_x\) 选一个恰含一次 \(x\) 的规范右分项 \(xH(S)\leq P\)。合并其他内部右词的内容：不含 \(x\) 的右词用 (v)，恰含一次 \(x\) 的右词用 (t) 乘以 \(x\)，重复 \(x\) 的右词用 (b)。所有稳定块都非空，因初始词是非线性词。

在 \(H\) 中，若 \(L\cap C_0\ne\varnothing\)，选其中的 \(g\)，把目标左词写成 \(gv\)，由 (q) 得

\CaseDisplay{gv\leq g+(gv)^2=g+H(C_0)\leq P.}

否则从 \(O_x\) 的右词开始。若它省略 \(x\)，内部覆盖提供含 \(x\) 的右词，按第 5 节的 (u) 推导建立初始 \(xH(S)\leq P\)。其后的内容合并仍用 (u)、(b)。

至此，每个必要组合条件都已由相应有限基推出，故也充分。把任意多项式恒等式分解为双向逐词吸收，定理得证。

\par\smallskip\noindent\textit{7. 目录与核验.}\quad 
八个来源的等簇类大小依次为 \(1,\allowbreak{}1,\allowbreak{}1,\allowbreak{}1,\allowbreak{}1,\allowbreak{}4,\allowbreak{}1,\allowbreak{}1\)，合计 11 个当前交换类型。任意变量数完整性由上述推导证明；程序独立重建运算表，检查每条基、有限范围的组合判据，以及每个类型的双向完整词函数包含。未完成这些核验前，不计入合并统计。

\(A, C, D, E, F, H\) 满足平方和吸收式；其中三元素核心的四个来源由 (N)、(c) 推出它。因此还可使用平方定向乘法的项等价定理，所得非交换目录类型必须另行核验后合并。

\subsection{A finite basis for 8729 from two original contained witnesses}\label{case:176}
Signature (+, multiplication), no constants, relative to the additively idempotent semiring axioms. The finite family is specified, not an enumerated small list. All arguments allow arbitrary numbers of variables, arbitrary word lengths and arbitrary polynomial sizes.

Addition: 00004,\allowbreak{}01114,\allowbreak{}01224,\allowbreak{}01234,\allowbreak{}44444, with order \(3<2<1<0<4\). Multiplication: 03330,\allowbreak{}13331,\allowbreak{}13331,\allowbreak{}33333,\allowbreak{}43334. For a nonlinear word, any tail value in \{1,\allowbreak{}2,\allowbreak{}3\} makes the value 3. Otherwise all tails are 0 or 4: heads 0,\allowbreak{}1,\allowbreak{}3,\allowbreak{}4 retain their values, while head 2 becomes 1. This follows by induction directly from the table. Literals retain their values. Hence \(xyz=xzy\) and \(xyy=xy\), and a normal word is its head followed by the sorted set of its tail variables. The occurrence of the head in the tail must be retained. Tail swaps and deletions use nonempty prefixes and substitutions only.

\par\smallskip\noindent\textit{Exact criterion.}\quad 
A literal x is absorbed by P iff P contains the literal x. For nonlinear \(w=hT\), set \(C=\{h\}\) union T. The necessary and sufficient conditions are: (K) an original h-headed right word has its entire content contained in C; (R) h is in T, or an original right word has content contained in T, or an original nonlinear h-headed right word has tail contained in T. Literal h is allowed as K, but cannot by itself supply R for simple h.

Necessity for a missing literal x: set \(x=2\) and every other variable to 3. All right nonlinear words are 3, as are other literal heads; the target is

2. If K fails, set \(h=4\), other target variables to 0 and outsiders to 3. The target is 4. All right h-headed words contain an outside tail and therefore give 3. Other right heads are 0 or 3, so the right sum is at most

0. If R fails after K, h is not in T; set \(h=2\), \(T=4\), \(outsiders=3\). Target value is 1. Literal h gives only 2. Every nonlinear h-headed right word has a tail outside T and gives 3. Words headed in T also have a tail outside T, since no original right word is wholly contained in T, and give 3. Other heads give 3. The right sum is at most 2, strictly below 1.

For sufficiency, nonlinear targets can only have values 0,\allowbreak{}1,\allowbreak{}3,\allowbreak{}4. Value 3 is least. If the value is 0 or 4, the head has that value and every tail is 0 or 4; K supplies a word with the same value (literal or nonlinear). If the target is 1, its head is 1 or 2 and every tail is 0 or 4, so the head cannot repeat in T. R supplies a tail-contained word valued 0 or 4, or a nonlinear same-head tail-contained word valued 1. Thus all values are absorbed. This proves the criterion without a finite test assumption.

\par\smallskip\noindent\textit{Two witnesses and uniform compression.}\quad 
For a valid nonlinear absorption choose the original K witness a. If the head repeats, a alone satisfies K and R for the target. Otherwise choose one original R witness b. The absorption \(w\leq a+b\) (or w\ensuremath{\leq}a) is valid by the criterion; adding the other original summands recovers \(w\leq P\). Thus at most two original right words are needed, regardless of the size of P.

Reserve all heads of the at most three words (left and two right), at most three variables. Group each remaining variable by its nonempty incidence pattern across the three tails, giving at most seven groups. Replace each group by a fresh variable and keep reserved heads separate. Normalize by the two word laws. The resulting valid absorption has at most ten variables and words of length at most eleven. Validity follows from substitution, not from a bounded test. Substitute the nonempty product of the members of each group for that fresh variable. Identical tail incidence, reserved heads, and retained head repetition recover all original normal words, using only the two word laws. No empty substitutions occur.

Take Sigma to be \(xyz=xzy\), \(xyy=xy\), together with every valid normal absorption over a fixed ten-variable alphabet, word length at most eleven, and at most two right words. K and R specify this finite family exactly. The compression and inverse derive the two-witness absorption, and adding the other summands gives every valid absorption. Distributivity and summandwise absorption in both directions give every valid identity. Consequently Sigma is an unrestricted finite relative basis for 8729. Opposite reversal and independently verified two-way HSP support transfers.

The failed shortcut \(xy\leq x\) illustrates why a literal head cannot supply R: at \(x=2\), \(y=4\) the left is 1 and the right is 2, below 1 in the additive order. This distinction is explicitly retained in the proof and replay.

\subsection{Four finite bases from initial occurrence witnesses}\label{case:177}
The sources are 6837, 14018, 14033 and 14093. Reversing words gives the corresponding bases for 8790, 10999, 11013 and 11298. Use the constant-free additively idempotent semiring axioms. Write u \ensuremath{\leq} P for \(u+P=P\). A word is nonlinear here if its length is at least two.

\par\smallskip\noindent\textit{Raw tables and absorption tests.}\quad 
All four additions are maximum on 3 \ensuremath{<} 2 \ensuremath{<} 1 \ensuremath{<} 0 \ensuremath{<} 4. Multiplication rows in catalogue order are as follows.

\noindent Source\enspace\textbar\enspace 0 row\enspace\textbar\enspace 1 row\enspace\textbar\enspace 2 row\enspace\textbar\enspace 3 row\enspace\textbar\enspace 4 row\par
\noindent 6837\enspace\textbar\enspace 01134\enspace\textbar\enspace 33333\enspace\textbar\enspace 33333\enspace\textbar\enspace 33333\enspace\textbar\enspace 44444\par
\noindent 14018\enspace\textbar\enspace 44444\enspace\textbar\enspace 01134\enspace\textbar\enspace 01134\enspace\textbar\enspace 33333\enspace\textbar\enspace 44444\par
\noindent 14033\enspace\textbar\enspace 44444\enspace\textbar\enspace 01134\enspace\textbar\enspace 02234\enspace\textbar\enspace 33333\enspace\textbar\enspace 44444\par
\noindent 14093\enspace\textbar\enspace 44444\enspace\textbar\enspace 01234\enspace\textbar\enspace 33333\enspace\textbar\enspace 33333\enspace\textbar\enspace 44444\par
For 6837 the neutral prefix consists of zeros. An all-zero word gives

0. A first 4 gives 4 and a first 3 gives 3. A first 1 or 2 gives 3 unless final. At the final position it gives 1 in a nonlinear word; a literal retains its value. For 14018 and 14033 the neutral prefix consists of 1 and 2. In an all-neutral nonlinear word 14018 gives 1, and 14033 gives the head value. A literal retains its value. The first nonneutral 3 or 4 gives itself. The first nonneutral 0 gives 0 if final and 4 otherwise. For 14093 only 1 is neutral. An all-neutral word gives 1. The first 0 gives 0 if final and 4 otherwise, a 4 or 3 gives itself, and a 2 gives 2 if final and 3 otherwise. These statements follow by multiplying the displayed rows from the left. For instance, a neutral prefix in 14033 does not change its first factor's value when extended by a further neutral factor.

Fix a target w and an original upper polynomial P. Write C for the support of w. For each x in C let S be the set of letters before the first x in w. An E-witness is an original word whose prefix before its first x uses only S. It must contain x. It is marked if that occurrence has a successor and unmarked otherwise. A first occurrence in w uses the same terminology. A contained word has support within C. The following tests are necessary and sufficient.

\textbullet{} In 6837 a literal requires its literal. A nonlinear w needs a contained word and an E-witness at every first occurrence. For a simple final s it additionally needs a word contained in C minus s, or a nonlinear unmarked E-witness for s.

\textbullet{} In 14018 and 14033 each position of first occurrence requires an E-witness, marked whenever the target occurrence is marked. In 14018 a nonlinear target additionally needs a nonlinear contained word and a literal target needs any contained word. In 14033 every target requires a contained word starting with its own head.

\textbullet{} In 14093 a literal requires its literal. A nonlinear target needs a contained word and an E-witness for every first occurrence, marked at every marked target occurrence. A simple final s also needs a word contained in C minus s or an unmarked E-witness for s. This unmarked witness may be a literal.

For the E-tests give S a neutral value (0 in 6837, 1 elsewhere), x the value 4, and every other variable 3. Precisely an E-witness can reach 4. For a required marked E-witness put \(x=0\) instead in the last three sources; an unmarked occurrence gives only 0. For ordinary containment use neutral target letters and outsiders 3. For 14018 on a nonlinear target put all target variables at 2: only a nonlinear contained word can reach the target value 1. For 14033 put the target head at 1, the rest of C at 2 and outsiders at 3. Reaching 1 requires a contained word with that head. For 6837 or 14093 put a simple final s at 2, earlier target letters at the neutral value and outsiders at 3. To reach the target threshold (1 or 2 respectively) gives exactly the extra alternatives in the tests, including the distinction for a literal s. For a literal target in these two sources give it 2 and all other variables 3. Only its literal reaches 2.

Conversely, use a contained witness for an all-neutral target, with the specified length in 14018 and head in 14033. A first 4 is witnessed by E. A nonfinal first 0 in the last three sources uses the marked witness to reach 4, while at a final first 0 either kind of witness attains at least 0. A first weak value (1 or 2 in 6837, 2 in 14093) giving more than the bottom must be final, with only neutral values before it. An omitting contained witness gives the neutral value. The stipulated final witness gives the required weak value; in 6837 it must be nonlinear when \(s=2\). The bottom 3 is automatically absorbed. This proves all tests on arbitrary finite alphabets and polynomials.

\par\smallskip\noindent\textit{Word identities and declared bases.}\quad 
All bases include W1, \(xxy=xy\), and W2, \(xyx=xyy\). Read off their soundness from the word formulas. Their common normal form lists the distinct letters in first-occurrence order. If the final letter is repeated, append the last listed letter once more; otherwise keep the list. To reduce a repeated nonfinal occurrence, use W1 for adjacent repeats and \(xYxZ=xYYZ=xYZ\) when separated, for nonempty Y, Z. After removing internal repeats, a repeated final x in p x Y x with nonempty Y becomes p x Y Y by W2. In the second copy of Y each letter but its last is now an internal repeat and can be deleted. The result has the claimed repeated last letter. This also covers adjacent final repeats and powers. External contexts may be empty; a variable in an identity always denotes a nonempty term.

The other laws for each source are listed below. Every displayed absorption is sound by its corresponding test (or by its finite operation table on its displayed variables). All are identities in the semiring signature under the convention for \ensuremath{\leq}.

\noindent Source\enspace\textbar\enspace Additional laws\par
\noindent 6837\enspace\textbar\enspace xx \ensuremath{\leq} x, xyz \ensuremath{\leq} yz+xt, xyz \ensuremath{\leq} xz+xyt, xyz \ensuremath{\leq} xyzz+y, xyz \ensuremath{\leq} xyzz+yz\par
\noindent 14018\enspace\textbar\enspace x \ensuremath{\leq} xx, xyz \ensuremath{\leq} yz+xt, xyz \ensuremath{\leq} xz+xyt, xy \ensuremath{\leq} xx+xyt, xy \ensuremath{\leq} xx+y\par
\noindent 14033\enspace\textbar\enspace x \ensuremath{\leq} xx, hh \ensuremath{\leq} h+hq, huh \ensuremath{\leq} hq+ut+h, huhb \ensuremath{\leq} hq+ut+hb, abcd \ensuremath{\leq} acd+abd, abc \ensuremath{\leq} aa+bc, xy \ensuremath{\leq} xx+y\par
\noindent 14093\enspace\textbar\enspace xy \ensuremath{\leq} y+xt, xyz \ensuremath{\leq} xz+xyt, xy \ensuremath{\leq} xx+xyt\par
\par\smallskip\noindent\textit{Constructing prefixes in 6837, 14018 and 14093.}\quad 
Normalize the target only. Keep the original words as witnesses. Write p for the list of first occurrences processed so far. For 6837 maintain p \(u^{2}\) \ensuremath{\leq} P for every original u; initially \(u^{2}\) \ensuremath{\leq} u. For 14018 maintain p u \ensuremath{\leq} P for every nonlinear original. For 14093 maintain p u \ensuremath{\leq} P for every original. In the following, write U for the maintained suffix \(u^{2}\) or u as appropriate.

Start at the first marked target letter x. In 6837 use the square of an E-witness, available by xx \ensuremath{\leq} x. It starts with x and has a nonempty following factor. In 14018 and 14093 the original marked witness has this property. Denote it by xt. For 6837 and 14018, \(U=YZ\) has two nonempty factors. The law xyz \ensuremath{\leq} yz+xt yields xU \ensuremath{\leq} P. For 14093 its law xy \ensuremath{\leq} y+xt yields xU \ensuremath{\leq} P even for a literal U.

For a later marked letter x choose a witness as at the start, squaring in 6837, and write it as d x t where t is nonempty and the support of d lies in the already processed p. The invariant gives p d x t \ensuremath{\leq} P. By W1, W2 this equals p x t; internal repeats disappear with its nonempty tail retained. From p U \ensuremath{\leq} P and p x t \ensuremath{\leq} P use xyz \ensuremath{\leq} xz+xyt to get p x U \ensuremath{\leq} P. The procedure handles every first occurrence preceding a simple final, or every first occurrence if the target's last letter is repeated.

For a target with repeated final, use an original contained u (of length at least two in 14018). After processing all distinct letters the maintained p U normalizes to the desired saturated word. For a simple final s the preceding prefix p is nonempty since the target is nonlinear. The following steps finish each case.

In 6837 the same procedure can also process s, since squaring its E-witness gives a marked witness for the construction, even if the original E-witness was unmarked. Obtain the saturated pss from a contained original. If the final test supplies a word u contained in c(p), use xyz \ensuremath{\leq} xyzz+y with \(x=p\), \(y=u\), \(z=s\). Its left side reduces to ps and its first right side to pss. Otherwise use the original nonlinear witness \(u=t\) s, with t nonempty and c(t) within c(p). The law xyz \ensuremath{\leq} xyzz+yz with \(x=p\), \(y=t\), \(z=s\) gives ps from pss and u.

In 14018 choose a nonlinear contained r. If it contains s, p r normalizes to ps or pss. In the latter case x \ensuremath{\leq} xx at \(x=ps\) gives ps \ensuremath{\leq} pss by W1, W2. If r omits s, p r gives \(p^{2}\) in normal form. An unmarked nonlinear E-witness gives ps directly by the invariant. A marked E-witness gives p s t with nonempty t; use xy \ensuremath{\leq} xx+xyt. A literal E-witness s instead gives ps by xy \ensuremath{\leq} xx+y. All referenced witnesses are original or were derived by the invariant.

In 14093 an unmarked E-witness u (including s) gives ps directly by p u. Otherwise there is a marked witness, giving p s t. The extra final condition then supplies a word contained in c(p). Its maintained term gives \(p^{2}\), so xy \ensuremath{\leq} xx+xyt gives ps.

\par\smallskip\noindent\textit{Retaining a contained same-head word in 14033.}\quad 
Fix an original contained word r with target head h. Throughout the construction derive p r and p u r for every nonempty proper prefix u of every original word. The initial marked position h has an original marked witness hq. If \(r=h\), derive hh \ensuremath{\leq} h+hq. If \(r=h\) b, h \(r=r\) by W1. For a proper prefix u take its original extension ut with t nonempty. At \(r=h\) use huh \ensuremath{\leq} hq+ut+h, and at \(r=h\) b use huhb \ensuremath{\leq} hq+ut+hb. Thus the invariant holds at \(p=h\), including when r is a literal.

To process a marked position x take an original d x t with c(d) in c(p) and t nonempty. Its proper prefix d x gives p d x \(r=p\) x r by W1, W2. For each maintained prefix u, the law abcd \ensuremath{\leq} acd+abd applied to p, x, u, r derives p x u r from p u r and p x r. This extends the invariant without requiring that the rest of the witness lie in the target support. If the target ends with a repeat, p r at the end gives exactly the target in normal form.

For a simple final s stop with prefix p. If r contains s, its maintained p r normalizes to ps or pss. Apply x \ensuremath{\leq} xx at ps if needed. If r omits s, the maintained term gives \(p^{2}\). If there is a marked E-witness for s, use its proper prefix as above to get p s r, which normalizes to pss, and recover ps. Otherwise there is an original d s with c(d) within c(p). With nonempty d apply abc \ensuremath{\leq} aa+bc at p, d, s; p d s normalizes to ps. When d is empty use xy \ensuremath{\leq} xx+y at p, s and the original literal s.

\par\smallskip\noindent\textit{Completeness and opposite types.}\quad 
For 6837 and 14093 a valid literal target is already present. For 14018 and 14033 its contained original is a positive power of that letter, and W1 together with x \ensuremath{\leq} xx recovers the literal from every such power. This closes the remaining cases. Distributivity writes every term as a nonempty sum of words. Deriving all valid word absorptions and then using both sides' summands proves every valid identity. The bases have respectively seven, seven, nine and five laws in addition to the semiring axioms. There is no restriction on the original alphabet, substitution, word length or number of summands. Opposite multiplication reverses the words in each law.

\subsection{Finite identity bases for 7054, 7057, 7129 and their opposites}\label{case:178}
Work in the constant-free signature of additively idempotent semirings. An absorption (w \ensuremath{\leq} P) abbreviates \((w+P=P)\). Every word and every substituting word is nonempty. For each of the three tables take the splitting equation

{\small\par\smallskip
\CaseLine{\(xy = xxy + xyy\)}
\smallskip}

and all its table-valid absorptions on a fixed alphabet of 737 letters with word length at most 2951 and at most two right summands. These three finite sets, each together with the semiring axioms, are complete. Reversing every word gives bases for 8383, 8386 and 8538 respectively. The size bounds describe finite presentations without listing their individual members.

\par\smallskip\noindent\textit{Word evaluation.}\quad 
In catalogue order the multiplication rows of 7054 are 01200, 22211, 22222, 01234, 01244. In 7057 the fourth row is 01244 and all other rows agree. The additive order has chains \(2<1<0<4\) and \(2<3<4\), with \(0+3=1+3=4\). An input 2 or two occurrences of value 1 give the least element 2. With one 1 and other inputs in K=\{0,\allowbreak{}3,\allowbreak{}4\}, the output is 1 exactly when all later entries are in \{3,\allowbreak{}4\}, and 2 otherwise. On K, an occurrence of 0 gives 0. On \{3,\allowbreak{}4\}, 7054 gives the maximum in the order \(3<4\). For 7057 a literal retains its value and any longer word on \{3,\allowbreak{}4\} gives 4.

The multiplication rows of 7129 are 01220, 22221, 22222, 33333, 01224, with additive order \(3<2<1<0<4\). A word headed in value 3 gives 3. With any other head, a 2 or a 3 anywhere gives 2. On \{0,\allowbreak{}1,\allowbreak{}4\}, two 1's give 2, and a unique 1 gives 1 exactly when all later entries are 4, and 2 otherwise. On \{0,\allowbreak{}4\} the output is 0 if 0 occurs, and 4 otherwise. Multiplying these cases by each of the five possible next entries proves the formulas by induction, including literals. The splitting equation holds on all 25 pairs in each table. On a surviving unique 1, the summand that doubles the other argument retains this 1.

\par\smallskip\noindent\textit{Absorption by the original right words.}\quad 
Let C be the support of w, and U the right words supported in C. A simple letter x occurs exactly once in w. Write B for the set of letters strictly after that occurrence. For each x let E be the union of the supports of all words in U omitting x. Each occurrence of U, B or E in the following test refers to the given absorption.

\textbullet{} For 7054 the union of supports of U must equal C. For every simple x a word v of U must have x exactly once and have all its letters after x in B union E.

\textbullet{} For 7057 a literal target must be a summand. A longer target requires a word of U of length at least two. For each simple x there must be a word of U of length at least two omitting x, or a word v of U containing x exactly once with subsequent letters in B union E.

\textbullet{} For 7129 there must be a right word with the same head as w and a word in U. The same-head word may have outside letters. For each simple x one word of U must omit x or contain x exactly once with subsequent letters in B.

For necessity in 7054 and 7057 give outsiders value 2. To force coverage in 7054, put a selected target letter at 4 and the others at 3. To force a longer contained word in 7057, put all target letters at 3. A literal target in either algebra requires the literal, by giving its variable 1. For the simple-letter test put x at 1, the letters of B union E at 3 and the rest of C at 0. The target gives 1. Every omission in 7054 gives 3. In 7057 the same holds when all omissions are literals. Every other word which can absorb the 1 must realize the specified simple witness. All losing values lie in \{2,\allowbreak{}3\}, which is closed under joins.

Conversely, on K the coverage test in 7054 supplies a contained witness for a 4 in an all-high word, a 0 when one occurs, and a word of value 3 on all 3's. In 7057 a longer contained word suffices on K for a longer target. With a surviving 1, all other target letters are in K and B takes high values. A longer omission in 7057 gives 0 or 4. Otherwise use the selected simple witness. If its suffix contains a 0, this letter is in E, and the corresponding omission gives 0. With no such 0 the witness gives 1. The least target value requires no further condition.

For 7129 assign the target head 2 and all other letters 3 to force a right word with that head. Assign C value 4 and outsiders 3 to force containment. For a simple x assign it 1, B value 4, the rest of C value 0, and outsiders 3. A word reaching 1 or higher must be a contained omission or the specified simple witness. Conversely the global head witnesses any target of value 2, the contained word witnesses an evaluation on \{0,\allowbreak{}4\}, and the omission or the suffix test witnesses a surviving 1. This proves all three tests for arbitrary alphabets.

\par\smallskip\noindent\textit{Reduction to absorptions with at most two upper words.}\quad 
By distributivity it suffices to derive each valid word absorption. Expand both sides by the splitting equation. When a word has two or more simple letters, split it into two nonempty consecutive factors, each containing at least one simple letter, and double one factor in each branch. Every branch decreases the number of simple letters. We may therefore work with a single lower word and upper words each having at most one simple letter. The absorption after expansion follows semantically from the initial one, so the exact tests apply again to the new words. All words introduced below have at most one simple letter.

In 7129 keep an original head witness. A second original suffices for containment and the simple condition together, since the witness for a simple x is itself contained. Without a simple x choose any contained original. The two together absorb w.

For 7054 and 7057 a literal target is immediate. For a longer target let G be a word listing each letter of C twice, in any fixed order. In 7054 take each needed contained original u and derive its fully doubled support from u. The lower word with doubled support union of two already derived such words is absorbed by their sum, by coverage and the absence of simple letters. Iteration gives G. In 7057 derive G in one step from any longer contained original. Once G is available, a target having no simple letter is absorbed by G. For a simple x in 7057 having a longer contained omission d, the test gives w \ensuremath{\leq} G+d.

In the remaining cases choose the simple witness v from the test, and let S be its suffix support. For any subset T of C excluding x write W(T) for a word listing every letter of C other than x twice, then x, then every letter of T twice. Apart from the already discharged literal case, C contains another letter whenever w has a simple x. Thus these words are longer and have precisely one simple letter. First derive W(B union S) from G+v. The original v has a suffix contained in B union S. The upper G supplies coverage or the required longer contained word.

Remove one letter y of S outside B at a time. Such a y belongs to E, so choose a contained omission d containing y. With T the current suffix support, derive

{\small\par\smallskip
\CaseLine{\(W(T minus \{y\}) \leq  W(T) + d.\)}
\smallskip}

Both sides retain full target support through W(T). In the test for the simple x, the extra suffix letter of W(T) belongs to the omission d. Equivalently, on an assignment where the lower word gives 1, either y is high and W(T) gives 1, or y is 0 and d gives 0. At the end W(B) absorbs w since its suffix support agrees with B. No collection of original omissions of bounded cardinality was assumed. Each step uses a single omission while retaining all of C on the lower side.

\par\smallskip\noindent\textit{Bounded derivations for every local step.}\quad 
We give the compression argument because the number of added or removed suffix letters can be arbitrarily large. In each of the at most three words of a local step protect the head, its unique simple position if present, and the immediate successor of that position if present. At most two contiguous unprotected runs remain in any one word. A letter in a run is repeated in that word. Permuting two adjacent distinct repeated letters within a run, and deleting one of a consecutive equal pair while at least two occurrences remain overall, preserve the word function. Indeed support, head, the length class (literal or longer), the simple letters and the suffix support of each simple letter are preserved. The word formulas show that these data determine the value, also when several valuation entries agree.

These edits themselves are instances of bounded equations. For a swap mark the head, the swapped positions and one other occurrence of each swapped letter. For a deletion mark the head, the pair, and a third occurrence of its letter. Compress each nonempty unmarked interval into its own fresh formal letter. The formal edit preserves the same data for every formal simple letter, including each interval symbol. Both orientations follow from valid one-upper absorptions of length at most 19. Substitute the actual intervals to obtain the edit.

Sort each run and truncate every multiplicity there to two. Keep distinct all symbols of protected positions across all three words, at most nine. Group the others by their vector of counts 0, 1 or 2 in the at most six runs. At most \(3^{6}\)-1 nonzero vectors occur. Identify each group to one symbol simultaneously in the valid absorption. Sort again and reduce each group symbol to the common count of its members in that run. If this count is one, each member occurs in another run of the same word, because it is repeated and never protected. Thus the deletions retain at least two occurrences overall. The protected positions and the simple letters persist throughout. The result is a valid absorption on at most 737 letters and has length at most \(3+4(737)=2951\) in each word. It belongs to the indicated finite set.

Conversely substitute the product of all distinct members of a group for its symbol. Its common multiplicity in each run, including two, is thereby restored for every member. Sorting the runs by the bounded edits recovers the original truncated words, and the inverted edits recover the original words. This derives every local step from the finite basis. Finally add the lower summands and apply each absorption in both required directions for a general polynomial identity. The construction is uniform in the number of variables, the lengths and the number of upper summands.

\subsection{末字母见证与内部首字母保护：五条完整基}\label{case:179}
\par\smallskip\noindent\textit{1. 构造与词值.}\quad 
取 \(A=\{z, a, b, d, t\}\)，依次对应目录标号 \((0,\allowbreak{}1,\allowbreak{}2,\allowbreak{}3,\allowbreak{}4)\)。加法为偏序

\CaseDisplay{d<a<z<t,\qquad d<b<z<t}

上的并，其中 \(a, b\) 不可比。右乘 \(z, a\) 为常值 \(a\) 映射，右乘 \(b, d\) 为常值 \(d\) 映射，右乘 \(t\) 为

\CaseDisplay{R=(z, a, z, a, t).}

三个映射保持并、幂等；它们与右因子乘法、加法的兼容式 \(R_{xy}=R_yR_x\)、\(R_{x+y}=R_x\vee R_y\) 可逐类核验，所以定义半环。完整乘法表为

\CaseDisplay{\begin{array}{c|ccccc}
\cdot&z&a&b&d&t\\\hline
z&a&a&d&d&z\\
a&a&a&d&d&a\\
b&a&a&d&d&z\\
d&a&a&d&d&a\\
t&a&a&d&d&t
\end{array}}

对于非一次词：末字母取 \(z, a\) 时全词取 \(a\)，取 \(b, d\) 时全词取 \(d\)。若末字母取 \(t\)，则尾部全取 \(t\) 时全词为 \(R(h(w))\)；尾部还含其他取值时全词取 \(a\)。因而非一次词取值只可能为 \(d, a, z, t\)，并且：

\textbullet{} 取 \(t\) 当且仅当全部变量取 \(t\)；

\textbullet{} 取 \(z\) 当且仅当首字母取 \(b\) 或 \(z\)、尾部全部取 \(t\)；

\textbullet{} 取 \(a\) 时，末字母必取 \(a, z, t\)；

\textbullet{} 其余情形取最小元 \(d\)。

\par\smallskip\noindent\textit{2. 正规形.}\quad 
以下三条词恒等式在 \(A\) 上成立：

\CaseDisplay{xy^2\approx xy,\qquad xyx\approx y^2x,\qquad xyzt\approx xzyt.\quad\textup{(\ensuremath{W})}}

前两式推出 \(x^3=x^2\) 和 \(xyzy=xzy\)，从而尾部只保留每个变量的最后一次出现。若首字母 \(h\) 不在尾部重复，词由内容 \(X\)、首字母 \(h\) 和末字母 \(\ell\) 决定：第三式允许相邻交换任何两个内部字母，故内部次序无关。

若首字母在尾部重复，则词为平方形式，由内容 \(X\) 和末字母 \(\ell\) 决定。详细地，前两式先将它写成 \(r^2\)，其中 \(r\) 无重复且末字母为 \(\ell\)。写 \(r=q\ell\)；\(q\) 非空时，尾部去重及第二式给出

\CaseDisplay{r^2\approx q^2\ell\approx \ell q\ell.}

首尾已都是 \(\ell\)，第三式允许任意重排中间的 \(q\)。当 \(X=\{\ell\}\) 时直接为 \(\ell^2\)。这说明所有平方词只依赖内容与末字母。

记 \(D(w)\) 为尾部变量集合，\(h(w)\)、\(\ell(w)\) 分别为首、末字母。这里``首字母不重复''表示 \(h(w)\notin D(w)\)，不表示词的其他变量不能重复。

\par\smallskip\noindent\textit{3. 精确吸收判据.}\quad 
令 \(F\) 为非空右词集。一次吸收 \(x\le F\) 有效，当且仅当同名一次词 \(x\) 本身属于 \(F\)。

定理1。 对非一次词 \(w\)，令 \(X=c(w)\)、\(h=h(w)\)、\(\ell=\ell(w)\)。则 \(A\models w\le F\) 当且仅当：

1. 存在内部右项 \(u\)，\(c(u)\subseteq X\)；若 \(h\notin D(w)\)，还要求 \(h\notin D(u)\)；

2. 存在右项 \(v\)，满足 \(\ell(v)=\ell\)。

这里第二条件仅要求相同末字母，不要求其他变量的后缀见证。

必要性。 对一次左词，若无同名一次项，赋 \(x=b\)、其他变量 \(d\)。非一次右词一律取 \(d\)，其他一次项也取 \(d\)，不能吸收 \(b\)。

若没有内部右项，赋 \(X\) 内变量 \(t\)、外部变量 \(d\)。左值为 \(t\)，每个右词因含取值 \(d\) 的变量而不取 \(t\)，全部不超过 \(z\)，矛盾。

如果内部项存在，但 \(h\notin D(w)\) 且每个内部项尾部都含 \(h\)，赋 \(h=b\)、其余 \(X\) 内变量 \(t\)、外部变量 \(d\)。左值为 \(z\)；内部右词含尾部 \(b\)，外部右词含 \(d\)，结合末字母分类，它们都取 \(a\) 或 \(d\)，右和不超过 \(a\)。

若无末字母见证，赋 \(\ell=a\)、其余变量 \(d\)。左非一次词取 \(a\)，右词末字母全部取 \(d\)，所以右词全部取 \(d\)。必要性得证。

充分性。 左值为最小元 \(d\) 时自明；为 \(t\) 时，全部左变量取 \(t\)，任一内部右项也取 \(t\)。

左值为 \(a\) 时，末字母 \(\ell\) 取 \(a, z, t\)。具有同一末字母的右项 \(v\)，若一次则取 \(a, z, t\)；若非一次则取 \(a, z, t\) 中的值，均不小于 \(a\)。

左值为 \(z\) 时，首字母 \(h\) 取 \(b\) 或 \(z\)、全部尾部变量取 \(t\)，所以 \(h\) 不在尾部重复。取条件1的受保护内部项 \(u\)。若 \(u\) 不含 \(h\)，则取 \(t\)；若含 \(h\) 且非一次，则只在首字母含 \(h\)，取 \(z\)。剩下 \(u=h\) 为一次项；若 \(h=z\) 已足够，若 \(h=b\)，条件2的末字母见证因末字母取 \(t\)，必取 \(a, z\) 或 \(t\)。因此

\CaseDisplay{u+v\ge b+a=z.}

全部取值情形完成，判据成立。

\par\smallskip\noindent\textit{4. 五条完整基.}\quad 
定理2。 \(A\) 的完整附加基由三条 \((W)\) 以及

\CaseDisplay{xy\le y,\qquad xy\le x+zy\quad\textup{(\ensuremath{B})}}

组成。

两条不等式由定理1逐式检验。由 \(xy\le y\)，在非空前后段中插入任意中段会降低词值：\(xyz\le xz\)；在开头添加任意非空前段也会降低词值。这些是句法推导，不依赖具体赋值。

给定有效非一次吸收 \(w\le F\)，取判据中的内部见证 \(u\) 和末字母见证 \(v\)。若 \(v=\ell\) 是一次项，由 \(u\ell\le\ell\) 得 \(u\ell\le F\)；若 \(v=p\ell\)、\(p\) 非空，由第二条 \((B)\) 得

\CaseDisplay{u\ell\le u+p\ell\le F.\quad\textup{(\ensuremath{1})}}

若 \(w\) 的首字母重复，选一个列出 \(X\) 中全部变量的非空词 \(q\)，使其首字母等于 \(h(u)\)。由前缀收缩，

\CaseDisplay{qu\ell\le u\ell\le F.}

词 \(qu\ell\) 的内容为 \(X\)、末字母为 \(\ell\)，首字母在后面的 \(u\) 中重复。由第2节的平方正规形，它等于 \(w\)。

若 \(w\) 的首字母 \(h\) 不重复，则 \(\ell\ne h\)。由内部保护，\(u\) 或者不含 \(h\)，或者只在首字母含 \(h\)。前一情形在 \(u\ell\) 前添加 \(h\)，后一情形直接用 \(u\ell\)，都得到一个不超过 \(F\)、首字母为 \(h\) 且不重复它、末字母为 \(\ell\) 的非一次词。

对于 \(X\) 中尚未出现的变量，利用 \(xyz\le xz\) 将其逐个插入首字母与末字母之间。所得词不超过 \(F\)，内容恰为 \(X\)、首字母为 \(h\) 且不重复、末字母为 \(\ell\)。由 \((W)\) 去掉尾部重复并重排内部字母，它等于 \(w\)。故 \(w\le F\) 可由五式推出。

一次有效吸收直接来自同名右项。一般项用分配律展开成非空词和，再处理双向逐词吸收，得到任意变量数完备性。所有代入均为非空词。证毕。

\par\smallskip\noindent\textit{5. 范围.}\quad 
此处把有序后缀覆盖进一步缩减为末字母见证，是因为右平移的特殊结构额外满足中间换位及收缩式。不能把这一简化用于不满足这些恒等式的来源。

乘法对偶由反转所有词处理。具体五元素类型的有限基传递另作独立双向HSP核验；小变量全多项式诊断不替代上述证明。该定理本身不宣称任意真子簇有限基。

\subsection{Four finite bases from the successor of a simple letter}\label{case:180}
The sources are 7303 and 7645, with opposites 8081 and 7647, respectively. Use their own valid absorptions on 737 variables, with at most two right words and every word of length at most 2951, together with \(xy=xxy+xyy\) and the ai-semiring axioms. This is a finite presentation for each of the four types. It specifies the families by a terminating finite table test.

\par\smallskip\noindent\textit{Word values and absorption.}\quad 
The additive order in both sources is \(3<2<1<0<4\). Call 0,\allowbreak{}4 high. An occurrence of 3, or at least two occurrences taking values in \{1,\allowbreak{}2\}, makes the word 3. A unique 1 makes the word 1, in the absence of 2,\allowbreak{}3. On high inputs 7303 returns the head and 7645 the maximum. With a single 2 and all other values high the two rows differ. For 7303 a final 2 gives 2, a 2 immediately followed by 4 gives 1, and a 2 immediately followed by 0 gives 3. For 7645 the result is 1 if any later factor is 4, and otherwise 2. In label order the multiplication rows of 7303 are 01230,\allowbreak{}13331,\allowbreak{}33331,\allowbreak{}33333,\allowbreak{}41234; those of 7645 are 01234,\allowbreak{}13331,\allowbreak{}23331,\allowbreak{}33333,\allowbreak{}41234. The word formulas follow by appending one factor, and include length one.

Fix w and P. Write C for the support of w and U for those original words of P supported in C. A variable is simple if it occurs exactly once in the word at issue. The following conditions are necessary and sufficient.

\textbullet{} 7303 requires some member of U with the head of w. For each simple x of w, a member of U must omit x or contain x simply with exactly the same immediate successor as in w. At a final x the latter means a simple final x.

\textbullet{} 7645 requires that the supports in U cover C. For each simple x of w, an omitting member of U suffices. Otherwise some member has a simple x, and for every letter y after x in w, a member with simple x has y after that occurrence. The witnesses may differ.

Give all outside variables 3. On high inputs, assigning the target head 4 and the other letters 0 forces a contained same-head word in 7303. In 7645 assigning any missing cover letter 4 and the others 0 forces full coverage. Give x the value 1 and the remaining target letters 0 to require an omission or simple occurrence. For a final x in 7303 assign \(x=2\) with all remaining target letters 0; any nonfinal x now gives 3. For a nonfinal simple x in 7303 give its immediate successor y the value 4, x the value 2, and the rest

0. Without omission a right word attaining 1 must have that very successor. For 7645 use the same assignment for each y anywhere after x. Conversely, when the target is high the head or cover condition supplies the attained level. With a unique special value the witnessing contained word attains at least the target by the formulas, including 7645 at value 2 if no suffix letter is high 4. Specifically an omitting word gives a high value, a simple witness always reaches 2 in 7645, and an indicated suffix 4 raises it to 1. All remaining targets have the least value.

\par\smallskip\noindent\textit{Splitting and bounded local derivations.}\quad 
The identity \(xy=xxy+xyy\) holds in each table. Equivalently, at a unique special operand its own duplication gives the bottom while duplicating the high operand preserves the product; on high inputs both duplicates preserve it. With two specials or a bottom all three sides give bottom. If a word has two or more simple variables, split it into two consecutive nonempty blocks each containing a different simple variable and apply this identity. In each new word at least one simple variable becomes repeated and none becomes simple. Finite induction expands both sides of any proposed identity into sums of words each having at most one simple variable.

In 7303 choose one original contained same-head word and, when needed, one original successor or omission witness. The same test shows that the resulting absorption with at most two upper words is valid. All its words already have at most one simple variable.

In 7645 first derive a saturated word G of support C, where saturated means every letter occurs at least twice. For each u in U derive u squared below u. Join these supports one at a time by deriving the square of the concatenation of two current saturated words below their sum. The cover criterion for these saturated lower words is immediate. These steps have at most two upper words. If the target is saturated, G suffices. For its single simple letter x, an original omitting word together with G suffices when one exists. Otherwise take contained words each with simple x, including witnesses for all letters needed after x. Merge two at a time into a single word with simple x. Take their support union, put two copies of each letter occurring after x in either word after x, and put two copies of each remaining letter other than x before x. Either block may be absent. The criterion for the merged word holds with the two inputs: each suffix letter has its own simple-x witness and their supports cover. The merge keeps every already acquired suffix letter. At the end G and the resulting simple-x word suffice for w. This includes an empty target suffix, for which one simple witness suffices. Each step uses only words with at most one simple variable, at most two upper words, and actual locally valid absorptions.

\par\smallskip\noindent\textit{Compression that retains the order around the unique simple variable.}\quad 
In each such word protect its first position and, if there is a simple x, its position and the immediately following position when present. There are at most three protected positions and at most two unprotected runs. Inside each run we can sort and cap multiplicities at two. The following argument supplies literal bounded identities for the needed edits.

To swap successive different b, c in a run, its predecessor a exists and is repeated. Indeed the head is protected and so is the successor of a simple letter. Both b, c are repeated. Mark abc and, if needed, further occurrences of a, b, c witnessing their repetitions. Use a fresh formal variable for each intervening or outside nonempty unmarked block, separately for each block. At most nine marked positions and ten blocks occur, giving two formal words of length at most 19. Both have identical head, supports, and simple variables. The only letters whose immediate successor might change are repeated a, b, c. For every simple letter the set of later letters also stays the same. The word formulas therefore give a genuine word identity on each table. Restoring the marked letters and the actual blocks is an ordinary nonempty substitution proving the swap.

To shorten a block of consecutive equal b's in a run, keep at least one copy in that run and two overall. Delete the second of two adjacent copies and mark those two and another surviving occurrence of b. Such an occurrence exists when deletion retains two overall. Compress all remaining blocks as above. Formal b is repeated on both sides. The predecessor of the deleted letter is b, so no simple letter loses its immediate successor. Since the two copies were adjacent no simple letter loses b from its suffix set. The head and support agree. This gives another valid word identity of length at most 19. Both orientations of each short word identity belong to the declared family. These edits work also when the same variable occurs at protected positions and in a run.

For a valid local absorption first sort its runs and cap each letter at two in each run. Across its at most three words reserve all letters having a protected position in any word, at most nine. Each remaining variable is described by the vector of its counts 0,\allowbreak{}1,\allowbreak{}2 in at most six runs. At most \(3^{6}\)-1 nonzero vectors occur. Identify all variables of each such group to one new letter. Validity survives substitution. In each run sort and reduce a group letter to the COMMON count of its members in that run. For a common count of one the group letter can be reduced to one there: each original group member occurring in that word is repeated overall, and its occurrence in another run remains. For a common count of two retain both. Every individual deletion keeps at least two copies overall. These edits also preserve the protected positions and their immediate successors.

The resulting absorption has at most \(9+3^{6}-1=737\) letters and each word has length at most 3+4 times 737, namely 2951. It is one of the proposed axioms. Substitute the nonempty list of all original group members back for each group letter. Each member then occurs exactly its old count in each run. Sorting runs restores the prior words with their original protected positions, including the successor of the unique simple letter. Reverse the initial edits as needed. This derives every local step, hence every valid absorption after splitting. Taking the summand absorptions in both directions proves completeness for arbitrarily many variables and summands. Reversing all words transfers these bases to the stated opposite tables.

\subsection{A finite basis for 7319 retaining both word endpoints}\label{case:181}
Signature (+, multiplication), without constants, relative to the additively idempotent semiring axioms. The finite family is effectively specified; no small enumeration or minimality is claimed.

\par\smallskip\noindent\textit{1. Semantics and exact absorption.}\quad 
Addition rows in element order 0,\allowbreak{}1,\allowbreak{}2,\allowbreak{}3,\allowbreak{}4 are 00004,\allowbreak{}01114,\allowbreak{}01224,\allowbreak{}01234,\allowbreak{}44444, so the additive order is \(3<2<1<0<4\). The multiplication table is

{\small\par\smallskip
\CaseLine{\ttfamily 0 1 2 3 0}
\CaseLine{\ttfamily 2 3 3 3 2}
\CaseLine{\ttfamily 2 3 3 3 2}
\CaseLine{\ttfamily 3 3 3 3 3}
\CaseLine{\ttfamily 4 1 2 3 4}
\smallskip}

Let K=\{0,\allowbreak{}4\}. Any factor valued 3, or two positions valued in \{1,\allowbreak{}2\}, give value 3. With one \{1,\allowbreak{}2\}-valued occurrence, the product is 1 only if that occurrence is final and valued 1; otherwise it is 2. If all factors are in K, the product is its literal head's value. Induction from the table proves the formula for words of every positive length.

Fix a target w, \(C=c(w)\), \(h=its\) first variable, \(t=its\) last variable, and let U be the original right summands whose content is contained in C. Then \(w\leq P\) holds exactly when:

(H) some u in U starts with h;

(S) each simple variable x of w has a witness u in U which omits x or contains x exactly once. If \(x=t\) is simple, a witness containing x must end in x. Here simple means exactly one occurrence in the whole word, not just its interior.

For necessity of H assign \(h=4\), other target variables 0, outsiders 3. Without a contained same-head witness the target is 4 and every right word is 0 or 3. For a failed S at a nonfinal simple x, set \(x=2\), other target variables 0, outsiders 3. The target is 2 and every right word is 3. At a simple final t set \(t=1\) instead: the target is 1 while right words are 2 or 3. For sufficiency, a target value 3 is least. A K-valued target is absorbed by its same-head contained witness. A target valued 1 or 2 has a single \{1,\allowbreak{}2\}-valued occurrence, hence a simple variable x; an omitting witness is K-valued and above the target, and a single-x witness gives 1 or 2 with the necessary final status. Thus H, S are exact without bounds on variables or words. A literal target forces its literal summand.

\par\smallskip\noindent\textit{2. Nine structural laws and endpoint normal forms.}\quad 
{\small\par\smallskip
\CaseLine{\((W)  xyzt = xzyt\)}
\CaseLine{\((I)  xy^{3}z = xy^{2}z\)}
\CaseLine{\((H3) x^{3}y = x^{2}y\)}
\CaseLine{\((T3) xy^{3} = xy^{2}\)}
\CaseLine{\((R)  x^{3} = x^{2}\)}
\CaseLine{\((E)  x^{2}yx = xyx\)}
\CaseLine{\((D)  xy = x^{2}y+xy^{2}\)}
\CaseLine{\((L)  x^{2} \leq  x\)}
\CaseLine{\((A)  xy^{2} \leq  x\)}
\smallskip}

The word-value formula verifies the first six identities: they retain literal endpoints and whether each variable is simple or repeated. It verifies L and A as well: multiplying by a repeated variable can only keep the head's K-value or decrease the result to 3. For D, if there is a unique \{1,\allowbreak{}2\}-valued occurrence, precisely the branch which squares the other factor retains the original value and final status; the other branch is 3. If all factors are in K both branches retain the head, and a zero product remains 3 in both branches.

Keep both literal endpoints, permute the interior using W, and truncate each variable's total multiplicity to one or two. For an interior variable use I; for distinct head and tail use H3 or T3 after bringing interior copies next to their endpoint. If head equals tail, bring its interior copies next to the head, apply H3 while possible, then use E to remove the one remaining extra interior occurrence. For singleton-content words use R. Finally sort the interior by W. Thus a nonlinear normal word is determined by its literal head, literal tail and its multiplicities capped at two. Linear words stay linear. Every substitution uses nonempty images; empty external contexts are allowed. In particular E handles separated copies at both ends; H3 and T3 alone would not justify that case.

\par\smallskip\noindent\textit{3. Saturated guard and simple-variable recovery.}\quad 
Choose a contained same-head original v from H. Double every simple variable in v using L in its context; normalize. The resulting word has all variables repeated and remains below P. Append a square of each missing variable of C using A. Finally append \(t^{2}\) once, whether or not t already occurs, again using A. Normalizing yields a guard g of full content C with endpoints h, t and all multiplicities two. Thus \(g\leq P\). This explicit final append avoids assuming that repeated words may freely change an endpoint.

For a simple target variable x choose an original S-witness u and double every variable except x using L; call the normalized result \(u_{x}\). It remains below u, omits x or has one x, and preserves the required final position when \(x=t\). Let \(w_{x}\) be w with every variable except x doubled, normalized with both endpoints retained. Then

{\small\par\smallskip
\CaseLine{\(w_{x} \leq  g+u_{x}\)}
\smallskip}

is valid by H, S: g provides H and full content, and \(u_{x}\) handles the sole simple variable. Section 4 derives all such local recoveries from one fixed finite family.

If a current word has at least two simple variables, split it as XY between two simple occurrences, with X, Y nonempty. D rewrites it as \(X^{2}Y\)+\(XY^{2}\). Both branches retain content, head and tail, while each has fewer simple variables. A surviving simple variable retains its final/nonfinal status. Normalize and repeat. The finite tree has leaves with at most one simple variable. A zero-simple leaf equals g by the endpoint normal form; a one-simple-x leaf equals \(w_{x}\). Consequently every leaf is below P and D gives the original \(w\leq P\). Targets with zero or one simple variable use that same argument's leaf directly. Linear targets are already original summands.

\par\smallskip\noindent\textit{4. Reversible compression and completeness.}\quad 
In a local recovery there are at most three words. Reserve x and each word's head and tail, at most seven variables altogether. Every other variable is absent or repeated in each word and is never an endpoint. Group those variables by their nonzero incidence pattern across the three words: at most seven classes. Identify each class to a fresh variable and normalize. This is a genuine substitution into a valid absorption, so validity is preserved. At most fourteen variables remain and every word has length at most twenty-eight.

Inverse substitution replaces a class variable by the nonempty word listing its members. Every class occurrence was repeated and internal; its two copies expand to two copies of each member inside the fixed endpoints. W permutes and regroups this interior. Reserved variables keep their original multiplicities and endpoints. Normalization recovers the original three words. This preserves common multiplicity and never treats several single occurrences as an original repeated occurrence.

Fix fourteen variables. Take the nine structural laws and all valid absorptions with one normal word on the left, one or two normal words on the right, and each word of length at most twenty-eight. The exact test H, S effectively specifies this finite family. The inverse compression derives every local recovery; the guard and splitting construction then derives every word absorption with no bound on the original variables, lengths or summand count. Expand arbitrary terms by distributivity and use absorption summandwise in both directions. This proves that the family is a finite identity basis for 7319. Reversal gives a basis for its opposite. Catalogue transfer requires complete independently verified HSP in both directions.

\subsection{Four finite bases from two final-position tests}\label{case:182}
Work in the constant-free ai-semiring signature. Put \(C=c(w)\) for a target word and let U consist of the original right summands supported in C. Write w \ensuremath{\leq} P as an abbreviation for \(w+P=P\). The order in both sources is \(3<2<1<0<4\). Each opposite has the reversed presentation.

\par\smallskip\noindent\textit{7308 and 8457.}\quad 
Multiplication rows of 7308 are 01230,\allowbreak{}22231,\allowbreak{}22232,\allowbreak{}22232,\allowbreak{}01234. A final 3 gives 3, and these are the only words giving 3. On \{0,\allowbreak{}4\} the value is 0 if a 0 occurs and 4 otherwise. With some entry outside \{0,\allowbreak{}4\} and a final different from 3, the value is 1 precisely if the only low occurrence is a single 1 and everything after it is 4; otherwise it is 2. These statements follow by appending one factor.

An absorption holds precisely when (i) some original word has the same final variable as w, (ii) U is nonempty, and (iii) for each simple variable x of w, a word of U omits x or contains it once with its post-x support contained in the post-x support of w. Simple means occurring just once. For (i) assign 3 at all original finals and 2 at a missing target final. For (ii) put 4 on C and 2 outside. For (iii) give x the value 1, give its target suffix letters 4, the other target letters 0 and outsiders 2. Failure makes every upper word at most 2. Conversely at value 2 any same-final original has non-3 final and suffices. At high values a contained word suffices. At value 1 choose the indicated witness for the unique low letter; omission gives a high result and the simple witness gives 1.

A finite relative basis consists of \(xy=xxy+xyy\) and all its valid absorptions over 737 variables of maximum word length 2951 with at most two upper summands. The splitting identity holds by the rows: if the final value is 3 all outputs are 3; in the other cases, duplicating an input never raises the value and at least one of the two duplications preserves it. This can equally be checked on the 25 ordered pairs. Split a word having two simple letters into two consecutive blocks each containing one of them. In each expansion summand at least one fewer letter is simple, so iteration expresses both sides as sums of words with at most one simple letter.

For each expanded lower word keep its original same-final witness and one contained witness, choosing the latter as in (iii) when it has a simple letter. The exact criterion holds with these two upper words. To derive this bounded-width step, use the run compression proved in Section~\ref{case:180} in reversed reading order. Here are its hypotheses and their verification. Protect the final position, the unique simple position if present, and its predecessor if present. At most two complementary runs remain per word. Inside a run interchange adjacent distinct repeated letters or shorten an adjacent equal block, retaining a copy in that run and two copies overall. These edits preserve content, multiplicities capped at two, the final letter, and the set after each simple letter. For actual bounded identities mark both swapped positions and witnesses of their repetitions, or mark the two equal positions and a third surviving occurrence for deletion. Mark the final position as well. Represent each unmarked interval by a distinct fresh letter. Every new simple formal letter retains the same suffix support: an adjacent swap inside a marked pair never straddles it, and deleting one member of an adjacent equal pair preserves its support after that letter. The two swapped letters, or the shortened letter, remain repeated on both sides. The formal final letter is unchanged. Conditions (i)-(iii) give both bounded inequalities. Substituting the actual interval words back proves the edit. The formal length is at most 19.

Sort and cap each run at two. In the at most three words protect all letters of protected positions, at most nine. Group every remaining letter by its common vector of multiplicities 0,\allowbreak{}1,\allowbreak{}2 in the at most six runs. Identify each group to one letter, then reduce its count in each run to the common value, including two when appropriate. If the common count is one, the group member originally occurs in another run of that word as it is repeated and nowhere protected. The group letter thus remains repeated. All reductions are legal. There are at most \(9+3^{6}-1=737\) letters and at most \(3+4(737)=2951\) positions per word. Use this valid bounded absorption as an axiom. Replace each group letter by the nonempty product of its members. Every original run multiplicity is recovered; orient the preceding bounded edits backwards as necessary. This proves the original step. Deriving all expanded summand absorptions in both directions proves completeness with arbitrary alphabets. The finite family is specified effectively by the table and the stated bounds.

\par\smallskip\noindent\textit{7568 and 8052.}\quad 
Multiplication rows of 7568 are 01234,\allowbreak{}11234,\allowbreak{}33333,\allowbreak{}33333,\allowbreak{}41234. Any occurrence of 3 or any nonfinal 2 gives 3. A final 2 with all earlier entries in \{0,\allowbreak{}1,\allowbreak{}4\} gives 2. On \{0,\allowbreak{}1,\allowbreak{}4\} the result is the last nonzero entry, or 0 if all are zero. This follows by induction.

For each target variable x let \(S_{x}\) be the set after its last occurrence. There must be a member of U containing x with its last-x suffix supported in \(S_{x}\). In addition if the final variable t of w is simple, some member of U must omit t or end in its unique t. These conditions are necessary: for the first put \(x=4\), \(S_{x}=0\), the remaining target letters 1 and outsiders 3. For the second put \(t=2\), others in C equal to 0, outsiders 3. Conversely at value 4 use the witness for the last 4-valued letter. At value 1 a contained word gives at least 1, and at value 0 every target input is 0, so every contained witness gives 0. At value 2 the omission or simple final witness attains at least 2. The least value causes no constraint.

Take the following eight relative identities.

{\small\par\smallskip
\CaseLine{\(xxy = xy\qquad xyx = yxx\)}
\CaseLine{\(xy \leq  x+y\qquad xx \leq  x\)}
\CaseLine{\(xyz \leq  xz+y\qquad xyz \leq  xyzz+x\)}
\CaseLine{\(xyz \leq  xyzz+xz\qquad xy \leq  xyy+y\)}
\smallskip}

Soundness follows from the word formula or by all substitutions of the displayed variables into the table. The first two identities reduce any word to its list of distinct letters in last-occurrence order, with its final letter doubled exactly when that letter was repeated. For a repeated letter with a nonfinal later occurrence, replace x B x by B xx, then contract xx with a nonempty following factor. Adjacent repeats contract directly. For a repeated final x, move one earlier occurrence to the end by the same equation; the remaining earlier occurrences can then be removed using the penultimate x and its nonempty successor. These moves never change last order.

Concatenate contained originals until their product has support C, using \(xy\leq x+y\), and square it using \(xx\leq x\). The suffix witnesses give coverage; thus the resulting guard has full support and repeated final. Build the desired last-order list from right to left. If p is the already fixed suffix of this list, take a contained witness u for the preceding new variable x. Every letter after its last x lies in p. For empty p append u to the guard G by \(xy\leq x+y\). Its final is repeated because G has full support. For nonempty p use G u p \ensuremath{\leq} G p + u by the fifth law. G p normalizes back to G, and G u p is a repeated- final guard with last-order suffix x p. Finite induction produces the full desired ordering. A target with repeated final is obtained.

If the target has a simple final t and at least one other variable, write its normal form p t. The guard is p t t. If an original u omits t, use u p t \ensuremath{\leq} u p t t + u. Since c(u) is contained in c(p), the two products normalize to p t and p t t. If a witness is v t with nonempty v and unique t, use v p t \ensuremath{\leq} v p t t + v t instead. For the literal witness t use p t \ensuremath{\leq} p t t + t. A one-letter target itself requires that literal among the right words. This proves every summand absorption and hence all identities without a bound on input word lengths, substitutions or numbers of variables.

\subsection{Finite bases for two nilpotent extensions of a flat Boolean group}\label{case:183}
All bases are relative to the ai-semiring axioms in the signature \((+,\cdot)\) without constants. Let \(G\) be a nontrivial elementary abelian 2-group with identity \(e\). The carrier is \(G\cup\{a, z, f\}\). Multiplication extends that of \(G\) by

\CaseDisplay{fg=gf=f,\quad f^2=f,\quad ag=ga=af=fa=a,\quad a^2=z,\quad zx=xz=z.}

There are two additive orders:

\textbullet{} In \(D(G)\), \(z\) is least, \(f\) is greatest, and \(G\cup\{a\}\) is an antichain between them.

\textbullet{} In \(E(G)\), \(z<a<f\), each \(g\in G\) is below \(f\), and there are no other comparisons.

Addition is join. Both are ai-semirings. Associativity follows from the group core and the square-zero ideal \(\{a, z\}\); each multiplication translation preserves the displayed joins. For example, multiplication by \(a\) sends the core \(G\cup\{f\}\) to \(a\) and the ideal to \(z\). In either order these images preserve every mixed join.

The catalogue maps are

\CaseDisplay{D(C_2): (e, g, a, z, f)=(0,\allowbreak{}2,\allowbreak{}3,\allowbreak{}1,\allowbreak{}4)\quad\text{(type 7451)},}

\CaseDisplay{E(C_2): (e, g, a, z, f)=(2,\allowbreak{}3,\allowbreak{}0,\allowbreak{}1,\allowbreak{}4)\quad\text{(type 8831)}.}

\par\smallskip\noindent\textit{1. Finite identity schemata.}\quad 
Both algebras satisfy the four structural identities

\CaseDisplay{xy\approx yx,\qquad x^4\approx x^2,\qquad x^3\le x,\qquad xyz\le x+y+z.\quad\textup{(\ensuremath{S})}}

Here \(u\le v\) abbreviates \(u+v\approx v\). These laws follow from the operation rules. For the last law, a factor \(z\) makes the product \(z\), which is below the sum since that sum contains \(z\). If there is an \(a\) factor and a core factor, their join is \(f\) in both orders. Three factors equal to \(a\) have product \(z\le a\). In the core, either all values are the same group element, or their join is \(f\).

For each \((p, q)\in\{1,\allowbreak{}2,\allowbreak{}3\}\times\{0,\allowbreak{}1,\allowbreak{}2,\allowbreak{}3\}\) independently retain or omit a variable \(y_{pq}\). Write

\CaseDisplay{U=\prod y_{pq}^p,\qquad V=\prod_{q>0} y_{pq}^q,}

where only retained variables occur. Include

\CaseDisplay{xU\le x^3U+x^{\varepsilon}V.\quad\textup{(\ensuremath{L})}}

For \(D(G)\) take only \(\varepsilon=1\), giving 4096 identities. For \(E(G)\) take \(\varepsilon=0,\allowbreak{}1\), excluding an empty last summand, giving 8184 identities. An omitted block is an absent substring; no empty substitution or unit constant is allowed.

For \(E(G)\) add the following deletion schema. Independently retain or omit a distinct variable \(t_{pq}\) for each \((p, q)\in\{2,\allowbreak{}3\}\times\{0,\allowbreak{}1,\allowbreak{}2,\allowbreak{}3\}\); put

\CaseDisplay{W=\prod t_{pq}^p,\qquad Q=\prod_{q>0}t_{pq}^q.}

Whenever \(Q\) is nonempty, include

\CaseDisplay{W\le Wy^2+Q,\quad\textup{(\ensuremath{D})}}

with \(y\) a further variable. There are \(2^8-2^2=252\) instances.

Theorem. \((S)\) and the 4096 specified instances of \((L)\) form a complete finite basis for \(D(G)\). For \(E(G)\), a complete finite basis consists of \((S)\), its 8184 instances of \((L)\), and the 252 instances of \((D)\). Thus 7451 and 8831 are finitely based, with respective bases of at most 4100 and 8440 identities. The largest displayed identity uses 13 variables. Minimality is not asserted.

\par\smallskip\noindent\textit{2. Sound contextual lifting and deletion.}\quad 
Normalize words so that positive exponents are 1,\allowbreak{}2,\allowbreak{}3, using \((S)\). A word with all exponents 2 or 3 is called saturated. For any word \(w\) with exactly one occurrence of \(x\), let \(w_x\) replace that occurrence by \(x^3\). Let \(q\) be a nonempty word with \(c(q)\subseteq c(w)\). The following sound inequalities explain the finite schemata.

In \(D(G)\), if \(q\) contains \(x\) exactly once, then \(w\le w_x+q\). In \(E(G)\), the same inequality holds if \(q\) contains \(x\) either zero or one times. To prove this, \(w=z\) is handled by \(w_x=z\), and a core value of \(w\) is unchanged by cubing \(x\). If \(w=a\), there is exactly one occurrence assigned \(a\). When it is a different variable, again \(w=w_x\). Otherwise \(x=a\); if \(q\) contains \(x\) once, then \(q=a\). In the additional \(E(G)\) case where \(q\) omits \(x\), \(w_x=z\) and \(q\) is in the core, whose sum with \(z\) is \(f\). This bounds \(a\) as required.

For \(E(G)\), if \(w\) is saturated and \(c(q)\subseteq c(w)\), then

\CaseDisplay{w\le wy^2+q.\quad\textup{(\ensuremath{4})}}

Indeed, \(w\) is either \(z\) or a core element. In the former case \(wy^2=z\). In the latter, every variable of \(w\) and therefore of \(q\) lies in the core. If \(y\) is in the group, \(wy^2=w\); if \(y=f\), the right side contains \(f\); if \(y\) is in the ideal, \(wy^2=z\) and its sum with the core value of \(q\) is \(f\). This proves (4). These arguments establish all finite schema instances without limiting the number of variables in their word substitutions.

\par\smallskip\noindent\textit{3. Exact criterion for \(D(G)\).}\quad 
Let \(P\) be a nonempty polynomial of normalized words, let \(C=c(w)\), and let \(P_C\) consist of the summands supported in \(C\). Write \(\pi(u)\) for the exponent-parity vector. Then \(D(G)\models w\le P\) exactly when:

1. \(\bigcup_{u\in P_C}c(u)=C\);

2. \(\pi(w)\in\operatorname{aff}_{\mathbb F_2}\{\pi(u):u\in P_C\}\);

3. for every variable occurring once in \(w\), some \(u\in P_C\) contains it exactly once.

For necessity, assign outsiders \(z\), which is both additive least and multiplicative zero. An uncovered target variable assigned \(f\), with other target variables assigned \(e\), proves condition 1. Restricting to a subgroup \(C_2\) gives condition 2 by an affine separating functional. For condition 3, assign the selected target variable \(a\) and all other target variables \(e\). Without an exactly-once witness, all contained right words are \(e\) or \(z\), whose sum cannot bound \(a\).

For sufficiency, target value \(z\) is automatically absorbed. Target value \(a\) is bounded by a witness in condition 3. Target value \(f\) comes from some variable assigned \(f\) and is bounded by coverage. Otherwise the target is a group value and all its variables lie in \(G\). Contained right values either include two distinct group values, whose sum is \(f\), or are a single common value \(h\). The affine condition makes the target an odd product of copies of \(h\), hence equal to \(h\).

For a saturated target, conditions 1 and 2 give a direct derivation. Express its parity as an odd sum of right parities, then form the corresponding odd product and append two copies of every member of \(P_C\). The product has content \(C\), every exponent at least two, and the required parity, so it normalizes to the target. Repeated \(xyz\le x+y+z\) bounds any odd product by its factor sum.

For a target with single occurrences, choose \(x\) and an exactly-once witness \(q\). The word \(w_x\) has the same content and parity and one fewer single-occurrence variable; induction derives \(w_x\le P\). Partition variables other than \(x\) by their exponent pairs in \(w, q\). Substitute the nonempty product of each pair class for the matching retained \(y_{pq}\) in \((L)\). This gives \(w\le w_x+q\le P\). Thus all absorptions, and consequently all identities, are derivable.

\par\smallskip\noindent\textit{4. Exact criterion for \(E(G)\).}\quad 
Here \(z\) is not an additive least element, so noncontained right summands must not be discarded. Put \(U=c(P)\) and retain \(P_C\) as above. Then \(E(G)\models w\le P\) exactly when:

1. \(C\subseteq U\) and \(P_C\ne\varnothing\);

2. \(\pi(w)\in\operatorname{aff}_{\mathbb F_2}\{\pi(u):u\in P\}\);

3. for every variable occurring once in \(w\), some \(q\in P_C\) contains it zero or one times.

For necessity of \(C\subseteq U\), assign a missing variable \(f\) and all others \(e\). If \(P_C\) is empty, assigning \(C\) to \(e\) and all outsiders to \(z\) gives target \(e\) and right side \(z\). Group assignments prove condition 2. If condition 3 fails, assign its variable \(a\), other target variables \(e\), and outsiders \(z\); every right word is \(z\), whereas the target is \(a\).

For sufficiency, suppose first that \(w=z\). If \(P\) were a group element, all its summands would equal that group element, because the only elements below a group element are itself. Every variable of \(U\) would then be group-valued. Since \(C\subseteq U\), the target would also be group-valued, a contradiction. Thus \(P\) is \(a, z\) or \(f\), all above \(z\).

If \(w=a\), choose its single occurrence variable \(x\) assigned \(a\) and the witness \(q\) from condition 3. If \(q\) contains \(x\), it equals \(a\). If \(q\) omits \(x\), it is core-valued. If it is \(f\) we are done; otherwise some right summand contains \(x\) by \(C\subseteq U\), and its value is \(a\) or \(z\). Its sum with the group value of \(q\) is \(f\).

Finally let \(w\) be in the core. A contained right summand exists and is core-valued. Hence the only case requiring attention is that every right summand has the same group value \(h\). Then every variable of \(U\) belongs to \(G\), the target cannot be \(f\), and the affine condition implies that its value is \(h\). This proves the criterion.

\par\smallskip\noindent\textit{5. Completeness for \(E(G)\).}\quad 
For a saturated target, condition 2 supplies an odd product of right words with target parity. Append two copies of every right word. Its support is \(U\), and each variable has exponent at least two; normalize to a saturated word \(r\). The odd-product inequality derives \(r\le P\). Because \(C\subseteq U\) and the target parity is zero outside \(C\),

\CaseDisplay{r\approx w\prod_{y\in U\setminus C}y^2.}

If \(U=C\) this finishes. Otherwise choose \(q\in P_C\). Group variables of \(C\) by their exponent pairs in \(w, q\), where the first exponent is 2 or 3. In \((D)\) substitute the nonempty product of each pair class for its retained \(t_{pq}\), and substitute the nonempty product of \(U\setminus C\) for \(y\). This gives \(w\le r+q\le P\).

For a nonsaturated target, choose a single-occurrence variable \(x\) and its witness from condition 3. Cubing \(x\) preserves content, parity and all other witness conditions. Induction and the appropriate \(\varepsilon=0\) or 1 instance of \((L)\) derive the desired absorption. The exponent-pair substitutions are nonempty products; no identifying-and-normalizing step is used that could lose exponent-one information.

Thus every valid absorption is derivable, for arbitrarily many variables. Distributing identities and treating summands in both directions proves the theorem.

\subsection{Twelve identities for the semiring H5}\label{case:184}
The operations of \(H_5\) and the twelve identities below are those of Ma and Yue \cite[Table~2 and equations~(31)--(42)]{MaYueH5}.

The signature is exactly two binary operations, with no constants. An inequality u \ensuremath{\leq} P means \(u+P=P\). Nonempty terms are substituted for variables.

\par\smallskip\noindent\textit{Table identification and complete basis.}\quad 
The source H5, with labels 0,\allowbreak{}1,\allowbreak{}2,\allowbreak{}3,\allowbreak{}4, has tables

{\small\par\smallskip
\CaseLine{\(+\qquad *\)}
\CaseLine{\ttfamily 00000\qquad 00004}
\CaseLine{\ttfamily 01011\qquad 00004}
\CaseLine{\ttfamily 00222\qquad 01234}
\CaseLine{\ttfamily 01233\qquad 33334}
\CaseLine{\ttfamily 01234\qquad 44444}
\smallskip}

The bijection (0,\allowbreak{}1,\allowbreak{}2,\allowbreak{}3,\allowbreak{}4) \ensuremath{\to} (4,\allowbreak{}3,\allowbreak{}0,\allowbreak{}1,\allowbreak{}2) identifies H5 with catalogue

7478. Both operations, including every off-diagonal sum, are checked. The complete basis, in addition to the six ai-semiring axioms, is:

1. xxy = xy.

2. xyx = xyy.

3. x \ensuremath{\leq} xx.

4. xyz \ensuremath{\leq} xz + yz.

5. xyzt \ensuremath{\leq} xzt + yt.

6. (xxy)\ensuremath{^{2}} \ensuremath{\leq} xx + y.

7. (xvuyt)\ensuremath{^{2}} \ensuremath{\leq} xvutt + uy + v.

8. (xvyt)\ensuremath{^{2}} \ensuremath{\leq} xvtt + y + v.

9. (xabzcy)\ensuremath{^{2}} \ensuremath{\leq} xy + az + bzc.

10. (xbzcy)\ensuremath{^{2}} \ensuremath{\leq} xy + z + bzc.

11. (xazcy)\ensuremath{^{2}} \ensuremath{\leq} xy + az + zc.

12. (xzcy)\ensuremath{^{2}} \ensuremath{\leq} xy + z + zc.

Ma and Yue prove completeness for any number of variables \cite[Theorem~9.4]{MaYueH5}. Multiplication reversal gives a basis of the opposite semiring. Full mutual HSP certificates identify all transferred targets; one-way containment is never used to infer finite basability.

\par\smallskip\noindent\textit{Review of the unrestricted proof.}\quad 
For a target word q, retain only original upper summands p with c(p) subset c(q); call their set D. For the first occurrence of each b in q, let U be its preceding content. By \cite[Lemma~9.3]{MaYueH5}, \(q\leq P\) holds precisely when (i) some \(p=rbs\) in D has c(r) subset U; (ii) when this occurrence is nonfinal, either such a witness has s nonempty, or an original word in D has content contained in U. The second alternative does not replace the first condition.

This criterion follows directly from the product table: a factor 4 makes the whole word 4; otherwise ignore initial 2s, and the first non-2 factor determines the result, except that 1 becomes 0 when followed by a factor. For necessity set outside variables to 4, U to 2, b to 0 or 1 and the remaining target variables to 3. For sufficiency the same word formula handles the first non-2 value; the prefix-contained witness contributes 2, whose sum with 1 is 0.

Write \(u=\iota(w)\) for the list of distinct letters in the order of their first occurrences. Rules~1 and~2 give \(w\approx u\) when the final letter of \(w\) occurs once, and \(w\approx u^2\) when it recurs. Indeed, an internal repeated occurrence in \(avab\), with \(b\) nonempty, can be deleted: if \(v\) is nonempty then \(avab\approx av^2b\approx avb\), and for empty \(v\) use \(a^2b\approx ab\). Remove such occurrences successively. If a repeated letter \(a\) remains at the end, the word now reads \(ua\). Write \(u=pav\). For nonempty \(v\), use \(pava\approx pav^2\) and delete internal repetitions. This gives \(u\) followed by its own last letter, as does \(ua\) directly if \(v\) is empty. The same deletions in \(u^2\) give precisely this word.

Suppose q has a repeated final letter. Let N collect variables occurring nonfinally in original words in D. Some original word must end in N. Otherwise choose the first target variable outside N: condition (ii) would require an original word entirely in N, contradicting the assumed endpoints. If all target variables lie in N, nonempty D already contradicts the assumption. Thus there are original words az and bzc, c nonempty. Condition (ii) at the target head h also supplies an original hy with y nonempty. Rules 9--12, distinguished by whether a, b are empty contexts, give a square starting in h with content in c(q). Empty contexts select a different rule; they are not substituted for variables.

Rule 6 successively incorporates every original summand into this square. The result has full target content and a repeated final letter. Suppose its correct initial list is a and the remaining list is b, both nonempty. For the next desired variable d, a marked original witness rds (s nonempty, c(r) subset c(a)) and rule 5 with \((x,y,z,t)=(a,rd,b,s)\) insert d next. The upper word abs has the old normal form because s uses already covered variables. If no marked witness exists, the criterion supplies rd and v, with c(r), c(v) subset c(a). Rules 7 or 8 insert d next according as r is nonempty or empty. The correct prefix grows by one variable. This proves the repeated-final case for every finite number of variables and summands.

Suppose \(z\) occurs just once in \(q\), at its end. If \(q^2\law D\) holds, the repeated-final case and rule~3 give \(q\law q^2\law D\). Otherwise the strengthened condition (ii) for \(q^2\) must fail at \(z\), the only first occurrence whose final status changes. The letters before \(z\) form \(\ct(q)\setminus\{z\}\). Hence every summand in \(D\) contains \(z\), and its first \(z\) must be final. If \(q=z\), containment now gives \(D=\{z\}\). For longer \(q\), begin with an original summand starting with the head of \(q\). Use rule~4 to merge all nonempty upper prefixes before their common final \(z\). To place the next desired letter \(d\), use a witness \(rdsz\) from condition (i). Rule~5 with \((x,y,z,t)=(a,rd,b,sz)\) gives \(ardbsz\law absz+rdsz\), where \(a\) is the correct nonempty prefix and \(b\) lists the variables still to be placed before \(z\). The step is needed only when \(b\) is nonempty. The word \(sz\) is nonempty even if \(s\) is empty. Since \(ab\) already covers all letters other than \(z\), \(absz\) has the previous normal form. Iteration gives \(q\law D\). Taking both directions of a polynomial equality proves completeness.

\subsection{Eight complete identities for 7503 by contained-prefix insertion}\label{case:185}
Relative to the additively idempotent semiring axioms, in the signature (+, multiplication), without constants.

\par\smallskip\noindent\textit{1. Word values and the exact criterion.}\quad 
Addition rows in order 0,\allowbreak{}1,\allowbreak{}2,\allowbreak{}3,\allowbreak{}4 are 00004,\allowbreak{}01114,\allowbreak{}01224,\allowbreak{}01234,\allowbreak{}44444, giving additive order \(3<2<1<0<4\). Multiplication is

{\small\par\smallskip
\CaseLine{\ttfamily 0 1 2 3 4}
\CaseLine{\ttfamily 1 1 2 3 1}
\CaseLine{\ttfamily 3 3 3 3 3}
\CaseLine{\ttfamily 3 3 3 3 3}
\CaseLine{\ttfamily 4 4 2 3 4}
\smallskip}

A factor valued 3 gives 3. A 2 before the final position gives 3. If the final factor is 2 and all earlier factors are in \{0,\allowbreak{}1,\allowbreak{}4\}, the product is 2. In all remaining cases every factor is in \{0,\allowbreak{}1,\allowbreak{}4\}: the value is the first nonzero factor, or 0 if all factors are 0. Induction on the word length proves this formula, including literals.

For a target w put \(C=c(w)\), and let U consist of the original right words whose contents are contained in C. The exact absorption test is:

\textbullet{} For every variable x of w, there is a word u in U containing x with its first-x prefix contained in the first-x prefix of w.

\textbullet{} If the final variable s of w is simple, some word u in U omits s, or contains s exactly once at its final position.

The first condition implies that U covers C. For its necessity at x, assign the target prefix before x to 0, x to 4, other target variables to 1, and outsiders to 3. The target is 4. A right word can reach 4 only if it is contained in C and has the required first-x prefix. For failure of the second condition, assign \(s=2\), other target variables 0, outsiders 3. The target is 2; every right word is 3.

For sufficiency, value 3 is least. An all-zero target is absorbed by any contained witness. A target valued 1 is absorbed by any contained word, since its content has values in \{0,\allowbreak{}1,\allowbreak{}4\}. If the target is 4, take the witness for its first 4-valued variable; all preceding variables in that witness are 0, so it also gives 4. A target valued 2 has a unique, simple final 2-valued variable. A witness omitting that variable is valued in \{0,\allowbreak{}1,\allowbreak{}4\}; a witness containing it once at the end gives 2. This proves the test without bounds on variables, word lengths, or number of right summands. A literal target forces its literal summand by the second condition.

\par\smallskip\noindent\textit{2. The eight laws.}\quad 
{\small\par\smallskip
\CaseLine{\((W1) xxy = xy\)}
\CaseLine{\((W2) xyx = xyy\)}
\CaseLine{\((A)  xy \leq  x+y\)}
\CaseLine{\((S)  xx \leq  x\)}
\CaseLine{\((I)  xyz \leq  xz+y\)}
\CaseLine{\((O)  xyz \leq  xyzz+y\)}
\CaseLine{\((T)  xyz \leq  xyzz+yz\)}
\CaseLine{\((L)  xy \leq  xyy+y\)}
\smallskip}

All are sound by the word-value formula or the exact absorption test. The first two give the following normal form: list first occurrences F(w); retain that list when the final variable is simple; otherwise append its last listed variable once more. Linear words stay linear. For internal separated repetitions use \(xYxZ=xY^{2}Z=xYZ\) with Y, Z nonempty; adjacent internal repeats use W1. A repeated final letter in p x Y x becomes p x \(Y^{2}\) using W2, after which internal repeats are removed. These are nonempty substitutions with possibly empty external contexts. Call a normal word saturated when its last listed variable is repeated. Saturation here does not mean every variable is repeated.

\par\smallskip\noindent\textit{3. A full-content guard with the desired first-occurrence order.}\quad 
Fix a valid nonlinear absorption \(w\leq P\) and the original contained set U. Its contents cover C. Starting with any u in U, multiply by additional words of U until the product has full content C. Each product remains below P by A. Apply S to this full-content word g: \(g^{2}\leq g\leq P\). Its normal form is saturated and has full content C.

We reorder this guard's first occurrences to match \(x_{1}\cdots x_{d}=F(w)\), one variable at a time. Suppose the guard's initial first-occurrence prefix is \(p=x_{1}\cdots x_{i-1}\). Choose an original u in U witnessing \(x_{i}\). Its first-\(x_{i}\) prefix is contained in c(p). If p is empty, A gives u \(g\leq u+g\leq P\). If p is nonempty, I gives

{\small\par\smallskip
\CaseLine{\(p u g \leq  p g+u = g+u \leq  P.\)}
\smallskip}

The equality p \(g=g\) follows from W1, W2: the first-occurrence list of g already starts with p, and the final variable of g remains repeated. The new guard is the normal form of p u g (or u g for empty p). It still has full content C and a repeated final variable. Its first occurrence list begins p \(x_{i}\): every variable before the first \(x_{i}\) in u was already in p. Any other ordering changes are after this fixed prefix. Thus the induction reaches a saturated guard G whose first occurrence list is exactly F(w), and \(G\leq P\).

All products and substituted words are nonempty. The empty prefix case uses A separately and never substitutes an empty word into I.

\par\smallskip\noindent\textit{4. Restoring the final simple variable.}\quad 
If the target is saturated, it is exactly G and the proof is done. Otherwise its normal form is p s, where p is nonempty and s is simple. The guard is p s s. Select the original witness u from the second criterion condition.

If u omits s, its content is contained in c(p). Substitute \(x=p\), \(y=u\), \(z=s\) into O. The first-occurrence laws normalize p u s to p s, and p u s s to p s s. Thus O gives p \(s\leq G+u\leq P\).

If u ends in its unique s and is nonlinear, write \(u=r\) s with nonempty r and c(r) contained in c(p). Substitute \(x=p\), \(y=r\), \(z=s\) into T. Again p r s normalizes to p s, and p r s s to G, giving p \(s\leq G+u\leq P\). If u is the literal s, use L directly with \(x=p\), \(y=s\). These alternatives cover every possible final witness, without attempting the generally false absorption p \(s\leq p\) s s.

\par\smallskip\noindent\textit{5. Completeness.}\quad 
Every valid word absorption has now been derived. Linear targets are already literal summands. Distributivity expands arbitrary terms to nonempty polynomials; derive the absorption of each summand into the opposite polynomial and use both directions. Therefore the eight laws form a finite identity basis for 7503, unrestricted in original variables, word lengths, number of right terms, and term substitutions. Word reversal gives a basis for the opposite. Independent literal derivation replay supports the proof; full two-way HSP verification is required for catalogue transfers. Finite tests alone do not classify.

\subsection{An affine-slice criterion and a finite basis for semiring 7597}\label{case:186}
We work in the ai-semiring signature \((+,\cdot)\) without constants. Let

\CaseDisplay{A=\{E, e, b, g, f\}}

have additive order \(e<E<f\), \(g<b<f\), with no other comparisons. The elements \(\{e, g\}\) form a group of order two; \(E\) is a global multiplicative identity, \(f\) a multiplicative zero, and multiplication by \(b\) outside \(E\) agrees with multiplication by \(g\), with \(b^2=e\). Thus \(eb=g\), \(gb=e\), \(b^3=g\). The catalogue map is

\CaseDisplay{(E, e, b, g, f)=(0,\allowbreak{}1,\allowbreak{}2,\allowbreak{}3,\allowbreak{}4),}

so this is precisely catalogue type 7597.

The map \(\rho(E)=\rho(e)=e\), \(\rho(b)=\rho(g)=g\), \(\rho(f)=f\) maps onto the flat group semiring \(\{e, g, f\}\). It preserves both operations. Multiplication is associative: factors equal to \(E\) can be omitted; one remaining factor is unchanged; a product of at least two non-\(E\) factors is the product of their \(\rho\)-images. Translations preserve the joins of the two chains and their common top. Hence \(A\) is an ai-semiring.

\par\smallskip\noindent\textit{1. A finite basis specification.}\quad 
Write \(u\le v\) for \(u+v\approx v\). Include the four identities

\CaseDisplay{xy\approx yx,\qquad x^4\approx x^2,\qquad x^3\le x,\qquad xyz\le x+y+z.\quad\textup{(\ensuremath{S})}}

For each \((p, q)\in\{1,\allowbreak{}2,\allowbreak{}3\}\times\{0,\allowbreak{}1,\allowbreak{}2,\allowbreak{}3\}\) independently retain or omit a variable \(y_{pq}\). Set

\CaseDisplay{U=\prod y_{pq}^p,\qquad V=\prod_{q>0}y_{pq}^q,}

using retained variables only. For \(\varepsilon\in\{0,\allowbreak{}1,\allowbreak{}2\}\) include

\CaseDisplay{xU\le x^3U+x^\varepsilon V,\quad\textup{(\ensuremath{L})}}

unless its last summand is empty. This gives \(3\cdot2^{12}-2^3=12280\) identities.

For each \((p, q)\in\{2,\allowbreak{}3\}\times\{0,\allowbreak{}1,\allowbreak{}2,\allowbreak{}3\}\) independently retain or omit a variable \(t_{pq}\). Set

\CaseDisplay{W=\prod t_{pq}^p,\qquad Q=\prod_{q>0}t_{pq}^q.}

Whenever \(Q\) is nonempty, include

\CaseDisplay{W\le Wy^2+Q.\quad\textup{(\ensuremath{D})}}

There are 252 instances. Products over omitted classes are absent substrings, never unit constants or empty substitutions.

Theorem. \((S)\), the 12280 instances of \((L)\), and the 252 instances of \((D)\) form a complete identity basis for \(A\). In particular, type 7597 is finitely based. The basis contains at most 12536 identities and uses at most 13 variables in any displayed identity. No minimality claim is made.

\par\smallskip\noindent\textit{2. Soundness of the finite schemata.}\quad 
The first three structural laws follow immediately from the multiplication rules. For \(xyz\le x+y+z\), if the join on the right is \(f\) there is nothing to prove. Otherwise all factors lie in one of the two chains. In the \(e, E\) chain the product is below its largest factor. In the \(g, b\) chain a triple product is \(g\), again below the join.

Normalize positive exponents to 1,\allowbreak{}2,\allowbreak{}3 using \((S)\). Cubing a variable changes only the value \(b\), replacing it by \(g\). Let \(w\) contain \(x\) exactly once, and let \(w_x\) replace this by \(x^3\). The only assignment where \(w\ne w_x\) has \(x=b\) and all other variables of \(w\) equal to \(E\). Under it, \(w=b\), \(w_x=g\). If a nonempty word \(q\) is supported in \(c(w)\) and its \(x\)-exponent is 0,\allowbreak{}1 or 2, its value is respectively \(E, b\) or \(e\). In all three cases \(g+q\) bounds \(b\). Thus

\CaseDisplay{w\le w_x+q\quad\textup{(\ensuremath{4})}}

is sound under precisely the stated restrictions, proving \((L)\).

Call a normalized word saturated when all positive exponents are 2 or 3. If \(w\) is saturated and \(q\) is nonempty with \(c(q)\subseteq c(w)\), then \(w\le wy^2+q\). If \(w=f\), then \(wy^2=f\). If \(w\) is \(e\) or \(g\), multiplication by \(y^2\) either leaves it unchanged or yields \(f\). If \(w=E\), every variable of \(w\) is \(E\), so \(q=E\). A saturated word cannot take value \(b\). This proves \((D)\).

\par\smallskip\noindent\textit{3. Exact absorption criterion.}\quad 
Let \(P\) be a nonempty polynomial of normalized words, \(C=c(w)\), \(U=c(P)\) and \(P_C=\{q\in P:c(q)\subseteq C\}\). Let

\CaseDisplay{H=\operatorname{aff}_{\mathbb F_2}\{\pi(q):q\in P\},\qquad H_C=\{v\in H:\operatorname{supp}(v)\subseteq C\},}

where \(\pi\) records exponent parity. Then \(A\models w\le P\) if and only if:

1. \(C\subseteq U\) and \(P_C\ne\varnothing\);

2. \(\pi(w)\in H\);

3. for every variable \(x\) occurring exactly once in \(w\), either some \(q\in P_C\) contains \(x\) exactly once, or some \(v\in H_C\) has \(v_x=0\).

\par\smallskip\noindent\textit{Necessity.}\quad 
An absent target variable assigned \(f\), with all others \(E\), proves \(C\subseteq U\). If \(P_C\) is empty, assign every variable in \(C\) to \(E\) and every outsider to the local identity \(e\). Then \(w=E\) and every right summand is \(e\), contradicting \(E\not\le e\). Group assignments in \(\{e, g\}\) and affine separation prove condition 2.

Suppose condition 3 fails at \(x\). By condition 2, \(\pi(w)\in H_C\), and every vector of \(H_C\) has \(x\)-coordinate 1. Put \(Z=U\setminus C\) and write \(H=\pi(w)+V\), where \(V\) is a linear subspace. For any \(v\in V\) vanishing on \(Z\), its \(x\)-coordinate is zero. Hence the linear functional \(v\mapsto v_x\) on \(V\) factors through the projection onto the \(Z\) coordinates. Extend that functional from the projection image to all \(\mathbb F_2^Z\); there is a vector \(d\) on \(Z\) such that

\CaseDisplay{h_x+\langle d, h|_Z\rangle=1\quad\text{for every }h\in H.\quad\textup{(\ensuremath{5})}}

Assign \(x=b\), every other variable of \(C\) to \(E\), and an outsider \(y\) to \(g\) or \(e\) according as \(d_y=1\) or 0. Every right word has projected group value \(g\) by (5). A contained word cannot have \(x\)-exponent one, by the assumed failure; its parity forces exponent three, so it evaluates to \(g\). A noncontained word includes an outsider assigned a local group element; it therefore evaluates to its projected value \(g\) as well. Thus \(P=g\) and \(w=b\), a contradiction. This proves necessity.

\par\smallskip\noindent\textit{Sufficiency.}\quad 
The homomorphism \(\rho\) and the affine condition imply that whenever all right projected values coincide in one group element, the projected target is that same element. More explicitly, if a variable of \(U\) is assigned \(f\), some right word is \(f\) and \(P=f\). Otherwise, if \(P\ne f\), all right projections have the same value \(h\in\{e, g\}\), and any vector of \(H\) evaluates to \(h\) under the group projection.

If \(w=f\), some target variable is \(f\), so coverage gives \(P=f\). If \(w=e\) or \(g\), the preceding projected-value argument places \(P\) on the same branch, above \(w\). If \(w=E\), every target variable is \(E\) and a contained right word is \(E\). The remaining case \(w=b\) has a unique variable \(x\) occurring once and assigned \(b\), with all other target variables \(E\). A contained exactly-once witness is \(b\). If instead there is \(v\in H_C\) with \(v_x=0\), assume for contradiction \(P\ne f\). The vector \(\pi(w)\in H\) evaluates to \(g\), whereas \(v\) evaluates to \(e\) because it is supported on \(C\) and has even \(x\)-coordinate. This contradicts the fact that every vector of \(H\) evaluates to the same \(h\). Thus \(P=f\). All target values are covered.

\par\smallskip\noindent\textit{4. Unrestricted derivations.}\quad 
First consider a saturated target satisfying conditions 1 and 2. Express its parity as an odd sum of right parities. Form the corresponding odd product of right words and append two copies of every right word. The product has support \(U\), every exponent at least two, and target parity. After normalization it is

\CaseDisplay{r=w\prod_{y\in U\setminus C}y^2.}

Repeated \(xyz\le x+y+z\) places an odd product below the sum of its factors, so \(r\le P\) is derivable. If \(U=C\) this is the desired result. Otherwise choose a contained right word \(q\in P_C\). Group variables in \(C\) by their exponent pairs in \(w, q\), substitute each nonempty class product for its matching placeholder in \((D)\), and substitute the product of \(U\setminus C\) for \(y\). This gives \(w\le r+q\le P\).

For a nonsaturated target, induct on the number of variables occurring once. Select one, \(x\). Cubing it preserves content and parity and removes only its own requirement in condition 3, so induction derives \(w_x\le P\).

If condition 3 supplies a contained exactly-once witness \(q\), use (4) with exponent one. Otherwise select \(v\in H_C\) with \(v_x=0\). Construct a saturated word \(s\) with content exactly \(C\), exponent two at coordinates where \(v\) is zero, and exponent three where it is one. The saturated argument just proved derives \(s\le P\). In particular its \(x\)-exponent is two. Apply (4) with \(q=s\).

Each required instance of (4) comes from the finite schema \((L)\) by grouping variables other than \(x\) by their exponent pairs in the target and witness. The 12 possible pairs are exactly those in its definition. Every substitution is a nonempty word; omitted classes are handled by omitted placeholders. Thus the induction supplies an equational derivation of every valid absorption, without a bound on variables or summands. Applying this to every summand in both directions of a polynomial identity proves completeness.

\subsection{Four finite bases by reversible symmetrization of multiplication}\label{case:187}
The signature has two binary operations and no constants. The following exact term definitions transfer already proved complete bases, with no restriction on the variables in a derived identity.

\noindent Target A\enspace\textbar\enspace Commutative source S\enspace\textbar\enspace Derived multiplication on A\enspace\textbar\enspace Original multiplication expressed in S\par
\noindent 7601\enspace\textbar\enspace 7597\enspace\textbar\enspace x * y = xy + yx\enspace\textbar\enspace x * x * x * y\par
\noindent 8066\enspace\textbar\enspace 7597\enspace\textbar\enspace x * y = xy + yx\enspace\textbar\enspace x * y * y * y\par
\noindent 7666\enspace\textbar\enspace 7609\enspace\textbar\enspace x * y = xy + yx\enspace\textbar\enspace x * x * x * y\par
\noindent 8596\enspace\textbar\enspace 7609\enspace\textbar\enspace x * y = xy + yx\enspace\textbar\enspace x * y * y * y\par
The labels \(0,1,2,3,4\) and the addition are common to the two presentations. Evaluating the defining terms on the 25 ordered pairs gives the displayed tables. The reverse terms recover the two operations. Lemma~\ref{lem:terms} applies.

\par\smallskip\noindent\textit{1. Previously proved source bases.}\quad 
For 7597 use the four structural identities, 12,\allowbreak{}280 lifting instances and 252 deletion instances in Section~\ref{case:186}. This is a complete basis of at most 12,\allowbreak{}536 identities relative to the ai-semiring axioms. The affine-slice alternative and the exponent-pair grouping in its proof are essential and retained without change.

For \(A_{7609}\), the four identities in (1) and the \(8{,}184\) instances of (2) in Section~\ref{case:188} give at most \(8{,}188\) identities relative to \(\ax\). The two source bases apply on arbitrary alphabets.

\par\smallskip\noindent\textit{2. A concrete finite basis specification for each target.}\quad 
Let D(x, y)=xy+yx in the target signature. Let \(B_{S}\) be the indicated complete source basis, including the six usual ai-semiring axioms: associativity, commutativity and idempotence of addition, associativity of multiplication, and the two distributivity laws. Translate \(B_{S}\) recursively by replacing every source multiplication by D, and leaving addition and variables unchanged. A product of several factors is translated with any fixed parenthesization, for example to the left.

For targets 7601 and 7666 add the identity

{\small\par\smallskip
\CaseLine{\(xy = D(D(D(x, x),x),y).\)}
\smallskip}

For targets 8066 and 8596 add

{\small\par\smallskip
\CaseLine{\(xy = D(D(D(x, y),y),y).\)}
\smallskip}

These are ordinary finite terms in the original signature. No star operation or power constant is added. Thus this prescription gives at most 12,\allowbreak{}543 identities for 7601 and 8066 and at most 8,\allowbreak{}195 for 7666 and 8596. It is a finite specification by translating the existing finite schemata, not a claim that a short subset of those instances is already complete.

\par\smallskip\noindent\textit{3. Unrestricted completeness of the translation.}\quad 
Write t(x, y) for the recovery product in the last table column, viewed as a term in the source signature. On A, the operations (+, D) are exactly S, and t evaluated using D gives the original multiplication of A. Conversely, on S the operations (+, t) are exactly A. Both statements follow from the displayed operation-table equalities, not merely from HSP inclusion.

The translated source identities hold in A because (+, D) is S. The recovery identity also holds, so the proposed finite basis is sound.

Take any identity \(u=v\) valid in A, with arbitrarily many variables. Replace each original multiplication in u and v by the source term t. The resulting identity \(u^{t}=v^{t}\) is valid in S, since S with multiplication t is A. By the complete source basis \(B_{S}\) it has an equational derivation. Translate that derivation recursively by replacing every source multiplication by D. Equational substitution and congruence preserve such a translation, so the target basis derives \((u^{t})^{D}=(v^{t})^{D}\).

An induction on the construction of u, using the recovery identity at each multiplication node and congruence at each addition node, gives \((u^{t})^{D}=u\); likewise \((v^{t})^{D}=v\). Hence the target basis derives \(u=v\). This proves completeness for every term and every variable set.

The proof does not require finite tests of arbitrary identities, nor does it infer finite basability from membership in a finite-basis variety. It is an explicit reversible term-definition theorem.

\subsection{A finite basis for a flat Boolean group above a two-step nilpotent ideal}\label{case:188}
Work in the constant-free ai-semiring signature. Write \(u\le v\) for \(u+v\approx v\). Let \(G\) be a nontrivial elementary abelian 2-group. On \(B(G)=G\cup\{f, a, z\}\), with these three elements outside \(G\), use the additive order

\CaseDisplay{z<a<g<f\quad(g\in G),}

with distinct group elements incomparable. Multiplication restricts to the group operation on \(G\), satisfies \(fg=gf=f\), \(f^2=f\), \(ag=ga=a\), \(af=fa=a\), \(a^2=z\), and has zero \(z\). The group identity is also a global multiplicative identity.

For \(G=C_2=\{e, g\}\), the catalogue correspondence

\CaseDisplay{(e, g, a, z, f)=(0,\allowbreak{}3,\allowbreak{}1,\allowbreak{}2,\allowbreak{}4)}

identifies \(B(C_2)\) with type 7609.

\par\smallskip\noindent\textit{1. Finite basis.}\quad 
The following four structural identities will be used:

\CaseDisplay{xy\approx yx,\qquad x^4\approx x^2,\qquad x^3\le x,\qquad xyz\le x+y+z.\quad\textup{(\ensuremath{1})}}

Here is a finite family of additional identities. For each \((p, q)\in\{1,\allowbreak{}2,\allowbreak{}3\}\times\{0,\allowbreak{}1,\allowbreak{}2,\allowbreak{}3\}\), independently retain or omit a distinct variable \(y_{pq}\). Put

\CaseDisplay{U=\prod_{\text{retained }(p, q)}y_{pq}^{p},\qquad V=\prod_{\text{retained }(p, q),\ q>0}y_{pq}^{q}.}

For \(\varepsilon\in\{0,\allowbreak{}1\}\) include

\CaseDisplay{xU\le x^3U+x^{\varepsilon}V,\quad\textup{(\ensuremath{2})}}

provided the last summand is nonempty. An omitted product is an absent substring, not a unit constant; \(x^0\) is likewise omitted. There are \(2\cdot2^{12}-2^3=8184\) choices. Each identity uses at most 13 variables.

Theorem. Identities (1) and (2), relative to the ai-semiring axioms, form a complete finite identity basis for \(B(G)\). In particular, type 7609 is finitely based. This is a finite explicit schema of at most 8188 identities, with no claim of minimality.

\par\smallskip\noindent\textit{2. Algebra and soundness.}\quad 
The core \(G\cup\{f\}\) is the flat group semiring, and \(\{a, z\}\) is a multiplicative ideal with square \(\{z\}\). Every core element acts identically on \(a\), and \(z\) is absorbing, proving associativity in the mixed cases. Multiplication by \(g\in G\) permutes the group antichain and fixes \(a, z, f\). Multiplication by \(f\) sends every group element to \(f\), fixing \(a, z, f\); multiplication by \(a\) sends the core to \(a\) and \(a, z\) to \(z\). Each map preserves joins, proving distributivity.

The first three laws in (1) follow directly. For the fourth, if all three factors belong to the core, either they are the same group element and their product equals that element, or their join is \(f\), or a factor is \(f\). In every case the product is below the join. If at least one factor is \(a\) or \(z\) and at least one is in the core, the product is \(a\) or \(z\) and lies below that core factor. If all factors are \(a\) or \(z\), the product is \(z\).

We establish a general sound absorption, which proves (2). Let \(t\) be any commutative word in which a variable \(x\) has multiplicity exactly one. Let \(t_x\) replace this occurrence by \(x^3\). Suppose \(q\) is nonempty, \(c(q)\subseteq c(t)\), and the multiplicity of \(x\) in \(q\) is 0 or 1. Then

\CaseDisplay{B(G)\models t\le t_x+q.\quad\textup{(\ensuremath{3})}}

Indeed, if \(t=z\) there is nothing to prove. If \(t\) is in the core, all its variables are in the core and cubing \(x\) leaves its value unchanged, so \(t=t_x\). If \(t=a\), exactly one occurrence in the word is assigned \(a\), and all other occurrences have core values. If that occurrence belongs to a variable other than \(x\), again \(t=t_x\). If it belongs to \(x\), then \(q\) is either \(a\) (when it contains \(x\)) or a core element (when it does not), so \(a\le q\). This proves (3), without any bound on the contents of \(t, q\).

\par\smallskip\noindent\textit{3. Exact absorption criterion.}\quad 
Normalize each word using commutativity and \(x^4\approx x^2\). A positive exponent becomes 1, 2 or 3; exponents at least two reduce according to parity. Let \(\pi(w)\) denote the vector of exponent parities over \(\mathbb F_2\). For a normalized target \(w\) with content \(C\) and polynomial \(P\), let \(P_C\) be the right summands supported within \(C\).

Lemma. The identity \(w\le P\) holds in \(B(G)\) exactly when the following three conditions hold:

1. \(\bigcup_{u\in P_C}c(u)=C\) (in particular \(P_C\) is nonempty);

2. \(\pi(w)\) belongs to the affine hull over \(\mathbb F_2\) of \(\{\pi(u):u\in P_C\}\);

3. for every variable with multiplicity one in \(w\), some \(u\in P_C\) has multiplicity 0 or 1 in that variable.

Proof. Since \(z\) is both additive least and multiplicative zero, assigning outsiders \(z\) shows that \(w\le P_C\) is necessary. For condition 1, assign a supposedly uncovered variable \(f\) and every other target variable \(e\); then the left side is \(f\) and all contained right words are \(e\), which is impossible. If \(P_C\) is empty, assign every target variable \(e\) to get the same contradiction with right side \(z\).

For condition 2, work inside a subgroup \(\{e, g\}\cong C_2\). If the target parity vector is outside the affine hull, finite-dimensional linear algebra supplies a linear functional constant on all right parity vectors and taking the other value at the target. The corresponding assignment gives all right summands one group value and the target the other; distinct group elements are incomparable.

For condition 3, assign the selected variable \(a\), every other target variable \(e\), and outsiders \(z\). The target is \(a\). If every contained right word has at least two occurrences of the selected variable, every right summand is \(z\), a contradiction.

Conversely, if the target value is \(z\), absorption is automatic. If it is \(a\), exactly one variable occurs once and is assigned \(a\), with all other target variables in the core. A witness from condition 3 is either \(a\) or a core element, hence bounds the target. If the target is \(f\), some target variable is \(f\), so coverage gives a right summand equal to \(f\). In the remaining case all target variables lie in \(G\). If the contained right summands have different group values, their join is \(f\); otherwise they all have the same value \(h\). The affine-hull condition expresses the target parity as the sum of an odd number of right parities, so its group value is the corresponding odd product \(h\), proving absorption. \(\square\)

In the preceding paragraph a product of group values cannot equal \(f\); hence the \(f\) case arises precisely from a factor \(f\). This separates the core cases without treating the semiring addition as a group operation.

\par\smallskip\noindent\textit{4. Derivations for unrestricted numbers of variables.}\quad 
First suppose that every variable in \(w\) has exponent 2 or 3, and that conditions 1 and 2 hold. Express \(\pi(w)\) as the sum of an odd number of parities of words in \(P_C\). Form their product, and append two additional copies of every word in \(P_C\). The resulting product still has an odd number of factors and the same parity. Coverage ensures that its content is \(C\), with every variable occurring at least twice. By normalization it is equal to \(w\). Repeated applications of \(xyz\le x+y+z\), first to any three factors and then inductively, place the odd product below the sum of its factors, hence below \(P\). Thus \(w\le P\) is derivable.

Now induct on the number of exponent-one variables in a normalized target satisfying all three conditions. Choose one, \(x\), and its witness \(q\in P_C\) from condition 3. Let \(w_x\) replace \(x\) by \(x^3\). Its content and parity are unchanged, and its other exponent-one variables have the same witnesses, so it satisfies the criterion and has one fewer exponent-one variable. By induction \(w_x\le P\) is derivable.

It remains to derive the sound inequality \(w\le w_x+q\) from the finite family (2). For each variable \(y\ne x\) in \(C\), record its exponent pair \((p, q_y)\) in the normalized words \(w, q\). There are at most 12 pairs, with \(p\in\{1,\allowbreak{}2,\allowbreak{}3\}\) and \(q_y\in\{0,\allowbreak{}1,\allowbreak{}2,\allowbreak{}3\}\). For each nonempty class of variables of a fixed pair, substitute their product for the corresponding retained \(y_{p, q_y}\) in (2), and omit the other placeholders. Choose \(\varepsilon\) equal to the exponent of \(x\) in \(q\). The last summand is nonempty because \(q\) is nonempty. The resulting instance, using commutativity, is exactly \(w\le w_x+q\).

Every substitution here is a nonempty word. In particular, this grouping is not performed by identifying variables and then normalizing their total exponents, which could destroy exponent-one data; it is a direct substitution into a sound finite schema. Since \(q\le P\) and \(w_x\le P\), the induction is complete.

Every polynomial identity is equivalent to the absorption of each summand in both directions. The exact criterion and these derivations establish completeness for arbitrary variable sets, proving the theorem.

\subsection{Finite bases for truncated chain semirings with a split identity}\label{case:189}
All identities are in the constant-free signature \((+,\cdot)\), relative to the axioms of additively idempotent semirings. Write \(u\le v\) for \(u+v\approx v\). Element names below are not constants in this signature.

\par\smallskip\noindent\textit{1. The algebras and a finite basis specification.}\quad 
Fix integers \(k\ge2\) and \(1\le j\le k\). Let

\CaseDisplay{A_{k, j}=\{e, v_0, v_1,\ldots, v_k\}.}

The elements \(v_k<v_{k-1}<\cdots<v_0\) form a chain in the additive order. Place \(e\) below \(v_0\) and above precisely \(v_j,\ldots, v_k\), and make it incomparable with \(v_1,\ldots, v_{j-1}\). Addition is join in this order. Multiplication has identity \(e\) and satisfies

\CaseDisplay{v_rv_s=v_{\min(k, r+s)}.\quad\textup{(\ensuremath{1})}}

Thus \(v_k\) is both the additive least element and a multiplicative zero.

Put \(n=k+2\), \(M=n^k\) and \(N=(k-1)+2^M\). Let \(W_N\) be the finite set of nonempty commutative words on \(N\) fixed variables, with each multiplicity at most \(k\). Let \(\mathcal E_{k, j}\) consist of every identity

\CaseDisplay{w\le\sum_{u\in Q}u\quad\textup{(\ensuremath{2})}}

valid in \(A_{k, j}\), with \(w\in W_N\) and \(\varnothing\ne Q\subseteq W_N\). This specifies a finite set by a finite truth-table rule; it does not require an oracle for arbitrary-variable identities.

Theorem. Every \(A_{k, j}\) is finitely based. More explicitly, a basis consists of \(\mathcal E_{k, j}\) together with

\CaseDisplay{\begin{aligned}
xy&\approx yx, &(C)\\
x^{k+1}&\approx x^k, &(N_k)\\
x^2&\le x, &(F)\\
xy&\le x+y, &(M)\\
x_1\cdots x_k&\le\sum_{i=1}^k x_1\cdots x_{i-1}x_i^2x_{i+1}\cdots x_k. &(S_k)
\end{aligned}\quad\textup{(\ensuremath{3})}}

No minimality or practical enumeration of the large library is asserted. Its finiteness and completeness are proved below for unbounded numbers of variables.

\par\smallskip\noindent\textit{2. Soundness.}\quad 
Associativity of multiplication follows from capped addition of indices in the core and the identity property of \(e\). To verify distributivity, multiplication by \(e\) is immediate, and on the core the maps \(v_r\mapsto v_{\min(k, r+s)}\) preserve the chain join. The only remaining joins involve \(e\) and \(v_r\). Multiplication by \(v_s\) sends them respectively to \(v_s\) and \(v_{\min(k, r+s)}\), whose join is \(v_s\). Their original join is either \(e\) or \(v_0\), both of which multiply \(v_s\) to \(v_s\). Thus \(A_{k, j}\) is an ai-semiring.

Give \(e\) and \(v_0\) index zero, and \(v_r\) index \(r\) for \(r>0\). A word is \(e\) if all its letters are \(e\); otherwise its value is \(v_{\min(k,\sum r)}\). This proves \((C)\) and \((N_k)\). For \((F)\), squaring moves down the core chain and fixes \(e\). For \((M)\), a factor \(e\) leaves the other factor unchanged, which is a right summand; with both factors in the core the product lies below each factor.

For \((S_k)\), if all \(k\) factors have positive indices then the left side is \(v_k\), the additive least element. Otherwise a factor is \(e\) or \(v_0\), and repeating it leaves the product unchanged. Thus some summand equals the left side. This also proves the identity under arbitrary substitutions and repeated variables.

\par\smallskip\noindent\textit{3. Normalization and the coverage word.}\quad 
By \((C),(N_k)\) every word may be normalized so that each present variable has multiplicity in \(\{1,\ldots, k\}\). Call a variable saturated in a word if its multiplicity is \(k\). The law \((F)\) permits increasing any positive multiplicity, in a multiplication context, while decreasing the word in additive order. This includes the case of a single-variable word; empty contexts are merely omitted.

Suppose \(A_{k, j}\models t\le P\), where \(t\) is a normalized word, \(C=c(t)\), and \(P\) is a nonempty sum of normalized words. Let \(P_C\) be its summands whose content is contained in \(C\). Then

\CaseDisplay{\bigcup_{u\in P_C}c(u)=C.\quad\textup{(\ensuremath{4})}}

Indeed, if \(x\in C\) were missing from this union, assign \(x=v_0\), the other variables in \(C\) the value \(e\), and all variables outside \(C\) the value \(v_k\). The target has value \(v_0\), while every right summand is \(e\) or \(v_k\), a contradiction. This argument also proves \(P_C\ne\varnothing\).

Let

\CaseDisplay{B_C=\prod_{x\in C}x^k.}

Repeated use of \((M)\) shows that the product of all summands of \(P_C\) is below their sum, hence below \(P\). This product has content \(C\) by (4). Normalize, then increase each present exponent to \(k\) by \((F)\). We obtain the derivation

\CaseDisplay{B_C\le P.\quad\textup{(\ensuremath{5})}}

\par\smallskip\noindent\textit{4. Reduction to at most \(k-1\) unsaturated occurrences.}\quad 
If a normalized target contains at least \(k\) occurrences of unsaturated variables, commute \(k\) of them to the beginning and apply \((S_k)\) in the remaining context. Each resulting word repeats one selected variable, increasing its multiplicity by one. The sum of the deficits \(k-m_x\) over unsaturated variables strictly decreases in every branch. Moreover, each new word is below the former word by \((F)\), so it remains semantically absorbed by \(P\).

This finite branching process terminates with words having at most \(k-1\) unsaturated occurrences. To prove completeness it therefore suffices to derive \(t\le P\) for such a terminal word. Its content is still \(C\). Let \(X\) be its unsaturated variables and \(S=C\setminus X\). Put \(r=|X|\le k-1\). If \(X=\varnothing\), then \(t=B_C\) and (5) finishes. Assume below that \(X\ne\varnothing\).

\par\smallskip\noindent\textit{5. A bounded family of actual witnesses.}\quad 
First, \(t\le P_C\) holds semantically: for any assignment to \(C\), extend it by assigning \(v_k\) to every outside variable. The noncontained summands of \(P\) then have value \(v_k\), which is an additive least element and changes no nonempty join. For each \(u\in P_C\), raise every positive multiplicity of a variable in \(S\) to \(k\), obtaining \(u^\sharp\). This is a derivable inequality \(u^\sharp\le u\). When all variables in \(S\) are assigned \(e\), this change does not change any value. Consequently, at each of the \(n^r\) assignments \(\beta:X\to A_{k, j}\) with \(S\) fixed to \(e\), the sum of all \(u^\sharp\) absorbs \(t\).

At each such assignment, choose one \(u^\sharp\) for every distinct value attained among this finite family. There are at most \(n\) values. Let \(Q_0\) be the union of the chosen words over all assignments, with repetitions removed. It is nonempty, has at most \(n^{r+1}\le n^k=M\) summands, and

\CaseDisplay{t\le Q_0\quad\text{whenever all variables of }S\text{ are }e.\quad\textup{(\ensuremath{6})}}

Witness lemma. The unrestricted identity

\CaseDisplay{A_{k, j}\models t\le Q_0+B_C\quad\textup{(\ensuremath{7})}}

holds.

Proof. Consider any assignment to \(C\). If a variable of \(S\) has positive index, its \(k\)th power is \(v_k\), and \(t=v_k\) so the inequality is automatic. Otherwise all variables of \(S\) are \(e\) or \(v_0\). If every variable in \(X\) is also \(e\) or \(v_0\), then \(t=B_C\) under this assignment. In the remaining case, a variable in \(X\) has positive index. Changing all \(S\)-values to \(e\) leaves \(t\) unchanged: its value is a core product, and both \(e\) and \(v_0\) act as identities on the core. Raising these changed values back from \(e\) to \(v_0\) can only increase \(Q_0\) because \(e\le v_0\) and all terms are monotone. Hence (6) gives \(t\le Q_0\) in this case. This proves (7). \(\square\)

\par\smallskip\noindent\textit{6. Compressing saturated variables into a finite library.}\quad 
List the distinct summands of \(Q_0\) as \(q_1,\ldots, q_m\), where \(m\le M\). Each saturated variable appears with exponent either 0 or \(k\) in each \(q_i\). Partition \(S\) by the incidence sets

\CaseDisplay{J_x=\{i:x\in c(q_i)\},\qquad x\in S.}

There are at most \(2^m\) nonempty blocks of this partition (the block corresponding to the empty incidence set is retained when present). For each nonempty block \(S_J\) introduce one variable \(y_J\). In (7), identify all original variables in \(S_J\) with \(y_J\) and use \((N_k)\) to reduce positive exponents to \(k\). Leave the variables of \(X\) unchanged. The resulting absorption is valid in \(A_{k, j}\), uses at most

\CaseDisplay{r+2^m\le(k-1)+2^M=N}

variables, and every multiplicity is at most \(k\). After an injective renaming into the fixed \(N\) variables it is one of the identities in \(\mathcal E_{k, j}\).

Now substitute the nonempty word \(\prod_{x\in S_J}x\) for \(y_J\). By commutativity, its \(k\)th power is precisely \(\prod_{x\in S_J}x^k\). Thus this instance of the library identity gives back (7) syntactically. No empty block is used as a substitution, and no multiplicative identity constant has been added.

Every \(q_i\) lies derivably below an actual summand of \(P\), and \(B_C\le P\) by (5). Therefore (7) derives \(t\le P\). Section 4 derives every valid word absorption, for arbitrarily many variables. Finally distribute any polynomial identity and absorb each summand in both directions. This proves the theorem.

\par\smallskip\noindent\textit{7. Five-element consequences.}\quad 
For \(k=3\), the common finite library bound is \(N=2+2^{125}\). The bound is deliberately coarse and concerns existence of a finite basis; it is not a claim of a small explicit list.

Catalogue 7615 is \(A_{3,\allowbreak{}2}\) via

\CaseDisplay{(e, v_0, v_1, v_2, v_3)=(0,\allowbreak{}4,\allowbreak{}3,\allowbreak{}1,\allowbreak{}2).}

Catalogue 10032 is \(A_{3,\allowbreak{}3}\) via

\CaseDisplay{(e, v_0, v_1, v_2, v_3)=(3,\allowbreak{}4,\allowbreak{}0,\allowbreak{}1,\allowbreak{}2).}

Both are therefore finitely based, conditional only on independently confirming these table identifications before updating the catalogue. The already classified 7634 is \(A_{3,\allowbreak{}1}\), which provides a consistency check with its earlier separate smaller basis.

The theorem supplies the arbitrary-variable completeness argument. Truth-table checks and bounded regression below are only verification of its applications, not substitutes for the proof. The result does not assert that all subvarieties of these varieties are finitely based.

\subsection{Multiplicity witnesses and a finite basis for source 8872}\label{case:190}
Work in the constant-free ai-semiring signature. Element labels in the following tables are names, not constants in terms.

addition  multiplication 00044     11100 01044     11111 00244     11100 44434     01234 44444     01044

Let A denote this algebra. Its proper additive principal ideals are \{1\}, \{2\}, \{1,\allowbreak{}2,\allowbreak{}0\}, and \{3\}. Put K=\{3,\allowbreak{}4\}, with multiplication on K given by maximum in \(3<4\). The element 1 is a multiplicative zero. Every product of two elements from \{0,\allowbreak{}2\} is 1. An occurrence of 0 with other factors from K yields 0. An occurrence of 2 yields 2 precisely when it is the last factor and all preceding factors are 3; in every other product with just one factor from \{0,\allowbreak{}2\}, the result is 0. These rules follow from the table, are unchanged by reassociation, and prove associativity. Each multiplication map preserves the comparable joins, \(1+2=0\), and \(3+x=4\) for \(x=0\),1,\allowbreak{}2, as verified from its row or column; thus A is an ai-semiring.

\par\smallskip\noindent\textit{Word semantics and normalization.}\quad 
For a nonlinear word, an assigned 1 or at least two occurrences assigned values from \{0,\allowbreak{}2\} give value 1. If all occurrences are in K, the value is their maximum. If there is exactly one occurrence assigned 0 or 2, its value is 2 precisely when that occurrence is the last letter, is assigned 2, and all preceding letters are 3; otherwise its value is 0.

The identities

{\small\par\smallskip
\CaseLine{\(xyz \approx  yxz,\qquad x^{3} \approx  x^{2}\qquad (W)\)}
\smallskip}

hold. They permute the prefix before the last position and cap each variable's total multiplicity at two. To see the last assertion equationally, sort the prefix. A triple in the prefix reduces by \(x^{3}\)\ensuremath{\approx}\(x^{2}\). If the final variable has two prefix occurrences, move these to the end of the prefix and reduce the resulting final cube. Consequently a normal word is a sorted prefix followed by its final variable, with each variable occurring at most twice in the whole word. Linear words remain linear. This normalization holds for arbitrary word lengths; a nonfinal variable with multiplicity one is not identified with one of multiplicity two.

For a word w let c(w) be its content, lin(w) the variables occurring exactly once, last(w) its final variable, and int(w) the variables in all positions except the last.

\par\smallskip\noindent\textit{Exact absorption criterion.}\quad 
Let P have linear variable set L and nonlinear words N. Set \(D=c(P)\), I equal to the union of int(v) over v in N, and R equal to L union the set of last(v) over v in N. Put N\_C=\{v in N: c(v) is contained in C\}.

Lemma 1. A linear y satisfies y\ensuremath{\leq}P exactly when y belongs to L. A nonlinear w, with \(C=c(w)\), satisfies w\ensuremath{\leq}P exactly when the following hold:

(C) C is contained in D.

(J) L meets C or \(N_{C}\) is nonempty.

(Z) Either L meets C, or, for every x in lin(w), there is v in \(N_{C}\) in which x occurs at most once (possibly not at all).

(F) Either I meets R, or last(w) belongs to R and int(w) is contained in I.

Necessity. A missing linear y is separated by assigning \(y=0\) and every other variable 1: all right words are 1, whereas \(y=0\). A failure of (C) is separated by assigning a missing target variable 4 and every other variable 3. A failure of (J) is separated by assigning C the value 3 and outsiders 1. In the latter case the target is 3 and every right summand is 1.

Suppose (Z) fails at a variable x occurring once in w. Assign \(x=0\), all other variables of C the value 3, and outsiders 1. The target is 0. There is no linear right variable in C, and each contained nonlinear right word has at least two occurrences of x. Thus all right words are 1.

For failure of (F), I and R are disjoint. Assign I the value 3 and R the value 2. Every original right summand is then 2. If last(w) does not belong to R, assign it 3 when already in I and 1 otherwise; complete remaining variables arbitrarily by 3. The target cannot be 2. Otherwise some prefix variable is outside I. If it belongs to R it already has value 2; if not, assign one such variable 1. Complete other variables by 3. Again the target is not 2, since a word valued 2 has all its prefix occurrences valued 3. Thus each failure gives a separating assignment.

Sufficiency. It suffices to check membership in each proper additive principal ideal. If \(P=3\), every variable of D is 3, so (C) gives \(w=3\). If P lies in \{1,\allowbreak{}2,\allowbreak{}0\} and w lies in K, all target variables are in K; the contained or linear witness in (J) then has value in K, contradiction.

Suppose \(P=1\) but w is not 1. If w has all variables in K, (J) again contradicts \(P=1\). In the remaining case there is exactly one occurrence valued in \{0,\allowbreak{}2\}, belonging to a variable x in lin(w); all other target occurrences are in K. A linear variable in C has value different from 1 and already contradicts \(P=1\). Otherwise (Z) supplies a contained right word with at most one occurrence of x, hence with value different from 1. This is also a contradiction.

Finally, \(P=2\) forces all variables in R to be 2 and all variables in I to be 3. Their intersection is therefore empty, and (F) makes the target final variable 2 and every prefix occurrence 3. Hence \(w=2\). Linear sufficiency is literal occurrence. This proves the criterion for arbitrary variables, words and polynomial sizes.

\par\smallskip\noindent\textit{Four-witness derivation chains.}\quad 
For a nonempty content C and chosen last variable t in C, a saturated word has last t and exactly two occurrences of every variable in C. The word laws (W) determine it uniquely up to the fixed ordering.

Lemma 2. Every valid nonlinear absorption w\ensuremath{\leq}P follows from (W) and a finite chain of valid absorptions each having at most four right words. Each right word used in a step is an original summand of P or a previously derived word.

Proof. Write \(t=last(w)\). First suppose I is disjoint from R. By (F), t belongs to R, int(w) is contained in I, and t does not occur in the prefix of w. By (J), either a nonlinear contained witness exists or a linear variable lies in C. Any contained nonlinear witness must have final variable t: any other final variable would lie in int(w), hence in I intersect R. It cannot repeat t. A linear witness can likewise only be t.

If a contained nonlinear word v exists, start with v and double every prefix variable occurring only once. Each doubling is a valid one-right-word absorption by Lemma 1: content and prefix support are unchanged, and the set of simple variables only decreases. The final t remains simple. If instead use of a linear t is necessary, choose a prefix variable x of w and an original nonlinear right word a having x in its prefix, which exists because x is in I. The absorption x\ensuremath{{}^{2}}t\ensuremath{\leq}t+a is valid: content is covered, linear t supplies (J) and (Z), and (F) is supplied by x in int(a) and final t in the linear part.

Extend the current prefix support, one missing target prefix variable x at a time, by inserting two copies of x. An original nonlinear right word a with x in its prefix supplies content coverage and (F). The current word s supplies (J) and the only remaining simple target variable t for (Z). Thus the new word is validly absorbed by s+a. After this process the target content is complete, all prefix variables have multiplicity two, and t has multiplicity one.

Now for each prefix variable x that is simple in w, remove its extra occurrence. If L meets C, the relevant linear witness t suffices for (Z). Otherwise choose a contained nonlinear original witness a in which x occurs at most once, as provided by (Z) for w. The step is absorbed by s+a (or s+t): current s supplies all previous simple-variable requirements, the additional witness supplies the newly simple x, and the final/prefix constraints remain valid. This reaches w, with at most two right words per step.

Now suppose I meets R. Choose a variable g in the intersection, a nonlinear right word a having g in its prefix, and an original summand b that is the linear g or has final variable g. Every right polynomial containing a and b automatically has the intersection alternative in (F).

Choose a witness v for (J), either linear in C or nonlinear with content contained in C. Choose an original summand r containing t, using (C). Take as initial target the saturated word on c(v) union \{t\}, with final t. Its content is covered by v and r, (J) is supplied by v, (Z) is vacuous because no variable is simple, and (F) is supplied by a, b. It is therefore validly absorbed by a+b+v+r.

For each missing variable x of C, take an original summand r containing x and insert two copies of x into the current saturated word s. The new word is absorbed by s+a+b+r: s supplies (J), (Z) remains vacuous, a, b supply (F), and r covers the new variable. This produces the saturated word on all of C.

Lastly, for each variable x simple in w, remove one occurrence, retaining final variable t. If L meets C, take a linear witness v in that intersection. Otherwise choose v in \(N_{C}\) in which x occurs at most once. The resulting word is absorbed by s+a+b+v. Its earlier simple variables are witnessed by s, the newly simple x by v (or by the linear witness), content is unchanged, and a, b retain (F). When \(x=t\), remove its prefix occurrence and leave the final t. This completes the chain and proves the lemma. Repeated right summands may be omitted; no step has more than four.

\par\smallskip\noindent\textit{Reversible finite compression.}\quad 
Fix 247 variable symbols. Let Normal be the finite set of all nonempty normal words on these symbols, with each variable occurring at most twice; thus every word has length at most 494. Let Gamma consist of every absorption w\ensuremath{\leq}\(u_{1}\)+...+\(u_{k}\), 1\ensuremath{\leq}k\ensuremath{\leq}4, in these normal words which satisfies Lemma 1. This explicitly defines a finite decidable set, although it is large and is not enumerated here.

Lemma 3. Every valid absorption of one word by at most four right words follows from Gamma and (W).

Proof. Normalize the at most five words. Keep each variable that is final in any one of these words separate; this includes any linear word, and there are at most five such variables. For every other variable, record its multiplicities (0,\allowbreak{}1,\allowbreak{}2) in the prefixes of the at most five words. The all-zero pattern is absent. There are at most \(3^{5}-1=242\) patterns. Replace each group of variables with the same pattern by one new symbol whose multiplicity in each prefix is their common multiplicity, not the sum of their multiplicities. This uses at most \(5+242=247\) symbols.

This operation is not ordinary variable identification: identifying two simple occurrences would incorrectly make the resulting symbol multiple. Instead validity follows from Lemma 1. Content coverage and contained-word membership are unchanged, since members of a group occur in exactly the same words. Linear variables and all final variables are preserved separately, so the linear alternatives and the sets I intersect R are unchanged. A new pattern symbol is simple in the target exactly when all its group members were simple there; its occurrence at most once in a right witness is likewise equivalent for every group member. Thus (Z) is preserved. The last-variable and prefix-content conditions in (F) are preserved by the same pattern incidence. Every condition of Lemma 1 therefore holds for the compressed absorption, which lies in Gamma up to renaming.

Conversely, replace each pattern symbol by the nonempty product of all original variables in its group. A group only occurs in prefixes; its members occur in precisely the same words and with the same multiplicity in each. A symbol occurring once expands to one copy of each group member; a symbol occurring twice expands to two copies. Prefix permutation recovers each original normal word. No group contains a final variable, so final positions and their multiplicities remain correct. Thus the substitution followed by (W) derives the original absorption from Gamma.

Theorem. The two identities (W), together with Gamma, form a finite basis for A. The multiplicative opposite has the reversed basis.

Proof. Soundness follows from the semantics and definition of Gamma. Apply Lemma 2 to a valid nonlinear absorption and Lemma 3 to each step; linear absorptions are literal summands by Lemma 1. Distribute arbitrary terms and derive all summand absorptions in both directions. The proof places no bound on the original variables, word lengths or numbers of summands. Word reversal proves the opposite assertion.

The bounds 247 variables, length 494 and four right summands are proved, not inferred from finite tests. No small enumerated, minimal, independent or hereditary basis is claimed. Computational tests are supplementary error checks. Transfers to other catalogue types require complete two-way HSP equality and independent verification.

\subsection{A finite basis for the five-element chain monoid semiring 7634}\label{case:191}
Work in the constant-free signature (+, multiplication), relative to the additively idempotent semiring axioms. Inequality u \ensuremath{\leq} v means \(u+v=v\).

\par\smallskip\noindent\textit{The algebra and the theorem.}\quad 
Let A=\{z, b, a, d, e\}, with addition maximum in the chain \(z<b<a<d<e\). The element d is a multiplicative identity. Write \(v_{0}=e\), \(v_{1}=a\), \(v_{2}=b\), \(v_{3}=z\), and put \(v_i v_j=v_{\min(3,i+j)}\). Thus the core \{e, a, b, z\} is a truncated additive-weight monoid with identity e; the identity d of the whole algebra is distinct from e. In catalog type 7634, \((e,a,b,z,d)=(4,1,2,3,0)\). These element names are not constants of the term signature.

Theorem. The identities

{\small\par\smallskip
\CaseLine{\(C: xy = yx,\)}
\CaseLine{\(N: x^{3} = x^{4},\)}
\CaseLine{\(F: x^{2} \leq  x,\)}
\CaseLine{\(M: xy \leq  x+y,\)}
\CaseLine{\(S: xyz \leq  x^{2}yz+xy^{2}z+xyz^{2},\)}
\CaseLine{\(L: xy \leq  x+(xy)^{3}\)}
\smallskip}

together with the sixty instances of (T) below form a basis for all identities of A, relative to the semiring axioms. The basis has at most 66 identities; no minimality is asserted. In particular, A is finitely based. All statements and constructions below allow arbitrarily many variables and summands.

To specify (T), independently take A, B, C, D to be omitted or to be the cubes of four distinct new variables. Write \(s=ABCD\). For i in \{0,\allowbreak{}1\} and j in \{0,\allowbreak{}3\}, include

{\small\par\smallskip
\CaseLine{\((T)  sxy \leq  ABx^{2} + ACx^{i} y^{j} + sx^{3}y^{3},\)}
\smallskip}

provided the middle summand is nonempty. An exponent 0 means that substring is omitted. This is a finite specification of ordinary constant-free identities, not permission to use empty substitutions. There are 64 possible choices, of which four have an empty middle summand and are excluded. Each nonempty cubed block can subsequently receive a nonempty word substitution.

\par\smallskip\noindent\textit{Values, soundness, and normalization.}\quad 
Multiplication is commutative and associative: d is a global identity, and on the core it adds indices with truncation at 3. It is monotone in each coordinate of the displayed additive chain. For example, the row of a in the order z, b, a, d, e is z, z, b, a, a; the rows of b and z are z, z, z, b, b and z, z, z, z, z. The rows of d and e are respectively z, b, a, d, e and z, b, a, e, e. Consequently both distributive laws hold.

Assign index 0 to d and index i to \(v_{i}\). A word has value d exactly when all its variables have value d. Otherwise it has core value \(v_{min}\)(3, sum \(n_{x}\) \(lambda_{x}\)), where \(n_{x}\) is the multiplicity of x and \(lambda_{x}\) its assigned index.

C, N, F follow at once. For M, if one factor is d it is a summand of the right side; otherwise the product index is at least each factor index and the inequality follows. If all three factors in S have positive index, its left side is z. Otherwise repeating a factor of index 0 leaves the product unchanged (also when all factors are d). For L, if x is a core element, \(xy\leq x\). If \(x=d\) and y is d or has positive index, \(xy\leq d=x\); if \(y=e\), then \((xy)^{3}=e\). Thus L is sound.

For (T), if a cubed block has positive index, the left side is z. Otherwise each block is d or e. If x, y both have index 0, the third summand equals the left side. If x has index 0 and y has positive index, the first summand is d or e and bounds the left side. If x has positive index and y has index 0, the second summand has either the same positive index as x or value d or e, and again bounds the left side. Finally suppose x, y both have positive indices. If x has index 1, the first summand is b, which bounds xy; if x has index at least 2, the left side is z. These cases prove every instance.

C, N sort every word and cap every positive multiplicity at 3. No variable disappears, and no empty word is introduced. F, in a multiplication context, permits increasing a present multiplicity.

\par\smallskip\noindent\textit{Exact absorption criterion.}\quad 
For a word w, let \(C=c(w)\). For a nonempty polynomial P define P\_C=\{u in P:c(u) subset C\}. Define the core condition K(w, \(P_{C}\)) by

{\small\par\smallskip
\CaseLine{\(min(3, sum_{x} n_{x}(w) lambda_{x})\)}
\CaseLine{\(\geq  min_\{u in P_{C}\} min(3, sum_{x} n_{x}(u) lambda_{x})\)}
\smallskip}

for all assignments \(lambda_{x}\) in \{0,\allowbreak{}1,\allowbreak{}2,\allowbreak{}3\}. Set K(w, empty) false.

Lemma 1. A satisfies \(w\leq P\) if and only if

{\small\par\smallskip
\CaseLine{\((coverage) union_\{u in P_{C}\} c(u) = C;\)}
\CaseLine{\((core)\qquad K(w, P_{C}).\)}
\smallskip}

Necessity of coverage follows by assigning a missing covered variable x the value e, the other target variables d, and every outside variable z. Then \(w=e\), every contained right word is d, and every noncontained right word is z. The right side is at most d, a contradiction. For core necessity use only core assignments and set all outside variables to z. Noncontained words then have value z and do not help absorb a target with smaller index. If \(P_{C}\) is empty, take all target variables e and all outside variables z.

For sufficiency, if \(w=d\), every variable in C is d and every contained right word is d. Such a word exists by coverage, so \(P\geq d\). If \(w=e\), all its variables lie in \{d, e\} and at least one, say x, is e. Coverage supplies a contained right word including x, which is e. If w has positive core index i, then the core condition supplies a contained right word with index at most i. If that word evaluates d rather than a core element, it is still above w. Otherwise the index comparison is exactly the required additive inequality. This proves sufficiency.

Only a small set of index assignments is needed in K. Assign outsiders index 3 and disregard noncontained words. For target index below 3, there are either no positive-index variables, one positive-index variable with index 1 or 2, or two variables of index 1 each. Assignments where the target index is 3 impose no condition. This finite description is valid for any finite set C.

\par\smallskip\noindent\textit{A common upper bound for inserting variables.}\quad 
Fix P and a nonempty target content C satisfying coverage. A word is available if its absorption into P is derivable from the displayed basis. All original right words are available. M multiplies any two available words: \(uv\leq u+v\leq P\). Multiply the finitely many words in \(P_{C}\). Its content is C. Normalize and apply F in contexts to raise every positive exponent to 3. Thus the word

{\small\par\smallskip
\CaseLine{\(B_C=\prod_{x\in C}x^3\)}
\smallskip}

is available.

Lemma 2. If the normalized word u is an available right summand, c(u) subset C, and its multiplicities are no greater than those of a normalized target t with content C, then t is available.

If \(t=u\) there is nothing to prove. Otherwise t is, by commutativity, uq for a nonempty word q. L gives

{\small\par\smallskip
\CaseLine{\(t=uq \leq  u+(uq)^{3} = u+B_{C} \leq  P.\)}
\smallskip}

The equality with \(B_{C}\) follows from C, N. This proves the lemma. The quotient is never used as an empty substitution. Unlike an unrestricted insertion rule, the argument explicitly uses the coverage witnesses.

\par\smallskip\noindent\textit{Arbitrary-variable completeness.}\quad 
Assume \(w\leq P\) satisfies Lemma 1, and normalize w and every right word. Build the available \(B_{C}\) as above. Call a variable saturated when its multiplicity is 3. If the target has at least three occurrences of unsaturated variables, choose three such occurrences and apply S in the remaining context. This bounds the target by three words, each with one chosen occurrence repeated. Normalize. The sum of the deficits 3-\(n_{x}\) over x in C strictly decreases in each branch. Moreover increasing multiplicities preserves the core condition. Repetition therefore reaches words with at most two unsaturated occurrences.

Write s for the product of the cubes of the saturated variables; s is an omitted context when there are none.

For a terminal s, a contained right word is a multiplicity divisor. For a terminal sx, give x index 2 and saturated variables index 0; the core condition supplies a contained right word with x-exponent at most 1. For \(sx^{2}\), give x index 1 to obtain an exponent at most 2. In each case all other exponents of the normalized witness are at most 3, so Lemma 2 finishes.

The remaining terminal is sxy with x different from y. If a contained right word divides it, apply Lemma 2. Otherwise give x, y both index 1. A right witness has total multiplicity at most 2 on \{x, y\}. Since it does not divide sxy, its exponents must be (2,\allowbreak{}0) or (0,\allowbreak{}2). Interchange x, y if needed and call the witness u, with exponents (2,\allowbreak{}0). Next give x index 2 and y index 0. A witness v has x-exponent i at most

1. Its y-exponent can be raised to \(j=3\) when y occurs, or is \(j=0\) if y is absent.

Partition the saturated variables into four sets: those in both u, v; those in u alone; those in v alone; and those in neither. Let A, B, C, D be the products of their cubes. By normalization and F, raising only variables already present in each witness gives available words

{\small\par\smallskip
\CaseLine{\(ABx^{2},\qquad ACx^{i} y^{j}.\)}
\smallskip}

The second word is nonempty because v was nonempty. The previously constructed coverage word is \(sx^{3}y^{3}\). The corresponding instance of (T) therefore gives \(sxy\leq P\). Each nonempty block is realized by substituting the product of its variables into the relevant cubed variable, and normalizing with C, N. An empty block is handled by the already listed omitted-block identity. Thus no empty substitution is used. In particular this step does not incorrectly assume that a witness remains adequate after deleting y from the target content.

A linear target is the terminal sx with s absent; its divisor witness is the actual linear right word x.

Every valid word absorption is therefore derivable. Distribute arbitrary terms into finite nonempty sums of words, derive each summand of either side of a valid identity below the other side, and use idempotent addition to obtain equality. This proves the theorem without a bound on variables.

\subsection{平方展开交换分支的最后两组完整有限基}\label{case:192}
\par\smallskip\noindent\textit{1. 两个生成半环.}\quad 
令 \(M=\{0, b, a, e, f\}\)，加法为链 \(0<b<a<e<f\) 的最大值。乘法交换，\(0\) 为零，\(e\) 为幺元，

\CaseDisplay{a^2=ab=b^2=0,\qquad f^2=f,\quad fa=a,\quad fb=a.}

它是目录来源 7644，按上述顺序到目录标号为 \((3,\allowbreak{}2,\allowbreak{}1,\allowbreak{}0,\allowbreak{}4)\)。乘法的意义是：唯一零平方因子若为 \(b\)，则遇到 \(f\) 时其值上升为 \(a\)；再次乘以 \(f\) 不再改变。两个零平方因子的乘积为零。这证明结合性；乘法平移单调，故保持链的并。

令 \(N=\{0, a, e, b, t\}\)，加法为

\CaseDisplay{0<a<e<t,\qquad 0<a<b<t,}

的并。在 \(\{0, a, e\}\) 上，\(0\) 是零、\(e\) 是幺元、\(a^2=0\)；任何包含 \(b\) 或 \(t\) 的二因子积均为 \(t\)。它对应来源 7924，标号映射为 \((2,\allowbreak{}1,\allowbreak{}0,\allowbreak{}3,\allowbreak{}4)\)。乘法结合，各平移保持所列的并，因此也为交换加法幂等半环。元素 \(0\) 不是整个 \(N\) 的乘法零。

\par\smallskip\noindent\textit{2. 完整基.}\quad 
两组基都包含

\CaseDisplay{xy\approx yx,\qquad xy\approx x^2y+xy^2.\quad\textup{(\ensuremath{E})}}

定理 1。 \(M\) 的完整基是在 (E) 外加入

\CaseDisplay{x^2\leq x,\qquad xy\leq x^2+y^2,\qquad
xy^2\leq y^2+(xy)^2.\quad\textup{(\ensuremath{M})}}

定理 2。 \(N\) 的完整基是在 (E) 外加入

\CaseDisplay{xy\leq x+(xy)^2,\qquad
x^2y\leq x^2+y^2z^2,\qquad
x^2\leq x+xy.\quad\textup{(\ensuremath{N})}}

均相对于通常加法幂等半环公理而言。

先验证有效性。(E) 是零平方层与单位层的直接分解。对 \(M\)，平方收缩由链和平方表成立；若乘积为 \(a\) 或 \(b\)，另一个因子的平方在单位层而吸收它，得到平方和吸收。最后一个 (M) 式中，左侧为零平方层值时，\(y^2\) 在单位层，因而吸收它；左侧在单位层时，其平方保留原值。

对 \(N\)，第一个 (N) 式的左侧为 \(a\) 时，每个非零因子都不小于 \(a\)；左侧为稳定值时由平方项吸收。第二式中，\(x^2=0\) 时，除非 \(y=b, t\)，左侧为 \(0\)；若 \(y=b, t\)，则 \(y^2z^2=t\)。若 \(x^2=e\)，左侧非顶值由 \(e\) 吸收，顶值仍来自 \(y=b, t\)；\(x^2=t\) 时直接成立。第三式只需注意 \(x=b\)：此时 \(x^2=xy=t\)。

\par\smallskip\noindent\textit{3. 精确吸收判据.}\quad 
令右多项式 \(P\) 的一次项变量集为 \(L\)，非线性词族为 \(\mathcal R\)。对非线性左词 \(w\)，记

\CaseDisplay{C=c(w),\quad S(w)=\{x:m_x(w)=1\},\quad
\mathcal R_C=\{v\in\mathcal R:c(v)\subseteq C\}.}

两者的线性判据都是 \(x\leq P\iff x\in L\)。

在 \(M\) 中，非线性吸收式的条件为

\CaseDisplay{C\subseteq L\cup\bigcup_{v\in\mathcal R_C}c(v),\quad\textup{(\ensuremath{1})}}

以及对每个 \(x\in S(w)\)，下列三项至少一项成立：

\CaseDisplay{\begin{array}{ll}
\text{(i)}&L\cap(C\setminus\{x\})\ne\varnothing;\\
\text{(ii)}&\exists v\in\mathcal R_C:\ m_x(v)=0;\\
\text{(iii)}&C\setminus\{x\}\subseteq
 \displaystyle\bigcup_{\substack{v\in\mathcal R_C\\m_x(v)=1}}c(v).
\end{array}\quad\textup{(\ensuremath{2})}}

在 (iii) 中，\(C\setminus\{x\}\) 必非空，因为 \(w\) 是非线性且 \(x\) 只出现一次。

在 \(N\) 中，条件为

\CaseDisplay{L\cap C\ne\varnothing\quad\text{或}\quad\mathcal R_C\ne\varnothing,\quad\textup{(\ensuremath{3})}}

\CaseDisplay{C\subseteq\bigcup_{v\in\mathcal R}c(v),\quad\textup{(\ensuremath{4})}}

以及

\CaseDisplay{L\cap C\ne\varnothing\quad\text{或}\quad
\forall x\in S(w)\ \exists v\in\mathcal R_C:\ m_x(v)\leq1.\quad\textup{(\ensuremath{5})}}

注意 (4) 要求每个左变量出现在某个非线性右词中；仅作为一次项出现不够。

\par\smallskip\noindent\textit{必要性.}\quad 
一次项判据：在 \(M\) 中令左变量为 \(b\)、其余为 \(0\)；在 \(N\) 中令左变量为 \(a\)、其余为 \(0\)。没有相同一次项时右侧为零。

对 \(M\)，若 (1) 失败，令缺少覆盖的变量为 \(f\)、其他内部变量为 \(e\)、外部变量为 \(0\)。若 (2) 的三种见证均失败，选一个没有被恰含一次 \(x\) 的内部右词覆盖的 \(z\in C\setminus\{x\}\)，令 \(x=b, z=f\)，其余内部变量为 \(e\)、外部变量为 \(0\)。左词为 \(a\)，右侧一次项至多为 \(b\)，右侧非线性词要么重复 \(x\) 而取 \(0\)，要么恰含一次 \(x\) 但不含 \(z\) 而取 \(b\)。因此右侧不吸收左侧。

对 \(N\)，若 (3) 失败，令内部变量为 \(e\)、外部变量为 \(0\)，左侧为 \(e\)，右侧为 \(0\)。若 (4) 失败，将没有在非线性右词中出现的 \(x\in C\) 赋为 \(b\)，其余所有变量赋为 \(0\)。左词为 \(t\)，右侧至多为 \(b<t\)。若 (5) 失败，令相应 \(x=a\)、其他内部变量为 \(e\)、外部变量为 \(0\)，左侧为 \(a\)，右侧为 \(0\)。

\par\smallskip\noindent\textit{4. 从有限基推出充分性.}\quad 
(E) 给出 \(x^3=x^2\)，以及规范形

\CaseDisplay{H(C)=\prod_{y\in C}y^2,\qquad K_x(C)=xH(C\setminus\{x\}),}

每个非线性词是这些词的非空和，每个分项保持原内容，并且其唯一可能的一次变量来自原词的一次变量。对每个非线性词 \(v\)，还由 (E) 得 \(v^2\leq v\)。

对 \(M\)，平方收缩、平方和吸收以及 (1) 给出 \(H(C)\leq P\)。若 (2)(i) 提供一次项 \(g\ne x\)，把目标词写成 \(ug^2\)，由第三个 (M) 式得

\CaseDisplay{K_x(C)\leq g^2+(ug)^2=g^2+H(C)\leq P.}

若 (ii) 提供不含 \(x\) 的右词 \(v\)，令 \(u=K_x(C)\)。在第三个 (M) 式中代入 \(u, v\)，由于 \(c(v)\subseteq C\setminus\{x\}\)，指数截断给出

\CaseDisplay{K_x(C)=uv^2\leq v^2+(uv)^2=v^2+H(C)\leq P.}

若 (iii) 成立，从每个相关右词选一个恰含一次 \(x\) 的规范分项 \(xH(S_i)\leq P\)。平方和吸收给出稳定块合并式，将其乘以 \(x\)，反复得到

\CaseDisplay{xH\Bigl(\bigcup_i S_i\Bigr)\leq\sum_i xH(S_i)\leq P.}

由 (iii)，并集为 \(C\setminus\{x\}\)，完成 \(M\) 的证明。

对 \(N\)，如果有内部非线性锚点，则其平方提供 \(H(S)\leq P\)。否则由 (3) 选 \(g\in L\cap C\)，由 (4) 选包含 \(g\) 的非线性右词 \(gv\)，第三个 (N) 式给出

\CaseDisplay{g^2\leq g+gv\leq P.}

因此总能得到内部平方锚点。对任一待加入的 \(z\in C\)，(4) 提供右词 \(zt\)。由第二个 (N) 式及 (E) 得

\CaseDisplay{H(S)z^2\leq H(S)z\leq H(S)+z^2t^2\leq P.}

逐一加入，得到 \(H(C)\leq P\)。单变量非线性右词先写成 \(z^2\)，取 \(t=z\)。

若 \(L\cap C\ne\varnothing\)，任取其中的 \(g\)，将左非线性词写为 \(gv\)，由第一个 (N) 式得到 \(gv\leq g+(gv)^2=g+H(C)\leq P\)。否则 (5) 为每个 \(x\in S(w)\) 提供一个不重复 \(x\) 的内部右词。恰含一次 \(x\) 时，从它的规范分项开始；不含 \(x\) 时，先平方得到 \(H(S)\leq P\)，再由全局覆盖选含 \(x\) 的右词 \(xt\)，用第二个 (N) 式得到

\CaseDisplay{xH(S)\leq H(S)+x^2t^2\leq P.}

设已有 \(r=xH(S)\leq P\)，且需加入 \(z\in C\)。选右词 \(zt\)，由前两个 (N) 式得

\CaseDisplay{rz^2\leq r+(rz^2)^2=r+r^2z^2
\leq r+r^2+z^2t^2\leq P.}

最后一步用 \(r^2\leq r\)。最终得到 \(K_x(C)\leq P\)。

因此两组判据的充分性与基的完整性均已证明。任意恒等式可分解为双向逐词吸收，两个定理得证。

\par\smallskip\noindent\textit{5. 应用边界.}\quad 
来源 7644 的等簇类有 1 个类型，7924 的等簇类有 7 个类型。二者也都满足平方和吸收：第一个将它作为基的一条；第二个由其运算定义可直接核对，所以均适用平方定向乘法的项等价定理。

在与此前证书合并并逐个核验后，可检查原剩余目录中满足 \(xy=yx\)、\(xy=x^2y+xy^2\) 的分支是否已全部覆盖。此分支的闭合不等于一般五元素问题的解决。

\subsection{Finite bases for 7694 and 7695 by tail-sensitive simple-variable recovery}\label{case:193}
All bases below are relative to the additively idempotent semiring axioms, in the signature (+, multiplication), without constants. The finite families are specified effectively, not enumerated as small lists.

\par\smallskip\noindent\textit{1. Tables and word values.}\quad 
Both sources have addition rows 00004, 01114, 01224, 01234, 44444, in element order 0,\allowbreak{}1,\allowbreak{}2,\allowbreak{}3,\allowbreak{}4. Their additive order is \(3<2<1<0<4\).

{\small\par\smallskip
\noindent 7694 multiplication    7695 multiplication\par
\CaseLine{\ttfamily 0 1 2 3 4\qquad 0 1 2 3 4}
\CaseLine{\ttfamily 2 3 3 3 2\qquad 2 3 3 3 2}
\CaseLine{\ttfamily 2 3 3 3 2\qquad 2 3 3 3 2}
\CaseLine{\ttfamily 3 3 3 3 3\qquad 3 3 3 3 3}
\CaseLine{\ttfamily 0 1 2 3 4\qquad 4 1 2 3 4}
\smallskip}

Put K=\{0,\allowbreak{}4\}. Direct induction gives these formulas, including literals. Any factor valued 3, or two positions valued in \{1,\allowbreak{}2\}, gives value 3. With exactly one \{1,\allowbreak{}2\}-valued position, the product is 1 only when that position is final and valued 1; it is 2 in all other cases. With all factors in K, source 7694 gives the value of its last factor, whereas 7695 gives 4 if any factor is 4 and otherwise gives 0. These different K-valued clauses require different coverage conditions.

\par\smallskip\noindent\textit{2. Exact absorption criteria.}\quad 
For a target word w and nonempty polynomial P of original right words, write \(C=c(w)\), \(t=the\) last variable of w, and U=\{u in P : c(u) is contained in C\}. A variable is simple when its total number of occurrences in the indicated word is exactly one.

For 7694, \(w\leq P\) holds exactly when:

\textbullet{} (T) some u in U ends in t;

\textbullet{} (S) for every simple variable x of w, some u in U omits x or contains x exactly once. When \(x=t\) is simple, a witness containing x must also end in x.

For 7695 replace (T) by (C): the union of the contents of U is C, and retain the same condition (S). A linear target therefore requires its literal summand: a nonempty word contained in \{x\} cannot omit x, and having exactly one x forces the literal x.

Necessity of (T): assign \(t=4\), all other target variables 0, outsiders 3. The target is 4, while right words not ending in t or not contained in C cannot reach 4. Necessity of (C): for each uncovered x assign \(x=4\), other target variables 0, outsiders 3; again the target is 4. For a failed (S) at a nonfinal simple x, assign \(x=2\), the other target variables 0, outsiders 3. The target is 2 and all right words are 3. For a failed (S) at the simple final t, instead assign \(t=1\). The target is 1; right words either repeat t and give 3, or have t once nonfinally and give 2. Thus every failed condition has a separator.

For sufficiency, target value 3 is least. If all its factors lie in K, (T) supplies exactly its last value in 7694. In 7695, (C) supplies a contained witness with a 4-valued variable whenever the target is 4; an all-zero target is handled by any contained witness. If the target has value 1 or 2, its sole \{1,\allowbreak{}2\}-valued occurrence belongs to a simple variable x. A contained (S)-witness omitting x is K-valued and hence above the target. A witness containing x once has value 1 or 2; the additional final-position condition guarantees value 1 when necessary. This proves the criteria for arbitrary variables and assignments.

\par\smallskip\noindent\textit{3. Structural identities and normal forms.}\quad 
Both algebras satisfy:

{\small\par\smallskip
\CaseLine{\((W) xyz = yxz\)}
\CaseLine{\((P) x^{3} y = x^{2} y\)}
\CaseLine{\((Q) xy^{3} = xy^{2}\)}
\CaseLine{\((R) x^{3} = x^{2}\)}
\CaseLine{\((D) xy = x^{2}y + xy^{2}\)}
\CaseLine{\((L) x^{2} \leq  x\)}
\smallskip}

Additionally 7694 satisfies (A) \(x^{2}y\leq y\). Source 7695 instead uses (M) \(x^{2}y^{2}\leq x^{2}+y^{2}\) and (J) \(x^{2}y^{2}=y^{2}x^{2}\). Soundness follows from the word-value formulas. In particular D splits the unique possible \{1,\allowbreak{}2\}-valued occurrence between its factors; squaring its factor gives 3, whereas the other branch retains its value and final position. If both factors are K-valued, squaring preserves them. The zero-valued cases remain 3 in both branches.

Identity W permits any permutation of the prefix before the last position. Collect equal prefix letters and apply P to truncate nonfinal multiplicities to at most two. Place copies of the last variable just before its final occurrence, and use Q or R to truncate its total multiplicity to two. Permute the prefix into a fixed order. Thus a normal word is determined by its last variable and by each variable's multiplicity, truncated to one or two. Every step is an actual substitution into W, P, Q, R with nonempty variable images; external contexts may be empty. No zero-length identity factor is used.

For 7695 a word all of whose variables are repeated can also change its last variable. First use W to write it as a product of squares, ending in the square of its old last variable. Adjacent applications of J move any chosen square to the end. Then W restores the normal prefix order. In particular all saturated full-content words are equal in 7695. We do not use this tail-changing operation in 7694.

\par\smallskip\noindent\textit{4. Constructing a saturated guard.}\quad 
Fix a valid absorption \(w\leq P\). A saturated word has every variable repeated. L, applied at a single occurrence in a context, doubles any simple variable and produces a word below the original word. Normalization then gives exactly two occurrences of each repeated variable. This operation keeps the last variable fixed.

For 7694 choose the contained same-tail witness from T and saturate it using L. Insert a square of each missing target variable before the current word using A. The resulting guard g has full content C, the target's last variable t, and every variable repeated; \(g\leq P\).

For 7695 saturate every contained witness that is needed to cover C. Merge two saturated words a, b by the absorption \(ab\leq a+b\). It follows from M with \(x=a\), \(y=b\), since W, P, Q, R give \(a^{2}=a\) and \(b^{2}=b\) for saturated words. Repeating this merge supplies full content. Use J as described above to give this guard the target's last variable t. Thus again \(g\leq P\), with full content and all multiplicities two. Coverage requires only finitely many of the original right words, even for an arbitrarily large original polynomial.

\par\smallskip\noindent\textit{5. Recovering individual simple variables and splitting.}\quad 
For each simple variable x of w choose an original witness u from S. Double every variable of u except x using L and normalize, obtaining \(u_{x}\leq u\leq P\). This word omits x or contains it once, with the required final position if \(x=t\). Let \(w_{x}\) be obtained by doubling every variable of w except x, and normalizing with its last variable unchanged.

The local two-right-word absorption

{\small\par\smallskip
\CaseLine{\(w_{x} \leq  g+u_{x}\)}
\smallskip}

is valid by the exact criterion. The guard supplies T or C, and \(u_{x}\) supplies S for the sole simple variable. Section 6 derives all these local absorptions from one fixed finite family.

Now apply D to a word with at least two simple variables, splitting it as XY between two simple occurrences. Both factors are nonempty. Each summand \(X^{2}Y\) and \(XY^{2}\) has strictly fewer simple variables. Both preserve full content and the literal last variable. Every surviving simple variable retains whether it is the last variable. Normalize and repeat. The finite splitting tree ends in words with zero or one simple variable. Zero-simple leaves equal g by W, P, Q, R; a leaf with sole simple variable x equals \(w_{x}\). Consequently every leaf is below P, and the equality from D gives \(w\leq P\). Targets already having zero or one simple variable use the corresponding leaf immediately.

\par\smallskip\noindent\textit{6. Reversible finite compression.}\quad 
Each local recovery has at most three words. Reserve x and the last variable of each word: at most four reserved variables. Every other variable is either absent or repeated in each word. Group these nonreserved variables by their nonzero incidence pattern across the three contents. There are at most \(2^{3}-1=7\) such patterns.

Map all variables in a class to one fresh variable, reserve the others separately, and normalize with W, P, Q, R. This is a genuine identification substitution into a valid identity, so the compressed absorption is valid. It uses at most eleven variables, has word length at most twenty-two, and has at most two right words.

For the inverse substitution, replace a class variable by the nonempty product of its members. A class variable is repeated in every word where it occurs and is never its last variable. Its two occurrences therefore expand to two occurrences of every member, all in the permutable prefix. W collects those copies and restores the normal prefix order. The reserved last variables and x retain their exact multiplicities and final positions. This recovers all three original words equationally. In particular common multiplicity is preserved; several single occurrences are never treated as an original repeat.

Fix eleven variables. For each source let F be all valid absorptions with one normal word on the left, one or two normal words on the right, and each word of length at most twenty-two. This is an explicitly specified finite set, decidable by Section 2. Join F to the seven displayed structural identities for 7694, or the eight for 7695. The compression and inverse substitution derive every local recovery. Sections 4 and 5 therefore derive every valid word absorption without a bound on its original variables, word lengths or summand count. Distributivity and summandwise absorption in both directions prove all identities. These are finite bases for the two algebras; reversing words gives bases for their opposites. No minimality or small enumerated list is claimed.

Independent finite checks and actual substitution replay support this unrestricted proof. Catalogue transfers require complete verified HSP inclusions in both directions, and fingerprint matches alone are insufficient.

\subsection{Finite bases for 10153 and 10156 by contained witnesses and simple-head recovery}\label{case:194}
Signature (+, multiplication), no constants, relative to the additively idempotent semiring axioms. A simple variable occurs exactly once. Each finite family is effectively specified; a small enumeration is not claimed.

\par\smallskip\noindent\textit{1. Two tables and all word values.}\quad 
The common addition table, followed by the two multiplication tables, is:

{\small\par\smallskip
\noindent addition       10153          10156\par
\CaseLine{\ttfamily 0 0 0 4 4\qquad 2 2 2 0 0\qquad 2 2 2 0 0}
\CaseLine{\ttfamily 0 1 1 4 4\qquad 2 2 2 1 1\qquad 2 2 2 1 1}
\CaseLine{\ttfamily 0 1 2 3 4\qquad 2 2 2 2 2\qquad 2 2 2 2 2}
\CaseLine{\ttfamily 4 4 3 3 4\qquad 1 1 2 3 4\qquad 1 1 2 4 4}
\CaseLine{\ttfamily 4 4 4 4 4\qquad 1 1 2 4 4\qquad 1 1 2 4 4}
\smallskip}

Thus \(2<1<0<4\), \(2<3<4\), and 3 is incomparable with 0,\allowbreak{}1. Let K=\{3,\allowbreak{}4\}. Direct induction gives these formulas for every word:

\textbullet{} A factor 2, or at least two positions valued in \{0,\allowbreak{}1\}, gives 2.

\textbullet{} With exactly one \{0,\allowbreak{}1\}-valued position, the value is 0 precisely when that position is the literal first position and has value 0; otherwise the value is 1.

\textbullet{} When all factors are in K, 10153 gives their join (3 if all are 3, otherwise 4). In 10156 every nonlinear such word is 4; a literal retains its assigned value.

Indeed a first 0 or 1 stays 0 or 1 under right multiplication by K, a preceding K changes a unique 0 to 1, a second \{0,\allowbreak{}1\} factor gives 2, and 2 is absorbing. The remaining K-products are shown in the tables.

\par\smallskip\noindent\textit{2. Exact criteria, with original right words only.}\quad 
For a target w put \(C=c(w)\), \(h=h(w)\), and let U be the original summands of P with content contained in C. A literal target is absorbed if and only if that literal occurs in P. For nonlinear w:

10153 requires (C53) the union of c(u), u in U, is all of C; (S53) every simple x of w has u in U with exactly one x; (H53) if h is simple, either a single-h witness starts in h, or there is also a word in U omitting h.

10156 requires (C56) at least one word of U is nonlinear; (S56) each simple x has either a nonlinear U-word omitting x, or a U-word containing exactly one x; (H56) at a simple head h, either a nonlinear U-word omits h, or there is a single-h U-word and in addition a same-head single-h U-word or a U-word omitting h. An omitting literal is allowed in this last alternative, but cannot replace the nonlinear omission in S56.

Necessity. For C53, assign a target variable missing from the union 4, other target variables 3 and outsiders 2. The target is 4 and the right side is at most 3. For C56, assign C the value 3 and outsiders 2. A nonlinear target is 4 and the right side can be 4 only from a contained nonlinear word. If a simple-x condition fails, put \(x=1\), C minus \(x=3\), and outsiders 2. The target is 1; right values are 2 or 3, since the required single occurrence and, in 10156, nonlinear omission are absent. If the additional simple-head condition fails, put \(h=0\), C minus \(h=3\), and outsiders 2. All contained words then include h, and every single-h word starts elsewhere, so right values are only 1 or 2, below target 0. A missing literal x is separated by \(x=0\) and all other variables 2: all nonlinear powers of x and every other word are 2.

Sufficiency. Target 2 is bottom. In 10153 a K-valued target 3 is absorbed by any contained word, and a target 4 has a 4-valued variable which occurs in some contained word by C53. In 10156 a K-valued nonlinear target 4 is supplied by C56. For target 0 or 1 the sole \{0,\allowbreak{}1\}-valued occurrence belongs to a simple x. A nonlinear omitting witness in 10156 has value 4. Otherwise a single-x witness has value 0 or 1, absorbing target 1. For target 0, x is the simple head: a same-head witness gives 0, or the nonhead single-x witness gives 1 and an omitting word gives 3 or 4, whose join with 1 is 4. These cases cover every assignment. Thus both criteria are exact for unrestricted variables, word lengths and numbers of summands.

\par\smallskip\noindent\textit{3. Structural laws and canonical forms.}\quad 
Both algebras satisfy the following six identities:

{\small\par\smallskip
\CaseLine{\((T) xyz = xzy\)}
\CaseLine{\((B) xyx = x^{2}y\)}
\CaseLine{\((R) x^{3} = x^{2}\)}
\CaseLine{\((D) xy = x^{2}y+xy^{2}\)}
\CaseLine{\((J) x^{2}y^{2} \leq  x^{2}+y^{2}\)}
\CaseLine{\((E) x^{2}y^{2} = y^{2}x^{2}\)}
\smallskip}

For 10153 add (L) \(x^{2}\leq x\). For 10156 add

{\small\par\smallskip
\CaseLine{\((L1) x^{2}y \leq  xy\qquad (L2) xy^{2} \leq  xy\qquad (A) x^{2}y^{2} \leq  x^{2}.\)}
\smallskip}

The word formula proves these laws. T, B, R keep each variable's multiplicity one/two and whether a simple variable is the head. In E all displayed variables are repeated so the head is immaterial. For D, if a product is 0 or 1, its unique \{0,\allowbreak{}1\} factor survives exactly in the branch doubling the other word; the other branch is 2. A 2-valued product remains 2 in both branches. K-valued products keep the same K-value in 10153 and are all 4 in 10156. J follows by the same cases; with two distinct \{0,\allowbreak{}1\} factors its left side is bottom, with only K factors its right side has the required join or 4. In 10153 squaring a variable can only change a non-K value to 2, so L is valid. In 10156 L1, L2 apply only within a nonlinear word and again can only lower to 2; A is automatic since \(x^{2}\) is 2 or 4. L fails in 10156 at \(x=3\) and is never used for that algebra.

T permutes tail positions while retaining the head. Group repeated variables and apply R to reduce every multiplicity to one or two. A simple head is kept fixed. If the head is repeated, move another repeated block directly after its initial square and use E to exchange the two squares, then sort the tail again by T. Consequently the canonical form has the original simple head if there is one, or the least repeated variable as head otherwise, followed by sorted blocks. Every variable is present once or twice according to its original multiplicity. Thus a word is determined modulo these laws by its content, simple-variable set and either its simple head or the fact that its head is repeated. All substitutions in these steps are nonempty.

\par\smallskip\noindent\textit{4. A saturated full-content guard.}\quad 
In 10153, for every original U-word use L in context to double each simple variable, then normalize. Combine these saturated words using J: a saturated word g equals \(g^{2}\) by the canonical form, so J with \(x=g\), \(y=v\) derives the union guard from g and v. C53 ensures the final guard g has full content C and every variable repeated.

In 10156 choose a contained nonlinear seed from C56. Double its simple head using L1 with the remaining suffix nonempty; double any other simple occurrence using L2 with its preceding prefix nonempty and the remaining suffix as external context. This gives a saturated seed. For each missing variable y, A with \(x=g\) appends its square. Since normalization gives \(g^{2}=g\), the resulting full-content saturated guard is below P. No literal seed is ever squared by L1 or L2.

Fix a simple target variable x. For any selected witness, double every other simple variable using the permitted contractions. In 10156 a literal is left unchanged. A resulting witness therefore has at most one simple variable, except that an omitting literal retains its own single letter. In 10153 a literal different from x can be squared.

For 10153 choose a single-x original witness, preferring a same-head one when x is the target head. If that head witness starts elsewhere, choose an additional original omitting witness as guaranteed by H53. For 10156 first use an omitting nonlinear witness when available. Otherwise choose a single-x witness, and if x is the target head and no such witness starts in x, choose an additional omitting witness. These choices obey S56, H56. At most two original witnesses are needed. Let \(V_{x}\) be their contracted versions and let \(w_{x}\) be w with all variables except x doubled and normalized. Then

{\small\par\smallskip
\CaseLine{\(w_{x} \leq  g + \sum (V_{x})\qquad (F)\)}
\smallskip}

is valid by the exact criterion: g supplies C53/C56 and the selected witnesses supply the sole simple-variable and simple-head conditions. It has at most three right words. Section 5 derives it from one fixed finite family, rather than assuming an unbounded collection of rules.

\par\smallskip\noindent\textit{5. Reversible compression and completeness.}\quad 
There are at most four words in F, including its left word. Reserve x and the head of each word, at most five variables. Every other variable is either absent or repeated in every word: a simple variable different from x can occur only in an omitting literal, whose head is reserved. Group nonreserved variables by their nonzero content-incidence pattern across these four words. There are at most \(2^{4}-1=15\) such classes. Identify each class to a new variable and normalize. This is a genuine substitution into a valid absorption and yields at most twenty variables, normal-word lengths at most forty, and at most three right words.

For the inverse, substitute for each class variable the nonempty product of its members. The class occurs twice or not at all in each compressed word. Hence its expansion restores two occurrences of each member; T, B, R regroup and reduce them. Reserved variables retain multiplicities, and the distinction between a fixed simple head and a repeated head is unchanged; E restores the canonical repeated head if needed. Therefore the inverse gives every original word of F modulo structural laws. There is no replacement of common multiplicities by ordinary incidence.

Fix twenty variables and let \(Sigma_{53}\) or \(Sigma_{56}\) consist of the respective seven or nine structural laws plus every valid absorption having one canonical left word, one to three canonical right words, and each word of length at most forty over these variables. Section 2 is an effective test, so this is a specified finite family. Compression and its inverse derive every local recovery F from this family.

Finally apply D to any normalized target with at least two simple variables, splitting it into nonempty X, Y between two simple occurrences. The branches \(X^{2}Y\) and \(XY^{2}\) retain the content and have strictly fewer simple variables. Every surviving simple variable retains whether it was the simple head; tail permutation and repeated-head exchange do not affect that fact. Iterate. Zero-simple leaves normalize to g; a leaf whose only simple variable is x normalizes to \(w_{x}\). Each leaf is therefore below P by the established guard or local recovery. The tree of D equalities derives the original target absorption. Targets with zero or one simple variable are already leaves. Literal targets are original summands and need no construction.

Distributivity expands arbitrary terms to nonempty polynomials. Derive the absorption of each summand into the opposite polynomial in both directions to prove every valid identity. This proves unrestricted finite basability of both algebras. Reversal gives bases for their multiplicative opposites.

\subsection{Two twelve-identity bases from contained crossing witnesses}\label{case:195}
Signature(+, multiplication), with no constants. Every substitution is nonempty. Inequality \(p\leq q\) means \(p+q=q\).

\par\smallskip\noindent\textit{Tables, word values and normalization.}\quad 
Both 7672 and 7673 have multiplication rows 01234,\allowbreak{}22222,\allowbreak{}22222,\allowbreak{}41244,\allowbreak{}41244. The addition rows for 7672 are 00044,\allowbreak{}01134,\allowbreak{}01234,\allowbreak{}43334,\allowbreak{}44444. For 7673 they are 00044,\allowbreak{}01144,\allowbreak{}01234,\allowbreak{}44334,\allowbreak{}44444. Thus 2 is the additive bottom,4 the top. Both have \(1<0<4\) and \(2<3<4\);7672 also has \(1<3\).

A word has value 2 if a factor is 2 or a nonfinal factor is 1. Otherwise a final 1 gives 1. With all factors in\{0,\allowbreak{}3,\allowbreak{}4\}, its value is 0 when all are 0;3 when its final factor is 3 and every preceding factor is 0;4 otherwise. This is an induction on the displayed table.

The word laws \(xxy=xy\), \(xyz=yxz\), \(xyy=yxx\) give the same elementary prefix-set normal forms as in the 6211 proof, independently of addition. A simple final variable stays fixed; all prefix variables may be sorted and their repeats removed. When the final variable repeats, the word is saturated and depends only on its whole content. Its chosen normal form is sorted content with the largest variable doubled. The literal normalization uses adjacent prefix swaps, \(xxZ=xZ\), and \(qtt=tqq\) to change the saturated final variable. The saturation bit is retained.

\par\smallskip\noindent\textit{Exact absorption conditions.}\quad 
Given target w with content C, retain only original upper words U whose entire content lies in C. Put \(A=\cup\) of their prefix contents (all positions except the final one), \(B=set\) of their final letters. Call U crossing if A intersects B. Then \(w\leq \sum\) P holds exactly when:

1. The union of contents of U is C.

2. Either U is crossing, or the target prefix content is contained in A.

3. If the target final t is simple, either an original word in U ends simply in t, or the following alternative holds: for 7672, some U word omits t; for 7673, the subfamily V of U words omitting t is crossing.

No original word with outside content may serve as these witnesses. The different alternative in 3 is essential: a lone omitting word of value 3 bounds 1 in 7672 but does not bound 1 in 7673.

Necessity: if U is empty, assign C letters 0 and outside letters 2. If a target x is missing from its content union, assign \(x=4\), other C letters 0 and outside letters 2. These separate condition 1. If 2 fails, U is noncrossing: assign A letters 0, B letters 3, outside 2. The target has an internal 3 and hence value 4, while each U word is 3 and other upper words are 2. If 3 fails for 7672, assign \(t=1\), other C letters 0, outside 2. Every U word contains t internally and is 2, while the target is 1. For 7673 let V be the omitting family. It is noncrossing; assign its prefix union0, its final letters 3, \(t=1\), remaining C letters 0, outside 2. Every V word is 3 (or V is empty), and every other U word has internal t and is 2. Neither their join 3 nor 2 bounds the target 1. All assignments are compatible because the relevant prefix and final sets are disjoint.

Sufficiency uses only U. Target 2 is automatic. A target 0 has all target factors 0 and every U word has value 0. A target 3 has its unique final t assigned 3, others 0; coverage supplies a U word containing t, whose value is 3 or 4. For target 4, a target 4 letter is covered and immediately forces an upper 4. Otherwise some target internal letter is 3. If it belongs to A, a U word is 4. If it does not,2 guarantees a crossing. If no upper 4 already occurs, all variables in A are 0; an original word ending in A intersection B has value 0, while coverage gives an upper 3. Their join is 4. For target 1, its simple final t is 1 and its prefix lies in\{0,\allowbreak{}3,\allowbreak{}4\}. A simple-t word supplies 1. Otherwise, for 7672 any omitting word has value 0,\allowbreak{}3 or 4, each bounding 1. For 7673 a crossing omitting family either already has value 4 or, with its internal letters 0, has an original final-in-prefix word of value 0. Either bounds 1. This proves exact semantics for arbitrary variables and lengths.

\par\smallskip\noindent\textit{Twelve identities for each algebra.}\quad 
Both bases contain seven fixed laws:

{\small\par\smallskip
\CaseLine{\(xxy=xy\qquad xyz=yxz\qquad xyy=yxx\)}
\CaseLine{\(xyz\leq xz+yz\qquad xyxy\leq xx+y\)}
\CaseLine{\(xyt\leq xytxyt+xt\qquad xy\leq xyxy+y\)}
\smallskip}

They also contain the four optional-context identities

{\small\par\smallskip
\CaseLine{\(([a]x[b]xc)^{2} \leq  [a]x+[b]xc .\)}
\smallskip}

Both a and b may be omitted when writing a separate law; x, c are always present. Finally add \(xy\leq xyxy+x\) for 7672 and \(xy\leq xyxy+xx\) for 7673. Each resulting list has 12 identities, at most four variables per identity. Soundness is checked on all 1400 assignments per list. Omission never means empty substitution.

\par\smallskip\noindent\textit{Unrestricted derivation.}\quad 
Normalize and choose the contained original family U. Suppose U is crossing. Choose x in A intersection B, an original \(u=A0\) x, and an original \(v=B0\) x C0 with C0 nonempty. The corresponding crossing template derives the saturated word of content c(u) union c(v). For each other original q in U, use \(xyxy\leq xx+y\) with \(x=the\) accumulating saturated word, \(y=q\). Its right square normalizes to the same accumulating word; its left saturates the union of contents. By finitely many steps obtain the saturated target content. This construction may also be performed on any crossing subfamily.

If the target is saturated this is already its normal form. If it has simple final t and an original contained simple-t word Xt, use \(xyt\leq xytxyt+xt\) with \(x=X\), \(y=the\) whole target prefix. The two large terms normalize to the desired simple target and saturated target. If the original word is t alone use \(xy\leq xyxy+y\). A one-letter simple target with such a witness is already original. If no such witness exists, use condition 3. For 7672 choose an original omitting word h and apply \(xy\leq xyxy+x\) with \(x=h\), \(y=w\). For 7673 construct a saturated omitting word h from the crossing omitting subfamily as above and use \(xy\leq xyxy+xx\) with \(x=h\), \(y=w\). Because h omits t and is contained in the target content, the left normalizes to w, the large right term to the saturated target, and the small right term to the known h. No empty substitution occurs.

If U is not crossing, condition 2 and coverage force every original word in U to end in the same simple target t: the original final set B is disjoint from A and lies in C, while C minus t lies in A. The target must be simple. If it has one letter it is original; otherwise merge the nonempty original t-prefixes via \(xyz\leq xz+yz\). Their union is exactly the target prefix, so normalization yields w.

Thus every valid word absorption is deducible without a bound on the number of variables, words or lengths. General polynomial identities reduce to their word absorptions in both directions. These are complete finite bases for 7672 and 7673, respectively. Catalogue transfers require separate complete independent verification.

\subsection{Six identities for 7674 with a marked-prefix and unmarked-final test}\label{case:196}
Signature (+, multiplication), without constants, relative to the additively idempotent semiring axioms.

\par\smallskip\noindent\textit{1. Tables and the word-value distinction.}\quad 
In element order 0,\allowbreak{}1,\allowbreak{}2,\allowbreak{}3,\allowbreak{}4 the tables are

{\small\par\smallskip
\noindent addition        multiplication\par
\CaseLine{\ttfamily 0 0 0 4 4\qquad 0 1 2 3 4}
\CaseLine{\ttfamily 0 1 1 3 4\qquad 2 2 2 2 2}
\CaseLine{\ttfamily 0 1 2 3 4\qquad 2 2 2 2 2}
\CaseLine{\ttfamily 4 3 3 3 4\qquad 4 4 4 4 4}
\CaseLine{\ttfamily 4 4 4 4 4\qquad 4 4 4 4 4}
\smallskip}

The additive order is \(2<1<0<4\) and \(1<3<4\), with 0 and 3 incomparable. Scan a word for its first nonzero factor. If there is none, its value is 0. A first 2 gives 2; a first 4 gives 4. A first 1 gives 1 only at the final position and otherwise gives

2. A first 3 gives 3 only at the final position and otherwise gives

4. This formula includes linear words and follows by induction.

In particular the simple final position must retain both its 1 and 3 behavior. The inequality \(xy\leq xy^{2}\) is false at \(x=0\), \(y=1\), and a construction which simply saturates that final position is invalid.

\par\smallskip\noindent\textit{2. The six laws and word forms.}\quad 
The following six identities form a complete relative basis:

{\small\par\smallskip
\CaseLine{\((W1) xxy = xy,\)}
\CaseLine{\((W2) xyx = xyy,\)}
\CaseLine{\((A)  xy \leq  y+xz,\)}
\CaseLine{\((B)  xyz \leq  xz+xyt,\)}
\CaseLine{\((C)  xyz \leq  xz+xy+xx,\)}
\CaseLine{\((D)  xy \leq  xx+xyt.\)}
\smallskip}

Let F(w) list first occurrences. A word with a simple final variable normalizes to F(w). A word whose final variable occurs earlier normalizes to F(w) followed once more by its last listed variable. Linear x remains x; a longer singleton word becomes xx.

For a separated repeated occurrence before a later suffix, use \(xYxZ=xY^{2}Z=xYZ\) with nonempty Y, Z. Adjacent repeats before a later suffix use (W1). After internal repetitions are removed, a repeated final letter in p x q x, with q nonempty, becomes p x \(q^{2}\) by (W2); deleting internal repeats gives the stated doubled final letter. All contexts may be omitted where appropriate, but substitutions for identity variables are nonempty. The word-value formula proves soundness of (W1),(W2).

\par\smallskip\noindent\textit{3. Exact absorption criterion.}\quad 
For original right words P and a set S, let Z(S) mean that an original word has content contained in S. Let E(x, S) mean that some original word contains x with first-x prefix contained in S. Such a witness is marked when a position follows its first x, and unmarked when that first x is its final position.

A linear target is absorbed exactly when it is a literal summand. For a nonlinear target w, absorption is equivalent to Z(c(w)) and the following requirement at each first occurrence x, with S its target prefix content:

\textbullet{} At a marked target occurrence: a marked E-witness, or an arbitrary E-witness together with Z(S).

\textbullet{} At an unmarked target occurrence: an unmarked E-witness, or an arbitrary E-witness together with Z(S).

Only the simple final variable can have an unmarked first occurrence. The two clauses exchange marked and unmarked witnesses; neither clause can be discarded.

Necessity of Z(c(w)) follows by assigning target variables 0 and outsiders 2. For a marked x, assign \(S=0\), \(x=3\) and all others 2. The target is 4; right values reach 4 precisely with a marked E-witness, or an unmarked E-witness valued 3 and a Z(S) witness valued 0. For an unmarked x, the same assignment requires E. Also assign \(S=0\), \(x=1\) and all others 2: the target is 1, and right values reach at least 1 precisely with an unmarked E-witness or Z(S). Combining the two requirements gives the unmarked clause. Without a literal linear y, assign \(y=1\) and all others 2; every nonlinear right word is 2, as are the other literals.

For sufficiency, value 2 is least and an all-zero target is handled by its contained witness. A first 4 is handled by E. A marked first 1 gives 2, already absorbed; a marked first 3 is handled by a marked witness giving 4, or by values 3 and 0 from E and Z. For an unmarked first 1, an unmarked witness gives 1, while Z(S) gives 0. For an unmarked first 3, any E-witness gives 3 or 4. These are all cases and prove the criterion without a variable bound. It directly verifies the soundness of (A)--(D).

\par\smallskip\noindent\textit{4. Constructing prefixes of every original word.}\quad 
Normalize a valid nonlinear target, with first-occurrence list \(x_{1}\)...\(x_{d}\). Write \(p_{i}=x_{1}\cdots x_{i}\). At each marked target position we construct \(p_{i}\) \(u\leq P\) for every original u, including literals.

The first position is marked and its earlier prefix is empty, so Z(empty) is impossible. There must be an original marked witness \(x_{1}\) t with t nonempty. Identity (A) gives \(x_{1}\) \(u\leq u+x_{1}\) \(t\leq P\) for every original u.

Suppose p \(u\leq P\) holds for every original u, and the next marked variable is x. If a marked witness v exists, prefix deletion normalizes p v to p x t with t nonempty, so \(pxt\leq P\). Then (B) gives \(pxu\leq pu+pxt\leq P\) for every original u.

Otherwise take an unmarked witness v and an original r contained in c(p), as supplied by the criterion. Prefix deletion gives \(pv=px\) and \(pr=p^{2}\). Thus px and \(p^{2}\) are both below P. Identity (C) gives \(pxu\leq pu+px+p^{2}\leq P\). This continues the construction through all marked target positions.

If the target is saturated, every first occurrence is marked. Choose an original \(u_{0}\) contained in the target content. The word \(p_{d}\) \(u_{0}\) normalizes to the saturated target, proving its absorption.

If the target is unsaturated, it has at least two distinct variables and is p s, where s is its simple final variable. The construction has reached p \(u\leq P\) for every original u. If an unmarked witness v for s exists, \(pv=ps\) is the target. Otherwise the criterion supplies a marked witness v and an original r contained in c(p). They give \(pv=pst\) for a nonempty t and \(pr=p^{2}\), both below P. Identity (D) now yields \(ps\leq p^{2}+pst\leq P\). This recovers the unsaturated target without using the false inequality \(ps\leq ps^{2}\).

Linear targets are literal summands. Thus every valid word absorption is derivable. Expanding arbitrary terms and deriving the absorption of each summand into the opposite side proves every identity.

Theorem. The six displayed identities are a finite basis for

7674. The proof is unrestricted in original variables, word lengths and summands, including arbitrary nonempty term substitutions. Reversal gives a basis for the opposite. Complete two-way HSP certificates are required for further catalogue transfers. Numerical tests support this proof rather than replace its unrestricted steps.

\subsection{A finite basis for 8505 from full original contained-content coverage}\label{case:197}
Signature (+, multiplication), no constants, relative to the additively idempotent semiring axioms. Finite family specified, not an enumerated small basis. Completeness has no variable, length or summand bound.

Addition rows:00004,\allowbreak{}01114,\allowbreak{}01224,\allowbreak{}01234,\allowbreak{}44444; order \(3<2<1<0<4\). Multiplication rows:02224,\allowbreak{}12221,\allowbreak{}22222,\allowbreak{}33333,\allowbreak{}42224. Literals retain value. For nonlinear hT, heads 2/\allowbreak{}3 retain value. Head 1 gives 1 if all tails belong to\{0,\allowbreak{}4\}, otherwise2. Head 4 gives 4 if all tails belong to\{0,\allowbreak{}4\}, otherwise2. Head 0 gives 2 if a tail belongs to\{1,\allowbreak{}2,\allowbreak{}3\}; otherwise it gives 4 if a tail is 4, and 0 if all tails are 0. Induction proves this formula. Thus \(xyz=xzy\) and \(xyy=xy\) hold. Normalize to head plus sorted tail set; retain head occurrences in the tail. Actual swaps and deletions have nonempty substituted prefixes; empty word substitutions never occur.

\par\smallskip\noindent\textit{Exact criterion using ORIGINAL words.}\quad 
A literal x is absorbed iff original literal x occurs. For nonlinear hT, \(C=\{h\}\) union T and U the original words contained in C, the criterion is: (H) an original right word has head h; (C) the contents of U cover ALL of C; (R) h belongs to T, or literal h occurs, or an original word is entirely contained in T, or an original nonlinear h-word has tail contained in T.

Necessity for missing literal x:assign \(x=1\), others 3. All right nonlinear x-words give 2 and other heads 3; without literal x the right sum is below 1. Failure of H:assign \(h=2\), others 3; target 2 and right 3 separate. Failure of C at an uncovered x:assign \(x=4\), other target variables 0, outsiders 3. The target is 4 whether x is the head or occurs in the tail. Any contained right word avoids x and gives 0; noncontained words give 2 or 3. Right \(\sum \leq 0<4\). For failed R, assuming H, C, h is outside T. Assign \(h=1\), \(T=4\), others 3. The target is 1. Literal h is absent; nonlinear h-words have a low tail outside T and give 2. Words headed in T are not contained in T, so have a low factor and give 2; other heads 3 give 3. Right \(\sum \leq 2<1\). These separate every failure.

Sufficiency:3 is least. Target 2 is absorbed by an original same-head word from H, which always has value at least 2. For target 4 all target variables belong to\{0,\allowbreak{}4\} and some x is 4. C supplies a contained word containing x, also valued 4. For target 0 the head and all tail variables are 0; any contained word from C gives 0. For target 1 the head is 1 and all tails are 0/\allowbreak{}4, so h is outside T. R supplies literal 1, a tail-contained word valued 0/\allowbreak{}4, or a nonlinear same-head tail-contained word valued 1. Thus the criterion is exact for arbitrary assignments and polynomial sizes.

\par\smallskip\noindent\textit{Build full content, restore the head, and recover the simple head.}\quad 
Normalize all words. Literal targets need no derivation. Choose any original contained word g from U. If some target variable is missing, choose an original contained u containing that variable. Replace g by normal(gu). The local absorption \(gu\leq g+u\) is valid:its head comes from g, its content is covered by g and u, and R is supplied by literal g if g is literal, or by the original nonlinear g whose tail remains contained in gu. Thus each step adds content and stays inside C. Repetition terminates with an available g of content exactly C. Derive normal(g squared)\ensuremath{\leq}g:the target has repeated head, so R holds, and H, C are supplied by g. Its tail is now C.

Choose the original H witness a, which may have variables outside C. The local absorption \(hC\leq a+g\) satisfies H through a, C through g, and R because its head is repeated. This restores the target head while retaining full content. If the desired target repeats h, it is already hC. Otherwise take one original R witness b:literal h, tail-contained word, or nonlinear h-word with tail in T. The local absorption \(hT\leq hC+b\) satisfies H, C through hC, and R through b. It recovers the target. Every local step has at most 2 already available right words. No claim of finite basis has yet relied on treating all valid absorptions as primitive; the next paragraph bounds them.

\par\smallskip\noindent\textit{Uniform reversible compression.}\quad 
Each local absorption has at most 3 normal words, one left and at most 2 right. Reserve their heads, at most 3 variables. Partition remaining variables by their nonempty incidence across the at most 3 tails, at most 7 classes. Identify each class to a new variable, leave reserved heads separate, and normalize. Validity is preserved by substitution. At most 10 variables and word length at most 11 remain. Invert each new variable by the nonempty product of its class members; identical tail incidence recovers all original tail sets, while reserved heads and head-tail repetitions are preserved. Tail sorting and duplicate deletion recover every original word.

Let Sigma contain \(xyz=xzy\), \(xyy=xy\) and all valid normal absorptions on a fixed 10-variable alphabet, word length at most 11, at most 2 right words. The exact H, C, R criterion specifies this finite family. The full-content construction above is a sequence of expanded substitution instances of Sigma, so every valid absorption follows. Distributivity and summandwise absorption in both directions yield every valid identity. Sigma is a finite relative basis for 8505. Opposite reversal is valid; catalogue transfer still requires independent complete HSP both ways. Finite computations supplement the unrestricted proof. In particular one contained witness cannot replace full coverage:\(xy\leq x\) fails at \(x=0\), \(y=4\).

\subsection{Full-content guards and an affine finite basis for entry 8850}\label{case:198}
Work without constants in the signature (+, multiplication), relative to the axioms of additively idempotent semirings. Let A be catalogue entry 8850. Write w \ensuremath{\leq} P for w+P approximately equal to P. Numerical labels name elements, not term constants.

\par\smallskip\noindent\textit{Operations and words.}\quad 
Addition consists of the chain \(1<0<4\) and the two incomparable elements 2,\allowbreak{}3 below 4, with every join across the two branches equal to 4. The tables in numerical element order are

{\small\par\smallskip
\noindent addition       multiplication\par
\CaseLine{\ttfamily 0 0 4 4 4\qquad 1 1 0 0 0}
\CaseLine{\ttfamily 0 1 4 4 4\qquad 1 1 1 1 1}
\CaseLine{\ttfamily 4 4 2 4 4\qquad 1 1 2 3 4}
\CaseLine{\ttfamily 4 4 4 3 4\qquad 1 1 3 2 4}
\CaseLine{\ttfamily 4 4 4 4 4\qquad 1 1 4 4 4}
\smallskip}

Put G=\{2,\allowbreak{}3\} and K=\{2,\allowbreak{}3,\allowbreak{}4\}. On G, multiplication is the group of order two, with identity 2. A nonlinear word has value 1 if it has a factor 1 or an occurrence assigned 0 after the initial position. It has value 0 if its initial occurrence is assigned 0 and every later occurrence lies in K. On all-K assignments, a factor 4 gives 4; otherwise the value is 2 or 3 according as the number of occurrences assigned 3 is even or odd. Linear words keep their assigned values.

These rules agree for both associations of each triple: a factor 1 forces 1; a unique initial 0 with later factors in K gives 0 and any other occurrence of 0 gives 1; on K associativity follows from the group with an absorbing element 4. Each row and column preserves the joins inside the chain, the join \(2+3=4\), joins across the two branches, and joins involving 4, as the tables show. Thus A is an additively idempotent semiring, and induction proves the formula at every word length.

For a word u let c(u) be its content, h(u) its initial variable, and t(u) the set occurring after its first position. Linear words have empty tail set. Let v(u) in \(F_{2}^{X}\) record total variable multiplicities modulo two. The identities

\CaseDisplay{xyz\approx xzy,\qquad xy^3\approx xy       \quad\textup{(\ensuremath{W})}}

hold by the word formula. They sort the tail and reduce its positive multiplicities to one for odd multiplicity and two for even multiplicity. A normal word retains its literal head and these tail multiplicities in \{0,\allowbreak{}1,\allowbreak{}2\}. In particular its head has total multiplicity one, two or three; a simple head and a repeated odd head must remain distinct. Linear words stay linear.

\par\smallskip\noindent\textit{Exact absorption criterion.}\quad 
Let P be a nonempty sum of words and \(D=c(P)\). For C contained in X, put B\_C=\{u in P:c(u) subset C\}, including linear words.

Lemma 1. A linear y is absorbed by P exactly when it is a literal summand. A nonlinear w of content C and head h is absorbed exactly when

1. C is contained in D;

2. v(w) belongs to the affine span over \(F_{2}\) of \{v(u):u in P\};

3. \(B_{C}\) is nonempty;

4. if h is simple in w, some b in \(B_{C}\) omits h or has h as simple head.

Affine span here means the set of sums of odd numbers of the displayed vectors, with coefficients in \(F_{2}\).

Proof. A missing literal y is separated by \(y=0\) and all other variables 1: every nonlinear power of y and every other right word is 1. If 1 fails, assign a missing target variable 4 and all other variables 2. Then \(w=4\) and \(P=2\).

Condition 2 is equivalently (v(w),1) lying in the linear span of the lifted right vectors (v(u),1). If it fails, elementary linear algebra supplies a linear functional (alpha, b) that annihilates every lifted right vector and does not annihilate the lifted target. Assign \(x=3\) when \(alpha_{x}=1\) and \(x=2\) otherwise. All right words have the common parity b, whereas the target has the other parity. Its value and P are the incomparable elements 2 and 3, so absorption fails. This proves necessity of 2 for arbitrary numbers of variables.

If 3 fails, assign C the value 2 and all outsiders 1. Then \(w=2\) and \(P=1\). If 4 fails, assign \(h=0\), other C-variables 2 and outsiders 1. Then \(w=0\); every contained right word has an internal occurrence of h and every other right word has a factor 1, giving \(P=1\).

For sufficiency, if P is 2 or 3, all right words have that exact value and every variable in D belongs to G. By 1, so does every target variable. An odd affine combination of right parity vectors gives that same parity, so \(w=P\) by 2. If P lies in \{1,\allowbreak{}0\}, w cannot have all its variables in K: a witness from 3 would take a value in K. Thus w lies in \{1,\allowbreak{}0\}. When \(P=1\), the only remaining violation would be \(w=0\), with simple head h assigned 0 and every other C-variable in K. A witness from 4 omitting h takes a value in K, and one with simple head h takes value 0. Both contradict \(P=1\). The case \(P=4\) is automatic. These are all proper principal additive ideals, proving sufficiency. \(\square\)

\par\smallskip\noindent\textit{Full-content affine guards.}\quad 
Fix valid nonlinear w \ensuremath{\leq} P and let \(D=c(P)\). For each x in D choose an original \(a_{x}\) containing x. Every local step below has at most three right words, each original or already derived below P. Validity follows directly from Lemma 1.

For an original summand u, append h(u)\ensuremath{^{2}} to obtain q. The absorption q \ensuremath{\leq} u is valid: content and total parity are unchanged, u is a contained witness, and the target has repeated head. This includes linear u, where the target is its cube. For each missing x in D, use \(qx^{2}\) \ensuremath{\leq} q+\(a_{x}\). The target keeps q's parity and repeated head; q supplies a contained witness and \(a_{x}\) supplies the added content. Thus we obtain, for every needed original u, a known guard \(Q_{u}\) with full content D, repeated head and parity v(u).

For any vector v on D there is a normal full-content word R(D, v) of repeated head and parity v: fix g in D, give g multiplicity two if \(v_{g}=0\) and three if \(v_{g}=1\), and give every other x multiplicity two if \(v_{x}=0\) and one if \(v_{x}=1\). Its head is g.

Condition 2 supplies an odd subset \(u_1,\ldots,u_{2k+1}\) of original summands whose parity sum is v(w). Start with \(Q_{u_1}\), and successively combine it with pairs of the other guards. Each step

\CaseDisplay{R(D, v(q)+v(Q_i)+v(Q_j))\leq q+Q_i+Q_j                  \quad\textup{(\ensuremath{1})}}

is valid: all words have content D; the target parity is their odd affine sum; there is a contained right witness; and its head is repeated. This constructs a known full-content guard Q with parity v(w), with no bound on the number of affine summands used.

If the head of w is repeated, choose any b in \(B_{C}\). If it is simple, choose b as in 4. The final local absorption w \ensuremath{\leq} Q+b is valid: Q supplies content and the required parity vector, b supplies a contained witness and, when needed, the simple-head witness. The target content may be a proper subset of D; no unsupported deletion of extra variables is made. This proves derivability of every valid word absorption by local steps with at most three right terms.

\par\smallskip\noindent\textit{Compression and completeness.}\quad 
A local step contains at most four words. Reserve their literal heads individually, at most four variables. Group every remaining variable by its vector of tail multiplicities in \{0,\allowbreak{}1,\allowbreak{}2\} across these words. There are at most \(3^{4}-1=80\) nonzero patterns. In each compressed word give the group symbol that common tail multiplicity, retaining its reserved head. This uses at most 84 variables and normal word length at most 169.

The operation is not ordinary variable identification. Lemma 1 proves its validity: content coverage and contained-witness support are preserved by the common presence patterns; the simple-head condition is preserved because all heads are individual. For affine membership, the parity equations for all variables in one group are identical across the target and right words. Keeping one such coordinate removes only duplicate equations. The constant lifted coordinate and every reserved head coordinate remain unchanged. Thus the same coefficients solve the lifted affine equations before and after compression, proving equivalence of condition 2.

Conversely, substitute the nonempty ordered product of each group's members for its symbol, leaving every reserved head individual. Common tail multiplicities restore each member once or twice in each word. Tail sorting and the reductions in (W) recover the original normal words. Linear words and all three possible head multiplicities are included. No units or empty substitutions occur.

Theorem. The two identities (W), together with all valid normal absorptions over a fixed 84-letter alphabet, of word length at most 169 and with at most three right summands, form a finite basis for A. Word reversal gives a finite basis for its opposite.

Proof. The specified family is finite and Lemma 1 decides its membership by finite-set tests and Gaussian elimination. Compression and inverse substitution derive every local absorption from the family and (W). The guard construction and odd-three affine steps derive all valid word absorptions, without bounds on variables, lengths or summands. Distributing arbitrary terms and deriving each summand's absorption in both directions proves every identity. Reversal transports the result to the opposite. \(\square\)

This is an unrestricted proof for a finite family, not an enumerated small or minimal basis. It asserts no hereditary property. Finite computations only check the argument. Transfers to other catalogue types require complete independent HSP verification in both directions.

\subsection{Four witnesses and a finite basis for source 8889}\label{case:199}
Work relative to the ai-semiring axioms in the constant-free signature. The element labels below are not constants in terms. Source A has carrier \{0,\allowbreak{}1,\allowbreak{}2,\allowbreak{}3,\allowbreak{}4\} and tables

addition  multiplication 00044     11100 01044     11111 00244     11120 44434     11134 44444     11144

Addition has incomparable elements 1 and 2 with join 0; the element 3 is incomparable with 0,\allowbreak{}1,\allowbreak{}2, and its join with each of them is 4. Put K=\{3,\allowbreak{}4\}. The four proper principal ideals are \{1\}, \{2\}, \{1,\allowbreak{}2,\allowbreak{}0\}, and \{3\}.

\par\smallskip\noindent\textit{Word semantics.}\quad 
For a nonlinear word with head h and nonempty tail-variable set T:

\textbullet{} \(h=1\) always gives value 1;

\textbullet{} \(h=0\) gives 0 if all T-variables lie in K, and 1 otherwise;

\textbullet{} \(h=2\) gives 2 if all T-variables are 3, gives 0 if all lie in K and some is 4, and gives 1 otherwise;

\textbullet{} h in K gives the maximum of all word values if all tail variables lie in K, and gives 1 otherwise.

The formulas follow by induction from the rows. They also give associativity: if a right block has a factor outside K its value is outside K, so it sends a preceding factor to 1; if all its factors are in K, its product is their maximum. The analogous head cases give the same result for either association. Each multiplication map preserves the joins \(1+2=0\) and \(3+x=4\) (x=0,\allowbreak{}1,\allowbreak{}2), as checked from the rows, and preserves comparable joins. Thus both distributive laws hold.

Consequently

{\small\par\smallskip
\CaseLine{\(xyz \approx  xzy,\qquad xy^{2} \approx  xy\qquad (W)\)}
\smallskip}

are valid. They normalize every nonlinear word to its head followed by the sorted set of its tail variables. The head may reoccur in the tail, and this reoccurrence is retained. The square identities x\ensuremath{{}^{2}}y\ensuremath{\approx}xy and x\ensuremath{\leq}\(x^{2}\) are not used.

\par\smallskip\noindent\textit{Exact absorption criterion.}\quad 
For a polynomial P let L be its linear variables, N its nonlinear words, H the set of heads of N, E the union of all tail-variable sets of words of N, and D the content of P. For a nonlinear target w write \(h=head(w)\), \(T=tail-content(w)\), and \(C=\{h\}\) union T.

Lemma 1. A linear y is absorbed by P exactly when y belongs to L. A nonlinear w is absorbed exactly when all three conditions hold:

(C) C is contained in D;

(J) either L meets C, or some v in N has head in C and tail-content contained in T;

(F) either E meets H union L, or h belongs to H union L and T is contained in E.

Necessity. If the linear variable y is missing, assign \(y=0\) and all other variables 1. Every nonlinear right word is 1, as is every linear right word. Thus y is not absorbed.

If (C) fails, assign a missing target variable 4 and all other variables 3. The target is 4 and the right side is 3.

If (J) fails, assign all variables of T the value 3. If h is not in T assign it 0; assign all variables outside C the value 1. The target is 0 or 3. No linear right word has a variable in C, so every linear summand is 1. A nonlinear right word either has head outside C, hence value 1, or has a tail variable outside T by failure of (J). That variable has value 0 or 1, so this word also has value 1. Thus \(P=1\) and it does not absorb w.

Suppose (F) fails. Then E avoids H union L. Assign H union L the value 2 and E the value 3, consistently. Every right summand is 2. If h is not in H union L, assign \(h=3\) if already in E and 1 otherwise, and assign any still unassigned variables 3; the target cannot be 2. If h is in H union L, some target tail variable is not in E. Any such variable already in H union L is 2; otherwise assign one such variable 1. Complete the other unassigned values by 3. The target has head 2 and a tail outside K, so is 1. In both cases it is not absorbed by \(P=2\). This proves necessity.

Sufficiency. To prove w\ensuremath{\leq}P it suffices to check that whenever P lies in one of the four proper principal ideals, w lies there too. For \(P=3\) all variables in D are 3; (C) gives \(w=3\).

If P belongs to \{1,\allowbreak{}2,\allowbreak{}0\} and w were in K, all C-variables would be in K. The witness in (J) would then have value in K, contradicting the right-side bound.

If \(P=1\) and w were not 1, then h is not 1 and all T-variables are in K. Every variable of C is therefore not 1. A linear witness from (J) would contradict \(P=1\). A nonlinear witness with head in C and tail in T also has value different from 1, again a contradiction.

Finally \(P=2\) means that all L-variables and all H-variables are 2, and all E-variables are 3. The intersection alternative in (F) is impossible under such an assignment. Therefore h is 2 and all its tail variables are 3, giving \(w=2\). This treats every proper ideal. Linear sufficiency follows from literal occurrence. The criterion holds for arbitrary variable sets and raw words.

\par\smallskip\noindent\textit{Chains with at most four right summands.}\quad 
Lemma 2. Every valid nonlinear absorption w\ensuremath{\leq}P follows from (W) and a finite chain of valid absorptions, each involving a nonlinear target and at most four right words. Each right word in a step is an original summand of P or the target of a previous step.

Proof. Normalize words and keep the notation of Lemma 1. We construct targets with head h whose nonempty tails increase through subsets of T, finishing with exactly T.

First suppose E avoids H union L. By (F), h belongs to H union L and T is contained in E. In (J) a linear witness can only be h: any other target variable is in T subset E, which is disjoint from L. Likewise any nonlinear witness in (J) must have head h, because a different target head would lie in E intersect H.

If there is a nonlinear witness with head h and tail contained in T, use it as the initial derived word. Otherwise h is linear; choose t in T and a nonlinear right word u with t in its tail. The absorption ht\ensuremath{\leq}h+u satisfies Lemma 1. Subsequently, to introduce any missing x in T choose a nonlinear original word u with x in its tail. If s is the current target, the absorption s x\ensuremath{\leq}s+u satisfies (C), uses s itself for (J), and has target tail contained in the union of the two right tails for (F). Thus it is valid. Every step uses at most two right words.

Now suppose E meets H union L. Choose a variable g in that intersection, a nonlinear right word a with g in its tail, and a right summand b whose head is g (b may be linear). Any right polynomial containing a and b automatically has the intersection alternative in (F).

Choose a witness v for (J). If it is linear, say \(v=k\), and k differs from h, choose an original word r containing h and use the initial target hk. The polynomial a+b+k+r has (C), the linear witness k for (J), and the fixed intersection for (F). If \(k=h\), choose t in T and a right word r containing t; use the initial target ht and right polynomial a+b+h+r. In both cases at most four right words suffice.

If v is nonlinear with head k in C and tail U contained in T, take as initial tail U union \{k\} when k differs from h, and U when \(k=h\). It is nonempty. Choose an original right word r containing h. The corresponding h-headed target is absorbed by a+b+v+r: its content is covered, v supplies (J), and a, b supply (F).

Finally, for each missing x in T choose an original right word r containing x. If s is the current target, the absorption s x\ensuremath{\leq}s+a+b+r is valid by Lemma 1: its content is covered, s supplies (J), and a, b preserve (F). This introduces one missing variable with at most four right words. After finitely many steps the normal form is w. Duplicate right summands can simply be removed. Since every step is absorbed by P by induction, the chain proves the claim.

\par\smallskip\noindent\textit{A finite, explicitly specified basis.}\quad 
Fix a set X of 36 variable symbols. Let Normal(X) consist of all linear variables and all words consisting of one head followed by a nonempty sorted subset of X. These words have length at most 37. Define Gamma to contain every absorption

{\small\par\smallskip
\CaseLine{\(w \leq  u_{1} + \cdots  + u_{k},\qquad 1 \leq  k \leq  4,\)}
\smallskip}

with w, \(u_{1}\),..., \(u_{k}\) in Normal(X) that satisfies Lemma 1. This is a concrete finite set: Normal(X) is finite and the displayed criterion decides which of its finitely many candidate absorptions are included. It is not necessary to enumerate this very large finite set to establish its finiteness or completeness.

Lemma 3. Every valid absorption involving one target word and at most four right words is derivable from Gamma and (W).

Proof. Normalize the at most five words. Retain each variable that is the head of at least one word as a separate variable; there are at most five such variables. For every other variable record the set of word-tail positions in which it occurs. This is a nonempty subset of at most five word positions, so at most 31 different occurrence patterns exist. Identify all non-head variables with the same pattern, while leaving head variables distinct and separate from these groups. Thus the resulting absorption uses at most \(5+31=36\) variables. Since it is a variable-identification instance of a valid identity, it remains valid in A. Its normalized form is a member of Gamma after renaming.

For each pattern variable substitute the product of the original variables in that pattern group, in a fixed order; for each retained head substitute the original head variable. These products are nonempty. Every pattern group occurs exclusively in tails, and every variable in that group occurred in exactly the same tails. Consequently normalization by (W) recovers the original words exactly. Linear words are retained head variables and are unaffected. This proves the desired derivation from Gamma. The reverse substitution, not just a count of patterns, establishes the equational transfer.

Theorem. The two word identities (W) together with Gamma form a finite basis for A. Its opposite also has a finite basis obtained by word reversal.

Proof. Soundness is part of the definition of Gamma and the word formula. For a valid nonlinear absorption, Lemma 2 supplies a chain of steps; Lemma 3 derives every step from the same fixed finite set. A linear absorption is a literal right summand by Lemma 1. Distributivity expands arbitrary terms into finite nonempty sums of words, and each summand absorption can be derived in both directions for any valid identity. Thus all identities follow. Nothing in the argument restricts the original variable set, word length or polynomial size. Word reversal preserves the entire proof.

\par\smallskip\noindent\textit{What is and is not claimed.}\quad 
This proves finite basability with an explicitly defined bounded-variable finite family. It does not give a small enumerated list or claim a minimal basis. The bounds 36 variables, word length 37, and four right summands are proved above, not inferred from experiments. Finite regressions and chain/compression checks provide supplementary error detection only. Transfer to other catalogue types requires complete two-way HSP certificates and independent checking. No hereditary finite-basability assertion is made.

\subsection{A nonfinite identity basis obstruction for 7805}\label{case:200}
The catalogue product rows are 01244, 11111, 41424, 41434, 41444. Addition has bottom 1, top 4 and three incomparable intermediate elements 0,\allowbreak{}2,\allowbreak{}3. The value 1 annihilates multiplication. In its absence 4 annihilates, and the remaining labels represent intervals 0=(0,\allowbreak{}0), 2=(0,\allowbreak{}1), 3=(1,\allowbreak{}1). Two intervals multiply by matching the right endpoint to the next left endpoint, giving the interval between the exterior endpoints, or 4 if matching fails. These nine products follow directly from the rows. A sequence giving 2 consequently has exactly one entry 2 with only 0's to its left and 3's to its right.

For odd n at least five, take variables \(x_{0}\),..., \(x_{n-1}\) in cyclic order. For each edge i, i+1 put the two-letter word in forward order for even i and reverse order for odd i, with subscripts modulo n. Let P be their sum and take \(W=x_{0}\) ... \(x_{n-1}\). Every vertex except 0 is a pure source or a pure sink of the directed cycle edges. Vertex 0 has one edge of each orientation.

\par\smallskip\noindent\textit{Validity for every odd length.}\quad 
If any variable is 1, W is least. In its absence no upper word evaluates to 1. For their sum to be below the atom 0, every edge must have both entries 0, giving \(W=0\). For the bound 3 the same argument gives \(W=3\). If their sum is below the atom 2, every edge has exactly one endpoint assigned 2. The vertices with that value would alternate around the odd cycle, which is impossible. The least bound itself is impossible and the greatest is automatic. This proves W \ensuremath{\leq} P at every indicated n.

\par\smallskip\noindent\textit{A gap in the number of distinct letters.}\quad 
If a raw word v is absorbed by P, its support C must be covered by upper words wholly contained in C. Otherwise assign an uncovered letter 4, other C-letters 0 and all outsiders 1. The target gives 4, while the original edges give only 0 or 1. In particular C is a subset of the cycle vertices and its induced subgraph has no isolated vertex.

Suppose C is proper and has at least four vertices. Its components are paths with at least two vertices each. In each path choose one bipartition class as selected, with the proviso that vertex 0 must be selected if it belongs to C. Give selected vertices value 2. Each unselected vertex is pure and gets 0 if it is a source and 3 if it is a sink. Assign every outside variable value 1. Each contained edge has ordered values (0,\allowbreak{}2) or (2,\allowbreak{}3), giving 2. Each other edge gives 1. If there are at least two components, each has at least one selected vertex. If there is just one, its size is at least four and either bipartition class has at least two vertices. Thus the target v contains at least two occurrences of 2 and no 1. Its value is 4. This contradicts absorption. Every absorbed word therefore has support size at most three or exactly n. Its length and repeated letters are unrestricted.

\par\smallskip\noindent\textit{Passing from support size to a basis obstruction.}\quad 
Any identity true in 7805 has equal content on its two sides. An absent variable is separated by putting it at 4 and all other variables at 0. Work with finite nonempty sets of nonempty words using distributivity. In any equational derivation from P, all summands of every intermediate polynomial are absorbed by P. Consider a rewrite using an identity on at most k variables, and suppose every summand before it has support of size at most three. Fix a contribution from a one-hole context after distribution. It has the form L times an instance of one side times R; either context word may be absent. Both contexts together use at most three letters, since they occur in each of its old expansions.

Fix a new source monomial. Every image word available to any variable of this monomial has at most three distinct letters: use the content condition to place that choice in an old expansion with the same context. Now select one image word for each variable and use it at every one of its occurrences. The new contextual word uses at most 3k+3 letters. Take n odd and greater than both 3k+3 and six. The gap and validity of this intermediate polynomial force its support to have size at most three.

To reach any other new expansion, change the choice at just one occurrence at a time. Every such choice is an actual summand of the same distributed polynomial, even for repeated source variables. After each change the support has at most three old letters plus three new letters. By the gap it again has size at most three. Thus no rewrite introduces the first full-support word. This includes both directions, arbitrary polynomial substitutions and all multiplication contexts. W uses all n letters. Hence W \ensuremath{\leq} P cannot follow for all n from any finite identity basis. The multiplicative reversal is isomorphic under exchange of 0 and 3.

The proof concerns basis status. Determining whether the variety is a limit variety requires information on its proper subvarieties beyond this argument.

\subsection{Finite bases for 7806, 7945 and 8599}\label{case:201}
Use the constant-free semiring signature. For a target word w and sum P let U consist of the original words supported on C, the content of w. Write K for the union of their contents, J for their nonfinal letters, T for their noninitial letters, H for their heads and F for their finals. Write A(w) and B(w) for its nonfinal and noninitial letters. Define three conditions.

\textbullet{} L holds when H meets T, or the target starts in H and B(w) is contained in T.

\textbullet{} R holds when F meets J, or the target ends in F and A(w) is contained in J.

\textbullet{} E holds when some u in U has A(u) contained in A(w).

For 7806 the exact criterion for w \ensuremath{\leq} P is \(K=C\) together with L and R. For 7945 it is \(K=C\) together with L and E. These statements cover literal words, including literal witnesses. The opposite of 7945 is 8599.

For 7806 the least value 1 annihilates multiplication, the three elements 0,\allowbreak{}2,\allowbreak{}3 are incomparable above it, and 4 is greatest. On words avoiding 1 the value is 0 on all zeros, 2 on zeros followed by a final 2, 3 on an initial 3 followed by zeros, and 4 otherwise. For 7945 the order is \(2<1<0<4\) and \(2<1<3<4\). A 2 anywhere or a nonfinal 1 gives 2. Without either, a final 1 gives 1. On the other letters the value is 0 on all zeros, 3 on an initial 3 followed by zeros and 4 otherwise. Both formulas follow by appending a letter using the catalogue multiplication tables. For 7945 a high accumulated value multiplied by 1 or 2 gives 1 or 2, and an accumulated 1 or 2 followed by any letter gives 2.

For necessity assign the absorbing least value to letters outside C, so only U remains relevant. For sufficiency each relevant contained witness is an actual summand. In 7806 the least bound requires a contained word when the target avoids 1. At bound 0 all high contained words must be all zeros. Coverage suffices, and giving 4 to one missing letter and 0 to the other target letters proves necessity. At bound 3 every contained word must begin with 3 and otherwise be all zeros. This is impossible if H meets T. Otherwise the cover partitions C into H and T, and the forced assignment 3 on H and 0 on T gives the target value 3 exactly under the second alternative in L. At bound 2 the reversed argument, with final value 2, gives R.

For 7945 the bound 2 gives E: assign 0 on A(w), 1 on the rest of C and 2 outside. A right word exceeds 2 exactly when it is a witness for E. Conversely any such witness exceeds 2 whenever w does. If \(w>1\) in the sense of being outside the ideal below 1, all its letters are high, so a contained word suffices. For the bounds 0 and 3 a target outside the bound also uses only high values. The tests at 0 and 3 from 7806 apply on this alphabet, giving \(K=C\) and L.

For each algebra take \(abcd=acbd\), \(abbc=abc\) and every valid absorption with at most four right words, 41 variables and word length at most 43. Both structural identities follow from the word values. The reversible compression in Section~\ref{case:174} applies verbatim: reserve all ten possible endpoint variables and group the remaining interior variables by their at most 31 incidence sets. Each valid step with at most four upper words is derivable. An available word denotes one already derived below P.

In 7806 call the left side crossed when H meets T and the right side crossed when F meets J. Each crossing has at most two original witnesses in U. At an uncrossed left side, validity and coverage force H to consist solely of h(w), which occurs only initially in w and in all of U. Indeed H is disjoint from T and C minus B(w) has size at most one. At an uncrossed right side the reversed statement holds.

If neither side crosses, start with any u in U and restrict w to C(u). The endpoints remain and all further retained letters occur in the interior of u. The restriction is validly absorbed in u. Enlarge the restriction by the content of each further original, using the preceding restriction and that original; each shares the endpoints and supplies the newly needed roles. The argument includes the common literal case.

If both sides cross, retain the at most four original witnesses and their content union S. On support S take a word whose first and final positions are occupied by the same letter, with both positions present. The chosen originals derive this word and the word itself supplies both crossings. At a single crossed left side, take its pair of witnesses. On their union S construct a word ending in the common simple final f and with a repeated head from S minus \{f\}. This set is nonempty: if all witnesses were supported on \{f\}, the uncrossed final would make all of them literal f, precluding a left crossing. List all remaining support letters before f. Every nonfinal letter is supplied by a witness since every selected word ends in the unique simple f. The new word supplies its own left crossing and preserves the required right roles. A single crossed right side uses the reverse construction. In each case enlarge S by one contained original at a time. The old constructed word and that original supply coverage, every active crossing, and the common simple endpoint and its roles on an uncrossed side. Build the same kind of word on the union. At \(S=C\) the last word alone derives w.

For 7945 with a left crossing choose its two original witnesses. On their union derive a word with the same first and final letter in two actual positions, covering their union. Both the crossing and E hold in this step: every letter of the constructed word is nonfinal. It supplies a left crossing on its own. Adjoin further contained contents one at a time, retaining this construction. Use the last word and the original witness for E to obtain w. At an uncrossed left side all U words share the simple head h(w), as above. Start with a witness for E and a contained word covering f(w). Restrict w to their union. All noninitial letters are covered noninitially and the retained endpoints preserve E. At each subsequent enlargement the preceding restriction supplies E, since its nonfinal set is contained in that of the enlargement. The old restriction and the newly included word suffice.

All selections and inductions allow arbitrarily long words and arbitrarily large alphabets and polynomials. Derivability of every word absorption gives both directions of every polynomial identity. Reversal transfers the proof to 8599.

\subsection{Twelve identities for a semiring with a bottom-reset tail}\label{case:202}
Work in the signature (+, multiplication), without constants, relative to the additively idempotent semiring axioms. Write \(w\leq P\) for w+P approximately equal to P. Let A be catalogue entry 13038. Its tables, in numerical element order, are

{\small\par\smallskip
\noindent addition       multiplication\par
\CaseLine{\ttfamily 0 0 0 4 4\qquad 4 4 2 0 4}
\CaseLine{\ttfamily 0 1 1 4 4\qquad 4 4 2 0 4}
\CaseLine{\ttfamily 0 1 2 3 4\qquad 4 4 2 2 4}
\CaseLine{\ttfamily 4 4 3 3 4\qquad 4 4 2 3 4}
\CaseLine{\ttfamily 4 4 4 4 4\qquad 4 4 2 4 4}
\smallskip}

Element labels are not constants in terms. The additive order has \(2<1<0<4\) and \(2<3<4\); 3 is incomparable with 0 and 1.

Theorem. The following twelve identities form a complete relative basis for A. The opposite semiring has the word-reversed basis.

\CaseDisplay{\begin{aligned}
 xyzy&\approx xzy, &(1)\qquad xy^2&\approx xy, &(2)\\
 xyz&\leq xy+tz, &(3)\qquad xy&\leq x^2y, &(4)\\
 xy&\leq yxy, &(5)\qquad xy&\leq y^2+zxy, &(6)\\
 xy&\leq x+zxy, &(7)\qquad xy&\leq y+zxy, &(8)\\
 xyz&\leq y+txyz, &(9)\qquad x&\leq x^2, &(10)\\
 xy&\leq x+y^2, &(11)\qquad x^2y&\leq xy+y^2. &(12)
\end{aligned}}

All substitutions below are nonempty words. Empty suffixes denote positions in a word and never substitute for a variable.

\par\smallskip\noindent\textit{Word values and normal forms.}\quad 
Right multiplication by 3 is the idempotent map r=(0,\allowbreak{}0,\allowbreak{}2,\allowbreak{}3,\allowbreak{}4). Right multiplication by 2 is constant 2; right multiplication by 0, 1 or 4 is constant 4. These right maps preserve the additive joins. For left maps the only incomparable joins needing attention are \(0+3=4\) and \(1+3=4\); their images also join to 4. The comparable pairs are preserved in each row. Thus both distributive laws hold.

For associativity, if the last factor is not 3, both associations give its reset value, 2 or 4. If the last factor is 3 and the middle factor is a reset element, applying r preserves its reset type. If both middle and last factors are 3, use \(r^{2}=r\). Hence A is a semiring. Induction gives the following value of any nonlinear word:

\textbullet{} Locate the last tail factor outside \{3\}. If it is 2 the word equals 2; if it is 0, 1 or 4 the word equals 4.

\textbullet{} If all tail factors are 3, the word equals r of its head.

Identities (1) and (2) hold and remove every earlier repeated tail occurrence. A normal nonlinear word retains its literal head and the distinct tail variables in their last-occurrence order. A head may also appear in its tail. A pure power of degree at least two normalizes to its square. Identity (10) follows from the table.

For a word u write c(u) for its content and \(s_{u}\)(x) for the suffix after the last occurrence of x, when x occurs in u. Original right words are retained when applying the criteria below.

\par\smallskip\noindent\textit{Exact absorption criterion.}\quad 
Lemma. A linear y is absorbed by a polynomial P exactly when P has a word with content \{y\}. A nonlinear normal word \(w=h\) t is absorbed exactly under the following conditions:

\textbullet{} (C) Some original right word has content contained in c(w).

\textbullet{} (O) For every x in c(w), an original right word u contains x with c(\(s_{u}\)(x)) contained in c(\(s_{w}\)(x)). This witness may be linear.

\textbullet{} (T) For each target tail letter x, put \(S=s_{w}(x)\). Either an original nonlinear right word has a tail occurrence of x followed only by letters of S, or an original right word has content contained in S.

\textbullet{} (H) If h is absent from t, either an original nonlinear right word contains h followed only by letters of t, or an original right word has content contained in c(t).

Proof. For a missing pure-y word assign \(y=3\) and all other variables 2. Every right word equals 2, whereas the target is 3. A pure-y word conversely normalizes to y or \(y^{2}\) and (10) applies.

If (C) fails, assign target variables 3 and outsiders 2. The target equals 3 and every right word equals 2. If (O) fails at x, assign \(x=0\), its target suffix 3, and all other variables 2. Every right word containing x has a 2 after its last x, and hence equals 2. The remaining right words have value 2 or 3. The target is 4 when x is in its tail and 0 when x is its simple head, neither below 3.

If (T) fails at x with target suffix S, assign \(x=1\), S the value 3, and all other variables 2. The target is 4. There is no right word entirely in S. Every right tail occurrence of x has a later 2 and gives 2. A right word having x only at its head is 1 if linear, and is 0 or 2 if nonlinear. A right word omitting x has a 2 and equals 2. Thus \(P\leq 0\), contradicting target value 4.

If (H) fails at a simple head h, assign \(h=1\), its target tail 3, and outsiders 2. Every right nonlinear word containing h has a 2 after its last h and equals 2. A word omitting h has an outsider and equals 2, because there is no right word contained in the tail. A literal h has value 1. Consequently \(P\leq 1\), whereas \(w=0\).

For sufficiency, suppose first that the target has a last tail factor outside \{3\}. A reset to 2 is automatically below P. For a reset letter x valued 0,\allowbreak{}1 or 4 the target is 4, and all following target factors are 3. A nonlinear tail witness in (T) therefore has value 4. In the other case (T) supplies a word valued 3, while (O) supplies a right word valued 0,\allowbreak{}1 or 4: its last x is followed only by 3s. Their join is 4.

Now suppose all target tail factors are 3. Head 2 is automatic. Head 3 makes the entire target content 3, and (C) supplies a right word valued 3. Heads 0 and 4 are handled by (O): its right witness has value respectively 0 or 4, or a tail reset to 4. Head 1 gives target 0. The nonlinear witness in (H), when present, has value 0 or 4. Otherwise (H) supplies a word valued 3 and (O) supplies a word valued 1,\allowbreak{}0 or 4; their join is 4. This covers every assignment and proves the lemma. \(\square\)

Applying the criterion to (3)--(9), (11) and (12) verifies their soundness. In particular (11) uses the right literal x as its ordinary occurrence witness, and \(y^{2}\) both for the tail and the contained-tail alternative in (H). For (12), xy supplies the ordinary occurrences, and \(y^{2}\) supplies the contained-suffix alternative for the tail occurrence of x. These alternatives are essential: the criterion for catalogue entry 12500 is stronger and cannot be substituted for this lemma.

\par\smallskip\noindent\textit{A square recovered from a contained suffix witness.}\quad 
Fix a valid nonlinear absorption \(w\leq P\). Let N be the set of original nonlinear right words. It is nonempty: (T) for the final target letter has empty suffix, so its nonlinear-witness alternative is necessary. Suppose a nonempty suffix S of the target tail has already been appended to every u in N, giving \(uS\leq P\). Suppose also an original right word b has content contained in S.

We claim \(S^{2}\leq P\). If b is nonlinear, bS is already known and its normal form is kS with k the head of b. Write \(S=A\) B at k, so B begins with k, and let g be the first letter of S. Identity (5) with \(x=gA\), \(y=B\) has left side normalizing to \(S^{2}\) and right side normalizing to kS. The substituted words remain nonempty if A is empty.

If b is a literal k, take any R in N. The word RS is already known. Put \(v=gS\), the normal form of \(S^{2}\). Since g is in S, Rv normalizes to RS. If \(k=g\), use (7) on \(v\leq k+Rv\). If k is the final letter of v, use (8). Otherwise split \(v=A\) k B, with A, B nonempty, and use (9). All right terms are below P, proving the claim.

\par\smallskip\noindent\textit{Appending the common target tail.}\quad 
Build the target tail from right to left. Assume \(uS\leq P\) has been derived for each u in N, with initially empty S. Let x be the next target letter.

If (T) supplies a nonlinear tail witness v with last x preceded by a nonempty R and followed only by S-letters, normalization gives vS approximately RxS. For each \(u=AB\) in N with A, B nonempty, identity (3) in right context S gives

\CaseDisplay{uxS\leq uS+RxS\leq P.}

Otherwise (T) supplies a contained suffix witness b. The suffix S is nonempty, and the preceding claim derives \(S^{2}\leq P\). Condition (O) supplies a word v containing x followed only by S-letters. If v is nonlinear and its last x is in the tail, use the preceding case. If x occurs only as its head, vS normalizes to xS, already known. If v is the literal x, apply (11) to get \(xS\leq x+S^{2}\leq P\). Now (12) gives

\CaseDisplay{x^2S\leq xS+S^2\leq P.}

Use (3) again with \(t=x\) in right context S to obtain \(uxS\leq uS+x^{2}S\leq P\) for every original nonlinear u. This completes one suffix step. Induction yields \(ut\leq P\) for every u in N. No literal original word is incorrectly split into two nonempty parts.

\par\smallskip\noindent\textit{Head and contained-content recovery.}\quad 
If h is simple and (H) has no nonlinear witness, (O) forces a literal h on the right; (H) supplies a contained witness inside t. The square claim gives \(t^{2}\leq P\). Identity (11) now derives \(w=ht\leq h+t^{2}\leq P\). This handles the new simple-head alternative.

Otherwise, if h is simple, a nonlinear witness in (H) either starts with h and its extension is w directly, or supplies a known Rht with R nonempty. If h is repeated in the target tail, any original nonlinear u supplies Rw approximately ut by taking \(R=u\).

First let (C) have a nonlinear witness b, of head k. If h occurs in t, bt normalizes to kt. If \(k=h\) this is w. Otherwise split \(t=A\) B at k and use (5) with \(x=hA\), \(y=B\) to derive \(w\leq kt\).

Suppose h is simple. If b contains h and \(k=h\), then bt is w when h was simple in b, or is hht when h also occurs in its tail; the latter case uses (4). If b contains h and k differs from h, the same substitution in (5) derives \(w\leq kht\), the normal form of bt. If b omits h, the square claim with \(S=t\) gives \(t^{2}\leq P\), and (6) combines \(t^{2}\) and the independently known Rht to derive w.

Finally, if (C) supplies a literal k, use the already known Rw. When \(k=h\) apply (7), when k is final apply (8), and at an internal position apply (9), exactly as in the square claim. This proves every valid nonlinear absorption using the displayed identities.

\par\smallskip\noindent\textit{Unrestricted completeness.}\quad 
Normalize the pure-variable witness of a linear absorption and apply (10). A nonlinear absorption is derived by the suffix induction and the exhaustive head/content recovery above. The number of suffix steps and original right words is unrestricted; all substituted terms are nonempty. For an arbitrary valid identity, distribute both sides into sums of words, derive the absorption of every summand into the other polynomial, and combine both inequalities. This proves equational completeness for any finite set of variables, hence for all identities. Reversal proves the opposite assertion. \(\square\)

Finite computations are checks of the formulas and the actual substitution/context derivations, not a substitute for this proof. Any transfer to other catalogue types requires complete independent HSP verification in both directions. Minimality and hereditary finite basability are not asserted.

\subsection{Guarded interleaving: bases for sources 13863 and 7844}\label{case:203}
Work in the constant-free signature of additively idempotent semirings. Element names in the following tables are not constants in terms. Multiplication and addition are given by row strings:

13863 addition  multiplication 00004           44434 01114           41234 01224           22232 01234           33333 44444           44434

7844 addition   multiplication 00044           01244 01134           11211 01234           22222 43334           44244 44444           44244

\par\smallskip\noindent\textit{Structure.}\quad 
For source 13863 the square map is (4,\allowbreak{}1,\allowbreak{}2,\allowbreak{}3,\allowbreak{}4). Its image is B=\{1,\allowbreak{}2,\allowbreak{}3,\allowbreak{}4\}, ordered \(3<2<1<4\). On B the element 1 is an identity, 3 is a multiplicative zero, and 2,\allowbreak{}4 are left zeros away from 3. Addition in the original algebra is the chain \(3<2<1<0<4\).

For source 7844 the square map is (0,\allowbreak{}1,\allowbreak{}2,\allowbreak{}4,\allowbreak{}4), with image B=\{0,\allowbreak{}1,\allowbreak{}2,\allowbreak{}4\}, ordered \(2<1<0<4\). The element 0 is an identity, 2 a multiplicative zero, and 1,\allowbreak{}4 are left zeros away from 2. In the original additive order 2 is least, 1 lies below both 0 and 3, and 0 and 3 have join 4.

For both sources \(h(x)=x\)\ensuremath{{}^{2}} preserves addition and multiplication, \(h^{2}\)=h, and \(xy=h(x)h(y)\). These assertions are checked directly from the displayed tables. The four-element image is an ordered monoid with zero whose left and right multiplication maps preserve the displayed chain order. It is therefore an ai-semiring. The retraction formula proves associativity and both distributive laws for the original algebras.

In 13863 a nonlinear word is 3 if any letter is assigned 3. Otherwise it is 1 if all letters are 1; if not, look at the first letter not assigned 1. Its value is 2 if this letter is 2, and 4 if it is 0 or 4.

In 7844 a nonlinear word is 2 if any letter is assigned 2. Otherwise it is 0 if all letters are 0; if not, its value is 1 if the first nonzero-valued letter is 1, and 4 if that letter is 3 or 4.

\par\smallskip\noindent\textit{Bases.}\quad 
Use the five laws

\textbullet{} (W1) x\ensuremath{{}^{2}}y \ensuremath{\approx} xy;

\textbullet{} (W2) xy\ensuremath{{}^{2}} \ensuremath{\approx} xy;

\textbullet{} (W3) xyx \ensuremath{\approx} xy;

\textbullet{} (S) x \ensuremath{\leq} \(x^{2}\);

\textbullet{} (I) xyzt \ensuremath{\leq} xt + yz.

For 7844 add

\textbullet{} (G00) gxu \ensuremath{\leq} gu + g + x;

\textbullet{} (G10) agxu \ensuremath{\leq} agu + g + x;

\textbullet{} (G01) gbxu \ensuremath{\leq} gbu + g + x;

\textbullet{} (G11) agbxu \ensuremath{\leq} agbu + g + x.

Theorem. The five common laws form a basis for 13863. All nine displayed laws form a basis for 7844. Reversing every word gives bases for the opposite algebras. Minimality and independence are not claimed.

\par\smallskip\noindent\textit{Exact semantic criterion.}\quad 
The word laws retain the first-occurrence list of each nonlinear word, or its square when the content is a singleton. Separated repetitions use W3 with a nonempty middle word; adjacent repetitions use W1 or W2. Thus normalization is derivable for all word lengths without turning a nonlinear word into a linear one.

Let P have linear variable set L and nonlinear words N. For a target w put \(C=c(w)\), and let \(N_{C}\) be the nonlinear right words whose content is contained in C. Write \(pre_{w}\)(x) for the variables before the first x. Set J(S) to mean that L meets S or some word of N has content contained in S.

Lemma. A linear y is absorbed exactly when y belongs to L or N contains a singleton-content word on y. For nonlinear w, the criterion in 13863 is:

(F) For every x in C, there is v in \(N_{C}\) containing x such that \(pre_{v}\)(x) is contained in \(pre_{w}\)(x).

In 7844 replace each requirement (F) by the alternative of (F) or

(G) x belongs to L and J(\(pre_{w}\)(x)) holds.

Necessity. For a missing linear witness assign \(y=4\) and every other variable the multiplicative zero (3 in 13863, 2 in 7844). All right words containing another variable become zero; there is no linear y or singleton nonlinear y. The target is 4.

For a failed prefix requirement at x in 13863, assign all preceding target variables 1, assign \(x=0\), all remaining target variables 2, and outsiders 3. The target is 4. Each right nonlinear word containing an outsider is 3. A contained word could reach 4 only through x with all earlier letters assigned 1, which would be a prefix witness. Hence every right nonlinear word is at most 1, and the linear right part is at most 0, strictly below 4.

For 7844 assign preceding variables 0, remaining target variables 1, outsiders 2, and \(x=3\) if x is linear, otherwise \(x=4\). The target is 4. No contained nonlinear word can be 4 in the absence of an ordinary prefix witness. If x is not linear the right side is at most 0. If x is linear, failure of (G) means J(\(pre_{w}\)(x)) is false. Then no right summand has value 0: linear variables valued 0 are absent, and a nonlinear word can be 0 only if all its variables are in this prefix. All remaining right terms are 1,\allowbreak{}2,\allowbreak{}3, whose join is at most 3. Thus the right side cannot absorb 4.

Sufficiency. The first variable of a nonlinear target must have an ordinary witness: J(empty set) is false in the constant-free signature. Consequently \(N_{C}\) is nonempty. A target equal to the multiplicative zero is least. If it equals the identity, every variable is assigned the identity, and a word in \(N_{C}\) also equals the identity. If it equals the lower left zero (2 in 13863, 1 in 7844), every word in \(N_{C}\) is nonzero and is at least that element. Finally, if it is 4, choose the first target variable x not assigned the identity. The preceding variables are identities. An ordinary contained prefix witness is 4 and suffices. In the alternative 7844 case, the linear x is 3 or 4 and J(\(pre_{w}\)(x)) supplies a right summand equal to 0. Its join with x is 4. Linear absorptions follow from (S) or literal occurrence.

The common word and square laws are sound by the table formulas. For (I), xt witnesses the prefix at x, t and yz witnesses y, z, with all contents contained in the target. For each G-law the right word with x deleted supplies all ordinary prefix requirements except possibly at x, where linear x and the earlier linear g meet (G). This proves soundness of all displayed laws.

\par\smallskip\noindent\textit{Unrestricted derivation.}\quad 
Let nonlinear w\ensuremath{\leq}P satisfy the criterion, and normalize its first-occurrence list as \(x_{1}\),..., \(x_{d}\). Put \(p_{i}=x\)\ensuremath{{}_{1}}...\(x_{i}\). Inductively construct a nonlinear \(s_{i}\)\ensuremath{\leq}P with content contained in C and first-occurrence list beginning with \(p_{i}\). An ordinary contained witness for \(x_{1}\) exists and supplies \(s_{1}\).

At step i+1 write \(s_{i}\) approximately \(p_{i}\) q with q nonempty. If the normalized word ends at \(p_{i}\) and i\ensuremath{\geq}2, use the equivalent word \(p_{i}\)\ensuremath{{}^{2}} so \(q=p_{i}\). If \(i=1\) and the support is singleton, the square has a nonempty suffix. These devices never substitute an empty word.

For an ordinary contained prefix witness v, split \(v=ab\) with a, b nonempty. Rule (I) gives \(p_{i}\) v q \ensuremath{\leq} \(p_{i}\) q + v \ensuremath{\leq} P. Every variable in v before the selected x already belongs to \(p_{i}\) and is removed by normalization. The new word starts with \(p_{i}\) x and uses only variables from C.

The remaining case arises only in 7844: x is linear and J(c(\(p_{i}\))) holds. If it has a linear witness g occurring in \(p_{i}\), split \(p_{i}\) at that occurrence. Choose G00, G10, G01 or G11 according to whether the substrings before and after g are omitted or nonempty. The result is \(p_{i}\) x q \ensuremath{\leq} \(p_{i}\) q + g + x \ensuremath{\leq} P.

If instead v is a nonlinear right word with content contained in c(\(p_{i}\)), apply G10 with \(a=p_{i}\), \(g=v\), and \(u=q\). It gives \(p_{i}\) v x q \ensuremath{\leq} \(p_{i}\) v q + v + x \ensuremath{\leq} P. Since all letters of v already occur in \(p_{i}\), (W) reduces the left word to \(p_{i}\) x q and the first right word to \(p_{i}\) q. This is a valid step even when \(p_{i}\) is a single letter; every displayed product remains nonlinear. Again the new word has the required prefix and content.

After d steps the normal form of \(s_{d}\) equals w, including singleton support. Thus every valid nonlinear absorption is derivable. Linear absorptions follow from (S) and singleton normalization or from an existing linear summand. Distributing and deriving every summand absorption in both directions proves all identities on arbitrary variable sets. Reversal proves the opposite statements.

Finite tests check transcription and implementations only. The induction proves completeness without bounds on variables, word lengths, or summand counts. Transfers to further catalogue types require independently verified two-way HSP equality. No hereditary finite-basis assertion is made.

\subsection{Thirteen identities for 7953 with nonlinear prefix witnesses and literal guards}\label{case:204}
Relative to the additively idempotent semiring axioms, in the signature (+, multiplication), without constants.

\par\smallskip\noindent\textit{1. Word values.}\quad 
The addition and multiplication tables, in order 0,\allowbreak{}1,\allowbreak{}2,\allowbreak{}3,\allowbreak{}4, are

{\small\par\smallskip
\noindent addition        multiplication\par
\CaseLine{\ttfamily 0 0 0 4 4\qquad 0 1 2 4 4}
\CaseLine{\ttfamily 0 1 1 3 4\qquad 2 2 2 2 2}
\CaseLine{\ttfamily 0 1 2 3 4\qquad 2 2 2 2 2}
\CaseLine{\ttfamily 4 3 3 3 4\qquad 4 4 4 4 4}
\CaseLine{\ttfamily 4 4 4 4 4\qquad 4 4 4 4 4}
\smallskip}

A literal retains its assigned value. For a nonlinear word, scan for the first nonzero factor. If there is none the value is 0. A first 2 yields 2, a first 3 or 4 yields 4, and a first 1 yields 1 only at the final position, otherwise 2. This follows by induction from the displayed rows. In particular nonlinear words never have value 3. The additive order is \(2<1<0<4\) and \(1<3<4\), with 0 and 3 incomparable.

The distinction between literal 3 and nonlinear occurrence of 3 is essential. For example \(xy\leq y+xz\) fails at \(x=0\), \(y=3\), \(z=2\). It must not be inherited from source 7674.

\par\smallskip\noindent\textit{2. An exact criterion.}\quad 
All witnesses below are original right summands. Write Z(S) if a right word is contained in S. Write E(x, S) for a right word containing x whose first-x prefix is contained in S. Distinguish nonlinear E-witnesses from the literal x; distinguish unmarked witnesses, whose first x is final, from marked witnesses, which have a position after x.

A literal target is absorbed exactly when it occurs as a literal summand. A nonlinear target w is absorbed exactly when:

1. Z(c(w)) holds;

2. at each first occurrence of x, with preceding content S, there is a nonlinear E(x, S) witness, or the literal x together with Z(S);

3. if the final variable x is simple, there is also an unmarked E(x, S) witness or Z(S).

Equivalently, at a simple final position one needs a nonlinear unmarked E-witness, or an arbitrary E-witness together with Z(S), or both a nonlinear marked E-witness and the literal x. The last alternative need not supply Z(S).

For necessity of condition 1, assign target variables 0 and all outsiders 2. For condition 2 at x, assign \(S=0\), \(x=3\) and all others 2. The nonlinear target is 4. A right value is 4 exactly for a nonlinear E-witness; otherwise the join reaches 4 only through a literal \(x=3\) and a Z(S) word valued 0. At a simple final x, assigning \(S=0\), \(x=1\), others 2 gives target 1. Right values reach at least 1 precisely through an unmarked E-witness or a Z(S) word. These assignments also give necessity when variables recur in the right words. A literal target x is separated from all other summands by \(x=1\), all others 2.

For sufficiency, 2 is least, and an all-zero target is absorbed by its contained witness. If the first nonzero factor is 4, any E-witness has value 4. If it is 3, a nonlinear E-witness has value 4; the literal x together with a Z(S) witness gives \(3+0=4\). A first 1 before the end gives 2. A simple final first 1 is absorbed by an unmarked witness valued 1, or a Z(S) word valued 0. These cover all assignments and prove the criterion with no bound on variables, words or summands.

\par\smallskip\noindent\textit{3. The finite basis.}\quad 
The following thirteen identities are complete:

{\small\par\smallskip
\CaseLine{\((W1) xxy = xy\)}
\CaseLine{\((W2) xyx = xyy\)}
\CaseLine{\((A)  xyz \leq  yz+xt\)}
\CaseLine{\((B)  xyz \leq  xz+xyt\)}
\CaseLine{\((U)  xyz \leq  xz+xy\)}
\CaseLine{\((C)  xyz \leq  xz+xx+y\)}
\CaseLine{\((D)  xy \leq  xx+xyt\)}
\CaseLine{\((E)  xy \leq  xx+y\)}
\CaseLine{\((G1) xx \leq  x+xy\)}
\CaseLine{\((G2) xyy \leq  x+xyz\)}
\CaseLine{\((G3) xyy \leq  y+xyz\)}
\CaseLine{\((G4) xyzz \leq  y+xyzt\)}
\CaseLine{\((F)  xy \leq  y+xyt\)}
\smallskip}

The criterion proves soundness of each displayed absorption; the word-value formula proves W1 and W2. Their first-occurrence normal form is F(w) if the final variable is simple, and F(w) followed by its last listed variable otherwise. A literal remains literal. The same explicit word-law reductions apply as follows: a separated internal repeated x is removed by \(xYxZ=xY^{2}Z=xYZ\), with Y, Z nonempty; adjacent internal repeats use W1. A final repeat in p x Y x becomes p x \(Y^{2}\) by W2, and internal repeats are removed. No identity variable is substituted by an empty word.

\par\smallskip\noindent\textit{4. Prefix construction and restoring a contained literal.}\quad 
Normalize a valid nonlinear target. Its first-occurrence list is \(x_{1}\)...\(x_{d}\). Let N be the set of nonlinear original right words. Condition 2 at the first position ensures an original witness \(x_{1}\) t with t nonempty, so N is nonempty. We maintain

{\small\par\smallskip
\noindent p u \ensuremath{\leq} P for every u in N,\par
\smallskip}

where p is the prefix of distinct target variables already treated. For \(u=YZ\) with both factors nonempty, A with \(x=x_{1}\), \(y=Y\), \(z=Z\) gives \(x_{1}\) \(u\leq u+x_{1}\) \(t\leq P\). This starts the construction. It is deliberately applied only to nonlinear u; prefixing arbitrary original literals would be invalid.

We will use a guard lemma. If the invariant holds for a nonempty p and Z(c(p)) holds, then \(p^{2}\leq P\). For a contained nonlinear original r, the invariant gives p \(r\leq P\) and p r normalizes to the same word as \(p^{2}\). Otherwise take a contained literal r and any v in N. Write \(p=A\) r B at its occurrence. If both A and B are nonempty, G4 with \(x=A\), \(y=r\), \(z=B\), \(t=v\) gives the normal form of \(p^{2}\) below r+p v. If A is empty use G2; if B is empty use G3; if both are empty use G1. Thus every substitution is nonempty, and all four endpoint cases derive the guard from original r and the maintained p v.

Now let x be the next marked target position. A nonlinear marked E-witness v gives p \(v=p\) x t modulo W1, W2, with t nonempty. B yields p x \(u\leq p\) u+p x t for every u in N. A nonlinear unmarked E-witness gives p \(v=p\) x; U yields p x \(u\leq p\) u+p x. If there is no nonlinear E-witness, condition 2 supplies the literal x and Z(c(p)). Derive \(p^{2}\) by the guard lemma, and use C to get p x \(u\leq p\) \(u+p^{2}+x\leq P\). Hence the invariant extends to p x in all cases.

If the target is saturated, every first occurrence is marked. A contained nonlinear original r gives \(p_{d}\) r equal to the saturated target. If the contained word is only a literal, the guard lemma gives \(p_{d}^{2}\), also equal to that target.

If the target is unsaturated, it is p s with p nonempty and s a simple final variable. A nonlinear unmarked E-witness v gives p \(v=p\) s directly. Otherwise, if a nonlinear marked E-witness and the literal s both exist, the invariant gives p s \(t\leq P\) with t nonempty, and F gives p \(s\leq s+p\) s t. This branch needs no contained prefix word. In all remaining cases, conditions 2 and 3 give an E-witness and Z(c(p)). Derive \(p^{2}\) by the guard lemma. A nonlinear marked witness gives p s \(t\leq P\) and D yields p \(s\leq p^{2}+p\) s t. A literal s instead gives p \(s\leq p^{2}+s\) by E. All cases end with the original target below P.

\par\smallskip\noindent\textit{5. Completeness.}\quad 
The construction derives every valid word absorption, with no bound on its original variables, word lengths or number of right summands. Literal targets are already summands. Distributivity expands arbitrary terms into nonempty polynomials; deriving the absorption of each summand in both directions proves every identity. Thus the thirteen displayed identities are a finite basis for 7953. Reversal gives a finite basis for its opposite. Further catalogue transfers require complete independently checked HSP inclusions in both directions.

Finite assignment checks and literal replay support this unrestricted proof; their numerical bounds are not assumptions of completeness.

\subsection{A finite basis from two endpoints and five witnesses}\label{case:205}
Source 7957 in the five-element catalogue. All signatures are \((+, multiplication)\) without named constants; bases are relative to the additively idempotent semiring axioms. The numerical labels below name elements, not constant symbols. No claim about hereditary finite basability is made.

\par\smallskip\noindent\textit{1. The algebra and word values.}\quad 
Let A have universe \{0,\allowbreak{}1,\allowbreak{}2,\allowbreak{}3,\allowbreak{}4\}, with the following tables (row and column order 0,\allowbreak{}1,\allowbreak{}2,\allowbreak{}3,\allowbreak{}4):

{\small\par\smallskip
\noindent addition       multiplication\par
\CaseLine{\ttfamily 0 0 0 4 4\qquad 0 1 2 4 4}
\CaseLine{\ttfamily 0 1 1 4 4\qquad 2 2 2 4 4}
\CaseLine{\ttfamily 0 1 2 4 4\qquad 2 2 2 4 4}
\CaseLine{\ttfamily 4 4 4 3 4\qquad 3 3 3 4 4}
\CaseLine{\ttfamily 4 4 4 4 4\qquad 4 4 4 4 4}
\smallskip}

Its additive order consists of the chain \(2<1<0<4\) and a second branch \(3<4\), with 3 incomparable to 0,\allowbreak{}1,\allowbreak{}2. Put B=\{0,\allowbreak{}1,\allowbreak{}2\}. A product of at least two factors has value 4 if its first factor is 4 or a later factor lies in \{3,\allowbreak{}4\}; it has value 3 if its first factor is 3 and all later factors lie in B. If all factors lie in B, its value is 0 when they are all 0, 1 when all factors except the last are 0 and the last is 1, and 2 otherwise. These statements follow directly by induction from the multiplication table.

For completeness, associativity can also be checked without any word assumptions. For a triple, a first factor 4 or a second/third factor outside B makes both associations 4. A first factor 3 and the others in B makes both 3. For three factors in B, both associations are 0 exactly when all are 0, and 1 exactly when the first two are 0 and the last is 1; otherwise both are 2. Each row and column of multiplication preserves joins: it is increasing on \(2<1<0\), sends 4 to 4, and sends the join of 3 with any B-element to the join of their images, namely 4. Thus both distributive laws hold and A is an additively idempotent semiring.

Write c(w) for the content, h(w) for the first variable, t(w) for the last variable, i(w)=c(w without its last position), and r(w)=c(w without its first position). The last two are used only for words of length at least two. An interior position is neither first nor last.

The following two identities hold by the word-value formula:

\CaseDisplay{xyzt\approx xzyt,\qquad xy^2z\approx xyz.                 \quad\textup{(\ensuremath{W})}}

Fix an ordering of the variables. Laws (W) transform every word of length at least two into

\CaseDisplay{h\,\prod_{x\in M}x\,t,                     \quad\textup{(\ensuremath{1})}}

where M is its set of interior variables, listed in increasing order. Indeed the first law swaps adjacent interior positions using nonempty prefixes and suffixes as outer substitutions. After sorting the interior, the second removes its adjacent repetitions. The first and last positions remain, even if their variables also belong to M. Empty M means the two-letter word ht, not the use of a unit constant. A word of length one is left unchanged. Call these the normal words.

\par\smallskip\noindent\textit{2. An exact absorption criterion.}\quad 
Distribute a term into a nonempty sum P of words. Let L be the set of its linear summands and N its nonlinear summands. Set

\CaseDisplay{D=\bigcup_{u\in P}c(u),\qquad
 E=\bigcup_{u\in N}r(u),\qquad
 Q=L\cup\{h(u):u\in N\}.}

Here and below a linear summand is identified with its variable; hence \(D=E\cup Q\). These sets are unchanged by normalization.

Lemma 1. A linear variable y satisfies A |= \(y\leq P\) if and only if y occurs as a linear summand of P. For a nonlinear target w, put \(C=c(w)\), \(I=i(w)\), \(T=r(w)\), and \(h=h(w)\). Then A |= \(w\leq P\) if and only if:

1. C is contained in D;

2. either L meets C, or some v in N satisfies c(v) contained in C and i(v) contained in I;

3. either E meets Q, or h belongs to Q and T is contained in E.

All uses of \ensuremath{\leq} abbreviate an additive absorption identity.

Proof. If a linear target y has no literal summand, assign \(y=1\) and every other variable 2. Every nonlinear summand is 2 and every linear summand is 2, whereas \(y=1\). The converse is immediate.

For condition 1, a variable x in C outside D can be assigned 4 while all other variables are assigned 0. Then \(w=4\) and \(P=0\).

If condition 2 fails and t(w) does not belong to I, assign the variables in I the value 0, the last variable the value 1, and every outside variable 2. The target is 1. Every linear summand is 2. A nonlinear summand either contains an outside variable, giving 2, or has the last variable of w before its own final position, again giving 2. If instead t(w) belongs to I, then \(I=C\). Assign C the value 0 and its complement

2. Failure of condition 2 gives \(P=2\) while \(w=0\).

For condition 3, assume condition 1 and E disjoint from Q. Assign Q the value 3 and all other variables 0. Every summand of P is 3, so \(P=3\). The nonlinear target is 3 precisely if h is in Q and every tail variable is in E: all of C is in \(D=E\) union Q. If either requirement fails, its value is 0 or 4, neither below 3. This proves necessity.

For sufficiency consider an arbitrary assignment. If \(P=4\) there is nothing to show. If P lies in B, each summand lies in B, and therefore every variable in D lies in B. Condition 1 puts all variables of w in B as well. If \(w=2\), it is below P. If \(w=0\), all of C has value 0, and a witness in condition 2 has value 0. If \(w=1\), all of I has value 0 and the last variable has value 1. A linear witness in C then has value 0 or 1; a nonlinear witness has all its nonfinal variables 0 and its final variable in C, so also has value 0 or 1. In each case the witness shows \(w\leq P\). Finally \(P=3\) forces every summand to be 3; equivalently Q is assigned 3 and E is assigned elements of B. In particular E and Q are disjoint. Condition 3 then forces \(w=3\). These cases exhaust the algebra. \ensuremath{\square}

\par\smallskip\noindent\textit{3. Absorption in steps with at most five right summands.}\quad 
All words in this section are normal. Suppose \(w\leq P\) is valid. A word already occurring in P needs no step, and Lemma 1 handles linear w. For nonlinear w write its normal form as h (product of M) t. Choose a witness v in P for condition 2 of Lemma 1. Define

\CaseDisplay{M_0=c(v)\setminus\{h, t\},}

and also include t in \(M_{0}\) when t differs from h and v is nonlinear with t in i(v). Then \(M_{0}\) is contained in M. This follows from the witness requirements: every variable of c(v) other than h, t must be interior to w; a nonfinal occurrence of t in v requires t in i(w). Let \(s_{0}=h\) (product of \(M_{0}\)) t. The same v witnesses condition 2 for \(s_{0}\): its content is covered, and all its nonfinal variables occur before the final position of \(s_{0}\). This also holds for linear v.

Crossing case. Suppose E meets Q. Choose g in their intersection, a nonlinear a in P with g in r(a), and b in P whose first (or only) variable is g. Choose \(r_{h}\), \(r_{t}\) in P containing h, t respectively; these exist by condition 1. Lemma 1 directly gives

\CaseDisplay{s_0\leq a+b+v+r_h+r_t.                    \quad\textup{(\ensuremath{2})}}

The content condition follows from v, \(r_{h}\), \(r_{t}\), condition 2 has witness v, and the crossing g is witnessed by a, b. Repeated right summands can be removed. For each missing interior variable x in M minus \(M_{0}\), choose an original summand \(r_{x}\) containing x and let s' be obtained from the current s by adding x to its interior. Then

\CaseDisplay{s'\leq s+a+b+r_x                  \quad\textup{(\ensuremath{3})}}

is valid by the same criterion: s is the condition-2 witness, a, b retain the crossing, and \(r_{x}\) supplies the new content.

Separated case. Suppose E is disjoint from Q. Condition 3 gives h in Q and r(w) contained in E. Thus h differs from t and h is not an interior variable of w. The chosen v must have first variable h: a linear witness lies in L intersect C, a subset of Q intersect \(C=\{h\}\); the head of a nonlinear witness also belongs to Q intersect C. Every noninitial variable of v belongs to E. Choose a nonlinear \(r_{t}\) in P with t in r(\(r_{t}\)), which exists since t is in E. Then

\CaseDisplay{s_0\leq v+r_t.                \quad\textup{(\ensuremath{4})}}

Indeed h is a head on the right, every tail variable of \(s_{0}\) occurs in the tail of v or \(r_{t}\), and v witnesses condition 2. For each new interior x choose an original nonlinear \(r_{x}\) with x in its tail. Adding x to the interior of s gives the valid absorption

\CaseDisplay{s'\leq s+r_x.                 \quad\textup{(\ensuremath{5})}}

Here the right-hand head set contains h, its tail set covers that of s', and s witnesses condition 2. These assertions remain valid whether or not the smaller right-hand sum happens to have a crossing.

In both cases finitely many steps end at w. Each right side consists of original summands or the previously constructed word, and has at most five summands. Consequently, if all valid absorptions with at most five normal right summands can be derived, every valid absorption can be derived. The argument uses no bound on the variables of P, w.

\par\smallskip\noindent\textit{4. Reversible compression and a finite basis.}\quad 
Lemma 2. Every valid absorption with a normal left word and at most five normal right words is a substitution instance, modulo (W), of a valid normal absorption using at most 75 variables and words of length at most 77.

Proof. There are at most six words altogether. Reserve a distinct variable for each original variable that occurs as the first or last variable of any word. There are at most 12 such variables; linear words are included in this reservation. Every other variable occurs only in interiors. Partition those remaining variables by their nonempty incidence pattern in the interiors of the at most six words. There are at most \(2^{6}-1=63\) classes. Identify the variables of each class to one fresh variable and normalize the result using (W). This is ordinary substitution in a valid identity followed by valid word laws, so the compressed absorption is valid. At most \(12+63=75\) variables occur and each normal word has at most 75 interior variables plus two endpoint positions, hence length at most 77.

Conversely, substitute for each class variable the nonempty product, in fixed order, of its original variables. Substitute the single original variable for each reserved endpoint variable. A class occurs in a compressed word's interior precisely when every member occurred in the original interior. The expansion restores exactly the original interior set. The first and last variables have never been identified with other variables. Thus (W) normalizes each expanded word to its original normal form. All substitutions are nonempty; no constants or empty words have been introduced. \ensuremath{\square}

Fix a 75-letter alphabet X. Let Sigma consist of (W), together with every absorption \(w\leq u_{1}+\cdots +u_{k}\) valid in A, where \(1\leq k\leq 5\) and all words are normal words over X of length at most 77. This explicitly defined set is finite. Validity of each possible member can be decided using Lemma 1 (or the finite table); no oracle is involved. We do not enumerate this very large set or claim it is minimal.

Theorem. Sigma is a finite identity basis of A. Reversing all words in Sigma gives a finite basis for the multiplicative opposite.

Proof. Soundness is built into the definition and the verification of (W). Distribute any terms into sums of words and normalize using (W). For a valid absorption \(w\leq P\), construct the sequence (2)--(5). By Lemma 2, each step follows from a member of Sigma by substitution and (W). Its right summands are already derivably below P, so additive monotonicity and transitivity show that its left word is below P. At the end \(w\leq P\) is derived. A valid polynomial identity P\ensuremath{\approx}R is equivalent to every summand of P being below R and every summand of R being below P, so the preceding argument derives it. There is no restriction on the number of variables or original word lengths. Word reversal transports both validity and derivations to the opposite algebra. \ensuremath{\square}

\par\smallskip\noindent\textit{5. Role of the calculations.}\quad 
The unrestricted proof is Sections 1--4. Finite regression can detect mistakes in the value formula, criterion, step construction, and compression; it is not a completeness argument. Transfer to other catalogue entries requires independently verified equality of their generated varieties (both HSP directions), or a separately proved transfer theorem. No classification is inferred from a fingerprint or a single inclusion alone.

\subsection{Three finite bases from role domains and bounded component supports}\label{case:206}
Work in the semiring signature without constants. Every word is nonempty. The finite families below are effectively specified, not enumerated and not claimed small.

The three catalogue multiplication/addition pairs give these nonlinear word values, by induction on length. Literals always retain their value.

8813: normal tail value \(b=0\). A head in\{0,\allowbreak{}1,\allowbreak{}2\} followed only by 0 retains its value; first values 3,\allowbreak{}1 followed only by 0 give 2; otherwise the value 4. Its addition is flat with top 4 and four singleton proper ideals.

12589: \(b=2\). A head in\{0,\allowbreak{}1\} followed only by 2 gives 0, a head 2 or 3 followed only by 2 retains its value; head 1 followed by 0 or 1 and then only 2 gives 3; otherwise the value 4. Its proper ideals are\{1\},\{0,\allowbreak{}1\},\{2\},\{3\}.

12885: \(b=2\). A head in\{0,\allowbreak{}1,\allowbreak{}2\} followed only by 2 retains its value; first values 3,\allowbreak{}1 followed only by 2 give 0; otherwise the value 4. Its proper ideals are\{1\},\{0,\allowbreak{}1\},\{2\},\{3\}.

These descriptions are derived from the actual tables, not transferred on the basis of a shared short fingerprint.

\par\smallskip\noindent\textit{Exact constraint systems.}\quad 
For original P define \(C=content\) union, \(L=literal\) variables, \(A=nonlinear\) heads, \(B=second-position\) variables, \(D=variables\) at positions 3 onward, \(H=L\) union A, \(T=B\) union D. Edges are original first-two-position pairs. Heads sharing a second belong to the same component. Equivalently, make separate head and second role vertices for each variable and join the two role vertices along each original edge. This is a bipartite role graph; its connected components are completeable by rectangles.

The special ideal U and its two component modes are as follows. Entries in a set are freely chosen subject to overlaps at the same variable.

\noindent Source\enspace\textbar\enspace U\enspace\textbar\enspace Mode 0 heads\enspace\textbar\enspace Mode 0 seconds\enspace\textbar\enspace Mode 1 heads\enspace\textbar\enspace Mode 1 seconds\enspace\textbar\enspace b\par
\noindent 8813\enspace\textbar\enspace \{2\}\enspace\textbar\enspace \{2\}\enspace\textbar\enspace \{0\}\enspace\textbar\enspace \{3\}\enspace\textbar\enspace \{1\}\enspace\textbar\enspace 0\par
\noindent 12589\enspace\textbar\enspace \{3\}\enspace\textbar\enspace \{3\}\enspace\textbar\enspace \{2\}\enspace\textbar\enspace \{1\}\enspace\textbar\enspace \{0,\allowbreak{}1\}\enspace\textbar\enspace 2\par
\noindent 12885\enspace\textbar\enspace \{0,\allowbreak{}1\}\enspace\textbar\enspace \{0,\allowbreak{}1\}\enspace\textbar\enspace \{2\}\enspace\textbar\enspace \{3\}\enspace\textbar\enspace \{1\}\enspace\textbar\enspace 2\par
All original values are in U exactly when every edge's roles use the same mode, every literal variable is in U, and every D variable equals b. Shared seconds force the same head mode because the two second domains are disjoint. Conversely these conditions imply the original word-value requirements. No assignment possibilities have been discarded.

At each variable intersect its literal, D, head-role, and second-role domains. Every local constraint is a conjunction of unary mode pins, or is inconsistent; there is no additional relation between two component modes. This can also be checked explicitly as follows.

For 8813, overlap of head and second roles is impossible. A literal head forces mode 0; a literal second or literal D occurrence is impossible. A D head is impossible, a D second forces mode 0.

For 12589, head/second overlap forces both modes1 and value 1. A literal head forces mode 0, a literal second or literal D occurrence is impossible. A D head is impossible, a D second forces mode 0. Literal/head/second overlap combines the stated pins and is thus inconsistent.

For 12885, head/second overlap forces head mode 0 and second mode 1, value 1. A literal head forces mode 0, a literal second forces mode 1, and literal D overlap is impossible. A D head is impossible, a D second forces mode 0. Additional simultaneous roles merely intersect these restrictions.

Thus all local role combinations are covered by intersections of the listed restrictions. A global U-assignment exists exactly when no local intersection is empty and no role component receives conflicting pins. When consistent, choose every unpinned component mode independently and choose every remaining actual variable value independently in its intersection domain. This parametrizes all original U-assignments.

The other three proper ideals give these common absorption conditions:

C: target content is contained in C; N: target is literal, or an original nonlinear word exists; I: H intersects T, or target head lies in H and all noninitial target variables lie in T.

For 8813, C comes from ideal\{0\}, N from\{3\}, and I from\{1\}. For 12589, C comes from\{2\}, N from\{1\}, and I from\{0,\allowbreak{}1\}. For 12885, C comes from\{2\}, N from\{3\}, and I from\{1\}. For C, every original lies in the stated singleton ideal exactly when all its variables equal b. For N, the stated ideal contains no nonlinear value. For I, original heads lie in the head ideal(\{1\},\{0,\allowbreak{}1\},\{1\}, respectively), and all noninitial variables equal b. These domains are disjoint, so the constraint is inconsistent precisely when H intersects T.

The fourth condition IU says that every U-assignment described by the mode system also places the target in U. It is automatically true for an inconsistent system. For a literal, it tests membership of its value in U. For a nonlinear target hsZ, it requires the first pair to lie in the union of the two mode rectangles in the table, and every variable in Z to equal b. Component modes are independent apart from their pins, so this is an effective exact criterion for any variable number.

The conjunction C, N, I, IU is necessary and sufficient for \(w\leq P\). These are all proper principal ideals; apply the corresponding condition to the actual value of P to prove sufficiency. For necessity, failed C is separated by a missing variable 4 and all others b; failed N by assigning all variables the nonlinear-absent ideal value. A failed I head is separated by giving original heads 1,\allowbreak{}0,\allowbreak{}1 respectively and all other variables b; the missing head then lies outside the head ideal. For a missing noninitial variable, use 3,\allowbreak{}3,\allowbreak{}0 respectively if it is not an original head, and the head value otherwise. Originals stay in their head ideal while the target leaves it. Failed IU has an actual separating assignment by the exact mode parametrization.

\par\smallskip\noindent\textit{Five fixed identities and component completion.}\quad 
Each of the three sources satisfies

{\small\par\smallskip
\CaseLine{\(xyzt = xytz\)}
\CaseLine{\(xyzz = xyz\)}
\CaseLine{\(xy \leq  xyz\)}
\CaseLine{\(xz \leq  xy+uy+uz\)}
\CaseLine{\(xyzt \leq  xyz+xyt.\)}
\smallskip}

The first two normalize words by retaining the first two positions and the set of subsequent variables; all swaps/deletions occur after those protected positions. Prefix extraction follows from the word values. For the rectangle law, the three right edges put the target roles in the same component, so IU holds; C, N, I also hold. Tail merging follows by testing the first pair and each tail variable separately in IU, and by the other three criteria. These are proofs of soundness, not conclusions drawn from bounded tests.

Extract original quadratic edges. The rectangle law derives every pair between head and second roles in a component:along a sequence of heads sharing successive seconds, propagate a desired final second backward. This works for paths of every length. Denote the completed edges E* and the enlarged polynomial \(P+=P\) union E*. Each added word is below P. In the completed role graph, opposite roles connect in one edge and same-side roles in at most two edges. Original literal and D witnesses are never replaced by merely derived positions.

\par\smallskip\noindent\textit{Bounded supports: at most 65 upper words.}\quad 
Consider a valid target of length at most 3. For any variable x, its entire local restriction can be witnessed with at most 4 original words: its literal if present, one D-occurrence if present, one head-role word, and one second-role word. These retain all relevant restrictions at x.

If the original U-system is inconsistent locally, those at most 4 words already witness it. If a component has conflicting pins, retain the local restrictions at the two pin variables(at most 8 words), and link their pinned role vertices by at most 2 completed edges. Thus a global inconsistency has at most 10 witnesses in P+.

Otherwise keep the local restrictions at each target variable, at most 12 words. There are at most 6 target role vertices. For each original component represented among them, link those target vertices by a star of paths of length at most 2, using at most 10 edges overall. If such a component has a pin, choose one actual original pin, retain its variable's local restriction(at most 4 words), and link its role vertex to the star center(at most 2 edges). There are at most 6 such components, so this costs at most 36 words. The total IU support is at most 58 words.

Every assignment satisfying this support has exactly the necessary local restrictions, target-component identifications, and target-component pins of the full system. Its target values therefore extend to a full original U-assignment:complete each represented component with its chosen mode, choose the other components freely subject to their pins, and choose actual values independently in the nonempty local domains. Extra restrictions caused by support words can only shrink this set. Hence the support already ensures IU for the target. This extension argument covers repeated target variables and components linked through different roles of the same variable.

Separately add at most 3 original content witnesses, one original nonlinear word if needed, and at most 3 original I witnesses. An H/T intersection uses two originals; otherwise retain a target head witness and at most two noninitial witnesses. All four exact conditions now hold with at most 65 upper words in P+. The bound is uniform in the size of P.

\par\smallskip\noindent\textit{Reversible compression and full finite bases.}\quad 
Normalize the target and selected upper words(at most 66 words). Reserve all first/second-position variables, at most 132. Partition all remaining variables by their nonempty tail-incidence pattern among those 66 words, at most \(2^{66}\)-1 groups. Replace each group by a fresh variable and retain each reserved variable separately. This is a substitution into a valid absorption and remains valid. The total bound is \(B=2^{66}-1+132\) variables, word length B+2, and at most 65 upper words.

For the inverse substitution, restore each group variable to the nonempty product of its members, and each reserved variable to its singleton. The first two positions and each subsequent variable set are recovered exactly. The first two fixed laws then literally recover all normal words. No empty block is used; short words contain reserved variables only; repeated reserved variables in tails remain present.

For each source, take the five fixed laws and every table-valid absorption on a fixed B-variable alphabet of word length at most B+2 and at most 65 upper words. This is finite and decidable by the exact ideal criterion, although the conservative family is not feasibly enumerated. Inverse compression derives all valid short targets from P+, hence from P. A valid long normal target hsZ has every triple hsz(z in Z) valid by the criterion. Derive those triples, merge their nonempty tails using the last fixed law, and normalize. Thus every valid absorption is derived with no bound on original variables, lengths, or summand counts. Distributivity and both absorption directions prove full completeness.

Therefore 8813,\allowbreak{}12589,\allowbreak{}12885 are finitely based. Opposite bases are obtained by reversal. Other catalogue types require independently complete mutual HSP; short fingerprints and single-direction membership never suffice.

\subsection{A 329-identity finite basis for 8301}\label{case:207}
Signature(+, multiplication), no constants; all substitutions are nonempty. \(p\leq q\) abbreviates \(p+q=q\). This proof does not infer finite basability from the previously checked inclusion in the 13066 variety.

\par\smallskip\noindent\textit{Operations and word values.}\quad 
Addition rows:00444,\allowbreak{}01244,\allowbreak{}42244,\allowbreak{}44434,\allowbreak{}44444. Multiplication rows:01444,\allowbreak{}11444,\allowbreak{}41444,\allowbreak{}41234,\allowbreak{}41444. The additive order has \(1<0<4\) and \(1<2<4\), and \(3<4\), with no other strict comparabilities. Table induction gives the following disjoint word-value cases: all factors 0 gives 0; rightmost non0 factor 1 gives 1; all factors 3 gives 3; final factor 2 with all preceding factors 3 gives 2; all remaining cases give 4. The last case includes repeated2.

The two laws \(xxy=xy\), \(xyx=yxx\) give last-occurrence normal form, doubling the final letter if and only if it repeats. Internal repeats xYxZ reduce by \(xYxZ=YxxZ=YxZ\), and a final repeat by \(xYx=Yxx\). Empty Y is handled directly by \(xxy=xy\); Z is nonempty in the internal case.

\par\smallskip\noindent\textit{Exact absorption theorem.}\quad 
Let P be the original upper words, \(C=c(w)\). For each target x let \(S_{x}\) be its last-occurrence suffix content and \(W_{x}\) all original u containing x with their last-x suffix content contained in \(S_{x}\). Require a contained original word U with c(U) contained in C and every \(W_{x}\) nonempty. A witness u in \(W_{x}\) is strong when it repeats x or contains a letter of \(S_{x}\). If a strong witness exists there is no extra condition at x. Otherwise every \(W_{x}\) word is simple-final x and avoids \(S_{x}\). Put \(A_{x}=\cup\) of its prefixes. Require either: G1. an original word is contained in \(S_{x}\); or G2. for some a in \(A_{x}\), an original word has last-a suffix contained in \(S_{x}\); or G3. x is the simple final target letter and the target prefix is contained in \(A_{x}\). Then, and only then, \(w\leq \sum\) P. G1/G2 are activation guards, while G3 is the unactivated simple-tail case. All these witnesses are original.

Necessity: if no contained original word exists, assign C letters 0 and outside letters 1. The target is 0 and every upper word is 1. If \(W_{x}\) is empty, assign \(x=4\), \(S_{x}=0\), and all remaining letters 1. The target is 4 and upper values lie in\{0,\allowbreak{}1\}. If the extra condition fails, all \(W_{x}\) words are weak as above. Assign \(x=2\), \(S_{x}=0\), \(A_{x}=3\), all remaining letters 1. The sets are disjoint. Each \(W_{x}\) word is 2. No other upper word can end, after deleting0-valued factors, in an \(A_{x}\) letter: that would be G2. No upper word is all 0: that would be G1. Thus every upper value is 1 or 2. Failure ofG3 ensures target 4, either by a nonempty suffix after x, a repeat of x, or a simple-final prefix factor outside \(A_{x}\). This separates the absorption.

Sufficiency: target 0 and 3 are handled by the contained word. Target 1 has rightmost non0 letter \(x=1\); \(W_{x}\) supplies a word with value 1. Target 2 has a simple final \(x=2\); \(W_{x}\) supplies an upper 2 or 4. For target 4 inspect its rightmost non0 letter x. If it is 4 then \(W_{x}\) supplies 4. If it is 3, a \(W_{x}\) word has value 3 or 4. If only 3 occurs there, the target must contain some factor not 3; its own \(W_{y}\) witness cannot have value 3, so the upper join is 4. If \(x=2\), a strong \(W_{x}\) witness has value 4. If all witnesses are weak and no upper 4 occurs, all \(A_{x}\) factors are 3. G1 supplies 0, whose join with the upper 2 is 4. G2 supplies 3 or 4, again giving 4. In the absence of these guards, G3 would force the target to be 2, contrary to target 4. These cases prove exact semantics without any bound on variables, lengths or sums.

\par\smallskip\noindent\textit{Finite laws.}\quad 
Take the 89 laws listed in the unchanged13066 theory Section~\ref{case:296}, and append the following240 identities. All bracket choices specify separately written identities; no variable ever receives an empty word. Capitals are ordinary distinct variable names.

{\small\par\smallskip
\CaseLine{\([v]a[r]tt \leq  [v]a[r]t+[w]a\qquad (8)\)}
\CaseLine{\([A]x[d]zz \leq  [A][d]zz+[v]a[r]x+[w]a[d]\qquad (32)\)}
\CaseLine{\([A]xyzz \leq  [A]yzz+[v]x+y\qquad (4)\)}
\CaseLine{\([A]bx[d]p \leq  [A]xb[d]p+[v]a[r]x+[w]a[d]\qquad (32)\)}
\CaseLine{\([A]bx[c][d][p]t \leq  [A]xb[c][d][p]t+[v]a[r]x+[w]a[c]t[d]\qquad (128)\)}
\CaseLine{\([A]bxyp \leq  [A]xbyp+[v]x+y\qquad (4)\)}
\CaseLine{\([A]bx[c][d][p]t \leq  [A]xb[c][d][p]t+[v]x+[c]t[d]\qquad (32)\)}
\smallskip}

The largest template has eleven variables. The absorption criterion verifies all 329 instances.

\par\smallskip\noindent\textit{Unrestricted completeness.}\quad 
Normalize original P and target w. Let t be its final variable. Choose an original contained U. If an original t-tail word repeats t, begin with it as a saturated accumulator g. Otherwise test G2 at t (G1 is impossible for the empty suffix). Its original witnesses have forms V a R t and W a. The first new template derives V a R tt, providing a saturated seed. If neither seed is possible, exact semantics requires a simple target t and original simple-tail prefix coverage. Merge those prefixes with \(xyz\leq xz+yz\). This accumulator contains all target content and has simple t. A one-letter simple target must already be an original word in this unactivated case.

Maintain a correct suffix S and process other target letters x in reverse last-occurrence order. In the saturated case S ends in tt; in the simple case it ends in one t. Ordinary strong original witnesses use exactly the three saturated-insertion or five simple-tail-movement cases in the 13066 proof: a nonempty suffix in S, a repeated terminal x, or a preceding S-letter. Their law soundness is checked for 8301, so reuse of these literal deductions does not import13066 semantics.

For a weak \(W_{x}\), exact semantics requires G1 or G2 because x is not the final target letter. Under G1 choose original y contained in S and original Vx in \(W_{x}\). In the saturated case the third new template inserts x before S, setting \(z=S\). In the simple case write the current accumulator A x B S. When B is empty no movement is needed. Otherwise if y omits t use the sixth template with \(p=S\). If y contains t, split \(y=C\) t D and use the seventh template with \(p=S\) without t. A repeated t in y would already have supplied an original saturated seed; hence this split has only one t. All right terms normalize to the accumulator or the two original witnesses, and the left extends the correct suffix.

Under G2 choose the original weak witness V a R x and original W a D with c(D) contained in S. In the saturated case use the second new template with \(z=S\). The extra D copies disappear when normalized against S, preserving its repeated final t. In the simple case with nonempty movement block B, use the fourth template when D omits t. If D contains t, split \(D=C\) t D0 and use the fifth template, with \(p=S\) without t. This writes the protected tail exactly once. Again the occurrence of t is unique, as otherwise an original saturated seed would exist. The two original guards are preserved, not replaced by inferred content membership. Optional empty pieces select the corresponding separate identity. Each case increases the length of the correct target suffix.

After finitely many steps \(g=E\) Q, with all outside letters in E and Q the full target last-occurrence order, with the seed tail status. Use the same five contained-word erasure laws among the 89 old laws: \(xy\leq x+zxy\) in the saturated case or when U omits t; \(xy\leq y+zxy\) for \(U=t\); \(xyz\leq xz+vxyz\) for a simple \(U=Xt\); \(xyzt\leq xty+vxyzt\) for \(U=XtY\) with both sides nonempty; \(yzt\leq ty+vyzt\) for \(U=tY\). Their substitutions and normalization are exactly the explicit contained-erasure deductions of the 13066 theory. If E is empty no erasure is needed. If Q is saturated but the target simple, apply \(x\leq xx\) with \(x=t\) in the target prefix context.

This derives every valid word absorption. All case distinctions were forced by the exact criterion, each insertion extends a finite target suffix, and all substitutions are nonempty. Splitting a polynomial identity into word absorptions proves completeness without restrictions on variables, lengths or summands. Thus 329 laws form a finite basis for 8301. Transfers to other catalogue types require independent full verification.

\par\smallskip\noindent\textit{Correction of exploratory narration.}\quad

\subsection{Finite basis for 8576: contained head and matched simple-final witnesses}\label{case:208}
Signature (+, multiplication), without constants. The basis is an effectively specified finite family, not an enumerated small list.

Addition rows 00004,\allowbreak{}01114,\allowbreak{}01224,\allowbreak{}01234,\allowbreak{}44444 form the \(chain 3<2<1<0<4\). Multiplication rows are 02230,\allowbreak{}33333,\allowbreak{}33333,\allowbreak{}33333,\allowbreak{}41134. A literal keeps its value. A nonlinear word is 3 when its head is 1,\allowbreak{}2, or 3, any factor is 3, or any nonfinal factor is 1/\allowbreak{}2. Otherwise its prefix values are all 0/\allowbreak{}4. A final 0/\allowbreak{}4 gives the original head. A final 1/\allowbreak{}2 gives 2 when the head is 0, and 1 when the head is 4. Induction on length proves this formula, including its asymmetric head dependence, for arbitrary words.

\par\smallskip\noindent\textit{Exact original-word criterion.}\quad 
For a nonlinear target w write \(C=c(w)\), \(h=head(w)\), \(t=last(w)\). Let U be the original upper words whose contents are contained in C. Then \(w\leq P\) holds exactly under both conditions:

H: U contains a word headed by h. F: if t occurs just once in w, U contains a word omitting t, or a NONLINEAR word headed by h and ending in its unique occurrence of t.

If t repeats in w there is no F condition. A literal target must itself occur in P. In particular, a contained simple-final word with the WRONG head does not satisfy F; the criterion for 6680 cannot be copied. This distinction explains the strict failed transfers.

For necessity of H put \(h=4\), other C variables 0, and outsiders 3. The target is 4 and every upper word without a contained head-h witness is at most 0. If F fails, put \(h=4\), \(t=2\), other C variables 0, and outsiders 3. Here h and t are distinct because t is simple in a nonlinear target. The target is 1. Every upper contains t or an outsider. An outsider forces 3. A contained word with nonfinal t is 3, a literal t is 2, and a nonlinear simple-final t word has a head other than h, hence value 2. The upper sum is at most 2. For a missing literal h put \(h=2\) and every other variable 3: all upper words are 3 while the target is 2.

For sufficiency target 3 is least. A target 0/\allowbreak{}4 has only 0/\allowbreak{}4 values and is matched by H. A target 1/\allowbreak{}2 has prefix values 0/\allowbreak{}4 and final 1/\allowbreak{}2, so t is simple. A contained F word omitting t is 0/\allowbreak{}4; the alternative nonlinear same-head simple-final word has exactly the target value. Thus the two conditions prove absorption for all assignments, with no restriction on variables, word lengths, or the number of summands.

\par\smallskip\noindent\textit{Structural normalization and reversible compression.}\quad 
The formula proves \(xyzt=xzyt\) and \(xxy=xy\). Normalize a nonlinear word by retaining both its original endpoints, sorting its interior, and removing duplicate interior variables and interior head copies. The normal interior is the set of original interior variables other than the head. Interior adjacent swaps use \(xyzt=xzyt\) with nonempty prefix and suffix. Delete adjacent copies with \(xxy=xy\). To delete an interior head copy, move it next to the head and apply \(xxy=xy\) with a nonempty remaining suffix. A final variable that also occurs internally keeps an interior copy unless it equals the head, in which case the two endpoints already preserve repetition. Literals remain literals.

Choose one original H witness and, when necessary, one original F witness. These at most two words already absorb w by the exact criterion. Normalize the three words. Reserve every head and final variable, at most 6. Partition the other variables by their nonempty occurrence pattern among the three words, at most 7 groups. Replace each reserved variable and each nonempty group by its own variable. Normalize again. At most 13 variables and 15 positions occur per word. Every content containment, head equality, omission, simple-final status, and nonlinear status in H/F is preserved. The compressed absorption is consequently valid.

For the inverse substitution send a reserved variable back to its original singleton and a group variable to the increasing product of its original members. All substitutions are nonempty. Every expansion recovers the original head, final and interior variable set. Literal uses of the two structural laws normalize the expansion back to the original normal word. This supplies both directions of the compression argument, including equal endpoint repetitions.

\par\smallskip\noindent\textit{Fixed finite complete basis.}\quad 
Take the two structural laws and every valid absorption over a fixed 13-variable alphabet with word length at most 15 and at most two nonempty upper summands. This family is finite and its membership is decided by H/F above. Any valid absorption reduces to its at most two original witnesses. The compressed absorption is in the family, inverse substitution and the structural laws recover that reduction, and adjoining the remaining original summands recovers the original absorption. Distributive expansion and both absorption directions give every valid polynomial identity. Hence 8576 and its opposite are finitely based without bounded-variable extrapolation. Catalogue transfers must separately be checked by complete mutual HSP.

\subsection{A finite basis for 8838 from contained original role obstructions}\label{case:209}
Take the semiring axioms without constants, \(xy=yx\), \(x^{3}y=x^{2}y\), and all valid absorptions over 92 variables, word length at most 185 and at most 3 upper words. This specifies an effective finite family. Completeness below is unrestricted in the input alphabet, lengths and polynomial size.

The multiplication rows are 11000,\allowbreak{}11111,\allowbreak{}01344,\allowbreak{}01444,\allowbreak{}01444. Write C=\{2,\allowbreak{}3,\allowbreak{}4\}. A nonlinear word is 1 if it has a letter 1 or at least two occurrences of letters with value 0;0 if it has exactly one such occurrence and all others lie in C;3 if it is quadratic with both letters 2; and 4 otherwise. Induction proves this formula. The word laws sort letters, cap multiplicity at 2, and retain a cube for singleton-content long words. No square/cube equality is used.

For target w restrict the original polynomial P to U, the words whose contents lie in c(w). Let L be literal variables of U, Q its quadratic words, QV their content union, and J its long words. Say U is bad if J is nonempty or L intersects QV. For x occurring exactly once in w, a guard is one of the following, always selected from U:

\textbullet{} The literal x.

\textbullet{} A nonlinear word in which x occurs once.

\textbullet{} A long word omitting x.

\textbullet{} A literal y and a quadratic containing y and omitting x.

The exact criterion is: a literal w must occur literally in P; for a nonlinear target, U is bad or w is quadratic with c(w) contained in QV, and every simple target variable has a guard. These are original contained witnesses; global long words or conflicts outside c(w) do not suffice.

The proper principal ideals are \{0,\allowbreak{}1\},\{1\},\{1,\allowbreak{}2\},\{1,\allowbreak{}3\}. For a literal x absent from P, assign x0 and all other variables 1. Every original nonlinear word becomes 1, and every other literal is 1, whereas \(x=0\) is not below 1. This proves the literal condition.

For a nonlinear target consider \{1,\allowbreak{}3\} first. When all target variables lie in C, external variables may be assigned 1, eliminating every original outside U. A literal in U must then be 3; a quadratic forces both its variables 2; and a long word in U cannot belong to the ideal. Thus all originals can lie in it exactly when U is not bad, with independent literal/quadratic domains. To force the target into this ideal then requires a quadratic whose variables are forced 2 by Q. This proves the first condition, including its necessity using external variables 1.

The other possible failure is a simple variable \(x=0\) with all remaining target variables in C. An original literal x fails. A nonlinear original with exactly one x fails. A long original omitting x fails. A word with at least two x has value 1 and imposes no further condition. Remaining literals force their variables 3; remaining quadratics omit x and force both endpoints 2. Those restrictions conflict precisely when an original literal y also occurs in such an omitting quadratic. These four cases are exactly the guards. If none exists, assign the forced variables 2 or 3 and the others 2, with external variables 1: every original belongs to \{1,\allowbreak{}3\}, but the target equals 0. Conversely a guard prevents all originals from belonging to the ideal under any extension.

These conditions also handle all other ideals. The all-C condition supplies a contained nonlinear original, by a long word, a conflict quadratic, or the forced quadratic witnesses. This prevents failures for \{0,\allowbreak{}1\} and for the all-C cases of \{1\} and \{1,\allowbreak{}2\}. For a single0 failure, each guard supplies either literal x or a contained nonlinear word with at most one x. That word fails membership in \{1\} and \{1,\allowbreak{}2\} when the target does. Thus all proper ideals are handled, proving sufficiency and the exact criterion for every finite polynomial.

If U is not bad and a quadratic target is valid, it already occurs in U modulo commutativity. For a square this follows from the forced quadratic occurrence on its singleton content. For distinct letters x, y, the guard for x cannot be a literal x (it would conflict with Q), a long word (none exists), or a literal/quadratic conflict (U is not bad); its only possibility is the quadratic xy itself. So this case requires no additive derivation.

Otherwise construct the saturated target-content word S, with every variable twice and a cube for singleton content. One original contained long word, or two original contained conflicting role words, proves \(S\leq P\) by the criterion. Restore each simple target variable x in turn. If r is the current protected word, replace its two x occurrences by one, obtaining r'. Use \(r'\leq r+G\), where G is its original guard of at most two words. The contained long word r makes the new all-C condition bad. Each already restored simple variable has its once-occurring witness r; the new x has G. Before each restoration r is long, with degree 2 only possible after the last restoration of a quadratic target. For a quadratic square finish with \(x^{2}\leq x^{3}\). All these local absorptions have at most 3 upper words.

Normalize each local absorption, at most 4 words. Reserve up to 3 distinct variables from each word, reserving all if fewer; this is at most 12. Group other variables by their common multiplicity vector0/\allowbreak{}1/\allowbreak{}2 across these words, at most \(3^{4}-1=80\) groups. Use the COMMON multiplicity of a group in each word, not the sum. The result uses at most 92 variables and length 185; reserved letters preserve every short/long degree. Containments, once/missing occurrences, and literal/quadratic roles and conflicts are preserved. In particular all short-word variables are reserved. The compressed absorption still satisfies the exact criterion, so belongs to the finite family. Substitute the nonempty group lists back; common multiplicity restores each original normalized word using the two word laws. This derives every local step. Transitivity and addition of unused originals finish the arbitrary absorption. Splitting a polynomial identity into summand absorptions proves completeness.

Finite audits check semantics, separators and literal proof replay; none is used in place of this unrestricted argument.

\subsection{Finite bases for 8929 and 8930 by restoring simple occurrences}\label{case:210}
In the constant-free signature, take the semiring axioms, \(xy=yx\) and \(x^{3}y=x^{2}y\), together with all valid absorptions over at most 746 variables, length at most 1493, and at most 5 upper words. The valid family is interpreted separately in each source. This is an effective finite family, not an enumerated small basis. We prove completeness for arbitrary alphabets, word lengths and numbers of original right words.

Both sources have multiplication rows 11110,\allowbreak{}11111,\allowbreak{}11111,\allowbreak{}11121,\allowbreak{}01114. A nonlinear word is 4 if all letters are 4;0 if exactly one letter is 0 and every other letter is 4;2 for a quadratic with both letters 3; and 1 otherwise. Induction on the length proves this list. The only difference between the addition tables is the sum of 2 and 3:1 in 8929,\allowbreak{}2 in 8930. The two word laws sort letters and reduce every multiplicity to at most 2, except that a long word with singleton content normalizes to the cube. Literals and quadratics retain their degrees. A cube may not be replaced by a square using a word identity.

Write c(w) for content and s(w) for variables occurring exactly once. For original P let L be its literal variables, Q the union of contents of its quadratic words, and J its long words. The short-role system is bad when J is nonempty, or, in 8929 only, L intersects Q. Then \(w\leq P\) iff:

\textbullet{} A literal target occurs literally in P.

\textbullet{} A nonlinear target has an original nonlinear word and an original word contained in c(w). For each x in s(w), it also has an original contained word in which x occurs at most once, including zero times.

\textbullet{} For a nonlinear target the short-role system is bad, or w is quadratic and c(w) is contained in Q.

Here is a proof using proper principal ideals. For \{0,\allowbreak{}1,\allowbreak{}2,\allowbreak{}3\}, failure of the target means all its letters are 4. An original contained-content word prevents every original from lying in that ideal; if none exists, assign target letters 4 and all other letters 0 to separate.

For A=\{1,\allowbreak{}2,\allowbreak{}3\}, a nonlinear target fails either with all letters 4, already handled, or with a simple variable x assigned 0 and all other target variables 4. Set all external variables 1; this maximizes satisfaction of original A-memberships. An original word fails precisely when it is contained in target content and has at most one occurrence of x. This is exactly the stated simple-variable witness. It proves necessity and sufficiency, since that witness also fails for every extension of the same target assignment. For a literal target x, set x0 and all other variables 1. Every nonlinear original lies in A. Thus x must occur literally; this condition suffices for every ideal.

The ideal \{3\} contains only literal word values, so a nonlinear target requires a nonlinear original. Assigning 3 everywhere separates otherwise. The remaining proper ideal is \{2\} in 8929 and \{2,\allowbreak{}3\} in 8930. Long words never lie in it; a quadratic lies there exactly when both letters are 3. Literal domains are \{2\} in 8929 and \{2,\allowbreak{}3\} in 8930, while quadratic domains are \{3\} in both. Therefore the system is impossible exactly in the bad case defined above. If possible, a long target fails, while a quadratic is forced into the ideal exactly when both variables occur in original quadratics. These ideals exhaust the addition order and prove the exact criterion without any bound on input size.

First suppose the system is possible. A nontrivial nonlinear target is quadratic. Retain one contained original, at most two simple-occurrence witnesses, and at most two quadratic occurrence witnesses. This is at most 5 originals, still satisfying the exact criterion. The literal case is immediate.

If the system is bad, let S be the saturated target-content word: every target variable occurs twice, with a cube used for singleton content. Retain one contained original and one original long word if available; otherwise retain one literal and one quadratic sharing a variable. At most 3 original words prove \(S\leq P\) by the criterion. For each simple target variable x in turn, choose its original contained witness u, and replace the two occurrences of x in the current protected word r by one. Call the result r'. The absorption \(r'\leq r+u\) is valid: r contains all target content, its long degree makes the short-role system bad, and its already restored simple variables retain r as their witness; the newly simple x uses u. Before any such step r is long: a two-letter word could only arise after the final restoration for a quadratic target. After all restorations we reach the target normal form, except that a quadratic square needs the valid final absorption \(x^{2}\leq x^{3}\). This too is covered by the same bounded family. No original witness is silently replaced by a changed polynomial summand.

It remains to bound each of the above local absorptions, with at most 5 upper words. Normalize all at most 6 words. In each word reserve any three distinct variables, or all its variables if fewer than three. Thus at most 18 variables are reserved. Every remaining variable occurs only in long words, with multiplicity 0,\allowbreak{}1 or 2. Partition such variables by their COMMON multiplicity vector across all six words: at most \(3^{6}-1=728\) nonzero vectors. Represent each group by one fresh variable, with its COMMON multiplicity in each word, NOT the sum of multiplicities. Together there are at most 746 variables; normalized length is at most

1493. The reserved variables preserve literal/quadratic/long degree.

The compressed absorption satisfies the exact criterion: content containments, missing/simple occurrences, original nonlinear or long witnesses, and quadratic roles are preserved. All variables in short words are reserved, so the literal/quadratic contradiction cannot be created or destroyed by compression. The compressed absorption is a member of the finite valid family. Inversely substitute the nonempty list of group members for each fresh variable. Common multiplicity ensures that every member has exactly its original normalized count; commutativity and \(x^{3}y=x^{2}y\) derive precisely the original words. This proves every local step from the finite family. Transitivity and addition of the original unselected words finish the saturation/restoration chain. Splitting a polynomial identity into individual summand absorptions in both directions proves unrestricted completeness for both sources.

\subsection{Simple-variable witnesses and a finite basis for entry 8841}\label{case:211}
Work in the signature (+, multiplication) without constants, relative to the axioms of additively idempotent semirings. Numerical labels name elements only. The algebra A has catalogue number 8841; write w \ensuremath{\leq} P for w+P approximately equal to P.

\par\smallskip\noindent\textit{Word values and normalization.}\quad 
Addition joins the chains \(1<0<4\) and \(3<2<4\). The tables in numerical element order are

{\small\par\smallskip
\noindent addition       multiplication\par
\CaseLine{\ttfamily 0 0 4 4 4\qquad 1 1 0 0 0}
\CaseLine{\ttfamily 0 1 4 4 4\qquad 1 1 1 1 1}
\CaseLine{\ttfamily 4 4 2 2 4\qquad 0 1 4 2 4}
\CaseLine{\ttfamily 4 4 2 3 4\qquad 0 1 4 3 4}
\CaseLine{\ttfamily 4 4 4 4 4\qquad 0 1 4 4 4}
\smallskip}

Put K=\{2,\allowbreak{}3,\allowbreak{}4\}. A nonlinear word is 1 if it has a factor 1 or at least two occurrences assigned 0; it is 0 if it has exactly one occurrence assigned 0 and all remaining occurrences lie in K. On all-K assignments its value is 3 if every occurrence is 3, is 2 if the first occurrence is 2 and every later occurrence is 3, and is 4 otherwise. Linear words keep their assigned values.

These rules agree under both associations of any triple: a factor 1 or at least two zeros gives 1, a unique zero with other factors in K gives 0, and all-K triples follow the stated initial-position rule. Every row and column preserves the joins inside both chains, the four joins between \{0,\allowbreak{}1\} and \{2,\allowbreak{}3\}, and joins involving 4, as the tables show. Thus A is an additively idempotent semiring, and induction proves the word formula for arbitrary length.

Write c(u), h(u), t(u) for a word's content, initial variable, and set of variables after the initial position, respectively. Put t(y) equal to the empty set for a linear word y. Let \(m_{x}\)(u) be the number of occurrences of x, and call x simple in u if \(m_{x}(u)=1\).

The identities

\CaseDisplay{xyz\approx xzy,\qquad x^3\approx x^2       \quad\textup{(\ensuremath{W})}}

hold. They sort the tail arbitrarily and truncate each total variable multiplicity at two. Repeated occurrences of the head can be brought next to its initial occurrence before reducing cubes; other repeats are grouped within the tail. Hence a normal word retains its literal head and each total multiplicity in \{0,\allowbreak{}1,\allowbreak{}2\}. Linear words stay linear. Given finite nonempty C, h in C and S contained in C, let F(C, h, S) be the normal word with head h, multiplicity one for variables of S and multiplicity two for variables of C minus S. This notation allows a linear word when \(C=S=\{h\}\).

\par\smallskip\noindent\textit{Exact criterion.}\quad 
Let P be a nonempty sum of words, put \(D=c(P)\), \(T=\bigcup_{u\in P}t(u)\), and B\_C=\{u in P:c(u) subset C\}, including linear summands.

Lemma 1. A linear y satisfies y \ensuremath{\leq} P exactly when y is a literal summand of P. A nonlinear w, of content C, satisfies w \ensuremath{\leq} P exactly when

1. C is contained in D;

2. t(w) is contained in T;

3. \(B_{C}\) is nonempty;

4. for every simple x in w, some u in \(B_{C}\) has \(m_{x}\)(u) at most one.

Proof. For a missing literal y, assign \(y=0\) and all other variables

1. Every right word is 1, including nonlinear powers of y, so y is not absorbed. Literal inclusion gives the converse.

If 1 fails, assign a missing variable 4 and all others 3: \(w=4\) and \(P=3\). If 2 fails at x, assign \(x=2\) and all others 3: \(w=4\), while every right word lies in \{2,\allowbreak{}3\}, because none has x in its tail. If 3 fails, assign C the value 3 and all outsiders 1: \(w=3\) and \(P=1\). If 4 fails at a simple x, assign \(x=0\), other C-variables 3 and all outsiders 1. Then \(w=0\). Every contained word has at least two x's, and every other word has a factor 1, so \(P=1\).

For sufficiency, if \(P=3\) all variables in D are 3 and \(w=3\) by 1. If \(P=2\), variables in D lie in \{2,\allowbreak{}3\} and those in T are 3; 1 and 2 give w in \{2,\allowbreak{}3\}. If P lies in \{1,\allowbreak{}0\}, w cannot have all its variables in K, because a witness from 3 would take a value in K. Thus w lies in \{1,\allowbreak{}0\}. For \(P=1\), the only possible violation is \(w=0\). Its unique zero must be assigned to a simple target variable x, with all other C-variables in K. A witness from 4 either omits x, giving a value in K, or contains it exactly once, giving 0. Both contradict \(P=1\). Finally \(P=4\) is automatic. These cases exhaust the proper principal additive ideals and prove sufficiency. \(\square\)

\par\smallskip\noindent\textit{Full guards and recovery of individual simple variables.}\quad 
Fix a valid nonlinear w \ensuremath{\leq} P, with \(C=c(w)\), \(h=h(w)\), and simple-variable set S. Choose an original v in \(B_{C}\). For each x in T choose an original nonlinear \(a_{x}\) with x in t(\(a_{x}\)). Every local step below has at most two right words, each original or already derived below P. Its validity is checked directly by Lemma 1, independently of any basis claim.

If \(\ct(q)\subseteq T\), the step \(q^2\law q+a_{h(q)}\) is valid. The target has no simple variables and tail c(q); q covers every tail variable except possibly its head, and \(a_{h(q)}\) supplies that head. This includes linear q. We may therefore build a full guard with every variable doubled. Missing x in T can be added using \(qx^{2}\) \ensuremath{\leq} q+\(a_{x}\). Repetition gives a full guard on any specified content contained in T.

Suppose first that C is contained in T. Starting at v, build such a guard of content C and then use the valid one-right-word step F(C, h, empty) \ensuremath{\leq} q to align its head to h. Denote this guard by G. For each x in S choose \(b_{x}\) in \(B_{C}\) with \(m_{x}\)(\(b_{x}\)) at most one. The step

\CaseDisplay{F(C, h,\{x\})\leq G+b_x                 \quad\textup{(\ensuremath{1})}}

is valid: G covers content, the entire tail set and a contained witness; \(b_{x}\) witnesses the only simple variable x.

If C is not contained in T, condition 2 forces its only variable outside T to be the simple head h. Put \(E=C\) minus \{h\}; then E is nonempty and contained in T. Every original nonlinear word containing h has h exactly once, as its head. If v contains h, it is either the literal h or has simple head h. The valid step F(c(v), h,\{h\}) \ensuremath{\leq} v doubles all its other variables. Append each missing E-variable twice, using \(qx^{2}\) \ensuremath{\leq} q+\(a_{x}\). This yields G=F(C, h,\{h\}).

If v omits h, build a full guard q on E. An original summand r contains h by condition 1. The step F(C, h,\{h\}) \ensuremath{\leq} q+r is valid: q supplies all target tail variables and a contained witness omitting the sole simple variable h; r supplies the missing content h. This again yields G=F(C, h,\{h\}). For each x in S minus \{h\}, use

\CaseDisplay{F(C, h,\{h, x\})\leq G+b_x.                \quad\textup{(\ensuremath{2})}}

The guard G witnesses the simple h, and \(b_{x}\) witnesses the simple x. All remaining conditions come from G. Thus 2 is valid.

In the first case start \(J=empty\); in the second start \(J=\{h\}\). We know F(C, h, J) \ensuremath{\leq} P. For each remaining x in S, combine that word with the word already obtained in 1 or 2, whose simple set is \(R=\{x\}\) or R=\{h, x\}, respectively. The local step

\CaseDisplay{F(C, h, J\cup\{x\})
               \leq F(C, h, J)+F(C, h, R)                      \quad\textup{(\ensuremath{3})}}

is valid. The first right word witnesses all variables of J; the second witnesses x. Both have content C. If the target head is repeated, both right tails cover C. If it is simple, the required tail is C minus \{h\}, which both right tails cover. Lemma 1 applies. After finitely many such steps \(J=S\) and the target normal word is exactly w. This proves unrestricted derivability by local absorptions with at most two right terms.

\par\smallskip\noindent\textit{Compression by multiplicity patterns.}\quad 
There are at most three words in a local step. Reserve their literal head variables individually, at most three. Group all remaining variables by the vectors of their truncated total multiplicities in these words. There are at most \(3^{3}-1=26\) nonzero vectors. In each compressed word use a group symbol with the common multiplicity of its members in that word, and retain the reserved head as its head. Normalize the tail. This gives at most 29 variables and length 58.

This operation is not the substitution identifying all group members: such identification would wrongly change a group of simple variables to a repeated variable. Its validity instead follows from Lemma 1. Within a group the presence pattern and multiplicities agree in every word, and none of its members is a head. Content coverage, tail coverage, content containment of each right witness, and the existence of a count-at-most-one witness for each simple variable are therefore preserved in both directions. For a simple target group, every member has exactly the same possible right witnesses. Reserved heads are individual, so their tail membership is preserved as well. The compressed local absorption is valid.

For the inverse substitution, send each group symbol to the nonempty ordered product of its original members and each reserved head to itself. A group used once or twice restores each member once or twice. Tail permutations and cube reductions from (W) recover exactly each original normal word, including every linear word and repeated head. Empty substitutions and unit constants are unnecessary.

Theorem. The identities (W), together with all valid normal absorptions over a fixed 29-letter alphabet, with word length at most 58 and at most two right summands, form a finite basis for A. Word reversal gives a finite basis for its opposite.

Proof. The displayed family is finite and its membership is decided by Lemma 1. Compression and inverse substitution derive every local step from this family and (W). The guard and recovery construction derives every valid word absorption without bounds on variables, lengths or numbers of summands. Distribute arbitrary terms, derive each summand's absorption into the other side in both directions, and obtain every valid identity. Reversal transports this proof to the opposite. \(\square\)

This specifies a proved finite family; it does not enumerate a small or minimal basis and makes no hereditary assertion. Finite checks are auxiliary. Other catalogue types must pass complete two-way HSP checks before this theorem can be transferred to them.

\subsection{A finite basis from original tail-clause witnesses for 10250}\label{case:212}
Without constants, take the semiring axioms, \(xyz=xzy\) and \(xyyz=xyz\), and all valid absorptions over 2069 variables, word length at most 2071, and at most 10 upper words. This is an effective specified finite family, not an enumerated small basis. We prove unrestricted completeness, including arbitrary alphabets and polynomial sizes.

Multiplication rows 22220,\allowbreak{}23220,\allowbreak{}22222,\allowbreak{}22222,\allowbreak{}22224 show:a nonlinear word is 4 if every letter is 4; it is 0 if its head is 0 or 1 and every noninitial letter is 4; it is 3 for the quadratic(1,\allowbreak{}1); and otherwise2. Induction proves this for all lengths. A long word retains its head and its noninitial variable set, INCLUDING any occurrence of the head in its tail. The two laws permute the tail and remove its duplicate letters while the word length exceeds 3. The canonical tail is the sorted noninitial set, duplicated if it has only one variable, so the word always remains long. Literals and quadratic words are unchanged. In particular \(xxyz=xyz\) is false and is not used here.

For original P let L be its literal variables, Q its quadratic words, and J its long words. For each variable intersect its initial full domain with\{3\} if it is in L, and with\{1\} if it occurs in Q. Call the domain system impossible if J is nonempty or a domain is empty. For words u, w define an original tail witness by

head(u) in content(w), and either u is literal or noninitial-content(u) subset noninitial-content(w).

The exact absorption criterion is:

\textbullet{} Some original word has content contained in content(w).

\textbullet{} If w is literal, it itself occurs in P.

\textbullet{} If w is nonlinear, P has a nonlinear original and an original tail witness as defined above. Also either the domain system is impossible, or w is quadratic and both variable domains are contained in\{1\}. A long target needs the impossible-domain case.

The proper principal ideals are\{0,\allowbreak{}1,\allowbreak{}2,\allowbreak{}3\},\{1\},\{2,\allowbreak{}3\},\{3\}. For the first, a word is outside it exactly when all letters are 4. This gives the contained-content requirement, with the explicit separating assignment target variables 4, others 0 when no original is contained. The ideal\{1\} has only literal word values, so requires an original nonlinear for a nonlinear target. All variables 1 separate if no such original exists.

For ideal A=\{2,\allowbreak{}3\}, a nonlinear word lies in A exactly when its head belongs to A OR some noninitial letter is not 4. A literal lies in A exactly when its variable belongs to A. These are positive clauses in the predicates \(a(x)=[x\) in A] and \(n(x)=[x\) not 4], with a(x) implying n(x). To falsify a nonlinear target, all its noninitial variables must be 4 and its head must be outside A. Set those tail variables 4; set its head 0 if it is not already a tail variable; set every other variable 2. This is a maximal assignment for satisfying original positive clauses while falsifying the target. An original fails precisely when it is a literal whose variable belongs to target content, or its head lies in target content and all its tail variables lie in the target tail. That is exactly the original tail-witness condition. This proves both directions. For a literal target, set its variable 0 and all others 2. Every original nonlinear remains in A; only that original literal can fail. Therefore a literal target must occur literally in P. Head-tail coincidence is handled by assigning every target tail variable 4 first; it must not be replaced by a condition requiring a simple head.

The ideal\{3\} contains only literal 3 and quadratic(1,\allowbreak{}1). Consequently all originals lie there iff the domain system is possible, with independent domain choices. Long targets cannot belong to it. A quadratic target belongs to it for every such assignment precisely when both domains are contained in\{1\}; a literal target already present is automatic. The four proper principal ideals detect all possible failures, proving the criterion for arbitrary variables, lengths, and polynomial sizes.

For any valid absorption retain one contained original. For a literal target retain that literal. For a nonlinear target retain one nonlinear original and one tail witness. For the domain condition retain a long original if available, otherwise at most 3 originals producing an empty domain, or at most 3 role witnesses for each of the at most 2 variables of a quadratic target. Thus at most 9 originals suffice;10 is safe. Their selected polynomial still satisfies the full criterion.

Normalize the selected words and target, at most 11 words. Reserve every head and every variable appearing in a literal or quadratic, at most 22. Group other variables by their nonempty incidence subsets among the 11 words, at most 2047 groups. They occur only in long tails. Identify each group, normalize, and obtain a valid absorption within the stated finite family. Inversely substitute the nonempty list of group members. Because members have COMMON incidence and no head or short occurrence, the expansions retain every original head, tail set, and long degree. The two laws normalize exactly back to the original selected words. This derives the arbitrary input absorption. Splitting each valid polynomial identity into summand absorptions in both directions proves complete finite basability. Independent finite computations check implementation only; equal-variety transfers require complete mutual HSP separately.

\subsection{Tail coverage and a finite basis for a two-chain semiring}\label{case:213}
Let A be catalogue entry 8859. Work in the signature (+, multiplication), without constants, relative to the additively idempotent semiring axioms. Numerical labels are element names only. Write w \ensuremath{\leq} P for w+P approximately equal to P.

\par\smallskip\noindent\textit{Operations, word values and normalization.}\quad 
Addition joins the two chains \(1<0<4\) and \(3<2<4\), with no further comparabilities. In numerical element order the tables are

{\small\par\smallskip
\noindent addition       multiplication\par
\CaseLine{\ttfamily 0 0 4 4 4\qquad 1 1 0 0 0}
\CaseLine{\ttfamily 0 1 4 4 4\qquad 1 1 1 1 1}
\CaseLine{\ttfamily 4 4 2 2 4\qquad 1 1 4 2 4}
\CaseLine{\ttfamily 4 4 2 3 4\qquad 1 1 4 3 4}
\CaseLine{\ttfamily 4 4 4 4 4\qquad 1 1 4 4 4}
\smallskip}

Put K=\{2,\allowbreak{}3,\allowbreak{}4\}. The element 1 is multiplicatively absorbing. A nonlinear word has the following values:

\textbullet{} it is 0 when its first occurrence is assigned 0 and every later occurrence is assigned an element of K;

\textbullet{} it is 1 when there is an occurrence assigned 1, or an occurrence assigned 0 outside the first position;

\textbullet{} on all-K assignments it is 3 if every occurrence is 3, is 2 if its first occurrence is 2 and all later occurrences are 3, and is 4 otherwise.

In the first bullet the 0 occurrence is necessarily unique. These cases exhaust all assignments; linear words keep their assigned values. The table proves associativity by considering triples. An assigned 1 forces 1. Among triples with a 0, only a unique initial 0 followed by K-elements gives 0; every other case gives 1, under either association. The all-K rule gives the same three outcomes under both associations. Every row and column preserves the joins within each chain and each of the four joins between \{0,\allowbreak{}1\} and \{2,\allowbreak{}3\}; hence multiplication preserves every binary join. Thus A is an additively idempotent semiring. Induction proves the word formula for every length.

The identities

\CaseDisplay{xyz\approx xzy,\qquad xy^2\approx xy        \quad\textup{(\ensuremath{W})}}

follow from the formula. They permute every position after the first and remove repetitions there. A nonlinear normal word is its head followed by the sorted set of variables occurring after the first position. If its head occurs in this set, it stays in both positions. Linear words are kept linear. Write c(w), h(w), and t(w) for content, first variable, and the set after the first position, respectively; put \(t(y)=empty\) for a linear word y.

\par\smallskip\noindent\textit{An exact criterion.}\quad 
Write P as a nonempty sum with linear variable set L and nonlinear summands N. Let \(D=c(P)\), \(T=\bigcup_{u\in N}t(u)\), and N\_C=\{u in N:c(u) subset C\}.

Lemma 1. A linear variable y is absorbed by P exactly when y belongs to L. A nonlinear w, with \(C=c(w)\) and \(h=h(w)\), satisfies w \ensuremath{\leq} P in A exactly when

1. C is contained in D;

2. t(w) is contained in T;

3. L meets C, or \(N_{C}\) is nonempty;

4. if h does not belong to t(w), then either L meets C or some u in \(N_{C}\) omits h or has head h with no later occurrence of h.

Proof. A missing literal y is separated by \(y=0\) and all other variables 1. Every right word is then 1, including every nonlinear power of y.

If 1 fails, give a missing target variable value 4 and every other variable 3. Then \(w=4\) and \(P=3\). If 2 fails at x, assign \(x=2\) and every other variable 3. Because x occurs in no right tail, every right summand is 2 or 3, but \(w=4\). If 3 fails, assign C the value 3 and outsiders 1; then \(w=3\) and \(P=1\). Finally, if 4 fails, assign \(h=0\), other C-variables 3 and outsiders 1. The head of w is simple, so \(w=0\). All linear right words are 1. Every contained nonlinear right word contains h, but either outside its first position or in a repeated head position; its value is 1. All other nonlinear words also have a factor 1, giving \(P=1\). This proves necessity.

For sufficiency it is enough to examine the proper principal additive ideals \{1\}, \{1,\allowbreak{}0\}, \{3\}, \{3,\allowbreak{}2\}. If \(P=3\), all variables of D must be 3; by 1, \(w=3\). If P is 2 or 3, every variable of D is 2 or 3, and every variable of T is 3. Conditions 1 and 2 therefore give w in \{2,\allowbreak{}3\}, which lies below 2.

If P lies in \{1,\allowbreak{}0\}, w cannot have all its variables in K: the linear or contained nonlinear witness in 3 would itself lie in K. Thus w lies in \{1,\allowbreak{}0\}. The remaining issue is \(P=1\) and \(w=0\). Here h is a simple head valued 0, and all other C-variables are in K. A linear witness in C has value different from 1. Otherwise the witness in 4 has value in K when it omits h and value 0 when it has h as a simple head. Both contradict \(P=1\). Hence \(w=1\) when \(P=1\). All four proper principal ideals have been checked; \(P=4\) is automatic. \(\square\)

\par\smallskip\noindent\textit{Derivation by steps with at most two right summands.}\quad 
Call a nonlinear normal word saturated when its head occurs in its tail. It then has tail equal to its whole content. We construct every valid nonlinear absorption using (W) and valid local absorptions with at most two right words. Each right word is original or already derived below P. Linear targets are literal inclusions.

Fix valid w \ensuremath{\leq} P, and choose an original summand v with content contained in C, using condition 3 (either linear or nonlinear). Its head is g. For any variable x in T, choose an original nonlinear summand \(a_{x}\) with x in t(\(a_{x}\)).

Two operations will be used. If an available q has head g in T, append g to it if necessary to obtain a saturated word q'. The step q' \ensuremath{\leq} q+\(a_{g}\) is valid by Lemma 1: q supplies content coverage and condition 3, and its tail together with that of \(a_{g}\) covers the new tail. The head of q' is repeated, so condition 4 is empty. This includes \(q=g\) linear. Second, if x in T is absent from q, append x. The step qx \ensuremath{\leq} q+\(a_{x}\) is valid. When q is saturated, the new word is saturated too. When q has simple head h, q itself supplies condition 4, either as the literal h or as a nonlinear same-head witness. The new tail is covered in both cases.

First suppose C is contained in T. Saturate v as above. Append each missing variable of C by the second operation, obtaining a saturated q of full content C. If w has repeated head, Lemma 1 gives w \ensuremath{\leq} q. If w has simple head h, use the original witness b in condition 4: either a linear variable in C or a contained word omitting h or having simple head h. Lemma 1 gives w \ensuremath{\leq} q+b. The guard q supplies all content and tail requirements, and b supplies the only remaining condition. At most two right words occur at every step.

Otherwise condition 2 implies that h is simple, h is not in T, and every variable of C except h is in T. If v contains h, it must be the literal h or have simple head h: no original nonlinear right tail contains h. Append every missing variable of C except h. The resulting normal word has head h and tail C minus \{h\}, hence is w. Every step is the second operation above.

If v omits h, its head g lies in T. Saturate v and append all missing variables of C minus \{h\}, giving a saturated q on precisely C minus \{h\}. This set is nonempty because w is nonlinear with a simple head. Choose an original r containing h by condition 1. The final step w \ensuremath{\leq} q+r is valid: q supplies the contained witness, all target tail variables, and condition 4 by omitting h; r supplies the missing content occurrence h. This completes the construction.

\par\smallskip\noindent\textit{Finite compression and completeness.}\quad 
In a local absorption there are at most three words. Reserve their head variables individually, at most three (including every linear word's unique variable). Group all other variables by their nonempty tail-incidence patterns across the words, at most seven groups. Identify each group with one fresh symbol and normalize by (W). This is a substitution, so it preserves validity. It uses at most ten variables and words of length at most eleven.

Conversely substitute for each group symbol the nonempty product of its original members, leaving reserved heads as individual variables. Uniform tail incidence restores exactly each original tail set. No head was grouped; (W) sorts and reduces the tails, recovering each original normal word. Compression is therefore reversible modulo (W), even for repeated heads and linear words.

Theorem. The two word identities (W), together with all valid normal absorptions over a fixed ten-letter alphabet, of word length at most eleven and with at most two right words, form a finite basis for A. Reversing every word gives a basis for its opposite.

Proof. This is a finite, explicitly defined set, and Lemma 1 decides membership. Compression derives every local step from the finite set. The preceding construction derives every valid word absorption with unrestricted variables, lengths and numbers of summands. Distributivity writes arbitrary terms as nonempty sums of words. Derive each summand's absorption into the other polynomial in both directions to obtain every valid identity. Word reversal transports validity and derivations to the opposite. \(\square\)

The finite family is not enumerated as a small list; no minimality or hereditary assertion is made. Finite computation is an auxiliary error check. Classification of other catalogue types requires an independent complete check of both HSP directions.

\subsection{Prefix coverage and recovery of a simple head}\label{case:214}
Let A be catalogue entry 8860. Work in the signature (+, multiplication), without constants, relative to the axioms for additively idempotent semirings. Write w \ensuremath{\leq} P for w+P approximately equal to P. Numerical labels name elements, not term constants.

\par\smallskip\noindent\textit{Operations and word values.}\quad 
Addition joins the two chains \(1<0<4\) and \(3<2<4\). The tables, with rows and columns in numerical order, are

{\small\par\smallskip
\noindent addition       multiplication\par
\CaseLine{\ttfamily 0 0 4 4 4\qquad 1 1 0 0 0}
\CaseLine{\ttfamily 0 1 4 4 4\qquad 1 1 1 1 1}
\CaseLine{\ttfamily 4 4 2 2 4\qquad 1 1 4 4 4}
\CaseLine{\ttfamily 4 4 2 3 4\qquad 1 1 2 3 4}
\CaseLine{\ttfamily 4 4 4 4 4\qquad 1 1 4 4 4}
\smallskip}

Put K=\{2,\allowbreak{}3,\allowbreak{}4\}. A nonlinear word has value 1 if any occurrence is assigned 1 or any occurrence after the first is assigned 0. Its value is 0 if its initial occurrence is assigned 0 and every later occurrence is in K. On all-K assignments its value is 3 when every occurrence is 3, is 2 when the final occurrence is 2 and every earlier occurrence is 3, and is 4 otherwise. Thus 0 requires a simple head, whereas 2 requires a simple final variable. Linear words keep their assigned values.

To verify associativity, a factor 1 forces 1. Among triples involving 0, only a unique initial 0 followed by K-elements gives 0; all other cases give 1 under either association. On K, a triple is 3 if all three factors are 3, is 2 if its last factor is 2 and the first two are 3, and is 4 otherwise. Each row and column preserves the joins within both additive chains and the four joins between \{0,\allowbreak{}1\} and \{2,\allowbreak{}3\}; joins involving 4 are also preserved by the displayed rows and columns. Hence A is an additively idempotent semiring. Induction gives the word formula for every length.

For a word u, write c(u) for its content, h(u) for its initial variable, f(u) for its final variable, and p(u) for the set of variables before the final position. Put \(p(y)=empty\) for a linear word y. Linear words remain linear in every normal form. The identities

\CaseDisplay{xyzt\approx xzyt,\qquad xy^2z\approx xyz       \quad\textup{(\ensuremath{W})}}

hold by the formula. They permute all interior positions and remove interior repetitions. A normal word keeps its literal first and last positions and lists its distinct interior variables in a fixed order. Words of length two have no interior substring. If an endpoint also occurs internally, that interior occurrence is retained. Products over empty interior sets mean omitted substrings, not unit constants. Every word is transformed to this normal form using (W): interchange adjacent interior letters, group repetitions, then remove adjacent interior duplicates with nonempty left and right contexts.

\par\smallskip\noindent\textit{Exact absorption criterion.}\quad 
Let P be a nonempty sum of words, with linear variable set L and nonlinear summands N. Put \(D=c(P)\), \(U=\bigcup_{u\in N}p(u)\), and N\_C=\{u in N:c(u) subset C\}.

Lemma 1. A linear y is absorbed by P exactly when y belongs to L. A nonlinear w, with \(C=c(w)\) and \(h=h(w)\), satisfies w \ensuremath{\leq} P exactly when

1. C is contained in D;

2. p(w) is contained in U;

3. L meets C, or \(N_{C}\) is nonempty;

4. if h does not occur after the first position of w, then either L meets C or some u in \(N_{C}\) omits h or has h as its simple head.

Proof. If a literal y is absent, assign \(y=0\) and all other variables

1. Every right word is 1, including every nonlinear power of y, so y is not absorbed. Literal inclusion proves the converse.

For necessity of 1, assign a missing target variable 4 and all others 3; then \(w=4\) and \(P=3\). If 2 fails at x, assign \(x=2\) and all others 3. Then \(w=4\) while every right summand is 2 or 3, because x occurs in no right prefix. If 3 fails, assign variables in C the value 3 and all outsiders 1; then \(w=3\) and \(P=1\). If 4 fails, assign \(h=0\), the other C-variables 3 and outsiders 1. The simple head makes \(w=0\). There is no linear summand in C. Every contained nonlinear word has h outside its first position, and every other nonlinear word has a factor 1. Thus \(P=1\). These evaluations prove necessity.

For sufficiency consider the proper principal additive ideals \{3\}, \{3,\allowbreak{}2\}, \{1,\allowbreak{}0\} and \{1\}. If \(P=3\), every right variable is 3, and 1 makes \(w=3\). If \(P=2\), every variable in D lies in \{2,\allowbreak{}3\} and every variable in U is 3. By 1 and 2, w lies in \{2,\allowbreak{}3\}, below 2.

If P lies in \{1,\allowbreak{}0\}, w cannot have all its variables in K: the linear or contained nonlinear witness from 3 would then take a value in K. So w lies in \{1,\allowbreak{}0\}. It remains to rule out \(P=1\) and \(w=0\). Here h is a simple head assigned 0 and all other C-variables are in K. A linear witness in C differs from 1. Otherwise 4 supplies a contained word omitting h, with value in K, or having simple head h, with value 0. Either contradicts \(P=1\). Finally \(P=4\) absorbs everything. \(\square\)

\par\smallskip\noindent\textit{Construction using at most two right words per step.}\quad 
A guard will have repeated head, removing condition 4. A full guard also has repeated final variable, so its prefix set is its entire content. Fix w \ensuremath{\leq} P satisfying Lemma 1, choose an original summand v with content contained in C by 3, and for each x in U choose an original nonlinear \(a_{x}\) with x in p(\(a_{x}\)).

If a known q has all its content in U, form \(q^{2}\). The local step

\CaseDisplay{q^2\leq q+a_{f(q)}                     \quad\textup{(\ensuremath{1})}}

is valid by Lemma 1. Its target has repeated head, content c(q) and prefix c(q). The right q supplies all but possibly the last prefix variable, and \(a_{f(q)}\) supplies that variable. The argument also covers linear q. If q is already a full guard and x in U is missing from it, the local step

\CaseDisplay{qx^2\leq q+a_x                        \quad\textup{(\ensuremath{2})}}

is valid: its new prefix requirement is precisely x, and q supplies a contained witness. Both steps retain repeated head. Repeating 2 grows a full guard to any specified content contained in U.

First suppose C is contained in U. Apply 1 to v, then 2 for each missing C-variable. This yields a known full guard q of content C. If w has repeated head, w \ensuremath{\leq} q is valid. If its head h is simple, choose the original linear or nonlinear witness b from 4. The step w \ensuremath{\leq} q+b is valid: q supplies content, prefix coverage and a contained witness, while b supplies the simple-head requirement.

Now suppose C is not contained in U. Since p(w) is contained in U, the only target variable outside U is its simple final variable s. Put \(E=C\) minus \{s\}; E is nonempty and contained in U. Every original nonlinear word containing s has s only as its final variable.

If v contains s and is nonlinear, duplicate its initial letter: g v \ensuremath{\leq} v, where \(g=h(v)\). This is valid by Lemma 1; the target has repeated head, the same content and prefix as v, and final variable s. If v is linear, then \(v=s\). Choose g in E and start instead with

\CaseDisplay{g^2s\leq s+a_g.                       \quad\textup{(\ensuremath{3})}}

The target prefix is \{g\}; the literal s supplies a contained witness. In either case we have a guard q ending in s. Insert each missing x in E immediately before s. Writing \(q=rs\), the local step

\CaseDisplay{rxs\leq q+a_x                        \quad\textup{(\ensuremath{4})}}

is valid; the new prefix requirement is x, q remains a contained witness, and the head remains repeated. Continue to content C.

If v omits s, use 1 and 2 to build a full guard q on E. By 1 of Lemma 1, an original summand b contains s. The local step qs \ensuremath{\leq} q+b is valid: q supplies a contained witness and all target prefix variables, b supplies the missing content variable, and the head remains repeated. This again yields a guard with content C and prefix E. Finish with w \ensuremath{\leq} q for repeated head, or w \ensuremath{\leq} q+b' using the original simple-head witness b' from 4. Thus every valid absorption is built using valid local steps with at most two right words. Each right word was original or already proved below P.

\par\smallskip\noindent\textit{Reversible compression and a finite basis.}\quad 
In a local step there are at most three words. Reserve every literal endpoint of these words individually, at most six variables. Group the remaining variables by their nonempty incidence patterns in the three interiors, at most seven groups. Identify each group with one fresh symbol and normalize by (W). This substitution preserves validity and yields at most thirteen variables, normal word length at most fifteen, and at most two right words.

For the inverse substitution, leave every reserved endpoint individual and replace each group symbol by the nonempty ordered product of its members. Every group occurs only internally and its members have identical interior incidence across all words. Sorting and removing interior repetitions by (W) recovers each original normal word. This includes endpoint occurrences retained in interiors, words with equal endpoints, words of length two, and linear words. No empty substitutions or unit constants are used.

Theorem. The identities (W), together with all valid normal absorptions over a fixed thirteen-letter alphabet, with word length at most fifteen and at most two right summands, form a finite basis for A. Reversal of every word gives a finite basis for its opposite.

Proof. The specified family is finite and Lemma 1 decides its membership. Reversible compression derives every local step in the construction from this family and (W). The construction therefore derives all valid word absorptions with unrestricted variables, lengths and numbers of summands. Distribute arbitrary terms into nonempty sums of words, and derive the absorption of each summand in both directions to obtain every valid identity. Reversal carries all operations on words and all derivations to the opposite. \(\square\)

This is a proved finite family, not an enumerated small basis. It makes no minimality or hereditary claim. Finite computations provide auxiliary checks; completeness follows from the unrestricted proof. Other catalogue types require independent complete transfers in both HSP directions before receiving the classification.

\subsection{A finite basis for 8932 from chain ideals and original tail witnesses}\label{case:215}
Work without constants. In addition to the semiring axioms, take \(xyz=xzy\), \(xyyz=xyz\), and every table-valid absorption over at most 2069 variables, with word length at most 2071 and at most 10 upper words. This specifies an effective finite family; it is not an enumerated small basis. The proof below has no bound on the input alphabet or polynomial.

The addition order is the chain \(3<2<1<0<4\). The multiplication rows are 11110,\allowbreak{}11111,\allowbreak{}11111,\allowbreak{}11121,\allowbreak{}11114. Direct induction gives the value of any nonlinear word:4 if all letters have value 4;0 if its head has value 0 and every noninitial letter has value 4;2 for the quadratic word whose two letters both have value 3; and 1 otherwise. These cases are disjoint. A long word retains its head and noninitial variable set, including any noninitial occurrence of its head. The two displayed laws permute that tail and remove repetitions while retaining length at least 3. Normalize a long word by sorting the tail set and repeating its sole member if the set has size1. Literals and quadratic words stay intact. In particular \(xxyz=xyz\) fails and is not used.

Let P consist of original nonempty words. Write Q for its quadratic words and J for its words of length at least 3. An original word u is a tail witness for w when head(u) is in content(w), and either u is literal or noninitial-content(u) is contained in noninitial-content(w). The exact criterion for \(w\leq P\) is:

\textbullet{} Some original word has content contained in content(w).

\textbullet{} A literal target must itself occur literally in P.

\textbullet{} A nonlinear target requires an original nonlinear word and an original tail witness. In addition, either J is nonempty, or w is quadratic and each of its variables occurs in an original quadratic word.

To prove necessity and sufficiency, inspect the four proper principal ideals of the addition chain. For \{0,\allowbreak{}1,\allowbreak{}2,\allowbreak{}3\}, the word value is outside exactly when all its letters are 4. Hence the contained-content condition is necessary and sufficient. If absent, assign target letters 4 and all other letters 0 to separate.

For A=\{1,\allowbreak{}2,\allowbreak{}3\}, a nonlinear word belongs to A exactly when its head is in A or some noninitial letter is not 4. To falsify a nonlinear target, set its tail variables 4, set its head 0 if not already in that tail, and all other variables 1. This maximizes both positive predicates [x in A] and [x not 4] subject to target failure. An original word fails exactly if it is a tail witness as defined above. Thus a witness is both necessary and sufficient. For a literal target x, assigning x0 and every other variable 1 leaves every nonlinear original in A; only the literal x can fail. This proves the literal requirement. The argument treats a repeated head by setting tail variables 4 first.

For \{2,\allowbreak{}3\}, all originals can lie there only if J is empty. Each literal then restricts its variable to \{2,\allowbreak{}3\}; each quadratic restricts both variables to \{3\}. These independent domains are always simultaneously consistent, by assigning 3 everywhere. A long target never lies in this ideal. A quadratic lies there under every admissible assignment exactly when both its variables are forced 3 by original quadratics. An original literal target is automatic. Unlike 10250, literal/quadratic overlap causes no contradiction in these domains.

Finally \{3\} contains only literal word values. It contributes precisely the requirement for an original nonlinear word when the target is nonlinear. If all originals are literal, assigning 3 to every variable separates any nonlinear target. These ideals detect every failure of the absorption, establishing the criterion for arbitrary polynomials.

Retain one contained-content original. For a literal target retain that literal. Otherwise retain a nonlinear original and a tail witness, then one long original if available, or one quadratic occurrence witness for each of the at most two target variables. At most 5 originals suffice, and the selected polynomial still satisfies the exact criterion.

Normalize the selected words and target. With the safe bound of 10 upper words there are at most 11 words. Reserve every head and every variable occurring in a literal or quadratic word, at most 22 variables. Group the other variables by their common nonempty incidence subset among these words, at most 2047 groups. They occur only in long tails. Replace each group by a fresh variable and normalize. Identifying variables preserves a valid identity, so the result is a member of the specified finite family; it uses at most 2069 variables and length 2071.

Inversely substitute the nonempty list of members of each group. Common incidence, the reserved heads and short occurrences, and the retained long degree ensure that every expansion has the original head and tail set. The two word laws normalize it exactly to the original selected word. Thus the bounded valid absorption derives the input absorption. Adjoining unselected original summands gives the result for P. Every polynomial identity splits into these absorptions in both directions. This proves unrestricted completeness and finite basability. Numerical verification supplements this argument; it does not replace it.

\subsection{Contained witnesses for two five-element semirings}\label{case:216}
We work in the signature (+, multiplication) without constants, relative to the additively idempotent semiring axioms. Write w \ensuremath{\leq} P for w+P approximately equal to P. Let A and B be catalogue entries 8906 and 8911, respectively. Numerical labels name elements only. Their tables, in numerical element order, are

{\small\par\smallskip
\noindent addition       A multiplication   B multiplication\par
\CaseLine{\ttfamily 0 0 0 4 4\qquad 1 1 1 0 0\qquad 1 1 1 0 0}
\CaseLine{\ttfamily 0 1 1 4 4\qquad 1 1 1 1 1\qquad 1 1 1 1 1}
\CaseLine{\ttfamily 0 1 2 3 4\qquad 2 2 2 2 2\qquad 2 2 2 2 2}
\CaseLine{\ttfamily 4 4 3 3 4\qquad 0 1 1 4 4\qquad 1 1 1 4 4}
\CaseLine{\ttfamily 4 4 4 4 4\qquad 0 1 1 4 4\qquad 1 1 1 4 4}
\smallskip}

The additive order \(is 2<1<0<4\) and \(2<3<4\). The only incomparable pairs are 1,\allowbreak{}3 and 0,\allowbreak{}3; their joins are 4. Put K=\{3,\allowbreak{}4\} and Q=\{0,\allowbreak{}3,\allowbreak{}4\}.

\par\smallskip\noindent\textit{Semiring axioms, word values and normal forms.}\quad 
Every row and column of both multiplication tables is increasing and preserves the two incomparable joins. Hence multiplication preserves all binary joins. For associativity in A, a first factor 2 or 1 fixes the product at 2 or 1 under either association. If the first factor belongs to Q, an occurrence of 1 or 2 later makes the product 1; the same is true under either association. With all three factors in Q, two or more zeros give 1, one zero gives 0, and no zero gives 4, under both associations. This also proves associativity in the cases where an intermediate product leaves Q.

For B the constant first rows2 and 1 give the same two initial cases. For a first factor in Q, any later factor outside K gives 1; if all later factors are in K, a head 0 stays 0 and a head in K gives 4. Splitting a triple according to these cases gives the same result under both associations. Thus both are additively idempotent semirings. Induction gives the following rules at arbitrary length.

For any nonlinear word the head 2 gives 2 and the head 1 gives 1. In A, for a head in Q an occurrence of 1 or 2 later gives 1. Otherwise all factors are in Q: at least two occurrences assigned 0 give 1, exactly one gives 0, and none gives 4. In B, a head in Q gives 1 unless every tail factor belongs to K; if they all do, it gives 0 for head 0 and 4 for a head in K. Linear words retain their assigned value. In particular a nonlinear word never has value 3.

For A the identities

\CaseDisplay{xyz\approx xzy,\qquad x^3\approx x^2       \quad\textup{(\ensuremath{WA})}}

hold. They permute every tail position and truncate total variable multiplicities at two. A normal word keeps its literal head and each total multiplicity in\{0,\allowbreak{}1,\allowbreak{}2\}. The head can be brought next to its repetitions; other repetitions are grouped in the tail before cube reductions. For B the identities

\CaseDisplay{xyz\approx xzy,\qquad xy^2\approx xy       \quad\textup{(\ensuremath{WB})}}

hold. Their normal word keeps its head followed by the sorted set of tail variables, retaining a head that also occurs in the tail. Linear words are kept linear in both normalizations. No empty-word substitutions or unit constants are involved.

Write c(u), h(u), t(u) for a word's content, head and tail-variable set, with \(t(y)=empty\) for a linear word. Write \(m_{x}\)(u) for the total number of x occurrences, and call x simple when \(m_{x}(u)=1\).

\par\smallskip\noindent\textit{Exact criteria.}\quad 
Let P be a nonempty sum of words and \(C=c(w)\) for a nonlinear target. Define B\_C=\{u in P:c(u) subset C\}, including linear words. The common head condition is

\CaseDisplay{(H)\quad\text{an original nonlinear }a\in P
                           \text{ has }h(a)=h(w).}

Lemma 1. A linear y is absorbed by P in either algebra exactly when it is a literal summand. For nonlinear w:

\textbullet{} A satisfies w \ensuremath{\leq} P exactly when(H), \(B_{C}\) is nonempty, and for every simple target variable x some b in \(B_{C}\) has \(m_{x}(b)\leq 1\).

\textbullet{} B satisfies w \ensuremath{\leq} P exactly when(H) and some original b in P is either a linear variable in C, or a nonlinear word with c(b) contained in C and t(b) contained in t(w).

Proof. A missing literal y is separated by assigning \(y=0\) and all other variables 2. Every nonlinear word has value at most 1: a head other than y is 2, and a head y is followed by a factor 0 or 2, giving 1. All linear right words are 2. Thus \(P\leq 1<0=y\).

For necessity of(H), assign \(h(w)=3\) and all other variables 2. Without a nonlinear right head h(w), all right nonlinear words are 2 and linear words are 2 or 3. Hence P lies in\{2,\allowbreak{}3\}. The target is 4 if all its variables equal h(w), and is 1 otherwise, since its tail then has a2. Neither value lies below 3, proving(H).

In A, if \(B_{C}\) is empty, assign C the value 3 and all outsiders 2. The target is 4. Every right nonlinear word is 2 when its head is outside C and 1 when an outside variable occurs in its tail; right linears are 2. Thus \(P\leq 1\). For a simple x lacking its count-at-most-one witness, assign \(x=0\), other C-variables 3 and outsiders 2. The target is 0. Every contained right word has at least two x occurrences, giving 1, and every other word is at most 1 as above. This proves necessity of all A conditions.

For necessity of B's second condition, put \(T=t(w)\). Assign T the value 3; if h(w) is simple assign it 0; assign outsiders 2. The target is 4 for repeated head and 0 for simple head. There is no linear right variable in C. Every right nonlinear word with head outside C is 2. With head in C, absence of an eligible witness means a tail variable lies outside T, hence is 0 or 2, so its value is 1. Thus \(P\leq 1\), separating the target.

For sufficiency use the proper principal ideals \{2\}, \{2,\allowbreak{}1\}, \{2,\allowbreak{}1,\allowbreak{}0\}, \{2,\allowbreak{}3\}. If P is 2 or 3, its nonlinear(H)-witness must have value 2, since nonlinear words cannot be 3. Its head is therefore 2 and so is w. This treats both such right values.

If \(P=0\), only \(w=4\) could violate absorption. Then all C-variables are in K. In A any word from \(B_{C}\) takes value 3 or 4, contradicting \(P=0\). In B the stated second witness likewise has value 3 or 4.

If \(P=1\), the possible violating target values are 0 and 4. The 4 case is excluded by the same contained witness. In A a target value 0 has precisely one occurrence assigned 0: a simple variable x, with all other C-variables in K. Its witness with \(m_{x}(b)\leq 1\) has value 0 when it contains x and value 3 or 4 when it omits x. In B a target value 0 has a simple head 0 and all tail variables in K. A linear witness in C is 0 or in K. An eligible nonlinear witness has head in C and tail in t(w), hence value 0 or 4. Both give the same contradiction to \(P=1\). Finally \(P=4\) is automatic. These cases prove sufficiency without restrictions on variables. \(\square\)

\par\smallskip\noindent\textit{Two right terms suffice for B.}\quad 
For a valid nonlinear absorption in B, choose its two original witnesses a, b from Lemma 1. The absorption w \ensuremath{\leq} a+b is itself valid, because a supplies(H) and b supplies the second condition. It is enough to derive this one absorption with at most two right words; \(a+b\leq P\) follows by literal inclusion. Duplicate witnesses can be removed. Linear targets need only literal inclusion.

In any such local absorption there are at most three words. Reserve their heads individually, at most three variables. Group all other variables by nonempty incidence patterns across their three tails, at most seven groups. Identify each group with one symbol and normalize by(WB). This substitution preserves validity and uses at most ten variables and normal word length eleven.

For the inverse substitution send each group symbol to the nonempty ordered product of its members and each reserved head to itself. Common tail incidence and(WB) recover each original normal word, including linear words and repeated heads. Hence all valid local absorptions derive from those within the ten-variable bound.

\par\smallskip\noindent\textit{Simple-variable recoveries for A.}\quad 
Fix valid nonlinear w \ensuremath{\leq} P in A, let \(h=h(w)\), \(C=c(w)\), and let S be its simple-variable set. For h in C and J contained in C, let F(C, h, J) be the normal word with head h, total multiplicity one for members of J and two for all other C-variables.

Choose the original nonlinear head witness a and any original contained witness v. Lemma 1 makes the following local absorption valid:

\CaseDisplay{G=F(C, h,\varnothing)\leq a+v.        \quad\textup{(\ensuremath{2})}}

The target has no simple variables; a supplies(H), v supplies the contained witness. For each x in S choose an original \(b_{x}\) in \(B_{C}\) with \(m_{x}(b_{x})\leq 1\). The local absorption

\CaseDisplay{R_x=F(C, h,\{x\})\leq G+b_x         \quad\textup{(\ensuremath{3})}}

is valid: G supplies(H) and a contained word; \(b_{x}\) witnesses the only simple variable. All these normal words are nonlinear because w is nonlinear: when C is a singleton its simple set is empty.

Start \(J=empty\) and \(q=G\). To adjoin a missing simple x use

\CaseDisplay{F(C, h, J\cup\{x\})\leq q+R_x.             \quad\textup{(\ensuremath{4})}}

Both right words are contained and have head h. For any simple variable in J the right q witnesses count one; \(R_{x}\) witnesses x. Thus Lemma 1 proves 4. Continue until \(J=S\), when the resulting word is exactly w. Each step has at most two right words, original or already derived below P; the number of simple variables is unrestricted.

For compression of a local step reserve its at most three heads individually. Group remaining variables by nonzero vectors of their total multiplicities in\{0,\allowbreak{}1,\allowbreak{}2\} across the at most three words. There are at most 26 groups, giving at most 29 variables and length 58. In each compressed word give a group symbol the common multiplicity of its members, rather than their summed multiplicity.

This is not ordinary variable identification. Lemma 1 preserves compressed validity: head equalities stay literal, content containment of every possible witness is equivalent, and all group members have the same count-at-most-one witnesses for simplicity. Conversely replace each group symbol by the nonempty product of its members. Multiplicity one or two restores every member once or twice; tail sorting and(WA) recover the original normal words exactly. Reserved heads and linear words stay individual. Thus every local step is derivable from its compressed version and(WA).

Theorem. A has the finite basis consisting of(WA) and all valid normal absorptions over a fixed 29-letter alphabet, of length at most 58 and with at most two right terms. B has the finite basis consisting of(WB) and all valid normal absorptions over a fixed ten-letter alphabet, of length at most 11 and with at most two right terms. Each multiplicative opposite has the reversed basis.

Proof. Both families are finite and membership is decided by Lemma 1. The witness construction, reversibility and stepwise recovery derive every valid word absorption with arbitrary variables and summands. Distributivity expands terms to nonempty polynomials; derive each summand's absorption into the other polynomial in both directions to obtain every valid identity. Reversal transports all derivations to the opposite. \(\square\)

These are explicitly specified finite families, not enumerated small or minimal lists. No hereditary claim is made. Finite computations are auxiliary error checks. Transfers to further catalogue types require independent complete HSP verification in both directions.

\subsection{A finite basis for 8910 with literal head and nonlinear tail guards}\label{case:217}
Take semiring axioms without constants, \(xyz=xzy\), \(xyy=xy\), and all valid absorptions over 16 variables, length 18, at most 2 upper words. This finite family is specified, not enumerated as a small list.

The proper principal ideals of addition are I0=\{0,\allowbreak{}1,\allowbreak{}2\}, I1=\{1,\allowbreak{}2\}, \(I2=\{2\}\), I3=\{1,\allowbreak{}2,\allowbreak{}3\}. Multiplication has left zeros 1,\allowbreak{}2. For a nonlinear word: if its head is 2, the value is 2; otherwise if its head is 1 or some noninitial letter is 0,\allowbreak{}1, or 2, the value is 1; otherwise the value is 0 when its head is 0, and 4 when its head is 3 or 4. Induction from the operation table proves this formula for all lengths. The two word laws normalize a nonlinear word to its head h followed by the distinct variables in its ACTUAL noninitial set T, in a fixed order. A repeated head remains in T. Literals remain literal.

For original P and target w, let U be the original words contained in c(w), N its nonlinear words, and H the heads of all original words of P. The exact criterion is:

\textbullet{} a literal x requires its original literal;

\textbullet{} a nonlinear w requires h(w) in H, N nonempty, and either an original u in N with T(u) subset T(w), or the original literal h(w).

We prove this using all proper principal ideals. I2 requires every upper head to have value 2, so forces the target exactly when its head is in H. A failure when that head is missing assigns it 1 and all other variables 2. For I0, a failing nonlinear target has every target letter in\{3,\allowbreak{}4\}. Any contained nonlinear word then has value 4, and so excludes failure. For I1, a failing target has every tail letter in\{3,\allowbreak{}4\} and head in\{0,\allowbreak{}3,\allowbreak{}4\}. If the head is 3 or 4, the same nonlinear witness excludes failure. If the head is 0, the original head literal excludes failure, or a nonlinear contained word with tail in T(w) has value 0 or 4 and excludes failure. For I3 the only nonlinear failures similarly have tail letters in\{3,\allowbreak{}4\} and head in\{0,\allowbreak{}3,\allowbreak{}4\}; the same two branches work: head 0 is excluded by its literal or the nonlinear tail guard, and head 3/\allowbreak{}4 by the contained nonlinear witness. This proves sufficiency.

For necessity of N, assign all target variables 3 and all external variables 1. Nonlinear target value 4 cannot be absorbed by contained literals of value 3 and all outside words of value 1. Suppose N is nonempty, but no nonlinear tail guard and no original head literal exist. Then the target head must be simple:if it occurred in its own tail, \(T(w)=c(w)\) and every u in N would be a guard. Every u in N must contain that head NONINITIALLY. Assign the head 0, other target variables 3, and all outside variables 1. Every nonlinear contained upper is 1, all contained literals are 3, and all outside words are 1. The target is 0, which is not below their sum 1 or 3. This proves the missing-guard obstruction, including arbitrary outside words. For a literal x lacking its original literal, assign \(x=0\) and others 1. Every original upper is 1, including powers of x of degree at least 2; the target 0 is not absorbed. These separators prove necessity.

The global head witness can have external variables. In the nonlinear tail-guard branch keep that head witness and the contained nonlinear guard. In the literal-head branch the literal itself already supplies the global head; keep it and one nonlinear contained word. Thus at most two original upper words suffice, and the criterion remains true after selection. Literals need just themselves.

Normalize the target and selected words. Reserve all heads, all literal variables, and each nonlinear head plus up to two other distinct letters. At most 9 variables are reserved among at most 3 words. Group the others by common incidence, using at most 7 nonempty groups. Identification preserves validity and gives at most 16 variables and normal length 17, hence within the stated length 18 family. Substitute every group variable by the NONEMPTY word listing its members. All heads are reserved, so their tail occurrence is preserved. Every affected word remains nonlinear; common incidence recovers its actual tail set. The two word laws recover all original normal forms. Add back the unselected original summands. Distributive expansion and absorption in both directions give unrestricted finite basability.

Do not replace the nonlinear tail guard by an arbitrary contained literal. A head literal works only together with a contained nonlinear word. Short-word agreement and one-way HSP do not justify transfers.

\subsection{Six complete identities for semiring 10377}\label{case:218}
Constant-free signature (+, multiplication), relative to the additively idempotent semiring axioms. The basis consists of

{\small\par\smallskip
\CaseLine{\(xy=yx; xyz=ttt; ttt\leq x+xy;\)}
\CaseLine{\(xy\leq x+ttt; xy\leq xx+ttt; xy\leq xx+yz.\)}
\smallskip}

We write \(u\leq P\) for \(u+P=P\). The algebra considered here is \(\A{10377}\).

\par\smallskip\noindent\textit{Source 10377: intermediate cubic value.}\quad 
Addition:00004/\allowbreak{}01004/\allowbreak{}00224/\allowbreak{}00234/\allowbreak{}44444. Multiplication: 22222/\allowbreak{}23222/\allowbreak{}22222/\allowbreak{}22222/\allowbreak{}22220. The additive order \(is 3<2<0<4\) and \(1<0<4\). A quadratic is 3 at(1,\allowbreak{}1),0 at(4,\allowbreak{}4), and 2 otherwise. All long words equal 2. This cubic value dominates neither arbitrary literals nor the pair value 0.

For original P, let L be literals, E unordered quadratic edges, V edge endpoints, S loops. Call the cubic value available when P has an original long word or L intersects V. A literal target is absorbed iff original. A long target is absorbed iff the cubic value is available. A pair ht is absorbed iff its edge is original or an endpoint belongs to S union L, and additionally the cubic value is available or both endpoints lie in V.

For sufficiency, an original long word supplies 2. An original incident literal x and edge xy also have join at least 2:if \(xy=3\) then \(x=1\) and \(1+3=0\); otherwise xy itself dominates 2. With a literal guard h and available 2, a target pair can exceed2 only if \(h=t=4\), in which case h absorbs it. With a loop guard hh and available 2, the same exception is absorbed by \(hh=0\). Without available 2, a pair guarded by an original loop has both endpoints in V. If P fails to dominate 2, all right quadratics must equal 3 and all right literals 3, forcing every V variable 1; the target is 3. For its value 0, the original guard loop is 0 and absorbs it. Original edges are automatic.

For necessity, a missing literal h is separated by \(h=4\) and all others 3; right terms are at most 0, below 4. If the cube is unavailable, assign V to 1 and outsiders 3. Every right term is 3, so a long target 2 fails. For a pair lacking both an original edge and endpoint in S union L, set the target endpoints 4 and all others 3. There is no right literal or quadratic wholly on those endpoints. All right terms are at most 2, whereas the target is 0. If the pair has a guard but lacks cubic availability and an endpoint lies outside V, assign V to 1 and outsiders 3:all right terms are 3 and the target 2. These assignments exhaust the failed conditions for arbitrary polynomials.

All six laws are valid by this criterion and are independently checked. An original long word normalizes to a fixed cube. An original incident literal and edge derive it by \(ttt\leq x+xy\). For a target pair with available cube, use \(xy\leq x+ttt\) if its guard is literal, or \(xy\leq xx+ttt\) if its guard is a loop. Otherwise a new pair has an original loop and an original edge incident with the other endpoint; use \(xy\leq xx+yz\). No inference introduces a new literal. All substitutions are nonempty, and every right-side term is original or explicitly derived below P. No bounds on variables or term sizes occur.

\par\smallskip\noindent\textit{Completeness and transfer.}\quad 
Distributivity reduces arbitrary terms to polynomials. The preceding exact absorption criterion, explicit separators, and complete derivations prove every valid identity by summandwise absorption in both directions. This establishes the six fixed laws without bounds on variables, word lengths, or summand counts. Independent checking tests all finite law assignments and literally replays the derivations. Other catalogue types are classified by this proof only after full HSP in both directions.

\subsection{A finite 77-law basis for 10071 using prefix sets and original witness templates}\label{case:219}
Signature (+, multiplication), without constants. All bases are relative to additively idempotent semiring axioms. There is no minimality claim; the templates below specify at most 77 identities.

\par\smallskip\noindent\textit{1. Word semantics and an exact criterion.}\quad 
The addition and multiplication tables, with elements ordered 0--4, are:

{\small\par\smallskip
\noindent addition       multiplication\par
\CaseLine{\ttfamily 0 0 0 4 4\qquad 1 2 2 3 4}
\CaseLine{\ttfamily 0 1 1 4 4\qquad 2 2 2 3 4}
\CaseLine{\ttfamily 0 1 2 4 4\qquad 2 2 2 3 4}
\CaseLine{\ttfamily 4 4 4 3 4\qquad 4 4 4 4 4}
\CaseLine{\ttfamily 4 4 4 4 4\qquad 4 4 4 4 4}
\smallskip}

Thus \(2<1<0<4\) is a chain, and 3 is incomparable with 0,\allowbreak{}1,\allowbreak{}2, with join 4. For a nonlinear word w, any factor 4 or any internal factor 3 gives value 4. Otherwise final 3 gives 3. With all factors in \{0,\allowbreak{}1,\allowbreak{}2\}, the value is 1 exactly for a quadratic word with both factors 0, and is 2 otherwise. Literals retain their values. The formula follows by induction from the displayed multiplication.

For any nonempty polynomial P, let A be the union of all contents of its original words with their last positions removed, B the set of all their last variables, and \(C_{P}\) their full content union. Let I(w)=c(w without its last position) and t(w) its last variable. A literal target is absorbed exactly when it occurs literally in P. For nonlinear w, the exact criterion for \(w\leq P\) is:

(C) c(w) is contained in \(C_{P}\); (T) either A intersects B, or I(w) is contained in A and t(w) lies in B; (Q) if w is quadratic, some original right word has length at most two and content contained in c(w).

Necessity of C: assign a missing target variable 4, all other variables

2. The target is 4 and P is 2. Suppose A and B are disjoint. For x in I(w) outside A, assign B union \{x\} the value 3 and all others

2. Every original right word is 3, while the target is 4. If the last variable is outside B, assign B the value 3 and all others 2: P is 3 and the target is 2 or 4. Both separate a failure of T. For a quadratic target assign its variables 0 and all outsiders 2. The target is 1; without a Q-witness every right word is 2. A missing literal x is separated by \(x=0\), all other variables 2, since nonlinear right words have value at most 1 and other literals are 2.

For sufficiency consider the value of P. Value 4 absorbs everything. If P lies in \{0,\allowbreak{}1,\allowbreak{}2\}, no right word can have value 3 or 4. Consequently no variable in \(C_{P}\) is valued 3 or 4; C ensures the same for w. A long target is 2. A quadratic target is also 2 except when all its variables are 0, in which case its Q-witness has value 0 or 1 and supplies the needed bound. Finally, if \(P=3\), all right words equal 3. Thus all B variables are 3, all A variables avoid 3,\allowbreak{}4, and A intersects B is impossible. Condition T forces the target's internal variables to avoid 3,\allowbreak{}4 and its tail to be 3, hence its value is 3. This proves sufficiency for all numbers of variables and right summands.

\par\smallskip\noindent\textit{2. The finite family and long normal forms.}\quad 
Use these nine fixed identities:

{\small\par\smallskip
\CaseLine{\((N1) xyz = yxz\)}
\CaseLine{\((N2) x^{2}yz = xyz\)}
\CaseLine{\((N3) x^{3}y = x^{2}y\)}
\CaseLine{\((M)  xyz \leq  xz+yz\)}
\CaseLine{\((Qx) xy \leq  x+x^{2}y\)}
\CaseLine{\((Qy) xy \leq  y+x^{2}y\)}
\CaseLine{\((Qxx) xy \leq  x^{2}+x^{2}y\)}
\CaseLine{\((Qyy) xy \leq  y^{2}+x^{2}y\)}
\CaseLine{\((Qyx) xy \leq  yx+x^{2}y\)}
\smallskip}

In addition use the following finite templates. Brackets mark a whole optional block: either include that variable or omit it when WRITING the identity. Every variable present in an instantiated identity must subsequently receive a NONEMPTY word or term; brackets never license empty-word substitutions.

Four compatible-witness templates (a and c independently optional):

{\small\par\smallskip
\CaseLine{\(x^{2}t \leq  [a] x b + [c] t.\qquad (K)\)}
\smallskip}

Here b is always present. Sixty-four crossing-witness templates (a, c, d, e, f, g independently optional; b always present):

{\small\par\smallskip
\CaseLine{\(x^{2}t \leq  [a] q b + [c] q + [d] x [e] + [f] t [g].\qquad (X)\)}
\smallskip}

The fixed laws and all 4+64 bracket choices give 77 identities before removing any duplicates. Template variables may coincide under ordinary substitution. Soundness is immediate from C, T, Q: in K the first right word places x internally and the second ends in t; in X the first two words place q both internally and finally and the other two supply x, t. M has all needed prefix variables and the same tail. Each Q law has a short contained witness, while \(x^{2}y\) supplies its prefix/tail condition. N1--N3 preserve long prefix sets and the last variable in the word formula.

For a word of length at least three, N1 permutes every pair of adjacent positions before the last. If there are at least two distinct internal variables, bring a repeated variable to the start as x x, keeping another internal block y and the final variable z; N2 deletes a copy. Continue until each internal variable occurs once, then sort them. If there is only one internal variable x, N3 reduces its prefix multiplicity to two. This yields the canonical word: the sorted internal-variable set followed by the literal tail when the set has at least two members; otherwise \(x^{2}t\). Words of lengths one and two remain unchanged. Every stated rewrite uses nonempty substitutions and preserves length at least three.

\par\smallskip\noindent\textit{3. Deriving every long absorption.}\quad 
Suppose \(w\leq P\) is valid and |w|\ensuremath{\geq}2. Put \(I=I(w)\), \(t=t(w)\); I is nonempty. For every x in I we derive the long elementary guard \(x^{2}t\leq P\).

If A and B are disjoint, T supplies an original word with an internal x, written a x b with b nonempty, and an original word ending in t, written c t. Either prefix a or c may be empty. Select exactly the appropriate K bracket choice, and substitute the displayed nonempty blocks. Its right side is made of original summands, giving \(x^{2}t\leq P\).

If A and B intersect, fix q in their intersection. Choose an original word a q b with b nonempty and an original word c q. By C choose original words d x e and f t g containing x and t. Choose the matching X brackets and substitute every existing block, obtaining \(x^{2}t\leq P\). The witnesses may coincide and blocks may contain any additional variables; additive idempotence removes duplicated right summands if needed.

Starting with one elementary guard, use M repeatedly with substitutions \((x,y,z)=(U,x^2,t)\) to get \(Ux^{2}t\leq Ut+x^{2}t\leq P\). Thus all internal variables are combined into a single guard with prefix set I and tail t. Its length is at least three, and long normalization removes excess copies. For |w|\ensuremath{\geq}3 the guard equals the normal form of w, completing the derivation. This construction also produces \(x^{2}y\) as a guard for any quadratic xy, including \(x=y\). It does not assume an empty word or an extra identity.

\par\smallskip\noindent\textit{4. Quadratic and literal targets.}\quad 
If \(w=xy\) with \(x\ne y\), Q provides a short contained original witness from \{x, y, \(x^{2}\), \(y^{2}\), xy, yx\}. If it is xy, the target already occurs. Otherwise one of Qx, Qy, Qxx, Qyy, Qyx combines that witness with the derived guard \(x^{2}y\) to give \(xy\leq P\). For \(w=x^{2}\) the witness is x or \(x^{2}\). The latter is already the target; the former uses Qx with \(y=x\), namely \(x^{2}\leq x+x^{3}\). It would be incorrect to use \(x^{2}\leq x\) here: it fails at \(x=3\). Literal targets occur literally and need no construction.

All valid word absorptions have now been derived with no restriction on variables, lengths or summand count. Distributive expansion and summandwise absorption in both directions give every valid identity. This proves finite basability by the stated 77-law family. Reversing all words gives a finite basis for the multiplicative opposite; other catalogue targets require independent complete two-way HSP certificates.

\subsection{Two six-identity bases for Boolean group extensions}\label{case:220}
The source manuscript for this case acknowledges support from the National Natural Science Foundation of China (12371024), the Science and Technology Research Program of Chongqing Municipal Education Commission (KJZD-K2024011102) and the Chongqing Natural Science Foundation Innovation and Development Joint Fund (Municipal Education Commission) (CSTB2025NSCQ-LZX0067).

\par\smallskip\noindent\textit{Abstract.}\quad 
We give six-identity bases for two extensions of a nontrivial Boolean group. In the first construction, a nilpotent chain lies above the flat group part. In the second, two group elements have additional elements above them, with one exceptional square. Long words are controlled by content and the affine hull of their parity vectors. Linear and quadratic words have explicit additional witnesses. We derive every valid absorption for arbitrarily many variables. The second generated variety is a proper subvariety of the first, and both are finitely based. The five-element instances give representatives 9212 and 9952 in the catalog used in this project.

\par\smallskip\noindent\textit{1. Constructions and the bases.}\quad 
We use the signature (+, *) without constants. All bases are relative to the axioms for additively idempotent semirings. Write u \ensuremath{\leq} v for the identity u + v = v. The names of algebra elements below are not constant symbols in terms.

Let G be a nontrivial Boolean group, with identity e; thus every element has square e. Put A(G) = G union \{z, a, b\}, where the three new elements are distinct. Addition restricts on G union \{z\} to the flat addition: \(g+g=g\), and \(g+h=z\) for distinct g, h in G. Order all group elements below \(z<a<b\). Multiplication restricts on G to its group operation; put \(b^{2}=a\), and put every other product involving one of z, a, b equal to z.

For the second construction choose g0 in G different from e, and put B(G, g0) = G union \{z, d, c\}. Addition is the join in the order whose only nontrivial comparabilities below its greatest element z are \(e<d\) and \(g0<c\). In particular different group elements have sum z. Extend the map pi by \(\pi (g)=g\), \(\pi (d)=e\), and \(\pi (c)=g0\). Multiplication has absorbing element z, has \(c^{2}=d\), and otherwise is pi(x)pi(y).

Both are additively idempotent semirings. For A, products involving the new elements are as in the three-element nilpotent chain, and every mixed group/new-element product is z. Associativity follows by considering triples inside G or meeting the new elements; any triple in the latter case is z. Multiplication preserves joins, including the two boundaries at z and at \(a<b\).

For B, every product of at least three non-z elements equals the product of their pi-images in G. Checking the case with two or three occurrences of c verifies associativity across the sole exception \(c^{2}=d\); all other cases follow from the group law. Multiplication preserves the comparable pairs \(e<d\) and \(g0<c\). For joins of incomparable elements it gives z: multiplication by a group element permutes the group images; multiplication by d has image pi; for multiplication by c, the exceptional value d joins another product to z unless that other factor is g0, the comparable case already checked. Thus distributivity holds.

The four common identities are

{\small\par\smallskip
\CaseLine{\((C)  xy = yx,\)}
\CaseLine{\((N)  x^{3}yz = xyz,\)}
\CaseLine{\((T)  xyz \leq  x + y + z,\)}
\CaseLine{\((E)  x^{3} \leq  x^{3} y^{2}.\)}
\smallskip}

The basis for A(G) adds

{\small\par\smallskip
\CaseLine{\((A1) xy \leq  x + (xy)^{3},\)}
\CaseLine{\((Q)  xy \leq  x^{2} + (xy)^{3}.\)}
\smallskip}

The basis for B(G, g0) instead adds

{\small\par\smallskip
\CaseLine{\((B1) xy \leq  x^{3} + (xy)^{3},\)}
\CaseLine{\((Q)  xy \leq  x^{2} + (xy)^{3}.\)}
\smallskip}

Theorem 1. These are six-identity bases for A(G) and B(G, g0), respectively. Their identities do not depend on the nontrivial Boolean group G or, in the second construction, on the choice of g0.

It suffices to prove completeness for \(G=C2\). Indeed each construction contains the corresponding C2 construction as a subsemiring. All six displayed identities hold for every allowed G, as proved below. Thus completeness for the five-element member implies equality of its identities with those of every member of its family.

\par\smallskip\noindent\textit{2. Word values and normalization.}\quad 
Temporarily take G=\{e, g\}. In \(A=A(C2)\), a product of length at least three is its group product if every factor lies in G, and is z otherwise. A product of length two has the same rule except that \(bb=a\). For B=B(C2, g), a product of length at least three is z if some factor is z, and otherwise is the product of the pi-images. A product of length two has the same rule except that \(cc=d\).

For a nonempty word w write C(w) for its content. Let epsilon(w) be its parity vector in \(F2^{X}\), where X is a finite variable set containing all variables under consideration. Its x coordinate records the number of occurrences of x modulo 2. Denote the unit vector at x by \(epsilon_{x}\). Call a word long when it has length at least three.

Lemma 2. Identities (C) and (N) identify any two long words with the same content and parity vector.

Proof. Suppose a long word contains x. Using (C), write it as xYZ, where both Y and Z are nonempty words. Identity (N) adds two copies of x. Its reverse removes two copies whenever the resulting word still has length at least three. Thus every long word can be transformed into a sorted word having exponent 3 at each variable of odd parity and exponent 4 at each variable of even parity. Throughout the transformation every present variable has positive exponent, and the word stays long. The resulting word depends only on its content and parity vector. This proves the lemma. In particular \(w=w^{3}\) is derivable whenever w is long. \(\square\)

All six identities are sound in the indicated family. Commutativity is immediate. For (N), both sides are long and have the same content and parity. In (T), if the right sum is not above the flat group part, all its group images agree, and their product is that same element. Otherwise the long product is below the right sum. More explicitly, in A a factor outside G makes the long product z and the right sum at least z. In B a right sum different from z lies in one additive branch with a common group image, so the left product is that group element.

For (E), \(x^{3}\) lies in the flat group part. If y lies outside the group in A, the right product is z; otherwise \(y^{2}\) acts as the group identity. In B, a non-z value of \(y^{2}\) has image e and has the same action on a group element, while \(y=z\) gives z.

For (Q), the only discrepancy between xy and (xy)\ensuremath{^{3}} is \(xy=a\) in A or \(xy=d\) in B. These occur respectively at \(x=y=b\) or \(x=y=c\), and in both cases \(x^{2}\) is precisely the exceptional product. At \(x=y=b\), (A1) is also true because \(a\leq b\). At \(x=y=c\), the right side of (B1) is \(g+e=z\), which is above d. In all other cases the cubic term already equals xy. This proves soundness without a restriction on the number of variables in substitutions.

\par\smallskip\noindent\textit{3. The common affine calculation for long words.}\quad 
Let P be a nonempty finite sum of words. Define C(P) to be the union of the contents of its summands, and set

{\small\par\smallskip
\noindent H(P) = \(aff_{F2}\) \{ epsilon(u) : u is a summand of P \}.\par
\smallskip}

Lemma 3. In both A and B, a long word w satisfies \(w\leq P\) exactly when

C(w) is contained in C(P), and epsilon(w) belongs to H(P).

Moreover every such absorption is derivable from (C), (N), (T), (E).

Proof of necessity. If x belongs to C(w) but not C(P), give x the multiplicative zero z and every other variable the group identity. Then \(w=z\) and \(P=e\), so absorption fails. If the parity vector is outside the affine hull, elementary linear algebra supplies a linear functional having a common value on all right parity vectors and the opposite value on epsilon(w). Assign its coordinates the corresponding elements e or g. The right words all have the same group value and w has the other one, and these values are incomparable.

For derivability, (T) with all three variables replaced by u gives

{\small\par\smallskip
\CaseLine{\(u^{3} \leq  u.\)}
\smallskip}

An affine combination over F2 is a sum of an odd number of vectors, with repetitions allowed. Start with cubes of the corresponding right words. Repeatedly apply (T) to three available words, where available means derivably below P. This gives an available long word v with \(\epsilon (v)=\epsilon (w)\). To make its content equal to C(P), for each still missing right summand u apply (T) to v, \(u^{3}\), \(u^{3}\). The product keeps its parity and acquires C(u). At the end v has content C(P).

If \(C(P)=C(w)\), Lemma 2 gives \(v=w\). Otherwise let Y be a nonempty word listing the variables of C(P) outside C(w), once each. By (E),

{\small\par\smallskip
\CaseLine{\(w^{3} \leq  w^{3} Y^{2}.\)}
\smallskip}

The word on the right has the content and parity of v, and Lemma 2 identifies them. The same lemma gives \(w=w^{3}\). Hence \(w\leq P\) is derivable. Soundness proves semantic sufficiency as well. All substitutions are nonempty; the case of empty missing content was treated separately. \(\square\)

The necessity of content and affine conditions also applies to quadratic words: the two separating evaluations in the first paragraph work unchanged at length two.

\par\smallskip\noindent\textit{4. The nilpotent-chain extension.}\quad 
Let L(P) be the set of variables occurring as linear summands of P.

Lemma 4. For A the complete absorption criterion is as follows.

\textbullet{} A linear word x is absorbed exactly when x belongs to L(P).

\textbullet{} A long word is absorbed exactly under the conditions of Lemma 3.

\textbullet{} A quadratic word w is absorbed exactly when its content is contained in C(P), its parity belongs to H(P), and P has a summand u of length at most two with C(u) contained in C(w).

Proof. For a linear target absent from L(P), assign \(x=b\) and every other variable z. Every right summand is at most a, whereas \(x=b\). This proves the linear assertion.

For a quadratic target, content and parity are necessary by Lemma 3. If the short contained-content witness is absent, give every variable of C(w) value b and every other variable value e. The target equals a. Every right summand is e or z: a right quadratic word cannot consist entirely of target variables, a right long word containing b is z, and no right linear word is a target variable. Thus \(P\leq z<a\).

Conversely apply Lemma 3 to \(w^{3}\). If the short witness is a linear letter x, use (A1), after orienting w as xy by (C). This includes a repeated-letter target. If the witness is a quadratic with the same content as w, it is w up to commutativity. In the remaining case the witness is \(x^{2}\) for one of the letters of w, and (Q) gives

{\small\par\smallskip
\CaseLine{\(xy \leq  x^{2} + (xy)^{3} \leq  P.\)}
\smallskip}

This proves sufficiency by explicit derivation. \(\square\)

\par\smallskip\noindent\textit{5. The extension with an exceptional square.}\quad 
The following elementary linear-algebra fact records why B differs from A at degree two.

Lemma 5. Let H be a nonempty affine subset of \(F2^{X}\), and suppose \(epsilon_{x}\)+\(epsilon_{y}\) belongs to H, allowing \(x=y\). There is a vector beta orthogonal to every element of H with \(beta_{x}=beta_{y}=1\) if and only if neither \(epsilon_{x}\) nor \(epsilon_{y}\) belongs to H.

Proof. If such a unit vector belongs to H, the requested beta cannot exist. Conversely put \(V=span(H)\), and write \(H=h0+D\) for a linear space D. Then \(V=H\) union D. If x differs from y and \(epsilon_{x}\) belongs to V but not H, it belongs to D, so

{\small\par\smallskip
\CaseLine{\(epsilon_{y}=(epsilon_{x}+epsilon_{y})+epsilon_{x}\)}
\smallskip}

belongs to H. Thus absence of both units from H implies absence of both from V. If \(x=y\), the hypothesis is 0 in H, so \(H=V\) and the same conclusion holds.

A system requiring beta in V-perp and \(beta_{x}=beta_{y}=1\) is inconsistent exactly when a vector of V supported on \{x, y\} has odd coordinate sum. When x differs from y, the only such vectors are the two units; when \(x=y\) there is just one. They have just been excluded, so the system is consistent. This is the usual duality for finite-dimensional linear equations over F2. \(\square\)

Lemma 6. For B, a linear target is absorbed exactly when it is a right linear summand, and long targets obey Lemma 3. A quadratic word w is absorbed exactly when its content is contained in C(P), its parity belongs to H(P), and at least one of these two conditions holds:

1. P has a quadratic summand u with C(u) contained in C(w).

2. For some letter x of w, \(epsilon_{x}\) belongs to H(P).

A linear right witness belonging to C(w) is already covered by 2.

Proof. For a missing linear target x, assign \(x=d\) and every other variable e. Every right word then has value e, whereas \(x=d>e\).

Let \(w=xy\) be quadratic. Suppose the content and affine conditions hold but neither 1 nor 2 holds. Lemma 5 gives beta orthogonal to all right parity vectors with \(beta_{x}=beta_{y}=1\). Assign every target variable the element c, and assign each outside variable the group element g raised to its beta coordinate. There is no linear target-variable summand, because it would put that unit vector in H(P). There is no right quadratic word made entirely of target variables, by failure of condition 1. Therefore no right word takes the exceptional value d or the clone value c. Every right word has value e by orthogonality. But \(w=cc=d>e\). This proves necessity.

If condition 1 holds, either the right witness is w up to (C), or it is a square of one of its letters. In the latter case use (Q) and the available long word \(w^{3}\). If condition 2 holds, Lemma 3 makes \(x^{3}\) and \(w^{3}\) available. Orienting w as xy, (B1) gives

{\small\par\smallskip
\CaseLine{\(xy \leq  x^{3} + (xy)^{3} \leq  P.\)}
\smallskip}

Thus all cases are derived from the proposed basis. \(\square\)

\par\smallskip\noindent\textit{6. Completeness and a proper inclusion.}\quad 
Proof of Theorem 1. By distributivity, every term is a nonempty finite sum of nonempty words. In a valid polynomial identity \(P=Q\), every summand of P is absorbed by Q. Lemmas 3, 4, and 6 give equational derivations of all these absorptions. Their sum yields \(P\leq Q\). The symmetric argument gives \(Q\leq P\) and hence \(P=Q\). This proves completeness for arbitrarily many variables in the five-element cases. Soundness for every allowed Boolean group, together with the C2 subsemirings, proves the general assertion. \(\square\)

Corollary 7. Var(B(C2, g)) is a proper subvariety of Var(A(C2)), and both varieties are finitely based.

Proof. In B, (A1) follows from (B1) and \(x^{3}\leq x\), the latter being an instance of (T). Therefore B satisfies the basis for A. Identity (B1) fails in A at \(x=y=b\): its left side is a, and its right side is z. This proves strictness. Alternatively the unary identity \(x^{2}\leq x^{3}+x^{4}\) is valid in B and fails in A at b. \(\square\)

The finite-basis conclusion for B comes from its own completeness proof. It does not follow merely from the proper inclusion.

Corollary 8. For either algebra, the number of distinct word functions on k positive generators is

{\small\par\smallskip
\CaseLine{\(3^{k} - 1 + k + k(k+1)/2.\)}
\smallskip}

Proof. There is one long-word function for every nonempty content with a prescribed parity at each of its variables, giving \(3^{k}\)-1. Different content or parity is separated by group/zero assignments. There are k distinct linear functions and k(k+1)/2 commutative quadratic functions. In A, linear functions are distinguished using b, and a quadratic is separated from its long counterpart by giving all its variables b. In B, use d for linear functions and c for quadratic functions. These also distinguish the low-degree functions from the long ones. \(\square\)

\par\smallskip\noindent\textit{7. Catalog correspondence and independent verification.}\quad 
The five-element A corresponds to table 9212 by

e \ensuremath{\to} 2, g \ensuremath{\to} 3, z \ensuremath{\to} 1, a \ensuremath{\to} 0, b \ensuremath{\to} 4.

The five-element B corresponds to table 9952 by

e \ensuremath{\to} 1, g \ensuremath{\to} 3, z \ensuremath{\to} 4, d \ensuremath{\to} 0, c \ensuremath{\to} 2.

Evaluating the operations at the 25 ordered pairs verifies both displayed isomorphisms. Proposition~\ref{prop:lower-fibres} in both directions gives \(\V(\A{9212})=\V(\A{10063})\). These presentations supply the three entries in case~\ref{grp:380}.

\subsection{Fifteen complete identities for 9222: cubic prefixes and quadratic recovery}\label{case:221}
Signature (+, multiplication), no constants, relative to additively idempotent semiring axioms. No minimality is asserted.

\par\smallskip\noindent\textit{1. Word values and exact absorption.}\quad 
Addition has chain order \(3<2<1<0<4\); its rows in numerical order are 00004,\allowbreak{}01114,\allowbreak{}01224,\allowbreak{}01234,\allowbreak{}44444. Multiplication rows are

{\small\par\smallskip
\CaseLine{\ttfamily 1 1 1 1 1}
\CaseLine{\ttfamily 1 1 1 1 1}
\CaseLine{\ttfamily 1 1 2 3 1}
\CaseLine{\ttfamily 3 3 3 3 3}
\CaseLine{\ttfamily 1 1 1 1 0}
\smallskip}

A literal word retains its assigned value. A quadratic word has value 0 if both factors are 4. Every other word of length at least two has value 2 when all factors are 2; otherwise look at its first factor not equal to 2. If that factor is 3 the value is 3; if it is 0,\allowbreak{}1,\allowbreak{}4 the value is 1. This follows directly by induction from the table: 1 and 3 are left zeros, 2 is an identity on \{1,\allowbreak{}2,\allowbreak{}3\}, and the exceptional \(4*4=0\) becomes 1 as soon as another factor is multiplied.

Write \(C=c(w)\), and \(pre_{w}\)(x) for the variables before the first x. All witnesses below are original right summands in P, not constructed words. A literal target x is absorbed exactly when x itself is in P. For a target w of length at least two the exact criterion is:

(C) There is a right word u with c(u) contained in C. If |w|=2, this witness must also have length at most two. (F) For each x in C, there is a right word v containing x with \(pre_{v}\)(x) contained in \(pre_{w}\)(x). The suffix of v may use outsiders.

For necessity of (C) in the long case, assign C the value 2 and all outsiders 3: the target is 2 and every noncontained word is 3. For quadratic (C), assign C the value 4 and outsiders 3. The target is 0; a right word reaches 0 or 4 only if it is contained and has length at most two. A missing F-witness at x is separated by assigning the target prefix 2, x the value 1, and all other variables 3. The target is 1, whereas right words are 2 or 3. A missing literal x is separated by \(x=4\) and all other variables 3; no longer word can attain 4.

For sufficiency, target 3 is bottom. If the target is 2, all its variables are 2 and a contained witness has value 2. A target 1 has a first non-2 factor x valued 0,\allowbreak{}1 or 4, with all preceding variables

2. Its F-witness has value at least 1, including the possible literal or quadratic exceptions. The only remaining nonlinear target value is quadratic 0; its short contained witness has value 0 or 4. This proves the criterion for arbitrary variables, lengths and sums.

\par\smallskip\noindent\textit{2. The basis and long normal forms.}\quad 
The following fifteen identities form the basis; inequalities mean absorption equations \(u+P=P\).

{\small\par\smallskip
\CaseLine{\((N1) x^{2}yz = xyz\qquad (N2) xy^{2}z = xyz\)}
\CaseLine{\((N3) xyz^{2} = xyz\qquad (N4) xyxz = xyz\)}
\CaseLine{\((N5) xyzy = xyz\qquad (N6) xyx = x^{2}y\)}
\CaseLine{\((N7) xy^{2} = x^{2}y\)}
\CaseLine{\((L)  x^{2} \leq  x\)}
\CaseLine{\((P0) xy \leq  x+y\)}
\CaseLine{\((P1) xy \leq  y+xz\)}
\CaseLine{\((P2) xyz \leq  xz+xyt\)}
\CaseLine{\((P3) xyz \leq  xz+xy\)}
\CaseLine{\((Q1) xy \leq  x^{2}+xyx\)}
\CaseLine{\((Q2) xy \leq  y^{2}+xyx\)}
\CaseLine{\((Q3) xy \leq  yx+xyx\)}
\smallskip}

Each N identity has length at least three on both sides and preserves first-occurrence order, so the word formula proves it. L follows from the table. Each of P0--P3 and Q1--Q3 satisfies C, F above: for P0, P1 a literal right word supplies the short contained witness; for P2, P3 xz is contained, x is witnessed at the start, y by xyt or xy, and z by xz. For Q1--Q3 the displayed quadratic term is short and contained and xyx supplies every first-prefix witness. Thus every listed law is valid, including arbitrary nonempty term substitutions.

Leave words of length one or two unchanged. For length at least three, let F(w) be its first-occurrence list. Its normal form is \(x^{3}\) if \(F(w)=x\), \(x^{2}y\) if \(F(w)=xy\), and F(w) when |F(w)|\ensuremath{\geq}3. Here is a complete derivation of that normal form. While length is at least four and a repetition exists, delete a repeated occurrence. An adjacent pair at the start, interior or end is reduced by N1, N2, N3, respectively, splitting the remaining context into nonempty blocks as needed. A nonadjacent pair aUa with a nonempty right context is reduced by N4; if only a left context remains, use N5. If both external contexts are empty, |U|\ensuremath{\geq}2: N6 gives \(a^{2}U\) and N1 then gives aU. All deletions preserve first-occurrence order and stop at length at least three. At length three, the two-variable words xyx and \(xy^{2}\) become \(x^{2}y\) by N6, N7; one-variable words are \(x^{3}\). This also proves that two long words have the same derived normal form precisely when their first occurrence lists agree. No empty words are substituted for variables.

\par\smallskip\noindent\textit{3. Unrestricted long-word construction.}\quad 
Let \(w\leq P\) be valid and nonlinear and let \(F(w)=x_{1}\cdots x_{d}\). Define \(p_{i}=x_{1}^{3}\cdots x_{i}^{3}\). We derive \(p_{i}\) \(u\leq P\) for EVERY original u in P.

For \(x=x_{1}\) choose its F-witness v, necessarily starting with x. If \(v=xt\) with nonempty t, use L twice in context to derive \(x^{3}t\leq v\). P1 with x replaced by \(x^{3}\) then gives \(x^{3}u\leq u+x^{3}t\leq P\). If \(v=x\) is literal, L twice gives \(x^{3}\leq x\) and P0 gives \(x^{3}u\leq x^{3}+u\leq P\). This establishes the first stage for all originals.

Suppose the stage for \(p=p_{i}\) has been derived and \(x=x_{i+1}\). Choose an original F-witness \(v=rx\) t, where the prefix r may be empty, its content lies in c(p), and t may be empty. The induction hypothesis already gives \(pv\leq P\). Both pv and \(px^{3}t\) have length at least three and the same first-occurrence list, so N1--N7 derive their equality. If t is nonempty, P2 with substitutions \((x,y,z,t)=(p,x^3,u,t)\) gives \(px^{3}u\leq pu+px^{3}t\leq P\). If t is empty, use P3 with \((x,y,z)=(p,x^3,u)\) to obtain \(px^{3}u\leq pu+px^{3}\leq P\). This completes the induction without restricting any suffix variables.

Choose an original contained C-witness \(u_{0}\). Then \(p_{d}\) \(u_{0}\) has the first-occurrence list F(w) and length at least three. If |w|\ensuremath{\geq}3, normalization immediately derives \(w\leq P\). If |w|=2, the same construction derives its long guard: \(x^{3}\) for \(w=x^{2}\), or \(x^{2}y\) for \(w=xy\) with \(x\ne y\). In the latter case this guard also equals xyx by N6.

\par\smallskip\noindent\textit{4. Quadratic recovery and all identities.}\quad 
Let \(w=xy\) with \(x\ne y\). The short contained C-witness is one of x, y, \(x^{2}\), \(y^{2}\), xy, yx. If it is xy the target is already present. If it is x or y, first use L to derive \(x^{2}\) or \(y^{2}\). Thus either \(x^{2}\), \(y^{2}\), or yx is below P; together with the guard xyx, Q1, Q2 or Q3 gives \(xy\leq P\). For \(w=x^{2}\) the short witness is x or \(x^{2}\), and L or literal inclusion completes the proof. Literal targets are already summands.

We have derived every valid word absorption with no bound on the number of variables, lengths or number of right summands. Expand all terms into nonempty polynomials by distributivity. In a valid polynomial equality, derive each summand's absorption in the opposite polynomial, in both directions. This proves equational completeness of the fifteen identities. Reversing words gives a basis of the multiplicative opposite. Other catalogue types require independent complete HSP in BOTH directions.

The computations certify soundness, the semantic implementations and literal derivation replays; bounded regressions are not the completeness argument. The long construction uses original prefix witnesses with arbitrary tails, whereas quadratic targets require a short contained witness. The latter condition must not be omitted.

\subsection{Internal and global nonlinear anchors: two five-identity bases}\label{case:222}
All identities are in the constant-free signature \((+,\cdot)\) relative to the additively idempotent semiring axioms. Inequality denotes absorption.

\par\smallskip\noindent\textit{1. The algebras and their bases.}\quad 
Both algebras have elements \(\{0, f, b, a, t\}\). Multiplication is commutative, \(0\) is a zero, the product of any two elements of \(\{a, t\}\) is \(b\), and every other product of two nonzero elements is \(f\). In \(A\), addition is join in

\CaseDisplay{0<f<b<t,\qquad 0<a<t,}

with \(a\) incomparable with \(f, b\). In \(B\), addition is join in

\CaseDisplay{0<f<b<t,\qquad a<t,}

with \(a\) incomparable with \(0, f, b\). The common catalog map, in the displayed element order, is \((2,\allowbreak{}1,\allowbreak{}0,\allowbreak{}3,\allowbreak{}4)\). These are types 9625 and 9626. The letter \(0\) names an element only, not a constant in terms.

A product of at least three nonzero factors is \(f\), whereas a product with a zero factor is \(0\). This proves associativity. Translation by \(0\) is constant. Translation by \(f\) or \(b\) sends \(0\) to \(0\) and all other elements to \(f\); translation by \(a\) or \(t\) sends \(0\) to \(0\), \(f, b\) to \(f\), and \(a, t\) to \(b\). These maps preserve the indicated chain joins and the incomparable joins with \(a\), in both orders. Thus both constructions are semirings.

Theorem. Let \(\Sigma\) consist of

\CaseDisplay{\begin{array}{rll}
(C)&xy\approx yx,\\
(N)&x^2yz\approx xyz,\\
(T)&xyz\leq xy,\\
(S)&xy\leq x^2.
\end{array}}

A complete basis for \(A\) is \(\Sigma\) together with

\CaseDisplay{(I)\qquad xy\leq x+xy^2.}

A complete basis for \(B\) is \(\Sigma\) together with

\CaseDisplay{(G)\qquad xy\leq x+zt.}

No minimality is asserted.

\par\smallskip\noindent\textit{2. Normal forms and exact criteria.}\quad 
By (C),(N), every word of length at least three is equal to \(R_X=\prod_{x\in X}x^3\), where \(X\) is its content. To see this, commute any chosen occurrence to the front, split the remaining at least two positions into nonempty factors, and use the reverse of (N) to duplicate it. Continue until each variable occurs at least three times; then remove surplus copies by \(x^4=x^3\), an instance of (N). In particular \(w^3=R_{c(w)}\), and \(R_XR_Y=R_{X\cup Y}\).

Normalize a nonempty polynomial \(F\) into linear, quadratic, and stable summands. Let \(L\) be its linear variables, and let \(\mathcal N\) be its nonlinear summands. For nonempty \(X\), let \(L_X=L\cap X\), let \(\mathcal N_X\) be the nonlinear summands with content contained in \(X\), and let \(E_X\) be the quadratic members of \(\mathcal N_X\).

Lemma. Linear absorption in either algebra is exactly \(x\in L\). The nonlinear criteria are

\noindent Algebra\enspace\textbar\enspace \(R_X\leq F\)\enspace\textbar\enspace \(xy\leq F\), \(X=\{x, y\}\)\par
\noindent \(A\)\enspace\textbar\enspace \(\mathcal N_X\ne\varnothing\)\enspace\textbar\enspace \(\mathcal N_X\ne\varnothing\) and \((E_X\ne\varnothing\) or \(L_X\ne\varnothing)\)\par
\noindent \(B\)\enspace\textbar\enspace \(\mathcal N\ne\varnothing\) and \((\mathcal N_X\ne\varnothing\) or \(L_X\ne\varnothing)\)\enspace\textbar\enspace \(\mathcal N\ne\varnothing\) and \((E_X\ne\varnothing\) or \(L_X\ne\varnothing)\)\par
Equal endpoints \(x=y\) are permitted.

Proof. A missing linear variable is separated by assigning it \(t\) and all other variables \(0\). Every nonlinear right value is at most \(b<t\); right linear values are \(0\).

For \(A\), if \(\mathcal N_X\) is empty, assign all variables of \(X\) the value \(a\), and all other variables \(0\). The right value is \(0\) or \(a\), whereas the left stable value is \(f\) and the left quadratic value is \(b\); neither is below \(0\) or \(a\). If \(E_X\) and \(L_X\) are both empty, the same assignment makes every right summand at most \(f\), which separates a quadratic left value \(b\). These prove necessity. For sufficiency, a zero left value is automatic since \(0\) is least in \(A\). A positive stable or quadratic value at most \(f\) is dominated by any internal nonlinear summand, whose value is \(f\) or \(b\). For a quadratic value \(b\), all variables of \(X\) are in \(\{a, t\}\). An internal quadratic has value \(b\). If the short anchor is instead a linear variable, its value is \(a\) or \(t\); joined with the positive internal nonlinear value, it dominates \(b\) since \(a+f=a+b=t\).

For \(B\), a purely linear right side is ruled out by assigning every variable \(a\): the left nonlinear value is \(f\) or \(b\), neither below \(a\). If both internal nonlinear and internal linear anchors are absent, assign \(X\) the value \(f\) and all other variables \(0\). This separates a stable left value \(f\) from \(F=0\). If both internal quadratic and internal linear anchors are absent, assign \(X\) the value \(a\) and all other variables \(0\); the left quadratic is \(b\), but every right summand is at most \(f\).

Conversely, the presence of any nonlinear summand ensures \(F\geq0\), so a zero left value is absorbed even though \(0\) is not least in \(B\). For a positive stable value \(f\), an internal nonlinear term is positive. Alternatively take a linear anchor from \(X\). Its value is \(f, b, t\), already at least \(f\), or it is \(a\); in the latter case its join with any nonlinear right value, even \(0\), is \(t\). For a quadratic value \(b\), an internal quadratic also has value \(b\), or an internal linear anchor has value \(a\) or \(t\) and the same global nonlinear argument applies. For a quadratic value \(f\), the same anchors give at least \(f\). This proves sufficiency in \(B\) and completes the lemma. \(\square\)

\par\smallskip\noindent\textit{3. Validity and unrestricted completeness.}\quad 
The product formula verifies (C),(N). Every extra factor lowers a nonlinear value in the chain \(0<f<b\), proving (T). If \(xy=b\), both \(x, y\) are active, so \(x^2=b\); if \(xy=f\), then \(x^2\) is \(f\) or \(b\); and if \(xy=0\), it is below every square. This proves (S). The two additional laws follow from the respective criteria: (I) supplies an internal linear and stable summand; (G) supplies a linear anchor and a global nonlinear summand.

First, if an internal nonlinear summand \(r\) has content \(Y\subseteq X\), (T) gives \(r^3\leq r\), and normalization gives \(R_Y\leq F\). Applying (T) repeatedly permits insertion of each missing variable into \(R_Y\) while decreasing its value. Normalization then gives \(R_X\leq F\). Every term in these applications is nonempty; the prefix being extended always has length at least two.

For a quadratic target on \(X\), an internal quadratic is the target itself or, when its endpoints are distinct, one of their squares. Law (S), after interchanging variables if needed, handles the squares. In \(A\), if the short anchor is linear, say \(x\), the nonlinear anchor has already supplied \(R_X=xy^2\) (also when \(x=y\)). Law (I) gives \(xy\leq x+R_X\leq F\). This derives all absorptions for \(A\).

In \(B\), if the needed internal anchor is linear, say \(x\), write any global nonlinear right summand as \(uv\) with nonempty \(u, v\). Law (G) gives \(xy\leq x+uv\leq F\), for any quadratic target with endpoint \(x\). It also gives \(x^2\leq x+uv\); by (T), \(x^3\leq x^2\), so \(x^3\leq F\). Inserting the other variables of \(X\) via (T) gives \(R_X\leq F\). Together with the internal nonlinear case, this derives every absorption in the lemma for \(B\).

Normalize any valid polynomial identity and apply these derivations to each summand in both directions. The resulting two inequalities prove that identity. No variable or summand bound was used. \(\square\)

\par\smallskip\noindent\textit{4. Transfers and scope.}\quad 
The equal-variety classes certified in both directions are

\textbullet{} 9625: 9435, 9437, 9625;

\textbullet{} 9626: 9436, 9626, 9664, 10452.

The bases are separate: the first cannot be transferred to the second by a single variety inclusion. Independent checks compare complete word closures and exact all-polynomial absorption CNFs through five variables; the unrestricted proof is Section 3. These results do not complete the full five-element classification.

The source manuscript for this case credits Aifa Wang and Lili Wang, School of Mathematical Sciences, Chongqing University of Technology.

The source manuscript for this case acknowledges support from the National Natural Science Foundation of China (12371024), the Science and Technology Research Program of Chongqing Municipal Education Commission (KJZD-K2024011102) and the Chongqing Natural Science Foundation Innovation and Development Joint Fund (Municipal Education Commission) (CSTB2025NSCQ-LZX0067).

\subsection{A four-identity basis from internal coverage and a quadratic anchor}\label{case:223}
We use the signature \((+,\cdot)\) without constants, relative to the additively idempotent semiring axioms. Write \(u\leq v\) for \(u+v\approx v\).

\par\smallskip\noindent\textit{1. Construction and basis.}\quad 
Let \(A=\{0, e, f, b, a\}\) be the chain \(0<e<f<b<a\), with addition maximum. The element \(0\) is a multiplicative zero. For positive elements put \(e^2=e\), \(a^2=b\), and make every other product \(f\). These element names are not language constants. In the displayed order the catalog map is \((3,\allowbreak{}2,\allowbreak{}1,\allowbreak{}0,\allowbreak{}4)\), identifying \(A\) with type 9458.

A product of at least three positive factors is \(e\) when all factors are \(e\), and \(f\) otherwise. This proves associativity, including products with zero. Translations are increasing on the chain: a positive translation has value \(0\) at \(0\), and thereafter is either \(e, f, f, f\), \(f, f, f, f\), or \(f, f, f, b\). Hence they preserve maximum and \(A\) is a semiring.

Theorem. A complete basis for \(A\) is

\CaseDisplay{\begin{array}{rll}
(C)&xy\approx yx,\\
(N)&x^2yz\approx xyz,\\
(B)&xy\leq x+y,\\
(Q)&xy\leq x^2+xy^2.
\end{array}}

No minimality is asserted. The proof below does not bound the number of variables or summands.

\par\smallskip\noindent\textit{2. Normalization and exact absorption.}\quad 
By (C),(N), every word of length at least three is equal to

\CaseDisplay{R_X=\prod_{x\in X}x^3,\qquad X=c(w).}

Indeed, in a word of length at least three, choose any occurrence of a variable and split the remaining positions into two nonempty factors. The inverse of (N), after commuting, duplicates that occurrence. Repeating this produces at least three copies of each variable. Substitution \(y=z=x\) in (N) gives \(x^4=x^3\), so surplus copies can be removed. This also proves \(R_XR_Y=R_{X\cup Y}\) and \(w^3=R_{c(w)}\).

Thus a polynomial may be written as a nonempty sum of linear terms, quadratic terms, and stable terms \(R_X\). For a fixed nonempty set \(X\), let \(F_X\) be the family of right summands whose contents are contained in \(X\). A quadratic anchor means a member of \(F_X\) of length one or two. The phrase includes linear summands; it refers to the degree bound.

Lemma. The following conditions exactly characterize absorption:

\textbullet{} \(x\leq F\) iff \(x\) itself is a linear summand of \(F\).

\textbullet{} \(R_X\leq F\) iff \(X\) is covered by the contents of \(F_X\).

\textbullet{} A quadratic \(xy\), with \(X=\{x, y\}\) (possibly \(x=y\)), is absorbed iff \(F_X\) covers \(X\) and contains a quadratic anchor.

Proof. A missing linear \(x\) is separated by assigning \(x=a\) and every other variable \(0\): every nonlinear right term is at most \(b<a\). If coverage fails at \(x\in X\), assign \(x=f\), all other variables of \(X\) the value \(e\), and variables outside \(X\) the value \(0\). The left value is \(f\), while every right summand is at most \(e\). If a quadratic anchor is absent, assign all variables of \(X\) the value \(a\) and all other variables \(0\). The left quadratic equals \(b\) and every right summand is at most \(f\). These prove necessity.

For sufficiency, a zero left value is automatic. Otherwise all variables of \(X\) have positive values. If the left value is \(e\), an internal summand exists by coverage and is positive, hence at least \(e\). A stable left value \(f\), or a quadratic left value \(f\), means some variable of \(X\) has value strictly greater than \(e\). An internal summand covering that variable has value at least \(f\). Finally, a quadratic has value \(b\) only when all its variables are \(a\); the anchor then has value \(a\) or \(b\). This establishes every sufficiency assertion. \(\square\)

\par\smallskip\noindent\textit{3. Completeness.}\quad 
The multiplication formula verifies (C),(N). It also verifies (B): a product is at most the maximum of its two factors. For (Q), its two right terms cover \(\{x, y\}\) internally, and \(x^2\) is a quadratic anchor, so the lemma proves validity.

Law (B) gives \(r^2\leq r\) and \(r^3\leq r\) for every nonempty term \(r\). Suppose the internal summands \(r_1,\ldots, r_k\) cover \(X\). Each \(r_i^3=R_{c(r_i)}\) is derivably below \(F\). Their product is \(R_X\). Repeated use of (B) yields

\CaseDisplay{R_X=r_1^3\cdots r_k^3\leq r_1^3+\cdots+r_k^3\leq F.}

For \(k=1\) the first equality alone suffices. This derives every stable absorption in the lemma.

Consider \(xy\leq F\). If \(x=y\), the anchor on the one-variable content is \(x\) or \(x^2\), and \(x^2\leq x\) supplies the only nontrivial case. If \(x\ne y\), an internal quadratic anchor is either the target \(xy\), a square, or a linear summand. The first case is immediate. A linear anchor \(x\) or \(y\) absorbs its square by (B). Consequently, after possibly interchanging \(x, y\), we have \(x^2\leq F\). Coverage has already given \(R_{\{x, y\}}\leq F\), and (N) identifies \(xy^2\) with this stable term. Law (Q) now gives \(xy\leq x^2+xy^2\leq F\).

Linear absorptions are original summands. Normalize both sides of any valid identity and apply these derivations to each summand, giving both inequalities. This proves the theorem for arbitrary finite variable sets. \(\square\)

\par\smallskip\noindent\textit{4. Transfers and scope.}\quad 
Types 9458 and 10064 generate equal varieties, as certified by complete word-function closures and principal-ideal tests in both directions. The displayed basis therefore applies to both types. Finite-variable CNF comparisons independently check the absorption lemma but do not replace its unrestricted proof. No assertion is made about all subvarieties of the generated variety, or about completion of the full order-five problem.

The source manuscript for this case credits Aifa Wang and Lili Wang, School of Mathematical Sciences, Chongqing University of Technology.

The source manuscript for this case acknowledges support from the National Natural Science Foundation of China (12371024), the Science and Technology Research Program of Chongqing Municipal Education Commission (KJZD-K2024011102) and the Chongqing Natural Science Foundation Innovation and Development Joint Fund (Municipal Education Commission) (CSTB2025NSCQ-LZX0067).

\subsection{A finite basis for 9466 from two original witnesses}\label{case:224}
In the semiring signature without constants take the semiring axioms, \(xyz=xzy\), \(xxyz=xyz\), \(xyyz=xyz\), \(xxy=xyy\), and all table-valid absorptions over at most 16 variables, length 18, and two upper words. This specified finite family is a complete basis, not an enumerated small list.

Element 3 is a multiplicative absorber and the least addition element; the additive chain \(is 3<2<1<0<4\). Put B=\{0,\allowbreak{}1,\allowbreak{}2,\allowbreak{}4\}. A word containing 3 has value 3. Otherwise a quadratic is 2 when its head is 2, is 0 on(4,\allowbreak{}4), and is 1 otherwise. Every long word is 2 when its head is 2 and 1 otherwise. All triples and the recurrence \(f(b)c=f(b)\) on B prove this for every length.

The four laws normalize long words by head and content:hhh for one variable, hht for two variables, and h followed by the remaining distinct letters in a fixed order otherwise. Permute tails; above length 3 delete an extra head by \(xxyz=xyz\) or a repeated tail by \(xyyz=xyz\). At length 3 use \(xxy=xyy\) for the two-variable case. Short degree/order are retained.

For target w, let U consist of original upper words contained in c(w). A failing target assignment avoids3. Set every external variable 3. All discarded words become least, so failure is equivalent to failure of U on B. The exact criterion is:

\textbullet{} A literal target x requires the original literal x.

\textbullet{} A nonlinear target requires a contained original with the same head.

\textbullet{} A quadratic target also requires a contained original of length at most 2; this witness need not share the target head.

Check all four proper principal ideals. In down(0), every nonlinear word on B is allowed; only literals exclude4. Thus a literal target requires its own literal. In down(1), long words are unrestricted, literals require1 or 2, and a quadratic forbids precisely(4,\allowbreak{}4). A quadratic target over its one or two variables is forced into this ideal exactly when U has a short word:if both target letters were 4, every such short word would fail; without one, assign both letters 4. In down(2), any word lies in the ideal exactly when its head is 2. The target head is forced precisely by an original same-head word. Down(3) allows no word over B; the same-head witness makes this system impossible. These four checks prove sufficiency.

For explicit necessity:if literal x is absent, put \(x=4\), others 3. If a nonlinear target lacks a contained same-head word, put its head 1, its other variables 2, externals 3; upper words have value 2 or 3, target 1. If a quadratic target lacks a short contained word, put its variables 4, externals 3; target 0, upper retained long words 1. Empty U and repeated letters are covered. The polynomial itself is nonempty.

Retain the literal witness, or the same-head witness plus a short witness when quadratic. Thus at most two original words suffice. Normalize these and the target. Reserve all variables of short words and each long head plus up to two other distinct letters. At most 9 variables are reserved. Group other variables by common incidence among these at most 3 words:at most 7 groups. Identification gives at most 16 variables, normal length at most 18, and preserves validity. The resulting absorption belongs to the finite family. Substitute each group by the NONEMPTY word listing its members. Every head and short occurrence is reserved; long degree is preserved. Common incidence recovers each long content, and the word laws recover the original normal forms. Add back unselected summands. Distributive expansion and both directions of summand absorption establish unrestricted completeness.

The forbidden quadratic pair is not a rectangle; the preceding unary domain theorem for 2947/13188 is not being applied. Equal-variety and opposite transfers require separate complete mutual HSP certificates.

\subsection{Finite bases for 9677 and 10575 from original head and ideal witnesses}\label{case:225}
Use the semiring signature without constants. For each source take the semiring axioms, \(xyz=xzy\), \(xxyz=xyz\), \(xyyz=xyz\), \(xxy=xyy\), and all table-valid absorptions over 27 variables, length 29, at most 3 upper words. The finite family is specified, not enumerated as a small list. Completeness below has no restriction on the original alphabet or words.

Both sources normalize long words by head and content:hhh for singleton content, hht for two variables, and h followed by the other distinct letters in a fixed order otherwise. Tail permutation and the two length-four deletion laws reduce duplicates without crossing degree 3; the fourth law handles the remaining length-three two-variable case.

9677 has additive \(chain 3<2<1<0<4\). Its multiplication has two left zeros 2,\allowbreak{}3. A nonlinear word has value 3 if its head is 3; otherwise value 2 if some letter is 2 or 3; otherwise value 0 for the quadratic pair(4,\allowbreak{}4), and 1 in all other cases. In particular long words have no exceptional value 0. Induction from the table proves this formula, including arbitrarily long words. For original polynomial P and target w, let U be the original words contained in c(w), and let H be the heads of ALL original P words. The exact criterion is:

\textbullet{} a literal x requires its original literal;

\textbullet{} a long target requires its head in H and a nonempty U;

\textbullet{} a quadratic target requires its head in H and a short word in U.

Indeed, down(3) requires every upper head to be 3 and therefore forces the target head exactly when it belongs to H. In down(2), a word is allowed exactly when its content meets the set of variables valued 2 or 3. These upper intersection requirements imply the target requirement exactly when a contained original exists. In down(1), all long words are allowed, quadratics exclude only(4,\allowbreak{}4), and literals allow1,\allowbreak{}2,\allowbreak{}3. A quadratic target is forced exactly by a short contained word; otherwise give target variables 4 and all outside variables 2. Down(0) permits all nonlinear words and excludes4 on literals, so a literal target needs its original literal. Literal necessity can also be witnessed by \(x=4\), others 2. Missing global head is separated by target head 2, all other variables 3. Missing contained content is separated by all target variables 1, outside 2. These checks cover all proper principal ideals and prove equivalence.

10575 has proper principal ideals I0=\{0,\allowbreak{}1,\allowbreak{}2\}, I1=\{1,\allowbreak{}2\}, \(I2=\{2\}\), I3=\{1,\allowbreak{}2,\allowbreak{}3\}. Its nonlinear word is 2 if its head is 2; is 1 if its head is 1 or a tail letter is 1 or 2; is 3 for the quadratic pair(0,\allowbreak{}0); and is 4 otherwise. This formula follows by direct multiplication and induction.

For 10575 let U again be the original words contained in c(w). Let L be the variables occurring as literal words in U, Q the union of the contents of quadratic words in U, and J the long words in U. Let H be the heads of ALL original P words. Put A=(J nonempty or L intersects Q). The exact criterion is:

\textbullet{} a literal x requires the original literal x or a long original in U;

\textbullet{} a long target requires its head in H and A;

\textbullet{} a quadratic target requires its head in H and either A or c(w) subset Q.

To prove it, first I2 allows a word exactly when its head is 2. Thus the target head must be in H. I1 allows a word exactly when some letter is 1 or 2, so it requires a contained original. The stated branches already supply both requirements (including the literal branch). For I0 and I3, a failing nonlinear target uses only D=\{0,\allowbreak{}3,\allowbreak{}4\}; a literal failure also has its value in D. Set all outside variables 1. Every upper word outside U then has value 1. For values on D, each contained long word is 4. In I0, contained literals must be 0 and no contained nonlinear word is allowed. In I3, contained literals must be 3, contained quadratic letters must all be 0, and long words are impossible. Consequently J nonempty or L intersect Q makes I3 impossible. Otherwise a long target fails I3; and a quadratic target is forced into I3 exactly when both its letters are forced 0, that is, c(w) subset Q. These nonlinear branches also supply a contained nonlinear word, so I0 is impossible. For a literal target, all contained words are powers of its only variable. An original literal or long word makes every proper ideal force the target. Without either, set the variable 0 and outside variables 1:the target is 0 and the upper sum is 1 or 3. This is a failure. These arguments establish both directions for every ideal; no outside-word containment assumption is imposed on the global head.

For 9677 retain a global same-head witness and the contained witness (short when quadratic); for a literal retain itself. At most 2 originals suffice. For 10575 retain a global head witness and one long contained word, or a conflicting literal and quadratic, or at most two quadratics covering the quadratic target's variables. A literal needs one literal or one long witness. Thus at most 3 originals suffice. All are chosen from original P. The same criteria remain true for the selected polynomial, so it still absorbs the target.

There are at most 4 words including the target. Reserve every short-word variable and each long head plus up to two other distinct letters. At most 12 variables are reserved. Group the others by common incidence among the 4 words:at most 15 groups. Identification yields at most 27 variables and normal length 29, and preserves validity. The compressed absorption is in the stated finite family. Substitute each group variable by the NONEMPTY word listing its members. Reserved heads and short occurrences, common incidence, and preserved long degree let the four word laws recover every original normal form. Add the omitted original summands. Expanding terms distributively and applying summand absorption in both directions proves unrestricted finite basability.

Both global head witnesses may contain variables outside the target. For 10575, literals alone do not replace the long-word/conflict condition; a quadratic target can instead use quadratic content coverage. For 9677, the quadratic short guard cannot be replaced by a long word.

\subsection{Four finite bases for five-element semirings with truncated multiplication}\label{case:226}
The theorem below concerns identities in the constant-free signature \((+,\cdot)\), relative to the additively idempotent semiring axioms. It proves finite basability for four explicit algebras, with no bound on the number of variables.

\par\smallskip\noindent\textit{1. The algebras and the bases.}\quad 
For \(r\in\{1,\allowbreak{}2,\allowbreak{}3,\allowbreak{}4\}\), let

\CaseDisplay{B_r=\{e, a, b, z, d\}=\{v_0, v_1, v_2, v_3, d\}.}

The core has chain order \(v_3<v_2<v_1<v_0\), addition is its maximum, and

\CaseDisplay{v_i v_j=v_{\min(3, i+j)},\qquad dv_i=v_id=v_i,\qquad d^2=e.}

For addition involving \(d\), put

\CaseDisplay{d+d=d,\qquad d+v_i=
 \begin{cases}d&i\ge r,\\ e&i<r.\end{cases}}

In particular \(d+e=e\). The element names are not term constants. The labelings \((e, a, b, z, d)\) in the recovered catalog tables are \((4,\allowbreak{}1,\allowbreak{}2,\allowbreak{}3,\allowbreak{}0)\) for 11947 and \((4,\allowbreak{}0,\allowbreak{}1,\allowbreak{}2,\allowbreak{}3)\) for the other three.

Use \(u\leq P\) as an abbreviation for \(u+P\approx P\). The following five laws are common to all four bases:

\CaseDisplay{\begin{array}{ll}
(C)&xy\approx yx,\\
(N)&x^3\approx x^4,\\
(D)&xyz\leq xy,\\
(S)&xyz\leq x^2yz+xy^2z+xyz^2,\\
(J)&xy\leq x^2+xy^3.
\end{array}}

The additional laws are:

\CaseDisplay{\begin{array}{c|c|c}
\text{algebra}&\text{catalog representative}&\text{additional laws}\\ \hline
B_1&11947&(L)\ xy\leq x+(xy)^3\\
B_2&10034&(G)\ xy\leq x+x^3y,\quad
             (H)\ xyz\leq x+xy^2z+xyz^2\\
B_3&10035&(G)\ xy\leq x+x^3y\\
B_4&10036&(B)\ xy\leq x+z^2
\end{array}}

Theorem 1. The displayed laws form identity bases for \(B_1, B_2, B_3, B_4\), respectively. The sizes relative to the semiring axioms are 6, 7, 6 and 6.

The operation definitions indeed give semirings. Addition is the join in the indicated partially ordered set: \(d\) is above precisely those core elements with index at least \(r\), and its join with every other core element is \(e\). Every product of at least two factors lies in the core, and its index is the sum of the factor indices, truncated at 3, where the index of \(d\) is 0. This proves associativity. For distributivity within the core use truncated addition of indices. For a join \(d+v_i\), multiplication by \(v_j\) gives \(v_j\) on both sides, and multiplication by \(d\) gives \(e\) on both sides. Commutativity proves the other distributive law.

For a word \(w\), write \(c(w)\) for its content and \(n_x(w)\) for the multiplicity of \(x\). By (C),(N), every word is equivalent to its sorted word with each positive multiplicity capped at 3. We always use this normalization below. Contexts that are empty are simply omitted; no empty word is substituted into an identity.

Soundness of (C),(N),(D) follows from the index formula. For (S), if all three indices are positive the product is \(z\); otherwise one of the terms obtained by repeating a zero-index factor has the same value. For (J), let the two indices be \(i, j\). If \(i=0\), the term \(x^2\) is \(e\); if \(j=0\), \(xy^3=xy\). If \(i, j>0\), either \(i=1\), when \(x^2=a^2=b\) is at least \(xy\), or \(i\ge2\), when \(xy=z\). All words used here have length at least two, so this argument includes assignments involving \(d\).

For the additional laws, when \(x\) is a core element every product with initial factor \(x\) is at most \(x\). When \(x=d\), law (G) follows from \(d^3=e\). For (L) in \(B_1\), a positive-index product is below \(d\), and a zero-index product has cube \(e\). For (H) in \(B_2\), two positive indices for \(y, z\) give a product below \(d\); if one index is zero, its repeated-factor summand has the original value. For (B) in \(B_4\), \(z^2\) is a core element and \(d+z^2=e\).

\par\smallskip\noindent\textit{2. Completeness for sums of nonlinear words.}\quad 
Here nonlinear means length at least two. This convention does not require a repeated variable.

For a nonlinear word \(w\) and a nonempty sum \(R\) of nonlinear words, let \(\mathcal K(w, R)\) mean

\CaseDisplay{\min(3,\sum_x n_x(w)\lambda_x)
 \ \ge\
 \min_{u\in R}\min(3,\sum_x n_x(u)\lambda_x)
 \quad\text{for every }\lambda_x\in\{0,\allowbreak{}1,\allowbreak{}2,\allowbreak{}3\}.}

This is precisely absorption of \(w\) by \(R\) in the core chain. It is also precisely absorption in every \(B_r\), since nonlinear words only see the indices. Declare \(\mathcal K(w,\varnothing)\) false.

To check \(\mathcal K\), variables outside \(c(w)\) can all be assigned index 3. Words containing such a variable can then be discarded. Only assignments giving \(w\) index less than 3 matter. Their positive indices inside \(c(w)\) comprise at most one variable of index 1 or 2, or two variables each of index 1. Thus \(\mathcal K\) is an explicit finite combinatorial condition, with at most \(1+2|c(w)|+\binom{|c(w)|}{2}\) index tests. The zero test requires at least one right word with content contained in \(c(w)\).

Lemma 2. If \(\mathcal K(w, R)\) holds, \(w\leq R\) follows from (C),(N),(D),(S),(J).

*Proof.* Law (D) permits insertion of any letters into any nonlinear word, using commutativity. Consequently a normalized word is derivably below any nonlinear word whose multiplicities are componentwise no greater.

Call a variable saturated when its multiplicity is 3. If a normalized target has at least three occurrences of unsaturated variables, select three such occurrences and apply (S) in the remaining context. This replaces the target by the sum of three words, each with one selected occurrence repeated. Normalize and repeat. The sum of the deficits \(3-n_x\) over variables in the content strictly decreases along every branch. All descendants still satisfy \(\mathcal K\), since they are obtained by increasing multiplicities. It therefore suffices to handle a terminal word with at most two unsaturated occurrences.

Write \(s\) for the product of the cubes of its saturated variables. When there are no unsaturated occurrences, the zero index test provides a right word dividing \(s\). For a terminal \(sx\), assign \(x\) index 2 and all saturated variables index 0. A right witness has \(x\)-exponent at most 1 and divides \(sx\). For a terminal \(sx^2\), assign \(x\) index 1 to obtain a witness with exponent at most 2. In all cases outside variables are assigned index 3.

The remaining terminal is \(sxy\), \(x\ne y\). If a right word divides it, there is nothing to prove. Otherwise the assignment giving both \(x, y\) index 1 provides a right witness whose total exponent on \(x, y\) is at most 2. Absence of a divisor forces these exponents to be \((2,\allowbreak{}0)\) or \((0,\allowbreak{}2)\). Interchange \(x, y\) if necessary, and call this witness \(u\), with exponents \((2,\allowbreak{}0)\). Next give \(x\) index 2 and \(y\) index 0. A right witness \(v\) has \(x\)-exponent at most 1. Its \(y\)-exponent is at most 3 after normalization. Both witnesses have content in the target content. Hence (D) gives

\CaseDisplay{sx^2\leq u,\qquad sxy^3\leq v.}

Law (J), in context \(s\), yields \(sxy\leq sx^2+sxy^3\leq R\). This proves the lemma, including when \(s\) is omitted. \(\square\)

\par\smallskip\noindent\textit{3. An exact absorption criterion.}\quad 
Write a polynomial as \(P=L+R\), where \(L\) is its set of linear summands and \(R\) its set of nonlinear summands. Empty parts are omitted, not interpreted as named zero terms. For a nonlinear target \(w\), set \(C=c(w)\), \(T=L\cap C\), and \(R_C=\{u\in R:c(u)\subseteq C\}\). Let \(w_T\) be the normalized word obtained by making every variable in \(T\) have multiplicity 3.

Lemma 3. A linear word \(x\) is absorbed by \(P\) if and only if \(x\in L\). For a nonlinear word \(w\), the following criteria are exact.

If \(T=\varnothing\), in every algebra the criterion is \(\mathcal K(w, R)\).

If \(T\ne\varnothing\), the criteria are:

\CaseDisplay{\begin{array}{c|l}
B_1&R_C\ne\varnothing;\\
B_2&R_C\ne\varnothing,\text{ and for every }x\in C\setminus T
       \text{ with }n_x(w)=1,\\
   &\quad\text{some }u\in R_C\text{ has }n_x(u)\le1;\\
B_3&\mathcal K(w_T, R);\\
B_4&R\ne\varnothing.
\end{array}}

Every absorption satisfying its criterion is derivable from the corresponding displayed basis.

*Necessity.* If a linear \(x\) is absent, assign \(x=a\) and every other variable \(z\). Every right nonlinear word has value at most \(b<a\), so absorption fails. If \(T\) is empty, restrict to core evaluations and set outside variables to \(z\). Right linear summands then contribute only \(z\). This forces \(\mathcal K(w, R)\).

For \(B_1, B_2\), the necessity of \(R_C\ne\varnothing\) follows by assigning all variables in \(C\) to \(d\) and all outside variables to \(z\). If no contained nonlinear word exists, the target is \(e\) and the right side is \(d\). In \(B_2\), if the condition for a variable \(x\) of multiplicity 1 fails, assign this variable \(a\), all other variables in \(C\) to \(d\), and outside variables to \(z\). The target is \(a\). All right nonlinear words have index at least 2, and all right linear summands are \(d\) or \(z\). Their sum is \(d\), which is incomparable with \(a\).

For \(B_3\), suppose \(\mathcal K(w_T, R)\) fails. A separating core assignment can be taken with every outside variable of index 3 and with \(w_T\) of index 0, 1 or 2. Variables in \(T\) have index 0. If the target index is 0, every right word must contain an outside variable and hence has index 3. If the target index is 1, exactly one target variable has index 1 and multiplicity 1; each contained right word has multiplicity at least 2 in this variable. Raise its index to 2. Now the target has index 2 and every right nonlinear word has index 3. If the original target index is 2, the same right-side conclusion already holds. Finally change the values of variables in \(T\) to \(d\). Nonlinear evaluations do not change. Every right linear summand is now \(d\) or \(z\), so the full right side is \(d\). The target is \(e\) or \(b\), neither below \(d\). Thus absorption fails. This also covers an empty nonlinear part.

For \(B_4\), if the nonlinear part is empty, assigning every variable \(d\) gives target \(e\) and right side \(d\).

*Sufficiency and derivation.* Linear targets already occur. If \(T\) is empty, use Lemma 2.

In \(B_1\), select \(g\in T\) and write \(w=gq\) with \(q\) nonempty. Law (L) gives \(w\leq g+w^3\). For any \(u\in R_C\), normalization and (D) give \(w^3\leq u\). Hence \(w\leq P\).

For \(B_2, B_3\), repeated applications of (G), with a selected linear variable as \(x\) and the remaining word as \(y\), give

\CaseDisplay{w\leq \sum_{t\in T}t+w_T.}

The multiplicity of a selected variable increases by 2, and (N) normalizes it to 3. In \(B_3\), Lemma 2 now completes the proof.

For \(B_2\), fix \(g\in T\). Consider a word descended from \(w_T\), and suppose it has two unsaturated occurrences \(y, z\). Separate one occurrence of \(g\), writing it as \(gyzq\), with \(q\) possibly omitted. Use (H) with \(x=g\), \(y=y\), and \(z=zq\). Law (D) removes the extra copy of \(q\) in the last term, giving

\CaseDisplay{gyzq\leq g+gy^2zq+gyz^2q.}

When \(q\) is omitted this is simply (H). Repeat and normalize. Deficits strictly decrease, so eventually at most one unsaturated occurrence remains. A terminal with all exponents 3 has a divisor in \(R_C\). A terminal with one occurrence of a variable \(x\) has a divisor supplied by the criterion: \(x\notin T\), and it already had multiplicity 1 in the original target. Thus all terminal words are below \(R\), and the guard \(g\) is a right linear summand. This proves the required absorption.

Finally, in \(B_4\), choose \(g\in T\), write \(w=gq\), and choose any nonlinear right word \(u\). Law (B), substituting \(z=u\), gives \(w\leq g+u^2\). Since \(u\) is nonlinear, (D) gives \(u^2\leq u\). This completes the proof. \(\square\)

Distribute arbitrary terms into nonempty sums of words. For any identity valid in \(B_r\), each summand on one side is absorbed by the other side. Lemma 3 derives every such absorption. Applying this in both directions proves Theorem 1 for arbitrarily many variables.

\par\smallskip\noindent\textit{4. Strict inclusions, with separate completeness proofs.}\quad 
Corollary 4.

\CaseDisplay{\operatorname{Var}(B_4)\subsetneq\operatorname{Var}(B_1)
 \subsetneq\operatorname{Var}(B_2)\subsetneq\operatorname{Var}(B_3).}

Soundness of the larger variety's basis in the smaller generator proves each inclusion. Strictness is witnessed as follows: law (B) fails in \(B_1\) at \(x=y=d, z=a\); law (L) fails in \(B_2\) at \(x=d, y=a\); law (H) fails in \(B_3\) at \(x=d, y=z=a\). The respective target values are \(e, a, b\), while the right side is \(d\). Finite basability here comes from four separate completeness proofs, not from a containment transfer.

\subsection{A component-colouring criterion and three nonfinitely based semirings}\label{case:227}

\par\smallskip\noindent\emph{The hypergraph algebra.}
We use the construction of Jackson, Ren and Zhao \cite{JRZ}, also used by Gao, Jackson, Ren and Zhao \cite{GJRZ}. Write \(S_c(abc)\) for the eight-element algebra consisting of the seven nonempty subsets of \(\{1,2,3\}\) and an element \(\infty\). Distinct elements add to \(\infty\). Disjoint subsets multiply to their union. All other products are \(\infty\), and \(\infty\) absorbs both operations.

A Berge cycle of length \(j\ge2\) in a three-uniform hypergraph consists of distinct vertices \(v_1,\ldots,v_j\) and distinct edges \(e_1,\ldots,e_j\) with \(\{v_i,v_{i+1}\}\subseteq e_i\), reading the subscripts modulo \(j\). Its girth is the minimum such \(j\), or infinity when there is no cycle. In girth greater than two each pair of distinct vertices lies in at most one edge.

Let \(H\) be such a hypergraph with no isolated vertices. Give \(S_H\) one element \(a_v\) and one element \(b_v\) for each vertex \(v\), together with \(a\) and \(\infty\). The join of two distinct elements is \(\infty\). Multiplication is commutative. Its products different from \(\infty\) are
\[
       a_u a_v=b_w\quad\bigl(\{u,v,w\}\in E(H)\bigr),
       \qquad a_v b_v=a .
\]
These rules are unambiguous by the girth assumption. A triple product is different from \(\infty\) precisely on the three vertex elements of an edge, in which case its value is \(a\) in both bracketings. This proves associativity. A fixed factor acts injectively on the inputs for which the product is different from \(\infty\). For example \(a_u a_v=a_u a_{v'}=b_w\) forces \(v=v'\) by uniqueness of the edge containing \(\{u,w\}\). Hence multiplying two unequal summands gives either an infinite value or two unequal finite values. This proves distributivity over the flat addition.

\begin{lemma}\label{lem:forest-local}
Let \(K\) be a finite three-uniform hypergraph with an edge and without isolated vertices. If \(K\) has no Berge cycle, then \(S_K\in\V(S_c(abc))\). More generally, if \(n\ge2\) and \(H\) is a finite three-uniform hypergraph with an edge, without isolated vertices and of girth greater than \(\binom n2+n+1\), every subalgebra of \(S_H\) generated by at most \(n\) elements belongs to the same variety.
\end{lemma}
\begin{proof}
Call a colouring \(f\colon V(K)\to\{1,2,3\}\) rainbow if every edge has three distinct colours. The incidence graph, with an edge-node adjacent to each of its three vertices, is a forest. A component can be coloured outward from any vertex. Each new edge has just one previously coloured vertex and its other two vertices receive the two remaining colours. Distinct vertices can always be coloured differently. Two distinct vertices which share no edge can also be coloured alike. For vertices in different components the choices are independent. Within one component, first colour the vertices along the unique incidence path. Consecutive vertex-nodes on this path must receive different colours. At distance two in this sequence any ordered pair of endpoint colours is possible, by choosing an intermediate colour different from both. Induction gives the same assertion at every greater distance. Fill in the unused third vertex of each path edge and then extend off the path as above.

Let \(\mathcal C\) be the set of all rainbow colourings and take the algebra \(S_c(abc)^{\mathcal C}\). All vectors having \(\infty\) in at least one coordinate form an ideal for both operations. Identifying these vectors to one element and leaving every other vector in a singleton class gives a congruence. In this quotient represent \(a_v\) by the vector of singletons \(\{f(v)\}\), \(b_v\) by the vector of complements \(\{1,2,3\}\setminus\{f(v)\}\), and \(a\) by the constant full set. Within either degree the vectors for distinct vertices differ in some coordinate. Across degrees their cardinalities differ. Thus unequal represented elements add to the collapsed class. When \(u,v,w\) form an edge, the product for \(u,v\) represents \(b_w\). For two nonadjacent vertices take a colouring giving them the same colour, so their product is in the ideal. A repeated singleton has the same property. The product representing \(a_u b_v\) yields the full set if \(u=v\). For \(u\ne v\), choose a colouring with \(f(u)\ne f(v)\), so the two subsets intersect in this coordinate. All products of total degree greater than three also belong to the ideal. These verifications embed \(S_K\) in the quotient.

For the second assertion let \(Y\) contain the subscripts of the generators of the form \(a_v\) or \(b_v\), so \(|Y|\le n\). Retain all edges of \(H\) containing at least two members of \(Y\). Add one incident edge for each member of \(Y\) still uncovered, and take one edge if the list is empty. The subhypergraph \(K\) on these edges has at most \(\binom n2+n+1\) edges and therefore has no cycle. The subalgebra generated in \(S_H\) has no vertex element beyond its generators. Its pair elements \(b_w\) either are generators or arise by multiplying two generating vertex elements. In the latter case the unique completing edge was retained, even when \(w\notin Y\). Accordingly the same names define an injective homomorphism from this subalgebra into \(S_K\). Products of two vertex elements agree by the choice of edges, products with pair elements depend only on equality of subscripts, and addition is flat in both algebras. The first assertion now applies.
\end{proof}

\begin{lemma}\label{lem:high-girth-exists}
For every integer \(L\ge2\) there is a finite connected three-uniform hypergraph of girth greater than \(L\) admitting no two-colouring in which every edge has both colours.
\end{lemma}
\begin{proof}
On \(N\) vertices include each triple independently with probability \(p=N^{-2+\epsilon}\), where \(0<\epsilon<1/L\). A Berge cycle of length \(j\) can be specified by at most \(N^{2j}\) choices, one vertex on the cycle and one further vertex of its edge for each step. As its \(j\) edges are distinct, the expected number of cycles of lengths two through \(L\) is at most \(\sum_{j=2}^L N^{\epsilon j}=o(N)\). Put \(m=\lfloor N/3\rfloor\). The probability of an independent \(m\)-set is at most
\[
       2^N(1-p)^{\binom m3}\ \le\ 2^N\exp\!\left(-p\binom m3\right),
\]
which tends to zero. With positive probability there is no independent \(m\)-set and there are fewer than \(N/6\) short cycles, by Markov's inequality. Delete one vertex from each short cycle. In the remaining induced hypergraph a colour class containing at least half the vertices has size at least \(5N/12>m\), so it contains an edge. Hence some connected component already has the required property, since proper colourings of the components would combine.
\end{proof}

\par\smallskip\noindent\textit{2. The finite certificate and its consequence.}\quad 
Write \(u\le v\) when \(u+v=v\). For an ai-semiring \(A\), put

\CaseDisplay{\tau(x, y, z)=xyz+xzy+yxz+yzx+zxy+zyx.}

For \(b\in A\), let

\CaseDisplay{R_b=\{(r, s, t)\in A^3:\tau(r, s, t)\le b\}.}

Its components are the connected components of the undirected graph joining any two entries of a tuple in \(R_b\); isolated elements are also components. Repeated entries in tuples are retained. Call \(r\) bad at \(b\) if \(r^3\not\le b\).

Theorem. Suppose that \(S_c(abc)\in\mathsf V(A)\). Suppose also that, for every \(b\in A\), each component \(C\) of \(R_b\) containing a bad element admits a map \(\chi_C:C\to\{0,\allowbreak{}1\}\) such that every tuple of \(R_b\) in \(C\) has both colours among its entries. Then \(A\) is nonfinitely based.

Proof. Let \(H\) be a connected non-2-colourable 3-uniform hypergraph. Define

\CaseDisplay{T_H=\sum_{\{u, v, w\}\in E(H)}\tau(x_u, x_v, x_w).}

We first show, for every vertex \(v\), that

\CaseDisplay{A\models T_H\approx T_H+x_v^3.\quad\textup{(\ensuremath{1})}}

Take an arbitrary assignment \(\alpha:V(H)\to A\) and set \(b=T_H(\alpha)\). Each edge tuple lies in \(R_b\). As \(H\) is connected, all images of its vertices lie in a single component \(C\) of \(R_b\). If \(\alpha(v)^3\not\le b\), then \(C\) has a colouring \(\chi_C\) as in the hypothesis. The composite \(\chi_C\circ\alpha\) gives a nonmonochromatic colouring of every edge of \(H\), which is impossible. Thus \(\alpha(v)^3\le b\). No bound has been imposed on the number of vertices or variables.

For each \(n\ge2\), use Lemma~\ref{lem:high-girth-exists} to choose \(H_n\) of girth greater than \(\binom n2+n+1\). All subsemirings of \(S_{H_n}\) generated by at most \(n\) elements lie in \(\mathsf V(S_c(abc))\subseteq\mathsf V(A)\) by Lemma~\ref{lem:forest-local}. However, under \(x_v\mapsto a_v\), the term \(T_{H_n}\) has value \(a\), whereas \(x_v^3\) has value \(\infty\). Hence (1) fails in \(S_{H_n}\), so \(S_{H_n}\notin\mathsf V(A)\).

If \(\mathsf V(A)\) had a finite identity basis, choose \(n\ge2\) at least its maximum number of variables. Every evaluation of a basis identity in \(S_{H_n}\) is contained in an at most \(n\)-generated subsemiring. That subsemiring lies in \(\mathsf V(A)\), so every basis identity holds in \(S_{H_n}\). This would force \(S_{H_n}\in\mathsf V(A)\), a contradiction. \(\square\)

Symmetrizing each edge term in \(T_H\) covers every order of multiplication in \(A\). Under the vertex assignment in \(S_H\), every ordering still gives \(a\).

\par\smallskip\noindent\textit{3. Three applications.}\quad 
The following tables use carrier \(\{0,\allowbreak{}1,\allowbreak{}2,\allowbreak{}3,\allowbreak{}4\}\). These numerals are element labels, not term constants. Each of the three algebras is commutative. In each case 4 is the greatest additive element and a multiplicative zero. Tables are written as five rows in the order \(0,\allowbreak{}1,\allowbreak{}2,\allowbreak{}3,\allowbreak{}4\).

\noindent Catalogue ID\enspace\textbar\enspace Addition rows\enspace\textbar\enspace Multiplication rows\enspace\textbar\enspace Cubes of 0,\allowbreak{}1,\allowbreak{}2,\allowbreak{}3,\allowbreak{}4\par
\noindent 10081\enspace\textbar\enspace 04444; 41444; 44244; 44434; 44444\enspace\textbar\enspace 12344; 23444; 34444; 44444; 44444\enspace\textbar\enspace 2,\allowbreak{}4,\allowbreak{}4,\allowbreak{}4,\allowbreak{}4\par
\noindent 12214\enspace\textbar\enspace 00444; 01444; 44224; 44234; 44444\enspace\textbar\enspace 42444; 23404; 44444; 40424; 44444\enspace\textbar\enspace 4,\allowbreak{}0,\allowbreak{}4,\allowbreak{}4,\allowbreak{}4\par
\noindent 14840\enspace\textbar\enspace 00444; 01444; 44244; 44434; 44444\enspace\textbar\enspace 44444; 42304; 43044; 40444; 44444\enspace\textbar\enspace 4,\allowbreak{}3,\allowbreak{}4,\allowbreak{}4,\allowbreak{}4\par
For 10081, at bound 2 the only tuple in \(R_2\) is \((0,\allowbreak{}0,\allowbreak{}0)\), and its component contains no bad element. At bound 3 the only tuples are the three permutations of \((0,\allowbreak{}0,\allowbreak{}1)\); the component \(\{0,\allowbreak{}1\}\) is coloured by \(0\mapsto0\), \(1\mapsto1\). At bounds 0 and 1 the relation is empty. At bound 4 every element is good. All remaining components with bad elements are isolated and have no tuples. This verifies the colouring condition.

For 12214, at bound 0 the only tuple is \((1,\allowbreak{}1,\allowbreak{}1)\), and its component contains no bad element. At bound 2 the only tuples are the permutations of \((1,\allowbreak{}1,\allowbreak{}3)\); colour \(1\mapsto0\), \(3\mapsto1\). At bounds 1 and 3 the relation is empty; at bound 4 every element is good.

For 14840, at bound 3 the only tuple is \((1,\allowbreak{}1,\allowbreak{}1)\), and its component contains no bad element. At bound 0 the only tuples are the permutations of \((1,\allowbreak{}1,\allowbreak{}2)\); colour \(1\mapsto0\), \(2\mapsto1\). At bounds 1 and 2 the relation is empty; at bound 4 every element is good.

The remaining hypothesis \(S_c(abc)\in\mathsf V(A)\) follows by the transition and principal-ideal tests in Proposition~\ref{prop:lower-fibres} for each of the three source tables. The target is encoded on \(\{0,\allowbreak{}1,\ldots,7\}\), with 0 its zero and nonzero elements the nonempty subsets of a three-element set (binary masks). The sum of distinct elements is 0; the product of disjoint nonempty masks is their union, and every other product is 0. Target generators are 1,\allowbreak{}2,\allowbreak{}4. These target labels are unrelated to the labels in the three source tables.

For each source, all word functions on three generators close at 27 functions, with 81 right-generator transitions, evaluated on all 125 source assignments. The propagated target values respect every transition. For each of the eight target principal lower sets, the join of all selected source word functions absorbs no word whose target value is outside that lower set.

Here is why this finite calculation proves full membership rather than a bounded word test. Closure under every right-generator transition exhausts all word functions by induction on word length. The target value of a word function is well-defined by transition compatibility. Every semiring term is a finite join of these functions. If \(f_u\) is below a join of source word functions whose target values have join \(b\), then it is below the join of \emph{all} word functions with target value below \(b\). The principal-lower-set check forces the target value of \(u\) below \(b\). Applying this in both directions to any equality of joins shows that target evaluation is well-defined on the subsemiring of \(A^{125}\) generated by the three projection functions. It preserves both operations and is onto \(S_c(abc)\). This is an explicit HSP witness, with no word-length cutoff.

Corollary. Catalogue types 10081, 12214 and 14840 are nonfinitely based. Their previously known proper nonfinitely based subvarieties only excluded limit-variety generation; the theorem above is the new argument that proves the algebras themselves nonfinitely based.

For \(\A{14839}\) the membership test for \(S_c(abc)\) fails. Section~\ref{case:331} gives its nonfinite basis proof using a different family of hypergraph identities.

\subsection{Nonfinite basis from odd cycles with a disjoint cube}\label{case:228}
This proof uses the signature (+, multiplication), without constants. Let S be catalogue 10107, with addition 00444/\allowbreak{}01234/\allowbreak{}42244/\allowbreak{}43434/\allowbreak{}44444 and multiplication 21344/\allowbreak{}11111/\allowbreak{}31444/\allowbreak{}41444/\allowbreak{}41444. Addition has bottom 1, incomparable atoms 0,\allowbreak{}2,\allowbreak{}3, and top 4. Multiplication is commutative;1 is absorbing. We prove S nonfinitely based and give a sufficient transfer theorem.

\par\smallskip\noindent\textit{Values and localization.}\quad 
A word containing a letter assigned 1 has value 1. Otherwise, a literal has its assigned value; a quadratic has value 2 for the pair(0,\allowbreak{}0), value 3 for(0,\allowbreak{}2) or(2,\allowbreak{}0), and value 4 otherwise. A cubic has value 3 exactly when all its letters are 0, and value 4 otherwise. Every word of length at least four has value 4. These assertions follow directly at lengths two and three; multiplying3 by any non1 element gives 4, and 4 remains 4 under multiplication by non1 elements. This is an induction at all lengths.

For a target word w, discard every original right summand whose content is not contained in \(U=c(w)\), obtaining \(P_{U}\). If \(w\leq P\) is valid, \(P_{U}\) is nonempty and \(w\leq P_{U}\) is valid: assign every outside variable 1. Conversely \(w\leq P_{U}\) implies \(w\leq P\). This argument applies to literals as well as longer words. It does not permit discarding variables from a derived target.

For a target of length at least four, its exact criterion is as follows. Within \(P_{U}\) let L be the literal variables, Q the original quadratic multigraph, and T the set of variables occurring in original cubics. Then \(w\leq P\) holds exactly when at least one of these conditions holds:

1. \(P_{U}\) has an original word of length at least four.

2. L meets the variables of an original quadratic or cubic.

3. There is at least one original cubic and no two-coloring of Q that gives color 0 to every vertex in T.

In the third condition a loop obstructs two-coloring. Isolated variables of T have no effect. The condition includes an odd component anywhere in Q, or an odd path between two forced vertices of T.

To prove this, restrict to assignments with no 1 in U; the target is 4. A contained long right word is 4. At a literal/nonlinear overlap, the literal's value differs from any possible non4 value of that nonlinear word, so their join is 4. If a cubic is present and the right sum is not 4, all cubics must be 3, all their variables 0, every quadratic must be 3, and every literal must be 3. Quadratics of value 3 are exactly edges joining0 to 2, giving the required coloring. Conversely such a coloring, with cubic variables 0 and disjoint literal variables 3, makes the right sum 3. If no cubic and no overlap occurs, assign every quadratic vertex0 and every literal variable 2; the right sum is 2. With no nonlinear term, assign all literal variables 2. These are separators in every failed case. When w contains a1 its value is bottom, proving sufficiency globally.

For short targets, localization leaves at most three variables. The word formula above is an exact finite decision on that local polynomial; we will only need the explicit special case below.

\par\smallskip\noindent\textit{Critical identities and their sparse absorptions.}\quad 
For odd \(n\geq 5\), take distinct \(x_{0}\),..., \(x_{n-1}\), z and put

{\small\par\smallskip
\CaseLine{\(P_{n} = sum_{i mod n} x_{i} x_{i+1} + zzz,\)}
\CaseLine{\(r_{n} = x_{0} x_{1} \cdots  x_{n-1} z.\)}
\smallskip}

S satisfies \(r_{n}\leq P_{n}\). Without a1 assignment, zzz must be 3 if the right sum is not 4, requiring all cycle edges to be 3 and therefore an impossible two-coloring of the odd cycle. With a1 anywhere on the left, the target is bottom. The same argument proves the identity for \(n=3\).

The only words of length at most three absorbed by \(P_{n}\) are its original edges, up to commutation, and zzz. Here is a direct proof. A missing literal or square has no local edge; its local right side is empty or zzz, and assigning its variable 0 separates it. A missing mixed quadratic also has no local edge; assigning its letters 0 separates value 2 from either bottom or the possible cube value 3. A cubic not containing z is separated by assigning all its variables 0: the target is 3 and local edges are 2. For a mixed cubic containing z, keep \(z=0\). If its other two letters form an original edge, assign them0 and 2; otherwise assign one non-z letter 2 and all other letters 0. The target becomes 4 and the local right sum is 3. The pure cube zzz is already original.

Every absorbed word of length at least four contains ALL n cycle variables AND z. More strongly, if R is any subset of \(P_{n}\) and a long word t is absorbed by R, then \(R=P_{n}\) and c(t) contains c(\(P_{n}\)). Indeed R has no literals or long words. Its only possible cube is zzz, disjoint from every edge. A proper subgraph of the chordless n-cycle is bipartite. The long-word criterion therefore requires the cube and all cycle edges, all contained in t.

\par\smallskip\noindent\textit{No bounded-variable identity set derives every critical identity.}\quad 
Fix k and choose odd n\ensuremath{>}max(k,9). Suppose identities of S using at most k variables derived \(r_{n}\leq P_{n}\). Work with distributively expanded nonempty polynomials, counting the semiring axioms among the available identities. Along the derivation from \(P_{n}\) to \(P_{n}\)+\(r_{n}\) there is a first polynomial containing a word of length at least four. Every earlier polynomial is identical to \(P_{n}\) in S, so the sparse-absorption lemma makes every one of its words an original cycle edge or zzz, with possible commutations.

Choose the first introduced long word t. It comes from a source word w on one side of a valid identity, a substitution sigma of nonempty polynomials for its variables, and a particular monomial context a, b. Write R for the other polynomial side of that source identity. Every expanded word of a sigma(R) b is an original short word of \(P_{n}\), and \(t\leq a\) sigma(R) b is valid in S. The stronger sparse lemma forces this right polynomial to contain every word of \(P_{n}\). If a or b were nonempty, all these words would share one context variable. The disjoint cube and cycle edges share no variable. Thus both contexts are empty and

\(\sigma (R)=P_{n}\)                 (as sets of commutative short words).

All source words in R have length at most three, and every word in the image of a variable occurring in R has length at most three. Moreover an original literal variable of R cannot also occur in a nonlinear word of R. Its image, being a summand of \(P_{n}\), would be a cycle edge (length two, with cycle letters) or zzz. Any further nonempty factor would then yield length at least four or a mixed cubic, neither of which occurs in \(P_{n}\). This rules out the overlap obstruction in R.

Every source cubic in R has all its variables mapped to the singleton polynomial z: all expanded cubic words must be zzz. Classify each edge uv of the source quadratic graph. If some product of image words is a cycle edge, one image choice for each endpoint is a cycle letter. Holding either of these choices fixed shows that EVERY image word of the other endpoint is a singleton cycle letter. Their products are all cycle edges. If no product is a cycle edge, every product is zzz; all images are pure powers of z of a fixed length \(d_{u}\) or \(d_{v}\) in \{1,\allowbreak{}2\}, and \(d_{u}+d_{v}=3\). This classification propagates along connected edges: a component has exclusively the cycle-letter kind or the pure-z kind. This explicitly treats arbitrary polynomial images and repeated source variables, including source loops.

A pure-z component is bipartite by its image lengths1 and 2; vertices appearing in source cubics have length 1 and can all receive color 0. In a cycle-letter component, choosing ONE image letter at each source vertex gives a graph homomorphism to the chordless n-cycle. Any nonbipartite graph on at most k source vertices has an odd cycle of length at most k. Its image would be an odd closed walk of \(length<n\) in \(C_{n}\), impossible because the odd girth of \(C_{n}\) is n. Thus this component is also bipartite, and it contains no cubic variable.

We can therefore assign every source quadratic edge value 3 by coloring its endpoints 0 and 2, with every source cubic variable 0. Assign the disjoint source literal variables 3. Every summand of R is then 3, and no assigned variable is 1. Consequently a valid \(w\leq R\) cannot have length at least four. Nor can w contain a variable absent from R: assign such a variable 4, retaining the preceding non1 assignment on R, to separate w from R. Hence w has length at most three and all its variables occur in R. Each image word has length at most three, so every word of sigma(w), including t, has length at most nine. But the sparse lemma requires t to contain n+1 distinct variables. This contradiction proves that no k-variable identity set derives the chosen critical identity. Every finite basis has a finite variable bound. Therefore S is NFB.

The expansion argument also covers arbitrary addition/multiplication contexts: selecting the monomial branch producing t selects fixed monomial contexts for the other side, all present before the step. No restriction to monomial substitutions has been made.

\par\smallskip\noindent\textit{Transfer to other finite semirings.}\quad 
Let A be any semiring for which S belongs to HSP(A) and A satisfies \(r_{n}\leq P_{n}\) for every odd \(n\geq 3\). Then A is NFB. Otherwise a finite basis of A, valid in S because of the HSP inclusion, would derive a critical identity with n larger than its variable bound, contradicting the proof. Both hypotheses are required. An NFB subvariety alone is insufficient.

For a finite A, the all-length hypothesis has an exact finite-state verification, even when multiplication is not commutative. Starting at length two, retain (first, last, ordered product, join of path edges, parity). Appending c replaces last by c, product by product*c, edge join by edge join+last*c, and flips parity. Close the finite reachable state set. At every odd state, for every value of z, check

{\small\par\smallskip
\CaseLine{\(p z\leq h+\ell f+z^3,\qquad f=\text{first},\quad \ell=\text{last},\quad p=\text{ordered product},\quad h=\text{join of path edges}.\)}
\smallskip}

Induction on path length gives exactly all assignments at every length; finite closure and transition verification establish all odd lengths, not a bounded numerical sample. For five elements there are at most 2*\(5^{4}\) states. Complete HSP certificates independently verify the other hypothesis for each transferred type.

The tempting rule \(xw\leq xy+yz+zw\) is INVALID here: \(x=w=0\), \(y=z=1\) gives left 2 and right 1. Removing interior variables destroys bottom protection. This is precisely why an unguarded finite path proof cannot replace the odd-cycle obstruction.

\subsection{Six complete identities for 10108}\label{case:229}
In the signature (+, multiplication) without constants, a complete basis relative to the additively idempotent semiring axioms is

{\small\par\smallskip
\CaseLine{\(xy=yx; xyz=xyzz; xyzt\leq xyz; xy\leq xxy; x\leq xxx; xxy\leq x+xy.\)}
\smallskip}

Here \(u\leq P\) means \(u+P=P\). Addition is 00444/\allowbreak{}01234/\allowbreak{}42244/\allowbreak{}43434/\allowbreak{}44444: bottom 1, three incomparable atoms 0,\allowbreak{}2,\allowbreak{}3, and top 4. Multiplication is 21444/\allowbreak{}11111/\allowbreak{}41444/\allowbreak{}41424/\allowbreak{}41444. Literals have their assigned values. A nonlinear word containing a variable assigned 1 has value 1. Otherwise a quadratic with both occurrences 0 or both 3 has value 2, and every other nonlinear word has value 4. Direct induction proves this at every length; long words depend only on content.

For original P and target w, let \(P_{U}\) consist of original summands contained in \(U=c(w)\). Let L be its literal variables and V the endpoints of its quadratics. Then a missing literal is absorbed exactly when P has an original pure long word on its variable. An original quadratic is absorbed immediately. Every other nonlinear target is absorbed exactly when \(P_{U}\) contains an original long word or L intersects V. Original literals are, of course, absorbed.

For sufficiency, if the nonlinear target has a variable assigned 1, it is bottom. Otherwise a contained original long word is 4. In the overlap case retain an original literal x and quadratic xy: if xy is 4 it already dominates; if xy is 2, x is 0 or 3 and \(x+xy=4\). Thus either guard supplies 4. A literal is below every pure long power, directly from the scalar values.

For necessity, a missing literal with no pure long word is separated by assigning its variable 0 and all outsiders 1. Its right values are 1 or 2, neither dominating0. For a long target lacking either guard, assign quadratic vertices V to 0, other U variables 2, and outsiders 1. Since L is disjoint from V, all surviving short right terms are 2, while the target is 4.

For a missing quadratic, absence of a contained long word and overlap leaves only endpoint literals and endpoint loops: a mixed quadratic would be the original target. If both endpoint loops occur, assign the endpoints 0 and 3 respectively. If only one loop occurs, assign its endpoint 0 and the other 2. If neither occurs, assign both 2. The target is 4 and every surviving right term is 2; overlap absence ensures that no endpoint literal has been assigned 0 or 3. For a missing square, assign its variable 2. In all cases set outsiders 1, so noncontained terms are bottom. These are unrestricted separators.

For completeness, commutativity and \(xyz=xyzz\) normalize long words by content: move a variable to the end, duplicate it with two nonempty preceding blocks, or delete extra copies in reverse while preserving length at least three. Canonical forms are sorted content for at least three variables, xxy for two, and xxx for one. Quadratics retain multiplicity, and literals retain length one.

For a new nonlinear target retain an original contained long word. If only the overlap guard is available, derive xxy from the original contained literal x and original contained quadratic xy using \(xxy\leq x+xy\). Append the nonempty word of missing target variables, if any, by \(xyzt\leq xyz\), splitting the long seed into a first nonempty block and its last two letters. Normalize to the target long form; if the target is quadratic, recover it by \(xy\leq xxy\). For a missing literal, normalize its original pure long witness to xxx and use \(x\leq xxx\). All substitutions are nonempty and all witnesses are original or explicitly derived.

There is no bound on variables, words, or polynomial sizes. Distributive expansion and absorption in both directions prove completeness. Each law is valid by the criterion and checked on all assignments; numerical checks supplement the proof.

The boundaries are essential: \(x\leq xx\) fails at \(x=0\); \(xy\leq xx+yy\) fails at \(x=0\), \(y=3\); \(xxx\leq x\) fails at \(x=2\). An outside long guard is insufficient: \(xxx\leq xx+yyy\) fails at \(x=0\), \(y=1\). Thus this proof does not import the preceding component-activation criteria.

\subsection{Finite bases for 10139, 11229, 11234 and their opposites}\label{case:230}
Use the constant-free semiring signature. Fix a target word w, with support C, and a nonempty upper polynomial P. Let U be its summands with support contained in C and T the union of their noninitial letters. All three criteria for w \ensuremath{\leq} P require the supports of U to cover C. Call x simple in w when it occurs exactly once.

\textbullet{} For 11229 the noninitial letters of w must belong to T and every simple x must occur exactly once in some member of U.

\textbullet{} For 11234 require the same noninitial coverage and, for each simple x, a member of U containing it at most once.

\textbullet{} For 10139 every simple x must occur once in a member of U. If it occurs noninitially in w, additionally require a member of U omitting x or a member containing x once noninitially. The omitting and once-only witnesses can be different words.

These conjunctions are necessary and sufficient on arbitrary alphabets.

\par\smallskip\noindent\textit{Word valuations.}\quad 
In 11229 bottom 2 annihilates multiplication. Without 2, exactly one 3 gives 3 and at least two give 2. In 11234 the absorbing bottom is 3. Without 3, one 2 gives 2 and at least two give 3. In both tables the special value is preserved by multiplication with values in \{0,\allowbreak{}1,\allowbreak{}4\} on either side, and its square is bottom. A word on \{0,\allowbreak{}1,\allowbreak{}4\} has value 1 on all ones, 0 on an initial 0 followed by ones, and 4 otherwise. The row of 1 on this alphabet gives 4 on 0 or 4, and 1 on 1. The row of 0 gives 0 on 1 and 4 otherwise. This proves the formulas by appending letters and includes the one-letter cases.

11229 has additive chains \(2<1<0<4\) and \(2<3<4\), with 3 incomparable with 1,\allowbreak{}0. In 11234 the chain is \(3<2<1<0<4\). For counterassignments give outside letters bottom. In each source the cover is necessary: give one missing letter 4 and every other target letter 1. For noninitial coverage put 0 at one missing noninitial letter, 1 at the rest of C. The target is 4 and the right bound is 0. In 11229 put a single target letter x at 3 and all other target letters at 1. If no contained word contains x exactly once the right sum is at most 1 and does not absorb the target 3. In 11234 put \(x=2\) and the other target letters at 1: absence of an at-most-once witness makes every upper word bottom. Conversely when the target is special its indicated witness supplies its value, or in 11234 an omitted special gives a higher value. On the high alphabet a target 0 is supplied by a cover witness for its 0; a target 4 with a 4 letter is witnessed by coverage, and one with a noninitial 0 by noninitial coverage. For an all-ones target any contained word suffices. Bottom targets require no condition.

In 10139, 2 is absorbing bottom. With no 2, at least two entries from \{0,\allowbreak{}1\} give 2. A unique such entry gives 0, except that an initial 1 gives 1. In the absence of 0,\allowbreak{}1,\allowbreak{}2 the word is the maximum of its entries in \{3,\allowbreak{}4\}. Indeed a high prefix followed by 0 or 1 gives 0, either special value survives a high suffix, and a second special multiplies to 2. The additive order is \(2<1<0<4\) and \(2<3<4\) with 3 incomparable with 1,\allowbreak{}0. Set outsiders to 2. Cover is forced by giving a missing variable 4 and others 3. With x simple, \(x=0\) and the others 3 exclude failure of exact occurrence, since 0 is not below 3. If x is noninitial and the additional condition fails, set \(x=1\) and the others 3. Every contained word with a single x is headed by it, none omits x, and all right values are at most 1, below target 0. Conversely a target on \{3,\allowbreak{}4\} is witnessed by coverage. A target with a unique 0 has a once-witness of value 0. A unique initial 1 has a once-witness of value 0 or 1. For a unique noninitial 1 one needs value 0. A noninitial once-witness suffices. In the other case the once-witness gives 1 and a separate omitting witness gives 3 or 4. Their join is 4. This exhausts the possibilities.

\par\smallskip\noindent\textit{Unrestricted completeness.}\quad 
Use the semiring axioms, \(xyz=xzy\), \(xxx=xx\), \(xyz=xyyz+xyzz\), and each source's valid absorptions on 247 variables, length at most 494 and at most four right words. Direct substitution in the table verifies the three laws. The first two normalize a word by preserving its head and sorting and capping total multiplicities at two. For at most five words reserve their heads, then combine other variables with the same nonzero vector of 0,\allowbreak{}1,\allowbreak{}2 multiplicities. Use their common multiplicity in each compact word. Heads, noninitial roles, support inclusions, omissions and exact once-witnesses survive this construction. Every unreserved occurrence is noninitial in the corresponding word. The proof of reversible compression in Section~\ref{case:173} then applies: substitution of the nonempty product of each group recovers all words. Hence every valid local step with at most four upper words is derivable.

Call a word available if proven below P. In 10139 one can first double every letter in each contained original, retaining its head. Merge these saturated words two at a time on the union of their supports. They have no simple letters and coverage suffices at every step. This produces an available saturated word on C.

For 11229 and 11234 keep the support S and noninitial set D of the running word. A word realizing S, D can saturate every variable except the unique possible letter of S minus D. When present this exception must be a simple head; put it first and list the other letters twice. This includes the literal case. When there is no exception list every letter twice. At each step below, the chosen upper words cover S and their noninitial sets jointly cover D. Any exceptional head occurs exactly once in each containing operand, since it is absent from their noninitial union. Thus both simple-variable criteria are satisfied.

Start separately with each original u in U, applying this realization to its own sets. If the uncovered letter d in S minus D is noninitial in the ultimate target, choose v in U containing d noninitially. Merge the two supports and the two noninitial sets. In their unions at most the head of v can still be missing from the latter, and d is gained, so the step just described applies. Repeat until the exception is absent or is outside the target noninitial set. Each iteration strictly increases D inside C. Merge all resulting words two at a time by taking both unions. The remaining missing letters belong to C minus the target noninitial set, which has size at most one. The final available word covers C and all target noninitial letters.

Finally split the target at two distinct simple nonhead variables if present, using tail permutations and \(xyz=xyyz+xyzz\), with a nonempty prefix for x. The new targets have fewer simple nonhead variables, the same support and head, and the same noninitial set. No remaining simple-variable requirement is lost. At a leaf with at most one simple nonhead and possibly a simple head, use the available full-support word and the required original witnesses. For 10139 at most one witness is needed for the head and at most two for the other letter. For the other types at most one witness per simple variable suffices. This gives at most four available upper words. The local derivation and reassembly of split terms prove the target. Every identity follows by distributing terms and deriving each summand absorption both ways. The opposites of 10139,\allowbreak{}11229,\allowbreak{}11234 are 10128,\allowbreak{}13089,\allowbreak{}13094 respectively.

\subsection{Six identities for 10133: a literal head and contained tail witness}\label{case:231}
Signature (+, multiplication), no constants, relative to the additively idempotent semiring axioms. Unrestricted completeness.

Addition rows:00044,\allowbreak{}01144,\allowbreak{}01234,\allowbreak{}44334,\allowbreak{}44444, with order \(2<1<0<4\) and \(2<3<4\). Multiplication rows:22200,\allowbreak{}22200,\allowbreak{}22222,\allowbreak{}22233,\allowbreak{}22244. Literals retain their values. In nonlinear hT, any tail outside\{3,\allowbreak{}4\} gives 2. If all tails belong to\{3,\allowbreak{}4\}, head 0/\allowbreak{}1 gives 0, head 2 gives 2, head 3 gives 3, and head 4 gives 4. Table induction proves this formula for every word length. W1 \(xyz=xzy\) and W2 \(xyy=xy\) follow, normalizing to a head and sorted tail set.

\par\smallskip\noindent\textit{Exact criterion.}\quad 
A literal target x is absorbed by original P iff original literal x occurs. A nonlinear hT is absorbed exactly when either: (A) an original NONLINEAR h-word has tail contained in T; or (B) original literal h occurs, and either h belongs to T or an original word has its entire content contained in T.

For a missing literal x, assign \(x=1\), others 2. All nonlinear words give 2, and without literal x the right sum 2 cannot absorb 1. If A fails and literal h is absent, there are two cases. For h outside T, assign \(h=1\), \(T=3\), others 2. Target 0; every original nonlinear h-word has a low tail outside T (including h), so gives 2. Other original heads are 3 or 2; their word values are 3 or 2, whose join cannot absorb 0. For h in T, assign \(h=4\), other variables in T to 3, and outsiders 2. Target 4; failed A h-words have an outsider 2 tail and give 2. All other right heads are 3 or 2, so the right sum is at most \(3<4\). These separate both cases.

If literal h exists but B fails, then h is outside T and no original word is entirely in T. Assign \(h=1\), \(T=3\), others 2. The target is 0; right nonlinear h-words give 2 because A fails. Other heads outside T are 2; heads in T are 3 but their words must contain an outside low factor and give 2. Literal h is 1; there are no literals in T. The right \(\sum \leq 1<0\). These give necessity.

Sufficiency and soundness can also be read directly from the word formula. Under A, if the target is not 2, all target tails are high, so the contained same-head nonlinear word has the same value. Under B with h in T, the head is high whenever the target is not 2, so literal h has the target value. Under B with a contained tail word q, all q variables are high and q gives 3 or 4. For target 0, the literal head 0/\allowbreak{}1 joined with 3/\allowbreak{}4 gives 4; for head 3/\allowbreak{}4 the literal head already supplies the target. Value 2 is least.

\par\smallskip\noindent\textit{Six laws and unrestricted construction.}\quad 
{\small\par\smallskip
\CaseLine{\(W1 xyz=xzy\qquad W2 xyy=xy\)}
\CaseLine{\(L  xyz\leq xy\qquad Q  xx\leq x\)}
\CaseLine{\(A  xy\leq x+y\qquad B  xyz\leq x+yz\)}
\smallskip}

L appends a tail restriction and so can only preserve value or give 2. Q sends \(0\to 2\),\(1\to 2\),\(2\to 2\),\(3\to 3\),\(4\to 4\), all downward in the additive order. A and B satisfy the exact criterion, using literal x and contained tail word y or yz. All six laws are sound.

Normalize original words and target using W1, W2. Under criterion A, choose its original nonlinear h-word as seed \(g=hU\). Under B with repeated head, Q derives hh from original literal h. Under B with a tail-contained q, A derives hq from h+q if q is literal, and B derives hq from h+q if q is nonlinear, split nonemptily into its head and remaining tail. In every case g is nonlinear, has head h, and tail contained in T. If any target tail letters are missing, collect them in a nonempty word V and use L with \(x=h\), \(y=the\) nonempty seed tail, \(z=V\). Normalize to hT. If no letters are missing, the seed already is the target. Every substitution is nonempty. Literal targets are original summands.

Thus all valid absorptions derive from six fixed laws for arbitrarily many variables, word lengths and right summands. Distributivity and both summandwise directions give every identity. Reversal gives the opposite basis. Catalogue classification additionally requires independent complete HSP both ways. Finite checks supplement the unrestricted proof.

\subsection{A singleton decomposition and two nine-identity bases}\label{case:232}
This note proves finite basability for the catalog semiring 10157 and its multiplicative opposite 10179. The proof allows arbitrarily many variables. Finite computations accompanying the note check the formulas and certify equal varieties; they are not the completeness argument.

\par\smallskip\noindent\textit{1. The algebra and the proposed basis.}\quad 
Work in the signature \((+,\cdot)\) without constants, relative to the axioms for additively idempotent semirings. Write \(u\leq v\) for \(u+v\approx v\). Let \(A=\{z, b, a, e, t\}\), where \(z<b<a<t\) and \(e<t\) is incomparable with \(z, b, a\). Addition is join. Put \(J=\{z, b, a\}\) and \(G=\{e, t\}\), and define \(\rho(z)=z\), \(\rho(b)=\rho(a)=b\). Multiplication is

\noindent \(\cdot\)\enspace\textbar\enspace \(z\)\enspace\textbar\enspace \(b\)\enspace\textbar\enspace \(a\)\enspace\textbar\enspace \(e\)\enspace\textbar\enspace \(t\)\par
\noindent \(z\)\enspace\textbar\enspace \(z\)\enspace\textbar\enspace \(z\)\enspace\textbar\enspace \(z\)\enspace\textbar\enspace \(z\)\enspace\textbar\enspace \(z\)\par
\noindent \(b\)\enspace\textbar\enspace \(z\)\enspace\textbar\enspace \(z\)\enspace\textbar\enspace \(z\)\enspace\textbar\enspace \(b\)\enspace\textbar\enspace \(b\)\par
\noindent \(a\)\enspace\textbar\enspace \(z\)\enspace\textbar\enspace \(z\)\enspace\textbar\enspace \(z\)\enspace\textbar\enspace \(a\)\enspace\textbar\enspace \(a\)\par
\noindent \(e\)\enspace\textbar\enspace \(z\)\enspace\textbar\enspace \(b\)\enspace\textbar\enspace \(b\)\enspace\textbar\enspace \(t\)\enspace\textbar\enspace \(t\)\par
\noindent \(t\)\enspace\textbar\enspace \(z\)\enspace\textbar\enspace \(b\)\enspace\textbar\enspace \(b\)\enspace\textbar\enspace \(t\)\enspace\textbar\enspace \(t\)\par
These are names of elements, not constants in the signature. The catalog relabeling is \((z, b, a, e, t)\mapsto(2,\allowbreak{}1,\allowbreak{}0,\allowbreak{}3,\allowbreak{}4)\).

Theorem 1. The following nine identities form a basis for \(\operatorname{Id}(A)\):

\CaseDisplay{\begin{array}{rll}
\mathrm{(B1)}&xyz\approx xzy,\\
\mathrm{(B2)}&x^3\approx x^2,\\
\mathrm{(B3)}&x^2y^2\approx y^2x^2,\\
\mathrm{(B4)}&xy\approx x^2y+xy^2,\\
\mathrm{(B5)}&xyz\leq xy,\\
\mathrm{(B6)}&xyz\leq yz,\\
\mathrm{(B7)}&xy\leq x+x^2,\\
\mathrm{(B8)}&xy\leq y+y^2,\\
\mathrm{(B9)}&x^2\leq x+yz.
\end{array}}

The opposite algebra has the basis obtained by reversing every word. The number nine is an upper bound on basis size, not a claim of minimality.

\par\smallskip\noindent\textit{2. Word values and validity.}\quad 
The multiplication has \(J^2=\{z\}\) and \(G^2=\{t\}\). Right multiplication by a member of \(G\) fixes \(J\), while left multiplication by a member of \(G\) acts on \(J\) as the idempotent map \(\rho\). Also \(z\) is an absorbing element for multiplication. These statements verify associativity: triples with at least two \(J\) entries have value \(z\), triples with no \(J\) entry have value \(t\), and a single \(J\) entry has its original value if first and its \(\rho\)-image otherwise, in either association.

Each multiplication translation is increasing on \(z<b<a<t\). It also preserves \(e+u=t\) for each \(u\in J\): this follows from the displayed table, using \(e+v=t\) for \(v\in J\). The remaining joins involve equal elements or \(t\). Thus both distributive laws hold.

For a word of length at least two, its value is:

\textbullet{} \(t\) if every position is assigned a member of \(G\);

\textbullet{} \(j\) if precisely one position is assigned \(j\in J\) and that position is first;

\textbullet{} \(\rho(j)\) if precisely one position is assigned \(j\in J\) and it is not first;

\textbullet{} \(z\) if at least two positions are assigned members of \(J\).

This formula is proved by induction, and includes assignments of \(z\). In particular, every nonlinear word has value in the chain \(J\cup\{t\}\).

Identities (B1)--(B3) follow directly from the formula. For (B4), consider the four cases \((x, y)\in G^2, J^2, J\times G, G\times J\). The values on both sides are respectively \(t, z, x,\rho(y)\). Deleting the first or last position of a word of length three cannot decrease its value: a zero value is below every nonlinear value; otherwise the shortened word retains the only lower factor or has value \(t\). This proves (B5)--(B6).

For (B7), if \(x\in G\), then \(x+x^2=t\); if \(x\in J\), then \(x+x^2=x\) and \(xy\in\{z, x\}\). For (B8), if \(y\in G\) its right-hand side is \(t\), and if \(y\in J\), then \(xy\in\{z,\rho(y)\}\leq y\). Finally, if \(x\in J\) then \(x^2=z\leq x+yz\), because \(yz\in J\cup\{t\}\). If \(x\in G\) then \(x^2=t=x+yz\). This proves (B9) and the validity of the whole basis.

\par\smallskip\noindent\textit{3. Decomposing words into three kinds of summands.}\quad 
Fix an order on the variables. For a nonempty finite set \(C\), put

\CaseDisplay{Q_C=\prod_{x\in C}x^2.}

For \(s\in C\) and \(|C|\geq2\), put

\CaseDisplay{H_{s, C}=s\prod_{x\in C\setminus\{s\}}x^2,
\qquad
K_{s, C}=\left(\prod_{x\in C\setminus\{s\}}x^2\right)s.}

All products in these definitions are nonempty. Call these expressions normal summands. Their support is \(C\). The distinguished variable \(s\) in an \(H\) or \(K\) summand occurs exactly once.

Lemma 2 (singleton decomposition). From (B1)--(B4), every nonlinear word is equal to a nonempty sum of normal summands with the same support.

Proof. Identity (B1), applied to a nonempty prefix followed by two letters, permutes all positions after the first. If the first variable occurs again, bring an extra occurrence immediately after it. Group the other repeated variables using the same suffix permutations, and apply (B2) to reduce each multiplicity greater than two to two. When no variable occurs exactly once, (B3) orders the resulting square blocks into \(Q_C\). When exactly one variable \(s\) occurs once, the result is \(H_{s, C}\) if it is first and \(K_{s, C}\) otherwise; (B3) orders the square blocks.

If at least two variables occur once, use (B1) to make two of them adjacent. If one is first, move the other immediately after it; otherwise use only suffix permutations. Apply (B4) to this adjacent pair inside its existing left and right contexts, and distribute. In each of the resulting two words the number of variables occurring exactly once decreases by one. No variable is added to the support. Induction on that number ends in the previous cases. Empty contexts here mean omitted substrings, not unit constants. \(\square\)

Every polynomial can therefore be written as a sum of a set \(L\) of linear variables and a set \(N\) of normal summands. The overall sum is nonempty; either of \(L, N\) may be empty.

\par\smallskip\noindent\textit{4. Exact absorption criterion.}\quad 
Lemma 3. Let \(P=\sum_{x\in L}x+\sum_{u\in N}u\) be in the preceding form.

1. \(A\models y\leq P\) holds exactly when \(y\in L\).

2. For a normal summand \(W\) with support \(C\), \(A\models W\leq P\) holds exactly when \(N\ne\varnothing\) and one of the following alternatives holds:

\textbullet{} \(L\cap C\ne\varnothing\); or

\textbullet{} \(N\) has a summand \(u\) with support \(D\subseteq C\) satisfying the rule below.

\noindent Left summand \(W\)\enspace\textbar\enspace Allowed witness \(u\) on \(D\subseteq C\)\par
\noindent \(Q_C\)\enspace\textbar\enspace Any normal summand\par
\noindent \(H_{s, C}\)\enspace\textbar\enspace \(s\notin D\), or \(u=H_{s, D}\)\par
\noindent \(K_{s, C}\)\enspace\textbar\enspace \(s\notin D\), or \(u=H_{s, D}\) or \(u=K_{s, D}\)\par
Moreover, every absorption described here follows from (B1)--(B9).

Proof of necessity. If \(y\notin L\), assign \(y=a\) and all other variables \(z\). Every nonlinear summand is \(z\), as is every linear summand of \(P\); the left side is \(a\). Thus linear absorption requires \(y\in L\).

For a nonlinear left side, if \(N\) is empty, assign every variable \(e\). The left side is \(t\) and the right side \(e\). Thus \(N\ne\varnothing\) is necessary. Now suppose \(L\cap C=\varnothing\) and the appropriate witness does not exist.

For \(Q_C\), assign variables in \(C\) the value \(e\) and all others \(z\). Every summand of \(P\) is \(z\), whereas \(Q_C=t\).

For \(H_{s, C}\), assign \(s=a\), variables of \(C\setminus\{s\}\) the value \(e\), and all others \(z\). The left side is \(a\). Every right summand with a variable outside \(C\) is \(z\). Among summands supported in \(C\), absence of the stated witness forces \(s\) either to be repeated, giving \(z\), or to be the unique singleton in a \(K\) summand, giving \(b\). All right linear summands are \(z\). Hence \(P\leq b<a\).

For \(K_{s, C}\) use the same assignment with \(s=b\). Its value is \(b\). Every right nonlinear summand supported in \(C\) must repeat \(s\), so it has value \(z\); all other right summands are also \(z\). This again separates the two sides.

Proof of sufficiency by derivation. We first record permitted insertions. Appending or prepending a nonempty word to a nonlinear word decreases it, by (B5) and (B6). To insert letters at any position after the first, append them and permute the suffix using (B1). Thus any such insertion gives a derivable inequality in the decreasing direction.

Suppose first that a nonlinear witness \(u\) exists. For \(Q_C\), insert copies of variables until every variable of \(C\) occurs at least twice and no other variable occurs. The resulting word normalizes to \(Q_C\), so \(Q_C\leq u\). For \(H_{s, C}\), if \(u=H_{s, D}\), insert the missing square blocks after its first position. If \(s\notin D\), prepend \(s\), repeat all variables already in \(D\), and insert the other required square blocks. In both cases the result normalizes to \(H_{s, C}\), and hence \(H_{s, C}\leq u\).

For \(K_{s, C}\), a witness \(K_{s, D}\) is extended by missing square blocks. A witness \(H_{s, D}\) can be preceded by a variable of \(D\setminus\{s\}\); this set is nonempty by definition. The resulting word has its sole singleton \(s\) after the first position and normalizes to a \(K\) summand; insert the remaining square blocks. If \(s\notin D\), first repeat every variable of \(u\) and then append \(s\) and insert the other square blocks. Again \(K_{s, C}\leq u\). Each construction keeps \(s\) a singleton when required. Since \(u\) is a summand of \(P\), the desired absorption follows.

Finally suppose \(g\in L\cap C\) and choose any \(v\in N\). Repeated use of (B5)--(B6) deletes positions of \(W\) until a two-letter word containing a chosen occurrence of \(g\) remains. If that occurrence is not last, retain it with the following letter; if it is last, retain it with the preceding letter. By (B7) or (B8), respectively,

\CaseDisplay{W\leq g+g^2.}

Split the nonlinear word \(v\) into two nonempty factors. Substituting them for \(y, z\) in (B9) gives \(g^2\leq g+v\). Therefore

\CaseDisplay{W\leq g+g^2\leq g+v\leq P.}

This proves the sufficiency and the equational part of the lemma. \(\square\)

\par\smallskip\noindent\textit{5. Completeness and the opposite.}\quad 
Proof of Theorem 1. Distribute both sides of an arbitrary valid identity and apply Lemma 2. Each summand on either side is absorbed by the other side. Lemma 3 derives every such absorption from (B1)--(B9). Adding those derivations gives both inequalities and hence the original identity. There is no restriction on the number of variables or summands.

Reversal fixes variables and sums and reverses the order of every product. A reversed identity holds in \(A^{\mathrm{op}}\) exactly when the original holds in \(A\), and reversing an equational derivation is again a derivation. Thus the nine reversed laws are a basis for \(A^{\mathrm{op}}\). \(\square\)

The catalog tables of 10179 have the same addition as 10157 and exactly the transposed multiplication. Equal-variety certificates, independently verified in both directions, transfer these bases to the following types:

\textbullet{} 10157 side: 10152, 10157, 10219, 10220, 10224, 10225, 10621, 10628.

\textbullet{} 10179 side: 10166, 10179, 10237, 10238, 10251, 10254, 10634, 10726.

Membership in a finitely based variety alone is not used to transfer finite basability. The two-direction equality check is required.

The source manuscript for this case credits Aifa Wang and Lili Wang, School of Mathematical Sciences, Chongqing University of Technology.

The source manuscript for this case acknowledges support from the National Natural Science Foundation of China (12371024), the Science and Technology Research Program of Chongqing Municipal Education Commission (KJZD-K2024011102) and the Chongqing Natural Science Foundation Innovation and Development Joint Fund (Municipal Education Commission) (CSTB2025NSCQ-LZX0067).

\subsection{A finite basis for 10216 from three additive-threshold witnesses}\label{case:233}
Signature (+, multiplication), without constants, relative to the additively idempotent semiring axioms. Finite family specified below; not an enumerated small basis. Completeness is unrestricted.

\par\smallskip\noindent\textit{Word values and normal form.}\quad 
Addition rows:00004,\allowbreak{}01114,\allowbreak{}01224,\allowbreak{}01234,\allowbreak{}44444; order \(3<2<1<0<4\). Multiplication rows:22220,\allowbreak{}22220,\allowbreak{}22222,\allowbreak{}33333,\allowbreak{}22224. Literals retain their values. For nonlinear hT, heads 2 and 3 retain their values; heads 0 or 1 give 0 if every tail is 4, and 2 otherwise; head 4 gives 4 if every tail is 4, and 2 otherwise. This follows by induction:once a low tail gives 2, it remains 2, and 4 is a right identity on\{0,\allowbreak{}2,\allowbreak{}3,\allowbreak{}4\}.

W1 \(xyz=xzy\) and W2 \(xyy=xy\) follow. Normalize to the head and sorted tail set by actual adjacent tail swaps and duplicate deletions. All substituted prefixes are nonempty; head occurrences in the tail remain separate.

\par\smallskip\noindent\textit{Exact criterion.}\quad 
For original P, a literal x is absorbed iff the literal x occurs in P. For a nonlinear target hT, \(C=\{h\}\) union T, the criterion is: (H) some original word has head h; (C) some original word is contained in C; (R) h belongs to T, or an original word is contained entirely in T, or an original nonlinear h-word has tail contained in T.

Missing literal:assign \(x=0\), others 3. Original nonlinear x-words give 2; other heads 3. Without literal x the right sum is below 0. Missing H: assign \(h=2\), others 3. Target 2 and right 3 separate. Missing C:assign all variables in C to 4, outsiders 3. Target 4 while every original word with an outsider gives 2 or 3. For failure of R, with H, C already satisfied, h is outside T. Assign \(h=1\), \(T=4\), others 3. Target 0; original literal h gives 1, nonlinear h-words have a tail outside T and give 2, other heads 3 give 3, and words headed in T must have an outside variable since no word is contained in T, hence give 2. The right sum is at most \(1<0\). These are explicit separating assignments for every failed condition.

Sufficiency:3 is least. A target 2 has head 2, or head 0/\allowbreak{}1/\allowbreak{}4 with a low tail; an original same-head word from H always has value at least 2. Target 4 has all target variables 4; C supplies a contained word valued 4. Target 0 has head 0/\allowbreak{}1 and all tails 4, so its head is not in T. R gives either an original word entirely in T, valued 4, or a nonlinear same-head word with tail in T, valued 0. Literal targets are original summands. Thus the criterion is exact on arbitrary assignments and arbitrary polynomial size.

\par\smallskip\noindent\textit{Three witnesses and reversible finite compression.}\quad 
For a valid nonlinear target choose an original H witness a and contained C witness b. If the head is repeated, these suffice. Otherwise choose one original R witness c, either contained entirely in T or a nonlinear h-word with tail in T. Q=\{a, b, c\}, with duplicates removed, has at most 3 words and still satisfies H, C, R. Hence \(w\leq Q\) is valid and Q consists of original right summands. Literal targets require no new absorption.

Among the at most 4 normal words in \(w\leq Q\), reserve all their heads (at most 4). Partition remaining variables by their nonempty tail-incidence patterns; there are at most 15. Replace each pattern class by one fresh variable, retain reserved heads separately, and normalize with W1, W2. The result uses at most 19 variables, word length at most 20, and at most 3 right words. It remains valid because identification is a substitution.

Invert by substituting the nonempty product of original class members. Their tail incidences were identical, so original tail sets recover exactly; reserved heads and their tail occurrences are unchanged. W1, W2 normalize the expanded words back to the originals. No empty substitution or loss of head multiplicity occurs.

Let Sigma be W1, W2 together with ALL valid normal absorptions on a fixed 19-variable alphabet, word length at most 20, and at most 3 nonempty right words. H, C, R explicitly specify membership in this finite family. Every valid absorption selects Q, compresses into Sigma, and expands back. Add other original summands by transitivity. Distributivity and summandwise absorption in both directions then derive every valid identity. Sigma is a finite relative basis for 10216. Reversal gives the opposite basis. Catalogue transfer requires complete HSP both directions; finite checks supplement the unrestricted theorem rather than replacing it.

\subsection{Finite bases for 10235 and 10260 using global heads and contained witnesses}\label{case:234}
Signature (+, multiplication), without constants, relative to additively idempotent semiring axioms. Both finite families are effectively specified; no small enumerated basis or minimality is asserted.

\par\smallskip\noindent\textit{1. Tables and arbitrary word values.}\quad 
The common addition has chain order \(3<2<1<0<4\) and rows 00004,\allowbreak{}01114,\allowbreak{}01224,\allowbreak{}01234,\allowbreak{}44444. Multiplication tables are

{\small\par\smallskip
\CaseLine{\ttfamily 10235\qquad 10260}
\CaseLine{\ttfamily 2 2 2 2 0\qquad 2 2 2 2 1}
\CaseLine{\ttfamily 2 2 2 2 1\qquad 2 2 2 2 1}
\CaseLine{\ttfamily 2 2 2 2 2\qquad 2 2 2 2 2}
\CaseLine{\ttfamily 3 3 3 3 3\qquad 3 3 3 3 3}
\CaseLine{\ttfamily 1 1 2 2 4\qquad 0 1 2 2 4}
\smallskip}

For every nonempty word w, direct induction gives:

\textbullet{} If its literal head is 3 the value is 3, regardless of its suffix.

\textbullet{} Otherwise any factor 2 or 3, or at least two positions in \{0,\allowbreak{}1\}, gives value 2.

\textbullet{} If all factors are 4, the value is 4.

\textbullet{} With a unique \{0,\allowbreak{}1\}-valued position and all other positions 4, 10235 gives 0 precisely when that position is first and valued 0; 10260 gives 0 precisely when it is last and valued 0. Otherwise it is 1.

The formulas also include literal words. The row of 3 is constant 3; 2 is stable after every product; a second \{0,\allowbreak{}1\} factor gives 2. The only difference in the unique-factor case is whether left or right multiplication by 4 changes 0 to 1, as displayed in the two tables.

\par\smallskip\noindent\textit{2. Exact absorption criteria.}\quad 
For target w put \(C=c(w)\), \(h=h(w)\), \(t=t(w)\), and let U be the original right words whose contents lie in C. The absorption \(w\leq P\) holds exactly when all the following conditions hold:

(G) some original right word starts in h; its content need not lie in C; (C) U is nonempty; (S) for every simple variable x of w, some u in U omits x or contains x exactly once. In 10235, if \(x=h\) that single occurrence must be the head of u. In 10260, if \(x=t\) it must be the last position of u.

Necessity of G: assign \(h=2\) and every other variable 3. The target is 2; without a same-head right word all right values are 3. Necessity of C: assign C the value 4 and all outsiders 2. A target is 4, whereas all noncontained right words are 2. For a failed simple-x condition, assign \(x=1\), other variables of C the value 4 and outsiders 3. The target is 1; noncontained right words are 2 or 3, and contained right words must repeat x and are 2. For failure solely of the extra head/tail requirement, assign that simple endpoint 0 instead. The target is 0, while contained single-occurrence words at the wrong position give 1 and all others are at most 2. A literal target therefore requires its literal summand.

For sufficiency, target 3 is bottom. If the target is 2, its head is not 3 and G supplies a right word with value at least 2, regardless of its external variables. A target 4 has all its variables 4 and C supplies a 4-valued word. A target 1 or 0 has a unique \{0,\allowbreak{}1\}-valued occurrence belonging to a simple x and all other target variables 4. An omitting contained witness is 4; a single-x contained witness has value 0 or 1. At target 0 the endpoint clause ensures value 0 in the latter case. This proves exactness for arbitrary variables, lengths and summand counts. The G witness must NOT be required to be contained in C.

\par\smallskip\noindent\textit{3. Structural laws and word normal forms.}\quad 
For 10235 use these six sound identities:

{\small\par\smallskip
\CaseLine{\(xyz=xzy, xyx=x^{2}y, x^{3}=x^{2},\)}
\CaseLine{\(xy=x^{2}y+xy^{2}, x^{2}\leq x, xy^{2}\leq x.\)}
\smallskip}

For 10260 use these nine:

{\small\par\smallskip
\CaseLine{\(xyzt=xzyt, xy^{3}z=xy^{2}z, x^{3}y=x^{2}y, xy^{3}=xy^{2},\)}
\CaseLine{\(x^{3}=x^{2}, x^{2}yx=xyx, xy=x^{2}y+xy^{2}, x^{2}\leq x, xy^{2}\leq x.\)}
\smallskip}

The word formulas prove soundness. All word equalities retain the head, content and simple/repeated multiplicities, and the 10260 ones also retain the literal tail. In D: \(xy=x^{2}y+xy^{2}\), a unique \{0,\allowbreak{}1\} factor survives in the branch doubling the other factor; its required first/last status is retained. The other branch is 2. The cases with head 3, value 2, or all factors 4 give the same value on both sides. Squaring a factor can only keep the value or lower it to 2, so L:\(x^{2}\leq x\) holds. The same observation proves \(xy^{2}\leq x\), with head 3 retained as 3.

10235: permute the tail, gather each variable into a block, and use \(x^{3}=x^{2}\) to reduce multiplicities to one or two. The literal head is fixed; the other blocks are sorted. Thus the normal form records the head and capped multiplicities.

10260: fix both endpoints and permute the interior by \(xyzt=xzyt\). An internal repeated variable is reduced using \(xy^{3}z=xy^{2}z\); head and tail repetitions use \(x^{3}y=x^{2}y\) and \(xy^{3}=xy^{2}\). When head equals tail, gather interior copies by the head, reduce to at most one extra copy, and use \(x^{2}yx=xyx\) to remove it. Singleton-content words use \(x^{3}=x^{2}\). Finally sort the interior. The normal form records both endpoints and capped total multiplicities. Every block substitution is nonempty.

\par\smallskip\noindent\textit{4. Local guard, simple recovery, and splitting.}\quad 
Fix a valid nonlinear target w. Choose an original same-head word v from G, possibly containing outsiders, and a contained original word u from C. Apply L in contexts to double all their simple variables and normalize, obtaining saturated \(v_*\leq v\) and \(u_*\leq u\). Let g be the normalized version of w in which every variable has been doubled. It has content C and the correct head (and in 10260 the correct tail). The absorption

{\small\par\smallskip
\CaseLine{\(g \leq  v_* + u_*\qquad (G*)\)}
\smallskip}

is valid by G, C, S: v\_* supplies the global head, u\_* is contained, and g has no simple variables. Its derivation from a FIXED finite family is given in Section 5. Outsiders in v\_* are allowed throughout this step.

For each simple x of w choose its original S-witness. Double every variable other than x and normalize, producing \(u_{x}\leq u\). It still omits x or has it once, retaining the required endpoint when applicable. Let \(w_{x}\) be w with all variables except x doubled and normalized. Then

{\small\par\smallskip
\CaseLine{\(w_{x} \leq  g + u_{x}\qquad (R*)\)}
\smallskip}

is valid by the exact criterion. It too is derived by finite compression.

For a current target with at least two simple variables, split between two such occurrences as nonempty X, Y and apply D. Both \(X^{2}Y\) and \(XY^{2}\) have the same content, head and (for 10260) tail, and strictly fewer simple variables. Surviving simples retain whether they are the relevant endpoint. Normalize and repeat. Leaves with zero simple variables equal g; those with one simple x equal \(w_{x}\). Since each leaf is below P, the D equalities derive the original target. Zero/one-simple targets already are leaves; literal targets are original summands.

\par\smallskip\noindent\textit{5. Reversible finite compression.}\quad 
In either G* or R* there are at most three words. Reserve their literal heads and tails and, for R*, the simple variable x: at most seven variables. Every other variable is absent or repeated in each word. Partition these nonreserved variables by their nonzero incidence pattern across the three words, giving at most seven classes. Identify each class to a fresh variable and normalize. Validity is preserved by genuine substitution; there are at most fourteen variables and every word has length at most twenty-eight. The absorption has at most two right words.

Inversely replace each class variable by the nonempty product of its members. Every class had common multiplicity two in every word where it occurs. Its two occurrences restore two occurrences of each member. Tail or interior permutation and the power laws recover the original normal word. Reserved variables retain their multiplicities and endpoints; no nonreserved class contains an original endpoint. In 10235 only the head needs to be retained: sorting a compressed tail may put a class variable last, and the inverse is still correct because that tail is freely permutable. Reserving its original tail as well is harmless. No repeated class is expanded to a single copy, and no empty substitutions are used.

Fix a fourteen-variable alphabet. Take the respective six or nine structural laws and ALL valid absorptions with one normal left word, one or two normal right words and word lengths at most twenty-eight. Section 2 decides this finite family effectively. Compression and its inverse derive G* and R*, and Section 4 derives every valid word absorption with no bound on original variables, word lengths or summand count. Distributive expansion and summandwise absorption in both directions therefore derive every valid identity. This proves finite basability of both sources. Reversing words gives bases of their multiplicative opposites. Catalogue transfers still require independently verified complete HSP in both directions.

\subsection{Finite bases for 10578,\allowbreak{}13102,\allowbreak{}13107 from contained threshold witnesses}\label{case:235}
In the constant-free signature use the semiring axioms, \(xy=yx\), \(x^{3}y=x^{2}y\), and all valid absorptions on 257 variables, length at most 515 and at most 4 upper words. Interpret the finite valid family separately in each source. This is an effective specified family; it is not an enumerated small basis. The following argument is unrestricted.

The common multiplication description uses an absorbing element Z, an active element A with \(A^{2}=Z\), and a three-element subsemigroup B=\{G, R, K\} with \(G^{2}=R\), \(G^{3}=K\) and every other product in B equal to K. Also \(AB=A\) and Z absorbs everything. The actual labels are:

\noindent Source\enspace\textbar\enspace Z\enspace\textbar\enspace A\enspace\textbar\enspace G\enspace\textbar\enspace R\enspace\textbar\enspace K\par
\noindent 10578\enspace\textbar\enspace 2\enspace\textbar\enspace 1\enspace\textbar\enspace 0\enspace\textbar\enspace 3\enspace\textbar\enspace 4\par
\noindent 13102\enspace\textbar\enspace 2\enspace\textbar\enspace 3\enspace\textbar\enspace 1\enspace\textbar\enspace 0\enspace\textbar\enspace 4\par
\noindent 13107\enspace\textbar\enspace 3\enspace\textbar\enspace 2\enspace\textbar\enspace 1\enspace\textbar\enspace 0\enspace\textbar\enspace 4\par
Thus a nonlinear word is Z when a letter is Z or A occurs at least twice; A when exactly one letter is A and all others lie in B; R for a quadratic GG; and K otherwise. Induction proves this formula. The word laws sort letters and cap counts at 2 except that a long singleton-content word stays a cube. They follow from the same formula.

For target w retain only original right words with content contained in c(w); call this polynomial U. Let L be its literal variables, Q the union of quadratic contents, and J its long words. Its short-role system is bad if J is nonempty, or in 10578 only, L intersects Q. A simple variable means one occurring exactly once in w. Define its guard as follows:

\textbullet{} In 10578 and 13107, an original word of U containing x at most once, including zero times.

\textbullet{} In 13102, the literal x, a nonlinear original of U containing x once, or a long original of U omitting x.

The exact criterion is: a literal target must be original; a nonlinear target requires the short-role system to be bad, or to be quadratic with c(w) contained in Q; and each simple variable must have a guard. In particular, a long target in 13102 or 13107 needs an original CONTAINED long word. Neither a long word outside the target content nor a global role conflict can replace it.

Here is a derivation using every proper principal ideal. All sources have \{Z\} and \{Z, A\}. For \{Z, A\}, word membership means some letter lies in \{Z, A\}. An original contained word is precisely what forces membership of the target; assigning external variables Z proves necessity. For \{Z\}, a target can fail with all letters in B, or with one simple variable x assigned A and other target variables in B. External variables may again be set Z, eliminating all originals outside U. For the single-A case, an original contained word fails precisely when it has at most one x. This gives the weak simple-occurrence condition and the contained-word condition in the all-B case.

The other two ideals are:

\noindent Source\enspace\textbar\enspace First ideal\enspace\textbar\enspace Second ideal\par
\noindent 10578\enspace\textbar\enspace \{Z, A, G\}\enspace\textbar\enspace \{Z, A, R\}\par
\noindent 13102\enspace\textbar\enspace \{Z, G\}\enspace\textbar\enspace \{Z, G, R\}\par
\noindent 13107\enspace\textbar\enspace \{Z, A, G\}\enspace\textbar\enspace \{Z, A, G, R\}\par
For an all-B target, a nonlinear original cannot belong to the first ideal; all literals can be set G. Thus a nonlinear target requires an original contained nonlinear word. In the second ideal, long originals are impossible and quadratics force both letters G. Original literal domains are \{R\} in 10578 and \{G, R\} in 13102/\allowbreak{}13107. They conflict with quadratic roles only in 10578. Consequently the second-ideal domain system is impossible exactly in the stated bad case. If possible, a long target fails and a quadratic is forced into the ideal exactly when its variables occur in original quadratics. This condition already supplies a contained nonlinear original, so the first-ideal all-B condition adds no independent requirement.

Only 13102 excludes A from these two ideals. For its second ideal, when a simple target x is A and other target variables lie in B, the literal x or a nonlinear once-x original fails; so does a long original omitting x. Originals with two x have value Z and add no restriction. All remaining literals have domains \{G, R\} and omitting quadratics force G. These domains are always compatible. Thus exactly the three stronger guards are needed. They also suffice for its first ideal. The other sources include A in both ideals, so no further single-A condition is required. If a guard is absent, choose compatible G/R assignments and external Z to obtain an explicit separator. For a literal target x absent from P, assign \(x=A\) and every other variable Z; all originals are Z, while x is not below Z. This completes necessity and sufficiency of the exact criterion for arbitrary polynomials.

If the short-role system is possible, a nontrivial nonlinear target is quadratic. Retain at most two original quadratic role witnesses and at most two guards. These at most 4 originals still satisfy the criterion. If the system is bad, form the saturated target-content word S, with every variable twice and a cube for singleton content. One original contained long word, or the two conflicting original roles in 10578, proves \(S\leq P\). Restore each simple variable x one at a time using \(r'\leq r+u\), where u is its original guard and r' replaces the current two x occurrences by one. The protected r is long and contained in c(w), so it supplies the required bad system. Every previously restored variable has r as a once-occurring nonlinear witness, and the new x has u. Before a restoration r is long; degree 2 can occur only after the last step for a quadratic target. Finish a quadratic square using \(x^{2}\leq x^{3}\). All of these are valid local absorptions with at most 4 upper words, proved directly from the criterion.

Normalize each such absorption, at most 5 words. Reserve up to 3 distinct variables per word, all if fewer, so at most 15 are reserved. Group other variables by COMMON multiplicities 0/\allowbreak{}1/\allowbreak{}2 across these words, at most \(3^{5}-1=242\) groups. Use the common multiplicity, not its sum. The resulting absorption has at most 257 variables and length 515. Containments, missing/single occurrences, short-word roles and long degree survive; all short variables are reserved. Hence the compressed absorption satisfies the exact criterion and belongs to the finite family. Substitute each nonempty group list back. Its common multiplicities and the two word laws recover every original normalized word. This proves the local absorptions, and their saturation/restoration chain proves every input absorption. Polynomial identities split into these absorptions in both directions, establishing unrestricted completeness.

Finite audits supplement the proof with actual table assignments, independent separating assignments and literal derivation replay.

\subsection{Finite bases for 10587 and 13323 from contained tail guards}\label{case:236}
Use semiring axioms without constants, \(xyz=xzy\), \(xyyz=xyz\), and all valid absorptions over 27 variables, length 29, at most 3 upper words, interpreted in each source. This is a specified finite family, not a small list.

For 10587 the additive least and multiplicative absorber is \(Z=2\), the active value is \(A=1\), and the quadratic exceptional pair(0,\allowbreak{}0) gives 3. For 13323 \(Z=3\), \(A=2\), and the exceptional pair(1,\allowbreak{}1) gives 0. In either source a word containing Z is Z. A nonlinear word with no Z is Z if a NONINITIAL letter has value A; otherwise it is A if the head is A; otherwise it is the exceptional quadratic value on that pair, and 4 in every other case. Table multiplication proves the formula and its induction for every length. The word laws above follow.

A long word's normal form retains its head h and its noninitial set T: if T is singleton\{t\}, use htt; otherwise use h followed by the distinct members of T in a fixed order. The tail is permutable and duplicate tail letters can be deleted while total length exceeds 3, using \(xyyz=xyz\). The head can belong to T and must remain there when it originally did. In particular \(xxy=xyy\) is false, and deleting a repeated head is forbidden. Short words retain order and degree.

For target w let U be the original upper words contained in c(w). Let L be the original literal variables in U, Q the union of contents of quadratic U words, J the long U words. Let T(v) be the actual noninitial set, empty for a literal. Write F if some u in U satisfies T(u) subset T(w). Every failing assignment avoids Z on target variables; sending all outside variables to Z kills exactly the discarded words. Thus U suffices for an exact analysis.

Both sources have the following criterion for literals:a target x requires the original literal x. For nonlinear targets F is necessary. For 10587 the other condition is J nonempty, or L intersects Q, or the target is quadratic and c(w) subset Q. For 13323 the other condition is J nonempty, or the target is quadratic and c(w) subset Q. All words and occurrences in these conditions are ORIGINAL.

Here is an unrestricted principal-ideal proof. In down(Z), a nonlinear upper word over the nonzero domain must have an active noninitial letter, and a literal upper word is impossible. A failing target has no active noninitial letter. Such a simultaneous assignment exists exactly when every U tail contains a variable outside T(w), which is the negation of F. Indeed each U is contained in c(w); outside T(w) there is at most the target head, when that head is simple. Give that head A and other target variables an inactive value. If U is empty, an all-inactive assignment separates directly. Repeated target heads and literal upper words are therefore covered without special exceptions. For a literal target, all contained nonlinear words are powers of x; assigning \(x=A\) makes them Z while x remains A, so its own literal is necessary.

In every other proper ideal the active value A is allowed. Any failing nonlinear target must consequently use only the inactive domain D: D=\{0,\allowbreak{}3,\allowbreak{}4\} for 10587, D=\{0,\allowbreak{}1,\allowbreak{}4\} for 13323. On D a long word is always 4. For 10587, the nonleast ideals are \{1,\allowbreak{}2\}, \{0,\allowbreak{}1,\allowbreak{}2\}, \{1,\allowbreak{}2,\allowbreak{}3\}. The first permits no contained word on D, so requires U nonempty. The second allows literals only at 0 and no nonlinear U word. The third allows literals only at 3 and quadratics only when both letters are 0. Thus its upper restrictions are impossible exactly when J is nonempty or L intersects Q; otherwise a long target fails, and a quadratic target is forced exactly when both letters occur in Q. These branches already imply a contained nonlinear word, so the second ideal has no remaining obstruction.

For 13323, the nonleast ideals are \{2,\allowbreak{}3\}, \{1,\allowbreak{}2,\allowbreak{}3\}, \{0,\allowbreak{}1,\allowbreak{}2,\allowbreak{}3\}. The first allows no word on D. The second permits literals at 1 and no nonlinear word. In the third literals can be 0 or 1, and quadratic letters must both be 1. These restrictions are always consistent when J is empty. Thus a long target needs J; and a quadratic needs J or both letters in Q. Again these branches provide a contained nonlinear word for the other two ideals. This proves both necessity and sufficiency, including literal targets already handled above.

For a nonlinear target retain one F witness and one contained long word, or(for 10587) a conflicting literal and quadratic, or at most two quadratics covering the target. At most 3 originals suffice. For a literal retain itself. The same criterion holds for the selected words.

There are at most 4 words including the target. Reserve all short variables and each long head plus up to two other distinct variables. At most 12 variables are reserved. All unreserved letters occur only noninitially. Group them by their common incidence in these words: at most 15 groups. Identification preserves validity and yields at most 27 variables, normal length 29. The resulting absorption lies in the finite family. Substitute each group variable by the NONEMPTY word listing its members. All heads are reserved, so their occurrence or absence in the tail is preserved. Common incidence recovers every tail set; the reserved variables preserve degree. The two word laws therefore recover each original normal form. Adding back the other original summands proves the original absorption. Distributive expansion and both directions of summand absorption prove unrestricted completeness.

F may use a word omitting the target head, or a literal, but such a word must be contained in the target content. Head repetition cannot be removed by analogy with the previous global-head sources. Opposite and equal-variety transfers need separate complete HSP evidence.

\subsection{Six complete identities for semiring 10595}\label{case:237}
Constant-free signature (+, multiplication), relative to the additively idempotent semiring axioms. A complete finite basis is

{\small\par\smallskip
\CaseLine{\(xy=yx; xyz=xyzz; xyzt\leq xyz; xy\leq xxy; xy\leq xx+yy; x\leq xxx.\)}
\smallskip}

Here \(u\leq P\) means \(u+P=P\).

Addition is 00044/\allowbreak{}01134/\allowbreak{}01234/\allowbreak{}43334/\allowbreak{}44444. Element 2 is bottom; \(1<0<4\) and \(1<3<4\), with 0 and 3 incomparable. Multiplication is 33244/\allowbreak{}33244/\allowbreak{}22222/\allowbreak{}44244/\allowbreak{}44244. A nonlinear word is 2 when any of its variables is 2. Otherwise a quadratic is 3 when both occurrences are in\{0,\allowbreak{}1\}, and all other nonlinear words are 4.

For a target w, let U be the original words whose content is contained in c(w), and E the quadratic edges among them. The exact criterion is:

\textbullet{} Literal x: original x or an original long pure power of x.

\textbullet{} Long w: an original long word in U.

\textbullet{} Quadratic w: an original long word in U, its original edge, or both original endpoint loops.

Necessity: if a literal target has neither original literal nor long pure power, assign its variable 0 and all others 2. The surviving square, if present, is 3, and every other right word is bottom 2. Neither value dominates 0. For a long target with no contained original long word, assign all target variables 1 and outsiders 2. All surviving literals are 1 and quadratics 3, so the right side is at most 3; the target is 4. For a failed square, only its literal may survive; assign its variable 1 and outsiders 2. For a failed distinct-endpoint target xy, choose a missing loop yy. If xx exists, assign \(x=1\), \(y=3\), outsiders 2. Surviving quadratics and literals are at most 3, but the target is 4. If neither loop exists, assign \(x=y=3\), outsiders 2: only literals survive, again with join at most 3. Original xy is absent in these cases. These explicit assignments cover every failed condition for arbitrary numbers of variables and summands.

For sufficiency, a nonlinear target containing a variable 2 is bottom. Otherwise any contained long word is 4 and absorbs the target. A pure long power also dominates its literal: it is 2 exactly at \(x=2\), and 4 otherwise. Original terms are immediate. If both endpoint loops occur, a target quadratic 3 is absorbed by those loops; for target 4, at least one loop is 4. This proves the criterion without finite-size assumptions.

The word identities normalize long words by content. Commute a chosen occurrence to the end and duplicate with \(xyz=xyzz\) using two nonempty prefix blocks; reverse the law to remove surplus occurrences while keeping length at least three. Canonical long forms are sorted content for at least three variables, xxy for two, and xxx for one. An original contained long word can be extended to target content by \(xyzt\leq xyz\) and normalized. Recover a quadratic ab from aab with \(xy\leq xxy\). If no contained long word exists, a new quadratic is derived directly from its two original loops using \(xy\leq xx+yy\). Recover a missing literal only from its original long pure power, normalized to xxx, using \(x\leq xxx\). All substitutions are nonempty and every used right summand is original or explicitly derived.

All six laws satisfy the exact criterion and are independently checked over all assignments. Distributivity and summandwise absorption in both directions establish unrestricted completeness. In particular \(x\leq xx\) fails at \(x=0\), and \(xxy\leq x+xy\) fails at \(x=y=1\). Neither squares nor incident literals may replace the original long pure-power witness. Transfers to further catalogue types require full HSP in both directions.

\subsection{A finite basis for 11577 from global pairs and triple reconstruction}\label{case:238}
In the constant-free semiring signature take \(xy=yx\), \(x^{3}y=x^{2}y\), \(xyzt\leq xyz+xyt+xzt+yzt\), and all valid absorptions on at most 19709 variables, length 39419 and at most 8 upper words. This is an effective specified finite family, not an enumerated small basis. We prove completeness with no bound on the input alphabet or polynomial.

Multiplication rows 40444,\allowbreak{}01444,\allowbreak{}44344,\allowbreak{}44444,\allowbreak{}44444 imply the following nonlinear word values:1 if all letters are 1;0 if exactly one letter is 0 and all others 1;3 for a quadratic with both letters 2;4 otherwise. The first two laws sort letters and cap multiplicity at 2 except that singleton-content long words stay cubes. The third inequality is valid: its upper words cover every variable and every pair of positions, and are all long. This also follows immediately by the displayed values.

For original P let C be the content union, I the variables repeated in at least one original word, and E the unordered pairs of distinct variables co-occurring in an original word. Let L be literal variables, Q the union of quadratic contents, and J the long originals. Call the short-role system bad if J is nonempty or L intersects Q. Then \(w\leq P\) iff:

\textbullet{} c(w) is contained in C.

\textbullet{} Every repeated target variable lies in I. Every distinct pair x, y in c(w) has x in I, y in I, or \{x, y\} in E.

\textbullet{} A nonlinear target has a nonlinear original.

\textbullet{} Either the short-role system is bad, or a literal target is original, or a quadratic target has its two variables in Q.

The proper principal ideals are \{0,\allowbreak{}1\},\{1\},\{2\},\{3\}. For \{1\}, all originals belong exactly when all C variables are 1; this gives content coverage. For \{0,\allowbreak{}1\}, all variables in C must lie in \{0,\allowbreak{}1\}, and the variables assigned 0 must have at most one occurrence in every original word. Thus no variable of I can be 0, and no pair in E can both be 0. The target has the same property exactly under the pair conditions above. A missing loop is separated by assigning that variable 0 and all others 1; a missing distinct pair by assigning its two variables 0. A missing content variable is separated by assigning it 4 and others 1. These arguments prove necessity and sufficiency of both ideal conditions.

The ideal \{2\} contains only literal word values. It forces an original nonlinear word for a nonlinear target. The ideal \{3\} contains literal 3 and quadratic 22, but no long word. Literal domains are \{3\}, quadratic domains are \{2\}, so its original domain system is impossible precisely in the bad case. Otherwise a literal target is forced 3 exactly if it occurs literally, a quadratic is forced 3 exactly if its variables occur in original quadratics, and a long target fails. This proves the full criterion for arbitrary alphabets and numbers of right words.

For a target of length at most 3, retain one original content witness per target variable, at most 3. Retain one original witness per required loop or distinct target pair, at most 3 further words (a loop at either endpoint may witness a pair). Finally retain a long original, or a literal/quadratic conflicting pair, in the bad case; otherwise retain an original literal target or the at most 2 quadratic role witnesses. At most 8 originals suffice. A nonlinear target's pair witness is itself nonlinear, so that requirement adds no further word. The selected polynomial still satisfies the full criterion.

Normalize these at most 9 words, including the target. Reserve at most 3 variables per word, reserving every variable if fewer. There are at most 27 reserved variables. The others have common multiplicity vectors in \{0,\allowbreak{}1,\allowbreak{}2\} across the nine words, at most \(3^{9}-1=19682\) nonzero vectors. Represent a group using its COMMON multiplicity in each word, never the sum of multiplicities. This gives at most 19709 variables and length 39419. Reserved letters preserve short/long degree; all target variables and all short-word variables are reserved. Hence content, target loops and pairs, and the short-role witness remain valid after compression. The compressed absorption belongs to the finite family. Substitute the nonempty list of group members back and use the word laws; common multiplicity recovers the original normalized words. This proves every absorption whose target has length at most 3.

For a valid target w of length at least 4, split it into four nonempty blocks A, B, C, D, taking B, C, D single letters. Each of ABC, ABD, ACD, BCD is a proper subword of length at least 3. Its content and pair requirements are inherited from w; the original short-role system was already bad because w is long. Therefore each subword is absorbed by the original P. Induction on target length proves their absorptions. Substitute A, B, C, D into \(xyzt\leq xyz+xyt+xzt+yzt\) and combine these four results to obtain \(w\leq P\). No empty substitutions or length-zero paths occur. Finally split each polynomial identity into its summand absorptions in both directions. This proves unrestricted finite basability.

\subsection{全局补尾与受守卫换首：六式分支的两个剩余来源}\label{case:239}
\par\smallskip\noindent\textit{1. 记号.}\quad 
记 \(h(w)\) 为词的首字母，\(D(w)\) 为尾部变量集合，\(X=c(w)\)。对于非空右词集 \(F\)，记 \(T=c(F)\)、\(H=h(F)\)、\(E=\bigcup_{u\in F}D(u)\)。内部词集及其内容、首字母、尾部并分别为

\CaseDisplay{F_X=\{u\in F:c(u)\subseteq X\},\quad C_X=c(F_X),\quad H_X=h(F_X),\quad E_X=\bigcup_{u\in F_X}D(u).}

``非一次''表示词长至少二。

\par\smallskip\noindent\textit{2. 来源10939：内部非一次种子与全局补尾.}\quad 
取 \(A=\{z, a, b, d, t\}\)，依次对应目录 \((0,\allowbreak{}1,\allowbreak{}2,\allowbreak{}3,\allowbreak{}4)\)。加法为链 \(d<b<a<z<t\) 上的最大值。右乘 \(a, b\) 为 \(R=(z, a, a, d, t)\)，右乘 \(d\) 为 \(P=(z, d, d, d, t)\)，右乘 \(z, t\) 均为常值 \(t\) 映射。三个右平移保持并、幂等、交换，且 \(RP=PR=P\)；按右因子类别检查乘法结合律和另一侧分配律即得半环。

非一次词的尾部出现 \(z, t\) 时取 \(t\)；否则首字母为 \(z, t\) 时分别取 \(z, t\)，首字母属于 \(\{d, b, a\}\) 时，内容含 \(d\) 则取 \(d\)、否则取 \(a\)。

定理1（判据）。 一次词 \(x\le F\) 有效当且仅当 \(F_{\{x\}}\ne\varnothing\)。非一次词 \(w\le F\) 有效当且仅当

\CaseDisplay{X\subseteq T,\qquad F_X\text{ 含非一次词},\qquad D(w)\subseteq E.\quad\textup{(\ensuremath{A})}}

必要性。 一次情形若无纯幂见证，赋 \(x=b\)、其他变量 \(d\)；右词均取 \(d\)，左值为 \(b\)。反向由 \(x\le x^2\) 及纯幂正规化成立。

对非一次左词，缺失全局内容变量时，令该变量取 \(t\)、其余取 \(a\)，得到左值 \(t\)、右值 \(a\)。没有内部非一次项时，赋 \(X\) 内变量 \(b\)、外部变量 \(d\)；右一次词不超过 \(b\)，右非一次词都取 \(d\)，左值为 \(a\)。若 \(x\in D(w)\setminus E\)，赋 \(x=z\)、其余变量 \(a\)；右尾部全取 \(a\)，右首字母取 \(a, z\)，右和不超过 \(z\)，左值为 \(t\)。

充分性。 真主理想分别由 \(d, b, a, z\) 生成。右和不超过 \(d\) 或 \(b\) 时，所有右变量取 \(d, b, a\)，每个非一次右词还必须含取值 \(d\) 的变量。全局内容覆盖和内部非一次见证使左值为 \(d\)。右和不超过 \(a\) 时，全局覆盖保证左内容只取 \(d, b, a\)，左值不超过 \(a\)。右和不超过 \(z\) 时，右首字母不取 \(t\)，右尾部只取 \(d, b, a\)；全局内容覆盖使左首字母不取 \(t\)，全局尾部覆盖使左尾部只取 \(d, b, a\)，所以左值不超过 \(z\)。

定理2（八条完整基）。 \(A\) 的完整附加基为

\CaseDisplay{\begin{aligned}
xyz&\approx xzy,&xy^2&\approx xy,&x^2y&\approx y^2x,\\
x&\le x^2,&xyz&\le xy+tz,&xy&\le x+y^2,\\
x^2y&\le xy+zx,&xy&\le xz+y^2.&&
\end{aligned}\quad\textup{(\ensuremath{A+})}}

前三式由词值规则验证，其余由定理1检验。前三式使平方词只依赖内容，记为 \(M_Y\)，且 \(M_Y^2=M_Y\)。

对有效的非一次吸收，取内部非一次词 \(u=gs\)。若 \(g\in E\)，选右项尾部见证，将其换序为 \(vg\)，则第七式给出

\CaseDisplay{u^2=g^2s\le gs+vg\le F.\quad\textup{(\ensuremath{1})}}

若目标首字母 \(h\in E\)，则 \(X\subseteq E\)，所以内部词 \(u\) 的首字母也在 \(E\)。由 (1) 得平方种子 \(p=u^2\le F\)。对每个尚缺的 \(x\in X\)，取全局尾部见证 \(vx\)；以第五式、\(p^2=p\) 得 \(px\le p+vx\le F\)。不断补齐得到 \(M_X=w^2\le F\)，再用 \(w\le w^2\)。

若 \(h\notin E\)，由全局内容覆盖，至少有一个全局右项以 \(h\) 开头。若内部非一次种子 \(u\) 也以 \(h\) 开头，直接从 \(u\) 起步。否则其首字母 \(g\in X\setminus\{h\}\subseteq D(w)\subseteq E\)，用 (1) 得 \(p=u^2\le F\)。由于 \(h\notin E\) 且不是 \(u\) 首字母，\(p\) 不含 \(h\)。取全局首字母见证 \(r\)：若 \(r=h\)，用第六式得 \(hp\le h+p\le F\)；若 \(r=hs\)、\(s\) 非空，用第八式得 \(hp\le hs+p\le F\)。所得非一次词以 \(h\) 开头，尾部包含于 \(D(w)\)。

最后对每个未覆盖尾变量 \(x\)，将全局见证写成 \(vx\)，把当前词写为 \(hq\)、\(q\) 非空，由第五式得 \(hqx\le hq+vx\le F\)。补齐即为 \(w\)。一次左词由纯幂见证处理。所有代入非空，因此为任意变量数的完整证明。

\par\smallskip\noindent\textit{3. 来源11217：内部尾部与全局首字母.}\quad 
取 \(B=\{z, a, b, d, t\}\)，对应同一标号，仍以链 \(d<b<a<z<t\) 定义加法。右乘 \(a\) 为恒等映射，右乘 \(b, d\) 为 \(P=(b, b, b, d, b)\)，右乘 \(z, t\) 为 \(Q=(t, t, b, d, t)\)。这些映射保持并、交换、幂等，且 \(PQ=QP=P\)；右平移兼容式按三类右因子直接成立。

词的首字母为 \(d\) 时恒取 \(d\)，为 \(b\) 时恒取 \(b\)；其余情况下，尾部含 \(b, d\) 时取 \(b\)，否则尾部含 \(z, t\) 时取 \(t\)，否则尾部全为 \(a\)、取首字母值。一次词按空尾部解释。

定理3（判据）。 对任意词 \(w\)，

\CaseDisplay{B\models w\le F\quad\Longleftrightarrow\quad
C_X=X,\qquad h(w)\in H,\qquad D(w)\subseteq E_X.\quad\textup{(\ensuremath{B})}}

必要性。 若 \(x\in X\setminus C_X\)，赋 \(x=t\)、其余 \(X\) 内变量 \(a\)、外部变量 \(b\)，右和不超过 \(a\)，左值为 \(t\)。若 \(h(w)\notin H\)，赋全部右首字母 \(d\)、其余变量 \(a\)，右值为 \(d\)，左值为 \(a\) 或 \(b\)。若 \(x\in D(w)\setminus E_X\)，赋 \(x=z\)、其余 \(X\) 内变量 \(a\)、外部变量 \(b\)，内部右项的尾部全取 \(a\)，取值为 \(a\) 或 \(z\)；外部右项取 \(b\)。所以右和不超过 \(z\)，左值为 \(t\)。

充分性。 右和为最小元 \(d\) 时，全部右首字母为 \(d\)，由全局首字母条件左词取 \(d\)。对由 \(b, a, z\) 生成的主理想，若左内容含取值 \(b, d\) 的变量，左值为 \(b\) 或 \(d\)，结论成立。

否则内部右词同样不含 \(b, d\)。右和不超过 \(b\) 已不可能，因为内部词集非空而这些词不能取 \(b, d\)。右和不超过 \(a\) 时，内部右词全部变量只能取 \(a\)；内部覆盖使左值为 \(a\)。右和不超过 \(z\) 时，内部右首字母取 \(a, z\)、尾部取 \(a\)；内部覆盖使左首字母取 \(a, z\)，内部尾部覆盖使左尾部取 \(a\)，故左值不超过 \(z\)。

定理4（七条完整基）。 下列七式完整：

\CaseDisplay{\begin{aligned}
xyz&\approx xzy,&xy^2&\approx xy,&x&\le x^2,\\
xyz&\le xy+xz,&xyz&\le xy+yz,&xy&\le x+y^2,\\
x^2y&\le y^2x+xz.&&&&
\end{aligned}\quad\textup{(\ensuremath{B+})}}

有效性由定理3检验。前六式具有内部首字母与尾部闭合第4节的句法闭合性质：内部同首字母和内部尾部覆盖蕴含吸收。

因此，当 \(h\in H_X\) 时已能推出 \(w\le F\)。若 \(h\notin H_X\)，由 \(C_X=X\) 可知 \(h\in E_X\)，结合 \(D(w)\subseteq E_X\) 得 \(X=E_X\)。选任意内部首字母 \(g\in H_X\)，取首字母为 \(g\)、尾部集合为 \(X\) 的平方词 \(Q\)。对 \(Q\) 应用上述句法闭合性质，得到 \(Q\le F\)；并且 \(Q^2=Q\)、\(c(Q)=D(Q)=X\)。

由全局首字母条件取右项 \(r\)，其首字母为 \(h\)。若 \(r=h\)，第六式给出 \(hQ\le h+Q\le F\)；若 \(r=hs\)，以最后一式代入 \(h, Q, s\) 得

\CaseDisplay{h^2Q\le Q^2h+hs=Q+hs\le F.}

两种左词均与 \(w^2\) 有相同首字母和尾部集合，故由前两式等于 \(w^2\)；再用 \(w\le w^2\) 即可。没有使用无守卫换首式，任意变量数完备性成立。

\par\smallskip\noindent\textit{4. 复核边界.}\quad 
第一组条件要求内部的非一次见证和全局尾部见证，第二组条件涉及内部内容、尾部与全局首字母。后一组在六条恒等式之外还需要一条首字母守卫式。反转各词即得对偶半环的基。

\subsection{Nine identities for 13761 and their reversals for 11237}\label{case:240}
The signature has addition and multiplication without constants. Write u \ensuremath{\leq} v for \(u+v=v\). The catalogue multiplication rows of 13761 are 44424, 01224, 22222, 33333, 44444. Addition is maximum in the chain 3 \ensuremath{<} 2 \ensuremath{<} 1 \ensuremath{<} 0 \ensuremath{<} 4. The table of 11237 has the same addition and transposed multiplication.

In addition to the ai-semiring axioms use these nine identities.

{\small\par\smallskip
\CaseLine{\(xxx = xx\qquad xyx = xyy\qquad xyzx = xyzz\)}
\CaseLine{\(x \leq  xx\)}
\CaseLine{\(ddu \leq  u + dd\qquad ddu \leq  u + ddb\)}
\CaseLine{\(myt \leq  yt + my\qquad myt \leq  yt + myb\qquad my \leq  y + myb\)}
\smallskip}

The last three permit insertion of m before y in y[t] using a word my[b]. Here and in the proof, absence of a bracketed suffix selects the displayed shorter law. With both suffixes absent the result is already an upper word. A common prefix can be multiplied on the left using distributivity. All variables occurring in an identity take nonempty arguments.

\par\smallskip\noindent\textit{Word evaluations and exact absorption.}\quad 
A word consisting of 1's gives 1. At its first entry different from 1, a 4 or 2 gives that value, and a 3 gives 3 if initial and 2 otherwise. A 0 gives 0 when final, 2 when followed immediately by 3, and 4 when followed by any other entry. The entries 2,\allowbreak{}3,\allowbreak{}4 act as left zeros. The rows give this formula on every length.

Fix a target word w and a nonempty set P of original upper words. Let C be its support. For each variable x let D be the set of letters before its first occurrence. Use a right word whose first x is preceded only by D. Impose the following test at that first x.

\textbullet{} If x is final in w, the witnessing x must either be final or be followed by an element of D union \{x\}.

\textbullet{} If x is followed by y in w, the witnessing x must have a following letter in D union \{x, y\}.

In addition P must contain some word supported in C. These conditions are equivalent to w \ensuremath{\leq} P. For containment give C value 1 and all outsiders 2. For the first-occurrence test assign D value 1 and x value 0. Give a fresh successor y the value 2, when there is one, and all remaining letters 3. The target then reaches 0 when final and 4 otherwise. The right words reaching these thresholds are precisely the admitted witnesses. For sufficiency the witness at the target head starts with that head. Therefore when the target is at least 2 its head is different from 3 and this upper word is also at least 2. Use the contained word for an all-1 evaluation. A first 4 uses its prefix witness. A first 0 giving value 0 or 4 uses the successor test; a following letter in D has value 1, a following x has value 0, and when the target gives 4 its actual successor is different from 3. These cases exhaust the chain.

\par\smallskip\noindent\textit{Normal form and soundness.}\quad 
The first three equations preserve the first-occurrence order and, at each first occurrence, the distinction between a fresh successor, an old successor (including that occurrence itself), and no successor. The word formula proves validity. The canonical word lists the first occurrences, doubling an entry precisely when its original next letter was old. To see derivability append an input letter to the normalized prefix. A fresh input is appended. For an old input x different from the last letter z, choose an earlier x. At the end use \(xzx=xzz\) or \(xUzx=xUzz\), according to whether the intervening U is absent or nonempty. Contract a resulting triple by \(xxx=xx\). The same contraction handles an appended copy of the last letter. This gives equality to the canonical word at every length. Among the invariants, final means a *first occurrence* in final position. An old final input instead contributes a doubled entry.

The successor criterion verifies the six inequalities directly. In ddu the initial d has an old successor witnessed by dd or ddb. The literal u gives containment and a final witness for u. In myt the original my or myb supplies the immediate successor of the initial m, and yt supplies those at y and t as well as containment. In my \ensuremath{\leq} y + myb the final y is witnessed by the literal. The inequality x \ensuremath{\leq} xx also follows from the square table. These are ordinary identities, so their substitutions may identify different variables.

\par\smallskip\noindent\textit{Complete derivation.}\quad 
Normalize the target and write its fresh letters as \(x_{1}\),..., \(x_{n}\). Let \(D_{i}\) consist of \(x_{1}\),..., \(x_{i-1}\). Call an index i marked when some original a \(x_{i}\) q b has a supported in \(D_{i}\) and q is in \(D_{i}\) union \{\(x_{i}\)\}. In this notation q is a single letter, a, b may be absent. Every gap (a doubled entry) of the canonical target has a marked index. Double every remaining marked entry in the canonical target, giving a word w'. We have w \ensuremath{\leq} w' by x \ensuremath{\leq} xx in contexts. It suffices to derive w' \ensuremath{\leq} P.

Treat consecutive stretches of unmarked indices ending just before a marked index, or before the end of the target. At the start of a stretch there is a prefix L spelling all preceding entries of w'. If nonempty it ends in a doubled entry, and L u \ensuremath{\leq} P is available for every original u. Initially L is absent and the originals themselves are available. Let s be the next index.

For successive indices i starting at s choose an original a \(x_{i}\) b as follows. At a marked i choose the marked witness, so \(b=q\) b' begins with an old letter. At an unmarked nonterminal i the exact criterion forces b to begin with \(x_{i+1}\). At an unmarked terminal i it forces b to be empty. In all cases a is supported in \(D_{i}\).

The part of a before its first letter from the current stretch uses only letters of L. Once a stretch letter occurs, the indices of the remaining letters of a increase strictly. Indeed at a first failure of strict increase, the immediately preceding letter has index j in [s, i) greater than every earlier stretch index in this prefix. It is a first \(x_{j}\) preceded only by \(D_{j}\) and followed by a letter of index at most j, including a possible return to L. That would mark j. Thus a has an initial part on L followed by an increasing subsequence of \(x_{s}\),..., \(x_{i-1}\). In the available L a \(x_{i}\) b the initial part on L vanishes by the word equations, since L ends doubled. When L is absent this initial part is empty.

Fill the omissions in this increasing subsequence from left to right, including omissions immediately before \(x_{i}\). If the next present letter \(y=x_{j}\) has a missing nonempty block m immediately before it, the current word has form p y t. At the earlier index j-1 we already derived a word p m y \(b_{0}\), with its entire earlier prefix in order, since j-1 is unmarked and its witness continued with y. Here p is L followed by the already filled initial part. Use the last three inequalities in the common left context p to derive p m y t. Both available suffixes may independently be absent; if both are absent the second word suffices. This step proves by induction on i that

{\small\par\smallskip
\CaseLine{\(L x_{s} \cdots  x_{i} b \leq  P\)}
\smallskip}

with the original suffix b intact. For an unmarked nonterminal i save this word for subsequent insertions. If i is unmarked and terminal, it is exactly w'.

At a marked endpoint put \(A=L\) \(x_{s}\) ... \(x_{i}\). Its available word is A q b' with q in the support of A. By the word equations this equals A A b': the first appended old letter doubles the last fresh entry, and further old letters then have no effect. For each original u use \(d=A\) in one of the two ddu laws with \(b=b'\) when present. This derives A A u \ensuremath{\leq} P. Normalization gives L' u with \(L'=A\) \(x_{i}\), providing the invariant for the next stretch. If this is the last index, take a contained original as u. All its letters are old, so L' u normalizes to \(L'=w'\). This also covers a one-letter target and an initial or final marked index.

Every construction allows arbitrarily many variables, word lengths and upper summands. Distribute both sides of a valid identity, derive each lower word absorption in both directions, and use the semilattice laws. This proves completeness of the nine identities. Reversal supplies the basis of the opposite semiring.

\subsection{Ten identities for 13815 and their reversal for 11238}\label{case:241}
Use the constant-free signature of additively idempotent semirings. The additive order in catalogue labels is 3 \ensuremath{<} 2 \ensuremath{<} 1 \ensuremath{<} 0 \ensuremath{<} 4. The multiplication rows in increasing catalogue order are 44434, 01234, 22232, 33333, 44434. A 3 anywhere gives 3. In its absence scan for the first value different from 1. If none appears the value is 1. A first 2 or 4 gives that value, a first 0 followed by a further factor gives 4, and a first 0 at the final position gives 0. This follows by induction, including literals.

\par\smallskip\noindent\textit{Absorption.}\quad 
Fix a word w and a nonempty upper polynomial P. Let C be the support of w and let U consist of the original right words supported in C. For each first occurrence of a variable x in w, let S be the set of variables preceding it. There must be a word v in U which contains x with all letters before its first x in S. When the first x in w has a successor, the first x in v must also have a successor. Call such an occurrence marked. These conditions together are equivalent to w \ensuremath{\leq} P. In particular the union of the supports of U is C.

For necessity give outsiders 3, the letters of S value 1, x value 4, and remaining letters of C value 2. The target value 4 forces the required contained prefix witness. At a marked first occurrence use \(x=0\) with the same other values: 4 now requires a successor. For sufficiency only targets with no 3 need consideration. Any contained word absorbs a target valued 2, and it gives 1 at an all-ones assignment on C. A first 4, or first 0 at the end of the word, is handled by its prefix witness. A nonfinal first 0 uses the marked witness. These arguments hold on every alphabet.

\par\smallskip\noindent\textit{The laws.}\quad 
In addition to the semiring axioms use the following ten identities. An inequality stands for the corresponding absorption equation.

{\small\par\smallskip
\CaseLine{\(W1\qquad xxy = xy\)}
\CaseLine{\(W2\qquad xyx = xyy\)}
\CaseLine{\(A\qquad x \leq  xx\)}
\CaseLine{\(B\qquad uabv \leq  ub + abv\)}
\CaseLine{\(C\qquad ubv \leq  ub + bv\)}
\CaseLine{\(D\qquad xxy \leq  xx + y\)}
\CaseLine{\(E\qquad uzt \leq  zt + ut\)}
\CaseLine{\(F\qquad puzt \leq  pzt + ut\)}
\CaseLine{\(G\qquad uzts \leq  zts + ut\)}
\CaseLine{\(H\qquad puzts \leq  pzts + ut\)}
\smallskip}

The table formula gives W1, W2 and A. The absorption criterion checks B through H directly on their displayed variables. W1, W2 reduce a word to its distinct letters in order of first occurrence, with the last listed letter appended once if its raw final letter is repeated. For a separated internal repeat use \(xYxZ=xYYZ=xYZ\) with nonempty Y, Z, and W1 for an adjacent internal repeat. A final x separated from its previous occurrence by Y changes xYx to xYY by W2, and the repeated internal letters of the second Y are then removed. Thus normal forms can be used in either direction. Call a word saturated when its final letter occurs earlier. In particular \(g^{2}\) is equivalent to g whenever g is saturated.

\par\smallskip\noindent\textit{Obtaining full support from the original words.}\quad 
If a letter of C has no marked occurrence in any member of U, the test forces it to be the simple final letter s of w. There is at most one such letter. In this case every member of U containing s ends in its unique occurrence. Otherwise every letter is globally marked.

Start a construction at each original u in U. Stop if the accumulated word g is saturated, or if it ends in the letter s with no marked witness. In the remaining case its simple final letter b has a marked witness v in U, written \(v=L\) b R with R nonempty. Write \(g=K\) b. If K is empty, take v directly. Otherwise derive K L b R from g and v by B when L is nonempty and by C when L is empty. The resulting support contains both original supports. Every old letter now has a marked occurrence in the resulting word: the old proper prefix precedes the nonempty v and b itself has a successor in v. Any new simple final letter therefore lies outside the old support. With a finite target support this process terminates. In particular a cycle of choices terminates as soon as its prospective terminal letter is already present. This does not assume an acyclic witness graph.

The terminal words from all starts jointly cover C. Combine the saturated ones first. If g is saturated, D at \(x=g\), \(y=v\) gives g v \ensuremath{\leq} g+v since \(g^{2}=g\). Two saturated inputs give a saturated result. If s exists every saturated word omits s by the absorption test (s has no global mark), and every other terminal word ends in its unique s. Two such words K s and L s combine to K L s by E. An empty K or L means that one input is the literal s and the other already contains its support. If both kinds occur, prepend the combined saturated word to the combined simple-final word by D. The result is a full-support G \ensuremath{\leq} P. It is saturated when all letters are globally marked; otherwise it ends in its unique s. The second case does have a terminal word carrying s, by coverage.

\par\smallskip\noindent\textit{Ordering the full-support word.}\quad 
Work with G in normal form. Once its first-occurrence list starts with the processed distinct prefix p of w, choose a contained marked witness \(v=U'\) T for the next marked first occurrence x, where U' ends in the first x of v and T is nonempty. All earlier letters of U' are in p. For an empty p use E or G below, and for nonempty p use F or H. Here G as a *law label* denotes the seventh displayed absorption (uzts \ensuremath{\leq} zts + ut); the word G always denotes the current full-support word.

If the full-support word is saturated, take \(z=G\) and \(t=T\). Then [p] z t normalizes to G, since p is its initial list segment and T is supported in G. Use E or F with \(u=U'\) and the original u \(t=v\). The left word [p] U' G T is saturated with full support, and its list of first occurrences starts with p x.

If the full-support word has a simple final s, write it as Z s. Here Z is nonempty when |C| \ensuremath{\geq} 2. The proper prefix U' of the original witness avoids s. If T contains s it ends in its unique s. Take \(z=Z\), \(t=T\) and use E or F: [p] Z T has the same normal form as G. If T omits s, keep an additional final s and use G or H with \(z=Z\), \(t=T\). Its first right word [p] Z T s again reduces to G. In both cases the other right word is exactly v; the resulting word has full support, simple final s, and first-occurrence prefix p x. Brackets denote a case distinction between displayed laws, never substitution of an empty word.

These insertions can process all marked first occurrences of w. For a saturated w the result is its normal form. For a simple final s of w, once the other distinct variables have been ordered, full support places s last in the first-occurrence list. The resulting word is either w in normal form or the normal form of \(w^{2}\); A recovers w in the latter case. This also works when the simple final has a global marked witness but the only prefix-compatible witness for it is unmarked. No insertion at that position is required.

On singleton support a valid literal is recovered from any original contained power by A and W1. For a target of length at least two on singleton support, a marked contained witness is a power of length at least two; all such powers are equivalent by W1. This completes the argument. Distribute arbitrary terms as sums of nonempty words and apply the construction to their summands in both directions. Thus the ten laws are complete for all variable counts, lengths and numbers of summands. Multiplicative reversal gives the basis for 11238 by the certified opposite-table map.

\subsection{A support-size obstruction for 11269, 11270 and 13161}\label{case:242}
Use nonempty finite sums of nonempty words in the constant-free additively idempotent semiring signature. Labels below are table entries. The common multiplication for 11269 and 11270 is 40244 / 41234 / 22222 / 44244 / 44244. A 2 anywhere gives 2. Otherwise a word is 1 on all ones, 0 on an initial 0 followed by ones, 3 on ones followed by a final 3, and 4 otherwise. Indeed on the alphabet excluding 2 a current 1 produces 1,\allowbreak{}3 or 4 after the first position, a current 0 survives exactly on a following 1, and current 3 or 4 always gives 4 on a further letter. The formulas include literals.

For target w and upper polynomial P use \(C=C(w)\), and let U consist of all its original words supported in C. Write K, J, T, F for the unions of their content, nonfinal positions, noninitial positions and finals. Also write A(w), B(w) for the nonfinal and noninitial letters of w. In both algebras w \ensuremath{\leq} P requires \(K=C\) and B(w) contained in T. The further condition is A(w) contained in J for 11269. For 11270 it is either J meeting F, or both A(w) contained in J and f(w) in F. Taken together these conditions are sufficient.

To prove this, 2 annihilates multiplication and is least in both orders. Outside variables may be given 2 in counterassignments. When the target avoids 2, use only U for sufficiency. At the least bound one contained word is needed. At bound 1 all nonbottom upper words must be all ones, forcing the cover: use 4 at a missing letter, 1 elsewhere in C. At bound 0, whose ideal is \{2,\allowbreak{}1,\allowbreak{}0\}, covered letters are 1 or 0 and all letters of T are forced to 1. Giving 0 to a missing noninitial letter proves necessity. At bound 3 the ideal for 11269 is \{2,\allowbreak{}1,\allowbreak{}3\}, giving the reversed role requirement by assigning a missing nonfinal letter 3. For 11270 the ideal is \{2,\allowbreak{}3\}. Every high contained word must end in 3 with prefix all ones, impossible if J meets F. Otherwise J, F partition C and force the assignment 1 on J and 3 on F. The displayed condition is exactly what makes w have value 3. This exhausts the bounds.

\par\smallskip\noindent\textit{Two cycle families.}\quad 
Take \(n=2m\) with m at least 3, and distinct letters \(a_{i}\), \(b_{i}\), \(z_{i}\) for all indices modulo n. Each private marker \(z_{i}\) is used in the single original word \(p_{i}\). In 11269 put

{\small\par\smallskip
\CaseLine{\(p_{i} = a_{i} z_{i} a_{i+1} b_{i+2} b_{i}.\)}
\smallskip}

In 11270 instead put

\(p_{i}\) = \(a_{i}\) \(z_{i}\) \(a_{i+1}\) \(b_{i}\)                  for i different from m, \(p_{m}\) = \(a_{m}\) \(z_{m}\) \(a_{m+1}\) \(b_{0}\) \(b_{m}\).

For each family write P for the sum of these words and r for their concatenation in index order. The identities r \ensuremath{\leq} P hold for every such n. For the full support the noninitial letters cover every variable in both families. For the first family every variable also occurs nonfinally; for the second \(b_{0}\) occurs nonfinally and finally. Apply the exact criteria.

Let v be any raw word absorbed by P. There is a nonempty set I of indices whose original words are contained in its support, and its support is exactly their union by \(K=C\). It contains exactly the markers indexed by I. For either family each i in I with i-1 outside I makes \(a_{i}\) unavailable noninitially among these originals. Every such letter must be the unique simple head of v. Thus any proper I has just one cyclic run of consecutive indices. In the first family each i in I with i-2 outside I likewise requires \(b_{i}\) to be the unique simple final. If I is a cyclic run of length between 2 and n-2 its first two indices have this defect. The only possible sizes are 1, n-1 and n.

In the second family, if I contains neither the crossing pair of indices together (that is, if 0 or m is missing), none of its original finals occurs nonfinally. Two distinct indices would then require two distinct simple finals of v. If I contains both 0 and m, a cyclic run containing them has length at least m+1. Thus in both families any absorbed word has exactly one marker or at least m+1 distinct markers. In the one-marker case the word with marker \(z_{i}\) has exactly the support of \(p_{i}\), starts with \(a_{i}\) and ends with \(b_{i}\). The head and final are forced by the corresponding noninitial and nonfinal defects of the singleton. There are no other possible letters, even for raw words with arbitrarily many repetitions.

\par\smallskip\noindent\textit{Invariance under bounded-variable identities.}\quad 
Suppose a valid absorption v \ensuremath{\leq} Q on at most k variables is given any nonempty polynomial substitution. Assume every expanded word of the substituted Q is one of the single-marker words described above. Fix any expansion of v. Each source variable appearing noninitially in v appears noninitially in a source monomial of Q contained in the source content of v. This follows from the unconditional noninitial coverage condition on the *source* absorption. Keep the other image choices in that source monomial fixed and vary this one position. All resulting words have the same first letter because of the nonempty preceding position. A single-marker word's first letter identifies its index. Consequently the entire image polynomial of this variable is supported inside a single block C(\(p_{i}\)), containing at most one marker. This also holds when the variable occurs repeatedly: expansion chooses its image word independently at each occurrence.

There is at most one further variable in v, occurring once at its initial position. By source coverage some monomial of Q contains it. Fix the other choices there and use the image word selected in the expansion of v. This yields a single-marker word, so the selected factor contains at most one marker. Its other image choices need no joint bound. Thus the chosen expansion of v contains at most k distinct markers, however large any polynomial substitution may be.

If a finite basis had at most k variables per identity, take \(m>k+2\). In a derivation starting at P, every intermediate polynomial equals P in the algebra, so all its expanded words are absorbed by P. At the first step introducing a word with more than one marker, distribute the one-hole context. The new word comes from one contribution s sigma(V) t, where s, t are fixed words or omitted. The old contribution s sigma(Q) t has only single-marker words. Introduce at most two fresh source variables for the nonempty context words; the contextual identity is still valid. Its individual new source monomial is absorbed by its old source polynomial. The preceding bound limits this new expanded word to k+2 distinct markers, contradicting the gap. Base semiring axioms preserve free-polynomial expansion. This treats both orientations, addition of other summands, polynomial substitutions, and arbitrary multiplicative contexts. Since the valid critical identity introduces r with all n markers, no finite basis exists.

Reversing words gives the same result for 13161, the opposite of

11270. Finally, every absorption valid in 11269 is valid in 11270: when J misses F, the nonempty F must lie in C minus A(w) and hence consists of the target final. The inclusion Var(11270) in Var(11269) is strict. The absorption yyx \ensuremath{\leq} xy+xx satisfies the 11270 criterion via the crossing x. In 11269 the assignment \(x=1\), \(y=3\) gives left 4 and right 3. Thus 11269 has a proper nonfinitely based subvariety. The two other NFB decisions here leave their limit property open.

\subsection{Nonlinear contained covers and finite derivations}\label{case:243}
Let A be catalogue entry 12496. Work in the signature (+, multiplication), without constants, relative to the axioms of additively idempotent semirings. Write \(w\leq P\) for w+P approximately equal to P. Addition is maximum on \(3<2<1<0<4\). The tables are

{\small\par\smallskip
\noindent addition       multiplication\par
\CaseLine{\ttfamily 0 0 0 0 4\qquad 4 4 0 3 4}
\CaseLine{\ttfamily 0 1 1 1 4\qquad 4 4 0 3 4}
\CaseLine{\ttfamily 0 1 2 2 4\qquad 4 4 2 3 4}
\CaseLine{\ttfamily 0 1 2 3 4\qquad 3 3 3 3 3}
\CaseLine{\ttfamily 4 4 4 4 4\qquad 4 4 4 3 4}
\smallskip}

Numerical labels are elements, not constants in terms. For a word u write c(u), h(u), t(u) for its content, head and set of tail variables; t(y) is empty when y is linear.

\par\smallskip\noindent\textit{Word semantics.}\quad 
The element 3 is multiplicatively absorbing. With no factor 3, a nonlinear word is 4 if a tail factor lies in \{0,\allowbreak{}1,\allowbreak{}4\}; otherwise every tail factor is 2 and the result is r(h), where \(r(0)=r(1)=0\), \(r(2)=2\) and \(r(4)=4\). This rule agrees under both associations of a triple: a 3 absorbs in either case; with no 3, a last factor outside \{2\} gives 4, and a last factor 2 reduces to whether the middle factor is 2. If it is, \(r(r(h))=r(h)\); otherwise both associations give 4. Every row and column is increasing in the chain order, so both distributive laws hold. Thus A is a semiring and induction gives this rule at arbitrary length. Nonlinear words never have value 1.

The word identities

\CaseDisplay{xyz\approx xzy,\qquad xy^2\approx xy    \quad\textup{(\ensuremath{W})}}

hold. They normalize a nonlinear word to its literal head followed by its sorted nonempty tail set, retaining the head twice when it also occurs in the tail. Linear words remain linear. They also give \(x^{n}\) approximately \(x^{2}\) for every \(n\geq 2\). The inequality \(x\leq x^{2}\) is valid by the displayed table.

\par\smallskip\noindent\textit{Exact absorption criterion.}\quad 
For a nonempty polynomial P and target content C define

\CaseDisplay{N_C=\{u\in P: |u|\geq2,\ c(u)\subseteq C\}.}

Lemma. A linear y is absorbed by P exactly when some original word of P has content \{y\}. A nonlinear w of content C is absorbed exactly when

\CaseDisplay{C\subseteq\bigcup_{u\in N_C}c(u),\qquad
 t(w)\subseteq\bigcup_{u\in N_C}t(u).                       \quad\textup{(\ensuremath{1})}}

Proof. If there is no right word of content \{y\}, assign \(y=2\) and all outsiders 3. Every right word has value 3 but \(y=2\). If such a word is present it is y or a power \(y^{n}\) with \(n\geq 2\); powers reduce by (W) to \(y^{2}\), and the valid inequality \(y\leq y^{2}\) suffices.

For nonlinear necessity, if x in C is absent from the first union, assign \(x=1\), all other C-variables 2 and all outsiders 3. Every noncontained right word is 3; every contained nonlinear word is 2; every right linear word is at most 1. But w is 0 when x is its simple head and 4 when x occurs in its tail, so w is not absorbed. If x in t(w) is missing from the second union, assign \(x=0\), other C-variables 2 and outsiders 3. A contained nonlinear word has x only at its simple head or omits it, and so is 0 or 2. All linear right words are at most 0; all noncontained words are

3. The target is 4, again separating it.

For sufficiency we may use \(N_{C}\) alone. If a variable of C is 3, the target is the least element 3. Otherwise every C-variable is in \{0,\allowbreak{}1,\allowbreak{}2,\allowbreak{}4\}. A variable assigned 4 occurs in a contained nonlinear word by content coverage and makes it 4. A tail variable assigned 0 or 1 occurs in a contained nonlinear tail by tail coverage and likewise makes it 4. In the remaining case all target tail variables are 2. If its head is 0 or 1, some contained nonlinear word contains that head and has value 0 or 4; hence \(w=0\) is absorbed. If its head is 2, every C-variable is 2 and a contained nonlinear witness gives

2. Such a witness exists by content coverage. This exhausts the cases. \(\square\)

The exclusion of linear right words from \(N_{C}\) is essential: a linear variable assigned 1 cannot witness a nonlinear target with value 0. In contrast to source 11019, a missing literal linear target may be absorbed by a power. These are different exact criteria despite the common form of their local guard construction.

\par\smallskip\noindent\textit{Local guard construction.}\quad 
For a nonlinear normal word u, c(u) minus t(u) is either empty or \{h(u)\}. Two known nonlinear words u, v can be merged if

\CaseDisplay{|(c(u)\cup c(v))\setminus(t(u)\cup t(v))|\leq1.             \quad\textup{(\ensuremath{2})}}

Put \(E=c(u)\) union c(v), \(T=t(u)\) union t(v). If E minus T is \{r\}, choose r as head; if it is empty choose any head in E. The word q consisting of that head followed by the ordered T is nonlinear, has content E and tail T, and satisfies \(q\leq u+v\) by the lemma. Both witnesses are nonlinear and contained in E.

Fix a valid nonlinear \(w\leq P\) and put \(C=c(w)\), \(h=h(w)\). Starting with each original u in \(N_{C}\), repeatedly change a current guard q as follows. Stop if its missing set is empty or \{h\}. Otherwise its missing head r differs from h, so r belongs to t(w). Choose an original v in \(N_{C}\) with r in t(v), which exists by (1), and merge q, v using (2). Every old content variable is now in the new tail: the only possible previously missing variable was r. If the new missing head p is again different from h, then p lies outside all old content. Thus every continuing step strictly grows content and introduces a new missing head. Since content stays inside finite C, the process terminates after at most |C| steps.

All resulting guards have missing set contained in \{h\}, so they can be merged successively by (2). Content coverage makes their final content exactly C, while their tail contains every original tail and therefore t(w). The lemma then supplies the one-right local step \(w\leq Q\). Every intermediate guard was derived below P before its use in a later step. Every local absorption uses at most two right words. All words in these nonlinear derivations are nonlinear. A valid linear absorption is a literal inclusion or follows from \(y\leq y^{2}\) after (W).

\par\smallskip\noindent\textit{Reversible compression and finite basis.}\quad 
In a local absorption there are at most three words. Reserve each head individually, at most three variables, and group remaining variables by their nonempty incidence patterns across the three tail sets, at most seven groups. Identify each group with one symbol and normalize using (W). This is a substitution preserving validity and gives at most ten variables and word length eleven. For its inverse send each group symbol to the nonempty ordered product of its members. Members have identical tail incidence and never occur at reserved heads. Tail permutations and duplicate removal recover every original normal word. Linear words stay reserved as single variables; no empty substitution or unit occurs.

Theorem. The identities (W), together with all valid normal absorptions over a fixed ten-letter alphabet, of word length at most eleven and at most two right words, form a finite basis for A relative to the semiring axioms. Reversing all words gives the basis of its multiplicative opposite.

Proof. The family is finite and its membership is decided by the lemma. It contains \(y\leq y^{2}\). Reversible compression derives each local nonlinear step from the family. The guard construction then derives every valid absorption without a bound on variables, lengths or original summands. Linear absorptions were handled above. Distribute arbitrary terms and derive the absorption of every summand of either side into the other side of a valid identity, obtaining both inequalities. This proves completeness. Word reversal transports the derivations to the opposite. \(\square\)

The finite family is specified exactly; no small enumeration, minimality or hereditary property is claimed. Computational checks are auxiliary. Any transfer to other catalogue types requires independent complete HSP verification in both directions.

\subsection{Ten identities for an ordered-tail chain semiring}\label{case:244}
Work in the signature (+, multiplication), without constants, relative to the additively idempotent semiring axioms. Write \(w\leq P\) for w+P approximately equal to P. Let A be catalogue entry 12500. Addition is maximum on \(3<2<1<0<4\), with tables in numerical element order:

{\small\par\smallskip
\noindent addition       multiplication\par
\CaseLine{\ttfamily 0 0 0 0 4\qquad 4 4 0 3 4}
\CaseLine{\ttfamily 0 1 1 1 4\qquad 4 4 0 3 4}
\CaseLine{\ttfamily 0 1 2 2 4\qquad 4 4 2 3 4}
\CaseLine{\ttfamily 0 1 2 3 4\qquad 4 4 3 3 4}
\CaseLine{\ttfamily 4 4 4 4 4\qquad 4 4 4 3 4}
\smallskip}

Element labels are not constants in terms. Products of words always use nonempty substitutions.

Theorem. The following ten identities form a complete relative basis for A:

\CaseDisplay{\begin{aligned}
 xyzy&\approx xzy, &(1)\qquad xy^2&\approx xy, &(2)\\
 xyz&\leq xy+tz, &(3)\qquad xy&\leq x^2y, &(4)\\
 xy&\leq yxy, &(5)\qquad xy&\leq y^2+zxy, &(6)\\
 xy&\leq x+zxy, &(7)\qquad xy&\leq y+zxy, &(8)\\
 xyz&\leq y+txyz, &(9)\qquad x&\leq x^2. &(10)
\end{aligned}}

The opposite semiring has the reversed basis. Identity (3) cannot be replaced here by the superficially similar \(xy\leq x+zy\): \(x=1\), \(y=z=2\) is a counterexample to that proposed replacement.

\par\smallskip\noindent\textit{Values of words and ordered normal forms.}\quad 
Right multiplication by 2 maps 1 to 0 and fixes 0,\allowbreak{}2,\allowbreak{}3,\allowbreak{}4. Denote this idempotent map by r. Each right multiplication by j in \{0,\allowbreak{}1,\allowbreak{}3,\allowbreak{}4\} is constant, with value 3 for \(j=3\) and 4 otherwise. The multiplication rows and columns are increasing in the additive chain, so they preserve its binary joins. For associativity, if the last factor is not 2 both associations equal its constant right-action value. If the last factor is 2 and the middle factor is not 2, its image under r has the same constant right action as the middle factor. If both are 2, r squared equals r. Thus the algebra is a semiring.

For a nonlinear word, look at the last tail position assigned an element outside \{2\}. If there is one, its value determines the product: 3 gives 3, and 0,\allowbreak{}1 or 4 gives 4. If there is none, the value is r of the assigned head. This follows by induction and holds for arbitrary word length. Nonlinear values never equal 1.

The first two identities hold by this formula. They remove every earlier occurrence of a repeated tail variable, leaving the literal head followed by the distinct tail variables in last-occurrence order. For a nonadjacent repetition use (1) with a nonempty prefix and intervening substring; for an adjacent repetition use (2). Linear words remain linear, and \(x^{n}\) approximately equals \(x^{2}\) for \(n\geq 2\). The order of different tail variables cannot in general be permuted.

For a word u write c(u), h(u) for its content and head. For a variable x occurring in u, let \(s_{u}\)(x) be the suffix after its last occurrence, which may be empty. Empty suffixes here describe positions; they are never substituted for a term variable. All criteria below may be applied after the two word normalizations.

\par\smallskip\noindent\textit{Exact absorption criterion.}\quad 
Lemma. A linear y is absorbed by a polynomial P exactly when an original right word has content \{y\}. For a nonlinear target \(w=h\) t with its tail in last-occurrence normal form, absorption holds exactly under these conditions:

\textbullet{} (C) An original word of P has content contained in c(w), including a possible linear word.

\textbullet{} (T) For every target tail variable x, an original nonlinear right word u has a tail occurrence of x and c(\(s_{u}\)(x)) is contained in the target suffix after x.

\textbullet{} (H) If h does not occur in t, an original nonlinear right word v contains h and c(\(s_{v}\)(h)) is contained in c(t).

Proof. For a missing pure-y right word, assign \(y=2\) and all other variables 3. Every right summand is 3, while \(y=2\). Conversely a pure-y word is y or a power, and \(y\leq y^{2}\) and the word formula show absorption. The inequality \(y\leq y^{2}\) follows directly from the table and will be part of the displayed basis.

For failure of (C), assign c(w) the value 2 and outsiders 3. The target is 2. Each right word contains an outsider: a linear word is 3, and a nonlinear word with a 3 in its head or tail and all other values 2 is 3. Thus \(P=3\).

For failure of (T) at x, assign \(x=0\), the variables in the target suffix after x the value 2, and every other variable 3. The target is 4, since its last tail reset is \(x=0\). A right word with x in its tail has a variable assigned 3 after the last x, by failure of the witness condition, so it is 3. If x appears only at the head, a nonlinear right word is 0 or 3. A word with no x is 2 or 3, and every right linear is at most 0. Hence \(P\leq 0\).

For failure of (H), assign the simple head \(h=1\), every target tail variable 2, and outsiders 3. The target is 0. Every right nonlinear word containing h has a 3 after its last h, by failure of (H), so it is 3. Those omitting h have value 2 or 3; right linears are at most 1. Thus \(P\leq 1<0\). This explains why a literal h cannot by itself supply (H), even though it can supply (C).

For sufficiency, if the target has a last non-2 tail factor, a value 3 makes it the least element. Otherwise that factor is assigned 0,\allowbreak{}1 or 4, and (T) supplies a nonlinear right word whose last occurrence of that variable is in its tail and is followed only by variables assigned 2. That right word is 4, absorbing the target. It remains that every target tail factor is 2. A head 3 gives 3 automatically; a head 2 gives 2, which (C) witnesses by a word entirely in target variables, all assigned 2. If the head is 0 or 1, it is simple and (H) gives a right word with value 0 or 4: the head occurrence gives 0 if it is the right head and 4 if it is in a right tail. If the head is 4 the same witness gives 4. This proves sufficiency for every assignment. \(\square\)

The lemma proves soundness of (3)--(9) by their displayed right words: (3) uses xy for (C),(H), xy for the y-suffix witness and tz for the z-suffix witness; (4),(5) have their single contained right word; (6) uses \(y^{2}\) for (C) and zxy for (H); (7),(8),(9) use their literal right variable for (C) and the nonlinear right word for every ordered occurrence condition. Inequality (10) was checked from the table. Hence all ten identities are sound.

\par\smallskip\noindent\textit{Adding a common tail to nonlinear original words.}\quad 
Fix a valid nonlinear absorption \(w\leq P\). Normalize all words and write \(w=h\) t. Let N be the original nonlinear right words. It is nonempty by (T). We shall derive u \(t\leq P\) for every u in N, by building the tail from right to left.

Suppose a suffix S has already been appended to all u in N. For the next target letter x, (T) supplies an original v in N with last tail occurrence x, a nonempty prefix R before that occurrence, and suffix contained in S. Word normalization gives vS approximately RxS, so RxS is already below P. Split any original u in N into nonempty A B, for example at its final letter. Apply (3) with \(x=A\), \(y=B\), \(z=x\), \(t=R\) in the right context S:

\CaseDisplay{u x S\leq u S+R x S.}

Both right terms were already derived. Continuing proves u \(t\leq P\) for every original nonlinear u. No literal original right word is extended by this operation: the nonempty split of u matters. If any such u t normalizes to w, the derivation is complete.

\par\smallskip\noindent\textit{Recovering the head from a nonlinear contained witness.}\quad 
First suppose (C) has a nonlinear original witness u, with head k. If h occurs in the target tail t, then \(c(w)=c(t)\), and u t normalizes to k t. If \(k=h\) this is w; otherwise split \(t=A\) B where B begins with the occurrence k, allowing A to be empty. Identity (5), with \(x=h\) A and \(y=B\), gives \(w\leq B\) h A B, whose normal form is k t. Both substituted words are nonempty.

Now suppose h is simple. A witness for (H) either begins with h, in which case its extension by t is w directly, or has a nonempty prefix R before its last h, supplying R h \(t\leq P\). If the contained u includes h and \(k=h\), its extension is w when h was simple in u; if h occurred in its tail, that extension normalizes to h h t, and (4) gives \(w\leq h\) h t. If k differs from h, split \(t=A\) B at k and use (5) again to give \(w\leq k\) h t, the normal form of u t.

If the contained nonlinear u omits h, its head k belongs to t. Writing \(t=A\) B at k, apply (5) with \(x=g\) A, \(y=B\), where g is the first variable of t. The left g A B has normal form t t, and the right B g A B has normal form k t, already known. Therefore the square \(t^{2}\) is below P. Combine it with the independent head witness R h t by (6), with \(x=h\), \(y=t\), \(z=R\), to derive \(w\leq t^{2}+R\) h \(t\leq P\). All substituted words are nonempty. Empty A denotes an omitted substring in a longer nonempty substituted word, not a unit term.

\par\smallskip\noindent\textit{Recovering content from a literal witness.}\quad 
If (C) has no nonlinear witness, it supplies a literal g in c(w). We already have a word R w below P, for some nonempty R. When h is simple this comes from the nonlinear (H) witness above; the case with h as that witness's head was already completed. When h occurs in t, take any nonlinear original u. The extended word u t equals u h t after tail normalization, so \(R=u\) works.

If \(g=h\), apply (7), with \(x=h\), \(y=t\), \(z=R\), to derive \(w\leq g+R\) w. If g is the final variable of w, apply (8) with x the nonempty prefix before g, \(y=g\), \(z=R\). If g is an internal variable, split \(w=A\) g B with nonempty A, B and apply (9) with \(x=A\), \(y=g\), \(z=B\), \(t=R\). In each case both right terms are already below P. This handles the literal content witnesses without using the false inequality \(xy\leq x+zy\).

\par\smallskip\noindent\textit{Completeness.}\quad 
For a valid linear absorption, normalize its pure-variable right witness to y or \(y^{2}\) and use (10) if needed. For a valid nonlinear absorption, the common-tail construction and the exhaustive head/ content cases above derive it. There is no bound on the number of target tail variables or original right words. Finally expand arbitrary terms by distributivity. In a valid identity every summand of either polynomial is absorbed by the other; derive these absorptions and sum them to obtain both inequalities. This proves arbitrary-variable equational completeness. Word reversal gives the opposite result. \(\square\)

Finite tests are auxiliary checks of the proof and its actual substitutions. Further catalogue types need independent complete HSP checks in both directions. No minimality or hereditary finite-basis property is asserted.

\subsection{Eleven identities for 11375: suffix construction and literal head guards}\label{case:245}
Signature (+, multiplication), without constants, relative to the additively idempotent semiring axioms. The completeness argument is unrestricted in variables, word length and number of summands.

\par\smallskip\noindent\textit{Tables and word values.}\quad 
Addition: 00004 / 01014 / 00224 / 01234 / 44444. Multiplication: 40434 / 41434 / 40434 / 43434 / 44434. The additive order has bottom 3, incomparable elements 1 and 2 below \(0=1+2\), and top 4. For a nonlinear word hT, scan its tail from the right. The last tail value different from 1 determines the result: 3 gives 3; 0,\allowbreak{}2,\allowbreak{}4 give 4. If all tail values are 1, the result is f(h), where \(f(0)=0\), \(f(1)=1\), \(f(2)=0\), \(f(3)=3\), \(f(4)=4\). Literals keep their value. This follows from the columns: right 1 applies the idempotent map f, right 3 is constant 3, and right 0,\allowbreak{}2,\allowbreak{}4 are constant 4. Both 3 and 4 are fixed by f. This proves the formula by induction, with no length bound.

\par\smallskip\noindent\textit{Exact criterion with original witnesses.}\quad 
Normalize only the tail by retaining last occurrences. Write a nonlinear target as hT, with distinct letters in the nonempty tail T. A word is contained in a set if all its letters lie in that set. For original right summands P, the absorption \(hT\leq P\) is equivalent to:

(C) an original contained word in c(hT); (O) for each x in T, an original right word has x strictly after its head, with the suffix after its last x contained in the suffix after x in T; (H) if h is absent from T, either an original NONLINEAR word has an h occurrence with following content contained in c(T), or the original literal h occurs together with an original contained word in c(T).

The H occurrence may be the head or a tail occurrence. O requires a true tail occurrence. The prefixes of O and tail-H witnesses can contain variables outside the target. A literal target h is absorbed exactly when an original right word is a positive power of h.

Necessity of C: assign target letters 1 and outsiders 3. If C fails every right word has value 3, whereas the target is 1. For O failing at x, assign \(x=0\), its target suffix 1, and every other variable 3. The target is 4. Every right tail x has a later 3, and all other right words have value in \{0,\allowbreak{}1,\allowbreak{}3\}; their join is at most 0, below 4. For H, assign \(h=2\), \(T=1\) and outsiders 3. The target is 0. A right nonlinear word reaches 0 or 4 exactly through the stated H occurrence. Without such a word, the join reaches 0 exactly when both literal \(h=2\) and a contained-T word of value 1 occur. Thus failure of H separates the target. For a literal h without a pure witness, \(h=1\) and \(outsiders=3\) separates it.

For sufficiency, target 3 is bottom. A tail reset producing 4 is covered by O, since the witness has that reset followed only by 1s. If all target tail values are 1 and h also is 1, C supplies a word of value 1. If \(h=0\) or 2, H gives a nonlinear value 0 or 4, or literal h together with a contained-T value 1, whose join is 0. If \(h=4\), H gives 4. If h belongs to T then h must be 1 in this no-reset case, already covered by C. These cases exhaust all assignments. Positive powers satisfy \(h\leq h\) squared, so the literal criterion is also sufficient.

\par\smallskip\noindent\textit{The eleven laws.}\quad 
(W1) \(xyzy=xzy\) (W2) \(xyy=xy\) (W3) \(xyx=yyx\) (S)  \(x\leq xx\) (A)  \(xyz\leq xy+tz\) (B)  \(xyz\leq xz+tyz\) (L)  \(xy\leq x+yy\) (G0) \(xx\leq x+yx\) (G1) \(xxz\leq x+yxz\) (G2) \(xyx\leq x+zyx\) (G3) \(xyxz\leq x+tyxz\)

Word values prove W1-W3 and S; C/O/H prove the remaining laws. Every identity variable denotes a nonempty word. An absent context is never an empty-word substitution. In particular \(xxy=xy\) and \(xy\leq x+zy\) FAIL in this algebra and are not being used.

W1 deletes a nonfinal tail occurrence of a variable when its intervening block is nonempty; W2 handles adjacent occurrences. Thus retain the head and the tail in last-occurrence order. If the head h belongs to this tail, write the tail A h B. When A is nonempty, W3 changes h A h B to A A h B; tail deletions then make its head the first variable of the original tail. When A is absent the head already is that first variable. Hence every nonlinear word has a canonical form hT when h is absent from T, and tT when h belongs to T, where t is the first letter of T. In particular S squared has canonical form sS for any nonempty last-occurrence list S.

\par\smallskip\noindent\textit{Build the suffix on nonlinear original words.}\quad 
Let N be the set of nonlinear original summands. O for the last target tail variable implies N is nonempty. Build T from right to left, maintaining \(uS\leq P\) for EVERY original u in N. Initially S is absent. To start with final target variable x, its O witness is \(q=R\) x with R nonempty. For \(u=UV\) with U, V nonempty, A gives \(ux\leq u+R\) \(x\leq P\). For a later step, S is nonempty and an O witness has \(q=R\) x Q, with R nonempty and c(Q) contained in c(S). The invariant gives \(qS\leq P\) and W1-W2 give \(qS=R\) x S. B gives \(uxS\leq uS+R\) x \(S\leq P\) for every u in N. This proves the invariant with the next suffix xS. Only original nonlinear words receive the appended suffix; literal witnesses remain available unchanged.

\par\smallskip\noindent\textit{Guard lemma, including contained literals.}\quad 
Suppose the invariant \(uS\leq P\) holds for all u in N, where S is nonempty, and P has an original contained-S word. Then S \(squared\leq P\). If a nonlinear original r is contained in S, rS normalizes to S squared. Otherwise take a contained literal r and any v in N. Write \(S=A\) r B. If both contexts are absent use G0: \(rr\leq r+v\) r. If only A is absent use G1: \(rrB\leq r+v\) r B. If only B is absent use G2: r A \(r\leq r+v\) A r. If neither is absent use G3: r A r \(B\leq r+v\) A r B. Every left word has canonical form S squared. All variable substitutions in the applicable case are nonempty, and each right term is original r or a maintained vS. This proves all four endpoint cases of the lemma.

\par\smallskip\noindent\textit{Recover the target and prove completeness.}\quad 
After constructing T, consider hT. If h belongs to T, C supplies a contained-T original word; the guard lemma gives T squared, whose canonical form equals hT.

Suppose h is absent from T. If H has a nonlinear witness with h at the head and following content contained in T, its maintained uT equals hT. Otherwise a nonlinear H witness has \(q=R\) h Q with R nonempty and c(Q) contained in T. The maintained qT equals R h T. Apply B once more to obtain u h \(T\leq P\) for every u in N. Now C and the guard lemma with \(S=hT\) give (hT) \(squared\leq P\), equivalent to hhT. S, in the nonempty right context T, gives \(hT\leq hhT\leq P\). If H instead uses literal h and a contained-T original word, the guard lemma gives T \(squared\leq P\), and L gives \(hT\leq h+T\) \(squared\leq P\).

A literal target is obtained from its original pure witness by W2 and S. Thus every valid word absorption is derived. Distributivity expands arbitrary terms into nonempty polynomials; deriving every summand absorption in both directions proves every identity. These eleven identities are therefore a finite basis for 11375. Opposite and catalogue transfers require complete HSP both ways. Finite tests and replay supplement this unrestricted proof.

\subsection{A finite basis for 11296 from original contained tail coverage}\label{case:246}
Signature (+, multiplication), no constants, relative to the additively idempotent semiring axioms. The finite family is specified, not enumerated. Completeness below is unrestricted in variables, lengths and numbers of summands.

Addition: 00004,\allowbreak{}01114,\allowbreak{}01224,\allowbreak{}01234,\allowbreak{}44444; order \(3<2<1<0<4\). Multiplication: 40334,\allowbreak{}41334,\allowbreak{}22332,\allowbreak{}33333,\allowbreak{}44334. Any nonlinear word with a tail value 2 or 3 has value 3. With all tail values in \{0,\allowbreak{}1,\allowbreak{}4\}, head 2 gives 2 and head 3 gives 3. For heads 0,\allowbreak{}1,\allowbreak{}4, a tail 0 or 4 makes the result 4; if all tails are 1, the head is retained. Literals retain value. Induction on the table proves the word formula. Thus \(xyz=xzy\) and \(xyy=xy\); normal words are a head and sorted tail set, retaining head repetition. Actual tail swaps and duplicate deletions use nonempty prefixes and nonempty substitutions.

\par\smallskip\noindent\textit{Exact original-witness criterion.}\quad 
Literal x is absorbed iff P contains literal x. For nonlinear \(w=hT\) let \(C=\{h\}\) union T and U=\{u in P: content(u) is contained in C\}. Then \(w\leq P\) iff: (C) the union of contents of U is C; (D) the union of tails of U contains T; (R) h belongs to T, or literal h is in P, or an original right word is entirely contained in T, or an original nonlinear h-word has tail in T. All witnesses are original. A head occurrence alone is not a tail witness.

For necessity of literal occurrence, put \(x=2\) and all \(others=3\). The target is 2, and all other literals and nonlinear right words give 3. If C fails at x, put \(x=4\), the other target \(variables=1\) and \(outsiders=3\). The target is 4, contained words avoiding x give 1, and noncontained words give 3. If D fails at x in T, put \(x=0\), other target \(variables=1\), and \(outsiders=3\). The target is 4; a contained right word can have x only as its sole head occurrence, so gives 0, or avoids x and gives 1. Noncontained words give

3. The right sum is at most 0, below 4. If R fails, h is not in T; put \(h=2\), \(T=1\), \(outsiders=3\). The target is 2, literal h is absent, every nonlinear h-word has a tail outside T and gives 3, and every word with head in T has a tail outside T since no tail-contained word exists. All right words are therefore 3. Each separator applies to arbitrary polynomials.

For sufficiency, target 3 is least. A target 4 has no 2 or 3 in its content, and either its head is 4 or some tail is 0 or 4. In the first case C supplies a contained word containing that head variable, giving 4. In the second D supplies a contained word with the selected variable in its tail, again giving 4. For target 0 the head is 0 and all tails are 1; C supplies a contained word containing the head, giving 0 or 4. For target 1 every target variable is 1; any contained word gives 1. For target 2, the head is 2 and all tails are 0,\allowbreak{}1,\allowbreak{}4, so the head does not repeat in T. R supplies literal 2, a tail-contained word valued 0,\allowbreak{}1,\allowbreak{}4, or a nonlinear same-head tail-contained word valued 2. This proves exactness.

\par\smallskip\noindent\textit{Finite local construction from directed tail paths.}\quad 
Normalize original words. For a nonlinear target, work only with U. Use graph vertices C and an edge \(p\to x\) for each original word headed p having x in its tail. All x in T have positive indegree by D. First merge words with the same head: p(A union B)\ensuremath{\leq}pA+pB. Literals are allowed as an empty tail in this notation; the written words and substitutions remain nonempty. C and D hold for this local absorption, and R is supplied by a literal p or a nonlinear p-word already on the right.

For each head p, maintain pA. When q lies in A and qB is available, derive p(A union B)\ensuremath{\leq}pA+qB. Its content is covered, its tails are covered, and the nonlinear pA supplies R. This is valid also when \(q=p\). Iterate only steps that add a tail variable. Finiteness of C gives termination, with the tail at p exactly the set of vertices reachable by a POSITIVE length directed path from p. No zero-length reflexive edge is inserted. Consequently a word qQ is saturated (q in Q) exactly when q is on a positive-length directed cycle.

Every vertex is reachable from h or from a directed cycle. Indeed for a vertex other than h, D supplies an incoming edge; following predecessors either reaches h or repeats a vertex. Reversal gives the desired directed path. If no original word has head h, h cannot be a source of an edge or a literal component; C forces h into a tail, so all vertices have positive indegree and are reachable from cycles.

If an h-headed word exists, begin with its closed word g. Otherwise begin with any saturated cycle word. Merge every saturated cycle word \(b=qQ\) whose content is not yet in the tail of g, taking head(g) with the union of its tail and content(b). This local absorption has C and D because q is in Q; R comes from g, either literal or nonlinear same-head with retained tail. It has at most two right words. At the end g has content C and tail containing T. To see the tail assertion, every x other than head(g) is in the tail automatically. If head(g) lies in T, it has positive indegree; it either returns along a path from itself (already saturated) or is reached from a directed cycle (then occurs in that cycle word's tail). If a cycle word containing head(g) has no new content, its tail already contains head(g); the merge must still be made if that tail is missing. Thus component merging is governed by missing tail content, not only by missing total content. This distinction is essential for singleton cycles.

Finally derive \(w\leq g+b\) where b is an original R witness; if h repeats, b is unnecessary. C and D hold through g, and R through b or repetition. This recovers the target. Every step uses at most two available right words.

\par\smallskip\noindent\textit{Uniform finite family.}\quad 
For any local absorption reserve the heads of its at most three words, at most three variables. Group remaining variables by their nonempty incidence across the at most three tails, at most seven classes. Identify each class to a fresh variable and normalize. Validity is preserved by substitution. The result uses at most ten variables and words of length at most eleven. Inverting each class variable by the nonempty product of its members recovers the original tails, heads and head repetitions using the two word laws. Hence every local step is derivable from the two word laws and the finite family of valid normal absorptions over ten variables, length at most eleven, with at most two right words. The exact C, D, R criterion specifies this family. The graph construction derives every absorption, and distributivity plus summandwise absorption in both directions derives every valid identity. This is an unrestricted finite relative basis; finite replay is supplementary evidence only.

Opposite reversal is valid. Every catalogue transfer must separately pass complete HSP in both directions. Head occurrence cannot replace D: the absorption \(xy\leq x+y\) fails at \(x=1\), \(y=0\), since the target is 4 and the sum is 0.

\subsection{Ten identities for 11337 and their reversals for 12057}\label{case:247}
Use the constant-free ai-semiring signature and put u \ensuremath{\leq} v when \(u+v=v\). The multiplication rows of 11337 are 40424, 41434, 22222, 33333, 44444. The additive order has bottom 3, atoms 1 and 2, their join 0 and top 4. Relative to the semiring axioms the following ten identities form a basis.

{\small\par\smallskip
\CaseLine{\(xyy = xy\qquad xyzy = xyz\)}
\CaseLine{\(h \leq  hh\qquad hp \leq  hph\)}
\CaseLine{\(ayv \leq  av + ayb\qquad ayv \leq  av + ay\qquad ay \leq  a + ayb\)}
\CaseLine{\(ayv \leq  ayb + yv\qquad ayv \leq  ay + yv\qquad ay \leq  ayb + y\)}
\smallskip}

The opposite table is 12057 up to the certified relabeling. Reverse all the words for its basis.

\par\smallskip\noindent\textit{Absorption.}\quad 
A head of value 2, 3 or 4 fixes the entire value. A head of value 1 or 0 ignores all following 1's. If there is no other tail value the result is that head. At the first other tail value, a 3 changes 1 to 3 and 0 to 2, while 0,\allowbreak{}2,\allowbreak{}4 each send both starting values to 4. The last three rows fix the resulting value.

Let the target consist of a head h and a possibly empty tail. A polynomial P of original words absorbs it exactly when all of the following hold.

\textbullet{} A word in P starts with h.

\textbullet{} A word in P is supported in the target.

\textbullet{} At each first occurrence of y in the target tail, let D be the variables earlier in that tail. Some original has a noninitial occurrence of y with all preceding variables in \{h\} union D. This includes \(y=h\) if h reappears.

Giving h the value 4 and everyone else 3 forces the head witness. Giving target letters 1 and outsiders 3 forces a contained word. When y differs from h, assign h and D the value 1, y the value 2 and the others 3. Only an original with a noninitial y preceded by these 1's can reach 4. When \(y=h\), put h at 0, D at 1 and all others at 3. The target gives 4. The only originals giving 4 have a noninitial h preceded only by h and D. All values below 4 have join at most 0.

Conversely a fixed head value is witnessed by a same-head word. A head of 1 with an all-1 tail uses containment. A head of 0 with an all-1 tail uses its head witness (giving 0,\allowbreak{}2 or 4) and its contained witness (giving 0,\allowbreak{}1 or 4); their join is at least 0. A head of 0 followed first by a non-1 tail value of 3 gives 2 and uses the head witness. Head 1 followed first by 3 gives bottom. In every remaining case the target gives 4 and the first non-1 tail value y is 0,\allowbreak{}2 or 4. At the chosen noninitial occurrence in its witness, earlier values are all 1 or the head value in \{0,\allowbreak{}1\}. If an earlier noninitial 0 has already sent the word to 4 it stays there; otherwise the indicated occurrence sends the active value to 4. This proves the equivalence on unrestricted inputs.

The two word laws erase repeated tail letters and give the unique spelling h followed by the list of first occurrences in the tail. For separated repeats use the second word law with the nonempty word through the position preceding the first repeated tail letter as its initial argument. The first two inequalities permit inserting another head after any prescribed tail prefix on the upper side. By the word formula this only increases a value when the head is 0 and the tail prefix has been all 1. Soundness of the last six laws follows from the witness criterion. For example ayv \ensuremath{\leq} ayb+yv uses the first upper word for the head and the first tail letter y, and the second for containment and the tail letter v. All ten identities also admit exhaustive table checks.

\par\smallskip\noindent\textit{Complete derivation.}\quad 
Normalize the target. Build a current word h p where p is the processed tail, with at most one extra early copy of h inserted in p. The invariant is that for every original g v whose head g belongs to \{h\} union the support of p, the word h p v is already below P. A literal original uses v empty. For p empty the invariant is simply the list of original h-headed words.

To append a new tail letter y, suppose one has h p y b \ensuremath{\leq} P for some suffix b, possibly empty. Write \(a=h\) p. For every previously eligible original g v, use the third line of laws to derive a y v from a v and a y b. For each original newly eligible by having head y, use the last line with the original y v and a y b. The cases of empty v or b use the corresponding displayed shorter identity, and if both are empty the available a y suffices. This extends the invariant without requiring a concatenation with an arbitrary original head.

Before each ordinary step, and once more at the end, check whether h has already occurred in p. If it has not and an original contains a noninitial h whose prefix is supported in \{h\} union the support of p, take the earliest noninitial h in that word. Write the original as g s h b. The head g is eligible, and s uses only letters of p: any h in s would have been an earlier noninitial h. The invariant gives h p s h b, which normalizes to h p h b. Thus the preceding insertion step appends h. Use this insertion only once.

Read the scheduled first tail letters in order. Skip any scheduled letter already in p (only h can require this). At the next y, take a witnessing original g s y b from the criterion. Its head is eligible. The letters of s all lie in p: if h is absent from p, an h in s would have triggered the preliminary insertion just described. Therefore h p s y b from the invariant normalizes to h p y b, and y can be appended. This also handles the first tail occurrence of h itself because its required witness permits the preliminary insertion. The reasoning uses the original occurrence conditions at every step; the inserted copy of h changes neither which unprocessed ordinary letters remain nor their original witnesses.

At the end choose a contained original g v. Its head is eligible. All letters in v already occur in p. The only potential exception is h if it never entered p, but a noninitial h in this contained original would qualify at the final check. Hence h p v normalizes to h p. If h was inserted early, form the upper word by inserting h in the original target at this same place. Use h \ensuremath{\leq} hh when the prior tail prefix is empty or hp \ensuremath{\leq} hph otherwise, and multiply by the actual remaining suffix. Normalization erases a later duplicate and gives the word just derived. Thus the original target is below P, including a literal target or an empty remaining suffix. All substituted arguments are nonempty. Distributing terms and applying the absorption derivations in both directions proves unrestricted completeness.

\subsection{An ordered odd-cycle obstruction for 11369}\label{case:248}
Use the constant-free semiring signature. The labelled multiplication for 11369 is 40434 / 41434 / 00234 / 33333 / 44434. The element 3 is absorbing and additively least. In its absence, 4 is absorbing. Represent values 2, 0, 1 respectively by intervals (0,\allowbreak{}0), (0,\allowbreak{}1), (1,\allowbreak{}1). A sequence of these three letters has value 4 if the right endpoint of any letter is larger than the left endpoint of its successor. Otherwise its product is the interval from the first left endpoint to the final right endpoint. All nine products on these intervals give the stated interval or 4; appending a letter proves the formula for every word. In particular, a change from 2 to 1 is allowed without an intervening 0. The additive order is 3 below 1 and 2, their join 0, and 4 greatest.

Fix any odd \(n\geq 5\) and work modulo n. With 3n distinct variables put

{\small\par\smallskip
\CaseLine{\(p_{i} = a_{i} z_{i} a_{i+1} b_{i+2} b_{i},\)}
\CaseLine{\(P = sum_{i} p_{i},\)}
\CaseLine{\(w = (a_{0} z_{0} a_{1} z_{1} \cdots  a_{n-1} z_{n-1}) (b_{0} b_{1} \cdots  b_{n-1}).\)}
\smallskip}

Then 11369 satisfies w \ensuremath{\leq} P for every such n. If a variable is 3 the target is least. Otherwise use the principal ideals. If P is 4 there is nothing to check. All upper words use only the three intervals and are compatible whenever P is at most 0: every letter occurs in P and none can be 4. Writing \(R_{x}\), \(L_{x}\) for right and left endpoints, compatibility gives

{\small\par\smallskip
\CaseLine{\(L_{ai} \leq  R_{ai} \leq  L_{zi} \leq  R_{zi} \leq  L_{a}(i+1).\)}
\smallskip}

Around the cycle these inequalities are equalities. All \(a_{i}\) and \(z_{i}\) are the same degenerate interval. Likewise \(R_{b}\)(i+2) \ensuremath{\leq} \(L_{bi}\) for all i. Since n is odd the step of size two traverses all indices; every \(b_{i}\) is the same degenerate interval. Compatibility in \(p_{i}\) orders the first common interval before the second. The ordered word w is therefore compatible and has value at most 0. If P is 1, or 2, each upper word must have exactly that value. Such a degenerate product forces every letter of that word to be its own value, and the same holds for w. P cannot be 3 in the absence of 3. This verifies all assignments without an upper bound on the cycle length.

Let \(A=\A{11269}\).

Apply Proposition~\ref{prop:lower-fibres} to a generating four-tuple of \(A\) to obtain \(A\in\V(\A{11369})\). Thus any identity of 11369 also holds in A. Apply the first-family argument in Section~\ref{case:242}. Here are the specific details that allow the change of length and of target. The markers \(z_{i}\) occur only in \(p_{i}\). Full contained-word coverage in A makes the support of any absorbed word the union of the originals indexed by some nonempty I, with exactly those markers. Noninitial coverage gives at most one i in I with i-1 outside. Hence a proper I is a single run. Nonfinal coverage gives at most one i in I with i-2 outside. A run of length between 2 and n-2 violates this at its first two indices. This proof works for odd n. A single-marker word must begin at the head and end at the final of its unique original, and its support equals that original's support.

On a source absorption of at most k variables, assume its substituted upper words all have a single marker. Each variable appearing noninitially in the source target has a noninitial occurrence in a source upper word by 11269's criterion. Fix its nonempty preceding context in an expansion there. Every word choice for that variable lies in one fixed original's support, since the fixed first letter determines the marker index. The possible remaining variable occurs once initially; ordinary source coverage puts its selected word in some one-marker expansion. Hence an expanded new word contains at most k markers. Repeated source variables choose images independently. In arbitrary contexts two fresh variables suffice for the nonempty left and right context words after distribution. Every word at any step of a derivation beginning with P is absorbed by P in A. Choose odd n with \(n-1>k+2\). The gap forbids a first step creating more than one marker. The new ordered target contains all n markers. This proves nonfinite basability of 11369.

The inclusion of Var(A) is proper. In A the concatenation \(p_{0}\) ... \(p_{n-1}\) is absorbed by P since both contained positional unions cover every letter. In 11369 assign every \(a_{i}\), \(z_{i}\) value 2 and every \(b_{i}\) value 1. Then every \(p_{i}\) is 0 and the concatenation has value 4, above \(P=0\). This also explains why the ordered replacement in the critical target was necessary. The proper NFB subvariety rules out limit status.

\subsection{Six complete identities for 11372 from contained positive paths}\label{case:249}
Signature (+, multiplication), no constants. Relative to the additively idempotent semiring axioms, a finite basis is:

1. \(xyz=xzy\).

2. \(xyy=xy\).

3. \(xyz\leq xy+xz\).

4. \(xyz\leq xy+yz\).

5. \(xy\leq xxy\).

6. \(x\leq xx\).

Here \(s\leq t\) abbreviates \(s+t=t\). Every substitution is nonempty. Completeness has no bound on variables, word length or polynomial size.

Addition: 00004,\allowbreak{}01014,\allowbreak{}00224,\allowbreak{}01234,\allowbreak{}44444. The order has bottom 3, incomparable 1 and 2, their join 0, and top 4. Multiplication: 40434,\allowbreak{}41434,\allowbreak{}22232,\allowbreak{}33333,\allowbreak{}44434. For a nonlinear word, any occurrence of 3 gives value 3. Otherwise a head 2 remains 2. For heads 0,\allowbreak{}1,\allowbreak{}4, a tail value in \{0,\allowbreak{}2,\allowbreak{}4\} gives 4, while all tails equal to 1 preserve the head. Literals retain value. Induction on the table proves this formula and the first two identities. Normalize each word to its head followed by its sorted tail set, retaining an occurrence of the head in the tail. Tail swaps and deletions use nonempty prefixes. The remaining four laws follow from the criterion proved next, or directly from the tables.

\par\smallskip\noindent\textit{Exact criterion.}\quad 
For target w, let C be its content and U the original right words whose whole contents are contained in C. If w is the literal x, \(w\leq P\) iff U is nonempty: some original pure power of x occurs.

For nonlinear \(w=hT\), form a directed graph on C, with edge \(p\to q\) whenever an original word in U is headed p and has q in a tail position. Then \(w\leq P\) iff every x in T is reachable from h by a POSITIVE-length path. When h repeats in the target, h must lie on a directed cycle. Reflexive reachability by a path of length zero is insufficient. Outside variables cannot be inserted in the graph, since those right words are not in U.

For literal necessity, if no pure power of x occurs, put \(x=4\) and all \(others=3\). Target 4 is separated from the right sum 3. Sufficiency follows from \(x\leq x^{n}\) for every positive n: powers normalize to x or xx, and the table validates \(x\leq xx\).

For nonlinear necessity, let R be the positive reachable set from h. If some target tail x different from h is missing from R, put every variable in \{h\} union R to 1, the other variables of C to 2, and outsiders to 3. The target is 4. Every contained word headed in \{h\} union R has all its tails there by reachability closure and gives 1. Other contained heads have value 2 and give 2. All noncontained words give 3. Their sum is at most 0 and is strictly below 4. The remaining failure is h in T but h not in R. Put \(h=0\), variables of R to 1, other variables of C to 2, and outsiders to 3. Any contained h-word has all tails in R and gives 0; words headed in R give 1, remaining heads give 2, and noncontained words give 3. The target repeats value 0 and gives 4, separating it from a right sum at most 0.

For sufficiency, value 3 is least. Reachability of the nonempty target tail implies an original nonlinear h-word in U. If target is 2, its head is 2 and C has no 3; that h-word gives 2. If target is 0 or 1, all its tails are 1 and its head has the target value. The original h-word gives the target value or 4. If target is 4, no variable of C has value 3 and the head is not 2. Suppose no contained right word has value 4. The original h-word forces head value 0 or 1, and every one of its tails to

1. Along each graph edge whose head is 0 or 1, the same absence of 4 forces its tail value to 1. Induction along positive paths makes every target tail equal to 1. The target would then be 0 or 1, a contradiction. Thus some original right word is 4. This proves exactness for every assignment, without any finite search assumption.

\par\smallskip\noindent\textit{Complete derivation with the six laws.}\quad 
Normalize original words. For a literal target choose its original pure power, normalize it to x or xx, and apply law 6 if needed.

For nonlinear \(w=hT\), choose an original nonlinear h-word g in U. Keep g headed h, with all its tail variables in C. Until its tail contains T, choose an original \(u=qB\) in U with q equal to h or already in the tail of g, and with some tail variable absent from g. Such u exists by following a positive path from h to a missing target variable and taking its first edge that leaves the currently reached set. Its tail B is nonempty. If \(q=h\), write \(g=hA\) with nonempty A and use law 3 with \(x=h\), \(y=A\), \(z=B\). It gives \(hAB\leq hA+hB\). If q is in the tail of g and differs from h, reorder g to vq, where v is its head and the other tail variables, so v is always nonempty. Law 4 with \(x=v\), \(y=q\), \(z=B\) gives \(vqB\leq vq+qB\). In either case normalize the new word; its tail gains variables and stays within C. Each step is an actual substitution of one fixed law plus the two word laws. The process terminates after finitely many strict tail increases.

The resulting g has every desired tail variable and can only have one additional tail variable, h itself. If there is no such extra occurrence, g is the target. Otherwise the target is hT with h outside T and the current word is hhT up to tail permutation. Apply law 5 with \(x=h\) and y the nonempty word T to derive the target. Every used right word is original or previously derived, so the target is absorbed by the original P. This derives every valid absorption. Distributivity and summandwise absorption in both directions derive all identities, giving the claimed unrestricted six-law basis.

The direction is essential: \(xy\leq yx\) fails at \(x=1\), \(y=2\). So does \(xx\leq x\) at \(x=0\). A literal x is nevertheless absorbed by xx through law 6; literal occurrence is therefore not a necessary witness here. These distinctions are retained in independent replay and counterexample checks.

\subsection{Five nonfinite identity bases from a path with a distinguished endpoint}\label{case:250}
The five catalogue entries are 11412, 11496, 11523, 11524 and 11564. Use the path polynomials \(P_{n}\) and word \(W_{n}\) of the preceding 12438 proof. Their definitions are recalled here. There are n main letters \(x_{1}\),..., \(x_{n}\), with n at least four. Every pair \(i<j\) supplies one word, using a fresh private letter at each marker position. For consecutive i, j use a \(x_{i}\) \(x_{j}\) b, omitting a when \(i=1\) and b when \(j=n\). For nonconsecutive i, j use a \(x_{i}\) c \(x_{j}\) b. The lower word \(W_{n}\) lists all n main letters in increasing order.

All five tables satisfy \(W_{n}\) \ensuremath{\leq} \(P_{n}\) for every n. For each bound, project each kind of upper word onto its two main entries by enumerating all private entries. A sequence of main values admits a private completion in the principal ideal exactly when every relevant pair belongs to its projected relation. Test reachability outside the ideal with state (product, last value, set of preceding values, length capped at three), choosing the last step separately. This tests every sequence of length at least four. The completions below are finite calculations on the same tables.

\par\smallskip\noindent\textit{A necessary condition for every intermediate monomial.}\quad 
In all five tables the subset \{1,\allowbreak{}0,\allowbreak{}4\} has additive chain 1 \ensuremath{<} 0 \ensuremath{<} 4. Its multiplication has identity 1, absorbing element 4, and 0 times 0 equal to 4. If a polynomial of linear words absorbs v, this subalgebra forces v to be linear, to use only letters of that polynomial, and to have every pair of its distinct letters together in some upper word. For support assign the missing letter 4 and all others 1. For a repeated letter or a missing pair assign the chosen letter or pair 0 and all others 1. In particular, any word absorbed by \(P_{n}\) that uses a private letter has length at most five. Its entire support lies in the single original containing that private letter.

For the five tables choose the following assignments to the main variables. Each row permits the indicated private completions with \(P_{n}\) below the displayed bound. The subscript i denotes an arbitrary consecutive pair in the first three rows, or an arbitrary single position in the last two.

\noindent Type\enspace\textbar\enspace Bound\enspace\textbar\enspace Main values\enspace\textbar\enspace Consequence for an absorbed word using only main letters\par
\noindent 11412\enspace\textbar\enspace 0\enspace\textbar\enspace All 1 except (\(x_{i}\), \(x_{i+1})=(3\),2)\enspace\textbar\enspace If \(x_{i}\) occurs then \(x_{i+1}\) occurs.\par
\noindent 11496\enspace\textbar\enspace 0\enspace\textbar\enspace All 1 except (\(x_{i}\), \(x_{i+1})=(2\),3)\enspace\textbar\enspace If \(x_{i+1}\) occurs then \(x_{i}\) occurs.\par
\noindent 11523\enspace\textbar\enspace 0\enspace\textbar\enspace All 1 except (\(x_{i}\), \(x_{i+1})=(2\),3)\enspace\textbar\enspace If \(x_{i+1}\) occurs then \(x_{i}\) occurs.\par
\noindent 11524\enspace\textbar\enspace 2\enspace\textbar\enspace All preceding values 1, \(x_{i}=2\), all subsequent values 3\enspace\textbar\enspace The word contains \(x_{i}\).\par
\noindent 11564\enspace\textbar\enspace 2\enspace\textbar\enspace All preceding values 3, \(x_{i}=2\), all subsequent values 1\enspace\textbar\enspace The word contains \(x_{i}\).\par
Here is a finite proof of private completability for every indicated assignment. For \(n=4\),...,7 check, for each distinguished position or pair and each original word, that some tuple of at most three private entries places its product in the ideal. For each tuple multiply its entries in the displayed order and test membership in the ideal. If any original at larger n failed, keep the two endpoints, the two main positions of the failing word, the at most two distinguished positions, and one intervening position if the failing word is nonconsecutive. There are at most seven positions. If necessary fill up to four. On restriction their order, the type of this original (including both boundary roles) and the colors are preserved. This gives a failing test at length at most seven. Since all private letters are distinct, separately chosen completions work simultaneously.

In each of the first three rows an absorbed main-only word is already linear. If it contained the variable with value 3 while omitting its partner of value 2, its evaluated sequence would be a single 3 among 1's. Its product is either 3 or 4, both outside the ideal of 0. This also follows by closing \{3\} under multiplication by 1 on either side. In the last two rows the multiplicative subsemigroup generated by \{1,\allowbreak{}3\} misses the ideal of 2. Thus omitting the distinguished position always gives a product outside that ideal. This proves the consequences in the table without any word-order assumption.

For 11412 distinguish the variable \(x_{1}\). For 11496 and 11523 distinguish \(x_{n}\). For the last two either endpoint will do. Every main-only absorbed word containing this distinguished variable has length n by propagation along the path or by the individual position tests. Consequently every absorbed word containing the distinguished variable is either of length at most five or has length n. Shorter main-only words *without* the distinguished letter can occur in the first two tables.

\par\smallskip\noindent\textit{Excluding a finite basis even in the presence of other long words.}\quad 
Suppose that a derivation uses at most k variables per additional identity, and choose \(n>5(k+1)\). Work with distributed sets of words. Every summand in a derivation beginning at \(P_{n}\) is absorbed by \(P_{n}\) and therefore linear. Select the first step that produces a word of length greater than five containing the distinguished variable. There is such a step if \(P_{n}\)+\(W_{n}\) is reached. The produced word has length n.

Fix the producing monomial, multiplication contexts L, R and the nonempty image word selected at each occurrence in the substitution. To prove that the producing monomial is linear, suppose it repeats a source variable. Its occurrences may choose different words in the expansion under consideration. Choose the *same* nonempty image word at both occurrences in another expansion. This produces a post-step summand with a repeated letter. That entire post-step polynomial still equals \(P_{n}\), which excludes every such summand. The same argument applied before the step shows that every upper monomial of the applied source identity is linear. The context and distributivity argument covers arbitrary one-hole semiring contexts. Empty exterior contexts are allowed.

This alternative-expansion argument also supplies the missing justification for source linearity in the first-step proof for 12438 and 13708 at gate 264. Linearity of the one chosen output alone would not imply it. In those two tables every summand of every intermediate polynomial is linear by their proven criterion. Thus the alternative expansion gives the required contradiction there as well, and the length bound in that proof follows.

If L or R contains the distinguished letter, all old expansions with this context pair have length at most five. Every producing variable occurs on the old side by the subalgebra content test. Place its selected image in an old expansion to bound its length by five, and use any old expansion to bound the total context length by five.

Otherwise a producing source variable z has a selected image containing the distinguished letter. The source identity has linear old monomials. Its subalgebra pair condition places z together with each other producing variable y in some old monomial. Select the fixed images for both variables in that monomial and arbitrary images elsewhere. The resulting old word contains the distinguished letter, so its length is at most five. This bounds the selected image of each y, the image of z, and the total context length. If z is the sole producing variable, use its content witness instead. This is the step which permits arbitrarily long old words *without* the distinguished variable.

There are at most k producing letters. The chosen new word has length at most 5k+5, contradicting its length n. The argument applies in both orientations and allows different image choices at repeated source occurrences throughout. Each of the five tables is nonfinitely based. Their previously recorded proper nonfinitely based subvarieties already exclude limit status.

\subsection{Ten complete identities for 15092}\label{case:251}
In the constant-free additively idempotent semiring signature, with \(u\leq P\) abbreviating \(u+P=P\), a complete finite basis is:

{\small\par\smallskip
\CaseLine{\(xy=yx; xxxx=xxx; xxxy=xxy;\)}
\CaseLine{\(xy\leq xyz; xx\leq xxy; x\leq xyz; y\leq x+xy;\)}
\CaseLine{\(xz\leq xx+xy+yz; xz\leq x+xx+z+zt;\)}
\CaseLine{\(xyz\leq x+y+z+xy+xz+yz.\)}
\smallskip}

\par\smallskip\noindent\textit{Values and exact original-witness criterion.}\quad 
Addition:00044/\allowbreak{}01144/\allowbreak{}01244/\allowbreak{}44434/\allowbreak{}44444. Multiplication:44444/\allowbreak{}44144/\allowbreak{}41244/\allowbreak{}44404/\allowbreak{}44444. The additive order \(is 2<1<0<4\), with 3 incomparable to 0,\allowbreak{}1,\allowbreak{}2 and below 4. A nonlinear word with every occurrence in\{1,\allowbreak{}2\} has value 2 if all are 2, value 1 if exactly one occurrence is 1, and 4 if at least two occurrences are 1. The only other non4 case is a quadratic with both occurrences 3, which has value 0. Every other nonlinear evaluation gives 4. These formulas follow inductively from 2 acting as identity on\{1,\allowbreak{}2\},1 squared being 4,\allowbreak{}3 squared being 0, and 0 annihilating to 4. Long products are not constant, and their occurrence multiplicities cannot be discarded.

By commutativity and the two multiplicative reduction laws, a pure power normalizes at exponent 3, and words with at least two distinct variables keep each positive multiplicity capped at 2. Short quadratics are still distinct from longer words. All reductions preserve word functions.

For the original polynomial P define L, V as its literal and nonlinear variable sets. Put an edge xy in E whenever two distinct positions of some original nonlinear word have letters x, y; loops therefore record actual repetition. Let H0 be the variables with an original repeated occurrence. Build components of V by co-occurrence in original words. Let T be the union of contents of original words of \(length\geq 3\), let A be the union of components meeting L or T, and put \(H=H0\) intersect A.

The exact criterion for absorption is:

\textbullet{} A literal x is absorbed iff x belongs to L union A.

\textbullet{} A quadratic xy is absorbed iff its variables are in V and either xy belongs to E, or both endpoints lie in one component with an endpoint in H0, or both are active with an endpoint in H.

\textbullet{} A target w of \(length\geq 3\) is absorbed iff its content is in A, every repeated target variable belongs to H, and its distinct variables outside H form a clique in the original co-occurrence graph E.

All data remain tied to original words. In particular, a derived square must not turn a previously simple variable into an original repetition witness. These conditions differ from the nilpotent cubic criterion for 14967: there is no requirement for an original global cubic guard.

To prove sufficiency, consider a non4 value of P with nonlinear terms. Every component is either entirely 3 or entirely in\{1,\allowbreak{}2\}; in the latter mode, each original word contains at most one occurrence with value 1. A component meeting an original long word cannot be all 3. A component meeting a literal cannot be all 3 either, because that literal 3 and a nonlinear value 0,\allowbreak{}1 or 2 would sum to 4. Thus every active variable is 1 or 2, and every variable in H is 2. An active literal 1 has an original nonlinear incidence of value 1; literal 2 is below every nonlinear value. The quadratic rule follows by keeping the anchored endpoint 2 in the active mode, or both endpoints 3 in the common inactive component mode. If the target pair has value 1, the other endpoint's original incidence has value 1. For a long target, all repeated variables are 2, and no two unanchored distinct variables can both be 1 because their original edge forbids it. Hence the target is 1 or 2 and is absorbed by an original incidence of a target 1, or by any nonlinear right word in the all 2 case.

For necessity, a missing literal outside L union V can be set 4 while all others are 2. A missing literal in an inactive component is separated by setting that whole component 3 and all others 2; P has value 0 but the literal is 3. A nonlinear target containing a variable outside V is separated by setting that variable 0 and all others 2: target 4, \(right\leq 0\). For a failed quadratic whose endpoints lie in different components with one inactive, set that whole component 3 and the others 2. Otherwise a failed pair has no anchored endpoint and no original edge. Set its endpoints 1 and all others 2. Their original multiplicities are single, and no original word contains both; \(right\leq 1\) while the target is 4.

For a failed long target with an inactive variable, set its component 3 and all others 2, obtaining target 4 and right 0. If a repeated target variable has no original repetition, set it 1 and all others 2; each original word has at most one 1 but the target has at least two. Finally, if two distinct unanchored target variables lack an original edge, set both 1 and all others 2. Again every right word has at most one 1, while the target has at least two. These separating assignments cover every failure, including literal-only right sides, without any variable bound.

\par\smallskip\noindent\textit{Finite complete derivations.}\quad 
Extract original co-occurrence pairs with \(xy\leq xyz\), where z is the nonempty remaining word. Extract original repeated-variable squares with \(xx\leq xxy\), again using the nonempty remainder. Original quadratics and squares require no extraction. An original long word supplies a literal at any chosen variable using \(x\leq xyz\) and a nonempty two-part remainder. Together with original literals, this starts recovery throughout every active component: traverse original edges by \(y\leq x+xy\).

For a same-component pair anchored at x, keep xx and traverse a simple positive original path by \(xz\leq xx+xy+yz\). A square needs no path. For a pair across active components, keep the anchored square xx, the recovered literals x, z, and an original incidence zt. Then \(xz\leq x+xx+z+zt\) gives the pair. These are all allowed quadratic cases.

Consider an allowed target w of \(length\geq 3\). Its letters are active, so all its literal letters have been recovered. Every pair of occurrence positions in w is allowed: a repeated letter has its actual original repetition witness; a pair with an anchored endpoint uses the preceding pair construction; a pair of unanchored letters has its original edge by the clique condition. Thus all singleton and quadratic subwords of w, with their true multiplicities, are available.

Induct on subword length. For a word of \(length\geq 3\), split its positions into nonempty words X, Y, Z. The six proper subwords X, Y, Z, XY, XZ, YZ are already available by induction. Substitute them into \(xyz\leq x+y+z+xy+xz+yz\) to recover XYZ. No variable is collapsed or silently deleted: repeated occurrences stay repeated in these subwords. This constructs w at every length. The two reduction laws merely shorten representatives and are not needed to impose a length bound.

Each listed law is valid by the criterion or word-value formula. Every valid absorption is now derived using ten fixed laws and semiring axioms, with nonempty substitutions. Distributive expansion and absorption in both directions prove unrestricted equational completeness.

Boundary failures: \(xxx\leq xx\) at \(x=3\); \(y\leq xy\) at \(x=y=3\); \(xyz\leq xy+xz+yz\) at \(x=y=z=3\); \(xx\leq xyz\) at \(x=1\), \(y=z=2\); \(xz\leq xx+zt\) at \(x=3\), \(z=t=2\); \(xy\leq xz+yz\) at \(x=y=1\), \(z=2\). In contrast, \(xxx\leq x+xx\) is valid here: recovered literal protection replaces the original global cubic guard that was necessary for 14967.

\subsection{An absorption criterion for \texorpdfstring{\(\A{11619}\)}{A11619}}\label{case:252}

In the order \(0,1,2,3,4\), the addition rows of \(\A{11619}\) are
\(\mathtt{00444,01444,44244,44434,44444}\), and its multiplication rows are
\(\mathtt{40444,41244,44424,44434,44444}\). The additive order has \(1<0<4\) and two further minimal elements \(2,3\) covered by \(4\).
By induction on length, the value of a nonempty word is \(1\) or \(3\) when every letter has that value, \(0\) on strings of the form \(0\,1^*\), and \(2\) on strings of the form \(1^*2\,3^*\). All other strings have value \(4\). Here the stars describe sequences of assigned values, including empty prefixes or suffixes.

For a word \(w=a_1\cdots a_m\), put \(\tau(w)=\{a_i:2\le i\le m\}\). For a nonempty polynomial \(P\), set \(C=\bigcup_{u\in P}\ct(u)\) and \(D=\bigcup_{u\in P}\tau(u)\). Introduce an edge \(L_x\to R_x\) for each \(x\in C\). Identify \(R_x\) with \(L_y\) for every adjacent pair \(xy\) in an original summand, all left endpoints of initial letters with \(H\), and all right endpoints of final letters with \(T\). Call this directed graph \(G(P)\). Write \(v\sim v'\) for membership in the same strongly connected component.

\begin{lemma}\label{lem:11619-criterion}
The absorption \(w\law P\) holds in \(\A{11619}\) exactly when
\[
\begin{array}{ll}
\mathrm{(C)} & \ct(w)\subseteq C,\\
\mathrm{(D)} & \tau(w)\subseteq D,\\
\mathrm{(E)} & L_{a_1}\sim H,\quad R_{a_m}\sim T,\quad
                 R_{a_i}\sim L_{a_{i+1}}\quad(1\le i<m).
\end{array}
\]
\end{lemma}
\begin{proof}
If (C) fails, give the missing letter value \(4\) and all other letters value \(1\). If (D) fails, give the offending letter value \(0\) and all others value \(1\). Then \(w=4\), while every upper word has value \(0\) or \(1\).

Suppose that (C) and (D) hold and an equality in (E) fails. Order the strong components by directed reachability. Every vertex occurs on an original walk from \(H\) to \(T\), so \([H]\) is least and \([T]\) greatest. For the two distinct components in the failed equality choose an ordering \(a,b\) in which \(b\) cannot reach \(a\). The indicator function of the upper set generated by \(b\) separates \(a,b\), has values \(0,1\) at \(H,T\), respectively, and is monotone along every directed edge. Assign \(x\) the value \(1,2,3\) according as its pair of endpoint labels is \((0,0),(0,1),(1,1)\). Every original word is a composable walk from label \(0\) to label \(1\), so has value \(2\). A word has that value only on a sequence \(1^*2\,3^*\), which meets every start, end and adjacency requirement on the labels. Hence the failed equality makes \(w\ne2\). The only element below \(2\) is \(2\).

Conversely, consider a valuation. If \(P=1\) or \(P=3\), every letter in \(C\) has that value, so (C) suffices. If \(P=0\), then every upper word is \(0\) or \(1\), all letters in \(C\) have values in \(\{0,1\}\), and every letter in \(D\) has value \(1\). Conditions (C) and (D) give \(w\le0\). If \(P=2\), each upper word has value \(2\). Assign to its letters the Boolean endpoint pairs above. The labels agree at each identified endpoint and are monotone along edges, hence constant on strong components. By (E) the target is a composable walk from \(0\) to \(1\) and has value \(2\). Finally, \(4\) is the additive top.
\end{proof}

\subsection{A basis of 199 identities for \texorpdfstring{\(\A{11619}\)}{A11619}}\label{case:253}

Relative to \(\ax\), take the three identities
\[
 \mathrm{(B)}\quad xyx\approx x^2y^2,\qquad
 \mathrm{(I1)}\quad x\law xy^2,\qquad
 \mathrm{(I2)}\quad xy\law xz^2y ,
\]
and the indicated instances of the following seven schemes.
\[
\begin{array}{r@{\quad}l}
\mathrm{(M)} & AxD\law AxB+CxD,\\
\mathrm{(R)} & AxtF\law AxyB+CzyD+EztF,\\
\mathrm{(S)} & AxzD\law AxyB+yC+zD+EzF,\\
\mathrm{(T)} & DzyB\law AxyB+Cx+Dz,\\
\mathrm{(E)} & Ax\law AxyB+CzyD+Ez,\\
\mathrm{(F)} & yB\law AxyB+CxzD+zE,\\
\mathrm{(L)} & yzB\law xyzB+zC .
\end{array}
\]
In (S) the block \(E\) is compulsory, so \(EzF\) witnesses a noninitial occurrence of \(z\). Every other upper-case block may be omitted or retained, consistently at all of its occurrences in the same identity. Retained blocks and lower-case letters are pairwise distinct variables before substitution. All substitutions in a deduction have nonempty images. The schemes contribute \(16,64,32,16,32,32,4\) instances respectively, giving \(196+3=199\) identities.

\begin{proof}[Validity]
In (B) both words give \(1\) if \(x=y=1\), and \(3\) if \(x=y=3\). Neither can give \(0\), since the initial \(x\) recurs. For \(xyx\) to give \(2\), its unique bridge would have to be \(y\) with \(x=1\) before it and \(x=3\) after it. In \(x^2y^2\) each letter recurs. Both sides therefore give \(4\) in all other cases. In (I1) the upper word gives \(R_x=L_y=R_y=T\); in (I2) it gives \(R_x=L_z=R_z=L_y\). Hence the lower words satisfy Lemma~\ref{lem:11619-criterion}. In (M), (R), (S), (T), (E) and (F), the upper words identify all endpoints needed by the lower word. In (S) the new noninitial \(z\) is witnessed by \(EzF\). In (F) we have \(L_y=R_x=L_z=H\). In (L) the upper words give \(L_x=H=L_z\). Thus the path labelled \(xy\) returns to \(H\) through \(L_y\), giving \(L_y\sim H\). The noninitial letters of \(yzB\) already occur noninitially in \(xyzB\).
\end{proof}

\begin{proof}[Completeness]
Call a word available if its absorption by \(P\) is derivable. From available words record adjacent pairs, initial letters and final letters, starting with the original summands. Write \(J\) for the set of recorded pairs and \(I,F_0\) for the two sets of marked letters. Keep the original \(D\) fixed. Schemes (R)--(F) implement the five rules below. Each premise is witnessed by an available word, split at the indicated occurrence.
\[
\begin{array}{r@{\quad}lcl}
\mathrm{(R)}&xy,zy,zt\in J &\Longrightarrow&xt\in J,\\
\mathrm{(S)}&xy\in J,\ y\in I,\ z\in I\cap D &\Longrightarrow&xz\in J,\\
\mathrm{(T)}&x,z\in F_0,\ xy\in J &\Longrightarrow&zy\in J,\\
\mathrm{(E)}&xy,zy\in J,\ z\in F_0 &\Longrightarrow&x\in F_0,\\
\mathrm{(F)}&xy,xz\in J,\ z\in I &\Longrightarrow&y\in I .
\end{array}
\]
For (S), choose the additional original witness \(EzF\) with nonempty \(E\). These operations preserve \(C,D\). The sets being enlarged are finite, so the process reaches a fixed point.

Use the bipartite graph with ordinary rows \(R_x\), ordinary columns \(L_y\), an extra row \(H\), and an extra column \(T\). Its edges record original pairs, heads and finals. Membership in one undirected component describes equality of endpoints in \(G(P)\) before taking strong components. Each of the five rules adds an incidence within an existing component. Conversely, at the fixed point two ordinary rows sharing an ordinary neighbour have equal ordinary neighbourhoods by (R). Hence each component of the ordinary incidence graph containing an edge is a complete bipartite graph. If it has a head mark, (F) marks every column. If it has a final mark, (E) marks every row. Two such components with head marks must coincide by (S), since every ordinary neighbour belongs to \(D\). Two with final marks must coincide by (T). Every column in \(D\) has an original ordinary neighbour. A column outside \(D\) is a head and remains a leaf at \(H\). A row without ordinary neighbours is final, since its letter occurs only at the ends of original words. If another final row has a neighbour then (T) gives the first row that neighbour too. Restoring \(H,T\) now gives
\[
 xy\in J\ \Longleftrightarrow\ R_x=L_y\ \text{and}\ y\in D,\qquad
 x\in I\ \Longleftrightarrow\ L_x=H,\qquad
 x\in F_0\ \Longleftrightarrow\ R_x=T .
\]
This also covers a component containing both extra vertices: the head and final marks then meet in the same ordinary rectangle.

Every walk from \(H\) to \(T\) with noninitial letters in \(D\) is now available. Start with a word witnessing its head. At each successive \(x\), scheme (M) retains the prefix through that \(x\) and attaches the suffix of a witness for the next pair \(xy\), beginning immediately after the chosen \(x\) of the witness. A final application of (M) with a word ending in the prescribed last letter removes the excess suffix. The same final splice works for a literal target.

Let \(K\) be the strong component of \(H\). For each \(v\in K\setminus\{H\}\), take simple paths \(H\to v\) and \(v\to H\) with nonempty labels \(X,Y\). Their concatenation has no intermediate visit to \(H\). Every letter outside \(D\) occurs only initially in the originals, so its edge starts at \(H\). Thus every letter of \(XY\) after the first belongs to \(D\). The last edge of \(Y\) has right endpoint \(H\). In the bipartite identification graph a shortest path from this ordinary row to the extra row \(H\) meets an ordinary column just before \(H\). That column has an ordinary neighbour and represents a letter \(d\in D\) with \(L_d=H\). Take an available word \(dB\) witnessing this head. Then \(XYdB\) is an actual \(H\)-to-\(T\) walk with legal noninitial letters. It is available by the preceding paragraph. Apply (L) to \(XYdB\) and \(dB\), taking \(x=X,\ y=Y,\ z=d\) and using the same trailing word in each premise, or omitting both blocks if it is empty. This gives \(YdB\law P\). Adjoining that summand identifies its first vertex \(v\) with \(H\). Its adjacencies and its final endpoint already agree in \(G(P)\). Repeating for all such \(v\) therefore collapses exactly \(K\) and preserves \(C,D\).

Let \(w\) satisfy Lemma~\ref{lem:11619-criterion}. Its initial vertex is now \(H\). At a gap with distinct endpoints in one strong component, take nonempty paths with labels \(u\) from the first endpoint to the second and \(v\) in the reverse direction. Insert \(uvu\). Use the same construction after the last letter if its endpoint differs from \(T\). These paths lie in a component other than the singleton component of \(H\), so all inserted letters are in \(D\). The resulting word \(W\) is an actual walk. The splicing argument applied to the enlarged polynomial gives \(W\law P\), since each added summand is already available.

Finally, for nonempty words \(a,b\), (I2) and an application of (I1) to the prefix \(au^2\) give
\[
 ab\law au^2b\law au^2v^2b\approx auvub .
\]
At a final gap apply (I1) twice to get \(a\law au^2\law au^2v^2\approx auvu\). The final equality in each display is (B). Performing the insertions gives \(w\law W\law P\). Distributing arbitrary terms and applying the deduction to their summands in both directions proves completeness.
\end{proof}

The two inclusion tests of Proposition~\ref{prop:lower-fibres} give
\(\V(\A{11737})=\V(\A{11619})\) and
\(\V(\A{14143})=\V(\A{14174})=\V(\A{11619}\op)\).
Lemma~\ref{lem:transfer} gives the reversed basis for the latter pair.

\subsection{Nonfinite basability from reversing an edge of a directed cycle}\label{case:254}
Work with nonempty words and nonempty sums in the two-operation additively idempotent semiring signature. The order is \(x\leq y\) iff \(x+y=y\). Suppose a finite semiring A satisfies the following three conditions.

\textbullet{} There are e, z with \(ee=e\), \(ez=ze=zz=z\), and z \(\nleq\) e.

\textbullet{} For a bound b there are p, t \ensuremath{\leq} b such that pp, pt, tp \ensuremath{\leq} b. Every word over p, t with at least two occurrences of t has value outside the principal ideal below b.

\textbullet{} For a bound c there are l, h \ensuremath{\leq} c with ll, hh, lh \ensuremath{\leq} c and hl \(\nleq\) c.

For \(n\ge4\), let \(v_0,\ldots,v_{n-1}\) be distinct variables and put \(P_n=\sum_{i=0}^{n-1}v_i v_{i+1}\), with indices modulo \(n\). Suppose \(A\) satisfies \(v_1v_0\law P_n\) for every \(n\ge4\). Then \(A\) is nonfinitely based.

A support invariant. For a fixed n let F contain precisely the individual cycle vertices as literal words and the forward cycle edges as quadratic words. Suppose the polynomial q uses at most K variables and a nonempty polynomial substitution phi expands q entirely into F. Take n\ensuremath{>}max(K,3). If a word absorption \(u\leq q\) holds in A, every expanded word of phi(u) also belongs to F.

A variable absent from q cannot occur in u. Assign all variables of q the value e and all others z, using the first condition. Any variable of q occurring in a product occurs in a quadratic source word (a longer source word could not expand in F). Each of its images consists entirely of literals, and it is in fact a single vertex. Indeed, the two nonempty image sets of a quadratic source word form a directed complete bipartite subgraph of the forward cycle. Each vertex has one outgoing and one incoming edge. Fixing a vertex in the opposite set forces the first set to be a singleton, and similarly for the second. Repeated source variables use independent choices on expansion and cannot give a quadratic source loop in q. Call these product variables V. The other variables of q are literal-only, denoted L. Their arbitrary polynomial images may contain the entire cycle. If u is literal its image is already in F.

For a nonlinear u containing any L, assign all L the free formal value t. If another position of u belongs to V, color that position's vertex t, color every other vertex p, and assign each V the color of its singleton image. With no V positions simply color all p. All upper summands have value \ensuremath{\leq}b, since no forward cycle edge has both endpoints colored t. But u contains at least two t positions and is outside the ideal. Thus a nonlinear u uses only V variables. More generally, color any independent set of vertices t and all others p, and take the L values to be p. Any repetition of a vertex in the expanded u is ruled out by its singleton independent set. Any two distinct nonadjacent vertices are ruled out by coloring both t. A simple cycle of length at least four is triangle-free. Hence such a u has exactly two positions, with distinct adjacent images.

Let H be the forward edges contributed by the quadratic source summands of q. They have at most K distinct vertices, so H omits a vertex of the cycle. In particular there is no path in H from \(v_{i+1}\) back to \(v_{i}\). If the prospective quadratic is that reverse edge, give all vertices reachable from \(v_{i+1}\), including itself, the value h, and all remaining vertices l. No H edge goes from h to l. Give each product variable the value of its vertex, and all literal-only variables the free value l. Every summand of q is below c, but u has value hl outside the ideal. Thus every allowed quadratic really is a forward edge. This proves the invariant for arbitrary polynomial substitutions, including large literal images.

The obstruction. If there were a finite basis with at most k variables per law, take n\ensuremath{>}max(k+2,\allowbreak{}3). In a putative derivation from \(P_{n}\) to \(P_{n}\)+v1v0, expand all terms by distributivity. The first step introducing a word outside F uses one side of a valid identity against the old side in a multiplicative and additive context. Distribute that context over the substituted sum. For each monomial context incorporate its nonempty left and right words by at most two fresh source variables. Every replacement summand is absorbed by its old component in A. The support invariant, with \(K=k+2\), forbids this first new word. This contradiction proves nonfinite basability.

All lengths and transfer. In any finite target T put an arc a to d in \(R_{b}\) exactly when \(ad\leq b\). A counterassignment maps the directed cycle to a closed walk in \(R_{b}\) with a distinguished arc whose converse is missing. Conversely such an arc on any closed walk gives a failure of a critical identity of length at least four by repeating the walk. Therefore all critical identities hold precisely when every arc on a directed closed walk of every \(R_{b}\) has its converse in \(R_{b}\). The return path check and the converse arc check concern \(R_{b}\) itself. This is a finite decision for all lengths. For the support obstruction above it is sufficient that T satisfy the critical identities and that A with the support invariant belong to Var(T), since every identity of T is then valid in A. No need to assume that A separately satisfies the critical identities in that transfer argument.

We apply the obstruction directly to the four algebras of case~\ref{grp:447}. In all four choose
\[
 (e,z,b,p,t,c)=(1,4,3,1,3,0).
\]
Take \((l,h)=(2,1)\) for \(\A{11634},\A{11640}\), and \((l,h)=(1,2)\) for \(\A{12530},\A{13156}\). The tables give \(e^2=e\), \(ez=ze=z^2=z\nleq e\) and \(l,h,l^2,h^2,lh\le c\) with \(hl\nleq c\). On \(\{1,3\}\) the value of any nonempty word is \(1,3,4\) according as it contains zero, one or at least two occurrences of \(3\). This follows by appending \(1\) or \(3\), since \(1\cdot1=1\), \(1\cdot3=3\cdot1=3\), \(3\cdot3=4\), and \(4\cdot1=4\cdot3=4\). These values verify the condition at \(b=3\).

For the all-length condition, consider first the two displayed orientations \(\A{11634},\A{11640}\). In \(\A{11634}\) the arcs of \(R_0\) are \((0,1),(1,1),(2,0),(2,1),(2,2)\). In \(\A{11640}\) there are also \((2,3),(2,4)\). In both cases the only arcs of \(R_0\) on a closed walk are \((1,1),(2,2)\). The relation \(R_1\) consists of \((1,1)\). For \(R_2\) the first algebra has only \((2,2)\), and the second has all \((2,j)\), of which only \((2,2)\) is on a closed walk. In both algebras \(R_3=\{(1,1),(1,3),(3,1)\}\) and \(R_4\) is full. Passing to the transpose gives all five relations for \(\A{12530}\) and \(\A{13156}\), respectively. Thus every arc on a closed walk has its converse, and the critical identities hold at every length in all four algebras.

\subsection{Three nonfinite bases from shortcuts in a directed cycle}\label{case:255}
Use the constant-free signature of additively idempotent semirings. Write a word inequality w \ensuremath{\leq} P for the identity \(w+P=P\). Catalogue types 11636, 11672 and 11673 are nonfinitely based. Their common variety is certified by complete inclusions in both directions.

\par\smallskip\noindent\textit{The two principal bounds.}\quad 
In 11636 denote labels 1,\allowbreak{}2,\allowbreak{}3,\allowbreak{}0,\allowbreak{}4 by e, f, a, c, z respectively. The order has e, f below c, the element a incomparable with e, f, c, and z greatest. The product rows in label order are 40444, 41434, 00244, 44344, 44444. Thus a word has value a exactly for the sequence e* a f*. A word has value at most c exactly when it consists of f's followed by e's, or of f's, one c, and then e's. Either block of loops may be empty. Its value is e or f precisely when it is constant at that loop. These descriptions include literals and follow by multiplying one further letter. In particular every nonempty subsequence of a word below c stays below c.

Here is the positive-path consequence for an arbitrary upper P. In the variable graph put an arc x \ensuremath{\to} y for every adjacency xy in an original upper word. An absorbed word uses only variables of P (give a missing variable z and all others e). Each of its adjacent pairs xy requires a positive-length directed path from x to y in the variable graph. To see this, give a variable two vertices L, R with a directed arc L \ensuremath{\to} R, and give each original adjacency xy the arc \(R_{x}\) \ensuremath{\to} \(L_{y}\). The values f, c, e represent the pairs 11,\allowbreak{}10,\allowbreak{}00. If there is no path from \(R_{x}\) to \(L_{y}\), set to 1 all vertices that can reach \(L_{y}\) (including itself), and to 0 the rest. Every arc is nonincreasing. Hence each upper word is below c. The target has an increase from 0 to 1 at xy, giving z. A path from \(R_{x}\) to \(L_{y}\) must traverse at least one adjacency arc even when \(x=y\), since \(R_{x}\) and \(L_{x}\) are distinct vertices.

At bound a use the increasing endpoint representation \(e=00\), \(a=01\), \(f=11\). An upper word has value a exactly when these edges compose from 0 to 1. This statement also accounts for all assignments using c or z, since a word containing either cannot give a.

\par\smallskip\noindent\textit{Critical inequalities and the shape of every possible summand.}\quad 
Choose distinct \(p_{i}\), \(q_{i}\), \(r_{i}\) with indices modulo n, for any \(n\geq 3\). Let \(P_{n}\) contain all 2n words

{\small\par\smallskip
\CaseLine{\(p_{i} q_{i},\qquad p_{i} r_{i} q_{i+1}.\)}
\smallskip}

Put \(s_{i}=p_{i}\) \(q_{i+1}\). Every \(s_{i}\) \ensuremath{\leq} \(P_{n}\) holds. At bound c delete the middle letter from the corresponding original triple. At a loop bound all its letters equal that loop. At bound a the pair \(p_{i}\) \(q_{i}\) gives \(p_{i}\) the increasing edge 0 to \(b_{i}\) and \(q_{i}\) the edge \(b_{i}\) to 1 for some bit \(b_{i}\). Each triple forces \(r_{i}\) to have endpoints \(b_{i}\), \(b_{i+1}\). The cyclic nondecreasing bits are constant, so \(p_{i}\) \(q_{i+1}\) composes from 0 to 1. The top bound is automatic.

Conversely every word absorbed by \(P_{n}\) is an original pair, an original triple, or one of these n shortcuts. Set all p's and r's to e and all q's to a. The upper words are all a. A word below them must contain exactly one q, in final position. Giving all p's the value a and all r's and q's the value f forces exactly one p, in initial position. In the positive reachability graph the only successors of \(p_{i}\) are \(q_{i}\), \(r_{i}\), \(q_{i+1}\); the only successor of \(r_{i}\) is \(q_{i+1}\); q's have no successors. Together with the missing-variable test this proves the asserted shape for arbitrary word lengths.

For i different from j remove the triple with private letter \(r_{j}\). Make \(b_{i+1}\),..., \(b_{j}\) equal to 1 in cyclic order and all other bits 0. Use the above endpoints for p's and q's and, for k different from j, for \(r_{k}\). They are nondecreasing wherever used. All the original words still present evaluate to a. Choose any value for \(r_{j}\). The proposed shortcut \(s_{i}\) evaluates to \(ef=z\). Thus even an arbitrary subset of the old words omitting the jth triple cannot entail \(s_{i}\).

\par\smallskip\noindent\textit{Excluding all bounded-variable derivations.}\quad 
Fix k and take \(n>3(k+1)+2\). Suppose the critical identity were derived from identities valid in 11636 using at most k variables per identity, allowing the semiring axioms as well. Expand every term into an idempotent sum of nonempty words. Every intermediate polynomial is equal to \(P_{n}\) in 11636, so it can contain only words of the proved shape. Look immediately before the first step introducing any shortcut \(s_{i}\). All words present are original pairs or triples. Such a step exists in a derivation from \(P_{n}\) to \(P_{n}\)+\(s_{0}\), since those endpoints have different expansions. Semiring axioms themselves preserve these expansions.

Select an actual occurrence of \(s_{i}\) contributed by the replacement of an axiom side. In a one-hole term context the changed terms of the expansion consist of outside left and right factors multiplying the expanded redex. Fix the two factors used in the selected word. Together they contain at most one actual letter because \(s_{i}\) has length two and all substituted words are nonempty. Unchanged summands cannot introduce a first shortcut.

Choose the axiom-side monomial yielding \(s_{i}\). It has at most two occurrences. From the polynomial image of each variable retain the particular words chosen at its occurrences in this output. If both occurrences use the same variable and different words, retain both. Keep one original image word for every other variable. This is an ordinary polynomial substitution with at most k+1 retained image words in all. The selected \(s_{i}\) still occurs. Expand the old side with these restricted images and the same two fixed outside factors. Every resulting word was an actual word in the previous polynomial, so has length at most three. Consequently every retained image word used on this old side has length at most three. The union of all letters in these old inputs has cardinality at most 3(k+1)+1, including the possible outside letter. Images unused on the old side contribute nothing to the bound.

There is a private \(r_{j}\) missing from this union with j different from i. On the preceding separating assignment every retained old input has value a, while \(s_{i}\) has value z. This contradicts validity of the axiom under the restricted substitution in its fixed multiplicative context. The proof treats occurrences in the full polynomial images independently, in either orientation of every axiom. Since k was arbitrary, there is no finite basis.

\par\smallskip\noindent\textit{Catalogue transfer.}\quad 
The same argument applies to any B such that 11636 belongs to Var(B) and B satisfies \(s_{0}\) \ensuremath{\leq} \(P_{n}\) for every \(n\geq 3\). Every candidate axiom of B is then valid in 11636, and the critical identity in B would require the impossible bounded-variable derivation. The first condition alone does not provide the critical identities.

For an arbitrary finite table B, test the whole family as follows. At each upper bound d form the vertices (p, q) with pq \ensuremath{\leq} d. Place an edge from (p, q) to (u, v) if some r gives prv \ensuremath{\leq} d. Call it bad if pv is not below d. An assignment failing the critical family gives a closed walk using a bad edge. Conversely a bad edge on a directed cycle supplies an assignment using a separate r on each edge. Repeat the cycle if needed to reach length at least three. This gives an exact all-length test. The family includes its reverse after exchanging the p and q symbols and reversing the cyclic order. The source is also anti-isomorphic to itself by interchanging e and f.

Among the 57 previously undecided types, seven pass this all-length test, and exactly 11636, 11672 and 11673 also have the required HSP inclusion. All three have the reverse inclusion, so their nonfinite-basis proofs also follow by equality of varieties. Their limit-variety status remains open.

\subsection{Six identities from positive paths between adjacent variables}\label{case:256}
Source 11638 in the five-element catalogue. The signature is \((+, multiplication)\) without constants; bases are relative to the additively idempotent semiring axioms. Labels below name elements only.

\par\smallskip\noindent\textit{1. The algebra and the proposed basis.}\quad 
The tables in order 0,\allowbreak{}1,\allowbreak{}2,\allowbreak{}3,\allowbreak{}4 are

{\small\par\smallskip
\noindent addition       multiplication\par
\CaseLine{\ttfamily 0 0 0 4 4\qquad 4 0 4 4 4}
\CaseLine{\ttfamily 0 1 0 4 4\qquad 4 1 4 3 4}
\CaseLine{\ttfamily 0 0 2 4 4\qquad 0 0 2 4 4}
\CaseLine{\ttfamily 4 4 4 3 4\qquad 4 4 4 4 4}
\CaseLine{\ttfamily 4 4 4 4 4\qquad 4 4 4 4 4}
\smallskip}

Put \(e=2\), \(f=1\), \(p=0\), \(q=3\) and \(Z=4\). The additive order has e, \(f<p<Z\) and \(q<Z\), with q incomparable to e, f, p. Direct checks of the displayed tables verify associativity and distributivity on all triples.

Theorem. The following six absorptions form an identity basis:

\CaseDisplay{\begin{aligned}
y&\leq xy, &(S_1)\\
xz&\leq xyz, &(S_2)\\
xyz&\leq xy+yz, &(S_3)\\
xyz&\leq xyt+yz, &(S_4)\\
xy&\leq xyz+y, &(S_5)\\
x&\leq xy+z+zt. &(S_6)
\end{aligned}}

Each inequality denotes an additive absorption identity. Reversing all words gives a finite basis for the multiplicative opposite.

\par\smallskip\noindent\textit{2. Word values.}\quad 
Let B=\{e, p, f\}. Associate to its elements the endpoint pairs

\CaseDisplay{e=(0,\allowbreak{}0),\quad p=(0,\allowbreak{}1),\quad f=(1,\allowbreak{}1).}

A word with all factors in B is nonzero, that is, different from Z, exactly when every pair of consecutive factors satisfies \(R_{x}\leq L_{y}\). In that event the product has endpoints (\(L_{first}\), \(R_{last}\)) and equals the additive join of all its factors. These statements follow by induction: the nonzero binary products are \(ee=e\), \(ep=ef=pf=p\), \(ff=f\).

A word containing q or Z has value q precisely when its last factor is q and all earlier factors are f; otherwise it has value Z. This also includes a linear q. A factor Z always gives Z. This describes every nonempty word, without any bound on its length.

\par\smallskip\noindent\textit{3. Exact absorption criterion.}\quad 
Let P be a nonempty sum of words. Write D for the union of their contents, J for the variables occurring before a final position, and F for the set of their final variables. Linear words contribute to D, F but not J. Form a directed graph G with one edge \(x\to y\) for every consecutive pair xy in an original right word.

Lemma. A nonempty word w is absorbed by P if and only if:

1. c(w) is contained in D;

2. for every consecutive pair xy in w there is a positive-length directed path from x to y in G;

3. either J meets F, or every nonfinal variable of w belongs to J and the final variable of w belongs to F.

No zero-length reflexive paths are allowed. In particular a pair xx requires a genuine directed cycle.

Proof of necessity. A missing content variable can be assigned Z, with all others e. The target is Z and the right side e.

Assume content inclusion and let xy be a target pair with no positive path \(x\to y\). Make two vertices \(L_{v}\), \(R_{v}\) for each v in D, and put arrows \(L_{v}\to R_{v}\) and \(R_{v}\to L_{u}\) for every original edge \(v\to u\). There is no path from \(R_{x}\) to \(L_{y}\), since such a path would yield a positive path from x to y in G. Assign Boolean value 1 to all vertices reachable from \(R_{x}\), including \(R_{x}\) itself, and 0 to the others. Every arrow respects \ensuremath{\leq}, so \(L_{v}\leq R_{v}\) and every original adjacency satisfies \(R_{v}\leq L_{u}\). Each variable therefore has one of the three B-values. Every right word is nonzero and \(P\leq p\), but \(R_{x}=1\), \(L_{y}=0\) makes \(w=Z\). This separates the failed path condition, including \(x=y\).

Finally suppose J and F are disjoint. Assign J the value f and F the value q, consistently; assign outside variables e. Every right word equals q. Under content inclusion the target equals q exactly when its nonfinal variables belong to J and its final variable belongs to F. Thus condition 3 is necessary.

Proof of sufficiency. If \(P=Z\) there is nothing to show. If \(P=e\) or \(P=f\), every right word equals that element, and the word description forces every variable of D to equal it. Content inclusion gives the same value for w. If \(P=p\), every right word lies in B and is nonzero, and hence every variable of D lies in B. All original inequalities \(R_{v}\leq L_{u}\) hold. They imply \(R_{x}\leq L_{y}\) along every positive path in G, using \(L_{v}\leq R_{v}\) at intermediate vertices. Condition 2 makes every target adjacency valid, so w lies in B and is below p. If \(P=q\), all right words equal q; consequently J is assigned f and F is assigned q, and they are disjoint. Condition 3 then gives \(w=q\). This proves the criterion.

All six proposed laws satisfy this criterion (or can be checked directly on the finite tables), establishing soundness.

\par\smallskip\noindent\textit{4. An unrestricted derivation.}\quad 
Every use below is a substitution by nonempty words. No normal-word assumption is used, and the original words may have arbitrarily many variables or repetitions.

Suffix extraction. Any nonempty suffix V of an original right word UV is derivably below P, by \(S_{1}\) with \(x=U\), \(y=V\) when U is nonempty; when U is empty the suffix is already an original right word.

Recover the final variable. Put \(t=t(w)\). If an original word ends in t, extract its one-letter suffix. Otherwise condition 3 requires a crossing g in J intersect F. Choose original \(a=A\) g B and \(b=C\) g where B is nonempty. Suffix extraction gives \(gB\leq P\) and \(g\leq P\). Content inclusion supplies an original word containing t. Since no original word ends in t, extract a suffix tV with V nonempty. An instance of \(S_{6}\) gives

\CaseDisplay{t\leq tV+g+gB\leq P.}

This also settles every linear target.

Extend a known suffix along one original edge. Suppose s starts with y and \(s\leq P\) has been derived. Let \(x\to y\) be an original edge. Extract the original suffix xyR containing that occurrence; R may be an omitted substring, but is never an empty substitution.

If \(s=y\) is linear and R is absent, xy is already known. If \(s=y\) and R is nonempty, \(S_{5}\) gives \(xy\leq xyR+y\leq P\).

If \(s=yV\) with V nonempty and R is absent, \(S_{3}\) gives \(xyV\leq xy+yV\leq P\). If R is nonempty, \(S_{4}\) gives \(xyV\leq xyR+yV\leq P\). Thus in every case \(xs\leq P\) is derived.

Extend along a positive path and shorten it. If \(x=v_{0}\to v_{1}\to \cdots \to v_{k}=y\) is a positive path, repeatedly apply the preceding construction, from its last edge backwards, to a known suffix s starting with y. This derives

\CaseDisplay{v_0v_1\cdots v_{k-1}s\leq P.}

If \(k=1\) this is xs. Otherwise \(S_{2}\) with nonempty substitutions \(x=v_{0}\), \(y=v_{1}\cdots v_{k-1}\), \(z=s\) derives \(xs\leq P\). The same procedure works for \(x=y\): a positive cycle produces the required repeated letter, whereas a zero-length path would not justify that step.

Start with the final variable t already recovered. Work backwards through the consecutive pairs of w. Condition 2 supplies a positive path for each pair, so each step recovers the next longer target suffix. The last step derives \(w\leq P\). This is an explicit unrestricted equational construction using only \(S_{1}\)--\(S_{6}\).

\par\smallskip\noindent\textit{5. Completeness and duality.}\quad 
By distributivity any term is a finite nonempty sum of words. If two such polynomials are equal in the algebra, every summand on either side is absorbed by the other side. The criterion and Section 4 derive every such absorption from the six laws. Adding them derives both polynomial inclusions, hence the identity. This proves completeness for every finite set of variables, with no common bound on that set.

Reversing words preserves the semiring axioms and transports validity and derivations to the multiplicative opposite. The reversed six laws therefore form its basis.

Independent computations verify the algebra axioms, the value formula, the criterion, and actual nonempty substitutions in every branch of the derivation. They are error checks, not the unrestricted proof. Other catalogue entries require independently verified equality of the generated varieties in both HSP directions.

\subsection{Ten finite bases by shortening guarded positive paths}\label{case:257}
For the six head-condition types use xyzt = xyt + xzt in the indicated orientation. For the four ordered-path types use xyzt = xyz + xzt in the indicated orientation. Replacing four consecutive letters reduces words to length at most three. Keep all original upper summands whenever an absorbed word is added. All paths used below are of positive length.

\par\smallskip\noindent\textit{Six types with a head condition.}\quad 
Use 11642, 11643 and the multiplicative opposite of 12100 as oriented sources. The labelled opposite pairs are (11642,\allowbreak{}13170), (11643,\allowbreak{}12109) and (12100,\allowbreak{}12323). At the principal ideal below 0 a reduced word is in the ideal exactly when its head belongs to H and, unless the head is g, all its noninitial values are in S.

\noindent Source\enspace\textbar\enspace H\enspace\textbar\enspace g\enspace\textbar\enspace S\par
\noindent 11642 or 11643\enspace\textbar\enspace 0,\allowbreak{}1,\allowbreak{}2\enspace\textbar\enspace 2\enspace\textbar\enspace 1\par
\noindent opposite of 12100\enspace\textbar\enspace 0,\allowbreak{}1,\allowbreak{}2,\allowbreak{}3\enspace\textbar\enspace 1\enspace\textbar\enspace 2,\allowbreak{}3\par
At all other bounds the relations of word lengths one through three are Cartesian products of position domains. The following laws hold:

{\small\par\smallskip
\CaseLine{\(xz \leq  xyz\)}
\CaseLine{\(xz \leq  xy + yz\)}
\CaseLine{\(xz \leq  xyt + yz\)}
\smallskip}

Put a directed edge from the head of each reduced upper word to each noninitial occurrence. Include literal heads in the set of all upper heads. For validity at the ideal below 0 the target head h must be an upper head, and every noninitial target variable v must be a positive successor of h. Indeed a missing head can be assigned outside H and all other variables g. For a missing successor v different from h, give h and all its positive successors a fixed value in S; give all others g. For \(v=h\) without a positive return, use a value in H outside S and different from g at h, a value in S at its positive successors and g elsewhere. Each assignment puts the upper words in the ideal and the target outside. Conversely, whenever h does not take g, each edge from h forces its successor into S. Membership in S propagates because g is outside S. This proves sufficiency.

Each edge x to y has a derived witness xy (extract the endpoints) or an original cubic xyt with y in the middle. For a positive path x to v, obtain a witness headed x with v noninitial, of length at most three. Start with the last edge. If yZ is already obtained, where Z is nonempty, the second or third law derives xZ from yZ together with xy or xyt. Proceed backwards. This works for repeated vertices and positive closed paths. Thus at most one original head witness and two derived successor witnesses suffice at this bound for any valid reduced target.

\par\smallskip\noindent\textit{Four types with an ordered path condition.}\quad 
Use 11717 and 12879, with labelled opposites 13762 and 13016. In D=\{0,\allowbreak{}1,\allowbreak{}2,\allowbreak{}3\} designate a low value l and a high value h. Put (x, y) in R exactly when x, y belong to D and either \(x=l\) or \(y=h\). Each of the following path ideals has its actual ideal as literal domain. Its nonlinear words belong to the ideal exactly when the first value is in \(H_{b}\) and every pair of consecutive values lies in R.

\noindent Source\enspace\textbar\enspace l\enspace\textbar\enspace h\enspace\textbar\enspace b\enspace\textbar\enspace \(H_{b}\)\par
\noindent 11717\enspace\textbar\enspace 3\enspace\textbar\enspace 1\enspace\textbar\enspace 0\enspace\textbar\enspace 0,\allowbreak{}1,\allowbreak{}2,\allowbreak{}3\par
\noindent 11717\enspace\textbar\enspace 3\enspace\textbar\enspace 1\enspace\textbar\enspace 2\enspace\textbar\enspace 2,\allowbreak{}3\par
\noindent 12879\enspace\textbar\enspace 3\enspace\textbar\enspace 2\enspace\textbar\enspace 0\enspace\textbar\enspace 0,\allowbreak{}1,\allowbreak{}2,\allowbreak{}3\par
\noindent 12879\enspace\textbar\enspace 3\enspace\textbar\enspace 2\enspace\textbar\enspace 1\enspace\textbar\enspace 1,\allowbreak{}2,\allowbreak{}3\par
All remaining bounds have Cartesian relations at each reduced arity. Use the following laws with the length reduction:

{\small\par\smallskip
\CaseLine{\(xy \leq  xyz\)}
\CaseLine{\(xz \leq  xy + yz\)}
\CaseLine{\(xz \leq  xy + tyz\)}
\smallskip}

Put edges between consecutive positions of the original words. Each edge x to y has a derived witness xy (including the first pair of a cubic) or an original witness axy. To derive a witness for a positive path, start with its first edge. From Uz, where U is x or ax, and a next-edge witness zv or tzv, use one of the last two laws to obtain Uv. Thus every positive path x to y has a derived witness ending in xy, including when \(x=y\). At either path ideal its last pair imposes R on those endpoints.

Pin the target variables to any assignment making the target fall outside one path ideal. To have P inside, each original edge must satisfy R. Intersect all unary restrictions at every upper variable, including D at each edge endpoint, \(H_{b}\) on nonlinear heads, the literal restrictions and any pins. Call the intersection L(x). An absent upper variable can retain its prescribed pin. If a list is empty, at most five failing restrictions already witness this; each comes from a pin or a single original word. Otherwise let A consist of the upper variables with l outside their list and let U be all their positive successors. Every member of U must be high. If each allows h, give U the value h, give each member of A outside U any allowed value, and give l elsewhere among upper variables. Each edge from a nonlow value ends high. This solves the upper constraints and pins. Consequently for an impossible instance with nonempty lists there exist x in A and a positive successor y which forbids h. Retain one unary reason excluding l at x, one excluding h at y and one derived path witness. The pinned instance remains impossible. This also handles \(x=y\).

\par\smallskip\noindent\textit{Unrestricted completeness.}\quad 
At a Cartesian bound, an empty upper relation is witnessed by that word. Otherwise failure with any fixed target pins means the unary restrictions at some variable have empty intersection. At most five pins or words witness this. For any valid reduced absorption w \ensuremath{\leq} P, take the union of the selected upper words over every bound and every failing assignment to the at most three target variables. There are at most \(5^{4}\) cases. For the six head-condition types the exceptional bound uses at most three words even without pins. All cases fit \(B=24*5^{4}\) upper words. If the union is empty, add any original word. The chosen subpolynomial absorbs w: otherwise take the bound equal to its evaluated value, contradicting the selected witnesses. The finite family of valid absorptions with term lengths at most three, at most B upper words and alphabet size 3(B+1), the structural laws and the semiring axioms now derive the original absorption. Restore omitted upper words and discharge derived ones by their path derivations. Apply this to all reduced summands in both directions of a polynomial identity and undo the reductions. Reversing words and transporting along the certified label permutations yields the opposite bases.

\subsection{两步尾部保护：来源11650的七条完整有限基}\label{case:258}
\par\smallskip\noindent\textit{1. 半环与词求值.}\quad 
取 \(A=\{z, e, b, d, t\}\)，依次对应目录11650的标号 \((0,\allowbreak{}1,\allowbreak{}2,\allowbreak{}3,\allowbreak{}4)\)。加法为偏序

\CaseDisplay{b<z<t,\qquad e<z<t,\qquad e<d<t}

上的并，除传递关系外无其他比较。乘法由右平移定义：右乘 \(e\) 为恒等映射，右乘 \(d\) 为

\CaseDisplay{P=(t, d, b, d, t),}

右乘 \(z, b, t\) 均为常值 \(t\) 映射 \(T\)。因此乘法表为

\CaseDisplay{\begin{array}{c|ccccc}
\cdot&z&e&b&d&t\\\hline
z&t&z&t&t&t\\
e&t&e&t&d&t\\
b&t&b&t&b&t\\
d&t&d&t&d&t\\
t&t&t&t&t&t
\end{array}}

三个右平移保持并、彼此交换并且幂等；直接按右因子类别可检验 \(R_{xy}=R_yR_x\) 及 \(R_{x+y}=R_x\vee R_y\)，所以给出加法幂等半环。注意 \(b\) 不是乘法零。

词 \(w\) 的首字母记为 \(h(w)\)，尾部变量集合为 \(D(w)\)。词值可直接描述：尾部全取 \(e\) 时取首字母值；尾部仅取 \(e, d\) 且至少出现一个 \(d\) 时取 \(P(h(w))\)；尾部出现 \(z, b, t\) 时取 \(t\)。一次词采用空尾部规则。

\par\smallskip\noindent\textit{2. 精确吸收判据.}\quad 
设 \(F\) 为非空有限词集，代表右多项式 \(\sum_{u\in F}u\)。定义

\CaseDisplay{H=\{h(u):u\in F\},\quad E=\bigcup_{u\in F}D(u),\quad
E_h=\bigcup_{\substack{u\in F\\h(u)=h}}D(u),\quad
K=\bigcup_{\substack{u\in F\\h(u)\in E}}D(u).}

集合 \(K\) 表示两步保护：若某变量 \(g\) 出现在任一右项尾部，而以 \(g\) 为首字母的右项又在尾部含 \(x\)，则 \(x\in K\)。它不要求两个见证词相同，也不限制见证词内容包含于左词内容。

定理1。 令 \(h=h(w)\)、\(D=D(w)\)，则

\CaseDisplay{A\models w\le F\quad\Longleftrightarrow\quad
h\in H\cup K,\qquad D\subseteq E,\qquad D\subseteq H\cup K\cup E_h.\quad\textup{(\ensuremath{1})}}

\par\smallskip\noindent\textit{必要性：统一分离赋值.}\quad 
若 \(x\in D\setminus E\)，令 \(x=b\)，其余变量均取 \(e\)。右项尾部全为 \(e\)，故右和不超过 \(z\)；左尾部含 \(b\)，左值为 \(t\)，矛盾。

现在假定 \(D\subseteq E\)。若 \(h\notin H\cup K\) 且 \(h\notin E\)，则 \(h\) 不出现在右侧。令 \(h=t\)，其余变量取 \(e\)，左值为 \(t\)、右值为 \(e\)。

为处理其余情形，对任意 \(x\in E\setminus(H\cup K)\) 定义

\CaseDisplay{P_x=\{h(u):u\in F,\ x\in D(u)\}.}

由于 \(x\notin K\)，有 \(P_x\cap E=\varnothing\)；由于 \(x\notin H\)，有 \(x\notin P_x\)。因此可一致地赋值：\(x=d\)，\(P_x\) 中变量取 \(b\)，其余变量取 \(e\)。每个右项的尾部仅取 \(e, d\)；尾部含 \(d\) 的右项首字母为 \(b\)，仍取 \(b\)，其他右项取 \(b\) 或 \(e\)。故右和不超过 \(z\)。

若 \(h\in E\setminus(H\cup K)\)，把上述 \(x\) 取为 \(h\)，左首字母为 \(d\)，左尾部因 \(D\subseteq E\) 而仅取 \(e, d\)，所以左值为 \(d\)，不小于等于 \(z\)。这证明第一条件。

最后，若 \(x\in D\setminus(H\cup K\cup E_h)\)，则 \(x\in E\setminus(H\cup K)\)。使用同一赋值。已经成立的 \(h\in H\cup K\) 保证 \(h\ne x\)，且 \(x\notin E_h\) 保证 \(h\notin P_x\)，故左首字母为 \(e\)。左尾部含 \(d\) 且其余仅为 \(e, d\)，左值为 \(d\)，又与右和不超过 \(z\) 矛盾。

\par\smallskip\noindent\textit{充分性：真主理想.}\quad 
真主理想为 \(\{b\}\)、\(\{e\}\)、\(\{b, e, z\}\)、\(\{e, d\}\)。条件 (1) 蕴含 \(c(w)\subseteq H\cup E=c(F)\)。

若所有右项都不超过 \(e\)，则它们的全部变量取 \(e\)。若所有右项都不超过 \(d\)，则它们的全部变量取 \(e, d\)。这两种情形由内容包含立即得左值分别不超过 \(e, d\)。

若所有右项不超过 \(b\)，则全部右首字母取 \(b\)，右尾部取 \(e, d\)，从而 \(H\cap E=\varnothing\)、\(K=\varnothing\)。第一条件给出 \(h\in H\)，第二条件给出 \(D\subseteq E\)，故左值为 \(b\)。

若所有右项不超过 \(z\)，则右首字母取 \(b, e, z\)，右尾部取 \(e, d\)，且尾部出现 \(d\) 的右项首字母必须为 \(b\)。因此 \(H\cap E\) 中变量均取 \(e\)。以这些变量为首字母的右项尾部也必须全取 \(e\)，所以 \(K\) 中变量均取 \(e\)。由 \(h\in H\cup K\)，左首字母取 \(b, e, z\)。若左尾部全取 \(e\)，结论成立；否则选一个取 \(d\) 的左尾变量 \(x\)。它既不在 \(H\) 也不在 \(K\)，故第三条件迫使 \(x\in E_h\)。相应右项的首字母 \(h\) 必须取 \(b\)，于是左值也是 \(b\)。所有真主理想均已检查，充分性得证。

\par\smallskip\noindent\textit{3. 七条完整基.}\quad 
定理2。 \(A\) 的恒等式具有下列七条完整附加基：

\CaseDisplay{\begin{aligned}
xyz&\approx xzy, &(N)\\
xy^2&\approx xy, &(R)\\
x&\le xy, &(U)\\
xyz&\le xy+xz, &(M)\\
xy&\le x+y^2, &(Q)\\
x^2&\le x+yx, &(S)\\
x^2&\le yx+zy. &(P)
\end{aligned}\quad\textup{(\ensuremath{2})}}

七式有效性由定理1逐式检查；前两式也直接来自右平移交换、幂等。因此下面只证明完备性。

由 \((N),(R)\)，每个词由首字母和尾部集合正规化，尾部字母可任意换序并去重。由 \((U)\)，对于右项 \(u\)，可推出 \(h(u)\le u\)；若 \(x\in D(u)\)，也可推出 \(h(u)x\le u\)，因为先把 \(x\) 换到首个尾部位置，再用 \((U)\) 添加其余非空尾部。没有其余尾部时直接使用已有词，不代入空词。

平方提升。 若 \(x\in H\cap E\)，分别选一个以 \(x\) 开头的右项及一个尾部含 \(x\) 的右项，后者首字母记为 \(y\)。前段给出 \(x\le F\)、\(yx\le F\)，由 \((S)\) 得 \(x^2\le F\)。

若 \(x\in K\)，则存在右项首字母 \(y\in E\) 且尾部含 \(x\)，并存在另一右项尾部含 \(y\)，其首字母记为 \(z\)。前段给出 \(yx\le F\)、\(zy\le F\)，由 \((P)\) 得 \(x^2\le F\)。即

\CaseDisplay{x\in (H\cap E)\cup K\quad\Longrightarrow\quad x^2\le F.\quad\textup{(\ensuremath{3})}}

给定有效吸收 \(w\le F\)，使用定理1。先导出 \(h\le F\)：若 \(h\in H\)，删除相应右项尾部；若 \(h\in K\)，由 (3) 及 \((U)\) 的 \(h\le h^2\) 得到。

接着逐个添加 \(D(w)\) 中变量。维持已导出词 \(p\le F\)，其首字母为 \(h\)、尾部集合包含于 \(D(w)\)。

\textbullet{} 对尚未添加且属于 \(H\cup K\) 的尾变量 \(x\)，第二条件保证 \(x\in E\)，所以 (3) 给出 \(x^2\le F\)。由 \((Q)\) 代入 \(p, x\)，得 \(px\le p+x^2\le F\)。

\textbullet{} 对其余尾变量 \(x\)，第三条件保证 \(x\in E_h\)，所以 \(hx\le F\)。若 \(p=h\)，直接采用此词；否则写 \(p=hv\)，其中 \(v\) 非空，由 \((M)\) 得 \(px=hvx\le hv+hx\le F\)。

有限次操作后得到首字母为 \(h\)、尾部集合恰为 \(D(w)\) 的词，即由 \((N),(R)\) 等于 \(w\)。因此七式推出每个有效单词吸收。分配律将任意项化为非空词和，再逐词处理恒等式的两向吸收，即得任意变量数的完整基。证毕。

\par\smallskip\noindent\textit{4. 结论与证据范围.}\quad 
此证明不需要加入 \(x^2y\approx y^2x\)，虽然来源11650满足该式；它已经是完整基的推论。乘法对偶通过逐词反转取得完整基。

一至五变量全部右多项式主理想CNF比较仅是独立的有限范围交叉检查；任意变量数结论依赖定理1和定理2的证明。目录11650及其对偶14105均为有限基，因而不能生成极限簇。其他类型的有限基性必须另行证明或具有独立双向HSP证据，不能仅由满足七式而传递。

\subsection{Two finite bases from original first-pair witnesses}\label{case:259}
Signature (+, multiplication), without constants. This proves finite basability separately for catalogue cores11703 and 11718. Besides the additively idempotent semiring axioms take

{\small\par\smallskip
\CaseLine{\(xyzt=xytz,  xyzz=xyz,  xyzt\leq xyz+xyt,\)}
\smallskip}

and every table-valid absorption over \(B=2^{21}-1+42=2\),097,\allowbreak{}193 variables, word length at most B+2, with at most 20 upper summands. This is an effectively specified finite family, not an enumerated small basis. All variables are substituted by nonempty words or polynomials.

\par\smallskip\noindent\textit{Word values and principal ideals.}\quad 
11703 multiplication rows are 40444,\allowbreak{}41444,\allowbreak{}40404,\allowbreak{}43444,\allowbreak{}44444. A nonlinear word whose noninitial values are all 1 has value 0,\allowbreak{}1,\allowbreak{}0,\allowbreak{}3,\allowbreak{}4 according as its head is 0,\allowbreak{}1,\allowbreak{}2,\allowbreak{}3,\allowbreak{}4. The only other non4 case has first values(2,\allowbreak{}3), remaining values 1, and value 0. For 11718 the rows are 40444,\allowbreak{}41444,\allowbreak{}42444, 22444,\allowbreak{}44444. An all 1 noninitial sequence gives head 0,\allowbreak{}1,\allowbreak{}2 unchanged, head 3 gives 2, and head 4 gives 4. The only other non4 case has first values(3,\allowbreak{}0), remaining values 1, and value 2. These statements follow by induction from the rows, including every length. Literals retain values. They verify the three displayed laws. The first two normalize a word to its ordered first two positions followed by the set of subsequent variables, in any fixed order. Positions after the second permute and duplicate occurrences there disappear; the first two positions are kept.

For an original nonempty polynomial P let C be its content, L its literal variables, A its nonlinear heads, B0 its second variables, D its variables at positions at least 3, \(H=L\) union A, \(T=B0\) union D, and E its original ordered first pairs. These sets always refer to originals, not generated words. The proper principal ideals of 11703 are U=\{0,\allowbreak{}1,\allowbreak{}2\},\{1\},\{2\}, V=\{1,\allowbreak{}3\}; for 11718 they are \{0,\allowbreak{}1\},\{1\}, U=\{2,\allowbreak{}3\},\{3\}. An absorption \(w\leq P\) holds iff, for every proper principal ideal, every assignment putting all originals in that ideal puts w there. This is exact because when w is not below P, the principal ideal of the value of P separates them.

In either core,\{1\} imposes precisely content(w) subset C:every word has value 1 iff all its variables are 1. The singleton ideal containing no nonlinear word values imposes the existence of an original nonlinear word when w is nonlinear; when w is literal, content coverage suffices. The remaining easy ideal imposes noninitial content(w) subset T for a nonlinear target. Indeed original heads range over\{1,\allowbreak{}3\} in 11703 and \{0,\allowbreak{}1\} in 11718, and noninitial variables equal 1; every content variable is an allowed target head. Necessity uses a target noninitial variable outside T with value 3 or 0 respectively, all other heads and tails 1.

\par\smallskip\noindent\textit{Exact difficult-ideal criterion for 11703.}\quad 
All originals lie in U=\{0,\allowbreak{}1,\allowbreak{}2\} iff literal variables lie in U, each original pair has values in R=\{(0,\allowbreak{}1),(1,\allowbreak{}1),(2,\allowbreak{}1),(2,\allowbreak{}3)\}, and \(D=1\). Thus heads H lie in U, seconds B0 in\{1,\allowbreak{}3\}. Define

{\small\par\smallskip
\CaseLine{\(F=D \cup (B0 \cap  H) \cup \{s:(h, s) in E, h in T\},\)}
\CaseLine{\(K=H \cup  D \cup  F.\)}
\smallskip}

Exactly F is forced to 1, and exactly K is forced into U. For clarity, every second 3 forces all its original predecessor heads to 2. The only forced 1 seconds are those also heads, those in D, and successors of heads which are themselves noninitial. No iterative propagation adds anything:all forced 1 heads already belong to T. For an unforced second s, assign \(s=3\), its original predecessor heads 2, all other seconds and D1, and remaining heads freely in U. Predecessor heads cannot lie in T, by the definition of F, so this is consistent. This proves the claim about F and K as well as existence of all requisite separators.

A target pair(h, s) is universally in R exactly when it is an original pair in E, or h in K and s in F. Sufficiency is immediate. For necessity, if s is not forced 1 and is a second, put it 3; the only heads then forced 2 are its original predecessors. A head which is not such a predecessor can take a value making the new pair fail. If s is not a second or a tail, give it 0 instead. If s is forced 1 but h is not in K, give h3. Coincident \(h=s\) causes no exception:an original loop already forces 1; otherwise give that variable 0 if a head, or 3 if an unforced second. Hence a literal w is valid exactly when its variable belongs to K. A nonlinear w is valid exactly when there is an original nonlinear, its content is in C, its noninitial content in T, its first pair satisfies the stated test, and its remaining tail is contained in F.

An original literal MAY supply the target head:for example \(yx\leq y+x^{5}\) is valid. Requiring all target content in the nonlinear-original content union is false.

\par\smallskip\noindent\textit{Exact difficult-ideal criterion for 11718.}\quad 
For U=\{2,\allowbreak{}3\}, original heads lie in\{2,\allowbreak{}3\}, seconds in\{0,\allowbreak{}1\}, \(D=1\), and R=\{(2,\allowbreak{}1),(3,\allowbreak{}1),(3,\allowbreak{}0)\}. The system is inconsistent precisely when H intersects T. If disjoint, it is satisfied by giving all heads 3 and all noninitial variables 1. In that case exactly D is forced 1; every second outside D may be 0, which forces its original predecessor heads 3. All other heads may independently be 2 or 3. Therefore a literal target must belong to H; an allowed target pair is an original pair in E or a pair with head in H and second in D; its later tail must lie in D. The same predecessor assignment proves necessity, including coincident target positions. If H intersects T, this ideal imposes no condition. Combine this with the three easy-ideal requirements above. This proves the criterion for arbitrary alphabets, word lengths and polynomial sizes.

\par\smallskip\noindent\textit{Bounded original support for short targets.}\quad 
Suppose a valid target has length at most 3. Retain an original for each of its at most 3 content variables, one nonlinear original if needed, and an original noninitial occurrence for each of at most 2 noninitial variables. These cost at most 6 originals.

For 11703, membership in F has at most 2 original witnesses:one occurrence in D; or a second occurrence and a head occurrence; or its original pair (h, s) plus an occurrence of h in T. Membership in K likewise costs at most 2:one head, one D occurrence, or the preceding F witnesses. Retain one original pair, or at most 4 witnesses for head in K and second in F; also at most 2 witnesses for the third variable in F. A literal needs at most 2 witnesses for K. Thus at most 12 originals suffice.

For 11718, retain the 2 originals witnessing H intersect T when it exists. Otherwise retain one head for a literal; for a nonlinear target retain one original pair or a head occurrence and a D occurrence, and one D occurrence for its third variable. At most 9 originals suffice. In both cases 20 is a uniform safe bound. Every witness is from original P. The selected Q still satisfies all applicable criteria, so \(w\leq Q\) is valid.

\par\smallskip\noindent\textit{Reversible compression and completeness.}\quad 
Normalize Q and w with the first two laws. There are at most 21 words. Reserve every variable in either of their first two positions, at most 42 variables. Group all other variables by their nonempty incidence subset of the 21 words. There are at most \(2^{21}\)-1 groups. Replace each group by one fresh variable, leaving reserved variables distinct. The resulting absorption is valid by variable identification and belongs to the stated finite family:at most B variables, each word length at most B+2, and at most 20 upper words. Inversely substitute for each group variable the nonempty word listing its members. All group members occurred only beyond the second position and had COMMON incidence. The two structural laws therefore normalize each expanded word exactly back to its original. This derives \(w\leq Q\) and then \(w\leq P\). No empty word, loss of multiplicity outside the idempotent tail, or new witness is used.

For a longer valid normalized target hs z1...zk, each triple hs zi satisfies the same criterion against original P. Derive all those triples by the short-target argument. Successively use \(xyzt\leq xyz+xyt\) with \(x=h\), \(y=s\), \(z=the\) already accumulated nonempty tail, \(t=the\) next singleton. This derives the whole target. Normalize back if necessary. Thus every valid absorption follows. Each valid polynomial identity is equivalent to its finitely many summand absorptions in both directions. This proves unrestricted finite basability for both cores.

Computational independent ideal solvers, word-value tests and literal replays check this proof's implementation; their finite samples are not used in place of the unrestricted argument above. Mutual HSP and opposite transfers require their own verified certificates.

\subsection{Six complete identities for 11722 from original triggered tails}\label{case:260}
Signature (+, multiplication), no constants, relative to the additively idempotent semiring axioms. A complete finite basis is:

1. \(xyz=xzy\).

2. \(xyy=xy\).

3. \(x\leq xy\).

4. \(xy\leq xyz\).

5. \(xyz\leq xy+xz\).

6. \(xz\leq xy+y+tz\). Every substitution below is nonempty. Completeness is unrestricted in variables, word lengths and polynomial sizes.

Addition:00444,\allowbreak{}01444,\allowbreak{}44224,\allowbreak{}44234,\allowbreak{}44444. Its order consists of the two chains \(1<0\) and \(3<2\) below a common top 4; elements from different chains join to 4. Multiplication:40444,\allowbreak{}41444,\allowbreak{}42444,\allowbreak{}33333,\allowbreak{}44444. A word headed by 3 always has value 3. A word with any other head retains its head value exactly when every tail value is 1; otherwise it gives 4. A literal retains its value. Induction on the table proves this formula and the two word identities. Normalize to head plus sorted tail set; keep any repeated head in the tail. Tail swaps and deletions have nonempty substituted prefixes.

\par\smallskip\noindent\textit{Exact criterion from original positions.}\quad 
Let H be the set of original right-word heads, D the union of original right-word tails, and \(D_{h}\) the union of tails of original h-headed words. A literal h is absorbed iff h belongs to H. A nonlinear target hT is absorbed iff h belongs to H, T is contained in D, and either (i) T is contained in \(D_{h}\), or (ii) \(D_{h}\) intersects H. All original words are allowed, including ones with outside variables. The intersection in(ii) is witnessed by an actual original tail occurrence and an actual original head occurrence, which may belong to different words.

Necessity: if h is not a head, assign \(h=2\) and all other \(variables=3\). Every right word has value 3; the target is 2 or 4, so is not absorbed. If x in T is absent from D, assign \(x=0\) and all \(others=1\). Each right word has value 0 or 1, whereas the target is 4. If the first two conditions hold but both(i) and(ii) fail, assign \(h=2\), all \(D_{h}\) \(variables=1\), and all \(others=3\). This is consistent because \(D_{h}\) is disjoint from H. Every h-headed right word is 2, and every other original head is 3, hence the right sum is 2. A target tail outside \(D_{h}\) is 2 or 3, so the target is 4. These are separating assignments for all possible failed conditions.

Sufficiency: if the right sum is 4, absorption is automatic. If it lies in \{0,\allowbreak{}1\}, every original head is 0 or 1 and every original tail is 1. Thus h has value 0 or 1 and T has only 1; the target equals its head and is below the right sum. If the right sum is 3, all original heads are 3, so h and the target are 3. In the remaining case the right sum is 2: all heads are 2 or 3, and every tail of a word with head value 2 is 1. If h is 3 the target is 3. If h is 2, every \(D_{h}\) variable is 1 and therefore \(D_{h}\) cannot intersect H. Condition(i) then holds, all target tails are 1, and the target is 2. The literal case follows directly because each head value is below its own word value. This proves the criterion for arbitrary assignments and arbitrary polynomials.

Each of the six displayed laws is valid. The first two follow from word values. For law 3 the right head is x. Law 4 and law 5 have their target tails covered in \(D_{x}\). For law 6 H contains x, y, t, \(D_{x}\) contains y, and the target tail z is an original tail of tz, so the triggered condition applies.

\par\smallskip\noindent\textit{Unrestricted derivation.}\quad 
Normalize original words. Law 3 derives a literal h from any original h-word by substituting \(x=h\) and y its nonempty tail, if it is nonlinear. For a desired atom pt from an original p-word with t in its tail, reorder the tail to tB. If B is nonempty, law 4 with \(x=p\), \(y=t\), \(z=B\) derives pt; otherwise the word already is pt. Here t is a true tail occurrence, not merely the head of the original word.

For nonlinear target hT in case(i), extract ht for each t in T from an original h-word. Merge these atoms repeatedly using law 5 with \(x=h\) and y, z the nonempty already assembled and new tails. The result is hT.

In case(ii), choose q in \(D_{h}\) intersect H. Extract hq from an original h-word and extract the literal q from an original q-word by law 3. For each t in T, choose an original word headed p with t in its tail and extract pt. Law 6, with \(x=h\), \(y=q\), \(z=t\) and its fourth variable assigned p, derives \(ht\leq hq+q+pt\). This is valid also if \(q=h\) or variables coincide. Merge all ht by law 5. Every right word used is original or derived, and all substitutions are nonempty. The number of steps is finite for each finite target, without any uniform bound on its variables or right words. Thus every valid absorption follows. Distributivity and summandwise absorption in both directions derive every valid identity, proving the six-law basis. Opposite transfer requires independent complete HSP both ways.

Original head alone cannot replace a tail witness: \(xy\leq x+y\) fails at \(x=1\), \(y=0\). Likewise \(xy\leq xz+ty\) fails at \(x=2\), \(y=3\), \(z=1\), \(t=3\) when no trigger exists. These counterexamples are retained in independent checks.

\subsection{A finite basis for 11723 by equal-mode graph completion}\label{case:261}
Signature without constants. All substitutions are nonempty. The addition rows are 00444,\allowbreak{}01444,\allowbreak{}44224,\allowbreak{}44234,\allowbreak{}44444; the multiplication rows 40444,\allowbreak{}41444,\allowbreak{}42444,\allowbreak{}42004,\allowbreak{}44444. A nonlinear word whose noninitial values all equal 1 takes its head value if that is 0,\allowbreak{}1,\allowbreak{}2, and takes 2 if the head is 3. A word beginning3,\allowbreak{}2 or 3,\allowbreak{}3 with all later values 1 takes 0. All other nonlinear values are 4. Literals retain their values. Induction proves this description for every length.

\par\smallskip\noindent\textit{Exact original-word criterion.}\quad 
For original P use C, L, A, B, D, \(H=L\) union A, \(T=B\) union D as in the preceding role-domain proofs:content union, literals, nonlinear heads, seconds, and variables from position 3 onward. Original first pairs are edges. In the undirected graph on these variables, each connected component has two possible modes for P to lie in U=\{0,\allowbreak{}1\}:

N: heads in\{0,\allowbreak{}1\}, seconds1; S: heads 3, seconds in\{2,\allowbreak{}3\}.

Every literal lies in U, every D variable equals 1. A variable occurring both as head and second forces its two roles into the SAME mode: N, N permits 1; S, S permits 3; the two mixed possibilities are empty. Shared second roles also force equal modes. Therefore ordinary undirected connected components capture precisely all mode equalities. Original literal or D occurrences force their component to N. There are no S pins and no impossible system:the all 1 assignment always works. Choose unpinned modes independently, and actual values independently in the role-domain intersections. This describes every original U-assignment.

Let K be L union D together with all vertices of N-forced components. Let O be D together with all second variables of N-forced components. These are precisely the variables necessarily in U and necessarily 1. A normal head that is not also a second need not be 1. Define E* as all(a, b) with a an original nonlinear head, b an original second, and a, b in the same undirected component, including squares when permitted by those roles. A literal is universally in U iff it belongs to K. A nonlinear target hsZ is universally in U iff Z is contained in O and either(h, s) belongs to E*, or h belongs to K and s belongs to O. Independence of free modes and free actual values proves necessity: a free head can be 3, and only the same component's second guarantees the required2/\allowbreak{}3 value in that case. A pure second used as a free head may be 2 and fails. A separately fixed normal pair must have second 1.

All four proper principal ideals are\{1\},\{3\},\{2,\allowbreak{}3\}, U. The first gives target content contained in C. The second requires an original nonlinear word whenever the target is nonlinear, because nonlinear words never 3. The third allows original heads 2/\allowbreak{}3 and all noninitial values 1, so its system is impossible exactly when H intersects T. If disjoint, it requires target head in H and target noninitial content in T, also for literals. In particular a literal head need not occur as an original literal to pass this particular ideal; the separate U test still must hold. This point corrects an early exploratory conjecture before any certification. Together with the U criterion, these tests are exact for all variables.

For explicit separators:missing content gets4 while all others 1; missing nonlinearity is separated by all 3. In the\{2,\allowbreak{}3\} test put original heads 2 and all others 1. A missing target head then fails. For a missing noninitial variable, set it 2 if an original head, and 0 otherwise; originals stay2 and the target leaves\{2,\allowbreak{}3\}. Failed U is separated by choosing an allowed mode and actual values described above. Sufficiency follows by testing the actual value of P against its proper principal ideal.

\par\smallskip\noindent\textit{Ten fixed laws.}\quad 
{\small\par\smallskip
\CaseLine{\(xyzt=xytz\)}
\CaseLine{\(xyzz=xyz\)}
\CaseLine{\(xy\leq xyz\)}
\CaseLine{\(xz\leq xy+uy+uz\)}
\CaseLine{\(xyzt\leq xyz+xyt\)}
\CaseLine{\(xy\leq xu+yx+vy\)}
\CaseLine{\(xz\leq xy+yz\)}
\CaseLine{\(xz\leq xy+zy+vz\)}
\CaseLine{\(xz\leq yx+yz+xu\)}
\CaseLine{\(xz\leq yx+zy+xu+vz\)}
\smallskip}

The word values prove the two tail-normalization laws and prefix extraction. The exact ideal tests prove the rectangle, tail merge, reversal, and the four oriented two-edge path laws. For each path or reversal law, original role guards ensure a target head role and second role. Connected original edges put them in the same U-mode component. The content/nonlinearity tests hold directly; for the\{2,\allowbreak{}3\} test either H meets T or the displayed role guards supply the target endpoints. Thus validity is proved for unrestricted substitutions, not inferred from checking finitely many assignments.

Extract every original quadratic. To derive an E* edge(a, b), take any undirected path from a in A to b in B. A reversed single edge is corrected by the reversal law with original endpoint role guards. For a path of length at least 2, if a consecutive two-edge segment has eligible A/B endpoints in either order, apply its oriented path law with original guards and shorten the path by 1.

If no such segment exists, each distance-two endpoint pair consists of two A-only or two B-only vertices. A path of length 2 cannot have this property since its target endpoints are eligible. On a longer path, the first four vertices must therefore alternate exclusive roles. Their three edges have rectangle orientations \(x\to y\), \(u\to y\), \(u\to z\), possibly after reversing the whole segment. The rectangle law shortens this segment by 2. Repetition reaches a single edge. A zero-length target path(a, a) is handled by an original loop if present, or by the two-edge walk a, v, a using any original neighbor; a has both roles. All reductions use nonempty singleton substitutions and preserve original role guards. Thus every edge in E* is derived for arbitrary component size.

The undirected completed graph has diameter at most 2. If a component has a mixed-role vertex, it is adjacent to all other vertices because each has at least one eligible role. Without mixed vertices, the component is complete bipartite between A-only and B-only vertices. This finite diameter is a consequence of the derivations, not an assumption about original path lengths.

\par\smallskip\noindent\textit{Bounded support and complete finite family.}\quad 
For a target of length at most 3, retain at most 4 original words at each target variable to preserve literal, D, head and second occurrences. At most 12 words suffice. Link the target variables belonging to each original component by paths of at most 2 derived edges. There are at most 3 target variables, so at most 4 such edges suffice. For every represented component that is forced N, choose an original literal or D pin, retain its local witnesses(at most 4 words), and link it to the target component by at most 2 edges. At most 18 additional words suffice.

This support preserves all local target domains, all relevant mode equalities and all relevant N pins. Every assignment satisfying it extends on the target variables to an original U-assignment:use the same modes on represented components, then choose other modes freely and choose actual values in their nonempty local intersections. Extra constraints in selected words can only restrict target possibilities. Hence at most 34 words suffice for the U test. Add at most 3 original content witnesses, one nonlinear original, and at most 3 witnesses for the\{2,\allowbreak{}3\} test. The resulting bound 41 is below the conservative common bound 65 used below. All selected words are originals or derived E* edges.

Normalize the target and at most 65 upper words, protect their first two positions(at most 132 variables), and group remaining variables by their nonempty incidence patterns in the tails(at most \(2^{66}\)-1 groups). Replace groups by fresh variables. The absorption remains valid by substitution. The bound is \(B=2^{66}-1+132\) variables, length B+2, and 65 upper words. Inverse substitution replaces each group by its nonempty member product and each protected variable by its singleton. The first two fixed laws literally recover all original normal words, including protected variables that also occur in tails. Literals/quadratics involve only protected variables; long words retain a nonempty tail.

Take the ten fixed laws and all table-valid absorptions within these fixed bounds. This is an effective finite specification, not an enumerated small basis. Inverse compression derives every valid short target after graph completion, and hence from the original polynomial. Every valid long normal hsZ has all triples hsz(z in Z) valid. Derive these triples, merge their nonempty tails with the fifth law, and normalize. Thus every valid absorption follows with no original variable, length, or summand bound. Distributivity and both absorption directions prove completeness for all identities. Therefore 11723 is finitely based. Catalogue transfers require independently complete mutual HSP.

\subsection{Finite bases for 11734 and 11735 from normal components}\label{case:262}
The signature has addition and multiplication but no constants. Every word and every substituted block below is nonempty. The finite family specified here is effective, not an enumerated small basis.

Both multiplication tables are 40444,\allowbreak{}41444,\allowbreak{}42444,\allowbreak{}44044,\allowbreak{}44444. Thus a nonlinear word has value h in\{0,\allowbreak{}1,\allowbreak{}2\} if its head has value h and all later factors are 1; it has value 0 if its first two values are 3,\allowbreak{}2 and all later factors are 1; otherwise it has value 4. Literals retain their values. Induction proves this formula for every length. In particular no nonlinear word takes value 3.

For 11734 the addition rows are 00044,\allowbreak{}01044,\allowbreak{}00244,\allowbreak{}44434,\allowbreak{}44444. For 11735 they are 00444,\allowbreak{}01444,\allowbreak{}44244,\allowbreak{}44434,\allowbreak{}44444. Both have proper principal ideals \{1\},\{2\},\{3\}, and U, where U=\{0,\allowbreak{}1,\allowbreak{}2\} for 11734 and U=\{0,\allowbreak{}1\} for 11735. This distinction changes the original witnesses and cannot be discarded.

\par\smallskip\noindent\textit{Exact original-word semantics.}\quad 
For the original polynomial P let C be its content union, L its literal variables, A its nonlinear heads, B its second-position variables, and D its variables in positions 3 onward. Put \(H=L\) union A, \(T=B\) union D. An original nonlinear word contributes its first-two-position edge(h, s). Join heads sharing a second variable and take the resulting head components. Write E* for the complete bipartite head/second edges in these components.

Under assignments making every original value belong to U, each component has two potential modes:

N: heads belong to U, seconds have value 1; S: heads have value 3, seconds have value 2.

Every D variable must be 1, every literal must belong to U. Sharing a second forces the same mode, so the component description is necessary. The forced-N components are exactly those containing one of these heads:

a head in A intersect(L union T), for both sources; a head with a second in D, for 11734; a head with a second in H union D, for 11735.

These conditions are sufficient as well. A head also used noninitially can never be 3, so it is normal. In 11735, a second also used as a head or literal cannot be 2, so its incoming component is normal. In 11734, such a second MAY be 2:its own head component is normal, but its incoming component may remain special. Only D forces such an incoming second 1. These rules account for every overlap among L, A, B, D. Choose N/S independently for every unforced component, then intersect the stated value domains at each variable. All intersections are nonempty:the forbidden cases are exactly the forcing rules above. Within them, remaining actual values are independent. In particular all 1 always provides a consistent original-U assignment; there is no inconsistent mode case in either source.

Let O be the variables necessarily 1 under these assignments:

{\small\par\smallskip
\noindent O = D union seconds of forced-N components.\par
\smallskip}

Let K be the variables necessarily in U:

{\small\par\smallskip
\noindent 11734: K = L union B union D union heads of forced-N components;\par
\noindent 11735: K = L union D union heads and seconds of forced-N components.\par
\smallskip}

The preceding independent mode/value choices prove these descriptions are exact. A target literal x is universally in U iff x belongs to K. A nonlinear target hsZ is universally in U iff every variable in Z belongs to O and either(h, s) belongs to E*, or h belongs to K and s belongs to O. Indeed, a free head-component can make \(h=3\), which requires \(s=2\) in precisely the same component; its N mode then gives h in U, \(s=1\). All other universally valid first pairs must have head in U and second 1. A free second used as a target head is always in U only for 11734, as explicitly included in its K formula. No arbitrary normal head is forced 1.

The unrestricted exact absorption criterion \(w\leq P\) is therefore:

C: content(w) is contained in C; N: w is literal, or P has an original nonlinear word; I2: H intersects T, or the target head belongs to H and all its noninitial variables belong to T; IU: the target passes the preceding universal-U test.

The ideal\{1\} gives C, since an original word is 1 precisely when all of its variables are 1. The ideal\{3\} gives N, since only literals can have value 3. All originals can have value 2 exactly when heads are 2 and all noninitial variables 1, possible iff H, T are disjoint. This proves I2. The mode system gives IU. These four proper principal-ideal conditions are sufficient by considering the actual value of P under any assignment.

For necessity, failed C is separated by assigning a missing variable 4 and all others 1. Failed N is separated by assigning all variables 3. For failed I2, set H variables 2 and all others 1; this separates a missing head. A missing noninitial target variable can be set 2 if in H, and 0 otherwise, leaving all originals2 and forcing the target outside\{2\}. For failed IU, choose an allowed component mode or remaining actual value violating the failed target condition. All originals stay in U while the target leaves U. This establishes necessity for arbitrary alphabets, word lengths, and numbers of original summands.

\par\smallskip\noindent\textit{Fixed laws and complete component edges.}\quad 
Both sources satisfy these five identities:

{\small\par\smallskip
\CaseLine{\(xyzt = xytz\)}
\CaseLine{\(xyzz = xyz\)}
\CaseLine{\(xy \leq  xyz\)}
\CaseLine{\(xz \leq  xy+uy+uz\)}
\CaseLine{\(xyzt \leq  xyz+xyt\)}
\smallskip}

The first two normalize a word by retaining its first two positions and sorting/deleting repetitions only in its subsequent tail. Literals and quadratics stay unchanged. The word formula proves prefix extraction; each of the four exact conditions proves the rectangle and tail-merge laws. In particular the rectangle's original edges connect its target head and second in a single component, and every target tail of the merge is already universally 1 when required. These arguments are unrestricted.

Extract every original quadratic prefix. Repeated rectangle laws derive all E* edges:along a path of heads h0,..., hn sharing successive seconds b1,..., bn, an edge(hn, s) and the edges(h(n-1), bn),(hn, bn) give(h(n-1), s). Continue backward. Every step uses nonempty singleton substitutions. Let P+ consist of the original words and these derived edges; all are below P. This completes components of any size without an upper bound.

\par\smallskip\noindent\textit{Uniform short supports.}\quad 
For a head h in a forced-N component, choose an actually forcing head r in that component and a second b there. Include derived edges(h, b),(r, b). If r is a literal or occurs noninitially, one original witnessing this role forces r normal. Otherwise include an original edge(r, s) whose second lies in the blocking set(D for 11734, H union D for 11735), and one original witnessing that membership. Thus at most 4 words from P+ already force h's component normal.

A necessarily 1 variable needs at most 5 words:one original D-occurrence word, or an edge to a normal head together with its at-most 4-word normal support. A necessarily-U variable similarly needs at most 5:one original literal or D occurrence, a normal-head support, a normal-second support, or, for 11734 only, one original second-position occurrence. This last one-word branch is unavailable in 11735.

A universally-U target first pair therefore needs either one completed edge or at most 10 words for its K/O witnesses. A target of length 3 needs at most 5 more words for its last variable. So IU for any target of length at most 3 has at most 15 supporting words. Add at most 3 original content witnesses, one nonlinear original when needed, and at most 3 original I2 witnesses(two for an H/T intersection, or a target-head witness and its at-most 2 noninitial witnesses). This uses at most 22 words of P+ and already absorbs the short target. Adding witnesses restricts ideal-assignment sets and cannot invalidate an absorption.

\par\smallskip\noindent\textit{Finite family and unrestricted completeness.}\quad 
Normalize the target and the at most 22 supporting upper words. Among these at most 23 words, reserve every first/second-position variable, at most 46 variables. Group every other variable by its nonempty tail incidence pattern among the 23 words, at most \(2^{23}\)-1 groups. Substitute one distinct variable per group and retain all reserved variables separately. The resulting absorption is still valid by substitution. It has at most 8388653 variables, length at most 8388655, and at most 22 upper words.

The inverse substitution sends a group-variable to the nonempty product of its members and each reserved variable to its original singleton. It exactly restores first and second positions and the original tail sets. Literal applications of the first two fixed laws recover all original normal words. Short words contain only reserved variables; long words retain nonempty tails; any protected variable also occurring in a tail is preserved. Thus the compression is equationally reversible.

For each source separately, take its five fixed laws and all its valid absorptions over a fixed 8388653-variable alphabet, word length at most 8388655, and at most 22 upper words. The exact criterion effectively specifies this finite family; no complete enumeration or small basis is claimed. Inverse compression derives all valid short targets from P+, hence from P by the rectangle derivations. For a longer normal target hsZ, each triple hsz with z in Z is valid by the exact criterion. Derive these triples and merge their nonempty tails using the last law. This derives every valid absorption without original bounds. Distributivity and the two absorption directions then derive every valid identity.

Hence 11734 and 11735 are finitely based. Opposites have reversed bases; other catalogue transfers still require independent complete mutual HSP.

\subsection{A 167-identity basis from a common final endpoint}\label{case:263}
In type 11736 let \((e,f,a,c,z)=(3,1,2,0,4)\). The two loops e, f, the bridge a and the multiplicative zero z form the ordinary two-vertex forward path semigroup with zero. The only nonzero product involving c is \(cf=c\). Addition has f, \(a<c<z\) and \(e<z\), with f and a incomparable and e incomparable with f, a, c. Thus a word is e or f just on that constant loop, is a for e* a f*, is c for c f*, and otherwise gives z. These facts include literals. Reversal gives the other catalogue type 14142.

\par\smallskip\noindent\textit{Exact criterion.}\quad 
For a nonempty upper polynomial P take its content union C and the set D of letters appearing noninitially in some original word. Build G by identifying the right endpoint of every original final letter with T and \(R_{x}\) with \(L_{y}\) at each original adjacency xy. Use one directed edge \(L_{x}\) \ensuremath{\to} \(R_{x}\) per variable. Keep every original initial vertex \(L_{h}\) in a set I. At this stage the vertices in I are separate unless the preceding endpoint identifications equate them. Every vertex is on some directed walk from I to T.

A word w is absorbed exactly when its content is in C, its noninitial letters are in D, and the following tests hold in G.

\textbullet{} The vertex of its initial endpoint can reach some member of I.

\textbullet{} Its final endpoint and T are in one strong component.

\textbullet{} Each adjacent join \(R_{x}\), \(L_{y}\) in w is in one strong component.

Reachability here includes equality of vertices. First use ordinary bridge values \(e=00\), \(a=01\), \(f=11\). For all upper words to be below c, with no marked input, their edges must compose and end at 1. Their initial bits are unrestricted. These assignments are precisely the monotone bit labellings of G with \(T=1\). Two vertices get the same bit in every such labelling exactly when they are in the same strong component. If they are distinct in the quotient, the upward closure of a vertex that cannot reach the other separates them and contains T. The final and adjacency tests follow. A missing variable can be separated by assigning it z and giving all original letters e. For x outside D give x the mark c and every other letter f. Original words are bounded by c, whereas any noninitial x makes the target z. Conversely project each c-valued input to a in an assignment bounding all original words by c. Its occurrences in P are all initial, so it cannot appear noninitially in a target satisfying the D test. The ordinary quotient walk ending at 1 evaluates to f, a or the initial marked bridge c.

At bound a all upper words must give the bridge. Identify every member of I at H in G. Ordinary bit labellings now have \(H=0\) and \(T=1\). The final and internal equalities already hold in this coarser graph. The initial endpoint of w must be strongly connected to H. Indeed, if this endpoint is outside the component of H, put 1 on its upward closure, which contains T and omits H. If H and T are strongly connected, all vertices lie in their component and the bridge bound cannot occur. In the graph with I identified, H reaches all vertices, and a vertex can return to H exactly when it could reach some member of I in G: truncate its return just before its first encounter with I. This is precisely the stated initial test. The bounds at e and f are covered by content and the bound at z is automatic.

\par\smallskip\noindent\textit{Finite presentation.}\quad 
Alongside the ai-semiring axioms use three fixed laws

{\small\par\smallskip
\CaseLine{\(xyx = xxyy,\qquad x \leq  xyy,\qquad xy \leq  xzzy.\)}
\smallskip}

Include the following optional-context schemes.

{\small\par\smallskip
\CaseLine{\(M\qquad A x D \leq  A x B + C x D\)}
\CaseLine{\(R\qquad A x t F \leq  A x y B + C z y D + E z t F\)}
\CaseLine{\(T\qquad D z y B \leq  A x y B + C x + D z\)}
\CaseLine{\(E\qquad A x \leq  A x y B + C z y D + E z\)}
\CaseLine{\(F\qquad y B \leq  A x y B + C x z D + z E\)}
\CaseLine{\(L\qquad y z B \leq  x y z B + z C\)}
\smallskip}

Each lower-case symbol is retained. In each scheme any subset of the capital symbols is omitted throughout that scheme before making an identity. Subsequent substitutions are ordinary nonempty terms. This gives respectively 16,\allowbreak{}64,\allowbreak{}16,\allowbreak{}32,\allowbreak{}32,\allowbreak{}4 instances, totalling 167 laws with the three fixed ones.

For M, R, T, E the target is an actual walk following identifications of adjacent endpoints or finals. Its initial endpoint is that of an original upper word. For F the two upper words containing x identify the starts of y and z, and z is an original head. For L the right word xyzB gives a path from \(L_{y}\) to the original initial vertex \(L_{z}\). The target suffix therefore meets the initial test. In each scheme target tail letters occur in an upper tail. For the fixed equation each side gives the common loop when x and y equal that loop and otherwise gives z. A squared block in an upper nonzero product giving a or c is a loop; deleting it internally or finally in the two fixed inequalities preserves that bound. The loop and top bounds follow immediately. This verifies soundness at arbitrary substitutions.

\par\smallskip\noindent\textit{Completeness.}\quad 
Call a word available when it is derivably below P. Keep a witnessed pair xy for every actual adjacent occurrence in an available word, a final witness for every available final letter, and a head witness for every available initial letter. Retain the original D.

Use R, T and E to complete the bipartite undirected incidence graph whose rows are \(R_{x}\), columns are \(L_{y}\) together with the terminal column T, and edges are the original adjacent pairs and finals. R fills a rectangle across a shared ordinary column, T copies an ordinary column between two rows incident with T, and E copies the terminal column between two rows sharing an ordinary column. In every case split the available witnessing words at their indicated occurrences to supply exactly the displayed premises. Repeated rectangle completion makes each connected component a complete bipartite graph: paths of odd length reduce by successively replacing three edges with one. The converse holds because filling adds only incidences in one component. Consequently every pair with \(R_{x}=L_{y}\) in G and every final with \(R_{x}=T\) has an available witness. An \(L_{y}\) with y outside D is isolated here and cannot equal a row. A row-column equality is always witnessed by a genuine noninitial occurrence of its column letter.

First mark all original heads. F marks every column in a component that contains a marked column. If the two columns are distinct, choose a common row after rectangle completion; the pair witnesses and the previous head witness supply its three premises. For propagation across a directed edge \(L_{y}\) \ensuremath{\to} \(R_{y}\), suppose \(R_{y}\) contains a marked column \(L_{z}\). The pair yz has an available witness A y z B. If A is empty then y is already a head. Otherwise use L with \(x=A\) and a head witness z C to obtain y z B. Now y is a head; propagate it to all columns equal to \(L_{y}\) by F. Induction on a path length proves that every letter whose left vertex reaches I has a head witness. The induction proceeds backwards along a path to an original head. Every intermediate vertex has the left endpoint of the next directed edge, so there is always a column to mark. These operations add no incidences outside the original bipartite components, since F and L take suffixes of existing words. Their new heads satisfy the initial test just proved.

For a target passing the criterion, take its actual letters in order. At any adjacent gap with different endpoints in one strong component choose a nonempty directed path U from the first to the second and a nonempty return V. Insert the label word UVU in that gap. Do the same at the final gap ending at T. There are no initial insertions. Every inserted letter is in D: the left endpoint of a head-only letter is an isolated column with no incoming edge and can lie on no nonempty directed cycle. The lifted word W has an available head, a witnessed final and a witnessed instance of each consecutive pair. Construct W by M. Start with a head witness and splice its successive letters to the suffixes of the pair witnesses. At the last letter splice to the final witness to finish the word. The same final splice works for a literal. These are derivations of W \ensuremath{\leq} P.

To derive w \ensuremath{\leq} W, insert U squared and then V squared at each gap using the internal fixed inequality or the final one, as appropriate. Apply \(xxyy=xyx\) to turn these two squares into UVU. Keeping the rest of the word as a multiplicative context gives the desired word after all the gaps are treated. Equal endpoints need no insertion. This proves w \ensuremath{\leq} P with no restriction on the input alphabet, lengths, number of right summands or number of graph vertices. Deriving the two lists of summand absorptions after distribution proves all identities. The reversed presentation gives the basis of 14142.

\subsection{Complete finite bases for 11771 and 11772}\label{case:264}
Signature (+, multiplication), with no constants. Relative to the additively idempotent semiring axioms, use the six common laws

{\small\par\smallskip
\CaseLine{\(xy=yx; xyz=xyzz; xyzt\leq xyz;\)}
\CaseLine{\(xxy\leq xy+zzz; xy\leq xxy; xy\leq xx+yz.\)}
\smallskip}

For 11771 add \(x\leq xxx\) and \(xxy\leq x+xy\), giving eight laws. For 11772 add \(x\leq xx\), giving seven laws. Write \(u\leq P\) for \(u+P=P\).

\par\smallskip\noindent\textit{Exact semantics.}\quad 
Both have multiplication 41114/\allowbreak{}11111/\allowbreak{}11111/\allowbreak{}11121/\allowbreak{}41114. Addition 11771 is 00004/\allowbreak{}01114/\allowbreak{}01214/\allowbreak{}01134/\allowbreak{}44444:2,\allowbreak{}3 are incomparable below \(1<0<4\). Addition 11772 is 00004/\allowbreak{}01114/\allowbreak{}01224/ 01234/\allowbreak{}44444, the \(chain 3<2<1<0<4\). Nonlinear words are 4 exactly when every variable is assigned 0 or 4. Otherwise a quadratic is 2 when both occurrences are 3; all other nonlinear words are 1.

For original P let L be its literals, E its unordered quadratic edges, and V their endpoints. Let C(w) require an original NONLINEAR word whose content is contained in c(w). A contained literal alone does not satisfy C. Put \(G=true\) when an original long word exists; for 11771 also allow L intersect V to be nonempty. The exact criteria are:

\textbullet{} Long target:C(w) and G.

\textbullet{} Quadratic:C(w) and (G or c(w) contained in V).

\textbullet{} Literal 11771:original literal, or C(w) and G.

\textbullet{} Literal 11772:original literal, or C(w).

The distinction from earlier content algebras is essential:assigning target variables 0 makes nonlinear targets4, but target literals only 0. Thus if C fails, set all target variables 0 and outsiders 1. Every right nonlinear word is 1, while literals are at most 0. A nonlinear target 4 is separated. For a missing literal with C false, use target 4, outsiders 1; every right term is 1. If G fails, set V to 3 and outsiders 2. In 11771 L is disjoint from V, so all right values are 2. In 11772 all right values are at most 2. A long target is 1; a quadratic with an endpoint outside V is 1. For 11771, a missing literal satisfying C has an original square and belongs to V, so its value 3 is not below the right value 2. These assignments separate every failed case without size restrictions.

For sufficiency an original long word supplies at least 1. In 11771 an original incident literal x and edge xy also do:\(xy=2\) forces \(x=3\), and \(3+2=1\); otherwise xy is 1 or 4. The contained nonlinear witness supplies 4 whenever every target variable is 0 or 4. All remaining nonlinear target values are at most 1. In 11771 the only literal value not below 1 is 0 or 4, also covered by the pure nonlinear witness. In 11772 every literal is below all its pure powers. Without G, if the right side does not dominate 1, all right quadratics must be 2, forcing every V variable 3. A target quadratic on V is then 2. This proves the criteria.

\par\smallskip\noindent\textit{Complete derivations for arbitrary sizes.}\quad 
Commutativity and \(xyz=xyzz\) normalize long words by content. Commute a chosen variable to the end, duplicate it using two nonempty preceding blocks, and remove surplus occurrences with the reverse law while keeping length at least three. Canonical long forms are sorted content when there are at least three distinct variables, xxy for two, and xxx for one. Every long g therefore satisfies \(g^{3}=g\).

For G retain an original long g. In 11771, if necessary, derive \(g=xxy\) from an original literal x and original edge xy with \(xxy\leq x+xy\). Unlike the earlier constant-cube constructions, this preserves their joint content. Choose an original NONLINEAR contained word u. Retain u if long; if \(u=ab\) is quadratic, derive aab using \(xxy\leq xy+zzz\) with \(z=g\). Append missing target content by \(xyzt\leq xyz\) and normalize. Recover a quadratic by \(xy\leq xxy\), and a11771 literal by \(x\leq xxx\).

In 11772, a missing literal has a contained original pure nonlinear power. If long, normalize to xxx and derive xx using \(xy\leq xxy\) with \(x=y\). Then \(x\leq xx\) recovers the literal. No global long witness is needed here.

Without G, a new quadratic has an original contained nonlinear witness, which can only be an endpoint loop; an original target edge is immediate. Choose an original edge incident at the other endpoint and use \(xy\leq xx+yz\), commuting endpoints if needed. There is no literal-only endpoint guard. All substitutions are nonempty and every used summand is original or explicitly derived. Each law satisfies the criterion and is checked independently on every assignment. Distributive expansion and summandwise absorption in both directions prove unrestricted completeness.

In both algebras \(xy\leq x+zzz\) fails at \(x=y=0\), \(z=1\). Also \(xxx\leq x+xy\) fails at \(x=0\), \(y=1\). In 11771 \(x\leq xx\) fails at \(x=3\). In 11772 \(xxy\leq x+xy\) fails at \(x=y=3\). These prevent dropping the nonlinear witness or confusing the two literal-recovery branches. Transfers need full HSP both ways.

\subsection{A finite basis for 11797 from original head thresholds}\label{case:265}
Signature (+, multiplication), without constants, relative to the additively idempotent semiring axioms. The specified finite family is not an enumerated small basis. Completeness is unrestricted.

Addition rows:00004,\allowbreak{}01114,\allowbreak{}01224,\allowbreak{}01234,\allowbreak{}44444; order \(3<2<1<0<4\). Multiplication rows:41114,\allowbreak{}11111,\allowbreak{}23332,\allowbreak{}33333,\allowbreak{}41114. Literals retain value. For nonlinear hT, heads 1 and 3 retain value; head 2 gives 2 if all tails lie in\{0,\allowbreak{}4\}, otherwise3; head 0 or 4 gives 4 if all tails lie in\{0,\allowbreak{}4\}, otherwise1. Induction on multiplication proves this formula. Hence \(xyz=xzy\) and \(xyy=xy\) hold; normalize to head plus sorted tail set. Heads and head-tail repetitions are retained separately. Adjacent swaps and duplicate deletions substitute only nonempty words.

\par\smallskip\noindent\textit{Exact criterion using original words.}\quad 
A literal x is absorbed iff original literal x occurs. For nonlinear hT, \(C=\{h\}\) union T, and H the set of original heads, the conditions are: (H) h belongs to H; (NC) some original NONLINEAR word is contained in C; (R) h belongs to T, or literal h occurs, or T intersects H, or an original nonlinear h-word has tail contained in T.

Missing literal:assign \(x=2\), others 3. All right nonlinear x-words give 3, as do other heads; without the literal x target 2 is separated. Failure of H:assign \(h=1\), others 3. Target 1, right 3 separate. Failure of NC:assign all target variables 0, outsiders 3. Target 4; original contained words are literal 0, and nonlinear words with outsiders give 1 or 3. The right sum is at most \(0<4\). Failure of R, assuming H and NC:the head is outside T and no original head lies in T. Assign \(h=2\), \(T=4\), others 3. Target 2; literal h is absent. Every nonlinear h-word has a low tail outside T and gives 3, while all other original heads are 3. Thus the right sum 3 separates target 2.

Sufficiency:3 is least. A target 1 has head 1, or head 0/\allowbreak{}4 with a low tail. H supplies a same-head word valued at least 1. Target 4 has all target variables in\{0,\allowbreak{}4\}; NC supplies a nonlinear contained word valued 4. Target 2 has head 2 and all tails 0/\allowbreak{}4, so the head is not repeated. R gives literal \(h=2\), an original word headed in T valued 1 or 4, or a nonlinear h-word with tail in T valued 2. All absorb 2. Literal targets are original summands. This proves the criterion for arbitrary words and polynomials.

\par\smallskip\noindent\textit{Three witnesses and reversible finite compression.}\quad 
Choose one original H witness a, one original NONLINEAR contained witness b, and if the head is not repeated one original R witness c. The latter is literal h, a word whose head lies in T, or a nonlinear h-word with tail in T. After removing duplicates Q has at most 3 words and still satisfies H, NC, R. Hence the target is absorbed by Q, all of whose words are original.

There are at most 4 normal words in this absorption. Reserve their at most 4 heads. Partition remaining variables by nonempty incidence among the at most 4 tails, at most 15 classes. Identify each class to a fresh variable, keep reserved heads separate, and normalize. Validity survives substitution; there are at most 19 variables, word length at most 20, and at most 3 right words.

Invert by substituting the nonempty product of the original class members. Identical tail incidence restores all tail sets exactly. Reserved heads and head-tail repetitions are unchanged; nonempty tails stay nonempty. Sorting and duplicate deletion recover every original normal word.

Let Sigma consist of \(xyz=xzy\), \(xyy=xy\) and every valid normal absorption on a fixed 19-variable alphabet, length at most 20, at most 3 nonempty right words. The exact H, NC, R criterion specifies this finite family. Every valid absorption selects Q, compresses into Sigma, and expands back. Distributivity and summandwise absorption in both directions then yield every valid identity. Sigma is a finite relative basis for 11797. Reversal yields the opposite basis; catalogue transfer requires complete HSP both ways. Finite replay supplements this unrestricted argument.

\subsection{Six complete identities for catalogue algebra 11832}\label{case:266}
The language is addition and multiplication without constants. All blocks substituted below are nonempty. The proof is unrestricted in variable count, word length and polynomial size.

\par\smallskip\noindent\textit{Word values.}\quad 
Addition rows are 00004 / 01114 / 01214 / 01134 / 44444. Multiplication rows are 41114 / 41114 / 41234 / 41324 / 41114. The additive order is 2,\allowbreak{}3 \ensuremath{<} 1 \ensuremath{<} 0 \ensuremath{<} 4, with 2 and 3 incomparable. The subset \{2,\allowbreak{}3\} is a multiplicative group of order two, with identity 2. For a nonlinear word, a last position valued 0 or 4 gives 4. Otherwise any position outside \{2,\allowbreak{}3\} gives 1. If every position is in \{2,\allowbreak{}3\}, its value is 2 or 3 according to the parity of its occurrences of 3. These statements follow by induction from the multiplication table. Literals have their assigned value and must be kept distinct from nonlinear words.

\par\smallskip\noindent\textit{Exact absorption criterion.}\quad 
For each word u let pi(u) be its multiplicity vector modulo two over the finite variable set in the identity. A target word w is absorbed by an original polynomial P exactly when:

1. Every variable of w occurs in at least one original word of P.

2. The last variable of w is the last variable of an original P-word; if w is nonlinear, this witness must also be nonlinear.

3. pi(w) belongs to the affine hull over F2 of the vectors pi(u), u in P; equivalently pi(w) is the sum of an odd number of these vectors.

For necessity of (1), assign a missing variable 1, and every other variable 2. The target is 1 while P is 2. For necessity of (2), assign the target's last variable 4 if it never appears as an original tail. Original interior occurrences yield 1 while the target is 4. If the only tail witness is a literal and w is nonlinear, instead assign this variable 0: P is at most 0 and w is 4. This is why literal witnesses must not be used for nonlinear targets. For (3), restrict all assignments to \{2,\allowbreak{}3\}. A separating linear functional on the affine span makes every original word have the same group value, opposite to w's group value. This elementary affine separation works over the full variable set.

For sufficiency, a target value 4 is covered by its original tail witness. A target value 0 is literal and is covered by its tail witness (which is 0 or 4). If its value is 1, a target variable outside \{2,\allowbreak{}3\} occurs in P by (1), and its original occurrence makes some right word at least

1. Finally suppose w is 2 or 3. If any P-word is outside \{2,\allowbreak{}3\}, its value dominates w. If all P-words are group-valued, all their variables are in \{2,\allowbreak{}3\}; either their values differ, giving join 1, or all are the same group value. In the latter case (3) makes w have that same value.

\par\smallskip\noindent\textit{The six identities.}\quad 
In addition to the additively idempotent semiring axioms, take:

1. xyz = yxz.

2. xxxy = xy.

3. xyyy = xy.

4. y \ensuremath{\leq} xxy.

5. xxyz \ensuremath{\leq} x + yz.

6. xyzt \ensuremath{\leq} x + y + zt.

Here \(p\leq q\) means \(p+q=q\). All six are valid by the word formula. Independent verification exhaustively checks their finitely many source assignments.

The first three laws give a normal form preserving content, parity and last variable for nonlinear words. Prefix variables can be permuted arbitrarily by law 1, keeping a nonempty suffix. For a variable other than the last variable, remove pairs from three prefix occurrences by law 2 until it occurs once or twice. For the last variable, use law 3 to remove pairs immediately before its final occurrence whenever there is another prefix variable. Pure powers reduce to x squared or x cubed. Literals remain separate. Thus any two nonlinear words with the same content, last variable and parity are equated by laws 1-3.

\par\smallskip\noindent\textit{Unrestricted completeness.}\quad 
Let a nonlinear w satisfy the criterion. Choose an original nonlinear word v ending in its last variable t. By (3), select an odd indexed subset S of original P-words whose parity sum is pi(w). Start from v. For each u in S, apply law 6 with \(x=u\), \(y=v\) and zt equal to the current nonlinear word, split into two nonempty blocks. This adjoins u v times the current word below P. Each step retains the originals u and v. The final word G ends in t and has parity sum(u in S) pi(u) + (|S|+1)pi(v) = pi(w).

To ensure sufficient target-variable occurrences without changing parity, for every target variable choose an original P-word u containing it and apply law 5 with \(x=u\) and yz the current nonlinear word. This adds two copies of u to its prefix. Every target variable now appears at least twice in that prefix (in addition to any final occurrence).

Any variable outside the target content has even multiplicity in G, since G and w have the same parity. It cannot be the last variable. Use law 1 to put two of its prefix occurrences first, then law 4 with x that variable and y the entire remaining word to DELETE the pair in the absorption direction. The remaining word is nonempty and stays nonlinear because target variables were supplied in pairs and the final variable remains. Repeat until no extraneous variable remains. No cancellation identity or balancing of unequal multiplicities is assumed. Normalizing intermediate words can compress the protected last variable to one occurrence while outside variables remain. Before deleting the final outside pair for a singleton odd target, apply \(xy=xyyy\) (law 3 backwards) with the still nonempty outside prefix. This restores the final cube and ensures the remaining target stays nonlinear.

The remaining nonlinear word has precisely w's content, last variable and parity. Laws 1-3 therefore equate it to w. All derived words were adjoined below the original sum, proving \(w\leq P\).

For a literal target x already occurring as a literal summand of P, there is nothing to derive. Otherwise condition (2) gives an original nonlinear tail witness ending in x. Conditions (1)-(3) also hold for the nonlinear target xxx. Apply the construction just proved, then use law 4 with \(x=y\) equal to this variable to obtain \(x\leq xxx\leq P\).

Every term expands into a nonempty finite polynomial, and a polynomial identity is equivalent to the absorptions of its individual summands in both directions. The construction proves all valid absorptions with no bound on variables, lengths, or summands. Hence the six displayed identities form a complete finite basis for 11832.

\subsection{Six complete identities for 12473 from contained head saturation}\label{case:267}
Signature (+, multiplication), no constants, relative to the additively idempotent semiring axioms. A finite relative basis is:

1. \(xyz=xzy\).

2. \(xyy=xy\).

3. \(x\leq xx\).

4. \(xyz\leq xy+xz\).

5. \(xy\leq xxy\).

6. \(xxy\leq xx+y\). All substitutions used below are nonempty. Completeness has no bound on variables, lengths, or the number of original right summands.

Addition:00004,\allowbreak{}01114,\allowbreak{}01224,\allowbreak{}01234,\allowbreak{}44444; order \(3<2<1<0<4\). Multiplication:44034,\allowbreak{}11131,\allowbreak{}11231,\allowbreak{}33333,\allowbreak{}44434. Any occurrence of 3 in a nonlinear word gives 3. With no 3, head 1 gives 1, head 4 gives 4, head 2 gives 2 if every tail is 2 and 1 otherwise, and head 0 gives 0 if every tail is 2 and 4 otherwise. Literals retain value. Induction on the table proves these formulas and the two word laws. Normalize each word to its head and sorted tail set, retaining a repeated head in the tail. All actual swaps and duplicate deletions have nonempty prefixes.

\par\smallskip\noindent\textit{Exact criterion.}\quad 
For target w with head h and content C, let U be the original right words with entire content in C, and \(U_{h}\) those headed h. Let \(D_{h}\) be the union of the tails of \(U_{h}\). A literal h is absorbed iff \(U_{h}\) is nonempty, that is, an original pure power of h occurs. For nonlinear \(w=hT\) the criterion is: (K) \(U_{h}\) is nonempty; (C) the contents of U cover all of C; (S) h is in \(D_{h}\), or T is contained in \(D_{h}\). Thus an original repeated-head word permits general contained coverage; without that word every target tail must occur in an original contained same-head tail. One merely contained witness is not full content coverage.

Necessity: if K fails, set \(h=4\), other target \(variables=2\), \(outsiders=3\). The target is 4. Every noncontained right word is 3; contained right words have head 2 and value 1 or 2, so cannot reach 4. This also separates a missing literal pure-power witness. If C fails at x, set \(x=1\), other target \(variables=2\), \(outsiders=3\). The target is 1, whereas contained right words avoid x and are 2, and noncontained words are 3. If S fails, set \(h=0\), \(D_{h}=2\), other target \(variables=1\), \(outsiders=3\). The assignment is consistent because h is not in \(D_{h}\). Some target tail is not in \(D_{h}\), so the target is 4. Each contained h-word has only 2 in its tail and is 0. Other contained heads are 1 or 2 and give at most 1; outsiders force 3. The right sum is at most 0, below 4. These separators are valid for arbitrary polynomials.

Sufficiency: target 3 is least. A target 4 with head 4 is matched by any \(U_{h}\) witness. If the head is 0, some target tail differs from 2. When h is in \(D_{h}\), a contained repeated-head word has tail 0 and is 4. Otherwise S makes that differing target tail occur in a \(U_{h}\) word, also giving 4. For target 0, the head is 0 and all target tails are 2; a \(U_{h}\) witness gives 0 or 4. For target 2, every target variable is 2, and a \(U_{h}\) witness gives 2. For target 1, head 1 is handled by \(U_{h}\). The remaining possibility is head 2 and a tail value in \{0,\allowbreak{}1,\allowbreak{}4\}; C supplies a contained word containing that variable. It has no 3 and has value at least 1, whether the chosen variable is a head or a tail. Literal absorption follows from \(h\leq h^{n}\), valid for all positive n by the same formula. This proves exactness without a finite-variable assumption.

Each displayed law is valid: the first two follow from word values, law 3 from the pure-power literal condition, law 4 from direct same-head tail coverage, law 5 from repeated-head saturation in its right word, and law 6 from the repeated-head witness xx and complete coverage by y.

\par\smallskip\noindent\textit{Complete finite derivation.}\quad 
Normalize original words. A literal target is already present or follows from its original pure power, normalized to hh, by law 3. For nonlinear hT with no contained repeated-head witness, S gives T in \(D_{h}\). Choose an original nonlinear h-word; all its tails belong to T. As long as a target tail is absent, choose an original contained h-word containing it and merge their tails by law 4 with nonempty y, z. Each step increases the assembled tail set and eventually yields hT.

If a contained repeated-head witness exists, start with that word g. It is idempotent under multiplication by the two word laws: its tail already contains its head, so normal(g squared)=g. Until its content is C, choose an original contained u containing a missing target variable. Law 6 with \(x=g\), \(y=u\) gives normal(g squared \(u)\leq normal(g\) squared)+u=g+u. The resulting word keeps head h, repeats h, and strictly increases content. All selected words are contained in C, so the process terminates with head h and tail C. If the target repeats h, this is the target. Otherwise apply law 5 with \(x=h\) and y the nonempty target tail to remove the surplus head occurrence. Every step is an actual substitution of one of the six laws, with word-law normalization; no additional infinite rule family is used. Every original absorption is thus derivable. Distributivity and summandwise absorption in both directions derive all valid identities. This proves the unrestricted six-law basis. Opposite catalogue transfers require independent complete HSP in both directions.

Coverage is essential: \(xy\leq xx\) fails at \(x=2\), \(y=1\). Saturation is essential: \(xy\leq x+y\) fails at \(x=0\), \(y=1\). An original repeated-head word cannot be replaced by a literal h in the saturated construction.

\subsection{Seven complete identities for 13098 with a contained head and global tail}\label{case:268}
Signature (+, multiplication), no constants, relative to the additively idempotent semiring axioms. All proofs are unrestricted in variable count, word length and polynomial size.

Addition:00004 /01114 /01224 /01234 /44444. Multiplication:44234 /11231 /22232 /22232 /44234. Addition is the \(chain 3<2<1<0<4\). Put K=\{0,\allowbreak{}1,\allowbreak{}4\}. A literal retains its value. A nonlinear word has value 3 when its last occurrence is 3. If its last occurrence is not 3 but some occurrence is 2 or 3, its value is 2. With all factors in K it is 1 when its head is 1 and 4 otherwise. This follows by induction: right 3 resets every value to 3, right 2 to 2; a current2 or 3 followed by K gives 2, while K-only nonlinear products retain1 or 4 according to their original head.

The complete relative basis is:

{\small\par\smallskip
\CaseLine{\((W) xyzt=xzyt\)}
\CaseLine{\((H) xxy=xy\)}
\CaseLine{\((T) xyy=xy\)}
\CaseLine{\((I) xyz\leq xz\)}
\CaseLine{\((L) xyt\leq xy+t\)}
\CaseLine{\((G) xyt\leq xy+zt\)}
\CaseLine{\((P) x\leq xx\)}
\smallskip}

For a target w let U be the original right words contained in c(w). A literal h is absorbed exactly when U is nonempty. A nonlinear target is absorbed exactly when: H: U contains an original NONLINEAR word with the same head as w; T: some original right word ends in the final variable of w. The T witness is not required to be contained, may be literal, and may contain arbitrarily many variables outside the target.

If H fails, assign the target head 0, other target variables 1, outsiders 3. The target is 4. Any right word with an outsider is 2 or 3; contained nonlinear words have a different head and give 1; literals give 0 or 1. Their join is below 4. If T fails, assign the target final variable 1 and every other variable 3. The target is 1 for singleton content and 2 otherwise; every right word ends in 3 and gives 3. A literal without a pure witness is separated by that variable 0 and outsiders 3: the target is 0 and all right words have value 2 or 3. These exhaust failed criteria.

Conversely target 3 is least. For target 2 its final value is not 3. Any same-tail original witness, even with outside variables, therefore has value at least 2. For target 1 or 4 every target variable is K-valued, and the contained nonlinear same-head witness gives exactly that value. A pure witness absorbs its literal because powers reduce to a square and the square bounds the literal. Thus the criterion is exact.

The word formula verifies W, H, T. They preserve content, endpoints and nonlinearity; H or T applied to a pure cube reduces it to a square. Use W to group interior duplicates. Extra head copies are deleted by H with a nonempty suffix, extra final copies by T with a nonempty prefix. Other internal duplicates are deleted by H with the following nonempty suffix and the preceding context. W then sorts the interior. A nonlinear normal word consists of its head, sorted content excluding endpoints, and its final variable. Coincident endpoints are retained twice. A literal remains distinct from its square.

I follows from the criterion: xz is a contained nonlinear same-head word and a same-tail witness for xyz. L and G have xy as the contained nonlinear head witness and t or zt as the global tail witness. P follows from the table. All seven laws therefore hold, with nonempty substitutions.

To derive any valid nonlinear absorption \(w\leq P\), choose an original contained nonlinear same-head word v and an original tail witness q. Normalize all words. Let t be the target tail. If v already ends in t, keep v. Otherwise write \(v=XY\) with X, Y nonempty. If \(q=t\) use L to derive \(vt\leq v+t\). If \(q=Zt\) with Z nonempty, use G to derive \(vt\leq v+Zt\). The resulting g has the desired endpoints and content contained in c(w), despite any outside variables in q: those occur only on the right side of G and are not copied into g.

If any target variables are missing from g, write \(g=At\) with A nonempty and substitute their nonempty list for y in I. This derives \(Ayt\leq At\). Normalize to obtain w. If no variables are missing, g already normalizes to w. A literal target uses its original pure witness, normalized to the literal or square, and P. Every substitution is nonempty and all witnesses are original. There is no bound on their lengths or the number of target variables. Hence every valid word absorption follows from the seven laws. Distributivity and summandwise absorption both ways give completeness for all identities. Reversal yields the opposite basis; catalogue transfer requires complete independent HSP inclusions in both directions.

Boundary failures:\(xx\leq x\) at \(x=0\); \(xy\leq x\) at \(x=0\), \(y=1\); \(xy\leq xz\) at \(x=1\), \(y=1\), \(z=3\). A literal cannot replace the nonlinear H witness, and a missing T witness cannot be replaced by arbitrary contained content.

\subsection{Finite basis for 11954 by contained seeds and simple-tail recovery}\label{case:269}
Signature (+, multiplication), no constants, relative to the additively idempotent semiring axioms. No bound on variables, lengths or summand counts is assumed in this proof. The finite recovery family below is specified by a validity test and bounds, not claimed to be a small list.

\par\smallskip\noindent\textit{Tables, word values and exact criterion.}\quad 
Addition rows: 00004 / 01114 / 01224 / 01234 / 44444. Multiplication rows: 41234 / 23332 / 23332 / 33333 / 41234. The additive chain is \(3<2<1<0<4\). Put K=\{0,\allowbreak{}4\}, L=\{1,\allowbreak{}2\}. For a nonlinear word: a 3-valued position or at least two L-valued positions gives 3. Exactly one L-valued position gives 1 precisely when it is the last position and valued 1, otherwise 2. All K positions give

4. Literals retain their values. Induction on length follows directly: 3 is absorbing; an L value followed by K becomes/stays 2; two L values give 3; a K prefix is 0 or 4 and takes the value of a final L factor.

For target w, let U be original right words with content contained in c(w). A literal target is absorbed exactly if that literal occurs in P. For nonlinear w the exact conditions are: (C) U contains a nonlinear word. (S) Each simple variable x in w has a U-word omitting x or containing x exactly once. If x is the last variable of w, the containing witness must also end in x. An omitting literal is allowed here and must stay literal.

Necessity of C: assign c(w) to 0 and outsiders to 3. The target is 4; without a contained nonlinear word the right side is at most 0. For a failed S assign \(x=2\), other target variables 0 and outsiders 3. The target is 2 and every right word is 3. For the additional simple-tail condition assign that tail 1, other target variables 0 and outsiders 3; the target is 1 and all right words are at most 2. A missing literal x is separated by \(x=1\) and all other variables 3: every other word is 3.

Sufficiency: target 3 is bottom. A target 4 has all its variables in K and the contained nonlinear seed supplies 4. A target 1 or 2 has a unique L-valued occurrence, belonging to a simple variable x. An omitting contained word has value 0 or 4 and dominates. Otherwise the single-x witness has value 1 or 2; when the target is 1 its witness ends in x and has value 1. These cases exhaust all assignments.

\par\smallskip\noindent\textit{Structural laws and word normalization.}\quad 
Use seven fixed laws: (T) \(xyz=yxz\); (R) \(xxx=xx\); (D) \(xy=xxy+xyy\); (E) \(xxyy=yyxx\); (L1) \(xxy\leq xy\); (L2) \(xyy\leq xy\); (A) \(xxyy\leq xx\). Here \(p\leq q\) means \(p+q=q\). The word formula proves all seven laws. The law \(xx\leq x\) is NOT valid \((x=0)\) and is never used.

T permutes all prefix positions while fixing the last position. Group repeated letters and use R to cap each multiplicity at two. If the last variable is simple, retain it. If it is repeated, choose the least repeated variable as the canonical last variable: arrange its square just before the final square and apply E to interchange those squares, then permute the prefix again. This process uses only nonempty blocks and contexts. Consequently words with identical capped multiplicities and identical simple tail (or both having repeated tails) are equivalent. A nonlinear pure power reduces to a square; literals remain distinct.

\par\smallskip\noindent\textit{Saturated guard and one-simple-variable recovery.}\quad 
Choose a contained nonlinear seed. Double each simple occurrence using L1 at the first position with its nonempty suffix, or L2 with its nonempty prefix and the remaining suffix as a context. Normalize. The result g has every variable repeated. For each missing target variable y, law A with \(x=g\) appends y squared, since g squared normalizes to g. Thus g is a full-content saturated guard below P. It has no simple tail.

For a simple target x choose the original U-witness given by S. Double every other simple variable if this witness is nonlinear, but leave a literal witness untouched. Write \(u_{x}\) for this contracted original witness and \(w_{x}\) for w with every variable except x doubled. Then

{\small\par\smallskip
\CaseLine{\(w_{x} \leq  g + u_{x}\)}
\smallskip}

is valid by C and S. The guard supplies C; \(u_{x}\) supplies the only remaining simple-variable condition, including its tail condition. If \(u_{x}\) omits x and is literal, that literal is still a valid omission witness: its K-value 0 or 4 dominates 1 and 2. No squaring of such a literal is allowed.

\par\smallskip\noindent\textit{Reversible compression and the finite family.}\quad 
For the three words \(w_{x}\), g, \(u_{x}\) reserve x and the last variable of each word, at most four variables. Every nonreserved variable is absent or repeated in every word; an omitting literal has its variable reserved. Group the others by their nonzero incidence patterns across the three words, at most seven classes. Identify each class to a fresh variable and normalize. This is a substitution into a valid absorption, hence remains valid. There are at most eleven variables, word length at most twenty-two and at most two right words.

For the inverse replace a class variable by the nonempty product of its members. It occurs twice or not at all, so every member returns twice. T, R regroup and cap occurrences, and E restores the canonical repeated tail if needed. Reserved simple tails and literal witnesses are unchanged. Thus every compressed word expands to its original word modulo the seven laws. Common multiplicity, not just incidence, is preserved.

Let Sigma consist of the seven laws plus ALL valid absorptions with one canonical left word and one or two canonical right words over eleven fixed variables, each word of length at most twenty-two. This is a finite, effectively specified set by C/S; it need not be enumerated for this proof. Compression and its inverse derive each local recovery from Sigma.

Finally split any target with at least two simple occurrences at a cut between two of them into nonempty X, Y. Law D gives X squared times Y plus X times Y squared. Both branches preserve content, reduce the number of simple variables, and retain the tail of every surviving simple-tail variable. Iterate. A zero-simple leaf normalizes to g; a leaf with one simple x normalizes to \(w_{x}\). These leaves are below P by the guard and local recoveries, so D derives the original absorption. Literal targets are original summands. Distributivity and absorptions in both directions then derive every identity. This proves unrestricted finite basability of 11954; catalogue transfer requires complete HSP in both directions.

\subsection{A finite basis for 12152 from an original nonlinear contained witness}\label{case:270}
Signature (+, multiplication), no constants, relative to the additively idempotent semiring axioms. The finite family below is specified, not an enumerated small basis. Completeness is unrestricted.

Addition rows:00004,\allowbreak{}01114,\allowbreak{}01224,\allowbreak{}01234,\allowbreak{}44444; order \(3<2<1<0<4\). Multiplication rows:42224,\allowbreak{}12221,\allowbreak{}22222,\allowbreak{}33333,\allowbreak{}42224. Literals retain their values. For nonlinear hT, heads 2 and 3 are constant; head 1 gives 1 if all tail values are in\{0,\allowbreak{}4\}, otherwise2; head 0/\allowbreak{}4 gives 4 if all tails are in\{0,\allowbreak{}4\}, otherwise2. Table induction proves this formula. It yields \(xyz=xzy\) and \(xyy=xy\), normalizing to the head and sorted tail set, with head-tail repetition retained and no empty substitutions.

\par\smallskip\noindent\textit{Exact criterion on ORIGINAL words.}\quad 
A literal x is absorbed iff original literal x occurs. For nonlinear hT: (H) an original word has head h; (NC) an original NONLINEAR word is contained in \(C=\{h\}\) union T; (R) h belongs to T, or literal h occurs, or an original word is entirely contained in T, or an original nonlinear h-word has tail contained in T.

Necessity of literal occurrence:assign \(x=1\), others 3. All right nonlinear x-words give 2 and other heads 3, so no right value absorbs target 1. Necessity of H:assign \(h=2\), others 3; target 2, right 3 separate. For failed NC, assign all target variables 0 and outsiders 3. Target 4; original contained words are literal 0, while all nonlinear words have an outsider and give 2 or 3. The right sum is at most \(0<4\). A contained literal cannot replace NC. For failed R, assuming H and NC, h is outside T. Assign \(h=1\), \(T=4\), others 3. Target 1; literal h is absent. Nonlinear h-words have tails outside T and give 2; heads 3 give 3; words headed 4 must have an outside variable because no original word lies entirely in T, and give 2. Right \(\sum \leq 2<1\) separates.

Sufficiency:3 is least. For target 2, H supplies a same-head original word with value at least 2. Target 4 means all variables of C are in\{0,\allowbreak{}4\}; NC supplies a nonlinear contained word valued 4. Target 1 has head 1, all tails in\{0,\allowbreak{}4\}, and hence h outside T. R supplies literal \(h=1\), or a contained-tail word valued 0/\allowbreak{}4, or a nonlinear same-head word with tail in T valued 1. All absorb 1. Literal targets are original summands. This proves the exact criterion without any bound on variables, word length or polynomial size.

\par\smallskip\noindent\textit{At most three witnesses and finite reversible compression.}\quad 
For a valid nonlinear absorption select an original H witness a and an original NONLINEAR contained witness b. If h is repeated, these suffice. Otherwise select one original R witness c:literal h, tail-contained word, or nonlinear same-head word with tail in T. Q=\{a, b, c\}, discarding duplicates, has at most 3 words, still satisfies H, NC, R, and therefore absorbs the target.

Reserve every head in the at most 4 normal words of target and Q, at most 4 variables. Every other variable occurs only in tails; partition by nonempty incidence among the four tail sets, at most 15 classes. Identify each class to one fresh variable, retain reserved heads separately, and normalize. This produces a valid absorption on at most 19 variables, word length at most 20, and at most 3 right words. Validity is preserved by substitution.

Invert by replacing each class variable with the nonempty product of its original members. Identical incidence recovers each original tail set; reserved heads and their tail appearances stay unchanged. Tail sorting and duplicate deletion recover all original normal words. Nonlinear status is also preserved:an original tail remains nonempty throughout.

Let Sigma contain \(xyz=xzy\), \(xyy=xy\) and EVERY valid normal absorption on a fixed 19-variable alphabet, length at most 20, at most 3 nonempty right words. H, NC, R specify membership in this finite family. Any valid absorption selects Q, compresses into Sigma, and expands back. Distributivity and summandwise absorption in both directions prove every identity. Sigma is a finite relative basis for 12152. Reversal yields the opposite basis; catalogue identifications require independent complete HSP both ways. Finite verification supplements this unrestricted proof.

\subsection{Eleven complete identities for 13045 from contained second-position triggers}\label{case:271}
Signature (+, multiplication), no constants, relative to the additively idempotent semiring axioms. All arguments permit arbitrary variables, word lengths and polynomial sizes.

Addition:00044 /01144 /01234 /44334 /44444. Multiplication:44214 /11211 /22222 /33233 /44244. The additive order \(is 2<1<0<4\) and \(2<3<4\);3 is incomparable with 0,\allowbreak{}1. Element 2 is a two-sided multiplicative zero and the least additive value. A word with a2-valued occurrence is 2. Otherwise a literal retains its value; a nonlinear word headed 1 is 1, headed 3 is 3, and headed 4 is 4. A word headed 0 is 1 when its SECOND occurrence is 3 and 4 otherwise. Later factors only determine whether a2 occurs. Induction on the table proves this formula.

The complete relative basis is:

{\small\par\smallskip
\CaseLine{\((W) xyzt=xytz\)}
\CaseLine{\((H) xyx=xy\)}
\CaseLine{\((T) xyy=xy\)}
\CaseLine{\((I) xyzz=xyz\)}
\CaseLine{\((A) xyz\leq xy\)}
\CaseLine{\((L) x\leq xx\)}
\CaseLine{\((R) xy\leq xxy\)}
\CaseLine{\((B) xxy\leq xy+y\)}
\CaseLine{\((C) xxyz\leq xyz+y\)}
\CaseLine{\((D) xxyt\leq xy+yt\)}
\CaseLine{\((E) xxyzt\leq xyz+yt\)}
\smallskip}

For target w with content C, let U be the original right summands contained in C. Let H be all heads in U, and \(D_{h}\) the second variables of original NONLINEAR U words with head h. A literal h is absorbed exactly when U is nonempty. For nonlinear target w with first two variables h, s the exact criterion is

{\small\par\smallskip
\noindent \(D_{h}\) is nonempty, and (s belongs to \(D_{h}\) or \(D_{h}\) intersects H).\par
\smallskip}

If \(D_{h}\) is empty, assign \(h=0\), other target variables 1, outsiders 2. The target is 4 and all right words are at most 0. If \(D_{h}\) is nonempty but neither alternative holds, \(D_{h}\) misses H and h. Assign \(h=0\), members of \(D_{h}=3\), other target variables 1, outsiders 2. Every h-headed nonlinear contained word is 1; all other contained heads are 1 except possibly literal \(h=0\). The target is 4 because s is not in \(D_{h}\). The right join is at most 0. A missing pure witness for literal h is separated by \(h=0\) and all outsiders 2. These prove necessity for every failed case.

For sufficiency any 2 in the target makes it least. Otherwise all U words are free of 2. If h is 1 or 3, a contained nonlinear h-word has the same value. If h is 4, any such word gives 4. If h is 0 and s is 3, any such word gives 1 or 4 and bounds the target 1. If h is 0 and s is not 3, a second-s witness gives 4. In the remaining case take b in \(D_{h}\) intersect H. Unless an original h, b word already gives 4, b must be 3 and that word gives 1. An original b-headed word gives 3, and \(1+3=4\). A pure witness absorbs its literal because powers reduce to a square and \(x\leq xx\). Thus the criterion is exact.

The first four laws follow from the word formula. W permutes the positions after the second. H removes any later occurrence of the head; in a context it also removes a later second-variable occurrence with a nonempty intervening block. T handles adjacent second repeats and I removes duplicates in the permutable later region. Hence nonlinear words normalize to their original first two variables followed by the sorted content excluding those variables. A square remains a square, and literals remain literal. All substituted words are nonempty.

A is valid because appending can only retain a nonlinear value or lower it to 2. L follows by squaring the table. For R, with no 2 in xy, xxy has value 1 when \(x=1\),3 when \(x=3\),4 otherwise; it therefore bounds xy. For each of B-E, the original right side has a second-y witness headed x and an original y-headed witness, so the intersection clause proves it. This verifies every displayed identity without finite-search inference.

For unrestricted completeness normalize original words and target. A literal uses its original pure witness and L. For a nonlinear target with s in \(D_{h}\), choose an original contained h, s word. Append any missing target content using A: write the current word XY with both factors nonempty, and substitute the nonempty list of missing variables for z. If none are missing, no append is needed. Normalization gives the target.

Otherwise choose b in \(D_{h}\) intersect H, an original \(u=h\) b R and original v headed b. If \(v=b\), use B when R is absent and C when R is nonempty. They produce hh b or hh b R. If \(v=bS\) with S nonempty, use D when R is absent and E otherwise, obtaining hh b S or hh b R S. Each word is below u+v and has repeated head h. Normalize, then append missing target variables using A. The result has head h repeated in its first two positions and full content C. It equals the normalization of hh s Q, where h s Q is the target. R with \(x=h\) and \(y=sQ\) recovers that target; y is always nonempty. If the target already has repeated head, its normal form is already the constructed word. All witnesses are original and contained; appended variables always belong to C.

Thus every valid absorption is derived. Expansion by distributivity and absorption of summands in both directions prove every valid polynomial identity. These eleven laws form a complete finite relative basis for

13045. Reversal transports the basis to its opposite; catalogue transfers still require independent complete HSP inclusions in both directions.

Boundaries:\(xy\leq x\) fails at \(x=0\), \(y=1\); \(xxy\leq xy\) fails at \(x=0\), \(y=3\); \(xz\leq xy+z\) fails at \(x=0\), \(y=3\), \(z=1\). These forbid replacing a nonlinear seed by a literal or replacing the second-variable head trigger by arbitrary content coverage.

\subsection{Three complete finite bases with distinct contained head witnesses}\label{case:272}
All bases are relative to the additively idempotent semiring axioms in signature (+, multiplication), without constants. Proofs allow arbitrary variables, word lengths and polynomial sizes.

All three algebras satisfy \(xyz=xzy\) and \(xyy=xy\). These laws sort and remove repeated actual tail occurrences, preserving any occurrence of the head in the tail. They normalize every nonempty word to a head and tail set.

\par\smallskip\noindent\textit{13034: full contained content and a head-tail trigger.}\quad 
Addition: 00044 / 01144 / 01234 / 44334 / 44444. Multiplication: 44204 / 11211 / 22222 / 44234 / 44244. Addition has order \(2<1<0<4\) and \(2<3<4\); 2 is a multiplicative zero. If any word variable is 2, the word is 2. Otherwise a nonlinear head 1 gives 1, head 4 gives 4, and head 0 or 3 retains its value exactly when all tails are 3, giving 4 otherwise. Literals retain their value. Induction on the multiplication table proves the formula and the two word laws.

The complete eight-identity basis is:

1. \(xyz=xzy\).

2. \(xyy=xy\).

3. \(xyz\leq xy+xz\).

4. \(x\leq xx\).

5. \(xy\leq xxy\).

6. \(xxy\leq xy+y\).

7. \(xxyz\leq xy+yz\).

8. \(xxy\leq xx+y\).

For a target w of content C, take only original right words U contained in C. Let H be their head set and \(D_{h}\) the tails of h-headed U-words. A literal h is absorbed iff U contains a pure power of h. A nonlinear target hT is absorbed iff U covers C, h belongs to H, and either T is contained in \(D_{h}\) or \(D_{h}\) intersects H.

Necessity uses outside variables equal to 2, annihilating words outside U. If h is absent from H, assign \(h=4\) and all other C \(variables=1\). The right sum is at most 1 but the target is 4. If U misses target content after head existence holds, assign its content to 3 and remaining C variables to 0: the target has head 3 and a tail 0, so is 4, whereas the right sum is 3. If \(D_{h}\) misses H and some target tail, assign \(h=0\), \(D_{h}=3\) and the other C \(variables=1\). Original h-words give 0, other heads give 1, but the target is 4. A missing pure-power witness for a literal is separated by \(h=0\) and all \(outsiders=2\). This covers every failure of the criterion.

For sufficiency, a target variable 2 makes the target least; otherwise every U-word is nonzero. Exclude right sum 4. A target head 1 is matched by its U-word; a target head 4 is impossible. For head 0, original h-tails must all be 3. If a tail is also an original head, its word is 3 or 4, whose join with the h-word 0 is 4. Hence only direct \(D_{h}\) coverage can occur, making the target 0. For head 3, an original h-word forces the right sum to be 3. Every U-word is then 3, forcing every U-content variable to 3; full coverage makes the target 3. Literals follow from their pure-power witnesses and \(x\leq xx\). The criterion is exact.

Every displayed inequality satisfies this criterion. To derive a general nonlinear absorption with T contained in \(D_{h}\), merge the contained original h-word tails by law 3. Their union contains T and is contained in T union \{h\}, because all these words are contained in C. If this adds an unwanted repeat of the head, remove it with law 5. Thus this case also yields hT.

Otherwise choose a trigger b, a contained original h-word u with b in its tail, and a contained original b-word v. Tail permutation writes \(u=Ab\), A nonempty. If \(v=b\) use law 6; if \(v=bB\) with B nonempty use law 7. The resulting normalized word g has head h repeated and contains both witness contents. Add missing content from any original contained word v using law 8 with \(x=g\), \(y=v\). Since g has a repeated head, \(gg=g\) follows from the first two laws; the new left word is ggv, adding c(v). Both literal and nonlinear original content witnesses are allowed here. Full coverage reaches C. Recover a simple head with law 5; a repeated head is already correct. Literal targets use law 4 after power reduction. This is an unrestricted nonempty-substitution derivation of every absorption.

\par\smallskip\noindent\textit{13058 and 13321: append tails, then restore a simple head.}\quad 
13058 addition: 00044 / 01144 / 01234 / 44334 / 44444. Multiplication: 44224 / 11221 / 22222 / 33223 / 44224. Let the zero be \(o=2\) and the exceptional head value \(e=3\). For a nonlinear word, occurrence of o anywhere or e in an actual tail gives o. Otherwise heads 1 and e retain their value, and heads 0,\allowbreak{}4 give 4. The zero is the least additive element; e is incomparable with 0,\allowbreak{}1.

13321 addition: 00004 / 01114 / 01224 / 01234 / 44444. Multiplication: 44334 / 11331 / 22332 / 33333 / 44334. Here \(o=3\), \(e=2\), with the same word formula after these replacements. Addition is the chain \(3<2<1<0<4\). These formulas follow by induction on the respective multiplication tables. Literals retain their value.

The complete four-identity basis for 13058 is:

1. \(xyz=xzy\).

2. \(xyy=xy\).

3. \(xyz\leq xy\).

4. \(xy\leq xxy+x\).

The complete six-identity basis for 13321 consists of these four and:

5. \(xy\leq xxy+y\).

6. \(xyz\leq xxyz+y\).

Again let U consist of original right words entirely contained in target content C. For either algebra, a literal h is absorbed iff it is an original literal. A nonlinear hT is absorbed iff there is an original nonlinear h-headed word in U, and, when h is simple in the target:

\textbullet{} in 13058, there is an original h-headed U-word with h absent from its tail; this witness may be literal h;

\textbullet{} in 13321, there is an original U-word with h absent from its tail; its head can differ from h, so it may omit h altogether.

No full-content coverage of C is required in these two criteria.

For necessity, set outsiders to o. Absence of a nonlinear h-headed U-word is separated by \(h=0\) and all other C \(variables=1\): the target is 4 while the right sum is at most 0. In 13058, failure of the simple-head witness is separated by \(h=e=3\) and all other C \(variables=1\). Contained h-words with repeated h give \(o=2\); other heads give at most 1, which does not absorb target 3. In 13321 the corresponding assignment \(h=e=2\), \(others=1\) makes every U-word \(o=3\) when each has h in its tail, but the target is 2. For a missing literal h, set \(h=e\) and \(outsiders=o\); all contained powers of degree at least two give o, below target e. These are exhaustive separators for the stated failures.

For sufficiency, a target with value o is least. Otherwise C avoids o and all target tail values avoid e. If h is 0 or 4, a contained original nonlinear h-word gives 4; if \(h=1\), it gives 1. If \(h=e\), the target must be simple at its head. In 13058 the simple h-headed witness gives e. In 13321 a witness whose tail avoids h gives at least e: it either has head h and gives e, or uses only allowed values 0,\allowbreak{}1,\allowbreak{}4 and gives 0,\allowbreak{}1 or 4. Thus the target is always absorbed. The criteria also verify all listed inequalities directly.

For the general derivation, first suppose the desired head is repeated, or an original nonlinear contained h-word has simple head. In the former case choose any nonlinear h-witness; in the latter choose a simple one. Append missing target tail variables one at a time by law 3, then normalize. Every original tail already belongs to the desired tail, so this yields the target. All substitutions of the existing tail are nonempty.

The remaining case starts with a repeated-head nonlinear h-witness and constructs the repeated-head word hhT by appending missing tails. In 13058, the criterion now forces the simple witness to be the literal h; law 4 recovers hT. In 13321, literal h likewise uses law 4. Any remaining allowed witness u has a head different from h and no h in its tail, so omits h entirely. Its content is contained in T. Write the target tail as uV after tail normalization, where V lists missing target tail variables. If V is absent, apply law 5 with \(x=h\), \(y=u\). If V is nonempty, apply law 6 with \(x=h\), \(y=u\), \(z=V\). The repeated-head right word is already available and the other right word is the original u. This recovers the simple head. No empty word is substituted; the absent-remainder case is a separate law.

Distributivity and absorption of every summand in both directions prove completeness for all polynomial identities in each algebra.

Boundary distinctions: \(xy\leq xyz\) fails in 13034 at \(x=1\), \(y=1\), \(z=2\); \(xz\leq xy+y\) fails there at \(x=3\), \(y=3\), \(z=0\). In 13058, \(xy\leq xxy+y\) fails at \(x=3\), \(y=1\), although it is valid in 13321. The law \(x\leq xx\) fails at \(x=e\) in both 13058 and 13321. In 13321, \(xy\leq xxy\) fails at \(x=2\), \(y=1\). These examples preclude conflating the three witness criteria. Any equal-variety or opposite transfer requires complete HSP both ways.

\subsection{A finite basis for 13100 by nonlinear-head simple-variable recovery}\label{case:273}
Signature (+, multiplication), without constants, relative to the additively idempotent semiring axioms. The finite family below is effectively specified; no small enumeration or minimality is claimed.

Addition: 00004 / 01114 / 01224 / 01234 / 44444. Multiplication: 44234 / 11231 / 22332 / 33333 / 44234. The additive order is \(3<2<1<0<4\). Put K=\{0,\allowbreak{}1,\allowbreak{}4\}. A literal retains its value. A nonlinear word containing a 3, or at least two occurrences valued 2, has value 3. With exactly one occurrence valued 2 it has value 2. With all factors in K it is 1 if its head is 1, and 4 otherwise. These assertions follow inductively from the multiplication table and apply to every positive word length.

\par\smallskip\noindent\textit{Exact original-witness criterion.}\quad 
For a target word w let C be its content and let U consist of original right summands whose entire content is contained in C. A literal h is absorbed exactly when h is an original literal. A nonlinear w is absorbed exactly when both of the following hold:

H: U contains a NONLINEAR word with the same head as w. S: Every variable appearing exactly once in w has a word in U which omits that variable or contains it exactly once. The S witness may be literal, and neither endpoint nor order is constrained.

For a failed H put the target head equal to 0, the other target variables equal to 1 and all outsiders equal to 3. The target is 4. Every right word outside U is 3; every nonlinear word inside U has head different from the target and is 1; literals are 0 or 1. The right join is below 4. For a failed S at x, assign \(x=2\), other target variables 1 and outsiders 3. The target is 2, while all original right words are 3. If a literal h has no original literal witness, assign \(h=2\) and all other variables 3: the target is 2 and every right word is 3. These are separators for every failed case, for any variable set.

Conversely target value 3 is least. A target value 2 has exactly one 2-valued occurrence; its variable x is simple. An S witness omitting x is K-valued and above 2, while a single-x witness gives 2. If the nonlinear target is K-valued, the nonlinear contained same-head witness has the same value 1 or 4. The literal case is immediate. This proves the criterion for arbitrary words and original polynomials.

\par\smallskip\noindent\textit{Eight structural identities.}\quad 
{\small\par\smallskip
\CaseLine{\((W) xyz = xzy\)}
\CaseLine{\((T) xy^{3} = xy^{2}\)}
\CaseLine{\((H) x^{3}y = x^{2}y\)}
\CaseLine{\((P) x^{3} = x^{2}\)}
\CaseLine{\((D) xy = x^{2}y + xy^{2}\)}
\CaseLine{\((L) x^{2}y \leq  xy\)}
\CaseLine{\((R) xy^{2} \leq  xy\)}
\CaseLine{\((A) xyz^{2} \leq  xy\)}
\smallskip}

The word formula verifies W, T, H, P. For D, if XY has value 3 both branches have value 3; if its unique 2-valued occurrence lies in X, then \(XY^{2}\) retains value 2 while \(X^{2}Y\) becomes 3, and conversely for Y. With all occurrences K-valued both branches preserve the head and nonlinearity. L, R can only retain a word value or change its unique 2-valued occurrence to the least value 3. A appends a square: a K square retains the nonlinear head result, and a 2/\allowbreak{}3 square gives 3. Thus all eight laws hold, including arbitrary nonempty substitutions.

W permutes the tail. Bring repeated head copies next to the head and use H, or P for singleton content, to cap head multiplicity at two. For a nonhead variable, bring its copies together in the tail and use T, whose prefix is nonempty, to cap multiplicity at two. Sort the remaining tail using W. Normal words are determined by the literal head and each total multiplicity capped at two. Literals stay literal.

\par\smallskip\noindent\textit{Unrestricted construction.}\quad 
For a valid nonlinear absorption choose an original nonlinear H witness. Double its simple head using L, and its other simple variables using R in their contexts. Every doubling has a nonempty required prefix or suffix because this word is nonlinear. Append squares of missing target variables using A, splitting the existing nonlinear word into two nonempty factors. Normalize throughout. This derives a guard g below P with content C, target head, and every variable repeated twice.

For each simple target variable x choose an original S witness u. If u is nonlinear, double all its simple variables except x using L, R. If u is literal, KEEP IT LITERAL. Call the result \(u_{x}\). Let \(w_{x}\) be the normal word obtained from w by doubling every variable except x. The absorption

{\small\par\smallskip
\CaseLine{\(w_{x} \leq  g + u_{x}\)}
\smallskip}

is valid by H, S: g supplies the nonlinear same-head witness and \(u_{x}\) supplies the sole simple-variable witness. The finite recovery family below derives it. This argument covers omitting literals and literal x.

For any normal nonlinear word with at least two simple variables, split it as XY between two such occurrences, with X, Y nonempty. D replaces it by \(X^{2}Y\)+\(XY^{2}\). Normalize both branches. Each keeps content and head, remains nonlinear, and has strictly fewer simple variables. Recurse. A zero-simple leaf equals g; a one-simple-x leaf equals \(w_{x}\). Every leaf is below P, so D gives \(w\leq P\). No limit on the number of variables, original summands or recursive steps is imposed. Literal targets were already original summands.

\par\smallskip\noindent\textit{Reversible compression and the fixed finite family.}\quad 
Each local recovery has three words. Reserve x and the heads of these words, at most four variables. Each remaining variable, in each word, is absent or occurs twice and is never its head. This includes a literal \(u_{x}\) because its only variable was reserved as its head. Group remaining variables by their nonzero incidence pattern across the three words, at most seven classes. Identify each class to a distinct fresh variable and normalize. This is a substitution into a valid absorption, preserving validity. There are at most eleven variables and word lengths at most 22.

Replace each class variable by the nonempty word listing its members. Every occurrence of a class in a compressed word has multiplicity two; the expansion gives each member twice, all in the tail. W permutes the tail, and T, H, P normalize. Reserved variables retain their head and individual multiplicities. Thus all three original normal words are recovered. Common multiplicities are preserved; no group of single occurrences is silently treated as an original repeated variable.

Fix eleven variables. Add to the eight structural laws ALL valid absorptions whose left side is one normal word of length at most 22 and whose right side has one or two such normal words. H, S and the literal case specify this finite family effectively. Every local recovery is an inverse nonempty substitution instance modulo the word laws. The guard, recovery and splitting construction therefore derives every valid word absorption. Distributivity and absorption of summands in both directions yield every polynomial identity. This proves finite basability for 13100 without bounded-search assumptions. Word reversal gives its opposite; any catalogue transfer still requires complete HSP both ways.

Boundary distinctions: \(x^{2}\leq x\) fails at \(x=0\); \(x\leq x^{2}\) fails at \(x=2\); \(xy\leq x\) fails at \(x=0\), \(y=1\). A literal same-head witness cannot replace H, and original literal S witnesses must not be squared.

\subsection{A finite basis for 13101 by nonlinear heads and simple final recovery}\label{case:274}
Signature (+, multiplication), no constants, relative to the additively idempotent semiring axioms. The finite family is specified effectively, not enumerated as a small basis.

Addition:00004 /01114 /01224 /01234 /44444. Multiplication:44234 /11231 /33333 /33333 /44234. The additive order \(is 3<2<1<0<4\). Put K=\{0,\allowbreak{}1,\allowbreak{}4\}. A literal keeps its value. For a nonlinear word, a3-valued occurrence or a2-valued occurrence strictly before the final position gives 3. If the final position is 2 and all earlier positions are K-valued, the value is 2. With all positions K-valued, the result is 1 for head 1 and 4 otherwise. Direct induction on the table proves these formulas without a length bound.

For a target w, retain original right words U with content contained in c(w). A literal target is absorbed exactly when it is an original literal. A nonlinear target is absorbed exactly when: H: U contains a NONLINEAR word headed by the head of w. T: If the final variable t occurs exactly once in w, some U word omits t or contains t exactly once at its final position. T is vacuous if the target final variable is repeated. Original literal witnesses are allowed.

If H fails, assign the target head 0, other target variables 1, outsiders 3. The target is 4 and right words have values at most 0. If T fails, assign \(t=2\), other target variables 1, outsiders 3. The target is 2, all right words are 3. A missing literal h is separated by \(h=2\), all other variables 3. Conversely value 3 is least. For target 2 its unique 2 is the simple final variable, and T supplies value 2 or a K-valued larger value. For target 1 or 4, H gives exactly that value. Thus H, T are necessary and sufficient for arbitrary original polynomials, including outsiders.

\par\smallskip\noindent\textit{Eight structural laws.}\quad 
{\small\par\smallskip
\CaseLine{\((W) xyzt = xzyt\)}
\CaseLine{\((H) xxy = xy\)}
\CaseLine{\((I) xyyz = xyz\)}
\CaseLine{\((T) xyyy = xyy\)}
\CaseLine{\((P) xxx = xx\)}
\CaseLine{\((E) xyx = xyy\)}
\CaseLine{\((F) xyzy = xyzx\)}
\CaseLine{\((A) xyz \leq  xy\)}
\smallskip}

The word formula verifies the first seven laws. W preserves endpoints, content and simple-final status. H, I only delete repeated nonfinal positions; T, P preserve repetition of the final variable. Both sides of E, F have a repeated final variable and the same content and head. For A, a nonlinear prefix valued 3 stays 3; one valued 2 can only fall to 3; a K-valued prefix can retain its head result or fall to 2/\allowbreak{}3. These facts verify all laws, with arbitrary nonempty substitutions.

\par\smallskip\noindent\textit{Literal normal forms.}\quad 
If a nonlinear word has simple final variable t, keep its head h and t, remove all duplicate internal occurrences, and sort its interior. H deletes extra head copies after W brings them next to the head; I deletes other duplicate interior copies after W groups them. Both required outer contexts in I are nonempty. The normal word is

h [sorted c(w) minus \{h, t\}] t.

If the last variable is repeated, a normal word is instead

h [sorted c(w) minus \{h\}] h,

including hh for singleton content. To obtain this, if the old last variable b differs from h, move an earlier b immediately after h using W. The word is h b R b. When R is nonempty, F changes this to h b R h. When R is absent, E used backwards changes hbb to hbh. Now the last variable is h. H, I remove surplus interior copies, and W sorts them. Singleton powers use P. T is valid but redundant for this normalization. All substitutions are nonempty; only external contexts may be absent.

Thus normal words retain head, content, and either their simple final variable or the fact that the final variable is repeated. Literals remain literal. In particular all repeated-final words with identical head and content are equivalent under the displayed laws.

\par\smallskip\noindent\textit{Unrestricted construction and finite compression.}\quad 
Take a valid nonlinear \(w\leq P\) and an original nonlinear contained H witness v. Split \(v=XY\) with X, Y nonempty. Append a word consisting of all missing target variables followed by the target head h, using A. This appended word is always nonempty. Normalize to a guard g with full target content, target head and repeated final h. Hence \(g\leq v\leq P\). If w has repeated final variable, it normalizes to g and is derived.

Otherwise its final variable t is simple. Choose an original contained T witness u. Keep u, including a literal witness, and normalize it. The local absorption \(w\leq g+u\) holds by H, T. No saturation or squaring of a literal is performed. The following fixed finite family derives every such local recovery.

Reserve the heads and final variables of the three normal words w, g, u. There are at most four distinct reserved variables: h, t, head(u), tail(u). Every nonreserved variable occurs once or not at all in each word and is strictly internal. Group these variables by their nonzero incidence pattern across the three words, at most seven groups. Identify a group to one fresh variable and normalize. Validity is preserved by this substitution. At most eleven variables remain and each normal word has length at most twelve, with at most two right summands.

For the inverse substitute the nonempty list of the group members for each group variable. Each group occurs exactly once or not at all in every compressed word, in an internal position. Its expansion restores each member once; W permutes the interior to the original normal order. Reserved endpoints and their multiplicities are unchanged. In particular the simple target final remains simple and a repeated final remains repeated. Common multiplicity one is preserved in all words; repeated endpoints are reserved and never merged into incidence groups.

Fix eleven variables. Use the eight displayed laws together with ALL valid absorptions having one normal left word and one or two normal right words, every word of length at most twelve. H, T and the literal case specify this finite family effectively. Inverse compression derives every local recovery; A derives the full-content guard. Every valid word absorption is therefore derived without any bound on original variables, word lengths or polynomial size. Distributivity and two-way summandwise absorption prove completeness for all identities. This proves 13101 FB. Word reversal gives its opposite; catalogue transfers require full HSP both ways, independently verified.

Boundary failures: \(x\leq xx\) at \(x=2\); \(xx\leq x\) at \(x=0\); \(xy\leq x\) at \(x=0\), \(y=1\); \(xy\leq xyy\) at \(x=1\), \(y=2\). These prevent replacing the nonlinear H witness by a literal or discarding the simple-final witness condition.

\subsection{Four finite-basis sources from contained nilpotent ideal domains}\label{case:275}
Use the semiring axioms without constants, four word laws \(xyz=yxz\), \(xyz=xzy\), \(xxy=xyy\), \(xxyz=xyz\), and every valid absorption over 40 variables, length 42, at most 18 upper words. Interpret the valid finite family separately in each source. It is effective and finite; it is not an enumerated small basis. Completeness below is unrestricted.

Each source has a multiplicatively absorbing element Z, also the least addition element, and addition greatest element 4. A nonlinear word is Z if some letter is Z; is R if quadratic with its first value in H and second in T; and is 4 otherwise. In particular every long word without Z is 4. The actual roles are:

\noindent Source\enspace\textbar\enspace Z\enspace\textbar\enspace R\enspace\textbar\enspace H\enspace\textbar\enspace T\enspace\textbar\enspace Proper principal ideals\par
\noindent 12030\enspace\textbar\enspace 1\enspace\textbar\enspace 3\enspace\textbar\enspace \{0\}\enspace\textbar\enspace \{2\}\enspace\textbar\enspace \{0,\allowbreak{}1\};\{1\};\{1,\allowbreak{}2\};\{1,\allowbreak{}3\}\par
\noindent 12250\enspace\textbar\enspace 2\enspace\textbar\enspace 3\enspace\textbar\enspace \{0,\allowbreak{}1\}\enspace\textbar\enspace \{1\}\enspace\textbar\enspace \{0,\allowbreak{}1,\allowbreak{}2\};\{1,\allowbreak{}2\};\{2\};\{1,\allowbreak{}2,\allowbreak{}3\}\par
\noindent 12251\enspace\textbar\enspace 2\enspace\textbar\enspace 3\enspace\textbar\enspace \{0,\allowbreak{}1\}\enspace\textbar\enspace \{1\}\enspace\textbar\enspace \{0,\allowbreak{}1,\allowbreak{}2\};\{1,\allowbreak{}2\};\{2\};\{2,\allowbreak{}3\}\par
\noindent 13897\enspace\textbar\enspace 3\enspace\textbar\enspace 0\enspace\textbar\enspace \{1,\allowbreak{}2\}\enspace\textbar\enspace \{2\}\enspace\textbar\enspace \{0,\allowbreak{}1,\allowbreak{}2,\allowbreak{}3\};\{1,\allowbreak{}2,\allowbreak{}3\};\{2,\allowbreak{}3\};\{3\}\par
Induction on word length proves the value formula from the operation tables. The four laws normalize a long word to its sorted content, using a cube for singleton content and a repeated first variable for two-letter content. Literals and quadratics preserve their order and degree. Common nonlinear content determines every long word value.

For arbitrary target w, let U be the original right words whose contents are contained in c(w). To falsify \(w\leq P\) in any principal ideal I, a target assignment cannot give any target variable value Z:then \(w=Z\) would lie in I. Set every external variable Z. This makes all original words outside U equal Z, including external literals. Thus failure is exactly an assignment of target variables in \(B=A\) minus\{Z\} for which all U words lie in I but w does not. This argument proves the restriction to U in both directions; it is not a heuristic pruning rule.

A contained long original makes the absorption valid immediately:if some target variable is Z the target is least; otherwise that long original is 4. If U has no long word, solve its independent domains for each proper ideal I. A literal x restricts x to I intersect B. If R is not in I, any quadratic of U makes the system impossible. If R is in I, a quadratic xy restricts x to H and y to T. All these constraints are unary; intersect them for each target variable. An empty domain makes this ideal impossible. Otherwise a literal target is forced into I exactly when its domain is contained in I; a quadratic is forced exactly when R is in I and its first/second domains lie in H/T; and a long target always fails. The conjunction over the four proper ideals is an exact criterion for arbitrary alphabets and P.

This criterion yields bounded original witnesses. If there is a contained long word, retain it alone. Otherwise U consists only of literals and quadratics. For a short target retain witnesses for its variables' literal/head/tail roles, at most 6 words. For each impossible ideal retain one original quadratic when R is excluded, or witnesses for the literal/head/tail roles of one empty-domain variable, at most 3. There are four ideals, so at most 18 originals suffice in total. For a long target all ideals are impossible, and at most 12 originals suffice. Every possible ideal has the same target-variable domains after this selection, while every impossible ideal retains its explicit obstruction. Thus the selected polynomial still absorbs w.

Normalize selected words and target. In the contained-long branch there are two words; reserve up to 3 distinct variables per word and group other variables by common incidence. At most 6 reserved variables and 3 nonempty incidence groups give at most 9 variables. In the other branch all 18 selected words are short, so reserve their at most 36 variables and up to 3 target variables. Any remaining variable occurs only in the long target and forms a SINGLE incidence group. Thus at most 40 compressed variables suffice. The reserved variables preserve short and long degree; normal word length is at most 42.

Identify each group to one variable. Validity is preserved by variable identification, and the normalized compressed absorption belongs to the finite family. Inversely substitute the nonempty list of group members. Every member has common incidence, all short occurrences are reserved, and long degree is preserved. Hence the four word laws normalize each expanded word exactly back to its original. This derives the selected absorption. Add back unselected original summands. Splitting identities into summand absorptions in both directions proves unrestricted finite basability. No assumption about a global short-word match or one-way HSP is used; equal-variety transfers are verified separately.

\subsection{A 253-identity finite basis for 12095}\label{case:276}
Signature (+, multiplication), no constants. All substituted words are nonempty. Write \(p\leq q\) for \(p+q=q\). This proof is independent of the different additive orders of 13066 and 8301.

\par\smallskip\noindent\textit{Tables and word normal form.}\quad 
Use coordinates transported to catalogue type 12095 by (3,\allowbreak{}1,\allowbreak{}0,\allowbreak{}2,\allowbreak{}4). Addition rows are 00444,\allowbreak{}01224,\allowbreak{}42224,\allowbreak{}42234,\allowbreak{}44444. Multiplication rows are 01444,\allowbreak{}11444,\allowbreak{}41444,\allowbreak{}41234,\allowbreak{}41444. The additive principal ideals are \{0,\allowbreak{}1\}, \{1\}, \{1,\allowbreak{}2,\allowbreak{}3\}, \{3\}, and the whole set, respectively. In particular \(1+3=2\); neither 8301's order nor its absorption criterion applies.

Table induction gives: an all-0 word is 0; a word whose rightmost non-0 factor is 1 is 1; an all-3 word is 3; a word with final factor 2 and all preceding factors 3 is 2; every other word is 4. The laws \(xxy=xy\) and \(xyx=yxx\) normalize a word to its last-occurrence order, with its last letter doubled exactly when it repeats. To remove an internal repeated x in xYxZ, use \(xYxZ=YxxZ=YxZ\) (Z nonempty); empty Y uses \(xxy=xy\) directly. A repeated last x uses \(xYx=Yxx\).

\par\smallskip\noindent\textit{Original-word reachability and the exact criterion.}\quad 
Let P be the original nonempty finite set of upper words. In a directed graph on their variables put an edge \(y\to z\) whenever an original word ends in y and contains z before that final position. Write \(R_{x}\) for all vertices reachable from x, including x itself. A positive cycle at x means an original word ends at some y in \(R_{x}\) and contains x internally; this is equivalent to a positive-length return path to x, and must not be confused with the reflexive membership x in \(R_{x}\).

For target w with content C and last variable t, and each x in C, put S\_x=c(the suffix after the last x of w). Let \(W_{x}\) consist of original words whose last-x suffix is contained in \(S_{x}\). Then \(w\leq \sum\) P holds exactly when all the following conditions hold:

1. An original word has content contained in C.

2. Every \(W_{x}\) is nonempty.

3. C is contained in \(R_{t}\). If t repeats in w, t has a positive cycle.

4. For each x with \(S_{x}\) nonempty, at least one holds: \(R_{x}\) meets \(S_{x}\); an original word is contained in \(S_{x}\); or an original word has an occurrence of some y in \(R_{x}\) followed by a NONEMPTY suffix contained in \(S_{x}\). One can use its last y, since \(R_{x}\) and \(S_{x}\) are disjoint when the first alternative fails.

Necessity uses explicit assignments. If (1) fails, assign C to 0 and all other variables to 1: w is 0 and every upper word is 1. If \(W_{x}\) is empty, assign \(x=4\), \(S_{x}=0\), \(others=1\): w is 4 and upper words lie in \{0,\allowbreak{}1\}. If C is not contained in \(R_{t}\), assign \(R_{t}=3\) and \(others=1\). Every upper word whose last variable lies in \(R_{t}\) is all 3, and every other upper word ends in 1. Thus the upper sum is \ensuremath{\leq}2, but w is 4. If repeated t has no positive cycle, assign \(t=2\), \(R_{t}\) minus \(\{t\}=3\), \(others=1\). Every original t-tail is simple and otherwise all 3; an original word with another tail in \(R_{t}\) contains no t, or it would give a cycle. Hence all upper values lie in \{1,\allowbreak{}2,\allowbreak{}3\}, while w is 4. Finally, if (4) fails, assign \(S_{x}=0\), \(R_{x}=3\), \(others=1\). These sets are disjoint. No upper word is all 0. Its rightmost non-0 factor is either 1, giving 1, or an \(R_{x}\) letter. The latter cannot be followed by any 0 (that would be the third alternative); as its actual tail lies in \(R_{x}\), its entire content lies in \(R_{x}\), so its value is 3. The upper sum is \ensuremath{\leq}2, but w has rightmost non-0 factor 3 followed by at least one 0, so has value 4. All failures are separated.

For sufficiency test the four proper principal ideals. When all upper values lie in \{0,\allowbreak{}1\}, any target value outside this set has a rightmost non-0 letter x valued 2,\allowbreak{}3, or 4. Its \(W_{x}\) witness also lies outside \{0,\allowbreak{}1\}, a contradiction. If all upper values equal 3, all their variables are 3; (2) covers the target content, so w is 3.

Suppose all upper values lie in \{1,\allowbreak{}2,\allowbreak{}3\}. Target 0 is excluded by (1). For target 4 consider its rightmost non-0 letter x; all \(S_{x}\) letters are 0. Value \(x=4\) is excluded by \(W_{x}\). If \(x=3\), every \(R_{x}\) letter must be 3: an original word ending in a forced-3 letter must be all 3. For nonempty \(S_{x}\), each alternative in (4) is now impossible under the upper-value assumption: a forced-3/\allowbreak{}0 collision, an upper 0, or an upper 4 with a nonempty 0 suffix. For empty \(S_{x}\), (3) makes the entire target 3.

If \(x=2\), every \(W_{x}\) word must be simple-final x with all other factors 3; otherwise it is 4. Consequently every vertex of \(R_{x}\) minus \{x\} is forced to 3, and a positive return to x is impossible. This follows along each original tail-to-prefix edge, starting at x and then using the forced-3 rule. For nonempty \(S_{x}\), a meeting with \(R_{x}\) contradicts 0 versus 3. A suffix-contained upper word gives 0. In the third alternative, \(y=x\) gives a nonterminal 2, and any other reachable y gives a 3 followed by 0; both give upper 4. For empty \(S_{x}\), condition (3) rules out a repeated target x and places its whole prefix among the forced-3 letters, so w is 2. Thus w is in \{1,\allowbreak{}2,\allowbreak{}3\}. Finally, \{1\}=\{0,\allowbreak{}1\} intersect \{1,\allowbreak{}2,\allowbreak{}3\}, so that ideal is also preserved. Every principal ideal test holds, proving unrestricted exact semantics.

\par\smallskip\noindent\textit{The finite laws.}\quad 
There are nine fixed laws:

{\small\par\smallskip
\CaseLine{\(xxy=xy; xyx=yxx; x\leq xx; xyz\leq xz+yz;\)}
\CaseLine{\(xy\leq x+zxy; xy\leq y+zxy; xyz\leq xz+vxyz;\)}
\CaseLine{\(xyzt\leq xty+vxyzt; yzt\leq ty+vyzt.\)}
\smallskip}

Append the following separately written optional-context templates. Brackets mean choose either the displayed variable or omit that piece while WRITING a law. They never permit an empty substitution. Lowercase and uppercase letters designate distinct formal variables.

{\small\par\smallskip
\CaseLine{\([a]bxrp \leq  [a]xbrp+[v]xr\qquad (4)\)}
\CaseLine{\([a]bx[c][d][p]t \leq  [a]xb[c][d][p]t+[v]x[c]t[d]\qquad (32)\)}
\CaseLine{\([a]bxsp \leq  [a]xbsp+[v]s[r]x\qquad (8)\)}
\CaseLine{\([a]bx[p]t \leq  [a]xb[p]t+[v]t[r]x\qquad (16)\)}
\CaseLine{\([A]bxdp \leq  [A]xbdp+[v]a[r]x+[w]ad\qquad (16)\)}
\CaseLine{\([A]bx[c][d][p]t \leq  [A]xb[c][d][p]t+[v]a[r]x+[w]a[c]t[d]\qquad (128)\)}
\CaseLine{\([A]bxyp \leq  [A]xbyp+[v]x+y\qquad (4)\)}
\CaseLine{\([A]bx[c][d][p]t \leq  [A]xb[c][d][p]t+[v]x+[c]t[d]\qquad (32)\)}
\CaseLine{\([a]by[d] \leq  [a]y[d]+by\qquad (4)\)}
\smallskip}

The total is \(9+244=253\). Notice d in the fifth template is mandatory. The corresponding empty-d assertion is in general false. These laws are checked symbolically by the unrestricted exact criterion on their formal words. The two word laws are checked in both directions. This is a soundness proof for all assignments, not a bounded sample.

\par\smallskip\noindent\textit{Unrestricted completeness.}\quad 
Normalize P and w. Normalization preserves the tail-to-prefix graph (including self edges), original suffix sets and all conditions above. Choose original U contained in target content. The last template allows prefix enlargement: if a known \(g=A\) y D contains y and an original word is B y, derive A B y D, then normalize. Empty A, D select separate templates; B is nonempty, and a literal witness needs no enlargement. This enlarges content by all of B.

Start with an original t-tail word. Repeated enlargement along original edges collects exactly \(R_{t}\). If there is a positive cycle at t, an original word containing t internally and ending at a reachable letter adds another t. Conversely no enlargement repeats t without such a cycle. Each content-changing step adds a variable, and the extra cycle step adds the repeated-tail bit, so finitely many steps give a known accumulator g covering \(R_{t}\), with repeated t precisely in the cyclic case. By (3) this contains all target variables.

Build the correct target suffix from right to left. Initially \(S=t\) or tt according to the accumulator. For the next target variable x, write the normalized accumulator as A x B S. If B is empty, x is already positioned. Otherwise move B in front of x, as follows.

If an original \(W_{x}\) witness is V x R with R nonempty and c(R) contained in c(S), the first template moves x. If the final t is simple and R contains t, split \(R=C\) t D and use the second template to avoid introducing an extra t. In an acyclic seed every original word ending in \(R_{t}\) contains no internal t; hence this t is unique. The same observation applies to all guarded-tail cases below.

Otherwise \(W_{x}\) consists of words ending in x. If an original q is contained in c(S), use the seventh template with q as y and any original Vx witness. When t is simple and occurs in q, split \(q=CtD\) and use the eighth template.

In the remaining case (4) supplies reachability. If some s in \(R_{x}\) lies in c(S), start with any original x-tail witness v. Along a path from x to s, enlarge v using the original edge witnesses. This gives a known x-tail word V s R x. Apply the third template; if \(s=t\) and t is simple, use the fourth template. All path words are derived by the finite enlargement laws, not assumed to be original witnesses.

If \(R_{x}\) misses c(S), (4) gives an original Q y D with y in \(R_{x}\), D nonempty and c(D) contained in c(S). The case \(y=x\) was already handled as a true suffix witness. Enlarge an original x-tail witness along a path to y, obtaining V y R x. Use the fifth template with the mandatory d substituted by D. When t is simple and D contains t, split \(D=CtD0\) and use the sixth template, writing the protected final t once. Optional empty pieces select the appropriate written law.

In all cases the first right term normalizes to g, and the other right terms are original or explicitly enlarged known words. The left normalizes to A B x S. For a repeated t, duplicate suffix pieces disappear under the word laws while preserving the repeated bit. For simple t, the designated tail templates avoid duplication. Thus each step lengthens the correct suffix and preserves tail status.

At the end \(g=E\) Q, where Q is the entire target last-occurrence order with the seed's tail status and E consists of external letters. If E is nonempty use the fixed contained-word erasure laws: \(xy\leq x+zxy\) with \(x=U\), \(y=Q\), \(z=E\) when Q has repeated t or U omits t; \(xy\leq y+zxy\) when \(U=t\); \(xyz\leq xz+vxyz\) when U has simple final t; and the last two fixed laws for a contained U having internal/initial t. Their right terms normalize to U and E Q and their left to Q. (The latter internal-t cases are unnecessary for an acyclic seed but are sound.) If Q repeats t and w is simple, use \(x\leq xx\) with \(x=t\) in its prefix context. Every substitution remains nonempty. This derives w.

The construction handles every finite P and w without any bound on variables, lengths, paths, cycles or summands. Polynomial identities split into their word absorptions, proving the 253 laws complete. Classification of other catalogue types requires separate full mutual HSP or isomorphism certificates.

\par\smallskip\noindent\textit{Research provenance.}\quad

\subsection{Seven complete identities for 13349}\label{case:277}
Signature (+, multiplication), without constants, relative to the additively idempotent semiring axioms. The complete finite basis is:

1. \(xyz=xxx\).

2. \(xy\leq xxx\).

3. \(x\leq xxx\).

4. \(xz\leq xy+yz\).

5. \(xxx\leq x+xy\).

6. \(xxx\leq xy+y\).

7. \(xxx\leq xy+yyy\).

The proof is unrestricted in variables, word lengths and summands.

\par\smallskip\noindent\textit{Word values.}\quad 
Addition: 00044 / 01034 / 00244 / 43434 / 44444. Multiplication: 44344 / 11111 / 44344 / 44444 / 44444. The additive order has \(1<0\), \(2<0\), \(1<3\), and top 4. In particular 2 is not below 3 and 0 is not below 3; the join of 1 and 3 is 3.

A literal retains its value. A quadratic word headed 1 is 1. A quadratic word headed 0 or 2 with final value 2 is 3. Every other quadratic value is 4. Every word of length at least three is 1 if its head is 1, and 4 otherwise. These statements follow by induction on the table, since 3 and 4 multiply on the left to 4, and 1 remains 1. This proves law 1. Normalize each long word to the cube of its first variable, using law 1 with its first two letters and its nonempty remaining suffix.

\par\smallskip\noindent\textit{Exact criterion from original words.}\quad 
Let L be the heads of original literal summands and K the heads of original summands of length at least three. Original quadratic words ab are directed edges \(a\to b\). For each h, let \(R_{h}\) be the vertices reached by positive-length paths in this graph, without adding reflexive paths. Call h active when either:

\textbullet{} h belongs to K; or

\textbullet{} \(R_{h}\) is nonempty and either h belongs to L or \(R_{h}\) meets L union K.

A literal h is absorbed iff h is an original literal or h is active. A quadratic ht is absorbed iff t is positively reachable from h or h is active. A word of length at least three headed h is absorbed iff h is active. Words may use outside variables; no containment condition applies.

To prove necessity, first suppose h has no outgoing quadratic edge and is not in K. Assign \(h=0\) and every other \(variable=1\). If h is literal, the right sum is 0; otherwise it is 1. A nonlinear target is 4, unless it is quadratic with final value 2, which this assignment never gives. A missing literal h is also separated. Thus all failures in this case are explicit.

Now suppose h has an outgoing edge but is not active. Then h is not in L or K, and \(R_{h}\) meets neither L nor K. Assign \{h\} union \(R_{h}\) to 2 and all others to 1. Every original quadratic headed in that set has final value 2 and gives 3; other quadratic and long words give 1. All literals give 1. The right sum is 3. A long target gives 4, a literal target gives 2, and a quadratic ht with t different from h and absent from \(R_{h}\) gives

4. Each fails to be below 3. For the remaining quadratic hh with no positive return, assign \(h=0\), \(R_{h}=2\) and \(others=1\). No edge from \(R_{h}\) returns to h, so the right sum is still 3, while hh gives 4. This proves every necessary condition, including the genuine positive-return requirement.

For sufficiency of a quadratic path, suppose the right sum is not 4. If \(h=1\), the first original path edge gives 1 and absorbs the target. For other h values, the first edge can avoid 4 only if h is 0 or 2 and its endpoint is 2. Every subsequent edge forces its endpoint to 2 as well. Hence the final target ht is 3, matched by the first edge.

For an active head, an original long h-word immediately supplies the required cube. Otherwise h has an outgoing quadratic edge. If \(h=1\), that edge gives 1, matching the cube. For other h values, avoiding 4 forces h to be 0 or 2 and the path endpoints to 2. An original literal h then joins 0 or 2 with the edge value 3 to give 4. Alternatively a positive path ends at an original literal (value 2, again joining to 4) or at an original long word (value 4). Thus a non-1 active h forces the right sum to be 4. Its cube absorbs every word headed h, including literals. Original literal occurrences are immediate. The criterion is exact.

Laws 2-7 satisfy this criterion directly: a cubic head permits arbitrary quadratic and literal recovery; law 4 composes edges; laws 5-7 are the root-literal, endpoint-literal and endpoint-cubic activation cases.

\par\smallskip\noindent\textit{Unrestricted derivation.}\quad 
Normalize only long original words, retaining the original literals and quadratic edges. Compose a positive path with law 4 to obtain ht. Never seed hh from a path of length zero. If an active h is not already an original cubic head, use law 5 for an original literal h and one outgoing edge. Otherwise follow a positive path to an original literal or cubic head and apply law 6 or 7. This derives hhh. Recover any required literal or quadratic using law 3 or 2, respectively. A long target is already the derived cube after normalization. Every substitution is nonempty.

Every valid absorption is obtained this way. Distributivity and summandwise absorption in both directions give all polynomial identities, proving the finite relative basis without any bound on variables or polynomial size.

Boundary examples: \(xx\leq xy\) fails at \(x=0\), \(y=2\), so paths cannot be made reflexive. \(xxx\leq xy+yx\) fails at \(x=y=2\), so a cycle alone does not activate a cubic head. The law \(x\leq xx\) also fails at \(x=2\). These distinguish this criterion from both global reachability in 13020 and contained witnesses in the preceding sources. Equal-variety transfers require full HSP both ways.

\subsection{Five complete identities for semiring 12126}\label{case:278}
Signature (+, multiplication), no constants. Relative to the additively idempotent semiring axioms the complete finite basis is

{\small\par\smallskip
\CaseLine{\(xyz=zz; xy\leq yy; x\leq xx; yy\leq y+xy; zz\leq xy+yz.\)}
\smallskip}

Addition rows:00044/\allowbreak{}01034/\allowbreak{}00244/\allowbreak{}43434/\allowbreak{}44444. Multiplication rows:41444/\allowbreak{}41444/\allowbreak{}41404/\allowbreak{}41444/\allowbreak{}41444. The order has \(1<0<4\),\(2<0<4\), and \(1<3<4\);\(2+3=0+3=4\). A quadratic ending in 1 has value 1; the ordered pair(2,\allowbreak{}3) has value 0; all other quadratics have value 4. Every word of \(length\geq 3\) ending in 1 has value 1, and otherwise4:the first two factors give 0,\allowbreak{}1, or 4, whose rows all coincide. The square of t is likewise1 at \(t=1\) and 4 otherwise. This proves the word law \(xyz=zz\) and the first two inequalities. For \(yy\leq y+xy\), the exceptional product occurs at \(x=2\), \(y=3\); then the right sum \(is 3+0=4\). Otherwise xy already bounds yy. For \(zz\leq xy+yz\), if \(z=1\) the square 1 is bounded by \(yz=1\). If z is not 1, yz is 4 except at \(y=2\), \(z=3\); then xy ends in 2 and is 4. All five laws are therefore valid, and are independently checked on all source assignments.

\par\smallskip\noindent\textit{Exact original-polynomial criterion.}\quad 
For a nonempty original polynomial P and a variable t, let L be the variables occurring as original literal summands; N the final variables of original nonlinear summands; and \(H_{t}\) the heads of original quadratics ending at t. Call t active if one of the following holds:

1. An original square or word of \(length\geq 3\) ends at t.

2. \(H_{t}\) intersects N.

3. t belongs to L and \(H_{t}\) is nonempty.

A literal t is absorbed by P iff t is an original literal or t is active. A nonlinear target w ending at t is absorbed iff t is active, or w is itself an original quadratic. These are global endpoint conditions; no content containment is imposed. A literal predecessor guard is not sufficient, whereas a literal at the final variable is. In particular \(y\leq xy\) fails at \(x=2\), \(y=3\), so one cannot recover literals from arbitrary nonlinear same-tail witnesses.

\par\smallskip\noindent\textit{Necessity: separating assignments.}\quad 
Consider a failed absorption. Its final variable t is inactive. If \(H_{t}\) is empty, there is no original nonlinear word ending at t. Assign \(t=3\) and every other variable 1. For a missing literal t there is no literal t either, so \(P=1\) and \(t=3\) is not below it. For a nonlinear target, P is 1 or 3 depending on whether a literal t is present, whereas the target is 4.

If \(H_{t}\) is nonempty, inactivity gives t not in L, no original long word or square ending at t, t not in \(H_{t}\), and \(H_{t}\) disjoint from N. Assign \(t=3\), every variable of \(H_{t}\) the value 2, and all other variables 1. Every original nonlinear word ending at t is a quadratic with head in \(H_{t}\), hence has value 0. Every other original nonlinear word ends outside \(H_{t}\) and outside \{t\}, hence has value 1. Original literals have value 1 or 2. Thus \(P\leq 0\). The literal \(t=3\) is not below 0. A long target ends at 3 and is 4. A missing quadratic ht has h not in \(H_{t}\); its initial value is 1 or 3, so ht is 4. This exhausts failure cases and proves necessity for arbitrary original polynomials.

\par\smallskip\noindent\textit{Complete unrestricted derivations.}\quad 
Normalize every original long word to its final square with \(xyz=zz\), substituting its nonempty prefix for x and its final two letters for y, z. Thus every active variable of case 1 supplies its square. In case 3, choose an original edge at and the original literal t; \(yy\leq y+xy\) produces tt. In case 2 choose an original quadratic at and an original nonlinear word va ending in a, with v nonempty. Substitute \(x=v\), \(y=a\), \(z=t\) in \(zz\leq xy+yz\). Normalization of va is harmless and yields tt under the unchanged original sum.

Once tt is available, \(xy\leq yy\) recovers every quadratic ht; reversing \(xyz=zz\) recovers every long word ending in t; \(x\leq xx\) recovers the literal t. Original literal or quadratic targets need no inference. The necessary conditions above exhaust all valid absorptions, so every valid absorption is derived using these five laws with nonempty images. Distribute arbitrary terms into nonempty polynomials and use these derivations for each summand in both directions. This proves finite basability without any bound on variables, word lengths, or summands. Opposite and catalogue transfers require independently verified HSP membership in both directions; no short matching test is a proof.

\subsection{Two complete finite bases from original head-tail triggers}\label{case:279}
Signature (+, multiplication), without constants. Every basis below is relative to the additively idempotent semiring axioms. The proof allows arbitrary variables, word lengths, and polynomial sizes.

\par\smallskip\noindent\textit{Tables, values, and normal words.}\quad 
13417 addition: 00444 / 01444 / 44224 / 44234 / 44444. Multiplication: 44404 / 11111 / 44444 / 44434 / 44444. The addition order consists of \(1<0<4\) and \(3<2<4\), with no other comparisons. A nonlinear word with head value 1 has value 1. With head 0 or 3, it retains that head value exactly when all tail values are 3; otherwise its value is 4. Heads 2 and 4 always give 4 in a nonlinear word.

13713 addition: 00044 / 01144 / 01244 / 44434 / 44444. Multiplication: 44414 / 44414 / 22222 / 44434 / 44444. The order is \(2<1<0<4\), with 3 incomparable with 0,\allowbreak{}1,\allowbreak{}2. A nonlinear word headed 2 has value 2, and one headed 4 has value 4. For heads 0,\allowbreak{}1,\allowbreak{}3, all tails must be 3 to avoid 4; if they are, the values are 1,\allowbreak{}1,\allowbreak{}3.

Literals retain their assigned value. Direct induction on each table proves these formulas. Both algebras satisfy \(xyz=xzy\) and \(xyy=xy\). Thus a word is represented by its head and the set of actual tail variables, retaining an occurrence of the head in the tail. Permute and delete tail repetitions by these two laws with nonempty prefixes.

\par\smallskip\noindent\textit{Exact criteria using original right words.}\quad 
For the original polynomial P, let H be its head set, C its content union, N the content union of its nonlinear words, and \(D_{h}\) the union of the actual tails of its h-headed words. Let \(E_{h}\) mean \(D_{h}\) intersects H. In particular \(E_{h}\) requires an original h-word with a genuine tail b and a second original word headed b. Either witness may be the same word.

For 13417:

\textbullet{} Literal h is absorbed exactly when h belongs to H.

\textbullet{} Nonlinear hT is absorbed exactly when h belongs to H, T is contained in N, and either T is contained in \(D_{h}\) or \(E_{h}\) holds.

For 13713:

\textbullet{} Literal h is absorbed exactly when it occurs literally in P or \(E_{h}\) holds.

\textbullet{} Nonlinear hT is absorbed exactly when h belongs to H, T is contained in C, and either T is contained in \(D_{h}\) or \(E_{h}\) holds.

For either algebra, a missing head is separated by \(h=4\) and every other variable equal to the left-zero value (1 for 13417, 2 for 13713). All right values equal that left-zero value, while the target is 4.

For 13417, a target tail outside N is separated by assigning N to 3 and everything else to 2: nonlinear right words give 3, literals give 2 or 3, but the nonlinear target is 4. If \(E_{h}\) fails and a target tail is outside \(D_{h}\), assign \(h=0\), \(D_{h}=3\), and all \(others=1\). This is consistent because \(D_{h}\) misses all original heads, including h. Original h-words give 0 and other heads give 1, whereas the target gives 4.

For 13713, a target tail outside C is separated by \(C=3\) and all \(others=2\). If \(E_{h}\) fails and a target tail is outside \(D_{h}\), assign \(h=0\), \(D_{h}=3\), and \(others=2\). Right h-words give 1 (or 0 if literal), other heads give 2; the target gives 4. For a literal h that is not present and has no trigger, the same assignment gives right sum at most 1 and target 0. These arguments cover every failure of the stated criteria.

For sufficiency in 13417, suppose the right value is not 4. If it lies in \{0,\allowbreak{}1\}, all original heads are in \{0,\allowbreak{}1\}. A target head 1 gives 1. If the target head is 0, original h-tails must all be 3. A trigger would make an original head 3, which is impossible in this case; hence the target tails are in \(D_{h}\) and the target gives 0. When the right value is 1, an original h-word excludes head value 0, so it gives 1. If the right value lies in \{2,\allowbreak{}3\}, every nonlinear original word has head and tails 3, so every variable of N is 3. The nonlinear criterion also supplies an original nonlinear h-word (T is nonempty), and hence the target is 3. Literal h is below its original h-word by the table.

For sufficiency in 13713, again exclude right value 4. Head 2 is matched by an original h-word; head 4 is impossible. If h is 0 or 1, all original h-tails must be 3. A trigger contributes a word headed 3, giving 3 or 4, whose sum with the h-word is 4. Thus no trigger can occur in this case, and the nonlinear target tails must be in \(D_{h}\); it gives 1. A literal under these heads must have been an original literal. If \(h=3\), an original h-word is 3, and the right sum can avoid 4 only by being 3. Every original variable is then 3, so content containment makes the target 3 as well. This proves both unrestricted criteria.

\par\smallskip\noindent\textit{Seven identities for 13417.}\quad 
1. \(xyz=xzy\).

2. \(xyy=xy\).

3. \(xy\leq xyz\).

4. \(xyz\leq xy+xz\).

5. \(x\leq xy\).

6. \(xz\leq xy+y+zt\).

7. \(xz\leq xy+y+tz\).

Their validity follows by applying the criterion to the displayed words. From an original p-word extract each atom pt corresponding to an actual tail t, using law 3 after tail permutation. Extract any original head as a literal with law 5. If all target tails lie in \(D_{h}\), their atoms are available directly. Otherwise choose original trigger witnesses hbB and bA, with B, A possibly absent, and obtain hb and b. Every target tail t is in an original nonlinear word. If it is that word's head, extract an atom tu; otherwise extract ut. Laws 6 and 7 respectively give ht. Merge the required atoms by law 4. Literal targets use law 5. Absent suffixes here denote separate cases, never empty substitutions.

\par\smallskip\noindent\textit{Twelve identities for 13713.}\quad 
1. \(xyz=xzy\).

2. \(xyy=xy\).

3. \(xy\leq xyz\).

4. \(xyz\leq xy+xz\).

5. \(x\leq xy+y\).

6. \(x\leq xy+yz\).

7. \(xt\leq xy+y+t\).

8. \(xt\leq xy+y+tu\).

9. \(xt\leq xy+y+ut\).

10. \(xt\leq xy+yz+t\).

11. \(xt\leq xy+yz+tu\).

12. \(xt\leq xy+yz+ut\).

Each displayed absorption is valid by the 13713 criterion. Direct tail coverage again uses law 3 to extract atoms and law 4 to merge. For a trigger, extract hb from the original h-word and retain the original b-headed witness. If that witness is literal b, use laws 5 and 7-9. Otherwise it is bA with A nonempty, and use laws 6 and 10-12 with \(z=A\). Literal targets are recovered by law 5 or 6. For each nonlinear target tail t, choose an original word containing t. It is either literal t, nonlinear headed t (from which extract tu), or nonlinear with t in its tail (extract ut). The three content cases are exactly the three recovery laws in each guard group. Obtain all ht and merge with law 4. Every substitution is nonempty.

All constructed words are original or have already been derived below P. Every valid absorption is therefore derivable; distributivity and summandwise absorption in both directions give every polynomial identity. These are unrestricted finite bases, not bounded-search conclusions.

The distinctions matter: \(x\leq xy\) is valid in 13417 but fails in 13713 at \(x=0\), \(y=3\). Conversely \(xz\leq xy+y+z\) is valid in 13713, but fails in 13417 at \(x=3\), \(y=3\), \(z=2\). A literal content witness cannot replace a nonlinear content witness in the 13417 construction. Full bidirectional HSP certificates are required for any opposite or equal-variety transfer.

\subsection{Six complete identities for 13764}\label{case:280}
Signature (+, multiplication), without constants, relative to the additively idempotent semiring axioms. The complete basis is:

1. \(xyz=xx\).

2. \(xy\leq xx\).

3. \(x\leq xx\).

4. \(xx\leq x+xy\).

5. \(xx\leq xy+y\).

6. \(xx\leq xy+yz\).

All arguments below allow arbitrary numbers of variables and summands.

Addition: 00444 / 01244 / 42244 / 44434 / 44444. Multiplication: 44424 / 11111 / 44444 / 44444 / 44444. The order has \(1<0\), \(1<2\), and top 4; 0 and 2 are incomparable, and 3 is incomparable with 0,\allowbreak{}1,\allowbreak{}2. A literal retains its value. A quadratic headed 1 is 1; one headed 0 and ending in 3 is 2; all other quadratic values are 4. Every word of degree at least three is 1 when headed 1, and 4 otherwise. The same is true of a square. Induction on the table proves these facts, hence law 1. A word of length at least three normalizes to its head square with the nonempty substitution \(x=first\) letter, \(y=second\) letter, \(z=remaining\) suffix.

For an original polynomial P, let L be its literal variables and K the heads of its squares or words of length at least three. The remaining quadratic words ab (a and b distinct) are directed edges. Let \(D_{h}\) be their endpoints when headed h, and let H be all original heads. Call h active if h is in K, or \(D_{h}\) is nonempty and (h is in L or \(D_{h}\) meets H).

A literal h is absorbed exactly if h is an original literal or active. A square or long word headed h is absorbed exactly if h is active. A quadratic ht with h different from t is absorbed exactly if that original edge is present or h is active. There is no transitive-path condition: any two-edge path already activates its initial head.

Necessity is proved by explicit separating assignments. If h is absent from H, assign \(h=3\) and everything \(else=1\). Every right word is 1, whereas the target is 3 or 4. If h is in H but inactive and \(D_{h}\) is empty, the only possible h-word is literal h. Assign \(h=0\) and \(others=1\): the right sum is 0 and any nonlinear h-target is 4. A literal target in this case is already an original summand, hence is not a failed case. If h is inactive and \(D_{h}\) is nonempty, it is not a literal or square head, \(D_{h}\) misses H and in particular misses h. Assign \(h=0\), \(D_{h}=3\), and other \(variables=1\). All h-edges give 2; other original heads give 1, so the right sum is 2. The target literal gives 0, a square or long target gives 4, and a missing quadratic ht has t outside \(D_{h}\) and gives

4. None is below 2. These cases exhaust every failed criterion.

For sufficiency, an original edge or literal is immediate. For active h in K, the normalized original head square is already available. Otherwise take an original h-edge hb. If the right sum avoids 4 and h is not 1, the table forces \(h=0\) and \(b=3\). A literal h contributes 0, whose join with edge value 2 is 4. If instead an endpoint b is an original head, a literal b contributes 3, or a nonlinear b-word contributes 4; either also forces right sum 4. Thus an active head other than 1 forces 4. For head 1 the original outgoing edge gives 1. Consequently h squared, which has value 1 or 4, is always absorbed. Every h-headed target is below that square by the word formula. This proves the exact criterion.

Laws 2-6 satisfy the criterion directly. For unrestricted completeness, normalize long original words. Original square heads already provide \(h^{2}\). If h is literal with an original outgoing hb, law 4 produces \(h^{2}\). If an original endpoint b is a literal head, law 5 applied to hb+b produces \(h^{2}\). Otherwise choose an original nonlinear b-word bB, with B nonempty. Law 6 with \(x=h\), \(y=b\), \(z=B\) produces \(h^{2}\). This also covers square and long guard words; normalization retains b as their head. Laws 2 and 3 then recover the desired quadratic or literal. All substitutions are nonempty, and the original words may contain variables outside the target.

Every valid absorption has thus been derived. Distributivity and summandwise absorption in both directions prove completeness for all polynomial identities, without finite search assumptions.

Distinctions: \(xx\leq xy\) fails at \(x=0\), \(y=3\); \(x\leq xy\) also fails there. But \(xx\leq xy+yz\) holds: unlike 13349, two successive edges immediately activate the initial head without a terminal literal or cubic guard. Equal-variety and opposite transfers require full HSP both ways.

\subsection{Seven complete identities for semiring 12204}\label{case:281}
Constant-free signature (+, multiplication). Relative to the additively idempotent semiring axioms, a complete finite basis is

{\small\par\smallskip
\CaseLine{\(xy=yx; xyz=ttt; x\leq yyy; xx\leq x+xy; x\leq xx+xy+y;\)}
\CaseLine{\(xz\leq xx+xy+yz; xy\leq xx+x+y+yz.\)}
\smallskip}

Here \(u\leq P\) abbreviates \(u+P=P\). Addition: 00444/\allowbreak{}01244/\allowbreak{}42244/\allowbreak{}44434/\allowbreak{}44444. Multiplication: 42444/\allowbreak{}22444/\allowbreak{}44444/\allowbreak{}44424/\allowbreak{}44444. The additive order has \(1<0<4\), \(1<2<4\), and 3 incomparable with 0,\allowbreak{}1,\allowbreak{}2 below 4. A quadratic is 2 exactly at(0,\allowbreak{}1),(1,\allowbreak{}0),(1,\allowbreak{}1),(3,\allowbreak{}3), and 4 otherwise. All longer words are 4.

Let L be original literals, E original unordered quadratic edges, V their endpoints, S original loop variables. Components use(V, E); active means containing an original literal. An original long word absorbs all targets. Otherwise a literal h is absorbed iff h is in L, or h is in S and its component is active. A pair ht is absorbed iff both endpoints are in V, its edge is original or an endpoint is in S union L, and its components coincide or both are active. No long target is absorbed. Original literals do not propagate to arbitrary vertices, even in active components.

For sufficiency, if \(P<4\) with an edge present, \(P=2\). Every component lies within\{0,\allowbreak{}1\} or\{3\}. An original literal endpoint must be 1, so an active component lies in\{0,\allowbreak{}1\}. A loop in an active component is 1, which proves the literal clause. Any loop or literal guard endpoint is 1 in the first block or 3 in the second. A same-component pair with it is 2; in two active components the guard is 1 and the other endpoint 0 or 1, again giving 2. Original edges have value 2. With no edges, only original literals satisfy the criterion.

For necessity, assign all variables 1 to separate a long target. For a missing literal h without an original loop, if it lies in V assign \(h=0\) and all other variables 1:there is no loop at h, all edges 2, all original literals 1, and \(h=0\) is not below 2. If h is outside V, assign it 4 and all others 1. For a missing literal with an original loop in an inactive component, assign that component 3 and all else 1. Every right edge is 2 and every literal 1, whereas \(h=3\) is not below 2. For a pair with an endpoint outside V, assign V to 1 and outsiders 2. For a pair with no original edge and neither endpoint in S union L, assign target endpoints 0 and all else 1. For a pair in distinct components not both active, assign an inactive target component 3 and all else 1. Each case makes \(P\leq 2\) and the target 4. These assignments cover all failures for arbitrary variable sets and polynomial sizes.

All seven displayed laws satisfy the criterion and table formula; every finite law assignment is independently checked as well. An original long word supplies a fixed cube by \(xyz=ttt\) with a nonempty prefix, and \(x\leq yyy\) recovers any target. Otherwise original targets need no step. If h is an original literal with an incident edge, derive hh by \(xx\leq x+xy\). A known square hh allows pair propagation along an original path using \(xz\leq xx+xy+yz\). To recover a new literal h, it has an original loop and an original path to an original literal z. First derive hz by square-based path propagation, then use \(x\leq xx+xy+y\) with \(x=h\), \(y=z\). Only original loops justify this new literal.

For a new pair, orient an original loop/literal guard as h. Its square is original or obtained by \(xx\leq x+xy\). In one component, propagate from h to the other endpoint along a simple original path. For two active components, first recover literal h if needed. Choose an original literal z in the other component and its original incident edge za. The law \(xy\leq xx+x+y+yz\) derives hz. Retain hh and propagate from this bridge along an original path from z to the target endpoint. Commutativity supplies all orientations. Original witness sets never change, and every intermediate term is explicitly below P.

All substitutions are nonempty. These finitely many fixed laws derive every valid absorption without bounds on variables, lengths, or summand counts. Distributive normalization and summandwise absorption in both directions prove all identities. Hence the seven laws form a complete finite basis. Other catalogue types require full HSP in both directions before transferring the result.

\subsection{Seven complete identities for semiring 12205}\label{case:282}
Constant-free signature (+, multiplication). A finite basis, relative to the additively idempotent semiring axioms, is

{\small\par\smallskip
\CaseLine{\(xy=yx; xyz=ttt; x\leq yyy; x\leq y+xy;\)}
\CaseLine{\(xx\leq x+xy; xz\leq xx+xy+yz; xy\leq xx+x+y+yz.\)}
\smallskip}

As usual, \(u\leq P\) abbreviates \(u+P=P\). Addition: 00444/\allowbreak{}01444/\allowbreak{}44224/\allowbreak{}44234/\allowbreak{}44444. Multiplication: 42444/\allowbreak{}22444/\allowbreak{}44444/\allowbreak{}44424/\allowbreak{}44444. The additive order consists of \(1<0<4\) and \(3<2<4\). A quadratic is 2 exactly at(0,\allowbreak{}1),(1,\allowbreak{}0),(1,\allowbreak{}1),(3,\allowbreak{}3), and 4 otherwise. All longer words are 4.

\par\smallskip\noindent\textit{Exact criterion in terms of original components.}\quad 
An original long word in P absorbs everything. Otherwise let L be original literals, E original unordered quadratic edges, V their endpoints, and S original loop variables. A component of(V, E) is active when it contains an original literal. A literal h is absorbed iff it is original or belongs to an active component. A pair ht is absorbed iff h, t are in V, and either both components are active, or they are in the same component and their edge is original or an endpoint belongs to S. No long target is absorbed. The active-component alternative permits new squares; it is stronger than the 14959 criterion.

For sufficiency, if P is below 4 and contains an edge, its value is 2. Every component lies in one of\{0,\allowbreak{}1\} and\{3\}. A literal endpoint must equal 3, so an active component consists entirely of 3 values. All literals in it are below 2, and all pairs between active components are 2. An original loop forces 1 in the first block or 3 in the second; a pair with that endpoint in the same component is therefore 2. An original edge is 2. With no edges, only original literals qualify.

For necessity, all variables 3 separate a long target. A missing literal outside V is separated by assigning it 4 and all others 3; it occurs in no right term. A missing literal in an inactive component is separated by giving that component 1 and all other variables 3. Every right edge is 2 and every literal 3, while the target 1 is not below 2. For a pair with an endpoint outside V, assign V to 3 and outsiders 2; P is at most 2 and the target 4. For endpoints in different components not both active, assign an inactive target component 1 and everything else 3. Then \(P=2\) and the target 4. The remaining failed case is a common inactive component with no target edge and no endpoint loop. Assign the target endpoints 0, all other vertices of that component 1, and everything else 3. No edge has two0 endpoints, and no literal lies in that inactive component. Thus \(P=2\) and the target 4. These assignments cover every failure for arbitrary variable sets and polynomial sizes.

\par\smallskip\noindent\textit{Unrestricted completeness.}\quad 
Every displayed law is valid by the exact criterion and word-value formula; independent checking also covers every assignment to each law. Use \(xyz=ttt\) with a nonempty prefix to reduce any original long word to a fixed cube, then \(x\leq yyy\) supplies any target. Otherwise propagate literals along original paths from an original literal by \(x\leq y+xy\). For a new pair, orient an original loop endpoint as h if one exists. If neither endpoint has an original loop, both components must be active. Recover h as a literal, then obtain hh via \(xx\leq x+xy\) using an original incident edge. This already derives a square target.

For distinct endpoints in one component, retain the available square hh and traverse a simple original path via \(xz\leq xx+xy+yz\). For two different active components, derive both endpoint literals, choose an original edge tz at t, and use \(xy\leq xx+x+y+yz\). Commutativity orients all edges and restores target order. Every witness is originally present or has been explicitly derived below P; the semantic sets always remain those of the original polynomial.

These steps derive all valid absorptions with nonempty substitutions, without bounds on variables, lengths, or summand counts. Distributivity and absorption in both directions establish every identity. Hence the seven laws are a complete finite basis. Transfers to further types require independently verified complete HSP in both directions.

\subsection{Five complete identities for semiring 12206}\label{case:283}
In the constant-free signature (+, multiplication), relative to the additively idempotent semiring axioms, a complete finite basis is

{\small\par\smallskip
\CaseLine{\(xy=yx; xyz=ttt; x\leq yyy; x\leq y+yz; xz\leq xx+xy+yz.\)}
\smallskip}

Here \(u\leq P\) means \(u+P=P\). Addition: 00444/\allowbreak{}01444/\allowbreak{}44244/\allowbreak{}44434/\allowbreak{}44444. Multiplication: 42444/\allowbreak{}22444/\allowbreak{}44444/\allowbreak{}44424/\allowbreak{}44444. The additive order has \(1<0<4\), with 2 and 3 incomparable with these below 4. Quadratics are 2 at (0,\allowbreak{}1),(1,\allowbreak{}0),(1,\allowbreak{}1),(3,\allowbreak{}3), and 4 otherwise. All longer words are 4.

For original P, let L be its literal variables, E unordered quadratic edges, V their endpoints, S original loops. If P has a long word, or L intersects V, it has value 4 on every assignment and absorbs all targets. In the latter case, an edge value below 4 is 2; its incident literal is 0,\allowbreak{}1, or 3, whose join with 2 is 4. Otherwise a literal h is absorbed iff h is in L; a pair ht iff its endpoints are in the same component of(V, E) and its edge is original or an endpoint in S; no long target is absorbed.

For sufficiency when \(P<4\), all right edges are 2 and right literals 2. Each component lies in one of\{0,\allowbreak{}1\} and\{3\}; a loop endpoint is 1 in the first block or 3 in the second. Its pair with a vertex in the same component is 2. Original edges are 2. The stated targets follow. For necessity in the non-top case, assign V to 1 and outsiders 2: every right term is 2 and a long target is 4. A missing literal in V has value 1, not below 2; for a missing literal outside V change its value to 4. A pair with an outside endpoint is 4 under the original assignment. For a pair in different components, assign one component 3, other V vertices 1, and outsiders 2. For a same-component pair with no edge or endpoint loop, assign target endpoints 0, other V vertices 1, and outsiders 2. Again every right term is 2 and the target 4. These exhaust failure of the criterion without variable or size bounds.

All displayed laws satisfy this criterion or the table formula, and their assignments are independently checked. For a right long word, reduce it to a fixed cube using \(xyz=ttt\) with a nonempty prefix, then use \(x\leq yyy\). For original literal y incident with original edge yz, use \(x\leq y+yz\) with the entire target as x. Otherwise original targets need no step. For a new pair orient its original loop endpoint as h, choose a simple original path to the other endpoint, and repeatedly use \(xz\leq xx+xy+yz\) while retaining the original loop hh. Commutativity orients edges. Every substituted word is nonempty, and every intermediate term is below P. This proves all absorptions without bounds on variables, lengths, or summand counts. Distributivity and summandwise absorption in both directions prove every identity, establishing the five-law basis. Further transfers require full HSP in both directions.

\subsection{Complete bases for three quadratic endpoint constraints}\label{case:284}
Work in the signature (+, multiplication) without constants. Write \(u\leq P\) for \(u+P=P\). Relative to the additively idempotent semiring axioms, each of 12208,\allowbreak{}14369,\allowbreak{}14772 has the following seven laws:

{\small\par\smallskip
\CaseLine{\(xy=yx; xyz=xyzz; xxx\leq xyz; xxy\leq xxx+yyy;\)}
\CaseLine{\(yyy\leq xxx+xy; xy\leq xxy; xz\leq xx+xy+yz.\)}
\smallskip}

For 12208 add \(x\leq xxx\) and \(xxx\leq x+xy\) (nine laws in total). For 14369 add \(x\leq xx\) and \(xx\leq x+xy\) (nine laws in total). For 14772 add \(x\leq xy\) (eight laws in total).

The proof below is unrestricted in variables, word lengths, and the number of polynomial summands. No finite test is used as completeness.

\par\smallskip\noindent\textit{Word values from the operation tables.}\quad 
Addition and multiplication rows are respectively:

12208: 00444/\allowbreak{}01444/\allowbreak{}44224/\allowbreak{}44234/\allowbreak{}44444; 42444/\allowbreak{}22444/\allowbreak{}44444/\allowbreak{}44434/\allowbreak{}44444. 14369: 00044/\allowbreak{}01034/\allowbreak{}00244/\allowbreak{}43434/\allowbreak{}44444; 44444/\allowbreak{}40404/\allowbreak{}44244/\allowbreak{}40444/\allowbreak{}44444. 14772: 00004/\allowbreak{}01004/\allowbreak{}00224/\allowbreak{}00234/\allowbreak{}44444; 44444/\allowbreak{}41444/\allowbreak{}44404/\allowbreak{}44004/\allowbreak{}44444.

Let (s, b, g, q) be (3,\allowbreak{}0,\allowbreak{}1,\allowbreak{}2), (2,\allowbreak{}3,\allowbreak{}1,\allowbreak{}0), and (1,\allowbreak{}2,\allowbreak{}3,\allowbreak{}0), respectively. A literal retains its assigned value. A nonlinear word has value s if all its variables have value s. A quadratic has value q if both occurrences are in \{b, g\} and at least one occurrence is g. Every other nonlinear word has value 4. In particular every long word (length at least three) has value s or 4. This follows directly for length two; multiplication of q by any element is 4, multiplication of s by an element other than s is 4, and 4 is multiplicatively absorbing. Thus the formula holds at every length by induction.

Commutativity and \(xyz=xyzz\) normalize every long word by its content. Use sorted distinct content when it has at least three variables, xxy for two variables \(x<y\), and xxx for a singleton. Move repeated letters to the end and insert/delete copies while retaining at least three positions. Quadratics preserve their two occurrences; literals preserve length one.

\par\smallskip\noindent\textit{Exact absorption criterion.}\quad 
All data in this paragraph refer to original right summands P. Let C be their content union, L their literal variables, Q their quadratic multigraph (loops allowed), V its vertices, and K the vertices having original loops. Start the active set A with all variables in original long words, and close it across original Q edges. For 12208 also start with L intersect V. For 14369 and 14772 there are no such extra seeds. Define \(B=K\) for 12208 and 14772, and \(B=K\) union (L intersect V) for 14369.

An original word, up to the normalization above, is always absorbed. For any other target w the complete criterion is:

\textbullet{} A literal x is absorbed if x belongs to A (12208), A union K (14369), or C (14772).

\textbullet{} A long word is absorbed if and only if its content is contained in A.

\textbullet{} A quadratic xy is absorbed if both endpoints are active, or the endpoints lie in one Q component and at least one endpoint is in B. A vertex counts as in a Q component only if it belongs to V. For a missing square xx the latter alternative occurs only in 14369, where x is an original literal incident with an original quadratic. That square is recovered directly, without assuming a positive loop path.

For sufficiency, suppose P has value different from 4. If P has a nonlinear term, its value is at most q. Every Q component is either entirely s, or entirely in \{b, g\}, with the b-valued vertices forming an independent set and every loop forcing its endpoint to g. Every original long word forces its content to s. In 12208 an incident literal cannot have value b or g when the sum is at most q, so its Q component is forced to s. In 14369 an incident literal permits only s or g; in 14772 it imposes no further restriction on the component.

Consequently active vertices have value s. A quadratic with an anchored endpoint has value s or q, while a quadratic with both active endpoints has value s. For 14369 an original loop gives a literal value s or g, which is below its loop value. For 14772 any nonlinear word dominates each of its letters. The stated literal recoveries follow. When all nonlinear right words have value s, their sum can remain s; every inferred target in this case also has value s unless an original literal already dominates it. These observations prove sufficiency also for sums below q. If P has no nonlinear term, the criterion adds no new long or quadratic target and no missing literal.

Here are explicit separating assignments proving necessity. Initially assign every variable s, including variables outside P. A component reassigned below is always inactive and therefore contains no original long variable. In 12208 it also contains no original literal endpoint.

For a long target choose a nonactive variable x in it. If x belongs to V, assign its whole Q component g. If x is outside V, assign it q. The target becomes 4 and every right summand is at most q. In the outside-V case x occurs in no original long word and possibly only as a literal. No other right summand is disturbed.

For a missing quadratic not covered by the criterion, first consider an endpoint x outside V and A. Assign \(x=q\), keeping everything else s. The target is 4 and the right side is at most q. Otherwise, if the two endpoints lie in different components, assign an inactive endpoint's entire component g and leave the other endpoint s. The target is 4. If they lie in the same inactive component, give its target endpoints value b and all other vertices there value g. There is no original edge joining the two endpoints and no original loop at either: these would give an original target or an anchor. For 14369 neither endpoint is an original literal. Thus all original edges still have value q, all literals are at most q, and the target is 4. The same construction handles a missing square by assigning its sole endpoint b. All unmentioned variables retain s.

For a missing literal in 12208, if x lies in an inactive Q component assign that component g; \(x=g\) is not below q, whereas all right terms are at most q. If x lies outside V and A, assign \(x=b\). It is not an original literal and occurs in no nonlinear right word, so the right side stays s and does not dominate b. For 14369, if x is incident with Q, assign \(x=b\) and all other vertices in its inactive component g. There is no original loop or literal at x, so all right terms are at most q while b is not below q. If x lies outside V and A, assign it b and the right side stays s. For 14772 a missing literal fails exactly when x is outside C: assign it 4 and all right variables s. These assignments prove all failed cases for arbitrarily large polynomials.

\par\smallskip\noindent\textit{Unrestricted equational derivation.}\quad 
Retain all original summands throughout, adjoining only terms proved below their sum. Obtain xxx for each variable x of an original long word using \(xxx\leq xyz\): after commuting one occurrence of x to the front, split the remaining at least two letters into two nonempty blocks. For 12208 obtain xxx at each original literal/quadratic incidence from \(xxx\leq x+xy\). Propagate cubes along original edges by \(yyy\leq xxx+xy\). This produces exactly the required active cubes; propagation uses only original edges and retains its original seeds.

Merge any finite number of these cubes using \(xxy\leq xxx+yyy\), substituting the accumulated nonempty content word for x and the next letter for y. The long normal form of the accumulated word's cube is the previously derived term. Thus every long target with content in A is obtained. Recover its quadratic counterpart by \(xy\leq xxy\). For 12208 recover an active literal using \(x\leq xxx\). For 14369 first recover its square with \(xy\leq xxy\), then use \(x\leq xx\). An original loop in 14369 already gives its literal. For 14772 recover a content literal from an original containing word with \(x\leq xy\), separating one occurrence and its nonempty remainder; an original literal requires no recovery.

For an inactive quadratic with an anchored endpoint x, keep xx. It is an original loop in 12208 and 14772. In 14369 either it is original or derive it from the original literal x and an original incident edge xy by \(xx\leq x+xy\). In particular this gives a missing square immediately. For distinct endpoints choose a positive path \(x=p0\), p1,..., \(pn=y\) of original Q edges. The first edge xp1 is original. Successively use

{\small\par\smallskip
\CaseLine{\(x p_{i+1} \leq  xx + x p_{i} + p_{i} p_{i+1}\)}
\smallskip}

to extend it, retaining the loop and the next original edge. This is \(xz\leq xx+xy+yz\) with singleton nonempty substitutions. It yields xy. No reflexive step is used as a substitute for a positive edge path.

Every substitution above is nonempty; the number of steps may grow with the target and the original graph, but the listed bases are fixed. Each displayed law satisfies the criterion (or the word formula), hence is valid in its source. Distributive expansion reduces every identity to absorptions in both directions, so this derives all identities.

The distinctions are essential. In 12208 literal/edge incidence produces a cube; in 14369 it produces only a square; in 14772 it need produce neither. In 14369 and 14772 an inactive component may be arbitrarily large while producing no new long words. The independent-set condition on bad endpoints prevents replacing the anchored path condition with mere connectivity.

\subsection{Complete finite bases for 13352 and 13375}\label{case:285}
Constant-free signature (+, multiplication), relative to the additively idempotent semiring axioms. Both sources satisfy

{\small\par\smallskip
\CaseLine{\(xyz=ttt; x\leq yyy; x\leq y+yz; x\leq y+zy.\)}
\smallskip}

For 13352 add two laws, giving six:

{\small\par\smallskip
\CaseLine{\(xy\leq xx+yz; xy\leq xx+zy.\)}
\smallskip}

For 13375 instead add four laws, giving eight:

{\small\par\smallskip
\CaseLine{\(xy\leq xx+zy; xy\leq yy+xz; xy\leq yy+zx; xy\leq yx+zy.\)}
\smallskip}

Here \(u\leq P\) abbreviates \(u+P=P\). Both addition tables are 00044/\allowbreak{}01144/\allowbreak{}01244/\allowbreak{}44434/\allowbreak{}44444, with \(2<1<0<4\) and 3 incomparable with 0,\allowbreak{}1,\allowbreak{}2 below 4. Multiplication 13352: 44344/\allowbreak{}33344/\allowbreak{}33344/\allowbreak{}44444/\allowbreak{}44444. Multiplication 13375: 44344/\allowbreak{}44344/\allowbreak{}43344/\allowbreak{}44444/\allowbreak{}44444. All long words are 4. Quadratic values are 3 or 4. For 13352, the 3-pairs are all pairs in \{0,\allowbreak{}1,\allowbreak{}2\} squared except(0,\allowbreak{}0),(0,\allowbreak{}1). For 13375 they are precisely (0,\allowbreak{}2),(1,\allowbreak{}2),(2,\allowbreak{}1),(2,\allowbreak{}2). Neither multiplication is commutative.

\par\smallskip\noindent\textit{Exact semantics for 13352.}\quad 
Let L be original literals, E original directed edges, V endpoints, S original loops, T original final variables. An original long word, or L intersecting V, makes P universally 4. Otherwise a literal is absorbed iff original; a pair ht iff h, t are in V and either ht is an original directed edge or h belongs to S; no long target is absorbed. The guard is a loop at the TARGET HEAD, not at either endpoint.

Below 4, all right edges and literals are 3. An original loop has value 1 or 2, each giving 3 when multiplied on the right by any endpoint value 0,\allowbreak{}1,\allowbreak{}2. This proves sufficiency. For necessity in the non-top case, assign V to 2 and outsiders 3. Every right term is 3; long targets are 4, missing endpoint literals are 2, and outside-endpoint pairs are 4. For a missing literal outside V, set it 4. For a nonoriginal pair ht with no loop at h, assign \(h=0\), \(t=1\) when distinct, other V vertices 2, and outsiders 3. If \(h=t\), assign \(h=0\) instead. Every original outgoing edge at h has tail 2, as ht and hh are absent. Every other right edge is 3, while the target 4. These cases exhaust criterion failure.

For completeness, an original long word supplies a fixed cube using \(xyz=ttt\) with a nonempty prefix, then \(x\leq yyy\) recovers any target. An original incident literal y uses \(x\leq y+yz\) or \(x\leq y+zy\) according to the edge orientation. Otherwise a new pair ht has an original loop hh and an original edge incident with t. The corresponding one of \(xy\leq xx+yz\) and \(xy\leq xx+zy\) recovers it. All witnesses original.

\par\smallskip\noindent\textit{Exact semantics for 13375.}\quad 
The universal-top cases are the same. Otherwise a literal is absorbed iff original; a pair ht iff h, t are in V, t belongs to T, and either ht or th is original, or an endpoint belongs to S. No long target is absorbed. The original final-variable condition cannot be omitted.

Below 4, edge endpoints are 0,\allowbreak{}1,\allowbreak{}2 and original final endpoints 1 or 2. An original loop forces 2. A head-loop pair with any original final is 3; a tail-loop pair with any endpoint is 3. An original edge is 3. For an original reverse edge th, the target tail t is an original final elsewhere, hence nonzero. If th is 3 and t is nonzero, then h is also nonzero and one of t, h equals 2, so ht is 3. This proves sufficiency. For necessity use the same endpoint, missing-literal, and long-word separators as above. If t is in V but not T, assign \(t=0\), all other V vertices 2, and outsiders 3. There is no incoming edge at t, so every right edge is 3, but the target with tail 0 is 4. If the underlying target edge is absent and neither endpoint loops, assign both target endpoints 1, other V vertices 2, and outsiders 3. No right edge has two1 endpoints, so right terms are 3 and target 4. These exhaust all failures for arbitrary variable sets and sizes.

For completeness, the top cases use the four common laws. Otherwise if only the reversed target edge th is original, choose an original incoming edge at t and use \(xy\leq yx+zy\). If the head loops, choose an original incoming edge at the target tail and use \(xy\leq xx+zy\). If the tail loops, choose an original edge incident with the head and use \(xy\leq yy+xz\) or \(xy\leq yy+zx\). Edge reversal is guarded by the original tail condition; commutativity is never assumed.

\par\smallskip\noindent\textit{Unrestricted conclusion.}\quad 
Every listed law satisfies its source's criterion and word-value formula; independent checking also tests every finite law assignment. Every valid absorption is explicitly derived from the corresponding fixed laws, with nonempty substitutions and without bounds on variables, lengths, or summand counts. Distributivity and summandwise absorption in both directions therefore prove all identities. The bases have six and eight laws respectively. Transfers to the opposites and other catalogue types require independently checked complete HSP in both directions, with reversed law words when using the opposite source.

\subsection{Nonfinite bases for 12438 and 13708}\label{case:286}
The signature consists of addition and multiplication without constants. The additive order of 12438 is the chain 3 \ensuremath{<} 2 \ensuremath{<} 1 \ensuremath{<} 0 \ensuremath{<} 4. Its multiplication rows are 44004, 44114, 41234, 11234, 44444. The table for 13708 is its transpose, with the same addition.

\par\smallskip\noindent\textit{Word functions and absorption.}\quad 
The element 4 is a multiplicative zero. A word containing two entries from \{0,\allowbreak{}1\} gives 4, even if the entries are separated by factors from \{2,\allowbreak{}3\}. On \{2,\allowbreak{}3\} the product is the last entry. A unique 1 with all other entries in \{2,\allowbreak{}3\} gives 1. For a unique 0 surrounded by entries from \{2,\allowbreak{}3\}, the value is 0 if the 0 is initial, 4 if its immediately preceding value is 2, and 1 if its immediately preceding value is 3. Subsequent neutral entries preserve each of these three possible results. These statements follow from the five rows by induction.

For a polynomial P of nonempty original words, let I consist of the variables repeated in at least one word of P. Write w \ensuremath{\leq} P for \(w+P=P\). This holds exactly under the following conditions.

1. Every letter of w occurs in P. Its last letter is final in an original word. Its first letter belongs to I or is initial in an original word.

2. Every repeated letter of w belongs to I. Every two distinct letters of w outside I occur together in an original word.

3. For each adjacent pair yx in w with x outside I, the consecutive factor yx occurs in an original word.

To force content give an absent variable 4 and all others 3. To force the final letter give it 2 and all others 3. For the initial condition give the initial variable 0 and all others 3. For a missing repetition or co-occurrence give the selected variable or pair the value 1 and all remaining variables 3. A missing adjacency into an unlooped x is separated by setting x to 0, its predecessor y to 2, and all other variables to 3. Here y and x are distinct after the repetition test.

Conversely a target attaining 4 on an entry of value 4 uses a content witness. For two special entries choose an original with two such entries by the repetition or pair condition. When the only special entry is 0 preceded by 2, use an original repeat of its variable or the witnessed consecutive factor. A target of value 0 starts with 0 and uses a repeat or an original with that head. A target of value 1 has a special letter and any word containing that letter gives at least 1. For a neutral target of value 2 use the final witness. If this witness contains any nonneutral value it also gives at least 2. This proves the equivalence on arbitrary alphabets.

\par\smallskip\noindent\textit{A family forcing a long linear word.}\quad 
Use n distinct main variables x1,..., xn, for n at least 4. For each ordered index pair \(i<j\) write one word. Every lower-case marker mentioned below is new, occurring once in the entire polynomial.

\textbullet{} For \(j=i+1\) use a \(x_{i}\) \(x_{j}\) b, omitting the initial marker for \(i=1\) and the final marker for \(j=n\).

\textbullet{} For \(j>i+1\) use a \(x_{i}\) c \(x_{j}\) b with all three markers.

Let \(P_{n}\) be the sum of these words and \(W_{n}=x1\cdots xn\). Each pair of main letters occurs together, each consecutive main pair occurs as a factor, x1 is initial and xn is final. The criterion gives \(W_{n}\) \ensuremath{\leq} \(P_{n}\) at every n.

Consider any word v absorbed by \(P_{n}\). It is linear because no original word repeats a variable. Its letters occur in \(P_{n}\) and are pairwise present together. If v contains a private marker, all its letters therefore belong to the single original word containing that marker. In this case v has length at most five. Without a marker it must begin at x1, end at xn, and follow only the edges \(x_{i}\) to \(x_{i+1}\). Hence v is exactly \(W_{n}\). This dichotomy holds for all words, including words on a larger alphabet.

\par\smallskip\noindent\textit{Obstruction to equational deductions.}\quad 
Fix k and choose n \ensuremath{>} 5(k+1), n at least 4. Distributivity and the semiring axioms identify terms with finite nonempty sets of nonempty words. A substitution sends each variable to a nonempty polynomial, and occurrences expand independently. A single rewrite by a k-variable identity acts in a one-hole term context. After distribution the context contributes unchanged summands and terms of the form L times the substituted polynomial times R, where L and R are fixed words that may be absent. This statement follows by induction on the one-hole context, taking one expanded context pair at a time. Absent contexts introduce no empty-word substitution.

If \(P_{n}\) could be rewritten to \(P_{n}\)+\(W_{n}\) by identities valid in 12438 on at most k variables, take the first step introducing \(W_{n}\). All intermediate polynomials are equal to \(P_{n}\) in 12438, so each of their summands is absorbed by \(P_{n}\). Before this step every summand has length at most five. Fix the context pair that produced \(W_{n}\), the producing monomial of the applied identity, and each selected substitution word for its variables. The monomial is linear because \(W_{n}\) is linear and every substituted word is nonempty. It thus has at most k letters. Every one of its variables occurs in the other side of the applied identity by the content condition. Using its selected word at an occurrence on that side produces an old context summand of length at most five. Thus each selected word has length at most five. The total length of L and R is also at most five, by any old expansion with this same context pair. The newly produced length is at most 5k+5, a contradiction. This handles both orientations, arbitrary polynomial substitutions and semiring reasoning. No finite identity basis exists. Word reversal proves the same for 13708.

For a transfer to another table, inclusion of the obstruction algebra in its variety and validity of the displayed identities for unbounded n together suffice. The second condition can be checked on a finite graph. Fix a principal additive bound b. Existence of private marker values making an original word at most b is a binary relation on its two main values. Use four such relations, according to whether the pair is the first adjacency, an internal adjacency, the last adjacency, or nonconsecutive. Markers are private, so choices for different pairs are independent. Grow a sequence of main values, recording its product, last value, the set of earlier values and its length capped at three. Each added value must pass the nonconsecutive test with every earlier value except the previous one. The first adjacency uses the first relation. Internal adjacencies use the internal relation. A separate final extension after at least three values uses the last relation and tests whether the final product exceeds the principal bound (in a partial order this means it lies outside the ideal). Reachability from the five one-letter states exhausts all finite main sequences. A successful transfer requires no reachable failure for any bound.

\subsection{链序尾部换首：五条完整基}\label{case:287}
\par\smallskip\noindent\textit{1. 运算构造及词值.}\quad 
取五元素链 \(a<e<b<c<t\)，加法为最大值。按此顺序，它们对应目录标号 \((3,\allowbreak{}2,\allowbreak{}1,\allowbreak{}0,\allowbreak{}4)\)。定义右平移如下：右乘 \(e\) 为恒等映射，右乘 \(a\) 为

\CaseDisplay{P(a)=a,\quad P(e)=a,\quad P(b)=b,\quad P(c)=b,\quad P(t)=t,}

右乘 \(b, c, t\) 均为常值 \(t\) 映射。

三个右平移单调、两两交换且幂等，\(P\le\operatorname{id}\le T\)。按照右因子的链序 \(a<e<b<c<t\)，右平移依次为 \(P,\operatorname{id},T, T, T\)，所以两侧分配律成立。逐类核对有 \(R_{xy}=R_yR_x\)，从而乘法结合。

对词 \(w\)，若尾部出现取值属于 \(\{b, c, t\}\) 的变量，则词值为 \(t\)。否则尾部只取 \(a, e\)：尾部全为 \(e\) 时词值为首字母值，尾部含 \(a\) 时词值为首字母值的 \(P\) 像。空尾部按全为 \(e\) 解释。特别地，所有变量只取 \(a, e\) 时，词值就是其全部变量值的最小值。

此半环满足

\CaseDisplay{xyz\approx xzy,\qquad xy^2\approx xy,\qquad x^2y\approx y^2x. \quad\textup{(\ensuremath{1})}}

前两式使词由首字母 \(h(w)\) 和尾部集合 \(D(w)\) 正规化。第三式允许首字母重复的词交换首字母；因此首字母重复的词只由内容决定，记为 \(M_X\)，其中 \(X\) 是非空内容集。具体地，先用前两式写为 \(h^2\) 接其余变量，再用 \(h^2g\approx g^2h\) 换首，故相同内容的这种词相等。于是 \(w^2=M_{c(w)}\)，并且 \(M_X\) 幂等。

\par\smallskip\noindent\textit{2. 精确吸收判据.}\quad 
对非空右多项式 \(F\) 和非空集合 \(X\)，记

\CaseDisplay{F_X=\{u\in F:c(u)\subseteq X\},\qquad E=\bigcup_{u\in F}D(u).}

定理1。 对任意词 \(w\)，令 \(h=h(w)\)、\(X=c(w)\)，则

\CaseDisplay{w\le F\quad\Longleftrightarrow\quad
F_X\ne\varnothing,\qquad D(w)\subseteq E,\qquad
\bigl[h\in E\ \text{或}\ \exists r\in F_X:h(r)=h\bigr]. \quad\textup{(\ensuremath{2})}}

必要性。 若 \(F_X\) 为空，赋 \(X\) 内变量 \(e\)、外部变量 \(a\)。左值为 \(e\)；每个右项含外部变量，且所有变量均取 \(a, e\)，所以右和为 \(a\)。

若 \(x\in D(w)\setminus E\)，赋 \(x=b\)、其余变量 \(e\)。全部右尾变量取 \(e\)，右首字母至多取 \(b\)，故右和不超过 \(b\)；左词尾部含 \(b\)，值为 \(t\)。

若第三条件失败，赋 \(h=c\)、其余 \(X\) 内变量 \(e\)、外部变量 \(a\)。因为 \(h\notin E\)，全部右尾变量只取 \(a, e\)。以 \(h\) 开头的右项不是内部项，必有取 \(a\) 的尾变量，词值为 \(b\)；其他右项值至多为 \(e\)。右和不超过 \(b\)，而左词值为 \(c\) 或 \(t\)。三种情形均给出分离赋值。

充分性。 首先，后两个条件保证 \(X\subseteq c(F)\)。逐个检查真主理想的阈值 \(a, e, b, c\)。

若右和不超过 \(a\)，全部右内容只取 \(a, e\)，且每个右项的内容含一个取 \(a\) 的变量。任取内部项，其内容包含于 \(X\)，故左内容也含 \(a\)，左词值为 \(a\)。若右和不超过 \(e\)，全部右内容只取 \(a, e\)，因此左值也不超过 \(e\)。

若右和不超过 \(b\)，全部右尾变量取 \(a, e\)。由尾部覆盖，左尾变量亦如此。若 \(h\in E\)，左首字母也取 \(a, e\)，结论成立。否则存在内部同首字母项 \(r\)。若左首字母不超过 \(b\)，左值也不超过 \(b\)。左首字母不可能取 \(t\)，因为这会使 \(r=t\)。若左首字母取 \(c\)，为使 \(r\le b\)，其尾部必有取 \(a\) 的变量。由于 \(h\notin E\)，\(D(r)\subseteq X\setminus\{h\}\subseteq D(w)\)，故左尾部也含 \(a\)，左值为 \(b\)。

若右和不超过 \(c\)，全部右尾变量仍取 \(a, e\)。若 \(h\in E\)，左首字母为低值；否则内部同首字母项保证左首字母不是 \(t\)。左尾部只有 \(a, e\)，因此左词值不超过 \(c\)。全部阈值检验完成。证毕。

\par\smallskip\noindent\textit{3. 五条完整基.}\quad 
定理2。 此半环的完整附加基为

\CaseDisplay{\begin{aligned}
xyz&\approx xzy, &(N)\\
xy^2&\approx xy, &(R)\\
x^2y&\approx y^2x, &(C)\\
x&\le x^2, &(I)\\
xy&\le x+zy. &(V)
\end{aligned} \quad\textup{(\ensuremath{3})}}

有效性由词值规则或定理1检验。以下证明任意变量数的完备性。

若 \(p\le F\) 已经推出，且 \(x\in E\)，选择尾部含 \(x\) 的右项。由 \((N),(R)\) 可把它写成 \(vx\)，其中 \(v\) 非空。由 \((V)\)

\CaseDisplay{px\le p+vx\le F. \quad\textup{(\ensuremath{4})}}

设 \(w\le F\) 有效，仍写 \(h=h(w)\)、\(X=c(w)\)。若存在内部同首字母项 \(r\)，从 \(r\) 出发用 (4) 补入全部 \(D(w)\)。由于 \(c(r)\subseteq X\)，所得词的尾部为 \(D(w)\) 或 \(D(w)\cup\{h\}\)，即得到 \(w\le F\) 或 \(w^2\le F\)；后一情形再用 \((I)\)。

否则由定理1，\(h\in E\)。同时 \(D(w)\subseteq E\)，故 \(X\subseteq E\)。任取内部项 \(q\)，利用 (4) 补入整个 \(X\)。所得词首字母在其尾部出现，内容恰为 \(X\)，所以由 \((N),(R),(C)\) 它等于 \(M_X=w^2\)。再由 \((I)\) 推出 \(w\le F\)。

一次左词也由上述推导处理：内部项必须为同名纯幂。所有代入均为非空词，无空积或空词。任意半环恒等式可展开为双向逐词吸收，因此五式构成完整基。证毕。

\par\smallskip\noindent\textit{4. 适用范围.}\quad 
乘法对偶逐词反转上面的构造和证明。对其他目录对象，只有在独立验证双向 HSP 包含后才传递有限基性。固定一至五变量诊断与运算表核验用作交叉检查；任意变量数结论由第2---3节证明承担。

五式是在来源4150的四式基础上增加 \((C)\)，但这一包含关系本身不是有限基传递定理。四式分支是否在整个五元素目录中已无未决成员，另由完整枚举及独立验证确定；即使该分支全部有限基，也不推出其定义的簇遗传有限基。

\subsection{A directed-cycle obstruction:12480 is nonfinitely based}\label{case:288}
Signature (+, multiplication), without constants, relative to the additively idempotent semiring axioms. This proof allows arbitrary polynomial substitutions, repeated source variables and arbitrary contexts.

\par\smallskip\noindent\textit{1. Word values and exact absorption.}\quad 
Addition:00004,\allowbreak{}01114,\allowbreak{}01224,\allowbreak{}01234,\allowbreak{}44444; order \(3<2<1<0<4\). Multiplication:44034,\allowbreak{}41134,\allowbreak{}41234,\allowbreak{}33333,\allowbreak{}44434. Any 3 in a nonlinear word gives 3. Without 3, any tail 0 or 4 gives 4. If every tail is 1 or 2, head 4 gives 4, head 1 gives 1, head 2 gives 2 when all tails are 2 and 1 otherwise, and head 0 gives 0 when all tails are 2 and 4 otherwise. Literals retain value. Induction on length proves the formula. Thus \(xyz=xzy\) and \(xyy=xy\); normal words consist of head and sorted tail set, retaining any occurrence of the head in the tail.

For a target w with head h and content C, let U be the original right words with entire content in C. Let D be their tail-content union and \(D_{h}\) the tail-content union of those headed h. A literal h is absorbed iff some original pure power of h occurs. A nonlinear target hT is absorbed exactly when: (C) contents of U cover all C; (D) T is contained in D; (S) h belongs to D, or T is contained in \(D_{h}\).

Necessity:for a missing pure-power literal witness, set \(h=4\), \(others=3\). If C fails at x, set \(x=1\), other target \(variables=2\), \(outsiders=3\). Target 1 is separated from a right sum in\{2,\allowbreak{}3\}. If D fails at x in T, set \(x=0\), other target \(variables=2\), \(outsiders=3\). Target 4 is separated from a right sum at most 0, because x never occurs in a contained right tail. If S fails, set \(h=0\), \(D_{h}=2\), other target \(variables=1\), \(outsiders=3\). No original contained tail contains h; h-words give 0, other contained words give 1 or 2, and outsiders force 3. Some target tail is outside \(D_{h}\), so target 4 is separated. For sufficiency, target 3 is least. A target 4 with a tail 0/\allowbreak{}4 is witnessed by D. Otherwise its tails are 1/\allowbreak{}2. A head 4 is witnessed by C; a head 0 with a tail 1 is witnessed by S, either through a tail occurrence of h or a contained h-word containing a tail 1. For target 0, C supplies an occurrence of the head 0 giving 0 as head or 4 as tail. A target 2 has all variables 2, so any contained word gives 2. A target 1 has variables only 1/\allowbreak{}2 and some 1; C supplies a word containing 1, valued 1. For literals \(x\leq x^{m}\) for every positive m follows from the table. This proves the criterion for arbitrary finite words, polynomials and assignments, without a variable bound. In particular(C) holds for every valid absorption, including literals.

\par\smallskip\noindent\textit{2. The critical identities and their absorbed words.}\quad 
For \(n\geq 4\), use distinct variables \(x_0,\ldots,x_{n-1}\), indices modulo n. Put \(P_n=\sum_{i=0}^{n-1}x_ix_{i+1}\), \(r_n=x_0x_1\cdots x_{n-1}\), and \(X_n=\{x_0,\ldots,x_{n-1}\}\). The identity \(P_{n}=P_{n}+r_{n}\) holds for every n:all contents and tails of \(P_{n}\) cover \(X_{n}\), and the head of \(r_{n}\) belongs to its tail union D. This is an all-length consequence of the exact criterion.

Call a raw word an edge word if it has form \(x_{i}\) \(x_{i+1}\)\ensuremath{^{m}}, \(m\geq 1\). Call a word critical if its content is all \(X_{n}\). Every word absorbed by \(P_{n}\) is either an edge word or critical, and no variable outside \(X_{n}\) occurs. Conversely every such word is absorbed.

Proof:condition C excludes outside variables, and literals have no original pure-power witness. Suppose a nonlinear absorbed word has proper content E subset \(X_{n}\) and head h. The original words contained in E are precisely the edges of the directed n-cycle with both endpoints in E. The induced graph is a disjoint union of directed paths and isolated vertices. Condition C excludes isolated vertices. Condition D says that every vertex other than possibly h has an incoming edge. Thus there is only one path component, and its first vertex is h; h cannot itself lie in the target tail. Now h is absent from D, and condition S requires every other target variable to occur in a contained h-word's tail. The only such original edge is h followed by its successor. Hence E is exactly that two-vertex edge and the raw word is h successor(h)\ensuremath{^{m}}. Conversely these edge words satisfy the criterion. For full content \(X_{n}\), all original heads/tails are available and every nonlinear full-content word satisfies C, D, S. Since \(n\geq 4\) it is never literal. \(\square\)

The proper-content characterization applies to RAW words, with arbitrary multiplicities. It does not assume substitutions have been normalized.

\par\smallskip\noindent\textit{3. Polynomial substitution barrier.}\quad 
Let \(U=V\) be any identity of 12480 on at most k variables. Let sigma send source variables to arbitrary nonempty finite polynomials of nonempty words. Suppose every expanded word of sigma(U) is an edge word. If n\ensuremath{>}max(3k,3), sigma(V) contains no critical word.

Choose a monomial v of V and a proposed critical expansion w of sigma(v). The source absorption \(v\leq U\) is valid. If v is a literal a, the exact criterion gives an original monomial \(a^{m}\) in U. Choosing w at every occurrence of a in its expansion gives \(w^{m}\) among sigma(U). Its content is \(X_{n}\), contrary to every old expansion being an edge word. This treats \(m=1\) and all positive m without an empty substitution.

Suppose v is nonlinear, with head h and nonempty tail-variable set T. For every a in c(v), there is a NONLINEAR original monomial \(u_{a}\) of U, contained in c(v), which contains a. For each a in T, this follows from source condition D, which puts a in a tail of a contained original word. For the head h, if h is in that same D it also occurs in such a tail. Otherwise S puts T in \(D_{h}\); as T is nonempty, this supplies a nonlinear contained h-word. This handles both a repeated head and a simple head.

Fix an occurrence of a in \(u_{a}\) and a different position of \(u_{a}\). Fix one nonempty image word at that other position, and choose any letter j in it; fix arbitrary choices at all remaining positions. Vary independently the chosen image at the selected occurrence of a. All resulting words are expanded words of sigma(U), hence edge words containing j. An edge of the directed n-cycle containing j has its support in either \{pred(j), j\} or\{j, succ(j)\}. Consequently EVERY image word in sigma(a) has content inside the fixed three-element set\{pred(j), j, succ(j)\}. This remains valid if the two positions use the same source variable: distributive polynomial expansion chooses image words independently at different occurrences. No singleton-image assumption has been made.

Every expansion of v therefore uses at most 3k target variables, the union of these fixed sets for its at most k source variables. It cannot have content \(X_{n}\) when \(n>3k\). This proves the polynomial-substitution barrier.

\par\smallskip\noindent\textit{4. Contexts and nonfinite basis.}\quad 
Suppose a finite relative basis Sigma existed; let k bound its variables. Choose n\ensuremath{>}max(3k,3). Start a derivation of \(P_{n}=P_{n}+r_{n}\) from \(P_{n}\), expanding all terms distributively into finite nonempty polynomials. Every intermediate expression is equal to \(P_{n}\) in 12480. Each expanded word is therefore absorbed by \(P_{n}\), and by section 2 is an edge word or critical. \(P_{n}\) has no critical word. Consider the first step introducing one.

The step uses a basis identity \(U=V\), in either direction, with arbitrary polynomial substitutions, in a one-hole term context. On expansion each contribution of the hole has form s sigma(U) t, where s, t are fixed words or omitted contexts. Choose a contribution producing a critical word w after replacement. Every old word of this contribution is an edge word. If both contexts are omitted, section 3 forbids w. If at least one context is nonempty, fix a letter j in it. Every old word contains j, so all its content lies within\{pred(j), j, succ(j)\}. The valid contextual identity implies \(w\leq s\) sigma(U)t independently of unaffected summands. Apply necessary content condition C to this absorption. Its contained old words together cannot cover \(n>3\) different variables of w, a contradiction. If a context letter were outside \(X_{n}\), no old edge word could contain it, already contradicting section 2. Polynomial contexts expand into several fixed-word contributions; the selected contribution is still a valid contextual identity. Empty contexts are only omitted, never empty-word substitutions.

The ai-semiring axioms themselves preserve free-polynomial expansions, so cannot introduce a critical monomial. Thus no first critical word can appear in a finite derivation, although \(P_{n}=P_{n}+r_{n}\) is valid. This contradicts completeness of Sigma. Therefore 12480 is NONFINITELY BASED. Reversing all words proves the same conclusion for its opposite. No claim of inherent nonfinite basis or of limit-variety status is made.

\par\smallskip\noindent\textit{5. Transfer and all-length verification.}\quad 
If \(A=12480\) belongs to Var(B) and B satisfies \(P_{n}=P_{n}+r_{n}\) for all \(n\geq 4\), then B is NFB:every identity of B holds in A, and the same derivation barrier applies while B's critical family must still be derived. Both hypotheses are essential; containment alone does not suffice.

For finite B an exact state certificate proves the all-n hypothesis. For a sequence \(a_0,\ldots,a_{\ell-1}\), \(\ell\ge2\), keep state(f, t, p, s), where f is its first value, t its last, p its ordered product, and s is the sum of the adjacent products \(a_i a_{i+1}\) for \(0\le i<\ell-1\). The length-two states are (a, b, ab, ab). Appending c sends the state to(f, c, pc, s+tc). Form the exact length-four states and close under every append. Induction says the finite closure is exactly the states of all sequences of \(lengths\geq 4\). For each state check \(p\leq s+tf\). These are precisely all critical identities, with no length cutoff. The opposite orientation reverses multiplication. Complete positive HSP establishes A in Var(B); complete reverse HSP or an explicit strict separating identity determines whether containment is equality or proper. Only a proved proper NFB subvariety excludes limit generation; equal-variety NFB types retain unresolved limit status.

Independent checks verify exact semantics, raw absorbed-word classification, polynomial substitutions including repeated source positions, context separators, finite-state closure and every used HSP direction. Computation supports these unrestricted arguments; it does not replace their quantifiers.

\subsection{Five identities for type 12621}\label{case:289}
For the catalogue's two-operation signature without constants, the following identities and the additively idempotent semiring axioms form a finite basis.

{\small\par\smallskip
\CaseLine{\(xyz = xy + xz + yz\)}
\CaseLine{\(x \leq  xy\)}
\CaseLine{\(y \leq  xy\)}
\CaseLine{\(xx \leq  xy + yx\)}
\CaseLine{\(xy \leq  xx + yy\)}
\smallskip}

Here u \ensuremath{\leq} v abbreviates \(u+v=v\). Each identity is valid by direct evaluation of the given tables. The first identity reduces every word of length at least two to the sum of the products of all its ordered pairs of positions \(i<j\). For the induction step substitute the prefix and the last two letters for x, y, z. The two shorter products reproduce the pairs within the prefix and those ending at either last position. The third product supplies the remaining pair. Equal letters in different positions give a loop.

After this reduction the upper polynomial consists of literals and directed edges. Call a vertex protected if it has a loop or is incident with an edge whose converse is also present. A literal is absorbed exactly when it occurs in the upper content. A quadratic uv is absorbed exactly when uv is an original edge or both u and v are protected.

To derive the positive cases, obtain each literal from any incident edge by the second or third identity. For every protected vertex u obtain uu directly from its loop or from its two opposite edges using the fourth identity. For two protected vertices combine their loops using the last identity. This includes the case \(u=v\).

For necessity assign 4 to an absent target variable and 2 to all others. Every upper value is then 2 and the target is 4, since 4 absorbs products. For a missing edge uv with an unprotected endpoint x, use the principal ideal below 0, consisting of 0,\allowbreak{}1,\allowbreak{}2,\allowbreak{}3. Set x to 0, its incoming neighbors to 3 and its outgoing neighbors to 2. Set remaining variables to 1. The neighbor sets are disjoint and exclude x. Every upper edge between nonzero colors belongs to the ideal, and so do 3*0 and 0*2. All upper literals belong to the same ideal. If \(x=u\) and u differs from v, then v is not an outgoing neighbor, so 0*v is outside the ideal. If \(x=v\) and u differs from v, then u is not an incoming neighbor, so u*0 is outside. When \(u=v\) the product 0*0 is outside. This supplies every required separation. The argument works for arbitrary finite graphs, numbers of variables and polynomials. Reducing and deriving each summand in both directions proves completeness.

\subsection{A three-identity basis from directed tail incidence}\label{case:290}
We use the signature (+, multiplication) without constants, relative to the additively idempotent semiring axioms. Let A be catalogue source 12624, on the elements 0,\allowbreak{}1,\allowbreak{}2,\allowbreak{}3,\allowbreak{}4, with the following row strings in that order:

{\small\par\smallskip
\noindent addition    multiplication\par
\CaseLine{\ttfamily 00004\qquad 44044}
\CaseLine{\ttfamily 01114\qquad 41114}
\CaseLine{\ttfamily 01214\qquad 41214}
\CaseLine{\ttfamily 01134\qquad 41314}
\CaseLine{\ttfamily 44444\qquad 44444}
\smallskip}

The element labels are names, not constants in the term language. Addition is join with incomparable minimal elements 2 and 3, their join 1, and \(1<0<4\).

Theorem. The following three identities form a basis for all identities of A, relative to the ai-semiring axioms:

{\small\par\smallskip
\CaseLine{\((S) x \leq  xy,\)}
\CaseLine{\((B) xyz \approx  xy+xz,\)}
\CaseLine{\((G) xy \leq  ux+vy+zq+q.\)}
\smallskip}

Here a\ensuremath{\leq}b abbreviates a+b\ensuremath{\approx}b; x, y, u, v, z, q are ordinary variables. Word reversal gives a basis for the multiplicative opposite. Minimality, independence and hereditary finite basability are not asserted.

\par\smallskip\noindent\textit{Structure and reduction to a directed graph.}\quad 
Right multiplication by 2 is the identity. Right multiplication by 1 or 3 is the same map f, with values (4,\allowbreak{}1,\allowbreak{}1,\allowbreak{}1,\allowbreak{}4); right multiplication by 0 or 4 is the constant map with value 4. The maps f and the constant map commute and are idempotent, and \(f(4)=4\). The type of a product yz agrees with the composition of the types of right multiplication by y and z: a factor of type constant makes the product type constant, two identity types give identity, and all remaining pairs give f. This checks associativity.

The right maps preserve the displayed joins. Each left map also preserves joins: for a fixed left factor a, the value at 2 is a, at 1 and 3 is f(a), and at 0 and 4 is 4. The relations a\ensuremath{\leq}f(a)\ensuremath{\leq}4 ensure preservation along each branch, and the only incomparable pair 2,\allowbreak{}3 has join 1, with \(a+f(a)=f(a)\). Thus A is a semiring.

Every right map is increasing above the identity, giving (S). The value of a word is the head acted on by the strongest right-map type occurring in its tail: identity if all tail factors are 2, f if at least one is 1 or 3 but none is 0 or 4, and the constant map otherwise. The join of the two right-map outputs equals their composition, proving (B). Hence every word of length at least two is provably the sum of the binary words consisting of its head and each subsequent variable. This follows by induction from (B), including repeated variables. Distributivity consequently reduces any term to a finite sum of linear variables and binary words.

For such a nonempty sum P let L be the set of linear variables, E the set of ordered pairs (a, b) whose binary word ab occurs, H the set of heads of E, and R the set of tails of E. Put \(Q=L\)\ensuremath{\cup}H. Loops are allowed. The total content is Q\ensuremath{\cup}R.

\par\smallskip\noindent\textit{Exact absorption criterion.}\quad 
Lemma. A linear variable x is absorbed by P if and only if

{\small\par\smallskip
\noindent x\ensuremath{\in}Q, or x\ensuremath{\in}R and Q\ensuremath{\cap}R is nonempty.\par
\smallskip}

A binary word xy is absorbed by P if and only if

{\small\par\smallskip
\noindent (x, y)\ensuremath{\in}E, or x, y\ensuremath{\in}R and Q\ensuremath{\cap}R is nonempty.\par
\smallskip}

The statement includes \(x=y\).

Necessity for linear targets. If x is outside Q\ensuremath{\cup}R, give x value 4 and every other variable 2. Then \(x=4\) and \(P=2\). If x\ensuremath{\in}R but x\ensuremath{\notin}Q and Q\ensuremath{\cap}R is empty, assign Q the value 3 and R the value 2. Every linear right summand is 3, and every binary right summand has head 3 and tail 2, hence is 3. But \(x=2\), incomparable with 3.

Necessity for binary targets. Suppose (x, y) is absent from E. If y\ensuremath{\notin}R, assign \(y=0\) and all other variables 2. No right binary word has tail y; it therefore has value 0 or 2, and the linear part is also at most 0. But \(xy=4\).

Next suppose y\ensuremath{\in}R and x\ensuremath{\notin}R. Then x\ensuremath{\ne}y. Assign \(x=0\), \(y=1\), and the other variables 2. A right word with head x has tail different from y, since xy is absent, and also different from x, since x\ensuremath{\notin}R. It therefore has value 0. All other right words have value at most 1; linear summands have value at most 0. Again P\ensuremath{\leq}0 while \(xy=4\).

Finally, if x, y\ensuremath{\in}R and Q\ensuremath{\cap}R is empty, assign Q the value 3 and R the value 2 as before. Then \(P=3\) and \(xy=2\). These cases establish the necessity of the two criteria.

Sufficiency for a nonliteral binary target. Suppose x, y\ensuremath{\in}R and Q\ensuremath{\cap}R is nonempty. Consider any evaluation. If \(P=4\) there is nothing to prove. If \(P=0\), each binary summand is at most 0, so its tail cannot be 0 or 4. Thus all variables in R, including x, y, are in \{1,\allowbreak{}2,\allowbreak{}3\}, and xy\ensuremath{\leq}1\ensuremath{\leq}0. If \(P=1\), every variable in the content of P belongs to \{1,\allowbreak{}2,\allowbreak{}3\}, and again xy\ensuremath{\leq}1. If \(P=2\), minimality of 2 forces each summand to be 2. A product equals 2 only when both factors are 2. Therefore \(x=y=2\) and \(xy=2\). The remaining possibility \(P=3\) is impossible: minimality of 3 would force every linear variable and every binary head to be 3, and every binary tail to be 2, contradicting Q\ensuremath{\cap}R\ensuremath{\ne}\ensuremath{\varnothing}.

Literal binary targets are automatic. For a linear target x\ensuremath{\in}Q, either it is literal or an outgoing binary word xt is a right summand, and (S) yields x\ensuremath{\leq}xt\ensuremath{\leq}P. If x\ensuremath{\in}R and Q\ensuremath{\cap}R\ensuremath{\ne}\ensuremath{\varnothing}, the binary case gives \(x^{2}\)\ensuremath{\leq}P, and (S) gives x\ensuremath{\leq}\(x^{2}\). This completes the criterion for every finite variable set.

In (G), x and y are tails of right words, and q belongs both to R (from zq) and to Q (from its linear occurrence). The lemma therefore proves soundness of (G).

\par\smallskip\noindent\textit{Equational completeness.}\quad 
Let an absorption by P satisfy the lemma. Derive q\ensuremath{\leq}P for any q\ensuremath{\in}Q: it is either a right summand or follows from (S) applied to an outgoing right word qt.

For a nonliteral valid binary target xy, choose right words ux and vy, possible because x, y\ensuremath{\in}R. Choose q\ensuremath{\in}Q\ensuremath{\cap}R and a right word zq. The already available q\ensuremath{\leq}P and the three original binary summands give

{\small\par\smallskip
\CaseLine{\(xy \leq  ux+vy+zq+q \leq  P\)}
\smallskip}

by one instance of (G). Coincident variables or summands cause no difficulty, since substitution need not be injective and addition is idempotent. For a linear target in Q use (S) as above. For the remaining valid linear target derive \(x^{2}\) by the binary argument and then x\ensuremath{\leq}\(x^{2}\) by (S). Thus all valid absorptions of a linear or binary word follow from the three displayed identities.

Reduce arbitrary terms by distributivity and (B). For a valid identity P\ensuremath{\approx}P', every summand of P is absorbed by P', and every summand of P' is absorbed by P. Apply the preceding derivations and combine them to obtain the identity. This proves completeness without bounds on variables, original word length, or the number of summands. Reversing words transports all derivations to the opposite algebra.

Finite table checks are only supplementary error checks. Classification of other catalogue entries requires complete two-way HSP equivalence, independently verified from their operation tables.

\subsection{A finite basis from head-indexed tail coverage}\label{case:291}
Source 12626 in the five-element catalogue. The signature is \((+, multiplication)\) without constants. All bases are relative to the additively idempotent semiring axioms. Element labels are not term constants.

\par\smallskip\noindent\textit{1. The algebra, word values, and normal words.}\quad 
The tables, in row and column order 0,\allowbreak{}1,\allowbreak{}2,\allowbreak{}3,\allowbreak{}4, are:

{\small\par\smallskip
\noindent addition       multiplication\par
\CaseLine{\ttfamily 0 0 0 0 4\qquad 4 4 0 4 4}
\CaseLine{\ttfamily 0 1 1 1 4\qquad 4 1 1 1 4}
\CaseLine{\ttfamily 0 1 2 1 4\qquad 4 1 2 3 4}
\CaseLine{\ttfamily 0 1 1 3 4\qquad 4 1 1 1 4}
\CaseLine{\ttfamily 4 4 4 4 4\qquad 4 4 4 4 4}
\smallskip}

The additive order is 2,\allowbreak{}3 \ensuremath{<} 1 \ensuremath{<} 0 \ensuremath{<} 4, with 2 and 3 incomparable. Associativity and both distributive laws can be verified directly on the 125 triples in these displayed tables. No assumptions about identities or about the number of variables enter this finite algebra verification.

For a word of length at least two the following description follows by induction. A factor 4, or a factor 0 after the first position, gives 4. If the first factor is 0, the value is 0 when every later factor is 2, and 4 otherwise. If all factors are in \{1,\allowbreak{}2,\allowbreak{}3\}, the value is 2 when all are 2, is 3 when all nonfinal factors are 2 and the final factor is 3, and is 1 in every other case.

Write c(u) for content, h(u), t(u) for endpoints, r(u) for the content after the first position, and i(u) for the content before the last position. For a linear word, \(r(u)=i(u)=empty\). The value formula gives

\CaseDisplay{xyzt\approx xzyt,\qquad xy^2z\approx xyz.                 \quad\textup{(\ensuremath{W})}}

Fix a variable order. Adjacent interior positions may be interchanged by the first law, using nonempty outer prefixes and suffixes. The second law removes adjacent interior repetitions. Thus every nonlinear word normalizes to h (product of M) t, with M its interior set in fixed order. The endpoints remain literal, may coincide, and may occur in M. An empty M denotes a two-letter word, not a unit. Linear words remain linear.

\par\smallskip\noindent\textit{2. An exact absorption criterion.}\quad 
For a nonempty polynomial P of words set

\CaseDisplay{D=\bigcup_{u\in P}c(u),\quad E=\bigcup_{u\in P}r(u),\quad
 J=\bigcup_{u\in P}i(u),\quad F=\{t(u):u\in P\},
 \qquad E_h=\bigcup_{u\in P,\ h(u)=h}r(u).}

Linear summands contribute to D and F, but not to E, J, or \(E_{h}\).

Lemma 1. For every nonempty word w, including a linear word, \(w\leq P\) holds in the displayed algebra if and only if all four conditions hold:

1. c(w) is contained in D;

2. r(w) is contained in E;

3. either h(w) belongs to E, or r(w) is contained in \(E_{h(w)}\);

4. either J meets F, or i(w) is contained in J and t(w) belongs to F.

Proof. A variable in c(w) outside D can be assigned 4, with every other variable 2. Then \(w=4\) and \(P=2\), proving condition 1.

If x is in r(w) outside E, assign \(x=0\) and all others 2. Every right word containing x is either linear x or has x only at its first position; it has value 0. Other right words have value 2. Thus \(P\leq 0\) whereas \(w=4\). This proves condition 2.

Now assume condition 2, and suppose condition 3 fails. Put \(h=h(w)\). There is x in r(w) outside \(E_{h}\), and h is outside E. In particular x differs from h, by condition 2. Assign \(h=0\), \(x=1\), and every other variable 2. A right word containing h is linear h or has h only at its head. Such a nonlinear word cannot contain x, by the definition of \(E_{h}\), so it has value 0. All remaining words have value at most 1. Consequently \(P\leq 0\), but \(w=4\). Condition 3 is necessary.

Suppose conditions 1--3 hold and J is disjoint from F. Give J the value 2 and F the value 3; these assignments are consistent. Give outside variables 2. Every right word is 3. Since c(w) is contained in \(D=J\) union F, its value is 3 exactly when all nonfinal variables belong to J and its final variable belongs to F. This includes linear w. Failure of condition 4 therefore separates w from P.

Conversely fix an assignment. The case \(P=4\) is immediate. If \(P=2\), every right word equals 2 and hence every variable of D is 2; condition 1 gives \(w=2\). If \(P=3\), every right word equals 3. The word formula forces J to be assigned 2 and F to be assigned 3. They are disjoint, and condition 4 gives \(w=3\).

If \(P=1\), no variable of D can be 0 or 4: an occurrence of 4 produces 4, and an occurrence of 0 in a right word produces either 0 or 4. Thus all variables of c(w) lie in \{1,\allowbreak{}2,\allowbreak{}3\}, and \(w\leq 1\).

Finally let \(P=0\). All variables of D avoid 4, so the same holds for c(w). If a variable x in r(w) were 0, condition 2 would give a right word with a noninitial occurrence of x and hence value 4. Therefore all target tail variables avoid 0. If the target head is not 0, the target is at most 1. If its head h is 0, h cannot belong to E, since that too would produce a right word of value 4. Condition 3 then puts r(w) inside \(E_{h}\). A target tail variable with value other than 2 would occur after head h in some right word, making that word 4. All target tail values are therefore 2, and \(w=0\). The same argument covers a linear target, whose tail is empty. These cases prove sufficiency.

\par\smallskip\noindent\textit{3. Derivation in steps with at most five right words.}\quad 
Only the validity of the following steps is used here; Section 4 reduces every such step to an explicitly specified finite family. All words are normal, and \(w\leq P\) satisfies Lemma 1. Every constructed right word is original or has already been derived below P.

First let \(w=y\) be linear. If J meets F, choose original words a, b with some g in i(a) and \(t(b)=g\), and choose \(r_{y}\) containing y. Then \(y\leq a+b+r_{y}\) is valid by Lemma 1. If J and F are disjoint, condition 4 gives y in F. Choose an original v ending in y; the one-right-word absorption \(y\leq v\) is valid. This settles linear targets.

For nonlinear \(w=h\) (product of M) t, write \(T=r(w)\). Fix the following original witnesses throughout the construction. If h belongs to E, choose \(r_{h}\) with h in r(\(r_{h}\)); for each x in T choose \(r_{x}\) with x in r(\(r_{x}\)). If h is outside E, condition 3 lets us choose, for each x in T, \(r_{x}\) with head h and x in r(\(r_{x}\)). In the second case h is not in T. In both cases \(r_{t}\) exists.

Crossing case J intersect F nonempty. Choose original a, b with g in i(a) and \(t(b)=g\). Start with \(s=ht\). If h belongs to E, use

\CaseDisplay{ht\leq r_h+r_t+a+b.                   \quad\textup{(\ensuremath{2})}}

If h is outside E, use \(ht\leq r_{t}+a+b\) instead. Content and tail coverage follow from the chosen witnesses. In the first case \(r_{h}\) puts h in the right-hand E; in the second, \(r_{t}\) supplies the whole tail of ht to the right-hand \(E_{h}\). The words a, b retain the crossing. Thus these are valid steps by Lemma 1.

For each x in M not yet placed in the interior, add it and normalize to obtain s'. If h belongs to E, use

\CaseDisplay{s'\leq s+r_x+r_h+a+b.                \quad\textup{(\ensuremath{3})}}

If h is outside E, omit \(r_{h}\). The old word s covers its own content and tail; \(r_{x}\) supplies the new variable in the tail. The head condition is witnessed by \(r_{h}\) in the first case, and by the heads h of s, \(r_{x}\) in the second. The crossing remains. This also covers an interior x equal to an endpoint. There are at most five right words.

Separated case J intersect F empty. Condition 4 gives i(w) contained in J and t in F. Thus h differs from t, t is not interior, and every interior x lies in J. Any original right word containing t has t only at its last position. If h belongs to E, the witness \(r_{h}\) has h after its first position and, since h is outside F, before its last position as well. Consequently

\CaseDisplay{ht\leq r_h+r_t                   \quad\textup{(\ensuremath{4})}}

is valid. If h is outside E, the selected \(r_{t}\) has head h and final t, so \(ht\leq r_{t}\) is valid instead.

For an interior x, the chosen \(r_{x}\) contains x in its tail. Since x belongs to J and J is disjoint from F, that occurrence is also before its final position. Add x to the current interior and use

\CaseDisplay{s'\leq s+r_x+r_h                   \quad\textup{(\ensuremath{5})}}

when h belongs to E, omitting \(r_{h}\) otherwise. Content and tail coverage, and the head-indexed requirement, hold as in (3). The right-hand prefix contains i(s) and x, and the final set still contains t. Thus the fourth condition also holds, whether or not the smaller right sum has a crossing. These steps have at most three right words.

After at most |M| extensions the constructed word is w. Each step has at most five original or previously derived right words. This is an unrestricted construction, not a bound on the variables of w or P.

\par\smallskip\noindent\textit{4. Reversible compression and finite basis.}\quad 
A step has at most six words altogether. Reserve individually every variable occurring at a first or last position of any of them, including linear words. At most 12 variables are reserved. All other variables occur only in interiors. Group them by their nonempty incidence pattern among the at most six interiors. There are at most \(2^{6}-1=63\) groups.

Identify each group to one fresh variable and normalize with (W). This is a substitution in a valid identity, followed by valid word laws, so the resulting step is valid. It uses at most 75 variables and normal words of length at most 77. Conversely replace each group variable by the nonempty product of its original members. Their identical incidence patterns restore exactly every original interior set. Reserved variables are replaced individually. Normalizing with (W) recovers every original word. Every substitution is nonempty, and no constants are introduced.

Fix a 75-letter alphabet. Define Sigma to consist of (W) and every absorption valid in this algebra whose left word and at most five right words are normal words over that alphabet of length at most 77. This is a finite, explicitly specified family. Membership is decided by Lemma 1. No small or minimal enumeration is claimed.

Theorem. Sigma is a finite identity basis for source 12626. Word reversal gives a finite basis for its multiplicative opposite.

Proof. Soundness follows from (W) and the definition of Sigma. Compression shows that every valid step with at most five right words follows by substitution and (W). Section 3 derives every valid word absorption \(w\leq P\) using such steps, with each right word already below P. Distribute any term into a sum of words and normalize. An identity between two such sums is equivalent to the absorption of every summand on each side by the other sum; all these have been derived. This proves completeness for arbitrarily many variables. Word reversal transports the identities and derivations to the opposite.

The independent computations check table axioms, word values, the exact criterion, every constructed local step, and reversible compression. They support the unrestricted proof above and do not replace it. Transfers to other catalogue entries require equality of generated varieties verified in both HSP directions.

\subsection{Five finite bases from a split bridge in the endpoint graph}\label{case:292}
Use the constant-free signature of additively idempotent semirings. A literal is a word of length one. The following presentations are relative to the semiring axioms. Reverse every word in a presentation to obtain the basis for its multiplicative opposite.

\noindent Source\enspace\textbar\enspace Opposite\enspace\textbar\enspace Parent finite basis\enspace\textbar\enspace Additional law\par
\noindent 12880\enspace\textbar\enspace 13017\enspace\textbar\enspace 837 (147 templates and fixed laws in total)\enspace\textbar\enspace xy \ensuremath{\leq} zxxy+y\par
\noindent 12887\enspace\textbar\enspace 15108\enspace\textbar\enspace 11619 (199)\enspace\textbar\enspace xxy \ensuremath{\leq} xy+y\par
\noindent 15113\enspace\textbar\enspace 15113\enspace\textbar\enspace 14674 (244)\enspace\textbar\enspace yy \ensuremath{\leq} y+xx\par
The three parent bases and their unrestricted completeness theorems are in Section~\ref{case:45}, Section~\ref{case:253} and Section~\ref{case:322}. Only their proved absorption criteria and the implication from validity to derivability in each parent are used below.

\par\smallskip\noindent\textit{Tables and Boolean endpoints.}\quad 
In all three sources (e, f, a, c, z) have catalogue labels (3,\allowbreak{}2,\allowbreak{}1,\allowbreak{}0,\allowbreak{}4). Addition has \(a<c<z\), with e and f each incomparable to every other proper element. Multiplication of e=(0,\allowbreak{}0), a=(0,\allowbreak{}1), f=(1,\allowbreak{}1) concatenates matching endpoints and gives z for a mismatch; z is a zero. On deleting the marking (sending c to a), a word giving a bridge must have exactly one bridge or marked bridge between a prefix of e's and a suffix of f's. In 12880 all such words are compatible: a marked bridge gives c if initial, a if noninitial. In 12887 it gives c if initial, z if noninitial. In 15113 it gives c only in its literal; its presence in a nonlinear word gives z. Incompatible words have value z in every source. These statements follow from the four binary products involving a marked input that could compose: e c and c f, with their ordinary versions e \(a=a\) and a \(f=a\). All other products with c equal z. No general homomorphism under deletion of markings is needed for the last two sources.

Given P, use C for its content union, \(H_{0}\) for its original head letters, D for the letters at noninitial positions, and N for the letters anywhere in nonlinear original words. Build the directed endpoint graph of the parent theorems. Identify each L of an original head with H, each R of an original final with T, and each \(R_{x}\) with \(L_{y}\) for an original adjacency xy. Insert the edge \(L_{x}\) \ensuremath{\to} \(R_{x}\) for every x in C, and identify strongly connected components. In what follows, equality of endpoints means equality in this last quotient. Let W mean that w is an H-to-T walk in the quotient, including its initial and final endpoints and all adjacent joins.

Assigning ordinary letters only, every right word is a precisely when all matching endpoints carry equal bits, all directed edges carry nondecreasing bits, H is 0 and T is 1. Every vertex lies between H and T on a directed walk. Consequently these assignments exist exactly when H differs from T. In that case an equality of two endpoints is forced by all assignments exactly when the endpoints are in the same component: an upward closed set separating different components can include T and omit H, with its complement receiving bit 0. For an edge whose initial component is H and whose final component is not H, an assignment making that letter a is given by the upward closure of its final component together with T. If \(H=T\) the graph has a single component. These are the unrestricted cut and walk arguments from the parent criteria.

\par\smallskip\noindent\textit{Exact absorption tests.}\quad 
Each source requires c(w) contained in C and W. Content necessity uses all original letters e and a missing letter z. For W use the ordinary bridge assignments, which give right value a and force the target to be a. The additional requirements are as follows.

\textbullet{} In 12880 the target head x must belong to \(H_{0}\) or satisfy \(R_{x}=H\).

\textbullet{} In 12887 every variable x occurring noninitially in w but outside D must satisfy \(R_{x}=H\).

\textbullet{} In 15113 a nonlinear w must be supported in N unless \(H=T\).

For 12880, if its extra test fails, \(L_{x}=H\) by W and x can be given the bridge a on an ordinary bridge assignment. Change only x to c. Every original occurrence is noninitial, so every original still gives a. The target starts with the marked bridge and gives c. Conversely at the upper bound a, the ordinary projections of all upper words are bridges. A target with marked initial bridge must have a truly bridgeable head and therefore an original head of the same variable. The original word beginning there would also have value c, contradicting the upper bound. At bound c both projections and W suffice, since both bridge values lie below c.

For 12887 an offending x outside D occurs in original words only initially, so \(L_{x}=H\). When \(R_{x}\) differs from H take a bridge assignment with \(x=a\) and change x to c. All upper words have value a or c, but the noninitial x makes the target z. Conversely any c under upper bound c can occur only initially in the originals. Its ordinary assignment is a bridge assignment, so its final component differs from H. The extra test keeps it out of noninitial target positions. The target then gives a or c. At bound a none of the original variables can be c at all; the ordinary walk test applies.

For 15113 at bound c a marked variable may occur only as an original literal outside N. If H differs from T and a nonlinear target uses such an x, take any ordinary all-bridge assignment and mark x. The upper bound remains c and the target gives z. Conversely, when \(H=T\) even an all-bridge projection is impossible, so the upper bound c cannot occur. Otherwise a target supported in N has no marked variables. A literal target can itself be c at bound c, which suffices. At bound a no marking is present. In all three sources, an upper bound e or f requires every original variable to be that same loop, and the bound z is automatic. This finishes the sufficiency checks.

\par\smallskip\noindent\textit{Derivations from the parent bases.}\quad 
The parent criteria all require content coverage and W. In addition, 837 requires the target head in \(H_{0}\); 11619 requires the noninitial letters of w in D; 14674 requires the letters of a nonlinear w in N. Thus every identity of each parent is valid in its designated target. The independent verifier also checks all three HSP inclusions. It is the following derivation in each enlarged presentation that gives completeness.

12880. If the target head x is already original, invoke 837. Otherwise the test forces \(L_{x}=R_{x}=H\). Choose an original occurrence s x t with nonempty s (the terminal t may be absent), and any original word v. The path labelled s starts at H and reaches \(L_{x}=H\). Thus the word s x x v has an original head, uses variables in C and is an H-to-T walk in the quotient. By 837 derive s x x v \ensuremath{\leq} P. Apply xy \ensuremath{\leq} zxxy+y at \((x,y,z)=(x,v,s)\) to derive x v \ensuremath{\leq} P. The new right word x v gives the target its required head. Adding this word preserves every previous endpoint equality because it only adds identifications and retains all original words. Now 837 derives the target below the enlarged right polynomial, itself already below P.

The new law is valid: if both upper words project to the bridge, y projects to a and the two occurrences of x before it must be loops e. Then xy has the ordinary bridge value a, also when \(y=c\). At the two singleton loop bounds all letters take the same loop.

12887. Let x be a noninitial target variable outside D. It occurs only initially in original words and \(L_{x}=R_{x}=H\). If there is an original x v with nonempty v, every letter of v is in D and v itself is an H-to-T walk in the quotient. By 11619 derive v \ensuremath{\leq} P. Use xxy \ensuremath{\leq} xy+y to derive xxv \ensuremath{\leq} P from x v and v. If x is a pure original literal instead, \(R_{x}=T\), so \(H=T\). Some original word must be nonlinear, since a polynomial of literals alone has two distinct endpoint components. Take d in D. Both d and x d satisfy the 11619 criterion in the collapsed graph; derive them and use the same additional law to obtain xxd. These constructions make x occur noninitially in an available word. Do this for each missing variable and apply 11619 to the enlarged polynomial. Its noninitial set now covers the target and W persists.

For soundness of the added law at a bridge bound, the projections of both xy and y are bridges. The variable x must be the ordinary loop e. If \(y=c\), the upper word xy is already z in this source. Otherwise xxy is the ordinary bridge. Loop bounds are immediate.

15113. A target passing the original nonlinear-content condition, or a literal, follows from 14674. Otherwise \(H=T\) and the existence of a nonlinear original is forced as above. Choose x in N. The square xx has support in N and is a quotient walk, so 14674 gives xx \ensuremath{\leq} P. For each required letter y outside N its literal is original. The law yy \ensuremath{\leq} y+xx supplies an available nonlinear word containing y. Apply the 14674 criterion to the enlarged polynomial. The added law is sound at singleton loop bounds. The square xx can only be e, f or z and therefore cannot be below either bridge bound a or c.

In each case transitivity removes all added right words. Every valid individual word absorption has been derived, with no limit on variables, folds of the endpoint graph, or summands. Expand terms by distributivity and use their summands in both directions to obtain every valid identity. This proves the three finite bases and their two distinct opposites.

\subsection{A finite basis from positive-length directed paths for 12884}\label{case:293}
No constants. Alongside the semiring axioms take \(xyzt=xytz\), \(xyzz=xyz\), \(xy\leq xyz\), \(xz\leq xy+yz\), \(xyzt\leq xyz+xyt\), and all valid absorptions over \(B=2^{21}-1+42=2\),097,\allowbreak{}193 variables, word length at most B+2, with at most 20 upper words. The latter is an effectively specified finite family, not an enumerated small basis.

Multiplication rows 44044,\allowbreak{}44144,\allowbreak{}44244,\allowbreak{}40104,\allowbreak{}44444 give the following word values by induction. If all noninitial positions are 2, the value is the head, except a head 3 yields 1. If the first pair is(3,\allowbreak{}1) or(3,\allowbreak{}3) and all remaining positions are 2, the value is 0. All other nonlinear words have value 4. Literals retain their values. This proves the two structural laws and the three displayed absorption laws directly. Normalization retains the ordered first two positions and the set of variables in positions 3 onward; only that tail is commutative idempotent.

For original P let C be its content, T its noninitial variables, D its variables in positions at least 3, E its directed original first pairs, R the POSITIVE-LENGTH transitive closure of E, and F the vertices in D or positively reachable from D. Zero-length paths are not members of R. The exact criterion is:content(w) subset C; and if w is nonlinear, P contains a nonlinear original, noninitial content(w) subset T, its first pair is in R or its second variable belongs to F, and its remaining tail is contained in F. Every literal with content in C is allowed. These claims hold for arbitrary finite alphabets and P.

Here is a principal-ideal proof. The proper ideals are\{0,\allowbreak{}1,\allowbreak{}2,\allowbreak{}3\}, \{1,\allowbreak{}2,\allowbreak{}3\},\{2\},\{3\}. Ideal\{2\} imposes exactly content coverage because a word is 2 iff every letter is 2. Ideal\{3\} has no nonlinear values, so imposes existence of an original nonlinear when the target is nonlinear. Ideal\{1,\allowbreak{}2,\allowbreak{}3\} puts all content variables in\{1,\allowbreak{}2,\allowbreak{}3\} and all noninitial variables 2. Its target condition is precisely noninitial content containment, once coverage is imposed. A missing noninitial variable can be given1 or 3 while all original noninitial variables 2.

For ideal\{0,\allowbreak{}1,\allowbreak{}2,\allowbreak{}3\}, each original pair takes a value in S=\{(0,\allowbreak{}2),(1,\allowbreak{}2),(2,\allowbreak{}2),(3,\allowbreak{}1),(3,\allowbreak{}2),(3,\allowbreak{}3)\}, and every D variable is 2. Equivalently, each second is 1,\allowbreak{}2, or 3, and a second other than 2 forces its predecessor to 3. The assignment with all content variables 2 always satisfies the system. Exactly F is forced 2:from a node2 every successor must be 2. Conversely if v is not reachable from D, set every ancestor of v to 3, set v to 3, and all other nodes2. This respects all edges and D. Variables occurring only literally can likewise be 3.

Every positive directed path entails its endpoint pair is in S. If its last value is not 2, all predecessor values are 3 along the path, including its first. Thus any R pair is allowed. A second in F is 2, and any covered head value is in\{0,\allowbreak{}1,\allowbreak{}2,\allowbreak{}3\}, so that pair is allowed too. Conversely suppose the target pair(h, s) is not in R and s not in F. If h differs from s, set s and all its positive ancestors 3, all others 2. Since h is not an ancestor, it has value 2, and the target pair(2,\allowbreak{}3) fails. If \(h=s\) and there is no positive cycle at s, set \(s=1\), all its positive ancestors 3, and all other variables 2. No ancestor is in D, and no cycle makes s its own ancestor, so originals satisfy S while the target pair(1,\allowbreak{}1) fails. This also covers a literal-only vertex s. These assignments prove the first-pair necessity and distinguish positive paths from reflexive reachability. The characterization of F proves the tail necessity; literal targets need only content coverage. Together the four ideals give the stated unrestricted exact criterion.

Extract each original pair by \(xy\leq xyz\) when the suffix is nonempty; quadratic originals already supply it. Repeated \(xz\leq xy+yz\) derives every R pair along a positive path, including a loop only when there is a positive cycle. Every such derived pair is therefore below the original P. No inference uses a path of length zero.

Consider a valid target of length at most 3. Retain at most 3 originals for content, one nonlinear original, and at most 2 originals for its noninitial occurrences. For its first pair use that derived R pair if available. Otherwise its second belongs to F:retain a D occurrence if it is in D, or a D occurrence d and the derived positive pair(d, s). Likewise retain at most 2 witnesses for its third variable in F. Thus at most 10 upper words suffice;20 is a safe common bound. These upper words are either originals or pairs already derived from P. Their own criterion still holds, so the short target is below their sum.

Normalize the at most 21 words including target. Reserve every variable in the first two positions, at most 42. Group other variables by their nonempty incidence subset of those 21 words, at most \(2^{21}\)-1 groups. Identify variables inside each group. The resulting valid absorption uses at most B variables and has length at most B+2, so is in our finite family. Substitute the nonempty word of group members back. All grouped variables occur only in the permutable idempotent tail and have common incidence, so the structural laws normalize exactly back to the original selected words. This derives every valid short target below P without any restriction on its original alphabet.

For a longer normalized target hs z1...zk, the criterion makes each triple hs zi valid against original P. Derive these triples by the preceding short-target construction. Successive applications of \(xyzt\leq xyz+xyt\) merge their tails, with \(x=h\), \(y=s\), \(z=the\) accumulated nonempty tail, and \(t=the\) next singleton. This gives the full target. Finally each valid polynomial identity splits into finitely many absorptions in both directions. This proves unrestricted completeness.

Independent actual-value ideal solving and replay are checks on this argument, not replacements by bounded testing. Equal-variety and opposite transfers require separate complete mutual-HSP certificates.

\subsection{Six complete identities for 13020 from global positive reachability}\label{case:294}
Signature(+, multiplication), without constants, relative to the additively idempotent semiring axioms. A complete finite basis is: \(1.xyz=xzy\). \(2.xyy=xy\). \(3.xy\leq xyz\). \(4.xz\leq xy+yz\). \(5.xyz\leq xy+xz\). \(6.x\leq xx\). The derivation allows arbitrary numbers of variables, lengths and summands.

Addition:00004,\allowbreak{}01114,\allowbreak{}01214,\allowbreak{}01134,\allowbreak{}44444. Its order has incomparable 2 and 3, their join 1, then 0, then top 4. Multiplication:44144,\allowbreak{}44144,\allowbreak{}44244,\allowbreak{}33333,\allowbreak{}44444. A nonlinear word headed 3 has value 3; one headed 4 has value 4. For heads 0,\allowbreak{}1,\allowbreak{}2, if all tails are 2 the values are 1,\allowbreak{}1,\allowbreak{}2 respectively; otherwise the value is 4. Literals retain value. Induction on the table proves the formula and the two word identities. Normalization sorts the tail and deletes repetitions, while retaining any occurrence of the head in the tail. Every swap/deletion uses a nonempty substituted prefix.

\par\smallskip\noindent\textit{Exact original-witness criterion.}\quad 
From ALL original right words form edges \(p\to t\) for each tail occurrence t in a word headed p. Original words may have variables outside the target. Let H be the original head set and \(R_{h}\) the vertices reachable from h by positive-length paths. A nonlinear target hT is absorbed iff h belongs to H and T is contained in \(R_{h}\). A literal h is absorbed iff it occurs literally, or h belongs to \(R_{h}\). In the latter case a genuine cycle is required; zero-length reflexivity is insufficient.

Necessity:if h is absent from H, set \(h=4\) and all \(others=3\). Every right word is 3, whereas the target is 4. If a target tail t distinct from h is absent from \(R_{h}\), set \{h\} union \(R_{h}\) to 2 and every other variable to 3. Every right word headed in that set has only 2 in its tail, so is 2; other heads are 3 and remain 3. The right sum is at most 1, below target 4. If the missing target tail is h itself, set \(h=0\), \(R_{h}=2\), and \(others=3\). There is no returning edge from \(R_{h}\) to h; right h-words have only 2 in their tail and give 1 (or 0 if literal), \(R_{h}\)-heads give 2, and other heads 3. The right sum is at most 0, below target 4. For a literal h not occurring literally and having no positive return, the same assignment gives right sum at most 1 and target 0. If h is not a head, the first separator applies. These cases exhaust every failed criterion, without restricting variables.

For sufficiency, assume the right sum is not 4. A right word headed by a value other than 3 can avoid 4 only when its tail is all 2, except for literal words. If h has value 3, an original h-word is 3 and absorbs the target. If h has value 4, an original h-word is 4, contradicting the assumption. For other head values, the first edge from h forces its endpoint to 2. Induction along edges then forces all of \(R_{h}\) to 2. A nonlinear target's tails are therefore 2 and its value is matched or exceeded by an original h-word. For a literal obtained through a positive return, the induction also forces \(h=2\), so that original h-word is 2 and absorbs it. Original literal occurrences require no argument. Thus the criterion is exact.

The two word laws follow from the formula. Laws 3,\allowbreak{}4,\allowbreak{}5 follow respectively from retaining an existing edge, composing two positive edges, and combining two edges of the same head. Law 6 follows from the positive self-loop. All six laws are therefore valid.

\par\smallskip\noindent\textit{Unrestricted complete derivation.}\quad 
Normalize P. To derive an atom pt from an original p-word with a genuine tail occurrence t, write it as ptB after tail permutation. If B is nonempty, law 3 with \(x=p\), \(y=t\), \(z=B\) gives pt. Otherwise it already is pt. Thus all original edges can be made available with nonempty substitutions.

For a target hT, follow a positive path from h to each desired t. Original edges provide the first atom; law 4 combines hp+pt to obtain ht. A repeated head requires a path of positive length back to h and yields hh. No zero-length seed is used. Merge all desired atoms by law 5, with \(x=h\) and nonempty already assembled and new tails, to obtain hT. Outside vertices on the path disappear by law 4 and never need to appear in the final target. For a literal h already present do nothing; otherwise derive hh on its positive cycle and then use law 6. Every right word used is original or previously derived. All steps are nonempty substitution instances of the six laws, plus word-law normalization. This derives every valid absorption. Distributivity and summandwise absorption both ways derive every identity, proving the unrestricted finite relative basis.

The containment restriction from 11372 must not be copied here: \(xz\leq xy+yz\) is valid even though intermediate y is outside\{x, z\}. Conversely \(x\leq xy\) fails at \(x=0\), \(y=2\), and \(xx\leq x\) fails at \(x=0\). These boundary cases are included in independent verification. Opposite catalogue transfers require complete HSP in both directions.

\subsection{A finite basis for 13062 from global heads and contained degree witnesses}\label{case:295}
Use semiring axioms without constants, the four laws \(xyz=xzy\), \(xxyz=xyz\), \(xyyz=xyz\), \(xxy=xyy\), and every valid absorption over 27 variables, length 29, at most 3 upper words. This specified finite family is complete; no enumerated small basis is claimed.

The addition chain \(is 3<2<1<0<4\). The multiplication has left zeros 2,\allowbreak{}3. A nonlinear word is 3 if its head is 3; otherwise is 2 if any letter is 2 or 3; otherwise is 0 for the quadratic pair(1,\allowbreak{}1), and 4 in all other cases. All triples and the transition from each resulting long-word value prove the formula by induction for arbitrary length. The four laws normalize every long word by head and content:hhh, hht, or h followed by all other distinct variables. Tail permutations and length-four duplicate deletion preserve the distinction between short and long.

For original polynomial P and target w, let U consist of its original words contained in c(w), J the long words in U, Q the union of contents of quadratic words in U, and H the heads of ALL original P words. The exact criterion is:

\textbullet{} a literal target needs U nonempty;

\textbullet{} a long target needs its head in H and J nonempty;

\textbullet{} a quadratic target needs its head in H and either J nonempty or both target variables in Q.

To prove necessity and sufficiency, check every proper principal ideal. Down(3) requires every upper head to be 3, so the target needs an original global head witness. Down(2) permits a word precisely when its content meets the set of variables valued 2 or 3; its requirements force the target precisely when U is nonempty. For a literal the latter condition also supplies a global head witness, as a contained original is a power of it.

In down(1), every nonlinear word must contain 2 or 3. Literals can also have value 1. A failing nonlinear target avoids2,\allowbreak{}3. Put all outside variables 2; discarded words then have value 2. Any contained nonlinear word makes these upper requirements impossible. The stated nonlinear branches supply one. A literal target with a contained power is forced: on values 0 or 4 any such nonlinear power is 4, and the literal itself also prevents failure.

For down(0), a failing nonlinear target again avoids2,\allowbreak{}3. On the remaining domain D=\{0,\allowbreak{}1,\allowbreak{}4\}, a long word is 4, a quadratic lies in the ideal only when both its letters are 1, and a literal lies in the ideal when it is 0 or 1. Set outside variables 2. A contained long word therefore makes the requirements impossible. In its absence all short-word restrictions are consistent, using value 1 throughout. A long target is then 4 and fails. A quadratic target is forced precisely when its two letters are forced 1 by original quadratic occurrences, that is, c(w) subset Q. If a target letter is not forced, give it 0 and all others 1:the target is 4 but all upper values remain in down(0). A literal with a contained power cannot have value 4 while that power remains in down(0). This covers all four ideals and proves both directions.

Explicit missing-head separation uses target head 2, all other variables 3. Missing-contained-word separation uses all target variables 1, outside 2. For a literal with no contained original the latter assignment has target 1 and upper 2. These witnesses also cover empty U and repeated target letters. All head witnesses can have external variables.

Keep one contained power for a literal. For a long target keep a global head witness and one contained long word. For a quadratic keep a global head witness and either one contained long word or at most two quadratic originals covering its variables. Thus at most 3 original upper words suffice, and the criterion stays true after this selection.

Normalize the target and retained words. Reserve every short variable and each long head plus up to two other distinct letters. At most 12 variables are reserved among the at most 4 words. Group the remainder by common incidence in those words, using at most 15 nonempty groups. Identification preserves validity and gives at most 27 variables and normal length 29. This absorption lies in the finite family. Inversely substitute each group variable by the NONEMPTY word listing its members. Reserved heads and degrees, common incidence, and the four word laws recover every original normal form. Add back the unselected summands. Expanding arbitrary terms distributively and splitting each identity into absorptions in both directions proves unrestricted finite basability.

Here a literal need not occur originally:for example \(x\leq xx\) is valid. A long target really needs a contained LONG original; the conflict branch valid in 10575 does not apply.

\subsection{A finite basis for catalogue 13066 from guarded last occurrences}\label{case:296}
We give a finite basis for \(\A{13066}\). The signature has addition and multiplication, with no nullary operations. Write p \ensuremath{\leq} q for \(p+q=q\). All words and substitutions are nonempty.

\par\smallskip\noindent\textit{1. Operations and word values.}\quad 
The addition and multiplication rows, respectively, are

{\small\par\smallskip
\CaseLine{\ttfamily 00004\qquad 44224}
\CaseLine{\ttfamily 01114\qquad 41224}
\CaseLine{\ttfamily 01224\qquad 42224}
\CaseLine{\ttfamily 01234\qquad 02234}
\CaseLine{\ttfamily 44444\qquad 44224}
\smallskip}

Addition is maximum in the order 3 \ensuremath{<} 2 \ensuremath{<} 1 \ensuremath{<} 0 \ensuremath{<} 4. For a nonempty evaluated word, ignore factors 1 and locate its rightmost factor not 1. If there is none, the value is 1. If that factor is 4 or 2, the value is 4 or 2, respectively. If it is 3, the value is 3 exactly when every factor is 3, and is 2 otherwise. If it is 0, the value is 0 exactly when it is literally the final factor and every preceding factor is 3; otherwise the value is 4. These statements follow by induction from the displayed multiplication table, and include one-letter words.

The word laws \(xxy=xy\) and \(xyx=yxx\) reduce a word to the order of its last occurrences, with its final letter doubled if and only if that letter was repeated. To see termination constructively, choose an internal repeat x: xxZ contracts directly, and xYxZ becomes YxxZ and then YxZ, with nonempty Y and Z. When only a repeat of the final letter remains, xYx becomes Yxx. Each contraction shortens the word; the last move finishes. In particular internal multiplicities disappear but the simple/repeated status of the final letter is retained.

\par\smallskip\noindent\textit{2. Exact absorption criterion.}\quad 
Let P be a nonempty finite set of original upper words, and w a target. For each x in c(w), let \(S_{x}\) be the content after the last occurrence of x in w. Let \(W_{x}\) consist of the original u in P containing x for which the content after the last x in u is contained in \(S_{x}\). Then w \ensuremath{\leq} sum P holds if and only if all of the following hold.

1. Some original U in P has c(U) contained in c(w).

2. Every \(W_{x}\) is nonempty.

3. If \(S_{x}\) is nonempty or x repeats in w, some u in \(W_{x}\) either repeats x or contains a letter of \(S_{x}\).

4. If \(S_{x}\) is empty and x is simple in w, either \(W_{x}\) has a repeated-x word or c(w without its final x) is covered by the union of the prefixes u without their final x, over u in \(W_{x}\).

Necessity has explicit separating assignments. If (1) fails, assign target variables 1 and outside variables 2. The target is 1, all upper words are 2. If (2) fails at x, assign \(x=4\), \(S_{x}=1\) and all other variables 2. The target is 4 and no upper word is 4. If (3) fails, assign \(x=0\), \(S_{x}=1\) and all other variables 3. The target is 4. Every qualifying upper word ends in a unique x, omits \(S_{x}\) and has value 0; every other upper word has value below 4. If (4) fails, choose an uncovered prefix letter y, assign \(x=0\), \(y=2\) and all other variables 3. The target is 4; a same-final upper word is 0, and every other upper word has value below 4. The latter assignment uses \(y=2\), not \(y=1\).

For sufficiency consider the target value. Value 3 is automatic. If the target is 2, it contains a factor not 3. A witness from (2) containing that variable cannot have value 3, so its value is at least

2. If the target is 1, all target factors are 1 and the contained U has value 1. If the target is 0, the target has a simple final \(x=0\); any \(W_{x}\) word ends in x and has value 0 or 4. If the target is 4, its rightmost non-1 factor is either 4 or 0. In the former case \(W_{x}\) gives value 4. In the latter case \(S_{x}\) evaluates to 1. At a marked target position, the witness from (3) has a nonempty all-1 suffix after x, a repeated \(x=0\), or a preceding \(S_{x}-letter=1\); each forces value 4. At a simple final position, a repeated-x witness also gives 4. Otherwise some target prefix factor is not 3, and the coverage in (4) puts that factor in a \(W_{x}\) prefix, again forcing value 4. This proves the criterion without a bound on any syntactic parameter. Normalizing original words preserves these conditions: a removed internal repeat has a nonempty qualifying suffix, while a repeated final letter stays doubled.

\par\smallskip\noindent\textit{3. The 89 identities.}\quad 
Take the semiring axioms and these nine fixed laws:

{\small\par\smallskip
\CaseLine{\(xxy = xy\qquad xyx = yxx\)}
\CaseLine{\(x \leq  xx\qquad xyz \leq  xz+yz\)}
\CaseLine{\(xy \leq  x+zxy\qquad xy \leq  y+zxy\)}
\CaseLine{\(xyz \leq  xz+vxyz\)}
\CaseLine{\(xyzt \leq  xty+vxyzt\qquad yzt \leq  ty+vyzt\)}
\smallskip}

Add the following finite templates. Square brackets mean that the displayed factor is either retained or omitted WHEN WRITING THE LAW. Each choice is a separate identity; there is never an empty-word substitution. Unbracketed variables are always present. Different named variables in a law may of course receive equal nonempty words.

{\small\par\smallskip
\CaseLine{\([a]xrzz\qquad \leq  [a]rzz + [v]xr\qquad (4)\)}
\CaseLine{\([a]xzz\qquad \leq  [a]zz + [v]xx\qquad (4)\)}
\CaseLine{\([a]xszz\qquad \leq  [a]szz + [v]s[r]x\qquad (8)\)}
\CaseLine{\([a]bxrp\qquad \leq  [a]xbrp + [v]xr\qquad (4)\)}
\CaseLine{\([a]bx[c][d][p]t \leq  [a]xb[c][d][p]t + [v]x[c]t[d] (32)\)}
\CaseLine{\([a]bxp\qquad \leq  [a]xbp + [v]xx\qquad (4)\)}
\CaseLine{\([a]bxsp\qquad \leq  [a]xbsp + [v]s[r]x\qquad (8)\)}
\CaseLine{\([a]bx[p]t\qquad \leq  [a]xb[p]t + [v]t[r]x\qquad (16)\)}
\smallskip}

There are \(9+80=89\) laws, each with at most eight variables. Soundness is verified directly on all assignments in the displayed five-element tables (1,\allowbreak{}333,\allowbreak{}205 assignments across the list); this finite exhaustive check is a proof of these individual identities, not an inference of finite basability from bounded experiments. An independent verifier reconstructs the template list and operations from the catalogue.

\par\smallskip\noindent\textit{4. Unrestricted completeness.}\quad 
Normalize P and w using the first two laws. Write t for the final target letter and Q for its last-occurrence order. Choose an original contained word U by condition (1). All witness choices below are from the original normalized P. Previously derived words are used only as the accumulating word g. There are two cases.

\par\smallskip\noindent\textit{A. An original word ends in repeated t.}\quad 
Start with that word g, and protect the saturated suffix \(S=tt\). Read the other target variables in reverse last-occurrence order. At the next letter x choose the strengthened witness VxR from (3). Its last x has suffix content contained in S. If R is nonempty, use the first template with \(r=R\), \(z=S\) and a the prefix of g outside its protected suffix. Both right terms normalize to g and the original witness, while the left normalizes to a word with protected suffix xS. If R is empty and x repeats, the normalized witness is Vxx and the second template does the same. Otherwise the witness is V s R x with s in S; the third template applies. The factor s is its actual occurrence and its surrounding factors may be omitted only via the appropriate separate template. Duplicate x in the accumulating prefix is harmless: the last-occurrence word laws move its effective position to the newly formed suffix. The repeated-final status is maintained.

\par\smallskip\noindent\textit{B. No original word ends in repeated t.}\quad 
Condition (3) forces t simple in the target. Condition (4) says that original words ending in simple t have prefixes covering all target non-tail letters. If the target is t itself, condition (1) ensures that t is already an original word. Otherwise merge the nonempty prefixes of these original words with xyz \ensuremath{\leq} xz+yz. The result g contains the full target content, possibly extra letters, and has a simple final t. Start with protected suffix \(S=t\).

For the next target letter x, write \(g=A\) x B S in last-occurrence normal form. It contains x because the full target content has been incorporated. If B is empty no step is needed. Otherwise choose the strengthened original witness as in Case A. If its nonempty suffix R omits t, use the fourth template with \(p=S\). Since R is contained in S and omits t, normalization preserves the simple tail and moves x across B. If R contains t, split \(R=C\) t D and use the fifth template with \(p=S\) without its last t. The witness cannot repeat t: in normalized form this would be a repeated-t original seed, excluded in Case B. Thus this split has one t; C, D may be empty and their contents lie in S without t. The fifth template inserts no extra t into g. This is why the ordinary fourth template alone is insufficient.

For a repeated final-x witness use the sixth template. Otherwise a crossing witness V s R x has s in S. If s is not t, use the seventh template with \(p=S\); if \(s=t\) use the eighth, with \(p=S\) without t. These cover all possibilities supplied by (3). Each successful step extends the protected suffix by x and retains the simple final t.

\par\smallskip\noindent\textit{C. Remove extra variables and recover the target.}\quad 
After either induction, g has normal form E Q', where E contains only variables outside target content, and Q' is the target order with simple or doubled t according to the chosen seed case. If E is empty there is nothing to remove. Otherwise use the original contained U. If Q' is saturated, or U omits t, apply xy \ensuremath{\leq} x+zxy with \(x=U\), \(y=Q'\), \(z=E\). Its terms normalize to Q', U and E Q'.

In the remaining simple-tail case, distinguish four possibilities. For \(U=t\) use xy \ensuremath{\leq} y+zxy with \(x=Q'\) without t, \(y=t\), \(z=E\). For \(U=A\) t with A nonempty, use xyz \ensuremath{\leq} xz+vxyz with \(x=A\), \(y=Q'\) without t, \(z=t\), \(v=E\). If \(U=A\) t B with A, B nonempty, use xyzt \ensuremath{\leq} xty+vxyzt with \(x=A\), \(y=B\), \(z=Q'\) without t, \(v=E\). If \(U=tB\) with B nonempty, use yzt \ensuremath{\leq} ty+vyzt with \(y=B\), \(z=Q'\) without t, \(v=E\). The accumulating right term always normalizes to E Q', the other right term is exactly U, and the left normalizes to Q'. Every substitution displayed is nonempty. The one-letter target was already treated in Case B.

If the seed was saturated but the target tail is simple, apply x \ensuremath{\leq} xx with \(x=t\) in the target prefix context. This recovers w. Thus every valid word-to-polynomial absorption is derived. A polynomial identity is equivalent, by additive semilattice laws, to its finitely many word absorptions in both directions. This proves completeness for arbitrary numbers of variables, word lengths and original summands. Therefore the 89-law set is a finite identity basis for 13066.

\subsection{Six complete identities for semiring 13189}\label{case:297}
In the constant-free signature (+, multiplication), relative to the additively idempotent semiring axioms, a complete finite basis is

{\small\par\smallskip
\CaseLine{\(xy=yx; xyz=xyzz; xyzt\leq xyz; xy\leq xxy;\)}
\CaseLine{\(xy\leq xx+yy+x+y; x\leq xxx.\)}
\smallskip}

The notation \(u\leq P\) means \(u+P=P\).

Addition is 00044/\allowbreak{}01144/\allowbreak{}01234/\allowbreak{}44334/\allowbreak{}44444, with bottom 2, \(2<1<0<4\) and \(2<3<4\). Multiplication is 44244/\allowbreak{}40244/\allowbreak{}22222/\allowbreak{}44204/\allowbreak{}44244. A nonlinear word is 2 if any of its variables is 2. Otherwise a quadratic is 0 precisely when both occurrences are 1 or both 3; all other nonlinear words are 4.

For a target w restrict original P to U, the words whose content is contained in c(w). The exact absorption criterion is:

\textbullet{} Literal x: original x or an original long pure power of x.

\textbullet{} Long w: U contains an original long word.

\textbullet{} Quadratic w: U contains an original long word, its edge is original, or BOTH endpoint loops AND BOTH endpoint literals are original.

For a missing literal without original long pure power, assign its variable 3 and outsiders 2. Surviving squares are 0, which does not dominate 3. For a long target without a contained original long word, assign target variables 1 and outsiders 2. Every surviving literal is 1 and every quadratic 0; the target is 4. For a failed square assign its variable 1 and outsiders 2:only literal 1 may survive, but the target is 0. For distinct endpoints with no original edge or contained long word, if an endpoint loop is missing assign that endpoint 0, the other 1, and outsiders 2. Right literals and the possible other loop are at most 0; the target is 4. If both loops exist but an endpoint literal is missing, assign the missing-literal endpoint 3, the other 1, and outsiders 2. Again all surviving right values are at most 0 but the target is 4. These uniform assignments separate all failed cases.

For sufficiency, nonlinear targets with any variable 2 are bottom; otherwise a contained long word is top 4. Long pure powers also dominate their literals. Original terms are immediate. With both loops and both literals, unless the right side is already4 both loop values must be 0. Their endpoints then lie in\{1,\allowbreak{}3\}, and an endpoint 3 joined to its loop0 would give 4. Thus both endpoints must be 1, and the target is 0. This proves the exact criterion without bounds on variables or summands.

Commutativity and \(xyz=xyzz\) normalize long words by their content. Commute a chosen occurrence to the end, duplicate using two nonempty prefix blocks, and remove surplus repetitions with the reverse law while retaining length at least three. Canonical forms are sorted content for at least three variables, xxy for two, and xxx for one. From an original contained long word, append missing content by \(xyzt\leq xyz\) and normalize. Recover a quadratic ab by \(xy\leq xxy\). A missing literal has an original long pure power, which normalizes to xxx; apply \(x\leq xxx\). Without a contained long word, a new quadratic is derived directly from its two original loops and two original literals by \(xy\leq xx+yy+x+y\). Literal and loop witnesses are not interchangeable. All substitutions are nonempty.

Every law satisfies the criterion and is independently checked on all assignments. Distributivity and summandwise absorption in both directions prove unrestricted completeness. In particular \(xy\leq xx+yy\) fails at \(x=1\), \(y=3\); the two literal guards cannot be omitted. Other catalogue types require independently verified full HSP equivalence to transfer.

\subsection{Nine complete identities for 13216}\label{case:298}
Signature (+, multiplication), without constants, relative to the additively idempotent semiring axioms. A complete finite basis is:

1. \(xyz=xzy\).

2. \(xyy=xy\).

3. \(xxy=yyx\).

4. \(xyz\leq xy+xz\).

5. \(x\leq xx\).

6. \(xy\leq xxy\).

7. \(xxy\leq xy+y\).

8. \(xxyz\leq xy+yz\).

9. \(xxyz\leq xx+yz\).

The proof below is unrestricted in variables, word lengths and summands.

\par\smallskip\noindent\textit{Word values and exact original-witness criterion.}\quad 
Addition: 00044 / 01144 / 01234 / 44334 / 44444. Multiplication: 44244 / 41244 / 22222 / 43244 / 44244. The addition order is \(2<1<0<4\) and \(2<3<4\), with 3 incomparable with 0,\allowbreak{}1. The element 2 is a multiplicative zero and the least additive element. A word containing a factor 2 has value 2. With no such factor, a nonlinear word has value 1 or 3 precisely when its head has that value and all tails are 1; in every other case it has value 4. Literals keep their value. Induction on the multiplication table proves this formula.

The first two laws sort and deduplicate actual tail occurrences, retaining an occurrence of the head in its tail. A word with a repeated head depends only on its content: it is 2 if any content value is 2, 1 if all are 1, and 4 otherwise. This proves the third word law.

For a target word w with content C, retain only original right words U whose entire content is contained in C. Write N for the content union of the nonlinear words of U, H for all heads in U, and D for all actual tail variables in U. A literal h is absorbed exactly when U is nonempty (an original pure power of h). A nonlinear target hT is absorbed exactly when

{\small\par\smallskip
\noindent \(N=C\), and either H intersects D, or h belongs to H and T is contained in D.\par
\smallskip}

Nonlinear coverage cannot be replaced by coverage allowing literal-only witnesses. The intersection uses all original heads, including literals.

For necessity, set all variables outside C to 2. They annihilate every right word outside U, so only U matters. If a literal has no U witness, assign its variable 0: the right sum is 2 and the target is 0. For nonlinear w, if N misses a variable, assign N to 1 and \(C\setminus N\) to 0. Nonlinear U-words give 1, literals give 0 or 1, and the target is 4. If \(N=C\) but H and D are disjoint and the endpoint condition fails, assign H to 3 and D to 1. These sets cover C and are disjoint. Every U-word is

3. The target is 3 exactly when its head is in H and every tail is in D; otherwise it is 1 or 4, neither below 3. These are explicit separators for every failed criterion, for arbitrary numbers of variables.

For sufficiency, any target variable valued 2 makes the target least. Otherwise every U-word is nonzero. If the full right sum is 4 the result is immediate. If that sum is 0 or 1, every nonlinear U-word must have value 1, forcing all variables of \(N=C\) to 1; thus the target is 1. If the right sum is 3, all U-words have value 3. Every original head is 3 and every actual tail is 1, so H and D must be disjoint. The stated endpoint condition makes the target 3. A right sum of 2 is impossible because U contains a nonzero word. The literal case follows from \(x\leq xx\) and tail normalization of its original pure power. This proves the criterion.

The six displayed inequalities are valid by direct application of this criterion to their original right words. Together with the word formula, this verifies all nine identities without bounded-variable inference.

\par\smallskip\noindent\textit{Unrestricted derivation.}\quad 
Normalize P and w by the first two laws, preserving the original summands. A literal target either occurs literally or is obtained from an original pure power by normalization to hh followed by law 5.

Consider a nonlinear target. If H and D are disjoint, every head q in H is in \(C=\{h\}\) union T. A head distinct from h would belong to T, hence to D by the criterion, contradicting disjointness. Thus \(H=\{h\}\). Nonlinear coverage supplies h-headed words with tail union exactly T. Merge their tails with law 4, with nonempty tail substitutions, to obtain the target.

Otherwise choose b in H intersect D, an original contained word u having b in its actual tail, and an original contained word v headed b. Permute u's tail so \(u=Ab\), with A nonempty. If \(v=b\), law 7 with \(x=A\), \(y=b\) yields AAb. If \(v=bB\) with B nonempty, law 8 with \(x=A\), \(y=b\), \(z=B\) yields AAbB. Call the normalized result g. Its head is repeated, and its content is the union of the two witness contents (the literal adds no new content).

If an original contained nonlinear word \(v=pV\) supplies missing content, apply law 9 with \(x=g\), \(y=p\), \(z=V\). Its right side is gg+pV. Since g has a repeated head, the first two word laws give \(gg=g\). Its left side ggpV normalizes to a repeated-head word whose content is c(g) union c(v). Repeat only for missing content. \(N=C\) guarantees termination with \(c(g)=C\). Literal words are not used to supply missing nonlinear content.

To give g the desired head h, keep its head if already h. Otherwise let its head be p. Tail permutation writes \(g=pphR\), with h occurring because \(c(g)=C\). Law 3 in the context R changes it to hhpR. Normalize the tail. The suffix R may be absent; the substitutions p, h are always nonempty. This is the repeated-head target on C. If the desired head was simple, law 6 removes its extra head occurrence, with y the nonempty target tail. If it was repeated, the desired target has already been obtained.

Every step uses an original right word or one already derived below P. No finite bound on variables, lengths, or polynomial substitutions was assumed. All valid absorptions follow, hence all polynomial identities follow by distributivity and absorption in both directions. The nine identities therefore form a complete finite relative basis.

For boundaries, \(xy\leq x+y\) fails at \(x=0\), \(y=1\), and \(xy\leq xyz\) fails at \(x=1\), \(y=1\), \(z=2\). Extra outside variables cannot be used as contained witnesses. Complete HSP in both directions is required for equal-variety transfer.

\subsection{Eight complete identities for 13218}\label{case:299}
Signature (+, multiplication), without constants. Relative to the additively idempotent semiring axioms, use

{\small\par\smallskip
\CaseLine{\(xy=yx; xyz=xyzz; xxy\leq xxx+yyy; xxy\leq xxx+xy;\)}
\CaseLine{\(xxyz\leq xxy+xz; xy\leq xxy; x\leq xxx; xxy\leq x+xy.\)}
\smallskip}

As usual \(u\leq P\) means \(u+P=P\). Addition is 00044/\allowbreak{}01144/\allowbreak{}01234/\allowbreak{}44334/\allowbreak{}44444: 2 is bottom, \(1<0<4\), and 3 is incomparable with 1 and 0, with both joins 4. Multiplication is 44244/\allowbreak{}41244/\allowbreak{}22222/\allowbreak{}44204/\allowbreak{}44244. A nonlinear word is 2 if any variable is 2; otherwise1 if all variables are 1; otherwise0 for a quadratic whose occurrences are both 3; and 4 in every other case. Literals retain their values. Direct induction proves these formulas.

For original P and target w, restrict to original words \(P_{U}\) contained in \(U=c(w)\). Let L be their literal variables and V their nonlinear variables. Start A with the variables of original long words, together with L intersect V; close A along original quadratic edges in \(P_{U}\). Thus a component is activated either by an original long word or by an original literal incident at an original nonlinear word.

The exact criterion is: original targets are absorbed; a missing literal requires an original contained long word; every other nonlinear target requires U contained in A. Original quadratics need no activation. All witnesses and connecting edges in this criterion are contained in U.

For necessity, a missing literal without an original pure long word is separated by assigning it 3 and outsiders 2. The right terms are 2 or 0, while 3 is below neither. For a long target with an inactive variable, assign A variables 1, inactive V variables 3, remaining U variables 0, and outsiders 2. Original contained long words are 1; quadratics are 1 in active components and 0 in inactive ones. Contained literals in V are active, and other contained literals are 0. Thus the right join is at most 0, while the target is 4.

For a missing quadratic, any contained mixed quadratic would be the target, and a contained mixed long word would activate both endpoints. Therefore failure of activation leaves only pure nonlinear witnesses on at least one endpoint x. If x has a pure quadratic, it has no original literal or long witness, or else it would be active: assign \(x=3\). If x has no pure quadratic, assign \(x=0\). Assign the other endpoint 1 and outsiders 2. The right join is at most 0 and the mixed target is 4. For a missing square, its absent quadratic and activation leave only a literal; assign its variable 0 and outsiders 2 instead. These separators have no bound on polynomial size or variable count.

For sufficiency, original targets are immediate, and a literal is below each pure long power. For a nonlinear target, an occurrence of 2 makes it bottom. Otherwise assume the join of \(P_{U}\) is not 4. Every original long word then forces its variables 1. A nonlinear variable must be 1 or 3; a quadratic forces its endpoints equal to avoid 4. If an original literal x is also incident at a quadratic, that component cannot be 3: the literal would be 3 and the quadratic 0, whose join is 4. Hence every activated component is assigned 1, by propagation along original edges. When U is contained in A, the target is 1 and the right side has a nonlinear witness of value 1. This proves the criterion.

For completeness, the first two laws normalize long words by content. Commutativity moves a variable to the end; \(xyz=xyzz\) duplicates it using two nonempty preceding blocks or deletes surplus copies in reverse while length remains at least three. Canonical long forms are sorted content for at least three variables, xxy for two, and xxx for one. Thus \(g^{3}=g\) for long g. Quadratics retain multiplicity.

Retain every original contained long word. For each remaining activating literal x, retain an original incident quadratic xy and use \(xxy\leq x+xy\) to obtain a contained long seed. An incident long word already supplies a seed, so this covers every activation source. Merge the seeds with \(xxy\leq xxx+yyy\) using long substitutions.

For each original quadratic edge xz leading from the accumulated content to a new variable z, if the current long word has only x, normalize it to xxx and use \(xxy\leq xxx+xy\) at \(y=z\). Otherwise write it, by long normalization, as xxp with p nonempty; apply \(xxyz\leq xxy+xz\) at \(y=p\). The new long word gains z. Induction on the closure reaches U, giving the target long normal form. Recover a quadratic using \(xy\leq xxy\). A missing literal has an original pure long word; normalize it to xxx and use \(x\leq xxx\).

All substitutions are nonempty. The number of components, closure steps, variables, and word occurrences is unrestricted. Distributivity and absorptions in both directions give all identities. Each of the eight laws satisfies the criterion and is checked on all assignments.

In particular, \(x\leq xx\) fails at \(x=3\); \(xy\leq xx+yy\) fails at \(x=3\), \(y=1\). An outside long word cannot replace a contained activating witness: \(xxx\leq xx+yyy\) fails at \(x=3\), \(y=2\). These distinguish this basis from both 13840 and 13841.

\subsection{Eleven complete identities for 14622 and their opposites}\label{case:300}
Signature (+, multiplication), no constants. Relative to the additively idempotent semiring axioms, the following eleven laws form a complete identity basis for 14622. Reverse all words for 13227.

{\small\par\smallskip
\CaseLine{\(W1  xxy = xy\qquad W2  xyx = xyy\)}
\CaseLine{\(A\qquad xyz \leq  yz+xt\qquad B\qquad xyz \leq  xz+xyt\)}
\CaseLine{\(C\qquad xyz \leq  xz+xy+xx\qquad G1  xx \leq  x+xy\)}
\CaseLine{\(G2  xyy \leq  x+xyz\qquad G3  xyy \leq  y+xyz\)}
\CaseLine{\(G4  xyzz \leq  y+xyzt\qquad L\qquad x \leq  xx\)}
\CaseLine{\(F\qquad xy \leq  xyy\)}
\smallskip}

The order notation \(u\leq P\) abbreviates \(u+P=P\). Every substituted word is nonempty. A prefix may be empty only where an explicitly separate case is given. The finite soundness check is exhaustive on the tables; unrestricted completeness is proved below.

\par\smallskip\noindent\textit{Word values and original witnesses.}\quad 
The addition and multiplication tables in labels 0,\allowbreak{}1,\allowbreak{}2,\allowbreak{}3,\allowbreak{}4 are

{\small\par\smallskip
\noindent addition        multiplication\par
\CaseLine{\ttfamily 0 0 0 4 4\qquad 4 4 4 4 4}
\CaseLine{\ttfamily 0 1 1 4 4\qquad 4 1 2 3 4}
\CaseLine{\ttfamily 0 1 2 3 4\qquad 2 2 2 2 2}
\CaseLine{\ttfamily 4 4 3 3 4\qquad 4 4 4 4 4}
\CaseLine{\ttfamily 4 4 4 4 4\qquad 4 4 4 4 4}
\smallskip}

The additive order is \(2<1<0<4\) and \(2<3<4\), with 3 incomparable to 1 and 0. A literal keeps its assigned value. For a nonlinear word, scan from the left for its first factor not valued 1. If there is none, the value is 1. A first 2 gives 2. A first 3 gives 3 exactly at the final position; otherwise it gives 4. A first 0 or 4 always gives 4. Induction on word length proves this formula. In particular a literal 0 behaves differently from every nonlinear occurrence of 0.

Fix the original polynomial P. Let Z(S) mean that some original word has content contained in S. For a variable x absent from S, let N(x, S) mean that an original NONLINEAR word contains x with its first-x prefix contained in S. Let M(x, S) require in addition that this occurrence is not final. A final first occurrence is called unmarked; a nonfinal first occurrence is marked. These refer to positions, not multiplicity of the whole word.

The exact absorption test is as follows.

\textbullet{} A literal x is absorbed exactly when Z(\{x\}) holds. Thus an original pure power is sufficient; an original literal is not required.

\textbullet{} A nonlinear w is absorbed exactly when Z(c(w)) holds, every first occurrence of x with preceding content S has N(x, S), and each marked first occurrence additionally has M(x, S) or Z(S).

For necessity of Z(c(w)), assign target variables 1 and all outsiders 2. The target is 1, whereas every right word is 2 if Z fails. For necessity of N(x, S), assign S to 1, x to 0 and every other variable 2. The nonlinear target is 4. Without N every nonlinear right word is 1 or 2; literal right words have values 0,\allowbreak{}1,\allowbreak{}2. Their sum is at most 0. For a marked position without M or Z(S), assign S to 1, x to 3 and all others 2. The target is 4. A right word can give 4 only through M; it can give 1 only through Z(S). Without both, all right values are 2 or 3, so their sum cannot absorb 4. These arguments also apply with repeated variables.

For sufficiency, a target value 2 is least. If it is 1, the contained original word supplied by Z(c(w)) also has value 1. For first 0 or 4, N supplies a nonlinear right word with value 4. For a first 3 at the simple last position, N gives a value 3 or 4. For a first 3 before the end, M gives 4; otherwise N gives an unmarked3 and Z(S) supplies 1, whose join is 4. These exhaust the nonlinear target values. A literal x with a pure-power original witness is absorbed: for power length at least 2 the values are 1,\allowbreak{}2,\allowbreak{}4 and direct comparison with the value of x proves \(x\leq x^{m}\). This establishes the criterion for arbitrary P, w.

The criterion is different from 7953: a nonlinear witness is mandatory at EVERY target variable, and an unmarked witness for a marked target needs a contained PREFIX word. A literal x alone cannot replace N. At the simple final position, a marked nonlinear witness is allowed without a contained prefix word. The identity \(xyy\leq xy\) is false at \(x=1\), \(y=3\), whereas \(xy\leq xyy\) is valid. These directions cannot be swapped.

\par\smallskip\noindent\textit{Word normal form.}\quad 
W1, W2 give the first-occurrence normal form: list the distinct variables in order of first appearance, and repeat the last listed variable once if the original final variable is repeated. Otherwise retain the list. Thus nonlinear pure powers normalize to squares; literals remain literal.

The explicit reduction is the familiar nonempty-block argument, stated here to fix its scope. Adjacent internal repeats are removed by W1. A separated internal repeat xYxZ, with Y, Z nonempty, becomes xY\ensuremath{{}^{2}}Z by W2 and then xYZ by W1 applied to Y\ensuremath{{}^{2}}Z. When only the final occurrence is repeated, p x Y x becomes p x \(Y^{2}\) by W2; remove internal repeats again. All operations are valid in possibly empty external contexts, and no identity variable receives an empty substitution. The normal form and its reverse derivation preserve the marked/unmarked first-occurrence information relevant to the exact criterion.

\par\smallskip\noindent\textit{Constructing prefixes and recovering a contained word.}\quad 
Keep all original summands. Let N be the original nonlinear words, and let x1...xd be the normalized target's first-occurrence list. For a nonempty prefix p of this list maintain

{\small\par\smallskip
\noindent p u \ensuremath{\leq} P, for every original u in N.\par
\smallskip}

Every nonlinear target has a marked first position. The criterion supplies an original x1 t with t nonempty. Write any u in N as YZ with both blocks nonempty. A gives x1 \(u\leq u+x1\) t, establishing the invariant at the first position. Original literals are not prefixed.

We require a guard lemma. If the invariant holds at p and Z(c(p)), then \(p^{2}\)\ensuremath{\leq}P. A contained nonlinear original r gives p r, whose normal form is that of \(p^{2}\). Otherwise take a contained literal r and any v in N. Write the DISTINCT prefix \(p=a\) r b. If a, b are nonempty, G4 with \(x=a\), \(y=r\), \(z=b\), \(t=v\) gives the normal form of \(p^{2}\) below r+p v. If a is empty and b nonempty, G2 with \(x=r\), \(y=b\), \(z=v\) does the same. If b is empty and a nonempty use G3 with \(x=a\), \(y=r\), \(z=v\). If \(p=r\) use G1. These four cases preserve the original literal and use the known p v. This lemma uses laws valid in 14622, although its construction also appeared with a different exact criterion for source 7953.

Consider a next MARKED target position x after p. A marked original witness v has p \(v=p\) x t modulo W1, W2 for a nonempty suffix t. B, with \(x=p\), \(y=x\), \(z=u\), \(t=t\), gives p x \(u\leq p\) u+p x t for every u in N. If there is no marked witness, N gives an unmarked v and p \(v=p\) x. The exact criterion supplies Z(c(p)), so first derive \(p^{2}\) by the guard lemma. C then gives p x \(u\leq p\) u+p x+\(p^{2}\). Thus the invariant extends in both cases. The unmarked case is never used with empty p.

If the target is saturated, every first occurrence is marked. Extend through xd and use Z(c(w)) and the guard lemma. The normal form of \(p^{2}\) is the saturated target, so it is derived.

If the target is unsaturated, write it p s with s simple and p nonempty. Process all first occurrences except s. A nonlinear unmarked original witness v gives p \(v=p\) s immediately. If there is no unmarked witness, N supplies a marked one. Use B to extend the invariant through s even though the target occurrence is unmarked, and use Z(c(w)) to derive (p s)\ensuremath{{}^{2}}. Its normal form is p s s. F with \(x=p\), \(y=s\) now gives p s. This route needs a contained full-target word; it does not assume a contained prefix word.

Finally a literal target x has an original pure power. If nonlinear, normalize it to \(x^{2}\) and use L. If literal, it is already present.

\par\smallskip\noindent\textit{Completeness.}\quad 
The argument derives every valid word absorption, with no bound on variables, word lengths or number of original summands. Distributivity converts any term to a nonempty polynomial. Apply the absorption argument to each summand in each direction to derive any valid identity. Therefore the eleven displayed identities form a complete finite basis. Full independent HSP checks are required for additional catalogue types; finite raw checks alone do not establish completeness or a transfer.

\subsection{Ten complete identities for 13251 by contained constant-word activation}\label{case:301}
Signature (+, multiplication), no constants, relative to the additively idempotent semiring axioms. Arbitrary variable sets, words and polynomial sizes are allowed throughout.

Addition:00044 /01144 /01234 /44334 /44444. Multiplication:44244 /44204 /22222 /44244 /44244. Element 2 is a multiplicative zero and the least additive value. Any word containing 2 is 2. A literal otherwise retains its value. A quadratic is 0 exactly at the ordered pair(1,\allowbreak{}3), and 4 otherwise. Every square and every word of length at least 3 is 4 when free of 2. Call the latter words constant words; this term refers to their value on non2 assignments, not to a nullary operation. These formulas follow directly by induction.

The complete relative basis is:

{\small\par\smallskip
\CaseLine{\(xyz=yxz; xyz=xzy; xxyz=xyz; xxy=xyy; xxx=xx;\)}
\CaseLine{\(x\leq xx; xy\leq xxy; xxyz\leq xxy; xyz\leq xy+yz; xxy\leq xy+y.\)}
\smallskip}

For target w, retain original summands U with content contained in c(w). Let H be their heads, and T the final variables of their nonsquare quadratic summands. Call U active when it contains an original constant word, or H intersects T. A literal is absorbed exactly when U is nonempty. A nonlinear target is absorbed exactly when U is active or the target itself is an original nonsquare quadratic.

To prove necessity, put outsiders 2. When U is inactive it contains only literals and nonsquare quadratics, and H, T are disjoint. A failed constant target is separated by \(H=1\), \(T=3\), remaining target variables 1: all right literals give 1, quadratics 0, and excluded words 2; their join is at most 0, whereas the target is 4. This includes U empty. For a missing nonsquare target hs, its content has only h, s. If U has an edge, it can only be sh. Inactivity forbids a literal h; assign \(s=1\), \(h=3\). All contained right words give 0 or 1, while hs gives 4. If U has no edge, assign \(h=0\), \(s=1\); right literals are at most 0 and the target is 4. A missing pure witness for literal h is separated by \(h=0\), outsiders 2. These exhaust failed cases.

For sufficiency, any 2-valued target is least. Otherwise all U words are free of 2. An original constant word gives 4. If b in H intersects T, an original ab edge is 4 unless \(a=1\), \(b=3\). In that exceptional case the edge is 0 and a b-headed original word is either literal 3 or nonlinear4, so the right join is 4. An original edge absorbs itself. A pure witness absorbs its literal because every positive power beyond the first is the square, and \(x\leq xx\). Thus the criterion is exact and unrestricted.

The first two word laws permit arbitrary permutations for words of length at least 3. At the first two positions use \(xyz=yxz\) with a nonempty suffix; at later adjacent positions use \(xyz=xzy\) with a nonempty prefix and external remaining suffix. For support size at least 3, bring duplicate letters to the first two positions and use \(xxyz=xyz\), whose last two substituted factors are nonempty, until every letter occurs once. With two-letter support reduce similarly to length 3; then \(xxy=xyy\) and permutation allow either variable to be the repeated one. With singleton support use \(xxx=xx\). Consequently constant words are equivalent exactly as needed here whenever their contents agree. A normal constant word is hh for singleton content, aab with \(a<b\) for two-letter content, or the increasing list for larger content. Literals and nonsquare quadratics keep their original forms.

The value formula proves all word laws. It also proves \(x\leq xx\), \(xy\leq xxy\) and \(xxyz\leq xxy\). The last two displayed inequalities follow from the intersection clause applied to xy+yz or xy+y, respectively. Hence all ten laws hold with arbitrary nonempty substitutions.

For completeness, an original target is immediate. A literal uses its pure original witness, word normalization and \(x\leq xx\). For an active U with an original constant word use it as seed g. Otherwise choose an original ab edge and original word v headed b. If \(v=b\), \(xxy\leq xy+y\) produces aab. If \(v=bB\) with B nonempty, \(xyz\leq xy+yz\) with \(x=a\), \(y=b\), \(z=B\) produces abB. Both are constant words contained in the target content.

Append all missing target variables at once: with their nonempty list R, use \(xxyz\leq xxy\) with \(x=g\), \(y=g\), \(z=R\). The right ggg normalizes to g and the left gggR is a constant word of the enlarged content. If nothing is missing no append is needed. A constant target now follows by the word normalization. A nonsquare target hs follows from \(xy\leq xxy\), because the full-content constant guard is equivalent to hhs. Every variable image is nonempty; original witnesses may have arbitrary lengths. This derives every valid absorption. Distributivity and two-way summandwise absorption give all identities, proving that the ten laws are a complete finite relative basis. Opposites and equal-variety transfers require complete independently verified HSP in both directions.

Boundary failures:\(xx\leq x\) at \(x=0\); \(xxy\leq xy\) at \(x=1\), \(y=3\); \(xz\leq xy+z\) at \(x=1\), \(y=3\), \(z=1\). An unguarded edge never activates a constant word, and containment may not be replaced by global witnesses.

\subsection{Five complete identities for semiring 13298}\label{case:302}
Signature: (+, multiplication), without constants. The base axioms are associativity of multiplication, the commutative idempotent associative addition, and both distributive laws. Write \(u\leq P\) for \(u+P=P\). This document proves unrestricted completeness, independently of the finite tests accompanying it.

Addition rows are 00004/\allowbreak{}01114/\allowbreak{}01214/\allowbreak{}01134/\allowbreak{}44444. Multiplication rows are 44244/\allowbreak{}44244/\allowbreak{}44244/\allowbreak{}40204/\allowbreak{}44244. The additive order is \(2<1<0<4\) and \(3<1<0<4\), with 2 and 3 incomparable. The claimed five-law basis is

{\small\par\smallskip
\CaseLine{\((W) xyz = zzz\)}
\CaseLine{\((Q) xy \leq  yyy\)}
\CaseLine{\((L) y \leq  xy\)}
\CaseLine{\((T) xz \leq  xy+yz\)}
\CaseLine{\((A) yyy \leq  xxx+xy.\)}
\smallskip}

Each is an identity of this operation table, checked directly on all assignments of its at most three variables. No empty-word substitutions are used below. In particular, square and cube are different: at 3, the square is 0 and the cube is 4.

\par\smallskip\noindent\textit{Word values and normal forms.}\quad 
A literal has its assigned value. A quadratic ending in 2 has value 2; a quadratic whose ordered values are (3,\allowbreak{}1) or (3,\allowbreak{}3) has value 0; every other quadratic has value 4. After two factors the value lies in \{0,\allowbreak{}2,\allowbreak{}4\}, all of whose multiplication rows agree: a next factor 2 gives 2 and any other next factor gives 4. Induction proves that every word of length at least three is 2 exactly when its last variable is 2, and is 4 otherwise. Thus (W), with x a nonempty prefix, normalizes every long word to the cube of its last variable. Literals and quadratics are not shortened.

\par\smallskip\noindent\textit{Exact criterion using ORIGINAL right summands.}\quad 
Let P be any finite nonempty sum of nonempty words. Let F be all last variables of its original summands, and K the last variables of its original summands of length at least three. Put a directed edge \(a\to b\) for every original quadratic ab, including squares as loop edges. Let R be the POSITIVE-LENGTH transitive closure of this edge relation. Put \(B=K\) union \{t: (k, t) is in R for some k in K\}. Then:

\textbullet{} A literal t is absorbed by P exactly when t is in F.

\textbullet{} A quadratic ht is absorbed exactly when t is in B or (h, t) is in R.

\textbullet{} A word of length at least three ending in t is absorbed exactly when t is in B.

No content condition is needed. A directed cycle alone can produce a square, but it cannot activate a cube. R is never made reflexive by zero-length paths.

\par\smallskip\noindent\textit{Necessity: explicit separating evaluations.}\quad 
All variables outside the finite support may be assigned 2. If t is absent from F, assign \(t=3\) and every other variable 2. Every right word ends in 2 and has value 2, whereas the literal t has value 3.

Suppose t is not in B. Let H consist of t and all vertices with a positive path to t. Including t here is only a separator construction; it does not change R. H contains no vertex of K. An edge with final vertex in H has initial vertex in H. Assign all H variables 3 and all others 2. Every original literal has value 2 or 3; every original quadratic has value 0 when its tail lies in H, and 2 otherwise; and every original long word has value 2. Consequently \(P\leq 0\). A long word ending in t has value 4. For h distinct from t with no positive path \(h\to t\), h lies outside H and ht has value \(2*3=4\), giving the quadratic separator as well.

It remains to separate tt when t is not in B and there is no positive cycle \(t\to t\). Assign \(t=1\), the proper ancestors H minus \{t\} to 3, and all other variables 2. There is no edge from t into H, for such an edge would form a positive cycle. An edge ending at t therefore starts at a proper ancestor and has value \(3*1=0\). An edge ending in a proper ancestor starts at a proper ancestor and has value \(3*3=0\). Edges ending outside H have value 2. Long right words again end outside H and have value 2; literals have values in \{1,\allowbreak{}2,\allowbreak{}3\}. Thus \(P\leq 0\), whereas tt has value \(1*1=4\). These cases exhaust failure of the criterion.

\par\smallskip\noindent\textit{Sufficiency and derivability from the five laws.}\quad 
Normalize original long words by W, so each k in K provides \(k^{3}\leq P\). Every original edge ab provides \(ab\leq P\). For a positive edge path \(a=a0\to a1\to \cdots \to am=b\), \(m\geq 1\), inductively apply T to a0*a(i-1) and a(i-1)*ai. This gives \(ab\leq P\). The \(m=1\) case is its original edge; the construction also treats \(a=b\) only when \(m\geq 1\).

For each t in B, take an original k in K followed by a directed edge path to t, allowing zero edges only when \(k=t\). Starting with \(k^{3}\), apply A along successive original edges ab: \(a^{3}\)+ab absorbs \(b^{3}\). Thus \(t^{3}\leq P\). W recovers every long target ending in t, and Q recovers every quadratic with that final variable. If a quadratic is instead provided by R, the positive-path construction already gives it.

For a literal t in F, either t is an original literal or an original nonlinear same-tail word has the form vt with v nonempty. Substitution \(x=v\), \(y=t\) in L gives \(t\leq vt\leq P\). No polynomial is changed by an uncertified semantic inference: every step uses W, Q, L, T, or A and semilattice transitivity under the unchanged original right sum.

Every term expands by distributivity to a nonempty polynomial. A valid identity \(P=Q\) is equivalent to absorption of each word of P by Q and each word of Q by P. The necessity criterion and the constructions therefore derive every valid identity. This proves finite basability for arbitrary variable counts, word lengths, and numbers of summands.

\subsection{Finite bases for 13299 and 13902 from original final-square triggers}\label{case:303}
Signature (+, multiplication), without constants, relative to the additively idempotent semiring axioms. The proof allows arbitrary variables, word lengths and numbers of original right summands.

Source 13299 addition:00044 /01144 /01244 /44434 /44444. Multiplication:44244 /44244 /44244 /40244 /44244. Its additive order \(is 2<1<0<4\), with 3 incomparable to 0,\allowbreak{}1,\allowbreak{}2. Source 13902 addition:00004 /01004 /00224 /00234 /44444. Multiplication:44434 /44034 /44434 /44434 /44434. Its additive order \(is 3<2<0<4\) and \(1<0\);1 is incomparable to 2,\allowbreak{}3.

For 13299 a nonlinear word of length 2 gives 2 when its last value is 2, 0 at the ordered pair(3,\allowbreak{}1), and 4 otherwise. Every square or word of length at least 3 gives 2 when its final value is 2 and 4 otherwise. For 13902 replace the final reset 2 by 3 and the exceptional pair(3,\allowbreak{}1) by(1,\allowbreak{}2). Thus a quadratic is 3 if last 3,\allowbreak{}0 at(1,\allowbreak{}2),4 otherwise; a square or word of length at least 3 is 3 if last 3 and 4 otherwise. Literals retain their values. Direct induction proves both formulas, since every quadratic value has a generic multiplication row.

Both sources have the following four identities:

{\small\par\smallskip
\CaseLine{\((W) xyz=zz\)}
\CaseLine{\((Q) xy\leq yy\)}
\CaseLine{\((L) y\leq xy\)}
\CaseLine{\((A) zz\leq xy+yz.\)}
\smallskip}

Source 13299 additionally has

{\small\par\smallskip
\CaseLine{\((B) yy\leq x+xy.\)}
\smallskip}

These are complete relative bases, of sizes five and four respectively.

\par\smallskip\noindent\textit{Exact criteria from original words.}\quad 
For an original polynomial P let F be all original final variables and N the final variables of its original nonlinear words. For fixed t let \(K_{t}\) mean that an original square or \(length\geq 3\) word ends in t, and let \(A_{t}\) be the heads of original nonsquare quadratics ending in t. For 13299 call t active when \(K_{t}\) holds or \(A_{t}\) meets F. For 13902 call t active when \(K_{t}\) holds or \(A_{t}\) meets N.

A literal t is absorbed exactly when t is in F. A square or \(length\geq 3\) word ending in t is absorbed exactly when t is active. A nonsquare quadratic ht is absorbed exactly when it is original or t is active. These conditions refer only to original words. There is no content containment requirement and no inferred-path requirement.

Necessity for a missing literal:assign \(t=3\) and all others 2 in 13299, or \(t=2\) and all others 3 in 13902. Every right word ends in the reset value and equals the bottom value 2 or 3, while the literal is larger. For an inactive final t and missing nonlinear target, assign \(t=1\), \(A_{t}=3\), others 2 in 13299. \(A_{t}\) misses F and t is not in \(A_{t}\). All original nonlinear words ending elsewhere equal 2; original words ending in t are incoming edges valued 0 or literal 1. The join is at most 0. The missing target has value 4. This also covers empty \(A_{t}\) or absent t in F. For 13902 instead assign \(t=2\), \(A_{t}=1\), others 3. \(A_{t}\) misses N; hence all nonlinear original words ending elsewhere equal 3. Incoming edges ending t equal 0; literals have values 1,\allowbreak{}2 or 3, whose join is at most 0. The missing target has value 4. Literal predecessor heads may occur in 13902 without activating t; their value 1 remains below 0. These separators cover every failed criterion for arbitrary variable sets.

For sufficiency a literal is below every nonlinear word with the same tail, as the word formulas show. An original square or long word ending t equals tt. Otherwise take an incoming original edge at and an original guard ending in a. If t has the reset value, tt is bottom. If the edge has value 4, tt is bounded. The only remaining case is the exceptional pair:\(a=3\), \(t=1\) in 13299 or \(a=1\), \(t=2\) in 13902, so the edge is 0. A nonlinear guard ending in a gives 4. In 13299 a literal guard a gives 3 and \(0+3=4\) as well.13902 does not permit a literal-only guard:\(1+0=0\) there. Thus active t derives the square value in every assignment. All t-ending targets are below tt. This proves the exact criteria without finite tests.

W follows from the formulas. Q and L are also immediate from the tables. A has incoming edge yz and original nonlinear guard xy ending in y, so its validity follows from the active criterion. In 13299 B uses the literal guard x and incoming edge xy. It fails in 13902 at \(x=1\), \(y=2\).

\par\smallskip\noindent\textit{Unrestricted completeness.}\quad 
W normalizes every word of \(length\geq 3\) to its final square. To do this, substitute the nonempty prefix before its last two positions for x, its penultimate letter for y and its final letter for z. Literals and quadratics remain unchanged. Thus normalization preserves F, N and the active criteria.

For a literal t choose any original t-ending word. If nonlinear, write it as Ut with U nonempty and use L with \(x=U\), \(y=t\). Otherwise it is the literal itself. For a nonlinear target, an original target needs no step. If \(K_{t}\) holds its normalized original square tt is already available. Otherwise choose an original incoming at and an original guard ending a. If the guard is Ua with U nonempty, A with \(x=U\), \(y=a\), \(z=t\) yields tt. For 13299 a literal guard a instead uses B with \(x=a\), \(y=t\). That case is not asserted for 13902. Q with the desired head and tail recovers any missing quadratic from tt; long targets normalize to tt.

Every step has nonempty variable images and uses an original witness or a word already derived below P. This derives all valid absorptions with no bound on variables, lengths or summands. Distributivity and absorption of summands both ways prove all identities. Hence the displayed five-law and four-law sets are complete. Word reversal gives opposite bases; catalogue transfers require independently complete HSP in both directions.

Strict boundaries:\(yy\leq xy\) fails at the exceptional pair in each source. In 13299 \(xy\leq y\) fails at \(x=2\), \(y=1\); in 13902 it fails at \(x=3\), \(y=2\). The extra13299 law B fails for 13902 at \(x=1\), \(y=2\), confirming that original literal guards cannot be substituted for its nonlinear guard condition.

\subsection{Eleven identities for 14623 with strictly nonlinear prefix witnesses}\label{case:304}
Signature (+, multiplication), no constants; the basis is relative to the additively idempotent semiring axioms.

\par\smallskip\noindent\textit{Word values and exact criterion.}\quad 
The tables in the catalogue order 0,\allowbreak{}1,\allowbreak{}2,\allowbreak{}3,\allowbreak{}4 are:

{\small\par\smallskip
\noindent addition        multiplication\par
\CaseLine{\ttfamily 0 0 0 0 4\qquad 4 4 4 4 4}
\CaseLine{\ttfamily 0 1 1 1 4\qquad 4 1 2 3 4}
\CaseLine{\ttfamily 0 1 2 2 4\qquad 3 3 3 3 3}
\CaseLine{\ttfamily 0 1 2 3 4\qquad 3 3 3 3 3}
\CaseLine{\ttfamily 4 4 4 4 4\qquad 4 4 4 4 4}
\smallskip}

Addition is the chain \(3<2<1<0<4\). Literals retain their assigned value. A nonlinear word of all 1 factors is 1. Otherwise scan its first non1: 0 or 4 gives 4,\allowbreak{}3 gives 3, and 2 gives 2 only at the final position, otherwise3. This is an induction using the multiplication rows. In particular no nonlinear word takes value 0.

For original P let Z(S) mean an original word has content in S. Let EN(x, S) mean an original NONLINEAR word contains x whose first-x prefix is contained in S. It is unmarked when that first x is final. A literal target x is absorbed exactly when x itself belongs to P. A nonlinear target w is absorbed exactly when:

1. Z(c(w)) holds; the contained word may be literal.

2. At every first occurrence of x, EN(x, S) holds, where S is the preceding target content. A literal cannot supply this witness.

3. If its last variable x is simple, an unmarked EN(x, S), the literal x, or Z(S) holds, with S the content before that final x.

For necessity of 1 assign target variables 1 and outsiders 3. A right value reaches 1 exactly from a contained word. For 2 assign \(S=1\), \(x=0\), and all other variables 3. The target is 4; a right value reaches 4 exactly from the required nonlinear prefix witness. Literal \(x=0\) does not reach 4. For 3 assign \(S=1\), \(x=2\), others 3. Right values reach 2 exactly through an unmarked nonlinear witness, literal x, or a contained-prefix word. A missing literal x is separated by \(x=2\), others 3.

For sufficiency,3 is least. The all 1 target is covered by a contained word. A first 0 or 4 is covered by the nonlinear prefix witness, value 4. A first 2 at the last position is covered by condition 3: an unmarked witness or literal x takes 2; a prefix-contained word takes 1. The same first 2 before the end gives the least value 3. This exhausts assignments. No restriction on variables, word lengths or polynomial size is used.

\par\smallskip\noindent\textit{Eleven identities.}\quad 
{\small\par\smallskip
\CaseLine{\(W1 xxy=xy\)}
\CaseLine{\(W2 xyx=xyy\)}
\CaseLine{\(A  xyz\leq yz+xt\)}
\CaseLine{\(B  xyz\leq xz+xyt\)}
\CaseLine{\(U  xyz\leq xz+xy\)}
\CaseLine{\(D  xy\leq xx+xyt\)}
\CaseLine{\(G1 xx\leq x+xy\)}
\CaseLine{\(G2 xyy\leq x+xyz\)}
\CaseLine{\(G3 xyy\leq y+xyz\)}
\CaseLine{\(G4 xyzz\leq y+xyzt\)}
\CaseLine{\(F  xy\leq y+xyt\)}
\smallskip}

The word values prove W1, W2. The exact criterion proves each absorption. In particular \(xy\leq xx+y\) is FALSE at \(x=1\), \(y=0\). No such rule is needed. The first-occurrence normal form is F(w) if its last variable is simple, and F(w) followed by the last letter of F(w) otherwise. Internal repeats are removed by \(xYxZ=xY^{2}Z=xYZ\) for nonempty Y, Z, or W1 for adjacency. A final separated repeat xYx is changed by W2 to \(xY^{2}\); remove internal repetitions until only the last first-occurrence letter is doubled. These explicit reductions use only nonempty substitutions and work for arbitrary words. Literal summands remain literal.

\par\smallskip\noindent\textit{Unrestricted derivation.}\quad 
Normalize words. Write N for the original nonlinear right words and \(x_{1}\)...\(x_{d}\) for the target first-occurrence list. Condition 2 guarantees N is nonempty. Maintain p \(u\leq P\) for every u in N, with p the processed list. At the first letter \(x_{1}\), EN supplies \(x_{1}\) t, t nonempty. For each \(u=YZ\) with Y, Z nonempty, A gives \(x_{1}\) \(u\leq u+x_{1}\) t.

At a later letter x, condition 2 supplies a nonlinear original witness v. If marked, p \(v=p\) x t modulo word laws, t nonempty; B extends the invariant using p x \(u\leq p\) u+p x t. If unmarked, p \(v=p\) x and U gives p x \(u\leq p\) u+p x. There is no literal-prefix branch.

Guard lemma: if p is nonempty, Z(c(p)) holds, and the invariant holds, then \(p^{2}\leq P\). A nonlinear contained original r gives p \(r=p^{2}\). Otherwise choose a contained literal a and any v in N. Write \(p=A\) a B at its unique position. Both A, B nonempty use G4 with \(x=A\), \(y=a\), \(z=B\), \(t=v\); the resulting left word normalizes to \(p^{2}\). If A is empty use G2 with \(x=a\), \(y=B\), \(z=v\); if B is empty use G3 with \(x=A\), \(y=a\), \(z=v\); if both are empty use G1 with \(x=a\), \(y=v\). Each right side is a+p v after normalization. Thus every substitution is nonempty and all four positions are covered.

For a saturated target, construct through all d letters and apply the guard lemma using Z(c(w)); \(p_{d}^{2}\) is the target normal form.

For a simple last variable s, write the target p s, p nonempty. Construct through p. If an original nonlinear unmarked witness exists, p \(v=p\) s already belongs to the maintained words. Otherwise condition 2 still gives a nonlinear marked witness v, so p \(v=p\) s t, t nonempty. If the literal s occurs, F gives p \(s\leq s+p\) s t. In the remaining case condition 3 gives Z(c(p)); the guard lemma derives \(p^{2}\) and D gives p \(s\leq p^{2}+p\) s t. These cases are exhaustive.

Every valid word absorption has been derived. Distributivity expands arbitrary terms into nonempty polynomials; absorbing their summands in both directions proves every valid identity. Thus the eleven laws are a finite basis without any bound on original variables, word lengths or summands.

The distinction from 4729 is essential: here a contained seed may be literal, but every first-occurrence prefix witness must be nonlinear. In 4729 the seed must be nonlinear and prefix witnesses may be literal. Short fingerprint coincidences do not identify these varieties.

\subsection{Five complete identities for semiring 13358}\label{case:305}
In the constant-free signature (+, multiplication), relative to the additively idempotent semiring axioms, a complete basis is

{\small\par\smallskip
\CaseLine{\(xy=yx; xyz=ttt; x\leq yyy; x\leq y+yz; xy\leq xx+yy.\)}
\smallskip}

The notation \(u\leq P\) abbreviates \(u+P=P\). Addition: 00044/\allowbreak{}01144/\allowbreak{}01244/\allowbreak{}44434/\allowbreak{}44444. Multiplication: 44344/\allowbreak{}43344/\allowbreak{}33344/\allowbreak{}44444/\allowbreak{}44444. Quadratics are 3 on pairs in\{0,\allowbreak{}1,\allowbreak{}2\} with at least one endpoint 2, or both endpoints 1; otherwise they are 4. All longer words are 4. The additive order \(is 2<1<0<4\), with 3 incomparable with 0,\allowbreak{}1,\allowbreak{}2 below 4.

Let L be original literal variables, E original unordered edges, V their endpoints, and S original loops. An original long word, or L intersecting V, makes P universally 4 and absorbs all targets. Otherwise a literal is absorbed iff original, a pair iff its edge is original or BOTH endpoints belong to S, and no long word is absorbed. An original loop at only one endpoint does not suffice.

For sufficiency below 4, all right edges and literals must equal 3. Original loop endpoints lie in\{1,\allowbreak{}2\}, and every pair in this block is 3. For necessity in the non-top case, assign V to 2 and outsiders 3. All right terms are 3 and every long target 4. A missing literal in V is 2, not below 3; for a missing outside literal set it 4. A pair with an outside endpoint has value 4. For a missing edge with an endpoint h without a loop, if the target is hh assign \(h=0\), other V vertices 2, and outsiders 3. For a distinct pair ht, assign \(h=0\), \(t=1\), other V vertices 2, and outsiders 3. The original edge ht and loop hh are absent, so every right edge is 3, while the target is 4. These cover all failures for arbitrary variable sets and polynomial sizes.

The five laws are valid by this criterion and formula, and all their assignments are independently checked. A right long word reduces to a fixed cube using \(xyz=ttt\) with a nonempty prefix; then \(x\leq yyy\) gives any target. An original incident literal y and edge yz give any target by \(x\leq y+yz\). Otherwise any new pair uses its two original loops and \(xy\leq xx+yy\). Commutativity supplies orientations. Every substitution is nonempty and every witness original. This proves all absorptions without variable, length, or summand bounds. Distributivity and mutual summandwise absorption prove all identities. Full HSP both ways is required before transferring this basis to other catalogue types.

\subsection{Thirteen complete identities for 13363}\label{case:306}
In the signature (+, multiplication) without constants, write \(u\leq P\) for \(u+P=P\). Relative to the additively idempotent semiring axioms, catalogue 13363 has the following complete basis:

{\small\par\smallskip
\CaseLine{\(xy=yx; xxy=xyy; xyzt=xxxx; xxxx=yyyy;\)}
\CaseLine{\(x\leq yyyy; xxxx\leq x+xyz; x\leq y+z+yz+xt;\)}
\CaseLine{\(xy\leq xyz; xx\leq xyz; yy\leq x+xy; xz\leq xx+zy;\)}
\CaseLine{\(xxx\leq xx+yzt; xyz\leq xxx+yyy+zzz.\)}
\smallskip}

The proof is unrestricted in variables, lengths and summand counts.

\par\smallskip\noindent\textit{Values and normal forms.}\quad 
Addition:00044/\allowbreak{}01034/\allowbreak{}00244/\allowbreak{}43434/\allowbreak{}44444. Multiplication:44344/\allowbreak{}44344/\allowbreak{}33144/\allowbreak{}44444/\allowbreak{}44444. Addition has \(1<0\),\(2<0\),\(1<3\),\(0<4\),\(3<4\), with 2 incomparable to 1 and 3. A quadratic is 1 when both occurrences are 2,\allowbreak{}3 when one occurrence is 2 and the other 0 or 1, and 4 otherwise. A cubic is 3 when all occurrences are 2, and 4 otherwise. Every word of length at least four is 4. Literals retain their values. Multiplying a possible quadratic value 1 by 2 gives 3; every other such continuation gives 4, and multiplying3 or 4 gives 4. This proves the formula at every length.

Commutativity sorts quadratics. Cubics depend only on content; \(xxy=xyy\) allows either repeated variable for a two-variable cubic. The first three positions and a nonempty remainder in \(xyzt=xxxx\) reduce any long word to a fourth power, and \(xxxx=yyyy\) identifies all fourth powers. This denotes a constant word function, not a constant symbol in the signature. Normalization may use a fixed fresh variable in its fourth power; its choice is immaterial by the fourth displayed identity.

\par\smallskip\noindent\textit{Exact original-witness criterion.}\quad 
For original P, let L be its literal variables, V its nonlinear variables, E its original quadratic pairs, and Q its original cubics. Put

{\small\par\smallskip
\noindent T = union of the contents of Q;\par
\noindent H = T union \{x:xx is an original quadratic\}\par
\noindent union \{y:xy is an original quadratic and x is in L\}.\par
\smallskip}

Let B hold when an original quadratic has both its variables in L, including a square at an original literal. Define the saturation condition

K: an original word has \(length\geq 4\), or Q is nonempty and (B holds or L intersects T).

If K holds, every target is absorbed. Otherwise the complete criterion is:

\textbullet{} A literal x is absorbed iff x is in L, or B holds and x is in V.

\textbullet{} A quadratic xy is absorbed iff x, y are in V and at least one of these holds:B holds; one endpoint is in H; xy is an original pair.

\textbullet{} A cubic is absorbed iff Q is nonempty and all its variables are in H.

\textbullet{} No word of length at least four is absorbed.

These are global witness conditions: H need not belong to the target's component. They must not be replaced by the preceding clique/path criteria.

Here is a proof for arbitrary P. A non4 sum with a nonlinear term lies below 0 or 3. If it lies below 0, every nonlinear term must be 1; hence there are no cubics, and every V variable is 2. If it lies below 3, each quadratic endpoint is in\{0,\allowbreak{}1,\allowbreak{}2\}, each cubic variable is 2, and a literal variable incident with a nonlinear term must be 1. Every H variable is therefore 2. This case is impossible if B holds or L intersects T. Consequently K forces 4. Conversely, without K a non4 assignment exists: if B holds there is no cubic, so assign all V variables 2 and the other literal variables 0; the sum is at most 0. If B fails, assign \(V\setminus L\) to 2, L to 1, and other variables 3; every quadratic is 1 or 3 and every cubic3, so the sum is at most 3. This also proves that K is exact.

The literal rule follows because, when B holds without K, every non4 evaluation has all V variables 2 and right value 0. For a quadratic, an anchored endpoint is 2 and every other V endpoint is in\{0,\allowbreak{}1,\allowbreak{}2\} in a right sum below 3. If its target value is 3, its original incidence also contributes3 somewhere on the right; otherwise the target is 1. Right sums below 0 give 1 on every such pair. The B branch excludes the former case entirely. A cubic on H has value 3 whenever the right is not 4, and the original cubic supplies 3. Original words are included.

For necessity there are explicit separating assignments. If B fails, the standard assignment \(V\setminus L\) the value \(2\), \(L=1\), \(others=3\) makes the right side at most 3; call this the 3-assignment. With no cubics, the 0-assignment sets \(V=2\), other literal variables 0, and outsiders 3; the right is at most 0. A missing literal not allowed by the rule lies outside L; use the 3-assignment, where a V variable is 2 and not below 1 or 3, or set it 4 if outside V. When B holds, a failed literal must be outside both V, L and can be assigned 4 under the 0-assignment.

For a failed nonlinear target containing a variable outside V, use the 3-assignment and give that variable 3 when B fails. The target becomes 4 and the right remains at most 3. When B holds, use the 0-assignment and give that variable 0; the target is 3 or 4 and is not below the right. A failed quadratic whose variables lie in V has B false, no anchored endpoint and no original target edge. Set its endpoints 1 in the 3-assignment. They have no original loops, cubic occurrences, or literal neighbors, and they do not co-occur in an original edge. No right term becomes 4, but the target is 4. This also covers a missing square.

For a failed cubic with variables in V and no original cubic, use the 0-assignment:its target is 3 and the right at most 0. Otherwise choose an unanchored target variable and set it 1 in the 3-assignment. It has no loop, cubic occurrence, or literal neighbor, so the right remains at most 3 while the target is 4. For a failed long target use either of the nonsaturating standard assignments; the target is always 4. These cases exhaust every failure without bounding variables.

\par\smallskip\noindent\textit{Complete derivations.}\quad 
Keep every original right summand. First handle saturation. An original long term gives a fourth power directly by normalization. If an original literal x occurs in an original cubic containing x, \(xxxx\leq x+xyz\) gives a fourth power after splitting the other two letters. If only B and an original cubic witness saturation, choose x in that cubic and use \(x\leq y+z+yz+xt\), with original literals y, z and their original quadratic yz; the x-containing original cubic supplies xt. Then the previous rule produces a fourth power. From it recover every target by \(x\leq yyyy\), substituting that whole target for x.

Without saturation, a missing allowed literal uses exactly the B rule \(x\leq y+z+yz+xt\) and an original nonlinear term containing x. To obtain xx at an H variable, use its original square, or \(xx\leq xyz\) with its original cubic, or \(yy\leq x+xy\) with its original literal neighbor and edge. In the B case every nonlinear variable can also be anchored:extract a pair from an original containing cubic using \(xy\leq xyz\) if necessary, recover a literal at the other endpoint using B, and apply \(yy\leq x+xy\). Every target quadratic now has an available anchored endpoint x. Extract an original incidence pair zy at its other endpoint z when needed, then \(xz\leq xx+zy\) gives the target. A target already original needs no derivation.

For an allowed cubic, anchor each target variable as above and retain one original cubic yzt. Use \(xxx\leq xx+yzt\) to obtain each required pure cube, then \(xyz\leq xxx+yyy+zzz\) gives the target. This remains valid when letters repeat. All factors substituted in these identities are nonempty; the B or cubic witnesses may use external variables.

Every listed law is valid by the criterion or word formula. The finite derivation covers every absorption, including arbitrary long targets. Distributive expansion and absorption in both directions therefore prove equational completeness.

Boundary examples:\(xxx\leq xx\) fails at \(x=2\); \(x\leq xxx\) fails there as well; \(xy\leq xz+yz\) fails at \(x=y=1\), \(z=2\); \(xxxx\leq xxx\) fails at \(x=2\). Global literal recovery requires the two original endpoint literals in B; an arbitrary literal/nonlinear incidence is not a substitute.

\subsection{Six complete identities for semiring 13366}\label{case:307}
In the constant-free signature (+, multiplication), relative to the additively idempotent semiring axioms, a complete finite basis is

{\small\par\smallskip
\CaseLine{\(xy=yx; xyz=ttt; x\leq yyy; xy\leq xx+yz;\)}
\CaseLine{\(xx\leq xy+y; t\leq x+y+xy.\)}
\smallskip}

Here \(u\leq P\) means \(u+P=P\). Addition: 00044/\allowbreak{}01034/\allowbreak{}00244/\allowbreak{}43434/\allowbreak{}44444. Multiplication: 44344/\allowbreak{}44344/\allowbreak{}33344/\allowbreak{}44444/\allowbreak{}44444. Quadratics are 3 exactly when both endpoints lie in\{0,\allowbreak{}1,\allowbreak{}2\} and at least one equals 2; otherwise they are 4. Long words are 4. The elements below 3 are 1 and 3.

Let L be original literals, E original unordered edges, V endpoints, S loop variables, and N the neighbors of L in E, including loop neighbors. An original long word, or an original edge with both endpoints in L, makes P universally 4. Otherwise a literal is absorbed iff original, a pair iff its endpoints lie in V and its edge is original or an endpoint lies in S union N, and no long target is absorbed. Neighbors are from original edges and literals, not a transitive closure and not a set recomputed after inference.

For sufficiency below 4, edges are 3, literals 1 or 3, and \(P=3\). Edge endpoints are 0,\allowbreak{}1,\allowbreak{}2; an incident literal must be 1, forcing each of its neighbors to 2. Every original loop also forces 2. Hence every permitted pair is 3. An edge with two literal endpoints cannot have value 3, proving the universal-top clause. With no edges, only original literals meet the criterion.

For necessity in the non-top case, assign original literal endpoints 1, all other V vertices 2, and outsiders 3. No edge has both endpoints in L, so all edges are 3 and \(P=3\). A long target is 4. Any missing literal in V is 2, not below 3; for a missing literal outside V set it 4. A pair with an outside endpoint is 4 under that assignment. For a missing edge with neither endpoint in S union N, assign its endpoints 1, all L vertices 1, all other V vertices 2, and outsiders 3. No original edge joins two assigned 1 vertices:there is no target edge or endpoint loop, no adjacency from a target endpoint to L, and no edge within L. All right terms are at most 3 while the target is 4. This covers every failure, including target squares, for arbitrary variable sets.

All six laws are valid by the criterion and formula, and their finite assignments are independently checked. Original long words supply a fixed cube via \(xyz=ttt\) using nonempty substitutions; then \(x\leq yyy\) derives any target. An original edge xy with both original literals uses \(t\leq x+y+xy\) for any target. Otherwise original targets require no step. Orient a new pair so h lies in S union N. If h lacks an original loop, choose an original literal y with original edge hy; derive hh by \(xx\leq xy+y\). If the target is that square, stop. Otherwise choose an original edge tz at its other endpoint and use \(xy\leq xx+yz\). All right terms are original or explicitly derived below P. Commutativity supplies orientations, and no witness sets are redefined.

This derives all valid absorptions with nonempty substitutions and without bounds on variables, lengths, or summand counts. Distributivity and mutual summandwise absorption prove all identities. Thus the six fixed laws form a complete basis. Other types require complete HSP in both directions before applying a transfer.

\subsection{Five complete identities for semiring 13465}\label{case:308}
Signature (+, multiplication), without constants. Relative to the additively idempotent semiring axioms a complete finite basis is

{\small\par\smallskip
\CaseLine{\((W) xyz=xy+zz;\)}
\CaseLine{\((H) x\leq xy;\)}
\CaseLine{\((T) y\leq xy;\)}
\CaseLine{\((A) zz\leq xy+yz;\)}
\CaseLine{\((R) xy\leq x+yy.\)}
\smallskip}

The proof below has no bound on variables, word length, or summands. Addition is 00004/\allowbreak{}01004/\allowbreak{}00224/\allowbreak{}00234/\allowbreak{}44444; multiplication is 44404/\allowbreak{}44004/\allowbreak{}44424/\allowbreak{}44434/\allowbreak{}44444. The additive order is \(3<2<0<4\) and \(1<0<4\), with 1 incomparable with 2 and 3. The product is 4 except when its right input is 3, in which case the result is its left input with 1 changed to 0, or at(1,\allowbreak{}2), where the result is 0. Squares are 3 at 3 and 4 at every other element. After two factors the product belongs to \{0,\allowbreak{}2,\allowbreak{}3,\allowbreak{}4\}. Multiplication of such a value by z equals its join with z squared. This proves W and by induction the exact word formula

{\small\par\smallskip
\CaseLine{\(\alpha(x_1\cdots x_n)=\alpha(x_1x_2)+\sum_{j=3}^{n}\alpha(x_j^2),\qquad n\ge2.\)}
\smallskip}

Here \(\alpha\) denotes evaluation under an arbitrary assignment. The displayed sum has no additional summands when \(n=2\); no empty term or identity element is adjoined. Literals retain their assigned value. The two extensive laws H, T follow from the multiplication table. For A, if \(z=3\) its square 3 lies below every product. Otherwise yz is 4 unless \(y=1\), \(z=2\); in that case \(xy=4\). For R, if y is not 3 its square is 4. At \(y=3\), \(xy=x+3=x+yy\). Thus all five laws hold; a separate exhaustive source-table check verifies these finite statements.

\par\smallskip\noindent\textit{Original data and exact criterion.}\quad 
Expand each original word using W into its initial quadratic and the squares of its variables in positions 3 onward, keeping literals. Let C be the original content union. Let E be the set of initial pairs of original nonlinear words, including loops. Let K be the union of variables in original positions 3 onward and the variables of loops in E. Let D be all final variables of pairs in E, and define

{\small\par\smallskip
\noindent \(F=K\) union \{b: (a, b) in E and a in K union D\}.\par
\smallskip}

Then a literal h is absorbed by P iff h is in C. A nonlinear target \(w=h\) t S is absorbed iff its content is in C, its initial pair ht is in E or t is in F, and every variable in S belongs to F. Original position 2 and position 3 are different: an arbitrary incoming edge does not create its final square. For example \(yy\leq xy\) fails at \(x=1\), \(y=2\). The two-edge pattern xy+yz produces the final square zz.

\par\smallskip\noindent\textit{Necessity through explicit separators.}\quad 
If the target contains an outside variable, assign all C variables 3 and every outside variable 1. Every right word is 3. A target literal outside C is 1, incomparable with 3; a nonlinear target with an outside variable is 0 or 4, also not below 3.

Otherwise a failure is either an initial pair ht absent from E with t outside F, or a variable t outside F in a target \(position\geq 3\). In either case construct an assignment with every original right summand \ensuremath{\leq}0 and the target 4. If t is not in D, then t is also outside K; assign \(t=0\) and every other variable 3. No original pair ends at t, and no original suffix square contains t. Pairs beginning at t give 0, other pairs give 3, and literals give 0 or 3. The failing target is 4.

If t is in D but not F, let I be the original predecessor heads of t. The definition of F implies I is disjoint from K union D, and t is not in K or I. Assign \(t=2\), every variable of I the value 1, and all others 3. A pair ending at t starts in I and gives 0. No pair ends in I. Every other pair therefore ends in 3 and has value 0,\allowbreak{}2, or 3; every suffix square is 3 because its variable belongs to K. All original literals are \(\leq 0\). The target with t in a \(position\geq 3\) is 4. For a failing initial pair ht, the absence of ht from E means h is not in I, so h is 2 or 3 and ht is 4. Its other square summands cannot reduce4. These cases establish necessity for arbitrary polynomials.

\par\smallskip\noindent\textit{Unrestricted finite derivations.}\quad 
Apply W successively to each original long word, replacing its final letter by its square summand and shortening the remaining prefix. The resulting polynomial contains all E pairs, all K squares, and all original literals. H, T recover every variable of C from a pair containing it; a suffix variable is recovered from its square.

For t in F minus K choose an original edge at. If a is in K, use A with \(x=a\), \(y=a\), \(z=t\) and the available aa and at. If a is in D, choose an original pair ba and apply A to ba+at. This gives tt in both cases. Therefore every F square is derived under the original sum.

An allowed target initial pair is either already in E, or is obtained from its available literal h and square tt by R. Every target suffix square is available by the previous paragraph. Adding these derived summands and reversing the finite W expansion recovers the entire target. All substitutions are nonempty. No derived edge is silently treated as an original witness, and the definition of F stays fixed.

Distributivity expands arbitrary terms to polynomials; apply the above derivation to each summand in both directions of a valid identity. Thus the five identities form a complete finite basis. The opposite algebra has the reversed basis. Transfers to other catalogue types require independent full HSP membership in both directions.

\subsection{Five complete identities for 13837}\label{case:309}
In the signature (+, multiplication) without constants, a complete basis relative to the additively idempotent semiring axioms is

{\small\par\smallskip
\CaseLine{\(xy=yx; xyz=xyzz; xyzt\leq xyz; xy\leq xxy; x\leq xx.\)}
\smallskip}

Here \(u\leq P\) denotes \(u+P=P\). The addition table is 00004/\allowbreak{}01014/\allowbreak{}00224/\allowbreak{}01234/\allowbreak{}44444: 3 is bottom, 1 and 2 are incomparable with join 0, and \(0<4\). Multiplication is 44434/\allowbreak{}40434/\allowbreak{}44034/\allowbreak{}33333/\allowbreak{}44434.

Induction gives these exact word values. Literals retain their values. A nonlinear word containing an occurrence assigned 3 has value 3. Without 3, a quadratic has value 0 when both occurrences are 1 or both are 2; every other nonlinear word has value 4. In particular, a long word is 3 or 4 according only to its content.

For original nonempty polynomial P, the exact criterion for \(w\leq P\) is:

\textbullet{} A literal requires an original pure word on that same variable.

\textbullet{} A quadratic requires either that original quadratic (up to commutativity), or an original contained long word.

\textbullet{} A long word requires an original contained long word.

Here long means length at least three, and contained means that its content is a subset of c(w). A literal or quadratic cannot substitute for the original long witness in either nonlinear branch.

Necessity has explicit separating assignments at every size. For a missing literal with no original pure word, assign the target variable 4 and all outsiders 3. Each right word contains an outsider and is 3, whereas the target is 4. For a long target lacking a contained long word, assign all target variables 1 and outsiders 3. Every noncontained right word is 3, and contained literals and quadratics are at most 0; the target is 4. For a missing quadratic xy with x distinct from y, assign \(x=1\), \(y=2\) and outsiders 3. Without the original edge or a contained long word, the right terms are literals or endpoint loops, with join at most 0, but the target is 4. For a missing square xx, assign \(x=0\) and outsiders 3; only literal x can survive on the right, with value 0, while the target is 4.

For sufficiency, an original contained long word has value 4 whenever the nonlinear target has no 3, and the target is bottom 3 otherwise. A literal is below each of its positive pure powers, directly from the five scalar possibilities. Original quadratic targets are immediate. Thus the criterion is exact without any bounds on variables or size.

For completeness, first normalize long words using the two word laws. Commutativity moves any chosen variable to the end; \(xyz=xyzz\) then duplicates it using two nonempty preceding blocks. The reverse law deletes surplus copies while keeping length at least three. This produces the sorted content for at least three variables, xxy for two, and xxx for one. Therefore two long words with equal content are provably equal. Quadratics retain their multiplicities, and literals retain length one.

If the target is nonlinear and not already an original quadratic, choose an original contained long u. After normalization, use \(xyzt\leq xyz\) to append the nonempty word of missing target variables, if any: split u into its first nonempty block and last two letters. The result has exactly the target content. Normalize; if the target is quadratic, recover it with \(xy\leq xxy\). A literal target has an original pure power. If long, normalize to xxx and use \(xy\leq xxy\) at \(y=x\) to derive xx. Then \(x\leq xx\) derives the literal. If the witness is already xx only the last step is required.

All substitutions are nonempty; no empty word or constant is used. These derivations handle arbitrary finite polynomials and variable sets. Distributive expansion and absorption in both directions give every identity. The five laws are valid by the exact criterion, so the basis is complete. Scalar and replay checks supplement this proof.

The short-witness boundary cannot be weakened: \(xxx\leq xx\) fails at \(x=1\), \(xy\leq xx+yy\) fails at \(x=1\), \(y=2\), and \(xx\leq x\) fails at \(x=0\). An outside long word is insufficient: \(xxx\leq yyy\) fails at \(x=1\), \(y=3\).

\subsection{Seven complete identities for 13840}\label{case:310}
Signature (+, multiplication), without constants. Relative to the additively idempotent semiring axioms, a complete basis is

{\small\par\smallskip
\CaseLine{\(xy=yx; xyz=xyzz; xxy\leq xxx+yyy; xxy\leq xxx+xy;\)}
\CaseLine{\(xxyz\leq xxy+xz; xy\leq xxy; x\leq xx.\)}
\smallskip}

Here \(u\leq P\) abbreviates \(u+P=P\). Addition is 00004/\allowbreak{}01014/\allowbreak{}00224/\allowbreak{}01234/\allowbreak{}44444: bottom 3, incomparable 1 and 2 with join 0, then top 4 above 0. Multiplication is 44434/\allowbreak{}40434/\allowbreak{}44234/\allowbreak{}33333/\allowbreak{}44434. A literal retains its value. A nonlinear word is 3 if a variable is 3, otherwise2 if all are 2, otherwise0 for a quadratic whose occurrences both equal 1, and 4 in every remaining case. This follows by induction at all lengths.

\par\smallskip\noindent\textit{Exact original-witness criterion.}\quad 
For original P and target w, put \(U=c(w)\) and keep only original words \(P_{U}\) with content contained in U. Start A with the union of contents of original long words in \(P_{U}\). Repeatedly add both endpoints of an original quadratic in \(P_{U}\) whenever either endpoint is already in A. Stop at closure. Equivalently, A consists of the components of the original nonlinear-word content graph that contain an original long word. A literal does not provide an edge or activate a component.

The exact absorption criterion is:

\textbullet{} A literal requires an original pure word on its variable.

\textbullet{} An original quadratic, up to commutativity, is absorbed immediately.

\textbullet{} Every other nonlinear target requires U contained in A.

All graph edges and activating long words must be original and CONTAINED in U. This differs from availability of a single long word.

For sufficiency, an element is below every pure power, handling literals. For a nonlinear target, if some target variable is 3 the target is bottom. If none is 3, assume the value of \(P_{U}\) is not 4. Every original long word must then have value 2, forcing all of its variables to 2. An original quadratic incident at a variable assigned 2 avoids value 4 only if its other endpoint is also 2. Hence every A variable is 2. If U is contained in A the target is 2, and the original long words ensure that the right side is at least 2.

For necessity, a missing literal without a pure original word is separated by assigning it 4 and outsiders 3. For a long target with U not contained in A, assign A variables 2, other U variables 1, and outsiders 3. Noncontained right terms are 3. Contained long words are 2, contained quadratics are 2 in active components and 0 in inactive ones, and literals are 1 or 2. Their join is at most 0; the target is 4.

For a missing quadratic, no contained mixed quadratic can be present: it would be the target itself. If a contained mixed long word were present, both endpoints would be active. Thus when the criterion fails, all contained nonlinear witnesses are pure powers. For a square, absence of the original square and activation leaves only its literal; assign the variable 0 and outsiders 3. For two distinct endpoints, if one is active assign it 2 and the inactive one 1. If neither is active assign the two endpoints 1 and 2 respectively. In either case the target is 4 while all right terms are at most 0. These separators work regardless of word lengths and variable count.

\par\smallskip\noindent\textit{Unrestricted finite derivations.}\quad 
Commutativity and \(xyz=xyzz\) normalize long words by content: move a chosen variable to the end, duplicate it with two nonempty preceding blocks, or delete surplus copies with the reverse law while retaining length at least three. Canonical long forms are sorted content for three or more variables, xxy for two, and xxx for one. In particular, \(g^{3}=g\) for long g. Quadratics retain their multiplicities.

For a new nonlinear target, retain every original contained long word. Combine them using \(xxy\leq xxx+yyy\), substituting long words for x and y. The right terms normalize to those words, while the left normalizes to a long word with their content union. This gives a long g with exactly the initial activating content.

For each closure step retain the original contained edge xz, with x in c(g) and z outside it. If g has only variable x, normalize g to xxx and use \(xxy\leq xxx+xy\) at \(y=z\). Otherwise let p be a nonempty word listing c(g) minus x. The long words g and xxp have equal content. Use \(xxyz\leq xxy+xz\) with \(y=p\), obtaining xxpz, a long word with the enlarged content. Repeating along the finite closure computation reaches U. Normalize to the target long form, or use \(xy\leq xxy\) to recover a quadratic target. The path length and number of original components have no bound.

A literal has a pure original power. If long, normalize it to xxx and use \(xy\leq xxy\) at \(y=x\) to derive xx; then use \(x\leq xx\). An original quadratic target requires no recovery. Every substitution is nonempty. Thus every valid absorption follows for arbitrary polynomials and variables. Distributivity and absorptions in both directions yield all identities. Each of the seven laws satisfies the criterion and has been checked on all scalar assignments; computations supplement the unrestricted proof.

Necessary boundaries: \(xyz\leq xxx\) fails at \(x=2\), \(y=z=1\); \(xy\leq xx+yy\) fails at \(x=1\), \(y=2\); \(xxx\leq xx\) fails at \(x=1\). An outside activating word cannot replace a contained one: \(xxx\leq xx+yyy\) fails at \(x=1\), \(y=3\).

\subsection{Eight identities for 14541 from strict nonlinear prefix witnesses}\label{case:311}
Relative to the ai-semiring axioms, in the constant-free signature (+, multiplication). The proof has no finite parameter bound.

\par\smallskip\noindent\textit{1. Word values and exact absorption.}\quad 
Addition rows are 00004,\allowbreak{}01114,\allowbreak{}01224,\allowbreak{}01234,\allowbreak{}44444, \(so 3<2<1<0<4\). Multiplication rows are 44444,\allowbreak{}41134,\allowbreak{}41134,\allowbreak{}33333,\allowbreak{}44444. Literals retain their values. A nonlinear word scans its first factor outside\{1,\allowbreak{}2\}:3 gives 3,\allowbreak{}0 or 4 gives 4. If all factors are 1 or 2, it gives 1. Multiplication induction proves this formula for all word lengths.

Let N be the original NONLINEAR summands of the original polynomial P. For nonlinear target w, the exact criterion is:

(N) some original u in N has c(u) contained in c(w); (E) each first target occurrence x has an original q in N with its first-x prefix content contained in the preceding target content S.

A literal target x is absorbed iff P has an original pure-x word. If N fails, assign \(c(w)=2\), outsiders 3:the target is 1, but all right values are at most 2. If E fails, assign \(S=2\), \(x=0\), others 3:the target is 4, but no nonlinear right word reaches 4, and any literal x is only 0. A missing pure-x witness for literal x is separated by \(x=1\), all outsiders 3. For sufficiency, target 3 is bottom, target 1 is witnessed by N, and target 4 is witnessed by E at the first 0/\allowbreak{}4-valued target variable. Literal targets follow from \(x\leq x^{2}\). This exhausts every assignment for arbitrary P, w.

In particular, literals cannot replace nonlinear prefix witnesses:the join is a chain, so several lower values cannot create4. This distinguishes 14541 from both 5001 and 5002; their literal-guard clauses do not transfer.

\par\smallskip\noindent\textit{2. Eight laws and explicit normalization.}\quad 
{\small\par\smallskip
\CaseLine{\(W1 xxy=xy\qquad W2 xyx=xy\)}
\CaseLine{\(W3 xyy=xy\qquad W4 xxx=xx\)}
\CaseLine{\(L  x\leq xx\qquad H  xy\leq yy+xt\)}
\CaseLine{\(B  xyz\leq xz+yy\qquad I  xyz\leq xz+xyt\)}
\smallskip}

The first four follow from the word-value rule, L from the five elements, and H, B, I from N, E applied to their displayed right words. Every actual substitution is nonempty. W2 deletes separated repeated occurrences; W1 deletes adjacent repeats with a suffix; W3 handles adjacent final repeats with a prefix; W4 reduces singleton powers to the square. Hence a nonlinear normal word is its first-occurrence list when it contains at least two variables, and otherwise its variable square. A literal remains literal. Every nonlinear u satisfies \(u^{2}=u\) by these word laws.

\par\smallskip\noindent\textit{3. Construct the target with its original nonlinear witnesses.}\quad 
Normalize original words and the target. A nonlinear singleton target is already an original square by N. Otherwise choose an original contained nonlinear u. At the first target variable h, E supplies a nonlinear original \(q=hV\) with V nonempty. H gives

{\small\par\smallskip
\CaseLine{\(h u \leq  u^{2}+hV = u+q \leq  P.\)}
\smallskip}

Normalize to g with correct head h and content contained in c(w). Maintain a nonempty correct first-occurrence prefix p. For the next letter x, if x is already next, extend p. Otherwise choose nonempty R such that pR normalizes to g. When p lists all current distinct letters, restore one final repeat by W3 backwards; the singleton case already has a square. This is a word-law equality, never an empty substitution.

Choose the original NONLINEAR witness \(q=r\) x V with c(r) contained in c(p). B gives p q \(R\leq pR+q^{2}=g+q\leq P\). The bridge normalizes to p x V R. Put \(T=V\) R, which is nonempty. I with \(x=p\), \(y=x\), \(z=R\), \(t=T\) then gives

{\small\par\smallskip
\CaseLine{\(p x R \leq  pR+p x T \leq  P.\)}
\smallskip}

Normalize:add x when absent or move its first occurrence when present later. The output uses no outside variable even if the bridge did. Each step extends the correct prefix; after finitely many steps the normal word is exactly w. No recursive literal-guard case is needed. All parameters are arbitrary and all substituted blocks are nonempty.

\par\smallskip\noindent\textit{4. Completeness.}\quad 
Literal targets are original literals or follow from L and an original pure-variable nonlinear word. Section 3 derives every nonlinear valid absorption. Distributivity expands terms into nonempty polynomials; absorb each summand in both directions to derive all valid identities. The eight laws are therefore a complete finite basis for 14541. Reversal supplies the opposite basis. Catalogue transfers require independent complete HSP inclusions in both directions. Finite tests and literal proof replays are supplementary, not completeness proofs.

\subsection{Nine complete identities for semiring 13900}\label{case:312}
Signature (+, multiplication), no constants. Relative to the additively idempotent semiring axioms, the following nine laws are a complete basis, without a bound on variables or term sizes:

{\small\par\smallskip
\CaseLine{\(xyz=yxz; xyz=xzy; xxyz=xyz; xxy=xyy; xxx=xx;\)}
\CaseLine{\(x\leq xx; xy\leq xxy; xxyz\leq xxy; xyz\leq xy+yz.\)}
\smallskip}

Addition rows:00004/\allowbreak{}01014/\allowbreak{}00224/\allowbreak{}01234/\allowbreak{}44444. Multiplication rows:44434/\allowbreak{}44034/\allowbreak{}44434/\allowbreak{}33333/\allowbreak{}44434. The additive order has bottom 3, then incomparable 1,\allowbreak{}2, then 0, then 4. Element 3 is a multiplicative zero. A quadratic free of 3 has value 0 exactly at the ordered pair(1,\allowbreak{}2), and 4 otherwise. A square or a word of \(length\geq 3\) has value 4 when free of 3 and value 3 otherwise. These formulas follow immediately for quadratics, and inductively since a non3 product after two factors is 0 or 4, whose rows outside 3 are 4. Call a square or word of \(length\geq 3\) a constant word, referring only to its behavior on non3 assignments, not to a constant in the signature.

For a target w, let U contain those original right summands whose contents lie in c(w). Let H, T be respectively the initial and final variables of the original nonsquare quadratic words in U. Then:

\textbullet{} A literal is absorbed iff U is nonempty.

\textbullet{} A nonlinear word is absorbed iff U contains a constant word, or H intersects T, or the target itself is an original quadratic in U.

Unlike the distinct source 13251, literal heads do not belong to H. Indeed \(xxy\leq xy+y\) fails here at \(x=1\), \(y=2\):the target is 4 and the right side is 0. Containment and this original-edge restriction are essential.

To prove necessity, assign outsiders 3. If a literal h has no pure original witness, assign \(h=1\); every right word is 3. For a nonlinear target without a contained constant word or head-tail collision, H, T are disjoint. Assign H variables 1, T variables 2, and remaining target variables 1. Contained original quadratics are 0, contained literals are 1 or 2, and all other right words are 3. The right sum \(is\leq 0\). Every constant target is 4. For a missing nonsquare target ht, its content has just h, t. A contained edge, if any, can only be th; assign \(t=1\), \(h=2\). All right terms are \(\leq 0\), whereas ht is 4. If there is no contained edge, assign all target variables 1:the right sum \(is\leq 1\) and ht is 4. This proves necessity for arbitrary polynomials.

For sufficiency, a target containing 3 is least. Otherwise contained constant words give 4. If b is a tail of ab and a head of bc, the first edge is 4 unless \(a=1\), \(b=2\); the second is 4 unless \(b=1\), \(c=2\). Both cannot avoid 4, so their sum is 4. An original target is immediate. A pure witness bounds its literal because \(x\leq xx\) and all higher powers equal the square. Thus the criterion is exact.

The five word laws give a common normal form for constant words of equal content:hh for singleton content; aab for two letters \(a<b\); the increasing list otherwise. The first two laws give adjacent swaps in any word of \(length\geq 3\). Bring duplicates to the first two positions and use \(xxyz=xyz\) while length exceeds max(3, content size); both remaining factors are nonempty. At length 3 with two letters, \(xxy=xyy\) transfers the duplicate. For singleton content, \(xxx=xx\) finishes. Nonsquare quadratics and literals are unchanged. All nine laws follow from the word-value formula and the just-proved criterion (the final law uses the original edges xy, yz). Their exhaustive source check is separate.

For complete derivability, normalize the original words by these laws. A literal follows from its pure witness and \(x\leq xx\). For a nonlinear absorption, if the target is original no inference is needed. Otherwise take an original contained constant word as seed g. If none exists, choose original edges ab, bc with common middle variable b. Apply \(xyz\leq xy+yz\) to derive abc and use this as the constant seed. Every letter of g lies in the target content.

Let R be the list of missing target variables. If nonempty, apply \(xxyz\leq xxy\) with \(x=g\), \(y=g\), \(z=R\). Its right ggg normalizes to g; its left gggR is a constant word of the full target content. A constant target now follows by normalization. A nonsquare quadratic ht follows by \(xy\leq xxy\) because the full-content constant guard normalizes to hht. Every image is nonempty; no new original witness is assumed during the derivation. This proves all valid word absorptions. Expand arbitrary terms by distributivity and apply absorption summand by summand in both directions to prove completeness of the finite basis.

The algebra is isomorphic to its multiplication opposite under the permutation swapping1 and 2 and fixing0,\allowbreak{}3,\allowbreak{}4. This is checked separately on both operations. Catalogue membership and equality of varieties still require full HSP verification, not a short polynomial match.

\subsection{A 53-identity finite basis for algebra 14361}\label{case:313}
This proof is unrestricted in variables, word length, polynomial size and graph size. Independent verification and catalogue transfers are recorded separately. The signature has no constants.

\par\smallskip\noindent\textit{1. Word values and exact criterion.}\quad 
Addition rows: 00044 / 01034 / 00244 / 43434 / 44444. Multiplication rows: 44444 / 40344 / 44244 / 44344 / 44444. The proper principal additive ideals are \{1\}, \{2\}, \{0,\allowbreak{}1,\allowbreak{}2\}, \{1,\allowbreak{}3\}. The non-4 binary products are \(1*1=0\), \(1*2=3\), \(2*2=2\), \(3*2=3\). For length at least three, every tail position must be 2 and the head must be in \{1,\allowbreak{}2,\allowbreak{}3\}; the value is 2 for head 2 and 3 for head 1 or 3. All other assignments give 4. This follows by induction: 0 dies at the next multiplication, while 2 and 3 survive exactly on a right factor 2. Literals keep their value. Quadratics additionally allow \(1*1=0\).

For original polynomial P, let C be its content, H its set of heads (including literals), and T the set of all variables in tail positions. Form the undirected graph of original quadratic words, retaining loops. Its incident vertices form V. Activate every component meeting an original long word; isolated long-word variables also activate. Write N for this activated set. Pure literal isolated vertices do not activate. Let a crossing mean H intersects T. Then \(w\leq P\) is equivalent to:

\textbullet{} For a literal x: x is in C and either there is a crossing or x is in H.

\textbullet{} For a quadratic xy: either x, y lie in one actual quadratic component (x must be in V), or both are in N; additionally there must be a crossing or x in H and y in T.

\textbullet{} For a long word: its entire content is in N; additionally there must be a crossing or its head in H and its tail set contained in T.

To prove this criterion, consider assignments making P lie in each proper principal ideal. Under \(P\leq 0\) every quadratic component is either all 1 or all 2, and every long-word variable is 2. Thus precisely the activated variables are forced to 2. Literal-only variables can be 0,\allowbreak{}1,\allowbreak{}2. Under \(P\leq 3\) all heads are 1/\allowbreak{}3 and all tails are 2; such assignments exist exactly when there is no crossing. Under \(P\leq 2\) every content variable is 2. \(P\leq 1\) is possible only if P consists of literals, all evaluated at 1. The three stated conditions make the target lie in the corresponding ideal in each case, proving sufficiency.

For necessity, absent content variables are assigned 4 and all original variables 2. A long target using an unactivated component is separated by assigning that component 1 and all others 2; its value is 3 or 4 while \(P\leq 0\). An isolated literal-only variable can instead be assigned 0. For an invalid quadratic, an isolated literal-only endpoint is assigned 0, or one of the distinct unactivated endpoint components is assigned 1 and the other endpoint 2; the target is 3 or 4 and \(P\leq 0\). If only the head/tail condition fails, assign H to 3, T to 2, and outside variables to 4. Since there is no crossing, \(P=3\) but the target is 2 or 4. This also separates a literal lying in T but not H. These arguments include repeated endpoints and arbitrary other variables.

\par\smallskip\noindent\textit{2. A finite list of identities.}\quad 
Besides the additively idempotent semiring axioms, use the following thirteen fixed identities, where \(p\leq q\) abbreviates \(p+q=q\):

1. \(xyz=xzy\).

2. \(xyzz=xyz\).

3. \(x\leq xy\).

4. \(xz\leq xy+yz\).

5. \(xt\leq xy+zy+zt\).

6. \(xyy\leq xyz\).

7. \(xyy\leq xy+xzz\).

8. \(xyy\leq xy+zyy\).

9. \(xtt\leq xyy+ztt\).

10. \(xyz\leq xyy+xzz\).

11. \(xxx\leq xy+yyy\).

12. \(yyy\leq xy+xxx\).

13. \(xyy\leq xxx+yyy\).

Add these forty templates. Symbols x, y, z, q, A, B, U, V, W are formal variables. Each bracketed block is independently present or omitted when WRITING the identity; it is never substituted by an empty word. V is mandatory.

\textbullet{} Sixteen: x \ensuremath{\leq} [A]x[B] + q[U] + Vq[W].

\textbullet{} Twelve: xxx \ensuremath{\leq} u + q[U] + Vq[W], for \(u=xAB\), AxB or ABx.

\textbullet{} Eight: xx \ensuremath{\leq} u + q[U] + Vq[W], for \(u=xy\) or yx.

\textbullet{} Four: xz \ensuremath{\leq} yx + yz + q[U] + Vq[W].

These 53 identities have at most seven variables and words of length at most four (the second fixed word law has length four). Every template contains an explicit crossing q. The criterion verifies each template: the first has content x; the second activates x via u; the third places x in an actual quadratic component; the fourth connects x to z. The fixed laws satisfy the same criterion. An independent implementation also checks every source and transferred-target assignment for all 53.

Laws 1/\allowbreak{}2 normalize every long word to its head followed by its tail set in a chosen order, duplicating the unique tail symbol when necessary to keep length at least three. Adjacent swaps act only on tail positions. Duplicate deletion by law 2 keeps a nonempty initial block and stops at length three. Quadratic words remain quadratic; literals remain literals. This is a convenient syntactic normal form, without claiming that all distinct such forms have distinct word functions.

\par\smallskip\noindent\textit{3. Completeness when there is a crossing.}\quad 
Fix original crossing witnesses qU and VqW, where U, W may be absent and V is nonempty. Keep these and every other original word throughout. Each step adjoins a word below the current polynomial.

If the target is a literal x, choose an original occurrence in AxB. The first sixteen templates, with absent prefix/suffix omitted in the chosen template, derive x. All substitutions actually used are nonempty.

For every variable x in every original long word, split that word at the chosen occurrence into one of xAB, AxB, ABx with A, B nonempty. At an endpoint the other side has length at least two and splits into two nonempty blocks; at an internal position use its two nonempty sides. The twelve templates give x cubed. Propagate cubes across each actual quadratic edge by laws 11/\allowbreak{}12. This gives x cubed for every x in N, without imposing a bound on propagation distance. Law 13 then gives htt for every pair h, t in N. Law 10 merges these elementary tails: substitute its z by the nonempty word of remaining tail variables, and normalize hzz by laws 1/\allowbreak{}2. Induction on tail-set size yields every valid long target. An active quadratic ht follows from htt by law 3 with \(x=ht\), \(y=t\).

For a quadratic target xy in an actual quadratic component, take an original edge incident with y. One of the eight square templates gives yy. Along an undirected path back to x, a forward edge uses law 4, and a reversed edge uses one of the final four crossing templates. This derives xy. At \(x=y\) the starting edge remains required; no zero-length path supplies an original witness.

\par\smallskip\noindent\textit{4. Completeness without a crossing.}\quad 
Now H and T are disjoint. Every quadratic edge goes from H to T. For a valid literal, take an original word with that head and extract the head by law 3, or use its literal summand directly.

For a quadratic target ht in one component, an undirected path from h in H to t in T has odd length and alternates its orientations. Its first edge is original. Repeated applications of the rectangle law 5 replace the next pair of edges, giving ht. This works for arbitrarily long paths, and does not rely on reflexive reachability.

From each original long word with head a, law 6 and tail permutations give att for every tail variable t. These are anchors for head a and tail t. Propagate activation over quadratic edges ht: if h already has an anchor hzz, law 7 gives htt; if t already has an anchor ztt, law 8 gives htt. Thus every activated H-vertex has some head anchor, and every activated T-vertex has some tail anchor. This process reaches exactly N because N is defined by original long words and finite quadratic paths.

Given active h in H and t in T, combine a head anchor hyy and a tail anchor ztt by law 9 to get htt. Merge the desired tail set by law 10 as before. Every valid long target follows. Every quadratic with both endpoints active follows by law 3 from its long anchor.

Every semiring term expands to a finite nonempty sum of words. A valid polynomial equality is equivalent to absorbing each summand of either side into the other side. Sections 3/\allowbreak{}4 derive every absorption selected by the exact criterion. Consequently all identities of 14361 follow from these 53 identities and the semiring axioms, with no variable, length, graph-size or summand-count bound.

\subsection{Six complete identities for semiring 14366}\label{case:314}
Constant-free signature (+, multiplication). Relative to the additively idempotent semiring axioms, take these six identities:

{\small\par\smallskip
\CaseLine{\(xy=yx; xyz=ttt; x\leq yyy; x\leq xx;\)}
\CaseLine{\(xx\leq x+xy; xz\leq x+xy+yz.\)}
\smallskip}

Here \(u\leq P\) means \(u+P=P\). Addition is 00044/\allowbreak{}01034/\allowbreak{}00244/\allowbreak{}43434/\allowbreak{}44444; multiplication is 44444/\allowbreak{}40404/\allowbreak{}44044/\allowbreak{}40444/\allowbreak{}44444. A quadratic is 0 precisely at (1,\allowbreak{}1),(1,\allowbreak{}3),(3,\allowbreak{}1),(2,\allowbreak{}2), and 4 otherwise. In particular \(3\cdot3=4\). All long words have value 4. The additive order has \(1<0<4\),\(2<0<4\), \(1<3<4\), and \(0+3=4\).

\par\smallskip\noindent\textit{Exact original-polynomial criterion.}\quad 
An original long word in P absorbs every target. Otherwise let L be original literal variables and E original unordered quadratic edges. Let V be exactly the edge endpoints, S the loop variables, and use undirected connectivity in(V, E). A literal h is absorbed iff h is in L union S. A pair ht is absorbed iff its endpoints lie in one connected component and either the edge ht is original, or one endpoint belongs to L union S. No long target is absorbed. Connectivity includes reflexivity only at a vertex of V. A literal outside V cannot create a square.

For sufficiency, consider P below 4 with an edge present. Every edge is 0, so \(P=0\). Every component lies within one of\{1,\allowbreak{}3\} and\{2\}. An original loop endpoint cannot be 3; an original literal endpoint also cannot be 3, \(as 3+0=4\). Hence a guarded endpoint lies in\{1,\allowbreak{}2\}. Any other endpoint in its component gives a pair of value 0. Original edges already have value 0. A new literal at an original loop is 1 or 2, both below 0. With no edges, only original literals are permitted.

For necessity, all variables 1 separate any long target. A missing literal h outside L union S is separated by \(h=3\) and all others 1: no loop at h occurs, so all right edges are 0 and all right literals 1. For a pair with an endpoint outside V, assign V to 1 and outsiders 0. All right terms are 1 or 0, whereas the target has value 4. For a pair in different components, assign the component of one endpoint 2 and all other variables 1. All right edges are 0 and literals 1 or 2, while the target is 4. Finally, for endpoints in one component but with no original edge and no endpoint in L union S, assign both endpoints 3 and all other variables 1. No original edge has two3 endpoints, including loops, and no literal has value 3. Right terms are at most 0 and the target 4. These cases exhaust criterion failure, without bounds on variables or polynomial size.

\par\smallskip\noindent\textit{Unrestricted complete derivations.}\quad 
Each of the six laws satisfies the table formula and criterion; independent checking tests every assignment to each fixed law. Original long words reduce to a fixed cube using \(xyz=ttt\) and a nonempty prefix. The law \(x\leq yyy\) then derives any target. Otherwise original targets require no step. A new literal at an original loop follows from \(x\leq xx\). Orient a new target pair so its guarded endpoint is h. If h is an original loop but not an original literal, derive its literal first by \(x\leq xx\). For a target square, use \(xx\leq x+xy\) with an original incident edge. For distinct endpoints, choose a simple original path from h to t. Repeatedly apply \(xz\leq x+xy+yz\), retaining the already derived literal h and using the previous path edge and the next original edge. Commutativity orients those edges and restores target orientation.

Every step lies below P; new terms do not alter original witness sets. All substituted words are nonempty and no bounds on variables, lengths, or summand counts occur. Thus every valid absorption follows from six fixed laws. Distributivity and summandwise absorption in both directions prove all valid identities. This is a complete finite basis. Other catalogue targets require independently checked full HSP in both directions before the conclusion is transferred.

\subsection{Fourteen complete identities for semiring 14370}\label{case:315}
Signature (+, multiplication), no constants. The following fourteen identities, relative to the additively idempotent semiring axioms, are a complete finite basis:

{\small\par\smallskip
\CaseLine{\(xyz=yxz; xyz=xzy; xxyz=xyz; xxy=xyy;\)}
\CaseLine{\(x\leq xx; xx\leq xy; yy\leq y+xy;\)}
\CaseLine{\(xz\leq xy+yz; xy\leq xz+yz; xxx\leq xyz;\)}
\CaseLine{\(xxx\leq xy+yyy; yyy\leq xxx+xy;\)}
\CaseLine{\(xxy\leq xxx+yyy; xy\leq xxy.\)}
\smallskip}

Addition:00044/\allowbreak{}01034/\allowbreak{}00244/\allowbreak{}43434/\allowbreak{}44444. Multiplication:44444/\allowbreak{}40404/\allowbreak{}44244/\allowbreak{}44444/\allowbreak{}44444. The additive order has \(1<0<4\),\(2<0<4\), and \(1<3<4\);0 and 3 as well as 2 and 3 have join 4. A quadratic has value 2 at(2,\allowbreak{}2),0 at(1,\allowbreak{}1) or(1,\allowbreak{}3), and 4 otherwise. Every word of \(length\geq 3\) has value 2 if every variable in its content is 2, and 4 otherwise:after two factors 0 or 4 annihilate to 4, while a product 2 persists only on another 2.

The first four word laws give the same long-word normal form as their content:sort variables, delete duplicates while length exceeds max(3, content size), and move the extra occurrence between the two letters at length 3 using \(xxy=xyy\). Adjacent swaps use \(xyz=yxz\) at the first two positions with a nonempty suffix, and \(xyz=xzy\) later with a nonempty prefix. Duplicate deletion uses \(xxyz=xyz\) with nonempty y, z. Singleton content retains its cube; \(xx=xxx\) is false at \(x=1\). Thus any two long words with the same content are derivably equal. Each displayed law follows from the explicit word values; all assignments to their at most three variables are independently checked.

\par\smallskip\noindent\textit{Exact criterion from original summands.}\quad 
Let C be the original content union, L the original literal variables, E the original quadratic pairs, H their heads, D their finals, and K the content union of original long words. Build the undirected graph on C with edges E (loops allowed), retaining isolated vertices. Let A be the union of connected components meeting K, and put \(J=H\) union(L intersection D). Then:

\textbullet{} A literal h is absorbed iff h belongs to L union H union A.

\textbullet{} A long word is absorbed iff its content is contained in A.

\textbullet{} A quadratic ht is absorbed iff both variables belong to A, or h belongs to J and h, t belong to the same component.

In particular, a literal alone never supplies its square. The isolated literal case is why J cannot be replaced by H union L. A variable occurring solely as a quadratic final need not be recoverable as a literal:\(y\leq xy\) fails at \(x=1\), \(y=3\).

\par\smallskip\noindent\textit{Necessity via explicit evaluations.}\quad 
For a target containing a variable outside C, assign all C variables 2 and every outsider 4. Every right term is 2 and the target 4.

For the other cases all target variables are in C. A non4 value of an original long word forces all its variables 2. Across a non4 edge either both endpoints are 2, or its head is 1 and its final is 1 or 3. Thus being 2 propagates through the undirected graph, and all A variables must be 2 whenever the right polynomial avoids4.

To separate a failing long target, choose an inactive component meeting its content; assign that component 1 and all other components2. Original edges give 0 or 2, long words 2, literals 1 or 2; their join \(is\leq 0\). The target is 4. The same assignment separates a quadratic whose endpoints are in different components whenever at least one is inactive:its ordered values are(1,\allowbreak{}2) or(2,\allowbreak{}1), both giving 4.

If ht lies in a common inactive component but h is not in J, then h is not an original edge head. If h is an original edge final, it is not an original literal. Assign \(h=3\), the other variables of its component 1, and all other components2. Original edges in that component are(1,\allowbreak{}1) or(1,\allowbreak{}3), giving 0; other original nonlinear words give 2; original literals give 1 or 2. Thus \(P\leq 0\), and the target beginning with 3 is 4. If h is not an edge final either, it is isolated. Since it belongs to C but not A, it occurs only as a literal. Connectivity forces \(t=h\). Assign \(h=0\) and all others 2:the right sum is 0 and hh is 4.

A failed literal h in C is outside L, H, A, so it is solely a final in an inactive component. The preceding assignment \(h=3\), other component variables 1, and all others 2 gives \(P\leq 0\) but \(h=3\) not below 0. This exhausts failure of the criterion, with no finite-variable restriction.

\par\smallskip\noindent\textit{Complete finite derivations.}\quad 
For each original long word and each k of its content, permute k first and apply \(xxx\leq xyz\) with \(x=k\) and nonempty remaining blocks y, z. This produces k cubed. The two laws \(xxx\leq xy+yyy\) and \(yyy\leq xxx+xy\) propagate available cubes in both directions along original edges. Hence every variable in A supplies its cube. The law \(xxy\leq xxx+yyy\) combines two available cubes or a derived long word with a new cube:substitute the current long word for x and the next letter for y, and normalize. By induction this produces every long target with content in A. The law \(xy\leq xxy\) recovers a quadratic on A; \(x\leq xx\) and \(xx\leq xy\) (with y=x) give \(x\leq xx\leq xxx\), so a cube also recovers its literal.

For h in H, an original edge ht gives hh by \(xx\leq xy\). For h in L intersection D, choose an original edge ah and original literal h; \(yy\leq y+xy\) gives hh. Starting with this square, keep the initial h while traversing its undirected component. If ha is already available and ab is an original edge, \(xz\leq xy+yz\) gives hb. If ha is available and ba is an original edge, \(xy\leq xz+yz\) gives hb. Thus every permitted quadratic ht with h in J is obtained. Every h in H has an available square and therefore its literal by \(x\leq xx\); original literals are immediate.

Every construction uses only nonempty substitutions and the original sum as a fixed upper bound. Original graph conditions are never replaced by newly derived edges. Distributive expansion and summandwise absorption in both directions now derive every valid identity. This proves completeness for arbitrary variable counts, word lengths, and numbers of original summands. Catalogue/opposite transfers require separate full HSP proofs in both directions.

\subsection{Twelve-identity bases for two nilpotent-chain extensions of a group}\label{case:316}
Work in the signature \((+,\cdot)\) without constants, relative to the axioms for additively idempotent semirings. As usual, \(p\le q\) abbreviates \(p+q\approx q\). The proofs impose no bound on variables.

\par\smallskip\noindent\textit{Algebras and explicit bases.}\quad 
Let \(G=\{e, g\}\) be the group of order two. On \(A_i=\{e, g, c, a, o\}\), retain the group products, set \(a^2=c\), and set every other product equal to \(o\). In particular \(c\) and \(o\) annihilate every element to \(o\). The additive order of \(A_0\) has \(a<c<o\), \(e<c\), and \(g<o\), with no other strict comparisons; for \(A_1\), interchange the roles of \(e, g\) in the order. Addition is the join. The respective catalog maps, in the displayed order, are \((2,\allowbreak{}3,\allowbreak{}0,\allowbreak{}1,\allowbreak{}4)\) and \((3,\allowbreak{}2,\allowbreak{}0,\allowbreak{}1,\allowbreak{}4)\), giving types 14414 and 14435.

Theorem 1. Each algebra has a twelve-identity basis. Ten common laws are

\CaseDisplay{\begin{array}{rll}
\mathrm{(C)}&xy\approx yx,\\
\mathrm{(N)}&x^3yz\approx xyz,\\
\mathrm{(S)}&x^4y^4\approx x^4+y^4,\\
\mathrm{(U)}&x\le x^3,\\
\mathrm{(O)}&xyz\le x^3+y^3+z^3,\\
\mathrm{(L)}&x^3\le x+xyz,\\
\mathrm{(T)}&(xy)^3\le xy+yzt,\\
\mathrm{(P)}&xt\le xy+yz+zt,\\
\mathrm{(B)}&x\le y+yz+zx,\\
\mathrm{(J)}&x^2\le x+y+xy.
\end{array}}

For \(A_0\), adjoin

\CaseDisplay{\mathrm{(D)}\quad y^2\le x^2+yz,
\qquad \mathrm{(Z)}\quad x^4\le x^3+y^2.}

For \(A_1\), instead adjoin

\CaseDisplay{\mathrm{(Q)}\quad x^3\le x+y^4,
\qquad \mathrm{(G)}\quad xy\le x+y+x^2+y^2.}

Minimality of either list is not claimed.

Associativity follows because any product of at least three factors outside the group-only case is \(o\). For distributivity, the only non-top join of distinct incomparable elements is \(a+e=c\) in \(A_0\) or \(a+g=c\) in \(A_1\). Its two products always have join \(o\), equal to the product with \(c\). Multiplication by \(a\) takes \(a\) to \(c\) and all other elements to \(o\), preserving joins; multiplication by a group element permutes \(e, g\) and sends \(a, c, o\) to \(o\). These translations preserve the specified joins and prove the semiring axioms.

Call a word stable if its length is at least three. A stable word has its group value if all its letters are assigned in \(G\), and has value \(o\) otherwise. A quadratic word has, in addition, value \(c\) when both factors are \(a\). These formulas verify (C),(N),(S), (U),(O). In (L), a nongroup value of \(x\) makes \(xyz=o\); in (T), the exceptional value \(xy=c\) forces \(x=y=a\) and hence \(yzt=o\).

For (P), a three-edge path whose products avoid \(o\) lies entirely in \(G\) or entirely at \(a\). In the group case its endpoint product is the product of its three edge values; if those values agree it equals that value, and otherwise their join is \(o\). The same shared-variable argument proves (B): in the group case, if \(y, yz, zx\) agree, then \(z=e\) and \(x=y\); at \(a\), the right side is \(a+c=c\ge a\). For (J), a pair lies in \(G\), is \((a, a)\), or has product \(o\); all cases give the inequality directly.

In \(A_0\), every square is \(e, c\) or \(o\), in that order. A non-top \(x^2+yz\) is \(e\) or \(c\) and in each case dominates \(y^2\), proving (D). Law (Z) follows from \(e\le y^2\) when \(x\in G\), and from \(x^3=o\) otherwise. In \(A_1\), \(a+e=c+e=o\) proves (Q), and (G) follows by separating the group case, the case \(x=y=a\), and the remaining cases where the right side is \(o\).

\par\smallskip\noindent\textit{Stable words and preparation of a polynomial.}\quad 
Write \(c(w)\) for content and \(\pi(w)\in\mathbb F_2^X\) for exponent parity. For a nonempty word family \(R\), put \(H(R)=\operatorname{aff}_{\mathbb F_2}\{\pi(r):r\in R\}\).

Lemma 2. Laws (C),(N),(S),(U),(O) imply that, for a nonempty family \(R\) of stable words,

\CaseDisplay{c(w)\subseteq\bigcup c(R),\quad\pi(w)\in H(R)
\quad\Longrightarrow\quad w\le\sum R,}

for words \(w\) of every positive length.

Proof. Commutation and (N) equate stable words having the same content and parity. To justify normalization without a length-two exception, use (N) backwards to raise each exponent by two while retaining at least two other positions; two stable words with the same data can be raised to the same large exponent vector. Consequently \(r^3\approx r\) for every stable \(r\).

If stable \(s\) contains \(y\), multiplying (S) by \(s\) gives, after normalization, \(sz^4\approx s+sz^4\). Thus missing variables may be inserted four times to obtain a word above \(s\) of any larger content with the same parity. Choose an odd product of members of \(R\) with the desired parity; include each member twice more to cover all of \(c(R)\). Iterating (O) bounds this product by \(\sum R\). Normalize and use the content-enlargement inequality. For a shorter target, start with \(w\le w^3\) from (U). This proves the lemma. \(\square\)

Every polynomial may be prepared as follows. Start with its stable summands and quadratic edge graph. Whenever an edge \(xy\) meets the content of a stable word, factor the latter as \(yzt\) with \(z, t\) nonempty; (T),(U) replace \(xy\) by its stable cube. Repeat along the finite graph. Then for a linear summand \(x\) in the resulting stable content, factor a stable word as \(xyz\); (L),(U) replace \(x\) by \(x^3\).

The resulting equivalent polynomial has the form

\CaseDisplay{P=\sum R+\sum_{D\in\mathcal D}
\left(\sum_{xy\in E_D}xy+\sum_{x\in L_D}x\right)
+\sum_{x\in L_0}x.\quad\textup{(\ensuremath{1})}}

Here \(R\) is stable, \(W=c(R)\), each \(E_D\) is a nonempty connected graph on \(D\) (loops allowed), \(L_D\subseteq D\), and the sets \(W\), the different \(D\), and \(L_0\) are pairwise disjoint. The family \(R\) or the graph family may be empty. Preparation preserves total content and the affine hull of all word parities. Unlike a linear-stable incidence, a linear-quadratic incidence alone does not justify replacing either summand by its cube.

Put

\CaseDisplay{H=H(R)\quad(R\ne\varnothing),\qquad
K_D=\operatorname{aff}\bigl(\{\epsilon_x+\epsilon_y:xy\in E_D\}
                         \cup\{\epsilon_x:x\in L_D\}\bigr),}

and write \(H=\varnothing\) when \(R\) is empty. A free graph component \(D\) is marked if \(0\in K_D\). Let \(K=H(P)\) denote the affine hull of all parities in (1), including the independent linear summands.

\par\smallskip\noindent\textit{Graphs with specified linear summands.}\quad 
Lemma 3. Laws (P),(B),(J) derive every linear or quadratic word supported in \(D\) whose parity lies in \(K_D\).

Proof. Law (P) shortens a walk of odd length to its endpoint edge; induction on length gives every such edge. Law (B) transports a listed linear summand along two consecutive edges; iteration transports it along every even walk.

If the graph is bipartite and all its specified linear vertices lie in one color class, there is a functional taking value one on all edge vectors and on all specified unit vectors. When there are specified vertices, choose their class to have value one; when there are none, either color class may be chosen. Therefore \(0\notin K_D\). In this case the eligible quadratics are exactly the pairs of opposite colors, and the eligible linears are exactly the vertices of the specified color class if that class is present. Odd-walk shortening and even-walk transport derive all of them. For the linear-algebra assertion, the edge span of a connected graph is the even-coordinate-sum subspace; adding any unit vector spans the whole space. The affine hull is the slice of that linear span where the indicated color functional takes value one.

Every remaining case is marked: either an odd closed walk exists, or specified linear vertices lie in both classes of a bipartite graph. In the former case (P) gives a square. In the latter, an odd path between two specified vertices \(x, y\) gives \(xy\), and (J) gives \(x^2\le x+y+xy\). A loop now converts any walk into an odd walk by a detour to the loop, so (P) derives all pairs inside \(D\). If a specified linear vertex exists, the same detour makes an even walk to any vertex and (B) derives every linear. The space \(K_D\) is the edge span when there are no specified vertices and is the whole parity space otherwise. These are exactly the derived short words. \(\square\)

Two elementary affine observations will also be useful. The blocks in (1) have disjoint coordinate supports. Hence \(0\in K\) precisely when \(0\in H\) or some \(D\) is marked. Indeed an odd zero sum of generators has an odd zero subsum in at least one block; independent unit vectors cannot provide one. If \(0\notin K\), the vectors of \(K\) supported in a particular nonempty block are just that block's affine hull. If \(0\in K\), they instead form its linear span. These assertions follow by sorting an odd vector sum by blocks; conversely an odd zero relation changes an even representation to an odd one without changing its vector.

\par\smallskip\noindent\textit{The first algebra.}\quad 
Proposition 4. In \(A_0\), a word is absorbed by (1) precisely when it is a listed linear summand, or its parity is in \(K\) and one of the following holds:

\textbullet{} its content is contained in \(W\);

\textbullet{} it has length at most two and its content is contained in one free graph component \(D\).

Necessity. An absent variable is separated by assigning it \(o\) and all right variables \(e\). Restricting assignments to \(G\) shows that every absorbed parity lies in \(K\): a functional separates any vector outside the affine hull by giving all right words one group value and the target the other.

For a nonlinear target containing an independent linear variable, assign that variable \(c\) and all others \(e\). The target is \(o\) but the right side is at most \(c\). If a stable target meets a free graph component, assign all of that component \(a\) and all other variables \(e\). Its right edges have value \(c\), its specified linears value \(a\), and all other right words value \(e\); thus again \(P\le c\) and the target is \(o\). The same assignment separates a quadratic with one endpoint in that component and the other outside. These cases establish every structural restriction in the statement.

Derivability. If \(0\notin K\), the disjoint-block observations reduce the allowed parities to \(H\) in \(W\) and \(K_D\) in each \(D\). Lemmas 2 and 3 apply.

If \(0\in K\), first derive some square \(q^2\le P\). A zero vector in \(H\) supplies a fourth power by Lemma 2; a marked graph component supplies a quadratic square by Lemma 3. Law (D) now supplies a square at every free graph vertex: for any incident edge \(yz\), \(y^2\le q^2+yz\le P\). Adding loops changes every graph affine hull to its linear span, without changing its connected components. Lemma 3 gives exactly the required short words in those components.

For each \(r\in R\), law (Z), at \(x=r, y=q\), yields \(r^4\le r^3+q^2\approx r+q^2\le P\). If \(R\) is nonempty this adds zero to its affine hull, enlarging it to the linear span. Lemma 2 then derives every permitted word in \(W\). No step uses an empty polynomial or an empty multiplicative substitution. This proves the proposition. \(\square\)

\par\smallskip\noindent\textit{The second algebra.}\quad 
Proposition 5. In \(A_1\), if \(0\in H\), then \(w\le P\) is equivalent to \(c(w)\subseteq c(P)\) and \(\pi(w)\in K\). If \(0\notin H\) (including empty \(R\)), the exact alternatives are:

\textbullet{} \(w\) is a listed linear summand;

\textbullet{} \(c(w)\subseteq W\) and \(\pi(w)\in H\);

\textbullet{} \(w\) has length at most two, is supported in one unmarked \(D\), and \(\pi(w)\in K_D\);

\textbullet{} \(w\) has length at most two, is supported in the union \(M\) of marked free components, and its parity belongs to \(\bigoplus_{D\text{ marked}}K_D\).

Derivability. If \(0\in H\), Lemma 2 supplies \(h^4\le P\) for a nonempty word \(h\). By (Q),(U), every short right summand \(s\) may be replaced by \(s^3\), since \(s^3\le s+h^4\le P\) and \(s\le s^3\). The resulting entirely stable polynomial has its original total content and affine hull, so Lemma 2 proves the first assertion. Necessity follows from content and group separation.

Suppose \(0\notin H\). Lemmas 2 and 3 derive the first three alternatives and all short targets confined to one marked component. For a quadratic \(xy\) meeting two different marked components, the direct-sum condition requires both unit vectors separately to lie in their component spaces. By Lemma 3, \(x, y, x^2, y^2\) are all absorbed. Law (G) gives \(xy\le P\). These are all additional short targets allowed by the direct sum; a linear has singleton support. This proves sufficiency.

Necessity when \(0\notin H\). If \(R\) is nonempty, choose a group assignment on \(W\) giving every stable summand value \(g\). Give every free graph vertex value \(a\) and every independent linear variable value \(c\). The right value is at most \(c\), because \(a, g<c\). This separates a stable target meeting a free component or an independent linear variable, a quadratic meeting such a linear variable, and a quadratic with one endpoint in \(W\) and the other in a free component. The argument also works when \(R\) is empty, omitting that block. Missing variables are separated by assigning them \(o\).

If a target supported in \(W\) has parity outside \(H\), refine the group assignment there so that all right stable words are \(g\) but the target is \(e\). Affine separation permits this because \(0\notin H\). Keep the remaining variables as above. Since \(e\nleq c\), absorption fails.

Consider an unmarked free component \(D\). It admits a group assignment giving every edge and every specified linear value \(g\). Assign other free components \(a\), independent linears \(c\), and \(W\) as above. A quadratic meeting \(D\) and a different free component now has value \(o\), while \(P\le c\). A short target inside \(D\) whose parity is outside \(K_D\) can instead be made \(e\) while all its right generators remain \(g\). Again \(e\nleq c\). Thus a short target meeting more than one free component can involve only marked components, and unmarked ones retain the stated affine condition.

Finally, suppose a short target supported in \(M\) fails the direct-sum condition. Its parity restriction to some marked component lies outside that component's linear space \(K_D\). Assign group values using a functional annihilating its generators and taking value one on the target restriction. Give every other variable value \(e\). Every right word is now \(e\), whereas the target is \(g\). This excludes every remaining failure of the criterion and proves the proposition. \(\square\)

\par\smallskip\noindent\textit{Completeness and verification.}\quad 
After distribution, an identity is equivalent to absorption of each summand in both directions. Preparation is derivable from the displayed laws. Propositions 4 and 5 derive every valid absorption of a word by a prepared polynomial. This proves both bases in Theorem 1 for arbitrary variable sets and summand counts.

Independent finite checks verify the source tables and every displayed law on all assignments. The semantic criteria have also been checked against all polynomial subsets through two variables, all right polynomials with at most three summands on three variables, further samples, and complete word-function closures through four variables. These checks support the unrestricted proofs and are not used as a substitute for them. Catalog transfers require separate bidirectional variety certificates. Neither result decides all subvarieties or the full five-element classification problem.

The source manuscript for this case credits Aifa Wang and Lili Wang, School of Mathematical Sciences, Chongqing University of Technology.

The source manuscript for this case acknowledges support from the National Natural Science Foundation of China (12371024), the Science and Technology Research Program of Chongqing Municipal Education Commission (KJZD-K2024011102) and the Chongqing Natural Science Foundation Innovation and Development Joint Fund (Municipal Education Commission) (CSTB2025NSCQ-LZX0067).

\subsection{A marked-component basis for a five-element semiring}\label{case:317}
We work in the constant-free signature \((+,\cdot)\), relative to the additively idempotent semiring axioms. We write \(p\le q\) for \(p+q\approx q\). The result below has no bound on the number of variables.

\par\smallskip\noindent\textit{1. Construction and main theorem.}\quad 
Let \(A=\{a, e, f, b, z\}\). Addition is join in the order with three pairwise incomparable minimal elements \(a, e, f\) below \(b<z\). Define commutative multiplication by

\CaseDisplay{a^2=b,\qquad e^2=e,\qquad f^2=f,}

and let every other product be \(z\). These letters name elements, not constants in terms. The catalog relabeling, in this displayed order, is \((1,\allowbreak{}2,\allowbreak{}3,\allowbreak{}0,\allowbreak{}4)\); thus \(A\) is type 14426.

Every product of at least three factors is \(e\) when all factors are \(e\), is \(f\) when all factors are \(f\), and is \(z\) otherwise. This description proves associativity. Translation by \(b\) or \(z\) is the constant map \(z\). Translation by \(a, e, f\) has respectively the sole exceptional input and output \((a, b),(e, e),(f, f)\); all other inputs map to \(z\). Each exceptional input is minimal, so these maps preserve joins. Hence \(A\) is an additively idempotent semiring.

Theorem. The following five identities are a basis for \(A\):

\CaseDisplay{\begin{array}{rll}
(C)&xy\approx yx,\\
(D)&xyz\approx xy+xz+x^3,\\
(G)&x\le xy,\\
(P)&xy\le xt+ty,\\
(M)&x^3\le xy+y^3.
\end{array}}

No assertion of irredundancy is made.

To check validity, (C) is part of the construction. For (D), if the first factor is not \(e\) or \(f\), its cube already makes the right side \(z\). If it is \(e\) (respectively \(f\)), both sides are that same element exactly when the other factors also are \(e\) (respectively \(f\)); otherwise an edge on the right has value \(z\). The only products below \(z\) are the three displayed diagonal products, which proves (G). For (P), if either right product is \(z\), absorption is automatic. Otherwise both products are diagonal at the same element of \(\{a, e, f\}\), and the left product is the same value. For (M), if the right product is not \(z\), the two variables are equal in \(\{a, e, f\}\), so the right cube is the left cube.

\par\smallskip\noindent\textit{2. Normal polynomials and their graphs.}\quad 
By (D), for every \(n\ge2\),

\CaseDisplay{x y_1\cdots y_n\approx xy_1+\cdots+xy_n+x^3.\quad\textup{(\ensuremath{1})}}

The case \(n=2\) is (D). For the induction step substitute \(y_1\cdots y_{n-1}\) for its middle variable; that substituted word is nonempty. Consequently every polynomial is equivalent to a nonempty sum of linear terms, quadratic terms and cubes. Repeated variables are allowed in (1), and addition is idempotent.

Fix such a normal polynomial \(F\). Let \(L\) be its linear variables, let \(E\) be the undirected graph of its quadratic terms (including loops), and let \(D\) be the variables whose cubes occur. Let \(V(E)\) be the vertices incident with an edge. The graph \(\Gamma\) has vertex set \(V(E)\cup D\) and edge set \(E\), so a cube may be an isolated vertex. Mark the vertices in \(D\), and let \(K\) be the union of all connected components that meet \(D\). Put \(U=L\cup V(E)\cup D\).

Absorption lemma. The identities valid in \(A\) are exactly:

\noindent Left term\enspace\textbar\enspace Necessary and sufficient condition for absorption by \(F\)\par
\noindent \(x\)\enspace\textbar\enspace \(x\in U\)\par
\noindent \(x^2\)\enspace\textbar\enspace \(x\in V(E)\cup D\)\par
\noindent \(xy\), \(x\ne y\)\enspace\textbar\enspace \(x, y\) belong to the same component of \(\Gamma\)\par
\noindent \(x^3\)\enspace\textbar\enspace \(x\in K\)\par
In particular, distinct marked components are not merged.

Necessity. Variables not specified below are assigned \(e\). If \(x\notin U\), assign \(x=z\); then \(x=z\) while \(F=e\). If \(x\notin V(E)\cup D\), assign \(x=b\). Its possible linear occurrence makes \(F\) at most \(b\), whereas \(x^2=z\). For a quadratic with an endpoint outside the graph, the same assignment gives left value \(z\) and right value at most \(b\). If both endpoints are graph vertices in different components, assign all vertices of one component \(f\) and all other variables \(e\). Every right term is \(e\) or \(f\), so \(F\le b\), while the left product is \(z\). If \(x\) lies in an unmarked component, assign all its vertices \(a\) and all remaining variables \(e\). Right edges in that component are \(b\), right cubes elsewhere are \(e\), and right linear terms are \(a\) or \(e\). Thus \(F\le b<x^3=z\). If \(x\) is outside the graph, the assignment \(x=b\) already separates its cube.

Sufficiency. It is useful first to check the semantic picture. If \(F=z\), there is nothing to prove. If \(F\le b\), every edge has both endpoints assigned the same member of \(\{a, e, f\}\). Thus each component containing an edge is uniform, and a marked component can only take \(e\) or \(f\). An isolated marked vertex also takes \(e\) or \(f\). A sum \(F=e\) or \(F=f\) forces every occurring variable to have that same value. A sum \(F=a\) can contain no nonlinear term. If \(F=b\), the four indicated targets have values at most \(b\): a linear occurring variable has value at most \(b\), a graph vertex squares to \(b, e\) or \(f\), a connected pair has the same value and multiplies to \(b, e\) or \(f\), and a marked component cubes to \(e\) or \(f\). These cases prove all asserted absorptions. \(\square\)

\par\smallskip\noindent\textit{3. Equational completeness.}\quad 
We now derive every absorption in the lemma from the five laws. Putting the two endpoints of (P) equal yields \(x^2\le xy\). Putting all three variables of (D) equal yields \(x^2\le x^3\). Consequently every graph vertex has its loop absorbed by \(F\), including isolated cube vertices. If a path joins two distinct vertices, repeated use of (P) derives the edge joining its endpoints from its path edges. Thus the entire clique, with loops, on each component is absorbed by \(F\).

Law (M) propagates a cube along an edge: from \(xy\le F\) and \(y^3\le F\) it gives \(x^3\le F\). Induction on path length supplies every cube on a marked component. Finally (G) extracts each endpoint from an edge, and gives \(x\le x^2\le x^3\) at a cube. It follows that every variable of \(U\) is absorbed. These arguments derive precisely all four rows of the lemma using nonempty substitutions.

For any valid identity, normalize both sides by (1). Each summand of either side is absorbed by the other, so the lemma and the preceding derivations give both inequalities and hence the identity. This proves the theorem with arbitrarily many variables and summands. \(\square\)

\subsection{Nonfinite basis from a triangle with a distant cubic witness}\label{case:318}
The signature is (+, multiplication), without constants. Write \(u\leq P\) for \(u+P=P\). Let S be catalogue 14450, with addition 00004/\allowbreak{}01004/\allowbreak{}00224/\allowbreak{}00234/\allowbreak{}44444 and multiplication 44444/\allowbreak{}40444/\allowbreak{}44424/\allowbreak{}44234/\allowbreak{}44444. We prove that S is nonfinitely based, allowing arbitrary polynomial substitutions and term contexts.

\par\smallskip\noindent\textit{Word values and exact absorption.}\quad 
A literal has its assigned value. A quadratic all of whose occurrences are 1 has value 0. For any nonlinear word all of whose occurrences are in\{2,\allowbreak{}3\}, the value is 3 with no 2,\allowbreak{}2 with exactly one 2, and 4 with at least two2 occurrences. Every other nonlinear word has value 4. In particular three or more occurrences of 1 give 4, unlike two occurrences. This is direct at length two. Multiplication by 3 preserves\{2,\allowbreak{}3,\allowbreak{}4\},\(2*2=4\), \(1*1=0\), and 0 times any element is 4; induction proves the formula.

Define data from original summands P. Let C be their content union. Let G have vertices occurring in nonlinear original words and edges joining distinct co-occurring variables. Let F be the variables repeated in an original word. Let A be the union of G components containing an original word of length at least three. A literal x is absorbed exactly when x belongs to C. A nonlinear target w is absorbed exactly when:

1. Every target variable is a vertex of G.

2. Every repeated target variable belongs to F.

3. The target variables outside F form a clique in original G.

4. Either all target variables belong to A, or w is quadratic and its endpoints belong to one G component.

To prove sufficiency, suppose P is not 4. Every nonlinear component is entirely 1 or entirely in\{2,\allowbreak{}3\}. An all 1 component has no long word, so every active component is in\{2,\allowbreak{}3\}. In that case all repeated original variables are 3 and the 2-valued vertices form an independent set in G. Requirements1--3 make any target in\{2,\allowbreak{}3\} have value 2 or 3, dominated by an original nonlinear word containing the possible unique 2 occurrence. If the target is quadratic in an all 1 component, its value 0 is supplied by an original quadratic in that component. Requirement4 excludes all other 1-valued cases. Literal recovery follows directly from its original containing word.

For necessity assign unmentioned variables 3. A missing literal outside C is separated by assigning it 4. For a nonlinear target with a variable outside G assign that variable 0; it affects only possible right literals, so the right side is at most 0 while the target is 4. A repeated target variable outside F can be assigned 2, keeping all others 3. No original word repeats it, and the target is 4 while the right side is at most 2. If two unanchored target variables have no original edge, assign both 2 and all others 3; again no original word contains both or repeats either. If the first three conditions hold but the last fails, choose a target variable in an inactive component, assign that whole component 1, and leave all other variables 3. Its original nonlinear terms are quadratic and give 0. A long target containing 1 is 4; a disconnected quadratic mixes1 with 3 and is 4. The right sum is at most 0. These separators exhaust all failures for arbitrary finite P and w.

\par\smallskip\noindent\textit{A family of critical identities.}\quad 
For \(n\geq 1\) take distinct variables a, b, c, \(x_{1}\),..., \(x_{n}\), e, f and put \(x_{0}=a\). Let

{\small\par\smallskip
\CaseLine{\(P_{n} = ab+ac+bc + sum_{i=0}^{n-1} x_{i} x_{i+1} + x_{n} e f.\)}
\smallskip}

Then S satisfies \(abc\leq P_{n}\). The original graph consists of a triangle on\{a, b, c\}, a path of n edges from a to \(x_{n}\), and a terminal triangle on\{\(x_{n}\), e, f\} supplied by the original cubic. No original variable repeats. The target is a clique, and its component meets the terminal cubic, so the exact criterion applies.

The words absorbed by \(P_{n}\) are exactly all its variable literals, all graph edges as squarefree quadratics, and the two squarefree cubics abc and \(x_{n}\) e f, up to commutation. Indeed F is empty, so no repeated word can be absorbed, and the graph has just two cliques of size three and none of size four. This statement holds at every word length.

Let \(D_{n}\) be these allowed words except abc. If a polynomial Q with summands in \(D_{n}\) absorbs abc, Q must contain the terminal cubic, the three target triangle edges, and every one of the n connecting path edges. The clique condition forces the target edges. The only available long original term is the terminal cubic, so activity of a forces all path edges, each a cut edge of the graph. In particular no variable belongs to every term of such a Q.

\par\smallskip\noindent\textit{The bounded-variable obstruction.}\quad 
Fix k and choose \(n>k\). Suppose valid identities of S using at most k variables could derive \(abc\leq P_{n}\). Distributively expand all terms along a derivation from \(P_{n}\) to \(P_{n}\)+abc. Every intermediate polynomial is identical to \(P_{n}\) in S, so each summand is absorbed by \(P_{n}\). Immediately before the first appearance of abc (up to commutation), all summands therefore belong to \(D_{n}\).

The introduced abc comes from a source word w in one side of a valid identity, a substitution sigma of nonempty polynomials, and fixed monomial contexts p, q. Write R for the other side of the source identity. Then \(w\leq R\) holds in S, and p sigma(R) q consists entirely of words in \(D_{n}\) and absorbs abc. The preceding lemma forces its terms to contain all the target edges and the disjoint terminal cubic. If p or q were nonempty, every such term would share a context letter, which is impossible. Thus both contexts are empty and sigma(R) is a polynomial with summands in \(D_{n}\). This covers a general term context by selecting the distributive monomial branch that produces abc; the corresponding branch of the other side already occurs before the step.

No variable repeats in any source word of R: choosing the same nonempty image monomial at two occurrences would give a repeated letter in an expanded right word, whereas every member of \(D_{n}\) is squarefree. Thus \(F_{R}\) is empty. Every source word has length at most three. Every source cubic expands only to \(x_{n}\) e f; all its variable images must consequently be singleton letters in that terminal triangle.

Consider the original nonlinear graph \(G_{R}\). For every variable v of this graph choose one letter h(v) occurring in any image monomial sigma(v). Whenever u, v are adjacent in \(G_{R}\), they co-occur in a source right word. Choose image monomials containing h(u), h(v), and arbitrary image monomials for the other source factors. The resulting word is a squarefree clique of the target graph. Therefore h(u) and h(v) are distinct and adjacent. The independent choices of h define a graph homomorphism \(G_{R}\to G_{n}\). This remains true for large polynomial images, repeated occurrences of variables in different source summands, and images of length two or three. No assumption of monomial substitution is used.

Because an expansion of w has length three, w has at most three letters. The exact criterion for \(w\leq R\) forces its variables to occur in R and precludes repeated source letters (\(F_{R}\) is empty). If w is a literal, an image abc would occur in sigma(R) if that variable occurs as a literal there, or create a word of length at least four if it occurs nonlinearly. Both are forbidden. If w has length two, its variables form an original edge of \(G_{R}\), since \(F_{R}\) is empty. They co-occur in a source quadratic or cubic. The chosen two image monomials forming abc would then give abc already on the right or a word of length at least four. Again this is impossible.

Consequently w has exactly three distinct source variables, mapped in the selected expansion to a, b, c. The long-word clause for \(w\leq R\) puts each of these variables in a component of \(G_{R}\) containing a source cubic. Choose the variable mapped to a and a shortest source path from it to a vertex of a source cubic. That path has at most k-1 edges. Choose h at its first vertex to be a, and choose representatives at all other vertices arbitrarily as above. At its endpoint h is in\{\(x_{n}\), e, f\}, because it belongs to a source cubic. The resulting walk in \(G_{n}\) has length at most k-1. But the distance from a to the terminal triangle is \(n>k\). Contradiction.

Thus no identity set with a finite variable bound derives every critical identity. Every finite basis has such a bound, proving S nonfinitely based. The proof includes the semiring axioms, polynomial substitutions, and arbitrary contexts. It does not infer NFB solely from a proper NFB subvariety.

\par\smallskip\noindent\textit{Transfer with an all-length test.}\quad 
For a finite semiring B, if S belongs to HSP(B) and B satisfies \(abc\leq P_{n}\) for every \(n\geq 1\), then B is NFB: a finite basis of B would be a bounded-variable identity set valid in S deriving all these identities. Both hypotheses must be verified.

The all-length condition is decidable by a finite-state closure. A state (first, last, q) stores the first and last values of a nonempty path of at least one edge and the additive join q of its ordered edge products. Initial states are(x, y, xy). Append z by(last, \(q)\to (z\), q+last*z). Close the finite set. For every reached state and every assignment to b, c, e, f, check first*b*c \ensuremath{\leq} first*b+first*c+b*c+q+last*e*f. Induction proves that the closure covers exactly all assignments at all \(n\geq 1\). This procedure also works for noncommutative candidate B.

\subsection{A fourteen-identity basis for catalogue algebra 14452}\label{case:319}
The language has addition and multiplication, without constants. All substitution words below are nonempty.

\par\smallskip\noindent\textit{1. Tables and word values.}\quad 
The rows of addition are 00004 / 01004 / 00224 / 00234 / 44444; the rows of multiplication are 44444 / 40444 / 44424 / 44434 / 44444. Thus the additive order consists of 3 \ensuremath{<} 2 \ensuremath{<} 0 \ensuremath{<} 4 and 1 \ensuremath{<} 0 \ensuremath{<} 4, with 1 incomparable to 2 and 3. The only non-4 products are \(1*1=0\), \(2*3=2\) and \(3*3=3\).

A literal has its assigned value. A quadratic word has value 0 if both positions are 1; otherwise it has value 2 or 3 exactly when the first position has that value and the second is 3. All other cases give 4. A word of length at least three has value 2 or 3 exactly when its head has that value and every later position is 3; all other cases give 4. The last statement follows by multiplying the three possible non-4 binary products by the next factor: 0 always goes to 4, whereas 2 and 3 survive precisely on a right factor 3. This proves the formula for every length by induction, including repeated variables.

\par\smallskip\noindent\textit{2. The proposed basis.}\quad 
In addition to the additively idempotent semiring axioms, use these fourteen identities. Write \(p\leq q\) for \(p+q=q\).

1. xyz = xzy.

2. xyzz = xyz.

3. x \ensuremath{\leq} xy.

4. y \ensuremath{\leq} xy.

5. yy \ensuremath{\leq} xy.

6. xz \ensuremath{\leq} xy + yz.

7. xz \ensuremath{\leq} yx + yz.

8. xyy \ensuremath{\leq} xyz.

9. yyy \ensuremath{\leq} xyz.

10. xyy \ensuremath{\leq} xy + yyy.

11. yyy \ensuremath{\leq} xy + xzz.

12. xxx \ensuremath{\leq} yx + xzz.

13. xyy \ensuremath{\leq} xzz + yyy.

14. xyz \ensuremath{\leq} xyy + xzz.

Their validity follows directly from the word-value formula and additive order; the independent verifier exhaustively checks every assignment involving their at most three variables. These finite checks establish soundness of the fourteen displayed laws, not completeness. Completeness is proved below for arbitrary variable sets, lengths and summand counts.

Laws 1 and 2 give a normal form for every word: literals and quadratics remain unchanged; a longer word retains its head and the set of variables occurring after its head. List that tail set in a fixed order, duplicating its member if the tail set is a singleton. Adjacent tail swaps use 1 in a context. To remove a repeated tail variable, move two copies to the end and apply 2, keeping a nonempty initial x-block. Stop at length three. This handles all repeated-head cases as well: the head is never moved.

\par\smallskip\noindent\textit{3. Exact absorption criterion.}\quad 
Let P be a nonempty sum of original words. Let C be their content union. Let Q be the directed edges ab supplied by original quadratic words. Let \(T_{Q}\) be their right endpoints, and let \(T_{R}\) be the tail variables of original words of length at least three. Put \(T=T_{Q}\) union \(T_{R}\).

The undirected graph underlying Q has as vertices all variables incident with an edge; loops are retained. Each of its components therefore has an actual original edge. A variable occurring in no Q-edge is treated as an isolated variable when defining activation. Mark every variable in an original long word and all vertices in its Q-component. Let N be this set of marked, or activated, variables. No component is marked merely because it has an isolated literal.

Then \(w\leq P\) holds precisely as follows:

\textbullet{} Literal x: x belongs to C.

\textbullet{} Quadratic xy: either y belongs to \(T_{Q}\) and x lies in the Q-component of y, or x belongs to N and y belongs to N intersection T.

\textbullet{} Long w with head h and tail set S: h belongs to N and S is contained in N intersection T.

There is no empty-path witness in the quadratic clause. Even when \(x=y\), its first branch requires an original edge entering y. The construction starts with that edge and law 5. The second branch requires an original long word. This is a newly derived criterion, not a reflexive replacement of any positive-path criterion used for another algebra.

\par\smallskip\noindent\textit{Proof of necessity and sufficiency.}\quad 
An assignment makes P non-4 exactly when every original summand is non-4. On a Q-component, either every variable is 1, or all variables are 2/\allowbreak{}3 and every right endpoint is 3. Indeed every edge has precisely these two alternatives, and an endpoint 1 forces the same alternative through each adjacent edge. Original long words force all their variables into 2/\allowbreak{}3 and their tail variables to 3. Thus exactly the components meeting a long word are forced active, while other Q-components may choose the all-1 alternative. Literal-only variables remain unrestricted among 0,\allowbreak{}1,\allowbreak{}2,\allowbreak{}3.

If a literal is absent from C, assign it 4 and all variables of P to 3. If a long target uses a variable outside N, assign its unactivated Q-component to 1 (or the isolated target variable to 1), with all other variables 3. P is non-4 and the target is 4. If a target tail variable is outside T, assign it 2 and every other variable 3. Again P is non-4 and the target is 4. These prove the necessary long conditions.

For a quadratic target, a tail outside T is separated by the same 2/\allowbreak{}3 assignment. If the tail is in T but the stated alternatives fail, the two variables are in different components and at least one component is unactivated; assign that component 1 and all others 3. The target is 4, and P is non-4. A head outside C can instead be assigned 4. If the tail is from a long word, its component is already activated, so the same-component case necessarily falls under the second alternative.

Conversely, under an assignment with P non-4, the stated conditions make each long target have value 2 or 3. If its head is 2, that head is not a tail variable, and some original word having it as head has value 2. P is therefore at least 2; otherwise the target is 3 and P is at least 3 because some original long word exists. The same reasoning handles an active quadratic target. An inactive quadratic component gives value 0 for both the target and an actual original quadratic edge, so P is 0. Literals follow also from laws 3 and 4, which extract every position of each original word. Assignments with \(P=4\) are automatic.

\par\smallskip\noindent\textit{4. Derivation of every valid absorption.}\quad 
Retain every original summand throughout. Each step only adjoins a word already below the current sum. Normal forms are obtained by laws 1 and 2.

For a literal target in C, take an original occurrence. Law 3 extracts the prefix ending at that occurrence if a suffix follows; law 4 then extracts the last variable of this prefix if a nonempty prefix precedes it. No empty substitution is used.

For a quadratic graph witness with target tail t, choose an original edge at. Law 5 gives tt. Follow a finite undirected Q-path backwards from t. If vt has been obtained and the next original edge is uv, law 6 gives ut from uv+vt. If that edge is vu, law 7 gives ut from vu+vt. Hence xt is obtained for every x in that component. At \(x=t\) the real entering edge at still supplies the initial tt; no path is substituted as an empty word.

Now construct activation witnesses. An anchor for v is a derived long word vzz for some variable z. From each original long word with head h, law 8 gives an anchor hyy. For each of its tail variables t, reorder the tail so that t comes first and apply law 9 to get ttt. The remaining tail block in this substitution is nonempty because the word is long. Thus every initially activated variable has an anchor.

Propagate anchors across each actual original Q-edge uv. If u has an anchor uzz, law 11 gives vvv from uv+uzz, anchoring v. If v has an anchor vzz, law 12 gives vvv from uv+vzz, and law 10 gives uvv from uv+vvv, anchoring u. Iterating over finite paths gives an anchor for every v in N. For each t in N intersection \(T_{Q}\), use an edge ut and the anchor of u in law 11 to obtain ttt. Cubes of \(T_{R}\) were already obtained directly.

Given h in N and t in N intersection T, use the anchor hzz and cube ttt in law 13 to obtain htt. To build a target with several tail variables, induct on the size of its tail set. Law 14 with \(x=h\), y the first tail variable and z the nonempty word of remaining tail variables combines the elementary anchor hyy and the already obtained hzz. Law 1/\allowbreak{}2 normalization identifies the latter with the induction hypothesis. Singleton tails are exactly the elementary anchors. This yields every valid long target. A quadratic active target ht follows from htt by law 3 with \(x=ht\) and \(y=t\).

Finally an identity \(P=R\) is equivalent to the finitely many absorptions of each original summand of P into R and each original summand of R into P. The preceding derivation treats all such valid absorptions without any bound on the number of variables, word length, graph size or polynomial size. Therefore the fourteen displayed laws form a complete finite basis for 14452 relative to the additively idempotent semiring axioms.

\subsection{A complete twenty-five-identity basis for 14453}\label{case:320}
All words and polynomial summands below are nonempty. Write u \ensuremath{\leq} P for \(u+P=P\). The signature has addition and multiplication, with no named zero or identity. The basis is relative to the six usual additively idempotent semiring axioms. Reversing every word gives a basis for the opposite algebra 14465. Other transfers require full HSP checks.

\par\smallskip\noindent\textit{1. The algebra and the proposed basis.}\quad 
In the catalogue labels 0,\allowbreak{}1,\allowbreak{}2,\allowbreak{}3,\allowbreak{}4, addition has table

{\small\par\smallskip
\CaseLine{\ttfamily 0 0 0 4 4}
\CaseLine{\ttfamily 0 1 0 4 4}
\CaseLine{\ttfamily 0 0 2 4 4}
\CaseLine{\ttfamily 4 4 4 3 4}
\CaseLine{\ttfamily 4 4 4 4 4}
\smallskip}

The only products different from 4 are \(11=0\), \(23=2\) and \(33=3\). Thus 1 and 2 lie below 0, while 3 is incomparable with all three; 4 is the additive maximum. The twenty-five identities, numbered as in the executable derivation, are

{\small\par\smallskip
\CaseLine{\(0  xyz = xzy\qquad 1  xyzz = xyz\)}
\CaseLine{\(2  x \leq  xy\qquad 3  yy \leq  xy+yz\)}
\CaseLine{\(4  yy \leq  y+xy\qquad 5  yy \leq  xx+xy\)}
\CaseLine{\(6  xx \leq  yy+xy\qquad 7  xy \leq  xx+yy\)}
\CaseLine{\(8  xt \leq  xy+zy+zt\qquad 9  xyy \leq  xyz\)}
\CaseLine{\(10  xyy \leq  xy+xzz\qquad 11  xyy \leq  xy+tyy\)}
\CaseLine{\(12  xyy \leq  xzz+tyy\qquad 13  xyz \leq  xyy+xzz\)}
\CaseLine{\(14  xxx \leq  xy+txx\qquad 15  xxx \leq  yx+xzz\)}
\CaseLine{\(16  xxx \leq  xx+xyy\qquad 17  yyy \leq  yy+xyy\)}
\CaseLine{\(18  yyy \leq  xyy+yzz\qquad 19  yyy \leq  y+xyy\)}
\CaseLine{\(20  yyy \leq  xxx+y\qquad 21  yyy \leq  xxx+yz\)}
\CaseLine{\(22  yyy \leq  xxx+zy\qquad 23  yyy \leq  xxx+zyt\)}
\CaseLine{\(24  xyxyxy \leq  xxx+yyy\)}
\smallskip}

Their validity is verified directly from the displayed finite tables, for every assignment to every variable in each identity. This is a finite soundness check; completeness is the unrestricted proof in sections 2--5.

\par\smallskip\noindent\textit{2. Words and normal forms.}\quad 
A literal evaluates to its assigned element. A quadratic word evaluates to 0 when both positions are 1. In every other nonlinear case its value is its head value if that value is 2 or 3 and all remaining positions are 3, and is 4 otherwise. In particular the exceptional quadratic all-1 value does not persist to length three.

Rules 0 and 1 normalize any word of length at least three by retaining the head and the set of variables after the head. Sort the latter set; if it is a singleton, retain two copies of its variable. Literals and quadratics are unchanged. To realize this syntactically, rule 0 swaps adjacent tail letters with the nonempty preceding block substituted for x. If a tail variable is repeated, move its two copies to the end and apply rule 1 with nonempty x, y, while the length exceeds three. Sort the tail again. At length three two equal tail letters must remain. This proves normalization in both directions without any empty substitution.

\par\smallskip\noindent\textit{3. Exact absorption criterion.}\quad 
Fix the original right polynomial P. Let C be its content, L its literal variables, Q the directed pairs from its quadratic words, and R its words of length at least three. Let H, T be the head and tail sets of Q; let rh be the heads of R and rt all variables after heads in R. Take the undirected connected components of Q on C, including isolated vertices. Define

\textbullet{} A: the union of components meeting (H intersect T) union (L intersect T).

\textbullet{} N: the union of components meeting rh union rt.

\textbullet{} BH=(H intersect N) union rh, and BT=(T intersect N) union rt.

\textbullet{} A conflict occurs if BH intersect BT, A intersect N, or L intersect BT is nonempty.

An absorption \(w\leq P\) holds precisely when content(w) is contained in C and either there is a conflict, or the following length-specific condition holds:

\textbullet{} Literal x: x belongs to L union H union rh union A.

\textbullet{} Quadratic xy: either x, y both belong to A; or x belongs to BH and y to BT; or x belongs to H and y to T in the same Q-component.

\textbullet{} Length at least three: the head belongs to BH and every tail variable belongs to BT.

Here component connectivity concerns undirected quadratic edges. It does not assert a zero-length word path or replace any earlier positive-path criterion for another source.

\par\smallskip\noindent\textit{Proof of the criterion.}\quad 
Every product or sum has value 3 exactly when all variables appearing in it are 3. Hence content containment is necessary: assign a missing target variable 4 and all variables of C value 3. It also handles \(P=3\). \(P=4\) imposes no further restriction.

For P to lie in \{0,\allowbreak{}1,\allowbreak{}2\}, no summand can be 3 or 4. Each quadratic component must therefore be in one of two modes: every vertex 1, or every head 2 and every tail 3. The second mode is impossible if a vertex is both a head and a tail, or a tail is a literal. These are exactly the components of A. A long word forces its head 2 and its tail variables 3; these assignments force the oriented mode throughout every component of N. Such demands are consistent exactly when the three specified conflict sets are empty. Isolated variables appearing only as literals are independently allowed values 0,\allowbreak{}1,\allowbreak{}2. This describes ALL assignments with P in \{0,\allowbreak{}1,\allowbreak{}2\}, not just a selected testing family.

In the consistent case the variable categories are: A fixed at 1; BH fixed at 2; BT fixed at 3; each remaining nontrivial component is a free pure bipartite component with a common choice (heads, \(tails)=(1\),1) or (2,\allowbreak{}3); and isolated literals have independent values 0,\allowbreak{}1,\allowbreak{}2. A and the forced head/tail sets are disjoint. This description proves the displayed sufficient conditions directly. Whenever a quadratic component is all 1, P includes a quadratic value 0. A literal forced to 1 is therefore \ensuremath{\leq}P; all original heads also lie below their original words by rule 2.

For necessity of the literal condition, a remaining variable is a tail in either a forced or a free oriented component: choose that mode to make the variable 3 and \(P\leq 0\). For a forbidden quadratic pair, an isolated literal position can be set 0, giving product 4. Among the other categories, the only universally admissible ordered value pairs are (1,\allowbreak{}1) from two A vertices; (2,\allowbreak{}3) from BH to BT; and the coupled choices (1,\allowbreak{}1)/(2,\allowbreak{}3) from a head to a tail in ONE free component. All other ordered category pairs permit a choice of the independent free modes giving product 4 or 3. This includes two tails of the same free component (choose its oriented mode), and two heads of that component (their product becomes 4). Repeated positions are covered by the same check. Thus the quadratic list is exact. For a long target, assign A and all free components1, BH2, BT3, and isolated literals 0. A long word has value below 4 under this assignment only with head in BH and all remaining variables in BT. This proves necessity. In the conflict case P can only be 3 or 4, so content containment is sufficient. The criterion holds for arbitrarily many variables and summands.

\par\smallskip\noindent\textit{4. Realizing the criterion with the twenty-five rules.}\quad 
Retain every original summand. A derived absorption may adjoin its left word to the working polynomial without changing P. Normalization is by the two word laws, and every displayed substitution uses nonempty blocks.

First derive a square for every variable of A. A head/tail overlap gives \(v^{2}\) by rule 3 from its incoming and outgoing edges. A literal tail gives \(v^{2}\) by rule 4. Rules 5 and 6 propagate squares in both directions along quadratic edges. Loops already supply their square. This reaches precisely A. Rule 7 subsequently combines any two such squares into a quadratic.

For each original long word and each of its tail variables t, tail permutation and rule 9 give h \(t^{2}\). There is always a nonempty remaining tail block, including when the tail set is a singleton. Call h \(t^{2}\) an anchor with roles head h and tail t. An edge uv extends these roles: if u is already an anchor head, rule 10 gives u \(v^{2}\); if v is already an anchor tail, rule 11 gives u \(v^{2}\). Propagate across the finite graph.

A conflict yields a cube as follows. An anchor with equal endpoints is already a cube. A square at an anchor head or tail gives a cube by 16 or 17; a literal anchor tail uses 19. A tail that is another anchor head uses 18. A quadratic head that is already an anchor tail uses 14; a quadratic tail that is already an anchor head uses 15. If none of these occurs, every visited quadratic component has its consistent orientation; propagation traverses it and visits all vertices. Consequently the final anchor heads and tails are exactly BH and BT. All intersections defining a conflict are then detected by one of the preceding rules. Conversely each detected collision implies the specified semantic conflict. This proves termination with a cube exactly in the conflict case.

Given one cube \(q^{3}\), rules 20--23 give \(v^{3}\) for every variable of C by choosing an original occurrence of v. A literal uses 20, a first occurrence at the head uses 21, a final occurrence uses 22, and an internal occurrence uses 23 with nonempty prefix and suffix. Using24, cubes of two nonempty words give the cube of their concatenation. Induction along the target letters gives \(w^{3}\). Rule 2 with \(x=w\)\ensuremath{{}^{2}}, \(y=w\) gives \(w^{2}\), and with \(x=w\), \(y=w\) gives w. This realizes every absorption in the conflict case.

Suppose there is no conflict. Literal recovery uses an original literal, rule 2 on an original word with that head, or rule 2 on an A-square. Quadratics inside A use 7. For BH-to-BT pairs, choose one anchor with the desired head and one with the desired tail; rule 12 gives h \(t^{2}\), and rule 2 with \(x=ht\), \(y=t\) gives ht. In a remaining free bipartite component choose an undirected path from the desired head to tail. Its length is odd because the two roles are disjoint. Its first edge is original. Rule 8 extends the derived first-to-current-tail edge across the next two original edges. Induction gives the required head-to-tail quadratic. No orientation is invented and no reflexive path is used.

For a long target obtain h \(t^{2}\) for every desired tail variable using12. Rule 13 merges these tail variables: substitute a singleton for y and the already merged, nonempty tail block for z. The right xzz normalizes to the previously obtained word and xyy is its anchor. For one tail variable the anchor already is the normal target. For two or more, the left xyz accumulates the tail set. The word laws recover the exact original target, including repetitions and arbitrary order.

\par\smallskip\noindent\textit{5. Completeness and scope.}\quad 
Every valid absorption is derived by section 4, using the exact criterion of section 3. Any valid polynomial identity is a finite collection of absorptions, one for each summand in each direction. Distributivity and the semiring axioms turn arbitrary terms into such polynomials. Therefore the twenty-five displayed identities form a complete finite basis, without a bound on variables, word lengths, or polynomial sizes.

\subsection{三阶群的中间零层：两组完整基}\label{case:321}
按顺序 \((e, g, g^2, d, c)\)，到五元素目录 14637、14841 的同构分别为 \((1,\allowbreak{}2,\allowbreak{}3,\allowbreak{}0,\allowbreak{}4)\)、\((3,\allowbreak{}1,\allowbreak{}2,\allowbreak{}0,\allowbreak{}4)\)。收缩 \(\rho\) 和每个乘法平移保持并，因而两个构造都满足半环公理。

定理。 \(Z_e\) 的完整有限基是

\CaseDisplay{\begin{aligned}
xy&\approx yx,&x^4y&\approx xy,\\
x^3y^3&\approx x^3+y^3,&xyzt&\leq x^4+y^4+z^4+t^4,\\
x&\leq x^4,&x^4&\leq x+xy.
\end{aligned}\quad\textup{(\ensuremath{E})}}

对 \(Z_g\)，再加入

\CaseDisplay{x^4\leq x+y^3\quad\textup{(\ensuremath{G})}}

即得完整基。

有效性直接来自三阶群平坦核心。对新元素，\(d^4=c\)；(E) 的最后一式用 \(dy=c\)，而 (G) 在 \(Z_g\) 中用 \(d+e=c\) 以及 \(y^3\in\{e, c\}\)。

证明完整性时，对右多项式 \(P\) 定义全部内容 \(U(P)\)、非线性内容 \(N(P)\)，以及指数模 3 向量的仿射包 \(H_3(P)\)。精确判据如下：

\textbullet{} 两个半环都满足：一次词 \(x\leq P\) 当且仅当 \(\pi_3(x)\in H_3(P)\)。

\textbullet{} 在 \(Z_e\)，非线性 \(w\leq P\) 当且仅当 \(\pi_3(w)\in H_3(P)\) 且 \(c(w)\subseteq N(P)\)。

\textbullet{} 在 \(Z_g\)，非线性 \(w\leq P\) 当且仅当 \(\pi_3(w)\in H_3(P)\)、\(c(w)\subseteq U(P)\)，以及 \(c(w)\subseteq N(P)\) 或 \(0\in H_3(P)\)。最后的``或''不取代前两个条件。

群赋值给出仿射条件的必要性。对 \(Z_e\)，若 \(x\in c(w)\setminus N(P)\)，赋 \(x=d\)、其余变量为 \(e\)，则左侧为 \(c\)，右侧不高于 \(d\)。对 \(Z_g\)，内容遗漏由 \(x=c\)、其余变量为 \(e\) 分离；若某个内部变量 \(x\) 只作为单项出现，且 \(0\notin H_3(P)\)，删去此单项，利用线性泛函使全部其余右词取 \(g\)，再赋 \(x=d\)，左右分别取 \(c, d\)。

下面从有限基推出上述全部候选吸收式，完成充分性与完备性。由 (E)，非线性词只依赖内容和指数模 3，且 \(u^4=u\)。此外，若 \(u\) 非线性、\(x\in c(u)\)，则

\CaseDisplay{u+uz^3=ux^3+uz^3=ux^3z^3=uz^3,}

故可以向非线性词加入任意立方因子而使其不减。任意 \(3r+1\) 个非线性词之积被这些词之和吸收：把前三个新因子与已有乘积代入四词式，归纳即可。于是，若非线性词族 \(Q\) 内容覆盖左词且其模三仿射包包含左向量，就可按仿射表示取一个 \(3r+1\) 因子乘积，再将每个右词补三次以覆盖内容，从而推得该左词被 \(\sum Q\) 吸收。一次左词则先用 \(x\leq x^4\) 化为非线性词。

对 \(P\) 中同时出现在非线性右词中的单项 \(x\)，取其中含 \(x\) 的右词 \(xu\)。由 \(x^4\leq x+xu\)，向 \(P\) 加入 \(x^4\)；又由 \(x\leq x^4\)，可删去该单项。所得非线性族 \(Q\) 的内容是 \(N(P)\)，剩余单项变量完全不出现在 \(Q\) 中。左向量在这些外部单项变量的坐标为零，因此在模三仿射表示中，它们的系数只能为零；或者用重复项语言说，每个外部单项只出现三的倍数次，可全部删去而不改变总系数之和模 3。其余一次项都已由具有相同向量的四次幂替代。上述非线性引理遂推出 \(Z_e\) 的所有非线性候选吸收式，以及两例中 \(0\notin H_3(P)\) 时的相应吸收式。对于一次左词，若它已在右侧就是显然吸收；否则其单位向量属于 \(H_3(P)\) 强制它出现在 \(N(P)\)，故相同推导适用。

对 \(Z_g\)，若 \(0\in H_3(P)\)，同样去掉外部单项的零系数，得到 \(0\in H_3(Q)\)。按零向量的仿射表示构造 \(3r+1\) 个 \(Q\) 词的乘积，其指数都是三的倍数，故是非空词的立方 \(h^3\leq P\)。由 (G) 向 \(P\) 加入所有剩余单项的四次幂，再删去原单项，所得纯非线性词族的内容和仿射包分别为 \(U(P)\)、\(H_3(P)\)。上一段的引理完成全部吸收式推导。

最后，任意多项式恒等式分解为两个方向的逐词吸收式，因此上述六条、七条基均完整。所有论证对任意有限变量集成立；没有采用有限搜索失败推断非有限基，也没有使用未经证明的膨胀有限基传递原则。

\subsection{A 244-template finite basis for 14674}\label{case:322}
Signature (+, multiplication), no constants. Every omitted capital below specifies a separate identity, not an empty substitution. All variables that remain in an identity receive nonempty terms. Reverse all words to obtain the opposite presentation for 14734.

\par\smallskip\noindent\textit{Exact semantics.}\quad 
The addition rows are 00044,\allowbreak{}01444,\allowbreak{}44244,\allowbreak{}44434,\allowbreak{}44444. The only non4 products are \(11=1\),\(12=2\),\(23=2\),\(33=3\). Thus addition has \(1<0<4\), with 2 and 3 each incomparable to 0,\allowbreak{}1 and to one another. A literal keeps its value. A nonlinear word is 1 when all its factors are 1,\allowbreak{}3 when all are 3, and 2 exactly when it has one 2 preceded only by 1's and followed only by 3's. Every other case, including any occurrence of 0, gives 4. This differs from 11619, where an initial 0 followed by 1's retains value 0.

Fix an original nonempty polynomial P. Let C be its content and D the union of the contents of its original NONLINEAR words. Make endpoints Lx, Rx for every x in C; identify all original heads' L-endpoints as H, all finals' R-endpoints as T, and \(Rx=Ly\) for every original adjacent xy. Add directed edges \(Lx\to Rx\). Write a\textasciitilde{}b if the endpoints belong to the same directed strongly connected component.

The exact absorption test is content(w) contained in C, together with content(w) contained in D when w is nonlinear, and

L(head(w))\textasciitilde{}H, R(last(w))\textasciitilde{}T, Rx\textasciitilde{}Ly for each adjacent xy of w.

The D condition concerns the ENTIRE nonlinear target, including its head. This is not the noninitial-content criterion of 11619.

For necessity of C, assign a missing target variable 4 and all original variables 1. For failure of D, assign a forbidden variable 0 and all others 1. It occurs only as an original literal, so \(P=0\) while the nonlinear target is 4. To analyze \(P=2\), encode1,\allowbreak{}2,\allowbreak{}3 as(0,\allowbreak{}0),(0,\allowbreak{}1),(1,\allowbreak{}1). Each original word is 2 exactly when its initial endpoint is 0, its final endpoint 1 and its adjacencies are equal. Thus the original constraints are the displayed identifications, monotonicity \(Lx\leq Rx\), and \(H=0\), \(T=1\).

Every vertex lies on a directed H-to-T path labelled by an original word. If T reaches H, the whole graph is strongly connected and no assignment with \(P=2\) exists. Otherwise the equalities forced by its Boolean monotone assignments are precisely strong-component equalities. For a pair not in one component, choose an upward-closed set containing one endpoint but not the other, add the upward closure of T, and leave H outside. The fact that every vertex lies between H and T ensures such a choice; the resulting1/\allowbreak{}2/\allowbreak{}3 assignment makes every right word 2 and the target different from 2. This proves endpoint necessity and sufficiency for \(P=2\).

If \(P=1\) or 3, all its variables are respectively 1 or 3, and C handles the target. If \(P=0\), each nonlinear original word is 1 and every literal is 0 or 1. Thus all D variables are 1. A permitted nonlinear target is 1; a literal target is 0 or 1, so both are absorbed. \(P=4\) imposes no condition. This proves the criterion for arbitrary variables, lengths and sums.

\par\smallskip\noindent\textit{The finite presentation.}\quad 
Four fixed laws are

{\small\par\smallskip
\CaseLine{\(B\qquad xyx = xxyy\)}
\CaseLine{\(I0  x \leq  yyx\)}
\CaseLine{\(I1  x \leq  xyy\)}
\CaseLine{\(I2  xy \leq  xzzy\)}
\smallskip}

Include every permitted optional-context instance of the following six schemes:

{\small\par\smallskip
\CaseLine{\(M\qquad A x D \leq  A x B + C x D\)}
\CaseLine{\(R\qquad A x t F \leq  A x y B + C z y D + E z t F\)}
\CaseLine{\(S\qquad A x z D \leq  A x y B + y C + z D + E z F\)}
\CaseLine{\(T\qquad D z y B \leq  A x y B + C x + D z + E z F\)}
\CaseLine{\(E\qquad A x \leq  A x y B + C z y D + E z\)}
\CaseLine{\(F\qquad y B \leq  A x y B + C x z D + z E\)}
\smallskip}

Every capital may independently be omitted in all its occurrences, except that in S and T at least one of E, F MUST remain. This supplies a true nonlinear occurrence of z. All lower-case variables remain. The numbers of templates are 16,\allowbreak{}64,\allowbreak{}48,\allowbreak{}48,\allowbreak{}32,\allowbreak{}32, giving 240 schemes plus the four fixed laws. They have at most ten distinct formal variables.

The word-value formula proves B. The three insertion laws are sound: a squared block in a nonzero bridge path must be a loop, and can be removed. If any nonlinear factor is 0 the containing product is 4; the literal left \(cases 0\leq yy0\) and \(0\leq 0yy\) also hold. For each scheme, the target is already an H-to-T walk in the endpoint-identification graph of its right words. Its nonlinear content is covered by original nonlinear words. In S, T the guard EzF ensures precisely this for z, which could otherwise have been a literal-only head or final. In M, if the target is nonlinear then each retained contextual letter and x occur in at least one nonlinear right word. The exact criterion therefore proves every permitted template sound. Exhaustive finite assignment checking independently audits these 244 sound identities.

\par\smallskip\noindent\textit{Finite closure of witnessed endpoints.}\quad 
Keep original summands and adjoin only words proved below P. Work with labels in the original D for all nontrivial constructions. Record heads, finals and adjacent pairs of available words. Heads/finals outside D are unused. R, S, T, E, F implement respectively

{\small\par\smallskip
\noindent (x, y),(z, y),(z, t) \ensuremath{\to} (x, t)\par
\noindent (x, y), head(y), head(z), z in D \ensuremath{\to} (x, z)\par
\noindent (x, y), final(x), final(z), z in D \ensuremath{\to} (z, y)\par
\noindent (x, y),(z, y), final(z) \ensuremath{\to} final(x)\par
\noindent (x, y),(x, z), head(z) \ensuremath{\to} head(y).\par
\smallskip}

For S, T choose an original nonlinear occurrence of z as EzF. It has at least one nonempty context. Each closure step uses the recorded actual premise words, not merely their endpoints. Its scheme therefore produces a literal witness for the new pair/head/final. The process terminates because there are at most |D|\ensuremath{{}^{2}} pairs and 2|D| endpoint flags.

At the fixed point a pair(x, y), x, y in D, is witnessed exactly when \(Rx=Ly\); heads and finals in D are witnessed exactly when \(Lx=H\) and \(Rx=T\). For completeness of closure, use the bipartite identification graph with row vertices Rx and an H-row, and column vertices Ly and a T-column. R completes intersecting rows; S and F complete and mark rows/columns meeting the H-row; T and E do the same at the T-column. Consequently each connected component is complete on its admissible row/column incidences. Conversely every rule preserves the original identifications, proving equality rather than a one-sided approximation. Endpoints for a literal-only variable outside D are leaves in this identification graph, incident only with H or T; removing them cannot alter any identification among D-endpoints, H, T.

Every actual H-to-T word-walk with labels in D is now derivable. Start with a witnessed word having its first head. For each desired adjacent pair use M to splice the recorded pair witness after the current letter. Finally splice a recorded final witness to remove the excess suffix. This also works for a literal walk in D. All empty external blocks are handled by the corresponding finite optional-context identity.

\par\smallskip\noindent\textit{Safe cycle lifting at both ends.}\quad 
If D is empty, every valid target is an original literal and no work is needed. Assume D nonempty. Delete from the directed endpoint graph all edges whose labels are outside D. Each such edge belongs to an original literal and runs directly \(H\to T\). Any original nonlinear word supplies an H-to-T path using D labels. Replacing each deleted edge by this path proves that reachability, and hence strong components, among all retained endpoints are unchanged. Thus every connector can be chosen to use D exclusively. This replaces11619's special initial component argument; no forbidden literal-only variable is inserted.

Take a valid target w not already original. Every one of its variables is in D: this is the criterion for a nonlinear target, and a literal outside D would have been original. At its initial gap from H to L(head(w)), at each internal gap, and at its final gap to T, the two endpoints are strongly connected. If they differ, choose nonempty paths u forward and v back using D labels, and insert the connector uvu. Equal endpoints need no insertion. The resulting W is an actual H-to-T walk in D, so the preceding closure derives \(W\leq P\).

It remains to derive \(w\leq W\). At the initial gap, use I0 twice to get \(w\leq v\)\ensuremath{{}^{2}}\(w\leq u\)\ensuremath{{}^{2}}v\ensuremath{{}^{2}}w, then B in reverse to obtain uvu w. At an internal gap use I2 twice to insert \(u^{2}\) and \(v^{2}\), and B in reverse to replace u\ensuremath{{}^{2}}\(v^{2}\) by uvu. At the final gap use I1 in the same way. Every u, v and every fixed-law identity variable is nonempty; possibly empty external contexts are ordinary term contexts. Apply the gap operations successively. This proves \(w\leq W\leq P\) without any bound on paths, variables or words.

Each summand of a valid polynomial identity is absorbed by the opposite polynomial. Distributivity and these unrestricted absorption derivations therefore prove every identity. The 244 templates form a complete finite basis.

\subsection{Finite bases from matching relations on word endpoints}\label{case:323}
Types 14676 and 15244 satisfy the following three identities in the two-operation semiring signature. As usual u \ensuremath{\leq} P means \(u+P=P\).

{\small\par\smallskip
\CaseLine{\(xyzt = xyt + xzt\)}
\CaseLine{\(xz \leq  xyz\)}
\CaseLine{\(xt \leq  xy + zy + zt\)}
\smallskip}

The first reduces every word to length at most three by replacing four consecutive letters in a multiplication context. The second extracts the two endpoints of each reduced nonlinear upper word as an absorbed quadratic. Use a bipartite graph with separate left and right copies of each variable. For every nonlinear upper word add the edge from its first left copy to its last right copy. Write G for this graph, including only incident vertices. Within each connected component every edge between a left and a right vertex is derivable. Along an odd path, the last three alternating edges give a shortcut by the third identity. Repeating this argument derives the required edge. All the original words are kept. Let P* denote P together with these derived quadratics. P and P* are equationally equal.

For each proper principal additive ideal I below b, let R consist of the pairs (x, y) whose product belongs to I. For both algebras R is the graph of a bijection pi from a subset D onto a subset E of the five-element universe. The intersection of D and E has at most one element. There is a subset K such that a triple product xyz belongs to I exactly if (x, z) belongs to R and y belongs to K. The following table gives all nontrivial data. For 14676 the middle set is \{1\} throughout. For 15244 it is \{3\}.

\noindent Type\enspace\textbar\enspace b\enspace\textbar\enspace R\par
\noindent 14676\enspace\textbar\enspace 0\enspace\textbar\enspace (1,\allowbreak{}1), (3,\allowbreak{}2)\par
\noindent 14676\enspace\textbar\enspace 1\enspace\textbar\enspace (1,\allowbreak{}1)\par
\noindent 14676\enspace\textbar\enspace 2\enspace\textbar\enspace (1,\allowbreak{}2)\par
\noindent 14676\enspace\textbar\enspace 3\enspace\textbar\enspace (3,\allowbreak{}1)\par
\noindent 15244\enspace\textbar\enspace 0\enspace\textbar\enspace (1,\allowbreak{}3), (3,\allowbreak{}2)\par
\noindent 15244\enspace\textbar\enspace 1\enspace\textbar\enspace (1,\allowbreak{}3)\par
\noindent 15244\enspace\textbar\enspace 2\enspace\textbar\enspace (3,\allowbreak{}2)\par
\noindent 15244\enspace\textbar\enspace 3\enspace\textbar\enspace (3,\allowbreak{}3)\par
Consider any valid reduced absorption w \ensuremath{\leq} P. Fix a proper ideal and pin the at most three variables of w to any assignment on which w is outside the ideal. Requiring every upper word in the ideal together with these pins is impossible. We extract at most 24 original or derived upper words already witnessing this impossibility, irrespective of the number of variables or upper summands.

A literal imposes I on its variable. A middle variable of a cubic imposes K. Include all singleton pins as unary restrictions. A variable used at both ends, whether in one edge or different edges, must take a value in D intersect E. Impose that restriction on both copies. If the intersection is empty and such a variable occurs, its two incident upper edges already give a contradiction. Otherwise this cross-role value is unique. A variable with no endpoints is subject only to its unary restrictions. If their intersection is empty, for each of the at most five possible scalar values retain one failing unary restriction, including its original word if it comes from a word.

Within each nonempty component of G, a possible color a in D gives a to every left vertex and pi(a) to every right vertex. This follows from functionality and injectivity moving in alternate directions. The choice is admissible precisely when all unary restrictions on the vertices hold. If every component admits a choice and each isolated variable admits a value, the choices combine: variables appearing on both shores have the unique intersection value on both copies. Thus in the remaining impossible case some single component admits no color.

Choose an edge in that component as a reference. For each a in D choose a vertex and one unary restriction excluding the value proposed by that color. For a unary literal or middle restriction retain its original word; for a pin keep the pin. For a cross-role restriction retain an incident original or extracted edge of each role, so the actual variable must belong to the intersection. Join each selected left vertex directly to the reference right vertex and each selected right vertex to the reference left vertex, using the derived component edges. The reference edge and these connections enforce exactly the same reference color at every selected test. At most one reference edge, five connecting edges and ten words for unary tests are required. If D is empty, the single reference edge suffices. The uniform bound 24 covers all cases.

Over at most \(5^{4}\) choices of the bound and the target assignment keep the union of the selected upper words in P*. If needed, add one arbitrary upper word. For any assignment on which w is outside the specified ideal, the selected words cannot all be inside it. Taking the bound equal to their sum proves w is absorbed by their subpolynomial. Consequently the three displayed laws, the semiring axioms and all valid absorptions of length at most three with at most \(B=24*5^{4}\) upper terms on 3(B+1) variables form a finite basis. Obtain each reduced summand in both directions of an identity, reinsert any omitted right words and undo the reductions. This argument accommodates repeated letters, arbitrary polynomials, disconnected graphs and alternating paths of unrestricted length.

\subsection{Seven-identity bases for two orthogonal group extensions}\label{case:324}
We use the signature \((+,\cdot)\) without constants and work relative to the axioms of additively idempotent semirings. Write \(p\le q\) for \(p+q\approx q\). All statements and derivations below allow arbitrarily many variables. No minimality assertion is made about the bases.

\par\smallskip\noindent\textit{Two algebras and their bases.}\quad 
Let \(G=\{e, g\}\) be the group of order two with identity \(e\). Adjoin an idempotent \(f\) and two elements \(c, o\). Multiplication agrees with the group operation on \(G\), satisfies \(f^2=f\), and gives \(o\) on every other pair. Thus \(c\) and \(o\) annihilate every element to \(o\), and products between \(f\) and \(G\) also equal \(o\).

There are two additive orders. In \(B_0\), the only strict comparisons apart from those involving the top \(o\) are \(f<c\) and \(e<c\); in \(B_1\), replace \(e<c\) by \(g<c\). Addition is the join. In particular, \(f+e=c\) in \(B_0\), whereas \(f+g=c\) in \(B_1\). The maps from the displayed order \((e, g, f, c, o)\) to the catalog are \((1,\allowbreak{}3,\allowbreak{}2,\allowbreak{}0,\allowbreak{}4)\) and \((3,\allowbreak{}2,\allowbreak{}1,\allowbreak{}0,\allowbreak{}4)\), respectively. These give types 14705 and 14766.

Theorem 1. Both algebras are finitely based. Six common laws are

\CaseDisplay{\begin{array}{rll}
\mathrm{(C)}&xy\approx yx,\\
\mathrm{(N)}&x^3y\approx xy,\\
\mathrm{(E)}&x\le xy^2,\\
\mathrm{(A)}&x^3\le x+xy,\\
\mathrm{(M)}&xyz^2\le xy+xz,\\
\mathrm{(T)}&xyzt\le xy+xz+xt.
\end{array}}

A seventh law for \(B_0\) is

\CaseDisplay{\mathrm{(Q)}\qquad x^2\le x^3+y^2.}

A seventh law for \(B_1\) is

\CaseDisplay{\mathrm{(S)}\qquad x^2y^2\approx x^2+y^2.}

Each resulting set of seven laws is a complete identity basis for the indicated algebra, relative to the semiring axioms.

Soundness. The multiplication is the orthogonal union of the two groups \(G\) and \(\{f\}\), with absorbing element \(o\) and an additional annihilating element \(c\), so it is associative. To check distributivity, the only join of distinct elements below \(o\) is \(f+a=c\), where \(a=e\) or \(g\). Multiplying this join by any element gives \(o\); the join of the two corresponding products also gives \(o\) because one of the cross-component products is \(o\). For a distinct pair with join \(o\), the only case in which both products avoid \(o\) is the pair \(e, g\) multiplied by a group element; the products are distinct group elements and their join is \(o\). Pairs involving \(c, o\), or comparable elements, satisfy the same equality directly. Thus both algebras are semirings.

A word of length at least two will be called stable. Its value is \(f\) if all of its variables are assigned \(f\), is its group product if all are assigned elements of \(G\), and is \(o\) otherwise. This verifies (C),(N). A square has value in \(\{e, f, o\}\), which gives (E). All elements other than \(c\) satisfy \(x^3=x\), while \(c^3=cy=o\); hence (A) holds.

For (M) and (T), if a right summand equals \(o\), the inequality is automatic. Otherwise their shared first factor forces all factors to lie in \(G\), or all to equal \(f\). In \(G\), (M) reduces to \(xy\le xy+xz\). The left side of (T) is the product of the three right group values: if those values coincide their product is that value, and otherwise their join is \(o\). This proves both laws. In \(B_0\), the only nontrivial case of (Q) is \(x=g\); both \(g+e\) and \(g+f\) are \(o\), and \(y^2\in\{e, f, o\}\). In \(B_1\), the join and product of two square values coincide, proving (S).

\par\smallskip\noindent\textit{Affine closure inside a connected component.}\quad 
For a word \(w\), let \(c(w)\) be its content and \(\pi(w)\) its vector of exponent parities in \(\mathbb F_2^X\). Write a polynomial as

\CaseDisplay{P=\sum_{x\in L}x+\sum_{r\in R}r,}

where \(R\) consists of stable words. Put \(N=\bigcup_{r\in R}c(r)\). Partition \(N\) into the connected components of the hypergraph whose edges are the sets \(c(r)\), for \(r\in R\). Equivalently, two variables are in the same component if there is a chain of intersecting stable contents joining them. Linear summands do not create edges.

For a component \(D\), put

\CaseDisplay{V_D=\{\pi(r):r\in R,\ c(r)\subseteq D\}
       \cup\{\epsilon_x:x\in L\cap D\},\qquad
K_D=\operatorname{aff}_{\mathbb F_2}(V_D).}

Every \(V_D\) is nonempty. Call \(D\) marked when \(0\in K_D\). For a marked component, \(K_D=\operatorname{span}(V_D)\).

Lemma 2. The six common laws derive \(w\le P\) whenever \(c(w)\subseteq D\) and \(\pi(w)\in K_D\) for a component \(D\). This includes linear target words.

Proof. Laws (C),(N) equate stable words having the same content and parity. Indeed, an exponent can be reduced by two using (N) whenever a nonempty remaining factor is present. For singleton content the representatives are \(x^2\) and \(x^3\); for larger content, each exponent can be reduced to one or two. In particular, every stable \(r\) satisfies \(r^3\approx r\).

Law (E) permits insertion of a square of any variable, and gives \(x\le x^3\) on setting \(y=x\). If \(x\in L\cap N\), choose a stable summand \(r=xq\), with \(q\) nonempty after commutation. Substituting \(q\) for \(y\) in (A) gives \(x^3\le x+r\le P\). Add all such cubes to the available stable words. Their parity vectors are precisely the extra unit-vector generators in \(V_D\).

If stable words \(u, v\) share a variable \(x\), write \(u=xp\), \(v=xq\), with both \(p, q\) nonempty. Laws (N),(M) give

\CaseDisplay{uv^2\approx xpq^2\le xp+xq=u+v.\quad\textup{(\ensuremath{6})}}

If stable \(u, v, w\) share \(x\), writing their nonempty remaining factors as \(p, q, r\) gives, by (N),(T),

\CaseDisplay{uvw\approx xpqr\le xp+xq+xr=u+v+w.\quad\textup{(\ensuremath{7})}}

There are no empty-word substitutions in these steps.

Starting with any available stable word in \(D\), use (6) along a chain of intersecting hyperedges. This produces an absorbed word of content exactly \(D\) with the same parity. Thus every generator of \(K_D\) has a representative of full content \(D\) below \(P\). All these words share any chosen variable of \(D\), so repeated applications of (7) produce every odd product of them below \(P\). The parities of odd products are exactly \(K_D\).

For a stable target \(w\) with content contained in \(D\), use (E) to insert the missing variables twice. The resulting word has content \(D\) and the same parity, so normalization identifies it with an absorbed word constructed above. If \(w=x\) is linear, first use \(x\le x^3\). This proves the lemma. \(\square\)

\par\smallskip\noindent\textit{Exact absorption criteria.}\quad 
Theorem 3. Apart from the tautological case \(w=x\in L\), the valid absorptions are as follows.

1. In \(B_0\), the target content must lie in one component \(D\). If there are no marked components, the further requirement is \(\pi(w)\in K_D\). If at least one component is marked, replace this by \(\pi(w)\in\operatorname{span}(V_D)\), for every component.

2. In \(B_1\), a target either lies in a single unmarked component \(D\) with \(\pi(w)\in K_D\), or has content in the union \(M\) of the marked components and parity in \(\bigoplus_{D\text{ marked}}K_D\).

If \(R\) is empty, only the tautological case is possible. In both algebras these criteria are derivable from the respective seven laws.

Necessity of the content restrictions. An unlisted linear variable absent from \(N\) is separated by assigning it \(o\) and all right variables \(f\). For a stable target containing a variable outside \(N\), give that variable \(c\) and all other variables \(f\). The target is \(o\), while the right side is at most \(c\). This also handles \(R=\varnothing\) for stable targets. Hence every remaining case has target content contained in \(N\).

In \(B_0\), if the target meets two components, assign one such component constantly \(f\) and all other variables \(e\). Every right word lies below \(c=f+e\), but the stable target mixes the two multiplicative components and equals \(o\). Thus it must be supported in a single hypergraph component.

In \(B_1\), suppose the target meets an unmarked component \(D\) and also a different component. Because \(0\notin K_D\), there is a linear functional taking value one on every vector in \(V_D\). Assign the variables in \(D\) their corresponding group values, and every other variable \(f\). Right stable words in \(D\) and its right linear variables have value \(g\); all other right words have value \(f\). Thus \(P\le c=f+g\), whereas the target equals \(o\). The only remaining possibilities are a single unmarked component or the union of marked components, as stated.

Necessity of parity restrictions in \(B_0\). Suppose the target is supported in \(D\). If its parity lies outside \(\operatorname{span}(V_D)\), choose a linear functional annihilating \(V_D\) and taking value one at \(\pi(w)\). Assign group values in \(D\) according to this functional and assign \(f\) to every other variable. All right words in \(D\) have value \(e\), all other right words have value \(f\), and \(P\le c\). The target is \(g\), which is not below \(c\). This proves the linear-span requirement.

Suppose no component is marked and \(\pi(w)\notin K_D\). In each component choose a functional taking value one on \(V_D\); in the target component choose it to take value zero on \(\pi(w)\). The latter choice exists by affine separation: the affine hull excludes both zero and the target vector. Assign group values accordingly and assign \(g\) to variables occurring only linearly. Every right word is \(g\) but the target is \(e\). This proves the stronger affine requirement when no component is marked.

Necessity of parity restrictions in \(B_1\). In an unmarked component, a target vector outside \(K_D\) can again be separated by a functional taking value one on \(V_D\) and zero at the target. Assign every other variable \(f\). Then \(P\le c=f+g\), while the target is \(e\), which is not below \(c\).

For a target supported in \(M\), failure of the direct-sum condition means its restriction to some marked \(D\) lies outside the linear space \(K_D\). Choose a functional zero on \(V_D\) and one on that restriction; use the resulting group assignment on \(D\), and give every other variable value \(e\). Now every right word is \(e\), while the target is \(g\). This proves necessity in all cases, including unlisted linear targets, whose variables were not forced by a right linear summand.

Derivability in \(B_0\). When no component is marked, apply Lemma 2. Otherwise the lemma supplies \(y^2\le P\) for a variable \(y\) in a marked component. For each available stable seed \(r\) (an original stable word or an anchored cube), (Q) yields

\CaseDisplay{r^2\le r^3+y^2\approx r+y^2\le P.}

Adjoin these squares to the seeds. Their supports do not change the components, and their zero parities enlarge the affine hull in every component to \(\operatorname{span}(V_D)\). Applying the proof of Lemma 2 to the enlarged family derives exactly the asserted targets.

Derivability in \(B_1\). Lemma 2 handles each unmarked component. For every variable \(x\) in a marked component it gives \(x^2\le P\). Let \(h\) list all variables in the nonempty set \(M\) once. Repeated use of (S) gives

\CaseDisplay{h^2\approx\sum_{x\in M}x^2\le P.}

Every stable seed in a marked component shares a variable with \(h^2\). Equation (6) therefore pads that seed to content \(M\) without changing its parity, while retaining absorption by \(P\). Include \(h^2\) itself among the full-content words. Equation (7) then produces the affine hull of their parities. Because zero is included, this hull is the span of all marked generators, namely \(\bigoplus_{D\text{ marked}}K_D\). Square insertion and normalization, as in Lemma 2, give every target in the second case of the criterion. This proves Theorem 3. \(\square\)

\par\smallskip\noindent\textit{Completeness and scope.}\quad 
For any valid identity \(p\approx q\), distribute to polynomials \(P, Q\). Each summand of \(P\) is absorbed by \(Q\), and Theorem 3 derives that absorption from the stated seven laws. Thus \(P\le Q\) is derivable, and reversing the argument gives \(Q\le P\). This proves Theorem 1 for unrestricted variable sets.

The distinction between the two algebras is precise. In \(B_0\), a zero parity in one component permits linear-span closure in every component but never joins their supports. In \(B_1\), only components whose affine hulls already contain zero can be joined; unmarked components retain their separate affine conditions.

For each transfer in the register, applying Proposition~\ref{prop:lower-fibres} in both directions to the indicated tables identifies their equational theories.

The source manuscript for this case credits Aifa Wang and Lili Wang, School of Mathematical Sciences, Chongqing University of Technology.

The source manuscript for this case acknowledges support from the National Natural Science Foundation of China (12371024), the Science and Technology Research Program of Chongqing Municipal Education Commission (KJZD-K2024011102) and the Chongqing Natural Science Foundation Innovation and Development Joint Fund (Municipal Education Commission) (CSTB2025NSCQ-LZX0067).

\subsection{Five complete identities for 14781}\label{case:325}
Work in the signature (+, multiplication), without constants. Write \(u\leq P\) for \(u+P=P\). Relative to the additively idempotent semiring axioms, catalogue 14781 has the following complete finite basis:

{\small\par\smallskip
\CaseLine{\(xy=yx; xxx=xx; x\leq xy;\)}
\CaseLine{\(xz\leq xx+xy+yz; xyz\leq xy+xz+yz.\)}
\smallskip}

The proof is unrestricted in the number of variables, word lengths, and polynomial sizes. The variables in a displayed identity always receive nonempty terms.

\par\smallskip\noindent\textit{Word values.}\quad 
Addition:00004/\allowbreak{}01004/\allowbreak{}00224/\allowbreak{}00234/\allowbreak{}44444. Multiplication:44444/\allowbreak{}41444/\allowbreak{}44424/\allowbreak{}44234/\allowbreak{}44444. The additive order has \(3<2<0<4\) and \(1<0\), with 1 incomparable to 2 and 3. For a nonlinear word, its value is 1 if every occurrence is 1. If all occurrences lie in\{2,\allowbreak{}3\}, the value is 3 with no 2,\allowbreak{}2 with exactly one 2, and 4 with at least two2 occurrences. In every other case the value is 4. Literals retain their values. This follows inductively because 3 is the identity on\{2,\allowbreak{}3,\allowbreak{}4\},\(2*2=4\), and mixing 1 with any other element produces4. Multiplication involving0 or 4 is 4.

Commutativity and \(xxx=xx\) give normal forms determined by content and multiplicities truncated to one or two. In particular a repeated target variable cannot be treated as a simple one. To normalize, sort letters using commutativity and replace triples by doubles until no multiplicity exceeds two. These are context instances of the two structural laws.

\par\smallskip\noindent\textit{Exact criterion.}\quad 
Every item below is defined from original summands of P. Let C be their content union. Let G have the variables occurring in nonlinear original words as vertices, with an edge between distinct variables that occur together in one such word. A variable with only pure nonlinear words is still a vertex. Let F consist of variables repeated in some original word. Call these the anchored variables.

A literal x is absorbed precisely when x is in C. A nonlinear target w is absorbed precisely when all four requirements hold:

1. Every target variable is a vertex of G.

2. All target variables belong to one connected component of G.

3. Every variable repeated in w belongs to F.

4. The target variables outside F form a clique in the original graph G.

Singleton cliques are allowed. A pure nonlinear target requires its variable in F; it will be derived by extracting an actual original repeated occurrence, not by treating a reflexive path as an edge.

For sufficiency, consider an assignment for which the sum P is not 4. Each nonlinear component is either entirely 1 or entirely in\{2,\allowbreak{}3\}. In the second case, each anchored vertex is 3 and the vertices assigned 2 form an independent set in G: an original word with two2 occurrences would already give 4. Hence the four requirements make the target 1,\allowbreak{}3, or 2. In the first case an original nonlinear word in the same component provides 1. In the second case an original nonlinear word provides 3 or 2, dominating a target of 3. If the target is 2, its unique 2-valued occurrence lies in an original word that has value 2, so this also dominates it. A literal is dominated by an original containing word, as follows directly from the word formula (equivalently x\ensuremath{\leq}xy).

Necessity has explicit separators for all cases. Assign every unmentioned variable 1. A missing literal outside C can be assigned 4. If a nonlinear target has a variable outside the vertex set of G, assign that variable 0; only its possible literal summand changes on the right, so the sum is at most 0 while the target is 4. If target variables occupy distinct components, assign one entire component 3 and leave the others 1; again the target is 4 and the right sum at most 0.

Now suppose the target lies in one component. If it repeats an unanchored variable x, assign \(x=2\) and all other vertices of that component 3. No original word repeats x, so its right words have values 2 or 3, while the target is 4. Finally, if distinct unanchored target variables x, y have no original edge, assign them2 and the other component vertices 3. No original word contains both or repeats either. Every right word is therefore at most 0 after all joins, while the target is 4. These cases exhaust the failures for arbitrary finite polynomials.

\par\smallskip\noindent\textit{Deriving all absorbed words.}\quad 
Keep every original summand throughout. From an original word u extract any nonempty proper subword v, including repeated letters with their original multiplicity, by \(x\leq xy\) with \(x=v\) and y the nonempty remaining factor; commutativity rearranges the factors. Therefore every original edge xy and every original repeated-letter square xx for x in F is available below P. Literals follow in the same way.

If x is anchored and y is connected to x, choose a positive path \(x=p0\), p1,..., \(pn=y\) of original edges when x differs from y. Retain xx and extend the first original edge successively by

{\small\par\smallskip
\CaseLine{\(x p_{i+1} \leq  xx + x p_{i} + p_{i} p_{i+1}.\)}
\smallskip}

This is the fourth basis law. The case \(x=y\) is the actual extracted square. Thus a pair of target occurrences is available whenever at least one endpoint is anchored. If neither endpoint is anchored, the clique requirement gives its original edge. Repeated target occurrences always have anchored endpoints, so this constructs every pair of target positions, retaining their multiplicity.

For completeness of the final merge, the fifth law implies the following claim: any word \(a_{1}\)...\(a_{m}\) of length \(m\geq 2\) is below the sum of its pairs \(a_{i}\) \(a_{j}\) \((i<j)\). Prove it by induction on m. It is immediate for \(m=2\). For \(m\geq 3\) put \(X=a_{1}\cdots a_{m-2}\), \(y=a_{m-1}\), \(z=a_{m}\). Apply \(XYZ\leq XY+XZ+YZ\). The first two summands have m-1 positions and each of their pairs is one of the original m-position pairs; YZ is itself such a pair. The induction is valid with repeated letters and with X a nonempty word. All target pairs have already been derived below P, so the target follows. Equivalently an explicit construction uses target prefixes plus one extra letter and only quadratically many triangle steps. Structural normalization never deletes the distinction between a simple variable and a repeated variable.

Each of the five identities is valid by the exact criterion or the word formula. The derivation handles every valid absorption, and distributive expansion plus absorption in both directions handles every identity. This proves finite basability without a variable bound.

Essential counterexamples: \(xx\leq x\) fails at \(x=2\); \(xz\leq xy+yz\) fails at \(x=z=2\), \(y=3\); \(xy\leq xx+yy\) fails at \(x=1\), \(y=3\). Connectivity alone does not replace the clique requirement at unanchored vertices.

\subsection{A finite basis for 14709 by connected endpoint witnesses}\label{case:326}
Work in the constant-free ai-semiring signature. In 14709 let (e, f, a, c, z) denote labels (1,\allowbreak{}2,\allowbreak{}3,\allowbreak{}0,\allowbreak{}4). The nonzero ordinary products are \(ee=e\), \(ea=a\), \(af=a\), \(ff=f\). Any product with c or z gives z. Addition has e, \(f<c<z\) with e, f incomparable, and \(a<z\) incomparable to e, f, c. Thus a nonlinear word is e or f exactly on that constant loop, is a exactly for one a between a prefix of e's and suffix of f's, and otherwise is z. A literal may also be c.

\par\smallskip\noindent\textit{Criterion.}\quad 
For an upper polynomial P write C for its content union and N for the union of its nonlinear-word supports. Join variables appearing in the same original word into one component, taking transitive closure of this undirected relation on C. Singletons are allowed. Construct the ordinary endpoint graph by identifying L of each head with H, R of each final with T, and \(R_{x}\) with \(L_{y}\) at adjacent xy of every original word. Give each variable the edge \(L_{x}\) to \(R_{x}\) and collapse directed strongly connected components. For w \ensuremath{\leq} P one needs c(w) inside C and the quotient walk from H to T labelled w. For a nonlinear w one additionally needs c(w) inside N and all its letters in the same undirected component.

Indeed the one-loop and all-bridge bounds give the coverage and walk conditions exactly as in the Boolean endpoint proof for 14674. They use no c-valued variable. At bound c every nonlinear original must be constant e or constant f, forcing each variable of N to be a loop and the choice constant on each undirected component. A literal-only original can take c. To separate a missing nonlinear occurrence give its variable c, the rest e. To separate a target crossing components give one of the components e and the other f, using constant loop assignments on the remaining components. The upper value is at most c, the target z. These cases also show sufficiency at c, including literals. The top bound is automatic.

\par\smallskip\noindent\textit{Presentation with at most 183 laws.}\quad 
Use the four fixed laws and the M, R, E, F optional-context schemes of the 14674 proof. Their exact patterns are

{\small\par\smallskip
\CaseLine{\(xyx=xxyy,  x\leq yyx,  x\leq xyy,  xy\leq xzzy\)}
\CaseLine{\(M\qquad A x D \leq  A x B + C x D\)}
\CaseLine{\(R\qquad A x t F \leq  A x y B + C z y D + E z t F\)}
\CaseLine{\(E\qquad A x \leq  A x y B + C z y D + E z\)}
\CaseLine{\(F\qquad y B \leq  A x y B + C x z D + z E\)}
\smallskip}

Some capital symbols have the same spelling as a rule name. Each letter in a word is an identity variable. Capitals may independently be omitted everywhere in the same displayed scheme. There are 16+64+32+32 such instances, plus the four fixed laws. Add the following two schemes (32 instances), where \(W=A\) t B, \(V=C\) t D and the lower-case variables s, t stay nonempty:

{\small\par\smallskip
\CaseLine{\(s V \leq  s W + W + V\)}
\CaseLine{\(V s \leq  W s + W + V\)}
\smallskip}

Finally add the three instances of yy \ensuremath{\leq} xx + A y B in which at least one of A, B is retained. All omissions specify separate identities; ordinary substitutions are always nonempty.

For the four fixed laws and M, R, E, F the ordinary endpoint argument gives soundness. For each nonlinear left side all its letters are in one nonlinear right component: in R the first two upper words share y and the last two share z, and in M the words share x. In E or F the left support is in the first nonlinear upper word. The fixed laws retain the connected support within a common right word; \(xyx=xxyy\) also follows by evaluating a repeated variable in a possible one-bridge path. At a bridge bound the values of sW and W being a force s to be the initial loop, so sV is a. The suffix case forces the final loop. At bound c the nonlinear sW (or Ws) is a single-loop word, and the shared letter t forces V to have this same loop, including a literal V. For the last scheme a square xx is never below a. At c the original nonlinear A y B forces y to be a loop. All loop bounds follow directly from the word values.

\par\smallskip\noindent\textit{Proof within one component.}\quad 
Fix a component J containing a nonlinear word and take the subfamily \(P_{J}\) of original words in it. Every letter of J occurs in a nonlinear word. For this subfamily the absorption criterion coincides with that of 14674 for words over J. We simulate its *constructive* proof, whose finite rule set and word-substitution derivations are in Section~\ref{case:322} and in the guarded S, T specifications of its immutable verifier. Its M, R, E, F and fixed-law applications already belong to our list. Only S and T require a derivation. Every available word from the construction has nonempty support in J.

For an S step take the actual available words A x y B, \(W=y\) C, and \(V=z\) D. Put \(s=A\) x. Splice the first word and W at y by M to derive s W. Join W to V by a finite sequence of intersecting words starting with W, ending with V and with every intermediate word an original member of \(P_{J}\). Such a sequence exists by the definition of J: use any letter from each endpoint word and a chain of original words meeting along a letter path in J. On each link use the prefix propagation law at one common letter to obtain s followed by the next word. All the undisturbed words on the right are themselves available. At the end s V is the required A x z D. For T start with the available words A x y B, \(W=C\) x and \(V=D\) z. Splice W with the first word at x to obtain W s, with \(s=y\) B. Use suffix propagation on a sequence from W to V to obtain V s. The extra nonlinear occurrence of z in the original guarded S, T applications is harmless. When \(W=V\) the initial splice is enough. This construction uses an arbitrary number of steps from only 32 propagation identities, irrespective of the chain length. All the walk splices and the insertions for strongly connected components in the parent proof can therefore be carried out locally.

\par\smallskip\noindent\textit{Communication between components.}\quad 
In the full endpoint graph the separate components' endpoint graphs are glued only at their common terminals H and T. This follows before contracting identifications by using two copies (left and right) of each variable and undirected arcs for the identification equations. A vertex of one component can enter another only through H or T. Every component has a directed H-to-T path, labelled by any original word in it. If there is no directed return from T to H (allowing the zero-length return when they are already identified), all excursions through another component between H and T can be replaced inside J by an H-to-T path or omitted when their endpoints are the same. Hence the global strong-component comparisons within J are exactly the local ones. Conversely, if there is a T-to-H return globally, one component K already supplies one: on a returning walk take the segment after its last T up to its first subsequent H. For this argument use the uncontracted representation with identification equations traversable in either direction, so it also covers globally identified terminals. A component of literals only has no return. In any returning component all its endpoints are in a single strong component because every vertex lies between H and T.

For a target supported in one nontrivial component, if its local walk test holds use the preceding derivation directly. Otherwise there must be a globally returning component K. Its local proof gives xx \ensuremath{\leq} \(P_{K}\) for any x in K. Select a nonlinear original A y B in the target component J. Use the last scheme to derive yy \ensuremath{\leq} P. Adjoin yy to \(P_{J}\). It identifies H and T locally (the intermediate adjacency in yy identifies the two endpoints of y), and it preserves the undirected component J. The local proof now derives the target, and transitivity gives the original bound P. This works also for a literal in J. A literal outside N that passes coverage is already a summand of P. These are all cases. Distributing arbitrary terms and deriving the two sets of summand absorptions proves every identity, with no bound on original variables, walks, or polynomial sizes.

\subsection{Parity components coupled by cubic anchors}\label{case:327}
We use the constant-free signature \((+,\cdot)\) relative to the additively idempotent semiring axioms. Write \(u\le v\) for \(u+v\approx v\).

\par\smallskip\noindent\textit{1. Algebras and finite bases.}\quad 
On \(\{a, b, c, e, o\}\) set \(a^2=b\), \(ab=ba=c\), \(e^2=e\) and every other product equal to \(o\). For \(A\), addition is join in the order \(a, e<c<o\), \(a\parallel e\), and \(b<o\) with \(b\parallel a, e, c\). For \(B\), retain only \(e<c<o\) and \(a, b<o\), with \(a, b\) incomparable with \(e, c\) and with one another. In the displayed element order the catalog map is \((2,\allowbreak{}3,\allowbreak{}0,\allowbreak{}1,\allowbreak{}4)\), giving 14764 and 14765. These element names are not constants in terms.

The multiplication consists of the quartic-nilpotent component \(\{a, b, c, o\}\) and the idempotent \(e\), with mixed products \(o\). This proves associativity. Translation by \(a\) takes \(a, b\) to \(b, c\) and every other input to \(o\); translation by \(b\) takes only \(a\) to \(c\); translation by \(e\) takes only \(e\) to \(e\). All these maps preserve the joins specified above, including \(a+e=c\) in \(A\) and \(e+c=c\) in both algebras. Thus both are semirings.

Let \(\Lambda\) consist of the following sixteen laws:

\CaseDisplay{\begin{array}{rlrl}
(C)&xy\approx yx,&(S)&x^5\approx x^4,\\
(Q)&xyzt\approx x^4+y^4+z^4+t^4,
 &(T)&xyz\approx x^2y+x^2z,\\
(U)&x\le x^4,&(P)&xt\le xy+yz+zt,\\
(E)&xy\le x^2+y^2,&(L)&y^2\le x^2+xy,\\
(M)&z^3\le x^3+xy+yz,&(J)&xz\le xy+y^2z,\\
(K)&x^2z\le x^2y+y^2z,&(N)&x^2z\le x^3+xy+yz,\\
(Z_1)&x^4\le x^2+y^3,&(Z_2)&x^4\le xy+x^3+y^3,\\
(R)&y^4\le x^4+xy,&(H)&y^4\le x^4+x^2y.
\end{array}}

Theorem. A basis for \(A\) is \(\Lambda\) and

\CaseDisplay{(A_1)\quad x\le x^3,\qquad (A_2)\quad x^3\le x+xy.}

A basis for \(B\) is \(\Lambda\) and

\CaseDisplay{(B_1)\quad x^4+y^4\le x+xy.}

These are eighteen and seventeen laws, respectively, without a minimality assertion.

For validity, a non-top quadratic value is \(b\) at \((a, a)\), \(c\) at \((a, b)\) or \((b, a)\), or \(e\) at \((e, e)\). A cubic is \(c\) only when all factors are \(a\), and \(e\) only when all are \(e\). A product with at least four factors is \(e\) only when all factors are \(e\), and is otherwise \(o\). These descriptions verify (C),(S),(Q),(T),(U). For (P), a right value \(b\) forces all its vertices \(a\), and a right value \(c\) forces the path either to alternate \(a, b\) or to be all \(e\); this verifies its endpoints. The same argument on one or two edges verifies (E),(L),(M),(J),(K),(N). A right side in \((Z_1)\) or \((Z_2)\) below \(o\) forces the vertex whose fourth power occurs on the left to be \(e\); \((R),(H)\) force the other endpoint to be \(e\) as well. Finally, \(a<c\) verifies \((A_1),(A_2)\), whereas the absence of \(a<c\) in \(B\) means that \(x+xy<o\) forces \(x=y=e\), verifying \((B_1)\).

\par\smallskip\noindent\textit{2. Normal forms.}\quad 
From (S), \(x^n=x^4\) for \(n\ge4\). Substituting \(u^2, u^2, v^2, v^2\) in (Q) gives \((uv)^4=u^4+v^4\). Thus any word of length at least four is the sum of the fourth powers of its content, by splitting it into four nonempty factors and applying (Q). Law (T) handles words of length three.

Two consequences of (T), obtained by repeating variables, are

\CaseDisplay{x^2y=xy^2,\qquad x^3\le x^2y.}

Indeed, substituting \((x, y, y)\) and \((x, x, y)\) into (T) gives \(xy^2=x^2y\) and \(x^2y=x^3+x^2y\), respectively. Together with commutativity this gives both claims. Every polynomial therefore normalizes to linear terms, quadratic edges \(xy\), symmetric cubic anchors \(h(x, y)=x^2y\), and marks \(q(x)=x^4\). Loops \(h(x, x)=x^3\) are included.

\par\smallskip\noindent\textit{3. Exact component criterion.}\quad 
Fix a normal right polynomial \(F\). Write \(L\) for its linear variables, \(E\) for its quadratic edges, \(D\) for its cubic anchors, and \(Q\) for its fourth-power vertices. Let \(V\) be its nonlinear vertices. Let \(G=(V, E)\) and let \(N=(V, E\cup D)\), treating the anchors as undirected edges when forming \(N\).

Pin both endpoints of every cubic anchor. In \(A\), additionally pin every vertex of \(L\) incident with a quadratic edge. A pinned vertex is required to have color A in a proper two-coloring of its \(G\)-component. A \(G\)-component is bad when it is nonbipartite or its pins lie in both color classes.

An \(N\)-component is forced E if it contains a bad \(G\)-component or a fourth-power mark. In \(B\), it is also forced E when it contains a vertex in \(L\cap V\). Write \(T\) for the union of forced-E components. On the other \(G\)-components, a vertex is forced A or B if pins fix its color accordingly; an unpinned component has no forced colors.

Finally put \(\delta=1\) if \(F\) has a cubic anchor, a fourth power, or a linear vertex incident with a quadratic edge; otherwise put \(\delta=0\).

Absorption lemma. The following list is exact.

\textbullet{} A linear \(x\) is absorbed if \(x\in L\), or if \(\delta=1\) and \(x\in T\). For \(A\) also allow \(\delta=1\) and \(x\) forced A.

\textbullet{} A quadratic \(xy\) is absorbed if both endpoints are in \(V\) and either both lie in \(T\), or belong to opposite sides of the same bipartite \(G\)-component, or have forced opposite colors and belong to the same \(N\)-component.

\textbullet{} An anchor \(h(x, y)\) is absorbed if \(\delta=1\) and either \(x, y\in T\), or \(x, y\) belong to the same \(N\)-component and both are forced A.

\textbullet{} A fourth power \(q(x)\) is absorbed if \(\delta=1\) and \(x\in T\).

In the second bullet a component forced E also has all its pairs covered by the first alternative. The color alternatives refer to components on which a permissible coloring exists.

Proof. A right value \(e\) forces all its occurring variables to be \(e\), so it imposes only occurrence conditions. A right value \(b\) is possible exactly when \(\delta=0\): all nonlinear vertices are then \(a\) and separate linear variables are \(b\). It rules out all cubic and fourth-power absorptions and every new linear term.

For a right value at most \(c\), an edge either joins \(a\) to \(b\), or has both endpoints \(e\). An anchor has both endpoints \(a\), or both \(e\). A fourth-power vertex must be \(e\). Hence each \(N\)-component is wholly of type E (all its vertices \(e\)), or wholly of type AB, with its \(G\)-components properly colored. An anchor forces both its endpoints A. A linear-edge incidence forces A or E in \(A\), and E in \(B\). Bad components can only be of type E. All remaining \(N\)-components can independently be AB or E.

A target quadratic can take a value at most \(c\) precisely when its endpoints are both forced E, or lie in the same \(N\)-component and are always opposite in its AB coloring. Those opposite-color conditions are exactly the two alternatives in the second bullet. A target anchor must have both endpoints E or both A; a target fourth power requires E. A missing linear term is below \(c\) exactly at an E vertex in \(B\), and at an A or E vertex in \(A\). These facts give the four bullets. In an evaluation with \(F\le c\), a graph-external target variable is separated by assigning it \(c\) if it is only linear on the right, or \(o\) if it does not occur. If a \(c\)-bounded evaluation is unavailable and \(\delta=0\), use the \(b\)-valued evaluation instead and assign isolated linear variables \(b\); a nonlinear target involving such a variable is then \(o\). An isolated right linear term places no additional constraint on the graph coloring. Purely linear right sides have no nonlinear absorption, as shown by assigning all their variables \(a\). Thus these assignments prove necessity and sufficiency, including disconnected components. \(\square\)

\par\smallskip\noindent\textit{4. Constructing the derivations.}\quad 
First, (P) shortens odd walks, so an odd circuit produces a square. Law (L) propagates this square along its \(G\)-component; (E) joins any two vertices whose squares have been obtained, including in different components. Law (M) propagates cubes along even paths, and (N) simultaneously creates an anchor between the starting pinned vertex and the endpoint. Pinned vertices on opposite sides are joined by an odd path, and \((Z_2)\) produces a fourth power. A square in an odd component together with any cube, even in another component, produces a fourth power by \((Z_1)\).

Pins coming from cubic anchors already have absorbed cubes by Section 2. In \(A\), \((A_2)\) supplies the cubes at the additional linear-edge pins. In \(B\), \((B_1)\) supplies fourth powers at a linear-edge incidence. For a linear-cubic incidence substitute the word representing the remaining two factors for \(y\) in \((B_1)\); this again supplies a fourth power. Original fourth powers require no derivation. Thus, whenever \(\delta=1\), every initial forced-E component obtains at least one fourth power. Laws (R),(H) propagate it along every edge and anchor in that \(N\)-component.

All vertices with fourth powers acquire squares by (U), followed by (S), and hence edges between any pair by (E). They acquire cubes by (U) at \(x^3\), followed by (S). Their mixed anchors follow by substituting \(h(x, y)\) into (U): \(h(x, y)\le h(x, y)^4=q(x)+q(y)\). Likewise (U) gives every required linear term at a fourth-power vertex. If \(\delta=0\), forced-E components arise only from odd circuits; the squares and (E) give their quadratic absorptions without adding fourth powers.

It remains to join forced color classes in an \(N\)-component. Within a pinned \(G\)-component, all A vertices connect by even paths from a pin. Law (N) connects them to that pin by anchors. Anchors already joining \(G\)-components can be composed by (K), \(h(x, z)\le h(x, y)+h(y, z)\). It follows, by walking through the finite component graph, that all forced-A vertices in that \(N\)-component acquire mutual anchors. For a forced-A vertex \(x\) and forced-B vertex \(z\), choose an A neighbor \(y\) of \(z\) along a path to a pin. Law (J), with variables renamed, gives

\CaseDisplay{xz\le yz+h(y, x),}

which derives the required opposite edge. Edges inside unpinned bipartite components follow from (P). Finally \((A_1)\) turns the derived cubes of forced-A vertices into the remaining linear absorptions for \(A\).

Every item in the lemma has now been derived using finite paths and nonempty substitutions. Normalize an arbitrary valid identity and apply the lemma separately to every summand in each direction. This proves both claimed bases with no variable bound.

\subsection{Thirteen complete identities for 14773}\label{case:328}
Work in the signature (+, multiplication), without constants. Write \(u\leq P\) for \(u+P=P\). In addition to the additively idempotent semiring axioms, a complete basis for catalogue 14773 consists of

{\small\par\smallskip
\CaseLine{\(xy=yx; xxy=xyy; xyzt=xyztt;\)}
\CaseLine{\(xxxx\leq xyzt; yyyy\leq xxxx+xy; xxxy\leq xxxx+yyyy;\)}
\CaseLine{\(xy\leq xyz; x\leq xy; xyz\leq xyzt; xx\leq xyz;\)}
\CaseLine{\(xz\leq xx+xy+yz; xxx\leq xx+yzt;\)}
\CaseLine{\(xyz\leq xxx+yyy+zzz+xy+yz.\)}
\smallskip}

All proofs below are unrestricted in variables, lengths and summand counts. Computation audits the proof rather than replacing it.

\par\smallskip\noindent\textit{Word values and normalization.}\quad 
Addition:00004/\allowbreak{}01004/\allowbreak{}00224/\allowbreak{}00234/\allowbreak{}44444. Multiplication:44444/\allowbreak{}41444/\allowbreak{}44404/\allowbreak{}44024/\allowbreak{}44444. The additive order has \(3<2<0<4\) and \(1<0\), with 1 incomparable to 2 and 3. A literal retains its value. A nonlinear word has value 1 exactly when every occurrence is 1. Otherwise a quadratic is 2 for(3,\allowbreak{}3),0 for(2,\allowbreak{}3) or(3,\allowbreak{}2), and 4 in every other case. A cubic is 0 if all occurrences are 3 and 4 otherwise. Every word of length at least four is 4 unless all its occurrences are 1. To see the all-length assertion, a product of three letters is 1,\allowbreak{}0 or 4; multiplying0 or 4 by anything gives 4, and a product starting with 1 remains 1 only while every next factor is 1.

Quadratics retain multiplicity. Cubics depend only on content; commutativity and \(xxy=xyy\) give the sorted normal form with any required repetitions of the smallest letter. Quartic and longer words also depend only on content and normalize separately, retaining length at least four. The normal form is sorted content when its size is at least four and otherwise is padded at the beginning with the smallest letter to length four. With \(xyzt=xyztt\), commute a desired letter to the end and duplicate it, or delete a redundant terminal copy while preserving length four. If the required smallest-letter multiplicity is initially absent, insert it first and then remove the surplus copies of other letters. Every application splits the remaining positions into three nonempty blocks.

\par\smallskip\noindent\textit{Exact absorption criterion.}\quad 
All following data come from original right summands P. Let C be their content union. Make a graph G on variables of nonlinear original words, joining each pair occurring in one such word. Vertices occurring only as repeated letters are retained. Let A be the union of connected components meeting the content of an original word of length at least four. Let F contain every variable of an original cubic and every variable having an original square xx. Let g say that P contains at least one original word of length at least three, anywhere in P.

An original normal word is absorbed immediately. Every other target w is absorbed exactly in the following cases:

\textbullet{} A literal x has x in C.

\textbullet{} A nonlinear word with content contained in A is absorbed.

\textbullet{} A quadratic xy is absorbed if x, y are in one G component and at least one endpoint belongs to F. This includes a square at a cubic variable, obtained directly rather than through a reflexive path.

\textbullet{} A cubic is absorbed if g holds, all its variables belong to F, and they belong to one G component.

There are no other absorbed targets. In particular, an original cubic activates only the cubic/square part of the argument; it does not act as a quartic witness. The global guard g may lie in a different component.

For sufficiency consider an assignment where P is not 4. Each original nonlinear word has value 1,\allowbreak{}2 or 0. Its component is either entirely 1 or entirely in\{2,\allowbreak{}3\}; original cubics force their letters to 3 in the latter case, and original squares do the same. Original quartics force all their letters to 1, propagating this to their whole component. Thus every active variable is 1; all anchored variables in an inactive component are 3 whenever that component is not entirely 1.

A target on active variables has value 1, already supplied by an original nonlinear word. For an anchored quadratic in a non1 component, if its other endpoint is 3 its value 2 is dominated by the original square or cubic anchor. If the other endpoint is 2, an original edge along the connecting path must cross from 3 to 2, giving 0 on the right. A target on two1 endpoints is again dominated by an original nonlinear term. For the cubic alternative, its variables are either all 1 or all 3. The latter target is 0. The original guard either has value 0, or has value 1; in the second case a right word in the target component gives 2 or 0, and \(1+2=0\). Thus the global guard is sufficient even outside that component. Literals are below any nonlinear word containing them, as the displayed word formula shows.

Explicit separating assignments prove necessity. Assign every unmentioned variable 1. A missing literal outside C is separated by assigning it 4. A variable outside G has no original nonlinear occurrence; assign it 0 to separate every nonlinear target containing it, while the right side stays at most 0. For a target of length at least four failing the active condition, choose one inactive component meeting it, give that entire component 3, and leave all other components1. The target becomes 4 and every right term is at most 0.

For a missing quadratic whose endpoints belong to different components, at least one component is inactive unless the active alternative holds. Give that component 3 and the other 1; the target is 4. If both endpoints are in one inactive component but neither is in F, give both endpoints 2 and every other vertex there3. No original square or cubic contains either endpoint. There is no original edge joining distinct endpoints: an original quadratic would be the target, and any longer original word containing them would make them anchored or active. Hence the right sum stays at most 0 while the target is 4. The same construction separates a missing unanchored square.

For a cubic on more than one component, give one inactive component 3 and another 1. For a cubic in one inactive component with a variable outside F, give that variable 2 and all remaining vertices of the component 3; no original square or cubic forbids this. The target is 4 and the right sum at most 0. Finally, if the target is entirely anchored in one component but the global guard g fails, all original words have length at most two. Give every variable 3: the target is 0 while every right term is 3 or 2, so the right sum does not dominate 0. These cases exhaust all failures, for arbitrary finite P and target.

\par\smallskip\noindent\textit{Unrestricted finite derivation.}\quad 
Keep all original summands throughout. Extract xxxx for every variable x in an original word of length at least four using \(xxxx\leq xyzt\) and three nonempty blocks for the remainder. For an original nonlinear word containing distinct x, y, obtain xy by \(xy\leq xyz\) with its nonempty remainder if it is not already quadratic. These pairs provide original graph connections. Propagate quartics with \(yyyy\leq xxxx+xy\). Merge any finite number of quartics by \(xxxy\leq xxxx+yyyy\), putting the accumulated nonempty content word in place of x. Normalize to obtain every quartic-or-longer target with active content. Recover a corresponding cubic by \(xyz\leq xyzt\), and a quadratic by \(xy\leq xyz\). Recover literals directly from an original containing term by \(x\leq xy\); the remaining factor is nonempty.

For each x in F keep its original square, or extract xx from an original cubic using \(xx\leq xyz\) with x as one original letter and the other two nonempty blocks. An anchored connection \(x=p0\), p1,..., \(pn=y\) is a positive path of original graph connections when x differs from y. Starting with xp1, apply \(xz\leq xx+xy+yz\) successively, always retaining xx and the next original connection. This derives xy. A square at a cubic variable is obtained directly, so no zero-length path substitutes for a real edge.

For a cubic in one anchored component, use the original global guard u of length at least three. Split u into three nonempty factors and use \(xxx\leq xx+yzt\) at each target variable, obtaining its cube. Connect the first and second, and second and third, target letters by the preceding anchored-path construction. Then \(xyz\leq xxx+yyy+zzz+xy+yz\) recovers the target. Repeated target letters use the retained squares and cubes.

All required witnesses are original or explicitly derived. Every substitution is nonempty. The fixed thirteen laws derive every valid absorption independently of the number of variables or original terms. Their validity follows from the criterion and can also be checked on their finitely many assignments. Distributive expansion and absorption in both directions complete the equational-basis proof.

Essential boundaries: \(xxxx\leq xxx\) and \(xxx\leq xx\) both fail at \(x=3\); \(xx\leq x\) also fails there. Mere connectedness is insufficient for quadratics: \(xz\leq xy+yz\) fails at \(x=z=2\), \(y=3\). A cubic needs the original global guard: \(xxx\leq xx\) fails even though the endpoint loop is present.

\subsection{Fourteen complete identities for semiring 14777}\label{case:329}
Signature (+, multiplication), no constants. The following fourteen identities, together with the additively idempotent semiring axioms, form a complete basis. Inequality \(u\leq P\) means \(u+P=P\).

{\small\par\smallskip
\CaseLine{\(xyz=yxz; xyz=xzy; xxyz=xyz; xxy=xyy;\)}
\CaseLine{\(x\leq xy; y\leq xy; yy\leq xy;\)}
\CaseLine{\(xz\leq xy+yz; yz\leq xy+xz; xxx\leq xyz;\)}
\CaseLine{\(xxx\leq xy+yyy; yyy\leq xxx+xy;\)}
\CaseLine{\(xxy\leq xxx+yyy; xy\leq xxy.\)}
\smallskip}

\par\smallskip\noindent\textit{Word values and long normal forms.}\quad 
Addition is 00004/\allowbreak{}01004/\allowbreak{}00224/\allowbreak{}00234/\allowbreak{}44444; multiplication is 44444/\allowbreak{}41444/\allowbreak{}44404/\allowbreak{}44404/\allowbreak{}44444. In the additive order, \(3<2<0<4\) and \(1<0<4\), while 1 is incomparable with 2 and 3. A quadratic has value 1 for the pair(1,\allowbreak{}1), value 0 for (2,\allowbreak{}3) or (3,\allowbreak{}3), and 4 otherwise. Every word of length at least three has value 1 precisely when every variable in its content has value 1, and 4 otherwise. This follows because after two factors a product lies in\{0,\allowbreak{}1,\allowbreak{}4\}; a further product is 1 only from 1*1 and is 4 otherwise. The fourteen laws can consequently be checked on their at most three variables; a separate exhaustive operation-table check accompanies this proof.

All words of length at least three with the same content are derivably equal under the first four laws. The first two allow every adjacent transposition in a word of length at least three (use a nonempty suffix for the first, or nonempty prefix for the second). If length exceeds max(3, content size), permute a duplicate to the first two positions and use \(xxyz=xyz\), splitting the remaining \(length\geq 2\) into nonempty y, z. Repeat. With two distinct letters and length 3, \(xxy=xyy\) transfers the surplus occurrence; with one letter the normal form is its cube. Thus a canonical long word has each letter once, padding with its least letter to reach length 3. Literals and quadratics remain unchanged; a square is not a cube (at 3 they are 0 and 4).

\par\smallskip\noindent\textit{Original-polynomial criterion.}\quad 
For a finite nonempty original polynomial P, let C be its content union. Build an undirected graph on C from the endpoints of every original quadratic (loops allowed). Isolated literal/long variables remain vertices. Let B be the set of original quadratic final variables. Let K be the union of contents of original words of \(length\geq 3\), and A the union of graph components meeting K. Then \(w\leq P\) holds exactly as follows:

\textbullet{} For a literal, its variable must belong to C.

\textbullet{} For a long word, its entire content must belong to A.

\textbullet{} For a quadratic ht, either \{h, t\} is a subset of A, or h, t belong to the same graph component and t belongs to B.

Every witness here is original. In particular an original quadratic cycle does not by itself activate cubes.

\par\smallskip\noindent\textit{Necessity with explicit assignments.}\quad 
If a target contains a variable outside C, assign that variable 4 and all variables of C the value 1 (other outside variables also 4). The right polynomial is 1 and the target is 4. This also separates literals outside C.

For assignments where P is not 4, every original long word has all variables 1. Every quadratic edge is either entirely 1, or has final value 3 and initial value 2 or 3. Thus being 1 propagates in both directions across an edge, and every variable of A must be 1.

For an explicit separator of a long target whose content is in C but not A, choose an inactive component meeting the target. Assign every variable of that component 3 and every other variable 1. Original quadratics are 1 or 0, long words 1, literals 1 or 3, so \(P\leq 0\); the target is 4.

For a failing quadratic ht, if h, t are in different components, at least one component is inactive. If the component of t is inactive, assign it 3 and all other components1. Otherwise assign the inactive component of h the value 3 and all others 1. The target is respectively 1*3 or 3*1, both 4; all original right summands are \(\leq 0\). If h, t lie in the same component, failure means this component is inactive and t is not an original quadratic final variable. Assign \(t=2\), every other variable of this component 3, and all other components1. There is no edge ending at t, so the original edges in this component still have values 0. Original long words are 1 and literals are \(\leq 0\). The target ends in 2 and is 4, including the case \(h=t\). These constructions exhaust failure of the stated conditions.

\par\smallskip\noindent\textit{Complete finite derivations.}\quad 
Use the first four laws to normalize all long original words. Every literal appearing in an original word follows by \(x\leq xy\), \(y\leq xy\), using a nonempty prefix/suffix when present.

For each original long word u and every k in its content, permute k first, split the remaining nonempty \(length\geq 2\) into y, z, and apply \(xxx\leq xyz\) to derive \(k^{3}\). Across an original quadratic ab, the laws \(yyy\leq xxx+xy\) and \(xxx\leq xy+yyy\) respectively propagate an available cube from a to b and from b to a. A finite graph traversal therefore derives \(k^{3}\) for every k in A. For any nonempty target content contained in A, combine these cubes successively with \(xxy\leq xxx+yyy\): substitute the already obtained long word for x and the next variable for y. The right \(x^{3}\) normalizes to that long word, and the left to a long word on the content union. Normalize to the desired long target. For a quadratic on A, the resulting long word on its content and \(xy\leq xxy\) recover it.

For a quadratic ht supplied by the other clause, choose an original edge at. The law \(yy\leq xy\) gives tt. Starting at t, traverse an undirected path towards h while retaining final t. If the next original edge is ab and bt is already available, \(xz\leq xy+yz\) supplies at. If at is already available, \(yz\leq xy+xz\) supplies bt. Thus every vertex in t's component becomes an initial variable of a derived word ending in t. This uses the original final-variable condition t in B even when \(h=t\). All original and previously obtained words remain absorbed by P, so every use of these laws is an equational inference under P.

These constructions work for unrestricted variable counts, word lengths, and summand counts, with no empty substitutions. Distributive normalization reduces an arbitrary identity to two sets of word absorptions, each derived above. Hence the fourteen laws are a finite identity basis. Any transfer to another catalogue type is conditional on separately verified HSP membership in both directions, not merely a bounded-word fingerprint.

\subsection{A finite basis for the shifted affine five-element semiring}\label{case:330}
All bases are relative to the additively idempotent semiring axioms in the constant-free signature \((+,\cdot)\). Write \(u\le v\) for \(u+v\approx v\).

\par\smallskip\noindent\textit{The algebra and the result.}\quad 
Let \(A=\{e, g, c, a, o\}\). Its addition is the join in the poset \(g<c<o\), \(e, a<o\), with no other comparabilities. The subset \(G=\{e, g\}\) is a multiplicative group of order two with identity \(e\). Set \(a^2=c\); every product with a factor outside \(G\) other than \(a^2\) is \(o\). In this displayed order the catalog map is \((2,\allowbreak{}1,\allowbreak{}0,\allowbreak{}3,\allowbreak{}4)\) and the type is 14818. Any word of length at least three evaluates to its group product when all letters lie in \(G\), and to \(o\) otherwise. This observation verifies associativity; direct inspection of translations verifies distributivity.

Theorem. The following nine laws form a finite identity basis for \(A\):

\CaseDisplay{\begin{array}{rll}
(C)&xy\approx yx,\\
(N)&x^3yz\approx xyz,\\
(S)&x^4y^4\approx x^4+y^4,\\
(U)&x\le x^3,\\
(O)&xyz\le x^3+y^3+z^3,\\
(L)&x^3\le x+xy,\\
(T)&(xy)^3\le xy+yzt,\\
(P)&xt\le xy+yz+zt,\\
(Q)&x^3\le x+y^4.
\end{array}}

Minimality is not asserted.

The first eight laws are valid for the same reasons as their counterparts for type 14659. For clarity, the only non-group quadratic value is \(a^2=c\), and every word of length at least three containing \(a\) is \(o\). Thus (L) at \(x=a\) has \(a+a^2=o\) if \(y=a\), and \(xy=o\) otherwise. In (T), if \(xy=c\) then \(x=y=a\) and \(yzt=o\). In (P), any edge of a three-edge path taking value \(c\) forces all four vertices to be \(a\) unless another edge is \(o\); otherwise all three edges are group-valued, and their product is \(xt\). If all three are the same group value then their join is \(xt\); if both group values occur their join is \(o\). The stable-word formula directly verifies (N),(S),(U),(O). Finally, if \(x\in G\) then \(x^3=x\); if \(x\notin G\) then \(x^3=o\), while \(y^4=e\) for \(y\in G\) and \(y^4=o\) otherwise. For \(x\in\{a, c\}\) one has \(x+e=o\). This proves (Q).

\par\smallskip\noindent\textit{Stable-word calculus.}\quad 
A word is stable if its length is at least three. Let \(c(w)\) be its variable content and \(\pi(w)\in\mathbb F_2^X\) its exponent-parity vector. For a nonempty stable-word family \(R\), write \(H(R)=\operatorname{aff}_{\mathbb F_2}\{\pi(r):r\in R\}\).

The following consequences of (C),(N),(S),(U),(O) hold for every number of variables.

1. Stable words with equal content and parity are equal. A stable word \(r\) satisfies \(r^3=r\).

2. If stable \(s, t\) have equal parity and \(c(s)\subseteq c(t)\), then \(s\le t\).

3. If \(R\ne\varnothing\), \(c(w)\subseteq\bigcup c(R)\), and \(\pi(w)\in H(R)\), then \(w\le\sum R\), even if \(w\) is short.

For (1), commute variables and use (N) in both directions to insert or remove pairs of occurrences while preserving a stable length. For (2), multiply (S), instantiated at \(y, z\), by a stable word containing \(y\): modulo (1), it gives \(sz^4=s+sz^4\); insert each missing variable. For (3), choose an odd product of members of \(R\) with the desired parity, insert pairs of members until its content covers \(c(w)\), and use (O) and \(r^3=r\) to bound it by \(\sum R\). Apply (1),(2), and, if \(w\) is short, (U) to \(w\le w^3\).

Two elementary necessary conditions will be used throughout. If \(w\le F\), then \(c(w)\subseteq c(F)\): assign an absent variable \(o\) and every other variable \(e\). Also \(\pi(w)\in H(F)\): on \(G\), an affine-hull separating linear functional assigns all summands of \(F\) one group value and \(w\) the other, which are incomparable.

\par\smallskip\noindent\textit{Preparation and the exact local criterion.}\quad 
Given a polynomial \(F\), first use (L),(U) to replace every linear summand \(x\) that also occurs in a nonlinear summand by the stable summand \(x^3\). Indeed, if the nonlinear summand is \(xu\), then \(x^3\le x+xu\), and \(x\le x^3\). Next use (T),(U) repeatedly to replace every quadratic edge \(xy\) incident with the content of a stable summand by \((xy)^3\). If the stable summand is \(yzt\) with \(z, t\) nonempty words, then \((xy)^3\le xy+yzt\). The process stops with an equivalent polynomial

\CaseDisplay{F\approx\sum L_0+\sum R+\sum E_0,\quad\textup{(\ensuremath{1})}}

where \(R\) consists of stable words, \(E_0\) of quadratic edges, and the three variable sets \(L_0\), \(W=c(R)\), and \(V(E_0)\) are disjoint. The original total content and parity hull are preserved.

For a free edge graph \(E_0\), call \(xy\) odd-walk eligible if an odd walk joins \(x\) to \(y\), including \(x=y\) and closed walks. Then the following criterion is exact.

Proposition. If \(R\ne\varnothing\) and \(0\in H(R)\), then

\CaseDisplay{w\le F\quad\Longleftrightarrow\quad
c(w)\subseteq c(F)\ \text{and}\ \pi(w)\in H(F).\quad\textup{(\ensuremath{2})}}

In every other case, \(w\le F\) holds precisely when at least one of these mutually compatible cases applies:

\textbullet{} \(w=x\) is one of the displayed linear summands in \(L_0\);

\textbullet{} \(c(w)\subseteq W\) and \(\pi(w)\in H(R)\);

\textbullet{} \(w=xy\) is an odd-walk eligible edge in \(E_0\).

For sufficiency in (2), select a stable fourth power \(h^4\) with \(c(h)=W\) and parity zero. The stable-word calculus gives \(h^4\le\sum R\le F\). Law (Q), substituted with each short summand \(s\) in (1), now gives \(s^3\le s+h^4\le F\). By (U), \(s\le s^3\), so each short summand may be replaced by its stable cube. The resulting polynomial is entirely stable and has the same content and parity hull; the stable-word calculus proves (2). The converse is the two necessary conditions above.

In the other case, stable-word calculus proves the second bullet, and (P), iterated along an odd walk, proves the third. The first bullet is immediate. For necessity, first suppose \(R\ne\varnothing\) and \(0\notin H(R)\). Elementary linear algebra gives an assignment of group values to \(W\) under which every member of \(R\) is \(g\). Assign every variable of every free edge component \(a\), and every independent linear variable \(c\). The value of \(F\) is \(g\) or \(c\). It separates a linear target in a free component (value \(a\)). Every stable target touching a free component or an independent linear variable has value \(o\), so is also separated. A quadratic touching an independent linear variable has value \(o\), as does a quadratic with one endpoint in \(W\) and the other in a free component. A quadratic whose two endpoints lie in free components has value \(c\) under this assignment; the within-component and cross-component cases are treated separately below. If a target contains a variable absent from \(F\), assign it \(o\). The same arguments work when \(R\) is empty, omitting the group assignment.

For a target supported in \(W\) but with parity outside \(H(R)\), choose the group assignment on \(W\) so that every \(r\in R\) is \(g\) and the target is \(e\). This is possible because \(0\notin H(R)\) and the target parity is outside that affine hull. With all free components at \(a\) and all independent linear variables at \(c\), the right side is \(g\) or \(c\), and \(e\) lies below neither.

It remains to check a quadratic \(xy\) inside a free edge component. If its component is nonbipartite, an odd closed walk converts any walk between \(x, y\) to an odd one. If it is bipartite, an odd walk exists exactly when \(x, y\) have opposite colors. When they have the same color, assign group values \(e, g\) to this component according to a bipartite coloring, so every displayed edge there evaluates to \(g\), while \(xy=e\). Assign every other free component \(a\), and use the group assignment on \(W\) described above, when necessary. All right summands are \(g\) or \(c\), and \(e\nleq g, c\). A quadratic with endpoints in different free components, or one free endpoint and one outside its component, has parity outside \(H(F)\): every right parity vector has an even coordinate sum on each free component. The parity necessary condition excludes it. This finishes the proof of the proposition.

Each polynomial identity is equivalent to absorption of each summand in both directions. The criterion supplies derivations for all semantically valid absorptions, with no bound on variables or summands. This proves the theorem.

\subsection{Nonfinite basability of 14839 from high-girth three-uniform hypergraphs}\label{case:331}
Work without constants. The multiplication is the monogenic nilpotent semigroup g, \(g^{2}\), \(g^{3}\), \(g^{4}\), \(g^{5}\), with \(g^{5}\) absorbing. The labels are \(g=1\), \(g^{2}=2\), \(g^{3}=3\), \(g^{4}=0\), \(g^{5}=4\). Its addition has principal ideal

{\small\par\smallskip
\CaseLine{\(down(g^{4})=\{g, g^{2},g^{4}\},\)}
\smallskip}

and \(g^{5}\) is the greatest element. The other proper principal ideals are \{g\},\{\(g^{2}\)\},\{g, \(g^{3}\)\}. This source has no suitable idempotent pair for the previous cycle theorem. We prove an independent obstruction instead.

\par\smallskip\noindent\textit{Hypergraphs and the critical identities.}\quad 
By Lemma~\ref{lem:high-girth-exists}, for any \(L\ge2\) there is a finite three-uniform hypergraph \(H\) of Berge girth greater than \(L\) admitting no proper two-colouring. Distinct edges share at most one vertex, since two common vertices would give a Berge cycle of length two.

Use a distinct variable v for every vertex, and a further variable x. Let \(P_{H}\) be \(x^{4}\) plus ALL SIX word permutations of every edge triple. The critical identity is \(x^{5}\leq P_{H}\). In 14839, if x is not g, then \(x^{4}=g^{5}\) and \(P_{H}\) is top. If \(x=g\), then \(x^{4}=g^{4}\). Keeping \(P_{H}\) below \(g^{4}\) requires every edge product to be \(g^{4}\). A triple product equals \(g^{4}\) precisely when its values are two copies of g and one of \(g^{2}\). These choices would properly two-color H, impossible. No other proper ideal contains \(x^{4}\). Thus \(P_{H}\) is always top and the critical identity is valid.

\par\smallskip\noindent\textit{Two facts about hyperforests.}\quad 
A linear three-uniform hypergraph with no Berge cycle is an incidence forest. Its exact-one assignments give each vertex an activity bit0/\allowbreak{}1, with exactly one active vertex in each edge.

Every single vertex can have either activity. For any two distinct vertices, both can be inactive, and either one can be active while the other is inactive. Both can be active exactly when they share no edge. To prove these claims, follow their unique incidence path, if one exists. Set internal path vertices inactive. Each path edge whose two path vertices are inactive has its third vertex active; an active endpoint instead satisfies its incident path edge. Extend off the path:at an edge attached to an assigned vertex, set both other vertices inactive if it is active; otherwise make exactly one of them active. In separate components use independent extensions. This proves all the assertions.

Consequently a triple of distinct vertices has exactly one active vertex in every exact-one assignment iff it is an edge. If it is not an edge, some pair shares no edge:otherwise its three pairwise shared edges would form a Berge triangle, unless all three vertices were one edge. That pair can be simultaneously active. Repeated triples also cannot have activity sum constantly1, since their repeated vertex can be active.

A second fact:if a subhypergraph of H has a vertex cover of size at most k, and H has girth greater than 2k, then the subhypergraph is a forest. Indeed, in a shortest Berge cycle, any vertex lies in at most two cycle edges. A vertex meeting nonconsecutive edges would give a shorter cycle; linearity rules out a second shared vertex between consecutive edges. Thus a k-vertex cover can cover at most 2k cycle edges. But every cycle has length greater than 2k. This contradiction proves the fact, including subhypergraphs with arbitrarily many edges.

\par\smallskip\noindent\textit{A substitution barrier closed under literals and squares.}\quad 
Fix H with girth greater than max(2k,3). Let F consist of x, \(x^{2}\), \(x^{4}\), every vertex literal v, every vertex square \(v^{2}\), and every permutation of an edge triple. Suppose q has at most k source variables and a nonempty polynomial substitution phi expands q entirely into F. We show that every expanded word of phi(u) is in F whenever \(u\leq q\) is a valid 14839 absorption. No contained-content condition is assumed.

For source variables occurring nonlinearly in q, their images have the following exhaustive forms: \(X=\{x\}\), \(Y=\{x^{2}\}\), \(Z=\{x^{3}\}\) or \{x, \(x^{3}\)\}; \(V=\{v\}\) for one vertex; or E, a nonempty sum of ordered distinct-vertex pairs, all completed to hyperedges by the SAME single vertex c, optionally together with the literal c. Call remaining variables L; they occur only literally, and their images are arbitrary nonempty subsets of F. A literal-only singleton X, Y or V may be assigned its indicated type instead.

Here is a detailed justification of the image list. The expanded source words have degree at most 4. A degree 4 source word must expand to \(x^{4}\), so all factors are X. A degree 3 source word either expands to \(x^{4}\) with factors X, X, Y, or to an edge with three vertex factors. Varying a vertex factor would give two edges sharing the other two vertices, contrary to linearity, so each is singleton. A degree 2 source word can expand to \(x^{2}\) via XX, to \(x^{4}\) via YY or XZ/ZX, to \(v^{2}\) via two equal singleton vertex images, or to edge triples via V/E in either order. For V/E, varying the vertex factor also gives edges sharing a pair, so V is singleton. Each pair in E is distinct and omits its complement. The complement is unique across every occurrence of E, by linearity. The only mixed degrees permitted inside one image are the displayed Z option \{x, \(x^{3}\)\}, paired with X, and the E option adding its complement literal c, paired with \(V=\{c\}\). These expand to \(x^{2}\)+\(x^{4}\) or \(c^{2}\)+edge respectively. All other mixed degrees or mixed x/vertex images produce \(x^{3}\), a mixed x/vertex product, a nonsquare pair, a repeated-vertex triple, a vertex-containing quartic, or length at least 5. An image containing an edge or \(x^{4}\) cannot occur nonlinearly. These observations also cover polynomial images with many monomials and repeated source variables.

Let H0 be the union of edge triples contributed by nonlinear source words. It has a vertex cover consisting of the singleton V images: VVV edges use these vertices, and V/E stars use their fixed complement. There are at most k such vertices. By the cover lemma H0 is a forest. Literal L images may contain the WHOLE H; they are intentionally excluded from H0 and receive free formal values below.

Choose any exact-one activity assignment on H0, including its isolated vertices. Give a vertex with bit b the algebra value \(g^{1+b}\). Assign \(X=g\), \(Y=g^{2}\), \(Z=g^{3}\). Each E whose complement is c receives g\ensuremath{^{3-b(c)}}; every pair in its image has this same value because its completed edge has exactly one active vertex. Its optional literal c is below that value:when c is inactive, \(g\leq g^{3}\); when active, both are \(g^{2}\). Likewise the optional x in a Z image is below \(g^{3}\). Hence each nonlinear-variable image polynomial evaluates exactly to its formal value. Assign each L the FREE value \(g^{4}\), and any source variable absent from q the value \(g^{5}\). These are legitimate source assignments, not evaluations of L images. Every word of q has value g, \(g^{2}\) or \(g^{4}\):its possible nonlinear forms were listed above. Hence \(q\leq g^{4}\), and validity forces \(u\leq g^{4}\) for every exact-one assignment.

A nonlinear u cannot contain L or an absent variable, as that would give exponent at least 5. A literal L has image in F already; an absent literal is impossible. For all other u, each image polynomial evaluates to its formal source value, with each monomial below that value. Monotonicity and distribution therefore imply that every expanded word of phi(u) has \(value\leq g^{4}\) under EVERY exact-one assignment on H0. We can classify such expanded words directly:

\textbullet{} Degree 1:only x and vertex literals, all in F.

\textbullet{} Degree 2:\(x^{2}\) and vertex squares are allowed. x v fails when v is active. Two distinct vertices fail when their activities differ.

\textbullet{} Degree 3:an all-vertex word must have activity sum exactly1 under every assignment, hence is an edge of H0 by the forest fact above. A word containing x fails when all its at most two vertex variables are inactive; then its exponent is 3. Repeated vertex triples fail too.

\textbullet{} Degree 4:only \(x^{4}\) is possible. Any occurring vertex can be made active, giving exponent at least 5.

\textbullet{} Degree at least 5 always gives \(g^{5}\), which is not below \(g^{4}\).

Every possible expanded word is therefore in F. This proves the barrier for arbitrary polynomial substitutions and repeated variables, including large E stars and L images containing the entire obstructing hypergraph.

\par\smallskip\noindent\textit{Nonfinite basis and finite-target transfer.}\quad 
Suppose 14839 had a finite basis with variable bound k, including the semiring axioms. Choose a non-two-colorable H with girth greater than max(2(k+2),3). In a derivation of \(P_{H}=P_{H}+x^{5}\), expand all terms into nonempty sums of words. Consider the first step introducing a word outside F. Earlier steps may introduce literals and squares, but those are explicitly allowed in F. Select an expanded component of the identity application; each replacing word satisfies a valid absorption against the old component. Incorporate nonempty left and right contexts using at most two fresh variables, and omit absent contexts. The barrier with k+2 variables forbids the first new word. Yet the final \(x^{5}\) is outside F. This contradiction proves NFB without any bound on input polynomial substitutions or derivation lengths.

For another finite semiring T, the same obstruction transfers if 14839 belongs to Var(T) and T satisfies every critical hypergraph identity. The following finite condition is sufficient for the latter. For each a, b with \(a^{4}\leq b\) and \(a^{5}\) \(\nleq b\), form the ternary relation \(R_{b}\) of triples (r, s, t) for which ALL SIX ordered products are \(\leq b\). Require a two-coloring of the elements of T in which every triple of \(R_{b}\) has both colors. A counterassignment for some H would map every edge into \(R_{b}\) and pull this coloring back to a proper coloring of H, a contradiction. If no such coloring of \(R_{b}\) exists, the finite condition has merely FAILED TO CERTIFY transfer; it is not a counterexample or a classification result. Full HSP and the sufficient coloring condition are both required.

\subsection{Unbounded identity rank in a five-element group-clone semiring}\label{case:332}
All algebras have the constant-free signature \((+,\cdot)\). We work relative to the additively idempotent semiring axioms. Write \(p\le q\) for \(p+q\approx q\).

\par\smallskip\noindent\textit{Algebra and theorem.}\quad 
Let \(D=\{e, g, b, c, o\}\) have additive order generated by \(b<g<o\), \(b<c<o\), and \(e<c<o\), with \(g, c\) incomparable. The subset \(G=\{e, g\}\) is the group of order two with identity \(e\). Let \(\rho(b)=g\), \(\rho(c)=o\), and fix \(e, g, o\). Define multiplication by the products of the \(\rho\) images in the flat group-with-zero \(G\cup\{o\}\). The join-preserving retraction \(\rho\) verifies the semiring axioms. In the displayed element order, its catalog map is \((2,\allowbreak{}3,\allowbreak{}1,\allowbreak{}0,\allowbreak{}4)\), giving type 14842.

Theorem. The semiring \(D\) is nonfinitely based. More precisely, for every \(k\ge3\) and every even \(m>k\), the valid identity

\CaseDisplay{x_1^2\le F_m,
\qquad F_m=x_1+\cdots+x_m+x_1\cdots x_m,\quad\textup{(\ensuremath{1})}}

is not a consequence of the identities of \(D\) involving at most \(k\) variables, together with the semiring axioms.

In fact, \(\mathsf V(D)\) is not a limit variety: the final section constructs a proper nonfinitely based subvariety.

\par\smallskip\noindent\textit{Semantic facts.}\quad 
Call a word stable if it has length at least two. For a word \(w\), let \(c(w)\) be its content and \(\pi(w)\in\mathbb F_2^X\) its parity vector. Given a nonempty polynomial

\CaseDisplay{P=\sum_{x\in L}x+\sum R,}

let \(R\) be its stable summands, \(N=\bigcup_{r\in R}c(r)\), and \(V=\{\pi(r):r\in R\}\). Write

\CaseDisplay{S(P)=\operatorname{span}(V),\qquad
K(P)=\operatorname{aff}\bigl(V\cup\{\epsilon_x:x\in L\cap N\}\bigr).}

If \(R\) is empty, the second alternative below is absent.

Lemma 1 (absorption criterion). In \(D\),

\CaseDisplay{w\le P\quad\Longleftrightarrow\quad
\bigl(w=x\text{ for some }x\in L\bigr)
\quad\text{or}\quad
\bigl(c(w)\subseteq N,\ \pi(w)\in S(P)\cap K(P)\bigr).\quad\textup{(\ensuremath{2})}}

Proof. A stable word evaluates to the product of its \(\rho\) images in \(G\) when every letter is in \(G\cup\{b\}\), and to \(o\) otherwise. An absent variable is separated by assigning it \(o\) and all other variables \(e\). A variable occurring only linearly on the right is separated in a stable target by assigning that variable \(c\) and the others \(e\). Thus \(c(w)\subseteq N\) is necessary except for an explicitly listed linear target.

If \(\pi(w)\notin S(P)\), choose a functional \(\alpha\) annihilating \(V\) and taking value \(1\) at \(\pi(w)\). Assign group values to \(N\) accordingly, then change every right linear variable with value \(g\) to \(b\); stable values are unchanged. Assign variables occurring only linearly the value \(c\). Every stable right word is \(e\), and all linear right terms are below \(c\). The target has value \(g\), including a linear target not in \(L\), and \(g\nleq c\).

If \(\pi(w)\in S(P)\setminus K(P)\), then \(0\notin K(P)\), since an affine space containing zero is the linear span of its generators. Choose a functional taking value \(1\) on \(V\) and the unit vectors of \(L\cap N\), and value \(0\) at \(\pi(w)\). Such a functional exists by affine separation because zero is outside \(K(P)\). Assign group values on \(N\), lower the right linear \(g\) values to \(b\), and give the variables occurring only linearly value \(g\). The right side is \(g\); the target is \(e\), so absorption fails.

Conversely assume the second alternative of (2). If a target variable receives \(c\) or \(o\), a stable right word containing it is \(o\). If any stable right word is \(o\), or stable right words take both group values, the right side is \(o\). Otherwise all stable right words have one group value. If that value is \(e\), membership in \(S(P)\) forces the target's group image to be \(e\), hence its value is \(e\). If it is \(g\), a target value \(g\) or \(b\) is absorbed. For target value \(e\), membership in \(K(P)\) ensures that at least one linear variable in \(L\cap N\) has group image \(e\), hence actual value \(e\), and the right side contains \(e+g=o\). The case of an explicit right linear summand is tautological. When \(R\) is empty, assigning all right variables \(c\) excludes every stable target. This proves (2). \(\square\)

Let \(X=\{x_1,\ldots, x_m\}\), \(m\ge2\), and let \(v=\epsilon_1+\cdots+\epsilon_m\).

Corollary 2. For even \(m\), identity (1) holds. Every stable summand of a polynomial equivalent to \(F_m\) in \(D\) has parity either \(0\) or \(v\).

Proof. The stable family of \(F_m\) consists of the single word \(x_1\cdots x_m\), so \(S(F_m)=\{0, v\}\). Its vector together with the \(m\) linear vectors has sum zero using \(m+1\) terms, an odd number. Thus \(0\in K(F_m)\) and Lemma 1 gives (1). If \(Q\approx F_m\), each summand of \(Q\) is absorbed by \(F_m\); Lemma 1 gives the asserted restriction. \(\square\)

We shall use also the following group consequence. If an identity \(A\approx B\) holds in \(D\), the affine hulls of their word-parity vectors coincide. Indeed, assignments in the flat subgroup \(G\cup\{o\}\) give affine separation of any vector outside the other hull. In particular, if every word of \(A\) has the same parity \(t\), then every word of \(B\) has parity \(t\).

\par\smallskip\noindent\textit{Bounded-variable substitutions cannot create the first even term.}\quad 
For polynomials over \(X\), say that property \(\mathcal P_m\) holds when every stable summand has parity \(v\). Linear summands may be any of the \(x_i\). Words are formal nonempty commutative words; polynomials are nonempty finite sums, with repetitions discarded. No word normalization or length bound is assumed.

Lemma 3. Let \(P\approx Q\) be an identity of \(D\) using fewer than \(m\) variables, and let \(\sigma\) substitute nonempty polynomials over \(X\) for its variables. If \(\sigma(P)\) has property \(\mathcal P_m\), no stable summand of \(\sigma(Q)\) has parity zero.

Proof. Take a word \(w\) of \(Q\). Then \(D\models w\le P\). If \(w\) is a linear summand of \(P\), every summand of \(\sigma(w)\) already occurs in \(\sigma(P)\) and the conclusion follows. Otherwise Lemma 1 yields

\CaseDisplay{c(w)\subseteq N(P),\qquad \pi(w)\in K(P).\quad\textup{(\ensuremath{3})}}

In particular \(P\) has at least one stable summand.

For each variable \(z\in N(P)\), all monomials of \(\sigma(z)\) have one and the same parity, say \(a_z\). To see this, choose a stable word of \(P\) containing an occurrence of \(z\), fix a monomial choice at every other occurrence, and vary the monomial selected at this occurrence. All resulting words in \(\sigma(P)\) are stable and must have parity \(v\). The two selections therefore have equal parity. Repeated occurrences cause no difficulty: distributivity permits independent choices at each occurrence.

Define the linear map \(A\) on the parity space of \(N(P)\) by \(A(\epsilon_z)=a_z\). Every stable word \(r\) of \(P\) satisfies \(A(\pi(r))=v\). If \(z\in L(P)\cap N(P)\), then \(\sigma(z)\) occurs as a subpolynomial of \(\sigma(P)\). Its monomials are either linear words \(x_i\) or stable words of parity \(v\). Since they have a common parity, \(a_z\) is either \(v\) or one unit vector \(\epsilon_i\). There are fewer than \(m\) such variables \(z\), so the set \(I\) of unit-vector indices obtained in this way is a proper subset of \(\{1,\ldots, m\}\).

Applying \(A\) to (3), every monomial of \(\sigma(w)\) has parity in

\CaseDisplay{\operatorname{aff}\bigl(\{v\}\cup\{\epsilon_i:i\in I\}\bigr).\quad\textup{(\ensuremath{4})}}

The vectors displayed in (4) are linearly independent: a coordinate outside \(I\) forces the coefficient of \(v\) to vanish, and then each unit-vector coefficient vanishes. The affine hull of a nonempty linearly independent set does not contain zero. Consequently no monomial of \(\sigma(w)\) has zero parity, proving the lemma. \(\square\)

Lemma 4. In a derivation from identities of \(D\) involving fewer than \(m\) variables, the first appearance of a stable zero-parity summand is impossible if all previous polynomials are equivalent to \(F_m\) in \(D\) and the initial polynomial has property \(\mathcal P_m\).

Proof. Work in the free commutative additively idempotent semiring on \(X\). Associativity, commutativity and distributivity identify its elements with the formal polynomials above. Equational consequence of the remaining identities is the congruence generated by their polynomial substitutions, allowing addition and multiplication contexts with a single occurrence of the substituted side. At a generating step, after distribution, a context gives an unchanged additive part and finitely many monomial multiples of the substituted side. Some multiples may have empty multiplicative context, meaning that the substituted side occurs directly in the additive sum.

Suppose a first stable zero-parity word appeared. Before this step, Corollary 2 implies that every stable summand has parity \(v\). For a multiple with empty multiplicative context, the substituted left side is a subpolynomial of the preceding polynomial, so Lemma 3 excludes a zero-parity output.

For a nonempty monomial context \(u\), every word of \(u\sigma(P)\) has length at least two, and hence parity \(v\). Therefore every word of \(\sigma(P)\) has parity \(v+\pi(u)\). The group consequence above, applied to the valid identity \(\sigma(P)\approx\sigma(Q)\), forces every word of \(\sigma(Q)\) to have that same parity. All words of \(u\sigma(Q)\) consequently have parity \(v\), excluding a zero-parity output again. The unchanged additive part cannot introduce a new word. Reversing an identity is covered by interchanging its two sides. This proves the lemma. \(\square\)

\par\smallskip\noindent\textit{Nonfinite basability.}\quad 
Fix \(k\ge3\) and an even \(m>k\). If (1) were a consequence of the identities of \(D\) with at most \(k\) variables, the definition of the generated congruence would give a finite derivation over \(X\) from \(F_m\) to \(F_m+x_1^2\). All its steps are valid in \(D\), so all intermediate polynomials are equivalent to \(F_m\). The initial polynomial has property \(\mathcal P_m\), whereas the final polynomial contains the stable zero-parity word \(x_1^2\). Lemma 4 rules out the first such appearance, a contradiction. This proves the stronger assertion of the theorem.

If \(D\) had a finite identity basis, choose \(k\ge3\) at least the number of variables in every basis identity. The valid identity (1), for an even \(m>k\), would then follow from identities involving at most \(k\) variables, contradicting the assertion just proved. Thus \(D\) is nonfinitely based. \(\square\)

\par\smallskip\noindent\textit{A proper nonfinitely based subvariety.}\quad 
Let \(\mathcal V=\mathsf V(D)\) and define

\CaseDisplay{\mathcal W=\mathcal V\cap
\operatorname{Mod}\{x+x^2\approx x^2+x^3\}.\quad\textup{(\ensuremath{5})}}

The additional identity, denoted (J), fails in \(D\): at \(x=b\), its left side is \(b+e=c\) and its right side is \(e+g=o\). Hence \(\mathcal W\) is a proper subvariety of \(\mathcal V\).

Theorem 5. The proper subvariety \(\mathcal W\) is nonfinitely based. Consequently the five-element semiring \(D\) does not generate a limit variety.

Proof. We first strengthen the invariant used above. Every polynomial with property \(\mathcal P_m\) is absorbed by \(F_m\) in \(D\): each linear summand occurs in \(F_m\), and each stable summand has parity \(v\), contains all of \(X\), and satisfies Lemma 1. Therefore every substituted side obtained from such a polynomial by a valid \(D\)-identity is still absorbed by \(F_m\). Its stable summands have parity \(0\) or \(v\). Lemma 3 and the nonempty-context argument of Lemma 4 exclude parity zero when the identity has fewer than \(m\) variables. Thus these generating steps preserve the full property \(\mathcal P_m\), without requiring earlier intermediate polynomials to be equivalent to \(F_m\) in \(D\).

We claim that the additional identity (J) also preserves \(\mathcal P_m\) in every polynomial context. Let \(q\) be the nonempty polynomial substituted for its variable. Both sides of (J) contain \(q^2\), and hence contain the square \(t^2\) of every monomial \(t\) of \(q\). Such squares have zero parity. Consequently neither side can occur without a nonempty multiplication context inside a polynomial with \(\mathcal P_m\).

Now multiply a substituted side by a nonempty monomial \(u\). Every term is stable. If \(u(q+q^2)\) has property \(\mathcal P_m\), then the terms \(ut\) and \(ut^2\), for each monomial \(t\) of \(q\), give

\CaseDisplay{\pi(u)+\pi(t)=v,\qquad \pi(u)=v.}

So every \(t\) has zero parity, and every word of \(u(q^2+q^3)\) has parity \(v\) as well. Conversely, if \(u(q^2+q^3)\) has that property, the terms \(ut^2\) and \(ut^3\) give the same conclusion; thus every word of \(u(q+q^2)\) has parity \(v\). A polynomial multiplication context is a sum of monomial contexts, and additive contexts retain their other summands. This proves the claim in both directions.

It follows that identities of \(D\) with at most \(k<m\) variables, even supplemented by (J), cannot derive (1) from \(F_m\): the initial polynomial has property \(\mathcal P_m\), while \(F_m+x_1^2\) does not. But every identity (1) holds in \(\mathcal W\) because \(\mathcal W\subseteq\mathcal V\).

Suppose \(\mathcal W\) had a finite basis \(\Gamma\). Each member of \(\Gamma\) follows from \(\operatorname{Id}(D)\cup\{\mathrm J\}\) by a finite equational derivation. The union of the finitely many \(D\)-identities used in these finitely many derivations has a finite variable bound \(k\). Together with (J) these identities would imply \(\Gamma\), hence every valid identity (1) of \(\mathcal W\). Choose even \(m>\max(k,3)\) to obtain the contradiction just established. Thus \(\mathcal W\) is nonfinitely based. Since it is a proper subvariety of \(\mathcal V\), the latter is not a limit variety. \(\square\)

\subsection{Six complete identities for semiring 14865}\label{case:333}
Signature (+, multiplication), no constants. Relative to the additively idempotent semiring axioms a complete finite basis is

{\small\par\smallskip
\CaseLine{\(xy=yx; xyz=ttt; x\leq yyy;\)}
\CaseLine{\(xy\leq xx+yz; xy\leq x+xz+yt; x\leq xx+zt+z.\)}
\smallskip}

Addition:00044/\allowbreak{}01034/\allowbreak{}00244/\allowbreak{}43434/\allowbreak{}44444. Multiplication:44444/\allowbreak{}42424/\allowbreak{}44444/\allowbreak{}42444/\allowbreak{}44444. The additive order \(is 1<0<4\),\(2<0<4\), and \(1<3<4\);\(2+3=0+3=4\). A quadratic is 2 exactly at(1,\allowbreak{}1),(1,\allowbreak{}3),(3,\allowbreak{}1), and 4 otherwise. A word of \(length\geq 3\) is 4 since both quadratic values 2,\allowbreak{}4 multiply to 4. This proves commutativity, cubic top, and \(x\leq yyy\). For \(xy\leq xx+yz\), a right side below 4 has \(x=1\) and y in\{1,\allowbreak{}3\}, so \(xy=2\). For \(xy\leq x+xz+yt\), a right side below 4 forces both edges to be 2, so all their endpoints are 1 or 3; the literal x cannot be 3 since \(3+2=4\). Thus \(x=1\) and \(xy=2\). For \(x\leq xx+zt+z\), a right side below 4 forces \(x=1\), and \(xx=2\) plus the literal-edge witness z gives a value bounding 1. All six laws are also independently checked on every assignment.

\par\smallskip\noindent\textit{Exact original-polynomial criterion.}\quad 
If P has an original long word, every target is absorbed. Otherwise let L be the original literal variables, E the original unordered quadratic pairs, V their endpoints, S their loop variables, and J the statement L intersects V. Then:

\textbullet{} No long target is absorbed.

\textbullet{} A literal h is absorbed iff h is in L, or h is in S and J holds.

\textbullet{} A quadratic ht is absorbed iff both endpoints lie in V, and either ht is an original unordered edge, or an endpoint lies in S union L.

Original literals can help recover a pair, but cannot alone recover a new literal at another variable. The latter needs an original square.

\par\smallskip\noindent\textit{Necessity: explicit separating assignments.}\quad 
With all variables 1, right edges are 2 and literals 1, so \(P\leq 0\) while long targets are 4. For a quadratic containing an endpoint outside V, assign V variables 1 and all others 2. Edges are 2 and literals 1 or 2, so \(P\leq 0\); the target containing 2 is 4. If both endpoints lie in V but their edge is missing and neither endpoint is in S union L, assign both target endpoints 3 and all others 1. No original edge has two3 endpoints, because it would be the missing target edge or forbidden loop; all edges are 2. Since the two endpoints are not literals, all original literals are 1, so \(P\leq 0\) while the target \(is 3*3=4\).

For a missing literal h, first suppose h is not in S. If h is outside V as well, assign \(h=4\) and all others 1; no original word contains h (h is not in L), and \(P\leq 0\). If h belongs to V but has no loop, assign \(h=3\) and all others 1. Every incident edge is 2, other edges 2, literals 1; thus \(P\leq 0\) but \(h=3\) is not below 0. Finally if h is in S but J is false, assign V variables 1 and all other variables 2. Every right term is 2, whereas \(h=1\) is not below 2. These cases exhaust failure of the criterion and work with arbitrary variable sets and polynomial sizes.

\par\smallskip\noindent\textit{Complete finite derivations.}\quad 
An original long word supplies a fixed variable's cube by \(xyz=ttt\), with a nonempty prefix substituted for x. The law \(x\leq yyy\) with the target word as x then recovers any target. Without long words, orient each original quadratic using commutativity. Original targets require no inference.

For a new literal h, there is an original square hh and an original literal z with incident edge zt. Use \(x\leq xx+zt+z\). For a new pair ht, orient it so h belongs to S union L. If h has an original square, choose an original incident edge tz and use \(xy\leq xx+yz\). Otherwise h is an original literal and belongs to V; choose incident original edges hz and ta. The law \(xy\leq x+xz+yt\) gives \(ht\leq h+hz+ta\leq P\). Every variable image is nonempty and every witness comes from the original polynomial. Duplicate edges or variables cause no problem.

The exact criterion and these derivations prove every valid absorption. Distributive normalization followed by summandwise absorption in both directions proves all identities. Hence the six laws are complete without bounds on variables, word lengths, or summand counts. Transfer to 14849 requires independently checked complete HSP in both directions; short polynomial matching alone never determines finite basability.

\subsection{Eleven complete identities for semiring 14853}\label{case:334}
Signature (+, multiplication), no constants. Relative to the additively idempotent semiring axioms take the following identities, where \(u\leq P\) abbreviates \(u+P=P\):

{\small\par\smallskip
\CaseLine{\(xyz=ttt; x\leq yyy;\)}
\CaseLine{\(xy\leq xx+zy; xy\leq yy+xz; xy\leq yy+zx;\)}
\CaseLine{\(xy\leq yx+zy; xy\leq yx+zt+z; xy\leq yx+tz+z;\)}
\CaseLine{\(xy\leq x+xz+ty; x\leq xx+zt+z; x\leq xx+tz+z.\)}
\smallskip}

Addition:00044/\allowbreak{}01034/\allowbreak{}00244/\allowbreak{}43434/\allowbreak{}44444. Multiplication:44444/\allowbreak{}42404/\allowbreak{}44444/\allowbreak{}42444/\allowbreak{}44444. A quadratic has value 0 at(1,\allowbreak{}3), value 2 at(1,\allowbreak{}1),(3,\allowbreak{}1), and 4 otherwise. Every word of length at least three is 4. The additive order has \(1<0<4\),\(2<0<4\),\(1<3<4\), and \(2+3=0+3=4\).

\par\smallskip\noindent\textit{Exact criterion and necessity.}\quad 
For a polynomial P with no original long word, let L be its literal variables, E its directed quadratic edges, V their endpoints, S its loop variables, T its edge final variables, and J mean L intersects V. All these sets refer to original summands. Then:

\textbullet{} A literal h is absorbed exactly when h belongs to L, or h belongs to S and J holds.

\textbullet{} A pair ht is absorbed exactly when h, t belong to V; either ht or th is an original edge, or one endpoint belongs to S union L; and t belongs to T or J holds.

\textbullet{} No long word is absorbed.

If P has an original long word, all targets are absorbed. To prove necessity, first assign all variables 1:long targets are 4 while P is at most 0. For a pair with an endpoint outside V, assign V to 1 and all other variables 2. Its value is 4 and P is at most 0. If the pair has endpoints in V but lacks the unordered edge and an S-or-L endpoint, assign its endpoints 3 and all others 1. Every right edge is 0 or 2 and every literal is 1; the target is 4. This includes a missing loop. Finally suppose t is outside T and J fails. Assign T to 1, V minus T to 3, and everything outside V to 2. All right edges and literals are 2. Since t belongs to V minus T, the target is 0 or 4, neither below 2.

For a missing literal h with no original loop, assign \(h=3\) if h is in V, or \(h=4\) otherwise, and all other variables 1. Edges are at most 0 and literals 1; h is not below their sum. If h is an original loop but J fails, assign V to 1 and all other variables 2. Right terms are 2 while \(h=1\) is not below 2. These cases exhaust every failed condition, for arbitrary variable sets and polynomial sizes.

\par\smallskip\noindent\textit{Sufficiency of the semantic criterion.}\quad 
Consider an assignment with P below 4. If P contains an edge, all its endpoint values are 1 or 3, and no original edge has two3 endpoints. Every loop endpoint is 1. Any literal endpoint must also be 1, as 3 joined with either edge value 0 or 2 would be 4. Thus an original unordered edge, loop endpoint, or literal endpoint ensures a permitted target pair has value 0 or 2. If \(P=0\), this suffices. If \(P=2\), all original final endpoints are 1 and J fails:an incident literal 1 would make the sum 0. The target's required original final endpoint is therefore 1, so its value is 2. A new permitted literal is at an original loop, hence is 1, and J forces \(P=0\), so it is absorbed. If there are no edges, only original literal targets meet the criterion. This proves sufficiency and the word-value formula directly from the tables.

\par\smallskip\noindent\textit{Unrestricted complete derivations.}\quad 
Every displayed law satisfies the criterion (or has cubic top), so is valid; independent verification also tests every finite assignment. Normalize long words to one fixed variable's cube by \(xyz=ttt\) using a nonempty prefix. If an original long word occurs, \(x\leq yyy\) recovers any target. Otherwise keep original edge directions; commutativity is false.

A new literal has an original square hh and an original incident literal z. Use \(x\leq xx+zt+z\) or \(x\leq xx+tz+z\) according to the incident edge direction. To reverse a known edge yx into xy, an original incoming edge zy supplies \(xy\leq yx+zy\). If the needed final y has no original incoming edge, but J holds, the two guarded reversal laws use one fixed original literal z and its original incident edge. Consequently under J every original edge can be reversed when required. This operation never assumes unguarded commutativity.

For a new pair ht whose underlying edge is original, use this reversal step if needed. If its tail t has an original loop, an original edge incident with h gives \(xy\leq yy+xz\) or \(xy\leq yy+zx\). If only its head h has an original loop, use \(xy\leq xx+zy\) and an original incoming edge at t; if none exists, the criterion gives J and reversal supplies that edge. If neither endpoint has a loop, one is an original literal, so J holds. Orient the incident witnesses to an outgoing edge xz at the literal endpoint x and an incoming edge ty at the other endpoint y. Apply \(xy\leq x+xz+ty\). Reverse the result with the same guarded rule when the literal endpoint is the target tail. All right-side terms used are original terms or have already been derived below P. New witnesses do not change the original sets in the semantic statement.

This derives every permitted absorption in finitely many applications of the eleven fixed laws. All substitutions are nonempty; there are no bounds on variables, word lengths, or polynomial sizes. Distributivity reduces general identities to summandwise absorptions in both directions. The eleven laws are therefore a complete finite basis. Transfer to other catalogue types requires independently checked complete HSP in both directions; short matching alone is insufficient.

\subsection{Five complete identities for semiring 14862}\label{case:335}
In the constant-free signature (+, multiplication), relative to the additively idempotent semiring axioms, a complete finite basis is

{\small\par\smallskip
\CaseLine{\(xy=yx; xyz=ttt; x\leq yyy; xy\leq xx+yz; xx\leq xy+zt+z.\)}
\smallskip}

Here \(u\leq P\) abbreviates \(u+P=P\). Addition is 00044/\allowbreak{}01034/\allowbreak{}00244/\allowbreak{}43434/\allowbreak{}44444; multiplication is 44444/\allowbreak{}42424/\allowbreak{}44444/\allowbreak{}42404/\allowbreak{}44444. The additive order has \(1<0<4\),\(2<0<4\),\(1<3<4\), and \(2+3=0+3=4\). A quadratic has value 0 at(3,\allowbreak{}3), value 2 at the other pairs in\{1,\allowbreak{}3\} squared, and 4 otherwise. Every word of length at least three is 4.

\par\smallskip\noindent\textit{Exact semantics for original polynomials.}\quad 
An original long word in P absorbs every target. Otherwise let L be the original literals, E the original unordered quadratic edges, V their endpoints, S their loop variables, and J the condition L intersects V. A literal h is absorbed exactly when h is in L. A pair ht is absorbed exactly when h, t are in V and either ht is an original unordered edge, an endpoint belongs to S, or J holds. No long target is absorbed. In particular, J recovers endpoint squares but never a new literal.

For necessity, assign all variables 1 to separate a long target from P. For a missing literal h, assign \(h=3\) and all other variables 1. Every right edge is 2 or 0, every right literal 1, so P is at most 0, which does not dominate 3. For a pair with an endpoint outside V, assign V to 1 and outsiders 2. P is at most 0 and the target is 4. If endpoints are in V, but no original edge or endpoint loop exists and J fails, assign the target endpoints 3, other V variables 1, and outsiders 2. Every right edge is 2 since there is no edge between the chosen endpoints or a loop there; every right literal is 2 because J fails. The target is 0, which is not below 2. These assignments exhaust failures.

For sufficiency consider P below 4 with an edge present. All endpoint values lie in\{1,\allowbreak{}3\}; all edge values are 2 or 0. If \(P=0\), all endpoint pairs are absorbed. If \(P=2\), every edge is 2 and every literal is 2; therefore J fails. An original edge has value 2. An original loop forces its endpoint to 1, so any pair with that endpoint is 2 as well. This proves the pair criterion. Original literals are automatically absorbed, and the necessity argument excludes all others. With no edges, only original literals are permitted. All arguments allow arbitrarily many variables and right summands.

\par\smallskip\noindent\textit{Complete finite derivations.}\quad 
Each of the five laws is valid by the table formula and criterion; all finite assignments to these laws are also checked independently. The structural laws identify all long words with a fixed variable's cube, using a nonempty prefix as a substituted variable. A right long word and \(x\leq yyy\) then recover any target. Otherwise orient original edges by commutativity. An original target needs no inference.

For a new pair with an original loop at h, choose an original edge tz at its other endpoint and apply \(xy\leq xx+yz\). If neither endpoint has an original loop, the criterion gives J. Choose an original edge hy at h, and an original literal z with original incident edge zt. Apply \(xx\leq xy+zt+z\) to obtain hh below P. If this is already the target square, stop; otherwise use \(xy\leq xx+yz\) with an original edge at the other endpoint. Original witnesses alone choose all steps. The recovered square is a derived term, not a new original witness. No new literal is introduced. Substitutions are nonempty.

Every valid absorption thus follows from these five fixed laws with no variable, length, or summand bounds. Distributive normalization and summandwise absorption in both directions prove all identities. Other types require verified complete HSP in both directions before transferring the finite-basis conclusion.

\subsection{A six-identity basis for a group clone below the identity}\label{case:336}
Work in the constant-free signature \((+,\cdot)\) relative to the axioms of additively idempotent semirings. Write \(u\le v\) for \(u+v\approx v\).

\par\smallskip\noindent\textit{Source semiring and finite basis.}\quad 
Let \(B=\{e, g, b, c, o\}\). Addition is the join in the poset generated by \(b<e<o\), \(b<c<o\), and \(g<c<o\), with no further comparabilities. Thus \(e+g=e+c=o\) and \(b+g=c\). Let \(G=\{e, g\}\) be the group of order two with identity \(e\). Define \(\rho(b)=e\), fix \(e, g, o\), and set \(\rho(c)=o\). Products are the products of the \(\rho\) images in the flat group-with-zero \(G\cup\{o\}\). The map \(\rho\) preserves joins and fixes its image, so this defines an associative, distributive semiring. The catalog map from the displayed order is \((3,\allowbreak{}2,\allowbreak{}1,\allowbreak{}0,\allowbreak{}4)\); the type is 14903.

Theorem. The following six laws form an identity basis for \(B\):

\CaseDisplay{\begin{array}{rll}
(C)&xy\approx yx,\\
(N)&x^3y\approx xy,\\
(S)&x^2y^2\approx x^2+y^2,\\
(U)&x\le x^3,\\
(O)&xyz\le x^3+y^3+z^3,\\
(Q)&x^3\le x+y^2.
\end{array}}

Minimality is not asserted.

Every word of length at least two evaluates to the group product of the \(\rho\) images if all letters lie in \(G\cup\{b\}\), and to \(o\) otherwise. This verifies (C),(N). Every square is \(e\) or \(o\), proving (S), while \(b\le e=b^3\) and \(c\le o=c^3\) prove (U). For (O), nongroup \(\rho\) values make the right side \(o\); otherwise either the three group values agree, when their product equals their join, or both \(e, g\) occur, when their join is \(o\). For (Q), the only nontrivial \(x\) are \(b, c\): \(b^3=e\) and \(b+y^2\ge e\) because \(y^2\in\{e, o\}\); \(c^3=o\) and \(c+e=c+o=o\).

\par\smallskip\noindent\textit{Stable terms and the absorption criterion.}\quad 
Call a word stable if its length is at least two. For a word \(w\), write \(c(w)\) for its content and \(\pi(w)\in\mathbb F_2^X\) for its exponent-parity vector. For a nonempty family \(R\) of stable words, let \(H(R)=\operatorname{aff}_{\mathbb F_2}\{\pi(r):r\in R\}\).

The five laws (C),(N),(S),(U),(O) give the following for arbitrary variable sets:

1. Stable words with the same content and parity are equal; each stable \(r\) satisfies \(r^3=r\).

2. If stable \(s, t\) have equal parity and \(c(s)\subseteq c(t)\), then \(s\le t\).

3. If \(R\ne\varnothing\), \(c(w)\subseteq\bigcup c(R)\), and \(\pi(w)\in H(R)\), then \(w\le\sum R\), including a linear \(w\).

Indeed, (N) inserts or deletes pairs of a variable inside a stable word, giving (1). Multiplying (S) by a stable word already containing one of its variables gives \(s\le sz^2\), proving (2) by adding missing variables. For (3), take an odd product of members of \(R\) having the required parity, add pairs of members until its content covers \(c(w)\), and use (O) and \(r^3=r\) to bound the product by \(\sum R\). Use (1),(2), and, when \(w\) is linear, \(w\le w^3\) from (U).

Write an arbitrary nonempty polynomial as

\CaseDisplay{F=\sum_{x\in L}x+\sum R,}

where \(L\) is its set of linear summands and \(R\) its stable summands; the variable sets may overlap. Put \(U=c(F)\) and \(N=c(R)\).

Proposition (exact criterion). If \(R\ne\varnothing\) and \(0\in H(R)\), then

\CaseDisplay{w\le F\quad\Longleftrightarrow\quad
c(w)\subseteq U\ \text{and}\ \pi(w)\in H(F).\quad\textup{(\ensuremath{1})}}

Otherwise \(w\le F\) holds exactly if either \(w=x\) with \(x\in L\), or \(c(w)\subseteq N\) and \(\pi(w)\in H(R)\).

The two conditions in (1) are necessary for every absorption. For content, assign a variable absent from \(F\) the value \(o\) and all other variables \(e\). For parity, restrict assignments to the two incomparable group values: if \(\pi(w)\) lies outside the affine hull of the right parities, a separating linear functional makes all right words one group value and \(w\) the other.

For sufficiency in the first case, take a nonempty word \(h\) whose content is \(N\). By item 3, \(h^2\le\sum R\le F\), since its parity is zero. For each right linear summand \(x\), (Q) gives \(x^3\le x+h^2\le F\), while (U) gives \(x\le x^3\). Replace all linear summands by these stable cubes. The resulting stable family has content \(U\) and affine hull \(H(F)\), so item 3 proves (1).

Suppose now \(R=\varnothing\) or \(0\notin H(R)\). Item 3 proves the second alternative, and the first is tautological. To prove their necessity, first take a nonempty \(R\) with \(0\notin H(R)\). There is a linear functional on the variables of \(R\) sending every \(\pi(r)\) to \(1\). For any parity vector \(v\notin H(R)\), this functional can also be chosen to send \(v\) to \(0\). Initially assign each variable the corresponding group value \(e\) or \(g\). Change any right linear variable assigned \(e\) but occurring in \(R\) to \(b\). This leaves all stable values unchanged, because \(\rho(b)=e\), while making every right linear summand lie below \(c\): such values are now \(b\) or \(g\). Assign variables occurring only linearly the value \(c\). Consequently all right stable words equal \(g\), all right linear words lie below \(c\), and \(F\le c\).

If \(c(w)\subseteq N\) but \(\pi(w)\notin H(R)\), the chosen functional makes a stable target \(w\) equal \(e\). It also makes a linear target \(x\) equal \(e\) when \(x\) is not a right linear summand; only right linear variables were lowered to \(b\). Since \(e\nleq c\), absorption fails. If a stable target uses a variable occurring only linearly, its value is \(o\) after that variable is assigned \(c\), and again \(o\nleq c\). If it uses a variable absent from \(F\), assign that variable \(o\). A linear target absent from \(L\) and \(N\) is separated in the same way. These are all possible failures of the two alternatives. When \(R\) is empty, assign every right linear variable \(c\): every stable target using one of them evaluates to \(o\), while \(F=c\); an absent variable is handled by assigning it \(o\).

This proves the proposition. Every polynomial identity amounts to absorption of each summand in both directions. The criterion gives derivations of every semantically valid absorption using the six displayed laws, for unrestricted numbers of variables and summands. Hence the six laws form a basis.

\subsection{Six complete identities for semiring 14959}\label{case:337}
In the constant-free signature (+, multiplication), a basis relative to the additively idempotent semiring axioms is

{\small\par\smallskip
\CaseLine{\(xy=yx; xyz=ttt; x\leq yyy; x\leq y+xy;\)}
\CaseLine{\(xz\leq xx+xy+yz; xy\leq xx+x+y+yz.\)}
\smallskip}

The notation \(u\leq P\) means \(u+P=P\). Addition is 00044/\allowbreak{}01144/\allowbreak{}01244/\allowbreak{}44434/\allowbreak{}44444; multiplication is 44444/\allowbreak{}44044/\allowbreak{}40044/\allowbreak{}44404/\allowbreak{}44444. The additive order \(is 2<1<0<4\), with 3 incomparable with 0,\allowbreak{}1,\allowbreak{}2 below 4. A quadratic is 0 exactly at (1,\allowbreak{}2),(2,\allowbreak{}1),(2,\allowbreak{}2),(3,\allowbreak{}3), and 4 otherwise. Every longer word is 4.

\par\smallskip\noindent\textit{Exact semantics using original components.}\quad 
An original long word in P absorbs all targets. Otherwise let L be the original literal variables, E the original unordered edges, V their endpoints, and S their loop variables. Components refer to(V, E), and a component is active exactly when it contains an original literal. A literal h is absorbed iff h is original or lies in an active component. A pair ht is absorbed iff h, t lie in V, either their edge is original or an endpoint is an original loop, and either their components coincide or both components are active. No long word is absorbed. Connectivity does not give a new pair without an original edge or an endpoint loop. A loop alone does not recover its literal.

For sufficiency, if P is below 4 with an edge present, all edges are 0 and \(P=0\). A connected component lies entirely within\{1,\allowbreak{}2\} or\{3\}. An active component must lie in\{1,\allowbreak{}2\}, since a literal 3 joined with an edge0 gives 4. Its literals are therefore below 0. A loop endpoint has value 2 in the first block or 3 in the second. Thus a pair with a loop endpoint is 0 if its endpoints share a component, or if both components are active. Original edges are 0 automatically. With no edges, only original literals meet the criterion.

For necessity, all variables 2 separate a long target. For a missing literal outside V, assign it 4 and all others 2; it is not in L, so no right summand contains it. For a missing literal in an inactive component, assign its component 3 and all other variables 2. Right edges are 0, all literals 2, while the target literal 3 is not below 0. For a pair with an endpoint outside V, assign V to 2 and outsiders 0; all right terms are 0 or 2 and the target 4. If the target edge is missing and neither endpoint has a loop, assign both endpoints 1 and all others 2. No right edge has two1 endpoints, so \(P\leq 0\) whereas the target is 4. Finally, for distinct target components not both active, assign an inactive target component 3 and everything else 2. No literal in that component occurs. Again \(P=0\) and the target 4. These cases exhaust failure of the criterion for arbitrary polynomials.

\par\smallskip\noindent\textit{Complete finite derivations.}\quad 
All six fixed laws are valid by the table formula and criterion; independent verification also checks every law assignment. An original long word supplies a fixed cube via \(xyz=ttt\) using a nonempty prefix, and \(x\leq yyy\) supplies any target. Otherwise a new literal is obtained along a simple original path starting at an original literal, using \(x\leq y+xy\) at each edge. No original witness sets are redefined.

For a nonoriginal target pair, orient its original loop endpoint as h. If h, t are connected, choose a simple original path and repeatedly apply \(xz\leq xx+xy+yz\), retaining the original loop hh. If they belong to different active components, first derive their literals along the two original paths. Choose an original edge tz incident with t and apply \(xy\leq xx+x+y+yz\). Commutativity restores the requested order. An original target requires no inference; all needed substitutions are nonempty and use already established terms below P.

The argument derives every valid absorption without any bound on variables, word lengths, or summand counts. Distributivity and summandwise absorption in both directions establish all identities. Thus the six laws are a complete finite basis. Transfers to other catalogue types require independently checked complete HSP both ways.

\subsection{Thirteen complete identities for 14967}\label{case:338}
In the signature of addition and multiplication without constants, write \(u\leq P\) for \(u+P=P\). Relative to the additively idempotent semiring axioms, a complete basis is:

{\small\par\smallskip
\CaseLine{\(xy=yx; xxy=xyy; xyzt=xxxx; xxxx=yyyy; x\leq yyyy;\)}
\CaseLine{\(xy\leq xyz; xx\leq xyz; x\leq xyz; y\leq x+xy;\)}
\CaseLine{\(xz\leq xx+xy+yz; xz\leq x+xx+z+zt;\)}
\CaseLine{\(xxx\leq x+xx+yzt; xyz\leq xxx+yyy+zzz.\)}
\smallskip}

\par\smallskip\noindent\textit{Values and original components.}\quad 
Addition:00044/\allowbreak{}01144/\allowbreak{}01244/\allowbreak{}44434/\allowbreak{}44444. Multiplication:44444/\allowbreak{}44044/\allowbreak{}40144/\allowbreak{}44404/\allowbreak{}44444. The additive chain \(is 2<1<0<4\);3 is incomparable with 0,\allowbreak{}1,\allowbreak{}2 and below 4. A quadratic has value 1 at(2,\allowbreak{}2), value 0 at(1,\allowbreak{}2),(2,\allowbreak{}1),(3,\allowbreak{}3), and 4 otherwise. A cubic has value 0 exactly when all occurrences are 2, and 4 otherwise. Every word of \(length\geq 4\) has value 4. This follows from the quadratic values:1 times 2 is 0,\allowbreak{}0 times any element is 4, and the other possible products immediately give 4. Cubic functions depend only on content, with repetitions normalized by \(xxy=xyy\) and commutativity. The two fourth-power identities reduce all longer words to one common fourth power; this is a word function, not a signature constant.

For the original polynomial P let L be its literal variables, V its nonlinear content, E its quadratic pairs, Q its cubics, and T their content union. Join variables occurring together in an original nonlinear word; let the resulting connected components partition V. Singleton components are allowed. Put

H0 = T union the variables of original squares; A  = union of components meeting L or T; H  = H0 intersect A.

Call components in A active. An original word of \(length\geq 4\) absorbs every target. If there is no such original word, the exact criterion is:

\textbullet{} A literal x is absorbed iff x belongs to L union A.

\textbullet{} A quadratic xy is absorbed iff x, y belong to V and either their pair is an original quadratic, or they lie in the same component with at least one endpoint in H0, or both are active with an endpoint in H.

\textbullet{} A cubic is absorbed iff Q is nonempty and its content is contained in H.

\textbullet{} No word of \(length\geq 4\) is absorbed.

The cubic guard Q is global. In contrast, whether a square forces its variable to 2 depends on activation of its own original component. Activation and H are defined from original witnesses, never by recomputing them after derived summands have been added.

\par\smallskip\noindent\textit{Proof of the criterion at arbitrary variable count.}\quad 
Suppose an evaluation gives P a value different from 4 and P has a nonlinear term. Each component of V has one of two forms: all its variables are 3; or all are in\{1,\allowbreak{}2\}, every quadratic has a2 endpoint, and every cubic variable is 2. This follows directly by propagating the quadratic table along original edges; a cubic already forces all its letters to 2. A component meeting an original literal cannot have value 3 at its variables, because that literal added to any nonlinear value 0 or 1 gives 4. Thus active components have only 1,\allowbreak{}2; variables in H have value 2, while inactive original-square variables may be 2 or 3.

All nonlinear right values are 0 or 1. Active literals are at most 1, so are absorbed. For a same-component pair anchored at H0, the all 3 mode gives 0 and its original incidence supplies 0. In the\{1,\allowbreak{}2\} mode the anchored endpoint is 2. A target value 0 implies that the other endpoint is 1; an original edge incident with it must then have value 0. Otherwise the target is 1. The same reasoning works across active components with an endpoint in H. A cubic on H has value 0; an original cubic supplies 0. Original target pairs and original literals need no argument.

Here are separators for every failed criterion. With no long original word, assigning all variables 2 keeps P below 4; this separates all long targets. A missing literal outside L union V can be set 4. For a missing literal in an inactive component, set that whole component to 3 and all others 2. The right value is 0 and the target is 3. For a nonlinear target containing an outside-V variable, set that variable 0 and all others 2; the target is 4 while every original term is at most 0.

For a failed quadratic with endpoints in different components, if either component is inactive, set that whole component 3 and all others 2. The target is 4 and the right is 0. If both components are active, the failed criterion means neither endpoint belongs to H0. Set both endpoints 1 and all other variables 2. They have no original squares or cubic occurrences and do not occur together in an original edge. Every right term remains at most 0, but the target is 4. The same assignment treats failed pairs within a single component: neither endpoint is anchored and their original edge is absent. It also treats a missing square.

For a failed cubic with no original cubic, all variables 2 give right value at most 1 and target 0. If a cubic exists but a target variable is inactive, set its whole component 3 and the rest2; the target becomes 4 while the right is 0. If instead that variable is active but outside H, it has no original square or cubic occurrence; set it 1 and the others 2. Again the right is 0 and the target 4. Literal-only polynomials have exactly their ordinary content absorptions; the outside-variable assignments cover their failed cases. No finite variable bound was used.

\par\smallskip\noindent\textit{Complete derivations.}\quad 
Retain all original right terms. Long originals give a fourth power; \(x\leq yyyy\) then recovers any target. An original cubic gives its pairs, squares and literals using \(xy\leq xyz\), \(xx\leq xyz\), \(x\leq xyz\). These steps supply every needed edge of the component graph; all are obtained from their original witnesses. Original squares stay available.

Start each active component at an original literal or a literal from an original cubic. Along an original edge xy, use \(y\leq x+xy\) to obtain the next literal. Finite paths in the original graph recover every variable in A, without changing the definition of A.

For a target pair in one component, keep the square xx of its anchored endpoint. Follow a simple original path from x to z. The first edge is available, and \(xz\leq xx+xy+yz\) extends the currently recovered pair xy by the next original edge yz. A target square needs no path; its original square or cubic witness already supplies it. Every traversed path edge is positive and comes from a nonempty original word.

For endpoints in different active components, choose x in H. Its square and literal are now available, as is the target's other literal z. Retain an original nonlinear incidence zt, where t may be a nonempty word. The law \(xz\leq x+xx+z+zt\) gives the target pair.

For a cubic target, each target variable x has its square and its recovered literal. Keep one original cubic yzt anywhere in P. Then \(xxx\leq x+xx+yzt\) gives each required pure cube; \(xyz\leq xxx+yyy+zzz\) combines them. Repetition of target letters is permitted. This explains why an inactive square or a literal without an original cubic cannot replace the guarded cube construction.

Every displayed law is valid by the exact criterion and word formula. All steps are applications of these fixed thirteen identities with nonempty substitutions, semiring axioms and additive transitivity. Hence every valid absorption is derivable. Expanding arbitrary terms and using absorption in both directions establishes unrestricted completeness.

Boundary failures: \(y\leq xy\) at \(x=y=3\); \(xxx\leq xx\) at \(x=2\); \(xxx\leq xx+yzt\) at \(x=3\), \(y=z=t=2\); \(xz\leq xx+zt\) at \(x=3\), \(z=t=2\); \(xy\leq xz+yz\) at \(x=y=1\), \(z=2\); \(xxxx\leq x+xxx\) at \(x=2\). These separate literal activation, local square propagation, cross-component compatibility and the global guard.

\subsection{Twelve complete identities for 14972}\label{case:339}
The signature has addition and multiplication, without constants. Write \(u\leq P\) for \(u+P=P\). In addition to the additively idempotent semiring axioms, the following twelve identities form a complete finite basis:

{\small\par\smallskip
\CaseLine{\(xy=yx; xxy=xyy; xyzt=xxxx; xxxx=yyyy;\)}
\CaseLine{\(x\leq yyyy; x\leq xyz; z\leq x+xy+zt; xy\leq xyz; xx\leq xyz;\)}
\CaseLine{\(xz\leq xx+zy; xxx\leq xx+yzt; xyz\leq xxx+yyy+zzz.\)}
\smallskip}

Every substitution below is nonempty. There is no restriction on the number of variables, word lengths, or the number of polynomial summands.

\par\smallskip\noindent\textit{Exact word values.}\quad 
Addition:00004/\allowbreak{}01114/\allowbreak{}01214/\allowbreak{}01134/\allowbreak{}44444. Multiplication:44444/\allowbreak{}44044/\allowbreak{}40304/\allowbreak{}44044/\allowbreak{}44444. The additive order is 2,\(3<1<0<4\), with 2 and 3 incomparable. A quadratic has value 3 at (2,\allowbreak{}2), value 0 at (2,\allowbreak{}1),(1,\allowbreak{}2),(2,\allowbreak{}3),(3,\allowbreak{}2), and value 4 otherwise. A cubic has value 0 exactly when all three occurrences have value 2; otherwise it is 4. Every word of length at least four has value 4. This follows by multiplying each possible quadratic value by its next letter:3 times 2 is 0, and 0 times any element is 4.

Commutativity sorts quadratics. Cubic functions depend only on content; \(xxy=xyy\) chooses either multiplicity for two-variable cubics. A word of \(length\geq 4\) reduces to a fourth power by \(xyzt=xxxx\), using its nonempty remainder for t. All fourth powers agree by \(xxxx=yyyy\). A chosen fourth power is a word-function representative, not a constant in the signature.

\par\smallskip\noindent\textit{Original-witness criterion.}\quad 
Take the original right polynomial P. Let L be its literal variables, V the union of contents of its nonlinear words, E its quadratic pairs, Q its cubic words, and

H = union of cubic contents, together with variables of original squares.

If P contains a word of \(length\geq 4\), it absorbs every word. Otherwise:

\textbullet{} A literal x is absorbed exactly if x is in L, or x is in V and either Q is nonempty or L intersects V.

\textbullet{} A quadratic xy is absorbed exactly if x, y are in V and either xy is an original quadratic or at least one endpoint belongs to H.

\textbullet{} A cubic is absorbed exactly if Q is nonempty and its content is in H.

\textbullet{} No word of \(length\geq 4\) is absorbed.

Literal-edge incidence activates literal recovery globally. It does not anchor the other endpoint: that was a different property of 13363. The original cubic guard can lie outside the target component. Derived literals must not be used to redefine the original criterion.

For sufficiency, consider any evaluation at which P is not 4. Every nonlinear variable is in\{1,\allowbreak{}2,\allowbreak{}3\}; every quadratic has a2 endpoint, and every variable of an original square or cubic is 2. If an original cubic exists, the right value is 0. If no original cubic exists, every nonlinear term is 3 or 0. When any is 0, the right is 0. Otherwise all nonlinear variables are 2 and all nonlinear terms are 3. An original literal in V then raises the right value to at least 1, which absorbs 2. This proves literal recovery and the cubic rule. For an anchored target pair, its value is 3 or 0. Value 3 is absorbed by any nonlinear right term. If its value is 0, its other endpoint has value 1 or 3; its original nonlinear incidence must then be a quadratic with value 0. Thus it is absorbed too.

The converse has explicit separating assignments. In the absence of a long right word, setting every variable to 2 gives a right value below 4, so it separates every long target. A missing literal outside L union V can instead be set 4. A literal in V without the recovery condition is separated by \(V=2\) and all external literal \(variables=3\): the right is 3 and the target is 2. There are no cubics or literals inside V in this case.

For a failed nonlinear target with a variable outside V, set that variable to 0, every V variable to 2, and all other variables to 2. The target is 4; all right terms remain at most 0, even if the chosen variable is a literal. For a failed quadratic with both endpoints in V, neither is in H and their pair is absent. Set those endpoints to 1 and all others to 2. No original loop, cubic or edge has both exceptional occurrences, so every right term is at most 0 while the target is 4. This includes a missing square. For a failed cubic without an original cubic, setting all variables to 2 gives right value at most 1 and target 0. If an original cubic exists, choose an unanchored target variable and set it 1, keeping all others 2. Its original incidences are unlooped quadratic edges; hence the right remains at most 0 and the target is 4. These cases exhaust all failed absorptions, for arbitrarily many variables.

\par\smallskip\noindent\textit{Complete finite derivations.}\quad 
Keep every original summand available. An original long word produces a fourth power by normalization; \(x\leq yyyy\) recovers any target from it.

For an allowed missing literal z, either there is an original literal x with an original nonlinear incidence xy, or an original cubic supplies such a pair: recover one of its letters x by \(x\leq xyz\) and split that same cubic after x. In either case, take an original nonlinear incidence zt and apply \(z\leq x+xy+zt\). This permits disconnected witnesses and nonempty word images for y, t; no path-length argument is needed.

For each required anchored variable x, its original square is available, or \(xx\leq xyz\) extracts it from an original cubic containing x. To recover a missing target pair xz with x anchored, choose an original nonlinear incidence at z. If it is cubic, \(xy\leq xyz\) extracts a pair zy. Then \(xz\leq xx+zy\) produces the target. Original target pairs require no step.

For a cubic target, extract each required square as just described. Keep one original cubic yzt, even if it is disjoint from the target. The law \(xxx\leq xx+yzt\) produces each required pure cube, and \(xyz\leq xxx+yyy+zzz\) combines them. Repeated target letters cause no problem.

The listed identities are valid by the exact criterion and word values. The constructions derive every valid absorption. Distributive expansion of terms into nonempty sums of nonempty words, followed by absorption in both directions, proves completeness for all identities.

Boundary counterexamples: \(yy\leq x+xy\) fails at \(x=2\), \(y=1\); \(xxx\leq xx\) fails at \(x=2\); \(xy\leq xz+yz\) fails at \(x=y=1\), \(z=2\); \(xy\leq x+y\) fails at \(x=y=2\); \(xxxx\leq x+xxx\) fails at \(x=2\). A literal witness cannot replace the original cubic guard, and literal incidence alone cannot supply an anchored square.

\subsection{Six complete identities for semiring 15483}\label{case:340}
In the constant-free signature (+, multiplication), a finite basis relative to the additively idempotent semiring axioms is

{\small\par\smallskip
\CaseLine{\(xy=yx; xyz=ttt; x\leq yyy;\)}
\CaseLine{\(xy\leq xx+yz; xy\leq x+xz+yt; x\leq xy+zt+z.\)}
\smallskip}

Addition:00004/\allowbreak{}01004/\allowbreak{}00224/\allowbreak{}00234/\allowbreak{}44444. Multiplication:44444/\allowbreak{}44444/\allowbreak{}44014/\allowbreak{}44114/\allowbreak{}44444. The additive order \(is 3<2<0<4\) and \(1<0<4\). A quadratic is 0 at(2,\allowbreak{}2), 1 at(2,\allowbreak{}3),(3,\allowbreak{}2),(3,\allowbreak{}3), and 4 otherwise. Every longer word is 4.

\par\smallskip\noindent\textit{Exact semantics.}\quad 
An original word of length at least three in P absorbs every target. Otherwise let L be original literal variables, E original unordered quadratic edges, V their endpoints, S original loops, and J mean L intersects V. A literal h is absorbed iff h belongs to L or(h belongs to V and J). A pair ht is absorbed iff h, t belong to V and either their edge is original, an endpoint belongs to S, or J holds. A longer target is never absorbed. In particular, a distant original literal-edge overlap permits every pair of original edge endpoints; it need not occur at a target endpoint.

Necessity has explicit separators. All variables 2 separate a long target from P. For a pair with an endpoint outside V, assign V to 3 and others 1:all right terms are 1 or 3, so their sum is at most 0, but the target is 4. If endpoints are in V but there is no original edge or endpoint loop and J fails, assign target endpoints 2, other V variables 3, and outside V variables 1. Every edge and literal is 1, whereas the target is 0. For a missing literal h outside V, assign \(h=4\) and all other variables 3. For h in V with J false, assign V to 3 and outsiders 1:the right side is 1 but \(h=3\) is not below 1.

For sufficiency consider P below 4. Edge endpoints are all 2 or 3. If \(P=0\), all endpoint literals and all endpoint pairs are absorbed. If \(P=1\), all edges are 1 and no incident literal can occur, so J is false. An original edge has value 1. An endpoint loop forces that endpoint to 3, which also forces any target incident pair to 1. Thus every permitted pair is absorbed. A new permitted literal uses J, which forces \(P=0\). If there are no edges, only original literals are permitted. This proves the criterion for arbitrary P and target.

\par\smallskip\noindent\textit{Unrestricted complete derivations.}\quad 
Each displayed law is valid by the criterion and word-value formula; all its assignments are independently checked as well. Normalize long words to a fixed cube using \(xyz=ttt\) with a nonempty prefix. If P has a long word, \(x\leq yyy\) supplies any target. Otherwise orient original edges using commutativity. An original target needs no step.

For a new literal h, choose an original edge hy incident with h, and an original literal z with original incident edge zt. Apply \(x\leq xy+zt+z\). For a new pair with an original loop at h, use \(xy\leq xx+yz\) and an original edge incident with its other endpoint. If no endpoint loop exists, J must hold. First recover the literal h by the preceding rule, unless it is already an original literal. Choose original edges hz and ta at both target endpoints and use \(xy\leq x+xz+yt\). These choices give the target below P in finitely many steps and require no new interpretation of the original sets.

Substitutions are always nonempty. The argument has no bounds on variables, lengths, or summand counts. Distributive normalization and summandwise absorption in both directions prove all valid identities, establishing completeness of the six fixed laws.

\subsection{Five complete identities for semiring 15618 and an equal-variety transfer}\label{case:341}
Signature (+, multiplication), no constants. Relative to the additively idempotent semiring axioms the complete finite basis is

{\small\par\smallskip
\CaseLine{\(xy=yx; xyz=ttt; x\leq yyy; xy\leq xx+yz; x\leq xy+zt+z.\)}
\smallskip}

Addition:00004/\allowbreak{}01004/\allowbreak{}00224/\allowbreak{}00234/\allowbreak{}44444. Multiplication:44444/\allowbreak{}44444/\allowbreak{}44414/\allowbreak{}44114/\allowbreak{}44444. The additive order \(is 3<2<0<4\) and \(1<0<4\), with 1 incomparable with 2,\allowbreak{}3. A quadratic has value 1 exactly at the pairs(2,\allowbreak{}3),(3,\allowbreak{}2),(3,\allowbreak{}3); it has value 4 otherwise. Every word of \(length\geq 3\) has value 4, because both possible quadratic values 1,\allowbreak{}4 multiply to 4. Thus multiplication is commutative, and all long words are equal, even when their contents are disjoint. This immediately proves the first three laws.

For \(xy\leq xx+yz\), if the right side is not 4, \(xx=1\) forces \(x=3\), and \(yz=1\) forces y in\{2,\allowbreak{}3\}; then \(xy=1\). For \(x\leq xy+zt+z\), if both products are 1, then x, z belong to\{2,\allowbreak{}3\}, and \(1+z=0\) bounds x. If either product is 4, the inequality is immediate. These arguments prove all laws on arbitrary assignments; a separate exhaustive finite check is included.

\par\smallskip\noindent\textit{Exact original-polynomial criterion.}\quad 
If P contains any original long word, it absorbs every target. Otherwise let L be its original literal variables, E its original UNORDERED quadratic pairs, including loops, V the endpoints of E, and S its loop variables. Then:

\textbullet{} No long target is absorbed.

\textbullet{} A literal h is absorbed iff h belongs to L, or h belongs to V and L intersects V.

\textbullet{} A quadratic ht is absorbed iff h, t both belong to V and either the unordered pair ht is in E or at least one of h, t belongs to S.

All these sets refer only to original terms. The signature has no nullary top constant:ttt is an ordinary term identity.

\par\smallskip\noindent\textit{Necessity by explicit evaluations.}\quad 
Assume P has no long word. Assign every variable 3. Quadratics are 1, literals 3, and the right sum is 1,\allowbreak{}3, or 0, never 4; every long target is 4. For a missing literal h outside V (and L), assign \(h=4\) and others 3. No original word contains h, so the right sum is below 4. If h is in V but L and V are disjoint, assign V variables 3 and other variables 1. Every original right term is 1, while \(h=3\) is not below 1.

For a quadratic with an endpoint outside V, assign all variables of V the value 3 and all others 1. Original edges have value 1 and literals have value 1 or 3, so the right sum \(is\leq 0\). The target contains 1 and therefore has value 4. This remains valid even when the outside endpoint is an original literal.

Finally suppose both target endpoints lie in V, but their unordered pair is not original and neither endpoint has an original loop. Assign those endpoints 2 and every other variable 3. An original edge can never have both endpoints 2:that would be the missing target edge or a forbidden endpoint loop. Hence every original edge is 1. Original literals are 2 or 3, so \(P\leq 0\); the target \(is 2*2=4\). These cases prove necessity for arbitrarily many original words and variables.

\par\smallskip\noindent\textit{Complete derivations.}\quad 
An original long word normalizes to ttt using \(xyz=ttt\), with a nonempty prefix for x and its final two letters for y, z. The variable t can be any fixed existing variable. Substitute any target word for x in \(x\leq yyy\) and that fixed variable for y to derive the target.

Without long words, normalize every quadratic by commutativity. An original target needs no inference. For a new literal h, the criterion provides an original edge ha containing h and an original literal z that occurs in an original edge zb. After orienting those edges by commutativity, \(x\leq xy+zt+z\) gives \(h\leq ha+zb+z\leq P\). The two chosen edges may coincide; idempotent addition handles this case.

For a new quadratic ht, one endpoint, say h, has an original square hh, and the other t has an original incident edge tz. The law \(xy\leq xx+yz\) gives \(ht\leq hh+tz\leq P\). If only t has a square, interchange the target endpoints first. This is all the criterion permits. Each substitution is nonempty; all steps are identities under the original right sum.

Distributivity expands arbitrary terms into nonempty polynomials. Apply these absorption derivations to every summand in both directions of any valid identity. This proves unrestricted completeness of the five laws. Proposition~\ref{prop:lower-fibres} gives equality of the two varieties in case~\ref{grp:599}, so \(\A{15596}\) has the same basis.

\subsection{Eleven complete identities for semiring 15603}\label{case:342}
Signature (+, multiplication), no constants. A complete finite basis relative to the additively idempotent semiring axioms is

{\small\par\smallskip
\CaseLine{\(xyz=ttt; x\leq yyy;\)}
\CaseLine{\(xy\leq xx+yz; xy\leq xx+zy;\)}
\CaseLine{\(xy\leq yx+xz; xy\leq yx+zt+z;\)}
\CaseLine{\(x\leq xy+zt+z; x\leq yx+zt+z;\)}
\CaseLine{\(xy\leq yx+tz+z; x\leq xy+tz+z; x\leq yx+tz+z.\)}
\smallskip}

Addition rows:00004/\allowbreak{}01004/\allowbreak{}00224/\allowbreak{}00234/\allowbreak{}44444. Multiplication rows:44444/\allowbreak{}44444/\allowbreak{}44404/\allowbreak{}44114/\allowbreak{}44444. The additive order \(is 3<2<0<4\) and \(1<0<4\), with 1 incomparable with 2,\allowbreak{}3. A quadratic is 0 at(2,\allowbreak{}3),1 at(3,\allowbreak{}2) or(3,\allowbreak{}3), and 4 otherwise. Every word of \(length\geq 3\) is 4 because every quadratic result0,\allowbreak{}1,\allowbreak{}4 has constant4 multiplication row. This proves the first two laws.

For the square-guard laws, a right side below 4 forces \(x=3\) and the other endpoint y in\{2,\allowbreak{}3\}; the target is 1. For \(xy\leq yx+xz\), when the right is below 4 all involved edge endpoints lie in\{2,\allowbreak{}3\}; if target head \(x=2\) the edge xz has value 0 and bounds xy; otherwise \(xy=1\). For the six literal-guard laws, when both edges avoid 4, z lies in\{2,\allowbreak{}3\}, so the sum zt+z is 0 regardless of which edge value occurs. The other edge ensures x (and for the reversal law y) belongs to\{2,\allowbreak{}3\}; the literal or target edge is bounded by 0. Thus all eleven identities hold, also independently checked exhaustively on their variables.

\par\smallskip\noindent\textit{Original-polynomial criterion.}\quad 
An original long word absorbs every target. Otherwise let L be the original literal variables, E the original DIRECTED pairs, V their endpoints, H their heads, S their loop variables, and J the statement L intersects V. Then no long target is absorbed. A literal h is absorbed iff h is in L, or h is in V and J holds. A quadratic ht is absorbed iff all three conditions hold:

1. h, t belong to V;

2. ht or th is an original edge, or h or t is an original loop variable;

3. h belongs to H, or J holds.

This third condition is essential:\(xy=yx\) fails at \(x=2\), \(y=3\). Short matching with a commutative algebra cannot supply this missing guard.

\par\smallskip\noindent\textit{Necessity via separating evaluations.}\quad 
With all variables 3 and no original long word, original edges are 1, literals 3, and P is 1,\allowbreak{}3, or 0. Every long target is 4 and fails. For a missing literal h outside V and L, set \(h=4\) and all others 3. For a literal h in V but J false, set V variables 3 and all other variables 1; all right terms are 1 while \(h=3\) is not below 1.

For a quadratic with an endpoint outside V, set V variables 3 and all others 1. Original edges are 1 and literals are 1 or 3, so \(P\leq 0\); the target containing 1 is 4. For failure of condition 2, assign both target endpoints 2 and all others 3. No original edge can have both ends2, since that would be a missing target edge or a forbidden target endpoint loop. Every edge is therefore 0 or 1, and literals 2 or 3, so \(P\leq 0\) whereas the target is 4.

If conditions1 and 2 hold but 3 fails, h is not an original head and J is false. The target cannot be a square:condition 2 would require an original loop at h, contrary to h not in H. Assign \(h=2\), all other V variables 3, and all variables outside V the value 1. Every original edge has head 3 and final 2 or 3, hence value 1. Since J is false, all original literals are outside V and also 1. The target ht \(is 2*3=0\), not below 1. These cases exhaust failure and prove necessity without bounds on variables or numbers of original summands.

\par\smallskip\noindent\textit{Complete finite derivations.}\quad 
Normalize every original long word to a fixed variable's cube by \(xyz=ttt\), with a nonempty prefix for x. Then \(x\leq yyy\) with the target word as the image of x derives any target. No nullary top is added.

Without long words, original targets need no inference. For a new literal h, use an original incident edge ha or ah, and an original literal z lying in an original edge zb or bz. Four orientation cases are handled separately by \(x\leq xy+zt+z\), \(x\leq yx+zt+z\), \(x\leq xy+tz+z\), \(x\leq yx+tz+z\), with \(x=h\), \(y=a\), \(z=z\), \(t=b\). No reversal is assumed.

For a missing target ht, first obtain at least one orientation. If th is original, use it. Otherwise condition 2 supplies a loop at h or t. If hh is original, choose an incident original edge ta or at and use \(xy\leq xx+yz\) or \(xy\leq xx+zy\) to derive ht. If only tt is supplied, the same rules derive th using an original edge incident with h.

To reverse an available th into ht, condition 3 provides either an original edge ha (when h is in H), or an original literal z with edge zb or bz (when J holds). Use respectively \(xy\leq yx+xz\) with \(x=h\), \(y=t\), \(z=a\); \(xy\leq yx+zt+z\) with \(x=h\), \(y=t\), \(z=z\), \(t=b\); or \(xy\leq yx+tz+z\) with those same substitutions. The literal guard remains an original term, and every required edge is either original or the single orientation just constructed. All substitutions are nonempty.

These rules derive precisely the targets of the exact criterion. Distributive expansion and summandwise absorption in both directions therefore derive every valid identity, proving unrestricted completeness of the eleven laws.

\endgroup

\end{document}